\documentclass[12pt,openany]{book} 
\usepackage[english]{babel}
\usepackage[cp1251]{inputenc}
\usepackage{amsmath}
\usepackage{amssymb}
\usepackage{amsfonts}

\usepackage[linktocpage=true, colorlinks=true, linkcolor=blue, citecolor=blue, urlcolor=blue]{hyperref}

\begin{document}

\pagestyle{myheadings}
\markboth{\underline{D.F. Kuznetsov {\footnotesize{S}}{\tiny{trong}}
{\footnotesize{A}}{\tiny{pproximation of}} 
{\footnotesize{I}}{\tiny{terated}} {\footnotesize{I}}{\tiny{to and}} 
{\footnotesize{S}}{\tiny{tratonovich}} 
{\footnotesize{S}}{\tiny{tochastic}} 
{\footnotesize{I}}{\tiny{ntegrals:}} 
{\footnotesize{M}}{\tiny{ethod of}} 
{\footnotesize{G}}{\tiny{eneralized}} 
{\footnotesize{M}}\tiny{ultiple}
{\footnotesize{F}}{\tiny{ourier}} 
{\footnotesize{S}}{\tiny{eries}}}}
{\underline{D.F. Kuznetsov {\footnotesize{S}}{\tiny{trong}}
{\footnotesize{A}}{\tiny{pproximation of}} 
{\footnotesize{I}}{\tiny{terated}} {\footnotesize{I}}{\tiny{to and}} 
{\footnotesize{S}}{\tiny{tratonovich}} 
{\footnotesize{S}}{\tiny{tochastic}} 
{\footnotesize{I}}{\tiny{ntegrals:}} 
{\footnotesize{M}}{\tiny{ethod of}} 
{\footnotesize{G}}{\tiny{eneralized}} 
{\footnotesize{M}}\tiny{ultiple}
{\footnotesize{F}}{\tiny{ourier}} 
{\footnotesize{S}}{\tiny{eries}}}} 
\markright{\underline{D.F. Kuznetsov  {\footnotesize{S}}{\tiny{trong}}
{\footnotesize{A}}{\tiny{pproximation of}} 
{\footnotesize{I}}{\tiny{terated}} {\footnotesize{I}}{\tiny{to and}} 
{\footnotesize{S}}{\tiny{tratonovich}} 
{\footnotesize{S}}{\tiny{tochastic}} 
{\footnotesize{I}}{\tiny{ntegrals:}} 
{\footnotesize{M}}{\tiny{ethod of}} 
{\footnotesize{G}}{\tiny{eneralized}} 
{\footnotesize{M}}\tiny{ultiple}
{\footnotesize{F}}{\tiny{ourier}} 
{\footnotesize{S}}{\tiny{eries}}}}

\thispagestyle{empty}
$$
{}
$$
$$
{}
$$
~\\
{{\huge{\bf Strong Approximation of}}}\\ 
~\\
~\\
{{\huge{\bf Iterated It\^{o} and Stratonovich}}}\\ 
~\\
~\\
{{\huge{\bf Stochastic Integrals:}}}\\ 
~\\
~\\
{{\huge{\bf Method of Generalized Multiple}}}\\ 
~\\
~\\ 
{{\huge{\bf Fourier Series.}}}\\
~\\
~\\
{{\huge{\bf Application to Numerical Solution}}}\\
~\\
~\\
{{\huge{\bf of It\^{o} SDEs}}}\\
~\\
~\\
{{\huge{\bf and Semilinear SPDEs}}}\\
~\\
$$
{}
$$
$$
{}
$$
\centerline{\Large {\bf Dmitriy F. Kuznetsov}}
$$
{}
$$
\begin{center}
{\large Russia, 195251, Saint-Petersburg, Polytechnicheskaya st., 29\\
Peter the Great Saint-Petersburg Polytechnic University\\
Institute of Physics and Mathematics\\
Department of Higher Mathematics\\
~\\
e-mail: sde\_kuznetsov@inbox.ru}
\end{center}

\newpage

\noindent
$$
{}
$$
$$
{}
$$
$$
{}
$$
$$
{}
$$
$$
{}
$$
$$
{}
$$
$$
{}
$$
$$
{}
$$

\noindent
~~~~~~~~~~~~~~~~~~~~~~~~~~~~~~~~~~~~~~~~~~~~~~~~~~~~~~~~~~~~~~~~~~~~~~~~~~~~~~~~~~~~~~Dedicated to My Family

\newpage
\noindent
$$
{}
$$
$$
{}
$$
$$
{}
$$
$$
{}
$$
$$
{}
$$
$$
{}
$$
$$
{}
$$
$$
{}
$$
\noindent
The first, second, and third editions of this monograph were published in the Journal\\
"Differencialnie Uravnenia i Protsesy Upravlenia"\\
(Differential Equations and Control Processes),\\
no. 4 (2020), A.1-A.606, no. 4 (2021), A.1-A.788, and no. 1 (2023), A.1-A.947\\
Available at:\\
\href{http://diffjournal.spbu.ru/EN/numbers/2020.4/article.1.8.html}
{http://diffjournal.spbu.ru/EN/numbers/2020.4/article.1.8.html}\\
\href{http://diffjournal.spbu.ru/EN/numbers/2021.4/article.1.9.html}
{http://diffjournal.spbu.ru/EN/numbers/2021.4/article.1.9.html}\\
\href{http://diffjournal.spbu.ru/EN/numbers/2023.1/article.1.10.html}
{http://diffjournal.spbu.ru/EN/numbers/2023.1/article.1.10.html}

\newpage

\addcontentsline{toc}{chapter}{Preface}

\vspace{30mm}
$$
{}
$$
$$
{}
$$

\noindent
{\Huge {\bf Preface}}

\vspace{10mm}

\normalsize

\hspace{100mm}~~~~~~~{\it God~does~not~care~about~our}\\ 

\vspace{-7mm}

\hspace{100mm}~~~~~~~{\it mathematical~difficulties. \hspace{-0.9mm}He}\\

\vspace{-7mm}

\hspace{100mm}~~~~~~~{\it integrates~empirically}

\vspace{4mm}

\hspace{100mm}~~~~~~~~~~~~~~~~~~~~---~{\it Albert Einstein}

\vspace{3mm}

\large

\vspace{10mm}

The book is devoted to the problem of strong 
(mean-square) approximation of iterated It\^{o} and Stratonovich 
stochastic integrals in the context of numerical 
integration of It\^{o} stochastic differential equations (SDEs)
and non-commutative semilinear stochastic
partial differential equations (SPDEs) with nonlinear multiplicative 
trace class noise. The presented monograph 
opens up a new direction in researching of iterated
stochastic integrals
and summarizes the author's research on the mentioned problem 
carried out in the period 1994--2026.

The basis of this book composes on
the monographs \cite{1}-\cite{12aa} and recent 
author's results \cite{art-1}-\cite{new-art-1}.

This monograph (also see books \cite{6}-\cite{10}, \cite{12a}-\cite{12aa})
is the first monograph where the problem of strong 
(mean-square) approximation of iterated It\^{o} and Stratonovich 
stochastic integrals with respect to components of a multidimensional
Wiener process is systematically analyzed in application
to the numerical solution of SDEs.

For the first time we successfully use
the generalized multiple Fourier series 
converging in the sense of norm
in Hilbert space $L_2([t, T]^k)$ 
for the expansion and strong approximation of iterated
It\^{o} stochastic integrals of arbitrary multiplicity $k,$ $k\in{\bf N}$
as well as
for the expansion of some other types of iterated 
stochastic integrals (Chapter 1).

The above result has been adapted for iterated Stratonovich stochastic
integrals of multiplicities 1 to 8 for the following two cases (Chapter 2).

1. The case of continuously differentiable 
weight functions (multiplicities 1 to 5) and 
weight functions identically equal to one (multiplicities 6 to 8).
In this case, we use
a complete orthonormal system of Legendre polynomials or 
trigonometric functions in $L_2([t, T])$.

2. The case of continuous weight functions (multiplicities 1 and 2),
binomial weight functions (multiplicities 3 and 4)
and weight functions identically equal to one (multiplicities 5 and 6).
In this case, we use
an arbitrary complete orthonormal system of functions in $L_2([t, T])$.

Recently (in 2024), the mentioned adaptation has also been carried out for 
iterated Stratonovich stochastic
integrals of multiplicity $k,$ $k\in{\bf N}$ (Chapter 2, Theorems~2.59, 2.61) 
but under one additional condition
(the case of continuous weight functions 
and an arbitrary complete orthonormal system of functions in $L_2([t, T])$).

Very recently (2026), Theorems~2.87 and 2.88 (Chapter 2) were proved, 
which justify expansions of iterated Stratonovich stochastic integrals of 
arbitrary multiplicity $k,$ $k\in{\bf N}$ without the additional condition mentioned above. 
These theorems prove Hypotheses 2.1, 2.2, and 2.5 (Chapter 2). The proof of these 
theorems is based on the connection between the results obtained in Chapter 2 and 
rough path theory.

Two theorems on expansions of iterated
Stratonovich stochastic
integrals of multiplicity $k, k\in{\bf N}$ based on 
iterated Fourier series with the pointwise convergence are 
formulated and proved
(Chapter 2).

The integration order replacement technique for 
the class of iterated It\^{o} stochastic integrals has been introduced
(Chapter 3). This result is generalized 
for the class of iterated stochastic integrals with respect to 
martingales.

Four new forms of the Taylor--It\^{o} and Taylor--Stratonovich 
expansions (the so-called unified Taylor--It\^{o} and Taylor--Stratonovich 
expansions) are presented (Chapter~4).

Exact expressions are obtained for the mean-square  
approximation error of iterated It\^{o} stochastic integrals of 
arbitrary multiplicity $k$, $k\in{\bf N}$ (Chapter 1) and iterated Stratonovich stochastic 
integrals of multiplicities 1 to 4 (Chapter 5).
Furthermore, 
we provided a significant practical material (Chapter 5)
devoted to the expansions and approximations of specific 
iterated It\^{o} and Stratonovich stochastic
integrals of multiplicities 1 to 6 from the 
Taylor--It\^{o} and Taylor--Stratonovich 
expansions (Chapter 4) using the system of Legendre polynomials 
and the system of trigonometric functions.

The methods formulated in this book have been compared 
with some existing methods 
of strong appro\-ximation of iterated It\^{o} and Stratonovich
stochastic integrals (Chapter 6).

The results of Chapter 1 were applied (Chapter 7) to the approximation 
of iterated stochastic integrals 
with respect to the finite-dimensional approximation ${\bf W}_t^M$
of the infinite-dimensional $Q$-Wiener process ${\bf W}_t$ (for integrals of 
arbitrary multiplicity $k,$
$k\in{\bf N}$) and to the approximation 
of iterated stochastic integrals 
with respect to the infinite-dimensional 
$Q$-Wiener process ${\bf W}_t$ (for integrals of multiplicities 1 to 3).

This book will be interesting for specialists dealing with the theory 
of stochastic processes, applied and computational mathematics
as well as
senior students and postgraduates of technical institutes and 
universities.

Exact solutions of It\^{o} SDEs and semilinear SPDEs 
are known in rather rare cases. 
Therefore, the
need arises
to construct numerical procedures for solving these equations.

The importance of the problem of numerical integration of 
It\^{o} SDEs and semilinear SPDEs is explained by a wide range of
their applications related to the construction of adequate 
mathematical models of dynamic systems of various
physical nature under random disturbances and to the application 
of these equations for solving
various mathematical problems, among which we mention signal 
filtering in the background of random
noise, stochastic optimal control, stochastic stability, evaluating 
the parameters of stochastic systems,
etc.

It is well known that one of the effective and perspective approaches 
to the numerical integration of It\^{o} SDEs and semilinear SPDEs
is an 
approach based on the stochastic analogues of the Taylor formula 
for solutions of 
these
equations. This approach uses the finite
discretization of temporal variable and performs numerical modeling 
of solutions of It\^{o} SDEs and semilinear SPDEs in discrete moments
of time using stochastic analogues of the Taylor formula.

Speaking about It\^{o} SDEs, note that
the most important feature 
of the mentioned 
stochastic analogues of the Taylor 
formula for solutions of It\^{o} SDEs  
is a presence in them of the so-called iterated It\^{o}
and Stratonovich stochastic integrals 
which are the functionals of a complex structure 
with respect to components of a multidimensional 
Wiener process. These iterated It\^{o} and Stratonovich stochastic integrals
are subject for study in this book
and are 
defined by the following formulas
$$
\int\limits_t^T\psi_k(t_k) \ldots \int\limits_t^{t_{2}}
\psi_1(t_1) d{\bf w}_{t_1}^{(i_1)}\ldots
d{\bf w}_{t_k}^{(i_k)}\ \ \ \ {\rm (It\hat{o}\ \ integrals),}
$$
$$
{\int\limits_t^{*}}^T
\psi_k(t_k) \ldots 
{\int\limits_t^{*}}^{t_{2}}
\psi_1(t_1) d{\bf w}_{t_1}^{(i_1)}\ldots
d{\bf w}_{t_k}^{(i_k)}\ \ \ \ {\rm (Stratonovich\ \ integrals),}
$$
where $\psi_1(\tau),\ldots, \psi_k(\tau): [t, T]\to {\bf R}$ 
are nonrandom functions  
(as a rule, in the applications they are identically equal to 1 or have 
a binomial form (see Chapter 4)),  ${\bf w}_{\tau}$ is a random vector with 
an $m+1$ components: ${\bf w}_{\tau}^{(i)}$ $(i=1,\ldots,m)$
are independent standard Wiener processes and
${\bf w}_{\tau}^{(0)}=\tau,$
$i_1,\ldots,i_k = 0, 1,\ldots,m.$

Apparently, one of the first who began the study of such stochastic 
integrals (the case $k=2,$ $m=2,$ $\psi_1(\tau), \psi_2(\tau)\equiv 1,$ 
$i_1=1,$ $i_2=2$) was L{\'e}vy, who introduced the concept of the so-called 
L{\'e}vy stochastic area and studied its properties.

The above iterated stochastic integrals are the specific objects in the 
theory of stochastic processes. From the one side, nonrandomness of weight 
functions $\psi_l(\tau)$ $(l=1,\ldots,k)$ is the factor simplifying their 
structure. From the other side, nonscalarity of the Wiener process 
${\bf w}_{\tau}$ with independent 
components 
${\bf w}_{\tau}^{(i)}$ $(i=1,\ldots,m)$ and the fact that 
the functions $\psi_l(\tau)$ $(l=1,\ldots,k)$
are different for various $l$ $(l=1,\ldots,k)$ are essential complicating 
factors of the structure of iterated stochastic integrals. 
Taking into account features mentioned above, we suppose that the systems of iterated
It\^{o} and  
Stratonovich stochastic integrals play the extraordinary and perhaps 
the key role 
for solving the problem of numerical integration of It\^{o} SDEs.

A natural question arises: is it possible to construct a numerical 
scheme for It\^{o} SDE that includes only increments of the Wiener 
processes ${\bf w}_{\tau}^{(i)}$ $(i=1,\ldots,m)$, 
but has a higher order of convergence than the Euler method?
It is known that this is impossible for $m>1$ in the general case.
This fact is called the "Clark--Cameron paradox" \cite{Clark1}
and explains the need 
to use iterated stochastic integrals for constructing 
high-order numerical
methods for It\^{o} SDEs.

We want to mention in short that there are 
two main criteria of numerical methods convergence for It\^{o} SDEs: 
a strong or 
mean-square
criterion and a 
weak criterion  where the subject of approximation is not the solution 
of It\^{o} SDE, simply stated, but the 
distribution of It\^{o} SDE solution.
Both mentioned criteria are independent, i.e. in general  it is 
impossible to state that from the execution of strong criterion follows 
the execution of weak criterion and vice versa.
Each of two convergence criteria is oriented on the solution of specific 
classes of mathematical problems connected with It\^{o} SDEs.

Numerical
integration of It\^{o} SDEs based on the strong convergence criterion of 
approximation is widely used
for the numerical simulation of sample trajectories of solutions 
to It\^{o} SDEs (which is required for
constructing new mathematical models based on
such equations and for the numerical solution
of different mathematical problems connected with It\^{o} SDEs). 
Among these problems, we note the
following: signal filtering under influence of random
noises in 
various statements (linear Kalman--Bucy filtering, nonlinear optimal 
filtering, filtering of continuous time Markov chains with a finite
space of states, etc.), optimal stochastic control (including 
incomplete data control), testing estimation
procedures of parameters of stochastic systems, stochastic stability 
and bifurcations analysis.

The problem of effective jointly numerical modeling 
(with respect to the mean-square convergence criterion) of iterated
It\^{o} or Stratonovich stochastic integrals is very important and difficult from 
theoretical and computing point of view.

Seems that iterated stochastic integrals may be approximated by multiple 
integral sums. However, this approach implies the partitioning of the interval 
of integration $[t, T]$ for iterated stochastic integrals. The length
$T-t$ of this interval is already fairly
small (because it is a step of integration of numerical methods for 
It\^{o} SDEs) and does not need to be partitioned.
Computational experiments show that the application
of numerical simulation for iterated stochastic integrals 
(in which the interval of integration is
partitioned) leads to unacceptably high computational cost and 
accumulation of computation errors.

The problem of effective decreasing of the mentioned cost (in several 
times or even in several orders) is very difficult and 
requires new complex investigations. 
The only exception is connected with a narrow particular case, when 
$i_1=\ldots=i_k\ne 0$ and
$\psi_1(\tau),\ldots,\psi_k(\tau)\equiv \psi(\tau)$.
This case allows 
the investigation with using of the It\^{o} formula.
In the more general case, when not all numbers $i_1,\ldots,i_k$
are equal, the 
mentioned problem turns out to be more complex (it cannot be solved 
using the It\^{o} formula and requires more deep 
and complex investigation). Note
that even for the case $i_1=\ldots=i_k\ne 0$, but for different 
functions $\psi_1(\tau),\ldots,\psi_k(\tau)$ the mentioned difficulties persist and 
simple sets of 
iterated It\^{o} and Stratonovich stochastic integrals, which can be often 
met in the applications, cannot be expressed effectively in a finite 
form (with respect to the mean-square approximation) 
using the system of standard 
Gaussian random variables. The It\^{o} formula is also useless in this case 
and as a result we need to use more complex but effective expansions.

Why the problem of the mean-square approximation of iterated stochastic 
integrals is so complex?

Firstly, the mentioned stochastic integrals (in the case of fixed limits 
of integration) are the random variables, whose density functions are 
unknown in the general case. The exception is connected with
the narrow particular case which is
the simplest iterated It\^{o} stochastic integral with multiplicity 2
and $\psi_1(\tau), \psi_2(\tau)\equiv 1;$
$i_1, i_2=1,\ldots,m$.
Nevertheless, the knowledge of this density function
not gives a simple way for approximation of iterated
It\^{o} stochastic integral of multiplicity 2.

Secondly, we need to approximate not only one stochastic integral, 
but several iterated stochastic integrals that are complexly dependent
in a probabilistic sense.

Often, the problem of combined mean-square approximation of iterated 
It\^{o} and Stratonovich stochastic integrals occurs even in cases when the 
exact solution of It\^{o} SDE is known.  It means
that even if you know the solution of It\^{o} SDE exactly, 
you cannot model it numerically without the combined
numerical modeling of 
iterated stochastic integrals.

Note that for a number of special types of It\^{o} SDEs
the problem of approximation of iterated
stochastic integrals may be simplified but cannot be solved. 
Equations with additive vector noise, with non-additive scalar noise, 
with additive scalar noise, with a small parameter are related to such 
types of equations. In these cases, simplifications are connected to the fact 
that some members
from stochastic Taylor expansions are equal to zero 
or we may neglect some members from these expansions
due to the presence of a small 
parameter.

Furthermore, the problem of combined numerical modeling (with respect
to the mean-square convergence criterion) of iterated It\^{o} 
and Stratonovich stochastic integrals is rather new.  

One of the main and unexpected achievements of this book is the successful 
usage of functional analysis methods (more concretely, we mean
generalized multiple
Fourier series in various systems of 
basis functions that converge in the sense of the norm in $L_2([t, T]^k)$) 
in this scientific field.

The problem of combined numerical modeling (with respect to the mean-square 
convergence criterion) of systems of iterated It\^{o} and Stra\-to\-no\-vich 
stochastic integrals 
was 
analyzed in the context of the problem of numerical integration 
of It\^{o} SDEs in the following monographs:

\vspace{2mm}

[I] Milstein G.N. Numerical Integration of Stochastic Differential 
Equations. Kluwer Academic Publishers. Dordrecht. 1995 (Russian 
Ed. 1988).
    
[II] Kloeden P.E., Platen E. Numerical Solution of Stochastic Differential 
Equations. Springer-Verlag. Berlin. 1992 (2nd Ed. 1995, 3rd
Ed. 1999).

[III] Milstein G.N., Tretyakov M. V. Stochastic Numerics for Mathematical 
Physics. Springer-Verlag. Berlin. 2004 (2nd Ed. 2021).

[IV] Kuznetsov D.F.  Stochastic Differential Equations: Theory and Practice 
of Numerical Solution. Polytechnical 
University Publ. St.-Petersburg. 2007 \cite{2}  (2nd Ed. 2007 \cite{3},
3rd Ed. 2009 \cite{4}, 4th Ed. 2010 \cite{5}, 5th Ed. 2017 \cite{11}, 
6th Ed. 2018 \cite{12}, 7th Ed. \cite{123789000}).

Note that the initial version of the book [IV] has been published in 2006 
\cite{1}. Also we mention the books
\cite{6} (2010), \cite{7} (2011), \cite{8} (2011), \cite{8xxxxxx} (2012), 
\cite{9} (2013), \cite{10} (2017) and \cite{12a} (2020), \cite{12aa-after} (2021),
\cite{12aa-afterxxx} (2023), \cite{12aa} (2026).

The books [I] and [III] analyze the problem of the
mean-square approximation of iterated stochastic integrals only 
for two simplest iterated It\^{o} stochastic integrals of 1st and 
2nd  multiplicities 
($k=1$ and $2,$
$\psi_1(\tau)$ and $\psi_2(\tau)\equiv 1$) 
for the 
multidimensional
case:
$i_1, i_2=0, 1,\ldots,m$. 
In addition, the main idea is based on the 
expansion of the so-called Brownian bridge process into the  
trigonometric Fourier series (version of the so-called Karhunen--Lo\`{e}ve
expansion). 
This method is called in [I] and [III] as the Fourier method${}^{1}$.

\footnotetext[1]{To date, there is confusion in the literature about who first 
proposed the Fourier method
[I], [III]. As far as the author of this book knows, 
the mentioned method first appeared in the Russian edition of the monograph 
by G.N. Milstein \cite{Zapad-1} (pp.~121--135), which was published in 1988.}

In [II] using the Fourier method [I], the attempt was made to obtain
the mean-square approximation of elementary iterated Stratonovich stochastic integrals of 
multiplicities 1 to 3 ($k=1,\ldots, 3,$
$\psi_1(\tau),\ldots,\psi_3(\tau)\equiv 1$)
for the multidimensional case: 
$i_1,\ldots,i_3=0, 1,\ldots,m$. However, as we can see in the
presented book, the 
results of the monograph [II], related to the mean-square approximation of 
iterated Stratonovich stochastic integrals of 3rd multiplicity, cause a number of critical 
remarks (see discussions in Sect.~2.42, 2.43, 6.2).
   
The main purpose of this book is to construct and develop
newer 
and more effective 
methods (than presented in the books [I]--[III]) of 
combined mean-square approximation of iterated It\^{o} and 
Stratonovich stochastic integrals.

Talking about the history of solving the problem of combined mean-square 
approximation of iterated stochastic integrals, the idea to find a basis 
of random variables
using which we may represent iterated 
stochastic 
integrals turned out to be useful. This idea was transformed several times 
during last decades.
   
Attempts to approximate the iterated  
stochastic integrals using various integral 
sums were made until 1980s and later, i.e. the interval of integration 
$[t, T]$ of the stochastic 
integral was divided
into $n$ parts and the iterated stochastic integral was 
represented approximately by the multiple integral sum,
which included the system of independent standard Gaussian random variables, 
whose numerical modeling is not a problem.
   
However, as we noted above, it is obvious that the length $T-t$
of integration interval $[t, T]$
of the iterated  stochastic integrals is a step of integration of
numerical methods  
for It\^{o} SDEs, which is already a rather 
small value even without the additional splitting. Numerical experiments 
demonstrate that such approach results in drastic increasing of computational 
costs accompanied by the growth of multiplicity of the stochastic integrals 
(beginning from 2nd and 3rd multiplicity) that is necessary for construction
of high-order strong numerical methods for It\^{o} SDEs
or in the case of decrease of integration step of numerical methods, 
and thereby it is almost useless for practice.

The new step for solution of the problem of combined mean-square 
approximation of iterated  stochastic integrals 
was made by Milstein G.N. in his 
monograph [I] (1988). 
For the expansion of iterated stochastic integrals, 
he proposed to use the trigonometric Fourier expansion of the Brownian bridge process
(version of the so-called Karhunen--Lo\`{e}ve 
expansion). Using this method, expansions of two 
simplest iterated It\^{o} stochastic integrals of multiplicities 1 and 2
are obtained and their mean-square convergence is proved.
   
As we noted above, the attempt to develop this idea 
together with the Wong--Zakai approximation \cite{W-Z-1}-\cite{Watanabe}
was made in the monograph 
[II] (1992), where 
the expansions of simplest iterated 
Stratonovich 
stochastic integrals of multiplicities 1 to 3 were
obtained. However, due to a number 
of limitations and technical difficulties which are typical for the method 
[I], in [II] and following 
publications 
this problem was not solved more completely.   
In addition, the author has reasonable doubts about 
application of the Wong--Zakai results
\cite{W-Z-1}-\cite{Watanabe} to
approximation of iterated Stratonovich stochastic integrals of 3rd 
multiplicity in the monograph [II] (see discussions in Sect.~2.42, 2.43, 6.2).  

It is necessary to note that the 
computational cost
for the method [I] 
is significantly less than
for the method of multiple 
integral sums.

Regardless of the method [I] positive features, the number of its 
limitations are also outlined. Among them let us mention the following.

1.\;The absence of explicit formula for 
calculation of expansion coefficients for iterated  stochastic 
in\-teg\-rals.

2.\;The practical impossibility of exact calculation of the 
mean-square approximation error
of iterated  stochastic integ\-rals with the exception of 
simplest integ\-rals of 1st and 2nd multiplicity (as a result, it is 
necessary to consider redundant terms of expansions and it results 
to the growth
of computational cost and complication of the numerical methods 
for It\^{o} SDEs).

3.\;There is a hard limitation on the system of basis functions 
---
it may be only the trigonometric functions.

4.\;There are some technical 
problems if we use this method for iterated stochastic integrals whose 
multiplicity is greater than 2nd. 

Nevertheless, it should be noted that the analyzed method is a concrete 
step forward in this scientific field.
   
The author 
thinks that the method presented by him in [IV] (for the first time 
this method is appeared in the final form in \cite{1} (2006)) 
and in this book 
(hereafter this method is reffered to as 
the method of generalized multiple Fourier series) is a 
breakthrough 
in solution of the problem of combined mean-square approximation of 
iterated It\^{o} stochastic integrals.
   
The idea of this method is as follows: the iterated It\^{o} stochastic 
integral of multiplicity $k,$ $k\in{\bf N}$
is represented as the multiple stochastic 
integral from the certain nonrandom discontinuous function of $k$ variables 
defined on the hypercube $[t, T]^k,$ where $[t, T]$ is the interval of 
integration of the iterated It\^{o} stochastic integral. Then, the mentioned
nonrandom function of $k$ variables
is expanded in the hypercube $[t, T]^k$ into the generalized 
multiple Fourier series converging 
in the mean-square sense
in the space 
$L_2([t,T]^k)$. After a number of nontrivial transformations we come 
to the 
mean-square converging expansion of the iterated It\^{o} stochastic 
integral into the multiple 
series of products
of standard  Gaussian random 
variables. The coefficients of this 
series are the coefficients of 
generalized multiple Fourier series for the mentioned nonrandom function 
of $k$ variables, which can be calculated using the explicit formula 
regardless 
of the multiplicity $k$ of the iterated It\^{o} stochastic integral.

As a result, we obtain the following new possibilities and advantages 
in comparison with the Fourier method [I].

1.\;There is an explicit formula for calculation of expansion coefficients 
of iterated It\^{o} stochastic integral with any
fixed multiplicity $k$. 
In other words, we can calculate (without any preliminary and additional 
work) the expansion coefficient with any fixed number in the expansion 
of iterated It\^{o} stochastic integral of the preset fixed multiplicity.
At that, we do not need any knowledge about coefficients with other numbers or 
about other iterated It\^{o} stochastic integrals included in the 
considered set.

2.\;We have new possibilies for obtainment the exact and approximate 
expressions for the mean-square
approximation errors of iterated It\^{o} stochastic integrals. These 
possibilities are realized by the exact and estimate formulas 
for the mentioned mean-square approximation errors. 
As a result, we would not need to consider 
redundant terms of expansions that may complicate approximations 
of iterated It\^{o} stochastic integrals.

3.\;Since the used multiple Fourier series is a generalized in the sense
that it is built using various complete orthonormal
systems of functions in the space $L_2([t, T]^k)$, we have new possibilities 
for approximation --- we can 
use 
not only the trigonometric functions as in [I] but the Legendre polynomials 
as well as the systems of Haar and Rademacher--Walsh functions.

4.\;As it turned out, it is more convenient to work 
with Legendre polynomials for approxi\-mation of iterated
It\^{o} stochastic integrals.
The ap\-pro\-xi\-mations themselves 
are simpler than their analogues based on
the system of trigonometric functions.   
Probably for the systems of Haar and Rademacher--Walsh functions the expansions 
of iterated stochastic integrals become more complex and 
less effective for practice [IV]. Expansions based on Haar functions
for $k=2$ were also considered in \cite{Zapad-5}, \cite{Rybak1}, \cite{Minoo}. 
Note that the multiple Fourier--Walsh and Fourier--Haar series $(k\in{\bf N})$ were applied
to the mean-square approximation of multiple Stratonovich 
stochastic integrals (defined as in \cite{SU11}, \cite{bardina10}) in \cite{Rybakov3000xxx}. 
The convergence of these approximations
was proved with respect to the special subsequence $n_m=2^m$ $(m\to\infty)$ \cite{Rybakov3000xxx}.

5.\;The question about what kind of functions (polynomial or trigonometric) 
is more convenient in the context of computational costs for approximation 
turns out to be nontrivial, since it is necessary to compare 
approximations not for one stochastic integral but for several stochastic 
integrals at the same time. At that there is a possibility 
that computational costs for some integrals will be smaller for 
the system of Legendre polynomials and for others --- for the system 
of trigonometric functions.
The author proved \cite{art-4} (also see Sect.~5.3 in this book) that the
computational costs are significantly less for the system of Legendre 
polynomials at least in the case of approximation of the special set of 
iterated It\^{o} stochastic integrals, which are necessary for the
implementation of strong numerical methods for It\^{o} SDEs with
the order of convergence $\gamma=1.5$.
In addition, 
the author supposes that this effect will be more impressive when 
analyzing more complex sets of iterated It\^{o} stochastic
integrals ($\gamma=$ $2.0,$ $2.5,$ $3.0,$ $\ldots $). 
This supposition is based on the fact that the polynomial 
system of functions has a significant advantage (in comparison with 
the trigonometric system of functions) 
in the mean-square approximation of iterated It\^{o} stochastic 
integrals for which not all weight functions are equal to 1.

6.\;The Milstein approach [I] for approximation of iterated It\^{o} 
stochastic integrals leads to 
iterated applicaton of the
operation of li\-mit transition (in contrast with the method 
of generalized multiple Fourier series, for which the operation 
of limit transition is implemented only once) 
starting at least from the second or third multiplicity of 
iterated It\^{o} stochastic integrals (we mean at least double or 
triple integration with respect to components
of a multidimensional Wiener process).
Multiple series are more 
preferential for approximation than the iterated ones, since the
partial sums of multiple series converge for any possible case of joint 
converging to infinity of their upper limits of summation (let us denote 
them as $p_1,\ldots, p_k$). 
For example,
when $p_1=\ldots=p_k=p\to\infty$. 
For iterated series, the condition $p_1=\ldots=p_k=p\to\infty$ obviously 
does not guarantee the convergence of this series.
However, in [II] the authors use (without rigorous proof) the condition 
$p_1=p_2=p_3=p\to\infty$
within the frames of the Milstein approach [I]
together with the Wong--Zakai approximation \cite{W-Z-1}-\cite{Watanabe}
(see discussions in Sect.~2.42, 2.43, 6.2).

7.\;The convergence in the mean of degree 
$2n,$ $n\in {\bf N}$ as well as the convergence
with probability 1 of approximations from the 
method of ge\-neralized multiple Fourier series
are proved. The convergence rate for these two types
of convergence is estimated.

8.\;The method of generalized multiple Fourier series
has been applied for some other types of iterated stochastic integrals
(iterated stochastic integrals with respect to martingale 
Poisson random measures
and iterated stochastic integrals with respect to martingales)
as well as for approximation of iterated stochastic 
integrals with respect to the infinite-dimensional
$Q$-Wiener process.

9.\;Another modification of the method of generalized
multiple Fourier series is connected with the 
application of complete orthonormal with weight $r(t_1)\ldots r(t_k)\ge 0$
systems of functions in the space $L_2([t, T]^k)$.

10.\;As it turned out, the method of generalized multiple Fourier
series can be adapted 
for iterated Stratonovich
stochastic integrals. This adaptation is carried out in Chapter 2 
for the following two cases.

1).\ The case of continuously differentiable 
weight functions (multiplicities 1 to 5) and 
weight functions identically equal to one (multiplicities 6 to 8).
In this case, we use
a complete orthonormal system of Legendre polynomials or 
trigonometric functions in $L_2([t, T])$.

2).\ The case of continuous weight functions (multiplicities 1 and 2),
binomial weight functions (multiplicities 3 and 4)
and weight functions identically equal to one (multiplicities 5 and 6).
In this case, we use
an arbitrary complete orthonormal system of functions in $L_2([t, T])$.

Recently (in 2024), the mentioned adaptation has also been carried out for 
iterated Stratonovich stochastic
integrals of multiplicity $k,$ $k\in{\bf N}$ but under one additional condition
(the case of continuous weight functions
and an arbitrary complete orthonormal system of functions in $L_2([t, T])$
(Chapter 2, Theorems~2.59, 2.61)).
The rate of mean-square convergence of approximations
of iterated Stratonovich stochastic integrals 
is found (Sect.~2.8, 2.15, 2.16).

Very recently (2026), Theorems~2.87 and 2.88 (Sect.~2.50) were proved, 
which justify expansions of iterated Stratonovich stochastic integrals of 
arbitrary multiplicity $k,$ $k\in{\bf N}$ without the additional condition mentioned above
(the case of continuously differentiable weight functions
and an arbitrary complete orthonormal system of functions in $L_2([t, T])$). 
These theorems prove Hypotheses 2.1, 2.2, and 2.5 (Sect.~2.5, 2.28). The proof of these 
theorems is based on the connection between the results obtained in Chapter 2 and 
rough path theory.

11.\;The method of generalized multiple Fourier
series is reformulated using Hermite polynomials in Sect.~1.10 and generalized 
to the case of an arbitrary complete orthonormal system of functions in the space 
$L_2([t, T])$ and $\psi_1(\tau),\ldots,\psi_k(\tau)\in L_2([t, T])$
in Sect.~1.11, 1.12, 1.14, 1.15.
At that, in Sect.~1.11, 1.12 we use the multiple Wiener stochastic integral 
with respect to the components of a multidimensional Wiener process.

12.\;The results of Chapter 1 
(Theorems 1.1, 1.2, 1.14, 1.16) and Chapter 2 
(Theorems 2.1--2.10, 2.14, 2.17, 2.30, 2.32--2.36, 2.41--2.51, 2.53, 2.55, 2.57, 2.59, 2.61--2.65,
2.87, 2.88) 
can be considered 
from the point of view of the Wong--Zakai approximation 
\cite{W-Z-1}-\cite{Watanabe}
for the case of a multidimensional Wiener process and the Wiener 
process approximation based on its series expansion
using various complete orthonormal systems of functions 
in the space $L_2([t,T])$
(see discussions in Sect.~2.42, 2.43, 6.2).
These results overcome a number of difficulties that 
were noted above and relate to the Fourier method [I].

In connection with the mention of iterated Stratonovich stochastic integrals, 
it is worth noting the scientific direction that was formed 
in the 1990s and is dedicated to multiple Stratonovich stochastic 
integrals and their representations
\cite{bugh1}-\cite{bugh3}, \cite{HuHu} (Chapter~5).
Let us note two key points that distinguish Chapter~2 of this book from 
\cite{bugh1}-\cite{bugh3}, \cite{HuHu} (Chapter~5).
In \cite{bugh1}-\cite{bugh3}, \cite{HuHu} (Chapter~5) the Wiener process is scalar,
while in this book the Wiener process is a multidimensional
process with independent components.
In addition, the conditions for the validity of the expansions of multiple Stratonovich 
stochastic integrals are formulated in \cite{bugh1}-\cite{bugh3}, \cite{HuHu} (Chapter~5) in terms of the 
existence of so-called limiting traces, but these conditions are 
not verified in \cite{bugh1}-\cite{bugh3}, \cite{HuHu}.
On the other hand, in Chapter 2 of this book, along with finding sufficient conditions 
for expansions of iterated Stratonovich stochastic integrals, much attention is paid 
to verifying the indicated conditions.

The theory presented in this book was realized \cite{Kuz-Kuz}, \cite{Mikh-1}
in the form of a software package in the Python programming language.
The mentioned software package implements the 
strong numerical methods with convergence orders 
0.5, 1.0, 1.5, 2.0, 2.5, and 3.0 for It\^{o}
SDEs (with multidimensional non-commutative noise)
based on the unified Taylor--It\^{o} 
and Taylor--Stratonovich expansions (Chapter 4).
At that for the numerical simulation of iterated 
It\^{o} and Stratonovich stochastic integrals of multiplicities 
1 to 6 we applied the formulas 
based on multiple Fourier--Legendre series (Chapter 5).
Moreover, we used \cite{Kuz-Kuz}, \cite{Mikh-1}
the database with 270,000 exactly calculated
Fourier--Legendre coefficients.

Throughout the book, special attention is paid to two systems 
of basis functions in the space $L_2([t,T])$. Namely, we mainly use 
the complete orthonormal systems of Legendre polynomials and trigonometric 
functions in the space $L_2([t,T])$.
This is due to two reasons. 
The first of these is that the trigonometric basis system has already been 
used to approximate iterated stochastic integrals 
in the 1980s-1990s (see above), and the author needed to compare his results with the 
results of other authors.
The second reason is that the system of Legendre polynomials is 
optimal (see Sect.~5.3) for the implementation of strong numerical methods with 
convergence order 1.5 and higher for It\^{o} SDEs with 
multidimensional non-commutative noise.
The system of Legendre polynomials was first applied to the approximation of 
iterated stochastic integrals in the author's work \cite{old-art-1} 
in 1997 (also see \cite{old-art-2}-\cite{old-art-4}).
According to the author's opinion, other complete orthonormal systems of functions in the space 
$L_2([t,T])$ (for example, systems of Haar and Rademacher--Walsh functions) 
turn out to be less efficient for the mean-square approximation of 
iterated It\^{o} and Stratonovich stochastic integrals.

The attentive reader will notice that Chapters 1 and 2 of this book can be somewhat 
shortened since Theorem 1.16 is a generalization of Theorems~1.1, 1.2 and
Theorems 2.3, 2.33, 2.34, 2.41 are generalizations of Theorems 2.1, 2.2, 2.4--2.9.
Theorems~2.87 and 2.88 also supersede some results from the first part
of Chapter 2 (these theorems were obtained very recently (2026), when the book
was already considered complete).
However, the author did not make the appropriate changes 
in Chapters 1, 2 for a number of reasons.
In particular, the application of the multiple Wiener 
stochastic integral with respect to the components of a multidimensional 
Wiener process to the expansion of iterated It\^{o} stochastic integrals
(Theorem~1.16) and a new approach to the expansion of 
iterated Stratonovich stochastic integrals (Theorems 2.30--2.65) were obtained 
by the author recently (in 2021--2025), while 
Theorems~1.1, 1.2, 2.4--2.9 were obtained by the author in the period 
from 2005 to 2013. Moreover, Theorems 2.87 and 2.88 
were added to a virtually completed book.
In addition, the proof of each of the mentioned theorems contains 
some original ideas that the author would like to keep in Chapters 1 and 2.
Moreover, a significant part of Chapter 2 is devoted 
to the proof of Hypothesis~2.5 (Sect.~2.28)
for various special cases (Theorems~2.1--2.9,
2.30, 2.33--2.36, 2.41, 2.45--2.48, 2.50, 2.51, 2.53, 2.55, 2.57, 2.59, 2.61--2.65).
In order to prove these theorems, we developed 
a number of approaches to the expansion of iterated Stratonovich stochastic
integrals. However, Hypothesis 2.5 was proved by Theorems 2.87 and 2.88.
Finally, the results of Chapters 1, 2 are presented primarily in the order 
in which they were obtained by the author.

\vspace{10mm}

Dmitriy F. Kuznetsov~~~~~~~~~~~~~~~~~~~~~~~~~~~~~~~~~~~~~~~~~~~~November, 
2026

\newpage

\addcontentsline{toc}{chapter}{Acknowledgements}

$$
{}
$$
$$
{}
$$
\noindent
{\bf {\Huge Acknowledgements}}
$$
{}
$$
$$
{}
$$

I would like to thank the Deputy Editor of the Journal 
"Differencialnie Uravnenia i Protsesy Upravlenia"
Dr. Nataly B. Ampilova for her timely administrative support 
and encouragement and Dr. Konstantin A. Rybakov 
for useful discussion of some presented results.

\begin{table}
\vspace{-30mm}
\begin{tabular}{p{5cm}p{11cm}}

\addcontentsline{toc}{chapter}{Basic Notations}

\vspace{30mm}\\

\multicolumn{2}{c}{}\\

\begin{tabular}{p{11cm}p{5cm}}
{\Huge {\bf \hspace{-6mm} Basic Notations}}
\end{tabular}

\vspace{10mm}\\

${\bf N}$&set of natural numbers

~\\

${\bf R},\ {\bf R}^1$&set of real numbers

~\\

${\bf R}^n$&$n$-dimensional Euclidean space

~\\

$(a_1,\ldots,a_n)$&ordered set with elements $a_1,\ldots,a_n$

~\\

$\{a_1,\ldots,a_n\}$&unordered set with elements $a_1,\ldots,a_n$

~\\

$n!$&$1\cdot 2 \cdot \ldots \cdot n$ 
for $n\in {\bf N}$\ \ $(0!=1)$

~\\

$(2n-1)!!$&$1\cdot 3 \cdot \ldots \cdot (2n-1)$ 
for $n\in {\bf N}$

~\\

$\stackrel{\sf def}{=}$&equal by definition

~\\

$\equiv$&identically equal to

~\\

$C_n^m$&binomial coefficient 
$n!/(m!(n-m)!)$

~\\

$\emptyset$&empty set

~\\

${\bf 1}_{A}$&indicator of the set $A$

~\\

$x\in X$&$x$ is an element of the set $X$

~\\

$X \bigcup Y$&union of sets $X$ and $Y$

~\\

$X\times Y$&Cartesian product of sets $X$ and $Y$

~\\

$C(D)$&set of continuous functions in $D$

~\\          

$\varlimsup\limits_{n\to\infty}$&$\limsup\limits_{n\to\infty}$

~\\

$\varliminf\limits_{n\to\infty}$&$\liminf\limits_{n\to\infty}$

~\\

$[x]$&largest integer number not exceeding $x$

\end{tabular}
\end{table}

\noindent
\begin{table}
\begin{tabular}{p{5cm}p{11cm}}

$\vspace{2mm}|x|$&absolute value of the real number $x$

~\\

$\vspace{2mm}F:\ X\rightarrow Y$&function $F$ from $X$ into $Y$\\

~\\

$A^{(ij)}$&$ij$th element of the matrix $A$

~\\

$A_i$&$i$th colomn of the matrix $A$

~\\

${\bf x}^{(i)}$&$i$th component of the vector
${\bf x}\in{\bf R}^n$

~\\

$O(x)$&expression being divided by $x$ remains bounded as 
$x\to 0$

~\\

$\sum\limits_{(i_1,\ldots,i_k)}$&sum with respect to all 
possible permutations $(i_1,\ldots,i_k)$

~\\

\vspace{0.1mm}

${\sf M}\{\xi\}$&\vspace{0.35mm}expectation of $\xi$

~\\

${\sf M}\{\xi|{\rm F}\}$&conditional expectation of 
$\xi$ with respect to ${\rm F}$

~\\

$\xi\sim {\rm N}(m,\sigma^2)$&Gaussian random variable $\xi$
with expectation $m$ and  variance
$\sigma^2$

~\\

$\hbox{\vtop{\offinterlineskip\halign{
\hfil#\hfil\cr
{\rm l.i.m.}\cr
$\stackrel{}{{}_{n\to \infty}}$\cr
}} }$&limit in the mean-square sense

~\\

${\cal B}(X)$&$\sigma$-algebra 
of Borel subsets of $X$

~\\

$w_t$&scalar standard Wiener process

~\\

${\bf w}_t^{(i)},$ $i=1,\ldots,m$&independent standard Wiener processes

~\\

${\bf w}_t$&vector with components 
${\bf w}_t^{(i)},$ $i=0, 1,\ldots,m,$ where
${\bf w}_t^{(i)},$ $i=1,\ldots,m$ are independent standard Wiener processes
and ${\bf w}^{(0)}_t=t$

~\\

w.~p.~1&with probability 1

~\\

\vspace{-4.5mm}
$$
\hspace{-42mm}\frac{\partial F}{\partial {\bf x}^{(i)}}
$$
&\vspace{-2mm}partial
derivative of $F:$ ${\bf R}^n\to{\bf R}$

~\\

\vspace{-6.5mm}
$$
\hspace{-33mm}\frac{\partial^2 F}{\partial{\bf x}^{(i)}
\partial{\bf x}^{(j)}}
$$&\vspace{-4.5mm}2nd order partial derivative 
of $F:$ ${\bf R}^n\to{\bf R}$

~\\

\vspace{-5.5mm}
$$
\hspace{-31mm}\int\limits_t^T \ldots d{\bf w}_{\tau}^{(i)}
$$
&\vspace{-0.2mm}It\^{o} stochastic integral

\end{tabular}
\end{table}

\noindent
\begin{table}
\vspace{-2mm}
\begin{tabular}{p{5cm}p{11cm}}

\vspace{-6mm}
$$
\hspace{-29mm}{\int\limits_{t}^{*}}^T \ldots d{\bf w}_{\tau}^{(i)}
$$
&\vspace{0mm}Stratonovich stochastic integral

~\\

\vspace{-6mm}
$$
\hspace{-27mm}\int\limits_t^T \ldots \circ d{\bf w}_{\tau}^{(i)}
$$
&\vspace{-0.5mm}Stratonovich stochastic integral \cite{SU11}

~\\

\vspace{-3mm}
${\bf W}_t$&\vspace{-3mm}$Q$-Wiener process

~\\

$J[\psi^{(k)}]_{T,t},\ \
I_{(l_1\ldots l_k)T,t}^{(i_1\ldots i_k)}$&iterated
It\^{o} stochastic integrals

~\\

$J^{*}[\psi^{(k)}]_{T,t},\ \
I_{(l_1\ldots l_k)T,t}^{*(i_1\ldots i_k)}$&iterated
Stratonovich stochastic integrals

~\\

$J^S[\psi^{(k)}]_{T,t}^{(i_1\ldots i_k)}$&iterated
Stratonovich stochastic integral \cite{bardina10}

~\\

$\bar J^S[\psi^{(k)}]_{T,t}^{(i_1\ldots i_k)}$&multiple
Stratonovich stochastic integral \cite{bardina10}

~\\

$J[\psi^{(k)}]^{p_1,\ldots,p_k}_{T,t},\ \
I_{(l_1\ldots l_k)T,t}^{(i_1\ldots i_k)p}$&approximations of iterated
It\^{o} stochastic integrals

~\\

$J^{*}[\psi^{(k)}]^{p}_{T,t},\ \
I_{(l_1\ldots l_k)T,t}^{*(i_1\ldots i_k)p}$&approximations of iterated
Stratonovich stochastic integrals

~\\

$J[\Phi]_{T,t}^{(k)},\ \
J[\Phi]_{T,t}^{(i_1\ldots i_k)}$&multiple Stratonovich
stochastic integrals

~\\

$J'[\Phi]_{T,t}^{(k)},\ \
J'[\Phi]_{T,t}^{(i_1\ldots i_k)}$&multiple Wiener
stochastic integrals

~\\

$P_n(x)$&Legendre polynomials

~\\

$H_n(x),\ h_n(x)$&Hermite polynomials

~\\

$H_n(x,y)$&polynomials related to the Hermite polynomials

~\\

$L_2(D)$&Hilbert space of square integrable functions on $D$

~\\

$\left\Vert \cdot \right\Vert_{L_2(D)}$&norm in the 
Hilbert space $L_2(D)$

~\\

${\rm tr}\ A$&trace of the operator $A$

~\\

$\left\Vert \cdot \right\Vert_{H}$&norm in the Hilbert space $H$

~\\

$\langle u, v\rangle_H$&scalar product in the Hilbert space $H$

~\\

$L_{HS}(U,H)$&space of Hilbert--Schmidt operators 
from $U$ to $H$
\end{tabular}
\end{table}

\newpage
\noindent
\begin{table}
\vspace{-2mm}
\begin{tabular}{p{5cm}p{11cm}}

\vspace{1mm}

$\left\Vert \cdot \right\Vert_{L_{HS}(U,H)}$&\vspace{2mm}operator
norm in the space of Hilbert--Schmidt operators from $U$ to $H$

~\\

\vspace{-5mm}

$$
\hspace{-31mm}\int\limits_t^T \ldots d{\bf W}_{\tau}
$$
&\vspace{1mm}stochastic integral 
with respect to the $Q$-Wiener process

\end{tabular}
\end{table}

\newpage

{\normalsize 
\tableofcontents
}

\chapter{Method of Expansion and Mean-Square Approximation
of Iterated It\^{o} Sto\-chas\-tic Integrals
Based on Generalized Multiple Fourier Series}

This chapter is devoted to the expansions of iterated
It\^{o} stochastic integrals with respect to
components of the multidimensional Wiener process
based on generalized multiple Fourier
series converging in the sense of norm in the space
$L_2([t, T]^k),$ $k\in {\bf N}.$ The
method of generalized multiple Fourier
series for ex\-pansion and mean-square approximation
of iterated It\^{o} stochastic integrals of arbitrary
multiplicity $k,$ $k\in {\bf N}$ is proposed and developed.
The obtained expansions contain only one operation
of the limit transition in contrast to existing analogues.
In this chapter it is also obtained the generalization 
of the proposed method for the case of an arbitrary complete orthonormal
system of functions in the space 
$L_2([t, T]^k),$ $k\in {\bf N}$ as well as
for the case of complete orthonormal with weight 
$r(t_1)\ldots r(t_k)\ge 0$
systems of functions in the space 
$L_2([t, T]^k),$ $k\in {\bf N}.$
It is shown that in the case of a scalar Wiener process
the proposed method leads to
the well known
expansion of iterated It\^{o} stochastic integrals
based on the It\^{o} formula and Hermite polynomials.
The convergence in the mean of degree $2n,$ $n \in {\bf N}$
as well as the convergence
with probability 1 of the proposed method are proved. 
The exact and approximate 
expressions for the mean-square
approximation error of iterated It\^{o} stochastic integrals
of multiplicity $k,$ $k\in{\bf N}$ have been derived.
The considered method 
has been applied for other types of iterated stochastic integrals
(iterated stochastic integrals with respect to martingale 
Poisson random measures
and iterated stochastic integrals with respect to martingales).

\section{Expansion of Iterated It\^{o} Stochastic Integrals 
of Arbitrary Multiplicity
Based on 
Generalized Multiple Fourier Series Converging in the Mean}

\vspace{5mm}

\subsection{Introduction}

The idea of representing the iterated It\^{o} and Stratonovich stochastic 
integrals 
in the form of multiple stochastic integrals from 
specific discontinuous nonrandom functions of several variables and following 
expansion of these functions using multiple and iterated 
Fourier series in order 
to get effective mean-square approximations of the mentioned stochastic 
integrals was proposed and developed in a lot of 
author's publications \cite{1}-\cite{new-2023a} (also see
early publications \cite{old-art-1} (1997), \cite{old-art-2} (1998),
\cite{old-art-3} (2000), \cite{old-art-4} (2001),
\cite{very-old-1} (1994), \cite{very-old-2} (1996)). 
Note that another
approaches to the mean-square approximation 
of iterated It\^{o} and Stratonovich stochastic integrals can be found in 
\cite{new-art-1}, \cite{Zapad-1}-\cite{Zapad-12xxx}. 

Specifically, 
the approach \cite{1}-\cite{new-2023a} appeared for the first time in 
\cite{very-old-1}, \cite{very-old-2}. In these
works the mentioned idea is formulated more likely at the level of 
guess (without any satisfactory grounding), and as a result the 
papers \cite{very-old-1}, \cite{very-old-2} 
contain rather fuzzy formulations and a number of 
incorrect conclusions. 
Note that in \cite{very-old-1}, \cite{very-old-2} we 
used the trigonometric multiple Fourier series 
converging  
in the sense of norm in the space
$L_2([t, T]^k),$ $k=1, 2, 3.$
It should be noted that the results 
of \cite{very-old-1}, \cite{very-old-2} 
are correct for a sufficiently narrow particular case 
when the numbers $i_1,...,i_k$ are pairwise different,
$i_1,\ldots,i_k = 1,\ldots,m$ (see Theorem 1.1 below).
   
Usage of Fourier series with respect to the system of Legendre polynomials 
for approximation of iterated stochastic integrals took place for 
the first time in the publications of the author
\cite{old-art-1}-\cite{old-art-4}
(also see \cite{1}-\cite{new-art-1}).

The question about what integrals (It\^{o} or Stratonovich) are more 
suitable for expansions within the frames of distinguished direction 
of researches has turned out to be rather interesting and difficult.

On the one side, the results of Chapter 1 (see Theorems 1.1, 1.2, 1.16) 
conclusively 
demonstrate that the structure 
of iterated It\^{o} stochastic integrals is rather convenient for expansions into 
multiple series with respect to the system of standard Gaussian random 
variables regardless of the multiplicity $k$ of the iterated It\^{o} stochastic integral.
   
On the other side, the results of Chapter 2
\cite{6}-\cite{art-6}, \cite{art-8}, \cite{art-9},
\cite{arxiv-2}, \cite{arxiv-4}-\cite{arxiv-11}, \cite{arxiv-14},
\cite{arxiv-15}, \cite{arxiv-17}-\cite{arxiv-19}, \cite{arxiv-24},
\cite{new-art-1xxy}, 
\cite{new-art-1xxys},
\cite{old-art-1}-\cite{old-art-4}
convincingly demontrate
that the final formulas for 
expansions of iterated Stratonovich stochastic integrals of multiplicities 1 to 8
(the case of continuously differentiable weight functions and 
a complete orthonormal system of Legendre polynomials or 
trigonometric functions in $L_2([t, T])$) and
iterated Stratonovich stochastic integrals
of multiplicity $k,$ $k\in {\bf N}$ (the case of continuous weight functions 
and an arbitrary complete orthonormal system of functions in $L_2([t, T])$)
are more compact than their analogues for iterated
It\^{o} stochastic integrals.

\subsection{It\^{o} Stochastic Integral}

Let $(\Omega,{\rm F},{\sf P})$ be a complete probability
space and let $w(t,\omega):$ $[0, T]\times \Omega\rightarrow {\bf R}$
be the standard Wiener process
defined on the probability space $(\Omega,{\rm F},{\sf P}).$
Further, we will use the following notation:
$w(t,\omega)\stackrel{\rm def}{=}w_t$.

Let us consider the right-continous family of $\sigma$-algebras
$\left\{{\rm F}_t,\ t\in[0,T]\right\}$ defined
on the probability space $(\Omega,{\rm F},{\sf P})$ and
connected
with the Wiener process $w_t$ in such a way that

1.\ ${\rm F}_s\subset {\rm F}_t\subset {\rm F}$\ for
$s<t.$

2.\ The Wiener process $w_t$ is ${\rm F}_t$-measurable for all
$t\in[0,T].$

3.\ The process $w_{t+\Delta}-w_{t}$ for all
$t\ge 0,$ $\Delta>0$ is independent with
the events of $\sigma$-algebra
${\rm F}_{t}.$

Let us introduce the class ${\rm M}_2([0,T])$ of functions
$\xi:$ $[0,T]\times\Omega\rightarrow {\bf R},$ which satisfy the
conditions:

1. The function $\xi(t,\omega)$ is 
measurable
with respect to the pair of variables
$(t,\omega).$

2. The function $\xi(t,\omega)$ is ${\rm F}_t$-measurable 
for all $t\in[0,T]$ and $\xi(\tau,\omega)$ is independent 
with increments $w_{t+\Delta}-w_{t}$ 
for $t\ge \tau,$ $\Delta>0.$

3.\ The following relation is fulfilled
$$\int\limits_0^T{\sf M}\left\{\left(\xi(t,\omega)\right)^2\right\}dt
<\infty.
$$

4.\  ${\sf M}\left\{\left(\xi(t,\omega)\right)^2\right\}<\infty$
for all $t\in[0,T].$

For any partition 
$\tau_j^{(N)},$ $j=0, 1, \ldots, N$ of
the interval $[0,T]$ such that
\begin{equation}
\label{pr}
0=\tau_0^{(N)}<\tau_1^{(N)}<\ldots <\tau_N^{(N)}=T,\ \ \ \
\max\limits_{0\le j\le N-1}\left|\tau_{j+1}^{(N)}-\tau_j^{(N)}\right|\to 0\ \
\hbox{if}\ \ N\to \infty
\end{equation}
we will define the sequence of step functions 
$$
\xi^{(N)}(t,\omega)=\xi_j\left(\omega\right)\ \ \ 
\hbox{w.~p.~1}\ \ \ 
\hbox{for}\ \ \ t\in\left[\tau_j^{(N)},\tau_{j+1}^{(N)}\right),
$$
where $\xi^{(N)}(t,\omega)\in{\rm M}_2([0,T]),$ $j=0, 1,\ldots,N-1,$  $N=1, 2,\ldots$ 
Here and further, w.~p.~1 means with probability 1.

Let us define the It\^{o} stochastic integral for
$\xi(t,\omega)\in{\rm M}_2([0,T])$ as the 
following mean-square limit \cite{Gih1}, \cite{Gih} (also see \cite{Zapad-3})
\begin{equation}
\label{10.1}
~~~~~~\hbox{\vtop{\offinterlineskip\halign{
\hfil#\hfil\cr
{\rm l.i.m.}\cr
$\stackrel{}{{}_{N\to \infty}}$\cr
}} }\sum_{j=0}^{N-1}\xi^{(N)}\left(\tau_j^{(N)},\omega\right)
\left(w\left(\tau_{j+1}^{(N)},\omega\right)-
w\left(\tau_j^{(N)},\omega\right)\right)
\stackrel{\rm def}{=}\int\limits_0^T\xi_\tau dw_\tau,
\end{equation}
where $\xi^{(N)}(t,\omega)$ is any step function from the class ${\rm M}_2([0,T])$, 
which converges
to the function $\xi(t,\omega)$
in the following sense
\begin{equation}
\label{jjj}
\hbox{\vtop{\offinterlineskip\halign{
\hfil#\hfil\cr
{\rm lim}\cr
$\stackrel{}{{}_{N\to \infty}}$\cr
}} }\int\limits_0^T{\sf M}\left\{\left|
\xi^{(N)}(t,\omega)-\xi(t,\omega)\right|^2\right\}dt=0.
\end{equation}

Further, we will denote  
$\xi(\tau,\omega)$ as $\xi_{\tau}.$

It is well known \cite{Gih1} that the It\^{o} stochastic integral
exists as the limit (\ref{10.1}) and it does not depend on the selection 
of sequence 
$\xi^{(N)}(t,\omega)$. 
Furthermore, the It\^{o} stochastic integral
satisfies w.~p.~1
to the following properties \cite{Gih1}
$$
{\sf M}\left\{\int\limits_0^T
\xi_t dw_t\biggl|\biggr.{\rm F}_0\right\}=0,
$$
$$
{\sf M}\left\{\left|\int\limits_0^T
\xi_t dw_t\right|^2\biggl|\biggr.{\rm F}_0\right\}=
{\sf M}\left\{\int\limits_0^T\xi_t^2 dt\biggl|\biggr.{\rm F}_0\right\},
$$
$$
\int\limits_0^T(\alpha\xi_t+\beta\psi_t)dw_t=
\alpha\int\limits_0^T\xi_t dw_t+\beta
\int\limits_0^T\psi_t dw_t,
$$
where $\xi_t,$ $\phi_t\in {\rm M}_2([0, T]),$
$\alpha,\ \beta\in{\bf R}^1.$

Let us define the stochastic integral for 
$\xi_{\tau}\in{\rm M}_2([0,T])$ as the 
following mean-square limit
\begin{equation}
\label{10.1selo}
\hbox{\vtop{\offinterlineskip\halign{
\hfil#\hfil\cr
{\rm l.i.m.}\cr
$\stackrel{}{{}_{N\to \infty}}$\cr
}} }\sum_{j=0}^{N-1}\xi^{(N)}\left(\tau_j^{(N)},\omega\right)
\left(\tau_{j+1}^{(N)}-
\tau_j^{(N)}\right)\stackrel{\rm def}{=}\int\limits_{0}^T \xi_\tau d\tau,
\end{equation}
where $\xi^{(N)}(t,\omega)$ is any step function
from the class ${\rm M}_2([0,T])$,
which converges
in the sense (\ref{jjj}) 
to the function $\xi(t,\omega).$

\subsection{Theorem on Expansion of Iterated It\^{o} Stochastic Integrals of 
Multiplicity $k$ ($k\in{\bf N}$)}

Let $(\Omega,$ ${\rm F},$ ${\sf P})$ be a complete probability space, let 
$\{{\rm F}_t, t\in[0,T]\}$ be a nondecreasing right-continuous family of 
$\sigma$-algebras of ${\rm F},$
and let ${\bf w}_t$ be a standard $m$-dimensional Wiener 
stochastic process, which is
${\rm F}_t$-measurable for any $t\in[0, T].$ We assume that the components
${\bf w}_{t}^{(i)}$ $(i=1,\ldots,m)$ of this process are independent.

Let us consider 
the following iterated It\^{o} 
stochastic integrals
\begin{equation}
\label{ito}
J[\psi^{(k)}]_{T,t}=\int\limits_t^T\psi_k(t_k) \ldots \int\limits_t^{t_{2}}
\psi_1(t_1) d{\bf w}_{t_1}^{(i_1)}\ldots
d{\bf w}_{t_k}^{(i_k)},
\end{equation}
where $\psi_l(\tau)$ $(l=1,\ldots,k)$ are
nonrandom functions on $[t, T]$,
${\bf w}_{\tau}^{(i)}$ $(i=1,\ldots,m)$ are independent standard Wiener
processes,
${\bf w}_{\tau}^{(0)}=\tau,$ 
$i_1,\ldots,i_k=0, 1,\ldots,m.$

Let us consider the approach to expansion of the iterated
It\^{o} stochastic integrals (\ref{ito}) 
\cite{1}-\cite{new-2023a} (the so-called
method of generalized
multiple Fourier series). 
The idea of this method is as follows: 
the iterated It\^{o} stochastic 
integral (\ref{ito}) of multiplicity $k,$ $k\in{\bf N}$
is represented as the multiple stochastic 
integral from the certain discontinuous nonrandom function of $k$ variables 
defined on the hypercube $[t, T]^k.$ Here $[t, T]$ is the interval of 
integration of the iterated It\^{o} stochastic integral (\ref{ito}). 
Then, the mentioned
nonrandom function of $k$ variables
is expanded in the hypercube $[t, T]^k$ into the generalized 
multiple Fourier series converging 
in the mean-square sense
in the space 
$L_2([t,T]^k)$. After a number of nontrivial transformations we come 
to the 
mean-square converging expansion of the iterated It\^{o} stochastic 
integral (\ref{ito}) into the multiple 
series of products
of standard  Gaussian random 
variables. The coefficients of this 
series are the coefficients of 
generalized multiple Fourier series for the mentioned nonrandom function 
of $k$ variables, which can be calculated using the explicit formula 
regardless 
of the multiplicity $k$ of the
iterated It\^{o} stochastic integral (\ref{ito}).

Suppose that every $\psi_l(\tau)$ $(l=1,\ldots,k)$ is a continuous 
nonrandom
function on $[t, T]$ (we will also consider the case $\psi_1(\tau),\ldots,\psi_k(\tau)
\in L_2([t,T])$ in Sect.~1.11, 1.12). 
Define the following function on the hypercube $[t, T]^k$
\begin{equation}
\label{ppp}
K(t_1,\ldots,t_k)=
\left\{\begin{matrix}
\psi_1(t_1)\ldots \psi_k(t_k),\ &t_1<\ldots<t_k\cr\cr
0,\ &\hbox{\rm otherwise}
\end{matrix}
\right.\ \ \
=\ \  
\prod\limits_{l=1}^k
\psi_l(t_l)\ \prod\limits_{l=1}^{k-1}{\bf 1}_{\{t_l<t_{l+1}\}},
\end{equation}
where $t_1,\ldots,t_k\in [t, T]$ $(k\ge 2)$ and 
$K(t_1)\equiv\psi_1(t_1)$ for $t_1\in[t, T].$ Here 
${\bf 1}_A$ denotes the indicator of the set $A$.

Suppose that $\{\phi_j(x)\}_{j=0}^{\infty}$
is a complete orthonormal system of functions in 
the space $L_2([t, T])$.

The function $K(t_1,\ldots,t_k)$ is piecewise continuous in the 
hypercube $[t, T]^k.$
At this situation it is well known that the generalized 
multiple Fourier series 
of $K(t_1,\ldots,t_k)\in L_2([t, T]^k)$ is converging 
to $K(t_1,\ldots,t_k)$ in the hypercube $[t, T]^k$ in 
the mean-square sense, i.e.
\begin{equation}
\label{sos1z}
~~~~~\hbox{\vtop{\offinterlineskip\halign{
\hfil#\hfil\cr
{\rm lim}\cr
$\stackrel{}{{}_{p_1,\ldots,p_k\to \infty}}$\cr
}} }\Biggl\Vert
K(t_1,\ldots,t_k)-
\sum_{j_1=0}^{p_1}\ldots \sum_{j_k=0}^{p_k}
C_{j_k\ldots j_1}\prod_{l=1}^{k} \phi_{j_l}(t_l)\Biggr\Vert_{L_2([t, T]^k)}=0,
\end{equation}
where
\begin{equation}
\label{ppppa}
C_{j_k\ldots j_1}=\int\limits_{[t,T]^k}
K(t_1,\ldots,t_k)\prod_{l=1}^{k}\phi_{j_l}(t_l)dt_1\ldots dt_k
\end{equation}
is the Fourier coefficient, and
$$
\left\Vert f\right\Vert_{L_2([t, T]^k)}=\left(~\int\limits_{[t,T]^k}
f^2(t_1,\ldots,t_k)dt_1\ldots dt_k\right)^{1/2}.
$$

Consider the partition $\{\tau_j\}_{j=0}^N$ of $[t,T]$ such that
\begin{equation}
\label{1111}
t=\tau_0<\ldots <\tau_N=T,\ \ \
\Delta_N=
\hbox{\vtop{\offinterlineskip\halign{
\hfil#\hfil\cr
{\rm max}\cr
$\stackrel{}{{}_{0\le j\le N-1}}$\cr
}} }\Delta\tau_j\to 0\ \ \hbox{if}\ \ N\to \infty,\ \ \ 
\Delta\tau_j=\tau_{j+1}-\tau_j.
\end{equation}

{\bf Theorem 1.1}${}^2$ \cite{1} (2006) (also see \cite{2}-\cite{new-2023a}). 
{\it Suppose that
every $\psi_l(\tau)$ $(l=$ $1,\ldots, k)$ is a continuous 
nonrandom function on 
$[t, T]$ and
$\{\phi_j(x)\}_{j=0}^{\infty}$ is a complete orthonormal system  
of continuous functions in the space $L_2([t,T]).$ 
Then
$$
J[\psi^{(k)}]_{T,t} =
\hbox{\vtop{\offinterlineskip\halign{
\hfil#\hfil\cr
{\rm l.i.m.}\cr
$\stackrel{}{{}_{p_1,\ldots,p_k\to \infty}}$\cr
}} }\sum_{j_1=0}^{p_1}\ldots\sum_{j_k=0}^{p_k}
C_{j_k\ldots j_1}\Biggl(
\prod_{l=1}^k\zeta_{j_l}^{(i_l)} -
\Biggr.
$$
\begin{equation}
\label{tyyy}
-\Biggl.
\hbox{\vtop{\offinterlineskip\halign{
\hfil#\hfil\cr
{\rm l.i.m.}\cr
$\stackrel{}{{}_{N\to \infty}}$\cr
}} }\sum_{(l_1,\ldots,l_k)\in {\rm G}_k}
\phi_{j_{1}}(\tau_{l_1})
\Delta{\bf w}_{\tau_{l_1}}^{(i_1)}\ldots
\phi_{j_{k}}(\tau_{l_k})
\Delta{\bf w}_{\tau_{l_k}}^{(i_k)}\Biggr),
\end{equation}
where
$$
{\rm G}_k={\rm H}_k\backslash{\rm L}_k,\ \ \
{\rm H}_k=\bigl\{(l_1,\ldots,l_k):\ l_1,\ldots,l_k=0,\ 1,\ldots,N-1\bigr\},
$$
$$
{\rm L}_k=\bigl\{(l_1,\ldots,l_k):\ l_1,\ldots,l_k=0,\ 1,\ldots,N-1;\
l_g\ne l_r\ (g\ne r);\ g, r=1,\ldots,k\bigr\},
$$

\noindent
${\rm l.i.m.}$ is a limit in the mean-square sense$,$
$i_1,\ldots,i_k=0,1,\ldots,m,$ 
\begin{equation}
\label{rr23}
\zeta_{j}^{(i)}=
\int\limits_t^T \phi_{j}(s) d{\bf w}_s^{(i)}
\end{equation} 
are independent standard Gaussian random variables
for various
$i$ or $j$ {\rm(}in the case when $i\ne 0${\rm),}
$C_{j_k\ldots j_1}$ is the Fourier coefficient {\rm(\ref{ppppa}),}
$\Delta{\bf w}_{\tau_{j}}^{(i)}=
{\bf w}_{\tau_{j+1}}^{(i)}-{\bf w}_{\tau_{j}}^{(i)}$
$(i=0, 1,\ldots,m),$\
$\left\{\tau_{j}\right\}_{j=0}^{N}$ is a partition of
$[t,T],$ which satisfies the condition {\rm (\ref{1111})}.}

\footnotetext[2]{Theorem 1.1 will be generalized to the case of an arbitrary 
complete ortho\-nor\-mal system of functions $\{\phi_j(x)\}_{j=0}^{\infty}$ 
in the space $L_2([t, T])$
and $\psi_1(\tau),$ $\ldots,\psi_k(\tau) \in L_2([t, T])$ in Sect.~1.11
(see Theorem~1.16). 
Theorem~1.1 marked the beginning of a systematic
study of the problem of strong approximation 
of iterated It\^{o} and Stratonovich stochastic integrals 
that have been most fully studied to date in this book.}

{\bf Proof.}\ 
At first, let us prove preparatory lemmas.

{\bf Lemma 1.1.} {\it Suppose that
every $\psi_l(\tau)$ $(l=1,\ldots, k)$ is a continuous nonrandom
function on 
$[t, T]$. Then
\begin{equation}
\label{30.30}
J[\psi^{(k)}]_{T,t}=
\hbox{\vtop{\offinterlineskip\halign{
\hfil#\hfil\cr
{\rm l.i.m.}\cr
$\stackrel{}{{}_{N\to \infty}}$\cr
}} }
\sum_{j_k=0}^{N-1}
\ldots \sum_{j_1=0}^{j_{2}-1}
\prod_{l=1}^k \psi_l(\tau_{j_l})\Delta{\bf w}_{\tau_
{j_l}}^{(i_l)}\ \ \ \hbox{\rm w.~p.~1},
\end{equation}
where $\Delta{\bf w}_{\tau_{j}}^{(i)}=
{\bf w}_{\tau_{j+1}}^{(i)}-{\bf w}_{\tau_{j}}^{(i)}$
$(i=0, 1,\ldots,m)$,
$\left\{\tau_{j}\right\}_{j=0}^{N}$ is a partition 
of the interval $[t,T]$ satisfying the condition {\rm (\ref{1111})}.
}

{\bf Proof.}\ It is easy to notice that using the 
property of stochastic integrals additivity, we can write 
\begin{equation}
\label{toto}
J[\psi^{(k)}]_{T,t}=
\sum_{j_k=0}^{N-1}\ldots
\sum_{j_{1}=0}^{j_{2}-1}\prod_{l=1}^{k}
J[\psi_l]_{\tau_{j_l+1},\tau_{j_l}}+
\varepsilon_N\ \ \ \ \ \hbox{w.~p.~1},
\end{equation}
where
$$
J[\psi_l]_{s,\theta}=\int\limits_{\theta}^s\psi_l(\tau)d{\bf w}_{\tau}
^{(i_l)}
$$
and
$$
\varepsilon_N = 
\sum_{j_k=0}^{N-1}\int\limits_{\tau_{j_k}}^{\tau_{j_k+1}}
\psi_k(s)\int\limits_{\tau_{j_k}}^{s}
\psi_{k-1}(\tau)J[\psi^{(k-2)}]_{\tau,t}d{\bf w}_\tau^{(i_{k-1})}
d{\bf w}_s^{(i_k)} +
$$
$$
+ \sum_{r=1}^{k-3}
G[\psi_{k-r+1}^{(k)}]_{N}\times
$$
$$
\times
\sum_{j_{k-r}=0}^{j_{k-r+1}-1}\int\limits_{\tau_{j_{k-r}}}^{\tau_{j_{k-r}+1}}
\psi_{k-r}(s)\int\limits_{\tau_{j_{k-r}}}^{s}
\psi_{k-r-1}(\tau)J[\psi^{(k-r-2)}]_{\tau,t}
d{\bf w}_\tau^{(i_{k-r-1})}
d{\bf w}_s^{(i_{k-r})}+
$$
$$
+ G[\psi_3^{(k)}]_{N}
\sum_{j_{2}=0}^{j_{3}-1}J[\psi^{(2)}]_{\tau_{j_{2}+1},
\tau_{j_{2}}},
$$
where
$$
G[\psi_m^{(k)}]_{N}=\sum_{j_k=0}^{N-1}\sum_{j_{k-1}=0}^{j_k-1}
\ldots \sum_{j_{m}=0}^{j_{m+1}-1}
\prod_{l=m}^{k}J[\psi_l]_{\tau_{j_l+1},\tau_{j_l}},
$$
$$
(\psi_m,\psi_{m+1},\ldots,\psi_{k})\stackrel{\rm def}{=}\psi_m^{(k)},\ \ \
(\psi_1,\ldots,\psi_{k})=\psi_1^{(k)}
\stackrel{\rm def}{=}\psi^{(k)}.
$$

\vspace{3mm}

Using the standard estimates (\ref{99.010}),
(\ref{99.010a}) (see below) for the moments of stochastic 
integrals, we obtain w. p. 1
\begin{equation}
\label{999.0001}
\hbox{\vtop{\offinterlineskip\halign{
\hfil#\hfil\cr
{\rm l.i.m.}\cr
$\stackrel{}{{}_{N\to \infty}}$\cr
}} }\varepsilon_N =0.
\end{equation}

Comparing (\ref{toto}) and (\ref{999.0001}), we get
\begin{equation}
\label{toto1}
J[\psi^{(k)}]_{T,t}=
\hbox{\vtop{\offinterlineskip\halign{
\hfil#\hfil\cr
{\rm l.i.m.}\cr
$\stackrel{}{{}_{N\to \infty}}
$\cr
}} }
\sum_{j_k=0}^{N-1}\ldots
\sum_{j_{1}=0}^{j_{2}-1}\prod_{l=1}^{k}
J[\psi_l]_{\tau_{j_l+1},\tau_{j_l}}\ \ \ \hbox{w.~p.~1}.
\end{equation}

Let us rewrite $J[\psi_l]_{\tau_{j_l+1},\tau_{j_l}}$ in the form 
$$
J[\psi_l]_{\tau_{j_l+1},\tau_{j_l}}=
\psi_l(\tau_{j_l})\Delta{\bf w}_{\tau_{j_l}}^{(i_l)}+
\int\limits_{\tau_{j_l}}^{\tau_{j_l+1}}
(\psi_l(\tau)-\psi_l(\tau_{j_l}))d{\bf w}_{\tau}^{(i_l)}
$$
and substitute it into
(\ref{toto1}).
Then, due to the moment properties of stochastic integrals and
continuity (which means uniform continuity) of the functions 
$\psi_l(s)$ ($l=1,\ldots,k$)
it is easy to see that the 
prelimit
expression on the right-hand side of (\ref{toto1}) is a sum of 
the prelimit
expression on the right-hand side of (\ref{30.30}) and the value which 
tends to zero in the mean-square sense if 
$N\to\infty.$ Lemma 1.1 is proved.

{\bf Remark 1.1.} {\it It is easy to see that if
$\Delta{\bf w}_{\tau_{j_l}}^{(i_l)}$ in {\rm (\ref{30.30})}
for some $l\in\{1,\ldots,k\}$ is replaced with 
$\left(\Delta{\bf w}_{\tau_{j_l}}^{(i_l)}\right)^p$ $(p=2,$
$i_l\ne 0),$ then
the differential $d{\bf w}_{t_{l}}^{(i_l)}$
in the integral $J[\psi^{(k)}]_{T,t}$
will be replaced with $dt_l$.
If $p=3, 4,\ldots,$ then the
right-hand side 
of the formula {\rm (\ref{30.30})}
will become zero w.~p.~{\rm 1}.
If we replace $\Delta{\bf w}_{\tau_{j_l}}^{(i_l)}$ in {\rm (\ref{30.30})}
for some $l\in\{1,\ldots,k\}$
with $\left(\Delta \tau_{j_l}\right)^p$ $(p=2, 3,\ldots),$
then the right-hand side of the formula
{\rm (\ref{30.30})} also 
will be equal to zero w. p. {\rm 1}.}

Let us define the following
multiple stochastic integral
\begin{equation}
\label{30.34}
\hbox{\vtop{\offinterlineskip\halign{
\hfil#\hfil\cr
{\rm l.i.m.}\cr
$\stackrel{}{{}_{N\to \infty}}$\cr
}} }\sum_{j_1,\ldots,j_k=0}^{N-1}
\Phi\left(\tau_{j_1},\ldots,\tau_{j_k}\right)
\prod\limits_{l=1}^k\Delta{\bf w}_{\tau_{j_l}}^{(i_l)}
\stackrel{\rm def}{=}J[\Phi]_{T,t}^{(k)},
\end{equation}
where $\Phi(t_1,\ldots,t_k): [t, T]^k\to{\bf R},$
$\Phi(t_1,\ldots,t_k)\in C([t, T]^k)$, i.e. $\Phi(t_1,\ldots,t_k)$
is a continuous nonrandom function in $[t, T]^k.$

Denote
\begin{equation}
\label{dom1}
D_k=\{(t_1,\ldots,t_k):\ t\le t_1<\ldots <t_k\le T\}.
\end{equation}

We will use the same symbol $D_k$ to denote the open and closed 
domains corresponding to the domain $D_k$ defined by (\ref{dom1}).
However, we always specify what domain we consider (open or closed).

Also we will write $\Phi(t_1,\ldots,t_k)\in C(D_k)$
if 
$\Phi(t_1,\ldots,t_k)$ is a continuous nonrandom function of $k$ variables
in the closed domain $D_k$. 

Let us consider the iterated It\^{o} stochastic integral
\begin{equation}
\label{rrr29}
I[\Phi]_{T,t}^{(k)}\stackrel{\rm def}{=}
\int\limits_t^T\ldots \int\limits_t^{t_2}
\Phi(t_1,\ldots,t_k)d{\bf w}_{t_1}^{(i_1)}\ldots
d{\bf w}_{t_k}^{(i_k)},
\end{equation}
where $\Phi(t_1,\ldots,t_k)\in C(D_k).$

Using the arguments which similar to the arguments used in the 
proof of Lemma 1.1
it is easy to demonstrate that if
$\Phi(t_1,\ldots,t_k)\in C(D_k),$ then the following equality is fulfilled
\begin{equation}
\label{30.52}
~~~~~~~~ I[\Phi]_{T,t}^{(k)}=\hbox{\vtop{\offinterlineskip\halign{
\hfil#\hfil\cr
{\rm l.i.m.}\cr
$\stackrel{}{{}_{N\to \infty}}$\cr
}} }
\sum_{j_k=0}^{N-1}
\ldots \sum_{j_1=0}^{j_{2}-1}
\Phi(\tau_{j_1},\ldots,\tau_{j_k})
\prod\limits_{l=1}^k\Delta {\bf w}_{\tau_{j_l}}^{(i_l)}\ \ \ \hbox{w.~p.~1}.
\end{equation}

In order to explain this, let us check the correctness of the equality 
(\ref{30.52}) when $k=3$.
For definiteness we will suppose that 
$i_1,i_2,i_3=1,\ldots,m.$ We have
$$
I[\Phi]_{T,t}^{(3)}\stackrel{\rm def}{=}
\int\limits_t^T\int\limits_t^{t_3}\int\limits_t^{t_2}
\Phi(t_1,t_2,t_3)d{\bf w}_{t_1}^{(i_1)}d{\bf w}_{t_2}^{(i_2)}
d{\bf w}_{t_3}^{(i_3)}=
$$
$$
=\hbox{\vtop{\offinterlineskip\halign{
\hfil#\hfil\cr
{\rm l.i.m.}\cr
$\stackrel{}{{}_{N\to \infty}}$\cr
}} }
\sum_{j_3=0}^{N-1}
\int\limits_{t}^{\tau_{j_3}}\int\limits_t^{t_2}
\Phi(t_1,t_2,\tau_{j_3})d{\bf w}_{t_1}^{(i_1)}d{\bf w}_{t_2}^{(i_2)}
\Delta{\bf w}_{\tau_{j_3}}^{(i_3)}=
$$
$$
=\hbox{\vtop{\offinterlineskip\halign{
\hfil#\hfil\cr
{\rm l.i.m.}\cr
$\stackrel{}{{}_{N\to \infty}}$\cr
}} }
\sum_{j_3=0}^{N-1}\sum_{j_2=0}^{j_3-1}
\int\limits_{\tau_{j_2}}^{\tau_{j_2+1}}\int\limits_t^{t_2}
\Phi(t_1,t_2,\tau_{j_3})d{\bf w}_{t_1}^{(i_1)}d{\bf w}_{t_2}^{(i_2)}
\Delta{\bf w}_{\tau_{j_3}}^{(i_3)}=
$$
$$
=\hbox{\vtop{\offinterlineskip\halign{
\hfil#\hfil\cr
{\rm l.i.m.}\cr
$\stackrel{}{{}_{N\to \infty}}$\cr
}} }
\sum_{j_3=0}^{N-1}\sum_{j_2=0}^{j_3-1}
\int\limits_{\tau_{j_2}}^{\tau_{j_2+1}}
\left(\ \int\limits_t^{\tau_{j_2}}\
+\ \int\limits_{\tau_{j_2}}^{t_2}\ \right)
\Phi(t_1,t_2,\tau_{j_3})d{\bf w}_{t_1}^{(i_1)}d{\bf w}_{t_2}^{(i_2)}
\Delta{\bf w}_{\tau_{j_3}}^{(i_3)}=
$$
$$
=\hbox{\vtop{\offinterlineskip\halign{
\hfil#\hfil\cr
{\rm l.i.m.}\cr
$\stackrel{}{{}_{N\to \infty}}$\cr
}} }
\sum_{j_3=0}^{N-1}\sum_{j_2=0}^{j_3-1}\sum_{j_1=0}^{j_2-1}
\int\limits_{\tau_{j_2}}^{\tau_{j_2+1}}\int\limits_{\tau_{j_1}}^{\tau_{j_1+1}}
\Phi(t_1,t_2,\tau_{j_3})d{\bf w}_{t_1}^{(i_1)}d{\bf w}_{t_2}^{(i_2)}
\Delta{\bf w}_{\tau_{j_3}}^{(i_3)}+
$$
\begin{equation}
\label{44444.25}
+\hbox{\vtop{\offinterlineskip\halign{
\hfil#\hfil\cr
{\rm l.i.m.}\cr
$\stackrel{}{{}_{N\to \infty}}$\cr
}} }
\sum_{j_3=0}^{N-1}\sum_{j_2=0}^{j_3-1}
\int\limits_{\tau_{j_2}}^{\tau_{j_2+1}}\int\limits_{\tau_{j_2}}^{t_2}
\Phi(t_1,t_2,\tau_{j_3})d{\bf w}_{t_1}^{(i_1)}d{\bf w}_{t_2}^{(i_2)}
\Delta{\bf w}_{\tau_{j_3}}^{(i_3)}.
\end{equation}

Let us demonstrate that the second limit on the right-hand side 
of (\ref{44444.25}) 
equals to zero.

Actually, for the second moment of its 
prelimit
expression we get
$$
\sum_{j_3=0}^{N-1}\sum_{j_2=0}^{j_3-1}
\int\limits_{\tau_{j_2}}^{\tau_{j_2+1}}\int\limits_{\tau_{j_2}}^{t_2}
\Phi^2(t_1,t_2,\tau_{j_3})dt_1 dt_2
\Delta\tau_{j_3}
\le M^2  \sum_{j_3=0}^{N-1}\sum_{j_2=0}^{j_3-1}
\frac{1}{2}\left(\Delta\tau_{j_2}\right)^2\Delta\tau_{j_3}\to 0
$$
when $N\to\infty.$
Here $M$ is a constant, which restricts the module of the
function
$\Phi(t_1,t_2,t_3)$ due to its continuity, $\Delta\tau_j=
\tau_{j+1}-\tau_j.$

Considering the obtained conclusions, we have
$$
I[\Phi]_{T,t}^{(3)}\stackrel{\rm def}{=}
\int\limits_t^T\int\limits_t^{t_3}\int\limits_t^{t_2}
\Phi(t_1,t_2,t_3)d{\bf w}_{t_1}^{(i_1)}d{\bf w}_{t_2}^{(i_2)}
d{\bf w}_{t_3}^{(i_3)}=
$$
$$
=\hbox{\vtop{\offinterlineskip\halign{
\hfil#\hfil\cr
{\rm l.i.m.}\cr
$\stackrel{}{{}_{N\to \infty}}$\cr
}} }
\sum_{j_3=0}^{N-1}\sum_{j_2=0}^{j_3-1}\sum_{j_1=0}^{j_2-1}
\int\limits_{\tau_{j_2}}^{\tau_{j_2+1}}\int\limits_{\tau_{j_1}}^{\tau_{j_1+1}}
\Phi(t_1,t_2,\tau_{j_3})d{\bf w}_{t_1}^{(i_1)}d{\bf w}_{t_2}^{(i_2)}
\Delta{\bf w}_{\tau_{j_3}}^{(i_3)}=
$$
$$
=\hbox{\vtop{\offinterlineskip\halign{
\hfil#\hfil\cr
{\rm l.i.m.}\cr
$\stackrel{}{{}_{N\to \infty}}$\cr
}} }
\sum_{j_3=0}^{N-1}\sum_{j_2=0}^{j_3-1}\sum_{j_1=0}^{j_2-1}
\int\limits_{\tau_{j_2}}^{\tau_{j_2+1}}\int\limits_{\tau_{j_1}}^{\tau_{j_1+1}}
\left(\Phi(t_1,t_2,\tau_{j_3})-\Phi(t_1,\tau_{j_2},\tau_{j_3})\right)
d{\bf w}_{t_1}^{(i_1)}d{\bf w}_{t_2}^{(i_2)}
\Delta{\bf w}_{\tau_{j_3}}^{(i_3)}+
$$
$$
+\hbox{\vtop{\offinterlineskip\halign{
\hfil#\hfil\cr
{\rm l.i.m.}\cr
$\stackrel{}{{}_{N\to \infty}}$\cr
}} }
\sum_{j_3=0}^{N-1}\sum_{j_2=0}^{j_3-1}\sum_{j_1=0}^{j_2-1}
\int\limits_{\tau_{j_2}}^{\tau_{j_2+1}}\int\limits_{\tau_{j_1}}^{\tau_{j_1+1}}
\left(\Phi(t_1,\tau_{j_2},\tau_{j_3})-
\Phi(\tau_{j_1},\tau_{j_2},\tau_{j_3})\right)
d{\bf w}_{t_1}^{(i_1)}d{\bf w}_{t_2}^{(i_2)}
\Delta{\bf w}_{\tau_{j_3}}^{(i_3)}+
$$
\begin{equation}
\label{4444.1}
+\hbox{\vtop{\offinterlineskip\halign{
\hfil#\hfil\cr
{\rm l.i.m.}\cr
$\stackrel{}{{}_{N\to \infty}}$\cr
}} }
\sum_{j_3=0}^{N-1}\sum_{j_2=0}^{j_3-1}\sum_{j_1=0}^{j_2-1}
\Phi(\tau_{j_1},\tau_{j_2},\tau_{j_3})
\Delta{\bf w}_{\tau_{j_1}}^{(i_1)}
\Delta{\bf w}_{\tau_{j_2}}^{(i_2)}
\Delta{\bf w}_{\tau_{j_3}}^{(i_3)}.
\end{equation}

\vspace{1mm}

In order to get the sought result, we just have to demonstrate that 
the first
two limits on the right-hand side of (\ref{4444.1}) equal to zero. 
Let us prove 
that the first one of them equals to zero (proof for the second limit 
is similar).
   
The second moment of prelimit expression of the first limit on the 
right-hand side of (\ref{4444.1}) equals to the following expression
\begin{equation}
\label{4444.01}
~~~~~~~~~~ \sum_{j_3=0}^{N-1}\sum_{j_2=0}^{j_3-1}\sum_{j_1=0}^{j_2-1}
\int\limits_{\tau_{j_2}}^{\tau_{j_2+1}}\int\limits_{\tau_{j_1}}^{\tau_{j_1+1}}
\left(\Phi(t_1,t_2,\tau_{j_3})-\Phi(t_1,\tau_{j_2},\tau_{j_3})\right)^2
dt_1 dt_2
\Delta\tau_{j_3}.
\end{equation}

Since the function $\Phi(t_1,t_2,t_3)$ is continuous in 
the closed bo\-un\-ded domain
$D_3,$ then
it is uniformly continuous in this domain. Therefore, if the 
distance between two points of the domain $D_3$ is less than 
$\delta(\varepsilon)$ ($\delta(\varepsilon)>0$ 
exists
for any $\varepsilon>0$ and it does not depend 
on mentioned points), then the cor\-res\-pond\-ing oscillation of the function 
$\Phi(t_1,t_2,t_3)$ for these two points of the domain $D_3$ is less than
$\varepsilon.$

If we assume that $\Delta\tau_j<\delta(\varepsilon)$ ($j=0, 1,\ldots,N-1$),
then the distance between points 
$(t_1,t_2,\tau_{j_3})$,\ $(t_1,\tau_{j_2},\tau_{j_3})$
is obviously less than $\delta(\varepsilon).$ In this case 
$$
|\Phi(t_1,t_2,\tau_{j_3})-\Phi(t_1,\tau_{j_2},\tau_{j_3})|<\varepsilon.
$$

Consequently, when $\Delta\tau_j<\delta(\varepsilon)$ ($j=0,\ 1,\ldots,N-1$)
the expression (\ref{4444.01})  
is estimated by the following value
$$
\varepsilon^2
\sum_{j_3=0}^{N-1}\sum_{j_2=0}^{j_3-1}\sum_{j_1=0}^{j_2-1}
\Delta\tau_{j_1}\Delta\tau_{j_2}\Delta\tau_{j_3}<
\varepsilon^2\frac{(T-t)^3}{6}.
$$ 

Therefore, the first limit on the right-hand side 
of (\ref{4444.1}) equals to zero.
Similarly, we can prove that the second limit on the right-hand
side
of (\ref{4444.1}) equals to zero.

Consequently, the equality (\ref{30.52}) is proved for $k=3$. 
The cases $k=2$ and $k>3$ 
are analyzed absolutely similarly.

It is necessary to note that the proof of 
correctness of (\ref{30.52}) 
is similar when the nonrandom function $\Phi(t_1,\ldots,t_k)$ is 
continuous in 
the open domain $D_k$ and bounded at its boundary.

Let us consider the following multiple 
stochastic integral
\begin{equation}
\label{mult11}
~~~~~~\hbox{\vtop{\offinterlineskip\halign{
\hfil#\hfil\cr
{\rm l.i.m.}\cr
$\stackrel{}{{}_{N\to \infty}}$\cr
}} }
\sum\limits_{\stackrel{j_1,\ldots,j_k=0}{{}_{j_q\ne j_r;\ q\ne r;\ 
q, r=1,\ldots, k}}}^{N-1}
\Phi\left(\tau_{j_1},\ldots,\tau_{j_k}\right)
\prod\limits_{l=1}^k
\Delta{\bf w}_{\tau_{j_l}}^{(i_l)}
\stackrel{\rm def}{=}J'[\Phi]_{T,t}^{(k)},
\end{equation}
where $\Phi(t_1,\ldots,t_k): [t, T]^k\to{\bf R}$ is the same function as in
(\ref{30.34}).

According to (\ref{30.52}), we get the following equality
\begin{equation}
\label{pobeda}
~~~~~J'[\Phi]_{T,t}^{(k)}=
\int\limits_t^T\ldots \int\limits_t^{t_2}
\sum\limits_{(t_1,\ldots,t_k)}\biggl(
\Phi(t_1,\ldots,t_k)
d{\bf w}_{t_1}^{(i_1)}\ldots
d{\bf w}_{t_k}^{(i_k)}\biggr)\ \ \ \hbox{w.~p.~1},
\end{equation}
where
$$
\sum\limits_{(t_1,\ldots,t_k)}
$$ 
means the sum with respect to all
possible permutations
$(t_1,\ldots,t_k).$ 
At the same time permutations $(t_1,\ldots,t_k)$
when summing 
are performed in (\ref{pobeda}) only in the expression, which
is enclosed in pa\-ren\-the\-ses. Moreover,
the nonrandom function $\Phi(t_1,\ldots,t_k)$ is assumed 
to be continuous in the 
cor\-res\-pond\-ing closed domains of integration.
The case when the nonrandom function $\Phi(t_1,\ldots,t_k)$ is 
continuous in the open domains of integration and bounded at 
their boundaries is also possible.

It is not difficult to see that (\ref{pobeda})
can be rewritten in the form
\begin{equation}
\label{s2s}
~~~~~~~ J'[\Phi]_{T,t}^{(k)}=\sum_{(t_1,\ldots,t_k)}
\int\limits_{t}^{T}
\ldots
\int\limits_{t}^{t_2}
\Phi(t_1,\ldots,t_k)d{\bf w}_{t_1}^{(i_1)}
\ldots
d{\bf w}_{t_k}^{(i_k)}\ \ \ \hbox{w.~p.~1},
\end{equation}
where permutations $(t_1,\ldots,t_k)$ when summing are 
performed only in the values
$d{\bf w}_{t_1}^{(i_1)}
\ldots $
$d{\bf w}_{t_k}^{(i_k)}.$ At the same time the indices near 
upper 
limits of integration in the iterated stochastic integrals are changed 
correspondently and if $t_r$ swapped with $t_q$ in the  
permutation $(t_1,\ldots,t_k)$, then $i_r$ swapped with $i_q$ in 
the permutation $(i_1,\ldots,i_k)$.

{\bf Lemma 1.2.}\ {\it Suppose that $\Phi(t_1,\ldots,t_k)\in C(D_k)$ 
or $\Phi(t_1,\ldots,t_k)$ 
is a continuous nonrandom function in the open domain $D_k$ and bounded at its boundary.
Then
$$
{\sf M}\left\{\biggl|I[\Phi]_{T,t}^{(k)}\biggr|^{2}\right\}
\le C_{k}
\int\limits_t^T\ldots \int\limits_t^{t_2}
\Phi^{2}(t_1,\ldots,t_k)dt_1\ldots dt_k,\ \ \
C_{k}<\infty,
$$
where $I[\Phi]_{T,t}^{(k)}$ is defined by the formula {\rm (\ref{rrr29})}.}

{\bf Proof.}\
Using standard properties and estimates of stochastic integrals 
for $\xi_{\tau}\in{\rm M}_2([t,T])$, we have \cite{Gih}
\begin{equation}
\label{99.010}
{\sf M}\left\{\left|\int\limits_{t}^T \xi_\tau
dw_\tau\right|^{2}\right\} =
\int\limits_{t}^T {\sf M}\{|\xi_\tau|^{2}\}d\tau,
\end{equation}
\begin{equation}
\label{99.010a}
{\sf M}\left\{\left|\int\limits_{t}^T \xi_\tau
d\tau\right|^{2}\right\} \le (T-t)
\int\limits_{t}^T {\sf M}\{|\xi_\tau|^{2}\}d\tau.
\end{equation}

Let us denote
$$
\xi[\Phi]_{t_{l+1},\ldots,t_k,t}^{(l)}=
\int\limits_t^{t_{l+1}}\ldots \int\limits_t^{t_2}
\Phi(t_1,\ldots,t_k)
d{\bf w}_{t_1}^{(i_1)}\ldots
d{\bf w}_{t_{l}}^{(i_{l})},
$$ 
where $l=1,\ldots,$ $k-1$ and
$\xi[\Phi]_{t_{1},\ldots,t_k,t}^{(0)}\stackrel{\rm def}{=}
\Phi(t_1,\ldots,t_k).$

By induction it is easy to demonstrate that
$$
\xi[\Phi]_{t_{l+1},\ldots,t_k,t}^{(l)}\in{\rm M}_2([t,T])
$$
with respect to the variable $t_{l+1}.$
Further, using the estimates (\ref{99.010}), (\ref{99.010a})
repeatedly we obtain the statement of Lemma 1.2.

It is not difficult to see that in the case $i_1,\ldots,i_k=1,\dots,m$
from the proof of Lemma 1.2 we obtain
\begin{equation}
\label{dobav1}
{\sf M}\left\{\biggl|I[\Phi]_{T,t}^{(k)}\biggr|^{2}\right\}
=
\int\limits_t^T\ldots \int\limits_t^{t_2}
\Phi^{2}(t_1,\ldots,t_k)dt_1\ldots dt_k.
\end{equation}

{\bf Lemma 1.3.}\ {\it Suppose that every $\varphi_l(s)$
$(l=1,\ldots,k)$ is a continuous nonrandom function on $[t, T]$.
Then
\begin{equation}
\label{30.39}
\prod_{l=1}^k 
J[\varphi_l]_{T,t}=J[\Phi]_{T,t}^{(k)}\ \ \ \hbox{\rm w.~p.~1},
\end{equation}
where 
$$
J[\varphi_l]_{T,t}
=\int\limits_t^T \varphi_l(s) d{\bf w}_{s}^{(i_l)},\ \ \ \ 
\Phi(t_1,\ldots,t_k)=\prod\limits_{l=1}^k\varphi_l(t_l),
$$
and the integral $J[\Phi]_{T,t}^{(k)}$ 
is defined
by the equality
{\rm (\ref{30.34})}.
}

{\bf Proof.}\
Let at first $i_l\ne 0,\ l=1,\ldots,k.$
Let us denote
$$
J[\varphi_l]_{N}\stackrel{\rm def}{=}\sum\limits_{j=0}^{N-1}
\varphi_l(\tau_j)\Delta{\bf w}_{\tau_j}^{(i_l)}.
$$

Since
$$
\prod_{l=1}^k J[\varphi_l]_{N}-\prod_{l=1}^k J[\varphi_l]_{T,t}
=
$$
\begin{equation}
\label{df2}
~~~=\sum_{l=1}^k \left(\prod_{g=1}^{l-1} J[\varphi_g]_{T,t}\right)
\biggl(J[\varphi_l]_{N}-J[\varphi_l]_{T,t}
\biggr)\left(\prod_{g=l+1}^k J[\varphi_g]_{N}\right),
\end{equation}

\vspace{2mm}
\noindent
then due to the Minkowski inequality and the inequality 
of Cauchy--Bu\-ny\-a\-kov\-sky we obtain

\vspace{-2mm}
$$
\left({\sf M}\left\{\left|\prod_{l=1}^k J[\varphi_l]_{N}
-\prod_{l=1}^k J[\varphi_l]_{T,t}\right|^2
\right\}\right)^{1/2}\le 
$$

\newpage
\noindent
\begin{equation}
\label{30.42}
\le C_k
\sum_{l=1}^k
\left({\sf M}\left\{
\biggl|J[\varphi_l]_{N}-J[\varphi_l]_{T,t}\biggr|^4\right\}\right)
^{1/4},
\end{equation}
where $C_k$ is a constant.

Note that
$$
J[\varphi_l]_{N}-J[\varphi_l]_{T,t}=\sum\limits_{j=0}^{N-1}
J[\Delta\varphi_l]_{\tau_{j+1},\tau_j},
$$
$$
J[\Delta\varphi_l]_{\tau_{j+1},\tau_j}
=\int\limits_{\tau_j}^{\tau_{j+1}}\left(
\varphi_l(\tau_j)-\varphi_l(s)\right)d{\bf w}_{s}^{(i_l)}.
$$

Since $J[\Delta\varphi_l]_{\tau_{j+1},\tau_j}$
are independent for various $j,$ then \cite{Scor}

\vspace{2mm}
$$
{\sf M}\left\{\left|\sum_{j=0}^{N-1}J[\Delta\varphi_l]_{\tau_{j+1},
\tau_j}\right|^4
\right\}=
\sum_{j=0}^{N-1}{\sf M}\left\{\biggl|J[\Delta\varphi_l]_{\tau_{j+1},
\tau_j}\biggr|^4
\right\}+ 
$$

\vspace{2mm}
\begin{equation}
\label{30.43}
+6 \sum_{j=0}^{N-1}{\sf M}
\left\{\biggl|J[\Delta\varphi_l]_{\tau_{j+1},\tau_j}\biggr|^2
\right\}
\sum_{q=0}^{j-1}{\sf M}\left\{\biggl|
J[\Delta\varphi_l]_{\tau_{q+1},\tau_q}\biggr|^2
\right\}.
\end{equation}

\vspace{6mm}

Moreover, since
$J[\Delta\varphi_l]_{\tau_{j+1},\tau_j}$ is a Gaussian random variable,
we have
$$
{\sf M}\left\{\biggl|J[\Delta\varphi_l]_{\tau_{j+1},\tau_j}\biggr|^2\right\}=
\int\limits_{\tau_j}^{\tau_{j+1}}(\varphi_l(\tau_j)-\varphi_l(s))^2ds,
$$
$$
{\sf M}\left\{\biggl|J[\Delta\varphi_l]_{\tau_{j+1},\tau_j}\biggr|^4\right\}=
3\left(\int\limits_{\tau_j}^{\tau_{j+1}}(\varphi_l(\tau_j)-\varphi_l(s))^2ds
\right)^2.
$$

Using these relations and continuity (which means uniform 
continuity) of the functions $\varphi_l(s),$ we get 
$$
{\sf M}\left\{\left|\sum_{j=0}^{N-1}J[\Delta\varphi_l]_{\tau_{j+1},
\tau_j}\right|^4
\right\}\le \varepsilon^4\left(
3 \sum_{j=0}^{N-1}(\Delta\tau_{j})^2+
6 \sum_{j=0}^{N-1}\Delta\tau_{j}
\sum_{q=0}^{j-1}\Delta\tau_{q}\right)
<
$$
$$
<3\varepsilon^4\left(\delta(\varepsilon) (T-t)+(T-t)^2\right),
$$
where $\Delta\tau_{j}<\delta(\varepsilon),$
$j=0,1,\ldots,N-1$  ($\delta(\varepsilon)>0$ exists
for any $\varepsilon>0$ and it does not
depend on points of the interval $[t, T]$).
Then the right-hand side of the formula 
(\ref{30.43}) tends to zero when $N\to \infty.$ 
Considering this fact as well 
as (\ref{30.42}), we obtain (\ref{30.39}). 

If ${\bf w}_{t_l}^{(i_l)}=t_l$ for some $l\in\{1,\ldots,k\},$
then
the proof of Lemma 1.3 becomes obviously simpler and  
it is performed similarly. Lemma 1.3 is proved.

{\bf Remark 1.2.} {\it It is easy to see that if
$\Delta{\bf w}_{\tau_{j_l}}^{(i_l)}$ in {\rm (\ref{30.39})}
for some\ $l\in\{1,\ldots,k\}$ is replaced with 
$\left(\Delta{\bf w}_{\tau_{j_l}}^{(i_l)}\right)^p$ $(p=2,$
$i_l\ne 0),$ then
the differential $d{\bf w}_{t_{l}}^{(i_l)}$
in the integral $J[\Phi^{(k)}]_{T,t}$
will be replaced with $dt_l$.
If $p=3, 4,\ldots,$ then the
right-hand side 
of the formula {\rm (\ref{30.39})}
will become zero w. p. {\rm 1}.}

Let us consider the case $p=2$ in detail.
Let $\Delta{\bf w}_{\tau_{j_l}}^{(i_l)}$ in {\rm (\ref{30.39})}
for some $l\in\{1,\ldots,k\}$ is replaced with 
$\left(\Delta{\bf w}_{\tau_{j_l}}^{(i_l)}\right)^2$ $(i_l\ne 0)$ 
and
$$
J[\varphi_l]_{N}\stackrel{\rm def}{=}\sum\limits_{j=0}^{N-1}
\varphi_l(\tau_j)\left(\Delta{\bf w}_{\tau_j}^{(i_l)}\right)^2,\ \ \ 
J[\varphi_l]_{T,t}\stackrel{\rm def}{=}
\int\limits_t^T\varphi_l(s)ds.
$$

We have
$$
\left({\sf M}\left\{
\biggl|J[\varphi_l]_{N}-J[\varphi_l]_{T,t}\biggr|^4\right\}\right)
^{1/4}=
$$
$$
=\left({\sf M}\left\{
\left|\sum\limits_{j=0}^{N-1}
\varphi_l(\tau_j)\left(\Delta{\bf w}_{\tau_j}^{(i_l)}\right)^2-
\int\limits_t^T\varphi_l(s)ds
\right|^4\right\}\right)
^{1/4}=
$$
$$
=\left({\sf M}\left\{
\left|\sum\limits_{j=0}^{N-1}\left(
\varphi_l(\tau_j)\left(\Delta{\bf w}_{\tau_j}^{(i_l)}\right)^2-
\int\limits_{\tau_j}^{\tau_{j+1}}\left(\varphi_l(s)
-\varphi_l(\tau_j)+\varphi_l(\tau_j)\right)
ds\right)
\right|^4\right\}\right)
^{1/4}\le
$$
$$
\le \left({\sf M}\left\{
\left|\sum\limits_{j=0}^{N-1}
\varphi_l(\tau_j)\left(\left(\Delta{\bf w}_{\tau_j}^{(i_l)}\right)^2-
\Delta\tau_j\right)
\right|^4\right\}\right)
^{1/4}+
$$
\begin{equation}
\label{df1}
+
\left|\sum\limits_{j=0}^{N-1}\int\limits_{\tau_j}^{\tau_{j+1}}\left(
\varphi_l(\tau_j)-\varphi_l(s)\right)ds
\right|.
\end{equation}

From the relation, which is similar to {\rm (\ref{30.43})}, we obtain
$$
{\sf M}\left\{
\left|\sum\limits_{j=0}^{N-1}
\varphi_l(\tau_j)\left(\left(\Delta{\bf w}_{\tau_j}^{(i_l)}\right)^2-
\Delta\tau_j\right)
\right|^4\right\}=
$$
$$
=
\sum\limits_{j=0}^{N-1}
\left(\varphi_l(\tau_j)\right)^4
{\sf M}\left\{
\left(\left(\Delta{\bf w}_{\tau_j}^{(i_l)}\right)^2-
\Delta\tau_j\right)^4\right\}+
$$
$$
+6\sum\limits_{j=0}^{N-1}
\left(\varphi_l(\tau_j)\right)^2
{\sf M}\left\{
\left(\left(\Delta{\bf w}_{\tau_j}^{(i_l)}\right)^2-
\Delta\tau_j\right)^2\right\}\times
$$
$$
\times
\sum\limits_{q=0}^{j-1}
\left(\varphi_l(\tau_q)\right)^2
{\sf M}\left\{
\left(\left(\Delta{\bf w}_{\tau_q}^{(i_l)}\right)^2-
\Delta\tau_q\right)^2\right\}
=
60\sum\limits_{j=0}^{N-1}
\left(\varphi_l(\tau_j)\right)^4
\left(\Delta\tau_j\right)^4+
$$
\begin{equation}
\label{cloun1}
~~~~~ +
24
\sum\limits_{j=0}^{N-1}
\left(\varphi_l(\tau_j)\right)^2
\left(\Delta\tau_j\right)^2
\sum\limits_{q=0}^{j-1}\left(\varphi_l(\tau_q)\right)^2
\left(\Delta\tau_q\right)^2\le C\left(\Delta_N\right)^2\ \to 0
\end{equation}

\vspace{3mm}
\noindent
if $N\to \infty$, where constant $C$ does not depend on $N.$

The second term on the right-hand side of 
{\rm (\ref{df1})} tends to zero if $N\to \infty$ 
due to continuity (which means uniform continuity)
of the function $\varphi_l(s)$ at the interval $[t, T].$
Then, taking into account {\rm (\ref{df2})}, {\rm (\ref{30.42})},
{\rm (\ref{df1})}, {\rm (\ref{cloun1}),
we come to the affirmation of Remark {\rm 1.2.}

Let us prove Theorem 1.1.
According to Lemma 1.1, we have
$$
J[\psi^{(k)}]_{T,t}=
\hbox{\vtop{\offinterlineskip\halign{
\hfil#\hfil\cr
{\rm l.i.m.}\cr
$\stackrel{}{{}_{N\to \infty}}$\cr
}} }\sum_{l_k=0}^{N-1}\ldots\sum_{l_1=0}^{l_2-1}
\psi_1(\tau_{l_1})\ldots\psi_k(\tau_{l_k})
\Delta{\bf w}_{\tau_{l_1}}^{(i_1)}
\ldots
\Delta{\bf w}_{\tau_{l_k}}^{(i_k)}=
$$
$$
=\hbox{\vtop{\offinterlineskip\halign{
\hfil#\hfil\cr
{\rm l.i.m.}\cr
$\stackrel{}{{}_{N\to \infty}}$\cr
}} }\sum_{l_k=0}^{N-1}\ldots\sum_{l_1=0}^{l_2-1}
K(\tau_{l_1},\ldots,\tau_{l_k})
\Delta{\bf w}_{\tau_{l_1}}^{(i_1)}
\ldots
\Delta{\bf w}_{\tau_{l_k}}^{(i_k)}=
$$
$$
=
\hbox{\vtop{\offinterlineskip\halign{
\hfil#\hfil\cr
{\rm l.i.m.}\cr
$\stackrel{}{{}_{N\to \infty}}$\cr
}} }\sum_{l_k=0}^{N-1}\ldots\sum_{l_1=0}^{N-1}
K(\tau_{l_1},\ldots,\tau_{l_k})
\Delta{\bf w}_{\tau_{l_1}}^{(i_1)}
\ldots
\Delta{\bf w}_{\tau_{l_k}}^{(i_k)}=
$$
\begin{equation}
\label{zido1}
~~~~~~~=\hbox{\vtop{\offinterlineskip\halign{
\hfil#\hfil\cr
{\rm l.i.m.}\cr
$\stackrel{}{{}_{N\to \infty}}$\cr
}} }
\sum\limits_{\stackrel{l_1,\ldots,l_k=0}{{}_{l_q\ne l_r;\ 
q\ne r;\ q, r=1,\ldots, k}}}^{N-1}
K(\tau_{l_1},\ldots,\tau_{l_k})
\Delta{\bf w}_{\tau_{l_1}}^{(i_1)}
\ldots
\Delta{\bf w}_{\tau_{l_k}}^{(i_k)}=
\end{equation}
\begin{equation}
\label{hehe100}
~~~~~~=
\int\limits_{t}^{T}
\ldots
\int\limits_{t}^{t_2}
\sum_{(t_1,\ldots,t_k)}\left(
K(t_1,\ldots,t_k)d{\bf w}_{t_1}^{(i_1)}
\ldots
d{\bf w}_{t_k}^{(i_k)}\right)\ \ \ \hbox{w.~p.~1},
\end{equation}

\noindent
where permutations 
$(t_1,\ldots,t_k)$ when summing
are performed only 
in the expression enclosed in pa\-ren\-the\-ses.

It is easy to see that (\ref{hehe100})
can be rewritten in the form
\begin{equation}
\label{zab1}
~~~~~ J[\psi^{(k)}]_{T,t}=\sum_{(t_1,\ldots,t_k)}
\int\limits_{t}^{T}
\ldots
\int\limits_{t}^{t_2}
K(t_1,\ldots,t_k)d{\bf w}_{t_1}^{(i_1)}
\ldots
d{\bf w}_{t_k}^{(i_k)}\ \ \ \hbox{w.~p.~1},
\end{equation}

\noindent
where permutations
$(t_1,\ldots,t_k)$ when summing
are performed only in 
the values
$d{\bf w}_{t_1}^{(i_1)}
\ldots $
$d{\bf w}_{t_k}^{(i_k)}$. At the same time the indices near upper 
limits of integration in the iterated stochastic integrals are changed 
correspondently and if $t_r$ swapped with $t_q$ in the  
permutation $(t_1,\ldots,t_k)$, then $i_r$ swapped with $i_q$ in 
the permutation $(i_1,\ldots,i_k)$.

Since integration of a bounded function with respect to the
set with measure zero for Riemann or Lebesgue integrals gives zero result, then the 
following formula is correct for these integrals
\begin{equation}
\label{riemann}
~~ \int\limits_{[t, T]^k}|G(t_1,\ldots,t_k)|dt_1\ldots dt_k=
\sum_{(t_1,\ldots,t_k)}
\int\limits_{t}^{T}
\ldots
\int\limits_{t}^{t_2}
|G(t_1,\ldots,t_k)|dt_1\ldots dt_k,
\end{equation}

\noindent
where permutations $(t_1,\ldots,t_k)$ when summing
are performed only 
in the 
va\-lues $dt_1\ldots dt_k$. At the same time the indices near upper 
limits of integration in the
iterated integrals
are changed correspondently 
and $|G(t_1,\ldots,t_k)|$ is the integrable function on
the hypercube $[t, T]^k.$

According to Lemmas 1.1, 1.3 and (\ref{pobeda}), 
(\ref{s2s}), (\ref{hehe100}), (\ref{zab1}),
we get the following representation

\vspace{-2mm}
$$
J[\psi^{(k)}]_{T,t}=
$$

\vspace{-4mm}
$$
=
\sum_{j_1=0}^{p_1}\ldots
\sum_{j_k=0}^{p_k}
C_{j_k\ldots j_1}
\int\limits_{t}^{T}
\ldots
\int\limits_{t}^{t_2}
\sum_{(t_1,\ldots,t_k)}\left(
\phi_{j_1}(t_1)
\ldots
\phi_{j_k}(t_k)
d{\bf w}_{t_1}^{(i_1)}
\ldots
d{\bf w}_{t_k}^{(i_k)}\right)
+
$$

\vspace{2mm}
$$
+R_{T,t}^{p_1,\ldots,p_k}=
$$

\newpage
\noindent
$$
=\sum_{j_1=0}^{p_1}\ldots
\sum_{j_k=0}^{p_k}
C_{j_k\ldots j_1}\             
\hbox{\vtop{\offinterlineskip\halign{
\hfil#\hfil\cr
{\rm l.i.m.}\cr
$\stackrel{}{{}_{N\to \infty}}$\cr
}} }
\sum\limits_{\stackrel{l_1,\ldots,l_k=0}{{}_{l_q\ne l_r;\ 
q\ne r;\ q, r=1,\ldots, k}}}^{N-1}
\phi_{j_1}(\tau_{l_1})\ldots
\phi_{j_k}(\tau_{l_k})
\Delta{\bf w}_{\tau_{l_1}}^{(i_1)}
\ldots
\Delta{\bf w}_{\tau_{l_k}}^{(i_k)}+
$$

\vspace{4mm}
\begin{equation}
\label{novoe1}
+R_{T,t}^{p_1,\ldots,p_k}=
\end{equation}

\vspace{-2mm}
$$
=\sum_{j_1=0}^{p_1}\ldots
\sum_{j_k=0}^{p_k}
C_{j_k\ldots j_1}\left(
\hbox{\vtop{\offinterlineskip\halign{
\hfil#\hfil\cr
{\rm l.i.m.}\cr
$\stackrel{}{{}_{N\to \infty}}$\cr
}} }\sum_{l_1,\ldots,l_k=0}^{N-1}
\phi_{j_1}(\tau_{l_1})
\ldots
\phi_{j_k}(\tau_{l_k})
\Delta{\bf w}_{\tau_{l_1}}^{(i_1)}
\ldots
\Delta{\bf w}_{\tau_{l_k}}^{(i_k)}
-\right.
$$
$$
-\left.
\hbox{\vtop{\offinterlineskip\halign{
\hfil#\hfil\cr
{\rm l.i.m.}\cr
$\stackrel{}{{}_{N\to \infty}}$\cr
}} }\sum_{(l_1,\ldots,l_k)\in {\rm G}_k}
\phi_{j_{1}}(\tau_{l_1})
\Delta{\bf w}_{\tau_{l_1}}^{(i_1)}\ldots
\phi_{j_{k}}(\tau_{l_k})
\Delta{\bf w}_{\tau_{l_k}}^{(i_k)}\right)
+
$$

\vspace{1mm}
$$
+R_{T,t}^{p_1,\ldots,p_k}=
$$

\vspace{5mm}
$$
=\sum_{j_1=0}^{p_1}\ldots\sum_{j_k=0}^{p_k}
C_{j_k\ldots j_1}\times
$$
$$
\times
\left(
\prod_{l=1}^k\zeta_{j_l}^{(i_l)}-
\hbox{\vtop{\offinterlineskip\halign{
\hfil#\hfil\cr
{\rm l.i.m.}\cr
$\stackrel{}{{}_{N\to \infty}}$\cr
}} }\sum_{(l_1,\ldots,l_k)\in {\rm G}_k}
\phi_{j_{1}}(\tau_{l_1})
\Delta{\bf w}_{\tau_{l_1}}^{(i_1)}\ldots
\phi_{j_{k}}(\tau_{l_k})
\Delta{\bf w}_{\tau_{l_k}}^{(i_k)}\right)+
$$

\vspace{3mm}
\begin{equation}
\label{novoe2}
+R_{T,t}^{p_1,\ldots,p_k}\ \ \ \hbox{w.~p.~1},
\end{equation}

\vspace{3mm}
\noindent
where
$$
R_{T,t}^{p_1,\ldots,p_k}
=\sum_{(t_1,\ldots,t_k)}
\int\limits_{t}^{T}
\ldots
\int\limits_{t}^{t_2}
\left(K(t_1,\ldots,t_k)-
\sum_{j_1=0}^{p_1}\ldots
\sum_{j_k=0}^{p_k}
C_{j_k\ldots j_1}
\prod_{l=1}^k\phi_{j_l}(t_l)\right)\times
$$

\begin{equation}
\label{y007}
\times
d{\bf w}_{t_1}^{(i_1)}
\ldots
d{\bf w}_{t_k}^{(i_k)},
\end{equation}

\vspace{3mm}
\noindent
where permutations $(t_1,\ldots,t_k)$ when summing are performed only 
in the values $d{\bf w}_{t_1}^{(i_1)}
\ldots $
$d{\bf w}_{t_k}^{(i_k)}$. At the same time the indices near 
upper limits of integration in the iterated stochastic integrals 
are changed correspondently and if $t_r$ swapped with $t_q$ in the  
permutation $(t_1,\ldots,t_k)$, then $i_r$ swapped with $i_q$ in the 
permutation $(i_1,\ldots,i_k)$.

Let us estimate the remainder
$R_{T,t}^{p_1,\ldots,p_k}$ of the series.
According to Lemma 1.2 and (\ref{riemann}), we have

\vspace{-2mm}
$$
{\sf M}\left\{\left(R_{T,t}^{p_1,\ldots,p_k}\right)^2\right\}
\le 
$$

\vspace{-1mm}
$$
\le C_k
\hspace{-2mm}\sum_{(t_1,\ldots,t_k)}
\int\limits_{t}^{T}
\ldots
\int\limits_{t}^{t_2}
\left(K(t_1,\ldots,t_k)-
\sum_{j_1=0}^{p_1}\ldots
\sum_{j_k=0}^{p_k}
C_{j_k\ldots j_1}
\prod_{l=1}^k\phi_{j_l}(t_l)\right)^2
\hspace{-2mm}dt_1
\ldots
dt_k=
$$

\vspace{-1mm}
\begin{equation}
\label{obana1}
=C_k\int\limits_{[t,T]^k}
\left(K(t_1,\ldots,t_k)-
\sum_{j_1=0}^{p_1}\ldots
\sum_{j_k=0}^{p_k}
C_{j_k\ldots j_1}
\prod_{l=1}^k\phi_{j_l}(t_l)\right)^2
dt_1
\ldots
dt_k\to 0
\end{equation}

\vspace{1mm}
\noindent
if $p_1,\ldots,p_k\to\infty,$ where constant $C_k$ 
depends only
on the multiplicity $k$ of the iterated It\^{o} stochastic integral
$J[\psi^{(k)}]_{T,t}$. 
Theorem 1.1 is proved.

Note that from (\ref{novoe1}) and (\ref{obana1}) it follows that
\begin{equation}
\label{drdr1}
~~~~~J[\psi^{(k)}]_{T,t}=
\hbox{\vtop{\offinterlineskip\halign{
\hfil#\hfil\cr
{\rm l.i.m.}\cr
$\stackrel{}{{}_{p_1,\ldots,p_k\to \infty}}$\cr
}} }\sum_{j_1=0}^{p_1}\ldots\sum_{j_k=0}^{p_k}
C_{j_k\ldots j_1}J'[\phi_{j_1}\ldots \phi_{j_k}]_{T,t}^{(i_1\ldots i_k)}\ \ \ \hbox{w.~p.~1},
\end{equation}

\noindent
where $J'[\phi_{j_1}\ldots \phi_{j_k}]_{T,t}^{(i_1\ldots i_k)}$
is defined by (\ref{mult11}).

It is not difficult to see that for the case of pairwise different numbers
$i_1,\ldots,i_k=0, 1,\ldots,m$ from Theorem 1.1 we obtain
\begin{equation}
\label{ziko5000}
J[\psi^{(k)}]_{T,t}=
\hbox{\vtop{\offinterlineskip\halign{
\hfil#\hfil\cr
{\rm l.i.m.}\cr
$\stackrel{}{{}_{p_1,\ldots,p_k\to \infty}}$\cr
}} }\sum_{j_1=0}^{p_1}\ldots\sum_{j_k=0}^{p_k}
C_{j_k\ldots j_1}\zeta_{j_1}^{(i_1)}\ldots \zeta_{j_k}^{(i_k)}.
\end{equation}

\subsection{Expansions of Iterated It\^{o} Stochastic 
Integrals with Multiplicities 
1 to 7 Based on Theorem 1.1}

In order to evaluate the significance of Theorem 1.1 for practice we will
demonstrate its transformed particular cases (see Remark 1.2) for 
$k=1,\ldots,7$ \cite{1}-\cite{new-new-6}
\begin{equation}
\label{a1}
J[\psi^{(1)}]_{T,t}
=\hbox{\vtop{\offinterlineskip\halign{
\hfil#\hfil\cr
{\rm l.i.m.}\cr
$\stackrel{}{{}_{p_1\to \infty}}$\cr
}} }\sum_{j_1=0}^{p_1}
C_{j_1}\zeta_{j_1}^{(i_1)},
\end{equation}
\begin{equation}
\label{a2}
~~~~~~ J[\psi^{(2)}]_{T,t}
=\hbox{\vtop{\offinterlineskip\halign{
\hfil#\hfil\cr
{\rm l.i.m.}\cr
$\stackrel{}{{}_{p_1,p_2\to \infty}}$\cr
}} }\sum_{j_1=0}^{p_1}\sum_{j_2=0}^{p_2}
C_{j_2j_1}\Biggl(\zeta_{j_1}^{(i_1)}\zeta_{j_2}^{(i_2)}
-{\bf 1}_{\{i_1=i_2\ne 0\}}
{\bf 1}_{\{j_1=j_2\}}\Biggr),
\end{equation}

\vspace{3mm}
$$
J[\psi^{(3)}]_{T,t}=
\hbox{\vtop{\offinterlineskip\halign{
\hfil#\hfil\cr
{\rm l.i.m.}\cr
$\stackrel{}{{}_{p_1,p_2,p_3\to \infty}}$\cr
}} }\sum_{j_1=0}^{p_1}\sum_{j_2=0}^{p_2}\sum_{j_3=0}^{p_3}
C_{j_3j_2j_1}\Biggl(
\zeta_{j_1}^{(i_1)}\zeta_{j_2}^{(i_2)}\zeta_{j_3}^{(i_3)}
-\Biggr.
$$
\begin{equation}
\label{a3}
~~ -\Biggl.
{\bf 1}_{\{i_1=i_2\ne 0\}}
{\bf 1}_{\{j_1=j_2\}}
\zeta_{j_3}^{(i_3)}
-{\bf 1}_{\{i_2=i_3\ne 0\}}
{\bf 1}_{\{j_2=j_3\}}
\zeta_{j_1}^{(i_1)}-
{\bf 1}_{\{i_1=i_3\ne 0\}}
{\bf 1}_{\{j_1=j_3\}}
\zeta_{j_2}^{(i_2)}\Biggr),
\end{equation}

\vspace{2mm}

$$
J[\psi^{(4)}]_{T,t}
=
\hbox{\vtop{\offinterlineskip\halign{
\hfil#\hfil\cr
{\rm l.i.m.}\cr
$\stackrel{}{{}_{p_1,\ldots,p_4\to \infty}}$\cr
}} }\sum_{j_1=0}^{p_1}\ldots\sum_{j_4=0}^{p_4}
C_{j_4\ldots j_1}\Biggl(
\prod_{l=1}^4\zeta_{j_l}^{(i_l)}
\Biggr.
-
$$
$$
-
{\bf 1}_{\{i_1=i_2\ne 0\}}
{\bf 1}_{\{j_1=j_2\}}
\zeta_{j_3}^{(i_3)}
\zeta_{j_4}^{(i_4)}
-
{\bf 1}_{\{i_1=i_3\ne 0\}}
{\bf 1}_{\{j_1=j_3\}}
\zeta_{j_2}^{(i_2)}
\zeta_{j_4}^{(i_4)}-
$$
$$
-
{\bf 1}_{\{i_1=i_4\ne 0\}}
{\bf 1}_{\{j_1=j_4\}}
\zeta_{j_2}^{(i_2)}
\zeta_{j_3}^{(i_3)}
-
{\bf 1}_{\{i_2=i_3\ne 0\}}
{\bf 1}_{\{j_2=j_3\}}
\zeta_{j_1}^{(i_1)}
\zeta_{j_4}^{(i_4)}-
$$
$$
-
{\bf 1}_{\{i_2=i_4\ne 0\}}
{\bf 1}_{\{j_2=j_4\}}
\zeta_{j_1}^{(i_1)}
\zeta_{j_3}^{(i_3)}
-
{\bf 1}_{\{i_3=i_4\ne 0\}}
{\bf 1}_{\{j_3=j_4\}}
\zeta_{j_1}^{(i_1)}
\zeta_{j_2}^{(i_2)}+
$$
$$
+
{\bf 1}_{\{i_1=i_2\ne 0\}}
{\bf 1}_{\{j_1=j_2\}}
{\bf 1}_{\{i_3=i_4\ne 0\}}
{\bf 1}_{\{j_3=j_4\}}
+
{\bf 1}_{\{i_1=i_3\ne 0\}}
{\bf 1}_{\{j_1=j_3\}}
{\bf 1}_{\{i_2=i_4\ne 0\}}
{\bf 1}_{\{j_2=j_4\}}+
$$
\begin{equation}
\label{a4}
+\Biggl.
{\bf 1}_{\{i_1=i_4\ne 0\}}
{\bf 1}_{\{j_1=j_4\}}
{\bf 1}_{\{i_2=i_3\ne 0\}}
{\bf 1}_{\{j_2=j_3\}}\Biggr),
\end{equation}

\vspace{4mm}

$$
J[\psi^{(5)}]_{T,t}
=\hbox{\vtop{\offinterlineskip\halign{
\hfil#\hfil\cr
{\rm l.i.m.}\cr
$\stackrel{}{{}_{p_1,\ldots,p_5\to \infty}}$\cr
}} }\sum_{j_1=0}^{p_1}\ldots\sum_{j_5=0}^{p_5}
C_{j_5\ldots j_1}\Biggl(
\prod_{l=1}^5\zeta_{j_l}^{(i_l)}
-\Biggr.
$$
$$
-
{\bf 1}_{\{i_1=i_2\ne 0\}}
{\bf 1}_{\{j_1=j_2\}}
\zeta_{j_3}^{(i_3)}
\zeta_{j_4}^{(i_4)}
\zeta_{j_5}^{(i_5)}-
{\bf 1}_{\{i_1=i_3\ne 0\}}
{\bf 1}_{\{j_1=j_3\}}
\zeta_{j_2}^{(i_2)}
\zeta_{j_4}^{(i_4)}
\zeta_{j_5}^{(i_5)}-
$$
$$
-
{\bf 1}_{\{i_1=i_4\ne 0\}}
{\bf 1}_{\{j_1=j_4\}}
\zeta_{j_2}^{(i_2)}
\zeta_{j_3}^{(i_3)}
\zeta_{j_5}^{(i_5)}-
{\bf 1}_{\{i_1=i_5\ne 0\}}
{\bf 1}_{\{j_1=j_5\}}
\zeta_{j_2}^{(i_2)}
\zeta_{j_3}^{(i_3)}
\zeta_{j_4}^{(i_4)}-
$$
$$
-
{\bf 1}_{\{i_2=i_3\ne 0\}}
{\bf 1}_{\{j_2=j_3\}}
\zeta_{j_1}^{(i_1)}
\zeta_{j_4}^{(i_4)}
\zeta_{j_5}^{(i_5)}-
{\bf 1}_{\{i_2=i_4\ne 0\}}
{\bf 1}_{\{j_2=j_4\}}
\zeta_{j_1}^{(i_1)}
\zeta_{j_3}^{(i_3)}
\zeta_{j_5}^{(i_5)}-
$$
$$
-
{\bf 1}_{\{i_2=i_5\ne 0\}}
{\bf 1}_{\{j_2=j_5\}}
\zeta_{j_1}^{(i_1)}
\zeta_{j_3}^{(i_3)}
\zeta_{j_4}^{(i_4)}
-{\bf 1}_{\{i_3=i_4\ne 0\}}
{\bf 1}_{\{j_3=j_4\}}
\zeta_{j_1}^{(i_1)}
\zeta_{j_2}^{(i_2)}
\zeta_{j_5}^{(i_5)}-
$$
$$
-
{\bf 1}_{\{i_3=i_5\ne 0\}}
{\bf 1}_{\{j_3=j_5\}}
\zeta_{j_1}^{(i_1)}
\zeta_{j_2}^{(i_2)}
\zeta_{j_4}^{(i_4)}
-{\bf 1}_{\{i_4=i_5\ne 0\}}
{\bf 1}_{\{j_4=j_5\}}
\zeta_{j_1}^{(i_1)}
\zeta_{j_2}^{(i_2)}
\zeta_{j_3}^{(i_3)}+
$$
$$
+
{\bf 1}_{\{i_1=i_2\ne 0\}}
{\bf 1}_{\{j_1=j_2\}}
{\bf 1}_{\{i_3=i_4\ne 0\}}
{\bf 1}_{\{j_3=j_4\}}\zeta_{j_5}^{(i_5)}+
{\bf 1}_{\{i_1=i_2\ne 0\}}
{\bf 1}_{\{j_1=j_2\}}
{\bf 1}_{\{i_3=i_5\ne 0\}}
{\bf 1}_{\{j_3=j_5\}}\zeta_{j_4}^{(i_4)}+
$$
$$
+
{\bf 1}_{\{i_1=i_2\ne 0\}}
{\bf 1}_{\{j_1=j_2\}}
{\bf 1}_{\{i_4=i_5\ne 0\}}
{\bf 1}_{\{j_4=j_5\}}\zeta_{j_3}^{(i_3)}+
{\bf 1}_{\{i_1=i_3\ne 0\}}
{\bf 1}_{\{j_1=j_3\}}
{\bf 1}_{\{i_2=i_4\ne 0\}}
{\bf 1}_{\{j_2=j_4\}}\zeta_{j_5}^{(i_5)}+
$$
$$
+
{\bf 1}_{\{i_1=i_3\ne 0\}}
{\bf 1}_{\{j_1=j_3\}}
{\bf 1}_{\{i_2=i_5\ne 0\}}
{\bf 1}_{\{j_2=j_5\}}\zeta_{j_4}^{(i_4)}+
{\bf 1}_{\{i_1=i_3\ne 0\}}
{\bf 1}_{\{j_1=j_3\}}
{\bf 1}_{\{i_4=i_5\ne 0\}}
{\bf 1}_{\{j_4=j_5\}}\zeta_{j_2}^{(i_2)}+
$$
$$
+
{\bf 1}_{\{i_1=i_4\ne 0\}}
{\bf 1}_{\{j_1=j_4\}}
{\bf 1}_{\{i_2=i_3\ne 0\}}
{\bf 1}_{\{j_2=j_3\}}\zeta_{j_5}^{(i_5)}+
{\bf 1}_{\{i_1=i_4\ne 0\}}
{\bf 1}_{\{j_1=j_4\}}
{\bf 1}_{\{i_2=i_5\ne 0\}}
{\bf 1}_{\{j_2=j_5\}}\zeta_{j_3}^{(i_3)}+
$$
$$
+
{\bf 1}_{\{i_1=i_4\ne 0\}}
{\bf 1}_{\{j_1=j_4\}}
{\bf 1}_{\{i_3=i_5\ne 0\}}
{\bf 1}_{\{j_3=j_5\}}\zeta_{j_2}^{(i_2)}+
{\bf 1}_{\{i_1=i_5\ne 0\}}
{\bf 1}_{\{j_1=j_5\}}
{\bf 1}_{\{i_2=i_3\ne 0\}}
{\bf 1}_{\{j_2=j_3\}}\zeta_{j_4}^{(i_4)}+
$$
$$
+
{\bf 1}_{\{i_1=i_5\ne 0\}}
{\bf 1}_{\{j_1=j_5\}}
{\bf 1}_{\{i_2=i_4\ne 0\}}
{\bf 1}_{\{j_2=j_4\}}\zeta_{j_3}^{(i_3)}+
{\bf 1}_{\{i_1=i_5\ne 0\}}
{\bf 1}_{\{j_1=j_5\}}
{\bf 1}_{\{i_3=i_4\ne 0\}}
{\bf 1}_{\{j_3=j_4\}}\zeta_{j_2}^{(i_2)}+
$$
$$
+
{\bf 1}_{\{i_2=i_3\ne 0\}}
{\bf 1}_{\{j_2=j_3\}}
{\bf 1}_{\{i_4=i_5\ne 0\}}
{\bf 1}_{\{j_4=j_5\}}\zeta_{j_1}^{(i_1)}+
{\bf 1}_{\{i_2=i_4\ne 0\}}
{\bf 1}_{\{j_2=j_4\}}
{\bf 1}_{\{i_3=i_5\ne 0\}}
{\bf 1}_{\{j_3=j_5\}}\zeta_{j_1}^{(i_1)}+
$$
\begin{equation}
\label{a5}
+\Biggl.
{\bf 1}_{\{i_2=i_5\ne 0\}}
{\bf 1}_{\{j_2=j_5\}}
{\bf 1}_{\{i_3=i_4\ne 0\}}
{\bf 1}_{\{j_3=j_4\}}\zeta_{j_1}^{(i_1)}\Biggr),
\end{equation}

\vspace{4mm}

$$
J[\psi^{(6)}]_{T,t}
=\hbox{\vtop{\offinterlineskip\halign{
\hfil#\hfil\cr
{\rm l.i.m.}\cr
$\stackrel{}{{}_{p_1,\ldots,p_6\to \infty}}$\cr
}} }\sum_{j_1=0}^{p_1}\ldots\sum_{j_6=0}^{p_6}
C_{j_6\ldots j_1}\Biggl(
\prod_{l=1}^6
\zeta_{j_l}^{(i_l)}
-\Biggr.
$$
$$
-
{\bf 1}_{\{i_1=i_6\ne 0\}}
{\bf 1}_{\{j_1=j_6\}}
\zeta_{j_2}^{(i_2)}
\zeta_{j_3}^{(i_3)}
\zeta_{j_4}^{(i_4)}
\zeta_{j_5}^{(i_5)}-
{\bf 1}_{\{i_2=i_6\ne 0\}}
{\bf 1}_{\{j_2=j_6\}}
\zeta_{j_1}^{(i_1)}
\zeta_{j_3}^{(i_3)}
\zeta_{j_4}^{(i_4)}
\zeta_{j_5}^{(i_5)}-
$$
$$
-
{\bf 1}_{\{i_3=i_6\ne 0\}}
{\bf 1}_{\{j_3=j_6\}}
\zeta_{j_1}^{(i_1)}
\zeta_{j_2}^{(i_2)}
\zeta_{j_4}^{(i_4)}
\zeta_{j_5}^{(i_5)}-
{\bf 1}_{\{i_4=i_6\ne 0\}}
{\bf 1}_{\{j_4=j_6\}}
\zeta_{j_1}^{(i_1)}
\zeta_{j_2}^{(i_2)}
\zeta_{j_3}^{(i_3)}
\zeta_{j_5}^{(i_5)}-
$$
$$
-
{\bf 1}_{\{i_5=i_6\ne 0\}}
{\bf 1}_{\{j_5=j_6\}}
\zeta_{j_1}^{(i_1)}
\zeta_{j_2}^{(i_2)}
\zeta_{j_3}^{(i_3)}
\zeta_{j_4}^{(i_4)}-
{\bf 1}_{\{i_1=i_2\ne 0\}}
{\bf 1}_{\{j_1=j_2\}}
\zeta_{j_3}^{(i_3)}
\zeta_{j_4}^{(i_4)}
\zeta_{j_5}^{(i_5)}
\zeta_{j_6}^{(i_6)}-
$$
$$
-
{\bf 1}_{\{i_1=i_3\ne 0\}}
{\bf 1}_{\{j_1=j_3\}}
\zeta_{j_2}^{(i_2)}
\zeta_{j_4}^{(i_4)}
\zeta_{j_5}^{(i_5)}
\zeta_{j_6}^{(i_6)}-
{\bf 1}_{\{i_1=i_4\ne 0\}}
{\bf 1}_{\{j_1=j_4\}}
\zeta_{j_2}^{(i_2)}
\zeta_{j_3}^{(i_3)}
\zeta_{j_5}^{(i_5)}
\zeta_{j_6}^{(i_6)}-
$$
$$
-
{\bf 1}_{\{i_1=i_5\ne 0\}}
{\bf 1}_{\{j_1=j_5\}}
\zeta_{j_2}^{(i_2)}
\zeta_{j_3}^{(i_3)}
\zeta_{j_4}^{(i_4)}
\zeta_{j_6}^{(i_6)}-
{\bf 1}_{\{i_2=i_3\ne 0\}}
{\bf 1}_{\{j_2=j_3\}}
\zeta_{j_1}^{(i_1)}
\zeta_{j_4}^{(i_4)}
\zeta_{j_5}^{(i_5)}
\zeta_{j_6}^{(i_6)}-
$$
$$
-
{\bf 1}_{\{i_2=i_4\ne 0\}}
{\bf 1}_{\{j_2=j_4\}}
\zeta_{j_1}^{(i_1)}
\zeta_{j_3}^{(i_3)}
\zeta_{j_5}^{(i_5)}
\zeta_{j_6}^{(i_6)}-
{\bf 1}_{\{i_2=i_5\ne 0\}}
{\bf 1}_{\{j_2=j_5\}}
\zeta_{j_1}^{(i_1)}
\zeta_{j_3}^{(i_3)}
\zeta_{j_4}^{(i_4)}
\zeta_{j_6}^{(i_6)}-
$$
$$
-
{\bf 1}_{\{i_3=i_4\ne 0\}}
{\bf 1}_{\{j_3=j_4\}}
\zeta_{j_1}^{(i_1)}
\zeta_{j_2}^{(i_2)}
\zeta_{j_5}^{(i_5)}
\zeta_{j_6}^{(i_6)}-
{\bf 1}_{\{i_3=i_5\ne 0\}}
{\bf 1}_{\{j_3=j_5\}}
\zeta_{j_1}^{(i_1)}
\zeta_{j_2}^{(i_2)}
\zeta_{j_4}^{(i_4)}
\zeta_{j_6}^{(i_6)}-
$$
$$
-
{\bf 1}_{\{i_4=i_5\ne 0\}}
{\bf 1}_{\{j_4=j_5\}}
\zeta_{j_1}^{(i_1)}
\zeta_{j_2}^{(i_2)}
\zeta_{j_3}^{(i_3)}
\zeta_{j_6}^{(i_6)}+
$$
$$
+
{\bf 1}_{\{i_1=i_2\ne 0\}}
{\bf 1}_{\{j_1=j_2\}}
{\bf 1}_{\{i_3=i_4\ne 0\}}
{\bf 1}_{\{j_3=j_4\}}
\zeta_{j_5}^{(i_5)}
\zeta_{j_6}^{(i_6)}+
$$
$$
+
{\bf 1}_{\{i_1=i_2\ne 0\}}
{\bf 1}_{\{j_1=j_2\}}
{\bf 1}_{\{i_3=i_5\ne 0\}}
{\bf 1}_{\{j_3=j_5\}}
\zeta_{j_4}^{(i_4)}
\zeta_{j_6}^{(i_6)}+
$$
$$
+
{\bf 1}_{\{i_1=i_2\ne 0\}}
{\bf 1}_{\{j_1=j_2\}}
{\bf 1}_{\{i_4=i_5\ne 0\}}
{\bf 1}_{\{j_4=j_5\}}
\zeta_{j_3}^{(i_3)}
\zeta_{j_6}^{(i_6)}
+
$$
$$
+
{\bf 1}_{\{i_1=i_3\ne 0\}}
{\bf 1}_{\{j_1=j_3\}}
{\bf 1}_{\{i_2=i_4\ne 0\}}
{\bf 1}_{\{j_2=j_4\}}
\zeta_{j_5}^{(i_5)}
\zeta_{j_6}^{(i_6)}+
$$
$$
+
{\bf 1}_{\{i_1=i_3\ne 0\}}
{\bf 1}_{\{j_1=j_3\}}
{\bf 1}_{\{i_2=i_5\ne 0\}}
{\bf 1}_{\{j_2=j_5\}}
\zeta_{j_4}^{(i_4)}
\zeta_{j_6}^{(i_6)}
+
$$
$$
+{\bf 1}_{\{i_1=i_3\ne 0\}}
{\bf 1}_{\{j_1=j_3\}}
{\bf 1}_{\{i_4=i_5\ne 0\}}
{\bf 1}_{\{j_4=j_5\}}
\zeta_{j_2}^{(i_2)}
\zeta_{j_6}^{(i_6)}+
$$
$$
+
{\bf 1}_{\{i_1=i_4\ne 0\}}
{\bf 1}_{\{j_1=j_4\}}
{\bf 1}_{\{i_2=i_3\ne 0\}}
{\bf 1}_{\{j_2=j_3\}}
\zeta_{j_5}^{(i_5)}
\zeta_{j_6}^{(i_6)}
+
$$
$$
+
{\bf 1}_{\{i_1=i_4\ne 0\}}
{\bf 1}_{\{j_1=j_4\}}
{\bf 1}_{\{i_2=i_5\ne 0\}}
{\bf 1}_{\{j_2=j_5\}}
\zeta_{j_3}^{(i_3)}
\zeta_{j_6}^{(i_6)}+
$$
$$
+
{\bf 1}_{\{i_1=i_4\ne 0\}}
{\bf 1}_{\{j_1=j_4\}}
{\bf 1}_{\{i_3=i_5\ne 0\}}
{\bf 1}_{\{j_3=j_5\}}
\zeta_{j_2}^{(i_2)}
\zeta_{j_6}^{(i_6)}
+
$$
$$
+
{\bf 1}_{\{i_1=i_5\ne 0\}}
{\bf 1}_{\{j_1=j_5\}}
{\bf 1}_{\{i_2=i_3\ne 0\}}
{\bf 1}_{\{j_2=j_3\}}
\zeta_{j_4}^{(i_4)}
\zeta_{j_6}^{(i_6)}+
$$
$$
+
{\bf 1}_{\{i_1=i_5\ne 0\}}
{\bf 1}_{\{j_1=j_5\}}
{\bf 1}_{\{i_2=i_4\ne 0\}}
{\bf 1}_{\{j_2=j_4\}}
\zeta_{j_3}^{(i_3)}
\zeta_{j_6}^{(i_6)}
+
$$
$$
+
{\bf 1}_{\{i_1=i_5\ne 0\}}
{\bf 1}_{\{j_1=j_5\}}
{\bf 1}_{\{i_3=i_4\ne 0\}}
{\bf 1}_{\{j_3=j_4\}}
\zeta_{j_2}^{(i_2)}
\zeta_{j_6}^{(i_6)}+
$$
$$
+
{\bf 1}_{\{i_2=i_3\ne 0\}}
{\bf 1}_{\{j_2=j_3\}}
{\bf 1}_{\{i_4=i_5\ne 0\}}
{\bf 1}_{\{j_4=j_5\}}
\zeta_{j_1}^{(i_1)}
\zeta_{j_6}^{(i_6)}
+
$$
$$
+{\bf 1}_{\{i_2=i_4\ne 0\}}
{\bf 1}_{\{j_2=j_4\}}
{\bf 1}_{\{i_3=i_5\ne 0\}}
{\bf 1}_{\{j_3=j_5\}}
\zeta_{j_1}^{(i_1)}
\zeta_{j_6}^{(i_6)}+
$$
$$
+
{\bf 1}_{\{i_2=i_5\ne 0\}}
{\bf 1}_{\{j_2=j_5\}}
{\bf 1}_{\{i_3=i_4\ne 0\}}
{\bf 1}_{\{j_3=j_4\}}
\zeta_{j_1}^{(i_1)}
\zeta_{j_6}^{(i_6)}
+
$$
$$
+{\bf 1}_{\{i_6=i_1\ne 0\}}
{\bf 1}_{\{j_6=j_1\}}
{\bf 1}_{\{i_3=i_4\ne 0\}}
{\bf 1}_{\{j_3=j_4\}}
\zeta_{j_2}^{(i_2)}
\zeta_{j_5}^{(i_5)}+
$$
$$
+
{\bf 1}_{\{i_6=i_1\ne 0\}}
{\bf 1}_{\{j_6=j_1\}}
{\bf 1}_{\{i_3=i_5\ne 0\}}
{\bf 1}_{\{j_3=j_5\}}
\zeta_{j_2}^{(i_2)}
\zeta_{j_4}^{(i_4)}
+
$$
$$
+{\bf 1}_{\{i_6=i_1\ne 0\}}
{\bf 1}_{\{j_6=j_1\}}
{\bf 1}_{\{i_2=i_5\ne 0\}}
{\bf 1}_{\{j_2=j_5\}}
\zeta_{j_3}^{(i_3)}
\zeta_{j_4}^{(i_4)}+
$$
$$
+
{\bf 1}_{\{i_6=i_1\ne 0\}}
{\bf 1}_{\{j_6=j_1\}}
{\bf 1}_{\{i_2=i_4\ne 0\}}
{\bf 1}_{\{j_2=j_4\}}
\zeta_{j_3}^{(i_3)}
\zeta_{j_5}^{(i_5)}
+
$$
$$
+{\bf 1}_{\{i_6=i_1\ne 0\}}
{\bf 1}_{\{j_6=j_1\}}
{\bf 1}_{\{i_4=i_5\ne 0\}}
{\bf 1}_{\{j_4=j_5\}}
\zeta_{j_2}^{(i_2)}
\zeta_{j_3}^{(i_3)}+
$$
$$
+
{\bf 1}_{\{i_6=i_1\ne 0\}}
{\bf 1}_{\{j_6=j_1\}}
{\bf 1}_{\{i_2=i_3\ne 0\}}
{\bf 1}_{\{j_2=j_3\}}
\zeta_{j_4}^{(i_4)}
\zeta_{j_5}^{(i_5)}
+
$$
$$
+{\bf 1}_{\{i_6=i_2\ne 0\}}
{\bf 1}_{\{j_6=j_2\}}
{\bf 1}_{\{i_3=i_5\ne 0\}}
{\bf 1}_{\{j_3=j_5\}}
\zeta_{j_1}^{(i_1)}
\zeta_{j_4}^{(i_4)}+
$$
$$
+
{\bf 1}_{\{i_6=i_2\ne 0\}}
{\bf 1}_{\{j_6=j_2\}}
{\bf 1}_{\{i_4=i_5\ne 0\}}
{\bf 1}_{\{j_4=j_5\}}
\zeta_{j_1}^{(i_1)}
\zeta_{j_3}^{(i_3)}
+
$$
$$
+{\bf 1}_{\{i_6=i_2\ne 0\}}
{\bf 1}_{\{j_6=j_2\}}
{\bf 1}_{\{i_3=i_4\ne 0\}}
{\bf 1}_{\{j_3=j_4\}}
\zeta_{j_1}^{(i_1)}
\zeta_{j_5}^{(i_5)}+
$$
$$
+
{\bf 1}_{\{i_6=i_2\ne 0\}}
{\bf 1}_{\{j_6=j_2\}}
{\bf 1}_{\{i_1=i_5\ne 0\}}
{\bf 1}_{\{j_1=j_5\}}
\zeta_{j_3}^{(i_3)}
\zeta_{j_4}^{(i_4)}
+
$$
$$
+{\bf 1}_{\{i_6=i_2\ne 0\}}
{\bf 1}_{\{j_6=j_2\}}
{\bf 1}_{\{i_1=i_4\ne 0\}}
{\bf 1}_{\{j_1=j_4\}}
\zeta_{j_3}^{(i_3)}
\zeta_{j_5}^{(i_5)}+
$$
$$
+
{\bf 1}_{\{i_6=i_2\ne 0\}}
{\bf 1}_{\{j_6=j_2\}}
{\bf 1}_{\{i_1=i_3\ne 0\}}
{\bf 1}_{\{j_1=j_3\}}
\zeta_{j_4}^{(i_4)}
\zeta_{j_5}^{(i_5)}
+
$$
$$
+{\bf 1}_{\{i_6=i_3\ne 0\}}
{\bf 1}_{\{j_6=j_3\}}
{\bf 1}_{\{i_2=i_5\ne 0\}}
{\bf 1}_{\{j_2=j_5\}}
\zeta_{j_1}^{(i_1)}
\zeta_{j_4}^{(i_4)}+
$$
$$
+
{\bf 1}_{\{i_6=i_3\ne 0\}}
{\bf 1}_{\{j_6=j_3\}}
{\bf 1}_{\{i_4=i_5\ne 0\}}
{\bf 1}_{\{j_4=j_5\}}
\zeta_{j_1}^{(i_1)}
\zeta_{j_2}^{(i_2)}
+
$$
$$
+{\bf 1}_{\{i_6=i_3\ne 0\}}
{\bf 1}_{\{j_6=j_3\}}
{\bf 1}_{\{i_2=i_4\ne 0\}}
{\bf 1}_{\{j_2=j_4\}}
\zeta_{j_1}^{(i_1)}
\zeta_{j_5}^{(i_5)}+
$$
$$
+
{\bf 1}_{\{i_6=i_3\ne 0\}}
{\bf 1}_{\{j_6=j_3\}}
{\bf 1}_{\{i_1=i_5\ne 0\}}
{\bf 1}_{\{j_1=j_5\}}
\zeta_{j_2}^{(i_2)}
\zeta_{j_4}^{(i_4)}
+
$$
$$
+{\bf 1}_{\{i_6=i_3\ne 0\}}
{\bf 1}_{\{j_6=j_3\}}
{\bf 1}_{\{i_1=i_4\ne 0\}}
{\bf 1}_{\{j_1=j_4\}}
\zeta_{j_2}^{(i_2)}
\zeta_{j_5}^{(i_5)}+
$$
$$
+
{\bf 1}_{\{i_6=i_3\ne 0\}}
{\bf 1}_{\{j_6=j_3\}}
{\bf 1}_{\{i_1=i_2\ne 0\}}
{\bf 1}_{\{j_1=j_2\}}
\zeta_{j_4}^{(i_4)}
\zeta_{j_5}^{(i_5)}
+
$$
$$
+{\bf 1}_{\{i_6=i_4\ne 0\}}
{\bf 1}_{\{j_6=j_4\}}
{\bf 1}_{\{i_3=i_5\ne 0\}}
{\bf 1}_{\{j_3=j_5\}}
\zeta_{j_1}^{(i_1)}
\zeta_{j_2}^{(i_2)}+
$$
$$
+
{\bf 1}_{\{i_6=i_4\ne 0\}}
{\bf 1}_{\{j_6=j_4\}}
{\bf 1}_{\{i_2=i_5\ne 0\}}
{\bf 1}_{\{j_2=j_5\}}
\zeta_{j_1}^{(i_1)}
\zeta_{j_3}^{(i_3)}
+
$$
$$
+{\bf 1}_{\{i_6=i_4\ne 0\}}
{\bf 1}_{\{j_6=j_4\}}
{\bf 1}_{\{i_2=i_3\ne 0\}}
{\bf 1}_{\{j_2=j_3\}}
\zeta_{j_1}^{(i_1)}
\zeta_{j_5}^{(i_5)}+
$$
$$
+
{\bf 1}_{\{i_6=i_4\ne 0\}}
{\bf 1}_{\{j_6=j_4\}}
{\bf 1}_{\{i_1=i_5\ne 0\}}
{\bf 1}_{\{j_1=j_5\}}
\zeta_{j_2}^{(i_2)}
\zeta_{j_3}^{(i_3)}
+
$$
$$
+{\bf 1}_{\{i_6=i_4\ne 0\}}
{\bf 1}_{\{j_6=j_4\}}
{\bf 1}_{\{i_1=i_3\ne 0\}}
{\bf 1}_{\{j_1=j_3\}}
\zeta_{j_2}^{(i_2)}
\zeta_{j_5}^{(i_5)}+
$$
$$
+
{\bf 1}_{\{i_6=i_4\ne 0\}}
{\bf 1}_{\{j_6=j_4\}}
{\bf 1}_{\{i_1=i_2\ne 0\}}
{\bf 1}_{\{j_1=j_2\}}
\zeta_{j_3}^{(i_3)}
\zeta_{j_5}^{(i_5)}
+
$$
$$
+{\bf 1}_{\{i_6=i_5\ne 0\}}
{\bf 1}_{\{j_6=j_5\}}
{\bf 1}_{\{i_3=i_4\ne 0\}}
{\bf 1}_{\{j_3=j_4\}}
\zeta_{j_1}^{(i_1)}
\zeta_{j_2}^{(i_2)}+
$$
$$
+
{\bf 1}_{\{i_6=i_5\ne 0\}}
{\bf 1}_{\{j_6=j_5\}}
{\bf 1}_{\{i_2=i_4\ne 0\}}
{\bf 1}_{\{j_2=j_4\}}
\zeta_{j_1}^{(i_1)}
\zeta_{j_3}^{(i_3)}
+
$$
$$
+{\bf 1}_{\{i_6=i_5\ne 0\}}
{\bf 1}_{\{j_6=j_5\}}
{\bf 1}_{\{i_2=i_3\ne 0\}}
{\bf 1}_{\{j_2=j_3\}}
\zeta_{j_1}^{(i_1)}
\zeta_{j_4}^{(i_4)}+
$$
$$
+
{\bf 1}_{\{i_6=i_5\ne 0\}}
{\bf 1}_{\{j_6=j_5\}}
{\bf 1}_{\{i_1=i_4\ne 0\}}
{\bf 1}_{\{j_1=j_4\}}
\zeta_{j_2}^{(i_2)}
\zeta_{j_3}^{(i_3)}
+
$$
$$
+{\bf 1}_{\{i_6=i_5\ne 0\}}
{\bf 1}_{\{j_6=j_5\}}
{\bf 1}_{\{i_1=i_3\ne 0\}}
{\bf 1}_{\{j_1=j_3\}}
\zeta_{j_2}^{(i_2)}
\zeta_{j_4}^{(i_4)}+
$$
$$
+
{\bf 1}_{\{i_6=i_5\ne 0\}}
{\bf 1}_{\{j_6=j_5\}}
{\bf 1}_{\{i_1=i_2\ne 0\}}
{\bf 1}_{\{j_1=j_2\}}
\zeta_{j_3}^{(i_3)}
\zeta_{j_4}^{(i_4)}-
$$
$$
-
{\bf 1}_{\{i_6=i_1\ne 0\}}
{\bf 1}_{\{j_6=j_1\}}
{\bf 1}_{\{i_2=i_5\ne 0\}}
{\bf 1}_{\{j_2=j_5\}}
{\bf 1}_{\{i_3=i_4\ne 0\}}
{\bf 1}_{\{j_3=j_4\}}-
$$
$$
-
{\bf 1}_{\{i_6=i_1\ne 0\}}
{\bf 1}_{\{j_6=j_1\}}
{\bf 1}_{\{i_2=i_4\ne 0\}}
{\bf 1}_{\{j_2=j_4\}}
{\bf 1}_{\{i_3=i_5\ne 0\}}
{\bf 1}_{\{j_3=j_5\}}-
$$
$$
-
{\bf 1}_{\{i_6=i_1\ne 0\}}
{\bf 1}_{\{j_6=j_1\}}
{\bf 1}_{\{i_2=i_3\ne 0\}}
{\bf 1}_{\{j_2=j_3\}}
{\bf 1}_{\{i_4=i_5\ne 0\}}
{\bf 1}_{\{j_4=j_5\}}-
$$
$$
-
{\bf 1}_{\{i_6=i_2\ne 0\}}
{\bf 1}_{\{j_6=j_2\}}
{\bf 1}_{\{i_1=i_5\ne 0\}}
{\bf 1}_{\{j_1=j_5\}}
{\bf 1}_{\{i_3=i_4\ne 0\}}
{\bf 1}_{\{j_3=j_4\}}-
$$
$$
-
{\bf 1}_{\{i_6=i_2\ne 0\}}
{\bf 1}_{\{j_6=j_2\}}
{\bf 1}_{\{i_1=i_4\ne 0\}}
{\bf 1}_{\{j_1=j_4\}}
{\bf 1}_{\{i_3=i_5\ne 0\}}
{\bf 1}_{\{j_3=j_5\}}-
$$
$$
-
{\bf 1}_{\{i_6=i_2\ne 0\}}
{\bf 1}_{\{j_6=j_2\}}
{\bf 1}_{\{i_1=i_3\ne 0\}}
{\bf 1}_{\{j_1=j_3\}}
{\bf 1}_{\{i_4=i_5\ne 0\}}
{\bf 1}_{\{j_4=j_5\}}-
$$
$$
-
{\bf 1}_{\{i_6=i_3\ne 0\}}
{\bf 1}_{\{j_6=j_3\}}
{\bf 1}_{\{i_1=i_5\ne 0\}}
{\bf 1}_{\{j_1=j_5\}}
{\bf 1}_{\{i_2=i_4\ne 0\}}
{\bf 1}_{\{j_2=j_4\}}-
$$
$$
-
{\bf 1}_{\{i_6=i_3\ne 0\}}
{\bf 1}_{\{j_6=j_3\}}
{\bf 1}_{\{i_1=i_4\ne 0\}}
{\bf 1}_{\{j_1=j_4\}}
{\bf 1}_{\{i_2=i_5\ne 0\}}
{\bf 1}_{\{j_2=j_5\}}-
$$
$$
-
{\bf 1}_{\{i_3=i_6\ne 0\}}
{\bf 1}_{\{j_3=j_6\}}
{\bf 1}_{\{i_1=i_2\ne 0\}}
{\bf 1}_{\{j_1=j_2\}}
{\bf 1}_{\{i_4=i_5\ne 0\}}
{\bf 1}_{\{j_4=j_5\}}-
$$
$$
-
{\bf 1}_{\{i_6=i_4\ne 0\}}
{\bf 1}_{\{j_6=j_4\}}
{\bf 1}_{\{i_1=i_5\ne 0\}}
{\bf 1}_{\{j_1=j_5\}}
{\bf 1}_{\{i_2=i_3\ne 0\}}
{\bf 1}_{\{j_2=j_3\}}-
$$
$$
-
{\bf 1}_{\{i_6=i_4\ne 0\}}
{\bf 1}_{\{j_6=j_4\}}
{\bf 1}_{\{i_1=i_3\ne 0\}}
{\bf 1}_{\{j_1=j_3\}}
{\bf 1}_{\{i_2=i_5\ne 0\}}
{\bf 1}_{\{j_2=j_5\}}-
$$
$$
-
{\bf 1}_{\{i_6=i_4\ne 0\}}
{\bf 1}_{\{j_6=j_4\}}
{\bf 1}_{\{i_1=i_2\ne 0\}}
{\bf 1}_{\{j_1=j_2\}}
{\bf 1}_{\{i_3=i_5\ne 0\}}
{\bf 1}_{\{j_3=j_5\}}-
$$
$$
-
{\bf 1}_{\{i_6=i_5\ne 0\}}
{\bf 1}_{\{j_6=j_5\}}
{\bf 1}_{\{i_1=i_4\ne 0\}}
{\bf 1}_{\{j_1=j_4\}}
{\bf 1}_{\{i_2=i_3\ne 0\}}
{\bf 1}_{\{j_2=j_3\}}-
$$
$$
-
{\bf 1}_{\{i_6=i_5\ne 0\}}
{\bf 1}_{\{j_6=j_5\}}
{\bf 1}_{\{i_1=i_2\ne 0\}}
{\bf 1}_{\{j_1=j_2\}}
{\bf 1}_{\{i_3=i_4\ne 0\}}
{\bf 1}_{\{j_3=j_4\}}-
$$
\begin{equation}
\label{a6}
\Biggl.-
{\bf 1}_{\{i_6=i_5\ne 0\}}
{\bf 1}_{\{j_6=j_5\}}
{\bf 1}_{\{i_1=i_3\ne 0\}}
{\bf 1}_{\{j_1=j_3\}}
{\bf 1}_{\{i_2=i_4\ne 0\}}
{\bf 1}_{\{j_2=j_4\}}\Biggr),
\end{equation}

\newpage
\noindent
$$
J[\psi^{(7)}]_{T,t}
=
\hbox{\vtop{\offinterlineskip\halign{
\hfil#\hfil\cr
{\rm l.i.m.}\cr
$\stackrel{}{{}_{p_1,\ldots, p_7\to \infty}}$\cr
}} }
\sum_{j_1=0}^{p_1}\ldots\sum_{j_7=0}^{p_7}
C_{j_7\ldots j_1}\Biggl(
\prod_{l=1}^7
\zeta_{j_l}^{(i_l)}
-\Biggr.
$$
$$
-
{\bf 1}_{\{i_1=i_6\ne 0,j_1=j_6\}}
\prod_{\stackrel{l=1}{{}_{l\ne 1, 6}}}^7\zeta_{j_l}^{(i_l)}
-
{\bf 1}_{\{i_2=i_6\ne 0,j_2=j_6\}}
\prod_{\stackrel{l=1}{{}_{l\ne 2, 6}}}^7\zeta_{j_l}^{(i_l)}
-
{\bf 1}_{\{i_3=i_6\ne 0,j_3=j_6\}}
\prod_{\stackrel{l=1}{{}_{l\ne 3, 6}}}^7\zeta_{j_l}^{(i_l)}-
$$
$$
-
{\bf 1}_{\{i_4=i_6\ne 0,j_4=j_6\}}
\prod_{\stackrel{l=1}{{}_{l\ne 4, 6}}}^7\zeta_{j_l}^{(i_l)}
-
{\bf 1}_{\{i_5=i_6\ne 0,j_5=j_6\}}
\prod_{\stackrel{l=1}{{}_{l\ne 5, 6}}}^7\zeta_{j_l}^{(i_l)}
-
{\bf 1}_{\{i_1=i_2\ne 0,j_1=j_2\}}
\prod_{\stackrel{l=1}{{}_{l\ne 1, 2}}}^7\zeta_{j_l}^{(i_l)}-
$$
$$
-
{\bf 1}_{\{i_1=i_3\ne 0,j_1=j_3\}}
\prod_{\stackrel{l=1}{{}_{l\ne 1, 3}}}^7\zeta_{j_l}^{(i_l)}
-
{\bf 1}_{\{i_1=i_4\ne 0,j_1=j_4\}}
\prod_{\stackrel{l=1}{{}_{l\ne 1, 4}}}^7\zeta_{j_l}^{(i_l)}
-
{\bf 1}_{\{i_1=i_5\ne 0,j_1=j_5\}}
\prod_{\stackrel{l=1}{{}_{l\ne 1, 5}}}^7\zeta_{j_l}^{(i_l)}-
$$
$$
-
{\bf 1}_{\{i_2=i_3\ne 0,j_2=j_3\}}
\prod_{\stackrel{l=1}{{}_{l\ne 2, 3}}}^7\zeta_{j_l}^{(i_l)}
-
{\bf 1}_{\{i_2=i_4\ne 0,j_2=j_4\}}
\prod_{\stackrel{l=1}{{}_{l\ne 2, 4}}}^7\zeta_{j_l}^{(i_l)}
-
{\bf 1}_{\{i_2=i_5\ne 0,j_2=j_5\}}
\prod_{\stackrel{l=1}{{}_{l\ne 2, 5}}}^7\zeta_{j_l}^{(i_l)}-
$$
$$
-
{\bf 1}_{\{i_3=i_4\ne 0,j_3=j_4\}}
\prod_{\stackrel{l=1}{{}_{l\ne 3, 4}}}^7\zeta_{j_l}^{(i_l)}
-
{\bf 1}_{\{i_3=i_5\ne 0,j_3=j_5\}}
\prod_{\stackrel{l=1}{{}_{l\ne 3, 5}}}^7\zeta_{j_l}^{(i_l)}
-
{\bf 1}_{\{i_4=i_5\ne 0,j_4=j_5\}}
\prod_{\stackrel{l=1}{{}_{l\ne 4, 5}}}^7\zeta_{j_l}^{(i_l)}-
$$
$$
-
{\bf 1}_{\{i_7=i_1\ne 0,j_7=j_1\}}
\prod_{\stackrel{l=1}{{}_{l\ne 1, 7}}}^7\zeta_{j_l}^{(i_l)}
-
{\bf 1}_{\{i_7=i_2\ne 0,j_7=j_2\}}
\prod_{\stackrel{l=1}{{}_{l\ne 2, 7}}}^7\zeta_{j_l}^{(i_l)}
-
{\bf 1}_{\{i_7=i_3\ne 0,j_7=j_3\}}
\prod_{\stackrel{l=1}{{}_{l\ne 3, 7}}}^7\zeta_{j_l}^{(i_l)}-
$$
$$
-
{\bf 1}_{\{i_7=i_4\ne 0,j_7=j_4\}}
\prod_{\stackrel{l=1}{{}_{l\ne 4, 7}}}^7\zeta_{j_l}^{(i_l)}
-
{\bf 1}_{\{i_7=i_5\ne 0,j_7=j_5\}}
\prod_{\stackrel{l=1}{{}_{l\ne 7, 5}}}^7\zeta_{j_l}^{(i_l)}
-
{\bf 1}_{\{i_7=i_6\ne 0,j_7=j_6\}}
\prod_{\stackrel{l=1}{{}_{l\ne 7, 6}}}^7\zeta_{j_l}^{(i_l)}+
$$
$$
+
{\bf 1}_{\{i_1=i_2\ne 0,j_1=j_2,i_3=i_4\ne 0,j_3=j_4\}}
\prod_{l=5,6,7}\zeta_{j_l}^{(i_l)}
+
{\bf 1}_{\{i_1=i_2\ne 0,j_1=j_2,i_3=i_5\ne 0,j_3=j_5\}}
\prod_{l=4,6,7}\zeta_{j_l}^{(i_l)}+
$$
$$
+
{\bf 1}_{\{i_1=i_2\ne 0,j_1=j_2,i_4=i_5\ne 0,j_4=j_5\}}
\prod_{l=3,6,7}\zeta_{j_l}^{(i_l)}
+
{\bf 1}_{\{i_1=i_3\ne 0,j_1=j_3,i_2=i_4\ne 0,j_2=j_4\}}
\prod_{l=5,6,7}\zeta_{j_l}^{(i_l)}+
$$
$$
+
{\bf 1}_{\{i_1=i_3\ne 0,j_1=j_3,i_2=i_5\ne 0,j_2=j_5\}}
\prod_{l=4,6,7}\zeta_{j_l}^{(i_l)}
+
{\bf 1}_{\{i_1=i_3\ne 0,j_1=j_3,i_4=i_5\ne 0,j_4=j_5\}}
\prod_{l=2,6,7}\zeta_{j_l}^{(i_l)}+
$$
$$
+
{\bf 1}_{\{i_1=i_4\ne 0,j_1=j_4,i_2=i_3\ne 0,j_2=j_3\}}
\prod_{l=5,6,7}\zeta_{j_l}^{(i_l)}
+
{\bf 1}_{\{i_1=i_4\ne 0,j_1=j_4,i_2=i_5\ne 0,j_2=j_5\}}
\prod_{l=3,6,7}\zeta_{j_l}^{(i_l)}+
$$
$$
+
{\bf 1}_{\{i_1=i_4\ne 0,j_1=j_4,i_3=i_5\ne 0,j_3=j_5\}}
\prod_{l=2,6,7}\zeta_{j_l}^{(i_l)}
+
{\bf 1}_{\{i_1=i_5\ne 0,j_1=j_5,i_2=i_3\ne 0,j_2=j_3\}}
\prod_{l=4,6,7}\zeta_{j_l}^{(i_l)}+
$$
$$
+
{\bf 1}_{\{i_1=i_5\ne 0,j_1=j_5,i_2=i_4\ne 0,j_2=j_4\}}
\prod_{l=3,6,7}\zeta_{j_l}^{(i_l)}
+
{\bf 1}_{\{i_1=i_5\ne 0,j_1=j_5,i_3=i_4\ne 0,j_3=j_4\}}
\prod_{l=2,6,7}\zeta_{j_l}^{(i_l)}+
$$
$$
+
{\bf 1}_{\{i_2=i_3\ne 0,j_2=j_3,i_4=i_5\ne 0,j_4=j_5\}}
\prod_{l=1,6,7}\zeta_{j_l}^{(i_l)}
+
{\bf 1}_{\{i_2=i_4\ne 0,j_2=j_4,i_3=i_5\ne 0,j_3=j_5\}}
\prod_{l=1,6,7}\zeta_{j_l}^{(i_l)}+
$$
$$
+
{\bf 1}_{\{i_2=i_5\ne 0,j_2=j_5,i_3=i_4\ne 0,j_3=j_4\}}
\prod_{l=1,6,7}\zeta_{j_l}^{(i_l)}
+
{\bf 1}_{\{i_6=i_1\ne 0,j_6=j_1,i_3=i_4\ne 0,j_3=j_4\}}
\prod_{l=2,5,7}\zeta_{j_l}^{(i_l)}+
$$
$$
+
{\bf 1}_{\{i_6=i_1\ne 0,j_6=j_1,i_3=i_5\ne 0,j_3=j_5\}}
\prod_{l=2,4,7}\zeta_{j_l}^{(i_l)}
+
{\bf 1}_{\{i_6=i_1\ne 0,j_6=j_1,i_2=i_5\ne 0,j_2=j_5\}}
\prod_{l=3,4,7}\zeta_{j_l}^{(i_l)}+
$$
$$
+
{\bf 1}_{\{i_6=i_1\ne 0,j_6=j_1,i_2=i_4\ne 0,j_2=j_4\}}
\prod_{l=3,5,7}\zeta_{j_l}^{(i_l)}
+
{\bf 1}_{\{i_6=i_1\ne 0,j_6=j_1,i_4=i_5\ne 0,j_4=j_5\}}
\prod_{l=2,3,7}\zeta_{j_l}^{(i_l)}+
$$
$$
+
{\bf 1}_{\{i_6=i_1\ne 0,j_6=j_1,i_2=i_3\ne 0,j_2=j_3\}}
\prod_{l=4,5,7}\zeta_{j_l}^{(i_l)}
+
{\bf 1}_{\{i_6=i_2\ne 0,j_6=j_2,i_3=i_5\ne 0,j_3=j_5\}}
\prod_{l=1,4,7}\zeta_{j_l}^{(i_l)}+
$$
$$
+
{\bf 1}_{\{i_6=i_2\ne 0,j_6=j_2,i_4=i_5\ne 0,j_4=j_5\}}
\prod_{l=1,3,7}\zeta_{j_l}^{(i_l)}
+
{\bf 1}_{\{i_6=i_2\ne 0,j_6=j_2,i_3=i_4\ne 0,j_3=j_4\}}
\prod_{l=1,5,7}\zeta_{j_l}^{(i_l)}+
$$
$$
+
{\bf 1}_{\{i_6=i_2\ne 0,j_6=j_2,i_1=i_5\ne 0,j_1=j_5\}}
\prod_{l=3,4,7}\zeta_{j_l}^{(i_l)}
+
{\bf 1}_{\{i_6=i_2\ne 0,j_6=j_2,i_1=i_4\ne 0,j_1=j_4\}}
\prod_{l=3,5,7}\zeta_{j_l}^{(i_l)}+
$$
$$
+
{\bf 1}_{\{i_6=i_2\ne 0,j_6=j_2,i_1=i_3\ne 0,j_1=j_3\}}
\prod_{l=4,5,7}\zeta_{j_l}^{(i_l)}
+
{\bf 1}_{\{i_6=i_3\ne 0,j_6=j_3,i_2=i_5\ne 0,j_2=j_5\}}
\prod_{l=1,4,7}\zeta_{j_l}^{(i_l)}+
$$
$$
+
{\bf 1}_{\{i_6=i_3\ne 0,j_6=j_3,i_4=i_5\ne 0,j_4=j_5\}}
\prod_{l=1,2,7}\zeta_{j_l}^{(i_l)}
+
{\bf 1}_{\{i_6=i_3\ne 0,j_6=j_3,i_2=i_4\ne 0,j_2=j_4\}}
\prod_{l=1,5,7}\zeta_{j_l}^{(i_l)}+
$$
$$
+
{\bf 1}_{\{i_6=i_3\ne 0,j_6=j_3,i_1=i_5\ne 0,j_1=j_5\}}
\prod_{l=2,4,7}\zeta_{j_l}^{(i_l)}
+
{\bf 1}_{\{i_6=i_3\ne 0,j_6=j_3,i_1=i_4\ne 0,j_1=j_4\}}
\prod_{l=2,5,7}\zeta_{j_l}^{(i_l)}+
$$
$$
+
{\bf 1}_{\{i_6=i_3\ne 0,j_6=j_3,i_1=i_2\ne 0,j_1=j_2\}}
\prod_{l=4,5,7}\zeta_{j_l}^{(i_l)}
+
{\bf 1}_{\{i_6=i_4\ne 0,j_6=j_4,i_3=i_5\ne 0,j_3=j_5\}}
\prod_{l=1,2,7}\zeta_{j_l}^{(i_l)}+
$$
$$
+
{\bf 1}_{\{i_6=i_4\ne 0,j_6=j_4,i_2=i_5\ne 0,j_2=j_5\}}
\prod_{l=1,3,7}\zeta_{j_l}^{(i_l)}
+
{\bf 1}_{\{i_6=i_4\ne 0,j_6=j_4,i_2=i_3\ne 0,j_2=j_3\}}
\prod_{l=1,5,7}\zeta_{j_l}^{(i_l)}+
$$
$$
+
{\bf 1}_{\{i_6=i_4\ne 0,j_6=j_4,i_1=i_5\ne 0,j_1=j_5\}}
\prod_{l=2,3,7}\zeta_{j_l}^{(i_l)}
+
{\bf 1}_{\{i_6=i_4\ne 0,j_6=j_4,i_1=i_3\ne 0,j_1=j_3\}}
\prod_{l=2,5,7}\zeta_{j_l}^{(i_l)}+
$$
$$
+
{\bf 1}_{\{i_6=i_4\ne 0,j_6=j_4,i_1=i_2\ne 0,j_1=j_2\}}
\prod_{l=3,5,7}\zeta_{j_l}^{(i_l)}
+
{\bf 1}_{\{i_6=i_5\ne 0,j_6=j_5,i_3=i_4\ne 0,j_3=j_4\}}
\prod_{l=1,2,7}\zeta_{j_l}^{(i_l)}+
$$
$$
+
{\bf 1}_{\{i_6=i_5\ne 0,j_6=j_5,i_2=i_4\ne 0,j_2=j_4\}}
\prod_{l=1,3,7}\zeta_{j_l}^{(i_l)}
+
{\bf 1}_{\{i_6=i_5\ne 0,j_6=j_5,i_2=i_3\ne 0,j_2=j_3\}}
\prod_{l=1,4,7}\zeta_{j_l}^{(i_l)}+
$$
$$
+
{\bf 1}_{\{i_6=i_5\ne 0,j_6=j_5,i_1=i_4\ne 0,j_1=j_4\}}
\prod_{l=2,3,7}\zeta_{j_l}^{(i_l)}
+
{\bf 1}_{\{i_6=i_5\ne 0,j_6=j_5,i_1=i_3\ne 0,j_1=j_3\}}
\prod_{l=2,4,7}\zeta_{j_l}^{(i_l)}+
$$
$$
+
{\bf 1}_{\{i_6=i_5\ne 0,j_6=j_5,i_1=i_2\ne 0,j_1=j_2\}}
\prod_{l=3,4,7}\zeta_{j_l}^{(i_l)}
+
{\bf 1}_{\{i_7=i_1\ne 0,j_7=j_1,i_2=i_3\ne 0,j_2=j_3\}}
\prod_{l=4,5,6}\zeta_{j_l}^{(i_l)}+
$$
$$
+
{\bf 1}_{\{i_7=i_1\ne 0,j_7=j_1,i_2=i_4\ne 0,j_2=j_4\}}
\prod_{l=3,5,6}\zeta_{j_l}^{(i_l)}
+
{\bf 1}_{\{i_7=i_1\ne 0,j_7=j_1,i_2=i_5\ne 0,j_2=j_5\}}
\prod_{l=3,4,6}\zeta_{j_l}^{(i_l)}+
$$
$$
+
{\bf 1}_{\{i_7=i_1\ne 0,j_7=j_1,i_2=i_6\ne 0,j_2=j_6\}}
\prod_{l=3,4,5}\zeta_{j_l}^{(i_l)}
+
{\bf 1}_{\{i_7=i_1\ne 0,j_7=j_1,i_3=i_4\ne 0,j_3=j_4\}}
\prod_{l=2,5,6}\zeta_{j_l}^{(i_l)}+
$$
$$
+
{\bf 1}_{\{i_7=i_1\ne 0,j_7=j_1,i_3=i_5\ne 0,j_3=j_5\}}
\prod_{l=2,4,6}\zeta_{j_l}^{(i_l)}
+
{\bf 1}_{\{i_7=i_1\ne 0,j_7=j_1,i_3=i_6\ne 0,j_3=j_6\}}
\prod_{l=2,4,5}\zeta_{j_l}^{(i_l)}+
$$
$$
+
{\bf 1}_{\{i_7=i_1\ne 0,j_7=j_1,i_4=i_5\ne 0,j_4=j_5\}}
\prod_{l=2,3,6}\zeta_{j_l}^{(i_l)}
+
{\bf 1}_{\{i_7=i_1\ne 0,j_7=j_1,i_4=i_6\ne 0,j_4=j_6\}}
\prod_{l=2,3,5}\zeta_{j_l}^{(i_l)}+
$$
$$
+
{\bf 1}_{\{i_1=i_2\ne 0,j_7=j_1,i_7=i_1\ne 0,j_5=j_6\}}
\prod_{l=2,3,4}\zeta_{j_l}^{(i_l)}
+
{\bf 1}_{\{i_7=i_2\ne 0,j_7=j_2,i_1=i_3\ne 0,j_1=j_3\}}
\prod_{l=4,5,6}\zeta_{j_l}^{(i_l)}+
$$
$$
+
{\bf 1}_{\{i_7=i_2\ne 0,j_7=j_2,i_1=i_4\ne 0,j_1=j_4\}}
\prod_{l=3,5,6}\zeta_{j_l}^{(i_l)}
+
{\bf 1}_{\{i_7=i_2\ne 0,j_7=j_2,i_1=i_5\ne 0,j_1=j_5\}}
\prod_{l=3,4,6}\zeta_{j_l}^{(i_l)}+
$$
$$
+
{\bf 1}_{\{i_7=i_2\ne 0,j_7=j_2,i_1=i_6\ne 0,j_1=j_6\}}
\prod_{l=3,4,5}\zeta_{j_l}^{(i_l)}
+
{\bf 1}_{\{i_7=i_2\ne 0,j_7=j_2,i_3=i_4\ne 0,j_3=j_4\}}
\prod_{l=1,5,6}\zeta_{j_l}^{(i_l)}+
$$
$$
+
{\bf 1}_{\{i_7=i_2\ne 0,j_7=j_2,i_3=i_5\ne 0,j_3=j_5\}}
\prod_{l=1,4,6}\zeta_{j_l}^{(i_l)}
+
{\bf 1}_{\{i_7=i_2\ne 0,j_7=j_2,i_3=i_6\ne 0,j_3=j_6\}}
\prod_{l=1,4,5}\zeta_{j_l}^{(i_l)}+
$$
$$
+
{\bf 1}_{\{i_7=i_2\ne 0,j_7=j_2,i_4=i_5\ne 0,j_4=j_5\}}
\prod_{l=1,3,6}\zeta_{j_l}^{(i_l)}
+
{\bf 1}_{\{i_7=i_2\ne 0,j_7=j_2,i_4=i_6\ne 0,j_4=j_6\}}
\prod_{l=1,3,5}\zeta_{j_l}^{(i_l)}+
$$
$$
+
{\bf 1}_{\{i_7=i_2\ne 0,j_7=j_2,i_5=i_6\ne 0,j_5=j_6\}}
\prod_{l=1,3,4}\zeta_{j_l}^{(i_l)}
+
{\bf 1}_{\{i_7=i_3\ne 0,j_7=j_3,i_1=i_2\ne 0,j_1=j_2\}}
\prod_{l=4,5,6}\zeta_{j_l}^{(i_l)}+
$$
$$
+
{\bf 1}_{\{i_7=i_3\ne 0,j_7=j_3,i_1=i_4\ne 0,j_1=j_4\}}
\prod_{l=2,3,5}\zeta_{j_l}^{(i_l)}
+
{\bf 1}_{\{i_7=i_3\ne 0,j_7=j_3,i_1=i_5\ne 0,j_1=j_5\}}
\prod_{l=2,4,6}\zeta_{j_l}^{(i_l)}+
$$
$$
+
{\bf 1}_{\{i_7=i_3\ne 0,j_7=j_3,i_1=i_6\ne 0,j_1=j_6\}}
\prod_{l=4,2,5}\zeta_{j_l}^{(i_l)}
+
{\bf 1}_{\{i_7=i_3\ne 0,j_7=j_3,i_2=i_4\ne 0,j_2=j_4\}}
\prod_{l=3,5,6}\zeta_{j_l}^{(i_l)}+
$$
$$
+
{\bf 1}_{\{i_7=i_3\ne 0,j_7=j_3,i_2=i_5\ne 0,j_2=j_5\}}
\prod_{l=1,4,6}\zeta_{j_l}^{(i_l)}
+
{\bf 1}_{\{i_7=i_3\ne 0,j_7=j_3,i_2=i_6\ne 0,j_2=j_6\}}
\prod_{l=1,4,5}\zeta_{j_l}^{(i_l)}+
$$
$$
+
{\bf 1}_{\{i_7=i_3\ne 0,j_7=j_3,i_4=i_5\ne 0,j_4=j_5\}}
\prod_{l=1,2,6}\zeta_{j_l}^{(i_l)}
+
{\bf 1}_{\{i_7=i_3\ne 0,j_7=j_3,i_4=i_6\ne 0,j_4=j_6\}}
\prod_{l=1,2,5}\zeta_{j_l}^{(i_l)}+
$$
$$
+
{\bf 1}_{\{i_7=i_3\ne 0,j_7=j_3,i_5=i_6\ne 0,j_5=j_6\}}
\prod_{l=1,2,4}\zeta_{j_l}^{(i_l)}
+
{\bf 1}_{\{i_7=i_4\ne 0,j_7=j_4,i_1=i_2\ne 0,j_1=j_2\}}
\prod_{l=3,5,6}\zeta_{j_l}^{(i_l)}+
$$
$$
+
{\bf 1}_{\{i_7=i_4\ne 0,j_7=j_4,i_1=i_3\ne 0,j_1=j_3\}}
\prod_{l=2,5,6}\zeta_{j_l}^{(i_l)}
+
{\bf 1}_{\{i_7=i_4\ne 0,j_7=j_4,i_1=i_5\ne 0,j_1=j_5\}}
\prod_{l=2,3,6}\zeta_{j_l}^{(i_l)}+
$$
$$
+
{\bf 1}_{\{i_7=i_4\ne 0,j_7=j_4,i_1=i_6\ne 0,j_1=j_6\}}
\prod_{l=2,3,5}\zeta_{j_l}^{(i_l)}
+
{\bf 1}_{\{i_7=i_4\ne 0,j_7=j_4,i_2=i_3\ne 0,j_2=j_3\}}
\prod_{l=1,5,6}\zeta_{j_l}^{(i_l)}+
$$
$$
+
{\bf 1}_{\{i_7=i_4\ne 0,j_7=j_4,i_2=i_5\ne 0,j_2=j_5\}}
\prod_{l=1,3,6}\zeta_{j_l}^{(i_l)}
+
{\bf 1}_{\{i_7=i_4\ne 0,j_7=j_4,i_2=i_6\ne 0,j_2=j_6\}}
\prod_{l=1,3,5}\zeta_{j_l}^{(i_l)}+
$$
$$
+
{\bf 1}_{\{i_7=i_4\ne 0,j_7=j_4,i_3=i_5\ne 0,j_3=j_5\}}
\prod_{l=1,2,6}\zeta_{j_l}^{(i_l)}
+
{\bf 1}_{\{i_7=i_4\ne 0,j_7=j_4,i_3=i_6\ne 0,j_3=j_6\}}
\prod_{l=1,2,5}\zeta_{j_l}^{(i_l)}+
$$
$$
+
{\bf 1}_{\{i_7=i_4\ne 0,j_7=j_4,i_5=i_6\ne 0,j_5=j_6\}}
\prod_{l=1,2,3}\zeta_{j_l}^{(i_l)}
+
{\bf 1}_{\{i_7=i_5\ne 0,j_7=j_5,i_1=i_2\ne 0,j_1=j_2\}}
\prod_{l=3,4,6}\zeta_{j_l}^{(i_l)}+
$$
$$
+
{\bf 1}_{\{i_7=i_5\ne 0,j_7=j_5,i_1=i_3\ne 0,j_1=j_3\}}
\prod_{l=2,4,6}\zeta_{j_l}^{(i_l)}
+
{\bf 1}_{\{i_7=i_5\ne 0,j_7=j_5,i_1=i_4\ne 0,j_1=j_4\}}
\prod_{l=2,3,6}\zeta_{j_l}^{(i_l)}+
$$
$$
+
{\bf 1}_{\{i_7=i_5\ne 0,j_7=j_5,i_1=i_6\ne 0,j_1=j_6\}}
\prod_{l=2,3,4}\zeta_{j_l}^{(i_l)}
+
{\bf 1}_{\{i_7=i_5\ne 0,j_7=j_5,i_2=i_3\ne 0,j_2=j_3\}}
\prod_{l=1,4,6}\zeta_{j_l}^{(i_l)}+
$$
$$
+
{\bf 1}_{\{i_7=i_5\ne 0,j_7=j_5,i_2=i_4\ne 0,j_2=j_4\}}
\prod_{l=1,3,6}\zeta_{j_l}^{(i_l)}
+
{\bf 1}_{\{i_7=i_5\ne 0,j_7=j_5,i_2=i_6\ne 0,j_2=j_6\}}
\prod_{l=1,3,5}\zeta_{j_l}^{(i_l)}+
$$
$$
+
{\bf 1}_{\{i_7=i_5\ne 0,j_7=j_5,i_3=i_4\ne 0,j_3=j_4\}}
\prod_{l=1,2,6}\zeta_{j_l}^{(i_l)}
+
{\bf 1}_{\{i_7=i_5\ne 0,j_7=j_5,i_3=i_6\ne 0,j_3=j_6\}}
\prod_{l=1,2,4}\zeta_{j_l}^{(i_l)}+
$$
$$
+
{\bf 1}_{\{i_7=i_5\ne 0,j_7=j_5,i_4=i_6\ne 0,j_4=j_6\}}
\prod_{l=1,2,3}\zeta_{j_l}^{(i_l)}
+
{\bf 1}_{\{i_7=i_6\ne 0,j_7=j_6,i_1=i_2\ne 0,j_1=j_2\}}
\prod_{l=3,4,5}\zeta_{j_l}^{(i_l)}+
$$
$$
+
{\bf 1}_{\{i_7=i_6\ne 0,j_7=j_6,i_1=i_3\ne 0,j_1=j_3\}}
\prod_{l=2,4,5}\zeta_{j_l}^{(i_l)}
+
{\bf 1}_{\{i_7=i_6\ne 0,j_7=j_6,i_1=i_4\ne 0,j_1=j_4\}}
\prod_{l=2,3,5}\zeta_{j_l}^{(i_l)}+
$$
$$
+
{\bf 1}_{\{i_7=i_6\ne 0,j_7=j_6,i_1=i_5\ne 0,j_1=j_5\}}
\prod_{l=2,3,4}\zeta_{j_l}^{(i_l)}
+
{\bf 1}_{\{i_7=i_6\ne 0,j_7=j_6,i_2=i_3\ne 0,j_2=j_3\}}
\prod_{l=1,4,5}\zeta_{j_l}^{(i_l)}+
$$
$$
+
{\bf 1}_{\{i_7=i_6\ne 0,j_7=j_6,i_2=i_4\ne 0,j_2=j_4\}}
\prod_{l=1,3,5}\zeta_{j_l}^{(i_l)}
+
{\bf 1}_{\{i_7=i_6\ne 0,j_7=j_6,i_2=i_5\ne 0,j_2=j_5\}}
\prod_{l=1,3,4}\zeta_{j_l}^{(i_l)}+
$$
$$
+
{\bf 1}_{\{i_7=i_6\ne 0,j_7=j_6,i_3=i_5\ne 0,j_3=j_5\}}
\prod_{l=1,2,4}\zeta_{j_l}^{(i_l)}
+
{\bf 1}_{\{i_7=i_6\ne 0,j_7=j_6,i_4=i_5\ne 0,j_4=j_5\}}
\prod_{l=1,2,3}\zeta_{j_l}^{(i_l)}+
$$
$$
+
{\bf 1}_{\{i_7=i_6\ne 0,j_7=j_6,i_3=i_4\ne 0,j_3=j_4\}}
\prod_{l=1,2,5}\zeta_{j_l}^{(i_l)} -
$$
$$
-\biggl(  
{\bf 1}_{\{i_2=i_3\ne 0,j_2=j_3,i_4=i_5\ne 0,j_4=j_5,i_6=i_7\ne 0,j_6=j_7\}}
+\biggr.
{\bf 1}_{\{i_2=i_3\ne 0,j_2=j_3,i_4=i_6\ne 0,j_4=j_6,i_5=i_7\ne 0,j_5=j_7\}}+
$$
$$
+
{\bf 1}_{\{i_2=i_3\ne 0,j_2=j_3,i_4=i_7\ne 0,j_4=j_7,i_5=i_6\ne 0,j_5=j_6\}}
+
{\bf 1}_{\{i_2=i_4\ne 0,j_2=j_4,i_3=i_5\ne 0,j_3=j_5,i_6=i_7\ne 0,j_6=j_7\}}+
$$
$$
+
{\bf 1}_{\{i_2=i_4\ne 0,j_2=j_4,i_3=i_6\ne 0,j_3=j_6,i_5=i_7\ne 0,j_5=j_7\}}
+
{\bf 1}_{\{i_2=i_4\ne 0,j_2=j_4,i_3=i_7\ne 0,j_3=j_7,i_5=i_6\ne 0,j_5=j_6\}}+
$$
$$
+
{\bf 1}_{\{i_2=i_5\ne 0,j_2=j_5,i_3=i_4\ne 0,j_3=j_4,i_6=i_7\ne 0,j_6=j_7\}}
+
{\bf 1}_{\{i_2=i_5\ne 0,j_2=j_5,i_3=i_6\ne 0,j_3=j_6,i_4=i_7\ne 0,j_4=j_7\}}+
$$
$$
+
{\bf 1}_{\{i_2=i_5\ne 0,j_2=j_5,i_3=i_7\ne 0,j_3=j_7,i_4=i_6\ne 0,j_4=j_6\}}
+
{\bf 1}_{\{i_2=i_6\ne 0,j_2=j_6,i_3=i_4\ne 0,j_3=j_4,i_5=i_7\ne 0,j_5=j_7\}}+
$$
$$
+
{\bf 1}_{\{i_2=i_6\ne 0,j_2=j_6,i_3=i_5\ne 0,j_3=j_5,i_4=i_7\ne 0,j_4=j_7\}}
+
{\bf 1}_{\{i_2=i_6\ne 0,j_2=j_6,i_3=i_7\ne 0,j_3=j_7,i_4=i_5\ne 0,j_4=j_5\}}+
$$
$$
+
{\bf 1}_{\{i_2=i_7\ne 0,j_2=j_7,i_3=i_4\ne 0,j_3=j_4,i_5=i_6\ne 0,j_5=j_6\}}
+
{\bf 1}_{\{i_2=i_7\ne 0,j_2=j_7,i_3=i_5\ne 0,j_3=j_5,i_4=i_6\ne 0,j_4=j_6\}}+
$$
$$
\biggl.
+{\bf 1}_{\{i_2=i_7\ne 0,j_2=j_7,i_3=i_6\ne 0,j_3=j_6,i_4=i_5\ne 0,j_4=j_5\}}
\biggr)\zeta_{j_1}^{(i_1)} -
$$
$$
-\biggl(
{\bf 1}_{\{i_1=i_3\ne 0,j_1=j_3,i_4=i_7\ne 0,j_4=j_7,i_5=i_6\ne 0,j_5=j_6\}}
+
{\bf 1}_{\{i_1=i_3\ne 0,j_1=j_3,i_4=i_5\ne 0,j_4=j_5,i_6=i_7\ne 0,j_6=j_7\}}+
\biggr.
$$
$$
+
{\bf 1}_{\{i_1=i_3\ne 0,j_1=j_3,i_4=i_6\ne 0,j_4=j_6,i_5=i_7\ne 0,j_5=j_7\}}
+
{\bf 1}_{\{i_1=i_4\ne 0,j_1=j_4,i_3=i_5\ne 0,j_3=j_5,i_6=i_7\ne 0,j_6=j_7\}}+
$$
$$
+
{\bf 1}_{\{i_1=i_4\ne 0,j_1=j_4,i_3=i_6\ne 0,j_3=j_6,i_5=i_7\ne 0,j_5=j_7\}}
+
{\bf 1}_{\{i_1=i_4\ne 0,j_1=j_4,i_3=i_7\ne 0,j_3=j_7,i_5=i_6\ne 0,j_5=j_6\}}+
$$
$$
+
{\bf 1}_{\{i_1=i_5\ne 0,j_1=j_5,i_3=i_4\ne 0,j_3=j_4,i_6=i_7\ne 0,j_6=j_7\}}
+
{\bf 1}_{\{i_1=i_5\ne 0,j_1=j_5,i_3=i_6\ne 0,j_3=j_6,i_4=i_7\ne 0,j_4=j_7\}}+
$$
$$
+
{\bf 1}_{\{i_1=i_5\ne 0,j_1=j_5,i_3=i_7\ne 0,j_3=j_7,i_4=i_6\ne 0,j_4=j_6\}}
+
{\bf 1}_{\{i_1=i_6\ne 0,j_1=j_6,i_3=i_4\ne 0,j_3=j_4,i_5=i_7\ne 0,j_5=j_7\}}+
$$
$$
+
{\bf 1}_{\{i_6=i_1\ne 0,j_6=j_1,i_3=i_5\ne 0,j_3=j_5,i_4=i_7\ne 0,j_4=j_7\}}
+
{\bf 1}_{\{i_6=i_1\ne 0,j_6=j_1,i_3=i_7\ne 0,j_3=j_7,i_4=i_5\ne 0,j_4=j_5\}}+
$$
$$
+
{\bf 1}_{\{i_1=i_7\ne 0,j_1=j_7,i_3=i_4\ne 0,j_3=j_4,i_5=i_6\ne 0,j_5=j_6\}}
+
{\bf 1}_{\{i_1=i_7\ne 0,j_1=j_7,i_3=i_5\ne 0,j_3=j_5,i_4=i_6\ne 0,j_4=j_6\}}+
$$
$$
\biggl.
+{\bf 1}_{\{i_1=i_7\ne 0,j_1=j_7,i_3=i_6\ne 0,j_3=j_6,i_4=i_5\ne 0,j_4=j_5\}}
\biggr)\zeta_{j_2}^{(i_2)} -
$$
$$
-\biggl(
{\bf 1}_{\{i_1=i_2\ne 0,j_1=j_2,i_4=i_5\ne 0,j_4=j_5,i_6=i_7\ne 0,j_6=j_7\}}
+
{\bf 1}_{\{i_1=i_2\ne 0,j_1=j_2,i_4=i_6\ne 0,j_4=j_6,i_5=i_7\ne 0,j_5=j_7\}}+
\biggr.
$$
$$
+
{\bf 1}_{\{i_1=i_2\ne 0,j_1=j_2,i_4=i_7\ne 0,j_4=j_7,i_5=i_6\ne 0,j_5=j_6\}}
+
{\bf 1}_{\{i_1=i_4\ne 0,j_1=j_4,i_2=i_5\ne 0,j_2=j_5,i_6=i_7\ne 0,j_6=j_7\}}+
$$
$$
+
{\bf 1}_{\{i_1=i_4\ne 0,j_1=j_4,i_2=i_6\ne 0,j_2=j_6,i_5=i_7\ne 0,j_5=j_7\}}
+
{\bf 1}_{\{i_1=i_4\ne 0,j_1=j_4,i_2=i_7\ne 0,j_2=j_7,i_5=i_6\ne 0,j_5=j_6\}}+
$$
$$
+
{\bf 1}_{\{i_1=i_5\ne 0,j_1=j_5,i_2=i_4\ne 0,j_2=j_4,i_6=i_7\ne 0,j_6=j_7\}}
+
{\bf 1}_{\{i_1=i_5\ne 0,j_1=j_5,i_2=i_6\ne 0,j_2=j_6,i_4=i_7\ne 0,j_4=j_7\}}+
$$
$$
+
{\bf 1}_{\{i_1=i_5\ne 0,j_1=j_5,i_2=i_7\ne 0,j_2=j_7,i_4=i_6\ne 0,j_4=j_6\}}
+
{\bf 1}_{\{i_6=i_1\ne 0,j_6=j_1,i_2=i_4\ne 0,j_2=j_4,i_5=i_7\ne 0,j_5=j_7\}}+
$$
$$
+
{\bf 1}_{\{i_6=i_1\ne 0,j_6=j_1,i_2=i_5\ne 0,j_2=j_5,i_4=i_7\ne 0,j_4=j_7\}}
+
{\bf 1}_{\{i_6=i_1\ne 0,j_6=j_1,i_2=i_7\ne 0,j_2=j_7,i_4=i_5\ne 0,j_4=j_5\}}+
$$
$$
+
{\bf 1}_{\{i_1=i_7\ne 0,j_1=j_7,i_2=i_4\ne 0,j_2=j_4,i_5=i_6\ne 0,j_5=j_6\}}
+
{\bf 1}_{\{i_1=i_7\ne 0,j_1=j_7,i_2=i_5\ne 0,j_2=j_5,i_4=i_6\ne 0,j_4=j_6\}}+
$$
$$
\biggl.+
{\bf 1}_{\{i_1=i_7\ne 0,j_1=j_7,i_2=i_6\ne 0,j_2=j_6,i_4=i_5\ne 0,j_4=j_5\}}
\biggr)\zeta_{j_3}^{(i_3)} -
$$
$$
-\biggl(
{\bf 1}_{\{i_1=i_2\ne 0,j_1=j_2,i_3=i_5\ne 0,j_3=j_5,i_6=i_7\ne 0,j_6=j_7\}}
+\biggr.
{\bf 1}_{\{i_1=i_2\ne 0,j_1=j_2,i_3=i_6\ne 0,j_3=j_6,i_5=i_7\ne 0,j_5=j_7\}}+
$$
$$
+
{\bf 1}_{\{i_1=i_2\ne 0,j_1=j_2,i_3=i_7\ne 0,j_3=j_7,i_5=i_6\ne 0,j_5=j_6\}}
+
{\bf 1}_{\{i_1=i_3\ne 0,j_1=j_3,i_2=i_5\ne 0,j_2=j_5,i_6=i_7\ne 0,j_6=j_7\}}+
$$
$$
+
{\bf 1}_{\{i_1=i_3\ne 0,j_1=j_3,i_2=i_6\ne 0,j_2=j_6,i_5=i_7\ne 0,j_5=j_7\}}
+
{\bf 1}_{\{i_1=i_3\ne 0,j_1=j_3,i_2=i_7\ne 0,j_2=j_7,i_5=i_6\ne 0,j_5=j_6\}}+
$$
$$
+
{\bf 1}_{\{i_1=i_5\ne 0,j_1=j_5,i_2=i_3\ne 0,j_2=j_3,i_6=i_7\ne 0,j_6=j_7\}}
+
{\bf 1}_{\{i_1=i_5\ne 0,j_1=j_5,i_2=i_6\ne 0,j_2=j_6,i_3=i_7\ne 0,j_3=j_7\}}+
$$
$$
+
{\bf 1}_{\{i_1=i_5\ne 0,j_1=j_5,i_2=i_7\ne 0,j_2=j_7,i_3=i_6\ne 0,j_3=j_6\}}
+
{\bf 1}_{\{i_6=i_1\ne 0,j_6=j_1,i_2=i_3\ne 0,j_2=j_3,i_5=i_7\ne 0,j_5=j_7\}}+
$$
$$
+
{\bf 1}_{\{i_6=i_1\ne 0,j_6=j_1,i_2=i_5\ne 0,j_2=j_5,i_3=i_7\ne 0,j_3=j_7\}}
+
{\bf 1}_{\{i_6=i_1\ne 0,j_6=j_1,i_2=i_7\ne 0,j_2=j_7,i_3=i_5\ne 0,j_3=j_5\}}+
$$
$$
+
{\bf 1}_{\{i_7=i_1\ne 0,j_7=j_1,i_2=i_3\ne 0,j_2=j_3,i_5=i_6\ne 0,j_5=j_6\}}
+
{\bf 1}_{\{i_7=i_1\ne 0,j_7=j_1,i_2=i_5\ne 0,j_2=j_5,i_3=i_6\ne 0,j_3=j_6\}}+
$$
$$
\biggl.
+{\bf 1}_{\{i_7=i_1\ne 0,j_7=j_1,i_2=i_6\ne 0,j_2=j_6,i_3=i_5\ne 0,j_3=j_5\}}
\biggr)\zeta_{j_4}^{(i_4)} -
$$
$$
-\biggl(
{\bf 1}_{\{i_1=i_2\ne 0,j_1=j_2,i_3=i_4\ne 0,j_3=j_4,i_6=i_7\ne 0,j_6=j_7\}}
+\biggr.
{\bf 1}_{\{i_1=i_2\ne 0,j_1=j_2,i_3=i_6\ne 0,j_3=j_6,i_4=i_7\ne 0,j_4=j_7\}}+
$$
$$
+
{\bf 1}_{\{i_1=i_2\ne 0,j_1=j_2,i_3=i_7\ne 0,j_3=j_7,i_4=i_6\ne 0,j_4=j_6\}}
+
{\bf 1}_{\{i_1=i_3\ne 0,j_1=j_3,i_2=i_4\ne 0,j_2=j_4,i_6=i_7\ne 0,j_6=j_7\}}+
$$
$$
+
{\bf 1}_{\{i_1=i_3\ne 0,j_1=j_3,i_2=i_6\ne 0,j_2=j_6,i_4=i_7\ne 0,j_4=j_7\}}
+
{\bf 1}_{\{i_1=i_3\ne 0,j_1=j_3,i_2=i_7\ne 0,j_2=j_7,i_4=i_6\ne 0,j_4=j_6\}}+
$$
$$
+
{\bf 1}_{\{i_1=i_4\ne 0,j_1=j_4,i_2=i_3\ne 0,j_2=j_3,i_6=i_7\ne 0,j_6=j_7\}}
+
{\bf 1}_{\{i_1=i_4\ne 0,j_1=j_4,i_2=i_6\ne 0,j_2=j_6,i_3=i_7\ne 0,j_3=j_7\}}+
$$
$$
+
{\bf 1}_{\{i_1=i_4\ne 0,j_1=j_4,i_2=i_7\ne 0,j_2=j_7,i_3=i_6\ne 0,j_3=j_6\}}
+
{\bf 1}_{\{i_6=i_1\ne 0,j_6=j_1,i_2=i_3\ne 0,j_2=j_3,i_4=i_7\ne 0,j_4=j_7\}}+
$$
$$
+
{\bf 1}_{\{i_6=i_1\ne 0,j_6=j_1,i_2=i_4\ne 0,j_2=j_4,i_3=i_7\ne 0,j_3=j_7\}}
+
{\bf 1}_{\{i_6=i_1\ne 0,j_6=j_1,i_2=i_7\ne 0,j_2=j_7,i_3=i_4\ne 0,j_3=j_4\}}+
$$
$$
+
{\bf 1}_{\{i_1=i_7\ne 0,j_1=j_7,i_2=i_3\ne 0,j_2=j_3,i_4=i_6\ne 0,j_4=j_6\}}
+
{\bf 1}_{\{i_1=i_7\ne 0,j_1=j_7,i_2=i_4\ne 0,j_2=j_4,i_3=i_6\ne 0,j_3=j_6\}}+
$$
$$
\biggl.+
{\bf 1}_{\{i_7=i_1\ne 0,j_7=j_1,i_2=i_6\ne 0,j_2=j_6,i_3=i_4\ne 0,j_3=j_4\}}
\biggr)\zeta_{j_5}^{(i_5)} -
$$
$$
-\biggl(
{\bf 1}_{\{i_1=i_2\ne 0,j_1=j_2,i_3=i_4\ne 0,j_3=j_4,i_5=i_7\ne 0,j_5=j_7\}}
+\biggr.
{\bf 1}_{\{i_1=i_2\ne 0,j_1=j_2,i_3=i_5\ne 0,j_3=j_5,i_4=i_7\ne 0,j_4=j_7\}}+
$$
$$
+
{\bf 1}_{\{i_1=i_2\ne 0,j_1=j_2,i_3=i_7\ne 0,j_3=j_7,i_4=i_5\ne 0,j_4=j_5\}}
+
{\bf 1}_{\{i_1=i_3\ne 0,j_1=j_3,i_2=i_4\ne 0,j_2=j_4,i_5=i_7\ne 0,j_5=j_7\}}+
$$
$$
+
{\bf 1}_{\{i_1=i_3\ne 0,j_1=j_3,i_2=i_5\ne 0,j_2=j_5,i_4=i_7\ne 0,j_4=j_7\}}
+
{\bf 1}_{\{i_1=i_3\ne 0,j_1=j_3,i_2=i_7\ne 0,j_2=j_7,i_4=i_5\ne 0,j_4=j_5\}}+
$$
$$
+
{\bf 1}_{\{i_1=i_4\ne 0,j_1=j_4,i_2=i_3\ne 0,j_2=j_3,i_5=i_7\ne 0,j_5=j_7\}}
+
{\bf 1}_{\{i_1=i_4\ne 0,j_1=j_4,i_2=i_5\ne 0,j_2=j_5,i_3=i_7\ne 0,j_3=j_7\}}+
$$
$$
+
{\bf 1}_{\{i_1=i_4\ne 0,j_1=j_4,i_2=i_7\ne 0,j_2=j_7,i_3=i_5\ne 0,j_3=j_5\}}
+
{\bf 1}_{\{i_1=i_5\ne 0,j_1=j_5,i_2=i_3\ne 0,j_2=j_3,i_4=i_7\ne 0,j_4=j_7\}}+
$$
$$
+
{\bf 1}_{\{i_1=i_5\ne 0,j_1=j_5,i_2=i_4\ne 0,j_2=j_4,i_3=i_7\ne 0,j_3=j_7\}}
+
{\bf 1}_{\{i_1=i_5\ne 0,j_1=j_5,i_2=i_7\ne 0,j_2=j_7,i_3=i_4\ne 0,j_3=j_4\}}+
$$
$$
+
{\bf 1}_{\{i_7=i_1\ne 0,j_7=j_1,i_2=i_3\ne 0,j_2=j_3,i_4=i_5\ne 0,j_4=j_5\}}
+
{\bf 1}_{\{i_7=i_1\ne 0,j_7=j_1,i_2=i_4\ne 0,j_2=j_4,i_3=i_5\ne 0,j_3=j_5\}}+
$$
$$
\biggl.+
{\bf 1}_{\{i_7=i_1\ne 0,j_7=j_1,i_2=i_5\ne 0,j_2=j_5,i_3=i_4\ne 0,j_3=j_4\}}
\biggr)\zeta_{j_6}^{(i_6)} - 
$$
$$
-\biggl(
{\bf 1}_{\{i_1=i_2\ne 0,j_1=j_2,i_3=i_4\ne 0,j_3=j_4,i_5=i_6\ne 0,j_5=j_6\}}
+\biggr.
{\bf 1}_{\{i_1=i_2\ne 0,j_1=j_2,i_3=i_5\ne 0,j_3=j_5,i_4=i_6\ne 0,j_4=j_6\}}+
$$
$$
+
{\bf 1}_{\{i_1=i_2\ne 0,j_1=j_2,i_3=i_6\ne 0,j_3=j_6,i_4=i_5\ne 0,j_4=j_5\}}
+
{\bf 1}_{\{i_1=i_3\ne 0,j_1=j_3,i_2=i_4\ne 0,j_2=j_4,i_5=i_6\ne 0,j_5=j_6\}}+
$$
$$
+
{\bf 1}_{\{i_1=i_3\ne 0,j_1=j_3,i_2=i_5\ne 0,j_2=j_5,i_4=i_6\ne 0,j_4=j_6\}}
+
{\bf 1}_{\{i_1=i_3\ne 0,j_1=j_3,i_2=i_6\ne 0,j_2=j_6,i_4=i_5\ne 0,j_4=j_5\}}+
$$
$$
+
{\bf 1}_{\{i_4=i_1\ne 0,j_4=j_1,i_2=i_3\ne 0,j_2=j_3,i_5=i_6\ne 0,j_5=j_6\}}
+
{\bf 1}_{\{i_4=i_1\ne 0,j_4=j_1,i_2=i_5\ne 0,j_2=j_5,i_3=i_6\ne 0,j_3=j_6\}}+
$$
$$
+
{\bf 1}_{\{i_4=i_1\ne 0,j_4=j_1,i_2=i_6\ne 0,j_2=j_6,i_3=i_5\ne 0,j_3=j_5\}}
+
{\bf 1}_{\{i_5=i_1\ne 0,j_5=j_1,i_2=i_3\ne 0,j_2=j_3,i_4=i_6\ne 0,j_4=j_6\}}+
$$
$$
+
{\bf 1}_{\{i_5=i_1\ne 0,j_5=j_1,i_2=i_4\ne 0,j_2=j_4,i_3=i_6\ne 0,j_3=j_6\}}
+
{\bf 1}_{\{i_5=i_1\ne 0,j_5=j_1,i_2=i_6\ne 0,j_2=j_6,i_3=i_4\ne 0,j_3=j_4\}}+
$$
$$
+
{\bf 1}_{\{i_6=i_1\ne 0,j_6=j_1,i_2=i_3\ne 0,j_2=j_3,i_4=i_5\ne 0,j_4=j_5\}}
+
{\bf 1}_{\{i_6=i_1\ne 0,j_6=j_1,i_2=i_4\ne 0,j_2=j_4,i_3=i_5\ne 0,j_3=j_5\}}+
$$
\begin{equation}
\label{a7}
\Biggl.\biggl.
+{\bf 1}_{\{i_6=i_1\ne 0,j_6=j_1,i_2=i_5\ne 0,j_2=j_5,i_3=i_4\ne 0,j_3=j_4\}}
\biggr)\zeta_{j_7}^{(i_7)}\Biggr),
\end{equation}

\vspace{2mm}

\noindent
where ${\bf 1}_A$ is the indicator of the set $A$.

\subsection{Expansion of Iterated It\^{o} Stochastic Integrals 
of Multiplicity $k$ ($k\in{\bf N}$) Based on Theorem 1.1}

Consider a generalization of the formulas (\ref{a1})--(\ref{a7}) 
for the case of arbitrary multiplicity $k$ for 
$J[\psi^{(k)}]_{T,t}$.
In order to do this, let us
consider the unordered
set $\{1, 2, \ldots, k\}$ 
and separate it into two parts:
the first part consists of $r$ unordered 
pairs (sequence order of these pairs is also unimportant) and the 
second one consists of the 
remaining $k-2r$ numbers.
So, we have
\begin{equation}
\label{leto5007}
(\{
\underbrace{\{g_1, g_2\}, \ldots, 
\{g_{2r-1}, g_{2r}\}}_{\small{\hbox{part 1}}}
\},
\{\underbrace{q_1, \ldots, q_{k-2r}}_{\small{\hbox{part 2}}}
\}),
\end{equation}
where 
$\{g_1, g_2, \ldots, 
g_{2r-1}, g_{2r}, q_1, \ldots, q_{k-2r}\}=\{1, 2, \ldots, k\},$
braces   
mean an unordered 
set, and pa\-ren\-the\-ses mean an ordered set.

We will say that (\ref{leto5007}) is a partition 
and consider the sum with respect to all possible
partitions
\begin{equation}
\label{leto5008}
\sum_{\stackrel{(\{\{g_1, g_2\}, \ldots, 
\{g_{2r-1}, g_{2r}\}\}, \{q_1, \ldots, q_{k-2r}\})}
{{}_{\{g_1, g_2, \ldots, 
g_{2r-1}, g_{2r}, q_1, \ldots, q_{k-2r}\}=\{1, 2, \ldots, k\}}}}
a_{g_1 g_2, \ldots, 
g_{2r-1} g_{2r}, q_1 \ldots q_{k-2r}},
\end{equation}
where $a_{g_1 g_2, \ldots, 
g_{2r-1} g_{2r}, q_1 \ldots q_{k-2r}}\in {\bf R}$.

Below there are several examples of sums in the form (\ref{leto5008})
$$
\sum_{\stackrel{(\{g_1, g_2\})}{{}_{\{g_1, g_2\}=\{1, 2\}}}}
a_{g_1 g_2}=a_{12},
$$
$$
\sum_{\stackrel{(\{\{g_1, g_2\}, \{g_3, g_4\}\})}
{{}_{\{g_1, g_2, g_3, g_4\}=\{1, 2, 3, 4\}}}}
a_{g_1 g_2, g_3 g_4}=a_{12,34} + a_{13,24} + a_{23,14},
$$

\vspace{-2mm}
$$
\sum_{\stackrel{(\{g_1, g_2\}, \{q_1, q_{2}\})}
{{}_{\{g_1, g_2, q_1, q_{2}\}=\{1, 2, 3, 4\}}}}
a_{g_1 g_2, q_1 q_{2}}=
$$
$$
=a_{12,34}+a_{13,24}+a_{14,23}
+a_{23,14}+a_{24,13}+a_{34,12},
$$

\vspace{-2mm}
$$
\sum_{\stackrel{(\{g_1, g_2\}, \{q_1, q_{2}, q_3\})}
{{}_{\{g_1, g_2, q_1, q_{2}, q_3\}=\{1, 2, 3, 4, 5\}}}}
a_{g_1 g_2, q_1 q_{2}q_3}
=
$$
$$
=a_{12,345}+a_{13,245}+a_{14,235}
+a_{15,234}+a_{23,145}+a_{24,135}+
$$
$$
+a_{25,134}+a_{34,125}+a_{35,124}+a_{45,123},
$$

\vspace{-2mm}
$$
\sum_{\stackrel{(\{\{g_1, g_2\}, \{g_3, g_{4}\}\}, \{q_1\})}
{{}_{\{g_1, g_2, g_3, g_{4}, q_1\}=\{1, 2, 3, 4, 5\}}}}
a_{g_1 g_2, g_3 g_{4},q_1}
=
$$
$$
=
a_{12,34,5}+a_{13,24,5}+a_{14,23,5}+
a_{12,35,4}+a_{13,25,4}+a_{15,23,4}+
a_{12,54,3}+a_{15,24,3}+
$$
$$
+a_{14,25,3}+a_{15,34,2}+a_{13,54,2}+a_{14,53,2}+
a_{52,34,1}+a_{53,24,1}+a_{54,23,1}.
$$

\vspace{4mm}

Now we can formulate Theorem 1.1 
(see (\ref{tyyy})) 
using alternative form.

\footnotetext[3]{The connection of 
formulas (\ref{a1})--(\ref{a7}), (\ref{leto6000}) 
with Hermite polynomials is studied in Sect.~1.10, 1.11 (see Theorems~1.14--1.17).}

{\bf Theorem 1.2}${}^3$ \cite{4} (2009) (also see
\cite{5}-\cite{12aa}, \cite{art-7}, \cite{arxiv-1}, \cite{arxiv-11},
\cite{arxiv-20}, \cite{arxiv-21}). {\it Under the conditions of Theorem {\rm 1.1} 
the following expansion 
$$
J[\psi^{(k)}]_{T,t}=
\hbox{\vtop{\offinterlineskip\halign{
\hfil#\hfil\cr
{\rm l.i.m.}\cr
$\stackrel{}{{}_{p_1,\ldots,p_k\to \infty}}$\cr
}} }
\sum\limits_{j_1=0}^{p_1}\ldots
\sum\limits_{j_k=0}^{p_k}
C_{j_k\ldots j_1}\Biggl(
\prod_{l=1}^k\zeta_{j_l}^{(i_l)}+\sum\limits_{r=1}^{[k/2]}
(-1)^r \times
\Biggr.
$$

\vspace{-2mm}
\begin{equation}
\label{leto6000}
\times
\sum_{\stackrel{(\{\{g_1, g_2\}, \ldots, 
\{g_{2r-1}, g_{2r}\}\}, \{q_1, \ldots, q_{k-2r}\})}
{{}_{\{g_1, g_2, \ldots, 
g_{2r-1}, g_{2r}, q_1, \ldots, q_{k-2r}\}=\{1, 2, \ldots, k\}}}}
\prod\limits_{s=1}^r
{\bf 1}_{\{i_{g_{{}_{2s-1}}}=~i_{g_{{}_{2s}}}\ne 0\}}
\Biggl.{\bf 1}_{\{j_{g_{{}_{2s-1}}}=~j_{g_{{}_{2s}}}\}}
\prod_{l=1}^{k-2r}\zeta_{j_{q_l}}^{(i_{q_l})}\Biggr)
\end{equation}

\vspace{1mm}
\noindent
con\-verg\-ing in the mean-square sense is valid,
where $[x]$ is an integer part of a real number $x,$ $\prod\limits_{\emptyset}
\stackrel{\sf def}{=}1,$ $\sum\limits_{\emptyset}
\stackrel{\sf def}{=}0;$ another notations are the same as in Theorem~{\rm 1.1.}}

\vspace{1mm}

{\bf Proof.}\ The equality (\ref{leto6000}) will be proved by induction in Sect.~1.14
(see the proof of Theorem~1.23).

In particular, from (\ref{leto6000}) for $k=5$ we obtain
$$
J[\psi^{(5)}]_{T,t}=
\hbox{\vtop{\offinterlineskip\halign{
\hfil#\hfil\cr
{\rm l.i.m.}\cr
$\stackrel{}{{}_{p_1,\ldots,p_5\to \infty}}$\cr
}} }\sum_{j_1=0}^{p_1}\ldots\sum_{j_5=0}^{p_5}
C_{j_5\ldots j_1}\Biggl(
\prod_{l=1}^5\zeta_{j_l}^{(i_l)}-\Biggr.
$$
$$
-
\sum\limits_{\stackrel{(\{g_1, g_2\}, \{q_1, q_{2}, q_3\})}
{{}_{\{g_1, g_2, q_{1}, q_{2}, q_3\}=\{1, 2, 3, 4, 5\}}}}
{\bf 1}_{\{i_{g_{{}_{1}}}=~i_{g_{{}_{2}}}\ne 0\}}
{\bf 1}_{\{j_{g_{{}_{1}}}=~j_{g_{{}_{2}}}\}}
\prod_{l=1}^{3}\zeta_{j_{q_l}}^{(i_{q_l})}+
$$
$$
+
\sum_{\stackrel{(\{\{g_1, g_2\}, 
\{g_{3}, g_{4}\}\}, \{q_1\})}
{{}_{\{g_1, g_2, g_{3}, g_{4}, q_1\}=\{1, 2, 3, 4, 5\}}}}
{\bf 1}_{\{i_{g_{{}_{1}}}=~i_{g_{{}_{2}}}\ne 0\}}
{\bf 1}_{\{j_{g_{{}_{1}}}=~j_{g_{{}_{2}}}\}}
\Biggl.{\bf 1}_{\{i_{g_{{}_{3}}}=~i_{g_{{}_{4}}}\ne 0\}}
{\bf 1}_{\{j_{g_{{}_{3}}}=~j_{g_{{}_{4}}}\}}
\zeta_{j_{q_1}}^{(i_{q_1})}\Biggr).
$$

\vspace{2mm}
\noindent
The last equality obviously agrees with
(\ref{a5}).

It is now appropriate
to make a remark about the structure 
of the formulas (\ref{a1})--(\ref{a7}) and           
(\ref{leto6000}). Using (\ref{novoe1}), (\ref{drdr1}),
(\ref{a1})--(\ref{a7}), (\ref{leto6000}), we obtain
$$
J'[\phi_{j_1}\ldots \phi_{j_k}]_{T,t}^{(i_1\ldots i_k)}=
\prod_{l=1}^k\zeta_{j_l}^{(i_l)}+\sum\limits_{r=1}^{[k/2]}
(-1)^r \times
\Biggr.
$$

\vspace{-2mm}
\begin{equation}
\label{leto60001a1b}
\times
\sum_{\stackrel{(\{\{g_1, g_2\}, \ldots, 
\{g_{2r-1}, g_{2r}\}\}, \{q_1, \ldots, q_{k-2r}\})}
{{}_{\{g_1, g_2, \ldots, 
g_{2r-1}, g_{2r}, q_1, \ldots, q_{k-2r}\}=\{1, 2, \ldots, k\}}}}
\prod\limits_{s=1}^r
{\bf 1}_{\{i_{g_{{}_{2s-1}}}=~i_{g_{{}_{2s}}}\ne 0\}}
\Biggl.{\bf 1}_{\{j_{g_{{}_{2s-1}}}=~j_{g_{{}_{2s}}}\}}
\prod_{l=1}^{k-2r}\zeta_{j_{q_l}}^{(i_{q_l})}\Biggr)
\end{equation}

\vspace{1mm}
\noindent
w.~p.~1, where the multiple stochastic integral
$J'[\phi_{j_1}\ldots \phi_{j_k}]_{T,t}^{(i_1\ldots i_k)}$
is defined by (\ref{mult11}); another notations in 
(\ref{leto60001a1b}) are the same as in Theorem~1.2.

The stochastic integral with respect to the scalar standard Wiener process
($i_1=\ldots=i_k\ne 0$)
and similar to (\ref{mult11}) was considered in \cite{ito1951} (1951)
and is called the multiple Wiener stochastic integral \cite{ito1951}.
Note that $\Phi(t_1,\ldots,t_k)\in L_2([t, T]^k)$ in \cite{ito1951}
(this case will be considered in Sect.~1.11--1.14).

As we will see in Sect.~1.10, 1.11, 1.14, the expression on the right-hand side
of (\ref{leto60001a1b}) is the Wick polynomial with arguments
$\zeta_{j_1}^{(i_1)},\ldots,\zeta_{j_k}^{(i_k)}.$
Moreover, the given expression is an explicit representation
of the Wick polynomial, in contrast to its 
representation in the form of a product of Hermite
polynomials (see Sect.~1.10, 1.11, 1.14) or its another representation
(or definition) using a recurrence relation (see (\ref{recur1})).

To best of our knowledge, the representation of the multiple Wiener
stochastic integral in the form of a Wick polynomial 
(see (\ref{leto60001a1b})) for the case of a multidimensional
Wiener process ($i_1,\ldots,i_k=0,1,\ldots,m$)
and the case $j_1,\ldots,j_k=0,1,2,\ldots $ was first obtained
in our monographs \cite{1} (2006), \cite{3} (2007), and
\cite{4} (2009). More precisely, 
the formula (\ref{leto60001a1b}) is obtained
in our monograph \cite{4} (2009) as part of the formula
(5.30) (see \cite{4}, p.~220).
Moreover, partiular cases $k=1,\ldots,5$ (see (\ref{a1})--(\ref{a5})) 
of the formula (\ref{leto60001a1b})
were obtained in \cite{1} (2006) as parts of the formulas 
on the pages 243-244 and partiular cases $k=1,\ldots,7$ 
(see (\ref{a1})--(\ref{a7})) 
of the formula (\ref{leto60001a1b})
were obtained in \cite{3} (2007) as parts of the formulas 
on the pages 208-218.

The indicated formulas are obtained for the case
when $\psi_1(\tau),\ldots,\psi_k(\tau)$ 
are continuous nonrandom functions on the interval $[t, T]$
and 
$\{\phi_j(x)\}_{j=0}^{\infty}$ is a complete orthonormal system  
of piecewise continuous functions in the space $L_2([t,T])$ 
(see Sect.~1.1.7 and \cite{1} (2006), \cite{3} (2007), and
\cite{4} (2009)). Note that the generality
of the above results is even too great when
applied to the numerical integration
of It\^{o} stochastic differential equations.

It should be noted that in \cite{fox} (1987)
an $L_2$--version of the formula (\ref{leto60001a1b})
was obtained, but only for the special case
$j_1=\ldots=j_k$. The above result in \cite{fox} (Proposition~5.1)
is obtained using diagrams, i.e. (unlike our results)
in an implicit form
(see Sect.~1.14 (below Remark~1.18) for details).

Let us turn to the comparison
of the formula (\ref{leto60001a1b}) with another interesting work \cite{major2} (2019).
An $L_2$-version of (\ref{leto60001a1b}) was obtained in \cite{major2} in terms 
of Wick polynomials and for the case of vector valued random measures 
(see \cite{major2}, Theorem~7.2, p.~69). In earlier works of this author
(see for example \cite{major1}) only the case of scalar valued random measures 
was considered (see Sect.~1.14 (below Remark~1.18) for details).

In Sect.~1.14 (Theorems~1.22, 1.23) 
we consider $L_2$--versions of the formula 
(\ref{leto60001a1b}). At that, to prove 
Theorems~1.22 and 1.23 we use only
the It\^{o} formula, in contrast to the diagram method from 
\cite{major2}.

\subsection{Comparison of Theorem 1.2 with the Representations of Iterated
It\^{o} Stochastic Integrals Based on Hermite Polynomials}

Note that the correctness of the formulas (\ref{a1})--(\ref{a7}) 
can be 
verified 
in the following way.
If 
$i_1=\ldots=i_7=i=1,\ldots,m$
and $\psi_1(s),\ldots,\psi_7(s)\equiv \psi(s)$,
then we can derive from (\ref{a1})--(\ref{a7}) 
\cite{2}-\cite{12aa}, \cite{arxiv-1}
the well known
equalities

\vspace{-2mm}
$$
J[\psi^{(1)}]_{T,t}
=\frac{1}{1!}\delta_{T,t},
$$
$$
J[\psi^{(2)}]_{T,t}
=\frac{1}{2!}\left(\delta^2_{T,t}-\Delta_{T,t}\right),\
$$

\vspace{-2mm}
$$
J[\psi^{(3)}]_{T,t}
=\frac{1}{3!}\left(\delta_{T,t}^3-3\delta_{T,t}\Delta_{T,t}\right),
$$

\vspace{-2mm}
$$
J[\psi^{(4)}]_{T,t}
=\frac{1}{4!}\left(\delta^4_{T,t}-6\delta_{T,t}^2\Delta_{T,t}
+3\Delta^2_{T,t}\right),\
$$

\vspace{-2mm}
$$
J[\psi^{(5)}]_{T,t}
=\frac{1}{5!}\left(\delta^5_{T,t}-10\delta_{T,t}^3\Delta_{T,t}
+15\delta_{T,t}\Delta^2_{T,t}\right),
$$

\vspace{-2mm}
$$
J[\psi^{(6)}]_{T,t}
=\frac{1}{6!}\left(\delta^6_{T,t}-15\delta_{T,t}^4\Delta_{T,t}
+45\delta_{T,t}^2\Delta^2_{T,t}-15\Delta_{T,t}^3\right),
$$

\vspace{-2mm}
$$
J[\psi^{(7)}]_{T,t}
=\frac{1}{7!}\left(\delta^7_{T,t}-21\delta_{T,t}^5\Delta_{T,t}
+105\delta_{T,t}^3\Delta^2_{T,t}-105\delta_{T,t}\Delta_{T,t}^3\right)
$$

\vspace{3mm}
\noindent
w. p. 1, where
$$
\delta_{T,t}=\int\limits_t^T\psi(s)d{\bf w}_s^{(i)},\ \ \ 
\Delta_{T,t}=\int\limits_t^T\psi^2(s)ds,
$$

\vspace{2mm}
\noindent
which can be independently  
obtained using the It\^{o} formula and Hermite polynomials
\cite{Ch}.

When $k=1$ everything is evident. Let us consider the cases  
$k=2$ and $k=3$ in detail.
When $k=2$ and $p_1=p_2=p$ 
we have (see (\ref{a2})) \cite{2}-\cite{12aa}, \cite{arxiv-1}
$$
J[\psi^{(2)}]_{T,t}=
\hbox{\vtop{\offinterlineskip\halign{
\hfil#\hfil\cr
{\rm l.i.m.}\cr
$\stackrel{}{{}_{p\to \infty}}$\cr
}} }\left(
\sum_{j_1,j_2=0}^{p}
C_{j_2j_1}\zeta_{j_1}^{(i)}\zeta_{j_2}^{(i)}-
\sum_{j_1=0}^{p}
C_{j_1j_1}\right)=
$$
$$
=
\hbox{\vtop{\offinterlineskip\halign{
\hfil#\hfil\cr
{\rm l.i.m.}\cr
$\stackrel{}{{}_{p\to \infty}}$\cr
}} }\left(
\sum_{j_1=0}^{p}\sum_{j_2=0}^{j_1-1}\biggl(
C_{j_2j_1}+C_{j_1j_2}\biggr)
\zeta_{j_1}^{(i)}\zeta_{j_2}^{(i)}+
\sum_{j_1=0}^{p}
C_{j_1j_1}\left(\left(\zeta_{j_1}^{(i)}\right)^2-1\right)\right)
=
$$
$$
=\hbox{\vtop{\offinterlineskip\halign{
\hfil#\hfil\cr
{\rm l.i.m.}\cr
$\stackrel{}{{}_{p\to \infty}}$\cr
}} }\left(
\sum_{j_1=0}^{p}\sum_{j_2=0}^{j_1-1}
C_{j_1}C_{j_2}
\zeta_{j_1}^{(i)}\zeta_{j_2}^{(i)}+
\frac{1}{2}\sum_{j_1=0}^{p}
C_{j_1}^2\left(\left(\zeta_{j_1}^{(i)}\right)^2-1\right)\right)=
$$
$$
=\hbox{\vtop{\offinterlineskip\halign{
\hfil#\hfil\cr
{\rm l.i.m.}\cr
$\stackrel{}{{}_{p\to \infty}}$\cr
}} }\left(
\frac{1}{2}
\sum_{\stackrel{j_1,j_2=0}{{}_{j_1\ne j_2}}}^{p}
C_{j_1}C_{j_2}
\zeta_{j_1}^{(i)}\zeta_{j_2}^{(i)}+
\frac{1}{2}\sum_{j_1=0}^{p}
C_{j_1}^2\left(\left(\zeta_{j_1}^{(i)}\right)^2-1\right)\right)=
$$
$$
=
\hbox{\vtop{\offinterlineskip\halign{
\hfil#\hfil\cr
{\rm l.i.m.}\cr
$\stackrel{}{{}_{p\to \infty}}$\cr
}} }\left(
\frac{1}{2}
\left(\sum_{j_1=0}^{p}
C_{j_1}\zeta_{j_1}^{(i)}\right)^2-
\frac{1}{2}\sum_{j_1=0}^{p}
C_{j_1}^2\right)
$$
\begin{equation}
\label{pipi20}
=\frac{1}{2!}\left(\delta^2_{T,t}-\Delta_{T,t}\right)\ \ \ \hbox{w.~p.~1}.
\end{equation}

Let us explain the last 
step in (\ref{pipi20}). For the It\^{o} stochastic 
integral the following estimate \cite{Kor} is valid
\begin{equation}
\label{pipi}
{\sf M}\left\{\left|\int\limits_t^T \xi_\tau dw_\tau\right|^q\right\}
\le K_q {\sf M}\left\{\left(\int\limits_t^T|\xi_\tau|^2 d\tau
\right)^{q/2}\right\},
\end{equation}
where $q>0$ is a fixed number, $w_\tau$ is a scalar
standard Wiener process, $\xi_\tau\in{\rm M}_2([t, T])$,
$K_q$ is a constant depending only on $q$,
$$
\int\limits_t^T|\xi_\tau|^2 d\tau<\infty\ \ \ {\rm w.~p.~1},
$$
$$
{\sf M}\left\{\left(\int\limits_t^T|\xi_\tau|^2 d\tau
\right)^{q/2}\right\}<\infty.
$$

Since
$$
\delta_{T,t}-\sum\limits_{j_1=0}^{p}
C_{j_1}\zeta_{j_1}^{(i)}=
\int\limits_t^T\biggl(\psi(s)-
\sum\limits_{j_1=0}^{p}
C_{j_1}\phi_{j_1}(s)\biggr)d{\bf w}_s^{(i)},
$$
then 
applying the 
estimate
(\ref{pipi}) to the right-hand side of this expression and considering that
$$
\int\limits_t^T\biggl(\psi(s)-
\sum\limits_{j_1=0}^{p}
C_{j_1}\phi_{j_1}(s)\biggr)^2ds\ \to 0
$$
if $p \to \infty$, we obtain
\begin{equation}
\label{pipi21}
\int\limits_t^T\psi(s)d{\bf w}_s^{(i)}=
\hbox{\vtop{\offinterlineskip\halign{
\hfil#\hfil\cr
$q$~-~{\rm l.i.m.}\cr
$\stackrel{}{{}_{p\to \infty}}$\cr
}} }\sum_{j_1=0}^{p}
C_{j_1}\zeta_{j_1}^{(i)},\ \ \ q>0.
\end{equation}

\noindent
Here $\hbox{\vtop{\offinterlineskip\halign{
\hfil#\hfil\cr
$q$~-~{\rm l.i.m.}\cr
$\stackrel{}{{}_{p\to \infty}}$\cr
}} }$ is a limit in the mean of degree $q.$
Hence, if $q=4$, then it is easy to conclude that w.~p.~1
$$
\hbox{\vtop{\offinterlineskip\halign{
\hfil#\hfil\cr
{\rm l.i.m.}\cr
$\stackrel{}{{}_{p\to \infty}}$\cr
}} }
\left(\sum\limits_{j_1=0}^{p}
C_{j_1}\zeta_{j_1}^{(i)}\right)^2=\delta^2_{T,t}.
$$

This equality as well as 
Parseval's equality 
were used in the last step of the formula (\ref{pipi20}).

When $k=3$ and $p_1=p_2=p_3=p$ we obtain (see (\ref{a3}))
\cite{2}-\cite{12aa}, \cite{arxiv-1}
$$
J[\psi^{(3)}]_{T,t}
=
\hbox{\vtop{\offinterlineskip\halign{
\hfil#\hfil\cr
{\rm l.i.m.}\cr
$\stackrel{}{{}_{p\to \infty}}$\cr
}} }\left(
\sum_{j_1,j_2,j_3=0}^{p}
C_{j_3j_2j_1}\zeta_{j_1}^{(i)}
\zeta_{j_2}^{(i)}
\zeta_{j_3}^{(i)}
-\right.
$$
$$
\left.-
\sum_{j_1,j_3=0}^{p}
C_{j_3j_1j_1}\zeta_{j_3}^{(i)}-
\sum_{j_1,j_2=0}^{p}
C_{j_2j_2j_1}\zeta_{j_1}^{(i)}-
\sum_{j_1,j_2=0}^{p}
C_{j_1j_2j_1}\zeta_{j_2}^{(i)}\right)=
$$
$$
=
\hbox{\vtop{\offinterlineskip\halign{
\hfil#\hfil\cr
{\rm l.i.m.}\cr
$\stackrel{}{{}_{p\to \infty}}$\cr
}} }\left(
\sum_{j_1,j_2,j_3=0}^{p}
C_{j_3j_2j_1}\zeta_{j_1}^{(i)}
\zeta_{j_2}^{(i)}
\zeta_{j_3}^{(i)}
-\sum_{j_1,j_3=0}^{p}\biggl(
C_{j_3j_1j_1}+
C_{j_1j_1j_3}+
C_{j_1j_3j_1}\biggr)\zeta_{j_3}^{(i)}\right)=
$$
$$
=
\hbox{\vtop{\offinterlineskip\halign{
\hfil#\hfil\cr
{\rm l.i.m.}\cr
$\stackrel{}{{}_{p\to \infty}}$\cr
}} }\left(
\sum_{j_1=0}^{p}
\sum_{j_2=0}^{j_1-1}
\sum_{j_3=0}^{j_2-1}
\biggl(C_{j_3j_2j_1}+
C_{j_3j_1j_2}+C_{j_2j_1j_3}+
C_{j_2j_3j_1}+
C_{j_1j_2j_3}+
C_{j_1j_3j_2}\biggr)\times\right.
$$
$$
\times
\zeta_{j_1}^{(i)}
\zeta_{j_2}^{(i)}
\zeta_{j_3}^{(i)}+
$$
$$
+\sum_{j_1=0}^{p}\sum_{j_3=0}^{j_1-1}\biggl(
C_{j_3j_1j_3}+
C_{j_1j_3j_3}+
C_{j_3j_3j_1}\biggr)\left(\zeta_{j_3}^{(i)}\right)^2
\zeta_{j_1}^{(i)}+
$$
$$
+
\sum_{j_1=0}^{p}\sum_{j_3=0}^{j_1-1}\biggl(
C_{j_3j_1j_1}+
C_{j_1j_1j_3}+
C_{j_1j_3j_1}\biggr)\left(\zeta_{j_1}^{(i)}\right)^2
\zeta_{j_3}^{(i)}+
$$
$$
\left.+
\sum_{j_1=0}^p C_{j_1j_1j_1}
\left(\zeta_{j_1}^{(i)}\right)^3-
\sum_{j_1,j_3=0}^{p}\left(
C_{j_3j_1j_1}+
C_{j_1j_1j_3}+
C_{j_1j_3j_1}\right)\zeta_{j_3}^{(i)}\right)=
$$

$$
=
\hbox{\vtop{\offinterlineskip\halign{
\hfil#\hfil\cr
{\rm l.i.m.}\cr
$\stackrel{}{{}_{p\to \infty}}$\cr
}} }\left(
\sum_{j_1=0}^{p}
\sum_{j_2=0}^{j_1-1}
\sum_{j_3=0}^{j_2-1}
C_{j_1}
C_{j_2}C_{j_3}
\zeta_{j_1}^{(i)}
\zeta_{j_2}^{(i)}
\zeta_{j_3}^{(i)}+\right.
$$
$$
+
\frac{1}{2}\sum_{j_1=0}^{p}\sum_{j_3=0}^{j_1-1}
C_{j_3}^2C_{j_1}\left(\zeta_{j_3}^{(i)}\right)^2
\zeta_{j_1}^{(i)}
+\frac{1}{2}\sum_{j_1=0}^{p}\sum_{j_3=0}^{j_1-1}
C_{j_1}^2C_{j_3}\left(\zeta_{j_1}^{(i)}\right)^2
\zeta_{j_3}^{(i)}+
$$
$$
\left.+\frac{1}{6}\sum_{j_1=0}^p C_{j_1}^3
\left(\zeta_{j_1}^{(i)}\right)^3-
\frac{1}{2}\sum_{j_1,j_3=0}^{p}
C_{j_1}^2C_{j_3}\zeta_{j_3}^{(i)}\right)=
$$

$$
=
\hbox{\vtop{\offinterlineskip\halign{
\hfil#\hfil\cr
{\rm l.i.m.}\cr
$\stackrel{}{{}_{p\to \infty}}$\cr
}} }\Biggl(
\frac{1}{6}\sum_{\stackrel{j_1,j_2,j_3=0}{{}_{j_1\ne j_2, j_2\ne j_3,
j_1\ne j_3}}}^{p}
C_{j_1}
C_{j_2}C_{j_3}
\zeta_{j_1}^{(i)}
\zeta_{j_2}^{(i)}
\zeta_{j_3}^{(i)}+\Biggr.
$$
$$
+\frac{1}{2}\sum_{j_1=0}^{p}\sum_{j_3=0}^{j_1-1}
C_{j_3}^2C_{j_1}\left(\zeta_{j_3}^{(i)}\right)^2
\zeta_{j_1}^{(i)}
+\frac{1}{2}\sum_{j_1=0}^{p}\sum_{j_3=0}^{j_1-1}
C_{j_1}^2C_{j_3}\left(\zeta_{j_1}^{(i)}\right)^2
\zeta_{j_3}^{(i)}+
$$
$$
\Biggl.+\frac{1}{6}\sum_{j_1=0}^p C_{j_1}^3
\left(\zeta_{j_1}^{(i)}\right)^3-
\frac{1}{2}\sum_{j_1,j_3=0}^{p}
C_{j_1}^2C_{j_3}\zeta_{j_3}^{(i)}\Biggr)
=
$$
$$
=
\hbox{\vtop{\offinterlineskip\halign{
\hfil#\hfil\cr
{\rm l.i.m.}\cr
$\stackrel{}{{}_{p\to \infty}}$\cr
}} }\left(
\frac{1}{6}\sum_{j_1,j_2,j_3=0}^{p}
C_{j_1}C_{j_2}C_{j_3}\zeta_{j_1}^{(i)}
\zeta_{j_2}^{(i)}\zeta_{j_3}^{(i)}-\right.
$$
$$
-\frac{1}{6}\left(
3\sum_{j_1=0}^{p}\sum_{j_3=0}^{j_1-1}
C_{j_3}^2C_{j_1}\left(\zeta_{j_3}^{(i)}\right)^2
\zeta_{j_1}^{(i)}
+3\sum_{j_1=0}^{p}\sum_{j_3=0}^{j_1-1}
C_{j_1}^2C_{j_3}\left(\zeta_{j_1}^{(i)}\right)^2
\zeta_{j_3}^{(i)}+\right.
$$
$$
\left.+
\sum_{j_1=0}^p C_{j_1}^3
\left(\zeta_{j_1}^{(i)}\right)^3\right)+
$$
$$
+\frac{1}{2}\sum_{j_1=0}^{p}\sum_{j_3=0}^{j_1-1}
C_{j_3}^2C_{j_1}\left(\zeta_{j_3}^{(i)}\right)^2
\zeta_{j_1}^{(i)}
+\frac{1}{2}\sum_{j_1=0}^{p}\sum_{j_3=0}^{j_1-1}
C_{j_1}^2C_{j_3}\left(\zeta_{j_1}^{(i)}\right)^2
\zeta_{j_3}^{(i)}+
$$
$$
\left.+\frac{1}{6}\sum_{j_1=0}^p C_{j_1}^3
\left(\zeta_{j_1}^{(i)}\right)^3-
\frac{1}{2}\sum_{j_1,j_3=0}^{p}
C_{j_1}^2C_{j_3}\zeta_{j_3}^{(i)}\right)=
$$
$$
=
\hbox{\vtop{\offinterlineskip\halign{
\hfil#\hfil\cr
{\rm l.i.m.}\cr
$\stackrel{}{{}_{p\to \infty}}$\cr
}} }\left(
\frac{1}{6}\left(\sum_{j_1=0}^{p}
C_{j_1}
\zeta_{j_1}^{(i)}\right)^3-
\frac{1}{2}\sum_{j_1=0}^{p}
C_{j_1}^2 \sum_{j_3=0}^p C_{j_3}\zeta_{j_3}^{(i)}\right)
=
$$

\vspace{-1mm}
\begin{equation}
\label{pipi22}
=\frac{1}{3!}\left(\delta_{T,t}^3-3\delta_{T,t}\Delta_{T,t}\right)\ \ \ \hbox{w.~p.~1}.
\end{equation}

\vspace{2mm}

The last step in (\ref{pipi22}) follows from Parseval's equality,
Theorem 1.1 for $k=1$, and the equality
$$
\hbox{\vtop{\offinterlineskip\halign{
\hfil#\hfil\cr
{\rm l.i.m.}\cr
$\stackrel{}{{}_{p\to \infty}}$\cr
}} }
\left(\sum\limits_{j_1=0}^{p}
C_{j_1}\zeta_{j_1}^{(i)}\right)^3=\delta^3_{T,t}\ \ \ \hbox{w.~p.~1},
$$
which can be obtained easily when
$q=8$ (see (\ref{pipi21})).

In addition, we used the following relations between Fourier 
coefficients for the considered case
\begin{equation}
\label{nahod0}
C_{j_1j_2}+C_{j_2j_1}=C_{j_1}C_{j_2},\ \ \
2C_{j_1j_1}=C_{j_1}^2,
\end{equation}
\vspace{-8mm}
\begin{equation}
\label{nahod}
~~~~~~C_{j_1j_2j_3}+
C_{j_1j_3j_2}+
C_{j_2j_3j_1}+
C_{j_2j_1j_3}+
C_{j_3j_2j_1}+
C_{j_3j_1j_2}=C_{j_1}C_{j_2}C_{j_3},
\end{equation}
\begin{equation}
\label{nahod1}
2\left(C_{j_1j_1j_3}+
C_{j_1j_3j_1}+
C_{j_3j_1j_1}\right)=C_{j_1}^2C_{j_3},
\end{equation}
\begin{equation}
\label{nahod2}
6C_{j_1j_1j_1}=C_{j_1}^3.
\end{equation}

\subsection{On Usage of Discontinuous Complete Orthonormal Systems 
of Functions in Theorem 1.1}

\footnotetext[4]{The results of this section will be generalized to the case of an arbitrary 
complete ortho\-nor\-mal system of functions $\{\phi_j(x)\}_{j=0}^{\infty}$ 
in the space $L_2([t, T])$
and $\psi_1(\tau),$ $\ldots,\psi_k(\tau) \in L_2([t, T])$ in Sect.~1.11
(see Theorem~1.16).}

Analyzing the proof of Theorem 1.1, we can ask the question: can we 
weaken the  continuity condition for the functions 
$\phi_j(x),$ $j=1, 2,\ldots$?${}^{4}$

{\it We will say that the function $f(x):$ $[t, T]\to {\bf R}$ 
satisfies the condition {\rm (}$\star ${\rm )}, if it is 
continuous at the interval $[t, T]$ except may be
for the finite number of points 
of the finite discontinuity as well as it is 
right-continuous 
at the interval  $[t, T].$}

Furthermore, let us suppose that $\{\phi_j(x)\}_{j=0}^{\infty}$
is a complete orthonormal system of functions in the 
space $L_2([t, T])$, each function $\phi_j(x)$ 
of which for $j<\infty$ satisfies the 
condition $(\star)$.

It is easy to see that continuity of the functions $\phi_j(x)$ was used 
substantially for the proof of Theorem 1.1 in two places. More precisely, 
we mean Lemma 1.3 and the 
formula (\ref{30.52}).
It is clear that without the loss of generality the partition 
$\{\tau_j\}_{j=0}^N$ of the interval $[t, T]$ in Lemma 1.3 and 
(\ref{30.52}) can be taken so ``dense"\  that 
among the 
points $\tau_j$ of this partition there will be all points 
of jumps of the functions
$\varphi_1(\tau)=\phi_{j_1}(\tau),$ $\ldots,$ 
$\varphi_k(\tau)=\phi_{j_k}(\tau)$
($j_1,\ldots,j_k<\infty$)
and among the points $(\tau_{j_1},\ldots,\tau_{j_k})$ for which
$0\le j_1<\ldots<j_k\le N-1$
there will be all points of jumps
of the function
$\Phi(t_1,\ldots,t_k)$.

Let us demonstrate how to modify the proofs of Lemma 1.3 and  
the formula (\ref{30.52}) 
in the case when 
$\{\phi_j(x)\}_{j=0}^{\infty}$
is a complete orthonormal system of functions in the 
space $L_2([t, T])$, each function $\phi_j(x)$ 
of which for $j<\infty$ satisfies the 
condition $(\star)$.

At first, appeal to Lemma 1.3. From the proof of this lemma 
it follows that 
$$
{\sf M}\left\{\left|\sum_{j=0}^{N-1}J[\Delta\varphi_l]_{\tau_{j+1},
\tau_j}\right|^4
\right\}=
\sum_{j=0}^{N-1}{\sf M}\left\{\biggl|J[\Delta\varphi_l]_{\tau_{j+1},
\tau_j}\biggr|^4
\right\}+
$$

\vspace{-2mm}
\begin{equation}
\label{4444.02}
+ 6 \sum_{j=0}^{N-1}{\sf M}
\left\{\biggl|J[\Delta\varphi_l]_{\tau_{j+1},\tau_j}\biggr|^2
\right\}
\sum_{q=0}^{j-1}{\sf M}\left\{\biggl|
J[\Delta\varphi_l]_{\tau_{q+1},\tau_q}\biggr|^2
\right\},
\end{equation}

$$
{\sf M}\left\{\left|J[\Delta\varphi_l]_{\tau_{j+1},\tau_j}\right|^2\right\}=
\int\limits_{\tau_j}^{\tau_{j+1}}(\varphi_l(\tau_j)-\varphi_l(s))^2ds,\
$$
$$
{\sf M}\left\{\left|J[\Delta\varphi_l]_{\tau_{j+1},\tau_j}\right|^4\right\}=
3\left(\int\limits_{\tau_j}^{\tau_{j+1}}(\varphi_l(\tau_j)-\varphi_l(s))^2ds
\right)^2.
$$

\vspace{2mm}

Suppose that the functions $\varphi_l(s)$ ($l=1,\ldots,k$)  satisfy 
the condition $(\star)$ and the partition $\{\tau_j\}_{j=0}^{N}$ 
includes all points of 
jumps of the functions
$\varphi_l(s)$ ($l=1,\ldots,k$).  It means that for the integral 
$$
\int\limits_{\tau_j}^{\tau_{j+1}}(\varphi_l(\tau_j)-\varphi_l(s))^2ds
$$
the integrand function is 
continuous 
at the interval
$[\tau_j, \tau_{j+1}],$ except possibly the point $\tau_{j+1}$ of
finite 
discontinuity.

Let $\mu\in (0, \Delta\tau_j)$ be fixed. 
Due to continuity (which means uniform continuity)
of the functions
$\varphi_l(s)$ ($l=1,\ldots,k$) 
at the interval $[\tau_j, \tau_{j+1}-\mu]$ we have
$$
\int\limits_{\tau_j}^{\tau_{j+1}}(\varphi_l(\tau_j)-\varphi_l(s))^2ds=
$$
$$
=
\int\limits_{\tau_j}^{\tau_{j+1}-\mu}(\varphi_l(\tau_j)-\varphi_l(s))^2ds 
+ \int\limits_{\tau_{j+1}-\mu}^{\tau_{j+1}}(\varphi_l(\tau_j)-\varphi_l(s))^2ds
<
$$

\begin{equation}
\label{4444.90}
<\varepsilon^2 (\Delta\tau_j-\mu)+M^2\mu.
\end{equation}

\vspace{2mm}

When obtaining the inequality (\ref{4444.90}) we supposed that
$\Delta\tau_j<\delta(\varepsilon)$ for all $j=0, 1,\ldots,N-1$ 
(here $\delta(\varepsilon)>0$ exists for any $\varepsilon>0$ 
and it does not depend 
on $s$), 
$$
|\varphi_l(\tau_j)-\varphi_l(s)|<\varepsilon
$$

\vspace{2mm}
\noindent
for $s\in [\tau_{j}, \tau_{j+1}-\mu]$
(due to uniform continuity of the 
functions $\varphi_l(s),$ $l=1,\ldots,k$),
$$
|\varphi_l(\tau_j)-\varphi_l(s)|<M
$$

\vspace{2mm}
\noindent
for
$s\in [\tau_{j+1}-\mu, \tau_{j+1}]$, $M$ is a constant
(potential discontinuity point
of the function $\varphi_l(s)$ 
is the 
point $\tau_{j+1}$).

Performing the 
passage to the limit
in the inequality (\ref{4444.90}) when $\mu\to +0$, we get 
\begin{equation}
\label{tra1}
\int\limits_{\tau_j}^{\tau_{j+1}}(\varphi_l(\tau_j)-\varphi_l(s))^2ds\le
\varepsilon^2 \Delta\tau_j.
\end{equation}

\vspace{2mm}

Using (\ref{tra1}) to estimate the right-hand side of 
(\ref{4444.02}),
we obtain
$$
{\sf M}\left\{\left|\sum_{j=0}^{N-1}J[\Delta\varphi_l]_{\tau_{j+1},
\tau_j}\right|^4
\right\}\le 
\varepsilon^4\left(
3 \sum_{j=0}^{N-1}(\Delta\tau_{j})^2+
6 \sum_{j=0}^{N-1}\Delta\tau_{j}
\sum_{q=0}^{j-1}\Delta\tau_{q}\right)<
$$

\vspace{-1mm}
\begin{equation}
\label{4444.92}
<
3\varepsilon^4\left(\delta(\varepsilon) (T-t)+(T-t)^2\right).
\end{equation}

\vspace{2mm}

This implies that
$$
{\sf M}\left\{\left|\sum\limits_{j=0}^{N-1}J[\Delta\varphi_l]_{\tau_{j+1},
\tau_j}\right|^4\right\} \to 0
$$

\vspace{2mm}
\noindent
when $N\to\infty$ and Lemma 1.3 remains correct.

Now, let us present explanations concerning the correctness of 
(\ref{30.52}), when $\{\phi_j(x)\}_{j=0}^{\infty}$ is a 
complete
orthonormal system of functions in the space $L_2([t, T])$, 
each function $\phi_j(x)$ of which 
for $j<\infty$ satisfies the condition $(\star)$.

Consider the case $k=3$ and the representation (\ref{4444.1}). 
Let us demonstrate that in the studied case the first limit 
on the right-hand side of (\ref{4444.1}) equals to zero 
(similarly, we can demonstrate that the second limit on the right-hand side 
of (\ref{4444.1}) equals to zero; 
proof of the second limit 
equality to zero on the right-hand side
of the formula
(\ref{44444.25}) is the same as for the case of continuous functions
$\phi_j(x),$ $j=0, 1,\ldots).$

The second moment of the prelimit expression of first limit on the 
right-hand side of (\ref{4444.1}) looks as follows

\vspace{-1mm}
$$
\sum_{j_3=0}^{N-1}\sum_{j_2=0}^{j_3-1}\sum_{j_1=0}^{j_2-1}
\int\limits_{\tau_{j_2}}^{\tau_{j_2+1}}\int\limits_{\tau_{j_1}}^{\tau_{j_1+1}}
\left(\Phi(t_1,t_2,\tau_{j_3})-\Phi(t_1,\tau_{j_2},\tau_{j_3})\right)^2
dt_1 dt_2
\Delta\tau_{j_3}.
$$

\vspace{2mm}
       
Further, for the fixed $\mu\in(0, \Delta\tau_{j_2})$ 
and $\rho\in(0, \Delta\tau_{j_1})$ we have

\vspace{-1mm}
$$
\int\limits_{\tau_{j_2}}^{\tau_{j_2+1}}\int\limits_{\tau_{j_1}}^{\tau_{j_1+1}}
\left(\Phi(t_1,t_2,\tau_{j_3})-\Phi(t_1,\tau_{j_2},\tau_{j_3})\right)^2
dt_1 dt_2=
$$
$$
=\left(\int\limits_{\tau_{j_2}}^{\tau_{j_2+1}-\mu}
+\int\limits_{\tau_{j_2+1}-\mu}^{\tau_{j_2+1}}\right)
\left(
\int\limits_{\tau_{j_1}}^{\tau_{j_1+1}-\rho}+
\int\limits_{\tau_{j_1+1}-\rho}^{\tau_{j_1+1}}\right)
\left(\Phi(t_1,t_2,\tau_{j_3})-\Phi(t_1,\tau_{j_2},\tau_{j_3})\right)^2
dt_1 dt_2=
$$

\vspace{-1mm}
$$
=\left(\int\limits_{\tau_{j_2}}^{\tau_{j_2+1}-\mu}
\int\limits_{\tau_{j_1}}^{\tau_{j_1+1}-\rho}+
\int\limits_{\tau_{j_2}}^{\tau_{j_2+1}-\mu}
\int\limits_{\tau_{j_1+1}-\rho}^{\tau_{j_1+1}}+
\int\limits_{\tau_{j_2+1}-\mu}^{\tau_{j_2+1}}
\int\limits_{\tau_{j_1}}^{\tau_{j_1+1}-\rho}+
\int\limits_{\tau_{j_2+1}-\mu}^{\tau_{j_2+1}}
\int\limits_{\tau_{j_1+1}-\rho}^{\tau_{j_1+1}}\right)\times
$$

\vspace{2mm}
$$
\times
\left(\Phi(t_1,t_2,\tau_{j_3})-\Phi(t_1,\tau_{j_2},\tau_{j_3})\right)^2
dt_1 dt_2<
$$

\vspace{-1mm}
$$
<\varepsilon^2\left(\Delta\tau_{j_2}-\mu\right)
\left(\Delta\tau_{j_1}-\rho\right)
+M^2\rho\left(\Delta\tau_{j_2}-\mu\right)+
$$

\vspace{-3mm}
\begin{equation}
\label{4444.54}
+M^2\mu\left(\Delta\tau_{j_1}-\rho\right)
+M^2\mu\rho,
\end{equation}

\vspace{4mm}
\noindent
where $M$ is a constant,
$\Delta\tau_j<\delta(\varepsilon)$ for $j=0, 1,\ldots,N-1$
($\delta(\varepsilon)>0$ exists for any $\varepsilon>0$ and it does not  
depend on points
$(t_1,t_2,\tau_{j_3}),$ $(t_1,\tau_{j_2},\tau_{j_3})$).
We suppose here that the partition 
$\{\tau_j\}_{j=0}^{N}$ contains all discontinuity points
of the function
$\Phi(t_1,t_2,t_3)$ 
as points $\tau_j$ (for each variable with fixed remaining two variables).
When obtaining
the inequality (\ref{4444.54})  
we also supposed that potential discontinuity points
of this function (for each variable with fixed remaining two variables) 
are contained among the points
$\tau_{j_1+1}, \tau_{j_2+1}, \tau_{j_3+1}$.

Let us explain in detail how we obtained the inequality 
(\ref{4444.54}). 
Since the function $\Phi(t_1,t_2,t_3)$ is continuous at the closed 
bounded set  
$$
Q_3=\biggl\{(t_1, t_2, t_3): t_1\in[\tau_{j_1}, \tau_{j_1+1}-\rho],
t_2\in[\tau_{j_2}, \tau_{j_2+1}-\mu],
t_3\in [\tau_{j_3}, \tau_{j_3+1}-\nu]\biggr\},
$$

\noindent
where $\rho, \mu, \nu$
are fixed small positive numbers such that
$$
\nu\in(0, \Delta\tau_{j_3}),\ \ \ 
\mu\in(0, \Delta\tau_{j_2}),\ \ \ 
\rho\in(0, \Delta\tau_{j_1}),
$$

\noindent
then this function is also uniformly continous 
at this set. Moreover, the function $\Phi(t_1,t_2,t_3)$ 
is supposed to be bounded at the closed set 
$D_3$ (see the proof of Theorem 1.1).

Since the distance between points
$(t_1,t_2,\tau_{j_3})$,
$(t_1,\tau_{j_2},\tau_{j_3})$ $\in $ $Q_3$ 
is obviously less than
$\delta(\varepsilon)$ ($\Delta\tau_j<\delta(\varepsilon)$ for
$j=0, 1,\ldots,N-1$),
then 

\vspace{-1mm}
$$
|\Phi(t_1,t_2,\tau_{j_3})-\Phi(t_1,\tau_{j_2},\tau_{j_3})|
<\varepsilon.
$$

\vspace{2mm}
 
This inequality was used 
to estimate the first double integral in (\ref{4444.54}). 
Estimating 
the three remaining double integrals
in (\ref{4444.54}) we used the boundedness property for
the function 
$\Phi(t_1,t_2,t_3)$ in the form of inequality  

\vspace{-1mm}
$$
|\Phi(t_1,t_2,\tau_{j_3})-\Phi(t_1,\tau_{j_2},\tau_{j_3})|
<M.
$$

\vspace{2mm}

Performing the 
passage to the limit 
in the inequality (\ref{4444.54})
when $\mu, \rho\to +0,$ we obtain the estimate

\vspace{-1mm}
$$
\int\limits_{\tau_{j_2}}^{\tau_{j_2+1}}\int\limits_{\tau_{j_1}}^{\tau_{j_1+1}}
\left(\Phi(t_1,t_2,\tau_{j_3})-\Phi(t_1,\tau_{j_2},\tau_{j_3})\right)^2
dt_1 dt_2\le
\varepsilon^2\Delta\tau_{j_2}
\Delta\tau_{j_1}.
$$

\vspace{2mm}

This estimate provides  
$$
\sum_{j_3=0}^{N-1}\sum_{j_2=0}^{j_3-1}\sum_{j_1=0}^{j_2-1}
\int\limits_{\tau_{j_2}}^{\tau_{j_2+1}}\int\limits_{\tau_{j_1}}^{\tau_{j_1+1}}
\left(\Phi(t_1,t_2,\tau_{j_3})-\Phi(t_1,\tau_{j_2},\tau_{j_3})\right)^2
dt_1 dt_2
\Delta\tau_{j_3}\le
$$

\vspace{-1mm}
$$
\le\varepsilon^2
\sum_{j_3=0}^{N-1}\sum_{j_2=0}^{j_3-1}\sum_{j_1=0}^{j_2-1}
\Delta\tau_{j_1}\Delta\tau_{j_2}\Delta\tau_{j_3}<
\varepsilon^2 \frac{(T-t)^3}{6}.
$$

\vspace{2mm}

The last inequality means that in the considered case the first 
limit on the right-hand side of (\ref{4444.1}) equals to zero 
(similarly, we can
demonstrate that the second limit on the right-hand side of (\ref{4444.1})
equals to zero).

Consequently, the formula (\ref{30.52}) is correct when 
$k=3$ in the studied case. 
Similarly, we can perform the argumentation for the cases
$k=2$ and $k>3.$

Therefore, in Theorem 1.1 we can use complete orthonormal systems of 
functions $\{\phi_j(x)\}_{j=0}^{\infty}$ in the space 
$L_2([t, T]),$ each function
$\phi_j(x)$ of which for $j<\infty$ satisfies the condition
$(\star)$.

One of the examples of such systems of functions is a complete orthonormal 
system of Haar functions in the space 
$L_2([t, T])$ 
$$
\phi_0(x)=\frac{1}{\sqrt{T-t}},\ \ \ 
\phi_{nj}(x)=\frac{1}{\sqrt{T-t}}\ \varphi_{nj}\biggl(\frac{x-t}{T-t}\biggr),
$$
where
$n=0, 1,\ldots,$\ \  $j=1, 2,\ldots, 2^n,$
and the functions $\varphi_{nj}(x)$ are defined as 
$$
\varphi_{nj}(x)=
\left\{
\begin{matrix}
2^{n/2},\ &x\in[(j-1)/2^n,\ (j-1)/2^n+
1/2^{n+1})\cr\cr
-2^{n/2},\ &x\in[(j-1)/2^n+1/2^{n+1},\
j/2^n)\cr\cr
0,\  &\hbox{otherwise}
\end{matrix}\right.,
$$

\noindent
$n=0, 1,\ldots,$\  $j=1, 2,\ldots, 2^n$ 
(we choose the values of Haar functions 
in the points of discontinuity in such a way that these functions will be 
right-continuous).

The other example of similar system of functions is a complete orthonormal 
system of Rademacher--Walsh functions in the space  $L_2([t, T])$
$$
\phi_0(x)=\frac{1}{\sqrt{T-t}},
$$
$$
\phi_{m_1\ldots m_k}(x)=
\frac{1}{\sqrt{T-t}}\ \varphi_{m_1}\biggl(\frac{x-t}{T-t}\biggr)
\ldots \varphi_{m_k}
\biggl(\frac{x-t}{T-t}\biggr),
$$

\vspace{1mm}
\noindent
where $0<m_1<\ldots<m_k,$\ \  $m_1,\ldots,m_k=1, 2,\ldots,$\ \  $k=1, 2,\ldots,$
$$
\varphi_m(x)=(-1)^{[2^m x]},
$$
$x\in [0, 1]$,\ \ $ m=1, 2,\ldots,$\ \ $[y]$ is an integer part of 
a real number $y.$

\subsection{Remark on Usage of Complete Orthonormal Systems 
of Functions in Theorem 1.1}

Note that actually the 
functions $\phi_j(s)$ from the complete orthonormal system 
of functions $\{\phi_j(s)\}_{j=0}^{\infty}$ in the space 
$L_2([t, T])$ depend not only on $s$, but 
on $t$ and $T.$
   
For example, the complete orthonormal systems of Legendre polynomials 
and trigonometric functions in the space $L_2([t, T])$ have the 
following form
$$
\phi_j(s,t,T)=\sqrt{\frac{2j+1}{T-t}}P_j\left(\left(
s-\frac{T+t}{2}\right)\frac{2}{T-t}\right),
$$
$$
P_j(y)=\frac{1}{2^j j!} \frac{d^j}{dy^j}\left(y^2-1\right)^j,
$$
where
$P_j(y)$ $(j=0, 1, 2,\ldots)$ is the Legendre polynomial,
\begin{equation}
\label{trig11}
~~~\phi_j(s,t,T)=\frac{1}{\sqrt{T-t}}
\left\{
\begin{matrix}
1,\ & j=0\cr\cr
\sqrt{2}{\rm sin} \left(2\pi r(s-t)/(T-t)\right),\ & j=2r-1\cr\cr
\sqrt{2}{\rm cos} \left(2\pi r(s-t)/(T-t)\right),\ & j=2r
\end{matrix}
,\right.
\end{equation}
where $r=1, 2,\ldots $

Note that the specified systems of functions are assumed to be used in
the context of implementation of numerical methods for It\^{o} SDEs
(see Chapter 4)
for the sequences of time intervals 
$$
[T_0, T_1],\
[T_1, T_2],\ [T_2, T_3],\ \ldots\ 
$$
and Hilbert spaces 
$$
L_2([T_0, T_1]),\ L_2([T_1, T_2]),\
L_2([T_2, T_3]),\ \ldots
$$

We can explain that the dependence of functions 
$\phi_j(s,t,T)$ on $t$ and $T$ (hereinafter these constants will 
mean fixed moments of time) will not affect on the main properties
of independence of random variables
$$
\zeta_{(j)T,t}^{(i)}=
\int\limits_t^T \phi_{j}(s,t,T) d{\bf w}_s^{(i)},
$$
where $i=1,\ldots,m$ and $j=0, 1, 2,\ldots $

Indeed, for fixed  $t$ and $T$ due to orthonormality of the mentioned systems 
of functions we have

\vspace{-2mm}
$$
{\sf M}\left\{\zeta_{(j)T,t}^{(i)}\zeta_{(g)T,t}^{(r)}\right\}=
{\bf 1}_{\{i=r\}}
{\bf 1}_{\{j=g\}},
$$

\vspace{3mm}
\noindent
where 
$i, r=1,\ldots,m,$\ \ $j, g=0, 1, 2,\ldots $

This means that $\zeta_{(j)T,t}^{(i)}$ and $\zeta_{(g)T,t}^{(r)}$
are independent for $j\ne g$ or $i\ne r$
(since these random variables
are Gaussian). 

From the other side, 
the random
variables   

\vspace{-2mm}
$$
\zeta_{(j)T_1,t_1}^{(i)}=
\int\limits_{t_1}^{T_1}\phi_{j}(s,t_1,T_1) d{\bf w}_s^{(i)},\ \ \ \
\zeta_{(j)T_2,t_2}^{(i)}=
\int\limits_{t_2}^{T_2}\phi_{j}(s,t_2,T_2) d{\bf w}_s^{(i)}
$$

\vspace{3mm}
\noindent
are independent if
$[t_1, T_1]\cap [t_2, T_2]=\emptyset$
(the case $T_1=t_2$ is possible)
according to 
the properties
of the It\^{o} stochastic integral.

Therefore, the important properties of random variables
$\zeta_{(j)T,t}^{(i)}$, which are the basic motive of their 
usage, are saved.

\subsection{Convergence in the 
Mean of Degree $2n$ ($n\in {\bf N}$) of Expansions of Iterated
It\^{o} Stochastic Integrals from Theorem 1.1}

Constructing the expansions of iterated
It\^{o} stochastic integrals from Theorem 1.1 we 
saved all 
information about these integrals. That is why it is 
natural to expect that the mentioned expansions will converge
not only in the mean-square sense but in the stronger probabilistic senses.

We will obtain the general estimate which prove convergence in 
the mean
of degree $2n$ ($n\in {\bf N}$) of expansion from Theorem 1.1.

According to the notations of Theorem 1.1 (see (\ref{s2s}), (\ref{y007})), we have

\vspace{-1mm}
$$
R_{T,t}^{p_1,\ldots, p_k}
=J[\psi^{(k)}]_{T,t}-
J[\psi^{(k)}]_{T,t}^{p_1,\ldots,p_k}=J'[R_{p_1\ldots p_k}]_{T,t}^{(i_1\ldots i_k)}=
$$

\vspace{-1mm}
\begin{equation}
\label{jjjye}
=\sum_{(t_1,\ldots,t_k)}
\int\limits_{t}^{T}
\ldots
\int\limits_{t}^{t_2}
R_{p_1\ldots p_k}(t_1,\ldots,t_k)
d{\bf w}_{t_1}^{(i_1)}
\ldots
d{\bf w}_{t_k}^{(i_k)},
\end{equation}

\noindent
where

\newpage
\noindent
\begin{equation}
\label{ch12026ll12}
~~~~~~R_{p_1\ldots p_k}(t_1,\ldots,t_k)\stackrel{{\rm def}}{=}
K(t_1,\ldots,t_k)-
\sum_{j_1=0}^{p_1}\ldots
\sum_{j_k=0}^{p_k}
C_{j_k\ldots j_1}
\prod_{l=1}^k\phi_{j_l}(t_l),
\end{equation}

\vspace{2mm}
\noindent
$J[\psi^{(k)}]_{T,t}$ is the stochastic integral {\rm (\ref{ito}),}
$J[\psi^{(k)}]_{T,t}^{p_1,\ldots,p_k}$ is the 
expression on the right-hand side of {\rm (\ref{tyyy})} before
passing to the limit 
$$
\hbox{\vtop{\offinterlineskip\halign{
\hfil#\hfil\cr
{\rm l.i.m.}\cr
$\stackrel{}{{}_{p_1,\ldots,p_k\to \infty}}$\cr
}} }.
$$ 

Note that for definiteness we consider 
the case $i_1,\ldots,i_k=1,\ldots,m$ in this section.
Another 
notations from this section are the same as in the
formulation and proof of Theorem 1.1.

When proving Theorem 1.1 we 
obtained the following estimate (see (\ref{obana1})) 

$$
{\sf M}\left\{\left(J'[R_{p_1\ldots p_k}]_{T,t}^{(i_1\ldots i_k)}\right)^2\right\}\le C_k 
\int\limits_{[t, T]^k}
R_{p_1\ldots p_k}^2(t_1,\ldots,t_k)
dt_1
\ldots
dt_k,
$$

\vspace{1mm}
\noindent
where $C_k$ is a constant.
Obviously, $C_k=k!$ for the case $i_1,\ldots,i_k=1,\ldots,m$.

First, we note that the iterated It\^{o} stochastic integral (\ref{ito})
can be considered as a multiple Wiener stochastic integral 
with respect to the components of a multidimensional
Wiener process. The multiple Wiener stochastic integral
with respect to a scalar Wiener process was first considered
in \cite{ito1951} (1951). 
Multiple Wiener stochastic integrals, including integrals
with respect to the components of a 
multidimensional Wiener process, are discussed in detail
in Sect.~1.10, 1.11, 1.14 (also see Sect.~1.9).
In fact, we have already considered the 
multiple Wiener stochastic integral 
(see (\ref{mult11})).

Let $H(s): [t,T]\to {\bf R}$.
Let $(\Omega,{\rm F},{\sf P})$ is a probability space, where
${\rm F}$ is the smallest $\sigma$-algebra such that
the random variables
$$
\int\limits_t^T H(s)dw_s
$$
are ${\rm F}$-measurable for every $H(s)\in L_2([t,T]),$
where 
$w_s$ is a standard Wiener process,
$$
\int\limits_t^T H(s)dw_s
$$
is a usual Wiener--It\^{o} stochastic integral.

It is well known (see Theorem~9.7.1, Theorem~9.7.3 and Theorem~9.6.7 in \cite{Kuo}) that
$$
L_2(\Omega,{\rm F},{\sf P})=\bigoplus\limits_{k=0}^{\infty}{\cal H}_k,
$$

\vspace{2mm}
\noindent
where ${\cal H}_0$ contains only constants,
the space ${\cal H}_k$ is the so-called $k$th $(k\ge 1)$ homogeneous Wiener chaos
which consists of all random variables of the form

\vspace{-3mm}
\begin{equation}
\label{ch12026oooo1}
k!\int\limits_t^{T}\ldots \int\limits_t^{t_2}
\hat H(t_1,\ldots,t_k)dw_{t_1}\ldots dw_{t_k},
\end{equation}

\noindent
where (\ref{ch12026oooo1}) is a representation for the 
multiple Wiener stochastic integral with respect to 
the scalar Wiener process (see (\ref{mult11}) or (\ref{WiI}) for the case 
$i_1=\ldots=i_k=i\in\{1,\ldots,m\}$ and with $w_s$ instead of ${\bf w}^{(i)}_s$),
$\hat H(t_1,\ldots,t_k)$ is a symmetrization
of $H(t_1,\ldots,t_k),$ $H(t_1,\ldots,t_k)\in L_2(D_k),$
$D_k=\{(t_1,\ldots,t_k):  t<t_1<\ldots <t_k<T\}.$

Note that (see \cite{Nurd2}, Corollary~2.8.4)

\vspace{-3.5mm}
\begin{equation}
\label{ch12026777}
{\sf M}\left\{\xi^{2n}\right\}\le (2n-1)^{kn}\left({\sf M}\left\{\xi^2\right\}\right)^n,
\end{equation}

\vspace{1mm}
\noindent
where $\xi\in {\cal H}_k,$ $n, k\in {\bf N}.$

The following estimate 
for the multiple Wiener stochastic integral
with respect to a scalar Wiener process 
is a consequence of the inequality (\ref{ch12026777})
\cite{Nual1}, \cite{Nurd1}
\begin{equation}
\label{2026ch1001}
~~~~~~~~{\sf M}\left\{\left(
J'[\Phi]_{T,t}^{(i_1\ldots i_k)}\right)^{2n}\right\}
\le (2n-1)^{nk} \left({\sf M}\left\{\left(
J'[\Phi]_{T,t}^{(i_1\ldots i_k)}\right)^{2}\right\}\right)^{n},
\end{equation}

\vspace{1mm}
\noindent
where $n\in{\bf N},$ 
$J'[\Phi]_{T,t}^{(i_1\ldots i_k)}$ is a multiple
Wiener stochastic integral defined as in Sect.~1.1.3 (see (\ref{mult11})) or as in Sect.~1.11
(see (\ref{WiI}))
but for the case of a scalar Wiener process $(i_1=\ldots=i_k=i\in\{1,\ldots,m\})$,
$\Phi(t_1,\ldots,t_k)\in L_2([t, T]^k)$ in (\ref{WiI}) and  
$\Phi(t_1,\ldots,t_k)\in C([t, T]^k)\subset L_2([t, T]^k)$ in (\ref{mult11}),
$J'[\Phi]_{T,t}^{(i_1\ldots i_k)}\in {\cal H}_k$ 
(see (\ref{pobeda}), (\ref{s2s}), (\ref{Wi110})). We also note that

\vspace{-2mm}
$$
J'[\Phi]_{T,t}^{(i_1\ldots i_k)}=J'[\hat\Phi]_{T,t}^{(i_1\ldots i_k)}\ \ \ \hbox{w.~p.~1,}
$$

\vspace{1mm}
\noindent
where $i_1=\ldots=i_k=i\in\{1,\ldots,m\}$ and $\hat\Phi(t_1,\ldots,t_k)$ is a symmetrization
of the function $\Phi(t_1,\ldots,t_k)$.

Consider the elementary inequality
\begin{equation}
\label{ch12026m300}
~~~~~~~~~~~ \left(a_1+a_2+\ldots+a_p\right)^2 \le
p\left(a_1^2+a_2^2+\ldots+a_p^2\right),\ \ \ p\in {\bf N}.
\end{equation}

Using the inequality (\ref{ch12026m300}) and (\ref{jjjye}), we obtain 
$$
{\sf M}\left\{\left(
J'[R_{p_1\ldots p_k}]_{T,t}^{(i_1\ldots i_k)}\right)^{2}\right\}=
$$
$$
=
{\sf M}\left\{\left(
\sum_{(t_1,\ldots,t_k)}
\int\limits_{t}^{T}
\ldots
\int\limits_{t}^{t_2}
R_{p_1\ldots p_k}(t_1,\ldots,t_k)
d{\bf w}_{t_1}^{(i_1)}
\ldots
d{\bf w}_{t_k}^{(i_k)}\right)^2\right\}\le
$$
$$
\le k!
\sum_{(t_1,\ldots,t_k)}
{\sf M}\left\{\left(
\int\limits_{t}^{T}
\ldots
\int\limits_{t}^{t_2}
R_{p_1\ldots p_k}(t_1,\ldots,t_k)
d{\bf w}_{t_1}^{(i_1)}
\ldots
d{\bf w}_{t_k}^{(i_k)}\right)^2\right\}=
$$
$$
= k!
\sum_{(t_1,\ldots,t_k)}
\int\limits_{t}^{T}
\ldots
\int\limits_{t}^{t_2}
R^2_{p_1\ldots p_k}(t_1,\ldots,t_k)
dt_1
\ldots
dt_k=
$$

\vspace{1mm}
\begin{equation}
\label{ch12026102}
=k!
\int\limits_{[t, T]^k}
R^2_{p_1\ldots p_k}(t_1,\ldots,t_k)
dt_1
\ldots
dt_k,
\end{equation}

\noindent
where $i_1,\ldots,i_k=1,\ldots,m.$

Suppose that $\{\phi_j(x)\}_{j=0}^{\infty}$
is a complete orthonormal system of functions in 
the space $L_2([t, T])$. 
Using the orthonormality of the functions $\phi_j(x)$ $(j=0, 1, 2,\ldots),$ 
we obtain
$$
\int\limits_{[t, T]^k}
R^2_{p_1\ldots p_k}(t_1,\dots,t_k)dt_1\ldots dt_k
=
$$
$$
=\int\limits_{[t,T]^k}
\Biggl(K(t_1,\ldots,t_k)-
\sum_{j_1=0}^{p_1}\ldots
\sum_{j_k=0}^{p_k}
C_{j_k\ldots j_1}
\prod_{l=1}^k\phi_{j_l}(t_l)\Biggr)^2
dt_1
\ldots
dt_k=
$$

\vspace{1mm}
$$
=\int\limits_{[t,T]^k}
K^2(t_1,\ldots,t_k)
dt_1\ldots dt_k -
$$
$$
- 2\int\limits_{[t,T]^k}
K(t_1,\ldots,t_k)\sum_{j_1=0}^{p_1}\ldots
\sum_{j_k=0}^{p_k}
C_{j_k\ldots j_1}
\prod_{l=1}^k\phi_{j_l}(t_l)
dt_1\ldots dt_k+
$$
$$
+\int\limits_{[t,T]^k}
\Biggl(\sum_{j_1=0}^{p_1}\ldots
\sum_{j_k=0}^{p_k}
C_{j_k\ldots j_1}
\prod_{l=1}^k\phi_{j_l}(t_l)\Biggr)^2
dt_1\ldots dt_k=
$$
$$
=\int\limits_{[t,T]^k}
K^2(t_1,\ldots,t_k)
dt_1\ldots dt_k - 
$$
$$
- 2 \sum_{j_1=0}^{p_1}\ldots
\sum_{j_k=0}^{p_k}C_{j_k\ldots j_1}\int\limits_{[t,T]^k}
K(t_1,\ldots,t_k)
\prod_{l=1}^k\phi_{j_l}(t_l)
dt_1\ldots dt_k+
$$
$$
+
\sum_{j_1=0}^{p_1}\sum_{j_1'=0}^{p_1}\ldots
\sum_{j_k=0}^{p_k}\sum_{j_k'=0}^{p_k}
C_{j_k\ldots j_1}C_{j_k'\ldots j_1'}
\prod_{l=1}^k\int\limits_{\stackrel{~}{t}}^T\phi_{j_l}(t_l)
\phi_{j_l'}(t_l)dt_l=
$$

\vspace{-1mm}
$$
=\int\limits_{[t,T]^k}
K^2(t_1,\ldots,t_k)
dt_1\ldots dt_k 
-2 \sum_{j_1=0}^{p_1}\ldots
\sum_{j_k=0}^{p_k}C^2_{j_k\ldots j_1}
+\sum_{j_1=0}^{p_1}\ldots
\sum\limits_{j_k=0}^{p_k}
C^2_{j_k\ldots j_1}=
$$
\begin{equation}
\label{dobav100}
=\int\limits_{[t,T]^k}
K^2(t_1,\ldots,t_k)
dt_1\ldots dt_k -\sum_{j_1=0}^{p_1}\ldots
\sum_{j_k=0}^{p_k}C^2_{j_k\ldots j_1}.
\end{equation}

\vspace{2mm}

Let us substitute (\ref{dobav100}) into (\ref{ch12026102})

\vspace{-2mm}
$$
{\sf M}\left\{\left(
J'[R_{p_1\ldots p_k}]_{T,t}^{(i_1\ldots i_k)}\right)^{2}\right\}\le
$$
\begin{equation}
\label{ch12026103}
~~~~~\le  k!
\left(~\int\limits_{[t,T]^k}
K^2(t_1,\ldots,t_k)
dt_1\ldots dt_k -\sum_{j_1=0}^{p_1}\ldots
\sum_{j_k=0}^{p_k}C^2_{j_k\ldots j_1}\right),
\end{equation}

\vspace{2mm}
\noindent
where $i_1,\ldots,i_k=1,\ldots,m.$

Due to Parseval's equality

\vspace{-2mm}
$$
\int\limits_{[t, T]^k}
R^2_{p_1\ldots p_k}(t_1,\dots,t_k)dt_1\ldots dt_k=
$$
\begin{equation}
\label{ziko11000}
~~~~~~=\int\limits_{[t,T]^k}
K^2(t_1,\ldots,t_k)
dt_1\ldots dt_k -\sum_{j_1=0}^{p_1}\ldots
\sum_{j_k=0}^{p_k}C^2_{j_k\ldots j_1}\  \to \  0
\end{equation}

\vspace{2mm}
\noindent
if $p_1,\ldots,p_k\to\infty.$

Combining (\ref{2026ch1001}) and (\ref{ch12026103}), we get
$$
{\sf M}\left\{\left(J'[R_{p_1\ldots p_k}]_{T,t}^{(i_1\ldots i_k)}\right)^{2n}\right\}\le
$$

\vspace{-1mm}
$$
\le
(k!)^{n}(2n-1)^{nk}\ \times
$$

\vspace{-7mm}
\begin{equation}
\label{2026ch1001s}
~~~~~~\times\ 
\left(~
\int\limits_{[t,T]^k}
K^2(t_1,\ldots,t_k)
dt_1\ldots dt_k -\sum_{j_1=0}^{p_1}\ldots
\sum_{j_k=0}^{p_k}C^2_{j_k\ldots j_1}
\right)^n,
\end{equation}

\vspace{4mm}
\noindent
or
$$
{\sf M}\left\{\left(J[\psi^{(k)}]_{T,t}-
J[\psi^{(k)}]_{T,t}^{p_1,\ldots,p_k}\right)^{2n}\right\}\le
$$

$$
\le
(k!)^{n} (2n-1)^{nk}\ \times
$$

\vspace{-7mm}
\begin{equation}
\label{dima2ye100}
~~~~~~\times\ 
\left(~
\int\limits_{[t,T]^k}
K^2(t_1,\ldots,t_k)
dt_1\ldots dt_k -\sum_{j_1=0}^{p_1}\ldots
\sum_{j_k=0}^{p_k}C^2_{j_k\ldots j_1}
\right)^n,
\end{equation}

\vspace{2mm}
\noindent
where $n\in{\bf N}$ and $i_1=\ldots=i_k=i\in\{1,\ldots,m\}.$

The inequality (\ref{2026ch1001s}) (or (\ref{dima2ye100}))
means that the expansion of 
iterated It\^{o} stochastic integral obtained using Theorem 1.1
(the case $k\in{\bf N},$ $i_1=\ldots=i_k=i\in\{1,\ldots,m\}$)
converges in the 
mean
of degree $2n$ ($n\in {\bf N}$) to the appropriate 
iterated It\^{o} stochastic integral.

Now we consider the case of a multidimensional Wiener process
and obtain an estimate of type (\ref{dima2ye100}) for the case $k=2,$ $i_1,i_2=1,\ldots,m.$

Suppose that $\{\phi_j(x)\}_{j=0}^{\infty}$
is a complete orthonormal system of continuous functions in 
the space $L_2([t, T])$ and $\psi_1(\tau),\psi_2(\tau)$
are continuous functions on $[t, T].$

Applying the Minkowski inequality and (\ref{jjjye}), we obtain 
$$
{\sf M}\left\{\left(J[\psi^{(2)}]_{T,t}-
J[\psi^{(2)}]_{T,t}^{p_1,p_2}\right)^{2n}\right\}=
{\sf M}\left\{\left(
J'[R_{p_1 p_2}]_{T,t}^{(i_1 i_2)}\right)^{2n}\right\}=
$$
$$
=
{\sf M}\left\{\left(
\sum_{(t_1,t_2)}
\int\limits_{t}^{T}
\int\limits_{t}^{t_2}
R_{p_1 p_2}(t_1,t_2)
d{\bf w}_{t_1}^{(i_1)}
d{\bf w}_{t_2}^{(i_2)}\right)^{2n}\right\}\le
$$
\begin{equation}
\label{ch12026e201}
~~~~~~~~~\le
\left(\sum_{(t_1,t_2)}
\left({\sf M}\left\{\left(
\int\limits_{t}^{T}
\int\limits_{t}^{t_2}
R_{p_1 p_2}(t_1,t_2)
d{\bf w}_{t_1}^{(i_1)}
d{\bf w}_{t_2}^{(i_2)}\right)^{2n}\right\}\right)^{1/2n}\right)^{2n},
\end{equation}

\vspace{2mm}
\noindent
where $n\in{\bf N}.$

Let us evaluate 
$$
{\sf M}\left\{\left(
\int\limits_{t}^{T}
\int\limits_{t}^{t_2}
R_{p_1 p_2}(t_1,t_2)
d{\bf w}_{t_1}^{(i_1)}
d{\bf w}_{t_2}^{(i_2)}\right)^{2n}\right\}.
$$

Denote
\begin{equation}
\label{ch1w1004}
\eta_s=\int\limits_{t}^{s}
\int\limits_{t}^{t_2}
R_{p_1 p_2}(t_1,t_2)
d{\bf w}_{t_1}^{(i_1)}
d{\bf w}_{t_2}^{(i_2)},\ \ \ s\in[t, T].
\end{equation}

We have
$$
d\eta_s=\xi_s d{\bf w}_{s}^{(i_2)},
$$

\vspace{2mm}
\noindent
where
\begin{equation}
\label{ch12026ner11}
\xi_s=\int\limits_{t}^s
R_{p_1 p_2}(t_1,s)
d{\bf w}_{t_1}^{(i_1)}.
\end{equation}

Using the It\^{o} formula it is easy to demonstrate that \cite{Gih1}

\vspace{-2mm}
$$
{\sf M}\left\{(\eta_{\tau})^{2n}
\right\}=n(2n-1)
{\sf M}\left\{\int\limits_{t}^{\tau}
\left(\eta_s\right)^{2n-2}\xi_s^2 ds\right\}=
$$
\begin{equation}
\label{ch12026qqq1}
=n(2n-1)
\int\limits_{t}^{\tau}
{\sf M}\left\{\left(\eta_s\right)^{2n-2}\xi_s^2\right\}ds.
\end{equation}

\vspace{2mm}

The last step in (\ref{ch12026qqq1}) is carried out
on the basis of a consequence from Fubini's Theorem \cite{Shir999},
since (as we will see later) 
\begin{equation}
\label{ch1ww2001}
{\sf M}\left\{\int\limits_{t}^{\tau}\left(\eta_s\right)^{2n-2}\xi_s^2 ds\right\}<\infty\ \ \
\hbox{for}\ \ \ p_1,p_2<\infty.
\end{equation}

Using the H\"{o}lder inequality (under the integral sign 
on the right-hand side of (\ref{ch12026qqq1})) for 
$p=n/(n-1)$, $q=n$ $(n>1)$ and using the non-decreasing property of 
the value ${\sf M}\left\{(\eta_{\tau})^{2n}\right\}$
with the 
growth of $\tau$ (see (\ref{ch12026qqq1})), we get
$$
{\sf M}\left\{\left(\eta_{\tau}\right)^{2n}
\right\}\le
n(2n-1)
\left({\sf M}\left\{\left(\eta_{\tau}\right)^{2n}
\right\}\right)^{(n-1)/n}
\int\limits_{t}^{\tau}\left({\sf M}\left\{(\xi_s)^{2n}\right\}\right)^{1/n}ds.
$$

After raising to power $n$ the obtained inequality and 
dividing the result by 
$$
\left({\sf M}\left\{\left(\eta_{\tau}\right)^{2n}
\right\}\right)^{n-1},
$$
we get the following estimate
\begin{equation}
\label{neogidal}
{\sf M}\left\{\left(\eta_{\tau}\right)^{2n}
\right\}\le
(n(2n-1))^n
\left(
\int\limits_{t}^{\tau}\left({\sf M}\left\{(\xi_s)^{2n}\right\}\right)^{1/n}ds
\right)^n.
\end{equation}

Note that
\begin{equation}
\label{ch12026r67}
{\sf M}\left\{(\xi_s)^{2n}\right\}=
(2n-1)!!\left(
\int\limits_t^{s}
R^2_{p_1 p_2}(t_1,s)dt_1
\right)^n,
\end{equation} 

\noindent
since the randon variable $\xi_s$ has a Gaussian distribution
and
$$
{\sf M}\left\{\xi_s^{2}\right\}={\sf M}\left\{\left(
\int\limits_{t}^s
R_{p_1 p_2}(t_1,s)
d{\bf w}_{t_1}^{(i_1)}\right)^2
\right\}=
$$
\begin{equation}
\label{ch12026art34}
=\int\limits_{t}^s
R^2_{p_1 p_2}(t_1,s)
dt_1.
\end{equation}

Combining (\ref{neogidal}) and (\ref{ch12026r67}), we obtain
\begin{equation}
\label{ch12026tttt1}
~~~~~~~{\sf M}\left\{\left(\eta_{\tau}\right)^{2n}
\right\}\le
(n(2n-1))^n (2n-1)!!
\left(\int\limits_{t}^{\tau}
\int\limits_t^{s}
R^2_{p_1 p_2}(t_1,s)dt_1 ds
\right)^n.
\end{equation}

\vspace{1mm}

Then
\begin{equation}
\label{ch12026r67d}
~~~~~~~{\sf M}\left\{\left(\eta_{T}\right)^{2n}
\right\}\le
(n(2n-1))^n (2n-1)!!
\left(~\int\limits_{[t, T]^2}
R^2_{p_1 p_2}(t_1,t_2)dt_1 dt_2
\right)^n.
\end{equation}

\vspace{2mm}

Finally, using (\ref{ch12026e201}) and (\ref{ch12026r67d}), we have
$$
{\sf M}\left\{\left(J[\psi^{(2)}]_{T,t}-
J[\psi^{(2)}]_{T,t}^{p_1 p_2}\right)^{2n}\right\}\le
$$
$$
\le
\left(\sum_{(t_1,t_2)}
\left({\sf M}\left\{\left(\eta_T\right)^{2n}
\right\}\right)^{1/2n}\right)^{2n}\le
$$
$$
\le
2^{2n}(n(2n-1))^n (2n-1)!! 
\left(~\int\limits_{[t, T]^2}
R^2_{p_1 p_2}(t_1,t_2)dt_1 dt_2
\right)^{n}=
$$

\vspace{2mm}
$$
=
2^{2n}(n(2n-1))^n (2n-1)!!\times
$$

\vspace{-1mm}
\begin{equation}
\label{ch12026ri12}
\times
\left(~
\int\limits_{[t,T]^2}
K^2(t_1,t_2)
dt_1dt_2 -\sum_{j_1=0}^{p_1}
\sum_{j_2=0}^{p_2}C^2_{j_2 j_1}
\right)^{n},
\end{equation}

\vspace{2mm}
\noindent
where $n\in{\bf N}.$

Let us show that for $p_1,p_2<\infty$ the following inequality

\vspace{-3.5mm}
\begin{equation}
\label{ch120204567}
{\sf M}\left\{(\eta_{\tau})^{2n}\right\}\le C<\infty
\end{equation}

\vspace{2mm}
\noindent
is satisfied, where $\tau\in [t, T]$ and $C$ is a constant
that depends on $p_1,p_2,T,n.$

Consider the following well known estimate
for the moments of the It\^{o} stochastic integral \cite{Gih1}

\vspace{-4mm}
\begin{equation}
\label{ch12026pupol1}
~~~~~~{\sf M}\left\{\left|\int\limits_{t}^T \phi_\tau
dw_\tau\right|^{2n}\right\} \le (T-t)^{n-1}\left(n(2n-1)\right)^n
\int\limits_{t}^T {\sf M}\left\{(\phi_\tau)^{2n}\right\}d\tau,
\end{equation}

\vspace{2mm}
\noindent
where the process $\phi_{\tau}$ is such that
$\left(\phi_{\tau}\right)^n\in{\rm M}_2
([t,T])$ and $w_{\tau}$ is a scalar standard Wiener 
process,
$n=1, 2,\ldots$ (definition of the class 
${\rm M}_2([t,T])$ see in Sect.~1.1.2).

Applying (\ref{ch12026pupol1}), we obtain
\begin{equation}
\label{ch12026pupol12}
{\sf M}\left\{(\eta_{\tau})^{2n}\right\} \le (T-t)^{n-1}\left(n(2n-1)\right)^n
\int\limits_{t}^{\tau} {\sf M}\left\{(\xi_s)^{2n}\right\}ds,
\end{equation}
where $\xi_\tau$ is defined by (\ref{ch12026ner11}).

Combining (\ref{ch12026pupol12}) and (\ref{ch12026r67}), we get

\vspace{-1mm}
$$
{\sf M}\left\{(\eta_{\tau})^{2n}\right\} \le (T-t)^{n-1}\left(n(2n-1)\right)^n
(2n-1)!!\times
$$

\vspace{-2mm}
\begin{equation}
\label{ch12026pupol122}
\times
\int\limits_{t}^{\tau} \left(
\int\limits_t^{s}
R^2_{p_1 p_2}(t_1,s)dt_1
\right)^n ds.
\end{equation}

\vspace{2mm}

Under the conditions of Theorem~1.1, 
the integrals

\vspace{-2mm}
$$
\int\limits_t^{s}
R^2_{p_1 p_2}(t_1,s)dt_1,\ \ \ \int\limits_{t}^{\tau} \left(
\int\limits_t^{s}
R^2_{p_1 p_2}(t_1,s)dt_1
\right)^n ds
$$

\vspace{2mm}
\noindent
are continuous functions with respect to $s$ and $\tau$, respectively.

Thus, the estimate (\ref{ch120204567}) is proved (see (\ref{ch12026pupol122})).
Then, the inequality (\ref{ch1ww2001}) holds (see (\ref{ch12026qqq1})).
This means that for $p_1,p_2<\infty$
we can apply the consequence from Fubini's Theorem \cite{Shir999}
in (\ref{ch12026qqq1}) and therefore for $p_1,p_2<\infty$ the estimate (\ref{ch12026ri12})
will be true. The proof of the estimate (\ref{ch12026ri12}) is completed.

Let us explain why this approach cannot be generalized
to the case $k\ge 3.$
Let $k=3.$ Now
$$
\eta_s=\int\limits_{t}^{s}
\int\limits_{t}^{t_3}
\int\limits_{t}^{t_2}
R_{p_1 p_2 p_3}(t_1,t_2,t_3)
d{\bf w}_{t_1}^{(i_1)}
d{\bf w}_{t_2}^{(i_2)}d{\bf w}_{t_3}^{(i_3)},\ \ \ s\in[t, T],
$$
and
$$
d\eta_s=\xi_s
d{\bf w}_{s}^{(i_3)},
$$
where
$$
\xi_s=\int\limits_{t}^{s}
\int\limits_{t}^{t_2}
R_{p_1 p_2 p_3}(t_1,t_2,s)
d{\bf w}_{t_1}^{(i_1)}
d{\bf w}_{t_2}^{(i_2)}.
$$

In the next step, we need that the stochastic differential
$d\xi_s$ to have the form
\begin{equation}
\label{ch12020yyy1}
d\xi_s=\mu_s d{\bf w}_{s}^{(i_2)}.
\end{equation}

\vspace{2mm}

According to (\ref{ch12026ll12}), $\xi_s$ is a finite
linear combination of integrals of the form
$$
\rho_s=h(s)\int\limits_{t}^{s}g(t_2)
\int\limits_{t}^{t_2}
q(t_1)d{\bf w}_{t_1}^{(i_1)}
d{\bf w}_{t_2}^{(i_2)},
$$

\noindent
where $h(s), g(s), q(s)$ are some continuous functions on $[t, T].$

If we assume that the function $h(s)$ is continuously differentiable,
then according to the It\^{o} formula the stochastic differential
$d\rho_s$ will have a non-zero drift coefficient.
This means that the stochastic differential $d\xi_s$
will have a more compex form than (\ref{ch12020yyy1})
(with a non-zero drift coefficient), which makes
it impossible to generalize this approach to the case $k=3.$

Let us generalize the estimate (\ref{dima2ye100}) to the case $k\in{\bf N},$
$i_1,\ldots,i_k=1,\ldots,m.$
Let $H_1(s),\ldots,H_m(s): [0,\infty)\to {\bf R}$.
Let $(\Omega,{\rm F},{\sf P})$ is a probability space, where
${\rm F}$ is the smallest $\sigma$-algebra such that
the random variables

\vspace{-2mm}
$$
\sum\limits_{i=1}^m\int\limits_0^{\infty} H_i(s)d{\bf w}_s^{(i)}
$$

\vspace{1mm}
\noindent
are ${\rm F}$-measurable for every $H_1(s),\ldots,H_m(s)\in L_2([0,\infty)),$
where 
${\bf w}_s^{(i)}$ are independent standard Wiener processes, $i=1,\ldots,m,$

\vspace{-2mm}
$$
\int\limits_0^{\infty} H_i(s)d{\bf w}_s^{(i)}
$$

\vspace{1mm}
\noindent
is a usual Wiener--It\^{o} stochastic integral.

It is well known \cite{hairer1} that

\vspace{-2mm}
$$
L_2(\Omega,{\rm F},{\sf P})=\bigoplus\limits_{k=0}^{\infty}{\cal H}_k,
$$

\vspace{2mm}
\noindent
where ${\cal H}_0$ contains only constants,
the space ${\cal H}_k$ is the so-called $k$th $(k\ge 1)$ homogeneous Wiener chaos
which consists of all random variables of the form

\newpage
\noindent
$$
\sum\limits_{i_1,\ldots,i_k=1}^m 
\int\limits_0^{\infty}\int\limits_0^{t_k}\ldots \int\limits_0^{t_2}
H_{i_1\ldots i_k}(t_1,\ldots,t_k)d{\bf w}_{t_1}^{(i_1)}\ldots d{\bf w}_{t_k}^{(i_k)},
$$

\vspace{1mm}
\noindent
where $H_{i_1\ldots i_k}(t_1,\ldots,t_k)\in L_2(D_k),$
$D_k=\{(t_1,\ldots,t_k):\  0<t_1<\ldots <t_k\}.$

Let 
$$
H_{i_1\ldots i_k}(t_1,\ldots,t_k)=
\Phi(t_1,\ldots,t_k){\bf 1}_{\{0\le t<t_1<\ldots <t_k\le T\}}
{\bf 1}_{\{i_1=j_1,\ldots,i_k=j_k\}},
$$

\vspace{1mm}
\noindent
where ${\bf 1}_A$ denotes the indicator of the set $A,$
$\Phi(t_1,\ldots,t_k)\in L_2([t,T]^k),$
$j_1,\ldots,j_k$ are some fixed numbers from the set
$\{1,\ldots,m\}.$

Then, we have w.~p.~1
$$
\sum\limits_{i_1,\ldots,i_k=1}^m 
\int\limits_0^{\infty}\int\limits_0^{t_k}\ldots \int\limits_0^{t_2}
H_{i_1\ldots i_k}(t_1,\ldots,t_k)d{\bf w}_{t_1}^{(i_1)}\ldots {\bf w}_{t_k}^{(i_k)}=
$$
$$
=\int\limits_t^{T}\ldots \int\limits_t^{t_2}
\Phi(t_1,\ldots,t_k)d{\bf w}_{t_1}^{(j_1)}\ldots 
d{\bf w}_{t_k}^{(j_k)}\in {\cal H}_k.
$$

\vspace{1mm}

Obviously (see (\ref{pobeda}), (\ref{s2s}), (\ref{Wi110})),
$$
J'[\Phi]_{T,t}^{(i_1\ldots i_k)}=\sum_{(t_1,\ldots,t_k)}
\int\limits_{t}^{T}
\ldots
\int\limits_{t}^{t_2}
\Phi(t_1,\ldots,t_k)
d{\bf w}_{t_1}^{(i_1)}
\ldots
d{\bf w}_{t_k}^{(i_k)}=
$$
$$
=
\int\limits_t^T\ldots \int\limits_t^{t_2}
\sum\limits_{(t_1,\ldots,t_k)}\biggl(
\Phi(t_1,\ldots,t_k)
d{\bf w}_{t_1}^{(i_1)}\ldots
d{\bf w}_{t_k}^{(i_k)}\biggr)
\in {\cal H}_k,
$$

\noindent
where 
$J'[\Phi]_{T,t}^{(i_1\ldots i_k)}$ is a multiple
Wiener stochastic integral 
with respect to components of a multidimensional 
Wiener process
defined as in Sect.~1.1.3 (see (\ref{mult11})) or as in Sect.~1.11
(see (\ref{WiI})),
$\Phi(t_1,\ldots,t_k)\in L_2([t, T]^k)$ in (\ref{WiI}) and  
$\Phi(t_1,\ldots,t_k)\in C([t, T]^k)\subset L_2([t, T]^k)$ in (\ref{mult11}),
$i_1,\ldots,i_k=1,\ldots,m.$

It is well known that (see \cite{hairer1})

\vspace{-2mm}
$$
{\sf M}\left\{\xi^{2n}\right\}\le (2n-1)^{kn}\left({\sf M}\left\{\xi^2\right\}\right)^n,
$$

\vspace{2mm}
\noindent
where $\xi\in {\cal H}_k,$ $n, k\in {\bf N}.$

Then, we have the following estimate for 
$J'[\Phi]_{T,t}^{(i_1\ldots i_k)}\in {\cal H}_k$

\newpage
\noindent
\begin{equation}
\label{2026ch1001dd}
~~~~~~~~{\sf M}\left\{\left(
J'[\Phi]_{T,t}^{(i_1\ldots i_k)}\right)^{2n}\right\}
\le (2n-1)^{nk} \left({\sf M}\left\{\left(
J'[\Phi]_{T,t}^{(i_1\ldots i_k)}\right)^{2}\right\}\right)^{n},
\end{equation}

\noindent
where $n\in{\bf N}$ and $k\in{\bf N},$ $i_1,\ldots,i_k=1,\ldots,m.$

Combining (\ref{jjjye}), (\ref{ch12026103}), and (\ref{2026ch1001dd}), we get

\vspace{-1.5mm}
$$
{\sf M}\left\{\left(J[\psi^{(k)}]_{T,t}-
J[\psi^{(k)}]_{T,t}^{p_1,\ldots,p_k}\right)^{2n}\right\}\le
$$

\vspace{-1mm}
$$
\le
(k!)^{n} (2n-1)^{nk}\ \times
$$

\vspace{-7mm}
\begin{equation}
\label{2026ch1001s11}
~~~~~~\times\ 
\left(~
\int\limits_{[t,T]^k}
K^2(t_1,\ldots,t_k)
dt_1\ldots dt_k -\sum_{j_1=0}^{p_1}\ldots
\sum_{j_k=0}^{p_k}C^2_{j_k\ldots j_1}
\right)^n,
\end{equation}

\noindent
where $n\in{\bf N}$ and $k\in{\bf N},$ $i_1,\ldots,i_k=1,\ldots,m.$

The inequality (\ref{2026ch1001s11})
means that the expansion of 
iterated It\^{o} stochastic integral obtained using Theorem 1.1
(the case $k\in{\bf N},$ $i_1,\ldots,i_k=1,\ldots,m$)
converges in the 
mean
of degree $2n$ ($n\in {\bf N}$) to the appropriate 
iterated It\^{o} stochastic integral.

\subsection{Conclusions}

Thus, we obtain the following useful possibilities and modifications
of the approach based on Theorem 1.1.${}^5$

1.\;There is an explicit formula (see (\ref{ppppa})) for calculation 
of expansion coefficients 
of the iterated It\^{o} stochastic integral (\ref{ito}) with any
fixed multiplicity $k$ ($k\in{\bf N}$).

2.\;We have possibilities for exact calculation of the mean-square 
approximation error
of the iterated It\^{o} stochastic integral (\ref{ito}) \cite{12a}-\cite{art-1},
\cite{arxiv-3} (see Sect.~1.2).

3.\;Since the used
multiple Fourier series is a generalized in the sense
that it is built using various complete orthonormal
systems of functions in the space $L_2([t, T])$, then we 
have new possibilities 
for approximation --- we can 
use not only the trigonometric functions as in \cite{Zapad-1}-\cite{Zapad-4},
\cite{Zapad-8}, \cite{Zapad-9}, \cite{Zapad-11}, \cite{Zapad-12a},
but the Legendre polynomials.

\footnotetext[5]{Theorem~1.1 will be generalized to the case of an arbitrary 
complete ortho\-nor\-mal system of functions $\{\phi_j(x)\}_{j=0}^{\infty}$ 
in the space $L_2([t, T])$
and $\psi_1(\tau),$ $\ldots,\psi_k(\tau) \in L_2([t, T])$ in Sect.~1.11
(see Theorem~1.16).}

4.\;As it turned out \cite{1}-\cite{new-new-6}, 
it is more convenient to work 
with Legendre polynomials for approximation
of the iterated It\^{o} stochastic integrals (\ref{ito}) (see Chapter 5). 
Approximations based on Legendre polynomials essentially simpler 
than their analogues based on trigonometric functions
\cite{1}-\cite{new-new-6}.
Another advantages of the application of Legendre polynomials 
in the framework of the mentioned problem are considered
in \cite{art-4}, \cite{arxiv-12} (see Sect.~5.3).

5.\;The Milstein approach \cite{Zapad-1} (see Sect.~6.2 in this book)
to expansion of iterated 
stochastic integrals based on the Karhunen--Lo\`{e}ve expansion
of the Brownian bridge process (also see \cite{Zapad-2}-\cite{Zapad-4},
\cite{Zapad-8}, \cite{Zapad-9}, \cite{Zapad-11}, \cite{Zapad-12a})
leads to 
iterated application of the operation of limit
transition (the operation of limit transition 
is implemented only once in Theorem 1.1)
starting from the 
second or third multiplicity
of the iterated It\^{o} stochastic integral (\ref{ito}).
Multiple series (the operation of limit transition 
is implemented only once) are more convenient 
for approximation than the iterated ones
(iterated application of the operation of limit
transition), 
since partial sums of multiple series converge for any possible case of  
convergence to infinity of their upper limits of summation 
(let us denote them as $p_1,\ldots, p_k$). 
For example, when 
$p_1=\ldots=p_k=p\to\infty$. 
For iterated series, the condition $p_1=\ldots=p_k=p\to\infty$ obviously 
does not guarantee the convergence of this series.
However, in \cite{Zapad-2}-\cite{Zapad-4}, \cite{Zapad-9} 
the authors use (without rigorous proof)
the condition $p_1=p_2=p_3=p\to\infty$
within the frames of the Milstein approach \cite{Zapad-1}
together with the
Wong--Zakai approximation  \cite{W-Z-1}-\cite{Watanabe}
(see discussions in Sect.~2.42, 2.43, 6.2).

6.\;As we mentioned above,
constructing the expansions of iterated
It\^{o} stochastic integrals from Theorem 1.1 we 
saved all 
information about these integrals. That is why it is 
natural to expect that the mentioned expansions will converge
with probability 1. The convergence with probability 1
in Theorem 1.1 has been proved for some
particular cases in \cite{3}-\cite{12aa},
\cite{arxiv-4} (see Sect.~1.7.1) and for the general case
of iterated
It\^{o} stochastic integrals
of multiplicity $k$ $(k\in{\bf N})$
in \cite{12a}-\cite{12aa}, \cite{OK1000}, \cite{arxiv-1}, \cite{arxiv-3},
\cite{arxiv-4} (see Sect.~1.7.2).

7.\;The generalizations of Theorem 1.1 for an arbitrary 
complete orthonormal system of functions in 
$L_2([t,T]^k)$ \cite{arxiv-1}
and complete 
orthonormal with weight  
$r(t_1)\ldots r(t_k)\ge 0$ systems of functions in 
$L_2([t,T]^k)$ 
\cite{11}-\cite{12aa}, \cite{arxiv-13}
as well 
as for iterated stochastic integrals 
with respect to martingale Poisson measures and 
iterated stochastic integrals with respect 
to martingales \cite{1}-\cite{12aa}, \cite{arxiv-13}
are presented in 
Sect.~1.3--1.6, 1.11.

8.\;The adaptation of Theorem 1.1 for iterated Stratonovich
stochastic integrals was carried out in 
\cite{6}-\cite{art-6}, \cite{art-9}, \cite{arxiv-2}, 
\cite{arxiv-4}-\cite{arxiv-11}, \cite{arxiv-15},
\cite{arxiv-17}-\cite{arxiv-19}, \cite{arxiv-22}, \cite{arxiv-24},
\cite{new-art-1xxy}, 
\cite{new-art-1xxys} (see Chapter 2).

9.\;Application of Theorem 1.1 for the mean-square
approximation of iterated stochastic integrals 
with respect to the 
infinite-dimensional $Q$-Wiener process can be found
in \cite{12a}-\cite{12aa}, \cite{art-7}, \cite{OK}, \cite{arxiv-20}, \cite{arxiv-21}
(see Chapter 7).

\section{Exact Calculation of the Mean-Square Error in the
Method of Approximation
of Iterated It\^{o} Stochastic integrals Based 
on Generalized Multiple Fourier Series}

This section is devoted
to the obtainment of 
exact and approximate 
expressions for the mean-square
approximation error in Theorem 1.1 for 
iterated It\^{o} stochastic integrals of
arbitrary multiplicity $k$ ($k\in {\bf N}$).
As a result, we do not need to use redundant terms
of expansions of iterated It\^{o} stochastic integrals.

\subsection{Introduction}

Recall that we called the method 
of expansion
and mean-square approximation of iterated It\^{o} stochastic integrals
based on Theorem 1.1 as the method of generalized multiple Fourier series.
The question about how estimate or even calculate 
exactly the mean-square approximation error of iterated It\^{o} stochastic 
integrals
for the method of generalized 
multiple Fourier series composes the subject of Sect.~1.2.
From the one side the mentioned question is essentially difficult
in the case of a multidimensional Wiener process, because of we need to take
into account all possible combinations of the components of a
multidimensional Wiener
process. From the other side an effective solution of the 
mentioned problem allows to construct
more simple expansions 
of iterated It\^{o} stochastic 
integrals than in \cite{Zapad-1}-\cite{Zapad-5}, 
\cite{Zapad-8}-\cite{Zapad-10}, \cite{Zapad-11}, \cite{Zapad-12a}.

Sect.~1.2.2 is devoted to the formulation and proof of 
Theorem 1.3, which allows to calculate exacly
the mean-square approximation error of iterated
It\^{o} stochastic integrals of arbitrary multiplicity $k$ ($k\in{\bf N}$)
for the method of generalized multiple Fourier series.
The particular cases ($k=1,\dots,5$) of Theorem 1.3 
are considered
in detail in Sect.~1.2.3.
In Sect.~1.2.4 we prove an effective estimate 
for the mean-square approximation error of iterated
It\^{o} stochastic integrals of arbitrary multiplicity $k$ ($k\in{\bf N}$)
for the method of generalized multiple Fourier series.

\subsection{Theorem on 
Exact Calculation of the Mean-Square Approximation Error 
for Iterated It\^{o} Stochastic integrals}

\footnotetext[6]{Theorem 1.3 will be generalized to the case of an arbitrary 
complete ortho\-nor\-mal system of functions $\{\phi_j(x)\}_{j=0}^{\infty}$ 
in the space $L_2([t, T])$
and $\psi_1(\tau),$ $\ldots,\psi_k(\tau) \in L_2([t, T])$ in Sect.~1.12
(see Theorem~1.18).}

{\bf Theorem 1.3}${}^6$ \cite{11}-\cite{art-1}, \cite{arxiv-3}. 
{\it Suppose that
every $\psi_l(\tau)$ $(l=1,\ldots, k)$ is a continuous nonrandom function on 
$[t, T]$ and
$\{\phi_j(x)\}_{j=0}^{\infty}$ is a complete orthonormal system  
of functions in the space $L_2([t,T]),$ 
each function $\phi_j(x)$ of which 
for finite $j$ satisfies the condition 
$(\star)$ {\rm (}see Sect.~{\rm 1.1.7)}.
Then

\vspace{-3mm}
$$
{\sf M}\left\{\left(J[\psi^{(k)}]_{T,t}-
J[\psi^{(k)}]_{T,t}^p\right)^2\right\}
= \int\limits_{[t,T]^k} K^2(t_1,\ldots,t_k)
dt_1\ldots dt_k - 
$$

\vspace{-3mm}
\begin{equation}
\label{tttr11}
- \hspace{-0.5mm}\sum_{j_1=0}^{p}\ldots\sum_{j_k=0}^{p}
C_{j_k\ldots j_1}
{\sf M}\left\{\hspace{-0.5mm}J[\psi^{(k)}]_{T,t}
\sum\limits_{(j_1,\ldots,j_k)}
\int\limits_t^T \phi_{j_k}(t_k)
\ldots
\int\limits_t^{t_{2}}\phi_{j_{1}}(t_{1})
d{\bf w}_{t_1}^{(i_1)}\ldots
d{\bf w}_{t_k}^{(i_k)}\hspace{-0.5mm}\right\}\hspace{-0.5mm},
\end{equation}

\vspace{2mm}
\noindent
where $i_1,\ldots,i_k=1,\ldots,m,$
$$
J[\psi^{(k)}]_{T,t}=\int\limits_t^T\psi_k(t_k) \ldots \int\limits_t^{t_{2}}
\psi_1(t_1) d{\bf w}_{t_1}^{(i_1)}\ldots
d{\bf w}_{t_k}^{(i_k)},
$$

\vspace{-2mm}
\begin{equation}
\label{yeee2}
J[\psi^{(k)}]_{T,t}^p=
\sum_{j_1=0}^{p}\ldots\sum_{j_k=0}^{p}
C_{j_k\ldots j_1}\left(
\prod_{l=1}^k\zeta_{j_l}^{(i_l)}-S_{j_1,\ldots,j_k}^{(i_1\ldots i_k)}
\right),
\end{equation}

\vspace{2mm}

\begin{equation}
\label{ppp1}
~~~~~~~~S_{j_1,\ldots,j_k}^{(i_1\ldots i_k)}=
\hbox{\vtop{\offinterlineskip\halign{
\hfil#\hfil\cr
{\rm l.i.m.}\cr
$\stackrel{}{{}_{N\to \infty}}$\cr
}} }\sum_{(l_1,\ldots,l_k)\in {\rm G}_k}
\phi_{j_{1}}(\tau_{l_1})
\Delta{\bf w}_{\tau_{l_1}}^{(i_1)}\ldots
\phi_{j_{k}}(\tau_{l_k})
\Delta{\bf w}_{\tau_{l_k}}^{(i_k)},
\end{equation}

\vspace{5mm}
\noindent
the Fourier coefficient $C_{j_k\ldots j_1}$ has the form {\rm (\ref{ppppa})},
\begin{equation}
\label{rr232}
\zeta_{j}^{(i)}=
\int\limits_t^T \phi_{j}(s) d{\bf w}_s^{(i)}
\end{equation}
are independent standard Gaussian random variables
for various
$i$ or $j$ $(i=1,\ldots,m),$
$$
\sum\limits_{(j_1,\ldots,j_k)}
$$ 
means the sum with respect to all
possible permutations 
$(j_1,\ldots,j_k).$ At the same time if 
$j_r$ swapped with $j_q$ in the permutation $(j_1,\ldots,j_k)$,
then $i_r$ swapped with $i_q$ in the permutation
$(i_1,\ldots,i_k)$ {\rm (}see {\rm (\ref{tttr11}));}
another notations are the same as in Theorem {\rm 1.1.}}

{\bf Remark 1.3.}\ {\it Note that

\vspace{-3mm}
$$
{\sf M}\left\{J[\psi^{(k)}]_{T,t}
\int\limits_t^T \phi_{j_k}(t_k)
\ldots
\int\limits_t^{t_{2}}\phi_{j_{1}}(t_{1})
d{\bf w}_{t_1}^{(i_1)}\ldots
d{\bf w}_{t_k}^{(i_k)}\right\}=
$$
$$
\hspace{-60mm}={\sf M}\left\{\int\limits_t^T\psi_k(t_k) \ldots \int\limits_t^{t_{2}}
\psi_1(t_1) d{\bf w}_{t_1}^{(i_1)}\ldots
d{\bf w}_{t_k}^{(i_k)}\times\right.
$$
$$
~~~~~~~~~~~~~~~~~~~~~~~~~~~~~~~~~\left.\times
\int\limits_t^T \phi_{j_k}(t_k)
\ldots
\int\limits_t^{t_{2}}\phi_{j_{1}}(t_{1})
d{\bf w}_{t_1}^{(i_1)}\ldots
d{\bf w}_{t_k}^{(i_k)}\right\}=
$$
\begin{equation}
\label{rwwr2}
~~~~~=\int\limits_t^T\psi_k(t_k) \phi_{j_k}(t_k)\ldots \int\limits_t^{t_{2}}
\psi_1(t_1)\phi_{j_1}(t_1) dt_1\ldots dt_k=
C_{j_k\ldots j_1}.
\end{equation}

\vspace{1mm}

Therefore, in the case of pairwise different numbers
$i_1,\ldots,i_k$ from Theorem {\rm 1.3} we obtain
\begin{equation}
\label{pochti}
{\sf M}\left\{\left(J[\psi^{(k)}]_{T,t}-
J[\psi^{(k)}]_{T,t}^p\right)^2\right\}=\hspace{-0.5mm}
\int\limits_{[t,T]^k} \hspace{-1.5mm} K^2(t_1,\ldots,t_k)
dt_1\ldots dt_k - \sum_{j_1=0}^{p}\ldots\sum_{j_k=0}^{p}
C_{j_k\ldots j_1}^2.
\end{equation}

Moreover, if $i_1=\ldots=i_k,$ then from Theorem {\rm 1.3} we get

$$
{\sf M}\left\{\left(J[\psi^{(k)}]_{T,t}-
J[\psi^{(k)}]_{T,t}^p\right)^2\right\}=
$$

$$
=\int\limits_{[t,T]^k} K^2(t_1,\ldots,t_k)
dt_1\ldots dt_k - \sum_{j_1=0}^{p}\ldots\sum_{j_k=0}^{p}
C_{j_k\ldots j_1}\Biggl(\sum\limits_{(j_1,\ldots,j_k)}
C_{j_k\ldots j_1}\Biggr),
$$

\vspace{4mm}
\noindent
where

\vspace{-4mm}
$$
\sum\limits_{(j_1,\ldots,j_k)}
$$ 

\newpage
\noindent
means the sum with respect to all
possible permutations
$(j_1,\ldots,j_k).$

For example, for the case $k=3$ we have
$$
{\sf M}\left\{\left(J[\psi^{(3)}]_{T,t}-
J[\psi^{(3)}]_{T,t}^p\right)^2\right\}=
\int\limits_t^T\psi_3^2(t_3)\int\limits_t^{t_3}\psi_2^2(t_2)
\int\limits_t^{t_2}\psi_1^2(t_1)dt_1dt_2dt_3 -
$$
$$
- \sum_{j_1,j_2,j_3=0}^{p}
C_{j_3j_2j_1}\biggl(C_{j_3j_2j_1}+C_{j_3j_1j_2}+C_{j_2j_3j_1}+
C_{j_2j_1j_3}+C_{j_1j_2j_3}+C_{j_1j_3j_2}\biggr).
$$
}

\vspace{2mm}

{\bf Proof.}
Using Theorem 1.1 for the case $i_1,\ldots,i_k=1,\ldots,m$ 
and $p_1=\ldots=p_k=p$, we obtain
\begin{equation}
\label{yyye1}
~~~~~~~~ J[\psi^{(k)}]_{T,t}=\
\hbox{\vtop{\offinterlineskip\halign{
\hfil#\hfil\cr
{\rm l.i.m.}\cr
$\stackrel{}{{}_{p\to \infty}}$\cr
}} }\sum_{j_1=0}^{p}\ldots\sum_{j_k=0}^{p}
C_{j_k\ldots j_1}\left(
\prod_{l=1}^k\zeta_{j_l}^{(i_l)}-S_{j_1,\ldots,j_k}^{(i_1\ldots i_k)}
\right).
\end{equation}

For $n>p$ we can write 
$$
J[\psi^{(k)}]_{T,t}^n=
\left(\sum_{j_1=0}^{p}+\sum_{j_1=p+1}^n\right)\ldots
\left(\sum_{j_k=0}^{p}+\sum_{j_k=p+1}^n\right)
C_{j_k\ldots j_1}\left(
\prod_{l=1}^k\zeta_{j_l}^{(i_l)}-S_{j_1,\ldots,j_k}^{(i_1\ldots i_k)}
\right)=
$$

\begin{equation}
\label{yyye}
=J[\psi^{(k)}]_{T,t}^p + \xi[\psi^{(k)}]_{T,t}^{p+1,n}.
\end{equation}

\vspace{1mm}

Let us prove that due to the special structure of random variables 
$S_{j_1,\ldots,j_k}^{(i_1\ldots i_k)}$ (see (\ref{a1})--(\ref{a7}), (\ref{leto6000}),
(\ref{ppp1}))
the following relations 
are correct
\begin{equation}
\label{tyty}
{\sf M}\left\{
\prod_{l=1}^k\zeta_{j_l}^{(i_l)}-S_{j_1,\ldots,j_k}^{(i_1\ldots i_k)}
\right\}=0,
\end{equation}
\begin{equation}
\label{tyty1}
~~~~~~~~~~{\sf M}\left\{
\left(\prod_{l=1}^k\zeta_{j_l}^{(i_l)}-S_{j_1,\ldots,j_k}^{(i_1\ldots i_k)}
\right)
\left(\prod_{l=1}^k\zeta_{j_l'}^{(i_l)}-S_{j_1',\ldots,j_k'}^{(i_1\ldots i_k)}
\right)\right\}=0,
\end{equation}

\vspace{2mm}
\noindent
where
$$
(j_1,\ldots,j_k)\in{\rm K}_p,\ \ \ (j_1',\ldots,j_k')
\in{\rm K}_n\backslash {\rm K}_{p}
$$ 

\vspace{1mm}
\noindent
and
$$
{\rm K}_n=\left\{(j_1,\ldots,j_k):\ 0\le j_1,\ldots,j_k\le n\right\},
$$

\newpage
\noindent
$$
{\rm K}_p=\left\{(j_1,\ldots,j_k):\ 0\le j_1,\ldots,j_k\le p\right\}.
$$

\vspace{3mm}

For the case $i_1,\ldots,i_k=1,\ldots,m$ 
and $p_1=\ldots=p_k=p$ from (\ref{novoe1}),
(\ref{novoe2}) (see the proof of Theorem 1.1) we obtain

\vspace{-4mm}
$$
\prod_{l=1}^k\zeta_{j_l}^{(i_l)}-
S_{j_1,\ldots,j_k}^{(i_1\ldots i_k)} =
\hbox{\vtop{\offinterlineskip\halign{
\hfil#\hfil\cr
{\rm l.i.m.}\cr
$\stackrel{}{{}_{N\to \infty}}$\cr
}} }
\sum_{\stackrel{l_1,\ldots,l_k=0}
{{}_{l_q\ne l_r;\ q\ne r;\
q,r=1,\ldots,k}}}^{N-1}
\phi_{j_1}(\tau_{l_1})\ldots
\phi_{j_k}(\tau_{l_k})
\Delta{\bf w}_{\tau_{l_1}}^{(i_1)}
\ldots
\Delta{\bf w}_{\tau_{l_k}}^{(i_k)}=
$$

\begin{equation}
\label{ttt2}
~~~~~~=\sum\limits_{(j_1,\ldots,j_k)}
\int\limits_t^T \phi_{j_k}(t_k)
\ldots
\int\limits_t^{t_{2}}\phi_{j_{1}}(t_{1})
d{\bf w}_{t_1}^{(i_1)}\ldots
d{\bf w}_{t_k}^{(i_k)}\ \ \ {\rm w.~p.~1,}
\end{equation}

\vspace{3mm}
\noindent
where 
$$
\sum\limits_{(j_1,\ldots,j_k)}
$$ 
means the sum with respect to all
possible permutations
$(j_1,\ldots,j_k).$ At the same time if 
$j_r$ swapped  with $j_q$ in the permutation $(j_1,\ldots,j_k)$,
then $i_r$ swapped  with $i_q$ in the permutation
$(i_1,\ldots,i_k);$
another notations are the same as in Theorem 1.1.

So, we obtain (\ref{tyty})
from (\ref{ttt2}) due to the moment property of the It\^{o}
stochastic integral.

Let us prove (\ref{tyty1}).
From (\ref{ttt2}) we have
$$
0\le \left\vert{\sf M}\left\{
\left(\prod_{l=1}^k\zeta_{j_l}^{(i_l)}-S_{j_1,\ldots,j_k}^{(i_1\ldots i_k)}
\right)
\left(\prod_{l=1}^k\zeta_{j_l'}^{(i_l)}-S_{j_1',\ldots,j_k'}^{(i_1\ldots i_k)}
\right)\right\}\right\vert=
$$
$$
=\left\vert
{\sf M}\left\{\sum\limits_{(j_1,\ldots,j_k)}\sum\limits_{(j_1',\ldots,j_k')}
\int\limits_t^T \phi_{j_k}(t_k)
\ldots
\int\limits_t^{t_{2}}\phi_{j_{1}}(t_{1})
d{\bf w}_{t_1}^{(i_1)}\ldots
d{\bf w}_{t_k}^{(i_k)}\times\right.\right.
$$
$$
\left.\left.\times
\int\limits_t^T \phi_{j_k'}(t_k)
\ldots
\int\limits_t^{t_{2}}\phi_{j_1'}(t_{1})
d{\bf w}_{t_1}^{(i_1)}\ldots
d{\bf w}_{t_k}^{(i_k)}\right\}\right\vert\le
$$
$$
\le\sum\limits_{(j_1',\ldots,j_k')}\ 
\int\limits_t^T \phi_{j_k}(t_k)\phi_{j_k'}(t_k)dt_k
\ldots
\int\limits_t^{T}\phi_{j_1}(t_{1})\phi_{j_1'}(t_{1})
dt_1=
$$
\begin{equation}
\label{hq11}
=
\sum\limits_{(j_1',\ldots,j_k')}
{\bf 1}_{\{j_1=j_1'\}}\ldots {\bf 1}_{\{j_k=j_k'\}},
\end{equation}

\noindent
where ${\bf 1}_A$ is the indicator of the set $A$.
From (\ref{hq11}) we obtain (\ref{tyty1}).

First, let us prove (\ref{hq11}) for the cases $k=2$ and $k=3$.
We have

\vspace{-4mm}
$$
{\sf M}\left\{\hspace{-0.5mm}\sum\limits_{(j_1,j_2)}\sum\limits_{(j_1',j_2')}
\int\limits_t^T\phi_{j_{2}}(t_{2})
\int\limits_t^{t_{2}}\phi_{j_{1}}(t_{1})
d{\bf w}_{t_1}^{(i_1)}d{\bf w}_{t_2}^{(i_2)}
\int\limits_t^T\phi_{j_{2}'}(t_{2})
\int\limits_t^{t_{2}}\phi_{j_{1}'}(t_{1})
d{\bf w}_{t_1}^{(i_1)}d{\bf w}_{t_2}^{(i_2)}\hspace{-0.5mm}
\right\}\hspace{-0.5mm}=
$$
$$
=
\int\limits_t^{T}\phi_{j_2}(s)\phi_{j_2'}(s)ds
\int\limits_t^{T}\phi_{j_1}(s)\phi_{j_1'}(s)ds+
$$
$$
+{\bf 1}_{\{i_1=i_2\}}
\int\limits_t^{T}\phi_{j_2}(s)\phi_{j_1'}(s)ds
\int\limits_t^{T}\phi_{j_1}(s)\phi_{j_2'}(s)ds=
$$
\begin{equation}
\label{gt22}
=
{\bf 1}_{\{j_1=j_1'\}}{\bf 1}_{\{j_2=j_2'\}}+
{\bf 1}_{\{i_1=i_2\}} \cdot
{\bf 1}_{\{j_2=j_1'\}}{\bf 1}_{\{j_1=j_2'\}},
\end{equation}

\vspace{2mm}

$$
{\sf M}\left\{\sum\limits_{(j_1,j_2,j_3)}\sum\limits_{(j_1',j_2',j_3')}
\int\limits_t^T \phi_{j_3}(t_3)
\int\limits_t^{t_{3}}\phi_{j_{2}}(t_{2})
\int\limits_t^{t_{2}}\phi_{j_{1}}(t_{1})
d{\bf w}_{t_1}^{(i_1)}d{\bf w}_{t_2}^{(i_2)}
d{\bf w}_{t_3}^{(i_3)}\times\right.
$$
$$
\left.\times
\int\limits_t^T \phi_{j_3'}(t_3)
\int\limits_t^{t_{3}}\phi_{j_{2}'}(t_{2})
\int\limits_t^{t_{2}}\phi_{j_{1}'}(t_{1})
d{\bf w}_{t_1}^{(i_1)}d{\bf w}_{t_2}^{(i_2)}
d{\bf w}_{t_3}^{(i_3)}\right\}=
$$
$$
=\int\limits_t^T \phi_{j_3}(s)\phi_{j_3'}(s)ds
\int\limits_t^{T}\phi_{j_2}(s)\phi_{j_2'}(s)ds
\int\limits_t^{T}\phi_{j_1}(s)\phi_{j_1'}(s)ds+
$$
$$
+{\bf 1}_{\{i_1=i_2\}}\int\limits_t^T \phi_{j_3}(s)\phi_{j_3'}(s)ds
\int\limits_t^{T}\phi_{j_1}(s)\phi_{j_2'}(s)ds
\int\limits_t^{T}\phi_{j_2}(s)\phi_{j_1'}(s)ds+
$$
$$
+{\bf 1}_{\{i_2=i_3\}}\int\limits_t^T \phi_{j_1}(s)\phi_{j_1'}(s)ds
\int\limits_t^{T}\phi_{j_2}(s)\phi_{j_3'}(s)ds
\int\limits_t^{T}\phi_{j_3}(s)\phi_{j_2'}(s)ds+
$$
$$
+{\bf 1}_{\{i_1=i_3\}}\int\limits_t^T \phi_{j_1}(s)\phi_{j_3'}(s)ds
\int\limits_t^{T}\phi_{j_2}(s)\phi_{j_2'}(s)ds
\int\limits_t^{T}\phi_{j_3}(s)\phi_{j_1'}(s)ds+
$$
$$
+{\bf 1}_{\{i_1=i_2=i_3\}}\int\limits_t^T \phi_{j_2}(s)\phi_{j_3'}(s)ds
\int\limits_t^{T}\phi_{j_1}(s)\phi_{j_2'}(s)ds
\int\limits_t^{T}\phi_{j_3}(s)\phi_{j_1'}(s)ds+
$$
$$
+{\bf 1}_{\{i_1=i_2=i_3\}}\int\limits_t^T \phi_{j_1}(s)\phi_{j_3'}(s)ds
\int\limits_t^{T}\phi_{j_3}(s)\phi_{j_2'}(s)ds
\int\limits_t^{T}\phi_{j_2}(s)\phi_{j_1'}(s)ds=
$$
$$
={\bf 1}_{\{j_3=j_3'\}}{\bf 1}_{\{j_2=j_2'\}}{\bf 1}_{\{j_1=j_1'\}}+
{\bf 1}_{\{i_1=i_2\}} \cdot
{\bf 1}_{\{j_3=j_3'\}}{\bf 1}_{\{j_1=j_2'\}}{\bf 1}_{\{j_2=j_1'\}}+
$$
$$
+{\bf 1}_{\{i_2=i_3\}} \cdot
{\bf 1}_{\{j_1=j_1'\}}{\bf 1}_{\{j_2=j_3'\}}{\bf 1}_{\{j_3=j_2'\}}+
{\bf 1}_{\{i_1=i_3\}} \cdot
{\bf 1}_{\{j_1=j_3'\}}{\bf 1}_{\{j_2=j_2'\}}{\bf 1}_{\{j_3=j_1'\}}+
$$
$$
+{\bf 1}_{\{i_1=i_2=i_3\}}\cdot
{\bf 1}_{\{j_2=j_3'\}}{\bf 1}_{\{j_1=j_2'\}}{\bf 1}_{\{j_3=j_1'\}}+
$$
\begin{equation}
\label{gt23}
+
{\bf 1}_{\{i_1=i_2=i_3\}}\cdot
{\bf 1}_{\{j_1=j_3'\}}{\bf 1}_{\{j_3=j_2'\}}{\bf 1}_{\{j_2=j_1'\}}.
\end{equation}

\vspace{5mm}

From (\ref{gt22}) and (\ref{gt23}) we get
$$
\left\vert{\sf M}\left\{\sum\limits_{(j_1,j_2)}\sum\limits_{(j_1',j_2')}
\int\limits_t^T\phi_{j_{2}}(t_{2})
\int\limits_t^{t_{2}}\phi_{j_{1}}(t_{1})
d{\bf w}_{t_1}^{(i_1)}d{\bf w}_{t_2}^{(i_2)}
\times\right.\right.
$$
$$
\left.\left.\times
\int\limits_t^T\phi_{j_{2}'}(t_{2})
\int\limits_t^{t_{2}}\phi_{j_{1}'}(t_{1})
d{\bf w}_{t_1}^{(i_1)}d{\bf w}_{t_2}^{(i_2)}
\right\}\right\vert\le
$$
$$
\le
{\bf 1}_{\{j_1=j_1'\}}{\bf 1}_{\{j_2=j_2'\}}+
{\bf 1}_{\{j_2=j_1'\}}{\bf 1}_{\{j_1=j_2'\}}=
$$
$$
=
\sum\limits_{(j_1',j_2')}
{\bf 1}_{\{j_1=j_1'\}}{\bf 1}_{\{j_2=j_2'\}},
$$

$$
\left\vert{\sf M}\left\{\sum\limits_{(j_1,j_2,j_3)}\sum\limits_{(j_1',j_2',j_3')}
\int\limits_t^T \phi_{j_3}(t_3)
\int\limits_t^{t_{3}}\phi_{j_{2}}(t_{2})
\int\limits_t^{t_{2}}\phi_{j_{1}}(t_{1})
d{\bf w}_{t_1}^{(i_1)}d{\bf w}_{t_2}^{(i_2)}
d{\bf w}_{t_3}^{(i_3)}\times\right.\right.
$$
$$
\left.\left.\times
\int\limits_t^T \phi_{j_3'}(t_3)
\int\limits_t^{t_{3}}\phi_{j_{2}'}(t_{2})
\int\limits_t^{t_{2}}\phi_{j_{1}'}(t_{1})
d{\bf w}_{t_1}^{(i_1)}d{\bf w}_{t_2}^{(i_2)}
d{\bf w}_{t_3}^{(i_3)}\right\}\right\vert\le
$$
$$
\le{\bf 1}_{\{j_3=j_3'\}}{\bf 1}_{\{j_2=j_2'\}}{\bf 1}_{\{j_1=j_1'\}}+
{\bf 1}_{\{j_3=j_3'\}}{\bf 1}_{\{j_1=j_2'\}}{\bf 1}_{\{j_2=j_1'\}}+
$$
$$
+
{\bf 1}_{\{j_1=j_1'\}}{\bf 1}_{\{j_2=j_3'\}}{\bf 1}_{\{j_3=j_2'\}}+
{\bf 1}_{\{j_1=j_3'\}}{\bf 1}_{\{j_2=j_2'\}}{\bf 1}_{\{j_3=j_1'\}}+
$$
$$
+
{\bf 1}_{\{j_2=j_3'\}}{\bf 1}_{\{j_1=j_2'\}}{\bf 1}_{\{j_3=j_1'\}}+
{\bf 1}_{\{j_1=j_3'\}}{\bf 1}_{\{j_3=j_2'\}}{\bf 1}_{\{j_2=j_1'\}}
=
$$
$$
=
\sum\limits_{(j_1',j_2',j_3')}
{\bf 1}_{\{j_1=j_1'\}}{\bf 1}_{\{j_2=j_2'\}}{\bf 1}_{\{j_3=j_3'\}},
$$

\vspace{2mm}
\noindent
where we used the relation

\vspace{-2mm}
$$
\int\limits_t^T
\phi_{i}(\tau)\phi_{j}(\tau)d\tau={\bf 1}_{\{i=j\}},\ \ \ i,j=0, 1, 2\ldots
$$

\vspace{2mm}

Now consider the case of an arbitrary $k\in {\bf N}.$ We have
$$
{\sf M}\left\{\sum\limits_{(j_1,\ldots,j_k)}\sum\limits_{(j_1',\ldots,j_k')}
\int\limits_t^T \phi_{j_k}(t_k)
\ldots
\int\limits_t^{t_{2}}\phi_{j_{1}}(t_{1})
d{\bf w}_{t_1}^{(i_1)}\ldots
d{\bf w}_{t_k}^{(i_k)}\times\right.
$$
$$
\left.\times
\int\limits_t^T \phi_{j_k'}(t_k)
\ldots
\int\limits_t^{t_{2}}\phi_{j_1'}(t_{1})
d{\bf w}_{t_1}^{(i_1)}\ldots
d{\bf w}_{t_k}^{(i_k)}\right\}=
$$
$$
={\sf M}\left\{\sum\limits_{(j_1,\ldots,j_k)}\sum\limits_{(j_1',\ldots,j_k')}
\int\limits_t^T \phi_{j_k}(t_k)
\ldots
\int\limits_t^{t_{2}}\phi_{j_{1}}(t_{1})
d{\bf w}_{t_1}^{(i_1)}\ldots
d{\bf w}_{t_k}^{(i_k)}\times\right.
$$
$$
\left.\times
\int\limits_t^T \phi_{j_k'}(t_k)
\ldots
\int\limits_t^{t_{2}}\phi_{j_1'}(t_{1})
d{\bf w}_{t_1}^{(i_1')}\ldots
d{\bf w}_{t_k}^{(i_k')}\right\}=
$$

$$
=\sum\limits_{(j_1,\ldots,j_k)}\sum\limits_{(j_1',\ldots,j_k')}
{\bf 1}_{\{i_k=i_k'\}}\ldots {\bf 1}_{\{i_1=i_1'\}}\times
$$
$$
\times
\int\limits_t^T \phi_{j_k}(t_k)\phi_{j_k'}(t_k)
\ldots
\int\limits_t^{t_{2}}\phi_{j_1}(t_1)\phi_{j_1'}(t_1)
dt_1\ldots dt_k=
$$
$$
=\sum\limits_{(j_1',\ldots,j_k')}
{\bf 1}_{\{i_k=i_k'\}}\ldots {\bf 1}_{\{i_1=i_1'\}}
\int\limits_t^T \phi_{j_k}(t_k)\phi_{j_k'}(t_k)dt_k
\ldots
\int\limits_t^T\phi_{j_1}(t_1)\phi_{j_1'}(t_1)dt_1=
$$
\begin{equation}
\label{wen100}
=\sum\limits_{(j_1',\ldots,j_k')}
{\bf 1}_{\{i_k=i_k'\}}\ldots {\bf 1}_{\{i_1=i_1'\}}
{\bf 1}_{\{j_k=j_k'\}}\ldots {\bf 1}_{\{j_1=j_1'\}},
\end{equation}

\vspace{3mm}
\noindent
where $(i_1',\ldots,i_k')=(i_1,\ldots,i_k).$
However, if $j_r'$ swapped  with $j_q'$ in the permutation $(j_1',\ldots,j_k')$,
then $i_r'$ swapped  with $i_q'$ in the permutation
$(i_1',\ldots,i_k')$ and if
$j_r$ swapped  with $j_q$ in the permutation $(j_1,\ldots,j_k)$,
then $i_r$ swapped  with $i_q$ in the permutation
$(i_1,\ldots,i_k)$.

From (\ref{wen100}) we obtain (\ref{hq11}).
The equality (\ref{tyty1}) is proved.

Note that the formula (\ref{tyty1}) (in the light of the results of 
Sect.~1.10, 1.11) can be 
interpreted as a consequence of the orthogonality of two random 
variables that are Hermite polynomials of vector random arguments.

From (\ref{tyty1}) we obtain

\vspace{1mm}
$$
{\sf M}\left\{J[\psi^{(k)}]_{T,t}^p\xi[\psi^{(k)}]_{T,t}^{p+1,n}
\right\}=0.
$$

\vspace{6mm}

Due to (\ref{yeee2}), (\ref{yyye1}), and (\ref{yyye}) we can write 

\vspace{0.5mm}
$$
\xi[\psi^{(k)}]_{T,t}^{p+1,n}=J[\psi^{(k)}]_{T,t}^n-J[\psi^{(k)}]_{T,t}^p,
$$

\vspace{3mm}
$$
\hbox{\vtop{\offinterlineskip\halign{
\hfil#\hfil\cr
{\rm l.i.m.}\cr
$\stackrel{}{{}_{n\to \infty}}$\cr
}} }
\xi[\psi^{(k)}]_{T,t}^{p+1,n}=J[\psi^{(k)}]_{T,t}-J[\psi^{(k)}]_{T,t}^p
\stackrel{\rm def}{=}\xi[\psi^{(k)}]_{T,t}^{p+1}.
$$

\vspace{7.5mm}

We have

\vspace{-2mm}
$$
0\le \left|{\sf M}\left\{
\xi[\psi^{(k)}]_{T,t}^{p+1}J[\psi^{(k)}]_{T,t}^p\right\}\right|=
$$

\vspace{4mm}
$$
=
\left|{\sf M}\left\{\left(
\xi[\psi^{(k)}]_{T,t}^{p+1}-
\xi[\psi^{(k)}]_{T,t}^{p+1,n}+\xi[\psi^{(k)}]_{T,t}^{p+1,n}\right)
J[\psi^{(k)}]_{T,t}^p\right\}\right|\le
$$

\vspace{2mm}
$$
\le 
\left|{\sf M}\left\{\left(
\xi[\psi^{(k)}]_{T,t}^{p+1}-
\xi[\psi^{(k)}]_{T,t}^{p+1,n}\right)
J[\psi^{(k)}]_{T,t}^p\right\}\right|+
\left|{\sf M}\left\{
\xi[\psi^{(k)}]_{T,t}^{p+1,n}J[\psi^{(k)}]_{T,t}^p\right\}\right|=
$$

\vspace{3mm}
$$
=\left|{\sf M}\left\{\left(
J[\psi^{(k)}]_{T,t}-
J[\psi^{(k)}]_{T,t}^{n}\right)
J[\psi^{(k)}]_{T,t}^p\right\}\right|\le
$$

\vspace{7mm}
$$
\le \sqrt{{\sf M}\left\{\left(
J[\psi^{(k)}]_{T,t}-
J[\psi^{(k)}]_{T,t}^{n}\right)^2\right\}}
\sqrt{{\sf M}\left\{\left(
J[\psi^{(k)}]_{T,t}^p\right)^2\right\}}\le
$$

\newpage
\noindent
$$
\le \sqrt{{\sf M}\left\{\left(
J[\psi^{(k)}]_{T,t}-
J[\psi^{(k)}]_{T,t}^{n}\right)^2\right\}}\times
$$

\vspace{1mm}
$$
\times
\left(\sqrt{{\sf M}\left\{\left(
J[\psi^{(k)}]_{T,t}^p - J[\psi^{(k)}]_{T,t}\right)^2\right\}}
+ \sqrt{{\sf M}\left\{\left(
J[\psi^{(k)}]_{T,t}\right)^2\right\}}\right)\le
$$

\begin{equation}
\label{rrre}
~~~~~~~~\le K \sqrt{{\sf M}\left\{\left(
J[\psi^{(k)}]_{T,t}-
J[\psi^{(k)}]_{T,t}^{n}\right)^2\right\}} \to 0\ \ \ {\rm if}\ \ \ n\to\infty,
\end{equation}

\vspace{2mm}
\noindent
where $K$ is a constant.

From (\ref{rrre}) it follows that

\vspace{-2mm}
$$
{\sf M}\left\{
\xi[\psi^{(k)}]_{T,t}^{p+1}J[\psi^{(k)}]_{T,t}^p\right\}=0
$$

\noindent
or

\vspace{-3mm}
$$
{\sf M}\left\{
\left(J[\psi^{(k)}]_{T,t}-J[\psi^{(k)}]_{T,t}^p\right)
J[\psi^{(k)}]_{T,t}^p\right\}=0.
$$

\vspace{2mm}

The last equality means that
\begin{equation}
\label{yyyw}
{\sf M}\left\{
J[\psi^{(k)}]_{T,t}J[\psi^{(k)}]_{T,t}^p\right\}=
{\sf M}\left\{
\left(J[\psi^{(k)}]_{T,t}^p\right)^2\right\}.
\end{equation}

\vspace{1mm}

Taking into account (\ref{yyyw}), we obtain

\vspace{-1mm}
$$
{\sf M}\left\{\left(J[\psi^{(k)}]_{T,t}-
J[\psi^{(k)}]_{T,t}^p\right)^2\right\}=
{\sf M}\left\{\left(J[\psi^{(k)}]_{T,t}\right)^2\right\}+
$$

\vspace{-2mm}
$$
+
{\sf M}\left\{\left(J[\psi^{(k)}]_{T,t}^p\right)^2\right\}
-2{\sf M}\left\{J[\psi^{(k)}]_{T,t}J[\psi^{(k)}]_{T,t}^p\right\}=
{\sf M}\left\{\left(J[\psi^{(k)}]_{T,t}\right)^2\right\}-
$$

\vspace{2mm}
$$
-
{\sf M}\left\{J[\psi^{(k)}]_{T,t}J[\psi^{(k)}]_{T,t}^p\right\}=
$$

\vspace{-2mm}
\begin{equation}
\label{tttr}
~~~~~~~ =\int\limits_{[t,T]^k} K^2(t_1,\ldots,t_k)
dt_1\ldots dt_k - 
{\sf M}\left\{J[\psi^{(k)}]_{T,t}J[\psi^{(k)}]_{T,t}^p\right\}.
\end{equation}

\vspace{2mm}

Let us consider the value
$$
{\sf M}\left\{J[\psi^{(k)}]_{T,t}J[\psi^{(k)}]_{T,t}^p\right\}.
$$

\vspace{2mm}

The relations (\ref{yeee2}) and (\ref{ttt2}) imply that
\begin{equation}
\label{z9}
J[\psi^{(k)}]_{T,t}^p=
\sum_{j_1=0}^{p}\ldots\sum_{j_k=0}^{p}
C_{j_k\ldots j_1}
\sum\limits_{(j_1,\ldots,j_k)}
\int\limits_t^T \phi_{j_k}(t_k)
\ldots
\int\limits_t^{t_{2}}\phi_{j_{1}}(t_{1})
d{\bf w}_{t_1}^{(i_1)}\ldots
d{\bf w}_{t_k}^{(i_k)}.
\end{equation}

\vspace{3mm}

After substituting (\ref{z9}) into (\ref{tttr}), we finally get

\vspace{-3mm}
$$
{\sf M}\left\{\left(J[\psi^{(k)}]_{T,t}-
J[\psi^{(k)}]_{T,t}^p\right)^2\right\}
= \int\limits_{[t,T]^k} K^2(t_1,\ldots,t_k)
dt_1\ldots dt_k - 
$$
$$
- \hspace{-0.5mm}\sum_{j_1=0}^{p}\ldots\sum_{j_k=0}^{p}
C_{j_k\ldots j_1}
{\sf M}\left\{\hspace{-0.5mm}J[\psi^{(k)}]_{T,t}
\sum\limits_{(j_1,\ldots,j_k)}
\int\limits_t^T \phi_{j_k}(t_k)
\ldots
\int\limits_t^{t_{2}}\phi_{j_{1}}(t_{1})
d{\bf w}_{t_1}^{(i_1)}\ldots
d{\bf w}_{t_k}^{(i_k)}\hspace{-0.5mm}\right\}\hspace{-0.6mm}.
$$

\vspace{5mm}

Theorem 1.3 is proved.

\subsection{Exact Calculation of the Mean-Square 
Approximation Errors
for the Cases $k=1,\ldots,5$}

Let us denote

\vspace{-2mm}
$$
{\sf M}\left\{\left(J[\psi^{(k)}]_{T,t}-
J[\psi^{(k)}]_{T,t}^{p}\right)^2\right\}\stackrel{{\rm def}}
{=}E_k^{p},
$$

\vspace{1mm}
$$
\left\Vert K\right\Vert_{L_2([t,T]^k)}^2=
\int\limits_{[t,T]^k} K^2(t_1,\ldots,t_k)
dt_1\ldots dt_k\stackrel{{\rm def}}{=}I_k.
$$

\vspace{7mm}

\centerline{\bf The case $k=1$}

\vspace{2mm}

In this case from Theorem 1.3 we obtain

\vspace{-2mm}
$$
E_1^p
=I_1
-\sum_{j_1=0}^p
C_{j_1}^2.
$$

\vspace{7mm}

\centerline{\bf The case $k=2$}

\vspace{1mm}

In this case from Theorem 1.3 we have

(I).\ $i_1\ne i_2$:
\begin{equation}
\label{ee1}
E_2^p
=I_2
-\sum_{j_1,j_2=0}^p
C_{j_2j_1}^2,
\end{equation}

\vspace{3mm}

(II).\ $i_1=i_2:$
\begin{equation}
\label{uyes2}
E_2^p
=I_2
-\sum_{j_1,j_2=0}^p
C_{j_2j_1}^2-
\sum_{j_1,j_2=0}^p
C_{j_2j_1}C_{j_1j_2}.
\end{equation}

\vspace{4mm}

Note that from (\ref{yeee2}), (\ref{ttt2}),
(\ref{gt22}), (\ref{yyyw}), and (\ref{tttr}) we obtain

\vspace{-1mm}
\begin{equation}
\label{ee1uu1}
E_2^p
=I_2
-\sum_{j_1,j_2=0}^p
C_{j_2j_1}^2-{\bf 1}_{\{i_1=i_2\}}\sum_{j_1,j_2=0}^p
C_{j_2j_1}C_{j_1j_2}.
\end{equation}

\vspace{4mm}
\noindent
Obviously, the relation (\ref{ee1uu1}) 
is consistent with (\ref{ee1}) and (\ref{uyes2}).

{\bf Example 1.1.} Let us consider the following
iterated It\^{o} stochastic integral

\vspace{-1mm}
\begin{equation}
\label{k1000}
I_{(00)T,t}^{(i_1i_2)}
=\int\limits_t^T\int\limits_t^{t_{2}}
d{\bf w}_{t_1}^{(i_1)}
d{\bf w}_{t_2}^{(i_2)},
\end{equation}

\vspace{2mm}
\noindent
where $i_1, i_2=1,\ldots,m.$

Approximation of the iterated It\^{o} stochastic integral (\ref{k1000}) 
based on the expansion (\ref{tyyy}) (Theorem 1.1,
the case of Legendre poly\-no\-mi\-als)
has the following form 

\vspace{-4mm}
\begin{equation}
\label{4004}
I_{(00)T,t}^{(i_1 i_2)p}=
\frac{T-t}{2}\left(\zeta_0^{(i_1)}\zeta_0^{(i_2)}+\sum_{i=1}^{p}
\frac{1}{\sqrt{4i^2-1}}\left(
\zeta_{i-1}^{(i_1)}\zeta_{i}^{(i_2)}-
\zeta_i^{(i_1)}\zeta_{i-1}^{(i_2)}\right)-{\bf 1}_{\{i_1=i_2\}}\right).
\end{equation}

\vspace{3mm}

Note that (\ref{4004}) has been derived for the first time
in \cite{old-art-1} (1997) (also see \cite{old-art-2}-\cite{old-art-4})
with using the another approach. This approach will be considered 
in Sect.~2.4.
Later (\ref{4004}) was obtained \cite{1} (2006), \cite{2}-\cite{new-new-6}
on the base of Theorem 1.1.

Using (\ref{ee1}), we get
\begin{equation}
\label{4007}
~~~~~~~~ {\sf M}\biggl\{\left(I_{(00)T,t}^{(i_1 i_2)}-
I_{(00)T,t}^{(i_1 i_2)p}
\right)^2\biggr\}
=\frac{(T-t)^2}{2}\left(\frac{1}{2}-\sum_{i=1}^p
\frac{1}{4i^2-1}\right),
\end{equation}

\vspace{1mm}
\noindent
where $i_1\ne i_2$.

It should also be noted that the formula (\ref{4007}) 
has been obtained for the first time in \cite{old-art-1} (1997)
by direct calculation.

\vspace{7mm}

\centerline{\bf The case $k=3$}

\vspace{2mm}

In this case from Theorem 1.3 we obtain

\vspace{2mm}

(I).\ $i_1\ne i_2, i_1\ne i_3, i_2\ne i_3:$
\begin{equation}
\label{ravnou1}
E_3^p=I_3
-\sum_{j_1,j_2,j_3=0}^p C_{j_3j_2j_1}^2,
\end{equation}

\vspace{2mm}

(II).\ $i_1=i_2=i_3:$
\begin{equation}
\label{ravno}
E_3^p = I_3 - \sum_{j_1,j_2,j_3=0}^{p}
C_{j_3j_2j_1}\Biggl(\sum\limits_{(j_1,j_2,j_3)}
C_{j_3j_2j_1}\Biggr),
\end{equation}

\vspace{2mm}

(III).1.\ $i_1=i_2\ne i_3:$
\begin{equation}
\label{qq1}
E_3^p=I_3
-\sum_{j_1,j_2,j_3=0}^p C_{j_3j_2j_1}^2-
\sum_{j_1,j_2,j_3=0}^p C_{j_3j_1j_2}C_{j_3j_2j_1},
\end{equation}

\vspace{2mm}

(III).2.\ $i_1\ne i_2=i_3:$
\begin{equation}
\label{dest30}
E_3^p=I_3-
\sum_{j_1,j_2,j_3=0}^p C_{j_3j_2j_1}^2-
\sum_{j_1,j_2,j_3=0}^p C_{j_2j_3j_1}C_{j_3j_2j_1},
\end{equation}

\vspace{2mm}

(III).3.\ $i_1=i_3\ne i_2:$
\begin{equation}
\label{dest40}
E_3^p=I_3
-\sum_{j_1,j_2,j_3=0}^p C_{j_3j_2j_1}^2-
\sum_{j_1,j_2,j_3=0}^p C_{j_3j_2j_1}C_{j_1j_2j_3}.
\end{equation}

\vspace{3mm}

It is not difficult to see that from (\ref{yeee2}), (\ref{ttt2}),
(\ref{gt23}), (\ref{yyyw}), and (\ref{tttr}) we obtain

\vspace{-5mm}
$$
E_3^p=I_3
-\sum_{j_1,j_2,j_3=0}^p C_{j_3j_2j_1}^2-
$$

\vspace{-1mm}
$$
-{\bf 1}_{\{i_1=i_2\}}
\sum_{j_1,j_2,j_3=0}^p C_{j_3j_2j_1}C_{j_3j_1j_2}-
$$

\vspace{-1mm}
$$
-{\bf 1}_{\{i_2=i_3\}}
\sum_{j_1,j_2,j_3=0}^p C_{j_3j_2j_1}C_{j_2j_3j_1}-
$$

\vspace{-1mm}
$$
-{\bf 1}_{\{i_1=i_3\}}
\sum_{j_1,j_2,j_3=0}^p C_{j_3j_2j_1}C_{j_1j_2j_3}-
$$

\vspace{-1mm}
\begin{equation}
\label{dest40xxx}
-{\bf 1}_{\{i_1=i_2=i_3\}}
\sum_{j_1,j_2,j_3=0}^p C_{j_3j_2j_1}\left(C_{j_2j_1j_3}+C_{j_1j_3j_2}\right). 
\end{equation}

\vspace{5mm}
\noindent
Obviously, the relation (\ref{dest40xxx}) 
is consistent with (\ref{ravnou1})--(\ref{dest40}).

Note that the cases $k=2$ and $k=3$ (excepting the formula (\ref{ravno}))
were investigated for the first time in \cite{2} (2007)
using the direct calculation.

{\bf Example 1.2.} Let us consider the following
iterated It\^{o} stochastic integral 

\vspace{-2mm}
\begin{equation}
\label{k1001}
I_{(000)T,t}^{(i_1i_2i_3)}
=\int\limits_t^T\int\limits_t^{t_{3}}\int\limits_t^{t_{2}}
d{\bf w}_{t_1}^{(i_1)}
d{\bf w}_{t_2}^{(i_2)}
d{\bf w}_{t_3}^{(i_3)},
\end{equation}

\noindent
where $i_1, i_2, i_3=1,\ldots,m.$

Approximation of the iterated It\^{o} stochastic integral
(\ref{k1001}) based on Theorem 1.1 (the case
of Legendre polynomials and $p_1=p_2=p_3=p$) has the following form 
\cite{1} (2006), \cite{2}-\cite{new-new-6}
$$
I_{(000)T,t}^{(i_1i_2i_3)p}
=\sum_{j_1,j_2,j_3=0}^{p}
C_{j_3j_2j_1}\Biggl(
\zeta_{j_1}^{(i_1)}\zeta_{j_2}^{(i_2)}\zeta_{j_3}^{(i_3)}
-{\bf 1}_{\{i_1=i_2\}}
{\bf 1}_{\{j_1=j_2\}}
\zeta_{j_3}^{(i_3)}-
\Biggr.
$$
\begin{equation}
\label{sad001}
\Biggl.
-{\bf 1}_{\{i_2=i_3\}}
{\bf 1}_{\{j_2=j_3\}}
\zeta_{j_1}^{(i_1)}-
{\bf 1}_{\{i_1=i_3\}}
{\bf 1}_{\{j_1=j_3\}}
\zeta_{j_2}^{(i_2)}\Biggr),
\end{equation}

\vspace{2mm}
\noindent
where
\begin{equation}
\label{w1}
~~~~~~~~~ C_{j_3j_2j_1}
=\frac{\sqrt{(2j_1+1)(2j_2+1)(2j_3+1)}}{8}(T-t)^{3/2}\bar
C_{j_3j_2j_1},
\end{equation}

\vspace{-1mm}
$$
\bar C_{j_3j_2j_1}=\int\limits_{-1}^{1}P_{j_3}(z)
\int\limits_{-1}^{z}P_{j_2}(y)
\int\limits_{-1}^{y}
P_{j_1}(x)dx dy dz,
$$

\vspace{3mm}
\noindent
where $P_i(x)$ is the Legendre polynomial $(i= 0, 1, 2,\ldots ).$

For example, using (\ref{qq1}) and (\ref{dest30}), we obtain

\vspace{-2mm}
$$
{\sf M}\left\{\left(
I_{(000)T,t}^{(i_1i_2 i_3)}-
I_{(000)T,t}^{(i_1i_2 i_3)p}\right)^2\right\}=
\frac{(T-t)^{3}}{6}-\sum_{j_1,j_2,j_3=0}^{p}
C_{j_3j_2j_1}^2
-\sum_{j_1,j_2,J_3=0}^{p}
C_{j_3j_1j_2}C_{j_3j_2j_1},
$$

\vspace{1mm}
\noindent
where $i_1=i_2\ne i_3$,

\vspace{1mm}

$$
{\sf M}\left\{\left(
I_{(000)T,t}^{(i_1i_2 i_3)}-
I_{(000)T,t}^{(i_1i_2 i_3)p}\right)^2\right\}=
\frac{(T-t)^{3}}{6}-\sum_{j_1,j_2,j_3=0}^{p}
C_{j_3j_2j_1}^2
-\sum_{j_1,j_2,J_3=0}^{p}
C_{j_2j_3j_1}C_{j_3j_2j_1},
$$

\vspace{1mm}
\noindent
where $i_1\ne i_2=i_3$.

The exact 
values of Fourier--Legendre coefficients 
$\bar C_{j_3j_2j_1}$ 
can be calculated for example using computer
algebra system Derive
\cite{1}-\cite{12aa}, \cite{arxiv-4} (see Sect.~5.1, Tables 5.4--5.36).
For more details on calculating of $\bar C_{j_3j_2j_1}$ using Python
programming lan\-gu\-age see
\cite{Kuz-Kuz}, \cite{Mikh-1}.

For the case $i_1=i_2=i_3$ it is convenient to use 
the following well known formula 

\vspace{-2mm}
$$
I_{(000)T,t}^{(i_1 i_1 i_1)}=\frac{1}{6}(T-t)^{3/2}\left(
\left(\zeta_0^{(i_1)}\right)^3-3
\zeta_0^{(i_1)}\right)\ \ \ \hbox{w.~p.~1.}
$$

\vspace{8mm}

\centerline{\bf The case $k=4$}

\vspace{2mm}

In this case from Theorem 1.3 we have

\vspace{2mm}

(I).\ $i_1,\ldots,i_4$ are pairwise different:
$$
E_4^p= I_4 - \sum_{j_1,\ldots,j_4=0}^{p}C_{j_4\ldots j_1}^2,
$$

\newpage
\noindent
\par
(II).\ $i_1=i_2=i_3=i_4$:
$$
E_4^p = I_4 -
 \sum_{j_1,\ldots,j_4=0}^{p}
C_{j_4\ldots j_1}\Biggl(\sum\limits_{(j_1,\ldots,j_4)}
C_{j_4\ldots j_1}\Biggr),
$$

\vspace{4mm}

(III).1.\ $i_1=i_2\ne i_3, i_4;\ i_3\ne i_4:$
\begin{equation}
\label{usl1}
E^p_4 = I_4 - \sum_{j_1,\ldots,j_4=0}^{p}
C_{j_4\ldots j_1}\Biggl(\sum\limits_{(j_1,j_2)}
C_{j_4\ldots j_1}\Biggr),
\end{equation}

\vspace{4mm}

(III).2.\ $i_1=i_3\ne i_2, i_4;\ i_2\ne i_4:$
\begin{equation}
\label{usl2}
E^p_4 = I_4 - \sum_{j_1,\ldots,j_4=0}^{p}
C_{j_4\ldots j_1}\Biggl(\sum\limits_{(j_1,j_3)}
C_{j_4\ldots j_1}\Biggr),
\end{equation}

\vspace{4mm}

(III).3.\ $i_1=i_4\ne i_2, i_3;\ i_2\ne i_3:$
\begin{equation}
\label{usl3}
E^p_4 = I_4 - \sum_{j_1,\ldots,j_4=0}^{p}
C_{j_4\ldots j_1}\Biggl(\sum\limits_{(j_1,j_4)}
C_{j_4\ldots j_1}\Biggr),
\end{equation}

\vspace{4mm}

(III).4.\ $i_2=i_3\ne i_1, i_4;\ i_1\ne i_4:$
\begin{equation}
\label{usl4}
E^p_4 = I_4 - \sum_{j_1,\ldots,j_4=0}^{p}
C_{j_4\ldots j_1}\Biggl(\sum\limits_{(j_2,j_3)}
C_{j_4\ldots j_1}\Biggr),
\end{equation}

\vspace{4mm}

(III).5.\ $i_2=i_4\ne i_1, i_3;\ i_1\ne i_3:$
\begin{equation}
\label{usl5}
E^p_4 = I_4 - \sum_{j_1,\ldots,j_4=0}^{p}
C_{j_4\ldots j_1}\Biggl(\sum\limits_{(j_2,j_4)}
C_{j_4\ldots j_1}\Biggr),
\end{equation}

\vspace{4mm}

(III).6.\ $i_3=i_4\ne i_1, i_2;\ i_1\ne i_2:$
\begin{equation}
\label{usl6}
E^p_4 = I_4 - \sum_{j_1,\ldots,j_4=0}^{p}
C_{j_4\ldots j_1}\Biggl(\sum\limits_{(j_3,j_4)}
C_{j_4\ldots j_1}\Biggr),
\end{equation}

\vspace{4mm}

(IV).1.\ $i_1=i_2=i_3\ne i_4$:
\begin{equation}
\label{usl7}
E_4^p = I_4 -
\sum_{j_1,\ldots,j_4=0}^{p}
C_{j_4\ldots j_1}\Biggl(\sum\limits_{(j_1,j_2,j_3)}
C_{j_4\ldots j_1}\Biggr),
\end{equation}

\vspace{4mm}

(IV).2.\ $i_2=i_3=i_4\ne i_1$:
\begin{equation}
\label{usl8}
E_4^p = I_4 -
 \sum_{j_1,\ldots,j_4=0}^{p}
C_{j_4\ldots j_1}\Biggl(\sum\limits_{(j_2,j_3,j_4)}
C_{j_4\ldots j_1}\Biggr),
\end{equation}

\vspace{4mm}

(IV).3.\ $i_1=i_2=i_4\ne i_3$:
\begin{equation}
\label{usl9}
E_4^p = I_4 -
 \sum_{j_1,\ldots,j_4=0}^{p}
C_{j_4\ldots j_1}\Biggl(\sum\limits_{(j_1,j_2,j_4)}
C_{j_4\ldots j_1}\Biggr),
\end{equation}

\vspace{4mm}

(IV).4.\ $i_1=i_3=i_4\ne i_2$:
\begin{equation}
\label{usl10}
E_4^p = I_4 -
 \sum_{j_1,\ldots,j_4=0}^{p}
C_{j_4\ldots j_1}\Biggl(\sum\limits_{(j_1,j_3,j_4)}
C_{j_4\ldots j_1}\Biggr),
\end{equation}

\vspace{4mm}

(V).1.\ $i_1=i_2\ne i_3=i_4$:
\begin{equation}
\label{usl11}
E^p_4 = I_4 - \sum_{j_1,\ldots,j_4=0}^{p}
C_{j_4\ldots j_1}\Biggl(\sum\limits_{(j_1,j_2)}\Biggl(
\sum\limits_{(j_3,j_4)}
C_{j_4\ldots j_1}\Biggr)\Biggr),
\end{equation}

\vspace{4mm}

(V).2.\ $i_1=i_3\ne i_2=i_4$:
\begin{equation}
\label{usl12}
E^p_4 = I_4 - \sum_{j_1,\ldots,j_4=0}^{p}
C_{j_4\ldots j_1}\Biggl(\sum\limits_{(j_1,j_3)}\Biggl(
\sum\limits_{(j_2,j_4)}
C_{j_4\ldots j_1}\Biggr)\Biggr),
\end{equation}

\vspace{4mm}

(V).3.\ $i_1=i_4\ne i_2=i_3$:
\begin{equation}
\label{usl13}
E^p_4 = I_4 - \sum_{j_1,\ldots,j_4=0}^{p}
C_{j_4\ldots j_1}\Biggl(\sum\limits_{(j_1,j_4)}\Biggl(
\sum\limits_{(j_2,j_3)}
C_{j_4\ldots j_1}\Biggr)\Biggr).
\end{equation}

\newpage
\noindent
\centerline{\bf The case $k=5$}

\vspace{1mm}

In this case from Theorem 1.3 we obtain

(I).\ $i_1,\ldots,i_5$ are pairwise different:
$$
E_5^p= I_5 - \sum_{j_1,\ldots,j_5=0}^{p}C_{j_5\ldots j_1}^2,
$$

(II).\ $i_1=i_2=i_3=i_4=i_5$:
$$
E_5^p = I_5 - \sum_{j_1,\ldots,j_5=0}^{p}
C_{j_5\ldots j_1}\Biggl(\sum\limits_{(j_1,\ldots,j_5)}
C_{j_5\ldots j_1}\Biggr),
$$

\vspace{3mm}

(III).1.\ $i_1=i_2\ne i_3,i_4,i_5$\ ($i_3,i_4,i_5$ are pairwise different):
$$
E^p_5 = I_5 - \sum_{j_1,\ldots,j_5=0}^{p}
C_{j_5\ldots j_1}\Biggl(\sum\limits_{(j_1,j_2)}
C_{j_5\ldots j_1}\Biggr),
$$

\vspace{3mm}

(III).2.\ $i_1=i_3\ne i_2,i_4,i_5$\ ($i_2,i_4,i_5$ are pairwise different):
$$
E^p_5 = I_5 - \sum_{j_1,\ldots,j_5=0}^{p}
C_{j_5\ldots j_1}\Biggl(\sum\limits_{(j_1,j_3)}
C_{j_5\ldots j_1}\Biggr),
$$

\vspace{3mm}

(III).3.\ $i_1=i_4\ne i_2,i_3,i_5$\ ($i_2,i_3,i_5$ are pairwise different):
$$
E^p_5 = I_5 - \sum_{j_1,\ldots,j_5=0}^{p}
C_{j_5\ldots j_1}\Biggl(\sum\limits_{(j_1,j_4)}
C_{j_5\ldots j_1}\Biggr),
$$

\vspace{3mm}

(III).4.\ $i_1=i_5\ne i_2,i_3,i_4$\ ($i_2,i_3,i_4$  are pairwise different):
$$
E^p_5 = I_5 - \sum_{j_1,\ldots,j_5=0}^{p}
C_{j_5\ldots j_1}\Biggl(\sum\limits_{(j_1,j_5)}
C_{j_5\ldots j_1}\Biggr),
$$

\vspace{3mm}

(III).5.\ $i_2=i_3\ne i_1,i_4,i_5$\ ($i_1,i_4,i_5$ are pairwise different):
$$
E^p_5 = I_5 - \sum_{j_1,\ldots,j_5=0}^{p}
C_{j_5\ldots j_1}\Biggl(\sum\limits_{(j_2,j_3)}
C_{j_5\ldots j_1}\Biggr),
$$

\vspace{3mm}

(III).6.\ $i_2=i_4\ne i_1,i_3,i_5$\ ($i_1,i_3,i_5$ are pairwise different):
$$
E^p_5 = I_5 - \sum_{j_1,\ldots,j_5=0}^{p}
C_{j_5\ldots j_1}\Biggl(\sum\limits_{(j_2,j_4)}
C_{j_5\ldots j_1}\Biggr),
$$

\vspace{3mm}

(III).7.\ $i_2=i_5\ne i_1,i_3,i_4$\ ($i_1,i_3,i_4$ are pairwise different):
$$
E^p_5 = I_5 - \sum_{j_1,\ldots,j_5=0}^{p}
C_{j_5\ldots j_1}\Biggl(\sum\limits_{(j_2,j_5)}
C_{j_5\ldots j_1}\Biggr),
$$

\vspace{3mm}

(III).8.\ $i_3=i_4\ne i_1,i_2,i_5$\ ($i_1,i_2,i_5$  are pairwise different):
$$
E^p_5 = I_5 - \sum_{j_1,\ldots,j_5=0}^{p}
C_{j_5\ldots j_1}\Biggl(\sum\limits_{(j_3,j_4)}
C_{j_5\ldots j_1}\Biggr),
$$

\vspace{3mm}

(III).9.\ $i_3=i_5\ne i_1,i_2,i_4$\ ($i_1,i_2,i_4$  are pairwise different):
$$
E^p_5 = I_5 - \sum_{j_1,\ldots,j_5=0}^{p}
C_{j_5\ldots j_1}\Biggl(\sum\limits_{(j_3,j_5)}
C_{j_5\ldots j_1}\Biggr),
$$

\vspace{3mm}

(III).10.\ $i_4=i_5\ne i_1,i_2,i_3$\ ($i_1,i_2,i_3$  are pairwise different):
$$
E^p_5 = I_5 - \sum_{j_1,\ldots,j_5=0}^{p}
C_{j_5\ldots j_1}\Biggl(\sum\limits_{(j_4,j_5)}
C_{j_5\ldots j_1}\Biggr),
$$

\vspace{3mm}

(IV).1.\ $i_1=i_2=i_3\ne i_4, i_5$\ $(i_4\ne i_5$):
$$
E^p_5 = I_5 - \sum_{j_1,\ldots,j_5=0}^{p}
C_{j_5\ldots j_1}\Biggl(\sum\limits_{(j_1,j_2,j_3)}
C_{j_5\ldots j_1}\Biggr),
$$

\vspace{3mm}

(IV).2.\ $i_1=i_2=i_4\ne i_3, i_5$\  $(i_3\ne i_5$):
$$
E^p_5 = I_5 - \sum_{j_1,\ldots,j_5=0}^{p}
C_{j_5\ldots j_1}\Biggl(\sum\limits_{(j_1,j_2,j_4)}
C_{j_5\ldots j_1}\Biggr),
$$

\vspace{3mm}

(IV).3.\ $i_1=i_2=i_5\ne i_3, i_4$\  $(i_3\ne i_4$):
$$
E^p_5 = I_5 - \sum_{j_1,\ldots,j_5=0}^{p}
C_{j_5\ldots j_1}\Biggl(\sum\limits_{(j_1,j_2,j_5)}
C_{j_5\ldots j_1}\Biggr),
$$

\vspace{3mm}

(IV).4.\ $i_2=i_3=i_4\ne i_1, i_5$\  $(i_1\ne i_5$):
$$
E^p_5 = I_5 - \sum_{j_1,\ldots,j_5=0}^{p}
C_{j_5\ldots j_1}\Biggl(\sum\limits_{(j_2,j_3,j_4)}
C_{j_5\ldots j_1}\Biggr),
$$

\vspace{3mm}

(IV).5.\ $i_2=i_3=i_5\ne i_1, i_4$\  $(i_1\ne i_4$):
$$
E^p_5 = I_5 - \sum_{j_1,\ldots,j_5=0}^{p}
C_{j_5\ldots j_1}\Biggl(\sum\limits_{(j_2,j_3,j_5)}
C_{j_5\ldots j_1}\Biggr),
$$

\vspace{3mm}

(IV).6.\ $i_2=i_4=i_5\ne i_1, i_3$\  $(i_1\ne i_3$):
$$
E^p_5 = I_5 - \sum_{j_1,\ldots,j_5=0}^{p}
C_{j_5\ldots j_1}\Biggl(\sum\limits_{(j_2,j_4,j_5)}
C_{j_5\ldots j_1}\Biggr),
$$

\vspace{3mm}

(IV).7.\ $i_3=i_4=i_5\ne i_1, i_2$\  $(i_1\ne i_2$):
$$
E^p_5 = I_5 - \sum_{j_1,\ldots,j_5=0}^{p}
C_{j_5\ldots j_1}\Biggl(\sum\limits_{(j_3,j_4,j_5)}
C_{j_5\ldots j_1}\Biggr),
$$

\vspace{3mm}

(IV).8.\ $i_1=i_3=i_5\ne i_2, i_4$\  $(i_2\ne i_4$):
$$
E^p_5 = I_5 - \sum_{j_1,\ldots,j_5=0}^{p}
C_{j_5\ldots j_1}\Biggl(\sum\limits_{(j_1,j_3,j_5)}
C_{j_5\ldots j_1}\Biggr),
$$

\vspace{3mm}

(IV).9.\ $i_1=i_3=i_4\ne i_2, i_5$\  $(i_2\ne i_5$):
$$
E^p_5 = I_5 - \sum_{j_1,\ldots,j_5=0}^{p}
C_{j_5\ldots j_1}\Biggl(\sum\limits_{(j_1,j_3,j_4)}
C_{j_5\ldots j_1}\Biggr),
$$

\vspace{3mm}

(IV).10.\ $i_1=i_4=i_5\ne i_2, i_3$\  $(i_2\ne i_3$):
$$
E^p_5 = I_5 - \sum_{j_1,\ldots,j_5=0}^{p}
C_{j_5\ldots j_1}\Biggl(\sum\limits_{(j_1,j_4,j_5)}
C_{j_5\ldots j_1}\Biggr),
$$

\vspace{3mm}

(V).1.\ $i_1=i_2=i_3=i_4\ne i_5$:
$$
E^p_5 = I_5 - \sum_{j_1,\ldots,j_5=0}^{p}
C_{j_5\ldots j_1}\Biggl(\sum\limits_{(j_1,j_2,j_3,j_4)}
C_{j_5\ldots j_1}\Biggr),
$$

\vspace{3mm}

(V).2.\ $i_1=i_2=i_3=i_5\ne i_4$:
$$
E^p_5 = I_5 - \sum_{j_1,\ldots,j_5=0}^{p}
C_{j_5\ldots j_1}\Biggl(\sum\limits_{(j_1,j_2,j_3,j_5)}
C_{j_5\ldots j_1}\Biggr),
$$

\vspace{3mm}

(V).3.\ $i_1=i_2=i_4=i_5\ne i_3$:
$$
E^p_5 = I_5 - \sum_{j_1,\ldots,j_5=0}^{p}
C_{j_5\ldots j_1}\Biggl(\sum\limits_{(j_1,j_2,j_4,j_5)}
C_{j_5\ldots j_1}\Biggr),
$$

\vspace{3mm}

(V).4.\ $i_1=i_3=i_4=i_5\ne i_2$:
$$
E^p_5 = I_5 - \sum_{j_1,\ldots,j_5=0}^{p}
C_{j_5\ldots j_1}\Biggl(\sum\limits_{(j_1,j_3,j_4,j_5)}
C_{j_5\ldots j_1}\Biggr),
$$

\vspace{3mm}

(V).5.\ $i_2=i_3=i_4=i_5\ne i_1$:
$$
E^p_5 = I_5 - \sum_{j_1,\ldots,j_5=0}^{p}
C_{j_5\ldots j_1}\Biggl(\sum\limits_{(j_2,j_3,j_4,j_5)}
C_{j_5\ldots j_1}\Biggr),
$$

\vspace{3mm}

(VI).1.\ $i_5\ne i_1=i_2\ne i_3=i_4\ne i_5$:
$$
E^p_5 = I_5 - \sum_{j_1,\ldots,j_5=0}^{p}
C_{j_5\ldots j_1}\Biggl(\sum\limits_{(j_1,j_2)}\Biggl(
\sum\limits_{(j_3,j_4)}
C_{j_5\ldots j_1}\Biggr)\Biggr),
$$

\vspace{3mm}

(VI).2.\ $i_5\ne i_1=i_3\ne i_2=i_4\ne i_5$:
$$
E^p_5 = I_5 - \sum_{j_1,\ldots,j_5=0}^{p}
C_{j_5\ldots j_1}\Biggl(\sum\limits_{(j_1,j_3)}\Biggl(
\sum\limits_{(j_2,j_4)}
C_{j_5\ldots j_1}\Biggr)\Biggr),
$$

\vspace{3mm}

(VI).3.\ $i_5\ne i_1=i_4\ne i_2=i_3\ne i_5$:
$$
E^p_5 = I_5 - \sum_{j_1,\ldots,j_5=0}^{p}
C_{j_5\ldots j_1}\Biggl(\sum\limits_{(j_1,j_4)}\Biggl(
\sum\limits_{(j_2,j_3)}
C_{j_5\ldots j_1}\Biggr)\Biggr),
$$

\vspace{3mm}

(VI).4.\ $i_4\ne i_1=i_2\ne i_3=i_5\ne i_4$:
$$
E^p_5 = I_5 - \sum_{j_1,\ldots,j_5=0}^{p}
C_{j_5\ldots j_1}\Biggl(\sum\limits_{(j_1,j_2)}\Biggl(
\sum\limits_{(j_3,j_5)}
C_{j_5\ldots j_1}\Biggr)\Biggr),
$$

\vspace{3mm}

(VI).5.\ $i_4\ne i_1=i_5\ne i_2=i_3\ne i_4$:
$$
E^p_5 = I_5 - \sum_{j_1,\ldots,j_5=0}^{p}
C_{j_5\ldots j_1}\Biggl(\sum\limits_{(j_1,j_5)}\Biggl(
\sum\limits_{(j_2,j_3)}
C_{j_5\ldots j_1}\Biggr)\Biggr),
$$

\vspace{3mm}

(VI).6.\ $i_4\ne i_2=i_5\ne i_1=i_3\ne i_4$:
$$
E^p_5 = I_5 - \sum_{j_1,\ldots,j_5=0}^{p}
C_{j_5\ldots j_1}\Biggl(\sum\limits_{(j_2,j_5)}\Biggl(
\sum\limits_{(j_1,j_3)}
C_{j_5\ldots j_1}\Biggr)\Biggr),
$$

\vspace{3mm}

(VI).7.\ $i_3\ne i_2=i_5\ne i_1=i_4\ne i_3$:
$$
E^p_5 = I_5 - \sum_{j_1,\ldots,j_5=0}^{p}
C_{j_5\ldots j_1}\Biggl(\sum\limits_{(j_2,j_5)}\Biggl(
\sum\limits_{(j_1,j_4)}
C_{j_5\ldots j_1}\Biggr)\Biggr),
$$

\vspace{3mm}

(VI).8.\ $i_3\ne i_1=i_2\ne i_4=i_5\ne i_3$:
$$
E^p_5 = I_5 - \sum_{j_1,\ldots,j_5=0}^{p}
C_{j_5\ldots j_1}\Biggl(\sum\limits_{(j_1,j_2)}\Biggl(
\sum\limits_{(j_4,j_5)}
C_{j_5\ldots j_1}\Biggr)\Biggr),
$$

\vspace{3mm}

(VI).9.\ $i_3\ne i_2=i_4\ne i_1=i_5\ne i_3$:
$$
E^p_5 = I_5 - \sum_{j_1,\ldots,j_5=0}^{p}
C_{j_5\ldots j_1}\Biggl(\sum\limits_{(j_2,j_4)}\Biggl(
\sum\limits_{(j_1,j_5)}
C_{j_5\ldots j_1}\Biggr)\Biggr),
$$

\vspace{3mm}

(VI).10.\ $i_2\ne i_1=i_4\ne i_3=i_5\ne i_2$:
$$
E^p_5 = I_5 - \sum_{j_1,\ldots,j_5=0}^{p}
C_{j_5\ldots j_1}\Biggl(\sum\limits_{(j_1,j_4)}\Biggl(
\sum\limits_{(j_3,j_5)}
C_{j_5\ldots j_1}\Biggr)\Biggr),
$$

\vspace{3mm}

(VI).11.\ $i_2\ne i_1=i_3\ne i_4=i_5\ne i_2$:
$$
E^p_5 = I_5 - \sum_{j_1,\ldots,j_5=0}^{p}
C_{j_5\ldots j_1}\Biggl(\sum\limits_{(j_1,j_3)}\Biggl(
\sum\limits_{(j_4,j_5)}
C_{j_5\ldots j_1}\Biggr)\Biggr),
$$

\vspace{3mm}

(VI).12.\ $i_2\ne i_1=i_5\ne i_3=i_4\ne i_2$:
$$
E^p_5 = I_5 - \sum_{j_1,\ldots,j_5=0}^{p}
C_{j_5\ldots j_1}\Biggl(\sum\limits_{(j_1,j_5)}\Biggl(
\sum\limits_{(j_3,j_4)}
C_{j_5\ldots j_1}\Biggr)\Biggr),
$$

\vspace{3mm}

(VI).13.\ $i_1\ne i_2=i_3\ne i_4=i_5\ne i_1$:
$$
E^p_5 = I_5 - \sum_{j_1,\ldots,j_5=0}^{p}
C_{j_5\ldots j_1}\Biggl(\sum\limits_{(j_2,j_3)}\Biggl(
\sum\limits_{(j_4,j_5)}
C_{j_5\ldots j_1}\Biggr)\Biggr),
$$

\vspace{3mm}

(VI).14.\ $i_1\ne i_2=i_4\ne i_3=i_5\ne i_1$:
$$
E^p_5 = I_5 - \sum_{j_1,\ldots,j_5=0}^{p}
C_{j_5\ldots j_1}\Biggl(\sum\limits_{(j_2,j_4)}\Biggl(
\sum\limits_{(j_3,j_5)}
C_{j_5\ldots j_1}\Biggr)\Biggr),
$$

\vspace{3mm}

(VI).15.\ $i_1\ne i_2=i_5\ne i_3=i_4\ne i_1$:
$$
E^p_5 = I_5 - \sum_{j_1,\ldots,j_5=0}^{p}
C_{j_5\ldots j_1}\Biggl(\sum\limits_{(j_2,j_5)}\Biggl(
\sum\limits_{(j_3,j_4)}
C_{j_5\ldots j_1}\Biggr)\Biggr),
$$

\vspace{3mm}

(VII).1.\ $i_1=i_2=i_3\ne i_4=i_5$:
$$
E^p_5 = I_5 - \sum_{j_1,\ldots,j_5=0}^{p}
C_{j_5\ldots j_1}\Biggl(\sum\limits_{(j_4,j_5)}\Biggl(
\sum\limits_{(j_1,j_2,j_3)}
C_{j_5\ldots j_1}\Biggr)\Biggr),
$$

\vspace{3mm}

(VII).2.\ $i_1=i_2=i_4\ne i_3=i_5$:
$$
E^p_5 = I_5 - \sum_{j_1,\ldots,j_5=0}^{p}
C_{j_5\ldots j_1}\Biggl(\sum\limits_{(j_3,j_5)}\Biggl(
\sum\limits_{(j_1,j_2,j_4)}
C_{j_5\ldots j_1}\Biggr)\Biggr),
$$

\vspace{3mm}

(VII).3.\ $i_1=i_2=i_5\ne i_3=i_4$:
$$
E_p = I - \sum_{j_1,\ldots,j_5=0}^{p}
C_{j_5\ldots j_1}\Biggl(\sum\limits_{(j_3,j_4)}\Biggl(
\sum\limits_{(j_1,j_2,j_5)}
C_{j_5\ldots j_1}\Biggr)\Biggr),
$$

\vspace{3mm}

(VII).4.\ $i_2=i_3=i_4\ne i_1=i_5$:
$$
E^p_5 = I_5 - \sum_{j_1,\ldots,j_5=0}^{p}
C_{j_5\ldots j_1}\Biggl(\sum\limits_{(j_1,j_5)}\Biggl(
\sum\limits_{(j_2,j_3,j_4)}
C_{j_5\ldots j_1}\Biggr)\Biggr),
$$

\vspace{3mm}

(VII).5.\ $i_2=i_3=i_5\ne i_1=i_4$:
$$
E^p_5 = I_5 - \sum_{j_1,\ldots,j_5=0}^{p}
C_{j_5\ldots j_1}\Biggl(\sum\limits_{(j_1,j_4)}\Biggl(
\sum\limits_{(j_2,j_3,j_5)}
C_{j_5\ldots j_1}\Biggr)\Biggr),
$$

\vspace{3mm}

(VII).6.\ $i_2=i_4=i_5\ne i_1=i_3$:
$$
E^p_5 = I_5 - \sum_{j_1,\ldots,j_5=0}^{p}
C_{j_5\ldots j_1}\Biggl(\sum\limits_{(j_1,j_3)}\Biggl(
\sum\limits_{(j_2,j_4,j_5)}
C_{j_5\ldots j_1}\Biggr)\Biggr),
$$

\vspace{3mm}

(VII).7.\ $i_3=i_4=i_5\ne i_1=i_2$:
$$
E^p_5 = I_5 - \sum_{j_1,\ldots,j_5=0}^{p}
C_{j_5\ldots j_1}\Biggl(\sum\limits_{(j_1,j_2)}\Biggl(
\sum\limits_{(j_3,j_4,j_5)}
C_{j_5\ldots j_1}\Biggr)\Biggr),
$$

\vspace{3mm}

(VII).8.\ $i_1=i_3=i_5\ne i_2=i_4$:
$$
E^p_5 = I_5 - \sum_{j_1,\ldots,j_5=0}^{p}
C_{j_5\ldots j_1}\Biggl(\sum\limits_{(j_2,j_4)}\Biggl(
\sum\limits_{(j_1,j_3,j_5)}
C_{j_5\ldots j_1}\Biggr)\Biggr),
$$

\vspace{3mm}

(VII).9.\ $i_1=i_3=i_4\ne i_2=i_5$:
$$
E^p_5 = I_5 - \sum_{j_1,\ldots,j_5=0}^{p}
C_{j_5\ldots j_1}\Biggl(\sum\limits_{(j_2,j_5)}\Biggl(
\sum\limits_{(j_1,j_3,j_4)}
C_{j_5\ldots j_1}\Biggr)\Biggr),
$$

\vspace{3mm}

(VII).10.\ $i_1=i_4=i_5\ne i_2=i_3$:
$$
E^p_5 = I_5 - \sum_{j_1,\ldots,j_5=0}^{p}
C_{j_5\ldots j_1}\Biggl(\sum\limits_{(j_2,j_3)}\Biggl(
\sum\limits_{(j_1,j_4,j_5)}
C_{j_5\ldots j_1}\Biggr)\Biggr).
$$

\vspace{2mm}

Let us make a remark about Theorem~1.3.
It is easy to see that the right-hand side of the formula (\ref{tttr11})
consists of two parts.
The first part tends to zero when $p\to\infty$ by Parseval's equality.
At the same time the second part also tends to 
zero when $p\to\infty$, but due to the generalized Parseval equality.
Let us explain the above reasoning in more detail for the case $k=3.$

For the case $k=3$ we have (see (\ref{dest40xxx}))

\vspace{-5mm}
$$
E_3^p=I_3
-\sum_{j_1,j_2,j_3=0}^p C_{j_3j_2j_1}^2-
$$

\vspace{-1mm}
$$
-{\bf 1}_{\{i_1=i_2\}}
\sum_{j_1,j_2,j_3=0}^p C_{j_3j_2j_1}C_{j_3j_1j_2}-
$$

\vspace{-1mm}
$$
-{\bf 1}_{\{i_2=i_3\}}
\sum_{j_1,j_2,j_3=0}^p C_{j_3j_2j_1}C_{j_2j_3j_1}-
$$

\vspace{-1mm}
$$
-{\bf 1}_{\{i_1=i_3\}}
\sum_{j_1,j_2,j_3=0}^p C_{j_3j_2j_1}C_{j_1j_2j_3}-
$$

\vspace{-1mm}
\begin{equation}
\label{ten-1000}
-{\bf 1}_{\{i_1=i_2=i_3\}}
\sum_{j_1,j_2,j_3=0}^p C_{j_3j_2j_1}\left(C_{j_2j_1j_3}+C_{j_1j_3j_2}\right). 
\end{equation}

\vspace{3mm}

Applying the Parseval equality, we obtain
\begin{equation}
\label{ten-1001}
\lim\limits_{p\to\infty}\left(I_3
-\sum_{j_1,j_2,j_3=0}^p C_{j_3j_2j_1}^2\right)=0.
\end{equation}

\vspace{2mm}

The generalized Parseval equality gives
\begin{equation}
\label{ten-1002}
~~~~~~\lim\limits_{p\to\infty}\sum_{j_1,j_2,j_3=0}^p C_{j_3j_2j_1}C_{j_1j_2j_3}=0,\ \ \ 
\lim\limits_{p\to\infty}\sum_{j_1,j_2,j_3=0}^p C_{j_3j_2j_1}C_{j_3j_1j_2}=0,
\end{equation}

\vspace{-5mm}
\begin{equation}
\label{ten-1003}
~~~~~~\lim\limits_{p\to\infty}\sum_{j_1,j_2,j_3=0}^p C_{j_3j_2j_1}C_{j_1j_3j_2}=0,\ \ \
\lim\limits_{p\to\infty}\sum_{j_1,j_2,j_3=0}^p C_{j_3j_2j_1}C_{j_2j_1j_3}=0,
\end{equation}

\vspace{-1mm}
\begin{equation}
\label{ten-1004}
\lim\limits_{p\to\infty}\sum_{j_1,j_2,j_3=0}^p C_{j_3j_2j_1}C_{j_2j_3j_1}=0.
\end{equation}

\vspace{3mm}

Let us explain in more detail the first equality in (\ref{ten-1002}).
Using the generalized Parseval equality, we have

\vspace{-3mm}
$$
\lim\limits_{p\to\infty}\sum_{j_1,j_2,j_3=0}^{p}
C_{j_1 j_2 j_3}C_{j_3 j_2 j_1}=
$$

\vspace{-2mm}
$$
=\lim\limits_{p\to\infty}\sum_{j_1,j_2,j_3=0}^{p}
\int\limits_t^T\psi_3(t_3)\phi_{j_1}(t_3)\int\limits_t^{t_3}\psi_2(t_2)\phi_{j_2}(t_2)
\int\limits_t^{t_2}\psi_1(t_1)\phi_{j_3}(t_1)dt_1dt_2dt_3\times
$$

\vspace{-2mm}
$$
\times
\int\limits_t^T\psi_3(t_3)\phi_{j_3}(t_3)\int\limits_t^{t_3}\psi_2(t_2)\phi_{j_2}(t_2)
\int\limits_t^{t_2}\psi_1(t_1)\phi_{j_1}(t_1)dt_1dt_2dt_3=
$$

\vspace{-2mm}
$$
=\lim\limits_{p\to\infty}\sum_{j_1,j_2,j_3=0}^{p}
\int\limits_t^T\psi_1(t_3)\phi_{j_3}(t_3)\int\limits_{t_3}^T\psi_2(t_2)\phi_{j_2}(t_2)
\int\limits_{t_2}^T\psi_3(t_1)\phi_{j_1}(t_1)dt_1dt_2dt_3\times
$$

\vspace{-2mm}
$$
\times
\int\limits_t^T\psi_3(t_3)\phi_{j_3}(t_3)\int\limits_t^{t_3}\psi_2(t_2)\phi_{j_2}(t_2)
\int\limits_t^{t_2}\psi_1(t_1)\phi_{j_1}(t_1)dt_1dt_2dt_3=
$$
$$
=\lim\limits_{p\to\infty}\sum_{j_1,j_2,j_3=0}^{p}~
\int\limits_{[t,T]^3}{\bf 1}_{\{t_3<t_2<t_1\}}
\psi_3(t_1)\psi_2(t_2)\psi_1(t_3)\prod\limits_{l=1}^3
\phi_{j_l}(t_l)dt_1dt_2dt_3 \times
$$
$$
\times
\int\limits_{[t,T]^3}{\bf 1}_{\{t_1<t_2<t_3\}}
\psi_1(t_1)\psi_2(t_2)\psi_3(t_3)\prod\limits_{l=1}^3
\phi_{j_l}(t_l)dt_1dt_2dt_3=
$$

\vspace{-2mm}
\begin{equation}
\label{ten-1005}
=\int\limits_{[t,T]^3}{\bf 1}_{\{t_3<t_2<t_1\}}{\bf 1}_{\{t_1<t_2<t_3\}}
\psi_3(t_1)\psi_1(t_1)\left(\psi_2(t_2)\right)^2\psi_3(t_3)\psi_1(t_3)
dt_1dt_2dt_3=0.
\end{equation}

\vspace{2mm}

Applying (\ref{ten-1000}), (\ref{ten-1001}), and (\ref{ten-1002})--(\ref{ten-1004}), we get
(see (\ref{tttr11}) for the case $k=3$)
$$
\lim\limits_{p\to\infty}E_3^p=0.
$$

\subsection{Estimate for the Mean-Square Approximation Error
of Iterated It\^{o} Stochastic Integrals Based on Theorem 1.1}

In this section, we prove the useful estimate for the
mean-square approximation error in Theorem 1.1.

{\bf Theorem 1.4} {\cite{11}-\cite{12aa}, \cite{arxiv-3}. 
{\it Suppose that
every $\psi_l(\tau)$ $(l=1,\ldots, k)$ is a continuous nonrandom function on 
$[t, T]$ and
$\{\phi_j(x)\}_{j=0}^{\infty}$ is a complete orthonormal system  
of functions in the space $L_2([t,T]),$ each function $\phi_j(x)$ of which 
for finite $j$ satisfies the condition 
$(\star)$ {\rm (}see Sect.~{\rm 1.1.7)}.
Then the estimate
$$
{\sf M}\left\{\left(
J[\psi^{(k)}]_{T,t}-J[\psi^{(k)}]_{T,t}^{p_1,\ldots,p_k}
\right)^2\right\}
\le 
$$
\begin{equation}
\label{z1}
~~~~\le k!\left(~\int\limits_{[t,T]^k}
K^2(t_1,\ldots,t_k)
dt_1\ldots dt_k -\sum_{j_1=0}^{p_1}\ldots
\sum_{j_k=0}^{p_k}C^2_{j_k\ldots j_1}\right)
\end{equation}

\vspace{3mm}
\noindent
is valid for the following cases{\rm :}

{\rm 1.}\ $i_1,\ldots,i_k=1,\ldots,m$\ \ and\ \ $0<T-t<\infty,$

{\rm 2.}\ $i_1,\ldots,i_k=0, 1,\ldots,m,$\ \ $i_1^2+\ldots+i_k^2>0,$\ \
and\ \ $0<T-t<1,$

\noindent
where $J[\psi^{(k)}]_{T,t}$ is the iterated It\^{o} stochastic integral {\rm (\ref{ito}),}
$J[\psi^{(k)}]_{T,t}^{p_1,\ldots,p_k}$ is the 
expression on the right-hand side of {\rm (\ref{tyyy})} before
passing to the limit 
$\hbox{\vtop{\offinterlineskip\halign{
\hfil#\hfil\cr
{\rm l.i.m.}\cr
$\stackrel{}{{}_{p_1,\ldots,p_k\to \infty}}$\cr
}} };$ another 
notations are the same as in Theorem {\rm 1.1}.
}

{\bf Proof.}\ In the proof of Theorem 1.1 
we obtained w.~p.~1 the 
following re\-pre\-sen\-ta\-ti\-on (see (\ref{novoe2}))
$$
J[\psi^{(k)}]_{T,t}=J[\psi^{(k)}]_{T,t}^{p_1,\ldots,p_k}+
R_{T,t}^{p_1,\ldots,p_k},
$$
where $J[\psi^{(k)}]_{T,t}^{p_1,\ldots,p_k}$
is the 
expression on the right-hand side of {\rm (\ref{tyyy})} before
passing to the limit 
$\hbox{\vtop{\offinterlineskip\halign{
\hfil#\hfil\cr
{\rm l.i.m.}\cr
$\stackrel{}{{}_{p_1,\ldots,p_k\to \infty}}$\cr
}} }$ and
$$
R_{T,t}^{p_1,\ldots,p_k}=
\sum_{(t_1,\ldots,t_k)}
\int\limits_{t}^{T}
\ldots
\int\limits_{t}^{t_2}
\Biggl(K(t_1,\ldots,t_k)-
\sum_{j_1=0}^{p_1}\ldots
\sum_{j_k=0}^{p_k}
C_{j_k\ldots j_1}
\prod_{l=1}^k\phi_{j_l}(t_l)\Biggr)\times
$$
\begin{equation}
\label{y007a}
\times d{\bf w}_{t_1}^{(i_1)}
\ldots
d{\bf w}_{t_k}^{(i_k)},
\end{equation}
where
$$
\sum_{(t_1,\ldots,t_k)}
$$
means the sum with respect to all
possible permutations
$(t_1,\ldots,t_k),$ which are
performed only 
in the values $d{\bf w}_{t_1}^{(i_1)}
\ldots $
$d{\bf w}_{t_k}^{(i_k)}$. At the same time the indices near 
upper limits of integration in the iterated stochastic integrals 
are changed correspondently and if $t_r$ swapped with $t_q$ in the  
permutation $(t_1,\ldots,t_k)$, then $i_r$ swapped with $i_q$ in the 
permutation $(i_1,\ldots,i_k)$.

The stochastic integrals on the right-hand side of (\ref{y007a})
will be dependent in a stochastic sense
($i_1,\ldots,i_k=1,\ldots,m,$\ $k\in {\bf N}).$
Let us estimate the second moment of 
$$
J[\psi^{(k)}]_{T,t}-J[\psi^{(k)}]_{T,t}^{p_1,\ldots,p_k}.
$$
 
Using (\ref{99.010}), (\ref{riemann}), (\ref{y007a}), 
the orthonormality
of the system $\{\phi_j(x)\}_{j=0}^{\infty}$
(see the relation (\ref{dobav100})), 
and the elementary inequality
\begin{equation}
\label{y5}
~~~~~~~~~~~ \left(a_1+a_2+\ldots+a_p\right)^2 \le
p\left(a_1^2+a_2^2+\ldots+a_p^2\right),\ \ \ p\in {\bf N},
\end{equation}
we obtain 
the following estimate 
$$
{\sf M}\left\{\left(J[\psi^{(k)}]_{T,t}-J[\psi^{(k)}]_{T,t}^{p_1,\ldots,p_k}
\right)^2\right\}
\le 
$$
$$
\le k!
\sum_{(t_1,\ldots,t_k)}
\int\limits_{t}^{T}
\ldots
\int\limits_{t}^{t_2}
\Biggl(K(t_1,\ldots,t_k)-
\sum_{j_1=0}^{p_1}\ldots
\sum_{j_k=0}^{p_k}
C_{j_k\ldots j_1}
\prod_{l=1}^k\phi_{j_l}(t_l)\Biggr)^2
dt_1\ldots dt_k=
$$
$$
=k!\int\limits_{[t,T]^k}
\Biggl(K(t_1,\ldots,t_k)-
\sum_{j_1=0}^{p_1}\ldots
\sum_{j_k=0}^{p_k}
C_{j_k\ldots j_1}
\prod_{l=1}^k\phi_{j_l}(t_l)\Biggr)^2
dt_1
\ldots
dt_k=
$$
\begin{equation}
\label{star00011}
~~~~~~~= k!\left(~\int\limits_{[t,T]^k}
K^2(t_1,\ldots,t_k)
dt_1\ldots dt_k -\sum_{j_1=0}^{p_1}\ldots
\sum_{j_k=0}^{p_k}C^2_{j_k\ldots j_1}\right),
\end{equation}

\vspace{2mm}
\noindent
where $T-t\in(0,\infty)$ and
$i_1,\ldots,i_k=1,\dots,m$.

From (\ref{99.010}), (\ref{99.010a}), (\ref{riemann}), (\ref{y007a}), (\ref{y5}),
and the orthonormality
of the system $\{\phi_j(x)\}_{j=0}^{\infty}$
we obtain

\vspace{-4mm}
$$
{\sf M}\left\{\left(J[\psi^{(k)}]_{T,t}-J[\psi^{(k)}]_{T,t}^{p_1,\ldots,p_k}
\right)^2\right\}
\le 
$$

\vspace{-5mm}
$$
\le 
C_k
\sum_{(t_1,\ldots,t_k)}
\int\limits_{t}^{T}
\ldots
\int\limits_{t}^{t_2}
\Biggl(K(t_1,\ldots,t_k)-
\sum_{j_1=0}^{p_1}\ldots
\sum_{j_k=0}^{p_k}
C_{j_k\ldots j_1}
\prod_{l=1}^k\phi_{j_l}(t_l)\Biggr)^2
dt_1
\ldots
dt_k=
$$

\vspace{-2mm}
$$
=C_k\int\limits_{[t,T]^k}
\Biggl(K(t_1,\ldots,t_k)-
\sum_{j_1=0}^{p_1}\ldots
\sum_{j_k=0}^{p_k}
C_{j_k\ldots j_1}
\prod_{l=1}^k\phi_{j_l}(t_l)\Biggr)^2
dt_1
\ldots
dt_k=
$$
$$
=C_k\left(~\int\limits_{[t,T]^k}
K^2(t_1,\ldots,t_k)
dt_1\ldots dt_k -\sum_{j_1=0}^{p_1}\ldots
\sum_{j_k=0}^{p_k}C^2_{j_k\ldots j_1}\right),
$$

\vspace{2mm}
\noindent
where $i_1,\ldots,i_k=0, 1,\ldots,m$,\ \ 
$i_1^2+\ldots+i_k^2>0$, and $C_k$ is a constant.

It is not difficult to see that the constant $C_k$ depends on 
$k$ ($k$ is the multiplicity 
of the iterated It\^{o} stochastic integral) and $T-t$ ($T-t$ is 
the length
of integration interval of the iterated It\^{o} stochastic integral).
Moreover, $C_k$ has the following form
$$
C_k=k!\cdot{\rm max}\Bigl\{
(T-t)^{\alpha_1},\ (T-t)^{\alpha_2},\ \ldots,\ (T-t)^{\alpha_{k!}}
\Bigr\},
$$

\noindent
where $\alpha_1, \alpha_2, \ldots, \alpha_{k!}=0,\ 1,\ldots,\ k-1.$

However, $T-t$ is an integration step of
numerical procedures 
for It\^{o} SDEs (see Chapter 4), which is 
a rather small value. For example, $0<T-t<1.$ Then $C_k\le k!$

It means that for the case
$i_1,\ldots,i_k = 0, 1,\ldots,m$,\  $i_1^2+\ldots+i_k^2>0$, and 
$0<T-t<1$
we get (\ref{z1}). Theorem 1.4 is proved.

{\bf Example 1.3.}\ The particular case of
the estimate (\ref{z1})
for the iterated It\^{o} stochastic integral $I_{(000)T,t}^{(i_1i_2 i_3)}$ 
(see (\ref{k1001})) has the following form
$$
{\sf M}\left\{\left(
I_{(000)T,t}^{(i_1i_2 i_3)}-
I_{(000)T,t}^{(i_1i_2 i_3)p}\right)^2\right\}\le
6\left(\frac{(T-t)^{3}}{6}-\sum_{j_1,j_2,j_3=0}^{p}
C_{j_3j_2j_1}^2\right),
$$

\noindent
where $i_1, i_2, i_3=1,\ldots,m$ and 
$C_{j_3j_2j_1}$ is defined by the formula (\ref{w1}).

Let us consider the case 
of pairwise different 
$i_1,\ldots,i_k=1,\ldots,m$ and prove
the following equality
$$
{\sf M}\left\{\left(J[\psi^{(k)}]_{T,t}-
J[\psi^{(k)}]_{T,t}^{p_1,\ldots,p_k}\right)^2\right\}=
$$
\begin{equation}
\label{1020}
=
\int\limits_{[t,T]^k} K^2(t_1,\ldots,t_k)
dt_1\ldots dt_k
-\sum_{j_1=0}^{p_1}\ldots\sum_{j_k=0}^{p_k}C_{j_k\ldots j_1}^2,
\end{equation}

\noindent
where 
notations are the same as in Theorem 1.4.

The stochastic integrals on the
right-hand side of 
(\ref{y007a}) are uncorrelated 
for the case of pairwise different 
$i_1,\ldots,i_k=1,\ldots,m.$ Moreover, these integrals
have zero expectations. Then
$$
{\sf M}\left\{\left(
J[\psi^{(k)}]_{T,t}-
J[\psi^{(k)}]_{T,t}^{p_1,\ldots,p_k}\right)^2\right\}=
$$
$$
=
{\sf M}\left\{\Biggl(
\sum_{(t_1,\ldots,t_k)}
\int\limits_{t}^{T}
\ldots
\int\limits_{t}^{t_2}
\Biggl(K(t_1,\ldots,t_k)
-
\sum_{j_1=0}^{p_1}\ldots
\sum_{j_k=0}^{p_k}
C_{j_k\ldots j_1}
\prod_{l=1}^k\phi_{j_l}(t_l)\Biggr)\times\Biggr.\right.
$$
$$
\left.\Biggl.\times
d{\bf w}_{t_1}^{(i_1)}
\ldots
d{\bf w}_{t_k}^{(i_k)}\Biggr)^2\right\}=
$$
$$
=
\sum_{(t_1,\ldots,t_k)}
{\sf M}\left\{\Biggl(
\int\limits_{t}^{T}
\ldots
\int\limits_{t}^{t_2}
\Biggl(K(t_1,\ldots,t_k)-
\sum_{j_1=0}^{p_1}\ldots
\sum_{j_k=0}^{p_k}
C_{j_k\ldots j_1}
\prod_{l=1}^k\phi_{j_l}(t_l)\Biggr)\times\Biggr.\right.
$$
$$
\left.\Biggl.\times
d{\bf w}_{t_1}^{(i_1)}
\ldots
d{\bf w}_{t_k}^{(i_k)}\Biggr)^2\right\}=
$$
$$
=
\sum_{(t_1,\ldots,t_k)}
\int\limits_{t}^{T}
\ldots
\int\limits_{t}^{t_2}
\Biggl(K(t_1,\ldots,t_k)-
\sum_{j_1=0}^{p_1}\ldots
\sum_{j_k=0}^{p_k}
C_{j_k\ldots j_1}
\prod_{l=1}^k\phi_{j_l}(t_l)\Biggr)^2
dt_1\ldots dt_k=
$$
$$
=\int\limits_{[t,T]^k}
\Biggl(K(t_1,\ldots,t_k)-\sum_{j_1=0}^{p_1}\ldots
\sum_{j_k=0}^{p_k}
C_{j_k\ldots j_1}
\prod_{l=1}^k\phi_{j_l}(t_l)\Biggr)^2
dt_1\ldots dt_k=
$$
$$
=\int\limits_{[t,T]^k}
K^2(t_1,\ldots,t_k)
dt_1\ldots dt_k -\sum_{j_1=0}^{p_1}\ldots
\sum_{j_k=0}^{p_k}C^2_{j_k\ldots j_1}.
$$

\section{Expansion of Iterated It\^{o} Stochastic Integrals 
Based on Generalized Multiple Fourier Series.
The Case of Complete Orthonormal with We\-ight $r(t_1)\ldots r(t_k)$  
Systems of Functions 
in the Space $L_2([t, T]^k)$}

In this section, we consider a modification of Theorem 1.1 for 
the case of complete orthonormal with weight $r(t_1)\ldots r(t_k)\ge 0$ 
systems of functions 
in the space $L_2([t, T]^k)$, $k\in{\bf N}$.${}^7$

\footnotetext[7]{The results of this section 
are generalized to the case of an arbitrary 
complete ortho\-nor\-mal with weight $r(x)\ge 0$
system of functions $\{\Psi_j(x)\sqrt{r(x)}\}_{j=0}^{\infty}$ 
in the space $L_2([t, T])$
and $\psi_1(x)\sqrt{r(x)},$ $\ldots,\psi_k(x)\sqrt{r(x)} \in L_2([t, T])$ in Sect.~1.13
(see Theorems~1.20, 1.21).}

Let $\{\Psi_j(x)\}_{j=0}^{\infty}$ be a complete orthonormal 
with weight $r(x)\ge 0$ 
system of functions in the space $L_2([t, T]).$ It is well known that the
Fourier
series of the function $f(x)$ 
$\left(f(x)\sqrt{r(x)}\in L_2([t, T])\right)$ 
with respect to the system $\{\Psi_j(x)\}_{j=0}^{\infty}$
converges 
to the function $f(x)$ in the
mean-square sense with weight $r(x),$ i.e.
\begin{equation}
\label{g1}
\lim\limits_{p\to\infty}
\int\limits_t^T\biggl(f(x)-\sum\limits_{j=0}^p 
{\tilde C}_j \Psi_j(x)\biggr)^2 r(x)dx = 0,
\end{equation}
where
\begin{equation}
\label{h1}
{\tilde C}_j=\int\limits_t^T f(x)\Psi_j(x)r(x)dx
\end{equation}

\noindent
is the Fourier coefficient.

The relations (\ref{g1}), (\ref{h1}) can be obtained if we will 
expand the function
$f(x)\sqrt{r(x)}\in L_2([t, T])$ 
into a usual Fourier series with respect
to the complete orthonormal with weight $1$ system of functions
$$
\left\{\Psi_j(x)\sqrt{r(x)}\right\}_{j=0}^{\infty}
$$ 
in 
the space $L_2([t, T]).$ Then
$$
\lim\limits_{p\to\infty}
\int\limits_t^T\biggl(f(x)\sqrt{r(x)}-\sum\limits_{j=0}^p {\tilde C}_j 
\Psi_j(x)\sqrt{r(x)}\biggr)^2dx = 
$$
\begin{equation}
\label{g2}
=\lim\limits_{p\to\infty}
\int\limits_t^T\biggl(f(x)-
\sum\limits_{j=0}^p {\tilde C}_j \Psi_j(x)\biggr)^2 r(x)dx = 0,
\end{equation}

\vspace{2mm}
\noindent
where ${\tilde C}_j$ is defined by (\ref{h1}).

Let us consider an obvious generalization of this approach 
to the case of $k$
variables.
Let us expand the function $K(t_1,\ldots,t_k)$ such that
$$
K(t_1,\ldots,t_k)\prod\limits_{l=1}^k \sqrt{r(t_l)}\in L_2([t, T]^k)
$$
using the complete orthonormal system of functions 
$$
\prod\limits_{l=1}^k \Psi_{j_l}(t_l)\sqrt{r(t_l)},\ \ \ 
j_l=0, 1, 2, \ldots,\ \ \  l=1,\ldots,k
$$
in the space $L_2([t, T]^k)$ into the generalized multiple Fourier 
series. 

It is well known that the mentioned
generalized multiple Fourier series converges in the mean-square sense,
i.e.

\vspace{-7mm}
$$
\lim\limits_{p_1,\ldots,p_k\to\infty}
\int\limits_{[t,T]^k}
\hspace{-1.5mm}\left(K(t_1,\ldots,t_k)\prod\limits_{l=1}^k \sqrt{r(t_l)}-
\sum\limits_{j_1=0}^{p_1}\ldots\sum\limits_{j_k=0}^{p_k}
{\tilde C}_{j_k\ldots j_1}
\prod\limits_{l=1}^k \Psi_{j_l}(t_l)\sqrt{r(t_l)}\right)^2
\hspace{-1.5mm}\times
$$
$$
\times
dt_1\ldots dt_k=
$$

\vspace{-6mm}
$$
=\lim\limits_{p_1,\ldots,p_k\to\infty}
\int\limits_{[t,T]^k}
\left(K(t_1,\ldots,t_k)-
\sum\limits_{j_1=0}^{p_1}\ldots\sum\limits_{j_k=0}^{p_k}
{\tilde C}_{j_k\ldots j_1}\prod\limits_{l=1}^k \Psi_{j_l}(t_l)\right)^2 
\prod\limits_{l=1}^k r(t_l)
\times
$$
\begin{equation}
\label{z1aaa}
\times
dt_1\ldots dt_k=0,
\end{equation}

\vspace{2mm}
\noindent
where
$$
{\tilde C}_{j_k\ldots j_1}=\int\limits_{[t,T]^k}
K(t_1,\ldots,t_k)\prod\limits_{l=1}^k 
\biggl(\Psi_{j_l}(t_l)r(t_l)\biggr)dt_1\ldots dt_k.
$$

Let us consider 
the following iterated It\^{o} 
stochastic integrals
\begin{equation}
\label{ito-rr}
~~~~~~~ {\tilde J}[\psi^{(k)}]_{T,t}=\int\limits_t^T\psi_k(t_k)\sqrt{r(t_k)} 
\ldots \int\limits_t^{t_{2}}
\psi_1(t_1)\sqrt{r(t_1)} d{\bf w}_{t_1}^{(i_1)}\ldots
d{\bf w}_{t_k}^{(i_k)},
\end{equation}
where $\psi_l(\tau)$ $(l=1,\ldots,k)$ are
nonrandom functions on $[t, T]$,
${\bf w}_{\tau}^{(i)}$ $(i=1,\ldots,m)$ are independent
standard Wiener processes,
${\bf w}_{\tau}^{(0)}=\tau,$ 
$i_1,\ldots,i_k=0, 1,\ldots,m.$

So, we obtain the following version of Theorem 1.1.

{\bf Theorem 1.5}\ \cite{12}-\cite{12aa}, \cite{arxiv-1}, \cite{arxiv-13}.
{\it Suppose that every $\psi_l(\tau)$ $(l=$ $1,\ldots, k)$ 
is a continuous 
nonrandom function on 
$[t, T].$ Moreover$,$ let 
$\{\Psi_j(x)\sqrt{r(x)}\}_{j=0}^{\infty}$ $(r(x)\ge 0)$
is a complete orthonormal 
system of functions in the space $L_2([t,T]),$ each function 
$\Psi_j(x)\sqrt{r(x)}$
of which 
for finite $j$ satisfies the condition 
$(\star)$ {\rm (}see Sect.~{\rm 1.1.7)}.
Then
$$
{\tilde J}[\psi^{(k)}]_{T,t} =
\hbox{\vtop{\offinterlineskip\halign{
\hfil#\hfil\cr
{\rm l.i.m.}\cr
$\stackrel{}{{}_{p_1,\ldots,p_k\to \infty}}$\cr
}} }\sum_{j_1=0}^{p_1}\ldots\sum_{j_k=0}^{p_k}
{\tilde C}_{j_k\ldots j_1}\Biggl(
\prod_{l=1}^k{\tilde \zeta}_{j_l}^{(i_l)} -
\Biggr.
$$

\vspace{-2mm}
\begin{equation}
\label{tyyy-rr}
~~~~ -\Biggl.
\hbox{\vtop{\offinterlineskip\halign{
\hfil#\hfil\cr
{\rm l.i.m.}\cr
$\stackrel{}{{}_{N\to \infty}}$\cr
}} }\sum_{(l_1,\ldots,l_k)\in {\rm G}_k}
\Psi_{j_{1}}(\tau_{l_1})\sqrt{r(\tau_{l_1})}
\Delta{\bf w}_{\tau_{l_1}}^{(i_1)}\ldots
\Psi_{j_{k}}(\tau_{l_k})\sqrt{r(\tau_{l_k})}
\Delta{\bf w}_{\tau_{l_k}}^{(i_k)}\Biggr),
\end{equation}

\vspace{2mm}
\noindent
where
$$
{\rm G}_k={\rm H}_k\backslash{\rm L}_k,\ \ \
{\rm H}_k=\bigl\{(l_1,\ldots,l_k):\ l_1,\ldots,l_k=0,\ 1,\ldots,N-1\bigr\},
$$
$$
{\rm L}_k=\bigl\{(l_1,\ldots,l_k):\ l_1,\ldots,l_k=0,\ 1,\ldots,N-1;\
l_g\ne l_r\ (g\ne r);\ g, r=1,\ldots,k\bigr\},
$$

\noindent
${\rm l.i.m.}$ is a limit in the mean-square sense$,$
$i_1,\ldots,i_k=0,1,\ldots,m,$ 
$$
{\tilde \zeta}_{j}^{(i)}=
\int\limits_t^T \Psi_{j}(s)\sqrt{r(s)}d{\bf w}_s^{(i)}
$$
are independent standard Gaussian random variables
for various
$i$ or $j$ {\rm(}in the case when $i\ne 0${\rm),}
$\Delta{\bf w}_{\tau_{j}}^{(i)}=
{\bf w}_{\tau_{j+1}}^{(i)}-{\bf w}_{\tau_{j}}^{(i)}$
$(i=0, 1,\ldots,m),$
$\left\{\tau_{j}\right\}_{j=0}^{N}$ is a partition of
$[t,T],$ which satisfies the condition {\rm (\ref{1111}),}
\begin{equation}
\label{koef}
~~~~~~~ {\tilde C}_{j_k\ldots j_1}=\int\limits_{[t,T]^k}
K(t_1,\ldots,t_k)
\prod_{l=1}^{k}\biggl(\Psi_{j_l}(t_l)r(t_l)\biggr)dt_1\ldots dt_k
\end{equation}
is the Fourier coefficient$,$
$$
K(t_1,\ldots,t_k)=
\left\{
\begin{matrix}
\psi_1(t_1)\ldots \psi_k(t_k),\ &t_1<\ldots<t_k\cr\cr
0,\ &\hbox{otherwise}
\end{matrix}\right.,\ \  t_1,\ldots,t_k\in[t, T],\ \  k\ge 2,
$$

\vspace{1mm}
\noindent
and 
$K(t_1)\equiv\psi_1(t_1)$ for $t_1\in[t, T].$  
}

{\bf Proof.}
According to Lemmas 1.1, 1.3 and (\ref{pobeda}), 
(\ref{s2s}), (\ref{hehe100}), (\ref{zab1}),
we get the following representation

\vspace{-2mm}
$$
{\tilde J}[\psi^{(k)}]_{T,t}=
\sum_{(t_1,\ldots,t_k)}
\int\limits_{t}^{T}
\ldots
\int\limits_{t}^{t_2}
K(t_1,\ldots,t_k)\prod\limits_{l=1}^k \sqrt{r(t_l)}d{\bf w}_{t_1}^{(i_1)}
\ldots
d{\bf w}_{t_k}^{(i_k)}=
$$

\vspace{1mm}
$$
=
\sum_{j_1=0}^{p_1}\ldots
\sum_{j_k=0}^{p_k}
{\tilde C}_{j_k\ldots j_1}
\int\limits_{t}^{T}
\ldots
\int\limits_{t}^{t_2}
\sum_{(t_1,\ldots,t_k)}\left(
\prod\limits_{l=1}^k\left(\Psi_{j_l}(t_l)\sqrt{r(t_l)}\right)
d{\bf w}_{t_1}^{(i_1)}
\ldots
d{\bf w}_{t_k}^{(i_k)}\right)
+
$$

\vspace{1mm}
$$
+{\tilde R}_{T,t}^{p_1,\ldots,p_k}=
$$

\vspace{5mm}
$$
=\sum_{j_1=0}^{p_1}\ldots
\sum_{j_k=0}^{p_k}
{\tilde C}_{j_k\ldots j_1}\times
$$

\vspace{-2mm}
$$
\times             
\hbox{\vtop{\offinterlineskip\halign{
\hfil#\hfil\cr
{\rm l.i.m.}\cr
$\stackrel{}{{}_{N\to \infty}}$\cr
}} }
\sum\limits_{\stackrel{l_1,\ldots,l_k=0}{{}_{l_q\ne l_r;\ 
q\ne r;\ q, r=1,\ldots, k}}}^{N-1}
\Psi_{j_1}(\tau_{l_1})\sqrt{r(\tau_{l_1})}\Delta{\bf w}_{\tau_{l_1}}^{(i_1)}
\ldots
\Psi_{j_k}(\tau_{l_k})\sqrt{r(\tau_{l_k})}
\Delta{\bf w}_{\tau_{l_k}}^{(i_k)}+
$$

$$
+{\tilde R}_{T,t}^{p_1,\ldots,p_k}= 
$$

\newpage
\noindent
$$
=\sum_{j_1=0}^{p_1}\ldots
\sum_{j_k=0}^{p_k}
{\tilde C}_{j_k\ldots j_1}\times
$$

\vspace{-2mm}
$$
\times\left(
\hbox{\vtop{\offinterlineskip\halign{
\hfil#\hfil\cr
{\rm l.i.m.}\cr
$\stackrel{}{{}_{N\to \infty}}$\cr
}} }\sum_{l_1,\ldots,l_k=0}^{N-1}
\Psi_{j_1}(\tau_{l_1})\sqrt{r(\tau_{l_1})}\Delta{\bf w}_{\tau_{l_1}}^{(i_1)}
\ldots
\Psi_{j_k}(\tau_{l_k})\sqrt{r(\tau_{l_k})}
\Delta{\bf w}_{\tau_{l_k}}^{(i_k)}
-\right.
$$

\vspace{-2mm}
$$
-\left.
\hbox{\vtop{\offinterlineskip\halign{
\hfil#\hfil\cr
{\rm l.i.m.}\cr
$\stackrel{}{{}_{N\to \infty}}$\cr
}} }\sum_{(l_1,\ldots,l_k)\in {\rm G}_k}
\Psi_{j_1}(\tau_{l_1})\sqrt{r(\tau_{l_1})}\Delta{\bf w}_{\tau_{l_1}}^{(i_1)}
\ldots
\Psi_{j_k}(\tau_{l_k})\sqrt{r(\tau_{l_k})}
\Delta{\bf w}_{\tau_{l_k}}^{(i_k)}
\right)
+
$$

$$
+{\tilde R}_{T,t}^{p_1,\ldots,p_k}=
$$

\vspace{5mm}
$$
=\sum_{j_1=0}^{p_1}\ldots\sum_{j_k=0}^{p_k}
{\tilde C}_{j_k\ldots j_1}\times
$$

\vspace{-2mm}
$$
\times
\hspace{-1.5mm}\left(
\prod_{l=1}^k {\tilde \zeta}_{j_l}^{(i_l)}-
\hbox{\vtop{\offinterlineskip\halign{
\hfil#\hfil\cr
{\rm l.i.m.}\cr
$\stackrel{}{{}_{N\to \infty}}$\cr
}} }\hspace{-1.3mm}
\sum_{(l_1,\ldots,l_k)\in {\rm G}_k}\hspace{-1.4mm}
\Psi_{j_1}(\tau_{l_1})\sqrt{r(\tau_{l_1})}\Delta{\bf w}_{\tau_{l_1}}^{(i_1)}
\ldots
\Psi_{j_k}(\tau_{l_k})\sqrt{r(\tau_{l_k})}
\Delta{\bf w}_{\tau_{l_k}}^{(i_k)}
\right)\hspace{-1.3mm}+
$$

\vspace{1mm}
$$
+{\tilde R}_{T,t}^{p_1,\ldots,p_k}\ \ \ \hbox{w.~p.~1},
$$

\vspace{4mm}
\noindent
where

\vspace{-2mm}
$$
{\tilde R}_{T,t}^{p_1,\ldots,p_k}
=\sum_{(t_1,\ldots,t_k)}
\int\limits_{t}^{T}
\ldots
\int\limits_{t}^{t_2}
\left(K(t_1,\ldots,t_k)\prod_{l=1}^k\sqrt{r(t_l)}-\right.
$$

\vspace{1mm}
$$
\left.
-\sum_{j_1=0}^{p_1}\ldots
\sum_{j_k=0}^{p_k}
{\tilde C}_{j_k\ldots j_1}
\prod_{l=1}^k\left(\Psi_{j_l}(t_l)\sqrt{r(t_l)}\right)\right)
d{\bf w}_{t_1}^{(i_1)}
\ldots
d{\bf w}_{t_k}^{(i_k)},
$$

\vspace{3mm}
\noindent
where permutations $(t_1,\ldots,t_k)$ when summing are performed only 
in the values $d{\bf w}_{t_1}^{(i_1)}
\ldots $
$d{\bf w}_{t_k}^{(i_k)}$. At the same time the indices near 
upper limits of integration in the iterated stochastic integrals 
are changed correspondently and if $t_r$ swapped with $t_q$ in the  
permutation $(t_1,\ldots,t_k)$, then $i_r$ swapped with $i_q$ in the 
permutation $(i_1,\ldots,i_k)$.

Let us estimate the remainder
${\tilde R}_{T,t}^{p_1,\ldots,p_k}$ of the series.

According to Lemma 1.2 and (\ref{riemann}), we have
$$
{\sf M}\left\{\left({\tilde R}_{T,t}^{p_1,\ldots,p_k}\right)^2\right\}
\le C_k
\sum_{(t_1,\ldots,t_k)}
\int\limits_{t}^{T}
\ldots
\int\limits_{t}^{t_2}
\left(K(t_1,\ldots,t_k)\prod_{l=1}^k\sqrt{r(t_l)}\right.-
$$

\vspace{-2mm}
\begin{equation}
\label{obana1eee}
~~~~~~~\left.-
\sum_{j_1=0}^{p_1}\ldots
\sum_{j_k=0}^{p_k}
{\tilde C}_{j_k\ldots j_1}
\prod_{l=1}^k\left(\Psi_{j_l}(t_l)\sqrt{r(t_l)}\right)\right)^2
dt_1
\ldots
dt_k=
\end{equation}

$$
=C_k\int\limits_{[t,T]^k}
\left(K(t_1,\ldots,t_k)-
\sum_{j_1=0}^{p_1}\ldots
\sum_{j_k=0}^{p_k}
{\tilde C}_{j_k\ldots j_1}
\prod_{l=1}^k\Psi_{j_l}(t_l)\right)^2
\times
$$

\vspace{-1mm}
\begin{equation}
\label{arch1}
~~~~\times
\left(\prod_{l=1}^k r(t_l)\right)
dt_1 \ldots
dt_k\to 0
\end{equation}

\vspace{4mm}
\noindent
if $p_1,\ldots,p_k\to\infty,$ where constant $C_k$ 
depends only
on the multiplicity $k$ of the iterated It\^{o} stochastic integral
(\ref{ito-rr}). 
Theorem 1.5 is proved.

Let us formulate the version of Theorem 1.4.

{\bf Theorem 1.6} \cite{12a}-\cite{12aa}, \cite{arxiv-1}, \cite{arxiv-13}. 
{\it Suppose that every $\psi_l(\tau)$ $(l=$ $1,\ldots, k)$ 
is a continuous 
nonrandom function on 
$[t, T].$ Moreover$,$ let 
$\{\Psi_j(x)\sqrt{r(x)}\}_{j=0}^{\infty}$ $(r(x)\ge 0)$
is a complete orthonormal 
system of functions in the space $L_2([t,T]),$ each function 
$\Psi_j(x)\sqrt{r(x)}$
of which 
for finite $j$ satisfies the condition 
$(\star)$ {\rm (}see Sect.~{\rm 1.1.7)}.
Then the estimate

\vspace{-2mm}
$$
{\sf M}\left\{\left(
{\tilde J}[\psi^{(k)}]_{T,t}-{\tilde J}[\psi^{(k)}]_{T,t}^{p_1,\ldots,p_k}
\right)^2\right\}
\le 
$$

\vspace{-5mm}
\begin{equation}
\label{z1-ura}
~ \le k!\left(~\int\limits_{[t,T]^k}
K^2(t_1,\ldots,t_k)\left(\prod_{l=1}^k r(t_l)\right)
dt_1\ldots dt_k -\sum_{j_1=0}^{p_1}\ldots
\sum_{j_k=0}^{p_k}{\tilde C}^2_{j_k\ldots j_1}\right)
\end{equation}

\vspace{2mm}
\noindent
is valid for the following cases{\rm :}

{\rm 1.}\ $i_1,\ldots,i_k=1,\ldots,m$\ \ and\ \ $0<T-t<\infty,$

{\rm 2.}\ $i_1,\ldots,i_k=0, 1,\ldots,m,$\ \ $i_1^2+\ldots+i_k^2>0,$\ \
and\ \ $0<T-t<1,$

\noindent
where ${\tilde J}[\psi^{(k)}]_{T,t}$ is the 
stochastic integral {\rm (\ref{ito-rr}),}
${\tilde J}[\psi^{(k)}]_{T,t}^{p_1,\ldots,p_k}$ is the 
expression on the right-hand side of {\rm (\ref{tyyy-rr})} before
passing to the limit 
$\hbox{\vtop{\offinterlineskip\halign{
\hfil#\hfil\cr
{\rm l.i.m.}\cr
$\stackrel{}{{}_{p_1,\ldots,p_k\to \infty}}$\cr
}} };$ another 
notations are the same as in Theorem {\rm 1.5}.
}

\section{Expansion of Iterated Stochastic Integrals 
with Respect to Martingale Poisson Measures
Based on Generalized Multiple Fourier Series}

In this section, we consider the version 
of Theorem 1.1 connected with the expansion
of iterated stochastic integrals
with respect to martingale Poisson measures.

\subsection{Stochastic Integral with Respect to 
Martingale Poisson Measure}

Let us consider the Poisson random measure on the set
$[0,T]\times{\bf Y}$
$({\bf R}^n\stackrel{\rm def}{=}{\bf Y})$.
We will denote the value of this measure at the set 
$\Delta\times A$ ($\Delta\subseteq
[0,T],$ $A\subset{\bf Y}$) 
as 
$\nu(\Delta,A).$ Assume that 
$$
{\sf M}\left\{\nu(\Delta,A)\right\}=
|\Delta|\Pi(A),
$$ 
where $|\Delta|$ is the Lebesgue measure
of $\Delta,$ 
$\Pi(A)$ is a measure on $\sigma$-algebra $\mathcal{B}$ of Borel
subsets of ${\bf Y},$ and $\mathcal{B}_0$ is a subalgebra of $\mathcal{B}$ 
consisting of sets
$A\subset \mathcal{B}$ that satisfy the condition 
$\Pi(A)<\infty.$

Let us consider the martingale Poisson measure
$$
\tilde\nu(\Delta,A)=\nu(\Delta,A)-|\Delta|\Pi(A).
$$

Let $(\Omega, {\rm F},{\sf P})$ be a fixed probability
space, let $\{{\rm F}_t,$ $t\in[0,T]\}$
be a non-decreasing family of 
$\sigma$-algebras ${\rm F}_t\subset{\rm F}$.

Assume that the following conditions are fulfilled:

1. The random variables $\nu([0,t),A)$ are ${\rm F}_t$-measurable
for all $A\subseteq \mathcal{B}_0,$ $t\in[0, T].$

2. The random variables $\nu([t,t+h),A),$ $A\subseteq \mathcal{B}_0,$
$h>0$ do not depend on events of $\sigma$-algebra ${\rm F}_t.$

Let us define the class $H_l(\Pi,[0,T])$ 
of random functions  
$\varphi:$ $[0,T]\times{\bf Y}\times\Omega\to{\bf R}^1$
that are ${\rm F}_t$-measurable for all
$t\in[0,T],$ ${\bf y}\in{\bf Y}$ 
and satisfy the
following condition 
$$
\int\limits_0^T\int\limits_{{\bf Y}}
{\sf M}\Bigl\{|\varphi(t,{\bf y})|^l\Bigr\}\Pi(d{\bf y})dt<\infty.
$$

Consider the partition $\{\tau_j\}_{j=0}^N$ of the 
interval $[0,T],$
which 
satisfies
the condition (\ref{1111}), and  
define the stochastic 
integral with respect to the martingale Poisson measure
for $\varphi(t,{\bf y})\in H_2(\Pi,[0,T])$ as the following 
mean-square limit 
\cite{Gih1}
\begin{equation}
\label{1.10}
~~~~~~~~~~\int\limits_0^T\int\limits_{{\bf Y}}
\varphi(t,{\bf y})\tilde\nu(dt,d{\bf y})
\stackrel{\rm def}{=}
\hbox{\vtop{\offinterlineskip\halign{
\hfil#\hfil\cr
{\rm l.i.m.}\cr
$\stackrel{}{{}_{N\to\infty}}$\cr
}} }\int\limits_0^T\int\limits_{{\bf Y}}
\varphi^{(N)}(t,{\bf y})\tilde\nu(dt,d{\bf y}),
\end{equation}
where $\varphi^{(N)}(t,{\bf y})$ is any
sequence of step functions 
from the class $H_2(\Pi,[0,T])$ such that 
$$
\hbox{\vtop{\offinterlineskip\halign{
\hfil#\hfil\cr
{\rm lim}\cr
$\stackrel{}{{}_{N\to\infty}}$\cr
}} }
\int\limits_0^T\int\limits_{{\bf Y}}{\sf M}\left\{\left|\varphi(t,{\bf y})-
\varphi^{(N)}(t,{\bf y})\right|^2\right\}\Pi(d{\bf y})dt\to 0.
$$

It is well known \cite{Gih1} that the stochastic integral 
(\ref{1.10}) exists, it does not depend on selection of the
sequence $\varphi^{(N)}(t,{\bf y})$ and it 
satisfies
w.~p.~1 
to the fol\-lowing properties
$$
{\sf M}\left\{\int\limits_0^T
\int\limits_{{\bf Y}}\varphi(t,{\bf y})\tilde\nu(dt,d{\bf y})
\biggl|\biggr.{\rm F}_0\right\}=0,
$$
$$
\int\limits_0^T\int\limits_{{\bf Y}}
(\alpha\varphi_1(t,{\bf y})+\beta\varphi_2(t,{\bf y}))
\tilde\nu(dt,d{\bf y})=
$$
$$
=
\alpha\int\limits_0^T
\int\limits_{{\bf Y}}\varphi_1(t,{\bf y})\tilde\nu(dt,d{\bf y})
+\beta\int\limits_0^T\int\limits_{{\bf Y}}
\varphi_2(t,{\bf y})\tilde\nu(dt,d{\bf y}),
$$

$$
{\sf M}\left\{\left|\int\limits_0^T
\int\limits_{{\bf Y}}\varphi(t,{\bf y})\tilde\nu(dt,d{\bf y})
\right|^2 \biggl|\biggr.{\rm F}_0\right\}=
\int\limits_0^T\int\limits_{{\bf Y}}{\sf M}\Bigl\{\left|\varphi(t,{\bf y})
\right|^2 \Bigl|\Bigr.{\rm F}_0
\Bigr\}\Pi(d{\bf y})dt,
$$

\vspace{2mm}
\noindent
where $\alpha,$ $\beta \in {\bf R}^1$ and $\varphi_1(t,{\bf y}),$
$\varphi_2(t,{\bf y}),$
$\varphi(t,{\bf y})$ from the class 
$H_2(\Pi,[0,T]).$

The stochastic integral
$$
\int\limits_0^T\int\limits_{{\bf Y}}
\varphi(t,{\bf y})\nu(dt,d{\bf y})
$$
with respect to the Poisson measure will be defined as follows \cite{Gih1}
\begin{equation}
\label{faru}
\int\limits_0^T\int\limits_{{\bf Y}}
\varphi(t,{\bf y})\nu(dt,d{\bf y})=
\int\limits_0^T\int\limits_{{\bf Y}}
\varphi(t,{\bf y})\tilde\nu(dt,d{\bf y})+
\int\limits_0^T\int\limits_{{\bf Y}}
\varphi(t,{\bf y})\Pi(d{\bf y})dt,
\end{equation}
where we suppose that the right-hand side of the last equality
exists.

According to the It\^{o} formula for It\^{o} processes with 
jumps, we get \cite{Gih1}
\begin{equation}
\label{16.008}
~~~~~~~~ \left({z}_t\right)^p=\int\limits_0^t\int\limits_{{\bf Y}}
\Bigl(({z}_{\tau-}+\gamma(\tau,{\bf y}))^p-\left({z}_{\tau-}\right)^p\Bigr)
\nu(d\tau,d{\bf y})\ \ \ \hbox{w.~p.~1},
\end{equation}
where $p\in{\bf N}$ and ${z}_{\tau-}$ means the left-sided limit value of the process
${z}_{\tau}$ at the point $\tau$,
$$
{z}_t=\int\limits_0^t\int\limits_{{\bf Y}}
\gamma(\tau,{\bf y})
\nu(d\tau,d{\bf y}).
$$

We suppose that the function $\gamma(\tau,{\bf y})$
satisfies
the conditions of existence of the right-hand side 
of (\ref{16.008}) \cite{Gih1}.

Let us consider the useful estimate for moments
of stochastic 
integrals with respect to the Poisson measure \cite{Gih1}
\begin{equation}
\label{16.010}
~~~ a_p(T)\le
\max\limits_{j\in\{p,\ 1\}}
\left\{\left(\int\limits_0^T\int\limits_{{\bf Y}}\left(
\left(\left(b_p(\tau,{\bf y})\right)^{1/p}+1\right)^p-1\right)
\Pi(d{\bf y})d\tau\right)^j\right\},
\end{equation}
where 
$$
a_p(t)=\sup\limits_{0\le\tau\le t}
{\sf M}\Bigl\{|{z}_{\tau}|^p\Bigr\},\ \ \
b_p(\tau,{\bf y})={\sf M}\Bigl\{\left|
\gamma(\tau,{\bf y})\right|^p\Bigr\}.
$$

\vspace{2mm}

We suppose that the right-hand side 
of (\ref{16.010}) exists. According to (see (\ref{faru}))
$$
\int\limits_0^t\int\limits_{{\bf Y}}
\gamma(\tau,{\bf y})\tilde \nu(d\tau,d{\bf y})=
\int\limits_0^t\int\limits_{{\bf Y}}
\gamma(\tau,{\bf y})\nu(d\tau,d{\bf y})-
\int\limits_0^t\int\limits_{{\bf Y}}
\gamma(\tau,{\bf y})\Pi(d{\bf y})d\tau
$$
and the Minkowski inequality, we obtain

\vspace{-4mm}
\begin{equation}
\label{16.011}
~~~~~~~~~ \left({\sf M}\Bigl\{\left|\tilde{z}_{t}\right|^{2p}\Bigr\}
\right)^{1/2p}\le
\left({\sf M}\Bigl\{\left|{z}_{t}\right|^{2p}\Bigr\}\right)^{1/2p}+
\left({\sf M}\Bigl\{\left|\hat{z}_{t}\right|^{2p}\Bigr\}
\right)^{1/2p},
\end{equation}

\noindent
where
$$
\tilde{z}_t=\int\limits_0^t\int\limits_{{\bf Y}}
\gamma(\tau,{\bf y})
\tilde\nu(d\tau,d{\bf y})
$$
and
$$
\hat{z}_{t}\stackrel{\rm def}{=}\int\limits_0^t\int\limits_{{\bf Y}}
\gamma(\tau,{\bf y})\Pi(d{\bf y})d\tau.
$$

The value ${\sf M}\Bigl\{|\hat{z}_{\tau}|^{2p}\Bigr\}$
can be estimated using the 
well known inequality \cite{Gih1}

\vspace{-2mm}
\begin{equation}
\label{dur}
{\sf M}\left\{|\hat{z}_t|^{2p}\right\}\le t^{2p-1}
\int\limits_{0}^t 
{\sf M}\left\{\left|\int\limits_{{\bf Y}}
\gamma(\tau,{\bf y})\Pi(d{\bf y})\right|^{2p}\right\}d\tau,
\end{equation}

\noindent
where we suppose that 
$$
\int\limits_0^t{\sf M}\left\{\left|\int\limits_{{\bf Y}}
\gamma(\tau,{\bf y})\Pi(d{\bf y})\right|^{2p}\right\}d\tau<\infty.
$$

\subsection{Expansion of Iterated Stochastic Integrals 
with Respect to Martingale Poisson Measures}

Let us consider the following iterated stochastic integrals 

\vspace{-1mm}
$$
P[\chi^{(k)}]_{T,t}=
$$
\begin{equation}
\label{2000.2.1}
=\int\limits_t^T\int\limits_{\bf X}\chi_k(t_k,{\bf y}_k)\ldots
\int\limits_t^{t_2}\int\limits_{\bf X}\chi_1(t_1,{\bf y}_1)
\tilde\nu^{(i_1)}(dt_1,d{\bf y}_1)\ldots
\tilde\nu^{(i_k)}(dt_k,d{\bf y}_k),
\end{equation}

\vspace{1mm}
\noindent
where $i_1,\ldots,i_k=0, 1,\ldots,m,$ 
${\bf R}^n\stackrel{\rm def}{=}{\bf X},$
$\chi_l(\tau,{\bf y})=\psi_l(\tau)\varphi_l({\bf y})$ $(l=1,\ldots,k),$
every function
$\psi_l(\tau): [t,T]\to {\bf R}^1$ $(l=1,\ldots,k)$ and every function 
$\varphi_l({\bf y}): {\bf X}\to {\bf R}^1$ $(l=1,\ldots,k)$ such that
$$
\chi_l(\tau,{\bf y})\in
H_2(\Pi,[t,T])\ \ \ (l=1,\ldots,k),
$$

\vspace{1mm}
\noindent
where definition of the class $H_2(\Pi,[t,T])$ see above,
$\nu^{(i)}(dt,d{\bf y})$ $(i=1,\ldots,m)$ 
are independent Poisson measures for various $i$,
which are defined on
$[0,T]\times {\bf X}$,

\vspace{-1mm}
$$
\tilde\nu^{(i)}(dt,d{\bf y})=
\nu^{(i)}(dt,d{\bf y})-\Pi(d{\bf y})dt\ \ \ (i=1,\ldots,m)
$$ 

\vspace{2mm}
\noindent
are independent martingale Poisson measures for various $i$,
$\tilde\nu^{(0)}(dt,d{\bf y})\stackrel{\rm def}{=}\Pi(d{\bf y})dt,$
$\nu^{(0)}(dt,d{\bf y})\stackrel{\rm def}{=}\Pi(d{\bf y})dt$.

Let us formulate an analogue of Theorem 1.1 for the iterated
stochastic integrals (\ref{2000.2.1}).

{\bf Theorem 1.7} \cite{1}-\cite{12aa}, \cite{arxiv-13}.
{\it Suppose that the following 
conditions are hold{\rm :}

{\rm 1}.\ Every $\psi_l(\tau)$ $(l=1,\ldots,k)$ is a 
continuous nonrandom function at
the interval $[t, T]$.

{\rm 2}.\ $\{\phi_j(x)\}_{j=0}^{\infty}$ is a complete orthonormal
system of functions in the space 
$L_2([t,T]),$ 
each function $\phi_j(x)$ of which for finite $j$ satisfies the condition 
$(\star)$ {\rm(}see Sect.~{\rm 1.1.7)}.

{\rm 3}.\ For $l=1,\ldots,k$ and  $q=2^{k+1}$
the following condition is satisfied
$$
\int\limits_{\bf X}\left|\varphi_l({\bf y})\right|^q
\Pi(d{\bf y})<\infty.
$$

Then$,$ for the iterated stochastic integral with respect to 
martingale Poisson measures $P[\chi^{(k)}]_{T,t}$ defined by
{\rm (\ref{2000.2.1})}
the following expansion 

\vspace{-2mm}
$$
P[\chi^{(k)}]_{T,t}=
\hbox{\vtop{\offinterlineskip\halign{
\hfil#\hfil\cr
{\rm l.i.m.}\cr
$\stackrel{}{{}_{p_1,\ldots,p_k\to \infty}}$\cr
}} }\sum_{j_1=0}^{p_1}\ldots\sum_{j_k=0}^{p_k}
C_{j_k\ldots j_1}
\Biggl(
\prod_{g=1}^k\pi_{j_g}^{(g,i_g)}
\Biggr.-
$$

\begin{equation}
\label{tyyys}
~~~~~~~-\Biggl.
\hbox{\vtop{\offinterlineskip\halign{
\hfil#\hfil\cr
{\rm l.i.m.}\cr
$\stackrel{}{{}_{N\to \infty}}$\cr
}} }\sum_{(l_1,\ldots,l_k)\in {\rm G}_k}
\prod_{g=1}^k 
\phi_{j_{g}}(\tau_{l_g})
\int\limits_{\bf X}\varphi_g({\bf y})
\tilde \nu^{(i_g)}([\tau_{l_g},\tau_{l_g+1}),d{\bf y})
\Biggr)
\end{equation}

\vspace{2mm}
\noindent
that converges in the mean-square sense is valid$,$ 
where $\left\{\tau_{j}\right\}_{j=0}^{N}$ is a partition of
the interval $[t, T]$ satisfying the condition {\rm (\ref{1111}),}
$$
{\rm G}_k={\rm H}_k\backslash{\rm L}_k,\ \ \
{\rm H}_k=\bigl\{(l_1,\ldots,l_k):\ l_1,\ldots,l_k=0,\ 1,\ldots,N-1\bigr\},
$$
$$
{\rm L}_k=\bigl\{(l_1,\ldots,l_k):\ l_1,\ldots,l_k=0,\ 1,\ldots,N-1;\
l_g\ne l_r\ (g\ne r);\ g, r=1,\ldots,k\bigr\},
$$

\noindent
${\rm l.i.m.}$ is a limit in the mean-square sense$,$
$i_1,\ldots,i_k=0,1,\ldots,m,$ 
random variables

\newpage
\noindent
$$
\pi_{j}^{(g,i_g)}=
\int\limits_t^T \phi_j(\tau)\int\limits_{\bf X}\varphi_g({\bf y})
\tilde\nu^{(i_g)}(d\tau,d{\bf y})
$$ 

\noindent
are independent for various
$i_g$ {\rm (}if $i_g\ne 0${\rm )} and uncorrelated for various
$j,$
$$
C_{j_k\ldots j_1}=\int\limits_{[t,T]^k}
K(t_1,\ldots,t_k)\prod_{l=1}^{k}\phi_{j_l}(t_l)dt_1\ldots dt_k
$$
is the Fourier coefficient$,$
$$
K(t_1,\ldots,t_k)=
\left\{
\begin{matrix}
\psi_1(t_1)\ldots \psi_k(t_k),\ &t_1<\ldots<t_k\cr\cr
0,\ &\hbox{otherwise}
\end{matrix}\right.,\ \  t_1,\ldots,t_k\in[t, T],\ \  k\ge 2,
$$
and 
$K(t_1)\equiv\psi_1(t_1)$ for $t_1\in[t, T].$  
}

{\bf Proof.} The scheme of the proof of Theorem 1.7
is the same with the sche\-me
of the proof of Theorem 1.1.
Some differences will take place in 
the proof of Lemmas 1.4, 1.5 (see below) and in the final part 
of the proof of Theorem 1.7.

{\bf Lemma 1.4.}\ {\it Suppose that
every $\psi_l(\tau)$ $(l=1,\ldots,k)$ is a continuous 
function at the interval
$[t, T]$ and every function $\varphi_l({\bf y})$ $(l=1,\ldots,k)$ such that
$$
\int\limits_{\bf X}\left|\varphi_l({\bf y})\right|^2
\Pi(d{\bf y})<\infty.
$$ 

Then$,$ the following equality 
\begin{equation}
\label{2000.2.11}
~~~~~P[\bar \chi^{(k)}]_{T,t}=
\hbox{\vtop{\offinterlineskip\halign{
\hfil#\hfil\cr
{\rm l.i.m.}\cr
$\stackrel{}{{}_{N\to \infty}}$\cr
}} }
\sum_{j_k=0}^{N-1}
\ldots \sum_{j_1=0}^{j_{2}-1}
\prod_{l=1}^k\int\limits_{\bf X}\chi_l(\tau_{j_l},{\bf y})
\bar \nu^{(i_l)}([\tau_{j_l},\tau_{j_l+1}),d{\bf y})
\end{equation}
is valid w. p. {\rm 1,}
where $\left\{\tau_{j}\right\}_{j=0}^{N}$ is a partition of
the interval $[t,T]$ satisfying the condition {\rm (\ref{1111}),}
$$
\bar \nu^{(i)}([\tau,s),d{\bf y})=
\left\{
\begin{matrix}
\tilde \nu^{(i)}([\tau,s),d{\bf y})\cr\cr
\nu^{(i)}([\tau,s),d{\bf y})
\end{matrix}\right. \ \ \ (i=0, 1, \ldots,m).
$$

\vspace{1mm}
\noindent
In contrast to the integral
$P[\chi^{(k)}]_{T,t}$
defined by {\rm (\ref{2000.2.1}),}
$\bar\nu^{(i_l)}(dt_l,d{\bf y}_l)$ is used
in the integral $P[\bar \chi^{(k)}]_{T,t}$
instead of 
$\tilde\nu^{(i_l)}(dt_l,d{\bf y}_l)$
$(l=1,\ldots,k).$
}

{\bf Proof.} Using the moment properties of stochastic integrals 
with respect to
the Poisson measure (see above) and the conditions of Lemma 1.4, 
it is easy to notice that the integral sum 
of the integral 
$P[\bar\chi^{(k)}]_{T,t}$ 
can be represented as a sum of the prelimit 
expression from the right-hand side of (\ref{2000.2.11}) and  
the value, which 
converges
to zero 
in the mean-square sense if 
$N\to \infty.$ Lemma 1.4 is proved.

Note that in the case when the functions 
$\psi_l(\tau)$ $(l=1,\ldots,k)$ satisfy the condition
$(\star)$ {\rm (}see Sect.~{\rm 1.1.7)}
we can suppose that among the points
$\tau_j,$ $j=0,1,\ldots,N$ there are all points of 
jumps of the functions $\psi_l(\tau)$ 
$(l=1,\ldots,k)$. Further,
we can apply the argumentation as in Sect.~1.1.7.

Let us consider the following multiple and iterated
stochastic integrals 
$$
\hbox{\vtop{\offinterlineskip\halign{
\hfil#\hfil\cr
{\rm l.i.m.}\cr
$\stackrel{}{{}_{N\to \infty}}$\cr
}} }
\sum_{j_1,\ldots,j_k=0}^{N-1}
\Phi(\tau_{j_1},\ldots,\tau_{j_k})
\prod_{l=1}^k 
\int\limits_{\bf X}\varphi_l({\bf y})
\tilde \nu^{(i_l)}([\tau_{j_l},\tau_{j_l+1}),d{\bf y})
\stackrel{\rm def}{=}P[\Phi]_{T,t}^{(k)},
$$
$$
\int\limits_t^T
\ldots\int\limits_t^{t_2}
\Phi(t_{1},\ldots,t_{k})\int\limits_{\bf X}\varphi_1({\bf y})
\tilde\nu^{(i_1)}(dt_1,d{\bf y})\ldots
\int\limits_{\bf X}\varphi_k({\bf y})\tilde\nu^{(i_k)}(dt_k,d{\bf y})
\stackrel{\rm def}{=}
$$
$$
\stackrel{\rm def}{=}
\hat P[\Phi]_{T,t}^{(k)},
$$

\vspace{2mm}
\noindent
where $\Phi(t_1,\ldots,t_k):\ [t, T]^k\to{\bf R}^1$ is a bounded nonrandom
function and
the sense of notations of the formula 
(\ref{2000.2.11}) is remaining.

Note that if the functions $\varphi_l({\bf y})$ $(l=1,\ldots,k)$
satisfy the conditions of Lemma 1.4 and the function 
$\Phi(t_1,\ldots,t_k)$
is continuous in the domain of integration, then for the integral 
$\hat P[\Phi]_{T,t}^{(k)}$ the equality similar to
{\rm (\ref{2000.2.11})}
is valid w.~p.~1.

{\bf Lemma 1.5.}\ {\it Assume that
the following representation takes place{\rm :}
$$
g_l(\tau,{\bf y})=h_l(\tau)\varphi_l({\bf y})\ \ \ (l=1,\ldots,k),
$$

\noindent
where 
the functions $h_l(\tau):$ $[t, T]\to{\bf R}^1$ 
$(l=1,\ldots,k)$ satisfy the condition
$(\star)$ {\rm (}see Sect.~{\rm 1.1.7)} and 
the
functions 
$\varphi_l({\bf y}):$ ${\bf X}\to{\bf R}^1$ 
$(l=1,\ldots,k)$ satisfy the condition
$$
\int\limits_{\bf X}\left|\varphi_l({\bf y})
\right|^p
\Pi(d{\bf y})<\infty\ \ \ \hbox{for} \ \ \ p=2^{k+1}.
$$

Then
\begin{equation}
\label{dur1}
\prod_{l=1}^k \int\limits_t^T\int\limits_{\bf X} g_l(s,{\bf y}) 
\bar\nu^{(i_l)}(ds,d{\bf y})=
P[\Phi]_{T,t}^{(k)}\ \ \ \hbox{w. p. {\rm 1}},
\end{equation}
where $i_l=0, 1, \ldots,m$ $(l=1,\ldots,k)$ and
$$
\Phi(t_1,\ldots,t_k)
=\prod\limits_{l=1}^k h_l(t_l).
$$
}

\vspace{-4mm}
{\bf Proof.} Let us introduce the following notations
$$
J[\bar g_l]_N\stackrel{\rm def}{=}\sum\limits_{j=0}^{N-1}\int\limits_{\bf X}
g_l(\tau_j,{\bf y})\bar\nu^{(i_l)}([\tau_{j},\tau_{j+1}),d{\bf y}),
$$
$$
J[\bar g_l]_{T,t}
\stackrel{\rm def}{=}\int\limits_t^T
\int\limits_{\bf X}
g_l(s,{\bf y})
\bar\nu^{(i_l)}(ds,d{\bf y}),
$$
where $\{\tau_j\}_{j=0}^N$ is a partition of the interval 
$[t,T]$ satisfying the condition (\ref{1111}).

It is easy to see that 
$$
\prod_{l=1}^k J[\bar g_l]_N-\prod_{l=1}^k J[\bar g_l]_{T,t}=
$$
$$
=
\sum_{l=1}^k \left(\prod_{q=1}^{l-1} J[\bar g_q]_{T,t}\right)
\left(J[\bar g_l]_N-
J[\bar g_l]_{T,t}\right)\left(\prod_{q=l+1}^k J[\bar g_q]_N\right).
$$

Using the Minkowski inequality and the inequality of Cauchy--Bunyakovsky
together with the estimates 
of moments of stochastic integrals with respect to the Poisson
measure and the conditions 
of Lemma 1.5, we obtain
\begin{equation}
\label{2000.4.3}
\left({\sf M}\left\{\left|\prod_{l=1}^k J[\bar g_l]_N-
\prod_{l=1}^k J[\bar g_l]_{T,t}\right|^2\right\}
\right)^{1/2}\hspace{-1mm}\le C_k
\sum_{l=1}^k
\left({\sf M}
\left\{\biggl|J[\bar g_l]_N-J[\bar g_l]_{T,t}
\biggr|^4\right\}\right)^{1/4}\hspace{-0.5mm},
\end{equation}
where $C_k<\infty.$

We have
$$
J[\bar g_l]_N-J[\bar g_l]_{T,t}
=\sum\limits_{q=0}^{N-1}J[\Delta\bar g_{l}]_{\tau_{q+1},\tau_q},
$$
where
$$
J[\Delta\bar g_{l}]_{\tau_{q+1},\tau_q}
=\int\limits_{\tau_q}^{\tau_{q+1}}\left(h_l(\tau_q)-
h_l(s)\right)\int\limits_{\bf X}
\phi_l({\bf y})
\bar\nu^{(i_l)}(ds,d{\bf y}).
$$

Let us introduce the notation
$$
h_l^{(N)}(s)=h_l(\tau_q),\ \ \ s\in [\tau_q, \tau_{q+1}), \ \ \ q=0, 1, \ldots,
N-1.
$$

Then
$$
J[\bar g_l]_N-J[\bar g_l]_{T,t}
=\sum\limits_{q=0}^{N-1}J[\Delta\bar g_{l}]_{\tau_{q+1},\tau_q}=
$$
$$
=
\int\limits_{t}^{T}\left(h_l^{(N)}(s)-
h_l(s)\right)\int\limits_{\bf X}
\phi_l({\bf y})
\bar\nu^{(i_l)}(ds,d{\bf y}).
$$

Applying the estimates (\ref{16.010}) (for $p=4$) and
(\ref{16.011}), (\ref{dur}) (for $p=2$) to
the value
$$
{\sf M}\left\{\left|
\int\limits_{t}^{T}\left(h_l^{(N)}(s)-
h_l(s)\right)\int\limits_{\bf X}
\phi_l({\bf y})
\bar\nu^{(i_l)}(ds,d{\bf y})\right|^4\right\},
$$
taking into account (\ref{2000.4.3}), the conditions of Lemma 1.5, and
the estimate
\begin{equation}
\label{dur2}
~~~~~~~~~~ \left|h_l(\tau_q)-
h_l(s)\right|<\varepsilon,\ \ \ s\in[\tau_q,\tau_{q+1}],\ \ \ q=0,1,\ldots,N-1,
\end{equation}
where $\varepsilon$ is an arbitrary small positive real number
and $|\tau_{q+1}-\tau_q|<\delta(\varepsilon)$,
we obtain that the right-hand side of (\ref{2000.4.3}) 
converges
to zero when $N\to\infty.$
Therefore, we come to 
the affirmation of Lemma 1.5. 

It should be noted that (\ref{dur2}) is valid
if the functions $h_l(s)$ are continuous at the interval
$[t, T]$, i.e. these functions are uniformly continuous at this interval.
So, $\left|h_l(\tau_q)-h_l(s)\right|<\varepsilon$
if $s\in [\tau_q, \tau_{q+1}],$ where
$|\tau_{q+1}-\tau_q|<\delta(\varepsilon),$ $q=0, 1,\ldots,N-1$
($\delta(\varepsilon)>0$ exists
for any $\varepsilon>0$ and it does not
depend on points of the interval $[t, T]$).

In the case when the functions 
$h_l(s)$ 
$(l=1,\ldots,k)$ satisfy the condition
$(\star)$ {\rm (}see Sect.~{\rm 1.1.7)}
we can suppose that among the points
$\tau_q,$ $q=0,1,\ldots,N$ there are all points of 
jumps of the functions $h_l(s)$ 
$(l=1,\ldots,k)$. Further,
we can apply the argumentation as in Sect.~1.1.7.

Obviously, if $i_l=0$ for some $l=1,\ldots,k,$ then
we also come to the 
affirmation of Lemma 1.5. Lemma 1.5 is proved.

Proving Theorem 1.7 by the scheme of the proof
of Theorem 1.1 
using Lemmas 1.4, 1.5 and 
moment properties 
of stochastic integrals with respect to the martingale Poisson measures, we obtain
$$
{\sf M}\left\{\biggl(R_{T,t}^{p_1,\ldots,p_k}\biggr)^2\right\}
\le 
C_k  \prod\limits_{l=1}^k
\int\limits_{\bf X}\varphi_l^2({\bf y})\Pi(d{\bf y})\times
$$
$$
\times
\sum_{(t_1,\ldots,t_k)}
\int\limits_{t}^{T}
\ldots
\int\limits_{t}^{t_2}
\left(K(t_1,\ldots,t_k)-
\sum_{j_1=0}^{p_1}\ldots
\sum_{j_k=0}^{p_k}
C_{j_k\ldots j_1}
\prod_{l=1}^k\phi_{j_l}(t_l)\right)^2\times
$$

\vspace{-2mm}
$$
\times
dt_1
\ldots
dt_k=
$$

\vspace{-4mm}
$$
=
C_k  \prod\limits_{l=1}^k
\int\limits_{\bf X}\varphi_l^2({\bf y})\Pi(d{\bf y})
\int\limits_{[t,T]^k}
\Biggl(K(t_1,\ldots,t_k)-
\sum_{j_1=0}^{p_1}\ldots
\sum_{j_k=0}^{p_k}
C_{j_k\ldots j_1}
\prod_{l=1}^k\phi_{j_l}(t_l)\Biggr)^2\times
$$
$$
\times
dt_1
\ldots
dt_k\le
$$

\vspace{-4mm}
$$
\le \bar{C_k}
\int\limits_{[t,T]^k}
\Biggl(K(t_1,\ldots,t_k)-
\sum_{j_1=0}^{p_1}\ldots
\sum_{j_k=0}^{p_k}
C_{j_k\ldots j_1}
\prod_{l=1}^k\phi_{j_l}(t_l)\Biggr)^2
dt_1
\ldots
dt_k\to 0
$$

\noindent
if $p_1,\ldots,p_k\to\infty,$
where constant $\bar{C_k}$ depends only on $k$ ($k$ is the multiplicity
of the iterated stochastic integral with respect 
to the martingale Poisson measures).
Moreover, $R_{T,t}^{p_1,\ldots,p_k}$ has the following form
$$
R_{T,t}^{p_1,\ldots,p_k}
=\sum_{(t_1,\ldots,t_k)}
\int\limits_{t}^{T}
\ldots
\int\limits_{t}^{t_2}
\left(K(t_1,\ldots,t_k)-
\sum_{j_1=0}^{p_1}\ldots
\sum_{j_k=0}^{p_k}
C_{j_k\ldots j_1}
\prod_{l=1}^k\phi_{j_l}(t_l)\right)\times
$$
\begin{equation}
\label{jter}
\times
\int\limits_{\bf X}\varphi_1({\bf y})
\tilde \nu^{(i_1)}(dt_1,d{\bf y})\ldots
\int\limits_{\bf X}\varphi_k({\bf y})
\tilde \nu^{(i_k)}(dt_k,d{\bf y}),
\end{equation}

\vspace{2mm}
\noindent
where permutations $(t_1,\ldots,t_k)$ when summing in (\ref{jter})
are performed only in the values
$\varphi_1({\bf y})
\tilde \nu^{(i_1)}(dt_1,d{\bf y})\ldots
\varphi_k({\bf y})
\tilde \nu^{(i_k)}(dt_k,d{\bf y}).$
At the same time, the indices near 
upper limits of integration in the iterated stochastic integrals are changed 
correspondently and if $t_r$ swapped with $t_q$ in the  
permutation $(t_1,\ldots,t_k)$, then $i_r$ swapped with $i_q$ in 
the permutation $(i_1,\ldots,i_k)$. Moreover,
$\varphi_r({\bf y})$ swapped with $\varphi_q({\bf y})$
in the permutation $(\varphi_1({\bf y}),\ldots,\varphi_k({\bf y}))$.
Theorem 1.7 is proved.

Let us consider the application of Theorem 1.7.
Let $i_1\ne i_2$ and $i_1,i_2=1,\ldots,m.$ 
Using Theorem 1.7 and
the system of Legendre polynomials, we obtain
$$
\int\limits_t^T\int\limits_{\bf X}\varphi_2({\bf y}_2)
\int\limits_{t}^{t_2}\int\limits_{\bf X}\varphi_1({\bf y}_1)
\tilde\nu^{(i_1)}(dt_1,d{\bf y}_1)
\tilde\nu^{(i_2)}(dt_2,d{\bf y}_2)=
$$
$$
=\frac{T-t}{2}\Biggl(
\pi_{0}^{(1,i_1)}\pi_{0}^{(2,i_2)}
+\sum_{i=1}^{\infty}\frac{1}{\sqrt{4i^2-1}}
\left(\pi_{i-1}^{(1,i_1)}
\pi_{i}^{(2,i_2)}-\pi_{i}^{(1,i_1)}\pi_{i-1}^{(2,i_2)}
\right)\Biggr),
$$
$$
\int\limits_t^T\int\limits_{\bf X}\varphi_1({\bf y}_1)
\tilde\nu^{(i_1)}(dt_1,d{\bf y}_1)=\sqrt{T-t}\pi_{0}^{(1,i_1)},
$$

\noindent
where 
$$
\pi_{j}^{(l,i_l)}=
\int\limits_t^T\phi_j(\tau)\int\limits_{\bf X}\varphi_l({\bf y})
\tilde\nu^{(i_l)}(d\tau,d{\bf y})\ \ \ (l=1, 2)
$$

\noindent
and $\{\phi_j(\tau)\}_{j=0}^{\infty}$ is a complete orthonormal
system of Legendre polynomials in the space
$L_2([t, T])$.

\section{Expansion 
of Iterated Stochastic Integrals with Respect to Martingales
Based on Generalized Multiple Fourier Series}

\subsection{Stochastic Integral with Respect to 
Martingale}

Let $(\Omega, {\rm F},{\sf P})$ be a fixed probability
space, let $\{{\rm F}_t,$ $t\in[0,T]\}$
be a non-decreasing family of 
$\sigma$-algebras ${\rm F}_t\subset{\rm F}$, and let
${\rm M}_2(\rho,[0,T])$ be a class 
of ${\rm F}_t$-measurable for each $t\in[0, T]$
martingales $M_t$ satisfying the conditions 
\begin{equation}
\label{rikos1}
{\sf M}\Bigl\{\left(M_s-M_t\right)^2\Bigr\}=\int\limits_t^s
\rho(\tau)d\tau,
\end{equation}
$$
~~~{\sf M}\Bigl\{\left|M_s-M_t\right|^p\Bigr\}\le C_p|s-t|,\ \ \ p=3, 4,\ldots,
$$

\vspace{1mm}
\noindent
where $0\le t< s\le T,$ $\rho(\tau)$ is a non-negative and continuously 
differentiable nonrandom function at the interval
$[0, T]$, $C_p<\infty$ is a constant.

Let us define the class $H_2(\rho,[0, T])$ 
of stochastic 
processes $\xi_t,$ $t\in[0, T],$ which are 
${\rm F}_t$-measurable for all $t\in[0, T]$ 
and satisfy the condition
$$
\int\limits_0^T{\sf M}\left\{\left|\xi_t\right|^2\right\}\rho(t)dt<\infty.
$$

For any partition 
$\left\{\tau_j^{(N)}\right\}_{j=0}^{N}$ of
the interval $[0,T]$ such that
\begin{equation}
\label{w11}
0=\tau_0^{(N)}<\tau_1^{(N)}<\ldots <\tau_N^{(N)}=T,\ \ \
\max\limits_{0\le j\le N-1}\left|\tau_{j+1}^{(N)}-\tau_j^{(N)}\right|\to\  
0\ \ 
\hbox{if}\ \ N\to \infty
\end{equation}
we will define the sequence of step functions 
$\xi^{(N)}(t,\omega)$ by the following relation
$$
\xi^{(N)}(t,\omega)=\xi_j\left(\omega\right)\ \ \ 
\hbox{w. p. 1}\ \ \ \hbox{for}\ \ \ t\in\left[\tau_j^{(N)},
\tau_{j+1}^{(N)}\right),
$$
where $\xi^{(N)}(t,\omega)\in H_2(\rho,[0, T]),$ $j=0, 1,\ldots,N-1,$\ $N=1, 2,\ldots$

Let us define the stochastic integral with respect to martingale 
from the process $\xi(t,\omega)\in H_2(\rho,[0,T])$ as the following 
mean-square limit \cite{Gih1}
\begin{equation}
\label{fff}
\hbox{\vtop{\offinterlineskip\halign{
\hfil#\hfil\cr
{\rm l.i.m.}\cr
$\stackrel{}{{}_{N\to \infty}}$\cr
}} }\sum_{j=0}^{N-1}\xi^{(N)}\left(\tau_j^{(N)},\omega\right)
\biggl(M\left(\tau_{j+1}^{(N)},\omega\right)-
M\left(\tau_j^{(N)},\omega\right)\biggr)
\stackrel{\rm def}{=}\int\limits_0^T\xi_\tau dM_\tau,
\end{equation}
where $\xi^{(N)}(t,\omega)$ is any step function
from the class $H_2(\rho,[0,T])$,
which converges
to the function $\xi(t,\omega)$
in the following sense
$$
\hbox{\vtop{\offinterlineskip\halign{
\hfil#\hfil\cr
{\rm lim}\cr
$\stackrel{}{{}_{N\to \infty}}$\cr
}} }\int\limits_0^T{\sf M}\left\{\left|\xi^{(N)}(t,\omega)-
\xi(t,\omega)\right|^2\right\}\rho(t)dt=0.
$$

It is well known  \cite{Gih1} that the stochastic integral (\ref{fff})
exists, it does not depend on selection 
of the sequence 
$\xi^{(N)}(t,\omega)$ and 
it satisfies w.~p.~1
to the following properties 
$$
{\sf M}\left\{\int\limits_0^T
\xi_t dM_t\biggl|\biggr.{\rm F}_0\right\}=0,
$$
$$
{\sf M}\left\{\left|\int\limits_0^T
\xi_t dM_t\right|^2\biggl|\biggr.{\rm F}_0\right\}=
{\sf M}\left\{\int\limits_0^T\xi_t^2\rho(t) dt\biggl|\biggr.{\rm F}_0\right\},
$$
$$
\int\limits_0^T(\alpha\xi_t+\beta\psi_t)dM_t=
\alpha\int\limits_0^T\xi_t dM_t+\beta
\int\limits_0^T\psi_t dM_t,
$$
where $\xi_t,$ $\psi_t\in H_2(\rho,[0, T]),$\ 
$\alpha, \beta\in{\bf R}^1.$

\subsection{Expansion 
of Iterated Stochastic Integrals with Respect to Martingales}

Let $Q_4(\rho,[0,T])$ be the class 
of martingales $M_t,$ $t\in[0,T],$ 
which satisfy the following conditions:

1. $M_t,$ $t\in[0,T]$ belongs to the class
${\rm M}_2(\rho,[0,T]).$ 

2. For some 
$\alpha>0$ the following estimate is correct
\begin{equation}
\label{durra1}
{\sf M}\left\{\left|\int\limits_{t}^{\tau}
g(s) dM_s\right|^4\right\}\le K_4 \int\limits_{t}^{\tau}
|g(s)|^{\alpha}ds,
\end{equation}
where $0\le t<\tau\le T,$
$g(s)$ is a bounded nonrandom function 
at the interval $[0, T]$,
$K_4<\infty$ is a constant.

Let $G_n(\rho,[0,T])$ be the class 
of martingales $M_t,$ $t\in[0,T],$ 
which satisfy the following conditions:

1. $M_t,$ $t\in[0,T]$ belongs to the class
${\rm M}_2(\rho,[0,T]).$ 

2. The following estimate is correct
$$
{\sf M}\left\{\left|\int\limits_{t}^{\tau}
g(s) dM_s\right|^n\right\} < \infty,
$$
where $0\le t<\tau\le T,$\ $n\in {\bf N},$\ $g(s)$ 
is the same function as in the definition 
of the class $Q_4(\rho,[0,T])$.

Let us 
remind that if $\left(\xi_t\right)^n\in H_2(\rho,[0,T])$
with $\rho(t)\equiv 1,$ then the following estimate is correct
\cite{Gih1}
\begin{equation}
\label{1.5aa}
~~~~ {\sf M}\left\{\left|\int\limits_{t}^{\tau}
\xi_s ds\right|^{2n}\right\}\le (\tau-t)^{2n-1}
\int\limits_{t}^{\tau}
{\sf M}\left\{\left|\xi_s\right|^{2n}\right\}ds,\ \ \ 0\le t<\tau\le T.
\end{equation}

Let us consider the iterated stochastic integral with respect to martingales
\begin{equation}
\label{mart}
~~~~~~~~~ J[\psi^{(k)}]_{T,t}^M=
\int\limits_t^T\psi_k(t_k)\ldots
\int\limits_t^{t_2}\psi_1(t_1)
dM_{t_1}^{(1,i_1)}\ldots dM_{t_k}^{(k,i_k)},
\end{equation}
where $i_1,\ldots,i_k=0, 1,\ldots,m,$
every $\psi_l(\tau)$ $(l=1,\ldots, k)$ is a continuous 
nonrandom function
at the interval $[t, T],$
$M^{(r,i)}_s$ $(r=1,\ldots,k,\ i=1,\ldots,m)$ are independent martingales 
for various $i=1,\ldots,m,$ 
$M_{s}^{(r,0)}\stackrel{\rm def}{=}s.$

Now we can formulate the following theorem.

{\bf Theorem 1.8}\ \cite{1}-\cite{12aa}, \cite{arxiv-13}.
{\it Suppose that the following 
conditions are hold{\rm :}

{\rm 1}.\ Every $\psi_l(\tau)\ (l=1,\ldots,k)$ is a 
continuous nonrandom function at
the interval $[t, T]$.

{\rm 2}.\ $\{\phi_j(x)\}_{j=0}^{\infty}$ is a complete orthonormal
system of functions in the space 
$L_2([t,T]),$ 
each function $\phi_j(x)$ of which for finite $j$ satisfies the condition 
$(\star)$ {\rm(}see Sect.~{\rm 1.1.7)}.

{\rm 3}. $M_{s}^{(l,i_l)}\in Q_4(\rho,[t,T]),$ $G_n(\rho,[t,T])$
with $n=2^{k+1},$  
$i_l=1,\ldots,m,$ $l=1,\ldots,k$ $(k\in {\bf N}).$

Then$,$ for the iterated stochastic integral 
$J[\psi^{(k)}]_{T,t}^M$ with respect to martingales 
defined by
{\rm (\ref{mart})}
the following expansion 
$$
J[\psi^{(k)}]_{T,t}^M=
\hbox{\vtop{\offinterlineskip\halign{
\hfil#\hfil\cr
{\rm l.i.m.}\cr
$\stackrel{}{{}_{p_1,\ldots,p_k\to \infty}}$\cr
}} }\sum_{j_1=0}^{p_1}\ldots\sum_{j_k=0}^{p_k}
C_{j_k\ldots j_1}\Biggl(
\prod_{l=1}^k\xi_{j_l}^{(l,i_l)}-
\Biggr.
$$

$$
-\Biggl.
\hbox{\vtop{\offinterlineskip\halign{
\hfil#\hfil\cr
{\rm l.i.m.}\cr
$\stackrel{}{{}_{N\to \infty}}$\cr
}} }\sum_{(l_1,\ldots,l_k)\in {\rm G}_k}
\phi_{j_{1}}(\tau_{l_1})
\Delta{M}_{\tau_{l_1}}^{(1,i_1)}\ldots
\phi_{j_{k}}(\tau_{l_k})
\Delta{M}_{\tau_{l_k}}^{(k,i_k)}\Biggr)
$$

\vspace{2mm}
\noindent
that
converges in the mean-square sense 
is valid$,$ where $i_1,\ldots,i_k=0,1,\ldots,m,$
$\left\{\tau_{j}\right\}_{j=0}^{N}$ is a partition of
the interval $[t, T]$ satisfying the condition similar to
{\rm (\ref{w11}),}
$\Delta{M}_{\tau_{j}}^{(r,i)}=
M_{\tau_{j+1}}^{(r,i)}-M_{\tau_{j}}^{(r,i)}$
$(i=0, 1,\ldots,m,$\ \ $r=1,\ldots,k),$
$$
{\rm G}_k={\rm H}_k\backslash{\rm L}_k,\ \ 
{\rm H}_k=\{(l_1,\ldots,l_k):\ l_1,\ldots,l_k=0,\ 1,\ldots,N-1\},
$$
$$
{\rm L}_k=\bigl\{(l_1,\ldots,l_k):\ l_1,\ldots,l_k=0,\ 1,\ldots,N-1;\
l_g\ne l_r\ (g\ne r);\ g, r=1,\ldots,k\bigr\},
$$

\vspace{1mm}
\noindent
${\rm l.i.m.}$ is a limit in the mean-square sense$,$
$$
\xi_{j}^{(l,i_l)}=
\int\limits_t^T \phi_{j}(s) d{M}_s^{(l,i_l)}
$$
are independent for various
$i_l$ {\rm (}if $i_l\ne 0${\rm )}
and uncorrelated for various $j$
{\rm (}if $\rho(\tau)$ is a constant$,$ $i_l\ne 0${\rm )} random variables$,$
$$
C_{j_k\ldots j_1}=\int\limits_{[t,T]^k}
K(t_1,\ldots,t_k)\prod_{l=1}^{k}\phi_{j_l}(t_l)dt_1\ldots dt_k
$$
is the Fourier coefficient$,$
$$
K(t_1,\ldots,t_k)=
\left\{
\begin{matrix}
\psi_1(t_1)\ldots \psi_k(t_k),\ &t_1<\ldots<t_k\cr\cr
0,\ &\hbox{otherwise}
\end{matrix}\right.,\ \  t_1,\ldots,t_k\in[t, T],\ \  k\ge 2,
$$

\vspace{1mm}
\noindent
and 
$K(t_1)\equiv\psi_1(t_1)$ for $t_1\in[t, T].$  
}

{\bf Remark 1.4.}\ {\it Note that from 
Theorem {\rm 1.8} for the case $\rho(\tau)\equiv 1$ we obtain the variant of
Theorem {\rm 1.1}.}

{\bf Proof.} The proof of Theorem 1.8 is similar to
the proof of Theorem 1.1.
Some differences will take place in the
proof of Lemmas 1.6, 1.7 (see below) and in the final part 
of the proof of Theorem 1.8.

{\bf Lemma 1.6.} {\it Assume that\ \
$M_{s}^{(r,i)}\in{\rm M}_2(\rho,[t,T])$ $(i=1,\ldots,m),$ $M_{s}^{(r,0)}=s$\\  
$(r=1,\ldots,k),$ and
every $\psi_l(\tau)$ $(l=1,\ldots, k)$ is a continuous nonrandom
function at the interval $[t,T]$.
Then
\begin{equation}
\label{1.9aa}
~~~~~~~~~~~~~ J[\psi^{(k)}]_{T,t}^M=
\hbox{\vtop{\offinterlineskip\halign{
\hfil#\hfil\cr
{\rm l.i.m.}\cr
$\stackrel{}{{}_{N\to \infty}}$\cr
}} }
\sum_{j_k=0}^{N-1}
\ldots \sum_{j_1=0}^{j_{2}-1}
\prod_{l=1}^k \psi_{l}(\tau_{j_l})\Delta M_{\tau_{j_l}}^{(l,i_l)}\ \ \
\hbox{w. p. {\rm 1}},
\end{equation}
where 
$\{\tau_j\}_{j=0}^N$ is a partition of the interval
$[t,T]$ satisfying the condition similar to {\rm(\ref{w11})},\ \
$i_l=0,1,\ldots,m,$\ \ $l=1,\ldots,k;$ another notations are the same as in Theorem {\rm 1.8.}
}

\par
{\bf Proof.} According to the properties of the stochastic integral 
with respect to martingales, we have \cite{Gih1}
\begin{equation}
\label{u1}
{\sf M}\left\{\left(\int\limits_{t}^{\tau}
\xi_s dM_s^{(l,i_l)}\right)^2\right\}=
\int\limits_{t}^{\tau}{\sf M}\left\{\left|\xi_s\right|^2\right\}\rho(s)ds,
\end{equation}
\begin{equation}
\label{u2}
{\sf M}\left\{\left(\int\limits_{t}^{\tau}
\xi_s ds\right)^2\right\}\le (\tau-t)
\int\limits_{t}^{\tau}{\sf M}\left\{\left|\xi_s\right|^2\right\}ds,
\end{equation}

\noindent
where 
$\xi_s\in H_2(\rho,[0,T]),$\ $0\le t<\tau\le T,$
$i_l=1,\ldots,m,$ $l=1,\ldots,k$. 
Then the integral sum 
for the integral 
$J[\psi^{(k)}]_{T,t}^M$
under
the conditions of Lemma 1.6 can be represented as a sum of the prelimit 
expression from the right-hand side of (\ref{1.9aa}) and the
value, which 
converges
to zero 
in the mean-square sense if 
$N\to \infty.$ 
More detailed proof of the similar lemma for the case $\rho(\tau)\equiv 1$ 
can be found in Sect.~1.1.3  (see Lemma 1.1).

In the case when the functions 
$\psi_l(\tau)$ $(l=1,\ldots,k)$ satisfy the condition
$(\star)$ {\rm (}see Sect.~{\rm 1.1.7)}
we can suppose that among the points
$\tau_j,$ $j=0,1,\ldots,N$ there are all points of 
jumps of the functions $\psi_l(\tau)$ 
$(l=1,\ldots,k).$ So, we can
apply the argumentation as in Sect.~1.1.7.

Let us define the following multiple stochastic integral

\vspace{-2mm}           
\begin{equation}
\label{777666}
\hbox{\vtop{\offinterlineskip\halign{
\hfil#\hfil\cr
{\rm l.i.m.}\cr
$\stackrel{}{{}_{N\to \infty}}$\cr
}} }
\sum_{j_1,\ldots,j_k=0}^{N-1}
\Phi(\tau_{j_1},\ldots,\tau_{j_k})
\prod_{l=1}^k 
\Delta M_{\tau_{j_l}}^{(l,i_l)}
\stackrel{\rm def}{=}I[\Phi]_{T,t}^{(k)},
\end{equation}

\vspace{2mm}
\noindent
where $\{\tau_j\}_{j=0}^N$ is a partition of the interval 
$[t,T]$ satisfying the condition similar to (\ref{w11})
and
$\Phi(t_1,\ldots,t_k):\ [t, T]^k\to{\bf R}^1$ is a bounded nonrandom
function; another notations are the same as in Theorem 1.8.

{\bf Lemma 1.7.} {\it Let\ \
$M_s^{(l,i_l)}\in Q_4(\rho,[t, T]),$ $G_n(\rho,[t,T])$ with
$n=2^{k+1},$\ \ $k\in{\bf N}$\\  $(i_l=1,\ldots,m,$\ \
$l=1,\ldots,k)$
and the functions $g_1(s),\ldots, g_k(s)$ satisfy the condition 
$(\star)$ {\rm(}see Sect.~{\rm 1.1.7)}.
Then
$$
\prod_{l=1}^k \int\limits_t^T g_l(s) 
dM_s^{(l,i_l)}=
I[\Phi]_{T,t}^{(k)}\ \ \ \hbox{w. p. {\rm 1}},
$$
where\  $i_l=0, 1,\ldots,m,$\ \ $l=1,\ldots, k,$
$$
\Phi(t_1,\ldots,t_k)=\prod\limits_{l=1}^k g_l(t_l).
$$}

\vspace{-6mm}

{\bf Proof.} Let us denote

\vspace{-3mm}
$$
J[g_l]_N\stackrel{\rm def}{=}\sum\limits_{j=0}^{N-1}
g_l(\tau_j)\Delta M_{\tau_j}^{(l,i_l)},\ \ \
J[g_l]_{T,t}\stackrel{\rm def}{=}\int\limits_t^T g_l(s)
dM_s^{(l,i_l)},
$$

\vspace{1mm}
\noindent
where $\{\tau_j\}_{j=0}^N$ is a partition of the interval 
$[t,T]$ satisfying the condition similar to (\ref{w11}).

Note that
$$
\prod_{l=1}^k J[g_l]_N-\prod_{l=1}^k J[g_l]_{T,t}
=
$$

\vspace{-2mm}
$$
=\sum_{l=1}^k \left(\prod_{q=1}^{l-1} J[g_q]_{T,t}\right)
\left(J[g_l]_N-
J[g_l]_{T,t}\right)\left(\prod_{q=l+1}^k J[g_q]_N\right).
$$

\vspace{2mm}

Using the Minkowski inequality and the inequality of Cauchy-Bu\-nya\-kov\-sky
as well as the conditions of Lemma 1.7, we obtain
$$
\left({\sf M}\left\{\left|\prod_{l=1}^k J[g_l]_N-
\prod_{l=1}^k J[g_l]_{T,t}\right|^2\right\}
\right)^{1/2}\le 
$$
\begin{equation}
\label{2000.4.300}
\le C_k
\sum_{l=1}^k
\left({\sf M}
\left\{\biggl|J[g_l]_N-J[g_l]_{T,t}
\biggr|^4\right\}\right)^{1/4},
\end{equation}

\noindent
where $C_k<\infty$ is a constant.

We have
$$
J[g_l]_N-J[g_l]_{T,t}
=\sum\limits_{q=0}^{N-1}J[\Delta g_{l}]_{\tau_{q+1},\tau_q},
$$
$$
J[\Delta g_{l}]_{\tau_{q+1},\tau_q}
=\int\limits_{\tau_q}^{\tau_{q+1}}\left(
g_l(\tau_q)-g_l(s)\right)
dM_s^{(l,i_l)}.
$$

Let us introduce the notation
$$
g_l^{(N)}(s)=g_l(\tau_q),\ \ \ s\in [\tau_q, \tau_{q+1}), 
\ \ \ q=0, 1, \ldots,N-1.
$$

Then
$$
J[g_l]_N-J[g_l]_{T,t}
=\sum\limits_{q=0}^{N-1}J[\Delta g_{l}]_{\tau_{q+1},\tau_q}=
$$
$$
=
\int\limits_{t}^{T}\left(g_l^{(N)}(s)-
g_l(s)\right)dM_s^{(l,i_l)}.
$$

Applying the estimate (\ref{durra1}), 
we obtain
$$
{\sf M}\left\{\left|
\int\limits_{t}^{T}\left(g_l^{(N)}(s)-
g_l(s)\right)dM_s^{(l,i_l)}\right|^4\right\}\le
$$
$$
\le K_4\int\limits_t^T\left|
g_l^{(N)}(s)-
g_l(s)\right|^{\alpha}ds=
$$
$$
=K_4\sum_{q=0}^{N-1}\int\limits_{\tau_q}^{\tau_{q+1}}\left|
g_l(\tau_q)-
g_l(s)\right|^{\alpha}ds< 
$$
\begin{equation}
\label{durra2}
< K_4 \varepsilon^{\alpha}\ 
\sum_{q=0}^{N-1}\left(\tau_{q+1}-\tau_q\right)
=K_4 \varepsilon^{\alpha}(T-t).
\end{equation}

Note that we used the estimate
\begin{equation}
\label{durra3}
~~~~~~~~~~ \left|g_l(\tau_q)-
g_l(s)\right|<\varepsilon,\ \ \ s\in[\tau_q,\tau_{q+1}],\ \ \ 
q=0,1,\ldots,N-1
\end{equation}
to derive (\ref{durra2}),
where $|\tau_{q+1}-\tau_q|<\delta(\varepsilon)$ and
$\varepsilon$ is an arbitrary small positive real number.

The inequality (\ref{durra3}) is valid
if the functions $g_l(s)$ are continuous at the interval
$[t, T]$, i.e. these functions are uniformly continuous at this interval.
So, $\left|g_l(\tau_q)-g_l(s)\right|<\varepsilon$
if $s\in [\tau_q, \tau_{q+1}],$ where
$|\tau_{q+1}-\tau_q|<\delta(\varepsilon),$ $q=0, 1,\ldots,N-1$
($\delta(\varepsilon)>0$ exists
for any $\varepsilon>0$ and it does not
depend on points of the interval $[t, T]$).

Thus, taking into account (\ref{durra2}),
we obtain that the right-hand side of (\ref{2000.4.300}) 
converges
to zero when $N\to\infty.$
Hence, we come to 
the affirmation of Lemma 1.7.

In the case when the functions 
$g_l(s)$ 
$(l=1,\ldots,k)$ satisfy the condition
$(\star)$ {\rm (}see Sect.~{\rm 1.1.7)}
we can suppose that among the points
$\tau_q,$ $q=0,1,\ldots,N$ there are all points of 
jumps of the functions $g_l(s)$ 
$(l=1,\ldots,k)$. So, we
can apply the argumentation as in Sect.~1.1.7.

Obviously if $i_l=0$ for some $l=1,\ldots,k,$ then
we also come to the 
affirmation of Lemma 1.7.
Lemma 1.7 is proved.

Proving Theorem 1.8 similar to the proof
of Theorem 1.1 
using Lemmas 1.6, 1.7 and
moment properties
of stochastic integrals with respect to martingales 
(see (\ref{u1}), (\ref{u2})), we obtain
$$
{\sf M}\left\{\left(R_{T,t}^{p_1,\ldots,p_k}\right)^2\right\}
\le 
$$

\vspace{-5mm}
$$
\le
C_k
\sum_{(t_1,\ldots,t_k)}
\int\limits_{t}^{T}
\ldots
\int\limits_{t}^{t_2}
\left(K(t_1,\ldots,t_k)-
\sum_{j_1=0}^{p_1}\ldots
\sum_{j_k=0}^{p_k}
C_{j_k\ldots j_1}
\prod_{l=1}^k\phi_{j_l}(t_l)\right)^2\times
$$

\vspace{-1mm}
\begin{equation}
\label{z2}
\times{\tilde \rho}_1(t_1)dt_1
\ldots
{\tilde \rho}_k(t_k)dt_k\le
\end{equation}

\vspace{-3mm}
$$
\le
\bar{C_k}
\hspace{-1.5mm}
\sum_{(t_1,\ldots,t_k)}
\int\limits_{t}^{T}
\ldots
\int\limits_{t}^{t_2}
\left(K(t_1,\ldots,t_k)-
\sum_{j_1=0}^{p_1}\ldots
\sum_{j_k=0}^{p_k}
C_{j_k\ldots j_1}
\prod_{l=1}^k\phi_{j_l}(t_l)\right)^2\hspace{-1.5mm}dt_1\ldots dt_k=
$$

$$
=\bar{C_k}
\int\limits_{[t,T]^k}
\Biggl(K(t_1,\ldots,t_k)-
\sum_{j_1=0}^{p_1}\ldots
\sum_{j_k=0}^{p_k}
C_{j_k\ldots j_1}
\prod_{l=1}^k\phi_{j_l}(t_l)\Biggr)^2
dt_1
\ldots
dt_k \to 0
$$

\vspace{1mm}
\noindent
if $p_1,\ldots,p_k\to\infty,$
where constant $\bar{C_k}$ depends only on $k$ ($k$ is the multiplicity of 
the iterated stochastic integral with respect to martingales) and  
${\tilde \rho}_l(s)\equiv\rho(s)$ or ${\tilde \rho}_l(s)\equiv 1$
$(l=1,\ldots,k)$. Moreover, $R_{T,t}^{p_1,\ldots,p_k}$ has the following 
form
$$
R_{T,t}^{p_1,\ldots,p_k}
=\sum_{(t_1,\ldots,t_k)}
\int\limits_{t}^{T}
\ldots
\int\limits_{t}^{t_2}
\left(K(t_1,\ldots,t_k)-
\sum_{j_1=0}^{p_1}\ldots
\sum_{j_k=0}^{p_k}
C_{j_k\ldots j_1}
\prod_{l=1}^k\phi_{j_l}(t_l)\right)\times
$$
\begin{equation}
\label{jter1}
\times
dM_{t_1}^{(1,i_1)}\ldots
dM_{t_k}^{(k,i_k)},
\end{equation}

\vspace{3mm}
\noindent
where permutations $(t_1,\ldots,t_k)$ when summing in (\ref{jter1})
are performed only in the values
$dM_{t_1}^{(1,i_1)}\ldots
dM_{t_k}^{(k,i_k)}.$
At the same time the indices near 
upper limits of integration in the iterated stochastic integrals are changed 
correspondently and if $t_r$ swapped with $t_q$ in the  
permutation $(t_1,\ldots,t_k)$, then $i_r$ swapped with $i_q$ in 
the permutation $(i_1,\ldots,i_k)$. Moreover,
$r$ swapped with $q$
in the permutation $(1,\ldots,k)$.
Theorem 1.8 is proved.

\section{One Modification of Theorems 1.5 and 1.8}

\subsection{Expansion 
of Iterated Stochastic Integrals with Respect to Martingales
Based on Generalized Multiple Fourier Series. The Case
$\rho(x)/r(x)<\infty$}

Let us compare the expressions (\ref{obana1eee}) and (\ref{z2}). 
If we suppose that $r(x)\ge 0$ and 
$$
\frac{\rho(x)}{r(x)}\le C<\infty,
$$  
where $\rho(x)$ as in (\ref{rikos1}), then 
$$
\int\limits_{[t,T]^k}
\left(K(t_1,\ldots,t_k)
-\sum_{j_1=0}^{p_1}\ldots
\sum_{j_k=0}^{p_k}
{\tilde C}_{j_k\ldots j_1}
\prod_{l=1}^k\Psi_{j_l}(t_l)\right)^2\times
$$

\vspace{-2mm}
$$
\times
\rho(t_1)dt_1
\ldots
\rho(t_k)dt_k=
$$

\vspace{-2mm}
$$
=
\int\limits_{[t,T]^k}
\left(K(t_1,\ldots,t_k)
-\sum_{j_1=0}^{p_1}\ldots
\sum_{j_k=0}^{p_k}
{\tilde C}_{j_k\ldots j_1}
\prod_{l=1}^k\Psi_{j_l}(t_l)\right)^2\times
$$
$$
\times
\frac{\rho(t_1)}{r(t_1)}r(t_1)dt_1
\ldots
\frac{\rho(t_k)}{r(t_k)}r(t_k)dt_k\le
$$

\vspace{-2mm}
$$
\le
C_k'\int\limits_{[t,T]^k}
\left(K(t_1,\ldots,t_k)
-\sum_{j_1=0}^{p_1}\ldots
\sum_{j_k=0}^{p_k}
{\tilde C}_{j_k\ldots j_1}
\prod_{l=1}^k\Psi_{j_l}(t_l)\right)^2\times
$$

\vspace{-2mm}
$$
\times
\left(\prod\limits_{l=1}^k r(t_l)\right)dt_1
\ldots dt_k \to 0
$$

\noindent
if $p_1,\ldots,p_k\to \infty$ (see (\ref{arch1})),
where $C_k'$ is a constant,
$\{\Psi_j(x)\}_{j=0}^{\infty}$ is a complete orthonormal 
with weight $r(x)\ge 0$ 
system of functions in the space $L_2([t, T]),$
and the Fourier coefficient 
${\tilde C}_{j_k\ldots j_1}$
has the form (\ref{koef}).

So, we obtain the following modification of Theorems 1.5 and 1.8.

{\bf Theorem 1.9}\ \cite{12}-\cite{12aa}, \cite{arxiv-13}.
{\it Suppose that the following 
conditions are fulfilled{\rm :}

{\rm 1}. Every $\psi_l(\tau)\ (l=1,\ldots,k)$ is a
continuous nonrandom function at
the interval $[t, T]$.

{\rm 2}. $\{\Psi_j(x)\}_{j=0}^{\infty}$ is a complete orthonormal 
with weight $r(x)\ge 0$ 
system of functions in the space $L_2([t,T]),$ each function 
$\Psi_j(x)$ of which 
for finite $j$ satisfies the condition $(\star)$ {\rm (}see Sect.~{\rm 1.1.7)}.
Moreover$,$
$$
\frac{\rho(x)}{r(x)}\le C<\infty.
$$  

{\rm 3}. $M_{s}^{(l,i_l)}\in Q_4(\rho,[t,T]),$ $G_n(\rho,[t,T])$ with
$n=2^{k+1},$  
$i_l=1,\ldots,m,$ $l=1,\ldots,k$ $(k\in {\bf N}).$

Then$,$ for the iterated stochastic integral 
$J[\psi^{(k)}]_{T,t}^M$ with respect to martingales 
defined by {\rm (\ref{mart})}
the following expansion 
$$
J[\psi^{(k)}]_{T,t}^M=
\hbox{\vtop{\offinterlineskip\halign{
\hfil#\hfil\cr
{\rm l.i.m.}\cr
$\stackrel{}{{}_{p_1,\ldots,p_k\to \infty}}$\cr
}} }\sum_{j_1=0}^{p_1}\ldots\sum_{j_k=0}^{p_k}
{\tilde C}_{j_k\ldots j_1}\Biggl(
\prod_{l=1}^k\xi_{j_l}^{(l,i_l)}-
\Biggr.
$$

$$
-\Biggl.
\hbox{\vtop{\offinterlineskip\halign{
\hfil#\hfil\cr
{\rm l.i.m.}\cr
$\stackrel{}{{}_{N\to \infty}}$\cr
}} }\sum_{(l_1,\ldots,l_k)\in {\rm G}_k}
\Psi_{j_{1}}(\tau_{l_1})
\Delta{M}_{\tau_{l_1}}^{(1,i_1)}\ldots
\Psi_{j_{k}}(\tau_{l_k})
\Delta{M}_{\tau_{l_k}}^{(k,i_k)}\Biggr)
$$

\vspace{3mm}
\noindent
that converges in the mean-square sense
is valid$,$ where $i_1,\ldots,i_k=1,\ldots,m,$
$\left\{\tau_{j}\right\}_{j=0}^{N}$ is a partition of
the interval $[t, T]$ satisfying the condition similar to
{\rm (\ref{w11}),}
$\Delta{M}_{\tau_{j}}^{(r,i)}=
M_{\tau_{j+1}}^{(r,i)}-M_{\tau_{j}}^{(r,i)}\
(i=1,\ldots,m,$\ \ $r=1,\ldots,k),$

\vspace{-5mm}
$$
{\rm G}_k={\rm H}_k\backslash{\rm L}_k,\ \ 
{\rm H}_k=\{(l_1,\ldots,l_k):\ l_1,\ldots,l_k=0,\ 1,\ldots,N-1\},
$$
$$
{\rm L}_k=\bigl\{(l_1,\ldots,l_k):\ l_1,\ldots,l_k=0,\ 1,\ldots,N-1;\
l_g\ne l_r\ (g\ne r);\ g, r=1,\ldots,k\bigr\},
$$

\vspace{4mm}
\noindent
${\rm l.i.m.}$ is a limit in the mean-square sense$,$
$$
\xi_{j}^{(l,i_l)}=
\int\limits_t^T \Psi_{j}(s) d{M}_s^{(l,i_l)}
$$
are independent for various 
$i_l=1,\ldots,m$ $(l=1,\ldots,k)$
and uncorrelated for various $j$
$\left(\hbox{if}\ \rho(x)\equiv r(x)\right)$ random variables$,$

\vspace{-1mm}
$$
{\tilde C}_{j_k\ldots j_1}=\int\limits_{[t,T]^k}
K(t_1,\ldots,t_k)
\prod_{l=1}^{k}\biggl(\Psi_{j_l}(t_l)r(t_l)\biggr)dt_1\ldots dt_k
$$

\newpage
\noindent
is the Fourier coefficient$,$
$$
K(t_1,\ldots,t_k)=
\left\{
\begin{matrix}
\psi_1(t_1)\ldots \psi_k(t_k),\ &t_1<\ldots<t_k\cr\cr
0,\ &\hbox{otherwise}
\end{matrix}\right.,\ \  t_1,\ldots,t_k\in[t, T],\ \  k\ge 2,
$$

\vspace{1mm}
\noindent
and 
$K(t_1)\equiv\psi_1(t_1)$ for $t_1\in[t, T].$  
}

{\bf Remark 1.5.}\ {\it Note that if $\rho(x), r(x)\equiv 1$ in 
Theorem {\rm 1.9,}
then we obtain the variant of Theorem {\rm 1.1}.}

\subsection{Example on Application of Theorem 1.9 and the System
of Bessel Functions}

Let us consider the following boundary-value problem

\vspace{-5mm}
$$
\left(p(x)\Phi'(x)\right)'+q(x)\Phi(x)=-\lambda r(x)\Phi(x),
$$
$$
\alpha\Phi(a)+\beta\Phi'(a)=0,\ \ \ 
\gamma\Phi(b)+\delta\Phi'(b)=0,
$$

\vspace{1mm}
\noindent
where the functions $p(x)$, $q(x)$, $r(x)$ satisfy 
the well known conditions and 
$\alpha,$ $\beta,$ $\gamma,$ $\delta,$ $\lambda$ are real numbers.

It is well known (Steklov V.A.) that 
the eigenfunctions
$\Phi_0(x),$ $\Phi_1(x),$ $\ldots $ of this boundary-value problem 
form a complete
orthonormal with weight $r(x)$ system of functions in the space 
$L_2([a, b]).$ 
It means that 
the Fourier series of the function $\sqrt{r(x)}f(x)\in L_2([a, b])$
with respect to the system of functions 
$\sqrt{r(x)}\Phi_0(x),$ $\sqrt{r(x)}\Phi_1(x),\ldots $
converges in the mean-square sense to the function $\sqrt{r(x)}f(x)$
at the interval $[a, b]$. Moreover,
the Fourier coefficients are defined by the formula
\begin{equation}
\label{www.67}
{\tilde C}_j=\int\limits_a^b f(x)\Phi_j(x)r(x)dx.
\end{equation}

It is known that when solving the problem on oscillations of 
a circular membrane (general case), a boundary-value problem arises 
for the following Euler--Bessel equation
\begin{equation}
\label{www.45}
~~~~~~~~~ r^2 R''(r)+rR'(r)+\left(\lambda^2 r^2-n^2\right)R(r)=0\ \ \ (\lambda\in{\bf R},\ \ \ 
n\in{\bf N}).
\end{equation}
The eigenfunctions of this problem, taking into account 
specific boundary 
conditions, are the following functions 
\begin{equation}
\label{www.55}
J_n\biggl(\mu_j\frac{\tau}{L}\biggr),
\end{equation}
where $\tau\in[0, L]$ and $\mu_j$ $(j=0, 1, 2,\ldots )$ 
are positive roots of the Bessel 
function $J_n(\mu)$ ($n=0, 1, 2,\ldots$)
numbered in ascending order. 

The problem on radial oscillations of a circular membrane leads to
the boundary-value problem for 
the
equation (\ref{www.45}) 
for $n=0$, 
the eigenfunctions of which are the functions (\ref{www.55}) when 
$n=0$.

Let us consider the system of functions
\begin{equation}
\label{dmitri4}
\Psi_j(\tau)=\frac{\sqrt{2}}{T J_{n+1}(\mu_j)} J_n\left(
\frac{\mu_j}{T}\tau\right),\ \ \ j=0, 1, 2,\ldots,
\end{equation}
where 
$$
J_n(x)=\sum\limits_{m=0}^{\infty}(-1)^m\left(\frac{x}{2}\right)^{n+2m}
\frac{1}{\Gamma(m+1)\Gamma(m+n+1)}
$$

\noindent
is the 
Bessel function of the first 
kind,
$$
\Gamma(z)=\int\limits_0^{\infty}e^{-x} x^{z-1} dx
$$
is the gamma-function, $\mu_j$ 
are positive  
roots of the function $J_n(x)$ numbered in ascending order, and
$n$ is a natural number or zero.

Due to the well known properties of 
the Bessel functions, the system $\left\{\Psi_j(\tau)\right\}_{j=0}^{\infty}$
is a complete orthonormal with weight $\tau$ 
system of continuous functions 
in the space $L_2([0, T])$.

Let us use the system of functions (\ref{dmitri4}) in 
Theorem 1.9.

Consider the following iterated stochastic integral with respect to
martingales
$$
\int\limits_0^T\int\limits_0^s dM_{\tau}^{(1)}dM_s^{(2)},
$$
where 
$$
M_s^{(i)}=\int\limits_0^s\sqrt{\tau}d{\bf w}_{\tau}^{(i)}\ \ \ (i=1, 2), 
$$
${\bf w}_{\tau}^{(i)}$ $(i=1, 2)$
are independent standard Wiener 
processes, 
$M_s^{(i)}$ $(i=1, 2)$ are martingales (here
$\rho(\tau)\equiv\tau$), $0\le s\le T.$
In addition, $M_s^{(i)}$ has a Gaussian distribution. 

It is obvious
that the conditions of Theorem 1.9 are fulfilled for $k=2$.
Using Theorem 1.9, we obtain

\newpage
\noindent
$$
\int\limits_0^T\int\limits_0^s dM_{\tau}^{(1)}dM_s^{(2)}=
\hbox{\vtop{\offinterlineskip\halign{
\hfil#\hfil\cr
{\rm l.i.m.}\cr
$\stackrel{}{{}_{p_1,p_2\to \infty}}$\cr
}} }\sum_{j_1=0}^{p_1}\sum_{j_2=0}^{p_2}
{\tilde C}_{j_2j_1}
\zeta_{j_1}^{(1)}\zeta_{j_2}^{(2)},
$$
where 
$$
\zeta_j^{(i)}=\int\limits_0^T\Psi_j(\tau)dM_{\tau}^{(i)}
$$
are independent standard  Gaussian random variables
for various $i$ or $j$ $(i=1, 2,$ $j=0, 1, 2,\ldots),$
$$
{\tilde C}_{j_2 j_1}=
\int\limits_0^T s\Psi_{j_2}(s)
\int\limits_0^s \tau\Psi_{j_1}(\tau)d\tau ds
$$
is the Fourier coefficient.

It is obvious that we can get the same result using the another 
approach: 
we can use Theorem 1.1 for the iterated It\^{o} stochastic integral
$$
\int\limits_0^T \sqrt{s}\int\limits_0^s \sqrt{\tau}d{\bf w}_{\tau}^{(1)}
d{\bf w}_s^{(2)},
$$
and as a system of functions $\{\phi_j(s)\}_{j=0}^{\infty}$
in Theorem 1.1 we can take
$$
\phi_j(s)=
\frac{\sqrt{2s}}{T J_{n+1}(\mu_j)} J_n\left(
\frac{\mu_j}{T}s\right),\ \ \ j=0, 1, 2,\ldots
$$

As a result, we obtain
$$
\int\limits_0^T \sqrt{s}\int\limits_0^s \sqrt{\tau}d{\bf w}_{\tau}^{(1)}
d{\bf w}_s^{(2)}=
\hbox{\vtop{\offinterlineskip\halign{
\hfil#\hfil\cr
{\rm l.i.m.}\cr
$\stackrel{}{{}_{p_1,p_2\to \infty}}$\cr
}} }\sum_{j_1=0}^{p_1}\sum_{j_2=0}^{p_2}
{C}_{j_2j_1}
\zeta_{j_1}^{(1)}\zeta_{j_2}^{(2)},
$$
where 
$$
\zeta_j^{(i)}=\int\limits_0^T\phi_j(\tau)d{\bf w}_{\tau}^{(i)}
$$
are independent standard Gaussian random variables
for various $i$ or $j$ $(i=1, 2,$
$j=0, 1, 2,\ldots ),$
$$
C_{j_2 j_1}=
\int\limits_0^T {\sqrt s}\phi_{j_2}(s)
\int\limits_0^s {\sqrt \tau}\phi_{j_1}(\tau)d\tau ds
$$
is the Fourier coefficient. Obviously that 
$C_{j_2 j_1}={\tilde C}_{j_2 j_1}.$

Easy calculation demonstrates that
$$
\tilde\phi_j(s)=
\frac{\sqrt{2(s-t)}}{(T-t) J_{n+1}(\mu_j)} J_n\left(
\frac{\mu_j}{T-t}(s-t)\right),\ \ \ j=0, 1, 2,\ldots
$$ 
is a complete orthonormal system of functions in the space 
$L_2([t, T]).$

Then, using Theorem 1.1, we obtain
$$
\int\limits_t^T \sqrt{s-t}\int\limits_t^s \sqrt{\tau-t}d{\bf w}_{\tau}^{(1)}
d{\bf w}_s^{(2)}=
\hbox{\vtop{\offinterlineskip\halign{
\hfil#\hfil\cr
{\rm l.i.m.}\cr
$\stackrel{}{{}_{p_1,p_2\to \infty}}$\cr
}} }\sum_{j_1=0}^{p_1}\sum_{j_2=0}^{p_2}
{C}_{j_2j_1}
\tilde\zeta_{j_1}^{(1)}\tilde\zeta_{j_2}^{(2)},
$$
where 
$$
\tilde\zeta_j^{(i)}=\int\limits_t^T\tilde\phi_j(\tau)d{\bf w}_{\tau}^{(i)}
$$
are independent standard Gaussian random variables
for various $i$ or $j$ $(i=1, 2,$
$j=0, 1, 2,\ldots ),$
$$
C_{j_2 j_1}=
\int\limits_t^T \sqrt{s-t}\tilde\phi_{j_2}(s)
\int\limits_t^s \sqrt{\tau-t}\tilde\phi_{j_1}(\tau)d\tau ds
$$
is the Fourier coefficient.

\section{Convergence with Probability 1
of Expansions of Iterated
It\^{o} Stochastic Integrals in Theorem 1.1}

\subsection{Convergence with Probability 1
of Expansions of Some Iterated
It\^{o} Stochastic Integrals of Multiplicities 1 and 2}

Let us address now to the convergence with probability 1 (w.~p.~1). 
Consider in detail the iterated It\^{o} stochastic integral
(\ref{k1000}) and its expansion, which is corresponds to (\ref{4004})
for the case $i_1\ne i_2$
\begin{equation}
\label{hqye}
~~~~~~ I_{(00)T,t}^{(i_1 i_2)}=
\frac{T-t}{2}\Biggl(\zeta_0^{(i_1)}\zeta_0^{(i_2)}+\sum_{i=1}^{\infty}
\frac{1}{\sqrt{4i^2-1}}\left(
\zeta_{i-1}^{(i_1)}\zeta_{i}^{(i_2)}-
\zeta_i^{(i_1)}\zeta_{i-1}^{(i_2)}\right)\Biggr).
\end{equation}

First, note the well known fact \cite{Shir999}.

{\bf Lemma 1.8.}\ {\it If for the sequence of random variables
$\xi_p$ and for some
$\alpha>0$ the number series 
$$
\sum\limits_{p=1}^{\infty}{\sf M}\left\{\left|\xi_p\right|^{\alpha}\right\}
$$
converges$,$ then the sequence $\xi_p$ converges to zero w.~p.~{\rm 1}.}

In our specific case $(i_1\ne i_2)$
$$
I_{(00)T,t}^{(i_1 i_2)}=I_{(00)T,t}^{(i_1 i_2)p}+\xi_p,\ \ \
\xi_p
=\frac{T-t}{2}\sum_{i=p+1}^{\infty}
\frac{1}{\sqrt{4i^2-1}}\left(
\zeta_{i-1}^{(i_1)}\zeta_{i}^{(i_2)}-
\zeta_i^{(i_1)}\zeta_{i-1}^{(i_2)}\right),
$$
where
\begin{equation}
\label{90}
~~~~~~ I_{(00)T,t}^{(i_1 i_2)p}=
\frac{T-t}{2}\Biggl(\zeta_0^{(i_1)}\zeta_0^{(i_2)}+\sum_{i=1}^{p}
\frac{1}{\sqrt{4i^2-1}}\left(
\zeta_{i-1}^{(i_1)}\zeta_{i}^{(i_2)}-
\zeta_i^{(i_1)}\zeta_{i-1}^{(i_2)}\right)\Biggr).
\end{equation}

Furthermore,
$$
{\sf M}\left\{|\xi_p|^2\right\}=
\frac{(T-t)^2}{2}\sum_{i=p+1}^{\infty}
\frac{1}{4i^2-1}
\le \frac{(T-t)^2}{2}\int\limits_{p}^{\infty}
\frac{1}{4x^2-1}dx
=
$$
\begin{equation}
\label{2017rock1}
=-\frac{(T-t)^2}{2}\frac{1}{4}{\rm ln}\biggl|
1-\frac{2}{2p+1}\biggr|\le \frac{C}{p},
\end{equation}

\vspace{1mm}
\noindent
where constant $C$ is independent of $p.$

Therefore, taking $\alpha=2$ in Lemma 1.8, we 
cannot
prove the convergence of 
$\xi_p$ to zero w.~p.~1, 
since the series 
$$
\sum\limits_{p=1}^{\infty}{\sf M}\left\{|\xi_p|^{2}\right\}
$$
will be 
majorized by the 
divergent
Dirichlet series
with the index 1. Let us take $\alpha=4$ and estimate the value
${\sf M}\left\{|\xi_p|^4\right\}$.

From (\ref{2026ch1001s11}) 
for $k=2$, $n=2$ and
(\ref{2017rock1}) we obtain
\begin{equation}
\label{91}
{\sf M}\left\{|\xi_p|^4\right\}\le \frac{K}{p^2}
\end{equation}
and
\begin{equation}
\label{hhq}
\sum_{p=1}^{\infty}
{\sf M}\left\{|\xi_p|^4\right\}\le {K}
\sum_{p=1}^{\infty}\frac{1}{p^2}<\infty,
\end{equation}
where constant $K$ is independent of $p.$

Since the series on the right-hand side of 
(\ref{hhq}) converges, then according to Lemma 1.8,
we obtain that
$\xi_p \to 0$ when $p\to \infty$ w.~p.~1. Then  
$$
I_{(00)T,t}^{(i_1 i_2)p}\to
I_{(00)T,t}^{(i_1 i_2)}\ \ \hbox{when}\ \ p\to \infty\ \ \hbox{w.~p.~1.}
$$

Let us analyze the following iterated It\^{o} stochastic integrals 
$$
I_{(01)T,t}^{(i_1 i_2)}=\int\limits_t^T(t-s)\int\limits_t^s 
d{\bf w}_{\tau}^{(i_1)}d{\bf w}_{s}^{(i_2)},\ \ \ 
I_{(10)T,t}^{(i_1 i_2)}=\int\limits_t^T\int\limits_t^s(t-\tau)
d{\bf w}_{\tau}^{(i_1)}d{\bf w}_{s}^{(i_2)},
$$
whose expansions based on Theorem 1.1 and Legendre polynomials
have the following form (also see Chapter~5,  Sect.~5.1)
$$
I_{(01)T,t}^{(i_1 i_2)}=-\frac{T-t}{2}I_{(00)T,t}^{(i_1 i_2)p}
-\frac{(T-t)^2}{4}\left(
\frac{\zeta_0^{(i_1)}\zeta_1^{(i_2)}}{\sqrt{3}}
+\right.
$$
$$
+\left.\sum_{i=0}^{p}\left(
\frac{(i+2)\zeta_i^{(i_1)}\zeta_{i+2}^{(i_2)}
-(i+1)\zeta_{i+2}^{(i_1)}\zeta_{i}^{(i_2)}}
{\sqrt{(2i+1)(2i+5)}(2i+3)}-
\frac{\zeta_i^{(i_1)}\zeta_{i}^{(i_2)}}{(2i-1)(2i+3)}\right)\right)
+ \xi_p^{(01)},
$$

\vspace{3mm}
$$
I_{(10)T,t}^{(i_1 i_2)}=-\frac{T-t}{2}I_{(00)T,t}^{(i_1 i_2)p}
-\frac{(T-t)^2}{4}\left(
\frac{\zeta_0^{(i_2)}\zeta_1^{(i_1)}}
{\sqrt{3}}+\right.
$$
$$
+\left.\sum_{i=0}^{p}\left(
\frac{(i+1)\zeta_{i+2}^{(i_2)}\zeta_{i}^{(i_1)}
-(i+2)\zeta_{i}^{(i_2)}\zeta_{i+2}^{(i_1)}}
{\sqrt{(2i+1)(2i+5)}(2i+3)}+
\frac{\zeta_i^{(i_1)}\zeta_{i}^{(i_2)}}{(2i-1)(2i+3)}\right)\right)
+\xi_p^{(10)},
$$

\vspace{4mm}
\noindent
where
$$
\xi_p^{(01)}=
-\frac{(T-t)^2}{4}
\Biggl(\sum_{i=p+1}^{\infty}
\frac{1}{\sqrt{4i^2-1}}\left(
\zeta_{i-1}^{(i_1)}\zeta_{i}^{(i_2)}-
\zeta_i^{(i_1)}\zeta_{i-1}^{(i_2)}\right)\Biggr.+
$$

\vspace{-2mm}
$$
+
\Biggl.\sum_{i=p+1}^{\infty}\Biggl(
\frac{(i+2)\zeta_i^{(i_1)}\zeta_{i+2}^{(i_2)}
-(i+1)\zeta_{i+2}^{(i_1)}\zeta_{i}^{(i_2)}}
{\sqrt{(2i+1)(2i+5)}(2i+3)}-
\frac{\zeta_i^{(i_1)}\zeta_{i}^{(i_2)}}{(2i-1)(2i+3)}\Biggr)\Biggr),
$$

\vspace{4mm}

$$
\xi_p^{(10)}=
-\frac{(T-t)^2}{4}
\Biggl(\sum_{i=p+1}^{\infty}
\frac{1}{\sqrt{4i^2-1}}\left(
\zeta_{i-1}^{(i_1)}\zeta_{i}^{(i_2)}-
\zeta_i^{(i_1)}\zeta_{i-1}^{(i_2)}\right)\Biggr.+
$$                                     

\newpage
\noindent
$$
+
\Biggl.\sum_{i=p+1}^{\infty}\Biggl(
\frac{(i+1)\zeta_i^{(i_1)}\zeta_{i+2}^{(i_2)}
-(i+2 )\zeta_{i+2}^{(i_1)}\zeta_{i}^{(i_2)}}
{\sqrt{(2i+1)(2i+5)}(2i+3)}-
\frac{\zeta_i^{(i_1)}\zeta_{i}^{(i_2)}}{(2i-1)(2i+3)}\Biggr)\Biggr).
$$

\vspace{2mm}

Then for the case $i_1\ne i_2$ we obtain
$$
{\sf M}\left\{\left|\xi_p^{(01)}\right|^2\right\}=
\frac{(T-t)^4}{16}\times
$$
$$
\times
\sum\limits_{i=p+1}^{\infty}
\Biggl(\frac{2}{4i^2-1}+\frac{(i+2)^2+(i+1)^2}{(2i+1)(2i+5)(2i+3)^2}
+\frac{1}{(2i-1)^2(2i+3)^2}\Biggr) 
\le 
$$
\begin{equation}
\label{2017mas1}
\le 
K\sum\limits_{i=p+1}^{\infty}\frac{1}{i^2}\le \frac{K}{p},
\end{equation}

\vspace{2mm}
\noindent
where constant $K$ is independent of $p$.

Analogously, we get
\begin{equation}
\label{2017mas2}
{\sf M}\left\{\left|\xi_p^{(10)}\right|^2\right\}
\le \frac{K}{p},
\end{equation}

\noindent
where constant $K$ does not depend on $p$.

From (\ref{2026ch1001s11}) 
for $k=2$, $n=2$ 
and (\ref{2017mas1}), (\ref{2017mas2}) we have
$$
{\sf M}\left\{\left|\xi_p^{(01)}\right|^4\right\}+
{\sf M}\left\{\left|\xi_p^{(10)}\right|^4\right\}\le \frac{K_1}{p^2}
$$

\noindent
and
\begin{equation}
\label{hhq11}
~~~~~~~~~~~~\sum_{p=1}^{\infty}
\left({\sf M}\left\{\left|\xi_p^{(01)}\right|^4\right\}+
{\sf M}\left\{\left|\xi_p^{(10)}\right|^4\right\}\right)\le {K_1}
\sum_{p=1}^{\infty}\frac{1}{p^2}<\infty,
\end{equation}

\noindent
where constant $K_1$ is independent of $p.$

According to (\ref{hhq11}) and Lemma 1.8, we obtain that
$\xi_p^{(01)},\ \xi_p^{(10)}\to 0$ when $p\to \infty$\ w.~p.~1. 
Then  
$$
I_{(01)T,t}^{(i_1 i_2)p}\to
I_{(01)T,t}^{(i_1 i_2)},\ \ \
I_{(10)T,t}^{(i_1 i_2)p}\to
I_{(10)T,t}^{(i_1 i_2)}\ \ \ \hbox{when}\ \ \  p\to \infty\ \ \ 
\hbox{w.~p.~1,}
$$
where $i_1\ne i_2.$

Let us consider the case $i_1=i_2$

$$
I_{(01)T,t}^{(i_1 i_1)}=
\frac{(T-t)^2}{4}-\frac{(T-t)^2}{4}\left(
\left(\zeta_0^{(i_1)}\right)^2+
\frac{\zeta_0^{(i_1)}\zeta_1^{(i_1)}}{\sqrt{3}}+\right.
$$

\newpage
\noindent
$$
\left.
+\sum_{i=0}^{p}\left(\frac{\zeta_i^{(i_1)}\zeta_{i+2}^{(i_1)}}
{\sqrt{(2i+1)(2i+5)}(2i+3)}
-
\frac{\zeta_i^{(i_1)}\zeta_i^{(i_1)}}{(2i-1)(2i+3)}
\right)\right)+
\mu_p^{(01)},
$$

\vspace{4mm}

$$
I_{(10)T,t}^{(i_1 i_1)}=
\frac{(T-t)^2}{4}-\frac{(T-t)^2}{4}\left(
\left(\zeta_0^{(i_1)}\right)^2+
\frac{\zeta_0^{(i_1)}\zeta_1^{(i_1)}}{\sqrt{3}}
+\right.
$$
$$
\left.
+\sum_{i=0}^{p}\left(-\frac{\zeta_i^{(i_1)}\zeta_{i+2}^{(i_1)}}
{\sqrt{(2i+1)(2i+5)}(2i+3)}
+
\frac{\zeta_i^{(i_1)}\zeta_i^{(i_1)}}
{(2i-1)(2i+3)}
\right)\right)+
\mu_p^{(10)},
$$

\vspace{3mm}
\noindent
where
$$
\mu_p^{(01)}=-\frac{(T-t)^2}{4}
\sum_{i=p+1}^{\infty}
\left(\frac{\zeta_i^{(i_1)}\zeta_{i+2}^{(i_1)}}{\sqrt{(2i+1)(2i+5)}(2i+3)}
-
\frac{\zeta_i^{(i_1)}\zeta_i^{(i_1)}}
{(2i-1)(2i+3)}
\right),
$$

$$
\mu_p^{(10)}=-\frac{(T-t)^2}{4}
\sum_{i=p+1}^{\infty}\left(-\frac{\zeta_i^{(i_1)}\zeta_{i+2}^{(i_1)}}
{\sqrt{(2i+1)(2i+5)}(2i+3)}
+
\frac{\zeta_i^{(i_1)}\zeta_i^{(i_1)}}{(2i-1)(2i+3)}
\right).
$$

\vspace{3mm}
Then
$$
{\sf M}\biggl\{\left(\mu_p^{(01)}\right)^2\biggr\}=
{\sf M}\biggl\{\left(\mu_p^{(10)}\right)^2\biggr\}=\frac{(T-t)^4}{16}\times
$$
$$
\times\Biggl(
\sum_{i=p+1}^{\infty}
\frac{1}{(2i+1)(2i+5)(2i+3)^2}+
\sum_{i=p+1}^{\infty}
\frac{2}{(2i-1)^2(2i+3)^2} + \Biggr.
$$
$$
+\Biggl.
\Biggl(\sum_{i=p+1}^{\infty}
\frac{1}{(2i-1)(2i+3)}\Biggr)^2\Biggr)
\le \frac{K}{p^2}
$$

\noindent
and
\begin{equation}
\label{hhq112}
~~~~~~~~~~~~\sum_{p=1}^{\infty}
\left({\sf M}\left\{\left|\mu_p^{(01)}\right|^2\right\}+
{\sf M}\left\{\left|\mu_p^{(10)}\right|^2\right\}\right)\le {K}
\sum_{p=1}^{\infty}\frac{1}{p^2}<\infty,
\end{equation}

\vspace{1mm}
\noindent
where constant $K$ is independent of $p.$

According to Lemma 1.8 and (\ref{hhq112}), we obtain that
$\mu_p^{(01)},$ $\mu_p^{(10)}\to 0$ when $p\to \infty$\ w.~p.~1. 
Then  
$$
I_{(01)T,t}^{(i_1 i_1)p}\to
I_{(01)T,t}^{(i_1 i_1)},\ \ \
I_{(10)T,t}^{(i_1 i_1)p}\to
I_{(10)T,t}^{(i_1 i_1)}\ \ \ \hbox{when}\ \ \  
p\to \infty\ \ \ \hbox{w.~p.~1.}
$$

Analogously, we have

\vspace{-3mm}
$$
I_{(02)T,t}^{(i_1 i_2)p}\to
I_{(02)T,t}^{(i_1 i_2)},\ \ \
I_{(11)T,t}^{(i_1 i_2)p}\to
I_{(11)T,t}^{(i_1 i_2)},\ \ \
I_{(20)T,t}^{(i_1 i_2)p}\to
I_{(20)T,t}^{(i_1 i_2)}\ \ \ \hbox{when}\ \ \ 
p\to \infty\ \ \ \hbox{w.~p.~1,}
$$

\vspace{3mm}
\noindent
where 
$$
I_{(02)T,t}^{(i_1 i_2)}=\int\limits_t^T(t-s)^2\int\limits_t^s 
d{\bf w}_{\tau}^{(i_1)}d{\bf w}_{s}^{(i_2)},\ \ \ 
I_{(20)T,t}^{(i_1 i_2)}=\int\limits_t^T\int\limits_t^s(t-\tau)^2
d{\bf w}_{\tau}^{(i_1)}d{\bf w}_{s}^{(i_2)},
$$
$$
I_{(11)T,t}^{(i_1 i_2)}=\int\limits_t^T(t-s)\int\limits_t^s 
(t-\tau)
d{\bf w}_{\tau}^{(i_1)}d{\bf w}_{s}^{(i_2)},
$$

\vspace{2mm}
\noindent
$i_1, i_2=1,\ldots,m.$
This result is based on the 
expansions of stochastic integrals
$I_{(02)T,t}^{(i_1 i_2)},$
$I_{(20)T,t}^{(i_1 i_2)},$
$I_{(11)T,t}^{(i_1 i_2)}$
(see the formulas (\ref{quak1})--(\ref{quak3}) in Chapter 5).

Let us denote
$$
I_{(l)T,t}^{(i_1)}=
\int\limits_t^T(t-s)^l
d{\bf w}_{s}^{(i_1)},
$$
where $l=0, 1, 2\ldots $
                    
The expansions (\ref{4001})--(\ref{4003}), (\ref{uruk}), (\ref{uruk1}) 
(see Chapter 5)
for stochastic integrals
$I_{(0)T,t}^{(i_1)},$
$I_{(1)T,t}^{(i_1)},$
$I_{(2)T,t}^{(i_1)},$
$I_{(3)T,t}^{(i_1)},$
$I_{(l)T,t}^{(i_1)}$
are
correct w.~p.~1  (they include $1, 2, 3, 4$, and $l+1$ members of 
expansion, correspondently).

\subsection{Convergence with Probability 1
of Expansions of Iterated
It\^{o} Stochastic Integrals of Multiplicity $k$ $(k\in{\bf N})$}

In this section, we formulate and prove the theorem on 
convergence with probability 1 (w. p. 1) of expansions 
of iterated It\^{o} stochastic integrals in Theorem 1.1
for the case of multiplicity $k$ $(k\in{\bf N})$.
This section is written on the base of Sect.~1.7.2 
from \cite{12a}-\cite{12aa} as well as on Sect.~6 from \cite{arxiv-3} and
Sect.~9 from \cite{arxiv-1}.

Let us remind the well known fact from the mathematical analysis,
which is connected to existence
of iterated limits.

{\bf Proposition 1.1.} {\it Let $\bigl\{x_{n,m}\bigr\}_{n,m=1}^{\infty}$
be a double sequence and let there exists the limit
$$
\lim\limits_{n,m\to\infty}x_{n,m}=a<\infty.
$$

Moreover, let there exist the limits
$$
\lim\limits_{n\to\infty}x_{n,m}<\infty\ \ \ \hbox{for any m},\ \ \ \ \
\lim\limits_{m\to\infty}x_{n,m}<\infty\ \ \ \hbox{for any n}.
$$

Then there exist the iterated limits
$$
\lim\limits_{n\to\infty}\lim\limits_{m\to\infty}x_{n,m},\ \ \ 
\lim\limits_{m\to\infty}\lim\limits_{n\to\infty}x_{n,m}
$$
and moreover,
$$
\lim\limits_{n\to\infty}\lim\limits_{m\to\infty}x_{n,m}=
\lim\limits_{m\to\infty}\lim\limits_{n\to\infty}x_{n,m}=a.
$$
}

{\bf Theorem 1.10} \cite{12a}-\cite{12aa}, \cite{OK1000}, \cite{arxiv-1}, \cite{arxiv-3}, 
\cite{arxiv-4}.
{\it Let 
$\psi_l(\tau)$ $(l=1,\ldots, k)$ are 
continuously differentiable nonrandom functions on the interval
$[t, T]$ and $\{\phi_j(x)\}_{j=0}^{\infty}$ is a complete
orthonormal system of Legendre polynomials or 
trigonometric functions in the space $L_2([t, T]).$
Then 
$$
J[\psi^{(k)}]_{T,t}^{p,\ldots,p}\ \to \ J[\psi^{(k)}]_{T,t}\ \ \ 
\hbox{if}\ \ \ p\to \infty
$$

\noindent
w.~p.~{\rm 1,} where $J[\psi^{(k)}]_{T,t}^{p,\ldots,p}$
is the expression on the right-hand side of {\rm (\ref{tyyy})}
before passing to the limit 
$\hbox{\vtop{\offinterlineskip\halign{
\hfil#\hfil\cr
{\rm l.i.m.}\cr
$\stackrel{}{{}_{p_1,\ldots,p_k\to \infty}}$\cr
}} }$  for the case $p_1=\ldots=p_k=p,$ i.e. {\rm (}see Theorem {\rm 1.1)}
$$
J[\psi^{(k)}]_{T,t}^{p,\ldots,p}=
\sum_{j_1=0}^{p}\ldots\sum_{j_k=0}^{p}
C_{j_k\ldots j_1}\Biggl(
\prod_{l=1}^k\zeta_{j_l}^{(i_l)}\ -
\Biggr.
$$
$$
-\ \Biggl.
\hbox{\vtop{\offinterlineskip\halign{
\hfil#\hfil\cr
{\rm l.i.m.}\cr
$\stackrel{}{{}_{N\to \infty}}$\cr
}} }\sum_{(l_1,\ldots,l_k)\in {\rm G}_k}
\phi_{j_{1}}(\tau_{l_1})
\Delta{\bf w}_{\tau_{l_1}}^{(i_1)}\ldots
\phi_{j_{k}}(\tau_{l_k})
\Delta{\bf w}_{\tau_{l_k}}^{(i_k)}\Biggr),
$$

\vspace{1mm}
\noindent
where $i_1,\ldots,i_k=1,\ldots,m,$ rest notations are the same as in 
Theorem {\rm 1.1.}}

{\bf Proof.} Let us consider the Parseval equality
\begin{equation}
\label{par1}
~~~~~~~~ \int\limits_{[t,T]^k}K^2(t_1,\ldots,t_k)dt_1\ldots dt_k=
\lim\limits_{p_1,\ldots,p_k\to\infty}
\sum_{j_1=0}^{p_1}\ldots \sum_{j_k=0}^{p_k}
C_{j_k\ldots j_1}^2,
\end{equation}
where
\begin{equation}
\label{pppx}
K(t_1,\ldots,t_k)=
\left\{\begin{matrix}
\psi_1(t_1)\ldots \psi_k(t_k),\ &t_1<\ldots<t_k\cr\cr
0,\ &\hbox{\rm otherwise}
\end{matrix}
\right.\ \ \
=\ \ 
\prod\limits_{l=1}^k
\psi_l(t_l)\ \prod\limits_{l=1}^{k-1}{\bf 1}_{\{t_l<t_{l+1}\}},
\end{equation}

\noindent
where $t_1,\ldots,t_k\in [t, T]$ for $k\ge 2$ and 
$K(t_1)\equiv\psi_1(t_1)$ for $t_1\in[t, T],$ 
${\bf 1}_A$ denotes the indicator of the set $A$,
\begin{equation}
\label{ppppax}
C_{j_k\ldots j_1}=\int\limits_{[t,T]^k}
K(t_1,\ldots,t_k)\prod_{l=1}^{k}\phi_{j_l}(t_l)dt_1\ldots dt_k
\end{equation}
is the Fourier coefficient.

Using (\ref{pppx}), we obtain
$$
C_{j_k\ldots j_1}=
\int\limits_t^T
\phi_{j_k}(t_k)\psi_k(t_k)\ldots \int\limits_t^{t_2}
\phi_{j_1}(t_1)\psi_1(t_1)dt_1\ldots dt_k.
$$

Further, we denote
$$
\lim\limits_{p_1,\ldots,p_k\to\infty}
\sum_{j_1=0}^{p_1}\ldots \sum_{j_k=0}^{p_k}
C_{j_k\ldots j_1}^2\stackrel{\sf def}{=}
\sum_{j_1,\ldots,j_k=0}^{\infty}
C_{j_k\ldots j_1}^2.
$$

If $p_1=\ldots=p_k=p,$ then we also write
$$
\lim\limits_{p\to\infty}
\sum_{j_1=0}^{p}\ldots \sum_{j_k=0}^{p}
C_{j_k\ldots j_1}^2\stackrel{\sf def}{=}
\sum_{j_1,\ldots,j_k=0}^{\infty}
C_{j_k\ldots j_1}^2.
$$

From the other hand, for iterated limits we write
$$
\lim\limits_{p_1\to\infty}\ldots \lim\limits_{p_k\to\infty}
\sum_{j_1=0}^{p_1}\ldots \sum_{j_k=0}^{p_k}
C_{j_k\ldots j_1}^2\stackrel{\sf def}{=}
\sum_{j_1=0}^{\infty}\ldots
\sum_{j_k=0}^{\infty}
C_{j_k\ldots j_1}^2,
$$
$$
\lim\limits_{p_1\to\infty}\lim\limits_{p_2,\ldots,p_k\to\infty}
\sum_{j_1=0}^{p_1}\ldots \sum_{j_k=0}^{p_k}
C_{j_k\ldots j_1}^2\stackrel{\sf def}{=}
\sum_{j_1=0}^{\infty}
\sum_{j_2,\ldots,j_k=0}^{\infty}
C_{j_k\ldots j_1}^2
$$
and so on.

Let us consider the following lemma.

{\bf Lemma 1.9.}\ {\it The following equalities are fulfilled
$$
\sum_{j_1,\ldots,j_k=0}^{\infty}
C_{j_k\ldots j_1}^2=
\sum_{j_1=0}^{\infty}\ldots
\sum_{j_k=0}^{\infty}
C_{j_k\ldots j_1}^2=
$$
\begin{equation}
\label{lem1}
=\sum_{j_k=0}^{\infty}\ldots
\sum_{j_1=0}^{\infty}
C_{j_k\ldots j_1}^2=
\sum_{j_{q_1}=0}^{\infty}\ldots
\sum_{j_{q_k}=0}^{\infty}
C_{j_k\ldots j_1}^2
\end{equation}
for any permutation $(q_1,\ldots,q_k)$ such that
$\{q_1,\ldots,q_k\}=\{1,\ldots,k\}.$}

{\bf Proof.} Let us consider the value
\begin{equation}
\label{21}
\sum_{j_{q_l}=0}^{p}\ldots
\sum_{j_{q_k}=0}^{p}
C_{j_k\ldots j_1}^2
\end{equation}
for any permutation $(q_l,\ldots,q_k)$, where $l=1,2,\ldots,k$,
$\{q_1,\ldots,q_k\}=\{1,\ldots,k\}.$

Obviously, the expression (\ref{21}) 
defines the non-decreasing sequence with respect to $p$.
Moreover,
$$
\sum_{j_{q_l}=0}^{p}\ldots
\sum_{j_{q_k}=0}^{p}
C_{j_k\ldots j_1}^2\le 
\sum_{j_{q_1}=0}^{p}\sum_{j_{q_2}=0}^{p}\ldots
\sum_{j_{q_k}=0}^{p}
C_{j_k\ldots j_1}^2\le 
$$
$$
\le
\sum_{j_1,\ldots,j_k=0}^{\infty}
C_{j_k\ldots j_1}^2<\infty.
$$

Then the following limit
$$
\lim\limits_{p\to\infty}\sum\limits_{j_{q_l}=0}^p \ldots 
\sum\limits_{j_{q_k}=0}^{p}
C_{j_k\ldots j_1}^2=
\sum_{j_{q_l},\ldots,j_{q_k}=0}^{\infty}
C_{j_k\ldots j_1}^2
$$
exists.

Let $p_l,\ldots,p_k$ simultaneously tend to infinity.
Then $g, r\to \infty$, where $g=\min\{p_l,\ldots,p_k\}$ and
$r=\max\{p_l,\ldots,p_k\}$. Moreover,
$$
\sum_{j_{q_l}=0}^{g}\ldots
\sum_{j_{q_k}=0}^{g}
C_{j_k\ldots j_1}^2\le 
\sum_{j_{q_l}=0}^{p_l}\ldots
\sum_{j_{q_k}=0}^{p_k}
C_{j_k\ldots j_1}^2\le
\sum_{j_{q_l}=0}^{r}\ldots
\sum_{j_{q_k}=0}^{r}
C_{j_k\ldots j_1}^2.
$$

This means that the existence of the limit 
\begin{equation}
\label{1c1c}
\lim\limits_{p\to\infty}\sum_{j_{q_l}=0}^{p}\ldots
\sum_{j_{q_k}=0}^{p}
C_{j_k\ldots j_1}^2
\end{equation}
implies the existence of the limit 
\begin{equation}
\label{1d1d}
\lim\limits_{p_l,\ldots,p_k\to\infty}\sum_{j_{q_l}=0}^{p_l}\ldots
\sum_{j_{q_k}=0}^{p_k}
C_{j_k\ldots j_1}^2
\end{equation}
and equality of the limits (\ref{1c1c}) and (\ref{1d1d}).
 
Taking into account the above reasoning, we have 
$$
\lim\limits_{p,q\to\infty}\sum_{j_{q_l}=0}^{q}\sum_{j_{q_{l+1}}=0}^{p}\ldots
\sum_{j_{q_k}=0}^{p}
C_{j_k\ldots j_1}^2=
\lim\limits_{p\to\infty}\sum_{j_{q_l}=0}^{p}\ldots
\sum_{j_{q_k}=0}^{p}
C_{j_k\ldots j_1}^2=
$$
\begin{equation}
\label{1h1h}
=\lim\limits_{p_l,\ldots,p_k\to\infty}\sum_{j_{q_l}=0}^{p_l}\ldots
\sum_{j_{q_k}=0}^{p_k}
C_{j_k\ldots j_1}^2.
\end{equation}

Since the limit
$$
\sum_{j_1,\ldots,j_k=0}^{\infty}
C_{j_k\ldots j_1}^2
$$

\vspace{1mm}
\noindent
exists (see the Parseval equality (\ref{par1})), then from Proposition 1.1
we have
$$
\sum_{j_{q_1}=0}^{\infty}\sum_{j_{q_2},\ldots,j_{q_k}=0}^{\infty}
C_{j_k\ldots j_1}^2=
\lim\limits_{q\to\infty}
\lim\limits_{p\to\infty}
\sum_{j_{q_1}=0}^{q}\sum_{j_{q_2}=0}^p \ldots \sum_{j_{q_k}=0}^{p}
C_{j_k\ldots j_1}^2=
$$
\begin{equation}
\label{1b1b}
=\lim\limits_{q,p\to\infty}
\sum_{j_{q_1}=0}^{q}\sum_{j_{q_2}=0}^p \ldots \sum_{j_{q_k}=0}^{p}
C_{j_k\ldots j_1}^2=
\sum_{j_1,\ldots,j_k=0}^{\infty}
C_{j_k\ldots j_1}^2.
\end{equation}

Using (\ref{1h1h}) and Proposition 1.1, we get
$$
\sum_{j_{q_2}=0}^{\infty}\sum_{j_{q_3},\ldots,j_{q_k}=0}^{\infty}
C_{j_k\ldots j_1}^2=
\lim\limits_{q\to\infty}
\lim\limits_{p\to\infty}
\sum_{j_{q_2}=0}^{q}\sum_{j_{q_3}=0}^p \ldots \sum_{j_{q_k}=0}^{p}
C_{j_k\ldots j_1}^2=
$$
\begin{equation}
\label{1a1a}
=\lim\limits_{q,p\to\infty}
\sum_{j_{q_2}=0}^{q}\sum_{j_{q_3}=0}^p \ldots \sum_{j_{q_k}=0}^{p}
C_{j_k\ldots j_1}^2=
\sum_{j_{q_2},\ldots,j_{q_k}=0}^{\infty}
C_{j_k\ldots j_1}^2.
\end{equation}

Combining (\ref{1a1a}) and (\ref{1b1b}), we obtain
$$
\sum_{j_{q_1}=0}^{\infty}\sum_{j_{q_2}=0}^{\infty}
\sum_{j_{q_3},\ldots,j_{q_k}=0}^{\infty}
C_{j_k\ldots j_1}^2=
\sum_{j_{1},\ldots,j_{k}=0}^{\infty}
C_{j_k\ldots j_1}^2.
$$

Repeating the above steps, we complete the proof of Lemma
1.9.

Further, let us show that for $s=1,\ldots,k$
$$
\sum_{j_1=0}^{\infty}\ldots
\sum_{j_{s-1}=0}^{\infty}
\sum_{j_s=p+1}^{\infty}\sum_{j_{s+1}=0}^{\infty}\ldots \sum_{j_k=0}^{\infty}
C_{j_k\ldots j_1}^2=
$$
\begin{equation}
\label{d11oh}
=
\sum_{j_s=p+1}^{\infty}\sum_{j_{s-1}=0}^{\infty}\ldots
\sum_{j_{1}=0}^{\infty}
\sum_{j_{s+1}=0}^{\infty}\ldots \sum_{j_k=0}^{\infty}
C_{j_k\ldots j_1}^2.
\end{equation}

\vspace{2mm}

Using the arguments which we used in the proof of Lemma 1.9, we have
$$
\lim\limits_{n\to\infty}
\sum_{j_1=0}^{n}\ldots
\sum_{j_{s-1}=0}^{n}
\sum_{j_s=0}^{p}\sum_{j_{s+1}=0}^{n}\ldots \sum_{j_k=0}^{n}
C_{j_k\ldots j_1}^2=
$$
\begin{equation}
\label{ura0}
~~~~~~~ =\sum_{j_s=0}^{p}\ \sum_{j_{1},\ldots, j_{s-1}, 
j_{s+1},\ldots,j_k=0}^{\infty}
C_{j_k\ldots j_1}^2
=\sum_{j_s=0}^{p}\sum_{j_{q_1}=0}^{\infty}\ldots
\sum_{j_{q_{k-1}}=0}^{\infty}
C_{j_k\ldots j_1}^2
\end{equation}

\vspace{2mm}
\noindent
for any permutation $(q_1,\ldots,q_{k-1})$ such that
$\{q_1,\ldots,q_{k-1}\}=\{1,\ldots,s-1,s+1,\ldots,k\}$,
where $p$ is a fixed natural number.

Obviously, we obtain
$$
\sum_{j_s=0}^{p}\sum_{j_{q_1}=0}^{\infty}\ldots
\sum_{j_{q_{k-1}}=0}^{\infty}
C_{j_k\ldots j_1}^2=
\sum_{j_{q_1}=0}^{\infty}\ldots \sum_{j_s=0}^{p} \ldots
\sum_{j_{q_{k-1}}=0}^{\infty}C_{j_k\ldots j_1}^2= \ldots =
$$
\begin{equation}
\label{ura1}
=
\sum_{j_{q_1}=0}^{\infty}\ldots 
\sum_{j_{q_{k-1}}=0}^{\infty}
\sum_{j_s=0}^{p}
C_{j_k\ldots j_1}^2.
\end{equation}

\vspace{2mm}

Using (\ref{ura0}), (\ref{ura1}) and Lemma 1.9, we get
$$
\sum_{j_1=0}^{\infty}\ldots
\sum_{j_{s-1}=0}^{\infty}
\sum_{j_s=p+1}^{\infty}\sum_{j_{s+1}=0}^{\infty}\ldots \sum_{j_k=0}^{\infty}
C_{j_k\ldots j_1}^2=
\sum_{j_1=0}^{\infty}\ldots
\sum_{j_{s-1}=0}^{\infty}
\sum_{j_s=0}^{\infty}\sum_{j_{s+1}=0}^{\infty}\ldots \sum_{j_k=0}^{\infty}
C_{j_k\ldots j_1}^2-
$$
$$
-\sum_{j_1=0}^{\infty}\ldots
\sum_{j_{s-1}=0}^{\infty}
\sum_{j_s=0}^{p}\sum_{j_{s+1}=0}^{\infty}\ldots \sum_{j_k=0}^{\infty}
C_{j_k\ldots j_1}^2=
$$

$$
=
\sum_{j_s=0}^{\infty}
\sum_{j_{s-1}=0}^{\infty}\ldots
\sum_{j_1=0}^{\infty}\sum_{j_{s+1}=0}^{\infty}\ldots \sum_{j_k=0}^{\infty}
C_{j_k\ldots j_1}^2-
\sum_{j_s=0}^{p}
\sum_{j_{s-1}=0}^{\infty}\ldots
\sum_{j_1=0}^{\infty}\sum_{j_{s+1}=0}^{\infty}\ldots \sum_{j_k=0}^{\infty}
C_{j_k\ldots j_1}^2=
$$

$$
=
\sum_{j_s=p+1}^{\infty}
\sum_{j_{s-1}=0}^{\infty}\ldots
\sum_{j_1=0}^{\infty}\sum_{j_{s+1}=0}^{\infty}\ldots \sum_{j_k=0}^{\infty}
C_{j_k\ldots j_1}^2.
$$

\vspace{2mm}

So, the equality (\ref{d11oh}) is proved.

Using the Parseval equality and Lemma 1.9, we obtain
$$
\int\limits_{[t,T]^k}K^2(t_1,\ldots,t_k)dt_1\ldots dt_k-
\sum_{j_1=0}^{p}\ldots \sum_{j_k=0}^{p}
C_{j_k\ldots j_1}^2=
$$

\vspace{-6mm}
$$
=\sum_{j_1,\ldots,j_k=0}^{\infty}
C_{j_k\ldots j_1}^2-
\sum_{j_1=0}^{p}\ldots \sum_{j_k=0}^{p}
C_{j_k\ldots j_1}^2=
$$
$$
=
\sum_{j_1=0}^{\infty}\ldots \sum_{j_k=0}^{\infty}
C_{j_k\ldots j_1}^2-
\sum_{j_1=0}^{p}\ldots \sum_{j_k=0}^{p}
C_{j_k\ldots j_1}^2=
$$
$$
=\sum_{j_1=0}^{p}\sum_{j_2=0}^{\infty}\ldots \sum_{j_k=0}^{\infty}
C_{j_k\ldots j_1}^2+
\sum_{j_1=p+1}^{\infty}\sum_{j_2=0}^{\infty}\ldots \sum_{j_k=0}^{\infty}
C_{j_k\ldots j_1}^2-
\sum_{j_1=0}^{p}\ldots \sum_{j_k=0}^{p}
C_{j_k\ldots j_1}^2=
$$

$$
=\sum_{j_1=0}^{p}\sum_{j_2=0}^{p}\sum_{j_3=0}^{\infty}
\ldots \sum_{j_k=0}^{\infty}
C_{j_k\ldots j_1}^2+
\sum_{j_1=0}^{p}\sum_{j_2=p+1}^{\infty}
\sum_{j_3=0}^{\infty}
\ldots \sum_{j_k=0}^{\infty}+
$$

$$
+\sum_{j_1=p+1}^{\infty}\sum_{j_2=0}^{\infty}\ldots \sum_{j_k=0}^{\infty}
C_{j_k\ldots j_1}^2-
\sum_{j_1=0}^{p}\ldots \sum_{j_k=0}^{p}
C_{j_k\ldots j_1}^2=\ldots =
$$

$$
=\sum_{j_1=p+1}^{\infty}\sum_{j_2=0}^{\infty}\ldots \sum_{j_k=0}^{\infty}
C_{j_k\ldots j_1}^2+
\sum_{j_1=0}^p
\sum_{j_2=p+1}^{\infty}\sum_{j_2=0}^{\infty}\ldots \sum_{j_k=0}^{\infty}
C_{j_k\ldots j_1}^2+
$$

$$
+\sum_{j_1=0}^p\sum_{j_2=0}^p
\sum_{j_3=p+1}^{\infty}\sum_{j_4=0}^{\infty}\ldots \sum_{j_k=0}^{\infty}
C_{j_k\ldots j_1}^2+ \ldots +
\sum_{j_1=0}^p\ldots \sum_{j_{k-1}=0}^p
\sum_{j_k=p+1}^{\infty}C_{j_k\ldots j_1}^2\le
$$

$$
\le\sum_{j_1=p+1}^{\infty}\sum_{j_2=0}^{\infty}\ldots \sum_{j_k=0}^{\infty}
C_{j_k\ldots j_1}^2+
\sum_{j_1=0}^{\infty}
\sum_{j_2=p+1}^{\infty}\sum_{j_2=0}^{\infty}\ldots \sum_{j_k=0}^{\infty}
C_{j_k\ldots j_1}^2+
$$

$$
+\sum_{j_1=0}^{\infty}\sum_{j_2=0}^{\infty}
\sum_{j_3=p+1}^{\infty}\sum_{j_4=0}^{\infty}\ldots \sum_{j_k=0}^{\infty}
C_{j_k\ldots j_1}^2+ \ldots +
\sum_{j_1=0}^{\infty}\ldots \sum_{j_{k-1}=0}^{\infty}
\sum_{j_k=p+1}^{\infty}C_{j_k\ldots j_1}^2=
$$

\begin{equation}
\label{aaap}
=\sum\limits_{s=1}^k \left(\sum_{j_1=0}^{\infty}\ldots
\sum_{j_{s-1}=0}^{\infty}
\sum_{j_s=p+1}^{\infty}\sum_{j_{s+1}=0}^{\infty}\ldots \sum_{j_k=0}^{\infty}
C_{j_k\ldots j_1}^2\right).
\end{equation}

\vspace{2mm}

Note that we use the following 
$$
\sum_{j_1=0}^{p}\ldots
\sum_{j_{s-1}=0}^{p}
\sum_{j_s=p+1}^{\infty}\sum_{j_{s+1}=0}^{\infty}\ldots \sum_{j_k=0}^{\infty}
C_{j_k\ldots j_1}^2\le
$$
$$
\le
\sum_{j_1=0}^{m_1}\ldots
\sum_{j_{s-1}=0}^{m_{s-1}}
\sum_{j_s=p+1}^{\infty}\sum_{j_{s+1}=0}^{\infty}\ldots \sum_{j_k=0}^{\infty}
C_{j_k\ldots j_1}^2\le
$$
$$
\le
\lim\limits_{m_{s-1}\to\infty}
\sum_{j_1=0}^{m_1}\ldots
\sum_{j_{s-1}=0}^{m_{s-1}}
\sum_{j_s=p+1}^{\infty}\sum_{j_{s+1}=0}^{\infty}\ldots \sum_{j_k=0}^{\infty}
C_{j_k\ldots j_1}^2=
$$
$$
=
\sum_{j_1=0}^{m_1}\ldots
\sum_{j_{s-2}=0}^{m_{s-2}}\sum_{j_{s-1}=0}^{\infty}
\sum_{j_s=p+1}^{\infty}\sum_{j_{s+1}=0}^{\infty}\ldots \sum_{j_k=0}^{\infty}
C_{j_k\ldots j_1}^2\le
$$
$$
\le\ldots\le
$$
$$
\le\sum_{j_1=0}^{\infty}\ldots
\sum_{j_{s-1}=0}^{\infty}
\sum_{j_s=p+1}^{\infty}\sum_{j_{s+1}=0}^{\infty}\ldots \sum_{j_k=0}^{\infty}
C_{j_k\ldots j_1}^2
$$

\noindent
to derive (\ref{aaap}), where $m_1,\ldots,m_{s-1}>p.$

Denote
$$
C_{j_s\ldots j_1}(\tau)=
\int\limits_t^{\tau}
\phi_{j_s}(t_s)\psi_s(t_s)\ldots \int\limits_t^{t_2}
\phi_{j_1}(t_1)\psi_1(t_1)dt_1\ldots dt_s,
$$
where
$s=1,\ldots,k-1.$

Let us remind the Dini Theorem, which we will use further.

{\bf Theorem (Dini).} {\it 
Let the functional sequence $u_n(x)$ 
be non-decreasing at each point of the interval $[a, b]$.
In addition, all the functions $u_n(x)$
of this sequence and the limit function $u(x)$ are continuous on the interval
$[a, b].$ Then the convergence $u_n(x)$ to 
$u(x)$ is uniform on the interval $[a,b].$}

For $s<k$ due to the Parseval equality and Dini Theorem
as well as (\ref{d11oh}) we obtain
$$
\sum_{j_1=0}^{\infty}\ldots
\sum_{j_{s-1}=0}^{\infty}
\sum_{j_s=p+1}^{\infty}\sum_{j_{s+1}=0}^{\infty}\ldots \sum_{j_k=0}^{\infty}
C_{j_k\ldots j_1}^2=
$$
$$
=
\sum_{j_s=p+1}^{\infty}
\sum_{j_{s-1}=0}^{\infty}\ldots
\sum_{j_{1}=0}^{\infty}
\sum_{j_{s+1}=0}^{\infty}\ldots \sum_{j_k=0}^{\infty}
C_{j_k\ldots j_1}^2=
$$
$$
=
\sum_{j_s=p+1}^{\infty}
\sum_{j_{s-1}=0}^{\infty}\ldots
\sum_{j_{1}=0}^{\infty}
\sum_{j_{s+1}=0}^{\infty}\ldots 
\sum_{j_{k-1}=0}^{\infty}
\int\limits_t^T \psi_k^2(t_k) \left(C_{j_{k-1}\ldots j_1}(t_k)\right)^2 dt_k=
$$
$$
=
\sum_{j_s=p+1}^{\infty}
\sum_{j_{s-1}=0}^{\infty}\ldots
\sum_{j_{1}=0}^{\infty}
\sum_{j_{s+1}=0}^{\infty}\ldots 
\sum_{j_{k-2}=0}^{\infty}
\int\limits_t^T \psi_k^2(t_k) 
\sum_{j_{k-1}=0}^{\infty}\left(C_{j_{k-1}\ldots j_1}(t_k)\right)^2 dt_k=
$$
$$
=
\sum_{j_s=p+1}^{\infty}
\sum_{j_{s-1}=0}^{\infty}\ldots
\sum_{j_{1}=0}^{\infty}
\sum_{j_{s+1}=0}^{\infty}\ldots 
\sum_{j_{k-2}=0}^{\infty}
\int\limits_t^T \psi_k^2(t_k) \int\limits_t^{t_k} \psi_{k-1}^2(t_{k-1}) 
\left(C_{j_{k-2}\ldots j_1}(t_{k-1})\right)^2\times
$$

\vspace{-1mm}
$$
\times dt_{k-1}dt_k\le
$$

\vspace{-5mm}
$$
\le C\sum_{j_s=p+1}^{\infty}
\sum_{j_{s-1}=0}^{\infty}\ldots
\sum_{j_{1}=0}^{\infty}
\sum_{j_{s+1}=0}^{\infty}\ldots 
\sum_{j_{k-2}=0}^{\infty}
\int\limits_t^T 
\left(C_{j_{k-2}\ldots j_1}(\tau)\right)^2 d\tau =
$$
$$
=
C\sum_{j_s=p+1}^{\infty}
\sum_{j_{s-1}=0}^{\infty}\ldots
\sum_{j_{1}=0}^{\infty}
\sum_{j_{s+1}=0}^{\infty}\ldots 
\sum_{j_{k-3}=0}^{\infty}
\int\limits_t^T 
\sum_{j_{k-2}=0}^{\infty}
\left(C_{j_{k-2}\ldots j_1}(\tau)\right)^2 d\tau =
$$
$$
=C\sum_{j_s=p+1}^{\infty}
\sum_{j_{s-1}=0}^{\infty}\ldots
\sum_{j_{1}=0}^{\infty}
\sum_{j_{s+1}=0}^{\infty}\ldots 
\sum_{j_{k-3}=0}^{\infty}
\int\limits_t^T \int\limits_t^{\tau}
\psi_{k-2}^2(\theta)
\left(C_{j_{k-3}\ldots j_1}(\theta)\right)^2
d\theta d\tau\le 
$$
$$
\le K
\sum_{j_s=p+1}^{\infty}
\sum_{j_{s-1}=0}^{\infty}\ldots
\sum_{j_{1}=0}^{\infty}
\sum_{j_{s+1}=0}^{\infty}\ldots 
\sum_{j_{k-3}=0}^{\infty}
\int\limits_t^T
\left(C_{j_{k-3}\ldots j_1}(\tau)\right)^2
d\tau\le 
$$

$$
\le \ldots \le
$$

\vspace{-4mm}
$$
\le C_k
\sum_{j_s=p+1}^{\infty}
\sum_{j_{s-1}=0}^{\infty}\ldots
\sum_{j_{1}=0}^{\infty}
\int\limits_t^T 
\left(C_{j_{s}\ldots j_1}(\tau)\right)^2 d\tau=
$$
\begin{equation}
\label{d14}
=C_k
\sum_{j_s=p+1}^{\infty}
\sum_{j_{s-1}=0}^{\infty}\ldots
\sum_{j_{2}=0}^{\infty}
\int\limits_t^T  \sum_{j_{1}=0}^{\infty}
\left(C_{j_{s}\ldots j_1}(\tau)\right)^2 d\tau,
\end{equation}

\vspace{4mm}
\noindent
where constants $C,$ $K$ depend on $T-t$ and
constant $C_k$ depends on $k$ and $T-t.$

Let us explane more precisely how we obtain (\ref{d14}).
For any function $g(s)\in L_2([t,T])$ we have the following
Parseval equality

\newpage
\noindent
$$
\sum\limits_{j=0}^{\infty}\left(\int\limits_t^{\tau}
\phi_j(s)g(s)ds\right)^2=
\sum\limits_{j=0}^{\infty}\left(\int\limits_t^T
{\bf 1}_{\{s<\tau\}}\phi_j(s)g(s)ds\right)^2=
$$
\begin{equation}
\label{d15}
=\int\limits_t^T
\left({\bf 1}_{\{s<\tau\}}\right)^2 g^2(s)ds=
\int\limits_t^{\tau}
g^2(s)ds.
\end{equation}

The equality (\ref{d15}) has been applied repeatedly when we obtaining
(\ref{d14}).

Using the replacement of integration order in Riemann integrals, we have
$$
C_{j_s\ldots j_1}(\tau)=
\int\limits_t^{\tau}
\phi_{j_s}(t_s)\psi_s(t_s)\ldots \int\limits_t^{t_2}
\phi_{j_1}(t_1)\psi_1(t_1)dt_1\ldots dt_s=
$$
$$
=\int\limits_t^{\tau}
\phi_{j_1}(t_1)\psi_1(t_1)\int\limits_{t_1}^{\tau}
\phi_{j_2}(t_2)\psi_2(t_2)
\ldots
\int\limits_{t_{s-1}}^{\tau}
\phi_{j_s}(t_s)\psi_s(t_s)dt_s\ldots dt_2dt_1
\stackrel{\sf def}{=}
$$
$$
\stackrel{\sf def}{=}
{\tilde C}_{j_s\ldots j_1}(\tau).
$$

\vspace{2mm}

For $l=1,\ldots,s$ we will use the following notation
$$
{\tilde C}_{j_s\ldots j_l}(\tau,\theta)=
$$
$$
=
\int\limits_{\theta}^{\tau}
\phi_{j_l}(t_l)\psi_l(t_l)\int\limits_{t_l}^{\tau}
\phi_{j_{l+1}}(t_{l+1})\psi_{l+1}(t_{l+1})
\ldots
\int\limits_{t_{s-1}}^{\tau}
\phi_{j_s}(t_s)\psi_s(t_s)dt_s\ldots dt_{l+1}dt_l.
$$

Using the Parseval equality and Dini Theorem, from (\ref{d14}) we obtain
$$
\sum_{j_1=0}^{\infty}\ldots
\sum_{j_{s-1}=0}^{\infty}
\sum_{j_s=p+1}^{\infty}\sum_{j_{s+1}=0}^{\infty}\ldots \sum_{j_k=0}^{\infty}
C_{j_k\ldots j_1}^2\le
$$
$$
\le
C_k
\sum_{j_s=p+1}^{\infty}
\sum_{j_{s-1}=0}^{\infty}\ldots
\sum_{j_{2}=0}^{\infty}
\int\limits_t^T  \sum_{j_{1}=0}^{\infty}
\left(C_{j_{s}\ldots j_1}(\tau)\right)^2 d\tau=
$$
$$
=C_k
\sum_{j_s=p+1}^{\infty}
\sum_{j_{s-1}=0}^{\infty}\ldots
\sum_{j_{2}=0}^{\infty}
\int\limits_t^T  \sum_{j_{1}=0}^{\infty}
\left({\tilde C}_{j_{s}\ldots j_1}(\tau)\right)^2 d\tau=
$$
\begin{equation}
\label{molod1}
~~~~~ =C_k
\sum_{j_s=p+1}^{\infty}
\sum_{j_{s-1}=0}^{\infty}\ldots
\sum_{j_{2}=0}^{\infty}
\int\limits_t^T\int\limits_t^{\tau}\psi_1^2(t_1)  
\left({\tilde C}_{j_{s}\ldots j_2}(\tau,t_1)\right)^2 dt_1d\tau=
\end{equation}
\begin{equation}
\label{molod2}
~~~~~ =C_k
\sum_{j_s=p+1}^{\infty}
\sum_{j_{s-1}=0}^{\infty}\ldots
\sum_{j_{3}=0}^{\infty}
\int\limits_t^T\int\limits_t^{\tau}\psi_1^2(t_1)  
\sum_{j_{2}=0}^{\infty}
\left({\tilde C}_{j_{s}\ldots j_2}(\tau,t_1)\right)^2 dt_1d\tau=
\end{equation}
$$
=C_k
\sum_{j_s=p+1}^{\infty}
\sum_{j_{s-1}=0}^{\infty}\ldots
\sum_{j_{3}=0}^{\infty}
\int\limits_t^T\int\limits_t^{\tau}\psi_1^2(t_1)  
\int\limits_{t_1}^{\tau}\psi_2^2(t_2)  
\left({\tilde C}_{j_{s}\ldots j_3}(\tau,t_2)\right)^2 dt_2dt_1d\tau \le
$$
$$
\le C_k
\sum_{j_s=p+1}^{\infty}
\sum_{j_{s-1}=0}^{\infty}\ldots
\sum_{j_{3}=0}^{\infty}
\int\limits_t^T\int\limits_t^{\tau}\psi_1^2(t_1)  
\int\limits_{t}^{\tau}\psi_2^2(t_2)  
\left({\tilde C}_{j_{s}\ldots j_3}(\tau,t_2)\right)^2 dt_2dt_1d\tau\le
$$
$$
\le C^{'}_k
\sum_{j_s=p+1}^{\infty}
\sum_{j_{s-1}=0}^{\infty}\ldots
\sum_{j_{3}=0}^{\infty}
\int\limits_t^T
\int\limits_{t}^{\tau}\psi_2^2(t_2)  
\left({\tilde C}_{j_{s}\ldots j_3}(\tau,t_2)\right)^2 dt_2d\tau
\le 
$$

\vspace{-2mm}
$$
\le  \ldots \le
$$
$$
\le C^{''}_k
\sum_{j_s=p+1}^{\infty}
\int\limits_t^T\int\limits_t^{\tau}
\psi_{s-1}^2(t_{s-1})
\left({\tilde C}_{j_{s}}(\tau,t_{s-1})\right)^2 dt_{s-1} d\tau\le
$$
\begin{equation}
\label{la}
\le {\tilde C}_k
\sum_{j_s=p+1}^{\infty}
\int\limits_t^T\int\limits_t^{\tau}
\left(~\int\limits_{u}^{\tau}\phi_{j_s}(\theta)
\psi_s(\theta)d\theta\right)^2 du d\tau,
\end{equation}

\vspace{1mm}
\noindent
where constants $C^{'}_k,$ $C^{''}_k,$ $\tilde C_k$
depend on $k$ and $T-t.$

Let us explane more precisely how we obtain (\ref{la}).
For any function $g(s)\in L_2([t,T])$ we have the following
Parseval equality
$$
\sum\limits_{j=0}^{\infty}\left(\int\limits_{\theta}^{\tau}
\phi_j(s)g(s)ds\right)^2=
\sum\limits_{j=0}^{\infty}\left(\int\limits_t^T
{\bf 1}_{\{\theta<s<\tau\}}\phi_j(s)g(s)ds\right)^2=
$$
\begin{equation}
\label{d22}
=\int\limits_t^T
\left({\bf 1}_{\{\theta<s<\tau\}}\right)^2 g^2(s)ds=
\int\limits_{\theta}^{\tau}
g^2(s)ds.
\end{equation}

The equality (\ref{d22}) has been applied repeatedly when we obtaining
(\ref{la}).

Let us explane more precisely the passing from (\ref{molod1})
to (\ref{molod2}) (the same steps have been used when we 
derive (\ref{la})).

We have
$$
\int\limits_t^T\int\limits_t^{\tau}\psi_1^2(t_1)  
\sum_{j_{2}=0}^{\infty}
\left({\tilde C}_{j_{s}\ldots j_2}(\tau,t_1)\right)^2 dt_1d\tau -
$$
$$
-
\sum_{j_{2}=0}^{n}\int\limits_t^T\int\limits_t^{\tau}\psi_1^2(t_1)  
\left({\tilde C}_{j_{s}\ldots j_2}(\tau,t_1)\right)^2 dt_1d\tau =
$$
$$
=\int\limits_t^T\int\limits_t^{\tau}\psi_1^2(t_1)  
\sum_{j_{2}=n+1}^{\infty}
\left({\tilde C}_{j_{s}\ldots j_2}(\tau,t_1)\right)^2 dt_1d\tau =
$$
\begin{equation}
\label{molod3}
~~~~~~~~~=\lim\limits_{N\to\infty}
\sum\limits_{j=0}^{N-1}\int\limits_t^{\tau_j}\psi_1^2(t_1)  
\sum_{j_{2}=n+1}^{\infty}
\left({\tilde C}_{j_{s}\ldots j_2}(\tau_j,t_1)\right)^2 dt_1 \Delta\tau_j,
\end{equation}

\noindent
where $\{\tau_j\}_{j=0}^{N}$ is a partition of the 
interval $[t, T]$ satisfying the condition (\ref{1111}).

Since the non-decreasing functional sequence $u_n(\tau_j,t_1)$ and its
limit function $u(\tau_j,t_1)$ are continuous on the
interval $[t,\tau_j]\subseteq [t, T]$ with respect to $t_1$,
where
$$
u_n(\tau_j,t_1)=
\sum_{j_{2}=0}^{n}
\left({\tilde C}_{j_{s}\ldots j_2}(\tau_j,t_1)\right)^2,
$$
$$
u(\tau_j,t_1)=
\sum_{j_{2}=0}^{\infty}
\left({\tilde C}_{j_{s}\ldots j_2}(\tau_j,t_1)\right)^2=
\int\limits_{t_1}^{\tau_j}
\psi_2^2(t_2)
\left({\tilde C}_{j_{s}\ldots j_3}(\tau_j,t_2)\right)^2 dt_2,
$$

\noindent 
then by Dini Theorem we have the uniform convergence
of $u_n(\tau_j,t_1)$ to $u(\tau_j,t_1)$ at the interval $[t,\tau_j]\subseteq
[t, T]$
with respect to $t_1.$ As a result, we obtain
\begin{equation}
\label{molod4}
\sum_{j_{2}=n+1}^{\infty}
\left({\tilde C}_{j_{s}\ldots j_2}(\tau_j,t_1)\right)^2<\varepsilon,\ \ \ 
t_1\in [t,\tau_j]
\end{equation}

\noindent
for $n>N(\varepsilon)\in{\bf N}$ ($N(\varepsilon)$ exists
for any $\varepsilon>0$ and it does not depend on $t_1$).

From (\ref{molod3}) and (\ref{molod4}) we obtain
$$
\lim\limits_{N\to\infty}
\sum\limits_{j=0}^{N-1}\int\limits_t^{\tau_j}\psi_1^2(t_1)  
\sum_{j_{2}=n+1}^{\infty}
\left({\tilde C}_{j_{s}\ldots j_2}(\tau_j,t_1)\right)^2 dt_1 \Delta\tau_j
\le
$$
\begin{equation}
\label{molod6}
~~~~~~~~~~~\le \varepsilon 
\lim\limits_{N\to\infty}
\sum\limits_{j=0}^{N-1}\int\limits_t^{\tau_j}\psi_1^2(t_1)  
dt_1 \Delta\tau_j
=\varepsilon \int\limits_t^T
\int\limits_t^{\tau}\psi_1^2(t_1)  
dt_1 d\tau.
\end{equation}

From (\ref{molod6}) we get
$$
\lim\limits_{n\to\infty}\int\limits_t^T\int\limits_t^{\tau}\psi_1^2(t_1)  
\sum_{j_{2}=n+1}^{\infty}
\left({\tilde C}_{j_{s}\ldots j_2}(\tau,t_1)\right)^2 dt_1d\tau = 0.
$$

This fact completes the proof of passing 
from (\ref{molod1})
to (\ref{molod2}).

Let us estimate the integral 
\begin{equation}
\label{st1}
\int\limits_{u}^{\tau}\phi_{j_s}(\theta)
\psi_s(\theta)d\theta
\end{equation}
from (\ref{la}) for the cases when $\{\phi_j(s)\}_{j=0}^{\infty}$
is a complete orthonormal system of Legendre polynomials or
trigonometric functions in the space $L_2([t,T])$.

Note that the estimates for the integral
\begin{equation}
\label{st2}
\int\limits_{t}^{\tau}\phi_{j}(\theta)\psi(\theta)d\theta,\ \ \ j\ge p+1,
\end{equation}
where $\psi(\theta)$ is a continuously
differentiable function on the interval $[t, T]$,
have been obtained in \cite{6}-\cite{12aa},
\cite{art-5}, \cite{arxiv-5}
(also see Sect.~2.2.5).

Let us estimate the integral (\ref{st1}) using the approach from
\cite{art-5}, \cite{arxiv-5}.

First, consider the case of Legendre polynomials.
Then $\phi_j(s)$ is defined as follows

\vspace{-4mm}
$$
\phi_j(\theta)=\sqrt{\frac{2j+1}{T-t}}P_j\left(\left(
\theta-\frac{T+t}{2}\right)\frac{2}{T-t}\right),\ \ \ j\ge 0,
$$

\vspace{3mm}
\noindent
where $P_j(x)$ $(j=0, 1, 2\ldots)$ is the Legendre polynomial.

Further, we have 
$$
\int\limits_v^x\phi_{j}(\theta)\psi(\theta)d\theta=
\frac{\sqrt{T-t}\sqrt{2j+1}}{2}
\int\limits_{z(v)}^{z(x)}P_{j}(y)
\psi(u(y))dy=
$$
$$
=\frac{\sqrt{T-t}}{2\sqrt{2j+1}}\Biggl((P_{j+1}(z(x))-
P_{j-1}(z(x)))\psi(x)-
(P_{j+1}(z(v))-
P_{j-1}(z(v)))\psi(v)-
\Biggr.
$$
\begin{equation}
\label{6000}
\Biggl.-
\frac{T-t}{2}
\int\limits_{z(v)}^{z(x)}((P_{j+1}(y)-P_{j-1}(y))
{\psi}'(u(y))dy\Biggr),
\end{equation}

\noindent
where $x, v\in (t, T),$ $j\ge p+1,$ 
$u(y)$ and $z(x)$ are defined by the following relations
$$
u(y)=\frac{T-t}{2}y+\frac{T+t}{2},\ \ \
z(x)=\left(x-\frac{T+t}{2}\right)\frac{2}{T-t},
$$
${\psi}'$ is a derivative of the function $\psi(\theta)$
with respect to the variable $u(y).$

Note that in (\ref{6000}) we used the following well known property
of the Legendre polynomials
$$
\frac{dP_{j+1}}{dx}(x)-\frac{dP_{j-1}}{dx}(x)=(2j+1)P_j(x),\ \ \ 
j=1, 2,\ldots
$$

From (\ref{6000}) and the well known estimate for the Legendre
polynomials \cite{Gob} (also see \cite{suet})
$$
|P_j(y)| <\frac{K}{\sqrt{j+1}(1-y^2)^{1/4}},\ \ \ 
y\in (-1, 1),\ \ \ j\in {\bf N},
$$

\noindent
where constant $K$ does not depend on $y$ and $j$, it follows that
\begin{equation}
\label{101oh}
~~~~~\left|
\int\limits_v^x\phi_{j}(\theta)\psi(\theta)d\theta
\right| <
\frac{C}{j}\Biggl(\frac{1}{(1-(z(x))^2)^{1/4}}+
\frac{1}{(1-(z(v))^2)^{1/4}}+C_1\Biggr),
\end{equation}
where $j\in {\bf N},$ $z(x), z(v)\in (-1, 1),$ $x, v\in (t, T)$ and
constants $C, C_1$ do not depend on $j$.

From (\ref{101oh}) we obtain
\begin{equation}
\label{102oh}
\left(
\int\limits_v^x\phi_{j}(\theta)\psi(\theta)d\theta
\right)^2 <
\frac{C_2}{j^2}\Biggl(\frac{1}{(1-(z(x))^2)^{1/2}}+
\frac{1}{(1-(z(v))^2)^{1/2}}+C_3\Biggr),
\end{equation}
where $j\in {\bf N},$ constants $C_2, C_3$ do not depend on $j$.

Let us apply (\ref{102oh}) for the estimate of the right-hand side
of (\ref{la}). We have
$$
\int\limits_t^T\int\limits_t^{\tau}
\left(~\int\limits_{u}^{\tau}\phi_{j_s}(\theta)
\psi_s(\theta)d\theta\right)^2 du d\tau\le
$$
$$
\le \frac{K_1}{j_s^2}
\left(
\int\limits_{-1}^1
\frac{dy}{\left(1-y^2\right)^{1/2}}+
\int\limits_{-1}^1\int\limits_{-1}^x
\frac{dy}{\left(1-y^2\right)^{1/2}}dx + K_2\right)\le
$$
\begin{equation}
\label{103}
\le \frac{K_3}{j_s^2},
\end{equation}

\vspace{2mm}
\noindent
where $j_s\in {\bf N},$ constants $K_1, K_2, K_3$ are independent of $j_s.$

Now consider the trigonometric case.
The complete orthonormal system of trigonometric functions
in the space $L_2([t, T])$ has the following form
\begin{equation}
\label{trig11oh}
~~~~\phi_j(\theta)=\frac{1}{\sqrt{T-t}}
\left\{
\begin{matrix}
1,\ & j=0\cr\cr
\sqrt{2}{\rm sin} \left(2\pi r(\theta-t)/(T-t)\right),\ & j=2r-1\cr\cr
\sqrt{2}{\rm cos} \left(2\pi r(\theta-t)/(T-t)\right),\ & j=2r
\end{matrix}
,\right.
\end{equation}

\noindent
where $r=1, 2,\ldots $

Using the system of functions 
(\ref{trig11oh}), we have
$$
\int\limits_v^x\phi_{2r-1}(\theta)\psi(\theta)d\theta=
\sqrt{\frac{2}{T-t}}\int\limits_v^x
{\rm sin} \frac{2\pi r(\theta-t)}{T-t}\psi(\theta)d\theta=
$$

\vspace{-1mm}
$$
=-\sqrt{\frac{T-t}{2}}\frac{1}{\pi r}\Biggl(
\psi(x){\rm cos}\frac{2\pi r(x-t)}{T-t}-
\psi(v){\rm cos}\frac{2\pi r(v-t)}{T-t}-\Biggr.
$$

\vspace{-1mm}
\begin{equation}
\label{201}
\Biggl.-
\int\limits_v^x
{\rm cos} \frac{2\pi r(\theta-t)}{T-t}\psi'(\theta)d\theta\Biggr),
\end{equation}

\vspace{2mm}
$$
\int\limits_v^x\phi_{2r}(\theta)\psi(\theta)d\theta=
\sqrt{\frac{2}{T-t}}\int\limits_v^x
{\rm cos} \frac{2\pi r(\theta-t)}{T-t}\psi(\theta)d\theta=
$$

\newpage
\noindent
$$
=\sqrt{\frac{T-t}{2}}\frac{1}{\pi r}\Biggl(
\psi(x){\rm sin}\frac{2\pi r(x-t)}{T-t}-
\psi(v){\rm sin}\frac{2\pi r(v-t)}{T-t}-\Biggr.
$$

\vspace{-1mm}
\begin{equation}
\label{202}
\Biggl.-
\int\limits_v^x
{\rm sin} \frac{2\pi r(\theta-t)}{T-t}\psi'(\theta)d\theta\Biggr),
\end{equation}

\vspace{2mm}
\noindent
where $\psi'(\theta)$ is a derivative of the function $\psi(\theta)$
with respect to the variable $\theta.$

Combining (\ref{201}) and (\ref{202}), we obtain for the
trigonometric case
\begin{equation}
\label{203}
\left(
\int\limits_v^x\phi_{j}(\theta)\psi(\theta)d\theta
\right)^2 \le 
\frac{C_4}{j^2},
\end{equation}

\noindent
where $j\in {\bf N},$ constant $C_4$ is independent of $j.$

From (\ref{203}) we finally have
\begin{equation}
\label{103x}
\int\limits_t^T\int\limits_t^{\tau}
\left(~\int\limits_{u}^{\tau}\phi_{j_s}(\theta)
\psi_s(\theta)d\theta\right)^2 du d\tau
\le \frac{K_4}{j_s^2},
\end{equation}

\noindent
where $j_s\in {\bf N},$ constant $K_4$ is independent of $j_s.$

Combining (\ref{la}), (\ref{103}), and (\ref{103x}), we obtain
$$
\sum_{j_1=0}^{\infty}\ldots
\sum_{j_{s-1}=0}^{\infty}
\sum_{j_s=p+1}^{\infty}\sum_{j_{s+1}=0}^{\infty}\ldots \sum_{j_k=0}^{\infty}
C_{j_k\ldots j_1}^2\le
$$

\vspace{-3mm}

\begin{equation}
\label{fffoh}
\le L_k
\sum_{j_s=p+1}^{\infty}\frac{1}{j_s^2} \le L_k\int\limits_p^{\infty}
\frac{dx}{x^2}=
\frac{L_k}{p},
\end{equation}

\vspace{2mm}
\noindent
where constant $L_k$ depends on $k$ and $T-t.$

Obviously, the case $s=k$ can be considered absolutely analogously to the
case $s<k$. Then from (\ref{aaap}) and (\ref{fffoh})
we obtain
\begin{equation}
\label{ddd1}
~~~~~~~~~ \int\limits_{[t,T]^k}K^2(t_1,\ldots,t_k)dt_1\ldots dt_k-
\sum_{j_1=0}^{p}\ldots \sum_{j_k=0}^{p}
C_{j_k\ldots j_1}^2\le \frac{G_k}{p},
\end{equation}

\noindent
where constant $G_k$ depends on $k$ and $T-t.$

For the further consideration we will use the estimate (\ref{2026ch1001s11}).
Using (\ref{ddd1}) and the estimate (\ref{2026ch1001s11})
for the case $p_1=\ldots=p_k=p$ and $n=2$, 
we obtain
$$
{\sf M}\left\{\biggl(J[\psi^{(k)}]_{T,t}-
J[\psi^{(k)}]_{T,t}^{p,\ldots,p}\biggr)^{4}\right\}\le
$$
$$
\le C_{2,k}
\left(
\int\limits_{[t,T]^k}
K^2(t_1,\ldots,t_k)
dt_1\ldots dt_k -\sum_{j_1=0}^{p}\ldots
\sum_{j_k=0}^{p}C^2_{j_k\ldots j_1}
\right)^2\le 
$$
\begin{equation}
\label{fff5}
\le\frac{H_{2,k}}{p^2},
\end{equation}

\noindent
where 
$$
C_{n,k}=(k!)^{n} (2n-1)^{nk} 
$$

\vspace{3mm}
\noindent
and $H_{2,k}=G_k^2{C}_{2,k}.$

Let $\alpha$ and $\xi_p$ in
Lemma 1.8 be chosen as follows
$$
\alpha=4,\ \ \ \xi_p=\biggl|J[\psi^{(k)}]_{T,t}-
J[\psi^{(k)}]_{T,t}^{p,\ldots,p}\biggr|.
$$

From (\ref{fff5}) we obtain
\begin{equation}
\label{qqq1oh}
~~~~~~~~~~\sum\limits_{p=1}^{\infty}
{\sf M}\left\{\biggl(J[\psi^{(k)}]_{T,t}-
J[\psi^{(k)}]_{T,t}^{p,\ldots,p}\biggr)^{4}\right\}\le
H_{2,k}\sum\limits_{p=1}^{\infty}\frac{1}{p^2}<\infty.
\end{equation}

\vspace{2mm}

Using Lemma 1.8 and the estimate (\ref{qqq1oh}), we have

\vspace{-1mm}
$$
J[\psi^{(k)}]_{T,t}^{p,\ldots,p}\ \to \ J[\psi^{(k)}]_{T,t}\ \ \ 
\hbox{if}\ \ \ p\to \infty
$$

\vspace{3mm}
\noindent
w.~p.~1, where (see Theorem 1.1)
$$
J[\psi^{(k)}]_{T,t}^{p,\ldots,p}=
\sum_{j_1=0}^{p}\ldots\sum_{j_k=0}^{p}
C_{j_k\ldots j_1}\Biggl(
\prod_{l=1}^k\zeta_{j_l}^{(i_l)}\ -
\Biggr.
$$

\vspace{-5mm}
\begin{equation}
\label{kk0}
~~~~~~ -\ \Biggl.
\hbox{\vtop{\offinterlineskip\halign{
\hfil#\hfil\cr
{\rm l.i.m.}\cr
$\stackrel{}{{}_{N\to \infty}}$\cr
}} }\sum_{(l_1,\ldots,l_k)\in {\rm G}_k}
\phi_{j_{1}}(\tau_{l_1})
\Delta{\bf w}_{\tau_{l_1}}^{(i_1)}\ldots
\phi_{j_{k}}(\tau_{l_k})
\Delta{\bf w}_{\tau_{l_k}}^{(i_k)}\Biggr)
\end{equation}

\vspace{3mm}
\noindent
or (see Theorem 1.2)
$$
J[\psi^{(k)}]_{T,t}^{p,\ldots,p}=
\sum\limits_{j_1=0}^{p}\ldots
\sum\limits_{j_k=0}^{p}
C_{j_k\ldots j_1}\Biggl(
\prod_{l=1}^k\zeta_{j_l}^{(i_l)}+\sum\limits_{r=1}^{[k/2]}
(-1)^r \times
\Biggr.
$$
\begin{equation}
\label{kkohh}
\times
\sum_{\stackrel{(\{\{g_1, g_2\}, \ldots, 
\{g_{2r-1}, g_{2r}\}\}, \{q_1, \ldots, q_{k-2r}\})}
{{}_{\{g_1, g_2, \ldots, 
g_{2r-1}, g_{2r}, q_1, \ldots, q_{k-2r}\}=\{1, 2, \ldots, k\}}}}
\prod\limits_{s=1}^r
{\bf 1}_{\{i_{g_{{}_{2s-1}}}=~i_{g_{{}_{2s}}}\ne 0\}}
\Biggl.{\bf 1}_{\{j_{g_{{}_{2s-1}}}=~j_{g_{{}_{2s}}}\}}
\prod_{l=1}^{k-2r}\zeta_{j_{q_l}}^{(i_{q_l})}\Biggr),
\end{equation}

\vspace{1mm}
\noindent
where $i_1,\ldots,i_k=1,\ldots,m$ in (\ref{kk0}) and (\ref{kkohh}).
Theorem 1.10 is proved.

{\bf Remark 1.6.} {\it From Theorem {\rm 1.4} and Lemma {\rm 1.9}
we obtain
$$
\lim\limits_{p_{q_1}\to \infty}
\varlimsup\limits_{p_{q_2}\to \infty}
\ldots\varlimsup\limits_{p_{q_k}\to\infty}
{\sf M}\left\{\left(
J[\psi^{(k)}]_{T,t}-J[\psi^{(k)}]_{T,t}^{p_1,\ldots,p_k}
\right)^2\right\}
\le 
$$
$$
\le k! \cdot \lim\limits_{p_{q_1}\to 0}\ldots\lim\limits_{p_{q_k}\to\infty}
\left(~\int\limits_{[t,T]^k}
K^2(t_1,\ldots,t_k)
dt_1\ldots dt_k -\sum_{j_1=0}^{p_{1}}\ldots
\sum_{j_{k}=0}^{p_{k}}C^2_{j_k\ldots j_1}\right)=
$$
$$
= k! 
\left(~\int\limits_{[t,T]^k}
K^2(t_1,\ldots,t_k)
dt_1\ldots dt_k -\sum_{j_{q_1}=0}^{\infty}\ldots
\sum_{j_{q_k}=0}^{\infty}C^2_{j_k\ldots j_1}\right)=0
$$

\vspace{2mm}
\noindent
for the following cases{\rm :}

{\rm 1.}\ $i_1,\ldots,i_k=1,\ldots,m$\ \ and\ \ $0<T-t<\infty,$

{\rm 2.}\ $i_1,\ldots,i_k=0, 1,\ldots,m,$\ \ $i_1^2+\ldots+i_k^2>0,$\ \
and\ \ $0<T-t<1.$

\noindent
At that$,$ 
$(q_1,\ldots,q_k)$
is any permutation such that
$\{q_1,\ldots,q_k\}=\{1,\ldots,k\},$
$J[\psi^{(k)}]_{T,t}$ is the stochastic integral {\rm (\ref{ito}),}
$J[\psi^{(k)}]_{T,t}^{p_1,\ldots,p_k}$ is the 
expression on the right-hand side of {\rm (\ref{tyyy})} before
passing to the limit 
$\hbox{\vtop{\offinterlineskip\halign{
\hfil#\hfil\cr
{\rm l.i.m.}\cr
$\stackrel{}{{}_{p_1,\ldots,p_k\to \infty}}$\cr
}} },$ $\varlimsup$ means ${\rm lim\ sup};$ another 
notations are the same as in Theorem {\rm 1.1}.
}

{\bf Remark 1.7.} {\it Taking into account Theorem {\rm 1.4} and
the estimate {\rm (\ref{ddd1}),} we obtain the following
inequality
\begin{equation}
\label{zsel1}
{\sf M}\left\{\left(
J[\psi^{(k)}]_{T,t}-J[\psi^{(k)}]_{T,t}^{p,\ldots,p}
\right)^2\right\}\le \frac{k! P_k (T-t)^k}{p},
\end{equation}

\noindent
where $i_1,\ldots,i_k=1,\ldots,m$ and
constant $P_k$ depends only on $k.$ 

The estimate {\rm (\ref{zsel1})} can be written
in a slightly different form.
Let us consider this
question in more detail.

By analogy with {\rm (\ref{ten-1005})} we have
\begin{equation}
\label{rwwr1}
\lim\limits_{p\to\infty}\sum\limits_{j_1,\ldots,j_k=0}^p
C_{j_k\ldots j_1}C_{j_{m_k}\ldots j_{m_1}}=0,
\end{equation}

\noindent
where $(m_1,\ldots,m_k)$ is any permutation of the set
$\{1,\ldots,k\}$ such that $(m_k,\ldots,m_1)\ne (k,\ldots,1);$
braces mean an unordered 
set$,$ and pa\-ren\-the\-ses mean an ordered set.

Further$,$ using {\rm (\ref{rwwr1})} and the estimate {\rm (\ref{ddd1}),} we obtain
$$
\left\vert\sum\limits_{j_1,\ldots,j_k=0}^p
C_{j_k\ldots j_1}C_{j_{m_k}\ldots j_{m_1}}\right\vert=
\left\vert\sum\limits_{j_1,\ldots,j_k=0}^{\infty}
C_{j_k\ldots j_1}C_{j_{m_k}\ldots j_{m_1}}-
\sum\limits_{j_1,\ldots,j_k=0}^p
C_{j_k\ldots j_1}C_{j_{m_k}\ldots j_{m_1}}\right\vert\le
$$
$$
\le
\left(\sum\limits_{j_1,\ldots,j_k=0}^{\infty}-
\sum\limits_{j_1,\ldots,j_k=0}^p\right)
\left\vert C_{j_k\ldots j_1}C_{j_{m_k}\ldots j_{m_1}}\right\vert\le
$$
$$
\le
\frac{1}{2}\left(\sum\limits_{j_1,\ldots,j_k=0}^{\infty}-
\sum\limits_{j_1,\ldots,j_k=0}^p\right)
\left( C_{j_k\ldots j_1}^2+C_{j_{m_k}\ldots j_{m_1}}^2\right)=
\left(\sum\limits_{j_1,\ldots,j_k=0}^{\infty}-
\sum\limits_{j_1,\ldots,j_k=0}^p\right)
C_{j_k\ldots j_1}^2=
$$
\begin{equation}
\label{rwwr3}
~~~~~~~=\int\limits_{[t,T]^k}
K^2(t_1,\ldots,t_k)
dt_1\ldots dt_k-
\sum\limits_{j_1,\ldots,j_k=0}^p
C_{j_k\ldots j_1}^2\le \frac{G_k}{p},
\end{equation}

\noindent
where constant $G_k$ depends on $k$ and $T-t.$

Combining {\rm (\ref{tttr11}),} {\rm (\ref{rwwr2}),} {\rm (\ref{ddd1}),} 
and {\rm (\ref{rwwr3}),} 
we get
$$
{\sf M}\left\{\left(
J[\psi^{(k)}]_{T,t}-J[\psi^{(k)}]_{T,t}^{p,\ldots,p}
\right)^2\right\}\le \frac{\tilde P_k (T-t)^k}{p},
$$
where $i_1,\ldots,i_k=1,\ldots,m$ and
constant $\tilde P_k$ depends only on $k.$ 

It is easy to see that from the proof of Theorem~{\rm 1.4}
and {\rm (\ref{ddd1})} we obtain the estimate
\begin{equation}
\label{zsel1xyzuv}
{\sf M}\left\{\left(
J[\psi^{(k)}]_{T,t}-J[\psi^{(k)}]_{T,t}^{p,\ldots,p}
\right)^2\right\}\le \frac{Q_k}{p},
\end{equation}

\noindent
where $i_1,\ldots,i_k=0,1,\ldots,m$ and
constant $Q_k$ depends only on $k$ and $T-t.$
}

{\bf Remark 1.8.} {\it The estimates
{\rm (\ref{2026ch1001s11})} and 
{\rm (\ref{ddd1})} imply the following 
inequality
$$
{\sf M}\left\{\left(
J[\psi^{(k)}]_{T,t}-J[\psi^{(k)}]_{T,t}^{p,\ldots,p}
\right)^{2n}\right\}\le 
$$

\vspace{1mm}
\begin{equation}
\label{xyzyx1}
\le (k!)^{n} (2n-1)^{nk}\
\frac{\left(P_k\right)^n (T-t)^{nk}}{p^n},
\end{equation}

\vspace{4mm}
\noindent
where $i_1,\ldots,i_k=1,\ldots,m,$\  $n\in{\bf N},$ and 
constant $P_k$ depends only on $k$.}

\subsection{Rate of Convergence with Probability 1
of Expansions of Iterated
It\^{o} Stochastic Integrals of Multiplicity $k$ $(k\in{\bf N})$}

Consider the question
on the rate of convergence w.~p.~1 in Theorem 1.10.
Using the inequality (\ref{xyzyx1}), we obtain

\vspace{-1mm}
\begin{equation}
\label{xyzyx11}
\left({\sf M}\left\{\left(
J[\psi^{(k)}]_{T,t}-J[\psi^{(k)}]_{T,t}^{p,\ldots,p}
\right)^{2n}\right\}\right)^{1/2n}\le \frac{Q_{n,k}}{\sqrt{p}},
\end{equation}

\vspace{3mm}
\noindent
where $n\in {\bf N}$ and 

$$
Q_{n,k}=(2n-1)^{k/2}\ \sqrt{k!}\
\sqrt{P_k}\ (T-t)^{k/2}.
$$

\vspace{5mm}

According to the Lyapunov inequality, we have

\vspace{-1mm}
\begin{equation}
\label{xyzyx12}
\left({\sf M}\biggl\{\left(
J[\psi^{(k)}]_{T,t}-J[\psi^{(k)}]_{T,t}^{p,\ldots,p}
\right)^{n}\biggr\}\right)^{1/n}\le \frac{Q_{n,k}}{\sqrt{p}}
\end{equation}

\vspace{3mm}
\noindent
for all $n\in {\bf N}$. Following \cite{xyz1001} (Lemma 2.1), we get

$$
\biggl|J[\psi^{(k)}]_{T,t}-
J[\psi^{(k)}]_{T,t}^{p,\ldots,p}\biggr|=
\frac{p^{1/2 - \varepsilon}}{p^{1/2 - \varepsilon}}\biggl|J[\psi^{(k)}]_{T,t}-
J[\psi^{(k)}]_{T,t}^{p,\ldots,p}\biggr|\le
$$

\vspace{2mm}
\begin{equation}
\label{xyzyx13}
~~~~~~~~\le 
\frac{1}{p^{1/2 - \varepsilon}}
\sup\limits_{p\in {\bf N}}\left(p^{1/2 - \varepsilon}
\biggl|J[\psi^{(k)}]_{T,t}-
J[\psi^{(k)}]_{T,t}^{p,\ldots,p}\biggr|\right)=\frac{\eta_{\varepsilon}}{p^{1/2 - \varepsilon}}
\end{equation}

\vspace{3mm}
\noindent
w. p. 1, where
$$
\eta_{\varepsilon}=
\sup\limits_{p\in {\bf N}}\left(p^{1/2 - \varepsilon}
\biggl|J[\psi^{(k)}]_{T,t}-
J[\psi^{(k)}]_{T,t}^{p,\ldots,p}\biggr|\right)
$$

\vspace{3mm}
\noindent
and $\varepsilon>0$ is fixed.

For $q>1/\varepsilon,$ $q\in {\bf N}$ we obtain (see (\ref{xyzyx12})) \cite{xyz1001}

\vspace{1mm}
$$
{\sf M}\left\{\left|\eta_{\varepsilon}\right|^q\right\}=
{\sf M}\left\{\left(\sup\limits_{p\in {\bf N}}\left(p^{1/2 - \varepsilon}
\biggl|J[\psi^{(k)}]_{T,t}-
J[\psi^{(k)}]_{T,t}^{p,\ldots,p}\biggr|\right)\right)^q\right\}=
$$

\vspace{4mm}
$$
=
{\sf M}\left\{\sup\limits_{p\in {\bf N}}\left(p^{(1/2 - \varepsilon)q}
\biggl|J[\psi^{(k)}]_{T,t}-
J[\psi^{(k)}]_{T,t}^{p,\ldots,p}\biggr|^q\right)\right\}\le
$$

\vspace{4mm}
$$
\le {\sf M}\left\{\sum\limits_{p=1}^{\infty}p^{(1/2 - \varepsilon)q}
\biggl|J[\psi^{(k)}]_{T,t}-
J[\psi^{(k)}]_{T,t}^{p,\ldots,p}\biggr|^q\right\}=
$$

\vspace{4mm}
$$
= \sum\limits_{p=1}^{\infty}p^{(1/2 - \varepsilon)q}
{\sf M}\left\{\biggl|J[\psi^{(k)}]_{T,t}-
J[\psi^{(k)}]_{T,t}^{p,\ldots,p}\biggr|^q\right\}\le
$$

\vspace{4mm}
\begin{equation}
\label{xyzyx14}
\le
\sum\limits_{p=1}^{\infty}p^{(1/2 - \varepsilon)q}
\frac{\left(Q_{q,k}\right)^q}{p^{q/2}}=
\left(Q_{q,k}\right)^q\sum\limits_{p=1}^{\infty}\frac{1}{p^{\varepsilon q}}<\infty.
\end{equation}

\vspace{6mm}

From (\ref{xyzyx13}) we obtain that for all $\varepsilon>0$
there exists a random variable $\eta_{\varepsilon}$ such that 
the inequality (\ref{xyzyx13}) is fulfilled w.~p.~1 for all $p\in {\bf N}.$
Moreover, from the Lyapunov inequality and (\ref{xyzyx14}), we obtain
${\sf M}\left\{\left|\eta_{\varepsilon}\right|^q\right\}<\infty$
for all $q\ge 1.$

\section{Modification of Theorem 1.1 for the Case
of In\-teg\-ra\-tion Interval $[t, s]$ $(s\in (t, T])$ 
of Iterated It\^{o} Sto\-chas\-tic Integrals}

\subsection{Formulation and Proof of Theorem 1.1 Modification}

Suppose that every $\psi_l(\tau)$ $(l=1,\ldots,k)$ is a continuous 
nonrandom
function on $[t, T]$. 
Define the following function on the hypercube $[t, T]^k$

\newpage
\noindent
\begin{equation}
\label{road901}
\bar K(t_1,\ldots,t_k,s)={\bf 1}_{\{t_k<s\}}K(t_1,\ldots,t_k),
\end{equation}

\vspace{2mm}
\noindent
where the function $K(t_1,\ldots,t_k)$ is defined by 
(\ref{ppp}), $s\in (t, T]$ ($s$ is fixed), 
and ${\bf 1}_A$ is the indicator of the set $A.$
So, we have

\vspace{-2mm}
$$
\bar K(t_1,\ldots,t_k,s)=
{\bf 1}_{\{t_1<\ldots <t_k<s\}}\psi_1(t_1)\ldots \psi_k(t_k)=
$$

\begin{equation}
\label{pppxyz}
=
\left\{\begin{matrix}
\psi_1(t_1)\ldots \psi_k(t_k),\ &t_1<\ldots<t_k<s\cr\cr
0,\ &\hbox{\rm otherwise}
\end{matrix}
\right.,
\end{equation}

\vspace{4mm}
\noindent
where $k\ge 1, $ $t_1,\ldots,t_k\in [t, T],$ and 
$s\in (t, T]$.

Suppose that $\{\phi_j(x)\}_{j=0}^{\infty}$
is a complete orthonormal system of functions in 
the space $L_2([t, T])$.

The function $\bar K(t_1,\ldots,t_k,s)$ defined by
(\ref{pppxyz})
is piecewise continuous in the 
hypercube $[t, T]^k.$
At this situation it is well known that the generalized 
multiple Fourier series 
of $\bar K(t_1,\ldots,t_k,s)\in L_2([t, T]^k)$ is converging 
to this function in the hypercube $[t, T]^k$ in 
the mean-square sense, i.e.

\vspace{-3mm}
\begin{equation}
\label{sos1zxyz}
\hbox{\vtop{\offinterlineskip\halign{
\hfil#\hfil\cr
{\rm lim}\cr
$\stackrel{}{{}_{p_1,\ldots,p_k\to \infty}}$\cr
}} }\Biggl\Vert
\bar K(t_1,\ldots,t_k,s)-
\sum_{j_1=0}^{p_1}\ldots \sum_{j_k=0}^{p_k}
C_{j_k\ldots j_1}(s)
\prod_{l=1}^{k} \phi_{j_l}(t_l)\Biggr\Vert_{L_2([t, T]^k)}=0,
\end{equation}

\vspace{3mm}
\noindent
where
$$
C_{j_k\ldots j_1}(s)=\int\limits_{[t,T]^k}
\bar K(t_1,\ldots,t_k,s)\prod_{l=1}^{k}\phi_{j_l}(t_l)dt_1\ldots dt_k=
$$
\begin{equation}
\label{ppppaxyz}
=\int\limits_t^s\psi_k(t_k)\phi_{j_k}(t_k)\ldots
\int\limits_t^{t_2}
\psi_1(t_1)\phi_{j_1}(t_1)
dt_1\ldots dt_k
\end{equation}

\vspace{2mm}
\noindent
is the Fourier coefficient, and
$$
\left\Vert f\right\Vert_{L_2([t, T]^k)}=\left(~\int\limits_{[t,T]^k}
f^2(t_1,\ldots,t_k)dt_1\ldots dt_k\right)^{1/2}.
$$

Note that
\newpage
\noindent
\begin{equation}
\label{opr22}
J[\psi^{(k)}]_{s,t}=\int\limits_t^s\psi_k(t_k) \ldots \int\limits_t^{t_{2}}
\psi_1(t_1) d{\bf w}_{t_1}^{(i_1)}\ldots
d{\bf w}_{t_k}^{(i_k)}=
\end{equation}
$$
=
\int\limits_t^T {\bf 1}_{\{t_k<s\}}\psi_k(t_k) \ldots \int\limits_t^{t_{2}}
\psi_1(t_1) d{\bf w}_{t_1}^{(i_1)}\ldots
d{\bf w}_{t_k}^{(i_k)}\ \ \ \hbox{w.~p.~1},
$$

\vspace{2mm}
\noindent
where $s\in (t, T]$ ($s$ is fixed), $i_1,\ldots,i_k=0,1,\ldots,m.$

Consider the partition $\{\tau_j\}_{j=0}^N$ of $[t,T]$ such that

\vspace{-4mm}
\begin{equation}
\label{1111xxx}
t=\tau_0<\ldots <\tau_N=T,\ \ \
\Delta_N=
\hbox{\vtop{\offinterlineskip\halign{
\hfil#\hfil\cr
{\rm max}\cr
$\stackrel{}{{}_{0\le j\le N-1}}$\cr
}} }\Delta\tau_j\to 0\ \ \hbox{if}\ \ N\to \infty,\ \ \ 
\Delta\tau_j=\tau_{j+1}-\tau_j.
\end{equation}

\vspace{2mm}

{\bf Theorem 1.11}\ \cite{12aa-after}--\cite{12aa}, \cite{arxiv-1}.\
{\it Suppose that
every $\psi_l(\tau)$ $(l=$ $1,\ldots, k)$ is a continuous 
non\-ran\-dom function on 
$[t, T]$ and $\{\phi_j(x)\}_{j=0}^{\infty}$ is a complete orthonormal system  
of continuous functions in the space $L_2([t,T]).$ 
Then

$$
J[\psi^{(k)}]_{s,t} =
\hbox{\vtop{\offinterlineskip\halign{
\hfil#\hfil\cr
{\rm l.i.m.}\cr
$\stackrel{}{{}_{p_1,\ldots,p_k\to \infty}}$\cr
}} }\sum_{j_1=0}^{p_1}\ldots\sum_{j_k=0}^{p_k}
C_{j_k\ldots j_1}(s)\Biggl(
\prod_{l=1}^k\zeta_{j_l}^{(i_l)} -
\Biggr.
$$

\begin{equation}
\label{road777}
~~~~~~~-\Biggl.
\hbox{\vtop{\offinterlineskip\halign{
\hfil#\hfil\cr
{\rm l.i.m.}\cr
$\stackrel{}{{}_{N\to \infty}}$\cr
}} }\sum_{(l_1,\ldots,l_k)\in {\rm G}_k}
\phi_{j_{1}}(\tau_{l_1})
\Delta{\bf w}_{\tau_{l_1}}^{(i_1)}\ldots
\phi_{j_{k}}(\tau_{l_k})
\Delta{\bf w}_{\tau_{l_k}}^{(i_k)}\Biggr),
\end{equation}

\vspace{4mm}
\noindent
where $J[\psi^{(k)}]_{s,t}$ is defined by {\rm (\ref{opr22}),}
$s\in (t, T]$ {\rm ($s$} is fixed{\rm ),}

\vspace{-5mm}
$$
{\rm G}_k={\rm H}_k\backslash{\rm L}_k,\ \ \
{\rm H}_k=\bigl\{(l_1,\ldots,l_k):\ l_1,\ldots,l_k=0,\ 1,\ldots,N-1\bigr\},
$$
$$
{\rm L}_k=\bigl\{(l_1,\ldots,l_k):\ l_1,\ldots,l_k=0,\ 1,\ldots,N-1;\
l_g\ne l_r\ (g\ne r);\ g, r=1,\ldots,k\bigr\},
$$

\vspace{3mm}
\noindent
${\rm l.i.m.}$ is a limit in the mean-square sense$,$
$i_1,\ldots,i_k=0,1,\ldots,m,$ 
$$
\zeta_{j}^{(i)}=
\int\limits_t^T \phi_{j}(\tau) d{\bf w}_{\tau}^{(i)}
$$
are independent standard Gaussian random variables
for various
$i$ or $j$ {\rm(}in the case when $i\ne 0${\rm),}
$C_{j_k\ldots j_1}(s)$ is the Fourier coefficient {\rm(\ref{ppppaxyz}),}
$\Delta{\bf w}_{\tau_{j}}^{(i)}=
{\bf w}_{\tau_{j+1}}^{(i)}-{\bf w}_{\tau_{j}}^{(i)}$
$(i=0,\ 1,\ldots,m),$\
$\left\{\tau_{j}\right\}_{j=0}^{N}$ is a partition of
$[t,T],$ which satisfies the condition {\rm (\ref{1111xxx})}.}

{\bf Proof.}\ Let us consider the multiple 
stochastic integrals (\ref{30.34}), (\ref{mult11}).
We will write $J[\Phi]_{s,t}^{(k)}$ and $J'[\Phi]_{s,t}^{(k)}$ $(s\in (t, T],$ $s$ is fixed)
if the function
$\Phi(t_1,\ldots,t_k)$ in (\ref{30.34}) and (\ref{mult11})
is replaced by ${\bf 1}_{\{t_1,\ldots,t_k<s\}}\Phi(t_1,\ldots,t_k).$

By analogy with (\ref{pobeda}), we have
\begin{equation}
\label{pobedaxxx}
J'[\Phi]_{s,t}^{(k)}=
\int\limits_t^T\ldots \int\limits_t^{t_2}
{\bf 1}_{\{t_k<s\}} \sum\limits_{(t_1,\ldots,t_k)}\biggl(
\Phi(t_1,\ldots,t_k)
d{\bf w}_{t_1}^{(i_1)}\ldots
d{\bf w}_{t_k}^{(i_k)}\biggr)\ \ \ \hbox{w.~p.~1},
\end{equation}
where 
$$
\sum\limits_{(t_1,\ldots,t_k)}
$$ 
means the sum with respect to all
possible permutations
$(t_1,\ldots,t_k).$ 
At the same time permutations $(t_1,\ldots,t_k)$
when summing 
are performed in (\ref{pobedaxxx}) only in the expression, which
is enclosed in pa\-ren\-the\-ses.
Moreover, 
the nonrandom function $\Phi(t_1,\ldots,t_k)$ is assumed 
to be continuous in the 
cor\-res\-pond\-ing closed domains of integration. The case
when the nonrandom function $\Phi(t_1,\ldots,t_k)$ is 
continuous in 
the open domains of integration and bounded at their boundaries is also possible.

Let us write (\ref{pobedaxxx}) as
\begin{equation}
\label{pobedaxyz}
J'[\Phi]_{s,t}^{(k)}=
\int\limits_t^T\ldots \int\limits_t^{t_2}
\sum\limits_{(t_1,\ldots,t_k)}\biggl(
{\bf 1}_{\{t_k<s\}}
\Phi(t_1,\ldots,t_k)
d{\bf w}_{t_1}^{(i_1)}\ldots
d{\bf w}_{t_k}^{(i_k)}\biggr)\ \ \ \hbox{w.~p.~1},
\end{equation}
where permutations $(t_1,\ldots,t_k)$
when summing 
are performed in (\ref{pobedaxyz}) only in the expression
$\Phi(t_1,\ldots,t_k)
d{\bf w}_{t_1}^{(i_1)}\ldots
d{\bf w}_{t_k}^{(i_k)}$.

It is not difficult to notice that (\ref{pobedaxxx}),
(\ref{pobedaxyz}) can be rewritten in the form (see (\ref{s2s}))
\begin{equation}
\label{s2sxxx}
J'[\Phi]_{s,t}^{(k)}=\sum_{(t_1,\ldots,t_k)}
\int\limits_{t}^{T}
\ldots
\int\limits_{t}^{t_2}
\Phi(t_1,\ldots,t_k){\bf 1}_{\{t_k<s\}}d{\bf w}_{t_1}^{(i_1)}
\ldots
d{\bf w}_{t_k}^{(i_k)}\ \ \ \hbox{w.~p.~1},
\end{equation}

\noindent
where permutations $(t_1,\ldots,t_k)$ when summing are 
performed only in the values
${\bf 1}_{\{t_k<s\}}d{\bf w}_{t_1}^{(i_1)}
\ldots $
$d{\bf w}_{t_k}^{(i_k)}.$ At the same time the indices near 
upper 
limits of integration in the iterated stochastic integrals are changed 
correspondently and if $t_r$ swapped with $t_q$ in the  
permutation $(t_1,\ldots,t_k)$, then $i_r$ swapped with $i_q$ in 
the permutation $(i_1,\ldots,i_k)$.

According to Lemma 1.1, we have

\vspace{-4mm}
$$
J[\psi^{(k)}]_{s,t}=
\hbox{\vtop{\offinterlineskip\halign{
\hfil#\hfil\cr
{\rm l.i.m.}\cr
$\stackrel{}{{}_{N\to \infty}}$\cr
}} }\sum_{l_k=0}^{N-1}\ldots\sum_{l_1=0}^{l_2-1}
{\bf 1}_{\{\tau_{l_k}<s\}}\psi_1(\tau_{l_1})\ldots\psi_k(\tau_{l_k})
\Delta{\bf w}_{\tau_{l_1}}^{(i_1)}
\ldots
\Delta{\bf w}_{\tau_{l_k}}^{(i_k)}=
$$

$$
=
\hbox{\vtop{\offinterlineskip\halign{
\hfil#\hfil\cr
{\rm l.i.m.}\cr
$\stackrel{}{{}_{N\to \infty}}$\cr
}} }\sum_{l_k=0}^{N-1}\ldots\sum_{l_1=0}^{N-1}
{\bf 1}_{\{\tau_{l_k}<s\}}K(\tau_{l_1},\ldots,\tau_{l_k})
\Delta{\bf w}_{\tau_{l_1}}^{(i_1)}
\ldots
\Delta{\bf w}_{\tau_{l_k}}^{(i_k)}=
$$

\vspace{-1mm}
$$
=\hbox{\vtop{\offinterlineskip\halign{
\hfil#\hfil\cr
{\rm l.i.m.}\cr
$\stackrel{}{{}_{N\to \infty}}$\cr
}} }
\sum\limits_{\stackrel{l_1,\ldots,l_k=0}{{}_{l_q\ne l_r;\ 
q\ne r;\ q, r=1,\ldots, k}}}^{N-1}
{\bf 1}_{\{\tau_{l_k}<s\}}K(\tau_{l_1},\ldots,\tau_{l_k})
\Delta{\bf w}_{\tau_{l_1}}^{(i_1)}
\ldots
\Delta{\bf w}_{\tau_{l_k}}^{(i_k)}=
$$

\vspace{-2mm}
\begin{equation}
\label{hehe100xyz}
~~~~~~~~~=
\int\limits_{t}^{T}
\ldots
\int\limits_{t}^{t_2}
\sum_{(t_1,\ldots,t_k)}\left(
{\bf 1}_{\{t_k<s\}}K(t_1,\ldots,t_k)d{\bf w}_{t_1}^{(i_1)}
\ldots
d{\bf w}_{t_k}^{(i_k)}\right)\ \ \ \hbox{w.~p.~1},
\end{equation}

\vspace{2mm}
\noindent
where $K(t_1,\ldots,t_k)$ is defined by 
(\ref{ppp}) and permutations 
$(t_1,\ldots,t_k)$ when summing
are performed only 
in the expression 
$K(t_1,\ldots,t_k)d{\bf w}_{t_1}^{(i_1)}
\ldots d{\bf w}_{t_k}^{(i_k)}.$

According to Lemmas 1.1, 1.3 and (\ref{pobeda}), (\ref{s2s}), (\ref{pobedaxyz})--(\ref{hehe100xyz}),
we get the following representation

\vspace{-5mm}
$$
J[\psi^{(k)}]_{s,t}=
$$

\vspace{-4mm}
$$
=
\sum_{j_1=0}^{p_1}\ldots
\sum_{j_k=0}^{p_k}
C_{j_k\ldots j_1}(s)
\int\limits_{t}^{T}
\ldots
\int\limits_{t}^{t_2}
\sum_{(t_1,\ldots,t_k)}\left(
\phi_{j_1}(t_1)
\ldots
\phi_{j_k}(t_k)
d{\bf w}_{t_1}^{(i_1)}
\ldots
d{\bf w}_{t_k}^{(i_k)}\right)
+
$$

\vspace{2mm}
$$
+R_{T,t,s}^{p_1,\ldots,p_k}=
$$

\vspace{1mm}
$$
=\sum_{j_1=0}^{p_1}\ldots
\sum_{j_k=0}^{p_k}
C_{j_k\ldots j_1}(s)\times
$$

\vspace{-2mm}
$$
\times\ 
\hbox{\vtop{\offinterlineskip\halign{
\hfil#\hfil\cr
{\rm l.i.m.}\cr
$\stackrel{}{{}_{N\to \infty}}$\cr
}} }
\sum\limits_{\stackrel{l_1,\ldots,l_k=0}{{}_{l_q\ne l_r;\ 
q\ne r;\ q, r=1,\ldots, k}}}^{N-1}
\phi_{j_1}(\tau_{l_1})\ldots
\phi_{j_k}(\tau_{l_k})
\Delta{\bf w}_{\tau_{l_1}}^{(i_1)}
\ldots
\Delta{\bf w}_{\tau_{l_k}}^{(i_k)}\ \
+\ R_{T,t,s}^{p_1,\ldots,p_k}=
$$

\newpage
\noindent
$$
=\sum_{j_1=0}^{p_1}\ldots
\sum_{j_k=0}^{p_k}
C_{j_k\ldots j_1}(s)\left(
\hbox{\vtop{\offinterlineskip\halign{
\hfil#\hfil\cr
{\rm l.i.m.}\cr
$\stackrel{}{{}_{N\to \infty}}$\cr
}} }\sum_{l_1,\ldots,l_k=0}^{N-1}
\phi_{j_1}(\tau_{l_1})
\ldots
\phi_{j_k}(\tau_{l_k})
\Delta{\bf w}_{\tau_{l_1}}^{(i_1)}
\ldots
\Delta{\bf w}_{\tau_{l_k}}^{(i_k)}
-\right.
$$
$$
-\left.
\hbox{\vtop{\offinterlineskip\halign{
\hfil#\hfil\cr
{\rm l.i.m.}\cr
$\stackrel{}{{}_{N\to \infty}}$\cr
}} }\sum_{(l_1,\ldots,l_k)\in {\rm G}_k}
\phi_{j_{1}}(\tau_{l_1})
\Delta{\bf w}_{\tau_{l_1}}^{(i_1)}\ldots
\phi_{j_{k}}(\tau_{l_k})
\Delta{\bf w}_{\tau_{l_k}}^{(i_k)}\right)
+
$$

\vspace{1mm}
$$
+R_{T,t,s}^{p_1,\ldots,p_k}=
$$

\vspace{4mm}
$$
=\sum_{j_1=0}^{p_1}\ldots\sum_{j_k=0}^{p_k}
C_{j_k\ldots j_1}(s)\times
$$
$$
\times
\left(
\prod_{l=1}^k\zeta_{j_l}^{(i_l)}-
\hbox{\vtop{\offinterlineskip\halign{
\hfil#\hfil\cr
{\rm l.i.m.}\cr
$\stackrel{}{{}_{N\to \infty}}$\cr
}} }\sum_{(l_1,\ldots,l_k)\in {\rm G}_k}
\phi_{j_{1}}(\tau_{l_1})
\Delta{\bf w}_{\tau_{l_1}}^{(i_1)}\ldots
\phi_{j_{k}}(\tau_{l_k})
\Delta{\bf w}_{\tau_{l_k}}^{(i_k)}\right)+
$$

\vspace{3mm}
$$
+R_{T,t,s}^{p_1,\ldots,p_k}\ \ \ \hbox{w.~p.~1},
$$

\vspace{7mm}
\noindent
where
$$
R_{T,t,s}^{p_1,\ldots,p_k}
=
$$

\vspace{-5mm}
$$
=
\sum_{(t_1,\ldots,t_k)}
\int\limits_{t}^{T}
\ldots
\int\limits_{t}^{t_2}
\left({\bf 1}_{\{t_k<s\}} K(t_1,\ldots,t_k)-
\sum_{j_1=0}^{p_1}\ldots
\sum_{j_k=0}^{p_k}
C_{j_k\ldots j_1}(s)
\prod_{l=1}^k\phi_{j_l}(t_l)\right)\times
$$

$$
\times
d{\bf w}_{t_1}^{(i_1)}
\ldots
d{\bf w}_{t_k}^{(i_k)}=
$$

\begin{equation}
\label{hhhsx1}
=
\sum_{(t_1,\ldots,t_k)}
\int\limits_{t}^{T}
\ldots
\int\limits_{t}^{t_2}
K(t_1,\ldots,t_k){\bf 1}_{\{t_k<s\}}d{\bf w}_{t_1}^{(i_1)}
\ldots
d{\bf w}_{t_k}^{(i_k)}-             
\end{equation}

\vspace{-2mm}
\begin{equation}
\label{y007xxx}
-
\sum_{(t_1,\ldots,t_k)}
\int\limits_{t}^{T}
\ldots
\int\limits_{t}^{t_2}
\sum_{j_1=0}^{p_1}\ldots
\sum_{j_k=0}^{p_k}
C_{j_k\ldots j_1}(s)
\prod_{l=1}^k\phi_{j_l}(t_l)
d{\bf w}_{t_1}^{(i_1)}
\ldots
d{\bf w}_{t_k}^{(i_k)}
\end{equation}

\vspace{4mm}
\noindent
w.~p.~1, where permutations $(t_1,\ldots,t_k)$ when summing 
in (\ref{hhhsx1})
are performed only 
in the values ${\bf 1}_{\{t_k<s\}}d{\bf w}_{t_1}^{(i_1)}
\ldots $
$d{\bf w}_{t_k}^{(i_k)}$. 
At the same time
permutations $(t_1,\ldots,t_k)$ when summing 
in (\ref{y007xxx})
are performed only 
in the values
$d{\bf w}_{t_1}^{(i_1)}
\ldots
d{\bf w}_{t_k}^{(i_k)}.$
Moreover, the indices near 
upper limits of integration in the iterated stochastic integrals 
in (\ref{hhhsx1}), (\ref{y007xxx})
are changed correspondently and if $t_r$ swapped with $t_q$ in the  
permutation $(t_1,\ldots,t_k)$, then $i_r$ swapped with $i_q$ in the 
permutation $(i_1,\ldots,i_k)$.

Let us estimate the remainder
$R_{T,t,s}^{p_1,\ldots,p_k}$ of the series.
According to Lemma 1.2, we have

\vspace{-0.5mm}
$$
{\sf M}\left\{\left(R_{T,t,s}^{p_1,\ldots,p_k}\right)^2\right\}
\le 
$$

\vspace{-2mm}
$$
\le C_k
\sum_{(t_1,\ldots,t_k)}
\int\limits_{t}^{T}
\ldots
\int\limits_{t}^{t_2}
\left(K(t_1,\ldots,t_k){\bf 1}_{\{t_k<s\}}-
\sum_{j_1=0}^{p_1}\ldots
\sum_{j_k=0}^{p_k}
C_{j_k\ldots j_1}(s)
\prod_{l=1}^k\phi_{j_l}(t_l)\right)^2\times
$$

\vspace{0.5mm}
\begin{equation}
\label{yes999}
\times
dt_1
\ldots
dt_k,
\end{equation}

\vspace{2mm}
\noindent
where constant $C_k$ 
depends only
on the multiplicity $k$ of the iterated It\^{o} stochastic integral
$J[\psi^{(k)}]_{s,t}$ and permutations $(t_1,\ldots,t_k)$ when summing 
in (\ref{yes999})
are performed only 
in the values ${\bf 1}_{\{t_k<s\}}$ and $dt_1\ldots dt_k.$
At the same time the indices near 
upper limits of integration in the iterated integrals 
in (\ref{yes999})
are changed correspondently.

Since 
$K(t_1,\ldots,t_k)\equiv 0$
if the condition $t_1<\ldots <t_k$ is not fulfilled, then

$$
{\sf M}\left\{\left(R_{T,t,s}^{p_1,\ldots,p_k}\right)^2\right\}
\le 
$$

\vspace{-1mm}
$$
\le C_k
\sum_{(t_1,\ldots,t_k)}
\int\limits_{t}^{T}
\ldots
\int\limits_{t}^{t_2}
\left(K(t_1,\ldots,t_k){\bf 1}_{\{t_k<s\}}-
\sum_{j_1=0}^{p_1}\ldots
\sum_{j_k=0}^{p_k}
C_{j_k\ldots j_1}(s)
\prod_{l=1}^k\phi_{j_l}(t_l)\right)^2\times
$$

\vspace{-1mm}
\begin{equation}
\label{yes9991}
\times
dt_1
\ldots
dt_k,
\end{equation} 

\vspace{2mm}
\noindent
where permutations $(t_1,\ldots,t_k)$ when summing 
in (\ref{yes9991})
are performed only 
in the values $dt_1\ldots dt_k.$
At the same time the indices near 
upper limits of integration in the iterated integrals 
in (\ref{yes9991})
are changed correspondently.

Then from (\ref{riemann}), (\ref{sos1zxyz}), and (\ref{yes9991}) we obtain

\vspace{-1mm}
$$
{\sf M}\left\{\left(R_{T,t,s}^{p_1,\ldots,p_k}\right)^2\right\}
\le 
$$

\vspace{-1mm}
$$
\le C_k
\sum_{(t_1,\ldots,t_k)}
\int\limits_{t}^{T}
\ldots
\int\limits_{t}^{t_2}
\left(K(t_1,\ldots,t_k){\bf 1}_{\{t_k<s\}}-
\sum_{j_1=0}^{p_1}\ldots
\sum_{j_k=0}^{p_k}
C_{j_k\ldots j_1}(s)
\prod_{l=1}^k\phi_{j_l}(t_l)\right)^2\times
$$

$$
\times
dt_1
\ldots
dt_k=
$$

\vspace{-2mm}
$$
=C_k\int\limits_{[t,T]^k}
\left(\bar K(t_1,\ldots,t_k,s)-
\sum_{j_1=0}^{p_1}\ldots
\sum_{j_k=0}^{p_k}
C_{j_k\ldots j_1}(s)
\prod_{l=1}^k\phi_{j_l}(t_l)\right)^2
dt_1
\ldots
dt_k\to 0
$$

\vspace{2mm}
\noindent
if $p_1,\ldots,p_k\to\infty,$ where constant $C_k$ 
depends only
on the multiplicity $k$ of the iterated It\^{o} stochastic integral
$J[\psi^{(k)}]_{s,t}$. 
Theorem 1.11 is proved.

{\bf Remark 1.9.} {\it Obviously from Theorem {\rm 1.11} for the case $s=T$
we obtain The\-o\-rem {\rm 1.1}.}

{\bf Remark 1.10.} {\it It is not difficult to see that Theorem {\rm 1.11} is valid for the case
when $\{\phi_j(x)\}_{j=0}^{\infty}$ is a complete orthonormal system  
of functions in the space $L_2([t,T]),$ 
each function $\phi_j(x)$ of which 
for finite $j$ satisfies the condition 
$(\star)$ {\rm (}see Sect.~{\rm 1.1.7} for details{\rm )}.}

From Theorem 1.11 for the case of pairwise different numbers
$i_1,\ldots,i_k=1,\ldots,m$ we obtain
\begin{equation}
\label{xyzyx300}
~~~~~~~~ J[\psi^{(k)}]_{s,t}=
\hbox{\vtop{\offinterlineskip\halign{
\hfil#\hfil\cr
{\rm l.i.m.}\cr
$\stackrel{}{{}_{p_1,\ldots,p_k\to \infty}}$\cr
}} }\sum_{j_1=0}^{p_1}\ldots\sum_{j_k=0}^{p_k}
C_{j_k\ldots j_1}(s)\zeta_{j_1}^{(i_1)}\ldots \zeta_{j_k}^{(i_k)}.
\end{equation}

\vspace{1mm}

Note that the expression on the right-hand side of (\ref{xyzyx300})
coincides for the case $k=1,$ $\psi_1(t_1)\equiv 1$ with the
right-hand side of the formula (\ref{jq1}) (approximation of the increment of the Wiener
process based on its series expansion).

\vspace{2mm}

{\bf Remark 1.11.}\ {\it Note that
by analogy with the proof of estimate {\rm (\ref{ddd1})} we obtain the following
inequality
\begin{equation}
\label{road900}
~~~~~ \int\limits_{[t,T]^k}\bar K^2(t_1,\ldots,t_k,s)dt_1\ldots dt_k-
\sum_{j_1=0}^{p}\ldots \sum_{j_k=0}^{p}
C^2_{j_k\ldots j_1}(s)\le \frac{G_k(s)}{p},
\end{equation}

\vspace{2mm}
\noindent
where $\bar K(t_1,\ldots,t_k,s)$ and  $C_{j_k\ldots j_1}(s)$
are defined by the equalities {\rm (\ref{road901})} and {\rm (\ref{ppppaxyz})}, respectively{\rm ;}  
constant $G_k(s)$ depends on $k$ and $s-t$ $(s\in (t, T],$ $s$ is fixed{\rm ).}}

\vspace{2mm}

The following obvious modification of Theorem {\rm 1.4} takes place.

\vspace{2mm}

{\bf Theorem 1.12.}\ {\it Suppose that
every $\psi_l(\tau)$ $(l=1,\ldots, k)$ is a continuous nonrandom function on 
$[t, T]$ and
$\{\phi_j(x)\}_{j=0}^{\infty}$ is a complete orthonormal system  
of functions in the space $L_2([t,T]),$ each function $\phi_j(x)$ of which 
for finite $j$ satisfies the condition 
$(\star)$ {\rm (}see Sect.~{\rm 1.1.7)}.
Then 

\vspace{-1mm}
$$
{\sf M}\left\{\left(
J[\psi^{(k)}]_{s,t}-J[\psi^{(k)}]_{s,t}^{p_1,\ldots,p_k}
\right)^2\right\}
\le 
$$

\vspace{-2mm}
\begin{equation}
\label{road888}
~~~~\le C_k(s)\left(~\int\limits_{[t,T]^k}
\bar K^2(t_1,\ldots,t_k,s)
dt_1\ldots dt_k -\sum_{j_1=0}^{p_1}\ldots
\sum_{j_k=0}^{p_k}C^2_{j_k\ldots j_1}(s)\right),
\end{equation}

\vspace{2mm}
\noindent
where $i_1,\ldots,i_k=0, 1,\ldots,m,$ constant $C_k(s)$ depends only on $k$ and $s-t.$ 
Moreover$,$ $C_k(s)\le k!$ for the following cases{\rm :}

\vspace{1mm}

{\rm 1.}\ $i_1,\ldots,i_k=1,\ldots,m$\ \ and\ \ $0<T-t<\infty,$

{\rm 2.}\ $i_1,\ldots,i_k=0, 1,\ldots,m,$\ \ $i_1^2+\ldots+i_k^2>0,$\ \
and\ \ $0<T-t<1,$

\vspace{2mm}
\noindent
where $J[\psi^{(k)}]_{s,t}$ is the stochastic integral {\rm (\ref{opr22}),}
$J[\psi^{(k)}]_{s,t}^{p_1,\ldots,p_k}$ is the 
expression on the right-hand side of {\rm (\ref{road777})} before
passing to the limit 
$\hbox{\vtop{\offinterlineskip\halign{
\hfil#\hfil\cr
{\rm l.i.m.}\cr
$\stackrel{}{{}_{p_1,\ldots,p_k\to \infty}}$\cr
}} },$ 
$\bar K(t_1,\ldots,t_k,s)$ and  $C_{j_k\ldots j_1}(s)$
are defined by the equalities {\rm (\ref{road901})} 
and {\rm (\ref{ppppaxyz}),} respectively{\rm ;}   $s\in (t, T]$ {\rm (}$s$ is fixed{\rm );} another 
notations are the same as in Theorem {\rm 1.11}.
}

\vspace{2mm}

{\bf Remark 1.12.}\ {\it Combining the estimates {\rm (\ref{road900})} and {\rm (\ref{road888}),}
we obtain

\vspace{-1mm}
\begin{equation}
\label{road1888}
{\sf M}\left\{\left(
J[\psi^{(k)}]_{s,t}-J[\psi^{(k)}]_{s,t}^{p,\ldots,p}
\right)^2\right\}
\le \frac{Q_k(s)}{p},
\end{equation}

\vspace{3mm}
\noindent
where $i_1,\ldots,i_k=0, 1,\ldots,m,$ constant $Q_k(s)$ depends only on $k$ and $s-t;$
another notations are the same as in {\rm (\ref{road900})} and {\rm (\ref{road888})}.}

\vspace{2mm}

{\bf Remark 1.13.}\ {\it An analogue of the estimate {\rm (\ref{2026ch1001s11})} for
the iterated It\^{o} stochastic integral {\rm (\ref{opr22})} has the
following form

\newpage
\noindent
$$
{\sf M}\left\{\left(J[\psi^{(k)}]_{s,t}-
J[\psi^{(k)}]_{s,t}^{p_1,\ldots,p_k}\right)^{2n}\right\}\le
$$

\vspace{2mm}
$$
\le
(k!)^{n} (2n-1)^{nk}\ \times
$$

\vspace{-2mm}
\begin{equation}
\label{agent01000}
~~~~~~~\times\ 
\left(~
\int\limits_{[t,T]^k}
\bar K^2(t_1,\ldots,t_k,s)
dt_1\ldots dt_k -\sum_{j_1=0}^{p_1}\ldots
\sum_{j_k=0}^{p_k}C^2_{j_k\ldots j_1}(s)
\right)^n,
\end{equation}

\vspace{4mm}
\noindent
where $J[\psi^{(k)}]_{s,t}^{p_1,\ldots,p_k}$ is the 
expression on the right-hand side of {\rm (\ref{road777})} before
passing to the limit 
$$
\hbox{\vtop{\offinterlineskip\halign{
\hfil#\hfil\cr
{\rm l.i.m.}\cr
$\stackrel{}{{}_{p_1,\ldots,p_k\to \infty}}$\cr
}} },
$$ 

\vspace{2mm}
\noindent
$\bar K(t_1,\ldots,t_k,s)$ and  $C_{j_k\ldots j_1}(s)$
are defined by the equalities {\rm (\ref{road901})} and 
{\rm (\ref{ppppaxyz}),} respectively{\rm ;}   $s\in (t, T]$ {\rm (}$s$ is fixed{\rm );}
$i_1,\ldots,i_k=1,\ldots,m.$}

\vspace{2mm}

{\bf Remark 1.14.} {\it The estimates
{\rm (\ref{road900})} and 
{\rm (\ref{agent01000})} imply the following 
inequality

\vspace{-6mm}
$$
{\sf M}\left\{\left(
J[\psi^{(k)}]_{s,t}-J[\psi^{(k)}]_{s,t}^{p,\ldots,p}
\right)^{2n}\right\}\le 
$$

\vspace{2mm}
$$
\le (k!)^{n} (2n-1)^{nk}\
\frac{\left(P_k\right)^n (s-t)^{nk}}{p^n},
$$

\vspace{4mm}
\noindent
where $i_1,\ldots,i_k=1,\ldots,m,$\ $n\in{\bf N},$ and 
constant $P_k$ depends only on $k$.}

\subsection{Expansions of Iterated It\^{o} Stochastic 
Integrals with Multiplicities 
1 to 5 and Miltiplicity $k$ Based on Theorem 1.11}

Consider particular cases of Theorem 1.11 for 
$k=1,\ldots,5$

\begin{equation}
\label{a1uuu}
J[\psi^{(1)}]_{s,t}
=\hbox{\vtop{\offinterlineskip\halign{
\hfil#\hfil\cr
{\rm l.i.m.}\cr
$\stackrel{}{{}_{p_1\to \infty}}$\cr
}} }\sum_{j_1=0}^{p_1}
C_{j_1}(s)\zeta_{j_1}^{(i_1)},
\end{equation}

\begin{equation}
\label{a2xxx}
~~~~J[\psi^{(2)}]_{s,t}
=\hbox{\vtop{\offinterlineskip\halign{
\hfil#\hfil\cr
{\rm l.i.m.}\cr
$\stackrel{}{{}_{p_1,p_2\to \infty}}$\cr
}} }\sum_{j_1=0}^{p_1}\sum_{j_2=0}^{p_2}
C_{j_2j_1}(s)\Biggl(\zeta_{j_1}^{(i_1)}\zeta_{j_2}^{(i_2)}
-{\bf 1}_{\{i_1=i_2\ne 0\}}
{\bf 1}_{\{j_1=j_2\}}\Biggr),
\end{equation}

\newpage
\noindent
$$
J[\psi^{(3)}]_{s,t}=
\hbox{\vtop{\offinterlineskip\halign{
\hfil#\hfil\cr
{\rm l.i.m.}\cr
$\stackrel{}{{}_{p_1,p_2,p_3\to \infty}}$\cr
}} }\sum_{j_1=0}^{p_1}\sum_{j_2=0}^{p_2}\sum_{j_3=0}^{p_3}
C_{j_3j_2j_1}(s)\Biggl(
\zeta_{j_1}^{(i_1)}\zeta_{j_2}^{(i_2)}\zeta_{j_3}^{(i_3)}
-\Biggr.
$$
\begin{equation}
\label{result3}
~-\Biggl.
{\bf 1}_{\{i_1=i_2\ne 0\}}
{\bf 1}_{\{j_1=j_2\}}
\zeta_{j_3}^{(i_3)}
-{\bf 1}_{\{i_2=i_3\ne 0\}}
{\bf 1}_{\{j_2=j_3\}}
\zeta_{j_1}^{(i_1)}-
{\bf 1}_{\{i_1=i_3\ne 0\}}
{\bf 1}_{\{j_1=j_3\}}
\zeta_{j_2}^{(i_2)}\Biggr),
\end{equation}

\vspace{4mm}

$$
J[\psi^{(4)}]_{s,t}
=
\hbox{\vtop{\offinterlineskip\halign{
\hfil#\hfil\cr
{\rm l.i.m.}\cr
$\stackrel{}{{}_{p_1,\ldots,p_4\to \infty}}$\cr
}} }\sum_{j_1=0}^{p_1}\ldots\sum_{j_4=0}^{p_4}
C_{j_4\ldots j_1}(s)\Biggl(
\prod_{l=1}^4\zeta_{j_l}^{(i_l)}
\Biggr.
-
$$
$$
-
{\bf 1}_{\{i_1=i_2\ne 0\}}
{\bf 1}_{\{j_1=j_2\}}
\zeta_{j_3}^{(i_3)}
\zeta_{j_4}^{(i_4)}
-
{\bf 1}_{\{i_1=i_3\ne 0\}}
{\bf 1}_{\{j_1=j_3\}}
\zeta_{j_2}^{(i_2)}
\zeta_{j_4}^{(i_4)}-
$$
$$
-
{\bf 1}_{\{i_1=i_4\ne 0\}}
{\bf 1}_{\{j_1=j_4\}}
\zeta_{j_2}^{(i_2)}
\zeta_{j_3}^{(i_3)}
-
{\bf 1}_{\{i_2=i_3\ne 0\}}
{\bf 1}_{\{j_2=j_3\}}
\zeta_{j_1}^{(i_1)}
\zeta_{j_4}^{(i_4)}-
$$
$$
-
{\bf 1}_{\{i_2=i_4\ne 0\}}
{\bf 1}_{\{j_2=j_4\}}
\zeta_{j_1}^{(i_1)}
\zeta_{j_3}^{(i_3)}
-
{\bf 1}_{\{i_3=i_4\ne 0\}}
{\bf 1}_{\{j_3=j_4\}}
\zeta_{j_1}^{(i_1)}
\zeta_{j_2}^{(i_2)}+
$$
$$
+
{\bf 1}_{\{i_1=i_2\ne 0\}}
{\bf 1}_{\{j_1=j_2\}}
{\bf 1}_{\{i_3=i_4\ne 0\}}
{\bf 1}_{\{j_3=j_4\}}
+
{\bf 1}_{\{i_1=i_3\ne 0\}}
{\bf 1}_{\{j_1=j_3\}}
{\bf 1}_{\{i_2=i_4\ne 0\}}
{\bf 1}_{\{j_2=j_4\}}+
$$
\begin{equation}
\label{cas1}
+\Biggl.
{\bf 1}_{\{i_1=i_4\ne 0\}}
{\bf 1}_{\{j_1=j_4\}}
{\bf 1}_{\{i_2=i_3\ne 0\}}
{\bf 1}_{\{j_2=j_3\}}\Biggr),
\end{equation}

\vspace{6mm}

$$
J[\psi^{(5)}]_{s,t}
=\hbox{\vtop{\offinterlineskip\halign{
\hfil#\hfil\cr
{\rm l.i.m.}\cr
$\stackrel{}{{}_{p_1,\ldots,p_5\to \infty}}$\cr
}} }\sum_{j_1=0}^{p_1}\ldots\sum_{j_5=0}^{p_5}
C_{j_5\ldots j_1}(s)\Biggl(
\prod_{l=1}^5\zeta_{j_l}^{(i_l)}
-\Biggr.
$$
$$
-
{\bf 1}_{\{i_1=i_2\ne 0\}}
{\bf 1}_{\{j_1=j_2\}}
\zeta_{j_3}^{(i_3)}
\zeta_{j_4}^{(i_4)}
\zeta_{j_5}^{(i_5)}-
{\bf 1}_{\{i_1=i_3\ne 0\}}
{\bf 1}_{\{j_1=j_3\}}
\zeta_{j_2}^{(i_2)}
\zeta_{j_4}^{(i_4)}
\zeta_{j_5}^{(i_5)}-
$$
$$
-
{\bf 1}_{\{i_1=i_4\ne 0\}}
{\bf 1}_{\{j_1=j_4\}}
\zeta_{j_2}^{(i_2)}
\zeta_{j_3}^{(i_3)}
\zeta_{j_5}^{(i_5)}-
{\bf 1}_{\{i_1=i_5\ne 0\}}
{\bf 1}_{\{j_1=j_5\}}
\zeta_{j_2}^{(i_2)}
\zeta_{j_3}^{(i_3)}
\zeta_{j_4}^{(i_4)}-
$$
$$
-
{\bf 1}_{\{i_2=i_3\ne 0\}}
{\bf 1}_{\{j_2=j_3\}}
\zeta_{j_1}^{(i_1)}
\zeta_{j_4}^{(i_4)}
\zeta_{j_5}^{(i_5)}-
{\bf 1}_{\{i_2=i_4\ne 0\}}
{\bf 1}_{\{j_2=j_4\}}
\zeta_{j_1}^{(i_1)}
\zeta_{j_3}^{(i_3)}
\zeta_{j_5}^{(i_5)}-
$$
$$
-
{\bf 1}_{\{i_2=i_5\ne 0\}}
{\bf 1}_{\{j_2=j_5\}}
\zeta_{j_1}^{(i_1)}
\zeta_{j_3}^{(i_3)}
\zeta_{j_4}^{(i_4)}
-{\bf 1}_{\{i_3=i_4\ne 0\}}
{\bf 1}_{\{j_3=j_4\}}
\zeta_{j_1}^{(i_1)}
\zeta_{j_2}^{(i_2)}
\zeta_{j_5}^{(i_5)}-
$$
$$
-
{\bf 1}_{\{i_3=i_5\ne 0\}}
{\bf 1}_{\{j_3=j_5\}}
\zeta_{j_1}^{(i_1)}
\zeta_{j_2}^{(i_2)}
\zeta_{j_4}^{(i_4)}
-{\bf 1}_{\{i_4=i_5\ne 0\}}
{\bf 1}_{\{j_4=j_5\}}
\zeta_{j_1}^{(i_1)}
\zeta_{j_2}^{(i_2)}
\zeta_{j_3}^{(i_3)}+
$$
$$
+
{\bf 1}_{\{i_1=i_2\ne 0\}}
{\bf 1}_{\{j_1=j_2\}}
{\bf 1}_{\{i_3=i_4\ne 0\}}
{\bf 1}_{\{j_3=j_4\}}\zeta_{j_5}^{(i_5)}+
{\bf 1}_{\{i_1=i_2\ne 0\}}
{\bf 1}_{\{j_1=j_2\}}
{\bf 1}_{\{i_3=i_5\ne 0\}}
{\bf 1}_{\{j_3=j_5\}}\zeta_{j_4}^{(i_4)}+
$$
$$
+
{\bf 1}_{\{i_1=i_2\ne 0\}}
{\bf 1}_{\{j_1=j_2\}}
{\bf 1}_{\{i_4=i_5\ne 0\}}
{\bf 1}_{\{j_4=j_5\}}\zeta_{j_3}^{(i_3)}+
{\bf 1}_{\{i_1=i_3\ne 0\}}
{\bf 1}_{\{j_1=j_3\}}
{\bf 1}_{\{i_2=i_4\ne 0\}}
{\bf 1}_{\{j_2=j_4\}}\zeta_{j_5}^{(i_5)}+
$$
$$
+
{\bf 1}_{\{i_1=i_3\ne 0\}}
{\bf 1}_{\{j_1=j_3\}}
{\bf 1}_{\{i_2=i_5\ne 0\}}
{\bf 1}_{\{j_2=j_5\}}\zeta_{j_4}^{(i_4)}+
{\bf 1}_{\{i_1=i_3\ne 0\}}
{\bf 1}_{\{j_1=j_3\}}
{\bf 1}_{\{i_4=i_5\ne 0\}}
{\bf 1}_{\{j_4=j_5\}}\zeta_{j_2}^{(i_2)}+
$$
$$
+
{\bf 1}_{\{i_1=i_4\ne 0\}}
{\bf 1}_{\{j_1=j_4\}}
{\bf 1}_{\{i_2=i_3\ne 0\}}
{\bf 1}_{\{j_2=j_3\}}\zeta_{j_5}^{(i_5)}+
{\bf 1}_{\{i_1=i_4\ne 0\}}
{\bf 1}_{\{j_1=j_4\}}
{\bf 1}_{\{i_2=i_5\ne 0\}}
{\bf 1}_{\{j_2=j_5\}}\zeta_{j_3}^{(i_3)}+
$$
$$
+
{\bf 1}_{\{i_1=i_4\ne 0\}}
{\bf 1}_{\{j_1=j_4\}}
{\bf 1}_{\{i_3=i_5\ne 0\}}
{\bf 1}_{\{j_3=j_5\}}\zeta_{j_2}^{(i_2)}+
{\bf 1}_{\{i_1=i_5\ne 0\}}
{\bf 1}_{\{j_1=j_5\}}
{\bf 1}_{\{i_2=i_3\ne 0\}}
{\bf 1}_{\{j_2=j_3\}}\zeta_{j_4}^{(i_4)}+
$$
$$
+
{\bf 1}_{\{i_1=i_5\ne 0\}}
{\bf 1}_{\{j_1=j_5\}}
{\bf 1}_{\{i_2=i_4\ne 0\}}
{\bf 1}_{\{j_2=j_4\}}\zeta_{j_3}^{(i_3)}+
{\bf 1}_{\{i_1=i_5\ne 0\}}
{\bf 1}_{\{j_1=j_5\}}
{\bf 1}_{\{i_3=i_4\ne 0\}}
{\bf 1}_{\{j_3=j_4\}}\zeta_{j_2}^{(i_2)}+
$$
$$
+
{\bf 1}_{\{i_2=i_3\ne 0\}}
{\bf 1}_{\{j_2=j_3\}}
{\bf 1}_{\{i_4=i_5\ne 0\}}
{\bf 1}_{\{j_4=j_5\}}\zeta_{j_1}^{(i_1)}+
{\bf 1}_{\{i_2=i_4\ne 0\}}
{\bf 1}_{\{j_2=j_4\}}
{\bf 1}_{\{i_3=i_5\ne 0\}}
{\bf 1}_{\{j_3=j_5\}}\zeta_{j_1}^{(i_1)}+
$$
$$
+\Biggl.
{\bf 1}_{\{i_2=i_5\ne 0\}}
{\bf 1}_{\{j_2=j_5\}}
{\bf 1}_{\{i_3=i_4\ne 0\}}
{\bf 1}_{\{j_3=j_4\}}\zeta_{j_1}^{(i_1)}\Biggr),
$$

\vspace{4mm}
\noindent
where ${\bf 1}_A$ is the indicator of the set $A,$
$C_{j_k\ldots j_1}(s)$ $(k=1,\ldots,5)$ has the form (\ref{ppppaxyz}),
$s\in (t, T]$ ($s$ is fixed).

Consider a generalization of the above formulas 
$$
J[\psi^{(k)}]_{s,t}=
\hbox{\vtop{\offinterlineskip\halign{
\hfil#\hfil\cr
{\rm l.i.m.}\cr
$\stackrel{}{{}_{p_1,\ldots,p_k\to \infty}}$\cr
}} }
\sum\limits_{j_1=0}^{p_1}\ldots
\sum\limits_{j_k=0}^{p_k}
C_{j_k\ldots j_1}(s)\Biggl(
\prod_{l=1}^k\zeta_{j_l}^{(i_l)}+\sum\limits_{r=1}^{[k/2]}
(-1)^r \times
\Biggr.
$$

\vspace{-2mm}
$$
\times
\sum_{\stackrel{(\{\{g_1, g_2\}, \ldots, 
\{g_{2r-1}, g_{2r}\}\}, \{q_1, \ldots, q_{k-2r}\})}
{{}_{\{g_1, g_2, \ldots, 
g_{2r-1}, g_{2r}, q_1, \ldots, q_{k-2r}\}=\{1, 2, \ldots, k\}}}}
\prod\limits_{s=1}^r
{\bf 1}_{\{i_{g_{{}_{2s-1}}}=~i_{g_{{}_{2s}}}\ne 0\}}
\Biggl.{\bf 1}_{\{j_{g_{{}_{2s-1}}}=~j_{g_{{}_{2s}}}\}}
\prod_{l=1}^{k-2r}\zeta_{j_{q_l}}^{(i_{q_l})}\Biggr),
$$

\vspace{4mm}
\noindent
where $k\in {\bf N}$, $C_{j_k\ldots j_1}(s)$ has the form (\ref{ppppaxyz});
another notations are the same as in Theorem~1.2.

\section{Expansion of Multiple Wiener Stochastic Integral Based on 
Generalized Multiple Fourier Series}

Let us consider the 
multiple stochastic integral (\ref{mult11})
\begin{equation}
\label{mult11www}
~~~~~~~~ \hbox{\vtop{\offinterlineskip\halign{
\hfil#\hfil\cr
{\rm l.i.m.}\cr
$\stackrel{}{{}_{N\to \infty}}$\cr
}} }
\sum\limits_{\stackrel{j_1,\ldots,j_k=0}{{}_{j_q\ne j_r;\ q\ne r;\ 
q, r=1,\ldots, k}}}^{N-1}
\Phi\left(\tau_{j_1},\ldots,\tau_{j_k}\right)
\prod\limits_{l=1}^k
\Delta{\bf w}_{\tau_{j_l}}^{(i_l)}
\stackrel{\rm def}{=}J'[\Phi]_{T,t}^{(k)},
\end{equation}

\noindent
where for simplicity we assume that
$\Phi(t_1,\ldots,t_k):\ [t, T]^k\to{\bf R}^1$ is a 
continuous nonrandom
function on $[t, T]^k.$ Moreover, 
$\left\{\tau_{j}\right\}_{j=0}^{N}$ is a partition of
$[t,T],$ which satisfies the condition {\rm (\ref{1111})}.

The stochastic integral with respect to the scalar standard Wiener process
($i_1=\ldots=i_k\ne 0$)
and similar to (\ref{mult11www}) was considered in \cite{ito1951}
and is called the multiple Wiener stochastic integral \cite{ito1951}.
Note that $\Phi(t_1,\ldots,t_k)\in L_2([t, T]^k)$ in \cite{ito1951}
(this case will be considered in Sect.~1.11, 1.12).

Consider the following theorem on expansion of the multiple
Wiener stochastic integral (\ref{mult11www}) based on generalized multiple Fourier series.

\footnotetext[8]{Theorem~1.13 
will be generalized to the case of an arbitrary 
complete ortho\-nor\-mal system of functions $\{\phi_j(x)\}_{j=0}^{\infty}$ 
in the space $L_2([t, T])$
and $\Phi(t_1,\ldots,t_k) \in L_2([t, T]^k)$ (see Sect.~1.11,
Theorem~1.17).}

{\bf Theorem 1.13.}${}^8$\ 
{\it Suppose that $\Phi(t_1,\ldots,t_k):\ [t, T]^k\to{\bf R}^1$ is a 
continuous nonrandom
function on $[t, T]^k$ and
$\{\phi_j(x)\}_{j=0}^{\infty}$ is a complete orthonormal system  
of functions in the space $L_2([t,T]),$ each function $\phi_j(x)$ of which 
for finite $j$ satisfies the condition 
$(\star)$ {\rm (}see Sect.~{\rm 1.1.7)}.
Then the following expansions

\vspace{-2mm}
$$
J'[\Phi]_{T,t}^{(k)} =
\hbox{\vtop{\offinterlineskip\halign{
\hfil#\hfil\cr
{\rm l.i.m.}\cr
$\stackrel{}{{}_{p_1,\ldots,p_k\to \infty}}$\cr
}} }\sum_{j_1=0}^{p_1}\ldots\sum_{j_k=0}^{p_k}
C_{j_k\ldots j_1}\Biggl(
\prod_{l=1}^k\zeta_{j_l}^{(i_l)} -
\Biggr.
$$

\begin{equation}
\label{quq50}
~~~~~~~~~-\Biggl.
\hbox{\vtop{\offinterlineskip\halign{
\hfil#\hfil\cr
{\rm l.i.m.}\cr
$\stackrel{}{{}_{N\to \infty}}$\cr
}} }\sum_{(l_1,\ldots,l_k)\in {\rm G}_k}
\phi_{j_{1}}(\tau_{l_1})
\Delta{\bf w}_{\tau_{l_1}}^{(i_1)}\ldots
\phi_{j_{k}}(\tau_{l_k})
\Delta{\bf w}_{\tau_{l_k}}^{(i_k)}\Biggr),
\end{equation}

\vspace{3mm}
$$
J'[\Phi]_{T,t}^{(k)}=
\hbox{\vtop{\offinterlineskip\halign{
\hfil#\hfil\cr
{\rm l.i.m.}\cr
$\stackrel{}{{}_{p_1,\ldots,p_k\to \infty}}$\cr
}} }
\sum\limits_{j_1=0}^{p_1}\ldots
\sum\limits_{j_k=0}^{p_k}
C_{j_k\ldots j_1}\Biggl(
\prod_{l=1}^k\zeta_{j_l}^{(i_l)}+\sum\limits_{r=1}^{[k/2]}
(-1)^r \times
\Biggr.
$$

\vspace{2mm}
\begin{equation}
\label{quq11}
\times
\sum_{\stackrel{(\{\{g_1, g_2\}, \ldots, 
\{g_{2r-1}, g_{2r}\}\}, \{q_1, \ldots, q_{k-2r}\})}
{{}_{\{g_1, g_2, \ldots, 
g_{2r-1}, g_{2r}, q_1, \ldots, q_{k-2r}\}=\{1, 2, \ldots, k\}}}}
\prod\limits_{s=1}^r
{\bf 1}_{\{i_{g_{{}_{2s-1}}}=~i_{g_{{}_{2s}}}\ne 0\}}
\Biggl.{\bf 1}_{\{j_{g_{{}_{2s-1}}}=~j_{g_{{}_{2s}}}\}}
\prod_{l=1}^{k-2r}\zeta_{j_{q_l}}^{(i_{q_l})}\Biggr)
\end{equation}

\vspace{1mm}
\noindent
con\-verg\-ing in the mean-square sense are valid$,$
where
$$
{\rm G}_k={\rm H}_k\backslash{\rm L}_k,\ \ \
{\rm H}_k=\bigl\{(l_1,\ldots,l_k):\ l_1,\ldots,l_k=0,\ 1,\ldots,N-1\bigr\},
$$
$$
{\rm L}_k=\bigl\{(l_1,\ldots,l_k):\ l_1,\ldots,l_k=0,\ 1,\ldots,N-1;\
l_g\ne l_r\ (g\ne r);\ g, r=1,\ldots,k\bigr\},
$$

\noindent
${\rm l.i.m.}$ is a limit in the mean-square sense$,$
$i_1,\ldots,i_k=0,1,\ldots,m,$ 
$$
\zeta_{j}^{(i)}=
\int\limits_t^T \phi_{j}(s) d{\bf w}_s^{(i)}
$$
are independent standard Gaussian random variables
for various
$i$ or $j$ {\rm(}in the case when $i\ne 0${\rm),}
\begin{equation}
\label{quq12}
C_{j_k\ldots j_1}=\int\limits_{[t,T]^k}
\Phi(t_1,\ldots,t_k)\prod_{l=1}^{k}\phi_{j_l}(t_l)dt_1\ldots dt_k
\end{equation}
is the Fourier coefficient$,$
$\Delta{\bf w}_{\tau_{j}}^{(i)}=
{\bf w}_{\tau_{j+1}}^{(i)}-{\bf w}_{\tau_{j}}^{(i)}$
$(i=0,\ 1,\ldots,m),$\
$\left\{\tau_{j}\right\}_{j=0}^{N}$ is a partition of
$[t,T],$ which satisfies the condition {\rm (\ref{1111});}
$[x]$ is an integer part of a real number $x;$
another notations are the same as in Theorem {\rm 1.2}.}

{\bf Proof.}\ Using Lemma 1.3 and (\ref{pobeda}), 
(\ref{s2s}),
we get the following representation
$$
J'[\Phi]_{T,t}^{(k)}=
$$

\vspace{-4mm}
$$
=
\sum_{j_1=0}^{p_1}\ldots
\sum_{j_k=0}^{p_k}
C_{j_k\ldots j_1}
\int\limits_{t}^{T}
\ldots
\int\limits_{t}^{t_2}
\sum_{(t_1,\ldots,t_k)}\left(
\phi_{j_1}(t_1)
\ldots
\phi_{j_k}(t_k)
d{\bf w}_{t_1}^{(i_1)}
\ldots
d{\bf w}_{t_k}^{(i_k)}\right)
+
$$

\vspace{2mm}
$$
+R_{T,t}^{p_1,\ldots,p_k}=
$$

$$
=\sum_{j_1=0}^{p_1}\ldots
\sum_{j_k=0}^{p_k}
C_{j_k\ldots j_1}\             
\hbox{\vtop{\offinterlineskip\halign{
\hfil#\hfil\cr
{\rm l.i.m.}\cr
$\stackrel{}{{}_{N\to \infty}}$\cr
}} }
\sum\limits_{\stackrel{l_1,\ldots,l_k=0}{{}_{l_q\ne l_r;\ 
q\ne r;\ q, r=1,\ldots, k}}}^{N-1}
\phi_{j_1}(\tau_{l_1})\ldots
\phi_{j_k}(\tau_{l_k})
\Delta{\bf w}_{\tau_{l_1}}^{(i_1)}
\ldots
\Delta{\bf w}_{\tau_{l_k}}^{(i_k)}+
$$

\vspace{4mm}
$$
+R_{T,t}^{p_1,\ldots,p_k}=
$$

\vspace{-2mm}
$$
=\sum_{j_1=0}^{p_1}\ldots
\sum_{j_k=0}^{p_k}
C_{j_k\ldots j_1}\left(
\hbox{\vtop{\offinterlineskip\halign{
\hfil#\hfil\cr
{\rm l.i.m.}\cr
$\stackrel{}{{}_{N\to \infty}}$\cr
}} }\sum_{l_1,\ldots,l_k=0}^{N-1}
\phi_{j_1}(\tau_{l_1})
\ldots
\phi_{j_k}(\tau_{l_k})
\Delta{\bf w}_{\tau_{l_1}}^{(i_1)}
\ldots
\Delta{\bf w}_{\tau_{l_k}}^{(i_k)}
-\right.
$$
$$
-\left.
\hbox{\vtop{\offinterlineskip\halign{
\hfil#\hfil\cr
{\rm l.i.m.}\cr
$\stackrel{}{{}_{N\to \infty}}$\cr
}} }\sum_{(l_1,\ldots,l_k)\in {\rm G}_k}
\phi_{j_{1}}(\tau_{l_1})
\Delta{\bf w}_{\tau_{l_1}}^{(i_1)}\ldots
\phi_{j_{k}}(\tau_{l_k})
\Delta{\bf w}_{\tau_{l_k}}^{(i_k)}\right)
+
$$

\vspace{2mm}
$$
+R_{T,t}^{p_1,\ldots,p_k}=
$$

\vspace{3mm}
$$
=\sum_{j_1=0}^{p_1}\ldots\sum_{j_k=0}^{p_k}
C_{j_k\ldots j_1}\times
$$

\vspace{-3mm}
$$
\times
\left(
\prod_{l=1}^k\zeta_{j_l}^{(i_l)}-
\hbox{\vtop{\offinterlineskip\halign{
\hfil#\hfil\cr
{\rm l.i.m.}\cr
$\stackrel{}{{}_{N\to \infty}}$\cr
}} }\sum_{(l_1,\ldots,l_k)\in {\rm G}_k}
\phi_{j_{1}}(\tau_{l_1})
\Delta{\bf w}_{\tau_{l_1}}^{(i_1)}\ldots
\phi_{j_{k}}(\tau_{l_k})
\Delta{\bf w}_{\tau_{l_k}}^{(i_k)}\right)+
$$
$$
+R_{T,t}^{p_1,\ldots,p_k}\ \ \ \hbox{w.~p.~1},
$$

\vspace{3mm}
\noindent
where
$$
R_{T,t}^{p_1,\ldots,p_k}
=\sum_{(t_1,\ldots,t_k)}
\int\limits_{t}^{T}
\ldots
\int\limits_{t}^{t_2}
\left(\Phi(t_1,\ldots,t_k)-
\sum_{j_1=0}^{p_1}\ldots
\sum_{j_k=0}^{p_k}
C_{j_k\ldots j_1}
\prod_{l=1}^k\phi_{j_l}(t_l)\right)\times
$$

$$
\times
d{\bf w}_{t_1}^{(i_1)}
\ldots
d{\bf w}_{t_k}^{(i_k)},
$$

\vspace{3mm}
\noindent
where permutations $(t_1,\ldots,t_k)$ when summing are performed only 
in the values $d{\bf w}_{t_1}^{(i_1)}
\ldots $
$d{\bf w}_{t_k}^{(i_k)}$. At the same time the indices near 
upper limits of integration in the iterated stochastic integrals 
are changed correspondently and if $t_r$ swapped with $t_q$ in the  
permutation $(t_1,\ldots,t_k)$, then $i_r$ swapped with $i_q$ in the 
permutation $(i_1,\ldots,i_k)$.

Let us estimate the remainder
$R_{T,t}^{p_1,\ldots,p_k}$ of the series using Lemma 1.2 and (\ref{riemann}). We have

\vspace{-2mm}
$$
{\sf M}\left\{\left(R_{T,t}^{p_1,\ldots,p_k}\right)^2\right\}
\le 
$$

\vspace{-2mm}
$$
\le C_k
\sum_{(t_1,\ldots,t_k)}
\int\limits_{t}^{T}
\ldots
\int\limits_{t}^{t_2}
\left(\Phi(t_1,\ldots,t_k)-
\sum_{j_1=0}^{p_1}\ldots
\sum_{j_k=0}^{p_k}
C_{j_k\ldots j_1}
\prod_{l=1}^k\phi_{j_l}(t_l)\right)^2\times
$$

\vspace{-2mm}
$$
\times
dt_1
\ldots
dt_k=
$$

\vspace{-2mm}
$$
=C_k\int\limits_{[t,T]^k}
\left(\Phi(t_1,\ldots,t_k)-
\sum_{j_1=0}^{p_1}\ldots
\sum_{j_k=0}^{p_k}
C_{j_k\ldots j_1}
\prod_{l=1}^k\phi_{j_l}(t_l)\right)^2\times
$$

\vspace{-2mm}
$$
\times
dt_1
\ldots
dt_k\to 0
$$

\vspace{4mm}
\noindent
if $p_1,\ldots,p_k\to\infty,$ where constant $C_k$ 
depends only
on the multiplicity $k$ of the multiple Wiener stochastic integral
$J'[\Phi]_{T,t}^{(k)}$.
The expansion (\ref{quq50}) is proved. Using (\ref{quq50}) and Remark 1.2,
we get the expansion (\ref{quq11}) (see Theorem 1.2).
Theorem 1.13 is proved.

Note that particular cases of the expansion (\ref{quq11})
are determined by the equalities (\ref{a1})--(\ref{a7}), in which the Fourier
coefficient $C_{j_k\ldots j_1}$ $(k=1,\ldots,7)$
has the form (\ref{quq12}).

\section{Reformulation of Theorems 1.1, 1.2, and 1.13 Using Hermite Polynomials}

In \cite{Rybakov1000} it was noted that
Theorem 3.1 (\cite{ito1951}, p.~162) can be applied
to the case of multiple Wiener stochastic integral
with respect to components of the multidimensional
Wiener process. As a result, Theorems 1.1, 1.2, and 1.13
can be reformulated 
using Hermite polynomials. Consider this approach 
using our notations.
Note that we derive the formula (\ref{ziko20}) (see below)
in two different ways. One of
them is not based on Theorem 3.1 \cite{ito1951} (see the proof of
Theorem~1.22 below for details).

{\it We will say that the condition {\rm ($\star\star$)} is fulfilled
for the multi-index $(i_1\ldots i_k)$ $(i_1,\ldots,i_k=0, 1,\ldots, m)$ if
$m_1,\ldots,m_k$ are multiplicities of the elements $i_1,\ldots,i_k,$ respectively$,$ i.e.
$$
\{i_1,\ldots, i_k\}\hspace{-0.4mm}=\hspace{-0.4mm}\{\overbrace{{i_1, \ldots, i_1}}^{m_1},
\overbrace{{i_2, \ldots, i_2}}^{m_2},
\ldots, \overbrace{{i_r, \ldots, i_r}}^{m_r}\}\ \ \ (m_{r+1}=\ldots=m_k=0),
$$
where $r=1,\ldots, k,$ braces   
mean an unordered 
set, and pa\-ren\-the\-ses mean an ordered set. At that, 
$m_1+\ldots+m_k=k,$\ $m_1,\ldots, m_k=0,1,\ldots,k,$\ 
and all elements with nonzero multiplicities are pairwise different.}

In this section, we consider the case $i_1,\ldots,i_k=0, 1,\ldots, m.$
Note that in \cite{Rybakov1000} the case
$i_1,\ldots,i_k=1,\ldots, m$ was considered.

Let the condition {\rm ($\star\star$)} is fulfilled
for the mul\-ti-\-in\-dex $(i_1 \ldots i_k).$ Then
$$
J'\left[\phi_{j_1}\ldots \phi_{j_k}\right]_{T,t}^{(i_1\ldots i_k)}=
J'\biggl[\underbrace{\phi_{j_{g_1}}
\ldots \phi_{j_{g_{{}_{m_1}}}}}_{m_1}
\underbrace{\phi_{j_{g_{m_1+1}}}
\ldots \phi_{j_{g_{m_1+m_2}}}}_{m_2}\ldots \biggr.
$$
\begin{equation}
\label{newe11}
~~~~~~\biggl.\ldots
\underbrace{\phi_{j_{g_{m_1+m_2+\ldots+m_{k-1}+1}}}\ldots
\phi_{j_{g_{m_1+m_2+\ldots+m_k}}}}_{m_k}\biggr]_{T,t}^
{(\overbrace{{}_{i_1 \ldots i_1}}^{m_1}
\overbrace{{}_{i_2 \ldots i_2}}^{m_2}
\ldots \overbrace{{}_{i_k \ldots i_k}}^{m_k})}
\end{equation}

\noindent
w.~p.~1, where 
$J'\left[\phi_{j_1}\ldots \phi_{j_k}\right]_{T,t}^{(i_1\ldots i_k)}$
is defined by (\ref{mult11}) (also see (\ref{mult11www})),
$\Phi(t_1,\ldots,t_k)=\phi_{j_1}(t_1)\ldots \phi_{j_k}(t_k),$
$\{\phi_j(x)\}_{j=0}^{\infty}$ is a complete orthonormal system  
of functions in the space $L_2([t,T]),$ each function $\phi_j(x)$ of which 
for finite $j$ satisfies the condition 
$(\star)$ (see Sect.~1.1.7),
$\{j_{g_1},\ldots,j_{g_{m_1+m_2+\ldots+m_k}}\}=\{j_1,\ldots,j_k\}$.

From (\ref{newe11}) we have 
$$
J'\left[\phi_{j_{1}}\ldots \phi_{j_{k}}\right]_{T,t}^{(i_1\ldots i_k)}=
J'\left[\phi_{j_{g_1}}
\ldots \phi_{j_{g_{{}_{m_1}}}}\right]_{T,t}^{
(\hspace{0.5mm}\overbrace{{}_{i_1 \ldots i_1}}^{m_1}\hspace{0.5mm})}
\cdot 
J'\left[\phi_{j_{g_{m_1+1}}}
\ldots \phi_{j_{g_{m_1+m_2}}}\right]_{T,t}^{
(\hspace{0.5mm}\overbrace{{}_{i_2 \ldots i_2}}^{m_2}\hspace{0.5mm})}
\cdot \ldots 
$$
\begin{equation}
\label{ziko30}
\ldots \cdot 
J'\left[\phi_{j_{g_{m_1+m_2+\ldots+m_{k-1}+1}}}\ldots
\phi_{j_{g_{m_1+m_2+\ldots+m_k}}}\right]_{T,t}^{
(\hspace{0.5mm}\overbrace{{}_{i_k \ldots i_k}}^{m_k}\hspace{0.5mm})}
\end{equation}

\noindent
w.~p.~1, where
\begin{equation}
\label{ziko10}
~~~~~~J'\left[\phi_{j_{g_{m_1+m_2+\ldots+m_{l-1}+1}}}\ldots
\phi_{j_{g_{m_1+m_2+\ldots+m_l}}}\right]_{T,t}^{
(\hspace{0.5mm}\overbrace{{}_{i_l \ldots i_l}}^{m_l}\hspace{0.5mm})}
\stackrel{\sf def}{=}1\ \ \ \hbox{for}\ \ \ m_l=0.
\end{equation}

The detailed proof of the equality (\ref{ziko30}) will be given in Sect.~1.14
(see the proof of Theorem~1.22).

Let us consider the following multiple Wiener stochastic integral 
$$
J'\left[\phi_{j_{g_{m_1+m_2+\ldots+m_{l-1}+1}}}\ldots
\phi_{j_{g_{m_1+m_2+\ldots+m_l}}}\right]_{T,t}^{
(\hspace{0.5mm}\overbrace{{}_{i_l \ldots i_l}}^{m_l}\hspace{0.5mm})}\ \ \ (m_l>0),
$$

\noindent
where we suppose that 
$$
\bigl\{j_{g_{m_1+m_2+\ldots+m_{l-1}+1}}, \ldots, j_{g_{m_1+m_2+\ldots+m_{l}}}
\bigr\}=
$$
\begin{equation}
\label{ziko999}
=\bigl\{\underbrace{j_{h_{1,l}}, \ldots, j_{h_{1,l}}}_{n_{1,l}}\ \hspace{-1mm},
\underbrace{j_{h_{2,l}}, \ldots, j_{h_{2,l}}}_{n_{2,l}}\ \hspace{-1mm}, \ldots,
\underbrace{j_{h_{d_l,l}}, \ldots, j_{h_{d_l,l}}}_{n_{d_l,l}}\bigr\},
\end{equation}
where
$n_{1,l}+n_{2,l}+\ldots+n_{d_l,l}=m_l;$\ \ $n_{1,l}, n_{2,l}, \ldots, n_{d_l,l}=1,\ldots, m_l;$\ \ 
$d_l=1,\ldots,m_l;$\ \ $l=1,\ldots,k.$ Note that the numbers $m_1,\ldots,m_k,$\ $g_1,\ldots,g_k$
depend on $(i_1,\ldots ,i_k)$ and the numbers
$n_{1,l},\ldots,n_{d_l,l},$\ $h_{1,l},\ldots,h_{d_l,l},$\ $d_l$
depend on $\{j_1,\ldots,j_k\}$. Moreover, 
$\left\{j_{g_1},\ldots,j_{g_k}\right\}
=\{j_1,\ldots,j_k\}.$

Using Theorem 3.1 \cite{ito1951}, we get w.~p.~1
$$
J'\left[\phi_{j_{g_{m_1+m_2+\ldots+m_{l-1}+1}}}\ldots
\phi_{j_{g_{m_1+m_2+\ldots+m_l}}}
\right]_{T,t}^{(\hspace{0.5mm}\overbrace{{}_{i_l \ldots i_l}}^{m_l}\hspace{0.5mm})}
=
$$

\vspace{-3mm}
\begin{equation}
\label{ziko20}
~~~~~~~ =\left\{
\begin{matrix}
H_{n_{1,l}}\left(\zeta_{j_{h_{1,l}}}^{(i_l)}\right)\ldots 
H_{n_{d_l,l}}\left(\zeta_{j_{h_{d_l,l}}}^{(i_l)}\right),\ 
&\hbox{\rm if}\ \ \ 
i_l\ne 0\cr\cr
\left(\zeta_{j_{h_{1,l}}}^{(0)}\right)^{n_{1,l}}\ldots
\left(\zeta_{j_{h_{d_l,l}}}^{(0)}\right)^{n_{d_l,l}},\  &\hbox{\rm if}\ \ \ 
i_l=0
\end{matrix}\right.\ \ \ (m_l>0),
\end{equation}

\vspace{2mm}
\noindent
where $H_n(x)$ is the Hermite polynomial of degree $n$
\begin{equation}
\label{ziko500}
H_n(x)=(-1)^n e^{x^2/2} \frac{d^n}{dx^n}\left(e^{-x^2/2}\right)
=n!\sum\limits_{m=0}^{[n/2]}\frac{(-1)^m x^{n-2m}}{m!(n-2m)! 2^m}\ \ \ (n\in{\bf N}),
\end{equation}

\vspace{1mm}
\noindent
and $\zeta_j^{(i)}$ $(i=0,1,\ldots,m,\ j=0,1,\ldots)$ is defined by (\ref{rr23}). 

For example,
$$
H_0(x)=1,\ \ \ H_1(x)=x,\ \ \ 
H_2(x)=x^2-1,
$$
$$
H_3(x)=x^3-3x,\ \ \
H_4(x)=x^4-6x^2 + 3,
$$

\vspace{-6mm}
$$
H_5(x)=x^5-10x^3 + 15x.
$$

\vspace{3mm}

From (\ref{ziko10}) and (\ref{ziko20}) we obtain w.~p.~1

$$
J'\left[\phi_{j_{g_{m_1+m_2+\ldots+m_{l-1}+1}}}\ldots
\phi_{j_{g_{m_1+m_2+\ldots+m_l}}}
\right]_{T,t}^{(\hspace{0.5mm}\overbrace{{}_{i_l \ldots i_l}}^{m_l}\hspace{0.5mm})}
=
$$

\vspace{-2mm}
\begin{equation}
\label{ziko40}
={\bf 1}_{\{m_l=0\}}+{\bf 1}_{\{m_l>0\}}\left\{
\begin{matrix}
H_{n_{1,l}}\left(\zeta_{j_{h_{1,l}}}^{(i_l)}\right)\ldots 
H_{n_{d_l,l}}\left(\zeta_{j_{h_{d_l,l}}}^{(i_l)}\right),\ 
&\hbox{\rm if}\ \ \ 
i_l\ne 0\cr\cr
\left(\zeta_{j_{h_{1,l}}}^{(0)}\right)^{n_{1,l}}\ldots
\left(\zeta_{j_{h_{d_l,l}}}^{(0)}\right)^{n_{d_l,l}},\  &\hbox{\rm if}\ \ \ 
i_l=0
\end{matrix}\right.,
\end{equation}

\vspace{4mm}
\noindent
where ${\bf 1}_A$ denotes the indicator of the set $A$.

Using (\ref{ziko30}) and (\ref{ziko40}), we get w.~p.~1

$$
J'\left[\phi_{j_{1}}\ldots \phi_{j_{k}}\right]_{T,t}^{(i_1\ldots i_k)}=
$$

\vspace{-2mm}
\begin{equation}
\label{ziko50}
=\prod_{l=1}^k\left({\bf 1}_{\{m_l=0\}}+{\bf 1}_{\{m_l>0\}}\left\{
\begin{matrix}
H_{n_{1,l}}\left(\zeta_{j_{h_{1,l}}}^{(i_l)}\right)\ldots 
H_{n_{d_l,l}}\left(\zeta_{j_{h_{d_l,l}}}^{(i_l)}\right),\ 
&\hbox{\rm if}\ \ \ 
i_l\ne 0\cr\cr
\left(\zeta_{j_{h_{1,l}}}^{(0)}\right)^{n_{1,l}}\ldots
\left(\zeta_{j_{h_{d_l,l}}}^{(0)}\right)^{n_{d_l,l}},\  &\hbox{\rm if}\ \ \ 
i_l=0
\end{matrix}\right.\ \right),
\end{equation}

\vspace{3mm}
\noindent
where notations are the same as in (\ref{ziko999}) and (\ref{ziko20}).

The equality (\ref{ziko50}) allows us to reformulate Theorems 1.1, 1.2, and 1.13
using the Hermite polynomials.${}^9$

\footnotetext[9]{Theorems~1.14, 1.15 (see below)
will be generalized to the case of an arbitrary 
complete ortho\-nor\-mal system of functions $\{\phi_j(x)\}_{j=0}^{\infty}$ 
in the space $L_2([t, T])$
and $\psi_1(\tau),$ $\ldots,\psi_k(\tau)\in L_2([t, T]),$ 
$\Phi(t_1,\ldots,t_k) \in L_2([t, T]^k)$ in Sect.~1.11
(see Theorems~1.16, 1.17).}

{\bf Theorem 1.14}\ \cite{arxiv-1} (reformulation of Theorems 1.1 and 1.2).\ {\it Suppose that
the condition {\rm ($\star\star$)} is fulfilled
for the multi-index $(i_1 \ldots i_k)$ and the condition {\rm (\ref{ziko999})} is also 
fulfilled.
Furthermore$,$ let 
every $\psi_l(\tau)$ $(l=$ $1,\ldots, k)$ is a continuous 
nonrandom function on 
$[t, T]$ and
$\{\phi_j(x)\}_{j=0}^{\infty}$ is a complete orthonormal system  
of functions in the space $L_2([t,T]),$ each function $\phi_j(x)$ of which 
for finite $j$ satisfies the condition 
$(\star)$ {\rm (}see Sect.~{\rm 1.1.7)}.
Then the following expansion

\vspace{-2mm}
$$
J[\psi^{(k)}]_{T,t}^{(i_1\ldots i_k)}=
\hbox{\vtop{\offinterlineskip\halign{
\hfil#\hfil\cr
{\rm l.i.m.}\cr
$\stackrel{}{{}_{p_1,\ldots,p_k\to \infty}}$\cr
}} }
\sum\limits_{j_1=0}^{p_1}\ldots
\sum\limits_{j_k=0}^{p_k}
C_{j_k\ldots j_1}\times
$$

\vspace{1.5mm}
\begin{equation}
\label{ziko800}
\times
\prod_{l=1}^k\left({\bf 1}_{\{m_l=0\}}+{\bf 1}_{\{m_l>0\}}\left\{
\begin{matrix}
H_{n_{1,l}}\left(\zeta_{j_{h_{1,l}}}^{(i_l)}\right)\ldots 
H_{n_{d_l,l}}\left(\zeta_{j_{h_{d_l,l}}}^{(i_l)}\right),\ 
&\hbox{\rm if}\ \ \ 
i_l\ne 0\cr\cr
\left(\zeta_{j_{h_{1,l}}}^{(0)}\right)^{n_{1,l}}\ldots
\left(\zeta_{j_{h_{d_l,l}}}^{(0)}\right)^{n_{d_l,l}},\  &\hbox{\rm if}\ \ \ 
i_l=0
\end{matrix}\right.\ \right)
\end{equation}

\vspace{3mm}
\noindent
con\-verg\-ing in the mean-square sense is valid$,$
where we denote the stochastic integral {\rm (\ref{ito})} as
$J[\psi^{(k)}]_{T,t}^{(i_1\ldots i_k)};$ 
$n_{1,l}+n_{2,l}+\ldots+n_{d_l,l}=m_l;$\ \ $n_{1,l}, n_{2,l}, \ldots, n_{d_l,l}=1,\ldots, m_l;$\ \ 
$d_l=1,\ldots,m_l;$\ \ $l=1,\ldots,k;$\ \ $m_1+\ldots+m_k=k;$\ \ 
the numbers $m_1,\ldots,m_k,$\ $g_1,\ldots,g_k$
depend on $(i_1,\ldots,i_k)$ and 
the numbers $n_{1,l},\ldots,n_{d_l,l},$\ $h_{1,l},\ldots,h_{d_l,l},$\ $d_l$
depend on $\{j_1,\ldots,j_k\};$ moreover$,$ $\left\{j_{g_1},\ldots,j_{g_k}\right\}
=\{j_1,\ldots,j_k\};$ $H_n(x)$ is the Hermite polynomial {\rm (\ref{ziko500});}
another
notations are the same as in Theorem {\rm 1.1}.}

{\bf Theorem 1.15}\ \cite{arxiv-1} (reformulation of Theorem 1.13).\ {\it Suppose that
the condition {\rm ($\star\star$)} is fulfilled
for the multi-index $(i_1 \ldots i_k)$ 
and the condition {\rm (\ref{ziko999})} is also 
fulfilled.
Furthermore$,$ let
$\Phi(t_1,\ldots,t_k):\ [t, T]^k\to{\bf R}^1$ is a 
continuous nonrandom
function on $[t, T]^k$
and
$\{\phi_j(x)\}_{j=0}^{\infty}$ is a complete orthonormal system  
of functions in the space $L_2([t,T]),$ each function $\phi_j(x)$ of which 
for finite $j$ satisfies the condition 
$(\star)$ {\rm (}see Sect.~{\rm 1.1.7)}.
Then the following expansion

\vspace{-2mm}
$$
J'[\Phi]_{T,t}^{(i_1\ldots i_k)}=
\hbox{\vtop{\offinterlineskip\halign{
\hfil#\hfil\cr
{\rm l.i.m.}\cr
$\stackrel{}{{}_{p_1,\ldots,p_k\to \infty}}$\cr
}} }
\sum\limits_{j_1=0}^{p_1}\ldots
\sum\limits_{j_k=0}^{p_k}
C_{j_k\ldots j_1}\times
$$

\vspace{1.5mm}
$$
\times
\prod_{l=1}^k\left({\bf 1}_{\{m_l=0\}}+{\bf 1}_{\{m_l>0\}}\left\{
\begin{matrix}
H_{n_{1,l}}\left(\zeta_{j_{h_{1,l}}}^{(i_l)}\right)\ldots 
H_{n_{d_l,l}}\left(\zeta_{j_{h_{d_l,l}}}^{(i_l)}\right),\ 
&\hbox{\rm if}\ \ \ 
i_l\ne 0\cr\cr
\left(\zeta_{j_{h_{1,l}}}^{(0)}\right)^{n_{1,l}}\ldots
\left(\zeta_{j_{h_{d_l,l}}}^{(0)}\right)^{n_{d_l,l}},\  &\hbox{\rm if}\ \ \ 
i_l=0
\end{matrix}\right.\ \right)
$$

\vspace{3mm}
\noindent
con\-verg\-ing in the mean-square sense is valid$,$
where we denote the multiple Wiener stochastic integral 
{\rm (\ref{mult11www})}
as
$J'[\Phi]_{T,t}^{(i_1\ldots i_k)};$ 
$n_{1,l}+n_{2,l}+\ldots+n_{d_l,l}=m_l;$\ \ $n_{1,l}, n_{2,l}, \ldots, n_{d_l,l}=1,\ldots, m_l;$\ \ 
$d_l=1,\ldots,m_l;$\ \ $l=1,\ldots,k;$\ \ $m_1+\ldots+m_k=k;$\ \ 
the numbers $m_1,\ldots,m_k,$\ $g_1,\ldots,g_k$
depend on $(i_1,\ldots,i_k)$ and 
the numbers $n_{1,l},\ldots,n_{d_l,l},$\ $h_{1,l},\ldots,h_{d_l,l},$\ $d_l$
depend on $\{j_1,\ldots,j_k\};$ moreover$,$ $\left\{j_{g_1},\ldots,j_{g_k}\right\}
=\{j_1,\ldots,j_k\};$ $H_n(x)$ is the Hermite polynomial {\rm (\ref{ziko500});}
another
notations are the same as in Theorem {\rm 1.13}.}

From (\ref{ziko40}) we have w.~p.~1
\begin{equation}
\label{ziko80}
~~~~~~~~~~ J'[\hspace{0.5mm}\underbrace{\phi_{j_1}\ldots
\phi_{j_1}}_{k}]_{T,t}^{(\hspace{0.5mm}\overbrace{{}_{i_1 \ldots i_1}}^{k}\hspace{0.5mm})}
=
\left\{
\begin{matrix}
H_{k}\left(\zeta_{j_{1}}^{(i_1)}\right),\ 
&\hbox{\rm if}\ \ \ 
i_1\ne 0\cr\cr
\left(\zeta_{j_{1}}^{(0)}\right)^{k},\  &\hbox{\rm if}\ \ \ 
i_1=0
\end{matrix}\right.\ \ \ (k>0) .
\end{equation}

\vspace{2mm}

Let us show how the relation (\ref{ziko80}) can be obtained from Theorem 1.2.
To prove (\ref{ziko80}) using Theorem 1.2 
we choose $i_1=\ldots=i_k$ and $j_1=\ldots=j_k$ $(i_1=0, 1,\ldots,m)$
in the following formula (this formula follows from a 
comparison of (\ref{drdr1}) and (\ref{leto6000}) or can be obtained
using the recurrence relation (\ref{recur1})) 
$$
J'[\phi_{j_1}\ldots \phi_{j_k}]_{T,t}^{(i_1\ldots i_k)}=\prod_{l=1}^k\zeta_{j_l}^{(i_l)}+
\sum\limits_{r=1}^{[k/2]}
(-1)^r \times
$$

\vspace{-2mm}
\begin{equation}
\label{rezo7}
\times \sum_{\stackrel{(\{\{g_1, g_2\}, \ldots, 
\{g_{2r-1}, g_{2r}\}\}, \{q_1, \ldots, q_{k-2r}\})}
{{}_{\{g_1, g_2, \ldots, 
g_{2r-1}, g_{2r}, q_1, \ldots, q_{k-2r}\}=\{1, 2, \ldots, k\}}}}
\prod\limits_{s=1}^r
{\bf 1}_{\{i_{g_{{}_{2s-1}}}=~i_{g_{{}_{2s}}}\ne 0\}}
\Biggl.{\bf 1}_{\{j_{g_{{}_{2s-1}}}=~j_{g_{{}_{2s}}}\}}
\prod_{l=1}^{k-2r}\zeta_{j_{q_l}}^{(i_{q_l})}
\end{equation}

\vspace{1mm}
\noindent
w.~p.~1, where notations are the same as in Theorem 1.2.

The case $i_1=0$ of (\ref{ziko80}) is obvious.
Simple combinatorial reasoning shows that 
$$
\sum_{\stackrel{(\{\{g_1, g_2\}, \ldots, 
\{g_{2r-1}, g_{2r}\}\}, \{q_1, \ldots, q_{k-2r}\})}
{{}_{\{g_1, g_2, \ldots, 
g_{2r-1}, g_{2r}, q_1, \ldots, q_{k-2r}\}=\{1, 2, \ldots, k\}}}}
\prod\limits_{s=1}^r
{\bf 1}_{\{i_{g_{{}_{2s-1}}}=~i_{g_{{}_{2s}}}\ne 0\}}
\Biggl.{\bf 1}_{\{j_{g_{{}_{2s-1}}}=~j_{g_{{}_{2s}}}\}}
\prod_{l=1}^{k-2r}\zeta_{j_{q_l}}^{(i_{q_l})}=
$$

\begin{equation}
\label{rezo15}
=\frac{C_k^2 \cdot C_{k-2}^2 \cdot \ldots \cdot 
C_{k-(r-1)2}^2}{r!} \left(\zeta_{j_{1}}^{(i_{1})}\right)^{k-2r},
\end{equation}

\vspace{3mm}
\noindent
where $i_1=\ldots=i_k,$\ \ $j_1=\ldots=j_k$\ $(i_1=1,\ldots,m),$ and
$$
C_k^l=\frac{k!}{l!(k-l)!}
$$
is the binomial coefficient.

We have

\vspace{-2mm}
\begin{equation}
\label{rezo16}
\frac{C_k^2 \cdot C_{k-2}^2 \cdot \ldots \cdot 
C_{k-(r-1)2}^2}{r!}=\frac{k!}{r!(k-2r)!2^r}.
\end{equation}

\vspace{3mm}

Combining (\ref{rezo7}), (\ref{rezo15}), and (\ref{rezo16}), we get w.~p.~1

\vspace{-2mm}
$$
J'[\underbrace{\phi_{j_1}
\ldots \phi_{j_1}}_{k}]_{T,t}^{(\hspace{0.5mm}\overbrace{{}_{i_1 \ldots i_1}}^{k}\hspace{0.5mm})}
=\left(\zeta_{j_1}^{(i_1)}\right)^k 
+ k!\sum\limits_{r=1}^{[k/2]}
\frac{(-1)^r}{r!(k-2r)!2^r}
\left(\zeta_{j_{1}}^{(i_{1})}\right)^{k-2r}=
$$

$$
=k!\sum\limits_{r=0}^{[k/2]}
\frac{(-1)^r}{r!(k-2r)!2^r}
\left(\zeta_{j_{1}}^{(i_{1})}\right)^{k-2r}=H_k\left(\zeta_{j_{1}}^{(i_{1})}\right).
$$

\vspace{4mm}
\noindent
The relation (\ref{ziko80}) is proved using (\ref{rezo7}).

From (\ref{ziko50}) and (\ref{rezo7}) we obtain the following equalities
for multiple Wiener stochastic integral

$$
J'[\phi_{j_1}\ldots \phi_{j_k}]_{T,t}^{(i_1\ldots i_k)}=\prod_{l=1}^k\zeta_{j_l}^{(i_l)}
+\sum\limits_{r=1}^{[k/2]}
(-1)^r \times
$$

\vspace{-3mm}
$$
\times\sum_{\stackrel{(\{\{g_1, g_2\}, \ldots, 
\{g_{2r-1}, g_{2r}\}\}, \{q_1, \ldots, q_{k-2r}\})}
{{}_{\{g_1, g_2, \ldots, 
g_{2r-1}, g_{2r}, q_1, \ldots, q_{k-2r}\}=\{1, 2, \ldots, k\}}}}
\prod\limits_{s=1}^r
{\bf 1}_{\{i_{g_{{}_{2s-1}}}=~i_{g_{{}_{2s}}}\ne 0\}}
\Biggl.{\bf 1}_{\{j_{g_{{}_{2s-1}}}=~j_{g_{{}_{2s}}}\}}
\prod_{l=1}^{k-2r}\zeta_{j_{q_l}}^{(i_{q_l})}=
$$

\begin{equation}
\label{ziko100}
=\prod_{l=1}^k\left({\bf 1}_{\{m_l=0\}}+{\bf 1}_{\{m_l>0\}}\left\{
\begin{matrix}
H_{n_{1,l}}\left(\zeta_{j_{h_{1,l}}}^{(i_l)}\right)\ldots 
H_{n_{d_l,l}}\left(\zeta_{j_{h_{d_l,l}}}^{(i_l)}\right),\ 
&\hbox{\rm if}\ \ \ 
i_l\ne 0\cr\cr
\left(\zeta_{j_{h_{1,l}}}^{(0)}\right)^{n_{1,l}}\ldots
\left(\zeta_{j_{h_{d_l,l}}}^{(0)}\right)^{n_{d_l,l}},\  &\hbox{\rm if}\ \ \ 
i_l=0
\end{matrix}\right.\ \right)
\end{equation}

\vspace{3mm}
\noindent
w.~p.~1, where notations are the same as in Theorem 1.2 and (\ref{ziko999}), (\ref{ziko20}).

Let us make a remark about how it is possible to obtain the formula
(\ref{ziko20}) without using Theorem 3.1 \cite{ito1951}. 

Consider the set of 
polynomials 
$H_n(x,y),$ $n=0, 1,\ldots$ defined by \cite{Ch}
\begin{equation}
\label{new1090}
\Biggl.H_n(x,y)=\left(\frac{d^n}{d\alpha^n} 
e^{\alpha x-\alpha^2 y/2}\right)
\Biggr|_{\alpha=0}\ \ \ (H_0(x,y)\stackrel{\sf def}{=}1).
\end{equation}

\vspace{1mm}

It is well known that polynomials $H_n(x,y)$ are connected with 
the Hermite polynomials (\ref{ziko500}) by the formula \cite{Ch}
\begin{equation}
\label{ziko1000}
H_n(x,y)=y^{n/2}
H_n\left(\frac{x}{\sqrt{y}}\right)=
n!\sum\limits_{i=0}^{[n/2]}\frac{(-1)^i x^{n-2i} y^i}{i!(n-2i)! 2^i}.
\end{equation}

For example,
$$
H_1(x,y)
=x,
$$
$$
H_2(x,y)
=x^2-y,
$$
$$
H_3(x,y)
=x^3-3xy,
$$
$$
H_4(x,y)
=x^4-6x^2 y
+3y^2,
$$
$$
H_5(x,y)=x^5-10x^3 y+15xy^2.
$$

\vspace{2mm}

From (\ref{ziko500}) and (\ref{ziko1000}) we get
\begin{equation}
\label{ziko1001}
H_n(x,1)=H_n(x).
\end{equation}

Obviously, without loss of generality, we can write
\begin{equation}
\label{ziko1002}
\left(j_1\ldots j_k\right)=
\bigl(\underbrace{j_1 \ldots j_1}_{m_1}\
\underbrace{j_2 \ldots j_2}_{m_2}\ \ldots\ 
\underbrace{j_r \ldots j_r}_{m_r}\bigr),
\end{equation}

\noindent
where $m_1+\ldots+m_r=k,$\ $m_1,\ldots, m_r=1,\ldots,k,$\ \
$r=1,\ldots,k,$\ \ $k>0,$ and $j_1,\ldots,j_r$ are pairwise different.

Analyzing the proof of Theorem 1.1 and using (\ref{new1010}), (\ref{new1200})
(see the proof of Theorem~1.22 below),
we can notice that w.~p.~1
(we suppose that the condition (\ref{ziko1002}) is fulfilled)

$$
J'[\phi_{j_1}\ldots \phi_{j_k}]_{T,t}^{(i_1\ldots i_1)}=
$$

\vspace{-2mm}
$$
=
\hbox{\vtop{\offinterlineskip\halign{
\hfil#\hfil\cr
{\rm l.i.m.}\cr
$\stackrel{}{{}_{N\to \infty}}$\cr
}} }
\sum\limits_{\stackrel{l_1,\ldots,l_k=0}{{}_{l_q\ne l_g;\ 
q\ne g;\ q, g=1,\ldots, k}}}^{N-1}
\phi_{j_1}(\tau_{l_1})\ldots
\phi_{j_k}(\tau_{l_k})
\Delta{\bf w}_{\tau_{l_1}}^{(i_1)}
\ldots
\Delta{\bf w}_{\tau_{l_k}}^{(i_1)}=
$$
$$
=
\hbox{\vtop{\offinterlineskip\halign{
\hfil#\hfil\cr
{\rm l.i.m.}\cr
$\stackrel{}{{}_{N\to \infty}}$\cr
}} }
\sum\limits_{\stackrel{l_1,\ldots,l_{m_1}=0}{{}_{l_q\ne l_g;\ 
q\ne g;\ q, g=1,\ldots, m_1}}}^{N-1}
\phi_{j_1}(\tau_{l_1})\ldots
\phi_{j_1}(\tau_{l_{m_1}})
\Delta{\bf w}_{\tau_{l_1}}^{(i_1)}
\ldots
\Delta{\bf w}_{\tau_{l_{m_1}}}^{(i_1)} \times
$$

\vspace{-1mm}
$$
\times
\hbox{\vtop{\offinterlineskip\halign{
\hfil#\hfil\cr
{\rm l.i.m.}\cr
$\stackrel{}{{}_{N\to \infty}}$\cr
}} }\hspace{-1mm}
\sum\limits_{\stackrel{l_{m_1+1},\ldots,l_{m_1+m_2}=0}{{}_{l_q\ne l_g;\ 
q\ne g;\ q, g= m_1+1,\ldots, m_1+m_2}}}^{N-1}
\phi_{j_2}(\tau_{l_{m_1+1}})\ldots
\phi_{j_2}(\tau_{l_{m_1+m_2}})
\Delta{\bf w}_{\tau_{l_{m_1+1}}}^{(i_1)}
\ldots
\Delta{\bf w}_{\tau_{l_{m_1+m_2}}}^{(i_1)} \times 
$$

\vspace{-6mm}
$$
\ldots
$$

\vspace{-6mm}
$$
\times
\hbox{\vtop{\offinterlineskip\halign{
\hfil#\hfil\cr
{\rm l.i.m.}\cr
$\stackrel{}{{}_{N\to \infty}}$\cr
}} }
\sum\limits_{\stackrel{l_{k-m_r+1},\ldots,l_k=0}{{}_{l_q\ne l_g;\ 
q\ne g;\ q, g=k-m_r+1,\ldots, k}}}^{N-1}
\phi_{j_r}(\tau_{l_{k-m_r+1}})\ldots
\phi_{j_r}(\tau_{l_k})
\Delta{\bf w}_{\tau_{l_{k-m_r+1}}}^{(i_1)}
\ldots
\Delta{\bf w}_{\tau_{l_k}}^{(i_1)}=
$$

\vspace{4mm}
$$
=
\hbox{\vtop{\offinterlineskip\halign{
\hfil#\hfil\cr
{\rm l.i.m.}\cr
$\stackrel{}{{}_{N\to \infty}}$\cr
}} }
\left(
\sum\limits_{l_1=0}^{N-1}
\phi_{j_1}(\tau_{l_1})\Delta{\bf w}_{\tau_{l_1}}^{(i_1)}
\ldots
\sum\limits_{l_{m_1}=0}^{N-1}
\phi_{j_1}(\tau_{l_{m_1}})
\Delta{\bf w}_{\tau_{l_{m_1}}}^{(i_1)}
-\right.
$$

\vspace{-1mm}
$$
\left.-
\sum_{(l_1,\ldots,l_{m_1})\in {\rm G}_{1,m_1}'}
\phi_{j_1}(\tau_{l_1})\Delta{\bf w}_{\tau_{l_1}}^{(i_1)}
\ldots
\phi_{j_1}(\tau_{l_{m_1}})
\Delta{\bf w}_{\tau_{l_{m_1}}}^{(i_1)}\right)\times
$$

\vspace{3mm}

$$
\times\hbox{\vtop{\offinterlineskip\halign{
\hfil#\hfil\cr
{\rm l.i.m.}\cr
$\stackrel{}{{}_{N\to \infty}}$\cr
}} }
\left(
\sum\limits_{l_{m_1+1}=0}^{N-1}
\phi_{j_2}(\tau_{l_{m_1}+1})\Delta{\bf w}_{\tau_{l_{m_1+1}}}^{(i_1)}
\ldots
\sum\limits_{l_{m_1+m_2}=0}^{N-1}
\phi_{j_2}(\tau_{l_{m_1+m_2}})
\Delta{\bf w}_{\tau_{l_{m_1+m_2}}}^{(i_1)}
-\right.
$$

\vspace{-1mm}
$$
\left.-
\sum_{(l_{m_1+1},\ldots,l_{m_1+m_2})\in {\rm G}_{m_1+1,m_1+m_2}'}
\phi_{j_2}(\tau_{l_{m_1}+1})\Delta{\bf w}_{\tau_{l_{m_1+1}}}^{(i_1)}
\ldots
\phi_{j_2}(\tau_{l_{m_1+m_2}})
\Delta{\bf w}_{\tau_{l_{m_1+m_2}}}^{(i_1)}
\right)\times
$$

\vspace{-4mm}
$$
\ldots 
$$
\vspace{-9mm}

$$
\times\hbox{\vtop{\offinterlineskip\halign{
\hfil#\hfil\cr
{\rm l.i.m.}\cr
$\stackrel{}{{}_{N\to \infty}}$\cr
}} }
\left(
\sum\limits_{l_{k-m_r+1}=0}^{N-1}
\phi_{j_r}(\tau_{l_{k-m_r+1}})
\Delta{\bf w}_{\tau_{l_{k-m_r+1}}}^{(i_1)}
\ldots
\sum\limits_{l_{k}=0}^{N-1}
\phi_{j_r}(\tau_{l_{k}})
\Delta{\bf w}_{\tau_{l_{k}}}^{(i_1)}
-\right.
$$

\vspace{-3mm}
$$
\left.-
\sum_{(l_{k-m_r+1},\ldots,l_{k})\in {\rm G}_{k-m_r+1,k}'}
\phi_{j_r}(\tau_{l_{k-m_r+1}})\Delta{\bf w}_{\tau_{k-m_r+1}}^{(i_1)}
\ldots
\phi_{j_r}(\tau_{l_{k}})
\Delta{\bf w}_{\tau_{l_{k}}}^{(i_1)}
\right),
$$

\vspace{3mm}
\noindent
where the set 
${\rm G}_{m,n}'$ is defined according to the same rule 
as the set ${\rm G}_{k}$ in (\ref{tyyy}). 
However, the elements of the set ${\rm G}_{m,n}'$ are the numbers 
$l_m,\ldots, l_n$ $(n>m)$, while the elements of the set ${\rm G}_{k}$ are the numbers 
$l_1,\ldots, l_k$.

We have (see the proof of Theorem 1.1) w.~p.~1 $(i_1\ne 0)$

$$
\hbox{\vtop{\offinterlineskip\halign{
\hfil#\hfil\cr
{\rm l.i.m.}\cr
$\stackrel{}{{}_{N\to \infty}}$\cr
}} }\left(\sum\limits_{l_1=0}^{N-1}
\phi_{j_1}(\tau_{l_1})\Delta{\bf w}_{\tau_{l_1}}^{(i_1)}
\ldots
\sum\limits_{l_{m_1}=0}^{N-1}
\phi_{j_1}(\tau_{l_{m_1}})
\Delta{\bf w}_{\tau_{l_{m_1}}}^{(i_1)}
-\right.
$$

\vspace{-1mm}
$$
\left.-
\sum_{(l_1,\ldots,l_{m_1})\in {\rm G}_{1,m_1}'}
\phi_{j_1}(\tau_{l_1})\Delta{\bf w}_{\tau_{l_1}}^{(i_1)}
\ldots
\phi_{j_1}(\tau_{l_{m_1}})
\Delta{\bf w}_{\tau_{l_{m_1}}}^{(i_1)}\right)=
$$

\vspace{4mm}
$$
=\hbox{\vtop{\offinterlineskip\halign{
\hfil#\hfil\cr
{\rm l.i.m.}\cr
$\stackrel{}{{}_{N\to \infty}}$\cr
}} }
\left(\left(\sum\limits_{l_1=0}^{N-1}\phi_{j_1}(\tau_{l_1})
\Delta {\bf w}_{\tau_{l_1}}^{(i_1)}\right)^{m_1}
+\sum\limits_{r=1}^{[m_1/2]}
(-1)^r \times\right.
$$

\vspace{-1mm}
$$
\times\sum_{\stackrel{(\{\{g_1, g_2\}, \ldots, 
\{g_{2r-1}, g_{2r}\}\}, \{q_1, \ldots, q_{m_1-2r}\})}
{{}_{\{g_1, g_2, \ldots, 
g_{2r-1}, g_{2r}, q_1, \ldots, q_{m_1-2r}\}=\{1, 2, \ldots, m_1\}}}}
\left(\sum\limits_{l_1=0}^{N-1}\phi_{j_1}^2(\tau_{l_1})
\left(\Delta {\bf w}_{\tau_{l_1}}^{(i_1)}\right)^2\right)^r\times
$$

\vspace{-1mm}
$$
\left.\times
\left(\sum\limits_{l_1=0}^{N-1}\phi_{j_1}(\tau_{l_1})
\Delta {\bf w}_{\tau_{l_1}}^{(i_1)}\right)^{m_1-2r}\right)= 
$$

\vspace{4mm}
$$
=\hbox{\vtop{\offinterlineskip\halign{
\hfil#\hfil\cr
{\rm l.i.m.}\cr
$\stackrel{}{{}_{N\to \infty}}$\cr
}} }\left(
\left(\sum\limits_{l_1=0}^{N-1}\phi_{j_1}(\tau_{l_1})
\Delta {\bf w}_{\tau_{l_1}}^{(i_1)}\right)^{m_1}
+\sum\limits_{r=1}^{[m_1/2]}
\frac{(-1)^r m_1!}{r!(m_1-2r)!2^r}\times\right.
$$

\vspace{-1mm}

$$
\left.\times
\left(\sum\limits_{l_1=0}^{N-1}\phi_{j_1}^2(\tau_{l_1})
\left(\Delta {\bf w}_{\tau_{l_1}}^{(i_1)}\right)^2\right)^r
\left(\sum\limits_{l_1=0}^{N-1}\phi_{j_1}(\tau_{l_1})
\Delta {\bf w}_{\tau_{l_1}}^{(i_1)}\right)^{m_1-2r}\right)=
$$

\vspace{4mm}
$$
=\hbox{\vtop{\offinterlineskip\halign{
\hfil#\hfil\cr
{\rm l.i.m.}\cr
$\stackrel{}{{}_{N\to \infty}}$\cr
}} }
\sum\limits_{r=0}^{[m_1/2]}
\frac{(-1)^r m_1!}{r!(m_1-2r)!2^r}
\left(\sum\limits_{l_1=0}^{N-1}\phi_{j_1}^2(\tau_{l_1})
\left(\Delta {\bf w}_{\tau_{l_1}}^{(i_1)}\right)^2\right)^r
\times
$$

$$
\times\left(\sum\limits_{l_1=0}^{N-1}\phi_{j_1}(\tau_{l_1})
\Delta {\bf w}_{\tau_{l_1}}^{(i_1)}\right)^{m_1-2r}=
$$
$$
=
\hbox{\vtop{\offinterlineskip\halign{
\hfil#\hfil\cr
{\rm l.i.m.}\cr
$\stackrel{}{{}_{N\to \infty}}$\cr
}} }
H_{m_1}\left(\sum\limits_{l_1=0}^{N-1}\phi_{j_1}(\tau_{l_1})
\Delta {\bf w}_{\tau_{l_1}}^{(i_1)},\
\sum\limits_{l_1=0}^{N-1}\phi_{j_1}^2(\tau_{l_1})
\left(\Delta {\bf w}_{\tau_{l_1}}^{(i_1)}\right)^2\right),
$$

\vspace{3mm}
\noindent
where notations are the same as in Theorems~1.1, 1.2.

Similarly we get w.~p.~1

\vspace{-2mm}
$$
\hbox{\vtop{\offinterlineskip\halign{
\hfil#\hfil\cr
{\rm l.i.m.}\cr
$\stackrel{}{{}_{N\to \infty}}$\cr
}} }
\left(\sum\limits_{l_{m_1+1}=0}^{N-1}
\phi_{j_2}(\tau_{l_{m_1}+1})\Delta{\bf w}_{\tau_{l_{m_1+1}}}^{(i_1)}
\ldots
\sum\limits_{l_{m_1+m_2}=0}^{N-1}
\phi_{j_2}(\tau_{l_{m_1+m_2}})
\Delta{\bf w}_{\tau_{l_{m_1+m_2}}}^{(i_1)}
-\right.
$$

\vspace{-4mm}
$$
-\left.
\sum_{(l_{m_1+1},\ldots,l_{m_1+m_2})\in {\rm G}_{m_1+1,m_1+m_2}'}
\phi_{j_2}(\tau_{l_{m_1}+1})\Delta{\bf w}_{\tau_{l_{m_1+1}}}^{(i_1)}
\ldots
\phi_{j_2}(\tau_{l_{m_1+m_2}})
\Delta{\bf w}_{\tau_{l_{m_1+m_2}}}^{(i_1)}\right)=
$$

\vspace{2mm}
$$
=\hbox{\vtop{\offinterlineskip\halign{
\hfil#\hfil\cr
{\rm l.i.m.}\cr
$\stackrel{}{{}_{N\to \infty}}$\cr
}} }
H_{m_2}\left(\sum\limits_{l_1=0}^{N-1}\phi_{j_2}(\tau_{l_1})
\Delta {\bf w}_{\tau_{l_1}}^{(i_1)},\
\sum\limits_{l_1=0}^{N-1}\phi_{j_2}^2(\tau_{l_1})
\left(\Delta {\bf w}_{\tau_{l_1}}^{(i_1)}\right)^2\right),
$$

$$
\ldots
$$

\vspace{-4mm}
$$
\hbox{\vtop{\offinterlineskip\halign{
\hfil#\hfil\cr
{\rm l.i.m.}\cr
$\stackrel{}{{}_{N\to \infty}}$\cr
}} }
\left(\sum\limits_{l_{k-m_r+1}=0}^{N-1}
\phi_{j_r}(\tau_{l_{k-m_r+1}})
\Delta{\bf w}_{\tau_{l_{k-m_r+1}}}^{(i_1)}
\ldots
\sum\limits_{l_{k}=0}^{N-1}
\phi_{j_r}(\tau_{l_{k}})
\Delta{\bf w}_{\tau_{l_{k}}}^{(i_1)}
- \right.
$$

\vspace{-4mm}
$$
\left.-
\sum_{(l_{k-m_r+1},\ldots,l_{k})\in {\rm G}_{k-m_r+1,k}'}
\phi_{j_r}(\tau_{l_{k-m_r+1}})\Delta{\bf w}_{\tau_{k-m_r+1}}^{(i_1)}
\ldots
\phi_{j_r}(\tau_{l_{k}})
\Delta{\bf w}_{\tau_{l_{k}}}^{(i_1)}\right)
=
$$

\vspace{2mm}
$$
=\hbox{\vtop{\offinterlineskip\halign{
\hfil#\hfil\cr
{\rm l.i.m.}\cr
$\stackrel{}{{}_{N\to \infty}}$\cr
}} }
H_{m_r}\left(\sum\limits_{l_1=0}^{N-1}\phi_{j_r}(\tau_{l_1})
\Delta {\bf w}_{\tau_{l_1}}^{(i_1)},\ 
\sum\limits_{l_1=0}^{N-1}\phi_{j_r}^2(\tau_{l_1})
\left(\Delta {\bf w}_{\tau_{l_1}}^{(i_1)}\right)^2\right).
$$

\vspace{5mm}

Then

\vspace{-7mm}
$$
J'[\phi_{j_1}\ldots \phi_{j_k}]_{T,t}^{(i_1\ldots i_1)}=
$$

$$
=\hbox{\vtop{\offinterlineskip\halign{
\hfil#\hfil\cr
{\rm l.i.m.}\cr
$\stackrel{}{{}_{N\to \infty}}$\cr
}} }
H_{m_1}\left(\sum\limits_{l_1=0}^{N-1}\phi_{j_1}(\tau_{l_1})
\Delta {\bf w}_{\tau_{l_1}}^{(i_1)},\
\sum\limits_{l_1=0}^{N-1}\phi_{j_1}^2(\tau_{l_1})
\left(\Delta {\bf w}_{\tau_{l_1}}^{(i_1)}\right)^2\right)\times 
$$
$$
\times \hbox{\vtop{\offinterlineskip\halign{
\hfil#\hfil\cr
{\rm l.i.m.}\cr
$\stackrel{}{{}_{N\to \infty}}$\cr
}} }
H_{m_2}\left(\sum\limits_{l_1=0}^{N-1}\phi_{j_2}(\tau_{l_1})
\Delta {\bf w}_{\tau_{l_1}}^{(i_1)},\
\sum\limits_{l_1=0}^{N-1}\phi_{j_2}^2(\tau_{l_1})
\left(\Delta {\bf w}_{\tau_{l_1}}^{(i_1)}\right)^2\right)\times \ldots
$$

\vspace{-2mm}
\begin{equation}
\label{ziko2000}
\ldots \times \hbox{\vtop{\offinterlineskip\halign{
\hfil#\hfil\cr
{\rm l.i.m.}\cr
$\stackrel{}{{}_{N\to \infty}}$\cr
}} }
H_{m_r}\left(\sum\limits_{l_1=0}^{N-1}\phi_{j_r}(\tau_{l_1})
\Delta {\bf w}_{\tau_{l_1}}^{(i_1)},\ 
\sum\limits_{l_1=0}^{N-1}\phi_{j_r}^2(\tau_{l_1})
\left(\Delta {\bf w}_{\tau_{l_1}}^{(i_1)}\right)^2\right)
\end{equation}

\vspace{3mm}
\noindent
w.~p.~1 for $i_1\ne 0$ and

\vspace{-2mm}
$$
J'[\phi_{j_1}\ldots \phi_{j_k}]_{T,t}^{(0\ldots 0)}=
\lim\limits_{N\to\infty}
\left(\sum\limits_{l_1=0}^{N-1}\phi_{j_1}(\tau_{l_1})
\Delta \tau_{l_1}\right)^{m_1}
\ldots 
\left(\sum\limits_{l_r=0}^{N-1}\phi_{j_r}(\tau_{l_r})
\Delta \tau_{l_r}\right)^{m_r}=
$$

\vspace{-2mm}
\begin{equation}
\label{ziko2001}
~~~~=\left(\int\limits_t^T\phi_{j_1}(s)
ds\right)^{m_1}
\ldots 
\left(\int\limits_t^T\phi_{j_r}(s)
ds\right)^{m_r}=
\left(\zeta_{j_1}^{(0)}\right)^{m_1}
\ldots 
\left(\zeta_{j_r}^{(0)}\right)^{m_r}
\end{equation}

\vspace{3mm}
\noindent
for $i_1=0,$ where we suppose that the condition 
(\ref{ziko1002}) is fulfilled; also we use in (\ref{ziko2000}) and (\ref{ziko2001}) 
the same notations as in the proof of Theorem 1.1.

Applying (\ref{ziko1000}), (\ref{ziko1001}), Lemma~1.3, and Remark~1.2 to the right-hand side of
(\ref{ziko2000}), we finally obtain w.~p.~1

\vspace{-2mm}
$$
J'[\phi_{j_1}\ldots \phi_{j_k}]_{T,t}^{(i_1\ldots i_1)}=
H_{m_1}\left(\int\limits_t^T\phi_{j_1}(s)
d{\bf w}_s^{(i_1)},
\int\limits_t^T\phi_{j_1}^2(s)ds
\right)\times
$$

\vspace{-2mm}
$$
\times
H_{m_2}\left(\int\limits_t^T\phi_{j_2}(s)
d{\bf w}_s^{(i_1)},
\int\limits_t^T\phi_{j_2}^2(s)ds
\right)\ldots
H_{m_r}\left(\int\limits_t^T\phi_{j_r}(s)
d{\bf w}_s^{(i_1)},
\int\limits_t^T\phi_{j_r}^2(s)ds\right)=
$$

\vspace{2mm}
$$
=
H_{m_1}\left(\zeta_{j_1}^{(i_1)},
1\right)H_{m_2}\left(\zeta_{j_2}^{(i_1)},
1\right)
\ldots 
H_{m_r}\left(\zeta_{j_r}^{(i_1)},1\right)=
$$

\vspace{2mm}
$$
=
H_{m_1}\left(\zeta_{j_1}^{(i_1)}\right)
H_{m_2}\left(\zeta_{j_2}^{(i_1)}\right)
\ldots 
H_{m_r}\left(\zeta_{j_r}^{(i_1)}\right)
$$

\vspace{3mm}
\noindent
for $i_1\ne 0,$ 
where we suppose that the condition 
(\ref{ziko1002}) is fulfilled.
Thus, an equality similar to (\ref{ziko20}) is proved without using 
Theorem 3.1 \cite{ito1951}.

Consider particular cases of the equality (\ref{ziko100}) for
$k=1,\ldots,4$ and $i_1,\ldots,i_4=1,\ldots,m$ (see (\ref{a1})--(\ref{a4})).
We have w.~p.~1

\vspace{-2mm}
$$
J'[\phi_{j_1}]_{T,t}^{(i_1)}=
\zeta_{j_1}^{(i_1)}=H_1\left(\zeta_{j_1}^{(i_1)}\right);
$$

\vspace{2mm}
$$
J'[\phi_{j_1}\phi_{j_2}]_{T,t}^{(i_1i_2)}=\zeta_{j_1}^{(i_1)}\zeta_{j_2}^{(i_2)}
-{\bf 1}_{\{i_1=i_2\}}
{\bf 1}_{\{j_1=j_2\}}=
$$

\begin{equation}
\label{ziko301}
=\left\{
\begin{matrix}
H_2\left(\zeta_{j_1}^{(i_1)}\right)H_0\left(\zeta_{j_2}^{(i_2)}\right),\ 
&\hbox{\rm if}\ \ \ 
i_1=i_2,\ j_1=j_2\cr\cr
H_1\left(\zeta_{j_1}^{(i_1)}\right)H_1\left(\zeta_{j_2}^{(i_2)}\right),\  &\hbox{\rm otherwise}
\end{matrix}\right.;
\end{equation}

\vspace{4.5mm}
$$
J'[\phi_{j_1}\phi_{j_2}\phi_{j_3}]_{T,t}^{(i_1i_1i_1)}=
\zeta_{j_1}^{(i_1)}\zeta_{j_2}^{(i_1)}\zeta_{j_3}^{(i_1)}
-
{\bf 1}_{\{j_1=j_2\}}
\zeta_{j_3}^{(i_1)}
-
{\bf 1}_{\{j_2=j_3\}}
\zeta_{j_1}^{(i_1)}-
{\bf 1}_{\{j_1=j_3\}}
\zeta_{j_2}^{(i_1)}=
$$

\vspace{-2mm}
\begin{equation}
\label{ziko302}
=\left\{
\begin{matrix}
H_3\left(\zeta_{j_1}^{(i_1)}\right)H_0\left(\zeta_{j_2}^{(i_1)}\right)
H_0\left(\zeta_{j_3}^{(i_1)}\right),\ 
&\hbox{\rm if}\ \ \ 
j_1=j_2=j_3
\cr\cr
H_2\left(\zeta_{j_1}^{(i_1)}\right)H_0\left(\zeta_{j_2}^{(i_1)}\right)
H_1\left(\zeta_{j_3}^{(i_1)}\right),\ 
&\hbox{\rm if}\ \ \ 
j_1=j_2\ne j_3
\cr\cr
H_1\left(\zeta_{j_1}^{(i_1)}\right)H_2\left(\zeta_{j_2}^{(i_1)}\right)
H_0\left(\zeta_{j_3}^{(i_1)}\right),\ 
&\hbox{\rm if}\ \ \ 
j_2=j_3\ne j_1
\cr\cr
H_0\left(\zeta_{j_1}^{(i_1)}\right)H_1\left(\zeta_{j_2}^{(i_1)}\right)
H_2\left(\zeta_{j_3}^{(i_1)}\right),\ 
&\hbox{\rm if}\ \ \ 
j_1=j_3\ne j_2
\cr\cr
H_1\left(\zeta_{j_1}^{(i_1)}\right)H_1\left(\zeta_{j_2}^{(i_1)}\right)
H_1\left(\zeta_{j_3}^{(i_1)}\right),\ 
&\hbox{\rm if}\ \ \ 
j_1\ne j_2,\ 
j_2\ne j_3,\ j_1\ne j_3
\end{matrix}\right.;
\end{equation}

\vspace{4mm}
$$
J'[\phi_{j_1}\phi_{j_2}\phi_{j_3}]_{T,t}^{(i_1i_2i_3)}=
\zeta_{j_1}^{(i_1)}\zeta_{j_2}^{(i_2)}\zeta_{j_3}^{(i_2)}=
H_1\left(\zeta_{j_1}^{(i_1)}\right)H_1\left(\zeta_{j_2}^{(i_2)}\right)
H_1\left(\zeta_{j_3}^{(i_3)}\right),
$$

\vspace{2mm}
\noindent
where $i_1,i_2,i_3$ are pairwise different;

\vspace{1.5mm}
$$
J'[\phi_{j_1}\phi_{j_2}\phi_{j_3}]_{T,t}^{(i_1i_1i_3)}=
\zeta_{j_1}^{(i_1)}\zeta_{j_2}^{(i_1)}\zeta_{j_3}^{(i_3)}
-
{\bf 1}_{\{j_1=j_2\}}
\zeta_{j_3}^{(i_3)}=
$$

\vspace{-1mm}
$$
=\left(\zeta_{j_1}^{(i_1)}\zeta_{j_2}^{(i_1)}
-
{\bf 1}_{\{j_1=j_2\}}\right)
\zeta_{j_3}^{(i_3)}=J'[\phi_{j_1}\phi_{j_2}]_{T,t}^{(i_1i_1)}J'[\phi_{j_3}]_{T,t}^{(i_3)}=
$$

\newpage
\noindent
$$
=\left\{
\begin{matrix}
H_2\left(\zeta_{j_1}^{(i_1)}\right)H_0\left(\zeta_{j_2}^{(i_1)}\right)
H_1\left(\zeta_{j_3}^{(i_3)}\right),\ 
&\hbox{\rm if}\ \ \ 
j_1=j_2\cr\cr
H_1\left(\zeta_{j_1}^{(i_1)}\right)H_1\left(\zeta_{j_2}^{(i_1)}\right)
H_1\left(\zeta_{j_3}^{(i_3)}\right),\  &\hbox{\rm if}\ \ \ 
j_1\ne j_2
\end{matrix}\right.,
$$

\vspace{1.5mm}
\noindent
where $i_1=i_2\ne i_3$;

\vspace{1.5mm}
$$
J'[\phi_{j_1}\phi_{j_2}\phi_{j_3}]_{T,t}^{(i_1i_2i_2)}=
\zeta_{j_1}^{(i_1)}\zeta_{j_2}^{(i_2)}\zeta_{j_3}^{(i_2)}
-
{\bf 1}_{\{j_2=j_3\}}
\zeta_{j_1}^{(i_1)}=
$$

\vspace{-2mm}
$$
=\zeta_{j_1}^{(i_1)}\left(\zeta_{j_2}^{(i_2)}\zeta_{j_3}^{(i_2)}
-
{\bf 1}_{\{j_2=j_3\}}\right)
=J'[\phi_{j_1}]_{T,t}^{(i_1)}J'[\phi_{j_2}\phi_{j_3}]_{T,t}^{(i_2i_2)}=
$$

\vspace{-1mm}
$$
=\left\{
\begin{matrix}
H_1\left(\zeta_{j_1}^{(i_1)}\right)H_2\left(\zeta_{j_2}^{(i_2)}\right)
H_0\left(\zeta_{j_3}^{(i_2)}\right),\ 
&\hbox{\rm if}\ \ \ 
j_2=j_3\cr\cr
H_1\left(\zeta_{j_1}^{(i_1)}\right)H_1\left(\zeta_{j_2}^{(i_2)}\right)
H_1\left(\zeta_{j_3}^{(i_2)}\right),\  &\hbox{\rm if}\ \ \ 
j_1\ne j_2
\end{matrix}\right.,
$$

\vspace{1.5mm}
\noindent
where $i_1\ne i_2=i_3$;

\vspace{1.5mm}
$$
J'[\phi_{j_1}\phi_{j_2}\phi_{j_3}]_{T,t}^{(i_1i_2i_1)}=
\zeta_{j_1}^{(i_1)}\zeta_{j_2}^{(i_2)}\zeta_{j_3}^{(i_1)}
-
{\bf 1}_{\{j_1=j_3\}}
\zeta_{j_2}^{(i_2)}=
$$

\vspace{-2mm}
$$
=\zeta_{j_2}^{(i_2)}\left(\zeta_{j_1}^{(i_1)}\zeta_{j_3}^{(i_1)}
-
{\bf 1}_{\{j_1=j_3\}}\right)
=J'[\phi_{j_2}]_{T,t}^{(i_2)}J'[\phi_{j_1}\phi_{j_3}]_{T,t}^{(i_1i_1)}=
$$

\vspace{-1mm}
$$
=\left\{
\begin{matrix}
H_2\left(\zeta_{j_1}^{(i_1)}\right)H_1\left(\zeta_{j_2}^{(i_2)}\right)
H_0\left(\zeta_{j_3}^{(i_1)}\right),\ 
&\hbox{\rm if}\ \ \ 
j_1=j_3\cr\cr
H_1\left(\zeta_{j_1}^{(i_1)}\right)H_1\left(\zeta_{j_2}^{(i_2)}\right)
H_1\left(\zeta_{j_3}^{(i_1)}\right),\  &\hbox{\rm if}\ \ \ 
j_1\ne j_3
\end{matrix}\right.,
$$

\vspace{1.5mm}
\noindent
where $i_1=i_3\ne i_2$;

$$
J'[\phi_{j_1}\phi_{j_2}\phi_{j_3}\phi_{j_4}]_{T,t}^{(i_1i_1i_1i_1)}=
\zeta_{j_1}^{(i_1)}\zeta_{j_2}^{(i_1)}\zeta_{j_3}^{(i_1)}\zeta_{j_4}^{(i_1)}-
$$

\vspace{-5mm}
$$
-
{\bf 1}_{\{j_1=j_2\}}
\zeta_{j_3}^{(i_1)}
\zeta_{j_4}^{(i_1)}
-
{\bf 1}_{\{j_1=j_3\}}
\zeta_{j_2}^{(i_1)}
\zeta_{j_4}^{(i_1)}-
{\bf 1}_{\{j_1=j_4\}}
\zeta_{j_2}^{(i_1)}
\zeta_{j_3}^{(i_1)}
-
$$

\vspace{-4mm}
$$
-
{\bf 1}_{\{j_2=j_3\}}
\zeta_{j_1}^{(i_1)}
\zeta_{j_4}^{(i_1)}-
{\bf 1}_{\{j_2=j_4\}}
\zeta_{j_1}^{(i_1)}
\zeta_{j_3}^{(i_1)}
-
{\bf 1}_{\{j_3=j_4\}}
\zeta_{j_1}^{(i_1)}
\zeta_{j_2}^{(i_1)}+
$$

\vspace{-5mm}
$$
+
{\bf 1}_{\{j_1=j_2\}}
{\bf 1}_{\{j_3=j_4\}}
+
{\bf 1}_{\{j_1=j_3\}}
{\bf 1}_{\{j_2=j_4\}}+
{\bf 1}_{\{j_1=j_4\}}
{\bf 1}_{\{j_2=j_3\}}=
$$
$$
=\left\{
\begin{matrix}
H_4\left(\zeta_{j_1}^{(i_1)}\right)H_0\left(\zeta_{j_2}^{(i_1)}\right)
H_0\left(\zeta_{j_3}^{(i_1)}\right)H_0\left(\zeta_{j_4}^{(i_1)}\right),\ 
&\hbox{\rm if}\ \ \ {\rm (I)}
\cr\cr
H_1\left(\zeta_{j_1}^{(i_1)}\right)H_1\left(\zeta_{j_2}^{(i_1)}\right)
H_1\left(\zeta_{j_3}^{(i_1)}\right)H_1\left(\zeta_{j_4}^{(i_1)}\right),\ 
&\hbox{\rm if}\ \ \ {\rm (II)}
\cr\cr
H_2\left(\zeta_{j_1}^{(i_1)}\right)H_0\left(\zeta_{j_2}^{(i_1)}\right)
H_1\left(\zeta_{j_3}^{(i_1)}\right)H_1\left(\zeta_{j_4}^{(i_1)}\right),\ 
&\hbox{\rm if}\ \ \ {\rm (III)}
\cr\cr
H_0\left(\zeta_{j_1}^{(i_1)}\right)H_1\left(\zeta_{j_2}^{(i_1)}\right)
H_2\left(\zeta_{j_3}^{(i_1)}\right)H_1\left(\zeta_{j_4}^{(i_1)}\right),\ 
&\hbox{\rm if}\ \ \ {\rm (IV)}
\cr\cr
H_0\left(\zeta_{j_1}^{(i_1)}\right)H_1\left(\zeta_{j_2}^{(i_1)}\right)
H_1\left(\zeta_{j_3}^{(i_1)}\right)H_2\left(\zeta_{j_4}^{(i_1)}\right),\ 
&\hbox{\rm if}\ \ \ {\rm (V)}
\cr\cr
H_1\left(\zeta_{j_1}^{(i_1)}\right)H_0\left(\zeta_{j_2}^{(i_1)}\right)
H_2\left(\zeta_{j_3}^{(i_1)}\right)H_1\left(\zeta_{j_4}^{(i_1)}\right),\ 
&\hbox{\rm if}\ \ \ {\rm (VI)}
\cr\cr
H_1\left(\zeta_{j_1}^{(i_1)}\right)H_0\left(\zeta_{j_2}^{(i_1)}\right)
H_1\left(\zeta_{j_3}^{(i_1)}\right)H_2\left(\zeta_{j_4}^{(i_1)}\right),\ 
&\hbox{\rm if}\ \ \ {\rm (VII)}
\cr\cr
H_1\left(\zeta_{j_1}^{(i_1)}\right)H_1\left(\zeta_{j_2}^{(i_1)}\right)
H_0\left(\zeta_{j_3}^{(i_1)}\right)H_2\left(\zeta_{j_4}^{(i_1)}\right),\ 
&\hbox{\rm if}\ \ \ {\rm (VIII)}
\cr\cr
H_3\left(\zeta_{j_1}^{(i_1)}\right)H_0\left(\zeta_{j_2}^{(i_1)}\right)
H_0\left(\zeta_{j_3}^{(i_1)}\right)H_1\left(\zeta_{j_4}^{(i_1)}\right),\ 
&\hbox{\rm if}\ \ \ {\rm (IX)}
\cr\cr
H_1\left(\zeta_{j_1}^{(i_1)}\right)H_3\left(\zeta_{j_2}^{(i_1)}\right)
H_0\left(\zeta_{j_3}^{(i_1)}\right)H_0\left(\zeta_{j_4}^{(i_1)}\right),\ 
&\hbox{\rm if}\ \ \ {\rm (X)}
\cr\cr
H_0\left(\zeta_{j_1}^{(i_1)}\right)H_0\left(\zeta_{j_2}^{(i_1)}\right)
H_1\left(\zeta_{j_3}^{(i_1)}\right)H_3\left(\zeta_{j_4}^{(i_1)}\right),\ 
&\hbox{\rm if}\ \ \ {\rm (XI)}
\cr\cr
H_0\left(\zeta_{j_1}^{(i_1)}\right)H_1\left(\zeta_{j_2}^{(i_1)}\right)
H_0\left(\zeta_{j_3}^{(i_1)}\right)H_3\left(\zeta_{j_4}^{(i_1)}\right),\ 
&\hbox{\rm if}\ \ \ {\rm (XII)}
\cr\cr
H_2\left(\zeta_{j_1}^{(i_1)}\right)H_0\left(\zeta_{j_2}^{(i_1)}\right)
H_0\left(\zeta_{j_3}^{(i_1)}\right)H_2\left(\zeta_{j_4}^{(i_1)}\right),\ 
&\hbox{\rm if}\ \ \ {\rm (XIII)}
\cr\cr
H_2\left(\zeta_{j_1}^{(i_1)}\right)H_2\left(\zeta_{j_2}^{(i_1)}\right)
H_0\left(\zeta_{j_3}^{(i_1)}\right)H_0\left(\zeta_{j_4}^{(i_1)}\right),\ 
&\hbox{\rm if}\ \ \ {\rm (XIV)}
\cr\cr
H_2\left(\zeta_{j_1}^{(i_1)}\right)H_0\left(\zeta_{j_2}^{(i_1)}\right)
H_2\left(\zeta_{j_3}^{(i_1)}\right)H_0\left(\zeta_{j_4}^{(i_1)}\right),\ 
&\hbox{\rm if}\ \ \ {\rm (XV)}
\end{matrix}\right.,
$$
where $H_n(x)$ is the Hermite polynomial (\ref{ziko500}) of degree $n$  
and (I)--(XV) are the following conditions

(I).\ $j_1=j_2=j_3=j_4,$

(II).\ $j_1, j_2, j_3, j_4$ are pairwise different,

(III).\ $j_1=j_2\ne j_3, j_4;\ j_3\ne j_4,$

(IV).\ $j_1=j_3\ne j_2, j_4;\ j_2\ne j_4,$

(V).\ $j_1=j_4\ne j_2, j_3;\ j_2\ne j_3,$

(VI).\ $j_2=j_3\ne j_1, j_4;\ j_1\ne j_4,$

(VII).\ $j_2=j_4\ne j_1, j_3;\ j_1\ne j_3,$

(VIII).\ $j_3=j_4\ne j_1, j_2;\ j_1\ne j_2,$

(IX).\ $j_1=j_2=j_3\ne j_4,$

(X).\ $j_2=j_3=j_4\ne j_1,$

(XI).\ $j_1=j_2=j_4\ne j_3,$

(XII).\ $j_1=j_3=j_4\ne j_2,$

(XIII).\ $j_1=j_2\ne j_3=j_4,$

(XIV).\ $j_1=j_3\ne j_2=j_4,$

(XV).\ $j_1=j_4\ne j_2=j_3.$

\vspace{2mm}

\noindent
Moreover, from (\ref{ziko30}) we have w.~p.~1

\vspace{-4mm}
$$
J'[\phi_{j_1}\phi_{j_2}\phi_{j_3}\phi_{j_4}]_{T,t}^{(i_1i_2i_3i_4)}=
H_1\left(\zeta_{j_1}^{(i_1)}\right)H_1\left(\zeta_{j_2}^{(i_2)}\right)
H_1\left(\zeta_{j_3}^{(i_3)}\right)H_1\left(\zeta_{j_4}^{(i_4)}\right),
$$

\vspace{1mm}
\noindent
where $i_1, i_2, i_3, i_4$ are pairwise different;

\vspace{-3mm}
\begin{equation}
\label{ziko201}
~~~~~~ J'[\phi_{j_1}\phi_{j_2}\phi_{j_3}\phi_{j_4}]_{T,t}^{(i_1i_1i_3i_4)}=
J'[\phi_{j_1}\phi_{j_2}]_{T,t}^{(i_1i_1)}
H_1\left(\zeta_{j_3}^{(i_3)}\right)H_1\left(\zeta_{j_4}^{(i_4)}\right),
\end{equation}

\vspace{1mm}
\noindent
where $i_1=i_2\ne i_3, i_4;\ i_3\ne i_4;$

\vspace{-3mm}
\begin{equation}
\label{ziko202}
~~~~~~J'[\phi_{j_1}\phi_{j_2}\phi_{j_3}\phi_{j_4}]_{T,t}^{(i_1i_2i_1i_4)}=
J'[\phi_{j_1}\phi_{j_3}]_{T,t}^{(i_1i_1)}
H_1\left(\zeta_{j_2}^{(i_2)}\right)H_1\left(\zeta_{j_4}^{(i_4)}\right),
\end{equation}

\vspace{1mm}
\noindent
where $i_1=i_3\ne i_2, i_4;\ i_2\ne i_4;$

\vspace{-3mm}
\begin{equation}
\label{ziko203}
~~~~~~J'[\phi_{j_1}\phi_{j_2}\phi_{j_3}\phi_{j_4}]_{T,t}^{(i_1i_2i_3i_1)}=
J'[\phi_{j_1}\phi_{j_4}]_{T,t}^{(i_1i_1)}
H_1\left(\zeta_{j_2}^{(i_2)}\right)H_1\left(\zeta_{j_3}^{(i_3)}\right),
\end{equation}

\vspace{1mm}
\noindent
where $i_1=i_4\ne i_2, i_3;\ i_2\ne i_3;$

\vspace{-3mm}
\begin{equation}
\label{ziko204}
~~~~~~J'[\phi_{j_1}\phi_{j_2}\phi_{j_3}\phi_{j_4}]_{T,t}^{(i_1i_2i_2i_4)}=
J'[\phi_{j_2}\phi_{j_3}]_{T,t}^{(i_2i_2)}
H_1\left(\zeta_{j_1}^{(i_1)}\right)H_1\left(\zeta_{j_4}^{(i_4)}\right),
\end{equation}

\vspace{1mm}
\noindent
where $i_2=i_3\ne i_1, i_4;\ i_1\ne i_4;$

\vspace{-3mm}
\begin{equation}
\label{ziko205}
~~~~~~J'[\phi_{j_1}\phi_{j_2}\phi_{j_3}\phi_{j_4}]_{T,t}^{(i_1i_2i_3i_2)}=
J'[\phi_{j_2}\phi_{j_4}]_{T,t}^{(i_2i_2)}
H_1\left(\zeta_{j_1}^{(i_1)}\right)H_1\left(\zeta_{j_3}^{(i_3)}\right),
\end{equation}

\vspace{1mm}
\noindent
where $i_2=i_4\ne i_1, i_3;\ i_1\ne i_3;$

\vspace{-3mm}
\begin{equation}
\label{ziko206}
~~~~~~J'[\phi_{j_1}\phi_{j_2}\phi_{j_3}\phi_{j_4}]_{T,t}^{(i_1i_2i_3i_3)}=
J'[\phi_{j_3}\phi_{j_4}]_{T,t}^{(i_3i_3)}
H_1\left(\zeta_{j_1}^{(i_1)}\right)H_1\left(\zeta_{j_2}^{(i_2)}\right),
\end{equation}

\vspace{1mm}
\noindent
where $i_3=i_4\ne i_1, i_2;\ i_1\ne i_2;$

\vspace{-3mm}
\begin{equation}
\label{ziko207}
~~~~~~J'[\phi_{j_1}\phi_{j_2}\phi_{j_3}\phi_{j_4}]_{T,t}^{(i_1i_1i_1i_4)}=
J'[\phi_{j_1}\phi_{j_2}\phi_{j_3}]_{T,t}^{(i_1i_1i_1)}
H_1\left(\zeta_{j_4}^{(i_4)}\right),
\end{equation}

\vspace{1mm}
\noindent
where $i_1=i_2=i_3\ne i_4;$

\vspace{-3mm}
\begin{equation}
\label{ziko208}
~~~~~~J'[\phi_{j_1}\phi_{j_2}\phi_{j_3}\phi_{j_4}]_{T,t}^{(i_1i_2i_2i_2)}=
J'[\phi_{j_2}\phi_{j_3}\phi_{j_4}]_{T,t}^{(i_2i_2i_2)}
H_1\left(\zeta_{j_1}^{(i_1)}\right),
\end{equation}

\vspace{1mm}
\noindent
where $i_2=i_3=i_4\ne i_1;$

\vspace{-3mm}
\begin{equation}
\label{ziko209}
~~~~~~J'[\phi_{j_1}\phi_{j_2}\phi_{j_3}\phi_{j_4}]_{T,t}^{(i_1i_1i_3i_1)}=
J'[\phi_{j_1}\phi_{j_2}\phi_{j_4}]_{T,t}^{(i_1i_1i_1)}
H_1\left(\zeta_{j_3}^{(i_3)}\right),
\end{equation}

\vspace{1mm}
\noindent
where $i_1=i_2=i_4\ne i_3;$

\vspace{-3mm}
\begin{equation}
\label{ziko210}
~~~~~~J'[\phi_{j_1}\phi_{j_2}\phi_{j_3}\phi_{j_4}]_{T,t}^{(i_1i_2i_1i_1)}=
J'[\phi_{j_1}\phi_{j_3}\phi_{j_4}]_{T,t}^{(i_1i_1i_1)}
H_1\left(\zeta_{j_2}^{(i_2)}\right),
\end{equation}

\vspace{1mm}
\noindent
where $i_1=i_3=i_4\ne i_2;$

\vspace{-3mm}
\begin{equation}
\label{ziko211}
~~~~~~J'[\phi_{j_1}\phi_{j_2}\phi_{j_3}\phi_{j_4}]_{T,t}^{(i_1i_1i_3i_3)}=
J'[\phi_{j_1}\phi_{j_2}]_{T,t}^{(i_1i_1)}J'[\phi_{j_3}\phi_{j_4}]_{T,t}^{(i_3i_3)},
\end{equation}

\vspace{1mm}
\noindent
where $i_1=i_2\ne i_3=i_4;$

\newpage
\noindent
\begin{equation}
\label{ziko212}
~~~~~~J'[\phi_{j_1}\phi_{j_2}\phi_{j_3}\phi_{j_4}]_{T,t}^{(i_1i_2i_1i_2)}=
J'[\phi_{j_1}\phi_{j_3}]_{T,t}^{(i_1i_1)}J'[\phi_{j_2}\phi_{j_4}]_{T,t}^{(i_2i_2)},
\end{equation}

\vspace{1mm}
\noindent
where $i_1=i_3\ne i_2=i_4;$

\vspace{-3mm}
\begin{equation}
\label{ziko213}
~~~~~~J'[\phi_{j_1}\phi_{j_2}\phi_{j_3}\phi_{j_4}]_{T,t}^{(i_1i_2i_2i_1)}=
J'[\phi_{j_1}\phi_{j_4}]_{T,t}^{(i_1i_1)}J'[\phi_{j_2}\phi_{j_3}]_{T,t}^{(i_2i_2)},
\end{equation}

\vspace{1mm}
\noindent
where $i_1=i_4\ne i_2=i_3.$

Note that the right-hand sides of (\ref{ziko201})--(\ref{ziko213})
contain multiple Wiener stochastic integrals of multiplicities 2 and 3.
These integrals are considered in detail in (\ref{ziko301}), (\ref{ziko302}).

It should be noted that the formulas (\ref{leto6000}) (Theorem~1.2) and 
(\ref{ziko800}) (Theorem~1.14) are interesting from various points of view.
The formulas (\ref{a1})--(\ref{a6}) (these formulas are particular cases of 
(\ref{leto6000}) for $k=1,\ldots,6$)
are convenient for numerical modeling of iterated It\^{o} stochastic integrals 
of multiplicities 1 to 6 (see Chapter 5).
For example, in \cite{Kuz-Kuz} and \cite{Mikh-1}, approximations 
of iterated It\^{o} stochastic integrals of multiplicities 1 to 6 
in the Python programming 
language were successfully implemented using (\ref{a1})--(\ref{a6}) and 
Legendre polynomials.

On the other hand, the equality (\ref{ziko800}) is interesting 
by a number of reasons.
Firstly, this equality connects It\^{o}'s results on multiple Wiener stochastic 
integral (\cite{ito1951}, Theorem~3.1) with the theory of 
mean-square approximation of iterated It\^{o} stochastic integrals 
presented in this book.
Secondly, the equality (\ref{ziko800}) is based on the 
Hermite polynomials, which have the 
orthogonality property on $\bf {R}$ with a Gaussian weight. 
This feature opens up new possibilities in the study 
of iterated It\^{o} stochastic integrals.
Note that the indicated orthogonality property 
is indirectly reflected by the formula (\ref{tyty1}) (see the proof of Theorem 1.3).

\section{Generalization of Theorems 1.1, 1.2, 1.14, and 1.15 to the Case of an Arbitrary 
Complete Ortho\-nor\-mal System of Functions in the Space $L_2([t, T])$
and $\psi_1(\tau),$ $\ldots,\psi_k(\tau)\in L_2([t, T]),$
$\Phi(t_1,\ldots,t_k)\in L_2([t, T]^k)$}

In this section, we will use the definition of the multiple Wiener 
stochastic integral from \cite{ito1951}, \cite{Kuo} to generalize Theorems 
1.1, 1.2, 1.14, and 1.15 to the case of an arbitrary 
complete orthonormal system of functions in the space $L_2([t, T])$
and $\psi_1(\tau),$ $\ldots,\psi_k(\tau)\in L_2([t, T]),$
$\Phi(t_1,\ldots,t_k)\in L_2([t, T]^k)$.

Consider the following step function on the hypercube $[t, T]^k$
\begin{equation}
\label{chain3}
~~~~~~~~~~~\Phi_N(t_1,\ldots,t_k)=\sum\limits_{l_1,\ldots,l_k=0}^{N-1}
a_{l_1 \ldots l_k} {\bf 1}_{[\tau_{l_1},\tau_{l_1+1})}(t_1) \ldots
{\bf 1}_{[\tau_{l_k},\tau_{l_k+1})}(t_k),
\end{equation}

\noindent
where $a_{l_1 \ldots l_k}\in{\bf R}$ and such that 
$a_{l_1 \ldots l_k}=0$ if $l_p=l_q$ for some $p\ne q,$
$$
{\bf 1}_A (\tau)=\left\{
\begin{matrix}
1\ &{\rm if}\ \tau\in A \cr\cr
0\ &\hbox{\rm otherwise}
\end{matrix}\right.,
$$

\noindent
$N\in{\bf N},$ $\left\{\tau_{j}\right\}_{j=0}^{N}$ is a partition of
$[t,T],$ which satisfies the condition (\ref{1111}):
\begin{equation}
\label{1111xxx1}
t=\tau_0<\ldots <\tau_N=T,\ \ \
\Delta_N=
\hbox{\vtop{\offinterlineskip\halign{
\hfil#\hfil\cr
{\rm max}\cr
$\stackrel{}{{}_{0\le j\le N-1}}$\cr
}} }\Delta\tau_j\to 0\ \ \hbox{if}\ \ N\to \infty,\ \ \ 
\Delta\tau_j=\tau_{j+1}-\tau_j.
\end{equation}

Let us define the multiple Wiener stochastic integral for $\Phi_N(t_1,\ldots,t_k)$ 
\cite{ito1951}, \cite{Kuo}
\begin{equation}
\label{chain9}
J'[\Phi_N]_{T,t}^{(i_1\ldots i_k)}\stackrel{\sf def}{=}
\sum\limits_{l_1,\ldots,l_k=0}^{N-1}
a_{l_1 \ldots l_k}
\Delta{\bf w}_{\tau_{l_1}}^{(i_1)}\ldots \Delta{\bf w}_{\tau_{l_k}}^{(i_k)},
\end{equation}

\vspace{1mm}
\noindent
where $\Delta{\bf w}_{\tau_{j}}^{(i)}=
{\bf w}_{\tau_{j+1}}^{(i)}-{\bf w}_{\tau_{j}}^{(i)},$\
$i=0, 1,\ldots,m,$\ ${\bf w}_{\tau}^{(0)}=\tau.$

It is known (see \cite{Kuo}, Lemma~9.6.4)
that for any $\Phi(t_1,\ldots,t_k)\in L_2([t, T]^k)$ 
there exists a sequence of step functions $\Phi_N(t_1,\ldots,t_k)$ of the form (\ref{chain3})
such that
\begin{equation}
\label{chain15}
~~~~~~~~\lim\limits_{N\to\infty} \int\limits_{[t,T]^k}
\left(\Phi(t_1,\ldots,t_k)-\Phi_N(t_1,\ldots,t_k)\right)^2 dt_1\ldots dt_k=0.
\end{equation}

We have
$$
\Phi_N(t_1,\ldots,t_k)=\sum\limits_{l_1,\ldots,l_k=0}^{N-1}
a_{l_1 \ldots l_k} {\bf 1}_{[\tau_{l_1},\tau_{l_1+1})}(t_1) \ldots
{\bf 1}_{[\tau_{l_k},\tau_{l_k+1})}(t_k)=
$$
\begin{equation}
\label{chain5}
~~~~~~=\sum\limits_{(l_1,\ldots,l_k)}
\sum_{\stackrel{l_1,\ldots,l_k=0}{{}_{l_1<l_2<\ldots < l_k}}}^{N-1}
a_{l_1 \ldots l_k} {\bf 1}_{[\tau_{l_1},\tau_{l_1+1})}(t_1) \ldots
{\bf 1}_{[\tau_{l_k},\tau_{l_k+1})}(t_k),
\end{equation}
where permutations $(l_1,\ldots,l_k)$ when summing are 
performed only in the expression $l_1<l_2<\ldots < l_k$
(recall that $a_{l_1 \ldots l_k}=0$ if $l_p=l_q$ for some $p\ne q$).

Using (\ref{chain5}), we get
\begin{equation}
\label{chain30}
\sum_{(t_1,\ldots,t_k)}
\int\limits_{t}^{T}
\ldots
\int\limits_{t}^{t_2}
\Phi_N(t_1,\ldots,t_k)d{\bf w}_{t_1}^{(i_1)}
\ldots
d{\bf w}_{t_k}^{(i_k)}=
\end{equation}
$$
=\sum\limits_{(l_1,\ldots,l_k)}
\sum_{\stackrel{l_1,\ldots,l_k=0}{{}_{l_1<l_2<\ldots < l_k}}}^{N-1}
a_{l_1 \ldots l_k} 
\Delta{\bf w}_{\tau_{l_1}}^{(i_1)} \ldots \Delta{\bf w}_{\tau_{l_k}}^{(i_k)}=
$$
\begin{equation}
\label{chain10}
~~~~~~=\sum\limits_{\stackrel{l_1,\ldots,l_k=0}{{}_{l_q\ne l_r;\ q\ne r;\ 
q, r=1,\ldots, k}}}^{N-1}
a_{l_1 \ldots l_k} 
\Delta{\bf w}_{\tau_{l_1}}^{(i_1)} \ldots \Delta{\bf w}_{\tau_{l_k}}^{(i_k)}
=J'[\Phi_N]_{T,t}^{(i_1\ldots i_k)}\ \ \ \hbox{w.~p.~1},
\end{equation}

\vspace{2mm}
\noindent
where permutations $(t_1,\ldots,t_k)$ when summing are 
performed only in the values
$d{\bf w}_{t_1}^{(i_1)}
\ldots $
$d{\bf w}_{t_k}^{(i_k)}$ 
and permutations $(l_1,\ldots,l_k)$ when summing are 
performed only in the expression $l_1<l_2<\ldots < l_k.$
At the same time the indices near 
upper 
limits of integration in the iterated stochastic integrals in (\ref{chain30}) are changed 
correspondently and if $t_r$ swapped with $t_q$ in the  
permutation $(t_1,\ldots,t_k)$, then $i_r$ swapped with $i_q$ in 
the permutation $(i_1,\ldots,i_k)$ (see (\ref{chain30})).
In addition, the multiple Wiener stochastic integral 
$J'[\Phi_N]_{T,t}^{(i_1\ldots i_k)}$ is defined by (\ref{chain9})
and 
$$
\int\limits_{t}^{T}
\ldots
\int\limits_{t}^{t_2}
\Phi_N(t_1,\ldots,t_k)d{\bf w}_{t_1}^{(i_1)}
\ldots
d{\bf w}_{t_k}^{(i_k)}
$$
is the iterated It\^{o} stochastic integral.

Using (\ref{chain15}), (\ref{chain10}), Lemma 1.2 for $\Phi(t_1,\ldots,t_k)\in L_2(D_k)$, 
and (\ref{riemann})
for Lebesgue integrals, we have
$$
{\sf M}\left\{\left(J'[\Phi_N]_{T,t}^{(i_1\ldots i_k)}-
J'[\Phi_M]_{T,t}^{(i_1\ldots i_k)}\right)^2\right\}\le
$$
$$
\le C_k 
\sum_{(t_1,\ldots,t_k)}
\int\limits_{t}^{T}
\ldots
\int\limits_{t}^{t_2}
\left(\Phi_N(t_1,\ldots,t_k)-\Phi_M(t_1,\ldots,t_k)\right)^2 dt_1
\ldots dt_k=
$$

\vspace{-1mm}
$$
=C_k 
\int\limits_{[t,T]^k}
\left(\Phi_N(t_1,\ldots,t_k)-\Phi_M(t_1,\ldots,t_k)\right)^2 dt_1
\ldots dt_k=
$$
$$
=C_k\left\Vert \Phi_N-\Phi_M\right\Vert_{L_2([t, T]^k)}^2\le
$$
$$
\le 2 C_k \left(\left\Vert \Phi_N-\Phi\right\Vert_{L_2([t, T]^k)}^2+
\left\Vert \Phi-\Phi_M\right\Vert_{L_2([t, T]^k)}^2\right)^2\ \to 0
$$

\vspace{2mm}
\noindent
if $N,M\to\infty,$ 
where constant $C_k$ 
depends only
on the multiplicity $k$ of the multiple Wiener stochastic integral.

Thus, there exists the limit 
$$
\hbox{\vtop{\offinterlineskip\halign{
\hfil#\hfil\cr
{\rm l.i.m.}\cr
$\stackrel{}{{}_{N\to \infty}}$\cr
}} }J'[\Phi_N]_{T,t}^{(i_1\ldots i_k)}.
$$

We will define the multiple Wiener stochastic integral for $\Phi(t_1,\ldots,t_k)\in L_2([t, T]^k)$ 
by the formula \cite{ito1951}, \cite{Kuo}
\begin{equation}
\label{WiI}
J'[\Phi]_{T,t}^{(i_1\ldots i_k)}\stackrel{\sf def}{=}
\hbox{\vtop{\offinterlineskip\halign{
\hfil#\hfil\cr
{\rm l.i.m.}\cr
$\stackrel{}{{}_{N\to \infty}}$\cr
}} }J'[\Phi_N]_{T,t}^{(i_1\ldots i_k)}=
\hbox{\vtop{\offinterlineskip\halign{
\hfil#\hfil\cr
{\rm l.i.m.}\cr
$\stackrel{}{{}_{N\to \infty}}$\cr
}} }
\sum\limits_{l_1,\ldots,l_k=0}^{N-1}
a_{l_1 \ldots l_k}
\Delta{\bf w}_{\tau_{l_1}}^{(i_1)}\ldots \Delta{\bf w}_{\tau_{l_k}}^{(i_k)},
\end{equation}
where $\Phi_N(t_1,\ldots,t_k)$ is defined by 
(\ref{chain3}),
$\Delta{\bf w}_{\tau_{j}}^{(i)}=
{\bf w}_{\tau_{j+1}}^{(i)}-{\bf w}_{\tau_{j}}^{(i)},$\
$i=0, 1,\ldots,m,$\ ${\bf w}_{\tau}^{(0)}=\tau.$

It is easy to see that the above definition coincides with 
(\ref{mult11}) if the function 
$\Phi(t_1,\ldots,t_k):\ [t, T]^k\to{\bf R}$ is continuous in the hypercube
$[t, T]^k$.

Let us prove the following equality 
\begin{equation}
\label{Wi110}
~~~J'[\Phi]_{T,t}^{(i_1\ldots i_k)}=\sum_{(t_1,\ldots,t_k)}
\int\limits_{t}^{T}
\ldots
\int\limits_{t}^{t_2}
\Phi(t_1,\ldots,t_k)d{\bf w}_{t_1}^{(i_1)}
\ldots
d{\bf w}_{t_k}^{(i_k)}\ \ \ \hbox{w.~p.~1},
\end{equation}
where permutations $(t_1,\ldots,t_k)$ when summing are 
performed only in the values
$d{\bf w}_{t_1}^{(i_1)}
\ldots $
$d{\bf w}_{t_k}^{(i_k)}.$ At the same time the indices near 
upper 
limits of integration in the iterated stochastic integrals are changed 
correspondently and if $t_r$ swapped with $t_q$ in the  
permutation $(t_1,\ldots,t_k)$, then $i_r$ swapped with $i_q$ in 
the permutation $(i_1,\ldots,i_k).$ 
In addition, the multiple Wiener stochastic integral 
$J'[\Phi]_{T,t}^{(i_1\ldots i_k)}$ is defined by (\ref{WiI})
and 
$$
\int\limits_{t}^{T}
\ldots
\int\limits_{t}^{t_2}
\Phi(t_1,\ldots,t_k)d{\bf w}_{t_1}^{(i_1)}
\ldots
d{\bf w}_{t_k}^{(i_k)}
$$
is the iterated It\^{o} stochastic integral.

The equality (\ref{Wi110}) has already been proved for the case 
$\Phi(t_1,\ldots,t_k)=\Phi_N(t_1,\ldots,t_k)$ (see (\ref{chain10})).

From (\ref{chain10}) we have

\vspace{-2mm}
$$
J'[\Phi_N]_{T,t}^{(i_1\ldots i_k)}=
\sum_{(t_1,\ldots,t_k)}
\int\limits_{t}^{T}
\ldots
\int\limits_{t}^{t_2}
\Phi_N(t_1,\ldots,t_k)d{\bf w}_{t_1}^{(i_1)}
\ldots
d{\bf w}_{t_k}^{(i_k)}=
$$
$$
=\sum_{(t_1,\ldots,t_k)}
\int\limits_{t}^{T}
\ldots
\int\limits_{t}^{t_2}
\Phi(t_1,\ldots,t_k)d{\bf w}_{t_1}^{(i_1)}
\ldots
d{\bf w}_{t_k}^{(i_k)}+
$$
\begin{equation}
\label{chain11}
+\sum_{(t_1,\ldots,t_k)}
\int\limits_{t}^{T}
\ldots
\int\limits_{t}^{t_2}
\left(\Phi_N(t_1,\ldots,t_k)-\Phi(t_1,\ldots,t_k)\right)d{\bf w}_{t_1}^{(i_1)}
\ldots
d{\bf w}_{t_k}^{(i_k)}\ \ \ \hbox{w.~p.~1.}
\end{equation}

Passing to the limit $\hbox{\vtop{\offinterlineskip\halign{
\hfil#\hfil\cr
{\rm l.i.m.}\cr
$\stackrel{}{{}_{N\to \infty}}$\cr
}} }$ in the equality (\ref{chain11}), we obtain
$$
J'[\Phi]_{T,t}^{(i_1\ldots i_k)}=
\sum_{(t_1,\ldots,t_k)}
\int\limits_{t}^{T}
\ldots
\int\limits_{t}^{t_2}
\Phi(t_1,\ldots,t_k)d{\bf w}_{t_1}^{(i_1)}
\ldots
d{\bf w}_{t_k}^{(i_k)}+
$$
\begin{equation}
\label{chain12}
+\hbox{\vtop{\offinterlineskip\halign{
\hfil#\hfil\cr
{\rm l.i.m.}\cr
$\stackrel{}{{}_{N\to \infty}}$\cr
}} }\sum_{(t_1,\ldots,t_k)}
\int\limits_{t}^{T}
\ldots
\int\limits_{t}^{t_2}
\left(\Phi_N(t_1,\ldots,t_k)-\Phi(t_1,\ldots,t_k)\right)d{\bf w}_{t_1}^{(i_1)}
\ldots
d{\bf w}_{t_k}^{(i_k)}\ \ \ \hbox{w.~p.~1.}
\end{equation}

Using Lemma 1.2 for $\Phi(t_1,\ldots,t_k)\in L_2(D_k)$, 
(\ref{riemann}) for Lebesgue integrals, and (\ref{chain15}), we get
$$
{\sf M}\left\{\left(
\sum_{(t_1,\ldots,t_k)}
\int\limits_{t}^{T}
\ldots
\int\limits_{t}^{t_2}
\left(\Phi_N(t_1,\ldots,t_k)-\Phi(t_1,\ldots,t_k)\right)d{\bf w}_{t_1}^{(i_1)}
\ldots
d{\bf w}_{t_k}^{(i_k)}\right)^2\right\}\le
$$
$$
\le C_k 
\sum_{(t_1,\ldots,t_k)}
\int\limits_{t}^{T}
\ldots
\int\limits_{t}^{t_2}
\left(\Phi_N(t_1,\ldots,t_k)-\Phi(t_1,\ldots,t_k)\right)^2 dt_1
\ldots dt_k=
$$

\vspace{-1mm}
\begin{equation}
\label{chain20}
~~~~~~~~=C_k 
\int\limits_{[t,T]^k}
\left(\Phi_N(t_1,\ldots,t_k)-\Phi(t_1,\ldots,t_k)\right)^2 dt_1
\ldots dt_k\ \to 0
\end{equation}
if $N\to\infty,$ 
where constant $C_k$ 
depends only
on the multiplicity $k$ of the multiple Wiener stochastic integral.
The relations (\ref{chain12}) and (\ref{chain20}) prove the equality 
(\ref{Wi110}).

Using (\ref{Wi110}) and the isometry property of the It\^{o} stochastic integral, we have
\begin{equation}            
\label{wi1001}
J[\psi^{(k)}]_{T,t}^{(i_1\ldots i_k)}=\int\limits_t^T\psi_k(t_k) \ldots \int\limits_t^{t_{2}}
\psi_1(t_1) d{\bf w}_{t_1}^{(i_1)}\ldots
d{\bf w}_{t_k}^{(i_k)}=J'[K]_{T,t}^{(i_1\ldots i_k)}\ \ \ \hbox{w.~p.~1},
\end{equation}
where 
$K=K(t_1,\ldots,t_k)$ is defined by (\ref{ppp}), i.e.
\begin{equation}
\label{chain200}
K(t_1,\ldots,t_k)=
\left\{\begin{matrix}
\psi_1(t_1)\ldots \psi_k(t_k),\ &t_1<\ldots<t_k\cr\cr
0,\ &\hbox{\rm otherwise}
\end{matrix}
\right.\ \ \ \left(\psi_l(\tau)\in L_2([t,T])\right),
\end{equation}
where $l=1,\ldots,k,$\ $t_1,\ldots,t_k\in [t, T]$ $(k\ge 2)$ and 
$K(t_1)\equiv\psi_1(t_1)$ for $t_1\in[t, T].$

Applying (\ref{wi1001}) and the linearity property of the It\^{o} stochastic integral, we obtain

\vspace{-5mm}
$$
J[\psi^{(k)}]_{T,t}^{(i_1\ldots i_k)}=J'[K]_{T,t}^{(i_1\ldots i_k)}=
$$

\vspace{-3mm}
\begin{equation}
\label{chain102}
~~~~=\sum_{j_1=0}^{p_1}\ldots
\sum_{j_k=0}^{p_k}
C_{j_k\ldots j_1}
J'[\phi_{j_1}\ldots \phi_{j_k}]_{T,t}^{(i_1\ldots i_k)}+
J'[R_{p_1\ldots p_k}]_{T,t}^{(i_1\ldots i_k)}\ \ \ \hbox{w.~p.~1,}
\end{equation}

\noindent
where
\begin{equation}
\label{chain30001}
~~~~~~~R_{p_1\ldots p_k}(t_1,\ldots,t_k)\stackrel{{\rm def}}{=}
K(t_1,\ldots,t_k)-
\sum_{j_1=0}^{p_1}\ldots
\sum_{j_k=0}^{p_k}
C_{j_k\ldots j_1}
\prod_{l=1}^k\phi_{j_l}(t_l)
\end{equation}
and
\begin{equation}
\label{chain300}
C_{j_k\ldots j_1}=\int\limits_{[t,T]^k}
K(t_1,\ldots,t_k)\prod_{l=1}^{k}\phi_{j_l}(t_l)dt_1\ldots dt_k
\end{equation}

\noindent
is the Fourier coefficient corresponding to $K(t_1,\ldots,t_k).$

Using the It\^{o} formula, we have

\vspace{-3mm}
$$
\sum\limits_{(j_1,\ldots,j_q)}\int\limits_t^T \phi_{j_q}(t_q)\ldots
\int\limits_t^{t_2}\phi_{j_1}(t_1)
d{\bf w}_{t_1}^{(i_1)}\ldots {\bf w}_{t_q}^{(i_q)}\times
$$
$$
\times \sum\limits_{(j_1',\ldots,j_n')}
\int\limits_t^T 
\phi_{j_n'}(t_n')\ldots \int\limits_t^{t_2'}\phi_{j_1'}(t_1')d{\bf w}_{t_1'}^{(g)}\ldots {\bf w}_{t_n'}^{(g)}=
$$
$$
=\sum\limits_{(j_1,\ldots,j_q,j_1',\ldots,j_n')}
\int\limits_t^T \phi_{j_q}(t_q)\ldots \int\limits_t^{t_2}\phi_{j_1}(t_1)
\int\limits_t^{t_1}\phi_{j_n'}(t_n')\ldots \int\limits_t^{t_2'}
\phi_{j_1'}(t_1')\times
$$

\vspace{1mm}
\begin{equation}
\label{new1000}
\times d{\bf w}_{t_1'}^{(g)}\ldots d{\bf w}_{t_n'}^{(g)}d{\bf w}_{t_1}^{(i_1)}\ldots d{\bf w}_{t_q}^{(i_q)}
\end{equation}

\vspace{5mm}
\noindent
w.~p.~1, where $g=0$ or $g=1,$\ $n, q\in {\bf N},$\ $i_1,\ldots,i_q\ne 0,\ 1,$

\vspace{-2mm}
$$
\sum\limits_{(j_1,\ldots,j_k)}
$$

\vspace{2mm}
\noindent
means the sum with respect to all possible permutations $(j_1,\ldots,j_k)$.
At the same time if $j_r$ swapped with $j_d$ in the permutation $(j_1,\ldots,j_k)$, then
$i_r$ swapped with $i_d$ in the permutation $(i_1,\ldots,i_k).$

The detailed proof of (\ref{new1000}) will be given in Sect.~1.14 (see the proof of Theorem~1.22).
The equality (\ref{new1000}) means that (see (\ref{Wi110}))

$$
J'[\phi_{j_1}\ldots
\phi_{j_q}]_{T,t}^{(i_1\ldots i_q)}\cdot J'[\phi_{j_1'}\ldots\phi_{j_n'}]^{(g\ldots g)}_{T,t}=
$$

\begin{equation}
\label{2023abc111}
=J'[\phi_{j_1}\ldots\phi_{j_q}\phi_{j_1'}\ldots\phi_{j_n'}]_{T,t}^{(i_1\ldots i_q g\ldots g)}
\end{equation}

\vspace{4mm}
\noindent
w.~p.~1, where $g=0$ or $g=1,$\ $n, q\in \{0\}\cup {\bf N},$\ $i_1,\ldots,i_q\ne 0,\ 1,$ 
and $J'[\phi_{j_1}\ldots\phi_{j_q}]_{T,t}^{(i_1\ldots i_q)}\stackrel{\sf def}{=}1$ for $q=0.$

Using the equality (\ref{2023abc111}), we obtain (\ref{ziko30}) for the case
of an arbitrary complete orthonormal system $\left\{\phi_j(x)\right\}_{j=0}^{\infty}$
of functions in $L_2([t,T])$.

Suppose that the conditions ($\star\star$) (see Sect.~1.10) and  (\ref{ziko999}) are fulfilled. 
Applying Theorem~9.6.9 \cite{Kuo} (also see 
\cite{ito1951}, Theorem~3.1) and (\ref{ziko100}) (also see Theorem~1.23 below), we get

$$
J'[\phi_{j_1}\ldots \phi_{j_k}]_{T,t}^{(i_1\ldots i_k)}=
$$

\vspace{1mm}
$$
=\prod_{l=1}^k\left({\bf 1}_{\{m_l=0\}}+{\bf 1}_{\{m_l>0\}}\left\{
\begin{matrix}
H_{n_{1,l}}\left(\zeta_{j_{h_{1,l}}}^{(i_l)}\right)\ldots 
H_{n_{d_l,l}}\left(\zeta_{j_{h_{d_l,l}}}^{(i_l)}\right),\ 
&\hbox{\rm if}\ \ \ 
i_l\ne 0\cr\cr
\left(\zeta_{j_{h_{1,l}}}^{(0)}\right)^{n_{1,l}}\ldots
\left(\zeta_{j_{h_{d_l,l}}}^{(0)}\right)^{n_{d_l,l}},\  &\hbox{\rm if}\ \ \ 
i_l=0
\end{matrix}\right.\ \right)=
$$

\newpage
\noindent
$$
=\prod_{l=1}^k\zeta_{j_l}^{(i_l)}
+\sum\limits_{r=1}^{[k/2]}
(-1)^r \times
$$

\vspace{-3mm}
\begin{equation}
\label{chain401}
~~\times\sum_{\stackrel{(\{\{g_1, g_2\}, \ldots, 
\{g_{2r-1}, g_{2r}\}\}, \{q_1, \ldots, q_{k-2r}\})}
{{}_{\{g_1, g_2, \ldots, 
g_{2r-1}, g_{2r}, q_1, \ldots, q_{k-2r}\}=\{1, 2, \ldots, k\}}}}
\prod\limits_{s=1}^r
{\bf 1}_{\{i_{g_{{}_{2s-1}}}=~i_{g_{{}_{2s}}}\ne 0\}}
\Biggl.{\bf 1}_{\{j_{g_{{}_{2s-1}}}=~j_{g_{{}_{2s}}}\}}
\prod_{l=1}^{k-2r}\zeta_{j_{q_l}}^{(i_{q_l})}
\end{equation}

\vspace{2mm}
\noindent
w.~p.~1, where notations are the same as in Theorems 1.2 and 1.14;
the multiple Wiener stochastic integral 
$J'[\phi_{j_1}\ldots \phi_{j_k}]_{T,t}^{(i_1\ldots i_k)}$ is defined by 
(\ref{WiI}).

Again applying (\ref{Wi110}), we have

\vspace{-4mm}
$$
J'[R_{p_1\ldots p_k}]_{T,t}^{(i_1\ldots i_k)}
=
\hspace{-1.9mm}\sum_{(t_1,\ldots,t_k)}\hspace{-0.5mm}
\int\limits_{t}^{T}
\ldots
\int\limits_{t}^{t_2}\hspace{-0.9mm}
\left(\hspace{-0.5mm}K(t_1,\ldots,t_k)-
\sum_{j_1=0}^{p_1}\ldots
\sum_{j_k=0}^{p_k}
C_{j_k\ldots j_1}
\prod_{l=1}^k\phi_{j_l}(t_l)\hspace{-0.5mm}\right)\hspace{-0.9mm}\times
$$

\vspace{-1mm}
\begin{equation}
\label{wi2005}
\times
d{\bf w}_{t_1}^{(i_1)}
\ldots
d{\bf w}_{t_k}^{(i_k)},
\end{equation}

\vspace{1mm}
\noindent
where permutations $(t_1,\ldots,t_k)$ when summing are performed only 
in the values $d{\bf w}_{t_1}^{(i_1)}
\ldots $
$d{\bf w}_{t_k}^{(i_k)}$. At the same time the indices near 
upper limits of integration in the iterated stochastic integrals 
are changed correspondently and if $t_r$ swapped with $t_q$ in the  
permutation $(t_1,\ldots,t_k)$, then $i_r$ swapped with $i_q$ in the 
permutation $(i_1,\ldots,i_k).$
In addition, the multiple Wiener stochastic integral
$J'[R_{p_1\ldots p_k}]_{T,t}^{(i_1\ldots i_k)}$ is defined by 
(\ref{WiI}).

According to Lemma 1.2 for $\Phi(t_1,\ldots,t_k)\in L_2(D_k)$, 
(\ref{sos1z}), and (\ref{riemann}) for Lebesgue integrals, we have

\vspace{-3mm}
$$
{\sf M}\left\{\left(J'[R_{p_1\ldots p_k}]_{T,t}^{(i_1\ldots i_k)}\right)^2\right\}
\le 
$$
$$
\le C_k \hspace{-0.4mm}
\sum_{(t_1,\ldots,t_k)}
\int\limits_{t}^{T}
\ldots
\int\limits_{t}^{t_2}
\left(K(t_1,\ldots,t_k)-
\sum_{j_1=0}^{p_1}\ldots
\sum_{j_k=0}^{p_k}
C_{j_k\ldots j_1}
\prod_{l=1}^k\phi_{j_l}(t_l)\right)^2
\hspace{-1mm}dt_1
\ldots
dt_k\hspace{-0.2mm}=
$$
\begin{equation}
\label{chain7771}
=C_k\int\limits_{[t,T]^k}
\left(K(t_1,\ldots,t_k)-
\sum_{j_1=0}^{p_1}\ldots
\sum_{j_k=0}^{p_k}
C_{j_k\ldots j_1}
\prod_{l=1}^k\phi_{j_l}(t_l)\right)^2
dt_1
\ldots
dt_k\to 0
\end{equation}
if $p_1,\ldots,p_k\to\infty,$ where constant $C_k$ 
depends only
on the multiplicity $k$ of the
iterated  It\^{o} stochastic integral
$J[\psi^{(k)}]_{T,t}^{(i_1\ldots i_k)}$. 
Thus (see (\ref{chain102}) and (\ref{chain7771})), 
\newpage
\noindent
\begin{equation}
\label{febr5000}
~~~~~~~~~~J[\psi^{(k)}]_{T,t}^{(i_1\ldots i_k)}=
\hbox{\vtop{\offinterlineskip\halign{
\hfil#\hfil\cr
{\rm l.i.m.}\cr
$\stackrel{}{{}_{p_1,\ldots,p_k\to \infty}}$\cr
}} }
\sum\limits_{j_1=0}^{p_1}\ldots
\sum\limits_{j_k=0}^{p_k}
C_{j_k\ldots j_1}J'[\phi_{j_1}\ldots \phi_{j_k}]_{T,t}^{(i_1\ldots i_k)}
\end{equation}
and the following theorem is proved.

{\bf Theorem 1.16}\ \cite{arxiv-1}\ (generalization of Theorems 1.1, 1.2, and 1.14).\ 
{\it Suppose that
the condition {\rm ($\star\star$)} is fulfilled
for the multi-index $(i_1 \ldots i_k)$ {\rm (}see Sect.~{\rm 1.10)} 
and the condition {\rm (\ref{ziko999})} is also 
fulfilled.
Furthermore$,$ let 
$\psi_l(\tau)\in L_2([t, T])$ $(l=$ $1,\ldots, k)$ and
$\{\phi_j(x)\}_{j=0}^{\infty}$ is an arbitrary complete orthonormal system  
of functions in the space $L_2([t,T]).$
Then the following expansions

\vspace{-2mm}
$$
J[\psi^{(k)}]_{T,t}^{(i_1\ldots i_k)}=
\hbox{\vtop{\offinterlineskip\halign{
\hfil#\hfil\cr
{\rm l.i.m.}\cr
$\stackrel{}{{}_{p_1,\ldots,p_k\to \infty}}$\cr
}} }
\sum\limits_{j_1=0}^{p_1}\ldots
\sum\limits_{j_k=0}^{p_k}
C_{j_k\ldots j_1}\times
$$

\vspace{1.5mm}
\begin{equation}
\label{new9999}
\times
\prod_{l=1}^k\left({\bf 1}_{\{m_l=0\}}+{\bf 1}_{\{m_l>0\}}\left\{
\begin{matrix}
H_{n_{1,l}}\left(\zeta_{j_{h_{1,l}}}^{(i_l)}\right)\ldots 
H_{n_{d_l,l}}\left(\zeta_{j_{h_{d_l,l}}}^{(i_l)}\right),\ 
&\hbox{\rm if}\ \ \ 
i_l\ne 0\cr\cr
\left(\zeta_{j_{h_{1,l}}}^{(0)}\right)^{n_{1,l}}\ldots
\left(\zeta_{j_{h_{d_l,l}}}^{(0)}\right)^{n_{d_l,l}},\  &\hbox{\rm if}\ \ \ 
i_l=0
\end{matrix}\right.\ \right),
\end{equation}

\vspace{-3mm}
$$
J[\psi^{(k)}]_{T,t}^{(i_1\ldots i_k)}=
\hbox{\vtop{\offinterlineskip\halign{
\hfil#\hfil\cr
{\rm l.i.m.}\cr
$\stackrel{}{{}_{p_1,\ldots,p_k\to \infty}}$\cr
}} }
\sum\limits_{j_1=0}^{p_1}\ldots
\sum\limits_{j_k=0}^{p_k}
C_{j_k\ldots j_1}\Biggl(
\prod_{l=1}^k\zeta_{j_l}^{(i_l)}+\sum\limits_{r=1}^{[k/2]}
(-1)^r \times
\Biggr.
$$

\begin{equation}
\label{razzar1}
\times
\sum_{\stackrel{(\{\{g_1, g_2\}, \ldots, 
\{g_{2r-1}, g_{2r}\}\}, \{q_1, \ldots, q_{k-2r}\})}
{{}_{\{g_1, g_2, \ldots, 
g_{2r-1}, g_{2r}, q_1, \ldots, q_{k-2r}\}=\{1, 2, \ldots, k\}}}}
\prod\limits_{s=1}^r
{\bf 1}_{\{i_{g_{{}_{2s-1}}}=~i_{g_{{}_{2s}}}\ne 0\}}
\Biggl.{\bf 1}_{\{j_{g_{{}_{2s-1}}}=~j_{g_{{}_{2s}}}\}}
\prod_{l=1}^{k-2r}\zeta_{j_{q_l}}^{(i_{q_l})}\Biggr)
\end{equation}

\vspace{0.5mm}
\noindent
con\-verg\-ing in the mean-square sense are valid$,$
where $[x]$ is an integer part of a real number $x;$\ \
$n_{1,l}+n_{2,l}+\ldots+n_{d_l,l}=m_l;$\ \ $n_{1,l}, n_{2,l}, \ldots, n_{d_l,l}=1,\ldots, m_l;$\ \ 
$d_l=1,\ldots,m_l;$\ \ $l=1,\ldots,k;$\ \ $m_1+\ldots+m_k=k;$\ \  
the numbers $m_1,\ldots,m_k,$\ $g_1,\ldots,g_k$
depend on $(i_1,\ldots,i_k)$ and 
the numbers $n_{1,l},\ldots,n_{d_l,l},$\ $h_{1,l},\ldots,h_{d_l,l},$\ $d_l$
depend on $\{j_1,\ldots,j_k\};$ moreover$,$ $\left\{j_{g_1},\ldots,j_{g_k}\right\}
=\{j_1,\ldots,j_k\};$ $H_n(x)$ is the Hermite polynomial {\rm (\ref{ziko500});}
another
notations as in Theorems {\rm 1.1, 1.2, and 1.14}.}

Replacing the function $K(t_1,\ldots,t_k)$ by 
$\Phi(t_1,\ldots,t_k)$ we get the following theorem.

{\bf Theorem 1.17}\ \cite{arxiv-1}\ (generalization of Theorems 1.13, 1.15).\ {\it Suppose that
the condition {\rm ($\star\star$)} is fulfilled
for the multi-index $(i_1 \ldots i_k)$ {\rm (}see Sect.~{\rm 1.10)} 
and the condition {\rm (\ref{ziko999})} is also 
fulfilled.
Furthermore$,$ let 
$\Phi(t_1,\ldots,t_k)\in L_2([t, T]^k)$ 
and
$\{\phi_j(x)\}_{j=0}^{\infty}$ is an arbitrary complete orthonormal system  
of functions in the 
space $L_2([t,T]).$
Then the following expansions

\vspace{-2mm}
$$
J'[\Phi]_{T,t}^{(i_1\ldots i_k)}=
\hbox{\vtop{\offinterlineskip\halign{
\hfil#\hfil\cr
{\rm l.i.m.}\cr
$\stackrel{}{{}_{p_1,\ldots,p_k\to \infty}}$\cr
}} }
\sum\limits_{j_1=0}^{p_1}\ldots
\sum\limits_{j_k=0}^{p_k}
C_{j_k\ldots j_1}\times
$$

\vspace{1mm}
\begin{equation}
\label{new1001}
\times
\prod_{l=1}^k\left({\bf 1}_{\{m_l=0\}}+{\bf 1}_{\{m_l>0\}}\left\{
\begin{matrix}
H_{n_{1,l}}\left(\zeta_{j_{h_{1,l}}}^{(i_l)}\right)\ldots 
H_{n_{d_l,l}}\left(\zeta_{j_{h_{d_l,l}}}^{(i_l)}\right),\ 
&\hbox{\rm if}\ \ \ 
i_l\ne 0\cr\cr
\left(\zeta_{j_{h_{1,l}}}^{(0)}\right)^{n_{1,l}}\ldots
\left(\zeta_{j_{h_{d_l,l}}}^{(0)}\right)^{n_{d_l,l}},\  &\hbox{\rm if}\ \ \ 
i_l=0
\end{matrix}\right.\ \right),
\end{equation}

\vspace{-3mm}
$$
J'[\Phi]_{T,t}^{(i_1\ldots i_k)}=
\hbox{\vtop{\offinterlineskip\halign{
\hfil#\hfil\cr
{\rm l.i.m.}\cr
$\stackrel{}{{}_{p_1,\ldots,p_k\to \infty}}$\cr
}} }
\sum\limits_{j_1=0}^{p_1}\ldots
\sum\limits_{j_k=0}^{p_k}
C_{j_k\ldots j_1}\Biggl(
\prod_{l=1}^k\zeta_{j_l}^{(i_l)}+\sum\limits_{r=1}^{[k/2]}
(-1)^r \times
\Biggr.
$$

\begin{equation}
\label{january20}
\times
\sum_{\stackrel{(\{\{g_1, g_2\}, \ldots, 
\{g_{2r-1}, g_{2r}\}\}, \{q_1, \ldots, q_{k-2r}\})}
{{}_{\{g_1, g_2, \ldots, 
g_{2r-1}, g_{2r}, q_1, \ldots, q_{k-2r}\}=\{1, 2, \ldots, k\}}}}
\prod\limits_{s=1}^r
{\bf 1}_{\{i_{g_{{}_{2s-1}}}=~i_{g_{{}_{2s}}}\ne 0\}}
\Biggl.{\bf 1}_{\{j_{g_{{}_{2s-1}}}=~j_{g_{{}_{2s}}}\}}
\prod_{l=1}^{k-2r}\zeta_{j_{q_l}}^{(i_{q_l})}\Biggr)
\end{equation}

\noindent
con\-verg\-ing in the mean-square sense are valid$,$
where $[x]$ is an integer part of a real number $x;$\ \
$n_{1,l}+n_{2,l}+\ldots+n_{d_l,l}=m_l;$\ \ $n_{1,l}, n_{2,l}, \ldots, n_{d_l,l}=1,\ldots, m_l;$\ \ 
$d_l=1,\ldots,m_l;$\ \ $l=1,\ldots,k;$\ \ $m_1+\ldots+m_k=k;$\ \  
the numbers $m_1,\ldots,m_k,$\ $g_1,\ldots,g_k$
depend on $(i_1,\ldots,i_k)$ and 
the numbers $n_{1,l},\ldots,n_{d_l,l},$\ $h_{1,l},\ldots,h_{d_l,l},$\ $d_l$
depend on $\{j_1,\ldots,j_k\};$ moreover$,$ $\left\{j_{g_1},\ldots,j_{g_k}\right\}
=\{j_1,\ldots,j_k\};$
the multiple Wiener stochastic integral
$J'[\Phi]_{T,t}^{(i_1\ldots i_k)}$ is defined by 
{\rm (\ref{WiI});} $H_n(x)$ is the Hermite polynomial {\rm (\ref{ziko500});}
another
notations as in Theorems {\rm 1.13, 1.15}.}

\vspace{2mm}

It should be noted that an analogue of Theorem 1.17 (more precisely, the expansion
like (\ref{new1001}) for the case $i_1,\ldots,i_k\ne 0$) was considered 
in \cite{Rybakov1000}. Also note that the proof in \cite{Rybakov1000} is 
different from the proof given in this section.         
In \cite{Rybakov1000}, the author interprets the multiple Wiener
stochastic integral from a finite-dimensional kernel
as a linear operator and proves that this operator is bounded.
In our proof of Theorems~1.16, 1.17 we several times use
the representation (\ref{Wi110}) of the multiple Wiener stochastic 
integral as the sum (with respect to permutations)
of iterated It\^{o} stochastic integrals and then estimate the 
remainder of the series (see (\ref{chain7771}) for details).

Note that the results of work \cite{Rybakov1000}, as well as 
the results of Chapter~1 of this book, are based on our idea 
\cite{1} (2006) on the expansion of the kernel (\ref{ppp}) (or $\Phi(t_1,\ldots,t_k)\in L_2([t,T]^k)$)
into a generalized multiple Fourier series 
(see \cite{1}, Chapter~5, Theorem~5.1, pp.~235-245 
or Sect.~1.1.3 of this book for details).

\section{Generalization of Theorems 1.3, 1.4 to the Case of an Arbitrary 
Complete Ortho\-nor\-mal System of Functions in the Space $L_2([t, T])$
and $\psi_1(\tau),$ $\ldots,$ $\psi_k(\tau)$ $\in $ $L_2([t, T])$}

In this section, we will use the multiple Wiener 
stochastic integral 
with respect 
to the components of a multidimensional Wiener
process
to generalize Theorems 
1.3, 1.4 to the case of an arbitrary 
complete orthonormal system of functions in the space $L_2([t, T])$
and $\psi_1(\tau),$ $\ldots,\psi_k(\tau)\in L_2([t, T]).$

{\bf Theorem 1.18.}\
{\it Suppose that
$\psi_1(\tau),\ldots,\psi_k(\tau)\in L_2([t, T])$ 
and
$\{\phi_j(x)\}_{j=0}^{\infty}$ is an arbitrary complete orthonormal system  
of functions in the space $L_2([t,T]).$ 
Then
$$
{\sf M}\left\{\left(J[\psi^{(k)}]_{T,t}-
J[\psi^{(k)}]_{T,t}^p\right)^2\right\}
= \int\limits_{[t,T]^k} K^2(t_1,\ldots,t_k)
dt_1\ldots dt_k - 
$$

\vspace{-4mm}
\begin{equation}
\label{chain100}
- \hspace{-0.5mm}\sum_{j_1=0}^{p}\ldots\sum_{j_k=0}^{p}
C_{j_k\ldots j_1}
{\sf M}\left\{\hspace{-0.5mm}J[\psi^{(k)}]_{T,t}
\sum\limits_{(j_1,\ldots,j_k)}
\int\limits_t^T \phi_{j_k}(t_k)
\ldots
\int\limits_t^{t_{2}}\phi_{j_{1}}(t_{1})
d{\bf w}_{t_1}^{(i_1)}\ldots
d{\bf w}_{t_k}^{(i_k)}\hspace{-0.5mm}\right\}\hspace{-0.6mm},
\end{equation}

\vspace{2mm}
\noindent
where
$$
J[\psi^{(k)}]_{T,t}=\int\limits_t^T\psi_k(t_k) \ldots \int\limits_t^{t_{2}}
\psi_1(t_1) d{\bf w}_{t_1}^{(i_1)}\ldots
d{\bf w}_{t_k}^{(i_k)},
$$

\vspace{-2mm}
\begin{equation}
\label{chain101}
J[\psi^{(k)}]_{T,t}^p=
\sum_{j_1=0}^{p}\ldots\sum_{j_k=0}^{p}
C_{j_k\ldots j_1} J'[\phi_{j_1}\ldots \phi_{j_k}]_{T,t}^{(i_1\ldots i_k)},
\end{equation}

\vspace{2mm}

\noindent
$J'[\phi_{j_1}\ldots \phi_{j_k}]_{T,t}^{(i_1\ldots i_k)}$
is the multiple Wiener stochastic integral 
defined by {\rm (\ref{WiI}),}
the Fourier coefficient $C_{j_k\ldots j_1}$ has the form {\rm (\ref{chain300}),}
$K(t_1,\ldots,t_k)$ is defined by {\rm (\ref{chain200}),}
$$
\zeta_{j}^{(i)}=
\int\limits_t^T \phi_{j}(s) d{\bf w}_s^{(i)}
$$
are independent standard Gaussian random variables
for various
$i$ or $j$ $(i=1,\ldots,m),$
$$
\sum\limits_{(j_1,\ldots,j_k)}
$$ 

\noindent
means the sum with respect to all
possible permutations 
$(j_1,\ldots,j_k).$ At the same time if 
$j_r$ swapped with $j_q$ in the permutation $(j_1,\ldots,j_k)$,
then $i_r$ swapped with $i_q$ in the permutation
$(i_1,\ldots,i_k)$ {\rm (}see {\rm (\ref{chain100})).}}

{\bf Proof.}\ First, note that the formula (\ref{chain101}) 
appears due to the equality (\ref{chain102}).
Using the equality (\ref{Wi110}), we get
\begin{equation}
\label{chain103}
J'[\phi_{j_1}\ldots \phi_{j_k}]_{T,t}^{(i_1\ldots i_k)}=\sum_{(t_1,\ldots,t_k)}
\int\limits_t^T \phi_{j_k}(t_k)
\ldots
\int\limits_t^{t_{2}}\phi_{j_{1}}(t_{1})
d{\bf w}_{t_1}^{(i_1)}\ldots
d{\bf w}_{t_k}^{(i_k)}\ \ \ \hbox{w.~p.~1},
\end{equation}

\noindent
where permutations $(t_1,\ldots,t_k)$ when summing are 
performed only in the values
$d{\bf w}_{t_1}^{(i_1)}
\ldots $
$d{\bf w}_{t_k}^{(i_k)}.$ At the same time the indices near 
upper 
limits of integration in the iterated stochastic integrals are changed 
correspondently and if $t_r$ swapped with $t_q$ in the  
permutation $(t_1,\ldots,t_k)$, then $i_r$ swapped with $i_q$ in 
the permutation $(i_1,\ldots,i_k).$ 

It is easy to see that the equality (\ref{chain103}) can be written in the form
\begin{equation}
\label{chain104}
J'[\phi_{j_1}\ldots \phi_{j_k}]_{T,t}^{(i_1\ldots i_k)}=
\sum\limits_{(j_1,\ldots,j_k)}
\int\limits_t^T \phi_{j_k}(t_k)
\ldots
\int\limits_t^{t_{2}}\phi_{j_{1}}(t_{1})
d{\bf w}_{t_1}^{(i_1)}\ldots
d{\bf w}_{t_k}^{(i_k)}\ \ \ \hbox{w.~p.~1},
\end{equation}
where 
$$
\sum\limits_{(j_1,\ldots,j_k)}
$$

\noindent
means the sum with respect to all
possible permutations 
$(j_1,\ldots,j_k).$ At the same time if 
$j_r$ swapped with $j_q$ in the permutation $(j_1,\ldots,j_k)$,
then $i_r$ swapped with $i_q$ in the permutation $(i_1,\ldots,i_k).$

Thus, an analogue of the equality (\ref{ttt2}) is proved under the conditions of Theorem~1.18
(compare (\ref{yeee2}), (\ref{ttt2}) and (\ref{chain101}), (\ref{chain104})).
Further proof of Theorem~1.18 is similar to the proof of Theorem~1.3. 
Theorem~1.18 is proved.

Consider the following obvious generalization of Theorem 1.4.

{\bf Theorem 1.19.} 
{\it Suppose that
$\psi_1(\tau),\ldots,\psi_k(\tau)\in L_2([t, T])$ 
and
$\{\phi_j(x)\}_{j=0}^{\infty}$ is an arbitrary complete orthonormal system  
of functions in the space $L_2([t,T]).$ 
Then the estimate
$$
{\sf M}\left\{\left(
J[\psi^{(k)}]_{T,t}-J[\psi^{(k)}]_{T,t}^{p_1,\ldots,p_k}
\right)^2\right\}
\le 
$$

\vspace{-4mm}
\begin{equation}
\label{z1new}
~~~~\le k!\left(~\int\limits_{[t,T]^k}
K^2(t_1,\ldots,t_k)
dt_1\ldots dt_k -\sum_{j_1=0}^{p_1}\ldots
\sum_{j_k=0}^{p_k}C^2_{j_k\ldots j_1}\right)
\end{equation}

\vspace{1mm}
\noindent
is valid for the following cases{\rm :}

{\rm 1.}\ $i_1,\ldots,i_k=1,\ldots,m$\ \ and\ \ $0<T-t<\infty,$

{\rm 2.}\ $i_1,\ldots,i_k=0, 1,\ldots,m,$\ \ $i_1^2+\ldots+i_k^2>0,$\ \
and\ \ $0<T-t<1,$

\noindent
where $J[\psi^{(k)}]_{T,t}$ is the iterated It\^{o} stochastic integral {\rm (\ref{ito}),}
$J[\psi^{(k)}]_{T,t}^{p_1,\ldots,p_k}$ is the 
expression on the right-hand side of {\rm (\ref{razzar1})} before
passing to the limit 
$\hbox{\vtop{\offinterlineskip\halign{
\hfil#\hfil\cr
{\rm l.i.m.}\cr
$\stackrel{}{{}_{p_1,\ldots,p_k\to \infty}}$\cr
}} };$ another 
notations are the same as in Theorems {\rm 1.1, 1.2, 1.16}.
}

\vspace{2mm}

In addition, under the conditions of Theorem 1.19 we have the estimate
(also see (\ref{2026ch1001s11}))
$$
{\sf M}\left\{\left(J[\psi^{(k)}]_{T,t}-
J[\psi^{(k)}]_{T,t}^{p_1,\ldots,p_k}\right)^{2n}\right\}\le
$$

$$
\le
(k!)^{n}
(2n-1)^{nk}\ \times
$$

\vspace{-4mm}
\begin{equation}
\label{start1000ww}
~~~~~~~~~~\times\ 
\left(~
\int\limits_{[t,T]^k}
K^2(t_1,\ldots,t_k)
dt_1\ldots dt_k -\sum_{j_1=0}^{p_1}\ldots
\sum_{j_k=0}^{p_k}C^2_{j_k\ldots j_1}
\right)^n.
\end{equation}

\section{Generalization of Theorems 1.5, 1.6 to the Case of an Arbitrary 
Complete Ortho\-nor\-mal with Weight $r(x)\ge 0$ System of Functions in the Space $L_2([t, T])$
and $\psi_1(x)\sqrt{r(x)},$ $\ldots,$ $\psi_k(x)\sqrt{r(x)}$ $\in $ $L_2([t, T])$}

In this section, we will use the multiple Wiener 
stochastic integral 
with respect 
to the components of a multidimensional Wiener
process
to generalize Theorems 
1.5, 1.6 to the case of an arbitrary 
complete ortho\-nor\-mal with weight $r(x)\ge 0$ system of functions in the 
space $L_2([t, T])$
and $\psi_1(x)\sqrt{r(x)},$ $\ldots,$ $\psi_k(x)\sqrt{r(x)}$ $\in $ $L_2([t, T])$.
From the results of Sect.~1.3, 1.11 we obtain the following two theorems.

{\bf Theorem 1.20.}
{\it Suppose that $\psi_1(x)\sqrt{r(x)},\ldots,
\psi_k(x)\sqrt{r(x)}\in L_2([t, T]),$ 
where $r(x)\ge 0.$
Moreover$,$ let 
$$
\left\{\Psi_j(x)\sqrt{r(x)}\right\}_{j=0}^{\infty}
$$

\noindent 
is an arbitrary complete orthonormal 
system of functions in the space $L_2([t,T]).$
Then$,$ for the iterated It\^{o} stochastic integral
\begin{equation}
\label{fifi1}
~~~~~~~{\tilde J}[\psi^{(k)}]_{T,t}=\int\limits_t^T\psi_k(t_k)\sqrt{r(t_k)} 
\ldots \int\limits_t^{t_{2}}
\psi_1(t_1)\sqrt{r(t_1)} d{\bf w}_{t_1}^{(i_1)}\ldots
d{\bf w}_{t_k}^{(i_k)}
\end{equation}

\noindent
the following expansion 
$$
{\tilde J}[\psi^{(k)}]_{T,t}=
\hbox{\vtop{\offinterlineskip\halign{
\hfil#\hfil\cr
{\rm l.i.m.}\cr
$\stackrel{}{{}_{p_1,\ldots,p_k\to \infty}}$\cr
}} }
\sum\limits_{j_1=0}^{p_1}\ldots
\sum\limits_{j_k=0}^{p_k}
\tilde C_{j_k\ldots j_1}\Biggl(
\prod_{l=1}^k {\tilde \zeta}_{j_l}^{(i_l)}+\sum\limits_{r=1}^{[k/2]}
(-1)^r \times
\Biggr.
$$

\vspace{-1mm}
\begin{equation}
\label{fifi2}
\times
\sum_{\stackrel{(\{\{g_1, g_2\}, \ldots, 
\{g_{2r-1}, g_{2r}\}\}, \{q_1, \ldots, q_{k-2r}\})}
{{}_{\{g_1, g_2, \ldots, 
g_{2r-1}, g_{2r}, q_1, \ldots, q_{k-2r}\}=\{1, 2, \ldots, k\}}}}
\prod\limits_{s=1}^r
{\bf 1}_{\{i_{g_{{}_{2s-1}}}=~i_{g_{{}_{2s}}}\ne 0\}}
\Biggl.{\bf 1}_{\{j_{g_{{}_{2s-1}}}=~j_{g_{{}_{2s}}}\}}
\prod_{l=1}^{k-2r} {\tilde \zeta}_{j_{q_l}}^{(i_{q_l})}\Biggr)
\end{equation}

\vspace{1mm}
\noindent
that converges in the mean-square
sense   
is valid, where 
$i_1,\ldots,i_k=0,1,\ldots,m,$ 
$$
{\tilde \zeta}_{j}^{(i)}=
\int\limits_t^T \Psi_{j}(s)\sqrt{r(s)}d{\bf w}_s^{(i)}
$$
are independent standard Gaussian random variables
for various
$i$ or $j$ {\rm(}in the case when $i\ne 0${\rm),}
$$
{\tilde C}_{j_k\ldots j_1}=\int\limits_{[t,T]^k}
K(t_1,\ldots,t_k)
\prod_{l=1}^{k}\biggl(\Psi_{j_l}(t_l)r(t_l)\biggr)dt_1\ldots dt_k
$$

\noindent
is the Fourier coefficient$,$
$K(t_1,\ldots,t_k)$ is defined by {\rm (\ref{chain200});}
another notations are the same as in Theorems {\rm 1.1, 1.2, 1.5.}
}

{\bf Theorem 1.21.}\
{\it Under the conditions of Theorem~{\rm 1.20}
the following estimate
$$
{\sf M}\left\{\left(
{\tilde J}[\psi^{(k)}]_{T,t}-{\tilde J}[\psi^{(k)}]_{T,t}^{p_1,\ldots,p_k}
\right)^2\right\}
\le 
$$

\vspace{-2mm}
$$
~ \le k!\left(~\int\limits_{[t,T]^k}
K^2(t_1,\ldots,t_k)\left(\prod_{l=1}^k r(t_l)\right)
dt_1\ldots dt_k -\sum_{j_1=0}^{p_1}\ldots
\sum_{j_k=0}^{p_k}{\tilde C}^2_{j_k\ldots j_1}\right)
$$

\vspace{1mm}
\noindent
is valid for the following cases{\rm :}

{\rm 1.}\ $i_1,\ldots,i_k=1,\ldots,m$\ \ and\ \ $0<T-t<\infty,$

{\rm 2.}\ $i_1,\ldots,i_k=0, 1,\ldots,m,$\ \ $i_1^2+\ldots+i_k^2>0,$\ \
and\ \ $0<T-t<1,$
                             \

\vspace{1mm}
\noindent
where ${\tilde J}[\psi^{(k)}]_{T,t}$ is the 
stochastic integral {\rm (\ref{fifi1}),}
${\tilde J}[\psi^{(k)}]_{T,t}^{p_1,\ldots,p_k}$ is the 
expression on the right-hand side of {\rm (\ref{fifi2})} before
passing to the limit 
$\hbox{\vtop{\offinterlineskip\halign{
\hfil#\hfil\cr
{\rm l.i.m.}\cr
$\stackrel{}{{}_{p_1,\ldots,p_k\to \infty}}$\cr
}} };$ another 
notations are the same as in Theorems {\rm 1.6, 1.20}.
}

\section{Proof of Theorems~1.16 and 1.17 on the Base of the It\^{o} Formula and 
Without Explicit Use of the Multiple Wiener Stochastic Integral}

Note that Theorems~1.16 and 1.17 can also be proved without explicit 
use of the multiple Wiener stochastic integral.
To do this, we introduce the following sum 
of iterated It\^{o} stochastic integrals
\begin{equation}
\label{chain10100}
~~~~~~J''[\Phi]_{T,t}^{(i_1\ldots i_k)}\stackrel{\sf def}{=}\sum_{(t_1,\ldots,t_k)}
\int\limits_{t}^{T}
\ldots
\int\limits_{t}^{t_2}
\Phi(t_1,\ldots,t_k)d{\bf w}_{t_1}^{(i_1)}
\ldots
d{\bf w}_{t_k}^{(i_k)},
\end{equation}

\noindent
where $\Phi(t_1,\ldots,t_k)\in L_2([t, T]^k),$ $i_1,\ldots,i_k=0,1,\ldots,m,$\ $d{\bf w}_{\tau}^{(0)}
=d\tau;$
another notations are the same as in (\ref{Wi110}).

Further, using the isometry property of the Ito stochastic integral as well as
the linearity property of this integral, we have

\vspace{-2mm}
$$
J[\psi^{(k)}]_{T,t}^{(i_1\ldots i_k)}=J''[K]_{T,t}^{(i_1\ldots i_k)}=
$$

\vspace{-4mm}
\begin{equation}
\label{chain10200}
~~~~=\sum_{j_1=0}^{p_1}\ldots
\sum_{j_k=0}^{p_k}
C_{j_k\ldots j_1}
J''[\phi_{j_1}\ldots \phi_{j_k}]_{T,t}^{(i_1\ldots i_k)}+
J''[R_{p_1\ldots p_k}]_{T,t}^{(i_1\ldots i_k)}\ \ \ \hbox{w.~p.~1,}
\end{equation}

\noindent
where $K(t_1,\ldots,t_k)$ and $R_{p_1\ldots p_k}(t_1,\ldots,t_k)$
are defined by (\ref{chain200}) and (\ref{chain30001})
correspondingly. Moreover, 
$J''[\phi_{j_1}\ldots \phi_{j_k}]_{T,t}^{(i_1\ldots i_k)}$
and $J''[R_{p_1\ldots p_k}]_{T,t}^{(i_1\ldots i_k)}$
are defined by (\ref{chain10100}). Obviously, we can consider
an analogue of (\ref{chain10200}) for  $\Phi(t_1,\ldots,t_k)$
instead of $K(t_1,\ldots,t_k)$.

Passing to the limit $\hbox{\vtop{\offinterlineskip\halign{
\hfil#\hfil\cr
{\rm l.i.m.}\cr
$\stackrel{}{{}_{p_1,\ldots,p_k\to \infty}}$\cr
}} }$ in (\ref{chain10200}) and using (\ref{wi2005}), (\ref{chain7771}), (\ref{chain10100}), we obtain

\vspace{-5mm}
$$
J[\psi^{(k)}]_{T,t}^{(i_1\ldots i_k)}=
\hbox{\vtop{\offinterlineskip\halign{
\hfil#\hfil\cr
{\rm l.i.m.}\cr
$\stackrel{}{{}_{p_1,\ldots,p_k\to \infty}}$\cr
}} }
\sum\limits_{j_1=0}^{p_1}\ldots
\sum\limits_{j_k=0}^{p_k}
C_{j_k\ldots j_1}J''[\phi_{j_1}\ldots \phi_{j_k}]_{T,t}^{(i_1\ldots i_k)}=
$$

\vspace{-1mm}
\begin{equation}
\label{chain4002}
=\hbox{\vtop{\offinterlineskip\halign{
\hfil#\hfil\cr
{\rm l.i.m.}\cr
$\stackrel{}{{}_{p_1,\ldots,p_k\to \infty}}$\cr
}} }
\sum\limits_{j_1=0}^{p_1}\ldots
\sum\limits_{j_k=0}^{p_k}
C_{j_k\ldots j_1}
\sum\limits_{(t_1,\ldots,t_k)}
\int\limits_t^T \phi_{j_k}(t_k)
\ldots
\int\limits_t^{t_{2}}\phi_{j_{1}}(t_{1})
d{\bf w}_{t_1}^{(i_1)}\ldots
d{\bf w}_{t_k}^{(i_k)},
\end{equation}

\vspace{4mm}
\noindent
where permutations $(t_1,\ldots,t_k)$ when summing are 
performed only in the values
$d{\bf w}_{t_1}^{(i_1)}
\ldots $
$d{\bf w}_{t_k}^{(i_k)}.$ At the same time the indices near 
upper 
limits of integration in the iterated stochastic integrals are changed 
correspondently and if $t_r$ swapped with $t_q$ in the  
permutation $(t_1,\ldots,t_k)$, then $i_r$ swapped with $i_q$ in 
the permutation $(i_1,\ldots,i_k).$ 

It is easy to see that the equality (\ref{chain4002}) can be written as

$$
J[\psi^{(k)}]_{T,t}^{(i_1\ldots i_k)}=
\hbox{\vtop{\offinterlineskip\halign{
\hfil#\hfil\cr
{\rm l.i.m.}\cr
$\stackrel{}{{}_{p_1,\ldots,p_k\to \infty}}$\cr
}} }
\sum\limits_{j_1=0}^{p_1}\ldots
\sum\limits_{j_k=0}^{p_k}
C_{j_k\ldots j_1}\times
$$

\vspace{2mm}
\begin{equation}
\label{chain7878}
\times \sum\limits_{(j_1,\ldots,j_k)}
\int\limits_t^T \phi_{j_k}(t_k)
\ldots
\int\limits_t^{t_{2}}\phi_{j_{1}}(t_{1})
d{\bf w}_{t_1}^{(i_1)}\ldots
d{\bf w}_{t_k}^{(i_k)},
\end{equation}

\vspace{3mm}
\noindent
where 
$$
\sum\limits_{(j_1,\ldots,j_k)}
$$ 

\vspace{1mm}
\noindent
means the sum with respect to all
possible permutations 
$(j_1,\ldots,j_k).$ At the same time if 
$j_r$ swapped with $j_q$ in the permutation $(j_1,\ldots,j_k)$,
then $i_r$ swapped with $i_q$ in the permutation $(i_1,\ldots,i_k)$.

Further, using the It\^{o} formula, we can prove the following equality
$$
\sum\limits_{(j_1,\ldots,j_k)}
\int\limits_t^T \phi_{j_k}(t_k)
\ldots
\int\limits_t^{t_{2}}\phi_{j_{1}}(t_{1})
d{\bf w}_{t_1}^{(i_1)}\ldots
d{\bf w}_{t_k}^{(i_k)}
=\prod_{l=1}^k\zeta_{j_l}^{(i_l)}
+\sum\limits_{r=1}^{[k/2]}
(-1)^r \times
$$
\begin{equation}
\label{chain405xx}
~~\times\sum_{\stackrel{(\{\{g_1, g_2\}, \ldots, 
\{g_{2r-1}, g_{2r}\}\}, \{q_1, \ldots, q_{k-2r}\})}
{{}_{\{g_1, g_2, \ldots, 
g_{2r-1}, g_{2r}, q_1, \ldots, q_{k-2r}\}=\{1, 2, \ldots, k\}}}}
\prod\limits_{s=1}^r
{\bf 1}_{\{i_{g_{{}_{2s-1}}}=~i_{g_{{}_{2s}}}\ne 0\}}
\Biggl.{\bf 1}_{\{j_{g_{{}_{2s-1}}}=~j_{g_{{}_{2s}}}\}}
\prod_{l=1}^{k-2r}\zeta_{j_{q_l}}^{(i_{q_l})}
\end{equation}

\vspace{3mm}
\noindent
w.~p.~1, where notations are the same as in Theorem 1.2 and (\ref{chain7878}).

The main difficulty in proving (\ref{chain405xx}) 
using the It\^{o} formula is related to the need to 
take into account various combinations of indices $i_1,\ldots,i_k=0,1,\ldots,m.$
To avoid this difficulty, consider another approach, also based
on the It\^{o} formula.

First, we prove the following modification and generalization 
of Theorem~3.1 from \cite{ito1951} (1951) for the case $i_1,\ldots,i_k=0, 1,\ldots,m$
using the It\^{o} formula and without explicit use of the multiple Wiener 
stochastic integral.

{\bf Theorem~1.22}\ \cite{arxiv-1}.\ {\it Suppose that
the condition {\rm ($\star\star$)} is fulfilled
for the multi-index $(i_1 \ldots i_k)$ {\rm (}see Sect.~{\rm 1.10)} 
and the condition {\rm (\ref{ziko999})} is also 
fulfilled.
Furthermore$,$ let 
$\{\phi_j(x)\}_{j=0}^{\infty}$ is an arbitrary complete orthonormal system  
of functions in the space $L_2([t,T]).$
Then 

\vspace{-1mm}
$$
J''[\phi_{j_1}\ldots \phi_{j_k}]_{T,t}^{(i_1\ldots i_k)}=
$$

\vspace{-2mm}
\begin{equation}
\label{new6000}
=\prod_{l=1}^k\left({\bf 1}_{\{m_l=0\}}+{\bf 1}_{\{m_l>0\}}\left\{
\begin{matrix}
H_{n_{1,l}}\left(\zeta_{j_{h_{1,l}}}^{(i_l)}\right)\ldots 
H_{n_{d_l,l}}\left(\zeta_{j_{h_{d_l,l}}}^{(i_l)}\right),\ 
&\hbox{\rm if}\ \ \ 
i_l\ne 0\cr\cr
\left(\zeta_{j_{h_{1,l}}}^{(0)}\right)^{n_{1,l}}\ldots
\left(\zeta_{j_{h_{d_l,l}}}^{(0)}\right)^{n_{d_l,l}},\  &\hbox{\rm if}\ \ \ 
i_l=0
\end{matrix}\right.\ \right)
\end{equation}

\vspace{2mm}
\noindent
w.~p.~{\rm 1,} where $i_1,\ldots,i_k=0, 1,\ldots,m;$\ \
$n_{1,l}+n_{2,l}+\ldots+n_{d_l,l}=m_l;$\ \ $n_{1,l}, n_{2,l}, \ldots, n_{d_l,l}=1,\ldots, m_l;$\ \ 
$d_l=1,\ldots,m_l;$\ \ $l=1,\ldots,k;$\ \ $m_1+\ldots+m_k=k;$\ \  
the numbers $m_1,\ldots,m_k,$\ $g_1,\ldots,g_k$
depend on $(i_1,\ldots,i_k)$ and 
the numbers $n_{1,l},\ldots,n_{d_l,l},$\ $h_{1,l},\ldots,h_{d_l,l},$\ $d_l$
depend on $\{j_1,\ldots,j_k\};$ moreover$,$ $\left\{j_{g_1},\ldots,j_{g_k}\right\}
=\{j_1,\ldots,j_k\};$ $H_n(x)$ is the Hermite polynomial {\rm (\ref{ziko500});}
another
notations are the same as in Theorem {\rm 1.14}.}

{\bf Proof.}\ First, consider the case $i_1=\ldots=i_k=1,\ldots, m$ and $j_1,\ldots,j_k\in \{0\}\cup {\bf N}$.
By induction, we prove the following equality
$$
p! \int\limits_t^T \phi_l(t_p)\ldots \int\limits_t^{t_2}
\phi_l(t_1)d{\bf w}_{t_1}^{(1)}\ldots d{\bf w}_{t_p}^{(1)}\times
$$
$$
\times \sum\limits_{(j_1,\ldots,j_q)}
\int\limits_t^T \phi_{j_q}(t_q)\ldots \int\limits_t^{t_2}
\phi_{j_1}(t_1)d{\bf w}_{t_1}^{(1)}\ldots d{\bf w}_{t_q}^{(1)}=
$$
$$
=\sum\limits_{(j_1,\ldots,j_q, \underbrace{{}_{l, \ldots ,l}}_{p})}
\int\limits_t^T \phi_{j_q}(t_q)\ldots \int\limits_t^{t_2}
\phi_{j_1}(t_1)
\int\limits_t^{t_1} \phi_{l}(t_p')\ldots \int\limits_t^{t_2'}
\phi_{l}(t_1')\times
$$
\begin{equation}
\label{new1010}
\times
d{\bf w}_{t_1'}^{(1)}\ldots d{\bf w}_{t_p'}^{(1)}
d{\bf w}_{t_1}^{(1)}\ldots d{\bf w}_{t_q}^{(1)}
\end{equation}

\vspace{5mm}
\noindent
w.~p.~1, where $p\in{\bf N},$\ $l\ne j_1,\ldots,j_q,$ and
$$
\sum\limits_{(q_1,\ldots, q_n)}
$$

\noindent
means the sum with respect to all possible permutations
$(q_1,\ldots, q_n)$.

Consider the case $p=1.$ Using the It\^{o} formula, we get w.~p.~1 for
$s\in[t, T]$
$$
\int\limits_t^s \phi_l(\tau)
d{\bf w}_{\tau}^{(1)}
\int\limits_t^s \phi_{j_q}(t_q)\ldots \int\limits_t^{t_2}
\phi_{j_1}(t_1)d{\bf w}_{t_1}^{(1)}\ldots d{\bf w}_{t_q}^{(1)}=
$$
$$
=\int\limits_t^s \phi_l(\tau)\phi_{j_q}(\tau)
\int\limits_t^{\tau} \phi_{j_{q-1}}(t_{q-1})\ldots \int\limits_t^{t_2}
\phi_{j_1}(t_1)d{\bf w}_{t_1}^{(1)}\ldots d{\bf w}_{t_{q-1}}^{(1)}d\tau+
$$
$$
+\int\limits_t^s \phi_l(\tau)\
\int\limits_t^{\tau} \phi_{j_q}(t_q)\ldots \int\limits_t^{t_2}
\phi_{j_1}(t_1)d{\bf w}_{t_1}^{(1)}\ldots d{\bf w}_{t_q}^{(1)}d{\bf w}_{\tau}^{(1)}+
$$
\begin{equation}
\label{new1011}
+\int\limits_t^s \phi_{j_q}(\tau)\hspace{-0.5mm}
\left(\int\limits_t^{\tau} \phi_{l}(\theta)
d{\bf w}_{\theta}^{(1)}
\int\limits_t^{\tau} \phi_{j_{q-1}}(t_{q-1})\ldots \int\limits_t^{t_2}
\phi_{j_1}(t_1)d{\bf w}_{t_1}^{(1)}\ldots d{\bf w}_{t_{q-1}}^{(1)}\right)
\hspace{-0.5mm}d{\bf w}_{\tau}^{(1)}\hspace{-0.2mm}.
\end{equation}

Hereinafter in this section always $s\in [t, T].$
Differentiating by the It\^{o} formula the expression in parentheses
on the right-hand side of equality (\ref{new1011}) and combining the result 
of differentiation with (\ref{new1011}), we obtain w.~p.~1

\vspace{-4mm}
$$
J_{(l)s,t} J_{(j_q\ldots j_1)s,t}=
$$

\vspace{-4mm}
$$
=\int\limits_t^s \phi_l(\tau)\phi_{j_q}(\tau)
\int\limits_t^{\tau} \phi_{j_{q-1}}(t_{q-1})\ldots \int\limits_t^{t_2}
\phi_{j_1}(t_1)d{\bf w}_{t_1}^{(1)}\ldots d{\bf w}_{t_{q-1}}^{(1)}d\tau+
$$

\vspace{-3mm}
$$
+
J_{(l j_q\ldots j_1)s,t}+
$$

\vspace{-4mm}
$$
+
\int\limits_t^s \hspace{-0.2mm}\phi_{j_q}(\tau)
\hspace{-0.2mm}
\int\limits_t^{\tau}\hspace{-0.2mm}
\phi_{l}(\theta)\phi_{j_{q-1}}(\theta)\hspace{-0.1mm}\int\limits_t^{\theta}
\hspace{-0.2mm}\phi_{j_{q-2}}(t_{q-2})
\ldots \int\limits_t^{t_2}\hspace{-0.2mm}
\phi_{j_1}(t_1)d{\bf w}_{t_1}^{(1)}\ldots d{\bf w}_{t_{q-2}}^{(1)}d\theta
d{\bf w}_{\tau}^{(1)}+
$$

\vspace{-3mm}
$$
+
J_{(j_q l j_{q-1}\ldots j_1)s,t}+
$$

$$
+\int\limits_t^s \phi_{j_q}(\tau)
\int\limits_t^{\tau} \phi_{j_{q-1}}(\theta)\times
$$
\begin{equation}
\label{new1012}
\times\left(\int\limits_t^{\theta} \phi_l(u)\
d{\bf w}_{u}^{(1)}\int\limits_t^{\theta} 
\phi_{j_{q-2}}(t_{q-2})\ldots \int\limits_t^{t_2}
\phi_{j_1}(t_1)d{\bf w}_{t_1}^{(1)}\ldots d{\bf w}_{t_{q-2}}^{(1)}
\right)
d{\bf w}_{\theta}^{(1)}d{\bf w}_{\tau}^{(1)},
\end{equation}

\vspace{2mm}
\noindent
where
$$
\int\limits_t^s \phi_{j_q}(t_q)\ldots \int\limits_t^{t_2}
\phi_{j_1}(t_1)d{\bf w}_{t_1}^{(1)}\ldots d{\bf w}_{t_q}^{(1)}
\stackrel{\sf def}{=}J_{(j_q \ldots j_1)s,t}.
$$

\vspace{1mm}

Continuing the process of iterative application of the It\^{o} formula, we have w.~p.~1

\vspace{-5mm}
$$
J_{(l)s,t} J_{(j_q\ldots j_1)s,t}=
$$

\vspace{1mm}
$$
=J_{(l j_q\ldots j_1)s,t}+ J_{(j_q l j_{q-1}\ldots j_1)s,t}+\ldots + J_{(j_q\ldots j_1 l)s,t}+
$$

\vspace{-1mm}
$$
+\int\limits_t^s \phi_l(\tau)\phi_{j_q}(\tau)
\int\limits_t^{\tau} \phi_{j_{q-1}}(t_{q-1})\ldots \int\limits_t^{t_2}
\phi_{j_1}(t_1)d{\bf w}_{t_1}^{(1)}\ldots d{\bf w}_{t_{q-1}}^{(1)}d\tau+\ldots
$$
\begin{equation}
\label{new1040ss}
~~~~~\ldots +
\int\limits_t^{s} \phi_{j_{q}}(t_{q})\ldots \int\limits_t^{t_3}
\phi_{j_2}(t_2)\int\limits_t^{t_2} \phi_{l}(\tau)\phi_{j_1}(\tau)      
d\tau d{\bf w}_{t_2}^{(1)}\ldots d{\bf w}_{t_{q}}^{(1)}.
\end{equation}

\vspace{3mm}

Summing the equality (\ref{new1040ss}) over permutations $(j_1,\ldots, j_q)$, we get 

\vspace{-1mm}
\begin{equation}
\label{new1025}
\sum\limits_{(j_1,\ldots, j_q)}J_{(l)s,t} J_{(j_q\ldots j_1)s,t}=
\sum\limits_{(j_1,\ldots, j_q,l)}J_{(l j_q\ldots j_1)s,t}+ S(s)
\end{equation}

\vspace{3mm}
\noindent
w.~p.~1, where
$$
S(s)=
$$

\vspace{-6mm}
$$
=\sum\limits_{(j_1,\ldots, j_q)}\left(\int\limits_t^s \phi_l(\tau)\phi_{j_q}(\tau)
\int\limits_t^{\tau} \phi_{j_{q-1}}(t_{q-1})\ldots \int\limits_t^{t_2}
\phi_{j_1}(t_1)d{\bf w}_{t_1}^{(1)}\ldots d{\bf w}_{t_{q-1}}^{(1)}d\tau+\ldots\right.
$$

\vspace{-2mm}
\begin{equation}
\label{new1040}
~~~~~\left.\ldots +
\int\limits_t^{s} \phi_{j_{q}}(t_{q})\ldots \int\limits_t^{t_3}
\phi_{j_2}(t_2)\int\limits_t^{t_2} \phi_{l}(\tau)\phi_{j_1}(\tau)      
d\tau d{\bf w}_{t_2}^{(1)}\ldots d{\bf w}_{t_{q}}^{(1)}\right).
\end{equation}

\vspace{2mm}

Consider 
$$
\int\limits_t^s \phi_l(\tau)\phi_{j_q}(\tau)d\tau
\int\limits_t^{s} \phi_{j_{q-1}}(t_{q-1})\ldots \int\limits_t^{t_2}
\phi_{j_1}(t_1)d{\bf w}_{t_1}^{(1)}\ldots d{\bf w}_{t_{q-1}}^{(1)}.
$$

\vspace{2mm}

Applying the It\^{o} formula, we get w.~p.~1

\vspace{-2mm}
$$
\int\limits_t^s \phi_l(\tau)\phi_{j_q}(\tau)d\tau
\int\limits_t^{s} \phi_{j_{q-1}}(t_{q-1})\ldots \int\limits_t^{t_2}
\phi_{j_1}(t_1)d{\bf w}_{t_1}^{(1)}\ldots d{\bf w}_{t_{q-1}}^{(1)}=
$$

\vspace{-2mm}
$$
=\int\limits_t^s \phi_l(\tau)\phi_{j_q}(\tau)
\int\limits_t^{\tau} \phi_{j_{q-1}}(t_{q-1})\ldots \int\limits_t^{t_2}
\phi_{j_1}(t_1)d{\bf w}_{t_1}^{(1)}\ldots d{\bf w}_{t_{q-1}}^{(1)}d\tau+
$$

\vspace{-2mm}
$$
+\int\limits_t^s \phi_{j_{q-1}}(t_{q-1})\times
$$
$$
\times\left(\int\limits_t^{t_{q-1}}\phi_l(\tau)\phi_{j_q}(\tau)d\tau
\int\limits_t^{t_{q-1}} \phi_{j_{q-2}}(t_{q-2})\ldots \int\limits_t^{t_2}
\phi_{j_1}(t_1)d{\bf w}_{t_1}^{(1)}\ldots d{\bf w}_{t_{q-2}}^{(1)}\right)d{\bf w}_{t_{q-1}}^{(1)}.
$$

\vspace{4mm}

By iterative application of the It\^{o} formula (as above), we obtain w.~p.~1

\vspace{-2mm}
$$
\int\limits_t^s \phi_l(\tau)\phi_{j_q}(\tau)d\tau
\int\limits_t^{s} \phi_{j_{q-1}}(t_{q-1})\ldots \int\limits_t^{t_2}
\phi_{j_1}(t_1)d{\bf w}_{t_1}^{(1)}\ldots d{\bf w}_{t_{q-1}}^{(1)}=
$$

\vspace{-2mm}
$$
=\int\limits_t^s \phi_l(\tau)\phi_{j_q}(\tau)
\int\limits_t^{\tau} \phi_{j_{q-1}}(t_{q-1})\ldots \int\limits_t^{t_2}
\phi_{j_1}(t_1)d{\bf w}_{t_1}^{(1)}\ldots d{\bf w}_{t_{q-1}}^{(1)}d\tau+\ldots
$$

\vspace{-2mm}
\begin{equation}
\label{newx1020ss}
\ldots +
\int\limits_t^{s} \phi_{j_{q-1}}(t_{q-1})\ldots \int\limits_t^{t_2}
\phi_{j_1}(t_1)\int\limits_t^{t_1} \phi_{l}(\tau)\phi_{j_q}(\tau)      
d\tau d{\bf w}_{t_1}^{(1)}\ldots d{\bf w}_{t_{q-1}}^{(1)}.
\end{equation}

\vspace{3mm}

Summing the equality (\ref{newx1020ss}) over permutations $(j_1,\ldots, j_q)$, we get 

\vspace{-4mm}
\begin{equation}
\label{new1020}
\sum\limits_{(j_1,\ldots, j_q)}\int\limits_t^s \phi_l(\tau)\phi_{j_q}(\tau)d\tau
\int\limits_t^{s} \phi_{j_{q-1}}(t_{q-1})\ldots \int\limits_t^{t_2}
\phi_{j_1}(t_1)d{\bf w}_{t_1}^{(1)}\ldots d{\bf w}_{t_{q-1}}^{(1)}=S_1(s),
\end{equation}

\vspace{2mm}
\noindent
w.~p.~1, where
$$
S_1(s)=
$$

\vspace{-4mm}
$$
=\sum\limits_{(j_1,\ldots, j_q)}\left(\int\limits_t^s \phi_l(\tau)\phi_{j_q}(\tau)
\int\limits_t^{\tau} \phi_{j_{q-1}}(t_{q-1})\ldots \int\limits_t^{t_2}
\phi_{j_1}(t_1)d{\bf w}_{t_1}^{(1)}\ldots d{\bf w}_{t_{q-1}}^{(1)}d\tau+\ldots\right.
$$

\vspace{-2mm}
\begin{equation}
\label{newx1020}
~~~~~\left.\ldots +
\int\limits_t^{s} \phi_{j_{q-1}}(t_{q-1})\ldots \int\limits_t^{t_2}
\phi_{j_1}(t_1)\int\limits_t^{t_1} \phi_{l}(\tau)\phi_{j_q}(\tau)      
d\tau d{\bf w}_{t_1}^{(1)}\ldots d{\bf w}_{t_{q-1}}^{(1)}\right).
\end{equation}

\vspace{3mm}

It is not difficult to see that
\begin{equation}
\label{new1021}
S(s)=S_1(s)\ \ \ \hbox{w.~p.~1.}
\end{equation}

\vspace{2mm}

Moreover, due to the orthogonality of $\{\phi_j(x)\}_{j=0}^{\infty}$
and (\ref{new1020}), (\ref{new1021}), we have
\begin{equation}
\label{new1026}
S(T)=S_1(T)=0\ \ \ \hbox{w.~p.~1.}
\end{equation}

\vspace{3mm}

Thus (see (\ref{new1025}), (\ref{new1026})), the equality (\ref{new1010}) is proved for the case
$p=1.$
Let us assume that the equality (\ref{new1010}) is true for $p=2, 3, \ldots, k-1$, and prove
its validity for $p=k.$

From (\ref{new1025}) for the case $q=k-1,$ $j_1=\ldots=j_{k-1}=l$ we obtain

\vspace{-1mm}
\begin{equation}
\label{new1042}
\left(J_1\right)_{s,t} (k-1)! \left(J_{k-1}\right)_{s,t}=k! \left(J_{k}\right)_{s,t} + S_2(s)
\end{equation}

\vspace{3mm}
\noindent
w.~p.~1, where 

\vspace{-3mm}
$$
S_2(s)=S(s)\biggl|_{j_1=\ldots=j_q=l,\ q=k-1}\biggr.\ \ (k\ge 2)\ \ \ \hbox{and}\ \ \
S_2(s)\stackrel{\sf def}{=}0\ \ (q=k-1,\ k=1),
$$

\vspace{-2mm}
$$
\int\limits_t^s \phi_{l}(t_r)\ldots \int\limits_t^{t_2}
\phi_{l}(t_1)d{\bf w}_{t_1}^{(1)}\ldots d{\bf w}_{t_r}^{(1)}
\stackrel{\sf def}{=}\left(J_{r}\right)_{s,t}\ \ (r\in {\bf N})\ \ \ \hbox{and}\ \ \
\left(J_{0}\right)_{s,t}\stackrel{\sf def}{=}1.
$$

\vspace{3mm}

Taking into account (\ref{new1040}), (\ref{new1020})--(\ref{new1021})
and the orthonormality of $\{\phi_j(x)\}_{j=0}^{\infty}$, we have
\begin{equation}
\label{new1043}
S_2(T)=(k-1)!\left(J_{k-2}\right)_{T,t}.
\end{equation}

\vspace{3mm}

Combining (\ref{new1042}) and (\ref{new1043}), we obtain the following recurrence 
relation

\vspace{-4mm}
\begin{equation}
\label{new1044}
~~~~~~~~~k! \left(J_{k}\right)_{T,t}=\left(J_1\right)_{T,t} (k-1)! \left(J_{k-1}\right)_{T,t}-
(k-1)!\left(J_{k-2}\right)_{T,t}
\end{equation}

\noindent
w.~p.~1.

Using (\ref{new1044}) and the induction hypothesis, we get w.~p.~1
$$
k! \int\limits_t^T \phi_l(t_k)\ldots \int\limits_t^{t_2}
\phi_l(t_1)d{\bf w}_{t_1}^{(1)}\ldots d{\bf w}_{t_k}^{(1)}\times
$$
$$
\times \sum\limits_{(j_1,\ldots,j_q)}
\int\limits_t^T \phi_{j_q}(t_q)\ldots \int\limits_t^{t_2}
\phi_{j_1}(t_1)d{\bf w}_{t_1}^{(1)}\ldots d{\bf w}_{t_q}^{(1)}=
$$
$$
=\int\limits_t^T \phi_l(\tau)\
d{\bf w}_{\tau}^{(1)}\Biggl(
(k-1)!\int\limits_t^T \phi_l(t_{k-1})\ldots \int\limits_t^{t_2}
\phi_l(t_1)d{\bf w}_{t_1}^{(1)}\ldots d{\bf w}_{t_{k-1}}^{(1)}\times\Biggr.
$$
$$
\Biggl.\times \sum\limits_{(j_1,\ldots,j_q)}
\int\limits_t^T \phi_{j_q}(t_q)\ldots \int\limits_t^{t_2}
\phi_{j_1}(t_1)d{\bf w}_{t_1}^{(1)}\ldots d{\bf w}_{t_q}^{(1)}\Biggr)-
$$
$$
-
(k-1)!\int\limits_t^T \phi_l(t_{k-2})\ldots \int\limits_t^{t_2}
\phi_l(t_1)d{\bf w}_{t_1}^{(1)}\ldots d{\bf w}_{t_{k-2}}^{(1)}\times
$$
$$
\times \sum\limits_{(j_1,\ldots,j_q)}
\int\limits_t^T \phi_{j_q}(t_q)\ldots \int\limits_t^{t_2}
\phi_{j_1}(t_1)d{\bf w}_{t_1}^{(1)}\ldots d{\bf w}_{t_q}^{(1)}=
$$
$$
=\int\limits_t^T \phi_l(\tau)\
d{\bf w}_{\tau}^{(1)}
\sum\limits_{(j_1,\ldots,j_q, \underbrace{{}_{l, \ldots ,l}}_{k-1})}
\int\limits_t^T \phi_{j_q}(t_q)\ldots \int\limits_t^{t_2}
\phi_{j_1}(t_1)
\int\limits_t^{t_1} \phi_{l}(t_{k-1}')\ldots \int\limits_t^{t_2'}
\phi_{l}(t_1')\times
$$
$$
\times
d{\bf w}_{t_1'}^{(1)}\ldots d{\bf w}_{t_{k-1}'}^{(1)}
d{\bf w}_{t_1}^{(1)}\ldots d{\bf w}_{t_{q}}^{(1)}-
$$

\vspace{-3mm}
$$
-(k-1)\sum\limits_{(j_1,\ldots,j_q, \underbrace{{}_{l, \ldots ,l}}_{k-2})}
\int\limits_t^T \phi_{j_q}(t_q)\ldots \int\limits_t^{t_2}
\phi_{j_1}(t_1)
\int\limits_t^{t_1} \phi_{l}(t_{k-2}')\ldots \int\limits_t^{t_2'}
\phi_{l}(t_1')\times
$$
\begin{equation}
\label{new1050a}
\times
d{\bf w}_{t_1'}^{(1)}\ldots d{\bf w}_{t_{k-2}'}^{(1)}
d{\bf w}_{t_1}^{(1)}\ldots d{\bf w}_{t_{q}}^{(1)}.
\end{equation}

\vspace{4mm}

Let $\fbox{\it l}$ be the symbol $l$ which does not participate
in the following sum with respect to permutations
$$
\sum\limits_{(j_1,\ldots,j_q, \underbrace{{}_{l, \ldots ,l}}_{k-1})}.
$$

Using (\ref{new1025}), we have w.~p.~1
$$
\int\limits_t^s \phi_{l}(\tau)\
d{\bf w}_{\tau}^{(1)}
\sum\limits_{(j_1,\ldots,j_q, \underbrace{{}_{l, \ldots ,l}}_{k-1})}
\int\limits_t^s \phi_{j_q}(t_q)\ldots \int\limits_t^{t_2}
\phi_{j_1}(t_1)
\int\limits_t^{t_1} \phi_{l}(t_{k-1}')\ldots \int\limits_t^{t_2'}
\phi_{l}(t_1')\times
$$
$$
\times
d{\bf w}_{t_1'}^{(1)}\ldots d{\bf w}_{t_{k-1}'}^{(1)}
d{\bf w}_{t_1}^{(1)}\ldots d{\bf w}_{t_{q}}^{(1)}=
$$

\vspace{-4mm}
$$
=
\sum\limits_{(j_1,\ldots,j_q, \underbrace{{}_{l, \ldots ,l}}_{k-1})}
\int\limits_t^s \phi_{\small{\fbox{\it l}}}(\tau)\
d{\bf w}_{\tau}^{(1)}\int\limits_t^s \phi_{j_q}(t_q)\ldots \int\limits_t^{t_2}
\phi_{j_1}(t_1)
\int\limits_t^{t_1} \phi_{l}(t_{k-1}')\ldots \int\limits_t^{t_2'}
\phi_{l}(t_1')\times
$$
$$
\times
d{\bf w}_{t_1'}^{(1)}\ldots d{\bf w}_{t_{k-1}'}^{(1)}
d{\bf w}_{t_1}^{(1)}\ldots d{\bf w}_{t_{q}}^{(1)}=
$$

$$
=\sum\limits_{(j_1,\ldots,j_q, \underbrace{{}_{l, \ldots ,l}}_{k-1})}
\left(J_{(\small{\fbox{\it l}}j_q\ldots j_1 \underbrace{l \ldots l}_{k-1})s,t}+
J_{(\small{j_q\fbox{\it l}}j_{q-1}\ldots j_1 \underbrace{l \ldots l}_{k-1})s,t}+\ldots\right.
$$
$$
\left.\ldots +J_{(j_q\ldots j_1 \small{\fbox{\it l}}\underbrace{l \ldots l}_{k-1})s,t}+
J_{(j_q\ldots j_1 l\small{\fbox{\it l}}\underbrace{l \ldots l}_{k-2})s,t}+
\ldots + J_{(j_q\ldots j_1 \small{\underbrace{l \ldots l}_{k-1}}\small{\fbox{\it l}})s,t}\right)+
S_3(s)=
$$

\vspace{4mm}
\begin{equation}
\label{new1050b}
=\sum\limits_{(j_1,\ldots,j_q, \underbrace{{}_{l, \ldots ,l}}_{k})}
J_{(j_q\ldots j_1 \small{\underbrace{l \ldots l}_{k}})s,t}
+S_3(s),
\end{equation}

\vspace{4mm}
\noindent
where
$$
S_3(s)=
$$

\vspace{-2mm}
$$
=
\sum\limits_{(j_1,\ldots,j_q, \underbrace{{}_{l, \ldots ,l}}_{k-1})}
\Biggl(\int\limits_t^s \phi_{\small{\fbox{\it l}}}(\tau)\phi_{j_q}(\tau)
\int\limits_t^{\tau} \phi_{j_{q-1}}(t_{q-1})\ldots \int\limits_t^{t_2}
\phi_{j_1}(t_1)
\times\Biggr.
$$
$$
\times\int\limits_t^{t_1} \phi_{l}(t_{k-1}')\ldots \int\limits_t^{t_2'}
\phi_{l}(t_1')
d{\bf w}_{t_1'}^{(1)}\ldots d{\bf w}_{t_{k-1}'}^{(1)}
d{\bf w}_{t_1}^{(1)}\ldots d{\bf w}_{t_{q-1}}^{(1)}d\tau +  \ldots
$$
$$
+\ldots 
\int\limits_t^{s} \phi_{j_{q}}(t_{q})\ldots \int\limits_t^{t_3}
\phi_{j_2}(t_2)\int\limits_t^{t_2} \phi_{\small{\fbox{\it l}}}(\tau)\phi_{j_1}(\tau)
\times\Biggr.
$$
$$
\times\int\limits_t^{\tau} \phi_{l}(t_{k-1}')\ldots \int\limits_t^{t_2'}
\phi_{l}(t_1')
d{\bf w}_{t_1'}^{(1)}\ldots d{\bf w}_{t_{k-1}'}^{(1)}
d\tau d{\bf w}_{t_2}^{(1)}\ldots d{\bf w}_{t_{q}}^{(1)}+
$$
$$
+
\int\limits_t^{s} \phi_{j_{q}}(t_{q})\ldots \int\limits_t^{t_2}
\phi_{j_1}(t_1)\int\limits_t^{t_1} \phi_{\small{\fbox{\it l}}}(\tau)\phi_{l}(\tau)
\times\Biggr.
$$
$$
\times\int\limits_t^{\tau} \phi_{l}(t_{k-2}')\ldots \int\limits_t^{t_2'}
\phi_{l}(t_1')
d{\bf w}_{t_1'}^{(1)}\ldots d{\bf w}_{t_{k-2}'}^{(1)}d\tau
d{\bf w}_{t_1}^{(1)}\ldots d{\bf w}_{t_{q}}^{(1)}+ \ldots
$$
$$
\ldots +
\int\limits_t^{s} \phi_{j_{q}}(t_{q})\ldots \int\limits_t^{t_2}
\phi_{j_1}(t_1)
\times\Biggr.
$$
$$
\Biggl.\times\int\limits_t^{t_1} \phi_{l}(t_{k-1}')\ldots \int\limits_t^{t_3'}
\phi_{l}(t_2')\int\limits_t^{t_2'} \phi_{\small{\fbox{\it l}}}(\tau)\phi_{l}(\tau)
d\tau d{\bf w}_{t_2'}^{(1)}\ldots d{\bf w}_{t_{k-1}'}^{(1)}
d{\bf w}_{t_1}^{(1)}\ldots d{\bf w}_{t_{q}}^{(1)}\Biggr).
$$

\vspace{4mm}

Using (\ref{new1040}), (\ref{new1020})--(\ref{new1021}), we get w.~p.~1

$$
S_3(s)=
$$

\vspace{-2mm}
$$
=\sum\limits_{(j_1,\ldots,j_q, \underbrace{{}_{l, \ldots ,l}}_{k-1})}
\int\limits_t^{s} \phi_{\small{\fbox{\it l}}}(\tau)\phi_{l}(\tau)d\tau
\int\limits_t^{s} \phi_{j_{q}}(t_{q})\ldots \int\limits_t^{t_2}
\phi_{j_1}(t_1)\times
$$
$$
\times
\int\limits_t^{t_1}\phi_{l}(t_{k-2}')\ldots \int\limits_t^{t_2'}
\phi_{l}(t_1')
d{\bf w}_{t_1'}^{(1)}\ldots d{\bf w}_{t_{k-2}'}^{(1)}
d{\bf w}_{t_1}^{(1)}\ldots d{\bf w}_{t_{q}}^{(1)} =
$$

$$
=(k-1)\sum\limits_{(j_1,\ldots,j_q, \underbrace{{}_{l, \ldots ,l}}_{k-2})}
\int\limits_t^{s} \phi_{\small{\fbox{\it l}}}(\tau)\phi_{l}(\tau)d\tau
\int\limits_t^{s} \phi_{j_{q}}(t_{q})\ldots \int\limits_t^{t_2}
\phi_{j_1}(t_1)\times
$$
$$
\times
\int\limits_t^{t_1}\phi_{l}(t_{k-2}')\ldots \int\limits_t^{t_2'}
\phi_{l}(t_1')
d{\bf w}_{t_1'}^{(1)}\ldots d{\bf w}_{t_{k-2}'}^{(1)}
d{\bf w}_{t_1}^{(1)}\ldots d{\bf w}_{t_{q}}^{(1)} +
$$

$$
+\sum\limits_{(j_1,\ldots,j_{q-1}, \underbrace{{}_{l, \ldots ,l}}_{k-1})}
\int\limits_t^{s} \phi_{\small{\fbox{\it l}}}(\tau)\phi_{j_q}(\tau)d\tau
\int\limits_t^{s} \phi_{j_{q-1}}(t_{q-1})\ldots \int\limits_t^{t_2}
\phi_{j_1}(t_1)\times
$$
$$
\times
\int\limits_t^{t_1}\phi_{l}(t_{k-1}')\ldots \int\limits_t^{t_2'}
\phi_{l}(t_1')
d{\bf w}_{t_1'}^{(1)}\ldots d{\bf w}_{t_{k-1}'}^{(1)}
d{\bf w}_{t_1}^{(1)}\ldots d{\bf w}_{t_{q-1}}^{(1)} + 
$$
$$
+\sum\limits_{(j_1,\ldots,j_{q-2}, j_q \underbrace{{}_{l, \ldots ,l}}_{k-1})}
\int\limits_t^{s} \phi_{\small{\fbox{\it l}}}(\tau)\phi_{j_{q-1}}(\tau)d\tau
\int\limits_t^{s} \phi_{j_{q}}(t_{q}) \int\limits_t^{t_q} \phi_{j_{q-2}}(t_{q-2})\ldots \int\limits_t^{t_2}
\phi_{j_1}(t_1)\times
$$

$$
\times
\int\limits_t^{t_1}\phi_{l}(t_{k-1}')\ldots \int\limits_t^{t_2'}
\phi_{l}(t_1')
d{\bf w}_{t_1'}^{(1)}\ldots d{\bf w}_{t_{k-1}'}^{(1)}
d{\bf w}_{t_1}^{(1)}\ldots d{\bf w}_{t_{q-2}}^{(1)}d{\bf w}_{t_{q}}^{(1)} + 
$$

\vspace{-4mm}
$$
\ldots
$$

\vspace{-4mm}
$$
+\sum\limits_{(j_2,\ldots,j_q \underbrace{{}_{l, \ldots ,l}}_{k-1})}
\int\limits_t^{s} \phi_{\small{\fbox{\it l}}}(\tau)\phi_{j_{1}}(\tau)d\tau
\int\limits_t^{s} \phi_{j_{q}}(t_{q}) \ldots \int\limits_t^{t_3}
\phi_{j_2}(t_2)\times
$$
\begin{equation}
\label{new1047}
~~~~~~~~\times
\int\limits_t^{t_2}\phi_{l}(t_{k-1}')\ldots \int\limits_t^{t_2'}
\phi_{l}(t_1')
d{\bf w}_{t_1'}^{(1)}\ldots d{\bf w}_{t_{k-1}'}^{(1)}
d{\bf w}_{t_2}^{(1)}\ldots d{\bf w}_{t_{q}}^{(1)}.
\end{equation}

\vspace{4mm}

Applying (\ref{new1047}) and the orthonormality of $\{\phi_j(x)\}_{j=0}^{\infty}$, we finally have
$$
S_3(T)=(k-1)\sum\limits_{(j_1,\ldots,j_q, \underbrace{{}_{l, \ldots ,l}}_{k-2})}
\int\limits_t^{T} \phi_{j_{q}}(t_{q})\ldots \int\limits_t^{t_2}
\phi_{j_1}(t_1)
\times
$$
\begin{equation}
\label{new1048}
~~~~~~~\times
\int\limits_t^{t_1}\phi_{l}(t_{k-2}')\ldots \int\limits_t^{t_2'}
\phi_{l}(t_1')
d{\bf w}_{t_1'}^{(1)}\ldots d{\bf w}_{t_{k-2}'}^{(1)}
d{\bf w}_{t_1}^{(1)}\ldots d{\bf w}_{t_{q}}^{(1)}.
\end{equation}

\vspace{3mm}

Combining (\ref{new1050a}), (\ref{new1050b}), (\ref{new1048}), we obtain w.~p.~1
$$
k! \int\limits_t^T \phi_l(t_k)\ldots \int\limits_t^{t_2}
\phi_l(t_1)d{\bf w}_{t_1}^{(1)}\ldots d{\bf w}_{t_k}^{(1)}\times
$$
$$
\times \sum\limits_{(j_1,\ldots,j_q)}
\int\limits_t^T \phi_{j_q}(t_q)\ldots \int\limits_t^{t_2}
\phi_{j_1}(t_1)d{\bf w}_{t_1}^{(1)}\ldots d{\bf w}_{t_q}^{(1)}=
$$
$$
=\sum\limits_{(\underbrace{{}_{l, \ldots ,l}}_{k})}
\int\limits_t^T \phi_l(t_k)\ldots \int\limits_t^{t_2}
\phi_l(t_1)d{\bf w}_{t_1}^{(1)}\ldots d{\bf w}_{t_k}^{(1)}\times
$$
$$
\times \sum\limits_{(j_1,\ldots,j_q)}
\int\limits_t^T \phi_{j_q}(t_q)\ldots \int\limits_t^{t_2}
\phi_{j_1}(t_1)d{\bf w}_{t_1}^{(1)}\ldots d{\bf w}_{t_q}^{(1)}=
$$
$$
=\sum\limits_{(j_1,\ldots,j_q, \underbrace{{}_{l, \ldots ,l}}_{k})}
\int\limits_t^T \phi_{j_q}(t_q)\ldots \int\limits_t^{t_2}
\phi_{j_1}(t_1)
\int\limits_t^{t_1} \phi_{l}(t_k')\ldots \int\limits_t^{t_2'}
\phi_{l}(t_1')\times
$$
\begin{equation}
\label{new1060}
\times
d{\bf w}_{t_1'}^{(1)}\ldots d{\bf w}_{t_k'}^{(1)}
d{\bf w}_{t_1}^{(1)}\ldots d{\bf w}_{t_q}^{(1)},
\end{equation}

\vspace{3mm}
\noindent
where $l\ne j_1,\ldots, j_q.$

The equality (\ref{new1010}) is proved. From the other hand, (\ref{new1060}) means that
\begin{equation}
\label{new1061}
J''[\phi_{j_1}\ldots \phi_{j_q}\underbrace{\phi_{l}\ldots \phi_{l}}_{n}]^
{(\hspace{0.5mm}\small{\overbrace{1 \ldots 1}^{q+n}}\hspace{0.5mm})}_{T,t}=
J''[\underbrace{\phi_{l}\ldots \phi_{l}}_{n}]^
{(\hspace{0.5mm}\small{\overbrace{1 \ldots 1}^{n}}\hspace{0.5mm})}_{T,t}
\cdot J''[\phi_{j_1}\ldots \phi_{j_q}]^
{(\hspace{0.5mm}\small{\overbrace{1 \ldots 1}^{q}}\hspace{0.5mm})}_{T,t}
\end{equation}
w.~p.~1, where $n, q=0,1,2\ldots;$\ $l\ne j_1,\ldots, j_q$ and
$$
J''[\phi_{j_1}\ldots \phi_{j_q}]^
{(\hspace{0.5mm}\small{\overbrace{1 \ldots 1}^{q}}\hspace{0.5mm})}_{T,t}\stackrel{\sf def}{=}1
$$
for $q=0.$

Note that \cite{Ch} (see Chapter 6, Sect.~6.6 of this book for details)
$$
\int\limits_t^T \phi_l(t_n)\ldots \int\limits_t^{t_2}
\phi_l(t_1)d{\bf w}_{t_1}^{(1)}\ldots d{\bf w}_{t_n}^{(1)}=
$$
$$
=
\frac{1}{n!} H_n\left(\int\limits_t^T \phi_l(\tau)d{\bf w}_{\tau}^{(1)},
\int\limits_t^T \phi_l^2(\tau)d\tau\right)=
$$
\begin{equation}
\label{new1100}
~~~~~~~=
\frac{1}{n!} H_n\left(\int\limits_t^T \phi_l(\tau)d{\bf w}_{\tau}^{(1)},1\right)=
\frac{1}{n!} H_n\left(\int\limits_t^T \phi_l(\tau)d{\bf w}_{\tau}^{(1)}\right)
\end{equation}

\vspace{2mm}
\noindent
w.~p.~1, where $n\in {\bf N},$ $H_n(x, y)$ is defined by (\ref{new1090})
(also see (\ref{ziko1000})), and 
$H_n(x)$ is the Hermite polynomial (\ref{ziko500}).

From (\ref{new1100}) we have w.~p.~1
$$
J''[\underbrace{\phi_{l}\ldots \phi_{l}}_{n}]^
{(\hspace{0.5mm}\small{\overbrace{1 \ldots 1}^{n}}\hspace{0.5mm})}_{T,t}
=n! \int\limits_t^T \phi_l(t_n)\ldots \int\limits_t^{t_2}
\phi_l(t_1)d{\bf w}_{t_1}^{(1)}\ldots d{\bf w}_{t_n}^{(1)}=
$$
\begin{equation}
\label{new1101}
~~~~~~=n! \frac{1}{n!} H_n\left(\int\limits_t^T \phi_l(\tau)d{\bf w}_{\tau}^{(1)}\right)=
H_n\left(\int\limits_t^T \phi_l(\tau)d{\bf w}_{\tau}^{(1)}\right),
\end{equation}

\vspace{2mm}
\noindent
where $n\in{\bf N}.$

Combining (\ref{new1061}) and (\ref{new1101}), we obtain
\begin{equation}
\label{new1102}
J''[\phi_{j_1}\ldots \phi_{j_q}\underbrace{\phi_{l}\ldots \phi_{l}}_{n}]^
{(\hspace{0.5mm}\small{\overbrace{1 \ldots 1}^{q+n}}\hspace{0.5mm})}_{T,t}=
H_n\left(\int\limits_t^T \phi_l(\tau)d{\bf w}_{\tau}^{(1)}\right)
\cdot J''[\phi_{j_1}\ldots \phi_{j_q}]^
{(\hspace{0.5mm}\small{\overbrace{1 \ldots 1}^{q}}\hspace{0.5mm})}_{T,t}
\end{equation}
w.~p.~1, where $n, q=0,1,2\ldots;$\ $l\ne j_1,\ldots, j_q.$ 

The iterated application of the formula (\ref{new1102})
completes the proof of Theorem~1.22 for the case 
$i_1=\ldots=i_k=1,\ldots, m$ and $j_1,\ldots,j_k\in \{0\}\cup {\bf N}$.

To prove Theorem~1.22 for the case $i_1=\ldots=i_k=0, 1,\ldots, m$ and 
$j_1,\ldots,j_k\in \{0\}\cup {\bf N}$, we need to prove
the following formula in addition to the previous proof
$$
p! \int\limits_t^T \phi_l(t_p)\ldots \int\limits_t^{t_2}
\phi_l(t_1)dt_1\ldots dt_p
\sum\limits_{(j_1,\ldots,j_q)}
\int\limits_t^T \phi_{j_q}(t_q)\ldots \int\limits_t^{t_2}
\phi_{j_1}(t_1)dt_1\ldots dt_q=
$$
\begin{equation}
\label{new1200}
=\sum\limits_{(j_1,\ldots,j_q, \underbrace{{}_{l, \ldots ,l}}_{p})}
\int\limits_t^T \phi_{j_q}(t_q)\ldots \int\limits_t^{t_2}
\phi_{j_1}(t_1)
\int\limits_t^{t_1} \phi_{l}(t_p')\ldots \int\limits_t^{t_2'}
\phi_{l}(t_1')
dt_1'\ldots dt_p'
dt_1\ldots dt_q,
\end{equation}

\vspace{3mm}
\noindent
where $p\in{\bf N}$,
$$
\sum\limits_{(j_1,\ldots, j_{d})}
$$

\noindent
means the sum with respect to all possible permutations $(j_1,\ldots,j_{d}).$ 

First, consider the case $p=1.$ We have 
$$
d\left(\int\limits_t^s \phi_l(\theta)d\theta
\int\limits_t^s \phi_{j_q}(t_q)\ldots \int\limits_t^{t_2}
\phi_{j_1}(t_1)dt_1\ldots dt_q\right)=
$$
$$
=\phi_l(s)
\int\limits_t^s \phi_{j_q}(t_q)\ldots \int\limits_t^{t_2}
\phi_{j_1}(t_1)dt_1\ldots dt_q ds+
$$
$$
+\phi_{j_q}(s)\left(\int\limits_t^{s}\phi_{j_{q-1}}(t_{q-1})\ldots \int\limits_t^{t_2}
\phi_{j_1}(t_1)dt_{1}\ldots dt_{q-1}
\cdot
\int\limits_t^s \phi_l(\theta)d\theta\right)ds.
$$

\vspace{3mm}

Then
$$
\int\limits_t^s \phi_l(\theta)d\theta
\int\limits_t^s \phi_{j_q}(t_q)\ldots \int\limits_t^{t_2}
\phi_{j_1}(t_1)dt_1\ldots dt_q=
$$
$$
=I_{(l j_q \ldots j_1)s,t}+
$$
$$
+\int\limits_t^s
\phi_{j_q}(\tau)\left(\int\limits_t^{\tau}\phi_{j_{q-1}}(t_{q-1})\ldots \int\limits_t^{t_2}
\phi_{j_1}(t_1)dt_{1}\ldots dt_{q-1}
\cdot
\int\limits_t^{\tau} \phi_l(\theta)d\theta\right)d\tau,
$$

\vspace{2mm}
\noindent
where
\begin{equation}
\label{new1701}
~~~~~~~~\int\limits_t^s \phi_{j_r}(t_r)\ldots \int\limits_t^{t_2}
\phi_{j_1}(t_1)dt_1\ldots dt_r
\stackrel{\sf def}{=}I_{(j_r \ldots j_1)s,t}.
\end{equation}

\vspace{2mm}

Continuing this process, we get
\begin{equation}
\label{new1301}
\int\limits_t^s \phi_l(\theta)d\theta\sum\limits_{(j_1,\ldots,j_q)}
I_{(j_q \ldots j_1)s,t}=
\sum\limits_{(j_1,\ldots,j_q, l)}
I_{(l j_q \ldots j_1)s,t},
\end{equation}

\noindent
where
$$
\sum\limits_{(j_1,\ldots,j_d)}
$$

\vspace{1mm}
\noindent
means the sum with respect to all possible permutations $(j_1,\ldots,j_d)$.

The equality (\ref{new1200}) is proved for the case $p=1.$
Let us assume that the equality (\ref{new1200}) is true for $p=2, 3, \ldots, k-1$, and prove
its validity for $p=k.$

From (\ref{new1301}) for $j_1=\ldots=j_q=l,$ $q=k-1$ we have

\vspace{-1mm}
\begin{equation}
\label{new1400}
\left(I_1\right)_{s,t} (k-1)! \left(I_{k-1}\right)_{s,t}=k!\left(I_{k}\right)_{s,t},
\end{equation}

\vspace{1mm}
\noindent
where $k\in {\bf N}$ and
$$
\int\limits_t^s \phi_{l}(t_k)\ldots \int\limits_t^{t_2}
\phi_{l}(t_1)dt_1\ldots dt_k
\stackrel{\sf def}{=}\left(I_k\right)_{s,t},\ \ \ \left(I_0\right)_{s,t}\stackrel{\sf def}{=}1.
$$

\vspace{2mm}

Using (\ref{new1400}) and the induction hypothesis, we obtain 

\vspace{-2mm}
$$
k! \left(I_k\right)_{s,t}
\sum\limits_{(j_1,\ldots,j_q)}
I_{(j_q\ldots j_1)s,t}=
\left(I_1\right)_{s,t} (k-1)! \left(I_{k-1}\right)_{s,t}
\sum\limits_{(j_1,\ldots,j_q)}I_{(j_q\ldots j_1)s,t}=
$$

\vspace{-2mm}
\begin{equation}
\label{new1499}
=I_{(l)s,t} 
\sum\limits_{(j_1,\ldots,j_q, \underbrace{{}_{l, \ldots ,l}}_{k-1})}
I_{(j_q \ldots j_1 \underbrace{{}_{l, \ldots ,l}}_{k-1})s,t}=
\sum\limits_{(j_1,\ldots,j_q, \underbrace{{}_{l, \ldots ,l}}_{k-1})}
I_{(\small{\fbox{\it l}})s,t} 
I_{(j_q \ldots j_1 \underbrace{{}_{l, \ldots ,l}}_{k-1})s,t},
\end{equation}

\vspace{3mm}
\noindent
where $I_{(j_r \ldots j_1)s,t}$ is defined by (\ref{new1701})
and $\fbox{\it l}$ is the symbol $l$ which does not participate
in the following sum with respect to permutations
$$
\sum\limits_{(j_1,\ldots,j_q, \underbrace{{}_{l, \ldots ,l}}_{k-1})}.
$$

\vspace{1mm}

By analogy with (\ref{new1301}) we obtain

\vspace{-2mm}
$$
\sum\limits_{(j_1,\ldots,j_q, \underbrace{{}_{l, \ldots ,l}}_{k-1})}
I_{(\small{\fbox{\it l}})s,t}
I_{(j_q \ldots j_1 \underbrace{{}_{l, \ldots ,l}}_{k-1})s,t}=
$$

\vspace{-2mm}
$$
=\sum\limits_{(j_1,\ldots,j_q, \underbrace{{}_{l, \ldots ,l}}_{k-1})}
\left(I_{(\small{\fbox{\it l}}j_q\ldots j_1 \underbrace{l \ldots l}_{k-1})s,t}+
I_{(\small{j_q\fbox{\it l}}j_{q-1}\ldots j_1 \underbrace{l \ldots l}_{k-1})s,t}
+\ldots\right.
$$
$$
\left.\ldots +I_{(j_q\ldots j_1 \small{\fbox{\it l}}\underbrace{l \ldots l}_{k-1})s,t}
+
I_{(j_q\ldots j_1 l\small{\fbox{\it l}}\underbrace{l \ldots l}_{k-2})s,t}
+
\ldots + I_{(j_q\ldots j_1 \small{\underbrace{l \ldots l}_{k-1}}\small{\fbox{\it l}})s,t}
\right)=
$$

\vspace{-2mm}
\begin{equation}
\label{new1500}
=\sum\limits_{(j_1,\ldots,j_q, \underbrace{{}_{l, \ldots ,l}}_{k})}
I_{(j_q\ldots j_1 \small{\underbrace{l \ldots l}_{k}})s,t}.
\end{equation}

Substituting $s=T$ into (\ref{new1499}), (\ref{new1500})
and combining (\ref{new1499}), (\ref{new1500}), we conlude that the equality (\ref{new1200}) 
is proved for $p=k.$ The equality (\ref{new1200}) 
is proved.

Note that
$$
n! \int\limits_t^T \phi_l(t_n)\ldots \int\limits_t^{t_2}
\phi_l(t_1)dt_1\ldots dt_n = n! \frac{1}{n!}
\left(\int\limits_t^T \phi_l(\tau)d\tau\right)^n=
$$
\begin{equation}
\label{new1505}
=
\left(\int\limits_t^T \phi_l(\tau)d\tau\right)^n,
\end{equation}

\noindent
where $n\in{\bf N}.$

After substituting (\ref{new1505}) into (\ref{new1200}), we have for $p=n$
\begin{equation}
\label{new1506}
~~~~~~\left(\int\limits_t^T \phi_l(\tau)d\tau\right)^n
\sum\limits_{(j_1,\ldots,j_q)}
J_{(j_q\ldots j_1)T,t}
=\sum\limits_{(j_1,\ldots,j_q, \underbrace{{}_{l, \ldots ,l}}_{n})}
J_{(j_q\ldots j_1 \small{\underbrace{l \ldots l}_{n}})T,t}.
\end{equation}

\vspace{2mm}

The equality (\ref{new1506}) means that
\begin{equation}
\label{new1507b}
~~~J''[\phi_{j_1}\ldots \phi_{j_q}\underbrace{\phi_{l}\ldots \phi_{l}}_{n}]^
{(\hspace{0.5mm}\small{\overbrace{0 \ldots 0}^{q+n}}\hspace{0.5mm})}_{T,t}=
\left(\int\limits_t^T \phi_l(\tau)d\tau\right)^n
\cdot J''[\phi_{j_1}\ldots \phi_{j_q}]^
{(\hspace{0.5mm}\small{\overbrace{0 \ldots 0}^{q}}\hspace{0.5mm})}_{T,t},
\end{equation}
where $n, q=0,1,2\ldots $ 
and $J''[\phi_{j_1}\ldots \phi_{j_q}]^
{(0 \ldots 0)}_{T,t}\stackrel{\sf def}{=}1$
for $q=0.$

The relations (\ref{new1102}) and (\ref{new1507b}) prove 
Theorem~1.22 for the case $i_1=\ldots=i_k=0, 1,\ldots, m$ and $j_1,\ldots,j_k\in \{0\}\cup {\bf N}$.

\vspace{2mm}

{\bf Remark~1.15.}\ {\it Note that the equality
{\rm (\ref{new1200})} can be obtained in another way.
Let $D_q=\left\{(t_1,\ldots,t_q)\in [t, T]^q:\
\exists\ i\ne j\ \hbox{such that}\ t_i=t_j\right\}$ be the "diagonal set" of $[t,T]^q$
$(q=2,3,\ldots)$
{\rm \cite{Kuo}}. Since the Lebesgue meashure of the set $D_q$ is equal to zero {\rm \cite{Kuo}}$,$ 
then {\rm(}see {\rm (\ref{chain10100}))}
\begin{equation}
\label{new9000a}
~~~~~~~J''[\phi_{j_1}\ldots \phi_{j_q}]^
{(\hspace{0.5mm}\small{\overbrace{0 \ldots 0}^{q}}\hspace{0.5mm})}_{T,t}=
\int\limits_{[t,T]^q}\phi_{j_1}(t_1)\ldots \phi_{j_q}(t_q)dt_1\ldots dt_q.
\end{equation}

From {\rm (\ref{new9000a})} we have
$$
J''[\phi_{l}\ldots \phi_{l}]^
{(\hspace{0.5mm}\small{\overbrace{0 \ldots 0}^{p}}\hspace{0.5mm})}_{T,t}\cdot
J''[\phi_{j_1}\ldots \phi_{j_q}]^
{(\hspace{0.5mm}\small{\overbrace{0 \ldots 0}^{q}}\hspace{0.5mm})}_{T,t}=
$$

\vspace{-2mm}
$$
=\int\limits_{[t,T]^q}\phi_{j_1}(t_1)\ldots \phi_{j_q}(t_q)dt_1\ldots dt_q
\int\limits_{[t,T]^p}\phi_{l}(t_1)\ldots \phi_{l}(t_p)dt_1\ldots dt_p=
$$

$$
=\int\limits_{[t,T]^{p+q}}\phi_{j_1}(t_1)\ldots \phi_{j_q}(t_q)
\phi_l(t_1')\ldots \phi_l(t_p')dt_1'\ldots dt_p'dt_1\ldots dt_q=
$$
\begin{equation}
\label{new100000}
=J''[\phi_{j_1}\ldots \phi_{j_q} \phi_{l}\ldots \phi_{l}]^
{(\hspace{0.5mm}\small{\overbrace{0 \ldots 0}^{p+q}}\hspace{0.5mm})}_{T,t}.
\end{equation}

\vspace{4mm}

It is not difficult to see that the equality {\rm (\ref{new100000})} is nothing but
the equality {\rm (\ref{new1200})} written in another form.}

\vspace{2mm}

To complete the proof of Theorem~1.22, we need to consider the case
$i_1,\ldots,i_k=0, 1,\ldots, m$ and $j_1,\ldots,j_k\in \{0\}\cup {\bf N}$.

Obviously, the proof of Theorem~1.22 will be completed if we prove the following equalities

\vspace{-3mm}
$$
\sum\limits_{(j_1,\ldots,j_q)}\int\limits_t^T \phi_{j_q}(t_q)\ldots
\int\limits_t^{t_2}\phi_{j_1}(t_1)
d{\bf w}_{t_1}^{(i_1)}\ldots d{\bf w}_{t_q}^{(i_q)}\times
$$

\vspace{-2mm}
$$
\times \sum\limits_{(j_1',\ldots,j_n')}
\int\limits_t^T 
\phi_{j_n'}(t_n')\ldots \int\limits_t^{t_2'}\phi_{j_1'}(t_1')d{\bf w}_{t_1'}^{(1)}\ldots 
d{\bf w}_{t_n'}^{(1)}=
$$

\vspace{-2mm}
$$
=\sum\limits_{(j_1,\ldots,j_q,j_1',\ldots,j_n')}
\int\limits_t^T \phi_{j_q}(t_q)\ldots \int\limits_t^{t_2}\phi_{j_1}(t_1)
\int\limits_t^{t_1}\phi_{j_n'}(t_n')\ldots \int\limits_t^{t_2'}
\phi_{j_1'}(t_1')\times
$$
\begin{equation}
\label{new1600}
\times d{\bf w}_{t_1'}^{(1)}\ldots d{\bf w}_{t_n'}^{(1)}d{\bf w}_{t_1}^{(i_1)}\ldots d{\bf w}_{t_q}^{(i_q)},
\end{equation}

\vspace{3mm}
$$
\sum\limits_{(j_1,\ldots,j_q)}\int\limits_t^T \phi_{j_q}(t_q)\ldots
\int\limits_t^{t_2}\phi_{j_1}(t_1)
d{\bf w}_{t_1}^{(i_1)}\ldots d{\bf w}_{t_q}^{(i_q)}\times
$$

\vspace{-1mm}
$$
\times \sum\limits_{(j_1',\ldots,j_n')}
\int\limits_t^T \phi_{j_n'}(t_n')\ldots \int\limits_t^{t_2'}\phi_{j_1'}(t_1')
d{\bf w}_{t_1'}^{(0)}\ldots d{\bf w}_{t_n'}^{(0)}=
$$

\vspace{-1mm}
$$
=\sum\limits_{(j_1,\ldots,j_q,j_1',\ldots,j_n')}
\int\limits_t^T \phi_{j_q}(t_q)\ldots \int\limits_t^{t_2}\phi_{j_1}(t_1)
\int\limits_t^{t_1}\phi_{j_n'}(t_n')\ldots \int\limits_t^{t_2'}
\phi_{j_1'}(t_1')\times
$$

\vspace{3mm}
\begin{equation}
\label{new1600a}
\times d{\bf w}_{t_1'}^{(0)}\ldots d{\bf w}_{t_n'}^{(0)}d{\bf w}_{t_1}^{(i_1)}\ldots d{\bf w}_{t_q}^{(i_q)}
\end{equation}

\vspace{5mm}
\noindent
w.~p.~1, where $n, q\in {\bf N},$\ $d{\bf w}_{\tau}^{(0)}
\stackrel{\sf def}{=}d\tau,$\ $i_1,\ldots,i_q\ne 1$ in (\ref{new1600})
and $i_1,\ldots,i_q\ne 0$ in (\ref{new1600a}),
$$
\sum\limits_{(j_1,\ldots,j_g)}
$$

\vspace{1mm}
\noindent
means the sum with respect to all possible permutations $(j_1,\ldots,j_g)$.
At the same time if $j_r$ swapped with $j_d$ in the permutation $(j_1,\ldots,j_g)$, then
$i_r$ swapped with $i_d$ in the permutation $(i_1,\ldots,i_g).$

The equalities (\ref{new1600}) and (\ref{new1600a}) mean that

\vspace{-3mm}
\begin{equation}
\label{new1601}
J''[\phi_{j_1}\ldots\phi_{j_q}\phi_{j_1'}\ldots\phi_{j_n'}]_{T,t}^{(i_1\ldots i_q 1\ldots 1)}=
J''[\phi_{j_1}\ldots\phi_{j_q}]_{T,t}^{(i_1\ldots i_q)}
\cdot J''[\phi_{j_1'}\ldots\phi_{j_n'}]^{(1\ldots 1)}_{T,t},
\end{equation}

\begin{equation}
\label{new1601a}
J''[\phi_{j_1}\ldots\phi_{j_q}\phi_{j_1'}\ldots\phi_{j_n'}]_{T,t}^{(i_1\ldots i_q 0\ldots 0)}=
J''[\phi_{j_1}\ldots\phi_{j_q}]_{T,t}^{(i_1\ldots i_q)}
\cdot J''[\phi_{j_1'}\ldots\phi_{j_n'}]^{(0\ldots 0)}_{T,t}
\end{equation}

\vspace{6mm}
\noindent
w.~p.~1, where $i_1,\ldots,i_q\ne 1$ in (\ref{new1601}) and
$i_1,\ldots,i_q\ne 0$ in (\ref{new1601a}).

First, we prove the equality (\ref{new1600}). Consider the case $n=1.$
Using the It\^{o} formula, we get w.~p.~1
$$
\int\limits_t^{s}\phi_{j_1'}(\theta)d{\bf w}_{\theta}^{(1)}
\int\limits_t^s \phi_{j_q}(t_q)\ldots
\int\limits_t^{t_2}\phi_{j_1}(t_1)
d{\bf w}_{t_1}^{(i_1)}\ldots d{\bf w}_{t_q}^{(i_q)}=
$$
$$
=J_{(j_1' j_q\ldots j_1)s,t}^{(1 i_q\ldots i_1)}+
$$

\vspace{-5mm}
$$
+\hspace{-0.6mm}
\int\limits_t^s\hspace{-0.8mm} \phi_{j_q}(\tau)\hspace{-0.6mm}\left(\int\limits_t^{\tau}
\phi_{j_{q-1}}(t_{q-1})\ldots
\int\limits_t^{t_2}\phi_{j_1}(t_1)
d{\bf w}_{t_1}^{(i_1)}\ldots d{\bf w}_{t_{q-1}}^{(i_{q-1})}
\hspace{-0.6mm}\int\limits_t^{\tau}\phi_{j_1'}(\theta)d{\bf w}_{\theta}^{(1)}\right)\hspace{-0.6mm}
d{\bf w}_{\tau}^{(i_{q})}\hspace{-0.5mm}
$$
$$
=\ldots =
$$

\vspace{-3mm}
\begin{equation}
\label{new1700}
=J_{(j_1' j_q\ldots j_1)s,t}^{(1 i_q\ldots i_1)}+J_{(j_q j_1' j_{q-1}\ldots j_1)s,t}
^{(i_q 1 i_{q-1}\ldots i_1)}+\ldots +J_{(j_q\ldots j_1 j_1')s,t}^{(i_q\ldots i_1 1)},
\end{equation}

\vspace{5mm}
\noindent
where
\begin{equation}
\label{new100001}
\int\limits_t^s \phi_{j_r}(t_r)\ldots \int\limits_t^{t_2}
\phi_{j_1}(t_1)d{\bf w}_{t_1}^{(i_1)}\ldots d{\bf w}_{t_r}^{(i_r)}
\stackrel{\sf def}{=}J_{(j_r \ldots j_1)s,t}^{(i_r\ldots i_1)},
\end{equation}

\vspace{3mm}
\noindent
$i_1,\ldots,i_r=0,1,\ldots,m.$

From (\ref{new1700}) we obtain
$$
\int\limits_t^{s}\phi_{j_1'}(\theta)d{\bf w}_{\theta}^{(1)}
\sum\limits_{(j_1,\ldots,j_q)}\int\limits_t^s \phi_{j_q}(t_q)\ldots
\int\limits_t^{t_2}\phi_{j_1}(t_1)
d{\bf w}_{t_1}^{(i_1)}\ldots d{\bf w}_{t_q}^{(i_q)}=
$$

\vspace{-2mm}
$$
=\sum\limits_{(j_1,\ldots,j_q)}\int\limits_t^{s}\phi_{j_1'}(\theta)d{\bf w}_{\theta}^{(1)}
\int\limits_t^s \phi_{j_q}(t_q)\ldots
\int\limits_t^{t_2}\phi_{j_1}(t_1)
d{\bf w}_{t_1}^{(i_1)}\ldots d{\bf w}_{t_q}^{(i_q)}=
$$

\vspace{3mm}
$$
=
\sum\limits_{(j_1,\ldots,j_q)}\left(
J_{(j_1' j_q\ldots j_1)s,t}^{(1 i_q\ldots i_1)}+J_{(j_q j_1' j_{q-1}\ldots j_1)s,t}
^{(i_q 1 i_{q-1}\ldots i_1)}+\ldots +J_{(j_q\ldots j_1 j_1')s,t}^{(i_q\ldots i_1 1)}\right)=
$$

\vspace{1mm}
\begin{equation}
\label{new1700a}
=
\sum\limits_{(j_1,\ldots,j_q, j_1')}
J_{(j_q\ldots j_1 j_1')s,t}^{(i_q\ldots i_1 1)}
\end{equation}

\vspace{4mm}
\noindent
w.~p.~1, where $J_{(j_r\ldots j_1)s,t}^{(i_r\ldots i_1)}$
is defined by (\ref{new100001}).
The equality (\ref{new1600}) is proved for the case $n=1.$

Let us assume that the equality (\ref{new1600}) is true for $n=2, 3, \ldots, k-1$, and prove
its validity for $n=k.$

Applying (\ref{new1025}), (\ref{new1040}), (\ref{new1020})--(\ref{new1021}), we obtain w.~p.~1
$$
\sum\limits_{(j_1',\ldots,j_k')}
\int\limits_t^s 
\phi_{j_k'}(t_k')\ldots \int\limits_t^{t_2'}\phi_{j_1'}(t_1')d{\bf w}_{t_1'}^{(1)}\ldots 
d{\bf w}_{t_k'}^{(1)}=
$$

\vspace{-2mm}
$$
=\int\limits_t^{s}\phi_{j_k'}(\theta)d{\bf w}_{\theta}^{(1)}
\sum\limits_{(j_1',\ldots,j_{k-1}')}
\int\limits_t^s 
\phi_{j_{k-1}'}(t_{k-1})\ldots \int\limits_t^{t_2}\phi_{j_1'}(t_1)d{\bf w}_{t_1}^{(1)}\ldots 
d{\bf w}_{t_{k-1}}^{(1)}-
$$

\begin{equation}
\label{new1800}
-\sum\limits_{(j_1',\ldots,j_{k-1}')}
\int\limits_t^s \phi_{j_k'}(\theta)\phi_{j_{k-1}'}(\theta)d\theta
\int\limits_t^s 
\phi_{j_{k-2}'}(t_{k-2})\ldots \int\limits_t^{t_2}\phi_{j_1'}(t_1)d{\bf w}_{t_1}^{(1)}\ldots 
d{\bf w}_{t_{k-2}}^{(1)}.
\end{equation}

\vspace{3mm}

After substituting $s=T$ in (\ref{new1800})
and applying the orthonormality of 
$\{\phi_j(x)\}_{j=0}^{\infty}$, we get w.~p.~1
$$
\sum\limits_{(j_1',\ldots,j_k')}
\int\limits_t^T
\phi_{j_k'}(t_k')\ldots \int\limits_t^{t_2'}\phi_{j_1'}(t_1')d{\bf w}_{t_1'}^{(1)}\ldots 
d{\bf w}_{t_k'}^{(1)}=
$$

\vspace{-2mm}
$$
=\int\limits_t^{T}\phi_{j_k'}(\theta)d{\bf w}_{\theta}^{(1)}
\sum\limits_{(j_1',\ldots,j_{k-1}')}
\int\limits_t^T 
\phi_{j_{k-1}'}(t_{k-1})\ldots \int\limits_t^{t_2}\phi_{j_1'}(t_1)d{\bf w}_{t_1}^{(1)}\ldots 
d{\bf w}_{t_{k-1}}^{(1)}-
$$

\begin{equation}
\label{new1801}
~~~~~-\sum\limits_{(j_1',\ldots,j_{k-1}')}
{\bf 1}_{\{j_k'=j_{k-1}'\}}
\int\limits_t^T
\phi_{j_{k-2}'}(t_{k-2})\ldots \int\limits_t^{t_2}\phi_{j_1'}(t_1)d{\bf w}_{t_1}^{(1)}\ldots 
d{\bf w}_{t_{k-2}}^{(1)},
\end{equation}

\vspace{4mm}
\noindent
where ${\bf 1}_{A}$ is the indicator of the set $A$.

Using (\ref{new1801}) and the induction hypothesis, we obtain w.~p.~1

\newpage
\noindent
$$
\sum\limits_{(j_1',\ldots,j_k')}
\int\limits_t^T 
\phi_{j_k'}(t_k)\ldots \int\limits_t^{t_2}\phi_{j_1'}(t_1)d{\bf w}_{t_1}^{(1)}\ldots 
d{\bf w}_{t_k}^{(1)}\times
$$

$$
\times\sum\limits_{(j_1,\ldots,j_q)}\int\limits_t^T \phi_{j_q}(t_q)\ldots
\int\limits_t^{t_2}\phi_{j_1}(t_1)
d{\bf w}_{t_1}^{(i_1)}\ldots d{\bf w}_{t_q}^{(i_q)}=
$$

\vspace{-2mm}
$$
=\int\limits_t^{T}\phi_{j_k'}(\theta)d{\bf w}_{\theta}^{(1)}
\sum\limits_{(j_1',\ldots,j_{k-1}')}
\int\limits_t^T 
\phi_{j_{k-1}'}(t_{k-1})\ldots \int\limits_t^{t_2}\phi_{j_1'}(t_1)d{\bf w}_{t_1}^{(1)}\ldots 
d{\bf w}_{t_{k-1}}^{(1)}\times
$$

$$
\times\sum\limits_{(j_1,\ldots,j_q)}\int\limits_t^T \phi_{j_q}(t_q)\ldots
\int\limits_t^{t_2}\phi_{j_1}(t_1)
d{\bf w}_{t_1}^{(i_1)}\ldots d{\bf w}_{t_q}^{(i_q)}-
$$

\vspace{-2mm}
$$
-\sum\limits_{(j_1',\ldots,j_{k-1}')}
{\bf 1}_{\{j_k'=j_{k-1}'\}}
\int\limits_t^T
\phi_{j_{k-2}'}(t_{k-2})\ldots \int\limits_t^{t_2}\phi_{j_1'}(t_1)d{\bf w}_{t_1}^{(1)}\ldots 
d{\bf w}_{t_{k-2}}^{(1)}\times
$$

\vspace{-2mm}
$$
\times\sum\limits_{(j_1,\ldots,j_q)}\int\limits_t^T \phi_{j_q}(t_q)\ldots
\int\limits_t^{t_2}\phi_{j_1}(t_1)
d{\bf w}_{t_1}^{(i_1)}\ldots d{\bf w}_{t_q}^{(i_q)}=
$$

\vspace{-1mm}
$$
=\int\limits_t^{T}\phi_{j_k'}(\theta)d{\bf w}_{\theta}^{(1)}\times
$$
$$
\times
\sum\limits_{(j_1,\ldots,j_q,j_1',\ldots,j_{k-1}')}
\int\limits_t^T \phi_{j_q}(t_q)\ldots \int\limits_t^{t_2}\phi_{j_1}(t_1)
\int\limits_t^{t_1}\phi_{j_{k-1}'}(t_{k-1}')\ldots \int\limits_t^{t_2'}
\phi_{j_1'}(t_1')\times
$$

\vspace{3mm}
$$
\times d{\bf w}_{t_1'}^{(1)}\ldots 
d{\bf w}_{t_{k-1}'}^{(1)}d{\bf w}_{t_1}^{(i_1)}\ldots d{\bf w}_{t_q}^{(i_q)}-
$$

\vspace{-2mm}
$$
-\sum\limits_{(j_1',\ldots,j_{k-1}')}
{\bf 1}_{\{j_k'=j_{k-1}'\}}
\int\limits_t^T
\phi_{j_{k-2}'}(t_{k-2})\ldots \int\limits_t^{t_2}\phi_{j_1'}(t_1)d{\bf w}_{t_1}^{(1)}\ldots 
d{\bf w}_{t_{k-2}}^{(1)}\times
$$

\vspace{-2mm}
\begin{equation}
\label{newi9000}
\times\sum\limits_{(j_1,\ldots,j_q)}\int\limits_t^T \phi_{j_q}(t_q)\ldots
\int\limits_t^{t_2}\phi_{j_1}(t_1)
d{\bf w}_{t_1}^{(i_1)}\ldots d{\bf w}_{t_q}^{(i_q)}.
\end{equation}

\vspace{5mm}

Further, applying the induction hypothesis, we have w.~p.~1

\vspace{-2mm}
$$
\sum\limits_{(j_1',\ldots,j_{k-1}')}
{\bf 1}_{\{j_k'=j_{k-1}'\}}
\int\limits_t^T
\phi_{j_{k-2}'}(t_{k-2})\ldots \int\limits_t^{t_2}\phi_{j_1'}(t_1)d{\bf w}_{t_1}^{(1)}\ldots 
d{\bf w}_{t_{k-2}}^{(1)}\times
$$

$$
\times\sum\limits_{(j_1,\ldots,j_q)}\int\limits_t^T \phi_{j_q}(t_q)\ldots
\int\limits_t^{t_2}\phi_{j_1}(t_1)
d{\bf w}_{t_1}^{(i_1)}\ldots d{\bf w}_{t_q}^{(i_q)}=
$$

\vspace{-2mm}
$$
=\Biggl(
\sum\limits_{(j_1',\ldots,j_{k-2}')}
{\bf 1}_{\{j_k'=j_{k-1}'\}}
\int\limits_t^T
\phi_{j_{k-2}'}(t_{k-2})\ldots \int\limits_t^{t_2}\phi_{j_1'}(t_1)d{\bf w}_{t_1}^{(1)}\ldots 
d{\bf w}_{t_{k-2}}^{(1)}+
\Biggr.
$$

\vspace{-2mm}
$$
+\sum\limits_{(j_1',\ldots,j_{k-3}', j_{k-1}')}
{\bf 1}_{\{j_k'=j_{k-2}'\}}
\int\limits_t^T
\phi_{j_{k-1}'}(t_{k-2})
\int\limits_t^{t_{k-2}}
\phi_{j_{k-3}'}(t_{k-3})\ldots \int\limits_t^{t_2}\phi_{j_1'}(t_1)\times
$$

\vspace{3mm}
$$
\times d{\bf w}_{t_1}^{(1)}\ldots 
d{\bf w}_{t_{k-3}}^{(1)}d{\bf w}_{t_{k-2}}^{(1)}+\ldots
$$

$$
\ldots+\sum\limits_{(j_2',\ldots,j_{k-1}')}
{\bf 1}_{\{j_k'=j_{1}'\}}
\int\limits_t^T
\phi_{j_{k-2}'}(t_{k-2})\ldots \int\limits_t^{t_3}\phi_{j_2'}(t_2)
\int\limits_t^{t_2}
\phi_{j_{k-1}'}(t_{1})\times
$$

$$
\Biggl.\times
d{\bf w}_{t_1}^{(1)}d{\bf w}_{t_2}^{(1)}\ldots 
d{\bf w}_{t_{k-2}}^{(1)}\Biggr)\times
$$

\vspace{2mm}
$$
\times\sum\limits_{(j_1,\ldots,j_q)}\int\limits_t^T \phi_{j_q}(t_q)\ldots
\int\limits_t^{t_2}\phi_{j_1}(t_1)
d{\bf w}_{t_1}^{(i_1)}\ldots d{\bf w}_{t_q}^{(i_q)}=
$$

\vspace{2mm}
$$
=\Biggl(
{\bf 1}_{\{j_k'=j_{k-1}'\}}\sum\limits_{(j_1',\ldots,j_{k-2}')}
\int\limits_t^T
\phi_{j_{k-2}'}(t_{k-2})\ldots \int\limits_t^{t_2}\phi_{j_1'}(t_1)d{\bf w}_{t_1}^{(1)}\ldots 
d{\bf w}_{t_{k-2}}^{(1)}+
\Biggr.
$$

\vspace{-2mm}
$$
+{\bf 1}_{\{j_k'=j_{k-2}'\}}\sum\limits_{(j_1',\ldots,j_{k-3}', j_{k-1}')}
\int\limits_t^T
\phi_{j_{k-1}'}(t_{k-2})
\int\limits_t^{t_{k-2}}
\phi_{j_{k-3}'}(t_{k-3})\ldots \int\limits_t^{t_2}\phi_{j_1'}(t_1)\times
$$
$$
\times d{\bf w}_{t_1}^{(1)}\ldots 
d{\bf w}_{t_{k-3}}^{(1)}d{\bf w}_{t_{k-2}}^{(1)}+\ldots
$$

$$
\ldots +{\bf 1}_{\{j_k'=j_{1}'\}}\sum\limits_{(j_2',\ldots,j_{k-1}')}
\int\limits_t^T
\phi_{j_{k-2}'}(t_{k-2})\ldots \int\limits_t^{t_3}\phi_{j_2'}(t_2)
\int\limits_t^{t_2}
\phi_{j_{k-1}'}(t_{1})\times
$$

\vspace{2mm}
$$
\Biggl.\times
d{\bf w}_{t_1}^{(1)}d{\bf w}_{t_2}^{(1)}\ldots 
d{\bf w}_{t_{k-2}}^{(1)}\Biggr)\times
$$

\vspace{2mm}
$$
\times\sum\limits_{(j_1,\ldots,j_q)}\int\limits_t^T \phi_{j_q}(t_q)\ldots
\int\limits_t^{t_2}\phi_{j_1}(t_1)
d{\bf w}_{t_1}^{(i_1)}\ldots d{\bf w}_{t_q}^{(i_q)}=
$$

\vspace{-2mm}
$$
={\bf 1}_{\{j_k'=j_{k-1}'\}}\sum\limits_{(j_1,\ldots,j_q,j_1',\ldots,j_{k-2}')}
\int\limits_t^T \phi_{j_q}(t_q)\ldots \int\limits_t^{t_2}\phi_{j_1}(t_1)
\int\limits_t^{t_1}\phi_{j_{k-2}'}(t_{k-2}')\ldots \int\limits_t^{t_2'}
\phi_{j_1'}(t_1')\times
$$

\vspace{3mm}
$$
\times d{\bf w}_{t_1'}^{(1)}\ldots d{\bf w}_{t_{k-2}'}^{(1)}
d{\bf w}_{t_1}^{(i_1)}\ldots d{\bf w}_{t_q}^{(i_q)}+
$$

\vspace{1mm}
$$
+{\bf 1}_{\{j_k'=j_{k-2}'\}}\sum\limits_{(j_1,\ldots,j_q,j_1',\ldots,j_{k-3}',j_{k-1}')}
\int\limits_t^T \phi_{j_q}(t_q)\ldots \int\limits_t^{t_2}\phi_{j_1}(t_1)
\int\limits_t^{t_1}\phi_{j_{k-1}'}(t_{k-2}')\times
$$

$$
\times\int\limits_t^{t_{k-2}'}
\phi_{j_{k-3}'}(t_{k-3}')
\ldots \int\limits_t^{t_2'}
\phi_{j_1'}(t_1')d{\bf w}_{t_1'}^{(1)}\ldots d{\bf w}_{t_{k-3}'}^{(1)}d{\bf w}_{t_{k-2}'}^{(1)}
d{\bf w}_{t_1}^{(i_1)}\ldots d{\bf w}_{t_q}^{(i_q)}+\ldots
$$

\vspace{-2mm}
$$
\ldots
$$

\vspace{-2mm}
$$
\ldots +{\bf 1}_{\{j_k'=j_{1}'\}}\sum\limits_{(j_1,\ldots,j_q,j_2',\ldots,j_{k-1}')}
\int\limits_t^T \phi_{j_q}(t_q)\ldots \int\limits_t^{t_2}\phi_{j_1}(t_1)\times
$$

$$
\times
\int\limits_t^{t_1}\phi_{j_{k-2}'}(t_{k-2}')\ldots \int\limits_t^{t_{3}'}
\phi_{j_{2}'}(t_{2}')
\int\limits_t^{t_2'}
\phi_{j_{k-1}'}(t_1')d{\bf w}_{t_1'}^{(1)}d{\bf w}_{t_2'}^{(1)}\ldots d{\bf w}_{t_{k-2}'}^{(1)}
d{\bf w}_{t_1}^{(i_1)}\ldots d{\bf w}_{t_q}^{(i_q)}\stackrel{\sf def}{=}
$$
\begin{equation}
\label{new1900}
\stackrel{\sf def}{=}S_4(T).
\end{equation}

\vspace{5mm}

By analogy with (\ref{newx1020ss}) we obtain w.~p.~1

\vspace{-2mm}
$$
\int\limits_t^T \phi_l(\tau)\phi_{j_r}(\tau)d\tau
\int\limits_t^{T} \phi_{j_{r-1}}(t_{r-1})\ldots \int\limits_t^{t_2}
\phi_{j_1}(t_1)d{\bf w}_{t_1}^{(i_1)}\ldots d{\bf w}_{t_{r-1}}^{(i_{r-1})}=
$$

\vspace{-2mm}
$$
=\int\limits_t^T \phi_l(\tau)\phi_{j_r}(\tau)
\int\limits_t^{\tau} \phi_{j_{r-1}}(t_{r-1})\ldots \int\limits_t^{t_2}
\phi_{j_1}(t_1)d{\bf w}_{t_1}^{(i_1)}\ldots d{\bf w}_{t_{r-1}}^{(i_{r-1})}d\tau+\ldots
$$

\vspace{-2mm}
\begin{equation}
\label{new5001}
\ldots +
\int\limits_t^{T} \phi_{j_{r-1}}(t_{r-1})\ldots \int\limits_t^{t_2}
\phi_{j_1}(t_1)\int\limits_t^{t_1} \phi_{l}(\tau)\phi_{j_r}(\tau)      
d\tau d{\bf w}_{t_1}^{(i_1)}\ldots d{\bf w}_{t_{r-1}}^{(i_{r-1})},
\end{equation}

\vspace{3mm}
\noindent
where $i_1,\ldots,i_{r-1}=0,1,\ldots,m.$

Using iteratively the It\^{o} formula, as well as (\ref{new5001})
and combinatorial reasoning, we obtain w.~p.~1 (see Remark~1.16 below for details)

\vspace{-2mm}
$$
\int\limits_t^{T}\phi_{j_k'}(\theta)d{\bf w}_{\theta}^{(1)}\times
$$

\vspace{-4mm}
$$
\times
\sum\limits_{(j_1,\ldots,j_q,j_1',\ldots,j_{k-1}')}
\int\limits_t^T \phi_{j_q}(t_q)\ldots \int\limits_t^{t_2}\phi_{j_1}(t_1)
\int\limits_t^{t_1}\phi_{j_{k-1}'}(t_{k-1}')\ldots \int\limits_t^{t_2'}
\phi_{j_1'}(t_1')\times
$$

\vspace{2mm}
$$
\times d{\bf w}_{t_1'}^{(1)}\ldots 
d{\bf w}_{t_{k-1}'}^{(1)}d{\bf w}_{t_1}^{(i_1)}\ldots d{\bf w}_{t_q}^{(i_q)}=
$$

\vspace{2mm}
$$
=\sum\limits_{(j_1,\ldots,j_q,j_1',\ldots,j_{k}')}
\int\limits_t^T \phi_{j_q}(t_q)\ldots \int\limits_t^{t_2}\phi_{j_1}(t_1)
\int\limits_t^{t_1}\phi_{j_{k}'}(t_{k}')\ldots \int\limits_t^{t_2'}
\phi_{j_1'}(t_1')\times
$$

\vspace{2mm}
$$
\times d{\bf w}_{t_1'}^{(1)}\ldots 
d{\bf w}_{t_{k}'}^{(1)}d{\bf w}_{t_1}^{(i_1)}\ldots d{\bf w}_{t_q}^{(i_q)}+
$$
$$
+\sum\limits_{(j_1,\ldots,j_q,j_1',\ldots,j_{k-1}')}
\Biggl(\int\limits_t^T \phi_{j_q}(t_q)\ldots \int\limits_t^{t_2}\phi_{j_1}(t_1)
\int\limits_t^{t_1}\phi_{j_{k}'}(\theta)\phi_{j_{k-1}'}(\theta)\int\limits_t^{\theta}
\phi_{j_{k-2}'}(t_{k-2}')\ldots\Biggr.
$$
$$
\Biggl.\ldots \int\limits_t^{t_2'}\phi_{j_1'}(t_1')
d{\bf w}_{t_1'}^{(1)}\ldots 
d{\bf w}_{t_{k-2}'}^{(1)}d{\bf w}_{\theta}^{(0)}d{\bf w}_{t_1}^{(i_1)}\ldots d{\bf w}_{t_q}^{(i_q)}+
$$

$$
+
\int\limits_t^T \phi_{j_q}(t_q)\ldots \int\limits_t^{t_2}
\phi_{j_1}(t_1)\int\limits_t^{t_1}\phi_{j_{k-1}'}(t_{k-1}')
\int\limits_t^{t_{k-1}'}\phi_{j_{k}'}(\theta)\phi_{j_{k-2}'}(\theta)\int\limits_t^{\theta}
\phi_{j_{k-3}'}(t_{k-3}')\ldots
$$
$$
\ldots
\int\limits_t^{t_2'}
\phi_{j_{1}'}(t_{1}')d{\bf w}_{t_1'}^{(1)}\ldots 
d{\bf w}_{t_{k-3}'}^{(1)}d{\bf w}_{\theta}^{(0)}d{\bf w}_{t_{k-1}'}^{(1)}
d{\bf w}_{t_1}^{(i_1)}\ldots d{\bf w}_{t_q}^{(i_q)}+\ldots
$$

$$
\ldots +
\int\limits_t^T \phi_{j_q}(t_q)\ldots \int\limits_t^{t_2}
\phi_{j_1}(t_1)\int\limits_t^{t_1}\phi_{j_{k-1}'}(t_{k-1}')\ldots
\int\limits_t^{t_{3}'}\phi_{j_{2}'}(t_2')\int\limits_t^{t_{2}'}
\phi_{j_{k}'}(\theta)\phi_{j_{1}'}(\theta)d{\bf w}_{\theta}^{(0)}\times
$$
$$
\Biggl.\times
d{\bf w}_{t_2'}^{(1)}\ldots 
d{\bf w}_{t_{k-1}'}^{(1)}
d{\bf w}_{t_1}^{(i_1)}\ldots d{\bf w}_{t_q}^{(i_q)}\Biggr)=
$$

\vspace{2mm}
$$
=\sum\limits_{(j_1,\ldots,j_q,j_1',\ldots,j_{k}')}
\int\limits_t^T \phi_{j_q}(t_q)\ldots \int\limits_t^{t_2}\phi_{j_1}(t_1)
\int\limits_t^{t_1}\phi_{j_{k}'}(t_{k}')\ldots \int\limits_t^{t_2'}
\phi_{j_1'}(t_1')\times
$$

$$
\times d{\bf w}_{t_1'}^{(1)}\ldots 
d{\bf w}_{t_{k}'}^{(1)}d{\bf w}_{t_1}^{(i_1)}\ldots d{\bf w}_{t_q}^{(i_q)}+
$$

$$
+\sum\limits_{(j_1,\ldots,j_q,j_1',\ldots,j_{k-2}')}
\Biggl\{\int\limits_t^T 
\phi_{j_{k}'}(\theta)\phi_{j_{k-1}'}(\theta)\int\limits_t^{\theta}
\phi_{j_q}(t_q)\ldots \int\limits_t^{t_2}\phi_{j_1}(t_1)
\int\limits_t^{t_1}
\phi_{j_{k-2}'}(t_{k-2}')\ldots\Biggr.
$$
$$
\Biggl.\ldots \int\limits_t^{t_2'}\phi_{j_1'}(t_1')
d{\bf w}_{t_1'}^{(1)}\ldots 
d{\bf w}_{t_{k-2}'}^{(1)}d{\bf w}_{t_1}^{(i_1)}\ldots d{\bf w}_{t_q}^{(i_q)}d{\bf w}_{\theta}^{(0)}+\ldots
$$
$$
\ldots+\int\limits_t^T 
\phi_{j_q}(t_q)\ldots \int\limits_t^{t_2}\phi_{j_1}(t_1)
\int\limits_t^{t_1}
\phi_{j_{k-2}'}(t_{k-2}')\ldots \int\limits_t^{t_2'}\phi_{j_1'}(t_1')
\int\limits_t^{t_1'} 
\phi_{j_{k}'}(\theta)\phi_{j_{k-1}'}(\theta)d{\bf w}_{\theta}^{(0)}\times
$$
$$
\Biggl.
\times d{\bf w}_{t_1'}^{(1)}\ldots 
d{\bf w}_{t_{k-2}'}^{(1)}d{\bf w}_{t_1}^{(i_1)}\ldots d{\bf w}_{t_q}^{(i_q)}\Biggr\}+
$$

$$
+\sum\limits_{(j_1,\ldots,j_q,j_1',\ldots,j_{k-3}',j_{k-1}')}
\Biggl\{\int\limits_t^T 
\phi_{j_{k}'}(\theta)\phi_{j_{k-2}'}(\theta)\int\limits_t^{\theta}
\phi_{j_q}(t_q)\ldots \int\limits_t^{t_2}\phi_{j_1}(t_1)
\int\limits_t^{t_1}
\phi_{j_{k-1}'}(t_{k-1}')\times\Biggr.
$$
$$
\times\int\limits_t^{t_{k-1}'}
\phi_{j_{k-3}'}(t_{k-3}')
\ldots \int\limits_t^{t_2'}\phi_{j_1'}(t_1')
d{\bf w}_{t_1'}^{(1)}\ldots 
d{\bf w}_{t_{k-3}'}^{(1)}d{\bf w}_{t_{k-1}'}^{(1)}
d{\bf w}_{t_1}^{(i_1)}\ldots d{\bf w}_{t_q}^{(i_q)}
d{\bf w}_{\theta}^{(0)}+\ldots
$$

$$
\ldots+\int\limits_t^T 
\phi_{j_q}(t_q)\ldots \int\limits_t^{t_2}\phi_{j_1}(t_1)
\int\limits_t^{t_1}
\phi_{j_{k-1}'}(t_{k-1}')\int\limits_t^{t_{k-1}'}
\phi_{j_{k-3}'}(t_{k-3}')\ldots
\int\limits_t^{t_2'}\phi_{j_1'}(t_1')
\times
$$
$$
\Biggl.\times
\int\limits_t^{t_1'} 
\phi_{j_{k}'}(\theta)\phi_{j_{k-2}'}(\theta)d{\bf w}_{\theta}^{(0)}
d{\bf w}_{t_1'}^{(1)}\ldots 
d{\bf w}_{t_{k-3}'}^{(1)}d{\bf w}_{t_{k-1}'}^{(1)}
d{\bf w}_{t_1}^{(i_1)}\ldots d{\bf w}_{t_q}^{(i_q)}\Biggr\}+\ldots
$$

$$
+\sum\limits_{(j_1,\ldots,j_q,j_2',\ldots,j_{k-1}')}
\Biggl\{\int\limits_t^T 
\phi_{j_{k}'}(\theta)\phi_{j_{1}'}(\theta)\int\limits_t^{\theta}
\phi_{j_q}(t_q)\ldots \int\limits_t^{t_2}\phi_{j_1}(t_1)
\int\limits_t^{t_1}
\phi_{j_{k-1}'}(t_{k-1}')\ldots\Biggr.
$$
$$
\ldots\int\limits_t^{t_{3}'}
\phi_{j_2'}(t_2')
d{\bf w}_{t_2'}^{(1)}\ldots 
d{\bf w}_{t_{k-1}'}^{(1)}
d{\bf w}_{t_1}^{(i_1)}\ldots d{\bf w}_{t_q}^{(i_q)}
d{\bf w}_{\theta}^{(0)}+\ldots
$$

$$
\ldots+\int\limits_t^T 
\phi_{j_q}(t_q)\ldots \int\limits_t^{t_2}\phi_{j_1}(t_1)
\int\limits_t^{t_1}
\phi_{j_{k-1}'}(t_{k-1}')\ldots
\int\limits_t^{t_3'}\phi_{j_2'}(t_2')\int\limits_t^{t_2'}
\phi_{j_{k}'}(\theta)\phi_{j_{1}'}(\theta)d{\bf w}_{\theta}^{(0)}\times
$$
$$
\Biggl.\times
d{\bf w}_{t_2'}^{(1)}\ldots 
d{\bf w}_{t_{k-1}'}^{(1)}
d{\bf w}_{t_1}^{(i_1)}\ldots d{\bf w}_{t_q}^{(i_q)}\Biggr\}=
$$
$$
=\sum\limits_{(j_1,\ldots,j_q,j_1',\ldots,j_{k}')}
\int\limits_t^T \phi_{j_q}(t_q)\ldots \int\limits_t^{t_2}\phi_{j_1}(t_1)
\int\limits_t^{t_1}\phi_{j_{k}'}(t_{k}')\ldots \int\limits_t^{t_2'}
\phi_{j_1'}(t_1')\times
$$

$$
\times d{\bf w}_{t_1'}^{(1)}\ldots 
d{\bf w}_{t_{k}'}^{(1)}d{\bf w}_{t_1}^{(i_1)}\ldots d{\bf w}_{t_q}^{(i_q)}+
$$

$$
+\int\limits_t^T 
\phi_{j_{k}'}(\theta)\phi_{j_{k-1}'}(\theta)d\theta\sum\limits_{(j_1,\ldots,j_q,j_1',\ldots,j_{k-2}')}
\int\limits_t^{T}
\phi_{j_q}(t_q)\ldots \int\limits_t^{t_2}\phi_{j_1}(t_1)
\int\limits_t^{t_1}
\phi_{j_{k-2}'}(t_{k-2}')\ldots
$$
$$
\ldots \int\limits_t^{t_2'}\phi_{j_1'}(t_1')
d{\bf w}_{t_1'}^{(1)}\ldots 
d{\bf w}_{t_{k-2}'}^{(1)}d{\bf w}_{t_1}^{(i_1)}\ldots d{\bf w}_{t_q}^{(i_q)}+
$$

$$
+\int\limits_t^T 
\phi_{j_{k}'}(\theta)\phi_{j_{k-2}'}(\theta)d\theta
\sum\limits_{(j_1,\ldots,j_q,j_1',\ldots,j_{k-3}',j_{k-1}')}
\int\limits_t^{T}
\phi_{j_q}(t_q)\ldots \int\limits_t^{t_2}\phi_{j_1}(t_1)
\int\limits_t^{t_1}
\phi_{j_{k-1}'}(t_{k-1}')\times\Biggr.
$$
$$
\times\int\limits_t^{t_{k-1}'}
\phi_{j_{k-3}'}(t_{k-3}')
\ldots \int\limits_t^{t_2'}\phi_{j_1'}(t_1')
d{\bf w}_{t_1'}^{(1)}\ldots 
d{\bf w}_{t_{k-3}'}^{(1)}d{\bf w}_{t_{k-1}'}^{(1)}
d{\bf w}_{t_1}^{(i_1)}\ldots d{\bf w}_{t_q}^{(i_q)}
+\ldots
$$

$$
\ldots +\int\limits_t^T 
\phi_{j_{k}'}(\theta)\phi_{j_{1}'}(\theta)d\theta\sum\limits_{(j_1,\ldots,j_q,j_2',\ldots,j_{k-1}')}
\int\limits_t^{T}
\phi_{j_q}(t_q)\ldots \int\limits_t^{t_2}\phi_{j_1}(t_1)
\int\limits_t^{t_1}
\phi_{j_{k-1}'}(t_{k-1}')\ldots
$$
$$
\ldots\int\limits_t^{t_{3}'}
\phi_{j_2'}(t_2')
d{\bf w}_{t_2'}^{(1)}\ldots 
d{\bf w}_{t_{k-1}'}^{(1)}
d{\bf w}_{t_1}^{(i_1)}\ldots d{\bf w}_{t_q}^{(i_q)}=
$$

$$
=\sum\limits_{(j_1,\ldots,j_q,j_1',\ldots,j_{k}')}
\int\limits_t^T \phi_{j_q}(t_q)\ldots \int\limits_t^{t_2}\phi_{j_1}(t_1)
\int\limits_t^{t_1}\phi_{j_{k}'}(t_{k}')\ldots \int\limits_t^{t_2'}
\phi_{j_1'}(t_1')\times
$$

\begin{equation}
\label{new1901}
\times d{\bf w}_{t_1'}^{(1)}\ldots 
d{\bf w}_{t_{k}'}^{(1)}d{\bf w}_{t_1}^{(i_1)}\ldots d{\bf w}_{t_q}^{(i_q)}+S_4(T).
\end{equation}

\vspace{6mm}

From (\ref{newi9000}), (\ref{new1900}), and (\ref{new1901}) we conclude that
the equality (\ref{new1600}) is proved for $n=k.$
The equality (\ref{new1600}) is proved.

\vspace{2mm}

{\bf Remark~1.16.}\ {\it It should be noted that the sums with respect to 
permutations 
$$
\sum\limits_{(j_1,\ldots,j_q,j_1',\ldots,j_{k-1}')}
$$

\noindent
in {\rm (\ref{new1901}),} containing the expressions 
$\phi_{j_{k}'}(\theta)\phi_{j_{k-1}'}(\theta),\ldots,
\phi_{j_{k}'}(\theta)\phi_{j_{1}'}(\theta),$
should be understood in a special way.
Let us explain this rule on the basis of the sum
$$
\sum\limits_{(j_1,\ldots,j_q,j_1',\ldots,j_{k-1}')}
\int\limits_t^T \phi_{j_q}(t_q)\ldots \int\limits_t^{t_2}\phi_{j_1}(t_1)
\int\limits_t^{t_1}\phi_{j_{k}'}(\theta)\phi_{j_{k-1}'}(\theta)\int\limits_t^{\theta}
\phi_{j_{k-2}'}(t_{k-2}')\ldots
$$
\begin{equation}
\label{new777100}
\Biggl.\ldots \int\limits_t^{t_2'}\phi_{j_1'}(t_1')
d{\bf w}_{t_1'}^{(1)}\ldots 
d{\bf w}_{t_{k-2}'}^{(1)}d{\bf w}_{\theta}^{(0)}d{\bf w}_{t_1}^{(i_1)}\ldots d{\bf w}_{t_q}^{(i_q)}.
\end{equation}

\vspace{1mm}

More precisely, permutations $\left(j_1,\ldots,j_q,j_1',\ldots,j_{k-1}'\right)$ 
when summing in {\rm (\ref{new777100})}
are performed in such a way that if
$j_r^{*}$ swapped with $j_d^{*}$ in the  
permutation $\left(j_{q+k-1}^{*},\ldots,j_1^{*}\right)=
\left(j_q,\ldots,j_1,j_{k-1}',j_{k-2}',\ldots,j_{1}'\right),$ 
then $i_r^{*}$ swapped with $i_d^{*}$ in 
the permutation 
$$
\left(i_{q+k-1}^{*},\ldots,i_1^{*}\right)=
\bigl(i_q,\ldots,i_1,0,\underbrace{1, \ldots ,1}_{k-2}\bigr).
$$

\noindent
Moreover, 
$\bar \phi_{j_r^{*}}$ swapped with $\bar \phi_{j_d^{*}}$
in the permutation 
$$
\bigl(\bar \phi_{j_{q+k-1}^{*}},\ldots,\bar\phi_{j_1^{*}}\bigr)=
\bigl(\phi_{j_q},\ldots,\phi_{j_1},\hspace{1.5mm} \phi_{j_{k}'}\hspace{-0.5mm}\cdot\hspace{-0.5mm}
\phi_{j_{k-1}'},\hspace{1.5mm}
\phi_{j_{k-2}'},\ldots, \phi_{j_{1}'}\bigr).
$$

\vspace{1mm}
\noindent
A similar rule should be applied to all other sums with respect to permutations
$$
\sum\limits_{(j_1,\ldots,j_q,j_1',\ldots,j_{k-1}')}
$$

\noindent
in {\rm (\ref{new1901})} that contain the expressions
$\phi_{j_{k}'}(\theta)\phi_{j_{k-2}'}(\theta),\ldots,
\phi_{j_{k}'}(\theta)\phi_{j_{1}'}(\theta).$}

\vspace{2mm}

Let us prove the equality (\ref{new1600a}). Consider the case  $n=1.$
By analogy with (\ref{new1700}) and (\ref{new1700a}) we obtain 
$$
\int\limits_t^{s}\phi_{j_1'}(\theta)d{\bf w}_{\theta}^{(0)}
\sum\limits_{(j_1,\ldots,j_q)}\int\limits_t^s \phi_{j_q}(t_q)\ldots
\int\limits_t^{t_2}\phi_{j_1}(t_1)
d{\bf w}_{t_1}^{(i_1)}\ldots {\bf w}_{t_q}^{(i_q)}=
$$
$$
=
\sum\limits_{(j_1,\ldots,j_q, j_1')}
J_{(j_q\ldots j_1 j_1')s,t}^{(i_q\ldots i_1 0)}
$$

\vspace{1mm}
\noindent
w.~p.~1, where $J_{(j_r\ldots j_1)s,t}^{(i_r\ldots i_1)}$
is defined by (\ref{new100001}).
The equality (\ref{new1600a}) is proved for the case $n=1.$

Let us assume that the equality (\ref{new1600a}) is true for $n=2, 3, \ldots, k-1$, and prove
its validity for $n=k.$

In complete analogy with (\ref{new1301}) we get
$$
\int\limits_t^{s}\phi_{j_k'}(\theta)d\theta
\int\limits_t^s \phi_{j_{k-1}'}(t_{k-1})\ldots
\int\limits_t^{t_2}\phi_{j_1'}(t_1)
dt_1\ldots dt_{k-1}=
$$
\begin{equation}
\label{new5000}
~~~~~~=
J_{(j_k'j_{k-1}'\ldots j_1')s,t}^{(0\ldots  0)}
+J_{(j_{k-1}'j_k' j_{k-2}'\ldots j_1')s,t}^{(0\ldots 0)}+
\ldots + J_{(j_{k-1}'\ldots j_1' j_k')s,t}^{(0\ldots 0)}.
\end{equation}

\vspace{5mm}

Applying (\ref{new5000}), we have
$$
\sum\limits_{(j_1',\ldots,j_k')}
\int\limits_t^T \phi_{j_k'}(t_k')\ldots \int\limits_t^{t_2'}\phi_{j_1'}(t_1')
d{\bf w}_{t_1'}^{(0)}\ldots d{\bf w}_{t_k'}^{(0)}=
$$
$$
=\sum\limits_{(j_1',\ldots,j_{k-1}')}
\left(J_{(j_k'j_{k-1}'\ldots j_1')s,t}^{(0\ldots  0)}
+J_{(j_{k-1}'j_k' j_{k-2}'\ldots j_1')s,t}^{(0\ldots 0)}+
\ldots + J_{(j_{k-1}'\ldots j_1' j_k')s,t}^{(0\ldots 0)}\right)=
$$
\begin{equation}
\label{new3005}
~~~~~~~~=\int\limits_t^{T}\phi_{j_k'}(\theta)d\theta\sum\limits_{(j_1',\ldots,j_{k-1}')}
\int\limits_t^T \phi_{j_{k-1}'}(t_{k-1})\ldots \int\limits_t^{t_2'}\phi_{j_1'}(t_1)
d{\bf w}_{t_1}^{(0)}\ldots d{\bf w}_{t_{k-1}}^{(0)}.
\end{equation}

\vspace{3mm}

Using (\ref{new3005}) and the induction hypothesis, we obtain w.~p.~1
$$
\sum\limits_{(j_1',\ldots,j_k')}
\int\limits_t^T 
\phi_{j_k'}(t_k)\ldots \int\limits_t^{t_2}\phi_{j_1'}(t_1)d{\bf w}_{t_1}^{(0)}\ldots 
d{\bf w}_{t_k}^{(0)}\times
$$

\vspace{-2mm}
$$
\times\sum\limits_{(j_1,\ldots,j_q)}\int\limits_t^T \phi_{j_q}(t_q)\ldots
\int\limits_t^{t_2}\phi_{j_1}(t_1)
d{\bf w}_{t_1}^{(i_1)}\ldots d{\bf w}_{t_q}^{(i_q)}=
$$

\vspace{-3mm}
$$
=\int\limits_t^{T}\phi_{j_k'}(\theta)d\theta\sum\limits_{(j_1',\ldots,j_{k-1}')}
\int\limits_t^T \phi_{j_{k-1}'}(t_{k-1}')\ldots \int\limits_t^{t_2'}\phi_{j_1'}(t_1')
d{\bf w}_{t_1'}^{(0)}\ldots d{\bf w}_{t_{k-1}'}^{(0)}\times
$$
$$
\times\sum\limits_{(j_1,\ldots,j_q)}\int\limits_t^T \phi_{j_q}(t_q)\ldots
\int\limits_t^{t_2}\phi_{j_1}(t_1)
d{\bf w}_{t_1}^{(i_1)}\ldots d{\bf w}_{t_q}^{(i_q)}=
$$

$$
=\int\limits_t^{T}\phi_{j_k'}(\theta)d\theta
\sum\limits_{(j_1,\ldots,j_q,j_1',\ldots,j_{k-1}')}
\int\limits_t^T \phi_{j_q}(t_q)\ldots \int\limits_t^{t_2}\phi_{j_1}(t_1)\times
$$

\vspace{1mm}
$$
\times
\int\limits_t^{t_1}\phi_{j_{k-1}'}(t_{k-1}')\ldots \int\limits_t^{t_2'}
\phi_{j_1'}(t_1')d{\bf w}_{t_1'}^{(0)}\ldots d{\bf w}_{t_{k-1}'}^{(0)}
d{\bf w}_{t_1}^{(i_1)}\ldots d{\bf w}_{t_q}^{(i_q)}=
$$

$$
=
\sum\limits_{(j_1,\ldots,j_q,j_1',\ldots,j_{k-1}')}\int\limits_t^{T}\phi_{j_k'}(\theta)d\theta
\int\limits_t^T \phi_{j_q}(t_q)\ldots \int\limits_t^{t_2}\phi_{j_1}(t_1)\times
$$

\vspace{-2mm}
\begin{equation}
\label{new4000}
~~~~~~~~~\times
\int\limits_t^{t_1}\phi_{j_{k-1}'}(t_{k-1}')\ldots \int\limits_t^{t_2'}
\phi_{j_1'}(t_1')d{\bf w}_{t_1'}^{(0)}\ldots d{\bf w}_{t_{k-1}'}^{(0)}
d{\bf w}_{t_1}^{(i_1)}\ldots d{\bf w}_{t_q}^{(i_q)}.
\end{equation}

\vspace{4mm}

An iterative application of the It\^{o} formula leads to the following equality

\vspace{-2mm}
$$
\int\limits_t^{T}\phi_{j_k'}(\theta)d\theta
\int\limits_t^T \phi_{j_q}(t_q)\ldots \int\limits_t^{t_2}\phi_{j_1}(t_1)\times
$$
$$
\times
\int\limits_t^{t_1}\phi_{j_{k-1}'}(t_{k-1}')\ldots \int\limits_t^{t_2'}
\phi_{j_1'}(t_1')d{\bf w}_{t_1'}^{(0)}\ldots d{\bf w}_{t_{k-1}'}^{(0)}
d{\bf w}_{t_1}^{(i_1)}\ldots d{\bf w}_{t_q}^{(i_q)}=
$$

$$
=J_{(j_k'j_q \ldots j_1 j_{k-1}'\ldots j_1')T,t}^{(0 i_q\ldots i_1 0\ldots 0)}+
J_{(j_q j_k'j_{q-1} \ldots j_1 j_{k-1}'\ldots j_1')T,t}^{(i_q 0 i_{q-1}\ldots i_1 0\ldots 0)}+\ldots
J_{(j_q \ldots j_1 j_k' j_{k-1}'\ldots j_1')T,t}^{(i_q\ldots i_1 0\ldots 0)}+
$$

\vspace{1mm}
\begin{equation}
\label{new4001}
+J_{(j_q \ldots j_1 j_{k-1}' j_k' j_{k-2}'\ldots j_1')T,t}^{(i_q\ldots i_1 0\ldots 0)}+\ldots
+J_{(j_q \ldots j_1 j_{k-1}' \ldots j_1' j_k')T,t}^{(i_q\ldots i_1 0\ldots 0)}
\end{equation}

\vspace{5mm}
\noindent
w.~p.~1.

Combining (\ref{new4000}) and (\ref{new4001}), we finally obtain w.~p.~1
$$
\sum\limits_{(j_1,\ldots,j_q)}\int\limits_t^T \phi_{j_q}(t_q)\ldots
\int\limits_t^{t_2}\phi_{j_1}(t_1)
d{\bf w}_{t_1}^{(i_1)}\ldots d{\bf w}_{t_q}^{(i_q)}\times
$$

\vspace{-2mm}
$$
\times \sum\limits_{(j_1',\ldots,j_k')}
\int\limits_t^T \phi_{j_k'}(t_k')\ldots \int\limits_t^{t_2'}\phi_{j_1'}(t_1')
d{\bf w}_{t_1'}^{(0)}\ldots d{\bf w}_{t_k'}^{(0)}=
$$

\vspace{-1mm}
$$
=\sum\limits_{(j_1,\ldots,j_q,j_1',\ldots,j_k')}
\int\limits_t^T \phi_{j_q}(t_q)\ldots \int\limits_t^{t_2}\phi_{j_1}(t_1)
\int\limits_t^{t_1}\phi_{j_k'}(t_k')\ldots \int\limits_t^{t_2'}
\phi_{j_1'}(t_1')\times
$$

\vspace{1mm}
$$
\times d{\bf w}_{t_1'}^{(0)}\ldots d{\bf w}_{t_k'}^{(0)}d{\bf w}_{t_1}^{(i_1)}\ldots d{\bf w}_{t_q}^{(i_q)}.
$$

\vspace{4mm}

The equality (\ref{new1600a}) is proved for $n=k.$
The equality (\ref{new1600a}) is proved. Theorem~1.22 is proved.

To complete the proof of Theorems~1.16 and 1.17, we prove the following theorem.

{\bf Theorem 1.23.}\ {\it Suppose that
$\{\phi_j(x)\}_{j=0}^{\infty}$ is an arbitrary complete orthonormal system  
of functions in the space $L_2([t,T]).$
Then the following representation
$$
J''[\phi_{j_1}\ldots\phi_{j_k}]^{(i_1 \ldots i_k)}=
\prod_{l=1}^k\zeta_{j_l}^{(i_l)}+\sum\limits_{r=1}^{[k/2]}
(-1)^r \times \Biggr.
$$

\begin{equation}
\label{leto6000xxa}
\times
\sum_{\stackrel{(\{\{g_1, g_2\}, \ldots, 
\{g_{2r-1}, g_{2r}\}\}, \{q_1, \ldots, q_{k-2r}\})}
{{}_{\{g_1, g_2, \ldots, 
g_{2r-1}, g_{2r}, q_1, \ldots, q_{k-2r}\}=\{1,2, \ldots, k\}}}}
\prod\limits_{s=1}^r
{\bf 1}_{\{i_{g_{{}_{2s-1}}}=~i_{g_{{}_{2s}}}\ne 0\}}
\Biggl.{\bf 1}_{\{j_{g_{{}_{2s-1}}}=~j_{g_{{}_{2s}}}\}}
\prod_{l=1}^{k-2r}\zeta_{j_{q_l}}^{(i_{q_l})}
\end{equation}

\vspace{2mm}
\noindent 
is valid w.~p.~{\rm 1,} where $i_1,\ldots,i_k=0,1,\ldots,m,$
$[x]$ is an integer part of a real number $x,$
$\prod\limits_{\emptyset}\stackrel{\sf def}{=}1,$ $\sum\limits_{\emptyset}
\stackrel{\sf def}{=}0;$
the sum in the second line of the formula {\rm (\ref{leto6000xxa})} 
is the sum with respect to all possible
partitions {\rm (\ref{leto5008});} 
another notations are the same as in Theorems~{\rm 1.1, 1.2.}}

\vspace{1mm}

{\bf Remark~1.17.}\ {\it It should be noted that the formulas {\rm (\ref{new1010}),}
{\rm (\ref{new100000}),} {\rm (\ref{new1601}),} {\rm (\ref{new1601a})}
follow from {\rm (\ref{leto6000xxa}).}
It is only necessary to set the values
of the corresponding indicators of the form ${\bf 1}_A$ from the formula
{\rm (\ref{leto6000xxa})} equal to $0$ or $1.$}

{\bf Proof.}\ The proof of Theorem~1.23 is carried out
by induction using the following recurrence 
relation

\vspace{-1mm}
$$
J''[\phi_{j_1}\ldots\phi_{j_k}]^{(i_1 \ldots i_k)}_{T,t}=
J''[\phi_{j_k}]^{(i_k)}_{T,t}\cdot
J''[\phi_{j_1}\ldots\phi_{j_{k-1}}]^{(i_1 \ldots i_{k-1})}_{T,t}-
$$

\vspace{-3mm}
\begin{equation}
\label{recur1}
~~~~~~ -\sum\limits_{l=1}^{k-1}{\bf 1}_{\{i_l=i_k\ne 0\}}
{\bf 1}_{\{j_l=j_k\}}\cdot 
J''[\phi_{j_1}\ldots\phi_{j_{l-1}}\phi_{j_{l+1}}\ldots\phi_{j_{k-1}}]^{(i_1
\ldots  i_{l-1}i_{l+1}\ldots i_{k-1})}_{T,t}
\end{equation}

\vspace{2mm}
\noindent
w.~p.~1.

Let us prove the recurrence relation (\ref{recur1}).
Using iteratively the It\^{o} formula, the orthonormality of $\{\phi_j(x)\}_{j=0}^{\infty}$,
as well as (\ref{new5001}) and 
combinatorial reasoning, we obtain w.~p.~1 (see Remark~1.18 below for details)

\vspace{1mm}
$$
J''[\phi_{j_k}]^{(i_k)}_{T,t}\cdot
J''[\phi_{j_1}\ldots\phi_{j_{k-1}}]^{(i_1 \ldots i_{k-1})}_{T,t}=
$$

\vspace{-3mm}
$$
=\int\limits_t^{T}\phi_{j_k}(\theta)d{\bf w}_{\theta}^{(i_k)}
\sum\limits_{(j_1,\ldots,j_{k-1})}
\int\limits_t^T \phi_{j_{k-1}}(t_{k-1})\ldots \int\limits_t^{t_2}\phi_{j_1}(t_1)
d{\bf w}_{t_1}^{(i_1)}\ldots d{\bf w}_{t_{k-1}}^{(i_{k-1})}=
$$

\vspace{-2mm}
$$
=\sum\limits_{(j_1,\ldots,j_{k-1})}\int\limits_t^{T}\phi_{j_k}(\theta)d{\bf w}_{\theta}^{(i_k)}
\int\limits_t^T \phi_{j_{k-1}}(t_{k-1})\ldots \int\limits_t^{t_2}\phi_{j_1}(t_1)
d{\bf w}_{t_1}^{(i_1)}\ldots d{\bf w}_{t_{k-1}}^{(i_{k-1})}=
$$

\vspace{-2mm}
$$
=\sum\limits_{(j_1,\ldots,j_{k})}
\int\limits_t^T \phi_{j_{k}}(t_{k})\ldots \int\limits_t^{t_2}\phi_{j_1}(t_1)
d{\bf w}_{t_1}^{(i_1)}\ldots d{\bf w}_{t_{k}}^{(i_{k})}+
$$

\vspace{-2mm}
$$
+\sum\limits_{(j_1,\ldots,j_{k-1})}
\Biggl({\bf 1}_{\{i_k=i_{k-1}\ne 0\}}
\int\limits_t^T \phi_{j_{k}}(\theta)\phi_{j_{k-1}}(\theta)
\int\limits_t^{\theta}
\phi_{j_{k-2}}(t_{k-2})
\ldots \int\limits_t^{t_2}\phi_{j_1}(t_1)\times\Biggr.
$$

$$
\times
d{\bf w}_{t_1}^{(i_1)}\ldots d{\bf w}_{t_{k-2}}^{(i_{k-2})}d{\bf w}_{\theta}^{(0)}+
$$

\vspace{-2mm}
$$
+{\bf 1}_{\{i_k=i_{k-2}\ne 0\}}
\int\limits_t^T \phi_{j_{k-1}}(t_{k-1})
\int\limits_t^{t_{k-1}} \phi_{j_{k}}(\theta)\phi_{j_{k-2}}(\theta)
\int\limits_t^{\theta}
\phi_{j_{k-3}}(t_{k-3})
\ldots \int\limits_t^{t_2}\phi_{j_1}(t_1)\times\Biggr.
$$
$$
\times
d{\bf w}_{t_1}^{(i_1)}\ldots d{\bf w}_{t_{k-3}}^{(i_{k-3})}d{\bf w}_{\theta}^{(0)}
d{\bf w}_{t_{k-1}}^{(i_{k-1})}+ \ldots
$$

\vspace{-2mm}
$$
\ldots +{\bf 1}_{\{i_k=i_{1}\ne 0\}}
\int\limits_t^T \phi_{j_{k-1}}(t_{k-1})\ldots 
\int\limits_t^{t_3}
\phi_{j_{2}}(t_{2})
\int\limits_t^{t_{2}} \phi_{j_{k}}(\theta)\phi_{j_{1}}(\theta)\times
$$

$$
\Biggl.\times
d{\bf w}_{\theta}^{(0)}d{\bf w}_{t_2}^{(i_2)}\ldots d{\bf w}_{t_{k-1}}^{(i_{k-1})}
\Biggr)=
$$

\vspace{-1mm}
$$
=\sum\limits_{(j_1,\ldots,j_{k})}
\int\limits_t^T \phi_{j_{k}}(t_{k})\ldots \int\limits_t^{t_2}\phi_{j_1}(t_1)
d{\bf w}_{t_1}^{(i_1)}\ldots d{\bf w}_{t_{k}}^{(i_{k})}+
$$

\vspace{-2mm}
$$
+\sum\limits_{(j_1,\ldots,j_{k-2})}
{\bf 1}_{\{i_k=i_{k-1}\ne 0\}}
\Biggl\{\int\limits_t^T \phi_{j_{k}}(\theta)\phi_{j_{k-1}}(\theta)
\int\limits_t^{\theta}
\phi_{j_{k-2}}(t_{k-2})
\ldots \int\limits_t^{t_2}\phi_{j_1}(t_1)\times\Biggr.
$$

$$
\times
d{\bf w}_{t_1}^{(i_1)}\ldots d{\bf w}_{t_{k-2}}^{(i_{k-2})}d{\bf w}_{\theta}^{(0)}+\ldots
$$

\vspace{-2mm}
$$
\Biggl.\ldots +\int\limits_t^T 
\phi_{j_{k-2}}(t_{k-2})
\ldots \int\limits_t^{t_2}\phi_{j_1}(t_1)
\int\limits_t^{t_1}
\phi_{j_{k}}(\theta)\phi_{j_{k-1}}(\theta)
d{\bf w}_{\theta}^{(0)}
d{\bf w}_{t_1}^{(i_1)}\ldots d{\bf w}_{t_{k-2}}^{(i_{k-2})}\Biggr\}+
$$

\vspace{-2mm}
$$
+\sum\limits_{(j_1,\ldots,j_{k-3},j_{k-1})}
{\bf 1}_{\{i_k=i_{k-2}\ne 0\}}
\Biggl\{\int\limits_t^T \phi_{j_{k}}(\theta)\phi_{j_{k-2}}(\theta)
\int\limits_t^{\theta}
\phi_{j_{k-1}}(t_{k-1})
\int\limits_t^{t_{k-1}}
\phi_{j_{k-3}}(t_{k-3})
\ldots \Biggr.
$$

$$
\ldots \int\limits_t^{t_2}\phi_{j_1}(t_1)
d{\bf w}_{t_1}^{(i_1)}\ldots d{\bf w}_{t_{k-3}}^{(i_{k-3})}d{\bf w}_{t_{k-1}}^{(i_{k-1})}
d{\bf w}_{\theta}^{(0)}+\ldots
$$

\vspace{-2mm}
$$
\ldots +\int\limits_t^T 
\phi_{j_{k-1}}(t_{k-1})
\int\limits_t^{t_{k-1}}
\phi_{j_{k-3}}(t_{k-3})
\ldots \int\limits_t^{t_2}\phi_{j_1}(t_1)
\int\limits_t^{t_1}
\phi_{j_{k}}(\theta)\phi_{j_{k-2}}(\theta)\times
$$

$$
\Biggl.\times d{\bf w}_{\theta}^{(0)}
d{\bf w}_{t_1}^{(i_1)}\ldots d{\bf w}_{t_{k-3}}^{(i_{k-3})}d{\bf w}_{t_{k-1}}^{(i_{k-1})}\Biggr\}+\ldots
$$
$$
\ldots +\sum\limits_{(j_2,\ldots,j_{k-1})}
{\bf 1}_{\{i_k=i_{1}\ne 0\}}
\Biggl\{\int\limits_t^T \phi_{j_{k}}(\theta)\phi_{j_{1}}(\theta)
\int\limits_t^{\theta}
\phi_{j_{k-1}}(t_{k-1})
\ldots \int\limits_t^{t_3}\phi_{j_2}(t_2)\times\Biggr.
$$

$$
\times
d{\bf w}_{t_2}^{(i_2)}\ldots d{\bf w}_{t_{k-1}}^{(i_{k-1})}d{\bf w}_{\theta}^{(0)}+\ldots
$$

\vspace{-2mm}
$$
\Biggl.\ldots +\int\limits_t^T 
\phi_{j_{k-1}}(t_{k-1})
\ldots \int\limits_t^{t_3}\phi_{j_2}(t_2)
\int\limits_t^{t_2}
\phi_{j_{k}}(\theta)\phi_{j_{1}}(\theta)
d{\bf w}_{\theta}^{(0)}
d{\bf w}_{t_2}^{(i_2)}\ldots d{\bf w}_{t_{k-1}}^{(i_{k-1})}\Biggr\}=
$$

\vspace{-2mm}
$$
=\sum\limits_{(j_1,\ldots,j_{k})}
\int\limits_t^T \phi_{j_{k}}(t_{k})\ldots \int\limits_t^{t_2}\phi_{j_1}(t_1)
d{\bf w}_{t_1}^{(i_1)}\ldots d{\bf w}_{t_{k}}^{(i_{k})}+
$$

\vspace{-2mm}
$$
+\int\limits_t^T \phi_{j_{k}}(\theta)\phi_{j_{k-1}}(\theta)d\theta
\sum\limits_{(j_1,\ldots,j_{k-2})}
{\bf 1}_{\{i_k=i_{k-1}\ne 0\}}
\int\limits_t^{T}
\phi_{j_{k-2}}(t_{k-2})
\ldots \int\limits_t^{t_2}\phi_{j_1}(t_1)\times\Biggr.
$$

$$
\times
d{\bf w}_{t_1}^{(i_1)}\ldots d{\bf w}_{t_{k-2}}^{(i_{k-2})}+
$$

\vspace{-4mm}
$$
+\int\limits_t^T \phi_{j_{k}}(\theta)\phi_{j_{k-2}}(\theta)d\theta
\sum\limits_{(j_1,\ldots,j_{k-3},j_{k-1})}
{\bf 1}_{\{i_k=i_{k-2}\ne 0\}}
\int\limits_t^{T}
\phi_{j_{k-1}}(t_{k-1})
\int\limits_t^{t_{k-1}}
\phi_{j_{k-3}}(t_{k-3})
\ldots 
$$

$$
\ldots \int\limits_t^{t_2}\phi_{j_1}(t_1)
d{\bf w}_{t_1}^{(i_1)}\ldots d{\bf w}_{t_{k-3}}^{(i_{k-3})}d{\bf w}_{t_{k-1}}^{(i_{k-1})}
+\ldots
$$

\vspace{-2mm}

$$
\ldots +\int\limits_t^T \phi_{j_{k}}(\theta)\phi_{j_{1}}(\theta)d\theta
\sum\limits_{(j_2,\ldots,j_{k-1})}
{\bf 1}_{\{i_k=i_{1}\ne 0\}}
\int\limits_t^{T}
\phi_{j_{k-1}}(t_{k-1})
\ldots \int\limits_t^{t_3}\phi_{j_2}(t_2)\times
$$

$$
\times
d{\bf w}_{t_2}^{(i_2)}\ldots d{\bf w}_{t_{k-1}}^{(i_{k-1})}=
$$

\vspace{-1mm}

$$
=J''[\phi_{j_1}\ldots\phi_{j_k}]^{(i_1 \ldots i_k)}_{T,t}+
{\bf 1}_{\{i_k=i_{k-1}\ne 0\}}{\bf 1}_{\{j_k=j_{k-1}\}}
\cdot J''[\phi_{j_1}\ldots \phi_{j_{k-2}}]_{T,t}^{(i_1\ldots i_{k-2})}+
$$

\vspace{-1mm}
$$
+
{\bf 1}_{\{i_k=i_{k-2}\ne 0\}}{\bf 1}_{\{j_k=j_{k-2}\}}
\cdot J''[\phi_{j_1}\ldots \phi_{j_{k-3}}\phi_{j_{k-1}}]_{T,t}^{(i_1\ldots i_{k-3}i_{k-1})}+\ldots
$$
$$
\ldots +
{\bf 1}_{\{i_k=i_{1}\ne 0\}}{\bf 1}_{\{j_k=j_{1}\}}
\cdot J''[\phi_{j_2}\ldots \phi_{j_{k-1}}]_{T,t}^{(i_2\ldots i_{k-1})}=
$$

\vspace{1mm}
$$
=J''[\phi_{j_1}\ldots\phi_{j_k}]^{(i_1 \ldots i_k)}_{T,t}+
$$

\vspace{-6mm}
\begin{equation}
\label{new00002}
~~~~~~+\sum\limits_{l=1}^{k-1}{\bf 1}_{\{i_l=i_k\ne 0\}}
{\bf 1}_{\{j_l=j_k\}}\cdot 
J''[\phi_{j_1}\ldots\phi_{j_{l-1}}\phi_{j_{l+1}}\ldots\phi_{j_{k-1}}]^{(i_1
\ldots  i_{l-1}i_{l+1}\ldots i_{k-1})}_{T,t}.
\end{equation}

\vspace{2mm}

The equality (\ref{recur1}) is proved. Theorem~1.23 is proved.

\vspace{2mm}

{\bf Remark~1.18.}\ {\it It should be noted that the sums with respect to 
permutations 
$$
\sum\limits_{(j_1,\ldots,j_{k-1})}
$$

\noindent
in {\rm (\ref{new00002}),} containing the expressions 

\vspace{-2mm}
$$
{\bf 1}_{\{i_k=i_{k-1}\ne 0\}}\phi_{j_{k}}(\theta)\phi_{j_{k-1}}(\theta),\ldots,
{\bf 1}_{\{i_k=i_{1}\ne 0\}}\phi_{j_{k}}(\theta)\phi_{j_{1}}(\theta),
$$

\vspace{2mm}
\noindent
should be understood in a special way.
Let us explain this rule on the basis of the sum
$$
\sum\limits_{(j_1,\ldots,j_{k-1})}
{\bf 1}_{\{i_k=i_{k-1}\ne 0\}}
\int\limits_t^T \phi_{j_{k}}(\theta)\phi_{j_{k-1}}(\theta)
\int\limits_t^{\theta}
\phi_{j_{k-2}}(t_{k-2})
\ldots \int\limits_t^{t_2}\phi_{j_1}(t_1)\times
$$

\begin{equation}
\label{new00003}
\times
d{\bf w}_{t_1}^{(i_1)}\ldots d{\bf w}_{t_{k-2}}^{(i_{k-2})}d{\bf w}_{\theta}^{(0)}.
\end{equation}

\vspace{2mm}

More precisely, permutations $(j_1,\ldots,j_{k-1})$ 
when summing in {\rm (\ref{new00003})}
are performed in such a way that if
$j_r$ swapped with $j_d$ in the  
permutation $(j_1,\ldots,j_{k-1}),$ 
then $i_r$ swapped with $i_d$ in 
the permutation $(i_1,\ldots,i_{k-2}, i_{k-1})$\ {\rm (}note that $i_{k-1}=0${\rm )}.
Moreover, 
$\bar \phi_{j_r}$ swapped with $\bar \phi_{j_d}$
in the permutation 

\vspace{-2mm}
$$
\bigl(\bar \phi_{j_{1}},\ldots,\bar\phi_{j_{k-1}}\bigr)=
\bigl(\phi_{j_1},\ldots,\phi_{j_{k-2}},\hspace{1.5mm} 
{\bf 1}_{\{i_k=i_{k-1}\ne 0\}}\cdot \phi_{j_k}\cdot \phi_{j_{k-1}}\bigr),
$$

\vspace{2mm}
\noindent
where $\bar\phi_{j_{k-1}}(\tau)=
{\bf 1}_{\{i_k=i_{k-1}\ne 0\}}\phi_{j_k}(\tau)\phi_{j_{k-1}}(\tau).$

A similar rule should be applied to all other sums with respect to permutations
$$
\sum\limits_{(j_1,\ldots,j_{k-1})}
$$

\noindent
in {\rm (\ref{new00002})} that contain the expressions
$$
{\bf 1}_{\{i_k=i_{k-2}\ne 0\}}\phi_{j_{k}}(\theta)\phi_{j_{k-2}}(\theta),\ldots,
{\bf 1}_{\{i_k=i_{1}\ne 0\}}\phi_{j_{k}}(\theta)\phi_{j_{1}}(\theta).
$$
}

\vspace{-2mm}

The relations (\ref{chain4002}), (\ref{new6000}), (\ref{leto6000xxa})
prove Theorem~1.16. An analogue of the formula (\ref{chain4002})
for $\Phi(t_1,\ldots,t_k)$ instead of $K(t_1,\ldots,t_k)$ and
(\ref{new6000}), (\ref{leto6000xxa}) prove Theorem~1.17.

We note a number of works \cite{fox}-\cite{major2} in which the properties
of multiple Wiener stochastic integrals were studied using
measure theory, in particular, the formulas for the product
of such integrals were obtained.

First of all, let us compare Theorem~1.23 with Proposition~5.1 from \cite{fox}.
An analogue of the right-hand side of (\ref{leto6000xxa})
for nonrandom $x_1,\ldots,x_k$
is constructed in \cite{fox} using diagrams (see the formula (5.1) in \cite{fox}).
This means that the application of the formula (5.1) from \cite{fox},
unlike the formula (\ref{leto6000xxa}), is difficult when
performing algebraic transformations.

Further, we note that the formula (5.1) from \cite{fox}
was applied to the representation of the multiple Wiener stochastic integral
somewhat differently than the formula (\ref{leto6000xxa}).
Namely, using Proposition~5.1 \cite{fox}.
Let us expain this difference in more detail.

Proposition~5.1 from \cite{fox} in our degree of generality 
and in our notations can be written as

\vspace{-2mm}
$$
J''\left[\phi_{j_1}\ldots \phi_{j_k}\right]_{T,t}^{(i_1\ldots i_k)}=
$$

\vspace{-2mm}
$$
=
J''\biggl[\underbrace{\phi_{j_1}
\ldots \phi_{j_1}}_{m_1}
\underbrace{\phi_{j_2}
\ldots \phi_{j_2}}_{m_2}\ldots 
\underbrace{\phi_{j_p}\ldots
\phi_{j_p}}_{m_p}\biggr]_{T,t}^
{(\overbrace{\hspace{0.5mm}{}_{i_1 \ldots i_{m_1}}}^{m_1}
\overbrace{\hspace{0.3mm}{}_{i_{m_1+1} \ldots i_{m_2}}}^{m_2}
\ldots \overbrace{\hspace{0.3mm}{}_{i_{m_1+\ldots +m_{p-1}+1} \ldots i_k}}^{m_p})}=
$$

\vspace{-2mm}
\begin{equation}
\label{new54321}
=J''\left[\phi_{j_1}
\hspace{-0.3mm}\ldots \hspace{-0.3mm}\phi_{j_1}\right]
_{T,t}^
{(\overbrace{\hspace{0.5mm}{}_{i_1 \ldots i_{m_1}}}^{m_1})}\hspace{-0.3mm}\cdot
J''\left[\phi_{j_2}\hspace{-0.3mm}
\ldots \hspace{-0.3mm}\phi_{j_2}\right]_{T,t}^
{(\overbrace{\hspace{0.3mm}{}_{i_{m_1+1} \ldots i_{m_2}}}^{m_2})}\hspace{-0.3mm}\cdot \ldots
\cdot J''\left[\phi_{j_p}
\hspace{-0.3mm}\ldots \hspace{-0.3mm}\phi_{j_p}\right]_{T,t}^
{(\overbrace{\hspace{0.3mm}{}_{i_{m_1+\ldots +m_{p-1}+1} \ldots i_k}}^{m_p})}
\end{equation}

\vspace{5mm}
\noindent
w.~p.~1, where

\vspace{-2mm}
$$
J''\left[\phi_{j_1}
\hspace{-0.3mm}\ldots \hspace{-0.3mm}\phi_{j_1}\right]
_{T,t}^
{(\overbrace{\hspace{0.5mm}{}_{i_1 \ldots i_{m_1}}}^{m_1})}\hspace{-1.5mm},
J''\left[\phi_{j_2}\hspace{-0.3mm}
\ldots \hspace{-0.3mm}\phi_{j_2}\right]_{T,t}^
{(\overbrace{\hspace{0.3mm}{}_{i_{m_1+1} \ldots i_{m_2}}}^{m_2})}\hspace{-1mm},\ldots,
J''\left[\phi_{j_p}
\hspace{-0.3mm}\ldots \hspace{-0.3mm}\phi_{j_p}\right]_{T,t}^
{(\overbrace{\hspace{0.3mm}{}_{i_{m_1+\ldots +m_{p-1}+1} \ldots i_k}}^{m_p})}
$$

\vspace{5mm}
\noindent
are defined by the right-hand side of the formula (5.1) from \cite{fox},
$m_1+\ldots +m_p=k,$ $m_1,\ldots, m_p>0,$ $j_q\ne j_d$ $(q\ne d,\ q,d=1,\ldots,p),$
$i_1,\ldots,i_k=1,\ldots,m.$

This actually means that in \cite{fox} an analogue of the formula
(\ref{leto6000xxa}) is constructed for the special case
$j_1=\ldots=j_k$. Moreover, the specified analogue 
is based on the formula (5.1) \cite{fox} obtained using diagrams.

Comparing the formulas (\ref{leto6000xxa}) and (\ref{new54321}) (or (5.1) from \cite{fox}), it is easy
to understand that the transition from 
(\ref{leto6000xxa}) to (\ref{new54321}) is obvious.
It is only necessary to set the values
of the corresponding indicators of the form ${\bf 1}_A$ from the formula
(\ref{leto6000xxa}) equal to $0$ or $1.$
The reverse transition from the formula (\ref{new54321})
to the formula (\ref{leto6000xxa}) is not obvious.
Note that the formula 
(\ref{leto6000xxa}) (not the formula (\ref{new54321})) is convenient for the  numerical
integration of It\^{o} stochastic differential equations (see Chapter~5 of this book for details).

Let us turn to the comparison of Theorem~1.23 with another interesting work \cite{major2} (2019).
As it turned out, a version of Theorem~1.23 was obtained in terms 
of Wick polynomials and for the case of vector valued random measures 
in \cite{major2} (see Theorem~7.2, p.~69).
However, much earlier the formula (\ref{leto6000xxa}) (Theorem~1.23) is obtained
in our monograph \cite{4} (2009) as part of the formula
(5.30) (see \cite{4}, p.~220).
Moreover, particular cases of the formula (\ref{leto6000xxa}) were obtained
even earlier in our works \cite{1} (2006) and \cite{3} (2007).
More precisely, partiular cases $k=1,\ldots,5$ of the formula (\ref{leto6000xxa})
were obtained in \cite{1} (2006) as parts of the formulas 
on the pages 243-244 and partiular cases $k=1,\ldots,7$ of the formula (\ref{leto6000xxa})
were obtained in \cite{3} (2007) as parts of the formulas 
on the pages 208-218.

We also note that we have found an explicit expression for the 
Wick polynomial of degree $k$ of the arguments $\zeta_{j_1}^{(i_1)},\ldots,\zeta_{j_k}^{(i_k)}$ 
(see the formula (\ref{leto6000xxa})),
which is very convenient for the numerical simulation of
iterated It\^{o} stochastic integrals (\ref{ito})  \cite{Kuz-Kuz}, \cite{Mikh-1}.
Note that the representation of the Wick polynomial
of the arguments $\zeta_{j_1}^{(i_1)},\ldots,\zeta_{j_k}^{(i_k)}$ 
in terms of the product of Hermite polynomials
is less convenient for the numerical simulation of
iterated It\^{o} stochastic integrals (\ref{ito}).
For example, the expression for $J''[\phi_{j_1}\phi_{j_2}\phi_{j_3}\phi_{j_4}]^{(i_1 i_2 i_3 i_4)}_{T,t}$
in terms of the product of Hermite polynomials,
even under the condition $i_1=i_2=i_3=i_4$, already contains
15 different expressions (see Sect.~1.10).
At the same time, all these 15 expressions are contained 
in one formula (\ref{leto6000xxa}) provided that $k=4$ and $i_1=i_2=i_3=i_4$.
It is very convenient, since in computer simulation
using the formula (\ref{leto6000xxa}), in addition to
modeling of random variables $\zeta_{j_1}^{(i_1)},\ldots,\zeta_{j_k}^{(i_k)}$,
it remains only to set  
the values
of the corresponding indicators of the form ${\bf 1}_A$ from the formula
(\ref{leto6000xxa}) equal to $0$ or $1.$

It should be noted that in \cite{major} (Theorem~6.1)
a diagram formula was obtained for the product
of two multiple Wiener stochastic integrals
with respect to vector valued random measures. 
The formula (\ref{new1600}) can be derived from the diagram formula \cite{major}.
Although the proof of the diagram formula \cite{major}
is much more complicated than our proof of the formula (\ref{new1600}).

To conclude this section, we say a few words about expansions 
(\ref{new9999}) and (\ref{razzar1}).
The transition from the expansion (\ref{razzar1}) to the expansion 
(\ref{new9999}) is obvious. It is only necessary to set the values
of the corresponding indicators of the form ${\bf 1}_A$ from the formula
(\ref{razzar1}) equal to $0$ or $1.$
The reverse transition from the formula (\ref{new9999})
to the formula (\ref{razzar1}) is also possible but not obvious.
However, Theorems~1.22 and 1.23 provide a transition from 
(\ref{new9999}) to (\ref{razzar1}) and vice versa.
Note that the expansion (\ref{new9999}) is interesting from the point of 
view of studying the structure of the expansion of iterated It\^{o}
stochastic integrals. On the orther hand, 
the expansion (\ref{razzar1}) is exceptionally convenient 
for applications (see Chapter 5 of this book and 
\cite{Kuz-Kuz}, \cite{Mikh-1}).

\section{Generalization of Theorem 1.11 to the Case of an Arbitrary 
Complete Ortho\-nor\-mal System of Functions in the Space $L_2([t, T])$
and $\psi_1(\tau),$ $\ldots,\psi_k(\tau)\in L_2([t, T])$}

Suppose that $\psi_1(\tau),$ $\ldots,\psi_k(\tau)\in L_2([t, T])$.
Define the following function on the hypercube $[t, T]^k$

\vspace{-3mm}
$$
\bar K(t_1,\ldots,t_k,s)={\bf 1}_{\{t_k<s\}}K(t_1,\ldots,t_k),
$$

\vspace{1mm}
\noindent
where the function $K(t_1,\ldots,t_k)$ has the form
(\ref{ppp}), $s\in (t, T]$ ($s$ is fixed), 
and ${\bf 1}_A$ is the indicator of the set $A.$

Further, we have (see (\ref{ppp}))

\vspace{-3mm}
$$
\bar K(t_1,\ldots,t_k,s)=
{\bf 1}_{\{t_1<\ldots <t_k<s\}}\psi_1(t_1)\ldots \psi_k(t_k)=
$$

\vspace{-2mm}
$$
=
\left\{\begin{matrix}
\psi_1(t_1)\ldots \psi_k(t_k),\ &t_1<\ldots<t_k<s\cr\cr
0,\ &\hbox{\rm otherwise}
\end{matrix}
\right.,
$$

\vspace{2.5mm}
\noindent
where $\bar K(t_1,\ldots,t_k,s)\in L_2([t,T]^k),$ $k\ge 1, $ $t_1,\ldots,t_k\in [t, T],$ and 
$s\in (t, T]$.

Note that
$$
J[\psi^{(k)}]_{s,t}=\int\limits_t^s\psi_k(t_k) \ldots \int\limits_t^{t_{2}}
\psi_1(t_1) d{\bf w}_{t_1}^{(i_1)}\ldots
d{\bf w}_{t_k}^{(i_k)}=
$$
\begin{equation}
\label{strange700}
~~~~~~~=
\int\limits_t^T {\bf 1}_{\{t_k<s\}}\psi_k(t_k) \ldots \int\limits_t^{t_{2}}
\psi_1(t_1) d{\bf w}_{t_1}^{(i_1)}\ldots
d{\bf w}_{t_k}^{(i_k)}\ \ \ \hbox{w.~p.~1},
\end{equation}

\vspace{2mm}
\noindent
where $s\in (t, T]$ ($s$ is fixed), $i_1,\ldots,i_k=0,1,\ldots,m.$ 

Applying Theorem~1.16 to the iterated It\^{o} stochastic integral
(\ref{strange700}), we obtain the following generalization 
of Theorem~1.11 to the case of an arbitrary 
complete ortho\-nor\-mal system of functions in the space $L_2([t, T])$
and $\psi_1(\tau),$ $\ldots,\psi_k(\tau)\in L_2([t, T]).$

{\bf Theorem~1.24.}\ {\it Suppose that
$\psi_1(\tau),$ $\ldots,\psi_k(\tau)\in L_2([t, T])$ and
$\{\phi_j(x)\}_{j=0}^{\infty}$ is an arbitrary complete orthonormal system  
of functions in the space $L_2([t,T]).$
Then$,$ the following expansion

\vspace{-4mm}
$$
J[\psi^{(k)}]_{s,t}^{(i_1\ldots i_k)}=
\hbox{\vtop{\offinterlineskip\halign{
\hfil#\hfil\cr
{\rm l.i.m.}\cr
$\stackrel{}{{}_{p_1,\ldots,p_k\to \infty}}$\cr
}} }
\sum\limits_{j_1=0}^{p_1}\ldots
\sum\limits_{j_k=0}^{p_k}
C_{j_k\ldots j_1}(s)\Biggl(
\prod_{l=1}^k\zeta_{j_l}^{(i_l)}+\sum\limits_{r=1}^{[k/2]}
(-1)^r \times
\Biggr.
$$

\vspace{-2mm}
$$
\times
\sum_{\stackrel{(\{\{g_1, g_2\}, \ldots, 
\{g_{2r-1}, g_{2r}\}\}, \{q_1, \ldots, q_{k-2r}\})}
{{}_{\{g_1, g_2, \ldots, 
g_{2r-1}, g_{2r}, q_1, \ldots, q_{k-2r}\}=\{1, 2, \ldots, k\}}}}
\prod\limits_{s=1}^r
{\bf 1}_{\{i_{g_{{}_{2s-1}}}=~i_{g_{{}_{2s}}}\ne 0\}}
\Biggl.{\bf 1}_{\{j_{g_{{}_{2s-1}}}=~j_{g_{{}_{2s}}}\}}
\prod_{l=1}^{k-2r}\zeta_{j_{q_l}}^{(i_{q_l})}\Biggr)
$$

\vspace{2mm}
\noindent
con\-verg\-ing in the mean-square sense is valid$,$
where $[x]$ is an integer part of a real number $x,$
$$
C_{j_k\ldots j_1}(s)=\int\limits_{[t,T]^k}
\bar K(t_1,\ldots,t_k,s)\prod_{l=1}^{k}\phi_{j_l}(t_l)dt_1\ldots dt_k=
$$
$$
=\int\limits_t^s\psi_k(t_k)\phi_{j_k}(t_k)\ldots
\int\limits_t^{t_2}
\psi_1(t_1)\phi_{j_1}(t_1)
dt_1\ldots dt_k
$$

\vspace{2mm}
\noindent
is the Fourier coefficient$,$ $\prod\limits_{\emptyset}
\stackrel{\sf def}{=}1,$ $\sum\limits_{\emptyset}
\stackrel{\sf def}{=}0;$ another
notations are the same as in Theorem {\rm 1.2}.}

\vspace{2mm}

Note that the estimates (\ref{road888}) and (\ref{agent01000})
will also be valid under the conditions of Theorem~1.24.

\vspace{-5mm}

\chapter{Expansions of Iterated Stratonovich Stochastic Integrals 
Based on Generalized Multiple and Iterated Fourier Series}

\vspace{-5mm}

This chapter is devoted to the adaptation
of Theorems 1.1, 1.16 for 
iterated Stratonovich stochastic integrals.
The case of continuously differentiable 
weight functions (multiplicities 1 to 5) and 
weight functions identically equal to one (multiplicities 6 to 8)
is considered. In this case, we use
a complete orthonormal system of Legendre polynomials or 
trigonometric functions in $L_2([t, T])$.
In addition, the case of continuous weight functions (multiplicities 1 and 2),
binomial weight functions (multiplicities 3 and 4)
and weight functions identically equal to one (multiplicities 5 and 6)
is studied. In this case, we use
an arbitrary complete orthonormal system of functions in $L_2([t, T])$.
Recently, the above adaptation has also been carried out for 
iterated Stratonovich stochastic
integrals of multiplicity $k,$ $k\in{\bf N},$ under one additional condition
(Theorems~2.59, 2.61 (2024)) and without this condition
(Theorems 2.87 and 2.88 (2026)).

\section{Expansions of Iterated Stratonovich Stochastic Integrals of 
Multiplicity 2 Based
on Theorem 1.1. The case $p_1, p_2\to \infty$ and Smooth Weight 
Functions}

\subsection{Approach Based on Theorem 1.1 and Integration by Parts}

Let $(\Omega,{\rm F},{\sf P})$ be a complete probability
space and let $w(t,\omega)\stackrel{\sf def}{=}w_t:$ 
$[0, T]\times \Omega\rightarrow {\bf R}$
be the standard Wiener process
defined on the probability space $(\Omega,{\rm F},{\sf P}).$

Consider the family of $\sigma$-algebras
$\left\{{\rm F}_t,\ t\in[0,T]\right\}$ defined
on the probability space $(\Omega,{\rm F},{\sf P})$ and
connected
with the Wiener process $w_t$ in such a way that

1.\ ${\rm F}_s\subset {\rm F}_t\subset {\rm F}$\ for
$s<t.$

2.\ The Wiener process $w_t$ is ${\rm F}_t$-measurable for all
$t\in[0,T].$

3.\ The process $w_{t+\Delta}-w_{t}$ for all
$t\ge 0,$ $\Delta>0$ is independent with
the events of $\sigma$-algebra
${\rm F}_{t}.$

Let ${\rm M}_2([t, T])$ $(t\ge 0)$  be the class of random functions 
$\xi(\tau,\omega)\stackrel{\sf def}{=}
\xi_{\tau}: [t, T]\times$ $\Omega \to {\bf R}$ 
defined as in 
Sect.~1.1.2.

We introduce the class ${\rm Q}_m([t, T])$ $(t\ge 0)$ of It\^{o} processes 
$\eta_{\tau},$ $\tau \in [t, T]$ of the form
\begin{equation}
\label{z900}
\eta_{\tau} = \eta_t + \int\limits_t^{\tau} a_sds +
\int\limits_t^{\tau} b_sdw_s, 
\end{equation}
where $\left(a_{\tau}\right)^m, \left(b_{\tau}\right)^m 
\in {\rm M}_2([t, T])$ and
$\lim\limits_{s\to\tau} {\sf M}\bigl\{\left\vert b_s - b_{\tau}\right\vert ^4\bigr\}=0$
for all $\tau \in [t, T].$
The second integral on the right-hand side of (\ref{z900})
is the It\^{o} stochastic integral (see Sect.~1.1.2).

Let $C^{2,1}({\bf R}\times [t, T])$ $(t\ge 0)$ be the space of functions
$F(x,\tau): {\bf R} \times [t, T] \to {\bf R}$ such that
$$
\left|\frac{\partial F}{\partial x}(x,\tau)\right|\le K,\ \ 
\left|\frac{\partial^2 F}{\partial x^2}(x,\tau) \right|\le K,\ \ 
\left|\frac{\partial F}{\partial \tau}(x,\tau)\right|\le K,\ \ 
\left|\frac{\partial^2 F}{\partial \tau \partial x}(x,\tau) \right|\le K
$$ 

\noindent
for all $x\in {\bf R}$ and $\tau\in [t, T],$ where constant $K$ does not depend on $x,\tau.$

Let $\tau_j^{(N)},$ $j=0, 1, \ldots, N$ 
be a partition of the interval $[t, T],$ $t\ge 0$ such that
\begin{equation}
\label{usl}
t=\tau_0^{(N)}<\tau_1^{(N)}<\ldots <\tau_N^{(N)}=T,\ \ \ \
\max\limits_{0\le j\le N-1}\left|\tau_{j+1}^{(N)}-\tau_j^{(N)}\right|\to 0\ \
\hbox{if}\ \ N\to \infty.
\end{equation}
\par
The mean-square limit
\begin{equation}
\label{123321.2}
\hbox{\vtop{\offinterlineskip\halign{
\hfil#\hfil\cr
{\rm l.i.m}\cr
$\stackrel{}{{}_{N\to \infty}}$\cr
}} }\sum_{j=0}^{N-1}F\left(\frac{1}{2}\left(
\eta_{\tau_j^{(N)}}+\eta_{\tau_{j+1}^{(N)}}\right),\tau_j^{(N)}\right)
\left(w_{\tau_{j+1}^{(N)}}-
w_{\tau_j^{(N)}}\right)
\stackrel{\sf def}{=}{\int\limits_t^{*}}^T F(\eta_{\tau},\tau)dw_\tau
\end{equation}
is called \cite{str} the Stratonovich stochastic integral 
of the process $F(\eta_{\tau}, \tau)$, $\tau\in [t, T]$,
where $\tau_j^{(N)},$ $j=0, 1, \ldots, N$
is a partition of the interval $[t, T]$ 
satisfying the condition (\ref{usl}).

It is known \cite{str} (also see \cite{123789000} (Sect.~1.2.5)) 
that under proper conditions, the following 
relation between Stratonovich and It\^{o} stochastic integrals holds
\newpage
\noindent
\begin{equation}
\label{d11}
~~~~~{\int\limits_{t}^{*}}^T F(\eta_{\tau},\tau)dw_{\tau}=
\int\limits_t^T F(\eta_{\tau},\tau)dw_{\tau}+
\frac{1}{2}\int\limits_t^T \frac{\partial F}{\partial x}(\eta_{\tau},\tau)
b_{\tau}d\tau\ \ \ \hbox{w.~p.~1}.
\end{equation}

If the Wiener processes in (\ref{z900}) and 
(\ref{123321.2}) are independent, then
\begin{equation}
\label{d11a}
{\int\limits_{t}^{*}}^T F(\eta_{\tau},\tau)dw_{\tau}=
\int\limits_t^T F(\eta_{\tau},\tau)dw_{\tau}\ \ \ \hbox{w.~p.~1}.
\end{equation}

A possible variant of conditions under which the formulas
(\ref{d11}) and (\ref{d11a}) are correct, for example, consists
of the conditions:
$\eta_{\tau}\in {\rm Q}_4([t,T]),$
$F(\eta_{\tau},\tau)\in {\rm M}_2([t,T]),$
$F(x,\tau)\in C^{2,1}({\bf R}\times [t, T]).$

Note that if $F(x,\tau)=F_1(x)F_2(\tau),$ then it suffices to require
that $F(x,\tau)$ be twice differentiable with respect to $x$ 
$($with bounded derivatives$)$ and continuous with respect to $\tau$
$($instead of the condition $F(x,\tau)\in C^{2,1}({\bf R}\times [t, T])).$

In Sect.~2.1--2.17, in most cases, $\{\phi_j(x)\}_{j=0}^{\infty}$ 
is a complete orthonormal systems of Legendre polynomials or
trigonometric functions 
in $L_2([t, T])$.
Therefore, we will pay attention on the 
following well known facts about these two systems of functions \cite{Gob}.

{\it Suppose that the function $f(x)$ is 
bounded at the interval $[t, T].$ Moreover, its derivative
$f'(x)$ is continuous function at the interval $[t, T]$ except may be
the finite number of points 
of the finite discontinuity.
Then the Fourier series 
$$
\sum\limits_{j=0}^{\infty}
C_j\phi_j(x),\ \ \ C_j=\int\limits_t^{T}f(x)\phi_j(x)dx
$$
converges at any internal point $x$ of 
the interval $[t, T]$ to the value 
$\left(f(x+0)+f(x-0)\right)/2$ and converges  
uniformly to $f(x)$ on any closed interval {\rm (}of continuity 
of the function
$f(x)${\rm )} lying inside 
$[t, T]$. At the same time the Fourier--Legendre series 
converges 
if $x=t$ and $x=T$ to $f(t+0)$ and $f(T-0)$ 
correspondently, and the trigonometric Fourier series converges if   
$x=t$ and $x=T$ to $\left(f(t+0)+f(T-0)\right)/2$ 
in the case of periodic continuation 
of the function $f(x)$}.

In Sect.~2.1 we consider the case $k=2$ of the following iterated
Stratonovich and It\^{o} stochastic integrals
\begin{equation}
\label{strxx}
J^{*}[\psi^{(k)}]_{T,t}=
{\int\limits_t^{*}}^T
\psi_k(t_k)\ldots {\int\limits_t^{*}}^{t_2}
\psi_1(t_1) d{\bf w}_{t_1}^{(i_1)}\ldots d{\bf w}_{t_k}^{(i_k)},
\end{equation}
\begin{equation}
\label{itoxx}
J[\psi^{(k)}]_{T,t}=\int\limits_t^T\psi_k(t_k)\ldots \int\limits_t^{t_{2}}
\psi_1(t_1) d{\bf w}_{t_1}^{(i_1)}\ldots
d{\bf w}_{t_k}^{(i_k)},
\end{equation}
where every $\psi_l(\tau)$ $(l=1,\ldots,k)$ 
is a continuous nonrandom function
at the interval $[t,T],$ 
${\bf w}_{\tau}^{(i)}$ $(i=1,\ldots,m)$ are independent 
standard Wiener processes
and ${\bf w}_{\tau}^{(0)}=\tau$.

Let us formulate and prove the following
theorem on expansion of iterated Stratonovich stochastic integrals of 
multiplicity 2.

{\bf Theorem 2.1} \cite{8} (2011), \cite{9}-\cite{art-5}, \cite{arxiv-5}.
{\it Suppose that 
$\{\phi_j(x)\}_{j=0}^{\infty}$ is a complete orthonormal system of 
Legendre polynomials or trigonometric functions in the space $L_2([t, T]).$
At the same time $\psi_2(s)$ is a continuously dif\-ferentiable 
nonrandom function on $[t, T]$ and $\psi_1(s)$ is twice 
continuously differentiable nonrandom function on $[t, T]$. 
Then$,$ for the iterated Stratonovich stochastic integral
$$
J^{*}[\psi^{(2)}]_{T,t}={\int\limits_t^{*}}^T\psi_2(t_2)
{\int\limits_t^{*}}^{t_2}\psi_1(t_1)d{\bf w}_{t_1}^{(i_1)}
d{\bf w}_{t_2}^{(i_2)}\ \ \ (i_1, i_2=1,\ldots,m)
$$
the following expansion 
$$
J^{*}[\psi^{(2)}]_{T,t}=\hbox{\vtop{\offinterlineskip\halign{
\hfil#\hfil\cr
{\rm l.i.m.}\cr
$\stackrel{}{{}_{p_1,p_2\to \infty}}$\cr
}} }\sum_{j_1=0}^{p_1}\sum_{j_2=0}^{p_2}
C_{j_2j_1}\zeta_{j_1}^{(i_1)}\zeta_{j_2}^{(i_2)}
$$
that converges in the mean-square
sense   
is valid, where 
$$
C_{j_2 j_1}=\int\limits_t^T\psi_2(s_2)\phi_{j_2}(s_2)
\int\limits_t^{s_2}\psi_1(s_1)\phi_{j_1}(s_1)ds_1ds_2
$$
and
$$
\zeta_{j}^{(i)}=
\int\limits_t^T \phi_{j}(s) d{\bf w}_s^{(i)}
$$ 
are independent
standard Gaussian random variables for various 
$i$ or $j$.}

{\bf Proof.} In accordance to the standard relations between
Stra\-to\-no\-vich and It\^{o} stochastic integrals (see (\ref{d11})
and (\ref{d11a})) 
we have w.~p.~1 
\begin{equation}
\label{oop51}
J^{*}[\psi^{(2)}]_{T,t}=
J[\psi^{(2)}]_{T,t}+
\frac{1}{2}{\bf 1}_{\{i_1=i_2\}}
\int\limits_t^T\psi_1(t_1)\psi_2(t_1)dt_1,
\end{equation}
where here and further ${\bf 1}_A$ is the indicator of the set $A.$

From the other side according to (\ref{a2}), we have
$$
J[\psi^{(2)}]_{T,t}=
\hbox{\vtop{\offinterlineskip\halign{
\hfil#\hfil\cr
{\rm l.i.m.}\cr
$\stackrel{}{{}_{p_1,p_2\to \infty}}$\cr
}} }\sum_{j_1=0}^{p_1}\sum_{j_2=0}^{p_2}
C_{j_2j_1}\Biggl(\zeta_{j_1}^{(i_1)}\zeta_{j_2}^{(i_2)}
-{\bf 1}_{\{i_1=i_2\}}
{\bf 1}_{\{j_1=j_2\}}\Biggr)=
$$
\begin{equation}
\label{yes2001}
~~~~ =\hbox{\vtop{\offinterlineskip\halign{
\hfil#\hfil\cr
{\rm l.i.m.}\cr
$\stackrel{}{{}_{p_1,p_2\to \infty}}$\cr
}} }\sum_{j_1=0}^{p_1}\sum_{j_2=0}^{p_2}
C_{j_2j_1}\zeta_{j_1}^{(i_1)}\zeta_{j_2}^{(i_2)}
-{\bf 1}_{\{i_1=i_2\}}\lim\limits_{p_1,p_2\to\infty}\sum_{j_1=0}^{\min\{p_1,p_2\}}
C_{j_1j_1}.
\end{equation}

From (\ref{oop51}) and (\ref{yes2001}) it follows that
Theorem 2.1 will be proved if 
\begin{equation}
\label{5t}
\frac{1}{2}
\int\limits_t^T\psi_1(t_1)\psi_2(t_1)dt_1
=\sum_{j_1=0}^{\infty}
C_{j_1j_1}.
\end{equation}

Note that in this section and in Sect.~2.1.2 we 
present two different proofs (under different conditions) of the existence of a limit on 
the right-hand side of (\ref{5t}) for the polynomial and trigonometric cases.

Let us prove (\ref{5t}). Consider the function
\begin{equation}
\label{yes2002}
K^{*}(t_1,t_2)=K(t_1,t_2)+\frac{1}{2}{\bf 1}_{\{t_1=t_2\}}
\psi_1(t_1)\psi_2(t_1),
\end{equation}

\noindent
where $t_1, t_2\in[t, T]$ and $K(t_1,t_2)$ is defined by
(\ref{ppp}) for $k=2.$

Let us expand the function $K^{*}(t_1,t_2)$ defined by (\ref{yes2002})
using the variable 
$t_1$, when $t_2$ is fixed, into the generalized Fourier series 
at the interval $(t, T)$
\begin{equation}
\label{leto8001yes1}
K^{*}(t_1,t_2)=
\sum_{j_1=0}^{\infty}C_{j_1}(t_2)\phi_{j_1}(t_1)\ \ \ (t_1\ne t, T),
\end{equation}
where
\begin{equation}
\label{fur1}
~~~~~~~~C_{j_1}(t_2)=\int\limits_t^T
K^{*}(t_1,t_2)\phi_{j_1}(t_1)dt_1=\psi_2(t_2)
\int\limits_t^{t_2}\psi_1(t_1)\phi_{j_1}(t_1)dt_1.
\end{equation}

The equality (\ref{leto8001yes1}) is 
satisfied
pointwise in each point of the interval $(t, T)$ with respect to the 
variable $t_1$, when $t_2\in [t, T]$ is fixed, due to 
a piecewise
smoothness of the function $K^{*}(t_1,t_2)$ with respect to the variable 
$t_1\in [t, T]$ ($t_2$ is fixed). 

Note also that due to well known properties of the 
Fourier--Legendre series
and trigonometric Fourier series, 
the series (\ref{leto8001yes1}) converges when $t_1=t, T$. 

Obtaining (\ref{leto8001yes1}) we also used the fact that the right-hand
side 
of (\ref{leto8001yes1}) converges when $t_1=t_2$ (point of a finite 
discontinuity
of the function $K(t_1,t_2)$) to the value
$$
\frac{1}{2}\left(K(t_2-0,t_2)+K(t_2+0,t_2)\right)=
\frac{1}{2}\psi_1(t_2)\psi_2(t_2)=
K^{*}(t_2,t_2).
$$

The function $C_{j_1}(t_2)$ is a continuously differentiable
one at the interval $[t, T]$. 
Let us expand it into the generalized Fourier series at the interval $(t, T)$
\begin{equation}
\label{leto8002yes}
C_{j_1}(t_2)=
\sum_{j_2=0}^{\infty}C_{j_2 j_1}\phi_{j_2}(t_2)\ \ \ (t_2\ne t, T),
\end{equation}
where 
$$
C_{j_2 j_1}=\int\limits_t^T
C_{j_1}(t_2)\phi_{j_2}(t_2)dt_2=
\int\limits_t^T
\psi_2(t_2)\phi_{j_2}(t_2)\int\limits_t^{t_2}
\psi_1(t_1)\phi_{j_1}(t_1)dt_1 dt_2,
$$
and the equality (\ref{leto8002yes}) is satisfied pointwise at any point 
of the interval $(t, T)$ (the right-hand side 
of (\ref{leto8002yes}) converges 
when $t_2=t, T$).

Let us substitute (\ref{leto8002yes}) into (\ref{leto8001yes1})
\begin{equation}
\label{leto8003yes}
~~~~~~~~~~K^{*}(t_1,t_2)=
\sum_{j_1=0}^{\infty}\sum_{j_2=0}^{\infty}C_{j_2 j_1}
\phi_{j_1}(t_1)\phi_{j_2}(t_2),\ \ \ (t_1, t_2)\in (t, T)^2,
\end{equation}

\noindent
where the 
series on the right-hand side of (\ref{leto8003yes}) converges at the 
boundary
of the square  $[t, T]^2$.

It is easy to see that substituting $t_1=t_2$ in (\ref{leto8003yes}), we
obtain
\begin{equation}
\label{uiyes}
\frac{1}{2}\psi_1(t_1)\psi_2(t_1)=
\sum_{j_1=0}^{\infty}\sum_{j_2=0}^{\infty}
C_{j_2j_1}\phi_{j_1}(t_1)\phi_{j_2}(t_1).
\end{equation}
                                                           
From (\ref{uiyes}) we formally have
$$
\frac{1}{2}\int\limits_t^T\psi_1(t_1)\psi_2(t_1)dt_1=
\int\limits_t^T
\sum_{j_1=0}^{\infty}\sum_{j_2=0}^{\infty}
C_{j_2j_1}\phi_{j_1}(t_1)\phi_{j_2}(t_1)dt_1=
$$
$$
=
\sum_{j_1=0}^{\infty}\sum_{j_2=0}^{\infty}
\int\limits_t^T C_{j_2j_1}\phi_{j_1}(t_1)\phi_{j_2}(t_1)dt_1=
$$
$$
=\lim\limits_{p_1\to\infty}\lim\limits_{p_2\to\infty}
\sum_{j_1=0}^{p_1}\sum_{j_2=0}^{p_2}
C_{j_2j_1}\int\limits_t^T\phi_{j_1}(t_1)\phi_{j_2}(t_1)dt_1=
$$
\begin{equation}
\label{rozayes}
=\lim\limits_{p_1\to\infty}\lim\limits_{p_2\to\infty}
\sum_{j_1=0}^{p_1}\sum_{j_2=0}^{p_2}
C_{j_2j_1}{\bf 1}_{\{j_1=j_2\}}=
\lim\limits_{p_1\to\infty}\lim\limits_{p_2\to\infty}
\sum_{j_1=0}^{{\rm min}\{p_1,p_2\}}
C_{j_1j_1}=
\sum_{j_1=0}^{\infty}C_{j_1j_1}.
\end{equation}

Let us explain the second step in (\ref{rozayes})
(the fourth step in (\ref{rozayes}) follows from the orthonormality of 
functions $\phi_j(s)$ at the interval $[t, T]$).

We have
$$
\left|\int\limits_t^T \sum_{j_1=0}^{\infty}C_{j_1}(t_1)\phi_{j_1}(t_1)dt_1
-\sum_{j_1=0}^{p_1}\int\limits_t^TC_{j_1}(t_1)\phi_{j_1}(t_1)dt_1\right|
\le
$$
\begin{equation}
\label{otit2001}
\le\int\limits_t^T 
\left|
\psi_2(t_1)
G_{p_1}(t_1)
\right| dt_1
\le C\int\limits_t^T \left|G_{p_1}(t_1)\right| dt_1,
\end{equation}
where $C<\infty$ and
$$
\sum_{j=p+1}^{\infty}
\int\limits_t^{\tau}\psi_1(s)\phi_{j}(s)ds
\phi_{j}(\tau)\stackrel{\sf def}{=}G_{p}(\tau).
$$

Let us consider the case of Legendre polynomials. Then
\begin{equation}
\label{otit2003}
~~~~~~~~~ \left|G_{p_1}(t_1)\right|
=
\frac{1}{2}\left| \sum_{j_1=p_1+1}^{\infty}(2j_1+1)
\int\limits_{-1}^{z(t_1)}\psi_1(u(y))
P_{j_1}(y)dy
P_{j_1}(z(t_1))
\right|,
\end{equation}
where 
\begin{equation}
\label{zz1}
u(y)=\frac{T-t}{2}y+\frac{T+t}{2},\ \ \
z(s)=\left(s-\frac{T+t}{2}\right)\frac{2}{T-t},
\end{equation}
and $P_j(s)$ is the Legendre polynomial.

From (\ref{otit2003}) and the well known formula
\begin{equation}
\label{w1ggg}
\frac{dP_{j+1}}{dx}(x)-\frac{dP_{j-1}}{dx}(x)=(2j+1)P_j(x),\ \ \ 
j=1, 2,\ldots
\end{equation}
we obtain
$$
\left|G_{p_1}(t_1)\right|=\frac{1}{2}\Biggl|\sum_{j_1=p_1+1}^{\infty}
\Biggl\{\left(P_{j_1+1}(z(t_1))-P_{j_1-1}(z(t_1))\right)\psi_1(t_1)-
\Biggr.\Biggr.
$$
$$
\Biggl.\Biggl.-\frac{T-t}{2}\int\limits_{-1}^{z(t_1)}
\left(P_{j_1+1}(y)-P_{j_1-1}(y)\right)\psi_1'
(u(y))dy\Biggr\}P_{j_1}(z(t_1))
\Biggr|\le
$$
$$
\le C_0
\Biggl|\sum_{j_1=p_1+1}^{\infty}
(P_{j_1+1}(z(t_1))P_{j_1}(z(t_1))-P_{j_1-1}(z(t_1))P_{j_1}(z(t_1)))
\Biggr|+
$$
$$
+\frac{T-t}{4}
\Biggl|\sum_{j_1=p_1+1}^{\infty}
\Biggl\{\psi_1'(t_1)\Biggl(\frac{1}{2j_1+3}(P_{j_1+2}(z(t_1))
-P_{j_1}(z(t_1)))
- \Biggr.\Biggr.\Biggr.
$$
$$
\Biggl.-\frac{1}{2j_1-1}(P_{j_1}(z(t_1))
-P_{j_1-2}(z(t_1)))\Biggr) - 
$$
$$
-\frac{T-t}{2}
\int\limits_{-1}^{z(t_1)}
\Biggl(\frac{1}{2j_1+3}(P_{j_1+2}(y)-P_{j_1}(y))-\Biggr.
$$
\begin{equation}
\label{otit2005}
\Biggl.\Biggl.\Biggl.-
\frac{1}{2j_1-1}(P_{j_1}(y)-P_{j_1-2}(y))\Biggr)
\psi_1''(u(y))dy\Biggr\}P_{j_1}(z(t_1))
\Biggr|,
\end{equation}

\vspace{1mm}
\noindent
where $C_0$ is a constant, $\psi_1'$ and $\psi_1''$ are
derivatives of the function $\psi_1(s)$ with respect to the variable
$u(y)$. 

From (\ref{otit2005}) and the well known estimate for Legendre
polynomials \cite{Gob}
\begin{equation}
\label{otit987}
\left|P_n(y)\right| <\frac{K}{\sqrt{n+1}(1-y^2)^{1/4}},\ \ \ 
y\in (-1, 1),\ \ \ n\in {\bf N},
\end{equation}

\noindent
where constant $K$ does not depend on $y$ and $n$, we have
$$
\left|G_{p_1}(t_1)\right|<
$$
$$
<C_0\Biggl|\lim_{n\to\infty} 
\sum_{j_1=p_1+1}^{n}
(P_{j_1+1}(z(t_1))P_{j_1}(z(t_1))-P_{j_1-1}(z(t_1))P_{j_1}(z(t_1)))
\Biggr|+
$$
$$
+C_1\sum_{j_1=p_1+1}^{\infty}\frac{1}{j_1^2}\left(
\frac{1}{\left(1-(z(t_1))^2\right)^{1/2}}+
\int\limits_{-1}^{z(t_1)}\frac{dy}{\left(1-y^2\right)^{1/4}}
\frac{1}{\left(1-(z(t_1))^2\right)^{1/4}}\right)<
$$
$$
< C_0\Biggl|\lim_{n\to\infty} 
\left(P_{n+1}(z(t_1))P_{n}(z(t_1))-P_{p_1}(z(t_1))P_{p_1+1}(z(t_1))\right)
\Biggr|+
$$
$$
+C_1\sum_{j_1=p_1+1}^{\infty}\frac{1}{j_1^2}\left(
\frac{1}{\left(1-(z(t_1))^2\right)^{1/2}}+
C_2
\frac{1}{\left(1-(z(t_1))^2\right)^{1/4}}\right)<
$$
$$
<
C_3\lim_{n\to\infty} 
\Biggl(\frac{1}{n}+\frac{1}{p_1}\Biggr)
\frac{1}{\left(1-(z(t_1))^2\right)^{1/2}}+
$$
$$
+C_1\sum_{j_1=p_1+1}^{\infty}\frac{1}{j_1^2}\Biggl(
\frac{1}{\left(1-(z(t_1))^2\right)^{1/2}}+
C_2
\frac{1}{\left(1-(z(t_1))^2\right)^{1/4}}\Biggr)\le
$$
$$
\le C_4\Biggl(\Biggl(\frac{1}{p_1}+\sum_{j_1=p_1+1}^{\infty}\frac{1}{j_1^2}
\Biggr)
\frac{1}{\left(1-(z(t_1))^2\right)^{1/2}}+
\sum_{j_1=p_1+1}^{\infty}\frac{1}{j_1^2}
\frac{1}{\left(1-(z(t_1))^2\right)^{1/4}}\Biggr)\le
$$
\begin{equation}
\label{otit2007}
\le\frac{K}{p_1}\Biggl(
\frac{1}{\left(1-(z(t_1))^2\right)^{1/2}}+
\frac{1}{\left(1-(z(t_1))^2\right)^{1/4}}\Biggr),
\end{equation}

\vspace{2mm}
\noindent
where $C_0, C_1,\ldots, C_4, K$ are constants, $t_1\in (t, T)$, and 
\begin{equation}
\label{obana}
\sum\limits_{j_1=p_1+1}^{\infty}\frac{1}{j_1^2}
\le \int\limits_{p_1}^{\infty}\frac{dx}{x^2}=\frac{1}{p_1}.
\end{equation}

From (\ref{otit2001}) and (\ref{otit2007}) we get
$$
\left| \int\limits_t^T \sum_{j_1=0}^{\infty}C_{j_1}(t_1)\phi_{j_1}(t_1)dt_1
-\sum_{j_1=0}^{p_1}\int\limits_t^TC_{j_1}(t_1)\phi_{j_1}(t_1)dt_1\right|<
$$
$$
<
\frac{K}{p_1}\left(
\int\limits_{-1}^{1}\frac{dy}{\left(1-y^2\right)^{1/2}}+
\int\limits_{-1}^{1}\frac{dy}{\left(1-y^2\right)^{1/4}}
\right)\ \to 0
$$

\noindent
if $p_1\to\infty$. So, we obtain
$$
\frac{1}{2}\int\limits_t^T\psi_1(t_1)\psi_2(t_1)dt_1=
\int\limits_t^T
\sum_{j_1=0}^{\infty}C_{j_1}(t_1)\phi_{j_1}(t_1)dt_1=
$$
$$
=
\sum_{j_1=0}^{\infty}\int\limits_t^TC_{j_1}(t_1)\phi_{j_1}(t_1)dt_1
=
\sum_{j_1=0}^{\infty}\int\limits_t^T
\sum_{j_2=0}^{\infty}
C_{j_2j_1}\phi_{j_2}(t_1)\phi_{j_1}(t_1)dt_1=
$$
\begin{equation}
\label{otit2009}
=
\sum_{j_1=0}^{\infty}\sum_{j_2=0}^{\infty}
\int\limits_t^T C_{j_2j_1}\phi_{j_2}(t_1)\phi_{j_1}(t_1)dt_1=
\sum_{j_1=0}^{\infty}C_{j_1j_1}.
\end{equation}

In (\ref{otit2009}) we used the fact that the Fourier--Legendre series
$$
\sum\limits_{j_2=0}^{\infty}C_{j_2j_1}\phi_{j_2}(t_1)
$$
of the smooth function $C_{j_1}(t_1)$ converges uniformly to
this function at the interval $[t+\varepsilon, T-\varepsilon]$ 
for any $\varepsilon>0$, converges to this function 
at the any point $t_1\in(t,T),$ 
and converges
to $C_{j_1}(t+0)$ and $C_{j_1}(T-0)$ when $t_1=t,$ $T.$

More precisely, we have
$$
\int\limits_t^T
\sum_{j_2=0}^{\infty}
C_{j_2j_1}\phi_{j_2}(t_1)\phi_{j_1}(t_1)dt_1=
\int\limits_{t+\varepsilon}^{T-\varepsilon}
\sum_{j_2=0}^{\infty}
C_{j_2j_1}\phi_{j_2}(t_1)\phi_{j_1}(t_1)dt_1+A_{\varepsilon}+B_{\varepsilon}=
$$
$$
=\sum_{j_2=0}^{\infty}C_{j_2j_1}\int\limits_{t+\varepsilon}^{T-\varepsilon}
\phi_{j_2}(t_1)\phi_{j_1}(t_1)dt_1+A_{\varepsilon}+B_{\varepsilon}=
$$
$$
=\sum_{j_2=0}^{\infty}C_{j_2j_1}\left(\int\limits_{t}^{T}-
\int\limits_t^{t+\varepsilon}-\int\limits_{T-\varepsilon}^{T}\right)
\phi_{j_2}(t_1)\phi_{j_1}(t_1)dt_1+A_{\varepsilon}+B_{\varepsilon}=
$$
$$
=\sum_{j_2=0}^{\infty}C_{j_2j_1}\biggl({\bf 1}_{\{j_1=j_2\}}-
\varepsilon\bigl(\phi_{j_2}(\lambda)\phi_{j_1}(\lambda)+
\phi_{j_2}(\theta)\phi_{j_1}(\theta)\bigr)\biggr)
+A_{\varepsilon}+B_{\varepsilon}=
$$
\begin{equation}
\label{dwdw30}
=C_{j_1j_1}-\varepsilon \left(\sum_{j_2=0}^{\infty}
C_{j_2j_1}\phi_{j_2}(\lambda)\phi_{j_1}(\lambda)+
\sum_{j_2=0}^{\infty}
C_{j_2j_1}\phi_{j_2}(\theta)\phi_{j_1}(\theta)\right)
+A_{\varepsilon}+B_{\varepsilon},
\end{equation}

\vspace{2mm}
\noindent
where $\theta\in [t,t+\varepsilon],$ $\lambda\in [T-\varepsilon,T]$, and
$$
A_{\varepsilon}=
\int\limits_t^{t+\varepsilon}
\sum_{j_2=0}^{\infty}
C_{j_2j_1}\phi_{j_2}(t_1)\phi_{j_1}(t_1)dt_1,\ \ \ 
B_{\varepsilon}=
\int\limits_{T-\varepsilon}^T
\sum_{j_2=0}^{\infty}
C_{j_2j_1}\phi_{j_2}(t_1)\phi_{j_1}(t_1)dt_1.
$$

In obtaining (\ref{dwdw30}) we used the theorem
on the mean value for the Riemann  
integral and orthonormality of the functions
$\phi_{j}(x)$ for $j=0, 1, 2\ldots$

Further, we have
$\left|A_{\varepsilon}\right|+\left|B_{\varepsilon}\right|\le \varepsilon C,$
where $C<\infty$ is a constant.
Performing the 
passage to the limit $\lim\limits_{\varepsilon\to +0}$
in the equality (\ref{dwdw30}), we get 
$$
\int\limits_t^T
\sum_{j_2=0}^{\infty}
C_{j_2j_1}\phi_{j_2}(t_1)\phi_{j_1}(t_1)dt_1=
C_{j_1j_1}.
$$

Then (see (\ref{otit2009}))
$$
\sum_{j_1=0}^{\infty}\int\limits_t^T
\sum_{j_2=0}^{\infty}
C_{j_2j_1}\phi_{j_2}(t_1)\phi_{j_1}(t_1)dt_1=
\sum_{j_1=0}^{\infty}C_{j_1j_1}
$$

\noindent
and the relation (\ref{5t}) is proved for the case of Legendre polynomials.

Let us consider the trigonometric case
and suppose that $\{\phi_j(x)\}_{j=0}^{\infty}$ is a complete orthonormal
system of trigonometric functions in $L_2([t, T])$.

Denote
$$
S_{p_1}\stackrel{\sf def}{=}
\left| \int\limits_t^T \sum_{j_1=0}^{\infty}C_{j_1}(t_1)\phi_{j_1}(t_1)dt_1
-\sum_{j_1=0}^{p_1}\int\limits_t^TC_{j_1}(t_1)\phi_{j_1}(t_1)dt_1\right|=
$$
$$
=
\left|\int\limits_t^T\sum_{j_1=p_1+1}^{\infty}
\psi_2(t_1)\phi_{j_1}(t_1)\int\limits_t^{t_1}
\psi_1(\theta)\phi_{j_1}(\theta)d\theta dt_1\right|.
$$

We have
$$
S_{2p_1}=
\left| \int\limits_t^T \sum_{j_1=0}^{\infty}C_{j_1}(t_1)\phi_{j_1}(t_1)dt_1
-\sum_{j_1=0}^{2p_1}\int\limits_t^TC_{j_1}(t_1)\phi_{j_1}(t_1)dt_1\right|=
$$
$$
=
\left|\int\limits_t^T\sum_{j_1=2p_1+1}^{\infty}
\psi_2(t_1)\phi_{j_1}(t_1)\int\limits_t^{t_1}
\psi_1(\theta)\phi_{j_1}(\theta)d\theta dt_1\right|=
$$
$$
=\frac{2}{T-t}\left|
\int\limits_t^T \psi_2(t_1)\sum_{j_1=p_1+1}^{\infty}\left(
\int\limits_t^{t_1}\psi_1(s){\rm sin}\frac{2\pi j_1(s-t)}{T-t}ds\
{\rm sin}\frac{2\pi j_1(t_1-t)}{T-t}+\right.\right.
$$
$$
\left.\left.+\int\limits_t^{t_1}\psi_1(s)
{\rm cos}\frac{2\pi j_1(s-t)}{T-t}ds\
{\rm cos}\frac{2\pi j_1(t_1-t)}{T-t}\right)dt_1\right|=
$$
$$
=
\frac{1}{\pi}\left|
\int\limits_t^T \Biggl(\psi_1(t)\psi_2(t_1)\sum_{j_1=p_1+1}^{\infty}
\frac{1}{j_1}{\rm sin}\frac{2\pi j_1(t_1-t)}{T-t}+\Biggr.\right.
$$
$$
+\frac{T-t}{2\pi}\psi_2(t_1)\sum_{j_1=p_1+1}^{\infty}\frac{1}{j_1^2}\Biggl(
\psi_1'(t_1)-\psi_1'(t){\rm cos}\frac{2\pi j_1(t_1-t)}{T-t}-\Biggr.
$$
$$
-
\int\limits_t^{t_1}{\rm sin}\frac{2\pi j_1(s-t)}{T-t}\psi_1''(s)ds\
{\rm sin}\frac{2\pi j_1(t_1-t)}{T-t}-
$$
$$
\left.\Biggl.\Biggl.-
\int\limits_t^{t_1}{\rm cos}\frac{2\pi j_1(s-t)}{T-t}\psi_1''(s)ds\
{\rm cos}\frac{2\pi j_1(t_1-t)}{T-t}\Biggr)\Biggr)dt_1\right|\le
$$
$$
\le C_1 
\left|
\int\limits_t^T \psi_2(t_1)\sum_{j_1=p_1+1}^{\infty}
\frac{1}{j_1}{\rm sin}\frac{2\pi j_1(t_1-t)}{T-t}dt_1\right|
+\frac{C_2}{p_1}=
$$
\begin{equation}
\label{2017zzz1}
=
C_1\left|
\sum_{j_1=p_1+1}^{\infty}\frac{1}{j_1}
\int\limits_t^T \psi_2(t_1)
{\rm sin}\frac{2\pi j_1(t_1-t)}{T-t}dt_1\right|
+\frac{C_2}{p_1},
\end{equation}

\vspace{2mm}
\noindent
where constants $C_1, C_2$ do not depend on $p_1.$

Here we used the fact that the functional series
\begin{equation}
\label{1010}
\sum\limits_{j_1=1}^{\infty}
\frac{1}{j_1}{\rm sin}\frac{2\pi j_1(t_1-t)}{T-t}
\end{equation}
converges uniformly at the interval $[t+\varepsilon, T-\varepsilon]$
for any $\varepsilon>0$
due to Dirichlet--Abel Theorem, and converges to zero at the 
points $t$ and $T$.
Moreover, the series (\ref{1010}) (with accuracy to a
linear transformation) 
is the trigonometric Fourier
series of the smooth function $K(t_1)=t_1-t,$ $t_1\in [t, T].$ 
Thus, (\ref{1010}) converges to the smooth function at any point $t_1\in(t,T)$.

From (\ref{2017zzz1}) we obtain
\begin{equation}
\label{agentyq5}
S_{2p_1}\le C_3\left|
\sum_{j_1=p_1+1}^{\infty}
\frac{1}{j_1^2}\Biggl(
\psi_2(T)-\psi_2(t)-
\int\limits_t^{T}{\rm cos}\frac{2\pi j_1(s-t)}{T-t}\psi_2'(s)ds
\Biggr)\right|
+\frac{C_2}{p_1}\le \frac{C_4}{p_1},
\end{equation}

\noindent
where constants $C_2, C_3, C_4$ do not depend on $p_1.$

Further,
$$
S_{2p_1-1}
=
\left|\int\limits_t^T\sum_{j_1=2p_1}^{\infty}
\psi_2(t_1)\phi_{j_1}(t_1)\int\limits_t^{t_1}
\psi_1(\theta)\phi_{j_1}(\theta)d\theta dt_1\right|=
$$
$$
=\left|S_{2p_1}+
\int\limits_t^T
\psi_2(t_1)\phi_{2p_1}(t_1)\int\limits_t^{t_1}
\psi_1(\theta)\phi_{2p_1}(\theta)d\theta dt_1\right|\le
$$
\begin{equation}
\label{agentyq1}
\le S_{2p_1}+\frac{2}{T-t}
\left|\int\limits_t^T
\psi_2(t_1){\rm cos}\frac{2\pi p_1(t_1-t)}{T-t}\int\limits_t^{t_1}
\psi_1(\theta){\rm cos}\frac{2\pi p_1(\theta-t)}{T-t}d\theta dt_1\right|.
\end{equation}

\vspace{2mm}

Moreover, 
$$
\int\limits_t^T
\psi_2(t_1){\rm cos}\frac{2\pi p_1(t_1-t)}{T-t}\int\limits_t^{t_1}
\psi_1(\theta){\rm cos}\frac{2\pi p_1(\theta-t)}{T-t}d\theta dt_1=
$$
$$
=\frac{T-t}{2\pi p_1}
\int\limits_t^T
\psi_2(t_1){\rm cos}\frac{2\pi p_1(t_1-t)}{T-t}
\Biggl(\psi_1(t_1){\rm sin}\frac{2\pi p_1(t_1-t)}{T-t}-\Biggr.
$$
\begin{equation}
\label{agentyq2}
-\Biggl.\int\limits_t^T
\psi_1'(\theta){\rm sin}\frac{2\pi p_1(\theta-t)}{T-t}d\theta\Biggr)dt_1.
\end{equation}

\vspace{2mm}

The relations (\ref{agentyq5})--(\ref{agentyq2}) imply that
\begin{equation}
\label{agentyq4}
S_{2p_1-1}\le \frac{C_5}{p_1},
\end{equation}

\noindent
where constant $C_5$ is independent of $p_1$.

From (\ref{agentyq5}) and (\ref{agentyq4}) we obtain
\begin{equation}
\label{2017zzz111}
~~~~~~~~ S_{p_1}=\left|\int\limits_t^T\sum_{j_1=p_1+1}^{\infty}
\psi_2(t_1)\phi_{j_1}(t_1)\int\limits_t^{t_1}
\psi_1(\theta)\phi_{j_1}(\theta)d\theta dt_1\right|
\le \frac{K}{p_1}\to 0
\end{equation}

\noindent
if $p_1\to\infty,$ where constant $K$ does not depend on $p_1$ $(p_1\in {\bf N})$.

Further steps are similar to the proof
of (\ref{5t}) for the case of Legendre polynomials.
Theorem 2.1 is proved.

Note that the estimate (\ref{2017zzz111}) will be used further.

\subsection{Approach Based on Theorem 1.1 and Double Fourier--Le\-gen\-dre 
Series Summarized by Pringsheim Method}

In Sect.~2.1.1 we considered the proof of Theorem 2.1 
based on Theorem 1.1 and
double integration by parts (this procedure leads to 
the requirement of double continuous differentiability of the function
$\psi_1(\tau)$ 
at the interval $[t, T]$).
In this section, we formulate and prove an analogue
of Theorem 2.1 
but under the weakened conditions:
the functions
$\psi_1(\tau),$  $\psi_2(\tau)$ only one time
continuously differentiable 
at the interval $[t, T]$.
At that we will use 
the double Fourier series summarized by Pringsheim method.

{\bf Theorem 2.2}\ \cite{12}-\cite{12aa}, \cite{art-9}, 
\cite{arxiv-19}. {\it Suppose that 
$\{\phi_j(x)\}_{j=0}^{\infty}$ is a complete orthonormal system of 
Legendre polynomials or trigonometric functions in the space $L_2([t, T]).$
Moreover, $\psi_1(s),$ $\psi_2(s)$ are 
continuously differentiable functions on $[t, T]$. Then$,$ 
for the iterated Stratonovich stochastic integral
$$
J^{*}[\psi^{(2)}]_{T,t}={\int\limits_t^{*}}^T\psi_2(t_2)
{\int\limits_t^{*}}^{t_2}\psi_1(t_1)d{\bf w}_{t_1}^{(i_1)}
d{\bf w}_{t_2}^{(i_2)}\ \ \ (i_1, i_2=1,\ldots,m)
$$
the following expansion 
\begin{equation}
\label{jes}
J^{*}[\psi^{(2)}]_{T,t}=\hbox{\vtop{\offinterlineskip\halign{
\hfil#\hfil\cr
{\rm l.i.m.}\cr
$\stackrel{}{{}_{p_1,p_2\to \infty}}$\cr
}} }\sum_{j_1=0}^{p_1}\sum_{j_2=0}^{p_2}
C_{j_2j_1}\zeta_{j_1}^{(i_1)}\zeta_{j_2}^{(i_2)}
\end{equation}
that converges in the mean-square
sense 
is valid, where 
\begin{equation}
\label{tupo11}
C_{j_2 j_1}=\int\limits_t^T\psi_2(s_2)\phi_{j_2}(s_2)
\int\limits_t^{s_2}\psi_1(s_1)\phi_{j_1}(s_1)ds_1ds_2
\end{equation}
and
$$
\zeta_{j}^{(i)}=
\int\limits_t^T \phi_{j}(s) d{\bf w}_s^{(i)}
$$ 
are independent
standard Gaussian random variables for various 
$i$ or $j$.}

{\bf Proof.}\ Theorem 2.2 will be proved if we prove the equality
(see the proof of Theorem 2.1)
\begin{equation}
\label{may62021}
\frac{1}{2}
\int\limits_t^T\psi_1(t_1)\psi_2(t_1)dt_1
=\sum_{j_1=0}^{\infty}
C_{j_1j_1},
\end{equation}
where $C_{j_1j_1}$ is defined by the formula (\ref{ppppa})
for $k=2$ and $j_1=j_2.$ At that
$\{\phi_j(x)\}_{j=0}^{\infty}$ is a complete orthonormal system of 
Legendre polynomials or tri\-go\-no\-met\-ric functions in the space $L_2([t, T]).$

Firstly, consider the
sufficient conditions of convergence of double Fourier--Legendre 
series summarized by Pringsheim method.

Let $P_j(x)$ $(j=0, 1, 2,\ldots )$ be the Legendre polynomial.
Consider the function $f(x,y)$ defined for $(x,y)\in [-1,1]^2.$
Furthermore, consider the double Fourier--Legendre series summarized by
Pringsheim method and corresponding 
to the function 
$f(x,y)$
$$
\lim_{n,m\to\infty}
\sum_{j=0}^n\sum_{i=0}^m \frac{1}{2}\sqrt{(2j+1)(2i+1)}{C}_{ij}^{*}
P_i(x)P_j(y)
\stackrel{\sf{def}}{=}
$$
\begin{equation}
\label{555}
\stackrel{\sf{def}}{=}
\sum_{i,j=0}^{\infty}\frac{1}{2}\sqrt{(2j+1)(2i+1)}{C}_{ij}^{*}
P_i(x)P_j(y),
\end{equation}
where
\begin{equation}
\label{777}
{C}_{ij}^{*}=\frac{1}{2}\sqrt{(2j+1)(2i+1)}\int\limits_{[-1,1]^2}
f(x,y)P_i(x)P_j(y)dxdy.
\end{equation}

Consider the generalization for the case of two variables
\cite{star}
of the theorem on equiconvergence for the Fourier--Legendre
series
\cite{suet}.

{\bf Proposition 2.1} \cite{star}. {\it Let $f(x,y)\in L_2([-1,1]^2)$ and 
the function
$$
f(x,y)\left(1-x^2\right)^{-1/4}\left(1-y^2\right)^{-1/4}
$$
is integrable on
$[-1,1]^2.$  Moreover, let
$$
|f(x,y)-f(u,v)|\le G(y)|x-u|+H(x)|y-v|,
$$
where $G(y), H(x)$ are bounded functions on $[-1,1]^2$.
Then for all $(x,y)\in(-1,1)^2$ the following equality
is satisfied

$$
\lim\limits_{n,m\to\infty}
\Biggl(\sum_{j=0}^n\sum_{i=0}^m \frac{1}{2}\sqrt{(2j+1)(2i+1)}{C}_{ij}^{*}
P_i(x)P_j(y)-\Biggr.
$$
\begin{equation}
\label{888}
~~~~~~~-\Biggl.
(1-x^2)^{-1/4}(1-y^2)^{-1/4}S_{nm}(\mathop{\mathrm{arccos}} x,
\mathop{\mathrm{arccos}} y,F)
\Biggr)=0.
\end{equation}

\noindent
At that, the convergence in {\rm (\ref{888})} is uniform on the rectangle
$$
[-1+\varepsilon, 1-\varepsilon]\times[-1+\delta,1-\delta]\ \ \
\hbox{for any}\ \ \ \varepsilon, \delta>0,
$$ 
$S_{nm}(\theta,\varphi,F)$ is a partial sum of
the double trigonometric Fourier series of the
auxiliary function
$$
F(\theta,\varphi)=\sqrt{|\mathrm{sin} \theta|}\sqrt{|\mathrm{sin} \varphi|}
f({\rm cos}\theta,{\rm cos}\varphi),\ \
\theta,\varphi\in[0, \pi],
$$
and the Fourier coefficient ${C}_{ij}^{*}$ is defined by {\rm (\ref{777}).}
}

Proposition 2.1 implies that the following equality

\vspace{-4mm}
\begin{equation}
\label{8888}
~~~~~\lim\limits_{n,m\to\infty}
\Biggl(\sum_{j=0}^n\sum_{i=0}^m \frac{1}{2}\sqrt{(2j+1)(2i+1)}{C}_{ij}^{*}
P_i(x)P_j(y)-f(x,y)\Biggr)=0
\end{equation}

\vspace{2mm}
\noindent
is fulfilled for all $(x,y)\in(-1,1)^2,$
and convergence in (\ref{8888}) is uniform on the rectangle
$$
[-1+\varepsilon, 1-\varepsilon]\times[-1+\delta,1-\delta]\ \ \
\hbox{for any}\ \ \ \varepsilon, \delta>0
$$ 
if the corresponding conditions of convergence 
of the double trigonometric Fourier series of the auxiliary 
function 
\begin{equation}
\label{999}
g(x,y)=f(x,y)\left(1-x^2\right)^{1/4}\left(1-y^2\right)^{1/4}
\end{equation}
are satisfied.

Note also that Proposition 2.1 does not imply any
conclusions on the behavior of the double Fourier--Legendre series 
on the boundary of
the square 
$[-1,1]^2.$

{\it For each $\delta>0$ let us call the exact upper edge 
of difference $\left|f({\bf t}')-f({\bf t}'')\right|$ 
in the set
of all points ${\bf t}'$, ${\bf t}''$ which 
belong 
to the domain $D$ 
as the module of 
continuity of the function
$f({\bf t})$ $({\bf t}=(t_1,\ldots,t_k))$ in the 
$k$-dimentional domain
$D$ $(k\ge 1)$ 
if the distance between ${\bf t}',{\bf t}''$
satisfies the condition
$\rho\left({\bf t}',{\bf t}''\right)<\delta.$}

{\it We will say that the function 
of $k$ $(k\ge 1)$ variables  
$f({\bf t})$ $({\bf t}=(t_1,\ldots,t_k))$
belongs 
to the H\"{o}lder class with 
the parameter $\alpha\in (0, 1]$ $(f({\bf t})\in C^{\alpha}(D))$ 
in the domain $D$ 
if the module of 
continuity of the function
$f({\bf t})$ $({\bf t}=(t_1,\ldots,t_k))$ 
in the domain $D$ has orders $o(\delta^{\alpha})$ $(\alpha \in (0, 1))$
and $O(\delta)$ $(\alpha=1)$.}

In 1967, Zhizhiashvili L.V.
proved
that the rectangular sums of multiple trigonometric Fourier series 
of the function of $k$ variables  
in the hypercube $[t,T]^k$ 
converge uniformly to this function in the hypercube $[t,T]^k$ if 
the function
belongs
to $C^{\alpha}([t,T]^k),$ $\alpha>0$ (definition
of the H\"{o}lder class with any parameter $\alpha>0$ can be found in 
the well known mathematical analysis tutorials \cite{IP}).

More precisely, the following statement is correct.

{\bf Proposition 2.2} \cite{IP}. {\it If the function
$f(x_1,\ldots,x_n)$ is periodic with period $2\pi$ with respect to each
variable and belongs in ${\bf R}^n$ to the H\"{o}lder class 
$C^{\alpha}({\bf R}^n)$
for any $\alpha>0,$ then the rectangular partial
sums of multiple trigonometric Fourier series of the function
$f(x_1,\ldots,x_n)$ converge to this function uniformly in 
${\bf R}^n$.}

Let us back to the proof of Theorem 2.2 and consider the following Lemma.

{\bf Lemma 2.1}. {\it 
Let the function $f(x,y)$ satisfies
to the following condition 
$$
|f(x,y)-f(x_1,y_1)|\le C_1|x-x_1|+C_2|y-y_1|,
$$
where $C_1, C_2<\infty$ and $(x,y),$ $(x_1,y_1)\in [-1,1]^2.$
Then the following inequality is fulfilled 
\begin{equation}
\label{2222}
|g(x,y)-g(x_1,y_1)|\le K\rho^{1/4},
\end{equation}

\noindent
where $g(x,y)$ 
in defined by
{\rm(\ref{999})}{\rm,}
$$
\rho=\sqrt{(x-x_1)^2+(y-y_1)^2},
$$ 
$(x,y)$ and $(x_1,y_1)\in [-1,1]^2,$ $K<\infty.$}

{\bf Proof.}\ First, we assume that $x\ne x_1,$ $y\ne y_1.$ In this case
we have
$$
\left|g(x,y)-g(x_1,y_1)\right|=
$$
$$
=
\left|\left(1-x^2\right)^{1/4}
\left(1-y^2\right)^{1/4}\left(f(x,y)-f(x_1,y_1)\right)+\right.
$$
$$
\left.+
f(x_1,y_1)\left(\left(1-x^2\right)^{1/4}\left(1-y^2\right)^{1/4}
-\left(1-x_1^2\right)^{1/4}\left(1-y_1^2\right)^{1/4}\right)\right|\le 
$$
$$
\le C_1\left|x-x_1\right|+C_2\left|y-y_1\right|+
$$
\begin{equation}
\label{6666}
~~~~~+C_3\left|\left(1-x^2\right)^{1/4}\left(1-y^2\right)^{1/4}
-\left(1-x_1^2\right)^{1/4}\left(1-y_1^2\right)^{1/4}\right|,
\end{equation}

\noindent
where $C_3<\infty.$

Moreover,
$$
\left|\left(1-x^2\right)^{1/4}\left(1-y^2\right)^{1/4}
-\left(1-x_1^2\right)^{1/4}\left(1-y_1^2\right)^{1/4}\right|=
$$
$$
=\left|\left(1-x^2\right)^{1/4}\left(\left(1-y^2\right)^{1/4}-
\left(1-y_1^2\right)^{1/4}\right)
+\right.
$$
$$
\left.
+\left(1-y_1^2\right)^{1/4}\left(
\left(1-x^2\right)^{1/4}-\left(1-x_1^2\right)^{1/4}\right)\right|\le
$$
\begin{equation}
\label{5555}
~~~~~~~~ \le \left|\left(1-y^2\right)^{1/4}-\left(1-y_1^2\right)^{1/4}\right|+
\left|\left(1-x^2\right)^{1/4}-\left(1-x_1^2\right)^{1/4}\right|,
\end{equation}

\vspace{-2mm}
$$
\left|\left(1-x^2\right)^{1/4}-\left(1-x_1^2\right)^{1/4}\right|=
$$
$$
=\left|\left(\left(1-x\right)^{1/4}-
\left(1-x_1\right)^{1/4}\right)(1+x)^{1/4}+\right.
$$
$$
\left.
+\left(1-x_1\right)^{1/4}
\left((1+x)^{1/4}-\left(1+x_1\right)^{1/4}\right)\right|\le
$$
\begin{equation}
\label{4444}
~~~~~~~ \le K_1\left(\left|(1-x)^{1/4}-(1-x_1)^{1/4}\right|+
\left|(1+x)^{1/4}-(1+x_1)^{1/4}\right|\right),
\end{equation}

\noindent
where $K_1<\infty.$

It is not difficult to see that
$$
\left|(1 \pm x)^{1/4}-(1 \pm x_1)^{1/4}\right|=
$$
$$
=
\frac{|(1\pm x) - (1\pm x_1)|}
{\left((1\pm x)^{1/2}+(1 \pm x_1)^{1/2}\right)
\left((1\pm x)^{1/4}+(1 \pm x_1)^{1/4}\right)}=
$$
$$
=|x_1-x|^{1/4} \frac{|x_1-x|^{1/2}}{(1\pm x)^{1/2}+(1 \pm x_1)^{1/2}}
\cdot \frac{|x_1-x|^{1/4}}{(1\pm x)^{1/4}+(1 \pm x_1)^{1/4}}
\le 
$$
\begin{equation}
\label{3333}
\le|x_1-x|^{1/4}.
\end{equation}

The last inequality follows from the obvious inequalities
$$
\frac{|x_1-x|^{1/2}}{(1\pm x)^{1/2}+(1 \pm x_1)^{1/2}}\le 1,
$$
$$
\frac{|x_1-x|^{1/4}}{(1\pm x)^{1/4}+(1 \pm x_1)^{1/4}}\le 1.
$$

From (\ref{6666})--(\ref{3333}) we obtain 
$$
\left|g(x,y)-g(x_1,y_1)\right|\le
$$
$$
\le
C_1|x-x_1|+C_2|y-y_1|+C_4\left(|x_1-x|^{1/4}+|y_1-y|^{1/4}\right)\le
$$
$$
\le 
C_5 \rho + C_6 \rho^{1/4} \le K \rho^{1/4},
$$
where $C_5, C_6, K <\infty.$

The cases $x=x_1,$ $y\ne y_1$ and $x\ne x_1,$ $y=y_1$ 
can be considered analogously to the case
$x\ne x_1,$ $y\ne y_1$. At that, the consideration 
begins from the inequalities 
$$
|g(x,y)-g(x_1,y_1)|\le
K_2 \left|\left(1-y^2\right)^{1/4}f(x,y)-
\left(1-y_1^2\right)^{1/4}f(x_1,y_1)\right|
$$
($x=x_1,\ y\ne y_1$) and
$$
|g(x,y)-g(x_1,y_1)|\le  
K_2 \left|\left(1-x^2\right)^{1/4}f(x,y)-
\left(1-x_1^2\right)^{1/4}f(x_1,y_1)\right|
$$
($x\ne x_1,\ y=y_1$),
where $K_2<\infty$. Lemma 2.1 is proved. 

\vspace{2mm}

Lemma 2.1 and Proposition 2.2 imply
that rectangular sums of double trigonometric
Fourier series of the function $g (x,y)$ converge uniformly 
to the function $g (x,y)$
in 
the square $[-1,1]^2$. This means that
the equality (\ref{8888}) holds.

Consider the auxiliary function
\begin{equation}
\label{ziko5001}
K'(t_1,t_2)=\left\{
\begin{matrix}
\psi_2(t_1)\psi_1(t_2),\ \ t_1\ge t_2\cr\cr
\psi_1(t_1)\psi_2(t_2),\ \ t_1\le t_2
\end{matrix}
\right.,\ \ \ t_1,t_2\in[t,T]
\end{equation}
and prove that 
\begin{equation}
\label{9090}
\left|K'(t_1,t_2)-K'(t_1^{*},t_2^{*})\right|\le L
\left(|t_1-t_1^{*}|+|t_2-t_2^{*}|\right),
\end{equation}
where $L<\infty$ and $(t_1,t_2)$, $(t_1^{*},t_2^{*})\in [t, T]^2.$

By the Lagrange formula for
the functions $\psi_1(t_1^{*}),$ $\psi_2(t_1^{*})$ at the interval
$$
\left[{\rm min}\left\{t_1, t_1^{*}\right\}, 
{\rm max}\left\{t_1, t_1^{*}\right\}\right]
$$ 
and for the functions 
$\psi_1(t_2^{*}),$ $\psi_2(t_2^{*})$ at the interval
$$
\left[{\rm min}\left\{t_2, t_2^{*}\right\}, 
{\rm max}\left\{t_2, t_2^{*}\right\}\right]
$$ 
we obtain
$$
\left|K'(t_1,t_2)-K'(t_1^{*},t_2^{*})\right|\le 
$$

\vspace{-3mm}
$$
\le\left|\ \left\{
\begin{matrix}
\psi_2(t_1)\psi_1(t_2),\ \ t_1\ge t_2\cr\cr
\psi_1(t_1)\psi_2(t_2),\ \ t_1\le t_2
\end{matrix}\right.
-\left\{
\begin{matrix}
\psi_2(t_1)\psi_1(t_2),\ t_1^{*}\ge t_2^{*}\cr\cr
\psi_1(t_1)\psi_2(t_2),\ \ t_1^{*}\le t_2^{*}
\end{matrix}
\right.\ \right|+
$$
\begin{equation}
\label{uuuu1}
+
L_1\left|t_1-t_1^{*}\right|+L_2\left|t_2-t_2^{*}\right|,\ \ \ L_1, L_2<\infty.
\end{equation}

We have
$$
\left\{
\begin{matrix}
\psi_2(t_1)\psi_1(t_2),\ \ t_1\ge t_2\cr\cr
\psi_1(t_1)\psi_2(t_2),\ \ t_1\le t_2
\end{matrix}
\right.
-\left\{
\begin{matrix}
\psi_2(t_1)\psi_1(t_2),\ \ t_1^{*}\ge t_2^{*}\cr\cr
\psi_1(t_1)\psi_2(t_2),\ \ t_1^{*}\le t_2^{*}
\end{matrix}\right.
=
$$
\begin{equation}
\label{uuuu2}
=\left\{
\begin{matrix}
0,\ \ \ t_1\ge t_2,\  t_1^{*}\ge t_2^{*}\ \ \ \hbox{or}\ \ \ 
t_1\le t_2,\ t_1^{*}\le t_2^{*}\cr\cr
\psi_2(t_1)\psi_1(t_2)-\psi_1(t_1)\psi_2(t_2),\ \ \
t_1\ge t_2,\ t_1^{*}\le t_2^{*}\cr\cr
\psi_1(t_1)\psi_2(t_2)-\psi_2(t_1)\psi_1(t_2),\ \ \
t_1\le t_2,\ t_1^{*}\ge t_2^{*}
\end{matrix}.\right.
\end{equation}

\vspace{2mm}

By Lagrange formula for the functions
$\psi_1(t_2),$ $\psi_2(t_2)$ at the interval
$$
[{\rm min}\{t_1, t_2\}, {\rm max}\{t_1, t_2\}]
$$ 
we obtain the estimate
$$
\left|\ \left\{
\begin{matrix}
\psi_2(t_1)\psi_1(t_2),\ \ t_1\ge t_2\cr\cr
\psi_1(t_1)\psi_2(t_2),\ \ t_1\le t_2
\end{matrix}\right.
-\left\{
\begin{matrix}
\psi_2(t_1)\psi_1(t_2),\ \ t_1^{*}\ge t_2^{*}\cr\cr
\psi_1(t_1)\psi_2(t_2),\ \ t_1^{*}\le t_2^{*}
\end{matrix}
\right.\ \right|\le
$$
\begin{equation}
\label{ppp1ggg}
~~~~~~\le L_3 |t_2-t_1|
\left\{
\begin{matrix}
0,\ \ \ t_1\ge t_2,\ t_1^{*}\ge t_2^{*}\ \ \ \hbox{or}\ \ \
t_1\le t_2,\ t_1^{*}\le t_2^{*}\cr\cr
1,\ \ \ t_1\le t_2,\ t_1^{*}\ge t_2^{*}\ \ \ \hbox{or}\ \ \
t_1\ge t_2,\ t_1^{*}\le t_2^{*}
\end{matrix},\right.\ 
\end{equation}
where $L_3<\infty.$

Let us show that if $t_1\le t_2,$ $t_1^{*}\ge t_2^{*}$ or
$t_1\ge t_2,$ $t_1^{*}\le t_2^{*}$, then the following 
inequality is satisfied  
\begin{equation}
\label{r5r5}
|t_2-t_1|\le |t_1^{*}-t_1|+|t_2^{*}-t_2|.
\end{equation}

First, consider the case $t_1\ge t_2,$ $t_1^{*}\le t_2^{*}$.
For this case
$$
t_2+(t_1^{*}-t_2^{*})\le t_2\le t_1.
$$

Then 
$$
(t_1^{*}-t_1)-(t_2^{*}-t_2)\le t_2-t_1\le 0
$$ 
and (\ref{r5r5}) is satisfied.

For the case $t_1\le t_2,$ $t_1^{*}\ge t_2^{*}$ we obtain
$$
t_1+(t_2^{*}-t_1^{*})\le t_1\le t_2.
$$

Then 
$$
(t_1-t_1^{*})-(t_2-t_2^{*})\le t_1-t_2\le 0
$$ 
and also (\ref{r5r5}) is satisfied.

From (\ref{ppp1ggg}) and (\ref{r5r5}) we have
$$
\left|\ \left\{
\begin{matrix}
\psi_2(t_1)\psi_1(t_2),\ \ t_1\ge t_2\cr\cr
\psi_1(t_1)\psi_2(t_2),\ \ t_1\le t_2
\end{matrix}\right.
-\left\{
\begin{matrix}
\psi_2(t_1)\psi_1(t_2),\ \ t_1^{*}\ge t_2^{*}\cr\cr
\psi_1(t_1)\psi_2(t_2),\ \ t_1^{*}\le t_2^{*}
\end{matrix}
\right.\ \right|\le
$$

\vspace{-3mm}
$$
\le L_3 \left(|t_1^{*}-t_1|+|t_2^{*}-t_2|\right)
\left\{
\begin{matrix}
0,\ \ \ t_1\ge t_2,\ t_1^{*}\ge t_2^{*}\ \ \ \hbox{or}\ \ \ 
t_1\le t_2,\ t_1^{*}\le t_2^{*}\cr\cr
1,\ \ \ t_1\le t_2,\ t_1^{*}\ge t_2^{*}\ \ \ \hbox{or}\ \ \ 
t_1\ge t_2,\ t_1^{*}\le t_2^{*}
\end{matrix}\ \right. \le
$$

\vspace{-3mm}
$$
\le
L_3 \left(|t_1^{*}-t_1|+|t_2^{*}-t_2|\right)
\left\{
\begin{matrix}
1,\ \ \ t_1\ge t_2,\ t_1^{*}\ge t_2^{*}\ \ \ \hbox{or}\ \ \ 
t_1\le t_2,\ t_1^{*}\le t_2^{*}\cr\cr
1,\ \ \ t_1\le t_2,\ t_1^{*}\ge t_2^{*}\ \ \ \hbox{or}\ \ \ 
t_1\ge t_2,\ t_1^{*}\le t_2^{*}
\end{matrix}=\right. 
$$

\begin{equation}
\label{t2t2}
=L_3 \left(|t_1^{*}-t_1|+|t_2^{*}-t_2|\right).
\end{equation}

\vspace{2mm}

From (\ref{uuuu1}), (\ref{t2t2}) we obtain (\ref{9090}).
Let 
$$
t_1=\frac{T-t}{2}x+\frac{T+t}{2},\ \ \
t_2=\frac{T-t}{2}y+\frac{T+t}{2},
$$
where $x, y \in [-1,1].$ 
Then
$$
K'(t_1,t_2)\equiv K''(x,y)
=\left\{
\begin{matrix}
\psi_2\left(h(x)\right)
\psi_1\left(h(y)\right),\ \ x\ge y\cr\cr
\psi_1\left(h(x)\right)
\psi_2\left(h(y)\right),\ \ x\le y
\end{matrix}\right.,
$$
where $x, y\in[-1,1]$ and
\begin{equation}
\label{tupo12}
h(x)=\frac{T-t}{2}x+\frac{T+t}{2}.
\end{equation}

The inequality (\ref{9090}) can be rewritten in the form 
\begin{equation}
\label{4343}
\left|K''(x,y)- K''(x^{*},y^{*})\right|\le L^{*}\left(|x-x^{*}|+
|y-y^{*}|\right),
\end{equation}
where $L^{*}<\infty$ and $(x,y)$, $(x^{*},y^{*})\in [-1, 1]^2.$

Thus, the function $K'' (x,y)$ satisfies the conditions of 
Lemma 2.1. Hence, for the function 
$$
K''(x,y)\left(1-x^2\right)^{1/4}\left(1-y^2\right)^{1/4}
$$ 
the inequality (\ref{2222}) is correct.

Due to the continuous differentiability of the functions 
$\psi_1\left(h(x)\right)$ and
$\psi_2\left(h(x)\right)$ at the interval
$[-1,1]$ we have 
$K''(x,y)\in L_2([-1,1]^2)$. In addition
$$
\int\limits_{[-1,1]^2}
\frac{K''(x,y)dxdy}{(1-x^2)^{1/4}(1-y^2)^{1/4}}\le
C \left(\int\limits_{-1}^{1}\frac{1}{(1-x^2)^{1/4}}
\int\limits_{-1}^{x}\frac{1}{(1-y^2)^{1/4}}dydx+\right.
$$
$$
\left.+
\int\limits_{-1}^{1}\frac{1}{(1-x^2)^{1/4}}
\int\limits_{x}^{1}\frac{1}{(1-y^2)^{1/4}}dydx\right)<\infty,\ \ \ C<\infty.
$$

Thus, the conditions of Proposition 2.1 are fulfilled for the 
function $K''(x,y).$ 
Note that the mentioned properties of the function
$K''(x,y),$ $x,y \in [-1,1]$ also correct
for the function $K'(t_1,t_2),$ $t_1, t_2\in [t, T].$

{\bf Remark 2.1.} {\it On the basis of {\rm(\ref{9090})} it can be 
argued that the 
function $K'(t_1, t_2)$ belongs to the H\"{o}lder class 
with parameter $1$ in $[t, T]^2.$ Hence by Proposition $2.2$ this 
function can be expanded into the uniformly convergent double 
trigonometric Fourier series in the square $[t, T]^2,$ 
which summarized by Pringsheim method.
However, the expansions of iterated stochastic integrals 
obtained by using the system of
Legendre polynomials are essentially simpler than their 
analogues obtained by using the trigonometric system of functions
{\rm (}see Chapter {\rm 5} for details{\rm )}.}

Let us expand the function $K'(t_1,t_2)$ into a multiple 
(double) Fourier--Legendre series or trigonometric Fourier series 
in the square $[t, T]^2.$ This 
series is
summable
by the method of rectangular sums (Pringsheim method), i.e.
$$
K'(t_1,t_2)=
\lim_{n_1,n_2\to\infty}
\sum_{j_1=0}^{n_1}\sum_{j_2=0}^{n_2}
\int\limits_t^T\int\limits_t^TK'(t_1,t_2)
\phi_{j_1}(t_1)\phi_{j_2}(t_2)dt_1 dt_2\cdot
\phi_{j_1}(t_1)\phi_{j_2}(t_2)=
$$
$$
=\lim_{n_1,n_2\to\infty}
\sum_{j_1=0}^{n_1}\sum_{j_2=0}^{n_2}\left(
\int\limits_t^T\psi_2(t_2)\phi_{j_2}(t_2)
\int\limits_t^{t_2}\psi_1(t_1)\phi_{j_1}(t_1)dt_1dt_2+\right.
$$
$$
\left. +
\int\limits_t^T\psi_1(t_2)\phi_{j_2}(t_2)
\int\limits_{t_2}^{T}\psi_2(t_1)\phi_{j_1}(t_1)dt_1\right)dt_2
\phi_{j_1}(t_1)\phi_{j_2}(t_2)
=
$$
\begin{equation}
\label{334.ye}
=\lim_{n_1,n_2\to\infty}
\sum_{j_1=0}^{n_1}\sum_{j_2=0}^{n_2}\left(C_{j_2j_1}+
C_{j_1j_2}\right)
\phi_{j_1}(t_1)\phi_{j_2}(t_2),
\end{equation}

\vspace{3mm}
\noindent
where $(t_1, t_2)\in (t, T)^2.$
At that, the convergence of the series (\ref{334.ye}) is uniform on the
rectangle
$$
[t+\varepsilon, T-\varepsilon]\times[t+\delta,T-\delta]\ 
\hbox{for any}\ \varepsilon, \delta>0\ \hbox{(in particular, we can 
choose}\ \varepsilon=\delta).
$$

In addition, the series (\ref{334.ye}) converges to 
$K'(t_1,t_2)$ at any inner point of the square $[t, T]^2.$

Note that Proposition 2.1 does not answer the question of convergence
of the series
(\ref{334.ye}) on the boundary of the square
$[t, T]^2.$

In obtaining (\ref{334.ye}) we replaced the order of
integration in the second iterated integral.

Let us substitute  $t_1=t_2$ in (\ref{334.ye}).
After that, let us rewrite the limit on 
the right-hand side of (\ref{334.ye}) as two limits.
Let us replace $j_1$ with $j_2$, $j_2$ with $j_1$, 
$n_1$ with $n_2,$ and $n_2$ with $n_1$ in the second limit. Thus, we get

\vspace{-2mm}
\begin{equation}
\label{5656}
~~~~~~\lim\limits_{n_1,n_2\to\infty}
\sum_{j_1=0}^{n_1} \sum_{j_2=0}^{n_2} C_{j_2j_1}\phi_{j_1}(t_1)
\phi_{j_2}(t_1)=\frac{1}{2}\psi_1(t_1)\psi_2(t_1),\ \ \ t_1\in (t, T).
\end{equation}

\vspace{3mm}

According to the above reasoning, the convergence in 
(\ref{5656}) is uniform on the interval
$[t+\varepsilon, T-\varepsilon]$ for any $\varepsilon>0.$
Additionally, (\ref{5656}) holds at each interior point
of the interval $[t, T].$

Let us fix $\varepsilon>0$ and 
integrate the equality (\ref{5656}) at the interval $[t+\varepsilon,
T-\varepsilon]$. 
Due to the uniform convergence of the series (\ref{5656}) we can swap  
the series and the integral 
\begin{equation}
\label{5657}
~~~~~~~\lim\limits_{n_1,n_2\to\infty}
\sum_{j_1=0}^{n_1} \sum_{j_2=0}^{n_2} C_{j_2j_1}
\int\limits_{t+\varepsilon}^{T-\varepsilon}\phi_{j_1}(t_1)
\phi_{j_2}(t_1)dt_1=\frac{1}{2}\int\limits_{t+\varepsilon}^{T-\varepsilon}
\psi_1(t_1)\psi_2(t_1)dt_1.
\end{equation}

\vspace{2mm}

{\bf Lemma 2.2.}\ {\it 
Under the conditions of Theorem {\rm 2.2} the following limit
$$
\lim\limits_{n\to\infty}\sum\limits_{j_1=0}^{n}C_{j_1 j_1}
$$

\noindent
exists and is finite, where $C_{j_1 j_1}$ 
is defined by {\rm (\ref{tupo11})} if $j_1=j_2,$
i.e. 
$$
C_{j_1j_1}=\int\limits_t^T \psi_2(t_2)\phi_{j_1}(t_2)
\int\limits_t^{t_2} \psi_1(t_1)\phi_{j_1}(t_1)dt_1 dt_2.
$$
}

\vspace{-3mm}

Lemma~2.2 has already been proved in Sect.~2.1.1 under stronger conditions. 
Further, in this section, another proof of Lemma~2.2 is given. 
This will allow us to obtain useful 
estimates that will be used later in Chapter~2.

Applying the equality (\ref{5657})
for $n_1=n_2=n$ and Lemma~2.2, we get
$$
\frac{1}{2}\int\limits_{t+\varepsilon}^{T-\varepsilon}
\psi_1(t_1)\psi_2(t_1)dt_1=\lim\limits_{n\to\infty}
\sum_{j_1,j_2=0}^{n} C_{j_2j_1}
\int\limits_{t+\varepsilon}^{T-\varepsilon}\phi_{j_1}(t_1)
\phi_{j_2}(t_1)dt_1=
$$
$$
=\lim\limits_{n\to\infty}
\sum_{j_1,j_2=0}^{n} C_{j_2j_1}\left(
\int\limits_{t}^{T}\phi_{j_1}(t_1)
\phi_{j_2}(t_1)dt_1-\int\limits_{t}^{t+\varepsilon}\phi_{j_1}(t_1)
\phi_{j_2}(t_1)dt_1-\right.
$$
$$
\left.
-\int\limits_{T-\varepsilon}^{T}\phi_{j_1}(t_1)
\phi_{j_2}(t_1)dt_1\right)=
$$
$$
=\lim\limits_{n\to\infty}
\sum_{j_1,j_2=0}^{n} C_{j_2j_1}\left(
{\bf 1}_{\{j_1=j_2\}}-
\biggl(\phi_{j_1}(\theta)
\phi_{j_2}(\theta)+\phi_{j_1}(\lambda)\phi_{j_2}(\lambda)\biggr)
\varepsilon\right)=
$$

\vspace{-2mm}
\begin{equation}
\label{566}
~~~~~~~=
\lim\limits_{n\to\infty}\sum_{j_1=0}^{n}C_{j_1j_1}-
\varepsilon
\lim\limits_{n\to\infty}
\sum_{j_1,j_2=0}^{n}
C_{j_2j_1}\biggl(\phi_{j_1}(\theta)
\phi_{j_2}(\theta)+
\phi_{j_1}(\lambda)
\phi_{j_2}(\lambda)\biggr),
\end{equation}

\vspace{2mm}
\noindent
where $\theta\in [t,t+\varepsilon],$ $\lambda\in [T-\varepsilon,T]$.
In obtaining (\ref{566}) we used the theorem
on the mean value for the Riemann  
integral and orthonormality of the functions
$\phi_{j}(x)$ for $j=0, 1, 2\ldots$

Applying (\ref{566}), we obtain
$$
\varepsilon
\lim\limits_{n\to\infty}
\sum_{j_1,j_2=0}^{n}
C_{j_2j_1}\biggl(\phi_{j_1}(\theta)
\phi_{j_2}(\theta)+
\phi_{j_1}(\lambda)
\phi_{j_2}(\lambda)\biggr)=
$$
$$
=\lim\limits_{n\to\infty}\sum_{j_1=0}^{n}C_{j_1j_1}-
\lim\limits_{n\to\infty}
\sum_{j_1,j_2=0}^{n} C_{j_2j_1}
\int\limits_{t+\varepsilon}^{T-\varepsilon}\phi_{j_1}(t_1)
\phi_{j_2}(t_1)dt_1,
$$

\vspace{2mm}
\noindent
where the limits
$$
\lim\limits_{n\to\infty}\sum_{j_1=0}^{n}C_{j_1j_1},\ \ \ 
\lim\limits_{n\to\infty}
\sum_{j_1, j_2=0}^{n} C_{j_2j_1}
\int\limits_{t+\varepsilon}^{T-\varepsilon}\phi_{j_1}(t_1)
\phi_{j_2}(t_1)dt_1
$$

\noindent
exist and are finite (see Lemma~2.2 and the equality (\ref{5657})). 
This means that the limit 
$$
\varepsilon
\lim\limits_{n\to\infty}
\sum_{j_1,j_2=0}^{n}
C_{j_2j_1}\biggl(\phi_{j_1}(\theta)
\phi_{j_2}(\theta)+
\phi_{j_1}(\lambda)
\phi_{j_2}(\lambda)\biggr)
$$

\noindent
also exists and is finite. 

Suppose that 
the following relations
\begin{equation}
\label{rr1-six}
~~~\left\vert
\sum\limits_{j_1,j_2=0}^n
C_{j_2j_1}\phi_{j_2}(T)\phi_{j_1}(T)\right\vert\le K <\infty,\ \ \
\left\vert
\sum\limits_{j_1,j_2=0}^n
C_{j_2j_1}\phi_{j_2}(t)\phi_{j_1}(t)\right\vert\le K<\infty
\end{equation}
are satisfied for $n\in{\bf N}$ (the relations (\ref{rr1-six})
will be proved further in this sec\-tion); constant $K$
does not depend on $n$.

Note that
$$
\left\vert\varepsilon
\lim\limits_{n\to\infty}
\sum_{j_1,j_2=0}^{n}
C_{j_2j_1}\biggl(\phi_{j_1}(\theta)
\phi_{j_2}(\theta)+
\phi_{j_1}(\lambda)
\phi_{j_2}(\lambda)\biggr)\right\vert=
$$
\begin{equation}
\label{seven-1000}
~~~~~~=\lim\limits_{n\to\infty} \varepsilon  \left\vert
\sum_{j_1,j_2=0}^{n}
C_{j_2j_1}\phi_{j_1}(\theta)
\phi_{j_2}(\theta)+
\sum_{j_1,j_2=0}^{n}
C_{j_2j_1}
\phi_{j_1}(\lambda)
\phi_{j_2}(\lambda)\right\vert.
\end{equation}

\vspace{2mm}

Using (\ref{5656}) ($n_1=n_2=n$) and (\ref{rr1-six}), we obtain
$$
\varepsilon  \lim\limits_{n\to\infty} \left\vert
\sum_{j_1,j_2=0}^{n}
C_{j_2j_1}\phi_{j_1}(\theta)
\phi_{j_2}(\theta)+
\sum_{j_1,j_2=0}^{n}
C_{j_2j_1}
\phi_{j_1}(\lambda)
\phi_{j_2}(\lambda)\right\vert\le
$$
\begin{equation}
\label{sixth-1}
\le \varepsilon \lim\limits_{n\to\infty} \left(
\left\vert
\sum_{j_1,j_2=0}^{n}
C_{j_2j_1}\phi_{j_1}(\theta)
\phi_{j_2}(\theta)\right\vert+
\left\vert\sum_{j_1,j_2=0}^{n}
C_{j_2j_1}
\phi_{j_1}(\lambda)
\phi_{j_2}(\lambda)\right\vert\right)\le 2\varepsilon K_1\ \to\ 0
\end{equation}

\noindent
if $\varepsilon\to +0,$ where 
$\theta\in [t,t+\varepsilon],$ $\lambda\in [T-\varepsilon,T]$,
constant $K_1$ is independent on $n.$

Performing the 
passage to the limit $\lim\limits_{\varepsilon\to +0}$
in the equality (\ref{566}) and taking into account (\ref{seven-1000}), (\ref{sixth-1}), we get 
\begin{equation}
\label{seven-1}
\frac{1}{2}\int\limits_t^T\psi_1(t_1)\psi_2(t_1)dt_1
=\sum_{j_1=0}^{\infty}
C_{j_1j_1}.
\end{equation}

Thus, to complete the proof of Theorem 2.2, it is necessary to prove 
(\ref{rr1-six}).
To prove (\ref{rr1-six}), as well as for further consideration, we need 
some well known properties 
of the Legendre polynomials \cite{Gob}, \cite{suet}.

The complete orthonormal system
of Legendre polynomials in the space $L_2([t, T])$
looks as follows

\vspace{-5mm}
\begin{equation}
\label{4009}
~~~~~~\phi_j(x)=\sqrt{\frac{2j+1}{T-t}}P_j\biggl(\biggl(
x-\frac{T+t}{2}\biggr)\frac{2}{T-t}\biggr),\ \ \ j=0, 1, 2,\ldots,
\end{equation}

\vspace{2mm}
\noindent
where $P_j(x)$ is the Legendre polynomial. 

It is known that the Legendre polynomial $P_j (x)$ is represented as 
$$
P_j(x)=\frac{1}{2^j j!} \frac{d^j}{dx^j}\left(x^2-1\right)^j.
$$

At the boundary points of the orthogonality interval the Legendre 
polynomials satisfy the following relations
$$
P_j(1)=1,\ \ \ P_j(-1)=(-1)^j,
$$
$$
P_{j+1}(1)-P_j(1)=0,\ \ \ P_{j+1}(-1)+P_j(-1)=0,
$$

\noindent
where $j=0, 1, 2, \ldots$

Relation of the Legendre polynomial $P_j(x)$ with derivatives 
of the Legendre polynomials $P_{j+1}(x)$ and $P_{j-1}(x)$ is expressed 
by the following equality
\begin{equation}
\label{seven-3}
P_j(x)=\frac{1}{2j+1}\left(P_{j+1}^{'}(x)-
P_{j-1}^{'}(x)\right),\ \ \
j=1, 2,\ldots
\end{equation}

The recurrent relation has the form

\vspace{-2mm}
$$
xP_{j}(x)=\frac{(j+1)P_{j+1}(x)+jP_{j-1}(x)}{2j+1},\ \ \ j=1, 2,\ldots
$$

\vspace{2mm}

Orthogonality of the Legendre polynomial $P_j(x)$ to any polynomial
$Q_k (x)$ of lesser degree $k$ we write in the following form
$$
\int\limits_{-1}^1 Q_k(x) P_j(x)dx=0,\ \ \ k=0, 1, 2,\ldots,j-1.
$$

From the property
$$
\int\limits_{-1}^1 P_k(x)P_j(x)dx=\left\{
\begin{matrix}
0\ \ &\hbox{if}\ \
k\ne j\cr\cr
2/(2j+1)\ \ &\hbox{if}\ \ k=j
\end{matrix}
\right.
$$

\vspace{2mm}
\noindent
it follows that the orthonormal on the interval $[-1,1]$ Legendre polynomials 
determined by the relation

\vspace{-1mm}
$$
P_j^{*}(x)=\sqrt{\frac{2j+1}{2}}P_j(x),\ \ \ j=0, 1, 2,\ldots
$$

\vspace{3mm}

Remind that there is the following estimate \cite{Gob}

\vspace{-4mm}
\begin{equation}
\label{otit987ggg}
~~~~~~~\left|P_j(y)\right| <\frac{K}{\sqrt{j+1}(1-y^2)^{1/4}},\ \ \ 
y\in (-1, 1),\ \ \ j=1, 2,\ldots,
\end{equation}

\vspace{3mm}
\noindent
where constant $K$ does not depend on $y$ and $j$.

Moreover,

\vspace{-4mm}
\begin{equation}
\label{q1}
\left|P_j(x)\right|\le 1,\ \ \ x\in[-1,1],\ \ \ j=0,\ 1,\ldots
\end{equation}

\vspace{4mm}
The Christoffel--Darboux formula has the form

\vspace{-3mm}
\begin{equation}
\label{ty}
~~~~~~~~~~\sum\limits_{j=0}^n (2j+1)P_j(x)P_j(y)=
(n+1)\frac{P_n(x)P_{n+1}(y)-P_{n+1}(x)P_{n}(y)}{y-x}.
\end{equation}

\vspace{1mm}

Let us prove (\ref{rr1-six}) (see \cite{art-9})).
From (\ref{ty}) for $x=\pm 1$ we obtain
\begin{equation}
\label{yy3-six}
\sum\limits_{j=0}^n (2j+1)P_j(y)=
(n+1)\frac{P_{n+1}(y)-P_{n}(y)}{y-1},
\end{equation}
\begin{equation}
~~~~~~~~~~\label{yy4-six}
\sum\limits_{j=0}^n (2j+1)(-1)^{j}P_j(y)=
(n+1)(-1)^n\frac{P_{n+1}(y)+P_{n}(y)}{y+1}.
\end{equation}

\vspace{3mm}

From the other hand (see (\ref{seven-3}))
$$
\sum\limits_{j=0}^n (2j+1)P_j(y)=1+\sum\limits_{j=1}^n (2j+1)P_j(y)=
$$
$$
=1+\sum\limits_{j=1}^n (P_{j+1}^{'}(y)-P_{j-1}^{'}(y))=
1+\biggl(\sum\limits_{j=1}^n (P_{j+1}(y)-P_{j-1}(y))\biggr)'=
$$

\vspace{1mm}
\begin{equation}
\label{yy2-six}
=1+(P_{n+1}(x)+P_n(x)-x-1)'=(P_n(x)+P_{n+1}(x))'
\end{equation}

\vspace{3mm}
\noindent
and
$$
\sum\limits_{j=0}^n (2j+1)(-1)^j P_j(y)=1+\sum\limits_{j=1}^n (-1)^j
(2j+1)P_j(y)=
$$
$$
=1+\sum\limits_{j=1}^n (-1)^j (P_{j+1}^{'}(y)-P_{j-1}^{'}(y))=
1+\biggl(\sum\limits_{j=1}^n (-1)^j(P_{j+1}(y)-P_{j-1}(y))\biggr)'=
$$

\vspace{-3mm}
\begin{equation}
\label{yy1-six}
~~~~~=1+((-1)^n(P_{n+1}(x)-P_n(x))-x+1)'=(-1)^n(P_{n+1}(x)-P_{n}(x))'.
\end{equation}

\vspace{4mm}

Applying (\ref{yy3-six})--(\ref{yy1-six}), we get

\begin{equation}
\label{uu1-six}
(n+1)\frac{P_{n+1}(y)-P_{n}(y)}{y-1}=(P_n(x)+P_{n+1}(x))',
\end{equation}

\vspace{1mm}
\begin{equation}
\label{uu2-six}
(n+1)\frac{P_{n+1}(y)+P_{n}(y)}{y+1}
=(P_{n+1}(x)-P_{n}(x))'.
\end{equation}

\vspace{4mm}

Let us prove the boundedness of the first sum in (\ref{rr1-six}). 
We have
$$
\sum\limits_{j_1,j_2=0}^n
C_{j_2j_1}\phi_{j_2}(T)\phi_{j_1}(T)=
$$
$$
=\frac{1}{4}
\sum\limits_{j_2=0}^n
\sum\limits_{j_1=0}^n
(2j_2+1)(2j_1+1)\int\limits_{-1}^{1}
\psi_2(h(y))P_{j_2}(y)
\int\limits_{-1}^{y}
\psi_1(h(y_1))P_{j_1}(y_1)dy_1 dy=
$$
$$
=\frac{1}{4}
\int\limits_{-1}^{1}
\psi_2(h(y))
\sum\limits_{j_2=0}^n
(2j_2+1)P_{j_2}(y)
\int\limits_{-1}^{y}
\psi_1(h(y_1))\sum\limits_{j_1=0}^n
(2j_1+1)P_{j_1}(y_1)dy_1 dy=
$$
$$
=\frac{1}{4}
\int\limits_{-1}^{1}
\psi_2(h(y))
\left(\int\limits_{-1}^{y}
\psi_1(h(y_1))
d(P_{n+1}(y_1)+P_{n}(y_1))\right)d(P_{n+1}(y)+P_{n}(y))=
$$
$$
=\frac{1}{4}
\int\limits_{-1}^{1}
\psi_1(h(y))
\left(\int\limits_{-1}^{y}
\psi_1(h(y_1))
d(P_{n+1}(y_1)+P_{n}(y_1))\right)d(P_{n+1}(y)+P_{n}(y))+
$$
$$
+\frac{1}{4}
\int\limits_{-1}^{1}
\Delta(h(y))
\left(\int\limits_{-1}^{y}
\psi_1(h(y_1))
d(P_{n+1}(y_1)+P_{n}(y_1))\right)d(P_{n+1}(y)+P_{n}(y))=
$$
$$
=\frac{1}{4}I_1+\frac{1}{4}I_2,
$$

\vspace{2mm}
\noindent
where
\begin{equation}
\label{191-six}
~~~~~~~~\Delta(h(y))=\psi_2(h(y))-\psi_1(h(y)),\ \ \ h(y)=\frac{T-t}{2}y+\frac{T+t}{2}.
\end{equation}

\vspace{2mm}

Further,
$$
I_1=\frac{1}{2}
\biggl(\int\limits_{-1}^{1}
\psi_1(h(y))
d(P_{n+1}(y)+P_{n}(y))\biggr)^2=
$$
$$
=\frac{1}{2}\biggl(2\psi_1(T)-\int\limits_{-1}^{1}
(P_{n+1}(y)+P_{n}(y))\psi_1'(h(y))\frac{T-t}{2}dy\biggr)^2<C_1<\infty,
$$

\vspace{2mm}
\noindent
where $\psi_1'$ is a derivative 
of the function 
$\psi_1$ with respect to the variable $y,$
constant $C_1$ does not depend on $n$.

By the Lagrange formula we obtain
$$
\Delta(h(y))=\psi_2\biggl(\frac{1}{2}(T-t)(y-1)+T\biggr)-
\psi_1\biggl(\frac{1}{2}(T-t)(y-1)+T\biggr)=
$$
$$
=\psi_2(T)-\psi_1(T)+(y-1)\biggl(\psi_2'(\xi_y)-\psi_1'(\theta_y)\biggr)
\frac{1}{2}(T-t)=
$$

\vspace{-1mm}
\begin{equation}
\label{hu-six}
=C_1+\alpha_y(y-1),
\end{equation}

\vspace{3mm}
\noindent
where $|\alpha_y|<\infty$ and $C_1=\psi_2(T)-\psi_1(T).$

Let us substitute (\ref{hu-six}) into the integral $I_2$

\vspace{-3mm}
$$
I_2=I_3+I_4,
$$
where
$$
I_3=\int\limits_{-1}^{1}
\alpha_y(y-1)
\left(\int\limits_{-1}^{y}
\psi_1(h(y_1))
d(P_{n+1}(y_1)+P_{n}(y_1))\right)d(P_{n+1}(y)+P_{n}(y)),
$$
$$
I_4=
C_1\int\limits_{-1}^{1}
\left(\int\limits_{-1}^{y}
\psi_1(h(y_1))
d(P_{n+1}(y_1)+P_{n}(y_1))\right)d(P_{n+1}(y)+P_{n}(y)).
$$

\vspace{4mm}

Integrating by parts and using (\ref{uu1-six}), we obtain

\vspace{-4mm}
$$
I_3=
\int\limits_{-1}^{1}
\frac{\alpha_y(y-1)(n+1)(P_{n+1}(y)-P_{n}(y))}{y-1}
\biggl(
\psi_1(h(y))(P_{n+1}(y)+P_{n}(y))-\biggr.
$$
$$
\biggl.-
\int\limits_{-1}^{y}
(P_{n+1}(y_1)+P_{n}(y_1))\psi_1'(h(y_1))\frac{1}{2}(T-t)dy_1\biggr)dy.
$$

\vspace{2mm}

Applying the etimate
(\ref{otit987ggg}) and taking into account the boundedness
of $\alpha_y$ and $\psi_1'(h(y_1))$, we have
that $|I_3|<\infty.$

Using the integration order replacement in $I_4,$ we get 

\vspace{-4mm}
$$
I_4=
C_1\int\limits_{-1}^{1}\psi_1(h(y_1))
\left(\int\limits_{y_1}^{1}
d(P_{n+1}(y)+P_{n}(y))\right)d(P_{n+1}(y_1)+P_{n}(y_1))=
$$
$$
=C_1\int\limits_{-1}^{1}\psi_1(h(y_1))d(P_{n+1}(y_1)+P_{n}(y_1))
\int\limits_{-1}^{1}
d(P_{n+1}(y)+P_{n}(y))-
$$
$$
-C_1\int\limits_{-1}^{1}\psi_1(h(y_1))
\left(\int\limits_{-1}^{y_1}
d(P_{n+1}(y)+P_{n}(y))\right)d(P_{n+1}(y_1)+P_{n}(y_1))=
$$

\vspace{-1mm}
$$
=I_5-I_6.
$$

\vspace{2mm}

Consider $I_5$
$$
I_5=2C_1\int\limits_{-1}^{1}\psi_1(h(y_1))d(P_{n+1}(y_1)+P_{n}(y_1))=
$$
$$
=2C_1\left(2\psi_1(T)-
\int\limits_{-1}^{1}
(P_{n+1}(y_1)+P_{n}(y_1))
\psi_1'(h(y_1))\frac{1}{2}(T-t)dy_1\right).
$$

\vspace{3mm} 

Applying the estimate (\ref{q1}) and using the boundedness of
$\psi_1'(h(y_1))$, we obtain
that $|I_5|<\infty.$

Since (see (\ref{hu-six}))
$$
\psi_1(h(y))=
\psi_1\biggl(\frac{1}{2}(T-t)(y-1)+T\biggr)=
$$
$$
=\psi_1(T)+(y-1)\psi_1'(\theta_y)
\frac{1}{2}(T-t)=C_2+\beta_y(y-1),
$$

\vspace{4mm}
\noindent
where $|\beta_y|<\infty$ and $C_2=\psi_1(T),$ then
$$
I_6=
C_3\int\limits_{-1}^{1}
\left(\int\limits_{-1}^{y_1}
d(P_{n+1}(y)+P_{n}(y))\right)d(P_{n+1}(y_1)+P_{n}(y_1))+
$$
$$
+C_1\int\limits_{-1}^{1}\beta_{y_1}(y_1-1)
\left(\int\limits_{-1}^{y_1}
d(P_{n+1}(y)+P_{n}(y))\right)d(P_{n+1}(y_1)+P_{n}(y_1))=
$$
$$
=\frac{C_3}{2}\left(\int\limits_{-1}^{1}
d(P_{n+1}(y)+P_{n}(y))\right)^2+
$$
$$
+C_1\int\limits_{-1}^{1}\frac{\beta_{y_1}(y_1-1)
(n+1)(P_{n+1}(y_1)-P_{n}(y_1))}{y_1-1}
\left(\int\limits_{-1}^{y_1}
d(P_{n+1}(y)+P_{n}(y))\right)dy_1=
$$
$$
=2C_3+
C_1\int\limits_{-1}^{1}\beta_{y_1}
(n+1)(P_{n+1}(y_1)-P_{n}(y_1))
(P_{n+1}(y_1)+P_{n}(y_1))dy_1.
$$

Using the estimate (\ref{otit987ggg}) and 
taking into account the bounedness of $\beta_{y_1}$, we obtain
that $|I_6|<\infty.$
Thus, the boundedness of the first sum in (\ref{rr1-six}) is proved.

Let us prove the boundedness of the second sum in (\ref{rr1-six}). 
We have
$$
\sum\limits_{j_1,j_2=0}^n
C_{j_2j_1}\phi_{j_2}(t)\phi_{j_1}(t)=
$$
$$
=\frac{1}{4}
\sum\limits_{j_2=0}^n
\sum\limits_{j_1=0}^n
(2j_2+1)(2j_1+1)(-1)^{j_1+j_2}\int\limits_{-1}^{1}
\psi_2(h(y))P_{j_2}(y)
\int\limits_{-1}^{y}
\psi_1(h(y_1))P_{j_1}(y_1)\times
$$
$$
\times
dy_1 dy=
$$
$$
=\frac{1}{4}
\int\limits_{-1}^{1}
\psi_2(h(y))
\sum\limits_{j_2=0}^n
(2j_2+1)P_{j_2}(y)(-1)^{j_2}
\int\limits_{-1}^{y}
\psi_1(h(y_1))\times
$$
$$
\times
\sum\limits_{j_1=0}^n
(2j_1+1)P_{j_1}(y_1)(-1)^{j_1}
dy_1 dy=
$$
$$
=\frac{(-1)^{2n}}{4}
\int\limits_{-1}^{1}
\psi_2(h(y))
\left(\int\limits_{-1}^{y}
\psi_1(h(y_1))
d(P_{n+1}(y_1)-P_{n}(y_1))\right)\times
$$

\vspace{-1mm}
$$
\times 
d(P_{n+1}(y)-P_{n}(y))=
$$
$$
=\frac{1}{4}
\int\limits_{-1}^{1}
\psi_1(h(y))
\left(\int\limits_{-1}^{y}
\psi_1(h(y_1))
d(P_{n+1}(y_1)-P_{n}(y_1))\right)d(P_{n+1}(y)-P_{n}(y))+
$$
$$
+\frac{1}{4}
\int\limits_{-1}^{1}
\Delta(h(y))
\left(\int\limits_{-1}^{y}
\psi_1(h(y_1))
d(P_{n+1}(y_1)-P_{n}(y_1))\right)d(P_{n+1}(y)-P_{n}(y))=
$$
$$
=\frac{1}{4}J_1+\frac{1}{4}J_2,
$$

\vspace{2mm}
\noindent
where $\Delta(h(y)),$
$h(y)$ are defined by (\ref{191-six}).

Further,
$$
J_1=\frac{1}{2}
\left(\int\limits_{-1}^{1}
\psi_1(h(y))
d(P_{n+1}(y)-P_{n}(y))\right)^2=
$$
\begin{equation}
\label{hh1-six}
=\frac{1}{2}\left(2(-1)^{n}\psi_1(t)-\int\limits_{-1}^{1}
(P_{n+1}(y)-P_{n}(y))\psi_1'(h(y))\frac{T-t}{2}dy\right)^2<K_1<\infty,
\end{equation}
where $\psi_1'$ is a derivative of the function
$\psi_1$ with respect to the variable $y,$
constant $K_1$ is independent of $n$.

By the Lagrange formula we obtain
$$
\Delta(h(y))=\psi_2\biggl(\frac{1}{2}(T-t)(y+1)+t\biggr)-
\psi_1\biggl(\frac{1}{2}(T-t)(y+1)+t\biggr)=
$$
$$
=\psi_2(t)-\psi_1(t)+(y+1)\biggl(\psi_2'(\mu_y)-\psi_1'(\rho_y)\biggr)
\frac{1}{2}(T-t)=
$$

\vspace{-1mm}
\begin{equation}
\label{hu11-six}
=K_2+\gamma_y(y+1),
\end{equation}

\vspace{3mm}
\noindent
where $|\gamma_y|<\infty$ and $K_2=\psi_2(t)-\psi_1(t).$

Consider $J_2$
$$
J_2=
\int\limits_{-1}^{1}
\Delta(h(y))
d(P_{n+1}(y)-P_{n}(y))
\int\limits_{-1}^{1}
\psi_1(h(y_1))
d(P_{n+1}(y_1)-P_{n}(y_1))-
$$
$$
-
\int\limits_{-1}^{1}
\Delta(h(y))
\left(\int\limits_{y}^{1}
\psi_1(h(y_1))
d(P_{n+1}(y_1)-P_{n}(y_1))\right)d(P_{n+1}(y)-P_{n}(y))=
$$
$$
=J_3J_4-J_5.
$$

\vspace{2mm}

The integral $J_4$ was considered earlier (see $J_1$ and (\ref{hh1-six})), i.e.
it has already been shown that
$|J_4|<\infty$. Analogously, we have that $|J_3|<\infty$.

Let us substitute (\ref{hu11-six}) into the integral $J_5$

\vspace{-3mm}
$$
J_5=J_6+J_7,
$$
where
$$
J_6=\int\limits_{-1}^{1}
\gamma_y(y+1)
\left(\int\limits_{y}^{1}
\psi_1(h(y_1))
d(P_{n+1}(y_1)-P_{n}(y_1))\right)d(P_{n+1}(y)-P_{n}(y)),
$$
$$
J_7=
K_2\int\limits_{-1}^{1}
\left(\int\limits_{y}^{1}
\psi_1(h(y_1))
d(P_{n+1}(y_1)-P_{n}(y_1))\right)d(P_{n+1}(y)-P_{n}(y)).
$$

\vspace{3mm}

Integrating by parts and using (\ref{uu2-six}), we get
$$
J_6=
\int\limits_{-1}^{1}
\frac{\gamma_y(y+1)(n+1)(P_{n+1}(y)+P_{n}(y))}{y+1}
\biggl(
-\psi_1(h(y))(P_{n+1}(y)-P_{n}(y))-\biggr.
$$
$$
\biggl.-
\int\limits_{y}^{1}
(P_{n+1}(y_1)-P_{n}(y_1))\psi_1'(h(y_1))\frac{1}{2}(T-t)dy_1\biggr)dy.
$$

\vspace{1mm}

Applying the etimate
(\ref{otit987ggg}) and taking into account the boundedness
of $\gamma_y$ and $\psi_1'(h(y_1))$, we have
that $|J_6|<\infty.$

Using the integration order replacement in $J_7,$ we obtain
$$
J_7=
K_2\int\limits_{-1}^{1}
\psi_1(h(y_1))
\left(\int\limits_{-1}^{y_1}
d(P_{n+1}(y)-P_{n}(y))\right)d(P_{n+1}(y_1)-P_{n}(y_1))=
$$
$$
=K_2\int\limits_{-1}^{1}
\psi_1(h(y_1))d(P_{n+1}(y_1)-P_{n}(y_1))
\int\limits_{-1}^{1}
d(P_{n+1}(y)-P_{n}(y))-K_2J_8=
$$

\vspace{-1mm}
$$
=K_2 J_4 2(-1)^{n} - K_2J_8,
$$

\vspace{2mm}
\noindent
where
$$
J_8=\int\limits_{-1}^{1}
\psi_1(h(y_1))
\left(\int\limits_{y_1}^{1}
d(P_{n+1}(y)-P_{n}(y))\right)d(P_{n+1}(y_1)-P_{n}(y_1)).
$$

\vspace{2mm}

Since (see (\ref{hu11-six}))
$$
\psi_1(h(y))=\psi_1\biggl(\frac{1}{2}(T-t)(y+1)+t\biggr)=
$$
\begin{equation}
\label{ooo1}
=\psi_1(t)+(y+1)\psi_1'(\rho_y)
\frac{1}{2}(T-t)=K_3+\varepsilon_y(y+1),
\end{equation}

\vspace{3mm}
\noindent
where $|\varepsilon_y|<\infty$ and $K_3=\psi_1(t),$ then
$$
J_8=K_3\int\limits_{-1}^{1}
\left(\int\limits_{y_1}^{1}
d(P_{n+1}(y)-P_{n}(y))\right)d(P_{n+1}(y_1)-P_{n}(y_1))+
$$
$$
+\int\limits_{-1}^{1}
\varepsilon_y(y+1)
\left(\int\limits_{y_1}^{1}
d(P_{n+1}(y)-P_{n}(y))\right)d(P_{n+1}(y_1)-P_{n}(y_1))=
$$
$$
=\frac{K_3}{2}\left(\int\limits_{-1}^{1}
d(P_{n+1}(y)-P_{n}(y))\right)^2+
$$
$$
+\int\limits_{-1}^{1}
\frac{\varepsilon_{y_1}(y_1+1)(n+1)(P_{n+1}(y_1)+P_{n}(y_1))}{y_1+1}
(P_{n}(y_1)-P_{n+1}(y_1))dy=
$$
\begin{equation}
\label{qqq1}
~~~~~~~=2K_3+\int\limits_{-1}^{1}
\varepsilon_{y_1}(n+1)(P_{n+1}(y_1)+P_{n}(y_1))
(P_{n}(y_1)-P_{n+1}(y_1))dy.
\end{equation}

When obtaining the equality (\ref{qqq1}), we used (\ref{uu2-six}).
Applying the estimate (\ref{otit987ggg}) and 
taking into account the bounedness of $\varepsilon_{y_1}$, we obtain
that $|J_8|<\infty.$
Thus, the boundedness of the second sum in (\ref{rr1-six}) is proved.
The relations (\ref{rr1-six}) are proved. 
Theorem~2.2 is proved.

Let us consider the proof of Lemma~2.2 under the conditions of Theorem~2.2.
We will prove
that 
$$
\sum\limits_{j_1=0}^n C_{j_1j_1}
$$

\noindent
is the Cauchy sequence for 
the cases of Legendre polynomials and trigonometric functions.

Consider the case of Legendre polynomials. Below in this section
we write $\lim\limits_{n,m\to\infty}$ instead of 
$\lim\limits_{\stackrel{n,m\to\infty}{{}_{n>m}}}$.
Fix $n>m$ $(n, m\in {\bf N})$. We have
$$
\sum\limits_{j_1=m+1}^n
C_{j_1j_1}=
\sum\limits_{j_1=m+1}^n
\int\limits_t^T \psi_2(s)\phi_{j_1}(s)
\int\limits_t^{s} \psi_1(\tau)\phi_{j_1}(\tau)d\tau ds=
$$
$$
=
\frac{T-t}{4}
\sum\limits_{j_1=m+1}^n
(2j_1+1)\int\limits_{-1}^{1}
\psi_2(h(x))P_{j_1}(x)
\int\limits_{-1}^{x}
\psi_1(h(y))P_{j_1}(y)dy dx=
$$
$$
=
\frac{T-t}{4}
\sum\limits_{j_1=m+1}^n
\int\limits_{-1}^{1}
\psi_1(h(x))\psi_2(h(x))\left(P_{j_1+1}(x)P_{j_1}(x)
-P_{j_1}(x)P_{j_1-1}(x)\right)dx-
$$
$$
-\frac{(T-t)^2}{8}\hspace{-2mm}
\sum\limits_{j_1=m+1}^n
\int\limits_{-1}^{1}
\psi_2(h(x))P_{j_1}(x)
\int\limits_{-1}^{x}
\left(P_{j_1+1}(y)-P_{j_1-1}(y)\right)\psi_1'(h(y))dy dx=
$$
$$
=
\frac{T-t}{4}
\int\limits_{-1}^{1}
\psi_1(h(x))\psi_2(h(x))\sum\limits_{j_1=m+1}^n
\left(P_{j_1+1}(x)P_{j_1}(x)
-P_{j_1}(x)P_{j_1-1}(x)\right)dx-
$$
$$
-\frac{(T-t)^2}{8}\hspace{-2mm}
\sum\limits_{j_1=m+1}^n
\int\limits_{-1}^{1}
\left(P_{j_1+1}(y)-P_{j_1-1}(y)\right)\psi_1'(h(y))
\int\limits_{y}^{1}
P_{j_1}(x)\psi_2(h(x))dx dy=
$$
$$
=
\frac{T-t}{4}
\int\limits_{-1}^{1}
\psi_1(h(x))\psi_2(h(x))
\left(P_{n+1}(x)P_{n}(x)
-P_{m+1}(x)P_{m}(x)\right)dx+
$$
$$
+\frac{(T-t)^2}{8}
\sum\limits_{j_1=m+1}^n
\frac{1}{2j_1+1}
\int\limits_{-1}^{1}
\left(P_{j_1+1}(y)-P_{j_1-1}(y)\right)\psi_1'(h(y))\times
$$
$$
\times
\Biggl(
\left(P_{j_1+1}(y)-P_{j_1-1}(y)\right)\psi_2(h(y))+\Biggr.
$$
\begin{equation}
\label{tupo14}
\Biggl.
+
\frac{T-t}{2}
\int\limits_{y}^{1}
\left(P_{j_1+1}(x)-
P_{j_1-1}(x)\right)\psi_2'(h(x))dx\Biggr)dy,
\end{equation}

\vspace{2mm}
\noindent
where $\psi_1'$, $\psi_2'$ are
derivatives of the functions $\psi_1$, $\psi_2$ with respect 
to the variable
$h(y)$ (see (\ref{tupo12})).

Applying the estimate (\ref{otit987ggg}) and tak\-ing into account 
the boundedness of the functions $\psi_1(\tau)$, $\psi_2(\tau)$
and their derivatives, we finally obtain
$$
\left\vert\sum\limits_{j_1=m+1}^n
C_{j_1j_1}\right\vert\le
C_1\left(\frac{1}{n}+\frac{1}{m}\right)
\int\limits_{-1}^1 \frac{dx}{\left(1-x^2\right)^{1/2}}+
$$
$$
+C_2 \sum\limits_{j_1=m+1}^n \frac{1}{j_1^2}\left(
\int\limits_{-1}^1 \frac{dy}{\left(1-y^2\right)^{1/2}}+
\int\limits_{-1}^1 \frac{1}{\left(1-y^2\right)^{1/4}}
\int\limits_{y}^1 \frac{dx}{\left(1-x^2\right)^{1/4}}dy\right)\le
$$
\begin{equation}
\label{tupo15}
\le C_3\left(\frac{1}{n}+\frac{1}{m}+\sum\limits_{j_1=m+1}^n 
\frac{1}{j_1^2}\right) \to 0
\end{equation}

\vspace{2mm}
\noindent
if $n, m\to\infty\ (n>m),$
where constants $C_1, C_2, C_3$ do not depend on $n$ and $m$.

Now consider the 
trigonometric case. Fix $n>m$ $(n, m\in {\bf N}).$
Denote
$$
S_{n,m}\stackrel{\sf def}{=}
\sum\limits_{j_1=m+1}^n
C_{j_1j_1}=
\sum\limits_{j_1=m+1}^n
\int\limits_t^T \psi_2(t_2)\phi_{j_1}(t_2)
\int\limits_t^{t_2} \psi_1(t_1)\phi_{j_1}(t_1)dt_1 dt_2.
$$

By analogy with (\ref{tupo14})
we obtain
$$
S_{2n,2m}=
\sum\limits_{j_1=2m+1}^{2n}
\int\limits_t^T \psi_2(t_2)\phi_{j_1}(t_2)
\int\limits_t^{t_2} \psi_1(t_1)\phi_{j_1}(t_1)dt_1 dt_2=
$$
$$
=\frac{2}{T-t}\sum\limits_{j_1=m+1}^n\left(
\int\limits_t^T \psi_2(t_2){\rm sin}\frac{2\pi j_1(t_2-t)}{T-t}
\int\limits_t^{t_2}\psi_1(t_1){\rm sin}\frac{2\pi j_1(t_1-t)}{T-t}dt_1 dt_2+
\right.
$$
$$
\left.+
\int\limits_t^T \psi_2(t_2){\rm cos}\frac{2\pi j_1(t_2-t)}{T-t}
\int\limits_t^{t_2}\psi_1(t_1){\rm cos}\frac{2\pi j_1(t_1-t)}{T-t}dt_1 dt_2\right)=
$$
$$
=\frac{T-t}{2\pi^2}\sum\limits_{j_1=m+1}^n\frac{1}{j_1^2}
\left(\psi_1(t)\left(\psi_2(t)-\psi_2(T)+
\int\limits_t^{T}\psi_2'(t_2){\rm cos}\frac{2\pi j_1(t_2-t)}{T-t}dt_2\right)-\right.
$$
$$
-\int\limits_t^{T}\psi_1'(t_1){\rm cos}\frac{2\pi j_1(t_1-t)}{T-t}\Biggl(
\psi_2(T)-\psi_2(t_1){\rm cos}\frac{2\pi j_1(t_1-t)}{T-t}-\Biggr.
$$
$$
\Biggl.-\int\limits_{t_1}^{T}\psi_2'(t_2){\rm cos}\frac{2\pi j_1(t_2-t)}{T-t}dt_2\Biggr)dt_1+
$$
$$
+\int\limits_t^{T}\psi_1'(t_1){\rm sin}\frac{2\pi j_1(t_1-t)}{T-t}\Biggl(
\psi_2(t_1){\rm sin}\frac{2\pi j_1(t_1-t)}{T-t}+\Biggr.
$$
\begin{equation}
\label{agentaa1000}
\left.\Biggl.
+\int\limits_{t_1}^{T}\psi_2'(t_2){\rm sin}\frac{2\pi j_1(t_2-t)}{T-t}dt_2\Biggr)dt_1\right),
\end{equation}

\noindent
where $\psi_1'(\tau),$ $\psi_2'(\tau)$ are
derivatives of the functions $\psi_1(\tau),$ $\psi_2(\tau)$ with respect to the variable
$\tau$.

From (\ref{agentaa1000}) we get
\begin{equation}
\label{agentaa1001}
\left|S_{2n,2m}\right|\le
C \sum\limits_{j_1=m+1}^n 
\frac{1}{j_1^2} \to 0
\end{equation}

\noindent
if $n, m\to\infty\ (n>m),$
where constant $C$ does not depend on $n$ and $m$.

Further,
$$
S_{2n-1,2m}=S_{2n,2m}-
$$
\begin{equation}
\label{agentaa1002}
-\frac{2}{T-t}
\int\limits_t^T \psi_2(t_2){\rm cos}\frac{2\pi n(t_2-t)}{T-t}
\int\limits_t^{t_2}\psi_1(t_1){\rm cos}\frac{2\pi n(t_1-t)}{T-t}dt_1 dt_2,
\end{equation}

$$
S_{2n,2m-1}=S_{2n,2m}+
$$
\begin{equation}
\label{agentaa1003}
+\frac{2}{T-t}
\int\limits_t^T \psi_2(t_2){\rm cos}\frac{2\pi m(t_2-t)}{T-t}
\int\limits_t^{t_2}\psi_1(t_1){\rm cos}\frac{2\pi m(t_1-t)}{T-t}dt_1 dt_2,
\end{equation}

$$
S_{2n-1,2m-1}=S_{2n,2m-1}-
$$
$$
-\frac{2}{T-t}
\int\limits_t^T \psi_2(t_2){\rm cos}\frac{2\pi n(t_2-t)}{T-t}
\int\limits_t^{t_2}\psi_1(t_1){\rm cos}\frac{2\pi n(t_1-t)}{T-t}dt_1 dt_2=
$$
$$
=S_{2n,2m}+
\frac{2}{T-t}
\int\limits_t^T \psi_2(t_2){\rm cos}\frac{2\pi m(t_2-t)}{T-t}
\int\limits_t^{t_2}\psi_1(t_1){\rm cos}\frac{2\pi m(t_1-t)}{T-t}dt_1 dt_2-
$$
\begin{equation}
\label{agentaa1004}
~~~~-\frac{2}{T-t}
\int\limits_t^T \psi_2(t_2){\rm cos}\frac{2\pi n(t_2-t)}{T-t}
\int\limits_t^{t_2}\psi_1(t_1){\rm cos}\frac{2\pi n(t_1-t)}{T-t}dt_1 dt_2.
\end{equation}

\vspace{2mm}

Integrating by parts in (\ref{agentaa1002})--(\ref{agentaa1004}), we obtain
\begin{equation}
\label{agentaa1005}
\left|S_{2n-1,2m}\right|\le \left|S_{2n,2m}\right|+\frac{C_1}{n},
\end{equation}
\begin{equation}
\label{agentaa1006}
\left|S_{2n,2m-1}\right|\le \left|S_{2n,2m}\right|+\frac{C_1}{m},
\end{equation}

\vspace{-3mm}
\begin{equation}
\label{agentaa1007}
\left|S_{2n-1,2m-1}\right|\le \left|S_{2n,2m}\right|+C_1\left(\frac{1}{m}+\frac{1}{n}\right),
\end{equation}

\vspace{3mm}
\noindent
where constant $C_1$ does not depend on $n$ and $m$.

The relations (\ref{agentaa1001}), (\ref{agentaa1005})--(\ref{agentaa1007}) imply that
\begin{equation}
\label{agentaa1008}
\lim\limits_{n,m\to\infty}\left|S_{2n,2m}\right|=
\lim\limits_{n,m\to\infty}\left|S_{2n-1,2m}\right|=
\lim\limits_{n,m\to\infty}\left|S_{2n,2m-1}\right|=
\lim\limits_{n,m\to\infty}\left|S_{2n-1,2m-1}\right|=0.
\end{equation}

From (\ref{agentaa1008}) we get
\begin{equation}
\label{agentaa1009}
\lim\limits_{n,m\to\infty}\left|S_{n,m}\right|=0.
\end{equation}

The relation (\ref{agentaa1009}) completes the proof.

\subsection{Approach Based on Generalized Double Multiple and 
Iterated Fourier Series}

This section is devoted to the proof of Theorem 2.1 
using a slightly different method than the method proposed in Sect.~2.1.1.
We will consider two different parts of the expansion of iterated
Stra\-to\-no\-vich stochastic integrals of second multiplicity.
The mean-square convergence of the first part will be proved on the base
of generalized multiple Fourier series converging 
in the mean-square sense in the space 
$L_2([t, T]^2).$ The mean-square convergence
of the second part will be proved on the base of 
generalized iterated (double) Fourier
series converging pointwise.

{\bf Proof.} Let us consider
Lemma 1.1, definition of the multiple stochastic integral (\ref{30.34})
together with the
formula (\ref{30.52}) when the function $\Phi(t_1,\ldots,t_k)$ is 
continuous in the open domain $D_k$ and bounded at its boundary
as well as
Lemma 1.3 for the case $k=2$ (see Sect.~1.1.3).

In accordance to the standard relation between
Stratonovich and It\^{o} sto\-chas\-tic integrals (see (\ref{oop51}))
we have w.~p.~1
\begin{equation}
\label{oop51xyx}
~~~~~~~~~ J^{*}[\psi^{(2)}]_{T,t}=
J[\psi^{(2)}]_{T,t}+
\frac{1}{2}{\bf 1}_{\{i_1=i_2\}}
\int\limits_t^T\psi_1(t_1)\psi_2(t_1)dt_1.
\end{equation}

Let us consider the function $K^{*}(t_1,t_2)$ defined by (\ref{yes2002})
\begin{equation}
\label{1999.1}
K^{*}(t_1,t_2)
=K(t_1,t_2)+\frac{1}{2}{\bf 1}_{\{t_1=t_{2}\}}\psi_1(t_1)\psi_2(t_2),
\end{equation}
where
\begin{equation}
\label{ziko0909}
K(t_1,t_2)={\bf 1}_{\{t_1<t_{2}\}}\psi_1(t_1)\psi_2(t_2),\ \ \ t_1, t_2\in[t, T].
\end{equation}

\vspace{4mm}

{\bf Lemma 2.3.}\
{\it Under the conditions of Theorem {\rm 2.2}
the following relation

\vspace{-2mm}
\begin{equation}
\label{30.36}
J[{K^{*}}]_{T,t}^{(2)}=
J^{*}[\psi^{(2)}]_{T,t}
\end{equation}

\vspace{2mm}
\noindent
is valid w.~p.~{\rm 1}, 
where $J[{K^{*}}]_{T,t}^{(2)}$ is defined by the equality 
{\rm (\ref{30.34})}.}

{\bf Proof.} Substituting
(\ref{1999.1}) into (\ref{30.34}) (the case $k=2$) and using 
Lemma 1.1 together with
(\ref{30.52}) (the case $k=2$)
it is easy to see that w.~p.~1

\vspace{-2mm}
$$
J[{K^{*}}]_{T,t}^{(2)}
=
J[\psi^{(2)}]_{T,t}+
\frac{1}{2}{\bf 1}_{\{i_1=i_2\}}
\int\limits_t^T\psi_1(t_1)\psi_2(t_1)dt_1
=
$$

\vspace{-2mm}
\begin{equation}
\label{30.37}
=J^{*}[\psi^{(2)}]_{T,t}.
\end{equation}

\vspace{3mm}

Let us consider the following generalized double Fourier sum

\vspace{-2mm}
$$
\sum_{j_1=0}^{p_1}\sum_{j_2=0}^{p_2}
C_{j_2j_1}\phi_{j_1}(t_1)\phi_{j_2}(t_2),
$$

\vspace{2mm}
\noindent
where $C_{j_2j_1}$ is the Fourier coefficient defined as follows

\vspace{-2mm}
\begin{equation}
\label{1}
C_{j_2j_1}=\int\limits_{[t,T]^2}
K^{*}(t_1,t_2)\phi_{j_1}(t_1)\phi_{j_2}(t_2)dt_1dt_2.
\end{equation}

Further, substitute the relation

\vspace{-2mm}
$$
K^{*}(t_1,t_2)=
\sum_{j_1=0}^{p_1}\sum_{j_2=0}^{p_2}
C_{j_2j_1}\phi_{j_1}(t_1)\phi_{j_2}(t_2)+
K^{*}(t_1,t_2)
-
$$

\vspace{-2mm}
$$
-
\sum_{j_1=0}^{p_1}\sum_{j_2=0}^{p_2}
C_{j_2j_1}\phi_{j_1}(t_1)\phi_{j_2}(t_2)
$$

\vspace{3mm}
\noindent
into $J[{K^{*}}]_{T,t}^{(2)}.$ At that we suppose that $p_1,p_2<\infty.$

Then using Lemma 1.3 (the case $k=2$), we obtain

\vspace{-3mm}
\begin{equation}
\label{proof1}
~~~~~~~~~ J^{*}[\psi^{(2)}]_{T,t}=
\sum_{j_1=0}^{p_1}\sum_{j_2=0}^{p_2}
C_{j_2j_1}
\zeta_{j_1}^{(i_1)}\zeta_{j_2}^{(i_2)}+
J[R_{p_1p_2}]_{T,t}^{(2)}\ \ \ \hbox{w. p. {\rm 1}},
\end{equation}

\vspace{2mm}
\noindent
where the stochastic integral
$J[R_{p_1p_2}]_{T,t}^{(2)}$
is defined in accordance with (\ref{30.34}) and

\vspace{-2mm}
\begin{equation}
\label{30.46}
~~~~~~~R_{p_1p_2}(t_1,t_2)=
K^{*}(t_1,t_2)-
\sum_{j_1=0}^{p_1}\sum_{j_2=0}^{p_2}
C_{j_2j_1}\phi_{j_1}(t_1)\phi_{j_2}(t_2),
\end{equation}

$$
\zeta_{j}^{(i)}=\int\limits_t^T \phi_{j}(s) d{\bf w}_s^{(i)},
$$

\vspace{-4mm}
$$
J[R_{p_1p_2}]_{T,t}^{(2)}=\int\limits_t^T\int\limits_t^{t_2}
R_{p_1p_2}(t_1,t_2)d{\bf w}_{t_1}^{(i_1)}d{\bf w}_{t_2}^{(i_2)}
+\int\limits_t^T\int\limits_t^{t_1}
R_{p_1p_2}(t_1,t_2)d{\bf w}_{t_2}^{(i_2)}d{\bf w}_{t_1}^{(i_1)}+
$$
$$
+{\bf 1}_{\{i_1=i_2\}}
\int\limits_t^T R_{p_1p_2}(t_1,t_1)dt_1.
$$

\vspace{4mm}

Using standard moment properties 
of stochastic integrals \cite{Gih1}
(see (\ref{99.010}), (\ref{99.010a})), we get

\vspace{-2mm}
$$
{\sf M}\left\{\left(J[R_{p_1p_2}]_{T,t}^{(2)}\right)^{2}
\right\}=
$$

\vspace{-2mm}
$$
={\sf M}\left\{\left(\int\limits_t^T\int\limits_t^{t_2}
R_{p_1p_2}(t_1,t_2)d{\bf w}_{t_1}^{(i_1)}d{\bf w}_{t_2}^{(i_2)}
+\int\limits_t^T\int\limits_t^{t_1}
R_{p_1p_2}(t_1,t_2)d{\bf w}_{t_2}^{(i_2)}d{\bf w}_{t_1}^{(i_1)}
\right)^2\right\}+
$$

\vspace{-2mm}
$$
+
{\bf 1}_{\{i_1=i_2\}}
\left(\int\limits_t^T R_{p_1p_2}(t_1,t_1)dt_1\right)^2\le
$$
$$
\le 
2\left(\int\limits_t^T\int\limits_t^{t_2}
\left(R_{p_1p_2}(t_1,t_2)\right)^{2}dt_1 dt_2
+
\int\limits_t^T\int\limits_t^{t_1}
\left(R_{p_1p_2}(t_1,t_2)\right)^{2}dt_2 dt_1\right)+
$$

\vspace{-2mm}
$$
+
{\bf 1}_{\{i_1=i_2\}}
\left(\int\limits_t^T R_{p_1p_2}(t_1,t_1)dt_1\right)^2=
$$

\vspace{-4mm}
\begin{equation}
\label{newbegin1}
~~~~~~~=
2\int\limits_{[t, T]^2}
\left(R_{p_1p_2}(t_1,t_2)\right)^{2}dt_1 dt_2
+
{\bf 1}_{\{i_1=i_2\}}
\left(\int\limits_t^T R_{p_1p_2}(t_1,t_1)dt_1\right)^2.
\end{equation}

\vspace{5mm}

We have
$$
\int\limits_{[t, T]^2}
\left(R_{p_1p_2}(t_1,t_2)\right)^{2}dt_1 dt_2=
$$
$$
=
\int\limits_{[t, T]^2}
\Biggl(
K^{*}(t_1,t_2)-
\sum_{j_1=0}^{p_1}\sum_{j_2=0}^{p_2}C_{j_2 j_1}
\phi_{j_1}(t_1)\phi_{j_2}(t_2)\Biggr)^2 dt_1 dt_2=
$$
\begin{equation}
\label{strange11}
~~~~~~~=\int\limits_{[t, T]^2}
\Biggl(
K(t_1,t_2)-
\sum_{j_1=0}^{p_1}\sum_{j_2=0}^{p_2}C_{j_2 j_1}
\phi_{j_1}(t_1)\phi_{j_2}(t_2)\Biggr)^2 dt_1 dt_2.
\end{equation}

\vspace{2mm}

The function $K(t_1,t_2)$ is piecewise continuous in the 
square $[t, T]^2$.
At this situation it is well known that the generalized
multiple Fourier series 
of the function $K(t_1,t_2)\in L_2([t, T]^2)$ is converging 
to this function in the square $[t, T]^2$ in the mean-square sense, i.e.

\vspace{-2mm}
$$
\hbox{\vtop{\offinterlineskip\halign{
\hfil#\hfil\cr
{\rm lim}\cr
$\stackrel{}{{}_{p_1,p_2\to \infty}}$\cr
}} }\Biggl\Vert
K(t_1,t_2)-
\sum_{j_1=0}^{p_1}\sum_{j_2=0}^{p_2}
C_{j_2 j_1}\prod_{l=1}^{2} \phi_{j_l}(t_l)\Biggr\Vert_{L_2([t,T]^2)}=0,
$$

\vspace{2mm}
\noindent
where notations are the same as in (\ref{sos1z}).

So, we obtain
\begin{equation}
\label{newbegin2}
\hbox{\vtop{\offinterlineskip\halign{
\hfil#\hfil\cr
{\rm lim}\cr
$\stackrel{}{{}_{p_1,p_2\to \infty}}$\cr
}} }
\int\limits_{[t, T]^2}
\left(R_{p_1p_2}(t_1,t_2)\right)^{2}dt_1 dt_2=0.
\end{equation}

\vspace{3mm}

Note that
$$
\int\limits_t^T R_{p_1p_2}(t_1,t_1)dt_1=
$$

\vspace{-1mm}
$$
=
\int\limits_t^T
\left(
\frac{1}{2}\psi_1(t_1)\psi_2(t_1) - 
\sum_{j_1=0}^{p_1}\sum_{j_2=0}^{p_2}C_{j_2 j_1}
\phi_{j_1}(t_1)\phi_{j_2}(t_1)\right) dt_1=
$$

\vspace{-1mm}
$$
=
\frac{1}{2}\int\limits_t^T
\psi_1(t_1)\psi_2(t_1)dt_1 -
\sum_{j_1=0}^{p_1}\sum_{j_2=0}^{p_2}C_{j_2 j_1}
\int\limits_t^T\phi_{j_1}(t_1)\phi_{j_2}(t_1)dt_1=
$$

\vspace{-1mm}
$$
=
\frac{1}{2}\int\limits_t^T
\psi_1(t_1)\psi_2(t_1)dt_1 -
\sum_{j_1=0}^{p_1}\sum_{j_2=0}^{p_2}C_{j_2 j_1}
{\bf 1}_{\{j_1=j_2\}}=
$$

\vspace{-1mm}
\begin{equation}
\label{newbegin3}
=
\frac{1}{2}\int\limits_t^T
\psi_1(t_1)\psi_2(t_1)dt_1 -
\sum_{j_1=0}^{{\rm min}\{p_1,p_2\}}C_{j_1 j_1}.
\end{equation}

\vspace{2mm}

From (\ref{newbegin3}) and Lemma 2.2 we get

\vspace{-3mm}
$$
\lim\limits_{p_1\to\infty}
\lim\limits_{p_2\to\infty}\int\limits_t^T R_{p_1p_2}(t_1,t_1)dt_1=
$$

\vspace{-1mm}
$$
=\frac{1}{2}\int\limits_t^T
\psi_1(t_1)\psi_2(t_1)dt_1 -
\lim\limits_{p_1\to\infty}
\sum_{j_1=0}^{p_1}C_{j_1 j_1}
=
$$

\vspace{-1mm}
$$
=\frac{1}{2}\int\limits_t^T
\psi_1(t_1)\psi_2(t_1)dt_1 -
\sum_{j_1=0}^{\infty}C_{j_1 j_1}
=
$$

\vspace{-1mm}
\begin{equation}
\label{s1}
=\lim\limits_{p_1,p_2\to\infty}\int\limits_t^T R_{p_1p_2}(t_1,t_1)dt_1.
\end{equation}

\vspace{3mm}

If we prove the following relation
\begin{equation}
\label{s11}
\lim\limits_{p_1\to\infty}
\lim\limits_{p_2\to\infty}\int\limits_t^T R_{p_1p_2}(t_1,t_1)dt_1=0,
\end{equation}

\noindent
then from (\ref{s1}) we obtain

\vspace{-2mm}
\begin{equation}
\label{44}
\frac{1}{2}\int\limits_t^T
\psi_1(t_1)\psi_2(t_1)dt_1 =
\sum_{j_1=0}^{\infty}C_{j_1 j_1},
\end{equation}
\begin{equation}
\label{444xx}
\lim\limits_{p_1,p_2\to\infty}\int\limits_t^T R_{p_1p_2}(t_1,t_1)dt_1=0.
\end{equation}

\vspace{3mm}

From (\ref{newbegin1}), (\ref{newbegin2}), and (\ref{444xx}) we get
$$
\lim\limits_{p_1,p_2\to\infty}
{\sf M}\left\{\left(J[R_{p_1p_2}]_{T,t}^{(2)}\right)^{2}
\right\}=0
$$

\vspace{2mm}
\noindent
and Theorem 2.1 will be proved (see (\ref{proof1})).

The proof of the equality
(\ref{s11}) can be carried out in the same way as in the proof of Theorem~2.1 or,
under weaker conditions, as in the proof of Theorem~2.2.

\subsection{Approach Based on Arbitrary Complete Orthonormal System of 
Functions in the Space $L_2([t,T])$ and $\psi_1(\tau), \psi_2(\tau)\in L_2([t,T])$}

Let us prove the equality (\ref{5t}) under weaker restrictions.
Suppose that 
$\{\phi_j(x)\}_{j=0}^{\infty}$ is an arbitrary complete orthonormal system of 
functions in the space $L_2([t, T])$
and $\psi_1(\tau)\equiv \psi_2(\tau)$ 
or 
\begin{equation}
\label{trace2}
\psi_1(\tau)=\psi_2(\tau)\int\limits_t^{\tau} g(\theta)d\theta,
\end{equation}
where $\tau\in [t, T]$ and $\psi_1(\tau), \psi_2(\tau)\in L_2([t, T]),$ $g(\tau)\in L_1([t, T])$.

Thus, we will prove the equality
\begin{equation}
\label{trace1}
~~~~~~~~~\sum_{j=0}^{\infty}\int\limits_t^T
\psi_2(t_2)\phi_j(t_2)
\int\limits_t^{t_2}
\psi_1(t_1)\phi_j(t_1)dt_1 dt_2=
\frac{1}{2}\int\limits_t^T 
\psi_1(\tau) \psi_2(\tau) d\tau
\end{equation}
under the above conditions.

Using Fubini's Theorem, Lebesgue's Dominated Convergence Theorem
and Parseval's equality, we have (see (\ref{trace2}))
$$
\sum_{j=0}^{\infty}\int\limits_t^T
\psi_2(t_2)\phi_j(t_2)
\int\limits_t^{t_2}
\psi_1(t_1)\phi_j(t_1)dt_1 dt_2=
$$
\begin{equation}
\label{july10000}
~~~~~~~~~~=\sum_{j=0}^{\infty}\int\limits_t^T
\psi_2(t_2)\phi_j(t_2)
\int\limits_t^{t_2}
\psi_2(t_1)\phi_j(t_1)
\int\limits_t^{t_1} g(\tau)d\tau dt_1 dt_2=
\end{equation}
$$
=\sum_{j=0}^{\infty}\int\limits_t^T
g(\tau)\int\limits_{\tau}^{T} \psi_2(t_1)\phi_j(t_1)
\int\limits_{t_1}^{T}
\psi_2(t_2)\phi_j(t_2)
dt_2 dt_1 d\tau =
$$
\begin{equation}
\label{after001}
=\frac{1}{2}\sum_{j=0}^{\infty}\int\limits_t^T
g(\tau)\left(\int\limits_{\tau}^{T} \psi_2(t_1)\phi_j(t_1)
dt_1\right)^2 d\tau =
\end{equation}
\begin{equation}
\label{after002}
=\frac{1}{2}\int\limits_t^T 
g(\tau)\sum_{j=0}^{\infty} \left(\int\limits_{t}^{T} {\bf 1}_{\{\tau<t_1\}}\psi_2(t_1)\phi_j(t_1)
dt_1\right)^2 d\tau =
\end{equation}
$$
=
\frac{1}{2}\int\limits_t^T 
g(\tau)\int\limits_{t}^{T} {\bf 1}_{\{\tau<t_1\}}\psi_2^2(t_1)
dt_1 d\tau 
=
\frac{1}{2}\int\limits_t^T 
g(\tau)\int\limits_{\tau}^{T} \psi_2^2(t_1)
dt_1 d\tau =
$$
\begin{equation}
\label{after003}
=
\frac{1}{2}\int\limits_t^T 
\psi_2^2(t_1) \int\limits_{t}^{t_1} g(\tau)
d\tau dt_1 =
\end{equation}
\begin{equation}
\label{trace10}
=
\frac{1}{2}\int\limits_t^T 
\psi_1(t_1) \psi_2(t_1) dt_1,
\end{equation}

\noindent 
where the transition from (\ref{after001}) to (\ref{after002})
is based on Lebesgue's 
Dominated Convergence Theorem. The integrable majorant exists due to 
Parseval's equality 
$$
\left\vert g(\tau)\right\vert
\sum_{j=0}^{q} \left(\int\limits_{\tau}^{T}\psi_2(t_1)\phi_j(t_1)
dt_1\right)^2\le
\left\vert g(\tau)\right\vert
\sum_{j=0}^{\infty} \left(\int\limits_{t}^{T} {\bf 1}_{\{\tau<t_1\}}\psi_2(t_1)\phi_j(t_1)
dt_1\right)^2=
$$
$$
=\left\vert g(\tau)\right\vert
\int\limits_{t}^{T}\left({\bf 1}_{\{\tau<t_1\}}\right)^2 \psi_2^2(t_1)
dt_1\le
\left\vert g(\tau)\right\vert \left\Vert \psi_2 \right\Vert_{L_2([t,T])}^2
=C\left\vert g(\tau)\right\vert
$$

\noindent
almost everywhere on $[t,T]$ with respect to Lebesgue's measure,
where constant $C$ does not depend on $p.$

From the other hand, using Fubini's Theorem 
and the generalized Parseval equality as well as the transition from (\ref{july10000})
to (\ref{after003}), we get
$$
\sum_{j=0}^{\infty}\int\limits_t^T
\psi_1(t_2)\phi_j(t_2)
\int\limits_t^{t_2}
\psi_2(t_1)\phi_j(t_1)dt_1 dt_2=
$$
$$
=\sum_{j=0}^{\infty}\int\limits_t^T
\psi_2(t_2)\phi_j(t_2)
\int\limits_t^{t_2}
g(\tau)d\tau
\int\limits_t^{t_2}  
\psi_2(t_1)\phi_j(t_1)
dt_1 dt_2=
$$
$$
=\sum_{j=0}^{\infty}\int\limits_t^T
\psi_2(t_1)\phi_j(t_1)
\int\limits_{t_1}^{T}
\psi_2(t_2)\phi_j(t_2)
\int\limits_t^{t_2}  
g(\tau)d\tau
dt_2 dt_1=
$$
$$
=\sum_{j=0}^{\infty}\int\limits_t^T
\psi_2(t_1)\phi_j(t_1) dt_1
\int\limits_{t}^{T}
\psi_2(t_2)\phi_j(t_2)
\int\limits_t^{t_2}  
g(\tau)d\tau
dt_2 -
$$
$$
-\sum_{j=0}^{\infty}\int\limits_t^T
\psi_2(t_1)\phi_j(t_1)
\int\limits_{t}^{t_1}
\psi_2(t_2)\phi_j(t_2)
\int\limits_t^{t_2}  
g(\tau)d\tau
dt_2 dt_1=
$$
$$
=
\int\limits_t^T 
\psi_2(t_1) \cdot \psi_2(t_1) \int\limits_{t}^{t_1} g(\tau)
d\tau dt_1 - \frac{1}{2}\int\limits_t^T 
\psi_2^2(t_1) \int\limits_{t}^{t_1} g(\tau)
d\tau dt_1 =
$$
\begin{equation}
\label{trace11}
=
\frac{1}{2}\int\limits_t^T 
\psi_2^2(t_1) \int\limits_{t}^{t_1} g(\tau)
d\tau dt_1 =
\frac{1}{2}\int\limits_t^T 
\psi_1(t_1) \psi_2(t_1) dt_1.
\end{equation}

\vspace{2mm}

In addition, for the case $\psi_1(\tau)\equiv \psi_2(\tau)$, 
using the Parseval equality, we ob\-tain
$$
\sum_{j=0}^{\infty}\int\limits_t^T
\psi_1(t_2)\phi_j(t_2)
\int\limits_t^{t_2}
\psi_1(t_1)\phi_j(t_1) dt_1 dt_2=
$$
\begin{equation}
\label{trace12}
=
\frac{1}{2}\sum_{j=0}^{\infty}
\left(\int\limits_{t}^{T} \psi_1(t_1)\phi_j(t_1)
dt_1\right)^2 =
\frac{1}{2}
\int\limits_{t}^{T} \psi_1^2(t_1)dt_1
=\frac{1}{2}
\int\limits_{t}^{T} \psi_1(t_1)\psi_1(t_2)dt_1.
\end{equation}

\vspace{2mm}

The equality (\ref{trace1}) is proved for $\psi_1(\tau)\equiv \psi_2(\tau)$ 
or when the equality (\ref{trace2}) is satisfied.

Further, let us suppose that
$\psi_2(\tau)=(\tau-t)^l,$ $g(\tau)=k (\tau-t)^{k-1},$
where $l=0,1,2,\ldots$ and $k=1,2,\ldots $
Note that this case is important for applications (see Sect.~4.11 and 4.12).

From (\ref{trace2}) we obtain 
$$
\psi_1(\tau)=\psi_2(\tau)\int\limits_t^{\tau} g(\theta)d\theta=
k(\tau-t)^l \int\limits_t^{\tau} (\theta-t)^{k-1}d\theta=(\tau-t)^{l+k}.
$$

Taking into account (\ref{trace10})--(\ref{trace12}),
we get
$$
\sum_{j=0}^{\infty}\int\limits_t^T 
(t_2-t)^l \phi_j(t_2)
\int\limits_t^{t_2} 
(t_1-t)^{l+k} \phi_j(t_1)dt_1 dt_2=
$$
$$
=
\sum_{j=0}^{\infty}\int\limits_t^T 
(t_2-t)^{l+k} \phi_j(t_2)
\int\limits_t^{t_2} 
(t_1-t)^{l} \phi_j(t_1)dt_1 dt_2=
$$
\begin{equation}
\label{dsds1}
=\frac{1}{2}\int\limits_t^T (\tau-t)^{2l+k} d\tau,
\end{equation}

\noindent
where $k, l=0,1,2,\ldots $

Let us rewrite the equality (\ref{dsds1}) in the following form
\begin{equation}
\label{strange902}
\sum_{j=0}^{\infty}\int\limits_t^T 
(t_2-t)^l \phi_j(t_2)
\int\limits_t^{t_2} 
(t_1-t)^{m} \phi_j(t_1)dt_1 dt_2=
\frac{1}{2}\int\limits_t^T (\tau-t)^{l}(\tau-t)^{m} d\tau,
\end{equation}

\noindent
where $l, m=0,1,2,\ldots $

The equality similar to (\ref{strange902})
was obtained in \cite{rybakov5000}, \cite{rybakov7000x} using other arguments.
These arguments are based on trace class operators and the equality
of matrix and integral traces for such operators (see Sect.~2.27 for details).

In addition, the formula similar to (\ref{strange902}) was used 
in \cite{rybakov5000}, \cite{rybakov7000x} to generalize the equality (\ref{trace1}) 
to the case of an arbitrary 
complete ortho\-nor\-mal system of functions in the space $L_2([t, T])$
and $\psi_1(\tau),\psi_2(\tau)$ $\in $ $L_2([t, T]).$
This means that Theorems~2.1, 2.2 can be generalized to the case
of continuous functions $\psi_1(\tau),\psi_2(\tau)$ 
(this condition is related to the definition (\ref{123321.2})
of the Stratonovich stochastic integral (see Sect.~2.1.1 for details))
and an arbitrary complete orthonormal system of 
functions in the space $L_2([t,T])$.

Consider the mentioned approach \cite{rybakov5000}, \cite{rybakov7000x} in our
interpretation
(after this, we will consider an approach that is slightly
different from the approach in \cite{rybakov5000}, \cite{rybakov7000x}).
Since the equality (\ref{strange902}) is valid for monomials 
with respect to $\tau-t$ ($\tau\in [t, T]$), it will obviously
also be valid for Legendre polynomials that form a complete 
orthonormal system of functions in the space $L_2([t, T])$
and finite linear combinations of Legendre polynomials.

Let $\psi_1(\tau), \psi_2(\tau)\in L_2([t, T])$ and 
$\psi_1^{(p)}(\tau), \psi_2^{(q)}(\tau)$
be approximations of the functions $\psi_1(\tau),\psi_2(\tau)$,
respectively, which are partial sums of the corresponding
Fourier--Legendre series. Then we have (see (\ref{strange902}))
\begin{equation}
\label{strange903}
~~~~~\sum_{j=0}^{\infty}\int\limits_t^T 
\psi_2^{(q)}(t_2) \phi_j(t_2)
\int\limits_t^{t_2} 
\psi_1^{(p)}(t_1) \phi_j(t_1)dt_1 dt_2=
\frac{1}{2}\int\limits_t^T \psi_1^{(p)}(\tau) \psi_2^{(q)}(\tau) d\tau,
\end{equation}

\noindent
where $p, q\in {\bf N},$ the series converges absolutly and 
its sum does not depend on a basis system $\left\{\phi_j(x)\right\}_{j=0}^{\infty}$
(we mean permutation of the terms
of the series on the left-hand side of (\ref{strange903})
(any permutation of basis functions $\phi_j(x)$ forms a basis 
in $L_2([t,T])$ \cite{gohb})).

Using Fubini's Theorem, we rewrite (\ref{strange903}) in the form
$$
\sum_{j=0}^{\infty}\left(\int\limits_t^T 
\psi_2^{(q)}(t_2) \phi_j(t_2)
\int\limits_t^{t_2} 
\psi_1^{(p)}(t_1) \phi_j(t_1)dt_1 dt_2+\right.
$$
\begin{equation}
\label{strange903xxx}
~~~~~~~~~\left.+
\int\limits_t^T 
\psi_1^{(p)}(t_2) \phi_j(t_2)
\int\limits_{t_2}^{T} 
\psi_2^{(q)}(t_1) \phi_j(t_1)dt_1 dt_2\right)=
\int\limits_t^T \psi_1^{(p)}(\tau) \psi_2^{(q)}(\tau) d\tau.
\end{equation}

Let us fix $q$ in (\ref{strange903xxx}). 
The right-hand side of (\ref{strange903xxx}) for a fixed $q$ defines
(as a scalar product in $L_2([t, T])$) a linear bounded
(and therefore continuous) functional in $L_2([t, T]),$
which is given by the function $\psi_2^{(q)}$.
The integral operator (which corresponds to the matrix trace 
on the left-hand side of (\ref{strange903xxx})) is a trace class operator
(see \cite{rybakov7000x}).
The matrix trace of the mentioned operator (on the left-hand side of (\ref{strange903xxx}))
is also a linear bounded (and therefore continuous) 
functional (in the space of trace class operators \cite{gohb},
\cite{goldberg})
which can be extended to the space $L_2([t, T])$ by continuity \cite{Pugach}.

Let us implement the passage to the limit $\lim\limits_{p\to\infty}$
in (\ref{strange903xxx})
$$
\sum_{j=0}^{\infty}\left(\int\limits_t^T 
\psi_2^{(q)}(t_2) \phi_j(t_2)
\int\limits_t^{t_2} 
\psi_1(t_1) \phi_j(t_1)dt_1 dt_2+\right.
$$
\begin{equation}
\label{strange904}
~~~~~~~~~\left.+
\int\limits_t^T 
\psi_1(t_2) \phi_j(t_2)
\int\limits_{t_2}^{T} 
\psi_2^{(q)}(t_1) \phi_j(t_1)dt_1 dt_2\right)=
\int\limits_t^T \psi_1(\tau) \psi_2^{(q)}(\tau) d\tau,
\end{equation}

\noindent
where $q\in {\bf N}.$ 
Recall that 
$\psi_2^{(q)}(\tau)$ is a partial sum of the Fourier--Legendre series
of any function $\psi_2(\tau)\in L_2([t,T]),$ i.e. the equality (\ref{strange904}) holds
on a dense subset in $L_2([t,T]).$
The right-hand side of (\ref{strange904}) defines
(as a scalar product in $L_2([t, T])$) a linear bounded (and therefore continuous)
functional in $L_2([t, T]),$
which is given by the function $\psi_1$.
On the left-hand side of (\ref{strange904}) (by virtue of the equality (\ref{strange904}))
there is a linear continuous functional on a dense subset in 
$L_2([t,T]).$ This functional can be uniquely extended 
to a linear continuous functional in $L_2([t, T])$
(see \cite{reed}, Theorem~I.7, P.~9).

Let us implement the passage to the limit $\lim\limits_{q\to\infty}$
in (\ref{strange904})
$$
\sum_{j=0}^{\infty}\left(\int\limits_t^T 
\psi_2(t_2) \phi_j(t_2)
\int\limits_t^{t_2} 
\psi_1(t_1) \phi_j(t_1)dt_1 dt_2+\right.
$$
\begin{equation}
\label{strange904xxx}
~~~~~~~~~\left.+
\int\limits_t^T 
\psi_1(t_2) \phi_j(t_2)
\int\limits_{t_2}^{T} 
\psi_2(t_1) \phi_j(t_1)dt_1 dt_2\right)=
\int\limits_t^T \psi_1(\tau) \psi_2(\tau) d\tau.
\end{equation}

\vspace{2mm}

Applying Fubini's Theorem to the left-hand side of (\ref{strange904xxx}), we obtain 
\begin{equation}
\label{start1000}
~~~~~~~~~~\sum_{j=0}^{\infty}\int\limits_t^T 
\psi_2(t_2) \phi_j(t_2)
\int\limits_t^{t_2} 
\psi_1(t_1) \phi_j(t_1)dt_1 dt_2=
\frac{1}{2}\int\limits_t^T \psi_1(\tau) \psi_2(\tau) d\tau,
\end{equation}

\noindent
where $\left\{\phi_j(x)\right\}_{j=0}^{\infty}$
is an arbitrary complete orthonormal system of 
functions in the space $L_2([t,T])$ and
$\psi_1(\tau),\psi_2(\tau)\in $ $L_2([t, T]).$

However, the equality (\ref{start1000})
can be obtained somewhat more simply.
Now let us consider an approach that is
slightly different form the approach in 
\cite{rybakov5000}, \cite{rybakov7000x}.

Consider the equality (\ref{strange903})
and fix $q$ in it.
The right-hand side of (\ref{strange903}) for a fixed $q$ defines
(as a scalar product of $\psi_1^{(p)}$ and $\frac{1}{2}\psi_2^{(q)}$ in $L_2([t, T])$) a linear bounded
(and therefore continuous) functional in $L_2([t, T]),$
which is given by the function $\frac{1}{2}\psi_2^{(q)}$.

On the left-hand side of (\ref{strange903}) (by virtue of the equality (\ref{strange903}))
there is a linear continuous functional on a dense subset in 
$L_2([t,T])$ (recall that 
$\psi_1^{(p)}(\tau)$ is a partial sum of the Fourier--Legendre series
of any function $\psi_1(\tau)\in L_2([t,T])$).
This functional can be uniquely extended 
to a linear continuous functional in $L_2([t, T])$
(see \cite{reed}, Theorem~I.7, P.~9).

Let us implement the passage to the limit $\lim\limits_{p\to\infty}$
in the equality (\ref{strange903}) 
\begin{equation}
\label{strange903october2024}
~~~~~\sum_{j=0}^{\infty}\int\limits_t^T 
\psi_2^{(q)}(t_2) \phi_j(t_2)
\int\limits_t^{t_2} 
\psi_1(t_1) \phi_j(t_1)dt_1 dt_2=
\frac{1}{2}\int\limits_t^T \psi_1(\tau) \psi_2^{(q)}(\tau) d\tau,
\end{equation}

\noindent
where $q\in {\bf N}.$

Recall that 
$\psi_2^{(q)}(\tau)$ is a partial sum of the Fourier--Legendre series
of any function $\psi_2(\tau)\in L_2([t,T]),$ i.e. the equality (\ref{strange903october2024}) holds
on a dense subset in $L_2([t,T]).$
The right-hand side of (\ref{strange903october2024}) defines
(as a scalar product of $\psi_2^{(q)}$ and $\frac{1}{2}\psi_1$ in $L_2([t, T])$) 
a linear bounded (and therefore continuous)
functional in $L_2([t, T]),$
which is given by the function $\frac{1}{2}\psi_1$.
On the left-hand side of (\ref{strange903october2024}) (by virtue of the equality 
(\ref{strange903october2024}))
there is a linear continuous functional on a dense subset in 
$L_2([t,T]).$ This functional can be uniquely extended 
to a linear continuous functional in $L_2([t, T])$
(see \cite{reed}, Theorem~I.7, P.~9).

Let us implement the passage to the limit $\lim\limits_{q\to\infty}$
in (\ref{strange903october2024})
$$
\sum_{j=0}^{\infty}\int\limits_t^T 
\psi_2(t_2) \phi_j(t_2)
\int\limits_t^{t_2} 
\psi_1(t_1) \phi_j(t_1)dt_1 dt_2=
\frac{1}{2}\int\limits_t^T \psi_1(\tau) \psi_2(\tau) d\tau,
$$

\noindent
where $\left\{\phi_j(x)\right\}_{j=0}^{\infty}$
is an arbitrary complete orthonormal system of 
functions in the space $L_2([t,T])$ and
$\psi_1(\tau),\psi_2(\tau)\in $ $L_2([t, T]).$
As a result, we obtained the equality (\ref{start1000}).

Thus, we have the following theorem.

{\bf Theorem 2.3.}\ {\it Suppose that 
$\{\phi_j(x)\}_{j=0}^{\infty}$ is an arbitrary complete orthonormal system of 
functions in the space $L_2([t, T]).$
Moreover$,$ $\psi_1(\tau), \psi_2(\tau)$ are continuous 
functions on $[t, T].$ 
Then$,$ 
for the iterated Stra\-to\-novich stochastic integral
$$
J^{*}[\psi^{(2)}]_{T,t}={\int\limits_t^{*}}^T\psi_2(t_2)
{\int\limits_t^{*}}^{t_2}\psi_1(t_1)d{\bf w}_{t_1}^{(i_1)}
d{\bf w}_{t_2}^{(i_2)}\ \ \ (i_1, i_2=1,\ldots,m)
$$
the following expansion  
\begin{equation}
\label{trace20}
J^{*}[\psi^{(2)}]_{T,t}=\hbox{\vtop{\offinterlineskip\halign{
\hfil#\hfil\cr
{\rm l.i.m.}\cr
$\stackrel{}{{}_{p_1,p_2\to \infty}}$\cr
}} }\sum_{j_1=0}^{p_1}\sum_{j_2=0}^{p_2}
C_{j_2j_1}\zeta_{j_1}^{(i_1)}\zeta_{j_2}^{(i_2)}
\end{equation}

\noindent
that converges in the mean-square
sence is valid$,$ where the notations are the same as in Theorems {\rm 2.1, 2.2.}
}

The condition of continuity of the functions
$\psi_1(\tau), \psi_2(\tau)$ 
is related to the definition (\ref{123321.2}) 
of the Stratonovich stochastic integral that we use.

Theorem~2.3 can be generalized to the case $\psi_1(\tau), \psi_2(\tau)\in L_2([t, T])$
if instead of the definition (\ref{123321.2})
we use another definition of the Stratonovich stochastic integral (see (\ref{dsds5})
and Theorem~2.44 in Sect.~2.18 for details).

\subsection{Approach Based on Arbitrary Complete Orthonormal System of 
Functions in the Space $L_2([t,T])$ and Symmetrized Kernel $K'(t_1,t_2)$}

Let us list some useful facts that we will need further in this section.

{\bf Theorem~A} (\cite{goldberg}, Theorem~8.1).\ {\it 
Let $\mathbb{K}: L_2([t,T])\rightarrow L_2([t,T])$
be an integral operator defined by
$$
\left(\mathbb{K}f\right)(\tau)=\int\limits_t^T K(\tau,s)f(s)ds,
$$
where $K(\tau,s)$ is a continuous function on $[t, T]\times[t, T]$.
If$,$ in addition$,$ $\mathbb{K}$ is a trace class operator then
\begin{equation}
\label{july11004}
tr\mathbb{K}=\int\limits_t^T K(s,s)ds,
\end{equation}

\noindent
where trace $tr\mathbb{K}$ is defined as a series
of singular values $s_j(\mathbb{K})$ of $\mathbb{K}$.}

\vspace{1mm}

{\bf Theorem~B} (\cite{goldberg}, P.~71).\ {\it 
Let 
$$
\left(\mathbb{K}f\right)(\tau)=\int\limits_t^T K(\tau,s)f(s)ds,
$$

\noindent
the kernel $K(\tau,s)$ is continuous on $[t, T]\times[t, T]$
and satisfies the condition
\begin{equation}
\label{july11002}
~~~~\left\vert K(\tau,s_2)-K(\tau,s_1)\right\vert \le C \left\vert s_2-s_1\right\vert^{\alpha},
\end{equation}

\noindent
where $0<\alpha\le 1$. If$,$ in addition$,$ $\mathbb{K}$ is a Hermitian
operator and $\alpha>1/2,$ then
$$
\sum\limits_{j=0}^{\infty} s_j(\mathbb{K})<\infty
$$

\noindent
i.e.$,$ $\mathbb{K}$ is a trace class operator.}

\vspace{1mm}

Suppose that $\mathbb{A}: H \rightarrow H$ is a linear bounded operator. 
Recall \cite{gohb} that $\mathbb{A}$ has a finite matrix trace
if for any orthonormal basis $\left\{\phi_j(x)\right\}_{j=0}^{\infty}$
of the space $H$ the series
\begin{equation}
\label{july11205}
\sum_{j=0}^{\infty} \left\langle
\mathbb{A}\phi_j, \phi_j\right\rangle_H
\end{equation}

\noindent
converges, where $\left\langle
\cdot , \cdot \right\rangle_H$ is a scalar probuct in $H$.

Note that the series (\ref{july11205}) converges absolutely
since its sum does not depend on the permutation of the terms
of the series (\ref{july11205})
(any permutation of basis functions $\phi_j(x)$ forms a basis 
in $H)$ \cite{gohb}.

\vspace{1mm}

{\bf Theorem~C} (\cite{goldberg}, Theorem~5.6).\ {\it 
Let $\mathbb{K}: H\rightarrow H$ be a trace class operator.
Then
\begin{equation}
\label{july11201}
tr\mathbb{A}=\sum_{j=0}^{\infty} \left\langle
\mathbb{A}\phi_j, \phi_j\right\rangle_H
\end{equation}

\noindent
for any orthonormal basis $\left\{\phi_j(x)\right\}_{j=0}^{\infty}$ of $H$.}

\vspace{1mm}

Consider an integral operator $\mathbb{K'}: L_2([t,T])\rightarrow L_2([t,T])$
defined by the equality
$$
\left(\mathbb{K'}f\right)(\tau)=\int\limits_t^T K'(\tau,s)f(s)ds,
$$

\vspace{1mm}
\noindent
where the continuous kernel $K'(\tau,s)$ has the form {\rm (\ref{ziko5001})}$,$ i.e.
$$
K'(t_1,t_2)=\left\{
\begin{matrix}
\psi_2(t_1)\psi_1(t_2),\ \ t_1\ge t_2\cr\cr
\psi_1(t_1)\psi_2(t_2),\ \ t_1\le t_2
\end{matrix}
\right.\ \ \ (t_1,t_2\in[t,T])
$$

\noindent
and $\psi_1(\tau), \psi_2(\tau)$ are continuously differentiable functions on $[t, T].$

Recall that (see (\ref{9090}))

\vspace{-2.5mm}
\begin{equation}
\label{july11000}
\left|K'(t_2,s_2)-K'(t_1,s_1)\right|\le L
\left(|t_2-t_1|+|s_2-s_1|\right),
\end{equation}

\vspace{2mm}
\noindent
where $L<\infty$ and $(t_1,s_1)$, $(t_2,s_2)\in [t, T]^2.$

Let us substitute $t_1=t_2=\tau$ into (\ref{july11000})

\vspace{-2mm}
\begin{equation}
\label{july11001}
\left|K'(\tau,s_2)-K'(\tau,s_1)\right|\le L
|s_2-s_1|.
\end{equation}

\vspace{3mm}

Thus, the condition (\ref{july11002}) is fulfilled $(\alpha=1$).
Further, using Fubini's Theorem, we have
$$
\left\langle
\mathbb{K'}x, y\right\rangle_{L_2([t,T])}=
\int\limits_t^T\psi_2(t_2)y(t_2)\int\limits_t^{t_2}\psi_1(t_1)x(t_1)dt_1 dt_2+
$$
$$
+
\int\limits_t^T\psi_1(t_2)y(t_2)\int\limits_{t_2}^T\psi_2(t_1)x(t_1)dt_1 dt_2=
\int\limits_t^T\psi_1(t_1)x(t_1) \int\limits_{t_1}^T \psi_2(t_2)y(t_2) dt_2 dt_1+
$$
\begin{equation}
\label{july11207}
~~~~~~~~~~+
\int\limits_t^T\psi_2(t_1)x(t_1) \int\limits_{t}^{t_2} \psi_1(t_2)y(t_2) dt_2 dt_1=
\left\langle\mathbb{K'}y, x\right\rangle_{L_2([t,T])}.
\end{equation}

The conditions of Theorem~B are fulfilled. Then, $\mathbb{K'}$ is a trace class operator.
Since the kernel $K'(t_1,t_2)$ is continuous, then by Theorems~A and C
(see (\ref{july11004}) and (\ref{july11201})) we obtain
\begin{equation}
\label{july11208}
~~~~~~~~~~\sum_{j_1=0}^{\infty} \left\langle
\mathbb{K'}\phi_{j_1}, \phi_{j_1}\right\rangle_{L_2([t,T])}
=\int\limits_t^T K'(s,s)ds=\int\limits_t^T \psi_1(s)\psi_2(s)ds.
\end{equation}

Combining (\ref{july11207}), (\ref{july11208}) and applying Fubini's Theorem, we get
$$
\sum_{j_1=0}^{\infty}
\left(\int\limits_t^T\psi_2(t_2)\phi_{j_1}(t_2)\int\limits_t^{t_2}\psi_1(t_1)\phi_{j_1}(t_1)dt_1 dt_2+
\right.
$$
$$
\left.+
\int\limits_t^T\psi_1(t_2)\phi_{j_1}(t_2)\int\limits_{t_2}^T\psi_2(t_1)\phi_{j_1}(t_1)dt_1 dt_2\right)=
$$
$$
=\sum_{j_1=0}^{\infty}
\left(\int\limits_t^T\psi_2(t_2)\phi_{j_1}(t_2)\int\limits_t^{t_2}\psi_1(t_1)\phi_{j_1}(t_1)dt_1 dt_2+
\right.
$$
$$
\left.+
\int\limits_t^T \psi_2(t_1)\phi_{j_1}(t_1)\int\limits_t^{t_2} \psi_1(t_2)\phi_{j_1}(t_2) dt_2 dt_1\right)=
$$
\begin{equation}
\label{july11209}
~~~~~~~~=2\sum_{j_1=0}^{\infty}
\int\limits_t^T\psi_2(t_2)\phi_{j_1}(t_2)\int\limits_t^{t_2}\psi_1(t_1)\phi_{j_1}(t_1)dt_1 dt_2=
\int\limits_t^T \psi_1(s)\psi_2(s)ds.
\end{equation}

\vspace{2mm}

From (\ref{july11209}) we obtain
\begin{equation}
\label{july11220}
~~~~~~~~~~\sum_{j_1=0}^{\infty}
\int\limits_t^T\psi_2(t_2)\phi_{j_1}(t_2)\int\limits_t^{t_2}\psi_1(t_1)\phi_{j_1}(t_1)dt_1 dt_2=
\frac{1}{2}\int\limits_t^T \psi_1(s)\psi_2(s)ds,
\end{equation}

\noindent
where $\left\{\phi_j(x)\right\}_{j=0}^{\infty}$
is an arbitrary complete orthonormal system of 
functions in the space $L_2([t,T])$ and
$\psi_1(\tau),\psi_2(\tau)$ are continuously
differentiable functions on $[t, T].$

To further generalize of the equality (\ref{july11220})
to the case when 
$\psi_1(\tau),\psi_2(\tau)\in L_2([t,T])$
it is necessary to set $\psi_2(\tau)=(\tau-t)^l,$ $\psi_1(\tau)=(\tau-t)^m$
$(l,m=0,1,2,\ldots)$ and apply the reasoning
of the previous section after the formula (\ref{strange902}).

\section{Expansion of Iterated Stratonovich Stochastic Integrals of 
Multiplicity 3 Based on Theorem 1.1}

This section is devoted to the development of the method of expansion
and mean-square approximation of iterated It\^{o} stochastic integrals
based on generalized multiple Fourier series converging in the mean
(Theorem 1.1).
We adapt this method for the iterated Stratonovich stochastic
integrals of multiplicity 3.
The main results of this section have been derived 
with using triple Fourier--Legendre series as well as
triple trigonometric Fourier series for different
cases of series summation and different cases of weight
functions of iterated Stratonovich stochastic
integrals.

\newpage
\noindent
\subsection{The Case $p_1, p_2, p_3\to \infty$ and Constant 
Weight Functions (The Case of Legendre Polynomials)}

{\bf Theorem 2.4}\ \cite{6}-\cite{12aa}, \cite{arxiv-7}. 
{\it Suppose that
$\{\phi_j(x)\}_{j=0}^{\infty}$ is a complete orthonormal
system of Legendre polynomials
in the space $L_2([t, T])$.
Then$,$ for the iterated Stra\-to\-no\-vich stochastic integral of 
third multiplicity
$$
{\int\limits_t^{*}}^T
{\int\limits_t^{*}}^{t_3}
{\int\limits_t^{*}}^{t_2}
d{\bf w}_{t_1}^{(i_1)}
d{\bf w}_{t_2}^{(i_2)}d{\bf w}_{t_3}^{(i_3)}\ \ \ (i_1, i_2, i_3=1,\ldots,m)
$$
the following expansion 
\begin{equation}
\label{feto19001}
~~{\int\limits_t^{*}}^T
{\int\limits_t^{*}}^{t_3}
{\int\limits_t^{*}}^{t_2}
d{\bf w}_{t_1}^{(i_1)}
d{\bf w}_{t_2}^{(i_2)}d{\bf w}_{t_3}^{(i_3)}\ 
=
\hbox{\vtop{\offinterlineskip\halign{
\hfil#\hfil\cr
{\rm l.i.m.}\cr
$\stackrel{}{{}_{p_1,p_2,p_3\to \infty}}$\cr
}} }\sum_{j_1=0}^{p_1}\sum_{j_2=0}^{p_2}\sum_{j_3=0}^{p_3}
C_{j_3 j_2 j_1}\zeta_{j_1}^{(i_1)}\zeta_{j_2}^{(i_2)}\zeta_{j_3}^{(i_3)}
\end{equation}
that converges in the mean-square sense is valid, where
$$
C_{j_3 j_2 j_1}=\int\limits_t^T
\phi_{j_3}(s)\int\limits_t^s
\phi_{j_2}(s_1)
\int\limits_t^{s_1}
\phi_{j_1}(s_2)ds_2ds_1ds
$$
and
$$
\zeta_{j}^{(i)}=
\int\limits_t^T \phi_{j}(s) d{\bf w}_s^{(i)}
$$ 
are independent standard Gaussian random variables for various 
$i$ or $j$.}

{\bf Proof.} If we prove w.~p.~1 the following equalities
\begin{equation}
\label{ogo12}
~~~~~~~~\hbox{\vtop{\offinterlineskip\halign{
\hfil#\hfil\cr
{\rm l.i.m.}\cr
$\stackrel{}{{}_{p_1, p_3\to \infty}}$\cr
}} }
\sum\limits_{j_1=0}^{p_1}\sum\limits_{j_3=0}^{p_3}
C_{j_3 j_1 j_1}\zeta_{j_3}^{(i_3)}
=
\frac{1}{4}(T-t)^{3/2}\left(
\zeta_0^{(i_3)}+\frac{1}{\sqrt{3}}\zeta_1^{(i_3)}\right),
\end{equation}
\begin{equation}
\label{ogo13}
~~~~~~~~\hbox{\vtop{\offinterlineskip\halign{
\hfil#\hfil\cr
{\rm l.i.m.}\cr
$\stackrel{}{{}_{p_1, p_3\to \infty}}$\cr
}} }
\sum\limits_{j_1=0}^{p_1}\sum\limits_{j_3=0}^{p_3}
C_{j_3 j_3 j_1}\zeta_{j_1}^{(i_1)}
=
\frac{1}{4}(T-t)^{3/2}\left(
\zeta_0^{(i_1)}-\frac{1}{\sqrt{3}}\zeta_1^{(i_1)}\right),
\end{equation}
\begin{equation}
\label{ogo13a}
\hbox{\vtop{\offinterlineskip\halign{
\hfil#\hfil\cr
{\rm l.i.m.}\cr
$\stackrel{}{{}_{p_1, p_3\to \infty}}$\cr
}} }
\sum\limits_{j_1=0}^{p_1}\sum\limits_{j_3=0}^{p_3}
C_{j_1 j_3 j_1}\zeta_{j_3}^{(i_2)}
=0,
\end{equation}

\noindent
then in accordance with the
formulas (\ref{ogo12})--(\ref{ogo13a}), Theorem 1.1 (see (\ref{a3})), 
standard 
relations between iterated It\^{o} and Stratonovich 
stochastic integrals as well as in accordance with 
the formulas (they also follow
from Theorem 1.1)
$$
\frac{1}{2}\int\limits_t^T\int\limits_t^{\tau}dsd{\bf w}_{\tau}^{(i_3)}=
\frac{1}{4}(T-t)^{3/2}\left(
\zeta_0^{(i_3)}+\frac{1}{\sqrt{3}}\zeta_1^{(i_3)}\right)\ \ \ \hbox{w.~p.~1},
$$
$$
\frac{1}{2}\int\limits_t^T\int\limits_t^{\tau}d{\bf w}_{s}^{(i_1)}d\tau=
\frac{1}{4}(T-t)^{3/2}\left(
\zeta_0^{(i_1)}-\frac{1}{\sqrt{3}}\zeta_1^{(i_1)}\right)\ \ \ \hbox{w.~p.~1}
$$
we will have
$$
\int\limits_t^T\int\limits_t^{t_3}\int\limits_t^{t_2}
d{\bf w}_{t_1}^{(i_1)}d{\bf w}_{t_2}^{(i_2)}d{\bf w}_{t_3}^{(i_3)}=
\hbox{\vtop{\offinterlineskip\halign{
\hfil#\hfil\cr
{\rm l.i.m.}\cr
$\stackrel{}{{}_{p_1,p_2,p_3\to \infty}}$\cr
}} }\sum_{j_1=0}^{p_1}\sum_{j_2=0}^{p_2}\sum_{j_3=0}^{p_3}
C_{j_3 j_2 j_1}\zeta_{j_1}^{(i_1)}\zeta_{j_2}^{(i_2)}\zeta_{j_3}^{(i_3)}
-
$$
$$
-
{\bf 1}_{\{i_1=i_2\}}
\frac{1}{2}\int\limits_t^T\int\limits_t^{\tau}dsd{\bf w}_{\tau}^{(i_3)}-
{\bf 1}_{\{i_2=i_3\}}
\frac{1}{2}\int\limits_t^T\int\limits_t^{\tau}d{\bf w}_{s}^{(i_1)}d\tau\ \ \ 
\hbox{w.~p.~1}.
$$

It means that the expansion (\ref{feto19001}) will be proved.

Let us at first prove that
\begin{equation}
\label{ogo3}
\sum\limits_{j_1=0}^{\infty}C_{0 j_1 j_1}=\frac{1}{4}(T-t)^{3/2},
\end{equation}
\begin{equation}
\label{ogo4}
\sum\limits_{j_1=0}^{\infty}C_{1 j_1 j_1}=
\frac{1}{4\sqrt{3}}(T-t)^{3/2}.
\end{equation}

We have
$$
C_{000}=\frac{(T-t)^{3/2}}{6},
$$
$$
C_{0 j_1 j_1}=\int\limits_t^T\phi_0(s)\int\limits_t^s\phi_{j_1}(s_1)
\int\limits_t^{s_1}\phi_{j_1}(s_2)ds_2ds_1ds
=
$$
\begin{equation}
\label{ogo6}
=\frac{1}{2}\int\limits_t^T\phi_0(s)\left(
\int\limits_t^s\phi_{j_1}(s_1)ds_1\right)^2ds,\ \ \ j_1\ge 1,
\end{equation}

\noindent
where $\phi_j(s)$ looks as follows
\begin{equation}
\label{ogo7}
~~~~~~~~~~\phi_j(s)=\sqrt{\frac{2j+1}{T-t}}P_j\left(\left(
s-\frac{T+t}{2}\right)\frac{2}{T-t}\right),\ \ \ j\ge 0,
\end{equation}
where $P_j(x)$ is the Legendre polynomial.

Let us substitute (\ref{ogo7}) into (\ref{ogo6}) and
calculate $C_{0 j_1 j_1}$\  $(j_1\ge 1)$
$$
C_{0 j_1 j_1}=\frac{2j_1+1}{2(T-t)^{3/2}}
\int\limits_t^T
\left(\int\limits_{-1}^{z(s)}
P_{j_1}(y)\frac{T-t}{2}dy\right)^2ds=
$$
$$
=\frac{(2j_1+1)\sqrt{T-t}}{8}
\int\limits_t^T
\left(\int\limits_{-1}^{z(s)}
\frac{1}{2j_1+1}\left(P_{j_1+1}^{'}(y)-P_{j_1-1}^{'}(y)\right)dy
\right)^2ds=
$$
\begin{equation}
\label{ogo8}
=\frac{\sqrt{T-t}}{8(2j_1+1)}
\int\limits_t^T\left(P_{j_1+1}(z(s))-P_{j_1-1}(z(s))\right)^2ds,
\end{equation}

\noindent
where here and further
$$
z(s)=\left(s-\frac{T+t}{2}\right)\frac{2}{T-t}.
$$

In (\ref{ogo8})
we used the following well known properties of the Legendre polynomials
$$
P_j(y)=\frac{1}{2j+1}\left(P_{j+1}^{'}(y)-P_{j-1}^{'}(y)\right),\ \ \ 
P_j(-1)=(-1)^j,\ \ \ j\ge 1.
$$

Also, we denote 
$$
\frac{dP_j}{dy}(y)\stackrel{{\rm def}}{=}P_j^{'}(y).
$$

From (\ref{ogo8}) using the property of orthogonality of the Legendre 
polynomials, we get the following relation
$$
C_{0 j_1 j_1}=\frac{(T-t)^{3/2}}{16(2j_1+1)}
\int\limits_{-1}^1\left(P_{j_1+1}^2(y)+P_{j_1-1}^2(y)\right)dy=
$$
$$
=
\frac{(T-t)^{3/2}}{8(2j_1+1)}
\left(\frac{1}{2j_1+3}+\frac{1}{2j_1-1}\right),
$$

\vspace{1mm}
\noindent
where we used the property
$$
\int\limits_{-1}^1 P_j^2(y)dy=\frac{2}{2j+1},\ \ \ j\ge 0.
$$

Then
$$
\sum\limits_{j_1=0}^{\infty}C_{0 j_1 j_1}=
\frac{(T-t)^{3/2}}{6}+
$$
$$
+
\frac{(T-t)^{3/2}}{8}
\left(
\sum_{j_1=1}^{\infty}\frac{1}{(2j_1+1)(2j_1+3)}+
\sum_{j_1=1}^{\infty}\frac{1}{4j_1^2-1}\right)=
$$
$$
=\frac{(T-t)^{3/2}}{6}+\frac{(T-t)^{3/2}}{8}
\left(\sum_{j_1=1}^{\infty}\frac{1}{4j_1^2-1}-\frac{1}{3}
+\sum_{j_1=1}^{\infty}\frac{1}{4j_1^2-1}\right)=
$$
$$
=\frac{(T-t)^{3/2}}{6}+\frac{(T-t)^{3/2}}{8}
\left(\frac{1}{2}-\frac{1}{3}+\frac{1}{2}\right)=
\frac{(T-t)^{3/2}}{4}.
$$

The relation (\ref{ogo3}) is proved.

Let us check the correctness of (\ref{ogo4}). 
Let us represent $C_{1 j_1 j_1}$ in the form
$$
C_{1 j_1 j_1}=\frac{1}{2}\int\limits_t^T
\phi_1(s)\left(\int\limits_t^s\phi_{j_1}(s_1)ds_1\right)^2 ds=
$$
$$
=\frac{(T-t)^{3/2}(2j_1+1)\sqrt{3}}{16}
\int\limits_{-1}^{1}
P_1(y)\left(\int\limits_{-1}^y P_{j_1}(y_1)dy_1\right)^2 dy,\ \ \ j_1\ge 1.
$$

Since the functions
$$
\left(\int\limits_{-1}^y P_{j_1}(y_1)dy_1\right)^2,\ \ \ j_1\ge 1
$$ 
are even, then the functions
$$
P_1(y)\left(\int\limits_{-1}^y P_{j_1}(y_1)dy_1\right)^2 dy,\ \ \ j_1\ge 1
$$
are uneven. 
It means that $C_{1 j_1 j_1}=0$ $(j_1\ge 1).$ From the other side
$$
C_{100}=\frac{\sqrt{3}(T-t)^{3/2}}{16}
\int\limits_{-1}^1 y(y+1)^2 dy=\frac{(T-t)^{3/2}}{4\sqrt{3}}.
$$

Then 
$$
\sum\limits_{j_1=0}^{\infty}C_{1 j_1 j_1}=C_{100}+
\sum\limits_{j_1=1}^{\infty}C_{1 j_1 j_1}=
\frac{(T-t)^{3/2}}{4\sqrt{3}}.
$$

The relation (\ref{ogo4}) is proved. 

Let us prove the equality (\ref{ogo12}). Using (\ref{ogo4}), we get
$$
\sum\limits_{j_1=0}^{p_1}\sum\limits_{j_3=0}^{p_3}
C_{j_3 j_1 j_1}\zeta_{j_3}^{(i_3)}=
\sum\limits_{j_1=0}^{p_1}C_{0 j_1 j_1}\zeta_{0}^{(i_3)}+
\frac{(T-t)^{3/2}}{4\sqrt{3}}\zeta_{1}^{(i_3)}+
\sum\limits_{j_1=0}^{p_1}\sum\limits_{j_3=2}^{p_3}
C_{j_3 j_1 j_1}\zeta_{j_3}^{(i_3)}=
$$
\begin{equation}
\label{ogo15}
~~~~~~~~~=\sum\limits_{j_1=0}^{p_1}C_{0 j_1 j_1}\zeta_{0}^{(i_3)}+
\frac{(T-t)^{3/2}}{4\sqrt{3}}\zeta_{1}^{(i_3)}+
\sum\limits_{j_1=0}^{p_1}\ \ \sum\limits_{j_3=2, j_3 - {\rm even}}^{2j_1+2}
C_{j_3 j_1 j_1}\zeta_{j_3}^{(i_3)}.
\end{equation}

Since
$$
C_{j_3j_1j_1}=\frac{(T-t)^{3/2}(2j_1+1)\sqrt{2j_3+1}}{16}
\int\limits_{-1}^{1}
P_{j_3}(y)\left(\int\limits_{-1}^y P_{j_1}(y_1)dy_1\right)^2 dy
$$
and degree of the polynomial
$$
\left(\int\limits_{-1}^y P_{j_1}(y_1)dy_1\right)^2
$$ 
equals
to $2j_1+2$, then 
$C_{j_3j_1j_1}=0$ for $j_3>2j_1+2.$ It explains
that we put
$2j_1+2$ instead of $p_3$ on the right-hand side 
of the formula (\ref{ogo15}).

Moreover, the function 
$$
\left(\int\limits_{-1}^y P_{j_1}(y_1)dy_1\right)^2
$$
is even. It means that the function
$$
P_{j_3}(y)\left(\int\limits_{-1}^y P_{j_1}(y_1)dy_1\right)^2
$$
is uneven
for uneven
$j_3.$ It means that $C_{j_3 j_1j_1}=0$ for 
uneven
$j_3.$
That is why we 
summarize using even
$j_3$ on the right-hand side
of the formula (\ref{ogo15}).

Then we have
$$
\sum\limits_{j_1=0}^{p_1}\ \sum\limits_{j_3=2, j_3 - {\rm even}}^{2j_1+2}
C_{j_3 j_1 j_1}\zeta_{j_3}^{(i_3)}=
\sum\limits_{j_3=2, j_3 - {\rm even}}^{2p_1+2}\ \
\sum\limits_{j_1=(j_3-2)/2}^{p_1}
C_{j_3 j_1 j_1}\zeta_{j_3}^{(i_3)}=
$$
\begin{equation}
\label{ogo16}
=
\sum\limits_{j_3=2, j_3 - {\rm even}}^{2p_1+2}\ \
\sum\limits_{j_1=0}^{p_1}
C_{j_3 j_1 j_1}\zeta_{j_3}^{(i_3)}.
\end{equation}

We replaced $(j_3-2)/2$ by zero on the right-hand side
of the formula (\ref{ogo16}), since $C_{j_3j_1j_1}=0$ for 
$0\le j_1< (j_3-2)/2.$

Let us substitute (\ref{ogo16}) into (\ref{ogo15})
$$
\sum\limits_{j_1=0}^{p_1}\sum\limits_{j_3=0}^{p_3}
C_{j_3 j_1 j_1}\zeta_{j_3}^{(i_3)}=
\sum\limits_{j_1=0}^{p_1}C_{0 j_1 j_1}\zeta_{0}^{(i_3)}+
\frac{(T-t)^{3/2}}{4\sqrt{3}}\zeta_{1}^{(i_3)}+
$$
\begin{equation}
\label{ogo17}
+
\sum\limits_{j_3=2, j_3 - {\rm even}}^{2p_1+2}\ \
\sum\limits_{j_1=0}^{p_1}
C_{j_3 j_1 j_1}\zeta_{j_3}^{(i_3)}.
\end{equation}

It is easy to see that the right-hand side
of the formula (\ref{ogo17}) does not depend on $p_3.$ 

If we prove that
\begin{equation}
\label{ogo18}
\hbox{\vtop{\offinterlineskip\halign{
\hfil#\hfil\cr
{\rm lim}\cr
$\stackrel{}{{}_{p_1\to \infty}}$\cr
}} }
{\sf M}\left\{\left(
\sum\limits_{j_1=0}^{p_1}\sum\limits_{j_3=0}^{p_3}
C_{j_3 j_1 j_1}\zeta_{j_3}^{(i_3)}-
\frac{1}{4}(T-t)^{3/2}\left(
\zeta_0^{(i_3)}+\frac{1}{\sqrt{3}}\zeta_1^{(i_3)}\right)\right)^2\right\}=0,
\end{equation}
then the relaion (\ref{ogo12}) will be proved.

Using (\ref{ogo17}) and (\ref{ogo3}), we can rewrite the left-hand side 
of (\ref{ogo18})
in the following form
$$
\hbox{\vtop{\offinterlineskip\halign{
\hfil#\hfil\cr
{\rm lim}\cr
$\stackrel{}{{}_{p_1\to \infty}}$\cr
}} }
{\sf M}\left\{\left(
\left(\sum\limits_{j_1=0}^{p_1}C_{0j_1j_1}-
\frac{(T-t)^{3/2}}{4}\right)\zeta_0^{(i_3)}+
\sum\limits_{j_3=2, j_3 - {\rm even}}^{2p_1+2}\ \
\sum\limits_{j_1=0}^{p_1}
C_{j_3 j_1 j_1}\zeta_{j_3}^{(i_3)}\right)^2\right\}=
$$
$$
=\hbox{\vtop{\offinterlineskip\halign{
\hfil#\hfil\cr
{\rm lim}\cr
$\stackrel{}{{}_{p_1\to \infty}}$\cr
}} }\left(\sum\limits_{j_1=0}^{p_1}C_{0j_1j_1}-
\frac{(T-t)^{3/2}}{4}\right)^2+
\hbox{\vtop{\offinterlineskip\halign{
\hfil#\hfil\cr
{\rm lim}\cr
$\stackrel{}{{}_{p_1\to \infty}}$\cr
}} }
\sum\limits_{j_3=2, j_3 - {\rm even}}^{2p_1+2}
\left(\sum\limits_{j_1=0}^{p_1}
C_{j_3 j_1 j_1}\right)^2=
$$
\begin{equation}
\label{ogo19}
=\hbox{\vtop{\offinterlineskip\halign{
\hfil#\hfil\cr
{\rm lim}\cr
$\stackrel{}{{}_{p_1\to \infty}}$\cr
}} }
\sum\limits_{j_3=2, j_3 - {\rm even}}^{2p_1+2}
\left(\sum\limits_{j_1=0}^{p_1}
C_{j_3 j_1 j_1}\right)^2.
\end{equation}

If we prove that
\begin{equation}
\label{ogo20}
\hbox{\vtop{\offinterlineskip\halign{
\hfil#\hfil\cr
{\rm lim}\cr
$\stackrel{}{{}_{p_1\to \infty}}$\cr
}} }
\sum\limits_{j_3=2, j_3 - {\rm even}}^{2p_1+2}
\left(\sum\limits_{j_1=0}^{p_1}
C_{j_3 j_1 j_1}\right)^2=0,
\end{equation}

\noindent
then the relation (\ref{ogo12}) will be proved.

We have
$$
\sum\limits_{j_3=2, j_3 - {\rm even}}^{2p_1+2}
\left(\sum\limits_{j_1=0}^{p_1}
C_{j_3 j_1 j_1}\right)^2=
$$
$$
=
\frac{1}{4}
\sum\limits_{j_3=2, j_3 - {\rm even}}^{2p_1+2}
\left(\int\limits_t^T\phi_{j_3}(s)\sum\limits_{j_1=0}^{p_1}
\left(\int\limits_t^s\phi_{j_1}(s_1)ds_1\right)^2ds\right)^2=
$$
$$
=\frac{1}{4}
\sum\limits_{j_3=2, j_3 - {\rm even}}^{2p_1+2}
\left(\int\limits_t^T\phi_{j_3}(s)\left((s-t)-\sum\limits_{j_1=p_1+1}^{\infty}
\left(\int\limits_t^s\phi_{j_1}(s_1)ds_1\right)^2\right)ds\right)^2=
$$
$$
=\frac{1}{4}
\sum\limits_{j_3=2, j_3 - {\rm even}}^{2p_1+2}
\left(\int\limits_t^T\phi_{j_3}(s)\sum\limits_{j_1=p_1+1}^{\infty}
\left(\int\limits_t^s\phi_{j_1}(s_1)ds_1\right)^2 ds\right)^2\le
$$
\begin{equation}
\label{ogo21}
~~~~~~~\le\frac{1}{4}
\sum\limits_{j_3=2, j_3 - {\rm even}}^{2p_1+2}
\left(\int\limits_t^T|\phi_{j_3}(s)| \sum\limits_{j_1=p_1+1}^{\infty}
\left(\int\limits_t^s\phi_{j_1}(s_1)ds_1\right)^2 ds\right)^2.
\end{equation}

Obtaining (\ref{ogo21}), we used 
the Parseval equality
\begin{equation}
\label{ogo10}
\sum_{j_1=0}^{\infty}\left(\int\limits_t^s\phi_{j_1}(s_1)ds_1\right)^2=
\int\limits_t^T \left({\bf 1}_{\{s_1<s\}}\right)^2ds_1=s-t
\end{equation}

\noindent
and the orthogonality property of the Legendre polynomials 
\begin{equation}
\label{ogo11}
\int\limits_t^T\phi_{j_3}(s)(s-t)ds=0,\ \ \ j_3\ge 2.
\end{equation}

Then we have for $j_1\in{\bf N}$
$$
\left(\int\limits_t^s\phi_{j_1}(s_1)ds_1\right)^2=
\frac{(T-t)(2j_1+1)}{4}
\left(\int\limits_{-1}^{z(s)}
P_{j_1}(y)dy\right)^2=
$$
$$
=\frac{T-t}{4(2j_1+1)}
\left(\int\limits_{-1}^{z(s)}
\left(P_{j_1+1}^{'}(y)-P_{j_1-1}^{'}(y)\right)dy
\right)^2=
$$
$$
=\frac{T-t}{4(2j_1+1)}
\left(P_{j_1+1}\left(z(s)\right)-
P_{j_1-1}\left(z(s)\right)\right)^2
\le
$$
\begin{equation}
\label{ogo22}
\le
\frac{T-t}{2(2j_1+1)}
\left(P_{j_1+1}^2\left(z(s)\right)+
P_{j_1-1}^2\left(z(s)\right)\right).
\end{equation}

\vspace{3mm}

Remind that for the Legendre polynomials the following 
estimate is correct
\begin{equation}
\label{ogo23}
~~~~~~~~~\left|P_j(y)\right|<\frac{K}{\sqrt{j+1}(1-y^2)^{1/4}},\ \ \ 
y\in (-1, 1),\ \ \ j\in {\bf N},
\end{equation}
where constant $K$ does not depend on $y$ and $j.$

The estimate (\ref{ogo23}) can be rewritten for the 
function $\phi_j(s)$ in 
the following form
$$
|\phi_j(s)|< \sqrt{\frac{2j+1}{j+1}}\frac{K}{\sqrt{T-t}}
\frac{1}
{\left(1-z^2(s)\right)^{1/4}}
<
$$
\begin{equation}
\label{ogo24}
<\frac{K_1}{\sqrt{T-t}}
\frac{1}
{\left(1-z^2(s)\right)^{1/4}},
\end{equation}
where
$K_1=K\sqrt{2},$\  $s\in (t, T).$

Let us estimate the right-hand side of (\ref{ogo22}) using the estimate
(\ref{ogo23})
$$
\left(\int\limits_t^s\phi_{j_1}(s_1)ds_1\right)^2 <
\frac{T-t}{2(2j_1+1)}\left(\frac{K^2}{j_1+2}+\frac{K^2}{j_1}\right)
\frac{1}
{(1-(z(s))^2)^{1/2}} <
$$
\begin{equation}
\label{ogo25}
<
\frac{(T-t)K^2}{2j_1^2}
\frac{1}
{(1-(z(s))^2)^{1/2}},
\end{equation}

\noindent
where $s\in(t, T),$\  $j_1\in {\bf N}.$

Substituting the estimate (\ref{ogo25}) into the relation (\ref{ogo21})
and using in (\ref{ogo21}) the estimate (\ref{ogo24})
for $|\phi_{j_3}(s)|$, we obtain
$$
\sum\limits_{j_3=2, j_3 - {\rm even}}^{2p_1+2}
\left(\sum\limits_{j_1=0}^{p_1}
C_{j_3 j_1 j_1}\right)^2<
$$
$$
<
\frac{(T-t)K^4 K_1^2}{16}
\sum\limits_{j_3=2, j_3 - {\rm even}}^{2p_1+2}
\left(\int\limits_t^T
\frac{ds}
{\left(1-\left(z(s)\right)^2
\right)^{3/4}}\sum\limits_{j_1=p_1+1}^{\infty}\frac{1}{j_1^2}
\right)^2=
$$
\begin{equation}
\label{ogo26}
~~~~~~ =\frac{(T-t)^3K^4 K_1^2(p_1+1)}{64}
\left(\int\limits_{-1}^1
\frac{dy}
{\left(1-y^2\right)^{3/4}}\right)^2\left(
\sum\limits_{j_1=p_1+1}^{\infty}\frac{1}{j_1^2}
\right)^2.
\end{equation}

Since
\begin{equation}
\label{ogo27}
\int\limits_{-1}^1
\frac{dy}
{\left(1-y^2\right)^{3/4}}<\infty
\end{equation}
and
\begin{equation}
\label{ogo28}
\sum\limits_{j_1=p_1+1}^{\infty}\frac{1}{j_1^2}
\le \int\limits_{p_1}^{\infty}\frac{dx}{x^2}=\frac{1}{p_1},
\end{equation}

\noindent
then from (\ref{ogo26}) we find
\begin{equation}
\label{ogo29}
\sum\limits_{j_3=2, j_3 - {\rm even}}^{2p_1+2}
\left(\sum\limits_{j_1=0}^{p_1}
C_{j_3 j_1 j_1}\right)^2<\frac{C(T-t)^3 (p_1+1)}{p_1^2}\ \  \to\  0\ \ \ 
\hbox{if}\ \
p_1\to \infty,
\end{equation}

\noindent
where constant $C$ does not depend on $p_1$ and $T-t.$
The relation (\ref{ogo29}) implies (\ref{ogo20}), and the relation
(\ref{ogo20})
implies the correctness of the formula (\ref{ogo12}).

Let us prove the equaity (\ref{ogo13}). 
Let us at first prove that
\begin{equation}
\label{ogo30}
\sum\limits_{j_3=0}^{\infty}C_{j_3 j_3 0}=\frac{1}{4}(T-t)^{3/2},
\end{equation}
\begin{equation}
\label{ogo31}
\sum\limits_{j_3=0}^{\infty}C_{j_3 j_3 1}=
-\frac{1}{4\sqrt{3}}(T-t)^{3/2}.
\end{equation}

We have
$$
\sum_{j_3=0}^{\infty}C_{j_3 j_3 0}=C_{000}+
\sum_{j_3=1}^{\infty}C_{j_3 j_3 0},
$$
$$
C_{000}=\frac{(T-t)^{3/2}}{6},
$$
$$
C_{j_3 j_3 0}=\frac{(T-t)^{3/2}}{16(2j_3+1)}
\int\limits_{-1}^1\left(P_{j_3+1}^2(y)+P_{j_3-1}^2(y)\right)dy=
$$
$$
=
\frac{(T-t)^{3/2}}{8(2j_3+1)}
\left(\frac{1}{2j_3+3}+\frac{1}{2j_3-1}\right),\ \ \ j_3\ge 1.
$$

Then
$$
\sum\limits_{j_3=0}^{\infty}C_{j_3 j_3 0}=
\frac{(T-t)^{3/2}}{6}+
$$
$$
+
\frac{(T-t)^{3/2}}{8}
\left(
\sum_{j_3=1}^{\infty}\frac{1}{(2j_3+1)(2j_3+3)}+
\sum_{j_3=1}^{\infty}\frac{1}{4j_3^2-1}\right)=
$$
$$
=\frac{(T-t)^{3/2}}{6}+\frac{(T-t)^{3/2}}{8}
\left(\sum_{j_3=1}^{\infty}\frac{1}{4j_3^2-1}-\frac{1}{3}
+\sum_{j_3=1}^{\infty}\frac{1}{4j_3^2-1}\right)=
$$
$$
=\frac{(T-t)^{3/2}}{6}+\frac{(T-t)^{3/2}}{8}
\left(\frac{1}{2}-\frac{1}{3}+\frac{1}{2}\right)=
\frac{(T-t)^{3/2}}{4}.
$$

The relation (\ref{ogo30}) is proved.
Let us check the equality (\ref{ogo31}). We have
$$
C_{j_3 j_3 j_1}=\int\limits_t^T
\phi_{j_3}(s)\int\limits_t^s
\phi_{j_3}(s_1)\int\limits_t^{s_1}
\phi_{j_1}(s_2)ds_2ds_1ds=
$$
$$
=
\int\limits_t^T\phi_{j_1}(s_2)ds_2
\int\limits_{s_2}^T
\phi_{j_3}(s_1)ds_1\int\limits_{s_1}^T
\phi_{j_3}(s)ds=
$$
$$
=\frac{1}{2}\int\limits_t^T
\phi_{j_1}(s_2)\left(\int\limits_{s_2}^T\phi_{j_3}(s_1)ds_1\right)^2 ds_2=
$$
\begin{equation}
\label{ogo33}
=\frac{(T-t)^{3/2}(2j_3+1)\sqrt{2j_1+1}}{16}
\int\limits_{-1}^{1}
P_{j_1}(y)\left(\int\limits_{y}^1 P_{j_3}(y_1)dy_1\right)^2 dy,\ \ \ j_3\ge 1.
\end{equation}

Since the functions
$$
\left(\int\limits_{y}^1 P_{j_3}(y_1)dy_1\right)^2,\ \ \ j_3\ge 1
$$
are even, then the functions
$$
P_1(y)\left(\int\limits_{y}^1 P_{j_3}(y_1)dy_1\right)^2 dy,\ \ \ j_3\ge 1
$$
are uneven. It means that $C_{j_3 j_3 1}=0$ $(j_3\ge 1).$

Moreover,
$$
C_{001}=\frac{\sqrt{3}(T-t)^{3/2}}{16}
\int\limits_{-1}^1 y(1-y)^2 dy=-\frac{(T-t)^{3/2}}{4\sqrt{3}}.
$$

\noindent
\par
Then 
$$
\sum\limits_{j_3=0}^{\infty}C_{j_3 j_3 1}=C_{001}+
\sum\limits_{j_3=1}^{\infty}C_{j_3 j_3 1}=
-\frac{(T-t)^{3/2}}{4\sqrt{3}}.
$$

\noindent
\par
The relation (\ref{ogo31}) is proved.
Using the obtained results, we get
$$
\sum\limits_{j_1=0}^{p_1}\sum\limits_{j_3=0}^{p_3}
C_{j_3 j_3 j_1}\zeta_{j_1}^{(i_1)}=
\sum\limits_{j_3=0}^{p_3}C_{j_3 j_3 0}\zeta_{0}^{(i_1)}-
\frac{(T-t)^{3/2}}{4\sqrt{3}}\zeta_{1}^{(i_1)}+
\sum\limits_{j_3=0}^{p_3}\sum\limits_{j_1=2}^{p_1}
C_{j_3 j_3 j_1}\zeta_{j_1}^{(i_1)}=
$$
\begin{equation}
\label{ogoo5}
~~~~~~~~~=\sum\limits_{j_3=0}^{p_3}C_{j_3 j_3 0}\zeta_{0}^{(i_1)}-
\frac{(T-t)^{3/2}}{4\sqrt{3}}\zeta_{1}^{(i_1)}+
\sum\limits_{j_3=0}^{p_3}\ \ \sum\limits_{j_1=2, j_1 - {\rm even}}^{2j_3+2}
C_{j_3 j_3 j_1}\zeta_{j_1}^{(i_1)}.
\end{equation}

Since 
$$
C_{j_3j_3j_1}=
\frac{(T-t)^{3/2}(2j_3+1)\sqrt{2j_1+1}}{16}
\int\limits_{-1}^{1}
P_{j_1}(y)\left(\int\limits_{y}^1 P_{j_3}(y_1)dy_1\right)^2 dy,\ \ \ j_3\ge 1,
$$
and degree of the polynomial
$$
\left(\int\limits_{y}^1 P_{j_3}(y_1)dy_1\right)^2
$$ 
equals to 
$2j_3+2$, 
then
$C_{j_3j_3j_1}=0$ for $j_1>2j_3+2.$ It explains
that we put $2j_3+2$ instead of $p_1$ on the right-hand side 
of the formula (\ref{ogoo5}).

Moreover, the function 
$$
\left(\int\limits_{y}^1 P_{j_3}(y_1)dy_1\right)^2
$$
is even.
It means that the 
function  
$$
P_{j_1}(y)\left(\int\limits_{y}^1 P_{j_3}(y_1)dy_1\right)^2
$$ 
is uneven for uneven $j_1.$
It means that $C_{j_3 j_3j_1}=0$ for uneven  $j_1.$
It explains the summation with respect to
even $j_1$ on the right-hand side of (\ref{ogoo5}).

Then we have
$$
\sum\limits_{j_3=0}^{p_3}\ \ \sum\limits_{j_1=2, j_1 - {\rm even}}^{2j_3+2}
C_{j_3 j_3 j_1}\zeta_{j_1}^{(i_1)}=
\sum\limits_{j_1=2, j_1 - {\rm even}}^{2p_3+2}\ \ 
\sum\limits_{j_3=(j_1-2)/2}^{p_3}
C_{j_3 j_3 j_1}\zeta_{j_1}^{(i_1)}
=
$$
\begin{equation}
\label{ogoo11}
=\sum\limits_{j_1=2, j_1 - {\rm even}}^{2p_3+2}\ \ 
\sum\limits_{j_3=0}^{p_3}
C_{j_3 j_3 j_1}\zeta_{j_1}^{(i_1)}.
\end{equation}

We replaced $(j_1-2)/2$ by zero on the right-hand side
of (\ref{ogoo11}), since $C_{j_3j_3j_1}=0$ for
$0\le j_3< (j_1-2)/2.$

Let us substitute (\ref{ogoo11}) into (\ref{ogoo5})
$$
\sum\limits_{j_1=0}^{p_1}\sum\limits_{j_3=0}^{p_3}
C_{j_3 j_3 j_1}\zeta_{j_1}^{(i_1)}=
\sum\limits_{j_3=0}^{p_3}C_{j_3 j_3 0}\zeta_{0}^{(i_1)}-
\frac{(T-t)^{3/2}}{4\sqrt{3}}\zeta_{1}^{(i_1)}+
$$
\begin{equation}
\label{ogoo12}
+
\sum\limits_{j_1=2, j_1 - {\rm even}}^{2p_3+2}\ \ 
\sum\limits_{j_3=0}^{p_3}
C_{j_3 j_3 j_1}\zeta_{j_1}^{(i_1)}.
\end{equation}

It is easy to see that the right-hand side of the formula 
(\ref{ogoo12}) does not depend on $p_1.$

If we prove that
\begin{equation}
\label{ogoo13}
\hbox{\vtop{\offinterlineskip\halign{
\hfil#\hfil\cr
{\rm lim}\cr
$\stackrel{}{{}_{p_3\to \infty}}$\cr
}} }
{\sf M}\left\{\left(
\sum\limits_{j_1=0}^{p_1}\sum\limits_{j_3=0}^{p_3}
C_{j_3 j_3 j_1}\zeta_{j_1}^{(i_1)}-
\frac{1}{4}(T-t)^{3/2}\left(
\zeta_0^{(i_1)}-\frac{1}{\sqrt{3}}\zeta_1^{(i_1)}\right)\right)^2\right\}=0,
\end{equation}
then (\ref{ogo13}) will be proved.

Using (\ref{ogoo12}) and (\ref{ogo30}), (\ref{ogo31}), we can rewrite 
the left-hand side of the formula (\ref{ogoo13}) in the 
following form
$$
\hbox{\vtop{\offinterlineskip\halign{
\hfil#\hfil\cr
{\rm lim}\cr
$\stackrel{}{{}_{p_3\to \infty}}$\cr
}} }
{\sf M}\left\{\left(
\left(\sum\limits_{j_3=0}^{p_3}C_{j_3j_3 0}-
\frac{(T-t)^{3/2}}{4}\right)\zeta_0^{(i_1)}+
\sum\limits_{j_1=2, j_1 - {\rm even}}^{2p_3+2}\ \
\sum\limits_{j_3=0}^{p_3}
C_{j_3 j_3 j_1}\zeta_{j_1}^{(i_1)}\right)^2\right\}=
$$
$$
=\hbox{\vtop{\offinterlineskip\halign{
\hfil#\hfil\cr
{\rm lim}\cr
$\stackrel{}{{}_{p_3\to \infty}}$\cr
}} }\left(\sum\limits_{j_3=0}^{p_1}C_{j_3j_3 0}-
\frac{(T-t)^{3/2}}{4}\right)^2+
\hbox{\vtop{\offinterlineskip\halign{
\hfil#\hfil\cr
{\rm lim}\cr
$\stackrel{}{{}_{p_3\to \infty}}$\cr
}} }
\sum\limits_{j_1=2, j_1 - {\rm even}}^{2p_3+2}
\left(\sum\limits_{j_3=0}^{p_3}
C_{j_3 j_3 j_1}\right)^2=
$$
$$
=\hbox{\vtop{\offinterlineskip\halign{
\hfil#\hfil\cr
{\rm lim}\cr
$\stackrel{}{{}_{p_3\to \infty}}$\cr
}} }
\sum\limits_{j_1=2, j_1 - {\rm even}}^{2p_3+2}
\left(\sum\limits_{j_3=0}^{p_3}
C_{j_3 j_3 j_1}\right)^2.
$$

If we prove that
\begin{equation}
\label{ogoo15}
\hbox{\vtop{\offinterlineskip\halign{
\hfil#\hfil\cr
{\rm lim}\cr
$\stackrel{}{{}_{p_3\to \infty}}$\cr
}} }
\sum\limits_{j_1=2, j_1 - {\rm even}}^{2p_3+2}
\left(\sum\limits_{j_3=0}^{p_3}
C_{j_3 j_3 j_1}\right)^2=0,
\end{equation}
then the relation (\ref{ogo13}) will be proved.

From (\ref{ogo33}) we obtain
$$
\sum\limits_{j_1=2, j_1 - {\rm even}}^{2p_3+2}
\left(\sum\limits_{j_3=0}^{p_3}
C_{j_3 j_3 j_1}\right)^2=
$$
$$
=
\frac{1}{4}
\sum\limits_{j_1=2, j_1 - {\rm even}}^{2p_3+2}
\left(\int\limits_t^T\phi_{j_1}(s_2)\sum\limits_{j_3=0}^{p_3}
\left(\int\limits_{s_2}^T\phi_{j_3}(s_1)ds_1\right)^2ds_2\right)^2=
$$
$$
=\frac{1}{4}
\sum\limits_{j_1=2, j_1 - {\rm even}}^{2p_3+2}
\left(\int\limits_t^T\phi_{j_1}(s_2)\left((T-s_2)-
\sum\limits_{j_3=p_3+1}^{\infty}
\left(\int\limits_{s_2}^T\phi_{j_3}(s_1)ds_1\right)^2\right)ds_2\right)^2=
$$
$$
=\frac{1}{4}
\sum\limits_{j_1=2, j_1 - {\rm even}}^{2p_3+2}
\left(\int\limits_t^T\phi_{j_1}(s_2)\sum\limits_{j_3=p_3+1}^{\infty}
\left(\int\limits_{s_2}^T\phi_{j_3}(s_1)ds_1\right)^2 ds_2\right)^2\le
$$
\begin{equation}
\label{ogoo21}
~~~~~~~~\le\frac{1}{4}
\sum\limits_{j_1=2, j_1 - {\rm even}}^{2p_3+2}
\left(\int\limits_t^T|\phi_{j_1}(s_2)|\sum\limits_{j_3=p_3+1}^{\infty}
\left(\int\limits_{s_2}^T\phi_{j_3}(s_1)ds_1\right)^2 ds_2\right)^2.
\end{equation}

In order to get (\ref{ogoo21}) we used 
the Parseval equality 
\begin{equation}
\label{ogo10ee}
\sum_{j_1=0}^{\infty}\left(\int\limits_s^T\phi_{j_1}(s_1)ds_1\right)^2=
\int\limits_t^T \left({\bf 1}_{\{s<s_1\}}\right)^2ds_1=T-s
\end{equation}
and the orthogonality property of the Legendre polynomials
\begin{equation}
\label{ogo11e}
\int\limits_t^T\phi_{j_3}(s)(T-s)ds=0,\ \ \ j_3\ge 2.
\end{equation}

Then we have for $j_3\in{\bf N}$
$$
\left(\int\limits_{s_2}^T\phi_{j_3}(s_1)ds_1\right)^2=
\frac{(T-t)}{4(2j_3+1)}
\left(P_{j_3+1}\left(
z(s_2)\right)-
P_{j_3-1}\left(
z(s_2)\right)\right)^2\le
$$
$$
\le
\frac{T-t}{2(2j_3+1)}
\left(P_{j_3+1}^2\left(
z(s_2)\right)+
P_{j_3-1}^2\left(
z(s_2)\right)\right)
<
$$
$$
<\frac{T-t}{2(2j_3+1)}\left(\frac{K^2}{j_3+2}+\frac{K^2}{j_3}\right)
\frac{1}
{(1-(z(s_2))^2)^{1/2}} <
$$
\vspace{2mm}
\begin{equation}
\label{ogoo25}
< \frac{(T-t)K^2}{2j_3^2}
\frac{1}
{(1-(z(s_2))^2)^{1/2}},\ \ \ s_2\in(t, T).
\end{equation}

\vspace{1mm}

In order to get (\ref{ogoo25}) we used the estimate
(\ref{ogo23}). 

Substituting the estimate (\ref{ogoo25}) into the relation (\ref{ogoo21})
and using in (\ref{ogoo21}) the estimate (\ref{ogo24})
for $|\phi_{j_1}(s_2)|$, we obtain
$$
\sum\limits_{j_1=2, j_1 - {\rm even}}^{2p_3+2}
\left(\sum\limits_{j_3=0}^{p_3}
C_{j_3 j_3 j_1}\right)^2<
$$
$$
<
\frac{(T-t)K^4 K_1^2}{16}
\sum\limits_{j_1=2, j_1 - {\rm even}}^{2p_3+2}
\left(\int\limits_t^T
\frac{ds_2}
{(1-z^2(s_2))^{3/4}}\sum\limits_{j_3=p_3+1}^{\infty}\frac{1}{j_3^2}
\right)^2=
$$
\begin{equation}
\label{ogoo26}
~~~~~~~ =\frac{(T-t)^3K^4 K_1^2(p_3+1)}{64}
\left(\int\limits_{-1}^1
\frac{dy}
{\left(1-y^2\right)^{3/4}}\right)^2\left(
\sum\limits_{j_3=p_3+1}^{\infty}\frac{1}{j_3^2}
\right)^2.
\end{equation}

Using (\ref{ogo27}) and (\ref{ogo28}) in
(\ref{ogoo26}), we get
\begin{equation}
\label{ogoo29}
\sum\limits_{j_1=2, j_1 - {\rm even}}^{2p_3+2}
\left(\sum\limits_{j_3=0}^{p_3}
C_{j_3 j_3 j_1}\right)^2<\frac{C(T-t)^3 (p_3+1)}{p_3^2}\ \ \to\ 0\ \ \
\hbox{with}\ \ p_3\to \infty,
\end{equation}
where constant $C$ does not depend on $p_3$ and $T-t.$

The relation (\ref{ogoo29}) implies (\ref{ogoo15}), and the
relation
(\ref{ogoo15})
implies the correctness 
of the formula (\ref{ogo13}). The relation (\ref{ogo13}) is proved.

Let us prove the equality (\ref{ogo13a}).
Since $\psi_1(\tau),$ $\psi_2(\tau),$ $\psi_3(\tau)\equiv 1,$
then the following relation 
for the Fourier coefficients is correct
$$
C_{j_1 j_1 j_3}+C_{j_1 j_3 j_1}+C_{j_3 j_1 j_1}=\frac{1}{2}
C_{j_1}^2 C_{j_3},
$$ 
where $C_j=0$ for $j\ge 1$ and $C_0=\sqrt{T-t}.$
Then\ w.~p.~1
\begin{equation}
\label{sodom31}
\hbox{\vtop{\offinterlineskip\halign{
\hfil#\hfil\cr
{\rm l.i.m.}\cr
$\stackrel{}{{}_{p_1, p_3\to \infty}}$\cr
}} }
\sum\limits_{j_1=0}^{p_1}\sum\limits_{j_3=0}^{p_3}
C_{j_1 j_3 j_1}\zeta_{j_3}^{(i_2)}=
\hbox{\vtop{\offinterlineskip\halign{
\hfil#\hfil\cr
{\rm l.i.m.}\cr
$\stackrel{}{{}_{p_1, p_3\to \infty}}$\cr
}} }
\sum\limits_{j_1=0}^{p_1}\sum\limits_{j_3=0}^{p_3}
\left(\frac{1}{2}C_{j_1}^2 C_{j_3}-C_{j_1 j_1 j_3}-C_{j_3 j_1 j_1}
\right)\zeta_{j_3}^{(i_2)}.
\end{equation}

Therefore, considering 
(\ref{ogo12}) and (\ref{ogo13}), we can write 
w.~p.~1 
$$
\hbox{\vtop{\offinterlineskip\halign{
\hfil#\hfil\cr
{\rm l.i.m.}\cr
$\stackrel{}{{}_{p_1, p_3\to \infty}}$\cr
}} }
\sum\limits_{j_1=0}^{p_1}\sum\limits_{j_3=0}^{p_3}
C_{j_1 j_3 j_1}\zeta_{j_3}^{(i_2)}=\frac{1}{2}C_0^3\zeta_0^{(i_2)}-
$$
$$
-
\hbox{\vtop{\offinterlineskip\halign{
\hfil#\hfil\cr
{\rm l.i.m.}\cr
$\stackrel{}{{}_{p_1, p_3\to \infty}}$\cr
}} }
\sum\limits_{j_1=0}^{p_1}\sum\limits_{j_3=0}^{p_3}
C_{j_1 j_1 j_3}\zeta_{j_3}^{(i_2)}-
\hbox{\vtop{\offinterlineskip\halign{
\hfil#\hfil\cr
{\rm l.i.m.}\cr
$\stackrel{}{{}_{p_1, p_3\to \infty}}$\cr
}} }
\sum\limits_{j_1=0}^{p_1}\sum\limits_{j_3=0}^{p_3}
C_{j_3 j_1 j_1}\zeta_{j_3}^{(i_2)}
=
$$
$$
=\frac{1}{2}(T-t)^{3/2}
\zeta_0^{(i_2)}
-\frac{1}{4}(T-t)^{3/2}\left(
\zeta_0^{(i_2)}-\frac{1}{\sqrt{3}}\zeta_1^{(i_2)}\right)
-
$$
\begin{equation}
\label{sodom3}
-\frac{1}{4}(T-t)^{3/2}\left(
\zeta_0^{(i_2)}+\frac{1}{\sqrt{3}}\zeta_1^{(i_2)}\right)=0.
\end{equation}

\vspace{1mm}

The relation (\ref{ogo13a}) is proved. Theorem 2.4 is proved.

It is easy to see that the formula (\ref{feto19001})  
can be proved for the case $i_1=i_2=i_3$  
using the It\^{o} formula
$$
{\int\limits_t^{*}}^T
{\int\limits_t^{*}}^{t_3}
{\int\limits_t^{*}}^{t_2}
d{\bf w}_{t_1}^{(i_1)}d{\bf w}_{t_2}^{(i_1)}d{\bf w}_{t_3}^{(i_1)}=
\frac{1}{6}\left(\int\limits_t^T d{\bf w}_{s}^{(i_1)}\right)^3=
\frac{1}{6}\left(C_0\zeta_{0}^{(i_1)}\right)^3=
C_{000}\left(\zeta_{0}^{(i_1)}\right)^3,
$$

\noindent
where the equality is fulfilled w.~p.~1.

\subsection{The Case $p_1, p_2, p_3\to \infty,$ Binomial 
Weight Functions, and Additional
Restrictive Conditions (The Case of Legendre Polynomials)}

Let us consider the following generalization of Theorem 2.4.

{\bf Theorem 2.5}\ \cite{6}-\cite{12aa}, \cite{arxiv-7}. {\it Suppose that
$\{\phi_j(x)\}_{j=0}^{\infty}$ is a complete orthonormal
system of Legendre polynomials
in the space $L_2([t, T])$.
Then$,$ for the iterated Stratonovich stochastic integral of 
third multiplicity
$$
I_{{l_1l_2l_3}_{T,t}}^{*(i_1i_2i_3)}={\int\limits_t^{*}}^T(t-t_3)^{l_3}
{\int\limits_t^{*}}^{t_3}(t-t_2)^{l_2}
{\int\limits_t^{*}}^{t_2}(t-t_1)^{l_1}
d{\bf w}_{t_1}^{(i_1)}
d{\bf w}_{t_2}^{(i_2)}d{\bf w}_{t_3}^{(i_3)}
$$
the following expansion 
\begin{equation}
\label{feto1900}
I_{{l_1l_2l_3}_{T,t}}^{*(i_1i_2i_3)}=
\hbox{\vtop{\offinterlineskip\halign{
\hfil#\hfil\cr
{\rm l.i.m.}\cr
$\stackrel{}{{}_{p_1,p_2,p_3\to \infty}}$\cr
}} }\sum_{j_1=0}^{p_1}\sum_{j_2=0}^{p_2}\sum_{j_3=0}^{p_3}
C_{j_3 j_2 j_1}\zeta_{j_1}^{(i_1)}\zeta_{j_2}^{(i_2)}\zeta_{j_3}^{(i_3)}
\end{equation}
that converges in the mean-square sense 
is valid for each of the following cases

\vspace{2mm}
\noindent
{\rm 1}.\ $i_1\ne i_2,\ i_2\ne i_3,\ i_1\ne i_3$\ and
$l_1, l_2, l_3=0, 1, 2,\ldots $\\
{\rm 2}.\ $i_1=i_2\ne i_3$ and $l_1=l_2\ne l_3$\ and
$l_1, l_2, l_3=0, 1, 2,\ldots $\\
{\rm 3}.\ $i_1\ne i_2=i_3$ and $l_1\ne l_2=l_3$\ and
$l_1, l_2, l_3=0, 1, 2,\ldots $\\
{\rm 4}.\ $i_1, i_2, i_3=1,\ldots,m;$ $l_1=l_2=l_3=l$\ and $l=0, 1, 
2,\ldots,$\\

\vspace{-3mm}
\noindent
where $i_1, i_2, i_3=1,\ldots,m,$
$$
C_{j_3 j_2 j_1}=\int\limits_t^T(t-s)^{l_3}\phi_{j_3}(s)
\int\limits_t^s(t-s_1)^{l_2}\phi_{j_2}(s_1)
\int\limits_t^{s_1}(t-s_2)^{l_1}\phi_{j_1}(s_2)ds_2ds_1ds,
$$
and
$$
\zeta_{j}^{(i)}=
\int\limits_t^T \phi_{j}(s) d{\bf w}_s^{(i)}
$$ 
are independent standard Gaussian random variables for various 
$i$ or $j$.}

{\bf Proof.} Case 1 directly follows from (\ref{a3}).
Let us consider Case 2, i.e. $i_1=i_2\ne i_3$, $l_1=l_2=l\ne l_3$, and
$l_1, l_3 = 0, 1, 2,\ldots$ So, we prove 
the following expansion 
\begin{equation}
\label{ogo101}
I_{{l_1 l_1 l_3}_{T,t}}^{*(i_1i_1i_3)}=
\hbox{\vtop{\offinterlineskip\halign{
\hfil#\hfil\cr
{\rm l.i.m.}\cr
$\stackrel{}{{}_{p_1,p_2,p_3\to \infty}}$\cr
}} }\sum_{j_1=0}^{p_1}\sum_{j_2=0}^{p_2}\sum_{j_3=0}^{p_3}
C_{j_3 j_2 j_1}\zeta_{j_1}^{(i_1)}\zeta_{j_2}^{(i_1)}\zeta_{j_3}^{(i_3)}\ \ \
(i_1, i_2, i_3=1,\ldots,m),
\end{equation}

\noindent
where 
$l_1, l_3=0, 1, 2,\ldots$ $(l_1=l)$ and
\begin{equation}
\label{ogo199}
C_{j_3 j_2 j_1}=\int\limits_t^T
\phi_{j_3}(s)(t-s)^{l_3}\int\limits_t^s(t-s_1)^{l}
\phi_{j_2}(s_1)
\int\limits_t^{s_1}(t-s_2)^l
\phi_{j_1}(s_2)ds_2ds_1ds.
\end{equation}

If we prove w.~p.~1 the formula
\begin{equation}
\label{ogo200}
~~~~~~\hbox{\vtop{\offinterlineskip\halign{
\hfil#\hfil\cr
{\rm l.i.m.}\cr
$\stackrel{}{{}_{p_1, p_3\to \infty}}$\cr
}} }
\sum\limits_{j_1=0}^{p_1}\sum\limits_{j_3=0}^{p_3}
C_{j_3 j_1 j_1}\zeta_{j_3}^{(i_3)}=
\frac{1}{2}\int\limits_t^T(t-s)^{l_3}
\int\limits_t^s(t-s_1)^{2l}ds_1d{\bf w}_s^{(i_3)},
\end{equation}
where 
coefficients $C_{j_3 j_1 j_1}$ are defined by (\ref{ogo199}), 
then using Theorem 1.1 and 
standard relations between iterated
It\^{o} and Stratonovich stochastic integrals, we obtain 
the expansion (\ref{ogo101}).

Using Theorem 
1.1, we obtain
$$
\frac{1}{2}\int\limits_t^T(t-s)^{l_3}
\int\limits_t^s(t-s_1)^{2l}ds_1d{\bf w}_s^{(i_3)}=
\frac{1}{2}\sum\limits_{j_3=0}^{2l+l_3+1}
\tilde C_{j_3}\zeta_{j_3}^{(i_3)}\ \ \ \hbox{w.~p.~1},
$$
where
$$
\tilde C_{j_3}=
\int\limits_t^T
\phi_{j_3}(s)(t-s)^{l_3}\int\limits_t^s(t-s_1)^{2l}ds_1ds.
$$

Then
$$
\sum\limits_{j_3=0}^{p_3}\sum\limits_{j_1=0}^{p_1}
C_{j_3 j_1 j_1}\zeta_{j_3}^{(i_3)}-
\frac{1}{2}\sum\limits_{j_3=0}^{2l+l_3+1}
\tilde C_{j_3}\zeta_{j_3}^{(i_3)}=
$$
$$
=\sum\limits_{j_3=0}^{2l+l_3+1}
\left(\sum\limits_{j_1=0}^{p_1}
C_{j_3 j_1 j_1}-\frac{1}{2}\tilde C_{j_3}\right)
\zeta_{j_3}^{(i_3)}+
\sum\limits_{j_3=2l+l_3+2}^{p_3}
\sum\limits_{j_1=0}^{p_1}
C_{j_3 j_1 j_1}\zeta_{j_3}^{(i_3)}.
$$

\noindent
\par
Therefore,
$$
\hbox{\vtop{\offinterlineskip\halign{
\hfil#\hfil\cr
{\rm lim}\cr
$\stackrel{}{{}_{p_1,p_3\to \infty}}$\cr
}} }
{\sf M}\left\{\left(
\sum\limits_{j_3=0}^{p_3}\sum\limits_{j_1=0}^{p_1}C_{j_3j_1 j_1}
\zeta_{j_3}^{(i_3)}-
\frac{1}{2}\int\limits_t^T(t-s)^{l_3}
\int\limits_t^s(t-s_1)^{2l}ds_1d{\bf w}_s^{(i_3)}\right)^2\right\}=
$$
$$
=\hbox{\vtop{\offinterlineskip\halign{
\hfil#\hfil\cr
{\rm lim}\cr
$\stackrel{}{{}_{p_1\to \infty}}$\cr
}} }\sum\limits_{j_3=0}^{2l+l_3+1}
\left(\sum\limits_{j_1=0}^{p_1}C_{j_3j_1 j_1}-
\frac{1}{2}\tilde C_{j_3}\right)^2+
$$
\begin{equation}
\label{ogo210}
+
\hbox{\vtop{\offinterlineskip\halign{
\hfil#\hfil\cr
{\rm lim}\cr
$\stackrel{}{{}_{p_1,p_3\to \infty}}$\cr
}} }{\sf M}\left\{\left(
\sum\limits_{j_3=2l+l_3+2}^{p_3}
\sum\limits_{j_1=0}^{p_1}
C_{j_3 j_1 j_1}\zeta_{j_3}^{(i_3)}\right)^2\right\}.
\end{equation}

Let us prove that
\begin{equation}
\label{ogo211}
\hbox{\vtop{\offinterlineskip\halign{
\hfil#\hfil\cr
{\rm lim}\cr
$\stackrel{}{{}_{p_1\to \infty}}$\cr
}} }
\left(\sum\limits_{j_1=0}^{p_1}C_{j_3j_1 j_1}-
\frac{1}{2}\tilde C_{j_3}\right)^2=0.
\end{equation}

We have
$$
\left(\sum\limits_{j_1=0}^{p_1}C_{j_3j_1 j_1}-
\frac{1}{2}\tilde C_{j_3}\right)^2=
$$
$$
=\left(\frac{1}{2}\sum\limits_{j_1=0}^{p_1}
\int\limits_t^T\phi_{j_3}(s)(t-s)^{l_3}
\left(\int\limits_t^s\phi_{j_1}(s_1)(t-s_1)^{l}ds_1\right)^2ds-\right.
$$
$$
\left.-
\frac{1}{2}
\int\limits_t^T
\phi_{j_3}(s)(t-s)^{l_3}\int\limits_t^s(t-s_1)^{2l}ds_1ds\right)^2=
$$
$$
=\frac{1}{4}\left(
\int\limits_t^T\phi_{j_3}(s)(t-s)^{l_3}\left(
\sum\limits_{j_1=0}^{p_1}
\left(\int\limits_t^s\phi_{j_1}(s_1)(t-s_1)^{l}ds_1\right)^2
-\right.\right.
$$
$$
\left.\left.-
\int\limits_t^s(t-s_1)^{2l}ds_1\right)ds\right)^2=
$$
$$
=\frac{1}{4}\left(
\int\limits_t^T\phi_{j_3}(s)(t-s)^{l_3}\left(
\int\limits_t^s(t-s_1)^{2l}ds_1-\sum\limits_{j_1=p_1+1}^{\infty}
\left(\int\limits_t^s\phi_{j_1}(s_1)(t-s_1)^{l}ds_1\right)^2
-\right.\right.
$$
$$
-\left.\left.
\int\limits_t^s(t-s_1)^{2l}ds_1\right)ds\right)^2=
$$
\begin{equation}
\label{ogo300}
~~~~~=\frac{1}{4}\left(
\int\limits_t^T\phi_{j_3}(s)(t-s)^{l_3}
\sum\limits_{j_1=p_1+1}^{\infty}
\left(\int\limits_t^s\phi_{j_1}(s_1)(t-s_1)^{l}ds_1\right)^2
ds\right)^2.
\end{equation}

In order to get (\ref{ogo300}) we used the Parseval equality
\begin{equation}
\label{ogo301}
\sum_{j_1=0}^{\infty}\left(\int\limits_t^s\phi_{j_1}(s_1)
(t-s_1)^lds_1\right)^2=
\int\limits_t^T K^2(s,s_1)ds_1,
\end{equation}
where
$$
K(s,s_1)=(t-s_1)^l\ {\bf 1}_{\{s_1<s\}},\ \ \ s, s_1\in [t, T].
$$

Taking into account the nondecreasing
of the functional sequence
$$
u_n(s)=\sum_{j_1=0}^{n}\left(\int\limits_t^s\phi_{j_1}(s_1)
(t-s_1)^lds_1\right)^2,
$$
continuity of its members and continuity of the limit function
$$
u(s)=\int\limits_t^s(t-s_1)^{2l}ds_1
$$ 
at the interval $[t, T]$
in accordance with the Dini Theorem we
have uniform
convergence of the functional sequences $u_n(s)$ to the limit function
$u(s)$ at the interval $[t, T]$.

From (\ref{ogo300}) using the inequality 
of Cauchy--Bunyakovsky, we obtain
$$
\left(\sum\limits_{j_1=0}^{p_1}C_{j_3j_1 j_1}-
\frac{1}{2}\tilde C_{j_3}\right)^2\le
$$
$$
\le
\frac{1}{4}
\int\limits_t^T\phi_{j_3}^2(s)(t-s)^{2l_3}ds
\int\limits_t^T\left(\sum\limits_{j_1=p_1+1}^{\infty}
\left(\int\limits_t^s\phi_{j_1}(s_1)(t-s_1)^{l}ds_1\right)^2
\right)^2 ds\le
$$
\begin{equation}
\label{ogo302}
~~~~~~~~\le\frac{1}{4}\varepsilon^2 (T-t)^{2l_3}\int\limits_t^T\phi_{j_3}^2(s)ds
(T-t)=\frac{1}{4}(T-t)^{2l_3+1}\varepsilon^2
\end{equation}

\noindent
when $p_1>N(\varepsilon),$ where $N(\varepsilon)\in{\bf N}$
exists for any $\varepsilon>0.$
The relation (\ref{ogo302}) implies
(\ref{ogo211}).

Further, 
\begin{equation}
\label{ogo303}
\sum\limits_{j_1=0}^{p_1}
\sum\limits_{j_3=2l+l_3+2}^{p_3}
C_{j_3 j_1 j_1}\zeta_{j_3}^{(i_3)}=
\sum\limits_{j_1=0}^{p_1}
\sum\limits_{j_3=2l+l_3+2}^{2(j_1+l+1)+l_3}
C_{j_3 j_1 j_1}\zeta_{j_3}^{(i_3)}.
\end{equation}

We put  $2(j_1+l+1)+l_3$ instead of $p_3$, since
$C_{j_3j_1j_1}=0$ for $j_3>2(j_1+l+1)+l_3.$ This conclusion
follows from the relation
$$
C_{j_3j_1j_1}=
\frac{1}{2}
\int\limits_t^T\phi_{j_3}(s)(t-s)^{l_3}
\left(
\int\limits_t^s\phi_{j_1}(s_1)(t-s_1)^{l}ds_1\right)^2ds=
$$
$$
=
\frac{1}{2}\int\limits_t^T\phi_{j_3}(s)Q_{2(j_1+l+1)+l_3}(s)ds,
$$
where $Q_{2(j_1+l+1)+l_3}(s)$ is a polynomial of degree
$2(j_1+l+1)+l_3.$

It is easy to see that
\begin{equation}
\label{ogo304}
~~~~~~~~~~~~ \sum\limits_{j_1=0}^{p_1}
\sum\limits_{j_3=2l+l_3+2}^{2(j_1+l+1)+l_3}
C_{j_3 j_1 j_1}\zeta_{j_3}^{(i_3)}=
\sum\limits_{j_3=2l+l_3+2}^{2(p_1+l+1)+l_3}
\sum\limits_{j_1=0}^{p_1}
C_{j_3 j_1 j_1}\zeta_{j_3}^{(i_3)}.
\end{equation}

Note that we included some zero coefficients $C_{j_3 j_1 j_1}$ 
into the sum $\sum\limits_{j_1=0}^{p_1}$.
From (\ref{ogo303}) and (\ref{ogo304}) we have
$$
{\sf M}\left\{\left(\sum\limits_{j_1=0}^{p_1}
\sum\limits_{j_3=2l+l_3+2}^{p_3}
C_{j_3 j_1 j_1}\zeta_{j_3}^{(i_3)}\right)^2\right\}=
$$
$$
=
{\sf M}\left\{\left(
\sum\limits_{j_3=2l+l_3+2}^{2(p_1+l+1)+l_3}
\sum\limits_{j_1=0}^{p_1}
C_{j_3 j_1 j_1}\zeta_{j_3}^{(i_3)}\right)^2\right\}
=\sum\limits_{j_3=2l+l_3+2}^{2(p_1+l+1)+l_3}
\left(\sum\limits_{j_1=0}^{p_1}
C_{j_3 j_1 j_1}\right)^2=
$$
$$
=
\sum\limits_{j_3=2l+l_3+2}^{2(p_1+l+1)+l_3}
\left(\frac{1}{2}\sum\limits_{j_1=0}^{p_1}
\int\limits_t^T\phi_{j_3}(s)(t-s)^{l_3}
\left(\int\limits_t^s\phi_{j_1}(s_1)(t-s_1)^{l}ds_1\right)^2ds\right)^2=
$$
$$
=\frac{1}{4}\sum\limits_{j_3=2l+l_3+2}^{2(p_1+l+1)+l_3}
\left(
\int\limits_t^T\phi_{j_3}(s)(t-s)^{l_3}
\sum\limits_{j_1=0}^{p_1}
\left(\int\limits_t^s\phi_{j_1}(s_1)(t-s_1)^{l}ds_1\right)^2
ds\right)^2=
$$
$$
=\frac{1}{4}\sum\limits_{j_3=2l+l_3+2}^{2(p_1+l+1)+l_3}\left(
\int\limits_t^T\phi_{j_3}(s)(t-s)^{l_3}\left(
\int\limits_t^s(t-s_1)^{2l}ds_1-\right.\right.
$$
$$
\left.\left.
-\sum\limits_{j_1=p_1+1}^{\infty}
\left(\int\limits_t^s\phi_{j_1}(s_1)(t-s_1)^{l}ds_1\right)^2
\right)ds\right)^2=
$$
\begin{equation}
\label{ogo310}
=\frac{1}{4}\sum\limits_{j_3=2l+l_3+2}^{2(p_1+l+1)+l_3}\left(
\int\limits_t^T\phi_{j_3}(s)(t-s)^{l_3}
\sum\limits_{j_1=p_1+1}^{\infty}
\left(\int\limits_t^s\phi_{j_1}(s_1)(t-s_1)^{l}ds_1\right)^2
ds\right)^2.
\end{equation}

In 
order to 
get
(\ref{ogo310}) we used the Parseval equality 
(\ref{ogo301}) and the following 
relation
$$
\int\limits_t^T\phi_{j_3}(s)Q_{2l+1+l_3}(s)ds=0,\ \ \ j_3>2l+1+l_3,
$$
where $Q_{2l+1+l_3}(s)$ is a polynomial of degree
$2l+1+l_3.$

Further, we have
$$
\left(\int\limits_t^s\phi_{j_1}(s_1)(t-s_1)^lds_1\right)^2=
$$
$$
=
\frac{(T-t)^{2l+1}(2j_1+1)}{2^{2l+2}}
\left(\int\limits_{-1}^{z(s)}
P_{j_1}(y)(1+y)^ldy\right)^2=
$$
$$
=\frac{(T-t)^{2l+1}}{2^{2l+2}(2j_1+1)}\times
$$
$$
\times\left(
\left(1+z(s)\right)^l
R_{j_1}(s)
-l\int\limits_{-1}^{z(s)}
\left(P_{j_1+1}(y)-P_{j_1-1}(y)\right)\left(1+y\right)^{l-1}dy\right)^2
\le
$$
$$
\le\frac{(T-t)^{2l+1}2}{2^{2l+2}(2j_1+1)}\times
$$
$$
\times\left(
\left(\frac{2(s-t)}{T-t}\right)^{2l}
R_{j_1}^2(s)
+l^2
\left(
\int\limits_{-1}^{z(s)}
\left(P_{j_1+1}(y)-P_{j_1-1}(y)\right)\left(1+y\right)^{l-1}dy\right)^2
\right)\le
$$
$$
\le\frac{(T-t)^{2l+1}}{2^{2l+1}(2j_1+1)}\times
$$
$$
\times\left(
2^{2l+1}
Z_{j_1}(s)+
l^2
\int\limits_{-1}^{z(s)}
(1+y)^{2l-2}dy
\int\limits_{-1}^{z(s)}
\left(P_{j_1+1}(y)-P_{j_1-1}(y)\right)^2dy
\right)\le
$$
$$
\le\frac{(T-t)^{2l+1}}{2^{2l+1}(2j_1+1)}\times
$$
$$
\times\left(
2^{2l+1}
Z_{j_1}(s)
+\frac{2l^2}{2l-1}\left(\frac{2(s-t)}{T-t}\right)^{2l-1}
\int\limits_{-1}^{z(s)}
\left(P_{j_1+1}^2(y)+P_{j_1-1}^2(y)\right)dy
\right)\le
$$
\begin{equation}
\label{ogo400}
~~~~~\le\frac{(T-t)^{2l+1}}{2(2j_1+1)}\left(
2
Z_{j_1}(s)
+\frac{l^2}{2l-1}
\int\limits_{-1}^{z(s)}
\left(P_{j_1+1}^2(y)+P_{j_1-1}^2(y)\right)dy
\right),
\end{equation}

\noindent
where $j_1\in{\bf N},$
$$
R_{j_1}(s)=P_{j_1+1}(z(s))-P_{j_1-1}(z(s)),
$$
$$
Z_{j_1}(s)=P_{j_1+1}^2(z(s))+P_{j_1-1}^2(z(s)).
$$

\vspace{2mm}

Let us estimate the right-hand side 
of (\ref{ogo400}) using (\ref{ogo23}) $(j_1\in{\bf N})$
$$
\left(\int\limits_t^s\phi_{j_1}(s_1)(t-s_1)^lds_1\right)^2 <
$$
$$
<
\frac{(T-t)^{2l+1}}{2(2j_1+1)}\left(\frac{K^2}{j_1+2}+\frac{K^2}{j_1}
\right)\left(\frac{2}
{(1-\left(z(s))^2
\right)^{1/2}}
+\frac{l^2}{2l-1}
\int\limits_{-1}^{z(s)}
\frac{dy}{\left(1-y^2\right)^{1/2}}\right)<
$$
\begin{equation}
\label{ogo401}
~~~~~~~<\frac{(T-t)^{2l+1}K^2}{2j_1^2}\left(
\frac{2}
{(1-\left(z(s))^2
\right)^{1/2}}+
\frac{l^2\pi}{2l-1}\right),\ \ \ s\in(t, T).
\end{equation}

\vspace{2mm}

From (\ref{ogo310}) and (\ref{ogo401}) we obtain
$$
{\sf M}\left\{\left(\sum\limits_{j_1=0}^{p_1}
\sum\limits_{j_3=2l+l_3+2}^{p_3}
C_{j_3 j_1 j_1}\zeta_{j_3}^{(i_3)}\right)^2\right\}\le
$$
$$
\le
\frac{1}{4}\sum\limits_{j_3=2l+l_3+2}^{2(p_1+l+1)+l_3}\left(
\int\limits_t^T|\phi_{j_3}(s)|(t-s)^{l_3}
\sum\limits_{j_1=p_1+1}^{\infty}
\left(\int\limits_t^s\phi_{j_1}(s_1)(t-s_1)^{l}ds_1\right)^2
ds\right)^2\le
$$
$$
\le
\frac{1}{4}(T-t)^{2l_3}\sum\limits_{j_3=2l+l_3+2}^{2(p_1+l+1)+l_3}\left(
\int\limits_t^T|\phi_{j_3}(s)|
\sum\limits_{j_1=p_1+1}^{\infty}
\left(\int\limits_t^s\phi_{j_1}(s_1)(t-s_1)^{l}ds_1\right)^2
ds\right)^2<
$$
$$
<
\frac{(T-t)^{4l+2l_3+1}K^4 K_1^2}{16}\ \
\sum\limits_{j_3=2l+l_3+2}^{2(p_1+l+1)+l_3}
\left(\left(\int\limits_t^T
\frac{2ds}
{\left(1-\left(z(s)\right)^2
\right)^{3/4}}
+\right.\right.
$$
$$
\left.\left.
+\frac{l^2\pi}{2l-1}
\int\limits_t^T
\frac{ds}
{\left(1-\left(z(s)\right)^2
\right)^{1/4}}\right)
\sum\limits_{j_1=p_1+1}^{\infty}\frac{1}{j_1^2}
\right)^2\le 
$$
$$
\le
\frac{(T-t)^{4l+2l_3+3}K^4 K_1^2}{64}\ \ \frac{2p_1+1}{p_1^2}
\left(\int\limits_{-1}^1
\frac{2dy}
{(1-y^2)^{3/4}}
+\frac{l^2\pi}{2l-1}
\int\limits_{-1}^1
\frac{dy}
{(1-y^2)^{1/4}}\right)^2\le
$$

\begin{equation}
\label{ogo500}
~~~~~~\le C(T-t)^{4l+2l_3+3}\ \frac{2p_1+1}{p_1^2}\to 0\ \ \ 
\hbox{when}\ \ p_1\ \to\ \infty,
\end{equation}

\vspace{2mm}
\noindent
where constant $C$ does not depend on $p_1$ and $T-t$.

The relations (\ref{ogo210}), (\ref{ogo211}), and (\ref{ogo500}) imply 
(\ref{ogo200}), and the relation (\ref{ogo200}) 
implies the correctness of the formula
(\ref{ogo101}).

Let us consider Case 3, i.e. $i_2=i_3\ne i_1$, $l_2=l_3=l\ne l_1$, and
$l_1, l_3=0, 1, 2,\ldots$ So, we prove 
the following expansion 
\begin{equation}
\label{ogo101ee}
I_{{l_1 l_3 l_3}_{T,t}}^{*(i_1i_3i_3)}=
\hbox{\vtop{\offinterlineskip\halign{
\hfil#\hfil\cr
{\rm l.i.m.}\cr
$\stackrel{}{{}_{p_1,p_2,p_3\to \infty}}$\cr
}} }\sum_{j_1=0}^{p_1}\sum_{j_2=0}^{p_2}\sum_{j_3=0}^{p_3}
C_{j_3 j_2 j_1}\zeta_{j_1}^{(i_1)}\zeta_{j_2}^{(i_3)}\zeta_{j_3}^{(i_3)}\ \ \
(i_1, i_2, i_3=1,\ldots,m),
\end{equation}

\noindent
where $l_1, l_3=0, 1, 2,\ldots$ $(l_3=l)$ and
\begin{equation}
\label{ogo1991}
C_{j_3 j_2 j_1}=\int\limits_t^T
\phi_{j_3}(s)(t-s)^{l}\int\limits_t^s(t-s_1)^{l}
\phi_{j_2}(s_1)
\int\limits_t^{s_1}(t-s_2)^{l_1}
\phi_{j_1}(s_2)ds_2ds_1ds.
\end{equation}

If we prove w.~p.~1 the formula
\begin{equation}
\label{ogo2000}
~~~\hbox{\vtop{\offinterlineskip\halign{
\hfil#\hfil\cr
{\rm l.i.m.}\cr
$\stackrel{}{{}_{p_1, p_3\to \infty}}$\cr
}} }
\sum\limits_{j_1=0}^{p_1}\sum\limits_{j_3=0}^{p_3}
C_{j_3 j_3 j_1}\zeta_{j_1}^{(i_1)}=
\frac{1}{2}\int\limits_t^T(t-s)^{2l}
\int\limits_t^s(t-s_1)^{l_1}d{\bf w}_{s_1}^{(i_1)}ds,
\end{equation}
where the 
coefficients $C_{j_3 j_3 j_1}$ are defined by (\ref{ogo1991}), 
then using Theorem 1.1 
and standard relations between iterated It\^{o} and 
Stratonovich stochastic integrals, we obtain the expansion (\ref{ogo101ee}).

Using Theorem 1.1 and the It\^{o} formula, we 
have
$$
\frac{1}{2}\int\limits_t^T(t-s)^{2l}
\int\limits_t^s(t-s_1)^{l_1}d{\bf w}_{s_1}^{(i_1)}ds=
\frac{1}{2}\int\limits_t^T(t-s_1)^{l_1}
\int\limits_{s_1}^T(t-s)^{2l}dsd{\bf w}_{s_1}^{(i_1)}=
$$
$$
=
\frac{1}{2}\sum\limits_{j_1=0}^{2l+l_1+1}
\tilde C_{j_1}\zeta_{j_1}^{(i_1)}\ \ \ \hbox{w.~p.~1},
$$
where
$$
\tilde C_{j_1}=
\int\limits_t^T
\phi_{j_1}(s_1)(t-s_1)^{l_1}\int\limits_{s_1}^T(t-s)^{2l}dsds_1.
$$

Then
$$
\sum\limits_{j_1=0}^{p_1}\sum\limits_{j_3=0}^{p_3}
C_{j_3 j_3 j_1}\zeta_{j_1}^{(i_1)}-
\frac{1}{2}\sum\limits_{j_1=0}^{2l+l_1+1}
\tilde C_{j_1}\zeta_{j_1}^{(i_1)}=
$$
$$
=
\sum\limits_{j_1=0}^{2l+l_1+1}
\left(\sum\limits_{j_3=0}^{p_3}
C_{j_3 j_3 j_1}-\frac{1}{2}\tilde C_{j_1}\right)
\zeta_{j_1}^{(i_1)}+
\sum\limits_{j_1=2l+l_1+2}^{p_1}
\sum\limits_{j_3=0}^{p_3}
C_{j_3 j_3 j_1}\zeta_{j_1}^{(i_1)}.
$$

\vspace{1mm}
\noindent
\par
Therefore,
$$
\hbox{\vtop{\offinterlineskip\halign{
\hfil#\hfil\cr
{\rm lim}\cr
$\stackrel{}{{}_{p_1,p_3\to \infty}}$\cr
}} }
{\sf M}\left\{\left(
\sum\limits_{j_1=0}^{p_1}\sum\limits_{j_3=0}^{p_3}C_{j_3j_3 j_1}
\zeta_{j_1}^{(i_1)}-
\frac{1}{2}\int\limits_t^T(t-s)^{2l}
\int\limits_t^s(t-s_1)^{l_1}d{\bf w}_{s_1}^{(i_1)}ds\right)^2\right\}=
$$
$$
=\hbox{\vtop{\offinterlineskip\halign{
\hfil#\hfil\cr
{\rm lim}\cr
$\stackrel{}{{}_{p_3\to \infty}}$\cr
}} }\sum\limits_{j_1=0}^{2l+l_1+1}
\left(\sum\limits_{j_3=0}^{p_3}C_{j_3j_3 j_1}-
\frac{1}{2}\tilde C_{j_1}\right)^2+
$$
\begin{equation}
\label{ogo2100}
+
\hbox{\vtop{\offinterlineskip\halign{
\hfil#\hfil\cr
{\rm lim}\cr
$\stackrel{}{{}_{p_1,p_3\to \infty}}$\cr
}} }{\sf M}\left\{\left(
\sum\limits_{j_1=2l+l_1+2}^{p_1}
\sum\limits_{j_3=0}^{p_3}
C_{j_3 j_3 j_1}\zeta_{j_1}^{(i_1)}\right)^2\right\}.
\end{equation}

Let us prove that
\begin{equation}
\label{ogo2110}
\hbox{\vtop{\offinterlineskip\halign{
\hfil#\hfil\cr
{\rm lim}\cr
$\stackrel{}{{}_{p_3\to \infty}}$\cr
}} }
\left(\sum\limits_{j_3=0}^{p_3}C_{j_3j_3 j_1}-
\frac{1}{2}\tilde C_{j_1}\right)^2=0.
\end{equation}

We have
$$
\left(\sum\limits_{j_3=0}^{p_3}C_{j_3j_3 j_1}-
\frac{1}{2}\tilde C_{j_1}\right)^2=
$$
$$
=
\left(\sum\limits_{j_3=0}^{p_3}
\int\limits_t^T\phi_{j_1}(s_2)(t-s_2)^{l_1}ds_2
\int\limits_{s_2}^T\phi_{j_3}(s_1)(t-s_1)^{l}ds_1
\int\limits_{s_1}^T\phi_{j_3}(s)(t-s)^{l}ds-\right.
$$
$$
-\left.
\frac{1}{2}
\int\limits_t^T
\phi_{j_1}(s_1)(t-s_1)^{l_1}\int\limits_{s_1}^T(t-s)^{2l}dsds_1\right)^2=
$$
$$
=\left(\frac{1}{2}\sum\limits_{j_3=0}^{p_3}
\int\limits_t^T\phi_{j_1}(s_2)(t-s_2)^{l_1}
\left(\int\limits_{s_2}^T\phi_{j_3}(s_1)(t-s_1)^{l}ds_1\right)^2ds_2-
\right.
$$
$$
\left.-
\frac{1}{2}
\int\limits_t^T
\phi_{j_1}(s_1)(t-s_1)^{l_1}\int\limits_{s_1}^T(t-s)^{2l}dsds_1\right)^2=
$$
$$
=\frac{1}{4}\left(
\int\limits_t^T\phi_{j_1}(s_1)(t-s_1)^{l_1}\left(
\sum\limits_{j_3=0}^{p_3}
\left(\int\limits_{s_1}^T\phi_{j_3}(s)(t-s)^{l}ds\right)^2
-\right.\right.
$$
$$
\left.\left.
-\int\limits_{s_1}^T(t-s)^{2l}ds\right)ds_1\right)^2=
$$
$$
=\frac{1}{4}\left(
\int\limits_t^T\phi_{j_1}(s_1)(t-s_1)^{l_1}\left(
\int\limits_{s_1}^T(t-s)^{2l}ds-
\sum\limits_{j_3=p_3+1}^{\infty}
\left(\int\limits_{s_1}^T\phi_{j_3}(s)(t-s)^{l}ds\right)^2
-\right.\right.
$$
$$
\left.\left.
-\int\limits_{s_1}^T(t-s)^{2l}ds\right)ds_1\right)^2=
$$
\begin{equation}
\label{ogo3000}
~~~~~=\frac{1}{4}\left(
\int\limits_t^T\phi_{j_1}(s_1)(t-s_1)^{l_1}
\sum\limits_{j_3=p_3+1}^{\infty}
\left(\int\limits_{s_1}^T\phi_{j_3}(s)(t-s)^{l}ds\right)^2
ds_1\right)^2.
\end{equation}

In order to 
get (\ref{ogo3000}) we used the Parseval equality
\begin{equation}
\label{ogo3010}
\sum_{j_3=0}^{\infty}\left(\int\limits_{s_1}^T\phi_{j_3}(s)
(t-s)^lds\right)^2=
\int\limits_t^T K^2(s,s_1)ds,
\end{equation}
where
$$
K(s,s_1)=(t-s)^l\
{\bf 1}_{\{s_1<s\}},\ \ \ s, s_1\in [t, T].
$$

Taking into account the nondecreasing
of the functional sequence
$$
u_n(s_1)=\sum_{j_3=0}^{n}\left(\int\limits_{s_1}^T\phi_{j_3}(s)
(t-s)^lds\right)^2,
$$
continuity of its members and continuity of the limit
function 
$$
u(s_1)=\int\limits_{s_1}^T(t-s)^{2l}ds
$$ 
at the interval $[t, T]$ in accordance with the Dini Theorem we have 
uniform 
convergence
of the functional sequence $u_n(s_1)$ to the limit 
function $u(s_1)$ at the interval $[t, T]$.

From (\ref{ogo3000}) using the inequality of Cauchy--Bunyakovsky, we obtain
$$
\left(\sum\limits_{j_3=0}^{p_3}C_{j_3j_3 j_1}-
\frac{1}{2}\tilde C_{j_1}\right)^2\le
$$
$$
\le
\frac{1}{4}
\int\limits_t^T\phi_{j_1}^2(s_1)(t-s_1)^{2l_1}ds_1
\int\limits_t^T\left(\sum\limits_{j_3=p_3+1}^{\infty}
\left(\int\limits_{s_1}^T\phi_{j_3}(s)(t-s)^{l}ds\right)^2
\right)^2 ds_1\le
$$
\begin{equation}
\label{ogo3020}
~~~~~~~~\le\frac{1}{4}\varepsilon^2 (T-t)^{2l_1}\int\limits_t^T\phi_{j_1}^2(s_1)ds_1
(T-t)=\frac{1}{4}(T-t)^{2l_1+1}\varepsilon^2
\end{equation}

\noindent 
when $p_3>N(\varepsilon),$ where $N(\varepsilon)\in{\bf N}$
exists for any $\varepsilon>0.$
The relation (\ref{ogo2110}) follows from (\ref{ogo3020}).

We have
\begin{equation}
\label{ogo3030}
\sum\limits_{j_3=0}^{p_3}
\sum\limits_{j_1=2l+l_1+2}^{p_1}
C_{j_3 j_3 j_1}\zeta_{j_1}^{(i_1)}
=
\sum\limits_{j_3=0}^{p_3}
\sum\limits_{j_1=2l+l_1+2}^{2(j_3+l+1)+l_1}
C_{j_3 j_3 j_1}\zeta_{j_1}^{(i_1)}.
\end{equation}

We put  $2(j_3+l+1)+l_1$ instead of $p_1$, since
$C_{j_3j_3j_1}=0$ when $j_1>2(j_3+l+1)+l_1.$ 
This conclusion follows from the relation
$$
C_{j_3j_3j_1}=
\frac{1}{2}
\int\limits_t^T\phi_{j_1}(s_2)(t-s_2)^{l_1}
\left(
\int\limits_{s_2}^T\phi_{j_3}(s_1)(t-s_1)^{l}ds_1\right)^2ds_2=
$$
$$
=
\frac{1}{2}\int\limits_t^T\phi_{j_1}(s_2)Q_{2(j_3+l+1)+l_1}(s_2)ds_2,
$$
where $Q_{2(j_3+l+1)+l_1}(s)$ is a polynomial of degree
$2(j_3+l+1)+l_1.$

It is easy to see that
\begin{equation}
\label{ogo3040}
~~~~~~~~~~~~ \sum\limits_{j_3=0}^{p_3}
\sum\limits_{j_1=2l+l_1+2}^{2(j_3+l+1)+l_1}
C_{j_3 j_3 j_1}\zeta_{j_1}^{(i_1)}=
\sum\limits_{j_1=2l+l_1+2}^{2(p_3+l+1)+l_1}
\sum\limits_{j_3=0}^{p_3}
C_{j_3 j_3 j_1}\zeta_{j_1}^{(i_1)}.
\end{equation}

Note that we 
included
some zero coefficients $C_{j_3 j_3 j_1}$ 
into the sum 
$\sum\limits_{j_3=0}^{p_3}$. 

From (\ref{ogo3030}) and (\ref{ogo3040}) we have
$$
{\sf M}\left\{\left(\sum\limits_{j_3=0}^{p_3}
\sum\limits_{j_1=2l+l_1+2}^{p_1}
C_{j_3 j_3 j_1}\zeta_{j_1}^{(i_1)}\right)^2\right\}=
$$
$$
=
{\sf M}\left\{\left(
\sum\limits_{j_1=2l+l_1+2}^{2(p_3+l+1)+l_1}
\sum\limits_{j_3=0}^{p_3}
C_{j_3 j_3 j_1}\zeta_{j_1}^{(i_1)}\right)^2\right\}
=\sum\limits_{j_1=2l+l_1+2}^{2(p_3+l+1)+l_1}
\left(\sum\limits_{j_3=0}^{p_3}
C_{j_3 j_3 j_1}\right)^2=
$$
$$
=\sum\limits_{j_1=2l+l_1+2}^{2(p_3+l+1)+l_1}
\left(\frac{1}{2}\sum\limits_{j_3=0}^{p_3}
\int\limits_t^T\phi_{j_1}(s_2)(t-s_2)^{l_1}
\left(\int\limits_{s_2}^T\phi_{j_3}(s_1)(t-s_1)^{l}ds_1\right)^2ds_2\right)^2=
$$
$$
=\frac{1}{4}\sum\limits_{j_1=2l+l_1+2}^{2(p_3+l+1)+l_1}
\left(
\int\limits_t^T\phi_{j_1}(s_2)(t-s_2)^{l_1}
\sum\limits_{j_3=0}^{p_3}
\left(\int\limits_{s_2}^T\phi_{j_3}(s_1)(t-s_1)^{l}ds_1\right)^2
ds_2\right)^2=
$$
$$
=\frac{1}{4}\sum\limits_{j_1=2l+l_1+2}^{2(p_3+l+1)+l_1}\left(
\int\limits_t^T\phi_{j_1}(s_2)(t-s_2)^{l_1}\left(
\int\limits_{s_2}^T(t-s_1)^{2l}ds_1-\right.\right.
$$
$$
\left.\left.-
\sum\limits_{j_3=p_3+1}^{\infty}
\left(\int\limits_{s_2}^T\phi_{j_3}(s_1)(t-s_1)^{l}ds_1\right)^2
\right)ds_2\right)^2=
$$
\begin{equation}
\label{ogo3100}
=\frac{1}{4}\sum\limits_{j_1=2l+l_1+2}^{2(p_3+l+1)+l_1}\left(
\int\limits_t^T\phi_{j_1}(s_2)(t-s_2)^{l_1}
\sum\limits_{j_3=p_3+1}^{\infty}
\left(\int\limits_{s_2}^T\phi_{j_3}(s_1)(t-s_1)^{l}ds_1\right)^2
ds_2\right)^2.
\end{equation}

In order to 
get
(\ref{ogo3100}) we used the Parseval equality 
(\ref{ogo3010}) and the following relation
$$
\int\limits_t^T\phi_{j_1}(s)Q_{2l+1+l_1}(s)ds=0,\ \ \ j_1>2l+1+l_1,
$$
where $Q_{2l+1+l_1}(s)$ is a polynomial of degree
$2l+1+l_1.$

Further, we have
$$
\left(\int\limits_{s_2}^T\phi_{j_3}(s_1)(t-s_1)^lds_1\right)^2=
$$
$$
=
\frac{(T-t)^{2l+1}(2j_3+1)}{2^{2l+2}}
\left(\int\limits_{z(s_2)}^1
P_{j_3}(y)(1+y)^ldy\right)^2=
$$
$$
=\frac{(T-t)^{2l+1}}{2^{2l+2}(2j_3+1)}\times
$$
$$
\times\left(
\left(1+z(s_2)\right)^l
Q_{j_3}(s_2)-
l\int\limits_{z(s_2)}^1
\left(P_{j_3+1}(y)-P_{j_3-1}(y)\right)\left(1+y\right)^{l-1}dy\right)^2\le
$$
$$
\le\frac{(T-t)^{2l+1}2}{2^{2l+2}(2j_3+1)}\times
$$
$$
\times
\left(
\left(\frac{2(s_2-t)}{T-t}\right)^{2l}
Q_{j_3}^2(s_2)+
l^2
\left(
\int\limits_{z(s_2)}^1
\left(P_{j_3+1}(y)-P_{j_3-1}(y)\right)\left(1+y\right)^{l-1}dy\right)^2
\right)\le
$$
$$
\le\frac{(T-t)^{2l+1}}{2^{2l+1}(2j_3+1)}\times
$$
\newpage
\noindent
$$
\times\left(
2^{2l+1}
H_{j_3}(s_2)+
l^2
\int\limits_{z(s_2)}^1
(1+y)^{2l-2}dy
\int\limits_{z(s_2)}^1
\left(P_{j_3+1}(y)-P_{j_3-1}(y)\right)^2dy
\right)\le
$$
$$
\le\frac{(T-t)^{2l+1}}{2^{2l+1}(2j_3+1)}\times
$$
$$
\times\hspace{-1mm}\left(
2^{2l+1}
H_{j_3}(s_2)
+\frac{2^{2l}l^2}{2l-1}\left(1-\left(\frac{s_2-t}{T-t}\right)^{2l-1}\right)
\int\limits_{z(s_2)}^1
\left(P_{j_3+1}^2(y)+P_{j_3-1}^2(y)\right)dy
\right)\hspace{-1mm}\le
$$
\begin{equation}
\label{ogo4000}
~~~~~\le\frac{(T-t)^{2l+1}}{2(2j_3+1)}\left(
2
H_{j_3}(s_2)
\Biggl.
+\frac{l^2}{2l-1}
\int\limits_{z(s_2)}^1
\left(P_{j_3+1}^2(y)+P_{j_3-1}^2(y)\right)dy
\right),
\end{equation}

\noindent
where $j_3\in{\bf N},$
$$
Q_{j_3}(s_2)=P_{j_3-1}(z(s_2))-P_{j_3+1}(z(s_2)),
$$
$$
H_{j_3}(s_2)=P_{j_3-1}^2(z(s_2))+P_{j_3+1}^2(z(s_2)).
$$

\vspace{2mm}

Let us estimate the right-hand side
of (\ref{ogo4000}) using
(\ref{ogo23}) $(j_3\in{\bf N})$
$$
\left(\int\limits_{s_2}^T\phi_{j_3}(s_1)(t-s_1)^lds_1\right)^2 <
$$
$$
<
\frac{(T-t)^{2l+1}}{2(2j_3+1)}\left(\frac{K^2}{j_3+2}+\frac{K^2}{j_3}
\right)\left(\frac{2}
{(1-\left(z(s_2))^2
\right)^{1/2}}
+\frac{l^2}{2l-1}
\int\limits_{z(s_2)}^1
\frac{dy}{\left(1-y^2\right)^{1/2}}\right)<
$$
\begin{equation}
\label{ogo4010}
~~~~~~~<\frac{(T-t)^{2l+1}K^2}{2j_3^2}\left(
\frac{2}
{(1-\left(z(s_2))^2
\right)^{1/2}}+
\frac{l^2\pi}{2l-1}\right),\ \ \ s_2\in(t, T).
\end{equation}

\vspace{3mm}

From (\ref{ogo3100}) and (\ref{ogo4010}) we obtain
$$
{\rm M}\left\{\left(\sum\limits_{j_3=0}^{p_3}
\sum\limits_{j_1=2l+l_1+2}^{p_1}
C_{j_3 j_3 j_1}\zeta_{j_1}^{(i_1)}\right)^2\right\}\le
$$
$$
\le
\frac{1}{4}\sum\limits_{j_1=2l+l_1+2}^{2(p_3+l+1)+l_1}\left(
\int\limits_t^T|\phi_{j_1}(s_2)|(t-s_2)^{l_1}
\sum\limits_{j_3=p_3+1}^{\infty}
\left(\int\limits_{s_2}^T\phi_{j_3}(s_1)(t-s_1)^{l}ds_1\right)^2
ds_2\right)^2\le
$$
$$
\le
\frac{1}{4}(T-t)^{2l_1}\sum\limits_{j_1=2l+l_1+2}^{2(p_3+l+1)+l_1}\left(
\int\limits_t^T|\phi_{j_1}(s_2)|
\sum\limits_{j_3=p_3+1}^{\infty}
\left(\int\limits_{s_2}^T\phi_{j_3}(s_1)(t-s_1)^{l}ds_1\right)^2
ds_2\right)^2<
$$
$$
<
\frac{(T-t)^{4l+2l_1+1}K^4 K_1^2}{16}\ \
\sum\limits_{j_1=2l+l_1+2}^{2(p_3+l+1)+l_1}
\left(\left(\int\limits_t^T
\frac{2ds_2}
{(1-\left(z(s_2))^2
\right)^{3/4}}+\right.\right.
$$
$$
\left.\left.
+\frac{l^2\pi}{2l-1}
\int\limits_t^T
\frac{ds_2}
{(1-\left(z(s_2))^2
\right)^{1/4}}\right)
\sum\limits_{j_3=p_3+1}^{\infty}\frac{1}{j_3^2}
\right)^2\le 
$$
$$
\le
\frac{(T-t)^{4l+2l_1+3}K^4 K_1^2}{64}\ \ \frac{2p_3+1}{p_3^2}
\left(\int\limits_{-1}^1
\frac{2dy}
{(1-y^2)^{3/4}}
+\frac{l^2\pi}{2l-1}
\int\limits_{-1}^1
\frac{dy}
{(1-y^2)^{1/4}}\right)^2\le
$$

\begin{equation}
\label{ogo5000}
~~~~~\le C(T-t)^{4l+2l_1+3}\ \frac{2p_3+1}{p_3^2}\to 0\ \ \ 
\hbox{when}\ \ p_3\ \to\ \infty,
\end{equation}

\vspace{2mm}
\noindent
where constant $C$ does not depend on $p_3$ and $T-t$.

The relations (\ref{ogo2100}), (\ref{ogo2110}), and (\ref{ogo5000}) 
imply
(\ref{ogo2000}), and the relation (\ref{ogo2000}) 
implies the correctness of the expansion
(\ref{ogo101ee}).

Let us consider Case 4, i.e. $l_1=l_2=l_3=l=0, 1, 2,\ldots$ and 
$i_1, i_2, i_3=1,\ldots,m$. So, we will prove the following expansion
for iterated Stratonovich stochastic integral of third multiplicity
\begin{equation}
\label{ogo10100}
I_{{l l l}_{T,t}}^{*(i_1i_2i_3)}=
\hbox{\vtop{\offinterlineskip\halign{
\hfil#\hfil\cr
{\rm l.i.m.}\cr
$\stackrel{}{{}_{p_1,p_2,p_3\to \infty}}$\cr
}} }\sum_{j_1=0}^{p_1}\sum_{j_2=0}^{p_2}\sum_{j_3=0}^{p_3}
C_{j_3 j_2 j_1}\zeta_{j_1}^{(i_1)}\zeta_{j_2}^{(i_2)}\zeta_{j_3}^{(i_3)}\ \ \
(i_1, i_2, i_3=1,\ldots,m),
\end{equation}
where the series converges in the mean-square sense,
$l=0, 1, 2,\ldots$\ \hspace{-0.7mm}, and
\begin{equation}
\label{ogo19900}
C_{j_3 j_2 j_1}=\int\limits_t^T
\phi_{j_3}(s)(t-s)^{l}\int\limits_t^s(t-s_1)^{l}
\phi_{j_2}(s_1)
\int\limits_t^{s_1}(t-s_2)^{l}
\phi_{j_1}(s_2)ds_2ds_1ds.
\end{equation}

If we prove w.~p.~1 the following formula
\begin{equation}
\label{ogo20000}
\hbox{\vtop{\offinterlineskip\halign{
\hfil#\hfil\cr
{\rm l.i.m.}\cr
$\stackrel{}{{}_{p_1,p_3\to \infty}}$\cr
}} }\sum_{j_1=0}^{p_1}\sum_{j_3=0}^{p_3}
C_{j_1 j_3 j_1}\zeta_{j_3}^{(i_2)}=0,
\end{equation}
where 
the coefficients $C_{j_3 j_2 j_1}$ 
are defined by
(\ref{ogo19900}), then
using the formulas (\ref{ogo200}), (\ref{ogo2000})
when $l_1=l_3=l$, Theorem 1.1,
and standard relations
between 
iterated It\^{o} and Stratonovich stochastic integrals,
we obtain the expansion (\ref{ogo10100}).

Since 
$\psi_1(s),$ $\psi_2(s),$ $\psi_3(s)\equiv (t-s)^l$, then
the following equality for the Fourier coefficients takes place
$$
C_{j_1 j_1 j_3}+C_{j_1 j_3 j_1}+C_{j_3 j_1 j_1}=\frac{1}{2}
C_{j_1}^2 C_{j_3},
$$ 
where the coefficients
$C_{j_3 j_2 j_1}$ are defined by (\ref{ogo19900}) and 
$$
C_{j_1}=\int\limits_t^T
\phi_{j_1}(s)(t-s)^{l}ds.
$$

Then w.~p.~1
$$
\hbox{\vtop{\offinterlineskip\halign{
\hfil#\hfil\cr
{\rm l.i.m.}\cr
$\stackrel{}{{}_{p_1,p_3\to \infty}}$\cr
}} }\sum_{j_1=0}^{p_1}\sum_{j_3=0}^{p_3}
C_{j_1 j_3 j_1}\zeta_{j_3}^{(i_2)}=
$$
\begin{equation}
\label{sodom310}
~~~~~~~~~~ =
\hbox{\vtop{\offinterlineskip\halign{
\hfil#\hfil\cr
{\rm l.i.m.}\cr
$\stackrel{}{{}_{p_1,p_3\to \infty}}$\cr
}} }\sum_{j_1=0}^{p_1}\sum_{j_3=0}^{p_3}
\left(\frac{1}{2}C_{j_1}^2 C_{j_3}-C_{j_1 j_1 j_3}-C_{j_3 j_1 j_1}
\right)\zeta_{j_3}^{(i_2)}.
\end{equation}

Taking into account (\ref{ogo200}) and (\ref{ogo2000}) 
when $l_3=l_1=l$ as well as the It\^{o} formula,
we have w.~p.~1
$$
\hbox{\vtop{\offinterlineskip\halign{
\hfil#\hfil\cr
{\rm l.i.m.}\cr
$\stackrel{}{{}_{p_1,p_3\to \infty}}$\cr
}} }\sum_{j_1=0}^{p_1}\sum_{j_3=0}^{p_3}
C_{j_1 j_3 j_1}\zeta_{j_3}^{(i_2)}
=
\frac{1}{2}\sum\limits_{j_1=0}^{l}
C_{j_1}^2\sum\limits_{j_3=0}^{l}C_{j_3}\zeta_{j_3}^{(i_2)}-
$$
$$
-
\hbox{\vtop{\offinterlineskip\halign{
\hfil#\hfil\cr
{\rm l.i.m.}\cr
$\stackrel{}{{}_{p_1,p_3\to \infty}}$\cr
}} }\sum_{j_1=0}^{p_1}\sum_{j_3=0}^{p_3}
C_{j_1 j_1 j_3}\zeta_{j_3}^{(i_2)}-
\hbox{\vtop{\offinterlineskip\halign{
\hfil#\hfil\cr
{\rm l.i.m.}\cr
$\stackrel{}{{}_{p_1,p_3\to \infty}}$\cr
}} }\sum_{j_1=0}^{p_1}\sum_{j_3=0}^{p_3}
C_{j_3 j_1 j_1}\zeta_{j_3}^{(i_2)}
=
$$
$$
=\frac{1}{2}\sum\limits_{j_1=0}^{l}
C_{j_1}^2\int\limits_t^T(t-s)^ld{\bf w}_s^{(i_2)}
-\frac{1}{2}\int\limits_t^T(t-s)^{l}
\int\limits_t^s(t-s_1)^{2l}ds_1d{\bf w}_s^{(i_2)}-
$$
$$
-\frac{1}{2}\int\limits_t^T(t-s)^{2l}
\int\limits_t^s(t-s_1)^{l}d{\bf w}_{s_1}^{(i_2)}ds=
$$
$$
=\frac{1}{2}\sum\limits_{j_1=0}^{l}
C_{j_1}^2\int\limits_t^T(t-s)^ld{\bf w}_s^{(i_2)}
+\frac{1}{2(2l+1)}\int\limits_t^T(t-s)^{3l+1}d{\bf w}_s^{(i_2)}-
$$
$$
-\frac{1}{2}\int\limits_t^T(t-s_1)^{l}
\int\limits_{s_1}^T(t-s)^{2l}dsd{\bf w}_{s_1}^{(i_2)}=
$$
$$
=\frac{1}{2}\sum\limits_{j_1=0}^{l}
C_{j_1}^2\int\limits_t^T(t-s)^ld{\bf w}_s^{(i_2)}
+\frac{1}{2(2l+1)}\int\limits_t^T(t-s)^{3l+1}d{\bf w}_s^{(i_2)}-
$$
$$
-\frac{1}{2(2l+1)}\left((T-t)^{2l+1}\int\limits_t^T(t-s)^{l}
d{\bf w}_{s}^{(i_2)}+\int\limits_t^T(t-s)^{3l+1}
d{\bf w}_{s}^{(i_2)}\right)=
$$
$$
=\frac{1}{2}\sum\limits_{j_1=0}^{l}
C_{j_1}^2\int\limits_t^T(t-s)^ld{\bf w}_s^{(i_2)}-
\frac{(T-t)^{2l+1}}{2(2l+1)}\int\limits_t^T(t-s)^{l}
d{\bf w}_{s}^{(i_2)}=
$$
$$
=
\frac{1}{2}\left(\sum\limits_{j_1=0}^{l}
C_{j_1}^2-\int\limits_t^T(t-s)^{2l}ds\right)
\int\limits_t^T(t-s)^ld{\bf w}_s^{(i_2)}=0.
$$

\vspace{2mm}

Here the Parseval equality looks as follows
$$
\sum\limits_{j_1=0}^{\infty}
C_{j_1}^2=
\sum\limits_{j_1=0}^{l}
C_{j_1}^2=\int\limits_t^T(t-s)^{2l}ds=\frac{(T-t)^{2l+1}}{2l+1}
$$
and      
$$
\int\limits_t^{T}(t-s)^{l}
d{\bf w}_{s}^{(i_2)}=
\sum\limits_{j_3=0}^{l}C_{j_3}\zeta_{j_3}^{(i_2)}\ \ \ \hbox{w.~p.~1}.
$$

The expansion (\ref{ogo10100}) is proved. Theorem 2.5 is proved.

It is easy to see that using the It\^{o} formula if
$i_1=i_2=i_3$ 
we obtain (see (\ref{nahod}))
$$
{\int\limits_t^{*}}^T
(t-s)^{l}
{\int\limits_t^{*}}^s
(t-s_1)^l
{\int\limits_t^{*}}^{s_1}
(t-s_2)^{l}
d{\bf w}_{s_2}^{(i_1)}d{\bf w}_{s_1}^{(i_1)}d{\bf w}_{s}^{(i_1)}=
$$
$$
=\frac{1}{6}\left(
\int\limits_t^{T}(t-s)^{l}
d{\bf w}_{s}^{(i_1)}\right)^3
=\frac{1}{6}\left(
\sum\limits_{j_1=0}^{l}C_{j_1}\zeta_{j_1}^{(i_1)}\right)^3=
$$
\begin{equation}
\label{sodom400}
=
\sum\limits_{j_1, j_2, j_3=0}^{l}
C_{j_3 j_2 j_1}\zeta_{j_1}^{(i_1)}\zeta_{j_2}^{(i_1)}\zeta_{j_3}^{(i_1)}\ \ \
\hbox{w.~p.~1}.
\end{equation}

\subsection{The Case $p_1, p_2, p_3\to \infty$ and Constant 
Weight Functions (The Case of Trigonometric Functions)}

In this section, we will prove the following theorem.

{\bf Theorem 2.6}\ \cite{6}-\cite{12aa}, \cite{arxiv-7}.
{\it Suppose that
$\{\phi_j(x)\}_{j=0}^{\infty}$ is a complete orthonormal
system of trigonometric functions
in the space $L_2([t, T])$.
Then$,$ for the iterated Stratonovich stochastic integral 
of third multiplicity
$$
{\int\limits_t^{*}}^T
{\int\limits_t^{*}}^{t_3}
{\int\limits_t^{*}}^{t_2}
d{\bf w}_{t_1}^{(i_1)}
d{\bf w}_{t_2}^{(i_2)}d{\bf w}_{t_3}^{(i_3)}\ \ \ (i_1, i_2, i_3=1,\ldots,m)
$$
the following 
expansion 
\begin{equation}
\label{feto19001ee}
~{\int\limits_t^{*}}^T
{\int\limits_t^{*}}^{t_3}
{\int\limits_t^{*}}^{t_2}
d{\bf w}_{t_1}^{(i_1)}
d{\bf w}_{t_2}^{(i_2)}d{\bf w}_{t_3}^{(i_3)}\ 
=
\hbox{\vtop{\offinterlineskip\halign{
\hfil#\hfil\cr
{\rm l.i.m.}\cr
$\stackrel{}{{}_{p_1,p_2,p_3\to \infty}}$\cr
}} }\sum_{j_1=0}^{p_1}\sum_{j_2=0}^{p_2}\sum_{j_3=0}^{p_3}
C_{j_3 j_2 j_1}\zeta_{j_1}^{(i_1)}\zeta_{j_2}^{(i_2)}\zeta_{j_3}^{(i_3)}
\end{equation}
that converges in the mean-square sense is valid, where
$$
C_{j_3 j_2 j_1}=\int\limits_t^T
\phi_{j_3}(s)\int\limits_t^s
\phi_{j_2}(s_1)
\int\limits_t^{s_1}
\phi_{j_1}(s_2)ds_2ds_1ds
$$
and
$$
\zeta_{j}^{(i)}=
\int\limits_t^T \phi_{j}(s) d{\bf w}_s^{(i)}
$$ 
are independent standard Gaussian random variables for various 
$i$ or $j$.}

{\bf Proof.}\ If we prove w.~p.~1 the following formulas
\begin{equation}
\label{ogo1299}
\hbox{\vtop{\offinterlineskip\halign{
\hfil#\hfil\cr
{\rm l.i.m.}\cr
$\stackrel{}{{}_{p_1, p_3\to \infty}}$\cr
}} }
\sum\limits_{j_1=0}^{p_1}\sum\limits_{j_3=0}^{p_3}
C_{j_3 j_1 j_1}\zeta_{j_3}^{(i_3)}
=
\frac{1}{2}\int\limits_t^T\int\limits_t^{\tau}dsd{\bf w}_{\tau}^{(i_3)},
\end{equation}
\begin{equation}
\label{ogo1399}
\hbox{\vtop{\offinterlineskip\halign{
\hfil#\hfil\cr
{\rm l.i.m.}\cr
$\stackrel{}{{}_{p_1, p_3\to \infty}}$\cr
}} }
\sum\limits_{j_1=0}^{p_1}\sum\limits_{j_3=0}^{p_3}
C_{j_3 j_3 j_1}\zeta_{j_1}^{(i_1)}
=
\frac{1}{2}\int\limits_t^T\int\limits_t^{\tau}d{\bf w}_{s}^{(i_1)}d\tau,
\end{equation}
\begin{equation}
\label{ogo13a99}
\hbox{\vtop{\offinterlineskip\halign{
\hfil#\hfil\cr
{\rm l.i.m.}\cr
$\stackrel{}{{}_{p_1, p_3\to \infty}}$\cr
}} }
\sum\limits_{j_1=0}^{p_1}\sum\limits_{j_3=0}^{p_3}
C_{j_1 j_3 j_1}\zeta_{j_3}^{(i_2)}
=0,
\end{equation}

\noindent
then from the equalities (\ref{ogo1299})--(\ref{ogo13a99}),
Theorem 1.1, and
standard relations 
between 
iterated It\^{o} and
Stratonovich stochastic integrals we will obtain
the expansion (\ref{feto19001ee}).

We have
$$
S_{p_1,p_3}\stackrel{\sf def}{=}\sum\limits_{j_3=0}^{p_3}\sum\limits_{j_1=0}^{p_1}
C_{j_3 j_1 j_1}\zeta_{j_3}^{(i_3)}=
\frac{(T-t)^{3/2}}{6}\zeta_{0}^{(i_3)}+
$$
$$
+\sum\limits_{j_1=1}^{p_1}C_{0,2j_1,2j_1}
\zeta_0^{(i_3)}+\sum\limits_{j_1=1}^{p_1}C_{0,2j_1-1,2j_1-1}\zeta_0^{(i_3)}+
\sum\limits_{j_3=1}^{p_1}C_{2j_3,0,0}\zeta_{2j_3}^{(i_3)}+
$$
$$
+
\sum\limits_{j_3=1}^{p_3}\sum_{j_1=1}^{p_1}
C_{2j_3,2j_1,2j_1}\zeta_{2j_3}^{(i_3)}+
\sum\limits_{j_3=1}^{p_3}\sum_{j_1=1}^{p_1}
C_{2j_3,2j_1-1,2j_1-1}\zeta_{2j_3}^{(i_3)}+
\sum\limits_{j_3=1}^{p_3}
C_{2j_3-1,0,0}\zeta_{2j_3-1}^{(i_3)}+
$$
\begin{equation}
\label{ogo900}
~~~~~~~~ +\sum\limits_{j_3=1}^{p_3}\sum_{j_1=1}^{p_1}
C_{2j_3-1,2j_1,2j_1}\zeta_{2j_3-1}^{(i_3)}+
\sum\limits_{j_3=1}^{p_3}\sum_{j_1=1}^{p_1}
C_{2j_3-1,2j_1-1,2j_1-1}\zeta_{2j_3-1}^{(i_3)},
\end{equation}

\vspace{2mm}
\noindent
where the summation is stopped, when $2j_1,$ $2j_1-1> p_1$
or $2j_3,$ $2j_3-1> p_3$ and
\begin{equation}
\label{ogo901}
C_{0,2l,2l}=\frac{(T-t)^{3/2}}{8\pi^2l^2},\ \ \
C_{0,2l-1,2l-1}=\frac{3(T-t)^{3/2}}{8\pi^2l^2},\ \ \
C_{2l,0,0}=\frac{\sqrt{2}(T-t)^{3/2}}{4\pi^2l^2},
\end{equation}
\begin{equation}
\label{ogo903}
~~~~C_{2r-1,2l,2l}=0,\ \ \
C_{2l-1,0,0}=-\frac{\sqrt{2}(T-t)^{3/2}}{4\pi l},\ \ \
C_{2r-1,2l-1,2l-1}=0,
\end{equation}
\begin{equation}
\label{ogo902}
C_{2r,2l,2l}=
\begin{cases}
-\sqrt{2}(T-t)^{3/2}/(16\pi^2l^2),\ &r=2l
\cr
\cr
0,\ & r\ne 2l\
\end{cases},
\end{equation}
\begin{equation}
\label{ogo902a}
~~~~~~~ C_{2r,2l-1,2l-1}=
\begin{cases}
\sqrt{2}(T-t)^{3/2}/(16\pi^2l^2),\ &r=2l
\cr
\cr
-\sqrt{2}(T-t)^{3/2}/(4\pi^2l^2),\ &r=l
\cr
\cr
0,\ & r\ne l,\ r\ne 2l
\end{cases}.
\end{equation}

Let us show that
\begin{equation}
\label{agenty100}
\hbox{\vtop{\offinterlineskip\halign{
\hfil#\hfil\cr
{\rm l.i.m.}\cr
$\stackrel{}{{}_{p_1, p_3\to \infty}}$\cr
}} }S_{2p_1,2p_3}=
\hbox{\vtop{\offinterlineskip\halign{
\hfil#\hfil\cr
{\rm l.i.m.}\cr
$\stackrel{}{{}_{p_1, p_3\to \infty}}$\cr
}} }S_{2p_1,2p_3-1}=
\hbox{\vtop{\offinterlineskip\halign{
\hfil#\hfil\cr
{\rm l.i.m.}\cr
$\stackrel{}{{}_{p_1, p_3\to \infty}}$\cr
}} }S_{2p_1-1,2p_3-1}=
\hbox{\vtop{\offinterlineskip\halign{
\hfil#\hfil\cr
{\rm l.i.m.}\cr
$\stackrel{}{{}_{p_1, p_3\to \infty}}$\cr
}} }S_{2p_1-1,2p_3}.
\end{equation}

We have
\begin{equation}
\label{agenty101}
S_{2p_1,2p_3}=S_{2p_1,2p_3-1}+
\sum\limits_{j_1=0}^{2p_1}
C_{2p_3, j_1, j_1}\zeta_{2p_3}^{(i_3)}.
\end{equation}

Using the relations (\ref{ogo901}), (\ref{ogo902}), and (\ref{ogo902a}), we obtain
$$
\sum\limits_{j_1=0}^{2p_1}
C_{2p_3, j_1, j_1}=
C_{2p_3, 0, 0}+\sum\limits_{j_1=1}^{2p_1}
C_{2p_3, j_1, j_1}=
$$
$$
=C_{2p_3, 0, 0}+\sum\limits_{j_1=1}^{p_1}
\biggl(C_{2p_3, 2j_1-1, 2j_1-1}+C_{2p_3, 2j_1, 2j_1}\biggr)=
$$
\begin{equation}
\label{agenty102}
~~ =\frac{\sqrt{2}(T-t)^{3/2}}{4\pi^2 p_3^2}\bigl(1-{\bf 1}_{\{p_1\ge p_3\}}\bigr).
\end{equation}

\vspace{2mm}

From (\ref{agenty101}), (\ref{agenty102}) we obtain
\begin{equation}
\label{agenty103}
\hbox{\vtop{\offinterlineskip\halign{
\hfil#\hfil\cr
{\rm l.i.m.}\cr
$\stackrel{}{{}_{p_1, p_3\to \infty}}$\cr
}} }S_{2p_1,2p_3}=
\hbox{\vtop{\offinterlineskip\halign{
\hfil#\hfil\cr
{\rm l.i.m.}\cr
$\stackrel{}{{}_{p_1, p_3\to \infty}}$\cr
}} }S_{2p_1,2p_3-1}.
\end{equation}

\vspace{2mm}

Further, we get (see (\ref{ogo901})--(\ref{ogo902}))
\begin{equation}
\label{agenty104}
S_{2p_1,2p_3-1}=S_{2p_1-1,2p_3-1}+
\sum\limits_{j_3=0}^{2p_3-1}
C_{j_3, 2p_1, 2p_1}\zeta_{j_3}^{(i_3)},
\end{equation}
$$
\sum\limits_{j_3=0}^{2p_3-1}
C_{j_3, 2p_1, 2p_1}\zeta_{j_3}^{(i_3)}=
C_{0, 2p_1, 2p_1}\zeta_{0}^{(i_3)}+
\sum\limits_{j_3=1}^{2p_3}
C_{j_3, 2p_1, 2p_1}\zeta_{j_3}^{(i_3)}-
C_{2p_3, 2p_1, 2p_1}\zeta_{2p_3}^{(i_3)}=
$$
$$
=
C_{0, 2p_1, 2p_1}\zeta_{0}^{(i_3)}+
\sum\limits_{j_3=1}^{p_3}\biggl(
C_{2j_3-1, 2p_1, 2p_1}\zeta_{2j_3-1}^{(i_3)}
+C_{2j_3, 2p_1, 2p_1}\zeta_{2j_3}^{(i_3)}\biggr)
-
C_{2p_3, 2p_1, 2p_1}\zeta_{2p_3}^{(i_3)}=
$$
\begin{equation}
\label{agenty105}
~~~~~~=\frac{(T-t)^{3/2}}{8\pi^2 p_1^2}\zeta_{0}^{(i_3)}+\frac{\sqrt{2}(T-t)^{3/2}}{16\pi^2 p_1^2}
\bigl({\bf 1}_{\{p_3=2p_1\}}-{\bf 1}_{\{p_3\ge 2p_1\}}\bigr)\zeta_{4p_1}^{(i_3)}.
\end{equation}

\vspace{2mm}

From (\ref{agenty104}), (\ref{agenty105}) we obtain
\begin{equation}
\label{agenty106}
\hbox{\vtop{\offinterlineskip\halign{
\hfil#\hfil\cr
{\rm l.i.m.}\cr
$\stackrel{}{{}_{p_1, p_3\to \infty}}$\cr
}} }S_{2p_1,2p_3-1}=
\hbox{\vtop{\offinterlineskip\halign{
\hfil#\hfil\cr
{\rm l.i.m.}\cr
$\stackrel{}{{}_{p_1, p_3\to \infty}}$\cr
}} }S_{2p_1-1,2p_3-1}.
\end{equation}

\vspace{2mm}

Further, we have
\begin{equation}
\label{agenty107}
S_{2p_1,2p_3}=S_{2p_1-1,2p_3}+
\sum\limits_{j_3=0}^{2p_3}
C_{j_3, 2p_1, 2p_1}\zeta_{j_3}^{(i_3)},
\end{equation}
$$
\sum\limits_{j_3=0}^{2p_3}
C_{j_3, 2p_1, 2p_1}\zeta_{j_3}^{(i_3)}=
C_{0, 2p_1, 2p_1}\zeta_{0}^{(i_3)}+
\sum\limits_{j_3=1}^{2p_3}
C_{j_3, 2p_1, 2p_1}\zeta_{j_3}^{(i_3)}=
$$
\begin{equation}
\label{agenty108}
~~~~~~~=
C_{0, 2p_1, 2p_1}\zeta_{0}^{(i_3)}+
\sum\limits_{j_3=1}^{p_3}
\biggl(C_{2j_3-1, 2p_1, 2p_1}\zeta_{2j_3-1}^{(i_3)}+
C_{2j_3, 2p_1, 2p_1}\zeta_{2j_3}^{(i_3)}\biggr).
\end{equation}

\vspace{2mm}

From (\ref{agenty108}), (\ref{ogo901})--(\ref{ogo902}) we obtain
\begin{equation}
\label{agenty109}
~~~~~~ \sum\limits_{j_3=0}^{2p_3}
C_{j_3, 2p_1, 2p_1}\zeta_{j_3}^{(i_3)}=
\frac{(T-t)^{3/2}}{8\pi^2 p_1^2}
\zeta_{0}^{(i_3)}-
\frac{\sqrt{2}(T-t)^{3/2}}{16\pi^2 p_1^2}
{\bf 1}_{\{p_3\ge 2p_1\}}\zeta_{4p_1}^{(i_3)}.
\end{equation}

\vspace{2mm}

The relations (\ref{agenty107}), (\ref{agenty109}) mean that
\begin{equation}
\label{agenty110}
\hbox{\vtop{\offinterlineskip\halign{
\hfil#\hfil\cr
{\rm l.i.m.}\cr
$\stackrel{}{{}_{p_1, p_3\to \infty}}$\cr
}} }S_{2p_1,2p_3}=
\hbox{\vtop{\offinterlineskip\halign{
\hfil#\hfil\cr
{\rm l.i.m.}\cr
$\stackrel{}{{}_{p_1, p_3\to \infty}}$\cr
}} }S_{2p_1-1,2p_3}.
\end{equation}

\vspace{2mm}

The equalities (\ref{agenty103}), (\ref{agenty106}), and (\ref{agenty110})
imply (\ref{agenty100}). This means that instead of (\ref{ogo1299}) it is enough
to prove the following equality
\begin{equation}
\label{agenty111}
~~~~~~~~ \hbox{\vtop{\offinterlineskip\halign{
\hfil#\hfil\cr
{\rm l.i.m.}\cr
$\stackrel{}{{}_{p_1, p_3\to \infty}}$\cr
}} }
\sum\limits_{j_1=0}^{2p_1}\sum\limits_{j_3=0}^{2p_3}
C_{j_3 j_1 j_1}\zeta_{j_3}^{(i_3)}
=
\frac{1}{2}\int\limits_t^T\int\limits_t^{\tau}dsd{\bf w}_{\tau}^{(i_3)}\ \ \ \hbox{w.~p.~1}.
\end{equation}

\vspace{2mm}

We have
$$
S_{2p_1,2p_3}=\sum\limits_{j_3=0}^{2p_3}\sum\limits_{j_1=0}^{2p_1}
C_{j_3 j_1 j_1}\zeta_{j_3}^{(i_3)}=
\frac{(T-t)^{3/2}}{6}\zeta_{0}^{(i_3)}+
$$
$$
+\sum\limits_{j_1=1}^{p_1}C_{0,2j_1,2j_1}
\zeta_0^{(i_3)}+\sum\limits_{j_1=1}^{p_1}C_{0,2j_1-1,2j_1-1}\zeta_0^{(i_3)}+
\sum\limits_{j_3=1}^{p_1}C_{2j_3,0,0}\zeta_{2j_3}^{(i_3)}+
$$
$$
+
\sum\limits_{j_3=1}^{p_3}\sum_{j_1=1}^{p_1}
C_{2j_3,2j_1,2j_1}\zeta_{2j_3}^{(i_3)}+
\sum\limits_{j_3=1}^{p_3}\sum_{j_1=1}^{p_1}
C_{2j_3,2j_1-1,2j_1-1}\zeta_{2j_3}^{(i_3)}+
\sum\limits_{j_3=1}^{p_3}
C_{2j_3-1,0,0}\zeta_{2j_3-1}^{(i_3)}+
$$
\begin{equation}
\label{agenty112}
~~~~~~~~~ +\sum\limits_{j_3=1}^{p_3}\sum_{j_1=1}^{p_1}
C_{2j_3-1,2j_1,2j_1}\zeta_{2j_3-1}^{(i_3)}+
\sum\limits_{j_3=1}^{p_3}\sum_{j_1=1}^{p_1}
C_{2j_3-1,2j_1-1,2j_1-1}\zeta_{2j_3-1}^{(i_3)}.
\end{equation}

\vspace{2mm}

After 
substituting 
(\ref{ogo901})--(\ref{ogo902a})
into (\ref{agenty112}), we obtain
$$
\sum\limits_{j_3=0}^{2p_3}\sum\limits_{j_1=0}^{2p_1}
C_{j_3 j_1 j_1}\zeta_{j_3}^{(i_3)}=(T-t)^{3/2}
\left(\frac{1}{6}\zeta_{0}^{(i_3)} + \frac{1}{2\pi^2}
\sum\limits_{j_1=1}^{p_1}\frac{1}{j_1^2}\ \zeta_{0}^{(i_3)}-\right.
$$
\begin{equation}
\label{agenty113}
~~~~~~~ \left.-
\frac{\sqrt{2}}{4\pi}\sum\limits_{j_3=1}^{p_3}\frac{1}{j_3}
\zeta_{2j_3-1}^{(i_3)}-
\frac{\sqrt{2}}{4\pi^2}\sum\limits_{j_3=1}^{\min\{p_1,p_3\}}\frac{1}{j_3^2}
\zeta_{2j_3}^{(i_3)}+
\frac{\sqrt{2}}{4\pi^2}\sum\limits_{j_3=1}^{p_3}\frac{1}{j_3^2}
\zeta_{2j_3}^{(i_3)}\right).
\end{equation}

\vspace{2mm}

From (\ref{agenty113}) we have w.~p.~1
$$
\hbox{\vtop{\offinterlineskip\halign{
\hfil#\hfil\cr
{\rm l.i.m.}\cr
$\stackrel{}{{}_{p_1, p_3\to \infty}}$\cr
}} }
\sum\limits_{j_3=0}^{2p_3}\sum\limits_{j_1=0}^{2p_1}
C_{j_3 j_1 j_1}\zeta_{j_3}^{(i_3)}=(T-t)^{3/2}
\left(\frac{1}{6}\zeta_{0}^{(i_3)} + \frac{1}{2\pi^2}
\sum\limits_{j_1=1}^{\infty}\frac{1}{j_1^2}\ \zeta_{0}^{(i_3)}-\right.
$$
\begin{equation}
\label{ogo905}
\left.-\hbox{\vtop{\offinterlineskip\halign{
\hfil#\hfil\cr
{\rm l.i.m.}\cr
$\stackrel{}{{}_{p_3\to \infty}}$\cr
}} }
\frac{\sqrt{2}}{4\pi}\sum\limits_{j_3=1}^{p_3}\frac{1}{j_3}
\zeta_{2j_3-1}^{(i_3)}\right).
\end{equation}

Using Theorem 1.1 and the system of trigonometric functions, we get
w.~p.~1
$$
\frac{1}{2}\int\limits_t^T\int\limits_t^s d\tau d{\bf w}_{s}^{(i_3)}=
\frac{1}{2}\int\limits_t^T(s-t)d{\bf w}_{s}^{(i_3)}
=
$$
\begin{equation}
\label{ogo906}
=\frac{(T-t)^{3/2}}{4}\hbox{\vtop{\offinterlineskip\halign{
\hfil#\hfil\cr
{\rm l.i.m.}\cr
$\stackrel{}{{}_{p_3\to \infty}}$\cr
}} }
\Biggl(\zeta_0^{(i_3)}-\frac{\sqrt{2}}{\pi}\sum_{j_3=1}^{p_3}
\frac{1}{j_3}
\zeta_{2j_3-1}^{(i_3)}
\Biggr).
\end{equation}

\vspace{2mm}

From (\ref{ogo905}) and (\ref{ogo906}) it follows that
$$
\hbox{\vtop{\offinterlineskip\halign{
\hfil#\hfil\cr
{\rm l.i.m.}\cr
$\stackrel{}{{}_{p_1, p_3\to \infty}}$\cr
}} }
\sum\limits_{j_3=0}^{2p_3}\sum\limits_{j_1=0}^{2p_1}
C_{j_3 j_1 j_1}\zeta_{j_3}^{(i_3)}=
$$
$$
=(T-t)^{3/2}
\left(\frac{1}{6}\zeta_{0}^{(i_3)} + \frac{1}{12}
\zeta_{0}^{(i_3)}-\hbox{\vtop{\offinterlineskip\halign{
\hfil#\hfil\cr
{\rm l.i.m.}\cr
$\stackrel{}{{}_{p_3\to \infty}}$\cr
}} }
\frac{\sqrt{2}}{4\pi}\sum\limits_{j_3=1}^{p_3}\frac{1}{j_3}
\zeta_{2j_3-1}^{(i_3)}\right)=
$$
$$
=(T-t)^{3/2}
\left(\frac{1}{4}\zeta_{0}^{(i_3)}
-\hbox{\vtop{\offinterlineskip\halign{
\hfil#\hfil\cr
{\rm l.i.m.}\cr
$\stackrel{}{{}_{p_3\to \infty}}$\cr
}} }
\frac{\sqrt{2}}{4\pi}\sum\limits_{j_3=1}^{p_3}\frac{1}{j_3}
\zeta_{2j_3-1}^{(i_3)}\right)
=
$$
$$
=\frac{1}{2}\int\limits_t^T\int\limits_t^s d\tau d{\bf w}_{s}^{(i_3)},
$$

\vspace{2mm}
\noindent
where the equality is fulfilled w.~p.~1.

So, the relations (\ref{agenty111}) and (\ref{ogo1299}) are proved for the case of 
trigonometric 
system of functions.

Let us prove the relation (\ref{ogo1399}). We have
$$
S'_{p_1,p_3}\stackrel{\sf def}{=}\sum\limits_{j_1=0}^{p_1}\sum\limits_{j_3=0}^{p_3}
C_{j_3 j_3 j_1}\zeta_{j_1}^{(i_1)}=
\frac{(T-t)^{3/2}}{6}\zeta_0^{(i_1)}+
$$
$$
+
\sum\limits_{j_3=1}^{p_3}C_{2j_3,2j_3,0}
\zeta_0^{(i_1)}+\sum\limits_{j_3=1}^{p_3}C_{2j_3-1,2j_3-1,0}\zeta_0^{(i_1)}+
\sum\limits_{j_1=1}^{p_1}\sum\limits_{j_3=1}^{p_3}
C_{2j_3,2j_3,2j_1-1}\zeta_{2j_1-1}^{(i_1)}+
$$
$$
+
\sum\limits_{j_1=1}^{p_1}\sum_{j_3=1}^{p_3}
C_{2j_3-1,2j_3-1,2j_1-1}\zeta_{2j_1-1}^{(i_1)}+
\sum\limits_{j_1=1}^{p_1}
C_{0,0,2j_1-1}\zeta_{2j_1-1}^{(i_1)}+
\sum\limits_{j_1=1}^{p_1}\sum\limits_{j_3=1}^{p_3}
C_{2j_3,2j_3,2j_1}\zeta_{2j_1}^{(i_1)}+
$$
\begin{equation}
\label{ogo9000}
+\sum\limits_{j_1=1}^{p_1}\sum_{j_3=1}^{p_3}
C_{2j_3-1,2j_3-1,2j_1}\zeta_{2j_1}^{(i_1)}+
\sum_{j_1=1}^{p_1}
C_{0,0,2j_1}\zeta_{2j_1}^{(i_1)},
\end{equation}

\noindent
where the summation is stopped, when
$2j_3,$ $2j_3-1> p_3$
or $2j_1,$ $2j_1-1> p_1$ and
\begin{equation}
\label{ogo9010}
C_{2l,2l,0}=\frac{(T-t)^{3/2}}{8\pi^2l^2},\ \ \
C_{2l-1,2l-1,0}=\frac{3(T-t)^{3/2}}{8\pi^2l^2},\ \ \
C_{0,0,2r}=\frac{\sqrt{2}(T-t)^{3/2}}{4\pi^2r^2},
\end{equation}
\begin{equation}
\label{ogo9030}
~~~~C_{2l-1,2l-1,2r-1}=0,\ \ \
C_{0,0,2r-1}=\frac{\sqrt{2}(T-t)^{3/2}}{4\pi r},\ \ \
C_{2l,2l,2r-1}=0,
\end{equation}
\begin{equation}
\label{ogo9020}
C_{2l,2l,2r}=
\begin{cases}
-\sqrt{2}(T-t)^{3/2}/(16\pi^2l^2),\ &r=2l
\cr
\cr
0,\  &r\ne 2l
\end{cases},
\end{equation}
\begin{equation}
\label{ogo9020a}
~~~~~~~~ C_{2l-1,2l-1,2r}=
\begin{cases}
\sqrt{2}(T-t)^{3/2}/(16\pi^2l^2),\ &r=2l
\cr
\cr
-\sqrt{2}(T-t)^{3/2}/(4\pi^2l^2),\ &r=l
\cr
\cr
0,\  &r\ne l,\ r\ne 2l
\end{cases}.
\end{equation}

Let us show that
\begin{equation}
\label{agenty1000}
\hbox{\vtop{\offinterlineskip\halign{
\hfil#\hfil\cr
{\rm l.i.m.}\cr
$\stackrel{}{{}_{p_1, p_3\to \infty}}$\cr
}} }S'_{2p_1,2p_3}=
\hbox{\vtop{\offinterlineskip\halign{
\hfil#\hfil\cr
{\rm l.i.m.}\cr
$\stackrel{}{{}_{p_1, p_3\to \infty}}$\cr
}} }S'_{2p_1,2p_3-1}=
\hbox{\vtop{\offinterlineskip\halign{
\hfil#\hfil\cr
{\rm l.i.m.}\cr
$\stackrel{}{{}_{p_1, p_3\to \infty}}$\cr
}} }S'_{2p_1-1,2p_3-1}=
\hbox{\vtop{\offinterlineskip\halign{
\hfil#\hfil\cr
{\rm l.i.m.}\cr
$\stackrel{}{{}_{p_1, p_3\to \infty}}$\cr
}} }S'_{2p_1-1,2p_3}.
\end{equation}

We have
\begin{equation}
\label{agenty1010}
S'_{2p_1,2p_3}=S'_{2p_1-1,2p_3}+
\sum\limits_{j_3=0}^{2p_3}
C_{j_3, j_3, 2p_1}\zeta_{2p_1}^{(i_1)}.
\end{equation}

Using the relations (\ref{ogo9010}), (\ref{ogo9020}), and (\ref{ogo9020a}), we obtain
$$
\sum\limits_{j_1=0}^{2p_3}
C_{j_3, j_3, 2p_1}=
C_{0, 0, 2p_1}+\sum\limits_{j_3=1}^{2p_3}
C_{j_3, j_3, 2p_1}=
$$
$$
=C_{0, 0, 2p_1}+\sum\limits_{j_3=1}^{p_3}
\biggl(C_{2j_3-1, 2j_3-1, 2p_1}+C_{2j_3, 2j_3, 2p_1}\biggr)=
$$
\begin{equation}
\label{agenty1020}
=\frac{\sqrt{2}(T-t)^{3/2}}{4\pi^2 p_1^2}\bigl(1-{\bf 1}_{\{p_3\ge p_1\}}\bigr).
\end{equation}

\vspace{2mm}

From (\ref{agenty1010}), (\ref{agenty1020}) we obtain
\begin{equation}
\label{agenty1030}
\hbox{\vtop{\offinterlineskip\halign{
\hfil#\hfil\cr
{\rm l.i.m.}\cr
$\stackrel{}{{}_{p_1, p_3\to \infty}}$\cr
}} }S'_{2p_1,2p_3}=
\hbox{\vtop{\offinterlineskip\halign{
\hfil#\hfil\cr
{\rm l.i.m.}\cr
$\stackrel{}{{}_{p_1, p_3\to \infty}}$\cr
}} }S'_{2p_1-1,2p_3}.
\end{equation}

\vspace{2mm}

Further, we get (see (\ref{ogo9010})--(\ref{ogo9020}))
\begin{equation}
\label{agenty1040}
S'_{2p_1-1,2p_3}=S'_{2p_1-1,2p_3-1}+
\sum\limits_{j_1=0}^{2p_1-1}
C_{2p_3, 2p_3, j_1}\zeta_{j_1}^{(i_1)},
\end{equation}
$$
\sum\limits_{j_1=0}^{2p_1-1}
C_{2p_3, 2p_3, j_1}\zeta_{j_1}^{(i_1)}=
C_{2p_3, 2p_3, 0}\zeta_{0}^{(i_1)}+
\sum\limits_{j_1=1}^{2p_1}
C_{2p_3, 2p_3, j_1}\zeta_{j_1}^{(i_1)}-
C_{2p_3, 2p_3, 2p_1}\zeta_{2p_1}^{(i_1)}=
$$
$$
=
C_{2p_3, 2p_3, 0}\zeta_{0}^{(i_1)}+
\sum\limits_{j_1=1}^{p_1}\biggl(
C_{2p_3, 2p_3, 2j_1-1}\zeta_{2j_1-1}^{(i_1)}
+C_{2p_3, 2p_3, 2j_1}\zeta_{2j_1}^{(i_1)}\biggr)
-
C_{2p_3, 2p_3, 2p_1}\zeta_{2p_1}^{(i_1)}=
$$
\begin{equation}
\label{agenty1050}
~~~~~~~ =\frac{(T-t)^{3/2}}{8\pi^2 p_3^2}\zeta_{0}^{(i_1)}+\frac{\sqrt{2}(T-t)^{3/2}}{16\pi^2 p_3^2}
\bigl({\bf 1}_{\{p_1=2p_3\}}-{\bf 1}_{\{p_1\ge 2p_3\}}\bigr)\zeta_{4p_3}^{(i_1)}.
\end{equation}

\vspace{2mm}

From (\ref{agenty1040}), (\ref{agenty1050}) we obtain
\begin{equation}
\label{agenty1060}
\hbox{\vtop{\offinterlineskip\halign{
\hfil#\hfil\cr
{\rm l.i.m.}\cr
$\stackrel{}{{}_{p_1, p_3\to \infty}}$\cr
}} }S'_{2p_1-1,2p_3}=
\hbox{\vtop{\offinterlineskip\halign{
\hfil#\hfil\cr
{\rm l.i.m.}\cr
$\stackrel{}{{}_{p_1, p_3\to \infty}}$\cr
}} }S'_{2p_1-1,2p_3-1}.
\end{equation}

\vspace{2mm}

Further, we have
\begin{equation}
\label{agenty1070}
S'_{2p_1,2p_3}=S'_{2p_1,2p_3-1}+
\sum\limits_{j_1=0}^{2p_1}
C_{2p_3, 2p_3, j_1}\zeta_{j_1}^{(i_1)},
\end{equation}
$$
\sum\limits_{j_1=0}^{2p_1}
C_{2p_3, 2p_3, j_1}\zeta_{j_1}^{(i_1)}=
C_{2p_3, 2p_3, 0}\zeta_{0}^{(i_1)}+
\sum\limits_{j_1=1}^{2p_1}
C_{2p_3, 2p_3, j_1}\zeta_{j_1}^{(i_1)}=
$$
\begin{equation}
\label{agenty1080}
~~~~~ =
C_{2p_3, 2p_3, 0}\zeta_{0}^{(i_1)}+
\sum\limits_{j_1=1}^{p_1}
\biggl(C_{2p_3, 2p_3, 2j_1-1}\zeta_{2j_1-1}^{(i_1)}+
C_{2p_3, 2p_3, 2j_1}\zeta_{2j_1}^{(i_1)}\biggr).
\end{equation}

\vspace{2mm}

From (\ref{agenty1080}), (\ref{ogo9010})--(\ref{ogo9020}) we obtain
\begin{equation}
\label{agenty1090}
~~~~~ \sum\limits_{j_1=0}^{2p_1}
C_{2p_3, 2p_3, j_1}\zeta_{j_1}^{(i_1)}=
\frac{(T-t)^{3/2}}{8\pi^2 p_3^2}
\zeta_{0}^{(i_1)}-
\frac{\sqrt{2}(T-t)^{3/2}}{16\pi^2 p_3^2}
{\bf 1}_{\{p_1\ge 2p_3\}}\zeta_{4p_3}^{(i_1)}.
\end{equation}

\vspace{2mm}

The relations (\ref{agenty1070}), (\ref{agenty1090}) mean that
\begin{equation}
\label{agenty1100}
\hbox{\vtop{\offinterlineskip\halign{
\hfil#\hfil\cr
{\rm l.i.m.}\cr
$\stackrel{}{{}_{p_1, p_3\to \infty}}$\cr
}} }S'_{2p_1,2p_3}=
\hbox{\vtop{\offinterlineskip\halign{
\hfil#\hfil\cr
{\rm l.i.m.}\cr
$\stackrel{}{{}_{p_1, p_3\to \infty}}$\cr
}} }S'_{2p_1,2p_3-1}.
\end{equation}

\vspace{2mm}

The equalities (\ref{agenty1030}), (\ref{agenty1060}), and (\ref{agenty1100})
imply (\ref{agenty1000}). This means that instead of (\ref{ogo1399}) it is enough
to prove the following equality
\begin{equation}
\label{agenty1110}
~~~~~~~~ \hbox{\vtop{\offinterlineskip\halign{
\hfil#\hfil\cr
{\rm l.i.m.}\cr
$\stackrel{}{{}_{p_1, p_3\to \infty}}$\cr
}} }
\sum\limits_{j_1=0}^{2p_1}\sum\limits_{j_3=0}^{2p_3}
C_{j_3 j_3 j_1}\zeta_{j_1}^{(i_1)}
=
\frac{1}{2}\int\limits_t^T\int\limits_t^{\tau}d{\bf w}_{s}^{(i_1)}d\tau\ \ \ \hbox{w.~p.~1}.
\end{equation}

\vspace{2mm}

We have
$$
S'_{2p_1,2p_3}=\sum\limits_{j_1=0}^{2p_1}\sum\limits_{j_3=0}^{2p_3}
C_{j_3 j_3 j_1}\zeta_{j_1}^{(i_1)}=
\frac{(T-t)^{3/2}}{6}\zeta_0^{(i_1)}+
$$
$$
+
\sum\limits_{j_3=1}^{p_3}C_{2j_3,2j_3,0}
\zeta_0^{(i_1)}+\sum\limits_{j_3=1}^{p_3}C_{2j_3-1,2j_3-1,0}\zeta_0^{(i_1)}+
\sum\limits_{j_1=1}^{p_1}\sum\limits_{j_3=1}^{p_3}
C_{2j_3,2j_3,2j_1-1}\zeta_{2j_1-1}^{(i_1)}+
$$
$$
+
\sum\limits_{j_1=1}^{p_1}\sum_{j_3=1}^{p_3}
C_{2j_3-1,2j_3-1,2j_1-1}\zeta_{2j_1-1}^{(i_1)}+
\sum\limits_{j_1=1}^{p_1}
C_{0,0,2j_1-1}\zeta_{2j_1-1}^{(i_1)}+
\sum\limits_{j_1=1}^{p_1}\sum\limits_{j_3=1}^{p_3}
C_{2j_3,2j_3,2j_1}\zeta_{2j_1}^{(i_1)}+
$$
\begin{equation}
\label{agenty1120}
+\sum\limits_{j_1=1}^{p_1}\sum_{j_3=1}^{p_3}
C_{2j_3-1,2j_3-1,2j_1}\zeta_{2j_1}^{(i_1)}+
\sum_{j_1=1}^{p_1}
C_{0,0,2j_1}\zeta_{2j_1}^{(i_1)}.
\end{equation}

\vspace{2mm}

After 
substituting 
(\ref{ogo9010})--(\ref{ogo9020a})
into (\ref{agenty1120}), we obtain
$$
\sum\limits_{j_1=0}^{2p_1}\sum\limits_{j_3=0}^{2p_3}
C_{j_3 j_3 j_1}\zeta_{j_1}^{(i_1)}=(T-t)^{3/2}
\left(\frac{1}{6}\zeta_{0}^{(i_1)} + \frac{1}{2\pi^2}
\sum\limits_{j_3=1}^{p_3}\frac{1}{j_3^2}\ \zeta_{0}^{(i_1)}+\right.
$$
\begin{equation}
\label{agenty1130}
~~~~~~~ \left.+
\frac{\sqrt{2}}{4\pi}\sum\limits_{j_1=1}^{p_1}\frac{1}{j_1}
\zeta_{2j_1-1}^{(i_1)}-
\frac{\sqrt{2}}{4\pi^2}\sum\limits_{j_1=1}^{\min\{p_1,p_3\}}\frac{1}{j_1^2}
\zeta_{2j_1}^{(i_1)}+
\frac{\sqrt{2}}{4\pi^2}\sum\limits_{j_1=1}^{p_1}\frac{1}{j_1^2}
\zeta_{2j_1}^{(i_1)}\right).
\end{equation}

\vspace{2mm}

From (\ref{agenty1130}) we have w.~p.~1
$$
\hbox{\vtop{\offinterlineskip\halign{
\hfil#\hfil\cr
{\rm l.i.m.}\cr
$\stackrel{}{{}_{p_1, p_3\to \infty}}$\cr
}} }
\sum\limits_{j_1=0}^{2p_1}\sum\limits_{j_3=0}^{2p_3}
C_{j_3 j_3 j_1}\zeta_{j_1}^{(i_1)}=(T-t)^{3/2}
\left(\frac{1}{6}\zeta_{0}^{(i_3)} + \frac{1}{2\pi^2}
\sum\limits_{j_3=1}^{\infty}\frac{1}{j_3^2}\ \zeta_{0}^{(i_1)}+\right.
$$
\begin{equation}
\label{ogo9050}
\left.+\hbox{\vtop{\offinterlineskip\halign{
\hfil#\hfil\cr
{\rm l.i.m.}\cr
$\stackrel{}{{}_{p_1\to \infty}}$\cr
}} }
\frac{\sqrt{2}}{4\pi}\sum\limits_{j_1=1}^{p_1}\frac{1}{j_1}
\zeta_{2j_1-1}^{(i_1)}\right).
\end{equation}

\vspace{2mm}

Using the It\^{o} formula and Theorem 1.1 for the case of 
trigonometric system of 
functions, we obtain w. p. 1
$$
\frac{1}{2}\int\limits_t^T\int\limits_t^{\tau}d{\bf w}_{s}^{(i_1)}d\tau=
\frac{1}{2}\left((T-t)
\int\limits_t^Td{\bf w}_{s}^{(i_1)}+
\int\limits_t^T(t-s)d{\bf w}_{s}^{(i_1)}\right)=
$$
\begin{equation}
\label{ogo9060}
=\frac{1}{4}(T-t)^{3/2}
\Biggl(\zeta_0^{(i_1)}+\hbox{\vtop{\offinterlineskip\halign{
\hfil#\hfil\cr
{\rm l.i.m.}\cr
$\stackrel{}{{}_{p_1\to \infty}}$\cr
}} }
\frac{\sqrt{2}}{\pi}\sum_{j_1=1}^{p_1}
\frac{1}{j_1}
\zeta_{2j_1-1}^{(i_1)}
\Biggr).
\end{equation}

\vspace{2mm}

From (\ref{ogo9050}) and (\ref{ogo9060}) it follows that 
$$
\hbox{\vtop{\offinterlineskip\halign{
\hfil#\hfil\cr
{\rm l.i.m.}\cr
$\stackrel{}{{}_{p_1, p_3\to \infty}}$\cr
}} }
\sum\limits_{j_1=0}^{2p_1}\sum\limits_{j_3=0}^{2p_3}
C_{j_3 j_3 j_1}\zeta_{j_1}^{(i_1)}=
$$
$$
=(T-t)^{3/2}\left(
\frac{1}{6}\zeta_{0}^{(i_1)}+\frac{1}{12}
\zeta_{0}^{(i_1)}+
\hbox{\vtop{\offinterlineskip\halign{
\hfil#\hfil\cr
{\rm l.i.m.}\cr
$\stackrel{}{{}_{p_1\to \infty}}$\cr
}} }
\frac{\sqrt{2}}{4\pi}\sum\limits_{j_1=1}^{p_1}\frac{1}{j_1}
\zeta_{2j_1-1}^{(i_1)}\right)=
$$
$$
=(T-t)^{3/2}\left(
\frac{1}{4}\zeta_{0}^{(i_1)}+
\hbox{\vtop{\offinterlineskip\halign{
\hfil#\hfil\cr
{\rm l.i.m.}\cr
$\stackrel{}{{}_{p_1\to \infty}}$\cr
}} }
\frac{\sqrt{2}}{4\pi}\sum\limits_{j_1=1}^{p_1}\frac{1}{j_1}
\zeta_{2j_1-1}^{(i_1)}\right)=
$$
$$
=
\frac{1}{2}\int\limits_t^T\int\limits_t^{\tau}d{\bf w}_{s}^{(i_1)}d\tau,
$$

\vspace{2mm}
\noindent
where the equality is fulfilled w.~p.~1.

So, the relations (\ref{agenty1110}) and (\ref{ogo1399}) are proved for the case 
of trigonometric system of functions.

Let us prove the equality (\ref{ogo13a99}).
Since $\psi_1(\tau),$ $\psi_2(\tau),$ $\psi_3(\tau)\equiv 1,$
then the following 
relation 
for the Fourier coefficients is correct
$$
C_{j_1 j_1 j_3}+C_{j_1 j_3 j_1}+C_{j_3 j_1 j_1}=\frac{1}{2}
C_{j_1}^2 C_{j_3}.
$$ 

\vspace{2mm}

Then w.~p.~1 
$$
\hbox{\vtop{\offinterlineskip\halign{
\hfil#\hfil\cr
{\rm l.i.m.}\cr
$\stackrel{}{{}_{p_1,p_3\to \infty}}$\cr
}} }\sum_{j_1=0}^{p_1}\sum_{j_3=0}^{p_3}
C_{j_1 j_3 j_1}\zeta_{j_3}^{(i_2)}=
$$
\begin{equation}
\label{ogo2010}
~~~~~~=
\hbox{\vtop{\offinterlineskip\halign{
\hfil#\hfil\cr
{\rm l.i.m.}\cr
$\stackrel{}{{}_{p_1,p_3\to \infty}}$\cr
}} }\sum_{j_1=0}^{p_1}\sum_{j_3=0}^{p_3}
\left(\frac{1}{2}C_{j_1}^2 C_{j_3}-C_{j_1 j_1 j_3}-C_{j_3 j_1 j_1}
\right)\zeta_{j_3}^{(i_2)}.
\end{equation}

\vspace{2mm}

Taking into account (\ref{ogo1299}) and (\ref{ogo1399}), 
we can write w.~p.~1 
$$
\hbox{\vtop{\offinterlineskip\halign{
\hfil#\hfil\cr
{\rm l.i.m.}\cr
$\stackrel{}{{}_{p_1,p_3\to \infty}}$\cr
}} }\sum_{j_1=0}^{p_1}\sum_{j_3=0}^{p_3}
C_{j_1 j_3 j_1}\zeta_{j_3}^{(i_2)}=
$$
$$
=
\frac{1}{2}C_0^3\zeta_0^{(i_2)}
-
\hbox{\vtop{\offinterlineskip\halign{
\hfil#\hfil\cr
{\rm l.i.m.}\cr
$\stackrel{}{{}_{p_1,p_3\to \infty}}$\cr
}} }\sum_{j_1=0}^{p_1}\sum_{j_3=0}^{p_3}
C_{j_1 j_1 j_3}\zeta_{j_3}^{(i_2)}-
$$
$$
-
\hbox{\vtop{\offinterlineskip\halign{
\hfil#\hfil\cr
{\rm l.i.m.}\cr
$\stackrel{}{{}_{p_1,p_3\to \infty}}$\cr
}} }\sum_{j_1=0}^{p_1}\sum_{j_3=0}^{p_3}
C_{j_3 j_1 j_1}\zeta_{j_3}^{(i_2)}=
$$
$$
=
\frac{1}{2}(T-t)^{3/2}
\zeta_0^{(i_2)}-
\frac{1}{4}(T-t)^{3/2}
\Biggl(\zeta_0^{(i_2)}+
\hbox{\vtop{\offinterlineskip\halign{
\hfil#\hfil\cr
{\rm l.i.m.}\cr
$\stackrel{}{{}_{p_1\to \infty}}$\cr
}} }
\frac{\sqrt{2}}{\pi}\sum_{j_1=1}^{p_1}
\frac{1}{j_1}
\zeta_{2j_1-1}^{(i_2)}
\Biggr)
-
$$
$$
-\frac{1}{4}(T-t)^{3/2}
\Biggl(\zeta_0^{(i_2)}-\hbox{\vtop{\offinterlineskip\halign{
\hfil#\hfil\cr
{\rm l.i.m.}\cr
$\stackrel{}{{}_{p_1\to \infty}}$\cr
}} }\frac{\sqrt{2}}{\pi}
\sum_{j_1=1}^{p_1}
\frac{1}{j_1}
\zeta_{2j_1-1}^{(i_2)}
\Biggr)=0.
$$

\vspace{2mm}

From Theorem 1.1 and (\ref{ogo1299})--(\ref{ogo13a99}) 
we obtain the expansion
(\ref{feto19001ee}). Theorem 2.6 is proved.

\subsection{The Case $p_1=p_2=p_3\to \infty,$ Smooth 
Weight Functions, and Additional
Restrictive Conditions (The Cases of Legendre 
Polynomials and Trigonometric Functions)}

Let us consider the following modification of Theorem 2.5.

{\bf Theorem 2.7}\ \cite{9}-\cite{12aa}, \cite{arxiv-7}. 
{\it Assume that
$\{\phi_j(x)\}_{j=0}^{\infty}$ is a complete orthonormal
system of Legendre polynomials or trigonometric functions
in the space $L_2([t, T])$
and $\psi_1(s),$ $\psi_2(s),$ $\psi_3(s)$ are continuously
differentiable functions at the interval $[t, T]$.
Then$,$ for the 
iterated Stratonovich stochastic integral of third multiplicity
$$
J^{*}[\psi^{(3)}]_{T,t}={\int\limits_t^{*}}^T\psi_3(t_3)
{\int\limits_t^{*}}^{t_3}\psi_2(t_2)
{\int\limits_t^{*}}^{t_2}\psi_1(t_1)
d{\bf w}_{t_1}^{(i_1)}
d{\bf w}_{t_2}^{(i_2)}d{\bf w}_{t_3}^{(i_3)},
$$
where $i_1, i_2, i_3=1,\ldots,m,$ the following 
expansion 
\begin{equation}
\label{feto19000}
J^{*}[\psi^{(3)}]_{T,t}=
\hbox{\vtop{\offinterlineskip\halign{
\hfil#\hfil\cr
{\rm l.i.m.}\cr
$\stackrel{}{{}_{p\to \infty}}$\cr
}} }
\sum\limits_{j_1, j_2, j_3=0}^{p}
C_{j_3 j_2 j_1}\zeta_{j_1}^{(i_1)}\zeta_{j_2}^{(i_2)}\zeta_{j_3}^{(i_3)}
\end{equation}
that converges in the mean-square sense 
is valid for each of the following cases

\noindent
{\rm 1}.\ $i_1\ne i_2,\ i_2\ne i_3,\ i_1\ne i_3,$\\
{\rm 2}.\ $i_1=i_2\ne i_3$ and
$\psi_1(s)\equiv\psi_2(s),$\\
{\rm 3}.\ $i_1\ne i_2=i_3$ and
$\psi_2(s)\equiv\psi_3(s),$\\
{\rm 4}.\ $i_1, i_2, i_3=1,\ldots,m$
and $\psi_1(s)\equiv\psi_2(s)\equiv\psi_3(s),$\

\noindent
where
$$
C_{j_3 j_2 j_1}=\int\limits_t^T\psi_3(s)\phi_{j_3}(s)
\int\limits_t^s\psi_2(s_1)\phi_{j_2}(s_1)
\int\limits_t^{s_1}\psi_1(s_2)\phi_{j_1}(s_2)ds_2ds_1ds
$$
and
$$
\zeta_{j}^{(i)}=
\int\limits_t^T \phi_{j}(s) d{\bf w}_s^{(i)}
$$ 
are independent standard Gaussian random variables for various 
$i$ or $j$.}

{\bf Proof.}\ Let us consider at first the polynomial case. Case 1
directly follows from Theorem 1.1. Further, consider Case 2.
We will prove the following re\-la\-ti\-on
$$
\hbox{\vtop{\offinterlineskip\halign{
\hfil#\hfil\cr
{\rm l.i.m.}\cr
$\stackrel{}{{}_{p\to \infty}}$\cr
}} }
\sum\limits_{j_1=0}^{p}\sum\limits_{j_3=0}^{p}
C_{j_3 j_1 j_1}\zeta_{j_3}^{(i_3)}=
\frac{1}{2}\int\limits_t^T\psi_3(s)
\int\limits_t^s\psi^2(s_1)ds_1d{\bf w}_s^{(i_3)}\ \ \ \hbox{w.~p.~1},
$$
where
$$
C_{j_3 j_1 j_1}=\int\limits_t^T\psi_3(s)\phi_{j_3}(s)
\int\limits_t^s\psi(s_1)\phi_{j_1}(s_1)
\int\limits_t^{s_1}\psi(s_2)\phi_{j_1}(s_2)ds_2ds_1ds.
$$

Using Theorem 
1.1, we can write w.~p.~1 
$$
\frac{1}{2}\int\limits_t^T\psi_3(s)
\int\limits_t^s\psi^2(s_1)ds_1d{\bf w}_s^{(i_3)}=
\frac{1}{2}\ 
\hbox{\vtop{\offinterlineskip\halign{
\hfil#\hfil\cr
{\rm l.i.m.}\cr
$\stackrel{}{{}_{p_3\to \infty}}$\cr
}} }\sum_{j_3=0}^{p_3}
\tilde C_{j_3}\zeta_{j_3}^{(i_3)},
$$
where 
$$
\tilde C_{j_3}=
\int\limits_t^T
\phi_{j_3}(s)\psi_3(s)\int\limits_t^s\psi^2(s_1)ds_1ds.
$$

We have
$$
{\sf M}\left\{\left(\sum_{j_3=0}^p\left(\sum_{j_1=0}^p
C_{j_3j_1j_1}-\frac{1}{2}\tilde C_{j_3}\right)
\zeta_{j_3}^{(i_3)}\right)^2\right\}=
\sum_{j_3=0}^p\left(\sum\limits_{j_1=0}^{p}C_{j_3j_1 j_1}-
\frac{1}{2}\tilde C_{j_3}\right)^2=
$$
$$
=\sum_{j_3=0}^p\left(\frac{1}{2}\sum\limits_{j_1=0}^{p}
\int\limits_t^T\phi_{j_3}(s)\psi_3(s)
\left(\int\limits_t^s\phi_{j_1}(s_1)\psi(s_1)ds_1\right)^2ds-\right.
$$
$$
\left.-
\frac{1}{2}
\int\limits_t^T
\phi_{j_3}(s)\psi_3(s)\int\limits_t^s\psi^2(s_1)ds_1ds\right)^2=
$$
$$
=\frac{1}{4}\sum_{j_3=0}^p\left(
\int\limits_t^T\phi_{j_3}(s)\psi_3(s)\left(
\sum\limits_{j_1=0}^{p}
\left(\int\limits_t^s\phi_{j_1}(s_1)\psi(s_1)ds_1\right)^2
-\int\limits_t^s\psi^2(s_1)ds_1\right)ds\right)^2=
$$
\begin{equation}
\label{otit5000}
~~~~~=\frac{1}{4}\sum_{j_3=0}^p\left(
\int\limits_t^T\phi_{j_3}(s)\psi_3(s)
\sum\limits_{j_1=p+1}^{\infty}
\left(\int\limits_t^s\phi_{j_1}(s_1)\psi(s_1)ds_1\right)^2
ds\right)^2.
\end{equation}

In order to get (\ref{otit5000}) we used the Parseval equality
$$
\sum_{j_1=0}^{\infty}\left(\int\limits_t^s\phi_{j_1}(s_1)
\psi(s_1)ds_1\right)^2=
\int\limits_t^T K^2(s,s_1)ds_1=\int\limits_t^s\psi^2(s_1)ds_1,
$$
where
$$
K(s,s_1)=\psi(s_1){\bf 1}_{\{s_1< s\}},\ \ \ s, s_1\in [t, T].
$$

We have for $j_1\in{\bf N}$
$$
\left(\int\limits_t^s\psi(s_1)\phi_{j_1}(s_1)ds_1\right)^2=
$$
$$
=
\frac{(T-t)(2j_1+1)}{4}
\left(\int\limits_{-1}^{z(s)}P_{j_1}(y)
\psi\left(\frac{T-t}{2}y+\frac{T+t}{2}\right)dy\right)^2=
$$
$$
=\frac{T-t}{4(2j_1+1)}\Biggl(\left(P_{j_1+1}(z(s))-
P_{j_1-1}(z(s))\right)\psi(s)-\Biggr.
$$
\begin{equation}
\label{otit6000}
~~~~~\Biggl.-\frac{T-t}{2}
\int\limits_{-1}^{z(s)}(\left(P_{j_1+1}(y)-P_{j_1-1}(y)\right)
\psi'\left(\frac{T-t}{2}y+\frac{T+t}{2}\right)dy\Biggr)^2,
\end{equation}
where 
$$
z(s)=\left(s-\frac{T+t}{2}\right)\frac{2}{T-t},
$$
and $\psi'$ is a derivative of the function $\psi(s)$
with respect to the variable 
$$
\frac{T-t}{2}y+\frac{T+t}{2}.
$$

Further consideration is similar to
the proof of Case 2 from Theorem 2.5. 
Finally, from (\ref{otit5000}) and (\ref{otit6000})
we obtain
$$
{\sf M}\left\{\left(\sum_{j_3=0}^p\left(\sum_{j_1=0}^p
C_{j_3j_1j_1}-\frac{1}{2}\tilde C_{j_3}\right)
\zeta_{j_3}^{(i_3)}\right)^2\right\}<
$$
$$
< K\frac{p}{p^2}\left(
\int\limits_{-1}^1\frac{dy}{(1-y^2)^{3/4}}
+\int\limits_{-1}^1\frac{dy}{(1-y^2)^{1/4}}\right)^2\le
$$
$$
\le \frac{K_1}{p} \to 0\ \ \  \hbox{if}\ \ \ p \to \infty,
$$

\noindent
where constants $K, K_1$ do not depend on $p$. Case 2 is proved.

Let us consider Case 3. In this case we will prove the following
relation
$$
\hbox{\vtop{\offinterlineskip\halign{
\hfil#\hfil\cr
{\rm l.i.m.}\cr
$\stackrel{}{{}_{p\to \infty}}$\cr
}} }
\sum\limits_{j_1=0}^{p}\sum\limits_{j_3=0}^{p}
C_{j_3 j_3 j_1}\zeta_{j_1}^{(i_1)}=
\frac{1}{2}\int\limits_t^T\psi^2(s)
\int\limits_t^s\psi_1(s_1)d{\bf w}_{s_1}^{(i_1)}ds\ \ \ \hbox{w.~p.~1},
$$
where
$$
C_{j_3 j_3 j_1}=\int\limits_t^T\psi(s)\phi_{j_3}(s)
\int\limits_t^s\psi(s_1)\phi_{j_3}(s_1)
\int\limits_t^{s_1}\psi_1(s_2)\phi_{j_1}(s_2)ds_2ds_1ds.
$$

Using the It\^{o} formula, we obtain w.~p.~1
\begin{equation}
\label{otit7000}
~~~~~~~~ \frac{1}{2}\int\limits_t^T\psi^2(s)
\int\limits_t^s\psi_1(s_1)d{\bf w}_{s_1}^{(i_1)}ds=
\frac{1}{2}\int\limits_t^T\psi_1(s_1)
\int\limits_{s_1}^T\psi^2(s)dsd{\bf w}_{s_1}^{(i_1)}.
\end{equation}

Moreover, using Theorem 
1.1, we have w.~p.~1
\begin{equation}
\label{otit9000}
~~~~~~~\frac{1}{2}\int\limits_t^T\psi_1(s_1)
\int\limits_{s_1}^T\psi^2(s)dsd{\bf w}_{s_1}^{(i_1)}=
\frac{1}{2}\
\hbox{\vtop{\offinterlineskip\halign{
\hfil#\hfil\cr
{\rm l.i.m.}\cr
$\stackrel{}{{}_{p_1\to \infty}}$\cr
}} }\sum_{j_1=0}^{p_1}
C_{j_1}^{*}\zeta_{j_1}^{(i_1)},
\end{equation}
where
$$
C_{j_1}^{*}=
\int\limits_t^T
\phi_{j_1}(s_1)\psi_1(s_1)\int\limits_{s_1}^T\psi^2(s)dsds_1.
$$

Further, 
$$
C_{j_3 j_3 j_1}=\int\limits_t^T\psi(s)\phi_{j_3}(s)
\int\limits_t^s\psi(s_1)\phi_{j_3}(s_1)
\int\limits_t^{s_1}\psi_1(s_2)\phi_{j_1}(s_2)ds_2ds_1ds=
$$
$$
=\int\limits_t^T\psi_1(s_2)\phi_{j_1}(s_2)
\int\limits_{s_2}^T\psi(s_1)\phi_{j_3}(s_1)
\int\limits_{s_1}^{T}\psi(s)\phi_{j_3}(s)dsds_1ds_2=
$$
\begin{equation}
\label{otit8000}
=\frac{1}{2}\int\limits_t^T\psi_1(s_2)\phi_{j_1}(s_2)
\left(\int\limits_{s_2}^T\psi(s_1)\phi_{j_3}(s_1)ds_1\right)^2
ds_2.
\end{equation}

From (\ref{otit7000})--(\ref{otit8000})
we obtain
$$
{\sf M}\left\{\left(\sum_{j_1=0}^p\left(\sum_{j_3=0}^p
C_{j_3j_3j_1}-\frac{1}{2} C_{j_1}^{*}\right)
\zeta_{j_1}^{(i_1)}\right)^2\right\}=
\sum_{j_1=0}^p\left(\sum\limits_{j_3=0}^{p}C_{j_3j_3 j_1}-
\frac{1}{2} C_{j_1}^{*}\right)^2=
$$
$$
=\frac{1}{4}\sum_{j_1=0}^p\left(
\int\limits_t^T\phi_{j_1}(s_1)\psi_1(s_1)\left(
\sum\limits_{j_3=0}^{p}
\left(\int\limits_{s_1}^T\phi_{j_3}(s)\psi(s)ds_1\right)^2
-\right.\right.
$$
$$
\left.\left.-
\int\limits_{s_1}^T\psi^2(s)ds\right)ds_1\right)^2=
$$
\begin{equation}
\label{otit9500}
~~~~~=\frac{1}{4}\sum_{j_1=0}^p\left(
\int\limits_t^T\phi_{j_1}(s_1)\psi_1(s_1)
\sum\limits_{j_3=p+1}^{\infty}
\left(\int\limits_{s_1}^T\phi_{j_3}(s)\psi(s)ds\right)^2
ds_1\right)^2.
\end{equation}

In order to get (\ref{otit9500}) we used the Parseval equality
$$
\sum_{j_3=0}^{\infty}\left(\int\limits_{s_1}^T\phi_{j_3}(s)
\psi(s)ds\right)^2=
\int\limits_t^T K^2(s,s_1)ds=\int\limits_{s_1}^T\psi^2(s)ds,
$$
where
$$
K(s,s_1)=\psi(s){\bf 1}_{\{s> s_1\}},\ \ \ s, s_1\in [t, T].
$$

Further consideration is similar to the proof of Case 3
from Theorem 2.5. 
Finally, from (\ref{otit9500})
we get
$$
{\sf M}\left\{\left(\sum_{j_1=0}^p\left(\sum_{j_3=0}^p
C_{j_3j_3j_1}-\frac{1}{2} C_{j_1}^{*}\right)
\zeta_{j_1}^{(i_1)}\right)^2\right\}<
$$
$$
< K\frac{p}{p^2}\left(
\int\limits_{-1}^1\frac{dy}{(1-y^2)^{3/4}}
+\int\limits_{-1}^1\frac{dy}{(1-y^2)^{1/4}}\right)^2\le
$$
$$
\le \frac{K_1}{p} \to 0\ \ \ \hbox{if}\ \ \  p \to \infty,
$$

\noindent 
where constants $K, K_1$ do not depend on $p$. Case 3 is proved.

Let us consider Case 4. We will prove w.~p.~1 the following
relation
$$
\hbox{\vtop{\offinterlineskip\halign{
\hfil#\hfil\cr
{\rm l.i.m.}\cr
$\stackrel{}{{}_{p\to \infty}}$\cr
}} }
\sum\limits_{j_1=0}^{p}\sum\limits_{j_3=0}^{p}
C_{j_1 j_3 j_1}\zeta_{j_3}^{(i_2)}=0\ \ \ 
(\psi_1(s), \psi_2(s), \psi_3(s)\equiv \psi(s)).
$$

In Case 4 we obtain w.~p.~1
$$
\hbox{\vtop{\offinterlineskip\halign{
\hfil#\hfil\cr
{\rm l.i.m.}\cr
$\stackrel{}{{}_{p\to \infty}}$\cr
}} }
\sum\limits_{j_1, j_3=0}^{p}
C_{j_1 j_3 j_1}\zeta_{j_3}^{(i_2)}=
$$
$$
=\hbox{\vtop{\offinterlineskip\halign{
\hfil#\hfil\cr
{\rm l.i.m.}\cr
$\stackrel{}{{}_{p\to \infty}}$\cr
}} }
\sum\limits_{j_1, j_3=0}^{p}
\left(\frac{1}{2}C_{j_1}^2 C_{j_3}-C_{j_1 j_1 j_3}-C_{j_3 j_1 j_1}
\right)\zeta_{j_3}^{(i_2)}=
$$
$$
=
\hbox{\vtop{\offinterlineskip\halign{
\hfil#\hfil\cr
{\rm l.i.m.}\cr
$\stackrel{}{{}_{p\to \infty}}$\cr
}} }
\frac{1}{2}\sum\limits_{j_1=0}^{p}
C_{j_1}^2\sum\limits_{j_3=0}^{p}C_{j_3}\zeta_{j_3}^{(i_2)}-
\hbox{\vtop{\offinterlineskip\halign{
\hfil#\hfil\cr
{\rm l.i.m.}\cr
$\stackrel{}{{}_{p\to \infty}}$\cr
}} }
\sum\limits_{j_1, j_3=0}^{p}
C_{j_1 j_1 j_3}\zeta_{j_3}^{(i_2)}-
$$
$$
-
\hbox{\vtop{\offinterlineskip\halign{
\hfil#\hfil\cr
{\rm l.i.m.}\cr
$\stackrel{}{{}_{p\to \infty}}$\cr
}} }
\sum\limits_{j_1, j_3=0}^{p}
C_{j_3 j_1 j_1}\zeta_{j_3}^{(i_2)}
=
$$
$$
=\frac{1}{2}\sum\limits_{j_1=0}^{\infty}
C_{j_1}^2\int\limits_t^T\psi(s)d{\bf w}_s^{(i_2)}
-\frac{1}{2}\int\limits_t^T\psi^2(s)
\int\limits_t^s\psi(s_1)d{\bf w}_{s_1}^{(i_2)}ds-
$$
$$
-\frac{1}{2}\int\limits_t^T\psi(s)
\int\limits_t^s\psi^2(s_1)ds_1d{\bf w}_{s}^{(i_2)}=
\frac{1}{2}\int\limits_t^T\psi^2(s)ds
\int\limits_t^T\psi(s)d{\bf w}_s^{(i_2)}-
$$
$$
-\frac{1}{2}\int\limits_t^T\psi(s_1)
\int\limits_{s_1}^T\psi^2(s)dsd{\bf w}_{s_1}^{(i_2)}-
\frac{1}{2}\int\limits_t^T\psi(s_1)
\int\limits_t^{s_1}\psi^2(s)dsd{\bf w}_{s_1}^{(i_2)}=
$$
$$
=\frac{1}{2}\int\limits_t^T\psi^2(s)ds
\int\limits_t^T\psi(s)d{\bf w}_s^{(i_2)}-
\frac{1}{2}\int\limits_t^T\psi(s_1)
\int\limits_t^T\psi^2(s)dsd{\bf w}_{s_1}^{(i_2)}=0,
$$

\noindent
where we used the Parseval equality 
$$
\sum\limits_{j_1=0}^{\infty}
C_{j}^2=
\sum\limits_{j=0}^{\infty}
\left(\int\limits_t^T\psi(s)\phi_j(s)ds\right)^2
=\int\limits_t^T\psi^2(s)ds.
$$

Case 4 and Theorem 2.7 are proved for the case of Legendre polynomials.

Let us consider the trigonometric case.
The complete orthonormal system of trigonometric functions
in the space $L_2([t, T])$ has the following form
$$
\phi_j(\theta)=\frac{1}{\sqrt{T-t}}
\left\{
\begin{matrix}
1,\ & j=0\cr\cr
\sqrt{2}{\rm sin} \left(2\pi r(\theta-t)/(T-t)\right),\ & j=2r-1\cr\cr
\sqrt{2}{\rm cos} \left(2\pi r(\theta-t)/(T-t)\right),\ & j=2r
\end{matrix}
,\right.
$$

\noindent
where $r=1, 2,\ldots $

Integrating by parts, we have
$$
\int\limits_t^{s}\phi_{2r-1}(\theta)\psi(\theta)d\theta=
\frac{\sqrt{2}}{\sqrt{T-t}}
\int\limits_t^{s}
\psi(\theta)\ {\rm sin}\frac{2\pi r(\theta-t)}{T-t}
d\theta=
$$
$$
=\sqrt{\frac{T-t}{2}}\frac{1}{\pi r}\Biggl(
-\psi(s)\ {\rm cos}\frac{2\pi r(s-t)}{T-t}+\psi(t)+\Biggr.
$$
$$
\Biggl.
+\int\limits_t^{s}
\psi'(\theta)\ {\rm cos}\frac{2\pi r(\theta-t)}{T-t}
d\theta\Biggr),
$$

$$
\int\limits_t^{s}\phi_{2r}(\theta)\psi(\theta)d\theta=
\frac{\sqrt{2}}{\sqrt{T-t}}
\int\limits_t^{s}
\psi(\theta)\ {\rm cos}\frac{2\pi r(\theta-t)}{T-t}
d\theta=
$$
$$
=\sqrt{\frac{T-t}{2}}\frac{1}{\pi r}\Biggl(\psi(s)\ {\rm sin}\frac{2\pi r(s-t)}{T-t}
-\Biggr.
$$
$$
\Biggl.-\int\limits_t^{s}\psi'(\theta)\ {\rm sin}\frac{2\pi r(\theta-t)}{T-t}
d\theta\Biggr),
$$

\noindent
where $r=1,2,\ldots$
and $\psi'(\theta)$ is a derivative 
of the function $\psi(\theta)$ with respect to the variable $\theta.$

Then

\vspace{-3mm}
\begin{equation}
\label{agentr100}
\left|\int\limits_t^{s}\phi_{2r-1}(\theta)\psi(\theta)d\theta\right|\le
\frac{C}{r}=\frac{2C}{2r}<\frac{2C}{2r-1},
\end{equation}

\newpage
\noindent
\begin{equation}
\label{agentr101}
\left|\int\limits_t^{s}\phi_{2r}(\theta)\psi(\theta)d\theta\right|\le
\frac{C}{r}=\frac{2C}{2r},
\end{equation}
where constant $C$ does not depend on $r$ ($r=1, 2,\ldots$).

From (\ref{agentr100}), (\ref{agentr101}) we get
\begin{equation}
\label{2017x11}
\left|\int\limits_t^{s}\phi_{j_1}(\theta)\psi(\theta)d\theta\right|\le
\frac{K}{j_1},
\end{equation}
where constant $K$ is independent of $j_1$ ($j_1=1, 2,\ldots$).

Analogously, we obtain
\begin{equation}
\label{2017x12}
\left|\int\limits_s^{T}\phi_{j_1}(\theta)\psi(\theta)d\theta\right|\le
\frac{K}{j_1},
\end{equation}
where constant $K$ does not depend on $j_1$ ($j_1=1, 2,\ldots$).

Using (\ref{otit5000}), (\ref{otit9500}), (\ref{2017x11}), and (\ref{2017x12}),
we get
$$
{\sf M}\left\{\left(\sum_{j_3=0}^p\left(\sum_{j_1=0}^p
C_{j_3j_1j_1}-\frac{1}{2}\tilde C_{j_3}\right)
\zeta_{j_3}^{(i_3)}\right)^2\right\}\le \frac{K_1}{p} \to 0\ \ \
\hbox{if}\ \ \  p \to \infty,
$$
$$
{\sf M}\left\{\left(\sum_{j_1=0}^p\left(\sum_{j_3=0}^p
C_{j_3j_3j_1}-\frac{1}{2} C_{j_1}^{*}\right)
\zeta_{j_1}^{(i_1)}\right)^2\right\}
\le \frac{K_1}{p} \to 0\ \ \ \hbox{if}\ \ \ p \to \infty,
$$

\noindent
where constant $K_1$ is independent of $p.$

The consideration of Case 4 is similar to the case of Legendre 
polynomials. Theorem 2.7 is proved.

In the next section, an analogue of Theorem 2.7 will be proved 
without the restrictions 1--4 (see the formulation of
Theorem 2.7).

\subsection{The Case  $p_1=p_2=p_3\to \infty,$ Smooth 
Weight Functions, and without Additional
Restrictive Conditions (The Cases of Legendre 
Polynomials and Trigonometric Functions)}

{\bf Theorem 2.8}\ \cite{9}-\cite{12aa}, \cite{art-5}, \cite{arxiv-5}.
{\it Suppose that 
$\{\phi_j(x)\}_{j=0}^{\infty}$ is a complete orthonormal system of 
Legendre polynomials or trigonometric functions in the space $L_2([t, T]).$
At the same time $\psi_2(s)$ is a continuously dif\-ferentiable 
nonrandom function on $[t, T]$ and $\psi_1(s),$ $\psi_3(s)$ are twice
continuously differentiable nonrandom functions on $[t, T]$. 
Then$,$ for the 
iterated Stratonovich stochastic integral of third multiplicity
$$
J^{*}[\psi^{(3)}]_{T,t}={\int\limits_t^{*}}^T\psi_3(t_3)
{\int\limits_t^{*}}^{t_3}\psi_2(t_2)
{\int\limits_t^{*}}^{t_2}\psi_1(t_1)
d{\bf w}_{t_1}^{(i_1)}
d{\bf w}_{t_2}^{(i_2)}d{\bf w}_{t_3}^{(i_3)},
$$
where $i_1, i_2, i_3=1,\ldots,m,$ the following 
expansion 
\begin{equation}
\label{feto19000a}
J^{*}[\psi^{(3)}]_{T,t}
=\hbox{\vtop{\offinterlineskip\halign{
\hfil#\hfil\cr
{\rm l.i.m.}\cr
$\stackrel{}{{}_{p\to \infty}}$\cr
}} }
\sum\limits_{j_1, j_2, j_3=0}^{p}
C_{j_3 j_2 j_1}\zeta_{j_1}^{(i_1)}\zeta_{j_2}^{(i_2)}\zeta_{j_3}^{(i_3)}
\end{equation}
that converges in the mean-square sense is valid, where
$$
C_{j_3 j_2 j_1}=\int\limits_t^T\psi_3(s)\phi_{j_3}(s)
\int\limits_t^s\psi_2(s_1)\phi_{j_2}(s_1)
\int\limits_t^{s_1}\psi_1(s_2)\phi_{j_1}(s_2)ds_2ds_1ds
$$
and
$$
\zeta_{j}^{(i)}=
\int\limits_t^T \phi_{j}(s) d{\bf w}_s^{(i)}
$$ 
are independent standard Gaussian random variables for various 
$i$ or $j$.}

{\bf Proof.} Let us consider the case of
Legendre polynomials. 
From (\ref{a3}) for the case $p_1=p_2=p_3=p$ and standard
relations between It\^{o} and Stratonovich stochastic integrals 
we conclude that Theorem 2.8 will be proved if w.~p.~1
\begin{equation}
\label{1xx}
~~~~~~\hbox{\vtop{\offinterlineskip\halign{
\hfil#\hfil\cr
{\rm l.i.m.}\cr
$\stackrel{}{{}_{p\to \infty}}$\cr
}} }
\sum\limits_{j_1=0}^{p}\sum\limits_{j_3=0}^{p}
C_{j_3 j_1 j_1}\zeta_{j_3}^{(i_3)}=
\frac{1}{2}\int\limits_t^T\psi_3(s)
\int\limits_t^s\psi_2(s_1)\psi_1(s_1)ds_1d{\bf w}_s^{(i_3)},
\end{equation}
\begin{equation}
\label{2xx}
~~~~~~\hbox{\vtop{\offinterlineskip\halign{
\hfil#\hfil\cr
{\rm l.i.m.}\cr
$\stackrel{}{{}_{p\to \infty}}$\cr
}} }
\sum\limits_{j_1=0}^{p}\sum\limits_{j_3=0}^{p}
C_{j_3 j_3 j_1}\zeta_{j_1}^{(i_1)}=
\frac{1}{2}\int\limits_t^T\psi_3(s)\psi_2(s)
\int\limits_t^s\psi_1(s_1)d{\bf w}_{s_1}^{(i_1)}ds,
\end{equation}
\begin{equation}
\label{3xx}
\hbox{\vtop{\offinterlineskip\halign{
\hfil#\hfil\cr
{\rm l.i.m.}\cr
$\stackrel{}{{}_{p\to \infty}}$\cr
}} }
\sum\limits_{j_1=0}^{p}\sum\limits_{j_3=0}^{p}
C_{j_1 j_3 j_1}\zeta_{j_3}^{(i_2)}=0.
\end{equation}

Let us prove (\ref{1xx}).
Using Theorem 
1.1 for $k=1$ (also see (\ref{a1})), we can write w.~p.~1
$$
\frac{1}{2}\int\limits_t^T\psi_3(s)
\int\limits_t^s\psi_2(s_1)\psi_1(s_1)ds_1d{\bf w}_s^{(i_3)}=
\frac{1}{2}\
\hbox{\vtop{\offinterlineskip\halign{
\hfil#\hfil\cr
{\rm l.i.m.}\cr
$\stackrel{}{{}_{p\to \infty}}$\cr
}} }
\sum\limits_{j_3=0}^{p}
\tilde C_{j_3}\zeta_{j_3}^{(i_3)},
$$
where 
$$
\tilde C_{j_3}=
\int\limits_t^T
\phi_{j_3}(s)\psi_3(s)\int\limits_t^s\psi_2(s_1)\psi_1(s_1)ds_1ds.
$$

We have
$$
E_p\stackrel{\sf def}{=}{\sf M}\left\{\left(
\sum\limits_{j_1=0}^{p}\sum\limits_{j_3=0}^{p}
C_{j_3 j_1 j_1}\zeta_{j_3}^{(i_3)} - 
\frac{1}{2}\sum\limits_{j_3=0}^{p}
\tilde C_{j_3}\zeta_{j_3}^{(i_3)}\right)^2\right\}=
$$
$$
={\sf M}\left\{\left(\sum_{j_3=0}^p\left(\sum_{j_1=0}^p
C_{j_3j_1j_1}-\frac{1}{2}\tilde C_{j_3}\right)
\zeta_{j_3}^{(i_3)}\right)^2\right\}=
$$
$$
=
\sum_{j_3=0}^p\left(\sum\limits_{j_1=0}^{p}C_{j_3j_1 j_1}-
\frac{1}{2}\tilde C_{j_3}\right)^2=
$$
$$
=\sum_{j_3=0}^p\left(\sum\limits_{j_1=0}^{p}
\int\limits_t^T\psi_3(s)\phi_{j_3}(s)
\int\limits_t^s\psi_2(s_1)\phi_{j_1}(s_1)
\int\limits_t^{s_1}\psi_1(s_2)\phi_{j_1}(s_2)
ds_2 ds_1 ds -\right.
$$
$$
\left.-
\frac{1}{2}
\int\limits_t^T
\psi_3(s)\phi_{j_3}(s)\int\limits_t^s\psi_1(s_1)\psi_2(s_1)ds_1ds\right)^2=
$$
$$
=\sum_{j_3=0}^p\left(
\int\limits_t^T\psi_3(s)\phi_{j_3}(s)
\int\limits_t^s\left(
\sum\limits_{j_1=0}^{p}
\psi_2(s_1)\phi_{j_1}(s_1)\times 
\right.\right.
$$
\begin{equation}
\label{otit5000xxx}
\left.\left.
\times \int\limits_t^{s_1}\psi_1(s_2)\phi_{j_1}(s_2)
ds_2- \frac{1}{2}
\psi_1(s_1)\psi_2(s_1)\right)ds_1ds\right)^2.
\end{equation}

Let us substitute $t_1=t_2=s_1$ into (\ref{leto8001yes1}).
Then for all $s_1\in (t, T)$
\begin{equation}
\label{4xx}
~~~~~~~~~~ \sum\limits_{j_1=0}^{\infty}
\psi_2(s_1)\phi_{j_1}(s_1)
\int\limits_t^{s_1}\psi_1(s_2)\phi_{j_1}(s_2)ds_2=
\frac{1}{2}\psi_1(s_1)\psi_2(s_1).
\end{equation}

From (\ref{otit5000xxx}) and (\ref{4xx}) it follows that
\begin{equation}
\label{otit5000x}
E_p
=\sum_{j_3=0}^p\left(
\int\limits_t^T\psi_3(s)\phi_{j_3}(s)
\int\limits_t^s
\sum\limits_{j_1=p+1}^{\infty}
\psi_2(s_1)\phi_{j_1}(s_1)
\int\limits_t^{s_1}\psi_1(s_2)\phi_{j_1}(s_2)
ds_2 ds_1ds\right)^2.
\end{equation}

Applying (\ref{otit5000x}) and (\ref{otit2007}), we obtain
$$
E_p< C_1 \sum\limits_{j_3=0}^p \left(
\int\limits_t^T |\phi_{j_3}(s)| 
\frac{1}{p} \left(
\int\limits_{-1}^{z(s)}\frac{dy}{(1-y^2)^{1/2}}
+
\int\limits_{-1}^{z(s)}\frac{dy}{(1-y^2)^{1/4}}\right)ds\right)^2 \le
$$
$$
\le \frac{C_2}{p^2} \sum\limits_{j_3=0}^p
\left(\int\limits_t^T |\phi_{j_3}(s)| ds\right)^2 \le
\frac{C_2(T-t)}{p^2} \sum\limits_{j_3=0}^p
\int\limits_t^T \phi_{j_3}^2(s) ds =
\frac{C_3 p}{p^2} \to   0
$$

\noindent
if $p \to \infty$,
where constants $C_1, C_2, C_3$ do not depend on $p$.
The equality (\ref{1xx}) is proved. 

Let us prove (\ref{2xx}).
Using the It\^{o} formula, we have
$$
\frac{1}{2}\int\limits_t^T\psi_3(s)\psi_2(s)
\int\limits_t^s\psi_1(s_1)d{\bf w}_{s_1}^{(i_1)}ds=
\frac{1}{2}\int\limits_t^T\psi_1(s_1)
\int\limits_{s_1}^T\psi_3(s)\psi_2(s)dsd{\bf w}_{s_1}^{(i_1)}\ \ \
\hbox{\rm w.~p.~1}.
$$

Moreover, using Theorem 1.1 for $k=1$ (also see (\ref{a1})), 
we obtain w.~p.~1
$$
\frac{1}{2}\int\limits_t^T\psi_1(s)
\int\limits_{s}^T\psi_3(s_1)\psi_2(s_1)ds_1d{\bf w}_s^{(i_1)}=
\frac{1}{2}\
\hbox{\vtop{\offinterlineskip\halign{
\hfil#\hfil\cr
{\rm l.i.m.}\cr
$\stackrel{}{{}_{p\to \infty}}$\cr
}} }
\sum\limits_{j_1=0}^{p}
C_{j_1}^{*}\zeta_{j_1}^{(i_1)},
$$
where 
\begin{equation}
\label{19xx}
C_{j_1}^{*}=
\int\limits_t^T
\psi_1(s)\phi_{j_1}(s)\int\limits_{s}^T\psi_3(s_1)\psi_2(s_1)ds_1ds.
\end{equation}

We have
$$
E_p'\stackrel{\sf def}{=}{\sf M}\left\{\left(
\sum\limits_{j_1=0}^{p}\sum\limits_{j_3=0}^{p}
C_{j_3 j_3 j_1}\zeta_{j_1}^{(i_1)} - 
\frac{1}{2}\sum\limits_{j_1=0}^{p}
C_{j_1}^{*}\zeta_{j_1}^{(i_1)}\right)^2\right\}=
$$
$$
={\sf M}\left\{\left(\sum_{j_1=0}^p\left(\sum_{j_3=0}^p
C_{j_3j_3j_1}-\frac{1}{2}C_{j_1}^{*}\right)
\zeta_{j_1}^{(i_1)}\right)^2\right\}
=
$$
\begin{equation}
\label{20xx}
=\sum_{j_1=0}^p\left(\sum\limits_{j_3=0}^{p}C_{j_3j_3 j_1}-
\frac{1}{2}C_{j_1}^{*}\right)^2,
\end{equation}
$$
C_{j_3 j_3 j_1}=\int\limits_t^T\psi_3(s)\phi_{j_3}(s)
\int\limits_t^s\psi_2(s_1)\phi_{j_3}(s_1)
\int\limits_t^{s_1}\psi_1(s_2)\phi_{j_1}(s_2)ds_2ds_1ds=
$$
\begin{equation}
\label{21xx}
~~~~~~~~~=\int\limits_t^T\psi_1(s_2)\phi_{j_1}(s_2)
\int\limits_{s_2}^T\psi_2(s_1)\phi_{j_3}(s_1)
\int\limits_{s_1}^T\psi_3(s)\phi_{j_3}(s)dsds_1ds_2.
\end{equation}

From (\ref{19xx})--(\ref{21xx}) we obtain
$$
E_p'
=\sum_{j_1=0}^p\left(
\int\limits_t^T\psi_1(s_2)\phi_{j_1}(s_2)
\int\limits_{s_2}^T\left(
\sum\limits_{j_3=0}^{p}
\psi_2(s_1)\phi_{j_3}(s_1)\times 
\right.\right.
$$
\begin{equation}
\label{otit5000xx}
\left.\left.\times \int\limits_{s_1}^T\psi_3(s)\phi_{j_3}(s)ds- \frac{1}{2}
\psi_3(s_1)\psi_2(s_1)\right)ds_1ds_2\right)^2.
\end{equation}

We will prove the following equality for all $s_1\in (t, T)$
\begin{equation}
\label{44xx}
~~~~~~~~~~\sum\limits_{j_3=0}^{\infty}
\psi_2(s_1)\phi_{j_3}(s_1)
\int\limits_{s_1}^T\psi_3(s)\phi_{j_3}(s)ds=
\frac{1}{2}\psi_2(s_1)\psi_3(s_1).
\end{equation}

Let us denote 
\begin{equation}
\label{yes2002x}
K_1^{*}(t_1,t_2)=K_1(t_1,t_2)+\frac{1}{2}{\bf 1}_{\{t_1=t_2\}}
\psi_2(t_1)\psi_3(t_1),
\end{equation}
where
$$
K_1(t_1,t_2)=\psi_2(t_1)\psi_3(t_2){\bf 1}_{\{t_1<t_2\}},\ \ \
t_1, t_2\in[t, T].
$$

\vspace{2mm}

Let us expand the function $K_1^{*}(t_1,t_2)$ using the variable 
$t_2$, when $t_1$ is fixed, into the Fourier--Legendre series 
at the interval $(t, T)$
\begin{equation}
\label{leto8001yesxx}
~~~~~~~~ K_1^{*}(t_1,t_2)=
\sum_{j_3=0}^{\infty}
\psi_2(t_1)
\int\limits_{t_1}^T\psi_3(t_2)\phi_{j_3}(t_2)dt_2\cdot
\phi_{j_3}(t_2)\ \ \ (t_2\ne t, T).
\end{equation}

The equality (\ref{leto8001yesxx}) is 
fulfilled in each point
of the interval $(t, T)$ with respect to the 
variable $t_2$, when $t_1\in [t, T]$ is fixed, due to 
piecewise
smoothness of the function $K_1^{*}(t_1,t_2)$ with respect to the variable 
$t_2\in [t, T]$ ($t_1$ is fixed).

Obtaining (\ref{leto8001yesxx}), we also used the fact that the 
right-hand side 
of (\ref{leto8001yesxx}) converges when $t_1=t_2$ (point of a finite 
discontinuity
of the function $K_1(t_1,t_2)$) to the value
$$
\frac{1}{2}\left(K_1(t_1,t_1-0)+K_1(t_1,t_1+0)\right)=
\frac{1}{2}\psi_2(t_1)\psi_3(t_1)=
K_1^{*}(t_1,t_1).
$$

Let us substitute $t_1=t_2$ into (\ref{leto8001yesxx}).
Then we have (\ref{44xx}).
From (\ref{otit5000xx}) and (\ref{44xx}) we get
\begin{equation}
\label{otit5000xy}
E_p'=\sum_{j_1=0}^p\left(
\int\limits_t^T\psi_1(s_2)\phi_{j_1}(s_2)
\int\limits_{s_2}^T
\sum\limits_{j_3=p+1}^{\infty}
\psi_2(s_1)\phi_{j_3}(s_1)
\int\limits_{s_1}^T\psi_3(s)\phi_{j_3}(s)
ds ds_1ds_2\right)^2.
\end{equation}

Analogously with (\ref{otit2007}) we obtain for
the twice continuously differentiable func\-tion $\psi_3(s)$
the following estimate
$$
\Biggl|
\sum\limits_{j_3=p+1}^{\infty}
\phi_{j_3}(s_1)
\int\limits_{s_1}^T\psi_3(s)\phi_{j_3}(s)
ds\Biggr| 
< 
$$
\begin{equation}
\label{55xx}
<\frac{C}{p}\Biggl(
\frac{1}{(1-(z(s_1))^2)^{1/2}}+
\frac{1}{(1-(z(s_1))^2)^{1/4}}\Biggr),
\end{equation}

\noindent
where $s_1\in (t, T)$, $z(s_1)$ is defined by (\ref{zz1}), and
constant $C$ does not depend on $p$.

Further consideration is analogously to the proof of
(\ref{1xx}). The relation (\ref{2xx})
is proved. 

Let us prove (\ref{3xx}). We have
\begin{equation}
\label{66xx}
~~~~~~~E_p''\stackrel{\sf def}{=}{\sf M}\left\{\left(
\sum\limits_{j_1=0}^{p}\sum\limits_{j_3=0}^{p}
C_{j_1 j_3 j_1}\zeta_{j_3}^{(i_2)}\right)^2\right\} =
\sum\limits_{j_3=0}^{p}\left(\sum\limits_{j_1=0}^{p}
C_{j_1 j_3 j_1}\right)^2,
\end{equation}
$$
C_{j_1 j_3 j_1}=\int\limits_t^T\psi_3(s)\phi_{j_1}(s)
\int\limits_t^s\psi_2(s_1)\phi_{j_3}(s_1)
\int\limits_t^{s_1}\psi_1(s_2)\phi_{j_1}(s_2)ds_2ds_1ds=
$$
\begin{equation}
\label{22xxx}
~~~~~~=\int\limits_t^T\psi_2(s_1)\phi_{j_3}(s_1)
\int\limits_t^{s_1}\psi_1(s_2)\phi_{j_1}(s_2)ds_2
\int\limits_{s_1}^T\psi_3(s)\phi_{j_1}(s)dsds_1.
\end{equation}

After substituting (\ref{22xxx}) into (\ref{66xx}), we obtain
\begin{equation}
\label{otit5000xyz}
E_p''=\sum_{j_3=0}^p\left(
\int\limits_t^T\psi_2(s_1)\phi_{j_3}(s_1)
\sum\limits_{j_1=0}^{p}\int\limits_{t}^{s_1}
\psi_1(\theta)\phi_{j_1}(\theta)d\theta
\int\limits_{s_1}^T\psi_3(s)\phi_{j_1}(s)
dsds_1\right)^2.
\end{equation}

The generalized Parseval equality gives
$$
\sum\limits_{j_1=0}^{\infty}\int\limits_{t}^{s_1}
\psi_1(\theta)\phi_{j_1}(\theta)d\theta
\int\limits_{s_1}^T\psi_3(s)\phi_{j_1}(s)
ds=
$$
$$
=\sum\limits_{j_1=0}^{\infty}\int\limits_{t}^{T}
{\bf 1}_{\{\theta<s_1\}}
\psi_1(\theta)\phi_{j_1}(\theta)d\theta
\int\limits_{t}^T  {\bf 1}_{\{s>s_1\}}\psi_3(s)\phi_{j_1}(s)
ds=
$$
\begin{equation}
\label{dwdw1}
=
\int\limits_{t}^{T}{\bf 1}_{\{\tau<s_1\}}
\psi_1(\tau)
{\bf 1}_{\{\tau>s_1\}}\psi_3(\tau)d\tau=0.
\end{equation}

\vspace{1mm}

Using (\ref{otit5000xyz}) and (\ref{dwdw1}), we get
\begin{equation}
\label{dwdw2}
E_p''=\sum_{j_3=0}^p\left(
\int\limits_t^T\psi_2(s_1)\phi_{j_3}(s_1)
\sum\limits_{j_1=p+1}^{\infty}\int\limits_{t}^{s_1}
\psi_1(\theta)\phi_{j_1}(\theta)d\theta
\int\limits_{s_1}^T\psi_3(s)\phi_{j_1}(s)
dsds_1\right)^2.
\end{equation}

Let us write the following relation
$$
\int\limits_t^x\psi_1(s)\phi_{j_1}(s)ds=
\frac{\sqrt{T-t}\sqrt{2j_1+1}}{2}
\int\limits_{-1}^{z(x)}P_{j_1}(y)
\psi_1(u(y))dy=
$$
$$
=\frac{\sqrt{T-t}}{2\sqrt{2j_1+1}}\Biggl((P_{j_1+1}(z(x))-
P_{j_1-1}(z(x)))\psi_1(x)-\Biggr.
$$
\begin{equation}
\label{otit6000x}
\Biggl.-
\frac{T-t}{2}
\int\limits_{-1}^{z(x)}((P_{j_1+1}(y)-P_{j_1-1}(y))
{\psi_1}'(u(y))dy\Biggr),
\end{equation}
where $x\in (t, T),$ $j_1\ge p+1,$ 
$z(x)$ and $u(y)$ are defined by (\ref{zz1}),
${\psi_1}'$ is a derivative of the function $\psi_1(s)$
with respect to the variable $u(y).$

Note that in (\ref{otit6000x}) we used the following well known property
of the Legendre polynomials \cite{suet}
$$
P_{j+1}(-1)=-P_j(-1),\ \ \ j=0, 1, 2, \ldots
$$ 
and (\ref{w1ggg}).

From (\ref{ogo23}) and (\ref{otit6000x}) we obtain
\begin{equation}
\label{101xx}
\left|
\int\limits_t^x\psi_1(s)\phi_{j_1}(s)ds
\right| <
\frac{C}{j_1}\Biggl(\frac{1}{(1-(z(x))^2)^{1/4}}+C_1\Biggr),
\end{equation}
where $j_1\in{\bf N},$ $x\in (t, T),$ constants $C, C_1$ do not depend on $j_1$.

Similarly to (\ref{101xx}) and due to 
$$
P_j(1)=1,\ \ \ j=0, 1, 2,\ldots
$$ 
we obtain an analogue of (\ref{101xx}) for the integral, which
is similar to the integral 
on the left-hand side of (\ref{101xx}), but with integration limits 
$x$ and $T$.

From the formula (\ref{101xx}) and its analogue for 
the integral with integration limits
$x$ and $T$ we obtain
\begin{equation}
\label{103xx}
~~~~~~~~~ \left|
\int\limits_t^x\psi_1(s)\phi_{j_1}(s)ds
\int\limits_x^T\psi_3(s)\phi_{j_1}(s)ds
\right| <
\frac{K}{j_1^2}\Biggl(\frac{1}{(1-(z(x))^2)^{1/2}}+K_1\Biggr),\
\end{equation}
where  $j_1\in{\bf N},$ $x\in (t, T),$ and constants $K, K_1$ do not depend on $j_1$.

Let us estimate the right-hand side of (\ref{dwdw2}) using (\ref{103xx})

\vspace{-2mm}
$$
E_p''\le
$$
$$
\le
L\sum_{j_3=0}^p\left(
\int\limits_t^T|\phi_{j_3}(s_1)|
\sum\limits_{j_1=p+1}^{\infty}\left|\int\limits_{t}^{s_1}
\psi_1(\theta)\phi_{j_1}(\theta)d\theta
\int\limits_{s_1}^T\psi_3(s)\phi_{j_1}(s)
ds\right|ds_1\right)^2<
$$
$$
< L_1
\sum_{j_3=0}^p\left(
\int\limits_t^T|\phi_{j_3}(s_1)|
\sum\limits_{j_1=p+1}^{\infty}\frac{1}{j_1^2}
\Biggl(\frac{1}{(1-(z(s_1))^2)^{1/2}}+K_1\Biggr)
ds_1\right)^2<
$$
$$
<
\frac{L_2}{p^2}
\sum_{j_3=0}^p
\left(\int\limits_t^T\frac{ds_1}{(1-(z(s_1))^2)^{3/4}}+
K_1\int\limits_t^T\frac{ds_1}{(1-(z(s_1))^2)^{1/4}}
\right)^2=
$$
$$
=
\frac{L_2(T-t)^2}{4 p^2}
\sum_{j_3=0}^p
\left(\int\limits_{-1}^1\frac{dy}{(1-y^2)^{3/4}}+
K_1\int\limits_{-1}^1\frac{dy}{(1-y^2)^{1/4}}
\right)^2
\le 
$$
\begin{equation}
\label{104xx}
\le\frac{L_3 p}{p^2}=\frac{L_3}{p}\to 0
\end{equation}

\vspace{2mm}
\noindent
if $p\to\infty$,
where constants $L, L_1, L_2, L_3$ do not depend on $p$ and we
used (\ref{obana}), (\ref{ogo24}) in (\ref{104xx}).
The relation (\ref{3xx}) is proved. Theorem 2.8 is proved
for the case of Legendre polynomials.

Let us consider the trigonometric case.
Analogously to (\ref{2017zzz111}) we obtain
\begin{equation}
\label{2017zzz1113}
\left|\int\limits_{s_2}^T\sum_{j_3=p+1}^{\infty}
\psi_2(s_1)\phi_{j_3}(s_1)\int\limits_{s_1}^{T}
\psi_3(s)\phi_{j_3}(s)ds ds_1\right|
\le \frac{K_1}{p},
\end{equation}

\vspace{1mm}
\noindent
where $s_2\in (t, T)$ 
and constant $K_1$ does not depend on $p$.

Using (\ref{2017zzz111}) for $T=s$ and
(\ref{otit5000x}), we obtain
$$
E_p\le K
\sum_{j_3=0}^p\left(
\int\limits_t^T
\left|\int\limits_t^s
\sum\limits_{j_1=p+1}^{\infty}
\psi_2(s_1)\phi_{j_1}(s_1)
\int\limits_t^{s_1}\psi_1(s_2)\phi_{j_1}(s_2)
ds_2 ds_1\right| ds\right)^2\le
$$
\begin{equation}
\label{2017abc}
~~~~~~~~\le K
\sum_{j_3=0}^p\left(
(T-t)
\frac{K_1}{p}\right)^2
\le \frac{K_2}{p^2}
\sum_{j_3=0}^p(T-t)^2\le \frac{L}{p} \to 0
\end{equation}

\vspace{1mm}
\noindent
if $p\to\infty$, where constants $K, K_1, K_2, L$ do not depend on $p$.

Analogously, using (\ref{2017zzz1113}) and (\ref{otit5000xy}), we obtain
that $E_p' \to 0$
if $p\to\infty$. 
It is not difficult to see that in our case we have (see (\ref{2017x11}), (\ref{2017x12}))
\begin{equation}
\label{2017era}
\left|
\int\limits_t^x\psi_1(s)\phi_{j_1}(s)ds
\int\limits_x^T\psi_3(s)\phi_{j_1}(s)ds
\right| <
\frac{C_1}{j_1^2},
\end{equation}

\vspace{1mm}
\noindent
where $j_1\in{\bf N},$ constant $C_1$ does not depend on $j_1.$

Using (\ref{dwdw2}) and (\ref{2017era}), we obtain

\vspace{-2mm}
$$
E_p''\le 
$$

\newpage
\noindent
$$
\le L
\sum_{j_3=0}^p\left(
\int\limits_t^T|\phi_{j_3}(s_1)|
\sum\limits_{j_1=p+1}^{\infty}\left|\int\limits_{t}^{s_1}
\psi_1(\theta)\phi_{j_1}(\theta)d\theta
\int\limits_{s_1}^T\psi_3(s)\phi_{j_1}(s)
ds\right|ds_1\right)^2\le
$$
$$
\le L_1
\sum_{j_3=0}^p\left((T-t)
\sum\limits_{j_1=p+1}^{\infty}
\frac{1}{j_1^2}
\right)^2\le \frac{L_1}{p^2}
\sum_{j_3=0}^p(T-t)^2\le 
$$

\begin{equation}
\label{may6000}
\le\frac{L_2}{p} \to 0
\end{equation}

\vspace{1mm}
\noindent
if $p\to\infty$, where constants $L, L_1, L_2$ do not depend on $p.$

Theorem 2.8 is proved for the trigonometric case.
Theorem 2.8 is proved.

\section{Expansion of Iterated Stratonovich Stochastic Integrals of 
Multiplicity 4 Based on Theorem 1.1. 
The Case $p_1=\ldots=p_4\to\infty,\ \psi_1(\tau),\ldots,\psi_4(\tau)\equiv 1$ (Cases of Legendre 
Polynomials and Trigonometric Functions)}

In this section, we will develop the approach to expansion
of iterated Stra\-to\-no\-vich stochatic integrals based on Theorem 1.1
for the stochastic integrals of multiplicity 4.

{\bf Theorem 2.9}\ \cite{8}-\cite{12aa}, \cite{art-5}, \cite{arxiv-5}. 
{\it Suppose that
$\{\phi_j(x)\}_{j=0}^{\infty}$ is a complete orthonormal
system of Legendre polynomials or trigonometric functions
in the space $L_2([t, T])$.
Then$,$ for the iterated 
Stratonovich stochastic integral of fourth multiplicity
$$
J^{*}[\psi^{(4)}]_{T,t}=
{\int\limits_t^{*}}^T
{\int\limits_t^{*}}^{t_4}
{\int\limits_t^{*}}^{t_3}
{\int\limits_t^{*}}^{t_2}
d{\bf w}_{t_1}^{(i_1)}
d{\bf w}_{t_2}^{(i_2)}d{\bf w}_{t_3}^{(i_3)}d{\bf w}_{t_4}^{(i_4)}\ \ \ 
(i_1, i_2, i_3, i_4=0, 1,\ldots,m)
$$
the following 
expansion 
\begin{equation}
\label{feto1900otit}
~~~~~~~ J^{*}[\psi^{(4)}]_{T,t}=
\hbox{\vtop{\offinterlineskip\halign{
\hfil#\hfil\cr
{\rm l.i.m.}\cr
$\stackrel{}{{}_{p\to \infty}}$\cr
}} }
\sum\limits_{j_1, j_2, j_3, j_4=0}^{p}
C_{j_4 j_3 j_2 j_1}\zeta_{j_1}^{(i_1)}\zeta_{j_2}^{(i_2)}\zeta_{j_3}^{(i_3)}
\zeta_{j_4}^{(i_4)}
\end{equation}
that converges in the mean-square sense is valid, where
$$
C_{j_4 j_3 j_2 j_1}=\int\limits_t^T\phi_{j_4}(s_4)\int\limits_t^{s_4}
\phi_{j_3}(s_3)
\int\limits_t^{s_3}\phi_{j_2}(s_2)\int\limits_t^{s_2}\phi_{j_1}(s_1)
ds_1ds_2ds_3ds_4
$$
and
$$
\zeta_{j}^{(i)}=
\int\limits_t^T \phi_{j}(s) d{\bf w}_s^{(i)}
$$ 
are independent standard Gaussian random variables for various 
$i$ or $j$ {\rm (}in the case when $i\ne 0${\rm ),}
${\bf w}_{\tau}^{(i)}$ $(i=1,\ldots,m)$ are independent standard Wiener processes$,$
${\bf w}_{\tau}^{(0)}=\tau.$}

{\bf Proof.} The relation (\ref{a4}) (in the case 
when $p_1=\ldots=p_4=p\to \infty$) implies that
$$
\hbox{\vtop{\offinterlineskip\halign{
\hfil#\hfil\cr
{\rm l.i.m.}\cr
$\stackrel{}{{}_{p\to \infty}}$\cr
}} }
\sum\limits_{j_1, j_2, j_3, j_4=0}^{p}
C_{j_4 j_3 j_2 j_1}\zeta_{j_1}^{(i_1)}\zeta_{j_2}^{(i_2)}\zeta_{j_3}^{(i_3)}
\zeta_{j_4}^{(i_4)}=
J[\psi^{(4)}]_{T,t}+
$$

\vspace{-4mm}
$$
+{\bf 1}_{\{i_1=i_2\ne 0\}}A_1^{(i_3i_4)}
+{\bf 1}_{\{i_1=i_3\ne 0\}}A_2^{(i_2i_4)}+
{\bf 1}_{\{i_1=i_4\ne 0\}}A_3^{(i_2i_3)}+
{\bf 1}_{\{i_2=i_3\ne 0\}}A_4^{(i_1i_4)}+
$$

\vspace{-3mm}
$$
+
{\bf 1}_{\{i_2=i_4\ne 0\}}A_5^{(i_1i_3)}
+{\bf 1}_{\{i_3=i_4\ne 0\}}A_6^{(i_1i_2)}-
{\bf 1}_{\{i_1=i_2\ne 0\}}
{\bf 1}_{\{i_3=i_4\ne 0\}}B_1-
$$

\vspace{-2mm}
\begin{equation}
\label{otiteee}
-{\bf 1}_{\{i_1=i_3\ne 0\}}
{\bf 1}_{\{i_2=i_4\ne 0\}}B_2-
{\bf 1}_{\{i_1=i_4\ne 0\}}
{\bf 1}_{\{i_2=i_3\ne 0\}}B_3,
\end{equation}

\vspace{4mm}
\noindent
where
$J[\psi^{(4)}]_{T,t}$ has the form {\rm (\ref{itoxx})}
for $\psi_1(s),\ldots,\psi_4(s)\equiv 1$ and
$i_1,\ldots,i_4=0, 1,\ldots,m,$
$$
A_1^{(i_3i_4)}=
\hbox{\vtop{\offinterlineskip\halign{
\hfil#\hfil\cr
{\rm l.i.m.}\cr
$\stackrel{}{{}_{p\to \infty}}$\cr
}} }
\sum\limits_{j_4, j_3, j_1=0}^{p}
C_{j_4 j_3 j_1 j_1}\zeta_{j_3}^{(i_3)}
\zeta_{j_4}^{(i_4)},
$$
$$
A_2^{(i_2i_4)}=
\hbox{\vtop{\offinterlineskip\halign{
\hfil#\hfil\cr
{\rm l.i.m.}\cr
$\stackrel{}{{}_{p\to \infty}}$\cr
}} }
\sum\limits_{j_4, j_3, j_2=0}^{p}
C_{j_4 j_3 j_2 j_3}\zeta_{j_2}^{(i_2)}
\zeta_{j_4}^{(i_4)},
$$
$$
A_3^{(i_2i_3)}=
\hbox{\vtop{\offinterlineskip\halign{
\hfil#\hfil\cr
{\rm l.i.m.}\cr
$\stackrel{}{{}_{p\to \infty}}$\cr
}} }
\sum\limits_{j_4, j_3, j_2=0}^{p}
C_{j_4 j_3 j_2 j_4}\zeta_{j_2}^{(i_2)}
\zeta_{j_3}^{(i_3)},
$$
$$
A_4^{(i_1i_4)}=
\hbox{\vtop{\offinterlineskip\halign{
\hfil#\hfil\cr
{\rm l.i.m.}\cr
$\stackrel{}{{}_{p\to \infty}}$\cr
}} }
\sum\limits_{j_4, j_3, j_1=0}^{p}
C_{j_4 j_3 j_3 j_1}\zeta_{j_1}^{(i_1)}
\zeta_{j_4}^{(i_4)},
$$
$$
A_5^{(i_1i_3)}=
\hbox{\vtop{\offinterlineskip\halign{
\hfil#\hfil\cr
{\rm l.i.m.}\cr
$\stackrel{}{{}_{p\to \infty}}$\cr
}} }
\sum\limits_{j_4, j_3, j_1=0}^{p}
C_{j_4 j_3 j_4 j_1}\zeta_{j_1}^{(i_1)}
\zeta_{j_3}^{(i_3)},
$$
$$
A_6^{(i_1i_2)}=
\hbox{\vtop{\offinterlineskip\halign{
\hfil#\hfil\cr
{\rm l.i.m.}\cr
$\stackrel{}{{}_{p\to \infty}}$\cr
}} }
\sum\limits_{j_3, j_2, j_1=0}^{p}
C_{j_3 j_3 j_2 j_1}\zeta_{j_1}^{(i_1)}
\zeta_{j_2}^{(i_2)},
$$
$$
B_1=
\hbox{\vtop{\offinterlineskip\halign{
\hfil#\hfil\cr
{\rm lim}\cr
$\stackrel{}{{}_{p\to \infty}}$\cr
}} }
\sum\limits_{j_1, j_4=0}^{p}
C_{j_4 j_4 j_1 j_1},\ \ \
B_2=
\hbox{\vtop{\offinterlineskip\halign{
\hfil#\hfil\cr
{\rm lim}\cr
$\stackrel{}{{}_{p\to \infty}}$\cr
}} }
\sum\limits_{j_4, j_3=0}^{p}
C_{j_3 j_4 j_3 j_4},
$$
$$
B_3=
\hbox{\vtop{\offinterlineskip\halign{
\hfil#\hfil\cr
{\rm lim}\cr
$\stackrel{}{{}_{p\to \infty}}$\cr
}} }
\sum\limits_{j_4, j_3=0}^{p}
C_{j_4 j_3 j_3 j_4}.
$$

Using the integration order replacement in Riemann integrals,
Theorem 1.1 for $k=2$ (see (\ref{a2})) and
(\ref{5t}), 
Parseval's equality and the integration order replacement
technique for It\^{o} stochastic integrals (see Chapter 3) 
\cite{1}-\cite{12aa}, \cite{old-art-2}, \cite{vini},
\cite{arxiv-25} or It\^{o}'s formula, we obtain
$$
A_1^{(i_3i_4)}=
$$
$$
=
\hbox{\vtop{\offinterlineskip\halign{
\hfil#\hfil\cr
{\rm l.i.m.}\cr
$\stackrel{}{{}_{p\to \infty}}$\cr
}} }
\sum\limits_{j_4, j_3, j_1=0}^{p}
\frac{1}{2}\int\limits_t^T\phi_{j_4}(s)\int\limits_t^s\phi_{j_3}(s_1)
\left(\int\limits_t^{s_1}\phi_{j_1}(s_2)ds_2\right)^2ds_1ds
\zeta_{j_3}^{(i_3)}
\zeta_{j_4}^{(i_4)}=
$$
$$
=\hbox{\vtop{\offinterlineskip\halign{
\hfil#\hfil\cr
{\rm l.i.m.}\cr
$\stackrel{}{{}_{p\to \infty}}$\cr
}} }
\sum\limits_{j_4, j_3=0}^{p}
\frac{1}{2}\int\limits_t^T\phi_{j_4}(s)\int\limits_t^s\phi_{j_3}(s_1)
\sum\limits_{j_1=0}^{p}\left(\int\limits_t^{s_1}
\phi_{j_1}(s_2)ds_2\right)^2ds_1ds
\zeta_{j_3}^{(i_3)}
\zeta_{j_4}^{(i_4)}=
$$
$$
=\hbox{\vtop{\offinterlineskip\halign{
\hfil#\hfil\cr
{\rm l.i.m.}\cr
$\stackrel{}{{}_{p\to \infty}}$\cr
}} }
\hspace{-1mm}\sum\limits_{j_4, j_3=0}^{p}
\frac{1}{2}\int\limits_t^T\phi_{j_4}(s)\int\limits_t^s\phi_{j_3}(s_1)
\hspace{-1mm}\left((s_1-t)-
\sum\limits_{j_1=p+1}^{\infty}\left(\int\limits_t^{s_1}
\phi_{j_1}(s_2)ds_2\right)^{\hspace{-2mm}2}\right)ds_1ds \times
$$
$$
\times
\zeta_{j_3}^{(i_3)}
\zeta_{j_4}^{(i_4)}=
$$
$$
=\hbox{\vtop{\offinterlineskip\halign{
\hfil#\hfil\cr
{\rm l.i.m.}\cr
$\stackrel{}{{}_{p\to \infty}}$\cr
}} }
\sum\limits_{j_4, j_3=0}^{p}
\frac{1}{2}\int\limits_t^T\phi_{j_4}(s)\int\limits_t^s\phi_{j_3}(s_1)
(s_1-t)ds_1ds
\zeta_{j_3}^{(i_3)}
\zeta_{j_4}^{(i_4)} - \Delta_1^{(i_3i_4)}=
$$
$$
=\frac{1}{2}\int\limits_t^T\int\limits_t^s(s_1-t)d{\bf w}_{s_1}^{(i_3)}
d{\bf w}_{s}^{(i_4)}
+
$$
$$
+\frac{1}{2}{\bf 1}_{\{i_3=i_4\ne 0\}}
\lim_{p\to\infty}
\sum\limits_{j_3=0}^{p}
\int\limits_t^T\phi_{j_3}(s)\int\limits_t^s\phi_{j_3}(s_1)(s_1-t)ds_1ds
- \Delta_1^{(i_3i_4)}=
$$
\begin{equation}
\label{otiteee1}
=\frac{1}{2}\int\limits_t^T\int\limits_t^s\int\limits_t^{s_1}ds_2
d{\bf w}_{s_1}^{(i_3)}
d{\bf w}_{s}^{(i_4)}+
\frac{1}{4}{\bf 1}_{\{i_3=i_4\ne 0\}}
\int\limits_t^T(s_1-t)ds_1
- \Delta_1^{(i_3i_4)}\ \ \ \hbox{w.~p.~1,}
\end{equation}

\noindent
where
$$
\Delta_1^{(i_3i_4)}=
\hbox{\vtop{\offinterlineskip\halign{
\hfil#\hfil\cr
{\rm l.i.m.}\cr
$\stackrel{}{{}_{p\to \infty}}$\cr
}} }
\sum\limits_{j_3, j_4=0}^{p}
a_{j_4 j_3}^p \zeta_{j_3}^{(i_3)}
\zeta_{j_4}^{(i_4)},
$$
\begin{equation}
\label{rr1xx}
~~~~~~~~~~a_{j_4 j_3}^p=
\frac{1}{2}\int\limits_t^T\phi_{j_4}(s)\int\limits_t^s\phi_{j_3}(s_1)
\sum\limits_{j_1=p+1}^{\infty}\left(\int\limits_t^{s_1}
\phi_{j_1}(s_2)ds_2\right)^2ds_1ds.
\end{equation}

\vspace{4mm}

Let us consider $A_2^{(i_2i_4)}$ 
$$
A_2^{(i_2i_4)}=
$$
$$
=
\hbox{\vtop{\offinterlineskip\halign{
\hfil#\hfil\cr
{\rm l.i.m.}\cr
$\stackrel{}{{}_{p\to \infty}}$\cr
}} }
\sum\limits_{j_4, j_3, j_2=0}^{p}
\int\limits_t^T\phi_{j_4}(s)
\int\limits_t^s\phi_{j_2}(s_2)
\int\limits_t^{s_2}\phi_{j_3}(s_3)ds_3
\int\limits_{s_2}^s\phi_{j_3}(s_1)ds_1 ds_2 ds
\zeta_{j_2}^{(i_2)}
\zeta_{j_4}^{(i_4)}=
$$
$$
=\hbox{\vtop{\offinterlineskip\halign{
\hfil#\hfil\cr
{\rm l.i.m.}\cr
$\stackrel{}{{}_{p\to \infty}}$\cr
}} }
\sum\limits_{j_4, j_3, j_2=0}^{p}\left(
\frac{1}{2}\int\limits_t^T\phi_{j_4}(s)
\left(\int\limits_t^{s}\phi_{j_3}(s_3)ds_3\right)^2
\int\limits_t^s\phi_{j_2}(s_2)
ds_2ds-\right.
$$
$$
-\frac{1}{2}\int\limits_t^T\phi_{j_4}(s)
\int\limits_t^s\phi_{j_2}(s_2)
\left(\int\limits_t^{s_2}\phi_{j_3}(s_3)ds_3\right)^2
ds_2ds-
$$
$$
\left.-\frac{1}{2}\int\limits_t^T\phi_{j_4}(s)
\int\limits_t^s\phi_{j_2}(s_2)
\left(\int\limits_{s_2}^{s}\phi_{j_3}(s_1)ds_1\right)^2
ds_2ds\right)
\zeta_{j_2}^{(i_2)}
\zeta_{j_4}^{(i_4)}=
$$
$$
=
\hbox{\vtop{\offinterlineskip\halign{
\hfil#\hfil\cr
{\rm l.i.m.}\cr
$\stackrel{}{{}_{p\to \infty}}$\cr
}} }
\sum\limits_{j_4, j_2=0}^{p}\left(
\frac{1}{2}\int\limits_t^T\phi_{j_4}(s)
(s-t)
\int\limits_t^s\phi_{j_2}(s_2)
ds_2ds-\right.
$$
$$
-
\frac{1}{2}\int\limits_t^T\phi_{j_4}(s)
\int\limits_t^s\phi_{j_2}(s_2)
(s_2-t)
ds_2ds-
$$
$$
\left.-\frac{1}{2}\int\limits_t^T\phi_{j_4}(s)
\int\limits_t^s\phi_{j_2}(s_2)
(s-t+t-s_2)
ds_2ds\right)
\zeta_{j_2}^{(i_2)}
\zeta_{j_4}^{(i_4)}-
$$

\vspace{-4mm}
\begin{equation}
\label{otit999}
~~~~~ -\Delta_2^{(i_2i_4)}+\Delta_1^{(i_2i_4)}+\Delta_3^{(i_2i_4)}=
-\Delta_2^{(i_2i_4)}+\Delta_1^{(i_2i_4)}+\Delta_3^{(i_2i_4)}\ \ \
\hbox{w.~p.~1,}
\end{equation}

\noindent
where
$$
\Delta_2^{(i_2i_4)}=
\hbox{\vtop{\offinterlineskip\halign{
\hfil#\hfil\cr
{\rm l.i.m.}\cr
$\stackrel{}{{}_{p\to \infty}}$\cr
}} }
\sum\limits_{j_4, j_2=0}^{p}
b_{j_4 j_2}^p \zeta_{j_2}^{(i_2)}
\zeta_{j_4}^{(i_4)},
$$
$$
\Delta_3^{(i_2i_4)}=
\hbox{\vtop{\offinterlineskip\halign{
\hfil#\hfil\cr
{\rm l.i.m.}\cr
$\stackrel{}{{}_{p\to \infty}}$\cr
}} }
\sum\limits_{j_4, j_2=0}^{p}
c_{j_4 j_2}^p \zeta_{j_2}^{(i_2)}
\zeta_{j_4}^{(i_4)},
$$
\begin{equation}
\label{may001}
~~~~~~~b_{j_4 j_2}^p=
\frac{1}{2}\int\limits_t^T\phi_{j_4}(s)
\sum\limits_{j_3=p+1}^{\infty}\left(\int\limits_t^{s}
\phi_{j_3}(s_1)ds_1\right)^2\int\limits_t^s\phi_{j_2}(s_1)ds_1ds,
\end{equation}
\begin{equation}
\label{may002}
~~~~~~~c_{j_4 j_2}^p=
\frac{1}{2}\int\limits_t^T\phi_{j_4}(s)\int\limits_t^s\phi_{j_2}(s_3)
\sum\limits_{j_3=p+1}^{\infty}\left(\int\limits_{s_3}^{s}
\phi_{j_3}(s_1)ds_1\right)^2ds_3ds.
\end{equation}

\vspace{4mm}

Let us consider $A_5^{(i_1i_3)}$
$$
A_5^{(i_1i_3)}=
$$
$$
=
\hbox{\vtop{\offinterlineskip\halign{
\hfil#\hfil\cr
{\rm l.i.m.}\cr
$\stackrel{}{{}_{p\to \infty}}$\cr
}} }
\sum\limits_{j_4, j_3, j_1=0}^{p}
\int\limits_t^T\phi_{j_1}(s_3)\int\limits_{s_3}^T\phi_{j_4}(s_2)
\int\limits_{s_2}^T\phi_{j_3}(s_1)\int\limits_{s_1}^T\phi_{j_4}(s)
dsds_1ds_2ds_3
\zeta_{j_1}^{(i_1)}
\zeta_{j_3}^{(i_3)}=
$$
$$
=\hbox{\vtop{\offinterlineskip\halign{
\hfil#\hfil\cr
{\rm l.i.m.}\cr
$\stackrel{}{{}_{p\to \infty}}$\cr
}} }
\sum\limits_{j_4, j_3, j_1=0}^{p}
\int\limits_t^T\phi_{j_1}(s_3)\int\limits_{s_3}^T\phi_{j_3}(s_1)
\int\limits_{s_1}^T\phi_{j_4}(s)ds\int\limits_{s_3}^{s_1}\phi_{j_4}(s_2)
ds_2ds_1ds_3
\zeta_{j_1}^{(i_1)}
\zeta_{j_3}^{(i_3)}=
$$
$$
=\hbox{\vtop{\offinterlineskip\halign{
\hfil#\hfil\cr
{\rm l.i.m.}\cr
$\stackrel{}{{}_{p\to \infty}}$\cr
}} }
\sum\limits_{j_4, j_3, j_1=0}^{p}\left(
\frac{1}{2}\int\limits_t^T\phi_{j_1}(s_3)
\left(\int\limits_{s_3}^{T}\phi_{j_4}(s)ds\right)^2
\int\limits_{s_3}^T\phi_{j_3}(s_1)
ds_1ds_3-\right.
$$
$$
-\frac{1}{2}\int\limits_t^T\phi_{j_1}(s_3)
\int\limits_{s_3}^T\phi_{j_3}(s_1)
\left(\int\limits_{s_3}^{s_1}\phi_{j_4}(s_2)ds_2\right)^2
ds_1ds_3-
$$
$$
\left.-\frac{1}{2}\int\limits_t^T\phi_{j_1}(s_3)
\int\limits_{s_3}^T\phi_{j_3}(s_1)
\left(\int\limits_{s_1}^{T}\phi_{j_4}(s)ds\right)^2
ds_1ds_3\right)
\zeta_{j_1}^{(i_1)}
\zeta_{j_3}^{(i_3)}=
$$
$$
=
\hbox{\vtop{\offinterlineskip\halign{
\hfil#\hfil\cr
{\rm l.i.m.}\cr
$\stackrel{}{{}_{p\to \infty}}$\cr
}} }
\sum\limits_{j_3, j_1=0}^{p}\left(
\frac{1}{2}\int\limits_t^T\phi_{j_1}(s_3)
(T-s_3)
\int\limits_{s_3}^T\phi_{j_3}(s_1)
ds_1ds_3-\right.
$$
$$
-\frac{1}{2}\int\limits_t^T\phi_{j_1}(s_3)
\int\limits_{s_3}^T\phi_{j_3}(s_1)
(s_1-s_3)
ds_1ds_3
-
$$
$$
\left.
-\frac{1}{2}\int\limits_t^T\phi_{j_1}(s_3)
\int\limits_{s_3}^T\phi_{j_3}(s_1)
(T-s_1)
ds_1ds_3\right)
\zeta_{j_1}^{(i_1)}
\zeta_{j_3}^{(i_3)}-
$$
\begin{equation}
\label{otit9999}
~~~~ -\Delta_4^{(i_1i_3)}+\Delta_5^{(i_1i_3)}+\Delta_6^{(i_1i_3)}=
-\Delta_4^{(i_1i_3)}+\Delta_5^{(i_1i_3)}+\Delta_6^{(i_1i_3)}\ \ \
\hbox{w.~p.~1,}
\end{equation}

\noindent
where
$$
\Delta_4^{(i_1i_3)}=
\hbox{\vtop{\offinterlineskip\halign{
\hfil#\hfil\cr
{\rm l.i.m.}\cr
$\stackrel{}{{}_{p\to \infty}}$\cr
}} }
\sum\limits_{j_3, j_1=0}^{p}
d_{j_3 j_1}^p \zeta_{j_1}^{(i_1)}
\zeta_{j_3}^{(i_3)},
$$
$$
\Delta_5^{(i_1i_3)}=
\hbox{\vtop{\offinterlineskip\halign{
\hfil#\hfil\cr
{\rm l.i.m.}\cr
$\stackrel{}{{}_{p\to \infty}}$\cr
}} }
\sum\limits_{j_3, j_1=0}^{p}
e_{j_3 j_1}^p \zeta_{j_1}^{(i_1)}
\zeta_{j_3}^{(i_3)},
$$
$$
\Delta_6^{(i_1i_3)}=
\hbox{\vtop{\offinterlineskip\halign{
\hfil#\hfil\cr
{\rm l.i.m.}\cr
$\stackrel{}{{}_{p\to \infty}}$\cr
}} }
\sum\limits_{j_3, j_1=0}^{p}
f_{j_3 j_1}^p \zeta_{j_1}^{(i_1)}
\zeta_{j_3}^{(i_3)},
$$
\begin{equation}
\label{may003}
~~~~~~~~~d_{j_3 j_1}^p=
\frac{1}{2}\int\limits_t^T\phi_{j_1}(s_3)
\sum\limits_{j_4=p+1}^{\infty}\left(\int\limits_{s_3}^{T}
\phi_{j_4}(s)ds\right)^2\int\limits_{s_3}^T\phi_{j_3}(s)dsds_3,
\end{equation}
\begin{equation}
\label{may004}
~~~~~~~~~e_{j_3 j_1}^p=
\frac{1}{2}\int\limits_t^T\phi_{j_1}(s_3)\int\limits_{s_3}^T\phi_{j_3}(s)
\sum\limits_{j_4=p+1}^{\infty}\left(\int\limits_{s_3}^{s}
\phi_{j_4}(s_1)ds_1\right)^2dsds_3,
\end{equation}
$$
f_{j_3 j_1}^p=
\frac{1}{2}\int\limits_t^T\phi_{j_1}(s_3)\int\limits_{s_3}^T\phi_{j_3}(s_2)
\sum\limits_{j_4=p+1}^{\infty}\left(\int\limits_{s_2}^{T}
\phi_{j_4}(s_1)ds_1\right)^2ds_2ds_3=
$$
\begin{equation}
\label{may005}
~~~~~~~~=
\frac{1}{2}\int\limits_t^T\phi_{j_3}(s_2)
\sum\limits_{j_4=p+1}^{\infty}\left(\int\limits_{s_2}^{T}
\phi_{j_4}(s_1)ds_1\right)^2
\int\limits_t^{s_2}\phi_{j_1}(s_3)ds_3ds_2.
\end{equation}

\vspace{4mm}

Moreover, 
$$
A_3^{(i_2i_3)}+A_5^{(i_2i_3)}=
$$
$$
=
\hbox{\vtop{\offinterlineskip\halign{
\hfil#\hfil\cr
{\rm l.i.m.}\cr
$\stackrel{}{{}_{p\to \infty}}$\cr
}} }
\sum\limits_{j_4, j_3, j_2=0}^{p}
\left(C_{j_4j_3j_2j_4}+C_{j_4j_3j_4j_2}\right)
\zeta_{j_2}^{(i_2)}
\zeta_{j_3}^{(i_3)}=
$$
$$
=\hbox{\vtop{\offinterlineskip\halign{
\hfil#\hfil\cr
{\rm l.i.m.}\cr
$\stackrel{}{{}_{p\to \infty}}$\cr
}} }
\sum\limits_{j_4, j_3, j_2=0}^{p}
\int\limits_t^T\phi_{j_4}(s)\int\limits_{t}^s\phi_{j_3}(s_1)
\int\limits_t^{s_1}\phi_{j_2}(s_2)\int\limits_t^{s_1}\phi_{j_4}(s_3)
ds_3ds_2ds_1ds
\zeta_{j_2}^{(i_2)}
\zeta_{j_3}^{(i_3)}=
$$
$$
=\hbox{\vtop{\offinterlineskip\halign{
\hfil#\hfil\cr
{\rm l.i.m.}\cr
$\stackrel{}{{}_{p\to \infty}}$\cr
}} }
\sum\limits_{j_4, j_3, j_2=0}^{p}
\int\limits_t^T\phi_{j_3}(s_1)\int\limits_{t}^{s_1}\phi_{j_2}(s_2)
\int\limits_{t}^{s_1}\phi_{j_4}(s_3)ds_3ds_2\int\limits_{s_1}^{T}\phi_{j_4}(s)
dsds_1
\zeta_{j_2}^{(i_2)}
\zeta_{j_3}^{(i_3)}=
$$
$$
=\hbox{\vtop{\offinterlineskip\halign{
\hfil#\hfil\cr
{\rm l.i.m.}\cr
$\stackrel{}{{}_{p\to \infty}}$\cr
}} }
\sum\limits_{j_4, j_3, j_2=0}^{p}
\left(\int\limits_t^T\phi_{j_3}(s_1)\int\limits_{t}^{s_1}\phi_{j_2}(s_2)
\int\limits_{t}^{T}\phi_{j_4}(s_3)ds_3\int\limits_{s_1}^{T}\phi_{j_4}(s)
dsds_2ds_1-\right.
$$
$$
\left.-
\int\limits_t^T\phi_{j_3}(s_1)\int\limits_{t}^{s_1}\phi_{j_2}(s_2)
\left(\int\limits_{s_1}^{T}\phi_{j_4}(s)ds\right)^{2}
ds_2ds_1\right)
\zeta_{j_2}^{(i_2)}
\zeta_{j_3}^{(i_3)}=
$$
$$
=\hbox{\vtop{\offinterlineskip\halign{
\hfil#\hfil\cr
{\rm l.i.m.}\cr
$\stackrel{}{{}_{p\to \infty}}$\cr
}} }
\sum\limits_{j_3, j_2=0}^{p}
\int\limits_t^T\phi_{j_3}(s_1)\int\limits_{t}^{s_1}\phi_{j_2}(s_2)
\hspace{-1mm}\left((T-s_1)-\sum\limits_{j_4=0}^p
\left(\int\limits_{s_1}^{T}\phi_{j_4}(s_3)ds_3\right)^{\hspace{-1mm}2}\right)
ds_2ds_1\times
$$
\begin{equation}
\label{strange500}
\times\zeta_{j_2}^{(i_2)}
\zeta_{j_3}^{(i_3)}
=
2\Delta_6^{(i_2i_3)}\ \ \ \hbox{w.~p.~1.}
\end{equation}

\vspace{2mm}
Then
\begin{equation}
\label{otit222}
~~~~~~~ A_3^{(i_2i_3)}=2\Delta_6^{(i_2i_3)}-A_5^{(i_2i_3)}=
\Delta_4^{(i_2i_3)}-\Delta_5^{(i_2i_3)}+\Delta_6^{(i_2i_3)}\ \ \
\hbox{w.~p.~1}.
\end{equation}

\vspace{4mm}

Let us consider $A_4^{(i_1i_4)}$
$$
A_4^{(i_1i_4)}=
$$
$$
=
\hbox{\vtop{\offinterlineskip\halign{
\hfil#\hfil\cr
{\rm l.i.m.}\cr
$\stackrel{}{{}_{p\to \infty}}$\cr
}} }
\sum\limits_{j_4, j_3, j_1=0}^{p}
\int\limits_t^T\phi_{j_4}(s)\int\limits_t^s\phi_{j_1}(s_3)
\int\limits_{s_3}^{s}\phi_{j_3}(s_2)
\int\limits_{s_2}^{s}\phi_{j_3}(s_1)
ds_1ds_2ds_3ds
\zeta_{j_1}^{(i_1)}
\zeta_{j_4}^{(i_4)}=
$$
$$
=\hbox{\vtop{\offinterlineskip\halign{
\hfil#\hfil\cr
{\rm l.i.m.}\cr
$\stackrel{}{{}_{p\to \infty}}$\cr
}} }
\sum\limits_{j_4, j_1=0}^{p}
\frac{1}{2}\int\limits_t^T\phi_{j_4}(s)\int\limits_t^s\phi_{j_1}(s_3)
\sum\limits_{j_3=0}^{p}\left(\int\limits_{s_3}^{s}
\phi_{j_3}(s_2)ds_2\right)^2ds_3ds
\zeta_{j_1}^{(i_1)}
\zeta_{j_4}^{(i_4)}=
$$
$$
=\hbox{\vtop{\offinterlineskip\halign{
\hfil#\hfil\cr
{\rm l.i.m.}\cr
$\stackrel{}{{}_{p\to \infty}}$\cr
}} }
\sum\limits_{j_4, j_1=0}^{p}
\frac{1}{2}\int\limits_t^T\phi_{j_4}(s)\int\limits_t^s\phi_{j_1}(s_3)
(s-s_3)ds_3ds
\zeta_{j_1}^{(i_1)}
\zeta_{j_4}^{(i_4)} - \Delta_3^{(i_1i_4)}=
$$
$$
=
\frac{1}{2}\int\limits_t^T\int\limits_t^s(s-s_3)d{\bf w}_{s_3}^{(i_1)}
d{\bf w}_{s}^{(i_4)}+
$$
$$
+\frac{1}{2}{\bf 1}_{\{i_1=i_4\ne 0\}}
\lim_{p\to\infty}
\sum\limits_{j_4=0}^{p}
\int\limits_t^T\phi_{j_4}(s)\int\limits_t^s\phi_{j_4}(s_3)(s-s_3)ds_3ds
- \Delta_3^{(i_1i_4)}=
$$
$$
=
\frac{1}{2}\int\limits_t^T\int\limits_t^{s_2}\int\limits_t^{s_1}
d{\bf w}_{s}^{(i_1)}ds_1
d{\bf w}_{s_2}^{(i_4)}+
$$
$$
+
\frac{1}{2}{\bf 1}_{\{i_1=i_4\ne 0\}}
\left(
\sum\limits_{j_4=0}^{\infty}
\int\limits_t^T(s-t)\phi_{j_4}(s)\int\limits_t^s\phi_{j_4}(s_3)ds_3ds-\right.
$$
$$
\left.
-
\sum\limits_{j_4=0}^{\infty}
\int\limits_t^T\phi_{j_4}(s)\int\limits_t^s(s_3-t)\phi_{j_4}(s_3)ds_3ds
\right)
- \Delta_3^{(i_1i_4)}=
$$
\begin{equation}
\label{otit555}
=\frac{1}{2}\int\limits_t^T\int\limits_t^{s_2}\int\limits_t^{s_1}
d{\bf w}_{s}^{(i_1)}ds_1
d{\bf w}_{s_2}^{(i_4)} - \Delta_3^{(i_1i_4)}\ \ \ \hbox{w.~p.~1.}
\end{equation}

\vspace{4mm}

Let us consider $A_6^{(i_1i_2)}$
$$
A_6^{(i_1i_2)}=
$$
$$
=
\hbox{\vtop{\offinterlineskip\halign{
\hfil#\hfil\cr
{\rm l.i.m.}\cr
$\stackrel{}{{}_{p\to \infty}}$\cr
}} }
\sum\limits_{j_3, j_2, j_1=0}^{p}
\int\limits_t^T\phi_{j_1}(s_3)\int\limits_{s_3}^T\phi_{j_2}(s_2)
\int\limits_{s_2}^{T}\phi_{j_3}(s_1)
\int\limits_{s_1}^{T}\phi_{j_3}(s)dsds_1ds_2ds_3
\zeta_{j_1}^{(i_1)}
\zeta_{j_2}^{(i_2)}=
$$
$$
=\hbox{\vtop{\offinterlineskip\halign{
\hfil#\hfil\cr
{\rm l.i.m.}\cr
$\stackrel{}{{}_{p\to \infty}}$\cr
}} }
\sum\limits_{j_1, j_2=0}^{p}
\frac{1}{2}\int\limits_t^T\phi_{j_1}(s_3)\int\limits_{s_3}^T\phi_{j_2}(s_2)
\sum\limits_{j_3=0}^{p}\left(\int\limits_{s_2}^{T}
\phi_{j_3}(s)ds\right)^2ds_2ds_3
\zeta_{j_1}^{(i_1)}
\zeta_{j_2}^{(i_2)}=
$$
$$
=\hbox{\vtop{\offinterlineskip\halign{
\hfil#\hfil\cr
{\rm l.i.m.}\cr
$\stackrel{}{{}_{p\to \infty}}$\cr
}} }
\sum\limits_{j_1, j_2=0}^{p}
\frac{1}{2}\int\limits_t^T\phi_{j_1}(s_3)\int\limits_{s_3}^T\phi_{j_2}(s_2)
(T-s_2)ds_2ds_3
\zeta_{j_1}^{(i_1)}
\zeta_{j_2}^{(i_2)} - \Delta_6^{(i_1i_2)}=
$$
$$
=\hbox{\vtop{\offinterlineskip\halign{
\hfil#\hfil\cr
{\rm l.i.m.}\cr
$\stackrel{}{{}_{p\to \infty}}$\cr
}} }
\sum\limits_{j_1, j_2=0}^{p}
\frac{1}{2}\int\limits_t^T\phi_{j_2}(s_2)(T-s_2)
\int\limits_t^{s_2}\phi_{j_1}(s_3)
ds_3ds_2
\zeta_{j_1}^{(i_1)}
\zeta_{j_2}^{(i_2)} - \Delta_6^{(i_1i_2)}=
$$
$$
=\frac{1}{2}\int\limits_t^T(T-s_2)\int\limits_t^{s_2}d{\bf w}_{s_3}^{(i_1)}
d{\bf w}_{s_2}^{(i_2)}+
$$
$$
+
\frac{1}{2}{\bf 1}_{\{i_1=i_2\ne 0\}}
\sum\limits_{j_2=0}^{\infty}
\int\limits_t^T\phi_{j_2}(s_2)(T-s_2)\int\limits_t^{s_2}
\phi_{j_2}(s_3)ds_3ds_2
- \Delta_6^{(i_1i_2)}=
$$
\begin{equation}
\label{otit001}
=\frac{1}{2}\int\limits_t^T\int\limits_t^{s_1}\int\limits_t^{s_2}
d{\bf w}_{s}^{(i_1)}
d{\bf w}_{s_2}^{(i_2)}ds_1+
\frac{1}{4}{\bf 1}_{\{i_1=i_2\ne 0\}}
\int\limits_t^T(T-s_2)ds_2
- \Delta_6^{(i_1i_2)}\ \ \ \hbox{w.~p.~1.}
\end{equation}

\vspace{2mm}

Let us consider $B_1, B_2, B_3$
$$
B_1
=\hbox{\vtop{\offinterlineskip\halign{
\hfil#\hfil\cr
{\rm lim}\cr
$\stackrel{}{{}_{p\to \infty}}$\cr
}} }
\sum\limits_{j_1, j_4=0}^{p}
\frac{1}{2}\int\limits_t^T\phi_{j_4}(s)\int\limits_{t}^s\phi_{j_4}(s_1)
\left(\int\limits_{t}^{s_1}
\phi_{j_1}(s_2)ds_2\right)^2ds_1ds=
$$
$$
=\hbox{\vtop{\offinterlineskip\halign{
\hfil#\hfil\cr
{\rm lim}\cr
$\stackrel{}{{}_{p\to \infty}}$\cr
}} }
\sum\limits_{j_4=0}^{p}
\frac{1}{2}\int\limits_t^T\phi_{j_4}(s)\int\limits_{t}^s\phi_{j_4}(s_1)
(s_1-t)ds_1ds - \lim_{p\to\infty}\sum\limits_{j_4=0}^{p}
a_{j_4j_4}^p=
$$
\begin{equation}
\label{otit239}
=
\frac{1}{4}\int\limits_t^T(s_1-t)ds_1
-\lim_{p\to\infty}\sum\limits_{j_4=0}^{p}
a_{j_4j_4}^p,
\end{equation}

$$
B_2=
\hbox{\vtop{\offinterlineskip\halign{
\hfil#\hfil\cr
{\rm lim}\cr
$\stackrel{}{{}_{p\to \infty}}$\cr
}} }
\sum\limits_{j_4, j_3=0}^{p}
\int\limits_t^T\phi_{j_3}(s)
\int\limits_t^s\phi_{j_3}(s_2)
\int\limits_t^{s_2}\phi_{j_4}(s_3)ds_3
\int\limits_{s_2}^s\phi_{j_4}(s_1)ds_1 ds_2 ds=
$$
$$
=\hbox{\vtop{\offinterlineskip\halign{
\hfil#\hfil\cr
{\rm lim}\cr
$\stackrel{}{{}_{p\to \infty}}$\cr
}} }
\sum\limits_{j_4, j_3=0}^{p}\left(
\frac{1}{2}\int\limits_t^T\phi_{j_3}(s)
\left(\int\limits_{t}^{s}
\phi_{j_4}(s_3)ds_3\right)^2\int\limits_{t}^{s}\phi_{j_3}(s_2)ds_2ds-\right.
$$
$$
-
\frac{1}{2}\int\limits_t^T\phi_{j_3}(s)
\int\limits_{t}^{s}\phi_{j_3}(s_2)
\left(\int\limits_{t}^{s_2}
\phi_{j_4}(s_3)ds_3\right)^2ds_2ds-
$$
$$
-
\left.
\frac{1}{2}\int\limits_t^T\phi_{j_3}(s)
\int\limits_{t}^{s}\phi_{j_3}(s_2)
\left(\int\limits_{s_2}^{s}
\phi_{j_4}(s_1)ds_1\right)^2ds_2ds\right)=
$$
$$
=\sum\limits_{j_3=0}^{\infty}
\frac{1}{2}\int\limits_t^T\phi_{j_3}(s)(s-t)
\int\limits_{t}^{s}
\phi_{j_3}(s_2)ds_2ds -\lim_{p\to\infty}\sum\limits_{j_3=0}^p b_{j_3j_3}^p -
$$
$$
-
\sum\limits_{j_3=0}^{\infty}
\frac{1}{2}\int\limits_t^T\phi_{j_3}(s)
\int\limits_{t}^{s}(s_2-t)
\phi_{j_3}(s_2)ds_2ds +
\lim_{p\to\infty}\sum\limits_{j_3=0}^p a_{j_3j_3}^p-
$$
$$
-
\sum\limits_{j_3=0}^{\infty}
\frac{1}{2}\int\limits_t^T\phi_{j_3}(s)
\int\limits_{t}^{s}
\phi_{j_3}(s_2)(s-t+t-s_2)ds_2ds +
\lim_{p\to\infty}\sum\limits_{j_3=0}^p c_{j_3j_3}^p =
$$
\begin{equation}
\label{otit990}
=\lim_{p\to\infty}\sum\limits_{j_3=0}^p a_{j_3j_3}^p
+\lim_{p\to\infty}\sum\limits_{j_3=0}^p c_{j_3j_3}^p
-\lim_{p\to\infty}\sum\limits_{j_3=0}^p b_{j_3j_3}^p.
\end{equation}

\vspace{2mm}

Moreover,
$$
B_2+B_3=
\hbox{\vtop{\offinterlineskip\halign{
\hfil#\hfil\cr
{\rm lim}\cr
$\stackrel{}{{}_{p\to \infty}}$\cr
}} }
\sum\limits_{j_4, j_3=0}^{p}
\left(C_{j_3j_4j_3j_4}+C_{j_3j_4j_4j_3}\right)=
$$
$$
=\hbox{\vtop{\offinterlineskip\halign{
\hfil#\hfil\cr
{\rm lim}\cr
$\stackrel{}{{}_{p\to \infty}}$\cr
}} }
\sum\limits_{j_4, j_3=0}^{p}
\int\limits_t^T\phi_{j_3}(s)\int\limits_{t}^s\phi_{j_4}(s_1)
\int\limits_t^{s_1}\phi_{j_4}(s_2)\int\limits_t^{s_1}\phi_{j_3}(s_3)
ds_3ds_2ds_1ds
=
$$
$$
=\hbox{\vtop{\offinterlineskip\halign{
\hfil#\hfil\cr
{\rm lim}\cr
$\stackrel{}{{}_{p\to \infty}}$\cr
}} }
\sum\limits_{j_4, j_3=0}^{p}
\int\limits_t^T\phi_{j_4}(s_1)\int\limits_{t}^{s_1}\phi_{j_4}(s_2)
\int\limits_{t}^{s_1}\phi_{j_3}(s_3)ds_3ds_2\int\limits_{s_1}^{T}\phi_{j_3}(s)
dsds_1=
$$
$$
=\hbox{\vtop{\offinterlineskip\halign{
\hfil#\hfil\cr
{\rm lim}\cr
$\stackrel{}{{}_{p\to \infty}}$\cr
}} }
\sum\limits_{j_4, j_3=0}^{p}
\left(\int\limits_t^T\phi_{j_4}(s_1)\int\limits_{t}^{s_1}\phi_{j_4}(s_3)
\int\limits_{t}^{T}\phi_{j_3}(s_2)ds_2\int\limits_{s_1}^{T}\phi_{j_3}(s)
dsds_3ds_1-\right.
$$
$$
\left.-
\int\limits_t^T\phi_{j_4}(s_1)\int\limits_{t}^{s_1}\phi_{j_4}(s_3)
\left(\int\limits_{s_1}^{T}\phi_{j_3}(s)ds\right)^2
ds_3ds_1\right)=
$$
$$
=
\sum\limits_{j_4=0}^{\infty}
\int\limits_t^T\phi_{j_4}(s_1)(T-s_1)\int\limits_{t}^{s_1}\phi_{j_4}(s_3)
ds_3ds_1-
$$
$$
-\sum\limits_{j_4=0}^{\infty}
\int\limits_t^T\phi_{j_4}(s_1)(T-s_1)\int\limits_{t}^{s_1}\phi_{j_4}(s_3)
ds_3ds_1+
2\lim_{p\to\infty}\sum\limits_{j_4=0}^p f_{j_4j_4}^p
=
$$
\begin{equation}
\label{cas1000}
=2\lim_{p\to\infty}\sum\limits_{j_4=0}^p f_{j_4j_4}^p.
\end{equation}

Therefore,
\begin{equation}
\label{otit123}
~~~~B_3=2\lim_{p\to\infty}\sum\limits_{j_3=0}^p f_{j_3j_3}^p-
\lim_{p\to\infty}\sum\limits_{j_3=0}^p a_{j_3j_3}^p-
\lim_{p\to\infty}\sum\limits_{j_3=0}^p c_{j_3j_3}^p
+\lim_{p\to\infty}\sum\limits_{j_3=0}^p b_{j_3j_3}^p.
\end{equation}

After substituting the relations (\ref{otiteee1})--(\ref{otit123}) 
into (\ref{otiteee}), we obtain
$$
\hbox{\vtop{\offinterlineskip\halign{
\hfil#\hfil\cr
{\rm l.i.m.}\cr
$\stackrel{}{{}_{p\to \infty}}$\cr
}} }
\sum\limits_{j_1, j_2, j_3, j_4=0}^{p}
C_{j_4 j_3 j_2 j_1}\zeta_{j_1}^{(i_1)}\zeta_{j_2}^{(i_2)}\zeta_{j_3}^{(i_3)}
\zeta_{j_4}^{(i_4)}=
$$
$$
=
J[\psi^{(4)}]_{T,t}+
\frac{1}{2}{\bf 1}_{\{i_1=i_2\ne 0\}}
\int\limits_t^T\int\limits_t^s\int\limits_t^{s_1}ds_2
d{\bf w}_{s_1}^{(i_3)}
d{\bf w}_{s}^{(i_4)}+
$$
$$
+\frac{1}{2}{\bf 1}_{\{i_2=i_3\ne 0\}}
\int\limits_t^T\int\limits_t^{s_2}\int\limits_t^{s_1}
d{\bf w}_{s}^{(i_1)}ds_1
d{\bf w}_{s_2}^{(i_4)}
+\frac{1}{2}{\bf 1}_{\{i_3=i_4\ne 0\}}
\int\limits_t^T\int\limits_t^{s_1}\int\limits_t^{s_2}
d{\bf w}_{s}^{(i_1)}
d{\bf w}_{s_2}^{(i_2)}ds_1+
$$
\begin{equation}
\label{otit7776}
~~~~~+\frac{1}{4}{\bf 1}_{\{i_1=i_2\ne 0\}}
{\bf 1}_{\{i_3=i_4\ne 0\}}
\int\limits_t^T\int\limits_t^{s_1}ds_2
ds_1 + R = J^{*}[\psi^{(4)}]_{T,t}+R\ \ \  \hbox{w.~p.~1,}
\end{equation}

\noindent
where
$$
R=-{\bf 1}_{\{i_1=i_2\ne 0\}}\Delta_1^{(i_3i_4)}
+{\bf 1}_{\{i_1=i_3\ne 0\}}\left(
-\Delta_2^{(i_2i_4)}
+\Delta_1^{(i_2i_4)}
+\Delta_3^{(i_2i_4)}\right)+
$$
$$
+{\bf 1}_{\{i_1=i_4\ne 0\}}\left(
\Delta_4^{(i_2i_3)}-
\Delta_5^{(i_2i_3)}
+\Delta_6^{(i_2i_3)}\right)-
{\bf 1}_{\{i_2=i_3\ne 0\}}\Delta_3^{(i_1i_4)}+
$$
$$
+{\bf 1}_{\{i_2=i_4\ne 0\}}
\left(-\Delta_4^{(i_1i_3)}
+\Delta_5^{(i_1i_3)}
+\Delta_6^{(i_1i_3)}\right)-
{\bf 1}_{\{i_3=i_4\ne 0\}}\Delta_6^{(i_1i_2)}-
$$
$$
-
{\bf 1}_{\{i_1=i_3\ne 0\}}
{\bf 1}_{\{i_2=i_4\ne 0\}}\Biggl(
\lim_{p\to\infty}\sum\limits_{j_3=0}^p a_{j_3j_3}^p
+\lim_{p\to\infty}\sum\limits_{j_3=0}^p c_{j_3j_3}^p
-\lim_{p\to\infty}\sum\limits_{j_3=0}^p b_{j_3j_3}^p\Biggr)-
$$
$$
-{\bf 1}_{\{i_1=i_4\ne 0\}}
{\bf 1}_{\{i_2=i_3\ne 0\}}
\Biggl(2\lim_{p\to\infty}\sum\limits_{j_3=0}^p f_{j_3j_3}^p
-\lim_{p\to\infty}\sum\limits_{j_3=0}^p a_{j_3j_3}^p
-\Biggr.
$$
$$
\Biggl.-\lim_{p\to\infty}\sum\limits_{j_3=0}^p c_{j_3j_3}^p
+\lim_{p\to\infty}\sum\limits_{j_3=0}^p b_{j_3j_3}^p\Biggr)+
$$
\begin{equation}
\label{otiteee0}
+{\bf 1}_{\{i_1=i_2\ne 0\}}
{\bf 1}_{\{i_3=i_4\ne 0\}}\lim_{p\to\infty}\sum\limits_{j_3=0}^p a_{j_3j_3}^p.
\end{equation}

\vspace{2mm}

From (\ref{otit7776}) and (\ref{otiteee0})
it follows that Theorem 2.9 will be proved if
\begin{equation}
\label{otitgggh-1}
\Delta_k^{(ij)}=0\ \ \ \hbox{w.~p.~1},
\end{equation}
\begin{equation}
\label{otitgggh}
~~~~~~ \lim_{p\to\infty}\sum\limits_{j_3=0}^p a_{j_3j_3}^p=
\lim_{p\to\infty}\sum\limits_{j_3=0}^p b_{j_3j_3}^p=
\lim_{p\to\infty}\sum\limits_{j_3=0}^p c_{j_3j_3}^p=
\lim_{p\to\infty}\sum\limits_{j_3=0}^p f_{j_3j_3}^p=0,
\end{equation}

\noindent
where $k=1, 2,\ldots,6,$\ \  $i,j=0, 1,\ldots,m.$

Consider the case of Legendre polynomials.
Let us prove that $\Delta_1^{(i_3i_4)}=0$ w.~p.~1.
We have
$$
{\sf M}\left\{\left(\sum\limits_{j_3, j_4=0}^{p}
a_{j_4 j_3}^p \zeta_{j_3}^{(i_3)}
\zeta_{j_4}^{(i_4)}\right)^2\right\}=
$$
$$
=
\sum\limits_{j_3'=0}^p\sum\limits_{j_3=0}^{j_3'-1}
\Biggl(2a_{j_3j_3}^pa_{j_3'j_3'}^p
+\left(a_{j_3j_3'}^p\right)^2+
2a_{j_3j_3'}^p a_{j_3'j_3}^p+
\left(a_{j_3'j_3}^p\right)^2\Biggr)
+3\sum_{j_3'=0}^p\left(a_{j_3'j_3'}^p\right)^2=
$$
\begin{equation}
\label{otit321}
=\left(\sum_{j_3=0}^p a_{j_3j_3}^p\right)^2+
\sum\limits_{j_3'=0}^p\sum\limits_{j_3=0}^{j_3'-1}
\left(a_{j_3j_3'}^p + a_{j_3'j_3}^p\right)^2
+2\sum\limits_{j_3'=0}^p\left(a_{j_3'j_3'}^p\right)^2\ \ \ (i_3=i_4\ne 0),
\end{equation}

\vspace{-1mm}
\begin{equation}
\label{otit3210}
{\sf M}\left\{\left(\sum\limits_{j_3, j_4=0}^{p}
a_{j_4 j_3}^p \zeta_{j_3}^{(i_3)}
\zeta_{j_4}^{(i_4)}\right)^2\right\}=
\sum\limits_{j_3,j_4=0}^p
\left(a_{j_4j_3}^p\right)^2\ \ \ (i_3\ne i_4,\ i_3\ne 0,\ i_4\ne 0),
\end{equation}

\vspace{-3mm}
\begin{equation}
\label{otit32101}
{\sf M}\left\{\left(\sum\limits_{j_3, j_4=0}^{p}
a_{j_4 j_3}^p \zeta_{j_3}^{(i_3)}
\zeta_{j_4}^{(i_4)}\right)^2\right\}=
\left\{
\begin{matrix}
(T-t)\sum\limits_{j_4=0}^p\left(a_{j_4,0}^p\right)^2\ 
&\hbox{\rm if}\ \ \ 
i_3=0,\  i_4\ne 0\cr\cr
(T-t)\sum\limits_{j_3=0}^p\left(a_{0,j_3}^p\right)^2\  &\hbox{\rm if}\ \ \ 
i_4=0,\  i_3\ne 0\cr\cr
(T-t)^2\left(a_{00}^p\right)^2\  &\hbox{\rm if}\ \ \ 
i_3=i_4=0
\end{matrix}\right..
\end{equation}

Let us consider the case $i_3=i_4\ne 0$
$$
a_{j_4j_3}^p=\frac{(T-t)^2\sqrt{(2j_4+1)(2j_3+1)}}{32}\times
$$
$$
\times
\int\limits_{-1}^1 P_{j_4}(y) \int\limits_{-1}^y
P_{j_3}(y_1)\sum\limits_{j_1=p+1}^{\infty}(2j_1+1)
\left(\int\limits_{-1}^{y_1}P_{j_1}(y_2)dy_2\right)^2 dy_1dy=
$$
$$
=\frac{(T-t)^2\sqrt{(2j_4+1)(2j_3+1)}}{32}\times
$$
$$
\times
\int\limits_{-1}^1 P_{j_3}(y_1) 
\sum\limits_{j_1=p+1}^{\infty}\frac{1}{2j_1+1}
\left(P_{j_1+1}(y_1)-P_{j_1-1}(y_1)\right)^2
\int\limits_{y_1}^{1}P_{j_4}(y)dy dy_1=
$$
$$
=\frac{(T-t)^2\sqrt{2j_3+1}}{32\sqrt{2j_4+1}}\times
$$
$$
\times
\int\limits_{-1}^1 P_{j_3}(y_1) \left(P_{j_4-1}(y_1)-P_{j_4+1}(y_1)\right)
\sum\limits_{j_1=p+1}^{\infty}\frac{1}{2j_1+1}
\left(P_{j_1+1}(y_1)-P_{j_1-1}(y_1)\right)^2 dy_1
$$

\noindent
if $j_4\ne 0$ and
$$
a_{j_4j_3}^p=\frac{(T-t)^2\sqrt{2j_3+1}}{32}\times
$$
$$
\times
\int\limits_{-1}^1 P_{j_3}(y_1) (1-y_1)
\sum\limits_{j_1=p+1}^{\infty}\frac{1}{2j_1+1}
\left(P_{j_1+1}(y_1)-P_{j_1-1}(y_1)\right)^2
dy_1
$$

\noindent
if $j_4=0.$

From (\ref{ogo23}) and the estimate $\left| P_j(y) \right|\le 1$, $y\in [-1, 1]$
we obtain
\begin{equation}
\label{may2021}
\left|P_{j}(y)\right|=\sqrt{\left|P_{j}(y)\right|}\cdot \sqrt{\left|P_{j}(y)\right|}\le 
\frac{C}{j^{1/4}(1-y^2)^{1/8}},\ \ \  y\in (-1, 1),\ \ \ j\in{\bf N}.
\end{equation}

Using (\ref{ogo23}) and (\ref{may2021}), we get
\begin{equation}
\label{otitf11}
\left|a_{j_4j_3}^p\right|\le \frac{C_0}{\left(j_4\right)^{3/4}}
\sum\limits_{j_1=p+1}^{\infty}\frac{1}{j_1^2}\int\limits_{-1}^1
\frac{dy}{(1-y^2)^{7/8}}\le \frac{C_1}{p \left(j_4\right)^{3/4}}\ \ \ (j_3\ne 0,\ j_4\ge 2),
\end{equation}
\begin{equation}
\label{otitf112}
~~~~~~~\left|a_{0 j_3}^p\right|+\left|a_{1 j_3}^p\right|\le C_0
\sum\limits_{j_1=p+1}^{\infty}\frac{1}{j_1^2}\int\limits_{-1}^1
\frac{dy}{(1-y^2)^{3/4}}\le \frac{C_1}{p}\ \ \ (j_3\ne 0),
\end{equation}
\begin{equation}
\label{otitf113}
~~~~~~~\left|a_{j_4 0}^p\right|+ |a_{00}^p|\le C_0
\sum\limits_{j_1=p+1}^{\infty}\frac{1}{j_1^2}\int\limits_{-1}^1
\frac{dy}{(1-y^2)^{1/2}}\le \frac{C_1}{p}\ \ \ (j_4\ge 1),
\end{equation}

\vspace{1mm}
\noindent
where constants $C_0, C_1$ do not depend on $p$.

Taking into account (\ref{otit321}), (\ref{otitf11})--(\ref{otitf113}),
we have 
$$
{\sf M}\left\{\left(\sum\limits_{j_3, j_4=0}^{p}
a_{j_4 j_3}^p \zeta_{j_3}^{(i_3)}
\zeta_{j_4}^{(i_4)}\right)^2\right\}
=\left(a_{00}^p+\sum_{j_3=1}^p a_{j_3j_3}^p\right)^2+
\sum\limits_{j_3'=1}^p
\left(a_{0 j_3'}^p + a_{j_3' 0}^p\right)^2+
$$
$$
+
\sum\limits_{j_3'=1}^p\sum\limits_{j_3=1}^{j_3'-1}
\left(a_{j_3j_3'}^p + a_{j_3'j_3}^p\right)^2
+2\left(\sum\limits_{j_3'=1}^p\left(a_{j_3'j_3'}^p\right)^2+
\left(a_{00}\right)^2\right)\le
$$
$$
\le K_0\left(\frac{1}{p}+\frac{1}{p}\sum\limits_{j_3=1}^p
\frac{1}{\left(j_3\right)^{3/4}}\right)^2+\frac{K_1}{p}
+K_2\sum\limits_{j_3'=1}^p\sum\limits_{j_3=1}^{j_3'-1}
\frac{1}{p^2}\Biggl(\frac{1}{\left(j_3'\right)^{3/4}}+\frac{1}{\left(j_3\right)^{3/4}}\Biggr)^2\le
$$
$$
\le
K_0\left(\frac{1}{p}+\frac{1}{p}\int\limits_0^p\frac{dx}{x^{3/4}}
\right)^2 + \frac{K_1}{p} + \frac{K_3}{p}\sum\limits_{j_3=1}^p\frac{1}{\left(j_3\right)^{3/2}}\le
$$
$$
\le K_0\Biggl(\frac{1}{p}+\frac{4}{p^{3/4}}\Biggr)^2
+
\frac{K_1}{p} +\frac{K_3}{p}\left(1+ \int\limits_1^p\frac{dx}{x^{3/2}}\right)\le
$$
$$
\le \frac{K_4}{p} + \frac{K_3}{p}\left(3-\frac{2}{\sqrt{p}}\right)\le \frac{K_5}{p}\ \to 0\
$$

\noindent
if $p\to \infty$\  $(i_3=i_4\ne 0).$

The same result for the cases  
(\ref{otit3210}), (\ref{otit32101})  
also follows from the estimates 
(\ref{otitf11})--(\ref{otitf113}). Therefore,
\begin{equation}
\label{otitf14}
\Delta_1^{(i_3i_4)}=0\ \ \ \hbox{w.~p.~1}.
\end{equation}

It is not difficult to see that the formulas
\begin{equation}
\label{2017f}
\Delta_2^{(i_2i_4)}=0,\ \ \ \Delta_4^{(i_1i_3)}=0,\ \ \
\Delta_6^{(i_1i_3)}=0\ \ \ \hbox{w.~p.~1}
\end{equation}

\noindent
can be proved similarly with the 
proof of (\ref{otitf14}).

Moreover, from the estimates (\ref{otitf11})--(\ref{otitf113}) we obtain
\begin{equation}
\label{20177}
\lim\limits_{p\to\infty}
\sum\limits_{j_3=0}^p a_{j_3j_3}^p=0.
\end{equation}

The relations 
\begin{equation}
\label{cas79}
\lim\limits_{p\to\infty}
\sum\limits_{j_3=0}^p b_{j_3j_3}^p=0\ \ \ \ {\rm and}\ \ \ \ 
\lim\limits_{p\to\infty}
\sum\limits_{j_3=0}^p f_{j_3j_3}^p=0
\end{equation}

\noindent
can also be proved analogously with (\ref{20177}).

Let us consider $\Delta_3^{(i_2i_4)}$
\begin{equation}
\label{otitzzz}
~~~~~~~~~\Delta_3^{(i_2i_4)}=\Delta_4^{(i_2i_4)}+
\Delta_6^{(i_2i_4)}-\Delta_7^{(i_2i_4)}=
-\Delta_7^{(i_2i_4)}\ \ \ \hbox{\rm w.~p.~1},
\end{equation}

\noindent
where
$$
\Delta_7^{(i_2i_4)}=
\hbox{\vtop{\offinterlineskip\halign{
\hfil#\hfil\cr
{\rm l.i.m.}\cr
$\stackrel{}{{}_{p\to \infty}}$\cr
}} }
\sum\limits_{j_2, j_4=0}^{p}
g_{j_4 j_2}^p \zeta_{j_2}^{(i_2)}
\zeta_{j_4}^{(i_4)},
$$
$$
g_{j_4 j_2}^p=
\int\limits_t^T\phi_{j_4}(s)\int\limits_{t}^s\phi_{j_2}(s_1)
\sum\limits_{j_1=p+1}^{\infty}\left(\int\limits_{s_1}^{T}
\phi_{j_1}(s_2)ds_2\int\limits_{s}^{T}
\phi_{j_1}(s_2)ds_2\right)  
ds_1ds=
$$
\begin{equation}
\label{otitzz}
~~~~~~~~=\sum\limits_{j_1=p+1}^{\infty}
\int\limits_t^T\phi_{j_4}(s)\int\limits_{s}^T\phi_{j_1}(s_2)
ds_2\int\limits_{t}^{s}
\phi_{j_2}(s_1)\int\limits_{s_1}^{T}
\phi_{j_1}(s_2)ds_2ds_1ds.
\end{equation}

\vspace{1mm}

The last step in (\ref{otitzz}) follows from the estimate
$$
\left|g_{j_4 j_2}^p\right|
\le K
\sum\limits_{j_1=p+1}^{\infty}\frac{1}{j_1^2}
\int\limits_{-1}^1 \frac{1}{(1-y^2)^{1/2}}
\int\limits_{-1}^y \frac{1}{(1-x^2)^{1/2}}dx dy\
\le \frac{K_1}{p}.
$$

Note that
\begin{equation}
\label{otitddd}
g_{j_4 j_4}^p=\sum\limits_{j_1=p+1}^{\infty}
\frac{1}{2}\left(
\int\limits_t^T\phi_{j_4}(s)\int\limits_{s}^T\phi_{j_1}(s_2)
ds_2ds\right)^2,
\end{equation}
\begin{equation}
\label{otitddd1}
g_{j_4 j_2}^p+g_{j_2 j_4}^p=\sum\limits_{j_1=p+1}^{\infty}
\int\limits_t^T\phi_{j_4}(s)\int\limits_{s}^T\phi_{j_1}(s_2)
ds_2ds
\int\limits_t^T\phi_{j_2}(s)\int\limits_{s}^T\phi_{j_1}(s_2)
ds_2ds,
\end{equation}

\noindent
and
$$
g_{j_4j_2}^p=
\frac{(T-t)^2\sqrt{(2j_4+1)(2j_2+1)}}{16}\times
$$
$$
\times
\sum\limits_{j_1=p+1}^{\infty}\frac{1}{2j_1+1}
\int\limits_{-1}^1 P_{j_4}(y_1)\left(P_{j_1-1}(y_1)-P_{j_1+1}(y_1)\right)
\times
$$
$$
\times
\int\limits_{-1}^{y_1} P_{j_2}(y)\left(P_{j_1-1}(y)-P_{j_1+1}(y)\right)
dydy_1,\ \ \ j_4, j_2\le p.
$$

Due to orthogonality of the Legendre polynomials we obtain
$$
g_{j_4j_2}^p+g_{j_2j_4}^p=
\frac{(T-t)^2\sqrt{(2j_4+1)(2j_2+1)}}{16}\times
$$
$$
\times
\sum\limits_{j_1=p+1}^{\infty}\frac{1}{2j_1+1}
\int\limits_{-1}^1 P_{j_4}(y_1)\left(P_{j_1-1}(y_1)-P_{j_1+1}(y_1)\right)
dy_1
\times
$$
$$
\times
\int\limits_{-1}^{1} P_{j_2}(y)\left(P_{j_1-1}(y)-P_{j_1+1}(y)\right)
dy=
$$
$$
=\frac{(T-t)^2(2p+1)}{16} \frac{1}{2p+3}\left(
\int\limits_{-1}^1 P_p^2(y_1)dy_1\right)^2
\cdot \left\{
\begin{matrix}
1\ &{\rm if}\ j_2=j_4=p\cr\cr
0\ &\hbox{\rm otherwise}
\end{matrix}\right.
=
$$
\begin{equation}
\label{otitz}
=\frac{(T-t)^2}{4(2p+3)(2p+1)}
\cdot \left\{
\begin{matrix}
1\ &\hbox{\rm if}\ j_2=j_4=p\cr\cr
0\ &\hbox{\rm otherwise}
\end{matrix}\right.,
\end{equation}

\vspace{5mm}
$$
g_{j_4j_4}^p=
\frac{(T-t)^2(2j_4+1)}{16}\times
$$
$$
\times
\sum\limits_{j_1=p+1}^{\infty}\frac{1}{2j_1+1}\cdot \frac{1}{2}
\left(\int\limits_{-1}^1 P_{j_4}(y_1)
\left(P_{j_1-1}(y_1)-P_{j_1+1}(y_1)\right)
dy_1\right)^2=
$$
$$
=\frac{(T-t)^2(2p+1)}{32} \frac{1}{2p+3}\left(
\int\limits_{-1}^1 P_p^2(y_1)dy_1\right)^2
\cdot \left\{
\begin{matrix}
1\  &\hbox{\rm if}\ j_4=p\cr\cr
0\ &\hbox{\rm otherwise}
\end{matrix}\right.
=
$$
\begin{equation}
\label{otitz21}
=\frac{(T-t)^2}{8(2p+3)(2p+1)}
\cdot \left\{
\begin{matrix}
1\  &\hbox{\rm if}\ j_4=p\cr\cr
0\ &\hbox{\rm otherwise}
\end{matrix}\right..
\end{equation}

\vspace{3mm}

From (\ref{otit321}), (\ref{otitz}), and (\ref{otitz21}) it follows that
$$
{\sf M}\left\{\left(\sum\limits_{j_2, j_4=0}^{p}
g_{j_4 j_2}^p \zeta_{j_2}^{(i_2)}
\zeta_{j_4}^{(i_4)}\right)^2\right\}
=
$$
$$
=\Biggl(\sum_{j_3=0}^p g_{j_3j_3}^p\Biggr)^2+
\sum\limits_{j_3'=0}^p\sum\limits_{j_3=0}^{j_3'-1}
\left(g_{j_3 j_3'}^p + g_{j_3' j_3}^p\right)^2 +
2\sum\limits_{j_3'=0}^p\left(g_{j_3'j_3'}^p\right)^2=
$$
$$
=
\left(\frac{(T-t)^2}{8(2p+3)(2p+1)}\right)^2 +\ 0\ +
2\left(\frac{(T-t)^2}{8(2p+3)(2p+1)}\right)^2\ \to 0
$$

\vspace{3mm}
\noindent
if $p\to \infty$\ $(i_2=i_4\ne 0).$

Let us consider the case
$i_2\ne i_4,\ i_2\ne 0,\ i_4\ne 0$ (see (\ref{otit3210})).
It is not difficult to see that 
$$
g_{j_4 j_2}^p=
\int\limits_t^T\phi_{j_4}(s)\int\limits_{t}^s\phi_{j_2}(s_1)
F_p(s,s_1)  
ds_1ds=
\int\limits_{[t,T]^2}
K_p(s,s_1)\phi_{j_4}(s)\phi_{j_2}(s_1)ds_1 ds
$$
is a coefficient of the double Fourier--Legendre series of the function
\begin{equation}
\label{2017i}
K_p(s,s_1)={\bf 1}_{\{s_1<s\}}F_p(s,s_1),
\end{equation}

\noindent
where
$$
\sum\limits_{j_1=p+1}^{\infty}
\int\limits_{s_1}^{T}
\phi_{j_1}(s_2)ds_2\int\limits_{s}^{T}
\phi_{j_1}(s_2)ds_2\stackrel{\sf def}{=}F_p(s,s_1).
$$

The Parseval equality in this case looks as follows
\begin{equation}
\label{2017j}
\lim_{p_1\to\infty}\sum\limits_{j_4,j_2=0}^{p_1}
\left(g_{j_4 j_2}^p\right)^2=\int\limits_{[t,T]^2}
\left(K_p(s,s_1)\right)^2 ds_1 ds=
\int\limits_t^T\int\limits_t^s\left(
F_p(s,s_1)\right)^2 ds_1 ds.
\end{equation}

From (\ref{ogo23}) we obtain
$$
\left|\hspace{0.3mm}\int\limits_{s_1}^T\phi_{j_1}(\theta)d\theta\right|=
\frac{1}{2}\sqrt{2 j_1+1}\sqrt{T-t}\left|~\int\limits_{z(s_1)}^1
P_{j_1}(y)dy\right|=
$$
\begin{equation}
\label{2017c}
~~~~~~ =\frac{\sqrt{T-t}}{2\sqrt{2j_1+1}}\left|P_{j_1-1}(z(s_1))-P_{j_1+1}(z(s_1))
\right|\le \frac{K}{j_1}\frac{1}{\left(1-z^2(s_1)\right)^{1/4}},
\end{equation}

\vspace{2mm}
\noindent
where $z(s_1)$ is defined by (\ref{zz1}),\ \ $s_1\in (t, T).$

From (\ref{2017c}) we have
\begin{equation}
\label{2017u}
~~~~~~~ \left(F_p(s,s_1)\right)^2
\le \frac{C^2}{p^2}
\frac{1}{\left(1-z^2(s)\right)^{1/2}}
\frac{1}{\left(1-z^2(s_1)\right)^{1/2}},\ \ \ s, s_1\in (t, T).
\end{equation}

From (\ref{2017u}) it follows that
$\left|F_p(s,s_1)\right|
\le M_{\varepsilon}/p$\ \
in the domain 
$$
D_{\varepsilon}=\{(s, s_1):\ s\in[t+\varepsilon,
T-\varepsilon],\ s_1\in [t+\varepsilon, s]\}\ \ \ \hbox{for some small}\ \
\varepsilon >0,
$$
where constant $M_{\varepsilon}$ does not depend on $s, s_1.$
Then we have the uniform convergence 
\begin{equation}
\label{2017z}
~~~~~~~\sum\limits_{j_1=0}^{p}\int\limits_s^T\phi_{j_1}(\theta)d\theta
\int\limits_{s_1}^T\phi_{j_1}(\theta)d\theta \to
\sum\limits_{j_1=0}^{\infty}\int\limits_s^T\phi_{j_1}(\theta)d\theta
\int\limits_{s_1}^T\phi_{j_1}(\theta)d\theta 
\end{equation}
at the set $D_{\varepsilon}$
if $p\to \infty.$

Because of continuity of the function on the left-hand side of (\ref{2017z})
we obtain continuity of the limit function on the right-hand side 
of (\ref{2017z})
at the set $D_{\varepsilon}.$

Using this fact and (\ref{2017u}), we obtain
$$
\int\limits_t^T\int\limits_t^s\left(
F_p(s,s_1)\right)^2 ds_1 ds=
\lim\limits_{\varepsilon\to +0}
\int\limits_{t+\varepsilon}^{T-\varepsilon}\int\limits_{t+\varepsilon}^s
\left(F_p(s,s_1)\right)^2 ds_1 ds\le
$$
$$
\le\frac{C^2}{p^2}
\lim\limits_{\varepsilon\to +0}
\int\limits_{t+\varepsilon}^{T-\varepsilon}\int\limits_{t+\varepsilon}^s
\frac{ds_1}{\left(1-z^2(s_1)\right)^{1/2}}
\frac{ds}{\left(1-z^2(s)\right)^{1/2}}=
$$
$$
=
\frac{C^2}{p^2}
\int\limits_{t}^{T}\int\limits_{t}^s
\frac{ds_1}{\left(1-z^2(s_1)\right)^{1/2}}
\frac{ds}{\left(1-z^2(s)\right)^{1/2}}=
$$
\begin{equation}
\label{2017x}
=\frac{K}{p^2}
\int\limits_{-1}^{1}\int\limits_{-1}^y
\frac{dy_1}{\left(1-y_1^2\right)^{1/2}}
\frac{dy}{\left(1-y^2\right)^{1/2}}
\le \frac{K_1}{p^2},
\end{equation}

\noindent
where constant $K_1$ does not depend on $p.$

From (\ref{2017x}) and (\ref{2017j}) we get
\begin{equation}
\label{2017xx}
~~~0\le \sum\limits_{j_2,j_4=0}^{p}
\left(g_{j_4 j_2}^p\right)^2\le
\lim_{p_1\to\infty}\sum\limits_{j_2,j_4=0}^{p_1}
\left(g_{j_4 j_2}^p\right)^2=
\sum\limits_{j_2,j_4=0}^{\infty}
\left(g_{j_4 j_2}^p\right)^2
\le \frac{K_1}{p^2} \to 0
\end{equation}

\noindent
if $p\to\infty$. The case $i_2\ne i_4,$\ $i_2\ne 0,$\ $i_4\ne 0$ is proved.

The same result for the cases\\
\noindent
1)\ $i_2=0,$\ $i_4\ne 0,$\\ 
2)\ $i_4=0,$\ $i_2\ne 0,$\\ 
3)\ $i_2=0,$\ $i_4=0$\\
\noindent
can also be obtained. Then $\Delta_7^{(i_2i_4)}=0$
and $\Delta_3^{(i_2i_4)}=0$\ \  w.~p.~1.

Let us consider $\Delta_5^{(i_1 i_3)}$
$$
\Delta_5^{(i_1 i_3)}=\Delta_4^{(i_1 i_3)}+
\Delta_6^{(i_1 i_3)}-\Delta_8^{(i_1 i_3)}\ \ \ \hbox{\rm w.~p.~1,}
$$
where 
$$
\Delta_8^{(i_1i_3)}=
\hbox{\vtop{\offinterlineskip\halign{
\hfil#\hfil\cr
{\rm l.i.m.}\cr
$\stackrel{}{{}_{p\to \infty}}$\cr
}} }
\sum\limits_{j_3, j_1=0}^{p}
h_{j_3 j_1}^p \zeta_{j_1}^{(i_1)}
\zeta_{j_3}^{(i_3)},
$$
$$
h_{j_3 j_1}^p=
\int\limits_t^T\phi_{j_1}(s_3)\int\limits_{s_3}^T\phi_{j_3}(s)
F_p(s_3,s)
dsds_3.
$$

Analogously, we obtain that $\Delta_8^{(i_1i_3)}=0$\ \ w.~p.~1.
Here we consider the function
$$
K_p(s_3,s)={\bf 1}_{\{s_3<s\}}F_p(s_3,s)
$$
and the relation 
$$
h_{j_3j_1}^{p}=\int\limits_{[t,T]^2}
K_p(s_3,s)\phi_{j_1}(s_3)\phi_{j_3}(s)ds ds_3
$$
for the case $i_1\ne i_3,$\ $i_1\ne 0,$\ $i_3\ne 0.$

For the case $i_1=i_3\ne 0$ we use (see (\ref{otitddd}),
(\ref{otitddd1}))
$$
h_{j_1 j_1}^p=\sum\limits_{j_4=p+1}^{\infty}
\frac{1}{2}\left(
\int\limits_t^T\phi_{j_1}(s)\int\limits_{s}^T\phi_{j_4}(s_1)
ds_1ds\right)^2,
$$
$$
h_{j_3 j_1}^p+h_{j_1 j_3}^p=\sum\limits_{j_4=p+1}^{\infty}
\int\limits_t^T\phi_{j_1}(s)\int\limits_{s}^T\phi_{j_4}(s_2)
ds_2ds
\int\limits_t^T\phi_{j_3}(s)\int\limits_{s}^T\phi_{j_4}(s_2)
ds_2ds.
$$

Let us prove that
\begin{equation}
\label{uu1ggg}
\lim\limits_{p\to\infty}
\sum\limits_{j_3=0}^p c_{j_3j_3}^p=0.
\end{equation}

We have 
\begin{equation}
\label{otitbb}
c_{j_3j_3}^p=
f_{j_3j_3}^p+
d_{j_3j_3}^p-g_{j_3j_3}^p.
\end{equation}

Moreover, 
\begin{equation}
\label{rrr}
\lim\limits_{p\to\infty}\sum\limits_{j_3=0}^p f_{j_3j_3}^p=0,\ \ \
\lim\limits_{p\to\infty}\sum\limits_{j_3=0}^p d_{j_3j_3}^p=0,
\end{equation}

\noindent
where the first equality in (\ref{rrr}) has been proved earlier.  
Analogously, we can prove the second equality in (\ref{rrr}).

From (\ref{otitz21}) we obtain
$$
0\le \lim\limits_{p\to\infty}\sum\limits_{j_3=0}^p g_{j_3j_3}^p\le
\lim_{p\to\infty}\frac{(T-t)^2}{8(2p+3)(2p+1)}=0.
$$

So, (\ref{uu1ggg}) is proved.
The relations (\ref{otitgggh-1}), (\ref{otitgggh}) are proved for the polynomial case.
Theorem 2.9 is proved for the case of Legendre polynomials.

Let us consider the trigonometric case. 
According to (\ref{rr1xx}), we have
\begin{equation}
\label{agentu2}
~~~~~~~~ a_{j_4 j_3}^p=\frac{1}{2}\int\limits_t^T\phi_{j_3}(s_1)
\sum\limits_{j_1=p+1}^{\infty}\left(\int\limits_t^{s_1}
\phi_{j_1}(s_2)ds_2\right)^2\int\limits_{s_1}^T\phi_{j_4}(s)dsds_1.
\end{equation}

Moreover (see (\ref{2017x11}), (\ref{2017x12})), 
\begin{equation}
\label{agentu1}
\left|\int\limits_t^{s_1}\phi_{j}(s_2)ds_2\right|\le
\frac{K}{j},\ \ \ 
\left|\int\limits_{s_1}^T\phi_{j}(s_2)ds_2\right|\le
\frac{K}{j},
\end{equation}

\noindent
where constant $K$ does not depend on $j$ ($j=1,2,\ldots $).

Note that
$$
\int\limits_{s_1}^{T}\phi_{0}(s)ds=\frac{T-s_1}{\sqrt{T-t}}.
$$

Using (\ref{agentu2}) and (\ref{agentu1}), we obtain
\begin{equation}
\label{otithj}
~~~~~~~~~\left\vert a_{j_4j_3}^p\right\vert \le \frac{C_1}{j_4}\sum\limits_{j_1=p+1}^{\infty}
\frac{1}{j_1^2}\le \frac{C_1}{p j_4}\ \ \ (j_4\ne 0),\ \ \ \ \ 
\left\vert a_{0 j_3}^p\right\vert \le \frac{C_1}{p},
\end{equation}

\vspace{1mm}
\noindent
where constant $C_1$ does not depend on $p.$

Taking into account (\ref{otit321})--(\ref{otit32101}) and
(\ref{otithj}),
we obtain that
$\Delta_1^{(i_3i_4)}=0$\ \  w.~p.~1.
Analogously, we get
$\Delta_2^{(i_2i_4)}=0,$ $\Delta_4^{(i_1i_3)}=0,$
$\Delta_6^{(i_1i_3)}=0$\ \ w.~p.~1 and
$$
\lim\limits_{p\to\infty}\sum\limits_{j_3=0}^p a_{j_3j_3}^p=0,\ \ \
\lim\limits_{p\to\infty}\sum\limits_{j_3=0}^p b_{j_3j_3}^p=0,\ \ \ 
\lim\limits_{p\to\infty}\sum\limits_{j_3=0}^p f_{j_3j_3}^p=0.
$$

\vspace{1mm}

Let us consider $\Delta_3^{(i_2i_4)}$ for the case $i_2=i_4\ne 0.$
For the values $g_{j_4j_2}^{2m}+g_{j_2j_4}^{2m}$ and
$g_{j_4j_2}^{2m-1}+g_{j_2j_4}^{2m-1}$ $(m\in{\bf N})$ we have
(see (\ref{otitddd1}))

$$
g_{j_4 j_2}^{2m}+g_{j_2 j_4}^{2m}=
$$
$$
=
\sum\limits_{j_1=2m+1}^{\infty}
\int\limits_t^T\phi_{j_4}(s)\int\limits_{s}^T\phi_{j_1}(s_2)
ds_2ds
\int\limits_t^T\phi_{j_2}(s)\int\limits_{s}^T\phi_{j_1}(s_2)
ds_2ds=
$$
$$
=\sum\limits_{r=m+1}^{\infty}
\left(\int\limits_t^T\phi_{j_4}(s)\int\limits_{s}^T\phi_{2r-1}(s_2)
ds_2ds
\int\limits_t^T\phi_{j_2}(s)\int\limits_{s}^T\phi_{2r-1}(s_2)
ds_2ds+\right.
$$
\begin{equation}
\label{agentu3}
~~~~~~~~ +\left.\int\limits_t^T\phi_{j_4}(s)\int\limits_{s}^T\phi_{2r}(s_2)
ds_2ds
\int\limits_t^T\phi_{j_2}(s)\int\limits_{s}^T\phi_{2r}(s_2)
ds_2ds\right),
\end{equation}

\vspace{3mm}

$$
g_{j_4 j_2}^{2m-1}+g_{j_2 j_4}^{2m-1}=
$$
$$
=
\sum\limits_{j_1=2m}^{\infty}
\int\limits_t^T\phi_{j_4}(s)\int\limits_{s}^T\phi_{j_1}(s_2)
ds_2ds
\int\limits_t^T\phi_{j_2}(s)\int\limits_{s}^T\phi_{j_1}(s_2)
ds_2ds=
$$
$$
=g_{j_4 j_2}^{2m}+g_{j_2 j_4}^{2m}+
$$
\begin{equation}
\label{agentu4}
~~~~~~~~+
\int\limits_t^T\phi_{j_4}(s)\int\limits_{s}^T\phi_{2m}(s_2)
ds_2ds
\int\limits_t^T\phi_{j_2}(s)\int\limits_{s}^T\phi_{2m}(s_2)
ds_2ds,
\end{equation}

\vspace{2mm}
\noindent
where
$$
\int\limits_t^{T}\phi_{j_4}(s)
\int\limits_s^{T}\phi_{2r-1}(s_2)ds_2ds=
\sqrt{\frac{2}{T-t}}\int\limits_t^{T}\phi_{j_4}(s)
\int\limits_s^{T}{\rm sin}\frac{2\pi r(s_2-t)}{T-t}ds_2ds=
$$
$$
=\frac{\sqrt{2}\sqrt{T-t}}{2\pi r}
\int\limits_t^T
\phi_{j_4}(s)\biggl(
{\rm cos}\frac{2\pi r(s-t)}{T-t}-1\biggr)ds,
$$

\vspace{2mm}
$$
\int\limits_t^{T}\phi_{j_4}(s)
\int\limits_s^{T}\phi_{2r}(s_2)ds_2ds=
\sqrt{\frac{2}{T-t}}\int\limits_t^{T}\phi_{j_4}(s)
\int\limits_s^{T}{\rm cos}\frac{2\pi r(s_2-t)}{T-t}ds_2ds=
$$
$$
=\frac{\sqrt{2}\sqrt{T-t}}{2\pi r}
\int\limits_t^T
\phi_{j_4}(s)\biggl(
-{\rm sin}\frac{2\pi r(s-t)}{T-t}\biggr)ds,
$$

\vspace{3mm}
\noindent
where $2r-1,\ 2r\ge p+1,$\ and $j_2, j_4=0, 1,\ldots,p.$ 

Due to orthogonality of the trigonometric functions we have

\vspace{-4mm}
\begin{equation}
\label{otitjjj}
~~~~~~\int\limits_t^{T}\phi_{j_4}(s)
\int\limits_s^{T}\phi_{2r-1}(s_2)ds_2ds=
\frac{\sqrt{2}(T-t)}{2\pi r}
\cdot \left\{
\begin{matrix}
-1\ &\hbox{\rm if}\ \ j_4=0\cr\cr
0\ &\hbox{\rm otherwise}
\end{matrix}\right.,
\end{equation}

\vspace{-2mm}
\begin{equation}
\label{otitjjj11}
\int\limits_t^{T}\phi_{j_4}(s)
\int\limits_s^{T}\phi_{2r}(s_2)ds_2ds=0,
\end{equation}

\vspace{2mm}
\noindent
where $2r-1,\ 2r\ge p+1,$\ and $j_4=0, 1,\ldots,p.$

From (\ref{agentu3}), (\ref{otitjjj}), and (\ref{otitjjj11})
we obtain

\vspace{-1mm}
$$
g_{j_4j_2}^{2m}+g_{j_2j_4}^{2m}=
$$

\vspace{-1mm}
$$
=\sum\limits_{j_1=m+1}^{\infty}
\frac{(T-t)^2}{2\pi^2 j_1^2}
\cdot \left\{
\begin{matrix}
1\ &\hbox{\rm if}\ \ j_2=j_4=0\cr\cr
0\ &\hbox{\rm otherwise}
\end{matrix}\right.,
$$
$$
g_{j_4j_4}^{2m}=
\frac{1}{2}\left(g_{j_4j_2}^{2m}+g_{j_2j_4}^{2m}\right)\biggl.\biggr|_{j_2=j_4}=
$$

$$
=
\sum\limits_{j_1=m+1}^{\infty}
\frac{(T-t)^2}{4\pi^2 j_1^2}
\cdot \left\{
\begin{matrix}
1\ &\hbox{\rm if}\ \ j_4=0\cr\cr
0\ \ &\hbox{\rm otherwise}
\end{matrix}\right..
$$

\vspace{3mm}

Therefore (see (\ref{obana})),
\begin{equation}
\label{otitkkk}
\left\{
\begin{matrix}
\left|g_{j_4j_2}^{2m}+g_{j_2j_4}^{2m}\right|\le K_1/(2m)
&\hbox{\rm if}\ \ j_2=j_4=0\cr\cr
g_{j_4j_2}^{2m}+g_{j_2j_4}^{2m}=0\ \ &\hbox{\rm otherwise}
\end{matrix}\right.,
\end{equation}

\vspace{3mm}
\begin{equation}
\label{otitkkk11}
\left\{
\begin{matrix}
\left|g_{j_4j_4}^{2m}\right|\le K_1/(2m)
&\hbox{\rm if}\ \ j_4=0\cr\cr
g_{j_4j_4}^{2m}=0\ \ &\hbox{\rm otherwise}
\end{matrix}\right.,
\end{equation}

\vspace{4mm}
\noindent
where constant $K_1$ does not depend on $p=2m$.

For $p=2m-1$ from (\ref{agentu4}) and (\ref{otitjjj11}) we have

\vspace{-2mm}
$$
g_{j_4 j_2}^{2m-1}+g_{j_2 j_4}^{2m-1}=
$$

\vspace{-3mm}
\begin{equation}
\label{agentrr200}
=\sum\limits_{j_1=m+1}^{\infty}
\frac{(T-t)^2}{2\pi^2 j_1^2}
\cdot \left\{
\begin{matrix}
1\ \ \hbox{\rm or}\ \ 0\ &\hbox{\rm if}\ j_2=j_4=0\cr\cr
0\ &\hbox{\rm otherwise}
\end{matrix}\right..
\end{equation}

\vspace{3mm}

The relation (\ref{agentrr200}) implies that

\vspace{-2mm}
$$
g_{j_4j_4}^{2m-1}=
\frac{1}{2}\left(g_{j_4j_2}^{2m-1}+g_{j_2j_4}^{2m-1}\right)\biggl.\biggr|_{j_2=j_4}
=
$$

\begin{equation}
\label{agentrr201}
=\sum\limits_{j_1=m+1}^{\infty}
\frac{(T-t)^2}{4\pi^2 j_1^2}
\cdot \left\{
\begin{matrix}
1\ \ \hbox{\rm or}\ \ 0\ &\hbox{\rm if}\ j_4=0\cr\cr
0\ &\hbox{\rm otherwise}
\end{matrix}\right..
\end{equation}

\vspace{3mm}

Using (\ref{agentrr200}) and (\ref{agentrr201}), we obtain
\begin{equation}
\label{agentoo1}
~~~~~~~~ \left\{
\begin{matrix}
\left|g_{j_4j_2}^{2m-1}+g_{j_2j_4}^{2m-1}\right|\le\ K_2/(2m-1)\ \
&\hbox{\rm if}\ \ j_2=j_4=0\cr\cr
g_{j_4j_2}^{2m-1}+g_{j_2j_4}^{2m-1}=0\ \ &\hbox{\rm otherwise}
\end{matrix}\right.,
\end{equation}

\vspace{2mm}
\begin{equation}
\label{agentoo2}
\left\{
\begin{matrix}
\left|g_{j_4j_4}^{2m-1}\right|\le\ K_2/(2m-1)\ \
&\hbox{\rm if}\ \ j_4=0\cr\cr
g_{j_4j_4}^{2m-1}=0\ \ &\hbox{\rm otherwise}
\end{matrix}\right.,
\end{equation}

\vspace{3mm}
\noindent
where constant $K_2$ does not depend on $p=2m-1$.

The relations (\ref{otitkkk}), (\ref{otitkkk11}), (\ref{agentoo1}), and (\ref{agentoo2})
imply the following formulas

\vspace{-4mm}
\begin{equation}
\label{agentooo1}
\left\{
\begin{matrix}
\left|g_{j_4j_2}^{p}+g_{j_2j_4}^{p}\right|\le\ K_3/p\ \
&\hbox{\rm if}\ \ j_2=j_4=0\cr\cr
g_{j_4j_2}^{p}+g_{j_2j_4}^{p}=0\ \ &\hbox{\rm otherwise}
\end{matrix}\right.,
\end{equation}

\vspace{3mm}
\begin{equation}
\label{agentooo1rrr}
\left\{
\begin{matrix}
\left|g_{j_4j_4}^{p}\right|\le\ K_3/p\ \
&\hbox{\rm if}\ \ j_4=0\cr\cr
g_{j_4j_4}^{p}=0\ \ &\hbox{\rm otherwise}
\end{matrix}\right.,
\end{equation}

\vspace{2mm}
\noindent
where constant $K_3$ does not depend on $p$ $(p\in {\bf N}).$
Moreover, $g_{j_4j_4}^{p}\ge 0$ (see (\ref{otitddd})).

From (\ref{otit321}), (\ref{agentooo1}), and (\ref{agentooo1rrr})
it follows that $\Delta_7^{(i_2i_4)}=0$ and 
$\Delta_3^{(i_2i_4)}=0$\ \ w.~p.~1 for $i_2=i_4\ne 0.$
Analogously to the polynomial case, we obtain 
$\Delta_7^{(i_2i_4)}=0$ and 
$\Delta_3^{(i_2i_4)}=0$\ \ w.~p.~1 for $i_2\ne i_4,$\
$i_2\ne 0,$\ $i_4\ne 0.$
The similar arguments prove that $\Delta_5^{(i_1i_3)}=0$\ \ w.~p.~1.

Taking into account (\ref{otitbb}), (\ref{agentooo1}), (\ref{agentooo1rrr}) and the relations
$$
\lim\limits_{p\to\infty}\sum\limits_{j_3=0}^p f_{j_3j_3}^p=
\lim\limits_{p\to\infty}\sum\limits_{j_3=0}^p d_{j_3j_3}^p=0,
$$

\noindent
which follow from the estimates 
\begin{equation}
\label{agentys1}
~~~~~~~~ |f_{jj}^p|\le \frac{C_1}{pj},\ \ \ 
|d_{jj}^p|\le \frac{C_1}{pj}\ \ \ (j\ne 0),\ \ \ \ 
|f_{00}^p|\le \frac{C_1}{p},\ \ \ 
|d_{00}^p|\le \frac{C_1}{p},
\end{equation}

\noindent
we obtain
$$
\lim\limits_{p\to\infty}\sum\limits_{j_3=0}^p c_{j_3j_3}^p=
-\lim\limits_{p\to\infty}\sum\limits_{j_3=0}^p g_{j_3j_3}^p,
$$
$$
0\le \lim\limits_{p\to\infty}\sum\limits_{j_3=0}^p g_{j_3j_3}^p
\le\lim\limits_{p\to\infty}\frac{K_3}{p}=0.
$$

\vspace{1mm}

Note that the estimates (\ref{agentys1}) can be obtained by analogy with (\ref{otithj});
constant $C_1$ in (\ref{agentys1}) has the same meaning as constant $C_1$ in (\ref{otithj}).

Finally, we have 
$$
\lim\limits_{p\to\infty}\sum\limits_{j_3=0}^p c_{j_3j_3}^p=0.
$$

The relations (\ref{otitgggh-1}), (\ref{otitgggh}) are proved for the trigonometric case. 
Theorem 2.9 is proved
for the trigonometric case. Theorem 2.9 is proved. 

\vspace{2mm}

{\bf Remark~2.2.}\ {\it It should be noted that the proof of Theorem~{\rm 2.9} 
can be somewhat simplified.
More precisely, instead of {\rm (\ref{otit321})--(\ref{otit32101}),} 
we can use only one and rather simple estimate.

We have
$$
{\sf M}\left\{\left(\sum\limits_{j_3, j_4=0}^{p}
a_{j_4 j_3}^p \zeta_{j_3}^{(i_3)}
\zeta_{j_4}^{(i_4)}\right)^2\right\}=
$$
$$
={\sf M}\left\{\left(\sum\limits_{j_3, j_4=0}^{p}
a_{j_4 j_3}^p \biggl(\zeta_{j_3}^{(i_3)}
\zeta_{j_4}^{(i_4)}-{\bf 1}_{\{i_3=i_4\ne 0\}}{\bf 1}_{\{j_3=j_4\}}
+{\bf 1}_{\{i_3=i_4\ne 0\}}{\bf 1}_{\{j_3=j_4\}}\biggr)
\right)^2\right\}=
$$
$$
={\sf M}\left\{\left(\sum\limits_{j_3, j_4=0}^{p}
a_{j_4 j_3}^p \biggl(\zeta_{j_3}^{(i_3)}
\zeta_{j_4}^{(i_4)}-{\bf 1}_{\{i_3=i_4\ne 0\}}{\bf 1}_{\{j_3=j_4\}}\biggr)
+{\bf 1}_{\{i_3=i_4\ne 0\}}\sum\limits_{j_4=0}^{p}
a_{j_4 j_4}^p\right)^2\right\}=
$$
$$
={\sf M}\left\{\left(\sum\limits_{j_3, j_4=0}^{p}
a_{j_4 j_3}^p \biggl(\zeta_{j_3}^{(i_3)}
\zeta_{j_4}^{(i_4)}-{\bf 1}_{\{i_3=i_4\ne 0\}}{\bf 1}_{\{j_3=j_4\}}\biggr)
\right)^2\right\}
+
$$
\begin{equation}
\label{riss0}
+{\bf 1}_{\{i_3=i_4\ne 0\}}\left(\sum\limits_{j_4=0}^{p}
a_{j_4 j_4}^p\right)^2.
\end{equation}

The expression 
$$
\sum\limits_{j_3, j_4=0}^{p}
a_{j_4 j_3}^p \biggl(\zeta_{j_3}^{(i_3)}
\zeta_{j_4}^{(i_4)}-{\bf 1}_{\{i_3=i_4\ne 0\}}{\bf 1}_{\{j_3=j_4\}}\biggr)
$$
can be interpreted as the multiple Wiener stochastic integral 
{\rm (\ref{mult11www})} {\rm (}also see {\rm (\ref{mult11}))} 
of multiplicity {\rm 2}
with nonrandom integrand function
$$
\sum\limits_{j_3, j_4=0}^{p}
a_{j_4 j_3}^p \phi_{j_3}(t_3)\phi_{j_4}(t_4).
$$

From {\rm (\ref{s2s})} we obtain
$$
{\sf M}\left\{\left(J'[\Phi]_{T,t}^{(k)}\right)^2\right\}\le
C_k
\sum_{(t_1,\ldots,t_k)}
\int\limits_{t}^{T}
\ldots
\int\limits_{t}^{t_2}
\Phi^2(t_1,\ldots,t_k)dt_1\ldots dt_k=
$$
\begin{equation}
\label{wiener1}
= C_k
\int\limits_{[t,T]^k}
\Phi^2(t_1,\ldots,t_k)dt_1\ldots dt_k,
\end{equation}

\vspace{2mm}
\noindent
where $J'[\Phi]_{T,t}^{(k)}$ is defined by {\rm (\ref{mult11})}
and $C_k$ is a constant.

Then
$$
{\sf M}\left\{\left(\sum\limits_{j_3, j_4=0}^{p}
a_{j_4 j_3}^p \biggl(\zeta_{j_3}^{(i_3)}
\zeta_{j_4}^{(i_4)}-{\bf 1}_{\{i_3=i_4\ne 0\}}{\bf 1}_{\{j_3=j_4\}}\biggr)
\right)^2\right\}\le 
$$
$$
\le C_2 \int\limits_{[t,T]^2}
\left(\sum\limits_{j_3, j_4=0}^{p}
a_{j_4 j_3}^p \phi_{j_3}(t_3)\phi_{j_4}(t_4)\right)^2 dt_3 dt_4=
$$
\begin{equation}
\label{riss1}
=C_2 \sum\limits_{j_3, j_4=0}^{p}
\left(a_{j_4 j_3}^p\right)^2.
\end{equation}

\vspace{3mm}

From {\rm (\ref{riss0})} and {\rm (\ref{riss1})} we get
\begin{equation}
\label{riss2}
{\sf M}\left\{\left(\sum\limits_{j_3, j_4=0}^{p}
a_{j_4 j_3}^p \zeta_{j_3}^{(i_3)}
\zeta_{j_4}^{(i_4)}\right)^2\right\}\le 
C_2 \sum\limits_{j_3, j_4=0}^{p}
\left(a_{j_4 j_3}^p\right)^2 +
{\bf 1}_{\{i_3=i_4\ne 0\}}\left(\sum\limits_{j_4=0}^{p}
a_{j_4 j_4}^p\right)^2.
\end{equation}

\vspace{2mm}

Obviously, the estimate {\rm (\ref{riss2})} can be used in the proof of Theorem~{\rm 2.9}
instead of {\rm (\ref{otit321})--(\ref{otit32101}).}

The estimate {\rm (\ref{riss2})} can be refined. Using {\rm (\ref{gt22}),} we obtain

\vspace{-3mm}
$$
{\sf M}\left\{\left(\sum\limits_{j_3, j_4=0}^{p}
a_{j_4 j_3}^p \biggl(\zeta_{j_3}^{(i_3)}
\zeta_{j_4}^{(i_4)}-{\bf 1}_{\{i_3=i_4\ne 0\}}{\bf 1}_{\{j_3=j_4\}}\biggr)
\right)^2\right\}=
$$

\newpage
\noindent
$$
=\sum\limits_{j_3, j_4=0}^{p}
\left(a_{j_4 j_3}^p\right)^2 +
{\bf 1}_{\{i_3=i_4\ne 0\}}
\sum\limits_{j_3, j_4=0}^{p}
a_{j_4 j_3}^p a_{j_3 j_4}^p\le
$$

\vspace{-2mm}
$$
\le \sum\limits_{j_3, j_4=0}^{p}
\left(a_{j_4 j_3}^p\right)^2 +
{\bf 1}_{\{i_3=i_4\ne 0\}}\frac{1}{2}
\sum\limits_{j_3, j_4=0}^{p}
\left(\left(a_{j_4 j_3}^p\right)^2+ \left(a_{j_3 j_4}^p\right)^2\right)=
$$
\begin{equation}
\label{riss5}
=\left(1+{\bf 1}_{\{i_3=i_4\ne 0\}}\right) 
\sum\limits_{j_3, j_4=0}^{p}
\left(a_{j_4 j_3}^p\right)^2.
\end{equation}

\vspace{2mm}

Combining {\rm (\ref{riss0})} and {\rm (\ref{riss5}),} we have
$$
{\sf M}\left\{\left(\sum\limits_{j_3, j_4=0}^{p}
a_{j_4 j_3}^p \zeta_{j_3}^{(i_3)}
\zeta_{j_4}^{(i_4)} 
\right)^2\right\}\le\left(1+{\bf 1}_{\{i_3=i_4\ne 0\}}\right) 
\sum\limits_{j_3, j_4=0}^{p}
\left(a_{j_4 j_3}^p\right)^2+
$$
\begin{equation}
\label{riss7}
+
{\bf 1}_{\{i_3=i_4\ne 0\}}\left(\sum\limits_{j_4=0}^{p}
a_{j_4 j_4}^p\right)^2.
\end{equation}

}

\section{Expansion of Iterated Stratonovich Stochastic Integrals of 
Multiplicity $k$ $(k\in {\bf N})$ Based on Generalized 
Iterated Fourier Series Converging Pointwise}

This section is devoted to the 
expansion of iterated Stratonovich stochastic integrals of arbitrary 
multiplicity $k$ $(k\in{\bf N})$ based on 
generalized iterated Fourier series. The case of 
trigonometric Fourier series are considered in detail. 
The obtained expansion provides a possibility to 
represent
the iterated Stratonovich sto\-chas\-tic integral in the form of iterated 
series 
of products of standard Gaussian random variables. Convergence
in the mean of degree $q=2n$ $(n\in {\bf N})$ of the expansion is proved.
The case of iterated Fourier--Legendre series for $k=2$ and $q=2$ is also considered.

The idea of representing of iterated Stratonovich stochastic 
integrals in the form of multiple stochastic integrals from 
specific discontinuous nonrandom functions of several variables and following 
expansion of these functions using generalized iterated Fourier series in order 
to get effective mean-square approximations of the mentioned stochastic 
integrals was proposed and developed in a lot of 
author's publications \cite{old-art-1} (1997), \cite{old-art-2} (1998)
(also see \cite{5}-\cite{12aa}, \cite{arxiv-6}).
The results 
of this section
convincingly testify
that there is a doubtless relation between 
the multiplier factor $1/2$, 
which is typical for Stratonovich stochastic integral and included 
into the sum connecting Stratonovich and It\^{o} stochastic integrals, 
and the fact that in the point of finite discontinuity of piecewise
smooth function $f(x)$ its trigonometric Fourier series 
and Fourier--Legendre series
converge to the value $(f(x+0)+f(x-0))/2.$

\subsection{Theorem on Expansion of 
Iterated Stratonovich Stochastic Integrals of 
Multiplicity $k$ $(k\in {\bf N})$}

Consider the following iterated Stratonovich and It\^{o} stochastic integrals
\begin{equation}
\label{str}
~~~~~~~~~~~J^{*}[\psi^{(k)}]_{T,t}=
{\int\limits_t^{*}}^T
\psi_k(t_k) \ldots 
{\int\limits_t^{*}}^{t_2}
\psi_1(t_1) d{\bf w}_{t_1}^{(i_1)}\ldots
d{\bf w}_{t_k}^{(i_k)},
\end{equation}
\begin{equation}
\label{itoxxx}
J[\psi^{(k)}]_{T,t}=\int\limits_t^T\psi_k(t_k) \ldots \int\limits_t^{t_{2}}
\psi_1(t_1) d{\bf w}_{t_1}^{(i_1)}\ldots
d{\bf w}_{t_k}^{(i_k)},
\end{equation}
\noindent
where $\psi_l(\tau)$ $(l=1,\ldots,k)$ are
nonrandom functions 
on $[t,T],$ ${\bf w}_{\tau}^{(i)}$
$(i=1,\ldots,m)$ are independent standard Wiener processes,
${\bf w}_{\tau}^{(0)}=\tau,$
$i_1,\ldots,i_k = 0, 1,\ldots,m.$

Let us
denote as $\{\phi_j(x)\}_{j=0}^{\infty}$ the
complete orthonormal systems of Legendre polynomials or
trigonometric functions 
in the space $L_2([t, T])$.

In this section, we will pay attention 
on the 
well known facts about Fourier series with respect to
these two systems of functions \cite{Gob} (also see Sect.~2.1.1).

Define the following function on the hypercube $[t, T]^k$
\begin{equation}
\label{pppxx}
K(t_1,\ldots,t_k)=
\left\{
\begin{matrix}
\psi_1(t_1)\ldots \psi_k(t_k),\ &t_1<\ldots<t_k\cr\cr
0,\ &\hbox{\rm otherwise}
\end{matrix}\right.\ \ \
=\ \ 
\prod\limits_{l=1}^k
\psi_l(t_l)\ \prod\limits_{l=1}^{k-1}{\bf 1}_{\{t_l<t_{l+1}\}}
\end{equation}

\noindent
for $t_1,\ldots,t_k\in [t, T]$ $(k\ge 2)$ and 
$K(t_1)\equiv\psi_1(t_1)$ for $t_1\in[t, T],$ where
${\bf 1}_A$ denotes the indicator of the set $A$.

Let us formulate the following theorem.

{\bf Theorem 2.10} \cite{old-art-1} (1997), \cite{old-art-2} (1998)
(also see \cite{5}-\cite{12aa}, \cite{arxiv-6}).
{\it Suppose that every function $\psi_l(\tau)$ $(l=1,\ldots,k)$ is twice continuously
differentiable at the interval
$[t, T]$ and
$\{\phi_j(x)\}_{j=0}^{\infty}$ is a complete
orthonormal system of trigonometric functions in the space $L_2([t, T])$. 
Then$,$ the iterated Stratonovich stochastic integral 
$J^{*}[\psi^{(k)}]_{T,t}$ defined by {\rm(\ref{str})}
is expanded into the 
conver\-ging 
in the mean of degree $2n$ $(n\in {\bf N})$
iterated series
\begin{equation}
\label{1500}
J^{*}[\psi^{(k)}]_{T,t}=
\sum_{j_1=0}^{\infty}\ldots\sum_{j_k=0}^{\infty}
C_{j_k\ldots j_1}
\prod_{l=1}^k
\zeta^{(i_l)}_{j_l},
\end{equation}

\noindent
where 
$$
\zeta_{j}^{(i)}=
\int\limits_t^T \phi_{j}(s) d{\bf w}_s^{(i)}
$$ 
are independent standard Gaussian random variables
for
various
$i$ or $j$ {\rm(}in the case when $i\ne 0${\rm)} and 
\begin{equation}
\label{333.40}
C_{j_k\ldots j_1}=\int\limits_{[t,T]^k}
K(t_1,\ldots,t_k)\prod_{l=1}^{k}\phi_{j_l}(t_l)dt_1\ldots dt_k
\end{equation}

\noindent
is the Fourier coefficient.}

Note that (\ref{1500})  means the following 
\begin{equation}
\label{1500e}
\lim\limits_{p_1\to\infty}
\varlimsup\limits_{p_2\to\infty}
\ldots\varlimsup\limits_{p_k\to\infty}
{\sf M}\left\{\left(J^{*}[\psi^{(k)}]_{T,t}-
\sum_{j_1=0}^{p_1}\ldots\sum_{j_k=0}^{p_k}
C_{j_k\ldots j_1}
\prod_{l=1}^k
\zeta^{(i_l)}_{j_l}\right)^{2n}\right\}=0,
\end{equation}
where $\varlimsup$ means ${\rm lim\  sup.}$

{\bf Proof.} The proof of Theorem 2.10 is based
on Lemmas 1.1, 1.3 (see Sect.~1.1.3) and Theorems 2.11--2.13 (see below).

Define the function $K^{*}(t_1,\ldots,t_k)$ on the hypercube 
$[t,T]^k$ as follows

\vspace{-1mm}
$$
K^{*}(t_1,\ldots,t_k)=\prod\limits_{l=1}^k\psi_l(t_l)
\prod_{l=1}^{k-1}\biggl({\bf 1}_{\{t_l<t_{l+1}\}}+
\frac{1}{2}{\bf 1}_{\{t_l=t_{l+1}\}}\biggr)=
$$

\vspace{-1mm}
\begin{equation}
\label{1999.1xx}
=\prod_{l=1}^k \psi_l(t_l)\left(\prod_{l=1}^{k-1}
{\bf 1}_{\{t_l<t_{l+1}\}}+
\sum_{r=1}^{k-1}\frac{1}{2^r}
\sum_{\stackrel{s_r,\ldots,s_1=1}{{}_{s_r>\ldots>s_1}}}^{k-1}\ 
\prod_{l=1}^r {\bf 1}_{\{t_{s_l}=t_{s_l+1}\}}
\prod_{\stackrel{l=1}{{}_{l\ne s_1,\ldots, s_r}}}^{k-1}
{\bf 1}_{\{t_{l}<t_{l+1}\}}\right)
\end{equation}

\vspace{3mm}
\noindent
for $t_1,\ldots,t_k\in[t, T]$ $(k\ge 2)$ and 
$K^{*}(t_1)\equiv\psi_1(t_1)$ for $t_1\in[t, T],$ 
where ${\bf 1}_A$ is the indicator of the set $A$.

{\bf Theorem 2.11} \cite{old-art-1} (1997).
{\it Suppose that every function $\psi_l(\tau)$ $(l=1,\ldots,k)$ is continuously
differentiable at the interval
$[t, T]$ and
$\{\phi_j(x)\}_{j=0}^{\infty}$ is a complete
orthonormal system of Legendre 
polynomials or trigonometric functions in the space $L_2([t, T])$. 
Then$,$ 
the function
$K^{*}(t_1,\ldots,t_k)$
is represented in any internal point of the hypercube   
$[t,T]^k$ by the generalized iterated Fourier series
\begin{equation}
\label{30.18}
~~~~~~~~~ K^{*}(t_1,\ldots,t_k)=
\lim\limits_{p_1\to\infty}\ldots\lim\limits_{p_k\to\infty}
\sum_{j_1=0}^{p_1}\ldots\sum_{j_k=0}^{p_k}
C_{j_k\ldots j_1}\prod_{l=1}^{k} \phi_{j_l}(t_l),
\end{equation}

\noindent
where 
$(t_1,\ldots,t_k)\in (t, T)^k$ and
$C_{j_k\ldots j_1}$ is defined by {\rm (\ref{333.40})}.
At that, the iterated series {\rm (\ref{30.18})} converges at the 
boundary
of the hypercube $[t,T]^k$
{\rm (}not necessarily to the function $K^{*}(t_1,\ldots,t_k)${\rm )}.}

{\bf Proof.}~We will perform the proof using induction.~Consider 
the case $k=2.$ Let us expand the function 
$K^{*}(t_1,t_2)$ using the variable 
$t_1$, when $t_2$ is fixed, into the generalized Fourier series 
with respect to the system $\{\phi_j(x)\}_{j=0}^{\infty}$
at the interval $(t, T)$ 
\begin{equation}
\label{leto8001}
K^{*}(t_1,t_2)=
\sum_{j_1=0}^{\infty}C_{j_1}(t_2)\phi_{j_1}(t_1)\ \ \ (t_1\ne t, T),
\end{equation}
where
$$
C_{j_1}(t_2)=\int\limits_t^T
K^{*}(t_1,t_2)\phi_{j_1}(t_1)dt_1=
\psi_2(t_2)
\int\limits_t^{t_2}\psi_1(t_1)\phi_{j_1}(t_1)dt_1.
$$

The equality (\ref{leto8001}) is 
satisfied
pointwise at each point of the interval $(t, T)$ with respect to the 
variable $t_1$, when $t_2\in [t, T]$ is fixed, due to 
a piecewise
smoothness of the function $K^{*}(t_1,t_2)$ with respect to the variable 
$t_1\in [t, T]$ ($t_2$ is fixed).

Note also that due to the well known properties of the 
Fourier--Legendre series
and trigonometric Fourier series, 
the series (\ref{leto8001}) converges when $t_1=t, T$ 
{\rm (}not necessarily to the function $K^{*}(t_1,t_2)${\rm )}.

Obtaining (\ref{leto8001}), we also used the fact that the right-hand side 
of (\ref{leto8001}) converges when $t_1=t_2$ (point of a finite
discontinuity
of the function $K(t_1,t_2)$ defined by (\ref{pppxx})) to the value
$$
\frac{1}{2}\left(K(t_2-0,t_2)+K(t_2+0,t_2)\right)=
\frac{1}{2}\psi_1(t_2)\psi_2(t_2)=
K^{*}(t_2,t_2).
$$

The function $C_{j_1}(t_2)$ is continuously differentiable
at the interval $[t, T]$. 
Let us expand it into the generalized Fourier series at the interval $(t, T)$

\newpage
\noindent
\begin{equation}
\label{leto8002}
C_{j_1}(t_2)=
\sum_{j_2=0}^{\infty}C_{j_2 j_1}\phi_{j_2}(t_2)\ \ \ (t_2\ne t, T),
\end{equation}
where 
$$
C_{j_2 j_1}=\int\limits_t^T
C_{j_1}(t_2)\phi_{j_2}(t_2)dt_2=
\int\limits_t^T
\psi_2(t_2)\phi_{j_2}(t_2)\int\limits_t^{t_2}
\psi_1(t_1)\phi_{j_1}(t_1)dt_1 dt_2
$$
and the equality (\ref{leto8002}) is satisfied pointwise at any point 
of the interval $(t, T)$. Moreover, the right-hand side 
of
(\ref{leto8002}) converges 
when $t_2=t, T$ (not necessarily to $C_{j_1}(t_2)$).

Let us substitute (\ref{leto8002}) into (\ref{leto8001})
\begin{equation}
\label{leto8003}
~~~~~~~ K^{*}(t_1,t_2)=
\sum_{j_1=0}^{\infty}\sum_{j_2=0}^{\infty}C_{j_2 j_1}
\phi_{j_1}(t_1)\phi_{j_2}(t_2),\ \ \ (t_1, t_2)\in (t, T)^2.
\end{equation}

Note that 
the series on the right-hand side of (\ref{leto8003}) converges at the 
boundary
of the square  $[t, T]^2$ (not necessarily to $K^{*}(t_1,t_2)$).
Theorem 2.11 is proved for the case $k=2.$

Note that proving Theorem 2.11 for the case $k=2$ we obtained the 
following equality (see (\ref{leto8001}))
\begin{equation}
\label{oop1}
~~~~~~\psi_1(t_1)\left({\bf 1}_{\{t_1<t_2\}}+
\frac{1}{2}{\bf 1}_{\{t_1=t_2\}}\right)=
\sum\limits_{j_1=0}^{\infty}\int\limits_{t}^{t_2}\psi_1(t_1)
\phi_{j_1}(t_1)dt_1\phi_{j_1}(t_1),
\end{equation}
which is 
satisfied pointwise at the interval $(t, T),$
besides
the series on the right-hand side 
of (\ref{oop1}) converges when $t_1=t, T.$

Let us introduce the induction assumption 
$$
\sum\limits_{j_1=0}^{\infty}\sum\limits_{j_2=0}^{\infty}\ldots
\sum\limits_{j_{k-2}=0}^{\infty}\psi_{k-1}(t_{k-1})\times
$$
$$
\times
\int\limits_t^{t_{k-1}}\psi_{k-2}(t_{k-2})\phi_{j_{k-2}}(t_{k-2})
\ldots\int\limits_t^{t_{2}}
\psi_{1}(t_{1})\phi_{j_{1}}(t_{1})dt_1\ldots dt_{k-2}
\prod_{l=1}^{k-2}\phi_{j_{l}}(t_{l})=
$$
\begin{equation}
\label{oop22}
=\prod\limits_{l=1}^{k-1}\psi_l(t_l)
\prod_{l=1}^{k-2}\left({\bf 1}_{\{t_l<t_{l+1}\}}+
\frac{1}{2}{\bf 1}_{\{t_l=t_{l+1}\}}\right).
\end{equation}

Then
$$
\sum\limits_{j_1=0}^{\infty}\sum\limits_{j_2=0}^{\infty}\ldots
\sum\limits_{j_{k-1}=0}^{\infty}\psi_{k}(t_{k})\times
$$

\vspace{-4mm}
$$
\times
\int\limits_t^{t_{k}}\psi_{k-1}(t_{k-1})\phi_{j_{k-1}}(t_{k-1})\ldots
\int\limits_t^{t_{2}}\psi_{1}(t_{1})
\phi_{j_{1}}(t_{1})
dt_1\ldots dt_{k-1}\prod_{l=1}^{k-1}\phi_{j_{l}}(t_{l})=
$$

\vspace{-4mm}
$$
=\sum\limits_{j_1=0}^{\infty}\sum\limits_{j_2=0}^{\infty}\ldots
\sum\limits_{j_{k-2}=0}^{\infty}
\psi_k(t_k)\left({\bf 1}_{\{t_{k-1}<t_{k}\}}+
\frac{1}{2}{\bf 1}_{\{t_{k-1}=t_{k}\}}\right)\psi_{k-1}(t_{k-1})\times
$$

\vspace{-4mm}
$$
\times
\int\limits_t^{t_{k-1}}\psi_{k-2}(t_{k-2})\phi_{j_{k-2}}(t_{k-2})\ldots
\int\limits_t^{t_{2}}\psi_{1}(t_{1})
\phi_{j_{1}}(t_{1})dt_1\ldots dt_{k-2}
\prod_{l=1}^{k-2}\phi_{j_{l}}(t_{l})=
$$

\vspace{-4mm}
$$
=\psi_k(t_k)\left({\bf 1}_{\{t_{k-1}<t_{k}\}}+
\frac{1}{2}{\bf 1}_{\{t_{k-1}=t_{k}\}}\right)
\sum\limits_{j_1=0}^{\infty}\sum\limits_{j_2=0}^{\infty}\ldots
\sum\limits_{j_{k-2}=0}^{\infty}
\psi_{k-1}(t_{k-1})\times
$$

\vspace{-4mm}
$$
\times
\int\limits_t^{t_{k-1}}\psi_{k-2}(t_{k-2})\phi_{j_{k-2}}(t_{k-2})\ldots
\int\limits_t^{t_{2}}\psi_{1}(t_{1})
\phi_{j_{1}}(t_{1})dt_1\ldots dt_{k-2}
\prod_{l=1}^{k-2}\phi_{j_{l}}(t_{l})=
$$

\vspace{-4mm}
$$
=\psi_k(t_k)\left({\bf 1}_{\{t_{k-1}<t_{k}\}}+
\frac{1}{2}{\bf 1}_{\{t_{k-1}=t_{k}\}}\right)
\prod\limits_{l=1}^{k-1}\psi_l(t_l)
\prod_{l=1}^{k-2}\left({\bf 1}_{\{t_l<t_{l+1}\}}+
\frac{1}{2}{\bf 1}_{\{t_l=t_{l+1}\}}\right)=
$$

\vspace{-1mm}
\begin{equation}
\label{oop30}
=\prod\limits_{l=1}^{k}\psi_l(t_l)
\prod_{l=1}^{k-1}\left({\bf 1}_{\{t_l<t_{l+1}\}}+
\frac{1}{2}{\bf 1}_{\{t_l=t_{l+1}\}}\right).
\end{equation}

\vspace{2mm}

On the other hand, the left-hand side 
of (\ref{oop30}) can be represented 
in the following form
$$
\sum_{j_1=0}^{\infty}\ldots \sum_{j_k=0}^{\infty}
C_{j_k\ldots j_1}\prod_{l=1}^{k} \phi_{j_l}(t_l)
$$
by 
expanding the function
$$
\psi_{k}(t_{k})
\int\limits_t^{t_{k}}\psi_{k-1}(t_{k-1})\phi_{j_{k-1}}(t_{k-1})\ldots
\int\limits_t^{t_{2}}\psi_{1}(t_{1})
\phi_{j_{1}}(t_{1})
dt_1\ldots dt_{k-1}
$$

\noindent
into the generalized Fourier series at the interval $(t, T)$ 
using the variable 
$t_k$.
Theorem 2.11 is proved.

Let us introduce the following notations
$$
J[\psi^{(k)}]_{T,t}^{s_l,\ldots,s_1}\ \stackrel{\rm def}{=}\ 
\prod_{p=1}^l {\bf 1}_{\{i_{s_p}=
i_{s_{p}+1}\ne 0\}}\ \times
$$
$$
\times
\int\limits_t^T\psi_k(t_k)\ldots \int\limits_t^{t_{s_l+3}}
\psi_{s_l+2}(t_{s_l+2})
\int\limits_t^{t_{s_l+2}}\psi_{s_l}(t_{s_l+1})
\psi_{s_l+1}(t_{s_l+1}) \times
$$
$$
\times
\int\limits_t^{t_{s_l+1}}\psi_{s_l-1}(t_{s_l-1})
\ldots
\int\limits_t^{t_{s_1+3}}\psi_{s_1+2}(t_{s_1+2})
\int\limits_t^{t_{s_1+2}}\psi_{s_1}(t_{s_1+1})
\psi_{s_1+1}(t_{s_1+1}) \times
$$
$$
\times
\int\limits_t^{t_{s_1+1}}\psi_{s_1-1}(t_{s_1-1})
\ldots \int\limits_t^{t_2}\psi_1(t_1)
d{\bf w}_{t_1}^{(i_1)}\ldots d{\bf w}_{t_{s_1-1}}^{(i_{s_1-1})}
dt_{s_1+1}d{\bf w}_{t_{s_1+2}}^{(i_{s_1+2})}\ldots
$$

\begin{equation}
\label{30.1}
\ldots\
d{\bf w}_{t_{s_l-1}}^{(i_{s_l-1})}
dt_{s_l+1}d{\bf w}_{t_{s_l+2}}^{(i_{s_l+2})}\ldots d{\bf w}_{t_k}^{(i_k)},
\end{equation}

\noindent
where 
\begin{equation}
\label{30.5550001}
{\rm A}_{k,l}
=\bigl\{(s_l,\ldots,s_1):\
s_l>s_{l-1}+1,\ldots,s_2>s_1+1,\ s_l,\ldots,s_1=1,\ldots,k-1\bigr\},
\end{equation}
$$
(s_l,\ldots,s_1)\in{\rm A}_{k,l},\ \ \ 
l=1,\ldots,\left[k/2\right],\ \ \
i_s=0, 1,\ldots,m,\ \ \
s=1,\ldots,k,
$$

\vspace{2mm}
\noindent
$[x]$ is an
integer
part of a real number $x,$
and ${\bf 1}_A$ is the indicator of the set $A$.

Let us formulate the statement on connection 
between
iterated 
Stra\-to\-no\-vich and It\^{o} stochastic integrals 
$J^{*}[\psi^{(k)}]_{T,t},$ $J[\psi^{(k)}]_{T,t}$ 
of fixed multiplicity $k,$ $k\in{\bf N}$ (see (\ref{str}), (\ref{itoxxx})).

{\bf Theorem 2.12} \cite{old-art-1} (1997). {\it Suppose that
every $\psi_l(\tau)$ $(l=1,\ldots,k)$ is a continuous
function at the interval $[t, T]$.
Then$,$ the following relation between iterated
Stra\-to\-no\-vich and It\^{o} stochastic integrals 
\begin{equation}
\label{30.4}
~~~~J^{*}[\psi^{(k)}]_{T,t}=J[\psi^{(k)}]_{T,t}+
\sum_{r=1}^{\left[k/2\right]}\frac{1}{2^r}
\sum_{(s_r,\ldots,s_1)\in {\rm A}_{k,r}}
J[\psi^{(k)}]_{T,t}^{s_r,\ldots,s_1}\ \ \ \hbox{{\rm w.~p.~1}}
\end{equation}

\noindent
is correct, 
where $\sum\limits_{\emptyset}$ is supposed to be equal to zero{\rm .}
}

\newpage
\noindent
\par
{\bf Proof.} Let us prove the equality (\ref{30.4}) using induction. 
The case $k=1$ is obvious.
If $k=2,$ then from (\ref{30.4}) we get 

\vspace{-2mm}
\begin{equation}
\label{30.6}
J^{*}[\psi^{(2)}]_{T,t}=J[\psi^{(2)}]_{T,t}+
\frac{1}{2}J[\psi^{(2)}]_{T,t}^{1}\ \ \ \hbox{w.~p.~1}.
\end{equation}

\vspace{2mm}

Let us demonstrate that the equality (\ref{30.6}) is correct 
w.~p.~1. In order to do it let us consider the 
function $F(x,\tau)=x\psi_2(\tau)$ and the
process $F(\eta_{\tau,t},\tau),$ where
$\eta_{\tau,t}=J[\psi^{(1)}]_{\tau,t},$ 
$\tau\in[t, T]$. Then 
\begin{equation}
\label{30.7}
\frac{\partial F}{\partial x}(x,\tau)=\psi_2(\tau),\ \ \ 
d\eta_{\tau,t}=\psi_1(\tau)d{\bf w}_{\tau}^{(i_1)}.
\end{equation}

From (\ref{30.7}) we obtain that the diffusion 
coefficient of the process $\eta_{\tau,t},$ $\tau \in [t, T]$ equals 
to 
${\bf 1}_{\{i_1\ne 0\}}\psi_1(\tau).$
Further, using the standard relations between 
Stratonovich and It\^{o} stochastic integrals (see
(\ref{d11}), (\ref{d11a})), 
we obtain the relation (\ref{30.6}). 
Thus, the statement
of Theorem 2.12 is proved for $k=1$ and $k=2.$

Assume that the statement of Theorem 2.12 is correct
for some integer $k$ ($k>2$). Let us prove its correctness when 
the value $k$ is greater per unit. Using the 
induction assumption,
we have w.~p.~1 

\vspace{-2mm}
$$
J^{*}[\psi^{(k+1)}]_{T,t}=
$$

\vspace{-3mm}
$$
=
{\int\limits_t^{*}}^T \psi_{k+1}(\tau)
\left(J[\psi^{(k)}]_{\tau,t}+
\sum_{r=1}^{\left[k/2\right]}\frac{1}{2^r}
\sum_{(s_r,\ldots,s_1)\in {\rm A}_{k,r}}
J[\psi^{(k)}]_{\tau,t}^{s_r,\ldots,s_1}
\right)
d{\bf w}_{\tau}^{(i_{k+1})}=
$$
$$
={\int\limits_t^{*}}^T \psi_{k+1}(\tau)
J[\psi^{(k)}]_{\tau,t}d{\bf w}_{\tau}^{(i_{k+1})}+
$$
\begin{equation}
\label{30.8}
~~~~~~~~ +
\sum_{r=1}^{\left[k/2\right]}\frac{1}{2^r}
\sum_{(s_r,\ldots,s_1)\in {\rm A}_{k,r}}
{\int\limits_t^{*}}^T\psi_{k+1}(\tau)
J[\psi^{(k)}]_{\tau,t}^{s_r,\ldots,s_1}d{\bf w}_{\tau}^{(i_{k+1})}.
\end{equation}

\vspace{4mm}

Using the standard relations
between Stratonovich and It\^{o} stochastic integrals
(see (\ref{d11}), (\ref{d11a})), similarly to (\ref{30.6}), we get
w.~p.~1

\vspace{-5mm}
\begin{equation}
\label{30.9}
~~~~~~~~~ {\int\limits_t^{*}}^T \psi_{k+1}(\tau)
J[\psi^{(k)}]_{\tau,t}d{\bf w}_{\tau}^{(i_{k+1})}=J[\psi^{(k+1)}]_{T,t}+
\frac{1}{2}J[\psi^{(k+1)}]_{T,t}^{k},
\end{equation}

\newpage
\noindent
$$
{\int\limits_t^{*}}^T\psi_{k+1}(\tau)
J[\psi^{(k)}]_{\tau,t}^{s_r,\ldots,s_1}d{\bf w}_{\tau}^{(i_{k+1})}=
$$

\begin{equation}
\label{30.10}
~~~~=\left\{
\begin{matrix}
J[\psi^{(k+1)}]_{T,t}^{s_r,\ldots,s_1}\ &\hbox{if}\ \ s_r=k-1\cr\cr\cr
J[\psi^{(k+1)}]_{T,t}^{s_r,\ldots,s_1}+
J[\psi^{(k+1)}]_{T,t}^{k,s_r,
\ldots,s_1}/2\
&\hbox{if}\ \ s_r<k-1
\end{matrix}\ \right. .
\end{equation}

\vspace{4mm}

After substituting (\ref{30.9}) and (\ref{30.10}) into 
(\ref{30.8}) and 
regrouping 
of summands, we pass to the following relations, which are valid
w.~p.~1

\vspace{-5mm}
\begin{equation}
\label{30.11}
~~J^{*}[\psi^{(k+1)}]_{T,t}=J[\psi^{(k+1)}]_{T,t}+
\sum_{r=1}^{\left[k/2\right]}\frac{1}{2^r}
\sum_{(s_r,\ldots,s_1)\in {\rm A}_{k+1,r}}
J[\psi^{(k+1)}]_{T,t}^{s_r,\ldots,s_1} 
\end{equation}
when $k$ is even and

\vspace{-5mm}
\begin{equation}
\label{30.12}
~~J^{*}[\psi^{(k'+1)}]_{T,t}=J[\psi^{(k'+1)}]_{T,t}+
\sum_{r=1}^{\left[k'/2\right]+1}\frac{1}{2^r}
\sum_{(s_r,\ldots,s_1)\in {\rm A}_{k'+1,r}}
J[\psi^{(k'+1)}]_{T,t}^{s_r,\ldots,s_1} 
\end{equation}

\noindent
when $k'=k+1$ is 
uneven.

From (\ref{30.11}) and (\ref{30.12})
we have w.~p.~1 

\vspace{-3mm}
\begin{equation}
\label{30.13}
~~J^{*}[\psi^{(k+1)}]_{T,t}=J[\psi^{(k+1)}]_{T,t}+
\sum_{r=1}^{\left[(k+1)/2\right]}\frac{1}{2^r}
\sum_{(s_r,\ldots,s_1)\in {\rm A}_{k+1,r}}
J[\psi^{(k+1)}]_{T,t}^{s_r,\ldots,s_1}.
\end{equation}

Theorem 2.12 is proved.

For example, from Theorem 2.12 for $k=1, 2, 3, 4$ we obtain
the following well known equalities \cite{Zapad-3}, which are fulfilled
w.~p.~1

$$
{\int\limits_t^{*}}^T\psi_1(t_1)d{\bf w}_{t_1}^{(i_1)}=
\int\limits_t^T\psi_1(t_1)d{\bf w}_{t_1}^{(i_1)},
$$
$$
{\int\limits_t^{*}}^T\psi_2(t_2)
{\int\limits_t^{*}}^{t_2}\psi_1(t_1)d{\bf w}_{t_1}^{(i_1)}
d{\bf w}_{t_2}^{(i_2)}=
\int\limits_t^T\psi_2(t_2)
\int\limits_t^{t_2}\psi_1(t_1)d{\bf w}_{t_1}^{(i_1)}
d{\bf w}_{t_2}^{(i_2)}+
$$
\begin{equation}
\label{zido10x}
+ \frac{1}{2}{\bf 1}_{\{i_1=i_2\ne 0\}}
\int\limits_t^T\psi_2(t_2)\psi_1(t_2)dt_2,
\end{equation}

\vspace{-2mm}
$$
{\int\limits_t^{*}}^T\psi_3(t_3)\ldots
{\int\limits_t^{*}}^{t_2}\psi_1(t_1)d{\bf w}_{t_1}^{(i_1)}\ldots
d{\bf w}_{t_3}^{(i_3)}=\hspace{-1mm}
\int\limits_t^T\psi_3(t_3)\ldots\hspace{-0.5mm}
\int\limits_t^{t_2}\psi_1(t_1)d{\bf w}_{t_1}^{(i_1)}\ldots
d{\bf w}_{t_3}^{(i_3)}+
$$
$$
+ \frac{1}{2}{\bf 1}_{\{i_1=i_2\ne 0\}}
\int\limits_t^T\psi_3(t_3)
\int\limits_t^{t_3}\psi_2(t_2)\psi_1(t_2)dt_2
d{\bf w}_{t_3}^{(i_3)}+
$$
\begin{equation}
\label{uyes3}
+ \frac{1}{2}{\bf 1}_{\{i_2=i_3\ne 0\}}
\int\limits_t^T\psi_3(t_3)\psi_2(t_3)
\int\limits_t^{t_3}\psi_1(t_1)d{\bf w}_{t_1}^{(i_1)}
dt_3,
\end{equation}

\vspace{-1mm}
$$
{\int\limits_t^{*}}^T\psi_4(t_4)\ldots
{\int\limits_t^{*}}^{t_2}\psi_1(t_1)d{\bf w}_{t_1}^{(i_1)}\ldots
d{\bf w}_{t_4}^{(i_4)}=\hspace{-1mm}
\int\limits_t^T\psi_4(t_4)\ldots\hspace{-0.5mm}
\int\limits_t^{t_2}\psi_1(t_1)d{\bf w}_{t_1}^{(i_1)}\ldots
d{\bf w}_{t_4}^{(i_4)}+
$$
$$
+ \frac{1}{2}{\bf 1}_{\{i_1=i_2\ne 0\}}
\int\limits_t^T\psi_4(t_4)\int\limits_t^{t_4}
\psi_3(t_3)\int\limits_t^{t_3}\psi_1(t_2)\psi_2(t_2)dt_2
d{\bf w}_{t_3}^{(i_3)}d{\bf w}_{t_4}^{(i_4)} +
$$
$$
+ \frac{1}{2}{\bf 1}_{\{i_2=i_3\ne 0\}}
\int\limits_t^T\psi_4(t_4)\int\limits_t^{t_4}
\psi_3(t_3)\psi_2(t_3)\int\limits_t^{t_3}\psi_1(t_1)
d{\bf w}_{t_1}^{(i_1)}dt_3d{\bf w}_{t_4}^{(i_4)} +
$$
$$
+ \frac{1}{2}{\bf 1}_{\{i_3=i_4\ne 0\}}
\int\limits_t^T\psi_4(t_4)\psi_3(t_4)\int\limits_t^{t_4}
\psi_2(t_2)\int\limits_t^{t_2}\psi_1(t_1)
d{\bf w}_{t_1}^{(i_1)}d{\bf w}_{t_2}^{(i_2)}dt_4 +
$$
\begin{equation}
\label{uyes3may}
~~~~~~~ + \frac{1}{4}{\bf 1}_{\{i_1=i_2\ne 0\}}{\bf 1}_{\{i_3=i_4\ne 0\}}
\int\limits_t^T\psi_4(t_4)\psi_3(t_4)\int\limits_t^{t_4}
\psi_2(t_2)\psi_1(t_2)dt_2 dt_4.
\end{equation}

\vspace{3mm}

Let us consider
Lemma 1.1, definition of the multiple stochastic integral (\ref{30.34})
together with the 
formula (\ref{30.52}) when the function $\Phi(t_1,\ldots,t_k)$ is 
continuous in the open domain $D_k$ and bounded at its boundary as well as
Lemma 1.3 (see Sect.~1.1.3).
Substituting
(\ref{1999.1xx}) into (\ref{30.34}) and using 
Lemma 1.1,
(\ref{30.52}), and Theorem 2.12
it is easy to see that w.~p.~1

\newpage
\noindent
\begin{equation}
\label{pp1}
~~~J^{*}[\psi^{(k)}]_{T,t}=
J[\psi^{(k)}]_{T,t}+
\sum_{r=1}^{\left[k/2\right]}\frac{1}{2^r}
\sum_{(s_r,\ldots,s_1)\in {\rm A}_{k,r}}
J[\psi^{(k)}]_{T,t}^{s_r,\ldots,s_1}=
J[{K^{*}}]_{T,t}^{(k)},
\end{equation}

\vspace{2mm}
\noindent
where 
$J[K^{*}]_{T,t}^{(k)}$
is defined by (\ref{30.34}) and $K^{*}(t_1,\ldots,t_k)$ has the form (\ref{1999.1xx}).

Let us substitute the relation

\vspace{-3mm}
$$
K^{*}(t_1,\ldots,t_k)=
\sum_{j_1=0}^{p_1}\ldots\sum_{j_k=0}^{p_k}
C_{j_k\ldots j_1} \prod_{l=1}^{k} \phi_{j_l}(t_l)
+K^{*}(t_1,\ldots,t_k)-
$$

$$
-\sum_{j_1=0}^{p_1}\ldots\sum_{j_k=0}^{p_k}
C_{j_k\ldots j_1} \prod_{l=1}^{k} \phi_{j_l}(t_l)
$$

\vspace{3mm}
\noindent
into the right-hand side of (\ref{pp1}) (here 
we suppose that $p_1,\ldots,p_k<\infty$).
Then using Lemma 1.3 (see Sect.~1.1.3), we obtain
\begin{equation}
\label{proof1ggg}
~~~~~J^{*}[\psi^{(k)}]_{T,t}=
\sum_{j_1=0}^{p_1}\ldots\sum_{j_k=0}^{p_k}
C_{j_k\ldots j_1}
\prod_{l=1}^k \zeta_{j_l}^{(i_l)}+
J[R_{p_1\ldots p_k}]_{T,t}^{(k)}\ \ \ \hbox{w.~p.~1,}
\end{equation}

\noindent
where the stochastic integral
$J[R_{p_1\ldots p_k}]_{T,t}^{(k)}$
is defined by (\ref{30.34}) and
\begin{equation}
\label{30.46xx}
~~~~~ R_{p_1\ldots p_k}(t_1,\ldots,t_k)
=
K^{*}(t_1,\ldots,t_k)-
\sum_{j_1=0}^{p_1}\ldots\sum_{j_k=0}^{p_k}
C_{j_k\ldots j_1} \prod_{l=1}^{k} \phi_{j_l}(t_l),
\end{equation}
$$
\zeta_{j_l}^{(i_l)}=\int\limits_t^T \phi_{j_l}(s) d{\bf w}_s^{(i_l)}.
$$

\vspace{2mm}

According to Theorem 2.11, we have

\vspace{-4mm}
\begin{equation}
\label{410}
~~~\hbox{\vtop{\offinterlineskip\halign{
\hfil#\hfil\cr
{\rm lim}\cr
$\stackrel{}{{}_{p_1\to \infty}}$\cr
}} }\ldots \hbox{\vtop{\offinterlineskip\halign{
\hfil#\hfil\cr
{\rm lim}\cr
$\stackrel{}{{}_{p_k\to \infty}}$\cr
}} }R_{p_1\ldots p_k}(t_1,\ldots,t_k)=0\ \ \ \hbox{when}\ \ \
(t_1,\ldots,t_k)\in (t,T)^k,
\end{equation}

\vspace{1mm}
\noindent
where the left-hand side of (\ref{410}) is bounded on the boundary of 
$[t, T]^k.$

{\bf Theorem 2.13.} {\it Under the conditions of Theorem {\rm 2.10} we have

\vspace{-1mm}
$$
\lim\limits_{p_1\to \infty}
\varlimsup\limits_{p_2\to \infty}
\ldots
\varlimsup\limits_{p_k\to \infty}
{\sf M}\left\{\left|J[R_{p_1\ldots p_k}]_{T,t}^{(k)}\right|^{2n}
\right\}=0,\ \ \ n\in {\bf N}.
$$
}

\newpage
\noindent
\par
{\bf Proof.} At first let us analize in detail the cases $k=2, 3, 4.$ 
Using (\ref{best1}) (see below) and 
(\ref{30.52}), we have w.~p.~1 
$$
J[R_{p_1p_2}]_{T,t}^{(2)}
=\hbox{\vtop{\offinterlineskip\halign{
\hfil#\hfil\cr
{\rm l.i.m.}\cr
$\stackrel{}{{}_{N\to \infty}}$\cr
}} }\sum_{l_2=0}^{N-1}
\sum_{l_1=0}^{N-1}
R_{p_1 p_2}(\tau_{l_1},\tau_{l_2})
\Delta{\bf w}_{\tau_{l_1}}^{(i_1)}
\Delta{\bf w}_{\tau_{l_2}}^{(i_2)}=
$$
$$
=\hbox{\vtop{\offinterlineskip\halign{
\hfil#\hfil\cr
{\rm l.i.m.}\cr
$\stackrel{}{{}_{N\to \infty}}$\cr
}} }\left(\sum_{l_2=0}^{N-1}
\sum_{l_1=0}^{l_2-1}+\sum_{l_1=0}^{N-1}
\sum_{l_2=0}^{l_1-1}
\right)
R_{p_1 p_2}(\tau_{l_1},\tau_{l_2})
\Delta{\bf w}_{\tau_{l_1}}^{(i_1)}
\Delta{\bf w}_{\tau_{l_2}}^{(i_2)}+
$$
$$
+
\hbox{\vtop{\offinterlineskip\halign{
\hfil#\hfil\cr
{\rm l.i.m.}\cr
$\stackrel{}{{}_{N\to \infty}}$\cr
}} }\sum_{l_1=0}^{N-1}
R_{p_1 p_2}(\tau_{l_1},\tau_{l_1})
\Delta{\bf w}_{\tau_{l_1}}^{(i_1)}
\Delta{\bf w}_{\tau_{l_1}}^{(i_2)}=
$$
$$
=\int\limits_t^T\int\limits_t^{t_2}
R_{p_1p_2}(t_1,t_2)d{\bf w}_{t_1}^{(i_1)}d{\bf w}_{t_2}^{(i_2)}
+\int\limits_t^T\int\limits_t^{t_1}
R_{p_1p_2}(t_1,t_2)d{\bf w}_{t_2}^{(i_2)}d{\bf w}_{t_1}^{(i_1)}+
$$
\begin{equation}
\label{nov800}
+{\bf 1}_{\{i_1=i_2\ne 0\}}
\int\limits_t^T R_{p_1p_2}(t_1,t_1)dt_1,
\end{equation}

\noindent
where we used the same notations as in the
formulas
(\ref{30.34}), (\ref{30.52}) and
Lemma 1.1 (see
Sect.~1.1.3). Moreover,
\begin{equation}
\label{zajaz}
~~~R_{p_1 p_2}(t_1,t_2)=
K^{*}(t_1,t_2)-
\sum_{j_1=0}^{p_1}\sum_{j_2=0}^{p_2}C_{j_2 j_1}
\phi_{j_1}(t_1)\phi_{j_2}(t_2),\ \ \ p_1, p_2<\infty.
\end{equation}

Let us consider the following well known estimates
for moments of stochastic integrals \cite{Gih1}
\begin{equation}
\label{pupol1}
~~~~~~{\sf M}\left\{\left|\int\limits_{t}^T \xi_\tau
dw_\tau\right|^{2n}\right\} \le (T-t)^{n-1}\left(n(2n-1)\right)^n
\int\limits_{t}^T {\sf M}\left\{\left|\xi_\tau \right|^{2n}\right\}d\tau,
\end{equation}
\begin{equation}
\label{pupol2}
{\sf M}\left\{\left|\int\limits_{t}^T \xi_\tau
d\tau\right|^{2n}\right\} \le (T-t)^{2n-1}
\int\limits_{t}^T {\sf M}\left\{\left|\xi_\tau\right|^{2n}\right\}d\tau,
\end{equation}

\noindent
where the process $\xi_{\tau}$ is such that
$\left(\xi_{\tau}\right)^n\in{\rm M}_2
([t,T])$ and $w_{\tau}$ is a scalar standard Wiener 
process,\
$n=1, 2,\ldots$ (definition of the class 
${\rm M}_2([t,T])$ see in Sect.~1.1.2).

Using (\ref{pupol1}) and (\ref{pupol2}),  we obtain 
$$
{\sf M}\left\{\left|J[R_{p_1p_2}]_{T,t}^{(2)}\right|^{2n}
\right\}\le C_n\left(\int\limits_t^T\int\limits_t^{t_2}
\left(R_{p_1p_2}(t_1,t_2)\right)^{2n}dt_1 dt_2
+\right.
$$
\begin{equation}
\label{leto80010}
~~~~\left.+
\int\limits_t^T\int\limits_t^{t_1}
\left(R_{p_1p_2}(t_1,t_2)\right)^{2n}dt_2 dt_1+
{\bf 1}_{\{i_1=i_2\ne 0\}}
\int\limits_t^T \left(R_{p_1p_2}(t_1,t_1)\right)^{2n}dt_1\right),
\end{equation}

\vspace{2mm}
\noindent
where constant $C_n<\infty$ depends on 
$n$ and $T-t$ $(n=1, 2,\ldots).$

Further, we have
$$
\int\limits_t^T\int\limits_t^{t_2}
\left(R_{p_1p_2}(t_1,t_2)\right)^{2n}dt_1 dt_2
+
\int\limits_t^T\int\limits_t^{t_1}
\left(R_{p_1p_2}(t_1,t_2)\right)^{2n}dt_2 dt_1=
$$
$$
=
\int\limits_t^T\int\limits_t^{t_2}
\left(R_{p_1p_2}(t_1,t_2)\right)^{2n}dt_1 dt_2
+
\int\limits_t^T\int\limits_{t_2}^{T}
\left(R_{p_1p_2}(t_1,t_2)\right)^{2n}dt_1 dt_2=
$$

\vspace{1mm}
\begin{equation}
\label{leb1}
=
\int\limits_{[t, T]^2}
\left(R_{p_1p_2}(t_1,t_2)\right)^{2n}dt_1 dt_2.
\end{equation}

\vspace{2mm}

Combining (\ref{leto80010}) and (\ref{leb1}), we obtain

\vspace{-1mm}
$$
{\sf M}\left\{\left|J[R_{p_1p_2}]_{T,t}^{(2)}\right|^{2n}
\right\}\le 
$$
\begin{equation}
\label{leto80010aaa}
\le C_n\left(
\int\limits_{[t, T]^2}
\left(R_{p_1p_2}(t_1,t_2)\right)^{2n}dt_1 dt_2
+{\bf 1}_{\{i_1=i_2\ne 0\}}
\int\limits_t^T \left(R_{p_1p_2}(t_1,t_1)\right)^{2n}dt_1\right),
\end{equation}

\vspace{2mm}
\noindent
where constant $C_n<\infty$ depends on 
$n$ and $T-t$ $(n=1, 2,\ldots).$

Since the integrals on the right-hand side of (\ref{leto80010aaa}) 
exist as Riemann integrals, then they are equal to the 
corresponding Lebesgue integrals. 
Moreover, 
\begin{equation}
\label{strange101}
\lim\limits_{p_1\to\infty}\lim\limits_{p_2\to\infty}
\left(R_{p_1p_2}(t_1,t_2)\right)^{2n}=0\ \ \ \hbox{when}\ \ \ 
(t_1,t_2)\in (t, T)^2,
\end{equation}
where $n\in \bf{N}$ and the left-hand side is bounded on the boundary of 
$[t, T]^2.$

According to (\ref{zajaz}), we have
$$
R_{p_1p_2}(t_1,t_2)=\left(K^{*}(t_1,t_2)-\sum\limits_{j_1=0}^{p_1}
C_{j_1}(t_2)\phi_{j_1}(t_1)\right)+
$$
\begin{equation}
\label{d2020}
+\left(
\sum\limits_{j_1=0}^{p_1}\left(C_{j_1}(t_2)-
\sum\limits_{j_2=0}^{p_2}
C_{j_2j_1}\phi_{j_2}(t_2)\right)
\phi_{j_1}(t_1)\right).
\end{equation}

\vspace{2mm}

Then, applying two times (we mean here an iterated passage to the limit
$\lim\limits_{p_1\to\infty}\varlimsup\limits_{p_2\to\infty}$)
the Lebesgue's 
Dominated Convergence Theorem and taking into account
(\ref{leto8001}), (\ref{leto8002}), and (\ref{d2020}),
we obtain
\begin{equation}
\label{leb2}
\lim\limits_{p_1\to\infty}\varlimsup\limits_{p_2\to\infty}
\int\limits_{[t, T]^2}
\left(R_{p_1p_2}(t_1,t_2)\right)^{2n}dt_1 dt_2=0,
\end{equation}
\begin{equation}
\label{leb2xx}
\lim\limits_{p_1\to\infty}\varlimsup\limits_{p_2\to\infty}
\int\limits_t^T
\left(R_{p_1p_2}(t_1,t_1)\right)^{2n}dt_1=0.
\end{equation}

\vspace{2mm}

We will discuss the choice of integrable majorants
when applying Lebesgue's 
Dominated Convergence Theorem when we consider the case 
of arbitrary $k\in{\bf N}$ later in this section.

From 
(\ref{leto80010aaa}), (\ref{leb2}), and (\ref{leb2xx}) we get
$$
\lim\limits_{p_1\to \infty}
\varlimsup\limits_{p_2\to \infty}
{\sf M}\left\{\left|J[R_{p_1p_2}]_{T,t}^{(2)}\right|^{2n}
\right\}=0,\ \ \ n\in {\bf N}.
$$

Recall that (\ref{leb2xx}) for $2n=1$
has also been proved in Sect.~2.1.1, 2.1.2.

Let us consider the case $k=3.$ Using (\ref{oop12}) (see below)
and (\ref{30.52}), we have w.~p.~1
$$
J[R_{p_1p_2p_3}]_{T,t}^{(3)}=
\hbox{\vtop{\offinterlineskip\halign{
\hfil#\hfil\cr
{\rm l.i.m.}\cr
$\stackrel{}{{}_{N\to \infty}}$\cr
}} }\sum_{l_3=0}^{N-1}\sum_{l_2=0}^{N-1}
\sum_{l_1=0}^{N-1}
R_{p_1 p_2 p_3}(\tau_{l_1},\tau_{l_2},\tau_{l_3})
\Delta{\bf w}_{\tau_{l_1}}^{(i_1)}
\Delta{\bf w}_{\tau_{l_2}}^{(i_2)}
\Delta{\bf w}_{\tau_{l_3}}^{(i_3)}=
$$
$$
=\hbox{\vtop{\offinterlineskip\halign{
\hfil#\hfil\cr
{\rm l.i.m.}\cr
$\stackrel{}{{}_{N\to \infty}}$\cr
}} }\sum_{l_3=0}^{N-1}\sum_{l_2=0}^{l_3-1}
\sum_{l_1=0}^{l_2-1}\Biggl(
R_{p_1 p_2 p_3}(\tau_{l_1},\tau_{l_2},\tau_{l_3})
\Delta{\bf w}_{\tau_{l_1}}^{(i_1)}
\Delta{\bf w}_{\tau_{l_2}}^{(i_2)}
\Delta{\bf w}_{\tau_{l_3}}^{(i_3)}+\Biggr.
$$

\vspace{-2mm}
$$
+
R_{p_1 p_2 p_3}(\tau_{l_1},\tau_{l_3},\tau_{l_2})
\Delta{\bf w}_{\tau_{l_1}}^{(i_1)}
\Delta{\bf w}_{\tau_{l_3}}^{(i_2)}
\Delta{\bf w}_{\tau_{l_2}}^{(i_3)}+
R_{p_1 p_2 p_3}(\tau_{l_2},\tau_{l_1},\tau_{l_3})
\Delta{\bf w}_{\tau_{l_2}}^{(i_1)}
\Delta{\bf w}_{\tau_{l_1}}^{(i_2)}
\Delta{\bf w}_{\tau_{l_3}}^{(i_3)}+
$$
$$
+
R_{p_1 p_2 p_3}(\tau_{l_2},\tau_{l_3},\tau_{l_1})
\Delta{\bf w}_{\tau_{l_2}}^{(i_1)}
\Delta{\bf w}_{\tau_{l_3}}^{(i_2)}
\Delta{\bf w}_{\tau_{l_1}}^{(i_3)}+
R_{p_1 p_2 p_3}(\tau_{l_3},\tau_{l_2},\tau_{l_1})
\Delta{\bf w}_{\tau_{l_3}}^{(i_1)}
\Delta{\bf w}_{\tau_{l_2}}^{(i_2)}
\Delta{\bf w}_{\tau_{l_1}}^{(i_3)}+
$$
$$
\Biggl.
+R_{p_1 p_2 p_3}
(\tau_{l_3},\tau_{l_1},\tau_{l_2})\Delta{\bf w}_{\tau_{l_3}}^{(i_1)}
\Delta{\bf w}_{\tau_{l_1}}^{(i_2)}
\Delta{\bf w}_{\tau_{l_2}}^{(i_3)}\Biggr)+
$$
$$
+
\hbox{\vtop{\offinterlineskip\halign{
\hfil#\hfil\cr
{\rm l.i.m.}\cr
$\stackrel{}{{}_{N\to \infty}}$\cr
}} }\sum_{l_3=0}^{N-1}\sum_{l_2=0}^{l_3-1}\Biggl(
R_{p_1 p_2 p_3}(\tau_{l_2},\tau_{l_2},\tau_{l_3})
\Delta{\bf w}_{\tau_{l_2}}^{(i_1)}
\Delta{\bf w}_{\tau_{l_2}}^{(i_2)}
\Delta{\bf w}_{\tau_{l_3}}^{(i_3)}+\Biggr.
$$
$$
+R_{p_1 p_2 p_3}(\tau_{l_2},\tau_{l_3},\tau_{l_2})
\Delta{\bf w}_{\tau_{l_2}}^{(i_1)}
\Delta{\bf w}_{\tau_{l_3}}^{(i_2)}
\Delta{\bf w}_{\tau_{l_2}}^{(i_3)}+
$$
$$
\Biggl.+R_{p_1 p_2 p_3}(\tau_{l_3},\tau_{l_2},\tau_{l_2})
\Delta{\bf w}_{\tau_{l_3}}^{(i_1)}
\Delta{\bf w}_{\tau_{l_2}}^{(i_2)}
\Delta{\bf w}_{\tau_{l_2}}^{(i_3)}\Biggr)+
$$
$$
+
\hbox{\vtop{\offinterlineskip\halign{
\hfil#\hfil\cr
{\rm l.i.m.}\cr
$\stackrel{}{{}_{N\to \infty}}$\cr
}} }\sum_{l_3=0}^{N-1}\sum_{l_1=0}^{l_3-1}\Biggl(
R_{p_1 p_2 p_3}(\tau_{l_1},\tau_{l_3},\tau_{l_3})
\Delta{\bf w}_{\tau_{l_1}}^{(i_1)}
\Delta{\bf w}_{\tau_{l_3}}^{(i_2)}
\Delta{\bf w}_{\tau_{l_3}}^{(i_3)}+\Biggr.
$$
$$
+R_{p_1 p_2 p_3}(\tau_{l_3},\tau_{l_1},\tau_{l_3})
\Delta{\bf w}_{\tau_{l_3}}^{(i_1)}
\Delta{\bf w}_{\tau_{l_1}}^{(i_2)}
\Delta{\bf w}_{\tau_{l_3}}^{(i_3)}+
$$
$$
\Biggl.+R_{p_1 p_2 p_3}(\tau_{l_3},\tau_{l_3},\tau_{l_1})
\Delta{\bf w}_{\tau_{l_3}}^{(i_1)}
\Delta{\bf w}_{\tau_{l_3}}^{(i_2)}
\Delta{\bf w}_{\tau_{l_1}}^{(i_3)}\Biggr)+
$$
$$
+\hbox{\vtop{\offinterlineskip\halign{
\hfil#\hfil\cr
{\rm l.i.m.}\cr
$\stackrel{}{{}_{N\to \infty}}$\cr
}} }\sum_{l_3=0}^{N-1}
R_{p_1 p_2 p_3}(\tau_{l_3},\tau_{l_3},\tau_{l_3})
\Delta{\bf w}_{\tau_{l_3}}^{(i_1)}
\Delta{\bf w}_{\tau_{l_3}}^{(i_2)}
\Delta{\bf w}_{\tau_{l_3}}^{(i_3)}=
$$
$$
=
\int\limits_t^T\int\limits_t^{t_3}\int\limits_t^{t_2}
R_{p_1 p_2 p_3}(t_1,t_2,t_3)
d{\bf w}_{t_1}^{(i_1)}
d{\bf w}_{t_2}^{(i_2)}
d{\bf w}_{t_3}^{(i_3)}+
$$
$$
+
\int\limits_t^T\int\limits_t^{t_3}\int\limits_t^{t_2}
R_{p_1 p_2 p_3}(t_1,t_3,t_2)
d{\bf w}_{t_1}^{(i_1)}
d{\bf w}_{t_2}^{(i_3)}
d{\bf w}_{t_3}^{(i_2)}+
$$
$$
+
\int\limits_t^T\int\limits_t^{t_3}\int\limits_t^{t_2}
R_{p_1 p_2 p_3}(t_2,t_1,t_3)
d{\bf w}_{t_1}^{(i_2)}
d{\bf w}_{t_2}^{(i_1)}
d{\bf w}_{t_3}^{(i_3)}+
$$
$$
+
\int\limits_t^T\int\limits_t^{t_3}\int\limits_t^{t_2}
R_{p_1 p_2 p_3}(t_2,t_3,t_1)
d{\bf w}_{t_1}^{(i_3)}
d{\bf w}_{t_2}^{(i_1)}
d{\bf w}_{t_3}^{(i_2)}+
$$
$$
+
\int\limits_t^T\int\limits_t^{t_3}\int\limits_t^{t_2}
R_{p_1 p_2 p_3}(t_3,t_2,t_1)
d{\bf w}_{t_1}^{(i_3)}
d{\bf w}_{t_2}^{(i_2)}
d{\bf w}_{t_3}^{(i_1)}+
$$
$$
+
\int\limits_t^T\int\limits_t^{t_3}\int\limits_t^{t_2}
R_{p_1 p_2 p_3}(t_3,t_1,t_2)
d{\bf w}_{t_1}^{(i_2)}
d{\bf w}_{t_2}^{(i_3)}
d{\bf w}_{t_3}^{(i_1)}+
$$
$$
+{\bf 1}_{\{i_1=i_2\ne 0\}}
\int\limits_t^T\int\limits_t^{t_3}
R_{p_1 p_2 p_3}(t_2,t_2,t_3)
dt_2
d{\bf w}_{t_3}^{(i_3)}+
$$
$$
+{\bf 1}_{\{i_1=i_3\ne 0\}}
\int\limits_t^T\int\limits_t^{t_3}
R_{p_1 p_2 p_3}(t_2,t_3,t_2)
dt_2
d{\bf w}_{t_3}^{(i_2)}+
$$
$$
+{\bf 1}_{\{i_2=i_3\ne 0\}}
\int\limits_t^T\int\limits_t^{t_3}
R_{p_1 p_2 p_3}(t_3,t_2,t_2)
dt_2
d{\bf w}_{t_3}^{(i_1)}+
$$
$$
+{\bf 1}_{\{i_2=i_3\ne 0\}}
\int\limits_t^T\int\limits_t^{t_3}
R_{p_1 p_2 p_3}(t_1,t_3,t_3)
d{\bf w}_{t_1}^{(i_1)}dt_3+
$$
$$
+{\bf 1}_{\{i_1=i_3\ne 0\}}
\int\limits_t^T\int\limits_t^{t_3}
R_{p_1 p_2 p_3}(t_3,t_1,t_3)
d{\bf w}_{t_1}^{(i_2)}dt_3
+
$$
\begin{equation}
\label{s1s}
+{\bf 1}_{\{i_1=i_2\ne 0\}}
\int\limits_t^T\int\limits_t^{t_3}
R_{p_1 p_2 p_3}(t_3,t_3,t_1)
d{\bf w}_{t_1}^{(i_3)}dt_3,
\end{equation}
where we used the same notations as in 
the formulas
(\ref{30.34}), (\ref{30.52}) and
Lemma 1.1 (see
Sect.~1.1.3).
Using (\ref{pupol1}) and (\ref{pupol2}), we obtain from (\ref{s1s}) 
$$
{\sf M}\left\{\left|J[R_{p_1p_2p_3}]_{T,t}^{(3)}\right|^{2n}
\right\}\le
$$
$$
\le
C_n\Biggl(\int\limits_t^T\int\limits_t^{t_3}\int\limits_t^{t_2}
\Biggl(
\left(R_{p_1 p_2 p_3}(t_1,t_2,t_3)\right)^{2n}+
\left(R_{p_1 p_2 p_3}(t_1,t_3,t_2)\right)^{2n}+\Biggr.\Biggr.
$$
$$
+\left(R_{p_1 p_2 p_3}(t_2,t_1,t_3)\right)^{2n}+
\left(R_{p_1 p_2 p_3}(t_2,t_3,t_1)\right)^{2n}+
\left(R_{p_1 p_2 p_3}(t_3,t_2,t_1)\right)^{2n}+
$$

\vspace{-4mm}
$$
\Biggl.
+\left(R_{p_1 p_2 p_3}(t_3,t_1,t_2)\right)^{2n}\Biggr)dt_1dt_2dt_3+
$$
$$
+
\int\limits_t^T\int\limits_t^{t_3}\Biggl(
{\bf 1}_{\{i_1=i_2\ne 0\}}\Biggl(
\left(R_{p_1 p_2 p_3}(t_2,t_2,t_3)\right)^{2n}+
\left(R_{p_1 p_2 p_3}(t_3,t_3,t_2)\right)^{2n}\Biggr)+\Biggr.
$$
$$
+{\bf 1}_{\{i_1=i_3\ne 0\}}\Biggl(
\left(R_{p_1 p_2 p_3}(t_2,t_3,t_2)\right)^{2n}+
\left(R_{p_1 p_2 p_3}(t_3,t_2,t_3)\right)^{2n}\Biggr)+
$$
\begin{equation}
\label{oop16}
\Biggl.
+{\bf 1}_{\{i_2=i_3\ne 0\}}\Biggl(
\left(R_{p_1 p_2 p_3}(t_3,t_2,t_2)\right)^{2n}+
\left(R_{p_1 p_2 p_3}(t_2,t_3,t_3)\right)^{2n}\Biggr)dt_2dt_3\Biggr),\ \ \
C_n<\infty.
\end{equation}

Due to (\ref{30.46xx}) and Theorem~2.11
the
function $R_{p_1 p_2 p_3}(t_1,t_2,t_3)$ 
is continuous in the 
open domains 
of integration of iterated integrals on the right-hand side of
(\ref{oop16})
and it is bounded
at the 
boundaries 
of these 
domains. Moreover, everywhere in 
$(t, T)^3$ the following formula takes place
\begin{equation}
\label{oop15}
\lim\limits_{p_1\to\infty}\lim\limits_{p_2\to\infty}\lim\limits_{p_3\to\infty}
R_{p_1 p_2 p_3}(t_1,t_2,t_3)=0.
\end{equation}

Further, we have
$$
\int\limits_t^T\int\limits_t^{t_3}\int\limits_t^{t_2}
\Biggl(
\left(R_{p_1 p_2 p_3}(t_1,t_2,t_3)\right)^{2n}+
\left(R_{p_1 p_2 p_3}(t_1,t_3,t_2)\right)^{2n}+
\left(R_{p_1 p_2 p_3}(t_2,t_1,t_3)\right)^{2n}+\Biggr.
$$
$$
+
\left(R_{p_1 p_2 p_3}(t_2,t_3,t_1)\right)^{2n}+
\left(R_{p_1 p_2 p_3}(t_3,t_2,t_1)\right)^{2n}
+\left(R_{p_1 p_2 p_3}(t_3,t_1,t_2)\right)^{2n}\Biggr)dt_1dt_2dt_3=
$$

\vspace{1mm}
\begin{equation}
\label{zero1}
=
\int\limits_{[t, T]^3}
\left(R_{p_1 p_2 p_3}(t_1,t_2,t_3)\right)^{2n}dt_1dt_2dt_3,
\end{equation}

\vspace{3mm}

$$
\int\limits_t^T\int\limits_t^{t_3}\biggl(
\left(R_{p_1 p_2 p_3}(t_2,t_2,t_3)\right)^{2n}+
\left(R_{p_1 p_2 p_3}(t_3,t_3,t_2)\right)^{2n}\biggr)dt_2dt_3=
$$
$$
=\int\limits_t^T\int\limits_t^{t_3}
\left(R_{p_1 p_2 p_3}(t_2,t_2,t_3)\right)^{2n}dt_2dt_3+
\int\limits_t^T\int\limits_{t_3}^{T}
\left(R_{p_1 p_2 p_3}(t_2,t_2,t_3)\right)^{2n}dt_2dt_3=
$$

\vspace{1mm}
\begin{equation}
\label{zero2}
=
\int\limits_{[t, T]^2}
\left(R_{p_1 p_2 p_3}(t_2,t_2,t_3)\right)^{2n}dt_2dt_3,
\end{equation}
$$
\int\limits_t^T\int\limits_t^{t_3}\biggl(
\left(R_{p_1 p_2 p_3}(t_2,t_3,t_2)\right)^{2n}+
\left(R_{p_1 p_2 p_3}(t_3,t_2,t_3)\right)^{2n}\biggr)dt_2dt_3=
$$
$$
=\int\limits_t^T\int\limits_t^{t_3}
\left(R_{p_1 p_2 p_3}(t_2,t_3,t_2)\right)^{2n}dt_2dt_3+
\int\limits_t^T\int\limits_{t_3}^{T}
\left(R_{p_1 p_2 p_3}(t_2,t_3,t_2)\right)^{2n}dt_2dt_3=
$$

\vspace{1mm}
\begin{equation}
\label{zero3}
=
\int\limits_{[t, T]^2}
\left(R_{p_1 p_2 p_3}(t_2,t_3,t_2)\right)^{2n}dt_2dt_3,
\end{equation}

\vspace{3mm}

$$
\int\limits_t^T\int\limits_t^{t_3}\biggl(
\left(R_{p_1 p_2 p_3}(t_3,t_2,t_2)\right)^{2n}+
\left(R_{p_1 p_2 p_3}(t_2,t_3,t_3)\right)^{2n}\biggr)dt_2dt_3=
$$
$$
=\int\limits_t^T\int\limits_t^{t_3}
\left(R_{p_1 p_2 p_3}(t_3,t_2,t_2)\right)^{2n}dt_2dt_3+
\int\limits_t^T\int\limits_{t_3}^{T}
\left(R_{p_1 p_2 p_3}(t_3,t_2,t_2)\right)^{2n}dt_2dt_3=
$$

\vspace{1mm}
\begin{equation}
\label{zero4}
=
\int\limits_{[t, T]^2}
\left(R_{p_1 p_2 p_3}(t_3,t_2,t_2)\right)^{2n}dt_2dt_3.
\end{equation}

\vspace{2mm}

Combining (\ref{oop16}) and (\ref{zero1})--(\ref{zero4}),
we obtain
$$
{\sf M}\left\{\left|J[R_{p_1p_2p_3}]_{T,t}^{(3)}\right|^{2n}
\right\}\le 
C_n\left(
\int\limits_{[t, T]^3}
\left(R_{p_1 p_2 p_3}(t_1,t_2,t_3)\right)^{2n}dt_1dt_2dt_3+\right.
$$
$$
+
{\bf 1}_{\{i_1=i_2\ne 0\}}
\int\limits_{[t, T]^2}
\left(R_{p_1 p_2 p_3}(t_2,t_2,t_3)\right)^{2n}dt_2dt_3+
$$
$$
+{\bf 1}_{\{i_1=i_3\ne 0\}}
\int\limits_{[t, T]^2}
\left(R_{p_1 p_2 p_3}(t_2,t_3,t_2)\right)^{2n}dt_2dt_3+
$$
\begin{equation}
\label{oop16xxx}
\left.
+{\bf 1}_{\{i_2=i_3\ne 0\}}
\int\limits_{[t, T]^2}
\left(R_{p_1 p_2 p_3}(t_3,t_2,t_2)\right)^{2n}dt_2dt_3\right).
\end{equation}

Since the integrals on the right-hand side of (\ref{oop16xxx}) 
exist as Riemann integrals, then they are equal to the 
corresponding Lebesgue integrals. 
Moreover, 
$$
\lim\limits_{p_1\to\infty}\lim\limits_{p_2\to\infty}\lim\limits_{p_3\to\infty}
R_{p_1 p_2 p_3}(t_1,t_2,t_3)=0\ \ \ \hbox{when}\ \ \ 
(t_1,t_2,t_3)\in (t, T)^3,
$$
where the left-hand side is bounded on the boundary of 
$[t, T]^3.$

According to the proof of Theorem 2.11 and (\ref{30.46xx}) for $k=3$, we have
$$
R_{p_1p_2p_3}(t_1,t_2,t_3)=\left(K^{*}(t_1,t_2,t_3)-\sum\limits_{j_1=0}^{p_1}
C_{j_1}(t_2,t_3)\phi_{j_1}(t_1)\right)+
$$
$$
+\left(
\sum\limits_{j_1=0}^{p_1}\left(C_{j_1}(t_2,t_3)-
\sum\limits_{j_2=0}^{p_2}
C_{j_2j_1}(t_3)\phi_{j_2}(t_2)\right)
\phi_{j_1}(t_1)\right)+
$$

\vspace{-3mm}
\begin{equation}
\label{novvv1}
~~~~~~~~~+\left(
\sum\limits_{j_1=0}^{p_1}\sum\limits_{j_2=0}^{p_2}\left(C_{j_2j_1}(t_3)-
\sum\limits_{j_3=0}^{p_3}
C_{j_3j_2j_1}\phi_{j_3}(t_3)\right)
\phi_{j_2}(t_2)\phi_{j_1}(t_1)\right),
\end{equation}

\vspace{2mm}
\noindent
where
$$
C_{j_1}(t_2,t_3)=\int\limits_t^T
K^{*}(t_1,t_2,t_3)\phi_{j_1}(t_1)dt_1,
$$
$$
C_{j_2j_1}(t_3)=\int\limits_{[t, T]^2}
K^{*}(t_1,t_2,t_3)\phi_{j_1}(t_1)\phi_{j_2}(t_2)dt_1 dt_2.
$$

\vspace{2mm}

Then, applying three times (we mean here an iterated passage to the limit
$\lim\limits_{p_1\to\infty}\varlimsup\limits_{p_2\to\infty}
\varlimsup\limits_{p_3\to\infty}$)
the Lebesgue's Dominated Convergence Theorem, 
we obtain
\begin{equation}
\label{leb21}
~~~~~~~~~~\lim\limits_{p_1\to\infty}\varlimsup\limits_{p_2\to\infty}
\varlimsup\limits_{p_3\to\infty}
\int\limits_{[t, T]^3}
\left(R_{p_1p_2p_3}(t_1,t_2,t_3)\right)^{2n}dt_1 dt_2 dt_3=0,
\end{equation}
\begin{equation}
\label{leb22}
\lim\limits_{p_1\to\infty}\varlimsup\limits_{p_2\to\infty}
\varlimsup\limits_{p_3\to\infty}
\int\limits_{[t, T]^2}
\left(R_{p_1 p_2 p_3}(t_2,t_2,t_3)\right)^{2n}dt_2dt_3=0,
\end{equation}
\begin{equation}
\label{leb23}
\lim\limits_{p_1\to\infty}\varlimsup\limits_{p_2\to\infty}
\varlimsup\limits_{p_3\to\infty}
\int\limits_{[t, T]^2}
\left(R_{p_1 p_2 p_3}(t_2,t_3,t_2)\right)^{2n}dt_2dt_3=0,
\end{equation}
\begin{equation}
\label{leb24}
\lim\limits_{p_1\to\infty}\varlimsup\limits_{p_2\to\infty}
\varlimsup\limits_{p_3\to\infty}
\int\limits_{[t, T]^2}
\left(R_{p_1 p_2 p_3}(t_3,t_2,t_2)\right)^{2n}dt_2dt_3=0.
\end{equation}

\vspace{2mm}

From 
(\ref{oop16xxx}) and (\ref{leb21})--(\ref{leb24}) we get
$$
\lim\limits_{p_1\to\infty}
\varlimsup\limits_{p_2\to\infty}\varlimsup\limits_{p_3\to\infty}
{\sf M}\left\{\left|J[R_{p_1p_2p_3}]_{T,t}^{(3)}\right|^{2n}
\right\}=0,\ \ \ n\in {\bf N}.
$$

\vspace{2mm}

Let us consider the case $k=4.$ Using (\ref{huh}) (see below) 
and (\ref{30.52}), we have
w.~p.~1 

\vspace{-3mm}
$$
J[R_{p_1p_2p_3p_4}]_{T,t}^{(4)}=
$$

\vspace{-5mm}
$$
=
\hbox{\vtop{\offinterlineskip\halign{
\hfil#\hfil\cr
{\rm l.i.m.}\cr
$\stackrel{}{{}_{N\to \infty}}$\cr
}} }\sum_{l_4=0}^{N-1}\sum_{l_3=0}^{N-1}\sum_{l_2=0}^{N-1}
\sum_{l_1=0}^{N-1}
R_{p_1 p_2 p_3 p_4}(\tau_{l_1},\tau_{l_2},\tau_{l_3}, \tau_{l_4})
\Delta{\bf w}_{\tau_{l_1}}^{(i_1)}
\Delta{\bf w}_{\tau_{l_2}}^{(i_2)}
\Delta{\bf w}_{\tau_{l_3}}^{(i_3)}\Delta{\bf w}_{\tau_{l_4}}^{(i_4)}=
$$

\vspace{-1mm}
$$
=\hbox{\vtop{\offinterlineskip\halign{
\hfil#\hfil\cr
{\rm l.i.m.}\cr
$\stackrel{}{{}_{N\to \infty}}$\cr
}} }\sum_{l_4=0}^{N-1}\sum_{l_3=0}^{l_4-1}\sum_{l_2=0}^{l_3-1}
\sum_{l_1=0}^{l_2-1}
\sum\limits_{(l_1,l_2,l_3,l_4)}
\biggl(R_{p_1 p_2 p_3 p_4}(\tau_{l_1},\tau_{l_2},\tau_{l_3}, \tau_{l_4})
\times\biggr.
$$

$$
\biggl.\times
\Delta{\bf w}_{\tau_{l_1}}^{(i_1)}
\Delta{\bf w}_{\tau_{l_2}}^{(i_2)}
\Delta{\bf w}_{\tau_{l_3}}^{(i_3)}\Delta{\bf w}_{\tau_{l_4}}^{(i_4)}\biggr)+
$$

\vspace{-5mm}
$$
+\hbox{\vtop{\offinterlineskip\halign{
\hfil#\hfil\cr
{\rm l.i.m.}\cr
$\stackrel{}{{}_{N\to \infty}}$\cr
}} }\hspace{-1.5mm}\sum_{l_4=0}^{N-1}\sum_{l_3=0}^{l_4-1}\sum_{l_2=0}^{l_3-1}
\sum\limits_{(l_2,l_2,l_3,l_4)}
\hspace{-1.5mm}
\biggl(\hspace{-0.4mm}
R_{p_1 p_2 p_3 p_4}(\tau_{l_2},\tau_{l_2},\tau_{l_3}, \tau_{l_4})
\Delta{\bf w}_{\tau_{l_2}}^{(i_1)}
\Delta{\bf w}_{\tau_{l_2}}^{(i_2)}
\Delta{\bf w}_{\tau_{l_3}}^{(i_3)}\Delta{\bf w}_{\tau_{l_4}}^{(i_4)}
\biggr)\hspace{-0.4mm}+
$$

\vspace{-3mm}
$$
+\hbox{\vtop{\offinterlineskip\halign{
\hfil#\hfil\cr
{\rm l.i.m.}\cr
$\stackrel{}{{}_{N\to \infty}}$\cr
}} }\hspace{-1.5mm}\sum_{l_4=0}^{N-1}\sum_{l_3=0}^{l_4-1}\sum_{l_1=0}^{l_3-1}
\sum\limits_{(l_1,l_3,l_3,l_4)}
\hspace{-1.5mm}
\biggl(\hspace{-0.4mm}
R_{p_1 p_2 p_3 p_4}(\tau_{l_1},\tau_{l_3},\tau_{l_3}, \tau_{l_4})
\Delta{\bf w}_{\tau_{l_1}}^{(i_1)}
\Delta{\bf w}_{\tau_{l_3}}^{(i_2)}
\Delta{\bf w}_{\tau_{l_3}}^{(i_3)}\Delta{\bf w}_{\tau_{l_4}}^{(i_4)}\biggr)
\hspace{-0.4mm}+
$$

\vspace{-3mm}
$$
+\hbox{\vtop{\offinterlineskip\halign{
\hfil#\hfil\cr
{\rm l.i.m.}\cr
$\stackrel{}{{}_{N\to \infty}}$\cr
}} }\hspace{-1.5mm}\sum_{l_4=0}^{N-1}\sum_{l_2=0}^{l_4-1}\sum_{l_1=0}^{l_2-1}
\sum\limits_{(l_1,l_2,l_4,l_4)}
\hspace{-1.5mm}
\biggl(\hspace{-0.4mm}
R_{p_1 p_2 p_3 p_4}(\tau_{l_1},\tau_{l_2},\tau_{l_4}, \tau_{l_4})
\Delta{\bf w}_{\tau_{l_1}}^{(i_1)}
\Delta{\bf w}_{\tau_{l_2}}^{(i_2)}
\Delta{\bf w}_{\tau_{l_4}}^{(i_3)}\Delta{\bf w}_{\tau_{l_4}}^{(i_4)}\biggr)
\hspace{-0.4mm}+
$$

\vspace{-3mm}
$$
+\hbox{\vtop{\offinterlineskip\halign{
\hfil#\hfil\cr
{\rm l.i.m.}\cr
$\stackrel{}{{}_{N\to \infty}}$\cr
}} }\sum_{l_4=0}^{N-1}\sum_{l_3=0}^{l_4-1}
\sum\limits_{(l_3,l_3,l_3,l_4)}
\biggl(R_{p_1 p_2 p_3 p_4}(\tau_{l_3},\tau_{l_3},\tau_{l_3}, \tau_{l_4})
\Delta{\bf w}_{\tau_{l_3}}^{(i_1)}
\Delta{\bf w}_{\tau_{l_3}}^{(i_2)}
\Delta{\bf w}_{\tau_{l_3}}^{(i_3)}\Delta{\bf w}_{\tau_{l_4}}^{(i_4)}\biggr)+
$$

\vspace{-3mm}
$$
+\hbox{\vtop{\offinterlineskip\halign{
\hfil#\hfil\cr
{\rm l.i.m.}\cr
$\stackrel{}{{}_{N\to \infty}}$\cr
}} }\sum_{l_4=0}^{N-1}\sum_{l_2=0}^{l_4-1}
\sum\limits_{(l_2,l_2,l_4,l_4)}
\biggl(R_{p_1 p_2 p_3 p_4}(\tau_{l_2},\tau_{l_2},\tau_{l_4}, \tau_{l_4})
\Delta{\bf w}_{\tau_{l_2}}^{(i_1)}
\Delta{\bf w}_{\tau_{l_2}}^{(i_2)}
\Delta{\bf w}_{\tau_{l_4}}^{(i_3)}\Delta{\bf w}_{\tau_{l_4}}^{(i_4)}\biggr)+
$$
$$
+\hbox{\vtop{\offinterlineskip\halign{
\hfil#\hfil\cr
{\rm l.i.m.}\cr
$\stackrel{}{{}_{N\to \infty}}$\cr
}} }\sum_{l_4=0}^{N-1}\sum_{l_1=0}^{l_4-1}
\sum\limits_{(l_1,l_4,l_4,l_4)}
\biggl(R_{p_1 p_2 p_3 p_4}(\tau_{l_1},\tau_{l_4},\tau_{l_4}, \tau_{l_4})
\Delta{\bf w}_{\tau_{l_1}}^{(i_1)}
\Delta{\bf w}_{\tau_{l_4}}^{(i_2)}
\Delta{\bf w}_{\tau_{l_4}}^{(i_3)}\Delta{\bf w}_{\tau_{l_4}}^{(i_4)}\biggr)+
$$

\vspace{-3mm}
$$
+\hbox{\vtop{\offinterlineskip\halign{
\hfil#\hfil\cr
{\rm l.i.m.}\cr
$\stackrel{}{{}_{N\to \infty}}$\cr
}} }\sum_{l_4=0}^{N-1}
R_{p_1 p_2 p_3 p_4}(\tau_{l_4},\tau_{l_4},\tau_{l_4}, \tau_{l_4})
\Delta{\bf w}_{\tau_{l_4}}^{(i_1)}
\Delta{\bf w}_{\tau_{l_4}}^{(i_2)}
\Delta{\bf w}_{\tau_{l_4}}^{(i_3)}\Delta{\bf w}_{\tau_{l_4}}^{(i_4)}=
$$

\vspace{-3mm}
$$
=
\int\limits_t^T\int\limits_t^{t_4}\int\limits_t^{t_3}\int\limits_t^{t_2}
\sum\limits_{(t_1,t_2,t_3,t_4)}\biggl(R_{p_1 p_2 p_3 p_4}(t_1,t_2,t_3,t_4)
d{\bf w}_{t_1}^{(i_1)}
d{\bf w}_{t_2}^{(i_2)}
d{\bf w}_{t_3}^{(i_3)}
d{\bf w}_{t_4}^{(i_4)}\biggr)+
$$

\vspace{-3mm}
$$
+
{\bf 1}_{\{i_1=i_2\ne 0\}}
\int\limits_t^T\int\limits_t^{t_4}\int\limits_t^{t_3}
\sum\limits_{(t_1,t_3,t_4)}\biggl(R_{p_1 p_2 p_3 p_4}(t_1,t_1,t_3,t_4)
dt_1
d{\bf w}_{t_3}^{(i_3)}
d{\bf w}_{t_4}^{(i_4)}\biggr)+
$$

\vspace{-3mm}
$$
+
{\bf 1}_{\{i_1=i_3\ne 0\}}
\int\limits_t^T\int\limits_t^{t_4}\int\limits_t^{t_2}
\sum\limits_{(t_1,t_2,t_4)}\biggl(
R_{p_1 p_2 p_3 p_4}(t_1,t_2,t_1,t_4)
dt_1
d{\bf w}_{t_2}^{(i_2)}
d{\bf w}_{t_4}^{(i_4)}\biggr)+
$$

\vspace{-3mm}
$$
+
{\bf 1}_{\{i_1=i_4\ne 0\}}
\int\limits_t^T\int\limits_t^{t_3}\int\limits_t^{t_2}
\sum\limits_{(t_1,t_2,t_3)}
\biggl(R_{p_1 p_2 p_3 p_4}(t_1,t_2,t_3,t_1)
dt_1 d{\bf w}_{t_2}^{(i_2)}
d{\bf w}_{t_3}^{(i_3)}\biggr)+
$$

\vspace{-3mm}
$$
+
{\bf 1}_{\{i_2=i_3\ne 0\}}
\int\limits_t^T\int\limits_t^{t_4}\int\limits_t^{t_2}
\sum\limits_{(t_1,t_2,t_4)}
\biggl(R_{p_1 p_2 p_3 p_4}(t_1,t_2,t_2,t_4)
d{\bf w}_{t_1}^{(i_1)}
dt_2
d{\bf w}_{t_4}^{(i_4)}\biggr)+
$$

\vspace{-3mm}
$$
+
{\bf 1}_{\{i_2=i_4\ne 0\}}
\int\limits_t^T\int\limits_t^{t_3}\int\limits_t^{t_2}
\sum\limits_{(t_1,t_2,t_3)}
\biggl(R_{p_1 p_2 p_3 p_4}(t_1,t_2,t_3,t_2)
d{\bf w}_{t_1}^{(i_1)}
dt_2
d{\bf w}_{t_3}^{(i_3)}\biggr)+
$$

\vspace{-3mm}
$$
+
{\bf 1}_{\{i_3=i_4\ne 0\}}
\int\limits_t^T\int\limits_t^{t_3}\int\limits_t^{t_2}
\sum\limits_{(t_1,t_2,t_3)}
\biggl(R_{p_1 p_2 p_3 p_4}(t_1,t_2,t_3,t_3)
d{\bf w}_{t_1}^{(i_1)}
d{\bf w}_{t_2}^{(i_2)}dt_3\biggr)+
$$

\vspace{-3mm}
$$
+
{\bf 1}_{\{i_1=i_2\ne 0\}}{\bf 1}_{\{i_3=i_4\ne 0\}}
\left(\int\limits_t^T\int\limits_t^{t_4}
R_{p_1 p_2 p_3 p_4}(t_2,t_2,t_4,t_4)dt_2 dt_4
+\right.
$$
$$
\left.+ \int\limits_t^T\int\limits_t^{t_4}
R_{p_1 p_2 p_3 p_4}(t_4,t_4,t_2,t_2)dt_2 dt_4\right)+
$$

\vspace{-3mm}
$$
+
{\bf 1}_{\{i_1=i_3\ne 0\}}{\bf 1}_{\{i_2=i_4\ne 0\}}
\left(\int\limits_t^T\int\limits_t^{t_4}
R_{p_1 p_2 p_3 p_4}(t_2,t_4,t_2,t_4)dt_2 dt_4
+\right.
$$

\vspace{-3mm}
$$
\left.+\int\limits_t^T\int\limits_t^{t_4}
R_{p_1 p_2 p_3 p_4}(t_4,t_2,t_4,t_2)dt_2 dt_4\right)+
$$

\vspace{-3mm}
$$
+
{\bf 1}_{\{i_1=i_4\ne 0\}}{\bf 1}_{\{i_2=i_3\ne 0\}}
\left(\int\limits_t^T\int\limits_t^{t_4}
R_{p_1 p_2 p_3 p_4}(t_2,t_4,t_4,t_2)dt_2 dt_4
+\right.
$$

\vspace{-3mm}
\begin{equation}
\label{sogl}
\left.+\int\limits_t^T\int\limits_t^{t_4}
R_{p_1 p_2 p_3 p_4}(t_4,t_2,t_2,t_4)dt_2 dt_4\right),
\end{equation}

\vspace{3mm}
\noindent
where
the expression
$$
\sum\limits_{(a_1, \ldots, a_k)}
$$
means the sum with respect to 
all possible  
permutations $(a_1, \ldots, a_k)$. 
Moreover, we used in (\ref{sogl}) 
the same notations as in the proof of Theorem 1.1
(see Sect.~1.1.3).
Note that an analogue of (\ref{sogl}) will be
obtained in Sect.~2.6 (also see \cite{9}-\cite{12aa}, \cite{arxiv-8})  
with using the another approach.

By analogy with
(\ref{oop16xxx}) we obtain

\vspace{-1mm}
$$
{\sf M}\left\{\left|J[R_{p_1p_2p_3p_4}]_{T,t}^{(4)}\right|^{2n}
\right\}\le
$$
$$
\le 
C_n\left(
\int\limits_{[t, T]^4}
\left(R_{p_1 p_2 p_3 p_4}(t_1,t_2,t_3,t_4)\right)^{2n}dt_1dt_2dt_3dt_4+\right.
$$
$$
+
{\bf 1}_{\{i_1=i_2\ne 0\}}
\int\limits_{[t, T]^3}
\left(R_{p_1 p_2 p_3 p_4}(t_2,t_2,t_3,t_4)\right)^{2n}dt_2dt_3dt_4+
$$
$$
+
{\bf 1}_{\{i_1=i_3\ne 0\}}
\int\limits_{[t, T]^3}
\left(R_{p_1 p_2 p_3 p_4}(t_2,t_3,t_2,t_4)\right)^{2n}dt_2dt_3dt_4+
$$
$$
+
{\bf 1}_{\{i_1=i_4\ne 0\}}
\int\limits_{[t, T]^3}
\left(R_{p_1 p_2 p_3 p_4}(t_2,t_3,t_4,t_2)\right)^{2n}dt_2dt_3dt_4+
$$
$$
+
{\bf 1}_{\{i_2=i_3\ne 0\}}
\int\limits_{[t, T]^3}
\left(R_{p_1 p_2 p_3 p_4}(t_3,t_2,t_2,t_4)\right)^{2n}dt_2dt_3dt_4+
$$
$$
+
{\bf 1}_{\{i_2=i_4\ne 0\}}
\int\limits_{[t, T]^3}
\left(R_{p_1 p_2 p_3 p_4}(t_3,t_2,t_4,t_2)\right)^{2n}dt_2dt_3dt_4+
$$
$$
+
{\bf 1}_{\{i_3=i_4\ne 0\}}
\int\limits_{[t, T]^3}
\left(R_{p_1 p_2 p_3 p_4}(t_3,t_4,t_2,t_2)\right)^{2n}dt_2dt_3dt_4+
$$
$$
+
{\bf 1}_{\{i_1=i_2\ne 0\}}{\bf 1}_{\{i_3=i_4\ne 0\}}
\int\limits_{[t, T]^2}
\left(R_{p_1 p_2 p_3 p_4}(t_2,t_2,t_4,t_4)\right)^{2n}dt_2dt_4+
$$
$$
+
{\bf 1}_{\{i_1=i_3\ne 0\}}{\bf 1}_{\{i_2=i_4\ne 0\}}
\int\limits_{[t, T]^2}
\left(R_{p_1 p_2 p_3 p_4}(t_2,t_4,t_2,t_4)\right)^{2n}dt_2dt_4+
$$
\begin{equation}
\label{udar}
\left.~~~+
{\bf 1}_{\{i_1=i_4\ne 0\}}{\bf 1}_{\{i_2=i_3\ne 0\}}
\int\limits_{[t, T]^2}
\left(R_{p_1 p_2 p_3 p_4}(t_2,t_4,t_4,t_2)\right)^{2n}dt_2dt_4\right),\ \
C_n<\infty.
\end{equation}

Since the integrals on the right-hand side of (\ref{udar}) 
exist as Riemann integrals, then they are equal to the 
corresponding Lebesgue integrals. 
Moreover, 

\vspace{-6mm}
$$
~~\lim\limits_{p_1\to\infty}\lim\limits_{p_2\to\infty}\lim\limits_{p_3\to\infty}
\lim\limits_{p_4\to\infty}
R_{p_1 p_2 p_3 p_4}(t_1,t_2,t_3,t_4)=0\ \ \ \hbox{when}\ \ \ 
(t_1,t_2,t_3,t_4)\in (t, T)^4,
$$

\vspace{1mm}
\noindent
where the left-hand side is bounded on the boundary of 
$[t, T]^4.$

According to the proof of Theorem 2.11 and (\ref{30.46xx}) for $k=4$, we have

\vspace{-2mm}
$$
R_{p_1p_2p_3p_4}
(t_1,t_2,t_3,t_4)=
$$

\vspace{-4mm}

$$
=\left(K^{*}(t_1,t_2,t_3,t_4)-
\sum\limits_{j_1=0}^{p_1}
C_{j_1}(t_2,t_3,t_4)\phi_{j_1}(t_1)\right)+
$$

\vspace{-3mm}
$$
+\left(
\sum\limits_{j_1=0}^{p_1}\left(C_{j_1}(t_2,t_3,t_4)-
\sum\limits_{j_2=0}^{p_2}
C_{j_2j_1}(t_3,t_4)\phi_{j_2}(t_2)\right)
\phi_{j_1}(t_1)\right)+
$$
$$
+\left(
\sum\limits_{j_1=0}^{p_1}\sum\limits_{j_2=0}^{p_2}\left(C_{j_2j_1}(t_3,t_4)-
\sum\limits_{j_3=0}^{p_3}
C_{j_3j_2j_1}(t_4)\phi_{j_3}(t_3)\right)
\phi_{j_2}(t_2)\phi_{j_1}(t_1)\right)+
$$

\vspace{-3mm}
$$
+\left(
\sum\limits_{j_1=0}^{p_1}\sum\limits_{j_2=0}^{p_2}\sum\limits_{j_3=0}^{p_3}
\left(C_{j_3j_2j_1}(t_4)-
\sum\limits_{j_4=0}^{p_4}
C_{j_4j_3j_2j_1}\phi_{j_4}(t_4)\right)
\phi_{j_3}(t_3)\phi_{j_2}(t_2)\phi_{j_1}(t_1)\right),
$$

\vspace{4mm}
\noindent
where
$$
C_{j_1}(t_2,t_3,t_4)=\int\limits_t^T
K^{*}(t_1,t_2,t_3,t_4)\phi_{j_1}(t_1)dt_1,
$$

\vspace{1mm}
$$
C_{j_2j_1}(t_3,t_4)=\int\limits_{[t, T]^2}
K^{*}(t_1,t_2,t_3,t_4)\phi_{j_1}(t_1)\phi_{j_2}(t_2)dt_1 dt_2,
$$

\vspace{1mm}
$$
C_{j_3j_2j_1}(t_4)=\int\limits_{[t, T]^3}
K^{*}(t_1,t_2,t_3,t_4)\phi_{j_1}(t_1)\phi_{j_2}(t_2)\phi_{j_3}(t_3)
dt_1 dt_2 dt_3.
$$

\vspace{2mm}

Then, applying  four times (we mean here an iterated passage to the limit
$\lim\limits_{p_1\to\infty}\varlimsup\limits_{p_2\to\infty}
\varlimsup\limits_{p_3\to\infty}\varlimsup\limits_{p_4\to\infty}$)
the Lebesgue's Dominated Convergence Theorem, 
we obtain
\begin{equation}
\label{final1}
~~~~\lim\limits_{p_1\to\infty}\varlimsup\limits_{p_2\to\infty}
\varlimsup\limits_{p_3\to\infty}\varlimsup\limits_{p_4\to\infty}\int\limits_{[t, T]^4}
\left(R_{p_1 p_2 p_3 p_4}(t_1,t_2,t_3,t_4)\right)^{2n}dt_1dt_2dt_3dt_4=0,
\end{equation}

\vspace{-3mm}
\begin{equation}
\label{final2}
~~~~\lim\limits_{p_1\to\infty}\varlimsup\limits_{p_2\to\infty}
\varlimsup\limits_{p_3\to\infty}\varlimsup\limits_{p_4\to\infty}\int\limits_{[t, T]^3}
\left(R_{p_1 p_2 p_3 p_4}(t_2,t_2,t_3,t_4)\right)^{2n}dt_2dt_3dt_4=0,
\end{equation}

\vspace{-3mm}
\begin{equation}
\label{final3}
~~~~\lim\limits_{p_1\to\infty}\varlimsup\limits_{p_2\to\infty}
\varlimsup\limits_{p_3\to\infty}\varlimsup\limits_{p_4\to\infty}\int\limits_{[t, T]^3}
\left(R_{p_1 p_2 p_3 p_4}(t_2,t_3,t_2,t_4)\right)^{2n}dt_2dt_3dt_4=0,
\end{equation}

\vspace{-3mm}
\begin{equation}
\label{final4}
~~~~\lim\limits_{p_1\to\infty}\varlimsup\limits_{p_2\to\infty}
\varlimsup\limits_{p_3\to\infty}\varlimsup\limits_{p_4\to\infty}\int\limits_{[t, T]^3}
\left(R_{p_1 p_2 p_3 p_4}(t_2,t_3,t_4,t_2)\right)^{2n}dt_2dt_3dt_4=0,
\end{equation}

\vspace{-3mm}
\begin{equation}
\label{final5}
~~~~\lim\limits_{p_1\to\infty}\varlimsup\limits_{p_2\to\infty}
\varlimsup\limits_{p_3\to\infty}\varlimsup\limits_{p_4\to\infty}\int\limits_{[t, T]^3}
\left(R_{p_1 p_2 p_3 p_4}(t_3,t_2,t_2,t_4)\right)^{2n}dt_2dt_3dt_4=0,
\end{equation}
\begin{equation}
\label{final6}
~~~~\lim\limits_{p_1\to\infty}\varlimsup\limits_{p_2\to\infty}
\varlimsup\limits_{p_3\to\infty}\varlimsup\limits_{p_4\to\infty}\int\limits_{[t, T]^3}
\left(R_{p_1 p_2 p_3 p_4}(t_3,t_2,t_4,t_2)\right)^{2n}dt_2dt_3dt_4=0,
\end{equation}

\vspace{-3.5mm}
\begin{equation}
\label{final7}
~~~~\lim\limits_{p_1\to\infty}\varlimsup\limits_{p_2\to\infty}
\varlimsup\limits_{p_3\to\infty}\varlimsup\limits_{p_4\to\infty}\int\limits_{[t, T]^3}
\left(R_{p_1 p_2 p_3 p_4}(t_3,t_4,t_2,t_2)\right)^{2n}dt_2dt_3dt_4=0,
\end{equation}

\vspace{-3.5mm}
\begin{equation}
\label{final8}
~~~~\lim\limits_{p_1\to\infty}\varlimsup\limits_{p_2\to\infty}
\varlimsup\limits_{p_3\to\infty}\varlimsup\limits_{p_4\to\infty}\int\limits_{[t, T]^2}
\left(R_{p_1 p_2 p_3 p_4}(t_2,t_2,t_4,t_4)\right)^{2n}dt_2dt_4=0,
\end{equation}

\vspace{-3.5mm}
\begin{equation}
\label{final9}
~~~~\lim\limits_{p_1\to\infty}\varlimsup\limits_{p_2\to\infty}
\varlimsup\limits_{p_3\to\infty}\varlimsup\limits_{p_4\to\infty}\int\limits_{[t, T]^2}
\left(R_{p_1 p_2 p_3 p_4}(t_2,t_4,t_2,t_4)\right)^{2n}dt_2dt_4=0,
\end{equation}

\vspace{-3.5mm}
\begin{equation}
\label{final10}
~~~~\lim\limits_{p_1\to\infty}\varlimsup\limits_{p_2\to\infty}
\varlimsup\limits_{p_3\to\infty}\varlimsup\limits_{p_4\to\infty}\int\limits_{[t, T]^2}
\left(R_{p_1 p_2 p_3 p_4}(t_2,t_4,t_4,t_2)\right)^{2n}dt_2dt_4=0.
\end{equation}

\vspace{2mm}

Combaining
(\ref{udar}) with (\ref{final1})--(\ref{final10}), we get 

\vspace{-2mm}
$$
\lim\limits_{p_1\to\infty}
\varlimsup\limits_{p_2\to\infty}\varlimsup\limits_{p_3\to\infty}
\varlimsup\limits_{p_4\to\infty}
{\sf M}\left\{\left|J[R_{p_1p_2p_3p_4}]_{T,t}^{(4)}\right|^{2n}
\right\}=0,\ \ \ n\in {\bf N}.
$$

\vspace{3mm}

Theorem 2.13 is proved for $k=4$.

Let us consider the case of arbitrary $k,$ $k\in{\bf N}$.
Let us analyze the stochastic integral defined by
(\ref{30.34}) and find its representation convenient  
for the following consideration. In order to do it we 
introduce several notations. Suppose that 

\vspace{-2mm}
$$
S_N^{(k)}(a)=
\sum\limits_{j_k=0}^{N-1}\ldots \sum\limits_{j_1=0}^{j_2-1}\ \
\sum\limits_{(j_1,\ldots,j_k)}a_{(j_1,\ldots,j_k)},
$$

\vspace{2mm}

$$
{\rm C}_{s_r}\ldots {\rm C}_{s_1}S_N^{(k)}(a)=
$$

\vspace{-6mm}
$$
=
\sum_{j_k=0}^{N-1}\ldots \sum_{j_{s_r+1}=0}^{j_{s_r+2}-1}
\sum_{j_{s_r-1}=0}^{j_{s_r+1}-1}\ldots \sum_{j_{s_1+1}=0}^{j_{s_1+2}-1}
\sum_{j_{s_1-1}=0}^{j_{s_1+1}-1}\ldots \sum_{j_1=0}^{j_2-1}
\sum\limits_{\prod\limits_{l=1}^r{\bf I}_{j_{s_l},j_{s_l+1}}
(j_1,\ldots,j_k)}
a_{\prod\limits_{l=1}^r{\bf I}_{j_{s_l},j_{s_l+1}}
(j_1,\ldots,j_k)},
$$
where
$$
\prod\limits_{l=1}^r{\bf I}_{j_{s_l},j_{s_l+1}}
(j_1,\ldots,j_k)\ \ \stackrel{\rm def}{=}\
{\bf I}_{j_{s_r},j_{s_r+1}}\ldots
{\bf I}_{j_{s_1},j_{s_1+1}}
(j_1,\ldots,j_k),
$$
$$
{\rm C}_{s_0}\ldots {\rm C}_{s_1}
S_N^{(k)}(a)=S_N^{(k)}(a),\ \ \ 
\prod\limits_{l=1}^0{\bf I}_{j_{s_l},j_{s_l+1}}
(j_1,\ldots,j_k)=(j_1,\ldots,j_k),
$$
$$
{\bf I}_{j_l,j_{l+1}}(j_{q_1},\ldots,j_{q_2},j_l,j_{q_3},\ldots,
j_{q_{k-2}},j_l,j_{q_{k-1}},\ldots,j_{q_{k}})\stackrel{\rm def}{=}
$$
$$
\stackrel{\rm def}{=}(j_{q_1},\ldots,j_{q_2},j_{l+1},j_{q_3},\ldots,
j_{q_{k-2}},j_{l+1},j_{q_{k-1}},\ldots,j_{q_{k}}),
$$

\vspace{2mm}
\noindent
where
$l\ne q_{1},\ldots,q_2,q_3,\ldots,q_{k-2},q_{k-1},\ldots,q_{k},$
$l\in{\bf N},$
$a_{(j_{q_1},\ldots,j_{q_k})}$ is a scalar value,
$s_1,\ldots,s_r = 1,\ldots, k-1$, $s_r>\ldots >s_1,$
$q_1,\ldots,q_k=1,\ldots,k,$
the expression
$$
\sum\limits_{(j_{q_1},\ldots,j_{q_k})}
$$
means the sum with respect to
all possible  
permutations
$(j_{q_1},\ldots,j_{q_k}).$

Using induction it is possible to prove the following equality
\begin{equation}
\label{979}
~~~~~~~~ \sum_{j_k=0}^{N-1}\ldots \sum_{j_1=0}^{N-1}
a_{(j_1,\ldots,j_k)}
=\sum_{r=0}^{k-1}\ \
\sum_{\stackrel{\Large{s_r,\ldots,s_1=1}}
{{}_{s_r>\ldots>s_1}}}^{k-1}  
{\rm C}_{s_r}\ldots {\rm C}_{s_1}S_N^{(k)}(a),
\end{equation}
where $k=2, 3,\ldots$

Hereinafter in this section, we will identify the following records
$$
a_{(j_1,\ldots,j_k)}
=a_{(j_1\ldots j_k)}=a_{j_1\ldots j_k}.
$$

In particular, from (\ref{979}) for $k=2, 3, 4$ 
we get the following formulas
$$
\sum\limits_{j_2=0}^{N-1}\sum\limits_{j_1=0}^{N-1}
a_{(j_1,j_2)}=S_N^{(2)}(a)
+{\rm C}_{1}S_N^{(2)}(a)=
$$
$$
=
\sum\limits_{j_2=0}^{N-1}\sum\limits_{j_1=0}^{j_2-1}
\sum\limits_{(j_1,j_2)}a_{(j_1j_2)}+
\sum\limits_{j_2=0}^{N-1}
a_{(j_2j_2)}
=\sum\limits_{j_2=0}^{N-1}\sum\limits_{j_1=0}^{j_2-1}(a_{j_1j_2}+
a_{j_2j_1})+
$$
\begin{equation}
\label{best1}
+
\sum\limits_{j_2=0}^{N-1}
a_{j_2j_2},
\end{equation}
$$
\sum\limits_{j_3=0}^{N-1}\sum\limits_{j_2=0}^{N-1}\sum\limits_{j_1=0}^{N-1}
a_{(j_1,j_2,j_3)}=S_N^{(3)}(a)
+{\rm C}_{1}S_N^{(3)}(a)+
{\rm C}_{2}S_N^{(3)}(a)+{\rm C}_{2}{\rm C}_{1}S_N^{(3)}(a)=
$$
$$
=\sum\limits_{j_3=0}^{N-1}\sum\limits_{j_2=0}^{j_3-1}
\sum\limits_{j_1=0}^{j_2-1}
\sum\limits_{(j_1,j_2,j_3)}a_{(j_1j_2j_3)}+
\sum\limits_{j_3=0}^{N-1}\sum\limits_{j_2=0}^{j_3-1}
\sum\limits_{(j_2,j_2,j_3)}a_{(j_2j_2j_3)}+
$$
$$
+\sum\limits_{j_3=0}^{N-1}
\sum\limits_{j_1=0}^{j_3-1}
\sum\limits_{(j_1,j_3,j_3)}a_{(j_1j_3j_3)}+
\sum\limits_{j_3=0}^{N-1}a_{(j_3j_3j_3)}=
$$
$$
=\sum\limits_{j_3=0}^{N-1}\sum\limits_{j_2=0}^{j_3-1}
\sum\limits_{j_1=0}^{j_2-1}
\left(a_{j_1j_2j_3}+a_{j_1j_3j_2}+a_{j_2j_1j_3}+
a_{j_2j_3j_1}+a_{j_3j_2j_1}+a_{j_3j_1j_2}\right)+
$$
$$
+\sum\limits_{j_3=0}^{N-1}\sum\limits_{j_2=0}^{j_3-1}
\left(a_{j_2j_2j_3}+a_{j_2j_3j_2}+a_{j_3j_2j_2}\right)+
\sum\limits_{j_3=0}^{N-1}
\sum\limits_{j_1=0}^{j_3-1}
\left(a_{j_1j_3j_3}+a_{j_3j_1j_3}+a_{j_3j_3j_1}\right)+
$$
\begin{equation}
\label{oop12}
+
\sum\limits_{j_3=0}^{N-1}a_{j_3j_3j_3},
\end{equation}

\vspace{6mm}

$$
\sum\limits_{j_4=0}^{N-1}\sum\limits_{j_3=0}^{N-1}
\sum\limits_{j_2=0}^{N-1}\sum\limits_{j_1=0}^{N-1}
a_{(j_1,j_2,j_3,j_4)}=
S_N^{(4)}(a)
+{\rm C}_{1}S_N^{(4)}(a)+{\rm C}_{2}S_N^{(4)}(a)+
$$
$$
+{\rm C}_{3}S_N^{(4)}(a)+
{\rm C}_{2}{\rm C}_{1}S_N^{(4)}(a)+
{\rm C}_{3}{\rm C}_{1}S_N^{(4)}(a)+
{\rm C}_{3}{\rm C}_{2}S_N^{(4)}(a)+
{\rm C}_{3}{\rm C}_{2}{\rm C}_{1}S_N^{(4)}(a)=
$$
$$
=\sum\limits_{j_4=0}^{N-1}\sum\limits_{j_3=0}^{j_4-1}
\sum\limits_{j_2=0}^{j_3-1}\sum\limits_{j_1=0}^{j_2-1}
\sum\limits_{(j_1,j_2,j_3,j_4)}a_{(j_1j_2j_3j_4)}
+\sum\limits_{j_4=0}^{N-1}\sum\limits_{j_3=0}^{j_4-1}
\sum\limits_{j_2=0}^{j_3-1}
\sum\limits_{(j_2,j_2,j_3,j_4)}a_{(j_2j_2j_3j_4)}
$$
$$
+\sum\limits_{j_4=0}^{N-1}\sum\limits_{j_3=0}^{j_4-1}
\sum\limits_{j_1=0}^{j_3-1}
\sum\limits_{(j_1,j_3,j_3,j_4)}a_{(j_1j_3j_3j_4)}
+\sum\limits_{j_4=0}^{N-1}\sum\limits_{j_2=0}^{j_4-1}
\sum\limits_{j_1=0}^{j_2-1}
\sum\limits_{(j_1,j_2,j_4,j_4)}a_{(j_1j_2j_4j_4)}+
$$
$$
+\sum\limits_{j_4=0}^{N-1}
\sum\limits_{j_3=0}^{j_4-1}
\sum\limits_{(j_3,j_3,j_3,j_4)}a_{(j_3j_3j_3j_4)}
+\sum\limits_{j_4=0}^{N-1}
\sum\limits_{j_2=0}^{j_4-1}
\sum\limits_{(j_2,j_2,j_4,j_4)}a_{(j_2j_2j_4j_4)}+
$$
$$
+\sum\limits_{j_4=0}^{N-1}
\sum\limits_{j_1=0}^{j_4-1}
\sum\limits_{(j_1,j_4,j_4,j_4)}a_{(j_1j_4j_4j_4)}+
\sum\limits_{j_4=0}^{N-1}a_{j_4j_4j_4j_4}=
$$
$$
=\sum\limits_{j_4=0}^{N-1}\sum\limits_{j_3=0}^{j_4-1}
\sum\limits_{j_2=0}^{j_3-1}\sum\limits_{j_1=0}^{j_2-1}
\left(a_{j_1j_2j_3j_4}+a_{j_1j_2j_4j_3}+
a_{j_1j_3j_2j_4}+a_{j_1j_3j_4j_2}+\right.
$$
$$
+a_{j_1j_4j_3j_2}+a_{j_1j_4j_2j_3}+a_{j_2j_1j_3j_4}+
a_{j_2j_1j_4j_3}+
a_{j_2j_4j_1j_3}+a_{j_2j_4j_3j_1}+a_{j_2j_3j_1j_4}+
$$
$$
+a_{j_2j_3j_4j_1}+a_{j_3j_1j_2j_4}+a_{j_3j_1j_4j_2}+
a_{j_3j_2j_1j_4}+
a_{j_3j_2j_4j_1}+a_{j_3j_4j_1j_2}+a_{j_3j_4j_2j_1}+
$$
$$
+a_{j_4j_1j_2j_3}+a_{j_4j_1j_3j_2}+a_{j_4j_2j_1j_3}+
\left.a_{j_4j_2j_3j_1}+
a_{j_4j_3j_1j_2}+a_{j_4j_3j_2j_1}\right)+
$$
$$
+\sum\limits_{j_4=0}^{N-1}\sum\limits_{j_3=0}^{j_4-1}
\sum\limits_{j_2=0}^{j_3-1}
\left(a_{j_2j_2j_3j_4}+a_{j_2j_2j_4j_3}+a_{j_2j_3j_2j_4}+\right.
a_{j_2j_4j_2j_3}+
a_{j_2j_3j_4j_2}+a_{j_2j_4j_3j_2}+
$$
$$
+a_{j_3j_2j_2j_4}+a_{j_4j_2j_2j_3}+a_{j_3j_2j_4j_2}
\left.+
a_{j_4j_2j_3j_2}+
a_{j_4j_3j_2j_2}+a_{j_3j_4j_2j_2}\right)+
$$
$$
+\sum\limits_{j_4=0}^{N-1}\sum\limits_{j_3=0}^{j_4-1}
\sum\limits_{j_1=0}^{j_3-1}
\left(a_{j_3j_3j_1j_4}+a_{j_3j_3j_4j_1}+a_{j_3j_1j_3j_4}+\right.
a_{j_3j_4j_3j_1}+
a_{j_3j_4j_1j_3}+a_{j_3j_1j_4j_3}+
$$
$$
+a_{j_1j_3j_3j_4}+a_{j_4j_3j_3j_1}+a_{j_4j_3j_1j_3}
\left.+a_{j_1j_3j_4j_3}+
a_{j_1j_4j_3j_3}+a_{j_4j_1j_3j_3}\right)+
$$
$$
+\sum\limits_{j_4=0}^{N-1}\sum\limits_{j_2=0}^{j_4-1}
\sum\limits_{j_1=0}^{j_2-1}
\left(a_{j_4j_4j_1j_2}+a_{j_4j_4j_2j_1}+a_{j_4j_1j_4j_2}+\right.
a_{j_4j_2j_4j_1}+
a_{j_4j_2j_1j_4}+a_{j_4j_1j_2j_4}+
$$
$$
+a_{j_1j_4j_4j_2}+a_{j_2j_4j_4j_1}+a_{j_2j_4j_1j_4}+
\left.a_{j_1j_4j_2j_4}+
a_{j_1j_2j_4j_4}+a_{j_2j_1j_4j_4}\right)+
$$
$$
+\sum\limits_{j_4=0}^{N-1}
\sum\limits_{j_3=0}^{j_4-1}
\left(a_{j_3j_3j_3j_4}+a_{j_3j_3j_4j_3}+a_{j_3j_4j_3j_3}+
a_{j_4j_3j_3j_3}\right)+
$$
$$
+\sum\limits_{j_4=0}^{N-1}
\sum\limits_{j_2=0}^{j_4-1}
\left(a_{j_2j_2j_4j_4}+a_{j_2j_4j_2j_4}+a_{j_2j_4j_4j_2}+\right.\left.
a_{j_4j_2j_2j_4}+
a_{j_4j_2j_4j_2}+a_{j_4j_4j_2j_2}\right)+
$$
$$
+\sum\limits_{j_4=0}^{N-1}
\sum\limits_{j_1=0}^{j_4-1}
\left(a_{j_1j_4j_4j_4}+a_{j_4j_1j_4j_4}+a_{j_4j_4j_1j_4}+
a_{j_4j_4j_4j_1}\right)+
$$
\begin{equation}
\label{huh}
+\sum\limits_{j_4=0}^{N-1}a_{j_4j_4j_4j_4}.
\end{equation}

\vspace{3mm}

Perhaps, the formula (\ref{979}) for any $k$ ($k\in {\bf N}$) was 
found
by the author for the first time \cite{old-art-1} (1997).

Assume that 
$$
a_{(j_1,\ldots,j_k)}=
\Phi\left(\tau_{j_1},\ldots,\tau_{j_k}\right)
\prod\limits_{l=1}^k\Delta{\bf w}_{\tau_{j_l}}^{(i_l)},
$$ 
where $\Phi\left(t_1,\ldots,t_k\right)$ is a nonrandom
function of $k$ variables.
Then from (\ref{30.34})
and (\ref{979}) we have
$$
J[\Phi]_{T,t}^{(k)}=\sum_{r=0}^{[k/2]}\ \ 
\sum_{(s_r,\ldots,s_1)\in{\rm A}_{k,r}}  \times
$$

$$
\times\ \ 
\hbox{\vtop{\offinterlineskip\halign{
\hfil#\hfil\cr
{\rm l.i.m.}\cr
$\stackrel{}{{}_{N\to \infty}}$\cr
}} }
\sum_{j_k=0}^{N-1}\ldots \sum_{j_{s_r+1}=0}^{j_{s_r+2}-1}
\sum_{j_{s_r-1}=0}^{j_{s_r+1}-1}\ldots \sum_{j_{s_1+1}=0}^{j_{s_1+2}-1}
\sum_{j_{s_1-1}=0}^{j_{s_1+1}-1}\ldots \sum_{j_1=0}^{j_2-1}\ \ 
\sum\limits_{\prod\limits_{l=1}^r{\bf I}_{j_{s_l},j_{s_l+1}}
(j_1,\ldots,j_k)}\times
$$
$$
\times
\Biggl[\Phi\biggl(\tau_{j_1},\ldots,\tau_{j_{s_1-1}},
\tau_{j_{s_1+1}},\tau_{j_{s_1+1}},\tau_{j_{s_1+2}},\ldots
,\tau_{j_{s_r-1}},\tau_{j_{s_r+1}},\tau_{j_{s_r+1}},\tau_{j_{s_r+2}},
\ldots,\tau_{j_k}\biggr)
\times\Biggr.
$$
$$
\times
\Delta{\bf w}_{\tau_{j_1}}^{(i_1)}
\ldots\Delta{\bf w}_{\tau_{j_{s_1-1}}}^{(i_{s_1-1})}
\Delta{\bf w}_{\tau_{j_{s_1+1}}}^{(i_{s_1})}
\Delta{\bf w}_{\tau_{j_{s_1+1}}}^{(i_{s_1+1})}
\Delta{\bf w}_{\tau_{j_{s_1+2}}}^{(i_{s_1+2})}
\ldots
$$
$$
\Biggl.
\ldots\Delta{\bf w}_{\tau_{j_{s_r-1}}}^{(i_{s_r-1})}
\Delta{\bf w}_{\tau_{j_{s_r+1}}}^{(i_{s_r})}
\Delta{\bf w}_{\tau_{j_{s_r+1}}}^{(i_{s_r+1})}
\Delta{\bf w}_{\tau_{j_{s_r+2}}}^{(i_{s_r+2})}\ldots
\Delta{\bf w}_{\tau_{j_{k}}}^{(i_{k})}\Biggr]=
$$
\begin{equation}
\label{30.51}
=\sum_{r=0}^{[k/2]}\sum_{(s_r,\ldots,s_1)\in{\rm A}_{k,r}}
I[\Phi]_{T,t}^{(k)s_1,\ldots,s_r}\ \ \ \hbox{w.~p.~1},
\end{equation}

\noindent
where
$$
I[\Phi]_{T,t}^{(k)s_1,\ldots,s_r}=
\int\limits_t^T\ldots\int\limits_t^{t_{s_r+3}}
\int\limits_t^{t_{s_r+2}}\int\limits_t^{t_{s_r}}\ldots
\int\limits_t^{t_{s_1+3}}
\int\limits_t^{t_{s_1+2}}\int\limits_t^{t_{s_1}}\ldots\int\limits_t^{t_2}
\sum\limits_{\prod\limits_{l=1}^r{\bf I}_{t_{s_l},t_{s_l+1}}
(t_1,\ldots,t_k)}\times
$$
$$
\times
\Biggl[\Phi\biggl(t_{1},\ldots,t_{s_1-1},t_{s_1+1},t_{s_1+1},t_{s_1+2},\ldots
,t_{s_r-1},t_{s_r+1},t_{s_r+1},t_{s_r+2},
\ldots,t_{k}\biggr)\times\Biggr.
$$
$$
\times
d{\bf w}_{t_1}^{(i_1)}\ldots d{\bf w}_{t_{s_1-1}}^{(i_{s_1-1})}
d{\bf w}_{t_{s_1+1}}^{(i_{s_1})}
d{\bf w}_{t_{s_1+1}}^{(i_{s_1+1})}
d{\bf w}_{t_{s_1+2}}^{(i_{s_1+2})}
\ldots 
$$
\begin{equation}
\label{99999}
\Biggl.
\ldots
d{\bf w}_{t_{s_r-1}}^{(i_{s_r-1})}
d{\bf w}_{t_{s_r+1}}^{(i_{s_r})}
d{\bf w}_{t_{s_r+1}}^{(i_{s_r+1})}
d{\bf w}_{t_{s_r+2}}^{(i_{s_r+2})}\ldots
d{\bf w}_{t_k}^{(i_k)}\Biggr],
\end{equation}

\noindent
where 
$k\ge 2,$ the set ${\rm A}_{k,r}$ is defined by the
relation (\ref{30.5550001}).
We suppose that
the right-hand side of (\ref{99999}) exists as the It\^{o} stochastic integral.

{\bf Remark 2.3.}\ {\it The summands on the right-hand side
of {\rm (\ref{99999})} should be 
understood
as follows{\rm :} for each 
permutation
from the set 
$$
\prod\limits_{l=1}^r{\bf I}_{t_{s_l},t_{s_l+1}}
(t_1,\ldots,t_k)
=
$$
$$
=
\biggl(t_{1},\ldots,t_{s_1-1},t_{s_1+1},t_{s_1+1},t_{s_1+2},\ldots
,t_{s_r-1},t_{s_r+1},t_{s_r+1},t_{s_r+2},
\ldots,t_{k}\biggr)
$$

\noindent
it is necessary to perform replacement on the right-hand side 
of {\rm (\ref{99999})} 
of all pairs {\rm (}their number is equal to 
$r${\rm )} of differentials $d{\bf w}_{t_p}^{(i)}d{\bf w}_{t_p}^{(j)}$
with similar 
lower indices
by the values ${\bf 1}_{\{i=j\ne 0\}}dt_p.$}

Note that the term in (\ref{30.51}) for $r=0$ should be understood as
follows
$$
\int\limits_t^T\ldots \int\limits_t^{t_2}
\sum\limits_{(t_1,\ldots,t_k)}\biggl(
\Phi\left(t_1,\ldots,t_k\right)
d{\bf w}_{t_{1}}^{(i_{1})}\ldots
d{\bf w}_{t_k}^{(i_k)}\biggr),
$$
where notations are the same as in (\ref{pobeda}).

Using (\ref{pupol1}), (\ref{pupol2}), (\ref{30.51}), and
(\ref{99999}), we get

$$
{\sf M}\left\{\left|J[\Phi]_{T,t}^{(k)}\right|^{2n}\right\}
\le 
$$

\vspace{-2mm}
\begin{equation}
\label{333.225}
\le C_{nk}\ \ \sum_{r=0}^{[k/2]}\sum_{(s_r,\ldots,s_1)\in{\rm A}_{k,r}}
{\sf M}\left\{\left|
I[\Phi]_{T,t}^{(k)s_1,\ldots,s_r}\right|^{2n}\right\},
\end{equation}

\vspace{2mm}
\noindent
where
$$
{\sf M}\left\{\left|
I[\Phi]_{T,t}^{(k)s_1,\ldots,s_r}\right|^{2n}\right\}\le
$$

\vspace{-2mm}
$$
\le C_{nk}^{s_1\ldots s_r}
\int\limits_t^T\ldots\int\limits_t^{t_{s_r+3}}
\int\limits_t^{t_{s_r+2}}\int\limits_t^{t_{s_r}}\ldots
\int\limits_t^{t_{s_1+3}}
\int\limits_t^{t_{s_1+2}}\int\limits_t^{t_{s_1}}\ldots\int\limits_t^{t_2}
\sum\limits_{\prod\limits_{l=1}^r{\bf I}_{t_{s_l},t_{s_l+1}}
(t_1,\ldots,t_k)}\times
$$
$$
\times
\Phi^{2n}\biggl(t_{1},\ldots,t_{s_1-1},t_{s_1+1},t_{s_1+1},t_{s_1+2},\ldots,
t_{s_r-1},t_{s_r+1},t_{s_r+1},t_{s_r+2},\ldots,t_{k}\biggr)\times
$$

\vspace{-6mm}
\begin{equation}
\label{333.226}
~~~~~~~~~~\times
dt_1\ldots dt_{s_1-1}
dt_{s_1+1}dt_{s_1+2}
\ldots dt_{s_r-1}dt_{s_r+1}
dt_{s_r+2}\ldots
dt_k,
\end{equation}

\vspace{2mm}
\noindent
where $C_{nk}$ and $C_{nk}^{s_1\ldots s_r}$ are constants
and
permutations
when summing are 
performed in
(\ref{333.226}) only in the values
$$
\Phi^{2n}\biggl(t_{1},\ldots,t_{s_1-1},t_{s_1+1},t_{s_1+1},t_{s_1+2},\ldots,
t_{s_r-1},t_{s_r+1},t_{s_r+1},t_{s_r+2},\ldots,t_{k}\biggr).
$$

Consider (\ref{333.225}) and (\ref{333.226}) for
$\Phi(t_1,\ldots,t_k)\equiv R_{p_1\ldots p_k}(t_1,\ldots,t_k)$

$$
{\sf M}\left\{\left|J[R_{p_1\ldots p_k}]_{T,t}^{(k)}\right|^{2n}\right\}
\le 
$$

\vspace{-1mm}
\begin{equation}
\label{333.225e}
\le
C_{nk}\ \ \sum_{r=0}^{[k/2]}\sum_{(s_r,\ldots,s_1)\in{\rm A}_{k,r}}
{\sf M}\left\{\left|
I[R_{p_1\ldots p_k}]_{T,t}^{(k)s_1,\ldots,s_r}\right|^{2n}\right\},
\end{equation}

\vspace{2mm}
\noindent
where
$$
{\sf M}\left\{\left|
I[R_{p_1\ldots p_k}]_{T,t}^{(k)s_1,\ldots,s_r}\right|^{2n}\right\}\le
$$

\vspace{-2mm}
$$
\le C_{nk}^{s_1\ldots s_r}
\int\limits_t^T\ldots\int\limits_t^{t_{s_r+3}}
\int\limits_t^{t_{s_r+2}}\int\limits_t^{t_{s_r}}\ldots
\int\limits_t^{t_{s_1+3}}
\int\limits_t^{t_{s_1+2}}\int\limits_t^{t_{s_1}}\ldots\int\limits_t^{t_2}
\sum\limits_{\prod\limits_{l=1}^r{\bf I}_{t_{s_l},t_{s_l+1}}
(t_1,\ldots,t_k)}\times
$$
$$
\times
R_{p_1\ldots p_k}^{2n}
\biggl(t_{1},\ldots,t_{s_1-1},t_{s_1+1},t_{s_1+1},t_{s_1+2},\ldots,
t_{s_r-1},t_{s_r+1},t_{s_r+1},t_{s_r+2},\ldots,t_{k}\biggr)\times
$$

\vspace{-5mm}
\begin{equation}
\label{333.226e}
~~~~~~~~~~\times
dt_1\ldots dt_{s_1-1}
dt_{s_1+1}dt_{s_1+2}
\ldots dt_{s_r-1}dt_{s_r+1}
dt_{s_r+2}\ldots
dt_k,
\end{equation}

\vspace{3mm}
\noindent
where $C_{nk}$ and $C_{nk}^{s_1\ldots s_r}$ are constants and
permutations
when summing are 
performed in
(\ref{333.226e}) only in the values
$$
R_{p_1\ldots p_k}^{2n}
\biggl(t_{1},\ldots,t_{s_1-1},t_{s_1+1},t_{s_1+1},t_{s_1+2},\ldots,
t_{s_r-1},t_{s_r+1},t_{s_r+1},t_{s_r+2},\ldots,t_{k}\biggr).
$$

From the other hand, we can consider
the generalization of the formulas (\ref{leto80010aaa}),
(\ref{oop16xxx}), (\ref{udar}) 
for the case of arbitrary $k$ ($k\in{\bf N}$).
In order to do this, let us
consider  
the sum with respect to all possible
partitions defined by (\ref{leto5008})
$$
\sum_{\stackrel{(\{\{g_1, g_2\}, \ldots, 
\{g_{2r-1}, g_{2r}\}\}, \{q_1, \ldots, q_{k-2r}\})}
{{}_{\{g_1, g_2, \ldots, 
g_{2r-1}, g_{2r}, q_1, \ldots, q_{k-2r}\}=\{1, 2, \ldots, k\}}}}
a_{g_1 g_2, \ldots, 
g_{2r-1} g_{2r}, q_1 \ldots q_{k-2r}}.
$$

Now we can 
generalize the formulas (\ref{leto80010aaa}),
(\ref{oop16xxx}), (\ref{udar}) 
for the case of arbitrary $k$ ($k\in{\bf N}$)
$$
{\sf M}\left\{\left|J[R_{p_1\ldots p_k}]_{T,t}^{(k)}\right|^{2n}\right\}
\le C_{nk}\left(~
\int\limits_{\stackrel{~}{[t, T]^k}}
\left(R_{p_1\ldots p_k}(t_1,\ldots,t_k)\right)^{2n}dt_1\ldots dt_k+
\right.
$$

\vspace{-2mm}
$$
+\sum\limits_{r=1}^{[k/2]}
\sum_{\stackrel{(\{\{g_1, g_2\}, \ldots, 
\{g_{2r-1}, g_{2r}\}\}, \{q_1, \ldots, q_{k-2r}\})}
{{}_{\{g_1, g_2, \ldots, 
g_{2r-1}, g_{2r}, q_1, \ldots, q_{k-2r}\}=\{1, 2, \ldots, k\}}}}
{\bf 1}_{\{i_{g_{{}_{1}}}=i_{g_{{}_{2}}}\ne 0\}}
\ldots 
{\bf 1}_{\{i_{g_{{}_{2r-1}}}=i_{g_{{}_{2r}}}\ne 0\}}\times
$$

\vspace{2mm}
$$
\times
\int\limits_{\stackrel{~}{[t, T]^{k-r}}}
\left(R_{p_1\ldots p_k}\biggl(
t_1,\ldots,t_k\biggr)\biggl.\biggr|_{t_{g_{{}_{1}}}=t_{g_{{}_{2}}},\ldots,
t_{g_{{}_{2r-1}}}=t_{g_{{}_{2r}}}}
\right)^{2n}\times
$$

\vspace{2mm}
\begin{equation}
\label{udar1}
~~~~~ \left.\times \biggl(dt_1\ldots dt_k\biggr)
\biggl|_{\left(dt_{g_{{}_{1}}}dt_{g_{{}_{2}}}\right)\biggr.
\curvearrowright dt_{g_1},\ldots,
\left(dt_{g_{{}_{2r-1}}}dt_{g_{{}_{2r}}}\right)\curvearrowright
dt_{g_{{}_{2r-1}}}}\right),
\end{equation}

\vspace{6mm}
\noindent
where $C_{nk}$ is a constant,
$$
\biggl(
t_1,\ldots,t_k\biggr)\biggl.\biggr|_{t_{g_{{}_{1}}}=t_{g_{{}_{2}}},\ldots,
t_{g_{{}_{2r-1}}}=t_{g_{{}_{2r}}}}
$$

\vspace{2mm}
\noindent
means the ordered set $(t_1,\ldots,t_k)$, where we put
$t_{g_{{}_{1}}}=t_{g_{{}_{2}}},$ $\ldots,$
$t_{g_{{}_{2r-1}}}=t_{g_{{}_{2r}}}.$

Moreover, 
$$
\biggl(dt_1\ldots dt_k\biggr)\biggl.
\biggr|_{\left(dt_{g_{{}_{1}}}dt_{g_{{}_{2}}}\right)
\curvearrowright dt_{g_1},\ldots,
\left(dt_{g_{{}_{2r-1}}}dt_{g_{{}_{2r}}}\right)\curvearrowright
dt_{g_{{}_{2r-1}}}}
$$

\vspace{2mm}
\noindent
means the product $dt_1\ldots dt_k$, where we replace
all pairs 
$dt_{g_{{}_{1}}}dt_{g_{{}_{2}}},$ $\ldots,$ 
$dt_{g_{{}_{2r-1}}}dt_{g_{{}_{2r}}}$ by 
$dt_{g_1},$ $\ldots,$ $dt_{g_{{}_{2r-1}}}$
correspondingly.

Note that the estimate like (\ref{udar1}), 
where all indicators ${\bf 1}_{\{\cdot\}}$ must be 
replaced with $1$, can be obtained from 
the estimates (\ref{333.225e}), (\ref{333.226e}).

The comparison of (\ref{udar1}) with the formula 
(\ref{leto6000}) (see Theorem 1.2)
shows their
similar structure.

Let us consider the particular case of (\ref{udar1}) for $k=4$

\newpage
\noindent
$$
{\sf M}\left\{\left|J[R_{p_1 p_2 p_3 p_4}]_{T,t}^{(4)}\right|^{2n}\right\}
\le C_{n4}\left(~
\int\limits_{\stackrel{~}{[t, T]^4}}
\left(R_{p_1 p_2 p_3 p_4}(t_1,t_2,t_3,t_4)\right)^{2n}dt_1 dt_2 dt_3 dt_4+
\right.
$$
$$
+
\sum_{\stackrel{(\{g_1, g_2\}, \{q_1, q_{2}\})}
{{}_{\{g_1, g_2, q_1, q_{2}\}=\{1, 2, 3, 4\}}}}
{\bf 1}_{\{i_{g_{{}_{1}}}=i_{g_{{}_{2}}}\ne 0\}}
\int\limits_{\stackrel{~}{[t, T]^3}}
\left(R_{p_1 p_2 p_3 p_4}\biggl(
t_1,t_2,t_3,t_4\biggr)\biggl.\biggr|_{t_{g_{{}_{1}}}=t_{g_{{}_{2}}}}
\right)^{2n}
\times
$$

\vspace{1mm}
$$
\times
\biggl(dt_1 dt_2 dt_3 dt_4\biggr)
\Biggl.\Biggr|_{\left(dt_{g_{{}_{1}}}dt_{g_{{}_{2}}}\right)
\curvearrowright dt_{g_1}}+
$$

\vspace{3mm}
$$
+
\sum_{\stackrel{(\{\{g_1, g_2\}, \{g_3, g_{4}\}\})}
{{}_{\{g_1, g_2, g_3, g_{4}\}=\{1, 2, 3, 4\}}}}
{\bf 1}_{\{i_{g_{{}_{1}}}=i_{g_{{}_{2}}}\ne 0\}}
{\bf 1}_{\{i_{g_{{}_{3}}}=i_{g_{{}_{4}}}\ne 0\}}
\times
$$

$$
\times
\int\limits_{\stackrel{~}{[t, T]^2}}
\left(R_{p_1 p_2 p_3 p_4}\biggl(
t_1,t_2,t_3,t_4\biggr)\biggl.\biggr|_{t_{g_{{}_{1}}}=t_{g_{{}_{2}}},
t_{g_{{}_{3}}}=t_{g_{{}_{4}}}}
\right)^{2n}
\times
$$

\begin{equation}
\label{f112}
\times\left.
\biggl(dt_1 dt_2 dt_3 dt_4\biggr)
\biggl|_{\left(dt_{g_{{}_{1}}}dt_{g_{{}_{2}}}\right)\biggr.
\curvearrowright dt_{g_1},
\left(dt_{g_{{}_{3}}}dt_{g_{{}_{4}}}\right)
\curvearrowright dt_{g_3}}\right).
\end{equation}

\vspace{5mm}

It is not difficult to notice that (\ref{f112}) is consistent with 
(\ref{udar}).

According to (\ref{1999.1xx}) and (\ref{30.46xx}), we have the following 
expression 

\vspace{-2mm}
$$
R_{p_1\ldots p_k}(t_1,\ldots,t_k)=
$$

\vspace{-7mm}
$$
=
\prod_{l=1}^k \psi_l(t_l)\left(\prod_{l=1}^{k-1}
{\bf 1}_{\{t_l<t_{l+1}\}}+
\sum_{r=1}^{k-1}\frac{1}{2^r}
\sum_{\stackrel{s_r,\ldots,s_1=1}{{}_{s_r>\ldots>s_1}}}^{k-1}
\prod_{l=1}^r {\bf 1}_{\{t_{s_l}=t_{s_l+1}\}}
\prod_{\stackrel{l=1}{{}_{l\ne s_1,\ldots, s_r}}}^{k-1}
{\bf 1}_{\{t_{l}<t_{l+1}\}}\right)-
$$

\vspace{-2mm}
\begin{equation}
\label{30.48}
-\sum_{j_1=0}^{p_1}\ldots\sum_{j_k=0}^{p_k}
C_{j_k\ldots j_1} \prod_{l=1}^{k} \phi_{j_l}(t_l).
\end{equation}

\vspace{2mm}

Due to (\ref{30.48}) 
the function $R_{p_1\ldots p_k}(t_1,\ldots,t_k)$
is continuous 
in the open
domains of integration of integrals on the right-hand side
of (\ref{333.226e}) 
and it is bounded
at the 
boundaries  
of these 
domains for $p_1,\ldots,p_k<\infty.$

Let us perform the iterated  passage to the limit
$\lim\limits_{p_1\to \infty}\varlimsup\limits_{p_2\to \infty}
\ldots\varlimsup\limits_{p_k\to \infty}$
under the  
integral signs on the right-hand side of the estimate (\ref{udar1})
(it was similarly performed for 
the 2-dimensional, 3-dimentional, and 4-dimensional cases (see above)).
Then, taking into account
(\ref{410}), we obtain the required 
result. More precisely,
since the integrals on the right-hand side of (\ref{udar1}) 
exist as Riemann integrals, then they are equal to the 
corresponding Lebesgue integrals. 
Moreover, 

\vspace{-3mm}
$$
\lim\limits_{p_1\to\infty}\ldots \lim\limits_{p_k\to\infty}
R_{p_1\ldots p_k}(t_1,\ldots,t_k)=0\ \ \ \hbox{when}\ \ \ 
(t_1,\ldots,t_k)\in (t, T)^k,
$$

\vspace{2mm}
\noindent
where the left-hand side is bounded on the boundary of 
$[t, T]^k.$

According to the proof of Theorem 2.11 and (\ref{30.46xx}), we have

$$
R_{p_1\ldots p_k}(t_1,\ldots,t_k)=
$$

\vspace{-1mm}
$$
=
\left(K^{*}(t_1,\ldots,t_k)-\sum\limits_{j_1=0}^{p_1}
C_{j_1}(t_2,\ldots,t_k)\phi_{j_1}(t_1)\right)+
$$

\vspace{-3mm}
$$
+\left(
\sum\limits_{j_1=0}^{p_1}\left(C_{j_1}(t_2,\ldots,t_k)-
\sum\limits_{j_2=0}^{p_2}
C_{j_2j_1}(t_3,\ldots,t_k)\phi_{j_2}(t_2)\right)
\phi_{j_1}(t_1)\right)+
$$

\vspace{-4.5mm}
$$
\ldots
$$

\vspace{-7.5mm}
\begin{equation}
\label{strange1}
+\left(
\sum\limits_{j_1=0}^{p_1}\ldots 
\sum\limits_{j_{k-1}=0}^{p_{k-1}}\left(C_{j_{k-1}\ldots j_1}(t_k)-
\sum\limits_{j_k=0}^{p_k}
C_{j_k\ldots j_1}\phi_{j_k}(t_k)\right)
\phi_{j_{k-1}}(t_{k-1})\ldots \phi_{j_1}(t_1)\right),
\end{equation}

\vspace{3mm}
\noindent
where
$$
C_{j_1}(t_2,\ldots,t_k)=\int\limits_t^T
K^{*}(t_1,\ldots,t_k)\phi_{j_1}(t_1)dt_1,
$$

\vspace{2mm}
$$
C_{j_2j_1}(t_3,\ldots,t_k)=\int\limits_{[t, T]^2}
K^{*}(t_1,\ldots,t_k)\phi_{j_1}(t_1)\phi_{j_2}(t_2)dt_1 dt_2,
$$

\newpage
\noindent
$$
\ldots
$$
$$
C_{j_{k-1}\ldots j_1}(t_k)=\int\limits_{[t, T]^{k-1}}
K^{*}(t_1,\ldots,t_k)\prod\limits_{l=1}^{k-1}
\phi_{j_l}(t_l)dt_1\ldots dt_{k-1}.
$$

\vspace{2mm}

Then, applying $k$ times (we mean here an iterated passage to the limit
$\lim\limits_{p_1\to\infty}\varlimsup\limits_{p_2\to\infty}\ldots \varlimsup\limits_{p_k\to\infty}$)
the Lebesgue's Dominated Convergence Theorem to the integrals
on the right-hand side of (\ref{udar1}),
we obtain

\vspace{-1.5mm}
$$
\lim\limits_{p_1\to\infty}
\varlimsup\limits_{p_2\to\infty}\ldots \varlimsup\limits_{p_k\to\infty}
{\sf M}\left\{\left|J[R_{p_1\ldots p_k}]_{T,t}^{(k)}\right|^{2n}
\right\}=0,\ \ \ n\in {\bf N}.
$$

\vspace{4mm}

Let us discuss the choice of integrable majorants
when applying Lebesgue's 
Dominated Convergence Theorem 
in (\ref{udar1}).

It is well known that \cite{bari}

\vspace{-3mm}
\begin{equation}
\label{strange2}
\left|\sum\limits_{k=1}^N \frac{\sin kx}{k}\right|\le C
\end{equation}

\vspace{3mm}
\noindent
for all $N$ and $x$, where constant $C$ does not depend on $N$ and $x$.

Moreover,

\vspace{-5.5mm}
\begin{equation}
\label{strange3}
\sum\limits_{j=1}^N \frac{1}{j^2}\le \sum\limits_{j=1}^{\infty} \frac{1}{j^2}=
\frac{\pi^2}{6}.
\end{equation}

\vspace{3mm}

Applying double integration by parts (as in (\ref{2017zzz1})), we estimate
the partial sums of one-dimensional trigonometric Fourier series

\vspace{-2mm}
$$
\sum\limits_{j_1=0}^{p_1}
C_{j_1}(t_2,\ldots,t_k)\phi_{j_1}(t_1),
$$

\vspace{-2mm}
$$
\sum\limits_{j_2=0}^{p_2}
C_{j_2j_1}(t_3,\ldots,t_k)\phi_{j_2}(t_2),
$$

\vspace{-8mm}
$$
\ldots
$$

\vspace{-7mm}
$$  
\sum\limits_{j_k=0}^{p_k}
C_{j_k\ldots j_1}\phi_{j_k}(t_k)
$$ 

\vspace{5mm}
\noindent
in (\ref{strange1}) using (\ref{strange3}) and (see (\ref{strange2}))

\newpage
\noindent
$$
\left|\sum\limits_{k=1}^N \frac{1}{k}\sin\frac{2\pi k (x-y)}{T-t}\right|\le C,
$$

$$
\left|\sum\limits_{k=1}^N \frac{1}{k}\sin\frac{2\pi k (x-t)}{T-t}\right|\le C
$$

\vspace{2mm}
\noindent 
(here $N\in {\bf N}$ and $x, y\in {\bf R}$, constant $C$ does not depend on $N$ and $x, y$) as follows

\vspace{-3mm}
$$
\left|\sum\limits_{j_1=0}^{p_1}
C_{j_1}(t_2,\ldots,t_k)\phi_{j_1}(t_1)\right|\le C_1,
$$

\vspace{-2mm}
$$
\left|\sum\limits_{j_1=0}^{p_1}
C_{j_1}(t_2,\ldots,t_k)\phi_{j_1}(t_1)\right|\le C_2,
$$

\vspace{-8mm}
$$
\ldots
$$

\vspace{-6mm}
$$
\left|\sum\limits_{j_k=0}^{p_k}
C_{j_k\ldots j_1}\phi_{j_k}(t_k)\right|\le C_k,
$$

\vspace{3mm}
\noindent
where constant $C_1$ does not depend on $p_1,$ constant $C_2$ does not depend on $p_2,$ etc.

Moreover, 

\vspace{-6mm}
$$
\left|K^{*}(t_1,\ldots,t_k)\right|\le \tilde C_1,\ \ \ 
\left|C_{j_1}(t_2,\ldots,t_k)\right|\le \tilde C_2,\ \ \ 
\ldots\ \ \
\left|C_{j_{k-1}\ldots j_1}(t_k)\right|\le \tilde C_k,
$$

\vspace{1.5mm}
\noindent
where constant $\tilde C_1$ does not depend on $p_1,$ 
constant $\tilde C_2$ does not depend on $p_2,$ etc.

Further, the construction of 
integrable majorants
when applying Lebesgue's 
Dominated Convergence Theorem 
in (\ref{udar1}) is obvious.

For example, to pass to the limit
$\varlimsup\limits_{p_k\to\infty},$
the integrable majorant has the form
(it is constructed on the base of (\ref{strange1}))

\vspace{-2mm}
$$
\biggl(R_{p_1\ldots p_k}(t_1,\ldots,t_k)\biggr)^{2n}\le 
$$
$$
\le\Biggl(\left(\tilde C_1 + C_1\right)+\Biggr.
$$
$$
+
\sum\limits_{j_1=0}^{p_1}\left(\tilde C_2 + C_2\right)
\left|\phi_{j_1}(t_1)\right|+
\ldots\ 
$$

\vspace{-1mm}
$$
\Biggl.\hspace{-10mm}
\ldots\ +
\sum\limits_{j_1=0}^{p_1}\ldots 
\sum\limits_{j_{k-1}=0}^{p_{k-1}}\left(\tilde C_k + C_k\right)
\left|\phi_{j_{k-1}}(t_{k-1})\ldots \phi_{j_1}(t_1)\right|\Biggr)^{2n}\le
$$

\vspace{-1mm}
$$
\le \Biggl(\left(\tilde C_1 + C_1\right)+\Biggr.
$$

\vspace{-1mm}
$$
+
\sqrt{\frac{2}{T-t}}\ (p_1+1)\left(\tilde C_2 + C_2\right)+
\ldots\
$$

\vspace{-1mm}
\begin{equation}
\label{strange6}
~~~\Biggl.\ldots\ +\Biggl(\sqrt{\frac{2}{T-t}}\Biggr)^{k-1}
(p_1+1)\ldots (p_{k-1}+1)
\left(\tilde C_k + C_k\right)\Biggr)^{2n},
\end{equation}

\vspace{3mm}
\noindent
where $n\in{\bf N},$ the numbers $p_1,\ldots,p_{k-1}$ are fixed and the right-hand side of
(\ref{strange6}) is independent of $p_k.$

Theorems 2.13 and 2.10 are proved.

It is easy to notice that if we expand the function
$K^{*}(t_1,\ldots,t_k)$ into the generalized 
Fourier series at the interval $(t, T)$
at first with respect to the variable $t_k$, after that
with respect to the variable $t_{k-1}$, etc., then
we will have the expansion

\vspace{-4mm}
\begin{equation}
\label{otit3333}
~~~~~~ K^{*}(t_1,\ldots,t_k)=
\lim\limits_{p_k\to\infty}\ldots \lim\limits_{p_1\to\infty}
\sum_{j_k=0}^{p_k}\ldots\sum_{j_1=0}^{p_1}
C_{j_k\ldots j_1}
\prod_{l=1}^k \phi_{j_l}(t_l)
\end{equation}

\vspace{3mm}
\noindent
instead of the expansion (\ref{30.18}).

Let us prove the expansion (\ref{otit3333}). Similarly 
with (\ref{oop1}) we have the following equality

\vspace{-5mm}
\begin{equation}
\label{oop1otit}
~~ \psi_k(t_k)\left({\bf 1}_{\{t_{k-1}<t_k\}}+
\frac{1}{2}{\bf 1}_{\{t_{k-1}=t_k\}}\right)=
\sum\limits_{j_k=0}^{\infty}\ \int\limits_{t_{k-1}}^{T}\psi_k(t_k)
\phi_{j_k}(t_k)dt_k\phi_{j_k}(t_k),
\end{equation}

\vspace{2mm}
\noindent
which is 
satisfied pointwise at the interval $(t, T),$
besides
the series on the right-hand side
of (\ref{oop1otit}) converges when $t_1=t, T.$

Let us introduce the induction assumption 

\vspace{-5mm}
$$
\sum\limits_{j_k=0}^{\infty}\ldots
\sum\limits_{j_{3}=0}^{\infty}\psi_{2}(t_{2})
\int\limits_{t_2}^T\psi_{3}(t_{3})\phi_{j_{3}}(t_{3})\ldots
\int\limits_{t_{k-1}}^T
\psi_{k}(t_{k})\phi_{j_{k}}(t_{k})dt_k\ldots dt_{3}
\prod_{l=3}^{k}\phi_{j_{l}}(t_{l})=
$$

\vspace{-1mm}
\begin{equation}
\label{oop22otit}
=\prod\limits_{l=2}^{k}\psi_l(t_l)
\prod_{l=2}^{k-1}\left({\bf 1}_{\{t_l<t_{l+1}\}}+
\frac{1}{2}{\bf 1}_{\{t_l=t_{l+1}\}}\right).
\end{equation}

\vspace{3mm}

Then

\vspace{-6mm}
$$
\sum\limits_{j_k=0}^{\infty}\ldots
\sum\limits_{j_{3}=0}^{\infty}\sum\limits_{j_{2}=0}^{\infty}
\psi_{1}(t_{1})
\int\limits_{t_1}^{T}\psi_{2}(t_{2})\phi_{j_{2}}(t_{2})\ldots
\int\limits_{t_{k-1}}^{T}\psi_{k}(t_{k})
\phi_{j_{k}}(t_{k})
dt_k\ldots dt_{2}\prod_{l=2}^{k}\phi_{j_{l}}(t_{l})=
$$

$$
=\sum\limits_{j_k=0}^{\infty}\ldots
\sum\limits_{j_{3}=0}^{\infty}
\psi_1(t_1)\left({\bf 1}_{\{t_{1}<t_{2}\}}+
\frac{1}{2}{\bf 1}_{\{t_{1}=t_{2}\}}\right)\psi_{2}(t_{2})\times
$$

$$
\times
\int\limits_{t_2}^{T}\psi_{3}(t_{3})\phi_{j_{3}}(t_{3})\ldots
\int\limits_{t_{k-1}}^T\psi_{k}(t_{k})
\phi_{j_{k}}(t_{k})dt_k\ldots dt_{3}
\prod_{l=3}^{k}\phi_{j_{l}}(t_{l})=
$$

$$
=\psi_1(t_1)\left({\bf 1}_{\{t_{1}<t_{2}\}}+
\frac{1}{2}{\bf 1}_{\{t_{1}=t_{2}\}}\right)
\sum\limits_{j_k=0}^{\infty}\ldots
\sum\limits_{j_{3}=0}^{\infty}
\psi_{2}(t_{2})\times
$$

$$
\times
\int\limits_{t_2}^T\psi_{3}(t_{3})\phi_{j_{3}}(t_{3})\ldots
\int\limits_{t_{k-1}}^T\psi_{k}(t_{k})
\phi_{j_{k}}(t_{k})dt_k\ldots dt_{3}
\prod_{l=3}^{k}\phi_{j_{l}}(t_{l})=
$$

\vspace{-1mm}
$$
=\psi_1(t_1)\left({\bf 1}_{\{t_{1}<t_{2}\}}+
\frac{1}{2}{\bf 1}_{\{t_{1}=t_{2}\}}\right)
\prod\limits_{l=2}^{k}\psi_l(t_l)
\prod_{l=2}^{k-1}\left({\bf 1}_{\{t_l<t_{l+1}\}}+
\frac{1}{2}{\bf 1}_{\{t_l=t_{l+1}\}}\right)=
$$

\vspace{3mm}
\begin{equation}
\label{oop30otit}
=\prod\limits_{l=1}^{k}\psi_l(t_l)
\prod_{l=1}^{k-1}\left({\bf 1}_{\{t_l<t_{l+1}\}}+
\frac{1}{2}{\bf 1}_{\{t_l=t_{l+1}\}}\right).
\end{equation}

\vspace{3mm}

From the other hand, the left-hand side
of (\ref{oop30otit}) can be represented 
in the following form
$$
\sum_{j_k=0}^{\infty}\ldots \sum_{j_1=0}^{\infty}
C_{j_k\ldots j_1}\prod_{l=1}^{k} \phi_{j_l}(t_l)
$$

\vspace{2mm}
\noindent
by 
expanding the function

\vspace{-2mm}
$$
\psi_{1}(t_{1})
\int\limits_{t_1}^{T}\psi_{2}(t_{2})\phi_{j_{2}}(t_{2})\ldots
\int\limits_{t_{k-1}}^{T}\psi_{k}(t_{k})
\phi_{j_{k}}(t_{k})
dt_k\ldots dt_{2}
$$

\vspace{2mm}
\noindent
into the generalized Fourier series at the interval $(t, T)$ 
using the variable 
$t_1.$
Here we applied the following replacement of 
integration order

\vspace{-3mm}
$$
\int\limits_t^T\psi_{1}(t_{1})\phi_{j_{1}}(t_{1})
\int\limits_{t_1}^{T}\psi_{2}(t_{2})\phi_{j_{2}}(t_{2})\ldots
\int\limits_{t_{k-1}}^{T}\psi_{k}(t_{k})
\phi_{j_{k}}(t_{k})
dt_k\ldots dt_2dt_{1}=
$$

\vspace{-3mm}
$$
=\int\limits_t^T\psi_{k}(t_{k})\phi_{j_{k}}(t_{k})
\ldots
\int\limits_t^{t_{3}}\psi_{2}(t_{2})
\phi_{j_{2}}(t_{2})\int\limits_t^{t_{2}}\psi_{1}(t_{1})
\phi_{j_{1}}(t_{1})
dt_1dt_2\ldots dt_k=
$$

$$
=C_{j_k\ldots j_1}.
$$

\vspace{4mm}

The expansion (\ref{otit3333}) is proved. So, we can formulate the 
following theorem.

{\bf Theorem 2.14} \cite{9} (2013) (also see
\cite{10}-\cite{12aa}, \cite{arxiv-6}). 
{\it Suppose that the conditions of 
Theorem {\rm 2.10} are fulfilled.
Then
\begin{equation}
\label{aq1}
J^{*}[\psi^{(k)}]_{T,t}=
\sum_{j_k=0}^{\infty}\ldots\sum_{j_1=0}^{\infty}
C_{j_k\ldots j_1}
\prod_{l=1}^k
\zeta^{(i_l)}_{j_l},
\end{equation}

\vspace{2mm}
\noindent
where notations are the same as in Theorem {\rm 2.10}.}

Note that (\ref{aq1})  means the following
$$
\lim\limits_{p_k\to\infty}
\varlimsup\limits_{p_{k-1}\to\infty}\ldots\varlimsup\limits_{p_1\to\infty}
{\sf M}\left\{\left(J^{*}[\psi^{(k)}]_{T,t}-
\sum_{j_k=0}^{p_k}\ldots\sum_{j_1=0}^{p_1}
C_{j_k\ldots j_1}
\prod_{l=1}^k
\zeta^{(i_l)}_{j_l}\right)^{2n}\right\}=0, 
$$

\vspace{1mm}
\noindent
where
$n\in{\bf N}.$

Let us make a remark about how one can obtain
an analogue of Theorem~2.10 for the 
complete orthonormal system of Legendre 
polynomials in the space $L_2([t, T])$
and $n=1$ (the case of mean-square convergence), $k=2.$

From (\ref{newbegin1}) we have

\vspace{-3mm}
$$
{\sf M}\left\{\left(J[R_{p_1p_2}]_{T,t}^{(2)}\right)^{2}
\right\}\le
$$
\begin{equation}
\label{strange10}
~~~~~~~\le
2\int\limits_{[t, T]^2}
\left(R_{p_1p_2}(t_1,t_2)\right)^{2}dt_1 dt_2
+
{\bf 1}_{\{i_1=i_2\ne 0\}}
\left(\int\limits_t^T R_{p_1p_2}(t_1,t_1)dt_1\right)^2.
\end{equation}

\vspace{3mm}

From Remark~1.6 and (\ref{dobav100}), (\ref{strange11}) we obtain
for the case of Legendre polynomials

\vspace{-3mm}
$$
\lim\limits_{p_1\to\infty}\varlimsup\limits_{p_2\to\infty}
\int\limits_{[t, T]^2}
\left(R_{p_1p_2}(t_1,t_2)\right)^{2}dt_1 dt_2=0.
$$

\vspace{4mm}

Further, we have (see (\ref{d2020}))

\vspace{-3mm}
$$
R_{p_1p_2}(t_1,t_1)=\left(K^{*}(t_1,t_1)-\sum\limits_{j_1=0}^{p_1}
C_{j_1}(t_1)\phi_{j_1}(t_1)\right)+
$$
\begin{equation}
\label{strange14}
+\left(
\sum\limits_{j_1=0}^{p_1}\left(C_{j_1}(t_1)-
\sum\limits_{j_2=0}^{p_2}
C_{j_2j_1}\phi_{j_2}(t_1)\right)
\phi_{j_1}(t_1)\right).
\end{equation}

\vspace{3mm}

Then, taking into account (\ref{strange101}),
(\ref{strange14}) and applying two times (we mean here an iterated passage to the limit
$\lim\limits_{p_1\to\infty}\varlimsup\limits_{p_2\to\infty}$)
the Lebesgue's 
Dominated Convergence Theorem,
we obtain

\vspace{-1mm}
$$
\lim\limits_{p_1\to\infty}\varlimsup\limits_{p_2\to\infty}
\int\limits_t^T
R_{p_1p_2}(t_1,t_1)dt_1=0.
$$

\vspace{2mm}

Let us discuss the choice of integrable majorants
when applying Lebesgue's 
Dominated Convergence Theorem in our case.

Using double integration by parts (as in (\ref{otit2005})), we estimate
the partial sums of one-dimensional Fourier--Legendre series
$$
\sum\limits_{j_1=0}^{p_1}
C_{j_1}(t_1)\phi_{j_1}(t_1),\ \ \ 
\sum\limits_{j_2=0}^{p_2}
C_{j_2j_1}\phi_{j_2}(t_1)
$$ 

\noindent
in (\ref{strange14}) using (\ref{otit987}) and (\ref{strange3})  as follows
\begin{equation}
\label{strange17}
~~~\left|\sum\limits_{j_1=0}^{p_1}
C_{j_1}(t_1)\phi_{j_1}(t_1)\right|
\le 
K_1\Biggl(1+
\frac{1}{\left(1-(z(t_1))^2\right)^{1/2}}+
\frac{1}{\left(1-(z(t_1))^2\right)^{1/4}}\Biggr),
\end{equation}
\begin{equation}
\label{strange18}
\left|\sum\limits_{j_2=0}^{p_2}
C_{j_2j_1}\phi_{j_2}(t_1)
\right|
\le K_2\left(1+\frac{1}{\left(1-(z(t_1))^2\right)^{1/4}}\right),
\end{equation}

\vspace{2mm}
\noindent
where $z(t_1)$ is defined by (\ref{zz1}),
constant $K_1$ does not depend on $p_1,$ and constant $K_2$ does not depend on $p_2.$ 

Thus, integrable majorants
in our case can be easily constracted 
using (\ref{strange14}), (\ref{strange17}) and (\ref{strange18}) (see the proof of Theorem~2.10
for details).

An analogue of Theorem~2.10 for the 
case of Legendre 
polynomials 
and $n=1$ (the case of mean-square convergence), $k=2$ is obtained.

\subsection{Further Remarks}

In this section, we consider some approaches on the base of Theorems 2.10 and 1.1
for the case $k=2.$ Moreover,
we explain the potential 
difficulties associated with the use 
of generalized multiple Fourier series converging 
almost everywhere  (with respect to Lebesgue's measure (here and further in this section))
on the hypercube 
$[t, T]^k$ in the proof of Theorem 2.10.

First, we show how iterated series can be replaced by 
multiple one in Theorem 2.10 (the case $k=2$ and $n=1$) 
and in analogue of Theorem~2.10 for the 
case of Legendre 
polynomials 
(the case $k=2$ and $n=1$).

We have

\vspace{-2mm}
$$
\lim\limits_{p\to\infty}
{\sf M}\left\{\left(
J^{*}[\psi^{(2)}]_{T,t}-
\sum\limits_{j_1=0}^{p}\sum\limits_{j_2=0}^{p}
C_{j_2 j_1}\zeta_{j_1}^{(i_1)}
\zeta_{j_2}^{(i_2)}\right)^{2}\right\}=
$$
$$
=\lim\limits_{p\to\infty}\varlimsup\limits_{q\to\infty}
{\sf M}\left\{\left(
J^{*}[\psi^{(2)}]_{T,t}-
\sum\limits_{j_1=0}^{p}\sum\limits_{j_2=0}^{p}
C_{j_2 j_1}\zeta_{j_1}^{(i_1)}
\zeta_{j_2}^{(i_2)}\right)^{2}\right\}\le
$$
$$
\le
\lim\limits_{p\to\infty}\varlimsup\limits_{q\to\infty}\left(
2{\sf M}\left\{\left(
J^{*}[\psi^{(2)}]_{T,t}-
\sum\limits_{j_1=0}^{p}\sum\limits_{j_2=0}^{q}
C_{j_2 j_1}\zeta_{j_1}^{(i_1)}
\zeta_{j_2}^{(i_2)}\right)^{2}\right\}+\right.
$$
$$
\left.+
2{\sf M}\left\{\left(
\sum\limits_{j_1=0}^{p}\sum\limits_{j_2=0}^{q}
C_{j_2 j_1}\zeta_{j_1}^{(i_1)}
\zeta_{j_2}^{(i_2)}
-
\sum\limits_{j_1=0}^{p}\sum\limits_{j_2=0}^{p}
C_{j_2 j_1}\zeta_{j_1}^{(i_1)}
\zeta_{j_2}^{(i_2)}\right)^{2}\right\}\right)=
$$
$$
=2\lim\limits_{p\to\infty}\varlimsup\limits_{q\to\infty}
{\sf M}\left\{\left(
\sum\limits_{j_1=0}^{p}\sum\limits_{j_2=p+1}^{q}
C_{j_2 j_1}\zeta_{j_1}^{(i_1)}
\zeta_{j_2}^{(i_2)}\right)^{2}\right\}=
$$

\vspace{-2mm}
$$
=2 \lim\limits_{p\to\infty}\varlimsup\limits_{q\to\infty}
\sum\limits_{j_1=0}^{p}\sum\limits_{j_1'=0}^{p}
\sum\limits_{j_2=p+1}^{q}\sum\limits_{j_2'=p+1}^{q}
C_{j_2 j_1}C_{j_2' j_1'}{\sf M}\left\{\zeta_{j_1}^{(i_1)}
\zeta_{j_1'}^{(i_1)}\right\}
{\sf M}\left\{\zeta_{j_2}^{(i_2)}\zeta_{j_2'}^{(i_2)}\right\}=
$$

\vspace{-2mm}
$$
=2\lim\limits_{p\to\infty}\lim\limits_{q\to\infty}
\sum\limits_{j_1=0}^{p}\sum\limits_{j_2=p+1}^{q}
C_{j_2 j_1}^2=
$$
\begin{equation}
\label{fd1}
=2\lim\limits_{p\to\infty}\lim\limits_{q\to\infty}
\left(\sum\limits_{j_1=0}^{p}\sum\limits_{j_2=0}^{q}
C_{j_2 j_1}^2-
\sum\limits_{j_1=0}^{p}\sum\limits_{j_2=0}^{p}
C_{j_2 j_1}^2\right)
=
\end{equation}

\begin{equation}
\label{fd2}
=2\left(\lim\limits_{p,q\to\infty}
\sum\limits_{j_1=0}^{p}\sum\limits_{j_2=0}^{q}
C_{j_2 j_1}^2-
\lim\limits_{p\to\infty}
\sum\limits_{j_1=0}^{p}\sum\limits_{j_2=0}^{p}
C_{j_2 j_1}^2\right)
=
\end{equation}

\begin{equation}
\label{fd3}
~~~~~~~~=2\int\limits_{[t,T]^2}K^2 (t_1,t_2)dt_1dt_2-
2\int\limits_{[t,T]^2}K^2 (t_1,t_2)dt_1dt_2=0,
\end{equation}

\vspace{2mm}
\noindent 
where the function $K(t_1,t_2)$ is defined by (\ref{ppp}) for $k=2.$

Note that the transition from (\ref{fd1}) to (\ref{fd2})
is based on the theorem on reducing 
of a limit to iterated one.
Moreover, the transition from (\ref{fd2}) to (\ref{fd3})
is based on the Parseval equality.

Thus, we obtain the following Theorem.

{\bf Theorem 2.15}\ \cite{12a}-\cite{12aa}, \cite{arxiv-6}. {\it Assume that
$\{\phi_j(x)\}_{j=0}^{\infty}$ is a complete orthonormal
system of Legendre polynomials or trigonometric functions
in the space $L_2([t, T])$. 
At the same time $\psi_2(\tau)$ is a continuously dif\-ferentiable 
nonrandom function on $[t, T]$ and $\psi_1(\tau)$ is twice 
continuously differentiable nonrandom function on $[t, T]$. 
Then$,$ for the iterated Stratonovich stochastic integral {\rm (\ref{str})}
of multiplicity $2$
$$
J^{*}[\psi^{(2)}]_{T,t}=
{\int\limits_t^{*}}^T
\psi_2(t_2)
{\int\limits_t^{*}}^{t_2}
\psi_1(t_1) d{\bf w}_{t_1}^{(i_1)}
d{\bf w}_{t_2}^{(i_2)}\ \ \ 
(i_1,i_2=0, 1,\ldots,m)
$$
the following 
expansion 
$$
J^{*}[\psi^{(2)}]_{T,t}=
\hbox{\vtop{\offinterlineskip\halign{
\hfil#\hfil\cr
{\rm l.i.m.}\cr
$\stackrel{}{{}_{p\to \infty}}$\cr
}} }
\sum\limits_{j_1,j_2=0}^{p}
C_{j_2 j_1}\zeta_{j_1}^{(i_1)}
\zeta_{j_2}^{(i_2)}
$$
that converges in the mean-square sense
is valid, where 
\begin{equation}
\label{xyzyx}
C_{j_2 j_1}=\int\limits_t^T\psi_2(t_2)\phi_{j_2}(t_2)
\int\limits_t^{t_2}
\psi_1(t_1)\phi_{j_1}(t_1)
dt_1dt_2
\end{equation}
is the Fourier coefficient and
$$
\zeta_{j}^{(i)}=
\int\limits_t^T \phi_{j}(s) d{\bf w}_s^{(i)}
$$ 
are independent standard Gaussian random variables for various 
$i$ or $j$ {\rm (}in the case when $i\ne 0${\rm )}$,$
${\bf w}_{\tau}^{(i)}$ $(i=1,\ldots,m)$ are independent 
standard Wiener processes$,$
${\bf w}_{\tau}^{(0)}=\tau.$}

Note that Theorem 2.15 is a modification (for the 
case $p_1=p_2=p$ of series summation) of Theorem 2.1.

Using Theorem 2.10, we get
$$
0\le
\left\vert \lim\limits_{p_1\to\infty}
\varlimsup\limits_{p_2\to\infty}\ldots 
\varlimsup\limits_{p_k\to\infty}{\sf M}\left\{
\sum_{j_1=0}^{p_1}\ldots\sum_{j_k=0}^{p_k}
C_{j_k\ldots j_1}
\prod_{l=1}^k
\zeta^{(i_l)}_{j_l}-J^{*}[\psi^{(k)}]_{T,t}\right\}\right\vert \le
$$
$$
\le \lim\limits_{p_1\to\infty}
\varlimsup\limits_{p_2\to\infty}\ldots 
\varlimsup\limits_{p_k\to\infty} \left\vert {\sf M}\left\{
\sum_{j_1=0}^{p_1}\ldots\sum_{j_k=0}^{p_k}
C_{j_k\ldots j_1}
\prod_{l=1}^k
\zeta^{(i_l)}_{j_l}-J^{*}[\psi^{(k)}]_{T,t}\right\}\right\vert \le
$$
$$
\le \lim\limits_{p_1\to\infty}
\varlimsup\limits_{p_2\to\infty}\ldots 
\varlimsup\limits_{p_k\to\infty} {\sf M}\left\{\left\vert 
J^{*}[\psi^{(k)}]_{T,t}-\sum_{j_1=0}^{p_1}\ldots\sum_{j_k=0}^{p_k}
C_{j_k\ldots j_1}
\prod_{l=1}^k
\zeta^{(i_l)}_{j_l}\right\vert\right\} \le
$$
\begin{equation}
\label{het}
\le 
\lim\limits_{p_1\to\infty}
\varlimsup\limits_{p_2\to\infty}\ldots 
\varlimsup\limits_{p_k\to\infty}
\left(\hspace{-1mm}{\sf M}\hspace{-1mm}
\left\{\left(J^{*}[\psi^{(k)}]_{T,t}-
\sum_{j_1=0}^{p_1}\ldots\sum_{j_k=0}^{p_k}
C_{j_k\ldots j_1}
\prod_{l=1}^k
\zeta^{(i_l)}_{j_l}\right)^{2}\right\}\right)^{\hspace{-2mm}1/2}
\hspace{-1mm}=0.
\end{equation}

From the other hand,
$$
\lim\limits_{p_1\to\infty}
\varlimsup\limits_{p_2\to\infty}\ldots 
\varlimsup\limits_{p_k\to\infty} 
\left(
\sum_{j_1=0}^{p_1}\ldots\sum_{j_k=0}^{p_k}
C_{j_k\ldots j_1}
{\sf M}\left\{\prod_{l=1}^k
\zeta^{(i_l)}_{j_l}\right\}-
{\sf M}\left\{J^{*}[\psi^{(k)}]_{T,t}\right\}\right)=
$$
\begin{equation}
\label{het100}
=
\lim\limits_{p_1\to\infty}
\varlimsup\limits_{p_2\to\infty}\ldots 
\varlimsup\limits_{p_k\to\infty}
\sum_{j_1=0}^{p_1}\ldots\sum_{j_k=0}^{p_k}C_{j_k\ldots j_1}
{\sf M}\left\{\prod_{l=1}^k
\zeta^{(i_l)}_{j_l}\right\}
-{\sf M}\left\{
J^{*}[\psi^{(k)}]_{T,t}\right\}.
\end{equation}

\vspace{3mm}

Combining (\ref{het}) and (\ref{het100}), we obtain
\begin{equation}
\label{het1}
~~{\sf M}\left\{
J^{*}[\psi^{(k)}]_{T,t}\right\}=
\lim\limits_{p_1\to\infty}
\varlimsup\limits_{p_2\to\infty}\ldots 
\varlimsup\limits_{p_k\to\infty}
\sum_{j_1=0}^{p_1}\ldots\sum_{j_k=0}^{p_k}
C_{j_k\ldots j_1}
{\sf M}\left\{\prod_{l=1}^k
\zeta^{(i_l)}_{j_l}\right\}.
\end{equation}

\vspace{2mm}

Note that the relations (\ref{het})--(\ref{het1})
are also valid for the case of Legendre polynomials and
$k=2.$

The formula (\ref{het1}) with $k=2$ implies the
following
$$
{\sf M}\left\{
J^{*}[\psi^{(2)}]_{T,t}\right\}= \frac{1}{2}{\bf 1}_{\{i_1=i_2\ne 0\}}
\int\limits_t^T
\psi_1(s)\psi_2(s)ds=
$$
\begin{equation}
\label{het2}
=\lim\limits_{p_1\to\infty}
\varlimsup\limits_{p_2\to\infty} 
\sum_{j_1=0}^{p_1}\sum_{j_2=0}^{p_2}
C_{j_2j_1}
{\sf M}\left\{
\zeta^{(i_1)}_{j_1}\zeta^{(i_2)}_{j_2}\right\},
\end{equation}

\vspace{4mm}
\noindent
where ${\bf 1}_A$ is the indicator of the set $A.$

Since
$$
{\sf M}\left\{
\zeta^{(i_1)}_{j_1}\zeta^{(i_2)}_{j_2}\right\}=
{\bf 1}_{\{i_1=i_2\ne 0\}}{\bf 1}_{\{j_1=j_2\}},
$$

\vspace{1mm}
\noindent
then from (\ref{het2}) we obtain
$$
{\sf M}\left\{
J^{*}[\psi^{(2)}]_{T,t}\right\}=\lim\limits_{p_1\to\infty}
\varlimsup\limits_{p_2\to\infty} 
\sum_{j_1=0}^{p_1}\sum_{j_2=0}^{p_2}
C_{j_2j_1}{\bf 1}_{\{j_1=j_2\}}{\bf 1}_{\{i_1=i_2\ne 0\}}=
$$
\begin{equation}
\label{het5}
~~~~~~~~~ ={\bf 1}_{\{i_1=i_2\ne 0\}}\lim\limits_{p_1\to\infty}
\varlimsup\limits_{p_2\to\infty} 
\sum_{j_1=0}^{{\rm min}\{p_1, p_2\}}
C_{j_1j_1}=
{\bf 1}_{\{i_1=i_2\ne 0\}}\sum_{j_1=0}^{\infty}
C_{j_1j_1},
\end{equation}
where $C_{j_1j_1}$ is defined by (\ref{xyzyx}) for $j_1=j_2,$
i.e. 
$$
C_{j_1j_1}=\int\limits_t^T \psi_2(t_2)\phi_{j_1}(t_2)
\int\limits_t^{t_2} \psi_1(t_1)\phi_{j_1}(t_1)dt_1 dt_2.
$$

From (\ref{het2}) and (\ref{het5}) we obtain the following relation
\begin{equation}
\label{het7}
\sum_{j_1=0}^{\infty}
C_{j_1j_1}=
\frac{1}{2}
\int\limits_t^T
\psi_1(s)\psi_2(s)ds.
\end{equation}

\vspace{3mm}

Combining (\ref{a2}) (also see (\ref{razzar1}) for $k=2$) and (\ref{het7}), we have
$$
J[\psi^{(2)}]_{T,t}
=\hbox{\vtop{\offinterlineskip\halign{
\hfil#\hfil\cr
{\rm l.i.m.}\cr
$\stackrel{}{{}_{p_1,p_2\to \infty}}$\cr
}} }\sum_{j_1=0}^{p_1}\sum_{j_2=0}^{p_2}
C_{j_2j_1}\Biggl(\zeta_{j_1}^{(i_1)}\zeta_{j_2}^{(i_2)}
-{\bf 1}_{\{i_1=i_2\ne 0\}}
{\bf 1}_{\{j_1=j_2\}}\Biggr)=
$$

\vspace{-3mm}
$$
=
\hbox{\vtop{\offinterlineskip\halign{
\hfil#\hfil\cr
{\rm l.i.m.}\cr
$\stackrel{}{{}_{p_1,p_2\to \infty}}$\cr
}} }\sum_{j_1=0}^{p_1}\sum_{j_2=0}^{p_2}
C_{j_2j_1}\zeta_{j_1}^{(i_1)}\zeta_{j_2}^{(i_2)}
-
{\bf 1}_{\{i_1=i_2\ne 0\}}\sum_{j_1=0}^{\infty}
C_{j_1j_1} =
$$

\vspace{-3mm}
\begin{equation}
\label{het10}
~~~~~~~~ =
\hbox{\vtop{\offinterlineskip\halign{
\hfil#\hfil\cr
{\rm l.i.m.}\cr
$\stackrel{}{{}_{p_1,p_2\to \infty}}$\cr
}} }\sum_{j_1=0}^{p_1}\sum_{j_2=0}^{p_2}
C_{j_2j_1}\zeta_{j_1}^{(i_1)}\zeta_{j_2}^{(i_2)}
-
\frac{1}{2}{\bf 1}_{\{i_1=i_2\ne 0\}}
\int\limits_t^T
\psi_1(s)\psi_2(s)ds.
\end{equation}

\vspace{3mm}

Since 
\begin{equation}
\label{uyes1}
~~~~~~~ J^{*}[\psi^{(2)}]_{T,t}=J[\psi^{(2)}]_{T,t}+
\frac{1}{2}{\bf 1}_{\{i_1=i_2\ne 0\}}
\int\limits_t^T
\psi_1(s)\psi_2(s)ds\ \ \ \hbox{w.~p.~1},
\end{equation}

\noindent
then from (\ref{het10}) we finally get
the following expansion
$$
J^{*}[\psi^{(2)}]_{T,t}
=\hbox{\vtop{\offinterlineskip\halign{
\hfil#\hfil\cr
{\rm l.i.m.}\cr
$\stackrel{}{{}_{p_1,p_2\to \infty}}$\cr
}} }\sum_{j_1=0}^{p_1}\sum_{j_2=0}^{p_2}
C_{j_2j_1}\zeta_{j_1}^{(i_1)}\zeta_{j_2}^{(i_2)}.
$$

\noindent
Thus, we obtain the statement of Theorem 2.1.

We have 
$$
J^{*}[\psi^{(2)}]_{T,t}^{p_1,p_2}\stackrel{\sf def}{=}J[\psi^{(2)}]_{T,t}^{p_1,p_2}+
\frac{1}{2}{\bf 1}_{\{i_1=i_2\ne 0\}}
\int\limits_t^T
\psi_1(s)\psi_2(s)ds=
$$
$$
=\sum_{j_1=0}^{p_1}\sum_{j_2=0}^{p_2}
C_{j_2j_1}\Biggl(\zeta_{j_1}^{(i_1)}\zeta_{j_2}^{(i_2)}
-{\bf 1}_{\{i_1=i_2\ne 0\}}
{\bf 1}_{\{j_1=j_2\}}\Biggr)+\frac{1}{2}{\bf 1}_{\{i_1=i_2\ne 0\}}
\int\limits_t^T
\psi_1(s)\psi_2(s)ds
=
$$
\begin{equation}
\label{ziko432}
=
\sum_{j_1=0}^{p_1}\sum_{j_2=0}^{p_2}
C_{j_2j_1}\zeta_{j_1}^{(i_1)}\zeta_{j_2}^{(i_2)}+{\bf 1}_{\{i_1=i_2\ne 0\}}
\left(\frac{1}{2}\int\limits_t^T
\psi_1(s)\psi_2(s)ds-\sum_{j_1=0}^{{\rm min}\{p_1,p_2\}}
C_{j_1j_1}\right),
\end{equation}

\noindent
where 
$$
J[\psi^{(2)}]_{T,t}^{p_1,p_2}=
\sum_{j_1=0}^{p_1}\sum_{j_2=0}^{p_2}
C_{j_2j_1}\Biggl(\zeta_{j_1}^{(i_1)}\zeta_{j_2}^{(i_2)}
-{\bf 1}_{\{i_1=i_2\ne 0\}}
{\bf 1}_{\{j_1=j_2\}}\Biggr)
$$

\noindent
is the approximation of iterated It\^{o} stochastic integral 
(\ref{itoxxx}) $(k=2)$ based on Theorem~1.1 (see (\ref{a2})).

Moreover, from (\ref{ziko11000}), (\ref{2026ch1001s11}), and (\ref{oop51}) we obtain
\begin{equation}
\label{ziko9991}
{\sf M}\left\{\left(J^{*}[\psi^{(2)}]_{T,t}-
J^{*}[\psi^{(2)}]_{T,t}^{p_1,p_2}\right)^{2n}\right\}=
{\sf M}\left\{\left(J[\psi^{(2)}]_{T,t}-
J[\psi^{(2)}]_{T,t}^{p_1,p_2}\right)^{2n}\right\}\ \to\ 0
\end{equation}

\noindent
if $p_1,p_2\to\infty$ $(n\in{\bf N})$.

Further, 
$$
{\sf M}\left\{\left(J^{*}[\psi^{(2)}]_{T,t}-
\sum\limits_{j_1=0}^{p_1}\sum\limits_{j_2=0}^{p_2}
C_{j_2 j_1}\zeta_{j_1}^{(i_1)}
\zeta_{j_2}^{(i_2)}
\right)^{2n}\right\}=
$$
$$
={\sf M}\Biggl\{\Biggl(J^{*}[\psi^{(2)}]_{T,t}-J^{*}[\psi^{(2)}]_{T,t}^{p_1,p_2}
+\Biggr.\Biggr.
$$
$$
\left.\left.+{\bf 1}_{\{i_1=i_2\ne 0\}}
\left(\frac{1}{2}\int\limits_t^T
\psi_1(s)\psi_2(s)ds-\sum_{j_1=0}^{{\rm min}\{p_1,p_2\}}
C_{j_1j_1}\right)\right)^{2n}\right\}\le
$$
$$
\le K_n \left( 
{\sf M}\left\{\Biggl(J^{*}[\psi^{(2)}]_{T,t}-
J^{*}[\psi^{(2)}]_{T,t}^{p_1,p_2}\Biggr)^{2n}\right\}+\right.
$$
\begin{equation}
\label{ziko9999}
~~~~~~~~+\left.
{\bf 1}_{\{i_1=i_2\ne 0\}}
\left(\frac{1}{2}\int\limits_t^T
\psi_1(s)\psi_2(s)ds-\sum_{j_1=0}^{{\rm min}\{p_1,p_2\}}
C_{j_1j_1}\right)^{2n}\right),
\end{equation}
where constant $K_n<\infty$ depends on 
$n$.

Taking into account (\ref{ziko9991}), (\ref{ziko9999}) and also that the equality
(\ref{het7}) is true under the conditions of Theorem~2.3 (see Sect.~2.1.4), we get
\begin{equation}
\label{ziko3210}
~~~~~~~~~~\lim\limits_{p_1,p_2\to\infty}
{\sf M}\left\{\left(J^{*}[\psi^{(2)}]_{T,t}-
\sum\limits_{j_1=0}^{p_1}\sum\limits_{j_2=0}^{p_2}
C_{j_2 j_1}\zeta_{j_1}^{(i_1)}
\zeta_{j_2}^{(i_2)}
\right)^{2n}\right\}=0.
\end{equation}

\noindent
Thus, we obtain the following theorem.

{\bf Theorem 2.16.}\ {\it Suppose that 
$\{\phi_j(x)\}_{j=0}^{\infty}$ is an arbitrary complete orthonormal system of 
functions in the space $L_2([t, T]).$
Moreover$,$ $\psi_1(\tau), \psi_2(\tau)$ are continuous 
functions on $[t, T].$ 
Then$,$ for the iterated Stratonovich stochastic integral {\rm (\ref{str})}
of multiplicity $2$
$$
J^{*}[\psi^{(2)}]_{T,t}=
{\int\limits_t^{*}}^T
\psi_2(t_2)
{\int\limits_t^{*}}^{t_2}
\psi_1(t_1) d{\bf w}_{t_1}^{(i_1)}
d{\bf w}_{t_2}^{(i_2)}\ \ \
(i_1,i_2=0, 1,\ldots,m)
$$
the following 
expansion 
$$
J^{*}[\psi^{(2)}]_{T,t}=
\sum\limits_{j_1,j_2=0}^{\infty}
C_{j_2 j_1}\zeta_{j_1}^{(i_1)}
\zeta_{j_2}^{(i_2)}
$$

\noindent
that converges in the mean of degree $2n,$ $n\in{\bf N}$ {\rm (see (\ref{ziko3210}))}
is valid, where the Fourier coefficient 
$C_{j_2 j_1}$ is defined by {\rm (\ref{xyzyx})}
and
$$
\zeta_{j}^{(i)}=
\int\limits_t^T \phi_{j}(s) d{\bf w}_s^{(i)}
$$ 
are independent standard Gaussian random variables for various 
$i$ or $j$ {\rm (}in the case when $i\ne 0${\rm );}
another notations are the same as in Theorem~{\rm 2.15.}}

Let us consider some other approaches
close to the approaches outlined in this section.

Now we turn to multiple trigonometric Fourier series converging 
almost everywhere. Let us formulate the well known result
from the theory of multiple
trigonometric Fourier series.

{\bf Proposition 2.3} \cite{dudu}. {\it Suppose that
$$
\int\limits_{[0,2\pi]^k}
|f(x_1,\ldots,x_k)|\left({\rm log}^{+}|f(x_1,\ldots,x_k)|\right)^k
{\rm log}^{+}{\rm log}^{+}|f(x_1,\ldots,x_k)|\times
$$
\begin{equation}
\label{vot}
\times dx_1\ldots dx_k<\infty.
\end{equation}

Then, for the square partial sums
$$
\sum_{j_1=0}^p\ldots \sum_{j_k=0}^p
C_{j_k\ldots j_1}\prod\limits_{l=1}^k\phi_{j_l}(x_l)
$$
of the multiple trigonometric Fourier series we have
$$
\lim\limits_{p\to\infty}
\sum_{j_1=0}^p\ldots \sum_{j_k=0}^p
C_{j_k\ldots j_1}\prod\limits_{l=1}^k\phi_{j_l}(x_l)=
f(x_1,\ldots,x_k)
$$
almost everywhere on $[0, 2\pi]^k,$ where 
$\{\phi_j(x)\}_{j=0}^{\infty}$ is a complete orthonormal
system of trigonometric functions
in the space $L_2([0, 2\pi]),$ ${\rm log}^{+}x={\rm log} \max\{1,\ x\},$
$$
C_{j_k\ldots j_1}=
\int\limits_{[0,2\pi]^k}
f(x_1,\ldots,x_k)\prod\limits_{l=1}^k
\phi_{j_l}(x_l)
dx_1\ldots dx_k
$$
is the Fourier coefficient of the function $f(x_1,\ldots,x_k).$}

Note that Proposition~2.3 can be reformulated for 
$[t, T]^k$ instead of 
$[0, 2\pi]^k.$
If we tried to apply Proposition~2.3 in the proof of Theorem 2.10, 
then we would encounter the following difficulties.
The right-hand side of (\ref{udar1}) contains
multiple integrals 
over hypercubes of various dimensions, namely over hypercubes
$[t, T]^k$, $[t, T]^{k-1},$ etc.
Obviously, the convergence almost everywhere on 
$[t, T]^k$ does not mean the convergence almost everywhere on 
$[t, T]^{k-1}$, $[t, T]^{k-2},$ etc.
This means that we could not apply the Lebesgue's Dominated 
Convergence Theorem 
in the proof of Theorem 2.13 
and thus we could not  complete the proof of Theorem 2.10.
Although multiple series are more convenient in terms of approximation 
than iterated series as in Theorem 2.10.

Suppose that the conditions of Theorem~2.10
are fulfilled. 
In the proof of Theorem~2.2 (see (\ref{5656})) we deduced for the trigonometric case that 
\begin{equation}
\label{za1}
~~~~~~\lim\limits_{p\to\infty}
\sum_{j_1=0}^{p} \sum_{j_2=0}^{p} C_{j_2j_1}\phi_{j_1}(t_1)
\phi_{j_2}(t_1)=\frac{1}{2}\psi_1(t_1)\psi_2(t_1)=K^{*}(t_1,t_1),
\end{equation}
where $t_1\in (t, T),$ $C_{j_2j_1}$ is defined by (\ref{xyzyx}).
This means that we can repeat the proof of Theorem 2.10 for 
the case $k=2$ and apply the Lebesgue's Dominated Convergence 
Theorem 
in the formula (\ref{udar1}), since Proposition~2.3 and (\ref{za1})
imply the convergence almost everywhere on $[t, T]^2$ and $[t, T]$ 
($t_1=t_2\in [t, T]$)
of the multiple trigonometric 
Fourier series 
\begin{equation}
\label{ziko456}
\lim\limits_{p\to\infty}
\sum_{j_1=0}^{p} \sum_{j_2=0}^{p} C_{j_2j_1}\phi_{j_1}(t_1)
\phi_{j_2}(t_2),\ \ \ t_1,t_2\in [t, T]^2
\end{equation}

\noindent
to the function 
$K^{*}(t_1,t_2)$ (the question of finding an integrable majorant
for Lebesgue's 
Dominated Convergence Theorem is omitted here).
So, we could obtain the particular case of Theorem~2.16.

Consider another possible way, which is based
on the function (\ref{ziko5001}) and Theorem 2.10.
The case $i_1\ne i_2$ follows from (\ref{ziko432}) and (\ref{ziko9991}).
Consider the case $i_1=i_2\ne 0.$ We have 
$K^{*}(t_1,t_2)+K^{*}(t_2,t_1)=
K'(t_1,t_2),$
where the functions $K'(t_1,t_2)$ and $K^{*}(t_1,t_2)$ are defined by (\ref{ziko5001})
and (\ref{1999.1}) 
correspondingly. Note that the function 
$K'(t_1,t_2)$ is symmetric, i.e. 
$K'(t_1,t_2)=K'(t_2,t_1).$

\vspace{1mm}

By analogy with (\ref{nov800}) we get w.~p.~1
$$
J[K'/2]_{T,t}^{(2)}=
\frac{1}{2}\hbox{\vtop{\offinterlineskip\halign{
\hfil#\hfil\cr
{\rm l.i.m.}\cr
$\stackrel{}{{}_{N\to \infty}}$\cr
}} }\sum_{l_2=0}^{N-1}
\sum_{l_1=0}^{N-1}
K'(\tau_{l_1},\tau_{l_2})
\Delta{\bf w}_{\tau_{l_1}}^{(i_1)}
\Delta{\bf w}_{\tau_{l_2}}^{(i_1)}=
$$
$$
=\frac{1}{2}\hbox{\vtop{\offinterlineskip\halign{
\hfil#\hfil\cr
{\rm l.i.m.}\cr
$\stackrel{}{{}_{N\to \infty}}$\cr
}} }\left(\sum_{l_2=0}^{N-1}
\sum_{l_1=0}^{l_2-1}+\sum_{l_1=0}^{N-1}
\sum_{l_2=0}^{l_1-1}
\right)
K'(\tau_{l_1},\tau_{l_2})
\Delta{\bf w}_{\tau_{l_1}}^{(i_1)}
\Delta{\bf w}_{\tau_{l_2}}^{(i_1)}+
$$
$$
+\frac{1}{2}
\hbox{\vtop{\offinterlineskip\halign{
\hfil#\hfil\cr
{\rm l.i.m.}\cr
$\stackrel{}{{}_{N\to \infty}}$\cr
}} }\sum_{l_1=0}^{N-1}
K'(\tau_{l_1},\tau_{l_1})
\left(\Delta{\bf w}_{\tau_{l_1}}^{(i_1)}\right)^2=
$$
$$
=\frac{1}{2}\hbox{\vtop{\offinterlineskip\halign{
\hfil#\hfil\cr
{\rm l.i.m.}\cr
$\stackrel{}{{}_{N\to \infty}}$\cr
}} }\sum_{l_2=0}^{N-1}
\sum_{l_1=0}^{l_2-1}
\left(K'(\tau_{l_1},\tau_{l_2})+K'(\tau_{l_2},\tau_{l_1})\right)
\Delta{\bf w}_{\tau_{l_1}}^{(i_1)}
\Delta{\bf w}_{\tau_{l_2}}^{(i_1)}+
$$
$$
+\frac{1}{2}
\hbox{\vtop{\offinterlineskip\halign{
\hfil#\hfil\cr
{\rm l.i.m.}\cr
$\stackrel{}{{}_{N\to \infty}}$\cr
}} }\sum_{l_1=0}^{N-1}
K'(\tau_{l_1},\tau_{l_1})
\left(\Delta{\bf w}_{\tau_{l_1}}^{(i_1)}\right)^2=
$$
$$
=\hbox{\vtop{\offinterlineskip\halign{
\hfil#\hfil\cr
{\rm l.i.m.}\cr
$\stackrel{}{{}_{N\to \infty}}$\cr
}} }\sum_{l_2=0}^{N-1}
\sum_{l_1=0}^{l_2-1}
K'(\tau_{l_1},\tau_{l_2})
\Delta{\bf w}_{\tau_{l_1}}^{(i_1)}
\Delta{\bf w}_{\tau_{l_2}}^{(i_1)}
+\frac{1}{2}
\hbox{\vtop{\offinterlineskip\halign{
\hfil#\hfil\cr
{\rm l.i.m.}\cr
$\stackrel{}{{}_{N\to \infty}}$\cr
}} }\sum_{l_1=0}^{N-1}
K'(\tau_{l_1},\tau_{l_1})
\left(\Delta{\bf w}_{\tau_{l_1}}^{(i_1)}\right)^2=
$$
$$
=\int\limits_t^T \psi_2(t_2)\int\limits_t^{t_2}\psi_1(t_1)
d{\bf w}_{t_1}^{(i_1)}d{\bf w}_{t_2}^{(i_1)}
+\frac{1}{2}
\int\limits_t^T \psi_1(t_1)\psi_2(t_1)dt_1= 
$$
\begin{equation}
\label{ziko6002}
={\int\limits_t^{*}}^T\psi_2(t_2)
{\int\limits_t^{*}}^{t_2}\psi_1(t_1)d{\bf w}_{t_1}^{(i_1)}
d{\bf w}_{t_2}^{(i_1)}\stackrel{\sf def}{=}
J^{*}[\psi^{(2)}]_{T,t},
\end{equation}

\noindent
where we used the same notations as in (\ref{nov800}).

Let us expand the function $K'(t_1,t_2)/2$ into a multiple 
(double) Fourier--Legendre series or trigonometric Fourier series 
in the square $[t, T]^2$ (see (\ref{334.ye}))
$$
\frac{1}{2}K'(t_1,t_2)=
$$
$$
=
\frac{1}{2}\lim_{p_1,p_2\to\infty}
\sum_{j_1=0}^{p_1}\sum_{j_2=0}^{p_2}
\int\limits_t^T\int\limits_t^T K'(t_1,t_2)
\phi_{j_1}(t_1)\phi_{j_2}(t_2)dt_1 dt_2\cdot
\phi_{j_1}(t_1)\phi_{j_2}(t_2)=
$$
$$
=\frac{1}{2}\lim_{p_1,p_2\to\infty}
\sum_{j_1=0}^{p_1}\sum_{j_2=0}^{p_2}\left(
\int\limits_t^T\psi_2(t_2)\phi_{j_2}(t_2)
\int\limits_t^{t_2}\psi_1(t_1)\phi_{j_1}(t_1)dt_1dt_2+\right.
$$
$$
\left. +
\int\limits_t^T\psi_1(t_2)\phi_{j_2}(t_2)
\int\limits_{t_2}^{T}\psi_2(t_1)\phi_{j_1}(t_1)dt_1\right)dt_2
\phi_{j_1}(t_1)\phi_{j_2}(t_2)
=
$$
\begin{equation}
\label{ziko6000}
=\frac{1}{2}\lim_{p_1,p_2\to\infty}
\sum_{j_1=0}^{p_1}\sum_{j_2=0}^{p_2}\left(C_{j_2j_1}+
C_{j_1j_2}\right)
\phi_{j_1}(t_1)\phi_{j_2}(t_2),
\end{equation}

\vspace{2mm}
\noindent
where the series (\ref{ziko6000}) converges to 
$K'(t_1,t_2)/2$ at any inner point of the square $[t, T]^2$
(see the proof of Theorem 2.2 for details).

In obtaining (\ref{ziko6000}) we replaced the order of
integration in the second iterated integral.

Using (\ref{ziko6002}), (\ref{ziko6000}), and the 
scheme of the proof of Theorem~2.10 $(k=2)$, we can obtain
the following relation
(the question of finding an integrable majorant
for Lebesgue's 
Dominated Convergence Theorem is omitted here)
\begin{equation}
\label{ziko7000}
\lim\limits_{p_1,p_2\to \infty}
{\sf M}\left\{\left(J^{*}[\psi^{(2)}]_{T,t}-\frac{1}{2}\sum_{j_1=0}^{p_1}\sum_{j_2=0}^{p_2}
\left(C_{j_2j_1}+C_{j_1j_2}\right)\zeta_{j_1}^{(i_1)}\zeta_{j_2}^{(i_1)}\right)^{2n}\right\}=0.
\end{equation}

Let us rewrite the sum on 
the left-hand side of (\ref{ziko7000}) as two  sums.
Let us replace $j_1$ with $j_2$, $j_2$ with $j_1$, 
$p_1$ with $p_2,$ and $p_2$ with $p_1$ in the second sum. 

Thus, we get 
$$
\lim\limits_{p_1,p_2\to \infty}
{\sf M}\left\{\left(J^{*}[\psi^{(2)}]_{T,t}-\sum_{j_1=0}^{p_1}\sum_{j_2=0}^{p_2}
C_{j_2j_1}\zeta_{j_1}^{(i_1)}\zeta_{j_2}^{(i_1)}\right)^{2n}\right\}=0.
$$

\vspace{1mm}

Let us consider another approach.
The following fact is well known \cite{IP}.

{\bf Proposition 2.4.} {\it Let 
$\bigl\{x_{n_1,\ldots,n_k}\bigr\}_{n_1,\ldots,n_k=1}^{\infty}$
be a multi-index sequence and let there exists the limit
$$
\lim\limits_{n_1,\ldots,n_k\to\infty}x_{n_1,\ldots,n_k}<\infty.
$$

\noindent
Moreover, let there exists the limit
$$
\lim\limits_{n_k\to\infty}x_{n_1,\ldots,n_k}=y_{n_1,\ldots,n_{k-1}}
<\infty\ \ \ \hbox{for any}\ \ \ n_1,\ldots,n_{k-1}.
$$

Then there exists the iterated limit
$$
\lim\limits_{n_1,\ldots,n_{k-1}\to\infty}\
\lim\limits_{n_k\to\infty}x_{n_1,\ldots,n_k}
$$
and moreover,
$$
\lim\limits_{n_1,\ldots,n_{k-1}\to\infty}\
\lim\limits_{n_k\to\infty}
x_{n_1,\ldots,n_k}=
\lim\limits_{n_1,\ldots,n_k\to\infty}x_{n_1,\ldots,n_k}.
$$
}

Denote
$$
C_{j_s\ldots j_1}(t_{s+1},\ldots,t_k)=
\int\limits_{[t,T]^s}K(t_1,\ldots,t_k)
\prod_{l=1}^s \phi_{j_l}(t_l)dt_1\ldots dt_s,
$$
where $s=1,\ldots,k-1$ and
$K(t_1,\ldots,t_k)$ is defined by (\ref{ppp}).
For $s=k$ we suppose that $C_{j_k\ldots j_1}$
is defined by
(\ref{ppppa}).

Consider the following Fourier series
\begin{equation}
\label{ww1}
\lim\limits_{p_1,p_2\to\infty}
\sum_{j_1=0}^{p_1} \sum_{j_2=0}^{p_2} C_{j_2j_1}(t_3,\ldots,t_k)
\phi_{j_1}(t_1)
\phi_{j_2}(t_2),
\end{equation}
\begin{equation}
\label{ww2}
~~~~~~~~\lim\limits_{p_1,p_2,p_3\to\infty}
\sum_{j_1=0}^{p_1}\sum_{j_2=0}^{p_2}
\sum_{j_3=0}^{p_3} C_{j_3j_2j_1}(t_4,\ldots,t_k)
\phi_{j_1}(t_1)
\phi_{j_2}(t_2)\phi_{j_3}(t_3),
\end{equation}
$$
~~~\ldots
$$
\begin{equation}
\label{ww3}
~~~~~\lim\limits_{p_1,\ldots,p_{k-1}\to\infty}
\sum_{j_1=0}^{p_1}\ldots
\sum_{j_{k-1}=0}^{p_{k-1}} C_{j_{k-1}\ldots j_1}(t_k)
\phi_{j_1}(t_1)\ldots
\phi_{j_{k-1}}(t_{k-1}),
\end{equation}
\begin{equation}
\label{ww4}
\lim\limits_{p_1,\ldots,p_{k}\to\infty}
\sum_{j_1=0}^{p_1}\ldots
\sum_{j_{k}=0}^{p_{k}} C_{j_{k}\ldots j_1}
\phi_{j_1}(t_1)\ldots
\phi_{j_{k}}(t_{k}),
\end{equation}

\vspace{2.5mm}
\noindent
where $t_1,\ldots,t_k\in[t,T],$ 
$\{\phi_j(x)\}_{j=0}^{\infty}$ is a complete orthonormal
system of Legendre polynomials or trigonometric functions
in the space $L_2([t, T])$. 

The author does not know the answer to the question
on the existence of limits (\ref{ww1})--(\ref{ww4}) even for
the case $p_1=\ldots=p_k$ and trigonometric Fourier series.
Obviously, at least 
for the case $k=2$ and $\psi_1(\tau),$ $\psi_2(\tau)\equiv 1$
the answere to the above question is positive
for the Fourier--Legendre series as well as 
for the trigonometric Fourier series.

If we suppose that the 
limits (\ref{ww1})--(\ref{ww4}) exist, then
combining Proposition 2.4 and the proof of Theorem 2.11, 
we obtain

\vspace{-1.5mm}
$$
K^{*}(t_1,\ldots,t_k)=
\sum_{j_1=0}^{\infty}C_{j_1}(t_2,\ldots,t_k)
\phi_{j_1}(t_1)=
$$

\vspace{-1mm}
\begin{equation}
\label{ww10}
=\sum_{j_1=0}^{\infty}\sum_{j_2=0}^{\infty}C_{j_2j_1}(t_3,\ldots,t_k)
\phi_{j_1}(t_1)\phi_{j_2}(t_2)=
\end{equation}

\vspace{-1mm}
$$
=
\lim\limits_{p_1,p_2\to\infty}
\sum_{j_1=0}^{p_1} \sum_{j_2=0}^{p_2} C_{j_2j_1}(t_3,\ldots,t_k)
\phi_{j_1}(t_1)
\phi_{j_2}(t_2)=
$$

\vspace{-1mm}
$$
=
\lim\limits_{p_1,p_2\to\infty}
\sum_{j_1=0}^{p_1} \sum_{j_2=0}^{p_2} 
\sum_{j_3=0}^{\infty}C_{j_3j_2j_1}(t_4,\ldots,t_k)
\phi_{j_1}(t_1)
\phi_{j_2}(t_2)\phi_{j_3}(t_3)=
$$

\vspace{-1mm}
\begin{equation}
\label{ww18}
~~~~~~ =
\lim\limits_{p_1,p_2,p_3\to\infty}
\sum_{j_1=0}^{p_1} \sum_{j_2=0}^{p_2} 
\sum_{j_3=0}^{p_3}C_{j_3j_2j_1}(t_4,\ldots,t_k)
\phi_{j_1}(t_1)
\phi_{j_2}(t_2)\phi_{j_3}(t_3)=
\end{equation}

\vspace{-2mm}
\begin{equation}
\label{ww181}
~~~~~~ =\sum_{j_1=0}^{\infty}\sum_{j_2=0}^{\infty}\sum_{j_3=0}^{\infty}
C_{j_3j_2j_1}(t_4,\ldots,t_k)
\phi_{j_1}(t_1)\phi_{j_2}(t_2)\phi_{j_3}(t_3)=
\end{equation}

\vspace{-2mm}
\begin{equation}
\label{ww182}
~~~~~~ =
\lim\limits_{p_1,p_2,p_3\to\infty}
\sum_{j_1=0}^{p_1} \sum_{j_2=0}^{p_2}\sum_{j_3=0}^{p_3} 
\sum_{j_4=0}^{\infty}C_{j_4\ldots j_1}(t_5,\ldots,t_k)
\phi_{j_1}(t_1)
\phi_{j_3}(t_4)= 
\end{equation}

\vspace{-2mm}
$$
~~~~~~=\ldots =
$$

\newpage
\noindent
\begin{equation}
\label{ww21}
=
\lim\limits_{p_1,\ldots,p_k\to\infty}
\sum_{j_1=0}^{p_1} \ldots
\sum_{j_k=0}^{p_k}C_{j_k\ldots j_1}
\phi_{j_1}(t_1)\ldots
\phi_{j_k}(t_k).
\end{equation}

Note that 
the transition from (\ref{ww18}) to (\ref{ww181})
is based on (\ref{ww10}) and the proof of Theorem 2.11.
The transition from (\ref{ww181}) to (\ref{ww182})
is based on (\ref{ww18}) and the proof of Theorem 2.11.

Using (\ref{ww21}), we could get the version of Theorem 2.10 
with multiple series instead of iterated ones
(see Hypothesis~2.3, Sect.~2.5).

\subsection{Refinement of Theorems 2.10 and 2.14 for Iterated 
Stra\-to\-no\-vich Stochastic Integrals of Multiplicities $2$ and $3$
$(i_1,i_2,i_3=1,\ldots,m).$
The Case of Mean-Square Convergence}

In this section, it will be shown that the upper limits
in Theorems 2.10 and 2.14 (the cases $k=2,$ $k=3$ and $n=1$)
can be replaced by the usual limits.

{\bf Theorem 2.17}\ \cite{arxiv-6}.
{\it Suppose that every $\psi_l(\tau)$ $(l=1,2,3)$ is twice continuously
differentiable function at the interval
$[t, T]$ and
$\{\phi_j(x)\}_{j=0}^{\infty}$ is a complete
orthonormal system of 
trigonometric functions in the space $L_2([t, T])$. 
Then$,$ the iterated Stratonovich stochastic integrals 
$J^{*}[\psi^{(2)}]_{T,t}$ and $J^{*}[\psi^{(3)}]_{T,t}$ $(i_1,i_2,i_3=1,\ldots,m)$
defined by {\rm(\ref{str})}
are expanded into the 
conver\-ging 
in the mean-square sense 
iterated series
\begin{equation}
\label{nov500}
~~~~~~~\lim\limits_{p_1\to\infty}
\lim\limits_{p_2\to\infty}
{\sf M}\left\{\left(J^{*}[\psi^{(2)}]_{T,t}-
\sum_{j_1=0}^{p_1}\sum_{j_2=0}^{p_2}
C_{j_2j_1}\zeta^{(i_1)}_{j_1}\zeta^{(i_2)}_{j_2}\right)^{2}\right\}=0,
\end{equation}
\begin{equation}
\label{nov501}
~~~~~~~\lim\limits_{p_2\to\infty}
\lim\limits_{p_1\to\infty}
{\sf M}\left\{\left(J^{*}[\psi^{(2)}]_{T,t}-
\sum_{j_2=0}^{p_2}\sum_{j_1=0}^{p_1}
C_{j_2j_1}\zeta^{(i_1)}_{j_1}\zeta^{(i_2)}_{j_2}\right)^{2}\right\}=0,
\end{equation}
\begin{equation}
\label{nov502}
\lim\limits_{p_1\to\infty}
\lim\limits_{p_2\to\infty}
\lim\limits_{p_3\to\infty}
{\sf M}\left\{\left(J^{*}[\psi^{(3)}]_{T,t}-
\sum_{j_1=0}^{p_1}\sum_{j_2=0}^{p_2}\sum_{j_3=0}^{p_3}
C_{j_3j_2j_1}\zeta^{(i_1)}_{j_1}\zeta^{(i_2)}_{j_2}\zeta^{(i_3)}_{j_3}\right)^{2}\right\}=0,
\end{equation}
\begin{equation}
\label{nov503}
\lim\limits_{p_3\to\infty}
\lim\limits_{p_2\to\infty}
\lim\limits_{p_1\to\infty}
{\sf M}\left\{\left(J^{*}[\psi^{(3)}]_{T,t}-
\sum_{j_3=0}^{p_3}\sum_{j_2=0}^{p_2}\sum_{j_1=0}^{p_1}
C_{j_3j_2j_1}\zeta^{(i_1)}_{j_1}\zeta^{(i_2)}_{j_2}\zeta^{(i_3)}_{j_3}\right)^{2}\right\}=0,
\end{equation}

\noindent
where 
$$
\zeta_{j}^{(i)}=
\int\limits_t^T \phi_{j}(s) d{\bf w}_s^{(i)}\ \ \ (i=1,\ldots,m,\ \ j=0, 1,\ldots )
$$ 
are independent standard Gaussian random variables
for
various
$i$ or $j$ and 
$C_{j_2j_1},$ $C_{j_3j_2j_1}$ are defined by 
{\rm (\ref{333.40})} and {\rm (\ref{pppxx})}.}

{\bf Proof.}\ We will prove the equalities (\ref{nov500}) and (\ref{nov502})
(the equalities (\ref{nov501}) and (\ref{nov503}) can be proved similarly
using the expansion (\ref{otit3333}) instead of the expansion
(\ref{30.18})).

From (\ref{nov800}) we have
w.~p.~1 
$$
J^{*}[\psi^{(2)}]_{T,t}-
\sum_{j_1=0}^{p_1}\sum_{j_2=0}^{p_2}
C_{j_2j_1}\zeta^{(i_1)}_{j_1}\zeta^{(i_2)}_{j_2}=J[R_{p_1p_2}]_{T,t}^{(2)}=
$$
$$
=\int\limits_t^T\int\limits_t^{t_2}
R_{p_1p_2}(t_1,t_2)d{\bf w}_{t_1}^{(i_1)}d{\bf w}_{t_2}^{(i_2)}
+\int\limits_t^T\int\limits_t^{t_1}
R_{p_1p_2}(t_1,t_2)d{\bf w}_{t_2}^{(i_2)}d{\bf w}_{t_1}^{(i_1)}+
$$
\begin{equation}
\label{nov801}
+{\bf 1}_{\{i_1=i_2\}}
\int\limits_t^T R_{p_1p_2}(t_1,t_1)dt_1,
\end{equation}

\noindent
where we used the same notations as in (\ref{nov800}).

Uning (\ref{nov801}), we obtain
$$
{\sf M}\left\{\left(J[R_{p_1p_2}]_{T,t}^{(2)}\right)^2\right\}=
\int\limits_t^T\int\limits_t^{t_2}
R_{p_1p_2}^2(t_1,t_2)dt_1 dt_2
+\int\limits_t^T\int\limits_t^{t_1}
R_{p_1p_2}^2(t_1,t_2)dt_2dt_1+
$$
$$
+{\bf 1}_{\{i_1=i_2\}} \left(2\int\limits_t^T\int\limits_t^{t_2}
R_{p_1p_2}(t_1,t_2)R_{p_1p_2}(t_2,t_1)dt_1dt_2+
\left(\int\limits_t^T R_{p_1p_2}(t_1,t_1)dt_1\right)^2\right)=
$$
$$
=\int\limits_t^T\int\limits_t^{t_2}
R_{p_1p_2}^2(t_1,t_2)dt_1 dt_2
+\int\limits_t^T\int\limits_{t_2}^{T}
R_{p_1p_2}^2(t_1,t_2)dt_1dt_2+
$$
$$
+{\bf 1}_{\{i_1=i_2\}} \left(\int\limits_t^T\int\limits_t^{t_2}
R_{p_1p_2}(t_1,t_2)R_{p_1p_2}(t_2,t_1)dt_1dt_2+\right.
$$
$$
+\left.
\int\limits_t^T\int\limits_{t_1}^T
R_{p_1p_2}(t_1,t_2)R_{p_1p_2}(t_2,t_1)dt_2dt_1\right)+
{\bf 1}_{\{i_1=i_2\}}\left(\int\limits_t^T R_{p_1p_2}(t_1,t_1)dt_1\right)^2=
$$
$$
=\int\limits_{[t,T]^2}
R_{p_1p_2}^2(t_1,t_2)dt_1 dt_2+
$$
$$
+{\bf 1}_{\{i_1=i_2\}} \left(\int\limits_t^T\int\limits_t^{t_2}
R_{p_1p_2}(t_1,t_2)R_{p_1p_2}(t_2,t_1)dt_1dt_2+\right.
$$
$$
+\left.
\int\limits_t^T\int\limits_{t_2}^T
R_{p_1p_2}(t_1,t_2)R_{p_1p_2}(t_2,t_1)dt_1dt_2\right)+
{\bf 1}_{\{i_1=i_2\}}\left(\int\limits_t^T R_{p_1p_2}(t_1,t_1)dt_1\right)^2=
$$
$$
=\int\limits_{[t,T]^2}
R_{p_1p_2}^2(t_1,t_2)dt_1 dt_2+
$$
\begin{equation}
\label{nov803}
+{\bf 1}_{\{i_1=i_2\}} \left(\int\limits_{[t,T]^2}
R_{p_1p_2}(t_1,t_2)R_{p_1p_2}(t_2,t_1)dt_1dt_2+
\left(\int\limits_t^T R_{p_1p_2}(t_1,t_1)dt_1\right)^2\right).
\end{equation}

\vspace{2mm}

Since the integrals on the right-hand side of (\ref{nov803}) 
exist as Riemann integrals, then they are equal to the 
corresponding Lebesgue integrals. 
Moreover,
$$
\lim\limits_{p_1\to\infty}\lim\limits_{p_2\to\infty}
R_{p_1p_2}(t_1,t_2)=0\ \ \ \hbox{when}\ \ \ 
(t_1,t_2)\in (t, T)^2,
$$
where the left-hand side is bounded on the boundary of 
$[t, T]^2$ (see (\ref{410})).

Then, applying two times (we mean here an iterated passage to the limit
$\lim\limits_{p_1\to\infty}\lim\limits_{p_2\to\infty}$)
the Lebesgue's 
Dominated Convergence Theorem (see the choice of integrable majorants in the proof of Theorem~2.10) 
and taking into account
(\ref{leto8001}), (\ref{leto8002}), and (\ref{d2020}),
we obtain
\begin{equation}
\label{nov805}
\lim\limits_{p_1\to\infty}\lim\limits_{p_2\to\infty}
\int\limits_{[t, T]^2}
R_{p_1p_2}^2(t_1,t_2)dt_1 dt_2=0,
\end{equation}
\begin{equation}
\label{nov806}
\lim\limits_{p_1\to\infty}\lim\limits_{p_2\to\infty}
\int\limits_{[t, T]^2}
R_{p_1p_2}(t_1,t_2)R_{p_1p_2}(t_2,t_1)dt_1 dt_2=0,
\end{equation}
\begin{equation}
\label{nov807}
\lim\limits_{p_1\to\infty}\lim\limits_{p_2\to\infty}
\int\limits_t^T
R_{p_1p_2}(t_1,t_1)dt_1=0.
\end{equation}

\vspace{2mm}

The relations (\ref{nov803})--(\ref{nov807}) imply the following equality
$$
\lim\limits_{p_1\to\infty}\lim\limits_{p_2\to\infty}
{\sf M}\left\{\left(J[R_{p_1p_2}]_{T,t}^{(2)}\right)^2\right\}=0.
$$

\vspace{1mm}
\noindent
The formula (\ref{nov500}) is proved.

Let us prove the relation (\ref{nov502}).
After replacement of the integration order 
in the iterated It\^{o} stochastic integrals 
from (\ref{s1s})
(see Chapter~3) \cite{1}-\cite{12aa}, \cite{old-art-2}, \cite{vini},
\cite{arxiv-25} 
we get w.~p.~1 
$$
J^{*}[\psi^{(3)}]_{T,t}-
\sum_{j_1=0}^{p_1}\sum_{j_2=0}^{p_2}\sum_{j_3=0}^{p_3}
C_{j_3j_2j_1}\zeta^{(i_1)}_{j_1}\zeta^{(i_2)}_{j_2}\zeta^{(i_3)}_{j_3}
=J[R_{p_1p_2p_3}]_{T,t}^{(3)}=
$$
$$
=
\int\limits_t^T\int\limits_t^{t_3}\int\limits_t^{t_2}
R_{p_1 p_2 p_3}(t_1,t_2,t_3)
d{\bf w}_{t_1}^{(i_1)}
d{\bf w}_{t_2}^{(i_2)}
d{\bf w}_{t_3}^{(i_3)}+
$$
$$
+
\int\limits_t^T\int\limits_t^{t_3}\int\limits_t^{t_2}
R_{p_1 p_2 p_3}(t_1,t_3,t_2)
d{\bf w}_{t_1}^{(i_1)}
d{\bf w}_{t_2}^{(i_3)}
d{\bf w}_{t_3}^{(i_2)}+
$$
$$
+
\int\limits_t^T\int\limits_t^{t_3}\int\limits_t^{t_2}
R_{p_1 p_2 p_3}(t_2,t_1,t_3)
d{\bf w}_{t_1}^{(i_2)}
d{\bf w}_{t_2}^{(i_1)}
d{\bf w}_{t_3}^{(i_3)}+
$$
$$
+
\int\limits_t^T\int\limits_t^{t_3}\int\limits_t^{t_2}
R_{p_1 p_2 p_3}(t_2,t_3,t_1)
d{\bf w}_{t_1}^{(i_3)}
d{\bf w}_{t_2}^{(i_1)}
d{\bf w}_{t_3}^{(i_2)}+
$$
$$
+
\int\limits_t^T\int\limits_t^{t_3}\int\limits_t^{t_2}
R_{p_1 p_2 p_3}(t_3,t_2,t_1)
d{\bf w}_{t_1}^{(i_3)}
d{\bf w}_{t_2}^{(i_2)}
d{\bf w}_{t_3}^{(i_1)}+
$$
$$
+
\int\limits_t^T\int\limits_t^{t_3}\int\limits_t^{t_2}
R_{p_1 p_2 p_3}(t_3,t_1,t_2)
d{\bf w}_{t_1}^{(i_2)}
d{\bf w}_{t_2}^{(i_3)}
d{\bf w}_{t_3}^{(i_1)}+
$$
$$
+{\bf 1}_{\{i_1=i_2\}}
\int\limits_t^T\left(\int\limits_t^{T}
R_{p_1 p_2 p_3}(t_2,t_2,t_3)dt_2\right)
d{\bf w}_{t_3}^{(i_3)}+
$$

\newpage  
\noindent
$$
+
{\bf 1}_{\{i_2=i_3\}}
\int\limits_t^T\left(\int\limits_t^{T}
R_{p_1 p_2 p_3}(t_1,t_2,t_2)dt_2\right)
d{\bf w}_{t_1}^{(i_1)}+
$$
\begin{equation}
\label{nov901}
\hspace{-2mm}+{\bf 1}_{\{i_1=i_3\}}
\int\limits_t^T\left(\int\limits_t^{T}
R_{p_1 p_2 p_3}(t_3,t_2,t_3)dt_3\right)
d{\bf w}_{t_2}^{(i_2)}.
\end{equation}

\vspace{2mm}

Let us calculate the second moment of 
$J[R_{p_1p_2p_3}]_{T,t}^{(3)}$ using (\ref{nov901}).
We have
$$
{\sf M}\left\{\left(J[R_{p_1p_2p_3}]_{T,t}^{(3)}\right)^2\right\}=
$$
\begin{equation}
\label{novv1}
=
\int\limits_t^T\int\limits_t^{t_3}\int\limits_t^{t_2}
\left(\sum\limits_{(t_1,t_2,t_3)} R_{p_1 p_2 p_3}^2(t_1,t_2,t_3)\right)dt_1
dt_2dt_3+
\end{equation}
$$
+2\left({\bf 1}_{\{i_1=i_2\}}
\int\limits_t^T\int\limits_t^{t_3}\int\limits_t^{t_2}
G^{(1)}_{p_1 p_2 p_3}(t_1,t_2,t_3)dt_1dt_2dt_3+\right.
$$
$$
+{\bf 1}_{\{i_1=i_3\}}
\int\limits_t^T\int\limits_t^{t_3}\int\limits_t^{t_2}
G^{(2)}_{p_1 p_2 p_3}(t_1,t_2,t_3)dt_1dt_2dt_3+
$$
$$
+{\bf 1}_{\{i_2=i_3\}}
\int\limits_t^T\int\limits_t^{t_3}\int\limits_t^{t_2}
G^{(3)}_{p_1 p_2 p_3}(t_1,t_2,t_3)dt_1dt_2dt_3+
$$
$$
\left.+{\bf 1}_{\{i_1=i_2=i_3\}}
\int\limits_t^T\int\limits_t^{t_3}\int\limits_t^{t_2}
G^{(4)}_{p_1 p_2 p_3}(t_1,t_2,t_3)dt_1dt_2dt_3\right)+
$$
$$
+\int\limits_{[t,T]^3}
\biggl({\bf 1}_{\{i_1=i_2\}}
R_{p_1 p_2 p_3}(t_1,t_1,t_3)R_{p_1 p_2 p_3}(t_2,t_2,t_3)+\biggr.
$$
$$
+{\bf 1}_{\{i_2=i_3\}}
R_{p_1 p_2 p_3}(t_3,t_1,t_1)R_{p_1 p_2 p_3}(t_3,t_2,t_2)+
$$
$$
+{\bf 1}_{\{i_1=i_3\}}
R_{p_1 p_2 p_3}(t_1,t_3,t_1)R_{p_1 p_2 p_3}(t_2,t_3,t_2)+
$$
$$
+2\cdot {\bf 1}_{\{i_1=i_2=i_3\}}\biggl(
R_{p_1 p_2 p_3}(t_1,t_1,t_3)R_{p_1 p_2 p_3}(t_3,t_2,t_2)+\biggr.
$$
$$
+
R_{p_1 p_2 p_3}(t_1,t_1,t_3)R_{p_1 p_2 p_3}(t_2,t_3,t_2)+
$$
\begin{equation}
\label{nov980}
\biggl.\biggl.+
R_{p_1 p_2 p_3}(t_3,t_1,t_1)R_{p_1 p_2 p_3}(t_2,t_3,t_2)\biggr)\biggr)dt_1dt_2dt_3,
\end{equation}

\vspace{1mm}
\noindent
where permutation $(t_1,t_2,t_3)$ when summing in (\ref{novv1}) are 
performed only in the value $R_{p_1 p_2 p_3}^2(t_1,t_2,t_3)$ and 
the functions $G^{(i)}_{p_1 p_2 p_3}(t_1,t_2,t_3)$ $(i=1,\ldots,4)$
are defined by the following relations

\vspace{-4mm}
$$
G^{(1)}_{p_1 p_2 p_3}(t_1,t_2,t_3)=
R_{p_1 p_2 p_3}(t_1,t_2,t_3)R_{p_1 p_2 p_3}(t_2,t_1,t_3)+
$$

\vspace{-4mm}
$$
+R_{p_1 p_2 p_3}(t_1,t_3,t_2)R_{p_1 p_2 p_3}(t_3,t_1,t_2)+
$$

\vspace{-4mm}
$$
+R_{p_1 p_2 p_3}(t_2,t_3,t_1)R_{p_1 p_2 p_3}(t_3,t_2,t_1),
$$

$$
G^{(2)}_{p_1 p_2 p_3}(t_1,t_2,t_3)=
R_{p_1 p_2 p_3}(t_1,t_2,t_3)R_{p_1 p_2 p_3}(t_3,t_2,t_1)+
$$

\vspace{-4mm}
$$
+R_{p_1 p_2 p_3}(t_1,t_3,t_2)R_{p_1 p_2 p_3}(t_2,t_3,t_1)+
$$

\vspace{-4mm}
$$
+R_{p_1 p_2 p_3}(t_2,t_1,t_3)R_{p_1 p_2 p_3}(t_3,t_1,t_2),
$$

$$
G^{(3)}_{p_1 p_2 p_3}(t_1,t_2,t_3)=
R_{p_1 p_2 p_3}(t_1,t_2,t_3)R_{p_1 p_2 p_3}(t_1,t_3,t_2)+
$$

\vspace{-4mm}
$$
+R_{p_1 p_2 p_3}(t_2,t_1,t_3)R_{p_1 p_2 p_3}(t_2,t_3,t_1)+
$$

\vspace{-4mm}
$$
+R_{p_1 p_2 p_3}(t_3,t_2,t_1)R_{p_1 p_2 p_3}(t_3,t_1,t_2),
$$

$$
G^{(4)}_{p_1 p_2 p_3}(t_1,t_2,t_3)=
R_{p_1 p_2 p_3}(t_1,t_2,t_3)R_{p_1 p_2 p_3}(t_2,t_3,t_1)+
$$

\vspace{-4mm}
$$
+R_{p_1 p_2 p_3}(t_1,t_2,t_3)R_{p_1 p_2 p_3}(t_3,t_1,t_2)+
$$

\vspace{-4mm}
$$
+R_{p_1 p_2 p_3}(t_1,t_3,t_2)R_{p_1 p_2 p_3}(t_2,t_1,t_3)+
$$

\vspace{-4mm}
$$
+R_{p_1 p_2 p_3}(t_1,t_3,t_2)R_{p_1 p_2 p_3}(t_3,t_2,t_1)+
$$

\vspace{-4mm}
$$
+R_{p_1 p_2 p_3}(t_2,t_1,t_3)R_{p_1 p_2 p_3}(t_3,t_2,t_1)+
$$

\vspace{-4mm}
$$
+R_{p_1 p_2 p_3}(t_2,t_3,t_1)R_{p_1 p_2 p_3}(t_3,t_1,t_2).
$$

\vspace{3mm}

Further (see (\ref{riemann})),
$$
\int\limits_t^T\int\limits_t^{t_3}\int\limits_t^{t_2}
\left(\sum\limits_{(t_1,t_2,t_3)} R_{p_1 p_2 p_3}^2(t_1,t_2,t_3)\right)dt_1
dt_2dt_3=
$$
\begin{equation}
\label{novv60}
=\int\limits_{[t,T]^3}R_{p_1 p_2 p_3}^2(t_1,t_2,t_3)dt_1
dt_2dt_3.
\end{equation}

We will say that the function $\Phi(t_1,t_2,t_3)$ is symmetric if

\vspace{-3mm}
$$
\Phi(t_1,t_2,t_3)=\Phi(t_1,t_3,t_2)=\Phi(t_2,t_1,t_3)=
\Phi(t_2,t_3,t_1)=
$$
$$
=\Phi(t_3,t_1,t_2)=\Phi(t_3,t_2,t_1).
$$

\vspace{2mm}

For the symmetric function $\Phi(t_1,t_2,t_3)$, we have
$$
\int\limits_t^T\int\limits_t^{t_3}\int\limits_t^{t_2}
\left(\sum\limits_{(t_1,t_2,t_3)} \Phi(t_1,t_2,t_3) \right)dt_1
dt_2dt_3= 
$$
$$
=6\int\limits_t^T\int\limits_t^{t_3}\int\limits_t^{t_2}
\Phi(t_1,t_2,t_3)dt_1dt_2dt_3= 
$$
\begin{equation}
\label{nov900}
=
\int\limits_{[t,T]^3}\Phi(t_1,t_2,t_3)dt_1
dt_2dt_3.
\end{equation}

The relation (\ref{nov900}) implies that
\begin{equation}
\label{nov901a}
~~~~~~~~~~ \int\limits_t^T\int\limits_t^{t_3}\int\limits_t^{t_2}
\Phi(t_1,t_2,t_3)dt_1dt_2dt_3=
\frac{1}{6}\int\limits_{[t,T]^3}\Phi(t_1,t_2,t_3)dt_1
dt_2dt_3.
\end{equation}

\vspace{2mm}

It is easy to check that the functions 
$G^{(i)}_{p_1 p_2 p_3}(t_1,t_2,t_3)$ $(i=1,\ldots,4)$
are symmetric. Using this property as well as 
(\ref{nov980}), (\ref{novv60}), and (\ref{nov901a}), we obtain

\vspace{-4mm}
$$
{\sf M}\left\{\left(J[R_{p_1p_2p_3}]_{T,t}^{(3)}\right)^2\right\}=
\int\limits_{[t,T]^3}R_{p_1 p_2 p_3}^2(t_1,t_2,t_3)dt_1dt_2dt_3+
$$
$$
+\frac{1}{3}\int\limits_{[t,T]^3}\biggl({\bf 1}_{\{i_1=i_2\}}
G^{(1)}_{p_1 p_2 p_3}(t_1,t_2,t_3)dt_1dt_2dt_3+\biggr.
$$
$$
+{\bf 1}_{\{i_1=i_3\}}
G^{(2)}_{p_1 p_2 p_3}(t_1,t_2,t_3)dt_1dt_2dt_3+
$$

\vspace{-2mm}
$$
+{\bf 1}_{\{i_2=i_3\}}
G^{(3)}_{p_1 p_2 p_3}(t_1,t_2,t_3)dt_1dt_2dt_3+
$$

\vspace{-2mm}
$$
\biggl.+{\bf 1}_{\{i_1=i_2=i_3\}}
G^{(4)}_{p_1 p_2 p_3}(t_1,t_2,t_3)dt_1dt_2dt_3\biggr)dt_1dt_2dt_3+
$$
$$
+
\int\limits_{[t,T]^3}\biggl({\bf 1}_{\{i_1=i_2\}}
R_{p_1 p_2 p_3}(t_1,t_1,t_3)R_{p_1 p_2 p_3}(t_2,t_2,t_3)+\biggr.
$$
$$
+{\bf 1}_{\{i_2=i_3\}}
R_{p_1 p_2 p_3}(t_3,t_1,t_1)R_{p_1 p_2 p_3}(t_3,t_2,t_2)+
$$

\vspace{-2mm}
$$
+{\bf 1}_{\{i_1=i_3\}}
R_{p_1 p_2 p_3}(t_1,t_3,t_1)R_{p_1 p_2 p_3}(t_2,t_3,t_2)+
$$

\vspace{-2mm}
$$
+2\cdot {\bf 1}_{\{i_1=i_2=i_3\}}\biggl(
R_{p_1 p_2 p_3}(t_1,t_1,t_3)R_{p_1 p_2 p_3}(t_3,t_2,t_2)+\biggr.
$$

\vspace{-2mm}
$$
+
R_{p_1 p_2 p_3}(t_1,t_1,t_3)R_{p_1 p_2 p_3}(t_2,t_3,t_2)+
$$

\vspace{-2mm}
\begin{equation}
\label{nov9801}
\biggl.\biggl.+
R_{p_1 p_2 p_3}(t_3,t_1,t_1)R_{p_1 p_2 p_3}(t_2,t_3,t_2)\biggr)\biggr)dt_1dt_2dt_3.
\end{equation}

\vspace{3mm}

Since the integrals on the right-hand side of (\ref{nov9801}) 
exist as Riemann integrals, then they are equal to the 
corresponding Lebesgue integrals. 
Moreover, 
$$
\lim\limits_{p_1\to\infty}\lim\limits_{p_2\to\infty}\lim\limits_{p_3\to\infty}
R_{p_1 p_2 p_3}(t_1,t_2,t_3)=0\ \ \ \hbox{when}\ \ \ 
(t_1,t_2,t_3)\in (t, T)^3,
$$
where the left-hand side is bounded on the boundary of 
$[t, T]^3$ (see (\ref{410})).

Using (\ref{novvv1}) and applying three times (we mean here an iterated passage to the limit
$\lim\limits_{p_1\to\infty}\lim\limits_{p_2\to\infty}
\lim\limits_{p_3\to\infty}$)
the Lebesgue's Dominated Convergence Theorem 
(see the choice of integrable majorants in the proof of Theorem~2.10) 
in the equality (\ref{nov9801}), 
we obtain
$$
\lim\limits_{p_1\to\infty}\lim\limits_{p_2\to\infty}\lim\limits_{p_3\to\infty}
{\sf M}\left\{\left(J[R_{p_1p_2p_3}]_{T,t}^{(3)}\right)^2\right\}=0.
$$

\vspace{1mm}
\noindent
The relation (\ref{nov502}) is proved. Theorem 2.17 is proved.

Developing the approach used in the proof of 
Theorem 2.17, we can in principle prove 
the following formulas

\vspace{-5mm}
$$
\lim\limits_{p_1\to\infty}
\ldots
\lim\limits_{p_k\to\infty}
{\sf M}\left\{\left(J^{*}[\psi^{(k)}]_{T,t}-
\sum_{j_1=0}^{p_1}\ldots \sum_{j_k=0}^{p_k}
C_{j_k \ldots j_1}\zeta^{(i_1)}_{j_1}
\ldots \zeta^{(i_k)}_{j_k}\right)^{2}\right\}=0,
$$

\vspace{-3mm}
$$
\lim\limits_{p_k\to\infty}
\ldots
\lim\limits_{p_1\to\infty}
{\sf M}\left\{\left(J^{*}[\psi^{(k)}]_{T,t}-
\sum_{j_k=0}^{p_k}\ldots \sum_{j_1=0}^{p_1}
C_{j_k \ldots j_1}\zeta^{(i_1)}_{j_1}
\ldots \zeta^{(i_k)}_{j_k}\right)^{2}\right\}=0,
$$

\vspace{2mm}
\noindent
which are correct under the conditions of Theorem 2.10 for $i_1,\ldots,i_k=1,\ldots,m$.

\section{The Hypotheses on Expansion of Iterated Stra\-to\-no\-vich Stochastic
Integrals of Multiplicity $k$ $(k\in{\bf N})$ Based on Theorem 1.1}

In this section, on the base of the presented theorems (see Sect.~1.1.3, 2.1--2.4) 
we formulate 3 hypotheses
on expansions of iterated Stratonovich stochastic integrals of arbitrary 
multiplicity $k$ ($k\in{\bf N}$)
based on generalized multiple Fourier series converging in $L_2([t, T]^k).$
The considered expansions contain
only one operation of the limit transition
and substantially simpler than
their analogues for iterated It\^{o} stochastic integrals (Theorem 1.1).

Taking into account (\ref{ziko5000}) and Theorems 2.1--2.10, 2.14, and 2.17, let us 
formulate the following hypotheses on expansions of iterated 
Stratonovich stochastic
integrals of multiplicity $k$ ($k\in{\bf N}$).

{\bf Hypothesis 2.1}\ \cite{8}-\cite{12aa}, \cite{arxiv-11}. {\it Assume that
$\{\phi_j(x)\}_{j=0}^{\infty}$ is a complete orthonormal
system of Legendre polynomials or trigonometric functions
in the space $L_2([t, T])$.
Then$,$ for the iterated Stratonovich stochastic 
integral of multiplicity $k$
\begin{equation}
\label{otitgo10}
I_{(\lambda_1\ldots\lambda_k)T,t}^{*(i_1\ldots i_k)}=
{\int\limits_t^{*}}^T
\ldots
{\int\limits_t^{*}}^{t_2}
d{\bf w}_{t_1}^{(i_1)}
\ldots d{\bf w}_{t_k}^{(i_k)}
\end{equation}
the following 
expansion 
\begin{equation}
\label{feto1900otita}
I_{(\lambda_1\ldots\lambda_k)T,t}^{*(i_1\ldots i_k)}=
\hbox{\vtop{\offinterlineskip\halign{
\hfil#\hfil\cr
{\rm l.i.m.}\cr
$\stackrel{}{{}_{p\to \infty}}$\cr
}} }
\sum\limits_{j_1,\ldots j_k=0}^{p}
C_{j_k \ldots j_1}\prod_{l=1}^k\zeta_{j_l}^{(i_l)}
\end{equation}

\vspace{2mm}
\noindent
that converges in the mean-square sense is valid, 
where  
$$
C_{j_k \ldots j_1}=\int\limits_t^T\phi_{j_k}(t_k)\ldots
\int\limits_t^{t_2}
\phi_{j_1}(t_1)
dt_1\ldots dt_k
$$
is the Fourier coefficient,
${\rm l.i.m.}$ is a limit in the mean-square sense,
$i_1,\ldots, i_k=0, 1,\ldots,m,$
$$
\zeta_{j}^{(i)}=
\int\limits_t^T \phi_{j}(s) d{\bf w}_s^{(i)}
$$ 
are independent standard Gaussian random variables for various 
$i$ or $j$ {\rm (}in the case when $i\ne 0${\rm )},
${\bf w}_{\tau}^{(i)}$ $(i=1,\ldots,m)$ are independent 
standard Wiener processes and 
${\bf w}_{\tau}^{(0)}=\tau,$
$\lambda_l=0$ if $i_l=0$ and $\lambda_l=1$ if $i_l=1,\ldots,m$
$(l=1,\ldots,k).$}

Hypothesis 2.1 allows to approximate the iterated
Stratonovich stochastic integral 
$I_{(\lambda_1\ldots\lambda_k)T,t}^{*(i_1\ldots i_k)}$
by the sum
\begin{equation}
\label{otit567}
I_{(\lambda_1\ldots\lambda_k)T,t}^{*(i_1\ldots i_k)p}=
\sum\limits_{j_1,\ldots j_k=0}^{p}
C_{j_k \ldots j_1}\prod\limits_{l=1}^k
\zeta_{j_l}^{(i_l)},
\end{equation}

\noindent
where
$$
\lim_{p\to\infty}{\sf M}\left\{\Biggl(
I_{(\lambda_1\ldots\lambda_k)T,t}^{*(i_1\ldots i_k)}-
I_{(\lambda_1\ldots\lambda_k)T,t}
^{*(i_1\ldots i_k)p}\Biggr)^2\right\}=0.
$$

The integrals (\ref{otitgo10}) will be used in the 
Taylor--Stratonovich expansion (see Chapter 4). It means that
the approximations (\ref{otit567})
may be very useful for the construction of high-order 
strong numerical methods
for It\^{o} SDEs (see Chapter 4 for details).

The expansion (\ref{feto1900otita}) contains only one operation of the limit
transition and by this reason is convenient for approximation
of iterated Stratonovich stochastic integrals.

Let us consider the more general hypothesis than Hypothesis 2.1.

{\bf Hypothesis 2.2} \cite{12a}-\cite{12aa},  \cite{arxiv-11}. {\it Assume that
$\{\phi_j(x)\}_{j=0}^{\infty}$ is a complete orthonormal
system of Legendre polynomials or trigonometric functions
in the space $L_2([t, T])$. Moreover,
every $\psi_l(\tau)$ $(l=1,\ldots,k)$ is 
an enough smooth nonrandom function
on $[t,T].$
Then$,$ for the iterated Stratonovich stochastic integral 
of multiplicity $k$
$$
~~J^{*}[\psi^{(k)}]_{T,t}=
{\int\limits_t^{*}}^T
\psi_k(t_k) \ldots 
{\int\limits_t^{*}}^{t_2}
\psi_1(t_1) d{\bf w}_{t_1}^{(i_1)}\ldots
d{\bf w}_{t_k}^{(i_k)}
$$

\noindent
the following 
expansion 
\begin{equation}
\label{feto1900otitarrr}
J^{*}[\psi^{(k)}]_{T,t}=
\hbox{\vtop{\offinterlineskip\halign{
\hfil#\hfil\cr
{\rm l.i.m.}\cr
$\stackrel{}{{}_{p\to \infty}}$\cr
}} }
\sum\limits_{j_1,\ldots j_k=0}^{p}
C_{j_k \ldots j_1}\prod\limits_{l=1}^k \zeta_{j_l}^{(i_l)}
\end{equation}

\vspace{1mm}
\noindent
that converges in the mean-square sense is valid, where 
$$
C_{j_k \ldots j_1}=\int\limits_t^T\psi_k(t_k)\phi_{j_k}(t_k)\ldots
\int\limits_t^{t_2}
\psi_1(t_1)\phi_{j_1}(t_1)
dt_1\ldots dt_k
$$
is the Fourier coefficient, 
${\rm l.i.m.}$ is a limit in the mean-square sense,
$i_1, \ldots, i_k=0, 1,\ldots,m,$
$$
\zeta_{j}^{(i)}=
\int\limits_t^T \phi_{j}(s) d{\bf w}_s^{(i)}
$$ 

\noindent
are independent standard Gaussian random variables for various 
$i$ or $j$ {\rm (}in the case when $i\ne 0${\rm );}
another notations are the same as above in this section.}

Hypothesis 2.2 allows to approximate the iterated
Stratonovich stochastic integral 
$J^{*}[\psi^{(k)}]_{T,t}$
by the sum
\begin{equation}
\label{otit567r}
J^{*}[\psi^{(k)}]_{T,t}^p=
\sum\limits_{j_1,\ldots j_k=0}^{p}
C_{j_k \ldots j_1}\prod\limits_{l=1}^k
\zeta_{j_l}^{(i_l)},
\end{equation}

\noindent
where
$$
\lim_{p\to\infty}{\sf M}\left\{\Biggl(
J^{*}[\psi^{(k)}]_{T,t}-
J^{*}[\psi^{(k)}]_{T,t}^p\Biggl)^2\right\}=0.
$$

Let us consider the more general hypothesis than 
Hypotheses 2.1 and 2.2.

{\bf Hypothesis 2.3} \cite{12a}-\cite{12aa}, \cite{arxiv-11}. {\it Assume that
$\{\phi_j(x)\}_{j=0}^{\infty}$ is a complete orthonormal
system of Legendre polynomials or trigonometric functions
in the space $L_2([t, T])$. Moreover,
every $\psi_l(\tau)$ $(l=1,\ldots,k)$ is 
an enough smooth nonrandom function
on $[t,T].$
Then$,$ for the iterated Stratonovich stochastic integral 
of multiplicity $k$
$$
~~J^{*}[\psi^{(k)}]_{T,t}=
{\int\limits_t^{*}}^T
\psi_k(t_k) \ldots 
{\int\limits_t^{*}}^{t_2}
\psi_1(t_1) d{\bf w}_{t_1}^{(i_1)}\ldots
d{\bf w}_{t_k}^{(i_k)}
$$

\noindent
the following 
expansion 
\begin{equation}
\label{feto1900otitarrri}
J^{*}[\psi^{(k)}]_{T,t}=
\hbox{\vtop{\offinterlineskip\halign{
\hfil#\hfil\cr
{\rm l.i.m.}\cr
$\stackrel{}{{}_{p_1,\ldots,p_k\to \infty}}$\cr
}} }
\sum\limits_{j_1=0}^{p_1}
\ldots\sum\limits_{j_k=0}^{p_k}
C_{j_k \ldots j_1}\prod\limits_{l=1}^k \zeta_{j_l}^{(i_l)}
\end{equation}

\vspace{1mm}
\noindent
that converges in the mean-square sense is valid, where 
$$
C_{j_k \ldots j_1}=\int\limits_t^T\psi_k(t_k)\phi_{j_k}(t_k)\ldots
\int\limits_t^{t_2}
\psi_1(t_1)\phi_{j_1}(t_1)
dt_1\ldots dt_k
$$
is the Fourier coefficient, 
${\rm l.i.m.}$ is a limit in the mean-square sense,
$i_1,\ldots, i_k=0, 1,\ldots,m,$
$$
\zeta_{j}^{(i)}=
\int\limits_t^T \phi_{j}(s) d{\bf w}_s^{(i)}
$$ 

\noindent
are independent standard Gaussian random variables for various 
$i$ or $j$ {\rm (}in the case when $i\ne 0${\rm );}
another notations are the same as above in this section.}

Let us consider the idea of the proof of Hypotheses 2.1--2.3.

According to (\ref{tyyy}), we have

\vspace{-1mm}
$$
\hbox{\vtop{\offinterlineskip\halign{
\hfil#\hfil\cr
{\rm l.i.m.}\cr
$\stackrel{}{{}_{p_1,\ldots,p_k\to \infty}}$\cr
}} }\sum_{j_1=0}^{p_1}\ldots\sum_{j_k=0}^{p_k}
C_{j_k\ldots j_1}\prod_{g=1}^k\zeta_{j_g}^{(i_g)}=
J[\psi^{(k)}]_{T,t}+
$$

\vspace{-1mm}
\begin{equation}
\label{tyyy1}
+
\hbox{\vtop{\offinterlineskip\halign{
\hfil#\hfil\cr
{\rm l.i.m.}\cr
$\stackrel{}{{}_{p_1,\ldots,p_k\to \infty}}$\cr
}} }\sum_{j_1=0}^{p_1}\ldots\sum_{j_k=0}^{p_k}
C_{j_k\ldots j_1}
~\hbox{\vtop{\offinterlineskip\halign{
\hfil#\hfil\cr
{\rm l.i.m.}\cr
$\stackrel{}{{}_{N\to \infty}}$\cr
}} }\sum_{(l_1,\ldots,l_k)\in {\rm G}_k}
\prod_{g=1}^k
\phi_{j_{g}}(\tau_{l_g})
\Delta{\bf w}_{\tau_{l_g}}^{(i_g)}\ \ \ \hbox{w.~p.~1},
\end{equation}

\vspace{4mm}
\noindent
where notations are the same as in (\ref{tyyy}).

From (\ref{tyyy1}) and Theorem 2.12 it follows that

\newpage
\noindent
\begin{equation}
\label{tt1}
J^{*}[\psi^{(k)}]_{T,t}=\hbox{\vtop{\offinterlineskip\halign{
\hfil#\hfil\cr
{\rm l.i.m.}\cr
$\stackrel{}{{}_{p_1,\ldots,p_k\to \infty}}$\cr
}} }\sum_{j_1=0}^{p_1}\ldots\sum_{j_k=0}^{p_k}
C_{j_k\ldots j_1}\prod_{g=1}^k\zeta_{j_g}^{(i_g)}
\end{equation}

\vspace{4mm}
\noindent
if 
$$
\sum_{r=1}^{\left[k/2\right]}\frac{1}{2^r}
\sum_{(s_r,\ldots,s_1)\in {\rm A}_{k,r}}
J[\psi^{(k)}]_{T,t}^{s_r,\ldots,s_1}=
$$

\vspace{-1mm}
$$
=
\hbox{\vtop{\offinterlineskip\halign{
\hfil#\hfil\cr
{\rm l.i.m.}\cr
$\stackrel{}{{}_{p_1,\ldots,p_k\to \infty}}$\cr
}} }\sum_{j_1=0}^{p_1}\ldots\sum_{j_k=0}^{p_k}
C_{j_k\ldots j_1}
~\hbox{\vtop{\offinterlineskip\halign{
\hfil#\hfil\cr
{\rm l.i.m.}\cr
$\stackrel{}{{}_{N\to \infty}}$\cr
}} }\sum_{(l_1,\ldots,l_k)\in {\rm G}_k}
\prod_{g=1}^k
\phi_{j_{g}}(\tau_{l_g})
\Delta{\bf w}_{\tau_{l_g}}^{(i_g)}\ \ \ \hbox{w.~p.~1,}
$$

\vspace{4mm}
\noindent
where notations are the same as in Theorems 1.1 and  2.12. 

Note that from Theorem~1.1 
for pairwise different $i_1,\ldots,i_k$
($i_1,\ldots,i_k=0,1,\ldots,m$) we obtain (\ref{tt1})
(compare (\ref{ziko5000}) and (\ref{tt1})).

In the case $p_1=\ldots=p_k=p$ and $\psi_l(s)\equiv 1$ $(l=1,\ldots,k)$
we obtain from (\ref{tt1}) the statement of Hypothesis 2.1 
(see (\ref{feto1900otita})).

If 
$p_1=\ldots=p_k=p$ and every $\psi_l(s)$ $(l=1,\ldots,k)$
is an enough smooth nonrandom function
on $[t,T],$ then
we obtain from (\ref{tt1}) the statement of Hypothesis 2.2 
(see (\ref{feto1900otitarrr})).

In the case when every $\psi_l(s)$ $(l=1,\ldots,k)$
is an enough smooth nonrandom function
on $[t,T]$ 
we obtain from (\ref{tt1}) the statement of Hypothesis 2.3 
(see (\ref{feto1900otitarrri})).

\section{Expansions of Iterated Stratonovich Stochastic Integrals 
of Multiplicities 3 and 4. Combained Approach Based on
Generalized Multiple and Iterated Fourier series.
Another Proof of Theorems
2.8 and 2.9}

In this section, we develop the approach from Sect.~2.1.3 for 
iterated Stra\-to\-no\-vich stochastic 
integrals of multiplicities 3 and 4. We 
call this approach 
the combined approach of generalized multiple and iterated
Fourier series. We consider two different parts of the expansion 
of iterated Stratonovich stochastic integrals. The mean-square 
convergence of the first part is proved on the base of generalized 
multiple Fourier series converging in the sense of 
norm in $L_2([t,T]^k)$, $k=3, 4$. The mean-square convergence 
of the second part is proved on the base 
of (\ref{a2}), (\ref{5t}), 
Parseval's equality, and generalized 
Fourier series converging pointwise. At that, we do 
not use iterated It\^{o} stochastic integrals as a tool of the 
proof and directly consider iterated Stratonovich stochastic 
integrals.

\subsection{Another Proof of Theorem 2.8}

Let us consider (\ref{proof1ggg}) for
$k=3$, $p_1=p_2=p_3=p$, and $i_1,i_2,i_3=1,\ldots,m$
\begin{equation}
\label{newbegin97}
~~~~~ J^{*}[\psi^{(3)}]_{T,t}=
\sum_{j_1=0}^{p}\sum_{j_2=0}^{p}\sum_{j_3=0}^{p}C_{j_3j_2 j_1}
\zeta_{j_1}^{(i_1)}\zeta_{j_2}^{(i_2)}\zeta_{j_3}^{(i_3)}
+J[R_{ppp}]_{T,t}^{(3)}\ \ \ \hbox{w.~p.~1,}
\end{equation}

\vspace{2mm}
\noindent
where
$$
J[R_{ppp}]_{T,t}^{(3)}=
\hbox{\vtop{\offinterlineskip\halign{
\hfil#\hfil\cr
{\rm l.i.m.}\cr
$\stackrel{}{{}_{N\to \infty}}$\cr
}} }\sum_{l_3=0}^{N-1}\sum_{l_2=0}^{N-1}
\sum_{l_1=0}^{N-1}
R_{ppp}(\tau_{l_1},\tau_{l_2},\tau_{l_3})
\Delta{\bf w}_{\tau_{l_1}}^{(i_1)}
\Delta{\bf w}_{\tau_{l_2}}^{(i_2)}
\Delta{\bf w}_{\tau_{l_3}}^{(i_3)},
$$

\vspace{-3mm}
$$
R_{ppp}(t_1,t_2,t_3)
\stackrel{\small{\sf def}}{=}K^{*}(t_1,t_2,t_3)-
\sum_{j_1=0}^{p}\sum_{j_2=0}^{p}\sum_{j_3=0}^{p}
C_{j_3j_2 j_1}\phi_{j_1}(t_1)\phi_{j_2}(t_2)\phi_{j_3}(t_3),
$$

\vspace{-3mm}
$$
K^{*}(t_1,t_2,t_3)
=\prod_{l=1}^{3}\psi_l(t_l)\Biggl(
{\bf 1}_{\{t_1<t_2\}}{\bf 1}_{\{t_2<t_3\}}+
\frac{1}{2}{\bf 1}_{\{t_1=t_2\}}{\bf 1}_{\{t_2<t_3\}}+\Biggr.
$$

\vspace{-3mm}
$$
\Biggl.
+\frac{1}{2}{\bf 1}_{\{t_1<t_2\}}{\bf 1}_{\{t_2=t_3\}}+
\frac{1}{4}{\bf 1}_{\{t_1=t_2\}}{\bf 1}_{\{t_2=t_3\}}\Biggr).
$$

\vspace{2mm}

Using (\ref{s1s}), we obtain w.~p.~1 
$$
J[R_{ppp}]_{T,t}^{(3)}=
R_{T,t}^{(1)ppp}+R_{T,t}^{(2)ppp},
$$
where
$$
R_{T,t}^{(1)ppp}
=\int\limits_t^T\int\limits_t^{t_3}\int\limits_t^{t_2}
R_{ppp}(t_1,t_2,t_3)
d{\bf w}_{t_1}^{(i_1)}
d{\bf w}_{t_2}^{(i_2)}
d{\bf w}_{t_3}^{(i_3)}+
$$
$$
+
\int\limits_t^T\int\limits_t^{t_3}\int\limits_t^{t_2}
R_{ppp}(t_1,t_3,t_2)
d{\bf w}_{t_1}^{(i_1)}
d{\bf w}_{t_2}^{(i_3)}
d{\bf w}_{t_3}^{(i_2)}+
$$
$$
+
\int\limits_t^T\int\limits_t^{t_3}\int\limits_t^{t_2}
R_{ppp}(t_2,t_1,t_3)
d{\bf w}_{t_1}^{(i_2)}
d{\bf w}_{t_2}^{(i_1)}
d{\bf w}_{t_3}^{(i_3)}+
$$
$$
+
\int\limits_t^T\int\limits_t^{t_3}\int\limits_t^{t_2}
R_{ppp}(t_2,t_3,t_1)
d{\bf w}_{t_1}^{(i_3)}
d{\bf w}_{t_2}^{(i_1)}
d{\bf w}_{t_3}^{(i_2)}+
$$
$$
+
\int\limits_t^T\int\limits_t^{t_3}\int\limits_t^{t_2}
R_{ppp}(t_3,t_2,t_1)
d{\bf w}_{t_1}^{(i_3)}
d{\bf w}_{t_2}^{(i_2)}
d{\bf w}_{t_3}^{(i_1)}+
$$
$$
+
\int\limits_t^T\int\limits_t^{t_3}\int\limits_t^{t_2}
R_{ppp}(t_3,t_1,t_2)
d{\bf w}_{t_1}^{(i_2)}
d{\bf w}_{t_2}^{(i_3)}
d{\bf w}_{t_3}^{(i_1)},
$$

\vspace{3mm}
$$
R_{T,t}^{(2)ppp}
={\bf 1}_{\{i_1=i_2\ne 0\}}
\int\limits_t^T\int\limits_t^{t_3}
R_{ppp}(t_2,t_2,t_3)
dt_2
d{\bf w}_{t_3}^{(i_3)}+
$$
$$
+
{\bf 1}_{\{i_1=i_3\ne 0\}}
\int\limits_t^T\int\limits_t^{t_3}
R_{ppp}(t_2,t_3,t_2)
dt_2
d{\bf w}_{t_3}^{(i_2)}+
$$
$$
+{\bf 1}_{\{i_2=i_3\ne 0\}}
\int\limits_t^T\int\limits_t^{t_3}
R_{ppp}(t_3,t_2,t_2)
dt_2
d{\bf w}_{t_3}^{(i_1)}+
$$
$$
+
{\bf 1}_{\{i_2=i_3\ne 0\}}
\int\limits_t^T\int\limits_t^{t_3}
R_{ppp}(t_1,t_3,t_3)
d{\bf w}_{t_1}^{(i_1)}dt_3+
$$
$$
+{\bf 1}_{\{i_1=i_3\ne 0\}}
\int\limits_t^T\int\limits_t^{t_3}
R_{ppp}(t_3,t_1,t_3)
d{\bf w}_{t_1}^{(i_2)}dt_3+
$$
$$
+
{\bf 1}_{\{i_1=i_2\ne 0\}}
\int\limits_t^T\int\limits_t^{t_3}
R_{ppp}(t_3,t_3,t_1)
d{\bf w}_{t_1}^{(i_3)}dt_3.
$$

We have
\begin{equation}
\label{newbegin98}
~~~~~~~~~~{\sf M}\left\{\left(J[R_{ppp}]_{T,t}^{(3)}\right)^2\right\}\le
2{\sf M}\left\{\left(R_{T,t}^{(1)ppp}\right)^2\right\}
+2{\sf M}\left\{\left(R_{T,t}^{(2)ppp}\right)^2\right\}.
\end{equation}

Now, using standard estimates for moments of stochastic 
integrals \cite{Gih1}, we obtain the following inequality

\vspace{-2mm}
$$
{\sf M}\left\{\left(R_{T,t}^{(1)ppp}\right)^{2}
\right\}\le 
$$
$$
\le 6 \int\limits_t^T\int\limits_t^{t_3}\int\limits_t^{t_2}
\Biggl(
\left(R_{p_1 p_2 p_3}(t_1,t_2,t_3)\right)^{2}+
\left(R_{p_1 p_2 p_3}(t_1,t_3,t_2)\right)^{2}+\Biggr.
\left(R_{p_1 p_2 p_3}(t_2,t_1,t_3)\right)^{2}+
$$
$$
+\left(R_{p_1 p_2 p_3}(t_2,t_3,t_1)\right)^{2}+
\left(R_{p_1 p_2 p_3}(t_3,t_2,t_1)\right)^{2}+
\Biggl.\left(R_{p_1 p_2 p_3}(t_3,t_1,t_2)\right)^{2}\Biggr)dt_1dt_2dt_3=
$$

\vspace{2mm}
$$
=6\int\limits_{[t, T]^3}
\left(R_{ppp}(t_1,t_2,t_3)\right)^{2}dt_1 dt_2 dt_3.
$$

\vspace{3mm}

We have

\vspace{-3mm}
$$
\int\limits_{[t, T]^3}
\left(R_{ppp}(t_1,t_2,t_3)\right)^{2}dt_1 dt_2 dt_3=
$$
$$
=
\int\limits_{[t, T]^3}
\Biggl(
K^{*}(t_1,t_2,t_3)-
\sum_{j_1=0}^{p}\sum_{j_2=0}^{p}\sum_{j_3=0}^{p}C_{j_3 j_2 j_1}
\phi_{j_1}(t_1)\phi_{j_2}(t_2)\phi_{j_3}(t_3)\Biggr)^2 
dt_1 dt_2 dt_3=
$$
$$
=
\int\limits_{[t, T]^3}
\Biggl(
K(t_1,t_2,t_3)-
\sum_{j_1=0}^{p}\sum_{j_2=0}^{p}\sum_{j_3=0}^{p}C_{j_3 j_2 j_1}
\phi_{j_1}(t_1)\phi_{j_2}(t_2)\phi_{j_3}(t_3)\Biggr)^2 dt_1 dt_2 dt_3,
$$

\vspace{4mm}
\noindent
where
$$
K(t_1,t_2,t_3)=
\left\{
\begin{matrix}
\psi_1(t_1)\psi_2(t_2)\psi_3(t_3),\ &t_1<t_2<t_3\cr\cr
0,\ &\hbox{\rm otherwise}
\end{matrix},\right. \ \ \ t_1, t_2, t_3\in[t, T].
$$

\vspace{3mm}

So, we get

\vspace{-2mm}
$$
\hbox{\vtop{\offinterlineskip\halign{
\hfil#\hfil\cr
{\rm lim}\cr
$\stackrel{}{{}_{p\to \infty}}$\cr
}} }
{\sf M}\left\{\left(R_{T,t}^{(1)ppp}\right)^{2}\right\}\le
$$
$$
\le 6\hbox{\vtop{\offinterlineskip\halign{
\hfil#\hfil\cr
{\rm lim}\cr
$\stackrel{}{{}_{p\to \infty}}$\cr
}} }
\int\limits_{[t, T]^3}\hspace{-1mm}
\Biggl(
K(t_1,t_2,t_3)-
\sum_{j_1=0}^{p}\sum_{j_2=0}^{p}\sum_{j_3=0}^{p}C_{j_3 j_2 j_1}
\phi_{j_1}(t_1)\phi_{j_2}(t_2)\phi_{j_3}(t_3)\Biggr)^2 dt_1 dt_2 dt_3=
$$

\vspace{-2mm}
\begin{equation}
\label{newbegin99}
=0,
\end{equation}

\vspace{1mm}
\noindent
where
$K(t_1,t_2,t_3) \in L_2([t, T]^3)$.

After replacement of the integration order 
in the iterated It\^{o} stochastic integrals 
from $R_{T,t}^{(2)ppp}$
\cite{1}-\cite{12aa}, \cite{old-art-2}, \cite{vini},
\cite{arxiv-25} (see Chapter 3)
we obtain  w.~p.~1 

$$
R_{T,t}^{(2)ppp}=
$$
$$
={\bf 1}_{\{i_1=i_2\ne 0\}}\left(
\int\limits_t^T\int\limits_t^{t_3}
R_{ppp}(t_2,t_2,t_3)
dt_2
d{\bf w}_{t_3}^{(i_3)}+\Biggr.
\Biggl.
\int\limits_t^T\int\limits_t^{t_3}
R_{ppp}(t_3,t_3,t_1)
d{\bf w}_{t_1}^{(i_3)}dt_3\right)+
$$
$$
+{\bf 1}_{\{i_2=i_3\ne 0\}}\left(
\int\limits_t^T\int\limits_t^{t_3}
R_{ppp}(t_3,t_2,t_2)
dt_2
d{\bf w}_{t_3}^{(i_1)}+\Biggr.
\Biggl.
\int\limits_t^T\int\limits_t^{t_3}
R_{ppp}(t_1,t_3,t_3)
d{\bf w}_{t_1}^{(i_1)}dt_3\right)+
$$
$$
+{\bf 1}_{\{i_1=i_3\ne 0\}}\left(
\int\limits_t^T\int\limits_t^{t_3}
R_{ppp}(t_2,t_3,t_2)dt_2
d{\bf w}_{t_3}^{(i_2)}+\Biggr.
\Biggl.
\int\limits_t^T\int\limits_t^{t_3}
R_{ppp}(t_3,t_1,t_3)
d{\bf w}_{t_1}^{(i_2)}dt_3\right)=
$$
$$
={\bf 1}_{\{i_1=i_2\ne 0\}}\left(
\int\limits_t^T\int\limits_t^{t_1}
R_{ppp}(t_2,t_2,t_1)
dt_2
d{\bf w}_{t_1}^{(i_3)}+\Biggr.
\Biggl.
\int\limits_t^T\int\limits_{t_1}^{T}
R_{ppp}(t_2,t_2,t_1)
dt_2d{\bf w}_{t_1}^{(i_3)}\right)+
$$
$$
+{\bf 1}_{\{i_2=i_3\ne 0\}}\left(
\int\limits_t^T\int\limits_t^{t_1}
R_{ppp}(t_1,t_2,t_2)
dt_2
d{\bf w}_{t_1}^{(i_1)}+\Biggr.
\Biggl.
\int\limits_t^T\int\limits_{t_1}^{T}
R_{ppp}(t_1,t_2,t_2)
dt_2d{\bf w}_{t_1}^{(i_1)}\right)+
$$
$$
+{\bf 1}_{\{i_1=i_3\ne 0\}}\left(
\int\limits_t^T\int\limits_t^{t_1}
R_{ppp}(t_2,t_1,t_2)dt_2
d{\bf w}_{t_1}^{(i_2)}+\Biggr.
\Biggl.
\int\limits_t^T\int\limits_{t_1}^{T}
R_{ppp}(t_2,t_1,t_2)
dt_2d{\bf w}_{t_1}^{(i_2)}\right)=
$$
$$
={\bf 1}_{\{i_1=i_2\ne 0\}}
\int\limits_t^T\left(\int\limits_t^{T}
R_{ppp}(t_2,t_2,t_3)dt_2\right)
d{\bf w}_{t_3}^{(i_3)}+
$$
$$
+
{\bf 1}_{\{i_2=i_3\ne 0\}}
\int\limits_t^T\left(\int\limits_t^{T}
R_{ppp}(t_1,t_2,t_2)dt_2\right)
d{\bf w}_{t_1}^{(i_1)}+
$$
$$
+{\bf 1}_{\{i_1=i_3\ne 0\}}
\int\limits_t^T\left(\int\limits_t^{T}
R_{ppp}(t_3,t_2,t_3)dt_3\right)
d{\bf w}_{t_2}^{(i_2)}=
$$
$$
={\bf 1}_{\{i_1=i_2\ne 0\}}
\int\limits_t^T\int\limits_t^{T}\Biggl(\Biggl(
\frac{1}{2}{\bf 1}_{\{t_2<t_3\}} + \frac{1}{4}{\bf 1}_{\{t_2=t_3\}}
\Biggr)
\psi_1(t_2)\psi_2(t_2)\psi_3(t_3)-\Biggr.
$$
$$
\Biggl.
-\sum_{j_1=0}^{p}\sum_{j_2=0}^{p}\sum_{j_3=0}^{p}
C_{j_3 j_2 j_1}\phi_{j_1}(t_2)\phi_{j_2}(t_2)\phi_{j_3}(t_3)
\Biggr)dt_2 d{\bf w}_{t_3}^{(i_3)}+
$$
$$
+{\bf 1}_{\{i_2=i_3\ne 0\}}
\int\limits_t^T\int\limits_t^{T}\Biggl(\Biggl(
\frac{1}{2}{\bf 1}_{\{t_1<t_2\}} + \frac{1}{4}{\bf 1}_{\{t_1=t_2\}}
\Biggr)
\psi_1(t_1)\psi_2(t_2)\psi_3(t_2)-\Biggr.
$$
$$
\Biggl.
-\sum_{j_1=0}^{p}\sum_{j_2=0}^{p}\sum_{j_3=0}^{p}
C_{j_3 j_2 j_1}\phi_{j_1}(t_1)\phi_{j_2}(t_2)\phi_{j_3}(t_2)
\Biggr)dt_2 d{\bf w}_{t_1}^{(i_1)}+
$$
$$
+{\bf 1}_{\{i_1=i_3\ne 0\}}
\int\limits_t^T\int\limits_t^{T}\Biggl(
\frac{1}{4}{\bf 1}_{\{t_2=t_3\}}
\psi_1(t_3)\psi_2(t_2)\psi_3(t_3)-\Biggr.
$$
$$
\Biggl.
-\sum_{j_1=0}^{p}\sum_{j_2=0}^{p}\sum_{j_3=0}^{p}
C_{j_3 j_2 j_1}\phi_{j_1}(t_3)\phi_{j_2}(t_2)\phi_{j_3}(t_3)
\Biggr)dt_3 d{\bf w}_{t_2}^{(i_2)}=
$$
$$
={\bf 1}_{\{i_1=i_2\ne 0\}}
\int\limits_t^T\left(\frac{1}{2}\psi_3(t_3)\int\limits_t^{t_3}
\psi_1(t_2)\psi_2(t_2)dt_2-\Biggr.
\Biggl.
\sum_{j_1=0}^{p}\sum_{j_3=0}^{p}
C_{j_3 j_1 j_1}\phi_{j_3}(t_3)
\right)d{\bf w}_{t_3}^{(i_3)}+
$$
$$
+{\bf 1}_{\{i_2=i_3\ne 0\}}
\int\limits_t^T\left(\frac{1}{2}\psi_1(t_1)\int\limits_{t_1}^{T}
\psi_2(t_2)\psi_3(t_2)dt_2-\Biggr.
\Biggl.
\sum_{j_1=0}^{p}\sum_{j_3=0}^{p}
C_{j_3 j_3 j_1}\phi_{j_1}(t_1)
\right)d{\bf w}_{t_1}^{(i_1)}+
$$
$$
+{\bf 1}_{\{i_1=i_3\ne 0\}}
\int\limits_t^T (-1)\sum_{j_1=0}^{p}\sum_{j_2=0}^{p}
C_{j_1 j_2 j_1}\phi_{j_2}(t_2)
d{\bf w}_{t_2}^{(i_2)}=
$$
$$
={\bf 1}_{\{i_1=i_2\ne 0\}}\left(
\frac{1}{2}\int\limits_t^T\psi_3(t_3)\int\limits_t^{t_3}
\psi_1(t_2)\psi_2(t_2)dt_2d{\bf w}_{t_3}^{(i_3)}-
\sum_{j_1=0}^{p}\sum_{j_3=0}^{p}
C_{j_3 j_1 j_1}\zeta_{j_3}^{(i_3)}\right)+
$$
$$
+{\bf 1}_{\{i_2=i_3\ne 0\}}\left(
\frac{1}{2}\int\limits_t^T\psi_1(t_1)\int\limits_{t_1}^{T}
\psi_2(t_2)\psi_3(t_2)dt_2d{\bf w}_{t_1}^{(i_1)}-\Biggr.
\Biggl.
\sum_{j_1=0}^{p}\sum_{j_3=0}^{p}
C_{j_3 j_3 j_1}\zeta_{j_1}^{(i_1)}\right)-
$$
$$
-{\bf 1}_{\{i_1=i_3\ne 0\}}
\sum_{j_1=0}^{p}\sum_{j_3=0}^{p}
C_{j_1 j_3 j_1}\zeta_{j_3}^{(i_2)}.
$$

\vspace{2mm}

From the proof of Theorem 2.8 we obtain

\vspace{-4mm}
$$
{\sf M}\left\{\left(R_{T,t}^{(2)ppp}\right)^2\right\}\le
3\left({\bf 1}_{\{i_1=i_2\ne 0\}}{\sf M}\left\{\Biggl(
\frac{1}{2}\int\limits_t^T\psi_3(t_3)\int\limits_t^{t_3}
\psi_1(t_2)\psi_2(t_2)dt_2d{\bf w}_{t_3}^{(i_3)}-\Biggr.\right.\right.
$$
$$
\Biggl.\left.-
\sum_{j_1=0}^{p}\sum_{j_3=0}^{p}
C_{j_3 j_1 j_1}\zeta_{j_3}^{(i_3)}\Biggr)^2\right\}
+{\bf 1}_{\{i_1=i_3\ne 0\}}
{\sf M}\left\{\left(\sum_{j_1=0}^{p}\sum_{j_3=0}^{p}
C_{j_1 j_3 j_1}\zeta_{j_3}^{(i_2)}\right)^2\right\}+
$$
$$
+{\bf 1}_{\{i_2=i_3\ne 0\}}{\sf M}\left\{\Biggl(
\frac{1}{2}\int\limits_t^T\psi_1(t_1)\int\limits_{t_1}^{T}
\psi_2(t_2)\psi_3(t_2)dt_2d{\bf w}_{t_1}^{(i_1)}-\Biggr.\right.
$$
\begin{equation}
\label{newbegin100}
\Biggl.\left.\left.-
\sum_{j_1=0}^{p}\sum_{j_3=0}^{p}
C_{j_3 j_3 j_1}\zeta_{j_1}^{(i_1)}\Biggr)^2\right\}\right)\ \to 0\
\end{equation}

\vspace{2mm}
\noindent
if $p \to \infty$. From (\ref{newbegin97})--(\ref{newbegin100})
we obtain the expansion (\ref{feto19000a}). 
Theorem 2.8 is proved.

\subsection{Another Proof of Theorem 2.9}

Let us consider (\ref{proof1ggg}) for $k=4,$
$p_1=\ldots=p_4=p$, and $\psi_1(s),\ldots,\psi_4(s)\equiv 1$

\vspace{-2mm}
$$
{\int\limits_t^{*}}^T
{\int\limits_t^{*}}^{t_4}
{\int\limits_t^{*}}^{t_3}
{\int\limits_t^{*}}^{t_2}
d{\bf w}_{t_1}^{(i_1)}
d{\bf w}_{t_2}^{(i_2)}d{\bf w}_{t_3}^{(i_3)}d{\bf w}_{t_4}^{(i_4)}
=
$$
\begin{equation}
\label{proof2}
~~~~~ =
\sum_{j_1=0}^{p}\sum_{j_2=0}^{p}\sum_{j_3=0}^{p}\sum_{j_4=0}^{p}
C_{j_4 j_3 j_2 j_1}
\zeta_{j_1}^{(i_1)}\zeta_{j_2}^{(i_2)}\zeta_{j_3}^{(i_3)}\zeta_{j_4}^{(i_4)}+
J[R_{pppp}]_{T,t}^{(4)}\ \ \ \hbox{w. p. 1},
\end{equation}

\vspace{4mm}
\noindent
where 

\vspace{-5mm}
$$
J[R_{pppp}]_{T,t}^{(4)}=
$$

\vspace{-4mm}
$$
=
\hbox{\vtop{\offinterlineskip\halign{
\hfil#\hfil\cr
{\rm l.i.m.}\cr
$\stackrel{}{{}_{N\to \infty}}$\cr
}} }\sum_{l_4=0}^{N-1}\sum_{l_3=0}^{N-1}\sum_{l_2=0}^{N-1}
\sum_{l_1=0}^{N-1}
R_{pppp}(\tau_{l_1},\tau_{l_2},\tau_{l_3},\tau_{l_4})
\Delta{\bf w}_{\tau_{l_1}}^{(i_1)}
\Delta{\bf w}_{\tau_{l_2}}^{(i_2)}
\Delta{\bf w}_{\tau_{l_3}}^{(i_3)}
\Delta{\bf w}_{\tau_{l_4}}^{(i_4)},
$$

\newpage
\noindent
$$
R_{pppp}(t_1,t_2,t_3,t_4)
\stackrel{\small{\sf def}}{=}
K^{*}(t_1,t_2,t_3,t_4)-
$$

\vspace{-5mm}
\begin{equation}
\label{que1}
~~~~~~ -
\sum_{j_1=0}^{p}\sum_{j_2=0}^{p}
\sum_{j_3=0}^{p}\sum_{j_4=0}^{p}
C_{j_4j_3j_2 j_1}\phi_{j_1}(t_1)\phi_{j_2}(t_2)\phi_{j_3}(t_3)\phi_{j_4}(t_4),
\end{equation}

\vspace{2mm}
$$
K^{*}(t_1,t_2,t_3,t_4)\stackrel{\small{\sf def}}{=}
\prod_{l=1}^3\Biggl({\bf 1}_{\{t_l<t_{l+1}\}}+
\frac{1}{2}{\bf 1}_{\{t_l=t_{l+1}\}}\Biggr)=
$$
$$
={\bf 1}_{\{t_1<t_2<t_3<t_4\}}+
\frac{1}{2}{\bf 1}_{\{t_1=t_2<t_3<t_4\}}+
\frac{1}{2}{\bf 1}_{\{t_1<t_2=t_3<t_4\}}+
$$
$$
+
\frac{1}{4}{\bf 1}_{\{t_1=t_2=t_3<t_4\}}+
\frac{1}{2}{\bf 1}_{\{t_1<t_2<t_3=t_4\}}+
\frac{1}{4}{\bf 1}_{\{t_1=t_2<t_3=t_4\}}+
$$
$$
+
\frac{1}{4}{\bf 1}_{\{t_1<t_2=t_3=t_4\}}+
\frac{1}{8}{\bf 1}_{\{t_1=t_2=t_3=t_4\}}.
$$

\vspace{5mm}

We have
\begin{equation}
\label{klor1}
J[R_{pppp}]_{T,t}^{(4)}=\sum\limits_{i=0}^7 R_{T,t}^{(i)pppp}\ \ \ 
\hbox{w.~p.~1},
\end{equation}

\vspace{1mm}
\noindent
where
$$
R_{T,t}^{(0)pppp}=
\hbox{\vtop{\offinterlineskip\halign{
\hfil#\hfil\cr
{\rm l.i.m.}\cr
$\stackrel{}{{}_{N\to \infty}}$\cr
}} }\sum_{l_4=0}^{N-1}\sum_{l_3=0}^{l_4-1}\sum_{l_2=0}^{l_3-1}
\sum_{l_1=0}^{l_2-1}
\sum\limits_{(l_1,l_2,l_3,l_4)}
\Biggl(
R_{pppp}(\tau_{l_1},\tau_{l_2},\tau_{l_3},\tau_{l_4})\times\Biggr.
$$
$$
\Biggl.
\times
\Delta{\bf w}_{\tau_{l_1}}^{(i_1)}
\Delta{\bf w}_{\tau_{l_2}}^{(i_2)}
\Delta{\bf w}_{\tau_{l_3}}^{(i_3)}
\Delta{\bf w}_{\tau_{l_4}}^{(i_4)}\Biggr),
$$

\noindent
where permutations
$(l_1,l_2,l_3,l_4)$ 
when summing
are performed only in the expression, which is enclosed in parentheses,
$$
R_{T,t}^{(1)pppp}={\bf 1}_{\{i_1=i_2\ne 0\}}
\hbox{\vtop{\offinterlineskip\halign{
\hfil#\hfil\cr
{\rm l.i.m.}\cr
$\stackrel{}{{}_{N\to \infty}}$\cr
}} }
\sum_{\stackrel{l_4,l_3,l_1=0}{{}_{l_1\ne l_3, l_1\ne l_4,
l_3\ne l_4}}}^{N-1}
R_{pppp}(\tau_{l_1},\tau_{l_1},\tau_{l_3},\tau_{l_4})
\Delta\tau_{l_1}
\Delta{\bf w}_{\tau_{l_3}}^{(i_3)}
\Delta{\bf w}_{\tau_{l_4}}^{(i_4)},
$$
$$
R_{T,t}^{(2)pppp}={\bf 1}_{\{i_1=i_3\ne 0\}}
\hbox{\vtop{\offinterlineskip\halign{
\hfil#\hfil\cr
{\rm l.i.m.}\cr
$\stackrel{}{{}_{N\to \infty}}$\cr
}} }
\sum_{\stackrel{l_4,l_2,l_1=0}{{}_{l_1\ne l_2, l_1\ne l_4,
l_2\ne l_4}}}^{N-1}
R_{pppp}(\tau_{l_1},\tau_{l_2},\tau_{l_1},\tau_{l_4})
\Delta\tau_{l_1}
\Delta{\bf w}_{\tau_{l_2}}^{(i_2)}
\Delta{\bf w}_{\tau_{l_4}}^{(i_4)},
$$
$$
R_{T,t}^{(3)pppp}={\bf 1}_{\{i_1=i_4\ne 0\}}
\hbox{\vtop{\offinterlineskip\halign{
\hfil#\hfil\cr
{\rm l.i.m.}\cr
$\stackrel{}{{}_{N\to \infty}}$\cr
}} }
\sum_{\stackrel{l_3,l_2,l_1=0}{{}_{l_1\ne l_2, l_1\ne l_3,
l_2\ne l_3}}}^{N-1}
R_{pppp}(\tau_{l_1},\tau_{l_2},\tau_{l_3},\tau_{l_1})
\Delta\tau_{l_1}
\Delta{\bf w}_{\tau_{l_2}}^{(i_2)}
\Delta{\bf w}_{\tau_{l_3}}^{(i_3)},
$$
$$
R_{T,t}^{(4)pppp}={\bf 1}_{\{i_2=i_3\ne 0\}}
\hbox{\vtop{\offinterlineskip\halign{
\hfil#\hfil\cr
{\rm l.i.m.}\cr
$\stackrel{}{{}_{N\to \infty}}$\cr
}} }
\sum_{\stackrel{l_4,l_2,l_1=0}{{}_{l_1\ne l_2, l_1\ne l_4,
l_2\ne l_4}}}^{N-1}
R_{pppp}(\tau_{l_1},\tau_{l_2},\tau_{l_2},\tau_{l_4})
\Delta{\bf w}_{\tau_{l_1}}^{(i_1)}
\Delta\tau_{l_2}
\Delta{\bf w}_{\tau_{l_4}}^{(i_4)},
$$
$$
R_{T,t}^{(5)pppp}={\bf 1}_{\{i_2=i_4\ne 0\}}
\hbox{\vtop{\offinterlineskip\halign{
\hfil#\hfil\cr
{\rm l.i.m.}\cr
$\stackrel{}{{}_{N\to \infty}}$\cr
}} }
\sum_{\stackrel{l_3,l_2,l_1=0}{{}_{l_1\ne l_2, l_1\ne l_3,
l_2\ne l_3}}}^{N-1}
R_{pppp}(\tau_{l_1},\tau_{l_2},\tau_{l_3},\tau_{l_2})
\Delta{\bf w}_{\tau_{l_1}}^{(i_1)}
\Delta\tau_{l_2}
\Delta{\bf w}_{\tau_{l_3}}^{(i_3)},
$$
$$
R_{T,t}^{(6)pppp}={\bf 1}_{\{i_3=i_4\ne 0\}}
\hbox{\vtop{\offinterlineskip\halign{
\hfil#\hfil\cr
{\rm l.i.m.}\cr
$\stackrel{}{{}_{N\to \infty}}$\cr
}} }
\sum_{\stackrel{l_3,l_2,l_1=0}{{}_{l_1\ne l_2, l_1\ne l_3,
l_2\ne l_3}}}^{N-1}
R_{pppp}(\tau_{l_1},\tau_{l_2},\tau_{l_3},\tau_{l_3})
\Delta{\bf w}_{\tau_{l_1}}^{(i_1)}
\Delta{\bf w}_{\tau_{l_2}}^{(i_2)}
\Delta\tau_{l_3},
$$
$$
R_{T,t}^{(7)pppp}={\bf 1}_{\{i_1=i_2\ne 0\}}{\bf 1}_{\{i_3=i_4\ne 0\}}
\hbox{\vtop{\offinterlineskip\halign{
\hfil#\hfil\cr
{\rm lim}\cr
$\stackrel{}{{}_{N\to \infty}}$\cr
}} }
\sum_{\stackrel{l_4,l_2=0}{{}_{l_2\ne l_4}}}^{N-1}
R_{pppp}(\tau_{l_2},\tau_{l_2},\tau_{l_4},\tau_{l_4})
\Delta\tau_{l_2}
\Delta\tau_{l_4}+
$$
$$
+{\bf 1}_{\{i_1=i_3\ne 0\}}{\bf 1}_{\{i_2=i_4\ne 0\}}
\hbox{\vtop{\offinterlineskip\halign{
\hfil#\hfil\cr
{\rm lim}\cr
$\stackrel{}{{}_{N\to \infty}}$\cr
}} }
\sum_{\stackrel{l_4,l_2=0}{{}_{l_2\ne l_4}}}^{N-1}
R_{pppp}(\tau_{l_2},\tau_{l_4},\tau_{l_2},\tau_{l_4})
\Delta\tau_{l_2}
\Delta\tau_{l_4}+
$$
$$
+{\bf 1}_{\{i_1=i_4\ne 0\}}{\bf 1}_{\{i_2=i_3\ne 0\}}
\hbox{\vtop{\offinterlineskip\halign{
\hfil#\hfil\cr
{\rm lim}\cr
$\stackrel{}{{}_{N\to \infty}}$\cr
}} }
\sum_{\stackrel{l_4,l_2=0}{{}_{l_2\ne l_4}}}^{N-1}
R_{pppp}(\tau_{l_2},\tau_{l_4},\tau_{l_4},\tau_{l_2})
\Delta\tau_{l_2}
\Delta\tau_{l_4}.
$$

\vspace{2mm}

From (\ref{proof2}) and (\ref{klor1}) it follows that Theorem 2.9
will be proved if 
$$
\hbox{\vtop{\offinterlineskip\halign{
\hfil#\hfil\cr
{\rm lim}\cr
$\stackrel{}{{}_{p\to \infty}}$\cr
}} }
{\sf M}\left\{\left(R_{T,t}^{(i)pppp}\right)^2\right\}=0,\ \ \ 
i=0, 1, \ldots,7.
$$

We have (see (\ref{30.52}), (\ref{pobeda}))
$$
R_{T,t}^{(0)pppp}=
\int\limits_t^T
\int\limits_t^{t_4}\int\limits_t^{t_3}\int\limits_t^{t_2}
\sum\limits_{(t_1,t_2,t_3,t_4)}
\Biggl(
R_{pppp}(t_1,t_2,t_3,t_4)
d{\bf w}_{t_1}^{(i_1)}
d{\bf w}_{t_2}^{(i_2)}
d{\bf w}_{t_3}^{(i_3)}
d{\bf w}_{t_4}^{(i_4)}\Biggr),
$$
where permutations
$(t_1,t_2,t_3,t_4)$ 
when summing
are performed only in the expression, which is enclosed in parentheses.

From the other hand (see (\ref{pobeda}), (\ref{s2s}))
$$
R_{T,t}^{(0)pppp}=
\sum\limits_{(t_1,t_2,t_3,t_4)}\int\limits_t^T
\int\limits_t^{t_4}\int\limits_t^{t_3}\int\limits_t^{t_2}
R_{pppp}(t_1,t_2,t_3,t_4)
d{\bf w}_{t_1}^{(i_1)}
d{\bf w}_{t_2}^{(i_2)}
d{\bf w}_{t_3}^{(i_3)}
d{\bf w}_{t_4}^{(i_4)},
$$

\noindent
where permutations
$(t_1,t_2,t_3,t_4)$ when summing
are performed only in 
the values
$d{\bf w}_{t_1}^{(i_1)}d{\bf w}_{t_2}^{(i_2)}d{\bf w}_{t_3}^{(i_3)}
d{\bf w}_{t_4}^{(i_4)}$. At the same time the indices near upper 
limits of integration in the iterated stochastic integrals are changed 
correspondently and if $t_r$ swapped with $t_q$ in the  
permutation $(t_1,t_2,t_3,t_4)$, then $i_r$ swapped with $i_q$ in 
the permutation $(i_1,i_2,i_3,i_4)$.

So, we obtain
$$
{\sf M}\left\{\left(R_{T,t}^{(0)pppp}\right)^2\right\}\le
24\sum\limits_{(t_1,t_2,t_3,t_4)}\int\limits_t^T
\int\limits_t^{t_4}\int\limits_t^{t_3}\int\limits_t^{t_2}
\left(R_{pppp}(t_1,t_2,t_3,t_4)\right)^2dt_1dt_2dt_3dt_4=
$$
$$
= 24\int\limits_{[t, T]^4}
\left(R_{pppp}(t_1,t_2,t_3,t_4))\right)^2
dt_1dt_2dt_3dt_4\ \to 0
$$

\noindent
if $p\to\infty,$ $K^{*}(t_1,t_2,t_3,t_4)\in L_2([t, T]^4)$ (see (\ref{que1})).

\vspace{2mm}

Let us consider $R_{T,t}^{(1)pppp}$
$$
R_{T,t}^{(1)pppp}={\bf 1}_{\{i_1=i_2\ne 0\}}
\hbox{\vtop{\offinterlineskip\halign{
\hfil#\hfil\cr
{\rm l.i.m.}\cr
$\stackrel{}{{}_{N\to \infty}}$\cr
}} }
\sum_{\stackrel{l_4,l_3,l_1=0}{{}_{l_1\ne l_3, l_1\ne l_4,
l_3\ne l_4}}}^{N-1}
R_{pppp}(\tau_{l_1},\tau_{l_1},\tau_{l_3},\tau_{l_4})
\Delta\tau_{l_1}
\Delta{\bf w}_{\tau_{l_3}}^{(i_3)}
\Delta{\bf w}_{\tau_{l_4}}^{(i_4)}=
$$
$$
={\bf 1}_{\{i_1=i_2\ne 0\}}
\hbox{\vtop{\offinterlineskip\halign{
\hfil#\hfil\cr
{\rm l.i.m.}\cr
$\stackrel{}{{}_{N\to \infty}}$\cr
}} }
\sum_{\stackrel{l_4,l_3,l_1=0}{{}_{l_3\ne l_4}}}^{N-1}
R_{pppp}(\tau_{l_1},\tau_{l_1},\tau_{l_3},\tau_{l_4})
\Delta\tau_{l_1}
\Delta{\bf w}_{\tau_{l_3}}^{(i_3)}
\Delta{\bf w}_{\tau_{l_4}}^{(i_4)}=
$$
$$
={\bf 1}_{\{i_1=i_2\ne 0\}}
\hbox{\vtop{\offinterlineskip\halign{
\hfil#\hfil\cr
{\rm l.i.m.}\cr
$\stackrel{}{{}_{N\to \infty}}$\cr
}} }
\sum_{\stackrel{l_4,l_3,l_1=0}{{}_{l_3\ne l_4}}}^{N-1}
\Biggl(
\frac{1}{2}{\bf 1}_{\{\tau_{l_1}<\tau_{l_3}<\tau_{l_4}\}}
+ 
\Biggr.
$$
$$
+\frac{1}{4}{\bf 1}_{\{\tau_{l_1}=\tau_{l_3}<\tau_{l_4}\}}+
\frac{1}{4}{\bf 1}_{\{\tau_{l_1}<\tau_{l_3}=\tau_{l_4}\}}+
\frac{1}{8}{\bf 1}_{\{\tau_{l_1}=\tau_{l_3}=\tau_{l_4}\}}-
$$
$$
\Biggl.
-\sum\limits_{j_4,j_3,j_2,j_1=0}^p C_{j_4j_3j_2j_1}
\phi_{j_1}(\tau_{l_1})
\phi_{j_2}(\tau_{l_1})
\phi_{j_3}(\tau_{l_3})
\phi_{j_4}(\tau_{l_4})\Biggr)
\Delta\tau_{l_1}
\Delta{\bf w}_{\tau_{l_3}}^{(i_3)}
\Delta{\bf w}_{\tau_{l_4}}^{(i_4)}=
$$
$$
={\bf 1}_{\{i_1=i_2\ne 0\}}
\hbox{\vtop{\offinterlineskip\halign{
\hfil#\hfil\cr
{\rm l.i.m.}\cr
$\stackrel{}{{}_{N\to \infty}}$\cr
}} }
\sum_{\stackrel{l_4,l_3,l_1=0}{{}_{l_3\ne l_4}}}^{N-1}
\Biggl(
\frac{1}{2}{\bf 1}_{\{\tau_{l_1}<\tau_{l_3}<\tau_{l_4}\}}
- 
\Biggr.
$$
$$
\Biggl.
-\sum\limits_{j_4,j_3,j_2,j_1=0}^p C_{j_4j_3j_2j_1}
\phi_{j_1}(\tau_{l_1})
\phi_{j_2}(\tau_{l_1})
\phi_{j_3}(\tau_{l_3})
\phi_{j_4}(\tau_{l_4})\Biggr)
\Delta\tau_{l_1}
\Delta{\bf w}_{\tau_{l_3}}^{(i_3)}
\Delta{\bf w}_{\tau_{l_4}}^{(i_4)}=
$$
$$
={\bf 1}_{\{i_1=i_2\ne 0\}}
\hbox{\vtop{\offinterlineskip\halign{
\hfil#\hfil\cr
{\rm l.i.m.}\cr
$\stackrel{}{{}_{N\to \infty}}$\cr
}} }
\sum\limits_{l_4=0}^{N-1}\sum\limits_{l_3=0}^{N-1}
\sum\limits_{l_1=0}^{N-1}
\Biggl(
\frac{1}{2}{\bf 1}_{\{\tau_{l_1}<\tau_{l_3}<\tau_{l_4}\}}
- 
\Biggr.
$$
$$
\Biggl.
-\sum\limits_{j_4,j_3,j_2,j_1=0}^p C_{j_4j_3j_2j_1}
\phi_{j_1}(\tau_{l_1})
\phi_{j_2}(\tau_{l_1})
\phi_{j_3}(\tau_{l_3})
\phi_{j_4}(\tau_{l_4})\Biggr)
\Delta\tau_{l_1}
\Delta{\bf w}_{\tau_{l_3}}^{(i_3)}
\Delta{\bf w}_{\tau_{l_4}}^{(i_4)}-
$$
$$
-{\bf 1}_{\{i_1=i_2\ne 0\}}{\bf 1}_{\{i_3=i_4\ne 0\}}
\hbox{\vtop{\offinterlineskip\halign{
\hfil#\hfil\cr
{\rm l.i.m.}\cr
$\stackrel{}{{}_{N\to \infty}}$\cr
}} }
\sum\limits_{l_4=0}^{N-1}
\sum\limits_{l_1=0}^{N-1}
\Biggl(
0
- 
\Biggr.
$$
$$
\Biggl.
-\sum\limits_{j_4,j_3,j_2,j_1=0}^p C_{j_4j_3j_2j_1}
\phi_{j_1}(\tau_{l_1})
\phi_{j_2}(\tau_{l_1})
\phi_{j_3}(\tau_{l_4})
\phi_{j_4}(\tau_{l_4})\Biggr)
\Delta\tau_{l_1}
\Delta\tau_{l_4}=
$$
$$
=
{\bf 1}_{\{i_1=i_2\ne 0\}}\left(\frac{1}{2}
\int\limits_t^T
\int\limits_t^{t_4}
\int\limits_t^{t_3}dt_1d{\bf w}_{t_3}^{(i_3)}
d{\bf w}_{t_4}^{(i_4)} - 
\sum\limits_{j_4,j_3,j_1=0}^p C_{j_4j_3j_1j_1}
\zeta_{j_3}^{(i_3)}
\zeta_{j_4}^{(i_4)}\right)+
$$
$$
+{\bf 1}_{\{i_1=i_2\ne 0\}}{\bf 1}_{\{i_3=i_4\ne 0\}}
\sum\limits_{j_4,j_1=0}^p C_{j_4j_4j_1j_1}\ \ \ \hbox{w. p. 1.}
$$

\vspace{2mm}

When proving Theorem 2.9
we have proved that
$$
\hbox{\vtop{\offinterlineskip\halign{
\hfil#\hfil\cr
{\rm lim}\cr
$\stackrel{}{{}_{p\to \infty}}$\cr
}} }
\sum\limits_{j_4,j_1=0}^{p} C_{j_4j_4j_1j_1}=
\frac{1}{4}
\int\limits_t^T
\int\limits_t^{t_2}dt_1dt_2,
$$
$$
\hbox{\vtop{\offinterlineskip\halign{
\hfil#\hfil\cr
{\rm l.i.m.}\cr
$\stackrel{}{{}_{p\to \infty}}$\cr
}} }
\sum\limits_{j_4,j_3,j_1=0}^p C_{j_4j_3j_1j_1}
\zeta_{j_3}^{(i_3)}
\zeta_{j_4}^{(i_4)}=
\frac{1}{2}
\int\limits_t^T
\int\limits_t^{t_4}
\int\limits_t^{t_3}dt_1d{\bf w}_{t_3}^{(i_3)}
d{\bf w}_{t_4}^{(i_4)}+
$$
$$
+{\bf 1}_{\{i_3=i_4\ne 0\}}
\frac{1}{4}
\int\limits_t^T
\int\limits_t^{t_2}dt_1dt_2\ \ \ \hbox{w.~p.~1.}
$$

\vspace{1mm}

Then
$$
\hbox{\vtop{\offinterlineskip\halign{
\hfil#\hfil\cr
{\rm lim}\cr
$\stackrel{}{{}_{p\to \infty}}$\cr
}} }
{\sf M}\left\{\left(R_{T,t}^{(1)pppp}\right)^2\right\}=0.
$$

\vspace{2mm}

Let us consider $R_{T,t}^{(2)pppp}$
$$
R_{T,t}^{(2)pppp}={\bf 1}_{\{i_1=i_3\ne 0\}}
\hbox{\vtop{\offinterlineskip\halign{
\hfil#\hfil\cr
{\rm l.i.m.}\cr
$\stackrel{}{{}_{N\to \infty}}$\cr
}} }
\sum_{\stackrel{l_4,l_2,l_1=0}{{}_{l_1\ne l_2, l_1\ne l_4,
l_2\ne l_4}}}^{N-1}
G_{pppp}(\tau_{l_1},\tau_{l_2},\tau_{l_1},\tau_{l_4})
\Delta\tau_{l_1}
\Delta{\bf w}_{\tau_{l_2}}^{(i_2)}
\Delta{\bf w}_{\tau_{l_4}}^{(i_4)}=
$$
$$
={\bf 1}_{\{i_1=i_3\ne 0\}}
\hbox{\vtop{\offinterlineskip\halign{
\hfil#\hfil\cr
{\rm l.i.m.}\cr
$\stackrel{}{{}_{N\to \infty}}$\cr
}} }
\sum_{\stackrel{l_4,l_2,l_1=0}{{}_{l_2\ne l_4}}}^{N-1}
G_{pppp}(\tau_{l_1},\tau_{l_2},\tau_{l_1},\tau_{l_4})
\Delta\tau_{l_1}
\Delta{\bf w}_{\tau_{l_2}}^{(i_2)}
\Delta{\bf w}_{\tau_{l_4}}^{(i_4)}=
$$
$$
={\bf 1}_{\{i_1=i_3\ne 0\}}
\hbox{\vtop{\offinterlineskip\halign{
\hfil#\hfil\cr
{\rm l.i.m.}\cr
$\stackrel{}{{}_{N\to \infty}}$\cr
}} }
\sum_{\stackrel{l_4,l_2,l_1=0}{{}_{l_2\ne l_4}}}^{N-1}
\Biggl(
\frac{1}{4}{\bf 1}_{\{\tau_{l_1}=\tau_{l_2}<\tau_{l_4}\}}
+ 
\Biggr.
\frac{1}{8}{\bf 1}_{\{\tau_{l_1}=\tau_{l_2}=\tau_{l_4}\}}-
$$
$$
\Biggl.
-\sum\limits_{j_4,j_3,j_2,j_1=0}^p C_{j_4j_3j_2j_1}
\phi_{j_1}(\tau_{l_1})
\phi_{j_2}(\tau_{l_2})
\phi_{j_3}(\tau_{l_1})
\phi_{j_4}(\tau_{l_4})\Biggr)
\Delta\tau_{l_1}
\Delta{\bf w}_{\tau_{l_2}}^{(i_2)}
\Delta{\bf w}_{\tau_{l_4}}^{(i_4)}=
$$
$$
={\bf 1}_{\{i_1=i_3\ne 0\}}
\hbox{\vtop{\offinterlineskip\halign{
\hfil#\hfil\cr
{\rm l.i.m.}\cr
$\stackrel{}{{}_{N\to \infty}}$\cr
}} }
\sum\limits_{l_4=0}^{N-1}\sum\limits_{l_2=0}^{N-1}
\sum\limits_{l_1=0}^{N-1}
(-1)\sum\limits_{j_4,j_3,j_2,j_1=0}^p C_{j_4j_3j_2j_1}\times
$$
$$
\times\phi_{j_1}(\tau_{l_1})
\phi_{j_2}(\tau_{l_2})
\phi_{j_3}(\tau_{l_1})
\phi_{j_4}(\tau_{l_4})
\Delta\tau_{l_1}
\Delta{\bf w}_{\tau_{l_2}}^{(i_2)}
\Delta{\bf w}_{\tau_{l_4}}^{(i_4)}-
$$
$$
-{\bf 1}_{\{i_1=i_3\ne 0\}}{\bf 1}_{\{i_2=i_4\ne 0\}}
\hbox{\vtop{\offinterlineskip\halign{
\hfil#\hfil\cr
{\rm l.i.m.}\cr
$\stackrel{}{{}_{N\to \infty}}$\cr
}} }
\sum\limits_{l_4=0}^{N-1}
\sum\limits_{l_1=0}^{N-1}
(-1)\sum\limits_{j_4,j_3,j_2,j_1=0}^p C_{j_4j_3j_2j_1}\times
$$
$$
\times
\phi_{j_1}(\tau_{l_1})
\phi_{j_2}(\tau_{l_4})
\phi_{j_3}(\tau_{l_1})
\phi_{j_4}(\tau_{l_4})
\Delta\tau_{l_1}
\Delta\tau_{l_4}=
$$
$$
=
-{\bf 1}_{\{i_1=i_3\ne 0\}}
\sum\limits_{j_4,j_2,j_1=0}^p C_{j_4j_1j_2j_1}
\zeta_{j_2}^{(i_2)}
\zeta_{j_4}^{(i_4)}+
$$
$$
+{\bf 1}_{\{i_1=i_3\ne 0\}}{\bf 1}_{\{i_2=i_4\ne 0\}}
\sum\limits_{j_4,j_1=0}^p C_{j_4j_1j_4j_1}\ \ \ \hbox{w.~p.~1.}
$$

\vspace{2mm}

When proving Theorem 2.9
we have proved that
$$
\hbox{\vtop{\offinterlineskip\halign{
\hfil#\hfil\cr
{\rm l.i.m.}\cr
$\stackrel{}{{}_{p\to \infty}}$\cr
}} }
\sum\limits_{j_4,j_2,j_1=0}^p C_{j_4j_1j_2j_1}
\zeta_{j_2}^{(i_2)}
\zeta_{j_4}^{(i_4)}=0\ \ \ \hbox{w. p. 1,}
$$
$$
\hbox{\vtop{\offinterlineskip\halign{
\hfil#\hfil\cr
{\rm lim}\cr
$\stackrel{}{{}_{p\to \infty}}$\cr
}} }
\sum\limits_{j_4,j_1=0}^p C_{j_4j_1j_4j_1}=0.
$$

\vspace{1mm}

Then
$$
\hbox{\vtop{\offinterlineskip\halign{
\hfil#\hfil\cr
{\rm lim}\cr
$\stackrel{}{{}_{p\to \infty}}$\cr
}} }
{\sf M}\left\{\left(R_{T,t}^{(2)pppp}\right)^2\right\}=0.
$$

\vspace{2mm}

Let us consider $R_{T,t}^{(3)pppp}$
$$
R_{T,t}^{(3)pppp}={\bf 1}_{\{i_1=i_4\ne 0\}}
\hbox{\vtop{\offinterlineskip\halign{
\hfil#\hfil\cr
{\rm l.i.m.}\cr
$\stackrel{}{{}_{N\to \infty}}$\cr
}} }
\sum_{\stackrel{l_3,l_2,l_1=0}{{}_{l_1\ne l_2, l_1\ne l_3,
l_2\ne l_3}}}^{N-1}
G_{pppp}(\tau_{l_1},\tau_{l_2},\tau_{l_3},\tau_{l_1})
\Delta\tau_{l_1}
\Delta{\bf w}_{\tau_{l_2}}^{(i_2)}
\Delta{\bf w}_{\tau_{l_3}}^{(i_3)}=
$$
$$
={\bf 1}_{\{i_1=i_4\ne 0\}}
\hbox{\vtop{\offinterlineskip\halign{
\hfil#\hfil\cr
{\rm l.i.m.}\cr
$\stackrel{}{{}_{N\to \infty}}$\cr
}} }
\sum_{\stackrel{l_3,l_2,l_1=0}{{}_{l_2\ne l_3}}}^{N-1}
G_{pppp}(\tau_{l_1},\tau_{l_2},\tau_{l_3},\tau_{l_1})
\Delta\tau_{l_1}
\Delta{\bf w}_{\tau_{l_2}}^{(i_2)}
\Delta{\bf w}_{\tau_{l_3}}^{(i_3)}=
$$
$$
={\bf 1}_{\{i_1=i_4\ne 0\}}
\hbox{\vtop{\offinterlineskip\halign{
\hfil#\hfil\cr
{\rm l.i.m.}\cr
$\stackrel{}{{}_{N\to \infty}}$\cr
}} }
\sum_{\stackrel{l_3,l_2,l_1=0}{{}_{l_2\ne l_3}}}^{N-1}
\Biggl(
\frac{1}{8}{\bf 1}_{\{\tau_{l_1}=\tau_{l_2}=\tau_{l_3}\}}-
$$
$$
\Biggl.
-\sum\limits_{j_4,j_3,j_2,j_1=0}^p C_{j_4j_3j_2j_1}
\phi_{j_1}(\tau_{l_1})
\phi_{j_2}(\tau_{l_2})
\phi_{j_3}(\tau_{l_3})
\phi_{j_4}(\tau_{l_1})\Biggr)
\Delta\tau_{l_1}
\Delta{\bf w}_{\tau_{l_2}}^{(i_2)}
\Delta{\bf w}_{\tau_{l_3}}^{(i_3)}=
$$
$$
={\bf 1}_{\{i_1=i_4\ne 0\}}
\hbox{\vtop{\offinterlineskip\halign{
\hfil#\hfil\cr
{\rm l.i.m.}\cr
$\stackrel{}{{}_{N\to \infty}}$\cr
}} }
\sum\limits_{l_3=0}^{N-1}\sum\limits_{l_2=0}^{N-1}
\sum\limits_{l_1=0}^{N-1}
(-1)\sum\limits_{j_4,j_3,j_2,j_1=0}^p C_{j_4j_3j_2j_1}\times
$$
$$
\times\phi_{j_1}(\tau_{l_1})
\phi_{j_2}(\tau_{l_2})
\phi_{j_3}(\tau_{l_3})
\phi_{j_4}(\tau_{l_1})
\Delta\tau_{l_1}
\Delta{\bf w}_{\tau_{l_2}}^{(i_2)}
\Delta{\bf w}_{\tau_{l_3}}^{(i_3)}-
$$
$$
-{\bf 1}_{\{i_1=i_4\ne 0\}}{\bf 1}_{\{i_2=i_3\ne 0\}}
\hbox{\vtop{\offinterlineskip\halign{
\hfil#\hfil\cr
{\rm l.i.m.}\cr
$\stackrel{}{{}_{N\to \infty}}$\cr
}} }
\sum\limits_{l_3=0}^{N-1}
\sum\limits_{l_1=0}^{N-1}
(-1)\sum\limits_{j_4,j_3,j_2,j_1=0}^p C_{j_4j_3j_2j_1}\times
$$
$$
\times
\phi_{j_1}(\tau_{l_1})
\phi_{j_2}(\tau_{l_3})
\phi_{j_3}(\tau_{l_3})
\phi_{j_4}(\tau_{l_1})
\Delta\tau_{l_1}
\Delta\tau_{l_3}=
$$
$$
=
-{\bf 1}_{\{i_1=i_4\ne 0\}}
\sum\limits_{j_4,j_3,j_2=0}^p C_{j_4j_3j_2j_4}
\zeta_{j_2}^{(i_2)}
\zeta_{j_3}^{(i_3)}+
$$
$$
+{\bf 1}_{\{i_1=i_4\ne 0\}}{\bf 1}_{\{i_2=i_3\ne 0\}}
\sum\limits_{j_4,j_2=0}^p C_{j_4j_2j_2j_4}\ \ \ \hbox{w. p. 1.}
$$

\vspace{2mm}

When proving Theorem 2.9
we have proved that
$$
\hbox{\vtop{\offinterlineskip\halign{
\hfil#\hfil\cr
{\rm l.i.m.}\cr
$\stackrel{}{{}_{p\to \infty}}$\cr
}} }
\sum\limits_{j_4,j_3,j_2=0}^p C_{j_4j_3j_2j_4}
\zeta_{j_2}^{(i_2)}
\zeta_{j_3}^{(i_3)}=0\ \ \ \hbox{w. p. 1,}
$$
$$
\hbox{\vtop{\offinterlineskip\halign{
\hfil#\hfil\cr
{\rm lim}\cr
$\stackrel{}{{}_{p\to \infty}}$\cr
}} }
\sum\limits_{j_4,j_2=0}^p C_{j_4j_2j_2j_4}=0.
$$

\vspace{1mm}

Then
$$
\hbox{\vtop{\offinterlineskip\halign{
\hfil#\hfil\cr
{\rm lim}\cr
$\stackrel{}{{}_{p\to \infty}}$\cr
}} }
{\sf M}\left\{\left(R_{T,t}^{(3)pppp}\right)^2\right\}=0.
$$

\vspace{2mm}

Let us consider $R_{T,t}^{(4)pppp}$
$$
R_{T,t}^{(4)pppp}={\bf 1}_{\{i_2=i_3\ne 0\}}
\hbox{\vtop{\offinterlineskip\halign{
\hfil#\hfil\cr
{\rm l.i.m.}\cr
$\stackrel{}{{}_{N\to \infty}}$\cr
}} }
\sum_{\stackrel{l_4,l_2,l_1=0}{{}_{l_1\ne l_2, l_1\ne l_4,
l_2\ne l_4}}}^{N-1}
G_{pppp}(\tau_{l_1},\tau_{l_2},\tau_{l_2},\tau_{l_4})
\Delta{\bf w}_{\tau_{l_1}}^{(i_1)}
\Delta\tau_{l_2}
\Delta{\bf w}_{\tau_{l_4}}^{(i_4)}=
$$
$$
={\bf 1}_{\{i_2=i_3\ne 0\}}
\hbox{\vtop{\offinterlineskip\halign{
\hfil#\hfil\cr
{\rm l.i.m.}\cr
$\stackrel{}{{}_{N\to \infty}}$\cr
}} }
\sum_{\stackrel{l_4,l_2,l_1=0}{{}_{l_1\ne l_4}}}^{N-1}
G_{pppp}(\tau_{l_1},\tau_{l_2},\tau_{l_2},\tau_{l_4})
\Delta{\bf w}_{\tau_{l_1}}^{(i_1)}
\Delta\tau_{l_2}
\Delta{\bf w}_{\tau_{l_4}}^{(i_4)}=
$$
$$
={\bf 1}_{\{i_2=i_3\ne 0\}}
\hbox{\vtop{\offinterlineskip\halign{
\hfil#\hfil\cr
{\rm l.i.m.}\cr
$\stackrel{}{{}_{N\to \infty}}$\cr
}} }
\sum_{\stackrel{l_4,l_2,l_1=0}{{}_{l_1\ne l_4}}}^{N-1}
\Biggl(
\frac{1}{2}{\bf 1}_{\{\tau_{l_1}<\tau_{l_2}<\tau_{l_4}\}}
+ 
\Biggr.
$$
$$
+\frac{1}{4}{\bf 1}_{\{\tau_{l_1}=\tau_{l_2}<\tau_{l_4}\}}+
\frac{1}{4}{\bf 1}_{\{\tau_{l_1}<\tau_{l_2}=\tau_{l_4}\}}+
\frac{1}{8}{\bf 1}_{\{\tau_{l_1}=\tau_{l_2}=\tau_{l_4}\}}-
$$
$$
\Biggl.
-\sum\limits_{j_4,j_3,j_2,j_1=0}^p C_{j_4j_3j_2j_1}
\phi_{j_1}(\tau_{l_1})
\phi_{j_2}(\tau_{l_2})
\phi_{j_3}(\tau_{l_2})
\phi_{j_4}(\tau_{l_4})\Biggr)
\Delta{\bf w}_{\tau_{l_1}}^{(i_1)}
\Delta\tau_{l_2}
\Delta{\bf w}_{\tau_{l_4}}^{(i_4)}=
$$
$$
={\bf 1}_{\{i_2=i_3\ne 0\}}
\hbox{\vtop{\offinterlineskip\halign{
\hfil#\hfil\cr
{\rm l.i.m.}\cr
$\stackrel{}{{}_{N\to \infty}}$\cr
}} }
\sum_{\stackrel{l_4,l_2,l_1=0}{{}_{l_1\ne l_4}}}^{N-1}
\Biggl(
\frac{1}{2}{\bf 1}_{\{\tau_{l_1}<\tau_{l_2}<\tau_{l_4}\}}
- 
\Biggr.
$$
$$
\Biggl.
-\sum\limits_{j_4,j_3,j_2,j_1=0}^p C_{j_4j_3j_2j_1}
\phi_{j_1}(\tau_{l_1})
\phi_{j_2}(\tau_{l_2})
\phi_{j_3}(\tau_{l_2})
\phi_{j_4}(\tau_{l_4})\Biggr)
\Delta{\bf w}_{\tau_{l_1}}^{(i_1)}
\Delta\tau_{l_2}
\Delta{\bf w}_{\tau_{l_4}}^{(i_4)}=
$$
$$
={\bf 1}_{\{i_2=i_3\ne 0\}}
\hbox{\vtop{\offinterlineskip\halign{
\hfil#\hfil\cr
{\rm l.i.m.}\cr
$\stackrel{}{{}_{N\to \infty}}$\cr
}} }
\sum\limits_{l_4=0}^{N-1}\sum\limits_{l_2=0}^{N-1}
\sum\limits_{l_1=0}^{N-1}
\Biggl(
\frac{1}{2}{\bf 1}_{\{\tau_{l_1}<\tau_{l_2}<\tau_{l_4}\}}
- 
\Biggr.
$$
$$
\Biggl.
-\sum\limits_{j_4,j_3,j_2,j_1=0}^p C_{j_4j_3j_2j_1}
\phi_{j_1}(\tau_{l_1})
\phi_{j_2}(\tau_{l_2})
\phi_{j_3}(\tau_{l_2})
\phi_{j_4}(\tau_{l_4})\Biggr)
\Delta{\bf w}_{\tau_{l_1}}^{(i_1)}
\Delta\tau_{l_2}
\Delta{\bf w}_{\tau_{l_4}}^{(i_4)}-
$$
$$
-{\bf 1}_{\{i_2=i_3\ne 0\}}{\bf 1}_{\{i_1=i_4\ne 0\}}
\hbox{\vtop{\offinterlineskip\halign{
\hfil#\hfil\cr
{\rm l.i.m.}\cr
$\stackrel{}{{}_{N\to \infty}}$\cr
}} }
\sum\limits_{l_4=0}^{N-1}
\sum\limits_{l_2=0}^{N-1}
(-1)
\sum\limits_{j_4,j_3,j_2,j_1=0}^p C_{j_4j_3j_2j_1}\times
$$

\vspace{-2mm}
$$
\times
\phi_{j_1}(\tau_{l_4})
\phi_{j_2}(\tau_{l_2})
\phi_{j_3}(\tau_{l_2})
\phi_{j_4}(\tau_{l_4})
\Delta\tau_{l_2}
\Delta\tau_{l_4}=
$$
$$
=
{\bf 1}_{\{i_2=i_3\ne 0\}}\left(\frac{1}{2}
\int\limits_t^T
\int\limits_t^{t_4}
\int\limits_t^{t_2}d{\bf w}_{t_1}^{(i_1)}dt_2
d{\bf w}_{t_4}^{(i_4)} - 
\sum\limits_{j_4,j_2,j_1=0}^p C_{j_4j_2j_2j_1}
\zeta_{j_1}^{(i_1)}
\zeta_{j_4}^{(i_4)}\right)+
$$
$$
+{\bf 1}_{\{i_2=i_3\ne 0\}}{\bf 1}_{\{i_1=i_4\ne 0\}}
\sum\limits_{j_4,j_2=0}^p C_{j_4j_2j_2j_4}\ \ \ \hbox{w. p. 1.}
$$

\vspace{3mm}

When proving Theorem 2.9
we have proved that
$$
\hbox{\vtop{\offinterlineskip\halign{
\hfil#\hfil\cr
{\rm lim}\cr
$\stackrel{}{{}_{p\to \infty}}$\cr
}} }
\sum\limits_{j_4,j_2=0}^{p} C_{j_4j_2j_2j_4}=0,
$$
$$
\hbox{\vtop{\offinterlineskip\halign{
\hfil#\hfil\cr
{\rm l.i.m.}\cr
$\stackrel{}{{}_{p\to \infty}}$\cr
}} }
\sum\limits_{j_4,j_2,j_1=0}^p C_{j_4j_2j_2j_1}
\zeta_{j_1}^{(i_1)}
\zeta_{j_4}^{(i_4)}=
\frac{1}{2}
\int\limits_t^T
\int\limits_t^{t_4}
\int\limits_t^{t_2}d{\bf w}_{t_1}^{(i_1)}dt_2
d{\bf w}_{t_4}^{(i_4)}\ \ \ \hbox{w. p. 1.}
$$

\vspace{4mm}

Then
$$
\hbox{\vtop{\offinterlineskip\halign{
\hfil#\hfil\cr
{\rm lim}\cr
$\stackrel{}{{}_{p\to \infty}}$\cr
}} }
{\sf M}\left\{\left(R_{T,t}^{(4)pppp}\right)^2\right\}=0.
$$

\vspace{2mm}

Let us consider $R_{T,t}^{(5)pppp}$
$$
R_{T,t}^{(5)pppp}={\bf 1}_{\{i_2=i_4\ne 0\}}
\hbox{\vtop{\offinterlineskip\halign{
\hfil#\hfil\cr
{\rm l.i.m.}\cr
$\stackrel{}{{}_{N\to \infty}}$\cr
}} }
\sum_{\stackrel{l_3,l_2,l_1=0}{{}_{l_1\ne l_2, l_1\ne l_3,
l_2\ne l_3}}}^{N-1}
G_{pppp}(\tau_{l_1},\tau_{l_2},\tau_{l_3},\tau_{l_2})
\Delta{\bf w}_{\tau_{l_1}}^{(i_1)}
\Delta\tau_{l_2}
\Delta{\bf w}_{\tau_{l_3}}^{(i_3)}=
$$
$$
={\bf 1}_{\{i_2=i_4\ne 0\}}
\hbox{\vtop{\offinterlineskip\halign{
\hfil#\hfil\cr
{\rm l.i.m.}\cr
$\stackrel{}{{}_{N\to \infty}}$\cr
}} }
\sum_{\stackrel{l_3,l_2,l_1=0}{{}_{l_1\ne l_3}}}^{N-1}
G_{pppp}(\tau_{l_1},\tau_{l_2},\tau_{l_3},\tau_{l_2})
\Delta{\bf w}_{\tau_{l_1}}^{(i_1)}
\Delta\tau_{l_2}
\Delta{\bf w}_{\tau_{l_3}}^{(i_3)}=
$$
$$
={\bf 1}_{\{i_2=i_4\ne 0\}}
\hbox{\vtop{\offinterlineskip\halign{
\hfil#\hfil\cr
{\rm l.i.m.}\cr
$\stackrel{}{{}_{N\to \infty}}$\cr
}} }
\sum_{\stackrel{l_3,l_2,l_1=0}{{}_{l_1\ne l_3}}}^{N-1}
\Biggl(
\frac{1}{4}{\bf 1}_{\{\tau_{l_1}<\tau_{l_2}=\tau_{l_3}\}}
+ 
\Biggr.
\frac{1}{8}{\bf 1}_{\{\tau_{l_1}=\tau_{l_2}=\tau_{l_3}\}}-
$$
$$
\Biggl.
-\sum\limits_{j_4,j_3,j_2,j_1=0}^p C_{j_4j_3j_2j_1}
\phi_{j_1}(\tau_{l_1})
\phi_{j_2}(\tau_{l_2})
\phi_{j_3}(\tau_{l_3})
\phi_{j_4}(\tau_{l_2})\Biggr)
\Delta{\bf w}_{\tau_{l_1}}^{(i_1)}
\Delta\tau_{l_2}
\Delta{\bf w}_{\tau_{l_3}}^{(i_3)}=
$$
$$
={\bf 1}_{\{i_2=i_4\ne 0\}}
\hbox{\vtop{\offinterlineskip\halign{
\hfil#\hfil\cr
{\rm l.i.m.}\cr
$\stackrel{}{{}_{N\to \infty}}$\cr
}} }
\sum_{\stackrel{l_3,l_2,l_1=0}{{}_{l_1\ne l_3}}}^{N-1}
(-1)\sum\limits_{j_4,j_3,j_2,j_1=0}^p C_{j_4j_3j_2j_1}\times
$$
$$
\times
\phi_{j_1}(\tau_{l_1})
\phi_{j_2}(\tau_{l_2})
\phi_{j_3}(\tau_{l_3})
\phi_{j_4}(\tau_{l_2})
\Delta{\bf w}_{\tau_{l_1}}^{(i_1)}
\Delta\tau_{l_2}
\Delta{\bf w}_{\tau_{l_3}}^{(i_3)}=
$$
$$
=
-{\bf 1}_{\{i_2=i_4\ne 0\}}
\sum\limits_{j_4,j_3,j_1=0}^p C_{j_4j_3j_4j_1}
\zeta_{j_1}^{(i_1)}
\zeta_{j_3}^{(i_3)}-
$$
$$
-{\bf 1}_{\{i_2=i_4\ne 0\}}{\bf 1}_{\{i_1=i_3\ne 0\}}
\hbox{\vtop{\offinterlineskip\halign{
\hfil#\hfil\cr
{\rm l.i.m.}\cr
$\stackrel{}{{}_{N\to \infty}}$\cr
}} }
\sum\limits_{l_3=0}^{N-1}
\sum\limits_{l_2=0}^{N-1}
(-1)\sum\limits_{j_4,j_3,j_2,j_1=0}^p C_{j_4j_3j_2j_1}\times
$$
$$
\times
\phi_{j_1}(\tau_{l_3})
\phi_{j_2}(\tau_{l_2})
\phi_{j_3}(\tau_{l_3})
\phi_{j_4}(\tau_{l_2})
\Delta\tau_{l_2}
\Delta\tau_{l_3}=
$$
$$
=
-{\bf 1}_{\{i_2=i_4\ne 0\}}
\sum\limits_{j_4,j_3,j_1=0}^p C_{j_4j_3j_4j_1}
\zeta_{j_1}^{(i_1)}
\zeta_{j_3}^{(i_3)}+
$$
$$
+{\bf 1}_{\{i_2=i_4\ne 0\}}{\bf 1}_{\{i_1=i_3\ne 0\}}
\sum\limits_{j_4,j_1=0}^p C_{j_4j_1j_4j_1}\ \ \ \hbox{w. p. 1.}
$$

\vspace{2mm}

When proving Theorem 2.9
we have proved that
$$
\hbox{\vtop{\offinterlineskip\halign{
\hfil#\hfil\cr
{\rm l.i.m.}\cr
$\stackrel{}{{}_{p\to \infty}}$\cr
}} }
\sum\limits_{j_4,j_3,j_1=0}^p C_{j_4j_3j_4j_1}
\zeta_{j_1}^{(i_1)}
\zeta_{j_3}^{(i_3)}=0\ \ \ \hbox{w. p. 1,}
$$
$$
\hbox{\vtop{\offinterlineskip\halign{
\hfil#\hfil\cr
{\rm lim}\cr
$\stackrel{}{{}_{p\to \infty}}$\cr
}} }
\sum\limits_{j_4,j_1=0}^p C_{j_4j_1j_4j_1}=0.
$$

\vspace{1mm}

Then
$$
\hbox{\vtop{\offinterlineskip\halign{
\hfil#\hfil\cr
{\rm lim}\cr
$\stackrel{}{{}_{p\to \infty}}$\cr
}} }
{\sf M}\left\{\left(R_{T,t}^{(5)pppp}\right)^2\right\}=0.
$$

\vspace{2mm}

Let us consider $R_{T,t}^{(6)pppp}$
$$
R_{T,t}^{(6)pppp}={\bf 1}_{\{i_3=i_4\ne 0\}}
\hbox{\vtop{\offinterlineskip\halign{
\hfil#\hfil\cr
{\rm l.i.m.}\cr
$\stackrel{}{{}_{N\to \infty}}$\cr
}} }
\sum_{\stackrel{l_3,l_2,l_1=0}{{}_{l_1\ne l_2, l_1\ne l_3,
l_2\ne l_3}}}^{N-1}
G_{pppp}(\tau_{l_1},\tau_{l_2},\tau_{l_3},\tau_{l_3})
\Delta{\bf w}_{\tau_{l_1}}^{(i_1)}
\Delta{\bf w}_{\tau_{l_2}}^{(i_2)}
\Delta\tau_{l_3}=
$$
$$
={\bf 1}_{\{i_3=i_4\ne 0\}}
\hbox{\vtop{\offinterlineskip\halign{
\hfil#\hfil\cr
{\rm l.i.m.}\cr
$\stackrel{}{{}_{N\to \infty}}$\cr
}} }
\sum_{\stackrel{l_3,l_2,l_1=0}{{}_{l_1\ne l_2}}}^{N-1}
G_{pppp}(\tau_{l_1},\tau_{l_2},\tau_{l_3},\tau_{l_3})
\Delta{\bf w}_{\tau_{l_1}}^{(i_1)}
\Delta{\bf w}_{\tau_{l_2}}^{(i_2)}
\Delta\tau_{l_3}=
$$
$$
={\bf 1}_{\{i_3=i_4\ne 0\}}
\hbox{\vtop{\offinterlineskip\halign{
\hfil#\hfil\cr
{\rm l.i.m.}\cr
$\stackrel{}{{}_{N\to \infty}}$\cr
}} }
\sum_{\stackrel{l_3,l_2,l_1=0}{{}_{l_1\ne l_2}}}^{N-1}
\Biggl(
\frac{1}{2}{\bf 1}_{\{\tau_{l_1}<\tau_{l_2}<\tau_{l_3}\}}
+ 
\Biggr.
$$
$$
+\frac{1}{4}{\bf 1}_{\{\tau_{l_1}=\tau_{l_2}<\tau_{l_3}\}}+
\frac{1}{4}{\bf 1}_{\{\tau_{l_1}<\tau_{l_2}=\tau_{l_3}\}}+
\frac{1}{8}{\bf 1}_{\{\tau_{l_1}=\tau_{l_2}=\tau_{l_3}\}}-
$$
$$
\Biggl.
-\sum\limits_{j_4,j_3,j_2,j_1=0}^p C_{j_4j_3j_2j_1}
\phi_{j_1}(\tau_{l_1})
\phi_{j_2}(\tau_{l_2})
\phi_{j_3}(\tau_{l_3})
\phi_{j_4}(\tau_{l_3})\Biggr)
\Delta{\bf w}_{\tau_{l_1}}^{(i_1)}
\Delta{\bf w}_{\tau_{l_2}}^{(i_2)}
\Delta\tau_{l_3}=
$$
$$
={\bf 1}_{\{i_3=i_4\ne 0\}}
\hbox{\vtop{\offinterlineskip\halign{
\hfil#\hfil\cr
{\rm l.i.m.}\cr
$\stackrel{}{{}_{N\to \infty}}$\cr
}} }
\sum_{\stackrel{l_3,l_2,l_1=0}{{}_{l_1\ne l_2}}}^{N-1}
\Biggl(
\frac{1}{2}{\bf 1}_{\{\tau_{l_1}<\tau_{l_2}<\tau_{l_3}\}}
- 
\Biggr.
$$
$$
\Biggl.
-\sum\limits_{j_4,j_3,j_2,j_1=0}^p C_{j_4j_3j_2j_1}
\phi_{j_1}(\tau_{l_1})
\phi_{j_2}(\tau_{l_2})
\phi_{j_3}(\tau_{l_3})
\phi_{j_4}(\tau_{l_3})\Biggr)
\Delta{\bf w}_{\tau_{l_1}}^{(i_1)}
\Delta{\bf w}_{\tau_{l_2}}^{(i_2)}
\Delta\tau_{l_3}=
$$
$$
={\bf 1}_{\{i_3=i_4\ne 0\}}
\left(\frac{1}{2}
\int\limits_t^T
\int\limits_t^{t_3}
\int\limits_t^{t_2}d{\bf w}_{t_1}^{(i_1)}
d{\bf w}_{t_2}^{(i_2)}dt_3 - 
\sum\limits_{j_4,j_2,j_1=0}^p C_{j_4j_4j_2j_1}
\zeta_{j_1}^{(i_1)}
\zeta_{j_2}^{(i_2)}\right)-
$$
$$
-{\bf 1}_{\{i_3=i_4\ne 0\}}{\bf 1}_{\{i_1=i_2\ne 0\}}
\hbox{\vtop{\offinterlineskip\halign{
\hfil#\hfil\cr
{\rm l.i.m.}\cr
$\stackrel{}{{}_{N\to \infty}}$\cr
}} }
\sum\limits_{l_3=0}^{N-1}
\sum\limits_{l_1=0}^{N-1}
(-1)
\sum\limits_{j_4,j_3,j_2,j_1=0}^p C_{j_4j_3j_2j_1}\times
$$
$$
\times
\phi_{j_1}(\tau_{l_1})
\phi_{j_2}(\tau_{l_1})
\phi_{j_3}(\tau_{l_3})
\phi_{j_4}(\tau_{l_3})
\Delta\tau_{l_1}
\Delta\tau_{l_3}=
$$
$$
={\bf 1}_{\{i_3=i_4\ne 0\}}
\left(\frac{1}{2}
\int\limits_t^T
\int\limits_t^{t_3}
\int\limits_t^{t_2}d{\bf w}_{t_1}^{(i_1)}
d{\bf w}_{t_2}^{(i_2)}dt_3 - 
\sum\limits_{j_4,j_2,j_1=0}^p C_{j_4j_4j_2j_1}
\zeta_{j_1}^{(i_1)}
\zeta_{j_2}^{(i_2)}\right)+
$$
$$
+{\bf 1}_{\{i_1=i_2\ne 0\}}{\bf 1}_{\{i_3=i_4\ne 0\}}
\sum\limits_{j_4,j_1=0}^p C_{j_4j_4j_1j_1}=
$$
$$
={\bf 1}_{\{i_3=i_4\ne 0\}}
\left(\frac{1}{2}
\int\limits_t^T
\int\limits_t^{t_3}
\int\limits_t^{t_2}d{\bf w}_{t_1}^{(i_1)}
d{\bf w}_{t_2}^{(i_2)}dt_3 
+\frac{1}{4}{\bf 1}_{\{i_1=i_2\ne 0\}}
\int\limits_t^T\int\limits_t^{t_3}dt_1dt_3-
\right.
$$
$$
\left.
-\sum\limits_{j_4,j_2,j_1=0}^p C_{j_4j_4j_2j_1}
\zeta_{j_1}^{(i_1)}
\zeta_{j_2}^{(i_2)}\right)+
$$
$$
+{\bf 1}_{\{i_1=i_2\ne 0\}}{\bf 1}_{\{i_3=i_4\ne 0\}}
\left(\sum\limits_{j_4,j_1=0}^p C_{j_4j_4j_1j_1}-
\frac{1}{4}
\int\limits_t^T\int\limits_t^{t_3}dt_1dt_3\right)\ \ \ \hbox{w.~p.~1.}
$$

\vspace{2mm}

When proving Theorem 2.9
we have proved that
$$
\hbox{\vtop{\offinterlineskip\halign{
\hfil#\hfil\cr
{\rm lim}\cr
$\stackrel{}{{}_{p\to \infty}}$\cr
}} }
\sum\limits_{j_4,j_1=0}^{p} C_{j_4j_4j_1j_1}=
\frac{1}{4}
\int\limits_t^T\int\limits_t^{t_3}dt_1dt_3,
$$
$$
\hbox{\vtop{\offinterlineskip\halign{
\hfil#\hfil\cr
{\rm l.i.m.}\cr
$\stackrel{}{{}_{p\to \infty}}$\cr
}} }
\sum\limits_{j_4,j_2,j_1=0}^p C_{j_4j_4j_2j_1}
\zeta_{j_1}^{(i_1)}
\zeta_{j_2}^{(i_2)}=
\frac{1}{2}
\int\limits_t^T
\int\limits_t^{t_3}
\int\limits_t^{t_2}{\bf w}_{t_1}^{(i_1)}
d{\bf w}_{t_2}^{(i_2)}dt_3+
$$
$$
+
{\bf 1}_{\{i_1=i_2\ne 0\}}\frac{1}{4}
\int\limits_t^T\int\limits_t^{t_3}dt_1dt_3\ \ \ \hbox{w. p. 1.}
$$

\vspace{1mm}

Then
$$
\hbox{\vtop{\offinterlineskip\halign{
\hfil#\hfil\cr
{\rm lim}\cr
$\stackrel{}{{}_{p\to \infty}}$\cr
}} }
{\sf M}\left\{\left(R_{T,t}^{(6)pppp}\right)^2\right\}=0.
$$

\vspace{2mm}

Finally, let us consider $R_{T,t}^{(7)pppp}$
$$
R_{T,t}^{(7)pppp}={\bf 1}_{\{i_1=i_2\ne 0\}}{\bf 1}_{\{i_3=i_4\ne 0\}}
\hbox{\vtop{\offinterlineskip\halign{
\hfil#\hfil\cr
{\rm l.i.m.}\cr
$\stackrel{}{{}_{N\to \infty}}$\cr
}} }
\sum_{\stackrel{l_4,l_2=0}{{}_{l_2\ne l_4}}}^{N-1}
G_{pppp}(\tau_{l_2},\tau_{l_2},\tau_{l_4},\tau_{l_4})
\Delta\tau_{l_2}
\Delta\tau_{l_4}
+
$$
$$
+{\bf 1}_{\{i_1=i_3\ne 0\}}{\bf 1}_{\{i_2=i_4\ne 0\}}
\hbox{\vtop{\offinterlineskip\halign{
\hfil#\hfil\cr
{\rm l.i.m.}\cr
$\stackrel{}{{}_{N\to \infty}}$\cr
}} }
\sum_{\stackrel{l_4,l_2=0}{{}_{l_2\ne l_4}}}^{N-1}
G_{pppp}(\tau_{l_2},\tau_{l_4},\tau_{l_2},\tau_{l_4})
\Delta\tau_{l_2}
\Delta\tau_{l_4}+
$$
$$
+{\bf 1}_{\{i_1=i_4\ne 0\}}{\bf 1}_{\{i_2=i_3\ne 0\}}
\hbox{\vtop{\offinterlineskip\halign{
\hfil#\hfil\cr
{\rm l.i.m.}\cr
$\stackrel{}{{}_{N\to \infty}}$\cr
}} }
\sum_{\stackrel{l_4,l_2=0}{{}_{l_2\ne l_4}}}^{N-1}
G_{pppp}(\tau_{l_2},\tau_{l_4},\tau_{l_4},\tau_{l_2})
\Delta\tau_{l_2}
\Delta\tau_{l_4}=
$$
$$
={\bf 1}_{\{i_1=i_2\ne 0\}}{\bf 1}_{\{i_3=i_4\ne 0\}}
\hbox{\vtop{\offinterlineskip\halign{
\hfil#\hfil\cr
{\rm l.i.m.}\cr
$\stackrel{}{{}_{N\to \infty}}$\cr
}} }
\sum\limits_{l_4=0}^{N-1}
\sum\limits_{l_2=0}^{N-1} 
G_{pppp}(\tau_{l_2},\tau_{l_2},\tau_{l_4},\tau_{l_4})
\Delta\tau_{l_2}
\Delta\tau_{l_4}+
$$
$$
+{\bf 1}_{\{i_1=i_3\ne 0\}}{\bf 1}_{\{i_2=i_4\ne 0\}}
\hbox{\vtop{\offinterlineskip\halign{
\hfil#\hfil\cr
{\rm l.i.m.}\cr
$\stackrel{}{{}_{N\to \infty}}$\cr
}} }
\sum\limits_{l_4=0}^{N-1}
\sum\limits_{l_2=0}^{N-1} 
G_{pppp}(\tau_{l_2},\tau_{l_4},\tau_{l_2},\tau_{l_4})
\Delta\tau_{l_2}
\Delta\tau_{l_4}+
$$
$$
+{\bf 1}_{\{i_1=i_4\ne 0\}}{\bf 1}_{\{i_2=i_3\ne 0\}}
\hbox{\vtop{\offinterlineskip\halign{
\hfil#\hfil\cr
{\rm l.i.m.}\cr
$\stackrel{}{{}_{N\to \infty}}$\cr
}} }
\sum\limits_{l_4=0}^{N-1}
\sum\limits_{l_2=0}^{N-1} 
G_{pppp}(\tau_{l_2},\tau_{l_4},\tau_{l_4},\tau_{l_2})
\Delta\tau_{l_2}
\Delta\tau_{l_4}=
$$
$$
={\bf 1}_{\{i_1=i_2\ne 0\}}{\bf 1}_{\{i_3=i_4\ne 0\}}
\hbox{\vtop{\offinterlineskip\halign{
\hfil#\hfil\cr
{\rm l.i.m.}\cr
$\stackrel{}{{}_{N\to \infty}}$\cr
}} }
\sum\limits_{l_4=0}^{N-1}
\sum\limits_{l_2=0}^{N-1} 
\Biggl(\frac{1}{4}{\bf 1}_{\{\tau_{l_2}<\tau_{l_4}\}}
+\frac{1}{8}{\bf 1}_{\{\tau_{l_2}=\tau_{l_4}\}}-\Biggr.
$$
$$
\Biggl.
-\sum\limits_{j_4,j_3,j_2,j_1=0}^p C_{j_4j_3j_2j_1}
\phi_{j_1}(\tau_{l_2})
\phi_{j_2}(\tau_{l_2})
\phi_{j_3}(\tau_{l_4})
\phi_{j_4}(\tau_{l_4})\Biggr)
\Delta\tau_{l_2}
\Delta\tau_{l_4}+
$$
$$
+{\bf 1}_{\{i_1=i_3\ne 0\}}{\bf 1}_{\{i_2=i_4\ne 0\}}
\hbox{\vtop{\offinterlineskip\halign{
\hfil#\hfil\cr
{\rm l.i.m.}\cr
$\stackrel{}{{}_{N\to \infty}}$\cr
}} }
\sum\limits_{l_4=0}^{N-1}
\sum\limits_{l_2=0}^{N-1}\Biggl(
\frac{1}{8}{\bf 1}_{\{\tau_{l_2}=\tau_{l_4}\}}-
\sum\limits_{j_4,j_3,j_2,j_1=0}^p C_{j_4j_3j_2j_1}\times\Biggr.
$$
$$
\Biggl.\times
\phi_{j_1}(\tau_{l_2})
\phi_{j_2}(\tau_{l_4})
\phi_{j_3}(\tau_{l_2})
\phi_{j_4}(\tau_{l_4})\Biggr)
\Delta\tau_{l_2}
\Delta\tau_{l_4}+
$$
$$
+{\bf 1}_{\{i_1=i_4\ne 0\}}{\bf 1}_{\{i_2=i_3\ne 0\}}
\hbox{\vtop{\offinterlineskip\halign{
\hfil#\hfil\cr
{\rm l.i.m.}\cr
$\stackrel{}{{}_{N\to \infty}}$\cr
}} }
\sum\limits_{l_4=0}^{N-1}
\sum\limits_{l_2=0}^{N-1} 
\Biggl(\frac{1}{8}{\bf 1}_{\{\tau_{l_2}=\tau_{l_4}\}}-
\sum\limits_{j_4,j_3,j_2,j_1=0}^p C_{j_4j_3j_2j_1}\times\Biggr.
$$
$$
\Biggl.\times
\phi_{j_1}(\tau_{l_2})
\phi_{j_2}(\tau_{l_4})
\phi_{j_3}(\tau_{l_4})
\phi_{j_4}(\tau_{l_2})\Biggr)
\Delta\tau_{l_2}
\Delta\tau_{l_4}=
$$
$$
=
{\bf 1}_{\{i_1=i_2\ne 0\}}{\bf 1}_{\{i_3=i_4\ne 0\}}
\left(
\frac{1}{4}
\int\limits_t^T\int\limits_t^{t_4}dt_2dt_4
-\sum\limits_{j_4,j_1=0}^p C_{j_4j_4j_1j_1}\right)-
$$
$$
-{\bf 1}_{\{i_1=i_3\ne 0\}}{\bf 1}_{\{i_2=i_4\ne 0\}}
\sum\limits_{j_4,j_1=0}^p C_{j_4j_1j_4j_1}-
$$
$$
-{\bf 1}_{\{i_1=i_4\ne 0\}}{\bf 1}_{\{i_2=i_3\ne 0\}}
\sum\limits_{j_4,j_2=0}^p C_{j_4j_2j_2j_4}.
$$

\vspace{2mm}

When proving Theorem 2.9
we have proved that
\begin{equation}
\label{july7010}
\hbox{\vtop{\offinterlineskip\halign{
\hfil#\hfil\cr
{\rm lim}\cr
$\stackrel{}{{}_{p\to \infty}}$\cr
}} }
\sum\limits_{j_4,j_1=0}^{p} C_{j_4j_4j_1j_1}=
\frac{1}{4}
\int\limits_t^T\int\limits_t^{t_4}dt_2dt_4,
\end{equation}
\begin{equation}
\label{july7011}
\hbox{\vtop{\offinterlineskip\halign{
\hfil#\hfil\cr
{\rm lim}\cr
$\stackrel{}{{}_{p\to \infty}}$\cr
}} }
\sum\limits_{j_4,j_1=0}^{p} C_{j_4j_1j_4j_1}=0,\ \ \ 
\hbox{\vtop{\offinterlineskip\halign{
\hfil#\hfil\cr
{\rm lim}\cr
$\stackrel{}{{}_{p\to \infty}}$\cr
}} }
\sum\limits_{j_4,j_2=0}^{p} C_{j_4j_2j_2j_4}=0.
\end{equation}

\vspace{1mm}

Then
$$
\hbox{\vtop{\offinterlineskip\halign{
\hfil#\hfil\cr
{\rm lim}\cr
$\stackrel{}{{}_{p\to \infty}}$\cr
}} }
R_{T,t}^{(7)pppp}=0.
$$

\vspace{2mm}

Theorem 2.9 is proved.

\section{Modification of Theorems 2.1, 2.8, and 2.9 for the Case
of In\-teg\-ra\-tion Interval $[t, s]$ $(s\in (t, T])$ 
of Iterated Stra\-to\-no\-vich Sto\-chas\-tic Integrals of Multiplicities
2 to 4 and Wong--Zakai Type Theorems}

\subsection{Modification of Theorem 2.1 for the Case
of In\-teg\-ra\-tion Interval $[t, s]$ $(s\in (t, T])$ 
of Iterated Stra\-to\-no\-vich Sto\-chas\-tic Integrals of Multiplicity
2}

Let us prove the following theorem.

\vspace{2mm}       

{\bf Theorem 2.18.}\ {\it Suppose that 
$\{\phi_j(x)\}_{j=0}^{\infty}$ is an arbitrary complete orthonormal system of 
functions in the space $L_2([t, T]).$
Moreover$,$ $\psi_1(\tau), \psi_2(\tau)$ are continuous functions on $[t, T]$.
Then$,$ 
for the iterated Stratonovich stochastic integral
$$
J^{*}[\psi^{(2)}]_{s,t}={\int\limits_t^{*}}^s\psi_2(t_2)
{\int\limits_t^{*}}^{t_2}\psi_1(t_1)d{\bf w}_{t_1}^{(i_1)}
d{\bf w}_{t_2}^{(i_2)}\ \ \ (i_1, i_2=1,\ldots,m)
$$

\noindent
the following expansion 
\begin{equation}
\label{jesxxx}
J^{*}[\psi^{(2)}]_{s,t}=\hbox{\vtop{\offinterlineskip\halign{
\hfil#\hfil\cr
{\rm l.i.m.}\cr
$\stackrel{}{{}_{p_1,p_2\to \infty}}$\cr
}} }\sum_{j_1=0}^{p_1}\sum_{j_2=0}^{p_2}
C_{j_2j_1}(s)\zeta_{j_1}^{(i_1)}\zeta_{j_2}^{(i_2)}
\end{equation}

\noindent
that converges in the mean-square
sense 
is valid$,$ where $s\in (t, T]$ $(s$ is fixed{\rm )},
\begin{equation}
\label{tupo11xxx}
C_{j_2 j_1}(s)=\int\limits_t^s\psi_2(t_2)\phi_{j_2}(t_2)
\int\limits_t^{t_2}\psi_1(t_1)\phi_{j_1}(t_1)dt_1dt_2,
\end{equation}
and
$$
\zeta_{j}^{(i)}=
\int\limits_t^T \phi_{j}(\tau) d{\bf w}_{\tau}^{(i)}
$$ 
are independent
standard Gaussian random variables for various 
$i$ or $j$.}

The condition of continuity of the functions
$\psi_1(\tau), \psi_2(\tau)$ 
is related to the definition (\ref{123321.2}) 
of the Stratonovich stochastic integral that we use.

{\bf Proof.}\ The case $s=T$ is considered in Theorems 2.1--2.3.
Below we consider the case $s\in (t, T).$ In accordance to the standard relations between
Stra\-to\-no\-vich and It\^{o} stochastic integrals (see (\ref{d11})
and (\ref{d11a})) 
we have w.~p.~1 
\begin{equation}
\label{oop51zzz}
J^{*}[\psi^{(2)}]_{s,t}=
J[\psi^{(2)}]_{s,t}+
\frac{1}{2}{\bf 1}_{\{i_1=i_2\}}
\int\limits_t^s\psi_1(t_1)\psi_2(t_1)dt_1,
\end{equation}
where $\psi_1(\tau), \psi_2(\tau)$ are continuous functions on $[t, T]$,
$s\in (t, T)$ ($s$ is fixed),
${\bf 1}_A$ is the indicator of the set $A.$

From the other side according to (\ref{a2xxx}), we obtain
$$
J[\psi^{(2)}]_{s,t}=
\hbox{\vtop{\offinterlineskip\halign{
\hfil#\hfil\cr
{\rm l.i.m.}\cr
$\stackrel{}{{}_{p_1,p_2\to \infty}}$\cr
}} }\sum_{j_1=0}^{p_1}\sum_{j_2=0}^{p_2}
C_{j_2j_1}(s)\Biggl(\zeta_{j_1}^{(i_1)}\zeta_{j_2}^{(i_2)}
-{\bf 1}_{\{i_1=i_2\}}
{\bf 1}_{\{j_1=j_2\}}\Biggr)=
$$
$$
=\hbox{\vtop{\offinterlineskip\halign{
\hfil#\hfil\cr
{\rm l.i.m.}\cr
$\stackrel{}{{}_{p_1,p_2\to \infty}}$\cr
}} }\sum_{j_1=0}^{p_1}\sum_{j_2=0}^{p_2}
C_{j_2j_1}(s)\zeta_{j_1}^{(i_1)}\zeta_{j_2}^{(i_2)}
-
$$
\begin{equation}
\label{yes2001zzz}
-
{\bf 1}_{\{i_1=i_2\}}\lim\limits_{p_1,p_2\to\infty}\sum_{j_1=0}^{\min\{p_1,p_2\}}
C_{j_1j_1}(s).
\end{equation}

From (\ref{oop51zzz}) and (\ref{yes2001zzz}) it follows that
Theorem 2.18 will be proved if 
\begin{equation}
\label{5tzzz}
\frac{1}{2}
\int\limits_t^s\psi_1(t_1)\psi_2(t_1)dt_1
=\sum_{j_1=0}^{\infty}
C_{j_1j_1}(s),
\end{equation}
\noindent
where $\psi_1(\tau), \psi_2(\tau)\in L_2([t, T]).$

We have (see Sect.~2.1.4)
\begin{equation}
\label{strange201}
~~~~~~~~~\frac{1}{2}\int\limits_t^T 
\bar \psi_1(\tau) \bar\psi_2(\tau) d\tau
=
\sum_{j=0}^{\infty}\int\limits_t^T
\bar \psi_2(t_2)\phi_j(t_2)
\int\limits_t^{t_2}
\bar \psi_1(t_1)\phi_j(t_1)dt_1 dt_2,
\end{equation}

\noindent
where $\bar \psi_1(\tau), \bar \psi_2(\tau)\in L_2([t, T]).$

Suppose that
\begin{equation}
\label{strange202}
\bar \psi_1(\tau)=\psi_1(\tau){\bf 1}_{\{\tau<s\}},\ \ \ 
\bar \psi_2(\tau)=\psi_2(\tau){\bf 1}_{\{\tau<s\}},
\end{equation}

\noindent
where $\psi_1(\tau), \psi_2(\tau)\in L_2([t,T]),$ $s\in (t, T)$ ($s$ is fixed).

Combining (\ref{strange201}) and (\ref{strange202}), we get
$$
\frac{1}{2}\int\limits_t^T
\hspace{-0.4mm}\psi_1(\tau)\psi_2(\tau){\bf 1}_{\{\tau<s\}} d\tau
=
\sum_{j=0}^{\infty}\int\limits_t^T\hspace{-0.4mm}
\psi_2(t_2){\bf 1}_{\{t_2<s\}}\phi_j(t_2)
\hspace{-0.4mm}\int\limits_t^{t_2}
\psi_1(t_1){\bf 1}_{\{t_1<s\}}\phi_j(t_1)dt_1 dt_2,
$$
i.e.
$$
\frac{1}{2}\int\limits_t^s
\psi_1(\tau)\psi_2(\tau)d\tau
=
\sum_{j=0}^{\infty}\int\limits_t^s
\psi_2(t_2)\phi_j(t_2)
\int\limits_t^{t_2}
\psi_1(t_1)\phi_j(t_1)dt_1 dt_2.
$$

The equality (\ref{5tzzz}) is proved.
Theorem 2.18 is proved.

Let us reformulate Theorem 2.18 in terms on the convergence 
of the solution of system of ordinary differential 
equations (ODEs) to the solution of system of Stratonovich 
SDEs (the so-called Wong--Zakai type theorem).

By analogy with (\ref{um1xxxx1}) for $k=2,$\
$i_1, i_2=1,\ldots,m$,  and $s\in (t, T]$ ($s$ is fixed) we
obtain
\begin{equation}
\label{rs222}
~~~~~~~\int\limits_t^s
\psi_2(t_2)\int\limits_t^{t_2}\psi_1(t_1)
d{\bf w}_{t_1}^{(i_1)p_1}d{\bf w}_{t_2}^{(i_2)p_2}=
\sum\limits_{j_1=0}^{p_1}\sum\limits_{j_2=0}^{p_2}
C_{j_2j_1}(s)\zeta_{j_1}^{(i_1)}\zeta_{j_2}^{(i_2)},
\end{equation}

\noindent
where $p_1, p_2\in{\bf N}$ and $d{\bf w}_{\tau}^{(i)p}$ is defined by
(\ref{um1xxx}); another notations are the same as in Theorem 2.18.

The iterated Riemann--Stiltjes integrals 
$$
Y_{s,t}^{(i_1i_2)p_1p_2}=\int\limits_t^s
\psi_2(t_2)\int\limits_t^{t_2}\psi_1(t_1)
d{\bf w}_{t_1}^{(i_1)p_1}d{\bf w}_{t_2}^{(i_2)p_2},
$$
$$
X_{s,t}^{(i_1)p_1}=\int\limits_t^{s}\psi_1(t_1)
d{\bf w}_{t_1}^{(i_1)p_1}
$$

\noindent
are the solution of the following system of ODEs

\vspace{-3mm}
$$
\left\{\begin{matrix}
dY_{s,t}^{(i_1i_2)p_1p_2}=\psi_2(s)X_{s,t}^{(i_1)p_1}
d{\bf w}_{s}^{(i_2)p_2},\ &Y_{t,t}^{(i_1i_2)p_1p_2}=0\cr\cr
dX_{s,t}^{(i_1)p_1}=\psi_1(s)d{\bf w}_{s}^{(i_1)p_1},\ 
&X_{t,t}^{(i_1)p_1}=0
\end{matrix}\right..
$$

\vspace{2mm}

From the other hand, the iterated Stratonovich
stochastic integrals 
$$
Y_{s,t}^{(i_1i_2)}={\int\limits_t^{*}}^s\psi_2(t_2)
{\int\limits_t^{*}}^{t_2}\psi_1(t_1)d{\bf w}_{t_1}^{(i_1)}
d{\bf w}_{t_2}^{(i_2)},
$$
$$
X_{s,t}^{(i_1)}
={\int\limits_t^{*}}^{s}\psi_1(t_1)d{\bf w}_{t_1}^{(i_1)}
$$

\vspace{1mm}
\noindent
are the solution of the following system of Stratonovich SDEs

\vspace{-3mm}
$$
\left\{\begin{matrix}
dY_{s,t}^{(i_1i_2)}=\psi_2(s)X_{s,t}^{(i_1)}
* d{\bf w}_{s}^{(i_2)},\ &Y_{t,t}^{(i_1i_2)}=0\cr\cr
dX_{s,t}^{(i_1)}=\psi_1(s) * d{\bf w}_{s}^{(i_1)},\ 
&X_{t,t}^{(i_1)}=0
\end{matrix}\right.,
$$

\vspace{2mm}
\noindent
where  $*~\hspace{-0.3mm}d{\bf w}_{s}^{(i)}$, $i=1,\ldots,m$ is the
Stratonovich differential.

Then from Theorem 2.18 and (\ref{a1uuu}) we obtain the following theorem.

{\bf Theorem 2.19}\ \cite{arxiv-5}. {\it Suppose that 
$\{\phi_j(x)\}_{j=0}^{\infty}$ is an arbitrary complete orthonormal system of 
functions in the space $L_2([t, T]).$
Moreover$,$ $\psi_1(\tau), \psi_2(\tau)$ are continuous functions on $[t, T]$.
Then
for any fixed $s$ $(s\in (t, T])$

\vspace{-3mm}
$$
\hbox{\vtop{\offinterlineskip\halign{
\hfil#\hfil\cr
{\rm l.i.m.}\cr
$\stackrel{}{{}_{p_1,p_2\to \infty}}$\cr
}} }Y_{s,t}^{(i_1i_2)p_1p_2}=Y_{s,t}^{(i_1i_2)},\ \ \
\hbox{\vtop{\offinterlineskip\halign{
\hfil#\hfil\cr
{\rm l.i.m.}\cr
$\stackrel{}{{}_{p_1\to \infty}}$\cr
}} }X_{s,t}^{(i_1)p_1}=X_{s,t}^{(i_1)}.
$$
}

\subsection{Modification of Theorem 2.8 for the Case
of In\-teg\-ra\-tion Interval $[t, s]$ $(s\in (t, T])$ 
of Iterated Stra\-to\-no\-vich Sto\-chas\-tic Integrals of Multiplicity 3}

Let us prove the following theorem.

{\bf Theorem 2.20}\ \cite{arxiv-5}.
{\it Suppose that 
$\{\phi_j(x)\}_{j=0}^{\infty}$ is a complete orthonormal system of 
Legendre polynomials or trigonometric functions in the space $L_2([t, T]).$
At the same time $\psi_2(\tau)$ is a continuously dif\-ferentiable 
nonrandom function on $[t, T]$ and $\psi_1(\tau),$ $\psi_3(\tau)$ are twice
continuously differentiable nonrandom functions on $[t, T]$. 
Then$,$ for the 
iterated Stratonovich stochastic integral of third multiplicity
$$
J^{*}[\psi^{(3)}]_{s,t}={\int\limits_t^{*}}^s\psi_3(t_3)
{\int\limits_t^{*}}^{t_3}\psi_2(t_2)
{\int\limits_t^{*}}^{t_2}\psi_1(t_1)
d{\bf w}_{t_1}^{(i_1)}
d{\bf w}_{t_2}^{(i_2)}d{\bf w}_{t_3}^{(i_3)},
$$
where $i_1, i_2, i_3=1,\ldots,m,$ the following 
expansion 
$$
J^{*}[\psi^{(3)}]_{s,t}
=\hbox{\vtop{\offinterlineskip\halign{
\hfil#\hfil\cr
{\rm l.i.m.}\cr
$\stackrel{}{{}_{p\to \infty}}$\cr
}} }
\sum\limits_{j_1, j_2, j_3=0}^{p}
C_{j_3 j_2 j_1}(s)\zeta_{j_1}^{(i_1)}\zeta_{j_2}^{(i_2)}\zeta_{j_3}^{(i_3)}
$$
that converges in the mean-square sense
is valid, where $s\in (t, T]$ $(s$ is fixed{\rm )},
$$
C_{j_3 j_2 j_1}(s)=\int\limits_t^s\psi_3(t_3)\phi_{j_3}(t_3)
\int\limits_t^{t_3}\psi_2(t_2)\phi_{j_2}(t_2)
\int\limits_t^{t_2}\psi_1(t_1)\phi_{j_1}(t_1)dt_1dt_2dt_3
$$
and
$$
\zeta_{j}^{(i)}=
\int\limits_t^T \phi_{j}(\tau) d{\bf w}_{\tau}^{(i)}
$$ 
are independent standard Gaussian random variables for various 
$i$ or $j$.}

{\bf Proof.} The case $s=T$ is considered in Theorem 2.8. 
Below we consider the case $s\in (t, T).$ First, let us consider the case of
Legendre polynomials. 
From (\ref{result3}) for the case $p_1=p_2=p_3=p$ and standard
relations between It\^{o} and Stratonovich stochastic integrals 
we conclude that Theorem 2.20 will be proved if w.~p.~1
\begin{equation}
\label{resul4}
~~~~\hbox{\vtop{\offinterlineskip\halign{
\hfil#\hfil\cr
{\rm l.i.m.}\cr
$\stackrel{}{{}_{p\to \infty}}$\cr
}} }
\sum\limits_{j_1=0}^{p}\sum\limits_{j_3=0}^{p}
C_{j_3 j_1 j_1}(s)\zeta_{j_3}^{(i_3)}=
\frac{1}{2}\int\limits_t^s\psi_3(\tau)
\int\limits_t^{\tau}\psi_2(s_1)\psi_1(s_1)ds_1d{\bf w}_{\tau}^{(i_3)},
\end{equation}
\begin{equation}
\label{resul5}
~~~~\hbox{\vtop{\offinterlineskip\halign{
\hfil#\hfil\cr
{\rm l.i.m.}\cr
$\stackrel{}{{}_{p\to \infty}}$\cr
}} }
\sum\limits_{j_1=0}^{p}\sum\limits_{j_3=0}^{p}
C_{j_3 j_3 j_1}(s)\zeta_{j_1}^{(i_1)}=
\frac{1}{2}\int\limits_t^s\psi_3(\tau)\psi_2(\tau)
\int\limits_t^{\tau}\psi_1(s_1)d{\bf w}_{s_1}^{(i_1)}d\tau,
\end{equation}
\begin{equation}
\label{resul6}
\hbox{\vtop{\offinterlineskip\halign{
\hfil#\hfil\cr
{\rm l.i.m.}\cr
$\stackrel{}{{}_{p\to \infty}}$\cr
}} }
\sum\limits_{j_1=0}^{p}\sum\limits_{j_3=0}^{p}
C_{j_1 j_3 j_1}(s)\zeta_{j_3}^{(i_2)}=0.
\end{equation}

\vspace{2mm}

The proof of the formulas (\ref{resul4}), (\ref{resul6}) 
is absolutely similar to the proof of 
the formulas (\ref{1xx}), (\ref{3xx}). It is only necessary 
to replace the 
interval of integration $[t,T]$ by $[t,s]$ 
in the proof of the formulas (\ref{1xx}), (\ref{3xx}) 
and use Theorem 1.11 instead of Theorem 1.1.
Also in the case (\ref{resul6}) it is 
necessary to use the estimate (\ref{101oh}).

Let us prove (\ref{resul5}).
Using Theorem 1.11 for $k=2$ (see (\ref{a2xxx}) for $i_1=1,\ldots,m,\ i_2=0$), 
we obtain w.~p.~1 (also see (\ref{dwdw21}), (\ref{dwdw22}))
$$
\frac{1}{2}\int\limits_t^s\psi_3(\tau)\psi_2(\tau)
\int\limits_t^{\tau}\psi_1(s_1)d{\bf w}_{s_1}^{(i_1)}d\tau=
\frac{1}{2}\
\hbox{\vtop{\offinterlineskip\halign{
\hfil#\hfil\cr
{\rm l.i.m.}\cr
$\stackrel{}{{}_{p\to \infty}}$\cr
}} }
\sum\limits_{j_1=0}^{p}
C_{j_1}^{*}(s)\zeta_{j_1}^{(i_1)},
$$
where 
$$
C_{j_1}^{*}(s)=
\int\limits_t^s \psi_3(\tau)\psi_2(\tau)
\int\limits_{t}^{\tau}\psi_1(s_1)\phi_{j_1}(s_1)
ds_1 d\tau=
$$
\begin{equation}
\label{resul11}
=
\int\limits_t^s
\psi_1(s_1)\phi_{j_1}(s_1)\int\limits_{s_1}^{s}
\psi_3(\tau)\psi_2(\tau)d\tau ds_1.
\end{equation}

We have
$$
E_p'(s)\stackrel{\sf def}{=}{\sf M}\left\{\left(
\sum\limits_{j_1=0}^{p}\sum\limits_{j_3=0}^{p}
C_{j_3 j_3 j_1}(s)\zeta_{j_1}^{(i_1)} - 
\frac{1}{2}\sum\limits_{j_1=0}^{p}
C_{j_1}^{*}(s)\zeta_{j_1}^{(i_1)}\right)^2\right\}=
$$
$$
={\sf M}\left\{\left(\sum_{j_1=0}^p\left(\sum_{j_3=0}^p
C_{j_3j_3j_1}(s)-\frac{1}{2}C_{j_1}^{*}(s)\right)
\zeta_{j_1}^{(i_1)}\right)^2\right\}
=
$$
\begin{equation}
\label{resul12}
=\sum_{j_1=0}^p\left(\sum\limits_{j_3=0}^{p}C_{j_3j_3 j_1}(s)-
\frac{1}{2}C_{j_1}^{*}(s)\right)^2,
\end{equation}
$$
C_{j_3 j_3 j_1}(s)=\int\limits_t^s\psi_3(\theta)\phi_{j_3}(\theta)
\int\limits_t^{\theta}\psi_2(\tau)\phi_{j_3}(\tau)
\int\limits_t^{\tau}\psi_1(s_1)\phi_{j_1}(s_1)ds_1d\tau d\theta=
$$
\begin{equation}
\label{resul15}
~~~~~~~~~=\int\limits_t^s\psi_1(s_1)\phi_{j_1}(s_1)
\int\limits_{s_1}^s\psi_2(\tau)\phi_{j_3}(\tau)
\int\limits_{\tau}^s\psi_3(\theta)\phi_{j_3}(\theta)d\theta
d\tau ds_1.
\end{equation}

From (\ref{resul11})--(\ref{resul15}) we obtain 
$$
E_p'(s)
=\sum_{j_1=0}^p\left(
\int\limits_t^s\psi_1(s_1)\phi_{j_1}(s_1)
\left(\sum\limits_{j_3=0}^{p}\int\limits_{s_1}^s
\psi_2(\tau)\phi_{j_3}(\tau)
\int\limits_{\tau}^s\psi_3(\theta)\phi_{j_3}
(\theta)d\theta d\tau- 
\right.\right.
$$
\begin{equation}
\label{resul20}
\left.\left.
-\frac{1}{2}
\int\limits_{s_1}^s \psi_3(\tau)\psi_2(\tau)d\tau\right) ds_1\right)^2.
\end{equation}

\vspace{2mm}

Let us show that
\begin{equation}
\label{dwdw6}
~~~~~~~~\sum\limits_{j_3=0}^{\infty}\int\limits_{s_1}^{s}\psi_2(\tau)\phi_{j_3}(\tau)
\int\limits_{\tau}^{s}\psi_3(\theta)\phi_{j_3}(\theta)d\theta d\tau
=\frac{1}{2}\int\limits_{s_1}^s \psi_3(\tau)\psi_2(\tau)d\tau.
\end{equation}

\vspace{1mm}

Using (\ref{strange201}) and Fubini's Theorem, we have
\begin{equation}
\label{strange300}
~~~~~~~~~\frac{1}{2}\int\limits_t^T 
\bar \psi_1(\tau) \bar\psi_2(\tau) d\tau
=
\sum_{j=0}^{\infty}\int\limits_t^T
\bar \psi_1(t_1)\phi_j(t_1)
\int\limits_{t_1}^{T}
\bar \psi_2(t_2)\phi_j(t_2)dt_2 dt_1,
\end{equation}

\vspace{2mm}
\noindent
where $\bar \psi_1(\tau), \bar \psi_2(\tau)\in L_2([t, T]).$

Suppose that
\begin{equation}
\label{strange402}
\bar \psi_1(\tau)=\psi_2(\tau){\bf 1}_{\{s_1<\tau<s\}},\ \ \ 
\bar \psi_2(\tau)=\psi_3(\tau){\bf 1}_{\{\tau<s\}}.
\end{equation}

\vspace{1mm}

Using (\ref{strange300}) and (\ref{strange402}), we get (\ref{dwdw6}).
Combining (\ref{resul20}) and (\ref{dwdw6}), we obtain
$$
E_p'(s)
=\sum_{j_1=0}^p\left(
\int\limits_t^s\psi_1(s_1)\phi_{j_1}(s_1)
\sum\limits_{j_3=p+1}^{\infty}\int\limits_{s_1}^s
\psi_2(\tau)\phi_{j_3}(\tau)
\int\limits_{\tau}^s\psi_3(\theta)\phi_{j_3}
(\theta)d\theta d\tau ds_1\right)^2\le
$$
\begin{equation}
\label{dwdw7}
\le K \sum_{j_1=0}^p\left(
\int\limits_t^s |\phi_{j_1}(s_1)|
\left|\sum\limits_{j_3=p+1}^{\infty}\int\limits_{s_1}^s
\psi_2(\tau)\phi_{j_3}(\tau)
\int\limits_{\tau}^s\psi_3(\theta)\phi_{j_3}
(\theta)d\theta d\tau \right| ds_1\right)^2,
\end{equation}

\noindent
where constant $K$ does not depend on $p$.

Let us estimate the value
$$
\left|\sum\limits_{j_3=p+1}^{\infty}\int\limits_{s_1}^s
\psi_2(\tau)\phi_{j_3}(\tau)
\int\limits_{\tau}^s\psi_3(\theta)\phi_{j_3}
(\theta)d\theta d\tau\right|.
$$

\vspace{1mm}

Note that, by virtue of the additivity property of the integral, we obtain
$$
\int\limits_{s_1}^{s}\psi_2(\tau)\phi_{j_3}(\tau)
\int\limits_{\tau}^{s}\psi_3(\theta)\phi_{j_3}(\theta)d\theta d\tau=
$$
$$
=\int\limits_{t}^{s}\psi_3(\theta)\phi_{j_3}(\theta)
\int\limits_{t}^{\theta}\psi_2(\tau)\phi_{j_3}(\tau)d\tau d\theta-
$$
$$
-\int\limits_{t}^{s_1}\psi_3(\theta)\phi_{j_3}(\theta)
\int\limits_{t}^{\theta}\psi_2(\tau)\phi_{j_3}(\tau)d\tau d\theta-
$$
$$
-\int\limits_{s_1}^{s}\psi_3(\theta)\phi_{j_3}(\theta)d\theta
\int\limits_{t}^{s_1}\psi_2(\tau)\phi_{j_3}(\tau)d\tau.
$$

\vspace{1mm}

Further, we have 
$$
\left|\sum\limits_{j_3=p+1}^{\infty}\int\limits_{s_1}^s
\psi_2(\tau)\phi_{j_3}(\tau)
\int\limits_{\tau}^s\psi_3(\theta)\phi_{j_3}
(\theta)d\theta d\tau\right|\le
$$
$$
\le\left|\sum\limits_{j_3=p+1}^{\infty}\int\limits_{t}^{s}\psi_3(\theta)\phi_{j_3}(\theta)
\int\limits_{t}^{\theta}\psi_2(\tau)\phi_{j_3}(\tau)d\tau d\theta\right|+
$$
$$
+\left|\sum\limits_{j_3=p+1}^{\infty}\int\limits_{t}^{s_1}\psi_3(\theta)\phi_{j_3}(\theta)
\int\limits_{t}^{\theta}\psi_2(\tau)\phi_{j_3}(\tau)d\tau d\theta\right|+
$$
\begin{equation}
\label{dwdw7a}
+\sum\limits_{j_3=p+1}^{\infty}\left|\int\limits_{s_1}^{s}\psi_3(\theta)\phi_{j_3}(\theta)d\theta
\int\limits_{t}^{s_1}\psi_2(\tau)\phi_{j_3}(\tau)d\tau\right|.
\end{equation}

\vspace{3mm}

Applying the estimate (\ref{fin1000}) (see Sect.~2.9), we can write
\begin{equation}
\label{dwdw9}
\left|\sum_{j_1=p+1}^{\infty}
C_{j_1j_1}(s)\right|\le \frac{C}{p}\left(
1+\frac{1}{\left(1-(z(s))^2\right)^{1/4}}\right),
\end{equation}

\noindent
where $s\in (t, T),$ constant $C$ does not depend on $p,$ $z(s)$ has the form (\ref{zz1}),
and $C_{j_1j_1}(s)$ is defined by (\ref{tupo11xxx}) for the case $j_1=j_2$.

Applying the estimates (\ref{101oh}), (\ref{101xx}), (\ref{dwdw9}) to the right-hand side 
of (\ref{dwdw7a}) gives
$$
\left|\sum\limits_{j_3=p+1}^{\infty}\int\limits_{s_1}^s
\psi_2(\tau)\phi_{j_3}(\tau)
\int\limits_{\tau}^s\psi_3(\theta)\phi_{j_3}
(\theta)d\theta d\tau\right|\le \frac{L}{p}
\left(1+\frac{1}{\left(1-(z(s_1))^2\right)^{1/4}}\right)\times
$$
\begin{equation}
\label{dwdw10}
\times
\left(1+\frac{1}{\left(1-(z(s))^2\right)^{1/4}}+
\frac{1}{\left(1-(z(s_1))^2\right)^{1/4}}
\right),
\end{equation}

\vspace{3mm}
\noindent
where $s, s_1\in (t, T)$ and constant $L$ is independent of $p$.

Combining the estimates (\ref{ogo24}), (\ref{dwdw7}), and (\ref{dwdw10}),
we finally obtain
$$
E_p'(s)\le
\frac{L(s)p}{p^2}=\frac{L(s)}{p}
$$
if $p\to\infty$,
where constant $L(s)$ ($s$ is fixed, $s\in (t, T)$) does not depend on $p$.
The relation (\ref{resul5}) is proved for the polynomial case. Theorem 2.20 is proved
for the case of Legendre polynomials.

For the trigonometric case, by analogy with the proof of Lemma~2.2 (Sect.~2.1.2), 
we obtain the following analog of 
(\ref{dwdw9}) 
\begin{equation}
\label{dwdw12}
\left|\sum_{j_1=p+1}^{\infty}
C_{j_1j_1}(s)\right|\le \frac{C}{p},
\end{equation}
where $s\in [t, T],$ constant $C$ does not depend on $p,$ 
and $C_{j_1j_1}(s)$ is defined by (\ref{tupo11xxx}) for the case $j_1=j_2$.

Note the following obvious estimates for the trigonometric case
\begin{equation}
\label{dwdw14}
~~~~~~~~\left|\int\limits_{s_1}^{s}\psi_3(\theta)\phi_{j}(\theta)d\theta\right|\le \frac{C}{j},\ \ \
\left|\int\limits_{t}^{s_1}\psi_2(\tau)\phi_{j}(\tau)d\tau\right|\le \frac{C}{j}\ \ \ (j\ne 0),
\end{equation}
where $s, s_1\in [t, T],$ constant $C$ does not depend on $p.$ 

Applying (\ref{dwdw7}), (\ref{dwdw7a}), (\ref{dwdw12}), and (\ref{dwdw14}),
we obtain the assertion of Theorem~2.20 for the trigonometric case.
Theorem 2.20 is proved.

Let us reformulate Theorem 2.20 in terms on the convergence 
of the solution of system of ODEs to the solution of  
system of Stratonovich 
SDEs (the so-called Wong--Zakai type theorem).

By analogy with (\ref{um1xxxx1}) for the case $k=3$,\ $p_1=p_2=p_3=p,$\
$i_1, i_2, i_3=1,\ldots,m$,  and $s\in (t, T]$ ($s$ is fixed) we
obtain
\begin{equation}
\label{resul80}
\int\limits_t^s
\psi_3(t_3)\int\limits_t^{t_3}\psi_2(t_2)
\int\limits_t^{t_2}\psi_1(t_1)
d{\bf w}_{t_1}^{(i_1)p}d{\bf w}_{t_2}^{(i_2)p}
d{\bf w}_{t_3}^{(i_3)p}=
\sum\limits_{j_1,j_2,j_3=0}^{p}
C_{j_3j_2j_1}(s)\zeta_{j_1}^{(i_1)}\zeta_{j_2}^{(i_2)}\zeta_{j_3}^{(i_3)},
\end{equation}

\vspace{1mm}
\noindent
where $p\in{\bf N}$ and $d{\bf w}_{\tau}^{(i)p}$ is defined by
(\ref{um1xxx}); another notations are the same as in Theorem 2.20.

The iterated Riemann--Stiltjes integrals 
$$
Z_{s,t}^{(i_1i_2i_3)p}=
\int\limits_t^s
\psi_3(t_3)\int\limits_t^{t_3}\psi_2(t_2)
\int\limits_t^{t_2}\psi_1(t_1)
d{\bf w}_{t_1}^{(i_1)p}d{\bf w}_{t_2}^{(i_2)p}
d{\bf w}_{t_3}^{(i_3)p},
$$
$$
Y_{s,t}^{(i_1i_2)p}=\int\limits_t^s
\psi_2(t_2)\int\limits_t^{t_2}\psi_1(t_1)
d{\bf w}_{t_1}^{(i_1)p}d{\bf w}_{t_2}^{(i_2)p},
$$
$$
X_{s,t}^{(i_1)p}=\int\limits_t^{s}\psi_1(t_1)
d{\bf w}_{t_1}^{(i_1)p}
$$

\noindent
are the solution of the following system of ODEs

$$
\left\{\begin{matrix}
dZ_{s,t}^{(i_1i_2i_3)p}=\psi_3(s)Y_{s,t}^{(i_1i_2)p}
d{\bf w}_{s}^{(i_3)p},\ &Z_{t,t}^{(i_1i_2i_3)p}=0\cr\cr
dY_{s,t}^{(i_1i_2)p}=\psi_2(s)X_{s,t}^{(i_1)p}
d{\bf w}_{s}^{(i_2)p},\ &Y_{t,t}^{(i_1i_2)p}=0\cr\cr
dX_{s,t}^{(i_1)p}=\psi_1(s)d{\bf w}_{s}^{(i_1)p},\ 
&X_{t,t}^{(i_1)p}=0
\end{matrix}\right..
$$

\vspace{4mm}

From the other hand, the iterated Stratonovich
stochastic integrals 

\vspace{-1mm}
$$
Z_{s,t}^{(i_1i_2i_3)}={\int\limits_t^{*}}^s\psi_3(t_3)
{\int\limits_t^{*}}^{t_3}\psi_2(t_2)
{\int\limits_t^{*}}^{t_2}\psi_1(t_1)d{\bf w}_{t_1}^{(i_1)}
d{\bf w}_{t_2}^{(i_2)}d{\bf w}_{t_3}^{(i_3)},
$$
$$
Y_{s,t}^{(i_1i_2)}={\int\limits_t^{*}}^s\psi_2(t_2)
{\int\limits_t^{*}}^{t_2}\psi_1(t_1)d{\bf w}_{t_1}^{(i_1)}
d{\bf w}_{t_2}^{(i_2)},
$$
$$
X_{s,t}^{(i_1)}
={\int\limits_t^{*}}^{s}\psi_1(t_1)d{\bf w}_{t_1}^{(i_1)}
$$

\vspace{3mm}
\noindent
are the solution of the following system of Stratonovich SDEs

$$
\left\{\begin{matrix}
dZ_{s,t}^{(i_1i_2i_3)}=\psi_3(s)Y_{s,t}^{(i_1i_2)}
* d{\bf w}_{s}^{(i_3)},\ &Z_{t,t}^{(i_1i_2i_3)}=0\cr\cr
dY_{s,t}^{(i_1i_2)}=\psi_2(s)X_{s,t}^{(i_1)}
* d{\bf w}_{s}^{(i_2)},\ &Y_{t,t}^{(i_1i_2)}=0\cr\cr
dX_{s,t}^{(i_1)}=\psi_1(s) * d{\bf w}_{s}^{(i_1)},\ 
&X_{t,t}^{(i_1)}=0
\end{matrix}\right.,
$$

\vspace{3mm}
\noindent
where  $*~\hspace{-0.3mm}d{\bf w}_{s}^{(i)}$, $i=1,\ldots,m$ is the
Stratonovich differential.

Then from Theorems 2.19 
and 2.20 we obtain the following theorem.

{\bf Theorem 2.21}\ \cite{arxiv-5}. {\it Suppose that 
$\{\phi_j(x)\}_{j=0}^{\infty}$ is a complete orthonormal system of 
Legendre polynomials or trigonometric 
functions in the space $L_2([t, T]).$
At the same time $\psi_2(\tau)$ is a continuously dif\-ferentiable 
nonrandom function on $[t, T]$ and $\psi_1(\tau),$ $\psi_3(\tau)$ are twice
continuously differentiable nonrandom functions on $[t, T]$. 
Then
for any fixed $s$ $(s\in (t, T])$

\vspace{-2mm}
$$
\hbox{\vtop{\offinterlineskip\halign{
\hfil#\hfil\cr
{\rm l.i.m.}\cr
$\stackrel{}{{}_{p\to \infty}}$\cr
}} }Z_{s,t}^{(i_1i_2i_3)p}=Z_{s,t}^{(i_1i_2i_3)},\ \ \ 
\hbox{\vtop{\offinterlineskip\halign{
\hfil#\hfil\cr
{\rm l.i.m.}\cr
$\stackrel{}{{}_{p\to \infty}}$\cr
}} }Y_{s,t}^{(i_1i_2)p}=Y_{s,t}^{(i_1i_2)},
$$
$$
\hbox{\vtop{\offinterlineskip\halign{
\hfil#\hfil\cr
{\rm l.i.m.}\cr
$\stackrel{}{{}_{p\to \infty}}$\cr
}} }X_{s,t}^{(i_1)p}=X_{s,t}^{(i_1)}.
$$
}

\subsection{Modification of Theorem 2.9 for the Case
of In\-teg\-ra\-tion Interval $[t, s]$ $(s\in (t, T])$ 
of Iterated Stra\-to\-no\-vich Sto\-chas\-tic Integrals of Multiplicity 4}

Let us prove the following theorem.

{\bf Theorem 2.22}\ \cite{arxiv-5}.
{\it Suppose that
$\{\phi_j(x)\}_{j=0}^{\infty}$ is a complete orthonormal
system of Legendre polynomials or trigonometric functions
in the space $L_2([t, T])$.
Then$,$ for the iterated 
Stratonovich stochastic integral of fourth multiplicity
$$
J^{*}[\psi^{(4)}]_{s,t}=
{\int\limits_t^{*}}^s
{\int\limits_t^{*}}^{t_4}
{\int\limits_t^{*}}^{t_3}
{\int\limits_t^{*}}^{t_2}
d{\bf w}_{t_1}^{(i_1)}
d{\bf w}_{t_2}^{(i_2)}d{\bf w}_{t_3}^{(i_3)}d{\bf w}_{t_4}^{(i_4)}\ \ \ 
(i_1, i_2, i_3, i_4=0, 1,\ldots,m)
$$
the following 
expansion 
$$
J^{*}[\psi^{(4)}]_{s,t}=
\hbox{\vtop{\offinterlineskip\halign{
\hfil#\hfil\cr
{\rm l.i.m.}\cr
$\stackrel{}{{}_{p\to \infty}}$\cr
}} }
\sum\limits_{j_1, j_2, j_3, j_4=0}^{p}
C_{j_4 j_3 j_2 j_1}(s)
\zeta_{j_1}^{(i_1)}\zeta_{j_2}^{(i_2)}\zeta_{j_3}^{(i_3)}
\zeta_{j_4}^{(i_4)}
$$
that converges in the mean-square sense is valid, where
$s\in (t,T]$ {\rm (}$s$ is fixed{\rm )},
$$
C_{j_4 j_3 j_2 j_1}(s)=\int\limits_t^s\phi_{j_4}(s_4)\int\limits_t^{s_4}
\phi_{j_3}(s_3)
\int\limits_t^{s_3}\phi_{j_2}(s_2)\int\limits_t^{s_2}\phi_{j_1}(s_1)
ds_1ds_2ds_3ds_4
$$
and
$$
\zeta_{j}^{(i)}=
\int\limits_t^T \phi_{j}(\tau) d{\bf w}_{\tau}^{(i)}
$$ 
are independent standard Gaussian random variables for various 
$i$ or $j$ {\rm (}in the case when $i\ne 0${\rm ),}
${\bf w}_{\tau}^{(i)}$
$(i=1,\ldots,m)$ are independent standard Wiener processes and 
${\bf w}_{\tau}^{(0)}=\tau.$}

{\bf Proof.} The case $s=T$ is considered in Theorem 2.9. 
Below we consider the case $s\in (t, T).$ The relation (\ref{cas1}) (in the case 
when $p_1=\ldots=p_4=p\to \infty$) implies that
$$
\hbox{\vtop{\offinterlineskip\halign{
\hfil#\hfil\cr
{\rm l.i.m.}\cr
$\stackrel{}{{}_{p\to \infty}}$\cr
}} }
\sum\limits_{j_1, j_2, j_3, j_4=0}^{p}
C_{j_4 j_3 j_2 j_1}(s)
\zeta_{j_1}^{(i_1)}\zeta_{j_2}^{(i_2)}\zeta_{j_3}^{(i_3)}
\zeta_{j_4}^{(i_4)}=
J[\psi^{(4)}]_{s,t}+
$$

\vspace{-5mm}
$$
+{\bf 1}_{\{i_1=i_2\ne 0\}}A_1^{(i_3i_4)}(s)
+{\bf 1}_{\{i_1=i_3\ne 0\}}A_2^{(i_2i_4)}(s)+
{\bf 1}_{\{i_1=i_4\ne 0\}}A_3^{(i_2i_3)}(s)+
{\bf 1}_{\{i_2=i_3\ne 0\}}A_4^{(i_1i_4)}(s)+
$$

\vspace{-7mm}
$$
+
{\bf 1}_{\{i_2=i_4\ne 0\}}A_5^{(i_1i_3)}(s)
+{\bf 1}_{\{i_3=i_4\ne 0\}}A_6^{(i_1i_2)}(s)-
{\bf 1}_{\{i_1=i_2\ne 0\}}
{\bf 1}_{\{i_3=i_4\ne 0\}}B_1(s)-
$$

\vspace{-7mm}
\begin{equation}
\label{cas2}
~~~~ -{\bf 1}_{\{i_1=i_3\ne 0\}}
{\bf 1}_{\{i_2=i_4\ne 0\}}B_2(s)-
{\bf 1}_{\{i_1=i_4\ne 0\}}
{\bf 1}_{\{i_2=i_3\ne 0\}}B_3(s),
\end{equation}

\vspace{3.5mm}
\noindent
where
$J[\psi^{(4)}]_{s,t}$ has the form {\rm (\ref{opr22})}
for $\psi_1(\tau),\ldots,\psi_4(\tau)\equiv 1$ and
$i_1,\ldots,i_4=0, 1,\ldots,m,$
$$
A_1^{(i_3i_4)}(s)=
\hbox{\vtop{\offinterlineskip\halign{
\hfil#\hfil\cr
{\rm l.i.m.}\cr
$\stackrel{}{{}_{p\to \infty}}$\cr
}} }
\sum\limits_{j_4, j_3, j_1=0}^{p}
C_{j_4 j_3 j_1 j_1}(s)\zeta_{j_3}^{(i_3)}
\zeta_{j_4}^{(i_4)},
$$
$$
A_2^{(i_2i_4)}(s)=
\hbox{\vtop{\offinterlineskip\halign{
\hfil#\hfil\cr
{\rm l.i.m.}\cr
$\stackrel{}{{}_{p\to \infty}}$\cr
}} }
\sum\limits_{j_4, j_3, j_2=0}^{p}
C_{j_4 j_3 j_2 j_3}(s)\zeta_{j_2}^{(i_2)}
\zeta_{j_4}^{(i_4)},
$$
$$
A_3^{(i_2i_3)}(s)=
\hbox{\vtop{\offinterlineskip\halign{
\hfil#\hfil\cr
{\rm l.i.m.}\cr
$\stackrel{}{{}_{p\to \infty}}$\cr
}} }
\sum\limits_{j_4, j_3, j_2=0}^{p}
C_{j_4 j_3 j_2 j_4}(s)\zeta_{j_2}^{(i_2)}
\zeta_{j_3}^{(i_3)},
$$
$$
A_4^{(i_1i_4)}(s)=
\hbox{\vtop{\offinterlineskip\halign{
\hfil#\hfil\cr
{\rm l.i.m.}\cr
$\stackrel{}{{}_{p\to \infty}}$\cr
}} }
\sum\limits_{j_4, j_3, j_1=0}^{p}
C_{j_4 j_3 j_3 j_1}(s)\zeta_{j_1}^{(i_1)}
\zeta_{j_4}^{(i_4)},
$$
$$
A_5^{(i_1i_3)}(s)=
\hbox{\vtop{\offinterlineskip\halign{
\hfil#\hfil\cr
{\rm l.i.m.}\cr
$\stackrel{}{{}_{p\to \infty}}$\cr
}} }
\sum\limits_{j_4, j_3, j_1=0}^{p}
C_{j_4 j_3 j_4 j_1}(s)\zeta_{j_1}^{(i_1)}
\zeta_{j_3}^{(i_3)},
$$
$$
A_6^{(i_1i_2)}(s)=
\hbox{\vtop{\offinterlineskip\halign{
\hfil#\hfil\cr
{\rm l.i.m.}\cr
$\stackrel{}{{}_{p\to \infty}}$\cr
}} }
\sum\limits_{j_3, j_2, j_1=0}^{p}
C_{j_3 j_3 j_2 j_1}(s)\zeta_{j_1}^{(i_1)}
\zeta_{j_2}^{(i_2)},
$$
$$
B_1(s)=
\hbox{\vtop{\offinterlineskip\halign{
\hfil#\hfil\cr
{\rm lim}\cr
$\stackrel{}{{}_{p\to \infty}}$\cr
}} }
\sum\limits_{j_1, j_4=0}^{p}
C_{j_4 j_4 j_1 j_1}(s),\ \ \
B_2(s)=
\hbox{\vtop{\offinterlineskip\halign{
\hfil#\hfil\cr
{\rm lim}\cr
$\stackrel{}{{}_{p\to \infty}}$\cr
}} }
\sum\limits_{j_4, j_3=0}^{p}
C_{j_3 j_4 j_3 j_4}(s),
$$
$$
B_3(s)=
\hbox{\vtop{\offinterlineskip\halign{
\hfil#\hfil\cr
{\rm lim}\cr
$\stackrel{}{{}_{p\to \infty}}$\cr
}} }
\sum\limits_{j_4, j_3=0}^{p}
C_{j_4 j_3 j_3 j_4}(s).
$$

Using the integration order replacement in Riemann integrals,
Theorem 1.11 for $k=2$ (see (\ref{a2xxx})) and
(\ref{5tzzz}), 
Parseval's equality and the integration order replacement
technique for It\^{o} stochastic integrals (see Chapter 3) 
\cite{1}-\cite{12aa}, \cite{old-art-2}, \cite{vini},
\cite{arxiv-25} or It\^{o}'s formula, we obtain
(see the derivation of the formula (\ref{otiteee1}))
$$
A_1^{(i_3i_4)}(s)=
\frac{1}{2}\int\limits_t^s\int\limits_t^{\tau}\int\limits_t^{s_1}ds_2
d{\bf w}_{s_1}^{(i_3)}
d{\bf w}_{\tau}^{(i_4)}+
$$
\begin{equation}
\label{cas0}
+
\frac{1}{4}{\bf 1}_{\{i_3=i_4\ne 0\}}
\int\limits_t^s(s_1-t)ds_1
- \Delta_1^{(i_3i_4)}(s)\ \ \ \hbox{w.~p.~1,}
\end{equation}

\noindent
where
$$
\Delta_1^{(i_3i_4)}(s)=
\hbox{\vtop{\offinterlineskip\halign{
\hfil#\hfil\cr
{\rm l.i.m.}\cr
$\stackrel{}{{}_{p\to \infty}}$\cr
}} }
\sum\limits_{j_3, j_4=0}^{p}
a_{j_4 j_3}^p (s) \zeta_{j_3}^{(i_3)}
\zeta_{j_4}^{(i_4)},
$$
\begin{equation}
\label{cas3}
~~~~~~~~~ a_{j_4 j_3}^p (s)=
\frac{1}{2}\int\limits_t^s\phi_{j_4}(\tau)\int\limits_t^{\tau}\phi_{j_3}(s_1)
\sum\limits_{j_1=p+1}^{\infty}\left(\int\limits_t^{s_1}
\phi_{j_1}(s_2)ds_2\right)^2ds_1d\tau.
\end{equation}

\vspace{4mm}

Let us consider $A_2^{(i_2i_4)}(s)$ 
(see the derivation of the formula (\ref{otit999}))
\begin{equation}
\label{cas800}
~~~~~~ A_2^{(i_2i_4)}(s)=
-\Delta_2^{(i_2i_4)}(s)+\Delta_1^{(i_2i_4)}(s)+\Delta_3^{(i_2i_4)}(s)\ \ \
\hbox{w.~p.~1,}
\end{equation}

\noindent
where
$$
\Delta_2^{(i_2i_4)}(s)=
\hbox{\vtop{\offinterlineskip\halign{
\hfil#\hfil\cr
{\rm l.i.m.}\cr
$\stackrel{}{{}_{p\to \infty}}$\cr
}} }
\sum\limits_{j_4, j_2=0}^{p}
b_{j_4 j_2}^p (s) \zeta_{j_2}^{(i_2)}
\zeta_{j_4}^{(i_4)},
$$
$$
\Delta_3^{(i_2i_4)}(s)=
\hbox{\vtop{\offinterlineskip\halign{
\hfil#\hfil\cr
{\rm l.i.m.}\cr
$\stackrel{}{{}_{p\to \infty}}$\cr
}} }
\sum\limits_{j_4, j_2=0}^{p}
c_{j_4 j_2}^p (s)\zeta_{j_2}^{(i_2)}
\zeta_{j_4}^{(i_4)},
$$
$$
b_{j_4 j_2}^p (s)=
\frac{1}{2}\int\limits_t^s\phi_{j_4}(\tau)
\sum\limits_{j_3=p+1}^{\infty}\left(\int\limits_t^{\tau}
\phi_{j_3}(s_1)ds_1\right)^2\int\limits_t^{\tau}\phi_{j_2}(s_1)ds_1d\tau,
$$
$$
c_{j_4 j_2}^p (s)=
\frac{1}{2}\int\limits_t^s\phi_{j_4}(\tau)\int\limits_t^{\tau}\phi_{j_2}(s_3)
\sum\limits_{j_3=p+1}^{\infty}\left(\int\limits_{s_3}^{\tau}
\phi_{j_3}(s_1)ds_1\right)^2ds_3d\tau.
$$

\vspace{4mm}

Further, we have w.~p.~1
(see the derivation of the formula (\ref{otit9999}))
\begin{equation}
\label{cas5}
\vspace{2mm}
~~~~~~~~~ A_5^{(i_1i_3)}(s)
=-\Delta_4^{(i_1i_3)}(s)+\Delta_5^{(i_1i_3)}(s)+\Delta_6^{(i_1i_3)}(s)\ \ \
\hbox{w.~p.~1,}
\end{equation}
where
$$
\Delta_4^{(i_1i_3)}(s)=
\hbox{\vtop{\offinterlineskip\halign{
\hfil#\hfil\cr
{\rm l.i.m.}\cr
$\stackrel{}{{}_{p\to \infty}}$\cr
}} }
\sum\limits_{j_3, j_1=0}^{p}
d_{j_3 j_1}^p (s)\zeta_{j_1}^{(i_1)}
\zeta_{j_3}^{(i_3)},
$$
$$
\Delta_5^{(i_1i_3)}(s)=
\hbox{\vtop{\offinterlineskip\halign{
\hfil#\hfil\cr
{\rm l.i.m.}\cr
$\stackrel{}{{}_{p\to \infty}}$\cr
}} }
\sum\limits_{j_3, j_1=0}^{p}
e_{j_3 j_1}^p (s)\zeta_{j_1}^{(i_1)}
\zeta_{j_3}^{(i_3)},
$$
$$
\Delta_6^{(i_1i_3)}(s)=
\hbox{\vtop{\offinterlineskip\halign{
\hfil#\hfil\cr
{\rm l.i.m.}\cr
$\stackrel{}{{}_{p\to \infty}}$\cr
}} }
\sum\limits_{j_3, j_1=0}^{p}
f_{j_3 j_1}^p (s)\zeta_{j_1}^{(i_1)}
\zeta_{j_3}^{(i_3)},
$$
$$
d_{j_3 j_1}^p (s)=
\frac{1}{2}\int\limits_t^s\phi_{j_1}(s_3)
\sum\limits_{j_4=p+1}^{\infty}\left(\int\limits_{s_3}^{s}
\phi_{j_4}(\tau)d\tau\right)^2\int\limits_{s_3}^s\phi_{j_3}(\tau)d\tau ds_3,
$$
$$
e_{j_3 j_1}^p (s)=
\frac{1}{2}\int\limits_t^s\phi_{j_1}(s_3)\int\limits_{s_3}^s\phi_{j_3}(\tau)
\sum\limits_{j_4=p+1}^{\infty}\left(\int\limits_{s_3}^{\tau}
\phi_{j_4}(s_1)ds_1\right)^2d\tau ds_3,
$$
$$
f_{j_3 j_1}^p (s)=
\frac{1}{2}\int\limits_t^s\phi_{j_1}(s_3)\int\limits_{s_3}^s\phi_{j_3}(s_2)
\sum\limits_{j_4=p+1}^{\infty}\left(\int\limits_{s_2}^{s}
\phi_{j_4}(s_1)ds_1\right)^2ds_2ds_3=
$$
$$
=
\frac{1}{2}\int\limits_t^s\phi_{j_3}(s_2)
\sum\limits_{j_4=p+1}^{\infty}\left(\int\limits_{s_2}^{s}
\phi_{j_4}(s_1)ds_1\right)^2
\int\limits_t^{s_2}\phi_{j_1}(s_3)ds_3ds_2.
$$

\vspace{2mm}

Let us consider $A_4^{(i_1i_4)}(s)$
(see the derivation of the formula (\ref{otit555}))
\begin{equation}
\label{cas7}
~~~~~~ A_4^{(i_1i_4)}(s)
=\frac{1}{2}\int\limits_t^s\int\limits_t^{s_2}\int\limits_t^{s_1}
d{\bf w}_{\tau}^{(i_1)}ds_1
d{\bf w}_{s_2}^{(i_4)} - \Delta_3^{(i_1i_4)}(s)\ \ \ \hbox{w.~p.~1.}
\end{equation}

\vspace{4mm}

Moreover
(see the derivation of the formula (\ref{otit001})),
$$
A_6^{(i_1i_2)}(s)=
\frac{1}{2}\int\limits_t^s\int\limits_t^{s_1}\int\limits_t^{s_2}
d{\bf w}_{\tau}^{(i_1)}
d{\bf w}_{s_2}^{(i_2)}ds_1+
$$
\begin{equation}
\label{cas8}
+
\frac{1}{4}{\bf 1}_{\{i_1=i_2\ne 0\}}
\int\limits_t^s(s-s_2)ds_2
- \Delta_6^{(i_1i_2)}(s)\ \ \ \hbox{w.~p.~1.}
\end{equation}

\vspace{3mm}

Further, 
we have w.~p.~1 (see the derivation of the formula (\ref{strange500}))

\vspace{-2mm}
$$
A_3^{(i_2i_3)}(s)+A_5^{(i_2i_3)}(s)=
$$
\begin{equation}
\label{strange501}
=\hbox{\vtop{\offinterlineskip\halign{
\hfil#\hfil\cr
{\rm l.i.m.}\cr
$\stackrel{}{{}_{p\to \infty}}$\cr
}} }
\sum\limits_{j_4, j_3, j_2=0}^{p}
\int\limits_t^s\phi_{j_3}(s_1)\int\limits_{t}^{s_1}\phi_{j_2}(s_2)ds_2
\int\limits_{t}^{s_1}\phi_{j_4}(s_3)ds_3\int\limits_{s_1}^{s}\phi_{j_4}(\tau)
d\tau ds_1
\zeta_{j_2}^{(i_2)}
\zeta_{j_3}^{(i_3)}.
\end{equation}

\vspace{3mm}

Using (\ref{strange501}) and the generalized Parseval equality, we obtain w.~p.~1

\vspace{-2mm}
$$
A_3^{(i_2i_3)}(s)+A_5^{(i_2i_3)}(s)=
$$
$$
=\hbox{\vtop{\offinterlineskip\halign{
\hfil#\hfil\cr
{\rm l.i.m.}\cr
$\stackrel{}{{}_{p\to \infty}}$\cr
}} }
\sum\limits_{j_3, j_2=0}^{p}
\int\limits_t^s\phi_{j_3}(s_1)\int\limits_{t}^{s_1}\phi_{j_2}(s_2)ds_2
\sum\limits_{j_4=0}^{p}\int\limits_{t}^{s_1}\phi_{j_4}(s_3)ds_3\int\limits_{s_1}^{s}\phi_{j_4}(\tau)
d\tau ds_1
\zeta_{j_2}^{(i_2)}
\zeta_{j_3}^{(i_3)}=
$$
$$
=-\hbox{\vtop{\offinterlineskip\halign{
\hfil#\hfil\cr
{\rm l.i.m.}\cr
$\stackrel{}{{}_{p\to \infty}}$\cr
}} }
\sum\limits_{j_3, j_2=0}^{p}
\int\limits_t^s\phi_{j_3}(s_1)\int\limits_{t}^{s_1}\phi_{j_2}(s_2)ds_2
\sum\limits_{j_4=p+1}^{\infty}\int\limits_{t}^{s_1}\phi_{j_4}(s_3)ds_3\int\limits_{s_1}^{s}\phi_{j_4}(\tau)
d\tau ds_1\times
$$
$$
\times
\zeta_{j_2}^{(i_2)}
\zeta_{j_3}^{(i_3)}=
$$

\vspace{-2mm}
\begin{equation}
\label{strange502}
=\Delta_6^{(i_2i_3)}(s)+\Delta_2^{(i_2i_3)}(s)-\Delta_9^{(i_2i_3)}(s),
\end{equation}

\vspace{2mm}
\noindent
where
$$
\Delta_9^{(i_2i_3)}(s)=
\hbox{\vtop{\offinterlineskip\halign{
\hfil#\hfil\cr
{\rm l.i.m.}\cr
$\stackrel{}{{}_{p\to \infty}}$\cr
}} }
\sum\limits_{j_3, j_2=0}^{p}
q_{j_2 j_3}^p (s)\zeta_{j_2}^{(i_3)}
\zeta_{j_3}^{(i_3)},
$$
$$
q_{j_2 j_3}^p (s)=
\frac{1}{2}\int\limits_t^s\phi_{j_3}(s_1)
\int\limits_{t}^{s_1}\phi_{j_2}(s_2)ds_2 ds_1
\sum\limits_{j_4=p+1}^{\infty}\left(\int\limits_{t}^{s}
\phi_{j_4}(\tau)d\tau\right)^2.
$$

\vspace{2mm}

From (\ref{cas5}) and (\ref{strange502}) we get
\begin{equation}
\label{strange505}
~~~~A_3^{(i_2i_3)}(s)
=\Delta_2^{(i_2i_3)}(s)+
\Delta_4^{(i_2i_3)}(s)-
\Delta_5^{(i_2i_3)}(s)-
\Delta_9^{(i_2i_3)}(s)
\ \ \ \hbox{w.~p.~1.}
\end{equation}

\vspace{3mm}

Let us consider $B_1(s), B_2(s), B_3(s).$
We have 
(see the derivation of the formulas (\ref{otit239}), (\ref{otit990}))
\begin{equation}
\label{cas9}
B_1(s)
=
\frac{1}{4}\int\limits_t^s(s_1-t)ds_1
-\lim_{p\to\infty}\sum\limits_{j_4=0}^{p}
a_{j_4j_4}^p (s),
\end{equation}
\begin{equation}
\label{cas10}
~~~~~~~~~~~B_2(s)
=\lim_{p\to\infty}\sum\limits_{j_3=0}^p a_{j_3j_3}^p (s)
+\lim_{p\to\infty}\sum\limits_{j_3=0}^p c_{j_3j_3}^p (s)
-\lim_{p\to\infty}\sum\limits_{j_3=0}^p b_{j_3j_3}^p (s).
\end{equation}

\vspace{4mm}

Moreover
(see the derivation of the formula (\ref{cas1000})),

\vspace{-2mm}
$$
B_2(s)+B_3(s)=
$$
\begin{equation}
\label{strange506}
~~~~=\hbox{\vtop{\offinterlineskip\halign{
\hfil#\hfil\cr
{\rm lim}\cr
$\stackrel{}{{}_{p\to \infty}}$\cr
}} }
\sum\limits_{j_4=0}^{p}
\int\limits_t^s\phi_{j_4}(s_1)\int\limits_{t}^{s_1}\phi_{j_4}(s_2)ds_2
\sum\limits_{j_3=0}^{p}
\int\limits_{t}^{s_1}\phi_{j_3}(s_3)ds_3\int\limits_{s_1}^{s}\phi_{j_3}(\tau)
d\tau ds_1.
\end{equation}

\vspace{2mm}

Using (\ref{strange506}) and the generalized Parseval equality, we obtain 

\vspace{-2mm}

$$
B_2(s)+B_3(s)=
$$
$$
=-\hbox{\vtop{\offinterlineskip\halign{
\hfil#\hfil\cr
{\rm lim}\cr
$\stackrel{}{{}_{p\to \infty}}$\cr
}} }
\sum\limits_{j_4=0}^{p}
\int\limits_t^s\phi_{j_4}(s_1)\int\limits_{t}^{s_1}\phi_{j_4}(s_2)ds_2
\sum\limits_{j_3=p+1}^{\infty}
\int\limits_{t}^{s_1}\phi_{j_3}(s_3)ds_3\int\limits_{s_1}^{s}\phi_{j_3}(\tau)
d\tau ds_1=
$$
\begin{equation}
\label{strange507}
~~~~~~~~=\lim_{p\to\infty}\sum\limits_{j_4=0}^p f_{j_4j_4}^p (s)+
\lim_{p\to\infty}\sum\limits_{j_4=0}^p b_{j_4j_4}^p (s)
-\lim_{p\to\infty}\sum\limits_{j_4=0}^p q_{j_4j_4}^p (s).
\end{equation}

\vspace{3mm}

Combining (\ref{cas10}) and (\ref{strange507}), we have
$$
B_3(s)=2\lim_{p\to\infty}\sum\limits_{j_4=0}^p b_{j_4j_4}^p (s)+
\lim_{p\to\infty}\sum\limits_{j_4=0}^p f_{j_4j_4}^p (s)-
\lim_{p\to\infty}\sum\limits_{j_4=0}^p c_{j_4j_4}^p (s)-
$$
\begin{equation}
\label{strange508}
-\lim_{p\to\infty}\sum\limits_{j_4=0}^p a_{j_4j_4}^p (s)-
\lim_{p\to\infty}\sum\limits_{j_4=0}^p q_{j_4j_4}^p (s).
\end{equation}

\vspace{2mm}

After substituting the relations (\ref{cas0}), (\ref{cas800})--(\ref{cas8}),
(\ref{strange505})--(\ref{cas10}), (\ref{strange508})
into (\ref{cas2}), we obtain
$$
\hbox{\vtop{\offinterlineskip\halign{
\hfil#\hfil\cr
{\rm l.i.m.}\cr
$\stackrel{}{{}_{p\to \infty}}$\cr
}} }
\sum\limits_{j_1, j_2, j_3, j_4=0}^{p}
C_{j_4 j_3 j_2 j_1}(s)
\zeta_{j_1}^{(i_1)}\zeta_{j_2}^{(i_2)}\zeta_{j_3}^{(i_3)}
\zeta_{j_4}^{(i_4)}=
$$
$$
=
J[\psi^{(4)}]_{s,t}+
\frac{1}{2}{\bf 1}_{\{i_1=i_2\ne 0\}}
\int\limits_t^s\int\limits_t^{\tau}\int\limits_t^{s_1}ds_2
d{\bf w}_{s_1}^{(i_3)}
d{\bf w}_{\tau}^{(i_4)}+
$$
$$
+\frac{1}{2}{\bf 1}_{\{i_2=i_3\ne 0\}}
\int\limits_t^s\int\limits_t^{s_2}\int\limits_t^{s_1}
d{\bf w}_{\tau}^{(i_1)}ds_1
d{\bf w}_{s_2}^{(i_4)}
+\frac{1}{2}{\bf 1}_{\{i_3=i_4\ne 0\}}
\int\limits_t^s\int\limits_t^{s_1}\int\limits_t^{s_2}
d{\bf w}_{\tau}^{(i_1)}
d{\bf w}_{s_2}^{(i_2)}ds_1+
$$
$$
+\frac{1}{4}{\bf 1}_{\{i_1=i_2\ne 0\}}
{\bf 1}_{\{i_3=i_4\ne 0\}}
\int\limits_t^T\int\limits_t^{s_1}ds_2
ds_1 + R(s) = J^{*}[\psi^{(4)}]_{s,t}+
$$

\begin{equation}
\label{cas99}
+R(s)\ \ \  \hbox{w.~p.~1,}
\end{equation}

\noindent
where
$$
R(s)=-{\bf 1}_{\{i_1=i_2\ne 0\}}\Delta_1^{(i_3i_4)}(s)+
{\bf 1}_{\{i_1=i_3\ne 0\}}\left(
-\Delta_2^{(i_2i_4)}(s)
+\Delta_1^{(i_2i_4)}(s)
+\Delta_3^{(i_2i_4)}(s)\right)+
$$
$$
+{\bf 1}_{\{i_1=i_4\ne 0\}}\left(
\Delta_2^{(i_2i_3)}(s)+
\Delta_4^{(i_2i_3)}(s)-
\Delta_5^{(i_2i_3)}(s)-
\Delta_9^{(i_2i_3)}(s)\right)-
{\bf 1}_{\{i_2=i_3\ne 0\}}\Delta_3^{(i_1i_4)}(s)+
$$
$$
+{\bf 1}_{\{i_2=i_4\ne 0\}}
\left(-\Delta_4^{(i_1i_3)}(s)
+\Delta_5^{(i_1i_3)}(s)
+\Delta_6^{(i_1i_3)}(s)\right)-
{\bf 1}_{\{i_3=i_4\ne 0\}}\Delta_6^{(i_1i_2)}(s)-
$$
$$
-
{\bf 1}_{\{i_1=i_3\ne 0\}}
{\bf 1}_{\{i_2=i_4\ne 0\}}\Biggl(
\lim_{p\to\infty}\sum\limits_{j_3=0}^p a_{j_3j_3}^p (s)
+\lim_{p\to\infty}\sum\limits_{j_3=0}^p c_{j_3j_3}^p (s)
-\lim_{p\to\infty}\sum\limits_{j_3=0}^p b_{j_3j_3}^p (s)\Biggr)-
$$
$$
-{\bf 1}_{\{i_1=i_4\ne 0\}}
{\bf 1}_{\{i_2=i_3\ne 0\}}
\Biggl(
2\lim_{p\to\infty}\sum\limits_{j_4=0}^p b_{j_4j_4}^p (s)+
\lim_{p\to\infty}\sum\limits_{j_4=0}^p f_{j_4j_4}^p (s)-
\lim_{p\to\infty}\sum\limits_{j_4=0}^p c_{j_4j_4}^p (s)-\Biggr.
$$
$$
\Biggl.-\lim_{p\to\infty}\sum\limits_{j_4=0}^p a_{j_4j_4}^p (s)-
\lim_{p\to\infty}\sum\limits_{j_4=0}^p q_{j_4j_4}^p (s)\Biggr)+
$$
\begin{equation}
\label{cas100}
+{\bf 1}_{\{i_1=i_2\ne 0\}}
{\bf 1}_{\{i_3=i_4\ne 0\}}
\lim_{p\to\infty}\sum\limits_{j_3=0}^p a_{j_3j_3}^p (s).
\end{equation}

\vspace{4mm}

Let us prove that
\begin{equation}
\label{cas101}
R(s)=0\ \ \ \hbox{w.~p.~1}.
\end{equation}

\vspace{2mm}

Consider the case of Legendre polynomials.
First, we prove that 
\begin{equation}
\label{cas103}
\Delta_1^{(i_3i_4)}(s)=0\ \ \  \hbox{w.~p.~1.}
\end{equation}

\vspace{3mm}

We have
$$
a_{j_4j_3}^p (s)=\frac{(T-t)^2\sqrt{(2j_4+1)(2j_3+1)}}{32}\times
$$
$$
\times
\int\limits_{-1}^{z(s)} P_{j_4}(y) \int\limits_{-1}^y
P_{j_3}(y_1)\sum\limits_{j_1=p+1}^{\infty}(2j_1+1)
\left(\int\limits_{-1}^{y_1}P_{j_1}(y_2)dy_2\right)^2 dy_1dy=
$$

\vspace{3mm}
$$
=\frac{(T-t)^2\sqrt{(2j_4+1)(2j_3+1)}}{32}\times
$$
$$
\times
\int\limits_{-1}^{z(s)} P_{j_3}(y_1) 
\sum\limits_{j_1=p+1}^{\infty}\frac{1}{2j_1+1}
\left(P_{j_1+1}(y_1)-P_{j_1-1}(y_1)\right)^2
\int\limits_{y_1}^{z(s)}P_{j_4}(y)dy dy_1=
$$

\vspace{3mm}
$$
=\frac{(T-t)^2\sqrt{2j_3+1}}{32\sqrt{2j_4+1}}\times
$$
$$
\times
\int\limits_{-1}^{z(s)}
P_{j_3}(y_1) \left(\left(
P_{j_4+1}(z(s))-P_{j_4-1}(z(s))\right)-
\left(P_{j_4+1}(y_1)-P_{j_4-1}(y_1)\right)\right)\times
$$
$$
\times
\sum\limits_{j_1=p+1}^{\infty}\frac{1}{2j_1+1}
\left(P_{j_1+1}(y_1)-P_{j_1-1}(y_1)\right)^2 dy_1
$$

\vspace{3mm}
\noindent
if $j_4\ne 0$ and
$$
a_{j_4j_3}^p (s)=\frac{(T-t)^2\sqrt{2j_3+1}}{32}\times
$$
$$
\times
\int\limits_{-1}^{z(s)} P_{j_3}(y_1) (z(s)-y_1)
\sum\limits_{j_1=p+1}^{\infty}\frac{1}{2j_1+1}
\left(P_{j_1+1}(y_1)-P_{j_1-1}(y_1)\right)^2
dy_1
$$

\vspace{3mm}
\noindent
if $j_4=0,$ where $z(s)$ is defined by (\ref{zz1}).

We assume that $s\in (t, T)$ $(z(s)\ne \pm 1)$ since the case
$s=T$ has already been considered in Theorem 2.9.
Now the further proof of the equality (\ref{cas103}) 
is completely analogous to the proof of the
equality (\ref{otitf14}).

It is not difficult to see that the formulas
\begin{equation}
\label{cas300}
\Delta_2^{(i_2i_4)}(s)=0,\ \ \ \Delta_4^{(i_1i_3)}(s)=0,\ \ \
\Delta_6^{(i_1i_3)}(s)=0,\ \ \ \Delta_9^{(i_2i_3)}(s)=0\ \ \ \hbox{w.~p.~1}
\end{equation}

\noindent
can be proved similarly with the 
proof of (\ref{cas103}).

Moreover, the relations
\begin{equation}
\label{cas301}
\lim\limits_{p\to\infty}
\sum\limits_{j_3=0}^p a_{j_3j_3}^p (s)=0,\ \ \ 
\lim\limits_{p\to\infty}
\sum\limits_{j_3=0}^p b_{j_3j_3}^p (s)=0,
\end{equation}
\begin{equation}
\label{cas301x}
\lim\limits_{p\to\infty}
\sum\limits_{j_3=0}^p f_{j_3j_3}^p (s)=0,\ \ \ 
\lim\limits_{p\to\infty}
\sum\limits_{j_3=0}^p q_{j_3j_3}^p (s)=0
\end{equation}

\vspace{2mm}
\noindent
can also be proved analogously with (\ref{20177}), (\ref{cas79}).

Let us consider $\Delta_3^{(i_2i_4)}(s)$ and prove that
\begin{equation}
\label{cas303a}
\Delta_3^{(i_2i_4)}(s)=0\ \ \ \hbox{w.~p.~1}.
\end{equation}

We have
\begin{equation}
\label{cas303}
~~~~~~~~~\Delta_3^{(i_2i_4)}(s)=\Delta_4^{(i_2i_4)}(s)+
\Delta_6^{(i_2i_4)}(s)-\Delta_7^{(i_2i_4)}(s)=
-\Delta_7^{(i_2i_4)}(s)
\end{equation}

\vspace{2mm}
\noindent
w.~p.~1, where
$$
\Delta_7^{(i_2i_4)}(s)=
\hbox{\vtop{\offinterlineskip\halign{
\hfil#\hfil\cr
{\rm l.i.m.}\cr
$\stackrel{}{{}_{p\to \infty}}$\cr
}} }
\sum\limits_{j_2, j_4=0}^{p}
g_{j_4 j_2}^p (s)\zeta_{j_2}^{(i_2)}
\zeta_{j_4}^{(i_4)},
$$

\vspace{-2mm}
$$
g_{j_4 j_2}^p (s)=
\int\limits_t^s\phi_{j_4}(\tau)\int\limits_{t}^{\tau}\phi_{j_2}(s_1)
\sum\limits_{j_1=p+1}^{\infty}\left(\int\limits_{s_1}^{s}
\phi_{j_1}(s_2)ds_2\int\limits_{\tau}^{s}
\phi_{j_1}(s_2)ds_2\right)  
ds_1d\tau.
$$

\vspace{3mm}

Note that (see (\ref{otitddd}))
\begin{equation}
\label{cas304}
g_{j_4 j_4}^p (s)=\sum\limits_{j_1=p+1}^{\infty}
\frac{1}{2}\left(
\int\limits_t^s\phi_{j_4}(\tau)\int\limits_{\tau}^s\phi_{j_1}(s_2)
ds_2d\tau\right)^2.
\end{equation}

\vspace{2mm}

The proof of (\ref{cas303a}) for the case $i_2=i_4\ne 0$
differs from the proof of the equality
$$
\Delta_3^{(i_2i_4)}=0\ \ \ \hbox{w.~p.~1}
$$
for the case $i_2=i_4\ne 0$ (see the proof of Theorem 2.9).
In our case we will use Parseval's equality 
instead of the orthogonality
property of the Legendre polynomials.

Using the Parseval equality, we obtain
$$
\sum\limits_{j_4=0}^p
g_{j_4 j_4}^p (s)=
\sum\limits_{j_4=0}^p
\sum\limits_{j_1=p+1}^{\infty}
\frac{1}{2}\left(
\int\limits_t^s\phi_{j_4}(\tau)\int\limits_{\tau}^s\phi_{j_1}(s_2)
ds_2d\tau\right)^2=
$$
$$
=
\sum\limits_{j_4=0}^p
\sum\limits_{j_1=p+1}^{\infty}
\frac{1}{2}\left(
\int\limits_t^s\phi_{j_4}(\tau)\left(
\int\limits_{t}^s\phi_{j_1}(s_2)
ds_2-
\int\limits_{t}^{\tau}\phi_{j_1}(s_2)
ds_2\right)
d\tau\right)^2\le
$$
$$
\le
\sum\limits_{j_4=0}^p
\left(
\int\limits_t^s\phi_{j_4}(\tau)d\tau\right)^2
\sum\limits_{j_1=p+1}^{\infty}\left(
\int\limits_t^s\phi_{j_1}(s_2)ds_2\right)^2+
$$
$$
+\sum\limits_{j_4=0}^p
\sum\limits_{j_1=p+1}^{\infty}
\left(
\int\limits_t^s\phi_{j_4}(\tau)
\int\limits_{t}^{\tau}\phi_{j_1}(s_2)
ds_2
d\tau\right)^2=
$$
$$
=
\sum\limits_{j_4=0}^p
\left(
\int\limits_t^T {\bf 1}_{\{\tau<s\}}\phi_{j_4}(\tau)d\tau\right)^2
\sum\limits_{j_1=p+1}^{\infty}\left(
\int\limits_t^s\phi_{j_1}(s_2)ds_2\right)^2+
$$
$$
+
\sum\limits_{j_1=p+1}^{\infty}\sum\limits_{j_4=0}^p
\left(
\int\limits_t^T {\bf 1}_{\{\tau<s\}}\phi_{j_4}(\tau)
\int\limits_{t}^{\tau}\phi_{j_1}(s_2)
ds_2
d\tau\right)^2\le
$$
$$
\le
\sum\limits_{j_4=0}^{\infty}
\left(
\int\limits_t^T {\bf 1}_{\{\tau<s\}}\phi_{j_4}(\tau)d\tau\right)^2
\sum\limits_{j_1=p+1}^{\infty}\left(
\int\limits_t^s\phi_{j_1}(s_2)ds_2\right)^2+
$$
$$
+
\sum\limits_{j_1=p+1}^{\infty}\sum\limits_{j_4=0}^{\infty}
\left(
\int\limits_t^T {\bf 1}_{\{\tau<s\}}\phi_{j_4}(\tau)
\int\limits_{t}^{\tau}\phi_{j_1}(s_2)
ds_2
d\tau\right)^2=
$$
$$
=
\int\limits_t^T \left({\bf 1}_{\{\tau<s\}}\right)^2
d\tau
\sum\limits_{j_1=p+1}^{\infty}\left(
\int\limits_t^s\phi_{j_1}(s_2)ds_2\right)^2+
$$
$$
+
\sum\limits_{j_1=p+1}^{\infty}
\int\limits_t^T \left({\bf 1}_{\{\tau<s\}}\right)^2
\left(\int\limits_{t}^{\tau}\phi_{j_1}(s_2)ds_2
\right)^2 d\tau=
$$
\begin{equation}
\label{cas400}
=
(s-t)
\sum\limits_{j_1=p+1}^{\infty}\left(
\int\limits_t^s\phi_{j_1}(s_2)ds_2\right)^2+
\sum\limits_{j_1=p+1}^{\infty}
\int\limits_t^s \
\left(\int\limits_{t}^{\tau}\phi_{j_1}(s_2)ds_2
\right)^2 d\tau.
\end{equation}

\vspace{2mm}

We assume that $s\in (t, T)$ $(z(s)\ne \pm 1)$ since the case
$s=T$ has already been considered in Theorem 2.9.
Then from (\ref{cas400}) and (\ref{ogo25}) we obtain 
\begin{equation}
\label{cas401}
0\le \sum\limits_{j_4=0}^p
g_{j_4 j_4}^p (s)\le \frac{C(s)}{p},
\end{equation}

\noindent
where constant $C(s)$ ($s$ is fixed) is independent of $p.$

Combining (\ref{otit987}) and (\ref{may2021}) with (\ref{6000}), we obtain for $j\in{\bf N}$
\begin{equation}
\label{after999x1}
~~~~~~~~~\left|\int\limits_{s_1}^{s}
\phi_{j}(\theta)d\theta\right|<\frac{K}{j^{1/2+m/4}} 
\Biggl(\frac{1}{(1-z^2(s))^{m/8}}+\frac{1}{(1-z^2(s_1))^{m/8}}\Biggr),
\end{equation}
where $s, s_1\in (t, T),$ $m=1$ or $m=2$, 
$z(s)$ is defined by (\ref{zz1}), constant $K$ does not depend on $j$.

Using the Parseval equality, we get
\begin{equation}
\label{after556}
\lim_{p_1\to\infty}\sum\limits_{j_4,j_2=0}^{p_1}
\left(g_{j_4 j_2}^p (s)\right)^2=\int\limits_{[t,T]^2}
\left(K_p(\tau,s_1,s)\right)^2 ds_1 d\tau=
\int\limits_t^s\int\limits_t^{\tau}\left(
F_p(\tau,s_1,s)\right)^2 ds_1 d\tau,
\end{equation}
where
$$
g_{j_4 j_2}^p (s)=
\int\limits_t^T {\bf 1}_{\{\tau<s\}}
\phi_{j_4}(\tau)\int\limits_{t}^{\tau}\phi_{j_2}(s_1)
F_p(\tau,s_1,s)  
ds_1d\tau=
$$
$$
=
\int\limits_{[t,T]^2}
K_p(\tau,s_1,s)\phi_{j_4}(\tau)\phi_{j_2}(s_1)ds_1 d\tau
$$

\vspace{1mm}
\noindent
is a coefficient of the double Fourier--Legendre series of the function
$$
K_p(\tau,s_1,s)={\bf 1}_{\{\tau<s\}}{\bf 1}_{\{s_1<\tau<s\}}F_p(\tau,s_1,s),
$$

\vspace{1mm}
\noindent
where
\begin{equation}
\label{cas700}
\sum\limits_{j_1=p+1}^{\infty}
\int\limits_{s_1}^{s}
\phi_{j_1}(s_2)ds_2\int\limits_{\tau}^{s}
\phi_{j_1}(s_2)ds_2\stackrel{\sf def}{=}F_p(\tau,s_1,s).
\end{equation}

\vspace{2mm}

From (\ref{after999x1}) for $m=1$ and $m=2$ we have
$$
\left|F_p(\tau,s_1,s)\right|<\sum\limits_{j_1=p+1}^{\infty}\frac{K_1}{(j_1)^{7/4}}
\Biggl(\frac{1}{(1-z^2(s))^{1/8}}+\frac{1}{(1-z^2(s_1))^{1/8}}\Biggr)\times
$$
$$
\times \Biggl(\frac{1}{(1-z^2(s))^{1/4}}+\frac{1}{(1-z^2(\tau))^{1/4}}\Biggr)\le
$$
\begin{equation}
\label{afterw34}
\le \frac{K_2}{p^{3/4}}
\Biggl(\frac{1}{(1-z^2(s))^{1/8}}+\frac{1}{(1-z^2(s_1))^{1/8}}\Biggr)
\Biggl(\frac{1}{(1-z^2(s))^{1/4}}+\frac{1}{(1-z^2(\tau))^{1/4}}\Biggr),
\end{equation}

\noindent
where $s, s_1, \tau\in(t, T),$ 
constant $K_2$ is independent of $p$
and we used the estimate (\ref{after1944}) in 
(\ref{afterw34}).

The relations (\ref{after556}) and (\ref{afterw34}) imply the estimate
\begin{equation}
\label{cas500}
\sum\limits_{j_2,j_4=0}^{p}
\left(g_{j_4 j_2}^p (s)\right)^2
\le \frac{C_1(s)}{p^{3/2}}
\end{equation}

\noindent
for the case
$s\in (t, T)$ or $z(s)\in (-1,1)$ (the case
$s=T$ has already been considered in Theorem 2.9),
where constant $C_1(s)$ ($s$ is fixed) does not depend on $p.$

Then from analogue of (\ref{riss7}) for $s\in (t, T)$
($s$ is fixed), (\ref{cas401}), and (\ref{cas500}) we have
$$
{\sf M}\left\{\left(\sum\limits_{j_2, j_4=0}^{p}
g_{j_4 j_2}^p (s)\zeta_{j_2}^{(i_2)}
\zeta_{j_4}^{(i_4)}\right)^2\right\}\le
\left(1+{\bf 1}_{\{i_2=i_4\ne 0\}}\right) 
\sum\limits_{j_2, j_4=0}^{p}
\left(g_{j_4 j_2}^p (s)\right)^2 +
$$

$$
+
{\bf 1}_{\{i_2=i_4\ne 0\}}\left(\sum\limits_{j_4=0}^{p}
g_{j_4 j_4}^p (s) \right)^2 \le
\frac{C_2(s)}{p^{3/2}}\ \to 0
$$

\vspace{2mm}
\noindent
if $p\to \infty,$
where constant $C_2(s)$ ($s$ is fixed) does not depend on $p.$
The equality (\ref{cas303a}) is proved.

Let us consider $\Delta_5^{(i_1 i_3)}(s)$
$$
\Delta_5^{(i_1 i_3)}(s)=\Delta_4^{(i_1 i_3)}(s)+
\Delta_6^{(i_1 i_3)}(s)-\Delta_8^{(i_1 i_3)}(s)\ \ \ \hbox{\rm w.~p.~1,}
$$
where 
$$
\Delta_8^{(i_1i_3)}(s)=
\hbox{\vtop{\offinterlineskip\halign{
\hfil#\hfil\cr
{\rm l.i.m.}\cr
$\stackrel{}{{}_{p\to \infty}}$\cr
}} }
\sum\limits_{j_3, j_1=0}^{p}
h_{j_3 j_1}^p (s)\zeta_{j_1}^{(i_1)}
\zeta_{j_3}^{(i_3)},
$$
$$
h_{j_3 j_1}^p (s) =
\int\limits_t^s\phi_{j_1}(s_3)\int\limits_{s_3}^s\phi_{j_3}(\tau)
F_p(s_3,\tau,s)
d\tau ds_3,
$$

\vspace{2mm}
\noindent
where $F_p(s_3,\tau,s)$ is defined by (\ref{cas700}).

Analogously to (\ref{cas303a}), 
we obtain that $\Delta_8^{(i_1i_3)}(s)=0$\ \ w.~p.~1.
In this case we consider the function
$$
K_p(s_3,\tau,s)={\bf 1}_{\{s_3<s\}}{\bf 1}_{\{s_3<\tau<s\}}F_p(s_3,\tau,s)
$$
and the relations (see (\ref{cas304}))
$$
h_{j_3j_1}^{p}(s)=\int\limits_{[t,T]^2}
K_p(s_3,\tau,s)\phi_{j_1}(s_3)\phi_{j_3}(\tau)d\tau ds_3,
$$
$$
h_{j_1 j_1}^p (s)=\sum\limits_{j_4=p+1}^{\infty}
\frac{1}{2}\left(
\int\limits_t^s\phi_{j_1}(\tau)\int\limits_{\tau}^s\phi_{j_4}(s_1)
ds_1d\tau\right)^2.
$$

\vspace{2mm}

Let us prove that
\begin{equation}
\label{cas900}
\lim\limits_{p\to\infty}
\sum\limits_{j_3=0}^p c_{j_3j_3}^p (s)=0.
\end{equation}

We have 
\begin{equation}
\label{cas901}
c_{j_3j_3}^p (s)=
f_{j_3j_3}^p (s)+
d_{j_3j_3}^p (s)-g_{j_3j_3}^p (s).
\end{equation}

\vspace{2mm}

Moreover, 
\begin{equation}
\label{cas902}
\lim\limits_{p\to\infty}\sum\limits_{j_3=0}^p f_{j_3j_3}^p(s)=0,\ \ \
\lim\limits_{p\to\infty}\sum\limits_{j_3=0}^p d_{j_3j_3}^p(s)=0,
\end{equation}

\vspace{1mm}
\noindent
where the first equality in (\ref{cas902}) has been proved earlier.  
Analogously, we can prove the second equality in (\ref{cas902}).

From (\ref{cas401}) we obtain
$$
\lim\limits_{p\to\infty}\sum\limits_{j_3=0}^p g_{j_3j_3}^p (s)=0.
$$

So, (\ref{cas900}) is proved.
The relation (\ref{cas101})
is proved for the polynomial case.
Theorem 2.22 is proved for the case of Legendre polynomials.

It is easy to see that the trigonometric case is considered
by analogy with the case of Legendre polynomials
using the estimates (\ref{dwdw14}).
Theorem 2.22 is proved.

Let us reformulate Theorem 2.22 in terms on the convergence 
of the solution of system of ODEs to the solution of 
system of Stratonovich SDEs.

By analogy with (\ref{um1xxxx1}) for the case $k=4$,\ $p_1=\ldots =p_4=p,$\
$i_1, \ldots , i_4=0, 1,\ldots,m$,  and $s\in (t, T]$ ($s$ is fixed) we
obtain
$$
\int\limits_t^s
\int\limits_t^{t_4}
\int\limits_t^{t_3}
\int\limits_t^{t_2}
d{\bf w}_{t_1}^{(i_1)p}d{\bf w}_{t_2}^{(i_2)p}
d{\bf w}_{t_3}^{(i_3)p}d{\bf w}_{t_4}^{(i_4)p}=
\sum\limits_{j_1,j_2,j_3,j_4=0}^{p}
C_{j_4j_3j_2j_1}(s)
\zeta_{j_1}^{(i_1)}\zeta_{j_2}^{(i_2)}\zeta_{j_3}^{(i_3)}
\zeta_{j_4}^{(i_4)},
$$
where $p\in{\bf N}$ and $d{\bf w}_{\tau}^{(i)p}$ is defined by
(\ref{um1xxx}); another notations are the same as in Theorem 2.22.

The iterated Riemann--Stiltjes integrals 

\vspace{-1mm}
\begin{equation}
\label{cas2000}
V_{s,t}^{(i_1i_2i_3i_4)p}=
\int\limits_t^s
\int\limits_t^{t_4}
\int\limits_t^{t_3}
\int\limits_t^{t_2}
d{\bf w}_{t_1}^{(i_1)p}d{\bf w}_{t_2}^{(i_2)p}
d{\bf w}_{t_3}^{(i_3)p}d{\bf w}_{t_4}^{(i_4)p},
\end{equation}
\begin{equation}
\label{cas2001}
Z_{s,t}^{(i_1i_2i_3)p}=
\int\limits_t^s
\int\limits_t^{t_3}
\int\limits_t^{t_2}
d{\bf w}_{t_1}^{(i_1)p}d{\bf w}_{t_2}^{(i_2)p}
d{\bf w}_{t_3}^{(i_3)p},
\end{equation}
\begin{equation}
\label{cas2002}
Y_{s,t}^{(i_1i_2)p}=\int\limits_t^s
\int\limits_t^{t_2}
d{\bf w}_{t_1}^{(i_1)p}d{\bf w}_{t_2}^{(i_2)p},
\end{equation}
\begin{equation}
\label{cas2003}
X_{s,t}^{(i_1)p}=\int\limits_t^{s}
d{\bf w}_{t_1}^{(i_1)p}
\end{equation}

\vspace{4mm}
\noindent
are the solution of the following system of ODEs

\vspace{1mm}
$$
\left\{\begin{matrix}
dV_{s,t}^{(i_1i_2i_3i_4)p}=Z_{s,t}^{(i_1i_2i_3)p}
d{\bf w}_{s}^{(i_4)p},\ &V_{t,t}^{(i_1i_2i_3i_4)p}=0\cr\cr
dZ_{s,t}^{(i_1i_2i_3)p}=Y_{s,t}^{(i_1i_2)p}
d{\bf w}_{s}^{(i_3)p},\ &Z_{t,t}^{(i_1i_2i_3)p}=0\cr\cr
dY_{s,t}^{(i_1i_2)p}=X_{s,t}^{(i_1)p}
d{\bf w}_{s}^{(i_2)p},\ &Y_{t,t}^{(i_1i_2)p}=0\cr\cr
dX_{s,t}^{(i_1)p}=1 \cdot d{\bf w}_{s}^{(i_1)p},\ 
&X_{t,t}^{(i_1)p}=0
\end{matrix}\right..
$$

\vspace{5mm}

From the other hand, the iterated Stratonovich
stochastic integrals\\
\begin{equation}
\label{cas2004}
V_{s,t}^{(i_1i_2i_3i_4)}={\int\limits_t^{*}}^s
{\int\limits_t^{*}}^{t_4}
{\int\limits_t^{*}}^{t_3}
{\int\limits_t^{*}}^{t_2}d{\bf w}_{t_1}^{(i_1)}
d{\bf w}_{t_2}^{(i_2)}d{\bf w}_{t_3}^{(i_3)}d{\bf w}_{t_4}^{(i_4)},
\end{equation}
\begin{equation}
\label{cas2005}
Z_{s,t}^{(i_1i_2i_3)}={\int\limits_t^{*}}^s
{\int\limits_t^{*}}^{t_3}
{\int\limits_t^{*}}^{t_2}d{\bf w}_{t_1}^{(i_1)}
d{\bf w}_{t_2}^{(i_2)}d{\bf w}_{t_3}^{(i_3)},
\end{equation}
\begin{equation}
\label{cas2006}
Y_{s,t}^{(i_1i_2)}={\int\limits_t^{*}}^s
{\int\limits_t^{*}}^{t_2}d{\bf w}_{t_1}^{(i_1)}
d{\bf w}_{t_2}^{(i_2)},
\end{equation}
\begin{equation}
\label{cas2007}
X_{s,t}^{(i_1)}
={\int\limits_t^{*}}^{s}d{\bf w}_{t_1}^{(i_1)}
\end{equation}

\vspace{4mm}
\noindent
are the solution of the following system of Stratonovich SDEs

$$
\left\{\begin{matrix}
dV_{s,t}^{(i_1i_2i_3i_4)}=Z_{s,t}^{(i_1i_2i_3)}
* d{\bf w}_{s}^{(i_4)},\ &V_{t,t}^{(i_1i_2i_3i_4)}=0\cr\cr
dZ_{s,t}^{(i_1i_2i_3)}=Y_{s,t}^{(i_1i_2)}
* d{\bf w}_{s}^{(i_3)},\ &Z_{t,t}^{(i_1i_2i_3)}=0\cr\cr
dY_{s,t}^{(i_1i_2)}=X_{s,t}^{(i_1)}
* d{\bf w}_{s}^{(i_2)},\ &Y_{t,t}^{(i_1i_2)}=0\cr\cr
dX_{s,t}^{(i_1)}=1 * d{\bf w}_{s}^{(i_1)},\ 
&X_{t,t}^{(i_1)}=0
\end{matrix}\right.,
$$

\vspace{3mm}
\noindent
where  $*~\hspace{-0.3mm}d{\bf w}_{s}^{(i)}$, $i=0,1,\ldots,m$ is the
Stratonovich differential, $*~\hspace{-0.3mm}d{\bf w}_{s}^{(0)}=ds.$

Then from Theorems 2.19,
2.21, and 2.22 we obtain the following theorem.

{\bf Theorem 2.23}\ \cite{arxiv-5}. {\it Suppose that 
$\{\phi_j(x)\}_{j=0}^{\infty}$ is a complete orthonormal system of 
Legendre polynomials or trigonometric 
functions in the space $L_2([t, T]).$
Then
for any fixed $s$ $(s\in (t, T])$

\vspace{-1mm}
$$
\hbox{\vtop{\offinterlineskip\halign{
\hfil#\hfil\cr
{\rm l.i.m.}\cr
$\stackrel{}{{}_{p\to \infty}}$\cr
}} }V_{s,t}^{(i_1i_2i_3i_4)p}=V_{s,t}^{(i_1i_2i_3i_4)},\ \ \ 
\hbox{\vtop{\offinterlineskip\halign{
\hfil#\hfil\cr
{\rm l.i.m.}\cr
$\stackrel{}{{}_{p\to \infty}}$\cr
}} }Z_{s,t}^{(i_1i_2i_3)p}=Z_{s,t}^{(i_1i_2i_3)},
$$

\vspace{-2mm}
$$
\hbox{\vtop{\offinterlineskip\halign{
\hfil#\hfil\cr
{\rm l.i.m.}\cr
$\stackrel{}{{}_{p\to \infty}}$\cr
}} }Y_{s,t}^{(i_1i_2)p}=Y_{s,t}^{(i_1i_2)},\ \ \ 
\hbox{\vtop{\offinterlineskip\halign{
\hfil#\hfil\cr
{\rm l.i.m.}\cr
$\stackrel{}{{}_{p\to \infty}}$\cr
}} }X_{s,t}^{(i_1)p}=X_{s,t}^{(i_1)}.
$$
}

\section{Rate of the Mean-Square Convergence of Expansions of Iterated
Stra\-to\-no\-vich Stochastic Integrals of Multiplicities 2 to 4 in Theorems 2.2, 2.8, and 2.9}

\subsection{Rate of the Mean-Square Convergence of Expansion of Iterated
Stra\-to\-no\-vich Stochastic Integrals of Multiplicity 2}

This section is devoted to the proof of the following theorem.

{\bf Theorem 2.24}\ \cite{arxiv-5}. {\it Suppose that
$\{\phi_j(x)\}_{j=0}^{\infty}$ is a complete orthonormal system of 
Legendre polynomials or trigonometric functions in the space $L_2([t, T]).$
Moreover$,$ $\psi_1(\tau),$ $\psi_2(\tau)$ are 
continuously differentiable functions on $[t, T]$. Then$,$ 
for the iterated Stratonovich stochastic integral
$$
J^{*}[\psi^{(2)}]_{T,t}={\int\limits_t^{*}}^T\psi_2(t_2)
{\int\limits_t^{*}}^{t_2}\psi_1(t_1)d{\bf w}_{t_1}^{(i_1)}
d{\bf w}_{t_2}^{(i_2)}\ \ \ (i_1, i_2=1,\ldots,m)
$$
the following estimate
\begin{equation}
\label{may2000}
{\sf M}\left\{\left(J^{*}[\psi^{(2)}]_{T,t}-\sum_{j_1, j_2=0}^{p}
C_{j_2j_1}\zeta_{j_1}^{(i_1)}\zeta_{j_2}^{(i_2)}\right)^2\right\}\le \frac{C}{p}
\end{equation}
is valid, where constant $C$ is independent of $p,$
$$
C_{j_2 j_1}=\int\limits_t^T\psi_2(s_2)\phi_{j_2}(s_2)
\int\limits_t^{s_2}\psi_1(s_1)\phi_{j_1}(s_1)ds_1ds_2,
$$
and
$$
\zeta_{j}^{(i)}=
\int\limits_t^T \phi_{j}(\tau) d{\bf w}_{\tau}^{(i)}
$$ 
are independent
standard Gaussian random variables for various 
$i$ or $j$.}

{\bf Proof.}\ Applying (\ref{oop51}), we obtain
$$
{\sf M}\left\{\left(J^{*}[\psi^{(2)}]_{T,t}-\sum_{j_1, j_2=0}^{p}
C_{j_2j_1}\zeta_{j_1}^{(i_1)}\zeta_{j_2}^{(i_2)}\right)^2\right\}=
$$
$$
={\sf M}\left\{\left(
J[\psi^{(2)}]_{T,t}+
\frac{1}{2}{\bf 1}_{\{i_1=i_2\}}
\int\limits_t^T\psi_1(t_1)\psi_2(t_1)dt_1
-\sum_{j_1, j_2=0}^{p}
C_{j_2j_1}\zeta_{j_1}^{(i_1)}\zeta_{j_2}^{(i_2)}\right)^2\right\}=
$$
$$
={\sf M}\Biggl\{\Biggl(
J[\psi^{(2)}]_{T,t}-
\sum_{j_1, j_2=0}^{p}
C_{j_2j_1}\Biggl(\zeta_{j_1}^{(i_1)}\zeta_{j_2}^{(i_2)}
-{\bf 1}_{\{i_1=i_2\}}
{\bf 1}_{\{j_1=j_2\}}\Biggr)
+\Biggr.\Biggr.
$$
$$
\Biggl.\Biggl.
+\frac{1}{2}{\bf 1}_{\{i_1=i_2\}}
\int\limits_t^T\psi_1(t_1)\psi_2(t_1)dt_1- 
{\bf 1}_{\{i_1=i_2\}}\sum_{j_1=0}^{p}
C_{j_1j_1}
\Biggr)^2\Biggr\}=
$$
$$
=
{\sf M}\left\{\left(
J[\psi^{(2)}]_{T,t}-
\sum_{j_1, j_2=0}^{p}
C_{j_2j_1}\Biggl(\zeta_{j_1}^{(i_1)}\zeta_{j_2}^{(i_2)}
-{\bf 1}_{\{i_1=i_2\}}
{\bf 1}_{\{j_1=j_2\}}\Biggr)\right)^2\right\}+
$$
$$
+\left(\frac{1}{2}{\bf 1}_{\{i_1=i_2\}}
\int\limits_t^T\psi_1(t_1)\psi_2(t_1)dt_1- 
{\bf 1}_{\{i_1=i_2\}}\sum_{j_1=0}^{p}
C_{j_1j_1}
\right)^2=
$$
$$
={\sf M}\left\{\biggl(
J[\psi^{(2)}]_{T,t}-J[\psi^{(2)}]_{T,t}^{p,p}
\biggr)^2\right\}+
$$
\begin{equation}
\label{may1002}
+ {\bf 1}_{\{i_1=i_2\}}\left(\frac{1}{2}
\int\limits_t^T\psi_1(t_1)\psi_2(t_1)dt_1- 
\sum_{j_1=0}^{p}
C_{j_1j_1}
\right)^2.
\end{equation}

From Remark 1.7 (see (\ref{zsel1})) we have
\begin{equation}
\label{may1003}
{\sf M}\left\{\biggl(
J[\psi^{(2)}]_{T,t}-J[\psi^{(2)}]_{T,t}^{p,p}
\biggr)^2\right\}\le \frac{C_1}{p},
\end{equation}

\noindent
where constant $C_1$ is independent of $p.$

From Theorem 2.2 (see (\ref{may62021})) we get
\begin{equation}
\label{may1000}
\frac{1}{2}
\int\limits_t^T\psi_1(t_1)\psi_2(t_1)dt_1- 
\sum_{j_1=0}^{p}
C_{j_1j_1}=\sum_{j_1=p+1}^{\infty}
C_{j_1j_1}.
\end{equation}

Let us consider the case of Legendre polynomials.
The estimate (\ref{tupo15}) implies that
\begin{equation}
\label{may1001}
\left\vert\sum\limits_{j_1=p+1}^{\infty}C_{j_1j_1}\right\vert
\le C_2\left(\frac{1}{p}+\sum\limits_{j_1=p+1}^{\infty}
\frac{1}{j_1^2}\right),
\end{equation}

\noindent
where constant $C_2$ does not depend on $p$.

Using (\ref{obana}) and (\ref{may1001}), we have
\begin{equation}
\label{may2001}
S_{p}\stackrel{\sf def}{=}\left\vert\sum\limits_{j_1=p+1}^{\infty}C_{j_1j_1}\right\vert
\le \frac{C_3}{p}, 
\end{equation}

\noindent
where constant $C_3$ is independent of $p$.

Applying the ideas that we used to obtain the relations
(\ref{agentaa1001}), (\ref{agentaa1005})--(\ref{agentaa1007}), we can prove
the following estimates for the trigonometric case
\begin{equation}
\label{agentt10}
S_{2p}=\left\vert\sum\limits_{j_1=2p+1}^{\infty}C_{j_1j_1}\right\vert
\le \frac{K_1}{p},
\end{equation}
\begin{equation}
\label{agentt11}
S_{2p-1}=\left\vert\sum\limits_{j_1=2p}^{\infty}C_{j_1j_1}\right\vert\le S_{2p}+\frac{K_2}{p},
\end{equation}

\noindent
where constants $K_1, K_2$ do not depend on $p.$

Using (\ref{agentt10}) and (\ref{agentt11}), we get the estimate (\ref{may2001})
for the trigonometric case.
Combining (\ref{may1002})--(\ref{may1000}), (\ref{may2001}), we obtain (\ref{may2000}).
Theorem 2.24 is proved.

\subsection{Rate of the Mean-Square Convergence of Expansion of Iterated
Stra\-to\-no\-vich Stochastic Integrals of Multiplicity 3}

In this section, we consider the following theorem.

{\bf Theorem 2.25}\ \cite{arxiv-5}.
{\it Suppose that 
$\{\phi_j(x)\}_{j=0}^{\infty}$ is a complete orthonormal system of 
Legendre polynomials or trigonometric functions in the space $L_2([t, T]).$
At the same time $\psi_2(\tau)$ is a continuously dif\-ferentiable 
nonrandom function on $[t, T]$ and $\psi_1(\tau),$ $\psi_3(\tau)$ are twice
continuously differentiable nonrandom functions on $[t, T]$. 
Then$,$ for the 
iterated Stratonovich stochastic integral of third multiplicity
$$
J^{*}[\psi^{(3)}]_{T,t}={\int\limits_t^{*}}^T\psi_3(t_3)
{\int\limits_t^{*}}^{t_3}\psi_2(t_2)
{\int\limits_t^{*}}^{t_2}\psi_1(t_1)
d{\bf w}_{t_1}^{(i_1)}
d{\bf w}_{t_2}^{(i_2)}d{\bf w}_{t_3}^{(i_3)},
$$
where $i_1, i_2, i_3=1,\ldots,m,$ the following estimate
\begin{equation}
\label{may3000}
~~~~~~~ {\sf M}\left\{\left(J^{*}[\psi^{(3)}]_{T,t}-
\sum\limits_{j_1, j_2, j_3=0}^{p}
C_{j_3 j_2 j_1}\zeta_{j_1}^{(i_1)}\zeta_{j_2}^{(i_2)}\zeta_{j_3}^{(i_3)}
\right)^2\right\}\le \frac{C}{p}
\end{equation}
is valid, where constant $C$ is independent of $p,$
$$
C_{j_3 j_2 j_1}=\int\limits_t^T\psi_3(s_3)\phi_{j_3}(s_3)
\int\limits_t^{s_3}\psi_2(s_2)\phi_{j_2}(s_2)
\int\limits_t^{s_2}\psi_1(s_1)\phi_{j_1}(s_1)ds_1ds_2ds_3,
$$
and
$$
\zeta_{j}^{(i)}=
\int\limits_t^T \phi_{j}(\tau) d{\bf w}_{\tau}^{(i)}
$$ 
are independent standard Gaussian random variables for various 
$i$ or $j$.}

{\bf Proof.}\ We have (see (\ref{uyes3}))
$$
{\sf M}\left\{\left(J^{*}[\psi^{(3)}]_{T,t}-
\sum\limits_{j_1, j_2, j_3=0}^{p}
C_{j_3 j_2 j_1}\zeta_{j_1}^{(i_1)}\zeta_{j_2}^{(i_2)}\zeta_{j_3}^{(i_3)}
\right)^2\right\}=
$$
$$
={\sf M}\left\{\left(
J[\psi^{(3)}]_{T,t}+
\frac{1}{2}{\bf 1}_{\{i_1=i_2\}}
\int\limits_t^T\psi_3(t_3)
\int\limits_t^{t_3}\psi_2(t_2)\psi_1(t_2)dt_2
d{\bf w}_{t_3}^{(i_3)}+\right.\right.
$$
$$
\left.\left.+ \frac{1}{2}{\bf 1}_{\{i_2=i_3\}}
\int\limits_t^T\psi_3(t_3)\psi_2(t_3)
\int\limits_t^{t_3}\psi_1(t_1)d{\bf w}_{t_1}^{(i_1)}
dt_3-
\hspace{-0.3mm}\sum\limits_{j_1, j_2, j_3=0}^{p}
C_{j_3 j_2 j_1}\zeta_{j_1}^{(i_1)}\zeta_{j_2}^{(i_2)}\zeta_{j_3}^{(i_3)}
\right)^2\right\}\hspace{-0.3mm}=
$$
$$
={\sf M}\left\{\Biggl(
J[\psi^{(3)}]_{T,t}-J[\psi^{(3)}]_{T,t}^{p,p,p}+\Biggr.\right.
$$
$$
+
{\bf 1}_{\{i_1=i_2\}}\left(\frac{1}{2}
\int\limits_t^T\psi_3(t_3)
\int\limits_t^{t_3}\psi_2(t_2)\psi_1(t_2)dt_2
d{\bf w}_{t_3}^{(i_3)}-\sum\limits_{j_1,j_3=0}^{p}
C_{j_3 j_1 j_1}\zeta_{j_3}^{(i_3)}\right)+
$$
$$
+ {\bf 1}_{\{i_2=i_3\}}\left(\frac{1}{2}
\int\limits_t^T\psi_3(t_3)\psi_2(t_3)
\int\limits_t^{t_3}\psi_1(t_1)d{\bf w}_{t_1}^{(i_1)}
dt_3 
-\sum\limits_{j_1,j_3=0}^{p}
C_{j_3 j_3 j_1}\zeta_{j_1}^{(i_1)}\right)-
$$
\begin{equation}
\label{may4000}
\left.\Biggl.-{\bf 1}_{\{i_1=i_3\}}\sum\limits_{j_1=0}^{p}\sum\limits_{j_3=0}^{p}
C_{j_1 j_3 j_1}\zeta_{j_3}^{(i_2)}\Biggr)^2\right\},
\end{equation}

\vspace{2mm}
\noindent
where (see (\ref{a3}))
$$
J[\psi^{(3)}]_{T,t}^{p,p,p}=
\sum_{j_1,j_2,j_3=0}^{p}
C_{j_3j_2j_1}\Biggl(
\zeta_{j_1}^{(i_1)}\zeta_{j_2}^{(i_2)}\zeta_{j_3}^{(i_3)}
-\Biggr.
$$
$$
-\Biggl.
{\bf 1}_{\{i_1=i_2\}}
{\bf 1}_{\{j_1=j_2\}}
\zeta_{j_3}^{(i_3)}
-{\bf 1}_{\{i_2=i_3\}}
{\bf 1}_{\{j_2=j_3\}}
\zeta_{j_1}^{(i_1)}-
{\bf 1}_{\{i_1=i_3\}}
{\bf 1}_{\{j_1=j_3\}}
\zeta_{j_2}^{(i_2)}\Biggr).
$$

\vspace{2mm}
  
Using (\ref{may4000}) and the elementary inequality 
$$
(a+b+c+d)^2\le 4\left(a^2+b^2+c^2+d^2\right),
$$
we obtain
$$
{\sf M}\left\{\left(J^{*}[\psi^{(3)}]_{T,t}-
\sum\limits_{j_1, j_2, j_3=0}^{p}
C_{j_3 j_2 j_1}\zeta_{j_1}^{(i_1)}\zeta_{j_2}^{(i_2)}\zeta_{j_3}^{(i_3)}
\right)^2\right\}\le 
$$
$$
\le 4 \biggl({\sf M}\left\{\left(
J[\psi^{(3)}]_{s,t}-J[\psi^{(3)}]_{T,t}^{p,p,p}\right)^2\right\}
+{\bf 1}_{\{i_1=i_2\}}E_p^{(1)}+
{\bf 1}_{\{i_2=i_3\}}E_p^{(2)}+\biggr.
$$
\begin{equation}
\label{teor100aaa}
\biggl.+{\bf 1}_{\{i_1=i_3\}}E_p^{(3)}\biggr),
\end{equation}

\vspace{2mm}
\noindent
where
$$
E_p^{(1)}=
{\sf M}\left\{\left(\frac{1}{2}
\int\limits_t^T\psi_3(t_3)
\int\limits_t^{t_3}\psi_2(t_2)\psi_1(t_2)dt_2
d{\bf w}_{t_3}^{(i_3)}-\sum\limits_{j_1,j_3=0}^{p}
C_{j_3 j_1 j_1}\zeta_{j_3}^{(i_3)}\right)^2\right\},
$$
$$
E_p^{(2)}=
{\sf M}\left\{\left(\frac{1}{2}
\int\limits_t^T\psi_3(t_3)\psi_2(t_3)
\int\limits_t^{t_3}\psi_1(t_1)d{\bf w}_{t_1}^{(i_1)}
dt_3 
-\sum\limits_{j_1,j_3=0}^{p}
C_{j_3 j_3 j_1}\zeta_{j_1}^{(i_1)}\right)^2\right\},
$$
$$
E_p^{(3)}={\sf M}\left\{\left(
\sum\limits_{j_1=0}^{p}\sum\limits_{j_3=0}^{p}
C_{j_1 j_3 j_1}\zeta_{j_3}^{(i_2)}\right)^2\right\}.
$$

\vspace{3mm}

From Remark 1.7 (see (\ref{zsel1})) we have
\begin{equation}
\label{may5000}
{\sf M}\left\{\biggl(
J[\psi^{(3)}]_{T,t}-J[\psi^{(3)}]_{T,t}^{p,p,p}
\biggr)^2\right\}\le \frac{C_1}{p},
\end{equation}

\noindent
where constant $C_1$ is independent of $p.$

Moreover, from (\ref{104xx}) and (\ref{may6000}) we have the following estimate
\begin{equation}
\label{teor200}
E_p^{(3)}\le\frac{C_2}{p}
\end{equation}

\noindent
for the polynomial and trigonometric cases,
where constant $C_2$ does not depend on $p$.

Using Theorem 
1.1 for $k=1$ (also see (\ref{a1})), we obtain w. p. 1
$$
\frac{1}{2}\int\limits_t^T\psi_3(s)
\int\limits_t^s\psi_2(s_1)\psi_1(s_1)ds_1d{\bf w}_s^{(i_3)}=
\frac{1}{2}\hbox{\vtop{\offinterlineskip\halign{
\hfil#\hfil\cr
{\rm l.i.m.}\cr
$\stackrel{}{{}_{p\to \infty}}$\cr
}} }
\sum\limits_{j_3=0}^{p}
\tilde C_{j_3}\zeta_{j_3}^{(i_3)},
$$

\noindent
where 
$$
\tilde C_{j_3}=
\int\limits_t^T
\phi_{j_3}(s)\psi_3(s)\int\limits_t^s\psi_2(s_1)\psi_1(s_1)ds_1ds.
$$

Applying the It\^{o} formula, we have
$$
\int\limits_t^T\psi_3(s)\psi_2(s)
\int\limits_t^s\psi_1(s_1)d{\bf w}_{s_1}^{(i_1)}ds=
\int\limits_t^T\psi_1(s_1)
\int\limits_{s_1}^T\psi_3(s)\psi_2(s)dsd{\bf w}_{s_1}^{(i_1)}\ \ \
\hbox{\rm w.~p.~1}.
$$

Using Theorem 
1.1 for $k=1$ (also see (\ref{a1})), we have w. p. 1
$$
\frac{1}{2}\int\limits_t^T\psi_1(s)
\int\limits_{s}^T\psi_3(s_1)\psi_2(s_1)ds_1d{\bf w}_s^{(i_1)}=
\frac{1}{2}\hbox{\vtop{\offinterlineskip\halign{
\hfil#\hfil\cr
{\rm l.i.m.}\cr
$\stackrel{}{{}_{p\to \infty}}$\cr
}} }
\sum\limits_{j_1=0}^{p}
C_{j_1}^{*}\zeta_{j_1}^{(i_1)},
$$

\noindent
where 
$$
C_{j_1}^{*}=
\int\limits_t^T
\psi_1(s)\phi_{j_1}(s)\int\limits_{s}^T\psi_3(s_1)\psi_2(s_1)ds_1ds.
$$

Further, we get
\begin{equation}
\label{teor300}
E_p^{(1)}\le 2 G_p^{(1)} + 2 G_p^{(2)},
\end{equation}
\begin{equation}
\label{teor400}
E_p^{(2)}\le 2 H_p^{(1)} + 2 H_p^{(2)},
\end{equation}
where
$$
G_p^{(1)}={\sf M}\left\{\frac{1}{4}\left(\int\limits_t^T\psi_3(t_3)
\int\limits_t^{t_3}\psi_2(t_2)\psi_1(t_2)dt_2
d{\bf w}_{t_3}^{(i_3)}-\sum\limits_{j_3=0}^{p}
\tilde C_{j_3}\zeta_{j_3}^{(i_3)}\right)^2\right\},
$$
$$
G_p^{(2)}={\sf M}\left\{\left(\frac{1}{2}\sum\limits_{j_3=0}^{p}
\tilde C_{j_3}\zeta_{j_3}^{(i_3)}-
\sum\limits_{j_1,j_3=0}^{p}
C_{j_3 j_1 j_1}\zeta_{j_3}^{(i_3)}\right)^2\right\},
$$
$$
H_p^{(1)}=
{\sf M}\left\{\frac{1}{4}\left(\int\limits_t^T\psi_3(t_3)\psi_2(t_3)
\int\limits_t^{t_3}\psi_1(t_1)d{\bf w}_{t_1}^{(i_1)}
dt_3 
-\sum\limits_{j_1=0}^{p}
C_{j_1}^{*}\zeta_{j_1}^{(i_1)}\right)^2\right\},
$$
$$
H_p^{(2)}={\sf M}\left\{\left(\frac{1}{2}\sum\limits_{j_1=0}^{p}
C_{j_1}^{*}\zeta_{j_1}^{(i_1)}-\sum\limits_{j_1,j_3=0}^{p}
C_{j_3 j_3 j_1}\zeta_{j_1}^{(i_1)}\right)^2\right\}.
$$

\vspace{3mm}

From Remark 1.7 (see (\ref{zsel1})) we have
\begin{equation}
\label{may7000}
G_p^{(1)}\le \frac{C_2}{p},\ \ \  H_p^{(1)}\le \frac{C_2}{p},
\end{equation}

\noindent
where constant $C_2$ is independent of $p.$

The estimates
\begin{equation}
\label{may8000}
G_p^{(2)}\le \frac{C_3}{p},\ \ \  H_p^{(2)}\le \frac{C_3}{p}
\end{equation}

\noindent
are proved in Sect.~2.2.5 (see the proof of Theorem 2.8)
for the polynomial and trigonometric cases; constant $C_3$ does not depend on $p.$

Combining the estimates (\ref{teor100aaa})--(\ref{may8000}), we obtain 
the inequality (\ref{may3000}). Theorem 2.25 is proved.

\subsection{Rate of the Mean-Square Convergence of Expansion of Iterated
Stra\-to\-no\-vich Stochastic Integrals of Multiplicity 4}

This section is devoted to the proof of the following theorem.

{\bf Theorem 2.26} \hspace{-1mm} \cite{arxiv-5}. {\it \hspace{-3mm} Suppose that
$\{\phi_j(x)\}_{j=0}^{\infty}$ is a complete orthonormal
system of Legendre polynomials or trigonometric functions
in the space $L_2([t, T])$.
Then$,$ for the iterated 
Stratonovich stochastic integral of fourth multiplicity
$$
J^{*}[\psi^{(4)}]_{T,t}=
{\int\limits_t^{*}}^T
{\int\limits_t^{*}}^{t_4}
{\int\limits_t^{*}}^{t_3}
{\int\limits_t^{*}}^{t_2}
d{\bf w}_{t_1}^{(i_1)}
d{\bf w}_{t_2}^{(i_2)}d{\bf w}_{t_3}^{(i_3)}d{\bf w}_{t_4}^{(i_4)}\ \ \ 
(i_1, i_2, i_3, i_4=0, 1,\ldots,m)
$$

\noindent
the following estimate
\begin{equation}
\label{may9000}
~~~~~~~~{\sf M}\left\{\left(J^{*}[\psi^{(4)}]_{T,t}-
\sum\limits_{j_1, j_2, j_3,j_4=0}^{p}
C_{j_4j_3 j_2 j_1}\zeta_{j_1}^{(i_1)}\zeta_{j_2}^{(i_2)}\zeta_{j_3}^{(i_3)}\zeta_{j_4}^{(i_4)}
\right)^2\right\}\le \frac{C}{p}
\end{equation}

\noindent
is valid, where constant $C$ is independent of $p,$
$$
C_{j_4 j_3 j_2 j_1}=\int\limits_t^T\phi_{j_4}(s_4)\int\limits_t^{s_4}
\phi_{j_3}(s_3)
\int\limits_t^{s_3}\phi_{j_2}(s_2)\int\limits_t^{s_2}\phi_{j_1}(s_1)
ds_1ds_2ds_3ds_4,
$$
and
$$
\zeta_{j}^{(i)}=
\int\limits_t^T \phi_{j}(\tau) d{\bf w}_{\tau}^{(i)}
$$ 
are independent standard Gaussian random variables for various 
$i$ or $j$ {\rm (}in the case when $i\ne 0${\rm ),}
${\bf w}_{\tau}^{(i)}$
$(i=1,\ldots,m)$ are independent standard Wiener processes$,$ 
${\bf w}_{\tau}^{(0)}=\tau.$}

{\bf Proof.}\ First, let us prove that Theorem 2.8 is valid for the case 
$i_1,i_2,i_3=0,1,\ldots,m.$ The case $i_1,i_2,i_3=1,\ldots,m$
has been proved in Theorem 2.8. 
From (\ref{a3}) and
the standard relation (\ref{uyes3}) between Stratonovich and It\^{o} stochastic
integrals of third multiplicity 
we have that Theorem 2.8 is valid for the following cases
$$
i_1=i_2=0,\ \ \ i_3=1,\ldots,m,
$$
$$
i_1=i_3=0,\ \ \ i_2=1,\ldots,m,
$$
$$
i_2=i_3=0,\ \ \ i_1=1,\ldots,m.
$$

Thus, it remains to consider the following three cases
\begin{equation}
\label{agentqq1}
i_1,\ i_2=1,\ldots,m,\ \ \ i_3=0,
\end{equation}
\begin{equation}
\label{agentqq2}
i_2,\ i_3=1,\ldots,m,\ \ \ i_1=0,
\end{equation}

\vspace{-5mm}
\begin{equation}
\label{agentqq3}
i_1,\ i_3=1,\ldots,m,\ \ \ i_2=0.
\end{equation}

\vspace{1mm}

The relations (\ref{a3}) and (\ref{uyes3}) imply that for the case
(\ref{agentqq1}) we need to prove the following equality
$$
\sum\limits_{j_1=0}^{\infty}
\int\limits_t^T
\psi_3(t_3)\int\limits_t^{t_3}\phi_{j_1}(t_2)\psi_2(t_2)
\int\limits_t^{t_2}\phi_{j_1}(t_1)\psi_1(t_1)dt_1 dt_2 dt_3=
$$
\begin{equation}
\label{agentqq4}
=
\frac{1}{2}\int\limits_t^T
\psi_3(t_3)\int\limits_t^{t_3}
\psi_1(t_1)\psi_2(t_1)dt_1dt_3.
\end{equation}

\vspace{1mm}

Using the relation (\ref{5t}), we get
$$
\sum\limits_{j_1=0}^{\infty}
\int\limits_t^T
\psi_3(t_3)\int\limits_t^{t_3}\phi_{j_1}(t_2)\psi_2(t_2)
\int\limits_t^{t_2}\phi_{j_1}(t_1)\psi_1(t_1)dt_1 dt_2 dt_3=
$$
$$
=\sum\limits_{j_1=0}^{\infty}
\int\limits_t^T
\phi_{j_1}(t_1)\psi_1(t_1)
\int\limits_{t_1}^T
\phi_{j_1}(t_2)\psi_2(t_2)
\int\limits_{t_2}^T \psi_3(t_3)dt_3 dt_2 dt_1=
$$
$$
=\sum\limits_{j_1=0}^{\infty}
\int\limits_t^T
\phi_{j_1}(t_1)\psi_1(t_1)
\int\limits_{t_1}^T
\phi_{j_1}(t_2)\tilde \psi_2(t_2)dt_2 dt_1=
$$
$$
=\sum\limits_{j_1=0}^{\infty}
\int\limits_t^T
\phi_{j_1}(t_2)\tilde \psi_2(t_2)
\int\limits_t^{t_2}
\phi_{j_1}(t_1)\psi_1(t_1)
dt_1 dt_2=
$$
\begin{equation}
\label{agentqq5}
=\frac{1}{2}\int\limits_t^T
\psi_1(t_2)\tilde \psi_2(t_2)dt_2,
\end{equation}

\noindent
where
\begin{equation}
\label{agent00x}
\tilde \psi_2(t_2)=
\psi_2(t_2)
\int\limits_{t_2}^T \psi_3(t_3)dt_3.
\end{equation}

From (\ref{agentqq5}) and (\ref{agent00x}) we obtain
$$
\sum\limits_{j_1=0}^{\infty}
\int\limits_t^T
\psi_3(t_3)\int\limits_t^{t_3}\phi_{j_1}(t_2)\psi_2(t_2)
\int\limits_t^{t_2}\phi_{j_1}(t_1)\psi_1(t_1)dt_1 dt_2 dt_3=
$$
$$
=\frac{1}{2}\int\limits_t^T
\psi_1(t_2)
\psi_2(t_2)\int\limits_{t_2}^T \psi_3(t_3)dt_3 dt_2=
$$
\begin{equation}
\label{agentqq6}
=\frac{1}{2}\int\limits_t^T
\psi_3(t_3)\int\limits_t^{t_3}
\psi_1(t_2)\psi_2(t_2)dt_2dt_3.
\end{equation}

The relation (\ref{agentqq4}) is proved.

From (\ref{a3}) and (\ref{uyes3}) it follows that for the case
(\ref{agentqq2}) we need to prove the following equality
$$
\sum\limits_{j_2=0}^{\infty}
\int\limits_t^T
\phi_{j_2}(t_3)\psi_3(t_3)\int\limits_t^{t_3}\phi_{j_2}(t_2)\psi_2(t_2)
\int\limits_t^{t_2}\psi_1(t_1)dt_1 dt_2 dt_3=
$$
\begin{equation}
\label{agentqq7}
=
\frac{1}{2}\int\limits_t^T
\psi_3(t_3)\psi_2(t_3)\int\limits_t^{t_3}
\psi_1(t_1)dt_1dt_3.
\end{equation}

\vspace{1mm}

Using the relation (\ref{5t}), we have
$$
\sum\limits_{j_2=0}^{\infty}
\int\limits_t^T
\phi_{j_2}(t_3)\psi_3(t_3)\int\limits_t^{t_3}\phi_{j_2}(t_2)\psi_2(t_2)
\int\limits_t^{t_2}\psi_1(t_1)dt_1 dt_2 dt_3=
$$
$$
=\sum\limits_{j_2=0}^{\infty}
\int\limits_t^T
\phi_{j_2}(t_3)\psi_3(t_3)\int\limits_t^{t_3}\phi_{j_2}(t_2)\bar \psi_2(t_2)
dt_2 dt_3=
$$
$$
=\frac{1}{2}\int\limits_t^T
\psi_3(t_3)\bar \psi_2(t_3)
dt_3,
$$

\noindent
where
\begin{equation}
\label{ag010}
\bar \psi_2(t_2)=
\psi_2(t_2)
\int\limits_t^{t_2}\psi_1(t_1)dt_1.
\end{equation}

The relation (\ref{agentqq7}) is proved.

The relations (\ref{a3}) and (\ref{uyes3}) imply that for the case
(\ref{agentqq3}) we need to prove the following equality
\begin{equation}
\label{agentqq8}
~~~~~~~~~ \sum\limits_{j_1=0}^{\infty}
\int\limits_t^T
\phi_{j_1}(t_3)\psi_3(t_3)\int\limits_t^{t_3}\psi_2(t_2)
\int\limits_t^{t_2}\phi_{j_1}(t_1)\psi_1(t_1)dt_1 dt_2 dt_3=0.
\end{equation}

\vspace{1mm}

We have
$$
\sum\limits_{j_1=0}^{\infty}
\int\limits_t^T
\phi_{j_1}(t_3)\psi_3(t_3)\int\limits_t^{t_3}\psi_2(t_2)
\int\limits_t^{t_2}\phi_{j_1}(t_1)\psi_1(t_1)dt_1 dt_2 dt_3=
$$
$$
=\sum\limits_{j_1=0}^{\infty}
\int\limits_t^T
\phi_{j_1}(t_3)\psi_3(t_3)\int\limits_t^{t_3}
\phi_{j_1}(t_1)\psi_1(t_1)
\int\limits_{t_1}^{t_3}
\psi_2(t_2)dt_2
dt_1dt_3=
$$
$$
=\sum\limits_{j_1=0}^{\infty}
\int\limits_t^T
\phi_{j_1}(t_3)\psi_3(t_3)\int\limits_t^{t_3}
\phi_{j_1}(t_1)\psi_1(t_1)
\left(\int\limits_{t_1}^{T}
\psi_2(t_2)dt_2-\int\limits_{t_3}^{T}
\psi_2(t_2)dt_2\right)
dt_1dt_3=
$$
$$
=\sum\limits_{j_1=0}^{\infty}
\int\limits_t^T
\phi_{j_1}(t_3)\psi_3(t_3)\int\limits_t^{t_3}
\phi_{j_1}(t_1)\psi_1(t_1)
\int\limits_{t_1}^{T}
\psi_2(t_2)dt_2
dt_1dt_3-
$$
$$
-\sum\limits_{j_1=0}^{\infty}
\int\limits_t^T
\phi_{j_1}(t_3)\psi_3(t_3)\int\limits_t^{t_3}
\phi_{j_1}(t_1)\psi_1(t_1)
\int\limits_{t_3}^{T}
\psi_2(t_2)dt_2
dt_1dt_3=
$$
$$
=\sum\limits_{j_1=0}^{\infty}
\int\limits_t^T
\phi_{j_1}(t_3)\psi_3(t_3)\int\limits_t^{t_3}
\phi_{j_1}(t_1)\tilde \psi_1(t_1)
dt_1dt_3-
$$
$$
-\sum\limits_{j_1=0}^{\infty}
\int\limits_t^T
\phi_{j_1}(t_3)\tilde \psi_3(t_3)\int\limits_t^{t_3}
\phi_{j_1}(t_1)\psi_1(t_1)
dt_1dt_3=
$$
$$
=\frac{1}{2}\int\limits_t^T
\psi_3(t_1)\tilde \psi_1(t_1)dt_1-
\frac{1}{2}\int\limits_t^T
\tilde \psi_3(t_1)\psi_1(t_1)dt_1=
$$
$$
=\frac{1}{2}\int\limits_t^T
\psi_3(t_1)\psi_1(t_1)\int\limits_{t_1}^{T}
\psi_2(t_2)dt_2dt_1-
\frac{1}{2}\int\limits_t^T
\psi_1(t_1)
\psi_3(t_1)
\int\limits_{t_1}^{T}
\psi_2(t_2)dt_2dt_1=0,
$$

\noindent
where
\begin{equation}
\label{ag200}
\tilde \psi_1(t_1)=
\psi_1(t_1)\int\limits_{t_1}^{T}
\psi_2(t_2)dt_2,
\end{equation}
\begin{equation}
\label{ag201}
\tilde \psi_3(t_3)=\psi_3(t_3)
\int\limits_{t_3}^{T}
\psi_2(t_2)dt_2.
\end{equation}

The relation (\ref{agentqq8}) is proved.
Theorem 2.8 is proved for the case 
$i_1,i_2,i_3=0,1,\ldots,m.$

Using
(\ref{uyes3}) and (\ref{uyes3may}), we obtain
$$
{\sf M}\left\{\left(J^{*}[\psi^{(4)}]_{T,t}-
\sum\limits_{j_1, j_2, j_3,j_4=0}^{p}
C_{j_4j_3 j_2 j_1}\zeta_{j_1}^{(i_1)}\zeta_{j_2}^{(i_2)}\zeta_{j_3}^{(i_3)}\zeta_{j_4}^{(i_4)}
\right)^2\right\}=
$$
$$
=
{\sf M}\left\{\left(
J[\psi^{(4)}]_{T,t}+
\frac{1}{2}{\bf 1}_{\{i_1=i_2\ne 0\}}
\int\limits_t^T\int\limits_t^s\int\limits_t^{s_1}ds_2
d{\bf w}_{s_1}^{(i_3)}
d{\bf w}_{s}^{(i_4)}+\right.\right.
$$
$$
+\frac{1}{2}{\bf 1}_{\{i_2=i_3\ne 0\}}
\int\limits_t^T\int\limits_t^{s_2}\int\limits_t^{s_1}
d{\bf w}_{s}^{(i_1)}ds_1
d{\bf w}_{s_2}^{(i_4)}
+\frac{1}{2}{\bf 1}_{\{i_3=i_4\ne 0\}}
\int\limits_t^T\int\limits_t^{s_1}\int\limits_t^{s_2}
d{\bf w}_{s}^{(i_1)}
d{\bf w}_{s_2}^{(i_2)}ds_1+
$$
$$
\left.\left.+\frac{1}{4}{\bf 1}_{\{i_1=i_2\ne 0\}}
{\bf 1}_{\{i_3=i_4\ne 0\}}
\int\limits_t^T\int\limits_t^{s_1}ds_2
ds_1-\sum\limits_{j_1, j_2, j_3,j_4=0}^{p}
C_{j_4j_3 j_2 j_1}\zeta_{j_1}^{(i_1)}\zeta_{j_2}^{(i_2)}\zeta_{j_3}^{(i_3)}\zeta_{j_4}^{(i_4)}
\right)^2\right\}=
$$
$$
=
{\sf M}\left\{\left(
J[\psi^{(4)}]_{T,t}+
\frac{1}{2}{\bf 1}_{\{i_1=i_2\ne 0\}}
{\int\limits_t^{*}}^T
{\int\limits_t^{*}}^s
{\int\limits_t^{*}}^{s_1}
ds_2
d{\bf w}_{s_1}^{(i_3)}
d{\bf w}_{s}^{(i_4)}-\right.\right.
$$
$$
-\frac{1}{4}{\bf 1}_{\{i_1=i_2\ne 0\}}{\bf 1}_{\{i_3=i_4\ne 0\}}
\int\limits_t^{T}\int\limits_t^{s_1}ds_2
ds_1 
+\frac{1}{2}{\bf 1}_{\{i_2=i_3\ne 0\}}
{\int\limits_t^{*}}^T
{\int\limits_t^{*}}^{s_2}
{\int\limits_t^{*}}^{s_1}
d{\bf w}_{s}^{(i_1)}ds_1
d{\bf w}_{s_2}^{(i_4)}
+
$$
$$
+\frac{1}{2}{\bf 1}_{\{i_3=i_4\ne 0\}}
{\int\limits_t^{*}}^T
{\int\limits_t^{*}}^{s_1}
{\int\limits_t^{*}}^{s_2}
d{\bf w}_{s}^{(i_1)}
d{\bf w}_{s_2}^{(i_2)}ds_1-\frac{1}{4}{\bf 1}_{\{i_1=i_2\ne 0\}}{\bf 1}_{\{i_3=i_4\ne 0\}}
\int\limits_t^{T}\int\limits_t^{s_1}ds_2ds_1+
$$
$$
\left.\left.+\frac{1}{4}{\bf 1}_{\{i_1=i_2\ne 0\}}
{\bf 1}_{\{i_3=i_4\ne 0\}}
\int\limits_t^T\int\limits_t^{s_1}ds_2
ds_1-\sum\limits_{j_1, j_2, j_3,j_4=0}^{p}
C_{j_4j_3 j_2 j_1}\zeta_{j_1}^{(i_1)}\zeta_{j_2}^{(i_2)}\zeta_{j_3}^{(i_3)}\zeta_{j_4}^{(i_4)}
\right)^2\right\}=
$$

$$
=
{\sf M}\left\{\left(
J[\psi^{(4)}]_{T,t}+
\frac{1}{2}{\bf 1}_{\{i_1=i_2\ne 0\}}
{\int\limits_t^{*}}^T
{\int\limits_t^{*}}^s
{\int\limits_t^{*}}^{s_1}
ds_2
d{\bf w}_{s_1}^{(i_3)}
d{\bf w}_{s}^{(i_4)}+\right.\right.
$$
$$
+\frac{1}{2}{\bf 1}_{\{i_2=i_3\ne 0\}}
{\int\limits_t^{*}}^T
{\int\limits_t^{*}}^{s_2}
{\int\limits_t^{*}}^{s_1}
d{\bf w}_{s}^{(i_1)}ds_1
d{\bf w}_{s_2}^{(i_4)}
+\frac{1}{2}{\bf 1}_{\{i_3=i_4\ne 0\}}
{\int\limits_t^{*}}^T
{\int\limits_t^{*}}^{s_1}
{\int\limits_t^{*}}^{s_2}
d{\bf w}_{s}^{(i_1)}
d{\bf w}_{s_2}^{(i_2)}ds_1-
$$
$$
\left.\left.-\frac{1}{4}{\bf 1}_{\{i_1=i_2\ne 0\}}
{\bf 1}_{\{i_3=i_4\ne 0\}}
\int\limits_t^T\int\limits_t^{s_1}ds_2
ds_1-\sum\limits_{j_1, j_2, j_3,j_4=0}^{p}
C_{j_4j_3 j_2 j_1}\zeta_{j_1}^{(i_1)}\zeta_{j_2}^{(i_2)}\zeta_{j_3}^{(i_3)}\zeta_{j_4}^{(i_4)}
\right)^2\right\}=
$$

$$
=
{\sf M}\Biggl\{\Biggl(
J[\psi^{(4)}]_{T,t}-J[\psi^{(4)}]_{T,t}^{p,p,p,p}+\Biggr.\Biggr.
$$
$$
+\frac{1}{2}{\bf 1}_{\{i_1=i_2\ne 0\}}\left(
{\int\limits_t^{*}}^T
{\int\limits_t^{*}}^s
{\int\limits_t^{*}}^{s_1}
ds_2
d{\bf w}_{s_1}^{(i_3)}
d{\bf w}_{s}^{(i_4)}-S_1^{(i_3i_4)p}\right) +
$$
$$
+\frac{1}{2}{\bf 1}_{\{i_2=i_3\ne 0\}}\left(
{\int\limits_t^{*}}^T
{\int\limits_t^{*}}^{s_2}
{\int\limits_t^{*}}^{s_1}
d{\bf w}_{s}^{(i_1)}ds_1
d{\bf w}_{s_2}^{(i_4)}-S_2^{(i_1i_4)p}\right) +
$$
$$
+\frac{1}{2}{\bf 1}_{\{i_3=i_4\ne 0\}}\left(
{\int\limits_t^{*}}^T
{\int\limits_t^{*}}^{s_1}
{\int\limits_t^{*}}^{s_2}
d{\bf w}_{s}^{(i_1)}
d{\bf w}_{s_2}^{(i_2)}ds_1-S_3^{(i_1i_2)p}\right) -
$$
$$
-{\bf 1}_{\{i_1=i_2\ne 0\}}{\bf 1}_{\{i_3=i_4\ne 0\}}\left(\frac{1}{4}
\int\limits_t^T\int\limits_t^{s_1}ds_2
ds_1-\right.
$$
\begin{equation}
\label{may10000}
~~~~~~~~\left.\left.\left.-\sum\limits_{j_4=0}^{p}
\frac{1}{2}\int\limits_t^T\phi_{j_4}(s)\int\limits_{t}^s\phi_{j_4}(s_1)
(s_1-t)ds_1ds\right)-R_p\right)^2\right\},
\end{equation}

\vspace{4mm}
\noindent
where $S_1^{(i_3i_4)p},$ $S_2^{(i_1i_4)p},$ $S_3^{(i_1i_2)p}$
are the approximations of the iterated Stratonovich stochastic integrals
$$
{\int\limits_t^{*}}^T
{\int\limits_t^{*}}^s
{\int\limits_t^{*}}^{s_1}
ds_2
d{\bf w}_{s_1}^{(i_3)}
d{\bf w}_{s}^{(i_4)},\ \ \ 
{\int\limits_t^{*}}^T
{\int\limits_t^{*}}^{s_2}
{\int\limits_t^{*}}^{s_1}
d{\bf w}_{s}^{(i_1)}ds_1
d{\bf w}_{s_2}^{(i_4)},
$$
$$
{\int\limits_t^{*}}^T
{\int\limits_t^{*}}^{s_1}
{\int\limits_t^{*}}^{s_2}
d{\bf w}_{s}^{(i_1)}
d{\bf w}_{s_2}^{(i_2)}ds_1,
$$

\vspace{1mm}
\noindent
respectively (these approximations are obtained by
the version of Theorem 2.8 for the case 
$i_1,i_2,i_3=0,1,\ldots,m$); $J[\psi^{(4)}]_{T,t}^{p,p,p,p}$
is the approximation of the iterated It\^{o} stochastic integral
$J[\psi^{(4)}]_{T,t}$ obtained by Theorem 1.1 (see (\ref{a4}))
$$
J[\psi^{(4)}]_{T,t}^{p,p,p,p}=
\sum\limits_{j_1, j_2, j_3, j_4=0}^{p}
C_{j_4 j_3 j_2 j_1}\zeta_{j_1}^{(i_1)}\zeta_{j_2}^{(i_2)}\zeta_{j_3}^{(i_3)}
\zeta_{j_4}^{(i_4)}-
$$
$$
-{\bf 1}_{\{i_1=i_2\ne 0\}}A_1^{(i_3i_4)p}
-{\bf 1}_{\{i_1=i_3\ne 0\}}A_2^{(i_2i_4)p}
-{\bf 1}_{\{i_1=i_4\ne 0\}}A_3^{(i_2i_3)p}
-{\bf 1}_{\{i_2=i_3\ne 0\}}A_4^{(i_1i_4)p}-
$$

\vspace{-2mm}
$$
-
{\bf 1}_{\{i_2=i_4\ne 0\}}A_5^{(i_1i_3)p}
-{\bf 1}_{\{i_3=i_4\ne 0\}}A_6^{(i_1i_2)p}+
{\bf 1}_{\{i_1=i_2\ne 0\}}
{\bf 1}_{\{i_3=i_4\ne 0\}}B_1^p+
$$

\vspace{-1mm}
$$
+{\bf 1}_{\{i_1=i_3\ne 0\}}
{\bf 1}_{\{i_2=i_4\ne 0\}}B_2^p+
{\bf 1}_{\{i_1=i_4\ne 0\}}
{\bf 1}_{\{i_2=i_3\ne 0\}}B_3^p,
$$

\vspace{4mm}
\noindent
where
$$
A_1^{(i_3i_4)p}=
\sum\limits_{j_4, j_3, j_1=0}^{p}
C_{j_4 j_3 j_1 j_1}\zeta_{j_3}^{(i_3)}
\zeta_{j_4}^{(i_4)},\ \ \ 
A_2^{(i_2i_4)p}=
\sum\limits_{j_4, j_3, j_2=0}^{p}
C_{j_4 j_3 j_2 j_3}\zeta_{j_2}^{(i_2)}
\zeta_{j_4}^{(i_4)},
$$
$$
A_3^{(i_2i_3)p}=
\sum\limits_{j_4, j_3, j_2=0}^{p}
C_{j_4 j_3 j_2 j_4}\zeta_{j_2}^{(i_2)}
\zeta_{j_3}^{(i_3)},\ \ \ 
A_4^{(i_1i_4)p}=
\sum\limits_{j_4, j_3, j_1=0}^{p}
C_{j_4 j_3 j_3 j_1}\zeta_{j_1}^{(i_1)}
\zeta_{j_4}^{(i_4)},
$$
$$
A_5^{(i_1i_3)p}=
\sum\limits_{j_4, j_3, j_1=0}^{p}
C_{j_4 j_3 j_4 j_1}\zeta_{j_1}^{(i_1)}
\zeta_{j_3}^{(i_3)},\ \ \ A_6^{(i_1i_2)p}=
\sum\limits_{j_3, j_2, j_1=0}^{p}
C_{j_3 j_3 j_2 j_1}\zeta_{j_1}^{(i_1)}
\zeta_{j_2}^{(i_2)},
$$
$$
B_1^p=
\sum\limits_{j_1, j_4=0}^{p}
C_{j_4 j_4 j_1 j_1},\ \ \
B_2^p=
\sum\limits_{j_4, j_3=0}^{p}
C_{j_3 j_4 j_3 j_4},
$$
$$
B_3^p=
\sum\limits_{j_4, j_3=0}^{p}
C_{j_4 j_3 j_3 j_4};
$$

\vspace{6mm}
\noindent
$R_p$ is the expression on the right-hand side of (\ref{otiteee0}) before passing to the limits, i.e.

\vspace{-5mm}
$$
R_p=
-{\bf 1}_{\{i_1=i_2\ne 0\}}\Delta_1^{(i_3i_4)p}
+{\bf 1}_{\{i_1=i_3\ne 0\}}\left(
-\Delta_2^{(i_2i_4)p}
+\Delta_1^{(i_2i_4)p}
+\Delta_3^{(i_2i_4)p}\right)+
$$

\vspace{-3mm}
$$
+{\bf 1}_{\{i_1=i_4\ne 0\}}\left(
\Delta_4^{(i_2i_3)p}-
\Delta_5^{(i_2i_3)p}
+\Delta_6^{(i_2i_3)p}\right)-
{\bf 1}_{\{i_2=i_3\ne 0\}}\Delta_3^{(i_1i_4)p}+
$$

\vspace{-3mm}
$$
+{\bf 1}_{\{i_2=i_4\ne 0\}}
\left(-\Delta_4^{(i_1i_3)p}
+\Delta_5^{(i_1i_3)p}
+\Delta_6^{(i_1i_3)p}\right)-
{\bf 1}_{\{i_3=i_4\ne 0\}}\Delta_6^{(i_1i_2)p}-
$$
$$
-
{\bf 1}_{\{i_1=i_3\ne 0\}}
{\bf 1}_{\{i_2=i_4\ne 0\}}\Biggl(
\sum\limits_{j_3=0}^p a_{j_3j_3}^p
+\sum\limits_{j_3=0}^p c_{j_3j_3}^p
-\sum\limits_{j_3=0}^p b_{j_3j_3}^p\Biggr)-
$$

\vspace{-4mm}
$$
-{\bf 1}_{\{i_1=i_4\ne 0\}}
{\bf 1}_{\{i_2=i_3\ne 0\}}
\Biggl(2\sum\limits_{j_3=0}^p f_{j_3j_3}^p
-\sum\limits_{j_3=0}^p a_{j_3j_3}^p
-\sum\limits_{j_3=0}^p c_{j_3j_3}^p
+\sum\limits_{j_3=0}^p b_{j_3j_3}^p\Biggr)+
$$

\vspace{-3mm}
$$
+{\bf 1}_{\{i_1=i_2\ne 0\}}
{\bf 1}_{\{i_3=i_4\ne 0\}}\sum\limits_{j_3=0}^p a_{j_3j_3}^p,
$$

\noindent
where

\vspace{-4mm}
$$
\Delta_1^{(i_3i_4)p}=
\sum\limits_{j_3, j_4=0}^{p}
a_{j_4 j_3}^p \zeta_{j_3}^{(i_3)}
\zeta_{j_4}^{(i_4)},\ \ \
\Delta_2^{(i_2i_4)p}=
\sum\limits_{j_4, j_2=0}^{p}
b_{j_4 j_2}^p \zeta_{j_2}^{(i_2)}
\zeta_{j_4}^{(i_4)},
$$

\vspace{-2mm}
$$
\Delta_3^{(i_2i_4)p}=
\sum\limits_{j_4, j_2=0}^{p}
c_{j_4 j_2}^p \zeta_{j_2}^{(i_2)}
\zeta_{j_4}^{(i_4)},\ \ \ \Delta_4^{(i_1i_3)p}=
\sum\limits_{j_3, j_1=0}^{p}
d_{j_3 j_1}^p \zeta_{j_1}^{(i_1)}
\zeta_{j_3}^{(i_3)},
$$

\vspace{-2mm}
$$
\Delta_5^{(i_1i_3)p}=
\sum\limits_{j_3, j_1=0}^{p}
e_{j_3 j_1}^p \zeta_{j_1}^{(i_1)}
\zeta_{j_3}^{(i_3)},\ \ \ \Delta_6^{(i_1i_3)p}=
\sum\limits_{j_3, j_1=0}^{p}
f_{j_3 j_1}^p \zeta_{j_1}^{(i_1)}
\zeta_{j_3}^{(i_3)},
$$

\vspace{3mm}
\noindent
where 

\vspace{-6mm}
$$
a_{j_4 j_3}^p,\ \ b_{j_4 j_2}^p,\ \ 
c_{j_4 j_2}^p,\ \  d_{j_3 j_1}^p,\ \ 
e_{j_3 j_1}^p,\ \  f_{j_3 j_1}^p
$$

\vspace{5mm}
\noindent
are defined by the relations
(\ref{rr1xx}), (\ref{may001}), (\ref{may002}), (\ref{may003})--(\ref{may005}).

Using (\ref{may10000}) and the elementary inequality 

\vspace{-2mm}
$$
(a_1+\ldots+a_6)^2\le 6\left(a_1^2+\ldots+a_6^2\right),
$$ 

\vspace{1mm}
\noindent
we get
$$
{\sf M}\left\{\left(J^{*}[\psi^{(4)}]_{T,t}-
\sum\limits_{j_1, j_2, j_3,j_4=0}^{p}
C_{j_4j_3 j_2 j_1}\zeta_{j_1}^{(i_1)}\zeta_{j_2}^{(i_2)}\zeta_{j_3}^{(i_3)}\zeta_{j_4}^{(i_4)}
\right)^2\right\}\le
$$

\begin{equation}
\label{mayx200}
\le
6\left(Q_p^{(1)}+Q_p^{(2)}+Q_p^{(3)}+Q_p^{(4)}+Q_p^{(5)}+Q_p^{(6)}\right),
\end{equation}

\vspace{6mm}
\noindent
where
$$
Q_p^{(1)}={\sf M}\left\{\biggl(
J[\psi^{(4)}]_{T,t}-J[\psi^{(4)}]_{T,t}^{p,p,p,p}\biggr)^2\right\},
$$
$$
Q_p^{(2)}=\frac{1}{4}{\bf 1}_{\{i_1=i_2\ne 0\}}{\sf M}\left\{\left(
{\int\limits_t^{*}}^T
{\int\limits_t^{*}}^s
{\int\limits_t^{*}}^{s_1}
ds_2
d{\bf w}_{s_1}^{(i_3)}
d{\bf w}_{s}^{(i_4)}-S_1^{(i_3i_4)p}\right)^2\right\},
$$

$$
Q_p^{(3)}=\frac{1}{4}{\bf 1}_{\{i_2=i_3\ne 0\}}{\sf M}\left\{\left(
{\int\limits_t^{*}}^T
{\int\limits_t^{*}}^{s_2}
{\int\limits_t^{*}}^{s_1}
d{\bf w}_{s}^{(i_1)}ds_1
d{\bf w}_{s_2}^{(i_4)}-S_2^{(i_1i_4)p}\right)^2\right\},
$$

$$
Q_p^{(4)}=\frac{1}{4}{\bf 1}_{\{i_3=i_4\ne 0\}}{\sf M}\left\{\left(
{\int\limits_t^{*}}^T
{\int\limits_t^{*}}^{s_1}
{\int\limits_t^{*}}^{s_2}
d{\bf w}_{s}^{(i_1)}
d{\bf w}_{s_2}^{(i_2)}ds_1-S_3^{(i_1i_2)p}\right)^2\right\},
$$

$$
Q_p^{(5)}={\bf 1}_{\{i_1=i_2\ne 0\}}{\bf 1}_{\{i_3=i_4\ne 0\}}\times
$$
$$
\times\left(\frac{1}{4}
\int\limits_t^T(s_1-t)ds_1-
\sum\limits_{j_4=0}^{p}
\frac{1}{2}\int\limits_t^T\phi_{j_4}(s)\int\limits_{t}^s\phi_{j_4}(s_1)
(s_1-t)ds_1ds\right)^2,
$$

\vspace{3mm}
$$
Q_p^{(6)}={\sf M}\left\{\left(R_p\right)^2\right\}.
$$

\vspace{4mm}

From Remark 1.7 (see (\ref{zsel1})) we have
\begin{equation}
\label{mayx101}
Q_p^{(1)}\le \frac{C_1}{p},
\end{equation}

\noindent
where constant $C_1$ is independent of $p.$

Let us prove the version of Theorem 2.25 for the case $i_1,i_2,i_3=0,1,\ldots,m$.
The case $i_1,i_2,i_3=1,\ldots,m$
has been proved in Theorem 2.25. 
It is easy to see that, in addition to the proof of 
Theorem 2.25, we need to prove the following inequalities
$$
\left|\frac{1}{2}\int\limits_t^T
\psi_3(t_3)\int\limits_t^{t_3}
\psi_1(t_1)\psi_2(t_1)dt_1dt_3
-\right.
$$
\begin{equation}
\label{ag001}
~~~~~~~\left.-\sum\limits_{j_1=0}^{p}
\int\limits_t^T
\psi_3(t_3)\int\limits_t^{t_3}\phi_{j_1}(t_2)\psi_2(t_2)
\int\limits_t^{t_2}\phi_{j_1}(t_1)\psi_1(t_1)dt_1 dt_2 dt_3\right|\le\frac{C}{p},
\end{equation}
$$
\left|\frac{1}{2}\int\limits_t^T
\psi_3(t_3)\psi_2(t_3)\int\limits_t^{t_3}
\psi_1(t_1)dt_1dt_3-\right.
$$
\begin{equation}
\label{ag002}
~~~~~~~\left.-
\sum\limits_{j_3=0}^{p}
\int\limits_t^T
\phi_{j_3}(t_3)\psi_3(t_3)\int\limits_t^{t_3}\phi_{j_3}(t_2)\psi_3(t_2)
\int\limits_t^{t_2}\psi_1(t_1)dt_1 dt_2 dt_3\right|\le\frac{C}{p},
\end{equation}

\vspace{2mm}
\begin{equation}
\label{ag0020}
~~~~~~~\left|\sum\limits_{j_1=0}^{p}
\int\limits_t^T
\phi_{j_1}(t_3)\psi_3(t_3)\int\limits_t^{t_3}\psi_2(t_2)
\int\limits_t^{t_2}\phi_{j_1}(t_1)\psi_1(t_1)dt_1 dt_2 dt_3\right|\le\frac{C}{p},
\end{equation}

\vspace{4mm}
\noindent
where constant $C$ is independent of $p.$

The inequalities (\ref{ag001}) and (\ref{ag002})
are equivalent to the following inequalities (see the proof of 
the cases (\ref{agentqq1}), (\ref{agentqq2})) 
\begin{equation}
\label{ag003}
\left|\frac{1}{2}\int\limits_t^T
\psi_1(t_2)
\tilde\psi_2(t_2)dt_2
-
\sum\limits_{j_1=0}^{p}
\int\limits_t^T
\phi_{j_1}(t_2)\tilde \psi_2(t_2)
\int\limits_t^{t_2}\phi_{j_1}(t_1)\psi_1(t_1)dt_1 dt_2\right|\le\frac{C}{p},
\end{equation}
\begin{equation}
\label{ag004}
\left|\frac{1}{2}\int\limits_t^T
\psi_3(t_3)
\bar\psi_2(t_3)dt_3
-
\sum\limits_{j_3=0}^{p}
\int\limits_t^T
\phi_{j_3}(t_3)\psi_3(t_3)
\int\limits_t^{t_3}\phi_{j_3}(t_2)\bar \psi_2(t_2)dt_2 dt_3\right|\le\frac{C}{p},
\end{equation}

\noindent
where $\tilde\psi_2(t_2),$ $\bar \psi_2(t_2)$ are defined by (\ref{agent00x}) and (\ref{ag010}),
respectively.
The inequalities (\ref{ag003}), (\ref{ag004}) follow from
(\ref{may1000}), (\ref{may2001})--(\ref{agentt11}).

Let us prove (\ref{ag0020}). By analogy with the proof of (\ref{agentqq8}) we have
$$
\sum\limits_{j_1=0}^{p}
\int\limits_t^T
\phi_{j_1}(t_3)\psi_3(t_3)\int\limits_t^{t_3}\psi_2(t_2)
\int\limits_t^{t_2}\phi_{j_1}(t_1)\psi_1(t_1)dt_1 dt_2 dt_3=
$$
$$
=\sum\limits_{j_1=0}^{p}
\int\limits_t^T
\phi_{j_1}(t_3)\psi_3(t_3)\int\limits_t^{t_3}
\phi_{j_1}(t_1)\tilde \psi_1(t_1)
dt_1dt_3-
$$
$$
-\sum\limits_{j_1=0}^{p}
\int\limits_t^T
\phi_{j_1}(t_3)\tilde \psi_3(t_3)\int\limits_t^{t_3}
\phi_{j_1}(t_1)\psi_1(t_1)
dt_1dt_3=
$$
$$
=\sum\limits_{j_1=0}^{\infty}
\int\limits_t^T
\phi_{j_1}(t_3)\psi_3(t_3)\int\limits_t^{t_3}
\phi_{j_1}(t_1)\tilde \psi_1(t_1)
dt_1dt_3-
$$
$$
-\sum\limits_{j_1=0}^{\infty}
\int\limits_t^T
\phi_{j_1}(t_3)\tilde \psi_3(t_3)\int\limits_t^{t_3}
\phi_{j_1}(t_1)\psi_1(t_1)
dt_1dt_3-
$$
$$
-\sum\limits_{j_1=p+1}^{\infty}
\int\limits_t^T
\phi_{j_1}(t_3)\psi_3(t_3)\int\limits_t^{t_3}
\phi_{j_1}(t_1)\tilde \psi_1(t_1)
dt_1dt_3+
$$
$$
+\sum\limits_{j_1=p+1}^{\infty}
\int\limits_t^T
\phi_{j_1}(t_3)\tilde \psi_3(t_3)\int\limits_t^{t_3}
\phi_{j_1}(t_1)\psi_1(t_1)
dt_1dt_3=
$$
$$
=-\sum\limits_{j_1=p+1}^{\infty}
\int\limits_t^T
\phi_{j_1}(t_3)\psi_3(t_3)\int\limits_t^{t_3}
\phi_{j_1}(t_1)\tilde \psi_1(t_1)
dt_1dt_3+
$$
\begin{equation}
\label{ag300}
+\sum\limits_{j_1=p+1}^{\infty}
\int\limits_t^T
\phi_{j_1}(t_3)\tilde \psi_3(t_3)\int\limits_t^{t_3}
\phi_{j_1}(t_1)\psi_1(t_1)
dt_1dt_3,
\end{equation}

\noindent
where
$\tilde \psi_1(t_1),$
$\tilde \psi_3(t_3)$ are defined by (\ref{ag200}), (\ref{ag201}), respectively.

Now the estimate (\ref{ag0020}) follows from (\ref{ag300}) and 
(\ref{may2001})--(\ref{agentt11}).
Theorem 2.25 is proved for the case $i_1,i_2,i_3=0,1,\ldots,m$.

Using the version of Theorem 2.25 for the case $i_1,i_2,i_3=0,1,\ldots,m$, we obtain
the following estimates
\begin{equation}
\label{mayx102}
Q_p^{(2)}\le \frac{C_2}{p},\ \ \ Q_p^{(3)}\le \frac{C_2}{p},\ \ \ Q_p^{(4)}\le \frac{C_2}{p},
\end{equation}

\noindent
where constant $C_2$ does not depend on $p.$

From Theorem 2.2 (see (\ref{may62021})) we get
$$
\frac{1}{2}
\int\limits_t^T(s_1-t)ds_1-
\sum\limits_{j_4=0}^{p}
\int\limits_t^T\phi_{j_4}(s)\int\limits_{t}^s\phi_{j_4}(s_1)
(s_1-t)ds_1ds=
$$
\begin{equation}
\label{mayx105}
=\sum_{j_4=p+1}^{\infty}
\int\limits_t^T\phi_{j_4}(s)\int\limits_{t}^s\phi_{j_4}(s_1)
(s_1-t)ds_1ds.
\end{equation}

\vspace{2mm}

Let us consider the case of Legendre polynomials.
From (\ref{may2001}) and (\ref{mayx105}) we have
\begin{equation}
\label{mayx106}
~ \left\vert\sum_{j_4=p+1}^{\infty}
\int\limits_t^T\phi_{j_4}(s)\int\limits_{t}^s\phi_{j_4}(s_1)
(s_1-t)ds_1ds\right\vert
\le \frac{C_3}{p}, 
\end{equation}

\noindent
where constant $C_3$ is independent of $p$.

By analogy with (\ref{agentt10}) and (\ref{agentt11}) we have the estimate (\ref{mayx106})
for the trigonometric case. Then
\begin{equation}
\label{mayx107}
Q_p^{(5)}\le \frac{C_4}{p^2},
\end{equation}

\noindent
where constant $C_4$ does not depend on $p.$

Analyzing the proof of Theorem 2.9, we conclude that
\begin{equation}
\label{mayx108}
Q_p^{(6)}\le \frac{C_5}{p}
\end{equation}

\noindent
for the polynomial and trigonometric cases; constant $C_5$ is independent of $p.$

Combining (\ref{mayx200}), (\ref{mayx101}), (\ref{mayx102}), (\ref{mayx107}), (\ref{mayx108}),
we get (\ref{may9000}). Theorem 2.26 is proved.

\section{Rate of the Mean-Square Convergence of Expansions of Iterated
Stra\-to\-no\-vich Stochastic Integrals of Multiplicities 2 to 4 in 
Theorems 2.18, 2.20, and 2.22 (The Case
of In\-teg\-ra\-tion Interval $[t, s]$ $(s\in (t, T])$)}

Let us prove the following theorem. 

{\bf Theorem 2.27} \cite{arxiv-5}. {\it Suppose that 
$\{\phi_j(x)\}_{j=0}^{\infty}$ is a complete orthonormal system of 
Legendre polynomials or trigonometric 
functions in the space $L_2([t, T]).$
Moreover$,$ $\psi_1(\tau),$ $\psi_2(\tau)$ are 
continuously differentiable functions on $[t, T]$. Then$,$ 
for the iterated Stratonovich stochastic integral
$$
J^{*}[\psi^{(2)}]_{s,t}=
{\int\limits_t^{*}}^s
\psi_2(t_2)
{\int\limits_t^{*}}^{t_2}
\psi_1(t_1)d{\bf w}_{t_1}^{(i_1)}
d{\bf w}_{t_2}^{(i_2)}\ \ \ (i_1, i_2=1,\ldots,m)
$$
the following estimate
\begin{equation}
\label{agent2000}
~~~~~~~~ {\sf M}\left\{\left(J^{*}[\psi^{(2)}]_{s,t}-\sum_{j_1, j_2=0}^{p}
C_{j_2j_1}(s)\zeta_{j_1}^{(i_1)}\zeta_{j_2}^{(i_2)}\right)^2\right\}\le \frac{C(s)}{p}
\end{equation}
is valid, where $s\in (t, T]$ $(s$ is fixed{\rm )},
constant $C(s)$ is independent of $p,$
$$
C_{j_2 j_1}(s)=\int\limits_t^s\psi_2(t_2)\phi_{j_2}(t_2)
\int\limits_t^{t_2}\psi_1(t_1)\phi_{j_1}(t_1)dt_1dt_2,
$$
and
$$
\zeta_{j}^{(i)}=
\int\limits_t^T \phi_{j}(\tau) d{\bf w}_{\tau}^{(i)}
$$ 
are independent
standard Gaussian random variables for various 
$i$ or $j$.}

{\bf Proof.} The case
$s=T$ has already been considered in Theorem 2.24.
Below we consider the case $s\in (t, T).$ 
By analogy with (\ref{may1002}) we obtain
$$
{\sf M}\left\{\left(J^{*}[\psi^{(2)}]_{s,t}-\sum_{j_1, j_2=0}^{p}
C_{j_2j_1}(s)\zeta_{j_1}^{(i_1)}\zeta_{j_2}^{(i_2)}\right)^2\right\}=
$$
$$
={\sf M}\left\{\biggl(
J[\psi^{(2)}]_{s,t}-J[\psi^{(2)}]_{s,t}^{p,p}
\biggr)^2\right\}+
$$
\begin{equation}
\label{agent10}
+ {\bf 1}_{\{i_1=i_2\}}\left(\frac{1}{2}
\int\limits_t^s\psi_1(t_1)\psi_2(t_1)dt_1- 
\sum_{j_1=0}^{p}
C_{j_1j_1}(s)
\right)^2,
\end{equation}

\noindent
where (see (\ref{a2xxx}))
$$
J[\psi^{(2)}]_{s,t}^{p,p}
=\sum_{j_1,j_2=0}^{p}
C_{j_2j_1}(s)\Biggl(\zeta_{j_1}^{(i_1)}\zeta_{j_2}^{(i_2)}
-{\bf 1}_{\{i_1=i_2\ne 0\}}
{\bf 1}_{\{j_1=j_2\}}\Biggr).
$$

\vspace{1mm}

From Remark 1.12 (see (\ref{road1888})) we have
\begin{equation}
\label{agent11}
{\sf M}\left\{\biggl(
J[\psi^{(2)}]_{s,t}-J[\psi^{(2)}]_{s,t}^{p,p}
\biggr)^2\right\}\le \frac{C_1(s)}{p},
\end{equation}

\noindent
where constant $C_1(s)$ is independent of $p.$

Using (\ref{5tzzz}), we obtain (the existence of a limit on the
right-hand side of (\ref{5tzzz}) and a useful estimate for this limit will be proved further in this section)
\begin{equation}
\label{agent0101}
\frac{1}{2}
\int\limits_t^s\psi_1(t_1)\psi_2(t_1)dt_1
-\sum_{j_1=0}^{p}
C_{j_1j_1}(s)=\sum_{j_1=p+1}^{\infty}
C_{j_1j_1}(s).
\end{equation}

Consider the case of Legendre polynomials.
By analogy with (\ref{tupo14}) we get
for $n>m$ $(n, m\in {\bf N})$
$$
\sum\limits_{j_1=m+1}^n
C_{j_1j_1}(s)=
\sum\limits_{j_1=m+1}^n
\int\limits_t^s \psi_2(\theta)\phi_{j_1}(\theta)
\int\limits_t^{\theta} \psi_1(\tau)\phi_{j_1}(\tau)d\tau d\theta=
$$
$$
=
\frac{T-t}{4}
\int\limits_{-1}^{z(s)}
\psi_1(h(x))\psi_2(h(x))
\left(P_{n+1}(x)P_{n}(x)
-P_{m+1}(x)P_{m}(x)\right)dx-
$$
$$
-\frac{(T-t)^2}{8}
\sum\limits_{j_1=m+1}^n
\frac{1}{2j_1+1}
\int\limits_{-1}^{z(s)}
\left(P_{j_1+1}(y)-P_{j_1-1}(y)\right)\psi_1'(h(y))\times
$$
$$
\times
\Biggl(\left(P_{j_1+1}(z(s))-P_{j_1-1}(z(s))\right)\psi_2(s)-
\left(P_{j_1+1}(y)-P_{j_1-1}(y)\right)\psi_2(h(y))-\Biggr.
$$
\begin{equation}
\label{agent14}
\Biggl.
-
\frac{T-t}{2}
\int\limits_{y}^{z(s)}
\left(P_{j_1+1}(x)-
P_{j_1-1}(x)\right)\psi_2'(h(x))dx\Biggr)dy,
\end{equation}

\vspace{2mm}
\noindent
where 
$$
h(y)=\frac{T-t}{2}y+\frac{T+t}{2},\ \ \
z(s)=\left(s-\frac{T+t}{2}\right)\frac{2}{T-t},
$$

\noindent
and $\psi_1'$, $\psi_2'$ are
derivatives of the functions $\psi_1(\tau)$, $\psi_2(\tau)$ with respect 
to the variable
$h(y)$ (see (\ref{tupo12})).

Applying the estimate (\ref{otit987ggg}) and tak\-ing into account 
the boundedness of the functions $\psi_1(\tau)$, $\psi_2(\tau)$
and their derivatives, we finally obtain
$$
\left\vert\sum\limits_{j_1=m+1}^n
C_{j_1j_1}(s)\right\vert\le
C_1\left(\frac{1}{n}+\frac{1}{m}\right)
\int\limits_{-1}^{z(s)} \frac{dx}{\left(1-x^2\right)^{1/2}}+
$$
$$
+C_2 \sum\limits_{j_1=m+1}^n \frac{1}{j_1^2}\left(
\int\limits_{-1}^{z(s)}\frac{dy}{\left(1-y^2\right)^{1/2}}
+\frac{1}{\left(1-z^2(s)\right)^{1/4}}
\int\limits_{-1}^{z(s)}\frac{dy}{\left(1-y^2\right)^{1/4}}+\right.
$$
\begin{equation}
\label{fin1000}
+
\left.\int\limits_{-1}^{z(s)}\frac{1}{\left(1-y^2\right)^{1/4}}
\int\limits_{y}^{z(s)}\frac{dx}{\left(1-x^2\right)^{1/4}}dy\right),
\end{equation}

\vspace{2mm}
\noindent
where constants $C_1, C_2$ do not depend on $n$ and $m$.

We assume that $s\in (t, T)$ $(z(s)\ne \pm 1)$ since the case
$s=T$ has already been considered in Theorem 2.24.
Then
\begin{equation}
\label{agent1515}
\left\vert\sum\limits_{j_1=m+1}^n
C_{j_1j_1}(s)\right\vert
\le
C_3(s)\left(\frac{1}{n}+\frac{1}{m}+\sum\limits_{j_1=m+1}^n 
\frac{1}{j_1^2}\right),
\end{equation}

\noindent
where constant $C_3(s)$ does not depend on $n$ and $m$.
Thus, the limit
\begin{equation}
\label{rum111}
\lim\limits_{p\to\infty}\sum\limits_{j_1=0}^p
C_{j_1j_1}(s)
\end{equation}

\noindent
exists for the polynomial case.
For the trigonometric case, the existence of the limit
(\ref{rum111}) can be proved by analogy with the proof of Lemma~2.2 
(Sect.~2.1.2). We also note that the existence of these limits
follows from Sect.~2.1.4.

The relations (\ref{agent1515}) and (\ref{obana}) imply that
\begin{equation}
\label{agent1516}
\left\vert\sum\limits_{j_1=p+1}^{\infty}
C_{j_1j_1}(s)\right\vert
\le C_3(s)\left(\frac{1}{p}+\sum\limits_{j_1=p+1}^{\infty}
\frac{1}{j_1^2}\right)\le   \frac{C_4(s)}{p},
\end{equation}

\noindent
where constant $C_4(s)$ is independent of $p$.

For the trigonometric case, the analog of the inequality (\ref{agent1516}) 
can be obtained by analogy with (\ref{agentt10}) and (\ref{agentt11})
(see the proof of Lemma~2.2).

Combining (\ref{agent10})--(\ref{agent0101}), (\ref{agent1516}), 
we obtain the estimate (\ref{agent2000}).
Theorem 2.27 is proved. 

The arguments given earlier in Chapters 1 and 2 of this book
allow us to formulate the following two theorems.

\vspace{2mm}

{\bf Theorem 2.28} \hspace{-1mm} \cite{arxiv-5}.
{\it \hspace{-3mm}  Suppose that 
$\{\phi_j(x)\}_{j=0}^{\infty}$ is a complete orthonormal system of 
Legendre polynomials or trigonometric functions in the space $L_2([t, T]).$
At the same time $\psi_2(\tau)$ is a continuously dif\-ferentiable 
nonrandom function on $[t, T]$ and $\psi_1(\tau),$ $\psi_3(\tau)$ are twice
continuously differentiable nonrandom functions on $[t, T]$. 
Then$,$ for the 
iterated Stratonovich stochastic integral of third multiplicity
$$
J^{*}[\psi^{(3)}]_{s,t}={\int\limits_t^{*}}^s\psi_3(t_3)
{\int\limits_t^{*}}^{t_3}\psi_2(t_2)
{\int\limits_t^{*}}^{t_2}\psi_1(t_1)
d{\bf w}_{t_1}^{(i_1)}
d{\bf w}_{t_2}^{(i_2)}d{\bf w}_{t_3}^{(i_3)},
$$

\noindent
where $i_1, i_2, i_3=0, 1,\ldots,m,$ the following estimate
$$
{\sf M}\left\{\left(J^{*}[\psi^{(3)}]_{s,t}-
\sum\limits_{j_1, j_2, j_3=0}^{p}
C_{j_3 j_2 j_1}(s)\zeta_{j_1}^{(i_1)}\zeta_{j_2}^{(i_2)}\zeta_{j_3}^{(i_3)}
\right)^2\right\}\le \frac{C(s)}{p}
$$

\vspace{1mm}
\noindent
is valid, where $s\in (t, T]$ $(s$ is fixed{\rm )}, constant $C(s)$ is independent of $p,$ 
$$
C_{j_3 j_2 j_1}(s)=\int\limits_t^s\psi_3(t_3)\phi_{j_3}(t_3)
\int\limits_t^{t_3}\psi_2(t_2)\phi_{j_2}(t_2)
\int\limits_t^{t_2}\psi_1(t_1)\phi_{j_1}(t_1)dt_1dt_2dt_3,
$$
and
$$
\zeta_{j}^{(i)}=
\int\limits_t^T \phi_{j}(\tau) d{\bf w}_{\tau}^{(i)}
$$ 
are independent standard Gaussian random variables for various 
$i$ or $j$.}

\vspace{2mm}

{\bf Theorem 2.29} \hspace{-1mm} \cite{arxiv-5}. {\it \hspace{-3mm}  Suppose that
$\{\phi_j(x)\}_{j=0}^{\infty}$ is a complete orthonormal
system of Legendre polynomials or trigonometric functions
in the space $L_2([t, T])$.
Then$,$ for the iterated 
Stratonovich stochastic integral of fourth multiplicity
$$
J^{*}[\psi^{(4)}]_{s,t}=
{\int\limits_t^{*}}^s
{\int\limits_t^{*}}^{t_4}
{\int\limits_t^{*}}^{t_3}
{\int\limits_t^{*}}^{t_2}
d{\bf w}_{t_1}^{(i_1)}
d{\bf w}_{t_2}^{(i_2)}d{\bf w}_{t_3}^{(i_3)}d{\bf w}_{t_4}^{(i_4)}\ \ \ 
(i_1, i_2, i_3, i_4=0, 1,\ldots,m)
$$
the following estimate
$$
{\sf M}\left\{\left(J^{*}[\psi^{(4)}]_{s,t}-
\sum\limits_{j_1, j_2, j_3,j_4=0}^{p}
C_{j_4j_3 j_2 j_1}(s)\zeta_{j_1}^{(i_1)}\zeta_{j_2}^{(i_2)}\zeta_{j_3}^{(i_3)}\zeta_{j_4}^{(i_4)}
\right)^2\right\}\le \frac{C(s)}{p}
$$

\vspace{1mm}
\noindent
is valid, where $s\in (t, T]$ $(s$ is fixed{\rm )}, constant $C(s)$ is independent of $p,$
$$
C_{j_4 j_3 j_2 j_1}(s)=\int\limits_t^s\phi_{j_4}(t_4)\int\limits_t^{t_4}
\phi_{j_3}(t_3)
\int\limits_t^{t_3}\phi_{j_2}(t_2)\int\limits_t^{t_2}\phi_{j_1}(t_1)
dt_1dt_2dt_3dt_4,
$$
and
$$
\zeta_{j}^{(i)}=
\int\limits_t^T \phi_{j}(\tau) d{\bf w}_{\tau}^{(i)}
$$ 
are independent standard Gaussian random variables for various 
$i$ or $j$ {\rm (}in the case when $i\ne 0${\rm ),}
${\bf w}_{\tau}^{(i)}$
$(i=1,\ldots,m)$ are independent standard Wiener processes$,$
${\bf w}_{\tau}^{(0)}=\tau.$}

\section{Expansion of Iterated 
Stratonovich Stochastic Integrals of Arbitrary 
Multiplicity $k$ $(k\in{\bf N})$. Proof of Hypotheses 2.2 and 2.3 Under the Condition of 
Convergence of Trace Series}

In this section, we 
prove the expansion of iterated 
Stratonovich stochastic integrals of arbitrary 
multiplicity $k$ ($k\in {\bf N}$)
under the condition of 
convergence of trace series.
Let us introduce some notations.

Consider the unordered
set $\{1, 2, \ldots, k\}$ 
and separate it into two parts:
the first part consists of $r$ unordered 
pairs (sequence order of these pairs is also unimportant) and the 
second one consists of the 
remaining $k-2r$ numbers.
So, we have

\vspace{-4mm}
\begin{equation}
\label{leto5007after}
(\{
\underbrace{\{g_1, g_2\}, \ldots, 
\{g_{2r-1}, g_{2r}\}}_{\small{\hbox{part 1}}}
\},
\{\underbrace{q_1, \ldots, q_{k-2r}}_{\small{\hbox{part 2}}}
\}),
\end{equation}

\vspace{2mm}
\noindent
where 
$$
\{g_1, g_2, \ldots, 
g_{2r-1}, g_{2r}, q_1, \ldots, q_{k-2r}\}=\{1, 2, \ldots, k\},
$$

\vspace{1mm}
\noindent
braces   
mean an unordered 
set, and pa\-ren\-the\-ses mean an ordered set.

Consider the sum (\ref{leto5008}) with respect to all possible
partitions (\ref{leto5007after})

\vspace{-3mm}
$$
\sum_{\stackrel{(\{\{g_1, g_2\}, \ldots, 
\{g_{2r-1}, g_{2r}\}\}, \{q_1, \ldots, q_{k-2r}\})}
{{}_{\{g_1, g_2, \ldots, 
g_{2r-1}, g_{2r}, q_1, \ldots, q_{k-2r}\}=\{1, 2, \ldots, k\}}}}
a_{g_1 g_2, \ldots, 
g_{2r-1} g_{2r}, q_1 \ldots q_{k-2r}}
$$

\vspace{2mm}
\noindent
and the Fourier coefficient
\begin{equation}
\label{after3000}
~~~~~~~~~~C_{j_k \ldots j_1}=\int\limits_t^T\psi_k(t_k)\phi_{j_k}(t_k)\ldots
\int\limits_t^{t_2}
\psi_1(t_1)\phi_{j_1}(t_1)
dt_1\ldots dt_k
\end{equation}
corresponding to the function {\rm (\ref{ppp}), where 
$\{\phi_j(x)\}_{j=0}^{\infty}$ is a complete orthonormal system
of functions 
in the space $L_2([t, T])$. At that we suppose $\phi_0(x)=1/\sqrt{T-t}.$ 

Denote
$$
C_{j_k \ldots j_{l+1}j_l j_l j_{l-2} \ldots j_1}\biggl|_{(j_l j_l)\curvearrowright (\cdot) }\biggr.
\stackrel{\sf def}{=}
$$
$$
\stackrel{\sf def}{=}\int\limits_t^T\psi_k(t_k)\phi_{j_k}(t_k)\ldots
\int\limits_t^{t_{l+2}}\psi_{l+1}(t_{l+1})\phi_{j_{l+1}}(t_{l+1})
\int\limits_t^{t_{l+1}}\psi_{l}(t_{l})\psi_{l-1}(t_{l})\times
$$
\begin{equation}
\label{after900}
\times
\int\limits_t^{t_{l}}\psi_{l-2}(t_{l-2})\phi_{j_{l-2}}(t_{l-2})\ldots 
\int\limits_t^{t_2}
\psi_1(t_1)\phi_{j_1}(t_1)
dt_1\ldots dt_{l-2}dt_{l}t_{l+1}\ldots dt_k=
\end{equation}
$$
=\sqrt{T-t}\int\limits_t^T\psi_k(t_k)\phi_{j_k}(t_k)\ldots
\int\limits_t^{t_{l+2}}\psi_{l+1}(t_{l+1})\phi_{j_{l+1}}(t_{l+1})
\int\limits_t^{t_{l+1}}\psi_{l}(t_{l})\psi_{l-1}(t_{l})\phi_0(t_l)\times
$$
$$
\times
\int\limits_t^{t_{l}}\psi_{l-2}(t_{l-2})\phi_{j_{l-2}}(t_{l-2})\ldots 
\int\limits_t^{t_2}
\psi_1(t_1)\phi_{j_1}(t_1)
dt_1\ldots dt_{l-2}dt_{l}t_{l+1}\ldots dt_k=
$$

$$
=\sqrt{T-t}\hat C_{j_k \ldots j_{l+1} 0 j_{l-2} \ldots j_1},
$$

\vspace{2mm}
\noindent
i.e. $\sqrt{T-t}\hat C_{j_k \ldots j_{l+1} 0 j_{l-2} \ldots j_1}$
is again the Fourier coefficient of type $C_{j_k \ldots j_1}$
but with a new shorter multi-index
$j_k \ldots j_{l+1} 0 j_{l-2} \ldots j_1$ 
and new weight functions $\psi_1(\tau),$ $\ldots,$ $\psi_{l-2}(\tau),$
$\sqrt{T-t}\psi_{l-1}(\tau)\psi_{l}(\tau),$ $\psi_{l+1}(\tau)$, $\ldots,$
$\psi_{k}(\tau)$ (also we suppose that $\{l, l-1\}$ is one of the pairs
$\{g_1, g_2\}, \ldots,\{g_{2r-1}, g_{2r}\}$).

Let
$$
C_{j_k \ldots j_{l+1}j_l j_l j_{l-2} \ldots j_1}\biggl|_{(j_l j_l)\curvearrowright j_m}\biggr.
\stackrel{\sf def}{=}
$$

\vspace{-3mm}
$$
\stackrel{\sf def}{=}\int\limits_t^T\psi_k(t_k)\phi_{j_k}(t_k)\ldots
\int\limits_t^{t_{l+2}}\psi_{l+1}(t_{l+1})\phi_{j_{l+1}}(t_{l+1})
\int\limits_t^{t_{l+1}}\psi_{l}(t_{l})\psi_{l-1}(t_{l})\phi_{j_m}(t_l)\times
$$
\begin{equation}
\label{after2000}
\times
\int\limits_t^{t_{l}}\psi_{l-2}(t_{l-2})\phi_{j_{l-2}}(t_{l-2})\ldots 
\int\limits_t^{t_2}
\psi_1(t_1)\phi_{j_1}(t_1)
dt_1\ldots dt_{l-2}dt_{l}t_{l+1}\ldots dt_k=
\end{equation}

\vspace{-1mm}
$$
= \bar C_{j_k \ldots j_{l+1} j_m j_{l-2} \ldots j_1},
$$

\vspace{2mm}
\noindent
i.e. $\bar C_{j_k \ldots j_{l+1} j_m j_{l-2} \ldots j_1}$
is again the Fourier coefficient of type $C_{j_k \ldots j_1}$
but with a new shorter multi-index
$j_k \ldots j_{l+1} j_m j_{l-2} \ldots j_1$ 
and new weight functions $\psi_1(\tau),$ $\ldots,$ $\psi_{l-2}(\tau),$
$\psi_{l-1}(\tau)\psi_{l}(\tau),$ $\psi_{l+1}(\tau)$, $\ldots,$
$\psi_{k}(\tau)$ (also we suppose that $\{l, l-1\}$ is one of the pairs
$\{g_1, g_2\}, \ldots,\{g_{2r-1}, g_{2r}\}$).

Denote
$$
\bar C^{(p)}_{j_k\ldots j_q \ldots j_1}\biggl|_{q\ne g_1,g_2,\ldots,g_{2r-1}, g_{2r}}
\stackrel{\sf def}{=}
$$

\vspace{-3mm}
\begin{equation}
\label{after5days}
~~~~~~~~~\stackrel{\sf def}{=}
\sum\limits_{j_{g_{2r-1}}=p+1}^{\infty}
\sum\limits_{j_{g_{2r-3}}=p+1}^{\infty}
\ldots \sum\limits_{j_{g_{3}}=p+1}^{\infty}
\sum\limits_{j_{g_{1}}=p+1}^{\infty}
C_{j_k \ldots j_1}\biggl|_{j_{g_1}=j_{g_2},\ldots, j_{g_{2r-1}}=j_{g_{2r}}}\biggl..
\end{equation}

Introduce the following notation

\vspace{-3mm}
$$
S_l \left\{\bar C^{(p)}_{j_k\ldots j_q \ldots j_1}\biggl|_{q\ne g_1,g_2,\ldots,g_{2r-1}, g_{2r}}
\right\}
\stackrel{\sf def}{=}
\frac{1}{2}{\bf 1}_{\{g_{2l}=g_{2l-1}+1\}}
\sum\limits_{j_{g_{2r-1}}=p+1}^{\infty}
\sum\limits_{j_{g_{2r-3}}=p+1}^{\infty}
\ldots 
$$

\vspace{-4mm}
\begin{equation}
\label{sars300}
\ldots
\sum\limits_{j_{g_{2l+1}}=p+1}^{\infty}
\sum\limits_{j_{g_{2l-3}}=p+1}^{\infty}
\ldots
\sum\limits_{j_{g_{3}}=p+1}^{\infty}
\sum\limits_{j_{g_{1}}=p+1}^{\infty}
C_{j_k \ldots j_1}\biggl|_{(j_{g_{2l}} j_{g_{2l-1}})\curvearrowright (\cdot),
j_{g_1}=j_{g_2},\ldots, j_{g_{2r-1}}=j_{g_{2r}}}\biggr..
\end{equation}

\vspace{4mm}

Note that the operation $S_l$ $(l=1,2,\ldots,r)$ acts on the value

\vspace{-3mm}
\begin{equation}
\label{after301}
\bar C^{(p)}_{j_k\ldots j_q \ldots j_1}\biggl|_{q\ne g_1,g_2,\ldots,g_{2r-1}, g_{2r}}
\end{equation}

\vspace{2mm}
\noindent
as follows: $S_l$ multiplies (\ref{after301})
by ${\bf 1}_{\{g_{2l}=g_{2l-1}+1\}}/2,$ removes the summation
$$
\sum\limits_{j_{g_{2l-1}}=p+1}^{\infty},
$$

\noindent
and replaces 
$$
C_{j_k \ldots j_1}\biggl|_{j_{g_1}=j_{g_2},\ldots, j_{g_{2r-1}}=j_{g_{2r}}}\biggl.
$$
with
\begin{equation}
\label{after300}
C_{j_k \ldots j_1}\biggl|_{(j_{g_{2l}} j_{g_{2l-1}})\curvearrowright (\cdot),
j_{g_1}=j_{g_2},\ldots, j_{g_{2r-1}}=j_{g_{2r}}}\biggr..
\end{equation}

\vspace{1mm}

Note that we write
$$
C_{j_k \ldots j_1}\biggl|_{(j_{g_{1}} j_{g_{2}})\curvearrowright (\cdot),
j_{g_{1}}=j_{g_{2}}}\biggr.
=
C_{j_k \ldots j_1}\biggl|_{(j_{g_{1}} j_{g_{1}})\curvearrowright (\cdot),
j_{g_{1}}=j_{g_{2}}},\biggr.
$$
$$                                                                   
C_{j_k \ldots j_1}\biggl|_{(j_{g_{1}} j_{g_{2}})\curvearrowright j_{m},
j_{g_{1}}=j_{g_{2}}}\biggr.
=
C_{j_k \ldots j_1}\biggl|_{(j_{g_{1}} j_{g_{1}})\curvearrowright j_{m},
j_{g_{1}}=j_{g_{2}}}\biggr.,\ \ \ 
$$

\vspace{2mm}
$$
C_{j_k \ldots j_1}\biggl|_{(j_{g_{1}} j_{g_{2}})\curvearrowright (\cdot),
(j_{g_{3}} j_{g_{4}})\curvearrowright (\cdot),
j_{g_{1}}=j_{g_{2}}, j_{g_{3}}=j_{g_{4}}}\biggr.
=
$$
$$
=
C_{j_k \ldots j_1}\biggl|_{(j_{g_{1}} j_{g_{1}})\curvearrowright (\cdot)
(j_{g_{3}} j_{g_{3}})\curvearrowright (\cdot),
j_{g_{1}}=j_{g_{2}}, j_{g_{3}}=j_{g_{4}}}\biggr.\ \ \ {\rm etc}.
$$

\vspace{4mm}
        
Since (\ref{after300}) is again the Fourier coefficient, 
then the action 
of superposition $S_l S_m$ on 
(\ref{after300}) is obvious. For example,
for $r=3$ 
$$
S_3S_2S_1
\left\{\bar C^{(p)}_{j_k\ldots j_q \ldots j_1}\biggl|_{q\ne g_1,g_2,\ldots,g_{5}, g_{6}}
\right\}
=
$$
$$
=\frac{1}{2^3}
\prod_{s=1}^3
{\bf 1}_{\{g_{2s}=g_{2s-1}+1\}}
C_{j_k \ldots j_1}\Biggl|_{(j_{g_2} j_{g_1})\curvearrowright (\cdot)
(j_{g_{4}} j_{g_{3}})\curvearrowright (\cdot)
(j_{g_{6}} j_{g_{5}})\curvearrowright (\cdot),
j_{g_{1}}=j_{g_{2}}, j_{g_{3}}=j_{g_{4}}, j_{g_{5}}=j_{g_{6}}}\Biggr.,
$$

\vspace{3mm}
$$
S_3S_1
\left\{\bar C^{(p)}_{j_k\ldots j_q \ldots j_1}\biggl|_{q\ne g_1,g_2,\ldots,g_{5}, g_{6}}
\right\}
=
$$

\vspace{-2mm}
$$
=\frac{1}{2^2}
{\bf 1}_{\{g_{6}=g_{5}+1\}}{\bf 1}_{\{g_{2}=g_{1}+1\}}
\sum_{j_{g_3}=p+1}^{\infty}
C_{j_k \ldots j_1}\Biggl|_{(j_{g_{2}} j_{g_{1}})\curvearrowright (\cdot)
(j_{g_{6}} j_{g_{5}})\curvearrowright (\cdot),
j_{g_{1}}=j_{g_{2}}, j_{g_{3}}=j_{g_{4}}, j_{g_{5}}=j_{g_{6}}}\Biggr.,
$$

\vspace{3mm}
$$
S_2
\left\{\bar C^{(p)}_{j_k\ldots j_q \ldots j_1}\biggl|_{q\ne g_1,g_2,\ldots,g_{5}, g_{6}}
\right\}
=
$$

$$
=\frac{1}{2}
{\bf 1}_{\{g_{4}=g_{3}+1\}}
\sum_{j_{g_1}=p+1}^{\infty}\sum_{j_{g_5}=p+1}^{\infty}
C_{j_k \ldots j_1}\Biggl|_{(j_{g_{4}} j_{g_{3}})\curvearrowright (\cdot),
j_{g_{1}}=j_{g_{2}}, j_{g_{3}}=j_{g_{4}}, j_{g_{5}}=j_{g_{6}}}\Biggr..
$$

\vspace{4mm}

{\bf Theorem~2.30}\ \cite{arxiv-5}, \cite{arxiv-10}, \cite{arxiv-11},
\cite{new-art-1xxy}.\ {\it Assume that
the continuously differentiable functions 
$\psi_l(\tau)$ $(l=1,\ldots,k)$ at the interval $[t, T]$ and 
the complete orthonormal system $\{\phi_j(x)\}_{j=0}^{\infty}$
of continuous functions $(\phi_0(x)=1/\sqrt{T-t})$ 
in the space $L_2([t, T])$ are such that the following 
conditions are satisfied{\rm :}

{\rm 1.}\ The equality 
\begin{equation}
\label{after200}
~~~~~\frac{1}{2}
\int\limits_t^s \Phi_1(t_1)\Phi_2(t_1)dt_1
=\sum_{j=0}^{\infty}
\int\limits_t^s
\Phi_2(t_2)\phi_{j}(t_2)\int\limits_t^{t_2}
\Phi_1(t_1)\phi_{j}(t_1)dt_1 dt_2
\end{equation}

\noindent
holds for all $s\in (t, T],$ where the nonrandom functions 
$\Phi_1(\tau),$ $\Phi_2(\tau)$
are continuously differentiable on $[t, T]$
and the series on the right-hand side of {\rm (\ref{after200})}
converges absolutely.

{\rm 2.}\ The estimates
$$
\left|\int\limits_t^s \phi_{j}(\tau)\Phi_1(\tau)d\tau\right|
\le \frac{\Psi_1(s)}{j^{1/2+\alpha}},\ \ \ 
\left|\int\limits_s^T \phi_{j}(\tau)\Phi_2(\tau)d\tau\right|\le
\frac{\Psi_1(s)}{j^{1/2+\alpha}},
$$
$$
\left|\sum_{j=p+1}^{\infty}\int\limits_t^s
\Phi_2(\tau)\phi_{j}(\tau)\int\limits_t^{\tau}
\Phi_1(\theta)\phi_{j}(\theta)d\theta d\tau\right|\le \frac{\Psi_2(s)}{p^{\hspace{0.5mm}\beta}}
$$

\vspace{1mm}
\noindent
hold for all $s\in (t, T)$ and for some $\alpha, \beta >0,$ where 
$\Phi_1(\tau),$ $\Phi_2(\tau)$
are continuously differentiable nonrandom functions on $[t, T],$\ $j, p\in {\bf N},$
and
$$
\int\limits_t^T \Psi_1^2(\tau) d\tau<\infty,\ \ \ 
\int\limits_t^T \left|\Psi_2(\tau)\right| d\tau<\infty.
$$

\vspace{-1mm}

{\rm 3.}\ The condition
$$
\lim\limits_{p\to\infty}
\sum\limits_{\stackrel{j_1,\ldots,j_q,\ldots,j_k=0}{{}_{q\ne g_1, g_2, \ldots, g_{2r-1},
g_{2r}}}}^p
\left(S_{l_1}S_{l_2}\ldots S_{l_{d}}
\left\{\bar C^{(p)}_{j_k\ldots j_q \ldots j_1}\biggl|_{q\ne g_1,g_2,\ldots,g_{2r-1}, g_{2r}}
\right\}\right)^2=0
$$

\vspace{1mm}
\noindent
holds for all possible $g_1,g_2,\ldots,g_{2r-1},g_{2r}$ {\rm (}see {\rm (\ref{leto5007after}))}
and $l_1, l_2, \ldots, l_{d}$ such that
$l_1, l_2, \ldots, l_{d}\in \{1,2,\ldots, r\},$\
$l_1>l_2>\ldots >l_{d},$\ $d=0, 1, 2,\ldots, r-1,$\ 
where $r=1, 2,\ldots,[k/2]$ and

\vspace{-5mm}
$$
S_{l_1}S_{l_2}\ldots S_{l_{d}}
\left\{\bar C^{(p)}_{j_k\ldots j_q \ldots j_1}\biggl|_{q\ne g_1,g_2,\ldots,g_{2r-1}, g_{2r}}
\right\}\stackrel{\sf def}{=}
\bar C^{(p)}_{j_k\ldots j_q \ldots j_1}\biggl|_{q\ne g_1,g_2,\ldots,g_{2r-1}, g_{2r}}
$$
for $d=0.$

Then, for the iterated Stratonovich stochastic integral 
of arbitrary multiplicity $k$

\newpage
\noindent
\begin{equation}
\label{afterstr}
~~~~~~~~~~~J^{*}[\psi^{(k)}]_{T,t}^{(i_1\ldots i_k)}=
{\int\limits_t^{*}}^T
\psi_k(t_k) \ldots 
{\int\limits_t^{*}}^{t_{2}}
\psi_1(t_1) d{\bf w}_{t_1}^{(i_1)}\ldots
d{\bf w}_{t_k}^{(i_k)}
\end{equation}
the following 
expansion 
\begin{equation}
\label{after1}
J^{*}[\psi^{(k)}]_{T,t}^{(i_1\ldots i_k)}=
\hbox{\vtop{\offinterlineskip\halign{
\hfil#\hfil\cr
{\rm l.i.m.}\cr
$\stackrel{}{{}_{p\to \infty}}$\cr
}} }
\sum\limits_{j_1,\ldots,j_k=0}^{p}
C_{j_k \ldots j_1}\prod\limits_{l=1}^k \zeta_{j_l}^{(i_l)}
\end{equation}

\vspace{1mm}
\noindent
that converges in the mean-square sense is valid, where 
\begin{equation}
\label{after1000}
~~~~~~~~C_{j_k \ldots j_1}=\int\limits_t^T\psi_k(t_k)\phi_{j_k}(t_k)\ldots
\int\limits_t^{t_2}
\psi_1(t_1)\phi_{j_1}(t_1)
dt_1\ldots dt_k
\end{equation}

\noindent
is the Fourier coefficient, 
${\rm l.i.m.}$ is a limit in the mean-square sense,
$i_1, \ldots, i_k=0, 1,\ldots,m,$
$$
\zeta_{j}^{(i)}=
\int\limits_t^T \phi_{j}(\tau) d{\bf w}_{\tau}^{(i)}
$$ 

\noindent
are independent standard Gaussian random variables for various 
$i$ or $j$ {\rm (}in the case when $i\ne 0${\rm ),}
${\bf w}_{\tau}^{(i)}$ 
$(i=1,\ldots,m)$ are independent standard Wiener processes$,$
${\bf w}_{\tau}^{(0)}=\tau.$}

\vspace{1mm}

{\bf Proof.} First note that (\ref{after200}) is true (see the proof of Theorem~2.18). 
The proof of Theorem~2.30 ($k\ge 2$) will consist 
of several steps (the case $k=1$ is obvious (see (\ref{a1})).

{\bf Step~1.}\ Let us find a representation of the quantity
$$
\sum\limits_{j_1,\ldots,j_k=0}^{p}
C_{j_k \ldots j_1}\prod\limits_{l=1}^k \zeta_{j_l}^{(i_l)}
$$

\noindent
that will be convenient for further consideration.

Recall the equality (\ref{rezo7}) (also see (\ref{chain401}), (\ref{leto6000xxa}))
$$
J'[\phi_{j_1}\ldots \phi_{j_k}]_{T,t}^{(i_1\ldots i_k)}=\prod_{l=1}^k\zeta_{j_l}^{(i_l)}+
$$
\begin{equation}
\label{2023abc300}
+\sum\limits_{r=1}^{[k/2]}
(-1)^r 
\hspace{-3mm}\sum_{\stackrel{(\{\{g_1, g_2\}, \ldots, 
\{g_{2r-1}, g_{2r}\}\}, \{q_1, \ldots, q_{k-2r}\})}
{{}_{\{g_1, g_2, \ldots, 
g_{2r-1}, g_{2r}, q_1, \ldots, q_{k-2r}\}=\{1, 2, \ldots, k\}}}}
\prod\limits_{s=1}^r
{\bf 1}_{\{i_{g_{{}_{2s-1}}}=~i_{g_{{}_{2s}}}\ne 0\}}
\Biggl.{\bf 1}_{\{j_{g_{{}_{2s-1}}}=~j_{g_{{}_{2s}}}\}}
\prod_{l=1}^{k-2r}\zeta_{j_{q_l}}^{(i_{q_l})}
\end{equation}

\vspace{1mm}
\noindent
w.~p.~1, where notations are the same as in Theorem 1.2
and
$J'[\phi_{j_1}\ldots\phi_{j_k}]_{T,t}^{(i_1\ldots i_k)}$
is the multiple Wiener stochastic
integral (\ref{mult11}) (also see (\ref{mult11www})).             

From (\ref{2023abc300}) we obtain

\vspace{-1mm}
$$
\prod_{l=1}^k\zeta_{j_l}^{(i_l)}=J'[\phi_{j_1}\ldots \phi_{j_k}]_{T,t}^{(i_1\ldots i_k)}-
$$
\begin{equation}
\label{2023abc300xz1}
-\sum\limits_{r=1}^{[k/2]}
(-1)^r 
\hspace{-3mm}\sum_{\stackrel{(\{\{g_1, g_2\}, \ldots, 
\{g_{2r-1}, g_{2r}\}\}, \{q_1, \ldots, q_{k-2r}\})}
{{}_{\{g_1, g_2, \ldots, 
g_{2r-1}, g_{2r}, q_1, \ldots, q_{k-2r}\}=\{1, 2, \ldots, k\}}}}
\prod\limits_{s=1}^r
{\bf 1}_{\{i_{g_{{}_{2s-1}}}=~i_{g_{{}_{2s}}}\ne 0\}}
\Biggl.{\bf 1}_{\{j_{g_{{}_{2s-1}}}=~j_{g_{{}_{2s}}}\}}
\prod_{l=1}^{k-2r}\zeta_{j_{q_l}}^{(i_{q_l})}
\end{equation}

\noindent
w.~p.~1.

By iteratively applying the formula (\ref{2023abc300xz1})
(also see (\ref{a2})--(\ref{a6})), we obtain the following
representation of the product
$$
\prod_{l=1}^k\zeta_{j_l}^{(i_l)}
$$

\noindent
as the sum of some constant value and multiple Wiener stochastic integrals 
of 
multiplicities not exceeding $k$ 

\vspace{-1mm}
$$
\prod_{l=1}^k\zeta_{j_l}^{(i_l)}=J'[\phi_{j_1}\ldots \phi_{j_k}]_{T,t}^{(i_1\ldots i_k)}+
$$
$$
+\sum\limits_{r=1}^{[k/2]}
\sum_{\stackrel{(\{\{g_1, g_2\}, \ldots, 
\{g_{2r-1}, g_{2r}\}\}, \{q_1, \ldots, q_{k-2r}\})}
{{}_{\{g_1, g_2, \ldots, 
g_{2r-1}, g_{2r}, q_1, \ldots, q_{k-2r}\}=\{1, 2, \ldots, k\}}}}
\prod\limits_{s=1}^r
{\bf 1}_{\{i_{g_{{}_{2s-1}}}=~i_{g_{{}_{2s}}}\ne 0\}}
\Biggl.{\bf 1}_{\{j_{g_{{}_{2s-1}}}=~j_{g_{{}_{2s}}}\}}\times
$$

\vspace{3mm}
\begin{equation}
\label{after8xx1}
\times
J'[\phi_{j_{q_1}}\ldots \phi_{j_{q_{k-2r}}}]_{T,t}^{(i_{q_1}\ldots i_{q_{k-2r}})}\ \ \ \hbox{w.~p.~1,}
\end{equation}

\vspace{6mm}
\noindent
where
$$
J'[\phi_{j_{q_1}}\ldots \phi_{j_{q_{k-2r}}}]_{T,t}^{(i_{q_1}\ldots i_{q_{k-2r}})}
\stackrel{\sf def}{=}1
$$

\vspace{2mm}
\noindent
for $k=2r$.

Multiplying both sides of the equality (\ref{after8xx1}) by $C_{j_k\ldots j_1}$
and summing over $j_1,\ldots,j_k,$ we get w.~p.~1
\newpage
\noindent
$$
\sum_{j_1=0}^{p_1}\ldots\sum_{j_k=0}^{p_k}
C_{j_k\ldots j_1}
\prod_{l=1}^k\zeta_{j_l}^{(i_l)}
=
\sum_{j_1=0}^{p_1}\ldots\sum_{j_k=0}^{p_k}
C_{j_k\ldots j_1}
J'[\phi_{j_1}\ldots \phi_{j_k}]_{T,t}^{(i_1\ldots i_k)}+
$$

\vspace{1mm}
$$
+\sum_{j_1=0}^{p_1}\ldots\sum_{j_k=0}^{p_k}
C_{j_k\ldots j_1}
\sum\limits_{r=1}^{[k/2]}
\sum_{\stackrel{(\{\{g_1, g_2\}, \ldots, 
\{g_{2r-1}, g_{2r}\}\}, \{q_1, \ldots, q_{k-2r}\})}
{{}_{\{g_1, g_2, \ldots, 
g_{2r-1}, g_{2r}, q_1, \ldots, q_{k-2r}\}=\{1, 2, \ldots, k\}}}}
\prod\limits_{s=1}^r
{\bf 1}_{\{i_{g_{{}_{2s-1}}}=~i_{g_{{}_{2s}}}\ne 0\}}\times
$$

\begin{equation}
\label{after8}
~~~~~~~~~~~~~~\times{\bf 1}_{\{j_{g_{{}_{2s-1}}}=~j_{g_{{}_{2s}}}\}}
J'[\phi_{j_{q_1}}\ldots \phi_{j_{q_{k-2r}}}]_{T,t}^{(i_{q_1}\ldots i_{q_{k-2r}})}\ \ \ \hbox{w.~p.~1.}
\end{equation}

\vspace{5mm}

Denote
\begin{equation}
\label{afterxx1}
K_{p_1\ldots p_k}(t_1,\ldots,t_k)=
\sum_{j_1=0}^{p_1}\ldots\sum_{j_k=0}^{p_k}
C_{j_k\ldots j_1}
\prod_{l=1}^k\phi_{j_l}(t_l),
\end{equation}

\begin{equation}
\label{afterxx2}
K_{p_1\ldots p_k}^{g_1\ldots g_{2r}, q_1\ldots q_{k-2r}}(t_{q_1},\ldots,t_{q_{k-2r}})=
\sum_{j_1=0}^{p_1}\ldots\sum_{j_k=0}^{p_k}
C_{j_k\ldots j_1}
\prod\limits_{s=1}^r
{\bf 1}_{\{j_{g_{{}_{2s-1}}}=~j_{g_{{}_{2s}}}\}}
\prod_{l=1}^{k-2r}\phi_{j_{q_l}}(t_{q_l}),
\end{equation}

\vspace{4mm}
\noindent
where $C_{j_k\ldots j_1}$ is defined by (\ref{after1000}) and
\begin{equation}
\label{2023abc450}
\prod_{l=1}^{0}\phi_{j_{q_l}}(t_{q_l})\stackrel{\sf def}{=}1,
\end{equation}

\noindent
i.e. $k=2r$ in (\ref{2023abc450}).

The equality (\ref{after8}) can be written as 

\vspace{-1mm}
$$
J[K_{p_1\ldots p_k}]_{T,t}^{(i_1\ldots i_k)}=
J'[K_{p_1\ldots p_k}]_{T,t}^{(i_1\ldots i_k)}+
$$

\vspace{1mm}
$$
+\sum\limits_{r=1}^{[k/2]}
\sum_{\stackrel{(\{\{g_1, g_2\}, \ldots, 
\{g_{2r-1}, g_{2r}\}\}, \{q_1, \ldots, q_{k-2r}\})}
{{}_{\{g_1, g_2, \ldots, 
g_{2r-1}, g_{2r}, q_1, \ldots, q_{k-2r}\}=\{1, 2, \ldots, k\}}}}
\prod\limits_{s=1}^r
{\bf 1}_{\{i_{g_{{}_{2s-1}}}=~i_{g_{{}_{2s}}}\ne 0\}}\times
$$

\vspace{6mm}
\begin{equation}
\label{after7}
\times 
J'[K_{p_1\ldots p_k}^{g_1\ldots g_{2r}, q_1\ldots q_{k-2r}}]_{T,t}^{(i_{q_1}\ldots i_{q_{k-2r}})}
\end{equation}

\vspace{5mm}
\noindent
w.~p.~{\rm 1}, where
$K_{p_1\ldots p_k}(t_1,\ldots,t_k)$ and 
$K_{p_1\ldots p_k}^{g_1\ldots g_{2r}, q_1\ldots q_{k-2r}}(t_{q_1},\ldots,t_{q_{k-2r}})$
are defined by the equalities (\ref{afterxx1}), (\ref{afterxx2}),
$J[K_{p_1\ldots p_k}]_{T,t}^{(i_1\ldots i_k)}$ 
is the multiple Stra\-to\-no\-vich stochastic integral (\ref{30.34})
(also see {\rm (\ref{30.34ququ})}) and
$J'[K_{p_1\ldots p_k}]_{T,t}^{(i_1\ldots i_k)},$ 
$J'[K_{p_1\ldots p_k}^{g_1\ldots g_{2r}, q_1\ldots q_{k-2r}}]_{T,t}^{(i_{q_1}\ldots i_{q_{k-2r}})}$
are multiple Wiener stochastic
integrals defined by (\ref{mult11}) (also see (\ref{mult11www})).             

Passing to the limit 
$\hbox{\vtop{\offinterlineskip\halign{
\hfil#\hfil\cr
{\rm l.i.m.}\cr
$\stackrel{}{{}_{p_1,\ldots,p_k\to \infty}}$\cr
}} }$ ($p_1=\ldots =p_k=p$) in (\ref{after8}) or (\ref{after7}), 
we get w.~p.~1 (see Theorems~1.1, 1.2 and (\ref{drdr1}))
$$
\hbox{\vtop{\offinterlineskip\halign{
\hfil#\hfil\cr
{\rm l.i.m.}\cr
$\stackrel{}{{}_{p\to \infty}}$\cr
}} }
\sum_{j_1,\ldots,j_k=0}^{p}
C_{j_k \ldots j_1}\prod\limits_{l=1}^k \zeta_{j_l}^{(i_l)}
=
J[\psi^{(k)}]_{T,t}^{(i_1\ldots i_k)}+
$$

$$
+
\hbox{\vtop{\offinterlineskip\halign{
\hfil#\hfil\cr
{\rm l.i.m.}\cr
$\stackrel{}{{}_{p\to \infty}}$\cr
}} }\sum_{j_1,\ldots,j_k=0}^{p}
C_{j_k \ldots j_1}
\sum\limits_{r=1}^{[k/2]}
\sum_{\stackrel{(\{\{g_1, g_2\}, \ldots, 
\{g_{2r-1}, g_{2r}\}\}, \{q_1, \ldots, q_{k-2r}\})}
{{}_{\{g_1, g_2, \ldots, 
g_{2r-1}, g_{2r}, q_1, \ldots, q_{k-2r}\}=\{1, 2, \ldots, k\}}}}
\prod\limits_{s=1}^r
{\bf 1}_{\{i_{g_{{}_{2s-1}}}=~i_{g_{{}_{2s}}}\ne 0\}}\times
$$

\vspace{4mm}
\begin{equation}
\label{after3}
\times
{\bf 1}_{\{j_{g_{{}_{2s-1}}}=~j_{g_{{}_{2s}}}\}}
J'[\phi_{j_{q_1}}\ldots \phi_{j_{q_{k-2r}}}]_{T,t}^{(i_{q_1}\ldots i_{q_{k-2r}})}=
\end{equation}

\vspace{6mm}
$$
=
J[\psi^{(k)}]_{T,t}^{(i_1\ldots i_k)}+
$$

$$
+
\hbox{\vtop{\offinterlineskip\halign{
\hfil#\hfil\cr
{\rm l.i.m.}\cr
$\stackrel{}{{}_{p\to \infty}}$\cr
}} }
\sum\limits_{r=1}^{[k/2]}
\sum_{\stackrel{(\{\{g_1, g_2\}, \ldots, 
\{g_{2r-1}, g_{2r}\}\}, \{q_1, \ldots, q_{k-2r}\})}
{{}_{\{g_1, g_2, \ldots, 
g_{2r-1}, g_{2r}, q_1, \ldots, q_{k-2r}\}=\{1, 2, \ldots, k\}}}}
\prod\limits_{s=1}^r
{\bf 1}_{\{i_{g_{{}_{2s-1}}}=~i_{g_{{}_{2s}}}\ne 0\}}\times
$$

\vspace{6mm}
\begin{equation}
\label{after7xx}
\times 
J'[K_{p_1\ldots p_k}^{g_1\ldots g_{2r}, q_1\ldots q_{k-2r}}]_{T,t}^{(i_{q_1}\ldots i_{q_{k-2r}})}
\end{equation}

\vspace{7mm}
\noindent
w.~p.~{\rm 1}, where
$J[\psi^{(k)}]_{T,t}^{(i_1\ldots i_k)}$ is the iterated It\^{o} stochastic
integral
\begin{equation}
\label{afterito}
~~~~~~~~~~J[\psi^{(k)}]_{T,t}^{(i_1\ldots i_k)}=
\int\limits_t^{T}\psi_k(t_k) \ldots 
\int\limits_t^{t_{2}}
\psi_1(t_1) d{\bf w}_{t_1}^{(i_1)}\ldots
d{\bf w}_{t_k}^{(i_k)}.
\end{equation}

If we prove that w.~p.~1
\newpage
\noindent
$$
\sum_{r=1}^{\left[k/2\right]}\frac{1}{2^r}
\sum_{(s_r,\ldots,s_1)\in {\rm A}_{k,r}}
J[\psi^{(k)}]_{T,t}^{s_r,\ldots,s_1}=
$$

$$
=
\hbox{\vtop{\offinterlineskip\halign{
\hfil#\hfil\cr
{\rm l.i.m.}\cr
$\stackrel{}{{}_{p\to \infty}}$\cr
}} }\sum_{j_1,\ldots,j_k=0}^{p}
C_{j_k \ldots j_1}
\sum\limits_{r=1}^{[k/2]}
\sum_{\stackrel{(\{\{g_1, g_2\}, \ldots, 
\{g_{2r-1}, g_{2r}\}\}, \{q_1, \ldots, q_{k-2r}\})}
{{}_{\{g_1, g_2, \ldots, 
g_{2r-1}, g_{2r}, q_1, \ldots, q_{k-2r}\}=\{1, 2, \ldots, k\}}}}
\prod\limits_{s=1}^r
{\bf 1}_{\{i_{g_{{}_{2s-1}}}=~i_{g_{{}_{2s}}}\ne 0\}}\times
$$

\vspace{4mm}
\begin{equation}
\label{after4}
\times
{\bf 1}_{\{j_{g_{{}_{2s-1}}}=~j_{g_{{}_{2s}}}\}}
J'[\phi_{j_{q_1}}\ldots \phi_{j_{q_{k-2r}}}]_{T,t}^{(i_{q_1}\ldots i_{q_{k-2r}})},
\end{equation}

\vspace{7mm}
\noindent
then (see (\ref{after3}), (\ref{after4}), and Theorem~2.12)
$$
\hbox{\vtop{\offinterlineskip\halign{
\hfil#\hfil\cr
{\rm l.i.m.}\cr
$\stackrel{}{{}_{p\to \infty}}$\cr
}} }
\sum_{j_1,\ldots,j_k=0}^{p}
C_{j_k \ldots j_1}\prod\limits_{l=1}^k \zeta_{j_l}^{(i_l)}
=
$$

\vspace{-2mm}
\begin{equation}
\label{after333}
~~~~~=
J[\psi^{(k)}]_{T,t}^{(i_1\ldots i_k)}+
\sum_{r=1}^{\left[k/2\right]}\frac{1}{2^r}
\sum_{(s_r,\ldots,s_1)\in {\rm A}_{k,r}}
J[\psi^{(k)}]_{T,t}^{s_r,\ldots,s_1}=
J^{*}[\psi^{(k)}]_{T,t}^{(i_1\ldots i_k)}
\end{equation}

\vspace{3mm}
\noindent
w.~p.~1, where notations in (\ref{after333}) are the same as in Theorem~2.12. 
Thus Theorem~2.30 will be proved.

From (\ref{after7}) we have that 
the multiple Stratonovich stochastic integral $J[K_{p_1\ldots p_k}]_{T,t}^{(i_1\ldots i_k)}$ 
of multiplicity $k$ is 
expressed as a sum of some constant value and
multiple Wiener stochastic 
integrals 

\vspace{-2mm}
$$
J'[K_{p_1\ldots p_k}]_{T,t}^{(i_1\ldots i_k)}
$$

\vspace{2mm}
\noindent
and 
$$
J'[K_{p_1\ldots p_k}^{g_1\ldots g_{2r}, q_1\ldots q_{k-2r}}]_{T,t}^{(i_{q_1}\ldots i_{q_{k-2r}})}
$$

\vspace{2.5mm}
\noindent 
of multiplicities $k,$ $k-2,$ $k-4$, $\ldots,$ $k-2[k/2]$  
($r=1,2,\ldots,[k/2]$).

The formulas (\ref{after8}), (\ref{after7}) can be considered
as new representations
of the Hu--Meyer formula
for the case of a multidimensional Wiener process 
\cite{Rybakov3000} (also see \cite{bugh1}, \cite{bugh3})
and kernel 
$K_{p_1\ldots p_k}(t_1,\ldots,t_k)$ (see (\ref{afterxx1})).

Further, we will use the representation (\ref{after8}) for $p_1=\ldots=p_k=p,$ i.e.
$$
\sum_{j_1,\ldots,j_k=0}^{p}
C_{j_k\ldots j_1}
\prod_{l=1}^k \zeta_{j_l}^{(i_l)}=
\sum_{j_1,\ldots,j_k=0}^{p}
C_{j_k\ldots j_1}
J'[\phi_{j_1}\ldots \phi_{j_k}]_{T,t}^{(i_1\ldots i_k)}+
$$

\vspace{-2mm}
$$
+\sum_{j_1,\ldots,j_k=0}^{p}
C_{j_k\ldots j_1}
\sum\limits_{r=1}^{[k/2]}
\sum_{\stackrel{(\{\{g_1, g_2\}, \ldots, 
\{g_{2r-1}, g_{2r}\}\}, \{q_1, \ldots, q_{k-2r}\})}
{{}_{\{g_1, g_2, \ldots, 
g_{2r-1}, g_{2r}, q_1, \ldots, q_{k-2r}\}=\{1, 2, \ldots, k\}}}}
\prod\limits_{s=1}^r
{\bf 1}_{\{i_{g_{{}_{2s-1}}}=~i_{g_{{}_{2s}}}\ne 0\}}\times
$$

\vspace{3mm}
\begin{equation}
\label{after8xx}
\times{\bf 1}_{\{j_{g_{{}_{2s-1}}}=~j_{g_{{}_{2s}}}\}}
J'[\phi_{j_{q_1}}\ldots \phi_{j_{q_{k-2r}}}]_{T,t}^{(i_{q_1}\ldots i_{q_{k-2r}})}\ \ \ \hbox{w.~p.~1.}
\end{equation}

\vspace{6mm}

For example, for $k=2,3,4,5,6$ from (\ref{after8}) we have w.~p.~1

\vspace{-2mm}
\begin{equation}
\label{after32}
~~\sum_{j_1=0}^{p_1}\sum_{j_2=0}^{p_2}
C_{j_2j_1}\zeta_{j_1}^{(i_1)}\zeta_{j_2}^{(i_2)}=
J'[K_{p_1p_2}]_{T,t}^{(i_1 i_2)}
+\sum_{j_1=0}^{p_1}\sum_{j_2=0}^{p_2}
C_{j_2j_1}
{\bf 1}_{\{i_1=i_2\ne 0\}}{\bf 1}_{\{j_1=j_2\}},
\end{equation}

\vspace{4mm}
$$
\sum_{j_1=0}^{p_1}\sum_{j_2=0}^{p_2}\sum_{j_3=0}^{p_3}
C_{j_3j_2j_1}
\zeta_{j_1}^{(i_1)}\zeta_{j_2}^{(i_2)}\zeta_{j_3}^{(i_3)}=
J'[K_{p_1p_2p_3}]_{T,t}^{(i_1 i_2 i_3)}+
$$
$$
+
\sum_{j_1=0}^{p_1}\sum_{j_2=0}^{p_2}\sum_{j_3=0}^{p_3}
C_{j_3j_2j_1}
\Biggl(
{\bf 1}_{\{i_1=i_2\ne 0\}}
{\bf 1}_{\{j_1=j_2\}}
J'[\phi_{j_3}]^{(i_3)}_{T,t}
+{\bf 1}_{\{i_2=i_3\ne 0\}}
{\bf 1}_{\{j_2=j_3\}}
J'[\phi_{j_1}]^{(i_1)}_{T,t}+\Biggr.
$$
\begin{equation}
\label{after33}
\Biggl.
+{\bf 1}_{\{i_1=i_3\ne 0\}}
{\bf 1}_{\{j_1=j_3\}}
J'[\phi_{j_2}]^{(i_2)}_{T,t}\Biggr),
\end{equation}

\vspace{4mm}
$$
\sum_{j_1=0}^{p_1}\ldots\sum_{j_4=0}^{p_4}
C_{j_4 j_3 j_2 j_1}
\zeta_{j_1}^{(i_1)}\zeta_{j_2}^{(i_2)}\zeta_{j_3}^{(i_3)}
\zeta_{j_4}^{(i_4)}=
J'[K_{p_1p_2p_3p_4}]_{T,t}^{(i_1 i_2 i_3 i_4)}+
$$
$$
+\sum_{j_1=0}^{p_1}\ldots\sum_{j_4=0}^{p_4}
C_{j_4 j_3 j_2 j_1}\Biggl(
\Biggr.
{\bf 1}_{\{i_1=i_2\ne 0\}}
{\bf 1}_{\{j_1=j_2\}}
J'[\phi_{j_3}\phi_{j_4}]^{(i_3 i_4)}_{T,t}
+
$$
$$
+
{\bf 1}_{\{i_1=i_3\ne 0\}}
{\bf 1}_{\{j_1=j_3\}}
J'[\phi_{j_2}\phi_{j_4}]^{(i_2 i_4)}_{T,t}
+
{\bf 1}_{\{i_1=i_4\ne 0\}}
{\bf 1}_{\{j_1=j_4\}}
J'[\phi_{j_2}\phi_{j_3}]^{(i_2 i_3)}_{T,t}
+
$$
$$
+
{\bf 1}_{\{i_2=i_3\ne 0\}}
{\bf 1}_{\{j_2=j_3\}}
J'[\phi_{j_1}\phi_{j_4}]^{(i_1 i_4)}_{T,t}
+
{\bf 1}_{\{i_2=i_4\ne 0\}}
{\bf 1}_{\{j_2=j_4\}}
J'[\phi_{j_1}\phi_{j_3}]^{(i_1 i_3)}_{T,t}
+
$$
$$
+
{\bf 1}_{\{i_3=i_4\ne 0\}}
{\bf 1}_{\{j_3=j_4\}}
J'[\phi_{j_1}\phi_{j_2}]^{(i_1 i_2)}_{T,t}
+
$$
$$
+
{\bf 1}_{\{i_1=i_2\ne 0\}}
{\bf 1}_{\{j_1=j_2\}}
{\bf 1}_{\{i_3=i_4\ne 0\}}
{\bf 1}_{\{j_3=j_4\}}
+
{\bf 1}_{\{i_1=i_3\ne 0\}}
{\bf 1}_{\{j_1=j_3\}}
{\bf 1}_{\{i_2=i_4\ne 0\}}
{\bf 1}_{\{j_2=j_4\}}+
$$
\begin{equation}
\label{after34}
+\Biggl.
{\bf 1}_{\{i_1=i_4\ne 0\}}
{\bf 1}_{\{j_1=j_4\}}
{\bf 1}_{\{i_2=i_3\ne 0\}}
{\bf 1}_{\{j_2=j_3\}}\Biggr),
\end{equation}

\vspace{3mm}
$$
\sum_{j_1=0}^{p_1}\ldots\sum_{j_5=0}^{p_5}
C_{j_5 j_4 j_3 j_2 j_1}
\zeta_{j_1}^{(i_1)}\zeta_{j_2}^{(i_2)}\zeta_{j_3}^{(i_3)}
\zeta_{j_4}^{(i_4)}\zeta_{j_5}^{(i_5)}
=
J'[K_{p_1p_2p_3p_4p_5}]_{T,t}^{(i_1 i_2 i_3 i_4 i_5)}
+
$$
$$
+\sum_{j_1=0}^{p_1}\ldots\sum_{j_5=0}^{p_5}
C_{j_5 j_4 j_3 j_2 j_1}\Biggl(
{\bf 1}_{\{i_1=i_2\ne 0\}}
{\bf 1}_{\{j_1=j_2\}}
J'[\phi_{j_3}\phi_{j_4}
\phi_{j_5}]^{(i_3i_4i_5)}_{T,t}+
$$
$$
+
{\bf 1}_{\{i_1=i_3\ne 0\}}
{\bf 1}_{\{j_1=j_3\}}
J'[\phi_{j_2}\phi_{j_4}\phi_{j_5}]^{(i_2i_4i_5)}_{T,t}+
{\bf 1}_{\{i_1=i_4\ne 0\}}
{\bf 1}_{\{j_1=j_4\}}
J'[\phi_{j_2}\phi_{j_3}\phi_{j_5}]^{(i_2i_3i_5)}_{T,t}
+
$$
$$
+
{\bf 1}_{\{i_1=i_5\ne 0\}}
{\bf 1}_{\{j_1=j_5\}}
J'[\phi_{j_2}\phi_{j_3}\phi_{j_4}]^{(i_2i_3i_4)}_{T,t}
+
{\bf 1}_{\{i_2=i_3\ne 0\}}
{\bf 1}_{\{j_2=j_3\}}
J'[\phi_{j_1}\phi_{j_4}\phi_{j_5}]^{(i_1i_4i_5)}_{T,t}
+
$$
$$
+
{\bf 1}_{\{i_2=i_4\ne 0\}}
{\bf 1}_{\{j_2=j_4\}}
J'[\phi_{j_1}\phi_{j_3}\phi_{j_5}]^{(i_1i_3i_5)}_{T,t}
+
{\bf 1}_{\{i_2=i_5\ne 0\}}
{\bf 1}_{\{j_2=j_5\}}
J'[\phi_{j_1}\phi_{j_3}\phi_{j_4}]^{(i_1i_3i_4)}_{T,t}
+
$$
$$
+
{\bf 1}_{\{i_3=i_4\ne 0\}}
{\bf 1}_{\{j_3=j_4\}}
J'[\phi_{j_1}\phi_{j_2}\phi_{j_5}]^{(i_1i_2i_5)}_{T,t}
+
{\bf 1}_{\{i_3=i_5\ne 0\}}
{\bf 1}_{\{j_3=j_5\}}
J'[\phi_{j_1}\phi_{j_2}\phi_{j_4}]^{(i_1i_2i_4)}_{T,t}
+
$$
$$
+{\bf 1}_{\{i_4=i_5\ne 0\}}
{\bf 1}_{\{j_4=j_5\}}
J'[\phi_{j_1}\phi_{j_2}\phi_{j_3}]^{(i_1i_2i_3)}_{T,t}
+
$$
$$
+
{\bf 1}_{\{i_1=i_2\ne 0\}}
{\bf 1}_{\{j_1=j_2\}}
{\bf 1}_{\{i_3=i_4\ne 0\}}
{\bf 1}_{\{j_3=j_4\}}J'[\phi_{j_5}]^{(i_5)}_{T,t}+
$$
$$
+
{\bf 1}_{\{i_1=i_2\ne 0\}}
{\bf 1}_{\{j_1=j_2\}}
{\bf 1}_{\{i_3=i_5\ne 0\}}
{\bf 1}_{\{j_3=j_5\}}
J'[\phi_{j_4}]^{(i_4)}_{T,t}
+
$$
$$
+
{\bf 1}_{\{i_1=i_2\ne 0\}}
{\bf 1}_{\{j_1=j_2\}}
{\bf 1}_{\{i_4=i_5\ne 0\}}
{\bf 1}_{\{j_4=j_5\}}
J'[\phi_{j_3}]^{(i_3)}_{T,t}+
$$
$$
+
{\bf 1}_{\{i_1=i_3\ne 0\}}
{\bf 1}_{\{j_1=j_3\}}
{\bf 1}_{\{i_2=i_4\ne 0\}}
{\bf 1}_{\{j_2=j_4\}}
J'[\phi_{j_5}]^{(i_5)}_{T,t}
+
$$
$$
+
{\bf 1}_{\{i_1=i_3\ne 0\}}
{\bf 1}_{\{j_1=j_3\}}
{\bf 1}_{\{i_2=i_5\ne 0\}}
{\bf 1}_{\{j_2=j_5\}}
J'[\phi_{j_4}]^{(i_4)}_{T,t}+
$$
$$
+
{\bf 1}_{\{i_1=i_3\ne 0\}}
{\bf 1}_{\{j_1=j_3\}}
{\bf 1}_{\{i_4=i_5\ne 0\}}
{\bf 1}_{\{j_4=j_5\}}
J'[\phi_{j_2}]^{(i_2)}_{T,t}+
$$
$$
+
{\bf 1}_{\{i_1=i_4\ne 0\}}
{\bf 1}_{\{j_1=j_4\}}
{\bf 1}_{\{i_2=i_3\ne 0\}}
{\bf 1}_{\{j_2=j_3\}}
J'[\phi_{j_5}]^{(i_5)}_{T,t}
+
$$
$$
+
{\bf 1}_{\{i_1=i_4\ne 0\}}
{\bf 1}_{\{j_1=j_4\}}
{\bf 1}_{\{i_2=i_5\ne 0\}}
{\bf 1}_{\{j_2=j_5\}}
J'[\phi_{j_3}]^{(i_3)}_{T,t}
+
$$
$$
+
{\bf 1}_{\{i_1=i_4\ne 0\}}
{\bf 1}_{\{j_1=j_4\}}
{\bf 1}_{\{i_3=i_5\ne 0\}}
{\bf 1}_{\{j_3=j_5\}}
J'[\phi_{j_2}]^{(i_2)}_{T,t}
+
$$
$$
+
{\bf 1}_{\{i_1=i_5\ne 0\}}
{\bf 1}_{\{j_1=j_5\}}
{\bf 1}_{\{i_2=i_3\ne 0\}}
{\bf 1}_{\{j_2=j_3\}}
J'[\phi_{j_4}]^{(i_4)}_{T,t}
+
$$
$$
+
{\bf 1}_{\{i_1=i_5\ne 0\}}
{\bf 1}_{\{j_1=j_5\}}
{\bf 1}_{\{i_2=i_4\ne 0\}}
{\bf 1}_{\{j_2=j_4\}}
J'[\phi_{j_3}]^{(i_3)}_{T,t}+
$$
$$
+
{\bf 1}_{\{i_1=i_5\ne 0\}}
{\bf 1}_{\{j_1=j_5\}}
{\bf 1}_{\{i_3=i_4\ne 0\}}
{\bf 1}_{\{j_3=j_4\}}
J'[\phi_{j_2}]^{(i_2)}_{T,t}+
$$
$$
+
{\bf 1}_{\{i_2=i_3\ne 0\}}
{\bf 1}_{\{j_2=j_3\}}
{\bf 1}_{\{i_4=i_5\ne 0\}}
{\bf 1}_{\{j_4=j_5\}}
J'[\phi_{j_1}]^{(i_1)}_{T,t}+
$$
$$
+
{\bf 1}_{\{i_2=i_4\ne 0\}}
{\bf 1}_{\{j_2=j_4\}}
{\bf 1}_{\{i_3=i_5\ne 0\}}
{\bf 1}_{\{j_3=j_5\}}
J'[\phi_{j_1}]^{(i_1)}_{T,t}+
$$
\begin{equation}
\label{after35}
+\Biggl.
{\bf 1}_{\{i_2=i_5\ne 0\}}
{\bf 1}_{\{j_2=j_5\}}
{\bf 1}_{\{i_3=i_4\ne 0\}}
{\bf 1}_{\{j_3=j_4\}}
J'[\phi_{j_1}]^{(i_1)}_{T,t}
\Biggr),
\end{equation}

\vspace{3mm}

$$
\sum_{j_1=0}^{p_1}\ldots\sum_{j_6=0}^{p_6}
C_{j_6 j_5 j_4 j_3 j_2 j_1}
\zeta_{j_1}^{(i_1)}\zeta_{j_2}^{(i_2)}\zeta_{j_3}^{(i_3)}
\zeta_{j_4}^{(i_4)}\zeta_{j_5}^{(i_5)}\zeta_{j_6}^{(i_6)}
=
J'[K_{p_1p_2p_3p_4p_5p_6}]_{T,t}^{(i_1 i_2 i_3 i_4 i_5 i_6)}+
$$
$$
+\sum_{j_1=0}^{p_1}\ldots\sum_{j_6=0}^{p_6}
C_{j_6j_5j_4j_3j_2j_1}\Biggl(
{\bf 1}_{\{i_1=i_6\ne 0\}}
{\bf 1}_{\{j_1=j_6\}}
J'[\phi_{j_2}\phi_{j_3}\phi_{j_4}\phi_{j_5}]^{(i_2i_3i_4i_5)}_{T,t}
+
$$
$$
+\hspace{-0.4mm}
{\bf 1}_{\{i_2=i_6\ne 0\}}
{\bf 1}_{\{j_2=j_6\}}
J'[\phi_{j_1}\phi_{j_3}\phi_{j_4}\phi_{j_5}]^{(i_1i_3i_4i_5)}_{T,t}
\hspace{-2.1mm}+\hspace{-1.5mm}
{\bf 1}_{\{i_3=i_6\ne 0\}}
{\bf 1}_{\{j_3=j_6\}}
J'[\phi_{j_1}\phi_{j_2}\phi_{j_4}\phi_{j_5}]^{(i_1i_2i_4i_5)}_{T,t}\hspace{-1.0mm}
+
$$
$$
+\hspace{-0.4mm}
{\bf 1}_{\{i_4=i_6\ne 0\}}
{\bf 1}_{\{j_4=j_6\}}
J'[\phi_{j_1}\phi_{j_2}\phi_{j_3}\phi_{j_5}]^{(i_1i_2i_3i_5)}_{T,t}
\hspace{-2.1mm}+\hspace{-1.5mm}
{\bf 1}_{\{i_5=i_6\ne 0\}}
{\bf 1}_{\{j_5=j_6\}}
J'[\phi_{j_1}\phi_{j_2}\phi_{j_3}\phi_{j_4}]^{(i_1i_2i_3i_4)}_{T,t}\hspace{-1.0mm}
+
$$
$$
+\hspace{-0.4mm}{\bf 1}_{\{i_1=i_2\ne 0\}}
{\bf 1}_{\{j_1=j_2\}}
J'[\phi_{j_3}\phi_{j_4}\phi_{j_5}\phi_{j_6}]^{(i_3i_4i_5i_6)}_{T,t}
\hspace{-2.1mm}+\hspace{-1.5mm}
{\bf 1}_{\{i_1=i_3\ne 0\}}
{\bf 1}_{\{j_1=j_3\}}
J'[\phi_{j_2}\phi_{j_4}\phi_{j_5}\phi_{j_6}]^{(i_2i_4i_5i_6)}_{T,t}\hspace{-1.0mm}
+
$$
$$
+\hspace{-0.4mm}{\bf 1}_{\{i_1=i_4\ne 0\}}
{\bf 1}_{\{j_1=j_4\}}
J'[\phi_{j_2}\phi_{j_3}\phi_{j_5}\phi_{j_6}]^{(i_2i_3i_5i_6)}_{T,t}
\hspace{-2.1mm}+\hspace{-1.5mm}
{\bf 1}_{\{i_1=i_5\ne 0\}}
{\bf 1}_{\{j_1=j_5\}}
J'[\phi_{j_2}\phi_{j_3}\phi_{j_4}\phi_{j_6}]^{(i_2i_3i_4i_6)}_{T,t}\hspace{-1.0mm}
+
$$
$$
+\hspace{-0.4mm}{\bf 1}_{\{i_2=i_3\ne 0\}}
{\bf 1}_{\{j_2=j_3\}}
J'[\phi_{j_1}\phi_{j_4}\phi_{j_5}\phi_{j_6}]^{(i_1i_4i_5i_6)}_{T,t}
\hspace{-2.1mm}+\hspace{-1.5mm}
{\bf 1}_{\{i_2=i_4\ne 0\}}
{\bf 1}_{\{j_2=j_4\}}
J'[\phi_{j_1}\phi_{j_3}\phi_{j_5}\phi_{j_6}]^{(i_1i_3i_5i_6)}_{T,t}\hspace{-1.0mm}
+
$$
$$
+
\hspace{-0.4mm}{\bf 1}_{\{i_2=i_5\ne 0\}}
{\bf 1}_{\{j_2=j_5\}}
J'[\phi_{j_1}\phi_{j_3}\phi_{j_4}\phi_{j_6}]^{(i_1i_3i_4i_6)}_{T,t}
\hspace{-2.1mm}+\hspace{-1.5mm}
{\bf 1}_{\{i_3=i_4\ne 0\}}
{\bf 1}_{\{j_3=j_4\}}
J'[\phi_{j_1}\phi_{j_2}\phi_{j_5}\phi_{j_6}]^{(i_1i_2i_5i_6)}_{T,t}\hspace{-1.0mm}
+
$$
$$
+\hspace{-0.4mm}
{\bf 1}_{\{i_3=i_5\ne 0\}}
{\bf 1}_{\{j_3=j_5\}}
J'[\phi_{j_1}\phi_{j_2}\phi_{j_4}\phi_{j_6}]^{(i_1i_2i_4i_6)}_{T,t}
\hspace{-2.1mm}+\hspace{-1.5mm}
{\bf 1}_{\{i_4=i_5\ne 0\}}
{\bf 1}_{\{j_4=j_5\}}
J'[\phi_{j_1}\phi_{j_2}\phi_{j_3}\phi_{j_6}]^{(i_1i_2i_3i_6)}_{T,t}\hspace{-1.0mm}
+
$$
$$
+
{\bf 1}_{\{i_1=i_2\ne 0\}}
{\bf 1}_{\{j_1=j_2\}}
{\bf 1}_{\{i_3=i_4\ne 0\}}
{\bf 1}_{\{j_3=j_4\}}
J'[\phi_{j_5}\phi_{j_6}]^{(i_5i_6)}_{T,t}+
$$
$$
+
{\bf 1}_{\{i_1=i_2\ne 0\}}
{\bf 1}_{\{j_1=j_2\}}
{\bf 1}_{\{i_3=i_5\ne 0\}}
{\bf 1}_{\{j_3=j_5\}}
J'[\phi_{j_4}\phi_{j_6}]^{(i_4i_6)}_{T,t}
+
$$
$$
+
{\bf 1}_{\{i_1=i_2\ne 0\}}
{\bf 1}_{\{j_1=j_2\}}
{\bf 1}_{\{i_4=i_5\ne 0\}}
{\bf 1}_{\{j_4=j_5\}}
J'[\phi_{j_3}\phi_{j_6}]^{(i_3i_6)}_{T,t}
+
$$
$$
+
{\bf 1}_{\{i_1=i_3\ne 0\}}
{\bf 1}_{\{j_1=j_3\}}
{\bf 1}_{\{i_2=i_4\ne 0\}}
{\bf 1}_{\{j_2=j_4\}}
J'[\phi_{j_5}\phi_{j_6}]^{(i_5i_6)}_{T,t}
+
$$
$$
+
{\bf 1}_{\{i_1=i_3\ne 0\}}
{\bf 1}_{\{j_1=j_3\}}
{\bf 1}_{\{i_2=i_5\ne 0\}}
{\bf 1}_{\{j_2=j_5\}}
J'[\phi_{j_4}\phi_{j_6}]^{(i_4i_6)}_{T,t}
+
$$
$$
+{\bf 1}_{\{i_1=i_3\ne 0\}}
{\bf 1}_{\{j_1=j_3\}}
{\bf 1}_{\{i_4=i_5\ne 0\}}
{\bf 1}_{\{j_4=j_5\}}
J'[\phi_{j_2}\phi_{j_6}]^{(i_2i_6)}_{T,t}
+
$$
$$
+
{\bf 1}_{\{i_1=i_4\ne 0\}}
{\bf 1}_{\{j_1=j_4\}}
{\bf 1}_{\{i_2=i_3\ne 0\}}
{\bf 1}_{\{j_2=j_3\}}
J'[\phi_{j_5}\phi_{j_6}]^{(i_5i_6)}_{T,t}
+
$$
$$
+
{\bf 1}_{\{i_1=i_4\ne 0\}}
{\bf 1}_{\{j_1=j_4\}}
{\bf 1}_{\{i_2=i_5\ne 0\}}
{\bf 1}_{\{j_2=j_5\}}
J'[\phi_{j_3}\phi_{j_6}]^{(i_3i_6)}_{T,t}
+
$$
$$
+
{\bf 1}_{\{i_1=i_4\ne 0\}}
{\bf 1}_{\{j_1=j_4\}}
{\bf 1}_{\{i_3=i_5\ne 0\}}
{\bf 1}_{\{j_3=j_5\}}
J'[\phi_{j_2}\phi_{j_6}]^{(i_2i_6)}_{T,t}
+
$$
$$
+
{\bf 1}_{\{i_1=i_5\ne 0\}}
{\bf 1}_{\{j_1=j_5\}}
{\bf 1}_{\{i_2=i_3\ne 0\}}
{\bf 1}_{\{j_2=j_3\}}
J'[\phi_{j_4}\phi_{j_6}]^{(i_4i_6)}_{T,t}
+
$$
$$
+
{\bf 1}_{\{i_1=i_5\ne 0\}}
{\bf 1}_{\{j_1=j_5\}}
{\bf 1}_{\{i_2=i_4\ne 0\}}
{\bf 1}_{\{j_2=j_4\}}
J'[\phi_{j_3}\phi_{j_6}]^{(i_3i_6)}_{T,t}
+
$$
$$
+
{\bf 1}_{\{i_1=i_5\ne 0\}}
{\bf 1}_{\{j_1=j_5\}}
{\bf 1}_{\{i_3=i_4\ne 0\}}
{\bf 1}_{\{j_3=j_4\}}
J'[\phi_{j_2}\phi_{j_6}]^{(i_2i_6)}_{T,t}
+
$$
$$
+
{\bf 1}_{\{i_2=i_3\ne 0\}}
{\bf 1}_{\{j_2=j_3\}}
{\bf 1}_{\{i_4=i_5\ne 0\}}
{\bf 1}_{\{j_4=j_5\}}
J'[\phi_{j_1}\phi_{j_6}]^{(i_1i_6)}_{T,t}
+
$$
$$
+{\bf 1}_{\{i_2=i_4\ne 0\}}
{\bf 1}_{\{j_2=j_4\}}
{\bf 1}_{\{i_3=i_5\ne 0\}}
{\bf 1}_{\{j_3=j_5\}}
J'[\phi_{j_1}\phi_{j_6}]^{(i_1i_6)}_{T,t}
+
$$
$$
+
{\bf 1}_{\{i_2=i_5\ne 0\}}
{\bf 1}_{\{j_2=j_5\}}
{\bf 1}_{\{i_3=i_4\ne 0\}}
{\bf 1}_{\{j_3=j_4\}}
J'[\phi_{j_1}\phi_{j_6}]^{(i_1i_6)}_{T,t}
+
$$
$$
+{\bf 1}_{\{i_6=i_1\ne 0\}}
{\bf 1}_{\{j_6=j_1\}}
{\bf 1}_{\{i_3=i_4\ne 0\}}
{\bf 1}_{\{j_3=j_4\}}
J'[\phi_{j_2}\phi_{j_5}]^{(i_2i_5)}_{T,t}
+
$$
$$
+
{\bf 1}_{\{i_6=i_1\ne 0\}}
{\bf 1}_{\{j_6=j_1\}}
{\bf 1}_{\{i_3=i_5\ne 0\}}
{\bf 1}_{\{j_3=j_5\}}
J'[\phi_{j_2}\phi_{j_4}]^{(i_2i_4)}_{T,t}
+
$$
$$
+{\bf 1}_{\{i_6=i_1\ne 0\}}
{\bf 1}_{\{j_6=j_1\}}
{\bf 1}_{\{i_2=i_5\ne 0\}}
{\bf 1}_{\{j_2=j_5\}}
J'[\phi_{j_3}\phi_{j_4}]^{(i_3i_4)}_{T,t}
+
$$
$$
+
{\bf 1}_{\{i_6=i_1\ne 0\}}
{\bf 1}_{\{j_6=j_1\}}
{\bf 1}_{\{i_2=i_4\ne 0\}}
{\bf 1}_{\{j_2=j_4\}}
J'[\phi_{j_3}\phi_{j_5}]^{(i_3i_5)}_{T,t}
+
$$
$$
+{\bf 1}_{\{i_6=i_1\ne 0\}}
{\bf 1}_{\{j_6=j_1\}}
{\bf 1}_{\{i_4=i_5\ne 0\}}
{\bf 1}_{\{j_4=j_5\}}
J'[\phi_{j_2}\phi_{j_3}]^{(i_2i_3)}_{T,t}
+
$$
$$
+
{\bf 1}_{\{i_6=i_1\ne 0\}}
{\bf 1}_{\{j_6=j_1\}}
{\bf 1}_{\{i_2=i_3\ne 0\}}
{\bf 1}_{\{j_2=j_3\}}
J'[\phi_{j_4}\phi_{j_5}]^{(i_4i_5)}_{T,t}
+
$$
$$
+{\bf 1}_{\{i_6=i_2\ne 0\}}
{\bf 1}_{\{j_6=j_2\}}
{\bf 1}_{\{i_3=i_5\ne 0\}}
{\bf 1}_{\{j_3=j_5\}}
J'[\phi_{j_1}\phi_{j_4}]^{(i_1i_4)}_{T,t}
+
$$
$$
+
{\bf 1}_{\{i_6=i_2\ne 0\}}
{\bf 1}_{\{j_6=j_2\}}
{\bf 1}_{\{i_4=i_5\ne 0\}}
{\bf 1}_{\{j_4=j_5\}}
J'[\phi_{j_1}\phi_{j_3}]^{(i_1i_3)}_{T,t}
+
$$
$$
+{\bf 1}_{\{i_6=i_2\ne 0\}}
{\bf 1}_{\{j_6=j_2\}}
{\bf 1}_{\{i_3=i_4\ne 0\}}
{\bf 1}_{\{j_3=j_4\}}
J'[\phi_{j_1}\phi_{j_5}]^{(i_1i_5)}_{T,t}
+
$$
$$
+
{\bf 1}_{\{i_6=i_2\ne 0\}}
{\bf 1}_{\{j_6=j_2\}}
{\bf 1}_{\{i_1=i_5\ne 0\}}
{\bf 1}_{\{j_1=j_5\}}
J'[\phi_{j_3}\phi_{j_4}]^{(i_3i_4)}_{T,t}
+
$$
$$
+{\bf 1}_{\{i_6=i_2\ne 0\}}
{\bf 1}_{\{j_6=j_2\}}
{\bf 1}_{\{i_1=i_4\ne 0\}}
{\bf 1}_{\{j_1=j_4\}}
J'[\phi_{j_3}\phi_{j_5}]^{(i_3i_5)}_{T,t}
+
$$
$$
+
{\bf 1}_{\{i_6=i_2\ne 0\}}
{\bf 1}_{\{j_6=j_2\}}
{\bf 1}_{\{i_1=i_3\ne 0\}}
{\bf 1}_{\{j_1=j_3\}}
J'[\phi_{j_4}\phi_{j_5}]^{(i_4i_5)}_{T,t}
+
$$
$$
+{\bf 1}_{\{i_6=i_3\ne 0\}}
{\bf 1}_{\{j_6=j_3\}}
{\bf 1}_{\{i_2=i_5\ne 0\}}
{\bf 1}_{\{j_2=j_5\}}
J'[\phi_{j_1}\phi_{j_4}]^{(i_1i_4)}_{T,t}
+
$$
$$
+
{\bf 1}_{\{i_6=i_3\ne 0\}}
{\bf 1}_{\{j_6=j_3\}}
{\bf 1}_{\{i_4=i_5\ne 0\}}
{\bf 1}_{\{j_4=j_5\}}
J'[\phi_{j_1}\phi_{j_2}]^{(i_1i_2)}_{T,t}
+
$$
$$
+{\bf 1}_{\{i_6=i_3\ne 0\}}
{\bf 1}_{\{j_6=j_3\}}
{\bf 1}_{\{i_2=i_4\ne 0\}}
{\bf 1}_{\{j_2=j_4\}}
J'[\phi_{j_1}\phi_{j_5}]^{(i_1i_5)}_{T,t}
+
$$
$$
+
{\bf 1}_{\{i_6=i_3\ne 0\}}
{\bf 1}_{\{j_6=j_3\}}
{\bf 1}_{\{i_1=i_5\ne 0\}}
{\bf 1}_{\{j_1=j_5\}}
J'[\phi_{j_2}\phi_{j_4}]^{(i_2i_4)}_{T,t}
+
$$
$$
+{\bf 1}_{\{i_6=i_3\ne 0\}}
{\bf 1}_{\{j_6=j_3\}}
{\bf 1}_{\{i_1=i_4\ne 0\}}
{\bf 1}_{\{j_1=j_4\}}
J'[\phi_{j_2}\phi_{j_5}]^{(i_2i_5)}_{T,t}
+
$$
$$
+
{\bf 1}_{\{i_6=i_3\ne 0\}}
{\bf 1}_{\{j_6=j_3\}}
{\bf 1}_{\{i_1=i_2\ne 0\}}
{\bf 1}_{\{j_1=j_2\}}
J'[\phi_{j_4}\phi_{j_5}]^{(i_4i_5)}_{T,t}
+
$$
$$
+{\bf 1}_{\{i_6=i_4\ne 0\}}
{\bf 1}_{\{j_6=j_4\}}
{\bf 1}_{\{i_3=i_5\ne 0\}}
{\bf 1}_{\{j_3=j_5\}}
J'[\phi_{j_1}\phi_{j_2}]^{(i_1i_2)}_{T,t}
+
$$
$$
+
{\bf 1}_{\{i_6=i_4\ne 0\}}
{\bf 1}_{\{j_6=j_4\}}
{\bf 1}_{\{i_2=i_5\ne 0\}}
{\bf 1}_{\{j_2=j_5\}}
J'[\phi_{j_1}\phi_{j_3}]^{(i_1i_3)}_{T,t}
+
$$
$$
+{\bf 1}_{\{i_6=i_4\ne 0\}}
{\bf 1}_{\{j_6=j_4\}}
{\bf 1}_{\{i_2=i_3\ne 0\}}
{\bf 1}_{\{j_2=j_3\}}
J'[\phi_{j_1}\phi_{j_5}]^{(i_1i_5)}_{T,t}
+
$$
$$
+
{\bf 1}_{\{i_6=i_4\ne 0\}}
{\bf 1}_{\{j_6=j_4\}}
{\bf 1}_{\{i_1=i_5\ne 0\}}
{\bf 1}_{\{j_1=j_5\}}
J'[\phi_{j_2}\phi_{j_3}]^{(i_2i_3)}_{T,t}
+
$$
$$
+{\bf 1}_{\{i_6=i_4\ne 0\}}
{\bf 1}_{\{j_6=j_4\}}
{\bf 1}_{\{i_1=i_3\ne 0\}}
{\bf 1}_{\{j_1=j_3\}}
J'[\phi_{j_2}\phi_{j_5}]^{(i_2i_5)}_{T,t}
+
$$
$$
+
{\bf 1}_{\{i_6=i_4\ne 0\}}
{\bf 1}_{\{j_6=j_4\}}
{\bf 1}_{\{i_1=i_2\ne 0\}}
{\bf 1}_{\{j_1=j_2\}}
J'[\phi_{j_3}\phi_{j_5}]^{(i_3i_5)}_{T,t}
+
$$
$$
+{\bf 1}_{\{i_6=i_5\ne 0\}}
{\bf 1}_{\{j_6=j_5\}}
{\bf 1}_{\{i_3=i_4\ne 0\}}
{\bf 1}_{\{j_3=j_4\}}
J'[\phi_{j_1}\phi_{j_2}]^{(i_1i_2)}_{T,t}
+
$$
$$
+
{\bf 1}_{\{i_6=i_5\ne 0\}}
{\bf 1}_{\{j_6=j_5\}}
{\bf 1}_{\{i_2=i_4\ne 0\}}
{\bf 1}_{\{j_2=j_4\}}
J'[\phi_{j_1}\phi_{j_3}]^{(i_1i_3)}_{T,t}
+
$$
$$
+{\bf 1}_{\{i_6=i_5\ne 0\}}
{\bf 1}_{\{j_6=j_5\}}
{\bf 1}_{\{i_2=i_3\ne 0\}}
{\bf 1}_{\{j_2=j_3\}}
J'[\phi_{j_1}\phi_{j_4}]^{(i_1i_4)}_{T,t}
+
$$
$$
+
{\bf 1}_{\{i_6=i_5\ne 0\}}
{\bf 1}_{\{j_6=j_5\}}
{\bf 1}_{\{i_1=i_4\ne 0\}}
{\bf 1}_{\{j_1=j_4\}}
J'[\phi_{j_2}\phi_{j_3}]^{(i_2i_3)}_{T,t}
+
$$
$$
+{\bf 1}_{\{i_6=i_5\ne 0\}}
{\bf 1}_{\{j_6=j_5\}}
{\bf 1}_{\{i_1=i_3\ne 0\}}
{\bf 1}_{\{j_1=j_3\}}
J'[\phi_{j_2}\phi_{j_4}]^{(i_2i_4)}_{T,t}
+
$$
$$
+
{\bf 1}_{\{i_6=i_5\ne 0\}}
{\bf 1}_{\{j_6=j_5\}}
{\bf 1}_{\{i_1=i_2\ne 0\}}
{\bf 1}_{\{j_1=j_2\}}
J'[\phi_{j_3}\phi_{j_4}]^{(i_3i_4)}_{T,t}
+
$$
$$
+
{\bf 1}_{\{i_6=i_1\ne 0\}}
{\bf 1}_{\{j_6=j_1\}}
{\bf 1}_{\{i_2=i_5\ne 0\}}
{\bf 1}_{\{j_2=j_5\}}
{\bf 1}_{\{i_3=i_4\ne 0\}}
{\bf 1}_{\{j_3=j_4\}}+
$$
$$
+
{\bf 1}_{\{i_6=i_1\ne 0\}}
{\bf 1}_{\{j_6=j_1\}}
{\bf 1}_{\{i_2=i_4\ne 0\}}
{\bf 1}_{\{j_2=j_4\}}
{\bf 1}_{\{i_3=i_5\ne 0\}}
{\bf 1}_{\{j_3=j_5\}}+
$$
$$
+
{\bf 1}_{\{i_6=i_1\ne 0\}}
{\bf 1}_{\{j_6=j_1\}}
{\bf 1}_{\{i_2=i_3\ne 0\}}
{\bf 1}_{\{j_2=j_3\}}
{\bf 1}_{\{i_4=i_5\ne 0\}}
{\bf 1}_{\{j_4=j_5\}}+
$$
$$
+
{\bf 1}_{\{i_6=i_2\ne 0\}}
{\bf 1}_{\{j_6=j_2\}}
{\bf 1}_{\{i_1=i_5\ne 0\}}
{\bf 1}_{\{j_1=j_5\}}
{\bf 1}_{\{i_3=i_4\ne 0\}}
{\bf 1}_{\{j_3=j_4\}}+
$$
$$
+
{\bf 1}_{\{i_6=i_2\ne 0\}}
{\bf 1}_{\{j_6=j_2\}}
{\bf 1}_{\{i_1=i_4\ne 0\}}
{\bf 1}_{\{j_1=j_4\}}
{\bf 1}_{\{i_3=i_5\ne 0\}}
{\bf 1}_{\{j_3=j_5\}}+
$$
$$
+
{\bf 1}_{\{i_6=i_2\ne 0\}}
{\bf 1}_{\{j_6=j_2\}}
{\bf 1}_{\{i_1=i_3\ne 0\}}
{\bf 1}_{\{j_1=j_3\}}
{\bf 1}_{\{i_4=i_5\ne 0\}}
{\bf 1}_{\{j_4=j_5\}}+
$$
$$
+
{\bf 1}_{\{i_6=i_3\ne 0\}}
{\bf 1}_{\{j_6=j_3\}}
{\bf 1}_{\{i_1=i_5\ne 0\}}
{\bf 1}_{\{j_1=j_5\}}
{\bf 1}_{\{i_2=i_4\ne 0\}}
{\bf 1}_{\{j_2=j_4\}}+
$$
$$
+
{\bf 1}_{\{i_6=i_3\ne 0\}}
{\bf 1}_{\{j_6=j_3\}}
{\bf 1}_{\{i_1=i_4\ne 0\}}
{\bf 1}_{\{j_1=j_4\}}
{\bf 1}_{\{i_2=i_5\ne 0\}}
{\bf 1}_{\{j_2=j_5\}}+
$$
$$
+
{\bf 1}_{\{i_3=i_6\ne 0\}}
{\bf 1}_{\{j_3=j_6\}}
{\bf 1}_{\{i_1=i_2\ne 0\}}
{\bf 1}_{\{j_1=j_2\}}
{\bf 1}_{\{i_4=i_5\ne 0\}}
{\bf 1}_{\{j_4=j_5\}}+
$$
$$
+
{\bf 1}_{\{i_6=i_4\ne 0\}}
{\bf 1}_{\{j_6=j_4\}}
{\bf 1}_{\{i_1=i_5\ne 0\}}
{\bf 1}_{\{j_1=j_5\}}
{\bf 1}_{\{i_2=i_3\ne 0\}}
{\bf 1}_{\{j_2=j_3\}}+
$$
$$
+
{\bf 1}_{\{i_6=i_4\ne 0\}}
{\bf 1}_{\{j_6=j_4\}}
{\bf 1}_{\{i_1=i_3\ne 0\}}
{\bf 1}_{\{j_1=j_3\}}
{\bf 1}_{\{i_2=i_5\ne 0\}}
{\bf 1}_{\{j_2=j_5\}}+
$$
$$
+
{\bf 1}_{\{i_6=i_4\ne 0\}}
{\bf 1}_{\{j_6=j_4\}}
{\bf 1}_{\{i_1=i_2\ne 0\}}
{\bf 1}_{\{j_1=j_2\}}
{\bf 1}_{\{i_3=i_5\ne 0\}}
{\bf 1}_{\{j_3=j_5\}}+
$$
$$
+
{\bf 1}_{\{i_6=i_5\ne 0\}}
{\bf 1}_{\{j_6=j_5\}}
{\bf 1}_{\{i_1=i_4\ne 0\}}
{\bf 1}_{\{j_1=j_4\}}
{\bf 1}_{\{i_2=i_3\ne 0\}}
{\bf 1}_{\{j_2=j_3\}}+
$$
$$
+
{\bf 1}_{\{i_6=i_5\ne 0\}}
{\bf 1}_{\{j_6=j_5\}}
{\bf 1}_{\{i_1=i_2\ne 0\}}
{\bf 1}_{\{j_1=j_2\}}
{\bf 1}_{\{i_3=i_4\ne 0\}}
{\bf 1}_{\{j_3=j_4\}}+
$$
\begin{equation}
\label{after36}
\Biggl.
~~~~~~~~~+
{\bf 1}_{\{i_6=i_5\ne 0\}}
{\bf 1}_{\{j_6=j_5\}}
{\bf 1}_{\{i_1=i_3\ne 0\}}
{\bf 1}_{\{j_1=j_3\}}
{\bf 1}_{\{i_2=i_4\ne 0\}}
{\bf 1}_{\{j_2=j_4\}}\Biggr).
\end{equation}

\vspace{4mm}

Note that the relation (\ref{after34})
can be written in the following form

\vspace{-3mm}
$$
\sum_{j_1=0}^{p_1}\ldots \sum_{j_4=0}^{p_4}
C_{j_4 j_3 j_2 j_1}\zeta_{j_1}^{(i_1)}\zeta_{j_2}^{(i_2)}\zeta_{j_3}^{(i_3)}\zeta_{j_4}^{(i_4)}=
\sum_{j_1=0}^{p_1}\ldots \sum_{j_4=0}^{p_4}C_{j_4 j_3 j_2 j_1}
J'[\phi_{j_1}\phi_{j_2}\phi_{j_3}\phi_{j_4}]_{T,t}^{(i_1i_2i_3i_4)}+
$$
$$
+{\bf 1}_{\{i_1=i_2\ne 0\}}
\sum_{j_3=0}^{p_3}\sum_{j_4=0}^{p_4}\left(\sum_{j_1=0}^{\min\{p_1,p_2\}} C_{j_4 j_3 j_1 j_1}\right)
J'[\phi_{j_3}\phi_{j_4}]_{T,t}^{(i_3i_4)}+
$$
$$
+{\bf 1}_{\{i_1=i_3\ne 0\}}
\sum_{j_2=0}^{p_2}\sum_{j_4=0}^{p_4}\left(\sum_{j_3=0}^{\min\{p_1,p_3\}} C_{j_4 j_3 j_2 j_3}\right)
J'[\phi_{j_2}\phi_{j_4}]_{T,t}^{(i_2i_4)}+
$$
$$
+{\bf 1}_{\{i_1=i_4\ne 0\}}
\sum_{j_2=0}^{p_2}\sum_{j_3=0}^{p_3}\left(\sum_{j_4=0}^{\min\{p_1,p_4\}} C_{j_4 j_3 j_2 j_4}\right)
J'[\phi_{j_2}\phi_{j_3}]_{T,t}^{(i_2i_3)}+
$$
$$
+{\bf 1}_{\{i_2=i_3\ne 0\}}
\sum_{j_1=0}^{p_1}\sum_{j_4=0}^{p_4}\left(\sum_{j_3=0}^{\min\{p_2,p_3\}} C_{j_4 j_3 j_3 j_1}\right)
J'[\phi_{j_1}\phi_{j_4}]_{T,t}^{(i_1i_4)}+
$$
$$
+{\bf 1}_{\{i_2=i_4\ne 0\}}
\sum_{j_1=0}^{p_1}\sum_{j_3=0}^{p_3}\left(\sum_{j_4=0}^{\min\{p_2,p_4\}} C_{j_4 j_3 j_4 j_1}\right)
J'[\phi_{j_1}\phi_{j_3}]_{T,t}^{(i_1i_3)}+
$$
$$
+{\bf 1}_{\{i_3=i_4\ne 0\}}
\sum_{j_1=0}^{p_1}\sum_{j_2=0}^{p_2}\left(\sum_{j_4=0}^{\min\{p_3,p_4\}} C_{j_4 j_4 j_2 j_1}\right)
J'[\phi_{j_1}\phi_{j_2}]_{T,t}^{(i_1i_2)}+
$$
$$
+{\bf 1}_{\{i_2=i_3\ne 0\}}{\bf 1}_{\{i_1=i_4\ne 0\}}
\sum_{j_2=0}^{\min\{p_2,p_3\}}\sum_{j_4=0}^{\min\{p_1,p_4\}}C_{j_4 j_2 j_2 j_4}+
$$
$$
+{\bf 1}_{\{i_2=i_4\ne 0\}}{\bf 1}_{\{i_1=i_3\ne 0\}}
\sum_{j_3=0}^{\min\{p_1,p_3\}}\sum_{j_4=0}^{\min\{p_2,p_4\}}C_{j_4 j_3 j_4 j_3}+
$$
$$
+{\bf 1}_{\{i_3=i_4\ne 0\}}{\bf 1}_{\{i_1=i_2\ne 0\}}
\sum_{j_2=0}^{\min\{p_1,p_2\}}\sum_{j_4=0}^{\min\{p_3,p_4\}}C_{j_4 j_4 j_2 j_2}\ \ \ \hbox{w.~p.~1}.
$$

\vspace{3mm}

{\bf Step~2.}\ Let us prove that
\begin{equation}
\label{after80}
\sum_{j_l=0}^{\infty} C_{j_k \ldots j_{l+1} j_l j_{l-1} \ldots j_{s+1} j_l j_{s-1} \ldots j_1}=0
\end{equation}
or
\begin{equation}
\label{after80xx}
~~~~~~~\sum_{j_l=0}^{p} C_{j_k \ldots j_{l+1} j_l j_{l-1} \ldots j_{s+1} j_l j_{s-1} \ldots j_1}=
-\sum_{j_l=p+1}^{\infty} C_{j_k \ldots j_{l+1} j_l j_{l-1} \ldots j_{s+1} j_l j_{s-1} \ldots j_1},
\end{equation}

\noindent
where
$l-1\ge s+1.$

Our further proof will not fundamentally depend on the weight
functions $\psi_1(\tau),\ldots,\psi_k(\tau).$
Therefore, sometimes in subsequent consideration we assume for simplicity
that $\psi_1(\tau),\ldots,\psi_k(\tau)\equiv 1.$

Using the integration order replacement, we have 

\vspace{-3mm}
$$
C_{j_k \ldots j_{l+1} j_l j_{l-1} \ldots j_{s+1} j_l j_{s-1} \ldots j_1}=
$$

\vspace{-2mm}
$$
=\int\limits_t^T \phi_{j_k}(t_k)\ldots \int\limits_t^{t_{l+2}} \phi_{j_{l+1}}(t_{l+1})
\int\limits_t^{t_{l+1}} \phi_{j_{l}}(t_{l})
\int\limits_t^{t_{l}} \phi_{j_{l-1}}(t_{l-1})\ldots
$$
$$
\ldots
\int\limits_t^{t_{s+2}} \phi_{j_{s+1}}(t_{s+1})
\int\limits_t^{t_{s+1}} \phi_{j_{l}}(t_{s})
\int\limits_t^{t_{s}} \phi_{j_{s-1}}(t_{s-1})\ldots
$$
$$
\ldots \int\limits_t^{t_{2}} \phi_{j_{1}}(t_{1})dt_1\ldots dt_{s-1}dt_{s}dt_{s+1}\ldots
dt_{l-1}dt_{l}dt_{l+1}\ldots dt_k=
$$
$$
=\int\limits_t^{T} \phi_{j_{s+1}}(t_{s+1})
\int\limits_t^{t_{s+1}} \phi_{j_{l}}(t_{s})
\int\limits_t^{t_{s}} \phi_{j_{s-1}}(t_{s-1})\ldots
\int\limits_t^{t_{2}} \phi_{j_{1}}(t_{1})dt_1\ldots dt_{s-1}dt_{s}\times
$$
$$
\times \left(~
\int\limits_{t_{s+1}}^T \phi_{j_{s+2}}(t_{s+2})
\ldots \int\limits_{t_{l-2}}^T \phi_{j_{l-1}}(t_{l-1})
\int\limits_{t_{l-1}}^T \phi_{j_{l}}(t_{l})
\int\limits_{t_{l}}^T \phi_{j_{l+1}}(t_{l+1})\ldots \right.
$$
$$
\left.\ldots
\int\limits_{t_{k-1}}^T \phi_{j_k}(t_k)dt_k\ldots 
dt_{l+1}dt_{l}dt_{l-1}\ldots dt_{s+2}\right)dt_{s+1}=
$$
$$
=\int\limits_t^{T} \phi_{j_{s+1}}(t_{s+1})
\int\limits_t^{t_{s+1}} \phi_{j_{l}}(t_{s})
\underbrace{
\int\limits_t^{t_{s}} \phi_{j_{s-1}}(t_{s-1})\ldots
\int\limits_t^{t_{2}} \phi_{j_{1}}(t_{1})dt_1\ldots dt_{s-1}}
_{G_{j_{s-1}\ldots j_1}(t_s)}
dt_{s}\times
$$

\vspace{2mm}
$$
\times
\int\limits_{t_{s+1}}^T \phi_{j_{l}}(t_{l})
\underbrace{\int\limits_{t_{l}}^T \phi_{j_{l+1}}(t_{l+1})
\ldots
\int\limits_{t_{k-1}}^T \phi_{j_k}(t_k)dt_k\ldots 
dt_{l+1}}_{H_{j_{k}\ldots j_{l+1}}(t_l)}\times
$$

\vspace{-1mm}
$$
\times \left(~
\underbrace{\int\limits_{t_{s+1}}^{t_l} \phi_{j_{l-1}}(t_{l-1})
\ldots \int\limits_{t_{s+1}}^{t_{s+3}} \phi_{j_{s+2}}(t_{s+2})
dt_{s+2}\ldots dt_{l-1}}_{Q_{j_{l-1}\ldots j_{s+2}}(t_l,t_{s+1})}
dt_{l}\right)dt_{s+1}=
$$

\vspace{4mm}
$$
=\int\limits_t^{T} \phi_{j_{s+1}}(t_{s+1})
\int\limits_t^{t_{s+1}} \phi_{j_{l}}(t_{s})
G_{j_{s-1}\ldots j_1}(t_s)
dt_{s}\times
$$

\begin{equation}
\label{after9}
\times
\int\limits_{t_{s+1}}^T \phi_{j_{l}}(t_{l})
H_{j_{k}\ldots j_{l+1}}(t_l)
Q_{j_{l-1}\ldots j_{s+2}}(t_l,t_{s+1})
dt_{l}dt_{s+1}.
\end{equation}

\vspace{4mm}

Applying the additive property of the integral, we obtain

\vspace{-2mm}
$$
Q_{j_{l-1}\ldots j_{s+2}}(t_l,t_{s+1})=
$$

$$
=
\int\limits_{t_{s+1}}^{t_l} \phi_{j_{l-1}}(t_{l-1})
\ldots \int\limits_{t_{s+1}}^{t_{s+3}} \phi_{j_{s+2}}(t_{s+2})
dt_{s+2}\ldots dt_{l-1}=
$$

\vspace{-2mm}
$$
=
\int\limits_{t_{s+1}}^{t_l} \phi_{j_{l-1}}(t_{l-1})
\ldots \int\limits_{t_{s+1}}^{t_{s+4}}\phi_{j_{s+3}}(t_{s+3})
\int\limits_{t}^{t_{s+3}} \phi_{j_{s+2}}(t_{s+2})
dt_{s+2}dt_{s+3}\ldots dt_{l-1}-
$$
$$
-\int\limits_{t_{s+1}}^{t_l} \phi_{j_{l-1}}(t_{l-1})
\ldots \int\limits_{t_{s+1}}^{t_{s+4}}\phi_{j_{s+3}}(t_{s+3})
dt_{s+3}\ldots dt_{l-1}\int\limits_{t}^{t_{s+1}} \phi_{j_{s+2}}(t_{s+2})dt_{s+2}=
$$

\vspace{-5mm}
$$
\ldots 
$$

\vspace{-5mm}
\begin{equation}
\label{after10}
=\sum_{m=1}^d h^{(m)}_{j_{l-1}\ldots j_{s+2}}(t_l)
q^{(m)}_{j_{l-1}\ldots j_{s+2}}(t_{s+1}),
\end{equation}

\vspace{3mm}
\noindent
where $d<\infty.$

Combining (\ref{after9}) and (\ref{after10}), we have

\vspace{-2mm}
$$
\sum_{j_l=0}^p C_{j_k \ldots j_{l+1} j_l j_{l-1} \ldots j_{s+1} j_l j_{s-1} \ldots j_1}=
$$

\vspace{-4mm}
$$
=\sum_{m=1}^d \left(
\int\limits_t^{T} \phi_{j_{s+1}}(t_{s+1})q^{(m)}_{j_{l-1}\ldots j_{s+2}}(t_{s+1})
\sum_{j_l=0}^p \int\limits_t^{t_{s+1}} \phi_{j_{l}}(t_{s})
G_{j_{s-1}\ldots j_1}(t_s)
dt_{s}\times\right.
$$

\vspace{-2mm}
\begin{equation}
\label{after11}
\left.\times
\int\limits_{t_{s+1}}^T \phi_{j_{l}}(t_{l})
H_{j_{k}\ldots j_{l+1}}(t_l)
h^{(m)}_{j_{l-1}\ldots j_{s+2}}(t_l)
dt_{l}dt_{s+1}\right).
\end{equation}

\vspace{4mm}

Using the generalized Parseval equality, we obtain

\vspace{-3mm}
$$
\sum_{j_l=0}^{\infty} \int\limits_t^{t_{s+1}} \phi_{j_{l}}(t_{s})
G_{j_{s-1}\ldots j_1}(t_s)
dt_{s}
\int\limits_{t_{s+1}}^T \phi_{j_{l}}(t_{l})
H_{j_{k}\ldots j_{l+1}}(t_l)
h^{(m)}_{j_{l-1}\ldots j_{s+2}}(t_l)
dt_{l}=
$$

\vspace{-5mm}
\begin{equation}
\label{after400}
~~~~~=\int\limits_t^T {\bf 1}_{\{\tau<t_{s+1}\}}
G_{j_{s-1}\ldots j_1}(\tau) \cdot
{\bf 1}_{\{\tau>t_{s+1}\}}
H_{j_{k}\ldots j_{l+1}}(\tau)
h^{(m)}_{j_{l-1}\ldots j_{s+2}}(\tau)
d\tau=0.
\end{equation}

\vspace{4mm}

From (\ref{after11}) and (\ref{after400}) we get

\vspace{-2mm}
$$
\sum_{j_l=0}^p C_{j_k \ldots j_{l+1} j_l j_{l-1} \ldots j_{s+1} j_l j_{s-1} \ldots j_1}=
$$
$$
=-\sum_{m=1}^d \left(
\int\limits_t^{T} \phi_{j_{s+1}}(t_{s+1})q^{(m)}_{j_{l-1}\ldots j_{s+2}}(t_{s+1})
\sum_{j_l=p+1}^{\infty} \int\limits_t^{t_{s+1}} \phi_{j_{l}}(t_{s})
G_{j_{s-1}\ldots j_1}(t_s)
dt_{s}\times\right.
$$

\begin{equation}
\label{after450}
\left.\times
\int\limits_{t_{s+1}}^T \phi_{j_{l}}(t_{l})
H_{j_{k}\ldots j_{l+1}}(t_l)
h^{(m)}_{j_{l-1}\ldots j_{s+2}}(t_l)
dt_{l}dt_{s+1}\right).
\end{equation}

\vspace{4mm}

Combining Condition~2 of Theorem~2.30 and (\ref{after9})--(\ref{after11}), (\ref{after450}), we have

\vspace{-5.5mm}
$$
\sum_{j_l=0}^p C_{j_k \ldots j_{l+1} j_l j_{l-1} \ldots j_{s+1} j_l j_{s-1} \ldots j_1}=
$$

$$
=-\sum_{j_l=p+1}^{\infty}\sum_{m=1}^d \left(
\int\limits_t^{T} \phi_{j_{s+1}}(t_{s+1})q^{(m)}_{j_{l-1}\ldots j_{s+2}}(t_{s+1})
\int\limits_t^{t_{s+1}} \phi_{j_{l}}(t_{s})
G_{j_{s-1}\ldots j_1}(t_s)
dt_{s}\times\right.
$$

$$
\left.\times
\int\limits_{t_{s+1}}^T \phi_{j_{l}}(t_{l})
H_{j_{k}\ldots j_{l+1}}(t_l)
h^{(m)}_{j_{l-1}\ldots j_{s+2}}(t_l)
dt_{l}dt_{s+1}\right)=
$$

\vspace{1.5mm}
$$
=-\sum_{j_l=p+1}^{\infty}
\int\limits_t^T \phi_{j_k}(t_k)\ldots \int\limits_t^{t_{l+2}} \phi_{j_{l+1}}(t_{l+1})
\int\limits_t^{t_{l+1}} \phi_{j_{l}}(t_{l})
\int\limits_t^{t_{l}} \phi_{j_{l-1}}(t_{l-1})\ldots
$$

\vspace{3mm}
$$
\ldots
\int\limits_t^{t_{s+2}} \phi_{j_{s+1}}(t_{s+1})
\int\limits_t^{t_{s+1}} \phi_{j_{l}}(t_{s})
\int\limits_t^{t_{s}} \phi_{j_{s-1}}(t_{s-1})\ldots
$$

$$
\ldots \int\limits_t^{t_{2}} \phi_{j_{1}}(t_{1})dt_1\ldots dt_{s-1}dt_{s}dt_{s+1}\ldots
dt_{l-1}dt_{l}dt_{l+1}\ldots dt_k=
$$

\begin{equation}
\label{after79}
=-\sum_{j_l=p+1}^{\infty} C_{j_k \ldots j_{l+1} j_l j_{l-1} \ldots j_{s+1} j_l j_{s-1} \ldots j_1}.
\end{equation}

\vspace{6.5mm}

The equality (\ref{after79}) implies (\ref{after80}), (\ref{after80xx}).

{\bf Step 3.}\ Under the conditions of Theorem~2.30 we prove that

\vspace{-2mm}
$$
\sum_{j_l=0}^{p} C_{j_k \ldots j_{l+1} j_l j_l j_{l-2} \ldots j_1}=
$$
\begin{equation}
\label{after500}
=
\frac{1}{2}
C_{j_k \ldots j_1}\biggl|_{(j_l j_l)\curvearrowright (\cdot)}\biggr.
-\sum_{j_l=p+1}^{\infty} C_{j_k \ldots j_{l+1} j_l j_l j_{l-2} \ldots j_1}
\end{equation}

\vspace{1mm}
or
\begin{equation}
\label{after500xx}
\sum_{j_l=0}^{\infty} C_{j_k \ldots j_{l+1} j_l j_l j_{l-2} \ldots j_1}=
\frac{1}{2}
C_{j_k \ldots j_1}\biggl|_{(j_l j_l)\curvearrowright (\cdot)}\biggr..
\end{equation}

\vspace{4mm}

Denote
$$
C_{j_{l-2}\ldots j_1}(t_{l-1})=
\int\limits_t^{t_{l-1}} \psi_{l-2}(t_{l-2})\phi_{j_{l-2}}(t_{l-2})\ldots
\int\limits_t^{t_{2}} \psi_1(t_1)\phi_{j_{1}}(t_{1})dt_1\ldots dt_{l-2}.
$$

Using the integration order replacement and 
Condition~1 of Theorem~2.30, we obtain

\vspace{-4.5mm}
$$
\sum_{j_l=0}^{\infty} C_{j_k \ldots j_{l+1} j_l j_l j_{l-2} \ldots j_1}=
$$

\vspace{1mm}
$$
=
\sum_{j_l=0}^{\infty} \int\limits_t^T \psi_k(t_k)\phi_{j_k}(t_k)\ldots 
\int\limits_t^{t_{l+2}} 
\psi_{l+1}(t_{l+1})\phi_{j_{l+1}}(t_{l+1})\times
$$

\vspace{-3mm}
$$
\times
\int\limits_t^{t_{l+1}} 
\psi_l(t_l)\phi_{j_{l}}(t_{l})
\int\limits_t^{t_{l}} \psi_{l-1}(t_{l-1}) \phi_{j_{l}}(t_{l-1})
C_{j_{l-2}\ldots j_1}(t_{l-1})
dt_{l-1}
dt_{l}dt_{l+1}\ldots dt_k=
$$

\vspace{-2mm}
$$
=
\sum_{j_l=0}^{\infty} \int\limits_t^T 
\psi_l(t_l)\phi_{j_{l}}(t_{l})
\int\limits_t^{t_{l}} \psi_{l-1}(t_{l-1}) \phi_{j_{l}}(t_{l-1})
C_{j_{l-2}\ldots j_1}(t_{l-1})
dt_{l-1}\times
$$

$$
\times
\int\limits_{t_l}^T
\psi_{l+1}(t_{l+1})\phi_{j_{l+1}}(t_{l+1})\ldots
\int\limits_{t_{k-1}}^T
\psi_{k}(t_{k})\phi_{j_{k}}(t_{k})dt_k\ldots dt_{l+1}dt_l=
$$
$$
=\frac{1}{2}
\sum_{j_l=0}^{\infty} \int\limits_t^T 
\psi_l(t_l)\psi_{l-1}(t_{l})
C_{j_{l-2}\ldots j_1}(t_{l})
\times
$$

\vspace{-2mm}
$$
\times
\int\limits_{t_l}^T
\psi_{l+1}(t_{l+1})\phi_{j_{l+1}}(t_{l+1})\ldots
\int\limits_{t_{k-1}}^T
\psi_{k}(t_{k})\phi_{j_{k}}(t_{k})dt_k\ldots dt_{l+1}dt_l=
$$

\vspace{-2mm}
$$
=
\frac{1}{2}\sum_{j_l=0}^{\infty} \int\limits_t^T \psi_k(t_k)\phi_{j_k}(t_k)\ldots 
\int\limits_t^{t_{l+2}} 
\psi_{l+1}(t_{l+1})\phi_{j_{l+1}}(t_{l+1})\times
$$

\vspace{-2mm}
$$
\times
\int\limits_t^{t_{l+1}} 
\psi_l(t_l)\psi_{l-1}(t_{l})
C_{j_{l-2}\ldots j_1}(t_{l})
dt_{l}dt_{l+1}\ldots dt_k=
$$

\begin{equation}
\label{r12345x}
=\frac{1}{2}
C_{j_k \ldots j_1}\biggl|_{(j_l j_l)\curvearrowright (\cdot)}\biggr..
\end{equation}

\vspace{5mm}

The equalities (\ref{after500}) and (\ref{after500xx}) are proved.

\vspace{4mm}

{\bf Step~4.}\ Passing to the limit 
$\hbox{\vtop{\offinterlineskip\halign{
\hfil#\hfil\cr
{\rm l.i.m.}\cr
$\stackrel{}{{}_{p\to \infty}}$\cr
}} }$ 
in (\ref{after8xx}), we have (see (\ref{drdr1}))

$$
\hbox{\vtop{\offinterlineskip\halign{
\hfil#\hfil\cr
{\rm l.i.m.}\cr
$\stackrel{}{{}_{p\to \infty}}$\cr
}} }\sum_{j_1,\ldots,j_k=0}^{p}
C_{j_k\ldots j_1}
\zeta_{j_1}^{(i_1)}\ldots \zeta_{j_k}^{(i_k)}
=J[\psi^{(k)}]_{T,t}^{(i_1\ldots i_k)}
+
$$

\vspace{2mm}
$$
+
\sum\limits_{r=1}^{[k/2]}
\sum_{\stackrel{(\{\{g_1, g_2\}, \ldots, 
\{g_{2r-1}, g_{2r}\}\}, \{q_1, \ldots, q_{k-2r}\})}
{{}_{\{g_1, g_2, \ldots, 
g_{2r-1}, g_{2r}, q_1, \ldots, q_{k-2r}\}=\{1, 2, \ldots, k\}}}}
\prod\limits_{s=1}^r
{\bf 1}_{\{i_{g_{{}_{2s-1}}}=~i_{g_{{}_{2s}}}\ne 0\}}\times
$$

\vspace{4mm}
\begin{equation}
\label{after501}
~~~~~~~~\times \hbox{\vtop{\offinterlineskip\halign{
\hfil#\hfil\cr
{\rm l.i.m.}\cr
$\stackrel{}{{}_{p\to \infty}}$\cr
}} }\sum_{j_1,\ldots,j_k=0}^{p}
C_{j_k\ldots j_1}
\prod\limits_{s=1}^r{\bf 1}_{\{j_{g_{{}_{2s-1}}}=~j_{g_{{}_{2s}}}\}}
J'[\phi_{j_{q_1}}\ldots \phi_{j_{q_{k-2r}}}]_{T,t}^{(i_{q_1}\ldots i_{q_{k-2r}})}
\end{equation}

\vspace{3mm}
\noindent
w.~p.~1.

Taking into account (\ref{after80xx}) and (\ref{after500}), we obtain for $r=1$

\vspace{-6mm}
$$
{\bf 1}_{\{i_{g_{{}_{1}}}=~i_{g_{{}_{2}}}\ne 0\}}\hbox{\vtop{\offinterlineskip\halign{
\hfil#\hfil\cr
{\rm l.i.m.}\cr
$\stackrel{}{{}_{p\to \infty}}$\cr
}} }\sum_{j_1,\ldots,j_k=0}^{p}
C_{j_k\ldots j_1}
{\bf 1}_{\{j_{g_{{}_{1}}}=~j_{g_{{}_{2}}}\}}
J'[\phi_{j_{q_1}}\ldots \phi_{j_{q_{k-2}}}]_{T,t}^{(i_{q_1}\ldots i_{q_{k-2}})}=
$$
$$
=-{\bf 1}_{\{i_{g_{{}_{1}}}=~i_{g_{{}_{2}}}\ne 0\}}\hbox{\vtop{\offinterlineskip\halign{
\hfil#\hfil\cr
{\rm l.i.m.}\cr
$\stackrel{}{{}_{p\to \infty}}$\cr
}} }\sum\limits_{j_{g_1}=p+1}^{\infty}
\sum\limits_{\stackrel{j_1,\ldots,j_q,\ldots,j_k=0}{{}_{q\ne g_1, g_2}}}^p
C_{j_k\ldots j_1}\biggl|_{j_{g_{{}_{1}}}=~j_{g_{{}_{2}}}}\biggr. 
{\bf 1}_{\{g_2>g_1+1\}}\times
$$

$$
\times
J'[\phi_{j_{q_1}}\ldots \phi_{j_{q_{k-2}}}]_{T,t}^{(i_{q_1}\ldots i_{q_{k-2}})}+
$$

\vspace{-3mm}
$$
+
{\bf 1}_{\{i_{g_{{}_{1}}}=~i_{g_{{}_{2}}}\ne 0\}}\hbox{\vtop{\offinterlineskip\halign{
\hfil#\hfil\cr
{\rm l.i.m.}\cr
$\stackrel{}{{}_{p\to \infty}}$\cr
}} }\sum\limits_{\stackrel{j_1,\ldots,j_q,\ldots,j_k=0}{{}_{q\ne g_1, g_2}}}^p
\frac{1}{2}
C_{j_k \ldots j_1}\biggl|_{(j_{g_2} j_{g_1})\curvearrowright (\cdot),
j_{g_{{}_{1}}}=~j_{g_{{}_{2}}}}\biggr.
{\bf 1}_{\{g_2=g_1+1\}}\times
$$

$$
\times
J'[\phi_{j_{q_1}}\ldots \phi_{j_{q_{k-2}}}]_{T,t}^{(i_{q_1}\ldots i_{q_{k-2}})}-
$$

\vspace{-1mm}
$$
-{\bf 1}_{\{i_{g_{{}_{1}}}=~i_{g_{{}_{2}}}\ne 0\}}\hbox{\vtop{\offinterlineskip\halign{
\hfil#\hfil\cr
{\rm l.i.m.}\cr
$\stackrel{}{{}_{p\to \infty}}$\cr
}} }\sum\limits_{j_{g_1}=p+1}^{\infty}
\sum\limits_{\stackrel{j_1,\ldots,j_q,\ldots,j_k=0}{{}_{q\ne g_1, g_2}}}^p
C_{j_k \ldots j_1}\biggl|_{j_{g_{{}_{1}}}=~j_{g_{{}_{2}}}}\biggr. 
{\bf 1}_{\{g_2=g_1+1\}}\times
$$

$$
\times
J'[\phi_{j_{q_1}}\ldots \phi_{j_{q_{k-2}}}]_{T,t}^{(i_{q_1}\ldots i_{q_{k-2}})}=
$$

$$
=-{\bf 1}_{\{i_{g_{{}_{1}}}=~i_{g_{{}_{2}}}\ne 0\}}\hbox{\vtop{\offinterlineskip\halign{
\hfil#\hfil\cr
{\rm l.i.m.}\cr
$\stackrel{}{{}_{p\to \infty}}$\cr
}} }\sum\limits_{j_{g_1}=p+1}^{\infty}
\sum\limits_{\stackrel{j_1,\ldots,j_q,\ldots,j_k=0}{{}_{q\ne g_1, g_2}}}^p
C_{j_k \ldots j_1}\biggl|_{j_{g_{{}_{1}}}=~j_{g_{{}_{2}}}}\biggr. 
\times
$$

\vspace{1mm}
$$
\times
J'[\phi_{j_{q_1}}\ldots \phi_{j_{q_{k-2}}}]_{T,t}^{(i_{q_1}\ldots i_{q_{k-2}})}+
$$

\vspace{-3mm}
$$
+
{\bf 1}_{\{i_{g_{{}_{1}}}=~i_{g_{{}_{2}}}\ne 0\}}
\hbox{\vtop{\offinterlineskip\halign{
\hfil#\hfil\cr
{\rm l.i.m.}\cr
$\stackrel{}{{}_{p\to \infty}}$\cr
}} }\sum\limits_{\stackrel{j_1,\ldots,j_q,\ldots,j_k=0}{{}_{q\ne g_1, g_2}}}^p
\frac{1}{2}
C_{j_k \ldots j_1}\biggl|_{(j_{g_2} j_{g_1})\curvearrowright (\cdot),
j_{g_{{}_{1}}}=~j_{g_{{}_{2}}}}\biggr.
{\bf 1}_{\{g_2=g_1+1\}}\times
$$

\begin{equation}
\label{after600}
\times
J'[\phi_{j_{q_1}}\ldots \phi_{j_{q_{k-2}}}]_{T,t}^{(i_{q_1}\ldots i_{q_{k-2}})}=
\end{equation}

\begin{equation}
\label{after607}
~~~~~~~~=\frac{1}{2}{\bf 1}_{\{g_2=g_1+1\}}
J[\psi^{(k)}]_{T,t}^{g_1}+{\bf 1}_{\{i_{g_{{}_{1}}}=~i_{g_{{}_{2}}}\ne 0\}}
\hbox{\vtop{\offinterlineskip\halign{
\hfil#\hfil\cr
{\rm l.i.m.}\cr
$\stackrel{}{{}_{p\to \infty}}$\cr
}} } R_{T,t}^{(p)1,g_1,g_2}\ \ \ \hbox{w.~p.~1},
\end{equation}

\vspace{6mm}
\noindent 
where $J[\psi^{(k)}]_{T,t}^{g_1}$ $(g_1=1,2,\ldots,k-1)$ is
defined by (\ref{30.1}), 
$$
R_{T,t}^{(p)1,g_1,g_2}=
-\sum\limits_{\stackrel{j_1,\ldots,j_q,\ldots,j_k=0}{{}_{q\ne g_1, g_2}}}^p
\bar C^{(p)}_{j_k\ldots j_q \ldots j_1}\biggl|_{q\ne g_1,g_2}
J'[\phi_{j_{q_1}}\ldots \phi_{j_{q_{k-2}}}]_{T,t}^{(i_{q_1}\ldots i_{q_{k-2}})}.
$$

\vspace{1mm}

Let us explain the transition from 
(\ref{after600}) to (\ref{after607}). We have for $g_2=g_1+1$
$$
{\bf 1}_{\{i_{g_{{}_{1}}}=~i_{g_{{}_{2}}}\ne 0\}}
\hbox{\vtop{\offinterlineskip\halign{
\hfil#\hfil\cr
{\rm l.i.m.}\cr
$\stackrel{}{{}_{p\to \infty}}$\cr
}} }\sum\limits_{\stackrel{j_1,\ldots,j_q,\ldots,j_k=0}{{}_{q\ne g_1, g_2}}}^p
\frac{1}{2}
C_{j_k \ldots j_1}\biggl|_{(j_{g_2} j_{g_1})\curvearrowright (\cdot),
j_{g_{{}_{1}}}=~j_{g_{{}_{2}}}}\biggr.
\times
$$

$$
\times
J'[\phi_{j_{q_1}}\ldots \phi_{j_{q_{k-2}}}]_{T,t}^{(i_{q_1}\ldots i_{q_{k-2}})}=
$$

\vspace{1mm}
$$
=\frac{1}{2}{\bf 1}_{\{i_{g_{{}_{1}}}=~i_{g_{{}_{2}}}\ne 0\}}
\hbox{\vtop{\offinterlineskip\halign{
\hfil#\hfil\cr
{\rm l.i.m.}\cr
$\stackrel{}{{}_{p\to \infty}}$\cr
}} }\sum\limits_{\stackrel{j_1,\ldots,j_q,\ldots,j_k=0}{{}_{q\ne g_1, g_2}}}^p
C_{j_k \ldots j_1}\biggl|_{(j_{g_2} j_{g_1})\curvearrowright 0,
j_{g_{{}_{1}}}=~j_{g_{{}_{2}}}}\biggr.
\times
$$

\vspace{1mm}
$$
\times
\zeta_{0}^{(0)} J'[\phi_{j_{q_1}}\ldots \phi_{j_{q_{k-2}}}]_{T,t}^{(i_{q_1}\ldots i_{q_{k-2}})}=
$$

\vspace{2mm}
$$
=\frac{1}{2}{\bf 1}_{\{i_{g_{{}_{1}}}=~i_{g_{{}_{2}}}\ne 0\}}
\hbox{\vtop{\offinterlineskip\halign{
\hfil#\hfil\cr
{\rm l.i.m.}\cr
$\stackrel{}{{}_{p\to \infty}}$\cr
}} }\sum\limits_{\stackrel{j_1,\ldots,j_q,\ldots,j_k=0}{{}_{q\ne g_1, g_2}}}^p
\sum_{j_{m_1}=0}^p
C_{j_k \ldots j_1}\biggl|_{(j_{g_2} j_{g_1})\curvearrowright  j_{m_1},
j_{g_{{}_{1}}}=~j_{g_{{}_{2}}}}\biggr.
\times
$$

\vspace{1mm}
$$
\times
\zeta_{j_{m_1}}^{(0)} 
J'[\phi_{j_{q_1}}\ldots \phi_{j_{q_{k-2}}}]_{T,t}^{(i_{q_1}\ldots i_{q_{k-2}})}=
$$

\vspace{2mm}
$$
=\frac{1}{2}{\bf 1}_{\{i_{g_{{}_{1}}}=~i_{g_{{}_{2}}}\ne 0\}}
\hbox{\vtop{\offinterlineskip\halign{
\hfil#\hfil\cr
{\rm l.i.m.}\cr
$\stackrel{}{{}_{p\to \infty}}$\cr
}} }\sum\limits_{\stackrel{j_1,\ldots,j_q,\ldots,j_k=0}{{}_{q\ne g_1, g_2}}}^p
\sum_{j_{m_1}=0}^p
C_{j_k \ldots j_1}\biggl|_{(j_{g_2} j_{g_1})\curvearrowright  j_{m_1},
j_{g_{{}_{1}}}=~j_{g_{{}_{2}}}}\biggr.
\times
$$

\vspace{1mm}
\begin{equation}
\label{after608}
\times
J'[\phi_{ j_{m_1}} \phi_{j_{q_1}}\ldots 
\phi_{j_{q_{k-2}}}]_{T,t}^{(0 i_{q_1}\ldots i_{q_{k-2}})}=
\end{equation}

\vspace{4mm}
\begin{equation}
\label{after609}
=
\frac{1}{2}J[\psi^{(k)}]_{T,t}^{g_1}\ \ \ \hbox{w.~p.~1},
\end{equation}

\vspace{4mm}
\noindent
where 
$$
C_{j_k \ldots j_1}\biggl|_{(j_{g_2} j_{g_1})\curvearrowright j_{m_1},
j_{g_{{}_{1}}}=~j_{g_{{}_{2}}}, g_2=g_1+1}\biggr.
=
$$
$$
=
\int\limits_t^T \psi_k(t_k)\phi_{j_k}(t_k)\ldots 
\int\limits_t^{t_{g_1+3}} 
\psi_l(t_{g_1+2})\phi_{j_{g_1+2}}(t_{g_1+2})
\int\limits_t^{t_{g_1+2}} 
\psi_{g_1+1}(t_{g_1})\psi_{g_1}(t_{g_1})\phi_{j_{m_1}}(t_{g_1})\times
$$
$$
\times
\int\limits_t^{t_{g_1}} \psi_l(t_{g_1-1})\phi_{j_{g_1-1}}(t_{g_1-1})\ldots
\int\limits_t^{t_{2}} \psi_1(t_1)\phi_{j_{1}}(t_{1})dt_1\ldots dt_{g_1-1}dt_{g_1}
dt_{g_1+2}\ldots dt_k,
$$

\vspace{-2mm}
\begin{equation}
\label{dwdw21}
\zeta_{j_{m_1}}^{(0)}=\int\limits_t^T\phi_{j_{m_1}}(\tau)d{\bf w}_{\tau}^{(0)}
=\int\limits_t^T\phi_{j_{m_1}}(\tau)d\tau=
\left\{
\begin{matrix}
\sqrt{T-t} &\hbox{if}\ j_{m_1}=0\cr\cr
0 &\hbox{if}\ j_{m_1}\ne 0
\end{matrix},\right.
\end{equation}

\vspace{-2mm}
\begin{equation}
\label{dwdw22}
\phi_0(\tau)=\frac{1}{\sqrt{T-t}}.
\end{equation}

\vspace{2mm}
\noindent
The transition from (\ref{after608}) to (\ref{after609}) is based
on (\ref{drdr1}) or (\ref{febr5000}).

By Condition~3 of Theorem~2.30 we have (also see the property (\ref{wiener1}) of
multiple Wiener stochastic integral)
$$
\lim\limits_{p\to\infty}
{\sf M}\left\{\left(R_{T,t}^{(p)1,g_1,g_2}\right)^2\right\}
\le 
K \lim\limits_{p\to\infty}
\sum\limits_{\stackrel{j_1,\ldots,j_q,\ldots,j_k=0}{{}_{q\ne g_1, g_2}}}^p
\left(\bar C^{(p)}_{j_k\ldots j_q \ldots j_1}\biggl|_{q\ne g_1,g_2}\right)^2 =0,
$$
  
\noindent
where constant $K$ does not depend on $p$.

Thus

\vspace{-7mm}
$$
{\bf 1}_{\{i_{g_{{}_{1}}}=~i_{g_{{}_{2}}}\ne 0\}}\hbox{\vtop{\offinterlineskip\halign{
\hfil#\hfil\cr
{\rm l.i.m.}\cr
$\stackrel{}{{}_{p\to \infty}}$\cr
}} }\sum_{j_1,\ldots,j_k=0}^{p}
C_{j_k\ldots j_1}
{\bf 1}_{\{j_{g_{{}_{1}}}=~j_{g_{{}_{2}}}\}}
J'[\phi_{j_{q_1}}\ldots \phi_{j_{q_{k-2}}}]_{T,t}^{(i_{q_1}\ldots i_{q_{k-2}})}=
$$

$$
=\frac{1}{2}{\bf 1}_{\{g_2=g_1+1\}}
J[\psi^{(k)}]_{T,t}^{g_1}\ \ \ \hbox{w.~p.~1}.
$$

\vspace{4mm}

Involving into consideration the second pair $\{g_3,g_4\}$
(the first pair is $\{g_1,g_2\}$), we obtain
from (\ref{after600}) for $r=2$
$$
\prod\limits_{s=1}^2
{\bf 1}_{\{i_{g_{{}_{2s-1}}}=~i_{g_{{}_{2s}}}\ne 0\}}\hbox{\vtop{\offinterlineskip\halign{
\hfil#\hfil\cr
{\rm l.i.m.}\cr
$\stackrel{}{{}_{p\to \infty}}$\cr
}} }\sum_{j_1,\ldots,j_k=0}^{p}
C_{j_k\ldots j_1}
\prod\limits_{s=1}^2{\bf 1}_{\{j_{g_{{}_{2s-1}}}=~j_{g_{{}_{2s}}}\}}\times
$$

$$
\times
J'[\phi_{j_{q_1}}\ldots \phi_{j_{q_{k-4}}}]_{T,t}^{(i_{q_1}\ldots i_{q_{k-4}})}=
$$

\vspace{-2mm}
$$
=\prod\limits_{s=1}^2
{\bf 1}_{\{i_{g_{{}_{2s-1}}}=~i_{g_{{}_{2s}}}\ne 0\}}
\times
$$
$$
\times\hbox{\vtop{\offinterlineskip\halign{
\hfil#\hfil\cr
{\rm l.i.m.}\cr
$\stackrel{}{{}_{p\to \infty}}$\cr
}} }
\sum\limits_{\stackrel{j_1,\ldots,j_q,\ldots,j_k=0}{{}_{q\ne g_1, g_2, g_3, g_4}}}^p
\Biggl(\frac{1}{4}
C_{j_k \ldots j_1}\biggl|_{(j_{g_2} j_{g_1})\curvearrowright (\cdot)
(j_{g_4} j_{g_3})\curvearrowright (\cdot),
j_{g_{{}_{1}}}=~j_{g_{{}_{2}}}, j_{g_{{}_{3}}}=~j_{g_{{}_{4}}}}\biggr.
\prod\limits_{s=1}^2 {\bf 1}_{\{g_{2s}=g_{2s-1}+1\}}-\Biggr.
$$

$$
-\frac{1}{2}\sum\limits_{j_{g_1}=p+1}^{\infty}
C_{j_k \ldots j_1}\biggl|_{(j_{g_4} j_{g_3})\curvearrowright (\cdot),
j_{g_{{}_{1}}}=~j_{g_{{}_{2}}}, j_{g_{{}_{3}}}=~j_{g_{{}_{4}}}}
{\bf 1}_{\{g_{4}=g_{3}+1\}}-
$$

\vspace{-1mm}
$$
-
\frac{1}{2}\sum\limits_{j_{g_3}=p+1}^{\infty}
C_{j_k \ldots j_1}\biggl|_{(j_{g_2} j_{g_1})\curvearrowright (\cdot),
j_{g_{{}_{1}}}=~j_{g_{{}_{2}}}, j_{g_{{}_{3}}}=~j_{g_{{}_{4}}}}\biggr.
{\bf 1}_{\{g_{2}=g_{1}+1\}}+
$$

\vspace{-2mm}
\begin{equation}
\label{after610}
~~~~~~\Biggl.+
\sum\limits_{j_{g_3}=p+1}^{\infty}\sum\limits_{j_{g_1}=p+1}^{\infty}
C_{j_k \ldots j_1}\biggl|_{j_{g_{{}_{1}}}=~j_{g_{{}_{2}}}, j_{g_{{}_{3}}}=~j_{g_{{}_{4}}}} 
\biggr.
\Biggr)J'[\phi_{j_{q_1}}\ldots 
\phi_{j_{q_{k-4}}}]_{T,t}^{(i_{q_1}\ldots i_{q_{k-4}})}=
\end{equation}

\vspace{-1mm}
\begin{equation}
\label{after610x}
~~~=
\frac{1}{4}\prod\limits_{s=1}^2 {\bf 1}_{\{g_{2s}=g_{2s-1}+1\}}
J[\psi^{(k)}]_{T,t}^{s_2,s_1}+  
\prod\limits_{s=1}^2
{\bf 1}_{\{i_{g_{{}_{2s-1}}}=~i_{g_{{}_{2s}}}\ne 0\}}
\hbox{\vtop{\offinterlineskip\halign{
\hfil#\hfil\cr
{\rm l.i.m.}\cr
$\stackrel{}{{}_{p\to \infty}}$\cr
}} } R_{T,t}^{(p)2,g_1,g_2,g_3,g_4}
\end{equation}

\vspace{2mm}
\noindent
w.~p.~1, where $g_3\stackrel{\sf def}{=}s_2,$ $g_1\stackrel{\sf def}{=}s_1,$
$(s_2,s_1)\in {\rm A}_{k,2},$ $J[\psi^{(k)}]_{T,t}^{s_2,s_1}$ is
defined by (\ref{30.1}) and ${\rm A}_{k,2}$ is defined by (\ref{30.5550001}),

\vspace{-2mm}
$$
R_{T,t}^{(p)2,g_1,g_2,g_3,g_4}=
\sum\limits_{\stackrel{j_1,\ldots,j_q,\ldots,j_k=0}{{}_{q\ne g_1, g_2, g_3, g_4}}}^p
\left(
\bar C^{(p)}_{j_k\ldots j_q \ldots j_1}\biggl|_{q\ne g_1,g_2,g_3,g_4}-\right.
$$

\vspace{-2mm}
$$
-S_1\left\{
\bar C^{(p)}_{j_k\ldots j_q \ldots j_1}\biggl|_{q\ne g_1,g_2,g_3,g_4}\right\}
\left.-S_2\left\{
\bar C^{(p)}_{j_k\ldots j_q \ldots j_1}\biggl|_{q\ne g_1,g_2,g_3,g_4}\right\}\right)\times
$$

\vspace{2mm}
$$
\times
J'[\phi_{j_{q_1}}\ldots \phi_{j_{q_{k-4}}}]_{T,t}^{(i_{q_1}\ldots i_{q_{k-4}})}.
$$

\vspace{4mm}

Let us explain the transition from (\ref{after610}) to (\ref{after610x}).
We have for $g_2=g_1+1,$ $g_4=g_3+1$

\vspace{-4mm}
$$
\hbox{\vtop{\offinterlineskip\halign{
\hfil#\hfil\cr
{\rm l.i.m.}\cr
$\stackrel{}{{}_{p\to \infty}}$\cr
}} }
\sum\limits_{\stackrel{j_1,\ldots,j_q,\ldots,j_k=0}{{}_{q\ne g_1, g_2, g_3, g_4}}}^p
\frac{1}{4}
C_{j_k \ldots j_1}\biggl|_{(j_{g_2} j_{g_1})\curvearrowright (\cdot)
(j_{g_4} j_{g_3})\curvearrowright (\cdot),
j_{g_{{}_{1}}}=~j_{g_{{}_{2}}}, j_{g_{{}_{3}}}=~j_{g_{{}_{4}}}}\biggr.
\times
$$

$$
\times 
\prod\limits_{s=1}^2
{\bf 1}_{\{i_{g_{{}_{2s-1}}}=~i_{g_{{}_{2s}}}\ne 0\}}
J'[\phi_{j_{q_1}}\ldots \phi_{j_{q_{k-4}}}]_{T,t}^{(i_{q_1}\ldots i_{q_{k-4}})}=
$$
$$
=\frac{1}{4}
\hbox{\vtop{\offinterlineskip\halign{
\hfil#\hfil\cr
{\rm l.i.m.}\cr
$\stackrel{}{{}_{p\to \infty}}$\cr
}} }
\sum\limits_{\stackrel{j_1,\ldots,j_q,\ldots,j_k=0}{{}_{q\ne g_1, g_2, g_3, g_4}}}^p
C_{j_k \ldots j_1}\biggl|_{(j_{g_2} j_{g_1})\curvearrowright 0
(j_{g_4} j_{g_3})\curvearrowright 0,
j_{g_{{}_{1}}}=~j_{g_{{}_{2}}}, j_{g_{{}_{3}}}=~j_{g_{{}_{4}}}}\biggr.
\times
$$

\vspace{1mm}
$$
\times
\prod\limits_{s=1}^2
{\bf 1}_{\{i_{g_{{}_{2s-1}}}=~i_{g_{{}_{2s}}}\ne 0\}}
\zeta_{0}^{(0)}\zeta_{0}^{(0)}
J'[\phi_{j_{q_1}}\ldots \phi_{j_{q_{k-4}}}]_{T,t}^{(i_{q_1}\ldots i_{q_{k-4}})}=
$$

\vspace{1.5mm}
$$
=\frac{1}{4}
\hbox{\vtop{\offinterlineskip\halign{
\hfil#\hfil\cr
{\rm l.i.m.}\cr
$\stackrel{}{{}_{p\to \infty}}$\cr
}} }
\sum\limits_{\stackrel{j_1,\ldots,j_q,\ldots,j_k=0}{{}_{q\ne g_1, g_2, g_3, g_4}}}^p
\sum\limits_{ j_{m_1}, j_{m_3}=0}^p
C_{j_k \ldots j_1}\biggl|_{(j_{g_2} j_{g_1})\curvearrowright j_{m_1}
(j_{g_4} j_{g_3})\curvearrowright j_{m_3},
j_{g_{{}_{1}}}=~j_{g_{{}_{2}}}, j_{g_{{}_{3}}}=~j_{g_{{}_{4}}}}\biggr.
\times
$$

$$
\times
\prod\limits_{s=1}^2
{\bf 1}_{\{i_{g_{{}_{2s-1}}}=~i_{g_{{}_{2s}}}\ne 0\}}
\zeta_{ j_{m_1}}^{(0)}\zeta_{j_{m_3}}^{(0)}
J'[\phi_{j_{q_1}}\ldots \phi_{j_{q_{k-4}}}]_{T,t}^{(i_{q_1}\ldots i_{q_{k-4}})}=
$$

\vspace{1.5mm}
$$
=\frac{1}{4}
\hbox{\vtop{\offinterlineskip\halign{
\hfil#\hfil\cr
{\rm l.i.m.}\cr
$\stackrel{}{{}_{p\to \infty}}$\cr
}} }
\sum\limits_{\stackrel{j_1,\ldots,j_q,\ldots,j_k=0}{{}_{q\ne g_1, g_2, g_3, g_4}}}^p
\sum\limits_{j_{m_1}, j_{m_3}=0}^p
C_{j_k \ldots j_1}\biggl|_{(j_{g_2} j_{g_1})\curvearrowright j_{m_1}
(j_{g_4} j_{g_3})\curvearrowright j_{m_3},
j_{g_{{}_{1}}}=~j_{g_{{}_{2}}}, j_{g_{{}_{3}}}=~j_{g_{{}_{4}}}}\biggr.
\times
$$

\vspace{-3mm}
\begin{equation}
\label{after900x}
~~~~~~~~~~\times
\prod\limits_{s=1}^2
{\bf 1}_{\{i_{g_{{}_{2s-1}}}=~i_{g_{{}_{2s}}}\ne 0\}}
J'[\phi_{ j_{m_1}}\phi_{j_{m_3}}
\phi_{j_{q_1}}\ldots \phi_{j_{q_{k-4}}}]_{T,t}^{(0 0 i_{q_1}\ldots i_{q_{k-4}})}=
\end{equation}

\vspace{3mm}
\begin{equation}
\label{after901}
=
\frac{1}{4}
J[\psi^{(k)}]_{T,t}^{s_2,s_1}\ \ \ \hbox{w.~p.~1}.
\end{equation}

\vspace{7mm}
\noindent
The transition from (\ref{after900x}) to (\ref{after901}) is based
on (\ref{drdr1}) or (\ref{febr5000}).

Note that

\vspace{-2mm}
$$
C_{j_k \ldots j_1}\biggl|_{(j_{g_2} j_{g_1})\curvearrowright j_{m_1},
j_{g_{{}_{1}}}=~j_{g_{{}_{2}}}}\biggr.
=C_{j_k \ldots j_1}\biggl|_{(j_{g_1} j_{g_1})\curvearrowright j_{m_1},
j_{g_{{}_{1}}}=~j_{g_{{}_{2}}}}\biggr.
$$

\vspace{4mm}
\noindent
is the Fourier coefficient, where $g_{2}=g_{1}+1.$
Therefore, the value
$$
C_{j_k \ldots j_1}\biggl|_{(j_{g_2} j_{g_1})\curvearrowright j_{m_1}
(j_{g_4} j_{g_3})\curvearrowright  j_{m_3},
j_{g_{{}_{1}}}=~j_{g_{{}_{2}}},
j_{g_{{}_{3}}}=~j_{g_{{}_{4}}}
}\biggr.=
$$
$$
=C_{j_k \ldots j_1}\biggl|_{(j_{g_1} j_{g_1})\curvearrowright j_{m_1}
(j_{g_3} j_{g_3})\curvearrowright  j_{m_3},
j_{g_{{}_{1}}}=~j_{g_{{}_{2}}},
j_{g_{{}_{3}}}=~j_{g_{{}_{4}}}
}\biggr.
$$

\vspace{3mm}
\noindent
is determined recursively using (\ref{after2000}) in an obvious way
for $g_2=g_1+1$ and $g_4=g_3+1$.

By Condition~3 of Theorem~2.30 we have (also see the property (\ref{wiener1}) of
multiple Wiener stochastic integral)

\vspace{-1mm}
$$
\lim\limits_{p\to\infty}
{\sf M}\left\{\left(R_{T,t}^{(p)2,g_1,g_2,g_3,g_4}\right)^2\right\}
\le
$$
$$
\le
K \lim\limits_{p\to\infty}
\sum\limits_{\stackrel{j_1,\ldots,j_q,\ldots,j_k=0}{{}_{q\ne g_1, g_2, g_3, g_4}}}^p
\left(
\left(
\bar C^{(p)}_{j_k\ldots j_q \ldots j_1}\biggl|_{q\ne g_1,g_2,g_3,g_4}\biggr.\right)^2+
\right.
$$

\vspace{-3mm}
$$
\left.+\left(S_1\left\{
\bar C^{(p)}_{j_k\ldots j_q \ldots j_1}\biggl|_{q\ne g_1,g_2,g_3,g_4}\biggr.
\right\}\right)^2+
\left(S_2\left\{
\bar C^{(p)}_{j_k\ldots j_q \ldots j_1}\biggl|_{q\ne g_1,g_2,g_3,g_4}\biggr.
\right\}\right)^2\right)=0,
$$

\vspace{2mm}
\noindent
where constant $K$ is independent of $p$.

Thus

\vspace{-2mm}
$$
\prod\limits_{s=1}^2
{\bf 1}_{\{i_{g_{{}_{2s-1}}}=~i_{g_{{}_{2s}}}\ne 0\}}\hbox{\vtop{\offinterlineskip\halign{
\hfil#\hfil\cr
{\rm l.i.m.}\cr
$\stackrel{}{{}_{p\to \infty}}$\cr
}} }\sum_{j_1,\ldots,j_k=0}^{p}
C_{j_k\ldots j_1}
\prod\limits_{s=1}^2{\bf 1}_{\{j_{g_{{}_{2s-1}}}=~j_{g_{{}_{2s}}}\}}
\times
$$

\vspace{-2mm}
$$
\times J'[\phi_{j_{q_1}}\ldots \phi_{j_{q_{k-4}}}]_{T,t}^{(i_{q_1}\ldots i_{q_{k-4}})}=
\frac{1}{4}\prod\limits_{s=1}^2 {\bf 1}_{\{g_{2s}=g_{2s-1}+1\}}
J[\psi^{(k)}]_{T,t}^{s_2,s_1}\ \ \ \hbox{w.~p.~1},
$$

\vspace{2mm}
\noindent
where $g_3\stackrel{\sf def}{=}s_2,$ $g_1\stackrel{\sf def}{=}s_1,$
$(s_2,s_1)\in {\rm A}_{k,2},$ $J[\psi^{(k)}]_{T,t}^{s_2,s_1}$ is
defined by (\ref{30.1}) and ${\rm A}_{k,2}$ is defined by (\ref{30.5550001}).

Involving into consideration the third pair $\{g_6,g_5\}$
($\{g_1,g_2\}$ is the first pair and $\{g_4,g_3\}$ is the second pair), we obtain
from (\ref{after610}) for $r=3$

\vspace{-1mm}
$$
\prod\limits_{s=1}^3
{\bf 1}_{\{i_{g_{{}_{2s-1}}}=~i_{g_{{}_{2s}}}\ne 0\}}\hbox{\vtop{\offinterlineskip\halign{
\hfil#\hfil\cr
{\rm l.i.m.}\cr
$\stackrel{}{{}_{p\to \infty}}$\cr
}} }\sum_{j_1,\ldots,j_k=0}^{p}
C_{j_k\ldots j_1}
\prod\limits_{s=1}^3{\bf 1}_{\{j_{g_{{}_{2s-1}}}=~j_{g_{{}_{2s}}}\}}
\times
$$

\vspace{2mm}
$$
\times
J'[\phi_{j_{q_1}}\ldots \phi_{j_{q_{k-6}}}]_{T,t}^{(i_{q_1}\ldots i_{q_{k-6}})}=
\prod\limits_{s=1}^3
{\bf 1}_{\{i_{g_{{}_{2s-1}}}=~i_{g_{{}_{2s}}}\ne 0\}}\times
$$
$$
\times\hbox{\vtop{\offinterlineskip\halign{
\hfil#\hfil\cr
{\rm l.i.m.}\cr
$\stackrel{}{{}_{p\to \infty}}$\cr
}} }\hspace{-3mm}
\sum\limits_{\stackrel{j_1,\ldots,j_q,\ldots,j_k=0}{{}_{q\ne g_1, g_2, g_3, g_4, g_5, g_6}}}^p
\hspace{-2mm}
\left(\frac{1}{2^3}
C_{j_k \ldots j_1}\biggl|_{(j_{g_2} j_{g_1})\curvearrowright (\cdot)
(j_{g_4} j_{g_3})\curvearrowright (\cdot)
(j_{g_6} j_{g_5})\curvearrowright (\cdot),
j_{g_{{}_{1}}}=~j_{g_{{}_{2}}}, j_{g_{{}_{3}}}=~j_{g_{{}_{4}}}, j_{g_{{}_{5}}}=~j_{g_{{}_{6}}}
}\biggr.\right.
\times
$$
$$
\times
\prod\limits_{s=1}^3
{\bf 1}_{\{g_{2s}=g_{2s-1}+1\}}-
$$

\vspace{-6mm}
$$
-\frac{1}{2^2}\sum\limits_{j_{g_1}=p+1}^{\infty}
C_{j_k \ldots j_1}\biggl|_{(j_{g_4} j_{g_3})\curvearrowright (\cdot)
(j_{g_6} j_{g_5})\curvearrowright (\cdot),
j_{g_{{}_{1}}}=~j_{g_{{}_{2}}}, j_{g_{{}_{3}}}=~j_{g_{{}_{4}}}, j_{g_{{}_{5}}}=~j_{g_{{}_{6}}}}
{\bf 1}_{\{g_{4}=g_{3}+1\}}{\bf 1}_{\{g_{6}=g_{5}+1\}}-
$$

\vspace{-1mm}
$$
-\frac{1}{2^2}\sum\limits_{j_{g_3}=p+1}^{\infty}
C_{j_k \ldots j_1}\biggl|_{(j_{g_2} j_{g_1})\curvearrowright (\cdot)
(j_{g_6} j_{g_5})\curvearrowright (\cdot),
j_{g_{{}_{1}}}=~j_{g_{{}_{2}}}, j_{g_{{}_{3}}}=~j_{g_{{}_{4}}}, j_{g_{{}_{5}}}=~j_{g_{{}_{6}}}}
{\bf 1}_{\{g_{2}=g_{1}+1\}}{\bf 1}_{\{g_{6}=g_{5}+1\}}-
$$

\vspace{-1mm}
$$
-\frac{1}{2^2}\sum\limits_{j_{g_5}=p+1}^{\infty}
C_{j_k \ldots j_1}\biggl|_{(j_{g_2} j_{g_1})\curvearrowright (\cdot)
(j_{g_4} j_{g_3})\curvearrowright (\cdot),
j_{g_{{}_{1}}}=~j_{g_{{}_{2}}}, j_{g_{{}_{3}}}=~j_{g_{{}_{4}}}, j_{g_{{}_{5}}}=~j_{g_{{}_{6}}}}
{\bf 1}_{\{g_{2}=g_{1}+1\}}{\bf 1}_{\{g_{4}=g_{3}+1\}}+
$$

\vspace{-1mm}
$$
+\frac{1}{2}\sum\limits_{j_{g_3}=p+1}^{\infty}\sum\limits_{j_{g_1}=p+1}^{\infty}
C_{j_k \ldots j_1}\biggl|_{(j_{g_6} j_{g_5})\curvearrowright (\cdot),
j_{g_{{}_{1}}}=~j_{g_{{}_{2}}}, j_{g_{{}_{3}}}=~j_{g_{{}_{4}}}, j_{g_{{}_{5}}}=~j_{g_{{}_{6}}}}
{\bf 1}_{\{g_{6}=g_{5}+1\}}+
$$

\vspace{-1mm}
$$
+\frac{1}{2}\sum\limits_{j_{g_5}=p+1}^{\infty}\sum\limits_{j_{g_1}=p+1}^{\infty}
C_{j_k \ldots j_1}\biggl|_{(j_{g_4} j_{g_3})\curvearrowright (\cdot),
j_{g_{{}_{1}}}=~j_{g_{{}_{2}}}, j_{g_{{}_{3}}}=~j_{g_{{}_{4}}}, j_{g_{{}_{5}}}=~j_{g_{{}_{6}}}}
{\bf 1}_{\{g_{4}=g_{3}+1\}}+
$$

\vspace{-1mm}
$$
+\frac{1}{2}\sum\limits_{j_{g_5}=p+1}^{\infty}\sum\limits_{j_{g_3}=p+1}^{\infty}
C_{j_k \ldots j_1}\biggl|_{(j_{g_2} j_{g_1})\curvearrowright (\cdot),
j_{g_{{}_{1}}}=~j_{g_{{}_{2}}}, j_{g_{{}_{3}}}=~j_{g_{{}_{4}}}, j_{g_{{}_{5}}}=~j_{g_{{}_{6}}}}
{\bf 1}_{\{g_{2}=g_{1}+1\}}-
$$

\vspace{-1mm}
$$
\left.-
\sum\limits_{j_{g_5}=p+1}^{\infty}
\sum\limits_{j_{g_3}=p+1}^{\infty}\sum\limits_{j_{g_1}=p+1}^{\infty}
C_{j_k \ldots j_1}\biggl|_{j_{g_{{}_{1}}}=~j_{g_{{}_{2}}}, 
j_{g_{{}_{3}}}=~j_{g_{{}_{4}}}, j_{g_{{}_{5}}}=~j_{g_{{}_{6}}}}
\right)\times
$$

\vspace{5mm}
$$
\times
J'[\phi_{j_{q_1}}\ldots \phi_{j_{q_{k-6}}}]_{T,t}^{(i_{q_1}\ldots i_{q_{k-6}})}=
$$

\vspace{-2mm}
$$
=
\frac{1}{2^3}\prod\limits_{s=1}^3
{\bf 1}_{\{g_{2s}=g_{2s-1}+1\}}
J[\psi^{(k)}]_{T,t}^{s_3,s_2,s_1}+
\prod\limits_{s=1}^3
{\bf 1}_{\{i_{g_{{}_{2s-1}}}=~i_{g_{{}_{2s}}}\ne 0\}}
\hbox{\vtop{\offinterlineskip\halign{
\hfil#\hfil\cr
{\rm l.i.m.}\cr
$\stackrel{}{{}_{p\to \infty}}$\cr
}} } R_{T,t}^{(p)3,g_1,g_2,\ldots,g_5,g_6}
$$

\vspace{3mm}
\noindent
w.~p.~1, where $g_{2i-1}\stackrel{\sf def}{=}s_i;$\ $i=1,2,3,$ 
$(s_3,s_2,s_1)\in {\rm A}_{k,3},$ $J[\psi^{(k)}]_{T,t}^{s_3,s_2,s_1}$ is
defined by (\ref{30.1}) and ${\rm A}_{k,3}$ is defined by (\ref{30.5550001}),

\newpage
\noindent
$$
R_{T,t}^{(p)3,g_1,g_2,\ldots,g_5,g_6}=
\sum\limits_{\stackrel{j_1,\ldots,j_q,\ldots,j_k=0}{{}_{q\ne g_1,g_2,\ldots,g_5,g_6}}}^p
\left(
-\bar C^{(p)}_{j_k\ldots j_q \ldots j_1}\biggl|_{q\ne g_1,g_2,\ldots,g_5,g_6}+\right.
$$

\vspace{-2mm}
$$
+S_1\left\{
\bar C^{(p)}_{j_k\ldots j_q \ldots j_1}\biggl|_{q\ne g_1,g_2,\ldots,g_5,g_6}\right\}
+S_2\left\{
\bar C^{(p)}_{j_k\ldots j_q \ldots j_1}\biggl|_{q\ne g_1,g_2,\ldots,g_5,g_6}\right\}+
$$

\vspace{-2mm}
$$       
+S_3\left\{
\bar C^{(p)}_{j_k\ldots j_q \ldots j_1}\biggl|_{q\ne g_1,g_2,\ldots,g_5,g_6}\right\}-
$$

\vspace{-2mm}
$$
-S_3S_1\left\{
\bar C^{(p)}_{j_k\ldots j_q \ldots j_1}\biggl|_{q\ne g_1,g_2,\ldots,g_5,g_6}\right\}
-S_3S_2\left\{
\bar C^{(p)}_{j_k\ldots j_q \ldots j_1}\biggl|_{q\ne g_1,g_2,\ldots,g_5,g_6}\right\}-
$$

$$
\left.-S_2S_1\left\{
\bar C^{(p)}_{j_k\ldots j_q \ldots j_1}\biggl|_{q\ne g_1,g_2,\ldots,g_5,g_6}\right\}\right)
J'[\phi_{j_{q_1}}\ldots \phi_{j_{q_{k-6}}}]_{T,t}^{(i_{q_1}\ldots i_{q_{k-6}})}.
$$

\vspace{4mm}

By Condition~3 of Theorem~2.30 we have (also see the property (\ref{wiener1}) of
multiple Wiener stochastic integral)
$$
\lim\limits_{p\to\infty}
{\sf M}\left\{\left(R_{T,t}^{(p)3,g_1,g_2,\ldots,g_5,g_6}\right)^2\right\}
\le K \lim\limits_{p\to\infty}
\sum\limits_{\stackrel{j_1,\ldots,j_q,\ldots,j_k=0}{{}_{q\ne g_1, g_2,\ldots, g_5, g_6}}}^p
\left(\left(
\bar C^{(p)}_{j_k\ldots j_q \ldots j_1}\biggl|_{q\ne g_1,g_2,\ldots, g_5, g_6}\right)^2+\right.
$$

\vspace{-2mm}
$$
+\left(S_1\left\{
\bar C^{(p)}_{j_k\ldots j_q \ldots j_1}\biggl|_{q\ne g_1,g_2,\ldots,g_5,g_6}\right\}\right)^2
+\left(S_2\left\{
\bar C^{(p)}_{j_k\ldots j_q \ldots j_1}\biggl|_{q\ne g_1,g_2,\ldots,g_5,g_6}\right\}\right)^2+
$$

\vspace{-2mm}
$$
+\left(S_3\left\{
\bar C^{(p)}_{j_k\ldots j_q \ldots j_1}\biggl|_{q\ne g_1,g_2,\ldots,g_5,g_6}\right\}\right)^2+
$$

\vspace{-2mm}
$$
+\left(S_3S_1\left\{
\bar C^{(p)}_{j_k\ldots j_q \ldots j_1}\biggl|_{q\ne g_1,g_2,\ldots,g_5,g_6}\right\}\right)^2
+\left(S_3S_2\left\{
\bar C^{(p)}_{j_k\ldots j_q \ldots j_1}\biggl|_{q\ne g_1,g_2,\ldots,g_5,g_6}\right\}\right)^2+
$$

\vspace{-2mm}
$$
\left.+\left(S_2S_1\left\{
\bar C^{(p)}_{j_k\ldots j_q \ldots j_1}\biggl|_{q\ne g_1,g_2,\ldots,g_5,g_6}\right\}\right)^2
\right)=0,
$$

\vspace{4mm}
\noindent
where constant $K$ does not depend on $p$.

Thus

\vspace{-1mm}
$$
\hbox{\vtop{\offinterlineskip\halign{
\hfil#\hfil\cr
{\rm l.i.m.}\cr
$\stackrel{}{{}_{p\to \infty}}$\cr
}} }\prod\limits_{s=1}^3
{\bf 1}_{\{i_{g_{{}_{2s-1}}}=~i_{g_{{}_{2s}}}\ne 0\}}\hbox{\vtop{\offinterlineskip\halign{
\hfil#\hfil\cr
{\rm l.i.m.}\cr
$\stackrel{}{{}_{p\to \infty}}$\cr
}} }\sum_{j_1,\ldots,j_k=0}^{p}
C_{j_k\ldots j_1}
\prod\limits_{s=1}^3{\bf 1}_{\{j_{g_{{}_{2s-1}}}=~j_{g_{{}_{2s}}}\}}\times
$$

\vspace{2mm}
$$
\times
J'[\phi_{j_{q_1}}\ldots \phi_{j_{q_{k-6}}}]_{T,t}^{(i_{q_1}\ldots i_{q_{k-6}})}=
\frac{1}{2^3}\prod\limits_{s=1}^3
{\bf 1}_{\{g_{2s}=g_{2s-1}+1\}}
J[\psi^{(k)}]_{T,t}^{s_3,s_2,s_1}
\ \ \ \hbox{w.~p.~1},
$$

\vspace{3mm}
\noindent
where $g_{2i-1}\stackrel{\sf def}{=}s_i;$\ $i=1,2,3,$ 
$(s_3,s_2,s_1)\in {\rm A}_{k,3},$ $J[\psi^{(k)}]_{T,t}^{s_3,s_2,s_1}$ is
defined by (\ref{30.1}) and ${\rm A}_{k,3}$ is defined by (\ref{30.5550001}).

Repeating the previous steps, we obtain for an arbitrary $r$ 
($r=1,2,$ $\ldots,$ $[k/2]$)

\vspace{-2mm}
$$
\prod\limits_{s=1}^r
{\bf 1}_{\{i_{g_{{}_{2s-1}}}=~i_{g_{{}_{2s}}}\ne 0\}}\hbox{\vtop{\offinterlineskip\halign{
\hfil#\hfil\cr
{\rm l.i.m.}\cr
$\stackrel{}{{}_{p\to \infty}}$\cr
}} }\sum_{j_1,\ldots,j_k=0}^{p}
C_{j_k\ldots j_1}
\prod\limits_{s=1}^r{\bf 1}_{\{j_{g_{{}_{2s-1}}}=~j_{g_{{}_{2s}}}\}}\times
$$

\vspace{4mm}
$$
\times
J'[\phi_{j_{q_1}}\ldots \phi_{j_{q_{k-2r}}}]_{T,t}^{(i_{q_1}\ldots i_{q_{k-2r}})}=
\prod\limits_{s=1}^r
{\bf 1}_{\{i_{g_{{}_{2s-1}}}=~i_{g_{{}_{2s}}}\ne 0\}}\times
$$

\vspace{-2mm}
$$
\times\hbox{\vtop{\offinterlineskip\halign{
\hfil#\hfil\cr
{\rm l.i.m.}\cr
$\stackrel{}{{}_{p\to \infty}}$\cr
}} }\hspace{-2.5mm}
\sum\limits_{\stackrel{j_1,\ldots,j_q,\ldots,j_k=0}{{}_{q\ne g_1, g_2,\ldots, g_{2r-1}, g_{2r}}}}^p
\hspace{-2.5mm}
\frac{1}{2^r}
C_{j_k \ldots j_1}\biggl|_{(j_{g_2} j_{g_1})\curvearrowright (\cdot)
\ldots (j_{g_{2r}} j_{g_{2r-1}})\curvearrowright (\cdot),
j_{g_{{}_{1}}}=~j_{g_{{}_{2}}},\ldots, j_{g_{{}_{2r-1}}}=~j_{g_{{}_{2r}}}
}\biggr.
\times
$$

\vspace{2mm}
$$
\times\prod\limits_{s=1}^r
{\bf 1}_{\{g_{2s}=g_{2s-1}+1\}}
J'[\phi_{j_{q_1}}\ldots \phi_{j_{q_{k-2r}}}]_{T,t}^{(i_{q_1}\ldots i_{q_{k-2r}})}
+
$$

\vspace{2mm}
\begin{equation}
\label{after903}
+\prod\limits_{s=1}^r
{\bf 1}_{\{i_{g_{{}_{2s-1}}}=~i_{g_{{}_{2s}}}\ne 0\}}
\hbox{\vtop{\offinterlineskip\halign{
\hfil#\hfil\cr
{\rm l.i.m.}\cr
$\stackrel{}{{}_{p\to \infty}}$\cr
}} } R_{T,t}^{(p)r,g_1,g_2,\ldots,g_{2r-1},g_{2r}}=
\end{equation}

\vspace{-1mm}
\begin{equation}
\label{after904}
=\frac{1}{2^r}\prod\limits_{s=1}^r
{\bf 1}_{\{g_{2s}=g_{2s-1}+1\}}
J[\psi^{(k)}]_{T,t}^{s_r, \ldots, s_1} + 
\prod\limits_{s=1}^r
{\bf 1}_{\{i_{g_{{}_{2s-1}}}=~i_{g_{{}_{2s}}}\ne 0\}}
\hbox{\vtop{\offinterlineskip\halign{
\hfil#\hfil\cr
{\rm l.i.m.}\cr
$\stackrel{}{{}_{p\to \infty}}$\cr
}} } R_{T,t}^{(p)r,g_1,g_2,\ldots,g_{2r-1},g_{2r}}
\end{equation}

\vspace{1mm}
\noindent
w.~p.~1, where $g_{2i-1}\stackrel{\sf def}{=}s_i;$\ $i=1,2,\ldots,r;$\
$r=1,2,\ldots,\left[k/2\right],$ 
$(s_r,\ldots,s_1)\in {\rm A}_{k,r},$ $J[\psi^{(k)}]_{T,t}^{s_r,\ldots,s_1}$ is
defined by (\ref{30.1}) and ${\rm A}_{k,r}$ is defined by (\ref{30.5550001}),

\newpage
\noindent
$$
R_{T,t}^{(p)r,g_1,g_2,\ldots,g_{2r-1},g_{2r}}=
$$

\vspace{-2mm}
$$
=
\sum\limits_{\stackrel{j_1,\ldots,j_q,\ldots,j_k=0}{{}_{q\ne g_1, g_2, \ldots, g_{2r-1}, g_{2r}}}}^p
\left(
(-1)^r \bar C^{(p)}_{j_k\ldots j_q \ldots j_1}\biggl|_{q\ne g_1,g_2, \ldots, g_{2r-1}, g_{2r}}+\right.
$$

\vspace{1mm}
$$
+(-1)^{r-1}\sum\limits_{l_1=1}^r S_{l_1}\left\{
\bar C^{(p)}_{j_k\ldots j_q \ldots j_1}\biggl|_{q\ne g_1,g_2, \ldots, g_{2r-1}, g_{2r}}\right\}+
$$

\vspace{3mm}
$$
+(-1)^{r-2}\sum\limits_{\stackrel{l_1,l_2=1}{{}_{l_1>l_2}}}^r
S_{l_1}S_{l_2}\left\{
\bar C^{(p)}_{j_k\ldots j_q \ldots j_1}\biggl|_{q\ne g_1,g_2, \ldots, g_{2r-1}, g_{2r}}\right\}+
$$

\vspace{-5mm}
$$
\ldots
$$
$$
\left.+(-1)^{1}\sum\limits_{\stackrel{l_1,l_2,\ldots, l_{r-1}=1}{{}_{l_1>l_2>\ldots > l_{r-1}}}}^r
S_{l_1}S_{l_2}\ldots S_{l_{r-1}}\left\{
\bar C^{(p)}_{j_k\ldots j_q \ldots j_1}\biggl|_{q\ne g_1,g_2, \ldots, g_{2r-1}, g_{2r}}\right\}
\right)\times
$$

\vspace{1mm}
\begin{equation}
\label{afterr1}
\times
J'[\phi_{j_{q_1}}\ldots \phi_{j_{q_{k-2r}}}]_{T,t}^{(i_{q_1}\ldots i_{q_{k-2r}})}.
\end{equation}

\vspace{7mm}

Let us explain the transition from (\ref{after903}) to (\ref{after904}).
We have for $g_{2}=g_{1}+1,\ldots, g_{2r}=g_{2r-1}+1$

\vspace{-3mm}
$$
\hbox{\vtop{\offinterlineskip\halign{
\hfil#\hfil\cr
{\rm l.i.m.}\cr
$\stackrel{}{{}_{p\to \infty}}$\cr
}} }
\sum\limits_{\stackrel{j_1,\ldots,j_q,\ldots,j_k=0}{{}_{q\ne g_1, g_2,\ldots, g_{2r-1}, g_{2r}}}}^p
\frac{1}{2^r}
C_{j_k \ldots j_1}\biggl|_{(j_{g_2} j_{g_1})\curvearrowright (\cdot)
\ldots (j_{g_{2r}} j_{g_{2r-1}})\curvearrowright (\cdot),
j_{g_{{}_{1}}}=~j_{g_{{}_{2}}},\ldots, j_{g_{{}_{2r-1}}}=~j_{g_{{}_{2r}}}}\biggr.
\times 
$$

\vspace{2mm}
$$
\times
\prod\limits_{s=1}^r
{\bf 1}_{\{i_{g_{{}_{2s-1}}}=~i_{g_{{}_{2s}}}\ne 0\}}
J'[\phi_{j_{q_1}}\ldots \phi_{j_{q_{k-2r}}}]_{T,t}^{(i_{q_1}\ldots i_{q_{k-2r}})}=
$$

\vspace{1mm}
$$
=\frac{1}{2^r}\hbox{\vtop{\offinterlineskip\halign{
\hfil#\hfil\cr
{\rm l.i.m.}\cr
$\stackrel{}{{}_{p\to \infty}}$\cr
}} }
\sum\limits_{\stackrel{j_1,\ldots,j_q,\ldots,j_k=0}{{}_{q\ne g_1, g_2,\ldots, g_{2r-1}, g_{2r}}}}^p
C_{j_k \ldots j_1}\biggl|_{(j_{g_2} j_{g_1})\curvearrowright 0
\ldots (j_{g_{2r}} j_{g_{2r-1}})\curvearrowright 0,
j_{g_{{}_{1}}}=~j_{g_{{}_{2}}},\ldots, j_{g_{{}_{2r-1}}}=~j_{g_{{}_{2r}}}}\biggr.
\times 
$$

\vspace{2mm}
$$
\times\prod\limits_{s=1}^r
{\bf 1}_{\{i_{g_{{}_{2s-1}}}=~i_{g_{{}_{2s}}}\ne 0\}}
\left(\zeta_{0}^{(0)}\right)^r 
J'[\phi_{j_{q_1}}\ldots \phi_{j_{q_{k-2r}}}]_{T,t}^{(i_{q_1}\ldots i_{q_{k-2r}})}=
$$

\newpage
\noindent
$$
=\frac{1}{2^r}\hbox{\vtop{\offinterlineskip\halign{
\hfil#\hfil\cr
{\rm l.i.m.}\cr
$\stackrel{}{{}_{p\to \infty}}$\cr
}} }
\sum\limits_{\stackrel{j_1,\ldots,j_q,\ldots,j_k=0}{{}_{q\ne g_1, g_2,\ldots, g_{2r-1}, g_{2r}}}}^p
\sum_{j_{m_1}, j_{m_3}\ldots,j_{m_{2r-1}}=0}^{p}
\prod\limits_{s=1}^r
{\bf 1}_{\{i_{g_{{}_{2s-1}}}=~i_{g_{{}_{2s}}}\ne 0\}}
\times
$$
$$
\times
C_{j_k \ldots j_1}\biggl|_{(j_{g_2} j_{g_1})\curvearrowright j_{m_1}
\ldots (j_{g_{2r}} j_{g_{2r-1}})\curvearrowright  j_{m_{2r-1}},
j_{g_{{}_{1}}}=~j_{g_{{}_{2}}},\ldots, j_{g_{{}_{2r-1}}}=~j_{g_{{}_{2r}}}}\biggr.
\times 
$$

\vspace{5mm}
$$
\times\zeta_{ j_{m_{1}}}^{(0)} \zeta_{ j_{m_{3}}}^{(0)}\ldots \zeta_{ j_{m_{2r-1}}}^{(0)}
J'[\phi_{j_{q_1}}\ldots \phi_{j_{q_{k-2r}}}]_{T,t}^{(i_{q_1}\ldots i_{q_{k-2r}})}=
$$

\vspace{1mm}
$$
=\frac{1}{2^r}\hbox{\vtop{\offinterlineskip\halign{
\hfil#\hfil\cr
{\rm l.i.m.}\cr
$\stackrel{}{{}_{p\to \infty}}$\cr
}} }
\sum\limits_{\stackrel{j_1,\ldots,j_q,\ldots,j_k=0}{{}_{q\ne g_1, g_2,\ldots, g_{2r-1}, g_{2r}}}}^p
\sum_{j_{m_1}, j_{m_3}\ldots,j_{m_{2r-1}}=0}^{p}
\prod\limits_{s=1}^r
{\bf 1}_{\{i_{g_{{}_{2s-1}}}=~i_{g_{{}_{2s}}}\ne 0\}}\times
$$

\vspace{1mm}
$$
\times
C_{j_k \ldots j_1}\biggl|_{(j_{g_2} j_{g_1})\curvearrowright j_{m_1}
\ldots (j_{g_{2r}} j_{g_{2r-1}})\curvearrowright  j_{m_{2r-1}},
j_{g_{{}_{1}}}=~j_{g_{{}_{2}}},\ldots, j_{g_{{}_{2r-1}}}=~j_{g_{{}_{2r}}}}\biggr.
\times
$$

\vspace{2mm}
\begin{equation}
\label{after905}
~~~~~~~~\times
J'[\phi_{ j_{m_1}}\phi_{ j_{m_3}}
\ldots \phi_{j_{m_{2r-1}}}
\phi_{j_{q_1}}\ldots \phi_{j_{q_{k-2r}}}]_{T,t}^{(00\ldots 0 i_{q_1}\ldots i_{q_{k-2r}})}=
\end{equation}

\begin{equation}
\label{after906}
=\frac{1}{2^r}
J[\psi^{(k)}]_{T,t}^{s_r, \ldots, s_1}\ \ \ \hbox{w.~p.~1}.
\end{equation}

\vspace{5mm}
\noindent
The transition from (\ref{after905}) to (\ref{after906}) is based
on (\ref{drdr1}) or (\ref{febr5000}).

Note that
$$
C_{j_k \ldots j_1}\biggl|_{(j_{g_2} j_{g_1})\curvearrowright j_{m_1},
j_{g_{{}_{1}}}=~j_{g_{{}_{2}}}}\biggr.
=
C_{j_k \ldots j_1}\biggl|_{(j_{g_1} j_{g_1})\curvearrowright j_{m_1},
j_{g_{{}_{1}}}=~j_{g_{{}_{2}}}}\biggr.
$$

\vspace{3mm}
\noindent
is the Fourier coefficient, where $g_{2}=g_{1}+1$. 
Therefore, the value
$$
C_{j_k \ldots j_1}\biggl|_{(j_{g_2} j_{g_1})\curvearrowright j_{m_1}
\ldots (j_{g_{2d}} j_{g_{2d-1}})\curvearrowright j_{m_{2d-1}},
j_{g_{{}_{1}}}=~j_{g_{{}_{2}}},\ldots,
j_{g_{{}_{2d-1}}}=~j_{g_{{}_{2d}}}}\biggr.=
$$
$$
=C_{j_k \ldots j_1}\biggl|_{(j_{g_1} j_{g_1})\curvearrowright j_{m_1}
\ldots (j_{g_{2d-1}} j_{g_{2d-1}})\curvearrowright j_{m_{2d-1}},
j_{g_{{}_{1}}}=~j_{g_{{}_{2}}},\ldots,
j_{g_{{}_{2d-1}}}=~j_{g_{{}_{2d}}}}\biggr.
$$

\vspace{4mm}
\noindent
is determined recursively using (\ref{after2000}) in an obvious way
for $g_2=g_1+1,$ $\ldots,$ $g_{2d}=g_{2d-1}+1$ and $d=2,\ldots,r$.

By Condition~3 of Theorem~2.30 we have (also see the property (\ref{wiener1}) of
multiple Wiener stochastic integral)

\vspace{-1mm}
$$
\lim\limits_{p\to\infty}
{\sf M}\left\{\left(R_{T,t}^{(p)r,g_1,g_2,\ldots,g_{2r-1},g_{2r}}\right)^2\right\}
\le 
$$

\vspace{-2mm}
$$
\le
K \lim\limits_{p\to\infty}
\sum\limits_{\stackrel{j_1,\ldots,j_q,\ldots,j_k=0}{{}_{q\ne g_1, g_2, \ldots, g_{2r-1}, g_{2r}}}}^p
\left(\left(
\bar C^{(p)}_{j_k\ldots j_q \ldots j_1}\biggl|_{q\ne g_1,g_2, \ldots, g_{2r-1}, g_{2r}}
\right)^2+\right.
$$

\vspace{3mm}
$$
+\sum\limits_{l_1=1}^r \left(S_{l_1}\left\{
\bar C^{(p)}_{j_k\ldots j_q \ldots j_1}\biggl|_{q\ne g_1,g_2, \ldots, g_{2r-1}, g_{2r}}\right\}
\right)^2+
$$

\vspace{4mm}
$$
+\sum\limits_{\stackrel{l_1,l_2=1}{{}_{l_1>l_2}}}^r
\left(S_{l_1}S_{l_2}\left\{
\bar C^{(p)}_{j_k\ldots j_q \ldots j_1}\biggl|_{q\ne g_1,g_2, \ldots, g_{2r-1}, g_{2r}}\right\}
\right)^2+
$$

\vspace{-4mm}
$$
\ldots
$$
$$
\left.+\sum\limits_{\stackrel{l_1,l_2,\ldots, l_{r-1}=1}{{}_{l_1>l_2>\ldots > l_{r-1}}}}^r
\left(S_{l_1}S_{l_2}\ldots S_{l_{r-1}}\left\{
\bar C^{(p)}_{j_k\ldots j_q \ldots j_1}\biggl|_{q\ne g_1,g_2, \ldots, g_{2r-1}, g_{2r}}\right\}
\right)^2
\right)=0,
$$

\vspace{4mm}
\noindent
where constant $K$ does not depend on $p$.

So we have

\vspace{-1mm}
$$
\prod\limits_{s=1}^r
{\bf 1}_{\{i_{g_{{}_{2s-1}}}=~i_{g_{{}_{2s}}}\ne 0\}}\hbox{\vtop{\offinterlineskip\halign{
\hfil#\hfil\cr
{\rm l.i.m.}\cr
$\stackrel{}{{}_{p\to \infty}}$\cr
}} }\sum_{j_1,\ldots,j_k=0}^{p}
C_{j_k\ldots j_1}
\prod\limits_{s=1}^r{\bf 1}_{\{j_{g_{{}_{2s-1}}}=~j_{g_{{}_{2s}}}\}}\times
$$

\vspace{2.5mm}
$$
\times
J'[\phi_{j_{q_1}}\ldots \phi_{j_{q_{k-2r}}}]_{T,t}^{(i_{q_1}\ldots i_{q_{k-2r}})}=
$$

\begin{equation}
\label{after801}
=\frac{1}{2^r}\prod\limits_{s=1}^r
{\bf 1}_{\{g_{2s}=g_{2s-1}+1\}}
J[\psi^{(k)}]_{T,t}^{s_r, \ldots, s_1}\ \ \ \hbox{w.~p.~1},
\end{equation}

\vspace{2mm}
\noindent
where $g_{2i-1}\stackrel{\sf def}{=}s_i;$\ $i=1,2,\ldots,r;$\
$r=1,2,\ldots,\left[k/2\right],$ 
$(s_r,\ldots,s_1)\in {\rm A}_{k,r},$ $J[\psi^{(k)}]_{T,t}^{s_r,\ldots,s_1}$ is
defined by (\ref{30.1}) and ${\rm A}_{k,r}$ is defined by (\ref{30.5550001}).

Note that

\vspace{-5mm}
$$
\sum_{\stackrel{(\{\{g_1, g_2\}, \ldots, 
\{g_{2r-1}, g_{2r}\}\}, \{q_1, \ldots, q_{k-2r}\})}
{{}_{\{g_1, g_2, \ldots, 
g_{2r-1}, g_{2r}, q_1, \ldots, q_{k-2r}\}=\{1, 2, \ldots, k\}}}}
\Biggl|_{g_2=g_1+1, g_3=g_2+1,\ldots, g_{2r}=g_{2r-1}+1}\Biggr.
A_{g_1,g_3,\ldots,g_{2r-1}}=
$$

\begin{equation}
\label{after800}
=\sum\limits_{(s_r,\ldots,s_1)\in {\rm A}_{k,r}}
A_{s_1,s_2,\ldots,s_r},
\end{equation}

\vspace{3mm}
\noindent
where $A_{g_1,g_3,\ldots,g_{2r-1}},$ 
$A_{s_1,s_2,\ldots,s_r}$ are scalar values,
$g_{2i-1}=s_i;$\ $i=1,2,\ldots,r;$\ $r=1,2,\ldots,\left[k/2\right],$
${\rm A}_{k,r}$ is defined by (\ref{30.5550001}):

\vspace{-3mm}
$$
{\rm A}_{k,r}
=\bigl\{(s_r,\ldots,s_1):\
s_r>s_{r-1}+1,\ldots,s_2>s_1+1,\ s_r,\ldots,s_1=1,\ldots,k-1\bigr\}.
$$

\vspace{4mm}

Using (\ref{after501}), (\ref{after801}), (\ref{after800}), and Theorem~2.12,
we finally get                                                           
$$
\hbox{\vtop{\offinterlineskip\halign{
\hfil#\hfil\cr
{\rm l.i.m.}\cr
$\stackrel{}{{}_{p\to \infty}}$\cr
}} }\sum_{j_1,\ldots,j_k=0}^{p}
C_{j_k\ldots j_1}
\prod\limits_{l=1}^k \zeta_{j_l}^{(i_l)}=
\hbox{\vtop{\offinterlineskip\halign{
\hfil#\hfil\cr
{\rm l.i.m.}\cr
$\stackrel{}{{}_{p\to \infty}}$\cr
}} }\sum_{j_1,\ldots,j_k=0}^{p}
C_{j_k\ldots j_1}
\zeta_{j_1}^{(i_1)}\ldots \zeta_{j_k}^{(i_k)}
=
$$
\begin{equation}
\label{afteru11}
~~~~~~~=J[\psi^{(k)}]_{T,t}^{(i_1\ldots i_k)}+
\sum_{r=1}^{\left[k/2\right]}\frac{1}{2^r}
\sum_{(s_r,\ldots,s_1)\in {\rm A}_{k,r}}
J[\psi^{(k)}]_{T,t}^{s_r, \ldots, s_1}=J^{*}[\psi^{(k)}]_{T,t}^{(i_1\ldots i_k)}
\end{equation}

\vspace{3mm}
\noindent
w.~p.~1, where (see (\ref{30.1}))
$$
J[\psi^{(k)}]_{T,t}^{s_r, \ldots, s_1} \stackrel{\rm def}{=}\ 
\prod_{q=1}^r {\bf 1}_{\{i_{s_q}=
i_{s_{q}+1}\ne 0\}}\ \times
$$
$$
\times
\int\limits_t^T\psi_k(t_k)\ldots \int\limits_t^{t_{s_r+3}}
\psi_{s_r+2}(t_{s_r+2})
\int\limits_t^{t_{s_r+2}}\psi_{s_r}(t_{s_r+1})
\psi_{s_r+1}(t_{s_r+1}) \times
$$
$$
\times
\int\limits_t^{t_{s_r+1}}\psi_{s_r-1}(t_{s_r-1})
\ldots
\int\limits_t^{t_{s_1+3}}\psi_{s_1+2}(t_{s_1+2})
\int\limits_t^{t_{s_1+2}}\psi_{s_1}(t_{s_1+1})
\psi_{s_1+1}(t_{s_1+1}) \times
$$
$$
\times
\int\limits_t^{t_{s_1+1}}\psi_{s_1-1}(t_{s_1-1})
\ldots \int\limits_t^{t_2}\psi_1(t_1)
d{\bf w}_{t_1}^{(i_1)}\ldots d{\bf w}_{t_{s_1-1}}^{(i_{s_1-1})}
dt_{s_1+1}d{\bf w}_{t_{s_1+2}}^{(i_{s_1+2})}\ldots
$$
\begin{equation}
\label{afterito1}
\ldots\
d{\bf w}_{t_{s_r-1}}^{(i_{s_r-1})}
dt_{s_r+1}d{\bf w}_{t_{s_r+2}}^{(i_{s_r+2})}\ldots d{\bf w}_{t_k}^{(i_k)}.
\end{equation}

\vspace{2mm}

Theorem~2.30 is proved.

{\bf Remark~2.4.}\ {\it Let us make a number of remarks about Theorem~{\rm 2.30}.
An expansion similar to {\rm (\ref{after1})}
was obtained in {\rm \cite{Rybakov3000}} for the multiple
Stratonovich stochastic integral
{\rm (\ref{january19})}.
The proof from {\rm \cite{Rybakov3000}} is somewhat simpler than the
proof proposed in this section. 
However$,$ the results from {\rm \cite{Rybakov3000}} 
were obtained under 
the condition of convergence of 
trace series. 
The verification of this condition for the kernel {\rm (\ref{ppp})} is a separate 
problem. 
In our proof
we essentially use the structure of the Fourier coefficients {\rm (\ref{after1000})}
corresponding to the kernel {\rm (\ref{ppp}).} 
This circumstance actually made it possible
to prove Theorem~{\rm 2.30} using not the condition of
finiteness of trace series$,$ but using the condition 
of convergence to zero of explicit expressions for the remainders
of the mentioned series. This leaves hope 
that it is possible to prove an analogue of Theorems 
{\rm 2.24--2.26, 2.37--2.39}
for the case of an arbitrary $k$ $(k\in{\bf N})$.

Note that under the conditions of Theorem~{\rm 2.30} 
{\rm (}also see {\rm (\ref{after80xx}), (\ref{after500}))}
the sequential order of the series
$$
\sum\limits_{j_{g_{2r-1}}=p+1}^{\infty}
\sum\limits_{j_{g_{2r-3}}=p+1}^{\infty}
\ldots \sum\limits_{j_{g_{3}}=p+1}^{\infty}
\sum\limits_{j_{g_{1}}=p+1}^{\infty}
$$

\noindent
in {\rm (\ref{after5days})} is not important.
We also note that Conditions $1, 2$
of Theorem~{\rm 2.30} are satisfied for complete
orthonormal systems of Legendre polynomials 
and trigonometric functions in the space
$L_2([t, T])$ {\rm (}see the proofs of 
Theorems~{\rm 2.1, 2.2, 2.8, 2.9, 2.27)}.
Moreover$,$ the equality  {\rm (\ref{after200})} is true for an arbitrary 
basis in $L_2([t, T])$ {\rm (}see {\rm (\ref{5tzzz})).}
It is easy to see that in the proofs of Theorems~{\rm 2.1--2.9}
the conditions of Theorem~{\rm 2.30}
are verified for various special cases
of iterated Stratonovich stochastic integrals
of multiplicities {\rm 2--4}.
}

Taking into account Theorem~1.11, 
we can formulate an analogue of Theorem~2.30
for the case of integration interval $[t, s]$ $(s\in (t, T);$ 
the case $s=T$ is considered in Theorem~2.30) 
of iterated Stratonovich stochastic integrals of multiplicity 
$k$ ($k\in{\bf N}$).

Denote
$$
\bar C^{(p)}_{j_k\ldots j_q \ldots j_1}(s)\biggl|_{q\ne g_1,g_2,\ldots,g_{2r-1}, g_{2r}}
\stackrel{\sf def}{=}
$$

\vspace{-1mm}
$$
\stackrel{\sf def}{=}
\sum\limits_{j_{g_{2r-1}}=p+1}^{\infty}
\sum\limits_{j_{g_{2r-3}}=p+1}^{\infty}
\ldots \sum\limits_{j_{g_{3}}=p+1}^{\infty}
\sum\limits_{j_{g_{1}}=p+1}^{\infty}
C_{j_k \ldots j_1}(s)\biggl|_{j_{g_1}=j_{g_2},\ldots, j_{g_{2r-1}}=j_{g_{2r}}}\biggl.
$$

\newpage
\noindent
and introduce the following notation

\vspace{-3mm}
$$
S_l \left\{\bar C^{(p)}_{j_k\ldots j_q \ldots j_1}(s)\biggl|_{q\ne g_1,g_2,\ldots,g_{2r-1}, g_{2r}}
\right\}
\stackrel{\sf def}{=}
\frac{1}{2}{\bf 1}_{\{g_{2l}=g_{2l-1}+1\}}
\sum\limits_{j_{g_{2r-1}}=p+1}^{\infty}
\sum\limits_{j_{g_{2r-3}}=p+1}^{\infty}
\ldots 
$$

\vspace{-4mm}
$$
\ldots
\sum\limits_{j_{g_{2l+1}}=p+1}^{\infty}
\sum\limits_{j_{g_{2l-3}}=p+1}^{\infty}
\ldots
\sum\limits_{j_{g_{3}}=p+1}^{\infty}
\sum\limits_{j_{g_{1}}=p+1}^{\infty}
C_{j_k \ldots j_1}(s)\biggl|_{(j_{g_{2l}} j_{g_{2l-1}})\curvearrowright (\cdot),
j_{g_1}=j_{g_2},\ldots, j_{g_{2r-1}}=j_{g_{2r}}}\biggr.,
$$

\vspace{4mm}
\noindent
where $l=1,2,\ldots,r,$ 
$$
C_{j_k \ldots j_1}(s)\Biggl|_{(j_{g_{2l}} j_{g_{2l-1}})\curvearrowright (\cdot)}\Biggr.
$$

\vspace{1mm}
\noindent
is defined by analogy with (\ref{after900}), 
\begin{equation}
\label{after1300}
~~~~~~~C_{j_k \ldots j_1}(s)=\int\limits_t^s\psi_k(t_k)\phi_{j_k}(t_k)\ldots
\int\limits_t^{t_2}
\psi_1(t_1)\phi_{j_1}(t_1)
dt_1\ldots dt_k.
\end{equation}

{\bf Theorem~2.31}\ \cite{arxiv-5}, 
\cite{arxiv-10}, \cite{arxiv-11}, \cite{new-art-1xxy}.\ {\it Assume that
the continuously differentiable functions 
$\psi_l(\tau)$ $(l=1,\ldots,k)$ at the interval $[t, T]$ and 
the complete orthonormal system $\{\phi_j(x)\}_{j=0}^{\infty}$
of continuous functions $(\phi_0(x)=1/\sqrt{T-t})$ 
in the space $L_2([t, T])$ are such that the following 
conditions are satisfied{\rm :}

{\rm 1.}\ The equality 
\begin{equation}
\label{after1301}
~~~~~\frac{1}{2}
\int\limits_t^s \Phi_1(t_1)\Phi_2(t_1)dt_1
=\sum_{j=0}^{\infty}
\int\limits_t^s
\Phi_2(t_2)\phi_{j}(t_2)\int\limits_t^{t_2}
\Phi_1(t_1)\phi_{j}(t_1)dt_1 dt_2
\end{equation}

\noindent
holds for all $s\in (t, T],$ where the nonrandom functions 
$\Phi_1(\tau),$ $\Phi_2(\tau)$
are continuously differentiable on $[t, T]$
and the series on the right-hand side of {\rm (\ref{after1301})}
converges absolutely.

{\rm 2.}\ The estimates
$$
\left|\int\limits_t^{s} \phi_{j}(\tau)\Phi_1(\tau)d\tau\right|
\le \frac{\Psi_1(s)}{j^{1/2+\alpha}},\ \ \ 
\left|\int\limits_{\tau}^s \phi_{j}(\theta)\Phi_2(\theta)d\theta\right|\le
\frac{\Psi_2(s,\tau)}{j^{1/2+\alpha}},
$$
$$
\left|\sum_{j=p+1}^{\infty}\int\limits_t^{s}
\Phi_2(\tau)\phi_{j}(\tau)\int\limits_t^{\tau}
\Phi_1(\theta)\phi_{j}(\theta)d\theta d\tau\right|\le \frac{\Psi_3(s)}{p^{\hspace{0.5mm}\beta}}
$$

\vspace{1mm}
\noindent
hold for all $s, \tau$ such that $t < \tau < s < T$ and for some $\alpha, \beta >0,$ where 
$\Phi_1(\tau),$ $\Phi_2(\tau)$
are continuously differentiable nonrandom functions on $[t, T],$\ $j, p\in {\bf N},$
and
$$
\int\limits_t^s \left|\Psi_1(\tau)
\Psi_2(s,\tau)\right| d \tau<\infty,\ \ \ 
\int\limits_t^s \left|\Psi_3(\tau)\right| d\tau<\infty
$$
for all $s\in (t, T).$

{\rm 3.}\ The condition
$$
\lim\limits_{p\to\infty}
\sum\limits_{\stackrel{j_1,\ldots,j_q,\ldots,j_k=0}{{}_{q\ne g_1, g_2, \ldots, g_{2r-1},
g_{2r}}}}^p
\left(S_{l_1}S_{l_2}\ldots S_{l_{d}}
\left\{\bar C^{(p)}_{j_k\ldots j_q \ldots j_1}(s)\biggl|_{q\ne g_1,g_2,\ldots,g_{2r-1}, g_{2r}}
\right\}\right)^2=0
$$

\vspace{2mm}
\noindent
holds for all possible $g_1,g_2,\ldots,g_{2r-1},g_{2r}$ {\rm (}see {\rm (\ref{leto5007after}))}
and $l_1, l_2, \ldots, l_{d}$ such that
$l_1, l_2, \ldots, l_{d}\in \{1,2,\ldots, r\},$\
$l_1>l_2>\ldots >l_{d},$\ $d=0, 1, 2,\ldots, r-1,$
where $r=1, 2,\ldots,[k/2]$ and
$$
S_{l_1}S_{l_2}\ldots S_{l_{d}}
\left\{\bar C^{(p)}_{j_k\ldots j_q \ldots j_1}(s)\biggl|_{q\ne g_1,g_2,\ldots,g_{2r-1}, g_{2r}}
\right\}\stackrel{\sf def}{=}
\bar C^{(p)}_{j_k\ldots j_q \ldots j_1}(s)\biggl|_{q\ne g_1,g_2,\ldots,g_{2r-1}, g_{2r}}
$$
for $d=0.$

Then$,$ for the iterated Stratonovich stochastic integral 
of arbitrary multiplicity $k$
\begin{equation}
\label{afterstr1}
~~~~~~~~~~~J^{*}[\psi^{(k)}]_{s,t}^{(i_1\ldots i_k)}=
{\int\limits_t^{*}}^s
\psi_k(t_k) \ldots 
{\int\limits_t^{*}}^{t_{2}}
\psi_1(t_1) d{\bf w}_{t_1}^{(i_1)}\ldots
d{\bf w}_{t_k}^{(i_k)}
\end{equation}
the following 
expansion 
$$
J^{*}[\psi^{(k)}]_{s,t}^{(i_1\ldots i_k)}=
\hbox{\vtop{\offinterlineskip\halign{
\hfil#\hfil\cr
{\rm l.i.m.}\cr
$\stackrel{}{{}_{p\to \infty}}$\cr
}} }
\sum\limits_{j_1,\ldots,j_k=0}^{p}
C_{j_k \ldots j_1}(s)\prod\limits_{l=1}^k \zeta_{j_l}^{(i_l)}
$$

\noindent
that converges in the mean-square sense is valid, where $C_{j_k \ldots j_1}(s)$
is the Fourier coefficient {\rm (\ref{after1300}),} 
${\rm l.i.m.}$ is a limit in the mean-square sense,
$i_1, \ldots, i_k=0, 1,\ldots,m,$ $s\in (t, T),$
$$
\zeta_{j}^{(i)}=
\int\limits_t^T \phi_{j}(\tau) d{\bf w}_{\tau}^{(i)}
$$ 
are independent standard Gaussian random variables for various 
$i$ or $j$ {\rm (}in the case when $i\ne 0${\rm ),}
${\bf w}_{\tau}^{(i)}$ 
$(i=1,\ldots,m)$ are independent standard Wiener processes$,$
${\bf w}_{\tau}^{(0)}=\tau.$}

\vspace{2mm}

It is easy to see that the estimates (\ref{101oh}), (\ref{203}), (\ref{2017x11}), (\ref{101xx}),
and the results of Sect.~2.9 imply 
the fulfillment of Condition~2 of Theorem~2.31 
for complete orthonormal systems of Legendre polynomials and
trigonometric functions 
in the space $L_2([t, T])$.

Also the equality (\ref{5tzzz}) guarantees the fulfillment of Condition 1 
of Theorem~2.31 for these two systems of functions.

It should be noted that (see (\ref{afterr1}))
$$
(-1)^r \bar C^{(p)}_{j_k\ldots j_q \ldots j_1}\biggl|_{q\ne g_1,g_2, \ldots, g_{2r-1}, g_{2r}}+
$$
$$
+(-1)^{r-1}\sum\limits_{l_1=1}^r S_{l_1}\left\{
\bar C^{(p)}_{j_k\ldots j_q \ldots j_1}\biggl|_{q\ne g_1,g_2, \ldots, g_{2r-1}, g_{2r}}\right\}+
$$

\vspace{3mm}
$$
+(-1)^{r-2}\sum\limits_{\stackrel{l_1,l_2=1}{{}_{l_1>l_2}}}^r
S_{l_1}S_{l_2}\left\{
\bar C^{(p)}_{j_k\ldots j_q \ldots j_1}\biggl|_{q\ne g_1,g_2, \ldots, g_{2r-1}, g_{2r}}\right\}+
$$

\vspace{-5mm}
$$
\ldots
$$
$$
+(-1)^{1}\sum\limits_{\stackrel{l_1,l_2,\ldots, l_{r-1}=1}{{}_{l_1>l_2>\ldots > l_{r-1}}}}^r
S_{l_1}S_{l_2}\ldots S_{l_{r-1}}\left\{
\bar C^{(p)}_{j_k\ldots j_q \ldots j_1}\biggl|_{q\ne g_1,g_2, \ldots, g_{2r-1}, g_{2r}}\right\}
=
$$

\vspace{4mm}
$$
=\sum\limits_{j_{g_1},j_{g_3},\ldots,j_{g_{2r-1}}=0}^p
C_{j_k\ldots j_1}\biggl|_{j_{g_1}=j_{g_2},\ldots, j_{g_{2r-1}}=j_{g_{2r}}}-
$$

\vspace{-2mm}
\begin{equation}
\label{drdr1000}
~~-\frac{1}{2^r} \prod\limits_{l=1}^r {\bf 1}_{\{g_{2l}=g_{2l-1}+1\}}
C_{j_k \ldots j_1}\biggl|_{(j_{g_2} j_{g_1})\curvearrowright (\cdot)
\ldots (j_{g_{2r}} j_{g_{2r-1}})\curvearrowright (\cdot),
j_{g_{{}_{1}}}=~j_{g_{{}_{2}}},\ldots, j_{g_{{}_{2r-1}}}=~j_{g_{{}_{2r}}}
}\biggr.,
\end{equation}

\vspace{3mm}
\noindent
where the meaning of the notations used in (\ref{afterr1}) is preserved.

For example, from (\ref{drdr1000}) for the case $r=2$ we get
$$
\sum\limits_{j_{g_3}=p+1}^{\infty}\sum\limits_{j_{g_1}=p+1}^{\infty}
C_{j_k \ldots j_1}\biggl|_{j_{g_{{}_{1}}}=~j_{g_{{}_{2}}}, j_{g_{{}_{3}}}=~j_{g_{{}_{4}}}}-
$$
$$
-
\frac{1}{2}{\bf 1}_{\{g_{4}=g_{3}+1\}}\sum\limits_{j_{g_1}=p+1}^{\infty}
C_{j_k \ldots j_1}\biggl|_{(j_{g_4} j_{g_3})\curvearrowright (\cdot),
j_{g_{{}_{1}}}=~j_{g_{{}_{2}}}, j_{g_{{}_{3}}}=~j_{g_{{}_{4}}}}
-
$$

\vspace{-2mm}
$$
-
\frac{1}{2}{\bf 1}_{\{g_{2}=g_{1}+1\}}\sum\limits_{j_{g_3}=p+1}^{\infty}
C_{j_k \ldots j_1}\biggl|_{(j_{g_2} j_{g_1})\curvearrowright (\cdot),
j_{g_{{}_{1}}}=~j_{g_{{}_{2}}}, j_{g_{{}_{3}}}=~j_{g_{{}_{4}}}}\biggr.
=
$$

$$
=\sum\limits_{j_{g_1}=0}^p \sum\limits_{j_{g_{3}}=0}^p
C_{j_k\ldots j_1}\biggl|_{j_{g_1}=j_{g_2},j_{g_{3}}=j_{g_{4}}}-
$$

\vspace{-2mm}
$$
-\frac{1}{4}{\bf 1}_{\{g_{2}=g_{1}+1\}}{\bf 1}_{\{g_{4}=g_{3}+1\}}
C_{j_k \ldots j_1}\biggl|_{(j_{g_2} j_{g_1})\curvearrowright (\cdot)
(j_{g_4} j_{g_3})\curvearrowright (\cdot),
j_{g_{{}_{1}}}=~j_{g_{{}_{2}}}, j_{g_{{}_{3}}}=~j_{g_{{}_{4}}}}\biggr..
$$

\vspace{3mm}

As a result, Condition 3 of Theorem~2.30 can be replaced by a weaker con\-di\-ti\-on
$$
\lim\limits_{p\to\infty}
\sum\limits_{\stackrel{j_1,\ldots,j_q,\ldots,j_k=0}{{}_{q\ne g_1, g_2, \ldots, g_{2r-1},
g_{2r}}}}^p
\Biggl(\sum\limits_{j_{g_1},j_{g_3},\ldots, j_{g_{2r-1}}=0}^p
C_{j_k\ldots j_1}\biggl|_{j_{g_1}=j_{g_2},\ldots, j_{g_{2r-1}}=j_{g_{2r}}}-\Biggr.
$$
\begin{equation}
\label{drdr1001}
\Biggl.-\frac{1}{2^r} \prod\limits_{l=1}^r {\bf 1}_{\{g_{2l}=g_{2l-1}+1\}}
C_{j_k \ldots j_1}\biggl|_{(j_{g_2} j_{g_1})\curvearrowright (\cdot)
\ldots (j_{g_{2r}} j_{g_{2r-1}})\curvearrowright (\cdot),
j_{g_{{}_{1}}}=~j_{g_{{}_{2}}},\ldots, j_{g_{{}_{2r-1}}}=~j_{g_{{}_{2r}}}
}\biggr.\Biggr)^2=0,
\end{equation}

\noindent
where $r=1, 2,\ldots,[k/2]$.

However, Condition~3 of Theorem~2.30 itself contains 
a way of proving of the condition (\ref{drdr1001}), which is partially 
realized in the proof of Theorems~2.33--2.36 (see below).

In fact, when proving Theorem~2.35 (the case $r=3$
is proved in Theorem~2.36 for $\psi_1(\tau),\ldots,\psi_6(\tau)\equiv 1$), 
we proved the following equality
$$
\lim\limits_{p\to\infty}
\sum\limits_{j_{g_1}=0}^p \sum\limits_{j_{g_{3}}=0}^p
C_{j_k\ldots j_1}\biggl|_{j_{g_1}=j_{g_2},j_{g_{3}}=j_{g_{4}}}=
$$

\vspace{-1mm}
$$
=\frac{1}{4}{\bf 1}_{\{g_{2}=g_{1}+1\}}{\bf 1}_{\{g_{4}=g_{3}+1\}}
C_{j_k \ldots j_1}\biggl|_{(j_{g_2} j_{g_1})\curvearrowright (\cdot)
(j_{g_4} j_{g_3})\curvearrowright (\cdot),
j_{g_{{}_{1}}}=~j_{g_{{}_{2}}}, j_{g_{{}_{3}}}=~j_{g_{{}_{4}}}}\biggr..
$$

\vspace{3mm}

On the other hand, iterative application of (\ref{after80xx}), 
(\ref{after500}) gives
$$
\sum\limits_{j_{g_1}=0}^{\infty} \sum\limits_{j_{g_3}=0}^{\infty}\ldots \sum\limits_{j_{g_{2r-1}}=0}^{\infty}
C_{j_k\ldots j_1}\biggl|_{j_{g_1}=j_{g_2},\ldots, j_{g_{2r-1}}=j_{g_{2r}}}=
$$
$$
=\frac{1}{2^r} \prod\limits_{l=1}^r {\bf 1}_{\{g_{2l}=g_{2l-1}+1\}}
C_{j_k \ldots j_1}\biggl|_{(j_{g_2} j_{g_1})\curvearrowright (\cdot)
\ldots (j_{g_{2r}} j_{g_{2r-1}})\curvearrowright (\cdot),
j_{g_{{}_{1}}}=~j_{g_{{}_{2}}},\ldots, j_{g_{{}_{2r-1}}}=~j_{g_{{}_{2r}}}
},
$$

\vspace{2mm}
\noindent
where $r=1, 2,\ldots,[k/2]$.

Moreover, we have (see (\ref{after906}))
$$
\hbox{\vtop{\offinterlineskip\halign{
\hfil#\hfil\cr
{\rm l.i.m.}\cr
$\stackrel{}{{}_{p\to \infty}}$\cr
}} }\sum_{j_1,\ldots,j_k=0}^{p}
C_{j_k\ldots j_1}
\prod\limits_{s=1}^r{\bf 1}_{\{j_{g_{{}_{2s-1}}}=~j_{g_{{}_{2s}}}\}}
{\bf 1}_{\{i_{g_{{}_{2s-1}}}=~i_{g_{{}_{2s}}}\ne 0\}}
J'[\phi_{j_{q_1}}\ldots \phi_{j_{q_{k-2r}}}]_{T,t}^{(i_{q_1}\ldots i_{q_{k-2r}})}=
$$

\vspace{4mm}
$$
=\hbox{\vtop{\offinterlineskip\halign{
\hfil#\hfil\cr
{\rm l.i.m.}\cr
$\stackrel{}{{}_{p\to \infty}}$\cr
}} }\hspace{-2.5mm}
\sum\limits_{\stackrel{j_1,\ldots,j_q,\ldots,j_k=0}{{}_{q\ne g_1, g_2,\ldots, g_{2r-1}, g_{2r}}}}^p
\sum\limits_{j_{g_1},j_{g_3},\ldots,j_{g_{2r-1}}=0}^p
\hspace{-2.5mm}
C_{j_k \ldots j_1}\biggl|_{j_{g_{{}_{1}}}=~j_{g_{{}_{2}}},\ldots, j_{g_{{}_{2r-1}}}=~j_{g_{{}_{2r}}}
}\biggr.
\times
$$

\vspace{4mm}
$$
\times
\prod\limits_{s=1}^r
{\bf 1}_{\{i_{g_{{}_{2s-1}}}=~i_{g_{{}_{2s}}}\ne 0\}}
J'[\phi_{j_{q_1}}\ldots \phi_{j_{q_{k-2r}}}]_{T,t}^{(i_{q_1}\ldots i_{q_{k-2r}})}=
$$

\vspace{4mm}
$$
=\hbox{\vtop{\offinterlineskip\halign{
\hfil#\hfil\cr
{\rm l.i.m.}\cr
$\stackrel{}{{}_{p\to \infty}}$\cr
}} }\hspace{-2.5mm}
\sum\limits_{\stackrel{j_1,\ldots,j_q,\ldots,j_k=0}{{}_{q\ne g_1, g_2,\ldots, g_{2r-1}, g_{2r}}}}^p
\Biggl(\sum\limits_{j_{g_1},j_{g_3},\ldots,j_{g_{2r-1}}=0}^p
\hspace{-2.5mm}
C_{j_k \ldots j_1}\biggl|_{j_{g_{{}_{1}}}=~j_{g_{{}_{2}}},\ldots, j_{g_{{}_{2r-1}}}=~j_{g_{{}_{2r}}}
}\biggr.-\Biggr.
$$

$$
\Biggl.-\frac{1}{2^r} \prod\limits_{l=1}^r {\bf 1}_{\{g_{2l}=g_{2l-1}+1\}}
C_{j_k \ldots j_1}\biggl|_{(j_{g_2} j_{g_1})\curvearrowright (\cdot)
\ldots (j_{g_{2r}} j_{g_{2r-1}})\curvearrowright (\cdot),
j_{g_{{}_{1}}}=~j_{g_{{}_{2}}},\ldots, j_{g_{{}_{2r-1}}}=~j_{g_{{}_{2r}}}}\Biggr)\times
$$

\vspace{6mm}
$$
\times
\prod\limits_{s=1}^r
{\bf 1}_{\{i_{g_{{}_{2s-1}}}=~i_{g_{{}_{2s}}}\ne 0\}}
J'[\phi_{j_{q_1}}\ldots \phi_{j_{q_{k-2r}}}]_{T,t}^{(i_{q_1}\ldots i_{q_{k-2r}})}+
$$

\vspace{4mm}
$$
+\hbox{\vtop{\offinterlineskip\halign{
\hfil#\hfil\cr
{\rm l.i.m.}\cr
$\stackrel{}{{}_{p\to \infty}}$\cr
}} }
\sum\limits_{\stackrel{j_1,\ldots,j_q,\ldots,j_k=0}{{}_{q\ne g_1, g_2,\ldots, g_{2r-1}, g_{2r}}}}^p
\frac{1}{2^r}
C_{j_k \ldots j_1}\biggl|_{(j_{g_2} j_{g_1})\curvearrowright (\cdot)
\ldots (j_{g_{2r}} j_{g_{2r-1}})\curvearrowright (\cdot),
j_{g_{{}_{1}}}=~j_{g_{{}_{2}}},\ldots, j_{g_{{}_{2r-1}}}=~j_{g_{{}_{2r}}}}\biggr.
\times 
$$

\vspace{2mm}
$$
\times
\prod\limits_{s=1}^r
{\bf 1}_{\{i_{g_{{}_{2s-1}}}=~i_{g_{{}_{2s}}}\ne 0\}}
\prod\limits_{s=1}^r
{\bf 1}_{\{g_{2s}=g_{2s-1}+1\}}
J'[\phi_{j_{q_1}}\ldots \phi_{j_{q_{k-2r}}}]_{T,t}^{(i_{q_1}\ldots i_{q_{k-2r}})}=
$$
$$
=\hbox{\vtop{\offinterlineskip\halign{
\hfil#\hfil\cr
{\rm l.i.m.}\cr
$\stackrel{}{{}_{p\to \infty}}$\cr
}} }\hspace{-2.5mm}
\sum\limits_{\stackrel{j_1,\ldots,j_q,\ldots,j_k=0}{{}_{q\ne g_1, g_2,\ldots, g_{2r-1}, g_{2r}}}}^p
\Biggl(\sum\limits_{j_{g_1},j_{g_3},\ldots,j_{g_{2r-1}}=0}^p
\hspace{-2.5mm}
C_{j_k \ldots j_1}\biggl|_{j_{g_{{}_{1}}}=~j_{g_{{}_{2}}},\ldots, j_{g_{{}_{2r-1}}}=~j_{g_{{}_{2r}}}
}\biggr.-\Biggr.
$$

$$
\Biggl.-\frac{1}{2^r} \prod\limits_{l=1}^r {\bf 1}_{\{g_{2l}=g_{2l-1}+1\}}
C_{j_k \ldots j_1}\biggl|_{(j_{g_2} j_{g_1})\curvearrowright (\cdot)
\ldots (j_{g_{2r}} j_{g_{2r-1}})\curvearrowright (\cdot),
j_{g_{{}_{1}}}=~j_{g_{{}_{2}}},\ldots, j_{g_{{}_{2r-1}}}=~j_{g_{{}_{2r}}}}\Biggr)\times
$$

\vspace{3mm}
$$
\times
\prod\limits_{s=1}^r
{\bf 1}_{\{i_{g_{{}_{2s-1}}}=~i_{g_{{}_{2s}}}\ne 0\}}
J'[\phi_{j_{q_1}}\ldots \phi_{j_{q_{k-2r}}}]_{T,t}^{(i_{q_1}\ldots i_{q_{k-2r}})}+
$$

\vspace{3mm}
\begin{equation}
\label{dydy11}
+
\frac{1}{2^r}\prod\limits_{s=1}^r
{\bf 1}_{\{g_{2s}=g_{2s-1}+1\}}
J[\psi^{(k)}]_{T,t}^{s_r, \ldots, s_1}\ \ \ \hbox{w.~p.~1.}
\end{equation}

\vspace{4mm}

Using (\ref{dydy11}) and the condition (\ref{drdr1001}), we 
obtain (\ref{after801}). This means that we get (\ref{afteru11}).
Thus the expansion (\ref{after1}) is proved.

Analyzing the proof of Theorems~2.30 and 2.12 and taking into account the above arguments,
it is easy to see that the following theorem is true.

{\bf Theorem~2.32}\ \cite{arxiv-10}, \cite{arxiv-11}.\ {\it Assume that
the continuous functions 
$\psi_1(\tau),\ldots,$ $\psi_k(\tau)$ at the interval $[t, T]$ and 
the complete orthonormal system $\{\phi_j(x)\}_{j=0}^{\infty}$
of functions $(\phi_0(x)=1/\sqrt{T-t})$ 
in the space $L_2([t, T])$ are such that the following 
condition 

\vspace{-3mm}
$$
\lim\limits_{p_1,\ldots,p_k\to\infty}~
\sum\limits_{j_1=0}^{p_1}\ldots \sum\limits_{j_q=0}^{p_q}\ldots \sum\limits_{j_k=0}^{p_k}~
\biggl|_{q\ne g_1, g_2, \ldots, g_{2r-1},g_{2r}}\times
$$

\vspace{2mm}
$$
\times
\Biggl(~\sum\limits_{j_{g_1}=0}^{\min\{p_{g_1}, p_{g_2}\}} \sum\limits_{j_{g_3}=0}^{\min\{p_{g_3}, p_{g_4}\}}\ldots \Biggr.
\sum\limits_{j_{g_{2r-1}}=0}^{\min\{p_{g_{2r-1}}, p_{g_{2r}}\}}
C_{j_k\ldots j_1}\biggl|_{j_{g_1}=j_{g_2},\ldots, j_{g_{2r-1}}=j_{g_{2r}}}-
$$

\vspace{-2mm}
\begin{equation}
\label{drdr1001xyx}
\Biggl.-\frac{1}{2^r} \prod\limits_{l=1}^r {\bf 1}_{\{g_{2l}=g_{2l-1}+1\}}
C_{j_k \ldots j_1}\biggl|_{(j_{g_2} j_{g_1})\curvearrowright (\cdot)
\ldots (j_{g_{2r}} j_{g_{2r-1}})\curvearrowright (\cdot),
j_{g_{{}_{1}}}=~j_{g_{{}_{2}}},\ldots, j_{g_{{}_{2r-1}}}=~j_{g_{{}_{2r}}}
}\biggr.\Biggr)^2=0
\end{equation}

\vspace{4mm}
\noindent
is satisfied for all $r=1, 2,\ldots,[k/2]$
and for all possible $g_1,g_2,\ldots,g_{2r-1},g_{2r}$ {\rm (}see {\rm (\ref{leto5007after}))}.
Then$,$ for the iterated Stratonovich sto\-chas\-tic integral 
of arbitrary multiplicity $k$
\newpage
\noindent
$$
J^{*}[\psi^{(k)}]_{T,t}^{(i_1\ldots i_k)}=
{\int\limits_t^{*}}^T
\psi_k(t_k) \ldots 
{\int\limits_t^{*}}^{t_{2}}
\psi_1(t_1) d{\bf w}_{t_1}^{(i_1)}\ldots
d{\bf w}_{t_k}^{(i_k)}
$$
the following 
expansion 
$$
J^{*}[\psi^{(k)}]_{T,t}^{(i_1\ldots i_k)}=
\hbox{\vtop{\offinterlineskip\halign{
\hfil#\hfil\cr
{\rm l.i.m.}\cr
$\stackrel{}{{}_{p_1,\ldots,p_k\to \infty}}$\cr
}} }
\sum\limits_{j_1=0}^{p_1}\ldots\sum\limits_{j_k=0}^{p_k}
C_{j_k \ldots j_1}\prod\limits_{l=1}^k \zeta_{j_l}^{(i_l)}
$$

\noindent
that converges in the mean-square sense is valid, where 
$$
C_{j_k \ldots j_1}=\int\limits_t^T\psi_k(t_k)\phi_{j_k}(t_k)\ldots
\int\limits_t^{t_2}
\psi_1(t_1)\phi_{j_1}(t_1)
dt_1\ldots dt_k
$$
is the Fourier coefficient, 
${\rm l.i.m.}$ is a limit in the mean-square sense,
$i_1, \ldots, i_k=0, 1,\ldots,m,$
$$
\zeta_{j}^{(i)}=
\int\limits_t^T \phi_{j}(\tau) d{\bf w}_{\tau}^{(i)}
$$ 
are independent standard Gaussian random variables for various 
$i$ or $j$ {\rm (}when $i\ne 0${\rm ),}
${\bf w}_{\tau}^{(i)}$ 
$(i=1,\ldots,m)$ are independent standard Wiener processes$,$
${\bf w}_{\tau}^{(0)}=\tau.$}

\section{Expansion of Iterated Stratonovich Stochastic Integrals
of Multiplicity 3. The Case $p_1=p_2=p_3\to \infty$ and 
Continuously Differentiable 
Weight Functions $\psi_1(\tau),$ $\psi_2(\tau),$ $\psi_3(\tau)$ 
(The Cases of Legendre 
Polynomials and Trigonometric Functions)}

In this section, we present a simple proof of Theorem~2.8  
based on Theorem~2.30. In this case, the conditions of Theorem~2.8 will be weakened.

First, consider the following equalities 
\begin{equation}
\label{after1400}
~~~~\frac{1}{2}
\int\limits_{t_1}^{t_2} \Phi_1(\tau)\Phi_2(\tau)d\tau
=\sum_{j=0}^{\infty}
\int\limits_{t_1}^{t_2}
\Phi_2(\tau)\phi_{j}(\tau)\int\limits_{t_1}^{\tau}
\Phi_1(\theta)\phi_{j}(\theta)d\theta d\tau,
\end{equation}
\begin{equation}
\label{after1401}
~~~~\frac{1}{2}
\int\limits_{t_1}^{t_2} \Phi_1(\tau)\Phi_2(\tau)d\tau
=\sum_{j=0}^{\infty}
\int\limits_{t_1}^{t_2}
\Phi_1(\theta)\phi_{j}(\theta)\int\limits_{\theta}^{t_2}
\Phi_2(\tau)\phi_{j}(\tau)d\tau d\theta
\end{equation}
that will be used further, where $t\le t_1<t_2\le T,$  $\Phi_1(\tau), \Phi_2(\tau)\in L_2([t,T]),$
$\{\phi_j(x)\}_{j=0}^{\infty}$ is the same 
as in the conditions of Theorem~2.8.

The equality (\ref{after1400}) follows from (\ref{start1000}).
Using (\ref{after1400}) and Fubini's Theorem, we get (\ref{after1401}).

{\bf Theorem 2.33}\ \cite{arxiv-5}, 
\cite{arxiv-10}, \cite{arxiv-11}, \cite{new-art-1xxy}.\
{\it Suppose that 
$\{\phi_j(x)\}_{j=0}^{\infty}$ is a complete orthonormal system of 
Legendre polynomials or trigonometric functions in the space $L_2([t, T]).$
Furthermore, let $\psi_1(\tau), \psi_2(\tau), \psi_3(\tau)$ are continuously dif\-ferentiable 
nonrandom functions on $[t, T].$ 
Then$,$ for the 
iterated Stratonovich stochastic integral of third multiplicity
\begin{equation}
\label{after1500}
~~~~~~J^{*}[\psi^{(3)}]_{T,t}={\int\limits_t^{*}}^T\psi_3(t_3)
{\int\limits_t^{*}}^{t_3}\psi_2(t_2)
{\int\limits_t^{*}}^{t_2}\psi_1(t_1)
d{\bf w}_{t_1}^{(i_1)}
d{\bf w}_{t_2}^{(i_2)}d{\bf w}_{t_3}^{(i_3)}
\end{equation}
the following 
expansion 
$$
J^{*}[\psi^{(3)}]_{T,t}
=\hbox{\vtop{\offinterlineskip\halign{
\hfil#\hfil\cr
{\rm l.i.m.}\cr
$\stackrel{}{{}_{p\to \infty}}$\cr
}} }
\sum\limits_{j_1, j_2, j_3=0}^{p}
C_{j_3 j_2 j_1}\zeta_{j_1}^{(i_1)}\zeta_{j_2}^{(i_2)}\zeta_{j_3}^{(i_3)}
$$
that converges in the mean-square sense is valid, where
$i_1, i_2, i_3=0, 1,\ldots,m,$
$$
C_{j_3 j_2 j_1}=\int\limits_t^T\psi_3(t_3)\phi_{j_3}(t_3)
\int\limits_t^{t_3}\psi_2(t_2)\phi_{j_2}(t_2)
\int\limits_t^{t_2}\psi_1(t_1)\phi_{j_1}(t_1)dt_1dt_2dt_3
$$
and
$$
\zeta_{j}^{(i)}=
\int\limits_t^T \phi_{j}(s) d{\bf w}_s^{(i)}
$$ 
are independent standard Gaussian random variables for various 
$i$ or $j$
{\rm (}in the case when $i\ne 0${\rm ),}
${\bf w}_{\tau}^{(i)}$
$(i=1,\ldots,m)$ are independent standard Wiener processes$,$
${\bf w}_{\tau}^{(0)}=\tau.$}

\vspace{2mm}

{\bf Proof.}\ As noted in Remark~2.4, Conditions 1 and 2
of Theorem~{\rm 2.30} are satisfied for complete
orthonormal systems of Legendre polynomials 
and trigonometric functions in the space
$L_2([t, T]).$ Let us verify Condition 3 of Theorem~2.30 for
the iterated Stratonovich stochastic integral (\ref{after1500}). 
Thus, we have to check the following conditions
\begin{equation}
\label{after1600}
\lim\limits_{p\to \infty}
\sum_{j_3=0}^p\left(\sum_{j_1=p+1}^{\infty}
C_{j_3 j_1 j_1}\right)^2=0,
\end{equation}
\begin{equation}
\label{after1601}
\lim\limits_{p\to \infty}
\sum_{j_1=0}^p\left(\sum_{j_3=p+1}^{\infty}
C_{j_3 j_3 j_1}\right)^2=0,
\end{equation}
\begin{equation}
\label{after1602}
\lim\limits_{p\to \infty}
\sum_{j_2=0}^p\left(\sum_{j_1=p+1}^{\infty}
C_{j_1 j_2 j_1}\right)^2=0.
\end{equation}

We have
$$
\sum_{j_3=0}^p\left(\sum_{j_1=p+1}^{\infty}
C_{j_3 j_1 j_1}\right)^2=
$$
\begin{equation}
\label{after1800}
=
\sum_{j_3=0}^p\left(\sum_{j_1=p+1}^{\infty}
\int\limits_t^T\psi_3(t_3)\phi_{j_3}(t_3)
\int\limits_t^{t_3}\psi_2(t_2)\phi_{j_1}(t_2)
\int\limits_t^{t_2}\psi_1(t_1)\phi_{j_1}(t_1)dt_1dt_2dt_3\right)^2=
\end{equation}
\begin{equation}
\label{after1801}
=
\sum_{j_3=0}^p\left(
\int\limits_t^T\psi_3(t_3)\phi_{j_3}(t_3)
\sum_{j_1=p+1}^{\infty}\int\limits_t^{t_3}\psi_2(t_2)\phi_{j_1}(t_2)
\int\limits_t^{t_2}\psi_1(t_1)\phi_{j_1}(t_1)dt_1dt_2dt_3\right)^2\le
\end{equation}
\begin{equation}
\label{after1802}
\le
\sum_{j_3=0}^{\infty}\left(
\int\limits_t^T\psi_3(t_3)\phi_{j_3}(t_3)
\sum_{j_1=p+1}^{\infty}\int\limits_t^{t_3}\psi_2(t_2)\phi_{j_1}(t_2)
\int\limits_t^{t_2}\psi_1(t_1)\phi_{j_1}(t_1)dt_1dt_2dt_3\right)^2=
\end{equation}
\begin{equation}
\label{after1803}
~~~~~=\int\limits_t^T\psi_3^2(t_3)
\left(\sum_{j_1=p+1}^{\infty}\int\limits_t^{t_3}\psi_2(t_2)\phi_{j_1}(t_2)
\int\limits_t^{t_2}\psi_1(t_1)\phi_{j_1}(t_1)dt_1dt_2\right)^2 dt_3\le
\end{equation}
\begin{equation}
\label{after1804}
\le \frac{K}{p^2}\ \to\  0
\end{equation}

\noindent
if $p\to\infty,$ where constant $K$ does not depend on $p.$

Note that the transition from (\ref{after1800}) to (\ref{after1801}) is based on 
the estimate (\ref{agent1516}) for the polynomial case and its analogue
for the
trigonometric case, 
the transition from (\ref{after1802}) to (\ref{after1803}) is based on the Parseval equality, 
and the transition from (\ref{after1803}) to (\ref{after1804}) 
is also based on the estimate (\ref{agent1516})
and its analogue
for the
trigonometric case.

By analogy with the previous case we have 
$$
\sum_{j_1=0}^p\left(\sum_{j_3=p+1}^{\infty}
C_{j_3 j_3 j_1}\right)^2=
$$
$$
=
\sum_{j_1=0}^p\left(\sum_{j_3=p+1}^{\infty}
\int\limits_t^T\psi_3(t_3)\phi_{j_3}(t_3)
\int\limits_t^{t_3}\psi_2(t_2)\phi_{j_3}(t_2)
\int\limits_t^{t_2}\psi_1(t_1)\phi_{j_1}(t_1)dt_1dt_2dt_3\right)^2=
$$
\begin{equation}
\label{after1900}
=
\sum_{j_1=0}^p\left(\sum_{j_3=p+1}^{\infty}
\int\limits_t^T \psi_1(t_1)\phi_{j_1}(t_1) \int\limits_{t_1}^{T} \psi_2(t_2)\phi_{j_3}(t_2)
\int\limits_{t_2}^T\psi_3(t_3)\phi_{j_3}(t_3)
dt_3dt_2dt_1\right)^2=
\end{equation}
\begin{equation}
\label{after1901}
=
\sum_{j_1=0}^p\left(
\int\limits_t^T \psi_1(t_1)\phi_{j_1}(t_1) 
\sum_{j_3=p+1}^{\infty}\int\limits_{t_1}^{T} \psi_2(t_2)\phi_{j_3}(t_2)
\int\limits_{t_2}^T\psi_3(t_3)\phi_{j_3}(t_3)
dt_3dt_2dt_1\right)^2\le
\end{equation}
$$
\le
\sum_{j_1=0}^{\infty}\left(
\int\limits_t^T \psi_1(t_1)\phi_{j_1}(t_1) 
\sum_{j_3=p+1}^{\infty}\int\limits_{t_1}^{T} \psi_2(t_2)\phi_{j_3}(t_2)
\int\limits_{t_2}^T\psi_3(t_3)\phi_{j_3}(t_3)
dt_3dt_2dt_1\right)^2=
$$
\begin{equation}
\label{after1902}
~~~~~=
\int\limits_t^T \psi_1^2(t_1)\left(
\sum_{j_3=p+1}^{\infty}\int\limits_{t_1}^{T} \psi_2(t_2)\phi_{j_3}(t_2)
\int\limits_{t_2}^T\psi_1(t_3)\phi_{j_3}(t_3)
dt_3dt_2\right)^2 dt_1\le
\end{equation}
\begin{equation}
\label{after1903}
\le \frac{K}{p^2}\ \to\  0
\end{equation}

\vspace{1mm}
\noindent
if $p\to\infty,$ where constant $K$ is independent of $p.$

The transition from (\ref{after1900}) to (\ref{after1901}) is based on 
an analogue of the estimate (\ref{agent1516}) 
for the value
$$
\left|\sum_{j_3=p+1}^{\infty}\int\limits_{t_1}^{T} \psi_2(t_2)\phi_{j_3}(t_2)
\int\limits_{t_2}^T\psi_3(t_3)\phi_{j_3}(t_3)
dt_3dt_2\right|
$$
for the polynomial and
trigonometric cases, 
the transition from (\ref{after1902}) to (\ref{after1903})
is also based on the mentioned analogue of the estimate (\ref{agent1516}).

Further, we have 
$$
\sum_{j_2=0}^p\left(\sum_{j_1=p+1}^{\infty}
C_{j_1 j_2 j_1}\right)^2=
$$
$$
=
\sum_{j_2=0}^p\left(\sum_{j_1=p+1}^{\infty}
\int\limits_t^T\psi_3(t_3)\phi_{j_1}(t_3)
\int\limits_t^{t_3}\psi_2(t_2)\phi_{j_2}(t_2)
\int\limits_t^{t_2}\psi_1(t_1)\phi_{j_1}(t_1)dt_1dt_2dt_3\right)^2=
$$
\begin{equation}
\label{after1920}
=
\sum_{j_2=0}^p\left(\sum_{j_1=p+1}^{\infty}
\int\limits_t^T \psi_2(t_2)\phi_{j_2}(t_2) \int\limits_{t}^{t_2} \psi_1(t_1)\phi_{j_1}(t_1)
dt_1 \int\limits_{t_2}^T\psi_3(t_3)\phi_{j_1}(t_3)
dt_3dt_2\right)^2=
\end{equation}
\begin{equation}
\label{after1921}
=
\sum_{j_2=0}^p\left(
\int\limits_t^T \psi_2(t_2)\phi_{j_2}(t_2) 
\sum_{j_1=p+1}^{\infty}\int\limits_{t}^{t_2} \psi_1(t_1)\phi_{j_1}(t_1)
dt_1 \int\limits_{t_2}^T\psi_3(t_3)\phi_{j_1}(t_3)
dt_3dt_2\right)^2\le
\end{equation}
$$
\le
\sum_{j_2=0}^{\infty}\left(
\int\limits_t^T \psi_2(t_2)\phi_{j_2}(t_2) \sum_{j_1=p+1}^{\infty}
\int\limits_{t}^{t_2} \psi_1(t_1)\phi_{j_1}(t_1)
dt_1 \int\limits_{t_2}^T\psi_3(t_3)\phi_{j_1}(t_3)
dt_3dt_2\right)^2=
$$
\begin{equation}
\label{after1922}
~~~~~=
\int\limits_t^T \psi_2^2(t_2)
\left(\sum_{j_1=p+1}^{\infty}\int\limits_{t}^{t_2} \psi_1(t_1)\phi_{j_1}(t_1)
dt_1 \int\limits_{t_2}^T\psi_3(t_3)\phi_{j_1}(t_3)
dt_3\right)^2 dt_2.
\end{equation}

\vspace{1mm}

The transition from (\ref{after1920}) to (\ref{after1921}) 
is based on the estimate (\ref{103xx}) and its obvious analogue
for the trigonometric case.
However, the estimate (\ref{103xx}) cannot be used to estimate
the right-hand side of (\ref{after1922}), 
since we get the divergent integral.
For this reason, we will obtain a new estimate based on the relation
(\ref{otit6000x}).

From (\ref{ogo23}) and the estimate $\left| P_j(y) \right|\le 1$, $y\in [-1, 1]$
we obtain

\vspace{-5mm}
\begin{equation}
\label{after5000}
\left|P_{j}(y)\right|=\left|P_{j}(y)\right|^{\varepsilon}
\cdot \left|P_{j}(y)\right|^{1-\varepsilon}\le \left|P_{j}(y)\right|^{1-\varepsilon}
< \frac{C}{j^{1/2-\varepsilon/2} (1-y^2)^{1/4-\varepsilon/4}},
\end{equation}

\vspace{1mm}
\noindent
where $y\in (-1, 1),$ $j\in{\bf N},$ $\varepsilon\in (0,1)$ is an arbitrary
small positive real number.

Combining (\ref{otit6000x}) and (\ref{after5000}), we have the following estimate

\vspace{-5mm}
\begin{equation}
\label{after1940}
~~~~~~\left|
\int\limits_t^s\psi_1(\tau)\phi_{j_1}(\tau)d\tau
\right| <
\frac{C}{(j_1)^{1-\varepsilon/2}}\Biggl(\frac{1}{(1-z^2(s))^{1/4-\varepsilon/4}}+1\Biggr),
\end{equation}

\vspace{1mm}
\noindent
where $j_1\in{\bf N},$ $s\in (t, T),$
$z(s)$ is defined by (\ref{zz1}), constant $C$ does not depend on $j_1.$

Similarly to (\ref{after1940}) we obtain 

\vspace{-5mm}
\begin{equation}
\label{after1941}
~~~~~~\left|\int\limits_s^T\psi_3(\tau)\phi_{j_1}(\tau)d\tau
\right| <
\frac{C}{(j_1)^{1-\varepsilon/2}}\Biggl(\frac{1}{(1-z^2(s))^{1/4-\varepsilon/4}}+1\Biggr),
\end{equation}

\vspace{1mm}
\noindent
where $j_1\in{\bf N},$ $s\in (t, T),$ constant $C$ does not depend on $j_1.$

Combining (\ref{101xx}) and (\ref{after1941}), we have

\vspace{-4mm}
$$
\left|
\int\limits_t^s\psi_1(\tau)\phi_{j_1}(\tau)d\tau
\int\limits_s^T\psi_3(\tau)\phi_{j_1}(\tau)d\tau\right|
<
$$
\begin{equation}
\label{after4000}
~~~~<\frac{L}{(j_1)^{2-\varepsilon/2}}\Biggl(\frac{1}{(1-z^2(s))^{1/4-\varepsilon/4}}+1\Biggr)
\Biggl(\frac{1}{(1-z^2(s))^{1/4}}+1\Biggr),
\end{equation}

\vspace{1mm}
\noindent
where $j_1\in{\bf N},$ $s\in (t, T),$ 
$z(s)$ is defined by (\ref{zz1}), constant $L$ does not depend on $j_1.$

Observe that

\vspace{-4mm}
\begin{equation}
\label{after1944}
\sum\limits_{j_1=p+1}^{\infty}\frac{1}{(j_1)^{2-\varepsilon/2}}
\le \int\limits_{p}^{\infty}\frac{dx}{x^{2-\varepsilon/2}}=
\frac{1}{(1-\varepsilon/2)p^{1-\varepsilon/2}}.
\end{equation}

\vspace{1mm}

Applying (\ref{after4000}) and (\ref{after1944})
to estimate the right-hand side of (\ref{after1922}) gives

\vspace{-4mm}
\begin{equation}
\label{after5001}
\sum_{j_2=0}^p\left(\sum_{j_1=p+1}^{\infty}
C_{j_1 j_2 j_1}\right)^2\le \frac{K}{p^{2-\varepsilon}}\ \to\  0
\end{equation}

\vspace{1mm}
\noindent
if $p\to\infty,$ where 
$\varepsilon$ is an arbitrary
small positive real number,
constant $K$ is independent of $p$.

The estimation of the right-hand side of (\ref{after1922})
for the trigonometric case is carried out using the estimates
(\ref{2017x11}), (\ref{2017x12}). At that we obtain the
estimate (\ref{after5001}) with $\varepsilon=0.$
Theorem~2.33 is proved.

\section{Expansion of Iterated Stratonovich Stochastic Integrals
of Multiplicity 4. The Case $p_1=\ldots =p_4\to \infty$ and 
Continuously Differentiable 
Weight Functions $\psi_1(\tau),$ $\ldots,$ $\psi_4(\tau)$ 
(The Cases of Legendre 
Polynomials and Trigonometric Functions)}

{\bf Theorem 2.34}\ \cite{arxiv-5}, 
\cite{arxiv-10}, \cite{arxiv-11}, \cite{new-art-1xxy}.\
{\it Suppose that 
$\{\phi_j(x)\}_{j=0}^{\infty}$ is a complete orthonormal system of 
Legendre polynomials or trigonometric functions in the space $L_2([t, T]).$
Furthermore, let $\psi_1(\tau), \ldots, \psi_4(\tau)$ are continuously dif\-ferentiable 
nonrandom functions on $[t, T].$ 
Then$,$ for the 
iterated Stratonovich stochastic integral of fourth multiplicity
\begin{equation}
\label{after2500}
J^{*}[\psi^{(4)}]_{T,t}={\int\limits_t^{*}}^T\psi_4(t_4)
{\int\limits_t^{*}}^{t_4}\psi_3(t_3)
{\int\limits_t^{*}}^{t_3}\psi_2(t_2)
{\int\limits_t^{*}}^{t_2}\psi_1(t_1)
d{\bf w}_{t_1}^{(i_1)}
d{\bf w}_{t_2}^{(i_2)}d{\bf w}_{t_3}^{(i_3)}d{\bf w}_{t_4}^{(i_4)}
\end{equation}
the following 
expansion 
$$
J^{*}[\psi^{(4)}]_{T,t}
=\hbox{\vtop{\offinterlineskip\halign{
\hfil#\hfil\cr
{\rm l.i.m.}\cr
$\stackrel{}{{}_{p\to \infty}}$\cr
}} }
\sum\limits_{j_1, j_2, j_3, j_4=0}^{p}
C_{j_4 j_3 j_2 j_1}\zeta_{j_1}^{(i_1)}\zeta_{j_2}^{(i_2)}\zeta_{j_3}^{(i_3)}
\zeta_{j_4}^{(i_4)}
$$
that converges in the mean-square sense is valid, where
$i_1, i_2, i_3, i_4=0, 1,\ldots,m,$
$$
C_{j_4 j_3 j_2 j_1}=
\int\limits_t^T\psi_4(t_4)\phi_{j_4}(t_4)
\int\limits_t^{t_4}\psi_3(t_3)\phi_{j_3}(t_3)
\int\limits_t^{t_3}\psi_2(t_2)\phi_{j_2}(t_2)
\int\limits_t^{t_2}\psi_1(t_1)\phi_{j_1}(t_1)dt_1\times
$$
$$
\times
dt_2dt_3dt_4
$$
and
$$
\zeta_{j}^{(i)}=
\int\limits_t^T \phi_{j}(s) d{\bf w}_s^{(i)}
$$ 
are independent standard Gaussian random variables for various 
$i$ or $j$
{\rm (}when $i\ne 0${\rm ),}
${\bf w}_{\tau}^{(i)}$ 
$(i=1,\ldots,m)$ are independent standard Wiener processes$,$
${\bf w}_{\tau}^{(0)}=\tau.$}

\vspace{2mm}

{\bf Proof.}\ As noted in Remark~2.4, Conditions 1 and 2
of Theorem~{\rm 2.30} are satisfied for complete
orthonormal systems of Legendre polynomials 
and trigonometric functions in the space
$L_2([t, T]).$ Let us verify Condition 3 of Theorem~2.30 for
the iterated Stratonovich stochastic integral (\ref{after2500}). 
Thus, we have to check the following conditions
\begin{equation}
\label{after2501}
\lim\limits_{p\to \infty}
\sum_{j_3,j_4=0}^p\left(\sum_{j_1=p+1}^{\infty}
C_{j_4 j_3 j_1 j_1}\right)^2=0,
\end{equation}
\begin{equation}
\label{after2502}
\lim\limits_{p\to \infty}
\sum_{j_2,j_4=0}^p\left(\sum_{j_1=p+1}^{\infty}
C_{j_4 j_1 j_2 j_1}\right)^2=0,
\end{equation}
\begin{equation}
\label{after2503}
\lim\limits_{p\to \infty}
\sum_{j_2,j_3=0}^p\left(\sum_{j_1=p+1}^{\infty}
C_{j_1 j_3 j_2 j_1}\right)^2=0,
\end{equation}
\begin{equation}
\label{after2504}
\lim\limits_{p\to \infty}
\sum_{j_1,j_4=0}^p\left(\sum_{j_2=p+1}^{\infty}
C_{j_4 j_2 j_2 j_1}\right)^2=0,
\end{equation}
\begin{equation}
\label{after2505}
\lim\limits_{p\to \infty}
\sum_{j_1,j_3=0}^p\left(\sum_{j_2=p+1}^{\infty}
C_{j_2 j_3 j_2 j_1}\right)^2=0,
\end{equation}
\begin{equation}
\label{after2506}
\lim\limits_{p\to \infty}
\sum_{j_1,j_2=0}^p\left(\sum_{j_3=p+1}^{\infty}
C_{j_3 j_3 j_2 j_1}\right)^2=0,
\end{equation}
\begin{equation}
\label{after2508}
\lim\limits_{p\to \infty}
\left(\sum_{j_2=p+1}^{\infty}\sum_{j_1=p+1}^{\infty}
C_{j_2 j_1 j_2 j_1}\right)^2=0,
\end{equation}
\begin{equation}
\label{after2509}
\lim\limits_{p\to \infty}
\left(\sum_{j_2=p+1}^{\infty}\sum_{j_1=p+1}^{\infty}
C_{j_1 j_2 j_2 j_1}\right)^2=0,
\end{equation}
\begin{equation}
\label{after2507}
\lim\limits_{p\to \infty}
\left(\sum_{j_3=p+1}^{\infty}\sum_{j_1=p+1}^{\infty}
C_{j_3 j_3 j_1 j_1}\right)^2=0,
\end{equation}
\begin{equation}
\label{after1602xxxx}
\lim\limits_{p\to \infty}
\left(\sum_{j_3=p+1}^{\infty}
C_{j_3 j_3 j_1 j_1}\biggl|_{(j_{1} j_{1})\curvearrowright (\cdot)}\biggr.\right)^2=0,
\end{equation}
\begin{equation}
\label{after1602xxxy}
\lim\limits_{p\to \infty}
\left(\sum_{j_1=p+1}^{\infty}
C_{j_3 j_3 j_1 j_1}\biggl|_{(j_{3} j_{3})\curvearrowright (\cdot)}\biggr.\right)^2=0,
\end{equation}
\begin{equation}
\label{after1602xxxz}
\lim\limits_{p\to \infty}
\left(\sum_{j_1=p+1}^{\infty}
C_{j_1 j_2 j_2 j_1}\biggl|_{(j_{2} j_{2})\curvearrowright (\cdot)}\biggr.\right)^2=0,
\end{equation}

\vspace{2mm}
\noindent
where in (\ref{after1602xxxx})--(\ref{after1602xxxz}) we use the notation (\ref{after900}).

Applying arguments similar to those we used in the proof of Theorem~2.33, we obtain
for (\ref{after2501})
$$
\sum_{j_3,j_4=0}^p\left(\sum_{j_1=p+1}^{\infty}
C_{j_4 j_3 j_1 j_1}\right)^2=
\sum_{j_3,j_4=0}^p\left(\sum_{j_1=p+1}^{\infty}
\int\limits_t^T\psi_4(t_4)\phi_{j_4}(t_4)
\int\limits_t^{t_4}\psi_3(t_3)\phi_{j_3}(t_3)\times\right.
$$
\begin{equation}
\label{after5010}
\left.\times
\int\limits_t^{t_3}\psi_2(t_2)\phi_{j_1}(t_2)
\int\limits_t^{t_2}\psi_1(t_1)\phi_{j_1}(t_1)dt_1dt_2dt_3dt_4\right)^2=
\end{equation}
$$
=
\sum_{j_3,j_4=0}^p\left(
\int\limits_t^T\psi_4(t_4)\phi_{j_4}(t_4)
\int\limits_t^{t_4}\psi_3(t_3)\phi_{j_3}(t_3)\times\right.
$$
\begin{equation}
\label{after5011}
\left.~~~~~~~~~~~~\times
\sum_{j_1=p+1}^{\infty}\int\limits_t^{t_3}\psi_2(t_2)\phi_{j_1}(t_2)
\int\limits_t^{t_2}\psi_1(t_1)\phi_{j_1}(t_1)dt_1dt_2dt_3dt_4\right)^2\le
\end{equation}
$$
\le
\sum_{j_3,j_4=0}^{\infty}\left(
\int\limits_t^T\psi_4(t_4)\phi_{j_4}(t_4)
\int\limits_t^{t_4}\psi_3(t_3)\phi_{j_3}(t_3)\times\right.
$$
\begin{equation}
\label{after5002}
\left.~~~~~~~~~~\times
\sum_{j_1=p+1}^{\infty}\int\limits_t^{t_3}\psi_2(t_2)\phi_{j_1}(t_2)
\int\limits_t^{t_2}\psi_1(t_1)\phi_{j_1}(t_1)dt_1dt_2dt_3dt_4\right)^2=
\end{equation}

\vspace{-1mm}
$$
=
\int\limits_{[t, T]^2}{\bf 1}_{\{t_3<t_4\}}
\psi_4^2(t_4)
\psi_3^2(t_3)\times
$$
\begin{equation}
\label{after5003}
~~~~~~~~~~~\times\left(
\sum_{j_1=p+1}^{\infty}\int\limits_t^{t_3}\psi_2(t_2)\phi_{j_1}(t_2)
\int\limits_t^{t_2}\psi_1(t_1)\phi_{j_1}(t_1)dt_1dt_2\right)^2 dt_3dt_4\le
\end{equation}
\begin{equation}
\label{after5004}
\le \frac{K}{p^2}\ \to\  0
\end{equation}

\noindent
if $p\to\infty,$ where constant $K$ is independent of $p.$

Note that the transition from (\ref{after5010}) to (\ref{after5011}) is based on 
the estimate (\ref{agent1516}) for the polynomial case and its analogue
for the
trigonometric case, 
the transition from (\ref{after5002}) to (\ref{after5003}) is based on the Parseval equality, 
and the transition from (\ref{after5003}) to (\ref{after5004}) 
is also based on the estimate (\ref{agent1516})
and its analogue
for the
trigonometric case.

Further, we have for (\ref{after2502})
$$
\sum_{j_2,j_4=0}^p\left(\sum_{j_1=p+1}^{\infty}
C_{j_4 j_1 j_2 j_1}\right)^2=
\sum_{j_2,j_4=0}^p\left(\sum_{j_1=p+1}^{\infty}
\int\limits_t^T\psi_4(t_4)\phi_{j_4}(t_4)
\int\limits_t^{t_4}\psi_3(t_3)\phi_{j_1}(t_3)\times\right.
$$
\begin{equation}
\label{after98}
\left.\times
\int\limits_t^{t_3}\psi_2(t_2)\phi_{j_2}(t_2)
\int\limits_t^{t_2}\psi_1(t_1)\phi_{j_1}(t_1)dt_1dt_2dt_3dt_4\right)^2=
\end{equation}
$$
=
\sum_{j_2,j_4=0}^p\left(\sum_{j_1=p+1}^{\infty}
\int\limits_t^T\psi_4(t_4)\phi_{j_4}(t_4)
\int\limits_t^{t_4}\psi_2(t_2)\phi_{j_2}(t_2)\times\right.
$$
\begin{equation}
\label{after99}
\left.\times
\int\limits_t^{t_2}\psi_1(t_1)\phi_{j_1}(t_1)dt_1
\int\limits_{t_2}^{t_4}\psi_3(t_3)\phi_{j_1}(t_3)dt_3dt_2dt_4\right)^2=
\end{equation}
$$
=
\sum_{j_2,j_4=0}^p\left(
\int\limits_t^T\psi_4(t_4)\phi_{j_4}(t_4)
\int\limits_t^{t_4}\psi_2(t_2)\phi_{j_2}(t_2)\times\right.
$$
$$
\left.\times
\sum_{j_1=p+1}^{\infty}\int\limits_t^{t_2}\psi_1(t_1)\phi_{j_1}(t_1)dt_1
\int\limits_{t_2}^{t_4}\psi_3(t_3)\phi_{j_1}(t_3)dt_3dt_2dt_4\right)^2\le
$$
$$
\le
\sum_{j_2,j_4=0}^{\infty}\left(
\int\limits_t^T\psi_4(t_4)\phi_{j_4}(t_4)
\int\limits_t^{t_4}\psi_2(t_2)\phi_{j_2}(t_2)\times\right.
$$
$$
\left.\times
\sum_{j_1=p+1}^{\infty}\int\limits_t^{t_2}\psi_1(t_1)\phi_{j_1}(t_1)dt_1
\int\limits_{t_2}^{t_4}\psi_3(t_3)\phi_{j_1}(t_3)dt_3dt_2dt_4\right)^2=
$$

\vspace{-1mm}
$$
=
\int\limits_{[t, T]^2}{\bf 1}_{\{t_2<t_4\}}
\psi_4^2(t_4)
\psi_2^2(t_2)\times
$$
$$
\times
\left(\sum_{j_1=p+1}^{\infty}\int\limits_t^{t_2}\psi_1(t_1)\phi_{j_1}(t_1)dt_1
\int\limits_{t_2}^{t_4}\psi_3(t_3)\phi_{j_1}(t_3)dt_3\right)^2dt_2dt_4\le 
$$
\begin{equation}
\label{after6009}
\le \frac{K}{p^{2-\varepsilon}}\ \to\  0
\end{equation}

\noindent
if $p\to\infty,$ where $\varepsilon$ is an arbitrary small positive real number
for the polynomial case and $\varepsilon=0$ for the 
trigonometric case, constant $K$ does not depend on $p.$

The relation (\ref{after6009}) was obtained by the same method as (\ref{after5004}). 
Note that in obtaining (\ref{after6009}) 
we used the estimates (\ref{101oh}) and (\ref{after1940}) for the polynomial 
case and (\ref{203}) and (\ref{2017x11}) for the trigonometric case. 
We also used the integration order replacement in the iterated Riemann integrals
(see (\ref{after98}), (\ref{after99})).

Repeating the previous steps for (\ref{after2503}) and (\ref{after2504}), we get
$$
\sum_{j_2,j_3=0}^p\left(\sum_{j_1=p+1}^{\infty}
C_{j_1 j_3 j_2 j_1}\right)^2=
\sum_{j_2,j_3=0}^p\left(\sum_{j_1=p+1}^{\infty}
\int\limits_t^T\psi_4(t_4)\phi_{j_1}(t_4)
\int\limits_t^{t_4}\psi_3(t_3)\phi_{j_3}(t_3)\times\right.
$$
$$
\left.\times
\int\limits_t^{t_3}\psi_2(t_2)\phi_{j_2}(t_2)
\int\limits_t^{t_2}\psi_1(t_1)\phi_{j_1}(t_1)dt_1dt_2dt_3dt_4\right)^2=
$$
$$
=
\sum_{j_2,j_3=0}^p\left(\sum_{j_1=p+1}^{\infty}
\int\limits_t^T\psi_3(t_3)\phi_{j_3}(t_3)
\int\limits_t^{t_3}\psi_2(t_2)\phi_{j_2}(t_2)\times\right.
$$
$$
\left.\times
\int\limits_t^{t_2}\psi_1(t_1)\phi_{j_1}(t_1)dt_1
\int\limits_{t_3}^{T}\psi_4(t_4)\phi_{j_1}(t_4)dt_4 dt_2 dt_3\right)^2=
$$
$$
=
\sum_{j_2,j_3=0}^p\left(
\int\limits_t^T\psi_3(t_3)\phi_{j_3}(t_3)
\int\limits_t^{t_3}\psi_2(t_2)\phi_{j_2}(t_2)\times\right.
$$
$$
\left.\times
\sum_{j_1=p+1}^{\infty}\int\limits_t^{t_2}\psi_1(t_1)\phi_{j_1}(t_1)dt_1
\int\limits_{t_3}^{T}\psi_4(t_4)\phi_{j_1}(t_4)dt_4 dt_2 dt_3\right)^2\le
$$
$$
\le
\sum_{j_2,j_3=0}^{\infty}\left(
\int\limits_t^T\psi_3(t_3)\phi_{j_3}(t_3)
\int\limits_t^{t_3}\psi_2(t_2)\phi_{j_2}(t_2)\times\right.
$$
$$
\left.\times
\sum_{j_1=p+1}^{\infty}\int\limits_t^{t_2}\psi_1(t_1)\phi_{j_1}(t_1)dt_1
\int\limits_{t_3}^{T}\psi_4(t_4)\phi_{j_1}(t_4)dt_4 dt_2 dt_3\right)^2=
$$

\vspace{-1mm}
$$
=
\int\limits_{[t, T]^2}{\bf 1}_{\{t_2<t_3\}}
\psi_3^2(t_3)
\psi_2^2(t_2)\times
$$
$$
\times
\left(\sum_{j_1=p+1}^{\infty}\int\limits_t^{t_2}\psi_1(t_1)\phi_{j_1}(t_1)dt_1
\int\limits_{t_3}^{T}\psi_4(t_4)\phi_{j_1}(t_4)dt_4\right)^2dt_2dt_3\le 
$$
\begin{equation}
\label{after59}
\le \frac{K}{p^2}\ \to\  0
\end{equation}

\noindent
if $p\to\infty,$ where constant $K$ does not depend on $p;$
$$
\sum_{j_1,j_4=0}^p\left(\sum_{j_2=p+1}^{\infty}
C_{j_4 j_2 j_2 j_1}\right)^2=
\sum_{j_1,j_4=0}^p\left(\sum_{j_2=p+1}^{\infty}
\int\limits_t^T\psi_4(t_4)\phi_{j_4}(t_4)
\int\limits_t^{t_4}\psi_3(t_3)\phi_{j_2}(t_3)\times\right.
$$
$$
\left.\times
\int\limits_t^{t_3}\psi_2(t_2)\phi_{j_2}(t_2)
\int\limits_t^{t_2}\psi_1(t_1)\phi_{j_1}(t_1)dt_1dt_2dt_3dt_4\right)^2=
$$
$$
=
\sum_{j_1,j_4=0}^p\left(\sum_{j_2=p+1}^{\infty}
\int\limits_t^T\psi_4(t_4)\phi_{j_4}(t_4)
\int\limits_t^{t_4}\psi_1(t_1)\phi_{j_1}(t_1)\times\right.
$$
$$
\left.\times
\int\limits_{t_1}^{t_4}\psi_2(t_2)\phi_{j_2}(t_2)
\int\limits_{t_2}^{t_4}\psi_3(t_3)\phi_{j_2}(t_3)dt_3dt_2 dt_1 dt_4\right)^2=
$$
$$
=
\sum_{j_1,j_4=0}^p\left(
\int\limits_t^T\psi_4(t_4)\phi_{j_4}(t_4)
\int\limits_t^{t_4}\psi_1(t_1)\phi_{j_1}(t_1)\times\right.
$$
$$
\left.\times
\sum_{j_2=p+1}^{\infty}\int\limits_{t_1}^{t_4}\psi_2(t_2)\phi_{j_2}(t_2)
\int\limits_{t_2}^{t_4}\psi_3(t_3)\phi_{j_2}(t_3)dt_3dt_2 dt_1 dt_4\right)^2\le
$$
$$
\le
\sum_{j_1,j_4=0}^{\infty}\left(
\int\limits_t^T\psi_4(t_4)\phi_{j_4}(t_4)
\int\limits_t^{t_4}\psi_1(t_1)\phi_{j_1}(t_1)\times\right.
$$
$$
\left.\times
\sum_{j_2=p+1}^{\infty}\int\limits_{t_1}^{t_4}\psi_2(t_2)\phi_{j_2}(t_2)
\int\limits_{t_2}^{t_4}\psi_3(t_3)\phi_{j_2}(t_3)dt_3dt_2 dt_1 dt_4\right)^2=
$$

\vspace{-1mm}
$$
=
\int\limits_{[t, T]^2}{\bf 1}_{\{t_1<t_4\}}
\psi_4^2(t_4)
\psi_1^2(t_1)\times
$$
\begin{equation}
\label{after72}
~~~~~~~~~\times
\left(\sum_{j_2=p+1}^{\infty}\int\limits_{t_1}^{t_4}\psi_2(t_2)\phi_{j_2}(t_2)
\int\limits_{t_2}^{t_4}\psi_3(t_3)\phi_{j_2}(t_3)dt_3 dt_2\right)^2 dt_1dt_4.
\end{equation}

\vspace{2mm}

Note that, by virtue of the additivity property of the integral, we have
\begin{equation}
\label{after8000}
\sum_{j_2=p+1}^{\infty}\int\limits_{t_1}^{t_4}\psi_2(t_2)\phi_{j_2}(t_2)
\int\limits_{t_2}^{t_4}\psi_3(t_3)\phi_{j_2}(t_3)dt_3 dt_2=
\end{equation}
$$
=\sum_{j_2=p+1}^{\infty}
\int\limits_{t}^{t_4}\psi_3(t_3)\phi_{j_2}(t_3)
\int\limits_{t}^{t_3}\psi_2(t_2)\phi_{j_2}(t_2)dt_2 dt_3-
$$
$$
-\sum_{j_2=p+1}^{\infty}\int\limits_{t}^{t_1}\psi_3(t_3)\phi_{j_2}(t_3)
\int\limits_{t}^{t_3}\psi_2(t_2)\phi_{j_2}(t_2)dt_2 dt_3-
$$
\begin{equation}
\label{after71}
-\sum_{j_2=p+1}^{\infty}
\int\limits_{t_1}^{t_4}\psi_3(t_3)\phi_{j_2}(t_3)dt_3
\int\limits_{t}^{t_1}\psi_2(t_2)\phi_{j_2}(t_2)dt_2.
\end{equation}

\vspace{2mm}

However, all three series on the right-hand side of (\ref{after71})
have already been evaluated in (\ref{after5004}) and (\ref{after6009}).
From (\ref{after72}) and (\ref{after71}) we finally obtain
\begin{equation}
\label{after8001}
\sum_{j_1,j_4=0}^p\left(\sum_{j_2=p+1}^{\infty}
C_{j_4 j_2 j_2 j_1}\right)^2\le \frac{K}{p^{2-\varepsilon}}\ \to\  0
\end{equation}

\noindent
if $p\to\infty,$ where $\varepsilon$ is an arbitrary small positive real number
for the polynomial case and $\varepsilon=0$ for the 
trigonometric case, constant $K$ does not depend on $p.$

In complete analogy with (\ref{after6009}), we have for (\ref{after2505})
$$
\sum_{j_1,j_3=0}^p\left(\sum_{j_2=p+1}^{\infty}
C_{j_2 j_3 j_2 j_1}\right)^2=
\sum_{j_1,j_3=0}^p\left(\sum_{j_2=p+1}^{\infty}
\int\limits_t^T\psi_4(t_4)\phi_{j_2}(t_4)
\int\limits_t^{t_4}\psi_3(t_3)\phi_{j_3}(t_3)\times\right.
$$
$$
\left.\times
\int\limits_t^{t_3}\psi_2(t_2)\phi_{j_2}(t_2)
\int\limits_t^{t_2}\psi_1(t_1)\phi_{j_1}(t_1)dt_1dt_2dt_3dt_4\right)^2=
$$
$$
=
\sum_{j_1,j_3=0}^p\left(\sum_{j_2=p+1}^{\infty}
\int\limits_t^T\psi_3(t_3)\phi_{j_3}(t_3)
\int\limits_t^{t_3}\psi_2(t_2)\phi_{j_2}(t_2)\times\right.
$$
$$
\left.\times
\int\limits_t^{t_2}\psi_1(t_1)\phi_{j_1}(t_1)dt_1 dt_2
\int\limits_{t_3}^{T}\psi_4(t_4)\phi_{j_2}(t_4)dt_4dt_3\right)^2=
$$
$$
=
\sum_{j_1,j_3=0}^p\left(\sum_{j_2=p+1}^{\infty}
\int\limits_t^T\psi_3(t_3)\phi_{j_3}(t_3)
\int\limits_t^{t_3}\psi_1(t_1)\phi_{j_1}(t_1)\times\right.
$$
$$
\left.\times
\int\limits_{t_1}^{t_3}\psi_2(t_2)\phi_{j_2}(t_2)dt_2 dt_1
\int\limits_{t_3}^{T}\psi_4(t_4)\phi_{j_2}(t_4)dt_4dt_3\right)^2=
$$
$$
=
\sum_{j_1,j_3=0}^p\left(
\int\limits_t^T\psi_3(t_3)\phi_{j_3}(t_3)
\int\limits_t^{t_3}\psi_1(t_1)\phi_{j_1}(t_1)\times\right.
$$
$$
\left.\times
\sum_{j_2=p+1}^{\infty}\int\limits_{t_1}^{t_3}\psi_2(t_2)\phi_{j_2}(t_2)dt_2 
\int\limits_{t_3}^{T}\psi_4(t_4)\phi_{j_2}(t_4)dt_4 dt_1 dt_3\right)^2\le
$$
$$
\le
\sum_{j_1,j_3=0}^{\infty}\left(
\int\limits_t^T\psi_3(t_3)\phi_{j_3}(t_3)
\int\limits_t^{t_3}\psi_1(t_1)\phi_{j_1}(t_1)\times\right.
$$
$$
\left.\times
\sum_{j_2=p+1}^{\infty}\int\limits_{t_1}^{t_3}\psi_2(t_2)\phi_{j_2}(t_2)dt_2 
\int\limits_{t_3}^{T}\psi_4(t_4)\phi_{j_2}(t_4)dt_4 dt_1 dt_3\right)^2=
$$

\vspace{-1mm}
$$
=
\int\limits_{[t,T]^2}{\bf 1}_{\{t_1<t_3\}}\psi_3^2(t_3)
\psi_1^2(t_1)\times
$$
\begin{equation}
\label{after9030}
\times
\left(\sum_{j_2=p+1}^{\infty}\int\limits_{t_1}^{t_3}\psi_2(t_2)\phi_{j_2}(t_2)dt_2 
\int\limits_{t_3}^{T}\psi_4(t_4)\phi_{j_2}(t_4)dt_4\right)^2
dt_1 dt_3
\le \frac{K}{p^{2-\varepsilon}}\ \to\  0
\end{equation}

\noindent
if $p\to\infty,$ where $\varepsilon$ is an arbitrary small positive real number
for the polynomial case and $\varepsilon=0$ for the 
trigonometric case, constant $K$ does not depend on $p.$

We have for (\ref{after2506})
$$
\sum_{j_1,j_2=0}^p\left(\sum_{j_3=p+1}^{\infty}
C_{j_3 j_3 j_2 j_1}\right)^2=
\sum_{j_1,j_2=0}^p\left(\sum_{j_3=p+1}^{\infty}
\int\limits_t^T\psi_4(t_4)\phi_{j_3}(t_4)
\int\limits_t^{t_4}\psi_3(t_3)\phi_{j_3}(t_3)\times\right.
$$
$$
\left.\times
\int\limits_t^{t_3}\psi_2(t_2)\phi_{j_2}(t_2)
\int\limits_t^{t_2}\psi_1(t_1)\phi_{j_1}(t_1)dt_1dt_2dt_3dt_4\right)^2=
$$
$$
=
\sum_{j_1,j_2=0}^p\left(\sum_{j_3=p+1}^{\infty}
\int\limits_t^T\psi_1(t_1)\phi_{j_1}(t_1)
\int\limits_{t_1}^{T}\psi_2(t_2)\phi_{j_2}(t_2)\times\right.
$$
$$
\left.\times
\int\limits_{t_2}^{T}\psi_3(t_3)\phi_{j_3}(t_3)
\int\limits_{t_3}^{T}\psi_4(t_4)\phi_{j_3}(t_4)dt_4dt_3dt_2dt_1\right)^2=
$$
$$
=
\sum_{j_1,j_2=0}^p\left(
\int\limits_t^T\psi_1(t_1)\phi_{j_1}(t_1)
\int\limits_{t_1}^{T}\psi_2(t_2)\phi_{j_2}(t_2)\times\right.
$$
$$
\left.\times
\sum_{j_3=p+1}^{\infty}\int\limits_{t_2}^{T}\psi_3(t_3)\phi_{j_3}(t_3)
\int\limits_{t_3}^{T}\psi_4(t_4)\phi_{j_3}(t_4)dt_4dt_3dt_2dt_1\right)^2\le
$$
$$
\le
\sum_{j_1,j_2=0}^{\infty}\left(
\int\limits_t^T\psi_1(t_1)\phi_{j_1}(t_1)
\int\limits_{t_1}^{T}\psi_2(t_2)\phi_{j_2}(t_2)\times\right.
$$
$$
\left.\times
\sum_{j_3=p+1}^{\infty}\int\limits_{t_2}^{T}\psi_3(t_3)\phi_{j_3}(t_3)
\int\limits_{t_3}^{T}\psi_4(t_4)\phi_{j_3}(t_4)dt_4dt_3dt_2dt_1\right)^2=
$$

\vspace{-1mm}
$$
=
\int\limits_{[t,T]^2}{\bf 1}_{\{t_1<t_2\}}\psi_1^2(t_1)
\psi_2^2(t_2)\times
$$
\begin{equation}
\label{after8002}
~~~~~~~~~~\times
\left(\sum_{j_3=p+1}^{\infty}\int\limits_{t_2}^{T}\psi_3(t_3)\phi_{j_3}(t_3)
\int\limits_{t_3}^{T}\psi_4(t_4)\phi_{j_3}(t_4)dt_4dt_3\right)^2 dt_2dt_1.
\end{equation}

\vspace{2mm}

It is easy to see that the integral (see (\ref{after8002}))
$$
\int\limits_{t_2}^{T}\psi_3(t_3)\phi_{j_3}(t_3)
\int\limits_{t_3}^{T}\psi_4(t_4)\phi_{j_3}(t_4)dt_4dt_3
$$

\noindent
is similar to the integral from the formula
(\ref{after8000}) if in the last integral we substitute $t_4=T.$
Therefore, by analogy with (\ref{after8001}), we obtain
\begin{equation}
\label{after9040}
\sum_{j_1,j_2=0}^p\left(\sum_{j_3=p+1}^{\infty}
C_{j_3 j_3 j_2 j_1}\right)^2
\le \frac{K}{p^{2-\varepsilon}}\ \to\  0
\end{equation}

\noindent
if $p\to\infty,$ where $\varepsilon$ is an arbitrary small positive real number
for the polynomial case and $\varepsilon=0$ for the 
trigonometric case, constant $K$ does not depend on $p.$

Now consider (\ref{after2508})--(\ref{after2507}).
We have for (\ref{after2508}) (see {\bf Step~2} in the proof
of Theorem~2.30)
$$
\left(\sum_{j_2=p+1}^{\infty}\sum_{j_1=p+1}^{\infty}
C_{j_2 j_1 j_2 j_1}\right)^2=
\left(\sum_{j_1=0}^{p}\sum_{j_2=p+1}^{\infty}
C_{j_2 j_1 j_2 j_1}\right)^2\le 
$$
\begin{equation}
\label{after8003}
\le (p+1)
\sum_{j_1=0}^{p}\left(\sum_{j_2=p+1}^{\infty}
C_{j_2 j_1 j_2 j_1}\right)^2.
\end{equation}

\vspace{2mm}

Consider (\ref{after2505}) and (\ref{after9030}). We have
$$
\sum_{j_1=0}^p\left(\sum_{j_2=p+1}^{\infty}
C_{j_2 j_1 j_2 j_1}\right)^2=
\sum_{j_1,j_3=0}^p\left(\sum_{j_2=p+1}^{\infty}
C_{j_2 j_3 j_2 j_1}\right)^2\Biggl|_{j_1=j_3}\Biggr.\le
$$
\begin{equation}
\label{after8009}
\le\sum_{j_1,j_3=0}^p\left(\sum_{j_2=p+1}^{\infty}
C_{j_2 j_3 j_2 j_1}\right)^2
 \le \frac{K}{p^{2-\varepsilon}},
\end{equation}

\vspace{2mm}
\noindent
where $\varepsilon$ is an arbitrary small positive real number
for the polynomial case and $\varepsilon=0$ for the 
trigonometric case, constant $K$ does not depend on $p.$
Combining (\ref{after8003}) and (\ref{after8009}), we obtain
$$
\left(\sum_{j_2=p+1}^{\infty}\sum_{j_1=p+1}^{\infty}
C_{j_2 j_1 j_2 j_1}\right)^2\le \frac{(p+1)K}{p^{2-\varepsilon}}\le
\frac{K_1}{p^{1-\varepsilon}}\ \to\  0
$$

\noindent
if $p\to\infty,$ where constant $K_1$ does not depend on $p.$

Similarly for (\ref{after2509}) we have (see (\ref{after2504}), (\ref{after8001}))
$$
\left(\sum_{j_2=p+1}^{\infty}\sum_{j_1=p+1}^{\infty}
C_{j_1 j_2 j_2 j_1}\right)^2=
\left(\sum_{j_1=0}^{p}\sum_{j_2=p+1}^{\infty}
C_{j_1 j_2 j_2 j_1}\right)^2\le 
$$
\begin{equation}
\label{after9002}
\le
(p+1)\sum_{j_1=0}^{p}
\left(\sum_{j_2=p+1}^{\infty}
C_{j_1 j_2 j_2 j_1}\right)^2,
\end{equation}

\vspace{1mm}
$$
\sum_{j_1=0}^p\left(\sum_{j_2=p+1}^{\infty}
C_{j_1 j_2 j_2 j_1}\right)^2=
\sum_{j_1,j_4=0}^p\left(\sum_{j_2=p+1}^{\infty}
C_{j_4 j_2 j_2 j_1}\right)^2\Biggl|_{j_1=j_4}\Biggr.\le
$$
\begin{equation}
\label{after9011}
\le
\sum_{j_1,j_4=0}^p\left(\sum_{j_2=p+1}^{\infty}
C_{j_4 j_2 j_2 j_1}\right)^2\le \frac{K}{p^{2-\varepsilon}},
\end{equation}

\vspace{2mm}
\noindent
where $\varepsilon$ is an arbitrary small positive real number
for the polynomial case and $\varepsilon=0$ for the 
trigonometric case, constant $K$ does not depend on $p.$
Combining (\ref{after9002}) and (\ref{after9011}), we obtain
$$
\left(\sum_{j_2=p+1}^{\infty}\sum_{j_1=p+1}^{\infty}
C_{j_1 j_2 j_2 j_1}\right)^2\le \frac{(p+1)K}{p^{2-\varepsilon}}\le
\frac{K_1}{p^{1-\varepsilon}}\ \to\  0
$$

\noindent
if $p\to\infty,$ where constant $K_1$ does not depend on $p.$

Consider (\ref{after2507}). Using (\ref{after500}), we obtain
$$
\sum_{j_3=p+1}^{\infty}\sum_{j_1=p+1}^{\infty}
C_{j_3 j_3 j_1 j_1}
=\sum_{j_3=p+1}^{\infty}\sum_{j_1=0}^{\infty}
C_{j_3 j_3 j_1 j_1}-\sum_{j_3=p+1}^{\infty}\sum_{j_1=0}^{p}
C_{j_3 j_3 j_1 j_1}=
$$

\vspace{-1mm}
\begin{equation}
\label{after9041}
=\frac{1}{2}\sum_{j_3=p+1}^{\infty}
C_{j_3 j_3 j_1 j_1}\biggl|_{(j_1 j_1)\curvearrowright (\cdot) }\biggr.
-\sum_{j_3=p+1}^{\infty}\sum_{j_1=0}^{p}
C_{j_3 j_3 j_1 j_1},
\end{equation}

\vspace{2mm}
\noindent
where (see (\ref{after900}))
$$
C_{j_3 j_3 j_1 j_1}\biggl|_{(j_1 j_1)\curvearrowright (\cdot) }\biggr.=
\int\limits_t^T\psi_4(t_4)\phi_{j_3}(t_4)
\int\limits_t^{t_4}\psi_3(t_3)\phi_{j_3}(t_3)
\int\limits_t^{t_3}\psi_2(t_2)\psi_1(t_2)dt_2 dt_3 dt_4.
$$

\vspace{2mm}

From the estimate (\ref{tupo15}) (polynomial case) and its
analogue for the trigonometric case (see the proof of Lemma~2.2, Sect.~2.1.2) we get
\begin{equation}
\label{after9042}
\left|\sum_{j_3=p+1}^{\infty}
C_{j_3 j_3 j_1 j_1}\biggl|_{(j_1 j_1)\curvearrowright (\cdot) }\biggr.\right|\le
\frac{C}{p},
\end{equation}

\vspace{2mm}
\noindent
where constant $C$ is independent of $p.$

Further, we have (see (\ref{after9040}))
$$
\left(\sum_{j_1=0}^{p}\sum_{j_3=p+1}^{\infty}
C_{j_3 j_3 j_1 j_1}\right)^2\le (p+1)\sum_{j_1=0}^{p}
\left(\sum_{j_3=p+1}^{\infty}C_{j_3 j_3 j_1 j_1}\right)^2
=$$
$$
=(p+1)\sum_{j_1,j_2=0}^{p}
\left(\sum_{j_3=p+1}^{\infty}C_{j_3 j_3 j_2 j_1}\right)^2\Biggl|_{j_1=j_2}\Biggr.\le
$$
\begin{equation}
\label{after9043}
~~~~~~~~~~~\le
(p+1)\sum_{j_1,j_2=0}^{p}
\left(\sum_{j_3=p+1}^{\infty}C_{j_3 j_3 j_2 j_1}\right)^2\le
\frac{(p+1)K}{p^{2-\varepsilon}}\le
\frac{K_1}{p^{1-\varepsilon}},
\end{equation}

\vspace{3mm}
\noindent
where constant $K_1$ does not depend on $p.$

Combining (\ref{after9041})--(\ref{after9043}), we obtain
$$
\left(\sum_{j_3=p+1}^{\infty}\sum_{j_1=p+1}^{\infty}
C_{j_3 j_3 j_1 j_1}\right)^2
\le \frac{K_2}{p^{1-\varepsilon}}\ \to\  0
$$

\vspace{2mm}
\noindent
if $p\to\infty,$ where constant $K_2$ does not depend on $p.$

Let us prove (\ref{after1602xxxx})--(\ref{after1602xxxz}).
It is not difficult to see that the estimate (\ref{after9042})
proves (\ref{after1602xxxx}).

Using the integration order replacement, we obtain
$$
\sum_{j_1=p+1}^{\infty}
C_{j_3 j_3 j_1 j_1}\biggl|_{(j_{3} j_{3})\curvearrowright (\cdot)}=
$$
$$
=
\sum_{j_1=p+1}^{\infty}\int\limits_t^T\psi_4(t_4)\psi_3(t_4)
\int\limits_t^{t_4}\psi_2(t_2)\phi_{j_1}(t_2)
\int\limits_t^{t_2}\psi_1(t_1)\phi_{j_1}(t_1)dt_1dt_2dt_4=
$$
\begin{equation}
\label{afterafter1}
~~~=
\sum_{j_1=p+1}^{\infty}
\int\limits_t^{T}\left(\psi_2(t_2)
\int\limits_{t_2}^T\psi_4(t_4)\psi_3(t_4)dt_4\right)\phi_{j_1}(t_2)
\int\limits_t^{t_2}\psi_1(t_1)\phi_{j_1}(t_1)dt_1dt_2,
\end{equation}

\vspace{2mm}
$$
\sum_{j_1=p+1}^{\infty}
C_{j_1 j_2 j_2 j_1}\biggl|_{(j_{2} j_{2})\curvearrowright (\cdot)}\biggr.=
$$
$$
=
\sum_{j_1=p+1}^{\infty}\int\limits_t^T\psi_4(t_4)\phi_{j_1}(t_4)
\int\limits_t^{t_4}\psi_3(t_3)\psi_2(t_3)
\int\limits_t^{t_3}\psi_1(t_1)\phi_{j_1}(t_1)dt_1dt_3dt_4=
$$
$$
=
\sum_{j_1=p+1}^{\infty}\int\limits_t^T\psi_4(t_4)\phi_{j_1}(t_4)
\int\limits_t^{t_4}\psi_1(t_1)\phi_{j_1}(t_1)\int\limits_{t_1}^{t_4}         
\psi_3(t_3)\psi_2(t_3)dt_3dt_1dt_4=
$$
$$
=
\sum_{j_1=p+1}^{\infty}\int\limits_t^T\psi_4(t_4)\phi_{j_1}(t_4)
\int\limits_t^{t_4}\psi_1(t_1)\phi_{j_1}(t_1)
\left(\int\limits_{t}^{t_4}-\int\limits_{t}^{t_1}\right)
\psi_3(t_3)\psi_2(t_3)dt_3dt_1dt_4=
$$

\vspace{-3mm}
\begin{equation}
\label{afterafter198}
=
\sum_{j_1=p+1}^{\infty}\int\limits_t^T\left(\psi_4(t_4)
\int\limits_{t}^{t_4}
\psi_3(t_3)\psi_2(t_3)dt_3\right)
\phi_{j_1}(t_4)
\int\limits_t^{t_4}\psi_1(t_1)\phi_{j_1}(t_1)
dt_1dt_4-
\end{equation}

\vspace{-3mm}
\begin{equation}
\label{afterafter199}
-
\sum_{j_1=p+1}^{\infty}\int\limits_t^T\psi_4(t_4)\phi_{j_1}(t_4)
\int\limits_t^{t_4}\left(\psi_1(t_1)
\int\limits_{t}^{t_1}\
\psi_3(t_3)\psi_2(t_3)dt_3\right)\phi_{j_1}(t_1)dt_1dt_4.
\end{equation}

\vspace{2mm}

Applying the estimate (\ref{tupo15}) (polynomial case) and its
analogue for the trigonometric case (see the proof of Lemma~2.2, Sect.~2.1.2) 
to the right-hand sides of (\ref{afterafter1})--(\ref{afterafter199}), we get
\begin{equation}
\label{after9042xxx}
\left|\sum_{j_3=p+1}^{\infty}
C_{j_3 j_3 j_1 j_1}\biggl|_{(j_3 j_3)\curvearrowright (\cdot) }\biggr.\right|\le
\frac{C}{p},
\end{equation}

\vspace{-1mm}
\begin{equation}
\label{after9042xxxe}
\left|\sum_{j_1=p+1}^{\infty}
C_{j_1 j_2 j_2 j_1}\biggl|_{(j_{2} j_{2})\curvearrowright (\cdot)}\biggr.\right|\le
\frac{C}{p},
\end{equation}

\vspace{3mm}
\noindent
where constant $C$ is independent of $p.$
The estimates (\ref{after9042xxx}), (\ref{after9042xxxe})
prove (\ref{after1602xxxy}), (\ref{after1602xxxz}).

The relations (\ref{after2501})--(\ref{after1602xxxz}) are proved. Theorem~2.34 is proved.

\section{Expansion of Iterated Stratonovich Stochastic Integrals
of Multiplicity 5. The Case $p_1=\ldots =p_5\to \infty$ and 
Continuously Differentiable 
Weight Functions $\psi_1(\tau),$ $\ldots,$ $\psi_5(\tau)$ 
(The Cases of Legendre 
Polynomials and Trigonometric Functions)}

{\bf Theorem 2.35}\ \cite{arxiv-5}, 
\cite{arxiv-10}, \cite{arxiv-11}, \cite{new-art-1xxy}.\
{\it Suppose that 
$\{\phi_j(x)\}_{j=0}^{\infty}$ is a complete orthonormal system of 
Legendre polynomials or trigonometric functions in the space $L_2([t, T]).$
Furthermore, let $\psi_1(\tau), \ldots, \psi_5(\tau)$ are continuously dif\-ferentiable 
nonrandom functions on $[t, T].$ 
Then$,$ for the 
iterated Stratonovich stochastic integral of fifth multiplicity
\begin{equation}
\label{after10001}
J^{*}[\psi^{(5)}]_{T,t}={\int\limits_t^{*}}^T\psi_5(t_5)
\ldots
{\int\limits_t^{*}}^{t_2}\psi_1(t_1)
d{\bf w}_{t_1}^{(i_1)}
\ldots d{\bf w}_{t_5}^{(i_5)}
\end{equation}
the following 
expansion 
$$
J^{*}[\psi^{(5)}]_{T,t}
=\hbox{\vtop{\offinterlineskip\halign{
\hfil#\hfil\cr
{\rm l.i.m.}\cr
$\stackrel{}{{}_{p\to \infty}}$\cr
}} }
\sum\limits_{j_1, \ldots, j_5=0}^{p}
C_{j_5 \ldots j_1}\zeta_{j_1}^{(i_1)}\ldots
\zeta_{j_5}^{(i_5)}
$$
that converges in the mean-square sense is valid, where
$i_1, \ldots, i_5=0, 1,\ldots,m,$
$$
C_{j_5 \ldots j_1}=
\int\limits_t^T\psi_5(t_5)\phi_{j_5}(t_5)\ldots
\int\limits_t^{t_2}\psi_1(t_1)\phi_{j_1}(t_1)dt_1\ldots dt_5
$$
and
$$
\zeta_{j}^{(i)}=
\int\limits_t^T \phi_{j}(s) d{\bf w}_s^{(i)}
$$ 
are independent standard Gaussian random variables for various 
$i$ or $j$
{\rm (}in the case when $i\ne 0${\rm ),}
${\bf w}_{\tau}^{(i)}$ 
$(i=1,\ldots,m)$ are independent standard Wiener processes$,$
${\bf w}_{\tau}^{(0)}=\tau.$}

\vspace{2mm}

{\bf Proof.}\ Note that in this proof we write $k$ instead of 5 when this is true for 
an arbitrary $k$ $(k\in {\bf N}).$
As noted in Remark~2.4, Conditions 1 and 2
of Theorem~{\rm 2.30} are satisfied for complete
orthonormal systems of Legendre polynomials 
and trigonometric functions in the space
$L_2([t, T]).$ Let us verify Condition 3 of Theorem~2.30 for
the iterated Stratonovich stochastic integral (\ref{after10001}). 
Thus, we have to check the following conditions
\begin{equation}
\label{after14000}
\lim\limits_{p\to\infty}
\sum\limits_{j_{q_1},j_{q_2},j_{q_3}=0}^p
\left(\sum_{j_{g_1}=p+1}^{\infty}
C_{j_5\ldots j_1}\biggl|_{j_{g_1}=j_{g_2}}\biggr.\right)^2=0,
\end{equation}
\begin{equation}
\label{after14001}
~~~~~~~~~\lim\limits_{p\to\infty}
\sum\limits_{j_{q_1}=0}^p
\left(\sum_{j_{g_1}=p+1}^{\infty}\sum_{j_{g_3}=p+1}^{\infty}
C_{j_5\ldots j_1}\biggl|_{j_{g_1}=j_{g_2},j_{g_3}=j_{g_4}}\biggr.\right)^2=0,
\end{equation}
\begin{equation}
\label{afterafter001}
~~~~~~~~\lim\limits_{p\to\infty}
\sum\limits_{j_{q_1}=0}^p
\left(\sum_{j_{g_3}=p+1}^{\infty}
C_{j_5\ldots j_1}\biggl|_{(j_{g_2}j_{g_1})\curvearrowright (\cdot),
j_{g_1}=j_{g_2},
j_{g_3}=j_{g_4}, g_2=g_1+1}\biggr.\right)^2=0,
\end{equation}

\newpage
\noindent
where 
$\left(\{g_1,g_2\},\{g_3,g_4\}, \{q_1\}\right)$ and 
$\left(\{g_1,g_2\}, \{q_1, q_2,q_3\}\right)$ 
are partitions of the set $\{1,2,\ldots,5\}$ that is
$\{g_1,g_2,g_3,g_4,q_1\}=\{g_1,g_2,q_1,q_2,q_3\}=\{1,2,\ldots,5\};$
braces mean an unordered 
set, and pa\-ren\-the\-ses mean an ordered set.

Let us find a representation for
$C_{j_k\ldots j_1}\bigl|_{j_{g_1}=j_{g_2},\ g_2>g_1+1}\biggr.$ 
that will be convenient for further consideration.

Using the integration order replacement in Riemann integrals, we obtain
$$
\int\limits_t^T h_{k}(t_k)\ldots \int\limits_t^{t_{l+2}} h_{l+1}(t_{l+1})
\int\limits_t^{t_{l+1}} h_{l}(t_{l})
\int\limits_t^{t_{l}} h_{l-1}(t_{l-1})\ldots
\int\limits_t^{t_2} h_{1}(t_1)
dt_1\ldots 
$$
$$
\ldots
dt_{l-1}dt_{l}dt_{l+1}\ldots dt_k=
$$

\vspace{-4mm}
$$
=\int\limits_t^T h_{k}(t_k)\ldots \int\limits_t^{t_{l+2}} h_{l+1}(t_{l+1})
\int\limits_t^{t_{l+1}} h_{1}(t_{1})
\int\limits_{t_1}^{t_{l+1}} h_{2}(t_{2})\ldots
\int\limits_{t_{l-2}}^{t_{l+1}} h_{l-1}(t_{l-1})
\int\limits_{t_{l-1}}^{t_{l+1}} h_{l}(t_{l})dt_l\times
$$
$$
\times dt_{l-1}\ldots dt_2dt_{1}dt_{l+1}\ldots dt_k=
$$

\vspace{-4mm}
$$
=\int\limits_t^T h_{k}(t_k)\ldots \int\limits_t^{t_{l+2}} h_{l+1}(t_{l+1})
\left(\int\limits_{t}^{t_{l+1}} h_{l}(t_{l})dt_l\right)\int\limits_t^{t_{l+1}} h_{1}(t_{1})
\int\limits_{t_1}^{t_{l+1}} h_{2}(t_{2})\ldots
\int\limits_{t_{l-2}}^{t_{l+1}} h_{l-1}(t_{l-1})
\times
$$
$$
\times dt_{l-1}\ldots dt_2dt_{1}dt_{l+1}\ldots dt_k-
$$

\vspace{-4mm}
$$
-\int\limits_t^T h_{k}(t_k)\ldots \int\limits_t^{t_{l+2}} h_{l+1}(t_{l+1})
\int\limits_t^{t_{l+1}} h_{1}(t_{1})
\int\limits_{t_1}^{t_{l+1}} h_{2}(t_{2})\ldots
\int\limits_{t_{l-2}}^{t_{l+1}} h_{l-1}(t_{l-1})
\left(\int\limits_{t}^{t_{l-1}} h_{l}(t_{l})dt_l\right)\times
$$
$$
\times dt_{l-1}\ldots dt_2dt_{1}dt_{l+1}\ldots dt_k=
$$

\vspace{-1mm}
$$
=\int\limits_t^T h_{k}(t_k)\ldots \int\limits_t^{t_{l+2}} h_{l+1}(t_{l+1})
\left(\int\limits_t^{t_{l+1}} h_{l}(t_{l})dt_l\right)
\int\limits_t^{t_{l+1}} h_{l-1}(t_{l-1})\ldots
$$
$$
\ldots 
\int\limits_t^{t_2} h_{1}(t_1)
dt_1\ldots dt_{l-1}dt_{l+1}\ldots dt_k-
$$
$$
-\int\limits_t^T h_{k}(t_k)\ldots \int\limits_t^{t_{l+2}} h_{l+1}(t_{l+1})
\int\limits_t^{t_{l+1}} h_{l-1}(t_{l-1})\left(\int\limits_t^{t_{l-1}} h_{l}(t_{l})dt_l\right)
\int\limits_t^{t_{l-1}} h_{l-2}(t_{l-2})
\ldots
$$
\begin{equation}
\label{after81}
\ldots
\int\limits_t^{t_2} h_{1}(t_1)
dt_1\ldots dt_{l-2}dt_{l-1}dt_{l+1}\ldots dt_k,
\end{equation}

\vspace{1mm}
\noindent
where $2<l<k-1$ and $h_1(\tau),\ldots,h_k(\tau)$ are continuous functions on the interval
$[t, T].$ The case $l=1$ is obvious. By analogy with (\ref{after81}) we have for $l=k$
$$
\int\limits_t^{T} h_{l}(t_{l})
\int\limits_t^{t_{l}} h_{l-1}(t_{l-1})\ldots
\int\limits_t^{t_2} h_{1}(t_1)
dt_1\ldots 
dt_{l-1}dt_{l}=
$$
$$
=\int\limits_{t}^{T} h_{1}(t_{1})
\int\limits_{t_1}^{T} h_{2}(t_{2})\ldots
\int\limits_{t_{l-2}}^{T} h_{l-1}(t_{l-1})\int\limits_{t_{l-1}}^{T} h_{l}(t_{l})
dt_ldt_{l-1}\ldots dt_2dt_{1}=
$$
$$
=\left(\int\limits_{t}^{T} h_{l}(t_{l})
dt_l\right)\int\limits_{t}^{T} h_{1}(t_{1})
\int\limits_{t_1}^{T} h_{2}(t_{2})\ldots
\int\limits_{t_{l-2}}^{T} h_{l-1}(t_{l-1})
dt_{l-1}\ldots dt_2dt_{1}-
$$
$$
-\int\limits_{t}^{T}h_{1}(t_{1})
\int\limits_{t_1}^{T} h_{2}(t_{2})\ldots
\int\limits_{t_{l-2}}^{T} h_{l-1}(t_{l-1})\left(
\int\limits_{t}^{t_{l-1}} h_{l}(t_{l})dt_l\right)
dt_{l-1}\ldots dt_2dt_{1}=
$$
$$
=\left(\int\limits_{t}^{T} h_{l}(t_{l})
dt_l\right)\int\limits_{t}^{T} h_{l-1}(t_{l-1})
\ldots
\int\limits_{t}^{t_2} h_{1}(t_{1})
dt_{1}\ldots dt_{l-1}-
$$
\begin{equation}
\label{after82}
-\int\limits_{t}^{T}
h_{l-1}(t_{l-1})\left(
\int\limits_{t}^{t_{l-1}} h_{l}(t_{l})dt_l\right)
\int\limits_{t}^{t_{l-1}} h_{l-2}(t_{l-2})\ldots
\int\limits_{t}^{t_2} 
h_{1}(t_{1})
dt_{1}\ldots dt_{l-1}.
\end{equation}

\vspace{2mm}

The formulas (\ref{after81}), (\ref{after82})
will be used further.

Our further proof will not fundamentally depend on the weight
functions $\psi_1(\tau),\ldots,\psi_k(\tau).$
Therefore, sometimes in subsequent consideration we assume for simplicity
that $\psi_1(\tau),\ldots,\psi_k(\tau)\equiv 1.$

Let us continue the proof. Applying (\ref{after81}) 
to $C_{j_k \ldots j_{l+1} j_l j_{l-1} \ldots j_{s+1} j_l j_{s-1} \ldots j_1}$ 
(more precisely to $h_s(t_s)=\psi_s(t_s)\phi_{j_{l}}(t_{s})$), we obtain
for $l+1\le k,$ $s-1\ge 1,$ $l-1\ge s+1$
\begin{equation}
\label{after90a}
\sum_{j_l=p+1}^{\infty}
C_{j_k \ldots j_{l+1} j_l j_{l-1} \ldots j_{s+1} j_l j_{s-1} \ldots j_1}=
\end{equation}
$$
=
\sum_{j_l=p+1}^{\infty}
\int\limits_t^T \phi_{j_k}(t_k)\ldots 
\int\limits_t^{t_{l+2}} \phi_{j_{l+1}}(t_{l+1})
\int\limits_t^{t_{l+1}} \phi_{j_{l}}(t_{l})
\int\limits_t^{t_{l}}\phi_{j_{l-1}}(t_{l-1})\ldots
$$
$$
\ldots \int\limits_t^{t_{s+2}} \phi_{j_{s+1}}(t_{s+1})
\int\limits_t^{t_{s+1}} \phi_{j_{l}}(t_{s})
\int\limits_t^{t_{s}} \phi_{j_{s-1}}(t_{s-1})\ldots 
$$
$$
\ldots\int\limits_t^{t_{2}} \phi_{j_{1}}(t_{1})dt_1
\ldots dt_{s-1}dt_{s}dt_{s+1}\ldots
dt_{l-1}dt_{l}dt_{l+1}\ldots dt_k=
$$
$$
=
\sum_{j_l=p+1}^{\infty}
\int\limits_t^T \phi_{j_k}(t_k)\ldots \int\limits_t^{t_{l+2}} \phi_{j_{l+1}}(t_{l+1})
\int\limits_t^{t_{l+1}} \phi_{j_{l}}(t_{l})
\int\limits_t^{t_{l}} \phi_{j_{l-1}}(t_{l-1})\ldots
$$
$$
\ldots
\int\limits_t^{t_{s+2}} \phi_{j_{s+1}}(t_{s+1})
\left(\int\limits_t^{t_{s+1}} \phi_{j_{l}}(t_{s})dt_s\right)
\int\limits_t^{t_{s+1}} \phi_{j_{s-1}}(t_{s-1})\ldots
$$

\vspace{-2mm}
$$
\ldots \int\limits_t^{t_{2}} \phi_{j_{1}}(t_{1})dt_1\ldots dt_{s-1}dt_{s+1}\ldots
dt_{l-1}dt_{l}dt_{l+1}\ldots dt_k-
$$
$$
-\sum_{j_l=p+1}^{\infty}
\int\limits_t^T \phi_{j_k}(t_k)\ldots \int\limits_t^{t_{l+2}} \phi_{j_{l+1}}(t_{l+1})
\int\limits_t^{t_{l+1}} \phi_{j_{l}}(t_{l})
\int\limits_t^{t_{l}} \phi_{j_{l-1}}(t_{l-1})\ldots
$$
$$
\ldots
\int\limits_t^{t_{s+2}} \phi_{j_{s+1}}(t_{s+1})
\int\limits_t^{t_{s+1}} \phi_{j_{s-1}}(t_{s-1})
\left(\int\limits_t^{t_{s-1}} \phi_{j_{l}}(t_{s})dt_s\right)
\int\limits_t^{t_{s-1}}\phi_{j_{s-2}}(t_{s-2})
\ldots 
$$
$$
\ldots \int\limits_t^{t_{2}} \phi_{j_{1}}(t_{1})dt_1\ldots dt_{s-2}dt_{s-1}dt_{s+1}\ldots
dt_{l-1}dt_{l}dt_{l+1}\ldots dt_k =
$$
$$
=\sum_{j_l=p+1}^{\infty} A_{j_k \ldots j_{l+1} j_l j_{l-1} \ldots j_{s+1} j_l j_{s-1} \ldots j_1}-
\sum_{j_l=p+1}^{\infty} B_{j_k \ldots j_{l+1} j_l j_{l-1} \ldots j_{s+1} j_l j_{s-1} \ldots j_1}.
$$

\vspace{2mm}

Now we apply the formula (\ref{after81}) 
to the quantities $A_{j_k \ldots j_{l+1} j_l j_{l-1} \ldots j_{s+1} j_l j_{s-1} \ldots j_1}$
and $B_{j_k \ldots j_{l+1} j_l j_{l-1} \ldots j_{s+1} j_l j_{s-1} \ldots j_1}$
(more precisely to $h_l(t_l)=\psi_l(t_l)\phi_{j_{l}}(t_{l})$). Then we have
for $l+1\le k,$ $s-1\ge 1,$ $l-1\ge s+1$
$$
\sum_{j_l=p+1}^{\infty}
C_{j_k \ldots j_{l+1} j_l j_{l-1} \ldots j_{s+1} j_l j_{s-1} \ldots j_1}=
$$
$$
=\int\limits_{[t, T]^{k-2}} \sum\limits_{d=1}^4
F_p^{(d)}(t_1,\ldots,t_{s-1},t_{s+1},\ldots ,t_{l-1},t_{l+1},\ldots,t_k)\times
$$
$$
\times 
\prod_{\stackrel{g=1}{{}_{g\ne l, s}}}^{k}
\psi_{g}(t_g) \phi_{j_g}(t_g)
dt_1\ldots dt_{s-1}dt_{s+1}\ldots dt_{l-1}dt_{l+1}\ldots dt_k=
$$
\begin{equation}
\label{after85a}
=\sum\limits_{d=1}^4
C_{j_k \ldots j_{l+1} j_{l-1} \ldots j_{s+1} j_{s-1} \ldots j_1}^{*(d)}
=\sum\limits_{d=1}^4
C_{j_k \ldots j_q \ldots j_1}^{*(d)}\biggl|_{q\ne l, s}\biggr.,
\end{equation}

\vspace{4mm}
\noindent
where

\vspace{-5mm}
$$
F_p^{(1)}(t_1,\ldots,t_{s-1},t_{s+1},\ldots ,t_{l-1},t_{l+1},\ldots,t_k)=
$$
\begin{equation}
\label{after90}
=
{\bf 1}_{\{t_1< \ldots <t_{s-1}<t_{s+1}< \ldots <t_{l-1}<t_{l+1}< \ldots <t_k\}}
\sum_{j_l=p+1}^{\infty}~
\int\limits_{t}^{t_{s+1}} \psi_s(\tau) \phi_{j_{l}}(\tau)d\tau
\int\limits_{t}^{t_{l+1}} \psi_l(\tau) \phi_{j_{l}}(\tau)d\tau,
\end{equation}

\vspace{-3mm}
$$
F_p^{(2)}(t_1,\ldots,t_{s-1},t_{s+1},\ldots ,t_{l-1},t_{l+1},\ldots,t_k)=
$$
\begin{equation}
\label{after91}
=
{\bf 1}_{\{t_1< \ldots <t_{s-1}<t_{s+1}< \ldots <t_{l-1}<t_{l+1}< \ldots <t_k\}}
\sum_{j_l=p+1}^{\infty}~
\int\limits_{t}^{t_{s-1}} \psi_s(\tau) \phi_{j_{l}}(\tau)d\tau
\int\limits_{t}^{t_{l-1}} \psi_l(\tau) \phi_{j_{l}}(\tau)d\tau,
\end{equation}

\vspace{-3mm}
$$
F_p^{(3)}(t_1,\ldots,t_{s-1},t_{s+1},\ldots ,t_{l-1},t_{l+1},\ldots,t_k)=
$$
\begin{equation}
\label{after92}
=
-{\bf 1}_{\{t_1< \ldots <t_{s-1}<t_{s+1}< \ldots <t_{l-1}<t_{l+1}< \ldots <t_k\}}
\sum_{j_l=p+1}^{\infty}~
\int\limits_{t}^{t_{s-1}} \psi_s(\tau)\phi_{j_{l}}(\tau)d\tau
\int\limits_{t}^{t_{l+1}} \psi_l(\tau)\phi_{j_{l}}(\tau)d\tau,
\end{equation}

\vspace{-3mm}
$$
F_p^{(4)}(t_1,\ldots,t_{s-1},t_{s+1},\ldots ,t_{l-1},t_{l+1},\ldots,t_k)=
$$
\begin{equation}
\label{after93}
=
-{\bf 1}_{\{t_1< \ldots <t_{s-1}<t_{s+1}< \ldots <t_{l-1}<t_{l+1}< \ldots <t_k\}}
\sum_{j_l=p+1}^{\infty}~
\int\limits_{t}^{t_{s+1}} \psi_s(\tau)\phi_{j_{l}}(\tau)d\tau
\int\limits_{t}^{t_{l-1}} \psi_l(\tau)\phi_{j_{l}}(\tau)d\tau.
\end{equation}

\vspace{1mm}

By analogy with (\ref{after85a}) we can consider the expressions
\begin{equation}
\label{afterx200}
\sum_{j_l=p+1}^{\infty}
C_{j_l j_{k-1}\ldots j_2 j_l},
\end{equation}
\begin{equation}
\label{afterx201}
\sum_{j_l=p+1}^{\infty}
C_{j_k \ldots j_{l+1} j_l j_{l-1} \ldots j_2 j_l}\ \ \ (l+1\le k),
\end{equation}
\begin{equation}
\label{afterx202}
\sum_{j_l=p+1}^{\infty}
C_{j_l j_{k-1} \ldots j_{s+1} j_l j_{s-1} \ldots j_1}\ \ \ (s-1\ge 1).
\end{equation}

\vspace{2mm}

Then we have for (\ref{afterx200})--(\ref{afterx202}) (see (\ref{after81}), (\ref{after82}))
\begin{equation}
\label{afterx204}
\sum_{j_l=p+1}^{\infty}
C_{j_l j_{k-1} \ldots j_2 j_l}=
\int\limits_{[t, T]^{k-2}} \sum\limits_{d=1}^2
G_p^{(d)}(t_2,\ldots,t_{k-1})
\prod_{g=2}^{k-1}
\psi_{g}(t_g) \phi_{j_g}(t_g)
dt_2\ldots dt_{k-1},
\end{equation}
$$
\sum_{j_l=p+1}^{\infty}
C_{j_k \ldots j_{l+1} j_l j_{l-1} \ldots j_2 j_l}=
\int\limits_{[t, T]^{k-2}} \sum\limits_{d=1}^2
E_p^{(d)}(t_2,\ldots,t_{l-1},t_{l+1}, \ldots, t_{k})\times
$$
\begin{equation}
\label{afterx205}
\times \prod_{\stackrel{g=2}{{}_{g\ne l}}}^{k}
\psi_{g}(t_g)\phi_{j_g}(t_g)
dt_2\ldots dt_{l-1}dt_{l+1} \ldots dt_{k},
\end{equation}
$$
\sum_{j_l=p+1}^{\infty}
C_{j_l j_{k-1} \ldots j_{s+1} j_l j_{s-1} \ldots j_1}=
\int\limits_{[t, T]^{k-2}} \sum\limits_{d=1}^4
D_p^{(d)}(t_1,\ldots,t_{s-1},t_{s+1}, \ldots, t_{k-1})\times
$$
\begin{equation}
\label{afterx206}
\times \prod_{\stackrel{g=1}{{}_{g\ne s}}}^{k-1}
\psi_{g}(t_g) \phi_{j_g}(t_g)
dt_1\ldots dt_{s-1}dt_{s+1} \ldots dt_{k-1},
\end{equation}

\noindent
where
$$
G_p^{(1)}(t_2,\ldots,t_{k-1})=
{\bf 1}_{\{t_2< \ldots <t_{k-1}\}}
\sum_{j_l=p+1}^{\infty}~
\int\limits_{t}^{T} \psi_k(\tau)\phi_{j_{l}}(\tau)d\tau
\int\limits_{t}^{t_{2}} \psi_1(\tau)\phi_{j_{l}}(\tau)d\tau,
$$
$$
G_p^{(2)}(t_2,\ldots,t_{k-1})=
-{\bf 1}_{\{t_2< \ldots <t_{k-1}\}}
\sum_{j_l=p+1}^{\infty}~
\int\limits_{t}^{t_{k-1}} \psi_k(\tau)\phi_{j_{l}}(\tau)d\tau
\int\limits_{t}^{t_{2}} \psi_1(\tau)\phi_{j_{l}}(\tau)d\tau,
$$

$$
E_p^{(1)}(t_2,\ldots,t_{l-1}, t_{l+1}, \ldots, t_k)=
$$
$$
=
{\bf 1}_{\{t_2< \ldots <t_{l-1}<t_{l+1}<\ldots < t_k\}}
\sum_{j_l=p+1}^{\infty}~
\int\limits_{t}^{t_{l+1}} \psi_l(\tau)\phi_{j_{l}}(\tau)d\tau
\int\limits_{t}^{t_{2}} \psi_1(\tau)\phi_{j_{l}}(\tau)d\tau,
$$

$$
E_p^{(2)}(t_2,\ldots,t_{l-1}, t_{l+1}, \ldots, t_k)=
$$
$$
=
-{\bf 1}_{\{t_2< \ldots <t_{l-1}<t_{l+1}<\ldots < t_k\}}
\sum_{j_l=p+1}^{\infty}~
\int\limits_{t}^{t_{l-1}} \psi_l(\tau)\phi_{j_{l}}(\tau)d\tau
\int\limits_{t}^{t_{2}} \psi_1(\tau)\phi_{j_{l}}(\tau)d\tau,
$$

$$
D_p^{(1)}(t_1,\ldots,t_{s-1}, t_{s+1}, \ldots, t_{k-1})=
$$
$$
=
{\bf 1}_{\{t_1< \ldots <t_{s-1}<t_{s+1}<\ldots < t_{k-1}\}}
\sum_{j_l=p+1}^{\infty}~
\int\limits_{t}^{T} \psi_k(\tau)\phi_{j_{l}}(\tau)d\tau
\int\limits_{t}^{t_{s+1}} \psi_s(\tau)\phi_{j_{l}}(\tau)d\tau,
$$

$$
D_p^{(2)}(t_1,\ldots,t_{s-1}, t_{s+1}, \ldots, t_{k-1})=
$$
$$
=
-{\bf 1}_{\{t_1< \ldots <t_{s-1}<t_{s+1}<\ldots < t_{k-1}\}}
\sum_{j_l=p+1}^{\infty}~
\int\limits_{t}^{T} \psi_k(\tau)\phi_{j_{l}}(\tau)d\tau
\int\limits_{t}^{t_{s-1}} \psi_s(\tau)\phi_{j_{l}}(\tau)d\tau,
$$

$$
D_p^{(3)}(t_1,\ldots,t_{s-1}, t_{s+1}, \ldots, t_{k-1})=
$$
$$
=
-{\bf 1}_{\{t_1< \ldots <t_{s-1}<t_{s+1}<\ldots < t_{k-1}\}}
\sum_{j_l=p+1}^{\infty}~
\int\limits_{t}^{t_{k-1}} \psi_k(\tau)\phi_{j_{l}}(\tau)d\tau
\int\limits_{t}^{t_{s+1}} \psi_s(\tau)\phi_{j_{l}}(\tau)d\tau,
$$

$$
D_p^{(4)}(t_1,\ldots,t_{s-1}, t_{s+1}, \ldots, t_{k-1})=
$$
$$
=
{\bf 1}_{\{t_1< \ldots <t_{s-1}<t_{s+1}<\ldots < t_{k-1}\}}
\sum_{j_l=p+1}^{\infty}~
\int\limits_{t}^{t_{k-1}} \psi_k(\tau)\phi_{j_{l}}(\tau)d\tau
\int\limits_{t}^{t_{s-1}} \psi_s(\tau)\phi_{j_{l}}(\tau)d\tau.
$$

\vspace{-6mm}

Now let us consider the value
$C_{j_k\ldots j_1}\bigl|_{j_{g_1}=j_{g_2},\ g_2=g_1+1}\biggr.$. 
To do this, we will make the following transformations
$$
\int\limits_t^T h_{k}(t_k)\ldots \int\limits_t^{t_{l+2}} h_{l+1}(t_{l+1})
\int\limits_t^{t_{l+1}} h_{l}(t_{l})
\int\limits_t^{t_{l}} h_{l}(t_{l-1})
\int\limits_t^{t_{l-1}} h_{l-2}(t_{l-2})\ldots
\int\limits_t^{t_2} h_{1}(t_1)
dt_1\ldots 
$$
$$
\ldots
dt_{l-2}dt_{l-1}dt_{l}dt_{l+1}\ldots dt_k=
$$

\vspace{-5mm}
$$
=\int\limits_t^T h_{k}(t_k)\ldots \int\limits_t^{t_{l+2}} h_{l+1}(t_{l+1})
\int\limits_t^{t_{l+1}} h_{1}(t_{1})
\int\limits_{t_1}^{t_{l+1}} h_{2}(t_{2})\ldots
\int\limits_{t_{l-3}}^{t_{l+1}} h_{l-2}(t_{l-2})\times
$$
$$
\times
\left(\int\limits_{t}^{t_{l+1}}-
\int\limits_{t}^{t_{l-2}}~\right)h_{l}(t_{l-1})
\left(\int\limits_{t}^{t_{l+1}}-
\int\limits_{t}^{t_{l-1}}~\right)
h_{l}(t_{l})dt_l
dt_{l-1}dt_{l-2}\ldots dt_2dt_{1}dt_{l+1}\ldots dt_k=
$$
$$
=\int\limits_t^T h_{k}(t_k)\ldots \int\limits_t^{t_{l+2}} h_{l+1}(t_{l+1})
\left(\int\limits_t^{t_{l+1}} h_{l}(t_{l})dt_l
\int\limits_t^{t_{l+1}} h_{l}(t_{l-1})dt_{l-1}\right)
\int\limits_t^{t_{l+1}} h_{1}(t_{1})
\times
$$
$$
\times
\int\limits_{t_1}^{t_{l+1}} h_{2}(t_{2})\ldots \int\limits_{t_{l-3}}^{t_{l+1}} h_{l-2}(t_{l-2})
dt_{l-2}\ldots dt_2dt_{1}dt_{l+1}\ldots dt_k-
$$
$$
-\int\limits_t^T h_{k}(t_k)\ldots \int\limits_t^{t_{l+2}} h_{l+1}(t_{l+1})
\left(\int\limits_t^{t_{l+1}} h_{l}(t_{l})dt_l\right)
\int\limits_t^{t_{l+1}} h_{1}(t_{1})
\int\limits_{t_1}^{t_{l+1}} h_{2}(t_{2})\ldots 
$$
$$
\ldots \int\limits_{t_{l-3}}^{t_{l+1}} h_{l-2}(t_{l-2})
\left(\int\limits_t^{t_{l-2}} h_{l}(t_{l-1})dt_{l-1}\right)
dt_{l-2}\ldots dt_2dt_{1}dt_{l+1}\ldots dt_k-
$$
$$
-\int\limits_t^T h_{k}(t_k)\ldots \int\limits_t^{t_{l+2}} h_{l+1}(t_{l+1})
\left(\int\limits_t^{t_{l+1}} h_{l}(t_{l-1})
\int\limits_t^{t_{l-1}} h_{l}(t_{l})dt_l dt_{l-1}\right)
\int\limits_t^{t_{l+1}} h_{1}(t_{1})\times
$$
$$
\times
\int\limits_{t_1}^{t_{l+1}} h_{2}(t_{2})\ldots 
\int\limits_{t_{l-3}}^{t_{l+1}} h_{l-2}(t_{l-2})
dt_{l-2}\ldots dt_2dt_{1}dt_{l+1}\ldots dt_k+
$$
$$
+\int\limits_t^T h_{k}(t_k)\ldots \int\limits_t^{t_{l+2}} h_{l+1}(t_{l+1})
\int\limits_t^{t_{l+1}} h_{1}(t_{1})
\int\limits_{t_1}^{t_{l+1}} h_{2}(t_{2})\ldots \int\limits_{t_{l-3}}^{t_{l+1}} h_{l-2}(t_{l-2})\times
$$
$$
\times
\left(\int\limits_t^{t_{l-2}} h_{l}(t_{l-1})
\int\limits_t^{t_{l-1}} h_{l}(t_{l})dt_ldt_{l-1}\right)
dt_{l-2}\ldots dt_2dt_{1}dt_{l+1}\ldots dt_k=
$$
$$
=\int\limits_t^T h_{k}(t_k)\ldots \int\limits_t^{t_{l+2}} h_{l+1}(t_{l+1})
\left(\int\limits_t^{t_{l+1}} h_{l}(t_{l})dt_l
\int\limits_t^{t_{l+1}} h_{l}(t_{l-1})dt_{l-1}\right)
\int\limits_t^{t_{l+1}} h_{l-2}(t_{l-2})
\times
$$
$$
\times
\int\limits_{t}^{t_{l-2}} h_{l-3}(t_{l-3})\ldots \int\limits_{t}^{t_2} h_{1}(t_{1})dt_1
\ldots dt_{l-3}dt_{l-2}dt_{l+1}\ldots dt_k-
$$
$$
-\int\limits_t^T h_{k}(t_k)\ldots \int\limits_t^{t_{l+2}} h_{l+1}(t_{l+1})
\left(\int\limits_t^{t_{l+1}} h_{l}(t_{l})dt_l\right)
\int\limits_t^{t_{l+1}} h_{l-2}(t_{l-2})
\times
$$
$$
\times\left(\int\limits_t^{t_{l-2}} h_{l}(t_{l-1})dt_{l-1}\right)
\int\limits_{t}^{t_{l-2}} h_{l-3}(t_{l-3})\ldots \int\limits_{t}^{t_2} h_{1}(t_{1})dt_1
\ldots dt_{l-3}dt_{l-2}dt_{l+1}\ldots dt_k-
$$
$$
-\int\limits_t^T h_{k}(t_k)\ldots \int\limits_t^{t_{l+2}} h_{l+1}(t_{l+1})
\left(\int\limits_t^{t_{l+1}} h_{l}(t_{l-1})
\int\limits_t^{t_{l-1}} h_{l}(t_{l})dt_l dt_{l-1}\right)
\times
$$
$$
\times
\int\limits_t^{t_{l+1}} h_{l-2}(t_{l-2})
\int\limits_{t}^{t_{l-2}} h_{l-3}(t_{l-3})\ldots \int\limits_{t}^{t_2} h_{1}(t_{1})dt_1
\ldots dt_{l-3}dt_{l-2}dt_{l+1}\ldots dt_k+
$$
$$
+\int\limits_t^T h_{k}(t_k)\ldots \int\limits_t^{t_{l+2}} h_{l+1}(t_{l+1})
\int\limits_t^{t_{l+1}}
h_{l-2}(t_{l-2})
\left(\int\limits_t^{t_{l-2}} h_{l}(t_{l-1})
\int\limits_t^{t_{l-1}} h_{l}(t_{l})dt_ldt_{l-1}\right)
\times
$$
\begin{equation}
\label{after9031}
~~~~~~\times
\int\limits_{t}^{t_{l-2}} h_{l-3}(t_{l-3})\ldots \int\limits_{t}^{t_2} h_{1}(t_{1})dt_1
\ldots dt_{l-3}dt_{l-2}dt_{l+1}\ldots dt_k,
\end{equation}

\vspace{1mm}
\noindent
where $l+1\le k,$ $l-2\ge 1,$ and
$h_1(\tau),\ldots,h_k(\tau)$ are continuous functions on the interval
$[t, T].$  The case $l=k$ follows from (\ref{after9031}) with $t_{l+1}=T,$ 
and the case $l=2$ is obvious.

Applying (\ref{after9031}) 
to $C_{j_k \ldots j_{l+1} j_l j_l j_{l-2} \ldots \ldots j_1}$, we obtain
for $l+1\le k,$ $l-2\ge 1$
$$
\sum_{j_l=p+1}^{\infty}
C_{j_k \ldots j_{l+1} j_l j_{l} j_{l-2} \ldots \ldots j_1}=
$$
$$
=\int\limits_{[t, T]^{k-2}} \sum\limits_{d=1}^4
H_p^{(d)}(t_1,\ldots,t_{l-2},t_{l+1},\ldots,t_k)
\prod_{\stackrel{g=1}{{}_{g\ne l-1, l}}}^{k}
\psi_{g}(t_g)\phi_{j_g}(t_g)\times
$$

\vspace{-2mm}
$$
\times dt_1\ldots dt_{l-2}dt_{l+1}\ldots dt_k=
$$

\vspace{-2mm}
\begin{equation}
\label{after85ayyy}
= \sum\limits_{d=1}^4
C_{j_k \ldots j_{l+1} j_{l-2} \ldots j_1}^{**(d)}
=\sum\limits_{d=1}^4
C_{j_k \ldots j_q \ldots j_1}^{**(d)}\biggl|_{q\ne l-1, l}\biggr.,
\end{equation}

\vspace{4mm}
\noindent
where
$$
H_p^{(1)}(t_1,\ldots,t_{l-2},t_{l+1},\ldots,t_k)=
$$
\begin{equation}
\label{after90yyy}
=
{\bf 1}_{\{t_1<\ldots<t_{l-2}<t_{l+1}<\ldots<t_k\}}
\sum_{j_l=p+1}^{\infty}~
\int\limits_{t}^{t_{l+1}} \psi_l(\tau) \phi_{j_{l}}(\tau)d\tau
\int\limits_{t}^{t_{l+1}} \psi_{l-1}(\tau) \phi_{j_{l}}(\tau)d\tau,
\end{equation}

$$
H_p^{(2)}(t_1,\ldots,t_{l-2},t_{l+1},\ldots,t_k)=
$$
\begin{equation}
\label{after91yyy}
=
-{\bf 1}_{\{t_1<\ldots<t_{l-2}<t_{l+1}<\ldots<t_k\}}
\sum_{j_l=p+1}^{\infty}~
\int\limits_{t}^{t_{l+1}} \psi_l(\tau) \phi_{j_{l}}(\tau)d\tau
\int\limits_{t}^{t_{l-2}} \psi_{l-1}(\tau) \phi_{j_{l}}(\tau)d\tau,
\end{equation}

$$
H_p^{(3)}(t_1,\ldots,t_{l-2},t_{l+1},\ldots,t_k)=
$$
\begin{equation}
\label{after92yyy}
=
-{\bf 1}_{\{t_1<\ldots<t_{l-2}<t_{l+1}<\ldots<t_k\}}
\sum_{j_l=p+1}^{\infty}~
\int\limits_{t}^{t_{l+1}} \psi_{l-1}(\tau)\phi_{j_{l}}(\tau)
\int\limits_{t}^{\tau} \psi_l(\theta)\phi_{j_{l}}(\theta)d\theta d\tau,
\end{equation}

$$
H_p^{(4)}(t_1,\ldots,t_{l-2},t_{l+1},\ldots,t_k)=
$$
\begin{equation}
\label{after93yyy}
=
{\bf 1}_{\{t_1<\ldots<t_{l-2}<t_{l+1}<\ldots<t_k\}}
\sum_{j_l=p+1}^{\infty}~
\int\limits_{t}^{t_{l-2}} \psi_{l-1}(\tau)\phi_{j_{l}}(\tau)
\int\limits_{t}^{\tau} \psi_l(\theta)\phi_{j_{l}}(\theta)d\theta d\tau.
\end{equation}

\vspace{1mm}

By analogy with (\ref{after85ayyy}) we can consider the expressions
\begin{equation}
\label{afterx300}
\sum_{j_l=p+1}^{\infty}
C_{j_k \ldots j_{l+1} j_l j_l},
\end{equation}
\begin{equation}
\label{afterx301}
\sum_{j_l=p+1}^{\infty}
C_{j_l j_l j_{k-2} \ldots j_1}.
\end{equation}

\vspace{2mm}

Then we have for (\ref{afterx300}), (\ref{afterx301}) 
(see (\ref{after9031}) and its analogue for $t_{l+1}=T$)
\begin{equation}
\label{afterx500}
~~~~~~\sum_{j_l=p+1}^{\infty}
C_{j_k \ldots j_{l+1} j_l j_{l}}=
\int\limits_{[t, T]^{k-2}} 
L_p(t_3,\ldots,t_k)\prod_{g=3}^{k}
\psi_{g}(t_g)\phi_{j_g}(t_g)
dt_3\ldots dt_k,
\end{equation}
\begin{equation}
\label{afterx501}
\sum_{j_l=p+1}^{\infty}
C_{j_l j_l j_{k-2} \ldots j_{1}}=
\int\limits_{[t, T]^{k-2}} 
\sum\limits_{d=1}^4 M_p^{(d)}(t_1,\ldots,t_{k-2})\prod_{g=1}^{k-2}
\psi_{g}(t_g)\phi_{j_g}(t_g)
dt_1\ldots dt_{k-2},
\end{equation}
where
$$
L_p(t_3,\ldots,t_k)
=
{\bf 1}_{\{t_3<\ldots<t_k\}}
\sum_{j_l=p+1}^{\infty}~
\int\limits_{t}^{t_{3}} \psi_2(\tau) \phi_{j_{l}}(\tau)
\int\limits_{t}^{\tau} \psi_{1}(\theta) \phi_{j_{l}}(\theta)d\theta d\tau,
$$

$$
M_p^{(1)}(t_1,\ldots,t_{k-2})=
$$
$$
=
{\bf 1}_{\{t_1<\ldots<t_{k-2}\}}
\sum_{j_l=p+1}^{\infty}~
\int\limits_{t}^{T} \psi_k(\tau) \phi_{j_{l}}(\tau)d\tau
\int\limits_{t}^{T} \psi_{k-1}(\tau) \phi_{j_{l}}(\tau)d\tau,
$$

$$
M_p^{(2)}(t_1,\ldots,t_{k-2})=
$$
$$
=
-{\bf 1}_{\{t_1<\ldots<t_{k-2}\}}
\sum_{j_l=p+1}^{\infty}~
\int\limits_{t}^{T} \psi_k(\tau) \phi_{j_{l}}(\tau)d\tau
\int\limits_{t}^{t_{k-2}} \psi_{k-1}(\tau) \phi_{j_{l}}(\tau)d\tau,
$$

$$
M_p^{(3)}(t_1,\ldots,t_{k-2})=
$$
$$
=
-{\bf 1}_{\{t_1<\ldots<t_{k-2}\}}
\sum_{j_l=p+1}^{\infty}~
\int\limits_{t}^{T} \psi_{k-1}(\tau) \phi_{j_{l}}(\tau)
\int\limits_{t}^{\tau} \psi_{k}(\theta) \phi_{j_{l}}(\theta)d\theta d\tau,
$$

$$
M_p^{(4)}(t_1,\ldots,t_{k-2})=
$$
$$
=
{\bf 1}_{\{t_1<\ldots<t_{k-2}\}}
\sum_{j_l=p+1}^{\infty}~
\int\limits_{t}^{t_{k-2}} \psi_{k-1}(\tau) \phi_{j_{l}}(\tau)
\int\limits_{t}^{\tau} \psi_{k}(\theta) \phi_{j_{l}}(\theta)d\theta d\tau.
$$

It is important to note that 
$C_{j_k \ldots j_{l+1} j_{l-2} \ldots j_1}^{*(d)},$
$C_{j_k \ldots j_{l+1} j_{l-2} \ldots j_1}^{**(d)}$
$(d=1,\ldots,4)$ 
are Fourier coefficients (see (\ref{after85a}), (\ref{after85ayyy})), 
that is, we can use Parseval's equality 
in the further proof.

Combining the equalities (\ref{after85a})--(\ref{after93})
(the case $g_2>g_1+1$), using Parseval's equality and applying the 
estimates for integrals from basis functions that we 
used in the proof of Theorems 2.33, 2.34, we obtain for (\ref{after85a}) 
$$
\sum\limits_{j_{q_1},\ldots,j_{q_{k-2}}=0}^p
\left(\sum_{j_{g_1}=p+1}^{\infty}
C_{j_k\ldots j_1}\biggl|_{j_{g_1}=j_{g_2}, g_2>g_1+1}\biggr.\right)^2=
$$
$$
=
\sum\limits_{\stackrel{j_{1},\ldots,j_q,\ldots,j_{k}=0}{{}_{q\ne g_1, g_2}}}^p
\left(\sum_{j_{g_1}=p+1}^{\infty}
C_{j_k\ldots j_1}\biggl|_{j_{g_1}=j_{g_2}, g_2>g_1+1}\biggr.\right)^2=
$$

\vspace{-3mm}
$$
= 
\sum\limits_{\stackrel{j_{1},\ldots,j_q,\ldots,j_{k}=0}{{}_{q\ne g_1, g_2}}}^p
\left(\sum\limits_{d=1}^4 
C_{j_k \ldots j_q \ldots j_1}^{*(d)}\biggl|_{q\ne g_1, g_2}\biggr.\right)^2
\le
\sum\limits_{\stackrel{j_{1},\ldots,j_q,\ldots,j_{k}=0}{{}_{q\ne g_1, g_2}}}^{\infty}
\left(\sum\limits_{d=1}^4
C_{j_k \ldots j_q \ldots j_1}^{*(d)}\biggl|_{q\ne g_1, g_2}\biggr.\right)^2
=
$$
$$
=\sum\limits_{\stackrel{j_{1},\ldots,j_q,\ldots,j_{k}=0}{{}_{q\ne g_1, g_2}}}^{\infty}
\left(~
\int\limits_{[t, T]^{k-2}} \sum\limits_{d=1}^4
F_p^{(d)}(t_1,\ldots,t_{g_1-1},t_{g_1+1},\ldots ,t_{g_2-1},t_{g_2+1},\ldots,t_k)\times\right.
$$
$$
\left.\times 
\prod_{\stackrel{q=1}{{}_{q\ne g_1, g_2}}}^{k}
\psi_{q}(t_q)\phi_{j_q}(t_q)
dt_1\ldots dt_{g_1-1}dt_{g_1+1}\ldots dt_{g_2-1}dt_{g_2+1}\ldots dt_k\right)^2=
$$
$$
=\int\limits_{[t, T]^{k-2}} 
\left(~
\sum\limits_{d=1}^4
F_p^{(d)}(t_1,\ldots,t_{g_1-1},t_{g_1+1},\ldots, t_{g_2-1},t_{g_2+1},\ldots,t_k)
\prod_{\stackrel{q=1}{{}_{q\ne g_1, g_2}}}^{k}
\psi_{q}(t_q)\right)^2\times
$$

$$
\times
dt_1\ldots dt_{g_1-1}dt_{g_1+1}\ldots dt_{g_2-1}dt_{g_2+1}\ldots dt_k\le
$$

\vspace{-5mm}
$$
\le 4\sum\limits_{d=1}^4 \int\limits_{[t, T]^{k-2}} 
\left(
F_p^{(d)}(t_1,\ldots,t_{g_1-1},t_{g_1+1},\ldots, t_{g_2-1},t_{g_2+1},\ldots,t_k)
\prod_{\stackrel{q=1}{{}_{q\ne g_1, g_2}}}^{k}
\psi_{q}(t_q)\right)^2\times
$$

\begin{equation}
\label{after12001}
~~~~~~~~\times dt_1\ldots dt_{g_1-1}dt_{g_1+1}\ldots dt_{g_2-1}dt_{g_2+1}\ldots dt_k\le
\frac{K}{p^{2-\varepsilon}}\ \to\  0
\end{equation}

\vspace{3mm}
\noindent
if $p\to\infty,$ where $\varepsilon$ is an arbitrary small positive real number
for the polynomial case and $\varepsilon=0$ for the 
trigonometric case, constant $K$ does not depend on $p.$
The cases (\ref{afterx200})--(\ref{afterx202}) are considered
analogously.

Absolutely similarly (see (\ref{after12001}))
combining the equalities (\ref{after85ayyy})--(\ref{after93yyy})
(the case $g_2=g_1+1$), using Parseval's equality and applying the 
estimates for integrals from basis functions that we 
used in the proof of Theorems 2.33, 2.34, we get for (\ref{after85ayyy})
$$
\sum\limits_{j_{q_1},\ldots,j_{q_{k-2}}=0}^p
\left(\sum_{j_{g_1}=p+1}^{\infty}
C_{j_k\ldots j_1}\biggl|_{j_{g_1}=j_{g_2}, g_2=g_1+1}\biggr.\right)^2=
$$
$$
=\sum\limits_{\stackrel{j_{1},\ldots,j_q,\ldots,j_{k}=0}{{}_{q\ne g_1, g_2}}}^p
\left(\sum_{j_{g_1}=p+1}^{\infty}
C_{j_k\ldots j_1}\biggl|_{j_{g_1}=j_{g_2}, g_2=g_1+1}\biggr.\right)^2=
$$
$$
= 
\sum\limits_{\stackrel{j_{1},\ldots,j_q,\ldots,j_{k}=0}{{}_{q\ne g_1, g_2}}}^p
\left(\sum\limits_{d=1}^4 
C_{j_k \ldots j_q \ldots j_1}^{**(d)}\biggl|_{q\ne g_1, g_2}\biggr.\right)^2
\le
\sum\limits_{\stackrel{j_{1},\ldots,j_q,\ldots,j_{k}=0}{{}_{q\ne g_1, g_2}}}^{\infty}
\left(\sum\limits_{d=1}^4
C_{j_k \ldots j_q \ldots j_1}^{**(d)}\biggl|_{q\ne g_1, g_2}\biggr.\right)^2
=
$$
$$
=\sum\limits_{\stackrel{j_{1},\ldots,j_q,\ldots,j_{k}=0}{{}_{q\ne g_1, g_2}}}^{\infty}
\left(~
\int\limits_{[t, T]^{k-2}} \sum\limits_{d=1}^4
H_p^{(d)}(t_1,\ldots,t_{g_1-1},t_{g_1+2},\ldots,t_k)\times\right.
$$
$$
\left.\times 
\prod_{\stackrel{q=1}{{}_{q\ne g_1, g_1+1}}}^{k}
\psi_{q}(t_q)\phi_{j_q}(t_q)
dt_1\ldots dt_{g_1-1}dt_{g_1+2}\ldots dt_k\right)^2=
$$
$$
=\int\limits_{[t, T]^{k-2}} 
\hspace{-1.5mm}\left(~\hspace{-1.5mm}
\sum\limits_{d=1}^4
H_p^{(d)}(t_1,\ldots,t_{g_1-1},t_{g_1+2},\ldots,t_k)\prod_{\stackrel{q=1}{{}_{q\ne g_1, g_1+1}}}^{k}
\psi_{q}(t_q)\right)^2\times
$$

$$
\times
dt_1\ldots dt_{g_1-1}dt_{g_1+2}\ldots \ldots dt_k\le
$$

\vspace{-5mm}
$$
\le 4\sum\limits_{d=1}^4 \int\limits_{[t, T]^{k-2}} 
\left(
H_p^{(d)}(t_1,\ldots,t_{g_1-1},t_{g_1+2},\ldots,t_k)\prod_{\stackrel{q=1}{{}_{q\ne g_1, g_1+1}}}^{k}
\psi_{q}(t_q)\right)^2\times
$$

\begin{equation}
\label{after12000}
\times
dt_1\ldots dt_{g_1-1}dt_{g_1+2}\ldots dt_k\le
\frac{K}{p^{2-\varepsilon}}\ \to\  0
\end{equation}

\vspace{3mm}
\noindent
if $p\to\infty,$ where $\varepsilon$ is an arbitrary small positive real number
for the polynomial case and $\varepsilon=0$ for the 
trigonometric case, constant $K$ does not depend on $p.$
The cases (\ref{afterx300}), (\ref{afterx301}) are considered
analogously.

From (\ref{after12001}), (\ref{after12000}) and their
analogues for the cases (\ref{afterx200})--(\ref{afterx202}),
(\ref{afterx300}), (\ref{afterx301})
we obtain 
\begin{equation}
\label{after17000}
\sum\limits_{j_{q_1},\ldots,j_{q_{k-2}}=0}^p
\left(\sum_{j_{g_1}=p+1}^{\infty}
C_{j_k\ldots j_1}\biggl|_{j_{g_1}=j_{g_2}}\biggr.\right)^2
\le \frac{K}{p^{2-\varepsilon}},
\end{equation}

\vspace{3mm}
\noindent
where constant $K$ is independent of $p.$ 
Thus the equality (\ref{after14000}) is proved.

Let us prove the equality (\ref{after14001}).
Consider the following cases

\vspace{4mm}
\centerline{1.\ $g_2>g_1+1$,\ $g_4=g_3+1$,\ \ \ 2.\ $g_2=g_1+1$,\ $g_4>g_3+1$,}

\vspace{2mm}
\centerline{\hspace{0.1mm}3.\ $g_2>g_1+1$,\ $g_4>g_3+1$,\ \ \ \ \hspace{-2.0mm}4.\ 
$g_2=g_1+1$,\ $g_4=g_3+1$.}

\vspace{4mm}
The proof for Cases 1--3 will be similar. Consider, for example, Case 2.
Using (\ref{after79}), we obtain
$$
\sum\limits_{j_{q_1}=0}^p
\left(\sum_{j_{g_1}=p+1}^{\infty}\sum_{j_{g_3}=p+1}^{\infty}
C_{j_5\ldots j_1}\biggl|_{j_{g_1}=j_{g_2},j_{g_3}=j_{g_4}, g_4>g_3+1, g_2=g_1+1}\biggr.\right)^2=
$$
$$
=
\sum\limits_{j_{q_1}=0}^p
\left(\sum_{j_{g_1}=p+1}^{\infty}\sum_{j_{g_3}=0}^{p}
C_{j_5\ldots j_1}\biggl|_{j_{g_1}=j_{g_2},j_{g_3}=j_{g_4}, g_4>g_3+1, g_2=g_1+1}\biggr.\right)^2=
$$
\begin{equation}
\label{afterpp1}
~~~~~~~~~=\sum\limits_{j_{q_1}=0}^p
\left(\sum_{j_{g_3}=0}^{p}\sum_{j_{g_1}=p+1}^{\infty}
C_{j_5\ldots j_1}\biggl|_{j_{g_1}=j_{g_2},j_{g_3}=j_{g_4}, g_4>g_3+1, g_2=g_1+1}\biggr.\right)^2\le
\end{equation}
$$
\le (p+1)\sum\limits_{j_{q_1}=0}^p\sum_{j_{g_3}=0}^{p}
\left(\sum_{j_{g_1}=p+1}^{\infty}
C_{j_5\ldots j_1}\biggl|_{j_{g_1}=j_{g_2},j_{g_3}=j_{g_4}, g_4>g_3+1, g_2=g_1+1}\biggr.\right)^2=
$$
$$
=(p+1)\sum\limits_{j_{q_1}=0}^p\sum_{j_{g_3}, j_{g_4}=0}^{p}
\left(\sum_{j_{g_1}=p+1}^{\infty}
C_{j_5\ldots j_1}\biggl|_{j_{g_1}=j_{g_2}, g_4>g_3+1, g_2=g_1+1}\biggr.
\right)^2\Biggl|_{j_{g_3}=j_{g_4}}\le
$$
\begin{equation}
\label{after15000}
~~~~~\le(p+1)\sum\limits_{j_{q_1}=0}^p\sum_{j_{g_3}, j_{g_4}=0}^{p}
\left(\sum_{j_{g_1}=p+1}^{\infty}
C_{j_5\ldots j_1}\biggl|_{j_{g_1}=j_{g_2}, g_4>g_3+1, g_2=g_1+1}\biggr.
\right)^2.
\end{equation}

\vspace{3mm}

It is easy to see that the expression 
(\ref{after15000}) (without the multiplier $p+1$) 
is a particular case $(k=5, g_4>g_3+1, g_2=g_1+1)$ of the left-hand side
of (\ref{after17000}).
Combining (\ref{after17000}) and (\ref{after15000}), we have
\begin{equation}
\label{after15000a}
\sum\limits_{j_{q_1}=0}^p
\left(\sum_{j_{g_1}=p+1}^{\infty}\sum_{j_{g_3}=p+1}^{\infty}
C_{j_5\ldots j_1}\biggl|_{j_{g_1}=j_{g_2},j_{g_3}=j_{g_4}, g_4>g_3+1, g_2=g_1+1}\biggr.
\right)^{\hspace{-0.7mm}2}
\hspace{-1.7mm}\le
\frac{(p+1)K}{p^{2-\varepsilon}}\le
\frac{K_1}{p^{1-\varepsilon}} \to  \hspace{-0.2mm}0
\end{equation}

\vspace{3mm}
\noindent
if $p\to\infty,$ where 
$\varepsilon$ is an arbitrary small positive real number
for the polynomial case and $\varepsilon=0$ for the 
trigonometric case,
constant $K_1$ does not depend on $p.$

Consider Case 4 ($g_2=g_1+1$,\ $g_4=g_3+1$). We have (see (\ref{after500}))
$$
\sum\limits_{j_{q_1}=0}^p
\left(\sum_{j_{g_1}=p+1}^{\infty}\sum_{j_{g_3}=p+1}^{\infty}
C_{j_5\ldots j_1}\biggl|_{j_{g_1}=j_{g_2},j_{g_3}=j_{g_4}}\biggr.\right)^2=
$$
$$
=\sum\limits_{j_{q_1}=0}^p
\left(\sum_{j_{g_1}=p+1}^{\infty}\left(\sum_{j_{g_3}=0}^{\infty}-\sum_{j_{g_3}=0}^{p}\right)
C_{j_5\ldots j_1}\biggl|_{j_{g_1}=j_{g_2},j_{g_3}=j_{g_4}}\biggr.\right)^2=
$$
$$
=\sum\limits_{j_{q_1}=0}^p
\left(\frac{1}{2}
\sum_{j_{g_1}=p+1}^{\infty}
C_{j_5\ldots j_1}\biggl|_{j_{g_1}=j_{g_2},
(j_{g_3} j_{g_3})\curvearrowright (\cdot)}\biggr.
-\sum_{j_{g_3}=0}^{p}\sum_{j_{g_1}=p+1}^{\infty}
C_{j_5\ldots j_1}\biggl|_{j_{g_1}=j_{g_2},j_{g_3}=j_{g_4}}\biggr.\right)^2\le
$$
\begin{equation}
\label{after19000}
\le
\frac{1}{2}\sum\limits_{j_{q_1}=0}^p
\left(
\sum_{j_{g_1}=p+1}^{\infty}
C_{j_5\ldots j_1}\biggl|_{j_{g_1}=j_{g_2},
(j_{g_3} j_{g_3})\curvearrowright (\cdot)}\biggr.\right)^2+
\end{equation}
\begin{equation}
\label{after19001}
+2\sum\limits_{j_{q_1}=0}^p
\left(\sum_{j_{g_3}=0}^{p}\sum_{j_{g_1}=p+1}^{\infty}
C_{j_5\ldots j_1}\biggl|_{j_{g_1}=j_{g_2},j_{g_3}=j_{g_4}}\biggr.\right)^2.
\end{equation}

\vspace{3mm}

An expression similar to (\ref{after19001}) was estimated 
(see (\ref{afterpp1})--(\ref{after15000a})). Let us estimate 
(\ref{after19000}). We have
$$
\sum\limits_{j_{q_1}=0}^p
\left(
\sum_{j_{g_1}=p+1}^{\infty}
C_{j_5\ldots j_1}\biggl|_{j_{g_1}=j_{g_2},
(j_{g_3} j_{g_3})\curvearrowright (\cdot)}\biggr.\right)^2=
$$
$$
=(T-t)\sum\limits_{j_{q_1}=0}^p
\left(
\sum_{j_{g_1}=p+1}^{\infty}
C_{j_5\ldots j_1}\biggl|_{j_{g_1}=j_{g_2},
(j_{g_3} j_{g_3})\curvearrowright 0}\biggr.\right)^2\le
$$
\begin{equation}
\label{after19005}
~~~~~~~\le (T-t)\sum\limits_{j_{q_1}=0}^p
\sum_{j_{g_3}=0}^{p}\left(
\sum_{j_{g_1}=p+1}^{\infty}
C_{j_5\ldots j_1}\biggl|_{j_{g_1}=j_{g_2},
(j_{g_3} j_{g_3})\curvearrowright j_{g_3}}\biggr.\right)^2,
\end{equation}

\vspace{3mm}
\noindent
where the notations are the same as in the proof of Theorem~2.30.

The expression (\ref{after19005}) 
without the multiplier $T-t$ is an expression of type 
(\ref{after2501})--(\ref{after2506})
before passing to the limit $\lim\limits_{p\to\infty}$ 
(the only difference is the replacement of one of the weight functions 
$\psi_1(\tau),\ldots, \psi_4(\tau)$ in 
(\ref{after2501})--(\ref{after2506}) by the product 
$\psi_{l+1}(\tau)\psi_l(\tau)$ $(l=1,\ldots,4).$
Therefore, for Case 4 ($g_2=g_1+1$,\ $g_4=g_3+1$), we obtain the estimate
\begin{equation}
\label{after20000}
~~~~~~~\sum\limits_{j_{q_1}=0}^p
\left(\sum_{j_{g_1}=p+1}^{\infty}\sum_{j_{g_3}=p+1}^{\infty}
C_{j_5\ldots j_1}\biggl|_{j_{g_1}=j_{g_2},j_{g_3}=j_{g_4}, g_4=g_3+1, g_2=g_1+1}\biggr.\right)^2\le
\frac{K}{p^{1-\varepsilon}},
\end{equation}

\vspace{3mm}
\noindent
where 
$\varepsilon$ is an arbitrary small positive real number
for the polynomial case and $\varepsilon=0$ for the 
trigonometric case,
constant $K$ is independent of $p.$

The estimates (\ref{after15000a}), (\ref{after20000}) 
prove (\ref{after14001}).

Let us prove (\ref{afterafter001}). By analogy with (\ref{after19005}) we have
$$
\sum\limits_{j_{q_1}=0}^p
\left(\sum_{j_{g_3}=p+1}^{\infty}
C_{j_5\ldots j_1}\biggl|_{(j_{g_2}j_{g_1})\curvearrowright (\cdot),
j_{g_1}=j_{g_2},
j_{g_3}=j_{g_4}, g_2=g_1+1}\biggr.\right)^2=
$$
$$
=\sum\limits_{j_{q_1}=0}^p
\left(\sum_{j_{g_3}=p+1}^{\infty}
C_{j_5\ldots j_1}\biggl|_{(j_{g_1}j_{g_1})\curvearrowright (\cdot),
j_{g_3}=j_{g_4}, g_2=g_1+1}\biggr.\right)^2=
$$
$$
=(T-t)\sum\limits_{j_{q_1}=0}^p
\left(\sum_{j_{g_3}=p+1}^{\infty}
C_{j_5\ldots j_1}\biggl|_{(j_{g_1}j_{g_1})\curvearrowright 0,
j_{g_3}=j_{g_4}, g_2=g_1+1}\biggr.\right)^2\le
$$
\begin{equation}
\label{afterafter120}
~~~~~~~~\le
(T-t)\sum\limits_{j_{q_1}=0}^p \sum\limits_{j_{g_1}=0}^p
\left(\sum_{j_{g_3}=p+1}^{\infty}
C_{j_5\ldots j_1}\biggl|_{(j_{g_1}j_{g_1})\curvearrowright j_{g_1},
j_{g_3}=j_{g_4}, g_2=g_1+1}\biggr.\right)^2.
\end{equation}

\vspace{3mm}

Thus, we obtain the estimate (see (\ref{after19005}) and the proof of Theorem~2.34)
\begin{equation}
\label{afterafter121}
~~~~~~~~~~\sum\limits_{j_{q_1}=0}^p
\left(\sum_{j_{g_3}=p+1}^{\infty}
C_{j_5\ldots j_1}\biggl|_{(j_{g_2}j_{g_1})\curvearrowright (\cdot),
j_{g_1}=j_{g_2},
j_{g_3}=j_{g_4}, g_2=g_1+1}\biggr.\right)^2
\le \frac{K}{p^{2-\varepsilon}},
\end{equation}

\vspace{3mm}
\noindent
where 
$\varepsilon$ is an arbitrary small positive real number
for the polynomial case and $\varepsilon=0$ for the 
trigonometric case,
constant $K$ does not depend on $p.$

The estimate (\ref{afterafter121}) proves (\ref{afterafter001}).
Theorem~2.35 is proved.

\section{Expansion of Iterated Stratonovich Stochastic Integrals
of Multiplicity 6. The Case $p_1=\ldots =p_6\to \infty$ and 
$\psi_1(\tau),\ldots,\psi_6(\tau)\equiv 1$ 
(The Cases of Legendre 
Polynomials and Trigonometric Functions)}

{\bf Theorem 2.36}\ \cite{arxiv-5}, \cite{arxiv-10}, \cite{arxiv-11}, \cite{new-art-1xxys}.\
{\it Suppose that 
$\{\phi_j(x)\}_{j=0}^{\infty}$ is a complete orthonormal system of 
Legendre polynomials or trigonometric functions in the space $L_2([t, T]).$
Then$,$ for the 
ite\-ra\-ted Stra\-to\-no\-vich stochastic integral of sixth multiplicity
\begin{equation}
\label{after10001qu1}
J_{T,t}^{*(i_1\ldots i_6)}={\int\limits_t^{*}}^T
\ldots
{\int\limits_t^{*}}^{t_2}
d{\bf w}_{t_1}^{(i_1)}
\ldots d{\bf w}_{t_6}^{(i_6)}
\end{equation}
the following 
expansion 
$$
J_{T,t}^{*(i_1\ldots i_6)}
=\hbox{\vtop{\offinterlineskip\halign{
\hfil#\hfil\cr
{\rm l.i.m.}\cr
$\stackrel{}{{}_{p\to \infty}}$\cr
}} }
\sum\limits_{j_1, \ldots, j_6=0}^{p}
C_{j_6 \ldots j_1}\zeta_{j_1}^{(i_1)}\ldots
\zeta_{j_6}^{(i_6)}
$$
that converges in the mean-square sense is valid, where
$i_1, \ldots, i_6=0, 1,\ldots,m,$
$$
C_{j_6 \ldots j_1}=
\int\limits_t^T\phi_{j_6}(t_6)\ldots
\int\limits_t^{t_2}\phi_{j_1}(t_1)dt_1\ldots dt_6
$$
and
$$
\zeta_{j}^{(i)}=
\int\limits_t^T \phi_{j}(s) d{\bf w}_s^{(i)}
$$ 
are independent standard Gaussian random variables for various 
$i$ or $j$
{\rm (}in the case when $i\ne 0${\rm ),}
${\bf w}_{\tau}^{(i)}$ 
$(i=1,\ldots,m)$ are independent standard Wiener processes$,$
${\bf w}_{\tau}^{(0)}=\tau.$}

\vspace{2mm}

{\bf Proof.}\ 
As noted in Remark~2.4, Conditions 1 and 2
of Theorem~{\rm 2.30} are satisfied for complete
orthonormal systems of Legendre polynomials 
and trigonometric functions in the space
$L_2([t, T]).$ Let us verify Condition 3 of Theorem~2.30 for
the iterated Stratonovich stochastic integral (\ref{after10001qu1}). 
Thus, we have to check the following conditions
\begin{equation}
\label{after14000qu1}
\lim\limits_{p\to\infty}
\sum\limits_{j_{q_1},j_{q_2},j_{q_3},j_{q_4}=0}^p
\left(\sum_{j_{g_1}=p+1}^{\infty}
C_{j_6\ldots j_1}\biggl|_{j_{g_1}=j_{g_2}}\biggr.\right)^2=0,
\end{equation}
\begin{equation}
\label{after14001qu2}
~~~~~~~~~\lim\limits_{p\to\infty}
\sum\limits_{j_{q_1},j_{q_2}=0}^p
\left(\sum_{j_{g_1}=p+1}^{\infty}\sum_{j_{g_3}=p+1}^{\infty}
C_{j_6\ldots j_1}\biggl|_{j_{g_1}=j_{g_2},j_{g_3}=j_{g_4}}\biggr.\right)^2=0,
\end{equation}
\begin{equation}
\label{afterafter001qu3}
~~~~~~~~~~\lim\limits_{p\to\infty}
\sum\limits_{j_{q_1},j_{q_2}=0}^p
\left(\sum_{j_{g_1}=p+1}^{\infty}
C_{j_6\ldots j_1}\biggl|_{(j_{g_4}j_{g_3})\curvearrowright (\cdot),
j_{g_1}=j_{g_2},
j_{g_3}=j_{g_4}, g_4=g_3+1}\biggr.\right)^2=0,
\end{equation}
\begin{equation}
\label{after14001qu10}
~~~~~~~~~\lim\limits_{p\to\infty}
\left(\sum_{j_{g_1}=p+1}^{\infty}\sum_{j_{g_3}=p+1}^{\infty}
\sum_{j_{g_5}=p+1}^{\infty}
C_{j_6\ldots j_1}\biggl|_{j_{g_1}=j_{g_2},j_{g_3}=j_{g_4},j_{g_5}=j_{g_6}}\biggr.\right)^2=0,
\end{equation}
\begin{equation}
\label{afterafter001qu33}
~~~\lim\limits_{p\to\infty}
\left(\sum_{j_{g_1}=p+1}^{\infty}\sum_{j_{g_3}=p+1}^{\infty}
C_{j_6\ldots j_1}\biggl|_{(j_{g_6}j_{g_5})\curvearrowright (\cdot),
j_{g_1}=j_{g_2},
j_{g_3}=j_{g_4}, j_{g_5}=j_{g_6}, g_6=g_5+1}\biggr.\right)^2=0,
\end{equation}
\begin{equation}
\label{afterafter001qu36}
\lim\limits_{p\to\infty}
\left(\sum_{j_{g_1}=p+1}^{\infty}
C_{j_6\ldots j_1}\biggl|_{(j_{g_4}j_{g_3})\curvearrowright (\cdot)
(j_{g_6}j_{g_5})\curvearrowright (\cdot),
j_{g_1}=j_{g_2},
j_{g_3}=j_{g_4}, j_{g_5}=j_{g_6}, g_4=g_3+1, g_6=g_5+1}\biggr.\right)^2=0,
\end{equation}

\noindent
where the expressions
$\left(\{g_1,g_2\},\{g_3,g_4\}, \{g_5,g_6\}\}\right),$
$\left(\{g_1,g_2\},\{g_3,g_4\}, \{q_1,q_2\}\}\right),$
$\left(\{g_1,g_2\}, \{q_1, q_2,q_3, q_4\}\right)$
are partitions of the set $\{1,2,\ldots,6\}$ that is
$\{g_1,g_2,g_3,$ $g_4,g_5,g_6\}=\{g_1,g_2,g_3,g_4,q_1,q_2\}=
\{g_1,g_2,q_1,q_2,q_3,q_4\}=
\{1,2,$ $\ldots,6\};$
braces mean an unordered 
set, and pa\-ren\-the\-ses mean an ordered set.

The equalities (\ref{after14000qu1}),
(\ref{afterafter001qu3}) were proved earlier (see the proof of equalities 
(\ref{after17000}), (\ref{after19005})).
The relation (\ref{afterafter001qu36}) follows 
from
the estimate (\ref{tupo15}) for the polynomial case and its analogue
for the
trigonometric case.
It is easy to see that the 
equalities (\ref{after14001qu2}) and 
(\ref{afterafter001qu33}) are proved in complete analogy with the proof of 
(\ref{after14001}), (\ref{after19005}).

Thus, we have to prove the relation (\ref{after14001qu10}).
The equality (\ref{after14001qu10}) is equivalent to the following equalities
\begin{equation}
\label{sixsix8}
\lim\limits_{p\to\infty}
\sum_{j_1=p+1}^{\infty}\sum_{j_2=p+1}^{\infty}
\sum_{j_3=p+1}^{\infty}
C_{j_3 j_2 j_1 j_3 j_2 j_1}=0,
\end{equation}
\begin{equation}
\label{sixsix9}
\lim\limits_{p\to\infty}
\sum_{j_1=p+1}^{\infty}\sum_{j_2=p+1}^{\infty}
\sum_{j_3=p+1}^{\infty}
C_{j_1 j_3 j_2 j_3 j_2 j_1}=0,
\end{equation}
\begin{equation}
\label{sixsix10}
\lim\limits_{p\to\infty}
\sum_{j_1=p+1}^{\infty}\sum_{j_2=p+1}^{\infty}
\sum_{j_3=p+1}^{\infty}
C_{j_3 j_2 j_3 j_1 j_2 j_1}=0,
\end{equation}
\begin{equation}
\label{sixsix4}
\lim\limits_{p\to\infty}
\sum_{j_1=p+1}^{\infty}\sum_{j_2=p+1}^{\infty}
\sum_{j_3=p+1}^{\infty}
C_{j_1 j_2 j_3 j_3 j_2 j_1}=0,
\end{equation}
\begin{equation}
\label{sixsix14}
\lim\limits_{p\to\infty}
\sum_{j_1=p+1}^{\infty}\sum_{j_2=p+1}^{\infty}
\sum_{j_3=p+1}^{\infty}
C_{j_1 j_2 j_2 j_3 j_3 j_1}=0,
\end{equation}
\begin{equation}
\label{sixsix3}
\lim\limits_{p\to\infty}
\sum_{j_1=p+1}^{\infty}\sum_{j_2=p+1}^{\infty}
\sum_{j_3=p+1}^{\infty}
C_{j_3 j_3 j_2 j_2 j_1 j_1}=0,
\end{equation}
\begin{equation}
\label{sixsix7}
\lim\limits_{p\to\infty}
\sum_{j_1=p+1}^{\infty}\sum_{j_2=p+1}^{\infty}
\sum_{j_3=p+1}^{\infty}
C_{j_2 j_3 j_3 j_2 j_1 j_1}=0,
\end{equation}
\begin{equation}
\label{sixsix6}
\lim\limits_{p\to\infty}
\sum_{j_1=p+1}^{\infty}\sum_{j_2=p+1}^{\infty}
\sum_{j_3=p+1}^{\infty}
C_{j_3 j_2 j_3 j_2 j_1 j_1}=0,
\end{equation}
\begin{equation}
\label{sixsix1}
\lim\limits_{p\to\infty}
\sum_{j_1=p+1}^{\infty}\sum_{j_2=p+1}^{\infty}
\sum_{j_3=p+1}^{\infty}
C_{j_3 j_3 j_2 j_1 j_2 j_1}=0,
\end{equation}
\begin{equation}
\label{sixsix2}
\lim\limits_{p\to\infty}
\sum_{j_1=p+1}^{\infty}\sum_{j_2=p+1}^{\infty}
\sum_{j_3=p+1}^{\infty}
C_{j_3 j_3 j_1 j_2 j_2 j_1}=0,
\end{equation}
\begin{equation}
\label{sixsix5}
\lim\limits_{p\to\infty}
\sum_{j_1=p+1}^{\infty}\sum_{j_2=p+1}^{\infty}
\sum_{j_3=p+1}^{\infty}
C_{j_2 j_1 j_3 j_3 j_2 j_1}=0,
\end{equation}
\begin{equation}
\label{sixsix12}
\lim\limits_{p\to\infty}
\sum_{j_1=p+1}^{\infty}\sum_{j_2=p+1}^{\infty}
\sum_{j_3=p+1}^{\infty}
C_{j_3 j_1 j_2 j_3 j_2 j_1}=0,
\end{equation}
\begin{equation}
\label{sixsix11}
\lim\limits_{p\to\infty}
\sum_{j_1=p+1}^{\infty}\sum_{j_2=p+1}^{\infty}
\sum_{j_3=p+1}^{\infty}
C_{j_2 j_3 j_1 j_3 j_2 j_1}=0,
\end{equation}
\begin{equation}
\label{sixsix13}
\lim\limits_{p\to\infty}
\sum_{j_1=p+1}^{\infty}\sum_{j_2=p+1}^{\infty}
\sum_{j_3=p+1}^{\infty}
C_{j_3 j_1 j_3 j_2 j_2 j_1}=0,
\end{equation}
\begin{equation}
\label{sixsix15}
\lim\limits_{p\to\infty}
\sum_{j_1=p+1}^{\infty}\sum_{j_2=p+1}^{\infty}
\sum_{j_3=p+1}^{\infty}
C_{j_2 j_3 j_3 j_1 j_2 j_1}=0.
\end{equation}

\vspace{2mm}

Consider in detail the case of Legendre polynomials 
(the case of trigonometric functions is considered in complete analogy).

First, we prove the following equality for the Fourier coefficients for the case 
$\psi_1(\tau),\ldots,\psi_6(\tau)\equiv 1$
$$
C_{j_6 j_5 j_4 j_3 j_2 j_1}+C_{j_1 j_2 j_3 j_4 j_5 j_6}=
C_{j_6}C_{j_5 j_4 j_3 j_2 j_1}-C_{j_5 j_6}C_{j_4 j_3 j_2 j_1}+
$$
\begin{equation}
\label{sixsix40}
+C_{j_4 j_5 j_6}C_{j_3 j_2 j_1}-C_{j_3 j_4 j_5 j_6}C_{j_2 j_1}+
C_{j_2 j_3 j_4 j_5 j_6}C_{j_1}.
\end{equation}

\vspace{2mm}

Using the integration order replacement, we have
$$
C_{j_6 j_5 j_4 j_3 j_2 j_1}=
$$
$$
=\int\limits_t^T\phi_{j_6}(t_6)\int\limits_t^{t_6}\phi_{j_5}(t_5)
\ldots
\int\limits_t^{t_2}\phi_{j_1}(t_1)dt_1 \ldots dt_5 dt_6=
$$
$$
=\int\limits_t^T\phi_{j_6}(t_6)\int\limits_t^{T}\phi_{j_5}(t_5)
\int\limits_t^{t_5}\phi_{j_4}(t_4)
\ldots \int\limits_t^{t_2}\phi_{j_1}(t_1)dt_1 \ldots  dt_4 dt_5 dt_6-
$$
$$
-\int\limits_t^T\phi_{j_6}(t_6)\int\limits_{t_6}^T\phi_{j_5}(t_5)
\int\limits_t^{t_5}\phi_{j_4}(t_4)
\ldots 
\int\limits_t^{t_2}\phi_{j_1}(t_1)dt_1 \ldots dt_4 dt_5 dt_6=
$$
$$
=C_{j_6}C_{j_5 j_4 j_3 j_2 j_1}-
$$

\vspace{-8mm}
$$
-\int\limits_t^T\phi_{j_6}(t_6)\int\limits_{t_6}^T\phi_{j_5}(t_5)
\int\limits_t^{T}\phi_{j_4}(t_4)
\int\limits_t^{t_4}\phi_{j_3}(t_3)
\ldots 
\int\limits_t^{t_2}\phi_{j_1}(t_1)dt_1 \ldots dt_3 dt_4 dt_5 dt_6+
$$
$$
+\int\limits_t^T\phi_{j_6}(t_6)\int\limits_{t_6}^T\phi_{j_5}(t_5)
\int\limits_{t_5}^{T}\phi_{j_4}(t_4)\int\limits_{t}^{t_4}\phi_{j_3}(t_3)
\ldots 
\int\limits_t^{t_2}\phi_{j_1}(t_1)dt_1 \ldots dt_3 dt_4 dt_5 dt_6=
$$
$$
=C_{j_6}C_{j_5 j_4 j_3 j_2 j_1}-
$$

\vspace{-3mm}
$$
-\int\limits_t^T\phi_{j_6}(t_6)\int\limits_{t_6}^T\phi_{j_5}(t_5)dt_5 dt_6\
C_{j_4 j_3 j_2 j_1}+
$$
$$
+\int\limits_t^T\phi_{j_6}(t_6)\int\limits_{t_6}^T\phi_{j_5}(t_5)
\int\limits_{t_5}^{T}\phi_{j_4}(t_4)\int\limits_{t}^{t_4}\phi_{j_3}(t_3)
\ldots 
\int\limits_t^{t_2}\phi_{j_1}(t_1)dt_1 \ldots dt_3 dt_4 dt_5 dt_6=
$$

\vspace{-2mm}
$$
=C_{j_6}C_{j_5 j_4 j_3 j_2 j_1}-
C_{j_5 j_6}C_{j_4 j_3 j_2 j_1}+
$$

\vspace{-5mm}
$$
+\int\limits_t^T\phi_{j_6}(t_6)\int\limits_{t_6}^T\phi_{j_5}(t_5)
\int\limits_{t_5}^{T}\phi_{j_4}(t_4)\int\limits_{t}^{t_4}\phi_{j_3}(t_3)
\ldots 
\int\limits_t^{t_2}\phi_{j_1}(t_1)dt_1 \ldots dt_3 dt_4 dt_5 dt_6=
$$
$$
\ldots
$$

\vspace{-12mm}
$$
=
C_{j_6}C_{j_5 j_4 j_3 j_2 j_1}-C_{j_5 j_6}C_{j_4 j_3 j_2 j_1}+
C_{j_4 j_5 j_6}C_{j_3 j_2 j_1}-C_{j_3 j_4 j_5 j_6}C_{j_2 j_1}+
C_{j_2 j_3 j_4 j_5 j_6}C_{j_1}-
$$

\vspace{-2mm}
$$
-\int\limits_t^T\phi_{j_6}(t_6)\int\limits_{t_6}^T\phi_{j_5}(t_5)
\ldots
\int\limits_{t_2}^T\phi_{j_1}(t_1)dt_1 \ldots dt_5 dt_6=
$$

\vspace{-2mm}
$$
=
C_{j_6}C_{j_5 j_4 j_3 j_2 j_1}-C_{j_5 j_6}C_{j_4 j_3 j_2 j_1}+
C_{j_4 j_5 j_6}C_{j_3 j_2 j_1}-
$$

\vspace{-5mm}
\begin{equation}
\label{sixsix41}
-C_{j_3 j_4 j_5 j_6}C_{j_2 j_1}+
C_{j_2 j_3 j_4 j_5 j_6}C_{j_1}-C_{j_1 j_2 j_3 j_4 j_5 j_6}.
\end{equation}

\vspace{5mm}

The equality (\ref{sixsix41}) completes the proof of the relation 
(\ref{sixsix40}).

Let us consider 
(\ref{sixsix8}). From (\ref{after80xx}) we obtain
\begin{equation}
\label{sixsix42}
~~~~~~~~~\sum_{j_1=p+1}^{\infty}\sum_{j_2=p+1}^{\infty}
\sum_{j_3=p+1}^{\infty}
C_{j_3 j_2 j_1 j_3 j_2 j_1}=
-\sum_{j_1=0}^{p}\sum_{j_2=0}^{p}
\sum_{j_3=0}^{p}
C_{j_3 j_2 j_1 j_3 j_2 j_1}.
\end{equation}

\vspace{2mm}

Applying (\ref{sixsix40}), we get
$$
\sum_{j_1,j_2,j_3=0}^{p}
C_{j_3 j_2 j_1 j_3 j_2 j_1}+
\sum_{j_1,j_2,j_3=0}^{p}
C_{j_1 j_2 j_3 j_1 j_2 j_3}=
2\sum_{j_1,j_2,j_3=0}^{p}
C_{j_3 j_2 j_1 j_3 j_2 j_1}=
$$
$$
=
\sum_{j_1,j_2,j_3=0}^{p}\biggl(
C_{j_3}C_{j_2 j_1 j_3 j_2 j_1}-C_{j_2 j_3}C_{j_1 j_3 j_2 j_1}+
C_{j_1 j_2 j_3}C_{j_3 j_2 j_1}-\biggr.
$$
\begin{equation}
\label{sixsix43}
\biggl.-C_{j_3 j_1 j_2 j_3}C_{j_2 j_1}+
C_{j_2 j_3 j_1 j_2 j_3}C_{j_1}\biggr).
\end{equation}

\vspace{1mm}

Note that
$$
C_{j_2j_1}=\int\limits_t^T\phi_{j_2}(\tau)\int\limits_t^{\tau}\phi_{j_1}(\theta)d\theta d\tau=
$$

\vspace{-3mm}
\begin{equation}
\label{sixsix50}
=
\frac{T-t}{2}\left\{
\begin{matrix}
1/\sqrt{(2j_1+1)(2j_1+3)} &\hbox{if}\ j_2=j_1+1,\ j_1=0,1,2,\ldots \cr\cr
-1/\sqrt{4j_1^2-1} &\hbox{if}\ j_2=j_1-1,\ j_1=1,2,\ldots \cr\cr
1 &\hbox{if}\ j_1=j_2=0 \cr\cr
0 &\hbox{otherwise}
\end{matrix},\right.
\end{equation}

\vspace{4mm}
\begin{equation}
\label{zero1s}
C_{j_1}=\int\limits_t^T\phi_{j_1}(\tau)d\tau=
\left\{
\begin{matrix}
\sqrt{T-t} &\hbox{if}\ j_1=0\cr\cr
0 &\hbox{if}\ j_1\ne 0
\end{matrix}.\right.
\end{equation}

\vspace{3mm}

Moreover, the generalized Parseval equality gives
$$
\lim\limits_{p\to\infty}\sum_{j_1,j_2,j_3=0}^{p}
C_{j_1 j_2 j_3}C_{j_3 j_2 j_1}=
$$
$$
=\lim\limits_{p\to\infty}\sum_{j_1,j_2,j_3=0}^{p}
\int\limits_t^T\phi_{j_1}(t_3)\int\limits_t^{t_3}\phi_{j_2}(t_2)
\int\limits_t^{t_2}\phi_{j_3}(t_1)dt_1dt_2dt_3\times
$$
$$
\times
\int\limits_t^T\phi_{j_3}(t_3)\int\limits_t^{t_3}\phi_{j_2}(t_2)
\int\limits_t^{t_2}\phi_{j_1}(t_1)dt_1dt_2dt_3=
$$
$$
=\lim\limits_{p\to\infty}\sum_{j_1,j_2,j_3=0}^{p}
\int\limits_t^T\phi_{j_3}(t_3)\int\limits_{t_3}^T\phi_{j_2}(t_2)
\int\limits_{t_2}^T\phi_{j_1}(t_1)dt_1dt_2dt_3\times
$$
$$
\times
\int\limits_t^T\phi_{j_3}(t_3)\int\limits_t^{t_3}\phi_{j_2}(t_2)
\int\limits_t^{t_2}\phi_{j_1}(t_1)dt_1dt_2dt_3=
$$
$$
=\lim\limits_{p\to\infty}\sum_{j_1,j_2,j_3=0}^{p}~
\int\limits_{[t,T]^3}{\bf 1}_{\{t_3<t_2<t_1\}}
\prod\limits_{l=1}^3
\phi_{j_l}(t_l)dt_1dt_2dt_3\times
$$
$$
\times
\int\limits_{[t,T]^3}{\bf 1}_{\{t_1<t_2<t_3\}}
\prod\limits_{l=1}^3
\phi_{j_l}(t_l)dt_1dt_2dt_3=
$$
\begin{equation}
\label{pars100}
=\int\limits_{[t,T]^3}{\bf 1}_{\{t_3<t_2<t_1\}}{\bf 1}_{\{t_1<t_2<t_3\}}
dt_1dt_2dt_3=0.
\end{equation}

Using the above arguments and also (\ref{after80xx}), (\ref{sixsix42}), and (\ref{sixsix43}), we get
$$
-\lim\limits_{p\to\infty}\sum_{j_1=p+1}^{\infty}\sum_{j_2=p+1}^{\infty}
\sum_{j_3=p+1}^{\infty}
C_{j_3 j_2 j_1 j_3 j_2 j_1}=
\lim\limits_{p\to\infty}\sum_{j_1,j_2,j_3=0}^{p}
C_{j_3 j_2 j_1 j_3 j_2 j_1}=
$$
$$
=
\frac{1}{2}\lim\limits_{p\to\infty}\sum_{j_1,j_2,j_3=0}^{p}\biggl(
C_{j_3}C_{j_2 j_1 j_3 j_2 j_1}-C_{j_2 j_3}C_{j_1 j_3 j_2 j_1}-\biggr.
$$
$$
\biggl.-C_{j_3 j_1 j_2 j_3}C_{j_2 j_1}+
C_{j_2 j_3 j_1 j_2 j_3}C_{j_1}\biggr)=
$$
$$
=
\lim\limits_{p\to\infty}\sum_{j_1,j_2,j_3=0}^{p}\biggl(
C_{j_3}C_{j_2 j_1 j_3 j_2 j_1}-C_{j_3 j_1 j_2 j_3}C_{j_2 j_1}\biggr)=
$$
$$
=
\sqrt{T-t}\lim\limits_{p\to\infty}\sum_{j_1,j_2=0}^{p}
C_{j_2 j_1 0 j_2 j_1}
-
\lim\limits_{p\to\infty}\sum_{j_1,j_2,j_3=0}^{p}
C_{j_3 j_1 j_2 j_3}C_{j_2 j_1}
=
$$
\begin{equation}
\label{sixsix71}
~~~~~~~~~~=\sqrt{T-t}\lim\limits_{p\to\infty}\sum_{j_1,j_2=0}^{p}
C_{j_2 j_1 0 j_2 j_1}+
\lim\limits_{p\to\infty}\sum_{j_1,j_2=0}^p\sum_{j_3=p+1}^{\infty}
C_{j_3 j_1 j_2 j_3}C_{j_2 j_1}.
\end{equation}

\vspace{2mm}

By analogy with the proof of (\ref{after2508}) (see the proof of Theorem~2.34)
we obtain
\begin{equation}
\label{sixsix72}
~~~~~~~\lim\limits_{p\to\infty}\sum_{j_1,j_2=0}^{p}
C_{j_2 j_1 0 j_2 j_1}=
\lim\limits_{p\to\infty}\sum_{j_1=p+1}^{\infty}\sum_{j_2=p+1}^{\infty} 
C_{j_2 j_1 0 j_2 j_1}=0,
\end{equation}

\vspace{2mm}
\noindent
where we used the following representation
$$
C_{j_2 j_1 0 j_2 j_1}=
$$
$$
=\frac{1}{\sqrt{T-t}}
\int\limits_t^T 
\phi_{j_2}(t_5)
\int\limits_t^{t_5} 
\phi_{j_1}(t_4)
\int\limits_t^{t_4} 
\int\limits_t^{t_3} 
\phi_{j_2}(t_2)
\int\limits_t^{t_2} 
\phi_{j_1}(t_1)
dt_1 dt_2 dt_3 dt_4 dt_5=
$$
$$
=\frac{1}{\sqrt{T-t}}
\int\limits_t^T 
\phi_{j_2}(t_5)
\int\limits_t^{t_5} 
\phi_{j_1}(t_4)
\int\limits_t^{t_4} 
\phi_{j_2}(t_2)
\int\limits_t^{t_2} 
\phi_{j_1}(t_1)
dt_1 
\int\limits_{t_2}^{t_4} 
dt_3 dt_2 dt_4 dt_5=
$$
$$
=\frac{1}{\sqrt{T-t}}
\int\limits_t^T 
\phi_{j_2}(t_5)
\int\limits_t^{t_5} 
\phi_{j_1}(t_4)(t_4-t)
\int\limits_t^{t_4} 
\phi_{j_2}(t_2)
\int\limits_t^{t_2} 
\phi_{j_1}(t_1)
dt_1 
dt_2 dt_4 dt_5+
$$
$$
+\frac{1}{\sqrt{T-t}}
\int\limits_t^T 
\phi_{j_2}(t_5)
\int\limits_t^{t_5} 
\phi_{j_1}(t_4)
\int\limits_t^{t_4} 
\phi_{j_2}(t_2)(t-t_2)
\int\limits_t^{t_2} 
\phi_{j_1}(t_1)
dt_1 
dt_2 dt_4 dt_5\stackrel{\sf def}{=}
$$
$$
\stackrel{\sf def}{=}
\bar C_{j_2 j_1 j_2 j_1}+\tilde C_{j_2 j_1 j_2 j_1}.
$$

\vspace{2mm}

Further, we have (see (\ref{sixsix50}))
$$
\lim\limits_{p\to\infty}\sum_{j_1,j_2=0}^p\sum_{j_3=p+1}^{\infty}
C_{j_3 j_1 j_2 j_3}C_{j_2 j_1}=
\lim\limits_{p\to\infty}\sum_{j_3=p+1}^{\infty}
\biggl(C_{00}C_{j_3 00 j_3}+\biggr.
$$
\begin{equation}
\label{sixsix70}
~~~~~~~\biggl.+\sum\limits_{j_1=1}^p C_{j_1-1,j_1}C_{j_3j_1,j_1-1,j_3}+
\sum\limits_{j_1=1}^{p-1} C_{j_1+1,j_1}C_{j_3j_1,j_1+1,j_3}+
C_{1,0}C_{j_301j_3}
\biggr).
\end{equation}

Observe that
\begin{equation}
\label{sixsix60}
|C_{j_1-1,j_1}|+|C_{j_1+1,j_1}|\le \frac{K}{j_1}\ \ \ (j_1=1,\ldots,p),
\end{equation}
$$
|C_{j_3 00 j_3}|+|C_{j_3j_1,j_1-1,j_3}|+|C_{j_3j_1,j_1+1,j_3}|+|C_{j_301j_3}|\le
$$
\begin{equation}
\label{sixsix61}
\le \frac{K_1}{j_3^2}\ \ \ (j_3\ge p+1),
\end{equation}

\noindent
where constants $K, K_1$ do not depend on $j_1, j_3.$

The estimate (\ref{sixsix60}) follows from (\ref{sixsix50}).
At the same time, the estimate (\ref{sixsix61}) can be obtained using the following reasoning.
First note that the integration order replacement gives
$$
C_{j_3 j_1 j_2 j_3}=\int\limits_t^T\phi_{j_3}(t_4)\int\limits_{t}^{t_4}\phi_{j_1}(t_3)
\int\limits_t^{t_3}\phi_{j_2}(t_2)
\int\limits_t^{t_2}\phi_{j_3}(t_1)
dt_1 dt_2 dt_3 dt_4=
$$
\begin{equation}
\label{sixsix62}
~~~~~~~~=\int\limits_t^T\phi_{j_1}(t_3)\int\limits_{t}^{t_3}\phi_{j_2}(t_2)
\left(\int\limits_t^{t_2}\phi_{j_3}(t_1)dt_1\right) dt_2
\left(\int\limits_{t_3}^T\phi_{j_3}(t_4)
dt_4\right) dt_3.
\end{equation}

\vspace{2mm}

Applying the estimates (\ref{ogo24}), (\ref{ogo25}), and (\ref{ogoo25}) to (\ref{sixsix62}) 
gives the estimate (\ref{sixsix61}).

Using (\ref{sixsix70}), (\ref{sixsix60}), and (\ref{sixsix61}),
we obtain
$$
\left\vert
\sum_{j_1,j_2=0}^p\sum_{j_3=p+1}^{\infty}
C_{j_3 j_1 j_2 j_3}C_{j_2 j_1}\right\vert\le
K\sum\limits_{j_3=p+1}^{\infty}\frac{1}{j_3^2}
\left(1+\sum\limits_{j_1=1}^p \frac{1}{j_1}\right)\le
$$
\begin{equation}
\label{sixsix74}
\le K\int\limits_p^{\infty} \frac{dx}{x^2}
\left(2+\int\limits_1^p \frac{dx}{x}\right)=
\frac{K(2+ln p)}{p}\to 0
\end{equation}

\noindent
if $p\to\infty$, where constant $K$ is independent of $p.$
Thus, the equality (\ref{sixsix8}) is proved (see (\ref{sixsix71}), (\ref{sixsix72}),
(\ref{sixsix74})).

The relation (\ref{sixsix9}) is proved in complete analogy 
with the proof of equality (\ref{sixsix8}).
For (\ref{sixsix9}) we have (see (\ref{sixsix40}))
$$
\lim\limits_{p\to\infty}
\left(\sum_{j_1,j_2,j_3=0}^{p}
C_{j_1 j_3 j_2 j_3 j_2 j_1}+
\sum_{j_1,j_2,j_3=0}^{p}
C_{j_1 j_2 j_3 j_2 j_3 j_1}\right)=
2\lim\limits_{p\to\infty}\sum_{j_1,j_2,j_3=0}^{p}
C_{j_1 j_3 j_2 j_3 j_2 j_1}=
$$
$$
=
\lim\limits_{p\to\infty}\sum_{j_1,j_2,j_3=0}^{p}\biggl(
C_{j_1}C_{j_3 j_2 j_3 j_2 j_1}-C_{j_3 j_1}C_{j_2 j_3 j_2 j_1}+
C_{j_2 j_3 j_1}C_{j_3 j_2 j_1}-\biggr.
$$
$$
\biggl.-C_{j_3 j_2 j_3 j_1}C_{j_2 j_1}+
C_{j_2 j_3 j_2 j_3 j_1}C_{j_1}\biggr)=
$$
$$
=2\lim\limits_{p\to\infty}\left(
\sqrt{T-t}\sum_{j_2,j_3=0}^{p}C_{j_3 j_2 j_3 j_2 0}-
\sum_{j_1,j_2,j_3=0}^{p}
C_{j_2 j_1}C_{j_3 j_2 j_3 j_1}\right)=
$$
$$
=-2\lim\limits_{p\to\infty}
\sum_{j_1,j_2,j_3=0}^{p}
C_{j_2 j_1}C_{j_3 j_2 j_3 j_1}.
$$

\vspace{2mm}

To estimate the Fourier coefficient $C_{j_3 j_2 j_3 j_1}$, 
we use the following 
(see the proof of (\ref{sixsix8}) for more details)
$$
C_{j_3 j_2 j_3 j_1}=\int\limits_t^T\phi_{j_3}(t_4)\int\limits_{t}^{t_4}\phi_{j_2}(t_3)
\int\limits_t^{t_3}\phi_{j_3}(t_2)
\int\limits_t^{t_2}\phi_{j_1}(t_1)
dt_1 dt_2 dt_3 dt_4=
$$
$$
=\int\limits_t^T\phi_{j_3}(t_4)\int\limits_{t}^{t_4}\phi_{j_2}(t_3)
\int\limits_t^{t_3}\phi_{j_1}(t_1)
\int\limits_{t_1}^{t_3}\phi_{j_3}(t_2)
dt_2 dt_1 dt_3 dt_4=
$$
$$
=\int\limits_t^T\phi_{j_3}(t_4)\int\limits_{t}^{t_4}\phi_{j_2}(t_3)
\left(\int\limits_{t}^{t_3}\phi_{j_3}(t_2)
dt_2\right)
\int\limits_t^{t_3}\phi_{j_1}(t_1)
dt_1 dt_3 dt_4-
$$
$$
-\int\limits_t^T\phi_{j_3}(t_4)\int\limits_{t}^{t_4}\phi_{j_2}(t_3)
\int\limits_t^{t_3}\phi_{j_1}(t_1)
\left(\int\limits_{t}^{t_1}\phi_{j_3}(t_2)
dt_2\right) dt_1 dt_3 dt_4=
$$
$$
=\int\limits_t^T
\phi_{j_2}(t_3)
\left(\int\limits_{t}^{t_3}\phi_{j_3}(t_2)
dt_2\right)
\int\limits_t^{t_3}\phi_{j_1}(t_1)
dt_1
\left(\int\limits_{t_3}^{T}
\phi_{j_3}(t_4)
dt_4\right) dt_3-
$$
$$
-\int\limits_t^T
\phi_{j_2}(t_3)
\int\limits_t^{t_3}\phi_{j_1}(t_1)
\left(\int\limits_{t}^{t_1}\phi_{j_3}(t_2)
dt_2\right) dt_1
\left(\int\limits_{t_3}^{T}
\phi_{j_3}(t_4) dt_4\right) dt_3.
$$

\vspace{2mm}

Let us prove (\ref{sixsix10}). 
From (\ref{after80xx}) we obtain
\begin{equation}
\label{sixsix80}
~~~~~~~~~\sum_{j_1=p+1}^{\infty}\sum_{j_2=p+1}^{\infty}
\sum_{j_3=p+1}^{\infty}
C_{j_3 j_2 j_3 j_1 j_2 j_1}=
-\sum_{j_1=0}^{p}\sum_{j_2=0}^{p}
\sum_{j_3=0}^{p}
C_{j_3 j_2 j_3 j_1 j_2 j_1}.
\end{equation}

\vspace{2mm}

Applying (\ref{sixsix40}) and (\ref{sixsix80}), we get (we replaced $j_3$ by $j_4$)
$$
\sum_{j_1,j_2,j_4=0}^{p}
C_{j_4 j_2 j_4 j_1 j_2 j_1}+
\sum_{j_1,j_2,j_4=0}^{p}
C_{j_1 j_2 j_1 j_4 j_2 j_4}=
2\sum_{j_1,j_2,j_4=0}^{p}
C_{j_4 j_2 j_4 j_1 j_2 j_1}=
$$
$$
=
\sum_{j_1,j_2,j_4=0}^{p}\biggl(
C_{j_4}C_{j_2 j_4 j_1 j_2 j_1}-C_{j_2 j_4}C_{j_4 j_1 j_2 j_1}+
C_{j_4 j_2 j_4}C_{j_1 j_2 j_1}-\biggr.
$$
$$
\biggl.-C_{j_1 j_4 j_2 j_4}C_{j_2 j_1}+
C_{j_2 j_1 j_4 j_2 j_4}C_{j_1}\biggr)=
$$
$$
=
2\sum_{j_1,j_2,j_4=0}^{p}\biggl(C_{j_2 j_1 j_4 j_2 j_4}C_{j_1}-C_{j_1 j_4 j_2 j_4}C_{j_2 j_1}
\biggr)+
$$
\begin{equation}
\label{sixsix83}
+
\sum_{j_1,j_2,j_4=0}^{p}
C_{j_4 j_2 j_4}C_{j_1 j_2 j_1}.
\end{equation}

\vspace{2mm}

Further, we have (see (\ref{after80xx}))
$$
\lim\limits_{p\to\infty}\sum_{j_1,j_2,j_4=0}^{p}
C_{j_4 j_2 j_4}C_{j_1 j_2 j_1}=
\lim\limits_{p\to\infty}\sum_{j_2=0}^{p}
\left(\sum_{j_1=0}^{p}C_{j_1 j_2 j_1}\right)^2=
$$
\begin{equation}
\label{sixsix84}
=\lim\limits_{p\to\infty}\sum_{j_2=0}^{p}
\left(\sum_{j_1=p+1}^{\infty}C_{j_1 j_2 j_1}\right)^2=0,
\end{equation}

\vspace{2mm}
\noindent
where we applied the equality (\ref{after1602}).

Furthermore, by analogy with the proof of (\ref{sixsix8}), we have
\begin{equation}
\label{sixsix85}
\lim\limits_{p\to\infty}
\sum_{j_1,j_2,j_4=0}^{p}\biggl(C_{j_2 j_1 j_4 j_2 j_4}C_{j_1}-C_{j_1 j_4 j_2 j_4}C_{j_2 j_1}
\biggr)=0.
\end{equation}

To estimate the Fourier coefficient $C_{j_1 j_4 j_2 j_4}$ in (\ref{sixsix85}),
we use the following 
(see the proof of (\ref{sixsix8}) for more details)
$$
C_{j_1 j_4 j_2 j_4}=\int\limits_t^T\phi_{j_1}(t_4)\int\limits_{t}^{t_4}\phi_{j_4}(t_3)
\int\limits_t^{t_3}\phi_{j_2}(t_2)
\left(\int\limits_t^{t_2}\phi_{j_4}(t_1)
dt_1\right) dt_2 dt_3 dt_4=
$$
$$
=\int\limits_t^T\phi_{j_1}(t_4)\int\limits_{t}^{t_4}
\phi_{j_2}(t_2)
\left(\int\limits_t^{t_2}\phi_{j_4}(t_1)
dt_1\right)
\int\limits_{t_2}^{t_4}
\phi_{j_4}(t_3)
dt_3 dt_2 dt_4=
$$
$$
=\int\limits_t^T\phi_{j_1}(t_4)
\left(\int\limits_{t}^{t_4}
\phi_{j_4}(t_3)
dt_3\right)
\int\limits_{t}^{t_4}
\phi_{j_2}(t_2)
\left(\int\limits_t^{t_2}\phi_{j_4}(t_1)
dt_1\right)
dt_2 dt_4-
$$
$$
-\int\limits_t^T\phi_{j_1}(t_4)\int\limits_{t}^{t_4}
\phi_{j_2}(t_2)\left(
\int\limits_{t}^{t_2}
\phi_{j_4}(t_3)
dt_3\right)
\left(\int\limits_t^{t_2}\phi_{j_4}(t_1)
dt_1\right)
dt_2 dt_4.
$$

\vspace{2mm}

The relations (\ref{sixsix80})--(\ref{sixsix85}) 
complete the proof of equality (\ref{sixsix10}).

Let us prove (\ref{sixsix4}).
Using (\ref{after80xx}), we get 
\begin{equation}
\label{sixsix90}
~~~~~~~~~\sum_{j_1=p+1}^{\infty}\sum_{j_2=p+1}^{\infty}
\sum_{j_3=p+1}^{\infty}
C_{j_1 j_2 j_3 j_3 j_2 j_1}=
\sum_{j_1=0}^{p}\sum_{j_2=0}^{p}
\sum_{j_3=p+1}^{\infty}
C_{j_1 j_2 j_3 j_3 j_2 j_1}.
\end{equation}

\vspace{2mm}

Applying (\ref{sixsix40}) and (\ref{sixsix90}), we obtain
$$
2\sum_{j_1,j_2=0}^{p}\sum_{j_3=p+1}^{\infty}
C_{j_1 j_2 j_3 j_3 j_2 j_1}=
$$
$$
=
\sum_{j_1,j_2=0}^{p}\sum_{j_3=p+1}^{\infty}\biggl(
C_{j_1}C_{j_2 j_3 j_3 j_2 j_1}-C_{j_2 j_1}C_{j_3 j_3 j_2 j_1}+
\left(C_{j_3 j_2 j_1}\right)^2-\biggr.
$$
$$
\biggl.-C_{j_3 j_3 j_2 j_1}C_{j_2 j_1}+
C_{j_2 j_3 j_3 j_2 j_1}C_{j_1}\biggr)=
$$
$$
=
2\sum_{j_1,j_2=0}^{p}\sum_{j_3=p+1}^{\infty}
\biggl(C_{j_1}C_{j_2 j_3 j_3 j_2 j_1}-C_{j_2 j_1}C_{j_3 j_3 j_2 j_1}
\biggr)+
$$
\begin{equation}
\label{sixsix91}
+
\sum_{j_1,j_2=0}^{p}\sum_{j_3=p+1}^{\infty}
\left(C_{j_3 j_2 j_1}\right)^2.
\end{equation}

\vspace{2mm}

Using the estimate (\ref{fffoh}),  we get
\begin{equation}
\label{sixsix92}
\lim\limits_{p\to\infty}\sum_{j_1,j_2=0}^{p}\sum_{j_3=p+1}^{\infty}
\left(C_{j_3 j_2 j_1}\right)^2=0.
\end{equation}

\vspace{2mm}

By analogy with the proof of (\ref{sixsix8}), we have
\begin{equation}
\label{sixsix93}
\lim\limits_{p\to\infty}
\sum_{j_1,j_2=0}^{p}\sum_{j_3=p+1}^{\infty}
\biggl(C_{j_1}C_{j_2 j_3 j_3 j_2 j_1}-C_{j_2 j_1}C_{j_3 j_3 j_2 j_1}
\biggr)=0,
\end{equation}

\noindent
where we applied
the equality (\ref{after2509}).
To estimate the Fourier coefficient $C_{j_3 j_3 j_2 j_1}$ in (\ref{sixsix93}),
we used the following 
(see the proof of (\ref{sixsix8}) for more details)
$$
C_{j_3 j_3 j_2 j_1}=\int\limits_t^T\phi_{j_3}(t_4)\int\limits_{t}^{t_4}\phi_{j_3}(t_3)
\int\limits_t^{t_3}\phi_{j_2}(t_2)
\int\limits_t^{t_2}\phi_{j_1}(t_1)
dt_1dt_2 dt_3 dt_4=
$$
$$
=\int\limits_t^T\phi_{j_1}(t_1)\int\limits_{t_1}^{T}
\phi_{j_2}(t_2)
\int\limits_{t_2}^{T}
\phi_{j_3}(t_3)
\int\limits_{t_3}^{T}\phi_{j_3}(t_4)
dt_4dt_3 dt_2 dt_1=
$$
\begin{equation}
\label{sept10}
=\frac{1}{2}\int\limits_t^T\phi_{j_1}(t_1)\int\limits_{t_1}^{T}
\phi_{j_2}(t_2)
\left(\int\limits_{t_2}^{T}
\phi_{j_3}(t_3)
dt_3\right)^2 dt_2 dt_1.
\end{equation}

\vspace{2mm}

Combining the equalities (\ref{sixsix90})--(\ref{sixsix93}), 
we obtain (\ref{sixsix4}).

Let us prove (\ref{sixsix14}) (we replace $j_2$ by $j_4$ and $j_3$ by $j_2$
in (\ref{sixsix14})).
As noted in Remark~2.4, the sequential order of the series
$$
\sum_{j_1=p+1}^{\infty}\sum_{j_2=p+1}^{\infty}
\sum_{j_4=p+1}^{\infty}
$$

\noindent
is not important. This follows directly from the formulas 
(\ref{after500}) and (\ref{after80xx}).

Applying the mentioned property and (\ref{after80xx}), we get
\begin{equation}
\label{sixsix94}
~~~~~~~~~~~~\sum_{j_1=p+1}^{\infty}\sum_{j_2=p+1}^{\infty}
\sum_{j_4=p+1}^{\infty}
C_{j_1 j_4 j_4 j_2 j_2 j_1}=
-\sum_{j_1=0}^{p}\sum_{j_2=p+1}^{\infty}
\sum_{j_4=p+1}^{\infty}
C_{j_1 j_4 j_4 j_2 j_2 j_1}.
\end{equation}

\vspace{2mm}

Observe that (see the above reasoning)
\begin{equation}
\label{sixsix95}
\sum_{j_2=p+1}^{\infty}
\sum_{j_4=p+1}^{\infty}
C_{j_1 j_4 j_4 j_2 j_2 j_1}=
\sum_{j_4=p+1}^{\infty}
\sum_{j_2=p+1}^{\infty}
C_{j_1 j_4 j_4 j_2 j_2 j_1}.
\end{equation}

\vspace{2mm}

Using (\ref{sixsix40}) and (\ref{sixsix95}), we obtain
$$
\sum_{j_1=0}^{p}\sum_{j_2=p+1}^{\infty}
\sum_{j_4=p+1}^{\infty}
\biggl(C_{j_1 j_4 j_4 j_2 j_2 j_1}+
C_{j_1 j_2 j_2 j_4 j_4 j_1}\biggr)=
2\sum_{j_1=0}^{p}\sum_{j_2=p+1}^{\infty}
\sum_{j_4=p+1}^{\infty}
C_{j_1 j_4 j_4 j_2 j_2 j_1}=
$$
$$
=
\sum_{j_1=0}^{p}\sum_{j_2=p+1}^{\infty}
\sum_{j_4=p+1}^{\infty}
\biggl(
C_{j_1}C_{j_4 j_4 j_2 j_2 j_1}-C_{j_4 j_1}C_{j_4 j_2 j_2 j_1}+
C_{j_4 j_4 j_1}C_{j_2 j_2 j_1}-\biggr.
$$
$$
\biggl.-C_{j_2 j_4 j_4 j_1}C_{j_2 j_1}+
C_{j_2 j_2 j_4 j_4 j_1}C_{j_1}\biggr)=
$$
$$
=
\sum_{j_1=0}^{p}\sum_{j_2=p+1}^{\infty}
\sum_{j_4=p+1}^{\infty}
\biggl(
C_{j_1}C_{j_4 j_4 j_2 j_2 j_1}-C_{j_4 j_1}C_{j_4 j_2 j_2 j_1}
-C_{j_2 j_4 j_4 j_1}C_{j_2 j_1}+
C_{j_2 j_2 j_4 j_4 j_1}C_{j_1}\biggr)+
$$
\begin{equation}
\label{sixsix96}
+\sum_{j_1=0}^{p}\left(\sum_{j_2=p+1}^{\infty}
C_{j_2 j_2 j_1}\right)^2.
\end{equation}

\vspace{2mm}

The equality 
\begin{equation}
\label{sixsix97}
\lim\limits_{p\to\infty}\sum_{j_1=0}^{p}\left(\sum_{j_2=p+1}^{\infty}
C_{j_2 j_2 j_1}\right)^2=0
\end{equation}

\vspace{3mm}
\noindent
follows from 
the relation (\ref{after1601}).

By analogy with the proof of equality (\ref{sixsix8}) we obtain
$$
\lim\limits_{p\to\infty}
\sum_{j_1=0}^{p}\sum_{j_2=p+1}^{\infty}
\sum_{j_4=p+1}^{\infty}
\biggl(
C_{j_1}C_{j_4 j_4 j_2 j_2 j_1}-C_{j_4 j_1}C_{j_4 j_2 j_2 j_1}
-\biggr.
$$
\begin{equation}
\label{sixsix98}
\biggl.-
C_{j_2 j_4 j_4 j_1}C_{j_2 j_1}+
C_{j_2 j_2 j_4 j_4 j_1}C_{j_1}\biggr)=0,
\end{equation}

\noindent
where we applied
the equality (\ref{after2507}).
To estimate the Fourier coefficient $C_{j_2 j_4 j_4 j_1}$ in (\ref{sixsix98}),
we used the following 
(see the proof of (\ref{sixsix8}) for more details)
$$
C_{j_2 j_4 j_4 j_1}=\int\limits_t^T\phi_{j_2}(t_4)\int\limits_{t}^{t_4}\phi_{j_4}(t_3)
\int\limits_t^{t_3}\phi_{j_4}(t_2)
\int\limits_t^{t_2}\phi_{j_1}(t_1)
dt_1dt_2 dt_3 dt_4=
$$
$$
=\int\limits_t^T\phi_{j_2}(t_4)\int\limits_{t}^{t_4}
\phi_{j_1}(t_1)
\int\limits_{t_1}^{t_4}
\phi_{j_4}(t_2)
\int\limits_{t_2}^{t_4}\phi_{j_4}(t_3)
dt_3dt_2 dt_1 dt_4=
$$
$$
=\frac{1}{2}\int\limits_t^T\phi_{j_2}(t_4)\int\limits_{t}^{t_4}
\phi_{j_1}(t_1)
\left(\int\limits_{t_1}^{t_4}
\phi_{j_4}(t_2)dt_2\right)^2 dt_1 dt_4=
$$
$$
=\frac{1}{2}\int\limits_t^T\phi_{j_2}(t_4)
\left(\int\limits_{t}^{t_4}
\phi_{j_4}(t_2)dt_2\right)^2\int\limits_{t}^{t_4}
\phi_{j_1}(t_1)
dt_1 dt_4+
$$
$$
+\frac{1}{2}\int\limits_t^T\phi_{j_2}(t_4)\int\limits_{t}^{t_4}
\phi_{j_1}(t_1)
\left(\int\limits_{t}^{t_1}
\phi_{j_4}(t_2)dt_2\right)^2 dt_1 dt_4-
$$
$$
-\int\limits_t^T\phi_{j_2}(t_4)
\left(\int\limits_{t}^{t_4}
\phi_{j_4}(t_2)dt_2\right)\int\limits_{t}^{t_4}
\phi_{j_1}(t_1)
\left(\int\limits_{t}^{t_1}
\phi_{j_4}(t_2)dt_2\right)
dt_1 dt_4.
$$

\vspace{2mm}

The relation (\ref{sixsix14}) follows from (\ref{sixsix94}),
(\ref{sixsix96})--(\ref{sixsix98}).

Consider (\ref{sixsix3}). Using the integration 
order replacement, we obtain
$$
C_{j_3j_3j_2j_2j_1j_1}=
$$
$$
=
\frac{1}{2}\int\limits_t^T \phi_{j_3}(t_6)
\int\limits_t^{t_6} \phi_{j_3}(t_5)
\int\limits_t^{t_5} \phi_{j_2}(t_4)
\int\limits_t^{t_4} \phi_{j_2}(t_3)
\left(\int\limits_t^{t_3} \phi_{j_1}(t_1)dt_1\right)^2 dt_3 dt_4 dt_5 dt_6=
$$
$$
=
\frac{1}{2}\int\limits_t^T 
\phi_{j_2}(t_3)
\left(\int\limits_t^{t_3} \phi_{j_1}(t_1)dt_1\right)^2
\int\limits_{t_3}^T 
\phi_{j_2}(t_4) 
\int\limits_{t_4}^T 
\phi_{j_3}(t_5) 
\int\limits_{t_5}^T 
\phi_{j_3}(t_6) 
dt_6 dt_5 dt_4 dt_3=
$$
\begin{equation}
\label{sixsix100}
~~~~~=
\frac{1}{4}\int\limits_t^T 
\phi_{j_2}(t_3)
\left(\int\limits_t^{t_3} \phi_{j_1}(t_1)dt_1\right)^2
\int\limits_{t_3}^T 
\phi_{j_2}(t_4) 
\left(\int\limits_{t_4}^T 
\phi_{j_3}(t_5) 
dt_5\right)^2 dt_4 dt_3.
\end{equation}

\vspace{3mm}

Applying the estimates (\ref{ogo24}), (\ref{ogo25}), and (\ref{ogoo25})  to 
(\ref{sixsix100}) 
gives the following estimate 
\begin{equation}
\label{sixsix101}
|C_{j_3j_3j_2j_2j_1j_1}|\le \frac{K}{j_1^2 j_3^2}\ \ \ (j_1, j_3>0,\ j_2\ge 0),
\end{equation}

\noindent
where constant $K$ does not depend on $j_1, j_2, j_3$.

Further, we get (see (\ref{after500}))
$$
\sum_{j_1=p+1}^{\infty}\sum_{j_2=p+1}^{\infty}
\sum_{j_3=p+1}^{\infty}
C_{j_3 j_3 j_2 j_2 j_1 j_1}=
\sum_{j_1=p+1}^{\infty}\sum_{j_3=p+1}^{\infty}
\sum_{j_2=p+1}^{\infty}
C_{j_3 j_3 j_2 j_2 j_1 j_1}=
$$
\begin{equation}
\label{sixsix102}
~~~~~=\frac{1}{2}
\sum_{j_1=p+1}^{\infty}\sum_{j_3=p+1}^{\infty}
C_{j_3j_3j_2j_2j_1j_1}\biggl|_{(j_2 j_2)\curvearrowright (\cdot)}\biggr. -
\sum_{j_2=0}^{p}\sum_{j_1=p+1}^{\infty}\sum_{j_3=p+1}^{\infty}
C_{j_3 j_3 j_2 j_2 j_1 j_1},
\end{equation}

\vspace{2mm}
\noindent
where
$$
C_{j_3j_3j_2j_2j_1j_1}\biggl|_{(j_2 j_2)\curvearrowright (\cdot)}\biggr.=
$$
$$
=
\int\limits_t^T \phi_{j_3}(t_6)
\int\limits_t^{t_6} \phi_{j_3}(t_5)
\int\limits_t^{t_5} 
\int\limits_t^{t_4} 
\phi_{j_1}(t_2)
\int\limits_t^{t_2} 
\phi_{j_1}(t_1)
dt_1 dt_2 dt_4 dt_5 dt_6=
$$
$$
=
\int\limits_t^T \phi_{j_3}(t_6)
\int\limits_t^{t_6} \phi_{j_3}(t_5)
\int\limits_t^{t_5} 
\phi_{j_1}(t_2)
\int\limits_t^{t_2} 
\phi_{j_1}(t_1)
dt_1
\int\limits_{t_2}^{t_5} 
dt_4 dt_2 dt_5 dt_6=
$$
$$
=
\int\limits_t^T \phi_{j_3}(t_6)
\int\limits_t^{t_6} \phi_{j_3}(t_5)(t_5-t)
\int\limits_t^{t_5} 
\phi_{j_1}(t_2)
\int\limits_t^{t_2} 
\phi_{j_1}(t_1)
dt_1dt_2 dt_5 dt_6+
$$
$$
+\int\limits_t^T \phi_{j_3}(t_6)
\int\limits_t^{t_6} \phi_{j_3}(t_5)
\int\limits_t^{t_5} 
\phi_{j_1}(t_2)(t-t_2)
\int\limits_t^{t_2} 
\phi_{j_1}(t_1)
dt_1dt_2 dt_5 dt_6\stackrel{\sf def}{=}
$$
\begin{equation}
\label{sixsix103}
\stackrel{\sf def}{=}
C'_{j_3 j_3 j_1 j_1}+C''_{j_3 j_3 j_1 j_1}.
\end{equation}

\vspace{2mm}

Let us substitute (\ref{sixsix103}) into (\ref{sixsix102})
$$
\sum_{j_1=p+1}^{\infty}\sum_{j_2=p+1}^{\infty}
\sum_{j_3=p+1}^{\infty}
C_{j_3 j_3 j_2 j_2 j_1 j_1}=
\frac{1}{2}
\sum_{j_1=p+1}^{\infty}\sum_{j_3=p+1}^{\infty}
C'_{j_3 j_3 j_1 j_1}+
$$
\begin{equation}
\label{sixsix104}
~~~~~~~+
\frac{1}{2}
\sum_{j_1=p+1}^{\infty}\sum_{j_3=p+1}^{\infty}
C''_{j_3 j_3 j_1 j_1}
-
\sum_{j_2=0}^{p}\sum_{j_1=p+1}^{\infty}\sum_{j_3=p+1}^{\infty}
C_{j_3 j_3 j_2 j_2 j_1 j_1}.
\end{equation}

\vspace{2mm}

The relation (\ref{after2507})
implies that
\begin{equation}
\label{sixsix105}
~~~~~~~\lim\limits_{p\to\infty}\sum_{j_1=p+1}^{\infty}\sum_{j_3=p+1}^{\infty}
C'_{j_3 j_3 j_1 j_1}=0,\ \ \ 
\lim\limits_{p\to\infty}\sum_{j_1=p+1}^{\infty}\sum_{j_3=p+1}^{\infty}
C''_{j_3 j_3 j_1 j_1}=0.
\end{equation}

From the estimate (\ref{sixsix101}) we get
$$
\left\vert\sum_{j_2=0}^{p}\sum_{j_1=p+1}^{\infty}\sum_{j_3=p+1}^{\infty}
C_{j_3 j_3 j_2 j_2 j_1 j_1}\right\vert\le
K(p+1) \sum_{j_1=p+1}^{\infty}\frac{1}{j_1^2}
\sum_{j_3=p+1}^{\infty}\frac{1}{j_3^2}\le 
$$
\begin{equation}
\label{sixsix106}
\le
K(p+1) \left(\int\limits_p^{\infty} \frac{dx}{x^2}\right)^2
\le \frac{K(p+1)}{p^2}\ \to\ 0
\end{equation}

\noindent
if $p\to\infty$, where constant $K$ is independent of $p$.

The relations (\ref{sixsix104})--(\ref{sixsix106})
complete the proof of (\ref{sixsix3}).

Let us prove (\ref{sixsix7}). Using the integration
order replacement, we get
$$
C_{j_2 j_3 j_3 j_2 j_1 j_1}=
$$
$$
=
\frac{1}{2}\int\limits_t^T 
\phi_{j_2}(t_6)
\int\limits_t^{t_6} 
\phi_{j_3}(t_5)
\int\limits_{t}^{t_5} 
\phi_{j_3}(t_4) 
\int\limits_t^{t_4}
\phi_{j_2}(t_3) 
\left(\int\limits_t^{t_3}
\phi_{j_1}(t_1)dt_1\right)^2 
dt_3 dt_4 dt_5 dt_6=
$$
$$
=
\frac{1}{2}\int\limits_t^T 
\phi_{j_2}(t_3) 
\left(\int\limits_t^{t_3}
\phi_{j_1}(t_1)dt_1\right)^2 
\int\limits_{t_3}^T
\phi_{j_3}(t_4) 
\int\limits_{t_4}^T
\phi_{j_3}(t_5)
\int\limits_{t_5}^T
\phi_{j_2}(t_6)dt_6
dt_5 dt_4 dt_3=
$$
$$
=
\frac{1}{2}\int\limits_t^T 
\phi_{j_2}(t_3) 
\left(\int\limits_t^{t_3}
\phi_{j_1}(t_1)dt_1\right)^2 
\int\limits_{t_3}^T
\phi_{j_3}(t_5)
\int\limits_{t_5}^T
\phi_{j_2}(t_6)dt_6
\int\limits_{t_3}^{t_5}
\phi_{j_3}(t_4) 
dt_4 dt_5 dt_3=
$$
$$
=
\frac{1}{2}\int\limits_t^T 
\phi_{j_2}(t_3) 
\left(\int\limits_t^{t_3}
\phi_{j_1}(t_1)dt_1\right)^2 
\int\limits_{t_3}^T
\phi_{j_3}(t_5)
\left(\int\limits_{t_5}^T
\phi_{j_2}(t_6)dt_6\right)
\left(\int\limits_{t}^{t_5}
\phi_{j_3}(t_4)dt_4\right)\times
$$
$$
\times 
dt_5 dt_3-
$$
$$
-
\frac{1}{2}\int\limits_t^T 
\phi_{j_2}(t_3) 
\left(\int\limits_t^{t_3}
\phi_{j_1}(t_1)dt_1\right)^2 
\left(\int\limits_{t}^{t_3}
\phi_{j_3}(t_4)dt_4\right)\int\limits_{t_3}^T
\phi_{j_3}(t_5)
\left(\int\limits_{t_5}^T
\phi_{j_2}(t_6)dt_6\right)
\times
$$
\begin{equation}
\label{sixsix107}
\times 
dt_5 dt_3.
\end{equation}

\vspace{2mm}

Applying (\ref{after80xx}) and (\ref{after500}), we obtain
$$
-\sum_{j_1=p+1}^{\infty}\sum_{j_2=p+1}^{\infty}
\sum_{j_3=p+1}^{\infty}
C_{j_2 j_3 j_3 j_2 j_1 j_1}=
-\sum_{j_1=p+1}^{\infty}\sum_{j_3=p+1}^{\infty}\sum_{j_2=p+1}^{\infty}
C_{j_2 j_3 j_3 j_2 j_1 j_1}=
$$
$$
=
\sum_{j_2=0}^{p}\sum_{j_1=p+1}^{\infty}\sum_{j_3=p+1}^{\infty}
C_{j_2 j_3 j_3 j_2 j_1 j_1}=
$$
$$
=\frac{1}{2}
\sum_{j_2=0}^{p}\sum_{j_1=p+1}^{\infty}
C_{j_2 j_3 j_3 j_2 j_1 j_1}\biggl|_{(j_3 j_3)\curvearrowright (\cdot)}\biggr. -
\sum_{j_2=0}^{p}\sum_{j_3=0}^{p}\sum_{j_1=p+1}^{\infty}
C_{j_2 j_3 j_3 j_2 j_1 j_1}=
$$
$$
=\frac{1}{2}
\sum_{j_2=0}^{p}\sum_{j_1=p+1}^{\infty}
C_{j_2 j_3 j_3 j_2 j_1 j_1}\biggl|_{(j_3 j_3)\curvearrowright (\cdot)}\biggr. -
\sum_{j_1=p+1}^{\infty}
C_{0000 j_1 j_1}-
$$
$$
-\sum_{j_3=1}^{p}\sum_{j_1=p+1}^{\infty}
C_{0 j_3 j_3 0 j_1 j_1}-
\sum_{j_2=1}^{p}\sum_{j_1=p+1}^{\infty}
C_{j_2 00 j_2 j_1 j_1}-
$$
\begin{equation}
\label{sixsix108}
-
\sum_{j_2=1}^{p}\sum_{j_3=1}^{p}\sum_{j_1=p+1}^{\infty}
C_{j_2 j_3 j_3 j_2 j_1 j_1}.
\end{equation}

\vspace{1mm}

The equality
\begin{equation}
\label{sixsix109a}
\lim\limits_{p\to\infty}
\frac{1}{2}
\sum_{j_2=0}^{p}\sum_{j_1=p+1}^{\infty}
C_{j_2 j_3 j_3 j_2 j_1 j_1}\biggl|_{(j_3 j_3)\curvearrowright (\cdot)}\biggr.=0
\end{equation}

\vspace{1mm}
\noindent
follows from the inequality similar to (\ref{after9043})
(see the proof of Theorem~2.34),
where we used the following representation
$$
C_{j_2j_3j_3j_2j_1j_1}\biggl|_{(j_3 j_3)\curvearrowright (\cdot)}\biggr.=
$$
$$
=
\int\limits_t^T \phi_{j_2}(t_6)
\int\limits_t^{t_6} 
\int\limits_t^{t_4} 
\phi_{j_2}(t_3)
\int\limits_t^{t_3} 
\phi_{j_1}(t_2)
\int\limits_t^{t_2} 
\phi_{j_1}(t_1)
dt_1 dt_2 dt_3 dt_4 dt_6=
$$
$$
=
\int\limits_t^T \phi_{j_2}(t_6)
\int\limits_t^{t_6} \phi_{j_2}(t_3)
\int\limits_t^{t_3} 
\phi_{j_1}(t_2)
\int\limits_t^{t_2} 
\phi_{j_1}(t_1)
dt_1 dt_2
\int\limits_{t_3}^{t_6} 
dt_4 dt_3 dt_6=
$$
$$
+\int\limits_t^T \phi_{j_2}(t_6)(t_6-t)
\int\limits_t^{t_6} \phi_{j_2}(t_3)
\int\limits_t^{t_3} 
\phi_{j_1}(t_2)
\int\limits_t^{t_2} 
\phi_{j_1}(t_1)
dt_1dt_2 dt_3 dt_6+
$$
$$
+\int\limits_t^T \phi_{j_2}(t_6)
\int\limits_t^{t_6} \phi_{j_2}(t_3)(t-t_3)
\int\limits_t^{t_3} 
\phi_{j_1}(t_2)
\int\limits_t^{t_2} 
\phi_{j_1}(t_1)
dt_1dt_2 dt_3 dt_6\stackrel{\sf def}{=}
$$
\begin{equation}
\label{sept12}
\stackrel{\sf def}{=}
C^{*}_{j_2 j_2 j_1 j_1}+C^{**}_{j_2 j_2 j_1 j_1}.
\end{equation}

\vspace{3mm}

Applying the estimates (\ref{ogo24}), (\ref{ogo25}), (\ref{ogoo25}),
and (\ref{after1940}) ($\varepsilon=1/2$) to 
(\ref{sixsix107}) 
gives the following estimates 
\begin{equation}
\label{sixsix109}
|C_{j_2 j_3 j_3 j_2 j_1 j_1}|\le \frac{K}{j_1^2 j_2 j_3^{3/4}}\ \ \ (j_1, j_2, j_3>0),
\end{equation}
\begin{equation}
\label{sixsix110}
|C_{j_2 00 j_2 j_1 j_1}|\le \frac{K}{j_1^2 j_2}\ \ \ (j_1, j_2> 0),
\end{equation}
\begin{equation}
\label{sixsix111}
|C_{0 j_3 j_3 0 j_1 j_1}|\le \frac{K}{j_1^2 j_3}\ \ \ (j_1, j_3>0),
\end{equation}
\begin{equation}
\label{sixsix112}
|C_{0000 j_1 j_1}|\le \frac{K}{j_1^2}\ \ \ (j_1> 0).
\end{equation}

\vspace{2mm}

Using the estimate (\ref{sixsix109}), we have
$$
\left\vert
\sum_{j_2=1}^{p}\sum_{j_3=1}^{p}\sum_{j_1=p+1}^{\infty}
C_{j_2 j_3 j_3 j_2 j_1 j_1}
\right\vert \le 
K\sum_{j_1=p+1}^{\infty}\frac{1}{j_1^2} \sum_{j_2=1}^{p} \frac{1}{j_2}
\sum_{j_3=1}^{p}\frac{1}{j_3^{3/4}}\le
$$
\begin{equation}
\label{sixsix113}
~~~~~~\le K \int\limits_p^{\infty} \frac{dx}{x^2}\left(1+\int\limits_1^p \frac{dx}{x}\right)
\left(1+\int\limits_1^p \frac{dx}{x^{3/4}}\right)\le
K_1 \frac{1+ ln p}{p^{3/4}}\ \to 0
\end{equation}

\vspace{2mm}
\noindent
if $p\to\infty,$ where constants $K, K_1$ do not depend on $p.$

Similarly we get (see (\ref{sixsix110})--(\ref{sixsix112}))
\begin{equation}
\label{sixsix114}
~~~~\left\vert\sum_{j_1=p+1}^{\infty}
C_{0000 j_1 j_1}\right\vert+
\left\vert\sum_{j_3=1}^{p}\sum_{j_1=p+1}^{\infty}
C_{0 j_3 j_3 0 j_1 j_1}\right\vert+
\left\vert\sum_{j_2=1}^{p}\sum_{j_1=p+1}^{\infty}
C_{j_2 00 j_2 j_1 j_1}\right\vert\ \to 0
\end{equation}

\noindent
if $p\to\infty.$

The relations 
(\ref{sixsix108}), (\ref{sixsix109a}), (\ref{sixsix113}), 
(\ref{sixsix114}) prove (\ref{sixsix7}).

Consider (\ref{sixsix6}). Using the integration
order replacement, we get
$$
C_{j_3 j_2 j_3 j_2 j_1 j_1}=
$$
$$
=
\frac{1}{2}\int\limits_t^T 
\phi_{j_3}(t_6)
\int\limits_t^{t_6} 
\phi_{j_2}(t_5)
\int\limits_{t}^{t_5} 
\phi_{j_3}(t_4) 
\int\limits_t^{t_4}
\phi_{j_2}(t_3) 
\left(\int\limits_t^{t_3}
\phi_{j_1}(t_1)dt_1\right)^2 
dt_3 dt_4 dt_5 dt_6=
$$
$$
=
\frac{1}{2}\int\limits_t^T 
\phi_{j_2}(t_3) 
\left(\int\limits_t^{t_3}
\phi_{j_1}(t_1)dt_1\right)^2 
\int\limits_{t_3}^T
\phi_{j_3}(t_4) 
\int\limits_{t_4}^T
\phi_{j_2}(t_5)
\int\limits_{t_5}^T
\phi_{j_3}(t_6)dt_6
dt_5 dt_4 dt_3=
$$
$$
=
\frac{1}{2}\int\limits_t^T 
\phi_{j_2}(t_3) 
\left(\int\limits_t^{t_3}
\phi_{j_1}(t_1)dt_1\right)^2 
\int\limits_{t_3}^T
\phi_{j_2}(t_5)
\int\limits_{t_5}^T
\phi_{j_3}(t_6)dt_6
\int\limits_{t_3}^{t_5}
\phi_{j_3}(t_4) 
dt_4 dt_5 dt_3=
$$
$$
=
\frac{1}{2}\int\limits_t^T 
\phi_{j_2}(t_3) 
\left(\int\limits_t^{t_3}
\phi_{j_1}(t_1)dt_1\right)^2 
\int\limits_{t_3}^T
\phi_{j_2}(t_5)
\left(\int\limits_{t}^{t_5}
\phi_{j_3}(t_4) 
dt_4\right)
\left(\int\limits_{t_5}^T
\phi_{j_3}(t_6)dt_6\right)\times
$$
$$
\times
dt_5 dt_3-
$$
$$
-
\frac{1}{2}\int\limits_t^T 
\phi_{j_2}(t_3) 
\left(\int\limits_t^{t_3}
\phi_{j_1}(t_1)dt_1\right)^2 
\left(\int\limits_{t}^{t_3}
\phi_{j_3}(t_4) 
dt_4 \right)
\int\limits_{t_3}^T
\phi_{j_2}(t_5)
\left(\int\limits_{t_5}^T
\phi_{j_3}(t_6)dt_6\right)\times
$$
\begin{equation}
\label{sixsix115}
\times
dt_5 dt_3.
\end{equation}

Applying (\ref{after80xx}), we obtain
$$
\sum_{j_1=p+1}^{\infty}\sum_{j_2=p+1}^{\infty}
\sum_{j_3=p+1}^{\infty}
C_{j_3 j_2 j_3 j_2 j_1 j_1}=
\sum_{j_1=p+1}^{\infty}\sum_{j_3=p+1}^{\infty}
\sum_{j_2=p+1}^{\infty}
C_{j_3 j_2 j_3 j_2 j_1 j_1}=
$$
\begin{equation}
\label{sixsix116}
=-
\sum_{j_2=0}^{p}
\sum_{j_1=p+1}^{\infty}\sum_{j_3=p+1}^{\infty}
C_{j_3 j_2 j_3 j_2 j_1 j_1}.
\end{equation}

\vspace{2mm}

Further proof of the equality (\ref{sixsix6})
is based on the relations (\ref{sixsix115}), (\ref{sixsix116}) and 
is similar to the proof of the formula (\ref{sixsix7}).

Let us prove (\ref{sixsix1}). Applying the integration 
order replacement, we obtain
$$
C_{j_3 j_3 j_2 j_1 j_2 j_1}=
$$
$$
=\hspace{-1mm}
\int\limits_t^T 
\phi_{j_3}(t_6)
\int\limits_t^{t_6} 
\phi_{j_3}(t_5)
\int\limits_{t}^{t_5} 
\phi_{j_2}(t_4) 
\int\limits_t^{t_4}
\phi_{j_1}(t_3) 
\int\limits_t^{t_3}
\phi_{j_2}(t_2)
\int\limits_t^{t_2}
\phi_{j_1}(t_1)
dt_1 dt_2 dt_3 dt_4 dt_5 dt_6=
$$
$$
=\hspace{-1mm}
\int\limits_t^T 
\phi_{j_1}(t_1)
\int\limits_{t_1}^T
\phi_{j_2}(t_2)
\int\limits_{t_2}^T
\phi_{j_1}(t_3) 
\int\limits_{t_3}^T
\phi_{j_2}(t_4) 
\int\limits_{t_4}^T
\phi_{j_3}(t_5)
\int\limits_{t_5}^T
\phi_{j_3}(t_6)
dt_6 dt_5 dt_4 dt_3 dt_2 dt_1=
$$
$$
=\frac{1}{2}
\int\limits_t^T 
\phi_{j_1}(t_1)
\int\limits_{t_1}^T
\phi_{j_2}(t_2)
\int\limits_{t_2}^T
\phi_{j_1}(t_3) 
\int\limits_{t_3}^T
\phi_{j_2}(t_4) 
\left(\int\limits_{t_4}^T
\phi_{j_3}(t_5)
dt_5\right)^2 dt_4 dt_3 dt_2 dt_1=
$$
$$
=\frac{1}{2}
\int\limits_t^T 
\phi_{j_2}(t_4) 
\left(\int\limits_{t_4}^T
\phi_{j_3}(t_5)
dt_5\right)^2
\int\limits_t^{t_4}
\phi_{j_1}(t_3) 
\int\limits_t^{t_3}
\phi_{j_2}(t_2)
\int\limits_t^{t_2}
\phi_{j_1}(t_1)
dt_1 dt_2 dt_3 dt_4=
$$
$$
=\frac{1}{2}
\int\limits_t^T 
\phi_{j_2}(t_4) 
\left(\int\limits_{t_4}^T
\phi_{j_3}(t_5)
dt_5\right)^2
\int\limits_t^{t_4}
\phi_{j_2}(t_2)
\int\limits_t^{t_2}
\phi_{j_1}(t_1)
dt_1
\int\limits_{t_2}^{t_4}
\phi_{j_1}(t_3) 
dt_3 dt_2 dt_4=
$$
$$
=\frac{1}{2}
\int\limits_t^T 
\phi_{j_2}(t_4) 
\left(\int\limits_{t_4}^T
\phi_{j_3}(t_5)
dt_5\right)^2
\left(\int\limits_{t}^{t_4}
\phi_{j_1}(t_3) 
dt_3\right)
\int\limits_t^{t_4}
\phi_{j_2}(t_2)
\left(\int\limits_t^{t_2}
\phi_{j_1}(t_1)
dt_1\right)\times
$$
$$
\times
dt_2 dt_4-
$$
$$
-\frac{1}{2}
\int\limits_t^T 
\phi_{j_2}(t_4) 
\left(\int\limits_{t_4}^T
\phi_{j_3}(t_5)
dt_5\right)^2
\int\limits_t^{t_4}
\phi_{j_2}(t_2)
\left(\int\limits_t^{t_2}
\phi_{j_1}(t_1)
dt_1\right)^2\times
$$
\begin{equation}
\label{sixsix120}
\times
dt_2 dt_4.
\end{equation}

\vspace{2mm}

Using (\ref{after80xx}), we get
$$
\sum_{j_1=p+1}^{\infty}\sum_{j_2=p+1}^{\infty}
\sum_{j_3=p+1}^{\infty}
C_{j_3 j_3 j_2 j_1 j_2 j_1}=
\sum_{j_1=p+1}^{\infty}\sum_{j_3=p+1}^{\infty}
\sum_{j_2=p+1}^{\infty}
C_{j_3 j_3 j_2 j_1 j_2 j_1}=
$$
\begin{equation}
\label{sixsix121}
=-
\sum_{j_2=0}^{p}
\sum_{j_1=p+1}^{\infty}\sum_{j_3=p+1}^{\infty}
C_{j_3 j_3 j_2 j_1 j_2 j_1}.
\end{equation}

\vspace{2mm}

Further proof of the equality (\ref{sixsix1})
is based on the relations (\ref{sixsix120}), (\ref{sixsix121}) and 
is similar to the proof of the relations (\ref{sixsix7}),
(\ref{sixsix6}).

Consider (\ref{sixsix2}). 
Using the integration 
order replacement, we have
$$
C_{j_3 j_3 j_1 j_2 j_2 j_1}=
$$
$$
=\hspace{-1mm}
\int\limits_t^T 
\phi_{j_3}(t_6)
\int\limits_t^{t_6} 
\phi_{j_3}(t_5)
\int\limits_{t}^{t_5} 
\phi_{j_1}(t_4) 
\int\limits_t^{t_4}
\phi_{j_2}(t_3) 
\int\limits_t^{t_3}
\phi_{j_2}(t_2)
\int\limits_t^{t_2}
\phi_{j_1}(t_1)
dt_1 dt_2 dt_3 dt_4 dt_5 dt_6=
$$
$$
=\hspace{-1mm}
\int\limits_t^T 
\phi_{j_1}(t_1)
\int\limits_{t_1}^T
\phi_{j_2}(t_2)
\int\limits_{t_2}^T
\phi_{j_2}(t_3) 
\int\limits_{t_3}^T
\phi_{j_1}(t_4) 
\int\limits_{t_4}^T
\phi_{j_3}(t_5)
\int\limits_{t_5}^T
\phi_{j_3}(t_6)
dt_6 dt_5 dt_4 dt_3 dt_2 dt_1=
$$
$$
=\frac{1}{2}
\int\limits_t^T 
\phi_{j_1}(t_1)
\int\limits_{t_1}^T
\phi_{j_2}(t_2)
\int\limits_{t_2}^T
\phi_{j_2}(t_3) 
\int\limits_{t_3}^T
\phi_{j_1}(t_4) 
\left(\int\limits_{t_4}^T
\phi_{j_3}(t_5)
dt_5\right)^2 dt_4 dt_3 dt_2 dt_1=
$$
$$
=\frac{1}{2}
\int\limits_t^T 
\phi_{j_1}(t_4) 
\left(\int\limits_{t_4}^T
\phi_{j_3}(t_5)
dt_5\right)^2
\int\limits_t^{t_4}
\phi_{j_2}(t_3) 
\int\limits_t^{t_3}
\phi_{j_2}(t_2)
\int\limits_t^{t_2}
\phi_{j_1}(t_1)
dt_1 dt_2 dt_3 dt_4=
$$
$$
=\frac{1}{2}
\int\limits_t^T 
\phi_{j_1}(t_4) 
\left(\int\limits_{t_4}^T
\phi_{j_3}(t_5)
dt_5\right)^2
\int\limits_t^{t_4}
\phi_{j_2}(t_2)
\int\limits_t^{t_2}
\phi_{j_1}(t_1)
dt_1
\int\limits_{t_2}^{t_4}
\phi_{j_2}(t_3) 
dt_3 dt_2 dt_4=
$$
$$
=\frac{1}{2}
\int\limits_t^T 
\phi_{j_1}(t_4) 
\left(\int\limits_{t_4}^T
\phi_{j_3}(t_5)
dt_5\right)^2
\left(\int\limits_{t}^{t_4}
\phi_{j_2}(t_3) 
dt_3\right)
\int\limits_t^{t_4}
\phi_{j_2}(t_2)
\left(\int\limits_t^{t_2}
\phi_{j_1}(t_1)
dt_1\right)\times
$$
$$
\times
dt_2 dt_4-
$$
$$
-\frac{1}{2}
\int\limits_t^T 
\phi_{j_1}(t_4) 
\left(\int\limits_{t_4}^T
\phi_{j_3}(t_5)
dt_5\right)^2
\int\limits_t^{t_4}
\phi_{j_2}(t_2)
\left(\int\limits_t^{t_2}
\phi_{j_1}(t_1)
dt_1\right)
\left(\int\limits_t^{t_2}
\phi_{j_2}(t_3)
dt_3\right)
\times
$$
\begin{equation}
\label{sixsix123}
\times
dt_2 dt_4.
\end{equation}

\vspace{2mm}

Applying (\ref{after80xx}) and (\ref{after500}), we obtain
$$
-\sum_{j_1=p+1}^{\infty}\sum_{j_2=p+1}^{\infty}
\sum_{j_3=p+1}^{\infty}
C_{j_3 j_3 j_1 j_2 j_2 j_1}=
-\sum_{j_2=p+1}^{\infty}\sum_{j_3=p+1}^{\infty}\sum_{j_1=p+1}^{\infty}
C_{j_2 j_3 j_1 j_2 j_2 j_1}=
$$
$$
=
\sum_{j_1=0}^{p}\sum_{j_2=p+1}^{\infty}\sum_{j_3=p+1}^{\infty}
C_{j_2 j_3 j_1 j_2 j_2 j_1}=
\sum_{j_1=0}^{p}\sum_{j_3=p+1}^{\infty}\sum_{j_2=p+1}^{\infty}
C_{j_2 j_3 j_1 j_2 j_2 j_1}=
$$
\begin{equation}
\label{sixsix124}
~~~~~~~~=\frac{1}{2}
\sum_{j_1=0}^{p}\sum_{j_3=p+1}^{\infty}
C_{j_3 j_3 j_1 j_2 j_2 j_1}\biggl|_{(j_2 j_2)\curvearrowright (\cdot)}\biggr. -
\sum_{j_1=0}^{p}\sum_{j_2=0}^{p}\sum_{j_3=p+1}^{\infty}
C_{j_3 j_3 j_1 j_2 j_2 j_1}.
\end{equation}

\vspace{4mm}

The equality
\begin{equation}
\label{sixsix125}
\lim\limits_{p\to\infty}\frac{1}{2}
\sum_{j_1=0}^{p}\sum_{j_3=p+1}^{\infty}
C_{j_3 j_3 j_1 j_2 j_2 j_1}\biggl|_{(j_2 j_2)\curvearrowright (\cdot)}\biggr. =0
\end{equation}

\vspace{1mm}
\noindent       
follows from the inequality (\ref{after9043}),
where we proceed similarly to the proof of equality (\ref{sixsix109a})
(see (\ref{sept12})).

The relation
\begin{equation}
\label{sixsix126}
\lim\limits_{p\to\infty}
\sum_{j_1=0}^{p}\sum_{j_2=0}^{p}\sum_{j_3=p+1}^{\infty}
C_{j_3 j_3 j_1 j_2 j_2 j_1}=0
\end{equation}

\vspace{1mm}
\noindent
is proved on the basis of (\ref{sixsix123}) and similarly with the proof 
of (\ref{sixsix7}).
The equalities (\ref{sixsix124})--(\ref{sixsix126}) prove
(\ref{sixsix2}). 

Let us prove (\ref{sixsix5}). Using (\ref{after80xx}) and (\ref{after500}), we get
$$
\sum_{j_1=p+1}^{\infty}\sum_{j_2=p+1}^{\infty}
\sum_{j_3=p+1}^{\infty}
C_{j_2 j_1 j_3 j_3 j_2 j_1}=
\sum_{j_3=p+1}^{\infty}\sum_{j_1, j_2= 0}^{p}
C_{j_2 j_1 j_3 j_3 j_2 j_1}=
$$
\begin{equation}
\label{sixsix127}
=\frac{1}{2}
\sum_{j_1,j_2=0}^{p}
C_{j_2 j_1 j_3 j_3 j_2 j_1}\biggl|_{(j_3 j_3)\curvearrowright (\cdot)}\biggr.
-
\sum_{j_1,j_2,j_3=0}^{p}
C_{j_2 j_1 j_3 j_3 j_2 j_1}.
\end{equation}

\vspace{2mm}

Using the equality (\ref{after2508}) we have
\begin{equation}
\label{sixsix128}
\lim\limits_{p\to\infty}
\frac{1}{2}
\sum_{j_1,j_2=0}^{p}
C_{j_2 j_1 j_3 j_3 j_2 j_1}\biggl|_{(j_3 j_3)\curvearrowright (\cdot)}\biggr.
=0,
\end{equation}

\vspace{2mm}
\noindent
where we proceed similarly to the proof of equality (\ref{sixsix109a})
(see (\ref{sept12})).

Further, we will prove the following relation
\begin{equation}
\label{sixsix129}
\lim\limits_{p\to\infty}
\sum_{j_1,j_2,j_3=0}^{p}
C_{j_2 j_1 j_3 j_3 j_2 j_1}=0
\end{equation}

\vspace{1mm}
\noindent
using the equality (\ref{sixsix40}). From (\ref{sixsix40}) we have
$$
\sum_{j_1,j_2,j_3=0}^{p}
C_{j_2 j_1 j_3 j_3 j_2 j_1}
=\frac{1}{2}\sum_{j_1,j_2,j_3=0}^{p}
\biggl(C_{j_2 j_1 j_3 j_3 j_2 j_1}+
C_{j_1 j_2 j_3 j_3 j_1 j_2}\biggr)=
$$
$$
=
\frac{1}{2}\sum_{j_1,j_2,j_3=0}^{p}\biggl(
C_{j_2}C_{j_1 j_3 j_3 j_2 j_1}-C_{j_1 j_2}C_{j_3 j_3 j_2 j_1}
+C_{j_3 j_1 j_2}C_{j_3 j_2 j_1}-\biggr.
$$
$$
\biggl.-C_{j_3 j_3 j_1 j_2}C_{j_2 j_1}+
C_{j_2 j_3 j_3 j_1 j_2}C_{j_1}\biggr)=
$$
$$
=
\sum_{j_1,j_2,j_3=0}^{p}\biggl(
C_{j_2 j_3 j_3 j_1 j_2}C_{j_1}-C_{j_3 j_3 j_1 j_2}C_{j_2 j_1}\biggr)+
$$
\begin{equation}
\label{sixsix130}
+
\frac{1}{2}\sum_{j_1,j_2,j_3=0}^{p}
C_{j_3 j_1 j_2}C_{j_3 j_2 j_1}.
\end{equation}

\vspace{2mm}

The generalized Parseval equality gives (by analogy with (\ref{pars100}))
\begin{equation}
\label{sixsix131}
\lim\limits_{p\to\infty}\frac{1}{2}\sum_{j_1,j_2,j_3=0}^{p}
C_{j_3 j_1 j_2}C_{j_3 j_2 j_1}=0.
\end{equation}

\vspace{2mm}

Let us prove the following equality
\begin{equation}
\label{sixsix132}
\lim\limits_{p\to\infty}\sum_{j_1,j_2,j_3=0}^{p}\biggl(
C_{j_2 j_3 j_3 j_1 j_2}C_{j_1}-C_{j_3 j_3 j_1 j_2}C_{j_2 j_1}\biggr)=0.
\end{equation}

\vspace{1mm}

The relation
\begin{equation}
\label{sixsix132s}
\lim\limits_{p\to\infty}\sum_{j_1,j_2,j_3=0}^{p}
C_{j_2 j_3 j_3 j_1 j_2}C_{j_1}=0
\end{equation}

\noindent
is proved by the same methods as 
in the proof of equality (\ref{sixsix8}) and also using 
Theorem~2.34 and (\ref{after500}).

Further, we have (see (\ref{after500}))
\begin{equation}
\label{sept2}
\sum_{j_3=0}^{p}
C_{j_3 j_3 j_1 j_2}=
\frac{1}{2}
C_{j_3 j_3 j_1 j_2}\biggl|_{(j_3 j_3)\curvearrowright (\cdot)}\biggr.-
\sum_{j_3=p+1}^{\infty}
C_{j_3 j_3 j_1 j_2}.
\end{equation}

Moreover,
$$
C_{j_3 j_3 j_1 j_2}\biggl|_{(j_3 j_3)\curvearrowright (\cdot)}\biggr.=
$$
$$
=\int\limits_t^T \int\limits_t^{t_3}\phi_{j_1}(t_2)
\int\limits_t^{t_2}\phi_{j_2}(t_1)dt_1 dt_2 dt_3=
$$
$$
=
\int\limits_t^T\phi_{j_1}(t_2)
\int\limits_t^{t_2}\phi_{j_2}(t_1)dt_1 \int\limits_{t_2}^T
dt_3 dt_2=
$$
$$
=
\int\limits_t^T(T-t_2)\phi_{j_1}(t_2)
\int\limits_t^{t_2}\phi_{j_2}(t_1)dt_1 
dt_2=
$$
$$
=\int\limits_t^T\phi_{j_2}(t_1)
\int\limits_{t_1}^T(T-t_2)\phi_{j_1}(t_2)dt_2 dt_1= 
$$
$$
=
\int\limits_t^T\phi_{j_2}(t_2)
\int\limits_{t_2}^T(T-t_1)\phi_{j_1}(t_1)dt_1 dt_2= 
$$
$$
=\int\limits_{[t,T]^2}
(T-t_1){\bf 1}_{\{t_2<t_1\}}\phi_{j_1}(t_1)\phi_{j_2}(t_2)dt_1 dt_2
\stackrel{\sf def}{=}
$$
\begin{equation}
\label{sept3}
\stackrel{\sf def}{=}\tilde C_{j_2 j_1}.
\end{equation}

\vspace{2mm}

Using (\ref{sept2}), (\ref{sept3}), and the generalized Parseval equality, we obtain 
$$
\lim\limits_{p\to\infty}\sum_{j_1,j_2,j_3=0}^{p}
C_{j_3 j_3 j_1 j_2}C_{j_2 j_1}=\frac{1}{2}\lim\limits_{p\to\infty}\sum_{j_1,j_2=0}^{p}
C_{j_2 j_1}\tilde C_{j_2 j_1}-
$$
\begin{equation}
\label{sept7}
~~~~~~~~~-\lim\limits_{p\to\infty}\sum_{j_1,j_2=0}^{p}\sum_{j_3=p+1}^{\infty}
C_{j_3 j_3 j_1 j_2}C_{j_2 j_1}=
-\lim\limits_{p\to\infty}\sum_{j_1,j_2=0}^{p}\sum_{j_3=p+1}^{\infty}
C_{j_3 j_3 j_1 j_2}C_{j_2 j_1}.
\end{equation}

\vspace{3mm}

We have (see (\ref{sept10}))
\begin{equation}
\label{sept5}
~~~~~~~C_{j_3 j_3 j_1 j_2}=
\frac{1}{2}\int\limits_t^T \phi_{j_2}(t_1)
\int\limits_{t_1}^T \phi_{j_1}(t_2)
\left(\int\limits_{t_2}^T \phi_{j_3}(t_3)
dt_3\right)^2 dt_2 dt_1.
\end{equation}

\vspace{2mm}
                        
By analogy with (\ref{sixsix74}) and also using (\ref{sept5}), we get
\begin{equation}
\label{sept6}
\lim\limits_{p\to\infty}\sum_{j_1,j_2=0}^{p}\sum_{j_3=p+1}^{\infty}
C_{j_3 j_3 j_1 j_2}C_{j_2 j_1}=0.
\end{equation}

\vspace{2mm}

Combining (\ref{sept7}) and (\ref{sept6}), we obtain
\begin{equation}
\label{sept9}
\lim\limits_{p\to\infty}\sum_{j_1,j_2,j_3=0}^{p}
C_{j_3 j_3 j_1 j_2}C_{j_2 j_1}=0.
\end{equation}

The relation (\ref{sixsix132}) follows from (\ref{sixsix132s}) and (\ref{sept9}).
From (\ref{sixsix130})--(\ref{sixsix132}) we get (\ref{sixsix129}).
The equalities (\ref{sixsix127})--(\ref{sixsix129})
complete the proof of (\ref{sixsix5}).

For the proof of (\ref{sixsix12})--(\ref{sixsix15})
we will use a new idea.  
More precisely, we will consider the sums of 
expressions (\ref{sixsix12})--(\ref{sixsix15}) with the expressions 
already studied throughout this proof.

Let us begin from (\ref{sixsix12}). 
Applying the integration order replacement, we obtain
$$
C_{j_3 j_1 j_2 j_3 j_2 j_1}+C_{j_3 j_1 j_2 j_3 j_1 j_2}=
$$
$$
=
\int\limits_t^T 
\phi_{j_3}(t_6)
\int\limits_t^{t_6} 
\phi_{j_1}(t_5)
\int\limits_{t}^{t_5} 
\phi_{j_2}(t_4) 
\int\limits_t^{t_4}
\phi_{j_3}(t_3) 
\left(\int\limits_t^{t_3}
\phi_{j_2}(t_2)dt_2\right)
\left(\int\limits_t^{t_3}
\phi_{j_1}(t_1)
dt_1\right)\times
$$
$$
\times 
dt_3 dt_4 dt_5 dt_6=
$$
$$
=
\int\limits_t^T 
\phi_{j_3}(t_6)
\int\limits_t^{t_6} 
\phi_{j_1}(t_5)
\int\limits_{t}^{t_5} 
\phi_{j_3}(t_3) 
\left(\int\limits_t^{t_3}
\phi_{j_2}(t_2)dt_2\right)
\left(\int\limits_t^{t_3}
\phi_{j_1}(t_1)
dt_1\right)
\int\limits_{t_3}^{t_5}
\phi_{j_2}(t_4) 
dt_4
\times
$$
$$
\times 
dt_3 dt_5 dt_6=
$$
$$
=
\int\limits_t^T 
\phi_{j_3}(t_6)
\int\limits_t^{t_6} 
\phi_{j_1}(t_5)
\left(\int\limits_{t}^{t_5}
\phi_{j_2}(t_4) 
dt_4
\right)
\int\limits_{t}^{t_5} 
\phi_{j_3}(t_3) 
\left(\int\limits_t^{t_3}
\phi_{j_2}(t_2)dt_2\right)
\times
$$
$$
\times\left(\int\limits_t^{t_3}
\phi_{j_1}(t_1)
dt_1\right) 
dt_3 dt_5 dt_6-
$$
$$
-
\int\limits_t^T 
\phi_{j_3}(t_6)
\int\limits_t^{t_6} 
\phi_{j_1}(t_5)
\int\limits_{t}^{t_5} 
\phi_{j_3}(t_3) 
\left(\int\limits_t^{t_3}
\phi_{j_2}(t_2)dt_2\right)^2
\left(\int\limits_t^{t_3}
\phi_{j_1}(t_1)
dt_1\right)
\times
$$
$$
\times 
dt_3 dt_5 dt_6=
$$
$$
=
\int\limits_t^T 
\phi_{j_1}(t_5)
\left(\int\limits_{t}^{t_5}
\phi_{j_2}(t_4) 
dt_4
\right)
\int\limits_{t}^{t_5} 
\phi_{j_3}(t_3) 
\left(\int\limits_t^{t_3}
\phi_{j_2}(t_2)dt_2\right)
\times
$$
$$
\times\left(\int\limits_t^{t_3}
\phi_{j_1}(t_1)
dt_1\right) 
dt_3
\left(\int\limits_{t_5}^T
\phi_{j_3}(t_6)dt_6\right) dt_5-
$$
$$
-
\int\limits_t^T 
\phi_{j_1}(t_5)
\int\limits_{t}^{t_5} 
\phi_{j_3}(t_3) 
\left(\int\limits_t^{t_3}
\phi_{j_2}(t_2)dt_2\right)^2
\left(\int\limits_t^{t_3}
\phi_{j_1}(t_1)
dt_1\right)
dt_3\times
$$
\begin{equation}
\label{sixsix133}
\times
\left(
\int\limits_{t_5}^T
\phi_{j_3}(t_6)
dt_6\right)dt_5.
\end{equation}

\vspace{2mm}

Using (\ref{after80xx}), we get
$$
\sum_{j_1=p+1}^{\infty}\sum_{j_2=p+1}^{\infty}
\sum_{j_3=p+1}^{\infty}
\biggl(
C_{j_3 j_1 j_2 j_3 j_2 j_1}+C_{j_3 j_1 j_2 j_3 j_1 j_2}\biggr)=
$$
\begin{equation}
\label{sixsix134}
=
\sum_{j_1=0}^{p}\sum_{j_3=0}^{p}
\sum_{j_2=p+1}^{\infty}
\biggl(
C_{j_3 j_1 j_2 j_3 j_2 j_1}+C_{j_3 j_1 j_2 j_3 j_1 j_2}\biggr).
\end{equation}

\vspace{2mm}

Further, by analogy with the proof of equality (\ref{sixsix7}) 
and using (\ref{sixsix133}), we obtain
\begin{equation}
\label{sixsix135}
\lim\limits_{p\to\infty}\sum_{j_1=0}^{p}\sum_{j_3=0}^{p}
\sum_{j_2=p+1}^{\infty}
\biggl(
C_{j_3 j_1 j_2 j_3 j_2 j_1}+C_{j_3 j_1 j_2 j_3 j_1 j_2}\biggr)=0.
\end{equation}

\vspace{2mm}

From (\ref{sixsix134}) and (\ref{sixsix135}) we get
\begin{equation}
\label{sixsix136}
~~~~~~~~\lim\limits_{p\to\infty}\sum_{j_1=p+1}^{\infty}\sum_{j_2=p+1}^{\infty}
\sum_{j_3=p+1}^{\infty}
\biggl(
C_{j_3 j_1 j_2 j_3 j_2 j_1}+C_{j_3 j_1 j_2 j_3 j_1 j_2}\biggr)=0.
\end{equation}

\vspace{2mm}

Moreover (see (\ref{sixsix8})),
\begin{equation}
\label{sixsix137}
\lim\limits_{p\to\infty}\sum_{j_1=p+1}^{\infty}\sum_{j_2=p+1}^{\infty}
\sum_{j_3=p+1}^{\infty}
C_{j_3 j_1 j_2 j_3 j_1 j_2}=0.
\end{equation}

\vspace{2mm}

Combining (\ref{sixsix136}) and (\ref{sixsix137}), we have
$$
\lim\limits_{p\to\infty}\sum_{j_1=p+1}^{\infty}\sum_{j_2=p+1}^{\infty}
\sum_{j_3=p+1}^{\infty}
C_{j_3 j_1 j_2 j_3 j_2 j_1}=0.
$$

\vspace{2mm}
\noindent
The equality (\ref{sixsix12}) is proved.

Consider (\ref{sixsix11}).
Using the integration order replacement, we have
$$
C_{j_2 j_3 j_1 j_3 j_2 j_1}+C_{j_2 j_3 j_1 j_3 j_1 j_2}=
$$
$$
=
\int\limits_t^T 
\phi_{j_2}(t_6)
\int\limits_t^{t_6} 
\phi_{j_3}(t_5)
\int\limits_{t}^{t_5} 
\phi_{j_1}(t_4) 
\int\limits_t^{t_4}
\phi_{j_3}(t_3) 
\left(\int\limits_t^{t_3}
\phi_{j_2}(t_2)dt_2\right)
\left(\int\limits_t^{t_3}
\phi_{j_1}(t_1)
dt_1\right)\times
$$
$$
\times 
dt_3 dt_4 dt_5 dt_6=
$$
$$
=
\int\limits_t^T 
\phi_{j_2}(t_6)
\int\limits_t^{t_6} 
\phi_{j_3}(t_5)
\int\limits_{t}^{t_5} 
\phi_{j_3}(t_3) 
\left(\int\limits_t^{t_3}
\phi_{j_2}(t_2)dt_2\right)
\left(\int\limits_t^{t_3}
\phi_{j_1}(t_1)
dt_1\right)
\int\limits_{t_3}^{t_5}
\phi_{j_1}(t_4) 
dt_4
\times
$$
$$
\times 
dt_3 dt_5 dt_6=
$$
$$
=
\int\limits_t^T 
\phi_{j_2}(t_6)
\int\limits_t^{t_6} 
\phi_{j_3}(t_5)
\left(\int\limits_{t}^{t_5}
\phi_{j_1}(t_4) 
dt_4
\right)
\int\limits_{t}^{t_5} 
\phi_{j_3}(t_3) 
\left(\int\limits_t^{t_3}
\phi_{j_2}(t_2)dt_2\right)
\times
$$
$$
\times\left(\int\limits_t^{t_3}
\phi_{j_1}(t_1)
dt_1\right) 
dt_3 dt_5 dt_6-
$$
$$
-
\int\limits_t^T 
\phi_{j_2}(t_6)
\int\limits_t^{t_6} 
\phi_{j_3}(t_5)
\int\limits_{t}^{t_5} 
\phi_{j_3}(t_3) 
\left(\int\limits_t^{t_3}
\phi_{j_2}(t_2)dt_2\right)
\left(\int\limits_t^{t_3}
\phi_{j_1}(t_1)
dt_1\right)^2
\times
$$
$$
\times 
dt_3 dt_5 dt_6=
$$
$$
=
\int\limits_t^T 
\phi_{j_3}(t_5)
\left(\int\limits_{t}^{t_5}
\phi_{j_1}(t_4) 
dt_4
\right)
\int\limits_{t}^{t_5} 
\phi_{j_3}(t_3) 
\left(\int\limits_t^{t_3}
\phi_{j_2}(t_2)dt_2\right)
\times
$$
$$
\times\left(\int\limits_t^{t_3}
\phi_{j_1}(t_1)
dt_1\right) 
dt_3
\left(\int\limits_{t_5}^T
\phi_{j_2}(t_6)dt_6\right) dt_5-
$$
$$
-
\int\limits_t^T 
\phi_{j_3}(t_5)
\int\limits_{t}^{t_5} 
\phi_{j_3}(t_3) 
\left(\int\limits_t^{t_3}
\phi_{j_2}(t_2)dt_2\right)
\left(\int\limits_t^{t_3}
\phi_{j_1}(t_1)
dt_1\right)^2
dt_3\times
$$
\begin{equation}
\label{sixsix150}
\times
\left(
\int\limits_{t_5}^T
\phi_{j_2}(t_6)
dt_6\right)dt_5.
\end{equation}

\vspace{2mm}

Using (\ref{after80xx}), we obtain
$$
-\sum_{j_1=p+1}^{\infty}\sum_{j_2=p+1}^{\infty}
\sum_{j_3=p+1}^{\infty}
\biggl(
C_{j_2 j_3 j_1 j_3 j_2 j_1}+C_{j_2 j_3 j_1 j_3 j_1 j_2}\biggr)=
$$
\begin{equation}
\label{sixsix151}
=
\sum_{j_3=0}^{p}
\sum_{j_1=p+1}^{\infty}\sum_{j_2=p+1}^{\infty}
\biggl(
C_{j_2 j_3 j_1 j_3 j_2 j_1}+C_{j_2 j_3 j_1 j_3 j_1 j_2}\biggr).
\end{equation}

\vspace{3mm}

By analogy with the proof of (\ref{sixsix7}) 
and applying (\ref{sixsix150}), we get
\begin{equation}
\label{sixsix152}
~~~~~~~~~~\lim\limits_{p\to\infty}\sum_{j_3=0}^{p}
\sum_{j_1=p+1}^{\infty}\sum_{j_2=p+1}^{\infty}
\biggl(
C_{j_2 j_3 j_1 j_3 j_2 j_1}+C_{j_2 j_3 j_1 j_3 j_1 j_2}\biggr)=0.
\end{equation}

\vspace{2mm}

From (\ref{sixsix151}) and (\ref{sixsix152}) we have
\begin{equation}
\label{sixsix153}
~~~~~~~~~~\lim\limits_{p\to\infty}
\sum_{j_1=p+1}^{\infty}\sum_{j_2=p+1}^{\infty}
\sum_{j_3=p+1}^{\infty}
\biggl(
C_{j_2 j_3 j_1 j_3 j_2 j_1}+C_{j_2 j_3 j_1 j_3 j_1 j_2}\biggr)=0.
\end{equation}

\vspace{2mm}

Moreover (see (\ref{sixsix9})),
\begin{equation}
\label{sixsix154}
\lim\limits_{p\to\infty}\sum_{j_1=p+1}^{\infty}\sum_{j_2=p+1}^{\infty}
\sum_{j_3=p+1}^{\infty}
C_{j_2 j_3 j_1 j_3 j_1 j_2}=0.
\end{equation}

\vspace{2mm}

Combining (\ref{sixsix153}) and (\ref{sixsix154}), we finally obtain
$$
\lim\limits_{p\to\infty}\sum_{j_1=p+1}^{\infty}\sum_{j_2=p+1}^{\infty}
\sum_{j_3=p+1}^{\infty}
C_{j_2 j_3 j_1 j_3 j_2 j_1}=0.
$$

\vspace{2mm}
\noindent
The equality (\ref{sixsix11}) is proved.

Now consider (\ref{sixsix13}).
Using the integration order replacement, we obtain
$$
C_{j_3 j_1 j_3 j_2 j_2 j_1}+C_{j_3 j_1 j_3 j_2 j_1 j_2}=
$$
$$
=
\int\limits_t^T 
\phi_{j_3}(t_6)
\int\limits_t^{t_6} 
\phi_{j_1}(t_5)
\int\limits_{t}^{t_5} 
\phi_{j_3}(t_4) 
\int\limits_t^{t_4}
\phi_{j_2}(t_3) 
\left(\int\limits_t^{t_3}
\phi_{j_2}(t_2)dt_2\right)
\left(\int\limits_t^{t_3}
\phi_{j_1}(t_1)
dt_1\right)\times
$$
$$
\times 
dt_3 dt_4 dt_5 dt_6=
$$
$$
=
\int\limits_t^T 
\phi_{j_3}(t_6)
\int\limits_t^{t_6} 
\phi_{j_1}(t_5)
\int\limits_{t}^{t_5} 
\phi_{j_2}(t_3) 
\left(\int\limits_t^{t_3}
\phi_{j_2}(t_2)dt_2\right)
\left(\int\limits_t^{t_3}
\phi_{j_1}(t_1)
dt_1\right)
\int\limits_{t_3}^{t_5}
\phi_{j_3}(t_4) 
dt_4
\times
$$
$$
\times 
dt_3 dt_5 dt_6=
$$
$$
=
\int\limits_t^T 
\phi_{j_3}(t_6)
\int\limits_t^{t_6} 
\phi_{j_1}(t_5)
\left(\int\limits_{t}^{t_5}
\phi_{j_3}(t_4) 
dt_4
\right)
\int\limits_{t}^{t_5} 
\phi_{j_2}(t_3) 
\left(\int\limits_t^{t_3}
\phi_{j_2}(t_2)dt_2\right)
\times
$$
$$
\times\left(\int\limits_t^{t_3}
\phi_{j_1}(t_1)
dt_1\right) 
dt_3 dt_5 dt_6-
$$
$$
-
\int\limits_t^T 
\phi_{j_3}(t_6)
\int\limits_t^{t_6} 
\phi_{j_1}(t_5)
\int\limits_{t}^{t_5} 
\phi_{j_2}(t_3) 
\left(\int\limits_t^{t_3}
\phi_{j_2}(t_2)dt_2\right)
\left(\int\limits_t^{t_3}
\phi_{j_1}(t_1)
dt_1\right)
\times
$$
$$
\times 
\left(\int\limits_t^{t_3}
\phi_{j_3}(t_4)
dt_4\right)
dt_3 dt_5 dt_6=
$$
$$
=
\int\limits_t^T 
\phi_{j_1}(t_5)
\left(\int\limits_{t}^{t_5}
\phi_{j_3}(t_4) 
dt_4
\right)
\int\limits_{t}^{t_5} 
\phi_{j_2}(t_3) 
\left(\int\limits_t^{t_3}
\phi_{j_2}(t_2)dt_2\right)
\times
$$
$$
\times\left(\int\limits_t^{t_3}
\phi_{j_1}(t_1)
dt_1\right) 
dt_3
\left(\int\limits_{t_5}^T
\phi_{j_3}(t_6)dt_6\right) dt_5-
$$
$$
-
\int\limits_t^T 
\phi_{j_1}(t_5)
\int\limits_{t}^{t_5} 
\phi_{j_2}(t_3) 
\left(\int\limits_t^{t_3}
\phi_{j_2}(t_2)dt_2\right)
\left(\int\limits_t^{t_3}
\phi_{j_1}(t_1)
dt_1\right)
\times
$$
\begin{equation}
\label{sixsix160}
\times
\left(\int\limits_t^{t_3}
\phi_{j_3}(t_4)
dt_4\right)dt_3
\left(
\int\limits_{t_5}^T
\phi_{j_3}(t_6)
dt_6\right)dt_5.
\end{equation}

\vspace{2mm}

Applying (\ref{after80xx}) and (\ref{after500}), we obtain
$$
\sum_{j_1=p+1}^{\infty}\sum_{j_2=p+1}^{\infty}
\sum_{j_3=p+1}^{\infty}
\biggl(C_{j_3 j_1 j_3 j_2 j_2 j_1}+C_{j_3 j_1 j_3 j_2 j_1 j_2}\biggr)=
$$
$$
=-\sum_{j_1=0}^{p}\sum_{j_3=p+1}^{\infty}\sum_{j_2=p+1}^{\infty}
\biggl(C_{j_3 j_1 j_3 j_2 j_2 j_1}+C_{j_3 j_1 j_3 j_2 j_1 j_2}\biggr)=
$$
$$
=
\sum_{j_1=0}^{p}\sum_{j_2=0}^{p}\sum_{j_3=p+1}^{\infty}
\biggl(C_{j_3 j_1 j_3 j_2 j_2 j_1}+C_{j_3 j_1 j_3 j_2 j_1 j_2}\biggr)-
$$
\begin{equation}
\label{sixsix170}
-\frac{1}{2}
\sum_{j_1=0}^{p}\sum_{j_3=p+1}^{\infty}
C_{j_3 j_1 j_3 j_2 j_2 j_1}\biggl|_{(j_2 j_2)\curvearrowright (\cdot)}\biggr..
\end{equation}

\vspace{2mm}

The equality
\begin{equation}
\label{sixsix171}
\lim\limits_{p\to\infty}\frac{1}{2}
\sum_{j_1=0}^{p}\sum_{j_3=p+1}^{\infty}
C_{j_3 j_1 j_3 j_2 j_2 j_1}\biggl|_{(j_2 j_2)\curvearrowright (\cdot)}\biggr. =0
\end{equation}

\vspace{1mm}
\noindent
follows from the 
equality (\ref{after2508}),
where we proceed similarly to the proof of equality (\ref{sixsix109a})
(see (\ref{sept12})).

By analogy with the proof of (\ref{sixsix7}) 
and applying (\ref{sixsix160}), we get
\begin{equation}
\label{sixsix171s}
\lim\limits_{p\to\infty}\sum_{j_1=0}^{p}\sum_{j_2=0}^{p}\sum_{j_3=p+1}^{\infty}
\biggl(C_{j_3 j_1 j_3 j_2 j_2 j_1}+C_{j_3 j_1 j_3 j_2 j_1 j_2}\biggr)=0.
\end{equation}

From (\ref{sixsix170})--(\ref{sixsix171s}) we have
\begin{equation}
\label{sixsix300}
~~~~~~~~~~\lim\limits_{p\to\infty}
\sum_{j_1=p+1}^{\infty}\sum_{j_2=p+1}^{\infty}
\sum_{j_3=p+1}^{\infty}
\biggl(C_{j_3 j_1 j_3 j_2 j_2 j_1}+C_{j_3 j_1 j_3 j_2 j_1 j_2}\biggr)=0.
\end{equation}

\vspace{2mm}

Moreover (see (\ref{sixsix10})),
\begin{equation}
\label{sixsix301}
\lim\limits_{p\to\infty}
\sum_{j_1=p+1}^{\infty}\sum_{j_2=p+1}^{\infty}
\sum_{j_3=p+1}^{\infty}
C_{j_3 j_1 j_3 j_2 j_1 j_2}=0.
\end{equation}

\vspace{2mm}

Combining (\ref{sixsix300}) and (\ref{sixsix301}), we finally obtain
$$
\lim\limits_{p\to\infty}
\sum_{j_1=p+1}^{\infty}\sum_{j_2=p+1}^{\infty}
\sum_{j_3=p+1}^{\infty}
C_{j_3 j_1 j_3 j_2 j_2 j_1}=0.
$$

\vspace{2mm}
\noindent
The equality (\ref{sixsix13}) is proved.

Finally consider (\ref{sixsix15}).
Using the integration order replacement, we have
$$
C_{j_2 j_3 j_3 j_1 j_2 j_1}+C_{j_2 j_3 j_3 j_1 j_1 j_2}=
$$
$$
=
\int\limits_t^T 
\phi_{j_2}(t_6)
\int\limits_t^{t_6} 
\phi_{j_3}(t_5)
\int\limits_{t}^{t_5} 
\phi_{j_3}(t_4) 
\int\limits_t^{t_4}
\phi_{j_1}(t_3) 
\left(\int\limits_t^{t_3}
\phi_{j_2}(t_2)dt_2\right)
\left(\int\limits_t^{t_3}
\phi_{j_1}(t_1)
dt_1\right)\times
$$
$$
\times 
dt_3 dt_4 dt_5 dt_6=
$$
$$
=
\int\limits_t^T 
\phi_{j_2}(t_6)
\int\limits_t^{t_6} 
\phi_{j_3}(t_5)
\int\limits_{t}^{t_5} 
\phi_{j_1}(t_3) 
\left(\int\limits_t^{t_3}
\phi_{j_2}(t_2)dt_2\right)
\left(\int\limits_t^{t_3}
\phi_{j_1}(t_1)
dt_1\right)
\int\limits_{t_3}^{t_5}
\phi_{j_3}(t_4) 
dt_4
\times
$$
$$
\times 
dt_3 dt_5 dt_6=
$$
$$
=
\int\limits_t^T 
\phi_{j_2}(t_6)
\int\limits_t^{t_6} 
\phi_{j_3}(t_5)
\left(\int\limits_{t}^{t_5}
\phi_{j_3}(t_4) 
dt_4
\right)
\int\limits_{t}^{t_5} 
\phi_{j_1}(t_3) 
\left(\int\limits_t^{t_3}
\phi_{j_2}(t_2)dt_2\right)
\times
$$
$$
\times\left(\int\limits_t^{t_3}
\phi_{j_1}(t_1)
dt_1\right) 
dt_3 dt_5 dt_6-
$$
$$
-
\int\limits_t^T 
\phi_{j_2}(t_6)
\int\limits_t^{t_6} 
\phi_{j_3}(t_5)
\int\limits_{t}^{t_5} 
\phi_{j_1}(t_3) 
\left(\int\limits_t^{t_3}
\phi_{j_2}(t_2)dt_2\right)
\left(\int\limits_t^{t_3}
\phi_{j_1}(t_1)
dt_1\right)
\times
$$
$$
\times 
\left(\int\limits_t^{t_3}
\phi_{j_3}(t_4)
dt_4\right)
dt_3 dt_5 dt_6=
$$
$$
=
\int\limits_t^T 
\phi_{j_3}(t_5)
\left(\int\limits_{t}^{t_5}
\phi_{j_3}(t_4) 
dt_4
\right)
\int\limits_{t}^{t_5} 
\phi_{j_1}(t_3) 
\left(\int\limits_t^{t_3}
\phi_{j_2}(t_2)dt_2\right)
\times
$$
$$
\times\left(\int\limits_t^{t_3}
\phi_{j_1}(t_1)
dt_1\right) 
dt_3
\left(\int\limits_{t_5}^T
\phi_{j_2}(t_6)dt_6\right) dt_5-
$$
$$
-
\int\limits_t^T 
\phi_{j_3}(t_5)
\int\limits_{t}^{t_5} 
\phi_{j_1}(t_3) 
\left(\int\limits_t^{t_3}
\phi_{j_2}(t_2)dt_2\right)
\left(\int\limits_t^{t_3}
\phi_{j_1}(t_1)
dt_1\right)
\times
$$
\begin{equation}
\label{sixsix190}
\times
\left(\int\limits_t^{t_3}
\phi_{j_3}(t_4)
dt_4\right)dt_3
\left(
\int\limits_{t_5}^T
\phi_{j_2}(t_6)
dt_6\right)dt_5.
\end{equation}

\vspace{2mm}

Using (\ref{after80xx}) and (\ref{after500}), we get
$$
\sum_{j_1=p+1}^{\infty}\sum_{j_2=p+1}^{\infty}
\sum_{j_3=p+1}^{\infty}
\biggl(C_{j_2 j_3 j_3 j_1 j_2 j_1}+
C_{j_2 j_3 j_3 j_1 j_1 j_2}\biggr)=
$$
$$
=\frac{1}{2}
\sum_{j_1=p+1}^{\infty}\sum_{j_2=p+1}^{\infty}
\left(C_{j_2 j_3 j_3 j_1 j_2 j_1}\biggl|_{(j_3 j_3)\curvearrowright (\cdot)}\biggr. +
C_{j_2 j_3 j_3 j_1 j_1 j_2}\biggl|_{(j_3 j_3)\curvearrowright (\cdot)}\biggr.\right)-
$$
$$
-\sum_{j_3=0}^{p}\sum_{j_1=p+1}^{\infty}
\sum_{j_2=p+1}^{\infty}
\biggl(C_{j_2 j_3 j_3 j_1 j_2 j_1}+
C_{j_2 j_3 j_3 j_1 j_1 j_2}\biggr)=
$$
$$
=\frac{1}{2}
\sum_{j_1=p+1}^{\infty}\sum_{j_2=p+1}^{\infty}
\left(C_{j_2 j_3 j_3 j_1 j_2 j_1}\biggl|_{(j_3 j_3)\curvearrowright (\cdot)}\biggr. +
C_{j_2 j_3 j_3 j_1 j_1 j_2}\biggl|_{(j_3 j_3)\curvearrowright (\cdot)}\biggr.\right)+
$$
$$
+\sum_{j_1=0}^{p}\sum_{j_3=0}^{p}
\sum_{j_2=p+1}^{\infty}
\biggl(C_{j_2 j_3 j_3 j_1 j_2 j_1}+
C_{j_2 j_3 j_3 j_1 j_1 j_2}\biggr)-
$$
\begin{equation}
\label{sixsix191}
-\frac{1}{2}
\sum_{j_3=0}^{p}\sum_{j_2=p+1}^{\infty}
C_{j_2 j_3 j_3 j_1 j_1 j_2}\biggl|_{(j_1 j_1)\curvearrowright (\cdot)}\biggr..
\end{equation}

\vspace{2mm}

The equalities
\begin{equation}
\label{sixsix192}
~~~~\lim\limits_{p\to\infty}\frac{1}{2}
\sum_{j_1=p+1}^{\infty}\sum_{j_2=p+1}^{\infty}
\left(C_{j_2 j_3 j_3 j_1 j_2 j_1}\biggl|_{(j_3 j_3)\curvearrowright (\cdot)}\biggr. +
C_{j_2 j_3 j_3 j_1 j_1 j_2}\biggl|_{(j_3 j_3)\curvearrowright (\cdot)}\biggr.\right)=0,
\end{equation}
$$
\lim\limits_{p\to\infty}\frac{1}{2}
\sum_{j_3=0}^{p}\sum_{j_2=p+1}^{\infty}
C_{j_2 j_3 j_3 j_1 j_1 j_2}\biggl|_{(j_1 j_1)\curvearrowright (\cdot)}\biggr.
=
$$
$$
=\lim\limits_{p\to\infty}\frac{1}{4}
\sum_{j_2=p+1}^{\infty}
C_{j_2 j_3 j_3 j_1 j_1 j_2}\biggl|_{(j_1 j_1)\curvearrowright (\cdot)
(j_3 j_3)\curvearrowright (\cdot)}\biggr.-
$$
\begin{equation}
\label{sixsix193}
-
\lim\limits_{p\to\infty}\frac{1}{2}
\sum_{j_3=p+1}^{\infty}\sum_{j_2=p+1}^{\infty}
C_{j_2 j_3 j_3 j_1 j_1 j_2}\biggl|_{(j_1 j_1)\curvearrowright (\cdot)}\biggr.
=0
\end{equation}

\vspace{2mm}
\noindent
follows from the 
equalities (\ref{after2508}), (\ref{after2509}),
where we used the same technique as in (\ref{sept12}).
When proving (\ref{sixsix193}), we also applied (\ref{after500})
and (\ref{tupo15}).

By analogy with the proof of (\ref{sixsix7}) and applying (\ref{sixsix190}), we obtain
\begin{equation}
\label{sixsix194}
\lim\limits_{p\to\infty}
\sum_{j_1=0}^{p}\sum_{j_3=0}^{p}
\sum_{j_2=p+1}^{\infty}
\biggl(C_{j_2 j_3 j_3 j_1 j_2 j_1}+
C_{j_2 j_3 j_3 j_1 j_1 j_2}\biggr)=0.
\end{equation}

\vspace{2mm}

From (\ref{sixsix191})--(\ref{sixsix194}) we have
\begin{equation}
\label{sixsix194a}
~~~~~~~~~\lim\limits_{p\to\infty}
\sum_{j_1=p+1}^{\infty}\sum_{j_2=p+1}^{\infty}
\sum_{j_3=p+1}^{\infty}
\biggl(C_{j_2 j_3 j_3 j_1 j_2 j_1}+
C_{j_2 j_3 j_3 j_1 j_1 j_2}\biggr)=0.
\end{equation}

\vspace{2mm}

Furthermore (see (\ref{sixsix14})),
\begin{equation}
\label{sixsix195}
\lim\limits_{p\to\infty}\sum_{j_1=p+1}^{\infty}\sum_{j_2=p+1}^{\infty}
\sum_{j_3=p+1}^{\infty}
C_{j_2 j_3 j_3 j_1 j_1 j_2}=0.
\end{equation}

\vspace{2mm}

Combining (\ref{sixsix194a}) and (\ref{sixsix195}), we finally obtain
$$
\lim\limits_{p\to\infty}\sum_{j_1=p+1}^{\infty}\sum_{j_2=p+1}^{\infty}
\sum_{j_3=p+1}^{\infty}
C_{j_2 j_3 j_3 j_1 j_2 j_1}=0.
$$

\vspace{2mm}
\noindent
The equality (\ref{sixsix15}) is proved. 
Theorem 2.36 is proved.

\section{Estimates for the Mean-Square Approximation Error of Iterated
Stra\-to\-no\-vich Stochastic Integrals of Multiplicity $k$
in Theorems 2.30, 2.31}

In this section, we estimate the mean-square approximation error
for iterated Stratonovich stochastic integrals of multiplicity
$k$ ($k\in{\bf N}$) in Theorems~2.30, 2.31.

{\bf Theorem~2.37}\ \cite{arxiv-5}, 
\cite{arxiv-10}, \cite{arxiv-11}, \cite{new-art-1xxy}.\
{\it Suppose that every $\psi_l(\tau)$ $(l=1,\ldots,k)$
is a continuously differentiable nonrandom function
at the interval $[t, T].$ Furthermore$,$ let
$\{\phi_j(x)\}_{j=0}^{\infty}$ is a complete orthonormal system of 
Legendre polynomials or trigonometric functions in the space $L_2([t, T]).$
Then the following estimates 

\vspace{-3mm}
$$
{\sf M}\left\{\left(
J^{*}[\psi^{(k)}]_{T,t}^{(i_1\ldots i_k)}-
\sum\limits_{j_1,\ldots,j_k=0}^{p}
C_{j_k \ldots j_1}\prod\limits_{l=1}^k \zeta_{j_l}^{(i_l)}
\right)^2\right\}\le
$$

\vspace{-2mm}
\begin{equation}
\label{after3407}
\le K_1 \left(\frac{1}{p} +\sum_{r=1}^{\left[k/2\right]}
\sum_{\stackrel{(\{\{g_1, g_2\}, \ldots, 
\{g_{2r-1}, g_{2r}\}\}, \{q_1, \ldots, q_{k-2r}\})}
{{}_{\{g_1, g_2, \ldots, 
g_{2r-1}, g_{2r}, q_1, \ldots, q_{k-2r}\}=\{1, 2, \ldots, k\}}}}
{\sf M}\left\{\left(R_{T,t}^{(p)r, g_1,g_2,\ldots,g_{2r-1}, g_{2r}}\right)^2\right\}\right),
\end{equation}

\vspace{2mm}
$$
{\sf M}\left\{\left(
J^{*}[\psi^{(k)}]_{s,t}^{(i_1\ldots i_k)}-
\sum\limits_{j_1,\ldots,j_k=0}^{p}
C_{j_k \ldots j_1}(s)\prod\limits_{l=1}^k \zeta_{j_l}^{(i_l)}
\right)^2\right\}\le 
$$

\vspace{-2mm}
\begin{equation}
\label{after3408}
\le K_2(s) \left(\frac{1}{p} +\sum_{r=1}^{\left[k/2\right]}
\sum_{\stackrel{(\{\{g_1, g_2\}, \ldots, 
\{g_{2r-1}, g_{2r}\}\}, \{q_1, \ldots, q_{k-2r}\})}
{{}_{\{g_1, g_2, \ldots, 
g_{2r-1}, g_{2r}, q_1, \ldots, q_{k-2r}\}=\{1, 2, \ldots, k\}}}}
{\sf M}\left\{\left(R_{s,t}^{(p)r, g_1,g_2,\ldots,g_{2r-1}, g_{2r}}\right)^2\right\}\right)
\end{equation}

\vspace{4mm}
\noindent
hold, where $s\in(t,T]$ $(s$ is fixed$),$ $i_1,\ldots,i_k=1,\ldots,m,$

$$
R_{s,t}^{(p)r, g_1,g_2,\ldots,g_{2r-1}, g_{2r}}=
R_{T,t}^{(p)r, g_1,g_2,\ldots,g_{2r-1}, g_{2r}}\biggl|_{T=s}\biggr.,
$$

\vspace{1mm}
\noindent
$R_{T,t}^{(p)r, g_1,g_2,\ldots,g_{2r-1}, g_{2r}}$ is defined by {\rm (\ref{afterr1}),}
$J^{*}[\psi^{(k)}]_{T,t}^{(i_1\ldots i_k)}$ and $J^{*}[\psi^{(k)}]_{s,t}^{(i_1\ldots i_k)}$
are iterated Stratonovich stochastic integrals {\rm (\ref{afterstr})} and {\rm (\ref{afterstr1}),}
$C_{j_k \ldots j_1}$ and $C_{j_k \ldots j_1}(s)$
are Fourier coefficients {\rm (\ref{after3000})} and {\rm (\ref{after1300}),} 
constants $K_1$ and $K_2(s)$ are independent of $p;$
another notations are the same as in Theorems {\rm 1.1,}
{\rm 2.30,} {\rm 2.31.}
}

{\bf Proof.}\ Note that Conditions {\rm 1} and {\rm 2} of Theorems
{\rm 2.30, 2.31} are satisfied under the conditions of Theorem~2.37
(see Remark~2.4). From the proof of Theorem~2.30 it follows that
the expression (\ref{afteru11}) ($i_1,\ldots,i_k=1,\ldots,m$)
before passing to the limit 
\hbox{\vtop{\offinterlineskip\halign{
\hfil#\hfil\cr
{\rm l.i.m.}\cr
$\stackrel{}{{}_{p\to \infty}}$\cr
}} } has the form

$$
\sum_{j_1,\ldots,j_k=0}^{p}
C_{j_k\ldots j_1}
\prod\limits_{l=1}^k \zeta_{j_l}^{(i_l)}=
J[\psi^{(k)}]_{T,t}^{(i_1\ldots i_k)p}+
$$

$$
+
\sum_{r=1}^{\left[k/2\right]}\Biggl(\frac{1}{2^r}
\sum_{(s_r,\ldots,s_1)\in {\rm A}_{k,r}}
I[\psi^{(k)}]_{T,t}^{(i_1\ldots i_{s_1-1}i_{s_1+2} \ldots i_{s_r-1}i_{s_r+2}
\ldots i_k)p}
+\Biggr.
$$

\vspace{-1mm}
\begin{equation}
\label{after3400}
\Biggl.
+\sum_{\stackrel{(\{\{g_1, g_2\}, \ldots, 
\{g_{2r-1}, g_{2r}\}\}, \{q_1, \ldots, q_{k-2r}\})}
{{}_{\{g_1, g_2, \ldots, 
g_{2r-1}, g_{2r}, q_1, \ldots, q_{k-2r}\}=\{1, 2, \ldots, k\}}}}
R_{T,t}^{(p)r, g_1,g_2,\ldots,g_{2r-1}, g_{2r}}
\Biggr)
\end{equation}

\vspace{3mm}
\noindent
w.~p.~1, where 
$J[\psi^{(k)}]_{T,t}^{(i_1\ldots i_k)p}$
is the approximation (\ref{kkohh}) of the iterated It\^{o}
stochastic integral (\ref{afterito}),
$I[\psi^{(k)}]_{T,t}^{(i_1\ldots i_{s_1-1}i_{s_1+2} \ldots i_{s_r-1}i_{s_r+2} \ldots i_k)p}$
is the approximation obtained using (\ref{kkohh}) for the 
iterated It\^{o}
stochastic integral 
$J[\psi^{(k)}]_{T,t}^{s_r, \ldots, s_1}$
(see (\ref{afterito1})).

Using (\ref{after3400}) and Theorem~2.12, we have

\vspace{-2mm}
$$
\sum_{j_1,\ldots,j_k=0}^{p}
C_{j_k\ldots j_1}
\prod\limits_{l=1}^k \zeta_{j_l}^{(i_l)}=
J[\psi^{(k)}]_{T,t}^{(i_1\ldots i_k)}+
$$

\vspace{-1mm}
$$
+
\sum_{r=1}^{\left[k/2\right]}\frac{1}{2^r}
\sum_{(s_r,\ldots,s_1)\in {\rm A}_{k,r}}
I[\psi^{(k)}]_{T,t}^{(i_1\ldots i_{s_1-1}i_{s_1+2} \ldots i_{s_r-1}i_{s_r+2}\ldots i_k)}
+
$$

\vspace{1mm}
$$
+\biggl(J[\psi^{(k)}]_{T,t}^{(i_1\ldots i_k)p}-J[\psi^{(k)}]_{T,t}^{(i_1\ldots i_k)}\biggr)+
$$

\vspace{-1mm}
$$
+
\sum_{r=1}^{\left[k/2\right]}\sum_{(s_r,\ldots,s_1)\in {\rm A}_{k,r}}
\frac{1}{2^r}
\Biggl(
I[\psi^{(k)}]_{T,t}^{(i_1\ldots i_{s_1-1}i_{s_1+2} \ldots i_{s_r-1}i_{s_r+2}\ldots i_k)p}-
\Biggr.
$$
$$
\Biggl.-
I[\psi^{(k)}]_{T,t}^{(i_1\ldots i_{s_1-1}i_{s_1+2} \ldots i_{s_r-1}i_{s_r+2}\ldots i_k)}
\Biggr)+
$$

$$
+\sum_{r=1}^{\left[k/2\right]}\sum_{\stackrel{(\{\{g_1, g_2\}, \ldots, 
\{g_{2r-1}, g_{2r}\}\}, \{q_1, \ldots, q_{k-2r}\})}
{{}_{\{g_1, g_2, \ldots, 
g_{2r-1}, g_{2r}, q_1, \ldots, q_{k-2r}\}=\{1, 2, \ldots, k\}}}}
R_{T,t}^{(p)r, g_1,g_2,\ldots,g_{2r-1}, g_{2r}}
=
$$

\vspace{4mm}
$$
=J^{*}[\psi^{(k)}]_{T,t}^{(i_1\ldots i_k)}+
\biggl(J[\psi^{(k)}]_{T,t}^{(i_1\ldots i_k)p}-J[\psi^{(k)}]_{T,t}^{(i_1\ldots i_k)}\biggr)+
$$

$$
+
\sum_{r=1}^{\left[k/2\right]}\sum_{(s_r,\ldots,s_1)\in {\rm A}_{k,r}}
\frac{1}{2^r}
\Biggl(
I[\psi^{(k)}]_{T,t}^{(i_1\ldots i_{s_1-1}i_{s_1+2} \ldots i_{s_r-1}i_{s_r+2}\ldots i_k)p}-
\Biggr.
$$
$$
\Biggl.-
I[\psi^{(k)}]_{T,t}^{(i_1\ldots i_{s_1-1}i_{s_1+2} \ldots i_{s_r-1}i_{s_r+2}\ldots i_k)}
\Biggr)+
$$

\begin{equation}
\label{after3401}
+\sum_{r=1}^{\left[k/2\right]}\sum_{\stackrel{(\{\{g_1, g_2\}, \ldots, 
\{g_{2r-1}, g_{2r}\}\}, \{q_1, \ldots, q_{k-2r}\})}
{{}_{\{g_1, g_2, \ldots, 
g_{2r-1}, g_{2r}, q_1, \ldots, q_{k-2r}\}=\{1, 2, \ldots, k\}}}}
R_{T,t}^{(p)r, g_1,g_2,\ldots,g_{2r-1}, g_{2r}}
\end{equation}

\vspace{4mm}
\noindent
w.~p.~1, where we denote $J[\psi^{(k)}]_{T,t}^{s_r, \ldots, s_1}$ as 
$I[\psi^{(k)}]_{T,t}^{(i_1\ldots i_{s_1-1}i_{s_1+2} \ldots i_{s_r-1}i_{s_r+2}\ldots i_k)}$.

Applying  (\ref{zsel1xyzuv}) (see Remark~1.7), we obtain the following
estimates

\vspace{-2mm}
\begin{equation}
\label{after3404}
{\sf M}\left\{\biggl(
J[\psi^{(k)}]_{T,t}^{(i_1\ldots i_k)p}-J[\psi^{(k)}]_{T,t}^{(i_1\ldots i_k)}\biggr)^2\right\}
\le \frac{C}{p},
\end{equation}
$$
{\sf M}\left\{\Biggl(
I[\psi^{(k)}]_{T,t}^{(i_1\ldots i_{s_1-1}i_{s_1+2} \ldots i_{s_r-1}i_{s_r+2}\ldots i_k)p}
-I[\psi^{(k)}]_{T,t}^{(i_1\ldots i_{s_1-1}i_{s_1+2} \ldots i_{s_r-1}i_{s_r+2}\ldots i_k)}
\Biggr)^2\right\}\le
$$
\begin{equation}
\label{after3405}
\le \frac{C}{p},
\end{equation}

\vspace{2mm}
\noindent
where constant $C$ does not depend on $p.$

From (\ref{after3401})--(\ref{after3405}) and
the elementary inequality
$$
\left(a_1+a_2+\ldots+a_n\right)^2 \le
n\left(a_1^2+a_2^2+\ldots+a_n^2\right),\ \ \ n\in {\bf N},
$$

\noindent
we obtain (\ref{after3407}).
The estimate 
(\ref{after3408}) is obtained similarly to the estimate 
(\ref{after3407})
using Theorems~1.11, 2.31 and (\ref{road1888}) (see Remark~1.12).
Theorem~2.37 is proved.

\section{Rate of the Mean-Square Convergence of Expansions of Iterated
Stra\-to\-no\-vich Stochastic Integrals of Multiplicities 3--5
in Theorems 2.33--2.35}

In this section, we consider the rate of convergence of 
approximations of iterated Stratonovich stochastic integrals in Theorems~2.33--2.35.
It is easy to see that in Theorems~2.33--2.35
the second term 
in parentheses
on the right-hand side of (\ref{after3407}) is estimated for $k=3, 4, 5$.
Combining these results with Theorem~2.37, we obtain the following theorems.

{\bf Theorem 2.38}\ \cite{arxiv-5}, 
\cite{arxiv-10}, \cite{arxiv-11}, \cite{new-art-1xxy}.\ {\it Suppose 
that $\{\phi_j(x)\}_{j=0}^{\infty}$ is a complete orthonormal system of 
Legendre polynomials or trigonometric functions in the space $L_2([t, T]).$
Furthermore$,$ let $\psi_1(\tau), \psi_2(\tau), \psi_3(\tau)$ are continuously dif\-ferentiable 
nonrandom functions on $[t, T].$ 
Then$,$ for the 
iterated Stratonovich stochastic integral of third multiplicity
$$
J^{*}[\psi^{(3)}]_{T,t}={\int\limits_t^{*}}^T\psi_3(t_3)
{\int\limits_t^{*}}^{t_3}\psi_2(t_2)
{\int\limits_t^{*}}^{t_2}\psi_1(t_1)
d{\bf w}_{t_1}^{(i_1)}
d{\bf w}_{t_2}^{(i_2)}d{\bf w}_{t_3}^{(i_3)}
$$

\noindent
the following 
estimate
$$
{\sf M}\left\{\left(
J^{*}[\psi^{(3)}]_{T,t}-
\sum\limits_{j_1, j_2, j_3=0}^{p}
C_{j_3 j_2 j_1}\zeta_{j_1}^{(i_1)}\zeta_{j_2}^{(i_2)}\zeta_{j_3}^{(i_3)}\right)^2\right\}
\le \frac{C}{p}
$$
is fulfilled, where $i_1, i_2, i_3=1,\ldots,m,$ 
constant $C$ is independent of $p,$
$$
C_{j_3 j_2 j_1}=\int\limits_t^T\psi_3(t_3)\phi_{j_3}(t_3)
\int\limits_t^{t_3}\psi_2(t_2)\phi_{j_2}(t_2)
\int\limits_t^{t_2}\psi_1(t_1)\phi_{j_1}(t_1)dt_1dt_2dt_3
$$
and
$$
\zeta_{j}^{(i)}=
\int\limits_t^T \phi_{j}(s) d{\bf w}_s^{(i)}
$$ 
are independent standard Gaussian random variables for various 
$i$ or $j.$}

{\bf Theorem 2.39}\ \cite{arxiv-5}, 
\cite{arxiv-10}, \cite{arxiv-11}, \cite{new-art-1xxy}.\ {\it Let
$\{\phi_j(x)\}_{j=0}^{\infty}$ be a complete orthonormal system of 
Legendre polynomials or trigonometric functions in the space $L_2([t, T]).$
Furthermore$,$ let $\psi_1(\tau), \ldots, \psi_4(\tau)$ be continuously dif\-ferentiable 
nonrandom functions on $[t, T].$ 
Then$,$ for the 
iterated Stratonovich stochastic integral of fourth multiplicity
$$
J^{*}[\psi^{(4)}]_{T,t}={\int\limits_t^{*}}^T\psi_4(t_4)
{\int\limits_t^{*}}^{t_4}\psi_3(t_3)
{\int\limits_t^{*}}^{t_3}\psi_2(t_2)
{\int\limits_t^{*}}^{t_2}\psi_1(t_1)
d{\bf w}_{t_1}^{(i_1)}
d{\bf w}_{t_2}^{(i_2)}d{\bf w}_{t_3}^{(i_3)}d{\bf w}_{t_4}^{(i_4)}
$$
the following 
estimate
$$
{\sf M}\left\{\left(
J^{*}[\psi^{(4)}]_{T,t}-
\sum\limits_{j_1, j_2, j_3, j_4=0}^{p}
C_{j_4 j_3 j_2 j_1}\zeta_{j_1}^{(i_1)}\zeta_{j_2}^{(i_2)}\zeta_{j_3}^{(i_3)}
\zeta_{j_4}^{(i_4)}
\right)^2\right\}
\le \frac{C}{p^{1-\varepsilon}}
$$
holds, where $i_1, i_2, i_3, i_4=1,\ldots,m,$ 
constant $C$ does not depend on $p,$
$\varepsilon$ is an arbitrary
small positive real number 
for the case of complete orthonormal system of 
Legendre polynomials in the space $L_2([t, T])$
and $\varepsilon=0$ for the case of
complete orthonormal system of 
trigonometric functions in the space $L_2([t, T]),$
$$
C_{j_4 j_3 j_2 j_1}=
\int\limits_t^T\psi_4(t_4)\phi_{j_4}(t_4)
\int\limits_t^{t_4}\psi_3(t_3)\phi_{j_3}(t_3)
\int\limits_t^{t_3}\psi_2(t_2)\phi_{j_2}(t_2)
\int\limits_t^{t_2}\psi_1(t_1)\phi_{j_1}(t_1)dt_1\times
$$
$$
\times dt_2dt_3dt_4;
$$

\vspace{2mm}
\noindent
another notations are the same as in Theorem~{\rm 2.38}.}

Note that Theorem~2.26 is an analog of Theorem~2.39. At that
$\varepsilon=0,$ 
$\psi_1(\tau),\ldots,\psi_4(\tau)\equiv 1,$
and $i_1,\ldots,i_k=0,1,\ldots,m$  in Theorem~2.26.

{\bf Theorem 2.40}\ \cite{arxiv-5}, 
\cite{arxiv-10}, \cite{arxiv-11}, \cite{new-art-1xxy}.\ {\it Assume 
that $\{\phi_j(x)\}_{j=0}^{\infty}$ is a complete orthonormal system of 
Legendre polynomials or trigonometric functions in the space $L_2([t, T])$
and $\psi_1(\tau), \ldots, \psi_5(\tau)$ are continuously dif\-ferentiable 
nonrandom functions on $[t, T].$ 
Then$,$ for the 
iterated Stratonovich stochastic integral of fifth multiplicity
$$
J^{*}[\psi^{(5)}]_{T,t}={\int\limits_t^{*}}^T\psi_5(t_5)
\ldots
{\int\limits_t^{*}}^{t_2}\psi_1(t_1)
d{\bf w}_{t_1}^{(i_1)}
\ldots d{\bf w}_{t_5}^{(i_5)}
$$
the following 
estimate
$$
{\sf M}\left\{\left(
J^{*}[\psi^{(5)}]_{T,t}-
\sum\limits_{j_1, \ldots, j_5=0}^{p}
C_{j_5 \ldots j_1}\zeta_{j_1}^{(i_1)}\ldots
\zeta_{j_5}^{(i_5)}
\right)^2\right\}
\le \frac{C}{p^{1-\varepsilon}}
$$
is valid, where $i_1, \ldots, i_5=1,\ldots,m,$ 
constant $C$ is independent of $p,$
$\varepsilon$ is an arbitrary
small positive real number 
for the case of complete orthonormal system of 
Legendre polynomials in the space $L_2([t, T])$
and $\varepsilon=0$ for the case of
complete orthonormal system of 
trigonometric functions in the space $L_2([t, T]),$
$$
C_{j_5 \ldots j_1}=
\int\limits_t^T\psi_5(t_5)\phi_{j_5}(t_5)\ldots
\int\limits_t^{t_2}\psi_1(t_1)\phi_{j_1}(t_1)dt_1\ldots dt_5;
$$
\noindent
another notations are the same as in Theorem~{\rm 2.38, 2.39}.}

\section{Generalization of Theorems~2.4--2.8. 
The Case $p_1,$ $p_2,$ $p_3\to \infty$ and Continuously
Differetiable Weight Functions (The Cases of Legendre Polynomials and 
Trigonometric Functions). Proof of Hypothesis~2.3 for the Case $k=3$}

This section is devoted to the following theorem.

{\bf Theorem~2.41}\ \cite{arxiv-5}, \cite{arxiv-10}, \cite{arxiv-11}.\ {\it Suppose that 
$\{\phi_j(x)\}_{j=0}^{\infty}$ is a complete orthonormal system of 
Legendre polynomials or trigonometric functions in the space $L_2([t, T]).$
Furthermore, let $\psi_1(\tau), \psi_2(\tau), \psi_3(\tau)$ are continuously dif\-ferentiable 
nonrandom functions on $[t, T]$.
Then$,$ for the 
iterated Stratonovich stochastic integral of third multiplicity
$$
J^{*}[\psi^{(3)}]_{T,t}^{(i_1 i_2 i_3)}={\int\limits_t^{*}}^T\psi_3(t_3)
{\int\limits_t^{*}}^{t_3}\psi_2(t_2)
{\int\limits_t^{*}}^{t_2}\psi_1(t_1)
d{\bf w}_{t_1}^{(i_1)}
d{\bf w}_{t_2}^{(i_2)}d{\bf w}_{t_3}^{(i_3)}
$$
the following 
expansion 
\begin{equation}
\label{20231}
~~~~~J^{*}[\psi^{(3)}]_{T,t}^{(i_1 i_2 i_3)}
=\hbox{\vtop{\offinterlineskip\halign{
\hfil#\hfil\cr
{\rm l.i.m.}\cr
$\stackrel{}{{}_{p_1,p_2,p_3\to \infty}}$\cr
}} }
\sum\limits_{j_1=0}^{p_1}\sum\limits_{j_2=0}^{p_2}\sum\limits_{j_3=0}^{p_3}
C_{j_3 j_2 j_1}\zeta_{j_1}^{(i_1)}\zeta_{j_2}^{(i_2)}\zeta_{j_3}^{(i_3)}
\end{equation}
that converges in the mean-square sense is valid, where
$i_1, i_2, i_3=0, 1,\ldots,m,$
$$
C_{j_3 j_2 j_1}=\int\limits_t^T\psi_3(t_3)\phi_{j_3}(t_3)
\int\limits_t^{t_3}\psi_2(t_2)\phi_{j_2}(t_2)
\int\limits_t^{t_2}\psi_1(t_1)\phi_{j_1}(t_1)dt_1dt_2dt_3
$$
and
$$
\zeta_{j}^{(i)}=
\int\limits_t^T \phi_{j}(s) d{\bf w}_s^{(i)}
$$ 
are independent standard Gaussian random variables for various 
$i$ or $j$
{\rm (}in the case when $i\ne 0${\rm ),}
${\bf w}_{\tau}^{(i)}$ 
$(i=1,\ldots,m)$ are independent standard Wiener processes$,$
${\bf w}_{\tau}^{(0)}=\tau.$}

\vspace{2mm}

{\bf Proof.} Let us consider the case of
Legendre polynomials (the trigonometric case is simpler
and can be considered similarly). Applying (\ref{after33}), we obtain
$$
\sum_{j_1=0}^{p_1}\sum_{j_2=0}^{p_2}\sum_{j_3=0}^{p_3}
C_{j_3j_2j_1}
\zeta_{j_1}^{(i_1)}\zeta_{j_2}^{(i_2)}\zeta_{j_3}^{(i_3)}=
J'[K_{p_1p_2p_3}]_{T,t}^{(i_1 i_2 i_3)}+
$$
$$
+{\bf 1}_{\{i_1=i_2\ne 0\}}
\sum_{j_3=0}^{p_3}\sum_{j_1=0}^{\min\{p_1,p_2\}}
C_{j_3j_1j_1}J'[\phi_{j_3}]^{(i_3)}_{T,t}+
$$
$$
+{\bf 1}_{\{i_2=i_3\ne 0\}}
\sum_{j_1=0}^{p_1}\sum_{j_3=0}^{\min\{p_2,p_3\}}
C_{j_3j_3j_1}J'[\phi_{j_1}]^{(i_1)}_{T,t}+
$$
\begin{equation}
\label{20232}
+{\bf 1}_{\{i_1=i_3\ne 0\}}
\sum_{j_2=0}^{p_2}\sum_{j_1=0}^{\min\{p_1,p_3\}}
C_{j_1j_2j_1}J'[\phi_{j_2}]^{(i_2)}_{T,t}
\end{equation}

\vspace{2mm}
\noindent
w.~p.~1, where notations are the same as in (\ref{after33}).

Using (\ref{uyes3}), Theorem~1.1 (see (\ref{drdr1})),
Theorem~2.12 (see (\ref{30.4})) as well as (\ref{after609}) (see the derivation
of (\ref{after609})) and (\ref{after500}), we get
$$
J^{*}[\psi^{(3)}]_{T,t}^{(i_1 i_2 i_3)}=J[\psi^{(3)}]_{T,t}^{(i_1 i_2 i_3)}
+ \frac{1}{2}{\bf 1}_{\{i_1=i_2\ne 0\}}
\int\limits_t^T\psi_3(t_3)
\int\limits_t^{t_3}\psi_2(t_2)\psi_1(t_2)dt_2
d{\bf w}_{t_3}^{(i_3)}+
$$
$$
+ \frac{1}{2}{\bf 1}_{\{i_2=i_3\ne 0\}}
\int\limits_t^T\psi_3(t_3)\psi_2(t_3)
\int\limits_t^{t_3}\psi_1(t_1)d{\bf w}_{t_1}^{(i_1)}
dt_3=
$$

\vspace{2.5mm}
$$
=J[\psi^{(3)}]_{T,t}^{(i_1 i_2 i_3)}
+ \frac{1}{2}J[\psi^{(3)}]_{T,t}^1+\frac{1}{2}J[\psi^{(3)}]_{T,t}^2
=
$$

\vspace{2.5mm}
$$
=
\hbox{\vtop{\offinterlineskip\halign{
\hfil#\hfil\cr
{\rm l.i.m.}\cr
$\stackrel{}{{}_{p_1,p_2,p_3\to \infty}}$\cr
}} }J'[K_{p_1p_2p_3}]_{T,t}^{(i_1 i_2 i_3)}+
$$
$$
+{\bf 1}_{\{i_1=i_2\ne 0\}}
\hbox{\vtop{\offinterlineskip\halign{
\hfil#\hfil\cr
{\rm l.i.m.}\cr
$\stackrel{}{{}_{p_3\to \infty}}$\cr
}} }\frac{1}{2}\sum_{j_3=0}^{p_3} 
C_{j_3 j_2 j_1}\biggl|_{(j_{2} j_{1})\curvearrowright (\cdot),
j_{1}=j_{2}}\biggr.
J'[\phi_{j_3}]^{(i_3)}_{T,t}+
$$
$$
+{\bf 1}_{\{i_2=i_3\ne 0\}}
\hbox{\vtop{\offinterlineskip\halign{
\hfil#\hfil\cr
{\rm l.i.m.}\cr
$\stackrel{}{{}_{p_1\to \infty}}$\cr
}} }
\frac{1}{2}\sum_{j_1=0}^{p_1} 
C_{j_3 j_2 j_1}\biggl|_{(j_{3} j_{2})\curvearrowright (\cdot),
j_{2}=j_{3}}\biggr.
J'[\phi_{j_1}]^{(i_1)}_{T,t}=
$$

\vspace{2.5mm}
$$
=
\hbox{\vtop{\offinterlineskip\halign{
\hfil#\hfil\cr
{\rm l.i.m.}\cr
$\stackrel{}{{}_{p_1,p_2,p_3\to \infty}}$\cr
}} }J'[K_{p_1p_2p_3}]_{T,t}^{(i_1 i_2 i_3)}+
$$

\vspace{1mm}
$$
+{\bf 1}_{\{i_1=i_2\ne 0\}}
\hbox{\vtop{\offinterlineskip\halign{
\hfil#\hfil\cr
{\rm l.i.m.}\cr
$\stackrel{}{{}_{p_3\to \infty}}$\cr
}} }\sum_{j_3=0}^{p_3}\sum_{j_1=0}^{\infty}  
C_{j_3 j_1 j_1}  
J'[\phi_{j_3}]^{(i_3)}_{T,t}+
$$

\begin{equation}
\label{20233}
+{\bf 1}_{\{i_2=i_3\ne 0\}}
\hbox{\vtop{\offinterlineskip\halign{
\hfil#\hfil\cr
{\rm l.i.m.}\cr
$\stackrel{}{{}_{p_1\to \infty}}$\cr
}} }
\sum_{j_1=0}^{p_1} \sum_{j_3=0}^{\infty}
C_{j_3 j_3 j_1}
J'[\phi_{j_1}]^{(i_1)}_{T,t}
\end{equation}

\vspace{1mm}
\noindent
w.~p.~1.

Using (\ref{20232}), (\ref{20233}) and the elementary inequality 
$$
(a+b+c+d)^2\le 4\left(a^2+b^2+c^2+d^2\right),
$$
we obtain
$$
{\sf M}\left\{\left(J^{*}[\psi^{(3)}]_{T,t}^{(i_1 i_2 i_3)}-
\sum\limits_{j_1=0}^{p_1}\sum\limits_{j_2=0}^{p_2}\sum\limits_{j_3=0}^{p_3}
C_{j_3 j_2 j_1}\zeta_{j_1}^{(i_1)}\zeta_{j_2}^{(i_2)}\zeta_{j_3}^{(i_3)}\right)^2\right\}\le
$$

\vspace{2mm}
$$
\le 4{\sf M}\left\{\biggl(J[\psi^{(3)}]_{T,t}^{(i_1 i_2 i_3)}-
J'[K_{p_1p_2p_3}]_{T,t}^{(i_1 i_2 i_3)}
\biggr)^2\right\}+
$$

\vspace{1mm}
$$
+4\cdot {\bf 1}_{\{i_1=i_2\ne 0\}}\times
$$
$$
\times
{\sf M}\left\{\left(
\hbox{\vtop{\offinterlineskip\halign{
\hfil#\hfil\cr
{\rm l.i.m.}\cr
$\stackrel{}{{}_{p_3\to \infty}}$\cr
}} }\sum_{j_3=0}^{p_3}\sum_{j_1=0}^{\infty}  
C_{j_3 j_1 j_1}  
J'[\phi_{j_3}]^{(i_3)}_{T,t}-
\sum_{j_3=0}^{p_3}\sum_{j_1=0}^{\min\{p_1,p_2\}}
C_{j_3j_1j_1}J'[\phi_{j_3}]^{(i_3)}_{T,t}\right)^2\right\}+
$$

\vspace{1mm}
$$
+4\cdot {\bf 1}_{\{i_2=i_3\ne 0\}}\times
$$
$$
\times
{\sf M}\left\{\left(
\hbox{\vtop{\offinterlineskip\halign{
\hfil#\hfil\cr
{\rm l.i.m.}\cr
$\stackrel{}{{}_{p_1\to \infty}}$\cr
}} }\sum_{j_1=0}^{p_1}\sum_{j_3=0}^{\infty}  
C_{j_3 j_3 j_1}  
J'[\phi_{j_1}]^{(i_1)}_{T,t}-
\sum_{j_1=0}^{p_1}\sum_{j_3=0}^{\min\{p_2,p_3\}}
C_{j_3j_3j_1}J'[\phi_{j_1}]^{(i_1)}_{T,t}\right)^2\right\}+
$$
$$
+4\cdot {\bf 1}_{\{i_1=i_3\ne 0\}}
{\sf M}\left\{\left(
\sum_{j_2=0}^{p_2}\sum_{j_1=0}^{\min\{p_1,p_3\}}
C_{j_1j_2j_1}J'[\phi_{j_2}]^{(i_2)}_{T,t}\right)^2\right\}=
$$

\vspace{-4mm}
\begin{equation}
\label{20234}
~~~~~~~= 4 A_{p_1 p_2 p_3} + 4\cdot {\bf 1}_{\{i_1=i_2\ne 0\}}B_{p_1 p_2 p_3}+
4\cdot {\bf 1}_{\{i_2=i_3\ne 0\}}C_{p_1 p_2 p_3}+
4\cdot {\bf 1}_{\{i_1=i_3\ne 0\}}D_{p_1 p_2 p_3}.
\end{equation}

\vspace{2mm}

Theorem~1.1 gives (see (\ref{drdr1}))

\vspace{-2mm}
\begin{equation}
\label{20235}
\lim\limits_{p_1,p_2,p_3\to\infty}
A_{p_1 p_2 p_3}=0.
\end{equation}

\vspace{2mm}

Further, in complete analogy with (\ref{after5001}) and using (\ref{after80xx}), we obtain
$$
D_{p_1 p_2 p_3}=
$$
$$
=
\sum_{j_2=0}^{p_2}\left(\sum_{j_1=0}^{\min\{p_1,p_3\}}
C_{j_1 j_2 j_1}\right)^2=
\sum_{j_2=0}^{p_2}\left(\sum_{j_1=\min\{p_1,p_3\}+1}^{\infty}
C_{j_1 j_2 j_1}\right)^2\le
$$
\begin{equation}
\label{20236}
~~~~~~~~~~\le
\sum_{j_2=0}^{\infty}\left(\sum_{j_1=\min\{p_1,p_3\}+1}^{\infty}
C_{j_1 j_2 j_1}\right)^2\le
\frac{K}{\left(\min\{p_1,p_3\}\right)^{2-\varepsilon}}\ \to\  0
\end{equation}

\vspace{1mm}
\noindent
if $p_1,p_2,p_3\to\infty,$ where 
$\varepsilon$ is an arbitrary
small positive real number,
constant $K$ is independent of $p$.

We have
$$
B_{p_1 p_2 p_3}=
$$
$$
=
{\sf M}\left\{\Biggl(\Biggl(
\hbox{\vtop{\offinterlineskip\halign{
\hfil#\hfil\cr
{\rm l.i.m.}\cr
$\stackrel{}{{}_{p_3\to \infty}}$\cr
}} }\sum_{j_3=0}^{p_3}\sum_{j_1=0}^{\infty}  
C_{j_3 j_1 j_1}  
J'[\phi_{j_3}]^{(i_3)}_{T,t}-
\sum_{j_3=0}^{p_3}\sum_{j_1=0}^{\infty}  
C_{j_3 j_1 j_1}  
J'[\phi_{j_3}]^{(i_3)}_{T,t}\Biggr)+\Biggr.\right.
$$
$$
\left.\Biggl.+\Biggl(\sum_{j_3=0}^{p_3}\sum_{j_1=0}^{\infty}  
C_{j_3 j_1 j_1}  
J'[\phi_{j_3}]^{(i_3)}_{T,t}-
\sum_{j_3=0}^{p_3}\sum_{j_1=0}^{\min\{p_1,p_2\}}
C_{j_3j_1j_1}J'[\phi_{j_3}]^{(i_3)}_{T,t}\Biggr)\Biggr)^2\right\}\le
$$

\vspace{-1mm}
\begin{equation}
\label{20237}
\le
2 E_{p_3}+ 2 F_{p_1 p_2 p_3},
\end{equation}

\vspace{2mm}
\noindent
where
$$
E_{p_3}=
$$
$$
=
{\sf M}\left\{\Biggl(
\hbox{\vtop{\offinterlineskip\halign{
\hfil#\hfil\cr
{\rm l.i.m.}\cr
$\stackrel{}{{}_{p_3\to \infty}}$\cr
}} }\sum_{j_3=0}^{p_3}\sum_{j_1=0}^{\infty}  
C_{j_3 j_1 j_1}  
J'[\phi_{j_3}]^{(i_3)}_{T,t}-
\sum_{j_3=0}^{p_3}\sum_{j_1=0}^{\infty}  
C_{j_3 j_1 j_1}  
J'[\phi_{j_3}]^{(i_3)}_{T,t}\Biggr)^2\right\},
$$

\vspace{2mm}
$$
F_{p_1 p_2 p_3}=
$$
$$
=
{\sf M}\left\{
\Biggl(\sum_{j_3=0}^{p_3}\sum_{j_1=0}^{\infty}  
C_{j_3 j_1 j_1}  
J'[\phi_{j_3}]^{(i_3)}_{T,t}-
\sum_{j_3=0}^{p_3}\sum_{j_1=0}^{\min\{p_1,p_2\}}
C_{j_3j_1j_1}J'[\phi_{j_3}]^{(i_3)}_{T,t}\Biggr)^2\right\}=
$$
$$
={\sf M}\left\{
\Biggl(
\sum_{j_3=0}^{p_3}\sum_{j_1=\min\{p_1,p_2\}+1}^{\infty}
C_{j_3j_1j_1}J'[\phi_{j_3}]^{(i_3)}_{T,t}\Biggr)^2\right\}=
$$
\begin{equation}
\label{20238}
=\sum_{j_3=0}^{p_3}\Biggl(\sum_{j_1=\min\{p_1,p_2\}+1}^{\infty}
C_{j_3j_1j_1}\Biggr)^2.
\end{equation}

\vspace{4mm}

By analogy with (\ref{after1804}) we get
$$
\sum_{j_3=0}^{p_3}\Biggl(\sum_{j_1=\min\{p_1,p_2\}+1}^{\infty}
C_{j_3j_1j_1}\Biggr)^2\le
$$
$$
\le
\sum_{j_3=0}^{\infty}\Biggl(\sum_{j_1=\min\{p_1,p_2\}+1}^{\infty}
C_{j_3j_1j_1}\Biggr)^2\le
$$
\begin{equation}
\label{20239}
\le \frac{K}{\left(\min\{p_1,p_2\}\right)^2}\ \to\  0
\end{equation}

\vspace{3mm}
\noindent
if $p_1,p_2,p_3\to\infty,$ where constant $K$ does not depend on $p.$

Moreover,
\begin{equation}
\label{202310}
\lim\limits_{p_3\to\infty}E_{p_3}=
\lim\limits_{p_1,p_2,p_3\to\infty}E_{p_3}=0.
\end{equation}

\vspace{2mm}

Combining (\ref{20237})--(\ref{202310}), we obtain
\begin{equation}
\label{202311}
\lim\limits_{p_1,p_2,p_3\to\infty}B_{p_1 p_2 p_3}=0.
\end{equation}

\vspace{2mm}

Consider $C_{p_1 p_2 p_3}$. We have
$$
C_{p_1 p_2 p_3}=
$$
$$
=
{\sf M}\left\{\Biggl(\Biggl(
\hbox{\vtop{\offinterlineskip\halign{
\hfil#\hfil\cr
{\rm l.i.m.}\cr
$\stackrel{}{{}_{p_1\to \infty}}$\cr
}} }\sum_{j_1=0}^{p_1}\sum_{j_3=0}^{\infty}  
C_{j_3 j_3 j_1}  
J'[\phi_{j_1}]^{(i_1)}_{T,t}-
\sum_{j_1=0}^{p_1}\sum_{j_3=0}^{\infty}  
C_{j_3 j_3 j_1}  
J'[\phi_{j_1}]^{(i_1)}_{T,t}\Biggr)+\Biggr.\right.
$$
$$
\left.\Biggl.+\Biggl(\sum_{j_1=0}^{p_1}\sum_{j_3=0}^{\infty}  
C_{j_3 j_3 j_1}  
J'[\phi_{j_1}]^{(i_1)}_{T,t}-
\sum_{j_1=0}^{p_1}\sum_{j_3=0}^{\min\{p_2,p_3\}}
C_{j_3j_3j_1}J'[\phi_{j_1}]^{(i_1)}_{T,t}\Biggr)\Biggr)^2\right\}\le
$$

\vspace{-1mm}
\begin{equation}
\label{202350}
\le
2 G_{p_1}+ 2 H_{p_1 p_2 p_3},
\end{equation}

\vspace{2mm}
\noindent
where
$$
G_{p_1}=
$$
$$
=
{\sf M}\left\{\Biggl(
\hbox{\vtop{\offinterlineskip\halign{
\hfil#\hfil\cr
{\rm l.i.m.}\cr
$\stackrel{}{{}_{p_1\to \infty}}$\cr
}} }\sum_{j_1=0}^{p_1}\sum_{j_3=0}^{\infty}  
C_{j_3 j_3 j_1}  
J'[\phi_{j_1}]^{(i_1)}_{T,t}-
\sum_{j_1=0}^{p_1}\sum_{j_3=0}^{\infty}  
C_{j_3 j_3 j_1}  
J'[\phi_{j_1}]^{(i_1)}_{T,t}\Biggr)^2\right\},
$$

\vspace{2mm}
$$
H_{p_1 p_2 p_3}=
$$
$$
=
{\sf M}\left\{
\Biggl(\sum_{j_1=0}^{p_1}\sum_{j_3=0}^{\infty}  
C_{j_3 j_3 j_1}  
J'[\phi_{j_1}]^{(i_1)}_{T,t}-
\sum_{j_1=0}^{p_1}\sum_{j_3=0}^{\min\{p_2,p_3\}}
C_{j_3j_3j_1}J'[\phi_{j_1}]^{(i_1)}_{T,t}\Biggr)^2\right\}=
$$
$$
={\sf M}\left\{
\Biggl(
\sum_{j_1=0}^{p_1}\sum_{j_3=\min\{p_2,p_3\}+1}^{\infty}
C_{j_3j_3j_1}J'[\phi_{j_1}]^{(i_1)}_{T,t}\Biggr)^2\right\}=
$$
\begin{equation}
\label{202351}
=\sum_{j_1=0}^{p_1}\Biggl(\sum_{j_3=\min\{p_2,p_3\}+1}^{\infty}
C_{j_3j_3j_1}\Biggr)^2.
\end{equation}

\vspace{4mm}

By analogy with (\ref{after1903}) we get
$$
\sum_{j_1=0}^{p_1}\Biggl(\sum_{j_3=\min\{p_2,p_3\}+1}^{\infty}
C_{j_3j_3j_1}\Biggr)^2\le
$$
$$
\le
\sum_{j_1=0}^{\infty}\Biggl(\sum_{j_3=\min\{p_2,p_3\}+1}^{\infty}
C_{j_3j_3j_1}\Biggr)^2\le
$$
\begin{equation}
\label{202352}
\le \frac{K}{\left(\min\{p_2,p_3\}\right)^2}\ \to\  0
\end{equation}

\vspace{3mm}
\noindent
if $p_1,p_2,p_3\to\infty,$ where constant $K$ does not depend on $p.$

Moreover,
\begin{equation}
\label{202353}
\lim\limits_{p_1\to\infty}G_{p_1}=
\lim\limits_{p_1,p_2,p_3\to\infty}G_{p_1}=0.
\end{equation}

\vspace{2mm}

Combining (\ref{202350})--(\ref{202353}), we obtain
\begin{equation}
\label{202354}
\lim\limits_{p_1,p_2,p_3\to\infty}C_{p_1 p_2 p_3}=0.
\end{equation}

\vspace{2mm}

The relations (\ref{20234})--(\ref{20236}),  (\ref{202311}),  
(\ref{202354}) complete the proof of Theorem~2.41.
Theorem~2.41 is proved.

\section{Generalization of Theorem~2.30 for 
Complete Or\-tho\-nor\-mal Sys\-tems of Functions in $L_2([t, T])$
and $\psi_1(\tau),\ldots,\psi_k(\tau)$ $\in $ $L_2([t, T])$
such that the Condition (\ref{drdr1001xyx}) is Satisfied}

In this section, we generalize Theorem~2.30 to the case of
complete ortho\-nor\-mal systems of functions in the space $L_2([t, T])$
and $\psi_1(\tau),$ $\ldots,$ $\psi_k(\tau)$ $\in $ $L_2([t, T])$
such that the condition (\ref{drdr1001xyx})  is satisfied.

Let $(\Omega,{\rm F},{\sf P})$ be a complete probability
space and let $w(t,\omega)\stackrel{\sf def}{=}w_t:$ 
$[0, T]\times \Omega\rightarrow {\bf R}$
be the standard Wiener process
defined on the probability space $(\Omega,{\rm F},{\sf P}).$

Let us consider the family of $\sigma$-algebras
$\left\{{\rm F}_t,\ t\in[0,T]\right\}$ defined
on the probability space $(\Omega,{\rm F},{\sf P})$ and
connected
with the Wiener process $w_t$ in such a way that

1.\ ${\rm F}_s\subset {\rm F}_t\subset {\rm F}$\ for
$s<t.$

2.\ The Wiener process $w_t$ is ${\rm F}_t$-measurable for all
$t\in[0,T].$

3.\ The process $w_{t+\Delta}-w_{t}$ for all
$t\ge 0,$ $\Delta>0$ is independent with
the events of $\sigma$-algebra
${\rm F}_{t}.$

Let $\xi(\tau,\omega)\stackrel{\sf def}{=}
\xi_{\tau}:$ $[0, T]\times\Omega \to {\bf R}$ 
be some random process, which is measurable
with respect to the pair of variables
$(\tau,\omega)$ and satisfies to the following
condition
$$
\int\limits_t^T |\xi_{\tau}| d\tau <\infty\ \ \ \hbox{w.~p.~1}\ \ \ (t\ge 0).
$$

Let $\tau_j^{(N)},$ $j=0, 1, \ldots, N$ 
be a partition of the interval $[t, T],$ $t\ge 0$ such that
\begin{equation}
\label{dsds4}
t=\tau_0^{(N)}<\tau_1^{(N)}<\ldots <\tau_N^{(N)}=T,\ \ \ \
\max\limits_{0\le j\le N-1}\left|\tau_{j+1}^{(N)}-\tau_j^{(N)}\right|\to 0\ \
\hbox{if}\ \ N\to \infty.
\end{equation}

\noindent
Further, for simplicity, we write $\tau_j$ instead of 
$\tau_j^{(N)}.$

Consider the definition of the Stratonovich stochastic integral, 
which differs from the definition given in Sect.~2.1.1.

The mean-square limit (if it exists)
\begin{equation}
\label{dsds5}
~~~~~~\hbox{\vtop{\offinterlineskip\halign{
\hfil#\hfil\cr
{\rm l.i.m.}\cr
$\stackrel{}{{}_{N\to \infty}}$\cr
}} }\sum_{j=0}^{N-1}
\frac{1}{\tau_{j+1}-\tau_j}
\int\limits_{\tau_j}^{\tau_{j+1}}\xi_s ds\
\left(w_{\tau_{j+1}}-
w_{\tau_j}\right)
\stackrel{\sf def}{=}\int\limits_t^T \xi_{\tau} \circ dw_\tau
\end{equation}
is called \cite{SU11}, \cite{bardina10} the Stratonovich stochastic integral 
of the process $\xi_{\tau}$, $\tau\in [t, T]$,
where $\tau_j,$ $j=0, 1, \ldots, N$
is a partition of the interval $[t, T]$ 
satisfying the condition (\ref{dsds4}).

We also denote by
$$
\int\limits_t^{\tau} \xi_{s} \circ dw_s
$$
the Stratonovich stochastic integral like (\ref{dsds5}) (if it exists) 
of $\xi_s {\bf 1}_{\{s\in [t, \tau]\}}$ for $\tau\in [t, T],$ $t\ge 0.$

It is known \cite{bardina10} (Lemma~A.2) that the following 
iterated Stratonovich stochastic integral 
\begin{equation}
\label{dsds7}
~~~~~~~J^{S}[\psi^{(k)}]_{\tau,t}^{(i_1\ldots i_k)}=
\int\limits_t^{\tau}
\psi_k(t_k)\ldots \int\limits_t^{t_2}
\psi_1(t_1) \circ d{\bf w}_{t_1}^{(i_1)}\ldots \circ d{\bf w}_{t_k}^{(i_k)}
\end{equation}
exists for the case $i_1=\ldots=i_k\ne 0$, 
where $\tau\in [t, T],$ $\psi_1(\tau),\ldots,\psi_k(\tau)\in L_2([t, T]),$
$i_1,\ldots,i_k=0,1,\ldots,m,$\ 
${\bf w}_{\tau}^{(i)}$ $(i=1,\ldots,m)$ 
are independent 
standard Wiener processes defined as above in this section,
and ${\bf w}_{\tau}^{(0)}=\tau$.

Note that in \cite{new-new-18} (2021) an analogue of Theorem~2.12 (1997) 
is proved for the integral $J^{S}[\psi^{(k)}]_{\tau,t}^{(i_1\ldots i_k)}$ 
($i_1=\ldots=i_k\ne 0,$ $\psi_1(\tau),\ldots,\psi_k(\tau)\in L_2([t, T])$).

Let us denote
\begin{equation}
\label{dsds9}
~~~~~~J[\psi^{(k)}]_{T,t}^{(i_1\ldots i_k)}+
\sum_{r=1}^{\left[k/2\right]}\frac{1}{2^r}
\sum_{(s_r,\ldots,s_1)\in {\rm A}_{k,r}}
J[\psi^{(k)}]_{T,t}^{s_r,\ldots,s_1}\stackrel{\sf def}{=}\bar J^{*}[\psi^{(k)}]_{T,t}^{(i_1\ldots i_k)},
\end{equation}
where $\psi_1(\tau),\ldots,\psi_k(\tau)\in L_2([t, T])$,
$\psi_l(\tau)\psi_{l-1}(\tau)\in L_2([t, T])$ $(l=2, 3,\ldots, k),$
$J[\psi^{(k)}]_{T,t}^{(i_1\ldots i_k)}$ is the iterated It\^{o} stochastic
integral
\begin{equation}
\label{dsds12}
~~~~~~~~~~~J[\psi^{(k)}]_{T,t}^{(i_1\ldots i_k)}=
\int\limits_t^{T}\psi_k(t_k) \ldots 
\int\limits_t^{t_{2}}
\psi_1(t_1) d{\bf w}_{t_1}^{(i_1)}\ldots
d{\bf w}_{t_k}^{(i_k)},
\end{equation}
$\sum\limits_{\emptyset}$ is supposed to be equal to zero;
another notations as in Theorem~2.12.

We will also notice that 
\begin{equation}
\label{april200}
J^{S}[\psi^{(1)}]_{T,t}^{(i_1)}=J[\psi^{(1)}]_{T,t}^{(i_1)}\ \ \ \hbox{w.~p.~1.}
\end{equation}

Further, by analogy with (\ref{after8}), (\ref{after7})
and using (\ref{chain401}) (also see Theorem~1.23)
instead of (\ref{rezo7}) 
we obtain 
the following generalization of (\ref{after8}) to the case 
of an arbitrary 
complete ortho\-nor\-mal system of functions in the space $L_2([t, T])$
and $\psi_1(\tau),$ $\ldots,$ $\psi_k(\tau)$ $\in $ $L_2([t, T])$

\vspace{-4mm}
$$
\sum_{j_1=0}^{p_1}\ldots\sum_{j_k=0}^{p_k}
C_{j_k\ldots j_1}
\prod_{l=1}^k \zeta_{j_l}^{(i_l)}=
\sum_{j_1=0}^{p_1}\ldots\sum_{j_k=0}^{p_k}
C_{j_k\ldots j_1}
J'[\phi_{j_1}\ldots \phi_{j_k}]_{T,t}^{(i_1\ldots i_k)}+
$$

$$
+\sum_{j_1=0}^{p_1}\ldots\sum_{j_k=0}^{p_k}
C_{j_k\ldots j_1}
\sum\limits_{r=1}^{[k/2]}
\sum_{\stackrel{(\{\{g_1, g_2\}, \ldots, 
\{g_{2r-1}, g_{2r}\}\}, \{q_1, \ldots, q_{k-2r}\})}
{{}_{\{g_1, g_2, \ldots, 
g_{2r-1}, g_{2r}, q_1, \ldots, q_{k-2r}\}=\{1, 2, \ldots, k\}}}}
\prod\limits_{s=1}^r
{\bf 1}_{\{i_{g_{{}_{2s-1}}}=~i_{g_{{}_{2s}}}\ne 0\}}\times
$$

\vspace{4mm} 
\begin{equation}
\label{after8xxds1}
\times{\bf 1}_{\{j_{g_{{}_{2s-1}}}=~j_{g_{{}_{2s}}}\}}
J'[\phi_{j_{q_1}}\ldots \phi_{j_{q_{k-2r}}}]_{T,t}^{(i_{q_1}\ldots i_{q_{k-2r}})}\ \ \ \hbox{w.~p.~1,}
\end{equation}

\vspace{3mm}
\noindent
where $k\ge 2,$ $J'[\phi_{j_1}\ldots \phi_{j_k}]_{T,t}^{(i_1\ldots i_k)},$
$J'[\phi_{j_{q_1}}\ldots \phi_{j_{q_{k-2r}}}]_{T,t}^{(i_{q_1}\ldots i_{q_{k-2r}})}$
are multiple Wie\-ner sto\-chas\-tic integrals 
(see (\ref{WiI})) and $J'[\phi_{j_{q_1}}\ldots \phi_{j_{q_{k-2r}}}]_{T,t}^{(i_{q_1}\ldots i_{q_{k-2r}})}
\stackrel{\sf def}{=}1$ for $k=2r.$

Using the equalities (\ref{chain401}) and (\ref{razzar1}), we can 
reformulate Theorem~1.16 as follows
\begin{equation}
\label{dsds11}
J[\psi^{(k)}]_{T,t}^{(i_1\ldots i_k)}=
\hbox{\vtop{\offinterlineskip\halign{
\hfil#\hfil\cr
{\rm l.i.m.}\cr
$\stackrel{}{{}_{p_1,\ldots,p_k\to \infty}}$\cr
}} }\sum_{j_1=0}^{p_1}\ldots\sum_{j_k=0}^{p_k}
C_{j_k\ldots j_1}J'[\phi_{j_1}\ldots \phi_{j_k}]_{T,t}^{(i_1\ldots i_k)}\ \ \ \hbox{w.~p.~1},
\end{equation}

\noindent
where 
$J'[\phi_{j_1}\ldots \phi_{j_k}]_{T,t}^{(i_1\ldots i_k)}$ is the 
multiple Wiener stochastic integral
defined by (\ref{WiI});
another notations are the same as in Theorem~1.16.

Passing to the limit 
$\hbox{\vtop{\offinterlineskip\halign{
\hfil#\hfil\cr
{\rm l.i.m.}\cr
$\stackrel{}{{}_{p_1,\ldots,p_k\to \infty}}$\cr
}} }$ in (\ref{after8xxds1}) and using the equality (\ref{dsds11}), we get w.~p.~1 

\newpage
\noindent
$$
\hbox{\vtop{\offinterlineskip\halign{
\hfil#\hfil\cr
{\rm l.i.m.}\cr
$\stackrel{}{{}_{p_1,\ldots,p_k\to \infty}}$\cr
}} }\sum_{j_1=0}^{p_1}\ldots\sum_{j_k=0}^{p_k}
C_{j_k\ldots j_1}
\zeta_{j_1}^{(i_1)}\ldots \zeta_{j_k}^{(i_k)}
=J[\psi^{(k)}]_{T,t}^{(i_1\ldots i_k)}
+
$$

$$
+
\sum\limits_{r=1}^{[k/2]}
\sum_{\stackrel{(\{\{g_1, g_2\}, \ldots, 
\{g_{2r-1}, g_{2r}\}\}, \{q_1, \ldots, q_{k-2r}\})}
{{}_{\{g_1, g_2, \ldots, 
g_{2r-1}, g_{2r}, q_1, \ldots, q_{k-2r}\}=\{1, 2, \ldots, k\}}}}
\prod\limits_{s=1}^r
{\bf 1}_{\{i_{g_{{}_{2s-1}}}=~i_{g_{{}_{2s}}}\ne 0\}}\times
$$

\vspace{1mm}
\begin{equation}
\label{after501ds1}
~~~~~~~~\times \hbox{\vtop{\offinterlineskip\halign{
\hfil#\hfil\cr
{\rm l.i.m.}\cr
$\stackrel{}{{}_{p_1,\ldots,p_k\to \infty}}$\cr
}} }\sum_{j_1=0}^{p_1}\ldots\sum_{j_k=0}^{p_k}
C_{j_k\ldots j_1}
\prod\limits_{s=1}^r{\bf 1}_{\{j_{g_{{}_{2s-1}}}=~j_{g_{{}_{2s}}}\}}
J'[\phi_{j_{q_1}}\ldots \phi_{j_{q_{k-2r}}}]_{T,t}^{(i_{q_1}\ldots i_{q_{k-2r}})},
\end{equation}

\vspace{1mm}
\noindent
where 
$J'[\phi_{j_{q_1}}\ldots \phi_{j_{q_{k-2r}}}]_{T,t}^{(i_{q_1}\ldots i_{q_{k-2r}})}$ is the 
multiple Wiener stochastic integral
defined by (\ref{WiI}),
$J[\psi^{(k)}]_{T,t}^{(i_1\ldots i_k)}$ is the iterated It\^{o} stochastic
integral (\ref{dsds12}).

Suppose that $\{\phi_j(x)\}_{j=0}^{\infty}$ is an arbitrary
complete orthonormal system of functions in $L_2([t, T])$
and $\Phi_1(\tau), \Phi_2(\tau)\in L_2([t, T])$.
Then we have
$$
\sum_{j=0}^{\infty}\left|\int\limits_t^s \phi_j(\tau)\Phi_1(\tau)d\tau 
\int\limits_s^T \phi_j(\tau)\Phi_2(\tau)d\tau\right|\le 
$$
$$
\le \frac{1}{2}\sum_{j=0}^{\infty}
\left(\left(\int\limits_t^T {\bf 1}_{\{\tau<s\}}\phi_j(\tau)\Phi_1(\tau)d\tau\right)^2+ 
\left(\int\limits_t^T {\bf 1}_{\{\tau>s\}}\phi_j(\tau)\Phi_2(\tau)d\tau\right)^2\right)=
$$
\begin{equation}
\label{dsds14}
=\frac{1}{2}\left(\int\limits_t^s \Phi_1^2(\tau)d\tau+
\int\limits_s^T\Phi_2^2(\tau)d\tau\right)\le
\frac{1}{2}\left(\left\Vert\Phi_1\right\Vert_{L_2([t,T])}^2+
\left\Vert\Phi_2\right\Vert_{L_2([t,T])}^2\right)<\infty,
\end{equation}
i.e. 
\begin{equation}
\label{dsds14fffff}
\left|\sum_{j=0}^{p}\int\limits_t^s \phi_j(\tau)\Phi_1(\tau)d\tau 
\int\limits_s^T \phi_j(\tau)\Phi_2(\tau)d\tau\right|\le C<\infty,
\end{equation}
where $p\in{\bf N}.$

By interpreting the integrals in (\ref{after9})--(\ref{after400}) as 
Lebesgue integrals, using Fubini's Theorem in (\ref{after9}) and
Lebesgue's 
Dominated Convergence Theorem in (\ref{after450}), we 
obtain (\ref{after80}) (see (\ref{dsds14fffff})) for
the case of an arbitrary complete 
orthonormal system of functions in the space $L_2([t, T])$
and $\psi_1(\tau),\ldots, \psi_k(\tau)\in L_2([t, T])$.

Using the equality (\ref{5tzzz}) 
for the case of an arbitrary complete 
orthonormal system of functions in the space $L_2([t, T])$
and $\psi_1(\tau),\psi_2(\tau)\in L_2([t, T])$ as well as
Fubini's Theorem when deriving (\ref{r12345x}), we obtain the generalization of
(\ref{after500}) for the case of an arbitrary complete 
orthonormal system of functions in the space $L_2([t, T])$
and $\psi_1(\tau),\ldots, \psi_k(\tau)\in L_2([t, T])$.

Repeating the steps of the proof of Theorem~2.30 below
the formula (\ref{after501}) using (\ref{dsds9}), (\ref{after501ds1}) or steps 
of the proof of Theorem~2.32 using (\ref{dsds9}), (\ref{after501ds1}), we obtain
for complete 
orthonormal systems $\{\phi_j(x)\}_{j=0}^{\infty}$
$(\phi_0(x)=1/\sqrt{T-t})$ 
in the space $L_2([t, T])$
and $\psi_1(\tau),\ldots, \psi_k(\tau)\in L_2([t, T]),$
$\psi_l(\tau)\psi_{l-1}(\tau)\in L_2([t, T])$ $(l=2, 3,\ldots, k)$
(for which the condition (\ref{drdr1001xyx}) is satisfied) the following equality
$$
\hbox{\vtop{\offinterlineskip\halign{
\hfil#\hfil\cr
{\rm l.i.m.}\cr
$\stackrel{}{{}_{p_1,\ldots,p_k\to \infty}}$\cr
}} }
\sum_{j_1=0}^{p_1}\ldots\sum_{j_k=0}^{p_k}
C_{j_k \ldots j_1}\prod\limits_{l=1}^k \zeta_{j_l}^{(i_l)}
=
$$

\vspace{-4mm}
\begin{equation}
\label{after333ds1}
~~~~~=
J[\psi^{(k)}]_{T,t}^{(i_1\ldots i_k)}+
\sum_{r=1}^{\left[k/2\right]}\frac{1}{2^r}
\sum_{(s_r,\ldots,s_1)\in {\rm A}_{k,r}}
J[\psi^{(k)}]_{T,t}^{s_r,\ldots,s_1}=
\bar J^{*}[\psi^{(k)}]_{T,t}^{(i_1\ldots i_k)}
\end{equation}

\vspace{1mm}
\noindent
w.~p.~1, where notations in (\ref{after333ds1}) are the same as in Theorem~2.12
and $\bar J^{*}[\psi^{(k)}]_{T,t}^{(i_1\ldots i_k)}$ is defined by
(\ref{dsds9}).

Thus the following two theorems are proved 
(compare with Theorem~1.16 (Sect. 1.11) on the expansion of iterated 
It\^{o} stochastic integrals).

\vspace{1mm}

{\bf Theorem~2.42}\ \cite{arxiv-5}, \cite{arxiv-10}, \cite{arxiv-11}.\ {\it Assume that
the complete orthonormal system $\{\phi_j(x)\}_{j=0}^{\infty}$
$(\phi_0(x)=1/\sqrt{T-t})$ 
in the space $L_2([t, T])$ and
$\psi_1(\tau),\ldots, \psi_k(\tau)\in L_2([t, T]),$
$\psi_l(\tau)\psi_{l-1}(\tau)\in L_2([t, T])$ $(l=2, 3,\ldots, k)$
are such that 

\vspace{-5mm}
$$
\lim\limits_{p_1,\ldots,p_k\to\infty}~
\sum\limits_{j_1=0}^{p_1}\ldots \sum\limits_{j_q=0}^{p_q}\ldots \sum\limits_{j_k=0}^{p_k}~
\biggl|_{q\ne g_1, g_2, \ldots, g_{2r-1},g_{2r}}\times
$$

\vspace{2mm}
$$
\times
\Biggl(~\sum\limits_{j_{g_1}=0}^{\min\{p_{g_1}, p_{g_2}\}} \sum\limits_{j_{g_3}=0}^{\min\{p_{g_3}, p_{g_4}\}}\ldots \Biggr.
\sum\limits_{j_{g_{2r-1}}=0}^{\min\{p_{g_{2r-1}}, p_{g_{2r}}\}}
C_{j_k\ldots j_1}\biggl|_{j_{g_1}=j_{g_2},\ldots, j_{g_{2r-1}}=j_{g_{2r}}}-
$$

\vspace{1mm}
\begin{equation}
\label{novorigin1}
\Biggl.-\frac{1}{2^r} \prod\limits_{l=1}^r {\bf 1}_{\{g_{2l}=g_{2l-1}+1\}}
C_{j_k \ldots j_1}\biggl|_{(j_{g_2} j_{g_1})\curvearrowright (\cdot)
\ldots (j_{g_{2r}} j_{g_{2r-1}})\curvearrowright (\cdot),
j_{g_{{}_{1}}}=~j_{g_{{}_{2}}},\ldots, j_{g_{{}_{2r-1}}}=~j_{g_{{}_{2r}}}
}\biggr.\Biggr)^2=0
\end{equation}

\vspace{3mm}
\noindent
for all $r=1, 2,\ldots,[k/2]$ 
and for all possible $g_1,g_2,\ldots,g_{2r-1},g_{2r}$ {\rm (}see {\rm (\ref{leto5007after}))}.
Then$,$ for the sum $\bar J^{*}[\psi^{(k)}]_{T,t}^{(i_1\ldots i_k)}$
of iterated It\^{o} stochastic integrals 
defined by {\rm (\ref{dsds9})}
the following 
expansion 
\begin{equation}
\label{january19c}
\bar J^{*}[\psi^{(k)}]_{T,t}^{(i_1\ldots i_k)}=
\hbox{\vtop{\offinterlineskip\halign{
\hfil#\hfil\cr
{\rm l.i.m.}\cr
$\stackrel{}{{}_{p_1,\ldots,p_k\to \infty}}$\cr
}} }
\sum_{j_1=0}^{p_1}\ldots\sum_{j_k=0}^{p_k}
C_{j_k \ldots j_1}\prod\limits_{l=1}^k \zeta_{j_l}^{(i_l)}
\end{equation}

\vspace{1mm}
\noindent
that converges in the mean-square sense is valid, where 

\vspace{-2mm}
$$
C_{j_k \ldots j_1}=\int\limits_t^T\psi_k(t_k)\phi_{j_k}(t_k)\ldots
\int\limits_t^{t_2}
\psi_1(t_1)\phi_{j_1}(t_1)
dt_1\ldots dt_k
$$

\vspace{1mm}
\noindent
is the Fourier coefficient, 
${\rm l.i.m.}$ is a limit in the mean-square sense,
$i_1, \ldots, i_k=0, 1,\ldots,m,$
$$
\zeta_{j}^{(i)}=
\int\limits_t^T \phi_{j}(\tau) d{\bf w}_{\tau}^{(i)}
$$ 

\noindent
are independent standard Gaussian random variables for various 
$i$ or $j$ {\rm (}when $i\ne 0${\rm ),}
${\bf w}_{\tau}^{(i)}$ 
$(i=1,\ldots,m)$ are independent standard Wiener processes$,$
${\bf w}_{\tau}^{(0)}=\tau.$}

\vspace{2mm}

{\bf Theorem~2.43}\ \cite{arxiv-5}, \cite{arxiv-10}, \cite{arxiv-11}.\ {\it Assume that
the complete orthonormal system $\{\phi_j(x)\}_{j=0}^{\infty}$
$(\phi_0(x)=1/\sqrt{T-t})$ 
in the space $L_2([t, T])$ and 
$\psi_1(\tau),\ldots,\psi_k(\tau)\in L_2([t, T]),$
$\psi_l(\tau)\psi_{l-1}(\tau)\in L_2([t, T])$ $(l=2, 3,\ldots, k)$
are such that 
the condition

\vspace{-4mm}
$$
\lim\limits_{p\to\infty}
\sum\limits_{\stackrel{j_1,\ldots,j_q,\ldots,j_k=0}{{}_{q\ne g_1, g_2, \ldots, g_{2r-1},
g_{2r}}}}^p
\left(S_{l_1}S_{l_2}\ldots S_{l_{d}}
\left\{\bar C^{(p)}_{j_k\ldots j_q \ldots j_1}\biggl|_{q\ne g_1,g_2,\ldots,g_{2r-1}, g_{2r}}
\right\}\right)^2=0
$$

\vspace{1mm}
\noindent
holds for all possible $g_1,g_2,\ldots,g_{2r-1},g_{2r}$ {\rm (}see {\rm (\ref{leto5007after}))}
and $l_1, l_2, \ldots, l_{d}$ such that
$l_1, l_2, \ldots, l_{d}\in \{1,2,\ldots, r\},$\
$l_1>l_2>\ldots >l_{d},$\ $d=0, 1, 2,\ldots, r-1,$\ 
where $r=1, 2,\ldots,[k/2]$ and

\vspace{-4mm}
$$
S_{l_1}S_{l_2}\ldots S_{l_{d}}
\left\{\bar C^{(p)}_{j_k\ldots j_q \ldots j_1}\biggl|_{q\ne g_1,g_2,\ldots,g_{2r-1}, g_{2r}}
\right\}\stackrel{\sf def}{=}
\bar C^{(p)}_{j_k\ldots j_q \ldots j_1}\biggl|_{q\ne g_1,g_2,\ldots,g_{2r-1}, g_{2r}}
$$

\vspace{2mm}
\noindent
for $d=0.$
Then$,$ for the sum $\bar J^{*}[\psi^{(k)}]_{T,t}^{(i_1\ldots i_k)}$
of iterated It\^{o} stochastic integrals 
defined by {\rm (\ref{dsds9})}
the following 
expansion 
$$
\bar J^{*}[\psi^{(k)}]_{T,t}^{(i_1\ldots i_k)}=
\hbox{\vtop{\offinterlineskip\halign{
\hfil#\hfil\cr
{\rm l.i.m.}\cr
$\stackrel{}{{}_{p\to \infty}}$\cr
}} }
\sum_{j_1,\ldots,j_k=0}^{p}
C_{j_k \ldots j_1}\prod\limits_{l=1}^k \zeta_{j_l}^{(i_l)}
$$

\vspace{2mm}
\noindent
that converges in the mean-square sense is valid, where 

\vspace{-2mm}
$$
C_{j_k \ldots j_1}=\int\limits_t^T\psi_k(t_k)\phi_{j_k}(t_k)\ldots
\int\limits_t^{t_2}
\psi_1(t_1)\phi_{j_1}(t_1)
dt_1\ldots dt_k
$$

\vspace{2mm}
\noindent
is the Fourier coefficient, 
${\rm l.i.m.}$ is a limit in the mean-square sense,
$i_1, \ldots, i_k=0, 1,\ldots,m,$

\vspace{-4mm}
$$
\zeta_{j}^{(i)}=
\int\limits_t^T \phi_{j}(\tau) d{\bf w}_{\tau}^{(i)}
$$ 

\vspace{1mm}
\noindent
are independent standard Gaussian random variables for various 
$i$ or $j$ {\rm (}in the case when $i\ne 0${\rm )},
${\bf w}_{\tau}^{(i)}$ 
$(i=1,\ldots,m)$ are independent standard Wiener processes$,$
${\bf w}_{\tau}^{(0)}=\tau.$}

Note that in Theorems~2.42, 2.43 (the case $k=2$)
the condition 
$\psi_1(\tau)\psi_{2}(\tau)$ $\in L_2([t, T])$ 
can be omitted.

Using Theorem 2.12 together with Proposition 3.1 \cite{new-new-18} and the proof of
Lemma A.2 \cite{bardina10}, we can write 
$\bar J^{*}[\psi^{(k)}]_{T,t}^{(i_1\ldots i_k)}=J^{S}[\psi^{(k)}]_{T,t}^{(i_1\ldots i_k)}$
w.~p.~1 and reformulate Theorems~2.42, 2.43 for 
$J^{S}[\psi^{(k)}]_{T,t}^{(i_1\ldots i_k)}$ defined by (\ref{dsds7}).

Let us consider the special case $k=2$ of Theorem~2.42 in more detail.
In this case, the condition (\ref{novorigin1}) takes the following form
(compare with (\ref{5t}))
\begin{equation}
\label{novorigin20}
\sum_{j_1=0}^{\infty}
C_{j_1j_1}=\frac{1}{2}
\int\limits_t^T\psi_1(t_1)\psi_2(t_1)dt_1.
\end{equation}

As follows from Sect.~2.1.4, the equality (\ref{novorigin20})
is valid for the case 
of an arbitrary complete orthonormal 
system of functions in $L_2([t, T])$ and
$\psi_1(\tau), \psi_2(\tau)\in L_2([t, T]).$

From Proposition 3.1 \cite{new-new-18} for the case $k=2$ we obtain
$$
\int\limits_t^{T}
\psi_2(t_2)\int\limits_t^{t_2}
\psi_1(t_1) \circ d{\bf w}_{t_1}^{(i)}\circ d{\bf w}_{t_2}^{(i)}
=\int\limits_t^{T}
\psi_2(t_2)\int\limits_t^{t_2}
\psi_1(t_1) d{\bf w}_{t_1}^{(i)} d{\bf w}_{t_2}^{(i)}+
$$
\begin{equation}
\label{novorigin3}
+
\frac{1}{2}
\int\limits_t^T\psi_1(t_1)\psi_2(t_1)dt_1
\end{equation}

\vspace{1mm}
\noindent
w.~p.~1, where $\psi_1(\tau), \psi_2(\tau)\in L_2([t, T]),$
$i=1,\ldots,m,$
$$
\int\limits_t^{T}
\psi_2(t_2)\int\limits_t^{t_2}
\psi_1(t_1) \circ d{\bf w}_{t_1}^{(i)}\circ d{\bf w}_{t_2}^{(i)}
$$
is defined by (\ref{dsds5}), (\ref{dsds7}) and 
$$
\int\limits_t^{T}
\psi_2(t_2)\int\limits_t^{t_2}
\psi_1(t_1) d{\bf w}_{t_1}^{(i)} d{\bf w}_{t_2}^{(i)}
$$
is the iterated It\^{o} stochastic integral of the form (\ref{itoxx})
($k=2$).

On the other hand, it is not difficult to show that
\begin{equation}
\label{novorigin4}
~~~~~~~~\int\limits_t^{T}
\psi_2(t_2)\int\limits_t^{t_2}
\psi_1(t_1) \circ d{\bf w}_{t_1}^{(i)}\circ d{\bf w}_{t_2}^{(j)}
=\int\limits_t^{T}
\psi_2(t_2)\int\limits_t^{t_2}
\psi_1(t_1) d{\bf w}_{t_1}^{(i)} d{\bf w}_{t_2}^{(j)}
\end{equation}

\noindent
w.~p.~1, where $\psi_1(\tau), \psi_2(\tau)\in L_2([t, T]),$
$i\ne j$ 
$(i, j=1,\ldots,m),$ another notations are the same as in (\ref{novorigin3}).

Combining (\ref{novorigin3}) and (\ref{novorigin4}), we get
(see (\ref{dsds9}))
$$
\int\limits_t^{T}
\psi_2(t_2)\int\limits_t^{t_2}
\psi_1(t_1) \circ d{\bf w}_{t_1}^{(i_1)}\circ d{\bf w}_{t_2}^{(i_2)}
=\int\limits_t^{T}
\psi_2(t_2)\int\limits_t^{t_2}
\psi_1(t_1) d{\bf w}_{t_1}^{(i_1)} d{\bf w}_{t_2}^{(i_2)}+
$$
\begin{equation}
\label{novorigin5}
+
\frac{1}{2}{\bf 1}_{\{i_1=i_2\}}
\int\limits_t^T\psi_1(t_1)\psi_2(t_1)dt_1
\stackrel{\sf def}{=}\bar J^{*}[\psi^{(2)}]_{T,t}^{(i_1 i_2)}
\end{equation}

\vspace{1mm}
\noindent
w.~p.~1, where $\psi_1(\tau),\psi_2(\tau)\in L_2([t, T]),$
$i_1, i_2=1,\ldots,m.$

It is easy to see that the condition 
$\phi_0(x)=1/\sqrt{T-t}$ 
can be omitted in Theorems~2.42, 2.43 for the case $k=2$
(see the proof of Theorem~2.30).

Summing up the above arguments, we obtain
the following generalization of Theorem~2.3 to the case 
of an arbitrary complete orthonormal 
system of functions in $L_2([t, T])$ and
$\psi_1(\tau), \psi_2(\tau)\in L_2([t, T]).$

{\bf Theorem 2.44}\ \cite{arxiv-5}, \cite{arxiv-10}, \cite{arxiv-11}.\ {\it Suppose that 
$\{\phi_j(x)\}_{j=0}^{\infty}$ is an arbitrary complete orthonormal system of 
functions in the space $L_2([t, T])$ and
$\psi_1(\tau), \psi_2(\tau)\in L_2([t, T])$.
Then$,$ 
for the iterated Stra\-to\-novich stochastic integral
$$
J^{S}[\psi^{(2)}]_{T,t}^{(i_1 i_2)}=
\int\limits_t^{T}
\psi_2(t_2)\int\limits_t^{t_2}
\psi_1(t_1) \circ d{\bf w}_{t_1}^{(i_1)}\circ d{\bf w}_{t_2}^{(i_2)}\ \ \ (i_1, i_2=1,\ldots,m)
$$
the following expansion  
\begin{equation}
\label{novorigin10}
J^{S}[\psi^{(2)}]_{T,t}^{(i_1 i_2)}=\hbox{\vtop{\offinterlineskip\halign{
\hfil#\hfil\cr
{\rm l.i.m.}\cr
$\stackrel{}{{}_{p_1,p_2\to \infty}}$\cr
}} }\sum_{j_1=0}^{p_1}\sum_{j_2=0}^{p_2}
C_{j_2j_1}\zeta_{j_1}^{(i_1)}\zeta_{j_2}^{(i_2)}
\end{equation}

\noindent
that converges in the mean-square
sence is valid$,$ where the notations are the same as in Theorems {\rm 2.1--2.3}
and $J^{S}[\psi^{(2)}]_{T,t}^{(i_1 i_2)}$ is defined by {\rm (\ref{dsds7})}. 
}

Note that the analog of (\ref{novorigin10}) for $k=1$ is also true
(see (\ref{a1}) and (\ref{april200})).

In this section, it is also appropriate to mention 
the so-called multiple Stratonovich stochastic integral
\cite{SU11}, \cite{bardina10} (also see \cite{bugh1}).

The mean-square limit (if it exists)

\vspace{-6mm}
$$
\hbox{\vtop{\offinterlineskip\halign{
\hfil#\hfil\cr
{\rm l.i.m.}\cr
$\stackrel{}{{}_{N\to \infty}}$\cr
}} }\sum_{l_1=0}^{N-1}\ldots \sum_{l_k=0}^{N-1}
\frac{1}{\Delta\tau_{l_1}\ldots \Delta\tau_{l_k}}
\int\limits_{[\tau_{l_1},\tau_{l_1+1}] \times \ldots \times [\tau_{l_k},\tau_{l_k+1}]}
K(t_1,\ldots,t_k)
dt_1\ldots dt_k\times 
$$

\begin{equation}
\label{january19}
\times \Delta{\bf w}_{\tau_{l_1}}^{(i_1)}\ldots
\Delta{\bf w}_{\tau_{l_k}}^{(i_k)}
\stackrel{\sf def}{=}\bar J^{S}[K]_{T,t}^{(i_1\ldots i_k)} 
\end{equation}

\vspace{1mm}
\noindent
is called \cite{SU11}, \cite{bardina10} the multiple Stratonovich stochastic integral 
of the function  $K(t_1,\ldots,t_k)\in L_2([t, T]^k)$, where
$\Delta {\bf w}_{\tau_j}^{(i)}={\bf w}_{\tau_{j+1}}^{(i)}-{\bf w}_{\tau_{j}}^{(i)}$
$(i=0, 1,\ldots,m),$ $\Delta\tau_j=\tau_{j+1}-\tau_{j},$
$\left\{\tau_j\right\}_{j=0}^N$ 
is a partition of the interval $[t, T]$ 
satisfying the condition (\ref{dsds4}),
$i_1,\ldots,i_k=0,1,\ldots,m,$\ 
${\bf w}_{\tau}^{(i)}$ $(i=1,\ldots,m)$
are independent 
standard Wiener processes defined as above in this section,
and
${\bf w}_{\tau}^{(0)}=\tau$.

Note that in \cite{bardina10} the case $i_1=\ldots=i_k\ne 0$
was considered.
We also denote by $\bar J^{S}[K]_{s,t}^{(i_1\ldots i_k)}$
the stochastic integral 
(\ref{january19}) (if it exists) 
of the function $K(t_1,\ldots,t_k){\bf 1}_{\{(t_1,\ldots,t_k)\in [t, s]^k\}},$ 
where $K(t_1,\ldots,t_k)\in L_2([t, T]^k),$
$s\in [t, T],$ $t\ge 0.$

Let the function $K(t_1,\ldots,t_k)$ be chosen as follows

\vspace{-3mm}
\begin{equation}
\label{january19a}
K(t_1,\ldots,t_k)=
\left\{\begin{matrix}
\psi_1(t_1)\ldots \psi_k(t_k),\ &t_1\le \ldots \le t_k\cr\cr
0,\ &\hbox{\rm otherwise}
\end{matrix}
\right.,
\end{equation}

\noindent
where $\psi_1(\tau),\ldots, \psi_k(\tau)\in L_2([t, T]),$
$t_1,\ldots,t_k\in [t, T]$ $(k\ge 2)$ and 
$K(t_1)\equiv\psi_1(t_1)$ for $t_1\in[t, T].$

We will denote the 
multiple Stratonovich stochastic integral (\ref{january19})
of the function (\ref{january19a}) as 
$\bar J^{S}[\psi^{(k)}]_{T,t}^{(i_1\ldots i_k)}$.

It is known \cite{bardina10} (Lemma~A.2) that the Stratonovich 
stochastic integrals
$J^{S}[\psi^{(k)}]_{T,t}^{(i_1\ldots i_k)}$ and
$\bar J^{S}[\psi^{(k)}]_{T,t}^{(i_1\ldots i_k)}$
exist for the case $i_1=\ldots=i_k\ne 0.$
Moreover,
$J^{S}[\psi^{(k)}]_{T,t}^{(i_1\ldots i_k)}=
\bar J^{S}[\psi^{(k)}]_{T,t}^{(i_1\ldots i_k)}$
w.~p.~1 for this case
\cite{bardina10} (Lemma~A.2).

Recall that an expansion similar to (\ref{january19c}) for $p_1=\ldots=p_k=p$
was obtained in {\rm \cite{Rybakov3000}} for the multiple
Stratonovich stochastic integral
(\ref{january19}) under the condition of convergence of trace series
(see Remarks~2.4, 2.7 for details).

Recently, 
another approach to the expansion of integral (\ref{january19})
has been proposed (assuming that the integral (\ref{january19}) exists), 
where multiple Fourier--Walsh and Fourier--Haar series $(k\in{\bf N})$ have been applied
\cite{Rybakov3000xxx}.
The convergence was proved with respect to the special 
subsequence ($p_1=\ldots=p_k=p=2^m,$ $m\to\infty$ in a formula
similar to (\ref{january19c}) \cite{Rybakov3000xxx}).

\section{Expansion of Iterated Stratonovich Stochastic Integrals
of Multiplicity 3. The Case of an Ar\-bit\-ra\-ry Complete Orthonormal System of 
Functions $(\phi_0(x)=1/\sqrt{T-t})$ in the Space $L_2([t,T])$ and $\psi_1(\tau), \psi_2(\tau), \psi_3(\tau)
\equiv 1$}

In this section, we will prove the following theorem.

{\bf Theorem~2.45}\ \cite{arxiv-5}, \cite{arxiv-10}, \cite{arxiv-11}.\ {\it Suppose that
$\{\phi_j(x)\}_{j=0}^{\infty}$ $(\phi_0(x)=1/\sqrt{T-t})$ is an arbitrary complete orthonormal system of 
functions in the space $L_2([t,T]).$
Then$,$ for the iterated Stra\-to\-no\-vich stochastic integral
of third multiplicity 
$$
{\int\limits_t^{*}}^T
{\int\limits_t^{*}}^{t_3}
{\int\limits_t^{*}}^{t_2}
d{\bf w}_{t_1}^{(i_1)}
d{\bf w}_{t_2}^{(i_2)}d{\bf w}_{t_3}^{(i_3)}\ \ \ (i_1,i_2,i_3=0,1,\ldots,m)
$$
the following expansion 
\begin{equation}
\label{2023novem1}
~~~~~~~{\int\limits_t^{*}}^T
{\int\limits_t^{*}}^{t_3}
{\int\limits_t^{*}}^{t_2}
d{\bf w}_{t_1}^{(i_1)}
d{\bf w}_{t_2}^{(i_2)}d{\bf w}_{t_3}^{(i_3)}=
\hbox{\vtop{\offinterlineskip\halign{
\hfil#\hfil\cr
{\rm l.i.m.}\cr
$\stackrel{}{{}_{p\to \infty}}$\cr
}} }\sum_{j_1,j_2,j_3=0}^{p}
C_{j_3 j_2 j_1}\zeta_{j_1}^{(i_1)}\zeta_{j_2}^{(i_2)}\zeta_{j_3}^{(i_3)}
\end{equation}
that converges in the mean-square sense is valid, where 
$$
C_{j_3 j_2 j_1}=\int\limits_t^T
\phi_{j_3}(t_3)\int\limits_t^{t_3}
\phi_{j_2}(t_2)
\int\limits_t^{t_2}
\phi_{j_1}(t_1)dt_1dt_2dt_3
$$
and
$$
\zeta_{j}^{(i)}=
\int\limits_t^T \phi_{j}(\tau) d{\bf w}_{\tau}^{(i)}
$$ 
are independent standard Gaussian random variables for various 
$i$ or $j$ {\rm (}when $i\ne 0${\rm ),}
${\bf w}_{\tau}^{(i)}$ 
$(i=1,\ldots,m)$ are independent 
standard Wiener processes$,$
${\bf w}_{\tau}^{(0)}=\tau.$}

{\bf Proof.} First, note that under the conditions of Theorem~2.45
the equality
$$
\bar J^{*}[\psi^{(3)}]_{T,t}^{(i_1 i_2 i_3)}=
{\int\limits_t^{*}}^T
{\int\limits_t^{*}}^{t_3}
{\int\limits_t^{*}}^{t_2}
d{\bf w}_{t_1}^{(i_1)}
d{\bf w}_{t_2}^{(i_2)}d{\bf w}_{t_3}^{(i_3)}
$$
is true w.~p.~1 (see Theorem~2.12), where $\bar J^{*}[\psi^{(3)}]_{T,t}^{(i_1 i_2 i_3)}$
is defined by (\ref{dsds9}).

According to Theorem~2.42, we come to the conclusion that 
Theorem~2.45 will be proved if we prove the following
equalities
\begin{equation}
\label{2023novem2}
~~~~~~~~~\lim\limits_{p\to\infty}
\sum\limits_{j_3=0}^{p}
\left(~\sum\limits_{j_1=0}^{p} 
C_{j_3 j_2 j_1}\biggl|_{j_{1}=j_{2}}-
\frac{1}{2} 
C_{j_3 j_2 j_1}\biggl|_{(j_{1} j_{2})\curvearrowright (\cdot),j_{1}=j_{2}}
\biggr.\right)^2=0,
\end{equation}
\begin{equation}
\label{2023novem3}
~~~~~~~~~\lim\limits_{p\to\infty}
\sum\limits_{j_1=0}^{p}
\left(~\sum\limits_{j_3=0}^{p} 
C_{j_3 j_2 j_1}\biggl|_{j_{2}=j_{3}}-
\frac{1}{2} 
C_{j_3 j_2 j_1}\biggl|_{(j_{2} j_{3})\curvearrowright (\cdot), j_{2}=j_{3}}
\biggr.\right)^2=0,
\end{equation}
\begin{equation}
\label{2023novem4}
\lim\limits_{p\to\infty}
\sum\limits_{j_2=0}^{p}
\left(~\sum\limits_{j_1=0}^{p} 
C_{j_3 j_2 j_1}\biggl|_{j_{1}=j_{3}}\right)^2=0.
\end{equation}

Note that using Theorem~2.43 (also see (\ref{start1000})), we can rewrite
the relations (\ref{2023novem2})--(\ref{2023novem4})) in the form
(compare with (\ref{after1600})--(\ref{after1602}))
$$
\lim\limits_{p\to\infty}
\sum\limits_{j_3=0}^{p}
\left(~\sum\limits_{j_1=p+1}^{\infty} 
C_{j_3 j_2 j_1}\biggl|_{j_{1}=j_{2}}\right)^2=0,\ \ \
\lim\limits_{p\to\infty}
\sum\limits_{j_1=0}^{p}
\left(~\sum\limits_{j_3=p+1}^{\infty} 
C_{j_3 j_2 j_1}\biggl|_{j_{2}=j_{3}}\right)^2=0,
$$
$$
\label{2023novem7}
\lim\limits_{p\to\infty}
\sum\limits_{j_2=0}^{p}
\left(~\sum\limits_{j_1=p+1}^{\infty} 
C_{j_3 j_2 j_1}\biggl|_{j_{1}=j_{3}}\right)^2=0.
$$

\vspace{2mm}

Let us prove (\ref{2023novem2}). Using Fubini's Theorem and Parseval's equality, we have
$$
\lim\limits_{p\to\infty}
\sum\limits_{j_3=0}^{p}
\left(~\sum\limits_{j_1=0}^{p} 
C_{j_3 j_2 j_1}\biggl|_{j_{1}=j_{2}}-
\frac{1}{2} 
C_{j_3 j_2 j_1}\biggl|_{(j_{1} j_{2})\curvearrowright (\cdot),j_{1}=j_{2}}
\biggr.\right)^2=
$$
$$
=\lim\limits_{p\to\infty}
\sum\limits_{j_3=0}^{p}
\left(
\frac{1}{2} 
C_{j_3 j_2 j_1}\biggl|_{(j_{1} j_{2})\curvearrowright (\cdot),j_{1}=j_{2}}
\biggr.-
\sum\limits_{j_1=0}^{p} C_{j_3 j_1 j_1}
\right)^2=
$$
$$
=\lim\limits_{p\to\infty}
\sum\limits_{j_3=0}^{p}
\left(\int\limits_t^T 
\phi_{j_3}(\tau)\left(\frac{1}{2}\int\limits_t^{\tau} ds
-
\sum\limits_{j_1=0}^p
\frac{1}{2}\left(\int\limits_t^{\tau}\phi_{j_1}(s)ds\right)^2
\right)d\tau\right)^2\le
$$
$$
\le\lim\limits_{p\to\infty}
\sum\limits_{j_3=0}^{\infty}
\left(\int\limits_t^T 
\phi_{j_3}(\tau)\left(\frac{1}{2}(\tau-t)
-
\sum\limits_{j_1=0}^p
\frac{1}{2}\left(\int\limits_t^{\tau}\phi_{j_1}(s)ds\right)^2
\right)d\tau\right)^2=
$$
\begin{equation}
\label{2023novem8}
=\lim\limits_{p\to\infty}
\int\limits_t^T 
\left(\frac{1}{2}(\tau-t)
-
\sum\limits_{j_1=0}^p
\frac{1}{2}\left(\int\limits_t^{\tau}\phi_{j_1}(s)ds\right)^2
\right)^2 d\tau.
\end{equation}

Applying the Parseval equality, we have
$$
\sum\limits_{j_1=0}^{\infty}
\frac{1}{2}\left(\int\limits_t^{\tau}\phi_{j_1}(s)ds\right)^2
=\sum\limits_{j_1=0}^{\infty}
\frac{1}{2}\left(\int\limits_t^T {\bf 1}_{\{s<\tau\}}\phi_{j_1}(s)ds\right)^2=
$$
\begin{equation}
\label{2023novem10}
=\frac{1}{2}\int\limits_t^T \left({\bf 1}_{\{s<\tau\}}\right)^2 ds=
\frac{1}{2}(\tau-t).
\end{equation}

Moreover, 
\begin{equation}
\label{2023novem11}
~~~~~~~~\left\vert
\frac{1}{2}(\tau-t)
-
\sum\limits_{j_1=0}^p
\frac{1}{2}\left(\int\limits_t^{\tau}\phi_{j_1}(s)ds\right)^2\right\vert\le
\frac{1}{2}(\tau-t)\le \frac{1}{2}(T-t)<\infty.
\end{equation}

Using (\ref{2023novem10}), (\ref{2023novem11}) and
applying Lebesgue's 
Dominated Convergence Theorem in (\ref{2023novem8}), we obtain
the equality (\ref{2023novem2}).

Note that we could  use Dini's Theorem 
instead of Lebesgue's 
Dominated Convergence Theorem. Using
the continuity of the functions $u_p(\tau)$ (see below),
the nondecreasing property of
the functional sequence
$$
u_p(\tau)=
\sum\limits_{j_1=0}^{p}
\frac{1}{2}\left(\int\limits_t^{\tau}\phi_{j_1}(s)ds\right)^2,
$$
and the continuity of the limit function
$u(\tau)=(\tau-t)/2$
according to Dini's 
Theorem,
we have the uniform convergence 
$u_p(\tau)$ to $u(\tau)$ at the interval $[t, T]$.
Then we can swap the limit and integral in (\ref{2023novem8})
and get (\ref{2023novem2}).

Let us prove (\ref{2023novem3}). Using Fubini's Theorem and Parseval's equality, we obtain
$$
\lim\limits_{p\to\infty}
\sum\limits_{j_1=0}^{p}
\left(~\sum\limits_{j_3=0}^{p} 
C_{j_3 j_2 j_1}\biggl|_{j_{2}=j_{3}}-
\frac{1}{2} 
C_{j_3 j_2 j_1}\biggl|_{(j_{2} j_{3})\curvearrowright (\cdot),j_{2}=j_{3}}
\biggr.\right)^2=
$$
$$
=\lim\limits_{p\to\infty}
\sum\limits_{j_1=0}^{p}
\left(
\frac{1}{2} 
C_{j_3 j_2 j_1}\biggl|_{(j_{2} j_{3})\curvearrowright (\cdot),j_{2}=j_{3}}
\biggr.-
\sum\limits_{j_3=0}^{p} C_{j_3 j_3 j_1}
\right)^2=
$$
$$
=\lim\limits_{p\to\infty}
\sum\limits_{j_1=0}^{p}
\left(\frac{1}{2}\int\limits_t^T 
\int\limits_t^{\tau} \phi_{j_1}(s)dsd\tau
-
\sum\limits_{j_3=0}^p
\int\limits_t^T \phi_{j_3}(\theta)\int\limits_t^{\theta} \phi_{j_3}(\tau)
\int\limits_t^{\tau} \phi_{j_1}(s)ds d\tau d\theta\right)^2=
$$
$$
=\lim\limits_{p\to\infty}
\sum\limits_{j_1=0}^{p}
\left(\frac{1}{2}\int\limits_t^T 
\phi_{j_1}(s)(T-s)ds
-
\sum\limits_{j_3=0}^p
\int\limits_t^T \phi_{j_1}(s)\int\limits_{s}^T \phi_{j_3}(\tau)
\int\limits_{\tau}^T \phi_{j_3}(\theta)d\theta d\tau ds\right)^2=
$$
$$
=\lim\limits_{p\to\infty}
\sum\limits_{j_1=0}^{p}
\left(\int\limits_t^T 
\phi_{j_1}(s)\left(\frac{1}{2}(T-s)
-
\sum\limits_{j_3=0}^p
\frac{1}{2}\left(\int\limits_{s}^T \phi_{j_3}(\tau)d\tau\right)^2\right) ds\right)^2\le
$$
$$
\le \lim\limits_{p\to\infty}
\sum\limits_{j_1=0}^{\infty}
\left(\int\limits_t^T 
\phi_{j_1}(s)\left(\frac{1}{2}(T-s)
-
\sum\limits_{j_3=0}^p
\frac{1}{2}\left(\int\limits_{s}^T \phi_{j_3}(\tau)d\tau\right)^2\right) ds\right)^2=
$$
\begin{equation}
\label{2023novem14}
~~~~=\lim\limits_{p\to\infty}
\int\limits_t^T 
\left(\frac{1}{2}(T-s)
-
\sum\limits_{j_3=0}^p
\frac{1}{2}\left(\int\limits_{s}^T \phi_{j_3}(\tau)d\tau\right)^2\right)^2 ds.
\end{equation}

Using the Parseval equality, we get
$$
\sum\limits_{j_3=0}^{\infty}
\frac{1}{2}\left(\int\limits_{s}^T \phi_{j_3}(\tau)d\tau\right)^2=
\sum\limits_{j_3=0}^{\infty}
\frac{1}{2}\left(\int\limits_t^T {\bf 1}_{\{s<\tau\}}\phi_{j_3}(\tau)d\tau\right)^2=
$$
\begin{equation}
\label{2023novem15}
=\frac{1}{2}\int\limits_t^T \left({\bf 1}_{\{s<\tau\}}\right)^2 d\tau=
\frac{1}{2}(T-s).
\end{equation}

Moreover, 
\begin{equation}
\label{2023novem16}
~~~~~~~~\left\vert
\frac{1}{2}(T-s)
-
\sum\limits_{j_3=0}^p
\frac{1}{2}\left(\int\limits_{s}^T \phi_{j_3}(\tau)d\tau\right)^2
\right\vert\le
\frac{1}{2}(T-s)\le \frac{1}{2}(T-t)<\infty.
\end{equation}

Combining (\ref{2023novem14})--(\ref{2023novem16}) and using
the same reasoning as in the proof of (\ref{2023novem2}), we obtain
$$
\lim\limits_{p\to\infty}
\int\limits_t^T 
\left(\frac{1}{2}(T-s)
-
\sum\limits_{j_3=0}^p
\frac{1}{2}\left(\int\limits_{s}^T \phi_{j_3}(\tau)d\tau\right)^2\right)^2 ds=0.
$$

The equality (\ref{2023novem3}) is proved.

Let us prove (\ref{2023novem4}). Applying Fubini's Theorem and Parseval's equality, we have
$$
\lim\limits_{p\to\infty}
\sum\limits_{j_2=0}^{p}
\left(~\sum\limits_{j_1=0}^{p} 
C_{j_1 j_2 j_1}\right)^2=
$$
$$
=\lim\limits_{p\to\infty}
\sum\limits_{j_2=0}^{p}
\left(~\sum\limits_{j_1=0}^{p} 
\int\limits_t^T \phi_{j_1}(\theta)\int\limits_t^{\theta}
\phi_{j_2}(\tau)\int\limits_t^{\tau} \phi_{j_1}(s)ds d\tau d\theta\right)^2=
$$
$$
=\lim\limits_{p\to\infty}
\sum\limits_{j_2=0}^{p}
\left(~\sum\limits_{j_1=0}^{p} 
\int\limits_t^T 
\phi_{j_2}(\tau)\int\limits_t^{\tau} \phi_{j_1}(s)ds
\int\limits_{\tau}^T
\phi_{j_1}(\theta)d\theta d\tau \right)^2\le
$$
$$
\le\lim\limits_{p\to\infty}
\sum\limits_{j_2=0}^{\infty}
\left(
\int\limits_t^T 
\phi_{j_2}(\tau)\sum\limits_{j_1=0}^{p}\int\limits_t^{\tau} \phi_{j_1}(s)ds
\int\limits_{\tau}^T
\phi_{j_1}(\theta)d\theta d\tau \right)^2=
$$
\begin{equation}
\label{2023novem18}
=\lim\limits_{p\to\infty}
\int\limits_t^T 
\left(\sum\limits_{j_1=0}^{p}\int\limits_t^{\tau} \phi_{j_1}(s)ds
\int\limits_{\tau}^T
\phi_{j_1}(\theta)d\theta\right)^2 d\tau.
\end{equation}

\vspace{2mm}

Applying (\ref{dsds14}), we obtain
$$
\left\vert
\sum\limits_{j_1=0}^{p}\int\limits_t^{\tau} \phi_{j_1}(s)ds
\int\limits_{\tau}^T
\phi_{j_1}(\theta)d\theta\right\vert\le
\sum\limits_{j_1=0}^{p}\left\vert
\int\limits_t^{\tau} \phi_{j_1}(s)ds
\int\limits_{\tau}^T
\phi_{j_1}(\theta)d\theta\right\vert\le
$$
\begin{equation}
\label{2023novem19}
\le 
\sum\limits_{j_1=0}^{\infty}\left\vert
\int\limits_t^{\tau} \phi_{j_1}(s)ds
\int\limits_{\tau}^T
\phi_{j_1}(\theta)d\theta\right\vert\le \frac{1}{2}(T-t)<\infty.
\end{equation}

\vspace{2mm}

Using the generalized Parseval equality, we get
$$
\lim\limits_{p\to\infty}
\sum\limits_{j_1=0}^{p}\int\limits_t^{\tau} \phi_{j_1}(s)ds
\int\limits_{\tau}^T
\phi_{j_1}(\theta)d\theta= 
\sum\limits_{j_1=0}^{\infty}\int\limits_t^T {\bf 1}_{\{s<\tau\}}\phi_{j_1}(s)ds
\int\limits_{t}^T
{\bf 1}_{\{s>\tau\}}
\phi_{j_1}(s)ds= 
$$
\begin{equation}
\label{2023novem20}
=
\int\limits_t^T {\bf 1}_{\{s<\tau\}}{\bf 1}_{\{s>\tau\}}ds=0.
\end{equation}

Taking into account (\ref{2023novem19}), (\ref{2023novem20}) and
applying Lebesgue's 
Dominated Convergence Theorem in (\ref{2023novem18}), we obtain
the equality (\ref{2023novem4}). Theorem~2.45 is proved.

\section{Expansion of Iterated Stratonovich Stochastic Integrals
of Multiplicity 4. The Case of an Ar\-bit\-ra\-ry Complete Orthonormal System of 
Functions $(\phi_0(x)=1/\sqrt{T-t})$ in the Space $L_2([t,T])$ and $\psi_1(\tau),\ldots, \psi_4(\tau)
\equiv 1$}

In this section, we will prove the following theorem.

{\bf Theorem~2.46}\ \cite{arxiv-5}, \cite{arxiv-10}, \cite{arxiv-11}.\ {\it Suppose that
$\{\phi_j(x)\}_{j=0}^{\infty}$ $(\phi_0(x)=1/\sqrt{T-t})$ is an arbitrary complete orthonormal system of 
functions in the space $L_2([t,T]).$
Then$,$ for the iterated Stra\-to\-no\-vich stochastic integral
of fourth multiplicity 
$$
J^{*}[\psi^{(4)}]_{T,t}=
{\int\limits_t^{*}}^T
{\int\limits_t^{*}}^{t_4}
{\int\limits_t^{*}}^{t_3}
{\int\limits_t^{*}}^{t_2}
d{\bf w}_{t_1}^{(i_1)}
d{\bf w}_{t_2}^{(i_2)}d{\bf w}_{t_3}^{(i_3)}d{\bf w}_{t_4}^{(i_4)}\ \ \ 
(i_1, i_2, i_3, i_4=0, 1,\ldots,m)
$$
the following 
expansion 
$$
J^{*}[\psi^{(4)}]_{T,t}=
\hbox{\vtop{\offinterlineskip\halign{
\hfil#\hfil\cr
{\rm l.i.m.}\cr
$\stackrel{}{{}_{p\to \infty}}$\cr
}} }
\sum\limits_{j_1, j_2, j_3, j_4=0}^{p}
C_{j_4 j_3 j_2 j_1}\zeta_{j_1}^{(i_1)}\zeta_{j_2}^{(i_2)}\zeta_{j_3}^{(i_3)}
\zeta_{j_4}^{(i_4)}
$$
that converges in the mean-square sense is valid, where 
$$
C_{j_4 j_3 j_2 j_1}=\int\limits_t^T
\phi_{j_4}(t_4)\int\limits_t^{t_4}
\phi_{j_3}(t_3)\int\limits_t^{t_3}
\phi_{j_2}(t_2)\int\limits_t^{t_2}
\phi_{j_1}(t_1)dt_1dt_2dt_3 dt_4
$$
and
$$
\zeta_{j}^{(i)}=
\int\limits_t^T \phi_{j}(\tau) d{\bf w}_{\tau}^{(i)}
$$ 
are independent standard Gaussian random variables for various 
$i$ or $j$ {\rm (}when $i\ne 0${\rm ),}
${\bf w}_{\tau}^{(i)}$ 
$(i=1,\ldots,m)$ are independent 
standard Wiener processes$,$
${\bf w}_{\tau}^{(0)}=\tau.$}

{\bf Proof.} First, note that under the conditions of Theorem~2.46
the equality
$$
\bar J^{*}[\psi^{(4)}]_{T,t}^{(i_1 i_2 i_3 i_4)}=
{\int\limits_t^{*}}^T
{\int\limits_t^{*}}^{t_4}
{\int\limits_t^{*}}^{t_3}
{\int\limits_t^{*}}^{t_2}
d{\bf w}_{t_1}^{(i_1)}
d{\bf w}_{t_2}^{(i_2)}d{\bf w}_{t_3}^{(i_3)}d{\bf w}_{t_4}^{(i_4)}
$$
is valid w.~p.~1 (see Theorem~2.12), where $\bar J^{*}[\psi^{(4)}]_{T,t}^{(i_1 i_2 i_3 i_4)}$
is defined by (\ref{dsds9}).

It is easy to see that Theorem~2.46 will be proved if we prove the following
equalities (see Theorem~2.42)
\begin{equation}
\label{2023novem200}
~~~~~~~\lim\limits_{p\to\infty}
\sum\limits_{j_3,j_4=0}^{p}
\left(~\sum\limits_{j_1=0}^{p} 
C_{j_4 j_3 j_1 j_1}-
\frac{1}{2} 
C_{j_4 j_3 j_1 j_1}\biggl|_{(j_{1} j_{1})\curvearrowright (\cdot)}
\biggr.\right)^2=0,
\end{equation}
\begin{equation}
\label{2023novem201}
\lim\limits_{p\to\infty}
\sum\limits_{j_2,j_4=0}^{p}
\left(~\sum\limits_{j_1=0}^{p} 
C_{j_4 j_1 j_2 j_1}\right)^2=0,
\end{equation}
\begin{equation}
\label{2023novem202}
\lim\limits_{p\to\infty}
\sum\limits_{j_2, j_3=0}^{p}
\left(~\sum\limits_{j_1=0}^{p} 
C_{j_1 j_3 j_2 j_1}\right)^2=0,
\end{equation}
\begin{equation}
\label{2023novem203}
~~~~~~~\lim\limits_{p\to\infty}
\sum\limits_{j_1,j_4=0}^{p}
\left(~\sum\limits_{j_2=0}^{p} 
C_{j_4 j_2 j_2 j_1}-
\frac{1}{2} 
C_{j_4 j_2 j_2 j_1}\biggl|_{(j_{2} j_{2})\curvearrowright (\cdot)}
\biggr.\right)^2=0,
\end{equation}
\begin{equation}
\label{2023novem204}
\lim\limits_{p\to\infty}
\sum\limits_{j_1, j_3=0}^{p}
\left(~\sum\limits_{j_2=0}^{p} 
C_{j_2 j_3 j_2 j_1}\right)^2=0,
\end{equation}
\begin{equation}
\label{2023novem205}
~~~~~~~\lim\limits_{p\to\infty}
\sum\limits_{j_1,j_2=0}^{p}
\left(~\sum\limits_{j_3=0}^{p} 
C_{j_3 j_3 j_2 j_1}-
\frac{1}{2} 
C_{j_3 j_3 j_2 j_1}\biggl|_{(j_{3} j_{3})\curvearrowright (\cdot)}
\biggr.\right)^2=0,
\end{equation}
\begin{equation}
\label{2023novem206}
~~~~~~~~\lim\limits_{p\to\infty}
\sum\limits_{j_1, j_3=0}^{p}
C_{j_3 j_3 j_1 j_1}=\frac{1}{4} 
C_{j_3 j_3 j_1 j_1}\biggl|_{(j_{3} j_{3})\curvearrowright (\cdot)
(j_{1} j_{1})\curvearrowright (\cdot)}=\frac{1}{8}(T-t)^2,
\biggr.
\end{equation}
\begin{equation}
\label{2023novem207}
\lim\limits_{p\to\infty}
\sum\limits_{j_1, j_3=0}^{p}
C_{j_1 j_3 j_3 j_1}=0,
\biggr.
\end{equation}
\begin{equation}
\label{2023novem208}
\lim\limits_{p\to\infty}
\sum\limits_{j_1, j_2=0}^{p}
C_{j_2 j_1 j_2 j_1}=0.
\biggr.
\end{equation}

Let us prove the equalities (\ref{2023novem200})--(\ref{2023novem205}).
Using Fubini's Theorem and Parseval's equality, we obtain
the following relations for the prelimit
expressions on the left-hand sides of (\ref{2023novem200})--(\ref{2023novem205})
$$
\sum\limits_{j_3,j_4=0}^{p}
\left(~\sum\limits_{j_1=0}^{p} 
C_{j_4 j_3 j_1 j_1}-
\frac{1}{2} 
C_{j_4 j_3 j_1 j_1}\biggl|_{(j_{1} j_{1})\curvearrowright (\cdot)}
\biggr.\right)^2=
$$
$$
=\sum\limits_{j_3,j_4=0}^{p}\left(
\frac{1}{2}\int\limits_t^T \phi_{j_4}(t_4)
\int\limits_t^{t_4}\phi_{j_3}(t_3)(t_3-t)dt_3 dt_4-\right.
$$
$$
\left.-\sum\limits_{j_1=0}^p\int\limits_t^T \phi_{j_4}(t_4)
\int\limits_t^{t_4}\phi_{j_3}(t_3)
\int\limits_t^{t_3}\phi_{j_1}(t_2)
\int\limits_t^{t_2}\phi_{j_1}(t_1)dt_1 dt_2 dt_3 dt_4\right)^2=
$$
$$
=\sum\limits_{j_3,j_4=0}^{p}\left(
\int\limits_t^T \phi_{j_4}(t_4)
\int\limits_t^{t_4}\phi_{j_3}(t_3)\Biggl(\frac{1}{2}(t_3-t)-\Biggr.\right.
$$
$$
\left.\Biggl.-\sum\limits_{j_1=0}^p
\int\limits_t^{t_3}\phi_{j_1}(t_2)
\int\limits_t^{t_2}\phi_{j_1}(t_1)dt_1 dt_2\Biggr) dt_3 dt_4\right)^2=
$$
$$
=\sum\limits_{j_3,j_4=0}^{p}\left(
\int\limits_t^T \phi_{j_4}(t_4)
\int\limits_t^{t_4}\phi_{j_3}(t_3)\left(\frac{1}{2}(t_3-t)
-\sum\limits_{j_1=0}^p
\frac{1}{2}\left(\int\limits_t^{t_3}\phi_{j_1}(s)ds\right)^2
\right) dt_3 dt_4\right)^2\le
$$
$$
\le\sum\limits_{j_3,j_4=0}^{\infty}\left(
\int\limits_t^T \phi_{j_4}(t_4)
\int\limits_t^{t_4}\phi_{j_3}(t_3)\left(\frac{1}{2}(t_3-t)
-\sum\limits_{j_1=0}^p
\frac{1}{2}\left(\int\limits_t^{t_3}\phi_{j_1}(s)ds\right)^2
\right) dt_3 dt_4\right)^2=
$$
\begin{equation}
\label{2023novem209}
~~~=
\int\limits_{[t, T]^2}{\bf 1}_{\{t_3<t_4\}}\left(\frac{1}{2}(t_3-t)
-\sum\limits_{j_1=0}^p
\frac{1}{2}\left(\int\limits_t^{t_3}\phi_{j_1}(s)ds\right)^2
\right)^2 dt_3 dt_4,
\end{equation}

\vspace{4mm}

$$
\sum\limits_{j_2,j_4=0}^{p}
\left(~\sum\limits_{j_1=0}^{p} 
C_{j_4 j_1 j_2 j_1}\right)^2=
$$
$$
=
\sum\limits_{j_2,j_4=0}^{p}\left(
\sum\limits_{j_1=0}^p\int\limits_t^T \phi_{j_4}(t_4)
\int\limits_t^{t_4}\phi_{j_1}(t_3)
\int\limits_t^{t_3}\phi_{j_2}(t_2)
\int\limits_t^{t_2}\phi_{j_1}(t_1)dt_1 dt_2 dt_3 dt_4\right)^2=
$$
$$
=
\sum\limits_{j_2,j_4=0}^{p}\left(
\sum\limits_{j_1=0}^p
\int\limits_t^T \phi_{j_4}(t_4)
\int\limits_t^{t_4}\phi_{j_2}(t_2)
\int\limits_t^{t_2}\phi_{j_1}(t_1)dt_1
\int\limits_{t_2}^{t_4}\phi_{j_1}(t_3)dt_3 dt_2 dt_4\right)^2=
$$
$$
=
\sum\limits_{j_2,j_4=0}^{p}\left(
\int\limits_t^T \phi_{j_4}(t_4)
\int\limits_t^{t_4}\phi_{j_2}(t_2)
\sum\limits_{j_1=0}^p
\int\limits_t^{t_2}\phi_{j_1}(t_1)dt_1
\int\limits_{t_2}^{t_4}\phi_{j_1}(t_3)dt_3 dt_2 dt_4\right)^2\le
$$
$$
\le
\sum\limits_{j_2,j_4=0}^{\infty}\left(
\int\limits_t^T \phi_{j_4}(t_4)
\int\limits_t^{t_4}\phi_{j_2}(t_2)
\sum\limits_{j_1=0}^p
\int\limits_t^{t_2}\phi_{j_1}(t_1)dt_1
\int\limits_{t_2}^{t_4}\phi_{j_1}(t_3)dt_3 dt_2 dt_4\right)^2=
$$
\begin{equation}
\label{2023novem210}
~~~=
\int\limits_{[t,T]^2} {\bf 1}_{\{t_2<t_4\}}
\left(\sum\limits_{j_1=0}^p
\int\limits_t^{t_2}\phi_{j_1}(t_1)dt_1
\int\limits_{t_2}^{t_4}\phi_{j_1}(t_3)dt_3\right)^2 dt_2 dt_4,
\end{equation}

\vspace{4mm}

$$
\sum\limits_{j_2,j_3=0}^{p}
\left(~\sum\limits_{j_1=0}^{p} 
C_{j_1 j_3 j_2 j_1}\right)^2=
$$
$$
=
\sum\limits_{j_2,j_3=0}^{p}\left(
\sum\limits_{j_1=0}^p
\int\limits_t^T \phi_{j_1}(t_4)
\int\limits_t^{t_4}\phi_{j_3}(t_3)
\int\limits_t^{t_3}\phi_{j_2}(t_2)
\int\limits_t^{t_2}\phi_{j_1}(t_1)dt_1 dt_2 dt_3 dt_4\right)^2=
$$
$$
=
\sum\limits_{j_2,j_3=0}^{p}\left(
\sum\limits_{j_1=0}^p
\int\limits_t^T \phi_{j_3}(t_3)
\int\limits_t^{t_3}\phi_{j_2}(t_2)
\int\limits_t^{t_2}\phi_{j_1}(t_1)dt_1
\int\limits_{t_3}^{T}\phi_{j_1}(t_4)dt_4 dt_2 dt_3\right)^2=
$$
$$
=
\sum\limits_{j_2,j_3=0}^{p}\left(
\int\limits_t^T \phi_{j_3}(t_3)
\int\limits_t^{t_3}\phi_{j_2}(t_2)
\sum\limits_{j_1=0}^p\int\limits_t^{t_2}\phi_{j_1}(t_1)dt_1
\int\limits_{t_3}^{T}\phi_{j_1}(t_4)dt_4 dt_2 dt_3\right)^2\le
$$
$$
\le
\sum\limits_{j_2,j_3=0}^{\infty}\left(
\int\limits_t^T \phi_{j_3}(t_3)
\int\limits_t^{t_3}\phi_{j_2}(t_2)
\sum\limits_{j_1=0}^p\int\limits_t^{t_2}\phi_{j_1}(t_1)dt_1
\int\limits_{t_3}^{T}\phi_{j_1}(t_4)dt_4 dt_2 dt_3\right)^2=
$$
\begin{equation}
\label{2023novem211}
~~~=
\int\limits_{[t, T]^2} {\bf 1}_{\{t_2<t_3\}}
\left(\sum\limits_{j_1=0}^p\int\limits_t^{t_2}\phi_{j_1}(t_1)dt_1
\int\limits_{t_3}^{T}\phi_{j_1}(t_4)dt_4\right)^2 dt_2 dt_3,
\end{equation}

\vspace{4mm}

$$
\sum\limits_{j_1,j_4=0}^{p}
\left(~\sum\limits_{j_2=0}^{p} 
C_{j_4 j_2 j_2 j_1}-
\frac{1}{2} 
C_{j_4 j_2 j_2 j_1}\biggl|_{(j_{2} j_{2})\curvearrowright (\cdot)}
\biggr.\right)^2=
$$
$$
=\sum\limits_{j_1,j_4=0}^{p}\left(
\frac{1}{2}\int\limits_t^T \phi_{j_4}(t_4)
\int\limits_t^{t_4} \int\limits_t^{t_2}\phi_{j_1}(t_1)dt_1 dt_2 dt_4-\right.
$$
$$
\left.-\sum\limits_{j_2=0}^p\int\limits_t^T \phi_{j_4}(t_4)
\int\limits_t^{t_4}\phi_{j_2}(t_3)
\int\limits_t^{t_3}\phi_{j_2}(t_2)
\int\limits_t^{t_2}\phi_{j_1}(t_1)dt_1 dt_2 dt_3 dt_4\right)^2=
$$
$$
=\sum\limits_{j_1,j_4=0}^{p}\left(
\frac{1}{2}\int\limits_t^T \phi_{j_4}(t_4)
\int\limits_t^{t_4}\phi_{j_1}(t_1)\int\limits_{t_1}^{t_4} dt_2 dt_1 dt_4-\right.
$$
$$
\left.-\sum\limits_{j_2=0}^p\int\limits_t^T \phi_{j_4}(t_4)
\int\limits_t^{t_4}\phi_{j_1}(t_1)
\int\limits_{t_1}^{t_4}\phi_{j_2}(t_2)
\int\limits_{t_2}^{t_4}\phi_{j_2}(t_3)dt_3 dt_2 dt_1 dt_4\right)^2=
$$
$$
=\sum\limits_{j_1,j_4=0}^{p}\left(
\int\limits_t^T \phi_{j_4}(t_4)
\int\limits_t^{t_4}\phi_{j_1}(t_1)\left(\frac{t_4-t_1}{2}
-\sum\limits_{j_2=0}^p
\frac{1}{2}\left(\int\limits_{t_1}^{t_4}\phi_{j_2}(s)ds\right)^2
\right) dt_1 dt_4\right)^2\le
$$
$$
\le\sum\limits_{j_1,j_4=0}^{\infty}\left(
\int\limits_t^T \phi_{j_4}(t_4)
\int\limits_t^{t_4}\phi_{j_1}(t_1)\left(\frac{t_4-t_1}{2}
-\sum\limits_{j_2=0}^p
\frac{1}{2}\left(\int\limits_{t_1}^{t_4}\phi_{j_2}(s)ds\right)^2
\right) dt_1 dt_4\right)^2=
$$
\begin{equation}
\label{2023novem212}
~~~=
\int\limits_{[t,T]^2} {\bf 1}_{\{t_1<t_4\}}\left(\frac{1}{2}(t_4-t_1)
-\sum\limits_{j_2=0}^p
\frac{1}{2}\left(\int\limits_{t_1}^{t_4}\phi_{j_2}(s)ds\right)^2
\right)^2 dt_1 dt_4,
\end{equation}

\vspace{4mm}

$$
\sum\limits_{j_1,j_3=0}^{p}
\left(~\sum\limits_{j_2=0}^{p} 
C_{j_2 j_3 j_2 j_1}\right)^2=
$$
$$
=
\sum\limits_{j_1,j_3=0}^{p}\left(
\sum\limits_{j_2=0}^p
\int\limits_t^T \phi_{j_2}(t_4)
\int\limits_t^{t_4}\phi_{j_3}(t_3)
\int\limits_t^{t_3}\phi_{j_2}(t_2)
\int\limits_t^{t_2}\phi_{j_1}(t_1)dt_1 dt_2 dt_3 dt_4\right)^2=
$$
$$
=
\sum\limits_{j_1,j_3=0}^{p}\left(
\sum\limits_{j_2=0}^p
\int\limits_t^T \phi_{j_3}(t_3)
\int\limits_t^{t_3}\phi_{j_2}(t_2)
\int\limits_t^{t_2}\phi_{j_1}(t_1)dt_1 dt_2
\int\limits_{t_3}^{T}\phi_{j_2}(t_4)dt_4 dt_3 \right)^2=
$$
$$
=
\sum\limits_{j_1,j_3=0}^{p}\left(
\sum\limits_{j_2=0}^p
\int\limits_t^T \phi_{j_3}(t_3)
\int\limits_t^{t_3}\phi_{j_1}(t_1)
\int\limits_{t_1}^{t_3}\phi_{j_2}(t_2)dt_2
\int\limits_{t_3}^{T}\phi_{j_2}(t_4)dt_4 dt_1 dt_3 \right)^2=
$$
$$
=
\sum\limits_{j_1,j_3=0}^{p}\left(
\int\limits_t^T \phi_{j_3}(t_3)
\int\limits_t^{t_3}\phi_{j_1}(t_1)
\sum\limits_{j_2=0}^p
\int\limits_{t_1}^{t_3}\phi_{j_2}(t_2)dt_2
\int\limits_{t_3}^{T}\phi_{j_2}(t_4)dt_4 dt_1 dt_3 \right)^2\le
$$
$$
\le
\sum\limits_{j_1,j_3=0}^{\infty}\left(
\int\limits_t^T \phi_{j_3}(t_3)
\int\limits_t^{t_3}\phi_{j_1}(t_1)
\sum\limits_{j_2=0}^p
\int\limits_{t_1}^{t_3}\phi_{j_2}(t_2)dt_2
\int\limits_{t_3}^{T}\phi_{j_2}(t_4)dt_4 dt_1 dt_3 \right)^2=
$$
\begin{equation}
\label{2023novem213}
~~~=
\int\limits_{[t,T]^2}{\bf 1}_{\{t_1<t_3\}}
\left(\sum\limits_{j_2=0}^p
\int\limits_{t_1}^{t_3}\phi_{j_2}(t_2)dt_2
\int\limits_{t_3}^{T}\phi_{j_2}(t_4)dt_4\right)^2 dt_1 dt_3,
\end{equation}

\vspace{4mm}

$$
\sum\limits_{j_1,j_2=0}^{p}
\left(~\sum\limits_{j_3=0}^{p} 
C_{j_3 j_3 j_2 j_1}-
\frac{1}{2} 
C_{j_3 j_3 j_2 j_1}\biggl|_{(j_{3} j_{3})\curvearrowright (\cdot)}
\biggr.\right)^2=
$$
$$
=\sum\limits_{j_1,j_2=0}^{p}\left(
\frac{1}{2}\int\limits_t^T
\int\limits_t^{t_3}\phi_{j_2}(t_2)\int\limits_t^{t_2}\phi_{j_1}(t_1)dt_1 dt_2 dt_3-\right.
$$
$$
\left.-\sum\limits_{j_3=0}^p\int\limits_t^T \phi_{j_3}(t_4)
\int\limits_t^{t_4}\phi_{j_3}(t_3)
\int\limits_t^{t_3}\phi_{j_2}(t_2)
\int\limits_t^{t_2}\phi_{j_1}(t_1)dt_1 dt_2 dt_3 dt_4\right)^2=
$$
$$
=\sum\limits_{j_1,j_2=0}^{p}\left(
\frac{1}{2}\int\limits_t^T\phi_{j_1}(t_1)
\int\limits_{t_1}^T\phi_{j_2}(t_2)\int\limits_{t_2}^T dt_3 dt_2 dt_1-\right.
$$
$$
\left.-\sum\limits_{j_3=0}^p\int\limits_t^T \phi_{j_1}(t_1)
\int\limits_{t_1}^T \phi_{j_2}(t_2)
\int\limits_{t_2}^T\phi_{j_3}(t_3)
\int\limits_{t_3}^T\phi_{j_3}(t_4)dt_4 dt_3 dt_2 dt_1\right)^2=
$$
$$
=\sum\limits_{j_1,j_2=0}^{p}\left(
\int\limits_t^T \phi_{j_1}(t_1)
\int\limits_{t_1}^T\phi_{j_2}(t_2)\left(\frac{T-t_2}{2}
-\sum\limits_{j_3=0}^p
\frac{1}{2}\left(\int\limits_{t_2}^T\phi_{j_3}(s)ds\right)^2
\right) dt_2 dt_1\right)^2\le
$$
$$
\le\sum\limits_{j_1,j_2=0}^{\infty}\left(
\int\limits_t^T \phi_{j_1}(t_1)
\int\limits_{t_1}^T\phi_{j_2}(t_2)\left(\frac{T-t_2}{2}
-\sum\limits_{j_3=0}^p
\frac{1}{2}\left(\int\limits_{t_2}^T\phi_{j_3}(s)ds\right)^2
\right) dt_2 dt_1\right)^2=
$$
\begin{equation}
\label{2023novem214}
~~~=
\int\limits_{[t, T]^2}{\bf 1}_{\{t_1<t_2\}}\left(\frac{1}{2}(T-t_2)
-\sum\limits_{j_3=0}^p
\frac{1}{2}\left(\int\limits_{t_2}^T\phi_{j_3}(s)ds\right)^2
\right)^2 dt_2 dt_1.
\end{equation}

\vspace{2mm}

Using Parseval's equality, generalized Parseval's equality and 
Lebesgue's Dominated Convergence Theorem, 
as well as applying the same reasoning as in the proof of Theorem~2.45, 
we obtain that the right-hand sides of (\ref{2023novem209})--(\ref{2023novem214}) 
tend to zero when $p\to\infty.$
The equalities (\ref{2023novem200})--(\ref{2023novem205}) are proved.

Let us prove the equalities (\ref{2023novem206})--(\ref{2023novem208}).
We will use our idea from Sect.~2.14. More precisely, we consider
the following analogue of the equality (\ref{sixsix40})
$$
C_{j_4 j_3 j_2 j_1}+C_{j_1 j_2 j_3 j_4}=
C_{j_4}C_{j_3 j_2 j_1}-C_{j_3 j_4}C_{j_2 j_1}+
$$
\begin{equation}
\label{2023novem215}
+C_{j_2 j_3 j_4}C_{j_1}.
\end{equation}

\vspace{2mm}

Using Fubini's Theorem, we have
$$
C_{j_4 j_3 j_2 j_1}=
$$
$$
=\int\limits_t^T\phi_{j_4}(t_4)\int\limits_t^{t_4}\phi_{j_3}(t_3)
\int\limits_t^{t_3}\phi_{j_2}(t_2)
\int\limits_t^{t_2}\phi_{j_1}(t_1)dt_1 dt_2 dt_3 dt_4=
$$
$$
=\int\limits_t^T\phi_{j_4}(t_4)\int\limits_t^{T}\phi_{j_3}(t_3)
\int\limits_t^{t_3}\phi_{j_2}(t_2)
\int\limits_t^{t_2}\phi_{j_1}(t_1)dt_1 dt_2 dt_3 dt_4-
$$
$$
-\int\limits_t^T\phi_{j_4}(t_4)\int\limits_{t_4}^T\phi_{j_3}(t_3)
\int\limits_t^{t_3}\phi_{j_2}(t_2)
\int\limits_t^{t_2}\phi_{j_1}(t_1)dt_1 dt_2 dt_3 dt_4=
$$
$$
=C_{j_4}C_{j_3 j_2 j_1}-
$$

\vspace{-5mm}
$$
-\int\limits_t^T\phi_{j_4}(t_4)\int\limits_{t_4}^T\phi_{j_3}(t_3)
\int\limits_t^{T}\phi_{j_2}(t_2)
\int\limits_t^{t_2}\phi_{j_1}(t_1)dt_1 dt_2 dt_3 dt_4+
$$
$$
+\int\limits_t^T\phi_{j_4}(t_4)\int\limits_{t_4}^T\phi_{j_3}(t_3)
\int\limits_{t_3}^{T}\phi_{j_2}(t_2)
\int\limits_t^{t_2}\phi_{j_1}(t_1)dt_1 dt_2 dt_3 dt_4=
$$

\vspace{-3mm}
$$
=C_{j_4}C_{j_3 j_2 j_1}-C_{j_3 j_4}C_{j_2 j_1}+
$$

\vspace{-5mm}
$$
+\int\limits_t^T\phi_{j_4}(t_4)\int\limits_{t_4}^T\phi_{j_3}(t_3)
\int\limits_{t_3}^{T}\phi_{j_2}(t_2)\int\limits_{t}^{T}\phi_{j_1}(t_1)
dt_1 dt_2 dt_3 dt_4-
$$
$$
-\int\limits_t^T\phi_{j_4}(t_4)\int\limits_{t_4}^T\phi_{j_3}(t_3)
\int\limits_{t_3}^{T}\phi_{j_2}(t_2)\int\limits_{t_2}^{T}\phi_{j_1}(t_1)
dt_1 dt_2 dt_3 dt_4=
$$

\vspace{-2mm}
\begin{equation}
\label{2023novem216}
=C_{j_4}C_{j_3 j_2 j_1}-C_{j_3 j_4}C_{j_2 j_1}+C_{j_2 j_3 j_4}C_{j_1}-
C_{j_1 j_2 j_3 j_4}.
\end{equation}

\vspace{5mm}

The equality (\ref{2023novem216}) completes the proof of the relation 
(\ref{2023novem215}).

Let us prove (\ref{2023novem206}). Substitute $j_4=j_3,$ $j_2=j_1$ into
(\ref{2023novem215})
$$
C_{j_3 j_3 j_1 j_1}+C_{j_1 j_1 j_3 j_3}=
C_{j_3}C_{j_3 j_1 j_1}-C_{j_3 j_3}C_{j_1 j_1}+
$$
\begin{equation}
\label{2023novem217}
+C_{j_1 j_3 j_3}C_{j_1}.
\end{equation}

\vspace{2mm}

From (\ref{2023novem217}) we obtain
$$
\sum\limits_{j_1,j_3=0}^p \bigl(C_{j_3 j_3 j_1 j_1}+C_{j_1 j_1 j_3 j_3}\bigr)=
\sum\limits_{j_1,j_3=0}^p C_{j_3}C_{j_3 j_1 j_1}-\sum\limits_{j_1,j_3=0}^p C_{j_3 j_3}C_{j_1 j_1}+
$$
$$
+\sum\limits_{j_1,j_3=0}^p C_{j_1 j_3 j_3}C_{j_1}.
$$

Then
\begin{equation}
\label{2023novem218}
~~~~~~~~~~~2\sum\limits_{j_1,j_3=0}^p C_{j_3 j_3 j_1 j_1}=
2\sum\limits_{j_1,j_3=0}^p C_{j_3}C_{j_3 j_1 j_1}-\left(\sum\limits_{j_1=0}^p C_{j_1 j_1}\right)^2.
\end{equation}

\vspace{2mm}

From (\ref{2023novem218}) we get
$$
\sum\limits_{j_1,j_3=0}^p C_{j_3 j_3 j_1 j_1}=
\sum\limits_{j_1,j_3=0}^p C_{j_3}C_{j_3 j_1 j_1}-\frac{1}{2}
\left(\sum\limits_{j_1=0}^p C_{j_1 j_1}\right)^2=
$$
\begin{equation}
\label{2023novem219}
=\sum\limits_{j_1,j_3=0}^p C_{j_3}C_{j_3 j_1 j_1}-\frac{1}{2}
\left(\sum\limits_{j_1=0}^p \frac{1}{2}\bigl(C_{j_1}\bigr)^2\right)^2=
\sum\limits_{j_1,j_3=0}^p C_{j_3}C_{j_3 j_1 j_1}-\frac{1}{8}
\left(\sum\limits_{j_1=0}^p \bigl(C_{j_1}\bigr)^2\right)^2.
\end{equation}

\vspace{1mm}

Recall that $\phi_0(\tau)=1/\sqrt{T-t}.$ Then 
\begin{equation}
\label{2023novem220}
C_j=\int\limits_t^T \phi_j(\tau)d\tau=
\left\{
\begin{matrix}
\sqrt{T-t} &\hbox{if}\ j=0\cr\cr
0 &\hbox{if}\ j\ne 0
\end{matrix}.\right.
\end{equation}

Combining (\ref{2023novem219}), (\ref{2023novem220}) and using Fubini's Theorem, we obtain
$$
\sum\limits_{j_1,j_3=0}^p C_{j_3 j_3 j_1 j_1}=
\sqrt{T-t}\sum\limits_{j_1=0}^p C_{0 j_1 j_1}-\frac{1}{8}(T-t)^2=
$$
$$
=\sum\limits_{j_1=0}^p \int\limits_t^T \int\limits_t^{t_3}
\phi_{j_1}(t_2)
\int\limits_t^{t_2}
\phi_{j_1}(t_1)dt_1 dt_2 dt_3-\frac{1}{8}(T-t)^2=
$$
$$
=\sum\limits_{j_1=0}^p \int\limits_t^T \phi_{j_1}(t_1) \int\limits_{t_1}^T
\phi_{j_1}(t_2)
\int\limits_{t_2}^T
dt_3 dt_2 dt_1-\frac{1}{8}(T-t)^2=
$$
$$
=\sum\limits_{j_1=0}^p \int\limits_t^T \phi_{j_1}(t_1) \int\limits_{t_1}^T
\phi_{j_1}(t_2)
(T-t_2)dt_2 dt_1-\frac{1}{8}(T-t)^2=
$$
\begin{equation}
\label{2023novem221}
~~~~~~~~=\sum\limits_{j_1=0}^p \int\limits_t^T 
\phi_{j_1}(t_2)
(T-t_2)\int\limits_{t}^{t_2}
\phi_{j_1}(t_1) dt_1 dt_2-\frac{1}{8}(T-t)^2.
\end{equation}

Finally applying (\ref{start1000}) and (\ref{2023novem221}), we have
$$
\lim\limits_{p\to\infty}\sum\limits_{j_1,j_3=0}^p C_{j_3 j_3 j_1 j_1}=
\frac{1}{2}\int\limits_t^T
(T-t_2)dt_2-\frac{1}{8}(T-t)^2=\frac{1}{8}(T-t)^2.
$$

The equality (\ref{2023novem206}) is proved.

Let us prove (\ref{2023novem207}). Substitute $j_4=j_1,$ $j_2=j_3$ into
(\ref{2023novem215})
$$
C_{j_1 j_3 j_3 j_1}+C_{j_1 j_3 j_3 j_1}=
C_{j_1}C_{j_3 j_3 j_1}-C_{j_3 j_1}C_{j_3 j_1}+
$$
\begin{equation}
\label{2023novem222}
+C_{j_3 j_3 j_1}C_{j_1}.
\end{equation}

\vspace{2mm}

Using (\ref{2023novem222}), we get
\begin{equation}
\label{2023novem223}
~~~~~~~~~2\sum\limits_{j_1,j_3=0}^p C_{j_1 j_3 j_3 j_1}=
2\sum\limits_{j_1,j_3=0}^p C_{j_1}C_{j_3 j_3 j_1}-\sum\limits_{j_1,j_3=0}^p \bigl(C_{j_3 j_1}\bigr)^2.
\end{equation}

Then applying (\ref{2023novem223}), (\ref{2023novem220}), Parseval's equality,
and (\ref{start1000}), we obtain
$$
\lim\limits_{p\to\infty}\sum\limits_{j_1,j_3=0}^p C_{j_1 j_3 j_3 j_1}=
\lim\limits_{p\to\infty}\sum\limits_{j_1,j_3=0}^p C_{j_1}C_{j_3 j_3 j_1}-
\frac{1}{2}\lim\limits_{p\to\infty}\sum\limits_{j_1,j_3=0}^p \bigl(C_{j_3 j_1}\bigr)^2=
$$
$$
=\sqrt{T-t}\sum\limits_{j_3=0}^{\infty} C_{j_3 j_3 0}-
\frac{1}{2}\sum\limits_{j_1,j_3=0}^{\infty}
\left(\int\limits_t^T \phi_{j_3}(t_2)\int\limits_t^{t_2}
\phi_{j_1}(t_1)dt_1 dt_2
\right)^2=
$$
$$
=\sum\limits_{j_3=0}^{\infty} \int\limits_t^T \phi_{j_3}(t_3)\int\limits_t^{t_3}
\phi_{j_3}(t_2)\int\limits_t^{t_2}dt_1 dt_2 dt_3-
$$
$$
-
\frac{1}{2}\sum\limits_{j_1,j_3=0}^{\infty}
\left(\int\limits_{[t,T]^2}{\bf 1}_{\{t_1<t_2\}} 
\phi_{j_1}(t_1) \phi_{j_3}(t_2) dt_1 dt_2
\right)^2=
$$
$$
=\sum\limits_{j_3=0}^{\infty} \int\limits_t^T \phi_{j_3}(t_3)\int\limits_t^{t_3}
\phi_{j_3}(t_2)(t_2-t) dt_2 dt_3-
\frac{1}{2}
\int\limits_{[t,T]^2}\left({\bf 1}_{\{t_1<t_2\}}\right)^2 dt_1 dt_2
=
$$
$$
=\frac{1}{2}\int\limits_t^T 
(t_2-t) dt_2-
\frac{1}{2}
\int\limits_t^T \int\limits_t^{t_2} dt_1 dt_2=0.
$$

The equality (\ref{2023novem207}) is proved.

Let us prove (\ref{2023novem208}). Substitute $j_3=j_1,$ $j_4=j_2$ into
(\ref{2023novem215})
$$
C_{j_2 j_1 j_2 j_1}+C_{j_1 j_2 j_1 j_2}=
C_{j_2}C_{j_1 j_2 j_1}-C_{j_1 j_2}C_{j_2 j_1}+
$$
\begin{equation}
\label{2023novem224}
+C_{j_2 j_1 j_2}C_{j_1}.
\end{equation}

Then
$$
\sum\limits_{j_1,j_2=0}^p \bigl(C_{j_2 j_1 j_2 j_1}+C_{j_1 j_2 j_1 j_2}\bigr)=
\sum\limits_{j_1,j_2=0}^p \bigl(C_{j_2}C_{j_1 j_2 j_1}+C_{j_2 j_1 j_2}C_{j_1}\bigr)-
$$
\begin{equation}
\label{2023novem225}
-\sum\limits_{j_1,j_2=0}^p C_{j_1 j_2}C_{j_2 j_1}.
\end{equation}

From (\ref{2023novem225}) we have
$$
2\sum\limits_{j_1,j_2=0}^p C_{j_2 j_1 j_2 j_1}=
2\sum\limits_{j_1,j_2=0}^p C_{j_1}C_{j_2 j_1 j_2}-
$$
$$
-\sum\limits_{j_1,j_2=0}^p \frac{1}{2} \left( \bigl(C_{j_1 j_2} + C_{j_2 j_1}\bigr)^2 
-\bigl(C_{j_1 j_2}\bigr)^2 - \bigl(C_{j_2 j_1}\bigr)^2\right)=
$$
$$
=
2\sum\limits_{j_1,j_2=0}^p C_{j_1}C_{j_2 j_1 j_2}
-\frac{1}{2}\sum\limits_{j_1,j_2=0}^p \bigl(C_{j_1 j_2} + C_{j_2 j_1}\bigr)^2 
+
$$
\begin{equation}
\label{2023novem226}
+\sum\limits_{j_1,j_2=0}^p \bigl(C_{j_2 j_1}\bigr)^2.
\end{equation}

\vspace{2mm}

Using Fubini's Theorem, we obtain (also see (\ref{nahod0}))
\begin{equation}
\label{2023novem227}
C_{j_1 j_2} + C_{j_2 j_1}=C_{j_1}C_{j_2}.
\end{equation}

Applying (\ref{2023novem226}), (\ref{2023novem227}), (\ref{2023novem220}),
Fubini's Theorem, Parseval's equality,
and (\ref{start1000}), we get
$$
\lim\limits_{p\to \infty}\sum\limits_{j_1,j_2=0}^p C_{j_2 j_1 j_2 j_1}=
\lim\limits_{p\to \infty}\sum\limits_{j_1,j_2=0}^p C_{j_1}C_{j_2 j_1 j_2}
-\frac{1}{4}\lim\limits_{p\to \infty}\sum\limits_{j_1,j_2=0}^p \bigl(C_{j_1 j_2} + C_{j_2 j_1}\bigr)^2 
+
$$
$$
+\frac{1}{2}\lim\limits_{p\to \infty}\sum\limits_{j_1,j_2=0}^p \bigl(C_{j_2 j_1}\bigr)^2=
$$

\vspace{2mm}
$$
=\sqrt{T-t}\sum\limits_{j_2=0}^{\infty} C_{j_2 0 j_2}
-\frac{1}{4}\sum\limits_{j_1,j_2=0}^{\infty} \bigl(C_{j_1}C_{j_2}\bigr)^2 
+\frac{1}{2}\sum\limits_{j_1,j_2=0}^{\infty} \bigl(C_{j_2 j_1}\bigr)^2=
$$

\vspace{-3mm}
$$
=\sum\limits_{j_2=0}^{\infty} \int\limits_t^T \phi_{j_2}(t_3)\int\limits_t^{t_3}\int\limits_t^{t_2}
\phi_{j_2}(t_1)dt_1 dt_2 dt_3
-\frac{1}{4}(T-t)^2 
+\frac{1}{2}\int\limits_{[t,T]^2} \bigl({\bf 1}_{\{t_1<t_2\}}\bigr)^2 dt_1 dt_2=
$$

\vspace{-3mm}
$$
=\sum\limits_{j_2=0}^{\infty} \int\limits_t^T \phi_{j_2}(t_3)\int\limits_t^{t_3}
\phi_{j_2}(t_1)\int\limits_{t_1}^{t_3} dt_2 dt_1 dt_3=
$$

\vspace{-3mm}
$$
=\sum\limits_{j_2=0}^{\infty} \int\limits_t^T \phi_{j_2}(t_3)(t_3-t)\int\limits_t^{t_3}
\phi_{j_2}(t_1)dt_1 dt_3+
\sum\limits_{j_2=0}^{\infty} \int\limits_t^T \phi_{j_2}(t_3)\int\limits_t^{t_3}
\phi_{j_2}(t_1)(t-t_1)dt_1 dt_3=
$$

\vspace{-3mm}
$$
=\frac{1}{2}\int\limits_t^T (t_3-t)dt_3+
\frac{1}{2}\int\limits_t^T(t-t_3)dt_3=0.
$$

\vspace{2mm}

The equality (\ref{2023novem208}) is proved. The equalities (\ref{2023novem200})--(\ref{2023novem208})
are proved. Theorem~2.46 is proved.

\section{Condition $\phi_0(x)=1/\sqrt{T-t}$ in Theorems~2.45 and 2.46
can be Omitted}

In this section, we will show that the 
condition $\phi_0(x)=1/\sqrt{T-t}$ in Theorems~2.45 and 2.46
can be omitted.

{\bf Theorem~2.47}\ \cite{arxiv-5}, \cite{arxiv-10}, \cite{arxiv-11}.\ {\it Suppose that
$\{\phi_j(x)\}_{j=0}^{\infty}$ is an arbitrary complete orthonormal system of 
functions in the space $L_2([t,T]).$
Then$,$ for the iterated Stra\-to\-no\-vich stochastic integral
of third multiplicity 
$$
{\int\limits_t^{*}}^T
{\int\limits_t^{*}}^{t_3}
{\int\limits_t^{*}}^{t_2}
d{\bf w}_{t_1}^{(i_1)}
d{\bf w}_{t_2}^{(i_2)}d{\bf w}_{t_3}^{(i_3)}\ \ \ (i_1,i_2,i_3=0,1,\ldots,m)
$$
the following expansion 
\begin{equation}
\label{2023novem1xyz}
~~~~{\int\limits_t^{*}}^T
{\int\limits_t^{*}}^{t_3}
{\int\limits_t^{*}}^{t_2}
d{\bf w}_{t_1}^{(i_1)}
d{\bf w}_{t_2}^{(i_2)}d{\bf w}_{t_3}^{(i_3)}=
\hbox{\vtop{\offinterlineskip\halign{
\hfil#\hfil\cr
{\rm l.i.m.}\cr
$\stackrel{}{{}_{p\to \infty}}$\cr
}} }\sum_{j_1,j_2,j_3=0}^{p}
C_{j_3 j_2 j_1}\zeta_{j_1}^{(i_1)}\zeta_{j_2}^{(i_2)}\zeta_{j_3}^{(i_3)}
\end{equation}
that converges in the mean-square sense is valid, where 
$$
C_{j_3 j_2 j_1}=\int\limits_t^T
\phi_{j_3}(t_3)\int\limits_t^{t_3}
\phi_{j_2}(t_2)
\int\limits_t^{t_2}
\phi_{j_1}(t_1)dt_1dt_2dt_3
$$
and
$$
\zeta_{j}^{(i)}=
\int\limits_t^T \phi_{j}(\tau) d{\bf w}_{\tau}^{(i)}
$$ 
are independent standard Gaussian random variables for various 
$i$ or $j$ {\rm (}in the case when $i\ne 0${\rm ),}
${\bf w}_{\tau}^{(i)}$ 
$(i=1,\ldots,m)$ are independent 
standard Wiener processes$,$ 
${\bf w}_{\tau}^{(0)}=\tau.$}

{\bf Proof.} Analyzing the proof of Theorems~2.42 and 2.45 (also see the derivation of (\ref{after906})
and (\ref{dydy11})),
we notice that Theorem~2.47 will be proved if we prove that
\begin{equation}
\label{febr1}
~~~~~~~~~\int\limits_t^T \int\limits_t^{t_3}
dt_2 d{\bf w}_{t_3}^{(i_3)}=
\hbox{\vtop{\offinterlineskip\halign{
\hfil#\hfil\cr
{\rm l.i.m.}\cr
$\stackrel{}{{}_{p\to \infty}}$\cr
}} }\sum_{j_3=0}^{p}\int\limits_t^T \phi_{j_3}(t_3)\int\limits_t^{t_3}
dt_2 dt_3\ \zeta_{j_3}^{(i_3)},
\end{equation}
\begin{equation}
\label{febr2}
~~~~~~~~~\int\limits_t^T \int\limits_t^{t_2}
d{\bf w}_{t_1}^{(i_1)}dt_2=
\hbox{\vtop{\offinterlineskip\halign{
\hfil#\hfil\cr
{\rm l.i.m.}\cr
$\stackrel{}{{}_{p\to \infty}}$\cr
}} }\sum_{j_1=0}^{p}\int\limits_t^T \int\limits_t^{t_2}\phi_{j_1}(t_1)dt_1 dt_2\
\zeta_{j_1}^{(i_1)}.
\end{equation}

The equality (\ref{febr1}) immediately follows from Theorem~1.16 (see (\ref{razzar1}) for $k=1$).

Let us prove (\ref{febr2}).
Using the theorem on replacement of the integration order in iterated 
It\^{o} stochastic integrals (see Theorem 3.1 and (\ref{febr3})) or 
the It\^{o} formula, Theorem~1.16 (see (\ref{razzar1}) for $k=1$) 
and Fubini's Theorem, we obtain~w.~p.~1
$$
\int\limits_t^T \int\limits_t^{t_2}
d{\bf w}_{t_1}^{(i_1)}dt_2=\int\limits_t^T \int\limits_{t_1}^{T}
dt_2 d{\bf w}_{t_1}^{(i_1)}=
\hbox{\vtop{\offinterlineskip\halign{
\hfil#\hfil\cr
{\rm l.i.m.}\cr
$\stackrel{}{{}_{p\to \infty}}$\cr
}} }\sum_{j_1=0}^{p}\int\limits_t^T \phi_{j_1}(t_1)\int\limits_{t_1}^{T}dt_2 dt_1\
\zeta_{j_1}^{(i_1)}=
$$
$$
=
\hbox{\vtop{\offinterlineskip\halign{
\hfil#\hfil\cr
{\rm l.i.m.}\cr
$\stackrel{}{{}_{p\to \infty}}$\cr
}} }\sum_{j_1=0}^{p}\int\limits_t^T \int\limits_{t}^{t_2}\phi_{j_1}(t_1)dt_1 dt_2\
\zeta_{j_1}^{(i_1)}.
$$

The equality (\ref{febr2}) is proved. Theorem~2.47 is proved.

Let us develop this approach and prove the following
generalization of Theorem~2.46.

{\bf Theorem~2.48}\ \cite{arxiv-5}, \cite{arxiv-10}, \cite{arxiv-11}.\ {\it Suppose that
$\{\phi_j(x)\}_{j=0}^{\infty}$ is an arbitrary complete orthonormal system of 
functions in the space $L_2([t,T]).$
Then$,$ for the iterated Stra\-to\-no\-vich stochastic integral
of fourth multiplicity 
$$
J^{*}[\psi^{(4)}]_{T,t}=
{\int\limits_t^{*}}^T
{\int\limits_t^{*}}^{t_4}
{\int\limits_t^{*}}^{t_3}
{\int\limits_t^{*}}^{t_2}
d{\bf w}_{t_1}^{(i_1)}
d{\bf w}_{t_2}^{(i_2)}d{\bf w}_{t_3}^{(i_3)}d{\bf w}_{t_4}^{(i_4)}\ \ \ 
(i_1, i_2, i_3, i_4=0, 1,\ldots,m)
$$
the following 
expansion 
$$
J^{*}[\psi^{(4)}]_{T,t}=
\hbox{\vtop{\offinterlineskip\halign{
\hfil#\hfil\cr
{\rm l.i.m.}\cr
$\stackrel{}{{}_{p\to \infty}}$\cr
}} }
\sum\limits_{j_1, j_2, j_3, j_4=0}^{p}
C_{j_4 j_3 j_2 j_1}\zeta_{j_1}^{(i_1)}\zeta_{j_2}^{(i_2)}\zeta_{j_3}^{(i_3)}
\zeta_{j_4}^{(i_4)}
$$
that converges in the mean-square sense is valid, where 
$$
C_{j_4 j_3 j_2 j_1}=\int\limits_t^T
\phi_{j_4}(t_4)\int\limits_t^{t_4}
\phi_{j_3}(t_3)\int\limits_t^{t_3}
\phi_{j_2}(t_2)\int\limits_t^{t_2}
\phi_{j_1}(t_1)dt_1dt_2dt_3 dt_4
$$
and
$$
\zeta_{j}^{(i)}=
\int\limits_t^T \phi_{j}(\tau) d{\bf w}_{\tau}^{(i)}
$$ 
are independent standard Gaussian random variables for various 
$i$ or $j$ {\rm (}when $i\ne 0${\rm ),}
${\bf w}_{\tau}^{(i)}$ 
$(i=1,\ldots,m)$ are independent 
standard Wiener processes$,$
${\bf w}_{\tau}^{(0)}=\tau.$}

{\bf Proof.}\ Considering the proof of Theorems~2.42 and 2.46 (also see the derivation of (\ref{after906})
and (\ref{dydy11})),
we conclude that Theorem~2.48 will be proved if we prove that
under the conditions of Theorem~2.48 the following equalities
$$
\int\limits_t^T \int\limits_t^{t_3} \int\limits_t^{t_2}
dt_1 d{\bf w}_{t_2}^{(i_2)}d{\bf w}_{t_3}^{(i_3)}=
\hbox{\vtop{\offinterlineskip\halign{
\hfil#\hfil\cr
{\rm l.i.m.}\cr
$\stackrel{}{{}_{p\to \infty}}$\cr
}} }\sum_{j_2,j_3=0}^{p}\int\limits_t^T \phi_{j_3}(t_3)\int\limits_t^{t_3} 
\phi_{j_2}(t_2)\int\limits_t^{t_2}
dt_1 dt_2 dt_3\times
$$
\begin{equation}
\label{febr5}
\times
J'[\phi_{j_2}\phi_{j_3}]_{T,t}^{(i_2 i_3)},
\end{equation}
$$
\int\limits_t^T \int\limits_t^{t_3} \int\limits_t^{t_2}
d{\bf w}_{t_1}^{(i_1)} dt_2 d{\bf w}_{t_3}^{(i_3)}=
\hbox{\vtop{\offinterlineskip\halign{
\hfil#\hfil\cr
{\rm l.i.m.}\cr
$\stackrel{}{{}_{p\to \infty}}$\cr
}} }\sum_{j_1,j_3=0}^{p}\int\limits_t^T \phi_{j_3}(t_3)\int\limits_t^{t_3} 
\int\limits_t^{t_2}
\phi_{j_1}(t_1)dt_1 
dt_2 dt_3\times
$$
\begin{equation}
\label{febr6}
\times 
J'[\phi_{j_1}\phi_{j_3}]_{T,t}^{(i_1 i_3)},
\end{equation}

\vspace{-6mm}
$$
\int\limits_t^T \int\limits_t^{t_3} \int\limits_t^{t_2}
d{\bf w}_{t_1}^{(i_1)}d{\bf w}_{t_2}^{(i_2)}dt_3=
\hbox{\vtop{\offinterlineskip\halign{
\hfil#\hfil\cr
{\rm l.i.m.}\cr
$\stackrel{}{{}_{p\to \infty}}$\cr
}} }\sum_{j_1,j_2=0}^{p}\int\limits_t^T \int\limits_t^{t_3} 
\phi_{j_2}(t_2)\int\limits_t^{t_2}
\phi_{j_1}(t_1)dt_1
dt_2 dt_3\times
$$
\begin{equation}
\label{febr7}
\times 
J'[\phi_{j_1}\phi_{j_2}]_{T,t}^{(i_1 i_2)},
\end{equation}

\vspace{-5mm}
\begin{equation}
\label{march10}
~~~~~~~~\lim\limits_{p\to\infty}
\sum\limits_{j_1, j_3=0}^{p}
C_{j_3 j_3 j_1 j_1}=\frac{1}{4} 
C_{j_3 j_3 j_1 j_1}\biggl|_{(j_{3} j_{3})\curvearrowright (\cdot)
(j_{1} j_{1})\curvearrowright (\cdot)}=\frac{1}{8}(T-t)^2,
\biggr.
\end{equation}
\begin{equation}
\label{march11}
\lim\limits_{p\to\infty}
\sum\limits_{j_1, j_2=0}^{p}
C_{j_2 j_1 j_2 j_1}=0
\biggr.,
\end{equation}
\begin{equation}
\label{march12}
\lim\limits_{p\to\infty}
\sum\limits_{j_1, j_3=0}^{p}
C_{j_1 j_3 j_3 j_1}=0
\biggr.
\end{equation}

\vspace{2mm}
\noindent 
holds, where we use (\ref{febr5000}), i.e. 
\begin{equation}
\label{febr9}
~~~~~~~~~J[\psi^{(k)}]_{T,t}^{(i_1\ldots i_k)}=
\hbox{\vtop{\offinterlineskip\halign{
\hfil#\hfil\cr
{\rm l.i.m.}\cr
$\stackrel{}{{}_{p_1,\ldots,p_k\to \infty}}$\cr
}} }
\sum\limits_{j_1=0}^{p_1}\ldots
\sum\limits_{j_k=0}^{p_k}
C_{j_k\ldots j_1}
J'[\phi_{j_1}\ldots \phi_{j_k}]_{T,t}^{(i_1\ldots i_k)},
\end{equation}

\noindent
where
$J'[\phi_{j_1}\ldots \phi_{j_k}]_{T,t}^{(i_1\ldots i_k)}$
is the multiple Wiener stochastic integral defined by (\ref{WiI}).

Moreover, for $k=4, r=2, g_1=1, g_2=2, g_3=3, g_4=4$ we can write 
(see the derivation of (\ref{after906}))
$$
\hbox{\vtop{\offinterlineskip\halign{
\hfil#\hfil\cr
{\rm l.i.m.}\cr
$\stackrel{}{{}_{p\to \infty}}$\cr
}} }
\sum\limits_{\stackrel{j_1,\ldots,j_q,\ldots,j_k=0}{{}_{q\ne g_1, g_2,\ldots, g_{2r-1}, g_{2r}}}}^p
\frac{1}{2^r}
C_{j_k \ldots j_1}\biggl|_{(j_{g_2} j_{g_1})\curvearrowright (\cdot)
\ldots (j_{g_{2r}} j_{g_{2r-1}})\curvearrowright (\cdot),
j_{g_{{}_{1}}}=~j_{g_{{}_{2}}},\ldots, j_{g_{{}_{2r-1}}}=~j_{g_{{}_{2r}}}}\biggr.
\times 
$$

\vspace{2mm}
$$
\times
\prod\limits_{s=1}^r
{\bf 1}_{\{i_{g_{{}_{2s-1}}}=~i_{g_{{}_{2s}}}\ne 0\}}
J'[\phi_{j_{q_1}}\ldots \phi_{j_{q_{k-2r}}}]_{T,t}^{(i_{q_1}\ldots i_{q_{k-2r}})}=
$$
$$
=\frac{1}{4}{\bf 1}_{\{i_{1}=i_{2}\ne 0\}}{\bf 1}_{\{i_{3}=i_{4}\ne 0\}}
C_{j_3 j_3 j_1 j_1}\biggl|_{(j_{3} j_{3})\curvearrowright (\cdot)
(j_{1} j_{1})\curvearrowright (\cdot)}\biggr.=
$$

\vspace{1mm}
$$
=\frac{1}{4}{\bf 1}_{\{i_{1}=i_{2}\ne 0\}}{\bf 1}_{\{i_{3}=i_{4}\ne 0\}}
\int\limits_t^T\int\limits_t^{t_2} dt_1 dt_2=
{\bf 1}_{\{i_{1}=i_{2}\ne 0\}}{\bf 1}_{\{i_{3}=i_{4}\ne 0\}}
\frac{(T-t)^2}{8},
$$

\vspace{2mm}
\noindent
where
$J'[\phi_{j_{q_1}}\ldots \phi_{j_{q_{k-2r}}}]_{T,t}^{(i_{q_1}\ldots i_{q_{k-2r}})}
\stackrel{\sf def}{=}1$
for $k=2r$.

The equality (\ref{febr5}) immediately follows from Theorem~1.16 (see (\ref{febr5000}) 
or (\ref{febr9}) for $k=2$).

Let us prove (\ref{febr7}).
Using the theorem on replacement of the integration order in iterated 
It\^{o} stochastic integrals (see Theorem 3.1 and (\ref{febr10})) or 
the It\^{o} formula, Theorem~1.16 (see (\ref{febr5000}) or (\ref{febr9}) for $k=2$) 
and Fubini's Theorem, we get~w.~p.~1
$$
\int\limits_t^T \int\limits_t^{t_3} \int\limits_t^{t_2}
d{\bf w}_{t_1}^{(i_1)}d{\bf w}_{t_2}^{(i_2)}dt_3=
\int\limits_t^T (T-t_2)\int\limits_t^{t_2} 
d{\bf w}_{t_1}^{(i_1)}d{\bf w}_{t_2}^{(i_2)}=
$$
$$
=\hbox{\vtop{\offinterlineskip\halign{
\hfil#\hfil\cr
{\rm l.i.m.}\cr
$\stackrel{}{{}_{p\to \infty}}$\cr
}} }\sum_{j_1,j_2=0}^{p}\int\limits_t^T (T-t_2)
\phi_{j_2}(t_2)\int\limits_t^{t_2}
\phi_{j_1}(t_1)dt_1 dt_2
J'[\phi_{j_1}\phi_{j_2}]_{T,t}^{(i_1 i_2)}=
$$
$$
=\hbox{\vtop{\offinterlineskip\halign{
\hfil#\hfil\cr
{\rm l.i.m.}\cr
$\stackrel{}{{}_{p\to \infty}}$\cr
}} }\sum_{j_1,j_2=0}^{p}\int\limits_t^T 
\phi_{j_1}(t_1)\int\limits_{t_1}^T
\phi_{j_2}(t_2)(T-t_2)dt_2 dt_1
J'[\phi_{j_1}\phi_{j_2}]_{T,t}^{(i_1 i_2)}=
$$
$$
=\hbox{\vtop{\offinterlineskip\halign{
\hfil#\hfil\cr
{\rm l.i.m.}\cr
$\stackrel{}{{}_{p\to \infty}}$\cr
}} }\sum_{j_1,j_2=0}^{p}\int\limits_t^T 
\phi_{j_1}(t_1)\int\limits_{t_1}^T
\phi_{j_2}(t_2)\int\limits_{t_2}^T dt_3 dt_2 dt_1
J'[\phi_{j_1}\phi_{j_2}]_{T,t}^{(i_1 i_2)}=
$$
$$
=
\hbox{\vtop{\offinterlineskip\halign{
\hfil#\hfil\cr
{\rm l.i.m.}\cr
$\stackrel{}{{}_{p\to \infty}}$\cr
}} }\sum_{j_1,j_2=0}^{p}\int\limits_t^T \int\limits_t^{t_3} 
\phi_{j_2}(t_2)\int\limits_t^{t_2}
\phi_{j_1}(t_1)dt_1
dt_2 dt_3
J'[\phi_{j_1}\phi_{j_2}]_{T,t}^{(i_1 i_2)}.
$$

\vspace{2mm}

The equality (\ref{febr7}) is proved. To prove (\ref{febr6}) we will use the above arguments
((\ref{febr12}) (see below) also directly follows from the It\^{o} formula)
$$
\int\limits_t^T \int\limits_t^{t_3} \int\limits_t^{t_2}
d{\bf w}_{t_1}^{(i_1)} dt_2 d{\bf w}_{t_3}^{(i_3)}=\hbox{[by Theorems~3.1, 3.3]}=
\int\limits_t^T \int\limits_t^{t_3} 
d{\bf w}_{t_1}^{(i_1)}\int\limits_{t_1}^{t_3} dt_2 d{\bf w}_{t_3}^{(i_3)}=
$$
$$
=\int\limits_t^T \int\limits_t^{t_3} (t_3-t_1)
d{\bf w}_{t_1}^{(i_1)}d{\bf w}_{t_3}^{(i_3)}=
$$
\begin{equation}
\label{febr12}
~~~~~~~~~=\int\limits_t^T (t_3-t)\int\limits_t^{t_3} 
d{\bf w}_{t_1}^{(i_1)}d{\bf w}_{t_3}^{(i_3)}-
\int\limits_t^T \int\limits_t^{t_3} (t_1-t)
d{\bf w}_{t_1}^{(i_1)}d{\bf w}_{t_3}^{(i_3)}=
\end{equation}
$$
=\hbox{\vtop{\offinterlineskip\halign{
\hfil#\hfil\cr
{\rm l.i.m.}\cr
$\stackrel{}{{}_{p\to \infty}}$\cr
}} }\sum_{j_1,j_3=0}^{p}\int\limits_t^T (t_3-t)\phi_{j_3}(t_3)\int\limits_t^{t_3} 
\phi_{j_1}(t_1)dt_1 
dt_3
J'[\phi_{j_1}\phi_{j_3}]_{T,t}^{(i_1 i_3)}-
$$
$$
-\hbox{\vtop{\offinterlineskip\halign{
\hfil#\hfil\cr
{\rm l.i.m.}\cr
$\stackrel{}{{}_{p\to \infty}}$\cr
}} }\sum_{j_1,j_3=0}^{p}\int\limits_t^T \phi_{j_3}(t_3)\int\limits_t^{t_3} 
(t_1-t)\phi_{j_1}(t_1)dt_1 
dt_3
J'[\phi_{j_1}\phi_{j_3}]_{T,t}^{(i_1 i_3)}=
$$
$$
=\hbox{\vtop{\offinterlineskip\halign{
\hfil#\hfil\cr
{\rm l.i.m.}\cr
$\stackrel{}{{}_{p\to \infty}}$\cr
}} }\sum_{j_1,j_3=0}^{p}\left(\int\limits_t^T (t_3-t)\phi_{j_3}(t_3)\int\limits_t^{t_3} 
\phi_{j_1}(t_1)dt_1 
dt_3-\right.
$$
$$
\left.-\int\limits_t^T \phi_{j_3}(t_3)\int\limits_t^{t_3} 
(t_1-t)\phi_{j_1}(t_1)dt_1 
dt_3\right)
J'[\phi_{j_1}\phi_{j_3}]_{T,t}^{(i_1 i_3)}=
$$
$$
=\hbox{\vtop{\offinterlineskip\halign{
\hfil#\hfil\cr
{\rm l.i.m.}\cr
$\stackrel{}{{}_{p\to \infty}}$\cr
}} }\sum_{j_1,j_3=0}^{p}\int\limits_t^T \phi_{j_3}(t_3)\int\limits_t^{t_3} (t_3-t+t-t_1)
\phi_{j_1}(t_1)dt_1 dt_3
J'[\phi_{j_1}\phi_{j_3}]_{T,t}^{(i_1 i_3)}=
$$
$$
=\hbox{\vtop{\offinterlineskip\halign{
\hfil#\hfil\cr
{\rm l.i.m.}\cr
$\stackrel{}{{}_{p\to \infty}}$\cr
}} }\sum_{j_1,j_3=0}^{p}\int\limits_t^T \phi_{j_3}(t_3)\int\limits_t^{t_3} 
\phi_{j_1}(t_1)\int\limits_{t_1}^{t_3}dt_2
dt_1 dt_3
J'[\phi_{j_1}\phi_{j_3}]_{T,t}^{(i_1 i_3)}=
$$
$$
=
\hbox{\vtop{\offinterlineskip\halign{
\hfil#\hfil\cr
{\rm l.i.m.}\cr
$\stackrel{}{{}_{p\to \infty}}$\cr
}} }\sum_{j_1,j_3=0}^{p}\int\limits_t^T \phi_{j_3}(t_3)\int\limits_t^{t_3} 
\int\limits_t^{t_2}
\phi_{j_1}(t_1)dt_1 
dt_2 dt_3
J'[\phi_{j_1}\phi_{j_3}]_{T,t}^{(i_1 i_3)}.
$$

\vspace{2mm}

The equality (\ref{febr6}) is proved. Let us prove (\ref{march10})--(\ref{march12}).
Using (\ref{2023novem219}), we obtain
\begin{equation}
\label{march13}
~~~~~~\sum\limits_{j_1,j_3=0}^p C_{j_3 j_3 j_1 j_1}=
\sum\limits_{j_1,j_3=0}^p C_{j_3}C_{j_3 j_1 j_1}-\frac{1}{8}
\left(\sum\limits_{j_1=0}^p \bigl(C_{j_1}\bigr)^2\right)^2.
\end{equation}

Applying Parseval's equality, we have
\begin{equation}
\label{march14}
\lim\limits_{p\to\infty}
\sum\limits_{j_1=0}^p \bigl(C_{j_1}\bigr)^2=\int\limits_t^T 1^2 d\tau=T-t.
\end{equation}

Combining (\ref{march13}) and (\ref{march14}), we get
\begin{equation}
\label{march15}
~~~~~~~~~~~~\lim\limits_{p\to\infty}\sum\limits_{j_1,j_3=0}^p C_{j_3 j_3 j_1 j_1}=
\lim\limits_{p\to\infty}\sum\limits_{j_1,j_3=0}^p C_{j_3}C_{j_3 j_1 j_1}-\frac{(T-t)^2}{8}.
\end{equation}

Further, we have
$$
\lim\limits_{p\to\infty}\sum\limits_{j_1,j_3=0}^p C_{j_3}C_{j_3 j_1 j_1}=
$$
\begin{equation}
\label{march16}
=
\frac{1}{2}\lim\limits_{p\to\infty}
\sum\limits_{j_3=0}^{p}C_{j_3}C_{j_3 j_1 j_1}\biggl|_{(j_{1} j_{1})\curvearrowright (\cdot)}-
\lim\limits_{p\to\infty}
\sum\limits_{j_3=0}^{p}C_{j_3}\left(\frac{1}{2}
C_{j_3 j_1 j_1}\biggl|_{(j_{1} j_{1})\curvearrowright (\cdot)}-
\sum\limits_{j_1=0}^p C_{j_3 j_1 j_1}\right).
\end{equation}

Applying the generalized Parseval equality, we obtain
$$
\lim\limits_{p\to\infty}
\sum\limits_{j_3=0}^{p}C_{j_3}C_{j_3 j_1 j_1}\biggl|_{(j_{1} j_{1})\curvearrowright (\cdot)}=
\lim\limits_{p\to\infty}
\sum\limits_{j_3=0}^{p}\int\limits_t^T \phi_{j_3}(\tau)d\tau
\int\limits_t^T \phi_{j_3}(\tau)\int\limits_t^{\tau}d\theta d\tau
=
$$
\begin{equation}
\label{march17}
=
\int\limits_t^T 1\cdot \int\limits_t^{\tau}d\theta d\tau=
\frac{(T-t)^2}{2}.
\end{equation}

From (\ref{march16}) and (\ref{march17}) we have
$$
\lim\limits_{p\to\infty}\sum\limits_{j_1,j_3=0}^p C_{j_3}C_{j_3 j_1 j_1}=
$$
\begin{equation}
\label{march18}
~~~~~~~~=
\frac{(T-t)^2}{4}-
\lim\limits_{p\to\infty}
\sum\limits_{j_3=0}^{p}C_{j_3}\left(\frac{1}{2}
C_{j_3 j_1 j_1}\biggl|_{(j_{1} j_{1})\curvearrowright (\cdot)}-
\sum\limits_{j_1=0}^p C_{j_3 j_1 j_1}\right).
\end{equation}

\vspace{3mm}

Combining (\ref{march15}) and (\ref{march18}), we obtain
\begin{equation}
\label{march19}
\lim\limits_{p\to\infty}\sum\limits_{j_1,j_3=0}^p C_{j_3 j_3 j_1 j_1}=
\frac{(T-t)^2}{8}-
\lim\limits_{p\to\infty}
\sum\limits_{j_3=0}^{p}C_{j_3}\left(\frac{1}{2}
C_{j_3 j_1 j_1}\biggl|_{(j_{1} j_{1})\curvearrowright (\cdot)}-
\sum\limits_{j_1=0}^p C_{j_3 j_1 j_1}\right).
\end{equation}

Due to the inequality of Cauchy--Bunyakovsky
and (\ref{2023novem2}), (\ref{march14}), we get
$$
\lim\limits_{p\to\infty}
\left(\sum\limits_{j_3=0}^{p}C_{j_3}\left(\frac{1}{2}
C_{j_3 j_1 j_1}\biggl|_{(j_{1} j_{1})\curvearrowright (\cdot)}-
\sum\limits_{j_1=0}^p C_{j_3 j_1 j_1}\right)\right)^2\le
$$
$$
\le \lim\limits_{p\to\infty}
\sum\limits_{j_3=0}^{p}\left(C_{j_3}\right)^2\
\sum\limits_{j_3=0}^{p}
\left(\frac{1}{2}
C_{j_3 j_1 j_1}\biggl|_{(j_{1} j_{1})\curvearrowright (\cdot)}-
\sum\limits_{j_1=0}^p C_{j_3 j_1 j_1}\right)^2\le
$$
$$
\le \lim\limits_{p\to\infty}
\sum\limits_{j_3=0}^{\infty}\left(C_{j_3}\right)^2\
\sum\limits_{j_3=0}^{p}
\left(\frac{1}{2}
C_{j_3 j_1 j_1}\biggl|_{(j_{1} j_{1})\curvearrowright (\cdot)}-
\sum\limits_{j_1=0}^p C_{j_3 j_1 j_1}\right)^2=
$$
\begin{equation}
\label{march20}
~~~~~~~~~~=(T-t)\lim\limits_{p\to\infty}
\sum\limits_{j_3=0}^{p}
\left(\frac{1}{2}
C_{j_3 j_1 j_1}\biggl|_{(j_{1} j_{1})\curvearrowright (\cdot)}-
\sum\limits_{j_1=0}^p C_{j_3 j_1 j_1}\right)^2=0.
\end{equation}
         
\vspace{3mm}

Taking into account (\ref{march19}) and (\ref{march20}), we obtain (\ref{march10}).
It is not difficult to see that by analogy with (\ref{march10}) we get
\begin{equation}
\label{march20a}
~~~~~~~~\lim\limits_{p\to\infty}
\sum\limits_{j_1, j_3=0}^{p}
C_{j_3 j_3 j_1 j_1}(s)=\frac{1}{8}(s-t)^2,
\biggr.
\end{equation}

\noindent
where $s\in (t, T]$ and
\begin{equation}
\label{march21a}
~~~~~~~~C_{j_4 j_3 j_2 j_1}(s)=\int\limits_t^s
\phi_{j_4}(t_4)\int\limits_t^{t_4}
\phi_{j_3}(t_3)\int\limits_t^{t_3}
\phi_{j_2}(t_2)\int\limits_t^{t_2}
\phi_{j_1}(t_1)dt_1dt_2dt_3 dt_4.
\end{equation}

Let us prove (\ref{march11}). Using (\ref{2023novem225}), we have
\begin{equation}
\label{march21}
~~~~~~~~~\sum\limits_{j_1,j_2=0}^p C_{j_2 j_1 j_2 j_1}=
\sum\limits_{j_1,j_2=0}^p C_{j_2}C_{j_1 j_2 j_1}
-\frac{1}{2}\sum\limits_{j_1,j_2=0}^p C_{j_1 j_2}C_{j_2 j_1}.
\end{equation}

Fubini's Theorem and the generalized Parseval equality give
$$
\lim_{p\to\infty}\sum\limits_{j_1,j_2=0}^p C_{j_1 j_2}C_{j_2 j_1}=
$$
$$
=
\lim_{p\to\infty}
\sum\limits_{j_1,j_2=0}^p \int\limits_t^T \phi_{j_2}(t_2)\int\limits_{t_2}^{T}\phi_{j_1}(t_1)dt_1 dt_2
\int\limits_t^T \phi_{j_2}(t_2)\int\limits_t^{t_2}\phi_{j_1}(t_1)dt_1 dt_2=
$$
$$
=\lim_{p\to\infty}\sum\limits_{j_1,j_2=0}^p 
\int\limits_{[t,T]^2}{\bf 1}_{\{t_2<t_1\}}\phi_{j_1}(t_1)\phi_{j_2}(t_2)dt_1 dt_2
\int\limits_{[t,T]^2} {\bf 1}_{\{t_1<t_2\}} \phi_{j_1}(t_1)\phi_{j_2}(t_2)dt_1 dt_2=
$$
\begin{equation}
\label{march22}
=\int\limits_{[t,T]^2}{\bf 1}_{\{t_2<t_1\}}{\bf 1}_{\{t_1<t_2\}}dt_1 dt_2=0.
\end{equation}

The equalities (\ref{march21}) and (\ref{march22}) imply the relation
\begin{equation}
\label{march23}
\lim_{p\to\infty}\sum\limits_{j_1,j_2=0}^p C_{j_2 j_1 j_2 j_1}=
\lim_{p\to\infty}\sum\limits_{j_1,j_2=0}^p C_{j_2}C_{j_1 j_2 j_1}.
\end{equation}

Further, we have (see the derivation of (\ref{march20}))
$$
\lim_{p\to\infty}\left(\sum\limits_{j_2=0}^p C_{j_2} \sum\limits_{j_1=0}^p  C_{j_1 j_2 j_1}\right)^2\le
\lim_{p\to\infty}\sum\limits_{j_2=0}^p \left(C_{j_2}\right)^2 
\sum\limits_{j_2=0}^p\left(\sum\limits_{j_1=0}^p  C_{j_1 j_2 j_1}\right)^2\le
$$
\begin{equation}
\label{march24}
\le \lim_{p\to\infty}\sum\limits_{j_2=0}^{\infty} \left(C_{j_2}\right)^2 
\sum\limits_{j_2=0}^p\left(\sum\limits_{j_1=0}^p  C_{j_1 j_2 j_1}\right)^2=
(T-t)\lim_{p\to\infty}\sum\limits_{j_2=0}^p\left(\sum\limits_{j_1=0}^p  C_{j_1 j_2 j_1}\right)^2=0,
\end{equation}
where (\ref{march24}) follows from (\ref{2023novem4}).

The relations (\ref{march23}) and (\ref{march24}) complete the proof of (\ref{march11}).
By analogy with the above reasoning, we obviously get
\begin{equation}
\label{march25}
\lim\limits_{p\to\infty}
\sum\limits_{j_1, j_2=0}^{p}
C_{j_2 j_1 j_2 j_1}(s)=0
\biggr.,
\end{equation}

\noindent
where $s\in (t, T]$ and $C_{j_2 j_1 j_2 j_1}(s)$ is defined by (\ref{march21a}).

Let us prove (\ref{march12}). Using (\ref{2023novem223}), we obtain
\begin{equation}
\label{march26}
~~~~~~~~~~\sum\limits_{j_1,j_3=0}^p C_{j_1 j_3 j_3 j_1}=
\sum\limits_{j_1,j_3=0}^p C_{j_1}C_{j_3 j_3 j_1}-\frac{1}{2}
\sum\limits_{j_1,j_3=0}^p \bigl(C_{j_3 j_1}\bigr)^2.
\end{equation}

Parseval's equality gives
$$
\lim\limits_{p\to\infty}\sum\limits_{j_1,j_3=0}^p \bigl(C_{j_3 j_1}\bigr)^2=
\lim\limits_{p\to\infty}\sum\limits_{j_1,j_3=0}^p
\left(~\int\limits_{[t,T]^2}{\bf 1}_{\{t_1<t_2\}}
\phi_{j_1}(t_1)\phi_{j_3}(t_2)dt_1 dt_2\right)^2=
$$
\begin{equation}
\label{march27}
=\int\limits_{[t,T]^2}\left({\bf 1}_{\{t_1<t_2\}}\right)^2 dt_1 dt_2=
\frac{(T-t)^2}{2}.
\end{equation}

Combining (\ref{march26}) and (\ref{march27}), we have
\begin{equation}
\label{march28}
~~~~~~~~~~~~\lim\limits_{p\to\infty}\sum\limits_{j_1,j_3=0}^p C_{j_1 j_3 j_3 j_1}=
\lim\limits_{p\to\infty}\sum\limits_{j_1,j_3=0}^p C_{j_1}C_{j_3 j_3 j_1}-\frac{(T-t)^2}{4}.
\end{equation}

Further, we have
$$
\lim\limits_{p\to\infty}\sum\limits_{j_1,j_3=0}^p C_{j_1}C_{j_3 j_3 j_1}=
$$
\begin{equation}
\label{march29}
=
\frac{1}{2}\lim\limits_{p\to\infty}
\sum\limits_{j_1=0}^{p}C_{j_1}C_{j_3 j_3 j_1}\biggl|_{(j_{3} j_{3})\curvearrowright (\cdot)}-
\lim\limits_{p\to\infty}
\sum\limits_{j_1=0}^{p}C_{j_1}\left(\frac{1}{2}
C_{j_3 j_3 j_1}\biggl|_{(j_{3} j_{3})\curvearrowright (\cdot)}-
\sum\limits_{j_3=0}^p C_{j_3 j_3 j_1}\right).
\end{equation}

Applying Fubini's Theorem and the generalized Parseval equality, we obtain
$$
\lim\limits_{p\to\infty}
\sum\limits_{j_1=0}^{p}C_{j_1}C_{j_3 j_3 j_1}\biggl|_{(j_{3} j_{3})\curvearrowright (\cdot)}=
\lim\limits_{p\to\infty}
\sum\limits_{j_1=0}^{p}\int\limits_t^T \phi_{j_1}(\tau)d\tau
\int\limits_t^T \int\limits_t^{t_2}\phi_{j_1}(\tau)d\tau dt_2
=
$$
\begin{equation}
\label{march30}
=\lim\limits_{p\to\infty}
\sum\limits_{j_1=0}^{p}\int\limits_t^T \phi_{j_1}(\tau)d\tau
\int\limits_t^T \phi_{j_1}(\tau) \int\limits_{\tau}^{T} dt_2 d\tau
=
\int\limits_t^T 1\cdot \int\limits_{\tau}^T d\theta d\tau=
\frac{(T-t)^2}{2}.
\end{equation}

From (\ref{march29}) and (\ref{march30}) we have
$$
\lim\limits_{p\to\infty}\sum\limits_{j_1,j_3=0}^p C_{j_1}C_{j_3 j_3 j_1}=
$$
\begin{equation}
\label{march31}
~~~~~~~~~~~=
\frac{(T-t)^2}{4}-
\lim\limits_{p\to\infty}
\sum\limits_{j_1=0}^{p}C_{j_1}\left(\frac{1}{2}
C_{j_3 j_3 j_1}\biggl|_{(j_{3} j_{3})\curvearrowright (\cdot)}-
\sum\limits_{j_3=0}^p C_{j_3 j_3 j_1}\right).
\end{equation}

\vspace{3mm}

Combining (\ref{march28}) and (\ref{march31}), we obtain
\begin{equation}
\label{march32}
~~~\lim\limits_{p\to\infty}\sum\limits_{j_1,j_3=0}^p C_{j_1 j_3 j_3 j_1}=
-
\lim\limits_{p\to\infty}
\sum\limits_{j_1=0}^{p}C_{j_1}\left(\frac{1}{2}
C_{j_3 j_3 j_1}\biggl|_{(j_{3} j_{3})\curvearrowright (\cdot)}-
\sum\limits_{j_3=0}^p C_{j_3 j_3 j_1}\right).
\end{equation}

Due to the inequality of Cauchy--Bunyakovsky
and (\ref{2023novem3}), (\ref{march14}), we get
$$
\lim\limits_{p\to\infty}
\left(\sum\limits_{j_1=0}^{p}C_{j_1}\left(\frac{1}{2}
C_{j_3 j_3 j_1}\biggl|_{(j_{3} j_{3})\curvearrowright (\cdot)}-
\sum\limits_{j_3=0}^p C_{j_3 j_3 j_1}\right)\right)^2\le
$$
$$
\le \lim\limits_{p\to\infty}
\sum\limits_{j_1=0}^{p}\left(C_{j_1}\right)^2\
\sum\limits_{j_1=0}^{p}
\left(\frac{1}{2}
C_{j_3 j_3 j_1}\biggl|_{(j_{3} j_{3})\curvearrowright (\cdot)}-
\sum\limits_{j_3=0}^p C_{j_3 j_3 j_1}\right)^2\le
$$
$$
\le \lim\limits_{p\to\infty}
\sum\limits_{j_1=0}^{\infty}\left(C_{j_1}\right)^2\
\sum\limits_{j_1=0}^{p}
\left(\frac{1}{2}
C_{j_3 j_3 j_1}\biggl|_{(j_{3} j_{3})\curvearrowright (\cdot)}-
\sum\limits_{j_3=0}^p C_{j_3 j_3 j_1}\right)^2=
$$
\begin{equation}
\label{march33}
~~~~~~~~~~=(T-t)\lim\limits_{p\to\infty}
\sum\limits_{j_1=0}^{p}
\left(\frac{1}{2}
C_{j_3 j_3 j_1}\biggl|_{(j_{3} j_{3})\curvearrowright (\cdot)}-
\sum\limits_{j_3=0}^p C_{j_3 j_3 j_1}\right)^2=0.
\end{equation}

The relations (\ref{march32}) and (\ref{march33}) complete the proof of (\ref{march12}).
By analogy with the above reasoning, we obviously have
\begin{equation}
\label{march34}
\lim\limits_{p\to\infty}
\sum\limits_{j_1, j_3=0}^{p}
C_{j_1 j_3 j_3 j_1}(s)=0
\biggr.,
\end{equation}

\noindent
where $s\in (t, T]$ and $C_{j_1 j_3 j_3 j_1}(s)$ is defined by (\ref{march21a}).

The equalities (\ref{febr5})--(\ref{march12}) are proved. Theorem~2.48 is proved.

Note that the equalities (\ref{march25}) and (\ref{march34}) can be proved by another way.
Using Fubini's Theorem, we obtain
\begin{equation}
\label{march35}
C_{j_2 j_1 j_2 j_1}(s)=\frac{1}{2}\left(C_{j_2 j_1}(s)\right)^2-2C_{j_2 j_2 j_1 j_1}(s),
\end{equation}
\begin{equation}
\label{march36}
~~~~~~~~\sum\limits_{(j_1,j_2,j_3,j_4)}C_{j_4 j_3 j_2 j_1}(s)=C_{j_1}(s)C_{j_2}(s)C_{j_3}(s)C_{j_4}(s),
\end{equation}

\noindent
where $s\in (t, T],$
$$
\sum\limits_{(j_1,j_2,j_3,j_4)}
$$
means the sum with respect to all
possible permutations
$(j_1,j_2,j_3,j_4)$ and
$$
C_{j_k \ldots j_1}(s)=\int\limits_t^s
\phi_{j_k}(t_k)\ldots
\int\limits_t^{t_2}
\phi_{j_1}(t_1)dt_1\ldots dt_k \ \ \ (k=1,\ldots,4).
$$

Taking into account (\ref{march20a}), (\ref{march27})
(for $s$ instead of $T$), (\ref{march35}), we get
$$
\lim\limits_{p\to\infty}\sum\limits_{j_1,j_2=0}^{p}
C_{j_2 j_1 j_2 j_1}(s)=\frac{1}{2}
\lim\limits_{p\to\infty}\sum\limits_{j_1,j_2=0}^{p}\left(C_{j_2 j_1}(s)\right)^2-2
\lim\limits_{p\to\infty}\sum\limits_{j_1,j_2=0}^{p}C_{j_2 j_2 j_1 j_1}(s)=
$$
$$
=\frac{1}{2}\cdot \frac{(s-t)^2}{2}-2\cdot \frac{(s-t)^2}{8}=0.
$$

\vspace{2mm}

The equality (\ref{march25}) is proved.
Let us substitute $j_2=j_1$ and $j_4=j_3$ into (\ref{march36}). Then we obtain
$$
4\biggl(C_{j_3 j_3 j_1 j_1}(s)+C_{j_1 j_1 j_3 j_3}(s)+C_{j_3 j_1 j_1 j_3}(s)
+C_{j_1 j_3 j_3 j_1}(s)+\biggr.
$$
\begin{equation}
\label{march40}
~~~~~~~~~~~\biggl.+C_{j_3 j_1 j_3 j_1}(s)+C_{j_1 j_3 j_1 j_3}(s)\biggr)=
\left(C_{j_1}(s)\right)^2\left(C_{j_3}(s)\right)^2.
\end{equation}

The equality (\ref{march40}) implies that
\begin{equation}
\label{march41}
8\sum\limits_{j_1,j_3=0}^{p}\biggl(C_{j_3 j_3 j_1 j_1}(s)+
C_{j_1 j_3 j_3 j_1}(s)+C_{j_3 j_1 j_3 j_1}(s)\biggr)=
\sum\limits_{j_1=0}^{p}\left(C_{j_1}(s)\right)^2\sum\limits_{j_3=0}^{p}\left(C_{j_3}(s)\right)^2.
\end{equation}

Passing to the limit $\lim\limits_{p\to\infty}$ in (\ref{march41})
and taking into account (\ref{march14}) (for $s$ instead of $T$), 
(\ref{march20a}), (\ref{march25}), we get
$$
8 \Biggl(\frac{(s-t)^2}{8} + \lim\limits_{p\to\infty}
\sum\limits_{j_1,j_3=0}^{p}C_{j_1 j_3 j_3 j_1}(s)+0\Biggr)=
(s-t)^2.
$$

The equality (\ref{march34}) is  proved.

\section{Generalization of Theorem~2.42 to the Case When the
Conditions $\phi_0(x)=1/\sqrt{T-t}$ and $\psi_l(\tau)\psi_{l-1}(\tau)\in L_2([t, T])$ 
$(l=2, 3,\ldots, k)$
are Omitted}

In this section, we will consider the following generalization of Theorem~2.42.

\vspace{1mm}

{\bf Theorem~2.49}\ \cite{arxiv-5}, \cite{arxiv-10}, \cite{arxiv-11}.\ {\it Assume that
the complete orthonormal system $\{\phi_j(x)\}_{j=0}^{\infty}$
in $L_2([t, T])$ and
$\psi_1(\tau),\ldots, \psi_k(\tau)\in L_2([t, T])$
are such that 

\vspace{-1mm}
$$
\lim\limits_{p_1,\ldots,p_k\to\infty}~
\sum\limits_{j_1=0}^{p_1}\ldots \sum\limits_{j_q=0}^{p_q}\ldots \sum\limits_{j_k=0}^{p_k}~
\biggl|_{q\ne g_1, g_2, \ldots, g_{2r-1},g_{2r}}\times
$$

\newpage
\noindent
$$
\times
\Biggl(~\sum\limits_{j_{g_1}=0}^{\min\{p_{g_1}, p_{g_2}\}} \sum\limits_{j_{g_3}=0}^{\min\{p_{g_3}, p_{g_4}\}}\ldots \Biggr.
\sum\limits_{j_{g_{2r-1}}=0}^{\min\{p_{g_{2r-1}}, p_{g_{2r}}\}}
C_{j_k\ldots j_1}\biggl|_{j_{g_1}=j_{g_2},\ldots, j_{g_{2r-1}}=j_{g_{2r}}}-
$$

\vspace{-2mm}
\begin{equation}
\label{april10}
\Biggl.-\frac{1}{2^r} \prod\limits_{l=1}^r {\bf 1}_{\{g_{2l}=g_{2l-1}+1\}}
C_{j_k \ldots j_1}\biggl|_{(j_{g_2} j_{g_1})\curvearrowright (\cdot)
\ldots (j_{g_{2r}} j_{g_{2r-1}})\curvearrowright (\cdot),
j_{g_{{}_{1}}}=~j_{g_{{}_{2}}},\ldots, j_{g_{{}_{2r-1}}}=~j_{g_{{}_{2r}}}
}\biggr.\Biggr)^2=0
\end{equation}

\vspace{2mm}
\noindent
for all $r=1, 2,\ldots,[k/2]$ 
and for all possible $g_1,g_2,\ldots,g_{2r-1},g_{2r}$ {\rm (}see {\rm (\ref{leto5007after}))}.
Then$,$ for the sum $\bar J^{*}[\psi^{(k)}]_{T,t}^{(i_1\ldots i_k)}$
of iterated It\^{o} stochastic integrals 
defined by {\rm (\ref{dsds9})}
the following 
expansion 
$$
\bar J^{*}[\psi^{(k)}]_{T,t}^{(i_1\ldots i_k)}=
\hbox{\vtop{\offinterlineskip\halign{
\hfil#\hfil\cr
{\rm l.i.m.}\cr
$\stackrel{}{{}_{p_1,\ldots,p_k\to \infty}}$\cr
}} }
\sum_{j_1=0}^{p_1}\ldots\sum_{j_k=0}^{p_k}
C_{j_k \ldots j_1}\prod\limits_{l=1}^k \zeta_{j_l}^{(i_l)}
$$

\noindent
that converges in the mean-square sense is valid, where 
\begin{equation}
\label{july15030}
~~~~~~~~C_{j_k \ldots j_1}=\int\limits_t^T\psi_k(t_k)\phi_{j_k}(t_k)\ldots
\int\limits_t^{t_2}
\psi_1(t_1)\phi_{j_1}(t_1)
dt_1\ldots dt_k
\end{equation}
is the Fourier coefficient, 
${\rm l.i.m.}$ is a limit in the mean-square sense,
$i_1, \ldots, i_k=0, 1,\ldots,m,$
$$
\zeta_{j}^{(i)}=
\int\limits_t^T \phi_{j}(\tau) d{\bf w}_{\tau}^{(i)}
$$ 

\noindent
are independent standard Gaussian random variables for various 
$i$ or $j$ {\rm (}when $i\ne 0${\rm )},
${\bf w}_{\tau}^{(i)}$ 
$(i=1,\ldots,m)$ are independent 
standard Wiener processes$,$
${\bf w}_{\tau}^{(0)}=\tau.$}

\vspace{2mm}

{\bf Proof.}\ To prove Theorem~2.49, we need to prove that under
the conditions of Theorem~2.49 the following equality 
$$
\hbox{\vtop{\offinterlineskip\halign{
\hfil#\hfil\cr
{\rm l.i.m.}\cr
$\stackrel{}{{}_{p\to \infty}}$\cr
}} }
\sum\limits_{\stackrel{j_1,\ldots,j_q,\ldots,j_k=0}{{}_{q\ne g_1, g_2,\ldots, g_{2r-1}, g_{2r}}}}^p
\frac{1}{2^r}
C_{j_k \ldots j_1}\biggl|_{(j_{g_2} j_{g_1})\curvearrowright (\cdot)
\ldots (j_{g_{2r}} j_{g_{2r-1}})\curvearrowright (\cdot),
j_{g_{{}_{1}}}=~j_{g_{{}_{2}}},\ldots, j_{g_{{}_{2r-1}}}=~j_{g_{{}_{2r}}}}\biggr.
\times 
$$

$$
\times
\prod\limits_{s=1}^r
{\bf 1}_{\{i_{g_{{}_{2s-1}}}=~i_{g_{{}_{2s}}}\ne 0\}}
J'[\phi_{j_{q_1}}\ldots \phi_{j_{q_{k-2r}}}]_{T,t}^{(i_{q_1}\ldots i_{q_{k-2r}})}=
$$

\vspace{3mm}
\begin{equation}
\label{febr14}
=\frac{1}{2^r}
J[\psi^{(k)}]_{T,t}^{s_r, \ldots, s_1}
\end{equation}

\vspace{5mm}
\noindent
holds w.~p.~1, where $g_{2}=g_{1}+1,\ldots, g_{2r}=g_{2r-1}+1,$
$g_{2i-1}\stackrel{\sf def}{=}s_i;$\ $i=1,2,\ldots,r;$\
$r=1,2,\ldots,\left[k/2\right],$ 
$(s_r,\ldots,s_1)\in {\rm A}_{k,r},$ $J[\psi^{(k)}]_{T,t}^{s_r,\ldots,s_1}$ is
defined by (\ref{30.1}) and ${\rm A}_{k,r}$ is defined by (\ref{30.5550001});
also we put $p_1=\ldots=p_k=p$ in (\ref{febr14}) to simplify the notation;
another notations in (\ref{febr14}) are the same as in Sect.~2.10.

Using the It\^{o} formula, we obtain w.~p.~1
$$
\int\limits_t^T \psi_k(t_k)\ldots \int\limits_t^{t_{l+2}}
\psi_{l+1}(t_{l+1}) \int\limits_t^{t_{l+1}}\psi_l(t_{l-1})\psi_{l-1}(t_{l-1})
\int\limits_t^{t_{l-1}}\psi_{l-2}(t_{l-2})\ldots
$$
$$
\ldots \int\limits_t^{t_2} \psi_1(t_1)d{\bf w}_{t_1}^{(i_1)}\ldots
d{\bf w}_{t_{l-2}}^{(i_{l-2})}dt_{l-1} d{\bf w}_{t_{l+1}}^{(i_{l+1})}\ldots
d{\bf w}_{t_k}^{(i_k)}=
$$
$$
=\int\limits_t^T \psi_k(t_k)\ldots \int\limits_t^{t_{l+2}}
\psi_{l+1}(t_{l+1}) \left(\int\limits_t^{t_{l+1}}\psi_l(t_{l-1})\psi_{l-1}(t_{l-1})dt_{l-1}\right)
\int\limits_t^{t_{l+1}}\psi_{l-2}(t_{l-2})\ldots
$$
$$
\ldots \int\limits_t^{t_2} \psi_1(t_1)d{\bf w}_{t_1}^{(i_1)}\ldots
d{\bf w}_{t_{l-2}}^{(i_{l-2})}d{\bf w}_{t_{l+1}}^{(i_{l+1})}\ldots
d{\bf w}_{t_k}^{(i_k)}-
$$
$$
-\int\limits_t^T \psi_k(t_k)\ldots \int\limits_t^{t_{l+2}}
\psi_{l+1}(t_{l+1}) 
\int\limits_t^{t_{l+1}}\psi_{l-2}(t_{l-2})
\left(\int\limits_t^{t_{l-2}}\psi_l(t_{l-1})\psi_{l-1}(t_{l-1})dt_{l-1}\right)\times
$$
\begin{equation}
\label{febr15}
~\times\int\limits_t^{t_{l-2}}\psi_{l-3}(t_{l-3})
\ldots
\int\limits_t^{t_2} \psi_1(t_1)d{\bf w}_{t_1}^{(i_1)}\ldots
d{\bf w}_{t_{l-3}}^{(i_{l-3})}d{\bf w}_{t_{l-2}}^{(i_{l-2})}d{\bf w}_{t_{l+1}}^{(i_{l+1})}\ldots
d{\bf w}_{t_k}^{(i_k)},
\end{equation}
where $l\ge 3.$ Note that the formula (\ref{febr15})
will change in an obvious way for the case $t_{l+1}=T.$
We will also assume that the transformation (\ref{febr15})
is not carried out for $l=2$ since the integral
$$
\int\limits_t^{t_{3}}\psi_2(t_{1})\psi_{1}(t_{1})dt_{1}
$$
is an internal integral on the left-hand side of 
(\ref{febr15}) for this case.

It is important to note that the transformation (\ref{febr15})
fully complies with the classical rules for replacing
the order of integration (Fubini's Theorem) if 
we replace
all differentials of the form
$d{\bf w}^{(i_j)}_{t_j}$ with $dt_j$
in (\ref{febr15}).

Indeed, formally changing the order of integration 
on the left-hand side of (\ref{febr15}) according 
to the classical rules, we have
\begin{equation}
\label{febr20}
~~~~~~~~\int\limits_t^T \psi_k(t_k)\ldots \int\limits_t^{t_{l+2}}
\psi_{l+1}(t_{l+1}) \int\limits_t^{t_{l+1}}\psi_l(t_{l-1})\psi_{l-1}(t_{l-1})
\int\limits_t^{t_{l-1}}\psi_{l-2}(t_{l-2})\ldots
\end{equation}
$$
\ldots \int\limits_t^{t_2} \psi_1(t_1)d{\bf w}_{t_1}^{(i_1)}\ldots
d{\bf w}_{t_{l-2}}^{(i_{l-2})}dt_{l-1} d{\bf w}_{t_{l+1}}^{(i_{l+1})}\ldots
d{\bf w}_{t_k}^{(i_k)}=
$$
$$
=\int\limits_t^T \psi_k(t_k)\ldots \int\limits_t^{t_{l+2}}
\psi_{l+1}(t_{l+1})\left(\int\limits_t^{t_{l+1}}\psi_1(t_{1})d{\bf w}_{t_1}^{(i_1)}
\ldots 
\int\limits_{t_{l-3}}^{t_{l+1}}\psi_{l-2}(t_{l-2})d{\bf w}_{t_{l-2}}^{(i_{l-2})}\times\right.
$$
$$
\left.\times \int\limits_{t_{l-2}}^{t_{l+1}}
\psi_{l}(t_{l-1})\psi_{l-1}(t_{l-1})dt_{l-1}\right)
d{\bf w}_{t_{l+1}}^{(i_{l+1})}\ldots
d{\bf w}_{t_k}^{(i_k)}=
$$
$$
=\int\limits_t^T \psi_k(t_k)\ldots \int\limits_t^{t_{l+2}}
\psi_{l+1}(t_{l+1})\left(\int\limits_t^{t_{l+1}}\psi_1(t_{1})d{\bf w}_{t_1}^{(i_1)}
\ldots 
\int\limits_{t_{l-3}}^{t_{l+1}}\psi_{l-2}(t_{l-2})d{\bf w}_{t_{l-2}}^{(i_{l-2})}\times\right.
$$
$$
\left.\times \left(\int\limits_t^{t_{l+1}}-\int\limits_t^{t_{l-2}}~\right)
\psi_{l}(t_{l-1})\psi_{l-1}(t_{l-1})dt_{l-1}\right)
d{\bf w}_{t_{l+1}}^{(i_{l+1})}\ldots
d{\bf w}_{t_k}^{(i_k)}=
$$
$$
=\int\limits_t^T \psi_k(t_k)\ldots \int\limits_t^{t_{l+2}}
\psi_{l+1}(t_{l+1})
\left(\int\limits_t^{t_{l+1}}
\psi_{l}(t_{l-1})\psi_{l-1}(t_{l-1})dt_{l-1}\right)
\int\limits_t^{t_{l+1}}\psi_1(t_{1})d{\bf w}_{t_1}^{(i_1)}
\ldots 
$$
$$
\ldots\int\limits_{t_{l-3}}^{t_{l+1}}\psi_{l-2}(t_{l-2})d{\bf w}_{t_{l-2}}^{(i_{l-2})}
d{\bf w}_{t_{l+1}}^{(i_{l+1})}\ldots
d{\bf w}_{t_k}^{(i_k)}-
$$
$$
-\int\limits_t^T \psi_k(t_k)\ldots \int\limits_t^{t_{l+2}}
\psi_{l+1}(t_{l+1})\int\limits_t^{t_{l+1}}\psi_1(t_{1})d{\bf w}_{t_1}^{(i_1)}
\ldots 
\int\limits_{t_{l-3}}^{t_{l+1}}\psi_{l-2}(t_{l-2})\times
$$
$$
\times \left(\int\limits_t^{t_{l-2}}
\psi_{l}(t_{l-1})\psi_{l-1}(t_{l-1})dt_{l-1}\right)
d{\bf w}_{t_{l-2}}^{(i_{l-2})}
d{\bf w}_{t_{l+1}}^{(i_{l+1})}\ldots
d{\bf w}_{t_k}^{(i_k)}=
$$
$$
=\int\limits_t^T \psi_k(t_k)\ldots \int\limits_t^{t_{l+2}}
\psi_{l+1}(t_{l+1}) \left(\int\limits_t^{t_{l+1}}\psi_l(t_{l-1})\psi_{l-1}(t_{l-1})dt_{l-1}\right)
\int\limits_t^{t_{l+1}}\psi_{l-2}(t_{l-2})\ldots
$$
$$
\ldots \int\limits_t^{t_2} \psi_1(t_1)d{\bf w}_{t_1}^{(i_1)}\ldots
d{\bf w}_{t_{l-2}}^{(i_{l-2})}d{\bf w}_{t_{l+1}}^{(i_{l+1})}\ldots
d{\bf w}_{t_k}^{(i_k)}-
$$
$$
-\int\limits_t^T \psi_k(t_k)\ldots \int\limits_t^{t_{l+2}}
\psi_{l+1}(t_{l+1}) 
\int\limits_t^{t_{l+1}}\psi_{l-2}(t_{l-2})
\left(\int\limits_t^{t_{l-2}}\psi_l(t_{l-1})\psi_{l-1}(t_{l-1})dt_{l-1}\right)\times
$$
\begin{equation}
\label{febr17}
~\times\int\limits_t^{t_{l-2}}\psi_{l-3}(t_{l-3})
\ldots
\int\limits_t^{t_2} \psi_1(t_1)d{\bf w}_{t_1}^{(i_1)}\ldots
d{\bf w}_{t_{l-3}}^{(i_{l-3})}d{\bf w}_{t_{l-2}}^{(i_{l-2})}d{\bf w}_{t_{l+1}}^{(i_{l+1})}\ldots
d{\bf w}_{t_k}^{(i_k)}.
\end{equation}

Comparing the right-hand sides of (\ref{febr15}) and (\ref{febr17})
we come to the conclusion that we got the same result.

The strict mathematical meaning of the transformations 
leading to (\ref{febr17}) is explained in Chapter~3,
at least for the case when $\psi_1(\tau),\ldots,\psi_k(\tau)$
are continuous functions on the interval $[t, T].$

Note that under the conditions of Theorem~2.49, the derivation of the 
formulas (\ref{febr15}) and (\ref{febr17}) will 
remain valid if in (\ref{febr15}) and (\ref{febr17}) we replace all 
differentials $d{\bf w}^{(i_j)}_{t_j}$ with $dt_j$
(this follows from Fubini's Theorem).

Recall that
$$
J[\psi^{(k)}]_{T,t}^{s_r, \ldots, s_1} \stackrel{\rm def}{=}\ 
\prod_{q=1}^r {\bf 1}_{\{i_{s_q}=
i_{s_{q}+1}\ne 0\}}\ \times
$$
$$
\times
\int\limits_t^T\psi_k(t_k)\ldots \int\limits_t^{t_{s_r+3}}
\psi_{s_r+2}(t_{s_r+2})
\int\limits_t^{t_{s_r+2}}\psi_{s_r}(t_{s_r+1})
\psi_{s_r+1}(t_{s_r+1}) \times
$$
$$
\times
\int\limits_t^{t_{s_r+1}}\psi_{s_r-1}(t_{s_r-1})
\ldots
\int\limits_t^{t_{s_1+3}}\psi_{s_1+2}(t_{s_1+2})
\int\limits_t^{t_{s_1+2}}\psi_{s_1}(t_{s_1+1})
\psi_{s_1+1}(t_{s_1+1}) \times
$$
$$
\times
\int\limits_t^{t_{s_1+1}}\psi_{s_1-1}(t_{s_1-1})
\ldots \int\limits_t^{t_2}\psi_1(t_1)
d{\bf w}_{t_1}^{(i_1)}\ldots d{\bf w}_{t_{s_1-1}}^{(i_{s_1-1})}
dt_{s_1+1}d{\bf w}_{t_{s_1+2}}^{(i_{s_1+2})}\ldots
$$

\vspace{-1mm}
\begin{equation}
\label{febr3000}
\ldots\
d{\bf w}_{t_{s_r-1}}^{(i_{s_r-1})}
dt_{s_r+1}d{\bf w}_{t_{s_r+2}}^{(i_{s_r+2})}\ldots d{\bf w}_{t_k}^{(i_k)},
\end{equation}

\vspace{3mm}
\noindent
where
${\rm A}_{k,r}$ is defined by (\ref{30.5550001}):

\vspace{-6mm}
$$
{\rm A}_{k,r}
=\bigl\{(s_r,\ldots,s_1):\
s_r>s_{r-1}+1,\ldots,s_2>s_1+1,\ s_r,\ldots,s_1=1,\ldots,k-1\bigr\}.
$$

\vspace{2mm}

Temporarily denote 

\vspace{-2mm}
$$
J[\psi^{(k)}]_{T,t}^{s_r, \ldots, s_1}\stackrel{\sf def}{=}
I[\psi^{(k)}]_{T,t}^{(i_1\ldots i_{s_1-1}i_{s_1+2} \ldots i_{s_r-1}i_{s_r+2}\ldots i_k)}.
$$

\vspace{3mm}

Let us carry out the transformation 
(\ref{febr15}) for the iterated It\^{o} stochastic integral 
$I[\psi^{(k)}]_{T,t}^{(i_1\ldots i_{s_1-1}i_{s_1+2} \ldots i_{s_r-1}i_{s_r+2}\ldots i_k)}$
iteratively for $s_1,\ldots,s_r.$ 
After this, apply (\ref{febr9}) to each of the obtained iterated It\^{o}
stochastic integrals.
As a result, we obtain w.~p.~1

\vspace{-2mm}
$$
I[\psi^{(k)}]_{T,t}^{(i_1\ldots i_{s_1-1}i_{s_1+2} \ldots i_{s_r-1}i_{s_r+2}\ldots i_k)}=
\prod_{q=1}^r {\bf 1}_{\{i_{s_q}=
i_{s_{q}+1}\ne 0\}}
\times
$$

\vspace{-2mm}
$$
\times\sum\limits_{d=1}^{2^r} \left(
\hat I[\psi^{(k)}]_{T,t}^{d(i_1\ldots i_{s_1-1}i_{s_1+2} \ldots i_{s_r-1}i_{s_r+2}\ldots i_k)}
-
\bar I[\psi^{(k)}]_{T,t}^{d(i_1\ldots i_{s_1-1}i_{s_1+2} \ldots i_{s_r-1}i_{s_r+2}\ldots i_k)}
\right)=
$$

\vspace{-2mm}
$$
=\prod_{q=1}^r {\bf 1}_{\{i_{s_q}=
i_{s_{q}+1}\ne 0\}}\times
$$

\vspace{-2mm}
$$
\times\hbox{\vtop{\offinterlineskip\halign{
\hfil#\hfil\cr
{\rm l.i.m.}\cr
$\stackrel{}{{}_{p\to \infty}}$\cr
}} }
\sum\limits_{j_1,\ldots, j_{s_1-1},j_{s_1+2}, \ldots, j_{s_r-1},j_{s_r+2},\ldots, j_k=0}^p
\ \sum\limits_{d=1}^{2^r}\ \Biggl(
\hat C_{j_1\ldots j_{s_1-1}j_{s_1+2}\ldots j_{s_r-1}j_{s_r+2}\ldots j_k}^{(d)}-\Biggr.
$$

\vspace{-1mm}
$$
\Biggl.-
\bar C_{j_1\ldots j_{s_1-1}j_{s_1+2}\ldots j_{s_r-1}j_{s_r+2}\ldots j_k}^{(d)}\Biggr)\times
$$

\vspace{-2mm}
\begin{equation}
\label{febr25}
~~~\times
J'[\phi_{j_1}\ldots \phi_{j_{s_1-1}}\phi_{j_{s_1+2}}\ldots \phi_{j_{s_r-1}}
\phi_{j_{s_r+2}}\ldots 
\phi_{j_k}]_{T,t}^{(i_1\ldots i_{s_1-1}i_{s_1+2} \ldots i_{s_r-1}i_{s_r+2}\ldots i_k)},
\end{equation}

\vspace{4mm}
\noindent
where some terms in the sum
$$
\sum\limits_{d=1}^{2^r}
$$

\noindent
can be identically equal to zero due to
the remark to (\ref{febr15}).

Taking into account that the integrals
$\hat I[\psi^{(k)}]_{T,t}^{d(i_1\ldots i_{s_1-1}i_{s_1+2} \ldots i_{s_r-1}i_{s_r+2}\ldots i_k)}$
and the Fourier coefficients
$\hat C_{j_1\ldots j_{s_1-1}j_{s_1+2}\ldots j_{s_r-1}j_{s_r+2}\ldots j_k}^{(d)}$
are 
formed on the basis of the same kernels (the same applies to 
the integrals
$\bar I[\psi^{(k)}]_{T,t}^{d(i_1\ldots i_{s_1-1}i_{s_1+2} \ldots i_{s_r-1}i_{s_r+2}\ldots i_k)}$
and 
the Fourier coefficients
$\bar C_{j_1\ldots j_{s_1-1}j_{s_1+2}\ldots j_{s_r-1}j_{s_r+2}\ldots j_k}^{(d)}$), as well 
as a remark about the relationship of the transformation (\ref{febr15}) 
based on the It\^{o} formula                  
and on the basis 
of classical rules for replacing
the order of integration (see the derivation of (\ref{febr17})), we obtain using Fubini's theorem 
(applying the inverse transformation from (\ref{febr17}) to (\ref{febr20})
in which all differentials of the form
$d{\bf w}^{(i_j)}_{t_j}$ are replaced with $dt_j$)
$$
\sum\limits_{d=1}^{2^r}\ \Biggl(
\hat C_{j_1\ldots j_{s_1-1}j_{s_1+2}\ldots j_{s_r-1}j_{s_r+2}\ldots j_k}^{(d)}-
\bar C_{j_1\ldots j_{s_1-1}j_{s_1+2}\ldots j_{s_r-1}j_{s_r+2}\ldots j_k}^{(d)}\Biggr)=
$$  
\begin{equation}
\label{febr26}
~~~~~~~~~~=C_{j_k \ldots j_1}\biggl|_{(j_{g_2} j_{g_1})\curvearrowright (\cdot)
\ldots (j_{g_{2r}} j_{g_{2r-1}})\curvearrowright (\cdot),
j_{g_{{}_{1}}}=~j_{g_{{}_{2}}},\ldots, j_{g_{{}_{2r-1}}}=~j_{g_{{}_{2r}}}}\biggr.,
\end{equation}

\vspace{5mm}
\noindent
where $g_{2}=g_{1}+1,\ldots, g_{2r}=g_{2r-1}+1.$
Combining (\ref{febr25}) and (\ref{febr26}), we get
for $g_{2}=g_{1}+1,\ldots, g_{2r}=g_{2r-1}+1$ w.~p.~1 

\vspace{1mm}
$$
I[\psi^{(k)}]_{T,t}^{(i_1\ldots i_{s_1-1}i_{s_1+2} \ldots i_{s_r-1}i_{s_r+2}\ldots i_k)}=
$$

\vspace{-4mm}
$$
=\hbox{\vtop{\offinterlineskip\halign{
\hfil#\hfil\cr
{\rm l.i.m.}\cr
$\stackrel{}{{}_{p\to \infty}}$\cr
}} }
\sum\limits_{\stackrel{j_1,\ldots,j_q,\ldots,j_k=0}{{}_{q\ne g_1, g_2,\ldots, g_{2r-1}, g_{2r}}}}^p
C_{j_k \ldots j_1}\biggl|_{(j_{g_2} j_{g_1})\curvearrowright (\cdot)
\ldots (j_{g_{2r}} j_{g_{2r-1}})\curvearrowright (\cdot),
j_{g_{{}_{1}}}=~j_{g_{{}_{2}}},\ldots, j_{g_{{}_{2r-1}}}=~j_{g_{{}_{2r}}}}\biggr.
\times 
$$

\vspace{2mm}
$$
\times
\prod\limits_{s=1}^r
{\bf 1}_{\{i_{g_{{}_{2s-1}}}=~i_{g_{{}_{2s}}}\ne 0\}}
J'[\phi_{j_{q_1}}\ldots \phi_{j_{q_{k-2r}}}]_{T,t}^{(i_{q_1}\ldots i_{q_{k-2r}})},
$$

\vspace{4mm}
\noindent
where 
we use the notations from Sect.~2.10. The equality (\ref{febr14}) is proved
for the case when $\{\phi_j(x)\}_{j=0}^{\infty}$ is an arbitrary
complete orthonormal system of functions
in the space $L_2([t, T]).$ 
Thus, the condition $\phi_0(x)=1/\sqrt{T-t}$ in Theorem~2.42 can be omitted.

Let us separately explain why the condition 
$\psi_l(\tau)\psi_{l-1}(\tau)\in L_2([t, T])$ 
$(l=2, 3,\ldots, k)$
in Theorem~2.42 can also be omitted.
Recall that this condition appeared due to the direct
application of (\ref{febr5000}) to the iterated
It\^{o} stochastic integral
$J[\psi^{(k)}]_{T,t}^{s_r, \ldots, s_1}$ defined by (\ref{febr3000})
(see the transition from (\ref{after905})
to (\ref{after906})).

It is easy to see that the kernels $\hat K_d(t_1,\ldots, t_{k-2r})$ and
$\bar K_d(t_1,\ldots, t_{k-2r})$ of the iterated It\^{o}
stochastic integrals
$\hat I[\psi^{(k)}]_{T,t}^{d(i_1\ldots i_{s_1-1}i_{s_1+2} \ldots i_{s_r-1}i_{s_r+2}\ldots i_k)}$
and 
$\bar I[\psi^{(k)}]_{T,t}^{d(i_1\ldots i_{s_1-1}i_{s_1+2} \ldots i_{s_r-1}i_{s_r+2}\ldots i_k)}$
have the same structure as (\ref{chain200}) but with new wight 
functions $\hat \psi_1(\tau),\ldots, \hat\psi_{k-2r}(\tau)$
and $\bar \psi_1(\tau),\ldots, \bar\psi_{k-2r}(\tau)$, some of which 
possibly coincide with $\psi_1(\tau),\ldots, \psi_k(\tau)$
(see (\ref{febr15})).
Moreover, the conditions $\psi_1(\tau),\ldots,\psi_{k}(\tau)\in L_2([t, T])$
and $\psi_l(\tau)\psi_{l-1}(\tau)\in L_1([t, T])$ 
$(l=2, 3,\ldots, k)$ guarantee that
$\hat K_d(t_1,\ldots, t_{k-2r}),
\bar K_d(t_1,\ldots, t_{k-2r})\in L_2([t, T])$ (see 
(\ref{febr15})).
This means that the formula (\ref{febr25}) is true if
$\psi_1(\tau),\ldots,\psi_{k}(\tau)\in L_2([t, T])$
and $\psi_l(\tau)\psi_{l-1}(\tau)\in L_1([t, T])$ 
$(l=2, 3,\ldots, k).$
Furthermore, the formula (\ref{febr26}) holds under
the conditions 
$\psi_1(\tau),\ldots,\psi_{k}(\tau)\in L_2([t, T])$
and $\psi_l(\tau)\psi_{l-1}(\tau)\in L_1([t, T])$ 
$(l=2, 3,\ldots, k).$

Since the condition $\psi_1(\tau),\ldots,\psi_{k}(\tau)\in L_2([t, T])$
implies the condition $\psi_l(\tau)\psi_{l-1}(\tau)\in L_1([t, T])$ 
$(l=2, 3,\ldots, k),$ then the condition $\psi_l(\tau)\psi_{l-1}(\tau)\in L_1([t, T])$ 
$(l=2, 3,\ldots, k)$ can be omitted in the above reasoning.

Thus, the equalities (\ref{febr25}) and (\ref{febr26})
are satisfied under the condition 
$\psi_1(\tau),\ldots,\psi_{k}(\tau)\in L_2([t, T])$
and the condition $\psi_l(\tau)\psi_{l-1}(\tau)\in L_2([t, T])$ 
$(l=2, 3,\ldots, k)$ can be omitted in 
Theorem~2.42.
Theorem~2.49 is proved.

\section{Expansion of Iterated Stratonovich Stochastic Integrals
of Multiplicity 5. The Case of an Ar\-bit\-ra\-ry Complete Orthonormal System of 
Functions in the Space $L_2([t,T])$ and $\psi_1(\tau),\ldots, \psi_5(\tau)
\equiv 1$}

{\bf Theorem~2.50}\ \cite{arxiv-5}, \cite{arxiv-10}, \cite{arxiv-11}.\ {\it Suppose that
$\{\phi_j(x)\}_{j=0}^{\infty}$ is an arbitrary complete orthonormal system of 
functions in the space $L_2([t,T]).$
Then$,$ for the iterated Stra\-to\-no\-vich stochastic integral
of fifth multiplicity 
$$
J^{*}[\psi^{(5)}]_{T,t}=
{\int\limits_t^{*}}^T
\ldots
{\int\limits_t^{*}}^{t_2}
d{\bf w}_{t_1}^{(i_1)}
\ldots d{\bf w}_{t_5}^{(i_5)}
$$
the following 
expansion 
$$
J^{*}[\psi^{(5)}]_{T,t}=
\hbox{\vtop{\offinterlineskip\halign{
\hfil#\hfil\cr
{\rm l.i.m.}\cr
$\stackrel{}{{}_{p\to \infty}}$\cr
}} }
\sum\limits_{j_1,\ldots, j_5=0}^{p}
C_{j_5 \ldots j_1}\zeta_{j_1}^{(i_1)}\ldots \zeta_{j_5}^{(i_5)}
$$
that converges in the mean-square sense is valid, where 
$i_1,\ldots,i_5=0, 1,\ldots,m,$
$$
C_{j_5\ldots j_1}=\int\limits_t^T
\phi_{j_5}(t_5)
\ldots
\int\limits_t^{t_2}
\phi_{j_1}(t_1)dt_1\ldots dt_5
$$
and
$$
\zeta_{j}^{(i)}=
\int\limits_t^T \phi_{j}(\tau) d{\bf w}_{\tau}^{(i)}
$$ 
are independent standard Gaussian random variables for various 
$i$ or $j$ {\rm (}when $i\ne 0${\rm ),}
${\bf w}_{\tau}^{(i)}$ 
$(i=1,\ldots,m)$ are independent 
standard Wiener processes$,$
${\bf w}_{\tau}^{(0)}=\tau.$}

{\bf Proof.}\ {\bf Step~1.}\ According to Theorem~2.49,
we conclude that Theorem~2.50 will be proved if we prove
the following equalities (see (\ref{april10}) for $k=5, r=1$ and $k=5, r=2$
($p_1=\ldots=p_5=p$))
under the conditions of Theorem~2.50 
\begin{equation}
\label{april11}
~~~~~~~~~~~\lim\limits_{p\to\infty}
\sum\limits_{j_3, j_4,j_5=0}^{p}
\left(\sum\limits_{j_1=0}^{p}
C_{j_5 j_4 j_3 j_1 j_1}-\frac{1}{2} 
C_{j_5 j_4 j_3 j_1 j_1}\biggl|_{(j_{1} j_{1})\curvearrowright (\cdot)}
\biggr.\right)^2=0,
\end{equation}
\begin{equation}
\label{april12}
\lim\limits_{p\to\infty}
\sum\limits_{j_2, j_4,j_5=0}^{p}
\left(\sum\limits_{j_1=0}^{p}
C_{j_5 j_4 j_1 j_2 j_1}\right)^2=0,
\end{equation}
\begin{equation}
\label{april13}
\lim\limits_{p\to\infty}
\sum\limits_{j_2, j_3,j_5=0}^{p}
\left(\sum\limits_{j_1=0}^{p}
C_{j_5 j_1 j_3 j_2 j_1}\right)^2=0,
\end{equation}
\begin{equation}
\label{april14}
\lim\limits_{p\to\infty}
\sum\limits_{j_2, j_3,j_4=0}^{p}
\left(\sum\limits_{j_1=0}^{p}
C_{j_1 j_4 j_3 j_2 j_1}\right)^2=0,
\end{equation}
\begin{equation}
\label{april15}
~~~~~~~~~~~\lim\limits_{p\to\infty}
\sum\limits_{j_1, j_4,j_5=0}^{p}
\left(\sum\limits_{j_2=0}^{p}
C_{j_5 j_4 j_2 j_2 j_1}-\frac{1}{2} 
C_{j_5 j_4 j_2 j_2 j_1}\biggl|_{(j_{2} j_{2})\curvearrowright (\cdot)}
\biggr.\right)^2=0,
\end{equation}
\begin{equation}
\label{april16}
\lim\limits_{p\to\infty}
\sum\limits_{j_1, j_3,j_5=0}^{p}
\left(\sum\limits_{j_2=0}^{p}
C_{j_5 j_2 j_3 j_2 j_1}\right)^2=0,
\end{equation}
\begin{equation}
\label{april17}
\lim\limits_{p\to\infty}
\sum\limits_{j_1, j_3,j_4=0}^{p}
\left(\sum\limits_{j_2=0}^{p}
C_{j_2 j_4 j_3 j_2 j_1}\right)^2=0,
\end{equation}
\begin{equation}
\label{april18}
~~~~~~~~~~~\lim\limits_{p\to\infty}
\sum\limits_{j_1, j_2,j_5=0}^{p}
\left(\sum\limits_{j_3=0}^{p}
C_{j_5 j_3 j_3 j_2 j_1}-\frac{1}{2} 
C_{j_5 j_3 j_3 j_2 j_1}\biggl|_{(j_{3} j_{3})\curvearrowright (\cdot)}
\biggr.\right)^2=0,
\end{equation}
\begin{equation}
\label{april19}
\lim\limits_{p\to\infty}
\sum\limits_{j_1, j_2,j_4=0}^{p}
\left(\sum\limits_{j_3=0}^{p}
C_{j_3 j_4 j_3 j_2 j_1}\right)^2=0,
\end{equation}
\begin{equation}
\label{april20}
~~~~~~~~~~~\lim\limits_{p\to\infty}
\sum\limits_{j_1, j_2,j_3=0}^{p}
\left(\sum\limits_{j_4=0}^{p}
C_{j_4 j_4 j_3 j_2 j_1}-\frac{1}{2} 
C_{j_4 j_4 j_3 j_2 j_1}\biggl|_{(j_{4} j_{4})\curvearrowright (\cdot)}
\biggr.\right)^2=0,
\end{equation}
\begin{equation}
\label{april21}
~~~~~~~~~~~\lim\limits_{p\to\infty}
\sum\limits_{j_5=0}^{p}
\left(\sum\limits_{j_1,j_3=0}^{p}
C_{j_5 j_3 j_3 j_1 j_1}-\frac{1}{4} 
C_{j_5 j_3 j_3 j_1 j_1}\biggl|_{(j_{1} j_{1})\curvearrowright (\cdot),
(j_{3} j_{3})\curvearrowright (\cdot)}
\biggr.\right)^2=0,
\end{equation}
\begin{equation}
\label{april22}
\lim\limits_{p\to\infty}
\sum\limits_{j_4=0}^{p}
\left(\sum\limits_{j_1,j_3=0}^{p}
C_{j_3 j_4 j_3 j_1 j_1}\right)^2=0,
\end{equation}
\begin{equation}
\label{april23}
~~~~~~~~~~~\lim\limits_{p\to\infty}
\sum\limits_{j_3=0}^{p}
\left(\sum\limits_{j_1,j_4=0}^{p}
C_{j_4 j_4 j_3 j_1 j_1}-\frac{1}{4} 
C_{j_4 j_4 j_3 j_1 j_1}\biggl|_{(j_{1} j_{1})\curvearrowright (\cdot),
(j_{4} j_{4})\curvearrowright (\cdot)}
\biggr.\right)^2=0,
\end{equation}
\begin{equation}
\label{april24}
\lim\limits_{p\to\infty}
\sum\limits_{j_5=0}^{p}
\left(\sum\limits_{j_1,j_2=0}^{p}
C_{j_5 j_2 j_1 j_2 j_1}\right)^2=0,
\end{equation}
\begin{equation}
\label{april25}
\lim\limits_{p\to\infty}
\sum\limits_{j_4=0}^{p}
\left(\sum\limits_{j_1,j_2=0}^{p}
C_{j_2 j_4 j_1 j_2 j_1}\right)^2=0,
\end{equation}
\begin{equation}
\label{april26}
\lim\limits_{p\to\infty}
\sum\limits_{j_2=0}^{p}
\left(\sum\limits_{j_1,j_4=0}^{p}
C_{j_4 j_4 j_1 j_2 j_1}\right)^2=0,
\end{equation}
\begin{equation}
\label{april27}
\lim\limits_{p\to\infty}
\sum\limits_{j_5=0}^{p}
\left(\sum\limits_{j_1,j_2=0}^{p}
C_{j_5 j_1 j_2 j_2 j_1}\right)^2=0,
\end{equation}
\begin{equation}
\label{april28}
\lim\limits_{p\to\infty}
\sum\limits_{j_3=0}^{p}
\left(\sum\limits_{j_1,j_2=0}^{p}
C_{j_2 j_1 j_3 j_2 j_1}\right)^2=0,
\end{equation}
\begin{equation}
\label{april29}
\lim\limits_{p\to\infty}
\sum\limits_{j_2=0}^{p}
\left(\sum\limits_{j_1,j_3=0}^{p}
C_{j_3 j_1 j_3 j_2 j_1}\right)^2=0,
\end{equation}
\begin{equation}
\label{april30}
\lim\limits_{p\to\infty}
\sum\limits_{j_4=0}^{p}
\left(\sum\limits_{j_1,j_2=0}^{p}
C_{j_1 j_4 j_2 j_2 j_1}\right)^2=0,
\end{equation}
\begin{equation}
\label{april31}
\lim\limits_{p\to\infty}
\sum\limits_{j_3=0}^{p}
\left(\sum\limits_{j_1,j_2=0}^{p}
C_{j_1 j_2 j_3 j_2 j_1}\right)^2=0,
\end{equation}
\begin{equation}
\label{april32}
\lim\limits_{p\to\infty}
\sum\limits_{j_2=0}^{p}
\left(\sum\limits_{j_1,j_3=0}^{p}
C_{j_1 j_3 j_3 j_2 j_1}\right)^2=0,
\end{equation}
\begin{equation}
\label{april33}
~~~~~~~~~~~\lim\limits_{p\to\infty}
\sum\limits_{j_1=0}^{p}
\left(\sum\limits_{j_2,j_4=0}^{p}
C_{j_4 j_4 j_2 j_2 j_1}-\frac{1}{4} 
C_{j_4 j_4 j_2 j_2 j_1}\biggl|_{(j_{2} j_{2})\curvearrowright (\cdot),
(j_{4} j_{4})\curvearrowright (\cdot)}
\biggr.\right)^2=0,
\end{equation}
\begin{equation}
\label{april34}
\lim\limits_{p\to\infty}
\sum\limits_{j_1=0}^{p}
\left(\sum\limits_{j_2,j_3=0}^{p}
C_{j_3 j_2 j_3 j_2 j_1}\right)^2=0,
\end{equation}
\begin{equation}
\label{april35}
\lim\limits_{p\to\infty}
\sum\limits_{j_1=0}^{p}
\left(\sum\limits_{j_2,j_3=0}^{p}
C_{j_2 j_3 j_3 j_2 j_1}\right)^2=0.
\end{equation}

\vspace{2mm}

{\bf Step~2.}\ Let us prove the equalities (\ref{april11})--(\ref{april20}).
Using Fubini's Theorem and Parseval's equality, we obtain
the following relations 
for the prelimit
expressions on the left-hand sides of (\ref{april11})--(\ref{april20})
$$
\sum\limits_{j_3, j_4,j_5=0}^{p}
\left(\sum\limits_{j_1=0}^{p}
C_{j_5 j_4 j_3 j_1 j_1}-\frac{1}{2} 
C_{j_5 j_4 j_3 j_1 j_1}\biggl|_{(j_{1} j_{1})\curvearrowright (\cdot)}
\biggr.\right)^2=
$$
$$
=\sum\limits_{j_3, j_4,j_5=0}^{p}
\Biggl(\int\limits_t^T\phi_{j_5}(t_5)\int\limits_t^{t_5}\phi_{j_4}(t_4)
\int\limits_t^{t_4}\phi_{j_3}(t_3)\left(\sum\limits_{j_1=0}^{p}
\frac{1}{2}\left(\int\limits_t^{t_3}\phi_{j_1}(\tau)d\tau\right)^2-\frac{t_3-t}{2}\right)\times
\Biggr.
$$
$$
\times
\Biggl.dt_3 dt_4 dt_5\Biggl)^2\le
$$
$$
\le\sum\limits_{j_3, j_4,j_5=0}^{\infty}
\Biggl(\int\limits_t^T\phi_{j_5}(t_5)\int\limits_t^{t_5}\phi_{j_4}(t_4)
\int\limits_t^{t_4}\phi_{j_3}(t_3)\left(\sum\limits_{j_1=0}^{p}
\frac{1}{2}\left(\int\limits_t^{t_3}\phi_{j_1}(\tau)d\tau\right)^2-\frac{t_3-t}{2}\right)\times
\Biggr.
$$
$$
\times
\Biggl.dt_3 dt_4 dt_5\Biggl)^2=
$$
\begin{equation}
\label{april36}
~~~=\int\limits_{[t,T]^3}\left({\bf 1}_{\{t_3<t_4<t_5\}}\right)^2
\left(\sum\limits_{j_1=0}^{p}
\frac{1}{2}\left(\int\limits_t^{t_3}\phi_{j_1}(\tau)d\tau\right)^2-\frac{t_3-t}{2}\right)^2
dt_3 dt_4 dt_5,
\end{equation}
$$
\sum\limits_{j_2, j_4,j_5=0}^{p}
\left(\sum\limits_{j_1=0}^{p}
C_{j_5 j_4 j_1 j_2 j_1}\right)^2=
$$
$$
=\sum\limits_{j_2, j_4,j_5=0}^{p}
\Biggl(\int\limits_t^T\phi_{j_5}(t_5)\int\limits_t^{t_5}\phi_{j_4}(t_4)
\int\limits_t^{t_4}\phi_{j_2}(t_2)\sum\limits_{j_1=0}^{p}
\int\limits_t^{t_2}\phi_{j_1}(t_1)dt_1
\int\limits_{t_2}^{t_4}\phi_{j_1}(t_3)dt_3\times\Biggr.
$$
$$
\Biggl.\times dt_2 dt_4 dt_5\Biggr)^2\le
$$
$$
\le\sum\limits_{j_2, j_4,j_5=0}^{\infty}
\Biggl(\int\limits_t^T\phi_{j_5}(t_5)\int\limits_t^{t_5}\phi_{j_4}(t_4)
\int\limits_t^{t_4}\phi_{j_2}(t_2)\sum\limits_{j_1=0}^{p}
\int\limits_t^{t_2}\phi_{j_1}(t_1)dt_1
\int\limits_{t_2}^{t_4}\phi_{j_1}(t_3)dt_3\times\Biggr.
$$
$$
\Biggl.\times dt_2 dt_4 dt_5\Biggr)^2=
$$
\begin{equation}
\label{april37}
~~~=\int\limits_{[t,T]^3}\left({\bf 1}_{\{t_2<t_4<t_5\}}\right)^2
\left(\sum\limits_{j_1=0}^{p}
\int\limits_t^{t_2}\phi_{j_1}(t_1)dt_1
\int\limits_{t_2}^{t_4}\phi_{j_1}(t_3)dt_3\right)^2
dt_2 dt_4 dt_5,
\end{equation}

\vspace{5mm}

$$
\sum\limits_{j_2, j_3,j_5=0}^{p}
\left(\sum\limits_{j_1=0}^{p}
C_{j_5 j_1 j_3 j_2 j_1}\right)^2=
$$
$$
=\sum\limits_{j_2, j_3,j_5=0}^{p}
\Biggl(\int\limits_t^T\phi_{j_5}(t_5)\int\limits_t^{t_5}\phi_{j_3}(t_3)
\int\limits_t^{t_3}\phi_{j_2}(t_2)\sum\limits_{j_1=0}^{p}
\int\limits_t^{t_2}\phi_{j_1}(t_1)dt_1
\int\limits_{t_3}^{t_5}\phi_{j_1}(t_4)dt_4\times\Biggr.
$$
$$
\Biggl.\times dt_2 dt_3 dt_5\Biggr)^2\le
$$
$$
\le\sum\limits_{j_2, j_3,j_5=0}^{\infty}
\Biggl(\int\limits_t^T\phi_{j_5}(t_5)\int\limits_t^{t_5}\phi_{j_3}(t_3)
\int\limits_t^{t_3}\phi_{j_2}(t_2)\sum\limits_{j_1=0}^{p}
\int\limits_t^{t_2}\phi_{j_1}(t_1)dt_1
\int\limits_{t_3}^{t_5}\phi_{j_1}(t_4)dt_4\times\Biggr.
$$
$$
\Biggl.\times dt_2 dt_3 dt_5\Biggr)^2=
$$
\begin{equation}
\label{april38}
~~~=\int\limits_{[t,T]^3}\left({\bf 1}_{\{t_2<t_3<t_5\}}\right)^2
\left(\sum\limits_{j_1=0}^{p}
\int\limits_t^{t_2}\phi_{j_1}(t_1)dt_1
\int\limits_{t_3}^{t_5}\phi_{j_1}(t_4)dt_4\right)^2
dt_2 dt_3 dt_5,
\end{equation}

\vspace{5mm}

$$
\sum\limits_{j_2, j_3,j_4=0}^{p}
\left(\sum\limits_{j_1=0}^{p}
C_{j_1 j_4 j_3 j_2 j_1}\right)^2=
$$
$$
=\sum\limits_{j_2, j_3,j_4=0}^{p}
\Biggl(\int\limits_t^T\phi_{j_4}(t_4)\int\limits_t^{t_4}\phi_{j_3}(t_3)
\int\limits_t^{t_3}\phi_{j_2}(t_2)\sum\limits_{j_1=0}^{p}
\int\limits_t^{t_2}\phi_{j_1}(t_1)dt_1
\int\limits_{t_4}^{T}\phi_{j_1}(t_5)dt_5\times\Biggr.
$$
$$
\Biggl.\times dt_2 dt_3 dt_4\Biggr)^2\le
$$
$$
\le\sum\limits_{j_2, j_3,j_4=0}^{\infty}
\Biggl(\int\limits_t^T\phi_{j_4}(t_4)\int\limits_t^{t_4}\phi_{j_3}(t_3)
\int\limits_t^{t_3}\phi_{j_2}(t_2)\sum\limits_{j_1=0}^{p}
\int\limits_t^{t_2}\phi_{j_1}(t_1)dt_1
\int\limits_{t_4}^{T}\phi_{j_1}(t_5)dt_5\times\Biggr.
$$
$$
\Biggl.\times dt_2 dt_3 dt_4\Biggr)^2=
$$
\begin{equation}
\label{april39}
~~~=\int\limits_{[t,T]^3}\left({\bf 1}_{\{t_2<t_3<t_4\}}\right)^2
\left(\sum\limits_{j_1=0}^{p}
\int\limits_t^{t_2}\phi_{j_1}(t_1)dt_1
\int\limits_{t_4}^{T}\phi_{j_1}(t_5)dt_5\right)^2
dt_2 dt_3 dt_4,
\end{equation}

\vspace{5mm}

$$
\sum\limits_{j_1, j_4,j_5=0}^{p}
\left(\sum\limits_{j_2=0}^{p}
C_{j_5 j_4 j_2 j_2 j_1}-\frac{1}{2} 
C_{j_5 j_4 j_2 j_2 j_1}\biggl|_{(j_{2} j_{2})\curvearrowright (\cdot)}
\biggr.\right)^2=
$$
$$
=\sum\limits_{j_1, j_4,j_5=0}^{p}
\Biggl(\int\limits_t^T\phi_{j_5}(t_5)\int\limits_t^{t_5}\phi_{j_4}(t_4)
\int\limits_t^{t_4}\phi_{j_1}(t_1)
\sum\limits_{j_2=0}^p \int\limits_{t_1}^{t_4}\phi_{j_2}(t_2)
\int\limits_{t_2}^{t_4}\phi_{j_2}(t_3)dt_3 dt_2\times
\Biggr.
$$
$$
\times
\Biggl. dt_1 dt_4 dt_5-
\frac{1}{2}\int\limits_t^T\phi_{j_5}(t_5)\int\limits_t^{t_5}\phi_{j_4}(t_4)
\int\limits_t^{t_4}\int\limits_t^{t_2}\phi_{j_1}(t_1)dt_1 dt_2 dt_4 dt_5
\Biggl)^2=
$$
$$
=\sum\limits_{j_1, j_4,j_5=0}^{p}
\Biggl(\int\limits_t^T\phi_{j_5}(t_5)\int\limits_t^{t_5}\phi_{j_4}(t_4)
\int\limits_t^{t_4}\phi_{j_1}(t_1)
\left(\sum\limits_{j_2=0}^p \frac{1}{2}\left(\int\limits_{t_1}^{t_4}\phi_{j_2}(t_2)
dt_2\right)^2-\frac{t_4-t_1}{2}\right)\times
\Biggr.
$$
$$
\times
\Biggl. dt_1 dt_4 dt_5
\Biggl)^2\le
$$
$$
\le\sum\limits_{j_1, j_4,j_5=0}^{\infty}
\Biggl(\int\limits_t^T\phi_{j_5}(t_5)\int\limits_t^{t_5}\phi_{j_4}(t_4)
\int\limits_t^{t_4}\phi_{j_1}(t_1)
\left(\sum\limits_{j_2=0}^p \frac{1}{2}\left(\int\limits_{t_1}^{t_4}\phi_{j_2}(t_2)
dt_2\right)^2-\frac{t_4-t_1}{2}\right)\times
\Biggr.
$$
$$
\times
\Biggl. dt_1 dt_4 dt_5
\Biggl)^2=
$$
\begin{equation}
\label{april40}
=\int\limits_{[t,T]^3}\left({\bf 1}_{\{t_1<t_4<t_5\}}\right)^2
\left(\sum\limits_{j_2=0}^p \frac{1}{2}\left(\int\limits_{t_1}^{t_4}\phi_{j_2}(t_2)
dt_2\right)^2-\frac{t_4-t_1}{2}\right)^2
dt_1 dt_4 dt_5,
\end{equation}

\vspace{5mm}

$$
\sum\limits_{j_1, j_3,j_5=0}^{p}
\left(\sum\limits_{j_2=0}^{p}
C_{j_5 j_2 j_3 j_2 j_1}\right)^2=
$$
$$
=\sum\limits_{j_1, j_3,j_5=0}^{p}
\Biggl(\sum\limits_{j_2=0}^p \int\limits_t^T\phi_{j_5}(t_5)\int\limits_t^{t_5}\phi_{j_3}(t_3)
\int\limits_t^{t_3}\phi_{j_2}(t_2)
\int\limits_t^{t_2}\phi_{j_1}(t_1)dt_1 dt_2
\int\limits_{t_3}^{t_5}\phi_{j_2}(t_4)dt_4\times
\Biggr.
$$
$$
\Biggl.\times dt_3 dt_5\Biggr)^2=
$$
$$
=\sum\limits_{j_1, j_3,j_5=0}^{p}
\Biggl(\int\limits_t^T\phi_{j_5}(t_5)\int\limits_t^{t_5}\phi_{j_3}(t_3)
\int\limits_t^{t_3}\phi_{j_1}(t_1)
\sum\limits_{j_2=0}^p \int\limits_{t_1}^{t_3}\phi_{j_2}(t_2)dt_2
\int\limits_{t_3}^{t_5}\phi_{j_2}(t_4)dt_4\times
\Biggr.
$$
$$
\Biggl.\times dt_1 dt_3 dt_5\Biggr)^2\le
$$
$$
\le\sum\limits_{j_1, j_3,j_5=0}^{\infty}
\Biggl(\int\limits_t^T\phi_{j_5}(t_5)\int\limits_t^{t_5}\phi_{j_3}(t_3)
\int\limits_t^{t_3}\phi_{j_1}(t_1)
\sum\limits_{j_2=0}^p \int\limits_{t_1}^{t_3}\phi_{j_2}(t_2)dt_2
\int\limits_{t_3}^{t_5}\phi_{j_2}(t_4)dt_4\times
\Biggr.
$$
$$
\Biggl.\times dt_1 dt_3 dt_5\Biggr)^2=
$$
\begin{equation}
\label{april41}
~~~=\int\limits_{[t,T]^3}\left({\bf 1}_{\{t_1<t_3<t_5\}}\right)^2
\left(\sum\limits_{j_2=0}^p \int\limits_{t_1}^{t_3}\phi_{j_2}(t_2)dt_2
\int\limits_{t_3}^{t_5}\phi_{j_2}(t_4)dt_4\right)^2
dt_1 dt_3 dt_5,
\end{equation}

\vspace{5mm}

$$
\sum\limits_{j_1, j_3,j_4=0}^{p}
\left(\sum\limits_{j_2=0}^{p}
C_{j_2 j_4 j_3 j_2 j_1}\right)^2=
$$
$$
=\sum\limits_{j_1, j_3,j_4=0}^{p}
\Biggl(\sum\limits_{j_2=0}^p \int\limits_t^T\phi_{j_4}(t_4)\int\limits_t^{t_4}\phi_{j_3}(t_3)
\int\limits_t^{t_3}\phi_{j_2}(t_2)
\int\limits_t^{t_2}\phi_{j_1}(t_1)dt_1 dt_2 dt_3
\int\limits_{t_4}^{T}\phi_{j_2}(t_5)dt_5\times
\Biggr.
$$
$$
\Biggl.\times dt_4\Biggr)^2=
$$
$$
=\sum\limits_{j_1, j_3,j_4=0}^{p}
\Biggl(\int\limits_t^T\phi_{j_4}(t_4)\int\limits_t^{t_4}\phi_{j_3}(t_3)
\int\limits_t^{t_3}\phi_{j_1}(t_1)
\sum\limits_{j_2=0}^p \int\limits_{t_1}^{t_3}\phi_{j_2}(t_2)dt_2
\int\limits_{t_4}^{T}\phi_{j_2}(t_5)dt_5\times
\Biggr.
$$
$$
\Biggl.\times dt_1 dt_3 dt_4\Biggr)^2\le
$$
$$
\le\sum\limits_{j_1, j_3,j_4=0}^{\infty}
\Biggl(\int\limits_t^T\phi_{j_4}(t_4)\int\limits_t^{t_4}\phi_{j_3}(t_3)
\int\limits_t^{t_3}\phi_{j_1}(t_1)
\sum\limits_{j_2=0}^p \int\limits_{t_1}^{t_3}\phi_{j_2}(t_2)dt_2
\int\limits_{t_4}^{T}\phi_{j_2}(t_5)dt_5\times
\Biggr.
$$
$$
\Biggl.\times dt_1 dt_3 dt_4\Biggr)^2=
$$
\begin{equation}
\label{april42}
~~~=\int\limits_{[t,T]^3}\left({\bf 1}_{\{t_1<t_3<t_4\}}\right)^2
\left(\sum\limits_{j_2=0}^p \int\limits_{t_1}^{t_3}\phi_{j_2}(t_2)dt_2
\int\limits_{t_4}^{T}\phi_{j_2}(t_5)dt_5\right)^2
dt_1 dt_3 dt_4,
\end{equation}

\vspace{5mm}

$$
\sum\limits_{j_1, j_2,j_5=0}^{p}
\left(\sum\limits_{j_3=0}^{p}
C_{j_5 j_3 j_3 j_2 j_1}-\frac{1}{2} 
C_{j_5 j_3 j_3 j_2 j_1}\biggl|_{(j_{3} j_{3})\curvearrowright (\cdot)}
\biggr.\right)^2=
$$
$$
=\sum\limits_{j_1, j_2,j_5=0}^{p}
\Biggl(\sum\limits_{j_3=0}^p
\int\limits_t^T\phi_{j_5}(t_5)\int\limits_t^{t_5}\phi_{j_1}(t_1)
\int\limits_{t_1}^{t_5}\phi_{j_2}(t_2)
\int\limits_{t_2}^{t_5}\phi_{j_3}(t_3)
\int\limits_{t_3}^{t_5}\phi_{j_3}(t_4)dt_4 dt_3\times
\Biggr.
$$
$$
\times
\Biggl. dt_2 dt_1 dt_5-
\frac{1}{2}\int\limits_t^T\phi_{j_5}(t_5)\int\limits_t^{t_5}\int\limits_t^{t_3}
\phi_{j_2}(t_2)
\int\limits_t^{t_2}\phi_{j_1}(t_1)dt_1 dt_2 dt_3 dt_5
\Biggl)^2=
$$
$$
=\sum\limits_{j_1, j_2,j_5=0}^{p}
\Biggl(\sum\limits_{j_3=0}^p
\int\limits_t^T\phi_{j_5}(t_5)\int\limits_t^{t_5}\phi_{j_1}(t_1)
\int\limits_{t_1}^{t_5}\phi_{j_2}(t_2)
\int\limits_{t_2}^{t_5}\phi_{j_3}(t_3)
\int\limits_{t_3}^{t_5}\phi_{j_3}(t_4)dt_4 dt_3\times
\Biggr.
$$
$$
\times
\Biggl. dt_2 dt_1 dt_5-
\frac{1}{2}\int\limits_t^T\phi_{j_5}(t_5)
\int\limits_{t}^{t_5}\phi_{j_1}(t_1)
\int\limits_{t_1}^{t_5}\phi_{j_2}(t_2)\int\limits_{t_2}^{t_5}dt_3 dt_2 dt_1 dt_5
\Biggl)^2=
$$
$$
=\sum\limits_{j_1, j_2,j_5=0}^{p}
\Biggl(\int\limits_t^T\phi_{j_5}(t_5)\int\limits_t^{t_5}\phi_{j_1}(t_1)
\int\limits_{t_1}^{t_5}\phi_{j_2}(t_2)
\left(\sum\limits_{j_3=0}^p \frac{1}{2}\left(\int\limits_{t_2}^{t_5}\phi_{j_3}(t_3)
dt_3\right)^2-\frac{t_5-t_2}{2}\right)\times
\Biggr.
$$
$$
\times
\Biggl. dt_2 dt_1 dt_5
\Biggl)^2\le
$$
$$
\le\sum\limits_{j_1, j_2,j_5=0}^{\infty}
\Biggl(\int\limits_t^T\phi_{j_5}(t_5)\int\limits_t^{t_5}\phi_{j_1}(t_1)
\int\limits_{t_1}^{t_5}\phi_{j_2}(t_2)
\left(\sum\limits_{j_3=0}^p \frac{1}{2}\left(\int\limits_{t_2}^{t_5}\phi_{j_3}(t_3)
dt_3\right)^2-\frac{t_5-t_2}{2}\right)\times
\Biggr.
$$
$$
\times
\Biggl. dt_2 dt_1 dt_5
\Biggl)^2=
$$
\begin{equation}
\label{april43}
=\int\limits_{[t,T]^3}\left({\bf 1}_{\{t_1<t_2<t_5\}}\right)^2
\left(\sum\limits_{j_3=0}^p \frac{1}{2}\left(\int\limits_{t_2}^{t_5}\phi_{j_3}(t_3)
dt_3\right)^2-\frac{t_5-t_2}{2}\right)^2
dt_2 dt_1 dt_5,
\end{equation}

\vspace{5mm}

$$
\sum\limits_{j_1, j_2,j_4=0}^{p}
\left(\sum\limits_{j_3=0}^{p}
C_{j_3 j_4 j_3 j_2 j_1}\right)^2=
$$
$$
=\sum\limits_{j_1, j_2,j_4=0}^{p}
\Biggl(\sum\limits_{j_3=0}^p \int\limits_t^T\phi_{j_1}(t_1)\int\limits_{t_1}^{T}\phi_{j_2}(t_2)
\int\limits_{t_2}^{T}\phi_{j_3}(t_3)
\int\limits_{t_3}^{T}\phi_{j_4}(t_4)
\int\limits_{t_4}^{T}\phi_{j_3}(t_5)dt_5 dt_4 dt_3\times
\Biggr.
$$
$$
\Biggl.\times dt_2 dt_1\Biggr)^2=
$$
$$
=\sum\limits_{j_1, j_2,j_4=0}^{p}
\Biggl(\int\limits_t^T\phi_{j_1}(t_1)\int\limits_{t_1}^{T}\phi_{j_2}(t_2)
\int\limits_{t_2}^{T}\phi_{j_4}(t_4)
\sum\limits_{j_3=0}^p\int\limits_{t_4}^{T}\phi_{j_3}(t_5)dt_5
\int\limits_{t_2}^{t_4}\phi_{j_3}(t_3)dt_3 dt_4 \times
\Biggr.
$$
$$
\Biggl.\times dt_2 dt_1\Biggr)^2\le
$$
$$
\le\sum\limits_{j_1, j_2,j_4=0}^{\infty}
\Biggl(\int\limits_t^T\phi_{j_1}(t_1)\int\limits_{t_1}^{T}\phi_{j_2}(t_2)
\int\limits_{t_2}^{T}\phi_{j_4}(t_4)
\sum\limits_{j_3=0}^p\int\limits_{t_4}^{T}\phi_{j_3}(t_5)dt_5
\int\limits_{t_2}^{t_4}\phi_{j_3}(t_3)dt_3 dt_4 \times
\Biggr.
$$
$$
\Biggl.\times dt_2 dt_1\Biggr)^2=
$$
\begin{equation}
\label{april44}
~~~=\int\limits_{[t,T]^3}\left({\bf 1}_{\{t_1<t_2<t_4\}}\right)^2
\left(
\sum\limits_{j_3=0}^p\int\limits_{t_4}^{T}\phi_{j_3}(t_5)dt_5
\int\limits_{t_2}^{t_4}\phi_{j_3}(t_3)dt_3
\right)^2
dt_4 dt_2 dt_1,
\end{equation}

\vspace{5mm}

$$
\sum\limits_{j_1, j_2,j_3=0}^{p}
\left(\sum\limits_{j_4=0}^{p}
C_{j_4 j_4 j_3 j_2 j_1}-\frac{1}{2} 
C_{j_4 j_4 j_3 j_2 j_1}\biggl|_{(j_{4} j_{4})\curvearrowright (\cdot)}
\biggr.\right)^2=
$$
$$
=\sum\limits_{j_1, j_2,j_3=0}^{p}
\Biggl(
\int\limits_t^T\phi_{j_3}(t_3)\int\limits_t^{t_3}\phi_{j_2}(t_2)
\int\limits_{t}^{t_2}\phi_{j_1}(t_1) dt_1 dt_2
\sum\limits_{j_4=0}^p\int\limits_{t_3}^{T}\phi_{j_4}(t_4)
\int\limits_{t_4}^{T}\phi_{j_4}(t_5)dt_5 dt_4\times
\Biggr.
$$
$$
\times
\Biggl. dt_3-
\frac{1}{2}\int\limits_t^T\int\limits_t^{t_4}\phi_{j_3}(t_3)\int\limits_t^{t_3}
\phi_{j_2}(t_2)
\int\limits_t^{t_2}\phi_{j_1}(t_1)dt_1 dt_2 dt_3 dt_4
\Biggl)^2=
$$
$$
=\sum\limits_{j_1, j_2,j_3=0}^{p}
\Biggl(
\int\limits_t^T\phi_{j_3}(t_3)\int\limits_t^{t_3}\phi_{j_2}(t_2)
\int\limits_{t}^{t_2}\phi_{j_1}(t_1) 
\left(\sum\limits_{j_4=0}^p\frac{1}{2}\left(\int\limits_{t_3}^{T}\phi_{j_4}(t_4)dt_4\right)^2
-\frac{T-t_3}{2}\right)\times
\Biggr.
$$
$$
\times
\Biggl. dt_1 dt_2dt_3
\Biggl)^2\le
$$
$$
\le\sum\limits_{j_1, j_2,j_3=0}^{\infty}
\Biggl(
\int\limits_t^T\phi_{j_3}(t_3)\int\limits_t^{t_3}\phi_{j_2}(t_2)
\int\limits_{t}^{t_2}\phi_{j_1}(t_1) 
\left(\sum\limits_{j_4=0}^p\frac{1}{2}\left(\int\limits_{t_3}^{T}\phi_{j_4}(t_4)dt_4\right)^2
-\frac{T-t_3}{2}\right)\times
\Biggr.
$$
$$
\times
\Biggl. dt_1 dt_2dt_3
\Biggl)^2=
$$
\begin{equation}
\label{april45}
=\int\limits_{[t,T]^3}\left({\bf 1}_{\{t_1<t_2<t_3\}}\right)^2
\left(
\sum\limits_{j_4=0}^p\frac{1}{2}\left(\int\limits_{t_3}^{T}\phi_{j_4}(t_4)dt_4\right)^2
-\frac{T-t_3}{2}
\right)^2
dt_1 dt_2 dt_3.
\end{equation}

\vspace{2mm}

Further, applying the Parseval equality and the generalized Parseval equality
as well as using the 
Cauchy--Bunyakovsky inequality,
we have (see the proof of Theorem~2.45)
\begin{equation}
\label{april46}
~~~~~~~~~~~~\sum\limits_{j=0}^{\infty}
\left(\int\limits_{t_1}^{t_2}\phi_{j}(s)ds\right)^2
=\int\limits_t^T \left({\bf 1}_{\{t_1<s<t_2\}}\right)^2 ds=
t_2-t_1,
\end{equation}
$$
\sum\limits_{j=0}^{\infty}\int\limits_{t_1}^{t_2}\phi_{j}(s)ds
\int\limits_{t_3}^{t_4}\phi_{j}(s)ds=
\sum\limits_{j=0}^{\infty}\int\limits_t^T {\bf 1}_{\{t_1<s<t_2\}}\phi_{j}(s)ds
\int\limits_t^T {\bf 1}_{\{t_3<s<t_4\}}\phi_{j}(s)ds=
$$
\begin{equation}
\label{april47}
=
\int\limits_t^T {\bf 1}_{\{t_1<s<t_2\}}{\bf 1}_{\{t_3<s<t_4\}}ds=0,
\end{equation}
\begin{equation}
\label{april48}
~~~~~~~~~~~~\left\vert
(t_2-t_1)
-\sum\limits_{j=0}^{p}
\left(\int\limits_{t_1}^{t_2}\phi_{j}(s)ds\right)^2
\right\vert\le
t_2-t_1\le T-t<\infty,
\end{equation}
$$
\left(\sum\limits_{j=0}^{p}\int\limits_{t_1}^{t_2}\phi_{j}(s)ds
\int\limits_{t_3}^{t_4}\phi_{j}(s)ds\right)^2\le
\sum\limits_{j=0}^{p}\left(\int\limits_{t_1}^{t_2}\phi_{j}(s)ds\right)^2
\sum\limits_{j=0}^{p}\left(\int\limits_{t_3}^{t_4}\phi_{j}(s)ds\right)^2\le
$$

\vspace{-2mm}
\begin{equation}
\label{april49}
\le
(t_2-t_1)(t_4-t_3)\le (T-t)^2<\infty,
\end{equation}

\vspace{3mm}
\noindent
where $t\le t_1<t_2\le t_3<t_4\le T.$

Using Lebesgue's Dominated Convergence Theorem and (\ref{april46})--(\ref{april49}), 
we obtain that the right-hand sides of (\ref{april36})--(\ref{april45}) 
tend to zero when $p\to\infty.$
The equalities (\ref{april11})--(\ref{april20}) are proved.

{\bf Step~3.}\ Before proving the equalities (\ref{april21})--(\ref{april35}), we show that
\begin{equation}
\label{april50}
\left\vert\sum\limits_{j_1, j_3=0}^p C_{j_3 j_3 j_1 j_1}(s,\tau)\right\vert\le K,
\end{equation}
\begin{equation}
\label{april51}
\left\vert\sum\limits_{j_1, j_3=0}^p C_{j_1 j_3 j_3 j_1}(s,\tau)\right\vert\le K,
\end{equation}
\begin{equation}
\label{april52}
\left\vert\sum\limits_{j_1, j_2=0}^p C_{j_2 j_1 j_2 j_1}(s,\tau)\right\vert\le K,
\end{equation}
\begin{equation}
\label{april53}
~~~\sum\limits_{j_2=0}^p
\left(\sum\limits_{j_1=0}^p C_{j_1 j_2 j_1}(s,\tau)\right)^2\le
\int\limits_{\tau}^s \left(\sum\limits_{j_1=0}^p
\int\limits_{\tau}^{t_2}\phi_{j_1}(t_1)dt_1\int\limits_{t_2}^{s}\phi_{j_1}(t_3)dt_3\right)^2 dt_2,
\end{equation}

\noindent
where constant $K$ does not depend on $p, s, \tau;$
here and further in this proof
$$
C_{j_k \ldots j_1}(s,\tau)=\int\limits_{\tau}^s
\phi_{j_k}(t_k)\ldots
\int\limits_{\tau}^{t_2}
\phi_{j_1}(t_1)dt_1\ldots dt_k \ \ \ (k=1,\ldots,4,\ t\le\tau<s\le T).
$$

Further, by $K, K_1, K_2$ we will denote contants
that can change from line to line.

By analogy with (\ref{march13}), (\ref{march21}), (\ref{march26}) 
and (\ref{march20a}), (\ref{march25}), (\ref{march34})
we get
\begin{equation}
\label{april54}
\sum\limits_{j_1,j_3=0}^p C_{j_3 j_3 j_1 j_1}(s,\tau)=
\sum\limits_{j_1,j_3=0}^p C_{j_3}(s,\tau)C_{j_3 j_1 j_1}(s,\tau)-\frac{1}{8}
\left(\sum\limits_{j_1=0}^p \bigl(C_{j_1}(s,\tau)\bigr)^2\right)^2,
\end{equation}
\begin{equation}
\label{april55}
\sum\limits_{j_1,j_2=0}^p C_{j_2 j_1 j_2 j_1}(s,\tau)=
\sum\limits_{j_1,j_2=0}^p C_{j_2}(s,\tau)C_{j_1 j_2 j_1}(s,\tau)
-\frac{1}{2}\sum\limits_{j_1,j_2=0}^p C_{j_1 j_2}(s,\tau)C_{j_2 j_1}(s,\tau),
\end{equation}
\begin{equation}
\label{april56}
\sum\limits_{j_1,j_3=0}^p C_{j_1 j_3 j_3 j_1}(s,\tau)=
\sum\limits_{j_1,j_3=0}^p C_{j_1}(s,\tau)C_{j_3 j_3 j_1}(s,\tau)-\frac{1}{2}
\sum\limits_{j_1,j_3=0}^p \bigl(C_{j_3 j_1}(s,\tau)\bigr)^2,
\end{equation}
\begin{equation}
\label{april56x}
\lim\limits_{p\to\infty}
\sum\limits_{j_1, j_3=0}^{p}
C_{j_3 j_3 j_1 j_1}(s,\tau)=\frac{1}{8}(s-\tau)^2,
\biggr.
\end{equation}
\begin{equation}
\label{april56xx}
\lim\limits_{p\to\infty}
\sum\limits_{j_1, j_2=0}^{p}
C_{j_2 j_1 j_2 j_1}(s,\tau)=0
\biggr.,
\end{equation}
\begin{equation}
\label{april56xxx}
\lim\limits_{p\to\infty}
\sum\limits_{j_1, j_3=0}^{p}
C_{j_1 j_3 j_3 j_1}(s,\tau)=0
\biggr..
\end{equation}

Using (\ref{april54}), Parseval's equality,
Cauchy--Bunyakovsky's inequality, as well as Fubini's Theorem and the elementary inequality
$(a+b)^2\le 2 a^2 + 2 b^2,$
we obtain
$$
\left(\sum\limits_{j_1,j_3=0}^p C_{j_3 j_3 j_1 j_1}(s,\tau)\right)^2\le
2\left(\sum\limits_{j_1,j_3=0}^p C_{j_3}(s,\tau)C_{j_3 j_1 j_1}(s,\tau)\right)^2+
$$
$$
+2\cdot \frac{1}{64}
\left(\sum\limits_{j_1=0}^p \bigl(C_{j_1}(s,\tau)\bigr)^2\right)^4\le
$$
$$
\le 2\sum\limits_{j_3=0}^p \left(C_{j_3}(s,\tau)\right)^2
\sum\limits_{j_3=0}^p \left(\sum\limits_{j_1=0}^p
C_{j_3 j_1 j_1}(s,\tau)\right)^2+K_1\le
$$
$$
\le K_2
\sum\limits_{j_3=0}^{\infty} \left(\sum\limits_{j_1=0}^p
C_{j_3 j_1 j_1}(s,\tau)\right)^2+K_1=
$$
$$
= K_2
\sum\limits_{j_3=0}^{\infty} \left(\int\limits_{\tau}^s
\phi_{j_3}(t_3)\sum\limits_{j_1=0}^p
\int\limits_{\tau}^{t_3}\phi_{j_1}(t_2)
\int\limits_{\tau}^{t_2}\phi_{j_1}(t_1)dt_1 dt_2 dt_3
\right)^2+K_1=
$$
$$
= K_2
\int\limits_{\tau}^s \left(\frac{1}{2}\sum\limits_{j_1=0}^p
\left(\int\limits_{\tau}^{t_3}\phi_{j_1}(t_2)dt_2\right)^2\right)^2 dt_3
+K_1\le 
$$
$$
\le K_2
\int\limits_{\tau}^s \left(\frac{1}{2}\sum\limits_{j_1=0}^{\infty}
\left(\int\limits_{\tau}^{t_3}\phi_{j_1}(t_2)dt_2\right)^2\right)^2 dt_3
+K_1=
$$
$$
=
K_2
\int\limits_{\tau}^s \left(\frac{1}{2}(t_3-\tau)\right)^2 dt_3
+K_1\le K<\infty,
$$
where constants $K, K_1, K_2$ do not depend on 
$p, s, \tau.$ The equality (\ref{april50}) is proved.

Let us prove (\ref{april51}). Using (\ref{april56}) and the
above reasoning, we get
$$
\left(\sum\limits_{j_1,j_3=0}^p C_{j_1 j_3 j_3 j_1}(s,\tau)\right)^2\le
2\left(\sum\limits_{j_1,j_3=0}^p C_{j_1}(s,\tau)C_{j_3 j_3 j_1}(s,\tau)\right)^2+
$$
$$
+2\cdot \frac{1}{4}\left(\sum\limits_{j_1,j_3=0}^p \bigl(C_{j_3 j_1}(s,\tau)\bigr)^2\right)^2\le
$$
$$
\le 2\sum\limits_{j_1=0}^p \left(C_{j_1}(s,\tau)\right)^2
\sum\limits_{j_1=0}^p \left(\sum\limits_{j_3=0}^p
C_{j_3 j_3 j_1}(s,\tau)\right)^2+K_1\le
$$
$$
\le K_2
\sum\limits_{j_1=0}^{\infty} \left(\sum\limits_{j_3=0}^p
C_{j_3 j_3 j_1}(s,\tau)\right)^2+K_1=
$$
$$
= K_2
\sum\limits_{j_1=0}^{\infty} \left(\int\limits_{\tau}^s
\phi_{j_1}(t_1)\sum\limits_{j_3=0}^p
\int\limits_{t_1}^{s}\phi_{j_3}(t_2)
\int\limits_{t_2}^{s}\phi_{j_3}(t_3)dt_3 dt_2 dt_1
\right)^2+K_1=
$$
$$
= K_2
\int\limits_{\tau}^s \left(\frac{1}{2}\sum\limits_{j_3=0}^p
\left(\int\limits_{t_1}^{s}\phi_{j_3}(t_2)dt_2\right)^2\right)^2 dt_1
+K_1\le 
$$
$$
\le K_2
\int\limits_{\tau}^s \left(\frac{1}{2}\sum\limits_{j_3=0}^{\infty}
\left(\int\limits_{t_1}^{s}\phi_{j_3}(t_2)dt_2\right)^2\right)^2 dt_1
+K_1=
$$
$$
=
K_2
\int\limits_{\tau}^s \left(\frac{1}{2}(s-t_1)\right)^2 dt_1
+K_1\le K<\infty,
$$
where constants $K, K_1, K_2$ do not depend on 
$p, s, \tau.$ The equality (\ref{april51}) is proved.

Let us prove (\ref{april52}), (\ref{april53}). Applying (\ref{april55}), (\ref{april49}) and the
above reasoning, we have
$$
\left(\sum\limits_{j_1,j_2=0}^p C_{j_2 j_1 j_2 j_1}(s,\tau)\right)^2\le
2\left(\sum\limits_{j_1,j_2=0}^p C_{j_2}(s,\tau)C_{j_1 j_2 j_1}(s,\tau)\right)^2+
$$
$$
+2\cdot \frac{1}{4}\left(\sum\limits_{j_1,j_2=0}^p C_{j_1 j_2}(s,\tau)C_{j_2 j_1}(s,\tau)\right)^2\le
$$
$$
\le 2\sum\limits_{j_2=0}^p \left(C_{j_2}(s,\tau)\right)^2
\sum\limits_{j_2=0}^p \left(\sum\limits_{j_1=0}^p
C_{j_1 j_2 j_1}(s,\tau)\right)^2+
$$
$$
+
\frac{1}{2}\sum\limits_{j_1,j_2=0}^p \left(C_{j_1 j_2}(s,\tau)\right)^2
\sum\limits_{j_1,j_2=0}^p \left(C_{j_2 j_1}(s,\tau)\right)^2
\le
$$
\begin{equation}
\label{april60}
\le K_2
\sum\limits_{j_2=0}^{p} \left(\sum\limits_{j_1=0}^p
C_{j_1 j_2 j_1}(s,\tau)\right)^2+K_1\le
K_2
\sum\limits_{j_2=0}^{\infty} \left(\sum\limits_{j_1=0}^p
C_{j_1 j_2 j_1}(s,\tau)\right)^2+K_1=
\end{equation}
$$
= K_2
\sum\limits_{j_2=0}^{\infty} \left(\int\limits_{\tau}^s
\phi_{j_2}(t_2)\sum\limits_{j_1=0}^p
\int\limits_{\tau}^{t_2}\phi_{j_1}(t_1)dt_1
\int\limits_{t_2}^{s}\phi_{j_1}(t_3)dt_3 dt_2 
\right)^2+K_1=
$$
\begin{equation}
\label{april61}
~~~~~~~~~~~= K_2
\int\limits_{\tau}^s \left(\sum\limits_{j_1=0}^p
\int\limits_{\tau}^{t_2}\phi_{j_1}(t_1)dt_1\int\limits_{t_2}^{s}\phi_{j_1}(t_3)dt_3\right)^2 dt_2
+K_1\le 
\end{equation}
$$
\le K_2
\int\limits_{\tau}^s \left((t_2-\tau)(s-t_2)\right)^2 dt_2
+K_1\le K<\infty,
$$
where constants $K, K_1, K_2$ do not depend on 
$p, s, \tau.$ The equalities (\ref{april52}) and (\ref{april53}) (see (\ref{april60}), (\ref{april61}))
are proved.

{\bf Step~4.}\ Let us start proving the equalities (\ref{april21})--(\ref{april35}).
Using Fubini's Theorem and Parseval's equality, we obtain
the following relations 
for the prelimit
expressions on the left-hand sides of (\ref{april21}), 
(\ref{april24}), (\ref{april27}),
(\ref{april33})--(\ref{april35})
$$
\sum\limits_{j_5=0}^{p}
\left(\sum\limits_{j_1,j_3=0}^{p}
C_{j_5 j_3 j_3 j_1 j_1}-\frac{1}{4} 
C_{j_5 j_3 j_3 j_1 j_1}\biggl|_{(j_{1} j_{1})\curvearrowright (\cdot),
(j_{3} j_{3})\curvearrowright (\cdot)}
\biggr.\right)^2=
$$
$$
=\sum\limits_{j_5=0}^{p}
\left(\int\limits_t^T \phi_{j_5}(t_5) \left(\sum\limits_{j_1,j_3=0}^{p}
C_{j_3 j_3 j_1 j_1} (t_5,t)-\frac{1}{4}\int\limits_t^{t_5}(\tau-t)d\tau\right)dt_5\right)^2\le
$$
$$
\le\sum\limits_{j_5=0}^{\infty}
\left(\int\limits_t^T \phi_{j_5}(t_5) \left(\sum\limits_{j_1,j_3=0}^{p}
C_{j_3 j_3 j_1 j_1} (t_5,t)-\frac{1}{4}\int\limits_t^{t_5}(\tau-t)d\tau\right)dt_5\right)^2=
$$
\begin{equation}
\label{april62}
=\int\limits_t^T
\left(\sum\limits_{j_1,j_3=0}^{p}
C_{j_3 j_3 j_1 j_1} (t_5,t)-\frac{1}{8}(t_5-t)^2\right)^2 dt_5,
\end{equation}
$$
\sum\limits_{j_5=0}^{p}
\left(\sum\limits_{j_1,j_2=0}^{p}
C_{j_5 j_2 j_1 j_2 j_1}\right)^2
=\sum\limits_{j_5=0}^{p}
\left(\int\limits_t^T \phi_{j_5}(t_5) \sum\limits_{j_1,j_2=0}^{p}
C_{j_2 j_1 j_2 j_1} (t_5,t)dt_5\right)^2\le
$$
\begin{equation}
\label{april63}
\le\sum\limits_{j_5=0}^{\infty}
\left(\int\limits_t^T \phi_{j_5}(t_5) \sum\limits_{j_1,j_2=0}^{p}
C_{j_2 j_1 j_2 j_1} (t_5,t)dt_5\right)^2
=\int\limits_t^T
\left(\sum\limits_{j_1,j_2=0}^{p}
C_{j_2 j_1 j_2 j_1} (t_5,t)\right)^2 dt_5,
\end{equation}

\vspace{-4mm}
$$
\sum\limits_{j_5=0}^{p}
\left(\sum\limits_{j_1,j_2=0}^{p}
C_{j_5 j_1 j_2 j_2 j_1}\right)^2=
\sum\limits_{j_5=0}^{p}
\left(\int\limits_t^T \phi_{j_5}(t_5) \sum\limits_{j_1,j_2=0}^{p}
C_{j_1 j_2 j_2 j_1} (t_5,t)dt_5\right)^2\le
$$
\begin{equation}
\label{april64}
\le\sum\limits_{j_5=0}^{\infty}
\left(\int\limits_t^T \phi_{j_5}(t_5) \sum\limits_{j_1,j_2=0}^{p}
C_{j_1 j_2 j_2 j_1} (t_5,t)dt_5\right)^2
=\int\limits_t^T
\left(\sum\limits_{j_1,j_2=0}^{p}
C_{j_1 j_2 j_2 j_1} (t_5,t)\right)^2 dt_5,
\end{equation}

$$
\sum\limits_{j_1=0}^{p}
\left(\sum\limits_{j_2,j_4=0}^{p}
C_{j_4 j_4 j_2 j_2 j_1}-\frac{1}{4} 
C_{j_4 j_4 j_2 j_2 j_1}\biggl|_{(j_{2} j_{2})\curvearrowright (\cdot),
(j_{4} j_{4})\curvearrowright (\cdot)}
\biggr.\right)^2=
$$
$$
=\sum\limits_{j_1=0}^{p}
\left(\int\limits_t^T  \phi_{j_1}(t_1)
\sum\limits_{j_2,j_4=0}^{p}
\int\limits_{t_1}^T  \phi_{j_2}(t_2)
\int\limits_{t_2}^T  \phi_{j_2}(t_3)
\int\limits_{t_3}^T  \phi_{j_4}(t_4)
\int\limits_{t_4}^T  \phi_{j_4}(t_5)
dt_5 dt_4 dt_3 dt_2 dt_1-\right.
$$
$$
\left.-\frac{1}{4} \int\limits_{t}^T\int\limits_{t}^{t_5} 
\int\limits_{t}^{t_3}\phi_{j_1}(t_1) dt_1 dt_3 dt_5\right)^2=
$$
$$
=\sum\limits_{j_1=0}^{p}
\left(\int\limits_t^T  \phi_{j_1}(t_1)
\left(\sum\limits_{j_2,j_4=0}^{p}
C_{j_4 j_4 j_2 j_2}(T,t_1)-
\frac{1}{4} \int\limits_{t_1}^T(T-t_3)dt_3
\right)dt_1\right)^2\le
$$
$$
\le\sum\limits_{j_1=0}^{\infty}
\left(\int\limits_t^T  \phi_{j_1}(t_1)
\left(\sum\limits_{j_2,j_4=0}^{p}
C_{j_4 j_4 j_2 j_2}(T,t_1)-
\frac{1}{8}(T-t_1)^2
\right)dt_1\right)^2=
$$
\begin{equation}
\label{april65}
~~~~~~~~=\int\limits_t^T
\left(\sum\limits_{j_2,j_4=0}^{p}
C_{j_4 j_4 j_2 j_2}(T,t_1)-
\frac{1}{8}(T-t_1)^2
\right)^2 dt_1,
\end{equation}
$$
\sum\limits_{j_1=0}^{p}
\left(\sum\limits_{j_2,j_3=0}^{p}
C_{j_3 j_2 j_3 j_2 j_1}\right)^2=
$$
$$
=\sum\limits_{j_1=0}^{p}
\Biggl(\int\limits_t^T  \phi_{j_1}(t_1)
\sum\limits_{j_2,j_3=0}^{p}
\int\limits_{t_1}^T  \phi_{j_2}(t_2)
\int\limits_{t_2}^T  \phi_{j_3}(t_3)
\int\limits_{t_3}^T  \phi_{j_2}(t_4)
\int\limits_{t_4}^T  \phi_{j_3}(t_5)
dt_5 dt_4\times\Biggr.
$$
$$
\Biggl.\times
dt_3 dt_2 dt_1\Biggr)^2=
$$
$$
=\sum\limits_{j_1=0}^{p}
\left(\int\limits_t^T  \phi_{j_1}(t_1)
\sum\limits_{j_2,j_3=0}^{p}
C_{j_3 j_2 j_3 j_2}(T,t_1)dt_1\right)^2\le
$$
\begin{equation}
\label{april66}
\le\sum\limits_{j_1=0}^{\infty}
\left(\int\limits_t^T  \phi_{j_1}(t_1)
\sum\limits_{j_2,j_3=0}^{p}
C_{j_3 j_2 j_3 j_2}(T,t_1)dt_1\right)^2=
\int\limits_t^T
\left(\sum\limits_{j_2,j_3=0}^{p}
C_{j_3 j_2 j_3 j_2}(T,t_1)
\right)^2 dt_1,
\end{equation}

\vspace{2mm}

$$
\sum\limits_{j_1=0}^{p}
\left(\sum\limits_{j_2,j_3=0}^{p}
C_{j_2 j_3 j_3 j_2 j_1}\right)^2=
$$
$$
=\sum\limits_{j_1=0}^{p}
\Biggl(\int\limits_t^T  \phi_{j_1}(t_1)
\sum\limits_{j_2,j_3=0}^{p}
\int\limits_{t_1}^T  \phi_{j_2}(t_2)
\int\limits_{t_2}^T  \phi_{j_3}(t_3)
\int\limits_{t_3}^T  \phi_{j_3}(t_4)
\int\limits_{t_4}^T  \phi_{j_2}(t_5)
dt_5 dt_4\times\Biggr.
$$
$$
\Biggl.\times
dt_3 dt_2 dt_1\Biggr)^2=
$$
$$
=\sum\limits_{j_1=0}^{p}
\left(\int\limits_t^T  \phi_{j_1}(t_1)
\sum\limits_{j_2,j_3=0}^{p}
C_{j_2 j_3 j_3 j_2}(T,t_1)dt_1\right)^2\le
$$
\begin{equation}
\label{april67}
\le\sum\limits_{j_1=0}^{\infty}
\left(\int\limits_t^T  \phi_{j_1}(t_1)
\sum\limits_{j_2,j_3=0}^{p}
C_{j_2 j_3 j_3 j_2}(T,t_1)dt_1\right)^2=
\int\limits_t^T
\left(\sum\limits_{j_2,j_3=0}^{p}
C_{j_2 j_3 j_3 j_2}(T,t_1)
\right)^2 dt_1.
\end{equation}

\vspace{3mm}

Using Lebesgue's Dominated Convergence Theorem and (\ref{april50})--(\ref{april52}), 
(\ref{april56x})--(\ref{april56xxx}),
we obtain that the right-hand sides of (\ref{april62})--(\ref{april67}) 
tend to zero when $p\to\infty.$
The equalities (\ref{april21}), 
(\ref{april24}), (\ref{april27}),
(\ref{april33})--(\ref{april35}) are proved.

Further, let us prove the equalities (\ref{april23}), (\ref{april25}), 
(\ref{april28}), (\ref{april29}), (\ref{april31}).

Using Fubini's Theorem, Parseval's equality
and Cauchy--Bunyakovsky's inequality, we have
the following relations 
for the prelimit
expressions on the left-hand sides of 
(\ref{april23}), (\ref{april25}), 
(\ref{april28}), (\ref{april29}), (\ref{april31})
$$
\sum\limits_{j_3=0}^{p}
\left(\sum\limits_{j_1,j_4=0}^{p}
C_{j_4 j_4 j_3 j_1 j_1}-\frac{1}{4} 
C_{j_4 j_4 j_3 j_1 j_1}\biggl|_{(j_{1} j_{1})\curvearrowright (\cdot),
(j_{4} j_{4})\curvearrowright (\cdot)}
\biggr.\right)^2=
$$
$$
=\sum\limits_{j_3=0}^{p}
\left(\int\limits_t^T  \phi_{j_3}(t_3)
\sum\limits_{j_1,j_4=0}^{p}
\int\limits_t^{t_3}  \phi_{j_1}(t_2)
\int\limits_t^{t_2}  \phi_{j_1}(t_1)dt_1 dt_2
\int\limits_{t_3}^{T}  \phi_{j_4}(t_4)
\int\limits_{t_4}^{T}  \phi_{j_4}(t_5)
dt_5 dt_4 dt_3-\right.
$$
$$
\left.-\frac{1}{4}
\int\limits_{t}^{T}\int\limits_{t}^{t_4}\phi_{j_3}(t_3)\int\limits_{t}^{t_3}dt_1
dt_3 dt_4\right)^2\le
$$
$$
\le\sum\limits_{j_3=0}^{\infty}
\left(\int\limits_t^T  \phi_{j_3}(t_3)
\left(\sum\limits_{j_1,j_4=0}^{p}
\frac{1}{4}\left(\int\limits_t^{t_3}  \phi_{j_1}(t_2)dt_2\right)^2
\left(\int\limits_{t_3}^{T}  \phi_{j_4}(t_4)dt_4\right)^2\right.\right.-
$$
$$
\left.\left.
-\frac{1}{4}(t_3-t)
\int\limits_{t_3}^{T}dt_4\right) dt_3\right)^2=
$$
\begin{equation}
\label{april100}
=\int\limits_t^T
\left(
\frac{1}{4}\sum\limits_{j_1=0}^{p}
\left(\int\limits_t^{t_3}  \phi_{j_1}(t_2)dt_2\right)^2
\sum\limits_{j_4=0}^{p}
\left(\int\limits_{t_3}^{T}  \phi_{j_4}(t_4)dt_4\right)^2
-\frac{1}{4}(t_3-t)(T-t_3)\right)^2 dt_3,
\end{equation}

$$
\sum\limits_{j_4=0}^{p}
\left(\sum\limits_{j_1,j_2=0}^{p}
C_{j_2 j_4 j_1 j_2 j_1}\right)^2=
$$
$$
=\hspace{-0.4mm}\sum\limits_{j_4=0}^{p}\hspace{-0.4mm}
\left(\int\limits_t^T  \hspace{-0.3mm}\phi_{j_4}(t_4)
\hspace{-1.8mm}\sum\limits_{j_1,j_2=0}^{p}
\int\limits_t^{t_4}  \hspace{-0.3mm}\phi_{j_1}(t_3)
\int\limits_t^{t_3}  \hspace{-0.3mm}\phi_{j_2}(t_2)
\int\limits_{t}^{t_2}  \hspace{-0.3mm}\phi_{j_1}(t_1)dt_1 dt_2 dt_3
\int\limits_{t_4}^{T}  \hspace{-0.3mm}\phi_{j_2}(t_5)
dt_5 dt_4\right)^2\hspace{-2.2mm}\le
$$
$$
\le\sum\limits_{j_4=0}^{\infty}
\left(\int\limits_t^T \phi_{j_4}(t_4)
\sum\limits_{j_1,j_2=0}^{p}C_{j_1 j_2 j_1}(t_4,t)
C_{j_2}(T,t_4) dt_4\right)^2=
$$
$$
=\int\limits_t^T
\left(
\sum\limits_{j_2=0}^{p} \sum\limits_{j_1=0}^{p}C_{j_1 j_2 j_1}(t_4,t)
C_{j_2}(T,t_4)\right)^2 dt_4\le
$$
$$
\le\int\limits_t^T
\sum\limits_{j_2=0}^{p}\left(C_{j_2}(T,t_4)\right)^2 
\sum\limits_{j_2=0}^{p}\left(\sum\limits_{j_1=0}^{p}C_{j_1 j_2 j_1}(t_4,t)
\right)^2 dt_4\le
$$
$$
\le\int\limits_t^T
\sum\limits_{j_2=0}^{\infty}\left(C_{j_2}(T,t_4)\right)^2 
\sum\limits_{j_2=0}^{p}\left(\sum\limits_{j_1=0}^{p}C_{j_1 j_2 j_1}(t_4,t)
\right)^2 dt_4\le
$$
\begin{equation}
\label{april69}
\le
K_1\int\limits_t^T
\sum\limits_{j_2=0}^{p}\left(\sum\limits_{j_1=0}^{p}C_{j_1 j_2 j_1}(t_4,t)
\right)^2 dt_4\le
\end{equation}
\begin{equation}
\label{april70}
~~~~~~~~~~~~\le
K_1\int\limits_t^T
\int\limits_{t}^{t_4} \left(\sum\limits_{j_1=0}^p
\int\limits_{t}^{t_2}\phi_{j_1}(t_1)dt_1\int\limits_{t_2}^{t_4}\phi_{j_1}(t_3)dt_3\right)^2 dt_2
dt_4=
\end{equation}
\begin{equation}
\label{april101}
~~~~~~~~~~~=
K_1\int\limits_{[t,T]^2}{\bf 1}_{\{t_2<t_4\}}
\left(\sum\limits_{j_1=0}^p
\int\limits_{t}^{t_2}\phi_{j_1}(t_1)dt_1\int\limits_{t_2}^{t_4}\phi_{j_1}(t_3)dt_3\right)^2 dt_2
dt_4,
\end{equation}

\noindent
where constant $K_1$ does not depend on $p$
and the transition from (\ref{april69}) to (\ref{april70}) is based on 
(\ref{april53});

$$
\sum\limits_{j_3=0}^{p}
\left(\sum\limits_{j_1,j_2=0}^{p}
C_{j_2 j_1 j_3 j_2 j_1}\right)^2=
$$
$$
=\sum\limits_{j_3=0}^{p}
\left(\int\limits_t^T  \hspace{-0.4mm}\phi_{j_3}(t_3)
\hspace{-1.2mm}\sum\limits_{j_1,j_2=0}^{p}
\int\limits_t^{t_3}  \hspace{-0.4mm}\phi_{j_2}(t_2)
\int\limits_t^{t_2}  \hspace{-0.4mm}\phi_{j_1}(t_1)dt_1 dt_2
\int\limits_{t_3}^{T} \hspace{-0.4mm} \phi_{j_1}(t_4)
\int\limits_{t_4}^{T}  \hspace{-0.4mm}\phi_{j_2}(t_5)dt_5 dt_4 dt_3
\hspace{-0.2mm}\right)^2
\hspace{-2.2mm}\le
$$
$$
\le\sum\limits_{j_3=0}^{\infty}
\left(\int\limits_t^T  \hspace{-0.4mm}\phi_{j_3}(t_3)
\hspace{-1.2mm}\sum\limits_{j_1,j_2=0}^{p}
\int\limits_t^{t_3}  \hspace{-0.4mm}\phi_{j_2}(t_2)
\int\limits_t^{t_2}  \hspace{-0.4mm}\phi_{j_1}(t_1)dt_1 dt_2
\int\limits_{t_3}^{T} \hspace{-0.4mm} \phi_{j_1}(t_1)
\int\limits_{t_1}^{T}  \hspace{-0.4mm}\phi_{j_2}(t_2)dt_2 dt_1 dt_3
\hspace{-0.2mm}\right)^2
\hspace{-2.2mm}=
$$
$$
=\int\limits_t^T
\left(
\sum\limits_{j_1,j_2=0}^{p}
\int\limits_t^{t_3}\phi_{j_2}(t_2)
\int\limits_t^{t_2} \phi_{j_1}(t_1)dt_1 dt_2
\int\limits_{t_3}^{T} \phi_{j_1}(t_1)
\int\limits_{t_1}^{T}\phi_{j_2}(t_2)dt_2 dt_1\right)^2 dt_3=
$$
$$
=\int\limits_t^T
\left(
\sum\limits_{j_1,j_2=0}^{p}~
\int\limits_{[t,T]^2}{\bf 1}_{\{t_1<t_2<t_3\}}
\phi_{j_2}(t_2)\phi_{j_1}(t_1)dt_1 dt_2\times\right.
$$
\begin{equation}
\label{april102}
\left.\times
\int\limits_{[t,T]^2}{\bf 1}_{\{t_2>t_1>t_3\}}
\phi_{j_2}(t_2)\phi_{j_1}(t_1)dt_1 dt_2
\right)^2 dt_3,
\end{equation}
where, using the generalized Parseval equality and the 
Cauchy--Bunyakovsky inequality, we obtain
$$
\lim\limits_{p\to\infty}\sum\limits_{j_1,j_2=0}^{p}~
\int\limits_{[t,T]^2}{\bf 1}_{\{t_1<t_2<t_3\}}
\phi_{j_2}(t_2)\phi_{j_1}(t_1)dt_1 dt_2
\int\limits_{[t,T]^2}{\bf 1}_{\{t_2>t_1>t_3\}}
\phi_{j_2}(t_2)\phi_{j_1}(t_1)dt_1 dt_2\hspace{-1mm}=
$$
$$
=
\int\limits_{[t,T]^2}{\bf 1}_{\{t_1<t_2<t_3\}}{\bf 1}_{\{t_2>t_1>t_3\}}dt_1 dt_2=0,
$$
$$
\left(\sum\limits_{j_1,j_2=0}^{p}~
\int\limits_{[t,T]^2}{\bf 1}_{\{t_1<t_2<t_3\}}
\phi_{j_2}(t_2)\phi_{j_1}(t_1)dt_1 dt_2
\int\limits_{[t,T]^2}{\bf 1}_{\{t_2>t_1>t_3\}}
\phi_{j_2}(t_2)\phi_{j_1}(t_1)dt_1 dt_2\right)^2\hspace{-1mm}\le
$$
$$
\le
\sum\limits_{j_1,j_2=0}^{p}~\left(\int\limits_{[t,T]^2}{\bf 1}_{\{t_1<t_2<t_3\}}
\phi_{j_2}(t_2)\phi_{j_1}(t_1)dt_1 dt_2
\right)^2\times
$$
$$
\times
\sum\limits_{j_1,j_2=0}^{p}~\left(
\int\limits_{[t,T]^2}{\bf 1}_{\{t_2>t_1>t_3\}}
\phi_{j_2}(t_2)\phi_{j_1}(t_1)dt_1 dt_2\right)^2\le K_1<\infty,
$$

\vspace{1mm}
\noindent
where constant $K_1$ does not depend on $p;$

$$
\sum\limits_{j_2=0}^{p}
\left(\sum\limits_{j_1,j_3=0}^{p}
C_{j_3 j_1 j_3 j_2 j_1}\right)^2=
$$
$$
=\sum\limits_{j_2=0}^{p}
\left(\int\limits_t^T  \hspace{-0.4mm}\phi_{j_2}(t_2)
\hspace{-1.2mm}\sum\limits_{j_1,j_3=0}^{p}
\int\limits_t^{t_2}  \hspace{-0.4mm}\phi_{j_1}(t_1)dt_1
\int\limits_{t_2}^T  \hspace{-0.4mm}\phi_{j_3}(t_3)
\int\limits_{t_3}^{T} \hspace{-0.4mm} \phi_{j_1}(t_4)
\int\limits_{t_4}^{T}  \hspace{-0.4mm}\phi_{j_3}(t_5)dt_5 dt_4 dt_3 dt_2
\hspace{-0.2mm}\right)^2
\hspace{-2.2mm}\le
$$
$$
\le\sum\limits_{j_2=0}^{\infty}
\left(\int\limits_t^T  \hspace{-0.4mm}\phi_{j_2}(t_2)
\hspace{-1.2mm}\sum\limits_{j_1,j_3=0}^{p}
\int\limits_t^{t_2}  \hspace{-0.4mm}\phi_{j_1}(t_1)dt_1
\int\limits_{t_2}^T  \hspace{-0.4mm}\phi_{j_3}(t_3)
\int\limits_{t_3}^{T} \hspace{-0.4mm} \phi_{j_1}(t_4)
\int\limits_{t_4}^{T}  \hspace{-0.4mm}\phi_{j_3}(t_5)dt_5 dt_4 dt_3 dt_2
\hspace{-0.2mm}\right)^2
\hspace{-2.2mm}=
$$
$$
=
\int\limits_t^T 
\left(\sum\limits_{j_1,j_3=0}^{p}
\int\limits_t^{t_2}\phi_{j_1}(t_1)dt_1
\int\limits_{t_2}^T\phi_{j_3}(t_3)
\int\limits_{t_3}^{T}\phi_{j_1}(t_4)
\int\limits_{t_4}^{T}\phi_{j_3}(t_5)dt_5 dt_4 dt_3\right)^2 dt_2
=
$$
$$
=
\int\limits_t^T 
\left(\sum\limits_{j_1=0}^{p}
C_{j_1}(t_2,t)
\sum\limits_{j_3=0}^{p}
\int\limits_{t_2}^T\phi_{j_3}(t_5)
\int\limits_{t_2}^{t_5}\phi_{j_1}(t_4)
\int\limits_{t_2}^{t_4}\phi_{j_3}(t_3)dt_3 dt_4 dt_5\right)^2 dt_2
=
$$
$$
=
\int\limits_t^T 
\left(\sum\limits_{j_1=0}^{p}
C_{j_1}(t_2,t)
\sum\limits_{j_3=0}^{p}
C_{j_3 j_1 j_3}(T,t_2)
\right)^2 dt_2
\le 
$$
$$
\le
\int\limits_t^T 
\sum\limits_{j_1=0}^{p}
\left(C_{j_1}(t_2,t)\right)^2
\sum\limits_{j_1=0}^{p}\left(\sum\limits_{j_3=0}^{p}
C_{j_3 j_1 j_3}(T,t_2)
\right)^2 dt_2\le
$$
\begin{equation}
\label{april72}
\le
K_1\int\limits_t^T 
\sum\limits_{j_1=0}^{p}\left(\sum\limits_{j_3=0}^{p}
C_{j_3 j_1 j_3}(T,t_2)
\right)^2 dt_2\le
\end{equation}
\begin{equation}
\label{april73}
~~~~~~~~~~~~\le
K_1\int\limits_t^T
\int\limits_{t_2}^{T} \left(\sum\limits_{j_3=0}^p
\int\limits_{t_2}^{\theta}\phi_{j_3}(t_1)dt_1\int\limits_{\theta}^{T}\phi_{j_3}(t_3)dt_3\right)^2 d\theta
dt_2=
\end{equation}
\begin{equation}
\label{april103}
~~~~~~~~~~~=
K_1\int\limits_{[t,T]^2}{\bf 1}_{\{t_2<\theta\}}
\left(\sum\limits_{j_3=0}^p
\int\limits_{t_2}^{\theta}\phi_{j_3}(t_1)dt_1\int\limits_{\theta}^{T}\phi_{j_3}(t_3)dt_3\right)^2 
d\theta dt_2,
\end{equation}

\noindent
where constant $K_1$ does not depend on $p$
and the transition from (\ref{april72}) to (\ref{april73}) is based on 
(\ref{april53});

$$
\lim\limits_{p\to\infty}
\sum\limits_{j_3=0}^{p}
\left(\sum\limits_{j_1,j_2=0}^{p}
C_{j_1 j_2 j_3 j_2 j_1}\right)^2=
$$
$$
=\sum\limits_{j_3=0}^{p}
\left(\int\limits_t^T  \hspace{-0.4mm}\phi_{j_3}(t_3)
\hspace{-1.2mm}\sum\limits_{j_1,j_2=0}^{p}
\int\limits_t^{t_3}  \hspace{-0.4mm}\phi_{j_2}(t_2)
\int\limits_t^{t_2}  \hspace{-0.4mm}\phi_{j_1}(t_1)dt_1 dt_2
\int\limits_{t_3}^{T} \hspace{-0.4mm} \phi_{j_2}(t_4)
\int\limits_{t_4}^{T}  \hspace{-0.4mm}\phi_{j_1}(t_5)dt_5 dt_4 dt_3
\hspace{-0.2mm}\right)^2
\hspace{-2.2mm}\le
$$
$$
\le\sum\limits_{j_3=0}^{\infty}
\left(\int\limits_t^T  \hspace{-0.4mm}\phi_{j_3}(t_3)
\hspace{-1.2mm}\sum\limits_{j_1,j_2=0}^{p}
\int\limits_t^{t_3}  \hspace{-0.4mm}\phi_{j_2}(t_2)
\int\limits_t^{t_2}  \hspace{-0.4mm}\phi_{j_1}(t_1)dt_1 dt_2
\int\limits_{t_3}^{T} \hspace{-0.4mm} \phi_{j_2}(t_2)
\int\limits_{t_2}^{T}  \hspace{-0.4mm}\phi_{j_1}(t_1)dt_1 dt_2 dt_3
\hspace{-0.2mm}\right)^2
\hspace{-2.2mm}=
$$
$$
=\int\limits_t^T
\left(
\sum\limits_{j_1,j_2=0}^{p}
\int\limits_t^{t_3}\phi_{j_2}(t_2)
\int\limits_t^{t_2} \phi_{j_1}(t_1)dt_1 dt_2
\int\limits_{t_3}^{T} \phi_{j_2}(t_2)
\int\limits_{t_2}^{T}\phi_{j_1}(t_1)dt_1 dt_2\right)^2 dt_3=
$$
$$
=\int\limits_t^T
\left(
\sum\limits_{j_1,j_2=0}^{p}~
\int\limits_{[t,T]^2}{\bf 1}_{\{t_1<t_2<t_3\}}
\phi_{j_2}(t_2)\phi_{j_1}(t_1)dt_1 dt_2\times\right.
$$
\begin{equation}
\label{april104}
\left.\times
\int\limits_{[t,T]^2}{\bf 1}_{\{t_1>t_2>t_3\}}
\phi_{j_2}(t_2)\phi_{j_1}(t_1)dt_1 dt_2
\right)^2 dt_3,
\end{equation}
where, using the generalized Parseval equality and the 
Cauchy--Bunyakovsky inequality, we obtain
$$
\lim\limits_{p\to\infty}\sum\limits_{j_1,j_2=0}^{p}~
\int\limits_{[t,T]^2}{\bf 1}_{\{t_1<t_2<t_3\}}
\phi_{j_2}(t_2)\phi_{j_1}(t_1)dt_1 dt_2
\int\limits_{[t,T]^2}{\bf 1}_{\{t_1>t_2>t_3\}}
\phi_{j_2}(t_2)\phi_{j_1}(t_1)dt_1 dt_2\hspace{-1mm}=
$$
$$
=
\int\limits_{[t,T]^2}{\bf 1}_{\{t_1<t_2<t_3\}}{\bf 1}_{\{t_1>t_2>t_3\}}dt_1 dt_2=0,
$$
$$
\left(\sum\limits_{j_1,j_2=0}^{p}~
\int\limits_{[t,T]^2}{\bf 1}_{\{t_1<t_2<t_3\}}
\phi_{j_2}(t_2)\phi_{j_1}(t_1)dt_1 dt_2
\int\limits_{[t,T]^2}{\bf 1}_{\{t_1>t_2>t_3\}}
\phi_{j_2}(t_2)\phi_{j_1}(t_1)dt_1 dt_2\right)^2\hspace{-1mm}\le
$$
$$
\le
\sum\limits_{j_1,j_2=0}^{p}~\left(\int\limits_{[t,T]^2}{\bf 1}_{\{t_1<t_2<t_3\}}
\phi_{j_2}(t_2)\phi_{j_1}(t_1)dt_1 dt_2
\right)^2\times
$$
$$
\times
\sum\limits_{j_1,j_2=0}^{p}~\left(
\int\limits_{[t,T]^2}{\bf 1}_{\{t_1>t_2>t_3\}}
\phi_{j_2}(t_2)\phi_{j_1}(t_1)dt_1 dt_2\right)^2\le K_1<\infty,
$$

\vspace{1mm}
\noindent
where constant $K_1$ does not depend on $p.$

Using Lebesgue's Dominated Convergence Theorem,
we obtain that the right-hand sides of 
(\ref{april100}), (\ref{april101}), 
(\ref{april102}), (\ref{april103}), (\ref{april104})
tend to zero when $p\to\infty.$
The equalities (\ref{april23}), (\ref{april25}), 
(\ref{april28}), (\ref{april29}), (\ref{april31}) are proved.

{\bf Step~5.}\ Finally, let us prove the equalities 
(\ref{april22}), (\ref{april26}), 
(\ref{april30}), (\ref{april32}).

Using Parseval's equality,
Cauchy--Bunyakovsky's inequality, as well as Fubini's Theorem and the elementary inequality
$(a+b)^2\le 2 a^2 + 2 b^2,$
we obtain for the prelimit expression on the left-hand side of 
(\ref{april22})
$$
\sum\limits_{j_4=0}^{p}
\left(\sum\limits_{j_1,j_3=0}^{p}
C_{j_3 j_4 j_3 j_1 j_1}\right)^2=
$$
$$
=\sum\limits_{j_4=0}^{p}
\left(\int\limits_t^T  \hspace{-0.4mm}\phi_{j_4}(t_4)
\hspace{-1.2mm}\sum\limits_{j_1,j_3=0}^{p}
\int\limits_t^{t_4}  \hspace{-0.4mm}\phi_{j_3}(t_3)
\int\limits_t^{t_3}  \hspace{-0.4mm}\phi_{j_1}(t_2)
\int\limits_{t}^{t_2} \hspace{-0.4mm} \phi_{j_1}(t_1)dt_1 dt_2dt_3
\int\limits_{t_4}^{T}  \hspace{-0.4mm}\phi_{j_3}(t_5)dt_5 dt_4
\hspace{-0.2mm}\right)^2
\hspace{-2.2mm}\le
$$
$$
\le\sum\limits_{j_4=0}^{\infty}
\left(\int\limits_t^T  \hspace{-0.4mm}\phi_{j_4}(t_4)
\hspace{-1.2mm}\sum\limits_{j_1,j_3=0}^{p}
\int\limits_t^{t_4}  \hspace{-0.4mm}\phi_{j_3}(t_3)
\int\limits_t^{t_3}  \hspace{-0.4mm}\phi_{j_1}(t_2)
\int\limits_{t}^{t_2} \hspace{-0.4mm} \phi_{j_1}(t_1)dt_1 dt_2dt_3
\int\limits_{t_4}^{T}  \hspace{-0.4mm}\phi_{j_3}(t_5)dt_5 dt_4
\hspace{-0.2mm}\right)^2
\hspace{-2.2mm}=
$$
$$
=
\int\limits_t^T  
\left(\sum\limits_{j_1,j_3=0}^{p}
\int\limits_t^{t_4}\phi_{j_3}(t_3)
\int\limits_t^{t_3}\phi_{j_1}(t_2)
\int\limits_{t}^{t_2}\phi_{j_1}(t_1)dt_1 dt_2dt_3
\int\limits_{t_4}^{T}\phi_{j_3}(t_5)dt_5 
\right)^2 dt_4=
$$
$$
=
\int\limits_t^T  
\left(\sum\limits_{j_3=0}^{p}
\int\limits_t^{t_4}\phi_{j_3}(t_3)
\left(\frac{1}{2}\sum\limits_{j_1=0}^{p}\left(
\int\limits_t^{t_3}\phi_{j_1}(t_2)dt_2\right)^2 \hspace{-2mm}\mp \frac{t_3-t}{2}\right)dt_3
\int\limits_{t_4}^{T}\phi_{j_3}(t_5)dt_5 
\right)^2 \hspace{-2mm}dt_4\hspace{-0.6mm}\le
$$
$$
\le 2
\int\limits_t^T  \hspace{-1mm}
\left(\sum\limits_{j_3=0}^{p}
\int\limits_t^{t_4}\hspace{-0.4mm}\phi_{j_3}(t_3)
\left(\frac{1}{2}\sum\limits_{j_1=0}^{p}\left(
\int\limits_t^{t_3}\hspace{-0.4mm}\phi_{j_1}(t_2)dt_2\right)^2 \hspace{-2mm} -\frac{t_3-t}{2}\right)dt_3
\int\limits_{t_4}^{T}\hspace{-0.4mm}\phi_{j_3}(t_5)dt_5 
\right)^2 \hspace{-2.3mm}dt_4\hspace{-0.1mm}+
$$
$$
+2
\int\limits_t^T  
\left(\sum\limits_{j_3=0}^{p}
\int\limits_t^{t_4}\phi_{j_3}(t_3)
\frac{t_3-t}{2}dt_3
\int\limits_{t_4}^{T}\phi_{j_3}(t_5)dt_5 
\right)^2 dt_4 \le
$$
$$
\le 2
\int\limits_t^T  
\sum\limits_{j_3=0}^{p} \left(C_{j_3}(T,t_4)\right)^2\times
$$
$$
\times
\sum\limits_{j_3=0}^{p} \left(
\int\limits_t^{t_4}\phi_{j_3}(t_3)
\left(\frac{1}{2}\sum\limits_{j_1=0}^{p}\left(
\int\limits_t^{t_3}\phi_{j_1}(t_2)dt_2\right)^2 -\frac{t_3-t}{2}\right)dt_3
\right)^2 dt_4 
+\varepsilon_p \le
$$
$$
\le K_1
\int\limits_t^T  
\sum\limits_{j_3=0}^{p} \left(
\int\limits_t^{t_4}\phi_{j_3}(t_3)
\left(\frac{1}{2}\sum\limits_{j_1=0}^{p}\left(
\int\limits_t^{t_3}\phi_{j_1}(t_2)dt_2\right)^2 -\frac{t_3-t}{2}\right)dt_3
\right)^2 dt_4 + \varepsilon_p\le
$$
$$
\le K_1
\int\limits_t^T  
\sum\limits_{j_3=0}^{\infty} \left(
\int\limits_t^{t_4}\phi_{j_3}(t_3)
\left(\frac{1}{2}\sum\limits_{j_1=0}^{p}\left(
\int\limits_t^{t_3}\phi_{j_1}(t_2)dt_2\right)^2 -\frac{t_3-t}{2}\right)dt_3
\right)^2 dt_4 +\varepsilon_p=
$$
$$
= K_1
\int\limits_t^T  
\int\limits_t^{t_4}
\left(\frac{1}{2}\sum\limits_{j_1=0}^{p}\left(
\int\limits_t^{t_3}\phi_{j_1}(t_2)dt_2\right)^2 -\frac{t_3-t}{2}\right)^2 
dt_3dt_4 + \varepsilon_p=
$$
\begin{equation}
\label{april105}
~~~~~= K_1
\int\limits_{[t,T]^2}
{\bf 1}_{\{t_3<t_4\}}
\left(\frac{1}{2}\sum\limits_{j_1=0}^{p}\left(
\int\limits_t^{t_3}\phi_{j_1}(t_2)dt_2\right)^2 -\frac{t_3-t}{2}\right)^2 
dt_3dt_4 + \varepsilon_p,
\end{equation}

\noindent
where constant $K_1$ does not depend on $p,$
$$
\varepsilon_p=2\int\limits_t^T  
\left(\sum\limits_{j_3=0}^{p}
\int\limits_t^{t_4}\phi_{j_3}(t_3)
\frac{t_3-t}{2}dt_3
\int\limits_{t_4}^{T}\phi_{j_3}(t_5)dt_5 
\right)^2 dt_4.
$$

By analogy with (\ref{april47}), (\ref{april49}) we get
\begin{equation}
\label{april106}
~~~~~~~~~~~\left(\sum\limits_{j_3=0}^{p}
\int\limits_t^{t_4}\phi_{j_3}(t_3)
\frac{t_3-t}{2}dt_3
\int\limits_{t_4}^{T}\phi_{j_3}(t_5)dt_5 
\right)^2\le K_2<\infty,
\end{equation}
\begin{equation}
\label{april107}
\sum\limits_{j_3=0}^{\infty}
\int\limits_t^{t_4}\phi_{j_3}(t_3)
\frac{t_3-t}{2}dt_3
\int\limits_{t_4}^{T}\phi_{j_3}(t_5)dt_5=0,
\end{equation}
where constant $K_2$ does not depend on $p.$

Using Lebesgue's Dominated Convergence Theorem and (\ref{april46}), (\ref{april48}), 
(\ref{april106}), (\ref{april107}), 
we obtain that the right-hand side of (\ref{april105})
tends to zero when $p\to\infty.$
The equality (\ref{april22}) is proved.

Let us prove the equality (\ref{april26}).
Using Parseval's equality,
Cauchy--Bunyakovsky's inequality, as well as Fubini's Theorem and the elementary inequality
$(a+b)^2\le 2 a^2 + 2 b^2,$
we obtain for the prelimit expression on the left-hand side of 
(\ref{april26})
$$
\sum\limits_{j_2=0}^{p}
\left(\sum\limits_{j_1,j_4=0}^{p}
C_{j_4 j_4 j_1 j_2 j_1}\right)^2=
$$
$$
=\sum\limits_{j_2=0}^{p}
\left(\int\limits_t^T  \hspace{-0.4mm}\phi_{j_2}(t_2)
\hspace{-1.2mm}\sum\limits_{j_1,j_4=0}^{p}
\int\limits_t^{t_2}  \hspace{-0.4mm}\phi_{j_1}(t_1)dt_1
\int\limits_{t_2}^{T}  \hspace{-0.4mm}\phi_{j_1}(t_3)
\int\limits_{t_3}^{T} \hspace{-0.4mm} \phi_{j_4}(t_4)
\int\limits_{t_4}^{T}  \hspace{-0.4mm}\phi_{j_4}(t_5)dt_5 dt_4 dt_3dt_2
\hspace{-0.2mm}\right)^2
\hspace{-2.2mm}\le
$$
$$
\le\sum\limits_{j_2=0}^{\infty}
\left(\int\limits_t^T  \hspace{-0.4mm}\phi_{j_2}(t_2)
\hspace{-1.2mm}\sum\limits_{j_1,j_4=0}^{p}
\int\limits_t^{t_2}  \hspace{-0.4mm}\phi_{j_1}(t_1)dt_1
\int\limits_{t_2}^{T}  \hspace{-0.4mm}\phi_{j_1}(t_3)
\int\limits_{t_3}^{T} \hspace{-0.4mm} \phi_{j_4}(t_4)
\int\limits_{t_4}^{T}  \hspace{-0.4mm}\phi_{j_4}(t_5)dt_5 dt_4 dt_3dt_2
\hspace{-0.2mm}\right)^2
\hspace{-2.2mm}=
$$
$$
=\int\limits_t^T \left(\sum\limits_{j_1,j_4=0}^{p}
\int\limits_t^{t_2} \phi_{j_1}(t_1)dt_1
\int\limits_{t_2}^{T} \phi_{j_1}(t_3)
\int\limits_{t_3}^{T} \phi_{j_4}(t_4)
\int\limits_{t_4}^{T}  \phi_{j_4}(t_5)dt_5 dt_4 dt_3\right)^2 dt_2
=
$$
$$
=
\int\limits_t^T  
\left(\sum\limits_{j_1=0}^{p}
\int\limits_t^{t_2}\phi_{j_1}(t_1)dt_1
\int\limits_{t_2}^{T}\phi_{j_1}(t_3)
\hspace{-0.4mm}\left(\hspace{-0.4mm}\frac{1}{2}\sum\limits_{j_4=0}^{p}\left(
\int\limits_{t_3}^T\phi_{j_4}(t_4)dt_4\right)^2 \hspace{-2mm}\mp \frac{T-t_3}{2}\hspace{-0.4mm}\right)
\hspace{-0.4mm}dt_3 
\hspace{-0.4mm}\right)^2 \hspace{-2mm}dt_2\hspace{-0.6mm}\le
$$
$$
\le 2
\int\limits_t^T  \hspace{-1mm}
\left(\sum\limits_{j_1=0}^{p}
\int\limits_t^{t_2}\hspace{-0.4mm}\phi_{j_1}(t_1)dt_1
\int\limits_{t_2}^{T}\phi_{j_1}(t_3)
\hspace{-0.4mm}\left(\hspace{-0.4mm}\frac{1}{2}\sum\limits_{j_4=0}^{p}\left(
\int\limits_{t_3}^T\hspace{-0.4mm}\phi_{j_4}(t_4)dt_4\right)^2 \hspace{-2mm} -\frac{T-t_3}{2}
\hspace{-0.4mm}\right)\hspace{-0.4mm}dt_3
\hspace{-0.4mm}\right)^2 \hspace{-2.3mm}dt_2\hspace{-0.1mm}+
$$
$$
+2
\int\limits_t^T  
\left(\sum\limits_{j_1=0}^{p}
\int\limits_t^{t_2}\phi_{j_1}(t_1)dt_1
\int\limits_{t_2}^{T}\phi_{j_1}(t_3)\frac{T-t_3}{2}dt_3
\right)^2 dt_2 \le
$$
$$
\le 2
\int\limits_t^T  
\sum\limits_{j_1=0}^{p}\left(C_{j_1}(t_2,t)\right)^2\times
$$
$$
\times
\sum\limits_{j_1=0}^{p}\left(
\int\limits_{t_2}^{T}\phi_{j_1}(t_3)
\left(\frac{1}{2}\sum\limits_{j_4=0}^{p}\left(
\int\limits_{t_3}^T\hspace{-0.4mm}\phi_{j_4}(t_4)dt_4\right)^2 \hspace{-2mm} -\frac{T-t_3}{2}
\right)dt_3
\right)^2 dt_2+\mu_p\le
$$
$$
\le K_1
\int\limits_t^T  
\sum\limits_{j_1=0}^{p}\left(
\int\limits_{t_2}^{T}\phi_{j_1}(t_3)
\left(\frac{1}{2}\sum\limits_{j_4=0}^{p}\left(
\int\limits_{t_3}^T\hspace{-0.4mm}\phi_{j_4}(t_4)dt_4\right)^2 \hspace{-2mm} -\frac{T-t_3}{2}
\right)dt_3
\right)^2 dt_2+\mu_p\le
$$
$$
\le K_1
\int\limits_t^T  
\sum\limits_{j_1=0}^{\infty}\left(
\int\limits_{t_2}^{T}\phi_{j_1}(t_3)
\left(\frac{1}{2}\sum\limits_{j_4=0}^{p}\left(
\int\limits_{t_3}^T\hspace{-0.4mm}\phi_{j_4}(t_4)dt_4\right)^2 \hspace{-2mm} -\frac{T-t_3}{2}
\right)dt_3
\right)^2 dt_2+\mu_p=
$$
$$
= K_1
\int\limits_t^T  
\int\limits_{t_2}^{T}
\left(\frac{1}{2}\sum\limits_{j_4=0}^{p}\left(
\int\limits_{t_3}^T\hspace{-0.4mm}\phi_{j_4}(t_4)dt_4\right)^2 \hspace{-2mm} -\frac{T-t_3}{2}
\right)^2 dt_3 dt_2+\mu_p=
$$
\begin{equation}
\label{april108}
~~~~~= K_1
\int\limits_{[t,T]^2}{\bf 1}_{\{t_2<t_3\}}
\left(\frac{1}{2}\sum\limits_{j_4=0}^{p}\left(
\int\limits_{t_3}^T\hspace{-0.4mm}\phi_{j_4}(t_4)dt_4\right)^2 \hspace{-2mm} -\frac{T-t_3}{2}
\right)^2 dt_3 dt_2+\mu_p,
\end{equation}

\noindent
where constant $K_1$ does not depend on $p,$
$$
\mu_p=2
\int\limits_t^T  
\left(\sum\limits_{j_1=0}^{p}
\int\limits_t^{t_2}\phi_{j_1}(t_1)dt_1
\int\limits_{t_2}^{T}\phi_{j_1}(t_3)\frac{T-t_3}{2}dt_3
\right)^2 dt_2.
$$

By analogy with (\ref{april47}), (\ref{april49}) we get
\begin{equation}
\label{april109}
~~~~~~~~~~~\left(\sum\limits_{j_1=0}^{p}
\int\limits_t^{t_2}\phi_{j_1}(t_1)dt_1
\int\limits_{t_2}^{T}\phi_{j_1}(t_3)\frac{T-t_3}{2}dt_3
\right)^2\le K_2<\infty,
\end{equation}
\begin{equation}
\label{april110}
\sum\limits_{j_1=0}^{\infty}
\int\limits_t^{t_2}\phi_{j_1}(t_1)dt_1
\int\limits_{t_2}^{T}\phi_{j_1}(t_3)\frac{T-t_3}{2}dt_3=0,
\end{equation}
where constant $K_2$ does not depend on $p.$

Using Lebesgue's Dominated Convergence Theorem and (\ref{april46}), (\ref{april48}), 
(\ref{april109}), (\ref{april110}), 
we obtain that the right-hand side of (\ref{april108})
tends to zero when $p\to\infty.$
The equality (\ref{april26}) is proved.

Let us prove the equality (\ref{april30}).
Using Parseval's equality,
Cauchy--Bunyakovsky's inequality, as well as Fubini's Theorem and the elementary inequality
$(a+b)^2\le 2 a^2 + 2 b^2,$
we obtain for the prelimit expression on the left-hand side of 
(\ref{april30})
$$
\sum\limits_{j_4=0}^{p}
\left(\sum\limits_{j_1,j_2=0}^{p}
C_{j_1 j_4 j_2 j_2 j_1}\right)^2=
$$
$$
=\sum\limits_{j_4=0}^{p}
\left(\int\limits_t^T  \hspace{-0.4mm}\phi_{j_4}(t_4)
\hspace{-1.2mm}\sum\limits_{j_1,j_2=0}^{p}
\int\limits_t^{t_4}  \hspace{-0.4mm}\phi_{j_2}(t_3)
\int\limits_t^{t_3} \hspace{-0.4mm}\phi_{j_2}(t_2)
\int\limits_t^{t_2} \hspace{-0.4mm} \phi_{j_1}(t_1)dt_1 dt_2 dt_3
\int\limits_{t_4}^{T}  \hspace{-0.4mm}\phi_{j_1}(t_5)dt_5 dt_4 
\hspace{-0.2mm}\right)^2
\hspace{-2.2mm}\le
$$
$$
\le\sum\limits_{j_4=0}^{\infty}
\left(\int\limits_t^T  \hspace{-0.4mm}\phi_{j_4}(t_4)
\hspace{-1.2mm}\sum\limits_{j_1,j_2=0}^{p}
\int\limits_t^{t_4}  \hspace{-0.4mm}\phi_{j_2}(t_3)
\int\limits_t^{t_3} \hspace{-0.4mm}\phi_{j_2}(t_2)
\int\limits_t^{t_2} \hspace{-0.4mm} \phi_{j_1}(t_1)dt_1 dt_2 dt_3
\int\limits_{t_4}^{T}  \hspace{-0.4mm}\phi_{j_1}(t_5)dt_5 dt_4 
\hspace{-0.2mm}\right)^2
\hspace{-2.2mm}=
$$
$$
=
\int\limits_t^T  
\left(\sum\limits_{j_1,j_2=0}^{p}
\int\limits_t^{t_4}\phi_{j_2}(t_3)
\int\limits_t^{t_3}\phi_{j_2}(t_2)
\int\limits_t^{t_2}\phi_{j_1}(t_1)dt_1 dt_2 dt_3
\int\limits_{t_4}^{T}\phi_{j_1}(t_5)dt_5 
\right)^2 dt_4 =
$$
$$
=
\int\limits_t^T  
\left(\sum\limits_{j_1,j_2=0}^{p}
\int\limits_t^{t_4}\phi_{j_1}(t_1)
\int\limits_{t_1}^{t_4}\phi_{j_2}(t_2)
\int\limits_{t_2}^{t_4}\phi_{j_2}(t_3)dt_3 dt_2 dt_1
\int\limits_{t_4}^{T}\phi_{j_1}(t_5)dt_5 
\right)^2 dt_4 =
$$
$$
=
\int\limits_t^T  
\left(\sum\limits_{j_1=0}^{p}
\int\limits_t^{t_4}\phi_{j_1}(t_1)
\hspace{-0.4mm}\left(\hspace{-0.4mm}\frac{1}{2}\sum\limits_{j_2=0}^{p}\left(
\int\limits_{t_1}^{t_4}\phi_{j_2}(t_2)dt_2\right)^2 \hspace{-2mm}\mp \frac{t_4-t_1}{2}\hspace{-0.4mm}\right)
\hspace{-0.4mm}dt_1\int\limits_{t_4}^T\phi_{j_1}(t_5)dt_5
\hspace{-0.4mm}\right)^2 \hspace{-2mm}dt_4\hspace{-0.6mm}\le
$$
$$
\le 2
\int\limits_t^T  
\left(\sum\limits_{j_1=0}^{p}
\int\limits_t^{t_4}\hspace{-0.4mm}\phi_{j_1}(t_1)
\hspace{-0.6mm}\left(\hspace{-0.6mm}\frac{1}{2}\sum\limits_{j_2=0}^{p}\left(
\int\limits_{t_1}^{t_4}\hspace{-0.4mm}
\phi_{j_2}(t_2)dt_2\right)^2 \hspace{-2.3mm}-\frac{t_4-t_1}{2}\hspace{-0.6mm}\right)
\hspace{-0.6mm}dt_1\int\limits_{t_4}^T\hspace{-0.4mm}\phi_{j_1}(t_5)dt_5
\hspace{-0.6mm}\right)^2 \hspace{-2.2mm}dt_4\hspace{-0.6mm}+
$$
$$
+2\int\limits_t^T  
\left(\sum\limits_{j_1=0}^{p}
\int\limits_t^{t_4}\phi_{j_1}(t_1)
\frac{t_4-t_1}{2}
dt_1\int\limits_{t_4}^T\phi_{j_1}(t_5)dt_5
\right)^2 dt_4\le
$$
$$
\le 2
\int\limits_t^T  
\sum\limits_{j_1=0}^{p}\left(C_{j_1}(T,t_4)\right)^2\times
$$
$$
\times
\sum\limits_{j_1=0}^{p}\left(
\int\limits_t^{t_4}\phi_{j_1}(t_1)
\left(\frac{1}{2}\sum\limits_{j_2=0}^{p}\left(
\int\limits_{t_1}^{t_4}
\phi_{j_2}(t_2)dt_2\right)^2-\frac{t_4-t_1}{2}\right)
dt_1
\right)^2 dt_4+\rho_p\le
$$
$$
\le K_1
\int\limits_t^T  
\sum\limits_{j_1=0}^{p}\left(
\int\limits_t^{t_4}\phi_{j_1}(t_1)
\left(\frac{1}{2}\sum\limits_{j_2=0}^{p}\left(
\int\limits_{t_1}^{t_4}
\phi_{j_2}(t_2)dt_2\right)^2-\frac{t_4-t_1}{2}\right)
dt_1
\right)^2 dt_4+\rho_p\le
$$
$$
\le K_1
\int\limits_t^T  
\sum\limits_{j_1=0}^{\infty}\left(
\int\limits_t^{t_4}\phi_{j_1}(t_1)
\left(\frac{1}{2}\sum\limits_{j_2=0}^{p}\left(
\int\limits_{t_1}^{t_4}
\phi_{j_2}(t_2)dt_2\right)^2-\frac{t_4-t_1}{2}\right)
dt_1
\right)^2 dt_4+\rho_p=
$$
$$
=K_1
\int\limits_t^T  
\int\limits_t^{t_4}
\left(\frac{1}{2}\sum\limits_{j_2=0}^{p}\left(
\int\limits_{t_1}^{t_4}
\phi_{j_2}(t_2)dt_2\right)^2-\frac{t_4-t_1}{2}
\right)^2 dt_1 dt_4+\rho_p=
$$
\begin{equation}
\label{april111}
~~~~~=K_1
\int\limits_{[t,T]^2}
{\bf 1}_{\{t_1<t_4\}}
\left(\frac{1}{2}\sum\limits_{j_2=0}^{p}\left(
\int\limits_{t_1}^{t_4}
\phi_{j_2}(t_2)dt_2\right)^2-\frac{t_4-t_1}{2}
\right)^2 dt_1 dt_4+\rho_p,
\end{equation}

\noindent
where constant $K_1$ does not depend on $p,$
$$
\rho_p=
2\int\limits_t^T  
\left(\sum\limits_{j_1=0}^{p}
\int\limits_t^{t_4}\phi_{j_1}(t_1)
\frac{t_4-t_1}{2}
dt_1\int\limits_{t_4}^T\phi_{j_1}(t_5)dt_5
\right)^2 dt_4.
$$

By analogy with (\ref{april47}), (\ref{april49}) we get $(t_4-t_1=(t_4-t)+(t-t_1))$
\begin{equation}
\label{april112}
~~~~~~~~~~~\left(\sum\limits_{j_1=0}^{p}
\int\limits_t^{t_4}\phi_{j_1}(t_1)
\frac{t_4-t_1}{2}
dt_1\int\limits_{t_4}^T\phi_{j_1}(t_5)dt_5
\right)^2\le K_2<\infty,
\end{equation}
\begin{equation}
\label{april113}
\sum\limits_{j_1=0}^{\infty}
\int\limits_t^{t_4}\phi_{j_1}(t_1)
\frac{t_4-t_1}{2}
dt_1\int\limits_{t_4}^T\phi_{j_1}(t_5)dt_5=0,
\end{equation}
where constant $K_2$ does not depend on $p.$

Using Lebesgue's Dominated Convergence Theorem and (\ref{april46}), (\ref{april48}), 
(\ref{april112}), (\ref{april113}), 
we obtain that the right-hand side of (\ref{april111})
tends to zero when $p\to\infty.$
The equality (\ref{april30}) is proved.

Let us prove the equality (\ref{april32}).
Using Parseval's equality,
Cauchy--Bunyakovsky's inequality, as well as Fubini's Theorem and the elementary inequality
$(a+b)^2\le 2 a^2 + 2 b^2,$
we obtain for the prelimit expression on the left-hand side of 
(\ref{april32})
$$
\sum\limits_{j_2=0}^{p}
\left(\sum\limits_{j_1,j_3=0}^{p}
C_{j_1 j_3 j_3 j_2 j_1}\right)^2=
$$
$$
=\sum\limits_{j_2=0}^{p}
\left(\int\limits_t^T  \hspace{-0.4mm}\phi_{j_2}(t_2)
\hspace{-1.2mm}\sum\limits_{j_1,j_3=0}^{p}
\int\limits_t^{t_2}  \hspace{-0.4mm}\phi_{j_1}(t_1)dt_1
\int\limits_{t_2}^T \hspace{-0.4mm}\phi_{j_3}(t_3)
\int\limits_{t_3}^T \hspace{-0.4mm} \phi_{j_3}(t_4)
\int\limits_{t_4}^{T}  \hspace{-0.4mm}\phi_{j_1}(t_5)dt_5 dt_4 dt_3 dt_2
\hspace{-0.2mm}\right)^2
\hspace{-2.2mm}\le
$$
$$
\le\sum\limits_{j_2=0}^{\infty}
\left(\int\limits_t^T  \hspace{-0.4mm}\phi_{j_2}(t_2)
\hspace{-1.2mm}\sum\limits_{j_1,j_3=0}^{p}
\int\limits_t^{t_2}  \hspace{-0.4mm}\phi_{j_1}(t_1)dt_1
\int\limits_{t_2}^T \hspace{-0.4mm}\phi_{j_3}(t_3)
\int\limits_{t_3}^T \hspace{-0.4mm} \phi_{j_3}(t_4)
\int\limits_{t_4}^{T}  \hspace{-0.4mm}\phi_{j_1}(t_5)dt_5 dt_4 dt_3 dt_2
\hspace{-0.2mm}\right)^2
\hspace{-2.2mm}=
$$
$$
=\int\limits_t^T  
\left(\sum\limits_{j_1,j_3=0}^{p}
\int\limits_t^{t_2}\phi_{j_1}(t_1)dt_1
\int\limits_{t_2}^T\phi_{j_3}(t_3)
\int\limits_{t_3}^T \phi_{j_3}(t_4)
\int\limits_{t_4}^{T} \phi_{j_1}(t_5)dt_5 dt_4 dt_3 
\right)^2 dt_2=
$$
$$
=\int\limits_t^T  
\left(\sum\limits_{j_1,j_3=0}^{p}
\int\limits_t^{t_2}\phi_{j_1}(t_1)dt_1
\int\limits_{t_2}^T\phi_{j_1}(t_5)
\int\limits_{t_2}^{t_5} \phi_{j_3}(t_4)
\int\limits_{t_2}^{t_4} \phi_{j_3}(t_3)dt_3 dt_4 dt_5 
\right)^2 dt_2=
$$
$$
=
\int\limits_t^T  
\left(\sum\limits_{j_1=0}^{p}
\int\limits_t^{t_2}\phi_{j_1}(t_1)dt_1\int\limits_{t_2}^T\phi_{j_1}(t_5)
\hspace{-0.4mm}\left(\hspace{-0.4mm}\frac{1}{2}\sum\limits_{j_3=0}^{p}\left(
\int\limits_{t_2}^{t_5}\phi_{j_3}(t_4)dt_4\right)^2 \hspace{-2mm}\mp \frac{t_5-t_2}{2}\hspace{-0.4mm}\right)
dt_5
\hspace{-0.4mm}\right)^2 \hspace{-2mm}dt_2\hspace{-0.6mm}\le
$$
$$
\le 2
\int\limits_t^T  \hspace{-0.8mm}
\left(\hspace{-0.2mm}\sum\limits_{j_1=0}^{p}
\int\limits_t^{t_2}\hspace{-0.4mm}\phi_{j_1}(t_1)dt_1\int\limits_{t_2}^T\hspace{-0.4mm}\phi_{j_1}(t_5)
\hspace{-0.6mm}\left(\hspace{-0.6mm}\frac{1}{2}\sum\limits_{j_3=0}^{p}\left(
\int\limits_{t_2}^{t_5}\hspace{-0.4mm}
\phi_{j_3}(t_4)dt_4\right)^2 \hspace{-2.2mm}- \frac{t_5-t_2}{2}\hspace{-0.6mm}\right)
dt_5
\hspace{-0.6mm}\right)^2 \hspace{-2.2mm}dt_2+
$$
$$
+2\int\limits_t^T  
\left(\sum\limits_{j_1=0}^{p}
\int\limits_t^{t_2}\phi_{j_1}(t_1)dt_1\int\limits_{t_2}^T\phi_{j_1}(t_5)
\frac{t_5-t_2}{2}
dt_5
\right)^2 dt_2\le
$$
$$
\le 
2\int\limits_t^T  
\sum\limits_{j_1=0}^{p}\left(C_{j_1}(t_2,t)\right)^2\times
$$
$$
\times
\sum\limits_{j_1=0}^{p}\left(
\int\limits_{t_2}^T\phi_{j_1}(t_5)
\left(\frac{1}{2}\sum\limits_{j_3=0}^{p}\left(
\int\limits_{t_2}^{t_5}
\phi_{j_3}(t_4)dt_4\right)^2 - \frac{t_5-t_2}{2}\right)
dt_5
\right)^2 dt_2+ \chi_p\le
$$
$$
\le K_1
\int\limits_t^T  
\sum\limits_{j_1=0}^{p}\left(
\int\limits_{t_2}^T\phi_{j_1}(t_5)
\left(\frac{1}{2}\sum\limits_{j_3=0}^{p}\left(
\int\limits_{t_2}^{t_5}
\phi_{j_3}(t_4)dt_4\right)^2 - \frac{t_5-t_2}{2}\right)
dt_5
\right)^2 dt_2+ \chi_p\le
$$
$$
\le K_1
\int\limits_t^T  
\sum\limits_{j_1=0}^{\infty}\left(
\int\limits_{t_2}^T\phi_{j_1}(t_5)
\left(\frac{1}{2}\sum\limits_{j_3=0}^{p}\left(
\int\limits_{t_2}^{t_5}
\phi_{j_3}(t_4)dt_4\right)^2 - \frac{t_5-t_2}{2}\right)
dt_5
\right)^2 dt_2+ \chi_p=
$$
$$
= K_1
\int\limits_t^T  
\int\limits_{t_2}^T
\left(\frac{1}{2}\sum\limits_{j_3=0}^{p}\left(
\int\limits_{t_2}^{t_5}
\phi_{j_3}(t_4)dt_4\right)^2 - \frac{t_5-t_2}{2}
\right)^2 dt_5 dt_2+ \chi_p=
$$
\begin{equation}
\label{april114}
~~~~~= K_1
\int\limits_{[t,T]^2}
{\bf 1}_{\{t_2<t_5\}}
\left(\frac{1}{2}\sum\limits_{j_3=0}^{p}\left(
\int\limits_{t_2}^{t_5}
\phi_{j_3}(t_4)dt_4\right)^2 - \frac{t_5-t_2}{2}
\right)^2 dt_5 dt_2+ \chi_p,
\end{equation}

\noindent
where constant $K_1$ does not depend on $p,$
$$
\chi_p=
2\int\limits_t^T  
\left(\sum\limits_{j_1=0}^{p}
\int\limits_t^{t_2}\phi_{j_1}(t_1)dt_1\int\limits_{t_2}^T\phi_{j_1}(t_5)
\frac{t_5-t_2}{2}
dt_5
\right)^2 dt_2.
$$

By analogy with (\ref{april47}), (\ref{april49}) we get $(t_5-t_2=(t_5-t)+(t-t_2))$
\begin{equation}
\label{april115}
~~~~~~~~~~~\left(\sum\limits_{j_1=0}^{p}
\int\limits_t^{t_2}\phi_{j_1}(t_1)dt_1\int\limits_{t_2}^T\phi_{j_1}(t_5)
\frac{t_5-t_2}{2}
dt_5
\right)^2\le K_2<\infty,
\end{equation}
\begin{equation}
\label{april116}
\sum\limits_{j_1=0}^{\infty}
\int\limits_t^{t_2}\phi_{j_1}(t_1)dt_1\int\limits_{t_2}^T\phi_{j_1}(t_5)
\frac{t_5-t_2}{2}
dt_5
=0,
\end{equation}
where constant $K_2$ does not depend on $p.$

Using Lebesgue's Dominated Convergence Theorem and (\ref{april46}), (\ref{april48}), 
(\ref{april115}), (\ref{april116}), 
we obtain that the right-hand side of (\ref{april114})
tends to zero when $p\to\infty.$
The equality (\ref{april32}) is proved.
The equalities (\ref{april11})--(\ref{april35})
are proved.
Theorem~2.50 is proved.

\section{Expansion of Iterated Stratonovich Stochastic Integrals
of Multiplicity 3. The Case of an Ar\-bit\-ra\-ry Complete Orthonormal System of 
Functions in the Space $L_2([t,T])$ and 
Binomial Weight Functions}

In this section, we will consider a generalization of Theorems~2.45, 2.47.
Namely, we will prove the following theorem.

{\bf Theorem~2.51}\ \cite{arxiv-5}, \cite{arxiv-10}, \cite{arxiv-11}.\ {\it Suppose that
$\{\phi_j(x)\}_{j=0}^{\infty}$ is an arbitrary complete orthonormal system of 
functions in the space $L_2([t,T]).$
Then$,$ for the iterated Stra\-to\-no\-vich stochastic integral
of third multiplicity 
\begin{equation}
\label{may99}
~~~~~I_{{l_1l_2l_3}_{T,t}}^{*(i_1i_2i_3)}={\int\limits_t^{*}}^T (t_3-t)^{l_3}
{\int\limits_t^{*}}^{t_3}(t_2-t)^{l_2}
{\int\limits_t^{*}}^{t_2}(t_1-t)^{l_1}
d{\bf w}_{t_1}^{(i_1)}
d{\bf w}_{t_2}^{(i_2)}d{\bf w}_{t_3}^{(i_3)}
\end{equation}
the following expansion 
\begin{equation}
\label{may100}
I_{{l_1l_2l_3}_{T,t}}^{*(i_1i_2i_3)}=
\hbox{\vtop{\offinterlineskip\halign{
\hfil#\hfil\cr
{\rm l.i.m.}\cr
$\stackrel{}{{}_{p\to \infty}}$\cr
}} }\sum_{j_1,j_2,j_3=0}^{p}
C_{j_3 j_2 j_1}\zeta_{j_1}^{(i_1)}\zeta_{j_2}^{(i_2)}\zeta_{j_3}^{(i_3)}
\end{equation}
that converges in the mean-square sense is valid, where 
$i_1,i_2,i_3=0,1,\ldots,m;$ $l_1,l_2,l_3=0,1,2,\ldots,$
$$
C_{j_3 j_2 j_1}=\int\limits_t^T
(t_3-t)^{l_3}\phi_{j_3}(t_3)\int\limits_t^{t_3}
(t_2-t)^{l_2}
\phi_{j_2}(t_2)
\int\limits_t^{t_2}
(t_1-t)^{l_1}\phi_{j_1}(t_1)dt_1dt_2dt_3
$$
and
$$
\zeta_{j}^{(i)}=
\int\limits_t^T \phi_{j}(\tau) d{\bf w}_{\tau}^{(i)}
$$ 
are independent standard Gaussian random variables for various 
$i$ or $j$ {\rm (}in the case when $i\ne 0${\rm ),}
${\bf w}_{\tau}^{(i)}$ 
$(i=1,\ldots,m)$ are independent 
standard Wiener processes$,$
${\bf w}_{\tau}^{(0)}=\tau.$}

Note that the iterated Stratonovich stochastic integrals (\ref{may99}) are important 
for applications (see Chapter~4).

{\bf Proof.}\ According to Theorems~2.49, 2.12, we come to the conclusion that 
Theorem~2.51 will be proved if we prove the following
equalities
\begin{equation}
\label{may101}
~~~~~~~~~~\lim\limits_{p\to\infty}
\sum\limits_{j_3=0}^{p}
\left(
\frac{1}{2} 
C_{j_3 j_1 j_1}\biggl|_{(j_{1} j_{1})\curvearrowright (\cdot)}
\biggr. - \sum\limits_{j_1=0}^{p} C_{j_3 j_1 j_1}
\right)^2=0,
\end{equation}
\begin{equation}
\label{may102}
~~~~~~~~~~~\lim\limits_{p\to\infty}
\sum\limits_{j_1=0}^{p}
\left(
\frac{1}{2} 
C_{j_2 j_2 j_1}\biggl|_{(j_{2} j_{2})\curvearrowright (\cdot)}
\biggr. - \sum\limits_{j_2=0}^{p}  C_{j_2 j_2 j_1}
\right)^2=0,
\end{equation}
\begin{equation}
\label{may103}
~~~~~\lim\limits_{p\to\infty}
\sum\limits_{j_2=0}^{p}
\left(~\sum\limits_{j_1=0}^{p} 
C_{j_1 j_2 j_1}\right)^2=0.
\end{equation}

First, we prove that
\begin{equation}
\label{may103x}
~~~~~~~~~~~~\left\vert\sum\limits_{j=0}^{p}\int\limits_{t_1}^{t_2}(s-t)^{l}\phi_{j}(s)
\int\limits_{t_1}^{s}(\tau-t)^{m}\phi_{j}(\tau)d\tau ds\right\vert\le K<\infty,
\end{equation}
where $l, m=0,1,2,\ldots,$ $t\le t_1<t_2\le T,$ constant $K$ does not depend on $p, t_1, t_2.$

Using Fubini's Theorem and Parseval's equality, we have for $m>l$ $(l, m=0,1,2,\ldots)$
$$
\sum\limits_{j=0}^{p}\int\limits_{t}^{t_2}(s-t)^{l}\phi_{j}(s)
\int\limits_{t}^{s}(\tau-t)^{m}\phi_{j}(\tau)d\tau ds=
$$
$$
=\sum\limits_{j=0}^{p}\int\limits_{t}^{t_2}(s-t)^{l}\phi_{j}(s)
\int\limits_{t}^{s}(\tau-t)^{l} (\tau-t)^{m-l}\phi_{j}(\tau)d\tau ds=
$$
$$
=\sum\limits_{j=0}^{p}\int\limits_{t}^{t_2}(s-t)^{l}\phi_{j}(s)
\int\limits_{t}^{s}(\tau-t)^{l}\phi_{j}(\tau)\int\limits_t^{\tau} (\theta-t)^{m-l-1}(m-l)d\theta
d\tau ds=
$$
$$
=(m-l)\sum\limits_{j=0}^{p}\int\limits_{t}^{t_2}(\theta-t)^{m-l-1}
\int\limits_{\theta}^{t_2}(\tau-t)^{l}\phi_{j}(\tau)\int\limits_{\tau}^{t_2} 
(s-t)^{l}\phi_j(s)ds d\tau d\theta=
$$
$$
=(m-l)\int\limits_{t}^{t_2}(\theta-t)^{m-l-1}
\frac{1}{2}\sum\limits_{j=0}^{p}\left(\int\limits_{\theta}^{t_2}(\tau-t)^{l}\phi_{j}(\tau)d\tau\right)^2
d\theta\le
$$
$$
\le\frac{m-l}{2}\int\limits_{t}^{t_2}(\theta-t)^{m-l-1}
\sum\limits_{j=0}^{\infty}\left(\int\limits_{\theta}^{t_2}(\tau-t)^{l}\phi_{j}(\tau)d\tau\right)^2
d\theta=
$$
\begin{equation}
\label{may104}
~~~~~~~~~=\frac{m-l}{2}\int\limits_{t}^{t_2}(\theta-t)^{m-l-1}
\int\limits_{\theta}^{t_2}(\tau-t)^{2l}d\tau
d\theta\le K_1 < \infty,
\end{equation}
where constant $K_1$ does not depend on $p, t_2.$

For $l>m$ $(l, m=0,1,2,\ldots)$ we get
$$
\sum\limits_{j=0}^{p}\int\limits_{t}^{t_2}(s-t)^{l}\phi_{j}(s)
\int\limits_{t}^{s}(\tau-t)^{m}\phi_{j}(\tau)d\tau ds=
$$
$$
=\sum\limits_{j=0}^{p}\int\limits_{t}^{t_2}(s-t)^{l}\phi_{j}(s)ds
\int\limits_{t}^{t_2}(\tau-t)^{m}\phi_{j}(\tau)d\tau-
$$
$$
-\sum\limits_{j=0}^{p}\int\limits_{t}^{t_2}(s-t)^{l}\phi_{j}(s)
\int\limits_{s}^{t_2}(\tau-t)^{m}\phi_{j}(\tau)d\tau ds=
$$
$$
=\sum\limits_{j=0}^{p}\int\limits_{t}^{t_2}(s-t)^{l}\phi_{j}(s)ds
\int\limits_{t}^{t_2}(\tau-t)^{m}\phi_{j}(\tau)d\tau-
$$
\begin{equation}
\label{may105}
-\sum\limits_{j=0}^{p}\int\limits_{t}^{t_2}(\tau-t)^{m}\phi_{j}(\tau)
\int\limits_{t}^{\tau}(s-t)^{l}\phi_{j}(s)ds d\tau.
\end{equation}

Applying 
Cauchy--Bunyakovsky's inequality  and Parseval's equality, we obtain
$$
\left(\sum\limits_{j=0}^{p}\int\limits_{t}^{t_2}(s-t)^{l}\phi_{j}(s)ds
\int\limits_{t}^{t_2}(\tau-t)^{m}\phi_{j}(\tau)d\tau\right)^2\le
$$
$$
\le \sum\limits_{j=0}^{p}\left(\int\limits_{t}^{t_2}(s-t)^{l}\phi_{j}(s)ds\right)^2
\sum\limits_{j=0}^{p}\left(\int\limits_{t}^{t_2}(\tau-t)^{m}\phi_{j}(\tau)d\tau\right)^2\le
$$
$$
\le \sum\limits_{j=0}^{\infty}\left(\int\limits_{t}^{t_2}(s-t)^{l}\phi_{j}(s)ds\right)^2
\sum\limits_{j=0}^{\infty}\left(\int\limits_{t}^{t_2}(\tau-t)^{m}\phi_{j}(\tau)d\tau\right)^2=
$$
\begin{equation}
\label{may106}
=\int\limits_{t}^{t_2}(s-t)^{2l}ds
\int\limits_{t}^{t_2}(\tau-t)^{2m}d\tau\le K_2 < \infty,
\end{equation}
where constant $K_2$ does not depend on $p, t_2.$

Using (\ref{may104})--(\ref{may106}), we obtain
\begin{equation}
\label{may107}
~~~~~~~~~~~~~~~\left\vert
\sum\limits_{j=0}^{p}\int\limits_{t}^{t_2}(s-t)^{l}\phi_{j}(s)
\int\limits_{t}^{s}(\tau-t)^{m}\phi_{j}(\tau)d\tau ds
\right\vert\le K_3<\infty,
\end{equation}
where $l>m$ $(l, m=0,1,2,\ldots),$ constant $K_3$ does not depend on $p, t_2.$

For the case $l=m$ we get
$$
\sum\limits_{j=0}^{p}\int\limits_{t}^{t_2}(s-t)^{l}\phi_{j}(s)
\int\limits_{t}^{s}(\tau-t)^{l}\phi_{j}(\tau)d\tau ds
=
$$
$$
=\sum\limits_{j=0}^{p}\frac{1}{2}\left(\int\limits_{t}^{t_2}(s-t)^{l}\phi_{j}(s)ds\right)^2
\le
\sum\limits_{j=0}^{\infty}\frac{1}{2}\left(\int\limits_{t}^{t_2}(s-t)^{l}\phi_{j}(s)ds\right)^2=
$$
\begin{equation}
\label{may108}
=
\frac{1}{2}\int\limits_{t}^{t_2}(s-t)^{2l}ds\le K_4<\infty,
\end{equation}

\noindent
where constant $K_4$ does not depend on $p, t_2.$

Combining (\ref{may104}), (\ref{may107}),  (\ref{may108}), we have 
\begin{equation}
\label{may109}
~~~~~~~~~~~~~\left\vert
\sum\limits_{j=0}^{p}\int\limits_{t}^{t_2}(s-t)^{l}\phi_{j}(s)
\int\limits_{t}^{s}(\tau-t)^{m}\phi_{j}(\tau)d\tau ds
\right\vert\le K_5<\infty,
\end{equation}
where $l, m=0,1,2,\ldots,$ constant $K_5$ does not depend on $p, t_2.$

Note that
$$
\sum\limits_{j=0}^{p}\int\limits_{t_1}^{t_2}(s-t)^{l}\phi_{j}(s)
\int\limits_{t_1}^{s}(\tau-t)^{m}\phi_{j}(\tau)d\tau ds=
$$
$$
=\sum\limits_{j=0}^{p}\int\limits_{t}^{t_2}(s-t)^{l}\phi_{j}(s)
\int\limits_{t}^{s}(\tau-t)^{m}\phi_{j}(\tau)d\tau ds-
$$
$$
-\sum\limits_{j=0}^{p}\int\limits_{t}^{t_1}(s-t)^{l}\phi_{j}(s)
\int\limits_{t}^{s}(\tau-t)^{m}\phi_{j}(\tau)d\tau ds-
$$
\begin{equation}
\label{may110}
-\sum\limits_{j=0}^{p}\int\limits_{t_1}^{t_2}(s-t)^{l}\phi_{j}(s)ds
\int\limits_{t}^{t_1}(\tau-t)^{m}\phi_{j}(\tau)d\tau,
\end{equation}
where $l, m=0,1,2,\ldots$ and $t\le t_1<t_2\le T.$

By analogy with (\ref{may106}) we get
\begin{equation}
\label{may111}
~~~~~~~~~~~\left\vert
\sum\limits_{j=0}^{p}\int\limits_{t_1}^{t_2}(s-t)^{l}\phi_{j}(s)ds
\int\limits_{t}^{t_1}(\tau-t)^{m}\phi_{j}(\tau)d\tau\right\vert
\le K_6<\infty,
\end{equation}
where $l, m=0,1,2,\ldots,$ constant $K_6$ does not depend on $p, t_2.$

Combining (\ref{may110}), (\ref{may109}), and (\ref{may111}), we obtain (\ref{may103x}).

Let us prove (\ref{may101}). Using Parseval's equality, we have
$$
\lim\limits_{p\to\infty}\sum\limits_{j_3=0}^{p}
\left(
\frac{1}{2} 
C_{j_3 j_1 j_1}\biggl|_{(j_{1} j_{1})\curvearrowright (\cdot)}
\biggr. - \sum\limits_{j_1=0}^{p} C_{j_3 j_1 j_1}
\right)^2=
$$
$$
=
\lim\limits_{p\to\infty}\sum\limits_{j_3=0}^{p}
\left(\int\limits_t^T (\tau-t)^{l_3}
\phi_{j_3}(\tau)\left(\frac{1}{2}\int\limits_t^{\tau}(s-t)^{l_1+l_2} ds
-\right.\right.
$$
$$
\left.\left.-
\sum\limits_{j_1=0}^p
\int\limits_t^{\tau}(s-t)^{l_2}\phi_{j_1}(s)\int\limits_t^s (\theta-t)^{l_1}
\phi_{j_1}(\theta)d\theta ds\right)
d\tau\right)^2\le
$$
$$
\le
\lim\limits_{p\to\infty}\sum\limits_{j_3=0}^{\infty}
\left(\int\limits_t^T (\tau-t)^{l_3}
\phi_{j_3}(\tau)\left(\frac{1}{2}\int\limits_t^{\tau}(s-t)^{l_1+l_2} ds
-\right.\right.
$$
$$
\left.\left.-
\sum\limits_{j_1=0}^p
\int\limits_t^{\tau}(s-t)^{l_2}\phi_{j_1}(s)\int\limits_t^s (\theta-t)^{l_1}
\phi_{j_1}(\theta)d\theta ds\right)
d\tau\right)^2=
$$
$$
=\lim\limits_{p\to\infty}
\int\limits_t^T (\tau-t)^{2 l_3}
\left(\frac{1}{2}\int\limits_t^{\tau}(s-t)^{l_1+l_2} ds
-\right.
$$
\begin{equation}
\label{may112}
~~~~~~~~~~\left.-
\sum\limits_{j_1=0}^p
\int\limits_t^{\tau}(s-t)^{l_2}\phi_{j_1}(s)\int\limits_t^s (\theta-t)^{l_1}
\phi_{j_1}(\theta)d\theta ds\right)^2 d\tau.
\end{equation}

Using (\ref{after1400}), (\ref{may103x}) and
applying Lebesgue's 
Dominated Convergence Theorem in (\ref{may112}), we obtain
the equality (\ref{may101}).

Let us prove (\ref{may102}). Using Fubini's Theorem and Parseval's equality, we obtain
$$
\lim\limits_{p\to\infty}
\sum\limits_{j_1=0}^{p}
\left(
\frac{1}{2} 
C_{j_2 j_2 j_1}\biggl|_{(j_{2} j_{2})\curvearrowright (\cdot)}
\biggr. - \sum\limits_{j_2=0}^{p}  C_{j_2 j_2 j_1}
\right)^2=
$$
$$
=
\lim\limits_{p\to\infty}\sum\limits_{j_1=0}^{p}
\left(\frac{1}{2}\int\limits_t^T (s-t)^{l_2+l_3}
\int\limits_t^{s}
(\theta-t)^{l_1}\phi_{j_1}(\theta)d\theta ds
-\right.
$$
$$
\left.-
\sum\limits_{j_2=0}^p
\int\limits_t^T (s-t)^{l_3}\phi_{j_2}(s)\int\limits_t^s (\tau-t)^{l_2}
\phi_{j_2}(\tau)\int\limits_t^{\tau} (\theta-t)^{l_1}
\phi_{j_1}(\theta) d\theta d\tau ds\right)^2=
$$
$$
=
\lim\limits_{p\to\infty}\sum\limits_{j_1=0}^{p}
\left(
\int\limits_t^T (\theta-t)^{l_1}\phi_{j_1}(\theta)
\left(\frac{1}{2}
\int\limits_{\theta}^{T}
(s-t)^{l_2+l_3}ds
-\right.\right.
$$
$$
\left.\left.-
\sum\limits_{j_2=0}^p
\int\limits_{\theta}^T
(\tau-t)^{l_2}
\phi_{j_2}(\tau)
\int\limits_{\tau}^T
(s-t)^{l_3}\phi_{j_2}(s)
ds d\tau \right)d\theta \right)^2\le
$$
$$
\le
\lim\limits_{p\to\infty}\sum\limits_{j_1=0}^{\infty}
\left(
\int\limits_t^T (\theta-t)^{l_1}\phi_{j_1}(\theta)
\left(\frac{1}{2}
\int\limits_{\theta}^{T}
(s-t)^{l_2+l_3}ds
-\right.\right.
$$
$$
\left.\left.-
\sum\limits_{j_2=0}^p
\int\limits_{\theta}^T
(\tau-t)^{l_2}
\phi_{j_2}(\tau)
\int\limits_{\tau}^T
(s-t)^{l_3}\phi_{j_2}(s)
ds d\tau \right)d\theta \right)^2=
$$
$$
=
\lim\limits_{p\to\infty}
\int\limits_t^T(\theta-t)^{2 l_1}
\left(\frac{1}{2}\int\limits_{\theta}^{T}
(s-t)^{l_2+l_3}ds
-\right.
$$
$$
\left.-
\sum\limits_{j_2=0}^p
\int\limits_{\theta}^T
(\tau-t)^{l_2}
\phi_{j_2}(\tau)
\int\limits_{\tau}^T
(s-t)^{l_3}\phi_{j_2}(s)
ds d\tau\right)^2d\theta =
$$
$$
=
\lim\limits_{p\to\infty}
\int\limits_t^T(\theta-t)^{2 l_1}
\left(\frac{1}{2}\int\limits_{\theta}^{T}
(s-t)^{l_2+l_3}ds
-\right.
$$
\begin{equation}
\label{may113}
~~~~~~~~~~\left.-
\sum\limits_{j_2=0}^p
\int\limits_{\theta}^T
(s-t)^{l_3}\phi_{j_2}(s)
\int\limits_{\theta}^s
(\tau-t)^{l_2}
\phi_{j_2}(\tau)
d\tau ds\right)^2d\theta.
\end{equation}

Applying (\ref{after1400}), (\ref{may103x}) and
using Lebesgue's 
Dominated Convergence Theorem in (\ref{may113}), we get
the equality (\ref{may102}).

Let us prove (\ref{may103}). Applying Fubini's Theorem and Parseval's equality, we have
$$
\lim\limits_{p\to\infty}
\sum\limits_{j_2=0}^{p}
\left(~\sum\limits_{j_1=0}^{p} 
C_{j_1 j_2 j_1}\right)^2=
$$
$$
=\lim\limits_{p\to\infty}
\sum\limits_{j_2=0}^{p}
\left(~\sum\limits_{j_1=0}^{p} 
\int\limits_t^T (\theta-t)^{l_3} \phi_{j_1}(\theta)\int\limits_t^{\theta}
(\tau-t)^{l_2} \phi_{j_2}(\tau)\int\limits_t^{\tau} (s-t)^{l_1}\phi_{j_1}(s)ds d\tau d\theta\right)^2
\hspace{-2mm}=
$$
$$
=\lim\limits_{p\to\infty}
\sum\limits_{j_2=0}^{p}
\left(~\sum\limits_{j_1=0}^{p} 
\int\limits_t^T 
(\tau-t)^{l_2}\phi_{j_2}(\tau)\int\limits_t^{\tau} (s-t)^{l_1}\phi_{j_1}(s)ds
\int\limits_{\tau}^T
(\theta-t)^{l_3}\phi_{j_1}(\theta)d\theta d\tau \right)^2
\hspace{-2mm}\le
$$
$$
\le\lim\limits_{p\to\infty}
\sum\limits_{j_2=0}^{\infty}
\left(~\int\limits_t^T 
(\tau-t)^{l_2}\phi_{j_2}(\tau)\sum\limits_{j_1=0}^{p} 
\int\limits_t^{\tau} (s-t)^{l_1}\phi_{j_1}(s)ds
\int\limits_{\tau}^T
(\theta-t)^{l_3}\phi_{j_1}(\theta)d\theta d\tau \right)^2
\hspace{-2mm}\le
$$
\begin{equation}
\label{may114}
=\lim\limits_{p\to\infty}
\int\limits_t^T 
(\tau-t)^{2 l_2}\left(\sum\limits_{j_1=0}^{p}\int\limits_t^{\tau} (s-t)^{l_1}\phi_{j_1}(s)ds
\int\limits_{\tau}^T
(\theta-t)^{l_3}\phi_{j_1}(\theta)d\theta\right)^2 d\tau.
\end{equation}

\vspace{2mm}

Applying (\ref{dsds14fffff}), we obtain
\begin{equation}
\label{may115}
~~~~~~~~~~~~~\left\vert
\sum\limits_{j_1=0}^{p}\int\limits_t^{\tau} (s-t)^{l_1}\phi_{j_1}(s)ds
\int\limits_{\tau}^T
(\theta-t)^{l_3}\phi_{j_1}(\theta)d\theta\right\vert\le C<\infty,
\end{equation}
where constant $C$ does not depend on $p, \tau.$

Using the generalized Parseval equality, we get
\begin{equation}
\label{may116}
\sum\limits_{j_1=0}^{\infty}\int\limits_t^{\tau} (s-t)^{l_1}\phi_{j_1}(s)ds
\int\limits_{\tau}^T
(\theta-t)^{l_3}\phi_{j_1}(\theta)d\theta= 
$$
$$
=
\int\limits_t^T (s-t)^{l_1+l_3}{\bf 1}_{\{s<\tau\}}{\bf 1}_{\{s>\tau\}}ds=0.
\end{equation}

Taking into account (\ref{may115}), (\ref{may116}) and
applying Lebesgue's 
Dominated Convergence Theorem in (\ref{may114}), we obtain
the equality (\ref{may103}). Theorem~2.51 is proved.

\section{Expansion of Iterated Stratonovich Stochastic Integrals
of Multiplicity 3 Under Additional Assumptions. 
The Case of an Ar\-bit\-ra\-ry Complete Orthonormal System of 
Functions in the Space $L_2([t,T])$ and 
$\psi_1(\tau), \psi_{2}(\tau), \psi_3(\tau)\in L_2([t, T])$}

In this section, we will prove the following two theorems. 

\vspace{2mm}

{\bf Theorem~2.52}\ \cite{arxiv-5}, \cite{arxiv-10}, \cite{arxiv-11}.\  {\it Suppose that
$\{\phi_j(x)\}_{j=0}^{\infty}$ is an arbitra\-ry complete ortho\-nor\-mal system of 
functions in the space $L_2([t,T])$ and $\psi_1(\tau),$ $\psi_2(\tau), \psi_3(\tau)
\in L_2([t,T])$ are such that 
\begin{equation}
\label{novemberxxx1}
~~~~~~~~~~~~~\Biggl|\sum\limits_{j_1=0}^{p}
\int\limits_t^{s}\psi_2(\tau)\phi_{j_1}(\tau)
\int\limits_t^{\tau}\psi_1(\theta)\phi_{j_1}(\theta)
d\theta d\tau\Biggr|^2\le K<\infty,
\end{equation}
\begin{equation}
\label{novemberxxx2}
~~~~~~~~~~~~~\Biggl|\sum\limits_{j_3=0}^{p}
\int\limits_{s}^T \psi_2(\tau)\phi_{j_3}(\tau)
\int\limits_{\tau}^T \psi_3(\theta)\phi_{j_3}(\theta)d\theta d\tau\Biggr|^2
\le K<\infty
\end{equation}

\noindent
$\forall p\in {\bf N},$ where constant $K$ does not depend on $p$ and $s$ $(t\le s\le T)$.
Then$,$ for the sum $\bar J^{*}[\psi^{(3)}]_{T,t}^{(i_1 i_2 i_3)}$
$(i_1,i_2,i_3=0,1,\ldots,m)$
of iterated It\^{o} stochastic integrals 
defined by {\rm (\ref{dsds9}) $(k=3)$}
the following 
expansion 
$$
\bar J^{*}[\psi^{(3)}]_{T,t}^{(i_1 i_2 i_3)}=
\hbox{\vtop{\offinterlineskip\halign{
\hfil#\hfil\cr
{\rm l.i.m.}\cr
$\stackrel{}{{}_{p\to \infty}}$\cr
}} }\sum_{j_1,j_2,j_3=0}^{p}
C_{j_3 j_2 j_1}\zeta_{j_1}^{(i_1)}\zeta_{j_2}^{(i_2)}\zeta_{j_3}^{(i_3)}
$$

\noindent
that converges in the mean-square sense is valid, where 
$$
C_{j_3 j_2 j_1}=\int\limits_t^T \psi_3(t_3)
\phi_{j_3}(t_3)\int\limits_t^{t_3}\psi_2(t_2)
\phi_{j_2}(t_2)
\int\limits_t^{t_2}\psi_1(t_1)
\phi_{j_1}(t_1)dt_1dt_2dt_3
$$
and
$$
\zeta_{j}^{(i)}=
\int\limits_t^T \phi_{j}(\tau) d{\bf w}_{\tau}^{(i)}
$$ 
are independent standard Gaussian random variables for various 
$i$ or $j$ {\rm (}in the case when $i\ne 0${\rm ),}
${\bf w}_{\tau}^{(i)}$ 
$(i=1,\ldots,m)$ are independent 
standard Wiener processes$,$
${\bf w}_{\tau}^{(0)}=\tau.$}

\vspace{2mm}

{\bf Theorem~2.53}\ \cite{arxiv-5}, \cite{arxiv-10}, \cite{arxiv-11}.\  {\it Suppose that
$\{\phi_j(x)\}_{j=0}^{\infty}$ is an arbitrary complete ortho\-nor\-mal system of 
functions in the space $L_2([t,T])$ and $\psi_1(\tau),$ $\psi_2(\tau), \psi_3(\tau)$
are continuous functions on $[t, T].$
Furthermore$,$ let the conditions {\rm (\ref{novemberxxx1}), (\ref{novemberxxx2})}
are satisfied.
Then$,$ for the iterated Stra\-to\-no\-vich stochastic integral
of third multiplicity 
$$
{\int\limits_t^{*}}^T \psi_3(t_3)
{\int\limits_t^{*}}^{t_3}\psi_2(t_2)
{\int\limits_t^{*}}^{t_2}\psi_1(t_1)
d{\bf w}_{t_1}^{(i_1)}
d{\bf w}_{t_2}^{(i_2)}d{\bf w}_{t_3}^{(i_3)}\ \ \ (i_1,i_2,i_3=0,1,\ldots,m)
$$
the following 
expansion 
$$
{\int\limits_t^{*}}^T \hspace{-1.5mm}\psi_3(t_3)
{\int\limits_t^{*}}^{t_3}\hspace{-1.5mm}\psi_2(t_2)
{\int\limits_t^{*}}^{t_2}\hspace{-1.5mm}\psi_1(t_1)
d{\bf w}_{t_1}^{(i_1)}
d{\bf w}_{t_2}^{(i_2)}d{\bf w}_{t_3}^{(i_3)}=
\hbox{\vtop{\offinterlineskip\halign{
\hfil#\hfil\cr
{\rm l.i.m.}\cr
$\stackrel{}{{}_{p\to \infty}}$\cr
}} }\hspace{-1.5mm}\sum_{j_1,j_2,j_3=0}^{p}\hspace{-1.5mm}
C_{j_3 j_2 j_1}\zeta_{j_1}^{(i_1)}\zeta_{j_2}^{(i_2)}\zeta_{j_3}^{(i_3)}
$$

\noindent
that converges in the mean-square sense is valid, where 
notations are the same as in Theorem~{\rm 2.52.}}

\vspace{2mm}

Note that Theorem~2.53 is a simple consequence of Theorem~2.52 and Theorem~2.12 $(k=3).$
Let us prove Theorem~2.52. 

{\bf Proof.}\ First, let us note
some facts that follows from Monotone Convergence Theorem
(\cite{Pugach}, Theorem~3.5.1) and Lebesgue's Dominated Convergence Theorem.
Suppose that $\left\{g_j(x)\right\}_{j=0}^{\infty}$ is an arbitrary
sequence of real-valued measurable functions such that
\begin{equation}
\label{july1000}
\sum\limits_{j=0}^{\infty}\left|g_j(x)\right|\le K<\infty
\end{equation} 
almost everywhere on $X$ (with respect to Lebesgue's measure),
where
constant $K$ does not depend on $x.$

It is easy to see that under the above conditions
the following equality 
\begin{equation}
\label{july999}
~~~~~~~~~~~\lim\limits_{p\to\infty}\int\limits_X h^2(x)\left(\sum\limits_{j=0}^{p}g_j(x)\right)^2 dx=
\int\limits_X h^2(x)\left(\sum\limits_{j=0}^{\infty} g_j(x)\right)^2 dx
\end{equation}

\noindent
is true, where $h(x)\in L_2(X)$ (further, we put $h(x)\equiv 1$ for simplicity). 
Indeed, we have 
$$
\left\vert g_j(x)\right\vert =g_j^{+}(x)+g_j^{-}(x),\ \ \ g_j(x)=g_j^{+}(x)-g_j^{-}(x),
$$
where 
$$
g_j^{+}(x)=\max\{g_j(x), 0\}\ge 0,\ \ \ 
g_j^{-}(x)=-\min\{g_j(x), 0\}\ge 0.
$$ 

\vspace{2mm}

Moreover,
\begin{equation}
\label{july1002}
\sum\limits_{j=0}^{\infty}\left\vert g_j(x)\right\vert =
\sum\limits_{j=0}^{\infty} g_j^{+}(x)+\sum\limits_{j=0}^{\infty} g_j^{-}(x),
\end{equation}
$$
\sum\limits_{j=0}^{\infty}g_j(x) =
\sum\limits_{j=0}^{\infty} g_j^{+}(x)-\sum\limits_{j=0}^{\infty} g_j^{-}(x).
$$

\vspace{1mm}

Using (\ref{july1000}), we obtain that
the series (with non-negative terms) on the right-hand side of 
(\ref{july1002}) satisfy the condition (\ref{july1000}).
Further, using Monotone Convergence Theorem, we obtain
$$
\lim\limits_{p\to\infty}\int\limits_X \left(\sum\limits_{j=0}^{p}g_j(x)\right)^2 dx=
\lim\limits_{p\to\infty}\int\limits_X \left(
\sum\limits_{j=0}^{p} g_j^{+}(x)-\sum\limits_{j=0}^{p} g_j^{-}(x)
\right)^2 dx=
$$
$$
=
\lim\limits_{p\to\infty} \int\limits_X  \left(
\sum\limits_{j=0}^{p} g_j^{+}(x)
\right)^2 dx-
\lim\limits_{p\to\infty} 2\int\limits_X 
\sum\limits_{j=0}^{p} g_j^{+}(x)\sum\limits_{j=0}^{p} g_j^{-}(x)
dx +
$$
$$
+\lim\limits_{p\to\infty} \int\limits_X \left(
\sum\limits_{j=0}^{p} g_j^{-}(x)
\right)^2 dx=
$$
$$
=
\int\limits_X \lim\limits_{p\to\infty} \left(
\sum\limits_{j=0}^{p} g_j^{+}(x)
\right)^2 dx-
2\int\limits_X \lim\limits_{p\to\infty}
\sum\limits_{j=0}^{p} g_j^{+}(x)\sum\limits_{j=0}^{p} g_j^{-}(x)
dx +
$$
$$
+\int\limits_X \lim\limits_{p\to\infty}\left(
\sum\limits_{j=0}^{p} g_j^{-}(x)
\right)^2 dx=
$$
\begin{equation}
\label{slovo4}
=
\int\limits_X \left(
\sum\limits_{j=0}^{\infty} g_j^{+}(x)
\right)^2 dx-
2\int\limits_X 
\sum\limits_{j=0}^{\infty} g_j^{+}(x)\sum\limits_{j=0}^{\infty} g_j^{-}(x)
dx +
\int\limits_X \left(
\sum\limits_{j=0}^{\infty} g_j^{-}(x)
\right)^2 dx=
\end{equation}
$$
=
\int\limits_X \left(
\sum\limits_{j=0}^{\infty} g_j^{+}(x)-
\sum\limits_{j=0}^{\infty} g_j^{-}(x)
\right)^2 dx=
\int\limits_X \left(\sum\limits_{j=0}^{\infty} g_j(x)\right)^2 dx.
$$

\vspace{1mm}

The equality (\ref{july999}) can be obtained under another
conditions. If we replace
the condition (\ref{july1000}) with
\begin{equation}
\label{july1000aaa1}
~~~~~~~~\left|\sum\limits_{j=0}^{p}g_j(x)\right| \le K<\infty\ \ \forall p\in {\bf N}\ \ \ 
\hbox{and}\ \ \ \lim\limits_{p\to\infty}\sum\limits_{j=0}^{p}g_j(x)\ \ \hbox{exists}
\end{equation} 
almost everywhere on $X$ (with respect to Lebesgue's measure),
then by Lebesgue's Dominated Convergence Theorem
we obtain (\ref{july999}). Here
constant $K$ does not depend on $x$ and $p.$

According to Theorem~2.49, we come to the conclusion that 
Theorem~2.52 will be proved if we prove the following
equalities
\begin{equation}
~~~~~~~~~\label{july1003}
\lim\limits_{p\to\infty}
\sum\limits_{j_3=0}^{p}
\left(\frac{1}{2} 
C_{j_3 j_1 j_1}\biggl|_{(j_{1} j_{1})\curvearrowright (\cdot)}
\biggr.-
\sum\limits_{j_1=0}^{p}  C_{j_3 j_1 j_1}
\right)^2=0,
\end{equation}
\begin{equation}
\label{july1004}
~~~~~~~~~\lim\limits_{p\to\infty}
\sum\limits_{j_1=0}^{p}
\left(\frac{1}{2} 
C_{j_3 j_3 j_1}\biggl|_{(j_{3} j_{3})\curvearrowright (\cdot)}
\biggr.- \sum\limits_{j_3=0}^{p} C_{j_3 j_3 j_1}
\right)^2=0,
\end{equation}

\begin{equation}
\label{july1005}
\lim\limits_{p\to\infty}
\sum\limits_{j_2=0}^{p}
\left(\sum\limits_{j_1=0}^{p} 
C_{j_1 j_2 j_1}\right)^2=0.
\end{equation}

\vspace{2mm}

Let us prove (\ref{july1003}). Using Parseval's equality, we have
$$
\lim\limits_{p\to\infty}
\sum\limits_{j_3=0}^{p}
\left(\frac{1}{2} 
C_{j_3 j_1 j_1}\biggl|_{(j_{1} j_{1})\curvearrowright (\cdot)}
\biggr.-
\sum\limits_{j_1=0}^{p}  C_{j_3 j_1 j_1}
\right)^2=
$$
$$
=\lim\limits_{p\to\infty}
\sum\limits_{j_3=0}^{p}
\left(\int\limits_t^T \psi_3(s)
\phi_{j_3}(s)\left(\frac{1}{2}\int\limits_t^{s} \psi_2(\tau)\psi_1(\tau)d\tau
-\right.\right.
$$
$$
\left.\left.-
\sum\limits_{j_1=0}^p
\int\limits_t^{s}\psi_2(\tau)\phi_{j_1}(\tau)
\int\limits_t^{\tau}\psi_1(\theta)\phi_{j_1}(\theta)
d\theta d\tau\right)ds\right)^2\le
$$
$$
\le\lim\limits_{p\to\infty}
\sum\limits_{j_3=0}^{\infty}
\left(\int\limits_t^T \psi_3(s)
\phi_{j_3}(s)\left(\frac{1}{2}\int\limits_t^{s} \psi_2(\tau)\psi_1(\tau)d\tau
-\right.\right.
$$
$$
\left.\left.-
\sum\limits_{j_1=0}^p
\int\limits_t^{s}\psi_2(\tau)\phi_{j_1}(\tau)
\int\limits_t^{\tau}\psi_1(\theta)\phi_{j_1}(\theta)
d\theta d\tau\right)ds\right)^2=
$$
$$
=\lim\limits_{p\to\infty}
\int\limits_t^T \psi_3^2(s)
\left(\frac{1}{2}\int\limits_t^{s} \psi_2(\tau)\psi_1(\tau)d\tau
-\right.
$$
\begin{equation}
\label{july1006}
~~~~~~~~~\left.
-\sum\limits_{j_1=0}^p
\int\limits_t^{s}\psi_2(\tau)\phi_{j_1}(\tau)
\int\limits_t^{\tau}\psi_1(\theta)\phi_{j_1}(\theta)
d\theta d\tau \right)^2 ds=
\end{equation}
$$
=
\int\limits_t^T 
\psi_3^2(s) \lim\limits_{p\to\infty}\left(\frac{1}{2}\int\limits_t^{s} \psi_2(\tau)\psi_1(\tau)d\tau
-\right.
$$
\begin{equation}
\label{july1007}
~~~~~~~~\left.-
\sum\limits_{j_1=0}^p
\int\limits_t^{s}\psi_2(\tau)\phi_{j_1}(\tau)
\int\limits_t^{\tau}\psi_1(\theta)\phi_{j_1}(\theta)
d\theta d\tau \right)^2 ds=0,
\end{equation}

\noindent
where (\ref{july1007}) follows from from (\ref{start1000}) (also see (\ref{after1400})) and
the transition from (\ref{july1006}) to (\ref{july1007}) 
is based on (\ref{july999}), (\ref{july1000aaa1}) and
Lebesgue's Dominated Convergence Theorem (see (\ref{novemberxxx1})).
The equality (\ref{july1003}) is proved.

Let us prove (\ref{july1004}). Using Fubini's Theorem and Parseval's equality, we obtain
$$
\lim\limits_{p\to\infty}
\sum\limits_{j_1=0}^{p}
\left(\frac{1}{2} 
C_{j_3 j_3 j_1}\biggl|_{(j_{3} j_{3})\curvearrowright (\cdot)}
\biggr.- \sum\limits_{j_3=0}^{p} C_{j_3 j_3 j_1}
\right)^2=
$$
$$
=\lim\limits_{p\to\infty}
\sum\limits_{j_1=0}^{p}
\left(\frac{1}{2}\int\limits_t^T \psi_3(\tau)\psi_2(\tau)
\int\limits_t^{\tau} \psi_1(s)\phi_{j_1}(s)dsd\tau
-\right.
$$
$$
-\left.
\sum\limits_{j_3=0}^p
\int\limits_t^T \psi_3(\theta)\phi_{j_3}(\theta)\int\limits_t^{\theta} \psi_2(\tau)\phi_{j_3}(\tau)
\int\limits_t^{\tau} \psi_1(s)\phi_{j_1}(s)ds d\tau d\theta\right)^2=
$$
$$
=\lim\limits_{p\to\infty}
\sum\limits_{j_1=0}^{p}
\left(\frac{1}{2}\int\limits_t^T 
\psi_1(s)\phi_{j_1}(s)\int\limits_s^T \psi_3(\tau)\psi_2(\tau) d\tau ds
-\right.
$$
$$
\left.
-\sum\limits_{j_3=0}^p
\int\limits_t^T \psi_1(s)\phi_{j_1}(s)\int\limits_{s}^T \psi_2(\tau)\phi_{j_3}(\tau)
\int\limits_{\tau}^T \psi_3(\theta)\phi_{j_3}(\theta)d\theta d\tau ds\right)^2=
$$
$$
=\lim\limits_{p\to\infty}
\sum\limits_{j_1=0}^{p}
\left(\int\limits_t^T 
\psi_1(s)\phi_{j_1}(s)\left(\frac{1}{2}\int\limits_s^T \psi_3(\tau)\psi_2(\tau) d\tau 
-\right.\right.
$$
$$
\left.\left.
-\sum\limits_{j_3=0}^p
\int\limits_{s}^T \psi_2(\tau)\phi_{j_3}(\tau)
\int\limits_{\tau}^T \psi_3(\theta)\phi_{j_3}(\theta)d\theta d\tau
\right)ds\right)^2\le
$$
$$
\le\lim\limits_{p\to\infty}
\sum\limits_{j_1=0}^{\infty}
\left(\int\limits_t^T 
\psi_1(s)\phi_{j_1}(s)\left(\frac{1}{2}\int\limits_s^T \psi_3(\tau)\psi_2(\tau) d\tau 
-\right.\right.
$$
$$
\left.\left.
-\sum\limits_{j_3=0}^p
\int\limits_{s}^T \psi_2(\tau)\phi_{j_3}(\tau)
\int\limits_{\tau}^T \psi_3(\theta)\phi_{j_3}(\theta)d\theta d\tau
\right)ds\right)^2=
$$
$$
=\lim\limits_{p\to\infty}
\int\limits_t^T 
\psi_1^2(s)\left(\frac{1}{2}\int\limits_s^T \psi_3(\tau)\psi_2(\tau) d\tau 
-\right.
$$
\begin{equation}
\label{july1008}
~~~~~~~~~\left.
-\sum\limits_{j_3=0}^p
\int\limits_{s}^T \psi_2(\tau)\phi_{j_3}(\tau)
\int\limits_{\tau}^T \psi_3(\theta)\phi_{j_3}(\theta)d\theta d\tau
\right)^2 ds=
\end{equation}
$$
=
\int\limits_t^T 
\psi_1^2(s) \lim\limits_{p\to\infty}\left(\frac{1}{2}\int\limits_s^T \psi_3(\tau)\psi_2(\tau) d\tau 
-\right.
$$
\begin{equation}
\label{july1009}
\left.
~~~~~~~~~~~~-\sum\limits_{j_3=0}^p
\int\limits_{s}^T \psi_2(\tau)\phi_{j_3}(\tau)
\int\limits_{\tau}^T \psi_3(\theta)\phi_{j_3}(\theta)d\theta d\tau
\right)^2 ds=0,
\end{equation}

\noindent
where (\ref{july1009}) follows from (\ref{start1000}) and
the transition from (\ref{july1008}) to (\ref{july1009}) 
is based on (\ref{july999}), (\ref{july1000aaa1}) and
Lebesgue's Dominated Convergence Theorem (see (\ref{novemberxxx2})).
The equality (\ref{july1004}) is proved.

Let us prove (\ref{july1005}).
Applying Fubini's Theorem and Parseval's equality, we have
$$
\lim\limits_{p\to\infty}
\sum\limits_{j_2=0}^{p}
\left(\sum\limits_{j_1=0}^{p} 
C_{j_1 j_2 j_1}\right)^2=
$$
$$
=\lim\limits_{p\to\infty}
\sum\limits_{j_2=0}^{p}
\left(\sum\limits_{j_1=0}^{p} 
\int\limits_t^T \psi_3(\theta)\phi_{j_1}(\theta)\int\limits_t^{\theta}
\psi_2(\tau)\phi_{j_2}(\tau)\int\limits_t^{\tau} \psi_1(s)\phi_{j_1}(s)ds d\tau d\theta\right)^2=
$$
$$
=\lim\limits_{p\to\infty}
\sum\limits_{j_2=0}^{p}
\left(\sum\limits_{j_1=0}^{p} 
\int\limits_t^T 
\psi_2(\tau)\phi_{j_2}(\tau)\int\limits_t^{\tau} \psi_1(s)\phi_{j_1}(s)ds
\int\limits_{\tau}^T
\psi_3(\theta)\phi_{j_1}(\theta)d\theta d\tau \right)^2\le
$$
$$
\le\lim\limits_{p\to\infty}
\sum\limits_{j_2=0}^{\infty}
\left(
\int\limits_t^T 
\psi_2(\tau)\phi_{j_2}(\tau) \sum\limits_{j_1=0}^{p} \int\limits_t^{\tau} \psi_1(s)\phi_{j_1}(s)ds
\int\limits_{\tau}^T
\psi_3(\theta)\phi_{j_1}(\theta)d\theta d\tau \right)^2=
$$
\begin{equation}
\label{july1010}
~~~~~~~=\lim\limits_{p\to\infty}
\int\limits_t^T 
\psi_2^2(\tau)\left(\sum\limits_{j_1=0}^{p} \int\limits_t^{\tau} \psi_1(s)\phi_{j_1}(s)ds
\int\limits_{\tau}^T
\psi_3(\theta)\phi_{j_1}(\theta)d\theta \right)^2 d\tau =
\end{equation}
\begin{equation}
\label{july1011}
~~~~~~~=
\int\limits_t^T 
\psi_2^2(\tau)\lim\limits_{p\to\infty}\left(\sum\limits_{j_1=0}^{p} 
\int\limits_t^{\tau} \psi_1(s)\phi_{j_1}(s)ds
\int\limits_{\tau}^T
\psi_3(\theta)\phi_{j_1}(\theta)d\theta \right)^2 d\tau = 0,
\end{equation}

\vspace{2mm}
\noindent
where (\ref{july1011}) follows from the equality
\begin{equation}
\label{july1012}
\sum\limits_{j_1=0}^{\infty} 
\int\limits_t^{\tau} \psi_1(s)\phi_{j_1}(s)ds
\int\limits_{\tau}^T
\psi_3(\theta)\phi_{j_1}(\theta)d\theta=
\int\limits_t^T \psi_1(s){\bf 1}_{\{s<\tau\}}\psi_3(s){\bf 1}_{\{s>\tau\}}ds=0
\end{equation}

\noindent
(the relation (\ref{july1012}) follows from the generalized Parseval equality)
and
the transition from (\ref{july1010}) to (\ref{july1011}) 
is based on (\ref{july999}), (\ref{july1000aaa1}) and
Lebesgue's Dominated Convergence Theorem (see (\ref{dsds14fffff})).
The equality (\ref{july1005}) is proved.
Theorem~2.52 is proved.

\section{Expansion of Iterated Stratonovich Stochastic Integrals
of Multiplicities 4 and 5 Under Additional Assumptions. 
The Case of an Ar\-bit\-ra\-ry Complete Orthonormal System of 
Functions in the Space $L_2([t,T])$ and 
$\psi_1(\tau),\ldots,\psi_5(\tau)\in L_2([t, T])$}

Let us develop the approach discussed in the previous section.
It is easy to see (according to Theorem~2.49) that analogues of Theorems~2.52 and 2.53
for the cases $k=4$ and $k=5$ ($\psi_1(\tau), \ldots, \psi_5(\tau)
\in L_2([t,T])$) will be true if the relations
(\ref{2023novem200})--(\ref{2023novem205}), 
(\ref{april11})--(\ref{april35}) 
as well as the equalities
\begin{equation}
\label{july13}
~~~~~~~~\lim\limits_{p\to\infty}
\sum\limits_{j_1, j_3=0}^{p}
C_{j_3 j_3 j_1 j_1}=
\frac{1}{4}\int\limits_{t}^{T} 
\psi_4(t_3)\psi_3(t_3)\int\limits_{t}^{t_3}
\psi_2(t_1)\psi_1(t_1)dt_1 dt_3,
\biggr.
\end{equation}
\begin{equation}
\label{july14}
\lim\limits_{p\to\infty}
\sum\limits_{j_1, j_3=0}^{p}
C_{j_1 j_3 j_3 j_1}=0,
\biggr.
\end{equation}
\begin{equation}
\label{july15}
\lim\limits_{p\to\infty}
\sum\limits_{j_1, j_2=0}^{p}
C_{j_2 j_1 j_2 j_1}=0,
\biggr.
\end{equation}
\begin{equation}
\label{july16}
~~~~~~~~\lim\limits_{p\to\infty}
\sum\limits_{j_1, j_3=0}^{p}
C_{j_3 j_3 j_1 j_1}(s,\tau)
=\frac{1}{4}\int\limits_{\tau}^{s} \psi_4(t_3)\psi_3(t_3)\int\limits_{\tau}^{t_3}
\psi_2(t_1)\psi_1(t_1)dt_1 dt_3,
\biggr.
\end{equation}
\begin{equation}
\label{july17}
\lim\limits_{p\to\infty}
\sum\limits_{j_1, j_3=0}^{p}
C_{j_1 j_3 j_3 j_1}(s,\tau)=0,
\biggr.
\end{equation}
\begin{equation}
\label{july18}
\lim\limits_{p\to\infty}
\sum\limits_{j_1, j_2=0}^{p}
C_{j_2 j_1 j_2 j_1}(s,\tau)=0
\biggr.
\end{equation}

\vspace{2mm}
\noindent
are satisfied, provided that 
$\{\phi_j(x)\}_{j=0}^{\infty}$ is an arbitrary complete ortho\-nor\-mal system of 
functions in the space $L_2([t,T]),$ $\psi_1(\tau), \ldots, \psi_5(\tau)
\in L_2([t,T]),$ the series on the left-hand sides of (\ref{july13})--(\ref{july18}) 
converge absolutely, and 
$$
C_{j_4 \ldots j_1}=\int\limits_t^T
\psi_4(t_4)\phi_{j_4}(t_4)\ldots 
\int\limits_t^{t_2}
\psi_1(t_1)\phi_{j_1}(t_1)dt_1\ldots dt_4,
$$
$$
C_{j_5 \ldots j_1}=\int\limits_t^T
\psi_5(t_5)\phi_{j_5}(t_5)\ldots 
\int\limits_t^{t_2}
\psi_1(t_1)\phi_{j_1}(t_1)dt_1\ldots dt_5,
$$
$$
C_{j_4 \ldots j_1}(s,\tau)=\int\limits_{\tau}^s
\psi_4(t_4)\phi_{j_4}(t_4)\ldots 
\int\limits_{\tau}^{t_2}
\psi_1(t_1)\phi_{j_1}(t_1)dt_1\ldots dt_4\ \ \ (t\le \tau<s\le T)
$$

\vspace{1mm}
\noindent
in (\ref{2023novem200})--(\ref{2023novem205}), 
(\ref{april11})--(\ref{april35}), (\ref{july13})--(\ref{july18}).

It is obvious that the equalities (\ref{july16})--(\ref{july18}) follow from the equalities 
(\ref{july13})--(\ref{july15})
if in (\ref{july13})--(\ref{july15}) we replace $\psi_4(t_4), 
\psi_3(t_3), \psi_2(t_2), \psi_1(t_1)$ with 
${\bf 1}_{\{\tau<t_4<s\}}\psi_4(t_4),$
${\bf 1}_{\{\tau<t_3\}}\psi_3(t_3),$
${\bf 1}_{\{\tau<t_2\}}\psi_2(t_2),$
${\bf 1}_{\{\tau<t_1\}}\psi_1(t_1),$
respectively.

Further, the proofs of Theorems~2.46 and 2.50 must be modified 
and carried out by analogy with the proof 
of Theorem~2.52, i.e. using the equalities (\ref{july999}), (\ref{july1000aaa1}) and 
Lebesgue's Dominated Convergence Theorem.
At that, the de\-ri\-va\-tion of formulas similar to (\ref{2023novem209})--(\ref{2023novem214}),
(\ref{april36})--(\ref{april45}), (\ref{april62})--(\ref{april100}),
(\ref{april101}), (\ref{april102}), (\ref{april103}),
(\ref{april104}), (\ref{april105}), (\ref{april108}), (\ref{april111}), (\ref{april114})
is carried out completely similarly to (\ref{2023novem209})--(\ref{2023novem214}),
(\ref{april36})--(\ref{april45}), (\ref{april62})--(\ref{april100}),
(\ref{april101}), (\ref{april102}), (\ref{april103}),
(\ref{april104}), (\ref{april105}), (\ref{april108}), (\ref{april111}), (\ref{april114}), adjusted
for the fact that in (\ref{2023novem209})--(\ref{2023novem214}),
(\ref{april36})--(\ref{april45}), (\ref{april62})--(\ref{april100}),
(\ref{april101}), (\ref{april102}), (\ref{april103}),
(\ref{april104}), (\ref{april105}), (\ref{april108}), (\ref{april111}), (\ref{april114}) 
the functions $\psi_1(\tau),\ldots,\psi_5(\tau)\equiv 1$ are replaced by 
$\psi_1(\tau),\ldots,\psi_5(\tau)\in L_2([t, T])$.

Furthermore, the following additional conditions 
\begin{equation}
\label{novemberxxx3}
\left|\sum\limits_{j=0}^{p}
\int\limits_{\tau}^{s}\psi_{m+1}(t_2)\phi_{j}(t_2)
\int\limits_{\tau}^{t_2}\psi_m(t_1)\phi_{j}(t_1)
dt_1 dt_2\right|^2\le K<\infty\ \ \ (m=1,2,3,4),
\end{equation}
\begin{equation}
\label{novemberxxx4}
~~~~~~~\left|\sum\limits_{j_1,j_2=0}^{p}
C_{j_2 j_2 j_1 j_1}^{\psi_{m+3}\psi_{m+2}\psi_{m+1}\psi_m}(s,\tau)\right|^2\le K<\infty\ \ \ (m=1,2),
\end{equation}

\vspace{-3mm}
\begin{equation}
\label{novemberxxx5}
~~~~~~~\left|\sum\limits_{j_1,j_2=0}^{p}
C_{j_2 j_1 j_2 j_1}^{\psi_{m+3}\psi_{m+2}\psi_{m+1}\psi_m}(s,\tau)\right|^2\le K<\infty\ \ \ (m=1,2),
\end{equation}

\vspace{-3mm}
\begin{equation}
\label{novemberxxx6}
~~~~~~~\left|\sum\limits_{j_1,j_2=0}^{p}
C_{j_1 j_2 j_2 j_1}^{\psi_{m+3}\psi_{m+2}\psi_{m+1}\psi_m}(s,\tau)\right|^2\le K<\infty\ \ \ (m=1,2),
\end{equation}

\vspace{1mm}
\noindent
must be satisfied $\forall p\in {\bf N},$ where constant $K$ does not depend on $p, \tau, s$, 
$$
C_{j_4 j_3 j_2 j_1}^{\psi_{m+3}\psi_{m+2}\psi_{m+1}\psi_m}(s,\tau)=
\int\limits_{\tau}^{s}\psi_{m+3}(t_4)\phi_{j_4}(t_4)\times
$$
$$
\times
\int\limits_{\tau}^{t_4}\psi_{m+2}(t_3)\phi_{j_3}(t_3)
\int\limits_{\tau}^{t_3}\psi_{m+1}(t_2)\phi_{j_2}(t_2)
\int\limits_{\tau}^{t_2}\psi_{m}(t_1)\phi_{j_1}(t_1)dt_1 dt_2 dt_3 dt_4,
$$

\noindent
where $m=1, 2$ and $t\le \tau<s\le T.$

The conditions (\ref{novemberxxx3})--(\ref{novemberxxx6}) are
required to perform the passage to the limit
using Lebesgue's Dominated Convergence Theorem
(see the proofs of Theorems~2.50, 2.52 for details).

The equality (\ref{july13}) is proved in \cite{rybakov7000x}
for the case when $\{\phi_j(x)\}_{j=0}^{\infty}$ is an arbitrary complete ortho\-nor\-mal system of 
functions in the space $L_2([t,T])$ and $\psi_1(\tau),$ $\ldots, \psi_4(\tau)
\in L_2([t,T]).$ 
The equalities (\ref{july14}), (\ref{july15}) can also be obtained
\cite{rybakov7000xa}
using the approach from \cite{rybakov7000x}. At that, the series on the left-hand
sides of (\ref{july13})--(\ref{july15}) converge absolutly.
We will return to these issues in Sect.~2.27.3 and 2.27.4. 
Sect.~2.27.3 will be devoted to the method from \cite{rybakov7000x}
based on trace class operators.
In Sect.~2.27.4 we will prove the equalities 
(\ref{july13})--(\ref{july15}) using an approach based on the 
generalized Parseval equality and (\ref{start1000})
(the case when $\{\phi_j(x)\}_{j=0}^{\infty}$ is an arbitrary complete ortho\-nor\-mal system of 
functions in the space $L_2([t,T])$ and $\psi_1(\tau),$ $\ldots, \psi_4(\tau)
\in L_2([t,T])).$

Taking into account everything said above in this section, 
we obtain the following four theorems.

{\bf Theorem~2.54}\ \cite{arxiv-5}, \cite{arxiv-10}, \cite{arxiv-11}.\ {\it Suppose that
$\{\phi_j(x)\}_{j=0}^{\infty}$ is an arbitrary complete ortho\-nor\-mal system of 
functions in the space $L_2([t,T])$ and $\psi_1(\tau),$ $\ldots, \psi_4(\tau)
\in L_2([t,T]).$ Furthermore, let the condition
{\rm (\ref{novemberxxx3})} $(m=1,2,3)$
is satisfied.
Then$,$ for the sum $\bar J^{*}[\psi^{(4)}]_{T,t}^{(i_1 \ldots i_4)}$
$(i_1,\ldots,i_4=0,1,\ldots,m)$
of iterated It\^{o} stochastic integrals 
defined by {\rm (\ref{dsds9}) $(k=4)$}
the following 
expansion 
$$
\bar J^{*}[\psi^{(4)}]_{T,t}^{(i_1 \ldots i_4)}=
\hbox{\vtop{\offinterlineskip\halign{
\hfil#\hfil\cr
{\rm l.i.m.}\cr
$\stackrel{}{{}_{p\to \infty}}$\cr
}} }\sum_{j_1,\ldots,j_4=0}^{p}
C_{j_4 \ldots j_1}\zeta_{j_1}^{(i_1)}\ldots \zeta_{j_4}^{(i_4)}
$$

\vspace{1mm}
\noindent
that converges in the mean-square sense is valid, where 
$$
C_{j_4 \ldots j_1}=\int\limits_t^T \psi_4(t_4)\phi_{j_4}(t_4)
\ldots
\int\limits_t^{t_2}\psi_1(t_1)
\phi_{j_1}(t_1)dt_1\ldots dt_4
$$
and
$$
\zeta_{j}^{(i)}=
\int\limits_t^T \phi_{j}(\tau) d{\bf w}_{\tau}^{(i)}
$$ 

\noindent
are independent standard Gaussian random variables for various 
$i$ or $j$ {\rm (}when $i\ne 0${\rm ),}
${\bf w}_{\tau}^{(i)}$ 
$(i=1,\ldots,m)$ are independent 
standard Wiener processes$,$
${\bf w}_{\tau}^{(0)}=\tau.$}

{\bf Theorem~2.55}\ \cite{arxiv-5}, \cite{arxiv-10}, \cite{arxiv-11}.\ {\it Suppose that
$\{\phi_j(x)\}_{j=0}^{\infty}$ is an arbitrary complete ortho\-nor\-mal system of 
functions in the space $L_2([t,T])$ and $\psi_1(\tau),$ $\ldots, \psi_4(\tau)$
are continuous functions on $[t, T].$
Furthermore, let the condition {\rm (\ref{novemberxxx3})} $(m=1,2,3)$
is satisfied.
Then$,$ for the iterated Stra\-to\-no\-vich stochastic integral
of fourth multiplicity 
$$
{\int\limits_t^{*}}^T \psi_4(t_4)
\ldots 
{\int\limits_t^{*}}^{t_2}\psi_1(t_1)
d{\bf w}_{t_1}^{(i_1)}
\ldots d{\bf w}_{t_4}^{(i_4)}\ \ \ (i_1,\ldots,i_4=0,1,\ldots,m)
$$

\noindent
the following 
expansion 
$$
{\int\limits_t^{*}}^T \psi_4(t_4)
\ldots 
{\int\limits_t^{*}}^{t_2}\psi_1(t_1)
d{\bf w}_{t_1}^{(i_1)}
\ldots d{\bf w}_{t_4}^{(i_4)}=
\hbox{\vtop{\offinterlineskip\halign{
\hfil#\hfil\cr
{\rm l.i.m.}\cr
$\stackrel{}{{}_{p\to \infty}}$\cr
}} }\sum_{j_1,\ldots,j_4=0}^{p}
C_{j_4 \ldots j_1}\zeta_{j_1}^{(i_1)}\ldots \zeta_{j_4}^{(i_4)}
$$

\noindent
that converges in the mean-square sense is valid, where 
notations are the same as in Theorem~{\rm 2.54.}}

{\bf Theorem~2.56}\ \cite{arxiv-5}, \cite{arxiv-10}, \cite{arxiv-11}.\ {\it Suppose that
$\{\phi_j(x)\}_{j=0}^{\infty}$ is an arbitrary complete ortho\-nor\-mal system of 
functions in the space $L_2([t,T])$ and $\psi_1(\tau),$ $\ldots, \psi_5(\tau)
\in L_2([t,T]).$ Furthermore, let the conditions
{\rm (\ref{novemberxxx3})--(\ref{novemberxxx6})} are satisfied.
Then$,$ for the sum $\bar J^{*}[\psi^{(5)}]_{T,t}^{(i_1 \ldots i_5)}$
$(i_1,\ldots,i_5=0,1,\ldots,m)$
of iterated It\^{o} stochastic integrals 
defined by {\rm (\ref{dsds9}) $(k=5)$}
the following 
expansion 
$$
\bar J^{*}[\psi^{(5)}]_{T,t}^{(i_1 \ldots i_5)}=
\hbox{\vtop{\offinterlineskip\halign{
\hfil#\hfil\cr
{\rm l.i.m.}\cr
$\stackrel{}{{}_{p\to \infty}}$\cr
}} }\sum_{j_1,\ldots,j_5=0}^{p}
C_{j_5 \ldots j_1}\zeta_{j_1}^{(i_1)}\ldots \zeta_{j_5}^{(i_5)}
$$

\noindent
that converges in the mean-square sense is valid, where 
$$
C_{j_5 \ldots j_1}=\int\limits_t^T \psi_5(t_5)\phi_{j_5}(t_5)
\ldots
\int\limits_t^{t_2}\psi_1(t_1)
\phi_{j_1}(t_1)dt_1\ldots dt_5
$$
and
$$
\zeta_{j}^{(i)}=
\int\limits_t^T \phi_{j}(\tau) d{\bf w}_{\tau}^{(i)}
$$ 

\noindent
are independent standard Gaussian random variables for various 
$i$ or $j$ {\rm (}when $i\ne 0${\rm ),}
${\bf w}_{\tau}^{(i)}$ 
$(i=1,\ldots,m)$ are independent 
standard Wiener processes$,$
${\bf w}_{\tau}^{(0)}=\tau.$}

{\bf Theorem~2.57}\ \cite{arxiv-5}, \cite{arxiv-10}, \cite{arxiv-11}.\ {\it Suppose that
$\{\phi_j(x)\}_{j=0}^{\infty}$ is an arbitrary complete ortho\-nor\-mal system of 
functions in the space $L_2([t,T])$ and $\psi_1(\tau),$ $\ldots, \psi_5(\tau)$
are continuous functions on $[t, T].$
Furthermore, let the conditions
{\rm (\ref{novemberxxx3})--(\ref{novemberxxx6})} are satisfied.
Then$,$ for the iterated Stra\-to\-no\-vich stochastic integral
of fifth multiplicity 
$$
{\int\limits_t^{*}}^T \psi_5(t_5)
\ldots 
{\int\limits_t^{*}}^{t_2}\psi_1(t_1)
d{\bf w}_{t_1}^{(i_1)}
\ldots d{\bf w}_{t_5}^{(i_5)}\ \ \ (i_1,\ldots,i_5=0,1,\ldots,m)
$$

\noindent
the following 
expansion 
$$
{\int\limits_t^{*}}^T \psi_5(t_5)
\ldots 
{\int\limits_t^{*}}^{t_2}\psi_1(t_1)
d{\bf w}_{t_1}^{(i_1)}
\ldots d{\bf w}_{t_5}^{(i_5)}=
\hbox{\vtop{\offinterlineskip\halign{
\hfil#\hfil\cr
{\rm l.i.m.}\cr
$\stackrel{}{{}_{p\to \infty}}$\cr
}} }\sum_{j_1,\ldots,j_5=0}^{p}
C_{j_5 \ldots j_1}\zeta_{j_1}^{(i_1)}\ldots \zeta_{j_5}^{(i_5)}
$$

\noindent
that converges in the mean-square sense is valid, where 
notations are the same as in Theorem~{\rm 2.56.}}

Note that Theorems~2.55 and 2.57 are simple consequences of Theorems~2.54 and 2.56, respectively
(see Theorem~2.12 $(k=4,\ 5)).$

\section{On the Calculation of Matrix Traces of Volterra--Type Integral Operators}

\subsection{Introduction}

It is easy to see that the function (\ref{chain200}) for even $k=2r$ $(r\in{\bf N})$
forms a family of integral operators
$\mathbb{K}: L_2([t, T]^r) \rightarrow L_2([t, T]^r)$
(with the kernel (\ref{chain200}))
of the form

\vspace{-3mm}
\begin{equation}
\label{july6999}
~~~~\left(\mathbb{K} f\right)(t_{g_1},\ldots,t_{g_r})=
\int\limits_{[t, T]^{r}}K(t_1,\ldots,t_k)f(t_{g_{r+1}}, \ldots, t_{g_k})
dt_{g_{r+1}} \ldots dt_{g_k},
\end{equation}

\vspace{2mm}
\noindent
where $\{g_1,\ldots,g_k\}=\{1,\ldots,k\},$
the kernel $K(t_1,\ldots,t_k)$ is defined by (\ref{chain200}), i.e. has
the form

\vspace{-3mm}
\begin{equation}
\label{july7000}
~~~~~~~~~K(t_1,\ldots,t_k)=
\left\{\begin{matrix}
\psi_1(t_1)\ldots \psi_k(t_k),\ &t_1<\ldots<t_k\cr\cr
0,\ &\hbox{\rm otherwise}
\end{matrix}
\right.,
\end{equation}

\vspace{2mm}
\noindent
where $\psi_1(\tau),\ldots,\psi_k(\tau)\in L_2([t,T])$,
$t_1,\ldots,t_k\in [t, T]$ $(k\ge 2)$ and 
$K(t_1)\equiv\psi_1(t_1)$ for $t_1\in[t, T].$

For example,
\begin{equation}
\label{july7013}
~~~~~~~~~~~\left(\mathbb{K} f\right)(t_2)=
\int\limits_t^T K(t_1,t_2)f(t_1)dt_1
=\psi_2(t_2)\int\limits_t^{t_2}\psi_1(t_1)f(t_1)dt_1,
\end{equation}

\vspace{1mm}
$$
\left(\mathbb{K} f\right)(t_3,t_4)=
\int\limits_{[t, T]^{2}}K(t_1,\ldots,t_4)f(t_1, t_2)
dt_1 dt_2=
$$
$$
=
{\bf 1}_{\{t_3<t_4\}}\psi_3(t_3)\psi_4(t_4)
\int\limits_t^{t_3} \psi_2(t_2)\int\limits_t^{t_2}\psi_1(t_1)
f(t_1, t_2)
dt_1 dt_2,
$$
$$
\left(\mathbb{K} f\right)(t_1,t_2)=
\int\limits_{[t, T]^{2}}K(t_1,\ldots,t_4)f(t_3, t_4)
dt_3 dt_4=
$$
$$
=
\psi_1(t_1)\psi_2(t_2){\bf 1}_{\{t_1<t_2\}}
\int\limits_{t_2}^{T} \psi_3(t_3)\int\limits_{t_3}^{T}\psi_4(t_4)
f(t_3, t_4)
dt_4 dt_3.
$$

\vspace{2mm}

The simplest representative of the family (\ref{july6999})
has the form
\begin{equation}
\label{july7003}
\left(\mathbb{V} f\right)(x)=\int\limits_0^x f(\tau)d\tau
\end{equation}
and is called the Volterra integral operator, where $\mathbb{V}: L_2([0,1]) \rightarrow L_2([0,1]),$
$f(\tau)\in L_2([0,1]).$
The kernel of the Volterra integral operator has the following form
$$
K(\tau,x)=
\left\{\begin{matrix}
1,\ &\tau < x\cr\cr
0,\ &\hbox{\rm otherwise}
\end{matrix}
\right.,\ \ \ \tau, x\in [0, 1].
$$

Suppose that $\mathbb{A}: H \rightarrow H$ is a linear bounded operator. 
Recall \cite{gohb} that $\mathbb{A}$ has a finite matrix trace
if for any orthonormal basis $\left\{\Psi_j(x)\right\}_{j=0}^{\infty}$
of the space $H$ the series
\begin{equation}
\label{july7002}
~~~~\sum_{j=0}^{\infty} \left\langle
\mathbb{A}\Psi_j, \Psi_j\right\rangle_H
\end{equation}

\noindent
converges, where $\left\langle
\cdot , \cdot \right\rangle_H$ is a scalar probuct in $H$.

Note that the series (\ref{july7002}) converges absolutely
since its sum does not depend on the permutation of the terms
of the series (\ref{july7002})
(any permutation of basis functions $\Psi_j(x)$ forms a basis 
in $H)$ \cite{gohb}.

It is well known that the Volterra integral operator 
(\ref{july7003}) is not a trace class operator since
its singular values are equal to 
$s_j(\mathbb{V})=2\left(\pi(2j+1)\right)^{-1}$ \cite{Brisl}.
On the other hand, it is known \cite{Brisl} that for trace class
operators the equality of matrix and integral traces holds.
It turns out that for the Volterra integral operator 
(\ref{july7003}) (although it is not a trace class operator),
the equality of matrix and integral traces is also true
\cite{Brisl}.
Thus, one cannot count on the fact that operators of the more
general form (\ref{july6999}) (from the same 
family of operators as the Volterra integral operator (\ref{july7003})) 
are operators of the trace class.
Nevertheless, the proof of the equalities 
of matrix and integral traces 
for Volterra--type integral operators (\ref{july6999}) (which is obviously a problem) 
provides a way
to calculate the matrix traces of these operators.

Why do we talk so much in this section about matrix
traces of operators from the family (\ref{july6999})?
The point is that matrix traces of operators of the form
(\ref{july6999}) are of great importance for obtaining
of expansions
of iterated Stratonovich stochastic integrals.

Throughout the Chapter~2, we have already 
calculated the matrix traces mentioned above many times
(see the formulas (\ref{5t}), (\ref{may62021}), (\ref{start1000}),
(\ref{otit239}),  (\ref{otit990}), (\ref{cas1000}), 
(\ref{july7010}), (\ref{july7011}), 
(\ref{5tzzz}), 
(\ref{cas9}), (\ref{cas10}), (\ref{strange508}), (\ref{after80}),
(\ref{after500xx}), (\ref{sixsix8})--(\ref{sixsix15}),
(\ref{sixsix129}), 
(\ref{2023novem206})--(\ref{2023novem208}), 
(\ref{march10})--(\ref{march12}), 
(\ref{july13})--(\ref{july18})).

Let us consider some illustrative examples.
We have
\begin{equation}
\label{july7015}
\sum_{j_1=0}^{\infty} \left\langle
\mathbb{K}\phi_{j_1}, \phi_{j_1}\right\rangle_{L_2([t,T])}=
\end{equation}
\begin{equation}
\label{july7020}
~~~~~~~~~~~=
\sum_{j_1=0}^{\infty}
\int\limits_t^T
\psi_2(t_2)\phi_{j_1}(t_2)\int\limits_t^{t_2}
\psi_1(t_1)\phi_{j_1}(t_1)dt_1 dt_2=
\sum_{j_1=0}^{\infty}C_{j_1 j_1},
\end{equation}

\begin{equation}
\label{july7016}
\sum_{j_1,j_2=0}^{\infty} \left\langle
\mathbb{K}\Psi_{j_1 j_2}, \Psi_{j_1 j_2}\right\rangle_{L_2([t,T]^2)}=
\end{equation}
$$
=
\sum_{j_1,j_2=0}^{\infty}
\int\limits_t^T
\psi_4(t_4)\phi_{j_2}(t_4)\int\limits_t^{t_4}
\psi_3(t_3)\phi_{j_2}(t_3)\int\limits_t^{t_3}
\psi_2(t_2)\phi_{j_1}(t_2)
\int\limits_t^{t_2}
\psi_1(t_1)\phi_{j_1}(t_1)\times
$$
\begin{equation}
\label{july7021}
\times
dt_1 dt_2 dt_3 dt_4=\sum_{j_1,j_2=0}^{\infty}
C_{j_2 j_2 j_1 j_1},
\end{equation}

\vspace{2mm}
\noindent
where $\left\{\Psi_{j_1 j_2}(x,y)\right\}_{j_1,j_2=0}^{\infty}=
\left\{\phi_{j_1}(x)\phi_{j_2}(y)\right\}_{j_1,j_2=0}^{\infty},$
$\left\{\phi_{j}(x)\right\}_{j=0}^{\infty}$
is an arbitrary complete orthonormal system of functions
in $L_2([t, T]),$ $\left(\mathbb{K}f\right)(t_2)$ in (\ref{july7015})
is defined by (\ref{july7013}),
and $\left(\mathbb{K}f\right)(t_2, t_3)$ in (\ref{july7016}) has the following form
$$
\left(\mathbb{K} f\right)(t_2,t_3)=
\int\limits_{[t, T]^{2}}K(t_1,\ldots,t_4)f(t_1, t_4)
dt_1 dt_4=
$$
$$
=\psi_2(t_2)\psi_3(t_3){\bf 1}_{\{t_2<t_3\}}
\int\limits_t^{t_2} \psi_1(t_1)\int\limits_{t_3}^{T}\psi_4(t_4)
f(t_1, t_4)
dt_4 dt_1,
$$

\noindent
where $K(t_1,\ldots,t_4)$ is defined by (\ref{july7000}).

The expressions on the right-hand sides of (\ref{july7020}) and (\ref{july7021})
were considered earlier in Chapter~2 under various assumptions
on $\left\{\phi_j(x)\right\}_{j=0}^{\infty}$ and $\psi_1(\tau),\ldots,\psi_4(\tau)$
(see the formulas (\ref{5t}), (\ref{may62021}), (\ref{start1000}),
(\ref{otit239}), (\ref{july7010}), (\ref{5tzzz}), 
(\ref{cas9}), (\ref{2023novem206}), (\ref{march10}), 
(\ref{july13})).

\subsection{Approach Based on Generalized Parseval's Equality and (\ref{start1000}).
Symmetrical Case When $\psi_1(\tau)=\psi_k(\tau),$ $\psi_2(\tau)=\psi_{k-1}(\tau),\ldots
,\psi_r(\tau)=\psi_{r+1}(\tau)$ $(k=2r, r=2,3,4,\ldots)$ and 
$\psi_1(\tau),$ $\ldots,\psi_k(\tau)\in L_2([t, T])$}

Let us consider one of the possible ways to calculate
matix traces of Volterra-type integral operators 
(\ref{july6999}) based Fubini's Theorem, Parseval's equality
and generalized Parseval's equality.

Recall the equalities (\ref{sixsix40}) and (\ref{2023novem215})

\vspace{-3mm}
$$
C_{j_6 j_5 j_4 j_3 j_2 j_1}+C_{j_1 j_2 j_3 j_4 j_5 j_6}=
C_{j_6}C_{j_5 j_4 j_3 j_2 j_1}-C_{j_5 j_6}C_{j_4 j_3 j_2 j_1}+
$$
\begin{equation}
\label{july7100}
+C_{j_4 j_5 j_6}C_{j_3 j_2 j_1}-C_{j_3 j_4 j_5 j_6}C_{j_2 j_1}+
C_{j_2 j_3 j_4 j_5 j_6}C_{j_1},
\end{equation}

\vspace{-4mm}
\begin{equation}
\label{july7101}
~~~~~~~~~~C_{j_4 j_3 j_2 j_1}+C_{j_1 j_2 j_3 j_4}=
C_{j_4}C_{j_3 j_2 j_1}-C_{j_3 j_4}C_{j_2 j_1}+
C_{j_2 j_3 j_4}C_{j_1},
\end{equation}

\vspace{3mm}
\noindent
where $C_{j_k\ldots j_1}$ is defined by
the formula
$$
C_{j_k \ldots j_1}=\int\limits_t^T\psi_k(t_k)\phi_{j_k}(t_k)\ldots
\int\limits_t^{t_2}
\psi_1(t_1)\phi_{j_1}(t_1)
dt_1\ldots dt_k\ \ \ (k\in{\bf N})
$$
for the case $\psi_1(\tau),\ldots,\psi_k(\tau)\equiv 1$.

It is easy to see (see the derivation of
(\ref{sixsix40}) and (\ref{2023novem215})) that analogues of the relations
(\ref{july7100}), (\ref{july7101}) (with appropriate changes) hold for 
$\psi_1(\tau),\ldots,\psi_6(\tau)\in L_2([t, T]).$

By analogy with (\ref{july7100}), (\ref{july7101})
(see the derivation of
(\ref{sixsix40}) and (\ref{2023novem215})) we obtain for $k=2r$ $(r=2,3,4,\ldots)$

\vspace{-5mm}
$$
C_{j_k j_{k-1}\ldots j_1}^{\psi_k \psi_{k-1} \ldots \psi_1} +
C_{j_1 j_2\ldots j_k}^{\psi_1 \psi_{2} \ldots \psi_k}=
C_{j_k}^{\psi_k} \cdot  C_{j_{k-1} j_{k-2}\ldots j_1}^{\psi_{k-1} \psi_{k-2} \ldots \psi_1}
-C_{j_{k-1} j_k}^{\psi_{k-1} \psi_{k}} \cdot C_{j_{k-2} j_{k-3}\ldots j_1}
^{\psi_{k-2} \psi_{k-3} \ldots \psi_1}+
$$
\begin{equation}
\label{july7028}
~~~~~~+C_{j_{k-2} j_{k-1} j_k}^{\psi_{k-2} \psi_{k-1} \psi_k} \cdot
C_{j_{k-3} j_{k-4}\ldots j_1}^{\psi_{k-3} \psi_{k-4} \ldots \psi_1}
- \ldots - C_{j_3 j_4 \ldots j_k}^{\psi_3 \psi_{4} \ldots \psi_k} \cdot C_{j_2 j_1}^{\psi_2 \psi_1}+
C_{j_2 j_3 \ldots j_k}^{\psi_2 \psi_{3} \ldots \psi_k} \cdot C_{j_1}^{\psi_1},
\end{equation}

\vspace{3mm}
\noindent
where 
\begin{equation}
\label{july40000}
~~~~~~~~C_{j_k \ldots j_1}^{\psi_k\ldots \psi_1}=\int\limits_t^T\psi_k(t_k)\phi_{j_k}(t_k)\ldots
\int\limits_t^{t_2}
\psi_1(t_1)\phi_{j_1}(t_1)
dt_1\ldots dt_k\ \ \ (k\in{\bf N}).
\end{equation}

When proving Theorems~2.46 and 2.48, 
using (\ref{july7028}) (the case $k=4,$ $\psi_1(\tau),\ldots,\psi_4(\tau)\equiv 1$),
we obtained the following formulas
$$
\lim\limits_{p\to\infty}
\sum\limits_{j_1, j_3=0}^{p}
C_{j_3 j_3 j_1 j_1}=\frac{1}{8}(T-t)^2,
$$
$$
\lim\limits_{p\to\infty}
\sum\limits_{j_1, j_3=0}^{p}
C_{j_1 j_3 j_3 j_1}=0,
$$
$$
\lim\limits_{p\to\infty}
\sum\limits_{j_1, j_2=0}^{p}
C_{j_2 j_1 j_2 j_1}=0,
$$

\vspace{2mm}
\noindent
where
$\{\phi_j(x)\}_{j=0}^{\infty}$ is an arbitrary complete orthonormal system of 
functions in the space $L_2([t,T])$ and we use the notation
$C_{j_k \ldots j_1}$ instead of 
$C_{j_k \ldots j_1}^{\psi_k \ldots \psi_1}$ for the case $\psi_1(\tau),\ldots,\psi_k(\tau)\equiv 1.$

In principle, using (\ref{july7028}), we can 
calculate some matrix traces for which the following
symmetry condition 
\begin{equation}
\label{july7030}
\psi_1(\tau)=\psi_k(\tau),\ \ \psi_2(\tau)=\psi_{k-1}(\tau),\ \ldots
,\ \ \psi_r(\tau)=\psi_{r+1}(\tau)\ \ (k=2r, r=2,3,4,\ldots)
\end{equation}

\noindent
is satisfied. Obviously, the case
$\psi_1(\tau),\ldots,\psi_k(\tau)\equiv 1$ 
is possible since it is a special case of (\ref{july7030}).
This case is important because it covers
the case of iterated Stratonovich
stochastic integrals from the classical Taylor--Stratonovich
expansions (see Chapter~4).

Consider the case $k=4$ of (\ref{july7028})

\vspace{-5mm}
\begin{equation}
\label{july5000}
~~~~~~~C_{j_4 j_3 j_2 j_1}^{\psi_4 \psi_3 \psi_2 \psi_1} + C_{j_1 j_2 j_3 j_4}^{\psi_1 \psi_2 \psi_3 \psi_4}=
C_{j_4}^{\psi_4}C_{j_3 j_2 j_1}^{\psi_3 \psi_2 \psi_1}-C_{j_3 j_4}^{\psi_3 \psi_4}
C_{j_2 j_1}^{\psi_2 \psi_1}+
C_{j_2 j_3 j_4}^{\psi_2 \psi_3 \psi_4}C_{j_1}^{\psi_1},
\end{equation}

\vspace{2mm}
\noindent
where $\psi_1(\tau), \ldots, \psi_4(\tau)\in L_2([t,T]).$

Substitute $j_4=j_3,$ $j_2=j_1$ into (\ref{july5000})

\vspace{-5mm}
\begin{equation}
\label{july5001}
~~~~~~~C_{j_3 j_3 j_1 j_1}^{\psi_4 \psi_3 \psi_2 \psi_1} + C_{j_1 j_1 j_3 j_3}^{\psi_1 \psi_2 \psi_3 \psi_4}=
C_{j_3}^{\psi_4}C_{j_3 j_1 j_1}^{\psi_3 \psi_2 \psi_1}-C_{j_3 j_3}^{\psi_3 \psi_4}
C_{j_1 j_1}^{\psi_2 \psi_1}+
C_{j_1 j_3 j_3}^{\psi_2 \psi_3 \psi_4}C_{j_1}^{\psi_1}.
\end{equation}

\vspace{2mm}

Applying (\ref{july5001}), we get
$$
\lim\limits_{p\to\infty}\sum\limits_{j_1,j_3=0}^{p}
\left(C_{j_3 j_3 j_1 j_1}^{\psi_4 \psi_3 \psi_2 \psi_1} +
C_{j_1 j_1 j_3 j_3}^{\psi_1 \psi_2 \psi_3 \psi_4}\right)=
\lim\limits_{p\to\infty}\sum\limits_{j_1,j_3=0}^{p}
C_{j_3}^{\psi_4}C_{j_3 j_1 j_1}^{\psi_3 \psi_2 \psi_1}-
$$
\begin{equation}
\label{july5002}
~~~~~~~~-
\lim\limits_{p\to\infty}\sum\limits_{j_1,j_3=0}^{p}
C_{j_3 j_3}^{\psi_3 \psi_4}
C_{j_1 j_1}^{\psi_2 \psi_1}+
\lim\limits_{p\to\infty}\sum\limits_{j_1,j_3=0}^{p}
C_{j_1 j_3 j_3}^{\psi_2 \psi_3 \psi_4}C_{j_1}^{\psi_1}.
\end{equation}

\vspace{2mm}

From (\ref{after1400}) we have
$$
\lim\limits_{p\to\infty}\sum\limits_{j_3=0}^{p} C_{j_3 j_3}^{\psi_3 \psi_4}
\sum\limits_{j_1=0}^{p}C_{j_1 j_1}^{\psi_2 \psi_1}=
\lim\limits_{p\to\infty}\sum\limits_{j_3=0}^{p} C_{j_3 j_3}^{\psi_3 \psi_4}
\lim\limits_{p\to\infty}\sum\limits_{j_1=0}^{p} C_{j_1 j_1}^{\psi_2 \psi_1}=
$$
\begin{equation}
\label{july5003}
=
\frac{1}{4}\int\limits_t^T \psi_4(s) \psi_3(s)ds
\int\limits_t^T \psi_2(s) \psi_1(s)ds.
\end{equation}

Further, we obtain
$$
\lim\limits_{p\to\infty}\sum\limits_{j_3=0}^{p} C_{j_3}^{\psi_4} 
\sum\limits_{j_1=0}^{p} C_{j_3 j_1 j_1}^{\psi_3 \psi_2 \psi_1}=
\frac{1}{2}\lim\limits_{p\to\infty}
\sum\limits_{j_3=0}^{p}C_{j_3}^{\psi_4} 
C_{j_3 j_1 j_1}^{\psi_3 \psi_2 \psi_1}\biggl|_{(j_{1} j_{1})\curvearrowright (\cdot)}-
$$
\begin{equation}
\label{july5004}
~~~~~~~~-
\lim\limits_{p\to\infty}
\sum\limits_{j_3=0}^{p}C_{j_3}^{\psi_4}
\left(\frac{1}{2}
C_{j_3 j_1 j_1}^{\psi_3 \psi_2 \psi_1}\biggl|_{(j_{1} j_{1})\curvearrowright (\cdot)}-
\sum\limits_{j_1=0}^{p} C_{j_3 j_1 j_1}^{\psi_3 \psi_2 \psi_1}\right).
\end{equation}

\vspace{2mm}

Applying the generalized Parseval equality, we have
$$
\lim\limits_{p\to\infty}
\sum\limits_{j_3=0}^{p}C_{j_3}^{\psi_4}
C_{j_3 j_1 j_1}^{\psi_3 \psi_2 \psi_1}\biggl|_{(j_{1} j_{1})\curvearrowright (\cdot)}=
$$
$$
=
\lim\limits_{p\to\infty}
\sum\limits_{j_3=0}^{p}\int\limits_t^T \psi_4(s)\phi_{j_3}(s)ds
\int\limits_t^T \psi_3(s)\phi_{j_3}(s)\int\limits_t^{s}
\psi_2(\tau)\psi_1(\tau)d\tau ds
=
$$
\begin{equation}
\label{july5005}
=
\int\limits_t^T \psi_4(s) \psi_3(s)
\int\limits_t^{s} \psi_2(\tau) \psi_1(\tau)d\tau ds.
\end{equation}

From (\ref{july5004}) and (\ref{july5005}) we obtain
$$
\lim\limits_{p\to\infty}\sum\limits_{j_3=0}^{p} C_{j_3}^{\psi_4} 
\sum\limits_{j_1=0}^{p} C_{j_3 j_1 j_1}^{\psi_3 \psi_2 \psi_1}=
\frac{1}{2}\int\limits_t^T \psi_4(s) \psi_3(s)
\int\limits_t^{s} \psi_2(\tau) \psi_1(\tau)d\tau ds
-
$$
\begin{equation}
\label{july5006}
~~~~~~~~~-\lim\limits_{p\to\infty}
\sum\limits_{j_3=0}^{p}C_{j_3}^{\psi_4}
\left(\frac{1}{2}
C_{j_3 j_1 j_1}^{\psi_3 \psi_2 \psi_1}\biggl|_{(j_{1} j_{1})\curvearrowright (\cdot)}-
\sum\limits_{j_1=0}^{p} C_{j_3 j_1 j_1}^{\psi_3 \psi_2 \psi_1}\right).
\end{equation}

\vspace{2mm}

Due to Cauchy--Bunyakovsky's inequality, Parseval's equality
and (\ref{july1003}), we get 
$$
\lim\limits_{p\to\infty}
\left(\sum\limits_{j_3=0}^{p}C_{j_3}^{\psi_4}
\left(\frac{1}{2}
C_{j_3 j_1 j_1}^{\psi_3 \psi_2 \psi_1}\biggl|_{(j_{1} j_{1})\curvearrowright (\cdot)}-
\sum\limits_{j_1=0}^{p} C_{j_3 j_1 j_1}^{\psi_3 \psi_2 \psi_1}\right)\right)^2\le
$$
$$
\le \lim\limits_{p\to\infty}
\sum\limits_{j_3=0}^{p}\left(C_{j_3}^{\psi_4}\right)^2\
\sum\limits_{j_3=0}^{p}
\left(\frac{1}{2}
C_{j_3 j_1 j_1}^{\psi_3 \psi_2 \psi_1}\biggl|_{(j_{1} j_{1})\curvearrowright (\cdot)}-
\sum\limits_{j_1=0}^{p} C_{j_3 j_1 j_1}^{\psi_3 \psi_2 \psi_1}\right)^2\le
$$
$$
\le \lim\limits_{p\to\infty}
\sum\limits_{j_3=0}^{\infty}\left(C_{j_3}^{\psi_4}\right)^2\
\sum\limits_{j_3=0}^{p}
\left(\frac{1}{2}
C_{j_3 j_1 j_1}^{\psi_3 \psi_2 \psi_1}\biggl|_{(j_{1} j_{1})\curvearrowright (\cdot)}-
\sum\limits_{j_1=0}^{p} C_{j_3 j_1 j_1}^{\psi_3 \psi_2 \psi_1}\right)^2=
$$
\begin{equation}
\label{july5007}
~~~~~~~~=\int\limits_t^T \psi_4^2(s)ds\lim\limits_{p\to\infty}
\sum\limits_{j_3=0}^{p}
\left(\frac{1}{2}
C_{j_3 j_1 j_1}^{\psi_3 \psi_2 \psi_1}\biggl|_{(j_{1} j_{1})\curvearrowright (\cdot)}-
\sum\limits_{j_1=0}^{p} C_{j_3 j_1 j_1}^{\psi_3 \psi_2 \psi_1}\right)^2=0.
\end{equation}

\vspace{2mm}         

Combining (\ref{july5006}) and (\ref{july5007}), we obtain
\begin{equation}
\label{july5008}
~~~~~~\lim\limits_{p\to\infty}\sum\limits_{j_3=0}^p C_{j_3}^{\psi_4} 
\sum\limits_{j_1=0}^p C_{j_3 j_1 j_1}^{\psi_3 \psi_2 \psi_1}=
\frac{1}{2}\int\limits_t^T \psi_4(s) \psi_3(s)
\int\limits_t^{s} \psi_2(\tau) \psi_1(\tau)d\tau ds.
\end{equation}

Absolutely similarly to (\ref{july5008}) we get
\begin{equation}
\label{july5009}
~~~~\lim\limits_{p\to\infty}\sum\limits_{j_1=0}^p C_{j_1}^{\psi_1}
\sum\limits_{j_3=0}^p C_{j_1 j_3 j_3}^{\psi_2 \psi_3 \psi_4}=
\frac{1}{2}\int\limits_t^T \psi_2(s) \psi_1(s)
\int\limits_t^{s} \psi_3(\tau) \psi_4(\tau)d\tau ds.
\end{equation}

Combining (\ref{july5002}), (\ref{july5003}), (\ref{july5008}), (\ref{july5009}) and 
applying Fubini's Theorem, we have
$$
\lim\limits_{p\to\infty}\sum\limits_{j_1,j_3=0}^p
\left(C_{j_3 j_3 j_1 j_1}^{\psi_4 \psi_3 \psi_2 \psi_1} +
C_{j_1 j_1 j_3 j_3}^{\psi_1 \psi_2 \psi_3 \psi_4}\right)=
\frac{1}{2}\int\limits_t^T \psi_4(s) \psi_3(s)
\int\limits_t^{s} \psi_2(\tau) \psi_1(\tau)d\tau ds+
$$
$$
+
\frac{1}{2}\int\limits_t^T \psi_2(s) \psi_1(s)
\int\limits_t^{s} \psi_3(\tau) \psi_4(\tau)d\tau ds
-\frac{1}{4}\int\limits_t^T \psi_4(s) \psi_3(s)ds
\int\limits_t^T \psi_2(s) \psi_1(s)ds=
$$
$$
=\frac{1}{4}\int\limits_t^T \psi_4(s) \psi_3(s)ds
\int\limits_t^T \psi_2(s) \psi_1(s)ds=
$$
\begin{equation}
\label{july5010}
=\frac{1}{4}\int\limits_t^T \psi_4(s) \psi_3(s)
\int\limits_t^{s} \psi_2(\tau) \psi_1(\tau)d\tau ds+
\frac{1}{4}\int\limits_t^T \psi_2(s) \psi_1(s)
\int\limits_t^{s} \psi_3(\tau) \psi_4(\tau)d\tau ds.
\end{equation}

Let us rewrite (\ref{july5010}) in the form
$$
\lim\limits_{p\to\infty}\sum\limits_{j_1,j_3=0}^p
\left(C_{j_3 j_3 j_1 j_1}^{\psi_4 \psi_3 \psi_2 \psi_1} +
C_{j_3 j_3 j_1 j_1}^{\psi_1 \psi_2 \psi_3 \psi_4}\right)=
$$
\begin{equation}
\label{july5010x}
=\frac{1}{4}\int\limits_t^T \psi_4(s) \psi_3(s)
\int\limits_t^{s} \psi_2(\tau) \psi_1(\tau)d\tau ds+
\frac{1}{4}\int\limits_t^T \psi_2(s) \psi_1(s)
\int\limits_t^{s} \psi_3(\tau) \psi_4(\tau)d\tau ds.
\end{equation}

It is easy to see that the left-hand side 
of (\ref{july5010x}) does not depend on 
the simultaneous rearrangement of $\psi_4$ with $\psi_1$
and $\psi_3$ with $\psi_2$.

Using the above arguments and using derivation method of (\ref{march11}) and (\ref{march12}), we get
\begin{equation}
\label{july8000}
\lim\limits_{p\to\infty}\sum\limits_{j_1,j_3=0}^p
\left(C_{j_3 j_1 j_3 j_1}^{\psi_4 \psi_3 \psi_2 \psi_1} +
C_{j_3 j_1 j_3 j_1}^{\psi_1 \psi_2 \psi_3 \psi_4}\right)=0,
\end{equation}
\begin{equation}
\label{july8001}
\lim\limits_{p\to\infty}\sum\limits_{j_1,j_3=0}^p
\left(C_{j_1 j_3 j_3 j_1}^{\psi_4 \psi_3 \psi_2 \psi_1} +
C_{j_1 j_3 j_3 j_1}^{\psi_1 \psi_2 \psi_3 \psi_4}\right)=0.
\end{equation}

\vspace{2mm}

Using (\ref{july5010x})--(\ref{july8001}) under the conditions
$\psi_1(\tau)=\psi_4(\tau),$ $\psi_2(\tau)=\psi_3(\tau),$ we obtain
$$
\lim\limits_{p\to\infty}\sum\limits_{j_1,j_3=0}^p
C_{j_3 j_3 j_1 j_1}^{\psi_1 \psi_2 \psi_2 \psi_1} 
=\frac{1}{4}\int\limits_t^T \psi_2(s) \psi_1(s)
\int\limits_t^{s} \psi_2(\tau) \psi_1(\tau)d\tau ds,
$$
$$
\lim\limits_{p\to\infty}\sum\limits_{j_1,j_3=0}^p
C_{j_3 j_1 j_3 j_1}^{\psi_1 \psi_2 \psi_2 \psi_1}=0,
$$
$$
\lim\limits_{p\to\infty}\sum\limits_{j_1,j_3=0}^p
C_{j_1 j_3 j_3 j_1}^{\psi_1 \psi_2 \psi_2 \psi_1}=0.
$$

\subsection{Approach Based on Trace Class Operators}

An efficient method for calculating of matrix traces 
of Volterra--type integral operators of the form (\ref{july6999})
was proposed in \cite{rybakov7000x}.
This method is based on Theorem~3.1 from \cite{Brisl}.
Theorem~3.1 \cite{Brisl} implies the following statement.

\vspace{2mm}

{\bf Theorem~D} (see \cite{Brisl} for details).\ {\it
Let $\mathbb{K}: L_2([t,T]^r)\rightarrow L_2([t,T]^r)$ $(r\in {\bf N})$
be a trace class operator with the kernel $K\in L_2([t,T]^{2r})$. Then
$\tilde K(t_1,\ldots,t_r,t_1,\ldots,t_r)$ exists
almost everywhere $[dt_1\ldots dt_r]$ and
\begin{equation}
\label{july15000}
tr\mathbb{K}=\int\limits_{[t,T]^r}\tilde K(t_1,\ldots,t_r,t_1,\ldots,t_r)
dt_1\ldots dt_r,
\end{equation}
where 
$$
\tilde F(x_1,\ldots,x_m)\stackrel{\sf def}{=}\lim\limits_{u\to 0}
A_u F(x_1,\ldots,x_m),
$$

\vspace{-3mm}
$$
A_u F(x_1,\ldots,x_m)\stackrel{\sf def}{=}
\frac{1}{\left(2u\right)^{m}}
\int\limits_{[-u,u]^{m}}
F(x_1+\tau_1,\ldots,x_m+\tau_m)
d\tau_1\ldots d\tau_m\ \ \ (m\in{\bf N}).
$$ 
}

Let us prove the equality (\ref{july13})
using the method from \cite{rybakov7000x} in our interpretation.
Consider two symmetric functions of the form (\ref{ziko5001}) which
we introduced in Sect.~2.1.2
\begin{equation}
\label{july15004}
~~~~~~~~~~~~K'(t_1,t_2)=\psi_1(t_1)f_2(t_2){\bf 1}_{\{t_1\le t_2\}}+
\psi_1(t_2)f_2(t_1){\bf 1}_{\{t_1\ge t_2\}},
\end{equation}
\begin{equation}
\label{july15005}
~~~~~~~~~~~~K''(t_3,t_4)=f_3(t_3)\psi_4(t_4){\bf 1}_{\{t_3\le t_4\}}+
f_3(t_4)\psi_4(t_3){\bf 1}_{\{t_3\ge t_4\}},
\end{equation}

\noindent
where we suppose that $\psi_1(\tau), \psi_4(\tau)$ are continuously differentiable
functions on $[t, T]$ (the case 
$\psi_1(\tau), \psi_4(\tau)\in L_2([t,T])$ will be considered further)
and $f_2(\tau), f_3(\tau)$  are polynomials of finite degrees.

By Theorem~B (see Sect.~2.1.5) and (\ref{july11001}) we have that the 
kernels $K'(t_1,t_2)$ and $K''(t_3,t_4)$ (see (\ref{july15004}), (\ref{july15005}))
correspond to the 
trace class integral operators.

It is known \cite{Brisl} that the integral operator $\mathbb{A}$ is a trace class operator
if and only if the kernel $K(x,y)$ of $\mathbb{A}$ has the following
representation
\begin{equation}
\label{july15007}
K(x,y)=\int\limits_{[t,T]^{2n}} K_1(x,\tau)K_2(\tau,y)d\tau
\end{equation}
almost everywhere $[dxdy]$,
where $K_1(x,y), K_2(x,y)$ are kernels of Hilbert--Schmidt operators,
$x, y\in {\bf R}^n$ $(n\ge 1).$

Since $K'(t_1,t_2)$ and $K''(t_3,t_4)$ are kernels of
the trace class integral operators, then (see (\ref{july15007}))
\begin{equation}
\label{july15008}
K'(t_1,t_2)=\int\limits_{[t,T]} K_1'(t_1,\tau)K_2'(\tau,t_2)d\tau,\ \ \
K''(t_1,t_2)=\int\limits_{[t,T]} K_1''(t_1,\tau)K_2''(\tau,t_2)d\tau
\end{equation}
almost everywhere $[dt_1 dt_2],$ where $K_1', K_2', K_1'', K_2''\in L_2([t, T]^2)$.
Then, we have
$$
K'(t_1,t_2)K''(t_3,t_4)
=\int\limits_{[t,T]} K_1'(t_1,\tau_1)K_2'(\tau_1,t_2)d\tau_1
\int\limits_{[t,T]} K_1''(t_3,\tau_2)K_2''(\tau_2,t_4)d\tau_2=
$$
\begin{equation}
\label{july15010}
~~~~~~~~~~=\int\limits_{[t,T]^2} K_1'(t_1,\tau_1)K_1''(t_3,\tau_2)
K_2'(\tau_1,t_2)
K_2''(\tau_2,t_4) d\tau_1 d\tau_2.
\end{equation}

The equality (\ref{july15010}) can be written as follows
$$
F(t_1,t_3,t_2,t_4)=\int\limits_{[t,T]^2} F_1(t_1,t_3,\tau_1,\tau_2)
F_2(\tau_1,\tau_2,t_2,t_4) d\tau_1 d\tau_2
$$
almost everywhere $[dt_1 dt_2 dt_3 dt_4]$, where
$F(t_1,t_3,t_2,t_4)=K'(t_1,t_2)K''(t_3,t_4),$
$F_1(t_1,t_3, \tau_1,\tau_2)=K_1'(t_1,\tau_1)K_1''(t_3,\tau_2),$
and $F_2(\tau_1,\tau_2,t_2,t_4)$ $=K_2'(\tau_1,t_2)K_2''(\tau_2,t_4).$

As a result, the product $K'(t_1,t_2)K''(t_3,t_4)$
is also the kernel of the trace class operator (see (\ref{july15007})).
Let us denote it by $\mathbb{K'}.$

Suppose that $\left\{\phi_{j}(x)\right\}_{j=0}^{\infty}$
is an arbitrary complete orthonormal system of functions
in $L_2([t, T]).$ Then 
$\left\{\Psi_{j_1 j_2}(x,y)\right\}_{j_1,j_2=0}^{\infty}=
\left\{\phi_{j_1}(x)\phi_{j_2}(y)\right\}_{j_1,j_2=0}^{\infty}$
is an orthonormal basis in $L_2([t, T]^2).$

Consider matrix trace of $\mathbb{K'}.$ Using Fubini's Theorem, we obtain
$$
\sum_{j_1,j_2=0}^{\infty} \left\langle
\Psi_{j_1 j_2}, \mathbb{K'}\Psi_{j_1 j_2}\right\rangle_{L_2([t,T]^2)}=
$$
$$
=\sum_{j_1,j_2=0}^{\infty}\int\limits_{[t,T]^2}
\phi_{j_2}(t_4)\phi_{j_1}(t_1)
\int\limits_{[t,T]^2}K'(t_1,t_2)K''(t_3,t_4)
\phi_{j_2}(t_3)\phi_{j_1}(t_2)dt_2 dt_3 dt_1 dt_4=
$$

\vspace{-2mm}
$$
=
\sum_{j_1,j_2=0}^{\infty}\Biggl(
\int\limits_t^T
\psi_4(t_4)\phi_{j_2}(t_4)\int\limits_t^{T}
\psi_1(t_1)\phi_{j_1}(t_1)\int\limits_t^{t_4}
f_3(t_3)\phi_{j_2}(t_3)
\int\limits_{t_1}^T
f_2(t_2)\phi_{j_1}(t_2)\times\Biggr.
$$
$$
\times
dt_2 dt_3 dt_1 dt_4+
$$

\vspace{-3mm}
$$
+
\int\limits_t^T
f_3(t_4)\phi_{j_2}(t_4)\int\limits_t^{T}
\psi_1(t_1)\phi_{j_1}(t_1)\int\limits_{t_4}^T
\psi_4(t_3)\phi_{j_2}(t_3)
\int\limits_{t_1}^T
f_2(t_2)\phi_{j_1}(t_2)\times
$$
$$
\times
dt_2 dt_3 dt_1 dt_4+
$$

\vspace{-3mm}
$$
+
\int\limits_t^T
\psi_4(t_4)\phi_{j_2}(t_4)\int\limits_t^{T}
f_2(t_1)\phi_{j_1}(t_1)\int\limits_t^{t_4}
f_3(t_3)\phi_{j_2}(t_3)
\int\limits_t^{t_1}
\psi_1(t_2)\phi_{j_1}(t_2)\times
$$
$$
\times
dt_2 dt_3 dt_1 dt_4+
$$

\vspace{-3mm}
$$
+
\int\limits_t^T
f_2(t_1)\phi_{j_1}(t_1)\int\limits_t^{T}
\psi_4(t_3)\phi_{j_2}(t_3)\int\limits_t^{t_3}
f_3(t_4)\phi_{j_2}(t_4)
\int\limits_t^{t_1}
\psi_1(t_2)\phi_{j_1}(t_2)\times
$$
$$
\Biggl.\times
dt_2 dt_4 dt_3 dt_1\Biggr)=
$$

\vspace{-3mm}
$$
=
\sum_{j_1,j_2=0}^{\infty}\Biggl(
\int\limits_t^T
\psi_4(t_4)\phi_{j_2}(t_4)\int\limits_t^{t_4}
f_3(t_3)\phi_{j_2}(t_3)\int\limits_t^{T}
f_2(t_2)\phi_{j_1}(t_2)
\int\limits_t^{t_2}
\psi_1(t_1)\phi_{j_1}(t_1)\times\Biggr.
$$
$$
\times
dt_1 dt_2 dt_3 dt_4+
$$

\vspace{-3mm}
$$
+
\int\limits_t^T
\psi_4(t_3)\phi_{j_2}(t_3)\int\limits_t^{t_3}
f_3(t_4)\phi_{j_2}(t_4)\int\limits_t^{T}
f_2(t_2)\phi_{j_1}(t_2)
\int\limits_t^{t_2}
\psi_1(t_1)\phi_{j_1}(t_1)\times
$$
$$
\times
dt_1 dt_2 dt_4 dt_3+
$$

\vspace{-3mm}
$$
+
\int\limits_t^T
\psi_4(t_4)\phi_{j_2}(t_4)\int\limits_t^{t_4}
f_3(t_3)\phi_{j_2}(t_3)\int\limits_t^{T}
f_2(t_1)\phi_{j_1}(t_1)
\int\limits_t^{t_1}
\psi_1(t_2)\phi_{j_1}(t_2)\times
$$
$$
\times
dt_2 dt_1 dt_3 dt_4+
$$
$$
+
\int\limits_t^T
\psi_4(t_3)\phi_{j_2}(t_3)\int\limits_t^{t_3}
f_3(t_4)\phi_{j_2}(t_4)\int\limits_t^{T}
f_2(t_1)\phi_{j_1}(t_1)
\int\limits_t^{t_1}
\psi_1(t_2)\phi_{j_1}(t_2)\times
$$
$$
\times
\Biggl. dt_2 dt_1 dt_4 dt_3\Biggr)=
$$

\vspace{-3mm}
$$
=4
\sum_{j_1,j_2=0}^{\infty}
\int\limits_t^T
\psi_4(t_4)\phi_{j_2}(t_4)\int\limits_t^{t_4}
f_3(t_3)\phi_{j_2}(t_3)\int\limits_t^{T}
f_2(t_2)\phi_{j_1}(t_2)
\int\limits_t^{t_2}
\psi_1(t_1)\phi_{j_1}(t_1)\times
$$
\begin{equation}
\label{july15012}
\times
dt_1 dt_2 dt_3 dt_4.
\end{equation}

\vspace{4mm}

According to (\ref{july15012}) and (\ref{july15000}), we get
$$
\sum_{j_1,j_2=0}^{\infty} \left\langle
\Psi_{j_1 j_2}, \mathbb{K'}\Psi_{j_1 j_2}\right\rangle_{L_2([t,T]^2)}=
$$
$$
=4\sum_{j_1,j_2=0}^{\infty}
\int\limits_t^T
\psi_4(t_4)\phi_{j_2}(t_4)\int\limits_t^{t_4}
f_3(t_3)\phi_{j_2}(t_3)\int\limits_t^{T}
f_2(t_2)\phi_{j_1}(t_2)
\int\limits_t^{t_2}
\psi_1(t_1)\phi_{j_1}(t_1)\times
$$
$$
\times
dt_1 dt_2 dt_3 dt_4=
$$

\vspace{-2mm}
$$
=
\int\limits_{[t,T]^2} 
\lim\limits_{u\to 0}
A_u K'(t_2,t_2)K''(t_4,t_4)dt_2 dt_4=
$$
$$
=
\int\limits_{[t,T]^2} 
\lim\limits_{u\to 0}
A_u K'(t_2,t_2) \lim\limits_{u\to 0}
A_u K''(t_4,t_4)dt_2 dt_4=
$$
$$
=
\int\limits_{[t,T]^2} 
K'(t_2,t_2)K''(t_4,t_4)dt_2 dt_4=
$$
\begin{equation}
\label{july15014}
=\int\limits_{[t,T]^2} 
\psi_4(t_4)f_3(t_4)f_2(t_2)\psi_1(t_2)dt_2 dt_4.
\end{equation}

\vspace{2mm}

Recall that 
$f_2(\tau)$ and $f_3(\tau)$ are polynomials of finite degrees.
For example, $f_2(\tau)$ and $f_3(\tau)$ can be Legendre polynomials
that form a complete orthonormal system of functions in $L_2([t,T]).$

Denote
\begin{equation}
\label{july15015}
s_q(t_2,t_3)=\sum\limits_{l_1,l_2=0}^q
C_{l_2 l_1}\bar \phi_{l_1}(t_2)\bar \phi_{l_2}(t_3),
\end{equation}

\noindent 
where $\left\{\bar \phi_j(x)\right\}_{j=0}^{\infty}$ 
is a complete orthonormal system of Legendre polynomials in $L_2([t,T])$
and
$C_{l_2 l_1}$ are Fourier--Legendre coefficients for the function
$g(t_2,t_3)=\psi_2(t_2)\psi_3(t_3){\bf 1}_{\{t_2<t_3\}}$ 
($\psi_2(\tau), \psi_3(\tau)\in L_2([t,T])),$ i.e.
$$
C_{l_2 l_1}=\int\limits_t^T \psi_3(t_3)\bar \phi_{l_2}(t_3)
\int\limits_t^{t_3}\psi_2(t_2)\bar \phi_{l_1}(t_2)dt_2 dt_3.
$$

Further, we have
$$
\lim\limits_{q\to\infty}
\int\limits_{[t,T]^2}
\left(s_q(t_2,t_3)-g(t_2,t_3)\right)^2 dt_2 dt_3=0\ \ \ \hbox{or}\ \ \ 
\lim\limits_{q\to\infty}\left\Vert s_q - g\right\Vert^2_{L_2([t,T]^2)}=0.
$$

From (\ref{july15014}) we obtain (the sum on the right-hand side of (\ref{july15015}) is finite)
$$
\sum_{j_1,j_2=0}^{\infty} \left\langle
\Psi_{j_1 j_2}, \mathbb{K'}_q\Psi_{j_1 j_2}\right\rangle_{L_2([t,T]^2)}=
$$
$$
=4\sum_{j_1,j_2=0}^{\infty}
\int\limits_{[t,T]^4}{\bf 1}_{\{t_1<t_2\}}{\bf 1}_{\{t_3<t_4\}}
\psi_4(t_4)\phi_{j_2}(t_4)
s_q(t_2,t_3)\phi_{j_2}(t_3)
\phi_{j_1}(t_2)
\psi_1(t_1)\phi_{j_1}(t_1)\times
$$
$$
\times
dt_1 dt_2 dt_3 dt_4=
$$

\vspace{-2mm}

\begin{equation}
\label{july15017}
=\int\limits_{[t,T]^2} 
\psi_4(t_4)s_q(t_2,t_4)\psi_1(t_2)dt_2 dt_4,
\end{equation}

\vspace{1mm}
\noindent
where the operator $\mathbb{K'}_q$ (more precisely, its kernel)
is obtained 
from the operator $\mathbb{K'}$ (more precisely, from its kernel) by replacing 
$f_2f_3$ with $s_q$.

Note that the equality (\ref{july15017}) remains true
when $s_q$ is a partial sum of the Fourier--Legendre series
of any function from $L_2([t,T]^2),$ i.e. the equality holds
on a dense subset in $L_2([t,T]^2).$

Trace class operators form a linear space. Therefore,
on the left-hand side of 
(\ref{july15017}) there is a matrix trace of a trace class
operator $\mathbb{K'}_q$. 
The mentioned matrix trace
is a linear bounded (and therefore continuous)
functional
in the space of trace class operators \cite{gohb},
\cite{goldberg}
(this functional can be extended to the space $L_2([t, T]^2)$ by continuity \cite{Pugach}).

From the other hand, the right-hand side of (\ref{july15017}) defines
(as a scalar product of $s_q(t_2,t_4)$ and $\psi_4(t_4)\psi_1(t_2)$
in $L_2([t, T]^2)$) a linear bounded (and therefore continuous)
functional in $L_2([t, T]^2),$
which is given by the function $\psi_4(t_4)\psi_1(t_2)$.
On the left-hand side of (\ref{july15017}) (by virtue of the equality (\ref{july15017}))
there is a linear continuous functional on a dense subset in 
$L_2([t,T]^2).$ This functional can be uniquely extended 
to a linear continuous functional in $L_2([t, T]^2)$
(see \cite{reed}, Theorem~I.7, P.~9).

Let us implement the passage to the limit $\lim\limits_{q\to\infty}$
in the equality (\ref{july15017}) (at that we suppose that $s_q$ is defined by (\ref{july15015}))
$$
\sum_{j_1,j_2=0}^{\infty} \left\langle
\Psi_{j_1 j_2}, \mathbb{K''}\Psi_{j_1 j_2}\right\rangle_{L_2([t,T]^2)}=
$$
$$
=4\sum_{j_1,j_2=0}^{\infty}
\int\limits_{[t,T]^4}{\bf 1}_{\{t_1<t_2<t_3<t_4\}}
\psi_4(t_4)\psi_3(t_3)\psi_2(t_2)\psi_1(t_1)
\phi_{j_2}(t_4)
\phi_{j_2}(t_3)
\phi_{j_1}(t_2)
\phi_{j_1}(t_1)\times
$$
$$
\times
dt_1 dt_2 dt_3 dt_4=
$$

\vspace{-2mm}

\begin{equation}
\label{july15019}
=\int\limits_t^T 
\psi_4(t_4) \psi_3(t_4) \int\limits_t^{t_4} \psi_2(t_2)\psi_1(t_2)dt_2 dt_4,
\end{equation}

\vspace{1mm}
\noindent
where the operator $\mathbb{K''}$ (more precisely, its kernel)
is obtained 
from the operator $\mathbb{K'}_q$ (more precisely, from its kernel) by replacing 
$s_q$ with $\lim\limits_{q\to\infty}s_q=g\in L_2([t,T]^2)$,
$\psi_2(\tau), \psi_3(\tau)\in L_2([t,T])$ and 
$\psi_1(\tau), \psi_4(\tau)$ are continuously differentiable
functions on $[t, T].$

Further, the formula (\ref{july15019}) will remain valid
if we choose
$$
\psi_1(\tau)=\bar \psi_1^{(p)}(\tau),\ \ \ \psi_4(\tau)=\bar \psi_4^{(p)}(\tau),
$$
where
\begin{equation}
\label{july15018}
\bar \psi_1^{(p)}(\tau)=\sum\limits_{j=0}^p \bar \phi_j(\tau)\int\limits_t^T \bar \psi_1(s) 
\bar \phi_j(s)ds,\ \ 
\bar \psi_4^{(p)}(\tau)=\sum\limits_{j=0}^p \bar \phi_j(\tau)
\int\limits_t^T \bar \psi_4(s) \bar \phi_j(s)ds,
\end{equation}
where $p\in {\bf N},$ 
$\bar \psi_1(\tau), \bar \psi_4(\tau)\in L_2([t,T]),$ and $\left\{\bar \phi_j(x)\right\}_{j=0}^{\infty}$ 
is a complete orthonormal system of Legendre polynomials in $L_2([t,T])$.

Substitute (\ref{july15018}) into (\ref{july15019})
$$
\sum_{j_1,j_2=0}^{\infty} \left\langle
\Psi_{j_1 j_2}, \mathbb{K''}_p\Psi_{j_1 j_2}\right\rangle_{L_2([t,T]^2)}=
$$
$$
=4\hspace{-1.3mm}\sum_{j_1,j_2=0}^{\infty}
\int\limits_{[t,T]^4}\hspace{-1.3mm}{\bf 1}_{\{t_1<t_2<t_3<t_4\}}
\bar \psi_4^{(p)}(t_4)\psi_3(t_3)\psi_2(t_2)\bar\psi_1^{(p)}(t_1)
\phi_{j_2}(t_4)
\phi_{j_2}(t_3)
\phi_{j_1}(t_2)
\phi_{j_1}(t_1)\times
$$
$$
\times
dt_1 dt_2 dt_3 dt_4=
$$

\vspace{-2mm}

\begin{equation}
\label{july15020}
=\int\limits_t^T 
\bar \psi_4^{(p)}(t_4) \psi_3(t_4) \int\limits_t^{t_4} \psi_2(t_2)\bar \psi_1^{(p)}(t_2)dt_2 dt_4.
\end{equation}

\noindent
where the operator $\mathbb{K''}_p$ (more precisely, its kernel)
is obtained 
from the operator $\mathbb{K''}$ (more precisely, from its kernel) by replacing 
$\psi_4$ and $\psi_1$ with $\bar \psi_4^{(p)}$ and $\bar \psi_1^{(p)},$
respectively.

Note that the equality (\ref{july15020}) will also remain true if
$\bar \psi_4^{(p)}\bar \psi_1^{(p)}$ is replaced by $s_p$
($s_p$ is the partial sum of the Fourier--Legendre series
of any function from $L_2([t, T]^2)$), i.e. 
the modified equality (\ref{july15020}) is true 
on a dense subset of $L_2([t,T]^2).$
Next, we can apply the reasoning below the formula 
(\ref{july15017}) and obtain the equality of two linear continuous functionals
in $L_2([t,T]^2).$
Let us implement the passage to the limit $\lim\limits_{p\to\infty}$
in the mentioned equality under the condition $s_p=\bar \psi_4^{(p)}\bar \psi_1^{(p)}$
$$
4\sum_{j_1,j_2=0}^{\infty}
\int\limits_{[t,T]^4}{\bf 1}_{\{t_1<t_2<t_3<t_4\}}
\bar \psi_4(t_4)\psi_3(t_3)\psi_2(t_2)\bar\psi_1(t_1)
\phi_{j_2}(t_4)
\phi_{j_2}(t_3)
\phi_{j_1}(t_2)
\phi_{j_1}(t_1)\times
$$
$$
\times
dt_1 dt_2 dt_3 dt_4=
$$

\vspace{-2mm}

\begin{equation}
\label{july15021}
=\int\limits_t^T 
\bar \psi_4(t_4) \psi_3(t_4) \int\limits_t^{t_4} \psi_2(t_2)\bar \psi_1(t_2)dt_2 dt_4,
\end{equation}

\noindent
where $\bar \psi_1(\tau), \psi_2(\tau), \psi_3(\tau), \bar \psi_4(\tau)\in L_2([t,T]).$

Rewrite the equality (\ref{july15021}) in the form
$$
\lim\limits_{p\to\infty}\sum_{j_1,j_2=0}^{p}C_{j_2 j_2 j_1 j_1}=
$$
$$
=\sum_{j_1,j_2=0}^{\infty}
\int\limits_t^T
\psi_4(t_4)\phi_{j_2}(t_4) 
\int\limits_t^{t_4}
\psi_3(t_3)\phi_{j_2}(t_3)
\int\limits_t^{t_3}
\psi_2(t_2)\phi_{j_1}(t_2)
\int\limits_t^{t_2}
\psi_1(t_1)\phi_{j_1}(t_1)\times
$$
$$
\times
dt_1 dt_2 dt_3 dt_4=
$$

\vspace{-2mm}

\begin{equation}
\label{july15022}
=\frac{1}{4}\int\limits_t^T 
\psi_4(t_4) \psi_3(t_4) \int\limits_t^{t_4} \psi_2(t_2)\psi_1(t_2)dt_2 dt_4,
\end{equation}

\noindent
where $\psi_1(\tau), \ldots, \psi_4(\tau)\in L_2([t,T]).$

Note that the series on the left-hand side of 
(\ref{july15022}) converges absolutly since
its sum does not depend 
on permutations of basis functions
(here the basis in $L_2([t,T]^2)$ is
$\left\{\phi_{j_1}(x)\phi_{j_2}(y)\right\}_{j_1,j_2=0}^{\infty}$).
The equality (\ref{july13}) is proved. 

In \cite{rybakov7000x}, the equality (\ref{july15022})
is generalized as follows
$$
\lim\limits_{p\to\infty}\sum_{j_k,j_{k-2},\ldots, j_2=0}^{p}
C_{j_k j_k j_{k-2} j_{k-2} \ldots j_2 j_2}=
$$
$$
=\frac{1}{2^r}\hspace{-0.4mm}\int\limits_t^T \hspace{-0.4mm}
\psi_k(t_k) \psi_{k-1}(t_k) 
\int\limits_t^{t_k} \psi_{k-2}(t_{k-2})\psi_{k-3}(t_{k-2})
\ldots
$$
\begin{equation}
\label{july15023}
\ldots
\int\limits_t^{t_4} \psi_{2}(t_{2})\psi_{1}(t_{2})dt_2 \ldots 
dt_{k-2} dt_k,
\end{equation}

\noindent
where $k=2r$ $(r=2,3,\ldots),$
$\psi_1(\tau),\ldots, \psi_k(\tau)\in L_2([t,T]).$

The equalities (\ref{july14}), (\ref{july15}) can also be obtained
\cite{rybakov7000xa}
using the approach from \cite{rybakov7000x} and the series
on the left-hand sides of (\ref{july14}), (\ref{july15})
converge absolutely.

In the notations of Theorem~2.49, the equality
(\ref{july15023}) can be written in the form
$$
\lim\limits_{p\to\infty}\sum\limits_{j_{1},j_{3},\ldots, j_{2r-1}=0}^p
C_{j_k\ldots j_1}\biggl|_{j_{1}=j_{2},\ldots, j_{2r-1}=j_{2r}}\biggr.=
$$
\begin{equation}
\label{july16000}
~~~~~~~~~=\frac{1}{2^r} 
C_{j_k \ldots j_1}\biggl|_{(j_{2} j_{1})\curvearrowright (\cdot)
(j_{4} j_{3})\curvearrowright (\cdot)
\ldots (j_{2r} j_{2r-1})\curvearrowright (\cdot),
j_1=j_2,j_3=j_4,\ldots, j_{2r-1}=j_{2r}
}\biggr.,
\end{equation}

\vspace{4mm}
\noindent
where $k=2r$ $(r=2,3,\ldots)$ and $C_{j_k\ldots j_1}$ is defined by (\ref{july15030}).

In principle, using the method from \cite{rybakov7000x}
the following equality can be obtained \cite{rybakov7000xa}
$$
\lim\limits_{p\to\infty}
\sum\limits_{j_{g_1},j_{g_3},\ldots,j_{g_{2r-1}}=0}^p
C_{j_k\ldots j_1}\biggl|_{j_{g_1}=j_{g_2},\ldots, j_{g_{2r-1}}=j_{g_{2r}}}=
$$
$$
=\frac{1}{2^r} \prod\limits_{l=1}^r {\bf 1}_{\{g_{2l}=g_{2l-1}+1\}}
C_{j_k \ldots j_1}\biggl|_{(j_{g_2} j_{g_1})\curvearrowright (\cdot)
\ldots (j_{g_{2r}} j_{g_{2r-1}})\curvearrowright (\cdot),
j_{g_{{}_{1}}}=~j_{g_{{}_{2}}},\ldots, j_{g_{{}_{2r-1}}}=~j_{g_{{}_{2r}}}
}\biggr.
$$

\vspace{4mm}
\noindent
for all possible $g_1,g_2,\ldots,g_{2r-1},g_{2r}$ (see {\rm (\ref{leto5007after})),
where $k=2r$ $(r=2,3,\ldots),$ $C_{j_k\ldots j_1}$ is defined by (\ref{july15030}),
another notations are the same as in Theorem~2.49.

\subsection{Approach Based on Generalized Parseval's Equality
and (\ref{start1000}).
General Case When $\psi_1(\tau),\ldots,\psi_k(\tau)\in L_2([t, T])$ 
$(k=2r, r=2,3,4,\ldots)$}

Let us prove the equalities (\ref{july13})--(\ref{july15})
using a method based on generalized Parseval's equality and (\ref{start1000}).

Consider (\ref{july13}). Using (\ref{start1000}), we have

\vspace{-4mm}
$$
\lim\limits_{p\to\infty}\sum_{j_1,j_2=0}^{p}
\int\limits_t^T
\psi_4(t_4)\phi_{j_2}(t_4)\int\limits_t^{t_4}
\psi_3(t_3)\phi_{j_2}(t_3)\int\limits_t^{T}
\psi_2(t_2)\phi_{j_1}(t_2)
\int\limits_{t}^{t_2}
\psi_1(t_1)\phi_{j_1}(t_1)\times\Biggr.
$$
$$
\times
dt_1 dt_2 dt_3 dt_4=
$$

\vspace{-4mm}
$$
=\lim\limits_{p\to\infty}\sum_{j_1,j_2=0}^{p}
\int\limits_t^T
\psi_4(t_4)\phi_{j_2}(t_4)\int\limits_t^{t_4}
\psi_3(t_3)\phi_{j_2}(t_3)dt_3 dt_4\times
$$
$$
\times
\int\limits_t^{T}
\psi_2(t_2)\phi_{j_1}(t_2)
\int\limits_{t}^{t_2}
\psi_1(t_1)\phi_{j_1}(t_1)
dt_1 dt_2=
$$
$$
=\lim\limits_{p\to\infty}\sum_{j_2=0}^{p}
\int\limits_t^T
\psi_4(t_4)\phi_{j_2}(t_4)\int\limits_t^{t_4}
\psi_3(t_3)\phi_{j_2}(t_3)dt_3 dt_4\times
$$
$$
\times
\lim\limits_{p\to\infty}\sum_{j_1=0}^{p}\int\limits_t^{T}
\psi_2(t_2)\phi_{j_1}(t_2)
\int\limits_{t}^{t_2}
\psi_1(t_1)\phi_{j_1}(t_1)
dt_1 dt_2=
$$
$$
=\frac{1}{4}\int\limits_t^T
\psi_4(t_4)\psi_3(t_4)dt_4
\int\limits_t^T
\psi_2(t_2)\psi_1(t_2)dt_2=
$$
\begin{equation}
\label{july29000}
=\frac{1}{4}\int\limits_{[t,T]^2}
\psi_4(t_4)\psi_3(t_4)
\psi_2(t_2)\psi_1(t_2)dt_2 dt_4,
\end{equation}

\noindent
where $\psi_1(\tau),\ldots,\psi_4(\tau)\in L_2([t, T]).$ 

Suppose that $\psi_2(\tau)$ and $\psi_3(\tau)$ are polynomials of finite degrees.
For example, $\psi_2(\tau)$ and $\psi_3(\tau)$ can be Legendre polynomials
that form a complete orthonormal system of functions in $L_2([t,T]).$

Denote
\begin{equation}
\label{july30000}
s_q(t_2,t_3)=\sum\limits_{l_1,l_2=0}^q
C_{l_2 l_1}\bar \phi_{l_1}(t_2)\bar \phi_{l_2}(t_3),
\end{equation}

\vspace{1mm}
\noindent 
where $\left\{\bar \phi_j(x)\right\}_{j=0}^{\infty}$ 
is a complete orthonormal system of Legendre polynomials in $L_2([t,T])$
and
$C_{l_2 l_1}$ are Fourier--Legendre coefficients for the function
$g(t_2,t_3)=\bar \psi_2(t_2)\bar \psi_3(t_3){\bf 1}_{\{t_2<t_3\}}$ 
($\bar \psi_2(\tau), \bar \psi_3(\tau)\in L_2([t,T])),$ i.e.
$$
C_{l_2 l_1}=\int\limits_t^T \bar \psi_3(t_3)\bar \phi_{l_2}(t_3)
\int\limits_t^{t_3}\bar \psi_2(t_2)\bar \phi_{l_1}(t_2)dt_2 dt_3
$$
and
$\lim\limits_{q\to\infty}\left\Vert s_q - g\right\Vert^2_{L_2([t,T]^2)}=0.$

From (\ref{july29000}) we obtain (the sum on the right-hand side of (\ref{july30000}) is finite)
$$
\sum_{j_1,j_2=0}^{\infty}
\int\limits_{[t,T]^4}{\bf 1}_{\{t_1<t_2\}}{\bf 1}_{\{t_3<t_4\}}
\psi_4(t_4)\phi_{j_2}(t_4)
s_q(t_2,t_3)\phi_{j_2}(t_3)
\phi_{j_1}(t_2)
\psi_1(t_1)\phi_{j_1}(t_1)\times
$$
$$
\times
dt_1 dt_2 dt_3 dt_4=
$$           
\begin{equation}
\label{july30001}
=\frac{1}{4}\int\limits_{[t,T]^2} 
\psi_4(t_4)s_q(t_2,t_4)\psi_1(t_2)dt_2 dt_4.
\end{equation}

\vspace{2mm}

Note that the equality (\ref{july30001}) remains true
when $s_q$ is a partial sum of the Fourier--Legendre series
of any function from $L_2([t,T]^2),$ i.e. the equality holds
on a dense subset in $L_2([t,T]^2).$

The right-hand side of (\ref{july30001}) defines
(as a scalar product of $s_q(t_2,t_4)$ and $\frac{1}{4}\psi_4(t_4)\psi_1(t_2)$
in $L_2([t, T]^2)$) a linear bounded (and therefore continuous)
functional in $L_2([t, T]^2),$
which is given by the function $\frac{1}{4}\psi_4(t_4)\psi_1(t_2)$.
On the left-hand side of (\ref{july30001}) (by virtue of the equality (\ref{july30001}))
there is a linear continuous functional on a dense subset in 
$L_2([t,T]^2).$ This functional can be uniquely extended 
to a linear continuous functional in $L_2([t, T]^2)$
(see \cite{reed}, Theorem~I.7, P.~9).

Let us implement the passage to the limit $\lim\limits_{q\to\infty}$
in (\ref{july30001}) (at that we suppose that $s_q$ is defined by (\ref{july30000}))
$$
\sum_{j_1,j_2=0}^{\infty}
\int\limits_{[t,T]^4}{\bf 1}_{\{t_1<t_2<t_3<t_4\}}
\psi_4(t_4)\bar \psi_3(t_3)\bar \psi_2(t_2)\psi_1(t_1)
\phi_{j_2}(t_4)
\phi_{j_2}(t_3)
\phi_{j_1}(t_2)
\phi_{j_1}(t_1)\times
$$
$$
\times
dt_1 dt_2 dt_3 dt_4=
$$

\vspace{-2mm}

\begin{equation}
\label{july30002}
=\frac{1}{4}\int\limits_t^T 
\psi_4(t_4) \bar\psi_3(t_4) \int\limits_t^{t_4} \bar\psi_2(t_2)\psi_1(t_2)dt_2 dt_4,
\end{equation}

\noindent
where $\psi_1(\tau), \bar \psi_2(\tau), \bar \psi_3(\tau), \psi_4(\tau)\in L_2([t,T]).$

Rewrite the equality (\ref{july30002}) in the form
$$
\lim\limits_{p\to\infty}\sum_{j_1,j_2=0}^{p}C_{j_2 j_2 j_1 j_1}=
$$
$$
=\sum_{j_1,j_2=0}^{\infty}
\int\limits_t^T
\psi_4(t_4)\phi_{j_2}(t_4) 
\int\limits_t^{t_4}
\psi_3(t_3)\phi_{j_2}(t_3)
\int\limits_t^{t_3}
\psi_2(t_2)\phi_{j_1}(t_2)
\int\limits_t^{t_2}
\psi_1(t_1)\phi_{j_1}(t_1)\times
$$
$$
\times
dt_1 dt_2 dt_3 dt_4=
$$

\vspace{-2mm}

\begin{equation}
\label{july30003}
=\frac{1}{4}\int\limits_t^T 
\psi_4(t_4) \psi_3(t_4) \int\limits_t^{t_4} \psi_2(t_2)\psi_1(t_2)dt_2 dt_4,
\end{equation}

\noindent
where $\psi_1(\tau), \ldots, \psi_4(\tau)\in L_2([t,T]).$

Note that the series on the left-hand side of 
(\ref{july30003}) converges absolutly since
its sum does not depend 
on permutations of basis functions
(here the basis in $L_2([t,T]^2)$ is
$\left\{\phi_{j_1}(x)\phi_{j_2}(y)\right\}_{j_1,j_2=0}^{\infty}$).
The equality (\ref{july13}) is proved.

Let us prove (\ref{july15}). Using the generalized Parseval equality, we obtain
$$
\lim\limits_{p\to\infty}\sum_{j_1,j_2=0}^{p}
\int\limits_t^T
\psi_4(t_4)\phi_{j_2}(t_4)\int\limits_t^{t_4}
\psi_3(t_3)\phi_{j_1}(t_3)\int\limits_t^{T}
\psi_2(t_2)\phi_{j_2}(t_2)
\int\limits_{t}^{t_2}
\psi_1(t_1)\phi_{j_1}(t_1)\times\Biggr.
$$
$$
\times
dt_1 dt_2 dt_3 dt_4=
$$

\vspace{-2mm}
$$
=\sum_{j_1,j_2=0}^{\infty}
\int\limits_t^T
\psi_4(t_4)\phi_{j_2}(t_4)\int\limits_t^{t_4}
\psi_3(t_3)\phi_{j_1}(t_3)dt_3 dt_4\times
$$
$$
\times
\int\limits_t^{T}
\psi_2(t_2)\phi_{j_2}(t_2)
\int\limits_{t}^{t_2}
\psi_1(t_1)\phi_{j_1}(t_1)
dt_1 dt_2=
$$
$$
=\sum_{j_1,j_2=0}^{\infty}
\int\limits_{[t,T]^2}
{\bf 1}_{\{t_3<t_4\}}\psi_3(t_3)\psi_4(t_4)\phi_{j_1}(t_3)
\phi_{j_2}(t_4)dt_3 dt_4\times
$$
$$
\times
\int\limits_{[t,T]^2}
{\bf 1}_{\{t_3<t_4\}}\psi_1(t_3)\psi_2(t_4)\phi_{j_1}(t_3)
\phi_{j_2}(t_4)dt_3 dt_4=
$$
$$
=
\int\limits_{[t,T]^2}
{\bf 1}_{\{t_3<t_4\}}\psi_3(t_3)\psi_2(t_4)\psi_4(t_4)\psi_1(t_3)
dt_3 dt_4=
$$
\begin{equation}
\label{july30004}
=
\int\limits_{[t,T]^2}
{\bf 1}_{\{t_3<t_2\}}\psi_3(t_3)\psi_2(t_2)\psi_4(t_2)\psi_1(t_3)
dt_3 dt_2,
\end{equation}

\noindent
where $\psi_1(\tau), \psi_2(\tau), \psi_3(\tau), \psi_4(\tau)\in L_2([t, T]).$

Suppose that $\psi_2(\tau)$ and $\psi_3(\tau)$ are Legendre polynomials of finite degrees.
Denote
\begin{equation}
\label{july30005}
s_q(t_2,t_3)=\sum\limits_{l_1,l_2=0}^q
C_{l_2 l_1}\bar \phi_{l_1}(t_2)\bar \phi_{l_2}(t_3),
\end{equation}

\noindent 
where $\left\{\bar \phi_j(x)\right\}_{j=0}^{\infty}$ 
is a complete orthonormal system of Legendre polynomials in $L_2([t,T])$
and
$C_{l_2 l_1}$ are Fourier--Legendre coefficients for the function
$g(t_2,t_3)=\bar \psi_2(t_2)\bar \psi_3(t_3){\bf 1}_{\{t_2<t_3\}}$ 
($\bar \psi_2(\tau), \bar \psi_3(\tau)\in L_2([t,T])),$ i.e.
$$
C_{l_2 l_1}=\int\limits_t^T \bar \psi_3(t_3)\bar \phi_{l_2}(t_3)
\int\limits_t^{t_3}\bar \psi_2(t_2)\bar \phi_{l_1}(t_2)dt_2 dt_3
$$
and 
$\lim\limits_{q\to\infty}\left\Vert s_q - g\right\Vert^2_{L_2([t,T]^2)}=0.$

\vspace{1mm}

From (\ref{july30004}) we obtain (the sum on the right-hand side of (\ref{july30005}) is finite)
$$
\sum_{j_1,j_2=0}^{\infty}
\int\limits_{[t,T]^4}{\bf 1}_{\{t_1<t_2\}}{\bf 1}_{\{t_3<t_4\}}
\psi_4(t_4)
s_q(t_2,t_3)\psi_1(t_1)\phi_{j_2}(t_4)\phi_{j_1}(t_3)
\phi_{j_2}(t_2)
\phi_{j_1}(t_1)\times
$$
$$
\times
dt_1 dt_2 dt_3 dt_4=
$$

\vspace{-2mm}
\begin{equation}
\label{july30006}
=
\int\limits_{[t,T]^2}
{\bf 1}_{\{t_3<t_2\}}s_q(t_2,t_3)\psi_1(t_3)\psi_4(t_2)
dt_3 dt_2.
\end{equation}

\vspace{1mm}

Note that the equality (\ref{july30006}) remains true
when $s_q$ is a partial sum of the Fourier--Legendre series
of any function from $L_2([t,T]^2),$ i.e. the equality holds
on a dense subset in $L_2([t,T]^2).$

The right-hand side of (\ref{july30006}) defines
(as a scalar product of $s_q(t_2,t_3)$ and ${\bf 1}_{\{t_3<t_2\}}\psi_1(t_3)\psi_4(t_2)$
in $L_2([t, T]^2)$) a linear bounded (and therefore continuous)
functional in $L_2([t, T]^2),$
which is given by the function ${\bf 1}_{\{t_3<t_2\}}\psi_1(t_3)\psi_4(t_2)$.
On the left-hand side of (\ref{july30006}) (by virtue of the equality (\ref{july30006}))
there is a linear continuous functional on a dense subset in 
$L_2([t,T]^2).$ This functional can be uniquely extended 
to a linear continuous functional in $L_2([t, T]^2)$
(see \cite{reed}, Theorem~I.7, P.~9).

Let us implement the passage to the limit $\lim\limits_{q\to\infty}$
in (\ref{july30006}) (at that we suppose that $s_q$ is defined by (\ref{july30005}))
$$
\sum_{j_1,j_2=0}^{\infty}
\int\limits_{[t,T]^4}{\bf 1}_{\{t_1<t_2<t_3<t_4\}}
\psi_4(t_4)
\bar\psi_3(t_3)\bar \psi_2(t_2)\psi_1(t_1)\phi_{j_2}(t_4)\phi_{j_1}(t_3)
\phi_{j_2}(t_2)
\phi_{j_1}(t_1)\times
$$
$$
\times
dt_1 dt_2 dt_3 dt_4=
$$
\begin{equation}
\label{july30007}
~~~~~~~~=
\int\limits_{[t,T]^2}
{\bf 1}_{\{t_2>t_3\}}{\bf 1}_{\{t_2<t_3\}} \bar\psi_3(t_3)\bar \psi_2(t_2)\psi_1(t_3)\psi_4(t_2)
dt_3 dt_2=0.
\end{equation}

\vspace{2mm}

Rewrite the equality (\ref{july30007}) in the form
$$
\lim\limits_{p\to\infty}\sum_{j_1,j_2=0}^{p}C_{j_2 j_1 j_2 j_1}=
$$
$$
=\sum_{j_1,j_2=0}^{\infty}
\int\limits_t^T
\psi_4(t_4)\phi_{j_2}(t_4) 
\int\limits_t^{t_4}
\psi_3(t_3)\phi_{j_1}(t_3)
\int\limits_t^{t_3}
\psi_2(t_2)\phi_{j_2}(t_2)
\int\limits_t^{t_2}
\psi_1(t_1)\phi_{j_1}(t_1)\times
$$
\begin{equation}
\label{july30008}
\times
dt_1 dt_2 dt_3 dt_4
=0,
\end{equation}

\vspace{3mm}
\noindent
where $\psi_1(\tau), \ldots, \psi_4(\tau)\in L_2([t,T]).$

Note that the series on the left-hand side of 
(\ref{july30008}) converges absolutly since
its sum does not depend 
on permutations of basis functions
(here the basis in $L_2([t,T]^2)$ is
$\left\{\phi_{j_1}(x)\phi_{j_2}(y)\right\}_{j_1,j_2=0}^{\infty}$).
The equality (\ref{july15}) is proved.

Let us prove (\ref{july14}). Using the generalized Parseval's
equality, we get
$$
\lim\limits_{p\to\infty}\sum_{j_1,j_2=0}^{p}
\int\limits_t^T
\psi_4(t_4)\phi_{j_2}(t_4)\int\limits_t^{t_4}
\psi_3(t_3)\phi_{j_1}(t_3)\int\limits_t^T
\psi_2(t_2)\phi_{j_1}(t_2)
\int\limits_{t}^{t_2}
\psi_1(t_1)\phi_{j_2}(t_1)\times\Biggr.
$$
$$
\times
dt_1 dt_2 dt_3 dt_4=
$$

\vspace{-2mm}
$$
=\sum_{j_1,j_2=0}^{\infty}
\int\limits_t^T
\psi_4(t_4)\phi_{j_2}(t_4)\int\limits_t^{t_4}
\psi_3(t_3)\phi_{j_1}(t_3)dt_3 dt_4\times
$$
$$
\times
\int\limits_t^{T}
\psi_2(t_2)\phi_{j_1}(t_2)
\int\limits_{t}^{t_2}
\psi_1(t_1)\phi_{j_2}(t_1)
dt_1 dt_2=
$$

\vspace{-2mm}
$$
=\sum_{j_1,j_2=0}^{\infty}
\int\limits_{[t,T]^2}
{\bf 1}_{\{t_3<t_4\}}\psi_3(t_3)\psi_4(t_4)\phi_{j_1}(t_3)
\phi_{j_2}(t_4)dt_3 dt_4\times
$$
$$
\times
\int\limits_{[t,T]^2}
{\bf 1}_{\{t_4<t_3\}}\psi_1(t_4)\psi_2(t_3)\phi_{j_1}(t_3)
\phi_{j_2}(t_4)dt_4 dt_3=
$$
\begin{equation}
\label{july30004fgf}
~~~~~~~~=
\int\limits_{[t,T]^2}
{\bf 1}_{\{t_3<t_4\}}{\bf 1}_{\{t_4<t_3\}}\psi_3(t_3)\psi_4(t_4)\psi_1(t_4)\psi_2(t_3)
dt_3 dt_4=0,
\end{equation}

\vspace{1mm}
\noindent
where $\psi_1(\tau), \psi_2(\tau), \psi_3(\tau), \psi_4(\tau)\in L_2([t, T]).$

Applying (\ref{july30004fgf}), we obtain
$$
\sum_{j_1,j_2=0}^{\infty}
\int\limits_{[t,T]^4}
{\bf 1}_{\{t_1<t_2\}}{\bf 1}_{\{t_3<t_4\}}
\psi_4(t_4)\phi_{j_2}(t_4)
\psi_3(t_3)\phi_{j_1}(t_3)
\psi_2(t_2)\phi_{j_1}(t_2)
\psi_1(t_1)\phi_{j_2}(t_1)\times\Biggr.
$$
\begin{equation}
\label{july30004fgf1}
\times
dt_1 dt_2 dt_3 dt_4=0.
\end{equation}

\vspace{3mm}

Suppose that $\psi_2(\tau)$ and $\psi_3(\tau)$ are Legendre polynomials of finite degrees.
Denote
\begin{equation}
\label{july30005sss}
s_q(t_2,t_3)=\sum\limits_{l_1,l_2=0}^q
C_{l_2 l_1}\bar \phi_{l_1}(t_2)\bar \phi_{l_2}(t_3),
\end{equation}

\noindent 
where $\left\{\bar \phi_j(x)\right\}_{j=0}^{\infty}$ 
is a complete orthonormal system of Legendre polynomials in $L_2([t,T])$
and
$C_{l_2 l_1}$ are Fourier--Legendre coefficients for the function
$g(t_2,t_3)=\bar \psi_2(t_2)\bar \psi_3(t_3){\bf 1}_{\{t_2<t_3\}}$ 
($\bar \psi_2(\tau), \bar \psi_3(\tau)\in L_2([t,T])),$ i.e.
$$
C_{l_2 l_1}=\int\limits_t^T \bar \psi_3(t_3)\bar \phi_{l_2}(t_3)
\int\limits_t^{t_3}\bar \psi_2(t_2)\bar \phi_{l_1}(t_2)dt_2 dt_3
$$
and 
$\lim\limits_{q\to\infty}\left\Vert s_q - g\right\Vert^2_{L_2([t,T]^2)}=0.$

\vspace{1mm}

From (\ref{july30004fgf1}) we obtain (the sum on the right-hand side of (\ref{july30005sss}) is finite)
$$
\sum_{j_1,j_2=0}^{\infty}
\int\limits_{[t,T]^4}
{\bf 1}_{\{t_1<t_2\}}{\bf 1}_{\{t_3<t_4\}}
\psi_4(t_4)s_q(t_2,t_3)\psi_1(t_1)
\phi_{j_2}(t_4)
\phi_{j_1}(t_3)
\phi_{j_1}(t_2)
\phi_{j_2}(t_1)\times\Biggr.
$$
\begin{equation}
\label{july30004fgf2}
\times
dt_1 dt_2 dt_3 dt_4=0.
\end{equation}

\vspace{3mm}

Note that the equality (\ref{july30004fgf2}) remains true
when $s_q$ is a partial sum of the Fourier--Legendre series
of any function from $L_2([t,T]^2),$ i.e. the equality holds
on a dense subset in $L_2([t,T]^2).$

The right-hand side of (\ref{july30004fgf2}) 
is the zero functional in $L_2([t, T]^2)$
(the zero functional in $L_2([t, T]^2)$ is a linear bounded (and therefore continuous)
functional in $L_2([t, T]^2)$).
On the left-hand side of (\ref{july30004fgf2}) (by virtue of the equality (\ref{july30004fgf2}))
there is a linear continuous functional on a dense subset in 
$L_2([t,T]^2).$ This functional can be uniquely extended 
to a linear continuous functional in $L_2([t, T]^2)$
(see \cite{reed}, Theorem~I.7, P.~9).

Let us implement the passage to the limit $\lim\limits_{q\to\infty}$
in (\ref{july30004fgf2}) (at that we suppose that $s_q$ is defined by (\ref{july30005sss}))
$$
\sum_{j_1,j_2=0}^{\infty}
\int\limits_{[t,T]^4}{\bf 1}_{\{t_1<t_2<t_3<t_4\}}
\psi_4(t_4)\phi_{j_2}(t_4)\bar \psi_3(t_3)\phi_{j_1}(t_3)
\bar \psi_2(t_2)\phi_{j_1}(t_2)
\psi_{1}(t_1)\phi_{j_2}(t_1)\times
$$
\begin{equation}
\label{july30014}
\times
dt_1 dt_2 dt_3 dt_4=0.
\end{equation}

\vspace{3mm}

Rewrite the equality (\ref{july30014}) in the form
$$
\lim\limits_{p\to\infty}\sum_{j_1,j_3=0}^{p}C_{j_1 j_3 j_3 j_1}=
$$
$$
=\sum_{j_1,j_3=0}^{\infty}
\int\limits_t^T
\psi_4(t_4)\phi_{j_1}(t_4) 
\int\limits_t^{t_4}
\psi_3(t_3)\phi_{j_3}(t_3)
\int\limits_t^{t_3}
\psi_2(t_2)\phi_{j_3}(t_2)
\int\limits_t^{t_2}
\psi_1(t_1)\phi_{j_1}(t_1)\times
$$
\begin{equation}
\label{july30015}
\times
dt_1 dt_2 dt_3 dt_4
=0,
\end{equation}

\vspace{3mm}
\noindent
where $\psi_1(\tau), \ldots, \psi_4(\tau)\in L_2([t,T]).$

Note that the series on the left-hand side of 
(\ref{july30015}) converges absolutly since
its sum does not depend 
on permutations of basis functions
(here the basis in $L_2([t,T]^2)$ is
$\left\{\phi_{j_1}(x)\phi_{j_2}(y)\right\}_{j_1,j_2=0}^{\infty}$).
The equality (\ref{july14}) is proved. The equalities
(\ref{july13})--(\ref{july15}) are proved.

By induction we prove the following equality (i.e. by a different method com\-pa\-red 
with \cite{rybakov7000x})
$$
\lim\limits_{p\to\infty}\sum_{j_{2r}, j_{2r-2}, \ldots, j_2=0}^{p}
C_{j_{2r}j_{2r} j_{2r-2}j_{2r-2} \ldots j_2 j_2}=
$$
$$
=
\frac{1}{2^{r}}
\int\limits_t^T
\psi_{2r}(t_{2r})
\psi_{2r-1}(t_{2r})\times
$$
\begin{equation}
\label{july30016}
~~~~~~~~\times
\int\limits_t^{t_{2r}}
\psi_{2r-2}(t_{2r-2})
\psi_{2r-3}(t_{2r-2})\ldots
\int\limits_t^{t_{4}}
\psi_{2}(t_{2})
\psi_{1}(t_{2})dt_2\ldots dt_{2r-2}dt_{2r},
\end{equation}

\noindent
where $r\in{\bf N},$ $C_{j_{2r}j_{2r} j_{2r-2}j_{2r-2} \ldots j_2 j_2}$
is defined by 
$$
C_{j_k \ldots j_1}=\int\limits_t^T\psi_k(t_k)\phi_{j_k}(t_k)\ldots
\int\limits_t^{t_2}
\psi_1(t_1)\phi_{j_1}(t_1)
dt_1\ldots dt_k\ \ \ (k\in{\bf N}),
$$
$\left\{\phi_j(x)\right\}_{j=0}^{\infty}$
is an arbitrary complete orthonormal system of 
functions in the space $L_2([t,T]),$ and
$\psi_1(\tau),\ldots ,\psi_{2r}(\tau)\in $ $L_2([t, T]).$

Note that the equality (\ref{july13}) is a particular case of 
(\ref{july30016}) for $r=2$ and the equality (\ref{start1000}) is a particular case of 
(\ref{july30016}) for $r=1$.
Thus, the equality
(\ref{july30016}) is true for $r=1, 2.$
Suppose that the equality (\ref{july30016}) is true for some 
$r>2.$ Then, using (\ref{start1000}), we get 
$$
\lim\limits_{p\to\infty}\sum_{j_{2r+2}, j_{2r}, \ldots, j_2=0}^{p}
\int\limits_t^T
\psi_{2r+2}(t_{2r+2})
\phi_{j_{2r+2}}(t_{2r+2})
\int\limits_t^{t_{2r+2}}
\psi_{2r+1}(t_{2r+1})
\phi_{j_{2r+2}}(t_{2r+1})\times
$$
$$
\times
\int\limits_t^T
\psi_{2r}(t_{2r})
\phi_{j_{2r}}(t_{2r})
\int\limits_t^{t_{2r}}
\psi_{2r-1}(t_{2r-1})
\phi_{j_{2r}}(t_{2r-1})\ldots
$$
$$
\ldots 
\int\limits_t^{t_3}
\psi_{2}(t_{2})
\phi_{j_{2}}(t_{2})
\int\limits_t^{t_{2}}
\psi_{1}(t_{1})
\phi_{j_{2}}(t_{1})
dt_1 dt_2\ldots dt_{2r-1}dt_{2r}dt_{2r+1}dt_{2r+2}=
$$
$$
=
\sum_{j_{2r+2}=0}^{\infty}
\int\limits_t^T
\psi_{2r+2}(t_{2r+2})
\phi_{j_{2r+2}}(t_{2r+2})
\int\limits_t^{t_{2r+2}}
\psi_{2r+1}(t_{2r+1})
\phi_{j_{2r+2}}(t_{2r+1})dt_{2r+1}dt_{2r+2}\times
$$
$$
\times\sum_{j_{2r}, j_{2r-2}, \ldots, j_2=0}^{\infty}
\int\limits_t^T
\psi_{2r}(t_{2r})
\phi_{j_{2r}}(t_{2r})
\int\limits_t^{t_{2r}}
\psi_{2r-1}(t_{2r-1})
\phi_{j_{2r}}(t_{2r-1})\times
$$
$$
\times
\int\limits_t^{t_{2r-1}}
\psi_{2r-2}(t_{2r-2})
\phi_{j_{2r-2}}(t_{2r-2})
\int\limits_t^{t_{2r-2}}
\psi_{2r-3}(t_{2r-3})
\phi_{j_{2r-2}}(t_{2r-3})\ldots
$$
$$
\ldots 
\int\limits_t^{t_3}
\psi_{2}(t_{2})
\phi_{j_{2}}(t_{2})
\int\limits_t^{t_{2}}
\psi_{1}(t_{1})
\phi_{j_{2}}(t_{1})
dt_1 dt_2\ldots dt_{2r-3}dt_{2r-2}dt_{2r-1}dt_{2r}=
$$
$$
=\frac{1}{2}
\int\limits_t^T
\psi_{2r+2}(t_{2r+2})
\psi_{2r+1}(t_{2r+2})
dt_{2r+2}\cdot
\frac{1}{2^{r}}
\int\limits_t^T
\psi_{2r}(t_{2r})
\psi_{2r-1}(t_{2r})\times
$$
\begin{equation}
\label{july30017}
~~~~~~~~~\times\int\limits_t^{t_{2r}}
\psi_{2r-2}(t_{2r-2})
\psi_{2r-3}(t_{2r-2})\ldots
\int\limits_t^{t_{4}}
\psi_{2}(t_{2})
\psi_{1}(t_{2})dt_2\ldots dt_{2r-2}dt_{2r}.
\end{equation}

\vspace{1mm}

Let us rewrite the equality (\ref{july30017}) in the form
$$
\lim\limits_{p\to\infty}\sum_{j_{2r+2}, j_{2r}, \ldots, j_2=0}^{p}
\int\limits_t^T
\psi_{2r+2}(t_{2r+2})
\phi_{j_{2r+2}}(t_{2r+2})
\int\limits_t^{t_{2r+2}}
\psi_{2r+1}(t_{2r+1})
\phi_{j_{2r+2}}(t_{2r+1})\times
$$
$$
\times
\int\limits_t^T
\psi_{2r}(t_{2r})
\phi_{j_{2r}}(t_{2r})
\int\limits_t^{t_{2r}}
\psi_{2r-1}(t_{2r-1})
\phi_{j_{2r}}(t_{2r-1})\ldots
$$
$$
\ldots 
\int\limits_t^{t_3}
\psi_{2}(t_{2})
\phi_{j_{2}}(t_{2})
\int\limits_t^{t_{2}}
\psi_{1}(t_{1})
\phi_{j_{2}}(t_{1})
dt_1 dt_2\ldots dt_{2r-1}dt_{2r}dt_{2r+1}dt_{2r+2}=
$$
$$
=\frac{1}{2^{r+1}}
\int\limits_t^T
\psi_{2r+2}(t_{2r+2})
\psi_{2r+1}(t_{2r+2})
\int\limits_t^T
\psi_{2r}(t_{2r})
\psi_{2r-1}(t_{2r})\times 
$$
\begin{equation}
\label{july30018}
\times 
\int\limits_t^{t_{2r}}
\psi_{2r-2}(t_{2r-2})
\psi_{2r-3}(t_{2r-2})\ldots
\int\limits_t^{t_{4}}
\psi_{2}(t_{2})
\psi_{1}(t_{2})dt_2\ldots dt_{2r-2}dt_{2r}dt_{2r+2},
\end{equation}
where $\psi_1(\tau),\ldots, \psi_{2r+2}(\tau)\in L_2([t, T]).$

Suppose that $\psi_{1}(\tau),\psi_{3}(\tau),\ldots, \psi_{2r-3}(\tau),
\psi_{2r}(\tau), \psi_{2r+1}(\tau)$ in (\ref{july30018}) are Legendre polynomials
of finite degrees.
Denote
$$
h(t_2,t_4,\ldots,t_{2r-2},t_{2r-1},t_{2r+2})=
$$
$$
=\psi_2(t_2)\psi_4(t_4)\ldots
\psi_{2r-2}(t_{2r-2})\psi_{2r-1}(t_{2r-1})\psi_{2r+2}(t_{2r+2}),
$$

\vspace{-2mm}
$$
g(t_1,t_3,\ldots,t_{2r-3},t_{2r},t_{2r+1})=
$$
\begin{equation}
\label{july50000}
~~~~~~~~~~=
\bar \psi_{1}(t_1)\bar \psi_{3}(t_3)\ldots \bar\psi_{2r-3}(t_{2r-3})
\bar \psi_{2r}(t_{2r})\bar \psi_{2r+1}(t_{2r+1})
{\bf 1}_{\{t_{2r}<t_{2r+1}\}},
\end{equation}

\vspace{-2mm}
$$
s_q(t_1,t_3,\ldots,t_{2r-3},t_{2r},t_{2r+1})=
$$
\begin{equation}
\label{july30019}
~~~~~~~~=\sum\limits_{l_1, \ldots, l_{r+1}=0}^q
C_{l_{r+1}\ldots l_1}\bar \phi_{l_1}(t_{1})\bar \phi_{l_2}(t_{3})\ldots \bar \phi_{l_{r-1}}(t_{2r-3})
\bar \phi_{l_r}(t_{2r})\bar \phi_{l_{r+1}}(t_{2r+1}),
\end{equation}
where $\left\{\bar \phi_j(x)\right\}_{j=0}^{\infty}$ 
is a complete orthonormal system of Legendre polynomials in $L_2([t,T]),$
$C_{l_{r+1}\ldots l_1}$ are Fourier--Legendre coefficients for the function
(\ref{july50000}), and
$\bar \psi_{1}(\tau),\bar \psi_{3}(\tau),\ldots ,\bar\psi_{2r-3}(\tau),
\bar \psi_{2r}(\tau), \bar \psi_{2r+1}(\tau)\in L_2([t,T]).$ 
Then we have
$$
\lim\limits_{q\to\infty}\left\Vert s_q - g\right\Vert^2_{L_2([t,T]^{r+1})}=0.
$$

From (\ref{july30018}) we obtain (the sum on the right-hand side of (\ref{july30019}) is finite)
$$
\lim\limits_{p\to\infty}\sum_{j_{2r+2}, j_{2r}, \ldots, j_2=0}^{p}
~\int\limits_{[t,T]^{2r+2}}
{\bf 1}_{\{t_1<t_2<\ldots <t_{2r}\}}{\bf 1}_{\{t_{2r+1}<t_{2r+2}\}}
s_q(t_1,t_3,\ldots,t_{2r-3},t_{2r},t_{2r+1})\times
$$
$$
\times
h(t_2,t_4,\ldots,t_{2r-2},t_{2r-1},t_{2r+2})
\times
$$

\vspace{-4mm}
$$
\times
\prod\limits_{d=1}^{r+1}\phi_{j_{2d}}(t_{2d-1})\phi_{j_{2d}}(t_{2d})
dt_1 dt_2\ldots dt_{2r-1}dt_{2r}dt_{2r+1}dt_{2r+2}=
$$

\vspace{1mm}
$$
=\frac{1}{2^{r+1}}
\int\limits_{[t,T]^{r+1}}
{\bf 1}_{\{t_2<t_4<\ldots <t_{2r}\}}s_q(t_2,t_4,\ldots,t_{2r-2},t_{2r},t_{2r+2})\times
$$

\vspace{-7mm}
\begin{equation}
\label{july30020}
~~~~~~~~~~~\times
h(t_2,t_4,\ldots,t_{2r-2},t_{2r},t_{2r+2})dt_2 dt_4\ldots dt_{2r-2}dt_{2r}dt_{2r+2}.
\end{equation}

\vspace{2mm}

The right-hand side of the equality (\ref{july30020}) defines
(as a scalar product of $s_q(t_2,t_4,\ldots,t_{2r-2},t_{2r},t_{2r+2})$
and 
$\frac{1}{2^{r+1}}{\bf 1}_{\{t_2<t_4<\ldots <t_{2r}\}}
h(t_2,t_4,\ldots,t_{2r-2},t_{2r},t_{2r+2})$
in the space $L_2([t, T]^{r+1})$) a linear bounded (and therefore continuous)
functional in $L_2([t, T]^{r+1}).$
The mentioned functional is given by the function 
$\frac{1}{2^{r+1}}{\bf 1}_{\{t_2<t_4<\ldots <t_{2r}\}}h(t_2,t_4,\ldots,t_{2r-2},t_{2r},t_{2r+2})$.

Note that the equality (\ref{july30020}) will also remain true if
$s_q$ is replaced 
by $\bar s_q$ ($\bar s_q$ is the partial sum of the Fourier--Legendre series
of any function from $L_2([t, T]^{r+1})$), i.e. 
the modified equality (\ref{july30020}) is true 
on a dense subset in $L_2([t,T]^{r+1}).$
On the left-hand side of (\ref{july30020}) (by virtue of the equality (\ref{july30020}))
there is a linear continuous functional on a dense subset in 
$L_2([t,T]^{r+1}).$ This functional can be uniquely extended 
to a linear continuous functional in $L_2([t, T]^{r+1})$
(see \cite{reed}, Theorem~I.7, P.~9).
Thus, we have the equality of two 
linear continuous functionals in $L_2([t, T]^{r+1}).$
Let us implement the passage to the limit $\lim\limits_{q\to\infty}$
in the mentioned equality if instead of $\bar s_q$
we choose $s_q$ of the form (\ref{july30019}) (i.e. passage to the limit $\lim\limits_{q\to\infty}$
in (\ref{july30020}))
$$
\lim\limits_{p\to\infty}\sum_{j_{2r+2}, j_{2r}, \ldots, j_2=0}^{p}
~\int\limits_{[t,T]^{2r+2}}
{\bf 1}_{\{t_1<t_2<\ldots <t_{2r}\}}{\bf 1}_{\{t_{2r+1}<t_{2r+2}\}}
g(t_1,t_3,\ldots,t_{2r-3},t_{2r},t_{2r+1})\times
$$
$$
\times
h(t_2,t_4,\ldots,t_{2r-2},t_{2r-1},t_{2r+2})
\times
$$

\vspace{-4mm}
$$
\times
\prod\limits_{d=1}^{r+1}\phi_{j_{2d}}(t_{2d-1})\phi_{j_{2d}}(t_{2d})
dt_1 dt_2\ldots dt_{2r-1}dt_{2r}dt_{2r+1}dt_{2r+2}=
$$

\vspace{2mm}
$$
=\frac{1}{2^{r+1}}
\int\limits_{[t,T]^{r+1}}
{\bf 1}_{\{t_2<t_4<\ldots <t_{2r}\}}g(t_2,t_4,\ldots,t_{2r-2},t_{2r},t_{2r+2})\times
$$

\vspace{-4mm}
\begin{equation}
\label{july30022}
~~~~~~~~~~~\times
h(t_2,t_4,\ldots,t_{2r-2},t_{2r},t_{2r+2})dt_2 dt_4\ldots dt_{2r-2}dt_{2r}dt_{2r+2},
\end{equation}

\vspace{4mm}
\noindent
where $\bar \psi_{1}(\tau),\bar \psi_{3}(\tau),\ldots ,\bar\psi_{2r-3}(\tau)
\bar \psi_{2r}(\tau), \bar \psi_{2r+1}(\tau)\in L_2([t,T]).$

It is easy to see that the equality (\ref{july30022}) (up to notations)
is the equality (\ref{july30016}) in which $r$ is replaced by $r+1.$
So, we proved the equality (\ref{july30016}) by induction.

Note that the series on the left-hand side of 
(\ref{july30016}) converges absolutly since
its sum does not depend 
on permutations of basis functions
(here the basis in $L_2([t,T]^{r})$ is
$\left\{\phi_{j_1}(x_1)\ldots \phi_{j_r}(x_r)\right\}_{j_1,\ldots,j_r=0}^{\infty}$).

Further, let us show that
$$
\lim\limits_{p\to\infty}
\sum\limits_{j_{g_1}, j_{g_3},\ldots ,j_{g_{2r-1}}=0}^p
C_{j_k\ldots j_1}\biggl|_{j_{g_1}=j_{g_2},\ldots, j_{g_{2r-1}}=j_{g_{2r}}}=
$$
\begin{equation}
\label{july90000}
=\frac{1}{2^r} \prod\limits_{l=1}^r {\bf 1}_{\{g_{2l}=g_{2l-1}+1\}}
C_{j_k \ldots j_1}\biggl|_{(j_{g_2} j_{g_1})\curvearrowright (\cdot)
\ldots (j_{g_{2r}} j_{g_{2r-1}})\curvearrowright (\cdot),
j_{g_{{}_{1}}}=~j_{g_{{}_{2}}},\ldots, j_{g_{{}_{2r-1}}}=~j_{g_{{}_{2r}}}
}\biggr.
\end{equation}

\vspace{4mm}
\noindent
for all possible $g_1,g_2,\ldots,g_{2r-1},g_{2r}$ (see {\rm (\ref{leto5007after})),
where $k=2r$ $(r=2,3,\ldots),$ $C_{j_k\ldots j_1}$ is defined by (\ref{july15030}),
another notations are the same as in Theorem~2.49.

The case
$$
\prod\limits_{l=1}^r {\bf 1}_{\{g_{2l}=g_{2l-1}+1\}}=1
$$

\noindent
corresponds to (\ref{july30016}).

Thus, it remains to prove that
\begin{equation}
\label{july80000}
\lim\limits_{p\to\infty}
\sum\limits_{j_{g_1}, j_{g_3},\ldots ,j_{g_{2r-1}}=0}^p
C_{j_k\ldots j_1}\biggl|_{j_{g_1}=j_{g_2},\ldots, j_{g_{2r-1}}=j_{g_{2r}}}=0
\end{equation}

\noindent
for the case
$$
\prod\limits_{l=1}^r {\bf 1}_{\{g_{2l}=g_{2l-1}+1\}}=0.
$$

Below we consider two examples that clearly explain 
the algorithm for the proof of equality (\ref{july80000}).
After this we will formulate the algorithm.

First, let us prove that
$$
\lim\limits_{p\to\infty}\sum_{j_1,j_3,j_4=0}^{p}
C_{j_3 j_4 j_4 j_3 j_1 j_1}=
$$
$$
=\lim\limits_{p\to\infty}\sum_{j_1,j_3,j_4=0}^{p}
\int\limits_t^T
\psi_6(t_6)\phi_{j_3}(t_6)\int\limits_t^{t_6}
\psi_5(t_5)\phi_{j_4}(t_5)\int\limits_t^{t_5}
\psi_4(t_4)\phi_{j_4}(t_4)
\int\limits_{t}^{t_4}
\psi_3(t_3)\phi_{j_3}(t_3)\times\Biggr.
$$
\begin{equation}
\label{july80013}
~~~~~~~~\times
\int\limits_t^{t_3}\psi_2(t_2)\phi_{j_1}(t_2)\int\limits_t^{t_2}
\psi_1(t_1)\phi_{j_1}(t_1)dt_1 dt_2 dt_3 dt_4 dt_5 dt_6=0,
\end{equation}

\noindent
where $\left\{\phi_j(x)\right\}_{j=0}^{\infty}$
is an arbitrary complete orthonormal system of 
functions in the space $L_2([t,T])$ and
$\psi_1(\tau),\ldots ,\psi_{6}(\tau)\in L_2([t, T]).$

{\bf Step~1.}\ Using (\ref{july30016}) ($r=1$) and generalized Parseval's equality, we obtain
$$
\lim\limits_{p\to\infty}\sum_{j_1,j_3,j_4=0}^{p}
\int\limits_t^T
\psi_6(t_6)\phi_{j_3}(t_6)\int\limits_t^{T}
\psi_5(t_5)\phi_{j_4}(t_5)\int\limits_t^{t_5}
\psi_4(t_4)\phi_{j_4}(t_4)
\int\limits_{t}^{T}
\psi_3(t_3)\phi_{j_3}(t_3)\times\Biggr.
$$
\begin{equation}
\label{july80100}
~~~~~~~~~~~~\times
\int\limits_t^{T}\psi_2(t_2)\phi_{j_1}(t_2)\int\limits_t^{t_2}
\psi_1(t_1)\phi_{j_1}(t_1)dt_1 dt_2 dt_3 dt_4 dt_5 dt_6=
\end{equation}
$$
=\lim\limits_{p\to\infty}\sum_{j_3=0}^{p}
\int\limits_t^T
\psi_6(t_6)\phi_{j_3}(t_6)dt_6 \int\limits_{t}^{T}
\psi_3(t_3)\phi_{j_3}(t_3)dt_3\times
$$
$$
\times
\lim\limits_{p\to\infty}\sum_{j_4=0}^{p}
\int\limits_t^{T}
\psi_5(t_5)\phi_{j_4}(t_5)\int\limits_t^{t_5}
\psi_4(t_4)\phi_{j_4}(t_4)dt_4 dt_5\times
$$
$$
\times
\lim\limits_{p\to\infty}\sum_{j_1=0}^{p}
\int\limits_t^{T}\psi_2(t_2)\phi_{j_1}(t_2)\int\limits_t^{t_2}
\psi_1(t_1)\phi_{j_1}(t_1)dt_1 dt_2=
$$
\begin{equation}
\label{july80001}
~~~~~~~~~~=
\int\limits_t^T
\psi_6(t_6)\psi_3(t_6)dt_6
\cdot \frac{1}{2}
\int\limits_t^{T}
\psi_5(t_4)\psi_4(t_4)dt_4
\cdot \frac{1}{2}
\int\limits_t^{T}\psi_2(t_2)\psi_1(t_2)dt_2.
\end{equation}

\vspace{2mm}

Let us rewrite (\ref{july80001}) in the form
$$
\sum_{j_1,j_3,j_4=0}^{\infty}
\int\limits_t^T
\psi_6(t_6)\phi_{j_3}(t_6)\int\limits_t^{T}
\psi_5(t_5)\phi_{j_4}(t_5)\int\limits_t^{t_5}
\psi_4(t_4)\phi_{j_4}(t_4)
\int\limits_{t}^{T}
\psi_3(t_3)\phi_{j_3}(t_3)\times\Biggr.
$$
$$
\times
\int\limits_t^{T}\psi_2(t_2)\phi_{j_1}(t_2)\int\limits_t^{t_2}
\psi_1(t_1)\phi_{j_1}(t_1)dt_1 dt_2 dt_3 dt_4 dt_5 dt_6=
$$
\begin{equation}
\label{july80002}
~~~~~~~~~~~=
\frac{1}{4}\int\limits_t^T
\psi_6(t_6)\psi_3(t_6)
\int\limits_t^{T}
\psi_5(t_4)\psi_4(t_4)
\int\limits_t^{T}\psi_2(t_2)\psi_1(t_2)dt_2 dt_4 dt_6.
\end{equation}

{\bf Step~2.}\ Suppose that $\psi_2(\tau),\psi_3(\tau),\psi_4(\tau)$ are Legendre
polynomials of finite degrees.
Denote
\begin{equation}
\label{july80003}
s_q(t_2,t_3,t_4)=\sum\limits_{l_1,l_2.l_3=0}^q
C_{l_3 l_2 l_1}\bar \phi_{l_1}(t_2)\bar \phi_{l_2}(t_3) \bar \phi_{l_3}(t_4),
\end{equation}

\noindent 
where $\left\{\bar \phi_j(x)\right\}_{j=0}^{\infty}$ 
is a complete orthonormal system of Legendre polynomials in $L_2([t,T])$
and
$C_{l_3 l_2 l_1}$ are Fourier--Legendre coefficients for the function
$g(t_2,t_3,t_4)=\bar \psi_2(t_2)\bar \psi_3(t_3)\bar \psi_4(t_4){\bf 1}_{\{t_2<t_3\}}$ 
($\bar \psi_2(\tau), \bar \psi_3(\tau), \bar \psi_4(\tau)\in L_2([t,T])),$ i.e.
$\lim\limits_{q\to\infty}\left\Vert s_q - g\right\Vert^2_{L_2([t,T]^3)}=0.$

From (\ref{july80002}) we obtain (the sum on the right-hand side of (\ref{july80003}) is finite)
$$
\sum_{j_1,j_3,j_4=0}^{\infty}
\int\limits_{[t,T]^6}{\bf 1}_{\{t_1<t_2\}}{\bf 1}_{\{t_4<t_5\}}
s_q(t_2,t_3,t_4)\psi_6(t_6)\psi_5(t_5)\psi_1(t_1)
\phi_{j_3}(t_6)
\phi_{j_3}(t_3)
\phi_{j_4}(t_5)
\times\Biggr.
$$
$$
\times
\phi_{j_4}(t_4)\phi_{j_1}(t_2)
\phi_{j_1}(t_1)dt_1 dt_2 dt_3 dt_4 dt_5 dt_6=
$$

\vspace{-5mm}
\begin{equation}
\label{july80004}
~~~~~~~~~~~=
\frac{1}{4}\int\limits_{[t,T]^3}
s_q(t_2,t_6,t_4)\psi_6(t_6)
\psi_5(t_4)\psi_1(t_2)dt_2 dt_4 dt_6.
\end{equation}

Note that the equality (\ref{july80004}) remains true
when $s_q$ is a partial sum of the Fourier--Legendre series
of any function from $L_2([t,T]^3),$ i.e. the equality holds
on a dense subset in $L_2([t,T]^3).$

The right-hand side of (\ref{july80004}) defines
(as a scalar product of $s_q(t_2,t_6,t_4)$ and $\frac{1}{4}\psi_6(t_6)\psi_5(t_4)\psi_1(t_2)$
in $L_2([t, T]^3)$) a linear bounded (and therefore continuous)
functional in $L_2([t, T]^3),$
which is given by the function $\frac{1}{4}\psi_6(t_6)\psi_5(t_4)\psi_1(t_2)$.
On the left-hand side of (\ref{july80004}) (by virtue of the equality (\ref{july80004}))
there is a linear continuous functional on a dense subset in 
$L_2([t,T]^3).$ This functional can be uniquely extended 
to a linear continuous functional in $L_2([t, T]^3)$
(see \cite{reed}, Theorem~I.7, P.~9).

Let us implement the passage to the limit $\lim\limits_{q\to\infty}$
in (\ref{july80004}) (at that we suppose that $s_q$ is defined by (\ref{july80003}))
$$
\sum_{j_1,j_3,j_4=0}^{\infty}
\int\limits_{[t,T]^6}{\bf 1}_{\{t_1<t_2<t_3\}}{\bf 1}_{\{t_4<t_5\}}
\psi_6(t_6)\psi_5(t_5)\bar \psi_4(t_4)\bar \psi_3(t_3)\bar \psi_2(t_2) \psi_1(t_1)
\phi_{j_3}(t_6)
\phi_{j_3}(t_3)
\times\Biggr.
$$
$$
\times
\phi_{j_4}(t_5)\phi_{j_4}(t_4)\phi_{j_1}(t_2)
\phi_{j_1}(t_1)dt_1 dt_2 dt_3 dt_4 dt_5 dt_6=
$$

\vspace{-5mm}
\begin{equation}
\label{july80005}
~~~~~~~~=
\frac{1}{4}\int\limits_{[t,T]^3}
{\bf 1}_{\{t_2<t_6\}}
\psi_6(t_6)\bar \psi_3(t_6)
\psi_5(t_4)\bar \psi_4(t_4)
\bar \psi_2(t_2)
\psi_1(t_2)dt_2 dt_4 dt_6.
\end{equation}

Rewrite the equality (\ref{july80005}) in the form
$$
\sum_{j_1,j_3,j_4=0}^{\infty}
\int\limits_{[t,T]^6}{\bf 1}_{\{t_1<t_2<t_3\}}{\bf 1}_{\{t_4<t_5\}}
\psi_6(t_6)\psi_5(t_5)\psi_4(t_4)\psi_3(t_3)\psi_2(t_2) \psi_1(t_1)
\phi_{j_3}(t_6)
\phi_{j_3}(t_3)
\times\Biggr.
$$
$$
\times
\phi_{j_4}(t_5)\phi_{j_4}(t_4)\phi_{j_1}(t_2)
\phi_{j_1}(t_1)dt_1 dt_2 dt_3 dt_4 dt_5 dt_6=
$$

\vspace{-5mm}
\begin{equation}
\label{july80006}
~~~~~~~~=
\frac{1}{4}\int\limits_{[t,T]^3}
{\bf 1}_{\{t_2<t_6\}}
\psi_6(t_6)\psi_3(t_6)
\psi_5(t_4)\psi_4(t_4)
\psi_2(t_2)
\psi_1(t_2)dt_2 dt_4 dt_6,
\end{equation}

\noindent
where $\psi_1(\tau),\ldots,\psi_6(\tau)\in L_2([t,T]).$

{\bf Step~3.}\ Suppose that $\psi_3(\tau),\psi_4(\tau),\psi_1(\tau)$ are 
Legendre polynomials of finite degrees.
Denote
\begin{equation}
\label{july80007}
s_q(t_3,t_4,t_1)=\sum\limits_{l_1,l_2.l_3=0}^q
C_{l_3 l_2 l_1}\bar \phi_{l_1}(t_3)\bar \phi_{l_2}(t_4) \bar \phi_{l_3}(t_1),
\end{equation}

\noindent 
where $\left\{\bar \phi_j(x)\right\}_{j=0}^{\infty}$ as in (\ref{july80003})
and
$C_{l_3 l_2 l_1}$ are Fourier--Legendre coefficients for the function
$g(t_3,t_4,t_1)=\bar \psi_3(t_3)\bar \psi_4(t_4)\bar \psi_1(t_1){\bf 1}_{\{t_3<t_4\}}$ 
($\bar \psi_3(\tau), \bar \psi_4(\tau), \bar \psi_1(\tau)\in L_2([t,T])),$ i.e.
$\lim\limits_{q\to\infty}\left\Vert s_q - g\right\Vert^2_{L_2([t,T]^3)}=0.$

From (\ref{july80006}) we obtain (the sum on the right-hand side of (\ref{july80007}) is finite)
$$
\sum_{j_1,j_3,j_4=0}^{\infty}
\int\limits_{[t,T]^6}{\bf 1}_{\{t_1<t_2<t_3\}}{\bf 1}_{\{t_4<t_5\}}
s_q(t_3,t_4,t_1)
\psi_6(t_6)\psi_5(t_5)\psi_2(t_2)
\phi_{j_3}(t_6)
\phi_{j_3}(t_3)
\times\Biggr.
$$
$$
\times
\phi_{j_4}(t_5)\phi_{j_4}(t_4)\phi_{j_1}(t_2)
\phi_{j_1}(t_1)dt_1 dt_2 dt_3 dt_4 dt_5 dt_6=
$$

\vspace{-5mm}
\begin{equation}
\label{july80008}
~~~~~~~~=
\frac{1}{4}\int\limits_{[t,T]^3}
{\bf 1}_{\{t_2<t_6\}}
s_q(t_6,t_4,t_2)
\psi_6(t_6)
\psi_5(t_4)
\psi_2(t_2)
dt_2 dt_4 dt_6.
\end{equation}

Note that the equality (\ref{july80008}) remains true
when $s_q$ is a partial sum of the Fourier--Legendre series
of any function from $L_2([t,T]^3),$ i.e. the equality holds
on a dense subset in $L_2([t,T]^3).$

The right-hand side of (\ref{july80008}) defines
(as a scalar product of $s_q(t_6,t_4,t_2)$ and $\frac{1}{4}{\bf 1}_{\{t_2<t_6\}}\psi_6(t_6)\psi_5(t_4)\psi_2(t_2)$
in $L_2([t, T]^3)$) a linear bounded (and therefore continuous)
functional in $L_2([t, T]^3),$
which is given by the function $\frac{1}{4}{\bf 1}_{\{t_2<t_6\}}\psi_6(t_6)\psi_5(t_4)\psi_2(t_2)$.
On the left-hand side of (\ref{july80008}) (by virtue of the equality (\ref{july80008}))
there is a linear continuous functional on a dense subset in 
$L_2([t,T]^3).$ This functional can be uniquely extended 
to a linear continuous functional in $L_2([t, T]^3)$
(see \cite{reed}, Theorem~I.7, P.~9).

Let us implement the passage to the limit $\lim\limits_{q\to\infty}$
in (\ref{july80008}) (at that we suppose that $s_q$ is defined by (\ref{july80007}))
$$
\sum_{j_1,j_3,j_4=0}^{\infty}
\int\limits_{[t,T]^6}{\bf 1}_{\{t_1<t_2<t_3<t_4<t_5\}}
\psi_6(t_6)\psi_5(t_5)\bar \psi_4(t_4)
\bar \psi_3(t_3)\psi_2(t_2)\bar \psi_1(t_1)
\phi_{j_3}(t_6)
\phi_{j_3}(t_3)
\times\Biggr.
$$
$$
\times
\phi_{j_4}(t_5)\phi_{j_4}(t_4)\phi_{j_1}(t_2)
\phi_{j_1}(t_1)dt_1 dt_2 dt_3 dt_4 dt_5 dt_6=
$$

\vspace{-5mm}
\begin{equation}
\label{july80009}
~~=
\frac{1}{4}\int\limits_{[t,T]^3}
{\bf 1}_{\{t_2<t_6\}}{\bf 1}_{\{t_6<t_4\}}
\psi_6(t_6)\bar \psi_3(t_6)
\psi_5(t_4)\bar \psi_4(t_4)
\psi_2(t_2)\bar \psi_1(t_2)
dt_2 dt_4 dt_6.
\end{equation}

Rewrite (\ref{july80009}) in the form
$$
\sum_{j_1,j_3,j_4=0}^{\infty}
\int\limits_{[t,T]^6}{\bf 1}_{\{t_1<t_2<t_3<t_4<t_5\}}
\psi_6(t_6)\psi_5(t_5)\psi_4(t_4)
\psi_3(t_3)\psi_2(t_2)\psi_1(t_1)
\phi_{j_3}(t_6)
\phi_{j_3}(t_3)
\times\Biggr.
$$
$$
\times
\phi_{j_4}(t_5)\phi_{j_4}(t_4)\phi_{j_1}(t_2)
\phi_{j_1}(t_1)dt_1 dt_2 dt_3 dt_4 dt_5 dt_6=
$$

\vspace{-5mm}
\begin{equation}
\label{july80009x}
~~=
\frac{1}{4}\int\limits_{[t,T]^3}
{\bf 1}_{\{t_2<t_6\}}{\bf 1}_{\{t_6<t_4\}}
\psi_6(t_6)\psi_3(t_6)
\psi_5(t_4)\psi_4(t_4)
\psi_2(t_2)\psi_1(t_2)
dt_2 dt_4 dt_6,
\end{equation}

\noindent
where $\psi_1(\tau),\ldots,\psi_6(\tau)\in L_2([t,T]).$

{\bf Step~4.}\ Suppose that $\psi_5(\tau),\psi_6(\tau),\psi_2(\tau)$ are 
Legendre polynomials of finite degrees.
Denote
\begin{equation}
\label{july80010}
s_q(t_5,t_6,t_2)=\sum\limits_{l_1,l_2.l_3=0}^q
C_{l_3 l_2 l_1}\bar \phi_{l_1}(t_5)\bar \phi_{l_2}(t_6) \bar \phi_{l_3}(t_2),
\end{equation}

\noindent 
where $\left\{\bar \phi_j(x)\right\}_{j=0}^{\infty}$ as in (\ref{july80003})
and
$C_{l_3 l_2 l_1}$ are Fourier--Legendre coefficients for the function
$g(t_5,t_6,t_2)=\bar \psi_5(t_5)\bar \psi_6(t_6)\bar \psi_2(t_2){\bf 1}_{\{t_5<t_6\}}$ 
($\bar \psi_5(\tau), \bar \psi_6(\tau), \bar \psi_2(\tau)\in L_2([t,T])),$ i.e.
$\lim\limits_{q\to\infty}\left\Vert s_q - g\right\Vert^2_{L_2([t,T]^3)}=0.$

From (\ref{july80009x}) we obtain (the sum on the right-hand side of (\ref{july80010}) is finite)
$$
\sum_{j_1,j_3,j_4=0}^{\infty}
\int\limits_{[t,T]^6}{\bf 1}_{\{t_1<t_2<t_3<t_4<t_5\}}
s_q(t_5,t_6,t_2)
\psi_4(t_4)
\psi_3(t_3)\psi_1(t_1)
\phi_{j_3}(t_6)
\phi_{j_3}(t_3)
\times\Biggr.
$$
$$
\times
\phi_{j_4}(t_5)\phi_{j_4}(t_4)\phi_{j_1}(t_2)
\phi_{j_1}(t_1)dt_1 dt_2 dt_3 dt_4 dt_5 dt_6=
$$

\vspace{-5mm}
\begin{equation}
\label{july80011}
~~~~~~~~=
\frac{1}{4}\int\limits_{[t,T]^3}
{\bf 1}_{\{t_2<t_6\}}{\bf 1}_{\{t_6<t_4\}}
s_q(t_4,t_6,t_2)
\psi_3(t_6)
\psi_4(t_4)
\psi_1(t_2)
dt_2 dt_4 dt_6.
\end{equation}

Note that the equality (\ref{july80011}) remains true
when $s_q$ is a partial sum of the Fourier--Legendre series
of any function from $L_2([t,T]^3),$ i.e. the equality holds
on a dense subset in $L_2([t,T]^3).$

The right-hand side of (\ref{july80011}) defines
(as a scalar product of $s_q(t_4,t_6,t_2)$ and $\frac{1}{4}{\bf 1}_{\{t_2<t_6\}}{\bf 1}_{\{t_6<t_4\}}
\psi_3(t_6)\psi_4(t_4)\psi_1(t_2)$
in $L_2([t, T]^3)$) a linear bounded (and therefore continuous)
functional in $L_2([t, T]^3),$
which is given by the function $\frac{1}{4}{\bf 1}_{\{t_2<t_6\}}{\bf 1}_{\{t_6<t_4\}}\psi_3(t_6)\psi_4(t_4)\psi_1(t_2)$.
On the left-hand side of (\ref{july80011}) (by virtue of the equality (\ref{july80011}))
there is a linear continuous functional on a dense subset in 
$L_2([t,T]^3).$ This functional can be uniquely extended 
to a linear continuous functional in $L_2([t, T]^3)$
(see \cite{reed}, Theorem~I.7, P.~9).

Let us implement the passage to the limit $\lim\limits_{q\to\infty}$
in (\ref{july80011}) (at that we suppose that $s_q$ is defined by (\ref{july80010}))
$$
\sum_{j_1,j_3,j_4=0}^{\infty}
\int\limits_{[t,T]^6}\hspace{-2mm}{\bf 1}_{\{t_1<t_2<t_3<t_4<t_5<t_6\}}
\bar \psi_6(t_6)\bar \psi_5(t_5)\psi_4(t_4)
\psi_3(t_3)\bar \psi_2(t_2)\psi_1(t_1)
\phi_{j_3}(t_6)
\phi_{j_3}(t_3)
\times\Biggr.
$$
$$
\times
\phi_{j_4}(t_5)\phi_{j_4}(t_4)\phi_{j_1}(t_2)
\phi_{j_1}(t_1)dt_1 dt_2 dt_3 dt_4 dt_5 dt_6=
$$

\vspace{-5mm}
\begin{equation}
\label{july80012}
=
\frac{1}{4}\int\limits_{[t,T]^3}
{\bf 1}_{\{t_2<t_6\}}{\bf 1}_{\{t_6<t_4\}}{\bf 1}_{\{t_4<t_6\}}
\bar \psi_6(t_6)\psi_3(t_6)
\bar \psi_5(t_4)\psi_4(t_4)
\bar \psi_2(t_2)\psi_1(t_2)
dt_2 dt_4 dt_6=0.
\end{equation}

It is obvious that the equality (\ref{july80012}) (up to notations)
is (\ref{july80013}). The equality (\ref{july80013}) is proved.

As a second example, we will prove the equality (\ref{july15}).
In this case, we will use the same approach as in the proof
of equality (\ref{july80013}). Thus, we prove that
\begin{equation}
\label{july80015}
\lim\limits_{p\to\infty}
\sum\limits_{j_1, j_2=0}^{p}
C_{j_2 j_1 j_2 j_1}=0.
\end{equation}

\vspace{2mm}

{\bf Step~1.}\ Using generalized Parseval's equality, we obtain
$$
\lim\limits_{p\to\infty}\sum_{j_1,j_2=0}^{p}
\int\limits_t^T
\psi_4(t_4)\phi_{j_2}(t_4)\int\limits_t^{T}
\psi_3(t_3)\phi_{j_1}(t_3)\int\limits_t^{T}
\psi_2(t_2)\phi_{j_2}(t_2)
\int\limits_{t}^{T}
\psi_1(t_1)\phi_{j_1}(t_1)\times
$$
\begin{equation}
\label{july80101}
\times
dt_1 dt_2 dt_3 dt_4=
\end{equation}
$$
=\lim\limits_{p\to\infty}\sum_{j_2=0}^{p}
\int\limits_t^T
\psi_4(t_4)\phi_{j_2}(t_4)dt_4
\int\limits_t^{T}
\psi_2(t_2)\phi_{j_2}(t_2)dt_2\times
$$
$$
\times
\lim\limits_{p\to\infty}\sum_{j_1=0}^{p}
\int\limits_t^{T}
\psi_3(t_3)\phi_{j_1}(t_3)dt_3
\int\limits_{t}^{T}
\psi_1(t_1)\phi_{j_1}(t_1)dt_1=
$$
\begin{equation}
\label{july80016}
=
\int\limits_t^T
\psi_4(t_4)\psi_2(t_4)dt_4 
\int\limits_t^{T}
\psi_3(t_3)\psi_1(t_3)dt_3.
\end{equation}

Rewrite the equality (\ref{july80016}) in the form
$$
\sum_{j_1,j_2=0}^{\infty}
\int\limits_{[t,T]^4}
\psi_4(t_4)\psi_3(t_3)\psi_2(t_2)\psi_1(t_1)
\phi_{j_2}(t_4)
\phi_{j_1}(t_3)
\phi_{j_2}(t_2)
\phi_{j_1}(t_1)
dt_1 dt_2 dt_3 dt_4=
$$
\begin{equation}
\label{july80017}
=\int\limits_{[t,T]^2}
\psi_4(t_4)\psi_2(t_4)
\psi_3(t_2)\psi_1(t_2)dt_2 dt_4.
\end{equation}

\vspace{2mm}

{\bf Step~2.}\ Suppose that $\psi_1(\tau), \psi_2(\tau)$ are Legendre polynomials of finite degrees.
Denote
$$
s_q(t_1,t_2)=\sum\limits_{l_1,l_2=0}^q
C_{l_2 l_1}\bar \phi_{l_1}(t_1)\bar \phi_{l_2}(t_2),
$$

\noindent 
where $\left\{\bar \phi_j(x)\right\}_{j=0}^{\infty}$ as in (\ref{july80003}),
$C_{l_2 l_1}$ are Fourier--Legendre coefficients for the function
$g(t_1,t_2)=\bar \psi_1(t_1)\bar \psi_2(t_2){\bf 1}_{\{t_1<t_2\}}$ 
($\bar \psi_1(\tau), \bar \psi_2(\tau)\in L_2([t,T])).$

From (\ref{july80017}) we obtain 
$$
\sum_{j_1,j_2=0}^{\infty}
\int\limits_{[t,T]^4}
s_q(t_1,t_2)\psi_4(t_4)\psi_3(t_3)
\phi_{j_2}(t_4)
\phi_{j_1}(t_3)
\phi_{j_2}(t_2)
\phi_{j_1}(t_1)
dt_1 dt_2 dt_3 dt_4=
$$
\begin{equation}
\label{july80019}
=\int\limits_{[t,T]^2}
s_q(t_2,t_4)\psi_4(t_4)
\psi_3(t_2)dt_2 dt_4.
\end{equation}

The left-hand and right-hand sides of (\ref{july80019}) define
linear continuous functionals in $L_2([t, T]^2)$ 
(see explanation earlier in this section).
Let us implement the passage to the limit $\lim\limits_{q\to\infty}$
in (\ref{july80019})
$$
\sum_{j_1,j_2=0}^{\infty}
\int\limits_{[t,T]^4}
\hspace{-2.2mm}
{\bf 1}_{\{t_1<t_2\}}\psi_4(t_4)\psi_3(t_3)\bar \psi_2(t_2)\bar \psi_1(t_1)
\phi_{j_2}(t_4)
\phi_{j_1}(t_3)
\phi_{j_2}(t_2)
\phi_{j_1}(t_1)
dt_1 dt_2 dt_3 dt_4\hspace{-0.8mm}=
$$
\begin{equation}
\label{july80020}
=\int\limits_{[t,T]^2}
{\bf 1}_{\{t_2<t_4\}}\psi_4(t_4)\bar \psi_2(t_4)
\psi_3(t_2)\bar \psi_1(t_2)dt_2 dt_4.
\end{equation}

\vspace{2mm}

Rewrite the equality (\ref{july80020}) in the form
$$
\sum_{j_1,j_2=0}^{\infty}
\int\limits_{[t,T]^4}
\hspace{-2.2mm}{\bf 1}_{\{t_1<t_2\}}\psi_4(t_4)\psi_3(t_3)\psi_2(t_2)\psi_1(t_1)
\phi_{j_2}(t_4)
\phi_{j_1}(t_3)
\phi_{j_2}(t_2)
\phi_{j_1}(t_1)
dt_1 dt_2 dt_3 dt_4\hspace{-0.8mm}=
$$
\begin{equation}
\label{july80021}
=\int\limits_{[t,T]^2}
{\bf 1}_{\{t_2<t_4\}}\psi_4(t_4)\psi_2(t_4)
\psi_3(t_2)\psi_1(t_2)dt_2 dt_4,
\end{equation}

\noindent
where $\psi_1(\tau),\ldots,\psi_4(\tau)\in L_2([t, T]).$

\vspace{2mm}

{\bf Step~3.}\ Suppose that $\psi_2(\tau), \psi_3(\tau)$ are Legendre polynomials of finite degrees.
Denote
$$
s_q(t_2,t_3)=\sum\limits_{l_1,l_2=0}^q
C_{l_2 l_1}\bar \phi_{l_1}(t_2)\bar \phi_{l_2}(t_3),
$$

\noindent 
where $\left\{\bar \phi_j(x)\right\}_{j=0}^{\infty}$ as in (\ref{july80003}),
$C_{l_2 l_1}$ are Fourier--Legendre coefficients for the function
$g(t_2,t_3)=\bar \psi_2(t_2)\bar \psi_3(t_3){\bf 1}_{\{t_2<t_3\}}$ 
($\bar \psi_2(\tau), \bar \psi_3(\tau)\in L_2([t,T])).$

From (\ref{july80021}) we obtain 
$$
\sum_{j_1,j_2=0}^{\infty}
\int\limits_{[t,T]^4}
{\bf 1}_{\{t_1<t_2\}}
s_q(t_2,t_3)
\psi_4(t_4)\psi_1(t_1)
\phi_{j_2}(t_4)
\phi_{j_1}(t_3)
\phi_{j_2}(t_2)
\phi_{j_1}(t_1)
dt_1 dt_2 dt_3 dt_4=
$$
\begin{equation}
\label{july80022}
=\int\limits_{[t,T]^2}
{\bf 1}_{\{t_2<t_4\}}s_q(t_4,t_2)\psi_4(t_4)
\psi_1(t_2)dt_2 dt_4.
\end{equation}

The left-hand and right-hand sides of (\ref{july80022}) define
linear continuous functionals in $L_2([t, T]^2)$.
Let us implement the passage to the limit $\lim\limits_{q\to\infty}$
in (\ref{july80022})
$$
\sum_{j_1,j_2=0}^{\infty}
\int\limits_{[t,T]^4}
{\bf 1}_{\{t_1<t_2<t_3\}}
\psi_4(t_4)\bar \psi_3(t_3)\bar \psi_2(t_2) \psi_1(t_1)
\phi_{j_2}(t_4)
\phi_{j_1}(t_3)
\phi_{j_2}(t_2)
\phi_{j_1}(t_1)\times
$$
$$
\times
dt_1 dt_2 dt_3 dt_4=
$$

\vspace{-4mm}
\begin{equation}
\label{july80023}
~~~~~~=\int\limits_{[t,T]^2}
{\bf 1}_{\{t_2<t_4\}}{\bf 1}_{\{t_4<t_2\}}\psi_4(t_4)\bar \psi_2(t_4)
\bar \psi_3(t_2)\psi_1(t_2)dt_2 dt_4=0.
\end{equation}

Rewrite the equality (\ref{july80023}) in the form
$$
\sum_{j_1,j_2=0}^{\infty}
\int\limits_{[t,T]^4}
{\bf 1}_{\{t_1<t_2<t_3\}}
\psi_4(t_4)\psi_3(t_3)\psi_2(t_2) \psi_1(t_1)
\phi_{j_2}(t_4)
\phi_{j_1}(t_3)
\phi_{j_2}(t_2)
\phi_{j_1}(t_1)\times
$$
\begin{equation}
\label{july80024}
\times
dt_1 dt_2 dt_3 dt_4=0.
\end{equation}

\vspace{2mm}

{\bf Step~4.}\ Suppose that $\psi_3(\tau), \psi_4(\tau)$ are Legendre polynomials of finite degrees.
Denote
$$
s_q(t_3,t_4)=\sum\limits_{l_1,l_2=0}^q
C_{l_2 l_1}\bar \phi_{l_1}(t_3)\bar \phi_{l_2}(t_4),
$$

\noindent 
where $\left\{\bar \phi_j(x)\right\}_{j=0}^{\infty}$ as in (\ref{july80003}),
$C_{l_2 l_1}$ are Fourier--Legendre coefficients for the function
$g(t_3,t_4)=\bar \psi_3(t_3)\bar \psi_4(t_4){\bf 1}_{\{t_3<t_4\}}$ 
($\bar \psi_3(\tau), \bar \psi_4(\tau)\in L_2([t,T])).$

From (\ref{july80024}) we obtain 
$$
\sum_{j_1,j_2=0}^{\infty}
\int\limits_{[t,T]^4}
{\bf 1}_{\{t_1<t_2<t_3\}}
s_q(t_3,t_4)\psi_2(t_2) \psi_1(t_1)
\phi_{j_2}(t_4)
\phi_{j_1}(t_3)
\phi_{j_2}(t_2)
\phi_{j_1}(t_1)\times
$$
\begin{equation}
\label{july80025}
\times
dt_1 dt_2 dt_3 dt_4=0.
\end{equation}

\vspace{2mm}

The left-hand and right-hand sides of (\ref{july80025}) define
linear continuous functionals in $L_2([t, T]^2)$
(we interpret the right-hand side of (\ref{july80025})
as the zero functional in $L_2([t, T]^2)$).
Let us implement the passage to the limit $\lim\limits_{q\to\infty}$
in (\ref{july80025})
$$
\sum_{j_1,j_2=0}^{\infty}
\int\limits_{[t,T]^4}
{\bf 1}_{\{t_1<t_2<t_3<t_4\}}
\bar \psi_4(t_4) \bar \psi_3(t_3)\psi_2(t_2) \psi_1(t_1)
\phi_{j_2}(t_4)
\phi_{j_1}(t_3)
\phi_{j_2}(t_2)
\phi_{j_1}(t_1)\times
$$
\begin{equation}
\label{july80026}
\times
dt_1 dt_2 dt_3 dt_4=0.
\end{equation}

\vspace{2mm}

It is easy to see that the equality (\ref{july80026}) (up to notations)
is the equality (\ref{july15}).
The equality (\ref{july15}) is proved.

Let us formulate the ideas used when considering
the two above examples in the form of an algorithm.

{\bf Step~1.}\ Suppose $k=2r$ $(r=2,3,4,\ldots)$, where $r$ is the number of pairs
$\{g_1, g_2\}, \ldots, 
\{g_{2r-1}, g_{2r}\}$ (see (\ref{leto5007after})). Let us select blocks
in the multi-index $j_k\ldots j_1$ that
correspond to the fulfillment of the condition
$$
\prod\limits_{l=1}^{r_d} {\bf 1}_{\{g_{2l}=g_{2l-1}+1\}}=1,
$$
where $r_d$ is the number of pairs (see (\ref{leto5007after}))
in the block with number $d.$

{\bf Step~2.}\ Let us write the Volterra--type kernel (\ref{july7000}) in the form
\begin{equation}
\label{july80027}
~~~~~~~~~~~~K(t_1,\ldots,t_k)=
\psi_1(t_1)\ldots \psi_k(t_k){\bf 1}_{\{t_1<t_{2}\}}{\bf 1}_{\{t_2<t_{3}\}}\ldots
{\bf 1}_{\{t_{k-1}<t_{k}\}},
\end{equation}

\noindent
where $\psi_1(\tau),\ldots,\psi_k(\tau)\in L_2([t,T])$,
$t_1,\ldots,t_k\in [t, T],$ $k\ge 4.$

Let us save multipliers of the form ${\bf 1}_{\{t_n<t_{n+1}\}}$
in the expression (\ref{july80027}) that correspond 
to the above blocks. At that, we remove the remaining 
multipliers of the form 
${\bf 1}_{\{t_n<t_{n+1}\}}$ from the expression (\ref{july80027}).
As a result, we get a modified kernel $\bar K(t_1,\ldots,t_k)$.
Let us write an analogue of the left-hand side
of equality (\ref{july80000}) for the modified kernel $\bar K(t_1,\ldots,t_k)$
(see (\ref{july80100}) and (\ref{july80101}) as examples).
For definiteness, let us denote this expression by
$({}^{-})$.

{\bf Step~3.}\ Using generalized Parseval's equality and (\ref{july30016}), we represent
the expression $({}^{-})$ as an integral over the hypercube $[t, T]^r$
(see the right-hand sides of (\ref{july80002}) and (\ref{july80017}) as examples).
For definiteness, let us denote the obtained equality by
$(\bar K)$ ((\ref{july80002}) and (\ref{july80017}) are examples of $(\bar K)$).

{\bf Step~4.}\ Further, transformations and passages to the limit
in the equality $(\bar K)$ are performed iteratively 
in such a way as to restore the removed multipliers ${\bf 1}_{\{t_n<t_{n+1}\}}$
on the left-hand side of $(\bar K)$
(for more details, see the proof of formulas (\ref{july80013}), (\ref{july80015})).
As a result, we obtain the equality (\ref{july80000}).
More precisely, we can move from left to right
along a multi-index corresponding to the left-hand side of $(\bar K)$.
Let us assume that at the $n$-th step we need to restore
the multiplier ${\bf 1}_{\{t_n<t_{n+1}\}}$. 
Then the function $g$ (see the proof of formulas (\ref{july80013}), (\ref{july80015})) 
will be the product
of ${\bf 1}_{\{t_n<t_{n+1}\}}\psi_n(t_n)\psi_{n+1}(t_{n+1})$
and $r-2$ weight functions that are chosen so that
on the right-hand side of the equality 
$(\bar K)$ there is a scalar product in $L_2([t, T]^r)$
involving $s_q$ ($s_q$ is an approximation of $g$).

Using the above algorithm, we prove the equality (\ref{july90000})
for the case $k=2r$ $(r=2,3,\ldots).$
The equality (\ref{july90000}) is proved.

Note that the series on the left-hand side of 
(\ref{july90000}) converges absolutly since
its sum does not depend 
on permutations of basis functions
(here the basis in $L_2([t,T]^{r})$ is
$\left\{\phi_{j_1}(x_1)\ldots \phi_{j_r}(x_r)\right\}_{j_1,\ldots,j_r=0}^{\infty}$).

\section{Revision of Hypotheses on Expansion of Iterated Stra\-to\-no\-vich Stochastic
Integrals of Multiplicity $k$ $(k\in{\bf N})$}

In Sect.~2.5, we formulated three 
hypotheses on expansion of iterated Stra\-to\-no\-vich stochastic 
integrals based on the results obtained by the author
in the 2010s. In light of recent results (Theorems~2.3, 2.42--2.57),
a new vision of the above problem 
has appeared.
In particular, it became clear that it is possible 
to methodically obtain results related to the expansion
of iterated Stratonovich stochastic 
integrals for the case of an 
arbitrary 
complete ortho\-nor\-mal system of functions in the space $L_2([t, T])$
and $\psi_1(\tau),$ $\ldots,$ $\psi_k(\tau)$ $\in $ $L_2([t, T])$.

Definition (\ref{123321.2}) of the Stratonovich stochastic 
integral, which we mainly use in this book, imposes
its own limitations. In particular, this definition assumes that
$\psi_1(\tau),$ $\ldots,$ $\psi_k(\tau)$ are continuous functions 
at the interval $[t, T]$.

Based on Theorems~2.3, 2.42--2.57, we formulate the following
hypothesis on expansion of the sum $\bar J^{*}[\psi^{(k)}]_{T,t}^{(i_1\ldots i_k)}$
of iterated It\^{o} stochastic integrals (see (\ref{dsds9})).

\vspace{1mm}

{\bf Hypothesis~2.4.}\ {\it Suppose that $\{\phi_j(x)\}_{j=0}^{\infty}$
is an arbitrary complete orthonormal system of functions
in $L_2([t, T])$ and
$\psi_1(\tau),\ldots, \psi_k(\tau)\in L_2([t, T]).$
Then$,$ for the sum $\bar J^{*}[\psi^{(k)}]_{T,t}^{(i_1\ldots i_k)}$
of iterated It\^{o} stochastic integrals 
$$
\bar J^{*}[\psi^{(k)}]_{T,t}^{(i_1\ldots i_k)}=J[\psi^{(k)}]_{T,t}^{(i_1\ldots i_k)}+
\sum_{r=1}^{\left[k/2\right]}\frac{1}{2^r}
\sum_{(s_r,\ldots,s_1)\in {\rm A}_{k,r}}
J[\psi^{(k)}]_{T,t}^{s_r,\ldots,s_1}
$$
the following 
expansion 
\begin{equation}
\label{july300000}
~~~~~~~~~~\bar J^{*}[\psi^{(k)}]_{T,t}^{(i_1\ldots i_k)}=
\hbox{\vtop{\offinterlineskip\halign{
\hfil#\hfil\cr
{\rm l.i.m.}\cr
$\stackrel{}{{}_{p_1,\ldots,p_k\to \infty}}$\cr
}} }
\sum_{j_1=0}^{p_1}\ldots\sum_{j_k=0}^{p_k}
C_{j_k \ldots j_1}\prod\limits_{l=1}^k \zeta_{j_l}^{(i_l)}
\end{equation}

\noindent
that converges in the mean-square sense is valid, where 
$$
C_{j_k \ldots j_1}=\int\limits_t^T\psi_k(t_k)\phi_{j_k}(t_k)\ldots
\int\limits_t^{t_2}
\psi_1(t_1)\phi_{j_1}(t_1)
dt_1\ldots dt_k
$$
is the Fourier coefficient, 
${\rm l.i.m.}$ is a limit in the mean-square sense,
$i_1, \ldots, i_k=0, 1,\ldots,m,$
$$
\zeta_{j}^{(i)}=
\int\limits_t^T \phi_{j}(\tau) d{\bf w}_{\tau}^{(i)}
$$ 
are independent standard Gaussian random variables for various 
$i$ or $j$ {\rm (}in the case when $i\ne 0${\rm )},
${\bf w}_{\tau}^{(i)}$ 
$(i=1,\ldots,m)$ are independent 
standard Wiener processes$,$
${\bf w}_{\tau}^{(0)}=\tau;$ another notations are the same as
in Theorem~{\rm 2.12}.}

\vspace{2mm}

Using Theorem~2.12, we obtain the following hypothesis.

\vspace{1mm}

{\bf Hypothesis~2.5.}\ {\it Suppose that $\{\phi_j(x)\}_{j=0}^{\infty}$
is an arbitrary complete orthonormal system of functions
in $L_2([t, T])$ and
$\psi_1(\tau),\ldots, \psi_k(\tau)$ are continuous functions
at the interval $[t, T].$
Then$,$ for the iterated Stratonovich sto\-chas\-tic integral 
of arbitrary multiplicity $k$
$$
J^{*}[\psi^{(k)}]_{T,t}^{(i_1\ldots i_k)}=
{\int\limits_t^{*}}^T
\psi_k(t_k) \ldots 
{\int\limits_t^{*}}^{t_{2}}
\psi_1(t_1) d{\bf w}_{t_1}^{(i_1)}\ldots
d{\bf w}_{t_k}^{(i_k)}
$$
the following 
expansion 
$$
J^{*}[\psi^{(k)}]_{T,t}^{(i_1\ldots i_k)}=
\hbox{\vtop{\offinterlineskip\halign{
\hfil#\hfil\cr
{\rm l.i.m.}\cr
$\stackrel{}{{}_{p_1,\ldots,p_k\to \infty}}$\cr
}} }
\sum\limits_{j_1=0}^{p_1}\ldots\sum\limits_{j_k=0}^{p_k}
C_{j_k \ldots j_1}\prod\limits_{l=1}^k \zeta_{j_l}^{(i_l)}
$$

\noindent
that converges in the mean-square sense is valid, where 
$$
C_{j_k \ldots j_1}=\int\limits_t^T\psi_k(t_k)\phi_{j_k}(t_k)\ldots
\int\limits_t^{t_2}
\psi_1(t_1)\phi_{j_1}(t_1)
dt_1\ldots dt_k
$$
is the Fourier coefficient, 
${\rm l.i.m.}$ is a limit in the mean-square sense,
$i_1, \ldots, i_k=0, 1,\ldots,m,$
$$
\zeta_{j}^{(i)}=
\int\limits_t^T \phi_{j}(\tau) d{\bf w}_{\tau}^{(i)}
$$ 
are independent standard Gaussian random variables for various 
$i$ or $j$ {\rm (}in the case when $i\ne 0${\rm )},
${\bf w}_{\tau}^{(i)}$ 
$(i=1,\ldots,m)$ are independent 
standard Wiener processes$,$
${\bf w}_{\tau}^{(0)}=\tau.$}

\section{Proof of Hypotheses~2.4, 2.5 Under the Condition (\ref{july90001})
for the Case $k\ge 2r,$ $p_1=\ldots=p_k=p$ and Under Some Additional
Assumptions}

Suppose that the equality

\vspace{-2mm}
$$
\lim\limits_{p\to\infty}
\sum\limits_{j_{g_1}, j_{g_3},\ldots ,j_{g_{2r-1}}=0}^p
C_{j_k\ldots j_1}\biggl|_{j_{g_1}=j_{g_2},\ldots, j_{g_{2r-1}}=j_{g_{2r}}}=
$$
\begin{equation}
\label{july90001}
=\frac{1}{2^r} \prod\limits_{l=1}^r {\bf 1}_{\{g_{2l}=g_{2l-1}+1\}}
C_{j_k \ldots j_1}\biggl|_{(j_{g_2} j_{g_1})\curvearrowright (\cdot)
\ldots (j_{g_{2r}} j_{g_{2r-1}})\curvearrowright (\cdot),
j_{g_{{}_{1}}}=~j_{g_{{}_{2}}},\ldots, j_{g_{{}_{2r-1}}}=~j_{g_{{}_{2r}}}
}\biggr.
\end{equation}

\vspace{2mm}
\noindent
is satisfied for all possible $g_1,g_2,\ldots,g_{2r-1},g_{2r}$ (see {\rm (\ref{leto5007after}))
and for any fixed $j_1,\ldots,j_q,\ldots,j_k$
($q\ne g_1, g_2, \ldots, g_{2r-1},$ $g_{2r}$),
where $k\ge 2r,$ $r=1,2,\ldots,[k/2],$ $C_{j_k\ldots j_1}$ is defined by (\ref{july15030}),
another notations are the same as in Theorem~2.49. Recall that the case
$k=2r$ is considered in Sect.~2.27.4.

Moreover, suppose that
the series 
$$
\lim\limits_{p\to\infty}
\sum\limits_{j_{g_1}, j_{g_3},\ldots ,j_{g_{2r-1}}=0}^p
C_{j_k\ldots j_1}\biggl|_{j_{g_1}=j_{g_2},\ldots, j_{g_{2r-1}}=j_{g_{2r}}}
$$

\vspace{2mm}
\noindent
converges absolutly for any fixed $j_1,\ldots,j_q,\ldots,j_k,$
where $q\ne g_1, g_2, \ldots, g_{2r-1},$ $g_{2r}$ and $k>2r.$

It should be noted that the above assumptions will be proved further (see Sect.~2.30).

Hypotheses~2.4 and 2.5 will be proved for the case $p_1=\ldots=p_k=p$
if we prove that (see 
Theorem~2.49 for the case $p_1=\ldots=p_k=p$)

\vspace{-2mm}
$$
\lim\limits_{p\to\infty}
\sum\limits_{\stackrel{j_1,\ldots,j_q,\ldots,j_k=0}{{}_{q\ne g_1, g_2, \ldots, g_{2r-1},
g_{2r}}}}^p
\Biggl(\sum\limits_{j_{g_1}, j_{g_3},\ldots ,j_{g_{2r-1}}=0}^p
C_{j_k\ldots j_1}\biggl|_{j_{g_1}=j_{g_2},\ldots, j_{g_{2r-1}}=j_{g_{2r}}}-\Biggr.
$$

\vspace{-3mm}
\begin{equation}
\label{july700000}
\Biggl.-\frac{1}{2^r} \prod\limits_{l=1}^r {\bf 1}_{\{g_{2l}=g_{2l-1}+1\}}
C_{j_k \ldots j_1}\biggl|_{(j_{g_2} j_{g_1})\curvearrowright (\cdot)
\ldots (j_{g_{2r}} j_{g_{2r-1}})\curvearrowright (\cdot),
j_{g_{{}_{1}}}=~j_{g_{{}_{2}}},\ldots, j_{g_{{}_{2r-1}}}=~j_{g_{{}_{2r}}}
}\biggr.\Biggr)^2=0
\end{equation}

\vspace{2mm}
\noindent 
for all $r=1,2,\ldots,[k/2]$ and 
for all possible $g_1,g_2,\ldots,g_{2r-1},g_{2r}$ (see {\rm (\ref{leto5007after})),
where notations are the same as in 
(\ref{july90001}).

We have
$$
\sum\limits_{\stackrel{j_1,\ldots,j_q,\ldots,j_k=0}{{}_{q\ne g_1, g_2, \ldots, g_{2r-1},
g_{2r}}}}^p
\Biggl(\sum\limits_{j_{g_1}, j_{g_3},\ldots ,j_{g_{2r-1}}=0}^p
C_{j_k\ldots j_1}\biggl|_{j_{g_1}=j_{g_2},\ldots, j_{g_{2r-1}}=j_{g_{2r}}}-\Biggr.
$$
$$
\Biggl.-\frac{1}{2^r} \prod\limits_{l=1}^r {\bf 1}_{\{g_{2l}=g_{2l-1}+1\}}
C_{j_k \ldots j_1}\biggl|_{(j_{g_2} j_{g_1})\curvearrowright (\cdot)
\ldots (j_{g_{2r}} j_{g_{2r-1}})\curvearrowright (\cdot),
j_{g_{{}_{1}}}=~j_{g_{{}_{2}}},\ldots, j_{g_{{}_{2r-1}}}=~j_{g_{{}_{2r}}}
}\biggr.\Biggr)^2\le
$$

\vspace{4mm}
$$
\le \sum\limits_{\stackrel{j_1,\ldots,j_q,\ldots,j_k=0}{{}_{q\ne g_1, g_2, \ldots, g_{2r-1},
g_{2r}}}}^{\infty}
\Biggl(\sum\limits_{j_{g_1}, j_{g_3},\ldots ,j_{g_{2r-1}}=0}^p
C_{j_k\ldots j_1}\biggl|_{j_{g_1}=j_{g_2},\ldots, j_{g_{2r-1}}=j_{g_{2r}}}-\Biggr.
$$

\vspace{-4mm}
\begin{equation}
\label{july90009xx}
\Biggl.-\frac{1}{2^r} \prod\limits_{l=1}^r {\bf 1}_{\{g_{2l}=g_{2l-1}+1\}}
C_{j_k \ldots j_1}\biggl|_{(j_{g_2} j_{g_1})\curvearrowright (\cdot)
\ldots (j_{g_{2r}} j_{g_{2r-1}})\curvearrowright (\cdot),
j_{g_{{}_{1}}}=~j_{g_{{}_{2}}},\ldots, j_{g_{{}_{2r-1}}}=~j_{g_{{}_{2r}}}
}\biggr.\Biggr)^2,
\end{equation}

\noindent
where
\begin{equation}
\label{july90009}
\sum\limits_{\stackrel{j_1,\ldots,j_q,\ldots,j_k=0}{{}_{q\ne g_1, g_2, \ldots, g_{2r-1},
g_{2r}}}}^{\infty}\stackrel{\sf def}{=}
\lim\limits_{q\to\infty}
\sum\limits_{\stackrel{j_1,\ldots,j_q,\ldots,j_k=0}{{}_{q\ne g_1, g_2, \ldots, g_{2r-1},
g_{2r}}}}^{q}.
\end{equation}

\vspace{2mm}

Consider the following analogue
of Monotone Convergence Theorem for infinite series.

\vspace{1mm}

{\bf Proposition 2.5.}\ {\it Suppose that $x_{m,n}\ge 0$ for all $m,n\in{\bf N},$
$$
\lim\limits_{m\to\infty}x_{m,n}=y_n\ \ \ (\hbox{for any fixed}~n\in {\bf N}),
$$

\noindent
and $x_{m,n}\le x_{m+1,n}$ for all $m\in{\bf N}$ and for any fixed $n\in {\bf N}.$
Then
\begin{equation}
\label{july90002x}
\lim_{m\to\infty}\sum\limits_{n=1}^{\infty}
x_{m,n}=\sum\limits_{n=1}^{\infty}
\lim_{m\to\infty}x_{m,n}=
\sum\limits_{n=1}^{\infty}y_n.
\end{equation}
}

\vspace{-1mm}

{\bf Proof.}\ Proposition 2.5 can be easily proved 
using the following version of Fatou's Lemma for infinite series
\begin{equation}
\label{july90002}
\sum\limits_{n=1}^{\infty}\liminf _{m\to\infty}
x_{m,n}\le \liminf _{m\to\infty}
\sum\limits_{n=1}^{\infty}x_{m,n},
\end{equation}

\noindent
where it is assumed that the conditions of Proposition 2.5
are fulfilled. Indeed, we have
$$
0\le x_{m,n}\le y_n.
$$

Then
$$
\sum\limits_{n=1}^{\infty}x_{m,n}\le 
\sum\limits_{n=1}^{\infty}y_n
$$
and (see (\ref{july90002}))
\begin{equation}
\label{july90003}
~~~~~~~~~~\limsup_{m\to\infty}\sum\limits_{n=1}^{\infty}x_{m,n}\le 
\sum\limits_{n=1}^{\infty}y_n=
\sum\limits_{n=1}^{\infty}\liminf_{m\to\infty}x_{m,n}
\le \liminf _{m\to\infty}
\sum\limits_{n=1}^{\infty}x_{m,n}.
\end{equation}

From (\ref{july90003}) we get
$$
\sum\limits_{n=1}^{\infty}y_n=\liminf_{m\to\infty}\sum\limits_{n=1}^{\infty}x_{m,n}=
\limsup_{m\to\infty}\sum\limits_{n=1}^{\infty}x_{m,n}=
\lim_{m\to\infty}\sum\limits_{n=1}^{\infty}x_{m,n},
$$

\noindent
i.e. the equality (\ref{july90002x}) is proved.

To prove (\ref{july90002}) we note that
$$
\inf_{j\ge m}x_{j,n}\le x_{k,n}\ \ \ (\forall k\ge m).
$$

Then
$$
\sum\limits_{n=1}^N \inf_{j\ge m}x_{j,n}\le \sum\limits_{n=1}^N x_{k,n}\ \ \ (\forall k\ge m)
$$
and
\begin{equation}
\label{july90004}
\sum\limits_{n=1}^N \inf_{j\ge m}x_{j,n}\le\inf\limits_{k\ge m}
\sum\limits_{n=1}^{N} x_{k,n}\le
\inf\limits_{k\ge m}
\sum\limits_{n=1}^{\infty} x_{k,n}.
\end{equation}

\vspace{2mm}

Passing to the limit $\lim\limits_{m\to\infty}$ in (\ref{july90004}), we obtain
\begin{equation}
\label{july90005}
\sum\limits_{n=1}^N \lim\limits_{m\to\infty} \inf_{j\ge m}x_{j,n}
\le
\lim\limits_{m\to\infty} \inf\limits_{k\ge m}
\sum\limits_{n=1}^{\infty} x_{k,n}.
\end{equation}

\vspace{2mm}

Passing to the limit $\lim\limits_{N\to\infty}$ in (\ref{july90005}), we get
$$
\sum\limits_{n=1}^{\infty} \lim\limits_{m\to\infty} \inf_{j\ge m}x_{j,n}
\le
\lim\limits_{m\to\infty} \inf\limits_{k\ge m}
\sum\limits_{n=1}^{\infty} x_{k,n},
$$

\noindent
i.e. the equality (\ref{july90002}) is satisfied.
Proposition~2.5 is proved.

\vspace{2mm}

{\bf Proposition~2.6.}\ {\it Suppose that 
\begin{equation}
\label{july90017}
\sum\limits_{j=1}^{\infty} g_{j,n}=0,
\end{equation}

\noindent
the series {\rm (\ref{july90017})} converges 
absolutely for any fixed $n\in{\bf N}$ and
$$
\sum\limits_{n=1}^{\infty}
\left(\sum\limits_{j=1}^{\infty} \left|g_{j,n}\right|\right)^2<\infty.
$$

Then
\begin{equation}
\label{july90015}
~~~~~~~~\lim_{m\to\infty}\sum\limits_{n=1}^{\infty}
\left(\sum\limits_{j=1}^m g_{j,n}\right)^2=\sum\limits_{n=1}^{\infty}
\lim_{m\to\infty}\left(\sum\limits_{j=1}^m g_{j,n}\right)^2=0.
\end{equation}
}

\vspace{2mm}

{\bf Proof.}\ We have 
$g_{j,n}=g_{j,n}^{+}-g_{j,n}^{-},$
$\left\vert g_{j,n}\right\vert =g_{j,n}^{+}+g_{j,n}^{-},$ 
where 
$$
g_{j,n}^{+}=\max\{g_{j,n}, 0\}=\frac{1}{2}\left(\left\vert g_{j,n}\right\vert + g_{j,n}\right)\ge 0,
$$
$$
g_{j,n}^{-}=-\min\{g_{j,n}, 0\}=\frac{1}{2}\left(\left\vert g_{j,n}\right\vert - g_{j,n}\right)
\ge 0.
$$ 

Moreover,
\begin{equation}
\label{july90016s}
\sum\limits_{j=1}^{\infty}g_{j,n}=
\sum\limits_{j=1}^{\infty} g_{j,n}^{+}-\sum\limits_{j=1}^{\infty} g_{j,n}^{-}=0,
\end{equation}
\begin{equation}
\label{july90016}
~~~~~~~~~~~\sum\limits_{j=1}^{\infty}\left\vert g_{j,n}\right\vert =
\sum\limits_{j=1}^{\infty} g_{j,n}^{+}+\sum\limits_{j=1}^{\infty} g_{j,n}^{-}=
2\sum\limits_{j=1}^{\infty} g_{j,n}^{+}=
2\sum\limits_{j=1}^{\infty} g_{j,n}^{-}.
\end{equation}

\vspace{2mm}

Since the series (\ref{july90017}) converges absolutely, then by virtue 
of the equality (\ref{july90016}) the series (with non-negative terms) on the right-hand side of 
(\ref{july90016}) and on the right-hand side of (\ref{july90016s}) converge. 

Further, using Proposition~2.5 and (\ref{july90016s}), (\ref{july90016}), we obtain
$$
\lim_{m\to\infty}\sum\limits_{n=1}^{\infty}
\left(\sum\limits_{j=1}^m g_{j,n}\right)^2=
\lim_{m\to\infty}\sum\limits_{n=1}^{\infty}
\left(\sum\limits_{j=1}^m g_{j,n}^{+}- \sum\limits_{j=1}^m g_{j,n}^{-}\right)^2=
$$
$$
=\lim_{m\to\infty}\sum\limits_{n=1}^{\infty}
\left(\sum\limits_{j=1}^m g_{j,n}^{+}\right)^2-
\lim_{m\to\infty}\sum\limits_{n=1}^{\infty}
\left(2\sum\limits_{j=1}^m g_{j,n}^{+}\sum\limits_{j=1}^m g_{j,n}^{-}\right)
+
$$
$$
+
\lim_{m\to\infty}\sum\limits_{n=1}^{\infty}
\left(\sum\limits_{j=1}^m g_{j,n}^{-}\right)^2=
$$

$$
=\sum\limits_{n=1}^{\infty}\lim_{m\to\infty}
\left(\sum\limits_{j=1}^m g_{j,n}^{+}\right)^2-
\sum\limits_{n=1}^{\infty}\lim_{m\to\infty}
\left(2\sum\limits_{j=1}^m g_{j,n}^{+}\sum\limits_{j=1}^m g_{j,n}^{-}\right)+
$$

$$
+
\sum\limits_{n=1}^{\infty}\lim_{m\to\infty}
\left(\sum\limits_{j=1}^m g_{j,n}^{-}\right)^2=
$$

$$
=\sum\limits_{n=1}^{\infty}
\left(\sum\limits_{j=1}^{\infty} g_{j,n}^{+}\right)^2-
2\sum\limits_{n=1}^{\infty}
\left(\sum\limits_{j=1}^{\infty} g_{j,n}^{+}\sum\limits_{j=1}^{\infty} g_{j,n}^{-}\right)+
$$

$$
+
\sum\limits_{n=1}^{\infty}
\left(\sum\limits_{j=1}^{\infty} g_{j,n}^{-}\right)^2=
$$

$$
=\frac{1}{4}\sum\limits_{n=1}^{\infty}
\left(\sum\limits_{j=1}^{\infty} \left|g_{j,n}\right|\right)^2-
\frac{1}{2}\sum\limits_{n=1}^{\infty}
\left(\sum\limits_{j=1}^{\infty} \left|g_{j,n}\right|\right)^2
+
$$

$$
+
\frac{1}{4}\sum\limits_{n=1}^{\infty}
\left(\sum\limits_{j=1}^{\infty} \left|g_{j,n}\right|\right)^2=0.
$$

\vspace{3mm}
\noindent
Proposition~2.6 is proved.

It is easy to see that by analogy with Propositions~2.5, 2.6 
the following statements can be proved.

\vspace{2mm}

{\bf Proposition 2.7.}\ {\it Suppose that $h_{p,k_1,\ldots,k_d}\ge 0$ 
for all $p\in {\bf N}$ and for any fixed $k_1,\ldots,k_d\in {\bf N},$
$$
\lim\limits_{p\to\infty}h_{p,k_1,\ldots,k_d}=u_{k_1,\ldots,k_d}\ \ \ (\hbox{for any fixed}~
k_1,\ldots,k_d\in {\bf N}),
$$

\noindent
and $h_{p,k_1,\ldots,k_d}\le h_{p+1,k_1,\ldots,k_d}$ for all $p\in {\bf N}$
and for any fixed $k_1,\ldots,k_d\in {\bf N}.$ 
Then
\begin{equation}
\label{july90006}
~~~~~~\lim_{p\to\infty}\sum\limits_{k_1,\ldots,k_d=1}^{\infty}
h_{p,k_1,\ldots,k_d}=\sum\limits_{k_1,\ldots,k_d=1}^{\infty}
\lim_{p\to\infty}h_{p,k_1,\ldots,k_d}=
\sum\limits_{k_1,\ldots,k_d=1}^{\infty}
u_{k_1,\ldots,k_d},
\end{equation}

\vspace{1mm}
\noindent
where $h_{p,k_1,\ldots,k_d}, u_{k_1,\ldots,k_d}\in {\bf R},$ 
$d\in {\bf N},$ the series on the left-hand side of {\rm (\ref{july90006})}
is understood in the same sense as in {\rm (\ref{july90009})}.
}

\vspace{2mm}

{\bf Proposition 2.8.}\ {\it Suppose that 
\begin{equation}
\label{july700010}
~~~~~~\lim\limits_{p\to\infty}
\sum\limits_{j_1,\ldots,j_q=1}^{p} h_{j_1,\ldots,j_q,k_1,\ldots,k_d}
\stackrel{\sf def}{=}\sum\limits_{j_1,\ldots,j_q=1}^{\infty} h_{j_1,\ldots,j_q,k_1,\ldots,k_d}=0,
\end{equation}

\noindent
the series {\rm (\ref{july700010})}
converges absolutely for any fixed $k_1,\ldots,k_d\in{\bf N}$ 
and
$$
\sum\limits_{k_1,\ldots,k_d=1}^{\infty} 
\left(\sum\limits_{j_1,\ldots,j_q=1}^{\infty} \left|h_{j_1,\ldots,j_q,k_1,\ldots,k_d}\right|
\right)^2<\infty.
$$

Then
$$
\lim\limits_{p\to\infty}
\sum\limits_{k_1,\ldots,k_d=1}^{\infty} 
\left(\sum\limits_{j_1,\ldots,j_q=1}^{p} h_{j_1,\ldots,j_q,k_1,\ldots,k_d}\right)^2
=
$$
$$
=
\sum\limits_{k_1,\ldots,k_d=1}^{\infty} 
\lim\limits_{p\to\infty}
\left(\sum\limits_{j_1,\ldots,j_q=1}^{p} h_{j_1,\ldots,j_q,k_1,\ldots,k_d}\right)^2=
0,
$$

\noindent
where 
$$
\lim\limits_{n\to\infty}
\sum\limits_{k_1,\ldots,k_d=1}^{n}
\stackrel{\sf def}{=}\sum\limits_{k_1,\ldots,k_d=1}^{\infty},
$$

\vspace{2mm}
\noindent
$h_{j_1,\ldots,j_q,k_1,\ldots,k_d}\in {\bf R}$ and
$d, q\in {\bf N}.$
}

Obviously, Proposition~2.8 follows from Proposition~2.7
in the same way as Proposition~2.6 follows from Proposition~2.5.
Applying Proposition~2.8 to the right-hand side of (\ref{july90009xx})
(using (\ref{july90001}) and the absolute convergence of the series on the left-hand side of 
(\ref{july90001})),
we obtain (\ref{july700000}). 
At that, we used the conditions
\begin{equation}
\label{slovo10}
~~~~\sum\limits_{\stackrel{j_1,\ldots,j_q,\ldots,j_k=0}{{}_{q\ne g_1, g_2, \ldots, g_{2r-1},
g_{2r}}}}^{\infty}
\left(\lim\limits_{p\to\infty}\sum\limits_{j_{g_1}, j_{g_3},\ldots ,j_{g_{2r-1}}=0}^p
\left\vert C_{j_k\ldots j_1}\biggl|_{j_{g_1}=j_{g_2},\ldots, j_{g_{2r-1}}=j_{g_{2r}}}
\right\vert\right)^2<\infty,
\end{equation}
\begin{equation}
\label{july700011}
~\sum\limits_{\stackrel{j_1,\ldots,j_q,\ldots,j_k=0}{{}_{q\ne g_1, g_2, \ldots, g_{2r-1},
g_{2r}}}}^{\infty}
\left(
C_{j_k \ldots j_1}\biggl|_{(j_{g_2} j_{g_1})\curvearrowright (\cdot)
\ldots (j_{g_{2r}} j_{g_{2r-1}})\curvearrowright (\cdot),
j_{g_{{}_{1}}}=~j_{g_{{}_{2}}},\ldots, j_{g_{{}_{2r-1}}}=~j_{g_{{}_{2r}}}
}\biggr.\right)^2<\infty.
\end{equation}

\vspace{2mm}

Note that (\ref{july700011}) follows from the Parseval equality
since the expression
$$
C_{j_k \ldots j_1}\biggl|_{(j_{g_2} j_{g_1})\curvearrowright (\cdot)
\ldots (j_{g_{2r}} j_{g_{2r-1}})\curvearrowright (\cdot),
j_{g_{{}_{1}}}=~j_{g_{{}_{2}}},\ldots, j_{g_{{}_{2r-1}}}=~j_{g_{{}_{2r}}}
}\biggr.\stackrel{\sf def}{=}H_{j_{q_1}\ldots j_{q_{k-2r}}}
$$ 

\noindent
is a finite linear combination 
of the Fourier coefficients
of $L_2([t,T]^{k-2r})$--func\-ti\-ons
after iteratively
applying
transformations (\ref{july100003}), (\ref{july100004}) (see Sect.~2.30) to
$H_{j_{q_1}\ldots j_{q_{k-2r}}}$
for integrations not involving the basis functions
$\phi_{j_{q_1}},\ldots, \phi_{j_{q_{k-2r}}}.$

Let us consider another sufficient condition under which the equality
(\ref{july700000}) is satisfied.
Suppose that $k>2r$ and
$$
\exists\ \ \ \lim\limits_{p,q\to\infty}
\sum\limits_{\stackrel{j_1,\ldots,j_m,\ldots,j_k=0}{{}_{m\ne g_1, g_2, \ldots, g_{2r-1},
g_{2r}}}}^q
\Biggl(\sum\limits_{j_{g_1}, j_{g_3},\ldots ,j_{g_{2r-1}}=0}^p
C_{j_k\ldots j_1}\biggl|_{j_{g_1}=j_{g_2},\ldots, j_{g_{2r-1}}=j_{g_{2r}}}-\Biggr.
$$

\vspace{-2mm}
\begin{equation}
\label{july700001xyz}
\Biggl.-\frac{1}{2^r} \prod\limits_{l=1}^r {\bf 1}_{\{g_{2l}=g_{2l-1}+1\}}
C_{j_k \ldots j_1}\biggl|_{(j_{g_2} j_{g_1})\curvearrowright (\cdot)
\ldots (j_{g_{2r}} j_{g_{2r-1}})\curvearrowright (\cdot),
j_{g_{{}_{1}}}=~j_{g_{{}_{2}}},\ldots, j_{g_{{}_{2r-1}}}=~j_{g_{{}_{2r}}}
}\biggr.\Biggr)^2<\infty
\end{equation}

\vspace{2mm}
\noindent 
for all $r=1,2,\ldots,[k/2],$
where notations are the same as in 
(\ref{july90001}).
Then, by Proposition~1.1 (see Sect.~1.7.2) and (\ref{july90001})
we obtain
$$
\lim\limits_{p,q\to\infty}
\sum\limits_{\stackrel{j_1,\ldots,j_m,\ldots,j_k=0}{{}_{m\ne g_1, g_2, \ldots, g_{2r-1},
g_{2r}}}}^q
\Biggl(\sum\limits_{j_{g_1}, j_{g_3},\ldots ,j_{g_{2r-1}}=0}^p
C_{j_k\ldots j_1}\biggl|_{j_{g_1}=j_{g_2},\ldots, j_{g_{2r-1}}=j_{g_{2r}}}-\Biggr.
$$

\vspace{-3mm}
$$
\Biggl.-\frac{1}{2^r} \prod\limits_{l=1}^r {\bf 1}_{\{g_{2l}=g_{2l-1}+1\}}
C_{j_k \ldots j_1}\biggl|_{(j_{g_2} j_{g_1})\curvearrowright (\cdot)
\ldots (j_{g_{2r}} j_{g_{2r-1}})\curvearrowright (\cdot),
j_{g_{{}_{1}}}=~j_{g_{{}_{2}}},\ldots, j_{g_{{}_{2r-1}}}=~j_{g_{{}_{2r}}}
}\biggr.\Biggr)^2=
$$

$$
=\lim\limits_{q\to\infty}
\sum\limits_{\stackrel{j_1,\ldots,j_m,\ldots,j_k=0}{{}_{m\ne g_1, g_2, \ldots, g_{2r-1},
g_{2r}}}}^q
\lim\limits_{p\to\infty}\Biggl(\sum\limits_{j_{g_1}, j_{g_3},\ldots ,j_{g_{2r-1}}=0}^p
C_{j_k\ldots j_1}\biggl|_{j_{g_1}=j_{g_2},\ldots, j_{g_{2r-1}}=j_{g_{2r}}}-\Biggr.
$$

\vspace{-3mm}
$$
\Biggl.-\frac{1}{2^r} \prod\limits_{l=1}^r {\bf 1}_{\{g_{2l}=g_{2l-1}+1\}}
C_{j_k \ldots j_1}\biggl|_{(j_{g_2} j_{g_1})\curvearrowright (\cdot)
\ldots (j_{g_{2r}} j_{g_{2r-1}})\curvearrowright (\cdot),
j_{g_{{}_{1}}}=~j_{g_{{}_{2}}},\ldots, j_{g_{{}_{2r-1}}}=~j_{g_{{}_{2r}}}
}\biggr.\Biggr)^2=0.
$$

\vspace{3mm}

Thus, we get
$$
\lim\limits_{p,q\to\infty}
\sum\limits_{\stackrel{j_1,\ldots,j_m,\ldots,j_k=0}{{}_{m\ne g_1, g_2, \ldots, g_{2r-1},
g_{2r}}}}^q
\Biggl(\sum\limits_{j_{g_1}, j_{g_3},\ldots ,j_{g_{2r-1}}=0}^p
C_{j_k\ldots j_1}\biggl|_{j_{g_1}=j_{g_2},\ldots, j_{g_{2r-1}}=j_{g_{2r}}}-\Biggr.
$$

\vspace{-3mm}
\begin{equation}
\label{july700002xyz}
\Biggl.-\frac{1}{2^r} \prod\limits_{l=1}^r {\bf 1}_{\{g_{2l}=g_{2l-1}+1\}}
C_{j_k \ldots j_1}\biggl|_{(j_{g_2} j_{g_1})\curvearrowright (\cdot)
\ldots (j_{g_{2r}} j_{g_{2r-1}})\curvearrowright (\cdot),
j_{g_{{}_{1}}}=~j_{g_{{}_{2}}},\ldots, j_{g_{{}_{2r-1}}}=~j_{g_{{}_{2r}}}
}\biggr.\Biggr)^2=0.
\end{equation}

\vspace{2mm}
\noindent
Substituting $p=q$ in (\ref{july700002xyz}), we obtain
(\ref{july700000}) (recall that the case
$k=2r$ is proved in Sect.~2.27.4).

As a result, 
Hypotheses~2.4 and 2.5 are proved under 
the conditions formulated above in this section.

\section{Expansion of Iterated Stratonovich Stochastic Integrals
of Arbitrary Multiplicity $k$ $(k\in{\bf N})$. 
The Case of an Ar\-bit\-ra\-ry Complete Orthonormal System of 
Functions in $L_2([t,T]),$ $\psi_1(\tau),\ldots, \psi_k(\tau)
\in L_2([t,T]).$
Proof of Hypotheses~2.4, 2.5 for the Case $p_1=\ldots=p_k=p$
and Under the Condition (\ref{july700001xyz})}

This section is devoted to the following two theorems.

\vspace{2mm}

{\bf Theorem~2.58.}\ {\it Suppose that 
the condition {\rm (\ref{july700001xyz})} 
is fulfilled$,$
$\{\phi_j(x)\}_{j=0}^{\infty}$
is an arbitrary complete orthonormal system of functions
in $L_2([t, T])$ and
$\psi_1(\tau),\ldots, \psi_k(\tau)\in L_2([t, T]).$
Then$,$ for the sum $\bar J^{*}[\psi^{(k)}]_{T,t}^{(i_1\ldots i_k)}$
of iterated It\^{o} stochastic integrals 
$$
\bar J^{*}[\psi^{(k)}]_{T,t}^{(i_1\ldots i_k)}=J[\psi^{(k)}]_{T,t}^{(i_1\ldots i_k)}+
\sum_{r=1}^{\left[k/2\right]}\frac{1}{2^r}
\sum_{(s_r,\ldots,s_1)\in {\rm A}_{k,r}}
J[\psi^{(k)}]_{T,t}^{s_r,\ldots,s_1}
$$
the following 
expansion 
$$
\bar J^{*}[\psi^{(k)}]_{T,t}^{(i_1\ldots i_k)}=
\hbox{\vtop{\offinterlineskip\halign{
\hfil#\hfil\cr
{\rm l.i.m.}\cr
$\stackrel{}{{}_{p\to \infty}}$\cr
}} }
\sum_{j_1,\ldots,j_k=0}^{p}
C_{j_k \ldots j_1}\prod\limits_{l=1}^k \zeta_{j_l}^{(i_l)}
$$

\noindent
that converges in the mean-square sense is valid, where 
\begin{equation}
\label{july99999}
~~~~~~~~~~C_{j_k \ldots j_1}=\int\limits_t^T\psi_k(t_k)\phi_{j_k}(t_k)\ldots
\int\limits_t^{t_2}
\psi_1(t_1)\phi_{j_1}(t_1)
dt_1\ldots dt_k
\end{equation}
is the Fourier coefficient, 
${\rm l.i.m.}$ is a limit in the mean-square sense,
$i_1, \ldots, i_k=0, 1,\ldots,m,$
$$
\zeta_{j}^{(i)}=
\int\limits_t^T \phi_{j}(\tau) d{\bf w}_{\tau}^{(i)}
$$ 
are independent standard Gaussian random variables for various 
$i$ or $j$ {\rm (}in the case when $i\ne 0${\rm )},
${\bf w}_{\tau}^{(i)}$ 
$(i=1,\ldots,m)$ are independent 
standard Wiener processes$,$
${\bf w}_{\tau}^{(0)}=\tau;$ another notations are the same as
in Theorem~{\rm 2.12}.}

\vspace{2mm}

Using Theorem~2.12, we obtain the following corollary of Theorem~2.58.

\vspace{2mm}

{\bf Theorem~2.59.}\ {\it Suppose that 
the condition {\rm (\ref{july700001xyz})} 
is satisfied$,$
$\{\phi_j(x)\}_{j=0}^{\infty}$
is an arbitrary complete orthonormal system of functions
in $L_2([t, T])$ and
$\psi_1(\tau),\ldots, \psi_k(\tau)$ are continuous functions
at the interval $[t, T].$
Then$,$ for the iterated Stratonovich sto\-chas\-tic integral 
of multiplicity $k$ $(k\in{\bf N})$
$$
J^{*}[\psi^{(k)}]_{T,t}^{(i_1\ldots i_k)}=
{\int\limits_t^{*}}^T
\psi_k(t_k) \ldots 
{\int\limits_t^{*}}^{t_{2}}
\psi_1(t_1) d{\bf w}_{t_1}^{(i_1)}\ldots
d{\bf w}_{t_k}^{(i_k)}
$$
the following 
expansion 
\begin{equation}
\label{july300001}
~~~~~~~~~~J^{*}[\psi^{(k)}]_{T,t}^{(i_1\ldots i_k)}=
\hbox{\vtop{\offinterlineskip\halign{
\hfil#\hfil\cr
{\rm l.i.m.}\cr
$\stackrel{}{{}_{p\to \infty}}$\cr
}} }
\sum\limits_{j_1,\ldots,j_k=0}^{p}
C_{j_k \ldots j_1}\prod\limits_{l=1}^k \zeta_{j_l}^{(i_l)}
\end{equation}

\noindent
that converges in the mean-square sense is valid, where 
$$
C_{j_k \ldots j_1}=\int\limits_t^T\psi_k(t_k)\phi_{j_k}(t_k)\ldots
\int\limits_t^{t_2}
\psi_1(t_1)\phi_{j_1}(t_1)
dt_1\ldots dt_k
$$
is the Fourier coefficient, 
${\rm l.i.m.}$ is a limit in the mean-square sense,
$i_1, \ldots, i_k=0, 1,\ldots,m,$
$$
\zeta_{j}^{(i)}=
\int\limits_t^T \phi_{j}(\tau) d{\bf w}_{\tau}^{(i)}
$$ 
are independent standard Gaussian random variables for various 
$i$ or $j$ {\rm (}when $i\ne 0${\rm )},
${\bf w}_{\tau}^{(i)}$ 
$(i=1,\ldots,m)$ are independent 
standard Wiener processes$,$
${\bf w}_{\tau}^{(0)}=\tau.$}

\vspace{2mm}

{\bf Proof of Theorem~2.58.}\ According to the results
of Sect.~2.29, Theorem~2.58 will be proved if we prove (see (\ref{july90001})) that
the equality

\vspace{-2mm}
$$
\lim\limits_{p\to\infty}
\sum\limits_{j_{g_1}, j_{g_3},\ldots ,j_{g_{2r-1}}=0}^p
C_{j_k\ldots j_1}\biggl|_{j_{g_1}=j_{g_2},\ldots, j_{g_{2r-1}}=j_{g_{2r}}}=
$$

\vspace{-2mm}
\begin{equation}
\label{july100000}
=\frac{1}{2^r} \prod\limits_{l=1}^r {\bf 1}_{\{g_{2l}=g_{2l-1}+1\}}
C_{j_k \ldots j_1}\biggl|_{(j_{g_2} j_{g_1})\curvearrowright (\cdot)
\ldots (j_{g_{2r}} j_{g_{2r-1}})\curvearrowright (\cdot),
j_{g_{{}_{1}}}=~j_{g_{{}_{2}}},\ldots, j_{g_{{}_{2r-1}}}=~j_{g_{{}_{2r}}}
}\biggr.
\end{equation}

\vspace{3mm}
\noindent
is satisfied for all possible $g_1,g_2,\ldots,g_{2r-1},g_{2r}$ (see {\rm (\ref{leto5007after})),
where $k\ge 2r,$ $r=1,2,\ldots,[k/2],$ $C_{j_k\ldots j_1}$ is defined by (\ref{july99999}),
another notations are the same as in Theorem~2.49.

Moreover (assuming that (\ref{july100000}) is proved), 
the series on the left-hand side of 
(\ref{july100000}) converges absolutly (the case $k=2r$ (see Sect.~2.27.4))
and converges absolutly
for any fixed $j_1,\ldots,j_q,\ldots,j_k$
and $q\ne g_1, g_2, \ldots, g_{2r-1},$ $g_{2r}$
(the case $k>2r$)
since
its sum does not depend 
on permutations of basis functions
(here the basis in $L_2([t,T]^{r})$ is
$\left\{\phi_{j_1}(x_1)\ldots \phi_{j_r}(x_r)\right\}_{j_1,\ldots,j_r=0}^{\infty}$).
Recall that any permutation of basis functions in a Hilbert space forms a basis 
in this Hilbert space \cite{gohb}.

The case
$k=2r$ of (\ref{july100000}) is considered in Sect.~2.27.4.
Consider the case $k>2r$ of (\ref{july100000}).

Using Fubini's Theorem, we obtain

\vspace{-5mm}
$$
\int\limits_t^T h_{k}(t_k)\ldots \int\limits_t^{t_{l+2}} h_{l+1}(t_{l+1})
\int\limits_t^{t_{l+1}} h_{l}(t_{l})
\int\limits_t^{t_{l}} h_{l-1}(t_{l-1})\ldots
\int\limits_t^{t_2} h_{1}(t_1)
dt_1\ldots 
$$

\vspace{-3mm}
$$
\ldots
dt_{l-1}dt_{l}dt_{l+1}\ldots dt_k=
$$

\vspace{-2mm}
$$
=\int\limits_t^T h_{k}(t_k)\ldots \int\limits_t^{t_{l+2}} h_{l+1}(t_{l+1})
\int\limits_t^{t_{l+1}} h_{1}(t_{1})
\int\limits_{t_1}^{t_{l+1}} h_{2}(t_{2})\ldots
\int\limits_{t_{l-2}}^{t_{l+1}} h_{l-1}(t_{l-1})
\int\limits_{t_{l-1}}^{t_{l+1}} h_{l}(t_{l})dt_l\times
$$

\vspace{-3mm}
$$
\times dt_{l-1}\ldots dt_2dt_{1}dt_{l+1}\ldots dt_k=
$$
$$
=\int\limits_t^T h_{k}(t_k)\ldots \int\limits_t^{t_{l+2}} h_{l+1}(t_{l+1})
\left(\int\limits_{t}^{t_{l+1}} h_{l}(t_{l})dt_l\right)\int\limits_t^{t_{l+1}} h_{1}(t_{1})
\int\limits_{t_1}^{t_{l+1}} h_{2}(t_{2})\ldots
\int\limits_{t_{l-2}}^{t_{l+1}} h_{l-1}(t_{l-1})
\times
$$
$$
\times dt_{l-1}\ldots dt_2dt_{1}dt_{l+1}\ldots dt_k-
$$

\vspace{-4mm}
$$
-\int\limits_t^T h_{k}(t_k)\ldots \int\limits_t^{t_{l+2}} h_{l+1}(t_{l+1})
\int\limits_t^{t_{l+1}} h_{1}(t_{1})
\int\limits_{t_1}^{t_{l+1}} h_{2}(t_{2})\ldots
\int\limits_{t_{l-2}}^{t_{l+1}} h_{l-1}(t_{l-1})
\left(\int\limits_{t}^{t_{l-1}} h_{l}(t_{l})dt_l\right)\times
$$
$$
\times dt_{l-1}\ldots dt_2dt_{1}dt_{l+1}\ldots dt_k=
$$

\vspace{-1mm}
$$
=\int\limits_t^T h_{k}(t_k)\ldots \int\limits_t^{t_{l+2}} h_{l+1}(t_{l+1})
\left(\int\limits_t^{t_{l+1}} h_{l}(t_{l})dt_l\right)
\int\limits_t^{t_{l+1}} h_{l-1}(t_{l-1})\ldots
$$
$$
\ldots 
\int\limits_t^{t_2} h_{1}(t_1)
dt_1\ldots dt_{l-1}dt_{l+1}\ldots dt_k-
$$
$$
-\int\limits_t^T h_{k}(t_k)\ldots \int\limits_t^{t_{l+2}} h_{l+1}(t_{l+1})
\int\limits_t^{t_{l+1}} h_{l-1}(t_{l-1})\left(\int\limits_t^{t_{l-1}} h_{l}(t_{l})dt_l\right)
\int\limits_t^{t_{l-1}} h_{l-2}(t_{l-2})
\ldots
$$
\begin{equation}
\label{july100003}
\ldots
\int\limits_t^{t_2} h_{1}(t_1)
dt_1\ldots dt_{l-2}dt_{l-1}dt_{l+1}\ldots dt_k,
\end{equation}

\vspace{1mm}
\noindent
where $2<l<k-1$ and $h_1(\tau),\ldots,h_k(\tau)\in L_2([t, T]).$ 

By analogy with (\ref{july100003}) we have for $l=k$
$$
\int\limits_t^{T} h_{l}(t_{l})
\int\limits_t^{t_{l}} h_{l-1}(t_{l-1})\ldots
\int\limits_t^{t_2} h_{1}(t_1)
dt_1\ldots 
dt_{l-1}dt_{l}=
$$
$$
=\int\limits_{t}^{T} h_{1}(t_{1})
\int\limits_{t_1}^{T} h_{2}(t_{2})\ldots
\int\limits_{t_{l-2}}^{T} h_{l-1}(t_{l-1})\int\limits_{t_{l-1}}^{T} h_{l}(t_{l})
dt_ldt_{l-1}\ldots dt_2dt_{1}=
$$

\vspace{-2mm}
$$
=\left(\int\limits_{t}^{T} h_{l}(t_{l})
dt_l\right)\int\limits_{t}^{T} h_{1}(t_{1})
\int\limits_{t_1}^{T} h_{2}(t_{2})\ldots
\int\limits_{t_{l-2}}^{T} h_{l-1}(t_{l-1})
dt_{l-1}\ldots dt_2dt_{1}-
$$
$$
-\int\limits_{t}^{T}h_{1}(t_{1})
\int\limits_{t_1}^{T} h_{2}(t_{2})\ldots
\int\limits_{t_{l-2}}^{T} h_{l-1}(t_{l-1})\left(
\int\limits_{t}^{t_{l-1}} h_{l}(t_{l})dt_l\right)
dt_{l-1}\ldots dt_2dt_{1}=
$$

\vspace{-1mm}
$$
=\left(\int\limits_{t}^{T} h_{l}(t_{l})
dt_l\right)\int\limits_{t}^{T} h_{l-1}(t_{l-1})
\ldots
\int\limits_{t}^{t_2} h_{1}(t_{1})
dt_{1}\ldots dt_{l-1}-
$$

\vspace{-3mm}
\begin{equation}
\label{july100004}
-\int\limits_{t}^{T}
h_{l-1}(t_{l-1})\left(
\int\limits_{t}^{t_{l-1}} h_{l}(t_{l})dt_l\right)
\int\limits_{t}^{t_{l-1}} h_{l-2}(t_{l-2})\ldots
\int\limits_{t}^{t_2} 
h_{1}(t_{1})
dt_{1}\ldots dt_{l-1}.
\end{equation}

\vspace{2mm}

We will assume that for $l=1$ the transformation 
(\ref{july100003}) is not carried out since
$$
\int\limits_t^{t_2} h_{1}(t_1)
dt_1
$$

\noindent
is the innermost integral on the left-hand side of (\ref{july100003}).
The formulas (\ref{july100003}), (\ref{july100004})
will be used further.

Let us carry out the transformations 
(\ref{july100003}), (\ref{july100004}) 
for 
$$
C_{j_k\ldots j_1}\biggl|_{j_{g_1}=j_{g_2},\ldots, j_{g_{2r-1}}=j_{g_{2r}}}
$$
iteratively for $j_1,\ldots,j_q,\ldots,j_k$
$(q\ne g_1, g_2, \ldots, g_{2r-1},$ $g_{2r}).$
As a result, we obtain 
$$
C_{j_k\ldots j_1}\biggl|_{j_{g_1}=j_{g_2},\ldots, j_{g_{2r-1}}=j_{g_{2r}}}=
$$
\begin{equation}
\label{july100005}
=\sum\limits_{d=1}^{2^{k-2r}}(-1)^{d-1}
\left(\hat C^{(d)}_{j_k\ldots j_1}\biggl|_{j_{g_1}=j_{g_2},\ldots, j_{g_{2r-1}}=j_{g_{2r}}}-~
\bar C^{(d)}_{j_k\ldots j_1}\biggl|_{j_{g_1}=j_{g_2},\ldots, j_{g_{2r-1}}=j_{g_{2r}}}\right),
\end{equation}

\vspace{3mm}
\noindent
where some terms in the sum
$$
\sum\limits_{d=1}^{2^{k-2r}}
$$

\noindent
can be identically equal to zero due to
the remark to (\ref{july100003}), (\ref{july100004}).

Using
(\ref{july100005}), we get
$$
\lim\limits_{p\to\infty}
\sum\limits_{j_{g_1}, j_{g_3},\ldots ,j_{g_{2r-1}}=0}^p
C_{j_k\ldots j_1}\biggl|_{j_{g_1}=j_{g_2},\ldots, j_{g_{2r-1}}=j_{g_{2r}}}=
$$

$$
=\lim\limits_{p\to\infty}
\sum\limits_{j_{g_1}, j_{g_3},\ldots ,j_{g_{2r-1}}=0}^p
\sum\limits_{d=1}^{2^{k-2r}}(-1)^{d-1}
\left(\hat C^{(d)}_{j_k\ldots j_1}\biggl|_{j_{g_1}=j_{g_2},\ldots, j_{g_{2r-1}}=j_{g_{2r}}}-
\right.
$$

$$
\left.
-~\bar C^{(d)}_{j_k\ldots j_1}\biggl|_{j_{g_1}=j_{g_2},\ldots, j_{g_{2r-1}}=j_{g_{2r}}}\right)=
$$

$$
=\sum\limits_{d=1}^{2^{k-2r}}(-1)^{d-1} \lim\limits_{p\to\infty}
\sum\limits_{j_{g_1}, j_{g_3},\ldots ,j_{g_{2r-1}}=0}^p
\left(\hat C^{(d)}_{j_k\ldots j_1}\biggl|_{j_{g_1}=j_{g_2},\ldots, j_{g_{2r-1}}=j_{g_{2r}}}-
\right.
$$

\begin{equation}
\label{july100010}
\left.
-~\bar C^{(d)}_{j_k\ldots j_1}\biggl|_{j_{g_1}=j_{g_2},\ldots, j_{g_{2r-1}}=j_{g_{2r}}}\right).
\end{equation}

\vspace{3mm}

Further, consider 3 possible cases.

{\bf Case~1.}\ The quantities
\begin{equation}
\label{july100010x}
~~~~~~~~~~~~\hat C^{(d)}_{j_k\ldots j_1}\biggl|_{j_{g_1}=j_{g_2},\ldots, j_{g_{2r-1}}=j_{g_{2r}}},\ \ \ 
\bar C^{(d)}_{j_k\ldots j_1}\biggl|_{j_{g_1}=j_{g_2},\ldots, j_{g_{2r-1}}=j_{g_{2r}}}
\end{equation}

\noindent
are such that
\begin{equation}
\label{july100010y}
\prod\limits_{l=1}^r {\bf 1}_{\{g_{2l}=g_{2l-1}+1\}}=1
\end{equation}

\noindent
for $d=1,2,\ldots, 2^{k-2r}$ and 
\begin{equation}
\label{july100010z}
C_{j_k\ldots j_1}\biggl|_{j_{g_1}=j_{g_2},\ldots, j_{g_{2r-1}}=j_{g_{2r}}}
\end{equation}

\noindent
is such that the condition (\ref{july100010y}) is fulfilled for
(\ref{july100010z}).

{\bf Case~2.}\ The quantities (\ref{july100010x}) 
are such that the condition (\ref{july100010y})
is satisfied for $d=1,2,\ldots, 2^{k-2r}$ and 
(\ref{july100010z}) is such that the condition 
\begin{equation}
\label{july100010v}
\prod\limits_{l=1}^r {\bf 1}_{\{g_{2l}=g_{2l-1}+1\}}=0
\end{equation}

\noindent
is fulfilled for (\ref{july100010z}).

{\bf Case~3.}\ The quantities (\ref{july100010x}) 
are such that the condition (\ref{july100010v}) 
is satisfied for $d=1,2,\ldots, 2^{k-2r}$ and 
(\ref{july100010z}) is such that the condition
(\ref{july100010v}) 
is fulfilled for (\ref{july100010z}).

For Case~1, applying 
(\ref{july100000}) for the case $k=2r$ and (\ref{july100010}),
we get
for any fixed $j_1,\ldots,j_q,\ldots,j_k$
$(q\ne g_1, g_2, \ldots, g_{2r-1},$ $g_{2r})$

\vspace{-2mm}
$$
\lim\limits_{p\to\infty}
\sum\limits_{j_{g_1}, j_{g_3},\ldots ,j_{g_{2r-1}}=0}^p
C_{j_k\ldots j_1}\biggl|_{j_{g_1}=j_{g_2},\ldots, j_{g_{2r-1}}=j_{g_{2r}}}=
$$

\vspace{3mm}
$$
=\sum\limits_{d=1}^{2^{k-2r}}(-1)^{d-1} \lim\limits_{p\to\infty}
\sum\limits_{j_{g_1}, j_{g_3},\ldots ,j_{g_{2r-1}}=0}^p
\left(\hat C^{(d)}_{j_k\ldots j_1}\biggl|_{j_{g_1}=j_{g_2},\ldots, j_{g_{2r-1}}=j_{g_{2r}}}-
\right.
$$

\vspace{1mm}
$$
\left.
-~\bar C^{(d)}_{j_k\ldots j_1}\biggl|_{j_{g_1}=j_{g_2},\ldots, j_{g_{2r-1}}=j_{g_{2r}}}\right)=
$$

\vspace{4mm}
$$
=
\sum\limits_{d=1}^{2^{k-2r}}(-1)^{d-1}\frac{1}{2^r} 
\prod\limits_{l=1}^r {\bf 1}_{\{g_{2l}=g_{2l-1}+1\}}\times
$$

\vspace{-1mm}
$$
\times
\left(\hat C^{(d)}_{j_k \ldots j_1}\biggl|_{(j_{g_2} j_{g_1})\curvearrowright (\cdot)
\ldots (j_{g_{2r}} j_{g_{2r-1}})\curvearrowright (\cdot),
j_{g_{{}_{1}}}=~j_{g_{{}_{2}}},\ldots, j_{g_{{}_{2r-1}}}=~j_{g_{{}_{2r}}}
}\biggr.-\right.
$$

\vspace{-1mm}
\begin{equation}
\label{july100006abc}
~~~~~~~~~~\left.-~\bar C^{(d)}_{j_k \ldots j_1}\biggl|_{(j_{g_2} j_{g_1})\curvearrowright (\cdot)
\ldots (j_{g_{2r}} j_{g_{2r-1}})\curvearrowright (\cdot),
j_{g_{{}_{1}}}=~j_{g_{{}_{2}}},\ldots, j_{g_{{}_{2r-1}}}=~j_{g_{{}_{2r}}}
}\biggr.
\right)=
\end{equation}

\vspace{0.5mm}
$$
=
\sum\limits_{d=1}^{2^{k-2r}}(-1)^{d-1}\frac{1}{2^r} 
\left(\hat C^{(d)}_{j_k \ldots j_1}\biggl|_{(j_{g_2} j_{g_1})\curvearrowright (\cdot)
\ldots (j_{g_{2r}} j_{g_{2r-1}})\curvearrowright (\cdot),
j_{g_{{}_{1}}}=~j_{g_{{}_{2}}},\ldots, j_{g_{{}_{2r-1}}}=~j_{g_{{}_{2r}}}
}\biggr.-\right.
$$

\vspace{-3mm}
\begin{equation}
\label{july100006}
~~~~~~~~\left.-~\bar C^{(d)}_{j_k \ldots j_1}\biggl|_{(j_{g_2} j_{g_1})\curvearrowright (\cdot)
\ldots (j_{g_{2r}} j_{g_{2r-1}})\curvearrowright (\cdot),
j_{g_{{}_{1}}}=~j_{g_{{}_{2}}},\ldots, j_{g_{{}_{2r-1}}}=~j_{g_{{}_{2r}}}
}\biggr.
\right),
\end{equation}

\vspace{5.5mm}
\noindent
where $g_1,g_2,\ldots,g_{2r-1},g_{2r}$ as in (\ref{leto5007after}),
$k>2r,$ $r=1,2,\ldots,[k/2].$

It is not difficult to see that 
the left-hand side of (\ref{july100010y}) is a constant for the
quantities (\ref{july100010x}) for all $d=1,2,\ldots, 2^{k-2r}$.

Using (\ref{july100003}), (\ref{july100004}), we obtain 

\vspace{-3mm}
$$
\sum\limits_{d=1}^{2^{k-2r}}(-1)^{d-1}\frac{1}{2^r} 
\left(\hat C^{(d)}_{j_k \ldots j_1}\biggl|_{(j_{g_2} j_{g_1})\curvearrowright (\cdot)
\ldots (j_{g_{2r}} j_{g_{2r-1}})\curvearrowright (\cdot),
j_{g_{{}_{1}}}=~j_{g_{{}_{2}}},\ldots, j_{g_{{}_{2r-1}}}=~j_{g_{{}_{2r}}}
}\biggr.-\right.
$$

\vspace{-1mm}
$$
\left.-~\bar C^{(d)}_{j_k \ldots j_1}\biggl|_{(j_{g_2} j_{g_1})\curvearrowright (\cdot)
\ldots (j_{g_{2r}} j_{g_{2r-1}})\curvearrowright (\cdot),
j_{g_{{}_{1}}}=~j_{g_{{}_{2}}},\ldots, j_{g_{{}_{2r-1}}}=~j_{g_{{}_{2r}}}
}\biggr.
\right)=
$$

\begin{equation}
\label{july100007}
~~~~~~~~=\frac{1}{2^r} C_{j_k \ldots j_1}\biggl|_{(j_{g_2} j_{g_1})\curvearrowright (\cdot)
\ldots (j_{g_{2r}} j_{g_{2r-1}})\curvearrowright (\cdot),
j_{g_{{}_{1}}}=~j_{g_{{}_{2}}},\ldots, j_{g_{{}_{2r-1}}}=~j_{g_{{}_{2r}}}
}\biggr..
\end{equation}

\vspace{6mm}

Combining (\ref{july100006}) and (\ref{july100007}), we have
for any fixed $j_1,\ldots,j_q,\ldots,j_k$
$(q\ne g_1, g_2, \ldots, g_{2r-1},$ $g_{2r})$

\vspace{-2mm}
$$
\lim\limits_{p\to\infty}
\sum\limits_{j_{g_1}, j_{g_3},\ldots ,j_{g_{2r-1}}=0}^p
C_{j_k\ldots j_1}\biggl|_{j_{g_1}=j_{g_2},\ldots, j_{g_{2r-1}}=j_{g_{2r}}}=
$$

\vspace{-1mm}
\begin{equation}
\label{july100008}
~~~~~~~~~~~=\frac{1}{2^r}
C_{j_k \ldots j_1}\biggl|_{(j_{g_2} j_{g_1})\curvearrowright (\cdot)
\ldots (j_{g_{2r}} j_{g_{2r-1}})\curvearrowright (\cdot),
j_{g_{{}_{1}}}=~j_{g_{{}_{2}}},\ldots, j_{g_{{}_{2r-1}}}=~j_{g_{{}_{2r}}}
}\biggr.,
\end{equation}

\vspace{4mm}
\noindent
where $g_1,g_2,\ldots,g_{2r-1},g_{2r}$ as in (\ref{leto5007after}),
$k>2r,$ $r=1,2,\ldots,[k/2].$

From (\ref{july100000}) for the case $k=2r$ and (\ref{july100008}) ($k>2r$) we obtain 
(\ref{july100000}) for the case $k\ge 2r$. The equality
(\ref{july100000}) is proved for Case~1.

For Case~2, applying 
(\ref{july100000}) for the case $k=2r$ and (\ref{july100010}),
we get (\ref{july100006})
for any fixed $j_1,\ldots,j_q,\ldots,j_k$
$(q\ne g_1, g_2, \ldots, g_{2r-1},$ $g_{2r}).$
Further, note that

\vspace{-4mm}
$$
\hat C^{(d)}_{j_k \ldots j_1}\biggl|_{(j_{g_2} j_{g_1})\curvearrowright (\cdot)
\ldots (j_{g_{2r}} j_{g_{2r-1}})\curvearrowright (\cdot),
j_{g_{{}_{1}}}=~j_{g_{{}_{2}}},\ldots, j_{g_{{}_{2r-1}}}=~j_{g_{{}_{2r}}}
}\biggr.=
$$

\begin{equation}
\label{july100008xyz}
~~~~~~~~=~\bar C^{(d)}_{j_k \ldots j_1}\biggl|_{(j_{g_2} j_{g_1})\curvearrowright (\cdot)
\ldots (j_{g_{2r}} j_{g_{2r-1}})\curvearrowright (\cdot),
j_{g_{{}_{1}}}=~j_{g_{{}_{2}}},\ldots, j_{g_{{}_{2r-1}}}=~j_{g_{{}_{2r}}}
}\biggr.
\end{equation}

\vspace{4mm}
\noindent
for Case~2. Combining (\ref{july100006}) and (\ref{july100008xyz}), we
obtain (Case~2)
for any fixed $j_1,\ldots,j_q,\ldots,j_k$
$(q\ne g_1, g_2, \ldots, g_{2r-1},$ $g_{2r})$

\vspace{-1mm}
\begin{equation}
\label{july100008xyz1}
~~~~~~~~\lim\limits_{p\to\infty}
\sum\limits_{j_{g_1}, j_{g_3},\ldots ,j_{g_{2r-1}}=0}^p
C_{j_k\ldots j_1}\biggl|_{j_{g_1}=j_{g_2},\ldots, j_{g_{2r-1}}=j_{g_{2r}}}=0.
\end{equation}

\vspace{3mm}

From (\ref{july100000}) for the case $k=2r$ and (\ref{july100008xyz1}) $(k>2r)$ we obtain 
(\ref{july100008xyz1}) for the case $k\ge 2r$. The equality
(\ref{july100000}) is proved for Case~2.

For Case~3, applying 
(\ref{july100000}) for the case $k=2r$ and (\ref{july100010}),
we get (\ref{july100006abc})
for any fixed $j_1,\ldots,j_q,\ldots,j_k$
$(q\ne g_1, g_2, \ldots, g_{2r-1},$ $g_{2r}).$

Since 

\vspace{-6mm}
\begin{equation}
\label{july400000}
\prod\limits_{l=1}^r {\bf 1}_{\{g_{2l}=g_{2l-1}+1\}}=0
\end{equation}

\vspace{3mm}
\noindent
for Case~3, then from (\ref{july100006abc})
we get (\ref{july100008xyz1}) for $k>2r$
(recall that the left-hand side of (\ref{july400000}) is a constant for the
quantities (\ref{july100010x}) for all $d=1,2,\ldots, 2^{k-2r}$).

From (\ref{july100000}) for $k=2r$ and (\ref{july100008xyz1}) for $k>2r$ (Case~3) we obtain 
(\ref{july100008xyz1}) for $k\ge 2r$ (Case~3). The equality
(\ref{july100000}) is proved for Case~3. 
The equality
(\ref{july100000}) is proved. 
Thus, Theorem~2.58 is proved. Theorem~2.59 is also proved.

In conclusion of this section, we will make a remark
about the condition (\ref{july700000}). It would seem
that according to (\ref{july100000}), we can write

$$
\lim\limits_{p\to\infty}
\sum\limits_{\stackrel{j_1,\ldots,j_q,\ldots,j_k=0}{{}_{q\ne g_1, g_2, \ldots, g_{2r-1},
g_{2r}}}}^p
\Biggl(\sum\limits_{j_{g_1}, j_{g_3},\ldots ,j_{g_{2r-1}}=0}^p
C_{j_k\ldots j_1}\biggl|_{j_{g_1}=j_{g_2},\ldots, j_{g_{2r-1}}=j_{g_{2r}}}-\Biggr.
$$

\vspace{-2mm}
$$
\Biggl.-\frac{1}{2^r} \prod\limits_{l=1}^r {\bf 1}_{\{g_{2l}=g_{2l-1}+1\}}
C_{j_k \ldots j_1}\biggl|_{(j_{g_2} j_{g_1})\curvearrowright (\cdot)
\ldots (j_{g_{2r}} j_{g_{2r-1}})\curvearrowright (\cdot),
j_{g_{{}_{1}}}=~j_{g_{{}_{2}}},\ldots, j_{g_{{}_{2r-1}}}=~j_{g_{{}_{2r}}}
}\biggr.\Biggr)^2=
$$

\vspace{5mm}
$$
=\sum\limits_{\stackrel{j_1,\ldots,j_q,\ldots,j_k=0}{{}_{q\ne g_1, g_2, \ldots, g_{2r-1},
g_{2r}}}}^{\infty}
\Biggl(\sum\limits_{j_{g_1}, j_{g_3},\ldots ,j_{g_{2r-1}}=0}^{\infty}
C_{j_k\ldots j_1}\biggl|_{j_{g_1}=j_{g_2},\ldots, j_{g_{2r-1}}=j_{g_{2r}}}-\Biggr.
$$
$$
\Biggl.-\frac{1}{2^r} \prod\limits_{l=1}^r {\bf 1}_{\{g_{2l}=g_{2l-1}+1\}}
C_{j_k \ldots j_1}\biggl|_{(j_{g_2} j_{g_1})\curvearrowright (\cdot)
\ldots (j_{g_{2r}} j_{g_{2r-1}})\curvearrowright (\cdot),
j_{g_{{}_{1}}}=~j_{g_{{}_{2}}},\ldots, j_{g_{{}_{2r-1}}}=~j_{g_{{}_{2r}}}
}\biggr.\Biggr)^2=
$$

\vspace{4mm}
$$
=\sum\limits_{\stackrel{j_1,\ldots,j_q,\ldots,j_k=0}{{}_{q\ne g_1, g_2, \ldots, g_{2r-1},
g_{2r}}}}^{\infty}
\Bigl( 0 \Bigr)^2= 0
$$

\vspace{3mm}
\noindent 
for all $r=1,2,\ldots,[k/2]$ and 
for all possible $g_1,g_2,\ldots,g_{2r-1},g_{2r}$ (see {\rm (\ref{leto5007after})),
where notations are the same as in 
(\ref{july700000}).

However, the above argument contains an error associated
with the replacement of the limit with an iterated one.
Let us consider this observation in more detail
using an example.

To begin, let us recall that the sum of an infinite number series
is defined as the limit of the partial sums of this series, i.e.
$$
\lim\limits_{n\to\infty}
\sum\limits_{i=1}^{n}a_i
\stackrel{\sf def}{=}
\sum\limits_{i=1}^{\infty}a_i.
$$

Let $k=3,$ $r=1$ and $g_1=1, g_2=3.$ Further, we have

\vspace{-1mm}
$$
\lim\limits_{p\to\infty}\sum\limits_{j_2=0}^p
\Biggl(
\sum\limits_{j_1=0}^{p}C_{j_1 j_2 j_1}\Biggr)^2=
\lim\limits_{p\to\infty}\sum\limits_{j_2=0}^p
\sum\limits_{j_1=0}^{p}C_{j_1 j_2 j_1}\sum\limits_{j_3=0}^{p}C_{j_3 j_2 j_3}=
$$

\vspace{1mm}

$$
=
\lim\limits_{p\to\infty}\sum\limits_{j_3, j_2, j_1=0}^p
C_{j_1 j_2 j_1}C_{j_3 j_2 j_3}\stackrel{\sf def}{=}
$$

\vspace{1mm}
\begin{equation}
\label{2025july2025}
\stackrel{\sf def}{=}\sum\limits_{j_3, j_2, j_1=0}^{\infty}
C_{j_1 j_2 j_1}C_{j_3 j_2 j_3},
\end{equation}

\vspace{4mm}
$$
\lim\limits_{q\to\infty}\lim\limits_{p\to\infty}\sum\limits_{j_2=0}^q
\Biggl(
\sum\limits_{j_1=0}^{p}C_{j_1 j_2 j_1}\Biggr)^2=
\lim\limits_{q\to\infty}\sum\limits_{j_2=0}^q
\lim\limits_{p\to\infty}\Biggl(
\sum\limits_{j_1=0}^{p}C_{j_1 j_2 j_1}\Biggr)^2=
$$

\vspace{1mm}
$$
=
\lim\limits_{q\to\infty}\sum\limits_{j_2=0}^q
\Biggl(
\sum\limits_{j_1=0}^{\infty}C_{j_1 j_2 j_1}\Biggr)^2=
$$
$$
=
\sum\limits_{j_2=0}^{\infty}
\Biggl(
\sum\limits_{j_1=0}^{\infty}C_{j_1 j_2 j_1}\Biggr)^2=
\sum\limits_{j_2=0}^{\infty}
\sum\limits_{j_1=0}^{\infty}C_{j_1 j_2 j_1}
\sum\limits_{j_3=0}^{\infty}C_{j_3 j_2 j_3}=
$$

\vspace{1mm}
\begin{equation}
\label{2025july2025a}
=
\sum\limits_{j_2=0}^{\infty}
\sum\limits_{j_1=0}^{\infty}
\sum\limits_{j_3=0}^{\infty}C_{j_1 j_2 j_1}C_{j_3 j_2 j_3}.
\end{equation}

\vspace{4mm}

It is obvious that the right-hand sides of equalities
(\ref{2025july2025}) and (\ref{2025july2025a})
are generally not equal. The equality 
of the mentioned expressions requires separate proof.

In the next section, we will consider a fairly
efficient approach to proving the equality 
(\ref{july700000}).

\section{Expansion of Iterated Stratonovich Stochastic Integrals
of Arbitrary Multiplicity $k$ $(k\in{\bf N})$. 
The Case of an Ar\-bit\-ra\-ry Complete Orthonormal System of 
Functions in $L_2([t,T]),$ $\psi_1(\tau),\ldots, \psi_k(\tau)
\in L_2([t,T]).$
Proof of Hypotheses~2.4, 2.5 for the Case $p_1=\ldots=p_k=p$
and Under the Condition (\ref{09091})}

We will start this section with an example. Let us assume that
$h_1(\tau),\ldots,$ $h_{12}(\tau)\in L_2([t, T])$ and consider the
following integral
$$
I\stackrel{\sf def}{=}\int\limits_t^T h_{12}(t_{12})
\int\limits_t^{t_{12}}h_{11}(t_{11})\ldots 
\int\limits_t^{t_{2}}h_{1}(t_{1})dt_1 \ldots dt_{11} dt_{12}.
$$

We want to transform the integral $I$ in such a way that
$$
I=\int\limits_t^T h_{10}(t_{10})
\int\limits_t^{t_{10}}h_{6}(t_{6})
\int\limits_t^{t_{6}}h_{4}(t_{4})
\int\limits_t^{t_{4}}h_{3}(t_{3})
\left(\ldots\right)
dt_3 dt_{4} dt_6 dt_{10},
$$
where $\left(\ldots\right)$ is some expression.

Using Fubini's Theorem, we obtain
$$
I=\int\limits_t^T h_{12}(t_{12})
\int\limits_t^{t_{12}}h_{11}(t_{11})
\int\limits_t^{t_{11}}\underline{h_{10}(t_{10})}
\int\limits_t^{t_{10}}h_{9}(t_{9})
\int\limits_t^{t_{9}}h_{8}(t_{8})
\int\limits_t^{t_{8}}h_{7}(t_{7})
\int\limits_t^{t_{7}}\underline{h_{6}(t_{6})}\times
$$
$$
\times
\int\limits_t^{t_{6}}h_{5}(t_{5})
\int\limits_t^{t_{5}}\underline{h_{4}(t_{4})}
\int\limits_t^{t_{4}}\underline{h_{3}(t_{3})}
\int\limits_t^{t_{3}}h_{2}(t_{2})
\int\limits_t^{t_{2}}h_{1}(t_{1})dt_1 dt_2 dt_3 dt_4
dt_5 dt_6 dt_7 dt_8 \times
$$

\vspace{-2mm}
$$
\times
dt_9 dt_{10} dt_{11} dt_{12}=
$$

\vspace{-2mm}
$$
=\int\limits_t^T
\underline{h_{10}(t_{10})}
\int\limits_t^{t_{10}}h_{9}(t_{9})
\int\limits_t^{t_{9}}h_{8}(t_{8})
\int\limits_t^{t_{8}}h_{7}(t_{7})
\int\limits_t^{t_{7}}\underline{h_{6}(t_{6})}\int\limits_t^{t_{6}}h_{5}(t_{5})
\times
$$

\vspace{-2mm}
$$
\times
\int\limits_t^{t_{5}}\underline{h_{4}(t_{4})}
\int\limits_t^{t_{4}}\underline{h_{3}(t_{3})}
\int\limits_t^{t_{3}}h_{2}(t_{2})
\int\limits_t^{t_{2}}h_{1}(t_{1})dt_1 dt_2 dt_3 dt_4
dt_5 dt_6 dt_7 dt_8 dt_9\times
$$

\vspace{-2mm}
$$
\times
\left(\int\limits_{t_{10}}^T h_{11}(t_{11})
\int\limits_{t_{11}}^T h_{12}(t_{12})
dt_{12} dt_{11}\right) dt_{10}=
$$

\vspace{-2mm}
$$
=\int\limits_t^T
\underline{h_{10}(t_{10})}
\int\limits_t^{t_{10}}\underline{h_{6}(t_{6})}\int\limits_t^{t_{6}}h_{5}(t_{5})
\int\limits_t^{t_{5}}\underline{h_{4}(t_{4})}
\int\limits_t^{t_{4}}\underline{h_{3}(t_{3})}
\int\limits_t^{t_{3}}h_{2}(t_{2})
\int\limits_t^{t_{2}}h_{1}(t_{1})\times
$$

\vspace{-1mm}
$$
\times
dt_1 dt_2 dt_3 dt_4 dt_5 \left(\int\limits_{t_6}^{t_{10}}h_7(t_7)
\int\limits_{t_7}^{t_{10}}h_8(t_8)
\int\limits_{t_8}^{t_{10}}h_9(t_9)dt_9 dt_8 dt_7
\right) dt_6 \times
$$

\vspace{1mm}
$$
\times
\left(\int\limits_{t_{10}}^T h_{11}(t_{11})
\int\limits_{t_{11}}^T h_{12}(t_{12})
dt_{12} dt_{11}\right) dt_{10}=
$$

\vspace{-2mm}
$$
=\int\limits_t^T
\underline{h_{10}(t_{10})}
\int\limits_t^{t_{10}}\underline{h_{6}(t_{6})}\int\limits_t^{t_{6}}
\underline{h_{4}(t_{4})}
\int\limits_t^{t_{4}}\underline{h_{3}(t_{3})}
\left(\int\limits_t^{t_{3}}h_{2}(t_{2})
\int\limits_t^{t_{2}}h_{1}(t_{1})
dt_1 dt_2\right) dt_3 \times
$$

\vspace{-2mm}
$$
\times \left(\int\limits_{t_4}^{t_6} h_5(t_5)dt_5 \right)dt_4 
\left(\int\limits_{t_6}^{t_{10}}h_7(t_7)
\int\limits_{t_7}^{t_{10}}h_8(t_8)
\int\limits_{t_8}^{t_{10}}h_9(t_9)dt_9 dt_8 dt_7
\right) dt_6 \times
$$

\vspace{-1mm}
$$
\times
\left(\int\limits_{t_{10}}^T h_{11}(t_{11})
\int\limits_{t_{11}}^T h_{12}(t_{12})
dt_{12} dt_{11}\right) dt_{10}=
$$
$$
=\int\limits_t^T
\underline{h_{10}(t_{10})}
\int\limits_t^{t_{10}}\underline{h_{6}(t_{6})}\int\limits_t^{t_{6}}
\underline{h_{4}(t_{4})}
\int\limits_t^{t_{4}}\underline{h_{3}(t_{3})}
\left(\int\limits_t^{t_{3}}h_{2}(t_{2})
\int\limits_t^{t_{2}}h_{1}(t_{1})
dt_1 dt_2\right) \times
$$

$$
\times \left(\int\limits_{t_4}^{t_6} h_5(t_5)dt_5 \right)
\left(\int\limits_{t_6}^{t_{10}}h_9(t_9)
\int\limits_{t_6}^{t_{9}}h_8(t_8)
\int\limits_{t_6}^{t_{8}}h_7(t_7)dt_7 dt_8 dt_9
\right)\times
$$

\vspace{-2mm}
\begin{equation}
\label{copa1}
~~~~~~~~~~\times
\left(\int\limits_{t_{10}}^T h_{12}(t_{12})
\int\limits_{t_{10}}^{t_{12}} h_{11}(t_{11})
dt_{11} dt_{12}\right) dt_3 dt_4 dt_6 dt_{10}.
\end{equation}

Further, suppose that $h_l(\tau)=\psi_l(\tau)\phi_{j_l}(\tau)$ $(l=1,\ldots,12)$ in 
(\ref{copa1}) (here $\left\{\phi_j(x)\right\}_{j=0}^{\infty}$
is an arbitrary complete orthonormal system of 
functions in the space $L_2([t,T])$ and
$\psi_1(\tau),\ldots ,\psi_{12}(\tau)\in $ $L_2([t, T])$).
Thus, we get
$$
C_{j_{12} j_{11} \underline{j_{10}} j_9 j_8 j_7 \underline{j_6} j_5 \underline{j_4} 
\underline{j_3} j_2 j_1}=
\int\limits_t^T
\psi_{10}(t_{10})\phi_{j_{10}}(t_{10})
\int\limits_t^{t_{10}}
\psi_{6}(t_{6})\phi_{j_{6}}(t_{6})
\int\limits_t^{t_6}
\psi_{4}(t_{4})\phi_{j_{4}}(t_{4})\times
$$
$$
\times
\int\limits_t^{t_4}
\psi_{3}(t_{3})\phi_{j_{3}}(t_{3})
C_{j_{12}j_{11}}^{\psi_{12}\psi_{11}}(T,t_{10})
C_{j_{9}j_{8}j_7}^{\psi_{9}\psi_{8}\psi_7}(t_{10},t_6)
C_{j_{5}}^{\psi_{5}}(t_6,t_{4})
C_{j_{2}j_{1}}^{\psi_{2}\psi_{1}}(t_{3},t)\times
$$

\vspace{-2mm}
\begin{equation}
\label{copa1a}
\times
dt_3 dt_4 dt_6 dt_{10},
\end{equation}

\vspace{2mm}
\noindent
where
$$
C_{j_m \ldots j_l}^{\psi_m\ldots \psi_l}(s,\tau)=\int\limits_{\tau}^s
\psi_m(t_m)\phi_{j_m}(t_m)\ldots
\int\limits_{\tau}^{t_{l+1}}
\psi_l(t_l)\phi_{j_l}(t_l)dt_l\ldots dt_m,
$$

\noindent
where $t\le\tau<s\le T$, $m\ge l$ $(m, l\in {\bf N}).$

Suppose that 
$g_1,g_2,\ldots,g_{2r-1},g_{2r}$ as in (\ref{leto5007after}) and $k>2r,$ $r\ge 1$
(the case $k=2r$ see in Sect.~2.27.4).
Consider $d_1,e_1,$ $\ldots,d_f,e_f, f \in{\bf N}$
such that 
$$
1\le d_1-e_1+1<\ldots < d_1-1 < d_1 <\ldots < d_f-e_f+1 <\ldots < d_f-1 < d_f \le k,
$$

\vspace{-6mm}
$$
\left\{g_1,g_2,\ldots,g_{2r-1},g_{2r}\right\}=
$$

\vspace{-9mm}
$$
=\left\{d_1-e_1+1,\ldots, d_1-1, d_1\right\}\cup\ldots
\cup  \left\{d_f-e_f+1,\ldots, d_f-1, d_f\right\},
$$
$$
e_1+e_2+\ldots+e_f=2r,\ \ \{1,\ldots, k\}~\backslash \left\{g_1,g_2,\ldots,g_{2r-1},g_{2r}\right\}=
\left\{q_1,\ldots,q_{k-2r}\right\}.
$$

\vspace{2mm}

{\it We will say that the condition $(A)$ is satisfied if~~$\forall$   
$\left\{g_{2l-1},g_{2l}\right\}$ $(l=1,\ldots,r)$ $\exists$ $h\in \left\{1,\ldots,f\right\}$
such that
\begin{equation}
\label{copa2}
\left\{g_{2l-1},g_{2l}\right\}\subset
\left\{d_h-e_h+1,\ldots, d_h-1, d_h\right\}.
\end{equation}
Moreover$,$ $\forall$ $h\in \left\{1,\ldots,f\right\}$
$\exists$ $\left\{g_{2l-1},g_{2l}\right\}$ $(l=1,\ldots,r)$
such that  {\rm (\ref{copa2})} is fulfilled.}

\vspace{2mm}

If the condition $(A)$ is satisfied, then $e_1, \ldots, e_f$ are even and we can write

$$
\left\{d_1-e_1+1,\ldots, d_1\right\}=
\left\{g_1^{(1)},g_2^{(1)},\ldots,g_{2r_1-1}^{(1)},g_{2r_1}^{(1)}\right\},
$$
$$
\ldots
$$

\vspace{-5mm}
$$
\left\{d_f-e_f+1,\ldots, d_f\right\}=
\left\{g_1^{(f)},g_2^{(f)},\ldots,g_{2r_f-1}^{(f)},g_{2r_f}^{(f)}\right\},
$$

\vspace{-1mm}
$$
\bigl\{g_1,g_2,\ldots,g_{2r-1},g_{2r}\bigr\}=
$$

\vspace{-2mm}
$$
=
\left\{g_1^{(1)},g_2^{(1)},\ldots,g_{2r_1-1}^{(1)},g_{2r_1}^{(1)},\ldots, 
g_1^{(f)},g_2^{(f)},\ldots,g_{2r_f-1}^{(f)},g_{2r_f}^{(f)}\right\}.
$$

\vspace{3mm}

\noindent
If the condition $(A)$ is not fulfilled, then some of $e_1, \ldots, e_f$ can be uneven.

Using (\ref{july90000}) and a modification of the algorithm from Sect.~2.27.4 
(see below for details) it can be proved that
$$
\lim\limits_{p\to\infty}
\sum\limits_{j_{g_1}, j_{g_3},\ldots ,j_{g_{2r-1}}=0}^p
\left(
C_{j_{d_f}\ldots j_{d_f-e_f+1}}^{\psi_{d_f}\ldots \psi_{d_f-e_f+1}}(t_{d_f+1},t_{d_f-e_f})
\ldots \right.
$$

$$
\left.
\ldots C_{j_{d_1}\ldots j_{d_1-e_1+1}}^{\psi_{d_1}\ldots \psi_{d_1-e_1+1}}(t_{d_1+1},t_{d_1-e_1})
\right)\biggl|_{j_{g_1}=j_{g_2},\ldots, j_{g_{2r-1}}=j_{g_{2r}}}= 
$$

\vspace{4mm}
$$
=\prod\limits_{h=1}^f
\frac{1}{2^{r_h}}\prod\limits_{l=1}^{r_h}
{\bf 1}_{\{g_{2l}^{(h)}=g_{2l-1}^{(h)}+1\}}\times
$$

\begin{equation}
\label{copa3}
\times C_{j_{d_h}\ldots j_{d_h-e_h+1}}^{\psi_{d_h}\ldots \psi_{d_h-e_h+1}}(t_{d_h+1},t_{d_h-e_h})
\biggl|_{(j_{g_2^{(h)}} j_{g_1^{(h)}})\curvearrowright (\cdot)
\ldots (j_{g_{2r_h}^{(h)}} j_{g_{2r_h-1}^{(h)}})\curvearrowright (\cdot),
j_{g_{{}_{1}}^{(h)}}=j_{g_{{}_{2}}^{(h)}},\ldots, j_{g_{{}_{2r_h-1}}^{(h)}}=j_{g_{{}_{2r_h}}^{(h)}}
}\biggr.
\end{equation}

\vspace{3mm}
\noindent
if the condition $(A)$ is satisfied, and 

\newpage
\noindent
$$
\lim\limits_{p\to\infty}
\sum\limits_{j_{g_1}, j_{g_3},\ldots ,j_{g_{2r-1}}=0}^p
\left(
C_{j_{d_f}\ldots j_{d_f-e_f+1}}^{\psi_{d_f}\ldots \psi_{d_f-e_f+1}}(t_{d_f+1},t_{d_f-e_f})
\ldots \right.
$$

\vspace{-1mm}
\begin{equation}
\label{copa4}
~~~~~~~~~\left.
\ldots C_{j_{d_1}\ldots j_{d_1-e_1+1}}^{\psi_{d_1}\ldots \psi_{d_1-e_1+1}}(t_{d_1+1},t_{d_1-e_1})
\right)\biggl|_{j_{g_1}=j_{g_2},\ldots, j_{g_{2r-1}}=j_{g_{2r}}}=0
\end{equation}

\vspace{3.5mm}
\noindent
if the condition $(A)$ is not fulfilled, where $t_{k+1}\stackrel{\sf def}{=}T,$
$t_{0}\stackrel{\sf def}{=}t,$
$e_1+\ldots+e_f=2r$ in (\ref{copa3}), (\ref{copa4})
and
$e_h=2 r_h$ $(h=1,\ldots, f),$ $r_1+\ldots+r_f=r$ in (\ref{copa3}). 

Note that the series on the left-hand sides of 
(\ref{copa3}) and (\ref{copa4}) converge absolutly 
since
their sums do not depend 
on permutations of basis functions
(here the basis in $L_2([t,T]^{r})$ is
$\left\{\phi_{j_1}(x_1)\ldots \phi_{j_r}(x_r)\right\}_{j_1,\ldots,j_r=0}^{\infty}$).
Recall that any permutation of basis functions in a Hilbert space forms a basis 
in this Hilbert space \cite{gohb}.

Let us prove the formulas (\ref{copa3}) and (\ref{copa4}).

1. Suppose that the condition $(A)$ is satisfied and 
\begin{equation}
\label{copa5}
\prod\limits_{l=1}^{r_h}
{\bf 1}_{\{g_{2l}^{(h)}=g_{2l-1}^{(h)}+1\}}=1
\end{equation}

\noindent
for all $h=1,\ldots,f.$ In this case we can use the results from Sect.~2.27.4.
We have (see (\ref{july90000}))

\vspace{-1.5mm}
$$
\lim\limits_{p\to\infty}
\sum\limits_{j_{g_1}, j_{g_3},\ldots ,j_{g_{2r-1}}=0}^p
\left(
C_{j_{d_f}\ldots j_{d_f-e_f+1}}^{\psi_{d_f}\ldots \psi_{d_f-e_f+1}}(t_{d_f+1},t_{d_f-e_f})
\ldots \right.
$$

\vspace{2.5mm}
$$
\left.
\ldots C_{j_{d_1}\ldots j_{d_1-e_1+1}}^{\psi_{d_1}\ldots \psi_{d_1-e_1+1}}(t_{d_1+1},t_{d_1-e_1})
\right)\biggl|_{j_{g_1}=j_{g_2},\ldots, j_{g_{2r-1}}=j_{g_{2r}}}= 
$$

\vspace{3mm}
$$
=\lim\limits_{p\to\infty}
\sum\limits_{j_{g_{{}_{1}}^{(1)}}, j_{g_{{}_{3}}^{(1)}},\ldots, j_{g_{{}_{2r_1-1}}^{(1)}}=0}^p
C_{j_{d_1}\ldots j_{d_1-e_1+1}}^{\psi_{d_1}\ldots \psi_{d_1-e_1+1}}(t_{d_1+1},t_{d_1-e_1})
\biggl|_{
j_{g_{{}_{1}}^{(1)}}=j_{g_{{}_{2}}^{(1)}},\ldots, j_{g_{{}_{2r_1-1}}^{(1)}}=j_{g_{{}_{2r_1}}^{(1)}}
}\times
$$
\vspace{-10mm}

$$
\ldots
$$
\vspace{-10mm}

$$
\times 
\lim\limits_{p\to\infty}
\sum\limits_{j_{g_{{}_{1}}^{(f)}}, j_{g_{{}_{3}}^{(f)}}, \ldots , j_{g_{{}_{2r_f-1}}^{(f)}}=0}^p
C_{j_{d_f}\ldots j_{d_f-e_f+1}}^{\psi_{d_f}\ldots \psi_{d_f-e_f+1}}(t_{d_f+1},t_{d_f-e_f})
\biggl|_{
j_{g_{{}_{1}}^{(f)}}=j_{g_{{}_{2}}^{(f)}},\ldots, j_{g_{{}_{2r_f-1}}^{(f)}}=j_{g_{{}_{2r_f}}^{(f)}}
}=
$$
\newpage
\noindent
$$
=\prod\limits_{h=1}^f
\frac{1}{2^{r_h}}\prod\limits_{l=1}^{r_h}
{\bf 1}_{\{g_{2l}^{(h)}=g_{2l-1}^{(h)}+1\}}\times
$$

$$
\times C_{j_{d_h}\ldots j_{d_h-e_h+1}}^{\psi_{d_h}\ldots \psi_{d_h-e_h+1}}(t_{d_h+1},t_{d_h-e_h})
\biggl|_{(j_{g_2^{(h)}} j_{g_1^{(h)}})\curvearrowright (\cdot)
\ldots (j_{g_{2r_h}^{(h)}} j_{g_{2r_h-1}^{(h)}})\curvearrowright (\cdot),
j_{g_{{}_{1}}^{(h)}}=j_{g_{{}_{2}}^{(h)}},\ldots, j_{g_{{}_{2r_h-1}}^{(h)}}=j_{g_{{}_{2r_h}}^{(h)}}
}\biggr..
$$

\vspace{3mm}

\noindent
Thus, we get the formula (\ref{copa3}).

2. Suppose that the condition $(A)$ is satisfied and for some $h=1,\ldots,f$
\begin{equation}
\label{copa6}
\prod\limits_{l=1}^{r_h}
{\bf 1}_{\{g_{2l}^{(h)}=g_{2l-1}^{(h)}+1\}}=0.
\end{equation}

In this case, we act the same as in the previous case.
Applying (\ref{july90000}), we obtain

\vspace{-2mm}
$$
\lim\limits_{p\to\infty}
\sum\limits_{j_{g_1}, j_{g_3}, \ldots,  j_{g_{2r-1}}=0}^p
\left(
C_{j_{d_f}\ldots j_{d_f-e_f+1}}^{\psi_{d_f}\ldots \psi_{d_f-e_f+1}}(t_{d_f+1},t_{d_f-e_f})
\ldots \right.
$$

\vspace{2mm}
$$
\left.
\ldots C_{j_{d_1}\ldots j_{d_1-e_1+1}}^{\psi_{d_1}\ldots \psi_{d_1-e_1+1}}(t_{d_1+1},t_{d_1-e_1})
\right)\biggl|_{j_{g_1}=j_{g_2},\ldots, j_{g_{2r-1}}=j_{g_{2r}}}= 
$$

\vspace{3mm}
$$
=\lim\limits_{p\to\infty}
\sum\limits_{j_{g_{{}_{1}}^{(1)}}, j_{g_{{}_{3}}^{(1)}}, \ldots ,
j_{g_{{}_{2r_1-1}}^{(1)}}=0}^p
C_{j_{d_1}\ldots j_{d_1-e_1+1}}^{\psi_{d_1}\ldots \psi_{d_1-e_1+1}}(t_{d_1+1},t_{d_1-e_1})
\biggl|_{
j_{g_{{}_{1}}^{(1)}}=j_{g_{{}_{2}}^{(1)}},\ldots, j_{g_{{}_{2r_1-1}}^{(1)}}=j_{g_{{}_{2r_1}}^{(1)}}
}\times
$$
\vspace{-10mm}

$$
\ldots
$$
\vspace{-10mm}

\begin{equation}
\label{2024october1}
\times 
\lim\limits_{p\to\infty}
\sum\limits_{j_{g_{{}_{1}}^{(f)}}, j_{g_{{}_{3}}^{(f)}},\ldots ,
j_{g_{{}_{2r_f-1}}^{(f)}}=0}^p
C_{j_{d_f}\ldots j_{d_f-e_f+1}}^{\psi_{d_f}\ldots \psi_{d_f-e_f+1}}(t_{d_f+1},t_{d_f-e_f})
\biggl|_{
j_{g_{{}_{1}}^{(f)}}=j_{g_{{}_{2}}^{(f)}},\ldots, j_{g_{{}_{2r_f-1}}^{(f)}}=j_{g_{{}_{2r_f}}^{(f)}}
}=0
\end{equation}

\vspace{4mm}
\noindent
(al least one of the multipliers is equal to zero on the right-hand side of (\ref{2024october1})).

The equality (\ref{copa3}) is proved in our case 
(the right-hand side of (\ref{copa3}) is equal to zero for the considered case (see (\ref{copa6}))).

3. Suppose that the condition $(A)$ is not satisfied. 
In this case, we act according to the algorithm
from Sect.~2.27.4.
More precisely, let us select blocks
in the multi-index $j_{d_h}\ldots j_{d_h-e_h+1}$ $(h=1,\ldots,f)$
that
correspond to the fulfillment of the condition 
$$
\prod\limits_{l=1}^{r_{m,h}} {\bf 1}_{\{g_{2l}^{(h)}=g_{2l-1}^{(h)}+1\}}=1,
$$

\noindent
where $r_{m,h}$ is the number of pairs $\{g_{2l-1}^{(h)}, g_{2l}^{(h)}\}$ (from the set 
$\{g_1,g_2,\ldots,$ $g_{2r-1},g_{2r}\})$
in the block with number $m$ that corresponds to 
the multi-index $j_{d_h}\ldots j_{d_h-e_h+1}$.

Let us save multipliers of the form 
${\bf 1}_{\{t_n<t_{n+1}\}}$
in the  Volterra--type kernels corresponding to the Fourier
coefficients
\begin{equation}
\label{copa7}
~~~~~~~~C_{j_{d_1}\ldots j_{d_1-e_1+1}}^{\psi_{d_1}\ldots \psi_{d_1-e_1+1}}(t_{d_1+1},t_{d_1-e_1}),
\ldots, 
C_{j_{d_f}\ldots j_{d_f-e_f+1}}^{\psi_{d_f}\ldots \psi_{d_f-e_f+1}}(t_{d_f+1},t_{d_f-e_f})
\end{equation}

\noindent
and corresponding to the above blocks.

At that, we remove the remaining 
multipliers of the form 
${\bf 1}_{\{t_n<t_{n+1}\}}$
in the  Volterra--type kernels corresponding to the Fourier
coefficients (\ref{copa7}).

As a result, we get a modified left-hand side
of the equality (\ref{copa4}).
For definiteness, let us denote this expression by
$({}^{-})$.

Using generalized Parseval's equality
(Parseval's equality for two functions)
and (\ref{july30016}), we represent
the expression $({}^{-})$ as an integral over the hypercube $[t, T]^r$.

It is not difficult to see that the indicated integral over the hypercube $[t, T]^r$ is represented as a product
of integrals over hypercubes of smaller dimentions.
At that, at least one of these integrals 
is equal to zero
due to the generalized Parseval equality (Parseval's equality for two functions)
and the fulfillment of the condition
$t\le t_{d_1-e_1}\le t_{d_1+1}\le \ldots \le  t_{d_f-e_f} \le  t_{d_f+1} \le T$
(see the above example and (\ref{copa1}) and (\ref{copa1a})).
For definiteness, let us denote the equality of $({}^{-})$ to zero by $(\bar K)$.
We interpret the above zero as the zero functional in $L_2([t,T]^r).$
Further, transformations and passages to the limit
in the equality $(\bar K)$ are performed iteratively 
in such a way as to restore the removed multipliers ${\bf 1}_{\{t_n<t_{n+1}\}}$
on the left-hand side of $(\bar K)$
(for more details, see Sect.~2.27.4).
As a result, we obtain the equality (\ref{copa4}).
The equalities (\ref{copa3}) and (\ref{copa4}) are proved.

For definiteness, suppose that $q_1<\ldots <q_{k-2r},$ $k>2r,$ $r\ge 1$
(recall that the case $k=2r$ is proved in Sect.~2.27.4).
Using Fubini's Theorem (as in the above example (see (\ref{copa1})), we 
obtain
$$
\sum\limits_{j_{g_1},j_{g_3},\ldots, j_{g_{2r-1}}=0}^p
C_{j_k\ldots j_1}\biggl|_{j_{g_1}=j_{g_2},\ldots, j_{g_{2r-1}}=j_{g_{2r}}}=
$$

\vspace{3mm}
$$
=
\int\limits_t^T \psi_{q_{k-2r}}(t_{q_{k-2r}})
\phi_{j_{q_{k-2r}}}(t_{q_{k-2r}})\ldots
\int\limits_t^{t_{q_1+1}} \psi_{q_{1}}(t_{q_{1}})
\phi_{j_{q_{1}}}(t_{q_{1}})\times
$$

\vspace{1.5mm}
$$
\times 
\sum\limits_{j_{g_1},j_{g_3},\ldots, j_{g_{2r-1}}=0}^p\left(
C_{j_{d_f}\ldots j_{d_f-e_f+1}}^{\psi_{d_f}\ldots \psi_{d_f-e_f+1}}(t_{d_f+1},t_{d_f-e_f})
\ldots \right.
$$

\vspace{4mm}
$$
\left.
\ldots C_{j_{d_1}\ldots j_{d_1-e_1+1}}^{\psi_{d_1}\ldots \psi_{d_1-e_1+1}}(t_{d_1+1},t_{d_1-e_1})
\right)\biggl|_{j_{g_1}=j_{g_2},\ldots, j_{g_{2r-1}}=j_{g_{2r}}}\times
$$

\vspace{1.5mm}
\begin{equation}
\label{copa9}
\times dt_{q_1}\ldots dt_{q_{k-2r}},
\end{equation}

\vspace{4mm}

$$
\frac{1}{2^r} \prod\limits_{l=1}^r {\bf 1}_{\{g_{2l}=g_{2l-1}+1\}}
C_{j_k \ldots j_1}\biggl|_{(j_{g_2} j_{g_1})\curvearrowright (\cdot)
\ldots (j_{g_{2r}} j_{g_{2r-1}})\curvearrowright (\cdot),
j_{g_{{}_{1}}}=~j_{g_{{}_{2}}},\ldots, j_{g_{{}_{2r-1}}}=~j_{g_{{}_{2r}}}
}\biggr.=
$$

\vspace{5mm}
$$
=
\int\limits_t^T \psi_{q_{k-2r}}(t_{q_{k-2r}})
\phi_{j_{q_{k-2r}}}(t_{q_{k-2r}})\ldots
\int\limits_t^{t_{q_1+1}} \psi_{q_{1}}(t_{q_{1}})
\phi_{j_{q_{1}}}(t_{q_{1}})\times
$$

\vspace{1.5mm}
$$
\times {\bf 1}_{\{the\ condition\ (A)\ is\ satisfied\}}\prod\limits_{h=1}^f
\frac{1}{2^{r_h}}\prod\limits_{l=1}^{r_h}
{\bf 1}_{\{g_{2l}^{(h)}=g_{2l-1}^{(h)}+1\}}\times
$$

$$
\times C_{j_{d_h}\ldots j_{d_h-e_h+1}}^{\psi_{d_h}\ldots \psi_{d_h-e_h+1}}(t_{d_h+1},t_{d_h-e_h})
\biggl|_{(j_{g_2^{(h)}} j_{g_1^{(h)}})\curvearrowright (\cdot)
\ldots (j_{g_{2r_h}^{(h)}} j_{g_{2r_h-1}^{(h)}})\curvearrowright (\cdot),
j_{g_{{}_{1}}^{(h)}}=j_{g_{{}_{2}}^{(h)}},\ldots, j_{g_{{}_{2r_h-1}}^{(h)}}=j_{g_{{}_{2r_h}}^{(h)}}
}\biggr.\hspace{-0.7mm}\times
$$

\vspace{1.5mm}
\begin{equation}
\label{copa10}
\times
dt_{q_1}\ldots dt_{q_{k-2r}}.
\end{equation}

\vspace{3mm}

Suppose that
$$
\Biggl|
\sum\limits_{j_{g_1}, j_{g_3},\ldots ,j_{g_{2r-1}}=0}^p
\left(
C_{j_{d_f}\ldots j_{d_f-e_f+1}}^{\psi_{d_f}\ldots \psi_{d_f-e_f+1}}(t_{d_f+1},t_{d_f-e_f})
\ldots \right.\Biggr.
$$
\begin{equation}
\label{09091}
~~~~~~~~~\Biggl.\left.
\ldots C_{j_{d_1}\ldots j_{d_1-e_1+1}}^{\psi_{d_1}\ldots \psi_{d_1-e_1+1}}(t_{d_1+1},t_{d_1-e_1})
\right)\biggl|_{j_{g_1}=j_{g_2},\ldots, j_{g_{2r-1}}=j_{g_{2r}}}\Biggr|\le K <\infty,
\end{equation}

\vspace{3mm}
\noindent
where constant $K$ does not depend on $p$ and 
$t_{d_1+1},t_{d_1-e_1},\ldots,t_{d_f+1},t_{d_f-e_f}$
(here $d_1-e_1\ge 1$ and $d_f+1\le k$).
In (\ref{09091}):\ 
$t_{k+1}\stackrel{\sf def}{=}T,$
$t_{0}\stackrel{\sf def}{=}t,$
$e_1+\ldots+e_f=2r;$ another notations as above in this section.

Applying (\ref{copa3}), (\ref{copa4}), (\ref{copa9}), (\ref{copa10}), we obtain 
($k>2r,$ $r\ge 1$)

$$
\lim\limits_{p\to\infty}
\sum\limits_{\stackrel{j_1,\ldots,j_q,\ldots,j_k=0}{{}_{q\ne g_1, g_2, \ldots, g_{2r-1},
g_{2r}}}}^p
\Biggl(\sum\limits_{j_{g_1},j_{g_3},\ldots, j_{g_{2r-1}}=0}^p
C_{j_k\ldots j_1}\biggl|_{j_{g_1}=j_{g_2},\ldots, j_{g_{2r-1}}=j_{g_{2r}}}-\Biggr.
$$

\vspace{-2mm}
$$
\Biggl.-\frac{1}{2^r} \prod\limits_{l=1}^r {\bf 1}_{\{g_{2l}=g_{2l-1}+1\}}
C_{j_k \ldots j_1}\biggl|_{(j_{g_2} j_{g_1})\curvearrowright (\cdot)
\ldots (j_{g_{2r}} j_{g_{2r-1}})\curvearrowright (\cdot),
j_{g_{{}_{1}}}=~j_{g_{{}_{2}}},\ldots, j_{g_{{}_{2r-1}}}=~j_{g_{{}_{2r}}}
}\biggr.\Biggr)^2\le
$$

\vspace{6mm}
$$
\le \lim\limits_{p\to\infty}
\sum\limits_{\stackrel{j_1,\ldots,j_q,\ldots,j_k=0}{{}_{q\ne g_1, g_2, \ldots, g_{2r-1},
g_{2r}}}}^{\infty}
\Biggl(\sum\limits_{j_{g_1},j_{g_3},\ldots, j_{g_{2r-1}}=0}^p
C_{j_k\ldots j_1}\biggl|_{j_{g_1}=j_{g_2},\ldots, j_{g_{2r-1}}=j_{g_{2r}}}-\Biggr.
$$

\vspace{-1mm}
$$
\Biggl.-\frac{1}{2^r} \prod\limits_{l=1}^r {\bf 1}_{\{g_{2l}=g_{2l-1}+1\}}
C_{j_k \ldots j_1}\biggl|_{(j_{g_2} j_{g_1})\curvearrowright (\cdot)
\ldots (j_{g_{2r}} j_{g_{2r-1}})\curvearrowright (\cdot),
j_{g_{{}_{1}}}=~j_{g_{{}_{2}}},\ldots, j_{g_{{}_{2r-1}}}=~j_{g_{{}_{2r}}}
}\biggr.\Biggr)^2=
$$

\vspace{3mm}
$$
=\lim\limits_{p\to\infty}
\sum\limits_{j_{q_1},\ldots,j_{q_{k-2r}}=0}^{\infty}
\Biggl(
\int\limits_t^T \psi_{q_{k-2r}}(t_{q_{k-2r}})
\phi_{j_{q_{k-2r}}}(t_{q_{k-2r}})\ldots
\int\limits_t^{t_{q_1+1}} \psi_{q_{1}}(t_{q_{1}})
\phi_{j_{q_{1}}}(t_{q_{1}})\times\Biggr.
$$

\vspace{3mm}
$$
\times 
\Biggl(\sum\limits_{j_{g_1},j_{g_3},\ldots, j_{g_{2r-1}}=0}^p\left(
C_{j_{d_f}\ldots j_{d_f-e_f+1}}^{\psi_{d_f}\ldots \psi_{d_f-e_f+1}}(t_{d_f+1},t_{d_f-e_f})
\ldots \right.\Biggr.
$$
$$
\left.
\ldots C_{j_{d_1}\ldots j_{d_1-e_1+1}}^{\psi_{d_1}\ldots \psi_{d_1-e_1+1}}(t_{d_1+1},t_{d_1-e_1})
\right)\biggl|_{j_{g_1}=j_{g_2},\ldots, j_{g_{2r-1}}=j_{g_{2r}}}-
$$

\vspace{3mm}
$$
-{\bf 1}_{\{the\ condition\ (A)\ is\ satisfied\}}\prod\limits_{h=1}^f
\frac{1}{2^{r_h}}\prod\limits_{l=1}^{r_h}
{\bf 1}_{\{g_{2l}^{(h)}=g_{2l-1}^{(h)}+1\}}\times
$$

\vspace{-2.5mm}
$$
\Biggl.\times C_{j_{d_h}\ldots j_{d_h-e_h+1}}^{\psi_{d_h}\ldots \psi_{d_h-e_h+1}}(t_{d_h+1},t_{d_h-e_h})
\biggl|_{(j_{g_2^{(h)}} j_{g_1^{(h)}})\curvearrowright (\cdot)
\ldots (j_{g_{2r_h}^{(h)}} j_{g_{2r_h-1}^{(h)}})\curvearrowright (\cdot),
j_{g_{{}_{1}}^{(h)}}=j_{g_{{}_{2}}^{(h)}},\ldots, j_{g_{{}_{2r_h-1}}^{(h)}}=j_{g_{{}_{2r_h}}^{(h)}}
}\biggr.\hspace{-0.5mm}\Biggr)\hspace{-0.5mm}\times
$$

\vspace{2mm}
\begin{equation}
\label{2024dec1}
\Biggl.\times
dt_{q_1}\ldots dt_{q_{k-2r}}\Biggr)^2=
\end{equation}

\vspace{7mm}
$$
=\lim\limits_{p\to\infty}
\int\limits_t^T \psi_{q_{k-2r}}^2(t_{q_{k-2r}})
\ldots
\int\limits_t^{t_{q_1+1}} \psi_{q_{1}}^2(t_{q_{1}})\times
$$

\vspace{2mm}
$$
\times
\Biggl(\sum\limits_{j_{g_1},j_{g_3},\ldots, j_{g_{2r-1}}=0}^p\left(
C_{j_{d_f}\ldots j_{d_f-e_f+1}}^{\psi_{d_f}\ldots \psi_{d_f-e_f+1}}(t_{d_f+1},t_{d_f-e_f})
\ldots \right.\Biggr.
$$

\vspace{3mm}
$$
\left.
\ldots C_{j_{d_1}\ldots j_{d_1-e_1+1}}^{\psi_{d_1}\ldots \psi_{d_1-e_1+1}}(t_{d_1+1},t_{d_1-e_1})
\right)\biggl|_{j_{g_1}=j_{g_2},\ldots, j_{g_{2r-1}}=j_{g_{2r}}}-
$$

\vspace{4mm}
$$
-{\bf 1}_{\{the\ condition\ (A)\ is\ satisfied\}}\prod\limits_{h=1}^f
\frac{1}{2^{r_h}}\prod\limits_{l=1}^{r_h}
{\bf 1}_{\{g_{2l}^{(h)}=g_{2l-1}^{(h)}+1\}}\times
$$

$$
\Biggl.\times C_{j_{d_h}\ldots j_{d_h-e_h+1}}^{\psi_{d_h}\ldots \psi_{d_h-e_h+1}}(t_{d_h+1},t_{d_h-e_h})
\biggl|_{(j_{g_2^{(h)}} j_{g_1^{(h)}})\curvearrowright (\cdot)
\ldots (j_{g_{2r_h}^{(h)}} j_{g_{2r_h-1}^{(h)}})\curvearrowright (\cdot),
j_{g_{{}_{1}}^{(h)}}=j_{g_{{}_{2}}^{(h)}},\ldots, j_{g_{{}_{2r_h-1}}^{(h)}}=j_{g_{{}_{2r_h}}^{(h)}}
}\biggr.\hspace{-0.5mm}\Biggr)^2\hspace{-1.7mm}\times
$$

\vspace{1mm}
\begin{equation}
\label{copa14}
\Biggl.\times
dt_{q_1}\ldots dt_{q_{k-2r}}=
\end{equation}
$$
=
\int\limits_t^T \psi_{q_{k-2r}}^2(t_{q_{k-2r}})
\ldots
\int\limits_t^{t_{q_1+1}} \psi_{q_{1}}^2(t_{q_{1}})\times
$$

\vspace{2mm}
$$
\times
\lim\limits_{p\to\infty}\Biggl(\sum\limits_{j_{g_1},j_{g_3},\ldots, j_{g_{2r-1}}=0}^p\left(
C_{j_{d_f}\ldots j_{d_f-e_f+1}}^{\psi_{d_f}\ldots \psi_{d_f-e_f+1}}(t_{d_f+1},t_{d_f-e_f})
\ldots \right.\Biggr.
$$

\vspace{4mm}
$$
\left.
\ldots C_{j_{d_1}\ldots j_{d_1-e_1+1}}^{\psi_{d_1}\ldots \psi_{d_1-e_1+1}}(t_{d_1+1},t_{d_1-e_1})
\right)\biggl|_{j_{g_1}=j_{g_2},\ldots, j_{g_{2r-1}}=j_{g_{2r}}}-
$$

\vspace{4mm}
$$
-{\bf 1}_{\{the\ condition\ (A)\ is\ satisfied\}}\prod\limits_{h=1}^f
\frac{1}{2^{r_h}}\prod\limits_{l=1}^{r_h}
{\bf 1}_{\{g_{2l}^{(h)}=g_{2l-1}^{(h)}+1\}}\times
$$

\vspace{-5mm}
$$
\Biggl.\times C_{j_{d_h}\ldots j_{d_h-e_h+1}}^{\psi_{d_h}\ldots \psi_{d_h-e_h+1}}(t_{d_h+1},t_{d_h-e_h})
\biggl|_{(j_{g_2^{(h)}} j_{g_1^{(h)}})\curvearrowright (\cdot)
\ldots (j_{g_{2r_h}^{(h)}} j_{g_{2r_h-1}^{(h)}})\curvearrowright (\cdot),
j_{g_{{}_{1}}^{(h)}}=j_{g_{{}_{2}}^{(h)}},\ldots, j_{g_{{}_{2r_h-1}}^{(h)}}=j_{g_{{}_{2r_h}}^{(h)}}
}\biggr.\hspace{-0.5mm}\Biggr)^2\hspace{-1.7mm}\times
$$

\vspace{-2mm}
\begin{equation}
\label{copa15}
\Biggl.\times
dt_{q_1}\ldots dt_{q_{k-2r}}=0,
\end{equation}

\vspace{1mm}
\noindent
where 
the transition from (\ref{2024dec1}) to (\ref{copa14})
is based on the Parseval equality
and the transition from (\ref{copa14}) to (\ref{copa15}) 
is based on Lebesgue's Dominated Convergence Theorem (see
(\ref{july999}), (\ref{july1000aaa1}), (\ref{copa3}), (\ref{copa4}), (\ref{09091})) 
and also on
convergence to zero (almost everywhere on 

\vspace{-3mm}
$$
X=\{(t_{q_1},\ldots ,t_{q_{k-2r}}):  t\le t_{q_1}\le \ldots \le t_{q_{k-2r}}\le T\}
$$

\vspace{3mm}
\noindent
with respect to Lebesgue's measure)
of the integrand function in (\ref{copa14}).

Thus, the equality (\ref{april10}) and Hypotheses~2.4, 2.5
are proved for the case $p_1=\ldots=p_k=p$ 
under the condition (\ref{09091})
and we have the following theorem.

\vspace{2mm}

{\bf Theorem~2.60.}\ {\it Suppose that 
the condition {\rm (\ref{09091})} 
is fulfilled$,$
$\{\phi_j(x)\}_{j=0}^{\infty}$
is an arbitrary complete orthonormal system of functions
in $L_2([t, T])$ and
$\psi_1(\tau),\ldots, \psi_k(\tau)\in L_2([t, T]).$
Then$,$ for the sum $\bar J^{*}[\psi^{(k)}]_{T,t}^{(i_1\ldots i_k)}$
of iterated It\^{o} stochastic integrals 
$$
\bar J^{*}[\psi^{(k)}]_{T,t}^{(i_1\ldots i_k)}=J[\psi^{(k)}]_{T,t}^{(i_1\ldots i_k)}+
\sum_{r=1}^{\left[k/2\right]}\frac{1}{2^r}
\sum_{(s_r,\ldots,s_1)\in {\rm A}_{k,r}}
J[\psi^{(k)}]_{T,t}^{s_r,\ldots,s_1}
$$
the following 
expansion 
$$
\bar J^{*}[\psi^{(k)}]_{T,t}^{(i_1\ldots i_k)}=
\hbox{\vtop{\offinterlineskip\halign{
\hfil#\hfil\cr
{\rm l.i.m.}\cr
$\stackrel{}{{}_{p\to \infty}}$\cr
}} }
\sum_{j_1,\ldots,j_k=0}^{p}
C_{j_k \ldots j_1}\prod\limits_{l=1}^k \zeta_{j_l}^{(i_l)}
$$

\noindent
that converges in the mean-square sense is valid, where 
\begin{equation}
\label{2025may21}
~~~~~~~~C_{j_k \ldots j_1}=\int\limits_t^T\psi_k(t_k)\phi_{j_k}(t_k)\ldots
\int\limits_t^{t_2}
\psi_1(t_1)\phi_{j_1}(t_1)
dt_1\ldots dt_k
\end{equation}
is the Fourier coefficient, 
${\rm l.i.m.}$ is a limit in the mean-square sense,
$i_1, \ldots, i_k=0, 1,\ldots,m,$
$$
\zeta_{j}^{(i)}=
\int\limits_t^T \phi_{j}(\tau) d{\bf w}_{\tau}^{(i)}
$$ 
are independent standard Gaussian random variables for various 
$i$ or $j$ {\rm (}in the case when $i\ne 0${\rm )},
${\bf w}_{\tau}^{(i)}$ 
$(i=1,\ldots,m)$ are independent 
standard Wiener processes$,$
${\bf w}_{\tau}^{(0)}=\tau;$ another notations are the same as
in Theorem~{\rm 2.12}.}

\vspace{2mm}

Using Theorem~2.12, we obtain the following corollary of Theorem~2.60.

\vspace{2mm}

{\bf Theorem~2.61.}\ {\it Suppose that 
the condition {\rm (\ref{09091})} 
is fulfilled$,$
$\{\phi_j(x)\}_{j=0}^{\infty}$
is an arbitrary complete orthonormal system of functions
in $L_2([t, T])$ and
$\psi_1(\tau),\ldots, \psi_k(\tau)$ are continuous functions
at the interval $[t, T].$
Then$,$ for the iterated Stratonovich sto\-chas\-tic integral 
of multiplicity $k$ $(k\in{\bf N})$
$$
J^{*}[\psi^{(k)}]_{T,t}^{(i_1\ldots i_k)}=
{\int\limits_t^{*}}^T
\psi_k(t_k) \ldots 
{\int\limits_t^{*}}^{t_{2}}
\psi_1(t_1) d{\bf w}_{t_1}^{(i_1)}\ldots
d{\bf w}_{t_k}^{(i_k)}
$$
the following 
expansion 
\begin{equation}
\label{march000195}
J^{*}[\psi^{(k)}]_{T,t}^{(i_1\ldots i_k)}=
\hbox{\vtop{\offinterlineskip\halign{
\hfil#\hfil\cr
{\rm l.i.m.}\cr
$\stackrel{}{{}_{p\to \infty}}$\cr
}} }
\sum\limits_{j_1,\ldots,j_k=0}^{p}
C_{j_k \ldots j_1}\prod\limits_{l=1}^k \zeta_{j_l}^{(i_l)}
\end{equation}

\noindent
that converges in the mean-square sense is valid$;$ another notations are the same as
in Theorem~{\rm 2.60}.}
 
\vspace{2mm}

Note that the condition (\ref{09091})
can be weakened. Namely, the constant $K^2$
can be replaced by the function $F$ such that
$\psi_{q_1}^2\ldots \psi_{q_{k-2r}}^2 F \in L_1([t, T]^{k-2r})$
(integrable majorant). More precisely, 
the condition (\ref{09091})
can be replaced by the following condition

\vspace{-1mm}
$$
\Biggl(
\sum\limits_{j_{g_1}, j_{g_3},\ldots ,j_{g_{2r-1}}=0}^p
\left(
C_{j_{d_f}\ldots j_{d_f-e_f+1}}^{\psi_{d_f}\ldots \psi_{d_f-e_f+1}}(t_{d_f+1},t_{d_f-e_f})
\ldots \right.\Biggr.
$$
\begin{equation}
\label{09091xxx}
\Biggl.\left.
\ldots C_{j_{d_1}\ldots j_{d_1-e_1+1}}^{\psi_{d_1}\ldots \psi_{d_1-e_1+1}}(t_{d_1+1},t_{d_1-e_1})
\right)\biggl|_{j_{g_1}=j_{g_2},\ldots, j_{g_{2r-1}}=j_{g_{2r}}}\Biggr)^2\le 
F(t_{q_1},\ldots,t_{q_{k-2r}})
\end{equation}

\noindent
almost everywhere
on 
$$
X=\{(t_{q_1},\ldots ,t_{q_{k-2r}}):  t\le t_{q_1}\le \ldots \le t_{q_{k-2r}}\le T\}
$$

\noindent
with respect to Lebesgue's measure,
where the function $F(t_{q_1},\ldots,t_{q_{k-2r}})$
is such that
$$
\psi_{q_1}^2(t_{q_1})\ldots \psi_{q_{k-2r}}^2(t_{q_{k-2r}}) 
F(t_{q_1},\ldots,t_{q_{k-2r}}) \in L_1([t, T]^{k-2r}),
$$

\noindent
where $F(t_{q_1},\ldots,t_{q_{k-2r}})$ does not depend on $p.$
In (\ref{09091xxx}):\ 
$t_{k+1}\stackrel{\sf def}{=}T,$
$t_{0}\stackrel{\sf def}{=}t,$
$e_1+\ldots+e_f=2r;$ 
another notations as above in this section.

\section{Expansion of Iterated Stratonovich Stochastic Integrals
of Multiplicity 6. The Case of an Ar\-bit\-ra\-ry Complete Orthonormal System of 
Functions in the Space $L_2([t,T])$ and $\psi_1(\tau),\ldots, \psi_6(\tau)
\equiv 1$}

This section is devoted to the following theorem.

\vspace{2mm}

{\bf Theorem~2.62.}\ {\it Suppose that
$\{\phi_j(x)\}_{j=0}^{\infty}$ is an arbitrary complete orthonormal system of 
functions in the space $L_2([t,T]).$
Then$,$ for the iterated Stra\-to\-no\-vich stochastic integral
of sixth multiplicity 
$$
J^{*}[\psi^{(6)}]_{T,t}=
{\int\limits_t^{*}}^T
\ldots
{\int\limits_t^{*}}^{t_2}
d{\bf w}_{t_1}^{(i_1)}
\ldots d{\bf w}_{t_6}^{(i_6)}
$$
the following 
expansion 
$$
J^{*}[\psi^{(6)}]_{T,t}=
\hbox{\vtop{\offinterlineskip\halign{
\hfil#\hfil\cr
{\rm l.i.m.}\cr
$\stackrel{}{{}_{p\to \infty}}$\cr
}} }
\sum\limits_{j_1,\ldots, j_6=0}^{p}
C_{j_6 \ldots j_1}\zeta_{j_1}^{(i_1)}\ldots \zeta_{j_6}^{(i_6)}
$$
that converges in the mean-square sense is valid, where 
$i_1,\ldots,i_6=0, 1,\ldots,m,$
\begin{equation}
\label{november1620241}
C_{j_6\ldots j_1}=\int\limits_t^T
\phi_{j_6}(t_6)
\ldots
\int\limits_t^{t_2}
\phi_{j_1}(t_1)dt_1\ldots dt_6
\end{equation}
and
$$
\zeta_{j}^{(i)}=
\int\limits_t^T \phi_{j}(\tau) d{\bf w}_{\tau}^{(i)}
$$ 
are independent standard Gaussian random variables for various 
$i$ or $j$ {\rm (}when $i\ne 0${\rm ),}
${\bf w}_{\tau}^{(i)}$ 
$(i=1,\ldots,m)$ are independent 
standard Wiener processes$,$
${\bf w}_{\tau}^{(0)}=\tau.$}

\vspace{2mm}

{\bf Proof.}\ Our proof will be based on Theorem~2.61 
and verification of the equality (\ref{09091}) under the conditions
of Theorem~2.62 (the case $k=6>2r$, where $r=1, 2$). Recall that the case
$k=2r$ is considered in Sect.~2.27.4 (see (\ref{july90000})).
Under the conditions of Theorem~2.62, this means that
$k=6=2r$, where $r=3$.

Let throughout this proof 
$$
C_{j_k \ldots j_1}(s,\tau)=\int\limits_{\tau}^s
\phi_{j_k}(t_k)\ldots
\int\limits_{\tau}^{t_2}
\phi_{j_1}(t_1)dt_1\ldots dt_k,
$$

\noindent
where $k=1,\ldots,4,$\ $t\le\tau<s\le T,$
and $C_{j_6\ldots j_1}$ is defined by (\ref{november1620241}).

Using Fubini's Theorem and the technique that leads to the formulas (\ref{copa1}),
(\ref{copa1a}),
we obtain (note that we find all possible
combinations of pairs using the equality (\ref{after36})):

\vspace{3mm}

1. $r=1$ (15 combinations)
$$
C_{j_1 j_5 j_4 j_3 j_2 j_1}=
\hspace{-0.7mm}\int\limits_t^T\hspace{-0.7mm}
\phi_{j_5}(t_5)
\hspace{-0.7mm}\int\limits_t^{t_5}\hspace{-0.7mm}
\phi_{j_4}(t_4)
\hspace{-0.7mm}\int\limits_t^{t_4}\hspace{-0.7mm}
\phi_{j_3}(t_3)
\hspace{-0.7mm}\int\limits_t^{t_3}\hspace{-0.7mm}
\phi_{j_2}(t_2)
C_{j_1}(t_2,t) C_{j_1}(T,t_5)dt_2 dt_3 dt_4 dt_5,
$$
$$
C_{j_2 j_5 j_4 j_3 j_2 j_1}=
\hspace{-0.7mm}\int\limits_t^T\hspace{-0.8mm}
\phi_{j_5}(t_5)
\hspace{-0.8mm}\int\limits_t^{t_5}\hspace{-0.8mm}
\phi_{j_4}(t_4)
\hspace{-0.8mm}\int\limits_t^{t_4}\hspace{-0.8mm}
\phi_{j_3}(t_3)
\hspace{-0.8mm}\int\limits_t^{t_3}\hspace{-0.8mm}
\phi_{j_1}(t_1)
C_{j_2}(t_3,t_1) C_{j_2}(T,t_5)dt_1 dt_3 dt_4 dt_5,
$$
$$
C_{j_3 j_5 j_4 j_3 j_2 j_1}=
\hspace{-0.8mm}\int\limits_t^T\hspace{-0.8mm}
\phi_{j_5}(t_5)
\hspace{-0.8mm}\int\limits_t^{t_5}\hspace{-0.8mm}
\phi_{j_4}(t_4)
\hspace{-0.8mm}\int\limits_t^{t_4}\hspace{-0.8mm}
\phi_{j_2}(t_2)
\hspace{-0.8mm}\int\limits_t^{t_2}\hspace{-0.8mm}
\phi_{j_1}(t_1)
C_{j_3}(t_4,t_2) C_{j_3}(T,t_5)dt_1 dt_2 dt_4 dt_5,
$$
$$
C_{j_4 j_5 j_4 j_3 j_2 j_1}=
\hspace{-0.8mm}\int\limits_t^T\hspace{-0.8mm}
\phi_{j_5}(t_5)
\hspace{-0.8mm}\int\limits_t^{t_5}\hspace{-0.8mm}
\phi_{j_3}(t_3)
\hspace{-0.8mm}\int\limits_t^{t_3}\hspace{-0.8mm}
\phi_{j_2}(t_2)
\hspace{-0.8mm}\int\limits_t^{t_2}\hspace{-0.8mm}
\phi_{j_1}(t_1)
C_{j_4}(t_5,t_3) C_{j_4}(T,t_5)dt_1 dt_2 dt_3 dt_5,
$$
$$
C_{j_5 j_5 j_4 j_3 j_2 j_1}=
\hspace{-0.8mm}\int\limits_t^T\hspace{-0.8mm}
\phi_{j_4}(t_4)
\hspace{-0.8mm}\int\limits_t^{t_4}\hspace{-0.8mm}
\phi_{j_3}(t_3)
\hspace{-0.8mm}\int\limits_t^{t_3}\hspace{-0.8mm}
\phi_{j_2}(t_2)
\hspace{-0.8mm}\int\limits_t^{t_2}\hspace{-0.8mm}
\phi_{j_1}(t_1)
C_{j_5 j_5}(T,t_4)dt_1 dt_2 dt_3 dt_4,
$$
$$
C_{j_6 j_5 j_4 j_3 j_1 j_1}=
\hspace{-0.8mm}\int\limits_t^T\hspace{-0.8mm}
\phi_{j_6}(t_6)
\hspace{-0.8mm}\int\limits_t^{t_6}\hspace{-0.8mm}
\phi_{j_5}(t_5)
\hspace{-0.8mm}\int\limits_t^{t_5}\hspace{-0.8mm}
\phi_{j_4}(t_4)
\hspace{-0.8mm}\int\limits_t^{t_4}\hspace{-0.8mm}
\phi_{j_3}(t_3)
C_{j_1 j_1}(t_3,t)dt_3 dt_4 dt_5 dt_6,
$$
$$
C_{j_6 j_5 j_4 j_1 j_2 j_1}=
\hspace{-0.8mm}\int\limits_t^T\hspace{-0.8mm}
\phi_{j_6}(t_6)
\hspace{-0.8mm}\int\limits_t^{t_6}\hspace{-0.8mm}
\phi_{j_5}(t_5)
\hspace{-0.8mm}\int\limits_t^{t_5}\hspace{-0.8mm}
\phi_{j_4}(t_4)
\hspace{-0.8mm}\int\limits_t^{t_4}\hspace{-0.8mm}
\phi_{j_2}(t_2)
C_{j_1}(t_2,t) C_{j_1}(t_4,t_2)dt_2 dt_4 dt_5 dt_6,
$$
$$
C_{j_6 j_5 j_1 j_3 j_2 j_1}=
\hspace{-0.8mm}\int\limits_t^T\hspace{-0.8mm}
\phi_{j_6}(t_6)
\hspace{-0.8mm}\int\limits_t^{t_6}\hspace{-0.8mm}
\phi_{j_5}(t_5)
\hspace{-0.8mm}\int\limits_t^{t_5}\hspace{-0.8mm}
\phi_{j_3}(t_3)
\hspace{-0.8mm}\int\limits_t^{t_3}\hspace{-0.8mm}
\phi_{j_2}(t_2)
C_{j_1}(t_2,t) C_{j_1}(t_5,t_3)dt_2 dt_3 dt_5 dt_6,
$$
$$
C_{j_6 j_1 j_4 j_3 j_2 j_1}=
\hspace{-0.8mm}\int\limits_t^T\hspace{-0.8mm}
\phi_{j_6}(t_6)
\hspace{-0.8mm}\int\limits_t^{t_6}\hspace{-0.8mm}
\phi_{j_4}(t_4)
\hspace{-0.8mm}\int\limits_t^{t_4}\hspace{-0.8mm}
\phi_{j_3}(t_3)
\hspace{-0.8mm}\int\limits_t^{t_3}\hspace{-0.8mm}
\phi_{j_2}(t_2)
C_{j_1}(t_2,t) C_{j_1}(t_6,t_4)dt_2 dt_3 dt_4 dt_6,
$$
$$
C_{j_6 j_5 j_4 j_2 j_2 j_1}=
\hspace{-0.8mm}\int\limits_t^T\hspace{-0.8mm}
\phi_{j_6}(t_6)
\hspace{-0.8mm}\int\limits_t^{t_6}\hspace{-0.8mm}
\phi_{j_5}(t_5)
\hspace{-0.8mm}\int\limits_t^{t_5}\hspace{-0.8mm}
\phi_{j_4}(t_4)
\hspace{-0.8mm}\int\limits_t^{t_4}\hspace{-0.8mm}
\phi_{j_1}(t_1)
C_{j_2 j_2}(t_4,t_1)dt_1 dt_4 dt_5 dt_6,
$$
$$
C_{j_6 j_5 j_2 j_3 j_2 j_1}=
\hspace{-0.9mm}\int\limits_t^T\hspace{-0.9mm}
\phi_{j_6}(t_6)
\hspace{-0.9mm}\int\limits_t^{t_6}\hspace{-0.9mm}
\phi_{j_5}(t_5)
\hspace{-0.9mm}\int\limits_t^{t_5}\hspace{-0.9mm}
\phi_{j_3}(t_3)
\hspace{-0.9mm}\int\limits_t^{t_3}\hspace{-0.9mm}
\phi_{j_1}(t_1)
C_{j_2}(t_3,t_1) C_{j_2}(t_5,t_3)dt_1 dt_3 dt_5 dt_6,
$$
$$
C_{j_6 j_2 j_4 j_3 j_2 j_1}=
\hspace{-0.9mm}\int\limits_t^T\hspace{-0.9mm}
\phi_{j_6}(t_6)
\hspace{-0.9mm}\int\limits_t^{t_6}\hspace{-0.9mm}
\phi_{j_4}(t_4)
\hspace{-0.9mm}\int\limits_t^{t_4}\hspace{-0.9mm}
\phi_{j_3}(t_3)
\hspace{-0.9mm}\int\limits_t^{t_3}\hspace{-0.9mm}
\phi_{j_1}(t_1)
C_{j_2}(t_3,t_1) C_{j_2}(t_6,t_4)dt_1 dt_3 dt_4 dt_6,
$$
$$
C_{j_6 j_5 j_3 j_3 j_2 j_1}=
\hspace{-0.8mm}\int\limits_t^T\hspace{-0.8mm}
\phi_{j_6}(t_6)
\hspace{-0.8mm}\int\limits_t^{t_6}\hspace{-0.8mm}
\phi_{j_5}(t_5)
\hspace{-0.8mm}\int\limits_t^{t_5}\hspace{-0.8mm}
\phi_{j_2}(t_2)
\hspace{-0.8mm}\int\limits_t^{t_2}\hspace{-0.8mm}
\phi_{j_1}(t_1)
C_{j_3 j_3}(t_5,t_2)dt_1 dt_2 dt_5 dt_6,
$$
$$
C_{j_6 j_3 j_4 j_3 j_2 j_1}=
\hspace{-0.9mm}\int\limits_t^T\hspace{-0.9mm}
\phi_{j_6}(t_6)
\hspace{-0.9mm}\int\limits_t^{t_6}\hspace{-0.9mm}
\phi_{j_4}(t_4)
\hspace{-0.9mm}\int\limits_t^{t_4}\hspace{-0.9mm}
\phi_{j_2}(t_2)
\hspace{-0.9mm}\int\limits_t^{t_2}\hspace{-0.9mm}
\phi_{j_1}(t_1)
C_{j_3}(t_4,t_2) C_{j_3}(t_6,t_4)dt_1 dt_2 dt_4 dt_6,
$$
$$
C_{j_6 j_4 j_4 j_3 j_2 j_1}=
\hspace{-0.8mm}\int\limits_t^T\hspace{-0.8mm}
\phi_{j_6}(t_6)
\hspace{-0.8mm}\int\limits_t^{t_6}\hspace{-0.8mm}
\phi_{j_3}(t_3)
\hspace{-0.8mm}\int\limits_t^{t_3}\hspace{-0.8mm}
\phi_{j_2}(t_2)
\hspace{-0.8mm}\int\limits_t^{t_2}\hspace{-0.8mm}
\phi_{j_1}(t_1)
C_{j_4 j_4}(t_6,t_3)dt_1 dt_2 dt_3 dt_6,
$$

\vspace{3mm}

2. $r=2$ (45 combinations)
$$
C_{j_6 j_5 j_3 j_3 j_1 j_1}=
\int\limits_t^T
\phi_{j_6}(t_6)
\int\limits_t^{t_6}
\phi_{j_5}(t_5)
C_{j_3 j_3 j_1 j_1}(t_5,t)dt_5 dt_6,
$$
$$
C_{j_6 j_3 j_4 j_3 j_1 j_1}=
\int\limits_t^T
\phi_{j_6}(t_6)
\int\limits_t^{t_6}
\phi_{j_4}(t_4)
C_{j_3 j_1 j_1}(t_4,t)C_{j_3}(t_6,t_4)dt_4 dt_6,
$$
$$
C_{j_6 j_4 j_4 j_3 j_1 j_1}=
\int\limits_t^T
\phi_{j_6}(t_6)
\int\limits_t^{t_6}
\phi_{j_3}(t_3)
C_{j_1 j_1}(t_3,t)C_{j_4 j_4}(t_6,t_3)dt_3 dt_6,
$$
$$
C_{j_6 j_5 j_2 j_1 j_2 j_1}=
\int\limits_t^T
\phi_{j_6}(t_6)
\int\limits_t^{t_6}
\phi_{j_5}(t_5)
C_{j_2 j_1 j_2 j_1}(t_5,t)dt_5 dt_6,
$$
$$
C_{j_6 j_2 j_4 j_1 j_2 j_1}=
\int\limits_t^T
\phi_{j_6}(t_6)
\int\limits_t^{t_6}
\phi_{j_4}(t_4)
C_{j_1 j_2 j_1}(t_4,t)C_{j_2}(t_6,t_4)dt_4 dt_6,
$$
$$
C_{j_6 j_4 j_4 j_1 j_2 j_1}=
\int\limits_t^T
\phi_{j_6}(t_6)
\int\limits_t^{t_6}
\phi_{j_2}(t_2)
C_{j_1}(t_2,t)C_{j_4 j_4 j_1}(t_6,t_2)dt_2 dt_6,
$$
$$
C_{j_6 j_5 j_1 j_2 j_2 j_1}=
\int\limits_t^T
\phi_{j_6}(t_6)
\int\limits_t^{t_6}
\phi_{j_5}(t_5)
C_{j_1 j_2 j_2 j_1}(t_5,t)dt_5 dt_6,
$$
$$
C_{j_6 j_2 j_1 j_3 j_2 j_1}=
\int\limits_t^T
\phi_{j_6}(t_6)
\int\limits_t^{t_6}
\phi_{j_3}(t_3)
C_{j_2 j_1}(t_3,t)C_{j_2 j_1}(t_6,t_3)dt_3 dt_6,
$$
$$
C_{j_6 j_3 j_1 j_3 j_2 j_1}=
\int\limits_t^T
\phi_{j_6}(t_6)
\int\limits_t^{t_6}
\phi_{j_2}(t_2)
C_{j_1}(t_2,t)C_{j_3 j_1 j_3}(t_6,t_2)dt_2 dt_6,
$$
$$
C_{j_6 j_1 j_4 j_2 j_2 j_1}=
\int\limits_t^T
\phi_{j_6}(t_6)
\int\limits_t^{t_6}
\phi_{j_4}(t_4)
C_{j_2 j_2 j_1}(t_4,t)C_{j_1}(t_6,t_4)dt_4 dt_6,
$$
$$
C_{j_6 j_1 j_2 j_3 j_2 j_1}=
\int\limits_t^T
\phi_{j_6}(t_6)
\int\limits_t^{t_6}
\phi_{j_3}(t_3)
C_{j_2 j_1}(t_3,t)C_{j_1 j_2}(t_6,t_3)dt_3 dt_6,
$$
$$
C_{j_6 j_1 j_3 j_3 j_2 j_1}=
\int\limits_t^T
\phi_{j_6}(t_6)
\int\limits_t^{t_6}
\phi_{j_2}(t_2)
C_{j_1}(t_2,t)C_{j_1 j_3 j_3}(t_6,t_2)dt_2 dt_6,
$$
$$
C_{j_6 j_4 j_4 j_2 j_2 j_1}=
\int\limits_t^T
\phi_{j_6}(t_6)
\int\limits_t^{t_6}
\phi_{j_1}(t_1)
C_{j_4 j_4 j_2 j_2}(t_6,t_1)dt_1 dt_6,
$$
$$
C_{j_6 j_3 j_2 j_3 j_2 j_1}=
\int\limits_t^T
\phi_{j_6}(t_6)
\int\limits_t^{t_6}
\phi_{j_1}(t_1)
C_{j_3 j_2 j_3 j_2}(t_6,t_1)dt_1 dt_6,
$$
$$
C_{j_6 j_2 j_3 j_3 j_2 j_1}=
\int\limits_t^T
\phi_{j_6}(t_6)
\int\limits_t^{t_6}
\phi_{j_1}(t_1)
C_{j_2 j_3 j_3 j_2}(t_6,t_1)dt_1 dt_6,
$$
$$
C_{j_1 j_5 j_3 j_3 j_2 j_1}=
\int\limits_t^T
\phi_{j_5}(t_5)
\int\limits_t^{t_5}
\phi_{j_2}(t_2)
C_{j_1}(t_2,t)C_{j_3 j_3}(t_5,t_2)C_{j_1}(T,t_5)dt_2 dt_5,
$$
$$
C_{j_1 j_3 j_4 j_3 j_2 j_1}=
\int\limits_t^T
\phi_{j_4}(t_4)
\int\limits_t^{t_4}
\phi_{j_2}(t_2)
C_{j_1}(t_2,t)C_{j_3}(t_4,t_2)C_{j_1 j_3}(T,t_4)dt_2 dt_4,
$$
$$
C_{j_1 j_2 j_4 j_3 j_2 j_1}=
\int\limits_t^T
\phi_{j_4}(t_4)
\int\limits_t^{t_4}
\phi_{j_3}(t_3)
C_{j_2 j_1}(t_3,t)C_{j_1 j_2}(T,t_4)dt_3 dt_4,
$$
$$
C_{j_1 j_5 j_2 j_3 j_2 j_1}=
\int\limits_t^T
\phi_{j_5}(t_5)
\int\limits_t^{t_5}
\phi_{j_3}(t_3)
C_{j_2 j_1}(t_3,t)C_{j_2}(t_5,t_3)C_{j_1}(T,t_5)dt_3 dt_5,
$$
$$
C_{j_1 j_4 j_4 j_3 j_2 j_1}=
\int\limits_t^T
\phi_{j_3}(t_3)
\int\limits_t^{t_3}
\phi_{j_2}(t_2)
C_{j_1}(t_2,t)C_{j_1 j_4 j_4}(T,t_3)dt_2 dt_3,
$$
$$
C_{j_1 j_5 j_4 j_2 j_2 j_1}=
\int\limits_t^T
\phi_{j_5}(t_5)
\int\limits_t^{t_5}
\phi_{j_4}(t_4)
C_{j_2 j_2 j_1}(t_4,t)C_{j_1}(T,t_5)dt_4 dt_5,
$$
$$
C_{j_2 j_3 j_4 j_3 j_2 j_1}=
\int\limits_t^T
\phi_{j_4}(t_4)
\int\limits_t^{t_4}
\phi_{j_1}(t_1)
C_{j_2 j_3}(t_4,t_1)C_{j_2 j_3}(T,t_4)dt_1 dt_4,
$$
$$
C_{j_2 j_4 j_4 j_3 j_2 j_1}=
\int\limits_t^T
\phi_{j_3}(t_3)
\int\limits_t^{t_3}
\phi_{j_1}(t_1)
C_{j_2}(t_3,t_1)C_{j_2 j_4 j_4}(T,t_3)dt_1 dt_3,
$$
$$
C_{j_2 j_5 j_3 j_3 j_2 j_1}=
\int\limits_t^T
\phi_{j_5}(t_5)
\int\limits_t^{t_5}
\phi_{j_1}(t_1)
C_{j_2 j_3 j_3}(t_5,t_1)C_{j_2}(T,t_5)dt_1 dt_5,
$$
$$
C_{j_2 j_1 j_4 j_3 j_2 j_1}=
\int\limits_t^T
\phi_{j_4}(t_4)
\int\limits_t^{t_4}
\phi_{j_3}(t_3)
C_{j_2 j_1}(t_3,t)C_{j_2 j_1}(T,t_4)dt_3 dt_4,
$$
$$
C_{j_2 j_5 j_1 j_3 j_2 j_1}=
\int\limits_t^T
\phi_{j_5}(t_5)
\int\limits_t^{t_5}
\phi_{j_3}(t_3)
C_{j_2 j_1}(t_3,t)C_{j_1}(t_5,t_3)C_{j_2}(T,t_5)dt_3 dt_5,
$$
$$
C_{j_2 j_5 j_4 j_1 j_2 j_1}=
\int\limits_t^T
\phi_{j_5}(t_5)
\int\limits_t^{t_5}
\phi_{j_4}(t_4)
C_{j_1 j_2 j_1}(t_4,t)C_{j_2}(T,t_5)dt_4 dt_5,
$$
$$
C_{j_3 j_2 j_4 j_3 j_2 j_1}=
\int\limits_t^T
\phi_{j_4}(t_4)
\int\limits_t^{t_4}
\phi_{j_1}(t_1)
C_{j_3 j_2}(t_4,t_1)C_{j_3 j_2}(T,t_4)dt_1 dt_4,
$$
$$
C_{j_3 j_4 j_4 j_3 j_2 j_1}=
\int\limits_t^T
\phi_{j_2}(t_2)
\int\limits_t^{t_2}
\phi_{j_1}(t_1)
C_{j_3 j_4 j_4 j_3}(T,t_2)dt_1 dt_2,
$$
$$
C_{j_3 j_5 j_2 j_3 j_2 j_1}=
\int\limits_t^T
\phi_{j_5}(t_5)
\int\limits_t^{t_5}
\phi_{j_1}(t_1)
C_{j_2 j_3 j_2}(t_5,t_1)C_{j_3}(T,t_5)dt_1 dt_5,
$$
$$
C_{j_3 j_1 j_4 j_3 j_2 j_1}=
\int\limits_t^T
\phi_{j_4}(t_4)
\int\limits_t^{t_4}
\phi_{j_2}(t_2)
C_{j_1}(t_2,t)C_{j_3}(t_4,t_2)C_{j_3 j_1}(T,t_4)dt_2 dt_4,
$$
$$
C_{j_3 j_5 j_1 j_3 j_2 j_1}=
\int\limits_t^T
\phi_{j_5}(t_5)
\int\limits_t^{t_5}
\phi_{j_2}(t_2)
C_{j_1}(t_2,t)C_{j_1 j_3}(t_5,t_2)C_{j_3}(T,t_5)dt_2 dt_5,
$$
$$
C_{j_3 j_5 j_4 j_3 j_1 j_1}=
\int\limits_t^T
\phi_{j_5}(t_5)
\int\limits_t^{t_5}
\phi_{j_4}(t_4)
C_{j_3 j_1 j_1}(t_4,t)C_{j_3}(T,t_5)dt_4 dt_5,
$$
$$
C_{j_4 j_3 j_4 j_3 j_2 j_1}=
\int\limits_t^T
\phi_{j_2}(t_2)
\int\limits_t^{t_2}
\phi_{j_1}(t_1)
C_{j_4 j_3 j_4 j_3}(T,t_2)dt_1 dt_2,
$$
$$
C_{j_4 j_2 j_4 j_3 j_2 j_1}=
\int\limits_t^T
\phi_{j_3}(t_3)
\int\limits_t^{t_3}
\phi_{j_1}(t_1)
C_{j_2}(t_3,t_1)C_{j_4 j_2 j_4}(T,t_3)dt_1 dt_3,
$$
$$
C_{j_4 j_5 j_4 j_2 j_2 j_1}=
\int\limits_t^T
\phi_{j_5}(t_5)
\int\limits_t^{t_5}
\phi_{j_1}(t_1)
C_{j_4 j_2 j_2}(t_5,t_1)C_{j_4}(T,t_5)dt_1 dt_5,
$$
$$
C_{j_4 j_1 j_4 j_3 j_2 j_1}=
\int\limits_t^T
\phi_{j_3}(t_3)
\int\limits_t^{t_3}
\phi_{j_2}(t_2)
C_{j_1}(t_2,t)C_{j_4 j_1 j_4}(T,t_3)dt_2 dt_3,
$$
$$
C_{j_4 j_5 j_4 j_1 j_2 j_1}=
\int\limits_t^T
\phi_{j_5}(t_5)
\int\limits_t^{t_5}
\phi_{j_2}(t_2)
C_{j_1}(t_2,t)C_{j_4 j_1}(t_5,t_2)C_{j_4}(T,t_5)dt_2 dt_5,
$$
$$
C_{j_4 j_5 j_4 j_3 j_1 j_1}=
\int\limits_t^T
\phi_{j_5}(t_5)
\int\limits_t^{t_5}
\phi_{j_3}(t_3)
C_{j_1 j_1}(t_3,t)C_{j_4}(t_5,t_3)C_{j_4}(T,t_5)dt_3 dt_5,
$$
$$
C_{j_5 j_5 j_3 j_3 j_2 j_1}=
\int\limits_t^T
\phi_{j_2}(t_2)
\int\limits_t^{t_2}
\phi_{j_1}(t_1)
C_{j_5 j_5 j_3 j_3}(T,t_2)dt_1 dt_2,
$$
$$
C_{j_5 j_5 j_2 j_3 j_2 j_1}=
\int\limits_t^T
\phi_{j_3}(t_3)
\int\limits_t^{t_3}
\phi_{j_1}(t_1)
C_{j_2}(t_3,t_1)C_{j_5 j_5 j_2}(T,t_3)dt_1 dt_3,
$$
$$
C_{j_5 j_5 j_4 j_2 j_2 j_1}=
\int\limits_t^T
\phi_{j_4}(t_4)
\int\limits_t^{t_4}
\phi_{j_1}(t_1)
C_{j_2 j_2}(t_4,t_1)C_{j_5 j_5}(T,t_4)dt_1 dt_4,
$$
$$
C_{j_5 j_5 j_1 j_3 j_2 j_1}=
\int\limits_t^T
\phi_{j_3}(t_3)
\int\limits_t^{t_3}
\phi_{j_2}(t_2)
C_{j_1}(t_2,t)C_{j_5 j_5 j_1}(T,t_3)dt_2 dt_3,
$$
$$
C_{j_5 j_5 j_4 j_1 j_2 j_1}=
\int\limits_t^T
\phi_{j_4}(t_4)
\int\limits_t^{t_4}
\phi_{j_2}(t_2)
C_{j_1}(t_2,t)C_{j_1}(t_4,t_2)C_{j_5 j_5}(T,t_4)dt_2 dt_4,
$$
$$
C_{j_5 j_5 j_4 j_3 j_1 j_1}=
\int\limits_t^T
\phi_{j_4}(t_4)
\int\limits_t^{t_4}
\phi_{j_3}(t_3)
C_{j_1 j_1}(t_3,t)C_{j_5 j_5}(T,t_4)dt_3 dt_4.
$$

It is not difficult to see (based on the above equalities)
that the condition (\ref{09091}) will be satisfied under the conditions
of Theorem~2.62 if
\begin{equation}
\label{2024december12}
\left\vert \sum\limits_{j_1=0}^p 
C_{j_1 j_1}(s,\tau)\right\vert\le K,
\end{equation}
\begin{equation}
\label{2024december13}
\left\vert \sum\limits_{j_1=0}^p 
C_{j_1}(s,\tau)C_{j_1}(\theta,u)\right\vert\le K,
\end{equation}
\begin{equation}
\label{2024december1}
\left\vert \sum\limits_{j_1, j_2=0}^p C_{j_2 j_2 j_1 j_1}(s,\tau)\right\vert\le K,
\end{equation}
\begin{equation}
\label{2024december2}
\left\vert \sum\limits_{j_1, j_2=0}^p C_{j_2 j_1 j_2 j_1}(s,\tau)\right\vert\le K,
\end{equation}
\begin{equation}
\label{2024december3}
\left\vert \sum\limits_{j_1, j_2=0}^p C_{j_1 j_2 j_2 j_1}(s,\tau)\right\vert\le K,
\end{equation}
\begin{equation}
\label{2024december4}
\left\vert \sum\limits_{j_1, j_2=0}^p C_{j_2 j_1 j_1}(s,\tau)C_{j_2}(\theta,u)\right\vert\le K,
\end{equation}
\begin{equation}
\label{2024december5}
\left\vert \sum\limits_{j_1, j_2=0}^p C_{j_1 j_2 j_1}(s,\tau)C_{j_2}(\theta,u)\right\vert\le K,
\end{equation}
\begin{equation}
\label{2024december6}
\left\vert \sum\limits_{j_1, j_2=0}^p C_{j_2 j_2 j_1}(s,\tau)C_{j_1}(\theta,u)\right\vert\le K,
\end{equation}
\begin{equation}
\label{2024december7}
\left\vert \sum\limits_{j_1, j_2=0}^p C_{j_1 j_1}(s,\tau)C_{j_2 j_2}(\theta,u)\right\vert\le K,
\end{equation}
\begin{equation}
\label{2024december8}
\left\vert \sum\limits_{j_1, j_2=0}^p C_{j_2 j_1}(s,\tau)C_{j_2 j_1}(\theta,u)\right\vert\le K,
\end{equation}
\begin{equation}
\label{2024december9}
\left\vert \sum\limits_{j_1, j_2=0}^p C_{j_2 j_1}(s,\tau)C_{j_1 j_2}(\theta,u)\right\vert\le K,
\end{equation}
\begin{equation}
\label{2024december10}
\left\vert \sum\limits_{j_1, j_2=0}^p C_{j_1}(s,\tau)C_{j_1}(\rho,v)
C_{j_2 j_2}(\theta,u)\right\vert\le K,
\end{equation}
\begin{equation}
\label{2024december11}
\left\vert \sum\limits_{j_1, j_2=0}^p C_{j_1}(s,\tau)C_{j_2}(\rho,v)
C_{j_1 j_2}(\theta,u)\right\vert\le K,
\end{equation}

\noindent
where $p\in{\bf N},$ $t\le \tau < s \le T,$ $t\le u<\theta \le T,$
$t\le v<\rho \le T,$ constant $K$ does not depend on 
$p, s, \tau, u, \theta, v, \rho$ (but only on $t, T$) and may differ from line to line.

The equalities (\ref{2024december1})--(\ref{2024december3}) 
have been proved earlier (see (\ref{april50})--(\ref{april52})).

Using Fubini's Theorem and Parseval's equality, we get
$$
\left\vert \sum\limits_{j_1=0}^p 
C_{j_1 j_1}(s,\tau)\right\vert=
\frac{1}{2}\sum\limits_{j_1=0}^p C_{j_1}^2(s,\tau)\le
$$
$$
\le
\frac{1}{2}\sum\limits_{j_1=0}^{\infty}
C_{j_1}^2(s,\tau)=
\frac{1}{2}(s-\tau) \le
$$
$$
\le \frac{1}{2}(T-t) \le K.
$$

\noindent
The equality (\ref{2024december12}) is proved.
Moreover, (\ref{2024december7}) follows from 
(\ref{2024december12}).

Using the inequality of Cauchy--Bunyakovsky and 
Parseval's equality, we obtain
$$
\left(\sum\limits_{j_1=0}^p 
C_{j_1}(s,\tau)C_{j_1}(\theta,u)\right)^2\le
$$
$$
\le
\sum\limits_{j_1=0}^p 
C_{j_1}^2(s,\tau)
\sum\limits_{j_1=0}^p 
C_{j_1}^2(\theta,u)\le
$$
$$
\le \sum\limits_{j_1=0}^{\infty} 
C_{j_1}^2(s,\tau)
\sum\limits_{j_1=0}^{\infty} 
C_{j_1}^2(\theta,u)=
$$
$$
=(s-\tau)(\theta-u)\le 
(T-t)^2 \le K^2,
$$

\vspace{-5mm}
$$
\left(\sum\limits_{j_1, j_2=0}^p C_{j_2 j_1}(s,\tau)C_{j_2 j_1}(\theta,u)\right)^2
\le
\sum\limits_{j_1,j_2=0}^p 
C_{j_2 j_1}^2(s,\tau)
\sum\limits_{j_1,j_2=0}^p 
C_{j_2 j_1}^2(\theta,u)\le
$$
$$
\le\sum\limits_{j_1,j_2=0}^{\infty}
C_{j_2 j_1}^2(s,\tau)
\sum\limits_{j_1,j_2=0}^{\infty} 
C_{j_2 j_1}^2(\theta,u)=
$$
$$
=\int\limits_{\tau}^s \int\limits_{\tau}^v dx dv
\int\limits_{u}^{\theta} \int\limits_{u}^v dx dv\le 
\frac{1}{4}(T-t)^4\le K^2.
$$

\vspace{1mm}

Thus, the inequalities (\ref{2024december13}), (\ref{2024december8}) are proved.
The inequalities (\ref{2024december9}), (\ref{2024december11}) are 
proved similarly to (\ref{2024december8}).
Moreover, (\ref{2024december10}) follows from
(\ref{2024december12}), (\ref{2024december13}).

Further, let us prove the equalities (\ref{2024december4})--(\ref{2024december6}).
Applying the 
Cau\-chy--Bunyakovsky inequality as well as 
Parseval's equality and (\ref{2024december12}), we have
$$
\left(\sum\limits_{j_1, j_2=0}^p C_{j_2 j_1 j_1}(s,\tau)C_{j_2}(\theta,u)\right)^2
\le 
\sum\limits_{j_2=0}^p \left(\sum\limits_{j_1=0}^p C_{j_2 j_1 j_1}(s,\tau)\right)^2
\sum\limits_{j_2=0}^p C_{j_2}^2(\theta,u)\le
$$
$$
\le \sum\limits_{j_2=0}^{\infty} \left(\sum\limits_{j_1=0}^p C_{j_2 j_1 j_1}(s,\tau)\right)^2
\sum\limits_{j_2=0}^{\infty} C_{j_2}^2(\theta,u)=
$$
$$
=\sum\limits_{j_2=0}^{\infty} \left(
\int\limits_{\tau}^s \phi_{j_2}(v) \sum\limits_{j_1=0}^p  C_{j_1 j_1}(v,\tau) dv
\right)^2 \cdot (\theta-u)=
$$
$$
=
(\theta-u)\int\limits_{\tau}^s \left(\sum\limits_{j_1=0}^p  C_{j_1 j_1}(v,\tau)\right)^2 dv
\le 
$$

\vspace{-2mm}
$$
\le K^2 (\theta-u)(s-\tau)\le K^2(T-t)^2=K_1.
$$

\vspace{1mm}
\noindent
The equality (\ref{2024december4}) is proved.

Using the 
Cau\-chy--Bunyakovsky inequality as well as 
Fubini's Theorem, Parseval's equality and (\ref{2024december13}), we have
$$
\left(\sum\limits_{j_1, j_2=0}^p C_{j_1 j_2 j_1}(s,\tau)C_{j_2}(\theta,u)\right)^2
\le 
\sum\limits_{j_2=0}^p \left(\sum\limits_{j_1=0}^p C_{j_1 j_2 j_1}(s,\tau)\right)^2
\sum\limits_{j_2=0}^p C_{j_2}^2(\theta,u)\le
$$
$$
\le \sum\limits_{j_2=0}^{\infty} \left(\sum\limits_{j_1=0}^p 
\int\limits_{\tau}^s \phi_{j_1}(z)\int\limits_{\tau}^{z} \phi_{j_2}(y)
\int\limits_{\tau}^{y} \phi_{j_1}(x)dx dy dz\right)^2
\sum\limits_{j_2=0}^{\infty} C_{j_2}^2(\theta,u)=
$$
$$
=\sum\limits_{j_2=0}^{\infty} \left(\sum\limits_{j_1=0}^p 
\int\limits_{\tau}^s 
\phi_{j_2}(y)
\int\limits_{\tau}^{y} \phi_{j_1}(x)dx
\int\limits_{y}^{s}
\phi_{j_1}(z) dz dy\right)^2
\cdot (\theta-u)=
$$
$$
=(\theta-u)\sum\limits_{j_2=0}^{\infty} \left(
\int\limits_{\tau}^s 
\phi_{j_2}(y)
\sum\limits_{j_1=0}^p C_{j_1}(y,\tau) C_{j_1}(s,y) dy\right)^2=
$$
$$
=(\theta-u)
\int\limits_{\tau}^s 
\left(\sum\limits_{j_1=0}^p C_{j_1}(y,\tau) C_{j_1}(s,y)\right)^2 dy
\le 
$$

\vspace{-1mm}
$$
\le K^2 (\theta-u)(s-\tau)\le K^2(T-t)^2=K_1.
$$

\vspace{1mm}
\noindent
The equality (\ref{2024december5}) is proved.

Using the 
Cau\-chy--Bunyakovsky inequality as well as 
Fubini's Theorem, Parseval's equality and (\ref{2024december12}), we have
$$
\left(\sum\limits_{j_1, j_2=0}^p C_{j_2 j_2 j_1}(s,\tau)C_{j_1}(\theta,u)\right)^2\le
\sum\limits_{j_1=0}^p \left(\sum\limits_{j_2=0}^p C_{j_2 j_2 j_1}(s,\tau)\right)^2
\sum\limits_{j_1=0}^p C_{j_1}^2(\theta,u)\le
$$
$$
\le \sum\limits_{j_1=0}^{\infty} \left(\sum\limits_{j_2=0}^p 
\int\limits_{\tau}^s \phi_{j_2}(z)\int\limits_{\tau}^{z} \phi_{j_2}(y)
\int\limits_{\tau}^{y} \phi_{j_1}(x)dx dy dz\right)^2
\sum\limits_{j_1=0}^{\infty} C_{j_1}^2(\theta,u)=
$$
$$
=\sum\limits_{j_1=0}^{\infty} \left(\sum\limits_{j_2=0}^p 
\int\limits_{\tau}^s \phi_{j_1}(x)\int\limits_{x}^{s} \phi_{j_2}(y)
\int\limits_{y}^{s} \phi_{j_2}(z)dz dy dx\right)^2
\cdot (\theta-u)=
$$
$$
=(\theta-u)\sum\limits_{j_1=0}^{\infty} \left(\sum\limits_{j_2=0}^p 
\int\limits_{\tau}^s \phi_{j_1}(x)\int\limits_{x}^{s} \phi_{j_2}(z)
\int\limits_{x}^{z} \phi_{j_2}(y)dy dz dx\right)^2=
$$
$$
=(\theta-u)\sum\limits_{j_1=0}^{\infty} \left(
\int\limits_{\tau}^s \phi_{j_1}(x)  
\sum\limits_{j_2=0}^p C_{j_2 j_2}(s,x)
dx\right)^2=
$$
$$
=(\theta-u)
\int\limits_{\tau}^s \left(
\sum\limits_{j_2=0}^p C_{j_2 j_2}(s,x)\right)^2
dx\le
$$

\vspace{-1mm}
$$
\le K^2 (\theta-u)(s-\tau)\le K^2(T-t)^2=K_1.
$$

\vspace{3mm}
\noindent
The equality (\ref{2024december6}) is proved.
The equalities (\ref{2024december12})--(\ref{2024december11})
are proved.

Thus, the condition (\ref{09091}) of Theorem~2.61 is satisfied 
under the conditions of Theorem~2.62. The assertion
of Theorem~2.62 now follows from Theorem~2.61. 
Theorem~2.62 is proved.

\section{Expansion of Iterated Stratonovich Stochastic Integrals
of Multiplicity 4. The Case of an Ar\-bit\-ra\-ry Complete Orthonormal System of 
Functions in the Space $L_2([t,T])$ and Binomial Weight Functions}

Let us prove the following theorem.

\vspace{2mm}

{\bf Theorem~2.63.}\ {\it Suppose that
$\{\phi_j(x)\}_{j=0}^{\infty}$ is an arbitrary complete orthonormal system of 
functions in the space $L_2([t,T]).$
Then$,$ for the iterated Stra\-to\-no\-vich stochastic integral
of fourth multiplicity 
$$
I_{{l_1l_2l_3 l_4}_{T,t}}^{*(i_1i_2i_3 i_4)}=
{\int\limits_t^{*}}^T (t_4-t)^{l_4}{\int\limits_t^{*}}^{t_4} (t_3-t)^{l_3}
{\int\limits_t^{*}}^{t_3}(t_2-t)^{l_2}
{\int\limits_t^{*}}^{t_2}(t_1-t)^{l_1}\times
$$
$$
\times
d{\bf w}_{t_1}^{(i_1)}
d{\bf w}_{t_2}^{(i_2)}d{\bf w}_{t_3}^{(i_3)}d{\bf w}_{t_4}^{(i_4)}
$$
the following expansion 
$$
I_{{l_1l_2l_3 l_4}_{T,t}}^{*(i_1i_2i_3i_4)}=
\hbox{\vtop{\offinterlineskip\halign{
\hfil#\hfil\cr
{\rm l.i.m.}\cr
$\stackrel{}{{}_{p\to \infty}}$\cr
}} }\sum_{j_1,j_2,j_3,j_4=0}^{p}
C_{j_4 j_3 j_2 j_1}\zeta_{j_1}^{(i_1)}\zeta_{j_2}^{(i_2)}\zeta_{j_3}^{(i_3)}\zeta_{j_4}^{(i_4)}
$$
that converges in the mean-square sense is valid, where 
$i_1,i_2,i_3,i_4=0,1,\ldots,m;$ $l_1,l_2,l_3,l_4=0,1,2,\ldots,$
$$
C_{j_4 j_3 j_2 j_1}=\int\limits_t^T
(t_4-t)^{l_4}\phi_{j_4}(t_4)\int\limits_t^{t_4}
(t_3-t)^{l_3}\phi_{j_3}(t_3)\int\limits_t^{t_3}
(t_2-t)^{l_2}
\phi_{j_2}(t_2)
\int\limits_t^{t_2}
(t_1-t)^{l_1}\phi_{j_1}(t_1)\times
$$
\begin{equation}
\label{2024decem1}
\times
dt_1dt_2dt_3dt_4
\end{equation}
and
$$
\zeta_{j}^{(i)}=
\int\limits_t^T \phi_{j}(\tau) d{\bf w}_{\tau}^{(i)}
$$ 
are independent standard Gaussian random variables for various 
$i$ or $j$ {\rm (}in the case when $i\ne 0${\rm ),}
${\bf w}_{\tau}^{(i)}$ 
$(i=1,\ldots,m)$ are independent 
standard Wiener processes$,$
${\bf w}_{\tau}^{(0)}=\tau.$}

\vspace{2mm}

{\bf Proof.}\ The following proof will be based on Theorem~2.61 
and verification of the equality (\ref{09091}) under the conditions
of Theorem~2.63 (the case $k=4>2r$, where $r=1$). Note that the case
$k=2r$ is proved in Sect.~2.27.4 (see (\ref{july90000})).
Under the conditions of Theorem~2.63, the equality $k=2r$ means that
$k=4$ and $r=2$.

Let throughout this proof 
$$
C_{j_1 j_1}^{\psi_{i+1}\psi_i}(s,\tau)\hspace{-0.3mm} =\hspace{-0.9mm}\int\limits_{\tau}^s \hspace{-0.8mm}
\psi_{i+1}(y)\phi_{j_1}(y)
\hspace{-0.8mm}\int\limits_{\tau}^{y}\hspace{-0.8mm}
\psi_i(x)\phi_{j_1}(x)dx dy,\ \ 
C_{j_1}^{\psi_q}(s,\tau)\hspace{-0.3mm}=\hspace{-0.9mm}\int\limits_{\tau}^s \hspace{-0.8mm}
\psi_{q}(x)\phi_{j_1}(x)dx,
$$
where $i=1,2,3,$ $t\le\tau<s\le T,$ $\psi_q(x)=(x-t)^{l_q},$ $l_q=0,1,2,\ldots,$ 
$q=1,\ldots,4,$ $x\in [t, T],$
and $C_{j_4 j_3 j_2 j_1}$ is defined by (\ref{2024decem1}).

Using Fubini's Theorem and the technique that leads to the formulas (\ref{copa1}),
(\ref{copa1a}),
we obtain (note that we find all possible
combinations of pairs using the equality (\ref{after34})):
$$
C_{j_4 j_3 j_1 j_1}=
\int\limits_t^T
\psi_4(t_4)\phi_{j_4}(t_4)\int\limits_t^{t_4}
\psi_3(t_3)\phi_{j_3}(t_3)C_{j_1 j_1}^{\psi_2 \psi_1}(t_3,t)dt_3 dt_4,
$$
$$
C_{j_4 j_1 j_2 j_1}=
\int\limits_t^T
\psi_4(t_4)\phi_{j_4}(t_4)\int\limits_t^{t_4}
\psi_2(t_2)\phi_{j_2}(t_2)C_{j_1}^{\psi_1}(t_2,t)C_{j_1}^{\psi_3}(t_4,t_2)dt_2 dt_4,
$$
$$
C_{j_1 j_3 j_2 j_1}=
\int\limits_t^T
\psi_3(t_3)\phi_{j_3}(t_3)\int\limits_t^{t_3}
\psi_2(t_2)\phi_{j_2}(t_2)C_{j_1}^{\psi_1}(t_2,t)C_{j_1}^{\psi_4}(T,t_3)dt_2 dt_3,
$$

$$
C_{j_4 j_2 j_2 j_1}=
\int\limits_t^T
\psi_4(t_4)\phi_{j_4}(t_4)\int\limits_t^{t_4}
\psi_1(t_1)\phi_{j_1}(t_1)C_{j_2 j_2}^{\psi_3 \psi_2}(t_4,t_1)dt_1 dt_4,
$$
$$
C_{j_2 j_3 j_2 j_1}=
\int\limits_t^T
\psi_3(t_3)\phi_{j_3}(t_3)\int\limits_t^{t_3}
\psi_1(t_1)\phi_{j_1}(t_1)C_{j_2}^{\psi_2}(t_3,t_1) C_{j_2}^{\psi_4}(T,t_3)dt_1 dt_3,
$$
$$
C_{j_3 j_3 j_1 j_1}=
\int\limits_t^T
\psi_2(t_2)\phi_{j_2}(t_2)\int\limits_t^{t_2}
\psi_1(t_1)\phi_{j_1}(t_1)C_{j_3 j_3}^{\psi_4 \psi_3}(T,t_2)dt_1 dt_2.
$$

\vspace{2mm}

It is easy to see (based on the above equalities)
that the condition (\ref{09091}) will be satisfied under the conditions
of Theorem~2.63 if

\vspace{-2mm}
\begin{equation}
\label{2024decem12}
\left\vert \sum\limits_{j_1=0}^p 
C_{j_1 j_1}^{\psi_{i+1} \psi_i}(s,\tau)\right\vert\le K,
\end{equation}

\vspace{-2mm}
\begin{equation}
\label{2024decem13}
\left\vert \sum\limits_{j_1=0}^p 
C_{j_1}^{\psi_k}(s,\tau)C_{j_1}^{\psi_q}(\theta,u)\right\vert\le K,
\end{equation}

\vspace{2mm}
\noindent
where $p\in{\bf N},$ $i=1,2,3,$ $k,q=1,\ldots,4,$ $t\le \tau < s \le T,$ $t\le u<\theta \le T,$
constant $K$ does not depend on 
$p, s, \tau, u, \theta$ (but only on $t, T$). 

The equality (\ref{2024decem12}) has been proved earlier (see (\ref{may103x})).
Obviously,  the relation (\ref{2024decem13}) is proved in complete
analogy with (\ref{may106}).

Thus, the condition (\ref{09091}) of Theorem~2.61 is fulfilled
under the conditions of Theorem~2.63. Then
Theorem~2.63 follows from Theorem~2.61. 
Theorem~2.63 is proved.

\section{Another Proof of Theorem~2.50 Based on Theorem~2.61}

The following proof will be based on Theorem~2.61 
and verification of the equality (\ref{09091}) under the conditions
of Theorem~2.50 (the case $k=5>2r$, where $r=1$ or $r=2$). 

Further, suppose that
$$
C_{j_k \ldots j_1}(s,\tau)=\int\limits_{\tau}^s
\phi_{j_k}(t_k)\ldots
\int\limits_{\tau}^{t_2}
\phi_{j_1}(t_1)dt_1\ldots dt_k, 
$$

\noindent
where $k=1,\ldots,4,$ $t\le\tau<s\le T$, and 
$$
C_{j_5\ldots j_1}=\int\limits_t^T
\phi_{j_5}(t_5)
\ldots
\int\limits_t^{t_2}
\phi_{j_1}(t_1)dt_1\ldots dt_5.
$$

Applying the technique that leads to (\ref{copa1}), we obtain 
(note that we find all possible
combinations of pairs using the equality (\ref{after35}))
$$
C_{j_5 j_4 j_3 j_1 j_1}
=\int\limits_t^T\phi_{j_5}(t_5)\int\limits_t^{t_5}\phi_{j_4}(t_4)
\int\limits_t^{t_4}\phi_{j_3}(t_3)C_{j_1 j_1}(t_3,t)
dt_3 dt_4 dt_5,
$$
$$
C_{j_5 j_4 j_1 j_2 j_1}=
\int\limits_t^T\phi_{j_5}(t_5)\int\limits_t^{t_5}\phi_{j_4}(t_4)
\int\limits_t^{t_4}\phi_{j_2}(t_2)
C_{j_1}(t_2,t)C_{j_1}(t_4,t_2)
dt_2 dt_4 dt_5,
$$
$$
C_{j_5 j_1 j_3 j_2 j_1}=
\int\limits_t^T\phi_{j_5}(t_5)\int\limits_t^{t_5}\phi_{j_3}(t_3)
\int\limits_t^{t_3}\phi_{j_2}(t_2)
C_{j_1}(t_2,t)C_{j_1}(t_5,t_3)
dt_2 dt_3 dt_5,
$$
$$
C_{j_1 j_4 j_3 j_2 j_1}=
\int\limits_t^T\phi_{j_4}(t_4)\int\limits_t^{t_4}\phi_{j_3}(t_3)
\int\limits_t^{t_3}\phi_{j_2}(t_2)
C_{j_1}(t_2,t)C_{j_1}(T,t_4)
dt_2 dt_3 dt_4,
$$
$$
C_{j_5 j_4 j_2 j_2 j_1}=
\int\limits_t^T\phi_{j_5}(t_5)\int\limits_t^{t_5}\phi_{j_4}(t_4)
\int\limits_t^{t_4}\phi_{j_1}(t_1)
C_{j_2 j_2}(t_4,t_1)
dt_1 dt_4 dt_5,
$$
$$
C_{j_5 j_2 j_3 j_2 j_1}=
\int\limits_t^T\phi_{j_5}(t_5)\int\limits_t^{t_5}\phi_{j_3}(t_3)
\int\limits_t^{t_3}\phi_{j_1}(t_1)
C_{j_2}(t_3,t_1)C_{j_2}(t_5,t_3)
dt_1 dt_3 dt_5,
$$
$$
C_{j_2 j_4 j_3 j_2 j_1}=
\int\limits_t^T\phi_{j_4}(t_4)\int\limits_t^{t_4}\phi_{j_3}(t_3)
\int\limits_t^{t_3}\phi_{j_1}(t_1)
C_{j_2}(t_3,t_1)C_{j_2}(T,t_4)
dt_1 dt_3 dt_4,
$$
$$
C_{j_5 j_3 j_3 j_2 j_1}
\int\limits_t^T\phi_{j_5}(t_5)\int\limits_t^{t_5}\phi_{j_2}(t_2)
\int\limits_{t}^{t_2}\phi_{j_1}(t_1)
C_{j_3 j_3}(t_5,t_2)
dt_1 dt_2 dt_5,
$$
$$
C_{j_3 j_4 j_3 j_2 j_1}=
\int\limits_t^T\phi_{j_4}(t_4)\int\limits_{t}^{t_4}\phi_{j_2}(t_2)
\int\limits_{t}^{t_2}\phi_{j_1}(t_1)
C_{j_3}(t_4,t_2)C_{j_3}(T,t_4)
dt_1 dt_2 dt_4,
$$
$$
C_{j_4 j_4 j_3 j_2 j_1}=
\int\limits_t^T\phi_{j_3}(t_3)\int\limits_t^{t_3}\phi_{j_2}(t_2)
\int\limits_{t}^{t_2}\phi_{j_1}(t_1) 
C_{j_4 j_4}(T,t_3)dt_1 dt_2dt_3,
$$
$$
C_{j_5 j_3 j_3 j_1 j_1}=
\int\limits_t^T \phi_{j_5}(t_5)
C_{j_3 j_3 j_1 j_1}(t_5,t)dt_5,
$$
$$
C_{j_5 j_2 j_1 j_2 j_1}=
\int\limits_t^T \phi_{j_5}(t_5) 
C_{j_2 j_1 j_2 j_1} (t_5,t)dt_5,
$$
$$
C_{j_5 j_1 j_2 j_2 j_1}=
\int\limits_t^T \phi_{j_5}(t_5)
C_{j_1 j_2 j_2 j_1} (t_5,t)dt_5,
$$
$$
C_{j_4 j_4 j_2 j_2 j_1}=
\int\limits_t^T  \phi_{j_1}(t_1)
C_{j_4 j_4 j_2 j_2}(T,t_1)dt_1,
$$
$$
C_{j_3 j_2 j_3 j_2 j_1}=
\int\limits_t^T  \phi_{j_1}(t_1)
C_{j_3 j_2 j_3 j_2}(T,t_1)dt_1,
$$
$$
C_{j_2 j_3 j_3 j_2 j_1}=
\int\limits_t^T  \phi_{j_1}(t_1)
C_{j_2 j_3 j_3 j_2}(T,t_1)dt_1,
$$
$$
C_{j_4 j_4 j_3 j_1 j_1}
=\int\limits_t^T  \phi_{j_3}(t_3)
C_{j_1 j_1}(t_3,t)C_{j_4 j_4}(T,t_3)
dt_3,
$$
$$
C_{j_2 j_4 j_1 j_2 j_1}=
\int\limits_t^T \phi_{j_4}(t_4)
C_{j_1 j_2 j_1}(t_4,t)
C_{j_2}(T,t_4) dt_4,
$$
$$
C_{j_2 j_1 j_3 j_2 j_1}=
\int\limits_t^T
\phi_{j_3}(t_3)C_{j_2 j_1}(t_3,t)
C_{j_2 j_1}(T,t_3)
dt_3,
$$
$$
C_{j_3 j_1 j_3 j_2 j_1}=\int\limits_t^T
\phi_{j_2}(t_2)C_{j_1}(t_2,t)
C_{j_3 j_1 j_3}(T,t_2)
dt_2,
$$
$$
C_{j_1 j_2 j_3 j_2 j_1}=
\int\limits_t^T
\phi_{j_3}(t_3)
C_{j_2 j_1}(t_3,t)C_{j_1 j_2}(T,t_3)dt_3,
$$

\vspace{-2.5mm}
$$
C_{j_3 j_4 j_3 j_1 j_1}=
\int\limits_t^T  
\phi_{j_4}(t_4)C_{j_3 j_1 j_1}(t_4,t)C_{j_3}(T,t_4)dt_4,
$$

\vspace{-2.5mm}
$$
C_{j_4 j_4 j_1 j_2 j_1}=
\int\limits_t^T  \phi_{j_2}(t_2)C_{j_1}(t_2,t)
C_{j_4 j_4 j_1}(T,t_2)dt_2,
$$

\vspace{-2.5mm}
$$
C_{j_1 j_4 j_2 j_2 j_1}=
\int\limits_t^T \phi_{j_4}(t_4)
C_{j_2 j_2 j_1}(t_4,t)C_{j_1}(T,t_4)dt_4,
$$

\vspace{-2.5mm}
$$
C_{j_1 j_3 j_3 j_2 j_1}=
\int\limits_t^T  
\phi_{j_2}(t_2)
C_{j_1}(t_2,t) C_{j_1 j_3 j_3}(T,t_2)
dt_2.
$$

It is easy to see (based on the above relations)
that (\ref{09091}) will be satisfied (under the conditions
of Theorem~2.50) if 
(\ref{2024december12})--(\ref{2024december9})
are fulfilled.
The equalities 
(\ref{2024december12})--(\ref{2024december9})
are proved in Sect.~2.32.
The assertion
of Theorem~2.50 now follows from Theorem~2.61. 
Theorem~2.50 is proved.

Recall that for the case $k=6$, together with
(\ref{2024december12})--(\ref{2024december9}), the conditions 
(\ref{2024december10}), (\ref{2024december11}) and the equality 
(\ref{july90000}) ($k=2r,$ $k=6,$ $r=3$)  
must be satisfied (see the proof of Theorem~2.62).

\section{Partial Proof of the Condition (\ref{09091})}

In this section, we will prove (\ref{09091})
for the case when the condition $(A)$
and the relation (\ref{copa5}) are satisfied (see Sect.~2.31).

Suppose that $\{\phi_j(x)\}_{j=0}^{\infty}$
is an arbitrary complete orthonormal system of functions
in $L_2([t, T])$ and $\psi_1(\tau),\ldots,\psi_k(\tau)\equiv 1.$

It is easy to see that (\ref{09091}) will be proved
for the above case if we prove that
\begin{equation}
\label{febr2025}
\left|\sum_{j_r,j_{r-2},\ldots, j_2=0}^{p}
C_{j_r j_r j_{r-2} j_{r-2} \ldots j_2 j_2}(s,\tau)\right|\le K <\infty,
\end{equation}

\noindent
where $p\in{\bf N},$ $r=2, 4, 6,\ldots,$ constant $K$ does not depend on $p, s, \tau$
(but only on $t, T$),

\newpage
\noindent
\begin{equation}
\label{utoch1}
C_{j_k \ldots j_1}(s,\tau)=\int\limits_{\tau}^s
\phi_{j_k}(t_k)\ldots
\int\limits_{\tau}^{t_2}
\phi_{j_1}(t_1)dt_1\ldots dt_k,
\end{equation}

\noindent
where $k\in {\bf N},$ $t\le\tau<s\le T$.

By analogy with (\ref{july7028}) we obtain

\vspace{-2mm}
$$
C_{j_r j_r j_{r-2}j_{r-2}\ldots j_2 j_2}(s,\tau) +
C_{j_2 j_2\ldots j_{r-2}j_{r-2} j_r j_r}(s,\tau)=
$$

\vspace{-2mm}
$$
=
C_{j_r}(s,\tau) \cdot  C_{j_{r} j_{r-2} j_{r-2} \ldots j_4 j_4 j_2 j_2}(s,\tau)
-C_{j_{r} j_r}(s,\tau) \cdot C_{j_{r-2} j_{r-2}\ldots j_4 j_4 j_2 j_2}(s,\tau)+
$$

\vspace{-2mm}
$$
+C_{j_{r-2} j_{r} j_r}(s,\tau) \cdot
C_{j_{r-2} j_{r-4}j_{r-4}\ldots j_4 j_4 j_2 j_2}(s,\tau)
- \ldots 
$$

\vspace{-2mm}
\begin{equation}
\label{febr2025a}
- C_{j_4 j_4 \ldots j_{r-2}j_{r-2} j_{r}j_{r}}(s,\tau) \cdot C_{j_2 j_2}(s,\tau)+
C_{j_2 j_4 j_4 \ldots j_{r-2}j_{r-2} j_r j_r}(s,\tau) \cdot C_{j_2}(s,\tau).
\end{equation}

\vspace{3mm}

Applying (\ref{febr2025a}), we get
$$
2 \sum_{j_r,j_{r-2},\ldots, j_4, j_2=0}^{p} C_{j_r j_r j_{r-2}j_{r-2}\ldots j_4 j_4 j_2 j_2}(s,\tau)=
$$

\vspace{0.5mm}
$$
=
\sum_{j_r=0}^p  C_{j_r}(s,\tau)\sum_{j_{r-2},\ldots, j_4, j_2=0}^{p}
 C_{j_{r} j_{r-2} j_{r-2} \ldots j_4 j_4 j_2 j_2}(s,\tau)
-
$$

\vspace{0.5mm}
$$
-\sum_{j_r=0}^{p} C_{j_{r} j_r}(s,\tau) \sum_{j_{r-2},\ldots, j_4, j_2=0}^{p}
C_{j_{r-2} j_{r-2}\ldots j_4 j_4 j_2 j_2}(s,\tau)+
$$

\vspace{0.5mm}
$$
+  \sum_{j_{r-2}=0}^{p} \sum_{j_{r}=0}^{p} C_{j_{r-2} j_{r} j_r}(s,\tau) 
\sum_{j_{r-4},\ldots, j_4, j_2=0}^{p} C_{j_{r-2} j_{r-4}j_{r-4}\ldots j_4 j_4 j_2 j_2}(s,\tau)
- \ldots 
$$

\vspace{0.5mm}
$$
- \sum_{j_r,j_{r-2},\ldots, j_4=0}^{p}
C_{j_4 j_4 \ldots j_{r-2}j_{r-2} j_{r}j_{r}}(s,\tau) \sum_{j_{2}=0}^{p}C_{j_2 j_2}(s,\tau)+
$$

\vspace{0.5mm}
\begin{equation}
\label{febr2025b}
~~~~~~~~+
\sum_{j_{2}=0}^{p} \sum_{j_r,j_{r-2},\ldots, j_4=0}^{p} 
C_{j_2 j_4 j_4\ldots j_{r-2}j_{r-2} j_r j_r}(s,\tau)\cdot  C_{j_2}(s,\tau).
\end{equation}

\vspace{3mm}

Let us prove (\ref{febr2025}) by induction.
The equality (\ref{febr2025}) is proved for $r=2, 4$ (see 
(\ref{nahod0}), (\ref{april48}), (\ref{april50})).
Suppose that
\begin{equation}
\label{febr2025b1}
\left|\sum_{j_6, j_4, j_2=0}^{p}
C_{j_6 j_6 j_{4} j_{4} j_2 j_2}(s,\tau)\right|\le K <\infty,
\end{equation}
\begin{equation}
\label{febr2025b2}
\left|\sum_{j_8, j_6, j_4, j_2=0}^{p}
C_{j_8 j_8 j_6 j_6 j_{4} j_{4} j_2 j_2}(s,\tau)\right|\le K <\infty,
\end{equation}

\vspace{-2mm}
$$
\ldots
$$
\begin{equation}
\label{febr2025b3}
~~~~~~~~~~~\left|\sum_{j_{r-2},j_{r-4},\ldots, j_2=0}^{p}
C_{j_{r-2} j_{r-2} j_{r-4} j_{r-4} \ldots j_2 j_2}(s,\tau)\right|\le K <\infty
\end{equation}

\vspace{2mm}
\noindent
and prove (\ref{febr2025}).

Using the induction hypothesis (see (\ref{febr2025b1})--(\ref{febr2025b3})), we obtain
\begin{equation}
\label{febr2025b4}
~~~~~~~~~~~\left|\sum_{j_r=0}^{p} C_{j_{r} j_r}(s,\tau) \sum_{j_{r-2},\ldots, j_4, j_2=0}^{p}
C_{j_{r-2} j_{r-2}\ldots j_4 j_4 j_2 j_2}(s,\tau)\right|\le K^2<\infty,
\end{equation}
\begin{equation}
\label{febr2025b5}
\left|
\sum_{j_r, j_{r-2}=0}^{p} C_{j_{r-2} j_{r-2} j_{r} j_{r}}(s,\tau) \sum_{j_{r-4},\ldots, j_4, j_2=0}^{p}
C_{j_{r-4} j_{r-4}\ldots j_4 j_4 j_2 j_2}(s,\tau)\right|\le K^2<\infty,
\end{equation}
$$
\ldots
$$

\vspace{-9mm}
\begin{equation}
\label{febr2025b6}
~~~~~~~~~~\left|\sum_{j_r,j_{r-2},\ldots, j_4=0}^{p}
C_{j_4 j_4 \ldots j_{r-2}j_{r-2} j_{r}j_{r}}(s,\tau) \sum_{j_{2}=0}^{p}C_{j_2 j_2}(s,\tau)
\right|\le K^2<\infty.
\end{equation}

\vspace{3mm}

Applying the inequality 
of Cauchy--Bunyakovsky, Parseval's equality and the induction hypothesis, we obtain
$$
\left(
\sum_{j_r=0}^p  C_{j_r}(s,\tau)\sum_{j_{r-2},\ldots, j_4, j_2=0}^{p}
C_{j_{r} j_{r-2} j_{r-2} \ldots j_4 j_4 j_2 j_2}(s,\tau)\right)^2\le
$$
$$
\le \sum_{j_r=0}^p  \left(C_{j_r}(s,\tau)\right)^2
\sum_{j_r=0}^p \left(
\sum_{j_{r-2},\ldots, j_4, j_2=0}^{p}
C_{j_{r} j_{r-2} j_{r-2} \ldots j_4 j_4 j_2 j_2}(s,\tau)\right)^2\le
$$
$$
\le \sum_{j_r=0}^{\infty}  \left(C_{j_r}(s,\tau)\right)^2
\sum_{j_r=0}^{\infty} \left(
\sum_{j_{r-2},\ldots, j_4, j_2=0}^{p}
C_{j_{r} j_{r-2} j_{r-2} \ldots j_4 j_4 j_2 j_2}(s,\tau)\right)^2\le
$$
$$
\le K_1 
\sum_{j_r=0}^{\infty} \left(
\sum_{j_{r-2},\ldots, j_4, j_2=0}^{p}
C_{j_{r} j_{r-2} j_{r-2} \ldots j_4 j_4 j_2 j_2}(s,\tau)\right)^2=
$$
$$
=K_1 
\sum_{j_r=0}^{\infty} \left(\int\limits_{\tau}^s \phi_{j_r}(u)
\sum_{j_{r-2},\ldots, j_4, j_2=0}^{p}
C_{j_{r-2} j_{r-2} \ldots j_4 j_4 j_2 j_2}(u,\tau)du\right)^2=
$$
$$
=K_1 
\int\limits_{\tau}^s \left(
\sum_{j_{r-2},\ldots, j_4, j_2=0}^{p}
C_{j_{r-2} j_{r-2} \ldots j_4 j_4 j_2 j_2}(u,\tau)\right)^2 du \le
$$
\begin{equation}
\label{febr2025b7}
\le K_1 K^2 
\int\limits_{\tau}^s du \le (T-t)K_1 K^2=K_2<\infty,
\end{equation}

\noindent
where constant $K_2$ does not depend on $p, s, \tau;$

\vspace{-3mm}

$$
\left(\sum_{j_{r-2}=0}^{p} \sum_{j_{r}=0}^{p} C_{j_{r-2} j_{r} j_r}(s,\tau) 
\sum_{j_{r-4},\ldots, j_4, j_2=0}^{p} C_{j_{r-2} j_{r-4}j_{r-4}\ldots 
j_4 j_4 j_2 j_2}(s,\tau)\right)^2\le
$$
$$
\le \sum_{j_{r-2}=0}^{p} \left(\sum_{j_{r}=0}^{p} C_{j_{r-2} j_{r} j_r}(s,\tau)\right)^2 
\sum_{j_{r-2}=0}^{p}
\left(\sum_{j_{r-4},\ldots, j_4, j_2=0}^{p} C_{j_{r-2} j_{r-4}j_{r-4}\ldots 
j_4 j_4 j_2 j_2}(s,\tau)\right)^2\le
$$
$$
\le \sum_{j_{r-2}=0}^{\infty} \left(\sum_{j_{r}=0}^{p} C_{j_{r-2} j_{r} j_r}(s,\tau)\right)^2 
\sum_{j_{r-2}=0}^{\infty}
\left(\sum_{j_{r-4},\ldots, j_4, j_2=0}^{p} C_{j_{r-2} j_{r-4}j_{r-4}\ldots 
j_4 j_4 j_2 j_2}(s,\tau)\right)^2=
$$
$$
=\sum_{j_{r-2}=0}^{\infty} \left(
\int\limits_{\tau}^s \phi_{j_{r-2}}(u)
\sum_{j_{r}=0}^{p} C_{j_{r} j_r}(u,\tau) du\right)^2 \times
$$
$$
\times
\sum_{j_{r-2}=0}^{\infty}
\left(\int\limits_{\tau}^s \phi_{j_{r-2}}(u)
\sum_{j_{r-4},\ldots, j_4, j_2=0}^{p} C_{j_{r-4}j_{r-4}\ldots j_4 j_4 j_2 j_2}(u,\tau) du\right)^2=
$$
$$
=
\int\limits_{\tau}^s \left(
\sum_{j_{r}=0}^{p} C_{j_{r} j_r}(u,\tau)\right)^2 du \times
$$
\begin{equation}
\label{febr2025b8}
\times
\int\limits_{\tau}^s \left(
\sum_{j_{r-4},\ldots, j_4, j_2=0}^{p} C_{j_{r-4}j_{r-4}\ldots j_4 j_4 j_2 j_2}(u,\tau)\right)^2 du\le
K^4 (T-t)^2 = K_3 <\infty.
\end{equation}

\vspace{3mm}

Similarly, we get
$$
\left(\sum_{j_{r-4}=0}^{p} \sum_{j_{r}, j_{r-2}=0}^{p} C_{j_{r-4} j_{r-2} j_{r-2} j_r j_r}(s,\tau) 
\sum_{j_{r-6},\ldots, j_4, j_2=0}^{p} C_{j_{r-4} j_{r-6}j_{r-6}\ldots j_4 j_4 j_2 j_2}(s,\tau)\right)^2
\le 
$$
\begin{equation}
\label{febr2025b9}
\le K_4 <
\infty,
\end{equation}

\vspace{-8mm}
$$
\ldots
$$
\begin{equation}
\label{febr2025b10}
\left(
\sum_{j_{4}=0}^{p} \sum_{j_r,j_{r-2},\ldots, j_6=0}^{p} 
C_{j_4 j_6 j_6\ldots j_{r-2}j_{r-2} j_r j_r}(s,\tau)
\sum_{j_{2}=0}^{p} C_{j_4 j_2 j_2}(s,\tau)\right)^2
\le K_4 <\infty,
\end{equation}
\begin{equation}
\label{febr2025b11}
~~~~~~~\left(
\sum_{j_{2}=0}^{p} \sum_{j_r,j_{r-2},\ldots, j_4=0}^{p} 
C_{j_2 j_4 j_4\ldots j_{r-2}j_{r-2} j_r j_r}(s,\tau)\cdot  C_{j_2}(s,\tau)\right)^2
\le K_4 <\infty,
\end{equation}

\vspace{3mm}
\noindent
where constant $K_4$ does not depend on $p, s, \tau.$

Combining (\ref{febr2025b}), (\ref{febr2025b4})--(\ref{febr2025b6}),
(\ref{febr2025b7}), (\ref{febr2025b8}),
(\ref{febr2025b9})--(\ref{febr2025b11}),
we obtain (\ref{febr2025}).
The equality (\ref{09091}) is proved
for the case when the condition $(A)$
and the relation (\ref{copa5}) are satisfied.

\section{Further Development of the Approach
Based on Theorem~2.61 for the Case $\psi_1(\tau),\ldots, \psi_7(\tau)
\equiv 1$. Expansion of Iterated Stratonovich Stochastic Integrals
of Multiplicity 7 (The Cases of Legendre 
Polynomials and Trigonometric Functions)}

Unfortunately, the approach from the previous section 
can be generalized only partially to the case
when the condition $(A)$
and the relation (\ref{copa6}) are satisfied (see Sect.~2.31).
In particular, the mentioned approach
is applicable to the proof of inequality
$$
\left|\sum\limits_{j_1,j_2,j_3=0}^p C_{j_3 j_2 j_1 j_3 j_2 j_1}(s,\tau)\right|\le K<\infty,
$$
but is not applicable to the proof of inequality
$$
\left|\sum\limits_{j_1,j_2,j_3=0}^p C_{j_2 j_3 j_3 j_1 j_2 j_1}(s,\tau)\right|\le K<\infty,
$$
where $C_{j_k \ldots j_1}(s,\tau)$ is defined by (\ref{utoch1}),
constant $K$ does not depend on $p, s, \tau$
$(p\in {\bf N},\ t\le \tau< s\le T).$

In this section, we will restrict ourselves to the case
$k=7,$ $r=1,2,3$ and we will also assume that 
$\{\phi_j(x)\}_{j=0}^{\infty}$ is a complete orthonormal
system of Legendre polynomials or trigonometric functions
in the space $L_2([t, T])$.

Recall that the condition (\ref{09091})
can be weakened. Namely, the constant $K^2$
can be replaced by the function $F$ such that
$\psi_{q_1}^2\ldots \psi_{q_{k-2r}}^2 F \in L_1([t, T]^{k-2r})$
(see (\ref{09091xxx})).
For the trigonometric case, we will prove (\ref{09091}) for $k=7,$ $r=1,2,3$.
For the polynomial case, we will prove a weakened version of
(\ref{09091}) for $k=7,$ $r=1,2,3$ (the constant $K$ and the above function 
$F$ will be used in the weakened version of 
(\ref{09091})).

Obviously, that the conditions
(\ref{2024december12})--(\ref{2024december11})
together with the following condition
\begin{equation}
\label{cc123}
~~~~~~~~~~\left\vert \sum\limits_{j_1, j_2=0}^p C_{j_1}(s,\tau)C_{j_2}(\rho,v)
C_{j_1}(\theta,u)
C_{j_2}(\mu,w)\right\vert\le K
\end{equation}

\noindent
cover the case $k=7,$ $r=1,2$ (see (\ref{09091})),
where $p\in{\bf N},$ $t\le \tau < s \le T,$ $t\le u<\theta \le T,$
$t\le v<\rho \le T,$ $t\le w<\mu \le T,$ constant $K$ does not depend on 
$p, s, \tau, u, \theta, v, \rho, w, \mu$ (but only on $t, T$).
The inequality (\ref{cc123}) is easily verified
using (\ref{dsds14fffff}).

Now let us focus on the proof of (\ref{09091})
for the case $k=7$ and $r=3$. So, we need to prove that
\begin{equation}
\label{march0001}
\left|\sum\limits_{j_{g_1},j_{g_3},j_{g_{5}}=0}^p
C_{j_{d_1} j_{d_1-1}j_{d_1-2}j_{d_1-3}j_{d_1-4}j_{d_1-5}}(s,\tau)
\biggl|_{j_{g_1}=j_{g_2},j_{g_3}=j_{g_4},j_{g_5}=j_{g_6}}\right|\le K<\infty,
\end{equation}
\begin{equation}
\label{march0002}
\left|\sum\limits_{j_{g_1},j_{g_3},j_{g_{5}}=0}^p
\bigl(C_{j_{d_2} j_{d_2-1}j_{d_2-2}j_{d_2-3}j_{d_2-4}}(s,\tau)
C_{j_{d_1}}(\theta,u)
\bigr)\biggl|_{j_{g_1}=j_{g_2},
j_{g_3}=j_{g_4},j_{g_5}=j_{g_6}}\right|\le K<\infty,
\end{equation}
\begin{equation}
\label{march0003}
\left|\sum\limits_{j_{g_1},j_{g_3},j_{g_{5}}=0}^p
\left(C_{j_{d_2} j_{d_2-1}j_{d_2-2}j_{d_2-3}}(s,\tau)
C_{j_{d_1}j_{d_1-1}}(\theta,u)\right)\biggl|_{j_{g_1}=j_{g_2},
j_{g_3}=j_{g_4},j_{g_5}=j_{g_6}}\right|\le K<\infty,
\end{equation}
\begin{equation}
\label{march0004}
\left|\sum\limits_{j_{g_1},j_{g_3},j_{g_{5}}=0}^p
\left(C_{j_{d_2} j_{d_2-1}j_{d_2-2}}(s,\tau)
C_{j_{d_1} j_{d_1-1}j_{d_1-2}}(\theta,u)\right)\biggl|_{j_{g_1}=j_{g_2},
j_{g_3}=j_{g_4},j_{g_5}=j_{g_6}}\right|\le K<\infty,
\end{equation}
where $p\in{\bf N},$ $t\le \tau < s \le T,$ $t\le u<\theta \le T,$
constant $K$ does not depend on 
$p, s, \tau, u, \theta$ (but only on $t, T$) and may differ from line to line;
another notations are the same as in Sect.~2.31, 2.35.

The inequalities (\ref{march0002})--(\ref{march0004})
are proved using the same technique as 
inequalities (\ref{2024december12})--(\ref{2024december11}) (see Sect.~2.32).
Here we will only prove as an example the following
special case of the inequality (\ref{march0003})
\begin{equation}
\label{march0005}
\left|\sum\limits_{j_1, j_2, j_3=0}^p C_{j_2 j_3 j_2 j_1}(s,\tau)C_{j_3 j_1}(\theta,u)\right|
\le K<\infty.
\end{equation}

Using the 
Cau\-chy--Bunyakovsky inequality as well as 
Fubini's Theorem, Parseval's equality and (\ref{2024december13}), we have
$$
\left(\sum\limits_{j_1, j_2, j_3=0}^p C_{j_2 j_3 j_2 j_1}(s,\tau)C_{j_3 j_1}(\theta,u)\right)^2\le
$$
$$
\le
\sum\limits_{j_1,j_3=0}^p \left(\sum\limits_{j_2=0}^p C_{j_2 j_3 j_2 j_1}(s,\tau)\right)^2
\sum\limits_{j_1, j_3=0}^p C_{j_3 j_1}^2(\theta,u)\le
$$
$$
\le \sum\limits_{j_1,j_3=0}^{\infty} \left(\sum\limits_{j_2=0}^p 
\int\limits_{\tau}^s \phi_{j_2}(u)
\int\limits_{\tau}^u \phi_{j_3}(z)
\int\limits_{\tau}^{z} \phi_{j_2}(y)
\int\limits_{\tau}^{y} \phi_{j_1}(x)dx dy dz du\right)^2 \times
$$
$$
\times
\sum\limits_{j_1,j_3=0}^{\infty} C_{j_3 j_1}^2(\theta,u)=
$$
$$
=\sum\limits_{j_1,j_3=0}^{\infty} \left(\sum\limits_{j_2=0}^p 
\int\limits_{\tau}^s 
\phi_{j_3}(z)
\int\limits_{\tau}^{z} \phi_{j_2}(y)
\int\limits_{\tau}^{y} \phi_{j_1}(x)dx dy
\int\limits_z^s 
\phi_{j_2}(u)
du dz\right)^2 
\cdot \frac{(\theta-u)^2}{2}=
$$
$$
=\frac{(\theta-u)^2}{2}\sum\limits_{j_1,j_3=0}^{\infty} \left(\sum\limits_{j_2=0}^p 
\int\limits_{\tau}^s 
\phi_{j_3}(z)
\int\limits_{\tau}^{z} \phi_{j_1}(x) \int\limits_{x}^{z} \phi_{j_2}(y)
dy dx
\int\limits_z^s 
\phi_{j_2}(u)
du dz\right)^2 =
$$
$$
=\frac{(\theta-u)^2}{2}\sum\limits_{j_1,j_3=0}^{\infty} \left(
\int\limits_{\tau}^s 
\phi_{j_3}(z)
\int\limits_{\tau}^{z} \phi_{j_1}(x) \sum\limits_{j_2=0}^p  C_{j_2}(z,x) 
C_{j_2}(s,z) dx dz\right)^2 
=
$$
$$
=
\frac{(\theta-u)^2}{2}\int\limits_{\tau}^s 
\int\limits_{\tau}^{z} \left(\sum\limits_{j_2=0}^p  C_{j_2}(z,x) 
C_{j_2}(s,z)\right)^2 dx dz \le
$$

\vspace{-1mm}
\begin{equation}
\label{marchto1}
\le K^2 \frac{(\theta-u)^2}{2}\frac{(s-\tau)^2}{2}\le K^2\frac{(T-t)^4}{4}=K_1.
\end{equation}

\vspace{3mm}
\noindent
The equality (\ref{march0005}) is proved.

The main difficulty is related to the proof
of the inequality (\ref{march0001}). Further, we prove (\ref{march0001})
for all 15 possible cases under the assumption
that 
$\{\phi_j(x)\}_{j=0}^{\infty}$ is a complete orthonormal
system of Legendre polynomials or trigonometric functions
in the space $L_2([t, T])$. As we noted above,
in some situations we will need a function $F\in L_1([t, T])$
instead of a constant $K^2$ for the polynomial case.

It is easy to see that (\ref{march0001}) reduces
to the following 15 inequalities
\begin{equation}
\label{marsixsix8}
\left|\sum_{j_1, j_2, j_3=0}^{p}
C_{j_3 j_2 j_1 j_3 j_2 j_1}(s,\tau)\right|\le K<\infty,
\end{equation}
\begin{equation}
\label{marsixsix9}
\left|\sum_{j_1, j_2, j_3=0}^{p}
C_{j_1 j_3 j_2 j_3 j_2 j_1}(s,\tau)\right|\le K<\infty,
\end{equation}
\begin{equation}
\label{marsixsix10}
\left|\sum_{j_1, j_2, j_3=0}^{p}
C_{j_3 j_2 j_3 j_1 j_2 j_1}(s,\tau)\right|\le K<\infty,
\end{equation}
\begin{equation}
\label{marsixsix4}
\left|\sum_{j_1, j_2, j_3=0}^{p}
C_{j_1 j_2 j_3 j_3 j_2 j_1}(s,\tau)\right|\le K<\infty,
\end{equation}
\begin{equation}
\label{marsixsix14}
\left|\sum_{j_1, j_2, j_3=0}^{p}
C_{j_1 j_2 j_2 j_3 j_3 j_1}(s,\tau)\right|\le K<\infty,
\end{equation}
\begin{equation}
\label{marsixsix3}
\left|\sum_{j_1, j_2, j_3=0}^{p}
C_{j_3 j_3 j_2 j_2 j_1 j_1}(s,\tau)\right|\le K<\infty,
\end{equation}
\begin{equation}
\label{marsixsix7}
\left|\sum_{j_1, j_2, j_3=0}^{p}
C_{j_2 j_3 j_3 j_2 j_1 j_1}(s,\tau)\right|\le K<\infty,
\end{equation}
\begin{equation}
\label{marsixsix6}
\left|\sum_{j_1, j_2, j_3=0}^{p}
C_{j_3 j_2 j_3 j_2 j_1 j_1}(s,\tau)\right|\le K<\infty,
\end{equation}
\begin{equation}
\label{marsixsix1}
\left|\sum_{j_1, j_2, j_3=0}^{p}
C_{j_3 j_3 j_2 j_1 j_2 j_1}(s,\tau)\right|\le K<\infty,
\end{equation}
\begin{equation}
\label{marsixsix2}
\left|\sum_{j_1, j_2, j_3=0}^{p}
C_{j_3 j_3 j_1 j_2 j_2 j_1}(s,\tau)\right|\le K<\infty,
\end{equation}
\begin{equation}
\label{marsixsix5}
\left|\sum_{j_1, j_2, j_3=0}^{p}
C_{j_2 j_1 j_3 j_3 j_2 j_1}(s,\tau)\right|\le K<\infty,
\end{equation}
\begin{equation}
\label{marsixsix12}
\left|\sum_{j_1, j_2, j_3=0}^{p}
C_{j_3 j_1 j_2 j_3 j_2 j_1}(s,\tau)\right|\le K<\infty,
\end{equation}
\begin{equation}
\label{marsixsix11}
\left|\sum_{j_1, j_2, j_3=0}^{p}
C_{j_2 j_3 j_1 j_3 j_2 j_1}(s,\tau)\right|\le K<\infty,
\end{equation}
\begin{equation}
\label{marsixsix13}
\left|\sum_{j_1, j_2, j_3=0}^{p}
C_{j_3 j_1 j_3 j_2 j_2 j_1}(s,\tau)\right|\le K<\infty,
\end{equation}
\begin{equation}
\label{marsixsix15}
\left|\sum_{j_1, j_2, j_3=0}^{p}
C_{j_2 j_3 j_3 j_1 j_2 j_1}(s,\tau)\right|\le K<\infty,
\end{equation}
where $p\in{\bf N},$ $t\le \tau < s \le T,$ 
constant $K$ does not depend on 
$p, s, \tau$ (but only on $t, T$) and may differ from line to line.

Using the technique that led to (\ref{copa1}) or (\ref{copa1a}), we obtain
$$
\int\limits_t^T h_{6}(t_{6})
\int\limits_t^{t_{6}}h_{5}(t_{5})\ldots
\int\limits_t^{t_{2}}h_{1}(t_{1})
\int\limits_t^{t_{1}}h_{7}(t_{7})dt_7 dt_1\ldots dt_6=
$$
$$
=
\int\limits_t^T h_{7}(t_{7})
\int\limits_{t_7}^T h_{1}(t_{1}) 
\int\limits_{t_1}^T h_{2}(t_{2})
\ldots
\int\limits_{t_5}^T h_{6}(t_{6})
dt_6\ldots dt_1 dt_7= 
$$
$$
=
\int\limits_t^T h_{7}(t_{7})
\left(\int\limits_{t_7}^T h_{6}(t_{6}) 
\int\limits_{t_7}^{t_6} h_{5}(t_{5})
\ldots
\int\limits_{t_7}^{t_2} h_{1}(t_{1})
dt_1\ldots dt_6\right) dt_7=
$$
\begin{equation}
\label{augu1}
~~~~~~~~~=
\int\limits_t^T h_{7}(\tau)
\left(\int\limits_{\tau}^T h_{6}(t_{6}) 
\int\limits_{\tau}^{t_6} h_{5}(t_{5})
\ldots
\int\limits_{\tau}^{t_2} h_{1}(t_{1})
dt_1\ldots dt_6\right) d\tau,
\end{equation}

\noindent
where $h_1(\tau),\ldots, h_7(\tau)\in L_2([t, T]).$

Moreover, 
$$
\int\limits_t^T h_{7}(t_{7})
\int\limits_t^{t_{7}}h_{6}(t_{6})\ldots
\int\limits_t^{t_{2}}h_{1}(t_{1})
dt_1\ldots dt_6 dt_7=
$$
\begin{equation}
\label{augu2}
~~~~~~~~~=\int\limits_t^T h_{7}(s)
\left(\int\limits_t^{s}h_{6}(t_{6})\ldots
\int\limits_t^{t_{2}}h_{1}(t_{1})
dt_1\ldots dt_6\right) ds.
\end{equation}

\vspace{1mm}

Taking into account (\ref{augu1}) and (\ref{augu2}), we note that
the conditions
(\ref{marsixsix8})--(\ref{marsixsix15})
need to be proved in two cases:
1.~$\tau=t,$\ \  2.~$s=T.$ Further, we will 
not carry out such a refinement if
some estimate from 
(\ref{marsixsix8})--(\ref{marsixsix15}) is true for all
$\tau, s\in [t, T]$ ($\tau<s$).
Looking ahead, we note that consideration
of Cases 1 and 2 will be required
only for some inequalities from (\ref{marsixsix8})--(\ref{marsixsix15})
for the polynomial case.

The relation (\ref{marsixsix3}) is a particular case of (\ref{febr2025}).
Let us prove (\ref{marsixsix8})--(\ref{marsixsix14}),
(\ref{marsixsix7})--(\ref{marsixsix15}) 
applying ideas from Sect.~2.11, 2.14, 2.32.

{\bf Step~1.}\ First, we prove (\ref{marsixsix8})--(\ref{marsixsix14}),
(\ref{marsixsix5}) using special
symmetry properties of the
Fourier coefficients. 

By analogy with (\ref{sixsix40})
we obtain
$$
C_{j_6 j_5 j_4 j_3 j_2 j_1}(s,\tau)+C_{j_1 j_2 j_3 j_4 j_5 j_6}(s,\tau)=
$$
$$
=C_{j_6}(s,\tau)C_{j_5 j_4 j_3 j_2 j_1}(s,\tau)-C_{j_5 j_6}(s,\tau)C_{j_4 j_3 j_2 j_1}(s,\tau)+
$$
$$
+C_{j_4 j_5 j_6}(s,\tau)C_{j_3 j_2 j_1}(s,\tau)-C_{j_3 j_4 j_5 j_6}(s,\tau)C_{j_2 j_1}(s,\tau)+
$$
\begin{equation}
\label{sixsix40eee}
+C_{j_2 j_3 j_4 j_5 j_6}(s,\tau)C_{j_1}(s,\tau).
\end{equation}

\vspace{2mm}

Substituting $j_1=j_4, j_2=j_5, j_3=j_6$ into (\ref{sixsix40eee}) and summing over
$j_1,j_2,$ $j_3,$ we have
$$
\sum_{j_1,j_2,j_3=0}^{p}
\hspace{-1mm}C_{j_3 j_2 j_1 j_3 j_2 j_1}(s,\tau)+
\sum_{j_1,j_2,j_3=0}^{p}
\hspace{-1mm}C_{j_1 j_2 j_3 j_1 j_2 j_3}(s,\tau)=
$$
$$
=
2\sum_{j_1,j_2,j_3=0}^{p}
\hspace{-1mm}C_{j_3 j_2 j_1 j_3 j_2 j_1}(s,\tau)=
\sum_{j_1,j_2,j_3=0}^{p}\biggl(C_{j_3}(s,\tau)C_{j_2 j_1 j_3 j_2 j_1}(s,\tau)-\biggr.
$$

\vspace{-3mm}
$$
-
C_{j_2 j_3}(s,\tau)C_{j_1 j_3 j_2 j_1}(s,\tau)+
C_{j_1 j_2 j_3}(s,\tau)C_{j_3 j_2 j_1}(s,\tau)-
$$
$$
\biggl.-C_{j_3 j_1 j_2 j_3}(s,\tau)C_{j_2 j_1}(s,\tau)+
C_{j_2 j_3 j_1 j_2 j_3}(s,\tau)C_{j_1}(s,\tau)\biggr).
$$

\vspace{2mm}

Then
$$
\sum_{j_1,j_2,j_3=0}^{p}
C_{j_3 j_2 j_1 j_3 j_2 j_1}(s,\tau)=
\frac{1}{2}\sum_{j_1,j_2,j_3=0}^{p}\biggl(
C_{j_3}(s,\tau)C_{j_2 j_1 j_3 j_2 j_1}(s,\tau)-\biggr.
$$
$$
-C_{j_2 j_3}(s,\tau)C_{j_1 j_3 j_2 j_1}(s,\tau)+
C_{j_1 j_2 j_3}(s,\tau)C_{j_3 j_2 j_1}(s,\tau)-
$$
\begin{equation}
\label{march0009}
~~~~~~~~~~\biggl.-C_{j_3 j_1 j_2 j_3}(s,\tau)C_{j_2 j_1}(s,\tau)+
C_{j_2 j_3 j_1 j_2 j_3}(s,\tau)C_{j_1}(s,\tau)\biggr).
\end{equation}

\vspace{3mm}

Similarly, we get
$$
\sum_{j_1,j_2,j_3=0}^{p}
C_{j_1 j_3 j_2 j_3 j_2 j_1}(s,\tau)=
\frac{1}{2}
\sum_{j_1,j_2,j_3=0}^{p}\biggl(
C_{j_1}(s,\tau)C_{j_3 j_2 j_3 j_2 j_1}(s,\tau)-\biggr.
$$
$$
-C_{j_3 j_1}(s,\tau)C_{j_2 j_3 j_2 j_1}(s,\tau)+
C_{j_2 j_3 j_1}(s,\tau)C_{j_3 j_2 j_1}(s,\tau)-
$$
\begin{equation}
\label{march00010}
~~~~~~~~~~\biggl.-C_{j_3 j_2 j_3 j_1}(s,\tau)C_{j_2 j_1}(s,\tau)+
C_{j_2 j_3 j_2 j_3 j_1}(s,\tau)C_{j_1}(s,\tau)\biggr),
\end{equation}

\vspace{3mm}

$$
\sum_{j_1,j_2,j_3=0}^{p}
C_{j_3 j_2 j_3 j_1 j_2 j_1}(s,\tau)=
\frac{1}{2}\sum_{j_1,j_2,j_3=0}^{p}\biggl(
C_{j_3}(s,\tau)C_{j_2 j_3 j_1 j_2 j_1}(s,\tau)-\biggr.
$$
$$
-C_{j_2 j_3}(s,\tau)C_{j_3 j_1 j_2 j_1}(s,\tau)+
C_{j_3 j_2 j_3}(s,\tau)C_{j_1 j_2 j_1}(s,\tau)-
$$
\begin{equation}
\label{march00011}
~~~~~~~~~~\biggl.-C_{j_1 j_3 j_2 j_3}(s,\tau)C_{j_2 j_1}(s,\tau)+
C_{j_2 j_1 j_3 j_2 j_3}(s,\tau)C_{j_1}(s,\tau)\biggr),
\end{equation}
$$
\sum_{j_1,j_2,j_3=0}^{p}
C_{j_1 j_2 j_3 j_3 j_2 j_1}(s,\tau)=
\frac{1}{2}
\sum_{j_1,j_2,j_3=0}^{p}\biggl(
C_{j_1}(s,\tau)C_{j_2 j_3 j_3 j_2 j_1}(s,\tau)-\biggr.
$$
$$
-C_{j_2 j_1}(s,\tau)C_{j_3 j_3 j_2 j_1}(s,\tau)+
\left(C_{j_3 j_2 j_1}(s,\tau)\right)^2-
$$
\begin{equation}
\label{march00012}
~~~~~~~~~~\biggl.-C_{j_3 j_3 j_2 j_1}(s,\tau)C_{j_2 j_1}(s,\tau)+
C_{j_2 j_3 j_3 j_2 j_1}(s,\tau)C_{j_1}(s,\tau)\biggr),
\end{equation}

\vspace{3mm}

$$
\sum_{j_1,j_2,j_3=0}^{p}
C_{j_1 j_3 j_3 j_2 j_2 j_1}(s,\tau)=
\frac{1}{2}
\sum_{j_1,j_2,j_3=0}^{p}
\biggl(
C_{j_1}(s,\tau)C_{j_3 j_3 j_2 j_2 j_1}(s,\tau)-\biggr.
$$
$$
-C_{j_3 j_1}(s,\tau)C_{j_3 j_2 j_2 j_1}(s,\tau)+
C_{j_3 j_3 j_1}(s,\tau)C_{j_2 j_2 j_1}(s,\tau)-
$$
\begin{equation}
\label{march00013}
~~~~~~~~~~\biggl.-C_{j_2 j_3 j_3 j_1}(s,\tau)C_{j_2 j_1}(s,\tau)+
C_{j_2 j_2 j_3 j_3 j_1}(s,\tau)C_{j_1}(s,\tau)\biggr),
\end{equation}

\vspace{3mm}

$$
\sum_{j_1,j_2,j_3=0}^{p}
C_{j_2 j_1 j_3 j_3 j_2 j_1}(s,\tau)=
\frac{1}{2}
\sum_{j_1,j_2,j_3=0}^{p}
\biggl(
C_{j_2}(s,\tau)C_{j_1 j_3 j_3 j_2 j_1}(s,\tau)-\biggr.
$$
$$
-C_{j_1 j_2}(s,\tau)C_{j_3 j_3 j_2 j_1}(s,\tau)+
C_{j_3 j_1 j_2}(s,\tau)C_{j_3 j_2 j_1}(s,\tau)-
$$
\begin{equation}
\label{march00014}
~~~~~~~~~~~~~\biggl.-C_{j_3 j_3 j_1 j_2}(s,\tau)C_{j_2 j_1}(s,\tau)+
C_{j_2 j_3 j_3 j_1 j_2}(s,\tau)C_{j_1}(s,\tau)\biggr).
\end{equation}

\vspace{3mm}

Using to the right-hand sides of (\ref{march0009})--(\ref{march00014})
the technique that led to the estimate (\ref{marchto1}), we obtain the inequalities
(\ref{marsixsix8})--(\ref{marsixsix14}), (\ref{marsixsix5}).

{\bf Step~2.}\ It is not difficult to see that
\begin{equation}
\label{march00015}
~~~~~~~~~~~\sum_{j_1, j_2, j_3=0}^{p}
C_{j_3 j_3 j_1 j_2 j_2 j_1}(s,\tau)=
\sum_{j_1, j_2, j_3=0}^{p}
C_{j_1 j_1 j_2 j_3 j_3 j_2}(s,\tau),
\end{equation}
\begin{equation}
\label{march00016}
~~~~~~~~~~~\sum_{j_1, j_2, j_3=0}^{p}
C_{j_3 j_3 j_2 j_1 j_2 j_1}(s,\tau)=
\sum_{j_1, j_2, j_3=0}^{p}
C_{j_1 j_1 j_2 j_3 j_2 j_3}(s,\tau),
\end{equation}
\begin{equation}
\label{march00017}
~~~~~~~~~~~\sum_{j_1, j_2, j_3=0}^{p}
C_{j_2 j_3 j_3 j_1 j_2 j_1}(s,\tau)=
\sum_{j_1, j_2, j_3=0}^{p}
C_{j_1 j_2 j_2 j_3 j_1 j_3}(s,\tau).
\end{equation}

\vspace{2mm}

Further, 
using (\ref{march00015})--(\ref{march00017}) and (\ref{sixsix40eee}), we get
$$
\sum_{j_1, j_2, j_3=0}^{p}
C_{j_2 j_3 j_3 j_2 j_1 j_1}(s,\tau)
+
\sum_{j_1, j_2, j_3=0}^{p}
C_{j_3 j_3 j_1 j_2 j_2 j_1}(s,\tau)
=
$$
$$
=\sum_{j_1, j_2, j_3=0}^{p}
C_{j_2 j_3 j_3 j_2 j_1 j_1}(s,\tau)
+
\sum_{j_1, j_2, j_3=0}^{p}
C_{j_1 j_1 j_2 j_3 j_3 j_2}(s,\tau)
=
$$
$$
=\sum_{j_1,j_2,j_3=0}^{p}
\biggl(
C_{j_2}(s,\tau)C_{j_3 j_3 j_2 j_1 j_1}(s,\tau)-\biggr.
$$
$$
-C_{j_3 j_2}(s,\tau)C_{j_3 j_2 j_1 j_1}(s,\tau)+
C_{j_3 j_3 j_2}(s,\tau)C_{j_2 j_1 j_1}(s,\tau)-
$$
\begin{equation}
\label{march00018}
~~~~~~~~~~\biggl.-C_{j_2 j_3 j_3 j_2}(s,\tau)C_{j_1 j_1}(s,\tau)+
C_{j_1 j_2 j_3 j_3 j_2}(s,\tau)C_{j_1}(s,\tau)\biggr),
\end{equation}

\vspace{2.5mm}

$$
\sum_{j_1, j_2, j_3=0}^{p}
C_{j_3 j_2 j_3 j_2 j_1 j_1}(s,\tau)+
\sum_{j_1, j_2, j_3=0}^{p}
C_{j_3 j_3 j_2 j_1 j_2 j_1}(s,\tau)=
$$
$$
=
\sum_{j_1, j_2, j_3=0}^{p}
C_{j_3 j_2 j_3 j_2 j_1 j_1}(s,\tau)+
\sum_{j_1, j_2, j_3=0}^{p}
C_{j_1 j_1 j_2 j_3 j_2 j_3}(s,\tau)=
$$
$$
=\sum_{j_1,j_2,j_3=0}^{p}
\biggl(
C_{j_3}(s,\tau)C_{j_2 j_3 j_2 j_1 j_1}(s,\tau)-\biggr.
$$
$$
-C_{j_2 j_3}(s,\tau)C_{j_3 j_2 j_1 j_1}(s,\tau)+
C_{j_3 j_2 j_3}(s,\tau)C_{j_2 j_1 j_1}(s,\tau)-
$$
\begin{equation}
\label{march00019}
~~~~~~~~~~\biggl.-C_{j_2 j_3 j_2 j_3}(s,\tau)C_{j_1 j_1}(s,\tau)+
C_{j_1 j_2 j_3 j_2 j_3}(s,\tau)C_{j_1}(s,\tau)\biggr),
\end{equation}

\vspace{2.5mm}

$$
\sum_{j_1, j_2, j_3=0}^{p}
C_{j_3 j_1 j_3 j_2 j_2 j_1}(s,\tau)+
\sum_{j_1, j_2, j_3=0}^{p}
C_{j_2 j_3 j_3 j_1 j_2 j_1}(s,\tau)=
$$
$$
=\sum_{j_1, j_2, j_3=0}^{p}
C_{j_3 j_1 j_3 j_2 j_2 j_1}(s,\tau)+
\sum_{j_1, j_2, j_3=0}^{p}
C_{j_1 j_2 j_2 j_3 j_1 j_3}(s,\tau)=
$$
$$
=\sum_{j_1,j_2,j_3=0}^{p}
\biggl(
C_{j_3}(s,\tau)C_{j_1 j_3 j_2 j_2 j_1}(s,\tau)-\biggr.
$$
$$
-C_{j_1 j_3}(s,\tau)C_{j_3 j_2 j_2 j_1}(s,\tau)+
C_{j_3 j_1 j_3}(s,\tau)C_{j_2 j_2 j_1}(s,\tau)-
$$
\begin{equation}
\label{march00020}
~~~~~~~~~~\biggl.-C_{j_2 j_3 j_1 j_3}(s,\tau)C_{j_2 j_1}(s,\tau)+
C_{j_2 j_2 j_3 j_1 j_3}(s,\tau)C_{j_1}(s,\tau)\biggr).
\end{equation}

\vspace{2mm}

Applying to the right-hand sides of (\ref{march00018})--(\ref{march00020})
the technique that led to the estimate (\ref{marchto1}), we obtain the inequalities
\begin{equation}
\label{march00021}
~~~~~~~~~~\left|\sum_{j_1, j_2, j_3=0}^{p}
C_{j_2 j_3 j_3 j_2 j_1 j_1}(s,\tau)
+
\sum_{j_1, j_2, j_3=0}^{p}
C_{j_3 j_3 j_1 j_2 j_2 j_1}(s,\tau)\right|\le K <\infty,
\end{equation}
\begin{equation}
\label{march00022}
~~~~~~~~~~\left|\sum_{j_1, j_2, j_3=0}^{p}
C_{j_3 j_2 j_3 j_2 j_1 j_1}(s,\tau)+
\sum_{j_1, j_2, j_3=0}^{p}
C_{j_3 j_3 j_2 j_1 j_2 j_1}(s,\tau)\right|\le K <\infty,
\end{equation}
\begin{equation}
\label{march00023}
~~~~~~~~~~\left|\sum_{j_1, j_2, j_3=0}^{p}
C_{j_3 j_1 j_3 j_2 j_2 j_1}(s,\tau)+
\sum_{j_1, j_2, j_3=0}^{p}
C_{j_2 j_3 j_3 j_1 j_2 j_1}(s,\tau)
\right|\le K <\infty,
\end{equation}

\vspace{3mm}
\noindent
where $p\in{\bf N},$ $t\le \tau < s \le T,$ 
constant $K$ does not depend on 
$p, s, \tau$ (but only on $t, T$) and may differ from line to line.

Note that $\left|a\right|\le K_1+K$ follows from
$\left|b\right|\le K$ and $\left|a+b\right|\le K_1,$
where $a, b, K, K_1\in {\bf R}.$
Indeed, we have $\left|a\right|=\left|a+b-b\right|\le
\left|a+b\right|+\left|b\right|\le K_1+K$. Then
from (\ref{march00021})--(\ref{march00023})
it follows that if we prove (\ref{marsixsix1}), (\ref{marsixsix2}), (\ref{marsixsix15}), then 
(\ref{marsixsix6}), (\ref{marsixsix7}), (\ref{marsixsix13}) will be proved.
Thus, it remains to prove 
(\ref{marsixsix1}), (\ref{marsixsix2}), (\ref{marsixsix12}), (\ref{marsixsix11}),
(\ref{marsixsix15}).

{\bf Step~3.}\ Let us prove 
(\ref{marsixsix1}), (\ref{marsixsix2}), (\ref{marsixsix12}), (\ref{marsixsix11}),
(\ref{marsixsix15}).
Consider (\ref{marsixsix11}). 
Using the 
Cau\-chy--Bunyakovsky inequality as well as 
Fubini's Theorem, Parseval's equality, (\ref{after80xx}), (\ref{dsds14fffff})
and Lebesgue's Dominated Convergence Theorem, we have
$$
\left(\sum_{j_1, j_2, j_3=0}^{p}
C_{j_2 j_3 j_1 j_3 j_2 j_1}(s,\tau)\right)^2=
$$
$$
=
\left(\sum_{j_2=0}^{p} 1\cdot \sum_{j_1, j_3=0}^{p}
C_{j_2 j_3 j_1 j_3 j_2 j_1}(s,\tau)\right)^2\le
$$
$$
\le \sum_{j_2=0}^{p} 1^2 \cdot \sum_{j_2=0}^{p}\left(\sum_{j_1, j_3=0}^{p}
C_{j_2 j_3 j_1 j_3 j_2 j_1}(s,\tau)\right)^2=
$$
$$
=
(p+1)\sum_{j_2=0}^{p}\left(\sum_{j_1, j_3=0}^{p}
C_{j_2 j_3 j_1 j_3 j_2 j_1}(s,\tau)\right)^2=
$$
$$
=
(p+1)\sum_{j_2=0}^{p}\left(\sum_{j_1, j_3=0}^{p}
\int\limits_{\tau}^s \phi_{j_2}(t_6)
\int\limits_{\tau}^{t_6} \phi_{j_2}(t_2)
C_{j_1}(t_2,\tau)C_{j_3 j_1 j_3}(t_6,t_2) dt_2 dt_6
\right)^2\le
$$
$$
\le
(p+1)\sum_{j_2, j_2'=0}^{p}\left(\sum_{j_1, j_3=0}^{p}
\int\limits_{\tau}^s \phi_{j_2}(t_6)
\int\limits_{\tau}^{t_6} \phi_{j_2'}(t_2)
C_{j_1}(t_2,\tau)C_{j_3 j_1 j_3}(t_6,t_2) dt_2 dt_6
\right)^2\le
$$
$$
\le
(p+1)\sum_{j_2, j_2'=0}^{\infty}\left(
\int\limits_{\tau}^s \phi_{j_2}(t_6)
\int\limits_{\tau}^{t_6} \phi_{j_2'}(t_2)
\sum_{j_1, j_3=0}^{p}C_{j_1}(t_2,\tau)C_{j_3 j_1 j_3}(t_6,t_2) dt_2 dt_6
\right)^2=
$$
$$
=
(p+1)
\int\limits_{\tau}^s 
\int\limits_{\tau}^{t_6} 
\left(\sum_{j_1=0}^{p}C_{j_1}(t_2,\tau)\sum_{j_3=0}^{p}C_{j_3 j_1 j_3}(t_6,t_2)\right)^2 dt_2 dt_6
=
$$
$$
=
(p+1)
\int\limits_{\tau}^s 
\int\limits_{\tau}^{t_6} 
\left(\sum_{j_1=0}^{p}C_{j_1}(t_2,\tau)\sum_{j_3=p+1}^{\infty}C_{j_3 j_1 j_3}(t_6,t_2)\right)^2 dt_2 dt_6
\le
$$
$$
\le
(p+1)
\int\limits_{\tau}^s 
\int\limits_{\tau}^{t_6} 
\sum_{j_1=0}^{p}C_{j_1}^2(t_2,\tau)
\sum_{j_1=0}^{p}\left(\sum_{j_3=p+1}^{\infty}C_{j_3 j_1 j_3}(t_6,t_2)\right)^2 dt_2 dt_6\le
$$
$$
\le
(p+1)
\int\limits_{\tau}^s 
\int\limits_{\tau}^{t_6} 
\sum_{j_1=0}^{\infty}C_{j_1}^2(t_2,\tau)
\sum_{j_1=0}^{p}\left(\sum_{j_3=p+1}^{\infty}C_{j_3 j_1 j_3}(t_6,t_2)\right)^2 dt_2 dt_6=
$$
$$
\le
(p+1)
\int\limits_{\tau}^s 
\int\limits_{\tau}^{t_6} 
(t_2-\tau)
\sum_{j_1=0}^{p}\left(\sum_{j_3=p+1}^{\infty}C_{j_3 j_1 j_3}(t_6,t_2)\right)^2 dt_2 dt_6=
$$
$$
=
(p+1)
\int\limits_{\tau}^s 
\int\limits_{\tau}^{t_6} 
(t_2-\tau)
\sum_{j_1=0}^{p}\left(\sum_{j_3=p+1}^{\infty}
\int\limits_{t_2}^{t_6}\phi_{j_1}(\theta)
C_{j_3}(\theta,t_2)C_{j_3}(t_6,\theta)d\theta\right)^2 dt_2 dt_6=
$$
$$
=
(p+1)
\int\limits_{\tau}^s 
\int\limits_{\tau}^{t_6} 
(t_2-\tau)
\sum_{j_1=0}^{p}\left(
\int\limits_{t_2}^{t_6}\phi_{j_1}(\theta)\sum_{j_3=p+1}^{\infty}
C_{j_3}(\theta,t_2)C_{j_3}(t_6,\theta)d\theta\right)^2 dt_2 dt_6\le
$$
$$
\le
(p+1)
\int\limits_{\tau}^s 
\int\limits_{\tau}^{t_6} 
(t_2-\tau)
\sum_{j_1=0}^{\infty}\left(
\int\limits_{t_2}^{t_6}\phi_{j_1}(\theta)\sum_{j_3=p+1}^{\infty}
C_{j_3}(\theta,t_2)C_{j_3}(t_6,\theta)d\theta\right)^2 dt_2 dt_6=
$$
\begin{equation}
\label{march00025}
~~~~~~~~=
(p+1)
\int\limits_{\tau}^s 
\int\limits_{\tau}^{t_6} 
(t_2-\tau)
\int\limits_{t_2}^{t_6}\left(\sum_{j_3=p+1}^{\infty}
C_{j_3}(\theta,t_2)C_{j_3}(t_6,\theta)\right)^2d\theta dt_2 dt_6.
\end{equation}

For the trigonometric case (see (\ref{trig11})), we have the following obvious estimate
\begin{equation}
\label{march00026}
\left|C_j(x,v)\right|=\left|
\int\limits_v^x\phi_{j}(\tau)d\tau
\right|< \frac{C}{j}\ \ \ (j>0),
\end{equation}

\noindent
where constant $C$ does not depend on $j, x, v.$

Recall that (see (\ref{obana}))
\begin{equation}
\label{march00027}
\sum\limits_{j=p+1}^{\infty}\frac{1}{j^2}
\le \int\limits_{p}^{\infty}\frac{dx}{x^2}=\frac{1}{p}.
\end{equation}

Combining (\ref{march00025})--(\ref{march00027}), we get
$$
\left(\sum_{j_1, j_2, j_3=0}^{p}
C_{j_2 j_3 j_1 j_3 j_2 j_1}(s,\tau)\right)^2\le \frac{K_1(p+1)}{p^2}\le K^2,
$$

\noindent
where constants $K, K_1$ depend only on  $t, T.$
The inequality (\ref{marsixsix11}) is proved for the trigonometric case.

For the polynomial case (see (\ref{4009})), by analogy with (\ref{101oh}) and (\ref{after1940})
we have 
$$
\left|C_j(x,v)\right|=\left|
\int\limits_v^x\phi_{j}(\tau)d\tau
\right| <
$$
\begin{equation}
\label{march00028}
~~~~~~~~<
\frac{C}{j^{1-\varepsilon/2}}\Biggl(\frac{1}{(1-z^2(x))^{1/4-\varepsilon/4}}+
\frac{1}{(1-z^2(v))^{1/4-\varepsilon/4}}\Biggr),
\end{equation}

\vspace{2mm}
\noindent
where 
$j\in {\bf N},$ $z(x), z(v)\in (-1, 1)$ 
($z(x)$ is defined by (\ref{zz1})),
$x, v\in (t, T),$ $\varepsilon\in (0,1)$ is an arbitrary
small positive real number,
constant $C$ does not depend on $j$.

Recall that
(see (\ref{after1944}))
\begin{equation}
\label{march00029}
\sum\limits_{j=p+1}^{\infty}\frac{1}{j^{2-\varepsilon}}
\le \int\limits_{p}^{\infty}\frac{dx}{x^{2-\varepsilon}}=
\frac{1}{(1-\varepsilon)p^{1-\varepsilon}}.
\end{equation}

\vspace{1mm}

Combining (\ref{march00025}), (\ref{march00028}), (\ref{march00029})
($\varepsilon=1/4$), we obtain
$$
\left(\sum_{j_1, j_2, j_3=0}^{p}
C_{j_2 j_3 j_1 j_3 j_2 j_1}(s,\tau)\right)^2\le \frac{K_1(p+1)}{p^{3/2}}\le K^2,
$$

\vspace{1mm}
\noindent
where constants $K, K_1$ depend only on  $t, T.$
The inequality (\ref{marsixsix11}) is proved for the polynomial case.

Let us prove (\ref{marsixsix12}). In complete analogy
with the proof of (\ref{marsixsix11}) we have
$$
\left(\sum_{j_1, j_2, j_3=0}^{p}
C_{j_3 j_1 j_2 j_3 j_2 j_1}(s,\tau)\right)^2\le
$$
$$
\le
(p+1)
\int\limits_{\tau}^s (s-t_5)
\int\limits_{\tau}^{t_5} 
\int\limits_{t_1}^{t_5}\left(\sum_{j_2=p+1}^{\infty}
C_{j_2}(\theta,t_1)C_{j_2}(t_5,\theta)\right)^2d\theta dt_1 dt_5.
$$

\vspace{1mm}

The further proof is the same as in the case of 
(\ref{marsixsix11}).
The inequality (\ref{marsixsix12}) is proved.

Let us prove (\ref{marsixsix15}). By analogy
with the proof of (\ref{marsixsix11}) (see (\ref{march00025}})) we get
$$
\left(\sum_{j_1, j_2, j_3=0}^{p}
C_{j_2 j_3 j_3 j_1 j_2 j_1}(s,\tau)\right)^2\le
$$
\begin{equation}
\label{march00030}
~~~~~~~~\le
(p+1)
\int\limits_{\tau}^s (s-t_5)
\int\limits_{\tau}^{t_5} 
\int\limits_{\tau}^{t_4}\left(\sum_{j_1=p+1}^{\infty}
C_{j_1}(\theta,\tau)C_{j_1}(t_4,\theta)\right)^2 d\theta dt_4 dt_5.
\end{equation}

\vspace{1mm}

The further proof for the trigonometric case is the same as for
the inequality (\ref{marsixsix11}).

Consider the polynomial case. In this case,
we note that it is actually
necessary to consider the following two cases of (\ref{march00030}) (see above)
\begin{equation}
\label{march00040}
1.\ \tau=t,\ \ \ 2.\ s=T\ \  (\tau=t\ \hbox{or}\ \tau>t).
\end{equation}

For Case~1, the estimate (\ref{march00028}) is simplified
as follows (see (\ref{otit6000x}), (\ref{after5000}) and (\ref{after1940}))
\begin{equation}
\label{march00031}
~~~~~~~~~\left|C_j(x,t)\right|=\left|
\int\limits_t^x\phi_{j}(\tau)d\tau
\right| <
\frac{C}{j^{1-\varepsilon/2}}\frac{1}{(1-z^2(x))^{1/4-\varepsilon/4}},
\end{equation}

\vspace{1mm}
\noindent
where notations are the same as in (\ref{march00028}).

Combining (\ref{march00030}), (\ref{march00028}), (\ref{march00029}), (\ref{march00031})
($\varepsilon=1/4$), we obtain
\begin{equation}
\label{march00042aa}
~~~~~~~~\left(\sum_{j_1, j_2, j_3=0}^{p}
C_{j_2 j_3 j_3 j_1 j_2 j_1}(s,t)\right)^2\le
\frac{K_1(p+1)}{p^{3/2}}\le K^2,
\end{equation}

\vspace{1mm}
\noindent
where constants $K, K_1$ depend only on  $t, T.$
The inequality (\ref{marsixsix15}) is proved for the polynomial case
(Case~1).

Consider Case~2.
Combining (\ref{march00030}), (\ref{march00028}), (\ref{march00029}), (\ref{march00031})
($\varepsilon=1/4$), we ob\-tain
$$
\left(\sum_{j_1, j_2, j_3=0}^{p}
C_{j_2 j_3 j_3 j_1 j_2 j_1}(T,\tau)\right)^2
\biggl|_{\tau>t}
\le
\frac{K_1(p+1)}{p^{3/2}}\frac{1}{(1-z^2(\tau))^{3/8}}\le 
$$
\begin{equation}
\label{augu3}
\le 
\frac{K^2}{(1-z^2(\tau))^{3/8}} \stackrel{\sf def}{=} F(\tau),
\end{equation}

\begin{equation}
\label{augu4}
\left(\sum_{j_1, j_2, j_3=0}^{p}
C_{j_2 j_3 j_3 j_1 j_2 j_1}(T,\tau)\right)^2
\biggl|_{\tau=t}
\le
\frac{K_1(p+1)}{p^{3/2}},
\end{equation}

\vspace{2mm}
\noindent
where constants $K, K_1$ depend only on  $t, T$
and $F(\tau)\in L_1([t, T])$ (integrable majorant (see above in this section)).

Combining (\ref{augu3}) and (\ref{augu4}), we get
the following weakened version of the inequality (\ref{marsixsix15})
\begin{equation}
\label{march00042aaa}
\left(\sum_{j_1, j_2, j_3=0}^{p}
C_{j_2 j_3 j_3 j_1 j_2 j_1}(T,\tau)\right)^2\le F(\tau)
\end{equation}

\vspace{1mm}
\noindent
for the polynomial case
(Case~2),
where
$$
F(\tau)=\frac{K^2}{(1-z^2(\tau))^{3/8}}.
$$

\vspace{2mm}

Let us prove (\ref{marsixsix2}).
Using the 
Cau\-chy--Bunyakovsky inequality as well as 
Fubini's Theorem and Parseval's equality, we have
$$
\left(\sum_{j_1, j_2, j_3=0}^{p}
C_{j_3 j_3 j_1 j_2 j_2 j_1}(s,\tau)\right)^2=
$$
$$
=
\left(\sum_{j_3=0}^{p} 1\cdot \sum_{j_1, j_2=0}^{p}
C_{j_3 j_3 j_1 j_2 j_2 j_1}(s,\tau)\right)^2\le
$$
$$
\le \sum_{j_3=0}^{p} 1^2 \cdot \sum_{j_3=0}^{p}\left(\sum_{j_1, j_2=0}^{p}
C_{j_3 j_3 j_1 j_2 j_2 j_1}(s,\tau)\right)^2=
$$
$$
=
(p+1)\sum_{j_3=0}^{p}\left(\sum_{j_1, j_2=0}^{p}
C_{j_3 j_3 j_1 j_2 j_2 j_1}(s,\tau)\right)^2=
$$
$$
=         
(p+1)\sum_{j_3=0}^{p}\left(\sum_{j_1, j_2=0}^{p}
\int\limits_{\tau}^s \phi_{j_3}(t_6)
\int\limits_{\tau}^{t_6} \phi_{j_3}(t_5)
C_{j_1 j_2 j_2 j_1}(t_5,\tau) dt_5 dt_6
\right)^2\le
$$
$$
\le
(p+1)\sum_{j_3, j_3'=0}^{p}\left(
\int\limits_{\tau}^s \phi_{j_3}(t_6)
\int\limits_{\tau}^{t_6} \phi_{j_3'}(t_5)
\sum_{j_1, j_2=0}^{p}C_{j_1 j_2 j_2 j_1}(t_5,\tau) dt_5 dt_6
\right)^2\le
$$
$$
\le
(p+1)\sum_{j_3, j_3'=0}^{\infty}\left(
\int\limits_{\tau}^s \phi_{j_3}(t_6)
\int\limits_{\tau}^{t_6} \phi_{j_3'}(t_5)
\sum_{j_1, j_2=0}^{p}C_{j_1 j_2 j_2 j_1}(t_5,\tau) dt_5 dt_6
\right)^2=
$$
$$
=
(p+1)
\int\limits_{\tau}^s 
\int\limits_{\tau}^{t_6}\left(
\sum_{j_1, j_2=0}^{p}C_{j_1 j_2 j_2 j_1}(t_5,\tau) \right)^2dt_5 dt_6
=
$$
$$
=
(p+1)
\int\limits_{\tau}^s 
\int\limits_{\tau}^{t_6}\left(
\sum_{j_2=0}^{p} 1\cdot \sum_{j_1=0}^{p}C_{j_1 j_2 j_2 j_1}(t_5,\tau) \right)^2dt_5 dt_6
\le
$$
$$
\le
(p+1)^2
\int\limits_{\tau}^s 
\int\limits_{\tau}^{t_6}
\sum_{j_2=0}^{p}\left(
\sum_{j_1=0}^{p}C_{j_1 j_2 j_2 j_1}(t_5,\tau)\right)^2 dt_5 dt_6
=
$$
$$
=
(p+1)^2
\int\limits_{\tau}^s 
\int\limits_{\tau}^{t_6}
\sum_{j_2=0}^{p}\left(
\sum_{j_1=0}^{p}
\int\limits_{\tau}^{t_5}\phi_{j_2}(t_3)
\int\limits_{\tau}^{t_3}\phi_{j_2}(t_2)
C_{j_1}(t_2,\tau)C_{j_1}(t_5,t_3)dt_2 dt_3\right)^2  \hspace{-1.5mm} dt_5 dt_6
\le
$$
$$
\le
(p+1)^2
\int\limits_{\tau}^s 
\int\limits_{\tau}^{t_6}
\sum_{j_2,j_2'=0}^{p}\left(
\int\limits_{\tau}^{t_5}\phi_{j_2}(t_3)\times \right.
$$
$$
\left. \times
\int\limits_{\tau}^{t_3}\phi_{j_2'}(t_2)
\sum_{j_1=0}^{p}C_{j_1}(t_2,\tau)C_{j_1}(t_5,t_3)dt_2 dt_3\right)^2 dt_5 dt_6
\le
$$
$$
\le
(p+1)^2
\int\limits_{\tau}^s 
\int\limits_{\tau}^{t_6}
\sum_{j_2,j_2'=0}^{\infty}\left(
\int\limits_{\tau}^{t_5}\phi_{j_2}(t_3)\times \right.
$$
$$
\left. \times
\int\limits_{\tau}^{t_3}\phi_{j_2'}(t_2)
\left(\sum_{j_1=0}^{\infty} - \sum_{j_1=p+1}^{\infty}\right)
C_{j_1}(t_2,\tau)C_{j_1}(t_5,t_3)dt_2 dt_3\right)^2 dt_5 dt_6
=
$$
\begin{equation}
\label{march00032}
~~~~~~~~=
(p+1)^2
\int\limits_{\tau}^s 
\int\limits_{\tau}^{t_6}
\int\limits_{\tau}^{t_5}
\int\limits_{\tau}^{t_3}
\left(\sum_{j_1=p+1}^{\infty}
C_{j_1}(t_2,\tau)C_{j_1}(t_5,t_3)\right)^2 dt_2 dt_3 dt_5 dt_6.
\end{equation}

\vspace{3mm}

Consider the trigonometric case.
Combining (\ref{march00032}), (\ref{march00026}), (\ref{march00027}}), we obtain
$$
\left(\sum_{j_1, j_2, j_3=0}^{p}
C_{j_3 j_3 j_1 j_2 j_2 j_1}(s,\tau)\right)^2
\le \frac{K_1(p+1)^2}{p^2}\le K^2,
$$

\vspace{1mm}
\noindent
where constants $K, K_1$ depend only on  $t, T.$
The inequality (\ref{marsixsix2}) is proved for the trigonometric case.

Consider the polynomial case for two cases (\ref{march00040}). 
Let $\tau=t.$ The modification of the estimate
(\ref{march00028}) for $\varepsilon =0$ is as follows
(see also (\ref{6000}), (\ref{101oh}))
$$
\left|C_j(x,v)\right|=\left|
\int\limits_v^x\phi_{j}(\tau)d\tau
\right| <
$$
\begin{equation}
\label{march00042}
<
\frac{C}{j}\Biggl(\frac{1}{(1-z^2(x))^{1/4}}+
\frac{1}{(1-z^2(v))^{1/4}}\Biggr),
\end{equation}

\vspace{3mm}
\noindent
where 
$j\in {\bf N},$ $z(x), z(v)\in (-1, 1)$ 
($z(x)$ is defined by (\ref{zz1})),
$x, v\in (t, T),$
constant $C$ does not depend on $j$.

For $v=t$, the estimate (\ref{march00042}) is simplified
as follows (see (\ref{otit6000x}), (\ref{101xx}))
\begin{equation}
\label{march00043}
\left|C_j(x,t)\right|=\left|
\int\limits_t^x\phi_{j}(\tau)d\tau
\right| <
\frac{C}{j(1-z^2(x))^{1/4}},
\end{equation}

\vspace{1mm}
\noindent
where notations are the same as in (\ref{march00042}).

Combining (\ref{march00032}), (\ref{march00042}), (\ref{march00043}), we get
$$
\left(\sum_{j_1, j_2, j_3=0}^{p}
C_{j_3 j_3 j_1 j_2 j_2 j_1}(s,t)\right)^2
\le \frac{K_1(p+1)^2}{p^2}\le K^2,
$$

\noindent
where constants $K, K_1$ depend only on  $t, T.$
The inequality (\ref{marsixsix2}) is proved for the polynomial case 
($\tau=t$).

Now let $s=T.$
Using (\ref{march00032}), (\ref{march00042}), (\ref{march00043}), we ob\-tain
$$
\left(\sum_{j_1, j_2, j_3=0}^{p}
C_{j_3 j_3 j_1 j_2 j_2 j_1}(T,\tau)\right)^2\biggl|_{\tau>t}\le
\frac{K_1(p+1)^2}{p^2}\frac{1}{(1-z^2(\tau))^{1/2}}\le 
$$
\begin{equation}
\label{augu8}
\le 
\frac{K^2}{(1-z^2(\tau))^{1/2}} \stackrel{\sf def}{=} F(\tau),
\end{equation}

\begin{equation}
\label{augu9}
\left(\sum_{j_1, j_2, j_3=0}^{p}
C_{j_3 j_3 j_1 j_2 j_2 j_1}(T,\tau)\right)^2\biggl|_{\tau=t}\le
\frac{K_1(p+1)^2}{p^2}, 
\end{equation}

\vspace{2mm}
\noindent
where constants $K, K_1$ depend only on  $t, T$
and $F(\tau)\in L_1([t, T])$ (integrable majorant (see above in this section)).

Combining (\ref{augu8}) and (\ref{augu9}), we get
the following weakened version of the inequality (\ref{marsixsix2})
\begin{equation}
\label{march00042aaaa}
\left(\sum_{j_1, j_2, j_3=0}^{p}
C_{j_3 j_3 j_1 j_2 j_2 j_1}(T,\tau)\right)^2\le F(\tau)
\end{equation}

\noindent
for the polynomial case
($s=T$),
where
$$
F(\tau)=\frac{K^2}{(1-z^2(\tau))^{1/2}}.
$$

\vspace{1mm}

Finally, we prove the inequality
(\ref{marsixsix1}). By analogy with (\ref{march00032}) we get
$$
\left(\sum_{j_1, j_2, j_3=0}^{p}
C_{j_3 j_3 j_2 j_1 j_2 j_1}(s,\tau)\right)^2\le
$$
$$
\le
(p+1)\sum_{j_3=0}^{p}\left(\sum_{j_1, j_2=0}^{p}
C_{j_3 j_3 j_2 j_1 j_2 j_1}(s,\tau)\right)^2=
$$
$$
=         
(p+1)\sum_{j_3=0}^{p}\left(\sum_{j_1, j_2=0}^{p}
\int\limits_{\tau}^s \phi_{j_3}(t_6)
\int\limits_{\tau}^{t_6} \phi_{j_3}(t_5)
C_{j_2 j_1 j_2 j_1}(t_5,\tau) dt_5 dt_6
\right)^2\le
$$
$$
\le
(p+1)\sum_{j_3, j_3'=0}^{\infty}\left(
\int\limits_{\tau}^s \phi_{j_3}(t_6)
\int\limits_{\tau}^{t_6} \phi_{j_3'}(t_5)
\sum_{j_1, j_2=0}^{p}C_{j_2 j_1 j_2 j_1}(t_5,\tau) dt_5 dt_6
\right)^2=
$$
$$
=
(p+1)
\int\limits_{\tau}^s 
\int\limits_{\tau}^{t_6}\left(
\sum_{j_1, j_2=0}^{p}C_{j_2 j_1 j_2 j_1}(t_5,\tau) \right)^2dt_5 dt_6
=
$$
$$
\le
(p+1)^2
\int\limits_{\tau}^s 
\int\limits_{\tau}^{t_6}
\sum_{j_2=0}^{p}\left(
\sum_{j_1=0}^{p}C_{j_2 j_1 j_2 j_1}(t_5,\tau)\right)^2 dt_5 dt_6
=
$$
$$
=
(p+1)^2
\int\limits_{\tau}^s 
\int\limits_{\tau}^{t_6}
\sum_{j_2=0}^{p}\left(
\sum_{j_1=0}^{p}
\int\limits_{\tau}^{t_5}\phi_{j_2}(t_4)
\int\limits_{\tau}^{t_4}\phi_{j_2}(t_2)
C_{j_1}(t_2,\tau)C_{j_1}(t_4,t_2)dt_2 dt_4\right)^2  \hspace{-1.5mm} dt_5 dt_6
\le
$$
$$
\le
(p+1)^2
\int\limits_{\tau}^s 
\int\limits_{\tau}^{t_6}
\sum_{j_2,j_2'=0}^{\infty}\left(
\int\limits_{\tau}^{t_5}\phi_{j_2}(t_4)
\times \right.
$$
$$
\left. \times
\int\limits_{\tau}^{t_4}
\phi_{j_2'}(t_2)
\left(\sum_{j_1=0}^{\infty} - \sum_{j_1=p+1}^{\infty}\right)
C_{j_1}(t_2,\tau)C_{j_1}(t_4,t_2)dt_2 dt_4\right)^2 dt_5 dt_6
=
$$
\begin{equation}
\label{march00044}
~~~~~~~~=
(p+1)^2
\int\limits_{\tau}^s 
\int\limits_{\tau}^{t_6}
\int\limits_{\tau}^{t_5}
\int\limits_{\tau}^{t_4}
\left(\sum_{j_1=p+1}^{\infty}
C_{j_1}(t_2,\tau)C_{j_1}(t_4,t_2)\right)^2 dt_2 dt_4 dt_5 dt_6.
\end{equation}

\vspace{2mm}

The further proof of inequality (\ref{marsixsix1}) for the trigonometric case and the 
weakened analogue of inequality (\ref{marsixsix1}) for the polynomial case is 
completely analogous to the proof of (\ref{marsixsix15}) and its weakened 
analogue (see (\ref{march00030}), (\ref{march00042aa}), (\ref{march00042aaa})).

Thus, the following theorem is proved.

\vspace{1mm}

{\bf Theorem~2.64.}\ {\it Suppose that
$\{\phi_j(x)\}_{j=0}^{\infty}$ is a complete orthonormal
system of Legendre polynomials or trigonometric functions
in the space $L_2([t, T]).$
Then$,$ for the iterated Stra\-to\-no\-vich stochastic integral
of seventh multiplicity 
$$
J^{*}[\psi^{(7)}]_{T,t}=
{\int\limits_t^{*}}^T
\ldots
{\int\limits_t^{*}}^{t_2}
d{\bf w}_{t_1}^{(i_1)}
\ldots d{\bf w}_{t_7}^{(i_7)}
$$
the following 
expansion 
$$
J^{*}[\psi^{(7)}]_{T,t}=
\hbox{\vtop{\offinterlineskip\halign{
\hfil#\hfil\cr
{\rm l.i.m.}\cr
$\stackrel{}{{}_{p\to \infty}}$\cr
}} }
\sum\limits_{j_1,\ldots, j_7=0}^{p}
C_{j_7 \ldots j_1}\zeta_{j_1}^{(i_1)}\ldots \zeta_{j_7}^{(i_7)}
$$
that converges in the mean-square sense is valid, where 
$i_1,\ldots,i_7=0, 1,\ldots,m,$
$$
C_{j_7\ldots j_1}=\int\limits_t^T
\phi_{j_7}(t_7)
\ldots
\int\limits_t^{t_2}
\phi_{j_1}(t_1)dt_1\ldots dt_7
$$
and
$$
\zeta_{j}^{(i)}=
\int\limits_t^T \phi_{j}(\tau) d{\bf w}_{\tau}^{(i)}
$$ 
are independent standard Gaussian random variables for various 
$i$ or $j$ {\rm (}in the case when $i\ne 0${\rm ),}
${\bf w}_{\tau}^{(i)}$
$(i=1,\ldots,m)$ are independent 
standard Wiener processes$,$
${\bf w}_{\tau}^{(0)}=\tau.$}

\section{Expansion of Iterated Stratonovich Stochastic Integrals
of Multiplicity 8 for the Case $\psi_1(\tau),\ldots, \psi_8(\tau)
\equiv 1$ (The Cases of Legendre 
Polynomials and Trigonometric Functions)}

This section is devoted to the following theorem.

{\bf Theorem~2.65.}\ {\it Suppose that
$\{\phi_j(x)\}_{j=0}^{\infty}$ is a complete orthonormal
system of Legendre polynomials or trigonometric functions
in the space $L_2([t, T]).$
Then$,$ for the iterated Stra\-to\-no\-vich stochastic integral
of eighth multiplicity 
$$
J^{*}[\psi^{(8)}]_{T,t}=
{\int\limits_t^{*}}^T
\ldots
{\int\limits_t^{*}}^{t_2}
d{\bf w}_{t_1}^{(i_1)}
\ldots d{\bf w}_{t_8}^{(i_8)}
$$
the following 
expansion 
$$
J^{*}[\psi^{(8)}]_{T,t}=
\hbox{\vtop{\offinterlineskip\halign{
\hfil#\hfil\cr
{\rm l.i.m.}\cr
$\stackrel{}{{}_{p\to \infty}}$\cr
}} }
\sum\limits_{j_1,\ldots, j_8=0}^{p}
C_{j_8 \ldots j_1}\zeta_{j_1}^{(i_1)}\ldots \zeta_{j_8}^{(i_8)}
$$
that converges in the mean-square sense is valid, where 
$i_1,\ldots,i_8=0, 1,\ldots,m,$
$$
C_{j_8\ldots j_1}=\int\limits_t^T
\phi_{j_8}(t_8)
\ldots
\int\limits_t^{t_2}
\phi_{j_1}(t_1)dt_1\ldots dt_8
$$
and
$$
\zeta_{j}^{(i)}=
\int\limits_t^T \phi_{j}(\tau) d{\bf w}_{\tau}^{(i)}
$$ 
are independent standard Gaussian random variables for various 
$i$ or $j$ {\rm (}when $i\ne 0${\rm ),}
${\bf w}_{\tau}^{(i)}$ 
$(i=1,\ldots,m)$ are independent 
standard Wiener processes$,$
${\bf w}_{\tau}^{(0)}=\tau.$}

{\bf Proof.}\ To prove the theorem, we need to check
the condition (\ref{09091}) (or the condition (\ref{09091xxx}))
for the case $k=8>2r$, where $r=1, 2, 3$  (see Theorem~2.61).
Recall that the case
$k=2r$ is considered in Sect.~2.27.4 (see (\ref{july90000})).
Under the conditions of Theorem~2.65, this means that
$k=8=2r$, where $r=4$.

The relations
(\ref{2024december12})--(\ref{2024december11}), (\ref{cc123})
cover the case $k=8,$ $r=1,2$ (see (\ref{09091})).

Thus, it remains to consider the case $k=8,$ $r=3.$
The case $k=7,$ $r=3$ was considered in the previous section.
Here we will focus on the differences
between these two cases.

Since now $k=8,$ then along with inequalities
(\ref{march0001})--(\ref{march0004}), it is necessary to prove the following 
inequalities
$$
\left|\sum\limits_{j_{g_1},j_{g_3}. j_{g_5}=0}^p
\bigl(C_{j_{d_3} j_{d_3-1}j_{d_3-2}j_{d_3-3}}(s,\tau)
C_{j_{d_2}}(\theta,u) C_{j_{d_1}}(\rho,v)
\bigr)\biggl|_{j_{g_1}=j_{g_2},j_{g_3}=j_{g_4},j_{g_5}=j_{g_6}}\right|\le 
$$
\begin{equation}
\label{march00045}
\le K<\infty,
\end{equation}
$$
\left|\sum\limits_{j_{g_1},j_{g_3}. j_{g_5}=0}^p
\bigl(C_{j_{d_3} j_{d_3-1}j_{d_3-2}}(s,\tau)
C_{j_{d_2}j_{d_2-1}}(\theta,u) C_{j_{d_1}}(\rho,v)
\bigr)\biggl|_{j_{g_1}=j_{g_2},j_{g_3}=j_{g_4},j_{g_5}=j_{g_6}}\right|\le 
$$
\begin{equation}
\label{march00046}
\le K<\infty,
\end{equation}
$$
\left|\sum\limits_{j_{g_1},j_{g_3}. j_{g_5}=0}^p
\bigl(C_{j_{d_3} j_{d_3-1}}(s,\tau)
C_{j_{d_2}j_{d_2-1}}(\theta,u) C_{j_{d_1} j_{d_1-1}}(\rho,v)
\bigr)\biggl|_{j_{g_1}=j_{g_2},j_{g_3}=j_{g_4},j_{g_5}=j_{g_6}}\right|\le 
$$
\begin{equation}
\label{march00047}
\le K<\infty,
\end{equation}

\noindent
where $p\in{\bf N},$ $t\le \tau < s \le T,$ $t\le u<\theta \le T,$
$t\le v<\rho \le T,$
constant $K$ does not depend on 
$p, s, \tau, \theta, u, \rho, v$ (but only on $t, T$) and may differ from line to line;
another notations are the same as in Sect.~2.31, 2.35.

The inequalities (\ref{march00045})--(\ref{march00047})
are proved using the same technique as 
inequalities (\ref{2024december12})--(\ref{2024december11}) (see Sect.~2.32).
Here we will only prove as an example the following
special case of the inequality (\ref{march00047})
\begin{equation}
\label{march00048}
~~~~~~~~~~~~~~~\left|\sum\limits_{j_1, j_2, j_3=0}^p C_{j_2 j_1}(s,\tau)C_{j_3 j_1}(\theta,u)
C_{j_2 j_3}(\rho,v)\right|
\le K<\infty.
\end{equation}

Using the 
Cau\-chy--Bunyakovsky inequality as well as 
Fubini's Theorem, Parseval's equality and (\ref{2024december13}), we have
$$
\left(
\sum\limits_{j_1, j_2, j_3=0}^p C_{j_2 j_1}(s,\tau)C_{j_3 j_1}(\theta,u)
C_{j_2 j_3}(\rho,v)\right)^2=
$$
$$
=\left(\sum\limits_{j_2, j_3=0}^p  C_{j_2 j_3}(\rho,v)
\sum\limits_{j_1=0}^p C_{j_2 j_1}(s,\tau)C_{j_3 j_1}(\theta,u)
\right)^2\le
$$
$$
\le \sum\limits_{j_2, j_3=0}^p  C_{j_2 j_3}^2(\rho,v)
\sum\limits_{j_2, j_3=0}^p \left(\sum\limits_{j_1=0}^p C_{j_2 j_1}(s,\tau)C_{j_3 j_1}(\theta,u)
\right)^2\le
$$
$$
\le \sum\limits_{j_2, j_3=0}^{\infty}  C_{j_2 j_3}^2(\rho,v)
\sum\limits_{j_2, j_3=0}^{\infty} \left(\sum\limits_{j_1=0}^p C_{j_2 j_1}(s,\tau)C_{j_3 j_1}(\theta,u)
\right)^2=
$$
$$
=\frac{(\rho-v)^2}{2}
\sum\limits_{j_2, j_3=0}^{\infty} 
\left(\sum\limits_{j_1=0}^p \int\limits_{\tau}^s \phi_{j_2}(t_2)\int\limits_{\tau}^{t_2} \phi_{j_1}(t_1)
dt_1 dt_2
\int\limits_{u}^{\theta} \phi_{j_3}(t_4)\int\limits_{u}^{t_4} \phi_{j_1}(t_3)
dt_3 dt_4
\right)^2=
$$
$$
=\frac{(\rho-v)^2}{2}
\sum\limits_{j_2, j_3=0}^{\infty} \left(
\int\limits_{\tau}^s \int\limits_{u}^{\theta} \phi_{j_2}(t_2)\phi_{j_3}(t_4)\times
\right.
$$
$$
\left.\times
\sum\limits_{j_1=0}^p 
\int\limits_{\tau}^{t_2} \phi_{j_1}(t_1)
dt_1 \int\limits_{u}^{t_4} \phi_{j_1}(t_3)
dt_1 dt_3 dt_4 dt_2
\right)^2=
$$
$$
=\frac{(\rho-v)^2}{2}
\int\limits_{\tau}^s \int\limits_{u}^{\theta} 
\left(
\sum\limits_{j_1=0}^p  C_{j_1}(t_2,\tau) C_{j_1}(t_4,u)
\right)^2 dt_4 dt_2\le
$$

\vspace{-1mm}
$$
\le
K_1^2 \frac{(\rho-v)^2}{2} (s-\tau)(\theta-u)\le
$$

\vspace{-1mm}
$$
\le 
K_1^2 \frac{(T-t)^4}{2}=K.
$$

\vspace{2mm}
\noindent
The inequality (\ref{march00048}) is proved.

The inequalities (\ref{march0001})--(\ref{march0004})
for the case $k=8$ 
are proved similarly to the 
inequalities (\ref{march0001})--(\ref{march0004})
for the case $k=7$ (see Sect.~2.36).
There will be minor differences only when proving
(\ref{march0001}) for the case $k=8$ (polynomial case).
The above differences will be due to 
the fact that along with the two cases (\ref{march00040}) the following third
case 
$$
\tau, s \in (t, T)
$$

\noindent
will now appear when proving 
(\ref{marsixsix1}), (\ref{marsixsix2}), (\ref{marsixsix15}).

Let us explain this in more detail. Cases~1 and 2 are
obtained similarly to (\ref{augu1}) and (\ref{augu2}).

Further, 
using the technique that led to (\ref{copa1}) or (\ref{copa1a}), we have
$$
\int\limits_t^T h_{7}(t_{7})
\int\limits_t^{t_7} h_{6}(t_{6})
\int\limits_t^{t_{6}}h_{5}(t_{5})\ldots
\int\limits_t^{t_{2}}h_{1}(t_{1})
\int\limits_t^{t_{1}}h_{8}(t_{8})dt_8 dt_1\ldots dt_6 dt_7=
$$

\vspace{-6mm}
$$
=\int\limits_t^T h_{7}(t_{7})
\int\limits_t^{t_7} h_{8}(t_{8})
\int\limits_{t_8}^{t_7}
h_{1}(t_{1})
\int\limits_{t_1}^{t_{7}}h_{2}(t_{2})\ldots
\int\limits_{t_5}^{t_{7}}h_{6}(t_{6})dt_6 \ldots dt_1 dt_8 dt_7=
$$

\vspace{-6mm}
$$
=\int\limits_t^T h_{7}(t_{7})
\int\limits_t^{t_7} h_{8}(t_{8})
\int\limits_{t_8}^{t_7}
h_{6}(t_{6})
\int\limits_{t_8}^{t_{6}}h_{5}(t_{5})\ldots
\int\limits_{t_8}^{t_{2}}h_{1}(t_{1})dt_1 \ldots dt_6 dt_8 dt_7=
$$
\begin{equation}
\label{augu11}
~~=\int\limits_t^T h_{7}(s)
\int\limits_t^{s} h_{8}(\tau)
\left(\int\limits_{\tau}^{s}
h_{6}(t_{6})
\int\limits_{\tau}^{t_{6}}h_{5}(t_{5})\ldots
\int\limits_{\tau}^{t_{2}}h_{1}(t_{1})dt_1 \ldots dt_6\right) d\tau ds,
\end{equation}

\noindent
where $h_1(\tau),\ldots, h_8(\tau)\in L_2([t, T]).$
The equality (\ref{augu11}) explains the occurrence
of the third case.

Applying the technique that led to the estimates (\ref{march00042aaa}),
(\ref{march00042aaaa}),
we obtain for Case~3
$$
\left(\sum_{j_1, j_2, j_3=0}^{p}
C_{j_3 j_3 j_2 j_1 j_2 j_1}(s,\tau)
\right)^2\le
\frac{K^2}{(1-z^2(\tau))^{3/8}}\stackrel{\sf def}{=}F(\tau)\ \ \ (\hbox{for}\ (\ref{marsixsix1})),
$$

\vspace{-1mm}
$$
\left(\sum_{j_1, j_2, j_3=0}^{p}
C_{j_3 j_3 j_1 j_2 j_2 j_1}(s,\tau)
\right)^2\le
\frac{K^2}{(1-z^2(\tau))^{1/2}}\stackrel{\sf def}{=}F(\tau)\ \ \ (\hbox{for}\ (\ref{marsixsix2})),
$$

\vspace{-1mm}
$$
\left(\sum_{j_1, j_2, j_3=0}^{p}
C_{j_2 j_3 j_3 j_1 j_2 j_1}(s,\tau)\right)^2\le 
\frac{K^2}{(1-z^2(\tau))^{3/8}}\stackrel{\sf def}{=}F(\tau)\ \ \ (\hbox{for}\ (\ref{marsixsix15})),
$$

\vspace{1mm}
\noindent
where constant $K$ depends only on  $t, T$
and $F(\tau)\in L_1([t, T])$ (integrable majorant (see (\ref{09091xxx})).
Theorem~2.65 is proved.

\section{Verification of the Conditions of Theorems~2.52--2.57
for the Case $\psi_1(\tau),\ldots,\psi_5(\tau)
\equiv \psi(\tau)$ and Generalization of Theorems
2.62, 2.64, 2.65 
to the Case $\psi_1(\tau),\ldots,\psi_8(\tau)
\equiv \psi(\tau)$}

It is easy to see that Theorems~2.52--2.57
will be true if
$\psi_1(\tau),\ldots,\psi_5(\tau)\equiv \psi(\tau),$ where $\psi(\tau)\in L_2([t, T])$
or $\psi(\tau)$ is a continuous function on $[t, T].$
Furthermore, Theorems~2.62, 2.64, 2.65 
can be generalized to the case 
$\psi_1(\tau),\ldots,\psi_8(\tau)\equiv \psi(\tau),$ where $\psi(\tau)$ is a continuous 
or continuously dif\-ferentiable 
function on $[t, T]$.
Let us provide the corresponding explanations.

Using Fubini's Theorem and Parseval's equality, we obtain for the case 
$\psi_1(\tau),\psi_2(\tau),\psi_3(\tau)\equiv \psi(\tau)\in L_2([t, T])$
$$
\sum\limits_{j_1=0}^{p}
\int\limits_t^{s}\psi_2(\tau)\phi_{j_1}(\tau)
\int\limits_t^{\tau}\psi_1(\theta)\phi_{j_1}(\theta)
d\theta d\tau=
\sum\limits_{j_1=0}^{p}
\int\limits_t^{s}\psi(\tau)\phi_{j_1}(\tau)
\int\limits_t^{\tau}\psi(\theta)\phi_{j_1}(\theta)
d\theta d\tau
=
$$
$$
=
\frac{1}{2}
\sum\limits_{j_1=0}^{p}
\left(\int\limits_t^{s}\psi(\tau)\phi_{j_1}(\tau)d\tau\right)^2\le
\frac{1}{2}
\sum\limits_{j_1=0}^{\infty}
\left(\int\limits_t^{s}\psi(\tau)\phi_{j_1}(\tau)d\tau\right)^2=
$$
\begin{equation}
\label{march00049}
=
\int\limits_t^s \psi^2(\tau) d\tau\le 
\int\limits_t^T \psi^2(\tau) d\tau =K < \infty,
\end{equation}
$$
\sum\limits_{j_3=0}^{p}
\int\limits_{s}^T \psi_2(\tau)\phi_{j_3}(\tau)
\int\limits_{\tau}^T \psi_3(\theta)\phi_{j_3}(\theta)d\theta d\tau
=\sum\limits_{j_3=0}^{p}
\int\limits_{s}^T \psi(\tau)\phi_{j_3}(\tau)
\int\limits_{\tau}^T \psi(\theta)\phi_{j_3}(\theta)d\theta d\tau=
$$
$$
=
\frac{1}{2}\sum\limits_{j_3=0}^{p}
\left(\int\limits_{s}^T \psi(\tau)\phi_{j_3}(\tau)d\tau\right)^2\le
\frac{1}{2}\sum\limits_{j_3=0}^{\infty}
\left(\int\limits_{s}^T \psi(\tau)\phi_{j_3}(\tau)d\tau\right)^2=
$$
$$
=
\int\limits_s^T \psi^2(\tau) d\tau\le 
\int\limits_t^T \psi^2(\tau) d\tau =K < \infty
$$

\vspace{2mm}
\noindent
$\forall p\in {\bf N},$ where constant $K$ does not depend on $p$ and $s$ $(t\le s\le T)$.

Thus, the conditions 
(\ref{novemberxxx1}), (\ref{novemberxxx2}) are fulfilled and Theorem~2.52
is true for the case $\psi_1(\tau),\psi_2(\tau),\psi_3(\tau)\equiv \psi(\tau)\in L_2([t, T])$.
Therefore, Theorem~2.53 is also true for the case
$\psi_1(\tau),\psi_2(\tau),\psi_3(\tau)\equiv \psi(\tau),$
where $\psi(\tau)$ is a
continuous function on $[t, T].$

By analogy with (\ref{march00049}) and (\ref{april50})--(\ref{april52}) we obtain the inequalities 
(\ref{novemberxxx3})--(\ref{novemberxxx6})
for the case $\psi_1(\tau),\dots,\psi_5(\tau)\equiv \psi(\tau)\in L_2([t, T])$.
Thus, Theorems~2.54, 2.56 are true for the above case.
Moreover, Theorems~2.55, 2.57
are true for the case 
$\psi_1(\tau),\dots,\psi_5(\tau)\equiv \psi(\tau),$
where $\psi(\tau)$ is a
continuous function on $[t, T].$

Generalizations 
(for the case $\psi_1(\tau),\dots,\psi_6(\tau)\equiv \psi(\tau)\in L_2([t, T])$)
of the relations (\ref{2024december12})--(\ref{2024december11})
(those that were not mentioned earlier in this section)
are proved similarly to (\ref{2024december12})--(\ref{2024december11}) 
(see Sect.~2.32).
This means that Theorem~2.62 is generalized to the case
$\psi_1(\tau),\dots,\psi_6(\tau)\equiv \psi(\tau),$
where $\psi(\tau)$ is a
continuous function on $[t, T].$

In addition to all that has been said, we note that the 
proofs of Theorems~2.64, 2.65 can be easily modified to the case 
$\psi_1(\tau),\dots,\psi_8(\tau)\equiv \psi(\tau),$
where $\psi(\tau)$ is a
continuously dif\-ferentiable 
function on $[t, T]$.
In this case, it is necessary to use (\ref{101oh}), (\ref{101xx}),
(\ref{after1940})
as well as the 
following 
estimate
for the case of Legendre polynomials
$$
\left|
\int\limits_v^x\psi(\tau)\phi_{j}(\tau)d\tau
\right| <
$$
$$
<
\frac{C}{j^{1-\varepsilon/2}}\Biggl(\frac{1}{(1-z^2(x))^{1/4-\varepsilon/4}}+
\frac{1}{(1-z^2(v))^{1/4-\varepsilon/4}}\Biggr),
$$

\vspace{1mm}
\noindent
where 
$j\in {\bf N},$ $z(x), z(v)\in (-1, 1)$ 
($z(x)$ is defined by (\ref{zz1})),
$x, v\in (t, T),$ $\varepsilon\in (0,1)$ is an arbitrary
small positive real number, $\psi(\tau)$ 
is a
continuously dif\-ferentiable 
function on $[t, T],$
constant $C$ does not depend on $j$.

\section{Convergence of the Expansion (\ref{march000195})
to the Iterated Stratonovich Stochastic Integrals
in the Sense of Mathematical Expectation}

In the previous sections, we actually proved
that the value
$$
\sum\limits_{j_1,\ldots,j_k=0}^{p}
C_{j_k \ldots j_1}\prod\limits_{l=1}^k \zeta_{j_l}^{(i_l)}
$$

\noindent
converges if $p\to\infty$ (under suitable conditions)
to the iterated Stratonovich stochastic integrals (\ref{strxx})
in the sense of mathematical expectation.
Let us explain this fact in more detail.

Suppose that $\psi_1(\tau),\ldots,\psi_k(\tau)$ ($k\in{\bf N}$) are
continuous functions on $[t, T]$ and consider Theorem~2.12.

First, let $k=2q+1,$ $q\in {\bf N}.$ We represent (w.~p.~1)
each stochastic integral $J[\psi^{(k)}]_{T,t}^{s_r,\ldots,s_1}$
from the right-hand side of (\ref{30.4})
using the transformation (\ref{febr15})
as a finite linear combination
of the iterated It\^{o} stochastic integrals. Thus, we have (see (\ref{30.4}))
\begin{equation}
\label{march000196}
{\sf M}\left\{J^{*}[\psi^{(k)}]_{T,t}\right\}=0,
\end{equation}

\vspace{1mm}
\noindent
where $J^{*}[\psi^{(k)}]_{T,t}$ is defined by
(\ref{strxx}). On the other hand,
\begin{equation}
\label{march000197}
{\sf M}\left\{\sum\limits_{j_1,\ldots,j_k=0}^{p}
C_{j_k \ldots j_1}\prod\limits_{l=1}^k \zeta_{j_l}^{(i_l)}
\right\}=0,
\end{equation}

\noindent
since $\zeta_{j_l}^{(i_l)}$ has Gaussian distribution and
$k=2q+1,$ $q\in {\bf N}.$

Combining (\ref{march000196}) and (\ref{march000197}), we obtain
\begin{equation}
\label{march000301}
~~~~~~~~~~~~\lim\limits_{p\to\infty}\left|{\sf M}\left\{
J^{*}[\psi^{(k)}]_{T,t}-
\sum\limits_{j_1,\ldots,j_k=0}^{p}
C_{j_k \ldots j_1}\prod\limits_{l=1}^k \zeta_{j_l}^{(i_l)}\right\}\right|=0.
\end{equation}

\vspace{1mm}

Now let $k=2q,$ $q\in {\bf N}.$ In this case, 
using the above reasoning, we get (see (\ref{30.4}))
$$
{\sf M}\left\{J^{*}[\psi^{(k)}]_{T,t}\right\}=
$$

\vspace{-2mm}
$$
=
\frac{1}{2^q}{\bf 1}_{\{i_1=i_2\ne 0\}}
{\bf 1}_{\{i_3=i_4\ne 0\}}\ldots {\bf 1}_{\{i_{2q-1}=i_{2q}\ne 0\}}\times
$$

\vspace{-5mm}
\begin{equation}
\label{march000199}
\times
\int\limits_t^T\psi_{2q}(t_{2q})\psi_{2q-1}(t_{2q}) \ldots 
\int\limits_t^{t_6}
\psi_4(t_4)\psi_3(t_4)
\int\limits_t^{t_4}
\psi_2(t_2)\psi_1(t_2)dt_2 dt_4\ldots  dt_{2q}.
\end{equation}

Recall that the multiple Wiener stochastic 
integral (\ref{WiI}) has zero expectation (see (\ref{Wi110})).
Then, using (\ref{after8xxds1}), (\ref{july30016}) and (\ref{march000199}), we have

\vspace{-3mm}
$$
\lim\limits_{p\to\infty}
{\sf M}\left\{\sum\limits_{j_1,\ldots,j_k=0}^{p}
C_{j_k \ldots j_1}\prod\limits_{l=1}^k \zeta_{j_l}^{(i_l)}\right\}=
$$

\vspace{-3mm}
$$
={\bf 1}_{\{i_1=i_2\ne 0\}}{\bf 1}_{\{i_3=i_4\ne 0\}}
\ldots {\bf 1}_{\{i_{2q-1}=i_{2q}\ne 0\}}\lim\limits_{p\to\infty}
\sum_{j_q,j_{q-2},\ldots, j_2=0}^{p}
C_{j_q j_q j_{q-2} j_{q-2} \ldots j_2 j_2}=
$$

\vspace{1mm}
$$
=
\frac{1}{2^q}{\bf 1}_{\{i_1=i_2\ne 0\}}
{\bf 1}_{\{i_3=i_4\ne 0\}}\ldots {\bf 1}_{\{i_{2q-1}=i_{2q}\ne 0\}}\times
$$

\vspace{-3.5mm}
$$
\times
\int\limits_t^T\psi_{2q}(t_{2q})\psi_{2q-1}(t_{2q}) \ldots 
\int\limits_t^{t_6}
\psi_4(t_4)\psi_3(t_4)
\int\limits_t^{t_4}
\psi_2(t_2)\psi_1(t_2)dt_2 dt_4\ldots  dt_{2q}=
$$

\vspace{1mm}

\begin{equation}
\label{march000300}
={\sf M}\left\{J^{*}[\psi^{(k)}]_{T,t}\right\}.
\end{equation}

\vspace{3mm}

Applying (\ref{march000300}), we obtain
$$
\lim\limits_{p\to\infty}
\left|{\sf M}\left\{
J^{*}[\psi^{(k)}]_{T,t}-\sum\limits_{j_1,\ldots,j_k=0}^{p}
C_{j_k \ldots j_1}\prod\limits_{l=1}^k \zeta_{j_l}^{(i_l)}\right\}\right|=
$$
$$
=\left|
{\sf M}\left\{
J^{*}[\psi^{(k)}]_{T,t}\right\}-
\lim\limits_{p\to\infty}
{\sf M}\left\{\sum\limits_{j_1,\ldots,j_k=0}^{p}
C_{j_k \ldots j_1}\prod\limits_{l=1}^k \zeta_{j_l}^{(i_l)}\right\}\right|=
$$

$$
=0.
$$

\vspace{3mm}

The equality (\ref{march000301}) is proved.

\section{Comparison of Theorems 2.2 and 2.7 with the Representations of Iterated
Stratonovich Stochastic Integrals With Respect to the Scalar Standard Wiener Process}

Note that the correctness of the formulas (\ref{jes}) and (\ref{feto19000}) 
can be 
verified 
in the following way.
If 
$i_1=i_2=i_3=i=1,\ldots,m$
and $\psi_1(\tau),\psi_2(\tau),\psi_3(\tau)\equiv \psi(\tau)$,
then we can derive from (\ref{jes}) and (\ref{feto19000}) 
the well known
equalities (see Sect.~6.7)
$$
{\int\limits_t^{*}}^T\psi(t_{2})
{\int\limits_t^{*}}^{t_{2}}\psi(t_1)
d{\bf w}_{t_1}^{(i)}d{\bf w}_{t_2}^{(i)}=
\frac{1}{2!}\left(\int\limits_t^T\psi(\tau)d{\bf w}_{\tau}^{(i)}\right)^2,
$$
$$
{\int\limits_t^{*}}^T\psi(t_{3})
{\int\limits_t^{*}}^{t_{3}}\psi(t_2)
{\int\limits_t^{*}}^{t_{2}}\psi(t_1)
d{\bf w}_{t_1}^{(i)}d{\bf w}_{t_2}^{(i)}d{\bf w}_{t_3}^{(i)}=
\frac{1}{3!}\left(\int\limits_t^T\psi(\tau)d{\bf w}_{\tau}^{(i)}\right)^3
$$

\noindent
w.~p.~1,
where $\psi(\tau)$ is a continuous
nonrandom function at the 
interval $[t, T]$.

From (\ref{jes}) (under the above assumptions and $p_1=p_2=p$)
we have (see (\ref{best1}) and (\ref{nahod0}))
$$
J^{*}[\psi^{(2)}]_{T,t}=
\hbox{\vtop{\offinterlineskip\halign{
\hfil#\hfil\cr
{\rm l.i.m.}\cr
$\stackrel{}{{}_{p\to \infty}}$\cr
}} }
\sum_{j_1,j_2=0}^{p}
C_{j_2j_1}\zeta_{j_1}^{(i)}\zeta_{j_2}^{(i)}=
$$
$$
=
\hbox{\vtop{\offinterlineskip\halign{
\hfil#\hfil\cr
{\rm l.i.m.}\cr
$\stackrel{}{{}_{p\to \infty}}$\cr
}} }\left(
\sum_{j_1=0}^{p}\sum_{j_2=0}^{j_1-1}\biggl(
C_{j_2j_1}+C_{j_1j_2}\biggr)
\zeta_{j_1}^{(i)}\zeta_{j_2}^{(i)}+
\sum_{j_1=0}^{p}
C_{j_1j_1}\left(\zeta_{j_1}^{(i)}\right)^2\right)
=
$$
$$
=\hbox{\vtop{\offinterlineskip\halign{
\hfil#\hfil\cr
{\rm l.i.m.}\cr
$\stackrel{}{{}_{p\to \infty}}$\cr
}} }\left(
\sum_{j_1=0}^{p}\sum_{j_2=0}^{j_1-1}
C_{j_1}C_{j_2}
\zeta_{j_1}^{(i)}\zeta_{j_2}^{(i)}+
\frac{1}{2}\sum_{j_1=0}^{p}
C_{j_1}^2\left(\zeta_{j_1}^{(i)}\right)^2\right)=
$$
$$
=\hbox{\vtop{\offinterlineskip\halign{
\hfil#\hfil\cr
{\rm l.i.m.}\cr
$\stackrel{}{{}_{p\to \infty}}$\cr
}} }\left(
\frac{1}{2}
\sum_{\stackrel{j_1,j_2=0}{{}_{j_1\ne j_2}}}^{p}
C_{j_1}C_{j_2}
\zeta_{j_1}^{(i)}\zeta_{j_2}^{(i)}+
\frac{1}{2}\sum_{j_1=0}^{p}
C_{j_1}^2\left(\zeta_{j_1}^{(i)}\right)^2\right)=
$$
\begin{equation}
\label{pipi20a}
~~~~~~~=
\hbox{\vtop{\offinterlineskip\halign{
\hfil#\hfil\cr
{\rm l.i.m.}\cr
$\stackrel{}{{}_{p\to \infty}}$\cr
}} }
\frac{1}{2}
\left(\sum_{j_1=0}^{p}
C_{j_1}\zeta_{j_1}^{(i)}\right)^2=
\frac{1}{2!}\left(\int\limits_t^T\psi(\tau)d{\bf w}_{\tau}^{(i)}\right)^2
\end{equation}

\noindent
w.~p.~1. Note that the last step in (\ref{pipi20a}) is performed
by analogy with (\ref{pipi20}).

From (\ref{feto19000}) (under the above assumptions) 
we obtain (see (\ref{oop12}) and (\ref{nahod})--(\ref{nahod2}))
$$
J^{*}[\psi^{(3)}]_{T,t}
=
\hbox{\vtop{\offinterlineskip\halign{
\hfil#\hfil\cr
{\rm l.i.m.}\cr
$\stackrel{}{{}_{p\to \infty}}$\cr
}} }
\sum_{j_1,j_2,j_3=0}^{p}
C_{j_3j_2j_1}\zeta_{j_1}^{(i)}
\zeta_{j_2}^{(i)}
\zeta_{j_3}^{(i)}=
$$
$$
=
\hbox{\vtop{\offinterlineskip\halign{
\hfil#\hfil\cr
{\rm l.i.m.}\cr
$\stackrel{}{{}_{p\to \infty}}$\cr
}} }\left(
\sum_{j_1=0}^{p}
\sum_{j_2=0}^{j_1-1}
\sum_{j_3=0}^{j_2-1}
\biggl(C_{j_3j_2j_1}+
C_{j_3j_1j_2}+C_{j_2j_1j_3}+
C_{j_2j_3j_1}+
C_{j_1j_2j_3}+
C_{j_1j_3j_2}\biggr)\times\right.
$$
$$
\times
\zeta_{j_1}^{(i)}
\zeta_{j_2}^{(i)}
\zeta_{j_3}^{(i)}+
$$
$$
+\sum_{j_1=0}^{p}\sum_{j_3=0}^{j_1-1}\biggl(
C_{j_3j_1j_3}+
C_{j_1j_3j_3}+
C_{j_3j_3j_1}\biggr)\left(\zeta_{j_3}^{(i)}\right)^2
\zeta_{j_1}^{(i)}+
$$
$$
+
\sum_{j_1=0}^{p}\sum_{j_3=0}^{j_1-1}\biggl(
C_{j_3j_1j_1}+
C_{j_1j_1j_3}+
C_{j_1j_3j_1}\biggr)\left(\zeta_{j_1}^{(i)}\right)^2
\zeta_{j_3}^{(i)}
\left.+
\sum_{j_1=0}^p C_{j_1j_1j_1}
\left(\zeta_{j_1}^{(i)}\right)^3\right)=
$$
$$
=
\hbox{\vtop{\offinterlineskip\halign{
\hfil#\hfil\cr
{\rm l.i.m.}\cr
$\stackrel{}{{}_{p\to \infty}}$\cr
}} }\left(
\sum_{j_1=0}^{p}
\sum_{j_2=0}^{j_1-1}
\sum_{j_3=0}^{j_2-1}
C_{j_1}
C_{j_2}C_{j_3}
\zeta_{j_1}^{(i)}
\zeta_{j_2}^{(i)}
\zeta_{j_3}^{(i)}+\right.
$$
$$
+
\frac{1}{2}\sum_{j_1=0}^{p}\sum_{j_3=0}^{j_1-1}
C_{j_3}^2C_{j_1}\left(\zeta_{j_3}^{(i)}\right)^2
\zeta_{j_1}^{(i)}
+\frac{1}{2}\sum_{j_1=0}^{p}\sum_{j_3=0}^{j_1-1}
C_{j_1}^2C_{j_3}\left(\zeta_{j_1}^{(i)}\right)^2
\zeta_{j_3}^{(i)}+
$$
$$
\left.+\frac{1}{6}\sum_{j_1=0}^p C_{j_1}^3
\left(\zeta_{j_1}^{(i)}\right)^3\right)=
$$
$$
=
\hbox{\vtop{\offinterlineskip\halign{
\hfil#\hfil\cr
{\rm l.i.m.}\cr
$\stackrel{}{{}_{p\to \infty}}$\cr
}} }\Biggl(
\frac{1}{6}\sum_{\stackrel{j_1,j_2,j_3=0}{{}_{j_1\ne j_2, j_2\ne j_3,
j_1\ne j_3}}}^{p}
C_{j_1}
C_{j_2}C_{j_3}
\zeta_{j_1}^{(i)}
\zeta_{j_2}^{(i)}
\zeta_{j_3}^{(i)}+\Biggr.
$$
$$
+\frac{1}{2}\sum_{j_1=0}^{p}\sum_{j_3=0}^{j_1-1}
C_{j_3}^2C_{j_1}\left(\zeta_{j_3}^{(i)}\right)^2
\zeta_{j_1}^{(i)}
+\frac{1}{2}\sum_{j_1=0}^{p}\sum_{j_3=0}^{j_1-1}
C_{j_1}^2C_{j_3}\left(\zeta_{j_1}^{(i)}\right)^2
\zeta_{j_3}^{(i)}+
$$
$$
\Biggl.+\frac{1}{6}\sum_{j_1=0}^p C_{j_1}^3
\left(\zeta_{j_1}^{(i)}\right)^3\Biggr)
=
$$
$$
=
\hbox{\vtop{\offinterlineskip\halign{
\hfil#\hfil\cr
{\rm l.i.m.}\cr
$\stackrel{}{{}_{p\to \infty}}$\cr
}} }\left(
\frac{1}{6}\sum_{j_1,j_2,j_3=0}^{p}
C_{j_1}C_{j_2}C_{j_3}\zeta_{j_1}^{(i)}
\zeta_{j_2}^{(i)}\zeta_{j_3}^{(i)}-\right.
$$
$$
-\frac{1}{6}\left(
3\sum_{j_1=0}^{p}\sum_{j_3=0}^{j_1-1}
C_{j_3}^2C_{j_1}\left(\zeta_{j_3}^{(i)}\right)^2
\zeta_{j_1}^{(i)}
+3\sum_{j_1=0}^{p}\sum_{j_3=0}^{j_1-1}
C_{j_1}^2C_{j_3}\left(\zeta_{j_1}^{(i)}\right)^2
\zeta_{j_3}^{(i)}+\right.
$$
$$
\left.+
\sum_{j_1=0}^p C_{j_1}^3
\left(\zeta_{j_1}^{(i)}\right)^3\right)+
$$
$$
+\frac{1}{2}\sum_{j_1=0}^{p}\sum_{j_3=0}^{j_1-1}
C_{j_3}^2C_{j_1}\left(\zeta_{j_3}^{(i)}\right)^2
\zeta_{j_1}^{(i)}
+\frac{1}{2}\sum_{j_1=0}^{p}\sum_{j_3=0}^{j_1-1}
C_{j_1}^2C_{j_3}\left(\zeta_{j_1}^{(i)}\right)^2
\zeta_{j_3}^{(i)}+
$$
$$
\left.+\frac{1}{6}\sum_{j_1=0}^p C_{j_1}^3
\left(\zeta_{j_1}^{(i)}\right)^3\right)=
$$
\begin{equation}
\label{pipi22a}
~~~~~~=
\hbox{\vtop{\offinterlineskip\halign{
\hfil#\hfil\cr
{\rm l.i.m.}\cr
$\stackrel{}{{}_{p\to \infty}}$\cr
}} }
\frac{1}{6}\left(\sum_{j_1=0}^{p}
C_{j_1}
\zeta_{j_1}^{(i)}\right)^3
=
\frac{1}{3!}\left(\int\limits_t^T\psi(\tau)d{\bf w}_{\tau}^{(i)}\right)^3
\end{equation}

\noindent
w.~p.~1. Note that the last step in (\ref{pipi22a}) is performed
by analogy with (\ref{pipi22}).

\section{Algorithm of the Proof of Hypothesis~2.2}

Sect.~2.27--2.39 were written recently, namely in 2024-2025.
At the same time, this section (Sect.~2.41) 
reflects the author's vision of the problem 
under consideration in 2021-2022.

Let us make some remarks about the development of the 
approach based on Theorem~2.30 and describe the algorithm
of the proof of Hypothesis~2.2 (see Sect.~2.5).
First, consider the case $k=2n+1,$ $n=3, 4, \ldots$
($k$ is the multiplicity 
of the iterated Stratonovich stochastic integral (\ref{afterstr})).
Let Condition 2 of Theorem~2.30 be satisfied
(Condition~1 of this theorem is satisfied automatically (see the proof of Theorem~2.18)).
Consider the equality (\ref{drdr1000}). The right-hand side of (\ref{drdr1000})
has the form

\newpage
\noindent
$$
\sum\limits_{j_{g_1},j_{g_3},\ldots, j_{g_{2r-1}}=0}^p
C_{j_k\ldots j_1}\biggl|_{j_{g_1}=j_{g_2},\ldots, j_{g_{2r-1}}=j_{g_{2r}}}-
$$
$$
-\frac{1}{2^r} \prod\limits_{l=1}^r {\bf 1}_{\{g_{2l}=g_{2l-1}+1\}}
C_{j_k \ldots j_1}\biggl|_{(j_{g_2} j_{g_1})\curvearrowright (\cdot)
\ldots (j_{g_{2r}} j_{g_{2r-1}})\curvearrowright (\cdot),
j_{g_{{}_{1}}}=~j_{g_{{}_{2}}},\ldots, j_{g_{{}_{2r-1}}}=~j_{g_{{}_{2r}}}
}\biggr..
$$

\vspace{3mm}

Iterated application
of the formulas (\ref{after81}), (\ref{after82}), (\ref{after9031})
separately to the values
$$
\sum\limits_{j_{g_1},j_{g_3},\ldots, j_{g_{2r-1}}=0}^p
C_{j_k\ldots j_1}\biggl|_{j_{g_1}=j_{g_2},\ldots, j_{g_{2r-1}}=j_{g_{2r}}}
$$

\vspace{-2mm}
\noindent
and 
$$
\frac{1}{2^r} \prod\limits_{l=1}^r {\bf 1}_{\{g_{2l}=g_{2l-1}+1\}}
C_{j_k \ldots j_1}\biggl|_{(j_{g_2} j_{g_1})\curvearrowright (\cdot)
\ldots (j_{g_{2r}} j_{g_{2r-1}})\curvearrowright (\cdot),
j_{g_{{}_{1}}}=~j_{g_{{}_{2}}},\ldots, j_{g_{{}_{2r-1}}}=~j_{g_{{}_{2r}}}
}\biggr.
$$

\vspace{3mm}
\noindent
($g_1, g_2,\ldots, g_{2r-1}, g_{2r}$ as in (\ref{leto5007after}), $r=1,2,\ldots,[k/2],$ $2r<k$)
gives the following representation (see (\ref{drdr1001}))

\vspace{-1mm}
$$
\sum\limits_{\stackrel{j_1,\ldots,j_q,\ldots,j_k=0}{{}_{q\ne g_1, g_2, \ldots, g_{2r-1},
g_{2r}}}}^p
\Biggl(\sum\limits_{j_{g_1},j_{g_3},\ldots, j_{g_{2r-1}}=0}^p
C_{j_k\ldots j_1}\biggl|_{j_{g_1}=j_{g_2},\ldots, j_{g_{2r-1}}=j_{g_{2r}}}-\Biggr.
$$

\vspace{-2mm}
$$
\Biggl.-\frac{1}{2^r} \prod\limits_{l=1}^r {\bf 1}_{\{g_{2l}=g_{2l-1}+1\}}
C_{j_k \ldots j_1}\biggl|_{(j_{g_2} j_{g_1})\curvearrowright (\cdot)
\ldots (j_{g_{2r}} j_{g_{2r-1}})\curvearrowright (\cdot),
j_{g_{{}_{1}}}=~j_{g_{{}_{2}}},\ldots, j_{g_{{}_{2r-1}}}=~j_{g_{{}_{2r}}}
}\biggr.\Biggr)^2\le
$$

\vspace{5mm}
$$
\le \sum\limits_{\stackrel{j_1,\ldots,j_q,\ldots,j_k=0}{{}_{q\ne g_1, g_2, \ldots, g_{2r-1},
g_{2r}}}}^{\infty}
\Biggl(\sum\limits_{j_{g_1},j_{g_3},\ldots, j_{g_{2r-1}}=0}^p
C_{j_k\ldots j_1}\biggl|_{j_{g_1}=j_{g_2},\ldots, j_{g_{2r-1}}=j_{g_{2r}}}-\Biggr.
$$

\vspace{-2mm}
$$
\Biggl.-\frac{1}{2^r} \prod\limits_{l=1}^r {\bf 1}_{\{g_{2l}=g_{2l-1}+1\}}
C_{j_k \ldots j_1}\biggl|_{(j_{g_2} j_{g_1})\curvearrowright (\cdot)
\ldots (j_{g_{2r}} j_{g_{2r-1}})\curvearrowright (\cdot),
j_{g_{{}_{1}}}=~j_{g_{{}_{2}}},\ldots, j_{g_{{}_{2r-1}}}=~j_{g_{{}_{2r}}}
}\biggr.\Biggr)^2=
$$

\vspace{1mm}
$$
=\sum\limits_{\stackrel{j_1,\ldots,j_q,\ldots,j_k=0}{{}_{q\ne g_1, g_2, \ldots, g_{2r-1},
g_{2r}}}}^{\infty}
\left(~
\int\limits_{[t, T]^{k-2r}} 
R_p(t_1,\ldots,t_{g_1-1},t_{g_1+1},\ldots, t_{g_{2r}-1},t_{g_{2r}+1},\ldots,t_k)\times\right.
$$
\begin{equation}
\label{pars0}
\left.\times 
\prod_{\stackrel{q=1}{{}_{q\ne g_1, g_2, \ldots, g_{2r-1},
g_{2r}}}}^{k}
\psi_{q}(t_q)\phi_{j_q}(t_q)\ 
dt_1\ldots dt_{g_1-1}dt_{g_1+1}\ldots dt_{g_{2r}-1}dt_{g_{2r}+1}\ldots dt_k\right)^2,
\end{equation}

\vspace{1mm}
\noindent
where
$$
\int\limits_{[t, T]^{k-2r}} 
R_p(t_1,\ldots,t_{g_1-1},t_{g_1+1},\ldots, 
t_{g_{2r}-1},t_{g_{2r}+1},\ldots,t_k)
\times
$$
$$
\times 
\prod_{\stackrel{q=1}{{}_{q\ne g_1, g_2, \ldots, g_{2r-1},
g_{2r}}}}^{k}
\psi_{q}(t_q)\phi_{j_q}(t_q)\ 
dt_1\ldots dt_{g_1-1}dt_{g_1+1}\ldots dt_{g_{2r}-1}dt_{g_{2r}+1}\ldots dt_k
$$

\vspace{1mm}
\noindent
is the Fourier coefficient of 
$$
\hat R_p(t_1,\ldots,t_{g_1-1},t_{g_1+1},\ldots, t_{g_{2r}-1},t_{g_{2r}+1},\ldots,t_k)=
$$
$$
=
R_p(t_1,\ldots,t_{g_1-1},t_{g_1+1},\ldots, t_{g_{2r}-1},t_{g_{2r}+1},\ldots,t_k)
\prod_{\stackrel{q=1}{{}_{q\ne g_1, g_2, \ldots, g_{2r-1},
g_{2r}}}}^{k}\psi_{q}(t_q),
$$

\noindent
where 
$$
R_p(t_1,\ldots,t_{g_1-1},t_{g_1+1},\ldots, t_{g_{2r}-1},t_{g_{2r}+1},\ldots,t_k)=
$$
$$
=\sum\limits_{d=1}^{4^r}
\bar R_p^{(d)}(t_1,\ldots,t_{g_1-1},t_{g_1+1},\ldots, t_{g_{2r}-1},t_{g_{2r}+1},\ldots,t_k)-
$$
$$
-
\sum\limits_{d=1}^{2^r}
\tilde R_p^{(d)}(t_1,\ldots,t_{g_1-1},t_{g_1+1},\ldots, t_{g_{2r}-1},t_{g_{2r}+1},\ldots,t_k)
\in L_2([t, T]^{k-2r})
$$

\noindent
and some of the func\-ti\-ons 
$\bar R_p^{(d)}(t_1,\ldots,t_{g_1-1},t_{g_1+1},\ldots, t_{g_{2r}-1},t_{g_{2r}+1},\ldots,t_k)$ and
$\tilde R_p^{(d)}(t_1,\ldots,t_{g_1-1},t_{g_1+1},\ldots, t_{g_{2r}-1},t_{g_{2r}+1},\ldots,t_k)$
can be iden\-ti\-cal\-ly equal to zero.
Obviously, we could use another representation for the function

\vspace{-3mm}
\begin{equation}
\label{de200}
R_p(t_1,\ldots,t_{g_1-1},t_{g_1+1},\ldots, t_{g_{2r}-1},t_{g_{2r}+1},\ldots,t_k)
\end{equation}

\vspace{3mm}
\noindent
based on the left-hand side of the equality (\ref{drdr1000})
and (\ref{after81}), (\ref{after82}), (\ref{after9031}) (see Sect.~2.13 for details).
In Sect.~2.13, we considered the function (\ref{de200}) in detail
for the case $k\ge 5,$ $r=1.$

Parseval's equality gives
$$
\sum\limits_{\stackrel{j_1,\ldots,j_q,\ldots,j_k=0}{{}_{q\ne g_1, g_2, \ldots, g_{2r-1},
g_{2r}}}}^{\infty}
\left(~
\int\limits_{[t, T]^{k-2r}} 
R_p(t_1,\ldots,t_{g_1-1},t_{g_1+1},\ldots, t_{g_{2r}-1},t_{g_{2r}+1},\ldots,t_k)\times\right.
$$
$$
\left.\times 
\prod_{\stackrel{q=1}{{}_{q\ne g_1, g_2, \ldots, g_{2r-1},
g_{2r}}}}^{k}
\psi_q(t_q)\phi_{j_q}(t_q)\ 
dt_1\ldots dt_{g_1-1}dt_{g_1+1}\ldots dt_{g_{2r}-1}dt_{g_{2r}+1}\ldots dt_k\right)^2=
$$

\vspace{1mm}
$$
=\int\limits_{[t, T]^{k-2r}} 
\hspace{-2mm}\left(
\hat R_p(t_1,\ldots,t_{g_1-1},t_{g_1+1},\ldots, t_{g_{2r}-1},t_{g_{2r}+1},\ldots,t_k)\right)^2
dt_1\ldots dt_{g_1-1}dt_{g_1+1}\ldots
$$

\vspace{-2mm}
\begin{equation}
\label{pars1}
\ldots dt_{g_{2r}-1}dt_{g_{2r}+1}\ldots dt_k 
=\bigl\Vert \hat R_p \bigr\Vert_{L_2([t, T]^{k-2r})}^2.
\end{equation}

\vspace{6mm}

Combining (\ref{pars0}) and (\ref{pars1}), we obtain
$$
\sum\limits_{\stackrel{j_1,\ldots,j_q,\ldots,j_k=0}{{}_{q\ne g_1, g_2, \ldots, g_{2r-1},
g_{2r}}}}^p
\Biggl(\sum\limits_{j_{g_1},j_{g_3},\ldots, j_{g_{2r-1}}=0}^p
C_{j_k\ldots j_1}\biggl|_{j_{g_1}=j_{g_2},\ldots, j_{g_{2r-1}}=j_{g_{2r}}}-\Biggr.
$$

\vspace{-6mm}
$$
\Biggl.-\frac{1}{2^r} \prod\limits_{l=1}^r {\bf 1}_{\{g_{2l}=g_{2l-1}+1\}}
C_{j_k \ldots j_1}\biggl|_{(j_{g_2} j_{g_1})\curvearrowright (\cdot)
\ldots (j_{g_{2r}} j_{g_{2r-1}})\curvearrowright (\cdot),
j_{g_{{}_{1}}}=~j_{g_{{}_{2}}},\ldots, j_{g_{{}_{2r-1}}}=~j_{g_{{}_{2r}}}
}\biggr.\Biggr)^2\le
$$

\vspace{3mm}
\begin{equation}
\label{pars2}
\le \bigl\Vert \hat R_p \bigr\Vert_{L_2([t, T]^{k-2r})}^2.
\end{equation}

\vspace{4mm}

Assume that we have succeeded in proving the following equality
\begin{equation}
\label{pars3s}
\lim\limits_{p\to\infty}
\bigl\Vert \hat R_p \bigr\Vert_{L_2([t, T]^{k-2r})}^2=0.
\end{equation}

Applying (\ref{pars2}) and (\ref{pars3s}), we get (compare with (\ref{drdr1001}))
$$
\lim\limits_{p\to\infty}
\sum\limits_{\stackrel{j_1,\ldots,j_q,\ldots,j_k=0}{{}_{q\ne g_1, g_2, \ldots, g_{2r-1},
g_{2r}}}}^p
\Biggl(\sum\limits_{j_{g_1},j_{g_3},\ldots, j_{g_{2r-1}}=0}^p
C_{j_k\ldots j_1}\biggl|_{j_{g_1}=j_{g_2},\ldots, j_{g_{2r-1}}=j_{g_{2r}}}-\Biggr.
$$
\begin{equation}
\label{for900}
\Biggl.-\frac{1}{2^r} \prod\limits_{l=1}^r {\bf 1}_{\{g_{2l}=g_{2l-1}+1\}}
C_{j_k \ldots j_1}\biggl|_{(j_{g_2} j_{g_1})\curvearrowright (\cdot)
\ldots (j_{g_{2r}} j_{g_{2r-1}})\curvearrowright (\cdot),
j_{g_{{}_{1}}}=~j_{g_{{}_{2}}},\ldots, j_{g_{{}_{2r-1}}}=~j_{g_{{}_{2r}}}
}\biggr.\Biggr)^2=0.
\end{equation}

\vspace{2mm}

As noted in Sect.~2.10, Condition 3 of Theorem~2.30 can be replaced by a weaker condition 
(\ref{drdr1001}) (or (\ref{for900})).
Also Condition 3 of Theorem~2.30 can be replaced by (\ref{pars3s}).
From (\ref{for900}) we obviously obtain for $2r<k$
$$
\lim\limits_{p\to\infty} \sum\limits_{j_{g_1},j_{g_3},\ldots, j_{g_{2r-1}}=0}^p
C_{j_k\ldots j_1}\biggl|_{j_{g_1}=j_{g_2},\ldots, j_{g_{2r-1}}=j_{g_{2r}}}=\Biggr.
$$

\vspace{-4mm}
\begin{equation}
\label{pars900}
\Biggl.=\frac{1}{2^r} \prod\limits_{l=1}^r {\bf 1}_{\{g_{2l}=g_{2l-1}+1\}}
C_{j_k \ldots j_1}\biggl|_{(j_{g_2} j_{g_1})\curvearrowright (\cdot)
\ldots (j_{g_{2r}} j_{g_{2r-1}})\curvearrowright (\cdot),
j_{g_{{}_{1}}}=~j_{g_{{}_{2}}},\ldots, j_{g_{{}_{2r-1}}}=~j_{g_{{}_{2r}}}
}\biggr..
\end{equation}

\vspace{4mm}

According to (\ref{drdr1000}), the equality (\ref{pars900}) will be satisfied 
if
\begin{equation}
\label{sars10}
\lim\limits_{p\to\infty}
S_{l_1}S_{l_2}\ldots S_{l_{d}}
\left\{\bar C^{(p)}_{j_k\ldots j_q \ldots j_1}\biggl|_{q\ne g_1,g_2,\ldots,g_{2r-1}, g_{2r}}
\right\}=0,
\end{equation}

\vspace{2mm}
\noindent
where $g_1,g_2,\ldots,g_{2r-1},g_{2r}$ as in (\ref{leto5007after}),
$l_1, l_2, \ldots, l_{d}$ such that
$l_1, l_2, \ldots, l_{d}\in \{1,2,\ldots, r\},$\
$l_1>l_2>\ldots >l_{d},$\ $d=0, 1, 2,\ldots, r-1,$\ 
$r=1, 2,\ldots,[k/2],$

\vspace{-4mm}
$$
S_{l_1}S_{l_2}\ldots S_{l_{d}}
\left\{\bar C^{(p)}_{j_k\ldots j_q \ldots j_1}\biggl|_{q\ne g_1,g_2,\ldots,g_{2r-1}, g_{2r}}
\right\}\stackrel{\sf def}{=}
\bar C^{(p)}_{j_k\ldots j_q \ldots j_1}\biggl|_{q\ne g_1,g_2,\ldots,g_{2r-1}, g_{2r}}
$$

\vspace{2mm}
\noindent
for $d=0,$ where 
$$
\bar C^{(p)}_{j_k\ldots j_q \ldots j_1}\biggl|_{q\ne g_1,g_2,\ldots,g_{2r-1}, g_{2r}},\ \ \
S_l \left\{\bar C^{(p)}_{j_k\ldots j_q \ldots j_1}\biggl|_{q\ne g_1,g_2,\ldots,g_{2r-1}, g_{2r}}
\right\}$$

\vspace{3mm}
\noindent
are defined by (\ref{after5days}), (\ref{sars300}), 
$l=1,2,\ldots,r$ (see Sect.~2.10 for details).

Let us make some remarks about the function
(\ref{de200})
for the case $k>5,$ $r=2.$
In this case, using the left-hand side of 
the equality (\ref{drdr1000})
and (\ref{after81}), (\ref{after82}), (\ref{after9031}), 
we represent the function (\ref{de200})
as the sum of several functions.
In particular, among these functions 
will be the following functions

\newpage
\noindent
$$
Q_p(t_1,\ldots,t_{s-1},t_{s+1},\ldots ,t_{l-1},t_{l+1},\ldots,
t_{q-1},t_{q+1},\ldots, t_{g-1},t_{g+1},\ldots,t_k)=
$$

\vspace{-4mm}
$$
=
{\bf 1}_{\{t_1< \ldots <t_{s-1}<t_{s+1}< \ldots <t_{l-1}<t_{l+1}< \ldots
<t_{q-1}<t_{q+1}< \ldots <t_{g-1}<t_{g+1}<\ldots<t_k\}}\times
$$

\vspace{-4mm}
$$
\times
\sum_{j_l=p+1}^{\infty}~
\int\limits_{t}^{t_{s+1}} \psi_s(\tau) \phi_{j_{l}}(\tau)d\tau
\int\limits_{t}^{t_{l-1}} \psi_l(\tau) \phi_{j_{l}}(\tau)d\tau\times
$$
\begin{equation}
\label{func500}
\times
\sum_{j_q=p+1}^{\infty}~
\int\limits_{t}^{t_{q+1}} \psi_q(\tau) \phi_{j_{q}}(\tau)d\tau
\int\limits_{t}^{t_{g-1}} \psi_g(\tau) \phi_{j_{q}}(\tau)d\tau,
\end{equation}

\vspace{6mm}

$$
\bar Q_p(t_1,\ldots,t_{l-2},t_{l+3},\ldots,t_k)=
$$

\vspace{-4mm}

$$
=
{\bf 1}_{\{t_1<\ldots<t_{l-2}<t_{l+3}<\ldots<t_k\}}\times
$$

\vspace{-4mm}
$$
\times
\sum_{j_l=p+1}^{\infty}~\left(\int\limits_{t}^{t_{l-2}} \psi_{l-1}(\theta)\phi_{j_{l}}(\theta)
\int\limits_{t}^{\theta} \psi_l(u)\phi_{j_{l}}(u)du d\theta\right)\times
$$
\begin{equation}
\label{func501}
~~~~~~~~~\times
\sum_{j_q=p+1}^{\infty}~\left(\int\limits_{t}^{t_{l-2}} \psi_{l+1}(\theta)\phi_{j_{q}}(\theta)
\int\limits_{t}^{\theta} \psi_{l+2}(u)\phi_{j_{q}}(u)du d\theta\right),
\end{equation}

\vspace{6mm}
$$
\tilde Q_p(t_1,\ldots,t_{l-2},t_{l+3},\ldots,t_k)=
$$

\vspace{-4mm}

$$
=
{\bf 1}_{\{t_1<\ldots<t_{l-2}<t_{l+3}<\ldots<t_k\}}\times
$$

\vspace{-4mm}
$$
\times
\sum_{j_l=p+1}^{\infty}\sum_{j_q=p+1}^{\infty}~
\int\limits_{t}^{t_{l+3}} \psi_{l+1}(\tau)\phi_{j_{q}}(\tau)
\left(\int\limits_{t}^{\tau} \psi_{l-1}(\theta)\phi_{j_{l}}(\theta)
\int\limits_{t}^{\theta} \psi_l(u)\phi_{j_{l}}(u)du d\theta\right)\times
$$
\begin{equation}
\label{func502}
\times
\int\limits_{t}^{\tau} \psi_{l+2}(u)\phi_{j_{q}}(u)du d\tau,
\end{equation}

\vspace{6mm}
$$
\hat Q_p(t_1,\ldots,t_{l-1},t_{l+2},\ldots, t_{q-1}, t_{q+2}, \ldots, t_k)=
$$

\vspace{-4mm}
$$
=
{\bf 1}_{\{t_1<\ldots<t_{l-1}<t_{l+2}<\ldots<t_{q-1}<t_{q+2}<\ldots <t_k\}}\times
$$

\newpage
\noindent
$$
\times
\sum_{j_{l}=p+1}^{\infty}
\sum_{j_{l+1}=p+1}^{\infty}~\left(\int\limits_{t}^{t_{l+2}} \psi_{l+1}(\theta)\phi_{j_{l+1}}(\theta)
\int\limits_{t}^{\theta} \psi_l(u)\phi_{j_{l}}(u)du d\theta\right)\times
$$
\begin{equation}
\label{func501qq}
~~~~~~~~~\times
\left(\int\limits_{t}^{t_{q+2}} \psi_{q+1}(\theta)\phi_{j_{l+1}}(\theta)
\int\limits_{t}^{\theta} \psi_{q}(u)\phi_{j_{l}}(u)du d\theta\right).
\end{equation}

\vspace{1mm}

Note that the pairs $(g_1,g_2),$ $(g_3,g_4)$ for the functions (\ref{func501}) and (\ref{func502}) 
have the property: $g_2=g_1+1,$ $g_4=g_3+1,$ $g_3=g_2+1.$
At the same time, the pairs $(g_1,g_2),$ $(g_3,g_4)$ for the function (\ref{func500})
have the following property: $g_2>g_1+1,$ $g_4>g_3+1,$ $g_3\ge g_2+1.$
For the function (\ref{func501qq}), the pairs $(g_1,g_2),$ $(g_3,g_4)$
chosen as follows: $g_2>g_1+1,$ $g_4>g_3+1,$ $g_4=g_2+1,$ $g_3=g_1+1.$
Generally speaking, all possible pairs $(g_1,g_2),$ $(g_3,g_4)$ must be considered.
We consider the functions (\ref{func500})--(\ref{func501qq}) only as an example.

Suppose that $s+1=l-1,$ $l+1=q-1,$ $q+1=g-1$ in (\ref{func500}). Let us show that
(we consider the case of Legendre polynomials; the trigonometric case is simpler
and can be considered similarly)

\vspace{-2mm}
\begin{equation}
\label{func503}
\lim\limits_{p\to\infty}\bigl\Vert Q_p \bigr\Vert_{L_2([t, T]^{k-4})}^2=0,
\end{equation}

\vspace{-1mm}
\begin{equation}
\label{func504}
\lim\limits_{p\to\infty}\bigl\Vert \bar Q_p \bigr\Vert_{L_2([t, T]^{k-4})}^2=0,
\end{equation}

\vspace{-1mm}
\begin{equation}
\label{func505}
\lim\limits_{p\to\infty}\bigl\Vert \tilde Q_p \bigr\Vert_{L_2([t, T]^{k-4})}^2=0,
\end{equation}

\vspace{-1mm}
\begin{equation}
\label{func505qq}
\lim\limits_{p\to\infty}\bigl\Vert \hat Q_p \bigr\Vert_{L_2([t, T]^{k-4})}^2=0.
\end{equation}

\vspace{4mm}

First consider the proof of (\ref{func503}). We have
($s+1=l-1,$ $l+1=q-1,$ $q+1=g-1$)

\vspace{-5mm}
$$
\left(Q_p(t_1,\ldots,t_{l-3},t_{l-1},t_{l+1},t_{l+3},t_{l+5},\ldots,t_k)\right)^2=
$$

\vspace{-5mm}
$$
=
{\bf 1}_{\{t_1< \ldots <t_{l-3}<t_{l-1}< t_{l+1}<t_{l+3}<t_{l+5}<\ldots<t_k\}}\times
$$
$$
\times
\left(\sum_{j_l=p+1}^{\infty}~
\int\limits_{t}^{t_{l-1}} \psi_{l-2}(\tau) \phi_{j_{l}}(\tau)d\tau
\int\limits_{t}^{t_{l-1}} \psi_l(\tau) \phi_{j_{l}}(\tau)d\tau\times\right.
$$
\begin{equation}
\label{func506}
~~~~~~~~~~~~~\left.\times \sum_{j_q=p+1}^{\infty}~
\int\limits_{t}^{t_{l+3}} \psi_{l+2}(\tau) \phi_{j_{q}}(\tau)d\tau
\int\limits_{t}^{t_{l+3}} \psi_{l+4}(\tau) \phi_{j_{q}}(\tau)d\tau\right)^2.
\end{equation}

\vspace{2mm}

Using the estimate (\ref{after1940}), we obtain
\begin{equation}
\label{func507}
\left|
\int\limits_t^s\psi(\tau)\phi_{j}(\tau)d\tau
\right| <
\frac{K}{j^{1-\varepsilon/2} (1-z^2(s))^{1/4-\varepsilon/4}},
\end{equation}

\vspace{1mm}
\noindent
where $j\in{\bf N},$ $s\in (t, T),$
$z(s)$ is defined by (\ref{zz1}), $\varepsilon\in (0,1),$ constant $K$ does not depend on $j,$
$\{\phi_j(x)\}_{j=0}^{\infty}$ is a complete orthonormal system of 
Legendre polynomials in the space $L_2([t, T]),$
$\psi(\tau)$ is a continuously dif\-ferentiable 
nonrandom function on $[t, T].$ 

Applying (\ref{func507}) and 
(\ref{after1944}) (we take $\varepsilon$ instead of $\varepsilon/2$ in (\ref{after1944})), we get
$$
\left(\sum_{j_l=p+1}^{\infty}~
\int\limits_{t}^{t_{l-1}} \psi_{l-2}(\tau) \phi_{j_{l}}(\tau)d\tau
\int\limits_{t}^{t_{l-1}} \psi_l(\tau) \phi_{j_{l}}(\tau)d\tau\times\right.
$$
$$
\left.\times \sum_{j_q=p+1}^{\infty}~
\int\limits_{t}^{t_{l+3}} \psi_{l+2}(\tau) \phi_{j_{q}}(\tau)d\tau
\int\limits_{t}^{t_{l+3}} \psi_{l+4}(\tau) \phi_{j_{q}}(\tau)d\tau\right)^2\le
$$

\begin{equation}
\label{func508}
\le \frac{K_1}{p^{4(1-\varepsilon)}(1-z^2(t_{l-1}))^{1-\varepsilon}
(1-z^2(t_{l+3}))^{1-\varepsilon}},
\end{equation}

\vspace{5mm}
\noindent
where $t_{l-1}, t_{l+3}\in (t, T),$ constant $K_1$ is independent of $p.$
Combining (\ref{func506}) and (\ref{func508}), we have (\ref{func503}).

Let us prove (\ref{func504}).
Applying the estimate (\ref{after5000}) in (\ref{agent14}) and tak\-ing into account 
the boundedness of the functions $\psi_1(\tau),$ $\psi_2(\tau)$ 
and their derivatives, we obtain
$$
\left\vert\sum\limits_{j=m+1}^n
C_{j j}(s)\right\vert\le
C_1\left(\frac{1}{n^{1-\varepsilon}}+\frac{1}{m^{1-\varepsilon}}\right)
\int\limits_{-1}^{z(s)} \frac{dx}{\left(1-x^2\right)^{1/2-\varepsilon/2}}+
$$
$$
+C_2 \sum\limits_{j=m+1}^n \frac{1}{j^{2-\varepsilon}}\left(
\int\limits_{-1}^{z(s)}\frac{dy}{\left(1-y^2\right)^{1/2-\varepsilon/2}}
+\frac{1}{\left(1-z^2(s)\right)^{1/4-\varepsilon/4}}
\int\limits_{-1}^{z(s)}\frac{dy}{\left(1-y^2\right)^{1/4-\varepsilon/4}}+\right.
$$
\begin{equation}
\label{func800}
+
\left.\int\limits_{-1}^{z(s)}\frac{1}{\left(1-y^2\right)^{1/4-\varepsilon/4}}
\int\limits_{y}^{z(s)}\frac{dx}{\left(1-x^2\right)^{1/4-\varepsilon/4}}dy\right),
\end{equation}

\vspace{1mm}
\noindent
where 
$$
C_{j j}(s)=\int\limits_t^s\psi_2(\tau)\phi_{j}(\tau)
\int\limits_t^{\tau}\psi_1(\theta)\phi_{j}(\theta)d\theta d\tau,
$$

\vspace{1mm}
\noindent
$s\in (t, T),$ constants $C_1, C_2$ do not depend on $n$ and $m$.

From (\ref{func800}) we have
\begin{equation}
\label{func800x}
~~~~~~~~\left\vert\sum\limits_{j=m+1}^{\infty}
C_{j j}(s)\right\vert\le
\frac{K_1}{m^{1-\varepsilon}}
+K_2 \sum\limits_{j=m+1}^{\infty} \frac{1}{j^{2-\varepsilon}}\left(
1+\frac{1}{\left(1-z^2(s)\right)^{1/4-\varepsilon/4}}\right),
\end{equation}

\vspace{2mm}
\noindent
where $s\in (t, T),$ constants $K_1, K_2$ do not depend on $m$.

Applying (\ref{after1944}) (we take $\varepsilon$ instead of $\varepsilon/2$ 
in (\ref{after1944})) in (\ref{func800x}),  we get
\begin{equation}
\label{func801}
\left\vert\sum\limits_{j=m+1}^{\infty}
C_{j j}(s)\right\vert\le
\frac{K}{m^{1-\varepsilon}\left(1-z^2(s)\right)^{1/4-\varepsilon/4}},
\end{equation}

\noindent
where $s\in (t, T),$ constant $K$ is independent of $m$.

Using the estimate (\ref{func801}), we obtain (see (\ref{func501}))
$$
\left(\bar Q_p(t_1,\ldots,t_{l-2},t_{l+3},\ldots,t_k)\right)^2=
{\bf 1}_{\{t_1<\ldots<t_{l-2}<t_{l+3}<\ldots<t_k\}}\times
$$
$$
\times
\left(\sum_{j_l=p+1}^{\infty}~\left(\int\limits_{t}^{t_{l-2}} \psi_{l-1}(\theta)\phi_{j_{l}}(\theta)
\int\limits_{t}^{\theta} \psi_l(u)\phi_{j_{l}}(u)du d\theta\right)\times\right.
$$
$$
\left.
\times
\sum_{j_q=p+1}^{\infty}~\left(\int\limits_{t}^{t_{l-2}} \psi_{l+1}(\theta)\phi_{j_{q}}(\theta)
\int\limits_{t}^{\theta} \psi_{l+2}(u)\phi_{j_{q}}(u)du d\theta\right)\right)^2\le
$$

\vspace{-2mm}
\begin{equation}
\label{func990}
\le \frac{K_1}{p^{4(1-\varepsilon)}(1-z^2(t_{l-2}))^{1-\varepsilon}},
\end{equation}

\vspace{3mm}
\noindent
where $t_{l-2}\in (t, T),$ constant $K_1$ is independent of $p.$
The inequality (\ref{func990}) completes the proof of (\ref{func504}).

Let us prove (\ref{func505}). Applying (\ref{agent0101}) in (\ref{func502}), we get
$$
\left(\tilde Q_p(t_1,\ldots,t_{l-2},t_{l+3},\ldots,t_k)\right)^2\le
$$
$$
\le
\left(\sum_{j_l=p+1}^{\infty}\sum_{j_q=p+1}^{\infty}~
\int\limits_{t}^{t_{l+3}} \psi_{l+1}(\tau)\phi_{j_{q}}(\tau)
\left(\int\limits_{t}^{\tau} \psi_{l-1}(\theta)\phi_{j_{l}}(\theta)
\int\limits_{t}^{\theta} \psi_l(u)\phi_{j_{l}}(u)du d\theta\right)\times\right.
$$
$$
\left.\times
\int\limits_{t}^{\tau} \psi_{l+2}(u)\phi_{j_{q}}(u)du d\tau\right)^2=
$$
$$
=\left(\frac{1}{2}\sum_{j_l=p+1}^{\infty}~
\int\limits_{t}^{t_{l+3}} \psi_{l+1}(\tau)
\left(\int\limits_{t}^{\tau} \psi_{l-1}(\theta)\phi_{j_{l}}(\theta)
\int\limits_{t}^{\theta} \psi_l(u)\phi_{j_{l}}(u)du d\theta\right)
\psi_{l+2}(\tau)d\tau-\right.
$$

\vspace{-3mm}
$$
-\sum_{j_q=0}^{p}~
\int\limits_{t}^{t_{l+3}}\psi_{l+1}(\tau)\phi_{j_{q}}(\tau)
\sum_{j_l=p+1}^{\infty}\left(\int\limits_{t}^{\tau} \psi_{l-1}(\theta)\phi_{j_{l}}(\theta)
\int\limits_{t}^{\theta}\psi_l(u)\phi_{j_{l}}(u)du d\theta\right)\times
$$
\begin{equation}
\label{func1000}
~~~~~~~~~\left.\times
\int\limits_{t}^{\tau} \psi_{l+2}(u)\phi_{j_{q}}(u)du d\tau\right)^2= 
(a-b)^2\le
2(|a|^2 + |b|^2).
\end{equation}

\vspace{2mm}

Further, we have
\begin{equation}
\label{func1001}
|a|\le 
\frac{1}{2}
\int\limits_{t}^{t_{l+3}} \left\vert\psi_{l+1}(\tau)\right\vert
\left\vert\sum_{j_l=p+1}^{\infty}~\int\limits_{t}^{\tau} \psi_{l-1}(\theta)\phi_{j_{l}}(\theta)
\int\limits_{t}^{\theta} \psi_l(u)\phi_{j_{l}}(u)du d\theta\right\vert
\left\vert\psi_{l+2}(\tau)\right\vert d\tau,
\end{equation}

\vspace{-5mm}
$$
|b|\le
\sum_{j_q=0}^{p}~
\int\limits_{t}^{t_{l+3}}\left\vert \psi_{l+1}(\tau)\phi_{j_{q}}(\tau)\right\vert
\left\vert\sum_{j_l=p+1}^{\infty}\int\limits_{t}^{\tau} \psi_{l-1}(\theta)\phi_{j_{l}}(\theta)
\int\limits_{t}^{\theta}\psi_l(u)\phi_{j_{l}}(u)du d\theta\right\vert \times
$$
\begin{equation}
\label{func1002}
\times
\left\vert\int\limits_{t}^{\tau} \psi_{l+2}(u)\phi_{j_{q}}(u)du \right\vert d\tau.
\end{equation}

\vspace{2mm}

Combining (\ref{func801}) and (\ref{func1001}), we obtain
\begin{equation}
\label{func1002x}
|a|\le \frac{C}{p^{1-\varepsilon}},
\end{equation}

\noindent
where constant $C$ is independent of $p.$

Separating in (\ref{func1002}) the term with the number $j_q=0$ and then applying 
(\ref{ogo24}), (\ref{101xx}), (\ref{func801}), we obtain
$$
|b|\le 
\frac{K}{p^{1-\varepsilon}}
\left(\int\limits_t^{t_{l+3}} \frac{d\tau}{\left(1-z^2(\tau)\right)^{1/2-\varepsilon/4}}+
\sum_{j_q=1}^{p}\frac{1}{j_q}
\int\limits_t^{t_{l+3}} \frac{d\tau}{\left(1-z^2(\tau)\right)^{3/4-\varepsilon/4}}\right)\le
$$

$$
\le \frac{K_1}{p^{1-\varepsilon}}\left(1+\sum_{j_q=1}^{p}\frac{1}{j_q}\right)\le
\frac{K_1}{p^{1-\varepsilon}}\left(2+\int\limits_1^p \frac{dx}{x}\right)=
$$

\vspace{3mm}
\begin{equation}
\label{func1003}
=\frac{K_1\left(2+ ln p\right)}{p^{1-\varepsilon}}\ \to\ 0
\end{equation}

\vspace{3mm}

\noindent
if $p\to\infty.$ The estimates (\ref{func1000}), (\ref{func1002x}), (\ref{func1003}) 
complete the proof of (\ref{func505}). 

Finally, consider the proof of (\ref{func505qq}).
Using the elementary inequality $|ab|\le (a^2+b^2)/2$ and Parseval's equality, we have
$$
\left(\hat Q_p(t_1,\ldots,t_{l-1},t_{l+2},\ldots, t_{q-1}, t_{q+2}, \ldots, t_k)\right)^2\le
$$
$$
\le 
\left(\sum_{j_{l}=p+1}^{\infty}
\sum_{j_{l+1}=p+1}^{\infty}~\left\vert\int\limits_{t}^{t_{l+2}} \psi_{l+1}(\theta)\phi_{j_{l+1}}(\theta)
\int\limits_{t}^{\theta} \psi_l(u)\phi_{j_{l}}(u)du d\theta\right\vert\times\right.
$$

$$
\times
\left.\left\vert\int\limits_{t}^{t_{q+2}} \psi_{q+1}(\theta)\phi_{j_{l+1}}(\theta)
\int\limits_{t}^{\theta} \psi_{q}(u)\phi_{j_{l}}(u)du d\theta\right\vert\right)^2\le
$$
$$
\le 
\frac{1}{4}\left(\sum_{j_{l}=p+1}^{\infty}
\sum_{j_{l+1}=p+1}^{\infty}~\left(\int\limits_{t}^{t_{l+2}} \psi_{l+1}(\theta)\phi_{j_{l+1}}(\theta)
\int\limits_{t}^{\theta} \psi_l(u)\phi_{j_{l}}(u)du d\theta\right)^2+\right.
$$
$$
+
\left.
\sum_{j_{l}=p+1}^{\infty}
\sum_{j_{l+1}=p+1}^{\infty}~
\left(\int\limits_{t}^{t_{q+2}} \psi_{q+1}(\theta)\phi_{j_{l+1}}(\theta)
\int\limits_{t}^{\theta} \psi_{q}(u)\phi_{j_{l}}(u)du d\theta\right)^2\right)^2\le
$$
$$
\le 
\frac{1}{4}\left(\sum_{j_{l}=p+1}^{\infty}
\sum_{j_{l+1}=0}^{\infty}~\left(\int\limits_{t}^{t_{l+2}} \psi_{l+1}(\theta)\phi_{j_{l+1}}(\theta)
\int\limits_{t}^{\theta} \psi_l(u)\phi_{j_{l}}(u)du d\theta\right)^2+\right.
$$
$$
+
\left.
\sum_{j_{l}=p+1}^{\infty}
\sum_{j_{l+1}=0}^{\infty}~
\left(\int\limits_{t}^{t_{q+2}} \psi_{q+1}(\theta)\phi_{j_{l+1}}(\theta)
\int\limits_{t}^{\theta} \psi_{q}(u)\phi_{j_{l}}(u)du d\theta\right)^2\right)^2\le
$$
$$
\le 
\frac{1}{4}\left(\sum_{j_{l}=p+1}^{\infty}~\int\limits_{t}^{t_{l+2}} 
\psi_{l+1}^2(\theta)
\left(\int\limits_{t}^{\theta} \psi_l(u)\phi_{j_{l}}(u)du\right)^2 d\theta+\right.
$$
\begin{equation}
\label{func651}
~~~~~~~~~~+
\left.
\sum_{j_{l}=p+1}^{\infty}~\int\limits_{t}^{t_{q+2}} 
\psi_{q+1}^2(\theta)
\left(\int\limits_{t}^{\theta} \psi_{q}(u)\phi_{j_{l}}(u)du\right)^2 d\theta\right)^2.
\end{equation}

\vspace{2mm}

From (\ref{func651}) and (\ref{obana}), (\ref{101xx}) we obtain
$$
\left(\hat Q_p(t_1,\ldots,t_{l-1},t_{l+2},\ldots, t_{q-1}, t_{q+2}, \ldots, t_k)\right)^2\le
\frac{K}{p^2}\ \to\ 0
$$

\noindent
if $p\to\infty,$ where constant $K$ does not depend on $p.$
Thus the equalities 
(\ref{func503})--(\ref{func505qq}) are proved.

Recall that the function 
(\ref{de200}) (this function is defined using the left-hand side of the equality (\ref{drdr1000}))
for the case $k > 5,$ $r=2$
is represented
as the sum of several functions. Four of them, namely $Q_p,$ $\bar Q_p,$ $\tilde Q_p,$
$\hat Q_p$
(these functions correspond to the particular case of choosing the pairs $(g_1,g_2),$ $(g_3,g_4)$;
generally speaking, all possible pairs $(g_1,g_2),$ $(g_3,g_4)$ must be considered),
have been studied above. Absolutely similarly, we can consider
the remaining functions (for all possible pairs $(g_1,g_2),$ $(g_3,g_4)$)
whose sum is the function 
(\ref{de200}) 
for the case $k > 5,$ $r=2.$ As a result, we will have
$$
\lim\limits_{p\to\infty}
\bigl\Vert \hat R_p \bigr\Vert_{L_2([t, T]^{k-2r})}^2=0\ \ \ (k > 5,\  r=2).
$$

After that, we can go to the function 
(\ref{de200}) 
for the case $k > 5,$ $r=3,$ $2r<k$
(this function is defined using the left-hand side of the equality (\ref{drdr1000}))
and follow the same steps as above. This will lead us to the following 
equality
$$
\lim\limits_{p\to\infty}
\bigl\Vert \hat R_p \bigr\Vert_{L_2([t, T]^{k-2r})}^2=0\ \ \ (k > 5,\  r=3,\ 2r<k).
$$

Then we can move on to the next step and so on.
As a result, we get the equality (\ref{pars3s}) ($r=1,2,\ldots,[k/2]$) and thus prove 
Hypothesis~2.2 for the case $k=2n+1,$ $n=3, 4, \ldots $ (see Sect.~2.5). 

For the case $k=2n,$ $n=3, 4, \ldots$ we follow the above steps
for $r=1,2,\ldots,[k/2]-1$ $(2r\le k-2$).
For $2r=k$ we use the same technique as in the proof of the equalities
(\ref{after2508})--(\ref{after2507}). Recall that we used 
(\ref{after80xx}), (\ref{after500}) and 
Parseval's equality in the proof of (\ref{after2508})--(\ref{after2507}).
For $2r=k$ we can also use the equality (\ref{july16000}).

The obvious disadvantage of the proposed algorithm is the drastic 
increase of complexity of the proof when moving from $r=1$ to $r=2,$
$r=2$ to $r=3$ and so on.

The proofs of Theorems~2.34 and 2.35 contain a rather simple trick
of passing from $r=1$ to $r=2.$
Unfortunately, this procedure cannot be 
applied already at the transition from $r=2$ to $r=3.$
Note that the case $k=6,$ $r=3$ was successfully 
considered in Theorem~2.36 under the following 
simplifying assumption:
$\psi_1(\tau),\ldots,\psi_6(\tau)\equiv 1.$

Nevertheless, the results obtained in the previous sections of Chapter~2
are quite sufficient for practical needs (see Chapters~4 and 5 for 
details).

\section{Theorems 2.1--2.9, 2.33--2.36, 2.41, 2.45--2.48, 2.50, 2.51, 2.53, 2.55, 2.57, 2.59, 2.61--2.65
on Expansion of Iterated Stra\-to\-no\-vich
Stochastic Integrals from Point
of View of the Wong\---Zakai Approximation}

The iterated It\^{o} stochastic integrals and solutions
of It\^{o} SDEs are complex and important functionals
from the independent components ${\bf w}_{\tau}^{(i)},$
$i=1,\ldots,m$ of the multidimensional
Wiener process ${\bf w}_{\tau},$ $s\in[0, T].$
Let ${\bf w}_{\tau}^{(i)p},$ $p\in{\bf N}$ 
be some approximation of
${\bf w}_{\tau}^{(i)},$
$i=1,\ldots,m$.
Suppose that 
${\bf w}_{\tau}^{(i)p}$
converges to
${\bf w}_{\tau}^{(i)},$
$i=1,\ldots,m$ if $p\to\infty$ in some sense and has
differentiable sample trajectories.

A natural question arises: if we replace 
${\bf w}_{\tau}^{(i)}$
by ${\bf w}_{\tau}^{(i)p},$
$i=1,\ldots,m$ in the functionals
mentioned above, will the resulting
functionals converge to the original
functionals from the components 
${\bf w}_{\tau}^{(i)},$
$i=1,\ldots,m$ of the multidimentional
Wiener process ${\bf w}_{\tau}$?

The answere to this question is negative 
in the general case. However, 
in the pioneering works of Wong E. and Zakai M. \cite{W-Z-1},
\cite{W-Z-2},
it was shown that under the special conditions and 
for some types of approximations 
of the Wiener process the answere is affirmative
with one peculiarity: the convergence takes place 
to the iterated Stratonovich stochastic integrals
and solutions of Stratonovich SDEs and not to the iterated 
It\^{o} stochastic integrals and solutions
of It\^{o} SDEs.

The piecewise 
linear approximation 
as well as the regularization by convolution 
\cite{W-Z-1}-\cite{Watanabe} relate to the 
mentioned types of approximations
of the Wiener process. The above approximation 
of stochastic integrals and solutions of SDEs 
is often called the Wong--Zakai approximation.

Let ${\bf w}_{\tau},$ $\tau\in[0, T]$ is a random vector with 
an $m+1$ components: ${\bf w}_{\tau}^{(i)}$ 
$(i=1,\ldots,m)$ are independent 
standard Wiener processes and 
${\bf w}_{\tau}^{(0)}=\tau$.

It is well known that the following representation 
takes place \cite{Lipt}, \cite{7e} 
(also see Sect.~6.1 of this book for detail)
\begin{equation}
\label{um1x}
~~~~~~~~~~ {\bf w}_{\tau}^{(i)}-{\bf w}_{t}^{(i)}=
\sum_{j=0}^{\infty}\int\limits_t^{\tau}
\phi_j(s)ds\ \zeta_j^{(i)},\ \ \ \zeta_j^{(i)}=
\int\limits_t^T \phi_j(s)d{\bf w}_s^{(i)},
\end{equation}
where $\tau\in[t, T],$ $t\ge 0,$
$\{\phi_j(x)\}_{j=0}^{\infty}$ is an arbitrary complete 
orthonormal system of functions in the space $L_2([t, T]),$ and
$\zeta_j^{(i)}$ are independent standard Gaussian 
random variables for various $i$ or $j$ (in the case when $i\ne 0$).
Moreover, the series (\ref{um1x}) converges for any $\tau\in [t, T]$
in the mean-square sense $(i=1,\ldots,m)$.

Let ${\bf w}_{\tau}^{(i)p}-{\bf w}_{t}^{(i)p}$ be 
the mean-square approximation of the process
${\bf w}_{\tau}^{(i)}-{\bf w}_{t}^{(i)},$
which has the following form

\vspace{-1mm}
\begin{equation}
\label{um1xx}
{\bf w}_{\tau}^{(i)p}-{\bf w}_{t}^{(i)p}=
\sum_{j=0}^{p}\int\limits_t^{\tau}
\phi_j(s)ds\ \zeta_j^{(i)}.
\end{equation}

\vspace{1mm}

From (\ref{um1xx}) we obtain
\begin{equation}
\label{um1xxx}
d{\bf w}_{\tau}^{(i)p}=
\sum_{j=0}^{p}
\phi_j(\tau)\zeta_j^{(i)} d\tau.
\end{equation}

\vspace{1mm}

Consider the following iterated Riemann--Stieltjes
integral
\begin{equation}
\label{um1xxxx}
\int\limits_t^T
\psi_k(t_k)\ldots \int\limits_t^{t_2}\psi_1(t_1)
d{\bf w}_{t_1}^{(i_1)p_1}\ldots d{\bf w}_{t_k}^{(i_k)p_k},
\end{equation}
where $p_1,\ldots,p_k\in{\bf N},$ $i_1,\ldots,i_k=0,1,\ldots,m,$ 
$d{\bf w}_{\tau}^{(i)p}$ is defined by the relation (\ref{um1xxx}).

\newpage
\noindent
\par
Let us substitute (\ref{um1xxx}) into (\ref{um1xxxx})
$$
\int\limits_t^T
\psi_k(t_k)\ldots \int\limits_t^{t_2}\psi_1(t_1)
d{\bf w}_{t_1}^{(i_1)p_1}\ldots d{\bf w}_{t_k}^{(i_k)p_k}
=
$$
\begin{equation}
\label{um1xxxx1}
=\sum\limits_{j_1=0}^{p_1}\ldots \sum\limits_{j_k=0}^{p_k}
C_{j_k \ldots j_1}\prod\limits_{l=1}^k \zeta_{j_l}^{(i_l)},
\end{equation}

\vspace{1mm}
\noindent
where $p_1,\ldots,p_k\in{\bf N},$
$$
\zeta_j^{(i)}=\int\limits_t^T \phi_j(\tau)d{\bf w}_{\tau}^{(i)}
$$ 
are independent standard Gaussian random variables for various 
$i$ or $j$ (in the case when $i\ne 0$),
${\bf w}_{\tau}^{(i)}$ 
$(i=1,\ldots,m)$ are independent standard Wiener processes,
${\bf w}_{\tau}^{(0)}=\tau,$
$$
C_{j_k \ldots j_1}=\int\limits_t^T\psi_k(t_k)\phi_{j_k}(t_k)\ldots
\int\limits_t^{t_2}
\psi_1(t_1)\phi_{j_1}(t_1)
dt_1\ldots dt_k
$$
is the Fourier coefficient.

To best of our knowledge \cite{W-Z-1}-\cite{Watanabe}
the approximations of the Wiener process
in the Wong--Zakai approximation must satisfy fairly strong
restrictions
\cite{Watanabe}
(see Definition 7.1, pp.~480--481).
Moreover, approximations of the Wiener process that are
similar to (\ref{um1xx})
were not considered in \cite{W-Z-1}, \cite{W-Z-2}
(also see \cite{Watanabe}, Theorems 7.1, 7.2).
Therefore, the proof of analogs of Theorems 7.1 and 7.2 \cite{Watanabe}
for approximations of the Wiener 
process based on its series expansion (\ref{um1x})
(also see (\ref{6.5.7})) should be carried out separately.
Thus, the mean-square convergence of the right-hand side
of (\ref{um1xxxx1}) to the iterated Stratonovich stochastic integral 
(\ref{strxx})
does not follow from the results of the papers
\cite{W-Z-1}, \cite{W-Z-2} (also see \cite{Watanabe},
Theorems 7.1, 7.2) even for the case $p_1=\ldots =p_k=p.$

From the other hand, Theorems 1.1, 1.16, 
2.1--2.9, 2.33--2.36, 2.41, 2.45--2.48, 2.50, 2.51, 2.53, 2.55, 2.57,
2.59, 2.61--2.65
from this 
monograph can be considered as the proof of the
Wong--Zakai approximation based on the iterated 
Riemann--Stieltjes integrals (\ref{um1xxxx}) 
as well as the Wiener process approximation (\ref{um1xx}) 
on the base of its series expansion.
At that, the mentioned Riemann--Stieltjes integrals converge
(according to Theorems 1.1, 1.16, 2.1--2.9, 2.33--2.36, 2.41, 2.45--2.48, 2.50, 2.51, 2.53,
2.55, 2.57, 2.59, 2.61--2.65)
to the appropriate Stratonovich 
stochastic integrals (\ref{strxx}). Recall that
$\{\phi_j(x)\}_{j=0}^{\infty}$ (see (\ref{um1x}), (\ref{um1xx}), and
Theorems 1.1, 2.1, 2.2, 2.4--2.9, 2.33--2.36, 2.41, 2.64, 2.65)
is a complete 
orthonormal system of Legendre polynomials or 
trigonometric functions 
in the space $L_2([t, T])$. In Theorems~1.16, 
2.3, 2.47, 2.48, 2.50, 2.51, 2.53, 2.55, 2.57, 2.59, 2.61--2.63
the system $\{\phi_j(x)\}_{j=0}^{\infty}$ 
can be arbitrary.

To illustrate the above reasoning, 
consider two examples for the case $k=2,$
$\psi_1(s),$ $\psi_2(s)\equiv 1;$ $i_1, i_2=1,\ldots,m.$

The first example relates to the piecewise linear approximation
of the multidimensional Wiener process (these approximations 
were considered in \cite{W-Z-1}-\cite{Watanabe}).

Let ${\bf b}_{\Delta}^{(i)}(t),$ $t\in[0, T]$ be the piecewise
linear approximation of the $i$th component ${\bf w}_t^{(i)}$
of the multidimensional standard Wiener process ${\bf w}_t,$
$t\in [0, T]$ with independent components
${\bf w}_t^{(i)},$ $i=1,\ldots,m,$ i.e.
$$
{\bf b}_{\Delta}^{(i)}(t)={\bf w}_{k\Delta}^{(i)}+
\frac{t-k\Delta}{\Delta}\Delta{\bf w}_{k\Delta}^{(i)},
$$

\noindent
where 
$\Delta{\bf w}_{k\Delta}^{(i)}={\bf w}_{(k+1)\Delta}^{(i)}-
{\bf w}_{k\Delta}^{(i)},$\ \
$t\in[k\Delta, (k+1)\Delta),$\ \ $k=0, 1,\ldots, N-1.$

Note that w.~p.~1
\begin{equation}
\label{pridum}
~~~~~~~~~\frac{d{\bf b}_{\Delta}^{(i)}}{dt}(t)=
\frac{\Delta{\bf w}_{k\Delta}^{(i)}}{\Delta},\ \ \
t\in[k\Delta, (k+1)\Delta),\ \ \ k=0, 1,\ldots, N-1.
\end{equation}

\vspace{2mm}

Consider the following iterated Riemann--Stieltjes
integral
$$
\int\limits_0^T
\int\limits_0^{s}
d{\bf b}_{\Delta}^{(i_1)}(\tau)d{\bf b}_{\Delta}^{(i_2)}(s),\ \ \ 
i_1,i_2=1,\ldots,m.
$$

Using (\ref{pridum}) and additive property of Riemann--Stieltjes integrals, 
we can write w.~p.~1
$$
\int\limits_0^T
\int\limits_0^{s}
d{\bf b}_{\Delta}^{(i_1)}(\tau)d{\bf b}_{\Delta}^{(i_2)}(s)=
\int\limits_0^T
\int\limits_0^{s}
\frac{d{\bf b}_{\Delta}^{(i_1)}}{d\tau}(\tau)d\tau
\frac{d {\bf b}_{\Delta}^{(i_2)}}{d s}(s)
ds =
$$

$$
=
\sum\limits_{l=0}^{N-1}\int\limits_{l\Delta}^{(l+1)\Delta}
\left(
\sum\limits_{q=0}^{l-1}\int\limits_{q\Delta}^{(q+1)\Delta}
\frac{\Delta{\bf w}_{q\Delta}^{(i_1)}}{\Delta}d\tau+
\int\limits_{l\Delta}^{s}
\frac{\Delta{\bf w}_{l\Delta}^{(i_1)}}{\Delta}d\tau\right)
\frac{\Delta{\bf w}_{l\Delta}^{(i_2)}}{\Delta}ds=
$$
$$
=\sum\limits_{l=0}^{N-1}\sum\limits_{q=0}^{l-1}
\Delta{\bf w}_{q\Delta}^{(i_1)}
\Delta{\bf w}_{l\Delta}^{(i_2)}+
\frac{1}{\Delta^2}\sum\limits_{l=0}^{N-1}
\Delta{\bf w}_{l\Delta}^{(i_1)}
\Delta{\bf w}_{l\Delta}^{(i_2)}
\int\limits_{l\Delta}^{(l+1)\Delta}
\int\limits_{l\Delta}^{s}d\tau ds=
$$
\begin{equation}
\label{oh-ty}
=\sum\limits_{l=0}^{N-1}\sum\limits_{q=0}^{l-1}
\Delta{\bf w}_{q\Delta}^{(i_1)}
\Delta{\bf w}_{l\Delta}^{(i_2)}+
\frac{1}{2}\sum\limits_{l=0}^{N-1}
\Delta{\bf w}_{l\Delta}^{(i_1)}
\Delta{\bf w}_{l\Delta}^{(i_2)}.
\end{equation}

\vspace{2mm}

Using (\ref{oh-ty}), it 
is not difficult to show (see Lemma 1.1, Remark 1.2,
and (\ref{oop51})) that
$$
\hbox{\vtop{\offinterlineskip\halign{
\hfil#\hfil\cr
{\rm l.i.m.}\cr
$\stackrel{}{{}_{N\to \infty}}$\cr
}} }
\int\limits_0^T
\int\limits_0^{s}
d{\bf b}_{\Delta}^{(i_1)}(\tau)d{\bf b}_{\Delta}^{(i_2)}(s)=
\int\limits_0^T
\int\limits_0^{s}
d{\bf w}_{\tau}^{(i_1)}d{\bf w}_{s}^{(i_2)}+
$$
\begin{equation}
\label{uh-111}
+
\frac{1}{2}{\bf 1}_{\{i_1=i_2\}}\int\limits_0^T ds=
{\int\limits_0^{*}}^T
{\int\limits_0^{*}}^{s}
d{\bf w}_{\tau}^{(i_1)}d{\bf w}_{s}^{(i_2)},
\end{equation}

\vspace{2mm}
\noindent
where $\Delta\to 0$ if $N\to\infty$ ($N\Delta=T$).
Obviously, (\ref{uh-111}) agrees with Theorem 7.1 (see \cite{Watanabe},
p.~486).

The next example relates to the approximation (\ref{um1xx})
of the Wiener process based on its series expansion (\ref{um1x}), where
$t=0$ and 
$\{\phi_j(x)\}_{j=0}^{\infty}$ 
is an arbitrary complete 
orthonormal system of functions
in the space $L_2([0, T]).$

Consider the following iterated Riemann--Stieltjes
integral
\begin{equation}
\label{abcd1}
\int\limits_0^T
\int\limits_0^{s}
d{\bf w}_{\tau}^{(i_1)p}d{\bf w}_{s}^{(i_2)p},\ \ \ 
i_1,i_2=1,\ldots,m,
\end{equation}

\noindent
where $d{\bf w}_{\tau}^{(i)p}$ is defined by the
relation
(\ref{um1xxx}).

Let us substitute (\ref{um1xxx}) into (\ref{abcd1}) 
\begin{equation}
\label{set18}
\int\limits_0^T
\int\limits_0^{s}
d{\bf w}_{\tau}^{(i_1)p}d{\bf w}_{s}^{(i_2)p}=
\sum\limits_{j_1,j_2=0}^p
C_{j_2 j_1} \zeta_{j_1}^{(i_1)}\zeta_{j_2}^{(i_2)},
\end{equation}

\noindent
where 
$$
C_{j_2 j_1}=
\int\limits_0^T \phi_{j_2}(s)\int\limits_0^s
\phi_{j_1}(\tau)d\tau ds
$$

\noindent
is the Fourier coefficient; another notations 
are the same as in (\ref{um1xxxx1}).

As we noted above, approximations of the Wiener process that are
similar to (\ref{um1xx})
were not considered in \cite{W-Z-1}, \cite{W-Z-2}
(also see Theorems 7.1, 7.2 in \cite{Watanabe}).
Furthermore, the extension of the results of Theorems 7.1 and 7.2
\cite{Watanabe} to the case under consideration is
not obvious.

On the other hand, we can apply the theory built in Chapters 1 and 2
of this book. More precisely, 
using 
Theorem 2.3, we obtain from (\ref{set18}) the desired result
$$
\hbox{\vtop{\offinterlineskip\halign{
\hfil#\hfil\cr
{\rm l.i.m.}\cr
$\stackrel{}{{}_{p\to \infty}}$\cr
}} }
\int\limits_0^T
\int\limits_0^{s}
d{\bf w}_{\tau}^{(i_1)p}d{\bf w}_{s}^{(i_2)p}=
\hbox{\vtop{\offinterlineskip\halign{
\hfil#\hfil\cr
{\rm l.i.m.}\cr
$\stackrel{}{{}_{p\to \infty}}$\cr
}} }
\sum\limits_{j_1,j_2=0}^p
C_{j_2 j_1} \zeta_{j_1}^{(i_1)}\zeta_{j_2}^{(i_2)}=
$$
\begin{equation}
\label{umen-bl}
=
{\int\limits_0^{*}}^T
{\int\limits_0^{*}}^{s}
d{\bf w}_{\tau}^{(i_1)}d{\bf w}_{s}^{(i_2)}.
\end{equation}

\vspace{2mm}

From the other hand, by Theorem 1.16 (see (\ref{razzar1}))
for the case
$k=2$ we obtain from (\ref{set18}) the following relation
$$
\hbox{\vtop{\offinterlineskip\halign{
\hfil#\hfil\cr
{\rm l.i.m.}\cr
$\stackrel{}{{}_{p\to \infty}}$\cr
}} }
\int\limits_0^T
\int\limits_0^{s}
d{\bf w}_{\tau}^{(i_1)p}d{\bf w}_{s}^{(i_2)p}=
\hbox{\vtop{\offinterlineskip\halign{
\hfil#\hfil\cr
{\rm l.i.m.}\cr
$\stackrel{}{{}_{p\to \infty}}$\cr
}} }
\sum\limits_{j_1,j_2=0}^p
C_{j_2 j_1} \zeta_{j_1}^{(i_1)}\zeta_{j_2}^{(i_2)}=
$$
$$
=
\hbox{\vtop{\offinterlineskip\halign{
\hfil#\hfil\cr
{\rm l.i.m.}\cr
$\stackrel{}{{}_{p\to \infty}}$\cr
}} }
\sum\limits_{j_1,j_2=0}^p
C_{j_2 j_1} \left(\zeta_{j_1}^{(i_1)}\zeta_{j_2}^{(i_2)}-
{\bf 1}_{\{i_1=i_2\}}{\bf 1}_{\{j_1=j_2\}}\right)+
{\bf 1}_{\{i_1=i_2\}}\sum\limits_{j_1=0}^{\infty}
C_{j_1 j_1}=
$$
\begin{equation}
\label{umen-blx}
=
\int\limits_0^T
\int\limits_0^{s}
d{\bf w}_{\tau}^{(i_1)}d{\bf w}_{s}^{(i_2)}+
{\bf 1}_{\{i_1=i_2\}}\sum\limits_{j_1=0}^{\infty}
C_{j_1 j_1}.
\end{equation}

\vspace{2mm}

Since
$$
\sum\limits_{j_1=0}^{\infty}
C_{j_1 j_1}=\frac{1}{2}\sum\limits_{j_1=0}^{\infty}
\left(\int\limits_0^T \phi_j(\tau)d\tau\right)^2
=\frac{1}{2}
\int\limits_0^T 1^2 ds=\frac{1}{2}
\int\limits_0^T ds,
$$

\vspace{3mm}
\noindent
then from (\ref{oop51}) and (\ref{umen-blx}) we obtain (\ref{umen-bl}).

\section{Wong--Zakai Type Theorems for Iterated Stra\-to\-no\-vich 
Sto\-chas\-tic 
Integrals. The Case of Approximation of the Multidimensional Wiener Process 
Based on its Series Expansion Using Legendre Polynomials and 
Trigonometric Functions}

As we mentioned above, there exists a lot of publications on the subject 
of Wong--Zakai approximation of stochastic integrals
and SDEs \cite{W-Z-1}-\cite{Watanabe}
(also see \cite{W-Z-3}-\cite{Gyon2}).
However, these works did not consider the approximation of 
iterated stochastic integrals and SDEs
for the case of approximation of the 
multidimensional Wiener process based on its series expansions.
Usually, as an approximation of the Wiener process in the theorems 
of the Wong--Zakai type, the authors \cite{W-Z-1}-\cite{Watanabe}
(also see \cite{W-Z-3}-\cite{Gyon2})
choose a piecewise linear approximation or an 
approximation based on the regularization by convolution.

The Wong--Zakai approximation is widely used to approximate 
stochastic integrals and SDEs. 
In particular, the Wong--Zakai approximation can be used to 
approximate the iterated Stratonovich stochastic integrals 
in the context of numerical integration of It\^{o} SDEs
in the framework of the approach 
based on the Taylor--Stratonovich expansion \cite{Zapad-3}, 
\cite{Zapad-4} (see Chapter 4).
It should be noted that the authors of the works
\cite{Zapad-2} (pp.~438-439),  \cite{Zapad-3}
(Sect.~5.8, pp.~202--204), \cite{Zapad-4} (pp.~82-84),
\cite{Zapad-9} (pp.~263-264) mention 
the Wong--Zakai approximation 
\cite{W-Z-1}-\cite{Watanabe} within the frames
of approximation of iterated Stratonovich stochastic integrals
based on the Karhunen--Loeve expansion of the Brownian bridge
process (see Sect.~6.2). However, in these works there is no rigorous proof 
of convergence for approximations of the mentioned stochastic integrals
of miltiplicity 3 and higher (see discussion in Sect.~6.2).

From the other hand, the theory constructed in Chapters 1 and 2
of this monograph
(also see \cite{12a}-\cite{12aa}) 
can be considered as the proof of the Wong--Zakai
approximation for iterated Stratonovich stochastic integrals.
At that, this approximation is based on the Wiener process series expansion 
using an arbitrary complete orthonormal system of functions in $L_2([t,T]).$

The subject of this section is to reformulate the main results of 
Chapter 2 of this book
in the form of theorems on convergence of iterated
Riemann--Stiltjes integrals 
to iterated Stratonovich stochastic integrals.

Let us reformulate Theorems 2.3--2.6, 2.8, 2.9, 
2.30, 2.32--2.36, 2.41, 2.49--2.51, 2.60--2.65
of this monograph as
statements on the convergence 
of the iterated Rie\-mann--Stiltjes integrals (\ref{um1xxxx})
to the iterated Stratonovich stochastic integrals (\ref{strxx}). 

\vspace{2mm}
                                                                  
{\bf Theorem 2.66}\ \cite{12aa}\ (reformulation of Theorem 2.3).\
{\it Let $\{\phi_j(x)\}_{j=0}^{\infty}$ be an arbitrary complete orthonormal system of 
functions in the space $L_2([t, T])$
and $\psi_1(\tau), \psi_2(\tau)$ are continuous 
functions at the interval $[t, T]$.
Then$,$ for the iterated 
Stratonovich stochastic integral of second multiplicity
$$
J^{*}[\psi^{(2)}]_{T,t}={\int\limits_t^{*}}^T\psi_2(t_2)
{\int\limits_t^{*}}^{t_2}\psi_1(t_1)d{\bf w}_{t_1}^{(i_1)}
d{\bf w}_{t_2}^{(i_2)}\ \ \ (i_1, i_2=1,\ldots,m)
$$
the following formula 
$$
J^{*}[\psi^{(2)}]_{T,t}=
\hbox{\vtop{\offinterlineskip\halign{
\hfil#\hfil\cr
{\rm l.i.m.}\cr
$\stackrel{}{{}_{p_1,p_2\to \infty}}$\cr
}} }\int\limits_t^T
\psi_2(t_2)\int\limits_t^{t_2}\psi_1(t_1)
d{\bf w}_{t_1}^{(i_1)p_1}
d{\bf w}_{t_2}^{(i_2)p_2}
$$
is valid.}

\vspace{2mm}

{\bf Theorem 2.67}\ \cite{arxiv-11}\ 
(reformulation of Theorems 2.4 and 2.6).\
{\it Suppose that
$\{\phi_j(x)\}_{j=0}^{\infty}$ is a complete orthonormal
system of Legendre polynomials or tri\-go\-no\-met\-ric functions
in the space $L_2([t, T])$.
Then$,$ for the iterated Stra\-to\-no\-vich stochastic integral of 
third multiplicity
$$
{\int\limits_t^{*}}^T
{\int\limits_t^{*}}^{t_3}
{\int\limits_t^{*}}^{t_2}
d{\bf w}_{t_1}^{(i_1)}
d{\bf w}_{t_2}^{(i_2)}d{\bf w}_{t_3}^{(i_3)}\ \ \ (i_1, i_2, i_3=1,\ldots,m)
$$
the following 
formula
$$
{\int\limits_t^{*}}^T
{\int\limits_t^{*}}^{t_3}
{\int\limits_t^{*}}^{t_2}
d{\bf w}_{t_1}^{(i_1)}
d{\bf w}_{t_2}^{(i_2)}d{\bf w}_{t_3}^{(i_3)}\ 
=
\hbox{\vtop{\offinterlineskip\halign{
\hfil#\hfil\cr
{\rm l.i.m.}\cr
$\stackrel{}{{}_{p_1,p_2,p_3\to \infty}}$\cr
}} }\int\limits_t^T
\int\limits_t^{t_3}
\int\limits_t^{t_2}
d{\bf w}_{t_1}^{(i_1)p_1}
d{\bf w}_{t_2}^{(i_2)p_2}
d{\bf w}_{t_3}^{(i_3)p_3}
$$

\noindent
is valid.}

\vspace{2mm}

{\bf Theorem 2.68}\ \cite{arxiv-11}\ 
(reformulation of Theorem 2.5).\ {\it Let
$\{\phi_j(x)\}_{j=0}^{\infty}$ be a complete orthonormal
system of Legendre poly\-no\-mi\-als
in the space $L_2([t, T])$.
Then$,$ for the iterated Stratonovich stochastic integral of 
third multiplicity
$$
I_{{l_1l_2l_3}_{T,t}}^{*(i_1i_2i_3)}={\int\limits_t^{*}}^T(t-t_3)^{l_3}
{\int\limits_t^{*}}^{t_3}(t-t_2)^{l_2}
{\int\limits_t^{*}}^{t_2}(t-t_1)^{l_1}
d{\bf w}_{t_1}^{(i_1)}
d{\bf w}_{t_2}^{(i_2)}d{\bf w}_{t_3}^{(i_3)}
$$

\noindent
the following 
formula
$$
I_{{l_1l_2l_3}_{T,t}}^{*(i_1i_2i_3)}=
\hbox{\vtop{\offinterlineskip\halign{
\hfil#\hfil\cr
{\rm l.i.m.}\cr
$\stackrel{}{{}_{p_1,p_2,p_3\to \infty}}$\cr
}} }
\int\limits_t^T
(t-t_3)^{l_3}
\int\limits_t^{t_3}
(t-t_2)^{l_2}
\int\limits_t^{t_2}
(t-t_1)^{l_1}
d{\bf w}_{t_1}^{(i_1)p_1}
d{\bf w}_{t_2}^{(i_2)p_2}
d{\bf w}_{t_3}^{(i_3)p_3},
$$

\noindent
where $i_1, i_2, i_3=1,\ldots,m,$ is valid 
for each of the following cases

\vspace{2mm}
\noindent
{\rm 1}.\ $i_1\ne i_2,\ i_2\ne i_3,\ i_1\ne i_3$\ and
$l_1, l_2, l_3=0, 1, 2,\ldots $\\
{\rm 2}.\ $i_1=i_2\ne i_3$ and $l_1=l_2\ne l_3$\ and
$l_1, l_2, l_3=0, 1, 2,\ldots $\\
{\rm 3}.\ $i_1\ne i_2=i_3$ and $l_1\ne l_2=l_3$\ and
$l_1, l_2, l_3=0, 1, 2,\ldots $\\
{\rm 4}.\ $i_1, i_2, i_3=1,\ldots,m;$ $l_1=l_2=l_3=l$\ and $l=0, 1, 
2,\ldots$
}

\vspace{2mm}

{\bf Theorem 2.69}\ \cite{arxiv-11}\ (reformulation of Theorem 2.8).\
{\it Let
$\{\phi_j(x)\}_{j=0}^{\infty}$ be a complete orthonormal
system of Legendre polynomials or trigonomertic functions
in the space $L_2([t, T])$. Furthermore, let
the function $\psi_2(\tau)$ is continuously
differentiable at the interval $[t, T]$ and
the functions $\psi_1(\tau),$ $\psi_3(\tau)$ are twice continuously
differentiable at the interval $[t, T]$.
Then$,$ for the iterated Stra\-to\-no\-vich stochastic integral 
of third multiplicity

\vspace{-2mm}
$$
J^{*}[\psi^{(3)}]_{T,t}={\int\limits_t^{*}}^T\psi_3(t_3)
{\int\limits_t^{*}}^{t_3}\psi_2(t_2)
{\int\limits_t^{*}}^{t_2}\psi_1(t_1)
d{\bf w}_{t_1}^{(i_1)}
d{\bf w}_{t_2}^{(i_2)}d{\bf w}_{t_3}^{(i_3)}
$$

\vspace{1mm}
\noindent
the following 
formula

\vspace{-5mm}
$$
J^{*}[\psi^{(3)}]_{T,t}=
\hbox{\vtop{\offinterlineskip\halign{
\hfil#\hfil\cr
{\rm l.i.m.}\cr
$\stackrel{}{{}_{p\to \infty}}$\cr
}} }
\int\limits_t^T
\psi_3(t_3)
\int\limits_t^{t_3}
\psi_2(t_2)
\int\limits_t^{t_2}
\psi_1(t_1)
d{\bf w}_{t_1}^{(i_1)p}
d{\bf w}_{t_2}^{(i_2)p}
d{\bf w}_{t_3}^{(i_3)p}
$$

\vspace{1mm}
\noindent
is valid$,$ where $i_1, i_2, i_3=1,\ldots,m.$}  

\vspace{2mm}

{\bf Theorem 2.70}\ \cite{arxiv-11}\ (reformulation of Theorem 2.9).\ 
{\it Let
$\{\phi_j(x)\}_{j=0}^{\infty}$ be a complete orthonormal
system of Legendre polynomials or trigonometric functions
in the space $L_2([t, T])$.
Then$,$ for the iterated Stra\-to\-no\-vich stochastic integral 
of fourth multiplicity

\vspace{-1mm}
$$
I_{T,t}^{*(i_1 i_2 i_3 i_4)}=
{\int\limits_t^{*}}^T
{\int\limits_t^{*}}^{t_4}
{\int\limits_t^{*}}^{t_3}
{\int\limits_t^{*}}^{t_2}
d{\bf w}_{t_1}^{(i_1)}
d{\bf w}_{t_2}^{(i_2)}d{\bf w}_{t_3}^{(i_3)}d{\bf w}_{t_4}^{(i_4)}
$$

\vspace{2mm}
\noindent
the following 
formula

\newpage
\noindent
$$
I_{T,t}^{*(i_1 i_2 i_3 i_4)}=
\hbox{\vtop{\offinterlineskip\halign{
\hfil#\hfil\cr
{\rm l.i.m.}\cr
$\stackrel{}{{}_{p\to \infty}}$\cr
}} }
\int\limits_t^T
\int\limits_t^{t_4}
\int\limits_t^{t_3}
\int\limits_t^{t_2}
d{\bf w}_{t_1}^{(i_1)p}
d{\bf w}_{t_2}^{(i_2)p}
d{\bf w}_{t_3}^{(i_3)p}
d{\bf w}_{t_4}^{(i_4)p}
$$

\vspace{1mm}
\noindent
is valid$,$ where $i_1, i_2, i_3, i_4=0, 1,\ldots,m$.}

\vspace{2mm}

{\bf Theorem 2.71}\ (reformulation of the modified Theorem 2.30 (see Sect. 2.22)).\
{\it Assume that
the continuously differentiable functions 
$\psi_l(\tau)$ $(l=1,\ldots,k)$ at the interval $[t, T]$ and 
the complete orthonormal system $\{\phi_j(x)\}_{j=0}^{\infty}$
of continuous functions  
in the space $L_2([t, T])$ are such that the following 
conditions are satisfied{\rm :}

{\rm 1.}\ The equality 
$$
\frac{1}{2}
\int\limits_t^s \Phi_1(t_1)\Phi_2(t_1)dt_1
=\sum_{j=0}^{\infty}
\int\limits_t^s
\Phi_2(t_2)\phi_{j}(t_2)\int\limits_t^{t_2}
\Phi_1(t_1)\phi_{j}(t_1)dt_1 dt_2
$$

\noindent
holds for all $s\in (t, T],$ where the nonrandom functions 
$\Phi_1(\tau),$ $\Phi_2(\tau)$
are continuously differentiable on $[t, T]$
and the series on the right-hand side of the above equality
converges absolutely.

{\rm 2.}\ The estimates
$$
\left|\int\limits_t^s \phi_{j}(\tau)\Phi_1(\tau)d\tau\right|
\le \frac{\Psi_1(s)}{j^{1/2+\alpha}},\ \ \ 
\left|\int\limits_s^T \phi_{j}(\tau)\Phi_2(\tau)d\tau\right|\le
\frac{\Psi_1(s)}{j^{1/2+\alpha}},
$$
$$
\left|\sum_{j=p+1}^{\infty}\int\limits_t^s
\Phi_2(\tau)\phi_{j}(\tau)\int\limits_t^{\tau}
\Phi_1(\theta)\phi_{j}(\theta)d\theta d\tau\right|\le \frac{\Psi_2(s)}{p^{\hspace{0.5mm}\beta}}
$$

\vspace{2mm}
\noindent
hold for all $s\in (t, T)$ and for some $\alpha, \beta >0,$ where 
$\Phi_1(\tau),$ $\Phi_2(\tau)$
are continuously differentiable nonrandom functions on $[t, T],$\ $j, p\in {\bf N},$
and
$$
\int\limits_t^T \Psi_1^2(\tau) d\tau<\infty,\ \ \ 
\int\limits_t^T \left|\Psi_2(\tau)\right| d\tau<\infty.
$$

\vspace{1mm}

{\rm 3.}\ The condition
$$
\lim\limits_{p\to\infty}
\sum\limits_{\stackrel{j_1,\ldots,j_q,\ldots,j_k=0}{{}_{q\ne g_1, g_2, \ldots, g_{2r-1},
g_{2r}}}}^p
\left(S_{l_1}S_{l_2}\ldots S_{l_{d}}
\left\{\bar C^{(p)}_{j_k\ldots j_q \ldots j_1}\biggl|_{q\ne g_1,g_2,\ldots,g_{2r-1}, g_{2r}}
\right\}\right)^2=0
$$

\vspace{2mm}
\noindent
holds for all possible $g_1,g_2,\ldots,g_{2r-1},g_{2r}$ {\rm (}see {\rm (\ref{leto5007after}))}
and $l_1, l_2, \ldots, l_{d}$ such that
$l_1, l_2, \ldots, l_{d}\in \{1,2,\ldots, r\},$\
$l_1>l_2>\ldots >l_{d},$\ $d=0, 1, 2,\ldots, r-1,$\ 
where $r=1, 2,\ldots,[k/2]$ and
$$
S_{l_1}S_{l_2}\ldots S_{l_{d}}
\left\{\bar C^{(p)}_{j_k\ldots j_q \ldots j_1}\biggl|_{q\ne g_1,g_2,\ldots,g_{2r-1}, g_{2r}}
\right\}\stackrel{\sf def}{=}
\bar C^{(p)}_{j_k\ldots j_q \ldots j_1}\biggl|_{q\ne g_1,g_2,\ldots,g_{2r-1}, g_{2r}}
$$

\noindent
for $d=0.$

Then$,$ for the iterated Stratonovich stochastic integral 
of arbitrary multiplicity $k$
$$
J^{*}[\psi^{(k)}]_{T,t}^{(i_1\ldots i_k)}=
{\int\limits_t^{*}}^T
\psi_k(t_k) \ldots 
{\int\limits_t^{*}}^{t_{2}}
\psi_1(t_1) d{\bf w}_{t_1}^{(i_1)}\ldots
d{\bf w}_{t_k}^{(i_k)}
$$

\vspace{1mm}
\noindent
the following 
formula
$$
J^{*}[\psi^{(k)}]_{T,t}^{(i_1\ldots i_k)}=
\hbox{\vtop{\offinterlineskip\halign{
\hfil#\hfil\cr
{\rm l.i.m.}\cr
$\stackrel{}{{}_{p\to \infty}}$\cr
}} }
\int\limits_t^{T}\psi_k(t_k) \ldots 
\int\limits_t^{t_{2}}
\psi_1(t_1) d{\bf w}_{t_1}^{(i_1)p}\ldots
d{\bf w}_{t_k}^{(i_k)p}
$$

\vspace{1mm}
\noindent
is valid$,$ where $i_1,\ldots, i_k=0, 1,\ldots,m$.}

\vspace{2mm}

{\bf Theorem 2.72}\ (reformulation of the modified Theorem 2.32 (see Sect. 2.22)).\
{\it Assume that
the continuous functions 
$\psi_1(\tau),\ldots,\psi_k(\tau)$ at the interval $[t, T]$ and 
the complete orthonormal system $\{\phi_j(x)\}_{j=0}^{\infty}$
of functions
in the space $L_2([t, T])$ are such that the following 
condition 

\vspace{-2mm}
$$
\lim\limits_{p_1,\ldots,p_k\to\infty}~
\sum\limits_{j_1=0}^{p_1}\ldots \sum\limits_{j_q=0}^{p_q}\ldots \sum\limits_{j_k=0}^{p_k}~
\biggl|_{q\ne g_1, g_2, \ldots, g_{2r-1},g_{2r}}\times
$$

$$
\times
\Biggl(~\sum\limits_{j_{g_1}=0}^{\min\{p_{g_1}, p_{g_2}\}} \sum\limits_{j_{g_3}=0}^{\min\{p_{g_3}, p_{g_4}\}}\ldots \Biggr.
\sum\limits_{j_{g_{2r-1}}=0}^{\min\{p_{g_{2r-1}}, p_{g_{2r}}\}}
C_{j_k\ldots j_1}\biggl|_{j_{g_1}=j_{g_2},\ldots, j_{g_{2r-1}}=j_{g_{2r}}}-
$$

\vspace{-2mm}
$$
\Biggl.-\frac{1}{2^r} \prod\limits_{l=1}^r {\bf 1}_{\{g_{2l}=g_{2l-1}+1\}}
C_{j_k \ldots j_1}\biggl|_{(j_{g_2} j_{g_1})\curvearrowright (\cdot)
\ldots (j_{g_{2r}} j_{g_{2r-1}})\curvearrowright (\cdot),
j_{g_{{}_{1}}}=~j_{g_{{}_{2}}},\ldots, j_{g_{{}_{2r-1}}}=~j_{g_{{}_{2r}}}
}\biggr.\Biggr)^2=0
$$

\vspace{2mm}
\noindent
is satisfied for all $r=1, 2,\ldots,[k/2]$
and for all possible $g_1,g_2,\ldots,g_{2r-1},g_{2r}$ {\rm (}see {\rm (\ref{leto5007after}))}.
Then$,$ for the iterated Stratonovich stochastic integral 
of arbitrary multiplicity $k$

\newpage
\noindent
$$
J^{*}[\psi^{(k)}]_{T,t}^{(i_1\ldots i_k)}=
{\int\limits_t^{*}}^T
\psi_k(t_k) \ldots 
{\int\limits_t^{*}}^{t_{2}}
\psi_1(t_1) d{\bf w}_{t_1}^{(i_1)}\ldots
d{\bf w}_{t_k}^{(i_k)}
$$
the following 
formula
$$
J^{*}[\psi^{(k)}]_{T,t}^{(i_1\ldots i_k)}=
\hbox{\vtop{\offinterlineskip\halign{
\hfil#\hfil\cr
{\rm l.i.m.}\cr
$\stackrel{}{{}_{p_1,\ldots,p_k\to \infty}}$\cr
}} }
\int\limits_t^{T}\psi_k(t_k) \ldots 
\int\limits_t^{t_{2}}
\psi_1(t_1) d{\bf w}_{t_1}^{(i_1)p_1}\ldots
d{\bf w}_{t_k}^{(i_k)p_k}
$$
is valid$,$ where $i_1,\ldots, i_k=0, 1,\ldots,m$.}

\vspace{2mm}

{\bf Theorem 2.73}\ (reformulation of Theorem 2.33).\
{\it Suppose that 
$\{\phi_j(x)\}_{j=0}^{\infty}$ is a complete orthonormal system of 
Legendre polynomials or trigonometric functions in the space $L_2([t, T]).$
Furthermore$,$ let $\psi_1(s), \psi_2(s), \psi_3(s)$ are continuously dif\-ferentiable 
nonrandom functions on $[t, T].$ 
Then$,$ for the 
iterated Stratonovich stochastic integral of third multiplicity
$$
J^{*}[\psi^{(3)}]_{T,t}={\int\limits_t^{*}}^T\psi_3(t_3)
{\int\limits_t^{*}}^{t_3}\psi_2(t_2)
{\int\limits_t^{*}}^{t_2}\psi_1(t_1)
d{\bf w}_{t_1}^{(i_1)}
d{\bf w}_{t_2}^{(i_2)}d{\bf w}_{t_3}^{(i_3)}
$$
the following 
formula
$$
J^{*}[\psi^{(3)}]_{T,t}
=\hbox{\vtop{\offinterlineskip\halign{
\hfil#\hfil\cr
{\rm l.i.m.}\cr
$\stackrel{}{{}_{p\to \infty}}$\cr
}} }
\int\limits_t^T\psi_3(t_3)
\int\limits_t^{t_3}\psi_2(t_2)
\int\limits_t^{t_2}\psi_1(t_1)
d{\bf w}_{t_1}^{(i_1)p}
d{\bf w}_{t_2}^{(i_2)p}d{\bf w}_{t_3}^{(i_3)p}
$$
is valid$,$ where $i_1, i_3, i_3=0, 1,\ldots,m$.
}

\vspace{2mm}

{\bf Theorem 2.74}\ (reformulation of Theorem 2.34).\
{\it Suppose that 
$\{\phi_j(x)\}_{j=0}^{\infty}$ is a complete orthonormal system of 
Legendre polynomials or trigonometric functions in the space $L_2([t, T]).$
Furthermore$,$ let $\psi_1(s), \ldots, \psi_4(s)$ are continuously dif\-ferentiable 
nonrandom functions on $[t, T].$ 
Then$,$ for the 
iterated Stratonovich stochastic integral of fourth multiplicity
$$
J^{*}[\psi^{(4)}]_{T,t}={\int\limits_t^{*}}^T\psi_4(t_4)
{\int\limits_t^{*}}^{t_4}\psi_3(t_3)
{\int\limits_t^{*}}^{t_3}\psi_2(t_2)
{\int\limits_t^{*}}^{t_2}\psi_1(t_1)
d{\bf w}_{t_1}^{(i_1)}
d{\bf w}_{t_2}^{(i_2)}d{\bf w}_{t_3}^{(i_3)}d{\bf w}_{t_4}^{(i_4)}
$$
the following 
formula
$$
J^{*}[\psi^{(4)}]_{T,t}
=\hbox{\vtop{\offinterlineskip\halign{
\hfil#\hfil\cr
{\rm l.i.m.}\cr
$\stackrel{}{{}_{p\to \infty}}$\cr
}} }
\int\limits_t^T\psi_4(t_4)
\int\limits_t^{t_4}\psi_3(t_3)
\int\limits_t^{t_3}\psi_2(t_2)
\int\limits_t^{t_2}\psi_1(t_1)
d{\bf w}_{t_1}^{(i_1)p}\times
$$
$$
\times
d{\bf w}_{t_2}^{(i_2)p}d{\bf w}_{t_3}^{(i_3)p}d{\bf w}_{t_4}^{(i_4)p}
$$

\vspace{2mm}
\noindent
is valid$,$ where $i_1, \ldots, i_4=0, 1,\ldots,m$.
}

\vspace{2mm}

{\bf Theorem 2.75}\ (reformulation of Theorem 2.35).\
{\it Suppose that 
$\{\phi_j(x)\}_{j=0}^{\infty}$ is a complete orthonormal system of 
Legendre polynomials or trigonometric functions in the space $L_2([t, T]).$
Furthermore$,$ let $\psi_1(s), \ldots, \psi_5(s)$ are continuously dif\-ferentiable 
nonrandom functions on $[t, T].$ 
Then$,$ for the 
iterated Stra\-to\-no\-vich stochastic integral of fifth multiplicity
$$
J^{*}[\psi^{(5)}]_{T,t}=
{\int\limits_t^{*}}^T\psi_5(t_5)
{\int\limits_t^{*}}^{t_5}\psi_4(t_4)
{\int\limits_t^{*}}^{t_4}\psi_3(t_3)
{\int\limits_t^{*}}^{t_3}\psi_2(t_2)
{\int\limits_t^{*}}^{t_2}\psi_1(t_1)
d{\bf w}_{t_1}^{(i_1)}\times
$$
$$
\times
d{\bf w}_{t_2}^{(i_2)}d{\bf w}_{t_3}^{(i_3)}d{\bf w}_{t_4}^{(i_4)}d{\bf w}_{t_5}^{(i_5)}
$$

\vspace{4mm}
\noindent
the following 
formula
$$
J^{*}[\psi^{(5)}]_{T,t}
=\hbox{\vtop{\offinterlineskip\halign{
\hfil#\hfil\cr
{\rm l.i.m.}\cr
$\stackrel{}{{}_{p\to \infty}}$\cr
}} }
\int\limits_t^T\psi_5(t_5)
\int\limits_t^{t_5}\psi_4(t_4)
\int\limits_t^{t_4}\psi_3(t_3)
\int\limits_t^{t_3}\psi_2(t_2)
\int\limits_t^{t_2}\psi_1(t_1)
d{\bf w}_{t_1}^{(i_1)p}\times
$$

\vspace{-2mm}
$$
\times
d{\bf w}_{t_2}^{(i_2)p}d{\bf w}_{t_3}^{(i_3)p}d{\bf w}_{t_4}^{(i_4)p}
d{\bf w}_{t_5}^{(i_5)p}
$$

\vspace{2mm}
\noindent
is valid$,$ where $i_1, \ldots, i_5=0, 1,\ldots,m$.
}

\vspace{2mm}

{\bf Theorem 2.76}\ (reformulation of Theorems 2.36, 2.64, 2.65).\
{\it Suppose that 
$\{\phi_j(x)\}_{j=0}^{\infty}$ is a complete orthonormal system of 
Legendre polynomials or trigonometric functions in the space $L_2([t, T]).$
Then$,$ for the 
iterated Stra\-to\-no\-vich stochastic integral

\vspace{-4mm}
$$
J_{T,t}^{*(i_1\ldots i_k)}={\int\limits_t^{*}}^T
{\int\limits_t^{*}}^{t_k}
\ldots
{\int\limits_t^{*}}^{t_2}
d{\bf w}_{t_1}^{(i_1)}
\ldots
d{\bf w}_{t_{k-1}}^{(i_{k-1})}
d{\bf w}_{t_k}^{(i_k)}\ \ \ (k=6, 7, 8)
$$

\noindent
the following 
formula
$$
J_{T,t}^{*(i_1\ldots i_k)}
=\hbox{\vtop{\offinterlineskip\halign{
\hfil#\hfil\cr
{\rm l.i.m.}\cr
$\stackrel{}{{}_{p\to \infty}}$\cr
}} }
\int\limits_t^T
\int\limits_t^{t_k}
\ldots
\int\limits_t^{t_2}
d{\bf w}_{t_1}^{(i_1)p}
\ldots
d{\bf w}_{t_{k-1}}^{(i_{k-1})p}d{\bf w}_{t_k}^{(i_k)p}\ \ \ (k=6, 7, 8)
$$

\noindent
is valid$,$ where $i_1, \ldots, i_8=0, 1,\ldots,m$.
}

\vspace{2mm}

{\bf Theorem~2.77}\ (reformulation of Theorem 2.41).\ {\it Suppose that 
$\{\phi_j(x)\}_{j=0}^{\infty}$ is a complete orthonormal system of 
Legendre polynomials or trigonometric functions in the space $L_2([t, T]).$
Furthermore, let $\psi_1(s), \psi_2(s), \psi_3(s)$ are continuously dif\-ferentiable 
nonrandom functions on $[t, T]$.
Then$,$ for the 
iterated Stratonovich stochastic integral of third multiplicity
$$
J^{*}[\psi^{(3)}]_{T,t}={\int\limits_t^{*}}^T\psi_3(t_3)
{\int\limits_t^{*}}^{t_3}\psi_2(t_2)
{\int\limits_t^{*}}^{t_2}\psi_1(t_1)
d{\bf w}_{t_1}^{(i_1)}
d{\bf w}_{t_2}^{(i_2)}d{\bf w}_{t_3}^{(i_3)}
$$
the following 
formula
$$
J^{*}[\psi^{(3)}]_{T,t}
=\hbox{\vtop{\offinterlineskip\halign{
\hfil#\hfil\cr
{\rm l.i.m.}\cr
$\stackrel{}{{}_{p_1,p_2,p_3\to \infty}}$\cr
}} }
\int\limits_t^T\psi_3(t_3)
\int\limits_t^{t_3}\psi_2(t_2)
\int\limits_t^{t_2}\psi_1(t_1)
d{\bf w}_{t_1}^{(i_1)p_1}
d{\bf w}_{t_2}^{(i_2)p_2}d{\bf w}_{t_3}^{(i_3)p_3}
$$
is valid$,$ where $i_1, i_2, i_3=0, 1,\ldots,m$.}

\vspace{2mm}

{\bf Theorem~2.78}\ (reformulation of Theorem 2.49).\ {\it Assume that
the complete orthonormal system $\{\phi_j(x)\}_{j=0}^{\infty}$
in the space $L_2([t, T])$ and
$\psi_1(\tau),\ldots, \psi_k(\tau)\in L_2([t, T])$
are such that 
$$
\lim\limits_{p_1,\ldots,p_k\to\infty}~
\sum\limits_{j_1=0}^{p_1}\ldots \sum\limits_{j_q=0}^{p_q}\ldots \sum\limits_{j_k=0}^{p_k}~
\biggl|_{q\ne g_1, g_2, \ldots, g_{2r-1},g_{2r}}\times
$$

$$
\times
\Biggl(~\sum\limits_{j_{g_1}=0}^{\min\{p_{g_1}, p_{g_2}\}} \sum\limits_{j_{g_3}=0}^{\min\{p_{g_3}, p_{g_4}\}}\ldots \Biggr.
\sum\limits_{j_{g_{2r-1}}=0}^{\min\{p_{g_{2r-1}}, p_{g_{2r}}\}}
C_{j_k\ldots j_1}\biggl|_{j_{g_1}=j_{g_2},\ldots, j_{g_{2r-1}}=j_{g_{2r}}}-
$$

$$
\Biggl.-\frac{1}{2^r} \prod\limits_{l=1}^r {\bf 1}_{\{g_{2l}=g_{2l-1}+1\}}
C_{j_k \ldots j_1}\biggl|_{(j_{g_2} j_{g_1})\curvearrowright (\cdot)
\ldots (j_{g_{2r}} j_{g_{2r-1}})\curvearrowright (\cdot),
j_{g_{{}_{1}}}=~j_{g_{{}_{2}}},\ldots, j_{g_{{}_{2r-1}}}=~j_{g_{{}_{2r}}}
}\biggr.\Biggr)^2=0
$$

\vspace{3mm}
\noindent
for all $r=1, 2,\ldots,[k/2]$
and for all possible $g_1,g_2,\ldots,g_{2r-1},g_{2r}$ {\rm (}see {\rm (\ref{leto5007after}))}.
Then$,$ for the sum 
of iterated It\^{o} sto\-chas\-tic integrals 
of the form
\begin{equation}
\label{july200000}
~~~~~~~~~J[\psi^{(k)}]_{T,t}^{(i_1\ldots i_k)}+
\sum_{r=1}^{\left[k/2\right]}\frac{1}{2^r}
\sum_{(s_r,\ldots,s_1)\in {\rm A}_{k,r}}
J[\psi^{(k)}]_{T,t}^{s_r,\ldots,s_1}\stackrel{\sf def}{=}\bar J^{*}[\psi^{(k)}]_{T,t}^{(i_1\ldots i_k)}
\end{equation}

\noindent
the following 
formula 

\newpage
\noindent
$$
\bar J^{*}[\psi^{(k)}]_{T,t}^{(i_1\ldots i_k)}=
\hbox{\vtop{\offinterlineskip\halign{
\hfil#\hfil\cr
{\rm l.i.m.}\cr
$\stackrel{}{{}_{p_1,\ldots,p_k\to \infty}}$\cr
}} }
\int\limits_t^{T}\psi_k(t_k) \ldots 
\int\limits_t^{t_{2}}
\psi_1(t_1) d{\bf w}_{t_1}^{(i_1)p_1}\ldots
d{\bf w}_{t_k}^{(i_k)p_k}
$$

\vspace{1mm}
\noindent
is valid$,$ where $i_1, \ldots, i_k=0, 1,\ldots,m$.}

\vspace{2mm}

{\bf Theorem 2.79}\ (reformulation of Theorem 2.51).\ {\it Suppose that
$\{\phi_j(x)\}_{j=0}^{\infty}$ is an arbitrary complete orthonormal system of 
functions in the space $L_2([t,T]).$
Then$,$ for the iterated Stra\-to\-no\-vich stochastic integral
of third multiplicity 
$$
I_{{l_1l_2l_3}_{T,t}}^{*(i_1i_2i_3)}={\int\limits_t^{*}}^T (t_3-t)^{l_3}
{\int\limits_t^{*}}^{t_3}(t_2-t)^{l_2}
{\int\limits_t^{*}}^{t_2}(t_1-t)^{l_1}
d{\bf w}_{t_1}^{(i_1)}
d{\bf w}_{t_2}^{(i_2)}d{\bf w}_{t_3}^{(i_3)}
$$

\noindent
the following formula
$$
I_{{l_1l_2l_3}_{T,t}}^{*(i_1i_2i_3)}=
\hbox{\vtop{\offinterlineskip\halign{
\hfil#\hfil\cr
{\rm l.i.m.}\cr
$\stackrel{}{{}_{p\to \infty}}$\cr
}} }
\int\limits_t^T (t_3-t)^{l_3}
\int\limits_t^{t_3} (t_2-t)^{l_2}
\int\limits_t^{t_2} (t_1-t)^{l_1}
d{\bf w}_{t_1}^{(i_1)p}
d{\bf w}_{t_2}^{(i_2)p}d{\bf w}_{t_3}^{(i_3)p}
$$

\noindent
is valid, where 
$i_1,i_2,i_3=0,1,\ldots,m;$ $l_1,l_2,l_3=0,1,2,\ldots$
}

\vspace{2mm}

{\bf Theorem~2.80}\ (reformulation of Theorem 2.63).\ {\it Suppose that
$\{\phi_j(x)\}_{j=0}^{\infty}$ is an arbitrary complete orthonormal system of 
functions in the space $L_2([t,T]).$
Then$,$ for the iterated Stra\-to\-no\-vich stochastic integral
of fourth multiplicity 
$$
I_{{l_1l_2l_3 l_4}_{T,t}}^{*(i_1i_2i_3 i_4)}=
{\int\limits_t^{*}}^T (t_4-t)^{l_4}{\int\limits_t^{*}}^{t_4} (t_3-t)^{l_3}
{\int\limits_t^{*}}^{t_3}(t_2-t)^{l_2}
{\int\limits_t^{*}}^{t_2}(t_1-t)^{l_1}\times
$$
$$
\times
d{\bf w}_{t_1}^{(i_1)}
d{\bf w}_{t_2}^{(i_2)}d{\bf w}_{t_3}^{(i_3)}d{\bf w}_{t_4}^{(i_4)}
$$

\noindent
the following formula
$$
I_{{l_1l_2l_3 l_4}_{T,t}}^{*(i_1i_2i_3i_4)}=
\hbox{\vtop{\offinterlineskip\halign{
\hfil#\hfil\cr
{\rm l.i.m.}\cr
$\stackrel{}{{}_{p\to \infty}}$\cr
}} }
\int\limits_t^T (t_4-t)^{l_4}
\int\limits_t^{t_4} (t_3-t)^{l_3}
\int\limits_t^{t_3} (t_2-t)^{l_2}
\int\limits_t^{t_2} (t_1-t)^{l_1}\times
$$
$$
\times
d{\bf w}_{t_1}^{(i_1)p}
d{\bf w}_{t_2}^{(i_2)p}d{\bf w}_{t_3}^{(i_3)p}d{\bf w}_{t_4}^{(i_4)p}
$$

\vspace{2mm}
\noindent
is valid, where 
$i_1,\ldots,i_4=0,1,\ldots,m;$ $l_1,\ldots,l_4=0,1,2,\ldots$
}

\vspace{2mm}

{\bf Theorem~2.81}\ (reformulation of Theorem 2.50).\ {\it Suppose that
$\{\phi_j(x)\}_{j=0}^{\infty}$ is an arbitrary complete orthonormal system of 
functions in the space $L_2([t,T]).$
Then$,$ for the iterated Stra\-to\-no\-vich stochastic integral
of fifth multiplicity 
$$
J_{T,t}^{*(i_1\ldots i_5)}={\int\limits_t^{*}}^T
{\int\limits_t^{*}}^{t_5}
{\int\limits_t^{*}}^{t_4}
{\int\limits_t^{*}}^{t_3}
{\int\limits_t^{*}}^{t_2}
d{\bf w}_{t_1}^{(i_1)}
d{\bf w}_{t_2}^{(i_2)}
d{\bf w}_{t_3}^{(i_3)}
d{\bf w}_{t_4}^{(i_4)}
d{\bf w}_{t_5}^{(i_5)}
$$

\noindent
the following 
formula
$$
J_{T,t}^{*(i_1\ldots i_5)}
=\hbox{\vtop{\offinterlineskip\halign{
\hfil#\hfil\cr
{\rm l.i.m.}\cr
$\stackrel{}{{}_{p\to \infty}}$\cr
}} }
\int\limits_t^T
\int\limits_t^{t_5}
\int\limits_t^{t_4}
\int\limits_t^{t_3}
\int\limits_t^{t_2}
d{\bf w}_{t_1}^{(i_1)p}
d{\bf w}_{t_2}^{(i_2)p}d{\bf w}_{t_3}^{(i_3)p}d{\bf w}_{t_4}^{(i_4)p}
d{\bf w}_{t_5}^{(i_5)p}
$$

\noindent
is valid$,$ where $i_1, \ldots, i_5=0, 1,\ldots,m$.
}

\vspace{2mm}                     

{\bf Theorem~2.82}\ (reformulation of Theorem 2.60).\ 
{\it Suppose that 
the condition {\rm (\ref{09091})} 
is fulfilled$,$
$\{\phi_j(x)\}_{j=0}^{\infty}$
is an arbitrary complete orthonormal system of functions
in the space $L_2([t, T])$ and
$\psi_1(\tau),\ldots, \psi_k(\tau)\in L_2([t, T]).$
Then$,$ for the sum $\bar J^{*}[\psi^{(k)}]_{T,t}^{(i_1\ldots i_k)}$
of iterated It\^{o} stochastic integrals
defined by {\rm (\ref{july200000}),}
we have
$$
\bar J^{*}[\psi^{(k)}]_{T,t}^{(i_1\ldots i_k)}=
\hbox{\vtop{\offinterlineskip\halign{
\hfil#\hfil\cr
{\rm l.i.m.}\cr
$\stackrel{}{{}_{p\to \infty}}$\cr
}} }
\int\limits_t^T
\psi_k(t_k)
\ldots
\int\limits_t^{t_2}
\psi_1(t_1)
d{\bf w}_{t_1}^{(i_1)p}
\ldots
d{\bf w}_{t_k}^{(i_k)p},
$$

\noindent
where $i_1,\ldots, i_k=0, 1, \ldots,m$.}

\vspace{2mm}

{\bf Theorem~2.83}\ (reformulation of Theorem 2.61).\
{\it Suppose that 
the condition {\rm (\ref{09091})} 
is fulfilled$,$
$\{\phi_j(x)\}_{j=0}^{\infty}$
is an arbitrary complete orthonormal system of functions
in the space $L_2([t, T])$ and
$\psi_1(\tau),\ldots, \psi_k(\tau)$ are continuous functions
at the interval $[t, T].$
Then$,$ for the iterated Stratonovich sto\-chas\-tic integral 
of arbitrary multiplicity $k$
$$
J^{*}[\psi^{(k)}]_{T,t}^{(i_1\ldots i_k)}=
{\int\limits_t^{*}}^T
\psi_k(t_k) \ldots 
{\int\limits_t^{*}}^{t_{2}}
\psi_1(t_1) d{\bf w}_{t_1}^{(i_1)}\ldots
d{\bf w}_{t_k}^{(i_k)}
$$

\noindent
the following 
formula
$$
J^{*}[\psi^{(k)}]_{T,t}^{(i_1\ldots i_k)}=
\hbox{\vtop{\offinterlineskip\halign{
\hfil#\hfil\cr
{\rm l.i.m.}\cr
$\stackrel{}{{}_{p\to \infty}}$\cr
}} }
\int\limits_t^T
\psi_k(t_k)
\ldots
\int\limits_t^{t_2}
\psi_1(t_1)
d{\bf w}_{t_1}^{(i_1)p}
\ldots
d{\bf w}_{t_k}^{(i_k)p}
$$

\noindent
is valid, where $i_1,\ldots, i_k=0, 1, \ldots,m$.}

{\bf Theorem~2.84}\ (reformulation of Theorem 2.62).\ {\it Suppose that
$\{\phi_j(x)\}_{j=0}^{\infty}$ is an arbitrary complete orthonormal system of 
functions in the space $L_2([t,T]).$
Then$,$ for the iterated Stra\-to\-no\-vich stochastic integral
of sixth multiplicity 
$$
J_{T,t}^{*(i_1\ldots i_6)}={\int\limits_t^{*}}^T
{\int\limits_t^{*}}^{t_6}
{\int\limits_t^{*}}^{t_5}
{\int\limits_t^{*}}^{t_4}
{\int\limits_t^{*}}^{t_3}
{\int\limits_t^{*}}^{t_2}
d{\bf w}_{t_1}^{(i_1)}
d{\bf w}_{t_2}^{(i_2)}
d{\bf w}_{t_3}^{(i_3)}
d{\bf w}_{t_4}^{(i_4)}
d{\bf w}_{t_5}^{(i_5)}
d{\bf w}_{t_6}^{(i_6)}
$$

\noindent
the following 
formula
$$
J_{T,t}^{*(i_1\ldots i_6)}
=\hbox{\vtop{\offinterlineskip\halign{
\hfil#\hfil\cr
{\rm l.i.m.}\cr
$\stackrel{}{{}_{p\to \infty}}$\cr
}} }
\int\limits_t^T
\int\limits_t^{t_6}
\int\limits_t^{t_5}
\int\limits_t^{t_4}
\int\limits_t^{t_3}
\int\limits_t^{t_2}
d{\bf w}_{t_1}^{(i_1)p}
d{\bf w}_{t_2}^{(i_2)p}d{\bf w}_{t_3}^{(i_3)p}d{\bf w}_{t_4}^{(i_4)p}
d{\bf w}_{t_5}^{(i_5)p}d{\bf w}_{t_6}^{(i_6)p}
$$

\noindent
is valid$,$ where $i_1, \ldots, i_6=0, 1,\ldots,m$.
}

\section{Expansion of Iterated Stratonovich Stochastic Integrals
of Multiplicity $k.$ The Case $i_1=\ldots=i_k\ne 0$ and
Different Continuously Differentiable 
Weight Functions $\psi_1(\tau),\ldots,\psi_k(\tau)$}

This section was written several years earlier
than Sect.~2.30--2.39. The results of Sect.~2.30--2.39
somewhat depreciate the results of Sect.~2.44, but we still
included it in this version of the book.

In this section, we generalize the approach considered in Sect.~2.1.2
to the case $i_1=\ldots=i_k\ne 0$ and 
different weight functions $\psi_1(\tau),\ldots,\psi_k(\tau)$ $(k>2).$
Let us formulate the following theorem.

\vspace{2mm}

{\bf Theorem 2.85} \cite{arxiv-6}. {\it Suppose that 
$\{\phi_j(x)\}_{j=0}^{\infty}$ is a complete orthonormal system of 
Legendre polynomials or trigonometric functions in the space $L_2([t, T]).$
Moreover, $\psi_1(\tau),\ldots,\psi_k(\tau)$ $(k\ge 2)$ are 
continuously differentiable nonrandom functions on $[t, T]$. Then$,$ 
for the iterated Stratonovich stochastic integral
$$
J^{*}[\psi^{(k)}]_{T,t}={\int\limits_t^{*}}^T\psi_k(t_k)\ldots
{\int\limits_t^{*}}^{t_2}\psi_1(t_1)d{\bf w}_{t_1}^{(i_1)}\ldots 
d{\bf w}_{t_k}^{(i_1)}\ \ \ (i_1=1,\ldots,m)
$$

\vspace{1mm}
\noindent
the following equality
$$
\lim\limits_{p\to\infty}
{\sf M}\left\{\left(
J^{*}[\psi^{(k)}]_{T,t}-\sum_{j_1=0}^{p}\ldots \sum_{j_k=0}^{p}
C_{j_k \ldots j_1}\zeta_{j_1}^{(i_1)}\ldots \zeta_{j_k}^{(i_1)}\right)^{2n}\right\}=0
$$

\vspace{1mm}
\noindent
is valid, where $n\in{\bf N},$
$$
C_{j_k \ldots j_1}=\int\limits_t^T\psi_k(t_k)\phi_{j_k}(t_k)
\ldots \int\limits_t^{t_2}\psi_1(t_1)\phi_{j_1}(t_1)dt_1\ldots dt_k
$$

\vspace{1mm}
\noindent
is the Fourier coefficient and
$$
\zeta_{j}^{(i_1)}=
\int\limits_t^T \phi_{j}(\tau) d{\bf w}_{\tau}^{(i_1)}\ \ \ (i_1=1,\ldots,m)
$$ 

\vspace{1mm}
\noindent
are independent
standard Gaussian random variables for various 
$j$.}

\vspace{2mm}

{\bf Proof.}\ The case $k=2$ is proved in Theorem 2.16.
Consider the case $k>2$. First, consider the case $k=3$ in detail.
Define the auxiliary function
$$
K'(t_1,t_2,t_3)=\frac{1}{6}\left\{
\begin{matrix}
\psi_1(t_1)\psi_2(t_2)\psi_3(t_3),
\ \ t_1\le t_2 \le t_3
\cr\cr
\psi_1(t_1)\psi_2(t_3)\psi_3(t_2),
\ \ t_1\le t_3 \le t_2
\cr\cr
\psi_1(t_2)\psi_2(t_1)\psi_3(t_3),
\ \ t_2\le t_1 \le t_3
\cr\cr
\psi_1(t_2)\psi_2(t_3)\psi_3(t_1),
\ \ t_2\le t_3 \le t_1
\cr\cr
\psi_1(t_3)\psi_2(t_2)\psi_3(t_1),
\ \ t_3\le t_2 \le t_1
\cr\cr
\psi_1(t_3)\psi_2(t_1)\psi_3(t_2),
\ \ t_3\le t_1 \le t_2
\end{matrix}\right.,\ \ \ \ t_1,t_2,t_3\in[t,T].
$$

\vspace{4mm}

Using Lemma~1.1, Remark~1.1 (see Sect.~1.1.3), and (\ref{uyes3}), we obtain w.~p.~1
$$
J[K']_{T,t}^{(3)}=
\hbox{\vtop{\offinterlineskip\halign{
\hfil#\hfil\cr
{\rm l.i.m.}\cr
$\stackrel{}{{}_{N\to \infty}}$\cr
}} }\sum_{l_3=0}^{N-1}\sum_{l_2=0}^{N-1}
\sum_{l_1=0}^{N-1}
K'(\tau_{l_1},\tau_{l_2},\tau_{l_3})
\Delta{\bf w}_{\tau_{l_1}}^{(i_1)}
\Delta{\bf w}_{\tau_{l_2}}^{(i_1)}
\Delta{\bf w}_{\tau_{l_3}}^{(i_1)}=
$$
$$
=\hbox{\vtop{\offinterlineskip\halign{
\hfil#\hfil\cr
{\rm l.i.m.}\cr
$\stackrel{}{{}_{N\to \infty}}$\cr
}} }\left(\sum_{l_3=0}^{N-1}\sum_{l_2=0}^{l_3-1}
\sum_{l_1=0}^{l_2-1}
K'(\tau_{l_1},\tau_{l_2},\tau_{l_3})
\Delta{\bf w}_{\tau_{l_1}}^{(i_1)}
\Delta{\bf w}_{\tau_{l_2}}^{(i_1)}
\Delta{\bf w}_{\tau_{l_3}}^{(i_1)}+\right.
$$
$$
+\sum_{l_3=0}^{N-1}\sum_{l_1=0}^{l_3-1}
\sum_{l_2=0}^{l_1-1}
K'(\tau_{l_1},\tau_{l_2},\tau_{l_3})
\Delta{\bf w}_{\tau_{l_1}}^{(i_1)}
\Delta{\bf w}_{\tau_{l_2}}^{(i_1)}
\Delta{\bf w}_{\tau_{l_3}}^{(i_1)}+
$$
$$
+\sum_{l_2=0}^{N-1}\sum_{l_1=0}^{l_2-1}
\sum_{l_3=0}^{l_1-1}
K'(\tau_{l_1},\tau_{l_2},\tau_{l_3})
\Delta{\bf w}_{\tau_{l_1}}^{(i_1)}
\Delta{\bf w}_{\tau_{l_2}}^{(i_1)}
\Delta{\bf w}_{\tau_{l_3}}^{(i_1)}+
$$
$$
+\sum_{l_2=0}^{N-1}\sum_{l_3=0}^{l_2-1}
\sum_{l_1=0}^{l_3-1}
K'(\tau_{l_1},\tau_{l_2},\tau_{l_3})
\Delta{\bf w}_{\tau_{l_1}}^{(i_1)}
\Delta{\bf w}_{\tau_{l_2}}^{(i_1)}
\Delta{\bf w}_{\tau_{l_3}}^{(i_1)}+
$$
$$
+\sum_{l_1=0}^{N-1}\sum_{l_2=0}^{l_1-1}
\sum_{l_3=0}^{l_2-1}
K'(\tau_{l_1},\tau_{l_2},\tau_{l_3})
\Delta{\bf w}_{\tau_{l_1}}^{(i_1)}
\Delta{\bf w}_{\tau_{l_2}}^{(i_1)}
\Delta{\bf w}_{\tau_{l_3}}^{(i_1)}+
$$
$$
+\sum_{l_1=0}^{N-1}\sum_{l_3=0}^{l_1-1}
\sum_{l_2=0}^{l_3-1}
K'(\tau_{l_1},\tau_{l_2},\tau_{l_3})
\Delta{\bf w}_{\tau_{l_1}}^{(i_1)}
\Delta{\bf w}_{\tau_{l_2}}^{(i_1)}
\Delta{\bf w}_{\tau_{l_3}}^{(i_1)}+
$$
$$
+\sum_{l_2=0}^{N-1}\sum_{l_1=0}^{l_2-1}
K'(\tau_{l_1},\tau_{l_2},\tau_{l_1})
\left(\Delta{\bf w}_{\tau_{l_1}}^{(i_1)}\right)^2
\Delta{\bf w}_{\tau_{l_2}}^{(i_1)}+
$$
$$
+
\sum_{l_3=0}^{N-1}\sum_{l_1=0}^{l_3-1}
K'(\tau_{l_1},\tau_{l_3},\tau_{l_3})
\left(\Delta{\bf w}_{\tau_{l_3}}^{(i_1)}\right)^2
\Delta{\bf w}_{\tau_{l_1}}^{(i_1)}+
$$
$$
+\sum_{l_1=0}^{N-1}\sum_{l_2=0}^{l_1-1}
K'(\tau_{l_1},\tau_{l_2},\tau_{l_2})
\left(\Delta{\bf w}_{\tau_{l_2}}^{(i_1)}\right)^2
\Delta{\bf w}_{\tau_{l_1}}^{(i_1)}+
$$
$$
+
\sum_{l_3=0}^{N-1}\sum_{l_2=0}^{l_3-1}
K'(\tau_{l_3},\tau_{l_2},\tau_{l_3})
\left(\Delta{\bf w}_{\tau_{l_3}}^{(i_1)}\right)^2
\Delta{\bf w}_{\tau_{l_2}}^{(i_1)}+
$$
$$
+\sum_{l_3=0}^{N-1}\sum_{l_2=0}^{l_3-1}
K'(\tau_{l_2},\tau_{l_2},\tau_{l_3})
\left(\Delta{\bf w}_{\tau_{l_2}}^{(i_1)}\right)^2
\Delta{\bf w}_{\tau_{l_3}}^{(i_1)}+
$$
$$
\left.+
\sum_{l_2=0}^{N-1}\sum_{l_3=0}^{l_2-1}
K'(\tau_{l_2},\tau_{l_2},\tau_{l_3})
\left(\Delta{\bf w}_{\tau_{l_2}}^{(i_1)}\right)^2
\Delta{\bf w}_{\tau_{l_3}}^{(i_1)}\right)=
$$
$$
=\frac{1}{6}\left(
\int\limits_t^T\psi_3(t_3)
\int\limits_t^{t_3}\psi_2(t_2)
\int\limits_t^{t_2}\psi_1(t_1)
d{\bf w}_{t_1}^{(i_1)}
d{\bf w}_{t_2}^{(i_1)}d{\bf w}_{t_3}^{(i_1)}+\right.
$$
$$
+
\int\limits_t^T\psi_3(t_2)
\int\limits_t^{t_2}\psi_2(t_1)
\int\limits_t^{t_1}\psi_1(t_3)
d{\bf w}_{t_3}^{(i_1)}
d{\bf w}_{t_1}^{(i_1)}d{\bf w}_{t_2}^{(i_1)}+
$$
$$
+\int\limits_t^T\psi_3(t_2)
\int\limits_t^{t_2}\psi_2(t_3)
\int\limits_t^{t_3}\psi_1(t_1)
d{\bf w}_{t_1}^{(i_1)}
d{\bf w}_{t_3}^{(i_1)}d{\bf w}_{t_2}^{(i_1)}+
$$
$$
+
\int\limits_t^T\psi_3(t_3)
\int\limits_t^{t_3}\psi_2(t_1)
\int\limits_t^{t_1}\psi_1(t_2)
d{\bf w}_{t_2}^{(i_1)}
d{\bf w}_{t_1}^{(i_1)}d{\bf w}_{t_3}^{(i_1)}+
$$
$$
+\int\limits_t^T\psi_3(t_1)
\int\limits_t^{t_1}\psi_2(t_2)
\int\limits_t^{t_2}\psi_1(t_3)
d{\bf w}_{t_3}^{(i_1)}
d{\bf w}_{t_2}^{(i_1)}d{\bf w}_{t_1}^{(i_1)}+
$$
$$
+
\int\limits_t^T\psi_3(t_1)
\int\limits_t^{t_1}\psi_2(t_3)
\int\limits_t^{t_3}\psi_1(t_2)
d{\bf w}_{t_2}^{(i_1)}
d{\bf w}_{t_3}^{(i_1)}d{\bf w}_{t_1}^{(i_1)}+
$$
$$
+
\int\limits_t^T\psi_3(t_2)
\int\limits_t^{t_2}\psi_2(t_1)\psi_1(t_1)
dt_1d{\bf w}_{t_2}^{(i_1)}+
\int\limits_t^T\psi_3(t_1)
\int\limits_t^{t_1}\psi_2(t_2)\psi_1(t_2)
dt_2d{\bf w}_{t_1}^{(i_1)}+
$$
$$
+
\int\limits_t^T\psi_3(t_3)
\int\limits_t^{t_3}\psi_2(t_1)\psi_1(t_1)
dt_1d{\bf w}_{t_3}^{(i_1)}+
\int\limits_t^T\psi_3(t_3)\psi_2(t_3)
\int\limits_t^{t_3}\psi_1(t_1)
d{\bf w}_{t_1}^{(i_1)}
dt_3+
$$
$$
\left.+\int\limits_t^T\psi_3(t_3)\psi_2(t_3)
\int\limits_t^{t_3}\psi_1(t_2)
d{\bf w}_{t_2}^{(i_1)}
dt_3+
\int\limits_t^T\psi_3(t_2)\psi_2(t_2)
\int\limits_t^{t_2}\psi_1(t_3)
d{\bf w}_{t_3}^{(i_1)}
dt_2\right)=
$$
$$
=\int\limits_t^T\psi_3(t_3)
\int\limits_t^{t_3}\psi_2(t_2)
\int\limits_t^{t_2}\psi_1(t_1)
d{\bf w}_{t_1}^{(i_1)}
d{\bf w}_{t_2}^{(i_1)}d{\bf w}_{t_3}^{(i_1)}+
$$
$$
+\frac{1}{2}
\int\limits_t^T\psi_3(t_3)
\int\limits_t^{t_3}\psi_2(t_1)\psi_1(t_1)
dt_1d{\bf w}_{t_3}^{(i_1)}+
\frac{1}{2}
\int\limits_t^T\psi_3(t_3)\psi_2(t_3)
\int\limits_t^{t_3}\psi_1(t_1)
d{\bf w}_{t_1}^{(i_1)}
dt_3=
$$
$$
=
{\int\limits_t^{*}}^T\psi_3(t_3)
{\int\limits_t^{*}}^{t_3}\psi_2(t_2)
{\int\limits_t^{*}}^{t_2}\psi_1(t_1)
d{\bf w}_{t_1}^{(i_1)}
d{\bf w}_{t_2}^{(i_1)}d{\bf w}_{t_3}^{(i_1)}\stackrel{\sf def}{=}
$$
\begin{equation}
\label{zikoxxx1}
\stackrel{\sf def}{=}
J^{*}[\psi^{(3)}]_{T,t},
\end{equation}

\vspace{2mm}
\noindent
where the multiple stochastic integral $J[K']_{T,t}^{(3)}$
is defined by (\ref{30.34}) and 
$\left\{\tau_{j}\right\}_{j=0}^{N}$ is a partition of
$[t,T],$ which satisfies the condition {\rm (\ref{1111}).

Using Proposition~2.2 for $n=3$ (see Sect.~2.1.2)
and generalizing the Fou\-rier--Legendre expansion (\ref{334.ye})
for the function $K'(t_1,t_2,t_3)$, we obtain 
$$
K'(t_1,t_2,t_3)=\lim_{p\to\infty}
\sum_{j_1=0}^{p}\sum_{j_2=0}^{p}
\sum_{j_3=0}^{p}\frac{1}{6}\biggl(C_{j_3j_2j_1}+
C_{j_3j_1j_2}+C_{j_2j_1j_3}+\biggr.
$$
\begin{equation}
\label{ziko9006}
\biggl.+C_{j_2j_3j_1}+C_{j_1j_2j_3}+C_{j_1j_3j_2}\biggr)
\phi_{j_1}(t_1)\phi_{j_2}(t_2)\phi_{j_3}(t_3),
\end{equation}

\noindent
where the multiple Fourier series (\ref{ziko9006}) converges to the function 
$K'(t_1,t_2,t_3)$ in $(t,T)^3$
and the partial sums of 
the series (\ref{ziko9006}) have an integrable majorant on $[t, T]^3$
that does not depend
on $p.$
For the trigonomertic case, the above statement follows from 
Proposition~2.2 (the proof that the function
$K'(t_1, t_2,t_3)$ belongs to the H\"{o}lder class 
with parameter $1$ in $[t, T]^3$ is omitted and 
can be carried out in the same way as for the function
$K'(t_1, t_2)$ in the two-dimensional case (see Sect.~2.1.2)).
The proof of generalization of the Fourier--Legendre expansion
(\ref{334.ye}) to the three-dimensional case (see (\ref{ziko9006}))
is omitted. The proof that the partial sums of 
the series (\ref{ziko9006}) have an integrable majorant on $[t, T]^3$
is also omitted.

Denote
$$
R'_{ppp}(t_1,t_2,t_3)=K'(t_1,t_2,t_3)-
\sum_{j_1=0}^{p}\sum_{j_2=0}^{p}
\sum_{j_3=0}^{p}\frac{1}{6}\biggl(C_{j_3j_2j_1}+
C_{j_3j_1j_2}+C_{j_2j_1j_3}+\biggr.
$$
$$
\biggl.+C_{j_2j_3j_1}+C_{j_1j_2j_3}+C_{j_1j_3j_2}\biggr)
\phi_{j_1}(t_1)\phi_{j_2}(t_2)\phi_{j_3}(t_3).
$$

\vspace{2mm}

Using Lemma~1.3 and (\ref{zikoxxx1}), we get w.~p.~1
$$
J^{*}[\psi^{(3)}]_{T,t}=J[K']_{T,t}^{(3)}=
\sum_{j_1=0}^{p}\sum_{j_2=0}^{p}
\sum_{j_3=0}^{p}\frac{1}{6}\biggl(C_{j_3j_2j_1}+
C_{j_3j_1j_2}+C_{j_2j_1j_3}+\biggr.
$$
$$
\biggl.+C_{j_2j_3j_1}+C_{j_1j_2j_3}+C_{j_1j_3j_2}\biggr)
\zeta_{j_1}^{(i_1)}\zeta_{j_2}^{(i_1)}\zeta_{j_3}^{(i_1)}+
J[R'_{ppp}]_{T,t}^{(3)}=
$$
$$
=\sum_{j_1=0}^{p}\sum_{j_2=0}^{p}
\sum_{j_3=0}^{p}C_{j_3j_2j_1}
\zeta_{j_1}^{(i_1)}\zeta_{j_2}^{(i_1)}\zeta_{j_3}^{(i_1)}+
J[R'_{ppp}]_{T,t}^{(3)}.
$$

\vspace{1mm}

Then 
$$
{\sf M}\left\{\left(J[R'_{ppp}]_{T,t}^{(3)}\right)^{2n}\right\}=
{\sf M}\left\{\left(J^{*}[\psi^{(3)}]_{T,t}-
\sum_{j_1=0}^{p}\sum_{j_2=0}^{p}
\sum_{j_3=0}^{p}C_{j_3j_2j_1}
\zeta_{j_1}^{(i_1)}\zeta_{j_2}^{(i_1)}\zeta_{j_3}^{(i_1)}\right)^{2n}\right\},
$$

\noindent
where $n\in {\bf N}.$

Applying (we mean here the passage to the limit
$\lim\limits_{p\to\infty}$)
the Lebesgue's Dominated Convergence Theorem to the integrals
on the right-hand side of (\ref{udar1}) for $k=3$ and 
$R'_{ppp}(t_1,t_2,t_3)$ instead of $R_{p_1p_2p_3}(t_1,t_2,t_3)$,
we obtain
$$
\lim\limits_{p\to\infty}
{\sf M}\left\{\left(J[R'_{ppp}]_{T,t}^{(3)}\right)^{2n}\right\}=0.
$$

\noindent
Theorem 2.85 is proved for the case $k=3.$ 

To prove Theorem 2.85 for the case $k>3$, consider the auxiliary function
\begin{equation}
\label{zikoyyy1}
K'(t_1,\ldots,t_k)=\frac{1}{k!}\left\{
\begin{matrix}
\psi_1(t_1)\ldots\psi_k(t_k),
\ \ t_1\le \ldots \le t_k
\cr\cr
\ldots 
\cr\cr
\psi_1(t_{g_1})\ldots\psi_k(t_{g_k}),
\ \ t_{g_1}\le \ldots \le t_{g_k}
\cr\cr
\ldots
\cr\cr
\psi_1(t_k)\ldots\psi_k(t_1),
\ \ t_k\le \ldots \le t_1
\end{matrix}\right.,\ \ \ \ t_1,\ldots,t_k\in[t,T],
\end{equation}
where $\{g_1,\ldots,g_k\}=\{1,\ldots,k\}$
and we take into account all
possible permutations $(g_1,\ldots,g_k)$ on the right-hand side of the formula
(\ref{zikoyyy1}).

Further, we have w.~p.~1
\begin{equation}
\label{zikozzz1a}
~~~~~~~~~~ J[K']_{T,t}^{(k)}=
J[\psi^{(k)}]_{T,t}+
\sum_{r=1}^{\left[k/2\right]}\frac{1}{2^r}
\sum_{(s_r,\ldots,s_1)\in {\rm A}_{k,r}}
J[\psi^{(k)}]_{T,t}^{s_r,\ldots,s_1},
\end{equation}

\noindent
where the function $K'(t_1,\ldots,t_k)$
is defined by (\ref{zikoyyy1}); another
notations are the same as in (\ref{30.1}) and 
Theorem~2.12 ($i_1=\ldots=i_k\ne 0$ in (\ref{30.1})).

From (\ref{zikozzz1a}) and Theorem~2.12 we obtain w.~p.~1
\begin{equation}
\label{zikozzz1}
J^{*}[\psi^{(k)}]_{T,t}=J[K']_{T,t}^{(k)}.
\end{equation}

Generalizing the above reasoning to the case $k>3$ and taking into account
(\ref{zikozzz1}), we get w.~p.~1
$$
J^{*}[\psi^{(k)}]_{T,t}=
\sum_{j_1=0}^{p}\ldots
\sum_{j_k=0}^{p}\frac{1}{k!}\left(\sum\limits_{(j_1,\ldots,j_k)}
C_{j_k\ldots j_1}\right)
\zeta_{j_1}^{(i_1)}\ldots\zeta_{j_k}^{(i_1)}+
J[R'_{p\ldots p}]_{T,t}^{(k)}=
$$
$$
=\sum_{j_1=0}^{p}\ldots
\sum_{j_k=0}^{p}C_{j_k\ldots j_1}
\zeta_{j_1}^{(i_1)}\ldots \zeta_{j_k}^{(i_1)}+
J[R'_{p\ldots p}]_{T,t}^{(k)},
$$

\noindent
where
$$
R'_{p\ldots p}(t_1,\ldots,t_k)\stackrel{\sf def}{=}K'(t_1,\ldots,t_k)-
$$
$$
-
\sum_{j_1=0}^{p}\ldots
\sum_{j_k=0}^{p}\frac{1}{k!}\left(\sum\limits_{(j_1,\ldots,j_k)}
C_{j_k\ldots j_1}\right)
\phi_{j_1}(t_1)\ldots\phi_{j_k}(t_k),
$$

\noindent
the expression
$$
\sum\limits_{(j_1, \ldots, j_k)}
$$

\noindent
means the sum with respect to 
all possible  
permutations $(j_1, \ldots, j_k)$.

Further,
$$
{\sf M}\left\{\left(J[R'_{p\ldots p}]_{T,t}^{(k)}\right)^{2n}\right\}=
{\sf M}\left\{\left(J^{*}[\psi^{(k)}]_{T,t}-
\sum_{j_1=0}^{p}\ldots
\sum_{j_k=0}^{p}C_{j_k\ldots j_1}
\zeta_{j_1}^{(i_1)}\ldots \zeta_{j_k}^{(i_1)}\right)^{2n}\right\},
$$
where $n\in {\bf N}.$

Applying (we mean here the passage to the limit
$\lim\limits_{p\to\infty}$)
the Lebesgue's Dominated Convergence Theorem to the integrals
on the right-hand side of (\ref{udar1}) for
$R'_{p\ldots p}(t_1,\ldots,t_k)$ instead of $R_{p_1\ldots,p_k}(t_1,\ldots,t_k)$,
we obtain
$$
\lim\limits_{p\to\infty}
{\sf M}\left\{\left(J[R'_{p\ldots p}]_{T,t}^{(k)}\right)^{2n}\right\}=0.
$$

\noindent
Theorem 2.85 is proved.

\section{Special Cases for Which the Expansion (\ref{march000195}) 
(Hypothesis~2.5 for $p_1=\ldots=p_k=p$) of the
Iterated Stratonovich Sto\-chastic Integrals of Multiplicity
$k$ $(k\in{\bf N})$ is Valid. Cases of Multidimensional
and Sca\-lar Wiener Process}

In this section, we will collect 
some special cases for which the expansion (\ref{march000195}) of the
iterated Stratonovich stochastic integrals of multiplicity
$k$ $(k\in{\bf N})$ is true.

{\bf Case~1.}\ The numbers $i_1,\ldots,i_k$ are pairwise different
$(i_1,\ldots,i_k=1,\ldots,m).$ It is also allowed
that some of the numbers $i_1,\ldots,i_k$ are equal to zero
(the remaining non-zero numbers from the set $\{i_1,\ldots,i_k\}$ are pairwise different).

In this case, according to Theorem~2.12, the iterated It\^{o}
and Stratonovich stochastic integrals (\ref{itoxx}) and (\ref{strxx}) coincide
w.~p.~1. At that $\psi_1(\tau),\ldots,\psi_k(\tau)$ are
continuous functions on $[t, T].$
On the other hand, the expansion (\ref{razzar1}) is valid for
the iterated It\^{o} stochastic integrals (\ref{itoxx})
when $\psi_1(\tau),\ldots,\psi_k(\tau)\in L_2([t, T]).$
Moreover, the right-hand sides of (\ref{razzar1}) ($p_1=\ldots=p_k=p$) and 
(\ref{march000195}) coincide for the case
when the non-zero numbers from the set $\{i_1,\ldots,i_k\}$ are pairwise different
$(i_1,\ldots,i_k=0,1,\ldots,m).$

{\bf Case~2.}\ The numbers $i_1,\ldots,i_k$ are such that $i_1=\ldots=i_k=i\in \{1,\ldots,m\}.$
Moreover, $\psi_1(\tau),\ldots,\psi_k(\tau)\equiv \psi(\tau),$
where $\psi(\tau)$ is a continuous function on $[t, T].$

It is well known that in this case there is an exact
representation for the iterated Stratonovich stochastic integral (\ref{strxx})
(see, for example, Sect.~6.7), which can be derived
from the expansion (\ref{march000195}) (see Sect.~2.40).

Consider here an approach that is somewhat different 
from the approach in Sect.~2.40. Using (\ref{ziko7771}), (\ref{a1}), and
Fubini's Theorem, we obtain w.~p.~1
$$
J^{*}[\psi^{(k)}]_{T,t}^{(\overbrace{{}_{i~ \ldots ~i}}^{k})}=
\frac{1}{k!}\left(\int\limits_t^T \psi(\tau) d{\bf w}^{(i)}_{\tau}\right)^k= 
$$
$$
=\frac{1}{k!}\left(\hbox{\vtop{\offinterlineskip\halign{
\hfil#\hfil\cr
{\rm l.i.m.}\cr
$\stackrel{}{{}_{p\to \infty}}$\cr
}} }\sum_{j_1=0}^{p}
C_{j_1}\zeta_{j_1}^{(i)}\right)^k=
\hbox{\vtop{\offinterlineskip\halign{
\hfil#\hfil\cr
{\rm l.i.m.}\cr
$\stackrel{}{{}_{p\to \infty}}$\cr
}} }
\frac{1}{k!}\left(\sum_{j_1=0}^{p}
C_{j_1}\zeta_{j_1}^{(i)}\right)^k= 
$$
$$
=\hbox{\vtop{\offinterlineskip\halign{
\hfil#\hfil\cr
{\rm l.i.m.}\cr
$\stackrel{}{{}_{p\to \infty}}$\cr
}} }
\frac{1}{k!}\sum_{j_1,\ldots,j_k=0}^{p}
C_{j_1}\ldots C_{j_k}\zeta_{j_1}^{(i)}\ldots \zeta_{j_k}^{(i)}= 
$$
$$
=\hbox{\vtop{\offinterlineskip\halign{
\hfil#\hfil\cr
{\rm l.i.m.}\cr
$\stackrel{}{{}_{p\to \infty}}$\cr
}} }
\frac{1}{k!}\sum_{j_1,\ldots,j_k=0}^{p}
\left(\sum\limits_{(j_1,\ldots,j_k)}
C_{j_k\ldots j_1}\right)\zeta_{j_1}^{(i)}\ldots \zeta_{j_k}^{(i)}= 
$$
$$
=\hbox{\vtop{\offinterlineskip\halign{
\hfil#\hfil\cr
{\rm l.i.m.}\cr
$\stackrel{}{{}_{p\to \infty}}$\cr
}} }
\frac{1}{k!}\sum\limits_{(j_1,\ldots,j_k)}\left(\sum_{j_1,\ldots,j_k=0}^{p}
C_{j_k\ldots j_1}\zeta_{j_1}^{(i)}\ldots \zeta_{j_k}^{(i)}\right)= 
$$
$$
=\hbox{\vtop{\offinterlineskip\halign{
\hfil#\hfil\cr
{\rm l.i.m.}\cr
$\stackrel{}{{}_{p\to \infty}}$\cr
}} }
\frac{1}{k!} \cdot k! \sum_{j_1,\ldots,j_k=0}^{p}
C_{j_k\ldots j_1}\zeta_{j_1}^{(i)}\ldots \zeta_{j_k}^{(i)}= 
$$
$$
=\hbox{\vtop{\offinterlineskip\halign{
\hfil#\hfil\cr
{\rm l.i.m.}\cr
$\stackrel{}{{}_{p\to \infty}}$\cr
}} }
\sum_{j_1,\ldots,j_k=0}^{p}
C_{j_k\ldots j_1}\zeta_{j_1}^{(i)}\ldots \zeta_{j_k}^{(i)},
$$

\vspace{2mm}
\noindent
where notations as in Theorems~2.61, 2.60, 
$$
\sum\limits_{(j_1,\ldots,j_k)}
$$ 
means the sum with respect to all
possible permutations 
$(j_1,\ldots,j_k),$
and
$\psi_1(\tau),\ldots,\psi_k(\tau)\equiv \psi(\tau),$
where $\psi(\tau)$ is a continuous function on $[t, T].$

{\bf Case~3.}\ Among the numbers $1, 2,\ldots, k$ there are
numbers $g_1,\ldots,g_r$ $(g_1+1 < g_2,\ldots, g_{r-1}+1 < g_r)$ 
such that $i_{g_1}=i_{g_1+1},\ldots,i_{g_r}=i_{g_r+1},$
where $r=1,2,\ldots,[k/2],$ $i_1,\ldots,i_k=1,\ldots,m,$ and
the set
$$
\bigl\{i_1,\ldots,i_k\bigr\}\backslash
\bigl\{i_{g_1},i_{g_1+1},\ldots,i_{g_r},i_{g_r+1}\bigr\}
$$                                               

\noindent
consists of pairwise different elements, none of which
is equal to $i_{g_1},i_{g_1+1},\ldots,$ $i_{g_r},i_{g_r+1}.$
Moreover, $\psi_1(\tau),\ldots,$ $\psi_k(\tau)\equiv 1$ or
$\psi_1(\tau),\ldots,\psi_k(\tau)\equiv \psi(\tau),$
where $\psi(\tau)$ is a continuous function on $[t, T].$

It is obvious that the left-hand side of condition
(\ref{09091}) for Case~3 will decompose into a product
of several expressions of the same type as the
left-hand side of inequality (\ref{febr2025}).
Thus, the condition (\ref{09091}) is satisfied
for Case~3 and the expansion (\ref{march000195}) 
is true for this case.

Consider the ordered set $(i_1,\ldots,i_k)$
that corresponds to Case~3. We have
\begin{equation}
\label{september20256}
~~~~~~~~~~~~(i_1,\ldots,i_{g_1-1},i_{g_1},i_{g_1},i_{g_1+2},\ldots,
i_{g_r-1},i_{g_r},i_{g_r},i_{g_r+2},\ldots,i_k),
\end{equation}

\vspace{1mm}
\noindent
where
$r=1,2,$ $\ldots,[k/2],$ $i_1,\ldots,i_k=1,\ldots,m,$ and
the set
$$
\bigl\{i_1,\ldots,i_k\bigr\}\backslash
\bigl\{i_{g_1},i_{g_1},\ldots,i_{g_r},i_{g_r}\bigr\}
$$

\noindent
consists of pairwise different elements, none of which
is equal to $i_{g_1},\ldots,i_{g_r}.$
Moreover, $i_{g_1},\ldots,i_{g_r}$ are pairwise different.

Consider some particular cases of (\ref{september20256}).

{\bf Case~3.1:} 
\begin{equation}
\label{september202510}
(i_1,i_1,i_3,i_3,i_5,i_6,\ldots,i_{k-1},i_{k}),
\end{equation}

\noindent
where $i_1,i_3,i_5,i_6,\ldots,i_{k-1},i_{k}=1,\ldots,m$ and
$i_1,i_3,i_5,i_6,\ldots,i_{k-1},i_{k}$
are pairwise different.

{\bf Case~3.2:} 
\begin{equation}
\label{september20259}
(i_1,i_1,i_3,i_4,\ldots,i_{k-2},i_{k-1},i_{k-1}),
\end{equation}

\noindent
where $i_1,i_3,i_4,\ldots,i_{k-2},i_{k-1}=1,\ldots,m$ and
$i_1,i_3,i_4,\ldots,i_{k-2},i_{k-1}$
are pairwise different.

{\bf Case~3.3:}
\begin{equation}
\label{september20257}
(i_1,i_1,i_2,i_2,\ldots,i_r,i_r),
\end{equation}

\noindent
where $k=2r,$ $i_1,\ldots,i_r=1,\ldots,m$
and $i_1,\ldots,i_r$ are pairwise different.

{\bf Case~3.4:} 
\begin{equation}
\label{september20258}
(i_1,i_1,i_3,i_4,i_4,i_6,\ldots,i_{3r-2},i_{3r-2}, i_{3r}),
\end{equation}

\noindent
where $k=3r,$ $i_1,i_3,i_4,i_6,\ldots,i_{3r-2}, i_{3r}=1,\ldots,m$
and $i_1,i_3,i_4,i_6,\ldots,i_{3r-2}, i_{3r}$ are pairwise different.

Also, we can consider a modification of (\ref{september20256})
in the sense that $i_1,\ldots,i_k=0, 1,\ldots,m$ in (\ref{september20256}), i.e.
some of the numbers $i_1,\ldots,i_k$ are equal to zero
in (\ref{september20256}).

In conclusion of this section, we will consider the expansion
(\ref{march000195}) of the iterated Stratonovich stochastic integral
(\ref{strxx}), which corresponds to
(\ref{september20257}) and (\ref{september20258}).
We have
$$
J^{*}[\psi^{(2r)}]_{T,t}^{(i_1 i_1 i_2 i_2\ldots i_r i_r)}
=
$$
$$
=
\hbox{\vtop{\offinterlineskip\halign{
\hfil#\hfil\cr
{\rm l.i.m.}\cr
$\stackrel{}{{}_{p\to \infty}}$\cr
}} }
\sum\limits_{j_1, j_2,\ldots,j_{2r}=0}^p
C_{j_{2r}\ldots j_2 j_1} \zeta_{j_1}^{(i_1)}\zeta_{j_2}^{(i_1)}\zeta_{j_3}^{(i_2)}\zeta_{j_4}^{(i_2)}
\ldots \zeta_{j_{2r-1}}^{(i_{2r})}\zeta_{j_{2r}}^{(i_{2r})},
$$

\vspace{2mm}

$$
J^{*}[\psi^{(3r)}]_{T,t}^{(i_1 i_1 i_3 i_4 i_4 i_6\ldots i_{3r-2} i_{3r-2}  i_{3r})}
=
$$
$$
=
\hbox{\vtop{\offinterlineskip\halign{
\hfil#\hfil\cr
{\rm l.i.m.}\cr
$\stackrel{}{{}_{p\to \infty}}$\cr
}} }
\sum\limits_{j_1,  j_2, \ldots, j_{3r}=0}^p
C_{j_{3r}\ldots j_2 j_1} \zeta_{j_1}^{(i_1)}\zeta_{j_2}^{(i_1)}\zeta_{j_3}^{(i_3)}
\zeta_{j_4}^{(i_4)}\zeta_{j_5}^{(i_4)}\zeta_{j_6}^{(i_6)}
\ldots \zeta_{j_{3r-2}}^{(i_{3r-2})}\zeta_{j_{3r-1}}^{(i_{3r-2})}\zeta_{j_{3r}}^{(i_{3r})},
$$

\noindent
where notations as in 
(\ref{march000195}) and $\psi_1(\tau),\ldots,\psi_{3r}(\tau)\equiv \psi(\tau)$,
where $\psi(\tau)$ is a continuous function on $[t, T]$ 
(in particular $\psi(\tau)\equiv 1$).

\section{One Result on the Expansion of Multiple Stra\-to\-no\-vich
Stochastic Integrals of Multiplicity $k.$
The Case $i_1=\ldots=i_k=1,\ldots,m$}

Let us consider the 
multiple stochastic integral (\ref{30.34})
\begin{equation}
\label{30.34ququ}
~~~~~~~~~~ \hbox{\vtop{\offinterlineskip\halign{
\hfil#\hfil\cr
{\rm l.i.m.}\cr
$\stackrel{}{{}_{N\to \infty}}$\cr
}} }\sum_{j_1,\ldots,j_k=0}^{N-1}
\Phi\left(\tau_{j_1},\ldots,\tau_{j_k}\right)
\prod\limits_{l=1}^k\Delta{\bf w}_{\tau_{j_l}}^{(i_l)}
\stackrel{\rm def}{=}J[\Phi]_{T,t}^{(i_1\ldots i_k)},
\end{equation}

\noindent
where we assume that
$\Phi(t_1,\ldots,t_k):\ [t, T]^k\to{\bf R}$ is a 
continuous nonrandom
function on $[t, T]^k,$ 
${\bf w}_{\tau}$ is a random vector with 
an $m+1$ components
(${\bf w}_{\tau}^{(i)} $ $(i=1,\ldots,m)$
are independent standard Wiener processes and
${\bf w}_{\tau}^{(0)}=\tau$),
$\Delta {\bf w}_{\tau_j}^{(i)}={\bf w}_{\tau_{j+1}}^{(i)}-{\bf w}_{\tau_j}^{(i)},$
$\left\{\tau_{j}\right\}_{j=0}^{N}$ is a partition of
$[t,T]$ which satisfies the condition (\ref{1111}),
and $i_1,\ldots,i_k=0,1,\ldots,m.$

The stochastic integral with respect to the scalar standard Wiener process
($i_1=\ldots=i_k\ne 0$)
and similar to (\ref{30.34ququ}) 
(the function $\Phi(t_1,\ldots,t_k)$ is assumed to be symmetric
on the hypercube $[t, T]^k$)
has been considered in literature
(see, for example, Remark~1.5.7 \cite{bugh1}). 
The integral (\ref{30.34ququ})
is sometimes called the multiple Stratonovich stochastic integral.
This is due to the fact that the following rule
of the classical integral calculus holds for this integral
(see~Lemma 1.3)
$$
J[\Phi]_{T,t}^{(i_1\ldots i_k)}=
J[\varphi_1]_{T,t}^{(i_1)}\ldots J[\varphi_k]_{T,t}^{(i_k)}\ \ \ \hbox{\rm w.~p.~1},
$$

\vspace{4mm}
\noindent
where $\Phi(t_1,\ldots,t_k)=\varphi_1(t_1)\ldots \varphi_k(t_k)$
and
$$
J[\varphi_l]_{T,t}^{(i_l)}
=\int\limits_t^T \varphi_l(s) d{\bf w}_{s}^{(i_l)}\ \ \ (l=1,\ldots,k).
$$

It is not difficult to see that for the case $i_1=\ldots=i_k\ne 0$
we have w.~p.~1
$$
\hbox{\vtop{\offinterlineskip\halign{
\hfil#\hfil\cr
{\rm l.i.m.}\cr
$\stackrel{}{{}_{N\to \infty}}$\cr
}} }\sum_{j_1,\ldots,j_k=0}^{N-1}
\Phi\left(\tau_{j_1},\ldots,\tau_{j_k}\right)
\Delta{\bf w}_{\tau_{j_1}}^{(i_1)}\ldots \Delta{\bf w}_{\tau_{j_k}}^{(i_1)}
=
$$
$$
=\hbox{\vtop{\offinterlineskip\halign{
\hfil#\hfil\cr
{\rm l.i.m.}\cr
$\stackrel{}{{}_{N\to \infty}}$\cr
}} }\sum_{j_1,\ldots,j_k=0}^{N-1}
\frac{1}{k!}\left(\sum\limits_{(j_1,\ldots,j_k)}
\Phi\left(\tau_{j_1},\ldots,\tau_{j_k}\right)\right)
\Delta{\bf w}_{\tau_{j_1}}^{(i_1)}\ldots \Delta{\bf w}_{\tau_{j_k}}^{(i_1)},
$$

\noindent
i.e.
\begin{equation}
\label{ziko6007}
J[\Phi]_{T,t}^{(\hspace{0.5mm}\overbrace{{}_{i_1 \ldots i_1}}^{k}\hspace{0.5mm})}
=J[\tilde \Phi]_{T,t}^{(\hspace{0.5mm}\overbrace{{}_{i_1 \ldots i_1}}^{k}\hspace{0.5mm})}
\ \ \ \hbox{w.~p.~1},
\end{equation}

\noindent
where 
$$
\tilde \Phi(t_1,\ldots,t_k)=
\frac{1}{k!}\left(\sum\limits_{(t_1,\ldots,t_k)}
\Phi(t_1,\ldots,t_k)\right)
$$

\vspace{1mm}
\noindent
is the symmetrization of the function $\Phi(t_1,\ldots,t_k)$;
the expression
$$
\sum\limits_{(a_1,\ldots,a_k)}
$$ 

\noindent
means the sum with respect to all
possible permutations
$(a_1,\ldots,a_k).$

Due to (\ref{ziko6007}) the condition of symmetry
of the function $\Phi(t_1,\ldots,t_k)$ need not be required
in the case $i_1=\ldots=i_k\ne 0.$

{\bf Definition 2.1} \cite{bugh2}.\ {\it Let $\Phi(t_1,\ldots,t_k)\in L_2([t, T]^k)$
is a symmetric function and $q=1, 2,\ldots, [k/2]$ {\rm (}$q$ is fixed{\rm )}.
Suppose that for every complete orthonormal system of functions
$\{\phi_j(x)\}_{j=0}^{\infty}$ in the space $L_2([t, T])$
the following sum
$$
\sum\limits_{j_1,\ldots,j_q=0}^p
\sum\limits_{j_{2q+1},\ldots,j_k=0}^p\
\int\limits_{[t,T]^{k}}
\Phi(t_1,\ldots,t_k)\phi_{j_1}(t_1)\phi_{j_1}(t_2)\ldots
\phi_{j_q}(t_{2q-1})\phi_{j_q}(t_{2q})\times
$$
\begin{equation}
\label{ziko6009}
~~~~~~~\times
\phi_{j_{2q+1}}(t_{2q+1})\ldots  \phi_{j_k}(t_k)dt_1\ldots dt_{k} \cdot
\phi_{j_{2q+1}}(t_{2q+1})\ldots  \phi_{j_k}(t_k)
\end{equation}

\vspace{3mm}
\noindent
converges in $L_2([t,T]^{k-2q})$ if $p\to\infty$
to a limit, which is independent of the choice of the
complete orthonormal system of functions
$\{\phi_j(x)\}_{j=0}^{\infty}$ in the space $L_2([t, T]).$
Then we say that the $q$th limiting trace for $\Phi(t_1,\ldots,t_k)$ exists,
which by definition is the limit of the sum {\rm (\ref{ziko6009})}
and is denoted as
$\overrightarrow{Tr}^{{}^{q}} \Phi$. Moreover, $\overrightarrow{Tr}^{{}^{0}} \Phi
\stackrel{\sf def}{=}\Phi$.}

Consider the following Theorem using our notations.

{\bf Theorem 2.86} \cite{bugh2}. {\it Let $\Phi(t_1,\ldots,t_k)\in L_2([t, T]^k)$
is a symmetric nonrandom function. Furthermore, let all limiting traces for
$\Phi(t_1,\ldots,t_k)$ {\rm (}see De\-finition {\rm 2.1}{\rm )} exist. Then
the following expansion
$$
J^{\circ}[\Phi]_{T,t}^{(\hspace{0.5mm}\overbrace{{}_{i_1 \ldots i_1}}^{k}\hspace{0.5mm})}=
\hbox{\vtop{\offinterlineskip\halign{
\hfil#\hfil\cr
{\rm l.i.m.}\cr
$\stackrel{}{{}_{p\to \infty}}$\cr
}} }
\sum\limits_{j_1,\ldots j_k=0}^{p}
C_{j_k \ldots j_1} \zeta_{j_1}^{(i_1)}\ldots \zeta_{j_k}^{(i_1)}
$$

\vspace{1mm}
\noindent
that converges in the mean-square sense is valid, where 
$$
J^{\circ}[\Phi]_{T,t}^{(\hspace{0.5mm}\overbrace{{}_{i_1 \ldots i_1}}^{k}\hspace{0.5mm})}
$$
is the multiple Stratonovich stochastic integral defined as in {\rm \cite{bugh3} (1993)}
{\rm (}also see {\rm \cite{bugh2},} pp.~{\rm 910--911),}
$$
C_{j_k\ldots j_1}=\int\limits_{[t,T]^k}
\Phi(t_1,\ldots,t_k)\prod_{l=1}^{k}\phi_{j_l}(t_l)dt_1\ldots dt_k
$$
is the Fourier coefficient, 
${\rm l.i.m.}$ is a limit in the mean-square sense,
$i_1=1,\ldots,m,$
$$
\zeta_{j}^{(i_1)}=
\int\limits_t^T \phi_{j}(s) d{\bf w}_s^{(i_1)}
$$ 

\noindent
are independent standard Gaussian random variables for various 
$j.$ 
}

In addition to the conditions of Theorem~2.86, we assume that
the function 
$\Phi(t_1,\ldots,t_k)$ is continuous on $[t, T]^k.$ Then 
\cite{bugh1}
\begin{equation}
\label{ziko09}
J^{\circ}[\Phi]_{T,t}^{(\hspace{0.5mm}\overbrace{{}_{i_1 \ldots i_1}}^{k}\hspace{0.5mm})}=
J[\Phi]_{T,t}^{(\hspace{0.5mm}\overbrace{{}_{i_1 \ldots i_1}}^{k}\hspace{0.5mm})}\ \ \ 
\hbox{w.~p.~1,}
\end{equation}

\noindent
where the multiple Stratonovich stochastic integral
$$
J[\Phi]_{T,t}^{(\hspace{0.5mm}\overbrace{{}_{i_1 \ldots i_1}}^{k}\hspace{0.5mm})}
$$

\noindent
is defined by (\ref{30.34ququ}).
As a result, we get the following expansion
\begin{equation}
\label{ziko777}
J[\Phi]_{T,t}^{(\hspace{0.5mm}\overbrace{{}_{i_1 \ldots i_1}}^{k}\hspace{0.5mm})}=
\hbox{\vtop{\offinterlineskip\halign{
\hfil#\hfil\cr
{\rm l.i.m.}\cr
$\stackrel{}{{}_{p\to \infty}}$\cr
}} }
\sum\limits_{j_1,\ldots j_k=0}^{p}
C_{j_k \ldots j_1} \zeta_{j_1}^{(i_1)}\ldots \zeta_{j_k}^{(i_1)}.
\end{equation}

It should be noted that the expansion (\ref{ziko777})
is valid provided that for the function 
$\Phi(t_1,\ldots,t_k)$ there exist all limiting traces
that do not depend on the choice of the 
complete orthonormal system of functions
$\{\phi_j(x)\}_{j=0}^{\infty}$ in the space $L_2([t, T]).$
The last condition is essential for the
proof of the equality (\ref{ziko09})
(this proof follows from Theorem~1.5.3, Remark~1.5.7, and Propositions~2.2.3, 2.2.5, 4.1.2
\cite{bugh1}).
More precisely, in \cite{bugh1}, to prove Proposition~4.1.2 (p.~65)
a special basis $\{\phi_j(x)\}_{j=0}^{\infty}$ was used.
This means that the existence of a limit of the sum (\ref{ziko6009})
for the function $\Phi(t_1,\ldots,t_k)$
in the case when $\{\phi_j(x)\}_{j=0}^{\infty}$ 
is an arbitrary complete orthonormal system of 
functions in the space $L_2([t, T])$
requires a separate proof.

It is not difficult to show that (see (\ref{zikozzz1}))
\begin{equation}
\label{ziko779}
~~J^{*}[\psi^{(k)}]_{T,t}^{(\hspace{0.5mm}\overbrace{{}_{i_1 \ldots i_1}}^{k}\hspace{0.5mm})}=
J[K']_{T,t}^{(\hspace{0.5mm}\overbrace{{}_{i_1 \ldots i_1}}^{k}\hspace{0.5mm})}
\ \ \ \hbox{w.~p.~1,}
\end{equation}

\vspace{2mm}
\noindent
where 
$$
J^{*}[\psi^{(k)}]_{T,t}^{(\hspace{0.5mm}\overbrace{{}_{i_1 \ldots i_1}}^{k}\hspace{0.5mm})}
$$

\vspace{1mm}
\noindent
is the iterated Stratonovich stochastic integral (\ref{str}) $(i_1=\ldots=i_k\ne 0),$ 
$$
J[K']_{T,t}^{(\hspace{0.5mm}\overbrace{{}_{i_1 \ldots i_1}}^{k}\hspace{0.5mm})}
$$

\vspace{1mm}
\noindent
is the multiple Stratonovich stochastic integral (\ref{30.34ququ}) $(i_1=\ldots=i_k\ne 0)$
for the continuous function
$K'(t_1,\ldots,t_k)$ defined by (\ref{zikoyyy1}).

If we assume that the limiting
traces from Theorem~2.86 exist, then we can write (see (\ref{ziko779}))
$$
J^{*}[\psi^{(k)}]_{T,t}^{(\overbrace{{}_{i_1 \ldots i_1}}^{k})}=
J[K']_{T,t}^{(\overbrace{{}_{i_1 \ldots i_1}}^{k})}=
$$
$$
=
\hbox{\vtop{\offinterlineskip\halign{
\hfil#\hfil\cr
{\rm l.i.m.}\cr
$\stackrel{}{{}_{p\to \infty}}$\cr
}} }
\sum\limits_{j_1,\ldots j_k=0}^{p}
\left(\frac{1}{k!}
\sum\limits_{(j_1,\ldots,j_k)}C_{j_k \ldots j_1}\right) 
\zeta_{j_1}^{(i_1)}\ldots \zeta_{j_k}^{(i_1)}=
$$
\begin{equation}
\label{ziko12000}
=\hbox{\vtop{\offinterlineskip\halign{
\hfil#\hfil\cr
{\rm l.i.m.}\cr
$\stackrel{}{{}_{p\to \infty}}$\cr
}} }
\sum\limits_{j_1,\ldots j_k=0}^{p}
C_{j_k \ldots j_1} \zeta_{j_1}^{(i_1)}\ldots \zeta_{j_k}^{(i_1)},
\end{equation}
where
\begin{equation}
\label{ziko9000}
~~~~~~~C_{j_k \ldots j_1}=\int\limits_t^T\psi_k(t_k)\phi_{j_k}(t_k)\ldots
\int\limits_t^{t_2}
\psi_1(t_1)\phi_{j_1}(t_1)
dt_1\ldots dt_k
\end{equation}

\vspace{1mm}
\noindent
is the Fourier coefficient,
$$
\sum\limits_{(j_1,\ldots,j_k)}
$$ 

\noindent
means the sum with respect to all
possible permutations
$(j_1,\ldots,j_k);$ another notations are the same as in Theorem~2.86.

The equality (\ref{ziko12000}) agrees with Theorems~2.49, 2.59, 2.61
for the particular case $i_1=\ldots=i_k\ne 0.$

From the other hand, the following expansion (see Theorem~2.49)
\begin{equation}
\label{ziko12001}
~~~~~~~J^{*}[\psi^{(k)}]_{T,t}^{(i_1\ldots i_k)}=
\hbox{\vtop{\offinterlineskip\halign{
\hfil#\hfil\cr
{\rm l.i.m.}\cr
$\stackrel{}{{}_{p\to \infty}}$\cr
}} }\sum_{j_1,\ldots,j_k=0}^{p}
C_{j_k\ldots j_1}\zeta_{j_1}^{(i_1)}\ldots \zeta_{j_k}^{(i_k)}
\end{equation}

\vspace{1mm}
\noindent
is valid, where $C_{j_k \ldots j_1}$ has the form (\ref{ziko9000}),
$i_1,\ldots,i_k=0, 1,\ldots,m$; another notations are the same as in Theorem~2.49.

\section{About Hypotheses~2.4 and 2.5}

In the previous section, we saw that in a number of papers (see, for example, 
\cite{bugh1}-\cite{bugh3}) 
the conditions of theorems related to multiple stochastic integrals 
(see Theorem~2.86 in Sect.~2.46)
are formulated 
in terms of limiting traces (see Definition~2.1).
In addition to limiting traces, the concept of Hilbert space 
valued traces (integral traces) is introduced in \cite{bugh1}.
The concepts of traces considered in \cite{bugh1}-\cite{bugh3} 
are close to some expressions 
that we used in this chapter.
For example, the following integral (see (\ref{5t}))
$$
\frac{1}{2}\int\limits_t^T\psi_1(t_1)\psi_2(t_1)dt_1=\int\limits_t^T K^{*}(t_1,t_1)dt_1
$$
is an example of a trace introduced in \cite{bugh1} (Definition~2.2.1).
However, the function $K^{*}(t_1,t_2)$ defined by (\ref{1999.1})
is not symmetric compared with \cite{bugh1} (De\-finition~2.2.1).
In addition, the expression (see (\ref{5t}))

\vspace{-1mm}
$$
\sum_{j_1=0}^{\infty}
C_{j_1j_1}=\sum_{j_1=0}^{\infty}~\int\limits_{[t,T]^2}K(t_1,t_2)\phi_{j_1}(t_1)
\phi_{j_1}(t_2)dt_1 dt_2
$$

\vspace{1mm}
\noindent
is an example of a limiting trace (see Definition~2.1)
for the function $K(t_1,t_2),$ which is not symmetric (see (\ref{ziko0909})).

In this section, we will talk about Hypotheses~2.4 and 2.5 again.
It should be noted that a significant part of Chapter 2 is devoted 
to the proof of Hypothesis~2.5 
for various special cases (Theorems~2.1--2.9,
2.33--2.36, 2.41, 2.45--2.48, 2.50, 2.51, 2.59, 2.61--2.65).
In order to prove these theorems, we developed 
a number of approaches for expansion of iterated Stratonovich stochastic integrals.

More precisely, in Theorems~2.1, 2.2, 2.4--2.9, 2.33--2.36, 2.41, 2.64, 2.65 we assume that
$\{\phi_j(x)\}_{j=0}^{\infty}$ is a complete orthonormal system of 
Legendre polynomials or trigonometric functions in the space $L_2([t, T]).$
The above systems of functions 
are most suitable for the expansion of iterated stochastic integrals
from the Taylor--It\^{o} and Taylor--Stra\-to\-no\-vich expansions
(see Chapter~5). In Theorems~2.3, 2.47, 2.48, 2.50, 2.51, 2.59, 2.61--2.63
the system $\{\phi_j(x)\}_{j=0}^{\infty}$ 
can be arbitrary.

Note that Theorems~2.1--2.3 
are special cases of Hypothesis~2.5
for $k=2$ and $p_1, p_2\to\infty$.
At that 
$\psi_2(\tau)$ is a continuously dif\-ferentiable 
nonrandom function on $[t, T]$ and $\psi_1(\tau)$ is twice 
continuously differentiable nonrandom function on $[t, T]$
(Theorem~2.1). In Theorem~2.2, the functions 
$\psi_1(\tau)$ and $\psi_2(\tau)$ are assumed to be
continuously differentiable only one time on $[t, T]$.
Theorem~2.3 is a ge\-ne\-ral\-ization of Theorems~2.1, 2.2.
More precisely, in Theorem~2.3 we assume that $\psi_1(\tau)$ and $\psi_2(\tau)$
are continuous functions on $[t, T]$.

Theorems~2.4--2.8, 2.33, 2.41, 2.47, 2.51 are special cases
of Hypothesis~2.5 for $k=3$.
In Theorems~2.4 and~2.6, the case
$\psi_1(\tau), \psi_2(\tau), \psi_3(\tau)\equiv 1$
and $p_1,p_2,p_3\to\infty$ is considered.
Theorem~2.8 is a special case of Hypothesis~2.5 for the case when
$\psi_2(\tau)$ is a continuously dif\-ferentiable 
nonrandom function on $[t, T]$ and $\psi_1(\tau), \psi_3(\tau)$ are twice
continuously differentiable nonrandom functions on $[t, T]$
$(p_1=p_2=p_3=p\to\infty)$.
Theorem~2.33 is an analogue
of Theorem~2.8 for continuously differentiable functions $\psi_1(\tau), 
\psi_2(\tau), \psi_3(\tau)$ and $p_1=p_2=p_3=p\to\infty$.
Theorem~2.41 is a generalization of Theorem~2.33 
for the case $p_1, p_2, p_3\to\infty$.
In Theorems~2.5, 2.7 and 2.51, we consider narrower particular
cases of the functions 
$\psi_1(\tau), \psi_2(\tau), \psi_3(\tau).$
For example, the functions
$\psi_1(\tau), \psi_2(\tau), \psi_3(\tau)$ 
have a binomial form (Theorems~2.5, 2.51).

Theorems~2.9, 2.34, 2.35, 2.48, 2.50,
2.63 are
special cases of Hypothesis~2.5 for $k=4$ and $k=5$.
The functions $\psi_1(\tau), \ldots, \psi_5(\tau)$
are continuously differentiable on $[t, T]$ 
in Theorems~2.34, 2.35,
$\psi_1(\tau), \ldots, \psi_5(\tau)\equiv 1$ in Theorems~2.9, 2.48, 2.50,
and 
$\psi_1(\tau), \ldots, \psi_4(\tau)$ have a binomial form in Theorem~2.63.

In Theorems~2.36, 2.62, 2.64, 2.65 the cases $k=6, 7, 8$ of Hypothesis~2.5 is considered. 
At that $\psi_1(\tau), \ldots, \psi_8(\tau)\equiv 1$
in these theorems.

Theorems~2.59, 2.61 prove Hypothesis~2.5 for $k\in {\bf N}$
but under one additional condition (see (\ref{july700001xyz}) or (\ref{09091})).

Let us conclude this section with a few remarks.

\vspace{2mm}          

{\bf Remark~2.5.}\ {\it The equalities 
{\rm (\ref{tyyy})} and {\rm (\ref{leto6000})} imply that the equality
{\rm (\ref{july300000})} is equivalent to the relation 
$$
\sum_{r=1}^{\left[k/2\right]}\frac{1}{2^r}
\sum_{(s_r,\ldots,s_1)\in {\rm A}_{k,r}}
J[\psi^{(k)}]_{T,t}^{s_r,\ldots,s_1}=
$$
$$
= -
\hbox{\vtop{\offinterlineskip\halign{
\hfil#\hfil\cr
{\rm l.i.m.}\cr
$\stackrel{}{{}_{p_1,\ldots,p_k\to \infty}}$\cr
}} }\sum_{j_1=0}^{p_1}\ldots\sum_{j_k=0}^{p_k}
C_{j_k\ldots j_1}
\sum\limits_{r=1}^{[k/2]}
(-1)^r \times
\Biggr.
$$

\vspace{-1mm}
\begin{equation}
\label{ziko0020ddd}
\times
\sum_{\stackrel{(\{\{g_1, g_2\}, \ldots, 
\{g_{2r-1}, g_{2r}\}\}, \{q_1, \ldots, q_{k-2r}\})}
{{}_{\{g_1, g_2, \ldots, 
g_{2r-1}, g_{2r}, q_1, \ldots, q_{k-2r}\}=\{1, 2, \ldots, k\}}}}
\prod\limits_{s=1}^r
{\bf 1}_{\{i_{g_{{}_{2s-1}}}=~i_{g_{{}_{2s}}}\ne 0\}}
\Biggl.{\bf 1}_{\{j_{g_{{}_{2s-1}}}=~j_{g_{{}_{2s}}}\}}
\prod_{l=1}^{k-2r}\zeta_{j_{q_l}}^{(i_{q_l})}\Biggr)
\end{equation}

\vspace{2mm}
\noindent
w.~p.~{\rm 1}$,$ where
notations are the same as in Theorems~{\rm 1.2, 1.16} and~{\rm 2.12.}}

\vspace{2mm}

{\bf Remark~2.6.}\ {\it Applying Theorems~{\rm 1.14,} {\rm 1.16,} we can reformulate
the equality
{\rm (\ref{ziko0020ddd})} as follows
$$
\sum_{r=1}^{\left[k/2\right]}\frac{1}{2^r}
\sum_{(s_r,\ldots,s_1)\in {\rm A}_{k,r}}
J[\psi^{(k)}]_{T,t}^{s_r,\ldots,s_1}= 
$$
$$
=
\hbox{\vtop{\offinterlineskip\halign{
\hfil#\hfil\cr
{\rm l.i.m.}\cr
$\stackrel{}{{}_{p_1,\ldots,p_k\to \infty}}$\cr
}} }\sum_{j_1=0}^{p_1}\ldots\sum_{j_k=0}^{p_k}
C_{j_k\ldots j_1}\times   
$$

\vspace{-3mm}
$$
\times \left(\zeta_{j_1}^{(i_1)}\ldots \zeta_{j_k}^{(i_k)} -
\prod_{l=1}^k\left({\bf 1}_{\{m_l=0\}}+{\bf 1}_{\{m_l>0\}}\left\{
\begin{matrix}
\prod\limits_{s=1}^{d_l} 
H_{n_{s,l}}\left(\zeta_{j_{h_{s,l}}}^{(i_l)}\right),\ 
&\hbox{\rm if}\ \ \ 
i_l\ne 0\cr\cr
\prod\limits_{s=1}^{d_l}  
\left(\zeta_{j_{h_{s,l}}}^{(0)}\right)^{n_{s,l}},\  &\hbox{\rm if}\ \ \ 
i_l=0
\end{matrix}\right.\ \right)\right)
$$

\vspace{1mm}
\noindent
w.~p.~{\rm 1}$,$ where
notations are the same as in Theorems~{\rm 1.14,} {\rm 1.16} and~{\rm 2.12}
{\rm(}$H_n(x)$ is the Hermite polynomial {\rm (\ref{ziko500}))}.}

\vspace{2mm}

{\bf Remark~2.7.}\ {\it Recently$,$ in {\rm \cite{Rybakov3000},} an approach to the proof 
of expansion similar to {\rm (\ref{july300001})} 
was proposed.
In particular$,$ this approach uses the representation of the multiple Stratonovich
stochastic integral {\rm (\ref{january19})}
as the sum of some constant value and multiple Wiener stochastic integrals of 
multiplicities not exceeding $k.$
Note that a similar representation in a different form is defined by 
the formula {\rm (\ref{after8xxds1}).}

It should be noted that an expansion similar to {\rm (\ref{july300001})} was considered 
in {\rm \cite{Rybakov3000}}
for an arbitrary $k.$ The system of basis functions $\{\phi_j(x)\}_{j=0}^{\infty}$
in the space $L_2([t, T])$
can also be arbitrary. However$,$ in {\rm \cite{Rybakov3000},} 
the condition on convergence of 
trace series is used as a sufficient condition for the validity of expansion similar to
{\rm (\ref{july300001})}
{\rm (}see {\rm \cite{Rybakov3000}} for details{\rm ).} 
Note that the verification of the above condition for the kernel {\rm (\ref{ppp})} is a separate
problem. 

In Theorems {\rm 2.24--2.29, 2.36--2.40} the rate of mean-square convergence of expansions 
of iterated Stratonovich stochastic integrals is found.
Determining the rate of mean-square convergence 
in the approach {\rm \cite{Rybakov3000}} is an open problem.}

\section{The Connection of Condition (\ref{july700000}) with the Concept
of Limiting Traces from the Work of G.W. Johnson and G. Kallianpur \cite{bugh3}}

Assume that $\{\phi_j(x)\}_{j=0}^{\infty}$
is an arbitrary complete orthonormal system of functions
in $L_2([t, T])$ and
$\psi_1(\tau),\ldots, \psi_k(\tau)\in L_2([t, T]).$

By analogy with Definition~2.1 (see Sect.~2.46), we define the 
$r$th limiting trace of the function $K(t_1,\ldots,t_k)\in L_2([t, T]^k)$
of type (\ref{chain200}) by the following expression
\begin{equation}
\label{2025may22}
~T^{k-2r}_{g_1,g_2,\ldots,g_{2r-1}, g_{2r}}
K (t_{q_1},\ldots,t_{q_{k-2r}})
\stackrel{\sf def}{=}
\lim\limits_{p\to\infty}T^{k-2r,\hspace{0.2mm}p}_{g_1,g_2,\ldots,g_{2r-1}, g_{2r}}
K(t_{q_1},\ldots,t_{q_{k-2r}})
\end{equation}
\noindent
in $L_2([t, T]^{k-2r})$ (here and further $L_2([t, T]^{0})\stackrel{\sf def}{=}{\bf R}$), where

\vspace{-2mm}
$$
T^{k-2r,\hspace{0.2mm}p}_{g_1,g_2,\ldots,g_{2r-1}, g_{2r}}
K(t_{q_1},\ldots,t_{q_{k-2r}})
=
$$
$$
=
\sum\limits_{j_{q_1},\ldots,j_{q_{k-2r}}=0}^p
\sum\limits_{j_{g_1}, j_{g_3},\ldots ,j_{g_{2r-1}}=0}^p
C_{j_k\ldots j_1}\biggl|_{j_{g_1}=j_{g_2},\ldots, j_{g_{2r-1}}=j_{g_{2r}}}
\cdot\hspace{0.7mm}\phi_{j_{q_1}}(t_{q_1})\ldots \phi_{j_{q_{k-2r}}}(t_{q_{k-2r}}),
$$

\vspace{3mm}
\noindent
where
$r=1,2,\ldots,[k/2]$ and 
$\{g_1,g_2,\ldots,g_{2r-1},g_{2r},
q_1,\ldots,q_{k-2r}\}=\{1,2,\ldots,k\}$ 
(see (\ref{leto5007after})),
$C_{j_k \ldots j_1}$ is the Fourier coefficient 
(\ref{2025may21}).
In addition we write
$T^{k-2r}_{g_1,g_2,\ldots,g_{2r-1}, g_{2r}}
K(t_{q_1},\ldots,t_{q_{k-2r}})
=K(t_1,\ldots,t_k)$ 
for $r=0$. 

Note that in \cite{bugh1}-\cite{bugh3} the Wiener process is scalar,
while in this book the Wiener process is a multidimensional
process with independent components.
One of the main results of work \cite{bugh3} (Theorem~5.1)
is obtained under the condition 
of existence of limiting traces (see Definition~2.1 in Sect.~2.46).

Further, we will show that the condition (\ref{july700000})
is a necessary and sufficient condition 
for the existence of limiting traces (\ref{2025may22})
(for all $r=1,2,\ldots,[k/2]$ and for all possible $g_1,g_2,\ldots,g_{2r-1},g_{2r}$ 
(see (\ref{leto5007after}))) for the case of a 
multidimensional Wiener process and $K(t_1,\ldots,t_k)$
defined by (\ref{chain200}).

Here it is also appropriate to recall the formula (\ref{after501ds1})

\vspace{-2mm}
$$
\hbox{\vtop{\offinterlineskip\halign{
\hfil#\hfil\cr
{\rm l.i.m.}\cr
$\stackrel{}{{}_{p_1,\ldots,p_k\to \infty}}$\cr
}} }\sum_{j_1=0}^{p_1}\ldots\sum_{j_k=0}^{p_k}
C_{j_k\ldots j_1}
\zeta_{j_1}^{(i_1)}\ldots \zeta_{j_k}^{(i_k)}
=J[\psi^{(k)}]_{T,t}^{(i_1\ldots i_k)}
+
$$

$$
+
\sum\limits_{r=1}^{[k/2]}
\sum_{\stackrel{(\{\{g_1, g_2\}, \ldots, 
\{g_{2r-1}, g_{2r}\}\}, \{q_1, \ldots, q_{k-2r}\})}
{{}_{\{g_1, g_2, \ldots, 
g_{2r-1}, g_{2r}, q_1, \ldots, q_{k-2r}\}=\{1, 2, \ldots, k\}}}}
\prod\limits_{s=1}^r
{\bf 1}_{\{i_{g_{{}_{2s-1}}}=~i_{g_{{}_{2s}}}\ne 0\}}\times
$$

\vspace{2mm}
\begin{equation}
\label{2025may23}
\times \hbox{\vtop{\offinterlineskip\halign{
\hfil#\hfil\cr
{\rm l.i.m.}\cr
$\stackrel{}{{}_{p_1,\ldots,p_k\to \infty}}$\cr
}} }\sum_{j_1=0}^{p_1}\ldots\sum_{j_k=0}^{p_k}
C_{j_k\ldots j_1}
\prod\limits_{s=1}^r{\bf 1}_{\{j_{g_{{}_{2s-1}}}=~j_{g_{{}_{2s}}}\}}
J'[\phi_{j_{q_1}}\ldots \phi_{j_{q_{k-2r}}}]_{T,t}^{(i_{q_1}\ldots i_{q_{k-2r}})}
\end{equation}

\vspace{1mm}
\noindent
w.~p.~1, where 
$J'[\phi_{j_{q_1}}\ldots \phi_{j_{q_{k-2r}}}]_{T,t}^{(i_{q_1}\ldots i_{q_{k-2r}})}$ is the 
multiple Wiener stochastic integral
defined by (\ref{WiI}),
$J[\psi^{(k)}]_{T,t}^{(i_1\ldots i_k)}$ is the iterated It\^{o} stochastic
integral (\ref{dsds12}),
$\{\phi_j(x)\}_{j=0}^{\infty}$
is an arbitrary complete orthonormal system of functions
in $L_2([t, T])$ and
$\psi_1(\tau),\ldots, \psi_k(\tau)\in L_2([t, T])$ (see Theorem~2.49).

The equality (\ref{2025may23}) is an analogue
of the formula (5.1) (see \cite{bugh3}, Theorem~5.1)
for the case of a multidimensional Wiener process.

We have

\vspace{-3mm}
$$
T^{k-2r,\hspace{0.2mm}p}_{g_1,g_2,\ldots,g_{2r-1}, g_{2r}}
K(t_{q_1},\ldots,t_{q_{k-2r}})
=
$$

$$
=
\sum\limits_{j_{q_1},\ldots,j_{q_{k-2r}}=0}^p
\Biggl(
\sum\limits_{j_{g_1}, j_{g_3},\ldots ,j_{g_{2r-1}}=0}^p
C_{j_k\ldots j_1}\biggl|_{j_{g_1}=j_{g_2},\ldots, j_{g_{2r-1}}=j_{g_{2r}}}-\Biggr.
$$
$$
\Biggl.-\frac{1}{2^r} \prod\limits_{l=1}^r {\bf 1}_{\{g_{2l}=g_{2l-1}+1\}}
C_{j_k \ldots j_1}\biggl|_{(j_{g_2} j_{g_1})\curvearrowright (\cdot)
\ldots (j_{g_{2r}} j_{g_{2r-1}})\curvearrowright (\cdot),
j_{g_{{}_{1}}}=~j_{g_{{}_{2}}},\ldots, j_{g_{{}_{2r-1}}}=~j_{g_{{}_{2r}}}
}\biggr.\Biggr)\times
$$

\vspace{4mm}

$$
\times\phi_{j_{q_1}}(t_{q_1})\ldots \phi_{j_{q_{k-2r}}}(t_{q_{k-2r}})+
$$

\vspace{-5mm}
$$
+\frac{1}{2^r} \prod\limits_{l=1}^r {\bf 1}_{\{g_{2l}=g_{2l-1}+1\}}
\hspace{-1.5mm}\sum\limits_{j_{q_1},\ldots,j_{q_{k-2r}}=0}^p\hspace{-1mm}
C_{j_k \ldots j_1}\biggl|_{(j_{g_2} j_{g_1})\curvearrowright (\cdot)
\ldots (j_{g_{2r}} j_{g_{2r-1}})\curvearrowright (\cdot),
j_{g_{{}_{1}}}=~j_{g_{{}_{2}}},\ldots, j_{g_{{}_{2r-1}}}=~j_{g_{{}_{2r}}}
}\biggr.
\hspace{-2mm}\times
$$

\vspace{2mm}
$$
\times\phi_{j_{q_1}}(t_{q_1})\ldots \phi_{j_{q_{k-2r}}}(t_{q_{k-2r}})\stackrel{\sf def}{=}
$$

\vspace{-2mm}
\begin{equation}
\label{2025may28}
~~~~~~~~~~\stackrel{\sf def}{=}F^{(p)}_{g_1,g_2,\ldots,g_{2r-1}, g_{2r}}(t_{q_1},\ldots,t_{q_{k-2r}})+
G^{(p)}_{g_1,g_2,\ldots,g_{2r-1}, g_{2r}}(t_{q_1},\ldots,t_{q_{k-2r}}).
\end{equation}

\vspace{3mm}

Denote
$$
C_{j_k \ldots j_1}\biggl|_{(j_{g_2} j_{g_1})\curvearrowright (\cdot)
\ldots (j_{g_{2r}} j_{g_{2r-1}})\curvearrowright (\cdot),
j_{g_{{}_{1}}}=~j_{g_{{}_{2}}},\ldots, j_{g_{{}_{2r-1}}}=~j_{g_{{}_{2r}}}
}\biggr.\stackrel{\sf def}{=}C_{j_{q_{k-2r}}\ldots j_{q_1}}^{g_1 g_2\ldots g_{2r-1} g_{2r}}.
$$ 

\vspace{2mm}

Applying transformations 
(\ref{july100003}), (\ref{july100004}) (see Sect.~2.30)
iteratively 
to $C_{j_{q_{k-2r}}\ldots j_{q_1}}^{g_1 g_2\ldots g_{2r-1} g_{2r}}$
for integrations not involving the basis functions
$\phi_{j_{q_1}},\ldots, \phi_{j_{q_{k-2r}}},$ 
we obtain
\begin{equation}
\label{2025may24}
~~~~~C_{j_{q_{k-2r}}\ldots j_{q_1}}^{g_1 g_2\ldots g_{2r-1}  g_{2r}}=
\sum\limits_{d=1}^{2^{r}}(-1)^{d-1}
\left(\hat C_{j_{q_{k-2r}}\ldots j_{q_1}}^{(d)g_1 g_2\ldots g_{2r-1} g_{2r}}
-
\bar C_{j_{q_{k-2r}}\ldots j_{q_1}}^{(d)g_1 g_2\ldots g_{2r-1} g_{2r}}
\right),
\end{equation}

\vspace{3mm}
\noindent
where some terms in the sum
$$
\sum\limits_{d=1}^{2^{r}}
$$

\noindent
can be identically equal to zero due to
the remark to (\ref{july100003}), (\ref{july100004}).

Using (\ref{2025may24}), we get
$$
G^{(p)}_{g_1,g_2,\ldots,g_{2r-1}, g_{2r}}(t_{q_1},\ldots,t_{q_{k-2r}})=
\frac{1}{2^r} \prod\limits_{l=1}^r {\bf 1}_{\{g_{2l}=g_{2l-1}+1\}}\times
$$
$$
\times
\sum\limits_{d=1}^{2^{r}}(-1)^{d-1}\Biggl(
\sum\limits_{j_{q_1},\ldots,j_{q_{k-2r}}=0}^p
\hat C_{j_{q_{k-2r}}\ldots j_{q_1}}^{(d)g_1 g_2\ldots g_{2r-1} g_{2r}}\cdot
\phi_{j_{q_1}}(t_{q_1})\ldots \phi_{j_{q_{k-2r}}}(t_{q_{k-2r}})-\Biggr.
$$
$$
\Biggl.-
\sum\limits_{j_{q_1},\ldots,j_{q_{k-2r}}=0}^p
\bar C_{j_{q_{k-2r}}\ldots j_{q_1}}^{(d)g_1 g_2\ldots g_{2r-1} g_{2r}}\cdot
\phi_{j_{q_1}}(t_{q_1})\ldots \phi_{j_{q_{k-2r}}}(t_{q_{k-2r}})\Biggr)\ \rightarrow
$$
$$
\rightarrow
\frac{1}{2^r} \prod\limits_{l=1}^r {\bf 1}_{\{g_{2l}=g_{2l-1}+1\}}\times
$$
$$
\times
\sum\limits_{d=1}^{2^{r}}(-1)^{d-1}\biggl(
\hat F^{(d)}_{g_1,g_2,\ldots,g_{2r-1}, g_{2r}}(t_{q_1},\ldots,t_{q_{k-2r}})-
\bar F^{(d)}_{g_1,g_2,\ldots,g_{2r-1}, g_{2r}}(t_{q_1},\ldots,t_{q_{k-2r}})\biggr)\stackrel{\sf def}{=}
$$

\vspace{-2mm}
\begin{equation}
\label{2025may25}
~~~~~~~~~~~\stackrel{\sf def}{=}
G_{g_1,g_2,\ldots,g_{2r-1}, g_{2r}}(t_{q_1},\ldots,t_{q_{k-2r}})\ \ \hbox{if}\ \ p\to\infty\ \ 
(\hbox{in}\ L_2([t, T]^{k-2r})), 
\end{equation}

\vspace{3mm}
\noindent
where $d=1,\ldots,2^{r},$ $G_{g_1,g_2,\ldots,g_{2r-1}, g_{2r}}(t_{q_1},\ldots,t_{q_{k-2r}}),$ 
$\hat F^{(d)}_{g_1,g_2,\ldots,g_{2r-1}, g_{2r}}(t_{q_1},\ldots,t_{q_{k-2r}}),$
$\bar F^{(d)}_{g_1,g_2,\ldots,g_{2r-1}, g_{2r}}(t_{q_1},\ldots,t_{q_{k-2r}})\in L_2([t, T]^{k-2r}).$

Futhermore, 

\vspace{-4mm}
$$
\bigl\Vert F^{(p)}_{g_1,g_2,\ldots,g_{2r-1}, g_{2r}} \bigr\Vert_{L_2([t, T]^{k-2r})}^2=
$$

\vspace{-2mm}
$$
=\int\limits_{[t, T]^{k-2r}}
\Biggl(
\sum\limits_{j_{q_1},\ldots,j_{q_{k-2r}}=0}^p
\Biggl(
\sum\limits_{j_{g_1}, j_{g_3},\ldots ,j_{g_{2r-1}}=0}^p
C_{j_k\ldots j_1}\biggl|_{j_{g_1}=j_{g_2},\ldots, j_{g_{2r-1}}=j_{g_{2r}}}-\Biggr.\Biggr.
$$

\vspace{-1mm}
$$
\Biggl.-\frac{1}{2^r} \prod\limits_{l=1}^r {\bf 1}_{\{g_{2l}=g_{2l-1}+1\}}
C_{j_k \ldots j_1}\biggl|_{(j_{g_2} j_{g_1})\curvearrowright (\cdot)
\ldots (j_{g_{2r}} j_{g_{2r-1}})\curvearrowright (\cdot),
j_{g_{{}_{1}}}=~j_{g_{{}_{2}}},\ldots, j_{g_{{}_{2r-1}}}=~j_{g_{{}_{2r}}}
}\biggr.\Biggr)\times
$$

\vspace{2mm}
$$
\Biggl.\times\phi_{j_{q_1}}(t_{q_1})\ldots \phi_{j_{q_{k-2r}}}(t_{q_{k-2r}})\Biggr)^2
dt_{q_1}\ldots dt_{q_{k-2r}}=
$$

\vspace{3mm}
$$
=
\sum\limits_{j_{q_1},\ldots,j_{q_{k-2r}}=0}^p
\Biggl(
\sum\limits_{j_{g_1}, j_{g_3},\ldots ,j_{g_{2r-1}}=0}^p
C_{j_k\ldots j_1}\biggl|_{j_{g_1}=j_{g_2},\ldots, j_{g_{2r-1}}=j_{g_{2r}}}-\Biggr.\Biggr.
$$

\vspace{-1mm}
$$
\Biggl.-\frac{1}{2^r} \prod\limits_{l=1}^r {\bf 1}_{\{g_{2l}=g_{2l-1}+1\}}
C_{j_k \ldots j_1}\biggl|_{(j_{g_2} j_{g_1})\curvearrowright (\cdot)
\ldots (j_{g_{2r}} j_{g_{2r-1}})\curvearrowright (\cdot),
j_{g_{{}_{1}}}=~j_{g_{{}_{2}}},\ldots, j_{g_{{}_{2r-1}}}=~j_{g_{{}_{2r}}}
}\biggr.\Biggr)^2.
$$

\vspace{5mm}

This means that the condition (\ref{july700000}) is equivalent to 
\begin{equation}
\label{2025may26}
\lim\limits_{p\to\infty}
\bigl\Vert F^{(p)}_{g_1,g_2,\ldots,g_{2r-1}, g_{2r}} \bigr\Vert_{L_2([t, T]^{k-2r})}^2=0.
\end{equation}

\vspace{1mm}

Suppose that the condition (\ref{july700000}) (or (\ref{2025may26}))
is fulfilled. 
Applying (\ref{2025may28}), (\ref{2025may25}) and (\ref{2025may26}), we obtain

\vspace{-1mm}
$$
\bigl\Vert T^{k-2r,\hspace{0.2mm}p}_{g_1,g_2,\ldots,g_{2r-1}, g_{2r}}
K- G_{g_1,g_2,\ldots,g_{2r-1}, g_{2r}} \bigr\Vert_{L_2([t, T]^{k-2r})}=
$$

\vspace{0.5mm}
$$
=\bigl\Vert F^{(p)}_{g_1,g_2,\ldots,g_{2r-1}, g_{2r}} + 
G^{(p)}_{g_1,g_2,\ldots,g_{2r-1}, g_{2r}} - G_{g_1,g_2,\ldots,g_{2r-1}, g_{2r}} 
\bigr\Vert_{L_2([t, T]^{k-2r})}\le
$$

\vspace{0.5mm}
$$
\le\bigl\Vert F^{(p)}_{g_1,g_2,\ldots,g_{2r-1}, g_{2r}} \bigr\Vert_{L_2([t, T]^{k-2r})} +
$$

\vspace{0.5mm}
$$
+ 
\bigl\Vert G^{(p)}_{g_1,g_2,\ldots,g_{2r-1}, g_{2r}} - 
G_{g_1,g_2,\ldots,g_{2r-1}, g_{2r}} \bigr\Vert_{L_2([t, T]^{k-2r})}\ \to 0
$$

\vspace{2mm}
\noindent
if $p\to\infty.$

Thus, the limiting trace $T^{k-2r}_{g_1,g_2,\ldots,g_{2r-1}, g_{2r}}
K(t_{q_1},\ldots,t_{q_{k-2r}})$
exists under the condition (\ref{july700000})
for all $r=1,2,\ldots,$ $[k/2]$ and 
for all possible $g_1,g_2,\ldots,g_{2r-1},g_{2r}$ 
(see (\ref{leto5007after})), i.e.

\vspace{-2mm}
$$
T^{k-2r}_{g_1,g_2,\ldots,g_{2r-1}, g_{2r}}K(t_{q_1},\ldots,t_{q_{k-2r}})
=
\lim\limits_{p\to\infty}T^{k-2r,\hspace{0.2mm}p}_{g_1,g_2,\ldots,g_{2r-1}, g_{2r}}
K(t_{q_1},\ldots,t_{q_{k-2r}})
=
$$

\begin{equation}
\label{october20251}
=G_{g_1,g_2,\ldots,g_{2r-1}, g_{2r}}(t_{q_1},\ldots,t_{q_{k-2r}})
\end{equation}

\vspace{3mm}
\noindent
in $L_2([t, T]^{k-2r}),$ where 
$G_{g_1,g_2,\ldots,g_{2r-1}, g_{2r}}(t_{q_1},\ldots,t_{q_{k-2r}})$ is defined by (\ref{2025may25}).

Further, we assume that the condition (\ref{october20251}) is fulfilled.
Using (\ref{october20251}) and (\ref{2025may28}), (\ref{2025may25}), we obtain 

\vspace{-2mm}
$$
\bigl\Vert F^{(p)}_{g_1,g_2,\ldots,g_{2r-1}, g_{2r}}\bigr\Vert_{L_2([t, T]^{k-2r})}=
$$

\vspace{0.5mm}
$$
=\bigl\Vert T^{k-2r,\hspace{0.2mm}p}_{g_1,g_2,\ldots,g_{2r-1}, g_{2r}}
K- G^{(p)}_{g_1,g_2,\ldots,g_{2r-1}, g_{2r}} \bigr\Vert_{L_2([t, T]^{k-2r})}=
$$

\vspace{0.5mm}
$$
=\bigl\Vert T^{k-2r,\hspace{0.2mm}p}_{g_1,g_2,\ldots,g_{2r-1}, g_{2r}}
K- G_{g_1,g_2,\ldots,g_{2r-1}, g_{2r}}+
G_{g_1,g_2,\ldots,g_{2r-1}, g_{2r}}
-
G^{(p)}_{g_1,g_2,\ldots,g_{2r-1}, g_{2r}} \bigr\Vert_{L_2([t, T]^{k-2r})}
$$

\newpage
\noindent
$$
\le\bigl\Vert T^{k-2r,\hspace{0.2mm}p}_{g_1,g_2,\ldots,g_{2r-1}, g_{2r}}
K- G_{g_1,g_2,\ldots,g_{2r-1}, g_{2r}} \bigr\Vert_{L_2([t, T]^{k-2r})} + 
$$

\vspace{-2mm}
$$
+
\bigl\Vert G_{g_1,g_2,\ldots,g_{2r-1}, g_{2r}}
-
G^{(p)}_{g_1,g_2,\ldots,g_{2r-1}, g_{2r}} \bigr\Vert_{L_2([t, T]^{k-2r})}\ \to 0
$$

\vspace{2mm}
\noindent
if $p\to\infty.$

Thus, the condition (\ref{2025may26}) is satisfied
under the condition (\ref{october20251}).
Note that if $T^{k-2r}_{g_1,g_2,\ldots,g_{2r-1}, g_{2r}}K(t_{q_1},\ldots,t_{q_{k-2r}})$
exists, then it will be equal to 
$G_{g_1,g_2,\ldots,g_{2r-1}, g_{2r}}(t_{q_1},\ldots,t_{q_{k-2r}}),$ which is defined by (\ref{2025may25})
(see Sect.~2.49 for details).

\section{New Representations of the Hu--Meyer Formulas for the Case
of a Multidimensional Wiener Process. Connection with Theorem~2.49}

Suppose that $\{\phi_j(x)\}_{j=0}^{\infty}$
is an arbitrary complete orthonormal system of functions
in $L_2([t, T])$ and
$\Phi(t_1,\ldots,t_k)\in L_2([t, T]^k).$

Let us generalize the definition of the limiting trace
from the previous section.

We define the 
$r$th limiting trace of the function $\Phi(t_1,\ldots,t_k)\in L_2([t, T]^k)$
by the following expression

\vspace{-3mm}
\begin{equation}
\label{novem2026xxx1122}
~~T^{k-2r}_{g_1,g_2,\ldots,g_{2r-1}, g_{2r}}\Phi (t_{q_1},\ldots,t_{q_{k-2r}})\stackrel{\sf def}{=}
\lim\limits_{p\to\infty}T^{k-2r,\hspace{0.2mm}p}_{g_1,g_2,\ldots,g_{2r-1}, g_{2r}}
\Phi (t_{q_1},\ldots,t_{q_{k-2r}})
\end{equation}

\vspace{1mm}
\noindent
in $L_2([t, T]^{k-2r})$ (here and further $L_2([t, T]^{0})\stackrel{\sf def}{=}{\bf R}$), where

\vspace{-3mm}
$$
T^{k-2r,\hspace{0.2mm}p}_{g_1,g_2,\ldots,g_{2r-1}, g_{2r}}\Phi (t_{q_1},\ldots,t_{q_{k-2r}})=
$$

\vspace{-5mm}
\begin{equation}
\label{trace11111111}
=
\sum\limits_{j_{q_1},\ldots,j_{q_{k-2r}}=0}^p
\sum\limits_{j_{g_1}, j_{g_3},\ldots ,j_{g_{2r-1}}=0}^p
C_{j_k\ldots j_1}\biggl|_{j_{g_1}=j_{g_2},\ldots, j_{g_{2r-1}}=j_{g_{2r}}}
\cdot\hspace{0.7mm}\phi_{j_{q_1}}(t_{q_1})\ldots \phi_{j_{q_{k-2r}}}(t_{q_{k-2r}}),
\end{equation}

\vspace{2mm}
\noindent
where
$r=1,2,\ldots,[k/2]$ and 
$\{g_1,g_2,\ldots,g_{2r-1},g_{2r},
q_1,\ldots,q_{k-2r}\}=\{1,2,\ldots,k\}$ 
(see (\ref{leto5007after})),

\vspace{-2mm}
\begin{equation}
\label{2025may30}
C_{j_k\ldots j_1}=\int\limits_{[t,T]^k}
\Phi(t_1,\ldots,t_k)\prod_{l=1}^{k}\phi_{j_l}(t_l)dt_1\ldots dt_k
\end{equation}

\vspace{3mm}
\noindent
is the Fourier coefficient.
Also we will write
$T^{k-2r}_{g_1,g_2,\ldots,g_{2r-1}, g_{2r}}\Phi (t_{q_1},\ldots,t_{q_{k-2r}})=\Phi(t_1,\ldots,t_k)$ 
and 
$T^{k-2r,\hspace{0.2mm}p}_{g_1,g_2,\ldots,g_{2r-1}, g_{2r}}
\Phi (t_{q_1},\ldots,t_{q_{k-2r}})=\Phi_p(t_1,\ldots,t_k)$
for $r=0$, where
$$
\Phi_p(t_1,\ldots,t_k)=
\sum\limits_{j_{1},\ldots,j_{k}=0}^p
C_{j_k\ldots j_1}
\phi_{j_{1}}(t_{1})\ldots \phi_{j_k}(t_k).
$$

\vspace{2mm}

Let us consider a variant of the formula (\ref{after8xxds1})
for the case $p_1=\ldots=p_k=p$ and replace 
$K(t_1,\ldots,t_k)$ with $\Phi(t_1,\ldots,t_k)$ in it. Thus, we have

\vspace{-1mm}
$$
\sum_{j_1,\ldots,j_k=0}^{p}
C_{j_k\ldots j_1}
\prod_{l=1}^k \zeta_{j_l}^{(i_l)}=
\sum_{j_1,\ldots,j_k=0}^{p}
C_{j_k\ldots j_1}
J'[\phi_{j_1}\ldots \phi_{j_k}]_{T,t}^{(i_1\ldots i_k)}+
$$

\vspace{2mm}
$$
+
\sum\limits_{r=1}^{[k/2]}
\sum_{\stackrel{(\{\{g_1, g_2\}, \ldots, 
\{g_{2r-1}, g_{2r}\}\}, \{q_1, \ldots, q_{k-2r}\})}
{{}_{\{g_1, g_2, \ldots, 
g_{2r-1}, g_{2r}, q_1, \ldots, q_{k-2r}\}=\{1, 2, \ldots, k\}}}}
\prod\limits_{s=1}^r
{\bf 1}_{\{i_{g_{{}_{2s-1}}}=~i_{g_{{}_{2s}}}\ne 0\}}
\sum_{j_1,\ldots,j_k=0}^{p}
C_{j_k\ldots j_1}\times
$$

\begin{equation}
\label{2025may31}
~~~~~~~~~~\times{\bf 1}_{\{j_{g_{{}_{2s-1}}}=~j_{g_{{}_{2s}}}\}}
J'[\phi_{j_{q_1}}\ldots \phi_{j_{q_{k-2r}}}]_{T,t}^{(i_{q_1}\ldots i_{q_{k-2r}})}\ \ \ \hbox{w.~p.~1,}
\end{equation}

\vspace{3mm}
\noindent
where $k\ge 2,$ $J'[\phi_{j_1}\ldots \phi_{j_k}]_{T,t}^{(i_1\ldots i_k)},$
$J'[\phi_{j_{q_1}}\ldots \phi_{j_{q_{k-2r}}}]_{T,t}^{(i_{q_1}\ldots i_{q_{k-2r}})}$
are multiple Wie\-ner sto\-chas\-tic integrals 
(see (\ref{WiI})), $J'[\phi_{j_{q_1}}\ldots \phi_{j_{q_{k-2r}}}]_{T,t}^{(i_{q_1}\ldots i_{q_{k-2r}})}
\stackrel{\sf def}{=}1$ for $k=2r$ and $C_{j_k\ldots j_1}$
is defined by (\ref{2025may30}).

Note that sometimes the multiple Stratonovich stochastic integral for
$\Phi(t_1,\ldots,t_k)\in L_2([t, T]^k)$
is defined as the limit in probability as $p\to\infty$ of the following expression
(the case of a scalar standard Wiener process)
$$
\sum_{j_1,\ldots,j_k=0}^{p}
C_{j_k\ldots j_1}
\prod_{l=1}^k \zeta_{j_l}^{(i)},
$$

\noindent
where $i_1=\ldots=i_k=i$ (see, for example,  Definition~5.9 in \cite{HuHu}).

By analogy with \cite{HuHu}, we define the multiple Stratonovich
stochastic integral
for $\Phi(t_1,\ldots,t_k)\in L_2([t, T]^k)$
(the case of a multidimensional Wiener process)
as the following mean-square limit

\vspace{-2mm}
\begin{equation}
\label{novem2026xxx4}
\hat J^{S}[\Phi]_{T,t}^{(i_1\ldots i_k)}=
\hbox{\vtop{\offinterlineskip\halign{
\hfil#\hfil\cr
{\rm l.i.m.}\cr
$\stackrel{}{{}_{p\to \infty}}$\cr
}} }
\hat J^{S}_p[\Phi]_{T,t}^{(i_1\ldots i_k)},
\end{equation}

\noindent
where
\begin{equation}
\label{novem2026xxx5}
\hat J^{S}_p[\Phi]_{T,t}^{(i_1\ldots i_k)}=
\sum_{j_1,\ldots,j_k=0}^{p}
C_{j_k\ldots j_1}
\prod_{l=1}^k \zeta_{j_l}^{(i_l)},
\end{equation}

\vspace{2mm}
\noindent
$C_{j_k\ldots j_1}$ is the Fourier coefficient
defined by (\ref{2025may30}).

Let us rewrite (\ref{2025may31}) in the form

\vspace{-2mm}
$$
\hat J^{S}_p[\Phi]_{T,t}^{(i_1\ldots i_k)}
=
J'[\Phi_p]_{T,t}^{(i_1\ldots i_k)}+
$$

\vspace{-2mm}
$$
+
\sum\limits_{r=1}^{[k/2]}
\sum_{\stackrel{(\{\{g_1, g_2\}, \ldots, 
\{g_{2r-1}, g_{2r}\}\}, \{q_1, \ldots, q_{k-2r}\})}
{{}_{\{g_1, g_2, \ldots, 
g_{2r-1}, g_{2r}, q_1, \ldots, q_{k-2r}\}=\{1, 2, \ldots, k\}}}}
\prod\limits_{s=1}^r
{\bf 1}_{\{i_{g_{{}_{2s-1}}}=~i_{g_{{}_{2s}}}\ne 0\}}\times
$$

\vspace{5mm}
\begin{equation}
\label{2026may1}
~~\times
J'\hspace{-1mm}\left[T^{k-2r,\hspace{0.2mm}p}_{g_1,g_2,\ldots,g_{2r-1}, 
g_{2r}}\Phi\right]_{T,t}^{(i_{q_1}\ldots i_{q_{k-2r}})}\ \ \ \hbox{w.~p.~1,}
\end{equation}

\vspace{5mm}
\noindent
where $J'\hspace{-1mm}\left[T^{k-2r,\hspace{0.2mm}p}_{g_1,g_2,\ldots,g_{2r-1}, 
g_{2r}}\Phi\right]_{T,t}^{(i_{q_1}\ldots i_{q_{k-2r}})}\stackrel{\sf def}{=}
T^{k-2r,\hspace{0.2mm}p}_{g_1,g_2,\ldots,g_{2r-1}, 
g_{2r}}\Phi$ for $k=2r.$

Assume that all limiting traces
$T^{k-2r}_{g_1,g_2,\ldots,g_{2r-1}, g_{2r}}\Phi (t_{q_1},\ldots,t_{q_{k-2r}})$
(for all $r=1,2,\ldots,$ $[k/2]$ and for all possible
$g_1,g_2,\ldots,g_{2r-1},g_{2r}$ (see {\rm (\ref{leto5007after})))
exist.

Recall the well-known property
of the multiple Wiener stochastic integral (\ref{WiI})

\vspace{-3mm}
\begin{equation}
\label{2026may2}
~~~~~~~~{\sf M}\left\{\left(J'[\Phi]_{T,t}^{(i_1\ldots i_k)}\right)^2\right\}\le
C_k
\int\limits_{[t,T]^k}
\Phi^2(t_1,\ldots,t_k)dt_1\ldots dt_k,
\end{equation}

\vspace{2mm}
\noindent
where $\Phi(t_1,\ldots,t_k)\in L_2([t, T]^k)$ and $C_k$ is a constant.

Using (\ref{2026may1}), (\ref{2026may2}) and the existence of limiting traces, we have

\vspace{-3mm}
$$
{\sf M}\Biggl\{\Biggl(\hat J^{S}_p[\Phi]_{T,t}^{(i_1\ldots i_k)}
-
J'[\Phi]_{T,t}^{(i_1\ldots i_k)}-\Biggr.\Biggr.
$$

\vspace{-2mm}

$$
-
\sum\limits_{r=1}^{[k/2]}
\sum_{\stackrel{(\{\{g_1, g_2\}, \ldots, 
\{g_{2r-1}, g_{2r}\}\}, \{q_1, \ldots, q_{k-2r}\})}
{{}_{\{g_1, g_2, \ldots, 
g_{2r-1}, g_{2r}, q_1, \ldots, q_{k-2r}\}=\{1, 2, \ldots, k\}}}}
\prod\limits_{s=1}^r
{\bf 1}_{\{i_{g_{{}_{2s-1}}}=~i_{g_{{}_{2s}}}\ne 0\}}\times
$$
$$
\Biggl.\Biggl.\times
J'\hspace{-1mm}\left[T^{k-2r}_{g_1,g_2,\ldots,g_{2r-1}, 
g_{2r}}\Phi\right]_{T,t}^{(i_{q_1}\ldots i_{q_{k-2r}})}\Biggr)^2\Biggr\}=
$$

\vspace{2mm}
$$
={\sf M}\Biggl\{\Biggl(
J'[\Phi_p-\Phi]_{T,t}^{(i_1\ldots i_k)}+\Biggr.\Biggr.
$$

\vspace{-2mm}
$$
+
\sum\limits_{r=1}^{[k/2]}
\sum_{\stackrel{(\{\{g_1, g_2\}, \ldots, 
\{g_{2r-1}, g_{2r}\}\}, \{q_1, \ldots, q_{k-2r}\})}
{{}_{\{g_1, g_2, \ldots, 
g_{2r-1}, g_{2r}, q_1, \ldots, q_{k-2r}\}=\{1, 2, \ldots, k\}}}}
\prod\limits_{s=1}^r
{\bf 1}_{\{i_{g_{{}_{2s-1}}}=~i_{g_{{}_{2s}}}\ne 0\}}\times
$$

$$
\Biggl.\Biggl.\times
J'\hspace{-1mm}\left[T^{k-2r,\hspace{0.2mm}p}_{g_1,g_2,\ldots,g_{2r-1}, 
g_{2r}}\Phi-
T^{k-2r}_{g_1,g_2,\ldots,g_{2r-1}, 
g_{2r}}\Phi
\right]_{T,t}^{(i_{q_1}\ldots i_{q_{k-2r}})}\Biggr)^2\Biggr\}\le
$$

\vspace{1mm}

$$
\le C_k'\Biggl(
\bigl\Vert\Phi_p-\Phi\bigr\Vert_{L_2([t, T]^k)}^2+\Biggr.
$$

\vspace{-3mm}
$$
+
\sum\limits_{r=1}^{[k/2]}
\sum_{\stackrel{(\{\{g_1, g_2\}, \ldots, 
\{g_{2r-1}, g_{2r}\}\}, \{q_1, \ldots, q_{k-2r}\})}
{{}_{\{g_1, g_2, \ldots, 
g_{2r-1}, g_{2r}, q_1, \ldots, q_{k-2r}\}=\{1, 2, \ldots, k\}}}}
\prod\limits_{s=1}^r
{\bf 1}_{\{i_{g_{{}_{2s-1}}}=~i_{g_{{}_{2s}}}\ne 0\}}\times
$$

\begin{equation}
\label{2026may5}
~~~~~~~~~~\Biggl.\times
\left\Vert T^{k-2r,\hspace{0.2mm}p}_{g_1,g_2,\ldots,g_{2r-1}, 
g_{2r}}\Phi-
T^{k-2r}_{g_1,g_2,\ldots,g_{2r-1}, 
g_{2r}}\Phi\right\Vert_{L_2([t, T]^{k-2r})}^2\Biggr)\ \to\ 0
\end{equation}

\vspace{2mm}
\noindent
if $p\to\infty,$ 
where $J'\hspace{-1mm}\left[T^{k-2r}_{g_1,g_2,\ldots,g_{2r-1}, 
g_{2r}}\Phi\right]_{T,t}^{(i_{q_1}\ldots i_{q_{k-2r}})}\stackrel{\sf def}{=}
T^{k-2r}_{g_1,g_2,\ldots,g_{2r-1}, 
g_{2r}}\Phi$ for $k=2r,$
and $C_k'$ is a constant.

Applying (\ref{2026may5}), we obtain the following 
new representation of the Hu--Meyer formula for the case
of a multidimensional Wiener process

\vspace{-2mm}
$$
\hat J^{S}[\Phi]_{T,t}^{(i_1\ldots i_k)}
=
J'[\Phi]_{T,t}^{(i_1\ldots i_k)}+
$$

\vspace{-3mm}
$$
+
\sum\limits_{r=1}^{[k/2]}
\sum_{\stackrel{(\{\{g_1, g_2\}, \ldots, 
\{g_{2r-1}, g_{2r}\}\}, \{q_1, \ldots, q_{k-2r}\})}
{{}_{\{g_1, g_2, \ldots, 
g_{2r-1}, g_{2r}, q_1, \ldots, q_{k-2r}\}=\{1, 2, \ldots, k\}}}}
\prod\limits_{s=1}^r
{\bf 1}_{\{i_{g_{{}_{2s-1}}}=~i_{g_{{}_{2s}}}\ne 0\}}\times
$$

\vspace{4mm}
\begin{equation}
\label{2026may10}
\times
J'\hspace{-1mm}\left[T^{k-2r}_{g_1,g_2,\ldots,g_{2r-1}, 
g_{2r}}\Phi\right]_{T,t}^{(i_{q_1}\ldots i_{q_{k-2r}})}\ \ \ \hbox{w.~p.~1.}
\end{equation}

\vspace{4mm}

The equality (\ref{2026may10}) is consistent with Theorem~5.1 \cite{bugh3}
(the case of a scalar Wiener process).

Further, let us obtain the inverse version of the Hu--Meyer
formula (\ref{2026may10}), i.e. a formula expressing
the multiple Wiener stochastic integral through the sum
of multiple Stratonovich stochastic integrals.

Let us give the following definition of the limiting trace.
We define the 
$r$th limiting trace of the function $\Phi(t_1,\ldots,t_k)\in L_2([t, T]^k)$
by the following expression
\begin{equation}
\label{soglas123sd}
~~\tilde T^{k-2r}_{g_1,g_2,\ldots,g_{2r-1}, g_{2r}}
\Phi (t_{q_1},\ldots,t_{q_{k-2r}})\stackrel{\sf def}{=}
\lim\limits_{p\to\infty}\tilde T^{k-2r,\hspace{0.2mm}p}_{g_1,g_2,\ldots,g_{2r-1}, g_{2r}}
\Phi (t_{q_1},\ldots,t_{q_{k-2r}})
\end{equation}

\noindent
in $L_2([t, T]^{k-2r})$, where

\vspace{-3mm}
$$
\tilde T^{k-2r,\hspace{0.2mm}p}_{g_1,g_2,\ldots,g_{2r-1}, g_{2r}}\Phi (t_{q_1},\ldots,t_{q_{k-2r}})=
$$

\vspace{-4mm}
$$
=\sum\limits_{j_{q_1},\ldots,j_{q_{k-2r}}=0}^p
\sum\limits_{j_{g_1}, j_{g_3},\ldots ,j_{g_{2r-1}}=0}^{\infty}
C_{j_k\ldots j_1}\biggl|_{j_{g_1}=j_{g_2},\ldots, j_{g_{2r-1}}=j_{g_{2r}}}
\cdot\hspace{0.7mm} \phi_{j_{q_1}}(t_{q_1})\ldots \phi_{j_{q_{k-2r}}}(t_{q_{k-2r}})=
$$

\vspace{-2mm}
\begin{equation}
\label{october202620}
~~~~~~=
\sum\limits_{j_{q_1},\ldots,j_{q_{k-2r}}=0}^p
\tilde C_{j_{q_{k-2r}}\ldots j_{q_1}}
\cdot\hspace{0.7mm}\phi_{j_{q_1}}(t_{q_1})\ldots \phi_{j_{q_{k-2r}}}(t_{q_{k-2r}}),
\end{equation}

\vspace{2mm}
\noindent
where
$$
\tilde C_{j_{q_{k-2r}}\ldots j_{q_1}}=
\sum\limits_{j_{g_1}, j_{g_3},\ldots ,j_{g_{2r-1}}=0}^{\infty}
C_{j_k\ldots j_1}\biggl|_{j_{g_1}=j_{g_2},\ldots, j_{g_{2r-1}}=j_{g_{2r}}}=
$$

\vspace{-1mm}
\begin{equation}
\label{october202621}
=\lim\limits_{p\to\infty}\sum\limits_{j_{g_1}, j_{g_3},\ldots ,j_{g_{2r-1}}=0}^{p}
C_{j_k\ldots j_1}\biggl|_{j_{g_1}=j_{g_2},\ldots, j_{g_{2r-1}}=j_{g_{2r}}},
\end{equation}

\vspace{3mm}
\noindent
where $r=1,2,\ldots,[k/2],$
$\{g_1,g_2,\ldots,g_{2r-1},g_{2r},
q_1,\ldots,q_{k-2r}\}=\{1,2,\ldots,k\}$ 
(see (\ref{leto5007after})),
$C_{j_k\ldots j_1}$ has the form
(\ref{2025may30}), 
$\tilde T^{k-2r}_{g_1,g_2,\ldots,g_{2r-1}, g_{2r}}\Phi (t_{q_1},\ldots,t_{q_{k-2r}})=\Phi(t_1,\ldots,t_k)$ 
and 
$\tilde T^{k-2r,\hspace{0.2mm}p}_{g_1,g_2,\ldots,g_{2r-1}, g_{2r}}\Phi (t_{q_1},\ldots,t_{q_{k-2r}})
=\Phi_p(t_1,\ldots,t_k)$
for $r=0$, where

\vspace{-3mm}
$$
\Phi_p(t_1,\ldots,t_k)=
\sum\limits_{j_{1},\ldots,j_{k}=0}^p
C_{j_k\ldots j_1}
\phi_{j_{1}}(t_{1})\ldots \phi_{j_k}(t_k).
$$

\vspace{2mm}

Suppose that all limiting traces
$T^{k-2r}_{g_1,g_2,\ldots,g_{2r-1}, g_{2r}}\Phi (t_{q_1},\ldots,t_{q_{k-2r}})$
exist
for all $r=1,2,\ldots,$ $[k/2]$ and for all possible
$g_1,g_2,$ $\ldots,g_{2r-1},g_{2r}$ (see {\rm (\ref{leto5007after})).

Let us fix
$j'_{q_1}, \ldots, j'_{q_{k-2r}}$
Applying 
the inequality 
of Cauchy--Bu\-ny\-a\-kov\-sky, we obtain
$$
\lim\limits_{p\to\infty}
\biggl\langle T^{k-2r,\hspace{0.2mm}p}_{g_1,g_2,\ldots,g_{2r-1}, g_{2r}}
\Phi (t_{q_1},\ldots,t_{q_{k-2r}})\ ,\  
\phi_{j'_{q_1}}(t_{q_1})\ldots \phi_{j'_{q_{k-2r}}}(t_{q_{k-2r}}) \biggr\rangle_{L_2([t, T]^{k-2r})}=
$$
$$
=
\biggl\langle T^{k-2r}_{g_1,g_2,\ldots,g_{2r-1}, g_{2r}}\Phi (t_{q_1},\ldots,t_{q_{k-2r}})\ ,\  
\phi_{j'_{q_1}}(t_{q_1})\ldots \phi_{j'_{q_{k-2r}}}(t_{q_{k-2r}}) \biggr\rangle_{L_2([t, T]^{k-2r})},
$$

\vspace{2mm}
\noindent
where 
$T^{k-2r,\hspace{0.2mm}p}_{g_1,g_2,\ldots,g_{2r-1}, g_{2r}}\Phi (t_{q_1},\ldots,t_{q_{k-2r}})$
is defined by (\ref{trace11111111}).

Let $p\ge \max\left\{j'_{q_1}, \ldots, j'_{q_{k-2r}}\right\}.$ Then
$$
\biggl\langle T^{k-2r,\hspace{0.2mm}p}_{g_1,g_2,\ldots,g_{2r-1}, g_{2r}}
\Phi (t_{q_1},\ldots,t_{q_{k-2r}})\ ,\  
\phi_{j'_{q_1}}(t_{q_1})\ldots \phi_{j'_{q_{k-2r}}}(t_{q_{k-2r}}) \biggr\rangle_{L_2([t, T]^{k-2r})}=
$$
$$
=\sum\limits_{j_{q_1},\ldots,j_{q_{k-2r}}=0}^p
\sum\limits_{j_{g_1}, j_{g_3},\ldots ,j_{g_{2r-1}}=0}^p
C_{j_k\ldots j_1}\biggl|_{j_{g_1}=j_{g_2},\ldots, j_{g_{2r-1}}=j_{g_{2r}}}\times
$$

\vspace{-1mm}
$$
\times
\biggl\langle\phi_{j_{q_1}}(t_{q_1})\ldots \phi_{j_{q_{k-2r}}}(t_{q_{k-2r}})\ ,\ 
\phi_{j'_{q_1}}(t_{q_1})\ldots \phi_{j'_{q_{k-2r}}}(t_{q_{k-2r}}) \biggr\rangle_{L_2([t, T]^{k-2r})}=
$$

\vspace{-2mm}
$$
=
\sum\limits_{j_{g_1}, j_{g_3},\ldots ,j_{g_{2r-1}}=0}^{p}
C_{j_k\ldots j_1}\biggl|_{j_{g_1}=j_{g_2},\ldots, j_{g_{2r-1}}=j_{g_{2r}},
j_{q_1}=j'_{q_1},\ldots, j_{q_{k-2r}}=j'_{q_{k-2r}}}.
$$

\vspace{4mm}

This means that
$$
\lim\limits_{p\to\infty}\sum\limits_{j_{g_1}, j_{g_3},\ldots ,j_{g_{2r-1}}=0}^{p}
C_{j_k\ldots j_1}\biggl|_{j_{g_1}=j_{g_2},\ldots, j_{g_{2r-1}}=j_{g_{2r}}}=
$$
$$
=\sum\limits_{j_{g_1}, j_{g_3},\ldots ,j_{g_{2r-1}}=0}^{\infty}
C_{j_k\ldots j_1}\biggl|_{j_{g_1}=j_{g_2},\ldots, j_{g_{2r-1}}=j_{g_{2r}}}
$$

\vspace{3mm}
\noindent
is a Fourier coefficient of
$T^{k-2r}_{g_1,g_2,\ldots,g_{2r-1}, g_{2r}}\Phi (t_{q_1},\ldots,t_{q_{k-2r}}),$
i.e. 

\vspace{-3mm}
\begin{equation}
\label{req1234}
~~~~~~~~~T^{k-2r}_{g_1,g_2,\ldots,g_{2r-1}, g_{2r}}\Phi (t_{q_1},\ldots,t_{q_{k-2r}})=
\tilde T^{k-2r}_{g_1,g_2,\ldots,g_{2r-1}, g_{2r}}\Phi (t_{q_1},\ldots,t_{q_{k-2r}})
\end{equation}

\vspace{3.5mm}
\noindent
in $L_2([t, T]^{k-2r})$.
Then, using the equality (\ref{req1234})
and (\ref{2026may10}), we obtain

\newpage
\noindent
$$
J'[\Phi]_{T,t}^{(i_1\ldots i_k)}=\hat J^{S}[\Phi]_{T,t}^{(i_1\ldots i_k)}-
$$

\vspace{-1mm}
$$
-
\sum\limits_{r=1}^{[k/2]}
\sum_{\stackrel{(\{\{g_1, g_2\}, \ldots, 
\{g_{2r-1}, g_{2r}\}\}, \{q_1, \ldots, q_{k-2r}\})}
{{}_{\{g_1, g_2, \ldots, 
g_{2r-1}, g_{2r}, q_1, \ldots, q_{k-2r}\}=\{1, 2, \ldots, k\}}}}
\prod\limits_{s=1}^r
{\bf 1}_{\{i_{g_{{}_{2s-1}}}=~i_{g_{{}_{2s}}}\ne 0\}}\times
$$

\vspace{5mm}
\begin{equation}
\label{2026may10sss}
\times
J'\hspace{-1mm}\left[\tilde T^{k-2r}_{g_1,g_2,\ldots,g_{2r-1}, 
g_{2r}}\Phi\right]_{T,t}^{(i_{q_1}\ldots i_{q_{k-2r}})}\ \ \ \hbox{w.~p.~1.}
\end{equation}

\vspace{5mm}

Since $\tilde T^{k-2r}_{g_1,g_2,\ldots,g_{2r-1}, 
g_{2r}}\Phi(t_{q_1},\ldots,t_{q_{k-2r}})\in L_2([t,T]^{k-2r})$,
then the multiple Stra\-to\-no\-vich stochastic integral 
$\hat J^{S}[\tilde T^{k-2r}_{g_1,g_2,\ldots,g_{2r-1}, g_{2r}}\Phi]_{T,t}^{(i_{q_1}\ldots i_{q_{k-2r}})}$
is defined according to (\ref{novem2026xxx4}), (\ref{novem2026xxx5}).
If we replace $\Phi$ with $\tilde T^{k-2r}_{g_1,g_2,\ldots,g_{2r-1}, g_{2r}}\Phi$
in (\ref{2026may10sss}), then the integral 
$\hat J^{S}[\tilde T^{k-2r}_{g_1,g_2,\ldots,g_{2r-1}, g_{2r}}\Phi]_{T,t}^{(i_{q_1}\ldots i_{q_{k-2r}})}$
will exist according the obtained version of the formula (\ref{2026may10sss})
under some additional conditions (see below).

When applying the formula (\ref{2026may10sss})
iteratively, we will obviously encounter with iterative 
application of traces $\tilde T^{k-2l}_{g_1,g_2,\ldots,g_{2l-1}, 
g_{2l}}\Phi(t_{q_1},\ldots,t_{q_{k-2l}})$ of different orders
$l.$ Let us assume that
$$
\tilde T^{k-2m}_{k_1,k_2,\ldots,k_{2m-1}, k_{2m}}
\tilde T^{k-2l}_{g_1,g_2,\ldots,g_{2l-1},g_{2l}}
\Phi(t_{q_1},\ldots,t_{q_{k-2(l+m)}})=
$$
\begin{equation}
\label{sogl111ds}
~~~~~~~~~~=
\tilde T^{k-2(l+m)}_{k_1,k_2,\ldots,k_{2m-1}, k_{2m},g_1,g_2,\ldots,g_{2l-1},g_{2l}}
\Phi(t_{q_1},\ldots,t_{q_{k-2(l+m)}})
\end{equation}

\vspace{1mm}
\noindent
in $L_2([t,T]^{k-2(l+m)}),$ where $l=0,1,2,\ldots,[k/2],$ $m=0,1,2,\ldots,[(k-2l)/2].$

We also suppose that conditions similar to (\ref{sogl111ds})
will be satisfied for the superposition of 3 or more traces
defined by the equality (\ref{soglas123sd}).

Note that the listed conditions will be satisfied if
the following condition is fulfilled.

{\it We will say that
condition {\rm $(\star\star\star)$} is satisfied if
for all $r=1,2,\ldots,[k/2]$ and for all possible $g_1,g_2,\ldots,$ $g_{2r-1},g_{2r}$ 
the sum of the series
$$
\sum\limits_{j_{g_1}, j_{g_3},\ldots ,j_{g_{2r-1}}=0}^{\infty}
C_{j_k\ldots j_1}\biggl|_{j_{g_1}=j_{g_2},\ldots, j_{g_{2r-1}}=j_{g_{2r}}}
$$

\noindent
does not depend on the method of summation$,$ i.e.$,$ for example$,$ for $r=3$}
$$
\sum\limits_{j_{g_1}, j_{g_3},j_{g_5}=0}^{\infty}
C_{j_k\ldots j_1}\biggl|_{j_{g_1}=j_{g_2},j_{g_3}=j_{g_4},j_{g_5}=j_{g_6}}=
\sum\limits_{j_{g_1}=0}^{\infty}
\sum\limits_{j_{g_3},j_{g_5}=0}^{\infty}
C_{j_k\ldots j_1}\biggl|_{j_{g_1}=j_{g_2},j_{g_3}=j_{g_4},j_{g_5}=j_{g_6}}=
$$
$$
=
\hspace{-1mm}\sum\limits_{j_{g_1},j_{g_3}=0}^{\infty}\sum\limits_{j_{g_5}=0}^{\infty}
C_{j_k\ldots j_1}\biggl|_{j_{g_1}=j_{g_2},j_{g_{3}}=j_{g_{4}},j_{g_5}=j_{g_6}}
=
\sum\limits_{j_{g_1}=0}^{\infty}\sum\limits_{j_{g_3}=0}^{\infty}\sum\limits_{j_{g_{5}}=0}^{\infty}
C_{j_k\ldots j_1}\biggl|_{j_{g_1}=j_{g_2},j_{g_{3}}=j_{g_{4}},j_{g_5}=j_{g_6}}\hspace{-0.5mm}.
$$

\vspace{2mm}

Suppose that the condition $(\star\star\star)$ is fulfilled.
Then,
by iteratively applying the formula (\ref{2026may10sss}), we obtain the following
inverse version of the Hu--Meyer
formula (\ref{2026may10})

\vspace{-3mm}
$$
J'[\Phi]_{T,t}^{(i_1\ldots i_k)}
=
\hat J^{S}[\Phi]_{T,t}^{(i_1\ldots i_k)}+
$$

\vspace{-3mm}
$$
+
\sum\limits_{r=1}^{[k/2]}
(-1)^r\sum_{\stackrel{(\{\{g_1, g_2\}, \ldots, 
\{g_{2r-1}, g_{2r}\}\}, \{q_1, \ldots, q_{k-2r}\})}
{{}_{\{g_1, g_2, \ldots, 
g_{2r-1}, g_{2r}, q_1, \ldots, q_{k-2r}\}=\{1, 2, \ldots, k\}}}}
\prod\limits_{s=1}^r
{\bf 1}_{\{i_{g_{{}_{2s-1}}}=~i_{g_{{}_{2s}}}\ne 0\}}\times
$$

\vspace{4mm}
\begin{equation}
\label{2026may38abcd}
\times
\hat J^{S}\hspace{-1mm}\left[\tilde T^{k-2r}_{g_1,g_2,\ldots,g_{2r-1}, 
g_{2r}}\Phi\right]_{T,t}^{(i_{q_1}\ldots i_{q_{k-2r}})}\ \ \ \hbox{w.~p.~1.}
\end{equation}

\vspace{3mm}

It is interesting to compare the formulas 
(\ref{2026may10}) and (\ref{2026may38abcd}) with similar
formulas (\ref{30.4}) and (\ref{2024str11}) that connect
the iterated Stratonovich and It\^{o} stochastic integrals.

Suppose that all limiting traces
$T^{k-2r}_{g_1,g_2,\ldots,g_{2r-1}, g_{2r}}\Phi (t_{q_1},\ldots,t_{q_{k-2r}})$
exist
for all $r=1,2,\ldots,$ $[k/2]$ and for all possible
$g_1,g_2,$ $\ldots,g_{2r-1},g_{2r}$ (see {\rm (\ref{leto5007after})).
Recall the equality $T^{k-2r}_{g_1,g_2,\ldots,g_{2r-1}, g_{2r}}\Phi (t_{q_1},\ldots,t_{q_{k-2r}})=
\tilde T^{k-2r}_{g_1,g_2,\ldots,g_{2r-1}, g_{2r}}\Phi (t_{q_1},\ldots,t_{q_{k-2r}})$ (in 
$L_2([t, T]^{k-2r})$), which is true in this case. 
Then
\begin{equation}
\label{2026may33}
~~~~~~~~~~~\lim\limits_{p\to\infty}
\left\Vert T^{k-2r,\hspace{0.2mm}p}_{g_1,g_2,\ldots,g_{2r-1}, 
g_{2r}}\Phi-
\tilde T^{k-2r,\hspace{0.2mm}p}_{g_1,g_2,\ldots,g_{2r-1}, 
g_{2r}}\Phi\right\Vert_{L_2([t, T]^{k-2r})}=0
\end{equation}

\noindent
for all $r=1,2,\ldots,[k/2]$ and for all possible
$g_1,g_2,\ldots,g_{2r-1},g_{2r}$.

Note that
$$
\left\Vert T^{k-2r,\hspace{0.2mm}p}_{g_1,g_2,\ldots,g_{2r-1}, 
g_{2r}}\Phi-
\tilde T^{k-2r,\hspace{0.2mm}p}_{g_1,g_2,\ldots,g_{2r-1}, 
g_{2r}}\Phi\right\Vert_{L_2([t, T]^{k-2r})}^2=
$$
\begin{equation}
\label{2026may32}
~~=\sum\limits_{j_{q_1},\ldots,j_{q_{k-2r}}=0}^p
\Biggl(
\sum\limits_{j_{g_1}, j_{g_3},\ldots ,j_{g_{2r-1}}=0}^p
C_{j_k\ldots j_1}\biggl|_{j_{g_1}=j_{g_2},\ldots, j_{g_{2r-1}}=j_{g_{2r}}}-
\tilde C_{j_{q_{k-2r}}\ldots j_{q_1}}
\Biggr)^2,
\end{equation}

\noindent
where $\tilde C_{j_{q_{k-2r}}\ldots j_{q_1}}$ is defined by (\ref{october202621}).

Now let us assume that 
all limiting traces
$\tilde T^{k-2r}_{g_1,g_2,\ldots,g_{2r-1}, g_{2r}}\Phi (t_{q_1},\ldots,t_{q_{k-2r}})$
exist
for all $r=1,2,\ldots,$ $[k/2]$ and for all possible
$g_1,g_2,$ $\ldots,g_{2r-1},g_{2r}$ (see {\rm (\ref{leto5007after})), and,
in addition, that the equality (\ref{2026may33}) is satisfied.

Then, by virtue of the equality
$$
T^{k-2r,\hspace{0.2mm}p}_{g_1,g_2,\ldots,g_{2r-1}, 
g_{2r}}\Phi=
$$

\vspace{-2mm}
$$
=\left(T^{k-2r,\hspace{0.2mm}p}_{g_1,g_2,\ldots,g_{2r-1}, 
g_{2r}}\Phi-
\tilde T^{k-2r,\hspace{0.2mm}p}_{g_1,g_2,\ldots,g_{2r-1}, 
g_{2r}}\Phi\right)+\tilde T^{k-2r,\hspace{0.2mm}p}_{g_1,g_2,\ldots,g_{2r-1}, 
g_{2r}}\Phi,
$$

\vspace{2mm}
\noindent
we obtain that 
all limiting traces
$T^{k-2r}_{g_1,g_2,\ldots,g_{2r-1}, g_{2r}}\Phi (t_{q_1},\ldots,t_{q_{k-2r}})$
exist
for all $r=1,2,\ldots,$ $[k/2]$ and for all possible
$g_1,g_2,$ $\ldots,g_{2r-1},g_{2r},$ and 
the equality $T^{k-2r}_{g_1,g_2,\ldots,g_{2r-1}, g_{2r}}\Phi (t_{q_1},\ldots,t_{q_{k-2r}})=
\tilde T^{k-2r}_{g_1,g_2,\ldots,g_{2r-1}, g_{2r}}\Phi (t_{q_1},\ldots,t_{q_{k-2r}})$ 
is fulfilled (in 
$L_2([t, T]^{k-2r})$).

Further, we will consider the connection of the above results
with Theorem~2.49. Suppose that the condition (\ref{2026may33})
is fulfilled
and all limiting traces 
$\tilde T^{k-2r}_{g_1,g_2,\ldots,g_{2r-1}, g_{2r}}\Phi(t_{q_1},\ldots,t_{q_{k-2r}})$ exist 
for all $r=1,2,\ldots,[k/2]$ and for all possible
$g_1,g_2,\ldots,g_{2r-1},g_{2r}.$

Using (\ref{2026may2}), we get
$$
{\sf M}\left\{\left(J'\hspace{-1mm}\left[T^{k-2r,\hspace{0.2mm}p}_{g_1,g_2,\ldots,g_{2r-1}, 
g_{2r}}\Phi\right]_{T,t}^{(i_{q_1}\ldots i_{q_{k-2r}})}-
J'\hspace{-1mm}\left[\tilde T^{k-2r,\hspace{0.2mm}p}_{g_1,g_2,\ldots,g_{2r-1}, 
g_{2r}}\Phi\right]_{T,t}^{(i_{q_1}\ldots i_{q_{k-2r}})}\right)^2\right\}\le
$$
$$
\le K_{k-2r} \left\Vert T^{k-2r,\hspace{0.2mm}p}_{g_1,g_2,\ldots,g_{2r-1}, 
g_{2r}}\Phi-
\tilde T^{k-2r,\hspace{0.2mm}p}_{g_1,g_2,\ldots,g_{2r-1}, 
g_{2r}}\Phi\right\Vert_{L_2([t, T]^{k-2r})}^2\ \to\ 0,
$$

\vspace{-2mm}
$$
{\sf M}\left\{\left(J'\hspace{-1mm}\left[\tilde T^{k-2r}_{g_1,g_2,\ldots,g_{2r-1}, 
g_{2r}}\Phi\right]_{T,t}^{(i_{q_1}\ldots i_{q_{k-2r}})}-
J'\hspace{-1mm}\left[\tilde T^{k-2r,\hspace{0.2mm}p}_{g_1,g_2,\ldots,g_{2r-1}, 
g_{2r}}\Phi\right]_{T,t}^{(i_{q_1}\ldots i_{q_{k-2r}})}\right)^2\right\}\le
$$
$$
\le K'_{k-2r} \left\Vert \tilde T^{k-2r}_{g_1,g_2,\ldots,g_{2r-1}, 
g_{2r}}\Phi-
\tilde T^{k-2r,\hspace{0.2mm}p}_{g_1,g_2,\ldots,g_{2r-1}, 
g_{2r}}\Phi\right\Vert_{L_2([t, T]^{k-2r})}^2\ \to\ 0
$$

\vspace{2mm}
\noindent
if $p\to\infty$ (for all $r=1,2,\ldots,[k/2]$ and for all possible
$g_1,g_2,\ldots,g_{2r-1},g_{2r}$
(see (\ref{leto5007after}))), where $K_{k-2r}, K'_{k-2r}$ are constants.
Then w.~p.~1
$$
\hbox{\vtop{\offinterlineskip\halign{
\hfil#\hfil\cr
{\rm l.i.m.}\cr
$\stackrel{}{{}_{p\to \infty}}$\cr
}} }
J'\hspace{-1mm}\left[T^{k-2r,\hspace{0.2mm}p}_{g_1,g_2,\ldots,g_{2r-1}, 
g_{2r}}\Phi\right]_{T,t}^{(i_{q_1}\ldots i_{q_{k-2r}})}=
$$

\vspace{-2mm}
$$
=
\hbox{\vtop{\offinterlineskip\halign{
\hfil#\hfil\cr
{\rm l.i.m.}\cr
$\stackrel{}{{}_{p\to \infty}}$\cr
}} }
J'\hspace{-1mm}\left[\tilde T^{k-2r,\hspace{0.2mm}p}_{g_1,g_2,\ldots,g_{2r-1}, 
g_{2r}}\Phi\right]_{T,t}^{(i_{q_1}\ldots i_{q_{k-2r}})}=
$$

\vspace{-2mm}
\begin{equation}
\label{october202615}
=
J'\hspace{-1mm}\left[\tilde T^{k-2r}_{g_1,g_2,\ldots,g_{2r-1}, 
g_{2r}}\Phi\right]_{T,t}^{(i_{q_1}\ldots i_{q_{k-2r}})}.
\end{equation}

\vspace{2mm}

Passing to the limit
$\hbox{\vtop{\offinterlineskip\halign{
\hfil#\hfil\cr
{\rm l.i.m.}\cr
$\stackrel{}{{}_{p\to \infty}}$\cr
}} }$ in (\ref{2026may1}) and applying (\ref{october202615}), we obtain
w.~p.~1

\newpage
\noindent
$$
\hat J^{S}[\Phi]_{T,t}^{(i_1\ldots i_k)}=\hbox{\vtop{\offinterlineskip\halign{
\hfil#\hfil\cr
{\rm l.i.m.}\cr
$\stackrel{}{{}_{p\to \infty}}$\cr
}} }\hat J^{S}_p[\Phi]_{T,t}^{(i_1\ldots i_k)}
=
\hbox{\vtop{\offinterlineskip\halign{
\hfil#\hfil\cr
{\rm l.i.m.}\cr
$\stackrel{}{{}_{p\to \infty}}$\cr
}} }J'[\Phi_p]_{T,t}^{(i_1\ldots i_k)}+
$$

\vspace{-3mm}
$$
+
\sum\limits_{r=1}^{[k/2]}
\sum_{\stackrel{(\{\{g_1, g_2\}, \ldots, 
\{g_{2r-1}, g_{2r}\}\}, \{q_1, \ldots, q_{k-2r}\})}
{{}_{\{g_1, g_2, \ldots, 
g_{2r-1}, g_{2r}, q_1, \ldots, q_{k-2r}\}=\{1, 2, \ldots, k\}}}}
\prod\limits_{s=1}^r
{\bf 1}_{\{i_{g_{{}_{2s-1}}}=~i_{g_{{}_{2s}}}\ne 0\}}\times
$$

\vspace{4mm}
\begin{equation}
\label{october202616}
~~\times
\hbox{\vtop{\offinterlineskip\halign{
\hfil#\hfil\cr
{\rm l.i.m.}\cr
$\stackrel{}{{}_{p\to \infty}}$\cr
}} }J'\hspace{-1mm}\left[\tilde T^{k-2r,\hspace{0.2mm}p}_{g_1,g_2,\ldots,g_{2r-1}, 
g_{2r}}\Phi\right]_{T,t}^{(i_{q_1}\ldots i_{q_{k-2r}})},
\end{equation}

\vspace{4mm}
\noindent
where notations are the same as in (\ref{2026may1}).

Replacing the function $\Phi(t_1,\ldots,t_k)$ in (\ref{october202616}) 
with the Volterra-type kernel (see (\ref{july80027}))
\begin{equation}
\label{october202633}
~~~~~~~~~~K(t_1,\ldots,t_k)=
\psi_1(t_1)\ldots \psi_k(t_k){\bf 1}_{\{t_1<\ldots<t_{k}\}}\ \ \ (k\ge 2),
\end{equation}

\noindent
where $\psi_1(\tau),\ldots,\psi_k(\tau)\in L_2([t,T])$ and
$t_1,\ldots,t_k\in [t, T],$ and using (\ref{febr5000}), we have

\vspace{-4mm}
$$
\hbox{\vtop{\offinterlineskip\halign{
\hfil#\hfil\cr
{\rm l.i.m.}\cr
$\stackrel{}{{}_{p\to \infty}}$\cr
}} }\sum_{j_1,\ldots,j_k=0}^{p}
C_{j_k\ldots j_1}
\prod_{l=1}^k \zeta_{j_l}^{(i_l)}
=
J[\psi^{(k)}]_{T,t}^{(i_1\ldots i_k)}+
$$

\vspace{-2mm}
$$
+
\sum\limits_{r=1}^{[k/2]}
\sum_{\stackrel{(\{\{g_1, g_2\}, \ldots, 
\{g_{2r-1}, g_{2r}\}\}, \{q_1, \ldots, q_{k-2r}\})}
{{}_{\{g_1, g_2, \ldots, 
g_{2r-1}, g_{2r}, q_1, \ldots, q_{k-2r}\}=\{1, 2, \ldots, k\}}}}
\prod\limits_{s=1}^r
{\bf 1}_{\{i_{g_{{}_{2s-1}}}=~i_{g_{{}_{2s}}}\ne 0\}}\times
$$

\vspace{5mm}
\begin{equation}
\label{october202617}
~~\times
\hbox{\vtop{\offinterlineskip\halign{
\hfil#\hfil\cr
{\rm l.i.m.}\cr
$\stackrel{}{{}_{p\to \infty}}$\cr
}} }J'\hspace{-1mm}\left[\tilde T^{k-2r,\hspace{0.2mm}p}_{g_1,g_2,\ldots,g_{2r-1}, 
g_{2r}}K\right]_{T,t}^{(i_{q_1}\ldots i_{q_{k-2r}})},
\end{equation}

\vspace{4mm}
\noindent
where $J[\psi^{(k)}]_{T,t}^{(i_1\ldots i_k)}$
is the iterated It\^{o} stochastic integral
(\ref{itoxx}).

Recall the equality (\ref{febr14}) (see the proof of Theorem~2.49)
$$
\hbox{\vtop{\offinterlineskip\halign{
\hfil#\hfil\cr
{\rm l.i.m.}\cr
$\stackrel{}{{}_{p\to \infty}}$\cr
}} }
\sum\limits_{j_{q_1},\ldots,j_{q_{k-2r}}=0}^p
\frac{1}{2^r}
C_{j_k \ldots j_1}\biggl|_{(j_{g_2} j_{g_1})\curvearrowright (\cdot)
\ldots (j_{g_{2r}} j_{g_{2r-1}})\curvearrowright (\cdot),
j_{g_{{}_{1}}}=~j_{g_{{}_{2}}},\ldots, j_{g_{{}_{2r-1}}}=~j_{g_{{}_{2r}}}}\biggr.
\times 
$$

\vspace{2mm}
$$
\times
\prod\limits_{s=1}^r
{\bf 1}_{\{i_{g_{{}_{2s-1}}}=~i_{g_{{}_{2s}}}\ne 0\}}
J'[\phi_{j_{q_1}}\ldots \phi_{j_{q_{k-2r}}}]_{T,t}^{(i_{q_1}\ldots i_{q_{k-2r}})}=
$$
\begin{equation}
\label{october202618}
=\frac{1}{2^r}
J[\psi^{(k)}]_{T,t}^{s_r, \ldots, s_1}
\end{equation}

\vspace{1mm}
\noindent
w.~p.~1, where $g_{2}=g_{1}+1,\ldots, g_{2r}=g_{2r-1}+1,$
$g_{2i-1}\stackrel{\sf def}{=}s_i;$\ $i=1,2,\ldots,r;$\
$r=1,2,\ldots,\left[k/2\right],$ 
$(s_r,\ldots,s_1)\in {\rm A}_{k,r},$ $J[\psi^{(k)}]_{T,t}^{s_r,\ldots,s_1}$ is
defined by (\ref{30.1}) and ${\rm A}_{k,r}$ is defined by (\ref{30.5550001});
another notations in (\ref{october202618}) are the same as in Sect.~2.10.

Also, recall that (see (\ref{july100000}))
$$
\lim\limits_{p\to\infty}
\sum\limits_{j_{g_1}, j_{g_3},\ldots ,j_{g_{2r-1}}=0}^p
C_{j_k\ldots j_1}\biggl|_{j_{g_1}=j_{g_2},\ldots, j_{g_{2r-1}}=j_{g_{2r}}}=
$$
\begin{equation}
\label{october202619}
=\frac{1}{2^r} \prod\limits_{l=1}^r {\bf 1}_{\{g_{2l}=g_{2l-1}+1\}}
C_{j_k \ldots j_1}\biggl|_{(j_{g_2} j_{g_1})\curvearrowright (\cdot)
\ldots (j_{g_{2r}} j_{g_{2r-1}})\curvearrowright (\cdot),
j_{g_{{}_{1}}}=~j_{g_{{}_{2}}},\ldots, j_{g_{{}_{2r-1}}}=~j_{g_{{}_{2r}}}
}\biggr.
\end{equation}

\vspace{3mm}
\noindent
for all possible $g_1,g_2,\ldots,g_{2r-1},g_{2r}$ (see {\rm (\ref{leto5007after})),
where $k\ge 2r,$ $r=1,2,\ldots,$ $[k/2],$ $C_{j_k\ldots j_1}$ is defined by (\ref{july99999});
another notations are the same as in Theorem~2.49.

Combining (\ref{october202620}), (\ref{october202621}), (\ref{october202618}),
(\ref{october202619}) and (\ref{after800}), we obtain w.~p.~1
$$
\prod\limits_{s=1}^r
{\bf 1}_{\{i_{g_{{}_{2s-1}}}=~i_{g_{{}_{2s}}}\ne 0\}}\
\hbox{\vtop{\offinterlineskip\halign{
\hfil#\hfil\cr
{\rm l.i.m.}\cr
$\stackrel{}{{}_{p\to \infty}}$\cr
}} }J'\hspace{-1mm}\left[\tilde T^{k-2r,\hspace{0.2mm}p}_{g_1,g_2,\ldots,g_{2r-1}, 
g_{2r}}K\right]_{T,t}^{(i_{q_1}\ldots i_{q_{k-2r}})}=
$$
$$
=\prod\limits_{l=1}^r {\bf 1}_{\{g_{2l}=g_{2l-1}+1\}}
\frac{1}{2^r}
J[\psi^{(k)}]_{T,t}^{s_r, \ldots, s_1},
$$

\vspace{1mm}
\noindent
and
$$
J[\psi^{(k)}]_{T,t}^{(i_1\ldots i_k)}+
$$

\vspace{-4mm}
$$
+\sum\limits_{r=1}^{[k/2]}
\sum_{\stackrel{(\{\{g_1, g_2\}, \ldots, 
\{g_{2r-1}, g_{2r}\}\}, \{q_1, \ldots, q_{k-2r}\})}
{{}_{\{g_1, g_2, \ldots, 
g_{2r-1}, g_{2r}, q_1, \ldots, q_{k-2r}\}=\{1, 2, \ldots, k\}}}}
\prod\limits_{s=1}^r
{\bf 1}_{\{i_{g_{{}_{2s-1}}}=~i_{g_{{}_{2s}}}\ne 0\}}\times
$$

\vspace{2mm}
$$
\times
\hbox{\vtop{\offinterlineskip\halign{
\hfil#\hfil\cr
{\rm l.i.m.}\cr
$\stackrel{}{{}_{p\to \infty}}$\cr
}} }J'\hspace{-1mm}\left[\tilde T^{k-2r,\hspace{0.2mm}p}_{g_1,g_2,\ldots,g_{2r-1}, 
g_{2r}}K\right]_{T,t}^{(i_{q_1}\ldots i_{q_{k-2r}})}=
$$

\begin{equation}
\label{october202623}
~~~~~~~~~=J[\psi^{(k)}]_{T,t}^{(i_1\ldots i_k)}+
\sum_{r=1}^{\left[k/2\right]}\frac{1}{2^r}
\sum_{(s_r,\ldots,s_1)\in {\rm A}_{k,r}}
J[\psi^{(k)}]_{T,t}^{s_r,\ldots,s_1}.
\end{equation}

\vspace{2mm}

Suppose that 
$\psi_1(\tau), \ldots, \psi_k(\tau)$ are continuous 
functions on $[t, T].$ Further, applying Theorem~2.12
to the right-hand side of (\ref{october202623})
and combining (\ref{october202617}), (\ref{october202623}), we get
the following expansion
$$
J^{*}[\psi^{(k)}]_{T,t}^{(i_1\ldots i_k)}=\hbox{\vtop{\offinterlineskip\halign{
\hfil#\hfil\cr
{\rm l.i.m.}\cr
$\stackrel{}{{}_{p\to \infty}}$\cr
}} }\sum_{j_1,\ldots,j_k=0}^{p}
C_{j_k\ldots j_1}
\prod_{l=1}^k \zeta_{j_l}^{(i_l)},
$$

\noindent
where $J^{*}[\psi^{(k)}]_{T,t}^{(i_1\ldots i_k)}$ is the iterated
Stratonovich stochastic integral (\ref{strxx}).

Note that in our case (the case of the Volterra-type kernel (\ref{october202633})), limiting traces
$\tilde T^{k-2r}_{g_1,g_2,\ldots,g_{2r-1}, g_{2r}}K$ exist 
for all $r=1,2,\ldots,[k/2]$ and for all possible
$g_1,g_2,\ldots,g_{2r-1},g_{2r},$ and are determined
by the equality (\ref{2025may25}) (also see (\ref{2025may28})).
Moreover, the condition (\ref{2026may33}) in this case
has the form (\ref{july700000}) (or (\ref{2025may26})).

Further, we will consider another approach to deriving
the Hu--Meyer formulas based on the so-called
Hilbert space valued traces \cite{bugh1}, \cite{bugh3}, \cite{HuHu}. At that we will still
consider the case of a multidimensional Wiener process.

Recall the following definition of multiple Stratonovich stochastic
integral (see (\ref{30.34ququ}) or (\ref{30.34}))
\begin{equation}
\label{30.34ququsss}
~~~~~~~~~~ \hbox{\vtop{\offinterlineskip\halign{
\hfil#\hfil\cr
{\rm l.i.m.}\cr
$\stackrel{}{{}_{N\to \infty}}$\cr
}} }\sum_{j_1,\ldots,j_k=0}^{N-1}
\Phi\left(\tau_{j_1},\ldots,\tau_{j_k}\right)
\prod\limits_{l=1}^k\Delta{\bf w}_{\tau_{j_l}}^{(i_l)}
\stackrel{\rm def}{=}J[\Phi]_{T,t}^{(i_1\ldots i_k)},
\end{equation}

\noindent
where 
$\Phi(t_1,\ldots,t_k):\ [t, T]^k\to{\bf R}$ is a 
continuous nonrandom
function, 
${\bf w}_{\tau}$ is a random vector with 
an $m+1$ components
(${\bf w}_{\tau}^{(i)} $ $(i=1,\ldots,m)$
are independent standard Wiener processes and
${\bf w}_{\tau}^{(0)}=\tau$),
$\Delta {\bf w}_{\tau_j}^{(i)}={\bf w}_{\tau_{j+1}}^{(i)}-{\bf w}_{\tau_j}^{(i)},$
$\left\{\tau_{j}\right\}_{j=0}^{N}$ is a partition of
$[t,T]$ which satisfies the condition (\ref{1111}),
and $i_1,\ldots,i_k=0,1,\ldots,m.$ 

Note that the function $\Phi(t_1,\ldots,t_k)$
in (\ref{30.34ququsss}) can be discontinuous.
For example, it can have the form (\ref{1999.1xx}).

Let us write the formula (\ref{sogl}), replacing  
$R_{p_1 p_2 p_3 p_4}(t_1,\ldots,t_4)$ with
$\Phi(t_1,\ldots,t_4)$ and using the folmula (\ref{Wi110}) for the
multiple Wiener stochastic integral $J'[\Phi]_{T,t}^{(i_1\ldots i_4)}$ defined by (\ref{WiI})
(also see (\ref{mult11}) and (\ref{pobeda}))

\vspace{-2mm}
$$
J[\Phi]_{T,t}^{(i_1\ldots i_4)}=J'[\Phi]_{T,t}^{(i_1\ldots i_4)}+
$$

\vspace{-3mm}
$$
+{\bf 1}_{\{i_1=i_2\ne 0\}}J'\left[\Phi(t_1,\ldots,t_4)\bigl.\bigr|_{t_1=t_2}
\right]_{T,t}^{(0 i_3 i_4)}+
{\bf 1}_{\{i_1=i_3\ne 0\}}J'\left[\Phi(t_1,\ldots,t_4)\bigl.\bigr|_{t_1=t_3}
\right]_{T,t}^{(0 i_2 i_4)}+
$$
$$
+{\bf 1}_{\{i_1=i_4\ne 0\}}J'\left[\Phi(t_1,\ldots,t_4)\bigl.\bigr|_{t_1=t_4}
\right]_{T,t}^{(0 i_2 i_3)}+
{\bf 1}_{\{i_2=i_3\ne 0\}}J'\left[\Phi(t_1,\ldots,t_4)\bigl.\bigr|_{t_2=t_3}
\right]_{T,t}^{(i_1 0 i_4)}+
$$

\vspace{-3mm}
$$
+{\bf 1}_{\{i_2=i_4\ne 0\}}J'\left[\Phi(t_1,\ldots,t_4)\bigl.\bigr|_{t_2=t_4}
\right]_{T,t}^{(i_1 0 i_3)}+
{\bf 1}_{\{i_3=i_4\ne 0\}}J'\left[\Phi(t_1,\ldots,t_4)\bigl.\bigr|_{t_3=t_4}
\right]_{T,t}^{(i_1 i_2 0)}+
$$

\vspace{-1mm}
$$
+{\bf 1}_{\{i_1=i_2\ne 0\}}{\bf 1}_{\{i_3=i_4\ne 0\}}
J'\left[\Phi(t_1,\ldots,t_4)\bigl.\bigr|_{t_1=t_2,t_3=t_4}
\right]_{T,t}^{(00)}+
$$

\vspace{-1mm}
$$
+{\bf 1}_{\{i_1=i_3\ne 0\}}{\bf 1}_{\{i_2=i_4\ne 0\}}
J'\left[\Phi(t_1,\ldots,t_4)\bigl.\bigr|_{t_1=t_3,t_2=t_4}
\right]_{T,t}^{(00)}+
$$

\vspace{-3mm}
\begin{equation}
\label{2025nov100}
~~~~~~~~+{\bf 1}_{\{i_1=i_4\ne 0\}}{\bf 1}_{\{i_2=i_3\ne 0\}}
J'\left[\Phi(t_1,\ldots,t_4)\bigl.\bigr|_{t_1=t_4,t_2=t_3}
\right]_{T,t}^{(00)}
\end{equation}

\vspace{2mm}
\noindent
w.~p.~1. Further, applying a Fubini-type theorem for stochastic and 
Lebesgue intgrals (as in \cite{farre}) to the right-hand side of
(\ref{2025nov100}), we obtain w.~p.~1

\vspace{-2mm}
$$
J[\Phi]_{T,t}^{(i_1\ldots i_4)}=
J'[\Phi]_{T,t}^{(i_1\ldots i_4)}
+
$$

\vspace{-2mm}
$$
+{\bf 1}_{\{i_1=i_2\ne 0\}}J'\Biggl[
~\int\limits_{[t, T]} \Phi(t_1,\ldots,t_4)\bigl.\bigr|_{t_1=t_2}
dt_1\Biggr]_{T,t}^{(i_3 i_4)}+
$$

\vspace{-2mm}
$$
+
{\bf 1}_{\{i_1=i_3\ne 0\}}J'\Biggl[~\int\limits_{[t, T]}\Phi(t_1,\ldots,t_4)\bigl.\bigr|_{t_1=t_3}
dt_1
\Biggr]_{T,t}^{(i_2 i_4)}+
$$

\vspace{-2mm}
$$
+{\bf 1}_{\{i_1=i_4\ne 0\}}J'\Biggl[~\int\limits_{[t, T]}
\Phi(t_1,\ldots,t_4)\bigl.\bigr|_{t_1=t_4}
dt_1\Biggr]_{T,t}^{(i_2 i_3)}+
$$

\vspace{-2mm}
$$
+
{\bf 1}_{\{i_2=i_3\ne 0\}}J'\Biggl[~\int\limits_{[t, T]}
\Phi(t_1,\ldots,t_4)\bigl.\bigr|_{t_2=t_3}
dt_2\Biggr]_{T,t}^{(i_1 i_4)}+
$$

\vspace{-2mm}
$$
+{\bf 1}_{\{i_2=i_4\ne 0\}}J'\Biggl[~\int\limits_{[t, T]}
\Phi(t_1,\ldots,t_4)\bigl.\bigr|_{t_2=t_4}
dt_2\Biggr]_{T,t}^{(i_1 i_3)}+
$$

\vspace{-2mm}
$$
+
{\bf 1}_{\{i_3=i_4\ne 0\}}J'\Biggl[~\int\limits_{[t, T]}
\Phi(t_1,\ldots,t_4)\bigl.\bigr|_{t_3=t_4}
dt_3\Biggr]_{T,t}^{(i_1 i_2)}+
$$
$$
+{\bf 1}_{\{i_1=i_2\ne 0\}}{\bf 1}_{\{i_3=i_4\ne 0\}}
\int\limits_{[t, T]^2}\Phi(t_1,\ldots,t_4)\bigl.\bigr|_{t_1=t_2,t_3=t_4}
dt_1 dt_3+
$$

$$
+{\bf 1}_{\{i_1=i_3\ne 0\}}{\bf 1}_{\{i_2=i_4\ne 0\}}
\int\limits_{[t, T]^2}\Phi(t_1,\ldots,t_4)\bigl.\bigr|_{t_1=t_3,t_2=t_4}
dt_1 dt_2+
$$

\vspace{-2mm}
\begin{equation}
\label{2025nov101}
~~~~~~~~~~+{\bf 1}_{\{i_1=i_4\ne 0\}}{\bf 1}_{\{i_2=i_3\ne 0\}}
\int\limits_{[t, T]^2}\Phi(t_1,\ldots,t_4)\bigl.\bigr|_{t_1=t_4,t_2=t_3}
dt_1 dt_2.
\end{equation}

\vspace{2mm}

Let us generalize the formula (\ref{2025nov101}) to the case $k\in{\bf N}$.
As a result, we obtain the following Hu--Meyer formula for the 
case of a multidimensional Wiener process

\vspace{-3mm}
$$
J[\Phi]_{T,t}^{(i_1\ldots i_k)}
=
J'[\Phi]_{T,t}^{(i_1\ldots i_k)}+
$$

\vspace{-3mm}
$$
+
\sum\limits_{r=1}^{[k/2]}
\sum_{\stackrel{(\{\{g_1, g_2\}, \ldots, 
\{g_{2r-1}, g_{2r}\}\}, \{q_1, \ldots, q_{k-2r}\})}
{{}_{\{g_1, g_2, \ldots, 
g_{2r-1}, g_{2r}, q_1, \ldots, q_{k-2r}\}=\{1, 2, \ldots, k\}}}}
\prod\limits_{s=1}^r
{\bf 1}_{\{i_{g_{{}_{2s-1}}}=~i_{g_{{}_{2s}}}\ne 0\}}\times
$$

\vspace{4mm}
\begin{equation}
\label{2026may10zyx}
\times
J'\hspace{-1mm}\left[\breve T^{k-2r}_{g_1,g_2,\ldots,g_{2r-1}, 
g_{2r}}\Phi\right]_{T,t}^{(i_{q_1}\ldots i_{q_{k-2r}})}\ \ \ \hbox{w.~p.~1,}
\end{equation}

\vspace{3mm}
\noindent
where the so-called 
Hilbert space valued trace 
$\breve T^{k-2r}_{g_1,g_2,\ldots,g_{2r-1}, 
g_{2r}}\Phi\in L_2([t, T]^{k-2r})$ 
has the form (we assume that this trace exists)

\vspace{-2mm}
$$
\breve T^{k-2r}_{g_1,g_2,\ldots,g_{2r-1}, 
g_{2r}}\Phi(t_{q_1},\ldots,t_{q_{k-2r}})\stackrel{\sf def}{=}
$$

$$
\stackrel{\sf def}{=}
\int\limits_{[t, T]^r}
\Phi(
t_1,\ldots,t_k)\biggl.\biggr|_{t_{g_{{}_{1}}}=t_{g_{{}_{2}}},\ldots,
t_{g_{{}_{2r-1}}}=t_{g_{{}_{2r}}}}dt_{g_{{}_{1}}}\ldots dt_{g_{{}_{2r-1}}},
$$

\vspace{2mm}
\noindent
where 
$$
J'\hspace{-1mm}\left[\breve T^{k-2r}_{g_1,g_2,\ldots,g_{2r-1}, 
g_{2r}}\Phi\right]_{T,t}^{(i_{q_1}\ldots i_{q_{k-2r}})}\stackrel{\sf def}{=}
\breve T^{k-2r}_{g_1,g_2,\ldots,g_{2r-1}, 
g_{2r}}\Phi
$$ 

\vspace{1mm}
\noindent
for $k=2r;$
another notations are the same as in (\ref{2026may10}).
Recall that
in \cite{bugh1}, the variant $\breve T^{k-2r}_{1,2\ldots,2r-1,2r}\Phi(t_{2r+1},\ldots,t_k)$
of $\breve T^{k-2r}_{g_1,g_2,\ldots,g_{2r-1}, 
g_{2r}}\Phi(t_{q_1},\ldots,t_{q_{k-2r}})$ was considered.

The inverse version of the Hu--Meyer
formula (\ref{2026may10zyx}) has the form

\vspace{-2mm}
$$
J'[\Phi]_{T,t}^{(i_1\ldots i_k)}
=
J[\Phi]_{T,t}^{(i_1\ldots i_k)}+
$$

\vspace{-3mm}
$$
+
\sum\limits_{r=1}^{[k/2]}
(-1)^r\sum_{\stackrel{(\{\{g_1, g_2\}, \ldots, 
\{g_{2r-1}, g_{2r}\}\}, \{q_1, \ldots, q_{k-2r}\})}
{{}_{\{g_1, g_2, \ldots, 
g_{2r-1}, g_{2r}, q_1, \ldots, q_{k-2r}\}=\{1, 2, \ldots, k\}}}}
\prod\limits_{s=1}^r
{\bf 1}_{\{i_{g_{{}_{2s-1}}}=~i_{g_{{}_{2s}}}\ne 0\}}\times
$$

\vspace{4mm}
$$
\times
J\hspace{-1mm}\left[\breve T^{k-2r}_{g_1,g_2,\ldots,g_{2r-1}, 
g_{2r}}\Phi\right]_{T,t}^{(i_{q_1}\ldots i_{q_{k-2r}})}\ \ \ \hbox{w.~p.~1.}
$$

\vspace{2mm}

Note that there are more general 
definitions of the multiple Stratonovich stochastic integral
(see, for example, (1.5.9) and Sect.~2.1 in \cite{bugh1})
that are consistent with definition (\ref{30.34ququsss}) on the class
of continuous functions. 

Further, consider one of these generalizations. 
But first, consider the following Riemann--Stieltjes integral
\begin{equation}
\label{novem2026xxx2}
\int\limits_{[t, T]^k}\Phi(t_1,\ldots,t_k)
d{\bf w}_{t_1}^{(i_1)p}\ldots d{\bf w}_{t_k}^{(i_k)p},
\end{equation}
where $p\in{\bf N}$, $i_1,\ldots,i_k=0,1,\ldots,m,$ $\Phi(t_1,\ldots,t_k):\ [t, T]^k\to{\bf R}$ is a 
nonrandom function, 
${\bf w}_{\tau}$ is a random vector with 
an $m+1$ components
(${\bf w}_{\tau}^{(i)} $ $(i=1,\ldots,m)$
are independent standard Wiener processes and
${\bf w}_{\tau}^{(0)}=\tau$),
${\bf w}_{\tau}^{(i)p}$ is the following 
mean-square approximation of 
${\bf w}_{\tau}^{(i)}$ \cite{Lipt}
\begin{equation}
\label{novem2026xxx1000}
{\bf w}_{\tau}^{(i)p}={\bf w}_{t}^{(i)p}+
\sum_{j=0}^{p}\int\limits_t^{\tau}
\phi_j(s)ds\ \zeta_j^{(i)},
\end{equation}
where
$$
\zeta_j^{(i)}=
\int\limits_t^T \phi_j(s)d{\bf w}_s^{(i)},
$$

\noindent
$\tau\in[t, T],$ $t\ge 0,$
$\{\phi_j(x)\}_{j=0}^{\infty}$ is an arbitrary complete 
orthonormal system of functions in the space $L_2([t, T]),$ and
$\zeta_j^{(i)}$ are independent standard Gaussian 
random variables for various $i$ or $j$ (in the case when $i\ne 0$).

From (\ref{novem2026xxx1000}) we have
\begin{equation}
\label{novem2026xxx1}
d{\bf w}_{\tau}^{(i)p}=
\sum_{j=0}^{p}
\phi_j(\tau)\zeta_j^{(i)} d\tau.
\end{equation}

Let us substitute (\ref{novem2026xxx1}) into (\ref{novem2026xxx2})
\begin{equation}
\label{novem2026xxx3}
~~~~~~~~~~~~\int\limits_{[t, T]^k}\Phi(t_1,\ldots,t_k)
d{\bf w}_{t_1}^{(i_1)p}\ldots d{\bf w}_{t_k}^{(i_k)p}
=\sum\limits_{j_1,\ldots,j_k=0}^{p}
C_{j_k \ldots j_1}\prod\limits_{l=1}^k \zeta_{j_l}^{(i_l)},
\end{equation}

\noindent
where
$$
C_{j_k \ldots j_1}=\int\limits_{[t,T]^k}
\Phi(t_1,\ldots,t_k)\prod\limits_{l=1}^k \phi_{j_l}(t_l)
dt_1\ldots dt_k
$$

\noindent
is the Fourier coefficient.

The formula (\ref{novem2026xxx3}) explains why
$\hat J^{S}[\Phi]_{T,t}^{(i_1\ldots i_k)}$ (see
(\ref{novem2026xxx4}), (\ref{novem2026xxx5}))
was called
the multiple Stratonovich stochastic integral.

Now consider another approximation of the Wiener process,
namely the so-called polygonal approximation (see \cite{HuHu}, Example~5.17)

\vspace{-3mm}
\begin{equation}
\label{novem2026xxx7}
~~~~~~~~~{\bf w}_{\tau}^{(i)N}=
\sum\limits_{j=0}^{N-1}\left({\bf w}_{\tau_{j}}^{(i)}+
\frac{1}{\Delta_j}\Delta {\bf w}_{\tau_j}^{(i)}
(\tau-\tau_j)\right)
{\bf 1}_{T_j}(\tau),
\end{equation}

\vspace{1mm}
\noindent
where $N\in{\bf N},$ $i=0,1,\ldots,m,$  $\tau\in [t, T]$, $t\ge 0,$
${\bf w}_{\tau}$ is a random vector with 
an $m+1$ components (${\bf w}_{\tau}^{(i)} $ $(i=1,\ldots,m)$
are independent standard Wiener processes and
${\bf w}_{\tau}^{(0)}=\tau$),
$\Delta {\bf w}_{\tau_j}^{(i)}={\bf w}_{\tau_{j+1}}^{(i)}-{\bf w}_{\tau_j}^{(i)},$
$\Delta_j=\tau_{j+1}-\tau_j,$
$T_j=[\tau_j, \tau_{j+1}),$ ${\bf 1}_A$ is the indicator of the set $A,$
$\{\tau_j\}_{j=0}^N$ is a partition of
the interval $[t,T]$ such that
\begin{equation}
\label{novem2026xxx6}
t=\tau_0<\tau_1<\ldots <\tau_N=T,\ \ \
\max\limits_{0\le j\le N-1}\left|\tau_{j+1}-\tau_j\right|\to 0\ \
\hbox{if}\ \ N\to \infty.
\end{equation}

From (\ref{novem2026xxx7}) we have
\begin{equation}
\label{novem2026xxx8}
d{\bf w}_{\tau}^{(i)N}=
\sum\limits_{j=0}^{N-1}
\frac{1}{\Delta_j}\Delta {\bf w}_{\tau_j}^{(i)}
d\tau
{\bf 1}_{T_j}(\tau).
\end{equation}

Using (\ref{novem2026xxx8}), we obtain
$$
\int\limits_{[t, T]^k}\Phi(t_1,\ldots,t_k)
d{\bf w}_{t_1}^{(i_1)N}\ldots d{\bf w}_{t_k}^{(i_k)N}
=\sum\limits_{j_1,\ldots,j_k=0}^{N-1}
\frac{1}{\Delta_{j_1}\ldots \Delta_{j_k}}\times
$$

\vspace{-1.5mm}
$$
\times
\int\limits_{[t,T]^k}\Phi(t_1,\ldots,t_k)
{\bf 1}_{T_{j_1}}(t_1)\ldots {\bf 1}_{T_{j_k}}(t_k)dt_1\ldots dt_k
\Delta {\bf w}_{\tau_{j_1}}^{(i_1)}\ldots \Delta {\bf w}_{\tau_{j_k}}^{(i_k)}=
$$
$$
=\sum\limits_{j_1,\ldots,j_k=0}^{N-1}
\frac{1}{\Delta_{j_1}\ldots \Delta_{j_k}}
\int\limits_{T_{j_1}\times \ldots \times T_{j_k}}\Phi(t_1,\ldots,t_k) dt_1\ldots dt_k
\Delta {\bf w}_{\tau_{j_1}}^{(i_1)}\ldots \Delta {\bf w}_{\tau_{j_k}}^{(i_k)}.
$$

\vspace{2mm}

We define the multiple Stratonovich
stochastic integral (see \cite{SU11}, \cite{HuHu})
for $\Phi(t_1,\ldots,t_k)\in L_2([t, T]^k)$ (but this function 
satisfies an additional condition (see below))
as the following mean-square limit (the case of a multidimensional Wiener process)

\vspace{-4mm}
$$
\bar J^{S}[\Phi]_{T,t}^{(i_1\ldots i_k)}=
\hbox{\vtop{\offinterlineskip\halign{
\hfil#\hfil\cr
{\rm l.i.m.}\cr
$\stackrel{}{{}_{N\to \infty}}$\cr
}} }
\sum\limits_{j_1,\ldots,j_k=0}^{N-1}
\frac{1}{\Delta_{j_1}\ldots \Delta_{j_k}}\times
$$

\begin{equation}
\label{novem2026xxx10}
~~~~~~~~\times
\int\limits_{T_{j_1}\times \ldots \times T_{j_k}}\Phi(t_1,\ldots,t_k) dt_1\ldots dt_k
\Delta {\bf w}_{\tau_{j_1}}^{(i_1)}\ldots \Delta {\bf w}_{\tau_{j_k}}^{(i_k)},
\end{equation}

\vspace{2mm}
\noindent
where $i_1,\ldots,i_k=0, 1,\ldots,m$ and $\{\tau_j\}_{j=0}^N$ 
as in (\ref{novem2026xxx6}).

We define the 
$r$th trace of the function $\Phi(t_1,\ldots,t_k)\in L_2([t, T]^k)$
by the following expression \cite{SU11}, \cite{HuHu}

\vspace{-1mm}
$$
\bar T^{k-2r}_{g_1,g_2,\ldots,g_{2r-1}, 
g_{2r}}\Phi(t_{q_1},\ldots,t_{q_{k-2r}})\stackrel{\sf def}{=}
$$

\begin{equation}
\label{novem2026xxx11}
\stackrel{\sf def}{=}
\lim\limits_{N\to\infty}\bar T^{k-2r, N}_{g_1,g_2,\ldots,g_{2r-1}, 
g_{2r}}\Phi(t_{q_1},\ldots,t_{q_{k-2r}})
\end{equation}

\vspace{3mm}
\noindent
in $L_2([t, T]^{k-2r})$, where 

\vspace{-4mm}
$$
\bar T^{k-2r, N}_{g_1,g_2,\ldots,g_{2r-1}, 
g_{2r}}\Phi(t_{q_1},\ldots,t_{q_{k-2r}})=
$$

$$
=
\sum\limits_{j_{g_1},j_{g_3},\ldots,j_{g_{2r-1}}=0}^{N-1}
\frac{1}{\Delta_{j_{g_1}}\Delta_{j_{g_3}},\ldots \Delta_{j_{g_{2r-1}}}}
\times
$$

\vspace{2mm}
$$
\times
\int\limits_{T_{j_{g_1}}^2\times T_{j_{g_3}}^2\times
\ldots \times T_{j_{g_{2r-1}}}^2}\Phi(t_1,\ldots,t_k) dt_{g_1}dt_{g_2}\ldots dt_{g_{2r-1}}dt_{g_{2r}},
$$

\vspace{3mm}
\noindent
$\{g_1,g_2,\ldots,g_{2r-1},g_{2r},
q_1,\ldots,q_{k-2r}\}=\{1,2,\ldots,k\}$ (see (\ref{leto5007after})),
$r=1,2,\ldots,$ $[k/2]$.

Assume that all traces
$\bar T^{k-2r}_{g_1,g_2,\ldots,g_{2r-1}, 
g_{2r}}\Phi(t_{q_1},\ldots,t_{q_{k-2r}})$
(for all $r=1,2,\ldots,$ $[k/2]$ and for all possible
$g_1,g_2,\ldots,g_{2r-1},g_{2r}$ (see {\rm (\ref{leto5007after})))
exist.

The analogue of the Hu--Meyer formula from \cite{SU11} for the multiple Stra\-to\-no\-vich 
stochastic integral (\ref{novem2026xxx10}) and the case of a multidimensional
Wiener process can be written as
$$
\bar J^{S}[\Phi]_{T,t}^{(i_1\ldots i_k)}
=
J'[\Phi]_{T,t}^{(i_1\ldots i_k)}+
$$

\vspace{-3mm}
$$
+
\sum\limits_{r=1}^{[k/2]}
\sum_{\stackrel{(\{\{g_1, g_2\}, \ldots, 
\{g_{2r-1}, g_{2r}\}\}, \{q_1, \ldots, q_{k-2r}\})}
{{}_{\{g_1, g_2, \ldots, 
g_{2r-1}, g_{2r}, q_1, \ldots, q_{k-2r}\}=\{1, 2, \ldots, k\}}}}
\prod\limits_{s=1}^r
{\bf 1}_{\{i_{g_{{}_{2s-1}}}=~i_{g_{{}_{2s}}}\ne 0\}}\times
$$

\vspace{4mm}
\begin{equation}
\label{novem2026xxx101}
\times
J'\hspace{-1mm}\left[\bar T^{k-2r}_{g_1,g_2,\ldots,g_{2r-1}, 
g_{2r}}\Phi\right]_{T,t}^{(i_{q_1}\ldots i_{q_{k-2r}})}\ \ \ \hbox{w.~p.~1,}
\end{equation}

\vspace{3mm}
\noindent
where $J'\hspace{-1mm}\left[\bar T^{k-2r}_{g_1,g_2,\ldots,g_{2r-1}, 
g_{2r}}\Phi\right]_{T,t}^{(i_{q_1}\ldots i_{q_{k-2r}})}\stackrel{\sf def}{=}
\bar T^{k-2r}_{g_1,g_2,\ldots,g_{2r-1},g_{2r}}\Phi$ for $k=2r,$ 
the existence of the right-hand side of (\ref{novem2026xxx101})
ensures the existence of the multiple Stratonovich stochastic
integral $\bar J^{S}[\Phi]_{T,t}^{(i_1\ldots i_k)};$
another notations are the same as in (\ref{2026may10}).

Assume that all limiting traces
$T^{k-2r}_{g_1,g_2,\ldots,g_{2r-1}, g_{2r}}\Phi (t_{q_1},\ldots,t_{q_{k-2r}})$
(see (\ref{novem2026xxx1122})) exist
for all $r=1,2,\ldots,$ $[k/2]$ and for all possible
$g_1,g_2,\ldots,g_{2r-1},g_{2r}$ (see {\rm (\ref{leto5007after})).
Moreover, suppose that

\vspace{-3.5mm}
$$
\bar T^{k-2r}_{g_1,g_2,\ldots,g_{2r-1}, 
g_{2r}}\Phi(t_{q_1},\ldots,t_{q_{k-2r}})=
T^{k-2r}_{g_1,g_2,\ldots,g_{2r-1}, 
g_{2r}}\Phi(t_{q_1},\ldots,t_{q_{k-2r}})
$$

\vspace{3mm}
\noindent
almost everywhere on $[t, T]^{k-2r}$ (with respect to Lebesgue's measure)
for all $r=1,2,\ldots,$ $[k/2]$ and for all possible
$g_1,g_2,\ldots,g_{2r-1},g_{2r}$ (see {\rm (\ref{leto5007after})),
where $\bar T^{k-2r}_{g_1,g_2,\ldots,g_{2r-1}, 
g_{2r}}\Phi(t_{q_1},\ldots,t_{q_{k-2r}})$
is defined by (\ref{novem2026xxx11}).

Then, using (\ref{2026may10}), (\ref{novem2026xxx101}), and (\ref{2026may2}), we obtain

\vspace{-3mm}
$$
{\sf M}\left\{\left(\bar J^{S}[\Phi]_{T,t}^{(i_1\ldots i_k)}-\hat J^{S}[\Phi]_{T,t}^{(i_1\ldots i_k)}  
\right)^2\right\}\le 
$$

$$
\le  
\bar C_k
\sum\limits_{r=1}^{[k/2]}
\sum_{\stackrel{(\{\{g_1, g_2\}, \ldots, 
\{g_{2r-1}, g_{2r}\}\}, \{q_1, \ldots, q_{k-2r}\})}
{{}_{\{g_1, g_2, \ldots, 
g_{2r-1}, g_{2r}, q_1, \ldots, q_{k-2r}\}=\{1, 2, \ldots, k\}}}}
\prod\limits_{s=1}^r
{\bf 1}_{\{i_{g_{{}_{2s-1}}}=~i_{g_{{}_{2s}}}\ne 0\}}\times
$$

$$
\Biggl.\times
\left\Vert \bar T^{k-2r}_{g_1,g_2,\ldots,g_{2r-1}, 
g_{2r}}\Phi-
T^{k-2r}_{g_1,g_2,\ldots,g_{2r-1}, 
g_{2r}}\Phi\right\Vert_{L_2([t, T]^{k-2r})}^2=0,
$$

\vspace{2mm}
\noindent
where $\bar C_k$ is a constant.

Then (see (\ref{novem2026xxx4}), (\ref{novem2026xxx5}))
$$
\bar J^{S}[\Phi]_{T,t}^{(i_1\ldots i_k)}=\hat J^{S}[\Phi]_{T,t}^{(i_1\ldots i_k)}=
\hbox{\vtop{\offinterlineskip\halign{
\hfil#\hfil\cr
{\rm l.i.m.}\cr
$\stackrel{}{{}_{p\to \infty}}$\cr
}} }
\sum_{j_1,\ldots,j_k=0}^{p}
C_{j_k\ldots j_1}
\prod_{l=1}^k \zeta_{j_l}^{(i_l)}
$$

\noindent
w.~p.~1, where 
$C_{j_k\ldots j_1}$ is the Fourier coefficient
defined by (\ref{2025may30}).

\section{Expansion of Iterated Stratonovich Stochastic Integrals of Arbitrary
Multiplicity. Proof of Hypotheses 2.1, 2.2, and 2.5}

{\bf Theorem~2.87.} {\it
Let $\{\phi_j(x)\}_{j=0}^{\infty}$ be an arbitrary 
CONS in $L^2([t, T]),$
and let $\psi_1(\tau) \equiv \dots \equiv \psi_k(\tau) \equiv 1$. 
For any multiplicity $k \in {\bf N}$ 
and any indices $i_1, \dots, i_k \in \{1, \dots, d\},$
the partial sum 
\begin{equation}
\label{ud11}
S_{T,t}^{(i_1, \dots, i_k)p}
=\sum_{j_1,\ldots,j_k=0}^{p}
C_{j_k\ldots j_1}
\prod_{l=1}^k \zeta_{j_l}^{(i_l)}
\end{equation}
with equal truncation limits $p_1 = \dots = p_k = p$ 
converges in the mean-square sense 
to the iterated Stratonovich stochastic integral 
$J^{*}[\psi^{(k)}]_{T,t}^{(i_1\ldots i_k)}$ as $p \to \infty:$
$$
\lim_{p \to \infty} 
{\sf M}\left\{
\left| S_{T,t}^{(i_1, \dots, i_k)p} 
- J^{*}[\psi^{(k)}]_{T,t}^{(i_1\ldots i_k)} \right|^2 
\right\} = 0,
$$
where $J^{*}[\psi^{(k)}]_{T,t}^{(i_1\ldots i_k)}$ is defined by {\rm (\ref{str})}
for $i_1,\ldots,i_k=1,\ldots,d,$ 
$C_{j_k \ldots j_1}$ has the form {\rm (\ref{after3000}),} 
and
$$
\zeta_{j}^{(i)}=
\int\limits_t^T \phi_{j}(\tau) d{\bf w}_{\tau}^{(i)}\ \ \ (i=1,\ldots,d,\ j=0,1,2,\ldots)
$$ 
are independent standard Gaussian random variables for various 
$i$ or $j.$ 
}

{\bf Proof.}\
Without loss of generality, 
by a linear time rescaling, 
we assume that $[t, T] = [0, 1]$, 
since all variation norms and estimates 
are invariant under linear time changes.

{\bf Step~1. Identification with the Karhunen--Lo\`eve lift.}
Let ${\bf w}_{\tau} = \bigl({\bf w}_{\tau}^{(1)}, \dots, {\bf w}_{\tau}^{(d)}\bigr)$, 
$\tau \in [0, 1]$, be a standard 
$d$-dimensional Wiener process with independent components and  ${\bf w}_0 = 0$ w.~p.~1. 
For each coordinate $\mathbf{w}^{(i)}$
($i = 1, \dots, d$), let $\mathcal{K}_i \subset L^2(\Omega)$ denote the
closed linear Gaussian subspace spanned by $\{\mathbf{w}_{\tau}^{(i)} : \tau \in [0,1]\}$.
Following Friz and Victoir (2010, Sect.~3.2) \cite{Friz}, the Cameron--Martin space
$\mathcal{H}_i$ associated with $\mathbf{w}^{(i)}$ is the Hilbert space of
paths $h : [0,1] \to \mathbf{R}$ of the form
$$
h(\tau) = \mathsf{M}\{Z \mathbf{w}_{\tau}^{(i)}\}, \quad \tau \in [0, 1]
$$
for some $Z \in \mathcal{K}_i$, equipped with the inner product
$$
\langle h, g \rangle_{\mathcal{H}_i} := \mathsf{M}\{Z Z'\},
$$
where $g(\tau) = \mathsf{M}\{Z' \mathbf{w}_{\tau}^{(i)}\}$ with $Z' \in \mathcal{K}_i$.

Since every element of $\mathcal{K}_i$ admits the Wiener integral representation
$$
Z = \int\limits_0^1 f_1(s)\, d\mathbf{w}_s^{(i)}, \ \ \ 
Z' = \int\limits_0^1 f_2(s)\, d\mathbf{w}_s^{(i)}
$$
for unique functions $f_1, f_2 \in L^2([0,1])$, the It\^o isometry yields
$$
h(\tau) = \mathsf{M}\left\{\int\limits_0^1 f_1(s)\, d\mathbf{w}_s^{(i)}
\int\limits_0^{\tau} 1\, d\mathbf{w}_s^{(i)}\right\} = \int\limits_0^{\tau} f_1(s)\, ds,
$$
and similarly 
$$
g(\tau) = \int\limits_0^{\tau} f_2(s) ds,
$$
with inner product
$$
\langle h, g \rangle_{\mathcal{H}_i}
= \mathsf{M}\left\{\int\limits_0^1 f_1(s)\, d\mathbf{w}_s^{(i)}
\int\limits_0^1 f_2(s)\, d\mathbf{w}_s^{(i)}\right\}
= \int\limits_0^1 f_1(s) f_2(s)ds.
$$

Defining the antiderivatives 
$$
h_j(\tau) = \int\limits_0^{\tau} \phi_j(s)ds,\ \ \
\tau \in [0,\,1],\ \ \  j \ge 0,
$$
we obtain 
$$
\langle h_j, h_m \rangle_{\mathcal{H}_i} 
= \int\limits_0^1 \phi_j(s)\phi_m(s)\,ds = \delta_{j,m}.
$$ 

Hence, $\{h_j\}_{j=0}^{\infty}$ forms 
an orthonormal basis of $\mathcal{H}_i$, 
and the vectors with $h_j$ in the $i$-th slot 
form an orthonormal basis of 
$$
\mathcal{H} = \bigoplus\limits_{i=1}^d \mathcal{H}_i.
$$

In Sect.~6 of Friz and Victoir (2010) \cite{Friz}, 
the general $L^2$-approximation of a Gaussian process $X_{\tau}$ 
for a subset $A \subset {\bf N}$ is defined by
$$
X_{\tau}^A := {\sf M}\bigl\{ X_{\tau} \mid \mathcal{F}_A \bigr\} 
= \sum_{j \in A} Z_j h_j(\tau),\ \ \
\mathcal{F}_A = \sigma(Z_j : j \in A).
$$

Setting $X_{\tau} = {\bf w}_{\tau}^{(i)}$, $A = A_p := \{0, 1, \dots, p\}$, 
and 
$$
Z_j = \zeta_j^{(i)} := \int\limits_0^1 \phi_j(s)d{\bf w}_{s}^{(i)},
$$ 
this specializes to the absolutely continuous 
componentwise approximation (also see (\ref{novem2026xxx1000}) for $t=0$)
$$
{\bf w}_{\tau}^{(i)p}
:= {\sf M}\bigl\{ {\bf w}_{\tau}^{(i)} \mid \mathcal{F}_{A_p}^{(i)} \bigr\}
= \sum_{j=0}^{p} \zeta_j^{(i)} h_j(\tau),\ \ \ 
i = 1, \dots, d,
$$
where $\mathcal{F}_{A_p}^{(i)} = \sigma\bigl(\zeta_j^{(i)} : 0 \le j \le p\bigr)$. 
We denote the resulting $d$-dimensional absolutely continuous approximation 
(having sample paths in $C^{1\text{-var}}([0,1],{\bf R}^d)$ w.~p.~1) by 
$$
{\bf w}_{\tau}^{p} := \bigl({\bf w}_{\tau}^{(1)p}, \dots, {\bf w}_{\tau}^{(d)p}\bigr),\ \ \
\tau \in [0, 1].
$$

Next, for any path $x \in C^{1\text{-var}}([0,\,1], {\bf R}^d)$, 
its canonical step-$N$ lift $S_N(x)$ 
into the truncated tensor algebra 
$$
T^N({\bf R}^d) = \bigoplus\limits_{m=0}^N ({\bf R}^d)^{\otimes m}
$$ 
is given by \cite{Friz} (Friz and Victoir, 2010, Sect.~1.1):
$$
S_N(x)_{\tau} 
= 1 + \sum_{m=1}^{N}\
\int\limits_{0 < s_1 < \dots < s_m < \tau} 
dx_{s_1} \otimes \dots \otimes dx_{s_m}.
$$

Let $\{e_1, \dots, e_d\}$ be the standard orthonormal basis of ${\bf R}^d$ 
given by $e_r = (\delta_{r,1}, \dots, \delta_{r,d})^{\top}$ ($r = 1, \dots, d$). 
Then, in the $k$-th tensor power 
$({\bf R}^d)^{\otimes k} = \underbrace{{\bf R}^d \otimes \dots \otimes {\bf R}^d}_{k \text{ times}}$ 
(of dimension $d^k$), 
the set $\{e_{j_1} \otimes \dots \otimes e_{j_k}\}_{j_1, \dots, j_k = 1}^d$ 
forms an orthonormal basis with respect to the canonical Euclidean structure, satisfying
$$
\bigl\langle e_{j_1} \otimes \dots \otimes e_{j_k},\, e_{i_1} \otimes \dots \otimes e_{i_k} \bigr\rangle 
= \prod_{m=1}^k \delta_{j_m, i_m}.
$$

Substituting $N = k$, $x = {\bf w}^p$, and $\tau = 1$ (noting that 
$S_k({\bf w}^p)_{0,1} := \left(S_k({\bf w}^p)_0\right)^{-1} \otimes S_k({\bf w}^p)_1 = S_k({\bf w}^p)_1$ 
since $S_k({\bf w}^p)_0 = 1$), the lift $S_k({\bf w}^p)_{0,1}$ in the direct sum 
$$
T^k({\bf R}^d) = \bigoplus\limits_{m=0}^k ({\bf R}^d)^{\otimes m}
$$ is written componentwise as
$$
\begin{aligned}
S_k({\bf w}^p)_{0,1} = 
\Biggl( 1,\; \int\limits_{0 < s_1 < 1} d{\bf w}^p_{s_1},\; 
&\int\limits_{0 < s_1 < s_2 < 1} d{\bf w}^p_{s_1} \otimes d{\bf w}^p_{s_2},\; \dots, \\
&\dots,\; \int\limits_{0 < s_1 < \dots < s_k < 1} d{\bf w}^p_{s_1} 
\otimes \dots \otimes d{\bf w}^p_{s_k} \Biggr).
\end{aligned}
$$

Applying the canonical projection $\pi_k : T^k({\bf R}^d) \to ({\bf R}^d)^{\otimes k}$ 
onto the $k$-th tensor level, defined explicitly 
for any tuple 
$$
a = (a_0, a_1, \dots, a_k) \in \bigoplus\limits_{m=0}^k ({\bf R}^d)^{\otimes m}
$$ 
by $\pi_k(a_0, a_1, \dots, a_k) = a_k$, we obtain the expansion
$$
\pi_k\bigl(S_k({\bf w}^p)_{0,1}\bigr) 
= \int\limits_{0 < s_1 < \dots < s_k < 1} 
d{\bf w}^p_{s_1} \otimes \dots \otimes d{\bf w}^p_{s_k}
= 
$$
$$
=\sum_{j_1, \dots, j_k = 1}^{d} 
S^{(j_1, \dots, j_k)p} \, e_{j_1} \otimes \dots \otimes e_{j_k},
$$
where $S^{(j_1, \dots, j_k)p}$ is defined by (\ref{ud11}) for $T=1,$ $t=0.$

Consequently, the partial sum $S^{(i_1, \dots, i_k)p}$ 
is the $(i_1, \dots, i_k)$-component of $\pi_k \bigl(S_k({\bf w}^p)_{0,1}\bigr)$:
$$
S^{(i_1, \dots, i_k)p} 
= \int\limits_{0 < s_1 < \dots < s_k < 1} 
d{\bf w}^{(i_1)p}_{s_1} \cdots d{\bf w}^{(i_k)p}_{s_k} 
= 
$$
$$
=\bigl\langle \pi_k\bigl(S_k({\bf w}^p)_{0,1}\bigr),\, 
e_{i_1} \otimes \dots \otimes e_{i_k} \bigr\rangle.
$$

\vspace{1mm}

\textbf{Step 2. Convergence at level~$2$ in $p'$-variation.} 
Fix a rough path exponent $p' \in (2, 3)$ (in the sense of Lyons, 1998, 
Sect.~1.2 and~2.2 \cite{Lyons}), so that $\lfloor p' \rfloor = 2$. 
Recall that for any $N \ge 1$ and geometric $p'$-rough paths 
$\mathbf{x}, \mathbf{y} \in C_0^{0,p'\text{-var}}\bigl(\lbrack 0, 1 \rbrack, G^N(\mathbf{R}^d)\bigr)$ 
\cite{Friz} (Friz and Victoir, 2010, Sect.~1.1, p.~372), the $p'$-variation 
distance is defined by
$$
d_{p'\text{-var};\lbrack 0, 1 \rbrack}(\mathbf{x}, \mathbf{y}) 
= \sup_{D=\{0=t_0<t_1<\dots<t_r=1\}} 
\left( \sum_{l=0}^{r-1} d(\mathbf{x}_{t_l, t_{l+1}}, \mathbf{y}_{t_l, t_{l+1}})^{p'} \right)^{1/p'},
$$
where $d(g, h) = \|g^{-1} \otimes h\|$ denotes the Carnot--Carath\'eodory 
metric on $G^N(\mathbf{R}^d)$ and $\mathbf{x}_{s,\tau} := \mathbf{x}_s^{-1} \otimes \mathbf{x}_\tau$.

For the $d$-dimensional Wiener process $\mathbf{w} = (\mathbf{w}^{(1)}, \dots, \mathbf{w}^{(d)})$, 
the covariance of each coordinate $R(s, \tau) = \min(s, \tau)$ has finite 2D $1$-variation 
($\rho = 1$, so $p' > 2\rho = 2$; see Friz and Victoir, 2010, Sect.~3.1.1 \cite{Friz}). 
By Theorem~35(i), Definition~36, and Theorem~51 (eq.~(19)) in Friz and Victoir (2010) \cite{Friz} 
applied to $X = \mathbf{w}$ and its projection approximations $X^{A_p} = \mathbf{w}^p$ 
($A_p = \{0, \dots, p\}$), projecting onto level~$2$ via the canonical contraction 
$\pi_{0,2} : G^3(\mathbf{R}^d) \to G^2(\mathbf{R}^d)$ (defined by 
$\pi_{0,2}(1, g_1, g_2, g_3) = (1, g_1, g_2)$) yields
$$
\lim_{p \to \infty} \mathsf{M}\Bigl\{ d_{p'\text{-var};\lbrack 0, 1 \rbrack}\bigl(S_2(\mathbf{w}^p), \mathbf{W}^{(2)}\bigr)^q \Bigr\} = 0 
\qquad \forall q \in \lbrack 1, \infty),
$$
where $\mathbf{W}^{(2)} := \pi_{0,2}(\mathbf{W}) \in C_0^{0,p'\text{-var}}\bigl(\lbrack 0, 1 \rbrack, G^2(\mathbf{R}^d)\bigr)$ 
and $\pi_{0,2}\bigl(S_3(\mathbf{w}^p)\bigr) = S_2(\mathbf{w}^p)$.

For any subinterval $\lbrack s, \tau \rbrack \subset \lbrack 0, 1 \rbrack$, choosing 
a partition $D$ containing $s$ and $\tau$ (e.g., $D = \{0, s, \tau, 1\}$) gives
$$
d\bigl(S_2(\mathbf{w}^p)_{s,\tau}, \mathbf{W}^{(2)}_{s,\tau}\bigr) 
\le d_{p'\text{-var};\lbrack 0, 1 \rbrack}\bigl(S_2(\mathbf{w}^p), \mathbf{W}^{(2)}\bigr) 
\xrightarrow[p \to \infty]{L^q(\Omega)} 0 
\qquad \forall q \in \lbrack 1, \infty).
$$

By the equivalence of homogeneous norms on $G^2(\mathbf{R}^d)$ \cite{Friz} (Friz and Victoir, 
2010, Sect.~1.1), 
it follows that for every $\lbrack s, \tau \rbrack \subset \lbrack 0, 1 \rbrack$, 
the second-level coordinates of $S_2(\mathbf{w}^p)_{s,\tau}$ converge in $L^2(\Omega)$ 
to those of $\mathbf{W}^{(2)}_{s,\tau}$. 

On the other hand, for any complete orthonormal 
system $\{\phi_j(x)\}_{j=0}^\infty$ in $L^2(\lbrack 0, 1 \rbrack)$ and every 
subinterval $[s, \tau]\subset [0, 1]$, the second-order 
iterated Riemann--Stieltjes integrals of the series approximations $\mathbf{w}^p$ 
(i.e. $S^{(i_1 i_2)p}_{\tau,s}$)
over $[s, \tau]$ $(0\le s<\tau\le 1)$ converge in $L^2(\Omega)$ to the $2$-fold iterated 
Stratonovich stochastic integrals over $[s, \tau].$ 

Indeed, from Theorem~1.16 for 
$\psi_1(t_1)={\bf 1}_{\{s<t_1\}},$ $\psi_2(t_2)={\bf 1}_{\{s<t_2<\tau\}},$ $p_1=p_2=p,$
we obtain 
\begin{equation}
\label{zor1q}
\int\limits_s^{\tau}
\int\limits_s^{t_2}
d{\bf w}_{t_1}^{(i_1)} d{\bf w}_{t_2}^{(i_2)}
=\hbox{\vtop{\offinterlineskip\halign{
\hfil#\hfil\cr
{\rm l.i.m.}\cr
$\stackrel{}{{}_{p\to \infty}}$\cr
}} }\sum_{j_1, j_2=0}^{p}
C_{j_2j_1}(\tau,s)\Biggl(\zeta_{j_1}^{(i_1)}\zeta_{j_2}^{(i_2)}
-{\bf 1}_{\{i_1=i_2\ne 0\}}
{\bf 1}_{\{j_1=j_2\}}\Biggr),
\end{equation}
where $i_1,i_2=1,\ldots,d,$
$$
C_{j_2j_1}(\tau,s)=
\int\limits_s^{\tau}\phi_{j_{2}}(t_{2})
\int\limits_s^{t_2}\phi_{j_{1}}(t_{1})dt_1 dt_2,
$$
$$
\zeta_{j}^{(i)}=
\int\limits_t^T \phi_{j}(s) d{\bf w}_s^{(i)} \ \ \ (i=1,\ldots,d,\ j=0,1,2,\ldots)
$$ 
are independent
standard Gaussian random variables for various 
$i$ or $j$.

Moreover, using (\ref{novorigin5}) for 
$\psi_1(t_1)={\bf 1}_{\{s<t_1\}}$ and $\psi_2(t_2)={\bf 1}_{\{s<t_2<\tau\}},$ we get w.~p.~1
$$
\int\limits_s^{\tau}
\int\limits_s^{t_2}
\circ\ d{\bf w}_{t_1}^{(i_1)}\circ d{\bf w}_{t_2}^{(i_2)}
=\int\limits_s^{\tau}
\int\limits_s^{t_2}
d{\bf w}_{t_1}^{(i_1)} d{\bf w}_{t_2}^{(i_2)}+
$$
\begin{equation}
\label{zor2q}
+
\frac{1}{2}{\bf 1}_{\{i_1=i_2\}}
\int\limits_s^{\tau}dt_1,
\end{equation}
where 
$i_1, i_2=1,\ldots,d.$

Furthermore, applying (\ref{start1000}) for
$\psi_1(t_1)={\bf 1}_{\{s<t_1\}},$ $\psi_2(t_2)={\bf 1}_{\{s<t_2<\tau\}},$ we have
\begin{equation}
\label{zor3q}
\sum_{j=0}^{\infty}\int\limits_s^{\tau}
\phi_j(t_2)
\int\limits_s^{t_2} 
\phi_j(t_1)dt_1 dt_2=
\frac{1}{2}\int\limits_s^{\tau}d\tau.
\end{equation}

Combining (\ref{zor1q}), (\ref{zor2q}), and (\ref{zor3q}), we obtain
$$
\int\limits_s^{\tau}
\int\limits_s^{t_2}
\circ\ d{\bf w}_{t_1}^{(i_1)}\circ d{\bf w}_{t_2}^{(i_2)}=
\hbox{\vtop{\offinterlineskip\halign{
\hfil#\hfil\cr
{\rm l.i.m.}\cr
$\stackrel{}{{}_{p\to \infty}}$\cr
}} }\sum_{j_1, j_2=0}^{p}
C_{j_2j_1}(\tau,s)\zeta_{j_1}^{(i_1)}\zeta_{j_2}^{(i_2)},
$$
i.e. 
the second-order 
iterated Riemann--Stieltjes integrals 
$$
S^{(i_1 i_2)p}_{\tau,s}=\sum_{j_1, j_2=0}^{p}
C_{j_2j_1}(\tau,s)\zeta_{j_1}^{(i_1)}\zeta_{j_2}^{(i_2)}
$$
of the series approximations $\mathbf{w}^p$ 
over $[s, \tau]$ $(0\le s<\tau\le 1)$ converge in $L^2(\Omega)$ to the $2$-fold iterated 
Stratonovich stochastic integrals over $[s, \tau].$ 

By uniqueness of the 
$L^2(\Omega)$-limit and path continuity, $\mathbf{W}^{(2)}$ coincides w.~p.~$1$ with 
the canonical Stratonovich Brownian rough path of step~$2$ on $\lbrack 0, 1 \rbrack$.

\vspace{1mm}

\textbf{Step 3. Extension to arbitrary multiplicity $k$.} 
By Step~2 and Markov's inequality, $d_{p'\text{-var};\lbrack 0, 1 \rbrack}\bigl(S_2(\mathbf{w}^p), \mathbf{W}^{(2)}\bigr) \xrightarrow{\mathsf{P}} 0$ 
as $p \to \infty$. Assume $k \ge 3$ (the case $k \le 2$ is already covered by Step~2). 
Since $p' \in (2, 3)$, we have $\lfloor p' \rfloor = 2$.

By Lyons' Extension Theorem \cite{Lyons} (Lyons, 1998, Theorem~2.2.1), 
the extension map
$$
\tilde{S}_k : C_0^{0,p'\text{-var}}\bigl(\lbrack 0, 1 \rbrack, G^2(\mathbf{R}^d)\bigr) 
\to C_0^{0,p'\text{-var}}\bigl(\lbrack 0, 1 \rbrack, G^k(\mathbf{R}^d)\bigr)
$$
is unique and continuous with respect to $d_{p'\text{-var};\lbrack 0, 1 \rbrack}$, 
and satisfies $\tilde{S}_k\bigl(S_2(\mathbf{w}^p)\bigr) = S_k(\mathbf{w}^p)$. 
Since the family of $k$-fold iterated Stratonovich stochastic integrals of $\mathbf{w}$
up to order $k$ has the same first two levels as $\mathbf{W}^{(2)}$, satisfies
Chen's identity, obeys the classical shuffle product rule (hence takes values in
$G^k(\mathbf{R}^d)$), and has finite $p'$-variation w.~p.~1, it is uniquely
determined on all higher levels $3, \dots, k$ by its second-level rough path
$\mathbf{W}^{(2)}$ (via Lyons' Extension Theorem
\cite{Lyons} (Lyons, 1998, Theorem~2.2.1)). Consequently,
$\pi_k(\tilde{S}_k(\mathbf{W}^{(2)})_{0,1}) \in (\mathbf{R}^d)^{\otimes k}$
coincides w.~p.~1 with the tensor $J^*$ of $k$-fold iterated Stratonovich
stochastic integrals of $\mathbf{w}$ over $[0,1]$:
$$
\pi_k\bigl(\tilde{S}_k(\mathbf{W}^{(2)})_{0,1}\bigr) 
= J^* 
= \sum_{j_1,\dots,j_k=1}^d J^{*(j_1,\dots,j_k)} e_{j_1} \otimes \dots \otimes e_{j_k},
$$
so that $J^{*(i_1,\dots,i_k)} = \bigl\langle \pi_k\bigl(\tilde{S}_k(\mathbf{W}^{(2)})_{0,1}\bigr), e_{i_1} \otimes \dots \otimes e_{i_k} \bigr\rangle$, 
where $J^*[\psi^{(k)}]_{T,t}^{(i_1 \dots i_k)}$ for $\psi_1 \equiv \dots \equiv \psi_k \equiv 1$ 
and $T = 1$, $t = 0$ is denoted by $J^{*(i_1,\dots,i_k)}$.

Applying the continuity of the extension map $\tilde{S}_k$ with respect to 
$d_{p'\text{-var};\lbrack 0, 1 \rbrack}$ (via the Continuous Mapping Theorem \cite{Mann}), we obtain
$$
d_{p'\text{-var};\lbrack 0, 1 \rbrack}\bigl(S_k(\mathbf{w}^p), \tilde{S}_k(\mathbf{W}^{(2)})\bigr) 
= d_{p'\text{-var};\lbrack 0, 1 \rbrack}
\Bigl(\tilde{S}_k\bigl(S_2(\mathbf{w}^p)\bigr), \tilde{S}_k\bigl(\mathbf{W}^{(2)}\bigr)\Bigr) 
\xrightarrow[p \to \infty]{\mathsf{P}} 0.
$$

To pass to individual tensor coordinates, recall from Sect.~1.1 of Friz and Victoir 
(2010) \cite{Friz} the equivalence of homogeneous norms on $G^k(\mathbf{R}^d)$ 
(where $\|g\| = d(1, g)$ is the Carnot--Carath\'eodory norm entering $d(g, h) = \|g^{-1} \otimes h\|$):
$$
\frac{1}{K_k} \max_{1 \le m \le k} \lvert\pi_m(g)\rvert^{1/m} 
\le \|g\| 
\le K_k \max_{1 \le m \le k} \lvert\pi_m(g)\rvert^{1/m}
$$
for some constant $K_k \ge 1$ and all $g \in G^k(\mathbf{R}^d)$. Choosing the coarse 
partition $D = \{0, 1\}$ in the definition of $d_{p'\text{-var};\lbrack 0, 1 \rbrack}$ for 
$\mathbf{x} = S_k(\mathbf{w}^p)$ and $\mathbf{y} = \tilde{S}_k(\mathbf{W}^{(2)})$ 
(since paths in $C_0^{0,p'\text{-var}}\bigl(\lbrack 0, 1 \rbrack, G^k(\mathbf{R}^d)\bigr)$ satisfy 
the initial condition $\mathbf{x}_0 = \mathbf{y}_0 = 1$, we have 
$\mathbf{x}_{0,1} = \mathbf{x}_0^{-1} \otimes \mathbf{x}_1 = 
\mathbf{x}_1$ and $\mathbf{y}_{0,1} = \mathbf{y}_1$), 
we obtain
$$
d(\mathbf{x}_{0,1}, \mathbf{y}_{0,1}) = 
\|\mathbf{x}_{0,1}^{-1} \otimes \mathbf{y}_{0,1}\| 
\le d_{p'\text{-var};\lbrack 0, 1 \rbrack}(\mathbf{x}, \mathbf{y}).
$$

Using the algebraic identity $\mathbf{y}_{0,1} - \mathbf{x}_{0,1} = 
\mathbf{x}_{0,1} \otimes (\mathbf{x}_{0,1}^{-1} \otimes \mathbf{y}_{0,1} - 1)$ 
and the multiplication rule
$$
\pi_k(\mathbf{a} \otimes \mathbf{b}) = \sum_{m=0}^k \pi_{k-m}(\mathbf{a}) \otimes \pi_m(\mathbf{b})
$$
in $T^k(\mathbf{R}^d)$ \cite{Lyons} (Lyons, 1998, eq.~(2.2)), we obtain
$$
\begin{aligned}
\pi_k(\mathbf{y}_{0,1}) - \pi_k(\mathbf{x}_{0,1}) 
&= \pi_k\Bigl(\mathbf{x}_{0,1} \otimes (\mathbf{x}_{0,1}^{-1} \otimes \mathbf{y}_{0,1} - 1)\Bigr) \\
&= \sum_{m=0}^k \pi_{k-m}(\mathbf{x}_{0,1}) \otimes \pi_m(\mathbf{x}_{0,1}^{-1} \otimes \mathbf{y}_{0,1} - 1) \\
&= \sum_{m=1}^k \pi_{k-m}(\mathbf{x}_{0,1}) \otimes \pi_m(\mathbf{x}_{0,1}^{-1} \otimes \mathbf{y}_{0,1}),
\end{aligned}
$$
since $\pi_0(\mathbf{x}_{0,1}^{-1} \otimes \mathbf{y}_{0,1} - 1) = 1 - 1 = 0$ and 
$\pi_m(g - 1) = \pi_m(g)$ for all $m \ge 1$.

Applying $\lvert\pi_m(g)\rvert \le K_k^m \|g\|^m$ (which follows from the lower bound in 
the norm equivalence above) to $g = \mathbf{x}_{0,1}^{-1} \otimes \mathbf{y}_{0,1}$ 
(so that $\|g\| = d(\mathbf{x}_{0,1}, \mathbf{y}_{0,1})$) and using the triangle inequality 
together with the submultiplicativity $\lvert u \otimes v \rvert \le \lvert u \rvert\,\lvert v \rvert$ of the tensor norm, 
we deduce the explicit estimate
$$
\bigl\lvert\pi_k(\mathbf{y}_{0,1}) - \pi_k(\mathbf{x}_{0,1})\bigr\rvert 
\le \sum_{m=1}^k K_k^m \lvert\pi_{k-m}(\mathbf{x}_{0,1})\rvert\, d(\mathbf{x}_{0,1}, \mathbf{y}_{0,1})^m 
\qquad \forall k \ge 1.
$$
Since $d_{p'\text{-var};\lbrack 0, 1 \rbrack}(\mathbf{x}, \mathbf{y}) = d_{p'\text{-var};\lbrack 0, 1 \rbrack}\bigl(S_k(\mathbf{w}^p), \tilde{S}_k(\mathbf{W}^{(2)})\bigr) \xrightarrow{\mathsf{P}} 0$, 
it follows that $d(\mathbf{x}_{0,1}, \mathbf{y}_{0,1}) \xrightarrow{\mathsf{P}} 0$ 
(and hence $d(\mathbf{x}_{0,1}, \mathbf{y}_{0,1})^m \xrightarrow{\mathsf{P}} 0$ for every $m = 1, \dots, k$). 
Since the estimate above holds at every level (with $\pi_0(\mathbf{x}_{0,1}) = 1$), 
it follows by induction on the tensor level that each lower-level projection 
$\pi_{k-m}(\mathbf{x}_{0,1})$ ($m \in \{1, \dots, k\}$) converges in probability 
to $\pi_{k-m}(\mathbf{y}_{0,1})$ and hence $\lvert\pi_{k-m}(\mathbf{x}_{0,1})\rvert$ is bounded 
in probability as $p \to \infty$. Since the product of a sequence bounded in 
probability and a sequence converging to $0$ in probability converges to $0$ 
in probability, each term in the sum and thus the entire right-hand side 
converges to $0$ in probability, so that $\lvert\pi_k(\mathbf{y}_{0,1}) - \pi_k(\mathbf{x}_{0,1})\rvert \xrightarrow{\mathsf{P}} 0$. 
Substituting back $\mathbf{x} = S_k(\mathbf{w}^p)$ and $\mathbf{y} = \tilde{S}_k(\mathbf{W}^{(2)})$, 
and noting that the $(i_1, \dots, i_k)$-coordinate of a tensor cannot exceed its 
Euclidean norm, we obtain for every multi-index $(i_1, \dots, i_k)$:
$$
\bigl\lvert S^{(i_1,\dots,i_k)p} - J^{*(i_1,\dots,i_k)} \bigr\rvert 
\le \bigl\lvert\pi_k(\mathbf{y}_{0,1}) - \pi_k(\mathbf{x}_{0,1})\bigr\rvert 
= 
$$
$$
=\bigl\lvert\pi_k\bigl(\tilde{S}_k(\mathbf{W}^{(2)})_{0,1}\bigr) - 
\pi_k\bigl(S_k(\mathbf{w}^p)_{0,1}\bigr)\bigr\rvert 
\xrightarrow[p \to \infty]{\mathsf{P}} 0.
$$

\medskip
\noindent\textbf{Step 4. Upgrading to mean-square convergence.}
By Proposition~21 in Friz and Victoir (2010) \cite{Friz}, 
the iterated integral $S^{(i_1, \dots, i_k)p}$ 
of the Gaussian process ${\bf w}^p$ (whose sample paths are absolutely continuous
and hence of finite $1$-variation w.~p.~1)
lies in the $k$-th (inhomogeneous) Wiener chaos $\mathcal{H}^{\le k}$. 
Since $\mathcal{H}^{\le k}$ is closed under convergence in probability, 
we also have $J^{*(i_1, \dots, i_k)} \in \mathcal{H}^{\le k}$, 
and thus the error 
$$
Z_p := S^{(i_1, \dots, i_k)p} - J^{*(i_1, \dots, i_k)} 
\in \mathcal{H}^{\le k}.
$$

Recall Lemma~22 in Friz and Victoir (2010) \cite{Friz}
(hypercontractivity on the inhomogeneous Wiener chaos $\mathcal{H}^{\le k}$): 
for any $Z \in \mathcal{H}^{\le k}$ and $1 < p_0 < q < \infty$,
$$ |Z|_{L^{p_0}(\Omega)} 
\le |Z|_{L^q(\Omega)} 
\le |Z|_{L^{p_0}(\Omega)} (k+1)(q-1)^{k/2} 
\max\bigl(1, (p_0 - 1)^{-k}\bigr).
$$

Furthermore, by Remark~24 in Friz and Victoir (2010) \cite{Friz}, 
on $\mathcal{H}^{\le k}$ convergence in probability 
is equivalent to convergence in $L^q(\Omega)$ 
for every $q \in (0, \infty)$:
$$
Z_p \xrightarrow[p \to \infty]{{\sf P}} 0 
\quad \Longleftrightarrow \quad 
\lim_{p \to \infty} |Z_p|_{L^q(\Omega)} = 0 
\qquad \forall\, q \in (0, \infty).
$$

Substituting $Z_p = S^{(i_1, \dots, i_k)p} - J^{*(i_1, \dots, i_k)}$ 
and $q = 2$ into this equivalence, 
and using Step~3 ($Z_p \xrightarrow{{\sf P}} 0$), 
we finally obtain
$$
\lim_{p \to \infty} |Z_p|_{L^2(\Omega)}^2 
= \lim_{p \to \infty} 
{\sf M}\left\{ 
\left| S^{(i_1, \dots, i_k)p} 
- J^{*(i_1, \dots, i_k)} \right|^2
\right\} = 0.
$$

Theorem~2.87 is proved. $\Box.$

\vspace{1mm}

{\bf Theorem~2.88\ (Extension to weighted iterated Stratonovich stochastic integrals).}
{\it Let $k \ge 1,$ 
let $\{\phi_j(x)\}_{j=0}^{\infty}$ be an arbitrary 
CONS in $L^2([0, 1]),$
$(i_1, \dots, i_k) \in \{1, \dots, d\}^k,$
and let $\psi_1(\tau), \dots,$ $\psi_k(\tau) \in C^1([0,\,1])$ be deterministic 
continuously differentiable weight functions. 
Define the weighted iterated stochastic integrals {\rm (}Riemann--Stiltjes and Stratonovich$,$ 
respectively{\rm )}
$$
S[\psi^{(k)}]^{(i_1, \dots, i_k)p}_{1,0}
:= \int\limits_{0 < t_1 < \dots < t_k < 1} 
\psi_1(t_1) \cdots \psi_k(t_k)\, 
d{\bf w}_{t_1}^{(i_1)p} \cdots d{\bf w}_{t_k}^{(i_k)p}, 
$$
$$
J^{*}[\psi^{(k)}]_{1,0}^{(i_1 \ldots i_k)}
:= \int\limits_{0 < t_1 < \dots < t_k < 1} 
\psi_1(t_1) \cdots \psi_k(t_k) 
\circ d{\bf w}_{t_1}^{(i_1)} \cdots \circ d{\bf w}_{t_k}^{(i_k)},
$$
where 
$$
S_{1,0}^{(i_1, \dots, i_k)p}
=\sum_{j_1,\ldots,j_k=0}^{p}
C_{j_k\ldots j_1}
\prod_{l=1}^k \zeta_{j_l}^{(i_l)};
$$

\noindent
another notations are the same as in Theorem~{\rm 2.87}. Then

\vspace{-3mm}
$$
\lim_{p \to \infty} 
{\sf M}\left\{
\left|S[\psi^{(k)}]^{(i_1, \dots, i_k)p}_{1,0}
- J^{*}[\psi^{(k)}]_{1,0}^{(i_1 \ldots i_k)} \right|^2 
\right\} = 0.
$$
}

{\bf Proof.} Recall that $\mathbf{w}^p = (\mathbf{w}^{(1)p}, \dots, \mathbf{w}^{(d)p})$ 
denotes the $p$-th partial sum expansion of the standard $d$-dimensional Wiener process 
$\mathbf{w}$ with independent components introduced in Theorem~2.87
(see (\ref{novem2026xxx1000}) for $t=0$). 
For each $i_m\in\{1,\ldots,d\},$ $m \in \{1, \dots, k\},$ and $t \in \lbrack 0, 1 \rbrack$, 
define the weighted process
$$
Y_t^{(i_m,m)} := \int\limits_0^t \psi_m(s)\, d\mathbf{w}_s^{(i_m)}
$$
and its absolutely continuous approximation
$$
Y_{t}^{(i_m,m)p} := \int\limits_0^t \psi_m(s)\, d\mathbf{w}_s^{(i_m)p},
$$
where $d{\bf w}_s^{(i)p}$ $(i=1,\ldots,d)$ is defined by (\ref{novem2026xxx1}).

Applying It\^o's formula to $F(t, {\bf w}_t^{(i_m)}) = \psi_m(t){\bf w}_t^{(i_m)}$ 
(where $F(t,x)=\psi_m(t)x,$ $\partial^2 F/\partial x^2 \equiv 0$) 
with $\mathbf{w}_0^{(i_m)} = 0$, we obtain w.~p.~1
$$
Y_t^{(i_m,m)} = \int\limits_0^t \psi_m(s)\, d\mathbf{w}_s^{(i_m)} 
= \psi_m(t)\mathbf{w}_t^{(i_m)} - \int\limits_0^t \psi_m'(s)\mathbf{w}_s^{(i_m)}\, ds,
$$
and similarly, by classical integration by parts for the absolutely continuous path
$\mathbf{w}_s^{(i_m)p}$ (of finite $1$-variation w.~p.~1),
$$
Y_{t}^{(i_m,m)p} = \psi_m(t)\mathbf{w}_t^{(i_m)p} - \int\limits_0^t \psi_m'(s)\mathbf{w}_s^{(i_m)p}\, ds\ \ \ 
\hbox{w.~p.~1}.
$$

Since $\psi_m(s)$ is deterministic and of bounded variation, the It\^o and Stratonovich 
differentials coincide 
($\psi_m(t)\,d\mathbf{w}_t^{(i_m)} = \psi_m(t) \circ d\mathbf{w}_t^{(i_m)} = \circ\, dY_t^{(i_m,m)}$), 
and $dY_{t}^{(i_m,m)p} = \psi_m(t)d{\bf w}_t^{(i_m,m)p}$. 
Hence the weighted iterated 
integrals rewrite as unweighted iterated integrals of the $k$-dimensional processes 
$Y^p_t = (Y^{(i_1,1)p}_t,$ $\ldots, Y^{(i_k,k)p}_t)$ and
$Y_t = (Y^{(i_1,1)}_t, \ldots, Y^{(i_k,k)}_t),$ $i_1,\ldots,i_k=1,\ldots,d$:
$$
S[\psi^{(k)}]^{(i_1, \dots, i_k)p}_{1,0}
= \int\limits_{0 < t_1 < \dots < t_k < 1} dY_{t_1}^{(i_1,1)p} \cdots dY_{t_k}^{(i_k,k)p}, 
$$
$$
J^{*}[\psi^{(k)}]_{1,0}^{(i_1 \ldots i_k)}
= \int\limits_{0 < t_1 < \dots < t_k < 1} \circ\, dY_{t_1}^{(i_1,1)} \cdots \circ dY_{t_k}^{(i_k,k)}.
$$

Note that the components $Y^{(i_m,m)}_t$ and $Y^{(i_l,l)}_t$ of
$Y_t = (Y^{(i_1,1)}_t,\ldots,Y^{(i_k,k)}_t)$
are independent whenever their first upper indices differ
($i_m \ne i_l$). Even though components with $i_m = i_l$
($m \ne l$) are correlated, the full covariance matrix
$R_Y$ has finite 2D $1$-variation ($\rho = 1$) on $[0,1]^2$ (see definition of 2D $1$-variation
in \cite{Friz}),
as shown below.

Now observe that $Y_t = (Y^{(i_1,1)}_t,\ldots,Y^{(i_k,k)}_t)$ is a centered $k$-dimensional 
Gaussian process whose covariance function $R_Y$ has finite 2D $1$-variation 
($\rho = 1$) on $[0, 1]^2$: indeed, by It\^o's isometry, for any 
rectangle $[s_1, s_2]\times [t_1, t_2] \subset [0, 1]^2$, 
the rectangular covariance increment of $(Y^{(i_m,m)}_t, Y^{(i_l,l)}_t)$ vanishes if $i_m \ne i_l$ 
and equals 
$$
\int\limits_{[s_1, s_2]\cap [t_1, t_2]} \psi_m(\tau)\psi_l(\tau) d\tau
$$ 
if $i_m = i_l$, so summing over any grid partition of $[0, 1]^2$ gives 
a total variation bounded by $\max\limits_{1 \le m \le k} \|\psi_m\|_\infty^2 < \infty$.

Next, let $\mathcal{F}_p := \sigma(\zeta_j^{(i)} : 0 \le j \le p,\, 1 \le i \le d)$ 
be the filtration generated by the first $p+1$ expansion coefficients of $\mathbf{w}_t$. 
Because $\{\phi_j(x)\}_{j=0}^\infty$ is a complete orthonormal basis of $L^2([0, 1])$, 
the coefficients $\zeta_j^{(i)}$ are independent standard Gaussian random variables 
and the $\sigma$-algebra 
$$
\mathcal{F}_\infty := \sigma\left(\bigcup\limits_{p=0}^\infty \mathcal{F}_p\right)
$$ 
coincides (up to ${\sf P}$-null sets)
with $\sigma(\mathbf{w}_t : t \in \lbrack 0, 1 \rbrack)$, so every $Y_t^{(i_m,m)}$ 
is $\mathcal{F}_\infty$-measurable. Moreover, since the tail series 
$$
\sum_{j=p+1}^\infty \zeta_j^{(i)} \int\limits_0^t \phi_j(s)ds
$$ 
is independent of $\mathcal{F}_p$ 
and has zero mean, we have w.~p.~1
$$
\mathbf{w}_t^{(i)p} = {\sf M}\left\{\mathbf{w}_t^{(i)} \mid \mathcal{F}_p\right\}
$$ 
for all $t \in \lbrack 0, 1 \rbrack$. By linearity of conditional expectation in the 
representations above, it follows that w.~p.~1
$$
Y_{t}^{(i_m,m)p} = {\sf M}\left\{Y_t^{(i_m,m)} \mid \mathcal{F}_p \right\}
$$
for all $m \in \{1, \dots, k\},$ $i_m\in\{1,\ldots,d\},$
and $t \in [0, 1]$.

In particular, $(Y_{t}^{(i_m,m)p})_{p \ge 0}$ is a closed centered Gaussian martingale with 
absolutely continuous sample paths (of finite $1$-variation w.~p.~1), 
and by L\'evy's Martingale Convergence Theorem, 
$Y_{t}^{(i_m,m)p} \to Y_t^{(i_m,m)}$ w.~p.~1 and in $L^2(\Omega)$ as $p \to \infty$ for every 
$t \in \lbrack 0, 1 \rbrack$.

Fix a rough path exponent $p' \in (2, 3)$ (so that $p' > 2\rho = 2$ and $\lfloor p' \rfloor = 2$). 
Applying Theorem~35(i), Definition~36, and Theorem~51 in Friz and Victoir (2010) \cite{Friz}
to the conditional expectation approximations $Y_{t}^p = {\sf M}\left\{Y_t \mid \mathcal{F}_p \right\}$ 
and projecting onto level~$2$ via the canonical 
contraction $\pi_{0,2} : G^3(\mathbf{R}^k) \to G^2(\mathbf{R}^k)$, 
we obtain
$$
\lim\limits_{p \to \infty} \mathsf{M}\Bigl\{ d_{p'\text{-var};\lbrack 0, 1 
\rbrack}\bigl(S_2(Y^p), \mathbf{Y}^{(2)}\bigr)^q \Bigr\} = 0 
\qquad \forall q \in \lbrack 1, \infty),
$$

\vspace{1mm}
\noindent
where $\mathbf{Y}^{(2)} := \pi_{0,2}(\mathbf{Y}) 
\in C_0^{0,p'\text{-var}}\bigl(\lbrack 0, 1 \rbrack, G^2(\mathbf{R}^k)\bigr)$. 

In particular, for every subinterval $\lbrack s, \tau \rbrack \subset \lbrack 0, 1 \rbrack$, 
the second-level coordinates of $S_2(Y^p)_{s,\tau}$ converge in $L^2(\Omega)$ to 
those of $\mathbf{Y}^{(2)}_{s,\tau}$. 

On the other hand, since the $(m, l)$-coordinate 
of $S_2(Y^p)_{s,\tau}$ is the weighted second-order iterated Riemann--Stieltjes integral 
$$
\int\limits_{s < t_1 < t_2 < \tau} 
\psi_l(t_2)\psi_m(t_1) d\mathbf{w}_{t_1}^{(i_m)p} d\mathbf{w}_{t_2}^{(i_l)p}
$$ 
with smooth weights $\psi_m(\tau), \psi_l(\tau) \in C^1(\lbrack 0, 1 \rbrack),$ then
we have
$L^2(\Omega)$-con\-ver\-gence of 
$$
\int\limits_{s < t_1 < t_2 < \tau} 
\psi_l(t_2)\psi_m(t_1)  d\mathbf{w}_{t_1}^{(i_m)p} d\mathbf{w}_{t_2}^{(i_l)p}
$$ 
to the $2$-fold iterated Stratonovich stochastic integral 
$$
\int\limits_{s < t_1 < t_2 < \tau} \psi_l(t_2)\psi_m(t_1) 
\circ d\mathbf{w}_{t_1}^{(i_m)} \circ d\mathbf{w}_{t_2}^{(i_l)} 
= \int\limits_{s < t_1 < t_2 < \tau} \circ\, dY_{t_1}^{(m)} \circ dY_{t_2}^{(l)}.
$$

Indeed, from Theorem~1.16 for 
$\tilde \psi_m(t_1)=\psi_m(t_1){\bf 1}_{\{s<t_1\}}$ and 
$\tilde \psi_l(t_2)=\psi_l(t_2){\bf 1}_{\{s<t_2<\tau\}}$ 
($\psi_m(\tau),\psi_l(\tau) \in C^1([0, 1]$) (instead of $\psi_1(t_1), \psi_2(t_2)$),
$p_1=p_2=p,$
we obtain 
$$
\int\limits_s^{\tau}\psi_l(t_2)
\int\limits_s^{t_2}\psi_m(t_1)
d{\bf w}_{t_1}^{(i_m)} d{\bf w}_{t_2}^{(i_l)}
=
$$
\begin{equation}
\label{zor1qx}
=\hbox{\vtop{\offinterlineskip\halign{
\hfil#\hfil\cr
{\rm l.i.m.}\cr
$\stackrel{}{{}_{p\to \infty}}$\cr
}} }\sum_{j_1, j_2=0}^{p}
C_{j_2j_1}(\tau,s)\Biggl(\zeta_{j_1}^{(i_m)}\zeta_{j_2}^{(i_l)}
-{\bf 1}_{\{i_m=i_l\ne 0\}}
{\bf 1}_{\{j_1=j_2\}}\Biggr),
\end{equation}
where $i_m,i_l=1,\ldots,d,$
$$
C_{j_2j_1}(\tau,s)=
\int\limits_s^{\tau}\psi_l(t_2)\phi_{j_{2}}(t_{2})
\int\limits_s^{t_2}\psi_m(t_1)\phi_{j_{1}}(t_{1})dt_1 dt_2,
$$
$$
\zeta_{j}^{(i)}=
\int\limits_t^T \phi_{j}(s) d{\bf w}_s^{(i)} \ \ \ (i=1,\ldots,d,\ j=0,1,2,\ldots)
$$ 
are independent
standard Gaussian random variables for various 
$i$ or $j$.

Moreover, using (\ref{novorigin5}) for 
$\tilde \psi_m(t_1)=\psi_m(t_1){\bf 1}_{\{s<t_1\}},$ 
$\tilde \psi_l(t_2)=\psi_l(t_2){\bf 1}_{\{s<t_2<\tau\}}$ 
($\psi_m(\tau), \psi_l(\tau) \in C^1([0, 1]$) (instead of $\psi_1(t_1), \psi_2(t_2)$), we get w.~p.~1
$$
\int\limits_s^{\tau}\psi_l(t_2)
\int\limits_s^{t_2}\psi_m(t_1)
\circ\ d{\bf w}_{t_1}^{(i_m)}\circ d{\bf w}_{t_2}^{(i_l)}
=\int\limits_s^{\tau}\psi_l(t_2)
\int\limits_s^{t_2}\psi_m(t_1)
d{\bf w}_{t_1}^{(i_m)} d{\bf w}_{t_2}^{(i_l)}+
$$
\begin{equation}
\label{zor2qx}
+
\frac{1}{2}{\bf 1}_{\{i_m=i_l\}}
\int\limits_s^{\tau} \psi_l(t_1)\psi_m(t_1)dt_1,
\end{equation}
where 
$i_m, i_l=1,\ldots,d.$

Furthermore, applying (\ref{start1000}) for
$\tilde \psi_m(t_1)=\psi_m(t_1){\bf 1}_{\{s<t_1\}}$ and
$\tilde \psi_l(t_2)=\psi_l(t_2){\bf 1}_{\{s<t_2<\tau\}}$ 
($\psi_m(\tau), \psi_l(\tau) \in C^1([0, 1]$) (instead of $\psi_1(t_1), \psi_2(t_2)$), we have
\begin{equation}
\label{zor3qx}
\sum_{j=0}^{\infty}\int\limits_s^{\tau}
\psi_l(t_2)\phi_j(t_2)
\int\limits_s^{t_2} 
\psi_m(t_1)\phi_j(t_1)dt_1 dt_2=
\frac{1}{2}\int\limits_s^{\tau} \psi_l(t_1)\psi_m(t_1) dt_1.
\end{equation}

Combining (\ref{zor1qx}), (\ref{zor2qx}), and (\ref{zor3qx}), we obtain
$$
\int\limits_s^{\tau}\psi_l(t_2)
\int\limits_s^{t_2}\psi_m(t_1)
\circ\ d{\bf w}_{t_1}^{(i_m)}\circ d{\bf w}_{t_2}^{(i_l)}=
\hbox{\vtop{\offinterlineskip\halign{
\hfil#\hfil\cr
{\rm l.i.m.}\cr
$\stackrel{}{{}_{p\to \infty}}$\cr
}} }\sum_{j_1, j_2=0}^{p}
C_{j_2j_1}(\tau,s)\zeta_{j_1}^{(i_m)}\zeta_{j_2}^{(i_l)},
$$
where
$$
\int\limits_{s < t_1 < t_2 < \tau} \psi_l(t_2) \psi_m(t_1)
\circ d\mathbf{w}_{t_1}^{(i_m)} \circ d\mathbf{w}_{t_2}^{(i_l)} 
= \int\limits_{s < t_1 < t_2 < \tau} \circ\, dY_{t_1}^{(m)} \circ dY_{t_2}^{(l)},
$$ 
i.e. 
the second-order 
iterated Riemann--Stieltjes integrals 
$$
S^{(i_1 i_2)p}_{\tau,s}=\sum_{j_1, j_2=0}^{p}
C_{j_2j_1}(\tau,s)\zeta_{j_1}^{(i_m)}\zeta_{j_2}^{(i_l)}
$$
of the series approximations $\mathbf{w}^p$ 
over $[s, \tau]$ $(0\le s<\tau\le 1)$ converge in $L^2(\Omega)$ to the $2$-fold iterated 
Stratonovich stochastic integrals over $[s, \tau].$

By uniqueness of the $L^2(\Omega)$-limit and path continuity, $\mathbf{Y}^{(2)}$ 
coincides w.~p.~$1$ with the canonical Stratonovich rough path of $Y$ of step~$2$ 
on $\lbrack 0, 1 \rbrack$.

Finally, repeating Steps~3 and~4 of the proof of Theorem~2.87 verbatim for the unweighted 
iterated integrals of $Y^p_t$ and $Y_t$ in $\mathbf{R}^k$ with multi-index $(1, 2, \dots, k)$ 
(using the uniqueness and continuity of the Lyons extension $\tilde{S}_k$ in Lyons, 
1998, Theorem~2.2.1~\cite{Lyons}, and hypercontractivity on the inhomogeneous 
Wiener chaos $\mathcal{H}^{\le k}$ 
in Friz and Victoir, 2010, Lemma~22 and Remark~24 \cite{Friz})
yields
$$
\lim_{p \to \infty} 
{\sf M}\left\{
\left|S[\psi^{(k)}]^{(i_1, \dots, i_k)p}_{1,0}
- J^{*}[\psi^{(k)}]_{1,0}^{(i_1 \ldots i_k)} \right|^2 
\right\} = 0.
$$

Theorem~2.88 is proved. $\Box$

\vspace{2mm}

\textbf{Remark~2.8.}\ {\it Since the approximation errors
$$
Z_p := S_{T,t}^{(i_1,\dots,i_k)p} - J^*[\psi^{(k)}]_{T,t}^{(i_1\dots i_k)}
$$
in Theorem~{\rm 2.87} {\rm (}for $\psi_1(\tau) \equiv \dots \equiv \psi_k(\tau) \equiv 1${\rm )} and
$$
Z_p[\psi^{(k)}] := S[\psi^{(k)}]_{1,0}^{(i_1,\dots,i_k)p}
- J^*[\psi^{(k)}]_{1,0}^{(i_1\dots i_k)}
$$
in Theorem~{\rm 2.88} belong to the $k$-th inhomogeneous Wiener chaos
$\mathcal{H}^{\le k}$ and converge to zero in probability as $p \to \infty,$
Remark~{\rm 24} in Friz and Victoir {\rm (2010) \cite{Friz}} immediately implies that both
Theorems~{\rm 9} and~{\rm 10} remain valid for convergence in $L^q(\Omega)$ for every
$q \in (0, \infty)$ {\rm (}rather than only for $q = 2${\rm )}. Namely$,$ under the
assumptions of Theorems~{\rm 2.87} and~{\rm 2.88,} respectively$,$ we have
$$
\lim_{p\to\infty} \mathsf{M}\left\{\left|
S_{T,t}^{(i_1,\dots,i_k)p} - J^*[\psi^{(k)}]_{T,t}^{(i_1\dots i_k)}
\right|^q\right\} = 0 \qquad \forall q \in (0, \infty)
$$
and
$$
\lim_{p\to\infty} \mathsf{M}\left\{\left|
S[\psi^{(k)}]_{1,0}^{(i_1,\dots,i_k)p} - J^*[\psi^{(k)}]_{1,0}^{(i_1\dots i_k)}
\right|^q\right\} = 0 \qquad \forall q \in (0, \infty).
$$
}

{\bf Remark~2.9.}\ {\it Theorems {\rm 2.87} and {\rm 2.88} prove Hypotheses {\rm 2.1, 2.2,} and {\rm 2.5}
for the case $p_1=\ldots=p_k=p,$ $i_1,\ldots,i_k=1,\ldots,d,$ 
and $\psi_1(\tau),\ldots,\psi_k(\tau)\in C^1([t, T])$
{\rm (}see Sect.~{\rm 2.5,} {\rm 2.28)}
on the mean-square convergence of expansions of iterated 
Stratonovich stochastic integrals of arbitrary multiplicity 
$(k \in {\bf N})$ with respect to components of a multidimensional 
Wiener process into generalized multiple Fourier series.
}

\section{Invariance of Expansions of Iterated It\^{o} and 
Stra\-to\-no\-vich Stochastic Integrals from Theorems 1.16 and
2.59}

In this section, we consider the invariance 
of expansions of iterated It\^{o} and 
Stra\-to\-no\-vich stochastic integrals from Theorems 1.16 and
2.59 (or 2.61).

Consider the multiple Wiener stochastic integral 
$J'[\phi_{j_1}\ldots \phi_{j_k}]_{T,t}^{(i_1\ldots i_k)}$
defined by (\ref{WiI})
$(\Phi(t_1,\ldots,t_k)=\phi_{j_1}(t_1)\ldots \phi_{j_k}(t_k))$, where
$\{\phi_j(x)\}_{j=0}^{\infty}$ is an arbitrary
complete orthonormal system  
of functions in the space $L_2([t,T]).$

Taking into account (\ref{wi1001}) and (\ref{febr5000}), we obtain
\begin{equation}
\label{zido2}
~~~~~~J'[K]_{T,t}^{(i_1\ldots i_k)} =
\hbox{\vtop{\offinterlineskip\halign{
\hfil#\hfil\cr
{\rm l.i.m.}\cr
$\stackrel{}{{}_{p\to \infty}}$\cr
}} }\sum_{j_1,\ldots,j_k=0}^{p}
C_{j_k\ldots j_1}J'[\phi_{j_1}\ldots \phi_{j_k}]_{T,t}^{(i_1\ldots i_k)}\ \ \ \hbox{w.~p.~1},
\end{equation}

\noindent
where $J'[K]_{T,t}^{(i_1\ldots i_k)}$ and 
$J'[\phi_{j_1}\ldots \phi_{j_k}]_{T,t}^{(i_1\ldots i_k)}$ are multiple Wiener stochastic integrals 
defined by (\ref{WiI}),
the function $K(t_1,\ldots,t_k)$ has the form (\ref{chain200}).

On the other hand, the expansion (\ref{july300001}) can be written as follows
\begin{equation}
\label{zido3}
~~~~~~\bar J^{S}[K]_{T,t}^{(i_1\ldots i_k)}=
\hbox{\vtop{\offinterlineskip\halign{
\hfil#\hfil\cr
{\rm l.i.m.}\cr
$\stackrel{}{{}_{p\to \infty}}$\cr
}} }
\sum\limits_{j_1,\ldots,j_k=0}^{p}
C_{j_k \ldots j_1}\bar J^{S}[\phi_{j_1}\ldots \phi_{j_k}]_{T,t}^{(i_1\ldots i_k)}\ \ \ \hbox{w.~p.~1},
\end{equation}

\noindent
where 
$\bar J^{S}[K]_{T,t}^{(i_1\ldots i_k)}$ and 
$\bar J^{S}[\phi_{j_1}\ldots \phi_{j_k}]_{T,t}^{(i_1\ldots i_k)}$
are multiple Stra\-to\-no\-vich stochastic integrals
defined by (\ref{january19}) or (\ref{novem2026xxx10}), the function
$K(t_1,\ldots,t_k)$ has the form (\ref{january19a}).

Therefore, the expansions 
(\ref{zido2}) and (\ref{zido3}) have the same form.
At that the expansion (\ref{zido2}) is formulated using 
multiple Wiener stochastic integrals
and the expansion (\ref{zido3}) is formulated using 
multiple Stra\-to\-no\-vich stochastic integrals.

The expansions 
(\ref{zido2}) and (\ref{zido3}) can be written in a
slightly different way. Using (\ref{chain104}),
we obtain
$$
J[\psi^{(k)}]_{T,t}^{(i_1\ldots i_k)} =
\hbox{\vtop{\offinterlineskip\halign{
\hfil#\hfil\cr
{\rm l.i.m.}\cr
$\stackrel{}{{}_{p\to \infty}}$\cr
}} }\sum_{j_1,\ldots,j_k=0}^{p}
C_{j_k\ldots j_1}\times
$$
\begin{equation}
\label{zido4}
~~~~~~~~~~\times
\sum\limits_{(j_1,\ldots,j_k)}
\int\limits_t^T \phi_{j_k}(t_k)
\ldots
\int\limits_t^{t_{2}}\phi_{j_{1}}(t_{1})
d{\bf w}_{t_1}^{(i_1)}\ldots
d{\bf w}_{t_k}^{(i_k)}\ \ \ {\rm w.~p.~1,}
\end{equation}

\vspace{3mm}
\noindent
where $J[\psi^{(k)}]_{T,t}^{(i_1\ldots i_k)}$ is the iterated
It\^{o} stochastic integral (\ref{wi1001}), 
$$
\sum\limits_{(j_1,\ldots,j_k)}
$$ 

\noindent
means the sum with respect to all
possible permutations
$(j_1,\ldots,j_k).$ At the same time if 
$j_r$ swapped  with $j_q$ in the permutation $(j_1,\ldots,j_k)$,
then $i_r$ swapped  with $i_q$ in the permutation
$(i_1,\ldots,i_k);$
another notations are the same as in Theorem 1.16.

The iterated Stratonovich stochastic integrals
$$
{\int\limits_t^{*}}^T \phi_{j_k}(t_k)
\ldots
{\int\limits_t^{*}}^{t_{2}}\phi_{j_{1}}(t_{1})
d{\bf w}_{t_1}^{(i_1)}\ldots
d{\bf w}_{t_k}^{(i_k)}
$$
satisfy the following equality
$$
\bar J^{S}[\phi_{j_1}\ldots \phi_{j_k}]_{T,t}^{(i_1\ldots i_k)}=
\zeta_{j_1}^{(i_1)}\ldots \zeta_{j_k}^{(i_k)}=
$$
\begin{equation}
\label{zido5}
~~~~~~=
\sum\limits_{(j_1,\ldots,j_k)}
{\int\limits_t^{*}}^T \phi_{j_k}(t_k)
\ldots
{\int\limits_t^{*}}^{t_{2}}\phi_{j_{1}}(t_{1})
d{\bf w}_{t_1}^{(i_1)}\ldots
d{\bf w}_{t_k}^{(i_k)}\ \ \ \hbox{w.~p.~1},
\end{equation}

\noindent
where $\phi_j(x)$ $(j=0,1,2,\ldots)$ are continuous functions on $[t, T],$
$$
\zeta_{j}^{(i)}=
\int\limits_t^T \phi_{j}(\tau) d{\bf w}_{\tau}^{(i)}\ \ \ (i=0, 1,\ldots,m,\ \ j=0, 1,\ldots)
$$
are independent standard Gaussian random variables
for various
$i$ or $j$ (in the case when $i\ne 0$), the expression
$$
\sum\limits_{(j_1,\ldots,j_k)}
$$ 

\noindent
has the same meaning as in (\ref{zido4}).

For the case $i_1=\ldots=i_k=0$ we obtain from (\ref{zido5})
the following well known formula from the classical
integral calculus (see (\ref{riemann}))
$$
\int\limits_{[t,T]^k}
\phi_{j_1}(t_1)\ldots \phi_{j_k}(t_k) 
dt_1\ldots dt_k=
\sum\limits_{(j_1,\ldots,j_k)}
\int\limits_t^T \phi_{j_k}(t_k)
\ldots
\int\limits_t^{t_{2}}\phi_{j_{1}}(t_{1})
dt_1\ldots
dt_k=
$$
\begin{equation}
\label{zido6}
=\sum\limits_{(t_1,\ldots,t_k)}
\int\limits_t^T\ldots \int\limits_t^{t_{2}} \phi_{j_1}(t_1)
\ldots
\phi_{j_{k}}(t_{k})
dt_1\ldots
dt_k,
\end{equation}
where
$$
\sum\limits_{(j_1,\ldots,j_k)}
$$ 

\noindent
means the sum with respect to all
possible permutations
$(j_1,\ldots,j_k)$ and
permutations $(t_1,\ldots,t_k)$ when summing
$$
\sum\limits_{(t_1,\ldots,t_k)}
$$ 

\noindent
(see (\ref{zido6})) are performed only 
in the 
va\-lues $dt_1\ldots dt_k$
(at the same time the indices near upper 
limits of integration in the
iterated integrals
are changed correspondently).

Let us check the formula (\ref{zido5})
for the cases $k=2$ and $k=3.$
Using (\ref{a2}), (\ref{zido10x}), and (\ref{zido4}) $(k=2)$, we have
$$
\sum\limits_{(j_1,j_2)}
{\int\limits_t^{*}}^T \phi_{j_2}(t_2)
{\int\limits_t^{*}}^{t_{2}}\phi_{j_{1}}(t_{1})
d{\bf w}_{t_1}^{(i_1)}d{\bf w}_{t_2}^{(i_2)}=
$$
$$
={\int\limits_t^{*}}^T \phi_{j_2}(t_2)
{\int\limits_t^{*}}^{t_{2}}\phi_{j_{1}}(t_{1})
d{\bf w}_{t_1}^{(i_1)}d{\bf w}_{t_2}^{(i_2)}+
{\int\limits_t^{*}}^T \phi_{j_1}(t_2)
{\int\limits_t^{*}}^{t_{2}}\phi_{j_{2}}(t_{1})
d{\bf w}_{t_1}^{(i_2)}d{\bf w}_{t_2}^{(i_1)}=
$$
$$
=\int\limits_t^T \phi_{j_2}(t_2)
\int\limits_t^{t_{2}}\phi_{j_{1}}(t_{1})
d{\bf w}_{t_1}^{(i_1)}d{\bf w}_{t_2}^{(i_2)}+
\int\limits_t^T \phi_{j_1}(t_2)
\int\limits_t^{t_{2}}\phi_{j_{2}}(t_{1})
d{\bf w}_{t_1}^{(i_2)}d{\bf w}_{t_2}^{(i_1)}+
$$
$$
+{\bf 1}_{\{i_1=i_2\ne 0\}}\int\limits_t^{T}\phi_{j_{1}}(t_{1})
\phi_{j_{2}}(t_{1})dt_1=
$$
$$
=\zeta_{j_1}^{(i_1)}\zeta_{j_2}^{(i_2)}-{\bf 1}_{\{i_1=i_2\ne 0\}}{\bf 1}_{\{j_1=j_2\}}+
{\bf 1}_{\{i_1=i_2\ne 0\}}{\bf 1}_{\{j_1=j_2\}}=
$$

$$
=
\zeta_{j_1}^{(i_1)}\zeta_{j_2}^{(i_2)}\ \ \ \hbox{w.~p.~1.}
$$

\vspace{4mm}

Applying (\ref{a3}), (\ref{uyes3}), (\ref{zido4}) $(k=3)$, 
and It\^{o}'s formula, we obtain w.~p.~1
$$
\sum\limits_{(j_1,j_2,j_3)}
{\int\limits_t^{*}}^T \phi_{j_3}(t_3)
{\int\limits_t^{*}}^{t_{3}}\phi_{j_{2}}(t_{2})
{\int\limits_t^{*}}^{t_{2}}\phi_{j_{1}}(t_{1})
d{\bf w}_{t_1}^{(i_1)}d{\bf w}_{t_2}^{(i_2)}
d{\bf w}_{t_3}^{(i_3)}=
$$
$$
=\sum\limits_{(j_1,j_2,j_3)}
\int\limits_t^T \phi_{j_3}(t_3)
\int\limits_t^{t_{3}}\phi_{j_{2}}(t_{2})
\int\limits_t^{t_{2}}\phi_{j_{1}}(t_{1})
d{\bf w}_{t_1}^{(i_1)}d{\bf w}_{t_2}^{(i_2)}
d{\bf w}_{t_3}^{(i_3)}+
$$
$$
+{\bf 1}_{\{i_1=i_2\ne 0\}}
\left(\int\limits_t^T
\phi_{j_3}(t_3)
\int\limits_t^{t_3}
\phi_{j_2}(t_1)\phi_{j_1}(t_1)dt_1 d{\bf w}_{t_3}^{(i_3)}+\right.
$$
$$
\left.+
\int\limits_t^T
\phi_{j_2}(t_1)\phi_{j_1}(t_1)
\int\limits_t^{t_1}
\phi_{j_3}(t_3)d{\bf w}_{t_3}^{(i_3)}dt_1\right)+
$$
$$
+{\bf 1}_{\{i_1=i_3\ne 0\}}
\left(\int\limits_t^T
\phi_{j_2}(t_2)
\int\limits_t^{t_2}
\phi_{j_1}(t_1)\phi_{j_3}(t_1)dt_1 d{\bf w}_{t_2}^{(i_2)}+\right.
$$
$$
\left.+
\int\limits_t^T
\phi_{j_3}(t_1)\phi_{j_1}(t_1)
\int\limits_t^{t_1}
\phi_{j_2}(t_2)d{\bf w}_{t_2}^{(i_2)}dt_1\right)+
$$
$$
+{\bf 1}_{\{i_2=i_3\ne 0\}}
\left(\int\limits_t^T
\phi_{j_1}(t_1)
\int\limits_t^{t_1}
\phi_{j_2}(t_3)\phi_{j_3}(t_3)dt_3 d{\bf w}_{t_1}^{(i_1)}+\right.
$$
$$
\left.+
\int\limits_t^T
\phi_{j_3}(t_3)\phi_{j_2}(t_3)
\int\limits_t^{t_3}
\phi_{j_1}(t_1)d{\bf w}_{t_1}^{(i_1)}dt_3\right)=
$$
$$
=\sum\limits_{(j_1,j_2,j_3)}
\int\limits_t^T \phi_{j_3}(t_3)
\int\limits_t^{t_{3}}\phi_{j_{2}}(t_{2})
\int\limits_t^{t_{2}}\phi_{j_{1}}(t_{1})
d{\bf w}_{t_1}^{(i_1)}d{\bf w}_{t_2}^{(i_2)}
d{\bf w}_{t_3}^{(i_3)}+
$$
$$
+{\bf 1}_{\{i_1=i_2\ne 0\}}
\left(\int\limits_t^T
\phi_{j_3}(t_3)
\int\limits_t^{t_3}
\phi_{j_2}(t_1)\phi_{j_1}(t_1)dt_1 d{\bf w}_{t_3}^{(i_3)}+\right.
$$
$$
\left.+
\int\limits_t^T
\phi_{j_3}(t_3)
\int\limits_{t_3}^{T}
\phi_{j_2}(t_1)\phi_{j_1}(t_1)dt_1 d{\bf w}_{t_3}^{(i_3)}\right)+
$$
$$
+{\bf 1}_{\{i_1=i_3\ne 0\}}
\left(\int\limits_t^T
\phi_{j_2}(t_2)
\int\limits_t^{t_2}
\phi_{j_1}(t_1)\phi_{j_3}(t_1)dt_1 d{\bf w}_{t_2}^{(i_2)}+\right.
$$
$$
\left.+
\int\limits_t^T
\phi_{j_2}(t_2)
\int\limits_{t_2}^{T}
\phi_{j_1}(t_1)\phi_{j_3}(t_1)dt_1 d{\bf w}_{t_2}^{(i_2)}\right)+
$$
$$
+{\bf 1}_{\{i_2=i_3\ne 0\}}
\left(\int\limits_t^T
\phi_{j_1}(t_1)
\int\limits_t^{t_1}
\phi_{j_2}(t_3)\phi_{j_3}(t_3)dt_3 d{\bf w}_{t_1}^{(i_1)}+\right.
$$
$$
\left.+
\int\limits_t^T
\phi_{j_1}(t_1)
\int\limits_{t_1}^{T}
\phi_{j_2}(t_3)\phi_{j_3}(t_3)dt_3 d{\bf w}_{t_1}^{(i_1)}\right)=
$$
$$
=\sum\limits_{(j_1,j_2,j_3)}
\int\limits_t^T \phi_{j_3}(t_3)
\int\limits_t^{t_{3}}\phi_{j_{2}}(t_{2})
\int\limits_t^{t_{2}}\phi_{j_{1}}(t_{1})
d{\bf w}_{t_1}^{(i_1)}d{\bf w}_{t_2}^{(i_2)}
d{\bf w}_{t_3}^{(i_3)}+
$$
$$
+{\bf 1}_{\{i_1=i_2\ne 0\}}
\int\limits_t^T
\phi_{j_3}(t_3)
\int\limits_t^{T}
\phi_{j_2}(t_1)\phi_{j_1}(t_1)dt_1 d{\bf w}_{t_3}^{(i_3)}+
$$
$$
+{\bf 1}_{\{i_1=i_3\ne 0\}}
\int\limits_t^T
\phi_{j_2}(t_2)
\int\limits_t^{T}
\phi_{j_1}(t_1)\phi_{j_3}(t_1)dt_1 d{\bf w}_{t_2}^{(i_2)}+
$$
$$
+{\bf 1}_{\{i_2=i_3\ne 0\}}
\int\limits_t^T
\phi_{j_1}(t_1)
\int\limits_t^{T}
\phi_{j_2}(t_3)\phi_{j_3}(t_3)dt_3 d{\bf w}_{t_1}^{(i_1)}=
$$
$$
=\sum\limits_{(j_1,j_2,j_3)}
\int\limits_t^T \phi_{j_3}(t_3)
\int\limits_t^{t_{3}}\phi_{j_{2}}(t_{2})
\int\limits_t^{t_{2}}\phi_{j_{1}}(t_{1})
d{\bf w}_{t_1}^{(i_1)}d{\bf w}_{t_2}^{(i_2)}
d{\bf w}_{t_3}^{(i_3)}+
$$
$$
+{\bf 1}_{\{i_1=i_2\ne 0\}}
{\bf 1}_{\{j_1=j_2\}}
\zeta_{j_3}^{(i_3)}
+
{\bf 1}_{\{i_1=i_3\ne 0\}}
{\bf 1}_{\{j_1=j_3\}}
\zeta_{j_2}^{(i_2)}
+
{\bf 1}_{\{i_2=i_3\ne 0\}}
{\bf 1}_{\{j_2=j_3\}}
\zeta_{j_1}^{(i_1)}
=
$$
$$
=\zeta_{j_1}^{(i_1)}\zeta_{j_2}^{(i_2)}\zeta_{j_3}^{(i_3)}-
$$
$$
-{\bf 1}_{\{i_1=i_2\ne 0\}}
{\bf 1}_{\{j_1=j_2\}}
\zeta_{j_3}^{(i_3)}
-
{\bf 1}_{\{i_1=i_3\ne 0\}}
{\bf 1}_{\{j_1=j_3\}}
\zeta_{j_2}^{(i_2)}
-
{\bf 1}_{\{i_2=i_3\ne 0\}}
{\bf 1}_{\{j_2=j_3\}}
\zeta_{j_1}^{(i_1)}+
$$
$$
+{\bf 1}_{\{i_1=i_2\ne 0\}}
{\bf 1}_{\{j_1=j_2\}}
\zeta_{j_3}^{(i_3)}
+
{\bf 1}_{\{i_1=i_3\ne 0\}}
{\bf 1}_{\{j_1=j_3\}}
\zeta_{j_2}^{(i_2)}
+
{\bf 1}_{\{i_2=i_3\ne 0\}}
{\bf 1}_{\{j_2=j_3\}}
\zeta_{j_1}^{(i_1)}=
$$
$$
=\zeta_{j_1}^{(i_1)}\zeta_{j_2}^{(i_2)}\zeta_{j_3}^{(i_3)}.
$$

\vspace{4mm}

Using (\ref{zido5}), we can write the expansion (\ref{zido3}) as follows
$$
J^{*}[\psi^{(k)}]_{T,t}^{(i_1\ldots i_k)}=
\hbox{\vtop{\offinterlineskip\halign{
\hfil#\hfil\cr
{\rm l.i.m.}\cr
$\stackrel{}{{}_{p\to \infty}}$\cr
}} }
\sum\limits_{j_1,\ldots,j_k=0}^{p}
C_{j_k \ldots j_1}\times
$$
\begin{equation}
\label{zido10}
~~~~~~~~~~\times
\sum\limits_{(j_1,\ldots,j_k)}
{\int\limits_t^{*}}^T \phi_{j_k}(t_k)
\ldots
{\int\limits_t^{*}}^{t_{2}}\phi_{j_{1}}(t_{1})
d{\bf w}_{t_1}^{(i_1)}\ldots
d{\bf w}_{t_k}^{(i_k)}\ \ \ \hbox{w.~p.~1},
\end{equation}

\vspace{3mm}
\noindent
where $J^{*}[\psi^{(k)}]_{T,t}^{(i_1\ldots i_k)}$ is the iterated
Stra\-to\-no\-vich stochastic integral;
another notations are the same as in (\ref{zido4}).

Obviously, the expansions 
(\ref{zido4}) and (\ref{zido10}) have the same form.
At that the expansion (\ref{zido4}) is formulated using 
iterated It\^{o} stochastic integrals
and the expansion (\ref{zido10}) is formulated using 
iterated Stra\-to\-no\-vich stochastic integrals.

\section{Expansion of Multiple Stratonovich Stochastic Integrals
of Arbitrary Multiplicity $k.$ 
The case of a multidimensional Wiener process and a smooth function 
$\Phi(t_1,\ldots,t_k)$}

As we have seen in this chapter, one of the main difficulties in 
obtaining expansions of iterated Stratonovich stochastic integrals is related to 
the properties of the kernel (\ref{ppp}).
The kernel (\ref{ppp}) is discontinuous, which causes difficulties in 
applying the theory of multiple Fourier series converging pointwise.
Moreover, the Volterra integral operator $\mathbb{V}: L_2([0,1]) \rightarrow L_2([0,1])$ 
of the form
$$
\left(\mathbb{V} f\right)(x)=\int\limits_0^x f(\tau)d\tau\ \ \ (f(\tau)\in L_2([0,1]))
$$
with the kernel 
$$
K(\tau,x)=
\left\{\begin{matrix}
1,\ &\tau < x\cr\cr
0,\ &\hbox{\rm otherwise}
\end{matrix}
\right.\ \ \ (\tau, x\in [0, 1])
$$

\vspace{1mm}
\noindent
is not a trace class operator \cite{Brisl} (see Sect.~2.27).

Thus, one cannot count on the fact that operators of the more
general form (\ref{july6999}) (from the same 
family of operators as the Volterra integral operator) 
with the kernel (\ref{ppp}) are operators of the trace class.
It is well known \cite{Brisl} that for trace class
operators the equality of matrix and integral traces holds
(this equality is very useful for expansion of iterated Stratonovich
stochastic integrals (see Sect.~2.27)).
As a result, the proof of the equalities 
of matrix and integral traces 
for Volterra--type integral operators (\ref{july6999}) 
provides a way
to calculate the matrix traces of these operators.
However, it is obvious that this is a separate problem.

Let us assume that the function $\Phi(t_1,\ldots,t_k): [t, T]^k\to {\bf R}$ satisfies 
sufficient conditions for 
its expansion into a multiple Fourier series 
\begin{equation}
\label{drdr1500}
\lim\limits_{p\to\infty}
\sum\limits_{j_1,\ldots,j_k=0}^p C_{j_k\ldots j_1}
\prod\limits_{l=1}^k \phi_{j_l}(t_l)
\end{equation}

\noindent
converging pointwise in $(t, T)^k$ to the function 
$\Phi(t_1,\ldots,t_k)$.
Also we suppose that the partial sums
$$
\sum\limits_{j_1,\ldots,j_k=0}^p C_{j_k\ldots j_1}
\prod\limits_{l=1}^k \phi_{j_l}(t_l)
$$

\vspace{1mm}
\noindent
of the multiple Fourier series (\ref{drdr1500})
have an integrable majorant on $[t, T]^k$ that does not depend
on $p.$
Here $\{\phi_j(x)\}_{j=0}^{\infty}$ is a complete orthonormal system of 
Legendre polynomials or trigonometric functions in the space $L_2([t, T])$
and
\begin{equation}
\label{drdr1200}
C_{j_k\ldots j_1}=\int\limits_{[t,T]^k}
\Phi(t_1,\ldots,t_k)\prod_{l=1}^{k}\phi_{j_l}(t_l)dt_1\ldots dt_k
\end{equation}

\vspace{1mm}
\noindent
is the Fourier coefficient.

The mentioned conditions for $k=1$ and $k\ge 2$ (trigonometric case) are given in 
Sect.~2.1.1, 2.1.2.

Consider the 
multiple Stratonovich stochastic integral (\ref{30.34}) (also see (\ref{30.34ququ}))
$$
\hbox{\vtop{\offinterlineskip\halign{
\hfil#\hfil\cr
{\rm l.i.m.}\cr
$\stackrel{}{{}_{N\to \infty}}$\cr
}} }\sum_{j_1,\ldots,j_k=0}^{N-1}
\Phi\left(\tau_{j_1},\ldots,\tau_{j_k}\right)
\prod\limits_{l=1}^k\Delta{\bf w}_{\tau_{j_l}}^{(i_l)}
\stackrel{\rm def}{=}J[\Phi]_{T,t}^{(i_1\ldots i_k)},
$$

\noindent
where we suppose that
the function $\Phi(t_1,\ldots,t_k)$ is the same as above,
${\bf w}_{\tau}$ is a random vector with 
an $m+1$ components
(${\bf w}_{\tau}^{(i)} $ $(i=1,\ldots,m)$
are independent standard Wiener processes and
${\bf w}_{\tau}^{(0)}=\tau$),
$\Delta {\bf w}_{\tau_j}^{(i)}={\bf w}_{\tau_{j+1}}^{(i)}-{\bf w}_{\tau_j}^{(i)},$
$\left\{\tau_{j}\right\}_{j=0}^{N}$ is a partition of
$[t,T]$ which satisfies the condition (\ref{1111}),
and $i_1,\ldots,i_k=0,1,\ldots,m.$

Denote
\begin{equation}
\label{drdr1201}
~~~~~~~~~~~~R_{p}(t_1,\ldots,t_k)
=
\Phi(t_1,\ldots,t_k)-
\sum_{j_1,\ldots,j_k=0}^{p}
C_{j_k\ldots j_1} \prod_{l=1}^{k} \phi_{j_l}(t_l),
\end{equation}

\noindent
where $C_{j_k\ldots j_1}$ has the form (\ref{drdr1200}).

Applying (we mean here the passage to the limit
$\lim\limits_{p\to\infty}$)
the Lebesgue's Dominated Convergence Theorem to the integrals
on the right-hand side of (\ref{udar1}) for
$R_{p}(t_1,\ldots,t_k)$ (see (\ref{drdr1201})) instead of $R_{p_1\ldots,p_k}(t_1,\ldots,t_k)$
(see (\ref{30.48})),
we obtain
$$
\lim\limits_{p\to\infty}
{\sf M}\left\{\left(J[\Phi]_{T,t}^{(i_1\ldots i_k)}-
\sum_{j_1,\ldots,j_k=0}^{p}
C_{j_k\ldots j_1} \prod_{l=1}^{k} \zeta_{j_l}^{(i_l)}
\right)^{2n}\right\}=
$$
\begin{equation}
\label{drdr1300}
=\lim\limits_{p\to\infty}
{\sf M}\left\{\left(J[R_{p}]_{T,t}^{(i_1\ldots i_k)}\right)^{2n}\right\}=
0,
\end{equation}

\vspace{2mm}
\noindent
where $n\in {\bf N}$ and
$$
\zeta_{j}^{(i)}=
\int\limits_t^T \phi_{j}(\tau) d{\bf w}_{\tau}^{(i)}
$$
are independent standard Gaussian random variables
for various
$i$ or $j$ {\rm(}in the case when $i\ne 0${\rm).}

Note that the equality (\ref{drdr1300}) will also be satisfied if the multiple Fourier 
series (\ref{drdr1500}) converges to the function $\Phi(t_1,\ldots,t_k)$ almost everywhere 
(with respect to Lebesgue's measure) on 
the hypercubes $[t, T]^{k-r}$ $(r=0,1,\ldots, [k/2])$ 
that are domains of integration for the integrals 
on the right-hand side of the inequality (\ref{udar1}).

\chapter{Integration Order Replacement Technique for 
Iterated It\^{o} Stochastic Integrals
and Iterated Stochastic Integrals with Respect to Martingales}

\vspace{6mm}

This chapter is devoted to 
the integration order replacement technique for iterated It\^{o} stochastic 
integrals
and iterated stochastic integrals with respect to martingales.
We consider the class of iterated It\^{o}
stochastic integrals, for which with probability 1 the 
formulas on integration order replacement corresponding to 
the rules of classical integral calculus are correct.
The theorems on integration order replacement for 
the class of iterated It\^{o} stochastic integrals are proved. 
Many examples of these theorems usage have been considered.
The mentioned results are generalized 
for the class of iterated stochastic integrals with respect to
martingales.

\vspace{6mm}

\section{Introduction}

\vspace{3mm}

In this chapter we performed rather laborious work connected with 
the theorems on integration order replacement for  
iterated It\^{o} stochastic integrals.
However, there may appear a
question about a practical usefulness of this theory, since the 
significant part of its conclusions directly follows from 
the It\^{o} formula.

It is not difficult to see that to obtain various relations 
for iterated It\^{o} stochastic integrals (see, for example, Sect.~3.6) 
using the It\^{o} formula, first 
of all these relations should be guessed. Then it is necessary 
to introduce corresponding It\^{o} processes and afterwards to use
the It\^{o} formula.
It is clear that this process requires 
intellectual expenses and it is not always trivial.

On the other hand, 
the technique on integration order replacement introduced in this
chapter is formally comply with the similar technique for 
Riemann integrals, although it is related to It\^{o} integrals, 
and it provides a possibility to perform transformations 
naturally (as with Riemann integrals) with iterated It\^{o} stochastic 
integrals and to obtain 
various relations for them.

So, in order to implementation of transformations of the specific 
class of It\^{o} processes, which is represented by iterated It\^{o} 
stochastic integrals, it is more naturally and easier to use 
the theorems on integration order replacement, than the It\^{o} formula.

Many examples of 
these theorems usage are presented in Sect.~3.6. 

Note that in Chapters 1, 2, and 4
the integration order replacement technique for iterated It\^{o} stochastic 
integrals has been successfully applied 
for the proof and development of the method of 
approximation of iterated It\^{o} and Stratonovich stochastic integrals
based on generalized multiple Fourier series (see Chapters 1 and 2)
as well as for the construction of the so-called
unified Taylor--It\^{o} and Taylor--Stratonovich expansions (see Chapter 4).

Let $(\Omega,{\rm F},{\sf P})$ be a complete probability
space and let $w(t,\omega):$ $[0, T]\times \Omega\rightarrow {\bf R}^1$
be the standard Wiener process
defined on the probability space $(\Omega,{\rm F},{\sf P}).$
Further, we will use the following notation:
$w(t,\omega)\stackrel{\rm def}{=}w_t$.

Let us consider the family of $\sigma$-algebras
$\left\{{\rm F}_t,\ t\in[0,T]\right\}$ defined
on the probability space $(\Omega,{\rm F},{\sf P})$ and
connected
with the Wiener process $w_t$ in such a way that

1.\ ${\rm F}_s\subset {\rm F}_t\subset {\rm F}$\ for
$s<t.$

2.\ The Wiener process $w_t$ is ${\rm F}_t$-measurable for all
$t\in[0,T].$

3.\ The process $w_{t+\Delta}-w_{t}$ for all
$t\ge 0,$ $\Delta>0$ is independent with
the events of $\sigma$-algebra
${\rm F}_{t}.$

Let us recall that the class ${\rm M}_2([0,T])$ (see Sect.~1.1.2)
consists of functions
$\xi:$ $[0,T]\times\Omega\rightarrow {\bf R}^1,$ which satisfy the
conditions:

1. The function $\xi(t,\omega)$ is 
measurable
with respect to the pair of variables
$(t,\omega).$

2. The function $\xi(t,\omega)$ is ${\rm F}_t$-measurable 
for all $t\in[0,T]$ and $\xi(\tau,\omega)$ is independent 
with increments $w_{t+\Delta}-w_{t}$ 
for $t\ge \tau,$ $\Delta>0.$

3.\ The following relation is fulfilled
$$\int\limits_0^T{\sf M}\left\{\left(\xi(t,\omega)\right)^2\right\}dt
<\infty.
$$

4.\  ${\sf M}\left\{\left(\xi(t,\omega)\right)^2\right\}<\infty$
for all $t\in[0,T].$

Let us recall (see Sect.~1.1.2) that the stochastic integrals 
\begin{equation}
\label{ball1}
\int\limits_0^T\xi_\tau dw_\tau\ \ \
\hbox{and}\ \ \  
\int\limits_0^T\xi_\tau d\tau,
\end{equation}
where $\xi_t\in {\rm M}_2([0,T])$ and the first integral in (\ref{ball1})
is the It\^{o} stochastic integral, can be defined in the mean-square
sense by the relations (\ref{10.1}) and (\ref{10.1selo}).

We will introduce the class ${\rm S}_2([0,T])$ of functions 
$\xi:$ $[0,T]\times\Omega\rightarrow
{\bf R}^1,$ which satisfy the conditions:

1.~$\xi(\tau,\omega) \in {\rm M}_2([0,T])$.

2.~$\xi(\tau,\omega)$
is the mean-square continuous random process at the interval
$[0,T].$

As we noted above,
the It\^{o} stochastic integral exists
in the mean-square sense (see (\ref{10.1})), if the random process
$\xi(\tau,\omega)\in {\rm M}_2([0,T]),$ 
i.e., perhaps this process does not satisfy 
the property of the mean-square continuity on the interval 
$[0,T].$ In this chapter we will formulate and prove
the theorems on integration order replacement for the special 
class of iterated It\^{o} stochastic integrals.
At the same time, the condition of the mean-square continuity 
of integrand in the innermost
stochastic integral will be significant.

Let us introduce the following class
of iterated stochastic integrals
$$
J[\phi,\psi^{(k)}]_{T,t}=\int\limits_{t}^{T}\psi_1(t_1)\ldots
\int\limits_t^{t_{k-1}}\psi_k(t_k)\int\limits_t^{t_k}
\phi_{\tau}dw^{(k+1)}_{\tau}dw^{(k)}
_{t_k}
\ldots dw^{(1)}_{t_1},
$$
where $\phi(\tau,\omega)\stackrel{\rm def}{=}\phi_{\tau},$ 
$\phi_\tau\in{\rm S}_2([t,T]),$ every
$\psi_l(\tau)$ $(l=1,\ldots,k)$ is a continous nonrandom function
at the interval $[t, T]$,
here and further 
$w_\tau^{(l)}=w_\tau$ or $w_\tau^{(l)}=\tau$
for $\tau\in[t,T]$
$(l=1,\ldots,k+1),$
$(\psi_1,\ldots,\psi_k)\stackrel{\rm def}{=}\psi^{(k)},$
$\psi^{(1)}\stackrel{\rm def}{=}\psi_1.$

We will call the stochastic integral $J[\phi,\psi^{(k)}]_{T,t}$
as the iterated It\^{o} stochastic integral.

It is well known that for the iterated Riemann integral 
in the case of specific conditions the formula on
integration order replacement is correct.
In particular, if the nonrandom functions
$f(x)$ and $g(x)$ are continuous at the interval $[a, b],$ then
\begin{equation}
\label{1.3000000}
\int\limits_a^b f(x)\int\limits_a^x g(y)dydx
=\int\limits_a^b g(y)\int\limits_y^b f(x)dxdy.
\end{equation}

A similar result follows from Fubini's Theorem
for $f(x),g(x)\in L_2([a,b]),$ where the integrals in 
(\ref{1.3000000}) are understood
as Lebesque integrals.

If we suppose that for the It\^{o} stochastic integral 
$$
J[\phi,\psi_1]_{T,t}=\int\limits_{t}^{T}\psi_1(s)
\int\limits_t^{s}\phi_{\tau}dw^{(2)}_{\tau}dw^{(1)}_{s}
$$ 
the formula on integration order replacement,   
which is similar to (\ref{1.3000000}), is valid, then we will have 
\begin{equation}
\label{1.4000000}
\int\limits_{t}^{T}\psi_1(s)
\int\limits_t^{s}\phi_{\tau}dw^{(2)}_{\tau}dw^{(1)}
_{s}=\int\limits_{t}^T \phi_{\tau}\int\limits_{\tau}^T\psi_1(s)
dw_s^{(1)}dw_{\tau}^{(2)}.
\end{equation}

If, in addition 
$w_s^{(1)},\ w_s^{(2)}=w_s$ $(s\in[t, T])$ in (\ref{1.4000000}), then
the stochastic process 
$$
\eta_\tau=
\phi_{\tau}\int\limits_{\tau}^T\psi_1(s)dw_s^{(1)}
$$ 
does not belong to the class ${\rm M}_2([t,T]),$
and, consequently, for the It\^{o} stochastic integral   
$$
\int\limits_t^T\eta_\tau dw_\tau^{(2)}
$$
on the right-hand side of (\ref{1.4000000})
the conditions of its existence are not fulfilled.

At the same time
\begin{equation}
\label{rrr111}
\int\limits_t^T dw_s\int\limits_t^T ds=
\int\limits_t^T (s-t)dw_s+\int\limits_t^T (w_s-w_t)ds\ \ \ \hbox{w.~p.~1},
\end{equation}
and we can obtain this equality, for example, using the It\^{o} formula,
but (\ref{rrr111}) can be considered as a result of 
integration order replacement (see below). 

Actually, we can demonstrate that 
$$
\int\limits_t^T (w_s-w_t)ds=
\int\limits_t^T\int\limits_t^s dw_\tau ds=
\int\limits_t^T\int\limits_{\tau}^T ds dw_{\tau}\ \ \ \hbox{w.~p.~1}.
$$

Then
$$
\int\limits_t^T (s-t)dw_s+\int\limits_t^T (w_s-w_t)ds=
\int\limits_t^T\int\limits_t^\tau ds dw_\tau+
\int\limits_t^T\int\limits_{\tau}^T ds dw_{\tau}=
\int\limits_t^T dw_s\int\limits_t^T ds\ \ \ \hbox{w.~p.~1}.
$$

The aim of this chapter is to establish the strict mathematical sense 
of the formula (\ref{1.4000000}) for the case
$w_s^{(1)},$ $w_s^{(2)}=w_s$ $(s\in[t, T])$ as well as its analogue
corresponding to the iterated It\^{o} stochastic integral 
$J[\phi,\psi^{(k)}]_{T,t},$ $k \ge 2.$
At that, we will use the definition of the It\^{o} stochastic integral 
which is more general than (\ref{10.1}).

Let us consider the partition 
$\tau_j^{(N)},$ $j=0, 1, \ldots, N$ of
the interval $[t,T]$ such that
\begin{equation}
\label{prgaba}
t=\tau_0^{(N)}<\tau_1^{(N)}<\ldots <\tau_N^{(N)}=T,\ \ \ \
\max\limits_{0\le j\le N-1}\left|\tau_{j+1}^{(N)}-\tau_j^{(N)}\right|\to 0\ \
\hbox{if}\ \ N\to \infty.
\end{equation}

In \cite{str} Stratonovich R.L. introduced the definition 
of the so-called 
combined stochastic integral for the specific class of 
integrated processes.
Taking this definition as a foundation, let us consider the 
following construction of stochastic integral 
\begin{equation}
\label{1.5000000}
\hbox{\vtop{\offinterlineskip\halign{
\hfil#\hfil\cr
{\rm l.i.m.}\cr
$\stackrel{}{{}_{N\to \infty}}$\cr
}} }
\sum^{N-1}_{j=0} \phi_{\tau_{j}}\left(
w_{\tau_{j+1}} - w_{\tau_{j}}
\right)\theta_{\tau_{j+1}} \stackrel {{\rm def}}{=}
\int\limits_{t}^{T}\phi_{\tau}dw_{\tau}\theta_{\tau},
\end{equation}
where $\phi_{\tau},$ $\theta_\tau\in{\rm S}_2([t,T]),$
$\{\tau_j\}_{j=0}^{N}$ is the partition 
of the interval $[t, T],$ which satisfies the condition (\ref{prgaba})
(here and sometimes further
for simplicity we write $\tau_j$ instead of $\tau_j^{(N)}$).

Further, we will prove existence of the integral
(\ref{1.5000000}) for 
$\phi_{\tau}\in{\rm S}_2([t,T])$
and $\theta_\tau$ from a little bit narrower class of processes
than ${\rm S}_2([t,T]).$ In addition, the integral defined by
(\ref{1.5000000})
will be used for the formulation and proof 
of the theorem on integration order replacement for the iterated It\^{o}
stochastic integrals $J[\phi,\psi^{(k)}]_{T,t},$ $k \ge 1.$

Note that under the appropriate conditions the following 
properties of stochastic integrals
defined by the formula (\ref{1.5000000}) can be proved
$$
\int\limits_{t}^{T} \phi_\tau dw_\tau g(\tau)=
\int\limits_{t}^{T} \phi_\tau g(\tau) dw_\tau\ \ \ \hbox{w.~p.~1}, 
$$
where $g(\tau)$ is a continuous nonrandom function at the
interval $[t,T]$,
$$
\int\limits_{t}^{T} 
\left(\alpha\phi_{\tau}+\beta\psi_{\tau}\right)dw_\tau\theta_\tau=
\alpha\int\limits_{t}^{T} \phi_\tau dw_\tau\theta_\tau+
\beta\int\limits_{t}^{T} \psi_{\tau} dw_\tau\theta_\tau\ \ \ \hbox{w.~p.~1},
$$
$$
\int\limits_{t}^{T} \phi_\tau dw_\tau
\left(\alpha\theta_\tau+\beta\psi_\tau\right)=
\alpha\int\limits_{t}^{T} \phi_\tau dw_\tau\theta_\tau+
\beta\int\limits_{t}^{T} \phi_\tau dw_\tau\psi_\tau\ \ \ \hbox{w.~p.~1},
$$
where $\alpha,$ $\beta\in{\bf R}^1.$

At that, we suppose that the stochastic processes 
$\phi_\tau,$ $\theta_\tau$, and $\psi_\tau$ are such that
the integrals included in the mentioned 
properties exist.

\section{Formulation of the Theorem on
Integration Order Replacement
for Iterated It\^{o} Stochastic Integrals of Multiplicity $k$ $(k\in{\bf N})$}

Let us define the stochastic integrals 
$\hat I[\psi^{(k)}]_{T,t},$ $k\ge 1$ of the form
$$
\hat I[\psi^{(k)}]_{T,t}=\int\limits_t^T\psi_k(t_k)dw_{t_k}^{(k)}
\int\limits_{t_k}^T\psi_{k-1}(t_{k-1})dw_{t_{k-1}}^{(k-1)}
\ldots \int\limits_{t_{2}}^T \psi_1(t_1)dw_{t_1}^{(1)}
$$
in accordance with the definition (\ref{1.5000000}) by the following 
recurrence relation
\begin{equation}
\label{2.1000000}
\hat I[\psi^{(k)}]_{T,t}
\stackrel{\rm def}{=}
\hbox{\vtop{\offinterlineskip\halign{
\hfil#\hfil\cr
{\rm l.i.m.}\cr
$\stackrel{}{{}_{N\to \infty}}$\cr
}} }
\sum^{N-1}_{l=0} \psi_k(\tau_{l})\Delta w_{\tau_l}^{(k)} 
\hat I[\psi^{(k-1)}]_{T,\tau_{l+1}},
\end{equation}
where $k\ge 1,$ 
$\hat I[\psi^{(0)}]_{T,s}\stackrel{\rm def}{=}1,$
$[s,T]\subseteq[t,T],$ here and further $\Delta w_{\tau_l}^{(i)}=
w_{\tau_{l+1}}^{(i)}-w_{\tau_l}^{(i)},$
$i=1,\ldots,k+1,$ $l=0, 1,\ldots,N-1.$

Then, we will define the iterated stochastic 
integral $\hat J[\phi,\psi^{(k)}]_{T,t},$ $k\ge 1$
$$
\hat J[\phi,\psi^{(k)}]_{T,t}=\int\limits_{t}^T \phi_s dw_s^{(k+1)}
\hat I[\psi^{(k)}]_{T,s}
$$
similarly in accordance with the definition (\ref{1.5000000})
$$
\hat J[\phi,\psi^{(k)}]_{T,t}
\stackrel{\rm def}{=}
\hbox{\vtop{\offinterlineskip\halign{
\hfil#\hfil\cr
{\rm l.i.m.}\cr
$\stackrel{}{{}_{N\to \infty}}$\cr
}} }
\sum^{N-1}_{l=0} \phi_{\tau_{l}}\Delta w_{\tau_l}^{(k+1)} 
\hat I[\psi^{(k)}]_{T,\tau_{l+1}}.
$$

Let us formulate the theorem on 
integration order replacement for
iterated It\^{o} stochastic integrals.

{\bf Theorem 3.1}\ \cite{vini} (1997), 
(also see \cite{1}-\cite{12aa}, \cite{old-art-2}, 
\cite{arxiv-25}).
{\it Suppose that 
$\phi_\tau\in{\rm S}_2([t,T])$ and every $\psi_l(\tau)$
$(l=1,\ldots,k)$ is a continuous nonrandom function
at the interval
$[t,T]$. Then, the stochastic integral
$\hat J[\phi,\psi^{(k)}]_{T,t},$ $k\ge 1$ exists and 
$$
J[\phi,\psi^{(k)}]_{T,t}=\hat J[\phi,\psi^{(k)}]_{T,t}\ \ \ 
\hbox{w.~p.~{\rm 1.}}
$$
}

\vspace{-6mm}

\section{Proof of Theorem 3.1 for the Case of 
Iterated It\^{o} Stochastic Integrals
of Multiplicity 2}

At first, let us prove Theorem 3.1 for the case $k=1.$ We have 
$$
J[\phi,\psi_1]_{T,t}\stackrel{\rm def}{=}
\hbox{\vtop{\offinterlineskip\halign{
\hfil#\hfil\cr
{\rm l.i.m.}\cr 
$\stackrel{}{{}_{N\to \infty}}$\cr 
}} }\sum_{l=0}^{N-1} \psi_1(\tau_l)\Delta w_{\tau_l}^{(1)}\int\limits_{t}
^{\tau_l}\phi_{\tau}dw_{\tau}^{(2)}=
$$
\begin{equation}
\label{2.2000000}
=
\hbox{\vtop{\offinterlineskip\halign{
\hfil#\hfil\cr
{\rm l.i.m.}\cr
$\stackrel{}{{}_{N\to \infty}}$\cr
}} }\sum_{l=0}^{N-1}\psi_1(\tau_l)\Delta w_{\tau_l}^{(1)}
\sum_{j=0}^{l-1}
\int\limits_{\tau_j}
^{\tau_{j+1}}\phi_{\tau}dw_{\tau}^{(2)},
\end{equation}
$$
\hat J[\phi,\psi_1]_{T,t}\stackrel{\rm def}{=}
\hbox{\vtop{\offinterlineskip\halign{
\hfil#\hfil\cr
{\rm l.i.m.}\cr
$\stackrel{}{{}_{N\to \infty}}$\cr
}} }\sum_{j=0}^{N-1} \phi_{\tau_j}\Delta w_{\tau_j}^{(2)}
\int\limits_{\tau_{j+1}}^T \psi_1(s)dw_{s}^{(1)}=
$$
$$
=
\hbox{\vtop{\offinterlineskip\halign{
\hfil#\hfil\cr
{\rm l.i.m.}\cr
$\stackrel{}{{}_{N\to \infty}}$\cr
}} }\sum_{j=0}^{N-1}\phi_{\tau_j}\Delta w_{\tau_j}^{(2)}
\sum_{l=j+1}^{N-1}\int\limits_{\tau_{l}}^{\tau_{l+1}}
\psi_1(s)dw_{s}^{(1)}=
$$
\begin{equation}
\label{2.3000000}
=\hbox{\vtop{\offinterlineskip\halign{
\hfil#\hfil\cr
{\rm l.i.m.}\cr
$\stackrel{}{{}_{N\to \infty}}$\cr
}} }\sum_{l=0}^{N-1} \int\limits_{\tau_{l}}^{\tau_{l+1}}
\psi_1(s)dw_{s}^{(1)}
\sum_{j=0}^{l-1}
\phi_{\tau_j}\Delta w_{\tau_j}^{(2)}.
\end{equation}

It is clear that if the difference $\varepsilon_N$ of
prelimit expressions  
on the right-hand sides of (\ref{2.2000000}) and  
(\ref{2.3000000}) tends
to zero when $N\to\infty$ in the mean-square sense, 
then the stochastic integral  
$\hat J[\phi,\psi_1]_{T,t}$ exists 
and 
$$
J[\phi,\psi_1]_{T,t}=\hat J[\phi,\psi_1]_{T,t}\ \ \ \hbox{\rm w.~p.~1.}
$$

\newpage
\noindent
\par
The difference $\varepsilon_N$ can be represented in the form 
$\varepsilon_N=\tilde\varepsilon_N+\hat\varepsilon_N$, where
$$
\tilde\varepsilon_N=\sum_{l=0}^{N-1} \psi_1(\tau_l)\Delta w_{\tau_l}^{(1)}
\sum_{j=0}^{l-1}
\int\limits_{\tau_j}
^{\tau_{j+1}}\left(\phi_{\tau}-\phi_{\tau_j}\right)dw_{\tau}^{(2)};
$$
$$
\hat\varepsilon_N=
\sum_{l=0}^{N-1} \int\limits_{\tau_{l}}^{\tau_{l+1}}
\left(\psi_1(\tau_l)-\psi_1(s)\right)dw_{s}^{(1)}
\sum_{j=0}^{l-1}
\phi_{\tau_j}\Delta w_{\tau_j}^{(2)}.
$$

We will demonstrate that w.~p.~1
$$
\hbox{\vtop{\offinterlineskip\halign{
\hfil#\hfil\cr
{\rm l.i.m.}\cr
$\stackrel{}{{}_{N\to \infty}}$\cr
}} }\varepsilon_N=0.
$$ 

In order to do it we will analyze four cases:

1.\ $w_{\tau}^{(2)}=w_{\tau},\ \Delta w_{\tau_l}^{(1)}=\Delta w_{\tau_l}.$

2.\ $w_{\tau}^{(2)}=\tau,\ \Delta w_{\tau_l}^{(1)}=\Delta w_{\tau_l}.$

3.\ $w_{\tau}^{(2)}=w_{\tau},\ \Delta w_{\tau_l}^{(1)}=\Delta \tau_l.$

4.\ $w_{\tau}^{(2)}=\tau,\ \Delta w_{\tau_l}^{(1)}=\Delta \tau_l.$

Let us recall the well known standard moment properties 
of stochastic integrals \cite{Gih1}
$$
{\sf M}\left\{\left|\int\limits_{t}^T \xi_\tau
dw_\tau\right|^{2}\right\}=
\int\limits_{t}^T {\sf M}\left\{|\xi_\tau|^{2}\right\}d\tau,
$$
\begin{equation}
\label{99.02}
{\sf M}\left\{\left|\int\limits_{t}^T \xi_\tau
d\tau\right|^{2}\right\} \le (T-t)
\int\limits_{t}^T {\sf M}\left\{|\xi_\tau|^{2}\right\}d\tau,
\end{equation}
where 
$\xi_{\tau}\in{\rm M}_2([t,T]).$

For Case 1 using 
standard moment properties for the It\^{o} 
stochastic integral as well as
mean-square continuity (which means
uniform mean-square continuity) 
of the process $\phi_\tau$ on the interval $[t, T]$, we obtain
$$
{\sf M}\left\{\left|\tilde\varepsilon_N\right|^2\right\}=
\sum_{k=0}^{N-1} \psi_1^2(\tau_k)\Delta \tau_k
\sum_{j=0}^{k-1}
\int\limits_{\tau_j}
^{\tau_{j+1}}{\sf M}\left\{\left|\phi_{\tau}-\phi_{\tau_j}\right|^2
\right\}d{\tau}<
$$
$$
<C^2\varepsilon
\sum_{k=0}^{N-1}\Delta \tau_k \sum_{j=0}^{k-1}\Delta \tau_j
<C^2\varepsilon\frac{(T-t)^2}{2},
$$

\noindent 
i.e. ${\sf M}\left\{\left|\tilde\varepsilon_N\right|^2\right\}
\to 0$ when
$N\to \infty.$ Here $\Delta\tau_j<\delta(\varepsilon),$
$j=0, 1,\ldots,N-1$
($\delta(\varepsilon)>0$ exists for any
$\varepsilon>0$ and it does not depend on $\tau$),
$|\psi_1(\tau)|<C.$

Let us consider Case 2. Using the Minkowski inequality,
uniform mean-square continuity of the process $\phi_\tau$
as well as the estimate (\ref{99.02}) for the stochastic integral,
we have
$$
{\sf M}\left\{\left|\tilde\varepsilon_N\right|^2\right\}=
\sum_{k=0}^{N-1} \psi_1^2(\tau_k)\Delta \tau_k
{\sf M}\left\{\left(\sum_{j=0}^{k-1}
\int\limits_{\tau_j}^{\tau_{j+1}}(\phi_{\tau}-\phi_{\tau_j})d\tau\right)^2
\right\}\le
$$
$$
\le 
\sum_{k=0}^{N-1} \psi_1^2(\tau_k)\Delta \tau_k
\left(\sum_{j=0}^{k-1}\left(
{\sf M}\left\{\left(
\int\limits_{\tau_j}^{\tau_{j+1}}(\phi_{\tau}-\phi_{\tau_j})d\tau\right)^2
\right\}\right)^{1/2}\ \right)^2<
$$

$$
< C^2\varepsilon
\sum_{k=0}^{N-1}\Delta \tau_k \left(\sum_{j=0}^{k-1}\Delta \tau_j\right)^2
<C^2\varepsilon \frac{(T-t)^3}{3},
$$

\noindent
i.e. ${\sf M}\left\{\left|\tilde\varepsilon_N\right|^2\right\}
\to 0$ 
when  
$N\to \infty.$ Here 
$\Delta\tau_j<\delta(\varepsilon),$ $j=0, 1,\ldots,N-1$
($\delta(\varepsilon)>0$ exists for any
$\varepsilon>0$ and it does not depend on $\tau$),
$|\psi_1(\tau)|< C.$

For Case 3 using the Minkowski inequality,
standard moment properties for the It\^{o} stochastic integral
as well as uniform mean-square continuity 
of the process $\phi_\tau$, we find
$$
{\sf M}\left\{\left|\tilde\varepsilon_N\right|^2\right\}\le
\left(
\sum_{k=0}^{N-1}\left| \psi_1(\tau_k)\right| \Delta \tau_k
\left({\sf M}\left\{\left(\sum_{j=0}^{k-1}
\int\limits_{\tau_j}^{\tau_{j+1}}(\phi_{\tau}-\phi_{\tau_j})dw_{\tau}\right)^2
\right\}\right)^{1/2}\ \right)^2=
$$
$$
=
\left(
\sum_{k=0}^{N-1}| \psi_1(\tau_k)| \Delta \tau_k
\left(\sum_{j=0}^{k-1}
\int\limits_{\tau_j}^{\tau_{j+1}}{\sf M}\left\{|\phi_{\tau}-\phi_{\tau_j}|
^2\right\}d\tau\right)^{1/2}\ \right)^2 <
$$
$$
< C^2\varepsilon
\left(\sum_{k=0}^{N-1}\Delta \tau_k \left(\sum_{j=0}^{k-1}\Delta \tau_j 
\right)^{1/2}\
\right)^2
<C^2\varepsilon\frac{4(T-t)^3}{9},
$$

\noindent
i.e. ${\sf M}\left\{\left|\tilde\varepsilon_N\right|^2\right\}
\to 0$ 
when  
$N\to \infty.$ Here 
$\Delta\tau_j<\delta(\varepsilon),$ $j=0, 1,\ldots,N-1$
($\delta(\varepsilon)>0$ exists for any
$\varepsilon>0$ and it does not depend on $\tau$),
$|\psi_1(\tau)|< C.$

Finally, for Case 4 using 
the Minkowski inequality, 
uniform mean-square continuity of the process $\phi_\tau$
as well as the estimate (\ref{99.02}) for the stochastic 
integral, 
we obtain
$$
{\sf M}\left\{\left|\tilde\varepsilon_N\right|^2\right\}\le
\left(
\sum_{k=0}^{N-1} \sum_{j=0}^{k-1}| \psi_1(\tau_k)| \Delta \tau_k
\left({\sf M}\left\{\left(
\int\limits_{\tau_j}^{\tau_{j+1}}(\phi_{\tau}-\phi_{\tau_j})d\tau\right)^2
\right\}\right)^{1/2}\ \right)^2 <
$$
$$
< C^2\varepsilon
\left(\sum_{k=0}^{N-1}\Delta \tau_k \sum_{j=0}^{k-1}\Delta \tau_j
\right)^2
<C^2\varepsilon\frac{(T-t)^4}{4},
$$

\noindent
i.e. ${\sf M}\left\{\left|\tilde\varepsilon_N\right|^2\right\} \to 0$
when  
$N\to \infty.$ Here 
$\Delta\tau_j<\delta(\varepsilon),$ $j=0, 1,\ldots,N-1$
($\delta(\varepsilon)>0$ exists for any
$\varepsilon>0$ and it does not depend on $\tau$),\
$|\psi_1(\tau)|< C.$

Thus, we have proved that w.~p.~1
$$
\hbox{\vtop{\offinterlineskip\halign{
\hfil#\hfil\cr
{\rm l.i.m.}\cr
$\stackrel{}{{}_{N\to \infty}}$\cr
}} }\tilde\varepsilon_N=0.
$$

Analogously, taking into account 
the uniform continuity of the function $\psi_1(\tau)$
on the interval $[t, T]$, we can demonstrate that w.~p.~1
$$
\hbox{\vtop{\offinterlineskip\halign{
\hfil#\hfil\cr
{\rm l.i.m.}\cr
$\stackrel{}{{}_{N\to \infty}}$\cr
}} }\hat\varepsilon_N=0.
$$ 

Consequently,
$$
\hbox{\vtop{\offinterlineskip\halign{
\hfil#\hfil\cr
{\rm l.i.m.}\cr
$\stackrel{}{{}_{N\to \infty}}$\cr
}} }\varepsilon_N=0\ \ \ \hbox{w.~p.~1}.
$$

Theorem 3.1 is proved for the case $k=1$.

{\bf Remark 3.1.}\ {\it Proving Theorem {\rm 3.1,} we used the fact that 
if the stochastic process 
$\phi_t$ is
mean-square continuous at the interval $[t, T],$ 
then it is uniformly mean-square continuous at
this interval, i.e. $\forall$ $\varepsilon>0$
$\exists$ $\delta(\varepsilon)>0$ such that
for all $t_1, t_2\in [t, T]$ satisfying 
the condition $|t_1-t_2|<\delta(\varepsilon)$ the inequality 
$$
{\sf M}\left\{\left|\phi_{t_1}-\phi_{t_2}\right|^2\right\}<\varepsilon
$$
is fulfilled
{\rm (}here $\delta(\varepsilon)$ does not depend on $t_1$ and $t_2${\rm )}.

{\bf Proof.}\ Suppose that the stochastic process $\phi_t$ is mean-square 
continuous at the interval $[t, T]$, but not 
uniformly mean-square continuous 
at this interval. Then for some $\varepsilon>0$
and $\forall$ $\delta(\varepsilon)>0$ $\exists$ $t_1, t_2\in[t, T]$ 
such that $|t_1-t_2|<\delta(\varepsilon),$ but
$$
{\sf M}\left\{\left|\phi_{t_1}-\phi_{t_2}\right|^2\right\}\ge\varepsilon.
$$

Consequently, for $\delta=\delta_n=1/n$ $(n\in{\bf N})$
$\exists$ $t_1^{(n)},$ $t_2^{(n)}\in[t, T]$
such that 
$$
\left|t_1^{(n)}-t_2^{(n)}\right|<\frac{1}{n},
$$
but 
$$
{\sf M}\left\{\left|\phi_{t_1^{(n)}}-\phi_{t_2^{(n)}}\right|^2\right\}
\ge\varepsilon.
$$

The sequence $t_1^{(n)}$ $(n\in{\bf N})$
is bounded, consequently, according to the Bol\-za\-no--Wei\-er\-strass 
Theorem,
we can choose from it the 
subsequence $t_1^{(k_n)}$ $(n\in{\bf N})$
that converges
to a certain number
$\tilde t$ {\rm (}it is simple to demonstrate that 
$\tilde t\in[t, T]${\rm )}.
Similarly to it and in virtue of the inequality 
$$
\left|t_1^{(n)}-t_2^{(n)}\right|<\frac{1}{n}
$$
we have $t_2^{(k_n)}\to\tilde t$ when $n\to\infty.$

According to the mean-square continuity of the process $\phi_t$
at the moment $\tilde t$ and the elementary inequality
$(a+b)^2\le 2(a^2+b^2)$, we obtain
$$
0\le {\sf M}\left\{\left|\phi_{t_1^{(k_n)}}-\phi_{t_2^{(k_n)}}
\right|^2\right\}\le
$$
$$
\le
2\left({\sf M}\left\{\left|\phi_{t_1^{(k_n)}}-\phi_{\tilde t}\right|^2\right\}+
{\sf M}\left\{\left|\phi_{t_2^{(k_n)}}-\phi_{\tilde t}\right|^2\right\}
\right)\to 0
$$
when $n\to\infty.$ Then 
$$
\lim\limits_{n\to\infty}{\sf M}
\left\{\left|\phi_{t_1^{(k_n)}}-\phi_{t_2^{(k_n)}}\right|^2\right\}=0.
$$

It is impossible by virtue of the fact that 
$$
{\sf M}\left\{\left|\phi_{t_1^{(k_n)}}-
\phi_{t_2^{(k_n)}}\right|^2\right\}\ge\varepsilon>0.
$$

The obtained contradiction proves the required statement.
}

\section{Proof of Theorem 3.1 for the 
Case of Iterated It\^{o} Stochastic Integrals
of Multiplicity $k$ $(k\in {\bf N})$}

Let us prove Theorem 3.1 for the case $k>1.$
In order to do it  
we will introduce the following notations
$$
I[\psi_q^{(r+1)}]_{\theta,s}
\stackrel{\rm def}{=}
\int\limits_{s}^{\theta}\psi_{q}(t_1)\ldots
\int\limits_{s}^{t_{r}} \psi_{q+r}(t_{r+1})
dw_{t_{r+1}}^{(q+r)}
\ldots dw_{t_1}^{(q)},
$$
$$
J[\phi,\psi_q^{(r+1)}]_{\theta,s}
\stackrel{\rm def}{=}
\int\limits_{s}^{\theta}\psi_{q}(t_1)\ldots
\int\limits_{s}^{t_{r}} \psi_{q+r}(t_{r+1})
\int\limits_s^{t_{r+1}}\phi_\tau dw_\tau^{(q+r+1)}
dw_{t_{r+1}}^{(q+r)}
\ldots dw_{t_1}^{(q)},
$$
$$
G[\psi_q^{(r+1)}]_{n,m}=\sum_{j_q=m}^{n-1}\sum_{j_{q+1}=m}^{j_q-1}
\ldots \sum_{j_{q+r}=m}^{j_{q+r-1}-1}
\prod_{l=q}^{r+q}I[\psi_l]_{\tau_{j_l+1},\tau_{j_l}},
$$

\vspace{2mm}
$$
(\psi_q,\ldots,\psi_{q+r})\stackrel{\rm def}{=}\psi_q^{(r+1)},\ \ \
\psi_q^{(1)}\stackrel{\rm def}{=}\psi_q,
$$

\vspace{-3mm}
$$
(\psi_1,\ldots,\psi_{r+1})\stackrel{\rm def}{=}\psi_1^{(r+1)},\ \ \ 
\psi_1^{(r+1)}\stackrel{\rm def}{=}\psi^{(r+1)}.
$$

\vspace{2mm}

Note that according to notations introduced above, we have
$$
I[\psi_l]_{s,\theta}=\int\limits_{\theta}^{s}
\psi_l(\tau)dw_{\tau}^{(l)}.
$$

To prove Theorem 3.1 for $k>1$
it is enough to show
that
\begin{equation}
\label{2.4000000}
~~~~~~~~ J[\phi,\psi^{(k)}]_{T,t}=
\hbox{\vtop{\offinterlineskip\halign{
\hfil#\hfil\cr
{\rm l.i.m.}\cr
$\stackrel{}{{}_{N\to \infty}}$\cr
}} }
S[\phi,\psi^{(k)}]_N=
\hat J[\phi,\psi^{(k)}]_{T,t}\ \ \ \hbox{w. p. 1},
\end{equation}
where
$$
S[\phi,\psi^{(k)}]_N=
G[\psi^{(k)}]_{N,0}
\sum_{l=0}^{j_k-1}\phi_{\tau_l}\Delta w_{\tau_l}^{(k+1)},
$$
where $\Delta w_{\tau_l}^{(k+1)}=w_{\tau_{l+1}}^{(k+1)}-w_{\tau_l}^{(k+1)}$.

At first, let us prove the right equality in (\ref{2.4000000}).
We have
\begin{equation}
\label{2.5000000}
\hat J[\phi,\psi^{(k)}]_{T,t}
\stackrel{\rm def}{=}
\hbox{\vtop{\offinterlineskip\halign{
\hfil#\hfil\cr
{\rm l.i.m.}\cr
$\stackrel{}{{}_{N\to \infty}}$\cr
}} }
\sum^{N-1}_{l=0} \phi_{\tau_l}\Delta w_{\tau_l}^{(k+1)}
\hat I[\psi^{(k)}]_{T,\tau_{l+1}}.
\end{equation}

On the basis of the inductive hypothesis we obtain that
\begin{equation}
\label{2.6000000}
I[\psi^{(k)}]_{T,\tau_{l+1}}=\hat I[\psi^{(k)}]_{T,\tau_{l+1}}\ \ \ 
\hbox{w. p. 1},
\end{equation}
where $\hat I[\psi^{(k)}]_{T,s}$ 
is defined in accordance with (\ref{2.1000000}) and
$$
I[\psi^{(k)}]_{T,s}=\int\limits_{s}^T \psi_1(t_1)\ldots
\int\limits_{s}^{t_{k-2}} \psi_{k-1}(t_{k-1})
\int\limits_{s}^{t_{k-1}} \psi_k(t_k)dw_{t_k}^{(k)}dw_{t_{k-1}}^{(k-1)}
\ldots dw_{t_1}^{(1)}.
$$

Let us note that when $k\ge 4$ (for $k=2,$ $3$ the arguments are similar)
due to additivity of the It\^{o} stochastic 
integral the following equalities are correct
$$
I[\psi^{(k)}]_{T,\tau_{l+1}}
=\sum_{j_1=l+1}^{N-1}\int\limits_{\tau_{j_1}}^{\tau_{j_1+1}}
\psi_1(t_1)\int\limits_{\tau_{l+1}}^{t_1}\psi_{2}(
t_2)I[\psi_3^{(k-2)}]_{t_2,\tau_{l+1}}
dw_{t_2}^{(2)}dw_{t_1}^{(1)}=
$$
$$
=\sum_{j_1=l+1}^{N-1}\int\limits_{\tau_{j_1}}^{\tau_{j_1+1}}
\psi_1(t_1)\left(\sum_{j_2=l+1}^{j_1-1}
\int\limits_{\tau_{j_2}}^{\tau_{j_2+1}}
+\int\limits_{\tau_{j_1}}^{t_1}\right) 
\psi_{2}(t_2)I[\psi_3^{(k-2)}]_{t_2,\tau_{l+1}}
dw_{t_2}^{(2)}dw_{t_1}^{(1)}=
$$

\begin{equation}
\label{2.7000000}
=\ldots=G[\psi^{(k)}]_{N,l+1}+H[\psi^{(k)}]_{N,l+1}\ \ \ \hbox{w. p. 1},
\end{equation}

\noindent
where
$$
H[\psi^{(k)}]_{N,l+1}=
\sum_{j_1=l+1}^{N-1}\int\limits_{\tau_{j_1}}^{\tau_{j_1+1}}
\psi_1(s)\int\limits_{\tau_{j_1}}^{s}
\psi_2(\tau)I[\psi_3^{(k-2)}]_{\tau,\tau_{l+1}}dw_\tau^{(2)}
dw_s^{(1)} +
$$
$$
+ \sum_{r=2}^{k-2}
G[\psi^{(r-1)}]_{N,l+1}
\sum_{j_r=l+1}^{j_{r-1}-1}\int\limits_{\tau_{j_r}}^{\tau_{j_r+1}}
\psi_r(s)\int\limits_{\tau_{j_r}}^{s}
\psi_{r+1}(\tau)I[\psi_{r+2}^{(k-r-1)}]_{\tau,\tau_{l+1}}
dw_\tau^{(r+1)}
dw_s^{(r)} +
$$
\begin{equation}
\label{2.8000000}
+ G[\psi^{(k-2)}]_{N,l+1}
\sum_{j_{k-1}=l+1}^{j_{k-2}-1}I[\psi_{k-1}^{(2)}]_{\tau_{j_{k-1}+1},
\tau_{j_{k-1}}}.
\end{equation}

\vspace{2mm}

Next, substitute (\ref{2.7000000}) into (\ref{2.6000000}) and
(\ref{2.6000000}) into (\ref{2.5000000}). Then w. p. 1
\begin{equation}
\label{2.9000000}
~~~~~ \hat J[\phi,\psi^{(k)}]_{T,t}=
\hbox{\vtop{\offinterlineskip\halign{
\hfil#\hfil\cr
{\rm l.i.m.}\cr
$\stackrel{}{{}_{N\to \infty}}$\cr
}} }
\sum_{l=0}^{N-1} \phi_{\tau_l}\Delta w_{\tau_l}^{(k+1)}
\biggl(G[\psi^{(k)}]_{N,l+1}+
H[\psi^{(k)}]_{N,l+1}\biggr).
\end{equation}

Since
\begin{equation}
\label{2.10000001}
\sum_{j_1=0}^{N-1} \sum_{j_2=0}^{j_1-1}\ldots
\sum_{j_k=0}^{j_{k-1}-1} 
a_{j_1\ldots j_k}=\sum_{j_k=0}^{N-1} 
\sum_{j_{k-1}=j_k+1}^{N-1}\ldots \sum_{j_1=j_2+1}^{N-1}
a_{j_1\ldots j_k},
\end{equation}

\noindent
where $a_{j_1\ldots j_k}$ are scalars,
then
\begin{equation}
\label{2.11000000}
G[\psi^{(k)}]_{N,l+1}=
\sum^{N-1}_{j_k=l+1}\ldots
\sum_{j_1=j_2+1}^{N-1}
\prod_{l=1}^k I[\psi_l]_{\tau_{j_l+1},\tau_{j_l}}.
\end{equation}

Let us substitute (\ref{2.11000000}) into 
$$
\sum\limits_{l=0}^{N-1}
\phi_{\tau_l}\Delta w_{\tau_l}^{(k+1)}
G[\psi^{(k)}]_{N,l+1}
$$
and use again 
the formula (\ref{2.10000001}). Then 
\begin{equation}
\label{2.12000000}
\sum\limits_{l=0}^{N-1}
\phi_{\tau_l}\Delta w_{\tau_l}^{(k+1)}
G[\psi^{(k)}]_{N,l+1}=S[\phi,\psi^{(k)}]_{N}.
\end{equation}

Suppose that the limit
\begin{equation}
\label{li}
\hbox{\vtop{\offinterlineskip\halign{
\hfil#\hfil\cr
{\rm l.i.m.}\cr
$\stackrel{}{{}_{N\to \infty}}$\cr
}} }S[\phi,\psi^{(k)}]_{N}
\end{equation} 

\noindent
exists (its existence will be proved further).

Then from (\ref{2.12000000}) and (\ref{2.9000000}) it follows that
for proof of the right equality in  
(\ref{2.4000000}) we have to demonstrate that w. p. 1
\begin{equation}
\label{2.13000000}
\hbox{\vtop{\offinterlineskip\halign{
\hfil#\hfil\cr
{\rm l.i.m.}\cr
$\stackrel{}{{}_{N\to \infty}}$\cr
}} }
\sum_{l=0}^{N-1}
\phi_{\tau_l}\Delta w_{\tau_l}^{(k+1)}
H[\psi^{(k)}]_{N,l+1}=0.
\end{equation}

Analyzing the second moment of the prelimit expression on the 
left-hand side of (\ref{2.13000000}) and taking into account
(\ref{2.8000000}), the independence of
$\phi_{\tau_l},$ $\Delta w_{\tau_l}^{(k+1)},$ and
$H[\psi^{(k)}]_{N,l+1}$ as well as the
standard estimates for second moments
of stochastic integrals and the Minkowski inequality,
we find that (\ref{2.13000000})
is correct.
Thus, by
the assumption of existence of the limit 
(\ref{li}) 
we obtain that the right equality in (\ref{2.4000000}) is fulfilled.

Let us demonstrate that the left equality in  
(\ref{2.4000000}) is also fulfilled.

We have 
\begin{equation}
\label{2.14000000}
J[\phi,\psi^{(k)}]_{T,t}\stackrel{\rm def}{=}
\hbox{\vtop{\offinterlineskip\halign{
\hfil#\hfil\cr
{\rm l.i.m.}\cr
$\stackrel{}{{}_{N\to \infty}}$\cr
}} }\sum_{l=0}^{N-1}\psi_1(\tau_l)\Delta w_{\tau_l}^{(1)} 
J[\phi,\psi_2^{(k-1)}]_{\tau_l,t}.
\end{equation}

Let us use for
the integral $J[\phi,\psi_2^{(k-1)}]_{\tau_l,t}$
in (\ref{2.14000000}) the same arguments, 
which resulted to the relation
(\ref{2.7000000})
for the integral
$I[\psi^{(k)}]_{T,\tau_{l+1}}$. 
After that let us substitute 
the expression obtained for
the integral $J[\phi,\psi_2^{(k-1)}]_{\tau_l,t}$
into (\ref{2.14000000}). 
Further, using 
the Minkowski inequality and standard 
estimates for second moments of stochastic integrals
it is easy to obtain that 
\begin{equation}
\label{2.15000000}
J[\phi,\psi^{(k)}]_{T,t}=
\hbox{\vtop{\offinterlineskip\halign{
\hfil#\hfil\cr
{\rm l.i.m.}\cr
$\stackrel{}{{}_{N\to \infty}}$\cr
}} }R[\phi,\psi^{(k)}]_N\ \ \ \hbox{w. p. 1},
\end{equation}
where
$$
R[\phi,\psi^{(k)}]_N=
\sum^{N-1}_{j_1=0}
\psi_1(\tau_{j_1})\Delta w_{\tau_{j_1}}^{(1)}
G[\psi_2^{(k-1)}]_{j_1,0}
\sum_{l=0}^{j_k-1}
\int\limits_{\tau_l}^{\tau_{l+1}}\phi_{\tau}dw_{\tau}^{(k+1)}.
$$

We will demonstrate that 
\begin{equation}
\label{2.16000000}
\hbox{\vtop{\offinterlineskip\halign{
\hfil#\hfil\cr
{\rm l.i.m.}\cr
$\stackrel{}{{}_{N\to \infty}}$\cr
}} }R[\phi,\psi^{(k)}]_N=
\hbox{\vtop{\offinterlineskip\halign{
\hfil#\hfil\cr
{\rm l.i.m.}\cr
$\stackrel{}{{}_{N\to \infty}}$\cr
}} }S[\phi,\psi^{(k)}]_N\ \ \ \hbox{w. p. 1}.
\end{equation}

It is easy to see that  
\begin{equation}
\label{2.17000000}
~~~~~~~~ R[\phi,\psi^{(k)}]_N=
U[\phi,\psi^{(k)}]_N+
V[\phi,\psi^{(k)}]_N+
S[\phi,\psi^{(k)}]_N\ \ \ \hbox{w. p. 1},
\end{equation}
where
$$
U[\phi,\psi^{(k)}]_N=
\sum^{N-1}_{j_1=0}
\psi_1(\tau_{j_1})\Delta w_{\tau_{j_1}}^{(1)}
G[\psi_2^{(k-1)}]_{j_1,0}
\sum_{l=0}^{j_k-1}
I[\Delta\phi]_{\tau_{l+1},\tau_l},
$$
$$
V[\phi,\psi^{(k)}]_N=
\sum^{N-1}_{j_1=0}
I[\Delta\psi_1]_{\tau_{j_1+1},\tau_{j_1}}
G[\psi_2^{(k-1)}]_{j_1,0}
\sum_{l=0}^{j_k-1}
\phi_{\tau_l}\Delta w_{\tau_l}^{(k+1)},
$$
$$
I[\Delta\psi_1]_{\tau_{j_1+1},\tau_{j_1}}
=\int\limits_{\tau_{j_1}}^{\tau_{j_1+1}}
(\psi_1(\tau_{j_1})-\psi_1(\tau))dw_{\tau}^{(1)},
$$
$$
I[\Delta\phi]_{\tau_{l+1},\tau_l}=\int\limits_{\tau_l}^{\tau_{l+1}}
(\phi_{\tau}-\phi_{\tau_l})dw_{\tau}^{(k+1)}.
$$

Using 
the Minkowski inequality, standard 
estimates for second moments of stochastic integrals,
the condition that the process $\phi_\tau$ belongs
to the class ${\rm S}_2([t,T])$
as well as   
continuity (which means uniform continuity)
of the function $\psi_1(\tau)$,
we obtain that
$$
\hbox{\vtop{\offinterlineskip\halign{
\hfil#\hfil\cr
{\rm l.i.m.}\cr
$\stackrel{}{{}_{N\to \infty}}$\cr
}} }V[\phi,\psi^{(k)}]_N=
\hbox{\vtop{\offinterlineskip\halign{
\hfil#\hfil\cr
{\rm l.i.m.}\cr
$\stackrel{}{{}_{N\to \infty}}$\cr
}} }U[\phi,\psi^{(k)}]_N=0\ \ \ \hbox{w. p. 1}.
$$

Then, considering (\ref{2.17000000}), we obtain
(\ref{2.16000000}). From (\ref{2.16000000}) and (\ref{2.15000000}) 
it follows that 
the left equality in (\ref{2.4000000}) is fulfilled.

Note that the limit  
(\ref{li})
exists because it is equal to the stochastic
integral $J[\phi,\psi^{(k)}]_{T,t}$,
which exists under the conditions of Theorem 3.1.
So, the chain of equalities (\ref{2.4000000}) is proved. 
Theorem 3.1 is proved.

\section{Corollaries and Generalizations of Theorem 3.1}

Assume that $D_k=\{(t_1,\ldots,t_k): t\le t_1<\ldots<t_k\le T\}$
and the following conditions are fulfilled:

AI.\ $\xi_{\tau}\in{\rm S}_2([t,T]).$

AII.\ $\Phi(t_1,\ldots,{t_{k-1}})$ is a continuous nonrandom function 
in the closed
domain $D_{k-1}$ (recall that 
we use the same symbol $D_{k-1}$ to denote the open and closed 
domains corresponding to the domain $D_{k-1}$).

Let us define the following stochastic integrals
$$
\hat J[\xi,\Phi]_{T,t}^{(k)}=\int\limits_t^T \xi_{t_k}d{\bf w}_{t_k}^{(i_k)}
\ldots 
\int\limits_{t_{3}}^T d{\bf w}_{t_{2}}^{(i_{2})}
\int\limits_{t_{2}}^{T}\Phi(t_1,t_2,\ldots,t_{k-1})
d{\bf w}_{t_1}^{(i_1)}
\stackrel{\rm def}{=}
$$
$$
\stackrel{\rm def}{=}\hbox{\vtop{\offinterlineskip\halign{
\hfil#\hfil\cr
{\rm l.i.m.}\cr
$\stackrel{}{{}_{N\to \infty}}$\cr
}} }\sum_{l=0}^{N-1}\xi_{\tau_l}\Delta{\bf w}_{\tau_{l}}^{(i_k)}
\int\limits_{\tau_{l+1}}^T
d{\bf w}_{t_{k-1}}^{(i_{k-1})}
\ldots\int\limits_{t_{3}}^T
d{\bf w}_{t_{2}}^{(i_{2})}
\int\limits_{t_{2}}^T
\Phi(t_1,t_2,\ldots,t_{k-1})d{\bf w}_{t_1}^{(i_1)}
$$

\vspace{1mm}
\noindent
for $k\ge 3$ and
$$
\hat J[\xi,\Phi]_{T,t}^{(2)}=\int\limits_t^T \xi_{t_2}d{\bf w}_{t_2}^{(i_2)}
\int\limits_{t_{2}}^{T}\Phi(t_1)d{\bf w}_{t_1}^{(i_1)}
\stackrel{\rm def}{=}
$$
$$
\stackrel{\rm def}{=}
\hbox{\vtop{\offinterlineskip\halign{
\hfil#\hfil\cr
{\rm l.i.m.}\cr
$\stackrel{}{{}_{N\to \infty}}$\cr
}} }\sum_{l=0}^{N-1}\xi_{\tau_l}\Delta{\bf w}_{\tau_{l}}^{(i_2)}
\int\limits_{\tau_{l+1}}^T
\Phi(t_1)d{\bf w}_{t_1}^{(i_1)}
$$

\noindent
for $k=2$. Here 
${\bf w}_{\tau}^{(i)}$ $(i=1,\ldots,m)$ 
are 
${\rm F}_{\tau}$-measurable for all $\tau\in[0,T]$ 
independent standard Wiener processes,
${\bf w}_\tau^{(0)}=\tau,$
$0\le t<T,$
$i_1,\ldots,i_k=0, 1,\ldots,m$.

Let us denote
\begin{equation}
\label{2.18000000}
~~~~~~~~ J[\xi,\Phi]_{T,t}^{(k)}=
\int\limits_t^T\ldots\int\limits_t^{t_{k-1}}
\Phi(t_1,\ldots,t_{k-1})\xi_{t_k}
d{\bf w}_{t_k}^{(i_k)}
\ldots d{\bf w}_{t_1}^{(i_1)},\ \ \ k\ge 2,
\end{equation}

\vspace{2mm}
\noindent
where the right-hand side of (\ref{2.18000000}) is the
iterated It\^{o} stochastic 
integral.

Let us introduce the following iterated stochastic integrals

\vspace{-3mm}
$$
\tilde J[\Phi]_{T,t}^{(k-1)}
=\int\limits_t^T d{\bf w}_{t_{k-1}}^{(i_{k-1})}\ldots 
\int\limits_{t_{3}}^T d{\bf w}_{t_{2}}^{(i_{2})}
\int\limits_{t_{2}}^{T}\Phi(t_1,t_2,\ldots,t_{k-1})
d{\bf w}_{t_1}^{(i_1)}
\stackrel{\rm def}{=}
$$

\vspace{-2mm}
$$
\stackrel{\rm def}{=}\hbox{\vtop{\offinterlineskip\halign{
\hfil#\hfil\cr
{\rm l.i.m.}\cr
$\stackrel{}{{}_{N\to \infty}}$\cr
}} }\sum_{l=0}^{N-1}\Delta{\bf w}_{\tau_{l}}^{(i_{k-1})}
\int\limits_{\tau_{l+1}}^T
d{\bf w}_{t_{k-2}}^{(i_{k-2})}
\ldots\int\limits_{t_{3}}^T
d{\bf w}_{t_{2}}^{(i_{2})}
\int\limits_{t_{2}}^T
\Phi(t_1,t_2,\ldots,t_{k-1})d{\bf w}_{t_1}^{(i_1)},
$$

$$
J'[\Phi]_{T,t}^{(k-1)}=
\int\limits_t^T\ldots\int\limits_t^{t_{k-2}}
\Phi(t_1,\ldots,t_{k-1})
d{\bf w}_{t_{k-1}}^{(i_{k-1})}
\ldots d{\bf w}_{t_1}^{(i_1)},\ \ \ k\ge 2.
$$

\vspace{3mm}

Similarly to the proof of 
Theorem 3.1 it is easy to demonstrate that under the condition  
AII the stochastic  
integral 
$\tilde J[\Phi]_{T,t}^{(k-1)}$ exists and 

\vspace{-1mm}
\begin{equation}
\label{432}
J'[\Phi]_{T,t}^{(k-1)}=\tilde J[\Phi]_{T,t}^{(k-1)}\ \ \ \hbox{w. p. 1.}
\end{equation}

\vspace{3mm}

Moreover, using (\ref{432}) 
the following generalization
of Theorem 3.1 can be proved
similarly to the proof of Theorem 3.1.

\vspace{2mm}

{\bf Theorem 3.2} \cite{vini} (1997)
(also see \cite{1}-\cite{12aa}, \cite{old-art-2}, 
\cite{arxiv-25}).
{\it Suppose that the conditions  
{\rm AI, AII} of this section are fulfilled.
Then, the stochastic integral
$\hat J[\xi,\Phi]_{T,t}^{(k)}$ exists and for $k\ge 2$
$$
J[\xi,\Phi]_{T,t}^{(k)}=\hat J[\xi,\Phi]_{T,t}^{(k)}\ \ \ 
\hbox{w. p. {\rm 1.}}
$$
}

Let us consider the following stochastic integrals
$$
I=\int\limits_t^{T}d{\bf w}_{t_2}^{(i_2)}\int\limits_{t_2}^T
\Phi_1(t_1,t_2)d{\bf w}_{t_1}^{(i_1)},\ \ \ 
J=\int\limits_t^{T}\int\limits_t^{t_2}
\Phi_2(t_1,t_2)d{\bf w}_{t_1}^{(i_1)}d{\bf w}_{t_2}^{(i_2)},
$$

\noindent
where $i_1,i_2=1,\ldots,m.$ If we consider 
$$\int\limits_{t_2}^T
\Phi_1(t_1,t_2)d{\bf w}_{t_1}^{(i_1)}
$$ 

\noindent
as the integrand of $I$ 
and 
$$
\int\limits_t^{t_2}
\Phi_2(t_1,t_2)d{\bf w}_{t_1}^{(i_1)}
$$ 

\noindent
as the integrand of $J$,
then, due to independence of these integrands we may mistakenly 
think that ${\sf M}\{IJ\}=0.$
But it is not the fact. Actually, using the 
integration order replacement technique in the stochastic
integral $I$, we have w.~p.~1 
$$
I=\int\limits_t^{T}\int\limits_t^{t_1}
\Phi_1(t_1,t_2)d{\bf w}_{t_2}^{(i_2)}d{\bf w}_{t_1}^{(i_1)}
=\int\limits_t^{T}\int\limits_t^{t_2}
\Phi_1(t_2,t_1)d{\bf w}_{t_1}^{(i_2)}d{\bf w}_{t_2}^{(i_1)}.
$$

\vspace{1mm}

So, using the standard properties of the It\^{o} stochastic 
integral \cite{Gih1}, we get
$$
{\sf M}\{IJ\}={\bf 1}_{\{i_1=i_2\}}\int\limits_t^{T}\int\limits_t^{t_2}
\Phi_1(t_2,t_1)\Phi_2(t_1,t_2)dt_1dt_2,
$$

\noindent
where ${\bf 1}_A$ is the indicator of the set $A$.

Let us consider the following statement.

\vspace{2mm}

{\bf Theorem 3.3} \cite{vini} (1997)
(also see \cite{1}-\cite{12aa}, \cite{old-art-2}, 
\cite{arxiv-25}).
{\it Let the 
conditions of Theorem {\rm 3.1} are fulfilled and $h(\tau)$  
is a continuous nonrandom function at the interval $[t,T]$.
Then
\begin{equation}
\label{2.19000000}
~~~~~ \int\limits_{t}^T \phi_\tau dw_\tau^{(k+1)}h(\tau) 
\hat I[\psi^{(k)}]_{T,\tau}=
\int\limits_{t}^T 
\phi_\tau h(\tau) dw_\tau^{(k+1)}\hat I[\psi^{(k)}]_{T,\tau}\ \ \
\hbox{{\rm w. p. 1}},
\end{equation}

\noindent
where stochastic integrals on the left-hand side of {\rm (\ref{2.19000000})} 
as well as on the right-hand side of {\rm (\ref{2.19000000})} 
exist.}

{\bf Proof.}\
According to Theorem 3.1, the iterated stochastic
integral on the right-hand side of (\ref{2.19000000}) exists.
In addition
$$
\int\limits_{t}^T 
\phi_\tau h(\tau) dw_\tau^{(k+1)} \hat I[\psi^{(k)}]_{T,\tau}=
\int\limits_{t}^T 
\phi_\tau dw_\tau^{(k+1)}h(\tau)\hat I[\psi^{(k)}]_{T,\tau} -
$$
$$
- \hbox{\vtop{\offinterlineskip\halign{
\hfil#\hfil\cr
{\rm l.i.m.}\cr
$\stackrel{}{{}_{N\to \infty}}$\cr
}} }\sum_{l=0}^{N-1} \phi_{\tau_l} \Delta h(\tau_{l})
\Delta w_{\tau_l}^{(k+1)} \hat I[\psi^{(k)}]_{T,\tau_{l+1}}\ \ \ 
\hbox{w. p. 1},
$$

\vspace{2mm}
\noindent
where $\Delta h(\tau_{l})=h(\tau_{l+1})-h(\tau_{l}).$

Using the arguments which resulted to the right equality
in (\ref{2.4000000}), we obtain
$$
\hbox{\vtop{\offinterlineskip\halign{
\hfil#\hfil\cr
{\rm l.i.m.}\cr
$\stackrel{}{{}_{N\to \infty}}$\cr
}} }\sum_{l=0}^{N-1}\phi_{\tau_l}\Delta h(\tau_{l})
\Delta w_{\tau_l}^{(k+1)} \hat I[\psi^{(k)}]_{T,\tau_{l+1}}=
$$
\begin{equation}
\label{2.20000001}
=\hbox{\vtop{\offinterlineskip\halign{
\hfil#\hfil\cr
{\rm l.i.m.}\cr
$\stackrel{}{{}_{N\to \infty}}$\cr
}} }G[\psi^{(k)}]_{N,0}
\sum_{l=0}^{j_k-1}\phi_{\tau_l}\Delta h(\tau_{l})
\Delta w_{\tau_l}^{(k+1)}\ \ \ \hbox{w. p. 1}.
\end{equation}

\vspace{2mm}

Using 
the Minkowski inequality, standard 
estimates for second moments of stochastic integrals as well as
continuity of the function $h(\tau)$, we obtain that
the second moment of the prelimit expression on  
the right-hand side of (\ref{2.20000001}) tends to zero 
when $N\to\infty.$ 
Theorem 3.3 is proved.

Let us consider one corollary of Theorem 3.1.

\vspace{2mm}

{\bf Theorem 3.4}\ \cite{vini} (1997)
(also see \cite{1}-\cite{12aa}, \cite{old-art-2}, 
\cite{arxiv-25}).
{\it Under
the conditions of Theorem {\rm 3.3}
the following equality 
$$
\int\limits_{t}^T h(t_1)\int\limits_{t}^{t_1}\phi_\tau dw_\tau^{(k+2)}
dw_{t_1}^{(k+1)} \hat I[\psi^{(k)}]_{T,t_1}=
$$
\begin{equation}
\label{2.21000000}
=\int\limits_{t}^T \phi_\tau dw_\tau^{(k+2)}\int\limits_{\tau}^T
h(t_1)dw_{t_1}^{(k+1)}\hat I[\psi^{(k)}]_{T,t_1}\ \ \ 
\hbox{w. p. {\rm 1}}
\end{equation}

\noindent
is fulfilled.
Moreover, the stochastic integrals in {\rm (\ref{2.21000000})} exist.}

{\bf Proof.}\ Using Theorem 3.1 two times,
we obtain
$$
\int\limits_{t}^T \phi_\tau dw_\tau^{(k+2)}\int\limits_{\tau}^T
h(t_1)dw_{t_1}^{(k+1)} \hat I[\psi^{(k)}]_{T,t_1}=
$$
$$
=\int\limits_{t}^T\psi_1(t_1)\ldots
\int\limits_{t}^{t_{k-1}}\psi_k(t_k)\int\limits_{t}^{t_k}
\rho_{\tau}dw_{\tau}^{(k+1)} dw_{t_k}^{(k)}\ldots dw_{t_1}^{(1)}=
$$
$$
=\int\limits_{t}^T \rho_\tau dw_\tau^{(k+1)}\int\limits_{\tau}^T\psi_k(t_k)
dw_{t_k}^{(k)}\ldots \int\limits_{t_{2}}^T\psi_1(t_1)dw_{t_1}^{(1)}\ \ \
\hbox{w. p. 1},
$$
where
$$
\rho_{\tau}\stackrel{\rm def}{=}h(\tau)\int\limits_{t}^{\tau}
\phi_s dw_s^{(k+2)}.
$$ 

\vspace{1mm}

Theorem 3.4 is proved.

\section{Examples of Integration Order Replacement Technique
for the Concrete Iterated It\^{o} Stochastic Integrals}

As we mentioned above, the formulas from this section
could be obtained using the It\^{o} formula. However,
the method based on Theorem 3.1 is more simple
and familiar, since it deals with usual 
rules of the integration order replacement for Riemann integrals.

Using the integration order replacement technique
for iterated It\^{o} stochastic integrals (Theorem 3.1), we obtain 
the following equalities which are fulfilled
w.~p.~1
\begin{equation}
\label{febr3}
\int\limits_t^T\int\limits_t^{t_2}dw_{t_1}dt_2=
\int\limits_t^T(T-t_1)dw_{t_1},
\end{equation}
$$
\int\limits_t^T {\rm cos}(t_2-T)\int\limits_t^{t_2}dw_{t_1}dt_2
=\int\limits_t^T {\rm sin}(T-t_1)dw_{t_1},
$$
$$
\int\limits_t^T {\rm sin}(t_2-T)\int\limits_t^{t_2}dw_{t_1}dt_2
=\int\limits_t^T 
\left({\rm cos}(T-t_1)-1\right)dw_{t_1},
$$
$$
\int\limits_t^T e^{\alpha(t_2-T)}\int\limits_t^{t_2}dw_{t_1}dt_2=
\frac{1}{\alpha}\int\limits_t^T\left(1-e^{\alpha(t_1-T)}\right)dw_{t_1},\ \ \
\alpha\ne 0,
$$
$$
\int\limits_t^T(t_2-T)^{\alpha}\int\limits_t^{t_2}dw_{t_1}dt_2=
-\frac{1}{\alpha+1}\int\limits_t^T(t_1-T)^{\alpha+1}dw_{t_1},\ \ \ 
\alpha\ne -1,
$$
$$
J_{(100)T,t}=\frac{1}{2}\int\limits_t^T(T-t_1)^2 dw_{t_1},\ \ \
J_{(010)T,t}=\int\limits_t^T(t_1-t)(T-t_1)dw_{t_1},
$$
\begin{equation}
\label{ex1}
J_{(110)T,t}=\int\limits_t^T(T-t_2)\int\limits_t^{t_2}dw_{t_1}dw_{t_2},\
\end{equation}
$$
J_{(101)T,t}=\int\limits_t^T\int\limits_t^{t_2}(t_2-t_1)dw_{t_1}dw_{t_2},\ \ \
J_{(1011)T,t}=\int\limits_t^T\int\limits_t^{t_3}\int\limits_t^{t_2}(t_2-t_1)
dw_{t_1}dw_{t_2}dw_{t_3},\
$$
$$
J_{(1101)T,t}=\int\limits_t^T\int\limits_t^{t_3}(t_3-t_2)\int\limits_t^{t_2}
dw_{t_1}dw_{t_2}dw_{t_3},
$$
$$
J_{(1110)T,t}=\int\limits_t^T(T-t_3)\int\limits_t^{t_3}\int\limits_t^{t_2}
dw_{t_1}dw_{t_2}dw_{t_3},\ \ \
J_{(1100)T,t}=\frac{1}{2}\int\limits_t^T(T-t_2)^2\int\limits_t^{t_2}
dw_{t_1}dw_{t_2},
$$
\begin{equation}
\label{ex2}
J_{(1010)T,t}=\int\limits_t^T(T-t_2)\int\limits_t^{t_2}(t_2-t_1)
dw_{t_1}dw_{t_2},
\end{equation}
$$
J_{(1001)T,t}=\frac{1}{2}\int\limits_t^T\int\limits_t^{t_2}(t_2-t_1)^2
dw_{t_1}dw_{t_2},\ \ \
J_{(0110)T,t}=\int\limits_t^T(T-t_2)\int\limits_t^{t_2}(t_1-t)
dw_{t_1}dw_{t_2},
$$
$$
J_{(0101)T,t}=\int\limits_t^T\int\limits_t^{t_2}(t_2-t_1)(t_1-t)
dw_{t_1}dw_{t_2},
$$
$$
J_{(0010)T,t}=\frac{1}{2}\int\limits_t^T(T-t_1)(t_1-t)^2
dw_{t_1},\ \ \
J_{(0100)T,t}=\frac{1}{2}\int\limits_t^T(T-t_1)^2(t_1-t)
dw_{t_1},
$$
$$
J_{(1000)T,t}=\frac{1}{3!}\int\limits_t^T(T-t_1)^3
dw_{t_1},\
$$
$$
J_{(1\underbrace{{}_{0\ldots 0}}_{k-1})T,t}=
\frac{1}{(k-1)!}\int\limits_t^T(T-t_1)^{k-1}dw_{t_1},
$$
$$
J_{(11\underbrace{{}_{0\ldots 0}}_{k-2})T,t}=
\frac{1}{(k-2)!}\int\limits_t^T(T-t_2)^{k-2}
\int\limits_t^{t_2}dw_{t_1}dw_{t_2},
$$
$$
J_{(\underbrace{{}_{1\ldots 1}}_{k-1}0)T,t}=
\int\limits_t^T(T-t_1)J_{(\underbrace{{}_{1\ldots 1}}_{k-2})t_1,t}dw_{t_1},
$$
$$
J_{(1\underbrace{{}_{0\ldots 0}}_{k-2}1)T,t}=
\frac{1}{(k-2)!}\int\limits_t^T\int\limits_t^{t_2}(t_2-t_1)^{k-2}
dw_{t_1}dw_{t_2},
$$
$$
J_{(10\underbrace{{}_{1\ldots 1}}_{k-2})T,t}=
\int\limits_t^T\ldots \int\limits_t^{t_3}\int\limits_t^{t_2}(t_2-t_1)
dw_{t_1}dw_{t_2}\ldots dw_{t_{k-1}},
$$
$$
J_{(\underbrace{{}_{1\ldots 1}}_{k-2}01)T,t}=
\int\limits_t^T\int\limits_t^{t_{k-1}}(t_{k-1}-t_{k-2})
\int\limits_t^{t_{k-2}}\ldots
\int\limits_t^{t_2}dw_{t_1}\ldots dw_{t_{k-3}}dw_{t_{k-2}}dw_{t_{k-1}},
$$

\vspace{1mm}
$$
J_{(10)T,t}+J_{(01)T,t}=(T-t)J_{(1)T,t},\
$$

\vspace{1mm}
$$
J_{(110)T,t}+J_{(101)T,t}+J_{(011)T,t}=(T-t)J_{(11)T,t},\
$$

\vspace{1mm}
$$
J_{(001)T,t}+J_{(010)T,t}+J_{(100)T,t}=\frac{(T-t)^2}{2}
J_{(1)T,t},
$$

$$
J_{(1100)T,t}+J_{(1010)T,t}+J_{(1001)T,t}+
J_{(0110)T,t}+
$$

\vspace{-2mm}
$$
+J_{(0101)T,t}+J_{(0011)T,t}=\frac{(T-t)^2}{2}
J_{(11)T,t},
$$

\vspace{1mm}

$$
J_{(1000)T,t}+J_{(0100)T,t}+J_{(0010)T,t}+J_{(0001)T,t}=\frac{(T-t)^3}{3!}
J_{(1)T,t},
$$

\vspace{1mm}
$$
J_{(1110)T,t}+J_{(1101)T,t}+J_{(1011)T,t}+J_{(0111)T,t}=
(T-t)J_{(111)T,t},
$$

\vspace{1mm}
$$
\sum_{l=1}^k J_{(\underbrace{{}_{0\ldots 0}}_{l-1}1 
\underbrace{{}_{0\ldots 0}}_{k-l})T,t}=\frac{1}{(k-1)!}(T-t)^{k-1}J_{(1)T,t},
$$
$$
\sum_{l=1}^k J_{(\underbrace{{}_{1\ldots 1}}_{l-1}0 
\underbrace{{}_{1\ldots 1}}_{k-l})T,t}=(T-t)
J_{(\underbrace{{}_{1\ldots 1}}_{k-1} )T,t},
$$
$$
\sum_{{}_{\stackrel{l_1+\ldots+l_k=m}{l_i\in\{0,\ 1\},\ i=1,\ldots,k}}}
J_{(l_1\ldots l_k)T,t}=\frac{(T-t)^{k-m}}{(k-m)!}
J_{(\underbrace{{}_{1\ldots 1}}_{m})T,t},
$$

\vspace{2mm}
\noindent
where
$$
J_{(l_1\ldots l_k)T,t}=\int\limits_t^T\ldots \int\limits_t^{t_2}dw_{t_1}
^{(1)}\ldots dw_{t_k}^{(k)},
$$

\vspace{1mm}
\noindent
$l_i=1$ when $w_{t_i}^{(i)}=w_{t_i}$ and $l_i=0$ when
$w_{t_i}^{(i)}=t_i$ $(i=1,\ldots,k),$
$w_\tau$ is a standard Wiener process.

Let us consider two examples and show explicitly the technique on 
integration
order replacement for iterated It\^{o} stochastic integrals.

{\bf Example 3.1.} {\it Let us prove the equality {\rm (\ref{ex1}).} 
Using Theorems {\rm 3.1} and {\rm 3.3,} 
we obtain{\rm :}}
$$
J_{(110)T,t}\stackrel{\rm def}{=}
\int\limits_t^T\int\limits_t^{t_3}\int\limits_t^{t_2} dw_{t_1}dw_{t_2}dt_3
=\int\limits_t^T dw_{t_1}\int\limits_{t_1}^T dw_{t_2} \int\limits_{t_2}^T dt_3=
$$
$$
=
\int\limits_t^T dw_{t_1}\int\limits_{t_1}^T dw_{t_2}(T-t_2)
=\int\limits_t^T dw_{t_1}\int\limits_{t_1}^T (T-t_2) dw_{t_2}=
$$
\begin{equation}
\label{febr10}
=
\int\limits_t^T(T-t_2)\int\limits_t^{t_2}dw_{t_1}dw_{t_2}\ \ \ 
\hbox{w. p. 1.}
\end{equation}

{\bf Example 3.2.} {\it Let us prove the equality {\rm (\ref{ex2}).} 
Using Theorems {\rm 3.1} and {\rm 3.3,} 
we obtain}
$$
J_{(1010)T,t}\stackrel{\rm def}{=}
\int\limits_t^T\int\limits_t^{t_4}
\int\limits_t^{t_3}\int\limits_t^{t_2}dw_{t_1}dt_2dw_{t_3}dt_4=
\int\limits_t^T dw_{t_1}\int\limits_{t_1}^T dt_2 \int\limits_{t_2}^T dw_{t_3}
\int\limits_{t_3}^T dt_4
=
$$
$$
=\int\limits_t^T dw_{t_1}\int\limits_{t_1}^T 
dt_2 \int\limits_{t_2}^T dw_{t_3}(T-t_3)
=\int\limits_t^T dw_{t_1}
\int\limits_{t_1}^T dt_2 \int\limits_{t_2}^T(T-t_3)dw_{t_3}=
$$
$$
=
\int\limits_t^T(T-t_3)
\int\limits_t^{t_3}\int\limits_t^{t_2} dw_{t_1}dt_2dw_{t_3}
=\int\limits_t^T(T-t_3)\left(\int\limits_t^{t_3}
\int\limits_t^{t_2} dw_{t_1}dt_2\right)dw_{t_3}=
$$
$$
=
\int\limits_t^T(T-t_3)\left(
\int\limits_t^{t_3}dw_{t_1}\int\limits_{t_1}^{t_3}dt_2\right)dw_{t_3}
=
$$
$$
=\int\limits_t^T(T-t_3)\left(
\int\limits_t^{t_3}dw_{t_1}(t_3-t_1)\right)dw_{t_3}=
$$
$$
=
\int\limits_t^T(T-t_3)\left(
\int\limits_t^{t_3}(t_3-t_1)dw_{t_1}\right)dw_{t_3}=
$$
$$
=
\int\limits_t^T(T-t_2)\int\limits_t^{t_2}(t_2-t_1)
dw_{t_1}dw_{t_2}\ \ \ \hbox{w. p. 1.}
$$

\section{Integration Order Replacement Technique for Iterated
Sto\-chas\-tic Integrals with Respect to Martingale}

In this section, we will generalize the theorems on
integration order replacement for iterated It\^{o} stochastic 
integrals to the class of iterated stochastic 
integrals with respect to martingale.

Let
$(\Omega,{\rm F},{\sf P})$ be a complete probability space 
and 
let $\{{\rm F}_t, t\in[0, T]\}$ be a nondecreasing 
family of $\sigma$-algebras defined on the probability space
$(\Omega,{\rm F},{\sf P}).$
Suppose that $M_t,$ $t\in[0,T]$ is an ${\rm F}_t$-measurable 
martingale for all 
$t\in[0,T]$, which 
satisfies
the condition
${\sf M}\left\{\left|M_t\right|\right\}<\infty$. Moreover,
for all $t\in[0,T]$ there exists
an ${\rm F}_t$- measurable and 
nonnegative w.~p.~1 stochastic process 
$\rho_t,$ $t\in[0,T]$ such that 
$$
{\sf M}\left\{\left(M_s-M_t\right)^2\left|\right.{\rm F}_t\right\}=
{\sf M}\left\{\int\limits_t^s\rho_{\tau}d\tau \biggl|\biggr.
{\rm F}_t\right\}\ \ \
\hbox{w.~p.~1},
$$
where $0\le t<s\le T.$

Let us consider the class $H_2(\rho,[0,T])$
of stochastic processes
$\varphi_t,$ $t\in[0,T]$, which are
${\rm F}_t$-measurable for all 
$t\in[0,T]$ and 
satisfy
the condition 
$$
{\sf M}\left\{\int\limits_0^T\varphi_t^2 \rho_t dt\right\}<\infty.
$$

For any partition 
$\tau_j^{(N)},$ $j=0, 1, \ldots, N$ of
the interval $[0,T]$ such that
\begin{equation}
\label{w11ggg}
0=\tau_0^{(N)}<\tau_1^{(N)}<\ldots <\tau_N^{(N)}=T,\ \ \ \
\max\limits_{0\le j\le N-1}\left|\tau_{j+1}^{(N)}-\tau_j^{(N)}\right|\to 0\ \
\hbox{if}\ \ N\to \infty
\end{equation}
we will define the sequence of step functions 

\vspace{-2mm}
$$
\varphi^{(N)}(t,\omega)=
\varphi_j\left(\omega\right)\ \ \ \hbox{w. p. 1}\ \ \ 
\hbox{for}\ \ \ t\in\left[\tau_j^{(N)},\tau_{j+1}^{(N)}\right),
$$

\vspace{1mm}
\noindent
where $\varphi^{(N)}(t,\omega)\in H_2(\rho,[0,T]),$ $j=0, 1,\ldots,N-1,$ $N=1, 2,\ldots$

Let us define the stochastic integral with respect to martingale
for
$\varphi(t,\omega)\in H_2(\rho,[0,T])$ 
as the 
following mean-square limit \cite{Gih1}
$$
\hbox{\vtop{\offinterlineskip\halign{
\hfil#\hfil\cr
{\rm l.i.m.}\cr
$\stackrel{}{{}_{N\to \infty}}$\cr
}} }\sum_{j=0}^{N-1}\varphi^{(N)}\left(\tau_j^{(N)},\omega\right)
\left(M\left(\tau_{j+1}^{(N)},\omega\right)-
M\left(\tau_j^{(N)},\omega\right)\right)
\stackrel{\rm def}{=}\int\limits_0^T\varphi_\tau dM_\tau,
$$
where $\varphi^{(N)}(t,\omega)$ is any step function
from the class $H_2(\rho,[0,T])$,
which converges
to the function $\varphi(t,\omega)$
in the following sense
$$
\hbox{\vtop{\offinterlineskip\halign{
\hfil#\hfil\cr
{\rm lim}\cr
$\stackrel{}{{}_{N\to \infty}}$\cr
}} }\int\limits_0^T{\sf M}\left\{\left|
\varphi^{(N)}(t,\omega)-\varphi(t,\omega)\right|^2\right\}\rho_tdt=0.
$$

It is well known  \cite{Gih1} that the stochastic integral
$$
\int\limits_0^T\varphi_{\tau} dM_{\tau}
$$
exists and it does not depend on the selection 
of sequence 
$\varphi^{(N)}(t,\omega)$.

Let $\tilde H_2(\rho,[0,T])$ be the class of 
stochastic processes  $\varphi_{\tau},$ $\tau\in[0,T],$
which are
mean-square
continuous for all $\tau\in[0,T]$ and belong to the 
class $H_2(\rho,[0,T])$.

Let us consider the following iterated stochastic integrals
\begin{equation}
\label{32.001}
~~~~S[\phi,\psi^{(k)}]_{T,t}=\int\limits_{t}^{T}\psi_1(t_1)\ldots
\int\limits_t^{t_{k-1}}\psi_k(t_k)\int\limits_t^{t_k}
\phi_{\tau}dM^{(k+1)}_{\tau}dM^{(k)}_{t_k}
\ldots dM^{(1)}_{t_1},
\end{equation}
\begin{equation}
\label{32.002}
S[\psi^{(k)}]_{T,t}=\int\limits_{t}^{T}\psi_1(t_1)\ldots
\int\limits_t^{t_{k-1}}\psi_k(t_k)dM^{(k)}
_{t_k}
\ldots dM^{(1)}_{t_1}.
\end{equation}

\noindent
Here $\phi_\tau\in \tilde H_2(\rho,[t,T])$ and
$\psi_1(\tau),\ldots,\psi_k(\tau)$ are continuous nonrandom 
functions at the interval
$[t,T]$,
$M_\tau^{(l)}=M_\tau$ or $M_\tau^{(l)}=\tau$
if $\tau\in[t,T],$
$l=1,\ldots,k+1,$\ $M_{\tau}$ is the martingale defined above.

Let us define the iterated stochastic 
integral 
$\hat S[\psi^{(k)}]_{T,s},$ $0\le t\le s\le T,$ $k\ge 1$
with respect to martingale
$$
\hat S[\psi^{(k)}]_{T,s}=\int\limits_s^T\psi_k(t_k)dM_{t_k}^{(k)}
\ldots \int\limits_{t_{2}}^T \psi_1(t_1)dM_{t_1}^{(1)}
$$

\noindent
by the following recurrence 
relation
\begin{equation}
\label{32.003}
\hat S[\psi^{(k)}]_{T,t}
\stackrel{\rm def}{=}
\hbox{\vtop{\offinterlineskip\halign{
\hfil#\hfil\cr
{\rm l.i.m.}\cr
$\stackrel{}{{}_{N\to \infty}}$\cr
}} }
\sum^{N-1}_{l=0} \psi_k(\tau_{l})\Delta M_{\tau_l}^{(k)} 
\hat S[\psi^{(k-1)}]_{T,\tau_{l+1}},
\end{equation}

\noindent
where $k\ge 1,$ 
$\hat S[\psi^{(0)}]_{T,s}\stackrel{\rm def}{=}1,$
$[s,T]\subseteq[t,T],$\ here and further $\Delta M_{\tau_l}^{(i)}=
M_{\tau_{l+1}}^{(i)}-M_{\tau_l}^{(i)},$
$i=1,\ldots,k+1,$ $l=0, 1,\ldots,N-1,$
$\{\tau_l\}_{l=0}^N$ is the partition of the interval $[t, T]$,
which satisfies the condition similar to (\ref{w11ggg}),
another notations are the same as in
(\ref{32.001}), (\ref{32.002}).

Further,
let us define the iterated stochastic 
integral 
$\hat S[\phi,\psi^{(k)}]_{T,t},$ $k\ge 1$ of the form
$$
\hat S[\phi,\psi^{(k)}]_{T,t}=\int\limits_{t}^T \phi_s dM_s^{(k+1)}
\hat S[\psi^{(k)}]_{T,s}
$$
by the equality
$$
\hat S[\phi,\psi^{(k)}]_{T,t}
\stackrel{\rm def}{=}
\hbox{\vtop{\offinterlineskip\halign{
\hfil#\hfil\cr
{\rm l.i.m.}\cr
$\stackrel{}{{}_{N\to \infty}}$\cr
}} }
\sum^{N-1}_{l=0} \phi_{\tau_{l}}\Delta M_{\tau_l}^{(k+1)} 
\hat S[\psi^{(k)}]_{T,\tau_{l+1}},
$$

\vspace{2mm}
\noindent
where the sense of notations included in 
(\ref{32.001})--(\ref{32.003}) is saved.

Let us formulate the theorem on integration order replacement 
for the iterated
stochastic integrals with respect to martingale,
which is the generalization of Theorem 3.1.

\vspace{1mm}

{\bf Theorem 3.5} \cite{old-preprint} (1999)
(also see \cite{1}-\cite{12aa}, \cite{arxiv-25}).
{\it Let 
$\phi_\tau\in\tilde H_2(\rho,[t,T]),$ every $\psi_l(\tau)$
$(l=1,\ldots,k)$ is a continuous nonrandom function at the interval 
$[t,T],$ and
$|\rho_{\tau}|\le K<\infty $ w. p. {\rm 1} for all
$\tau\in[t,T].$
Then, the stochastic integral 
$\hat S[\phi,\psi^{(k)}]_{T,t}$ exists and
$$
S[\phi,\psi^{(k)}]_{T,t}=\hat S[\phi,\psi^{(k)}]_{T,t}\ \ \ 
\hbox{w. p. {\rm 1}}.
$$ 
}

\newpage
\noindent
\par
The proof of Theorem 3.5 is similar to the proof of 
Theorem 3.1.

{\bf Remark 3.2.}\ {\it Let us note that we can propose 
another variant of the conditions in Theorem {\rm 3.5}. For example,
if we not require the boundedness of the process  
$\rho_{\tau}$, then it is necessary to require
the fulfillment of the following additional 
conditions{\rm :}

{\rm 1}. ${\sf M}\{|\rho_{\tau}|\}<\infty$ for all $\tau\in[t,T].$

{\rm 2}. The process $\rho_{\tau}$ is independent with 
the processes  $\phi_{\tau}$
and $M_{\tau}.$
}

{\bf Remark 3.3.}\ {\it Note that it is well known
the construction of stochastic integral
with respect to the Wiener process with integrable process, which is
not an ${\rm F}_{\tau}$-measurable 
stochastic  
process  --- the so-called 
Stratonovich stochastic integral {\rm \cite{str}}.

The stochastic integral $\hat S[\phi,\psi^{(k)}]_{T,t}$
is also the stochastic  
integral with 
integrable process, which is 
not an ${\rm F}_{\tau}$-measurable 
stochastic  
process. 
However, under the conditions of Theorem {\rm 3.5}
$$
S[\phi,\psi^{(k)}]_{T,t}=\hat S[\phi,\psi^{(k)}]_{T,t}\ \ \ 
\hbox{w. p. {\rm 1,}}
$$
where $S[\phi,\psi^{(k)}]_{T,t}$ is a usual
iterated stochastic integral with respect to martingale.
If, for example, $M_{\tau}, \tau\in[t,T]$
is the Wiener process, then the question on connection between stochastic 
integral $\hat S[\phi,\psi^{(k)}]_{T,t}$
and Stratonovich stochastic integral
is solving as a standard question on connection
between Stratonovich and It\^{o} stochastic integrals {\rm \cite{str}}.
}

Let us consider several statements, which are the generalizations of theorems
formulated in the previous sections.

Assume that $D_k=\{(t_1,\ldots,t_k):\ t\le t_1<\ldots<t_k\le T\}$
and the following conditions are fulfilled:

BI.\ $\xi_{\tau}\in \tilde H_2(\rho,[t,T]).$

BII.\ $\Phi(t_1,\ldots,{t_{k-1}})$ is a continuous nonrandom function
in the closed
domain $D_{k-1}$ (recall that 
we use the same symbol $D_{k-1}$ to denote the open and closed 
domains corresponding to the domain $D_{k-1}$).

Let us define the following stochastic integrals
with respect to martingale
$$
\hat S[\xi,\Phi]_{T,t}^{(k)}=\int\limits_t^T \xi_{t_k}dM_{t_k}^{(k)}\ldots 
\int\limits_{t_{3}}^T dM_{t_{2}}^{(2)}
\int\limits_{t_{2}}^{T}\Phi(t_1,t_2,\ldots,t_{k-1})
dM_{t_1}^{(1)}
\stackrel{\rm def}{=}
$$
$$
\stackrel{\rm def}{=}\hbox{\vtop{\offinterlineskip\halign{
\hfil#\hfil\cr
{\rm l.i.m.}\cr
$\stackrel{}{{}_{N\to \infty}}$\cr
}} }\sum_{l=0}^{N-1}\xi_{\tau_l}\Delta M_{\tau_{l}}^{(k)}
\int\limits_{\tau_{l+1}}^T
dM_{t_{k-1}}^{(k-1)}
\ldots\int\limits_{t_{3}}^T
dM_{t_{2}}^{(2)}
\int\limits_{t_{2}}^T
\Phi(t_1,t_2,\ldots,t_{k-1})dM_{t_1}^{(1)}
$$

\vspace{3mm}
\noindent
for $k\ge 3$ and
$$
\hat S[\xi,\Phi]_{T,t}^{(2)}=\int\limits_t^T \xi_{t_2}dM_{t_2}^{(2)}
\int\limits_{t_{2}}^{T}\Phi(t_1)dM_{t_1}^{(1)}
\stackrel{\rm def}{=}
$$
$$
\stackrel{\rm def}{=}
\hbox{\vtop{\offinterlineskip\halign{
\hfil#\hfil\cr
{\rm l.i.m.}\cr
$\stackrel{}{{}_{N\to \infty}}$\cr
}} }\sum_{l=0}^{N-1}\xi_{\tau_l}\Delta M_{\tau_{l}}^{(2)}
\int\limits_{\tau_{l+1}}^T
\Phi(t_1)dM_{t_1}^{(1)}
$$

\vspace{1mm}
\noindent
for $k=2$, where the sense of notations included 
in (\ref{32.001})--(\ref{32.003}) is saved.
Moreover,
the stochastic process  $\xi_{\tau},$ $\tau\in[t,T]$
belongs to the class $\tilde H_2(\rho,[t,T]).$

In addition, let
\begin{equation}
\label{32.006}
~~~~~~~~S[\xi,\Phi]_{T,t}^{(k)}=
\int\limits_t^T\ldots\int\limits_t^{t_{k-1}}
\Phi(t_1,\ldots,t_{k-1})\xi_{t_k}
dM_{t_k}^{(k)}
\ldots dM_{t_1}^{(1)},\ \ \ k\ge 2,
\end{equation}
where the right-hand side of (\ref{32.006}) is the iterated stochastic 
integral with respect to martingale.

Let us introduce the following iterated stochastic integrals
with respect to martingale
$$
\tilde S[\Phi]_{T,t}^{(k-1)}
=\int\limits_t^T d{M}_{t_{k-1}}^{(k-1)}
\ldots 
\int\limits_{t_{3}}^T d{M}_{t_{2}}^{(2)}
\int\limits_{t_{2}}^{T}\Phi(t_1,t_2,\ldots,t_{k-1})
d{M}_{t_1}^{(1)}
\stackrel{\rm def}{=}
$$
$$
\stackrel{\rm def}{=}\hbox{\vtop{\offinterlineskip\halign{
\hfil#\hfil\cr
{\rm l.i.m.}\cr
$\stackrel{}{{}_{N\to \infty}}$\cr
}} }\sum_{l=0}^{N-1}\Delta{M}_{\tau_{l}}^{(k-1)}
\int\limits_{\tau_{l+1}}^T
d{M}_{t_{k-2}}^{(k-2)}
\ldots\int\limits_{t_{3}}^T
d{M}_{t_{2}}^{(2)}
\int\limits_{t_{2}}^T
\Phi(t_1,t_2,\ldots,t_{k-1})d{M}_{t_1}^{(1)},
$$
$$
S'[\Phi]_{T,t}^{(k-1)}=
\int\limits_t^T\ldots\int\limits_t^{t_{k-2}}
\Phi(t_1,\ldots,t_{k-1})
d{M}_{t_{k-1}}^{(k-1)}
\ldots d{M}_{t_1}^{(1)},\ \ \ k\ge 2.
$$

\vspace{1mm}

It is easy to demonstrate similarly to the proof of 
Theorem 3.5 that 
under the condition  
BII the stochastic  
integral
$\tilde S[\Phi]_{T,t}^{(k-1)}$ exists
and

\vspace{-2mm}
$$
S'[\Phi]_{T,t}^{(k-1)}=\tilde S[\Phi]_{T,t}^{(k-1)}\ \ \ 
\hbox{w. p. {\rm 1}}.
$$ 

\vspace{3mm}

In its turn, using this fact we can
prove the following theorem
similarly to the proof of Theorem 3.5.

{\bf Theorem 3.6}\ \cite{old-preprint} (1999)
(also see \cite{1}-\cite{12aa}, \cite{arxiv-25}).
{\it Let the conditions {\rm BI, BII} of this section are fulfilled
and $|\rho_{\tau}|\le K<\infty $ w. p. {\rm 1} for all 
$\tau\in[t,T].$
Then, the stochastic integral
$\hat S[\xi,\Phi]_{T,t}^{(k)}$ exists and for $k\ge 2$
$$
S[\xi,\Phi]_{T,t}^{(k)}=\hat S[\xi,\Phi]_{T,t}^{(k)}\ \ \ 
\hbox{w. p. {\rm 1}}.
$$
}

Theorem 3.6 is the generalization of Theorem 3.2 for the case  
of iterated stochastic integrals with respect to martingale.

Let us  consider two statements.

{\bf Theorem 3.7} \cite{old-preprint} (1999)
(also see \cite{1}-\cite{12aa}, \cite{arxiv-25}). {\it
Let the conditions of Theorem {\rm 3.5} are fulfilled and $h(\tau)$ 
is a continuous nonrandom function at the interval $[t,T]$.
Then
\begin{equation}
\label{32.007}
~~~~~~\int\limits_{t}^T \phi_\tau dM_\tau^{(k+1)}h(\tau) 
\hat S[\psi^{(k)}]_{T,\tau}=
\int\limits_{t}^T \phi_\tau h(\tau) dM_\tau^{(k+1)}
\hat S[\psi^{(k)}]_{T,\tau}\ \ \ \hbox{w. p. {\rm 1}},
\end{equation}

\noindent
where the stochastic  
integrals in {\rm (\ref{32.007})}
exist.}

{\bf Theorem 3.8} \cite{old-preprint} (1999)
(also see \cite{1}-\cite{12aa}, \cite{arxiv-25}). {\it Under the 
conditions of Theorem {\rm 3.5}
$$
\int\limits_{t}^T h(t_1)\int\limits_{t}^{t_1}\phi_\tau dM_\tau^{(k+2)}
dM_{t_1}^{(k+1)} \hat S[\psi^{(k)}]_{T,t_1}=
$$
\begin{equation}
\label{32.008}
=\int\limits_{t}^T \phi_\tau dM_\tau^{(k+2)}\int\limits_{\tau}^T
h(t_1)dM_{t_1}^{(k+1)}\hat S[\psi^{(k)}]_{T,t_1}\ \hbox{w.~p.~{\rm 1}},
\end{equation}

\noindent
where the stochastic integrals in {\rm (\ref{32.008})} exist.}

The proofs of Theorems 3.7 and 3.8 are similar to the proofs of 
Theorems 3.3 and 3.4 correspondingly.

{\bf Remark 3.4.}\
{\it The integration order replacement technique for iterated It\^{o} stochastic 
integrals {\rm (}Theorems {\rm 3.1--3.4)} 
has been successfully applied for construction of the so-called
unified Taylor--It\^{o} and Taylor--Stratonovich expansions {\rm (}see Chapter 
{\rm 4)}
as well as
for proof and development of the mean-square 
approximation method for iterated 
It\^{o} and Stratonovich stochastic integrals
based on generalized multiple Fourier series {\rm (}see Chapters {\rm 1}
and {\rm 2)}.}

\chapter{Four New Forms of the Taylor--It\^{o} and Taylor--Stratonovich 
Expansions and its
Application to the High-Order Strong Numerical Methods
for It\^{o} Stochastic Differential Equations}

The problem of the Taylor--It\^{o} and Taylor--Stratonovich 
expansions of the It\^{o} stochastic 
processes in a neighborhood of a fixed time moment is considered
in this chapter.
The classical forms of the Taylor--It\^{o} and Taylor--Stratonovich 
expansions are transformed to 
four new representations, which include
the minimal sets of different types of iterated It\^{o} and Stratonovich
stochastic integrals. 
Therefore, these representations
(the so-called unified Taylor--It\^{o} and Taylor--Stratonovich expansions)
are more convenient for constructing 
of the high-order strong numerical 
methods for It\^{o} SDEs.
Explicit one-step strong numerical schemes with the 
convergence orders 1.0, 1.5, 2.0, 2.5, and 3.0
based on the unified Taylor--It\^{o} and Taylor--Stratonovich expansions 
are derived.

\section{Introduction}

Let $(\Omega,$ ${\rm F},$ ${\sf P})$ be a complete probability space, let 
$\{{\rm F}_t, t\in[0,T]\}$ be a nondecreasing right-con\-ti\-nu\-o\-us family 
of $\sigma$-algebras of ${\rm F},$
and let ${\bf w}_t$ be a standard $m$-dimensional Wiener
process, which is
${\rm F}_t$-measurable for any $t\in[0, T].$ We assume that the components
${\bf w}_{t}^{(i)}$ $(i=1,\ldots,m)$ of this process are independent. Consider
an It\^{o} SDE in the integral form
\begin{equation}
\label{1.5.2}
{\bf x}_t={\bf x}_0+\int\limits_0^t {\bf a}({\bf x}_{\tau},\tau)d\tau+
\int\limits_0^t B({\bf x}_{\tau},\tau)d{\bf w}_{\tau},\ \ \
{\bf x}_0={\bf x}(0,\omega).
\end{equation}
\noindent
Here ${\bf x}_t$ is some $n$-dimensional stochastic process 
satisfying to the It\^{o} SDE (\ref{1.5.2}). 
The nonrandom functions ${\bf a}: {\bf R}^n\times[0, T]\to{\bf R}^n$,
$B: {\bf R}^n\times[0, T]\to{\bf R}^{n\times m}$
guarantee the existence and uniqueness (up to stochastic 
equivalence) of a strong solution
to the equation (\ref{1.5.2}) \cite{Gih1}. The second integral on 
the right-hand side of (\ref{1.5.2}) is 
interpreted as an It\^{o} stochastic integral.
Let ${\bf x}_0$ be an $n$-dimensional random variable, which is 
${\rm F}_0$-measurable and 
${\sf M}\bigl\{\left|{\bf x}_0\right|^2\bigr\}<\infty$.
Also we assume that
${\bf x}_0$ and ${\bf w}_t-{\bf w}_0$ are independent when $t>0.$

It is well known \cite{Zapad-3}, \cite{Zapad-4}, \cite{Zapad-9},
\cite{Arato}, \cite{Shir}
(also see \cite{12}) that It\^{o} SDEs are 
adequate mathematical models of dynamic systems of 
different physical nature that are affected by random
perturbations. For example, It\^{o} SDEs are used as 
mathematical models in stochastic mathematical finance,
hydrology, seismology, geophysics, chemical kinetics, 
population dynamics, electrodynamics, medicine and other fields 
\cite{Zapad-3}, \cite{Zapad-4}, \cite{Zapad-9},
\cite{Arato}, \cite{Shir}
(also see \cite{12}).

Numerical integration of It\^{o} SDEs based on the 
strong convergence criterion of approximations \cite{Zapad-3} is widely 
used for the numerical simulation of sample trajectories of 
solutions to It\^{o} SDEs (which is required for constructing new 
mathematical models on the basis of such equations and for the 
numerical solution of different mathematical problems
connected with It\^{o} SDEs).
Among these problems, we note the following:
filtering of signals under influence of random noises 
in various statements (linear Kalman--Bucy filtering, 
nonlinear optimal filtering, filtering of continuous time Markov 
chains with a finite space of states, etc.), optimal stochastic 
control (including incomplete data control), 
testing estimation procedures of parameters
of stochastic systems, 
stochastic stability and bifurcations analysis
\cite{Zapad-1}, \cite{Zapad-3}, \cite{Zapad-4}, \cite{Zapad-8},
\cite{Zapad-9}, \cite{Lipt}, 
\cite{rybakov}-\cite{nasyrov2}.

Exact solutions of It\^{o} SDEs are known in rather rare cases. 
For this reason it is necessary 
to construct numerical procedures for solving these equations.

In this chapter, a promising approach 
\cite{Zapad-1}, \cite{Zapad-3}, \cite{Zapad-4},
\cite{Zapad-8}, \cite{Zapad-9} to the numerical 
integration of It\^{o} SDEs based on the stochastic analogues of the 
Taylor formula (Taylor--It\^{o} and Taylor--Stratonovich 
expansions) \cite{KlPl1},
\cite{PW1} (also see \cite{arxiv-24}, \cite{old-art-2}, 
\cite{Kul-Kuz}-\cite{Kuz-3}) is used. 
This approach uses a finite discretization of the time variable 
and the numerical simulation of the solution to the It\^{o} SDE at 
discrete time moments using the stochastic analogues of the Taylor 
formula mentioned above. 
A number of works (e.g., \cite{Zapad-1}-\cite{Zapad-4}, \cite{Zapad-8},
\cite{Zapad-9}) describe numerical schemes 
with the strong convergence orders $1.5, 2.0$, $2.5$, and $3.0$ for 
It\^{o} SDEs; however, they do not contain efficient procedures 
of the mean-square approximation of the iterated stochastic 
integrals for the case of multidimensional nonadditive noise. 

In this chapter, we consider the unified Taylor--It\^{o} 
and Taylor--Stratonovich expansions \cite{Kul-Kuz}, 
\cite{Kuz} (also see \cite{arxiv-24}, \cite{old-art-2})
which makes 
it possible (in contrast with its classical 
analogues \cite{Zapad-3}, \cite{KlPl1}) 
to use the 
minimal sets of iterated It\^{o} and Stratonovich stochastic integrals; this is a 
simplifying factor 
for the numerical methods implementation. 
We prove 
the unified Taylor--It\^{o} expansion \cite{Kul-Kuz}
with using of the slightly different approach (which is taken from \cite{Kuz})
in comparison with the approach from \cite{Kul-Kuz}.
Moreover, we obtain another (second) version of 
the unified Taylor--It\^{o} expansion \cite{very-old-2},
\cite{Kul-Kuz-1}.
In addition we construct two new forms of the Taylor--Stratonovich expansion
(the so-called unified Taylor--Stratonovich expansions \cite{Kuz}).

It should be noted that in Chapter 5 
on the base of the results of Chapters 1, 2
we study methods 
of numerical simulation of specific iterated It\^{o} and Stratonovich
stochastic integrals of multiplicities $1, 2, 3, 4, 5,$ and $6$
from the Taylor--It\^{o} and Taylor--Stratonovich expansions.
These stochastic integrals
are used in the 
strong numerical methods for It\^{o} SDEs 
\cite{Zapad-1}, \cite{Zapad-3}, \cite{Zapad-4}, \cite{Zapad-8}
(also see \cite{12}).
To approximate the iterated It\^{o} and Stratonovich
stochastic integrals 
appearing in the numerical schemes with the strong 
convergence orders $1.0, 1.5, 2.0, 2.5,$ and $3.0$, 
the method of generalized multiple Fourier series (see Chapter 1)
and especially method of multiple Fourier--Legendre 
series will be applied in Chapter 5.
It is important that the method of 
generalized multiple Fourier series (Theorems 1.1, 1.16)
does not lead to the partitioning of the 
integration interval of the iterated It\^{o} and Stratonovich 
stochastic integrals under consideration; 
this interval length is the integration step of the numerical methods used 
to solve It\^{o} SDEs; therefore, it is already fairly small and does not 
need to be partitioned. 
Computational experiments \cite{1} show that the numerical 
simulation for iterated stochastic integrals
(in which the interval of integration is partitioned) leads to 
unacceptably high computational cost and accumulation of 
computation errors. 
Also note that the Legendre polynomials have
essential advantage (in a number of aspects) over the trigonomentric functions
(see Chapter 5) 
constructing the mean-square approximations of iterated It\^{o} and
Stratonovich stochastic integrals
in the framework of the 
method of generalized multiple Fourier series (Theorems 1.1, 1.16).

Let us consider the following 
iterated It\^{o} and Stratonovich stochastic integrals:
\begin{equation}
\label{ito4}
J[\psi^{(k)}]_{s,t}=\int\limits_t^s\psi_k(t_k) \ldots \int\limits_t^{t_{2}}
\psi_1(t_1) d{\bf w}_{t_1}^{(i_1)}\ldots
d{\bf w}_{t_k}^{(i_k)},
\end{equation}
\begin{equation}
\label{str4}
J^{*}[\psi^{(k)}]_{s,t}=
{\int\limits_t^{*}}^s
\psi_k(t_k) 
\ldots 
{\int\limits_t^{*}}^{t_{2}}
\psi_1(t_1) d{\bf w}_{t_1}^{(i_1)}\ldots
d{\bf w}_{t_k}^{(i_k)},
\end{equation}

\noindent
where $0\le t <s \le T$, every $\psi_l(\tau)$ $(l=1,\ldots,k)$ is a continuous nonrandom
function 
at the interval $[0,T],$ ${\bf w}_{\tau}^{(i)}$
$(i=1,\ldots,m)$ are independent standard Wiener processes,
${\bf w}_{\tau}^{(0)}=\tau,$
$i_1,\ldots,i_k = 0, 1,\ldots,m.$

It should be noted that one of the main problems 
when constructing the high-order strong numerical methods
for It\^{o} SDEs on the base of the Taylor--It\^{o} and
Taylor--Stratonovich expansions is the mean-square approximation of
the iterated It\^{o} and Stratonovich stochastic integrals (\ref{ito4})
and (\ref{str4}).
Obviously, in the absence of procedures for the numerical simulation 
of stochastic integrals, the mentioned numerical methods 
are unrealizable in practice.
For this reason, in Chapter 5
we give 
the extensive practical material on expansions
and mean-square approximations 
of iterated It\^{o} and Stratonovich stochastic integrals
of multiplicities 1 to 6 from the Taylor--It\^{o} and 
Taylor--Stratonovich expansions.
In Chapter 5, the main focus is on approximations based on 
multiple Fourier--Legendre series. Such approximations 
is more effective in comparison with the trigonometric 
approximations (see Sect.~5.2) at least 
for the numerical methods with the
strong convergence order 1.5 and higher \cite{art-4}, 
\cite{arxiv-12}.

The rest of this Chapter is organized as follows.
In Sect.~4.1 (below) we 
consider a brief review of publications on the 
problem of construction of the Taylor--It\^{o} and Taylor--Stratonovich
expansions for the solutions of It\^{o} SDEs.
Sect.~4.2 is devoted to some auxiliary lemmas.
In Sect.~4.3 we consider the classical Taylor--It\^{o} expansion
while Sect.~4.4 and Sect.~4.5 are devoted to 
first and second forms of the so-called unified
Taylor--It\^{o} expansion correspondingly.
The classical Taylor--Stratonovich expansion is considered in Sect.~4.6.
First and second forms of the
unified Taylor--Stratonovich expansion 
are derived in Sect.~4.7 and Sect.~4.8.
In Sect.~4.9 we give a comparative analysis 
of the unified Taylor--It\^{o} and Taylor--Stratonovich expansions
with the classical Taylor--It\^{o} and Taylor--Stratonovich expansions.
Application of the first form of the unified 
Taylor--It\^{o} expansion to the high-order
strong numerical methods for It\^{o} SDEs
is considered in Sect.~4.10.
In Sect.~4.11 we construct 
the high-order
strong numerical methods for It\^{o} SDEs
on the base of the first form of the 
unified Taylor--Stratonovich expansion.
 
Let us give a brief review of publications on the 
problem of construction of the Taylor--It\^{o} and Taylor--Stratonovich
expansions for the solutions of It\^{o} SDEs.
A few variants of a stochastic analog of the Taylor formula 
have been obtained in \cite{KlPl1}, \cite{PW1} 
(also see \cite{Zapad-1}, \cite{Zapad-3})
for the stochastic processes
in the form $R({\bf x}_s, s)$, where ${\bf x}_s$ is a solution 
of the It\^{o} SDE (\ref{1.5.2}) and $R: {\bf R}^n\times [0, T]
\to {\bf R}$ is a 
sufficiently smooth nonrandom function.

The first result in this direction called the It\^{o}--Taylor expansion 
has been obtained in \cite{PW1} (also see \cite{KlPl1}). This result gives
an expansion of the process $R({\bf x}_s, s)$ into a series 
such that every term (if $k>0$) contains the iterated It\^{o} stochastic
integral 

\vspace{-2mm}
\begin{equation}
\label{re11}
\int\limits_t^s\ldots \int\limits_t^{t_{2}}
d{\bf w}_{t_1}^{(i_1)}\ldots
d{\bf w}_{t_k}^{(i_k)}
\end{equation}

\noindent
as a multiplier factor, where $0\le t<s\le T,$ $i_1,\ldots,i_k=0,1,\ldots,m.$
Obviously, the iterated It\^{o} stochastic integral (\ref{re11})
is a particular case of (\ref{ito4}) for
$\psi_1(\tau),\dots,\psi_k(\tau)\equiv 1.$

In \cite{KlPl1} another expansion of 
the stochastic process $R({\bf x}_s, s)$ into 
a series has been
derived. The iterated Stratonovich 
stochastic integrals

\vspace{-2mm}
\begin{equation}
\label{str11}
{\int\limits_t^{*}}^s
\ldots {\int\limits_t^{*}}^{t_2}
d{\bf w}_{t_1}^{(i_1)}\ldots
d{\bf w}_{t_k}^{(i_k)}
\end{equation}

\noindent
were used instead of the iterated It\^{o} stochastic integrals; 
the corresponding expansion was called the 
Stratonovich--Taylor expansion. In the formula (\ref{str11}) the indices
$i_1,\ldots,i_k$ take values $0, 1,\ldots,m$.

In \cite{Kul-Kuz}
the It\^{o}--Taylor expansion \cite{KlPl1} 
is reduced to the
interesting and unexpected
form
(the so-called unified Taylor--It\^{o} expansion)
by special transformations (see Chapter 3). 
Every term of this
expansion (if $k>0$) contains the iterated It\^{o} stochastic
integral 

\vspace{-2mm}
\begin{equation}
\label{ll1}
\int\limits^ {s} _ {t} (s-t_
{k}) ^ {l_ {k}} 
\ldots \int\limits^ {t _ {2}} _ {t} (s-t _ {1}) ^ {l_ {1}} d
{\bf w} ^ {(i_ {1})} _ {t_ {1}} \ldots 
d {\bf w} _ {t_ {k}} ^ {(i_ {k})},
\end{equation}

\noindent
where
$l_1,\ldots,l_k=0, 1, 2,\ldots $ and $i_1,\ldots,i_k = 1,\ldots,m$.

It is worth to mention another form of the 
unified Taylor--It\^{o} expansion \cite{very-old-2},
\cite{Kul-Kuz-1}
(also see \cite{1}-\cite{12aa}). 
Terms of the latter expansion contain iterated 
It\^{o} stochastic integrals of the form
\begin{equation}
\label{ll11}
\int\limits^ {s} _ {t} (t-t _{k}) ^ {l_ {k}} 
\ldots \int\limits^ {t _ {2}} _ {t} (t-t _ {1}) ^ {l_ {1}} 
d{\bf w} ^ {(i_ {1})} _ {t_ {1}} \ldots 
d {\bf w} _ {t_ {k}} ^ {(i_ {k})},
\end{equation}
where
$l_1,\ldots,l_k=0, 1, 2,\ldots $ and $i_1,\ldots,i_k = 1,\ldots,m$.

Obviously that some of the iterated It\^{o} stochastic integrals 
(\ref{re11}) or (\ref{str11}) are connected by linear relations, 
while
this is not the case for integrals defined by
(\ref{ll1}), (\ref{ll11}). In this sense, the total quantity
of stochastic integrals defined by
(\ref{ll1}) or (\ref{ll11}) is minimal. 
Futhermore, in this chapter we construct two new forms of the 
Taylor--Stratonovich expansion (the so-called unified
Taylor--Stratonovich expansions) \cite{Kuz-3} (also see \cite{Kuz})
such that every term
(if $k>0$)
contains as a multiplier 
the iterated Stratonovich stochastic integral
of one of two types

\vspace{-4mm}
\begin{equation}
\label{a111}
{\int\limits_t^{*}}^{s}(t-t_k)^{l_{k}}\ldots
{\int\limits_t^{*}}^{t _ {2}}(t-t _ {1}) ^ {l_ {1}} 
d{\bf w} ^ {(i_ {1})} _ {t_ {1}} \ldots
d{\bf w}_{t_ {k}}^{(i_ {k})},
\end{equation}
\begin{equation}
\label{a112}
{\int\limits_t^{*}}^{s}(s-t_k)^{l_{k}}\ldots
{\int\limits_t^{*}}^{t _ {2}}(s-t _ {1}) ^ {l_ {1}} 
d{\bf w} ^ {(i_ {1})} _ {t_ {1}} \ldots
d{\bf w}_{t_ {k}}^{(i_ {k})},
\end{equation}

\noindent
where 
$l_1,\ldots,l_k=0, 1, 2,\ldots $, $i_1,\ldots,i_k = 1,\ldots,m$,
and $k=1, 2,\ldots $

It is not difficult to see that for the sets of iterated Stratonovich
stochastic integrals
(\ref{a111}) and (\ref{a112}) 
the property of minimality (see above) 
also holds as for 
the sets of iterated It\^{o} 
stochastic integrals (\ref{ll1}), (\ref{ll11}).

As we noted above, the main problem in implementation
of high-order strong numerical methods 
for It\^{o} SDEs is the mean-square approximation 
of the iterated stochastic integrals
(\ref{re11})--(\ref{a112}). 
Obviously, these stochastic integrals are particular cases
of the stochastic integrals 
(\ref{ito4}), (\ref{str4}).

Taking into account the results of 
Chapters 1, 2, 3, 5
and the minimality
of the sets of stochastic integrals (\ref{ll1})--(\ref{a112}), 
we conclude that the unified Taylor--It\^{o} and Taylor--Stratonovich
expansions based on the iterated stochastic
integrals (\ref{ll1})--(\ref{a112}) can be useful for
constructing of high-order strong numerical
methods with the convergence orders 1.5, 
2.0, 2.5, 3.0, 3.5, 4.0, 4.5, $\ldots$ 
for It\^{o} SDEs.

\section{Auxiliary Lemmas}

Let $(\Omega,{\rm F},{\sf P})$ be a complete probability
space and let $w(t,\omega)\stackrel{\sf def}{=}w_t:$ 
$[0, T]\times \Omega\rightarrow {\bf R}$
be the standard Wiener process
defined on the probability space $(\Omega,{\rm F},{\sf P}).$

Consider the family of $\sigma$-algebras
$\left\{{\rm F}_t,\ t\in[0,T]\right\}$ defined
on the probability space $(\Omega,{\rm F},{\sf P})$ and
connected
with the Wiener process $w_t$ in such a way that

1.\ ${\rm F}_s\subset {\rm F}_t\subset {\rm F}$\ for
$s<t.$

2.\ The Wiener process $w_t$ is ${\rm F}_t$-measurable for all
$t\in[0,T].$

3.\ The process $w_{t+\Delta}-w_{t}$ for all
$t\ge 0,$ $\Delta>0$ is independent with
the events of $\sigma$-algebra
${\rm F}_{t}.$

Let us consider
the class ${\rm M}_2([t, T])$ $(t\ge 0)$ of random functions 
$\xi(\tau,\omega)\stackrel{\sf def}{=}
\xi_{\tau}:$ $[t, T]\times\Omega \to {\bf R}$ 
(see Sect.~1.1.2).

Recall (see Sect.~2.1.1) that the class ${\rm Q}_m([t, T])$ $(t\ge 0)$
consists of It\^{o} processes 
$\eta_{\tau},$ $\tau \in [t, T]$ of the form
\begin{equation}
\label{z900a}
\eta_{\tau} = \eta_t + \int\limits_t^{\tau} a_sds +
\int\limits_t^{\tau} b_sdw_s, 
\end{equation}
where $(a_{\tau})^m, (b_{\tau})^m \in {\rm M}_2([t, T])$ and
$\lim\limits_{s\to\tau} {\sf M}\bigl\{\left\vert b_s - b_{\tau}\right\vert ^4\bigr\}=0$
for all $\tau \in [t, T].$
The second integral on the right-hand side of (\ref{z900a})
is the It\^{o} stochastic integral.

Also note that the definition of the Stratonovich stochastic integral
in the mean-square sense is given by (\ref{123321.2}) (Sect.~2.1.1)
and the relation between Stratonovich and It\^{o}
stochastic integrals (see Sect.~2.1.1) has the following form
\cite{str} (also see \cite{123789000} (Sect.~1.2.5)) 
\begin{equation}
\label{d11kuz}
~~~~~~~~~ {\int\limits_{t}^{*}}^T F(\eta_{\tau},\tau)dw_{\tau}=
\int\limits_t^T F(\eta_{\tau},\tau)dw_{\tau}+
\frac{1}{2}\int\limits_t^T \frac{\partial F}{\partial x}(\eta_{\tau},\tau)
b_{\tau}d\tau\ \ \ \hbox{w.~p.~1}.
\end{equation}

If the Wiener processes in (\ref{z900a}) and (\ref{d11kuz}) 
are independent, then
\begin{equation}
\label{d11ahoh}
{\int\limits_{t}^{*}}^T F(\eta_{\tau},\tau)dw_{\tau}=
\int\limits_t^T F(\eta_{\tau},\tau)dw_{\tau}\ \ \ \hbox{w.~p.~1}.
\end{equation}

Recall that a possible variant of conditions providing
the correctness of the formulas 
(\ref{d11kuz}) and (\ref{d11ahoh}) consists of the 
conditions:
$\eta_{\tau}\in {\rm Q}_4([t,T]),$
$F(\eta_{\tau},\tau)\in {\rm M}_2([t,T]),$
$F(x,\tau)\in C^{2,1}({\bf R}\times [t, T]),$
where $C^{2,1}({\bf R}\times [t, T])$ $(t\ge 0)$ is the space of functions
$F(x,\tau): {\bf R} \times [t, T] \to {\bf R}$ such that
$$
\left|\frac{\partial F}{\partial x}(x,\tau)\right|\le K,\ \ 
\left|\frac{\partial^2 F}{\partial x^2}(x,\tau) \right|\le K,\ \ 
\left|\frac{\partial F}{\partial \tau}(x,\tau)\right|\le K,\ \ 
\left|\frac{\partial^2 F}{\partial \tau \partial x}(x,\tau) \right|\le K
$$ 
for all $x\in {\bf R}$ and $\tau\in [t, T],$ where constant $K$ does not depend on $x,\tau.$

{\bf Remark~4.1.}\ {\it Note that if $F(x,\tau)=F_1(x)F_2(\tau),$ then it suffices to require
that $F(x,\tau)$ be twice differentiable with respect to $x$ 
$($with bounded derivatives$)$ and continuous with respect to $\tau$
$($instead of the condition $F(x,\tau)\in C^{2,1}({\bf R}\times [t, T])).$}

Also remind that ${\rm S}_2([t,T])$ $(t\ge 0)$  is a subset of ${\rm M}_2([t,T])$
and ${\rm S}_2([t,T])$
consists of the mean-square continuous 
random functions 
(see Sect.~3.1).

Let us apply Theorem 3.1 (see Sect.~3.2) to derive one property for It\^{o}
stochastic integrals.

{\bf Lemma 4.1} \cite{12a}-\cite{12aa}, \cite{arxiv-24}. {\it Let $h(\tau),$ 
$g(\tau),$ $G(\tau):$ $[t,s]\to{\bf R}$ 
be continuous nonrandom functions at the interval $[t, s]$ and let $G(\tau)$         
be a antiderivative of the function
$g(\tau).$ Furthermore$,$ let $\xi_{\tau}\in {\rm S}_2([t,s]).$
Then
$$
\int\limits_t^s g(\tau)\int\limits_t^{\tau}h(\theta)
\int\limits_t^{\theta}\xi_u d{\bf w}_u^{(i)}
d{\bf w}_{\theta}^{(j)}
d\tau=\int\limits_t^s(G(s)-G(\theta))h(\theta)
\int\limits_t^{\theta}\xi_u d{\bf w}_u^{(i)}d{\bf w}_{\theta}^{(j)}
$$
w.~p.~{\rm 1}, 
where $i, j =1, 2$ and ${\bf w}_{\tau}^{(1)}$, ${\bf w}_{\tau}^{(2)}$
are independent standard Wiener processes that are
${\rm F}_{\tau}$--measurable for all $\tau\in[t,s]$.}

{\bf Proof.} Applying Theorem 3.1 two times and Theorem 3.3, 
we get the following relations
$$
\int\limits_t^s g(\tau)\int\limits_t^{\tau}h(\theta)
\int\limits_t^{\theta}\xi_u d{\bf w}_u^{(i)}
d{\bf w}_{\theta}^{(j)}
d\tau=\int\limits_t^s \xi_u d{\bf w}_u^{(i)}
\int\limits_u^s h(\theta)d{\bf w}_{\theta}^{(j)}\int\limits_{\theta}^s
g(\tau)d\tau=
$$
$$
=G(s)\int\limits_t^s \xi_u d{\bf w}_u^{(i)}
\int\limits_u^s h(\theta)d{\bf w}_{\theta}^{(j)}-
\int\limits_t^s \xi_u d{\bf w}_u^{(i)}
\int\limits_u^s G(\theta)h(\theta)d{\bf w}_{\theta}^{(j)}=
$$
$$
=G(s)\int\limits_t^{s}h(\theta)
\int\limits_t^{\theta}\xi_u d{\bf w}_u^{(i)}d{\bf w}_{\theta}^{(j)}
-
\int\limits_t^{s}G(\theta)h(\theta)
\int\limits_t^{\theta}\xi_u d{\bf w}_u^{(i)}d{\bf w}_{\theta}^{(j)}=
$$
$$
=\int\limits_t^s(G(s)-G(\theta))h(\theta)
\int\limits_t^{\theta}\xi_u d{\bf w}_u^{(i)}d{\bf w}_{\theta}^{(j)}\ \ \ 
\hbox{w.~p.~1.}
$$

The proof of Lemma 4.1 is completed.

Let us consider an analogue of Lemma 4.1 for
Stratonovich stochastic integrals.

{\bf Lemma 4.2} \cite{Kuz} (also see \cite{1}-\cite{5}, 
\cite{11}-\cite{12aa}, \cite{arxiv-24}). 
{\it Let $h(\tau),$ 
$g(\tau),$ $G(\tau):$ $[t,s]\to{\bf R}$ 
be continuous nonrandom functions at the interval $[t, s]$ and let $G(\tau)$         
be a antiderivative of the function
$g(\tau).$ Furthermore, let $\xi_{\tau}^{(l)}\in {\rm Q}_4([t,s])$ 
and
$$
\xi_{\tau}^{(l)}=\int\limits_t^{\tau}a_u du+
\int\limits_t^{\tau}b_u d{\bf w}_u^{(l)},\ \ \ l=1, 2.
$$

Then
\begin{equation}
\label{a1xxx}
\int\limits_t^s g(\tau){\int\limits_t^{*}}^{\tau}h(\theta)
{\int\limits_t^{*}}^{\theta}\xi^{(l)}_u d{\bf w}_u^{(i)}
d{\bf w}_{\theta}^{(j)}
d\tau={\int\limits_t^{*}}^s(G(s)-G(\theta))h(\theta)
{\int\limits_t^{*}}^{\theta}\xi_u^{(l)} d{\bf w}_u^{(i)}d{\bf w}_{\theta}^{(j)}
\end{equation}
w.~p.~{\rm 1}, 
where $i, j, l=1, 2$ and ${\bf w}_{\tau}^{(1)}$, ${\bf w}_{\tau}^{(2)}$
are independent standard Wiener processes that are
${\rm F}_{\tau}$--measurable for all $\tau\in[t,s]$.}

{\bf Proof.} Under the conditions of Lemma 4.2, we can 
apply the equalities (\ref{d11kuz}) and (\ref{d11ahoh}) with
$F(x,\theta)\equiv x h(\theta)$ and
$$
\eta_{\theta}={\int\limits_t^{*}}^{\theta}\xi_u^{(l)} d{\bf w}_u^{(i)},
$$
since the function $xh(\theta)$ is sufficiently smooth (see Remark~4.1)
and the following obvious inclusions hold:
$\eta_{\theta}\in {\rm Q}_4([t,s]),$ $\eta_{\theta}h(\theta)\in {\rm M}_2([t,s]).$

Thus, we have the equalities
\begin{equation}
\label{a222}
{\int\limits_t^{*}}^{\tau}h(\theta)
{\int\limits_t^{*}}^{\theta}\xi^{(l)}_u d{\bf w}_u^{(i)}d{\bf w}_{\theta}^{(j)}
=\int\limits_t^{\tau}h(\theta)
{\int\limits_t^{*}}^{\theta}\xi_u^{(l)} 
d{\bf w}_u^{(i)}d{\bf w}_{\theta}^{(j)}+
\frac{1}{2}{\bf 1}_{\{i=j\}}\int\limits_t^{\tau}h(\theta)\xi_{\theta}^{(l)}
d\theta,
\end{equation}
\begin{equation}
\label{a2.hoho}
{\int\limits_t^{*}}^{\theta}\xi^{(l)}_u d{\bf w}_u^{(i)}
=\int\limits_t^{\theta}\xi^{(l)}_u d{\bf w}_u^{(i)}+
\frac{1}{2}{\bf 1}_{\{l=i\}}\int\limits_t^{\theta}b_u du
\end{equation}
w.~p.~{\rm 1}, where  ${\bf 1}_A$ 
is the indicator of the set $A$. 

Substituting the 
formulas (\ref{a222}) and (\ref{a2.hoho}) into the left-hand
side of the equality (\ref{a1xxx}) and applying Theorem 3.1 twice
and Theorem 3.3, 
we get the following relations
$$
\int\limits_t^s g(\tau){\int\limits_t^{*}}^{\tau}h(\theta)
{\int\limits_t^{*}}^{\theta}\xi_u^{(l)} d{\bf w}_u^{(i)}d{\bf w}_{\theta}^{(j)}
d\tau=
\int\limits_t^s \xi_u^{(l)} d{\bf w}_u^{(i)}
\int\limits_u^s h(\theta)d{\bf w}_{\theta}^{(j)}\int\limits_{\theta}^s
g(\tau)d\tau+
$$
$$
+\frac{1}{2}{\bf 1}_{\{l=i\}}\int\limits_t^s b_u du
\int\limits_u^s h(\theta)d{\bf w}_{\theta}^{(j)}
\int\limits_{\theta}^s g(\tau)d\tau
+\frac{1}{2}{\bf 1}_{\{i=j\}}\int\limits_t^s h(\theta)\xi_{\theta}^{(l)}
d\theta \int\limits_{\theta}^s g(\tau)d\tau=
$$
$$
=G(s)\left(\int\limits_t^s \xi_u^{(l)} d{\bf w}_u^{(i)}
\int\limits_u^s h(\theta)d{\bf w}_{\theta}^{(j)}+
\frac{1}{2}{\bf 1}_{\{i=j\}}\int\limits_t^{s}h(\theta)\xi_{\theta}^{(l)}
d\theta
+\right.
$$
$$
\left.+\frac{1}{2}{\bf 1}_{\{l=i\}}\int\limits_t^s b_u du
\int\limits_u^s h(\theta)d{\bf w}_{\theta}^{(j)}
\right)-
$$
$$
-\left(\int\limits_t^s \xi_u^{(l)} d{\bf w}_u^{(i)}
\int\limits_u^s G(\theta)h(\theta)d{\bf w}_{\theta}^{(j)}+
\frac{1}{2}{\bf 1}_{\{i=j\}}\int\limits_t^{s}
G(\theta)h(\theta)\xi_{\theta}^{(l)}d\theta
+\right.
$$
$$
\left.+\frac{1}{2}{\bf 1}_{\{l=i\}}\int\limits_t^s b_u du
\int\limits_u^s h(\theta)G(\theta)d{\bf w}_{\theta}^{(j)}
\right)=
$$
$$
=G(s)\left(\int\limits_t^{s}h(\theta)
\int\limits_t^{\theta}\xi_u^{(l)} d{\bf w}_u^{(i)}d{\bf w}_{\theta}^{(j)}+
\frac{1}{2}{\bf 1}_{\{i=j\}}\int\limits_t^{s}h(\theta)\xi_{\theta}^{(l)}
d\theta+\right.
$$
$$
\left.+\frac{1}{2}{\bf 1}_{\{l=i\}}\int\limits_t^s
h(\theta)\int\limits_t^{\theta}b_u du d{\bf w}_{\theta}^{(j)}
\right)-
$$
$$
-\left(\int\limits_t^{s}G(\theta)h(\theta)
\int\limits_t^{\theta}\xi_u^{(l)} d{\bf w}_u^{(i)}d{\bf w}_{\theta}^{(j)}+
\frac{1}{2}{\bf 1}_{\{i=j\}}\int\limits_t^{s}
G(\theta)h(\theta)\xi_{\theta}^{(l)}d\theta
+\right.
$$
\begin{equation}
\label{a3xxx}
\left.+\frac{1}{2}{\bf 1}_{\{l=i\}}\int\limits_t^s
h(\theta)G(\theta)\int\limits_t^{\theta}b_u du d{\bf w}_{\theta}^{(j)}
\right)
\end{equation}

\noindent
w.~p.~1. Applying successively the formulas (\ref{a222}),
(\ref{a2.hoho})
together with the 
formula (\ref{a222}) in which $h(\theta)$ replaced 
by $G(\theta)h(\theta)$ as well as
the relation (\ref{a3xxx}), we obtain the equality (\ref{a1xxx}). 
The proof of Lemma 4.2 is completed.

\section{The Taylor--It\^{o} Expansion}

In this section, we use the Taylor--It\^{o} expansion \cite{KlPl1} 
and introduce some necessary notations.
At that we will use the original notations introduced by the author
of this book.

Let $C^{2,1}({\bf R}^n\times [0,T])\stackrel{{\sf def}}{=}{\rm L}$ 
be the space of functions $R({\bf x}, t): 
{\bf R}^n\times [0, T] \to {\bf R}$ with the 
following property: these functions are continuous and twice
continuously differentiable in ${\bf x}$ and have one continuous 
derivative in $t$. We consider the following operators on
the space ${\rm L}$
$$
L= {\partial \over \partial t}
+ \sum^ {n} _ {i=1} a^{(i)} ({\bf x},  t) 
{\partial \over  \partial  {\bf  x}^{(i)}} +
$$
\begin{equation}
\label{2.3}
+ {1\over 2} \sum^ {m} _ {j=1} \sum^ {n} _ {l,i=1}
B^ {(lj)} ({\bf x}, t) B^ {(ij)} ({\bf x}, t) {\partial
^{2} \over\partial{\bf x}^{(l)}\partial{\bf x}^{(i)}},
\end{equation}

\begin {equation}
\label{2.4}
G^ {(i)} _ {0} = \sum^ {n} _ {j=1} B^ {(ji)} ({\bf x}, t)
{\partial \over \partial {\bf x} ^ {(j)}},\ \ \
i=1,\ldots,m,
\end {equation}
where ${\bf x}^{(j)}$ is the $j$th component of
${\bf x}$, $a^{(j)} ({\bf x},  t)$ is the $j$th component
of $a({\bf x},  t)$, and $B^ {(ij)} ({\bf x}, t)$
is the $ij$th element of $B({\bf x}, t)$.

By the It\^{o} formula, we have the equality
\begin{equation}
\label{sa}
~~~~~~~~R({\bf x}_s,s)= R({\bf x}_t,t) + 
\int\limits_t^s LR({\bf x}_\tau,\tau)d\tau
+\sum^m_{i=1}\int\limits_t^s G_0^{(i)}R({\bf x}_\tau,\tau)d{\bf w}_\tau^{(i)}
\end{equation}
w.~p.~1, where  
${\bf x}_t$ is a strong solution of the It\^{o} SDE (\ref{1.5.2}),
$0 \le t<s \le T.$ In the formula (\ref{sa})
it is assumed that the functions ${\bf a}({\bf x}, t)$, 
$B({\bf x}, t)$, and
$R({\bf x}, t)$ satisfy the following condition: 
$LR({\bf x}_{\tau},\tau)$, $G_0^{(i)}R({\bf x}_{\tau},\tau)
\in {\rm M}_2([0, T])$
for $i = 1,\ldots,m$.

Introduce the following notations
\begin{equation}
\label{999.003}
{}^{(k)}A =\Biggl\| A^{(i_{1}\ldots  i_{k})}
\Biggr\|^{m_{1}~\ldots~~ m_{k}}_{i_{1}=1,\ldots ,i_{k}=1},\ \ \ 
m_1,\ldots,m_k\ge 1,
\end{equation}
$$
{ }^{(k+l)}A \stackrel{l}{\cdot}{}^{(l)}B
=\begin{cases} 
\Biggl\|\sum\limits^{m_{1}}_{i_1=1}\ldots\sum\limits_{i_l=1}^{m_l}
A^{(i_{1}\ldots i_{k+l})}B^{(i_{1}
\ldots  i_{l})}\Biggr\|^{m_{l+1}~\ldots ~~m_{l+k}}_{i_{l+1}=1,
\ldots ,i_{l+k}=1}\ \ \ &\hbox{for}\ \ \ k\ge 1\cr\cr
\sum\limits^{m_{1}}_{i_1=1}\ldots\sum\limits_{i_l=1}^{m_l}
A^{(i_{1}\ldots i_{l})}B^{(i_{1}
\ldots  i_{l})}\ \ \ &\hbox{for}\ \ \ k= 0
\end{cases},
$$
\begin{equation}
\label{sa901}
\Biggl\| A_{k+1}D_{k}^{(i_k)}A_k\ldots A_{2}D_{1}^{(i_1)}A_1
R({\bf x},t)
\Biggr\|^{m_{1}~\ldots~~ m_{k}}_{i_{1}=1,\ldots ,i_{k}=1}=
{}^{(k)}A_{k+1}D_{k}A_k\ldots A_{2}D_{1}A_1
R({\bf x},t),
\end{equation}
where $A_p$ and $D_q^{(i_q)}$ are operators defined on the 
space ${\rm L}$ for $p = 1,\ldots,k+1$, 
$q = 1,\ldots,k$, and $i_q = 1,\ldots,m_q$. It
is assumed that the left-hand side of (\ref{sa901}) exists. 
The symbol $\stackrel{0}{\cdot}$
is treated as the usual multiplication. If $m_l = 0$
in (\ref{999.003}) for some $l\in\{1,\ldots,k\}$, 
then the right-hand side of (\ref{999.003}) 
is treated as
$$
\Biggl\| A^{(i_{1}\ldots i_{l-1}i_{l+1}\ldots i_{k})}
\Biggr\|^{m_{1}~\ldots~~ m_{l-1}~~m_{l+1}~\ldots~~ m_{k}}
_{i_{1}=1,\ldots, i_{l-1}=1, i_{l+1}=1,\ldots,i_{k}=1},
$$
\noindent
(shortly, ${}^{(k-1)}A$).

We also introduce the following notations
$$
\Biggl\Vert Q_{\lambda_l}^{(i_l)}\ldots Q_{\lambda_1}^{(i_1)}
R({\bf x},t)\Biggr\Vert_{i_1=\lambda_1,\ldots,i_l=\lambda_l}^{m\lambda_1
~\ldots~m\lambda_l}\stackrel{\sf def}{=}{}^{(p_l)}Q_{\lambda_l}
\ldots Q_{\lambda_1}R({\bf x},t),
$$
$$
{}^{(p_{k})}J_{(\lambda_{k}\ldots \lambda_1)s,t}
=\Biggl\Vert J_{(\lambda_{k}\ldots \lambda_1)s,t}^{(i_k\ldots
i_1)}\Biggr\Vert_
{i_1=\lambda_1,\ldots,i_k=\lambda_k}^{m\lambda_1~\ldots~m\lambda_k},
$$
$$
{M}_k=\biggl\{(\lambda_k,\ldots,\lambda_1):
\lambda_l=1\ \hbox{or}\ \lambda_l=0;\ l=1,\ldots,k\biggr\},\ \ \ k\ge 1,
$$
$$
J_{(\lambda_{k}\ldots \lambda_1)s,t}^{(i_k\ldots
i_1)}=
\int\limits_t^s\ldots
\int\limits_t^{t_{2}}
d{\bf w}_{t_{1}}^{(i_k)}\ldots
d{\bf w}_{t_k}^{(i_1)},\ \ \ k\ge 1,
$$
where $\lambda_l=1$ or $\lambda_l=0$,
$Q_{\lambda_l}^{(i_l)}={L}$
and $i_l=0$ for $\lambda_l=0,$  $Q_{\lambda_l}^{(i_l)}=G_0^{(i_l)}$
and $i_l=1,\ldots,m$ for $\lambda_l=1,$ 
$$
p_l=\sum\limits_{j=1}^l \lambda_j\ \ \ \hbox{for}\ \ \ l=1,\ldots, r+1,\ \ \
r\in{\bf N},
$$
${\bf w}_{\tau}^{(i)}$ $(i=1,\ldots,m)$ are ${\rm F}_{\tau}$-measurable 
for all $\tau\in [0, T]$
independent standard Wiener processes and
${\bf w}_{\tau}^{(0)}=\tau.$

Applying (\ref{sa}) to the process 
$R({\bf x}_s, s)$ repeatedly, we obtain the following 
Taylor--It\^{o} expansion
\cite{KlPl1}
$$
R({\bf x}_s, s)=R({\bf x}_t, t)+\sum_{k=1}^r \sum_{(\lambda_{k},\ldots,
\lambda_1)\in {M}_k}
{}^{(p_{k})}Q_{\lambda_{k}}\ldots Q_{\lambda_1}
R({\bf x}_t, t)\stackrel{p_k}{\cdot}
{}^{(p_{k})}J_{(\lambda_{k}\ldots \lambda_1)s,t}
+
$$
\begin{equation}
\label{5.7.11}
+\left(D_{r+1}\right)_{s,t}
\end{equation}

\noindent
w.~p.~1, where $s, t\in [0, T],$ $s>t,$
$$
\left(D_{r+1}\right)_{s,t}
=
$$
\begin{equation}
\label{5.7.12}
=\sum_{(\lambda_{r+1},\ldots,\lambda_{1})\in{M}_{r+1}}
\int\limits_t^s\ldots
\left(
\int\limits_t^{t_2}
{}^{(p_{r+1})}Q_{\lambda_{r+1}}\ldots
Q_{\lambda_{1}}
R({\bf x}_{t_1}, t_1)\stackrel{\lambda_{r+1}}{\cdot}
d{\bf w}_{t_{1}}\right)
\ldots
\stackrel{\lambda_1}{\cdot}d{\bf w}_{t_{r+1}}.
\end{equation}

\noindent
It is assumed that the right-hand sides of (\ref{5.7.11}), 
(\ref{5.7.12}) exist.

A possible variant of the conditions, under which the right-hand 
sides of (\ref{5.7.11}), (\ref{5.7.12}) exist is as follows

(i)\ $Q_{\lambda_l}^{(i_l)}\ldots Q_{\lambda_1}^{(i_1)}R({\bf x},t)
\in {\rm L}$ for all 
$(\lambda_l,\ldots,\lambda_1)\in\bigcup\limits_{g=1}^{r} M_g$;

(ii)\ $Q_{\lambda_l}^{(i_l)}\ldots Q_{\lambda_1}^{(i_1)}R({\bf x}_{\tau},\tau)
\in {\rm M}_2([0,T])$
for all 
$(\lambda_l,\ldots,\lambda_1)\in\bigcup\limits_{g=1}^{r+1}
{M}_g.$

Let us rewrite 
the expansion (\ref{5.7.11}) in the another form 
$$
R({\bf x}_s,s)=R({\bf x}_t,t)
+\sum_{k=1}^r \sum_{(\lambda_{k},\ldots,\lambda_1)
\in{M}_k}
\sum_{i_1=\lambda_1}^{m\lambda_1}
\ldots 
\sum_{i_k=\lambda_k}^{m\lambda_k}
Q_{\lambda_k}^{(i_k)}\ldots Q_{\lambda_1}^{(i_1)}
R({\bf x}_t,t){J}_{(\lambda_{k}\ldots \lambda_1)s,t}^{(i_k\ldots i_1)}
+
$$
$$
+\left(D_{r+1}\right)_{s,t}\ \ \ \hbox{w.~p.~1}.
$$

Denote
$$
{G}_{rk}=\biggl\{(\lambda_k,\ldots,\lambda_1):\ r+1\le 
2k-\lambda_1-\ldots-\lambda_k\le 2r\biggr\},
$$
$$
{E}_{qk}=\biggl\{(\lambda_k,\ldots,\lambda_1):\ 2k-\lambda_1-\ldots-
\lambda_k=q\biggr\},
$$
where $\lambda_l=1$ or $\lambda_l=0$ $(l=1,\ldots,k).$

The Taylor--It\^{o} expansion ordered according to the 
order of smallness (in the mean-square sense when $s\downarrow  t$) 
of its terms has the form
$$
R({\bf x}_s,s)=
R({\bf x}_t,t)+
$$
\begin{equation}
\label{5.6.1rrr}
+
\sum_{q,k=1}^r \sum_{(\lambda_{k},\ldots,\lambda_1)
\in{\rm E}_{qk}}
\sum_{i_1=\lambda_1}^{m\lambda_1}
\ldots 
\sum_{i_k=\lambda_k}^{m\lambda_k}
Q_{\lambda_k}^{(i_k)}\ldots Q_{\lambda_1}^{(i_1)}
R({\bf x}_t,t){J}_{(\lambda_{k}\ldots \lambda_1)s,t}^{(i_k\ldots
i_1)}+
\left(H_{r+1}\right)_{s,t}
\end{equation}

\noindent
w.~p.~1, where
\begin{equation}
\label{5.6.1rrrh}
\left(H_{r+1}\right)_{s,t}=\hspace{-1mm}
\sum_{k=1}^r \sum_{(\lambda_{k},\ldots,\lambda_1)
\in{\rm G}_{rk}}
\sum_{i_1=\lambda_1}^{m\lambda_1}
\ldots 
\sum_{i_k=\lambda_k}^{m\lambda_k}
Q_{\lambda_k}^{(i_k)}\ldots Q_{\lambda_1}^{(i_1)}
R({\bf x}_t,t)
{J}_{(\lambda_{k}\ldots \lambda_1)s,t}^{(i_k\ldots
i_1)}
+\left(D_{r+1}\right)_{s,t}\hspace{-1mm}.
\end{equation}

\section{The First Form of the Unified Taylor--It\^{o} Expansion}

In this section, we transform the right-hand side of (\ref{5.7.11}) by
Theorem 3.1 and Lemma 4.1 to a
representation including the iterated 
It\^{o} stochastic integrals (\ref{ll11}).

Denote
\begin{equation}
\label{opr1}
~~~~~~I^{(i_1\ldots i_k)} _ {{l_1 \ldots l_k}_{s, t}} 
=\int\limits_t^{s}(t-t_k)^{l_{k}}\ldots\
\int\limits_t^{t _ {2}}(t-t _ {1}) ^ {l_ {1}} 
d{\bf w} ^ {(i_ {1})} _ {t_ {1}} \ldots
d{\bf w}_{t_ {k}}^{(i_ {k})}\ \ \ \hbox{for}\ \ \ k\ge 1
\end{equation}
and
$$
I^{(i_1\ldots i_k)} _ {{l_1 \ldots l_k}_{s, t}}=1\ \ \ 
\hbox{for}\ \ \ k=0,
$$

\vspace{2mm}
\noindent
where $i_1,\ldots,i_k=1,\ldots,m.$ Moreover, let
$$
{}^{(k)}I_{{l_1\ldots l_k}_{s,t}}=\Biggl\Vert 
I^{(i_1\ldots i_k)} _ {{l_1 \ldots l_k}_{s, t}}\Biggr\Vert
_{i_1,\ldots,i_k=1}^{m},
$$
\begin{equation}
\label{a9}
~~~~~~~G_p^{(i)}\stackrel{\sf def}{=}\frac{1}{p}\left(
G_{p-1}^{(i)}L-LG_{p-1}^{(i)}\right),\ \ \
p=1, 2,\ldots,\ \ \ i=1,\ldots,m,
\end{equation}

\vspace{1mm}
\noindent
where $L$ and $G_0^{(i)},$ $i=1,\ldots,m$
are determined by the equalities
(\ref{2.3}), (\ref{2.4}). 

Denote
$$
{A}_q\stackrel{\sf def}{=}
\Biggl\{
(k,j,l_1,\ldots,l_k):\ k+j+\sum_{p=1}^k l_p=q;\ k,j,l_1,\ldots,l_k=0, 1,\ldots
\Biggr\},
$$
$$
\Biggl\Vert 
G_{l_1}^{(i_1)}\ldots  G_{l_k}^{(i_k)} L^j
R({\bf x},t) \Biggr\Vert
_{i_1,\ldots,i_k=1}^
{m}\stackrel{\sf def}{=}{}^{(k)}
G_{l_1}\ldots G_{l_k} L^j
R({\bf x},t),
$$
$$
L^j R({\bf x},t)\stackrel{\sf def}{=}
\begin{cases}\underbrace{L\ldots L}_j
R({\bf x},t)\ &\hbox{for}\ j\ge 1\cr\cr
R({\bf x},t)\ &\hbox{for}\ j=0
\end{cases}.
$$

\vspace{2mm}

{\bf Theorem 4.1.}\ {\it Let conditions {\rm (i), (ii)}
be satisfied. Then for any $s, t \in [0, T]$ such that $s>t$ 
and for any positive
integer $r,$ the following expansion takes place w.~p.~{\rm 1}
$$
R({\bf x}_s,s)=
R({\bf x}_t,t)+
$$

\vspace{-5mm}
\begin{equation}
\label{razl4}
+
\sum_{q=1}^r
\sum_{(k,j,l_1,\ldots,l_k) \in {\rm A}_q}
\frac{(s-t)^j}{j!}
\sum_{i_1,\ldots,i_k=1}^m G_{l_1}^{(i_1)}\ldots
G_{l_k}^{(i_k)}L^j R({\bf x}_t,t) 
I_{{l_1\ldots l_k}_{s,t}}^{(i_1\ldots i_k)}+
\left(D_{r+1}\right)_{s,t},
\end{equation}

\vspace{1mm}
\noindent
where $\left(D_{r+1}\right)_{s,t}$ is defined by {\rm (\ref{5.7.12})}.
}

{\bf Proof.} We claim that
$$
\sum_{(\lambda_{q},\ldots,
\lambda_1)\in {M}_q}
{}^{(p_{q})}Q_{\lambda_{q}}\ldots Q_{\lambda_1}
R({\bf x}_t,t)\stackrel{p_q}{\cdot}
{}^{(p_{q})}J_{(\lambda_{q}\ldots \lambda_1)s,t}=
$$
\begin{equation}
\label{a22}
~~~~~=
\sum_{(k,j,l_1,\ldots,l_k) \in {\rm A}_q} 
\frac{(s-t)^j}{j!}
\sum_{i_1,\ldots,i_k=1}^m G_{l_1}^{(i_1)}\ldots
G_{l_k}^{(i_k)}L^j R({\bf x}_t,t)
I_{{l_1\ldots l_k}_{s,t}}^{(i_1\ldots i_k)}
\end{equation}

\noindent
w.~p.~1. The equality (\ref{a22}) is valid for $q = 1$. Assume 
that (\ref{a22}) is valid for some $q > 1$. 
In this case, using the induction hypothesis, we obtain
$$
\sum_{(\lambda_{q+1},\ldots,
\lambda_1)\in {M}_{q+1}}
{}^{(p_{q+1})}Q_{\lambda_{1}}\ldots Q_{\lambda_{q+1}}
R({\bf x}_t,t)\stackrel{p_{q+1}}{\cdot}
{}^{(p_{q+1})}J_{(\lambda_{1}\ldots \lambda_{q+1})s,t}=
$$
$$
=\sum_{\lambda_{q+1}\in\{1,\ 0\}}
\int\limits_t^s
\sum_{(\lambda_{q},\ldots,
\lambda_1)\in {M}_{q}}
\Biggl({}^{(p_{q+1})}Q_{\lambda_{1}}
\ldots
Q_{\lambda_{q+1}}
R({\bf x}_t,t)\stackrel{p_{q}}{\cdot}
{}^{(p_{q})}J_{(\lambda_{1}\ldots \lambda_{q})\theta,t}\Biggr)
\stackrel{\lambda_{q+1}}{\cdot}d{\bf w}_{\theta}=
$$
$$
=\sum_{\lambda_{q+1}\in\{1,\ 0\}}
\int\limits_t^s
\sum_{(k,j,l_1,\ldots,l_k)\in{A}_q}
\frac{(\theta-t)^j}{j!}\times
$$
$$
\times
\Biggl({}^{(k+\lambda_{q+1})}G_{l_1}\ldots G_{l_k} L^j
Q_{\lambda_{q+1}}R({\bf x}_t,t)
\stackrel{k}{\cdot}
{}^{(k)}I_{{l_1\ldots l_k}_{s,t}}\Biggr)
\stackrel{\lambda_{q+1}}{\cdot}d{\bf w}_{\theta}=
$$
$$
=\sum_{(k,j,l_1,\ldots,l_k)\in{A}_q}\left(
{}^{(k)}G_{l_1}\ldots G_{l_k} L^{j+1}
R({\bf x}_t,t)
\stackrel{k}{\cdot}
\int\limits_t^s\frac{(\theta-t)^j}{j!}
{}^{(k)}I_{{l_1\ldots l_k}_{\theta,t}}d\theta+\right.
$$
\begin{equation}
\label{a30}
~~~~~~\left.+\left(
{}^{(k+1)}G_{l_1}\ldots G_{l_k} L^{j} G_0
R({\bf x}_t,t)
\stackrel{k}{\cdot}
\int\limits_t^s\frac{(\theta-t)^j}{j!}
{}^{(k)}I_{{l_1\ldots l_k}_{\theta,t}}\right)
\stackrel{1}{\cdot}d{\bf w}_{\theta}\right)
\end{equation}
w.~p.~1.

Using Lemma 4.1, we obtain
$$
\int\limits_t^s\frac{(\theta-t)^j}{j!}
{}^{(k)}I_{{l_1\ldots l_k}_{\theta,t}}d\theta=
$$
\begin{equation}
\label{a31}
=\frac{1}{(j+1)!}
\begin{cases}(s-t)^{j+1}\ &\hbox{for}\ k=0 \cr \cr 
(s-t)^{j+1} \cdot {}^{(k)}I_{{l_1\ldots l_k}_{s,t}}-
(-1)^{j+1} \cdot
{}^{(k)}I_{{l_1\ldots l_{k-1}\ l_k+j+1}_{s,t}}\ &\hbox{for}\ k>0
\end{cases}
\end{equation}

\noindent
w.~p.~1. In addition (see (\ref{opr1})) we get
\begin{equation}
\label{a32}
\int\limits_t^s\frac{(\theta-t)^j}{j!}
I^{(i_1\ldots i_k)}_{{l_1\ldots l_k}_{\theta,t}}
d{\bf w}_{\theta}^{(i_{k+1})}=
\frac{(-1)^j}{j!} I_{{l_1\ldots l_k j}_{s,t}}^{(i_1\ldots i_k i_{k+1})}
\end{equation}

\vspace{-2mm}
\noindent
in the notations just introduced.
Substitute (\ref{a31}) and (\ref{a32}) 
into the formula (\ref{a30}). Grouping summands in the 
obtained expression with
equal lower indices at iterated It\^{o} stochastic
integrals and using (\ref{a9}) as well as the equality
\begin{equation}
\label{a33}
~~~~~~~G_p^{(i)}R({\bf x},t)=\frac{1}{p!}
\sum_{q=0}^p(-1)^q C_p^q L^q G_0^{(i)} L^{p-q}
R({\bf x},t),\ \ \ C_p^q=\frac{p!}{q!(p-q)!}
\end{equation}

\noindent
(this equality follows from (\ref{a9})), we note that the obtained 
expression equals to
$$
\sum_{(k,j,l_1,\ldots,l_k)\in{A}_{q+1}}
\frac{(s-t)^j}{j!}
{}^{(k)} G_{l_1}\ldots G_{l_k} L^j\{\eta_t\} 
\stackrel{k}{\cdot}
{}^{(k)}I_{{l_1\ldots l_k}_{s,t}}
$$

\noindent
w.~p.~1. Summing the equalities (\ref{a22}) for $q = 1, 2,\ldots,r$ 
and applying the formula (\ref{5.7.11}), we obtain the expression
(\ref{razl4}). The proof is completed.

Let us order terms of the expansion (\ref{razl4}) according to 
their smallness orders as $s \downarrow t$ in the mean-square sense
$$
R({\bf x}_s,s)=
R({\bf x}_t,t)+
$$

\vspace{-8mm}
\begin{equation}
\label{t100}
+
\sum_{q=1}^r
\sum_{(k,j,l_1,\ldots,l_k) \in {\rm D}_q}
\frac{(s-t)^j}{j!}
\sum_{i_1,\ldots,i_k=1}^m  G_{l_1}^{(i_1)}\ldots
G_{l_k}^{(i_k)} L^j R({\bf x}_t,t)
I_{{l_1\ldots l_k}_{s,t}}^{(i_1\ldots i_k)}+
\left(H_{r+1}\right)_{s,t}
\end{equation}

\noindent
w.~p.~1, where 
$$
\left(H_{r+1}\right)_{s,t}=
\sum_{(k,j,l_1,\ldots,l_k) \in {\rm U}_r}
\frac{(s-t)^j}{j!}
\sum_{i_1,\ldots,i_k=1}^m  G_{l_1}^{(i_1)}\ldots
G_{l_k}^{(i_k)} L^j R({\bf x}_t,t)
I_{{l_1\ldots l_k}_{s,t}}^{(i_1\ldots i_k)}+
$$
$$
+
\left(D_{r+1}\right)_{s,t},
$$

\vspace{-4mm}
\begin{equation}
\label{w1xx}
{\rm D}_{q}=\left\{
(k, j, l_ {1},\ldots, l_ {k}): k + 2\left(j +
\sum\limits_{p=1}^k l_p\right)= q;\ k, j, l_{1},\ldots, 
l_{k} =0,1,\ldots\right\},
\end{equation}
$$
{\rm U}_{r}=\left\{
(k, j, l_ {1},\ldots, l_ {k}):
k + j +
\sum\limits_{p=1}^k l_p\le r,\right.
$$
\begin{equation}
\label{w2xx}
\left.  
k + 2\left(j +
\sum\limits_{p=1}^k l_p\right)\ge r+1;\ k, j, l_{1},\ldots, 
l_{k} =0,1,\ldots\right\},
\end{equation}

\noindent
and $\left(D_{r+1}\right)_{s,t}$ is defined by (\ref{5.7.12}).
Note that the remainder term $\left(H_{r+1}\right)_{s,t}$ in 
(\ref{t100}) has a higher order of smallness in the mean-square sense as
$s \downarrow t$ than the terms of the main part of 
the expansion (\ref{t100}).

\section{The Second Form of the Unified Taylor--It\^{o} Expansion}

Consider iterated It\^{o} stochastic integrals of the form
$$
J^{(i_1\ldots i_k)} _ {{l_1 \ldots l_k}_{s, t}} 
=\int\limits_t^{s}(s-t_k)^{l_{k}}\ldots\
\int\limits_t^{t _ {2}}(s-t _ {1}) ^ {l_ {1}} 
d{\bf w} ^ {(i_ {1})} _ {t_ {1}} \ldots
d{\bf w}_{t_ {k}}^{(i_ {k})}\ \ \ \hbox{for}\ \ \ k\ge 1
$$
and
$$
J^{(i_1\ldots i_k)} _ {{l_1 \ldots l_k}_{s, t}}=1\ \ \ \hbox{for}\ \ \
k=0,
$$

\noindent
where $i_1,\ldots,i_k=1,\ldots,m.$ 

The additive property of stochastic integrals and the 
Newton binomial formula imply the following
equality
\begin{equation}
\label{70}
~~~~I^{(i_1\ldots i_k)} _ {{l_1 \ldots l_k}_{s, t}}=
\sum_{j_1=0}^{l_1}\ldots \sum_{j_k=0}^{l_k}
\prod_{g=1}^k C_{l_g}^{j_g}(t-s)^{l_1+\ldots+l_k-j_1-\ldots-j_k}\
J^{(i_1\ldots i_k)} _ {{j_1 \ldots j_k}_{s, t}}\ \ \ \hbox{w.~p.~1},
\end{equation}
where 
$$
C_l^k=\frac{l!}{k!(l-k)!}
$$

\noindent
is the binomial coefficient. Thus, the Taylor--It\^{o} expansion 
of the process $\eta_s = $ $R({\bf x}_s, s)$, $s\in [0, T]$ 
can be constructed either using the iterated 
stochastic integrals 
$I^{(i_1\ldots i_k)} _ {{l_1 \ldots l_k}_{s, t}}$
similarly to the
previous section or using the iterated stochastic integrals 
$J^{(i_1\ldots i_k)} _ {{l_1 \ldots l_k}_{s, t}}.$
This is the main subject of this section.

Denote
$$
\Biggl\Vert 
J^{(i_1\ldots i_k)} _ {{l_1 \ldots l_k}_{s, t}}
\Biggr\Vert_{i_1,\ldots,i_k=1}^{m}\stackrel{\sf def}
{=}{}^{(k)}J_{{l_1\ldots l_k}_{s,t}},
$$
$$
\Biggl\Vert 
L^j  G_{l_1}^{(i_1)}\ldots  G_{l_k}^{(i_k)}
R({\bf x},t) \Biggr\Vert
_{i_1,\ldots,i_k=1}^
{m}\stackrel{\sf def}{=}{}^{(k)}
L^j  G_{l_1}\ldots  G_{l_k}
R({\bf x},t).
$$

{\bf Theorem 4.2.}\ {\it Let conditions {\rm (i), (ii)} be satisfied. 
Then for any $s, t\in [0, T]$ such that $s>t$ and for any positive
integer $r,$ the following expansion is valid w.~p.~{\rm 1}
$$
R({\bf x}_s,s)=
R({\bf x}_t,t)+
$$

\vspace{-5mm}
\begin{equation}
\label{razl44}
+
\sum_{q=1}^r
\sum_{(k,j,l_1,\ldots,l_k) \in {\rm A}_q}
\frac{(s-t)^j}{j!}
\sum_{i_1,\ldots,i_k=1}^m  L^j G_{l_1}^{(i_1)}\ldots
G_{l_k}^{(i_k)}R({\bf x}_t,t)
J_{{l_1\ldots l_k}_{s,t}}^{(i_1\ldots i_k)}+
\left(D_{r+1}\right)_{s,t},
\end{equation}

\vspace{1mm}
\noindent
where $\left(D_{r+1}\right)_{s,t}$ is defined by {\rm (\ref{5.7.12})}.
}

{\bf Proof.} To prove the theorem, we check the equalities
$$
\sum_{(k,j,l_1,\ldots,l_k) \in {\rm A}_q}
\frac{(s-t)^j}{j!}
\sum_{i_1,\ldots,i_k=1}^m  L^j G_{l_1}^{(i_1)}\ldots
G_{l_k}^{(i_k)}R({\bf x}_t,t)
J_{{l_1\ldots l_k}_{s,t}}^{(i_1\ldots i_k)}=
$$
\begin{equation}
\label{f11}
=\sum_{(k,j,l_1,\ldots,l_k) \in {\rm A}_q}
\frac{(s-t)^j}{j!}
\sum_{i_1,\ldots,i_k=1}^m  G_{l_1}^{(i_1)}\ldots
G_{l_k}^{(i_k)} L^j R({\bf x}_t,t)
I_{{l_1\ldots l_k}_{s,t}}^{(i_1\ldots i_k)}\ \ \ \hbox{w.~p.~1}
\end{equation}

\noindent
for $q=1,2,\ldots,r.$ To check (\ref{f11}), substitute the expression 
(\ref{70}) 
into the right-hand side of (\ref{f11}) and then use the formulas
(\ref{a9}), (\ref{a33}).

Let us order terms of the expansion (\ref{razl44}) according to 
their smallness orders as $s \downarrow t$ in the mean-square sense
$$
R({\bf x}_s,s)=
R({\bf x}_t,t)+
$$

\vspace{-6mm}
$$
+\sum_{q=1}^r
\sum_{(k,j,l_1,\ldots,l_k) \in {\rm D}_q} 
\frac{(s-t)^j}{j!}
\sum_{i_1,\ldots,i_k=1}^m  L^j  G_{l_1}^{(i_1)}\ldots
G_{l_k}^{(i_k)}R({\bf x}_t,t)
J_{{l_1\ldots l_k}_{s,t}}^{(i_1\ldots i_k)}+
\left(H_{r+1}\right)_{s,t}
$$

\vspace{1mm}
\noindent
w.~p.~1, where
$$
\left(H_{r+1}\right)_{s,t}=
\sum_{(k,j,l_1,\ldots,l_k) \in {\rm U}_r}
\frac{(s-t)^j}{j!}
\sum_{i_1,\ldots,i_k=1}^m  L^j G_{l_1}^{(i_1)}\ldots
G_{l_k}^{(i_k)}R({\bf x}_t,t)
J_{{l_1\ldots l_k}_{s,t}}^{(i_1\ldots i_k)}+
$$
$$
+
\left(D_{r+1}\right)_{s,t}.
$$

\vspace{2mm}

The remainder term $\left(D_{r+1}\right)_{s,t}$ 
is defined by (\ref{5.7.12}); the sets ${\rm D}_q$ and ${\rm U}_r$ 
are defined by 
(\ref{w1xx}) and (\ref{w2xx}), respectively.
Finally, we note that the convergence w.~p.~1 of the 
truncated Taylor--It\^{o} expansion (\ref{5.7.11}) (without
the remainder term $\left(D_{r+1}\right)_{s,t}$) to the process 
$R({\bf x}_s, s)$ as $r\to\infty$  
for all $s, t \in [0, T]$ such that $s>t$ and $T < \infty$
has been proved in \cite{Zapad-3} (Proposition 5.9.2).
Since the expansions (\ref{razl4}) and 
(\ref{razl44}) are obtained from the Taylor--It\^{o}
expansion (\ref{5.7.11}) without any additional
conditions, the truncated expansions (\ref{razl4}) and 
(\ref{razl44}) (without the 
reminder term $\left(D_{r+1}\right)_{s,t}$) under the
conditions of Proposition 5.9.2 \cite{Zapad-3} converge to the 
process $R({\bf x}_s, s)$ w.~p.~1 as $r\to\infty$ for all
$s, t \in [0, T]$ such that $s>t$ and $T < \infty.$

\section{The Taylor--Stratonovich Expansion}

In this section, we use the Taylor--Stratonovich expansion \cite{KlPl1} 
and introduce some necessary notations.
At that we will use the original notations introduced by the author
of this book.

Let us consider two classic results.

{\bf Proposition~4.1}\ \cite{Gih1}.\ {\it Suppose that the following conditions
are satisfied.

{\rm AI.}\ The functions ${\bf a}({\bf x},t),\
B_j({\bf x},t):\ {\bf R}^n\times[0, T]\to {\bf R}^{n}$ $(j=1,\ldots,m)$ 
are measurable for all
$({\bf x},t)\in {\bf R}^n \times [0,T],$
where
$B_j({\bf x},t)$ is the $j$th column
of the matrix $B({\bf x},t)$ {\rm(}see {\rm (\ref{1.5.2}))}.

{\rm AII.}\ For all ${\bf x},\ {\bf y} \in {\bf R}^n$
there exists a constant $K<\infty$ such that
$$
\left|{\bf a}({\bf x},t)-{\bf a}({\bf y},t)\right| + 
\sum_{j=1}^m \left|B_{j}({\bf x},t)-B_{j}({\bf y},t)\right|
\le K \left|{\bf x}-{\bf y}\right|,
$$
$$
\left|{\bf a}({\bf x},t)\right|^2 + 
\sum_{j=1}^m\left|B_j({\bf x},t)\right|^2 \le 
K^2\bigl(1+\left|{\bf x}\right|^2 \bigr),
$$

\noindent
where $\left|\cdot\right|$ is the Euclidean norm of the vector.

{\rm AIII.}\ A random variable ${\bf x}_0$ is ${\rm F}_0$-measurable and
${\sf M}\bigl\{\left|{\bf x}_0\right|^2\bigr\}<\infty.$

Then there exists a unique {\rm(}up to stochastic equivalence{\rm)}
and continuous w.~p.~{\rm 1} strong solution of the It\^{o}
SDE {\rm (\ref{1.5.2})}.}

{\bf Proposition~4.2}\ \cite{Zapad-3}.\ {\it Suppose that the conditions
{\rm AI--AIII} {\rm(}see Proposition~{\rm 4.1)}
are satisfied and
${\sf M}\bigl\{\left|{\bf x}_{t_0}\right|^{2n}\bigr\}<\infty$ $(n\ge 1).$
Then
$$
{\sf M}\bigl\{\left|{\bf x}_t\right|^{2n}\bigr\}\le \bigl(1+{\sf M}\bigl\{\left|{\bf x}_{t_0}\right|
^{2n}\bigr\}\bigr)e^{C(t-t_0)},
$$
$$
{\sf M}\bigl\{\left|{\bf x}_t - {\bf x}_{t_0}\right|^{2n}\bigr\}\le
C_1\bigl(1+{\sf M}\bigl\{\left|{\bf x}_{t_0}\right|^{2n}\bigr\}\bigr)(t-t_0)^n e^{C(t-t_0)},
$$

\vspace{1mm}
\noindent
where ${\bf x}_t$ is the solution of the It\^{o}
SDE {\rm (\ref{1.5.2}),}
$t\in[t_0,T],$ $T<\infty,$
constant $C_1$ $(C_1\in (0, \infty))$ depends only on 
$n, K, T-t_0,$ $C=2n(2n+1)K^2,$ $K<\infty$
is a constant.}

Assume that
$R({\bf x},t)\in {\rm L},$
$LR({\bf x}_{\tau},\tau)$, $G_0^{(i)}R({\bf x}_{\tau},\tau)
\in {\rm M}_2([0, T])$
for $i = 1,\ldots,$ $m$ and
consider the It\^{o} formula (\ref{sa}).

In addition, suppose that the function $G_0^{(i)}R({\bf x},t)$
$(i = 1,\ldots,m)$ is such that the formulas
(\ref{d11kuz}) and (\ref{d11ahoh}) can be applied. 
For example, assume that

{\rm 1}.\ $G_0^{(i)}R({\bf x},t)\in {\rm L},$ $i=1,\ldots,m.$

{\rm 2}.\ For all ${\bf x}, {\bf y}\in{\bf R}^n,$ $t, s\in[0,T],$
$i_1, i_2=1,\ldots,m$ and for some $\nu>0$
$$
\left|G_0^{(i_2)}G_0^{(i_1)}R({\bf x},t)-G_0^{(i_2)}G_0^{(i_1)}R({\bf y},t)\right|
\le K_1\left|{\bf x}-{\bf y}\right|,
$$
$$
\left|G_0^{(i_2)}G_0^{(i_1)}R({\bf x},t)\right|+
\left|LG_0^{(i_1)}R({\bf x},t)\right|\le K_1\left(1+\left|{\bf x}\right|\right),
$$
$$
\left|G_0^{(i_2)}G_0^{(i_1)}R({\bf x},s)-
G_0^{(i_2)}G_0^{(i_1)}R({\bf x},t)\right|
\le K_1\left|s-t\right|^{\nu}
\left(1+\left|{\bf x}\right|\right),
$$
where $K_1<\infty$ is a constant.

{\rm 3}.\ Conditions {\rm AI, AII} are fulfilled {\rm (see Proposition~4.1)}.

{\rm 4}.\ ${\sf M}\{\left|{\bf x}_0\right|^8\}<\infty.$

\vspace{1mm}

Indeed, using the above conditions, Proposition~4.2 and the 
elementary inequality
$(a+b)^2\le 2a^2+2b^2,$ we obtain 
$$
{\sf M}\left\{\left|
G_0^{(i_2)}G_0^{(i_1)}R({\bf x}_{s},s)-
G_0^{(i_2)}G_0^{(i_1)}R({\bf x}_{t},t)
\right|^4\right\}\le
$$
$$
\le 8{\sf M}\left\{\left|
G_0^{(i_2)}G_0^{(i_1)}R({\bf x}_{s},s)-
G_0^{(i_2)}G_0^{(i_1)}R({\bf x}_{t},s)
\right|^4\right\}+
$$
$$
+ 8{\sf M}\left\{\left|
G_0^{(i_2)}G_0^{(i_1)}R({\bf x}_{t},s)
-G_0^{(i_2)}G_0^{(i_1)}R({\bf x}_{t},t)
\right|^4\right\}\le
$$

\vspace{-3mm}
$$
\le 8K_1^4{\sf M}\bigl\{\left|{\bf x}_{s}-{\bf x}_{t}\right|^4 \bigr\}+
8K_1^4 \left|s-t\right|^{4\nu}{\sf M}\bigl\{\bigl(1+\left|{\bf x}_{t}\right|\bigr)^4\bigr\}\le
$$

\vspace{-4mm}
\begin{equation}
\label{2024str1}
\le C_2\left|s-t\right|^2+C_3\left|s-t\right|^{4\nu}\to 0\ \ \ \hbox{if}\ \ \ 
s-t\to 0,
\end{equation}

$$
{\sf M}\left\{\left|
G_0^{(i_2)}G_0^{(i_1)}R({\bf x}_{s},s)\right|^8\right\}
\le K_1^8 {\sf M}\bigl\{\bigl(1+\left|{\bf x}_{s}\right|\bigr)^8\bigr\}
\le
$$
\begin{equation}
\label{2024str2}
~~~~~~~~\le C_4\bigl(1+{\sf M}\bigl\{\left|{\bf x}_{s}\right|^8\bigr\}\bigr)\le
C_5\bigl(1+\bigl(1+{\sf M}\bigl\{\left|{\bf x}_{0}\right|^{8}\bigr\}\bigr)e^{Cs}\bigr)<\infty,
\end{equation}

\vspace{2mm}
\noindent
where $C_2,\ldots,C_5<\infty$ are constants, $t, s\in[0,T].$

Analogously, we get
\begin{equation}
\label{2024str3}
{\sf M}\left\{\left|
LG_0^{(i_1)}R({\bf x}_s,s)\right|^8\right\}<\infty,\ \ \ s\in[0,T].
\end{equation}

Applying the It\^{o} formula, we obtain w.~p.~1
\begin{equation}
\label{8888.01}
~~~~~R({\bf x}_s,s)=R({\bf x}_t,t)+\int\limits_t^s
LR({\bf x}_{\tau},\tau)d\tau+\sum_{i_1=1}^m
\int\limits_t^s
G_0^{(i_1)}R({\bf x}_{\tau},\tau)d{\bf w}_{\tau}^{(i_1)},
\end{equation}
$$
G_0^{(i_1)}R({\bf x}_s,s)=G_0^{(i_1)}R({\bf x}_t,t)+\int\limits_t^s
LG_0^{(i_1)}R({\bf x}_{\tau},\tau)d\tau+
$$
\begin{equation}
\label{8888.02}
+\sum_{i_2=1}^m
\int\limits_t^s
G_0^{(i_2)}G_0^{(i_1)}R({\bf x}_{\tau},\tau)d{\bf w}_{\tau}^{(i_2)},
\end{equation}
where $i_1, i_2=1,\ldots,m.$ 

Thus, using (\ref{2024str1})--(\ref{2024str3}),
(\ref{d11kuz}) and (\ref{d11ahoh}),
we have 
\begin{equation}
\label{8888.03}
~~~~\int\limits_t^s G_0^{(i)}R({\bf x}_{\tau},\tau)d{\bf w}_{\tau}^{(i)}=
{\int\limits_t^{*}}^s 
G_0^{(i)}R({\bf x}_{\tau},\tau)d{\bf w}_{\tau}^{(i)}-
\frac{1}{2}\int\limits_t^s G_0^{(i)}G_0^{(i)}R({\bf x}_{\tau},\tau)d\tau
\end{equation}
w.~p.~1, where $s, t\in [0, T],$ $s>t,$ $i = 1,\ldots,m.$

Using the relation (\ref{8888.03}), let us write 
(\ref{8888.01}) in the following form
\begin{equation}
\label{sa1}
R({\bf x}_s,s)= R({\bf x}_t,t) + 
\int\limits_t^s {\bar L}R({\bf x}_\tau,\tau)d\tau
+\sum^m_{i=1}{\int\limits_t^{*}}^s 
G_0^{(i)}R({\bf x}_\tau,\tau)d{\bf w}_\tau^{(i)}\ \ \ \hbox{w.~p.~1,}
\end{equation}
where
\begin{equation}
\label{2.4a}
{\bar L}R({\bf x},t)=LR({\bf x},t)-
\frac{1}{2}\sum^m_{i=1}G_0^{(i)}G_0^{(i)}R({\bf x},t).
\end{equation}

Introduce the following notations
$$
\Biggl\Vert D_{\lambda_l}^{(i_l)}\ldots D_{\lambda_1}^{(i_1)}
R({\bf x},t)\Biggr\Vert_{i_1=\lambda_1,\ldots,i_l=\lambda_l}^{m\lambda_1
~\ldots~m\lambda_l}\stackrel{\sf def}{=}{}^{(p_l)}D_{\lambda_l}
\ldots D_{\lambda_1}R({\bf x},t),
$$
$$
{}^{(p_{k})}J^{*}_{(\lambda_{k}\ldots \lambda_1)s,t}
=\Biggl\Vert J_{(\lambda_{k}\ldots \lambda_1)s,t}^{*(i_k\ldots
i_1)}\Biggr\Vert_
{i_1=\lambda_1,\ldots,i_k=\lambda_k}^{m\lambda_1~\ldots~m\lambda_k},
$$
$$
{M}_k=\biggl\{(\lambda_k,\ldots,\lambda_1):
\lambda_l=1\ \hbox{or}\ \lambda_l=0;\ l=1,\ldots,k\biggr\},\ \ \ k\ge 1,
$$
$$
J_{(\lambda_{k}\ldots \lambda_1)s,t}^{*(i_k\ldots
i_1)}=
{\int\limits_t^{*}}^s\ldots
{\int\limits_t^{*}}^{t_{2}}
d{\bf w}_{t_{1}}^{(i_k)}\ldots
d{\bf w}_{t_k}^{(i_1)},\ \ \ k\ge 1,
$$
where $\lambda_l=1$ or $\lambda_l=0$,
$D_{\lambda_l}^{(i_l)}={\bar L}$
and $i_l=0$ for $\lambda_l=0,$  $D_{\lambda_l}^{(i_l)}=G_0^{(i_l)}$
and $i_l=1,\ldots,m$ for $\lambda_l=1,$ 
$$
p_l=\sum\limits_{j=1}^l \lambda_j\ \ \ \hbox{for}\ \ \ l=1,\ldots, r+1,\ \ \
r\in{\bf N},
$$
${\bf w}_{\tau}^{(i)}$ $(i=1,\ldots,m)$ are ${\rm F}_{\tau}$-measurable 
for all $\tau\in [0, T]$
independent standard Wiener processes and
${\bf w}_{\tau}^{(0)}=\tau.$

Applying the formula (\ref{sa1}) to the process 
$R({\bf x}_s, s)$ repeatedly, we obtain the following 
Taylor--Stra\-to\-no\-vich expansion
\cite{KlPl1}
$$
R({\bf x}_s, s)=R({\bf x}_t, t)+\sum_{k=1}^r \sum_{(\lambda_{k},\ldots,
\lambda_1)\in {M}_k}
{}^{(p_{k})}D_{\lambda_{k}}\ldots D_{\lambda_1}
R({\bf x}_t, t)\stackrel{p_k}{\cdot}
{}^{(p_{k})}J^{*}_{(\lambda_{k}\ldots \lambda_1)s,t}
+
$$
\begin{equation}
\label{5.7.11xxx}
+\left(D_{r+1}\right)_{s,t}
\end{equation}

\noindent
w.~p.~1, where $s, t\in [0, T],$ $s>t,$
$$
\left(D_{r+1}\right)_{s,t}
=
$$
\begin{equation}
\label{5.7.12xxx}
=\sum_{(\lambda_{r+1},\ldots,\lambda_1)\in{M}_{r+1}}
{\int\limits_t^{*}}^s\ldots
\Biggl(
{\int\limits_t^{*}}^{t_2}
\hspace{-2mm}{}^{(p_{r+1})}D_{\lambda_{r+1}}\ldots
D_{\lambda_1}
R({\bf x}_{t_1}, t_1)\stackrel{\lambda_{r+1}}{\cdot}
d{\bf w}_{t_{1}}\Biggr)
\ldots
\stackrel{\lambda_1}{\cdot}d{\bf w}_{t_{r+1}}.
\end{equation}

Let us rewrite 
the expansion (\ref{5.7.11xxx}) in another form w.~p.~1
$$
R({\bf x}_s,s)=
R({\bf x}_t,t)+
$$
\begin{equation}
\label{5.7.11xxxh}
+\sum_{k=1}^r \sum_{(\lambda_{k},\ldots,\lambda_1)
\in{M}_k}
\sum_{i_1=\lambda_1}^{m\lambda_1}
\ldots 
\sum_{i_k=\lambda_k}^{m\lambda_k}
D_{\lambda_k}^{(i_k)}\ldots D_{\lambda_1}^{(i_1)}
R({\bf x}_t,t){J}_{(\lambda_{k}\ldots \lambda_1)s,t}^{*(i_k\ldots i_1)}
+\left(D_{r+1}\right)_{s,t}.
\end{equation}

Denote
$$
{G}_{rk}=\bigl\{(\lambda_k,\ldots,\lambda_1):\ r+1\le 
2k-\lambda_1-\ldots-\lambda_k\le 2r\bigr\},
$$
$$
{E}_{qk}=\bigl\{(\lambda_k,\ldots,\lambda_1):\ 2k-\lambda_1-\ldots-
\lambda_k=q\bigr\},
$$

\vspace{1mm}
\noindent
where $\lambda_l=1$ or $\lambda_l=0$ $(l=1,\ldots,k).$

Let us order terms of the Taylor--Stratonovich 
expansion (\ref{5.7.11xxx}) or (\ref{5.7.11xxxh}) according to 
their smallness orders as $s \downarrow t$ in the mean-square sense
$$
R({\bf x}_s,s)=
R({\bf x}_t,t)+
$$
\begin{equation}
\label{5.6.1rrrx}
+
\sum_{q,k=1}^r \sum_{(\lambda_{k},\ldots,\lambda_1)
\in{\rm E}_{qk}}
\sum_{i_1=\lambda_1}^{m\lambda_1}
\ldots 
\sum_{i_k=\lambda_k}^{m\lambda_k}
D_{\lambda_k}^{(i_k)}\ldots D_{\lambda_1}^{(i_1)}
R({\bf x}_t,t)
{J}_{(\lambda_{k}\ldots \lambda_1)s,t}^{*(i_k\ldots
i_1)}+
\left(H_{r+1}\right)_{s,t}
\end{equation}

\noindent
w.~p.~1, where
$$
\left(H_{r+1}\right)_{s,t}=
\sum_{k=1}^r \sum_{(\lambda_{k},\ldots,\lambda_1)
\in{\rm G}_{rk}}
\sum_{i_1=\lambda_1}^{m\lambda_1}
\ldots 
\sum_{i_k=\lambda_k}^{m\lambda_k}
D_{\lambda_k}^{(i_k)}\ldots D_{\lambda_1}^{(i_1)}
R({\bf x}_t,t)
{J}_{(\lambda_{k}\ldots \lambda_1)s,t}^{*(i_k\ldots
i_1)}+
$$
\begin{equation}
\label{5.6.1rrrxh}
+
\left(D_{r+1}\right)_{s,t}.
\end{equation}

The following two questions seem interesting.

1. Under what conditions do the right-hand sides
of the formulas (\ref{5.6.1rrrx})
and (\ref{5.6.1rrrxh})
exist for $r\ge 2$?

2. Is it possible to obtain another
representation of the remainder term (\ref{5.6.1rrrxh})
for $r\ge 2$?

Below we will provide compelling 
arquments in favor of the following two facts.

(A). First, one can construct the Taylor--Stratonovich expansion
(\ref{5.6.1rrrx}) $(r\ge 2)$ in such a way that its remainder term 
will coincide w.~p.~1 with the remainder term (\ref{5.6.1rrrh}) $(r\ge 2)$
of the Taylor--It\^{o} expansion (\ref{5.6.1rrr}) $(r\ge 2)$.

(B). Second, the truncated Taylor--Stratonovich expansion (\ref{5.6.1rrrx}) $(r\ge 2)$
(without the remainder term (\ref{5.6.1rrrxh}) $(r\ge 2)$)
will coincide w.~p.~1 with the truncated Taylor--It\^{o}
expansion (\ref{5.6.1rrr}) $(r\ge 2)$
(without the remainder term (\ref{5.6.1rrrh}) $(r\ge 2)$).

This means that the right-hand side of (\ref{5.6.1rrrx}) $(r\ge 2)$
(in which the remainder term will have the form (\ref{5.6.1rrrh}) $(r\ge 2)$)
will exist under the conditions (i), (ii) (see Sect.~4.3).

Let us begin our reasoning with Theorem~2.12 (see Sect.~2.4.1).
This theorem allows us to represent 
the iterated Stratonovich stochastic integral
of multiplicity $k$ $(k\in{\bf N})$
as a sum of iterated It\^{o} stochastic integrals
and its mathematical expectation.
It is obvious that it is possible to obtain an inverse
formula that will express the iterated It\^{o} stochastic integral (\ref{itoxxx})
as a sum of iterated Stratonovich stochastic integrals (\ref{str}).
Below we present the corresponding proposition.

\vspace{2mm}

{\bf Proposition~4.3.}\ {\it Suppose that
every $\psi_l(\tau)$ $(l=1,\ldots,k)$ is a continuous
function at the interval $[t, T]$.
Then$,$ the following relation between iterated
It\^{o} and Stra\-to\-no\-vich stochastic integrals 
\begin{equation}
\label{2024str11}
~~~~J[\psi^{(k)}]_{T,t}=J^{*}[\psi^{(k)}]_{T,t}+
\sum_{r=1}^{\left[k/2\right]}\frac{(-1)^r}{2^r}
\sum_{(s_r,\ldots,s_1)\in {\rm A}_{k,r}}
J^{*}[\psi^{(k)}]_{T,t}^{s_r,\ldots,s_1}\ \ \ \hbox{{\rm w.~p.~1}}
\end{equation}
is correct, 
where $\sum\limits_{\emptyset}$ is supposed to be equal to zero$,$
$J[\psi^{(k)}]_{T,t}$ and $J^{*}[\psi^{(k)}]_{T,t}$
are defined by {\rm (\ref{itoxxx})} and {\rm (\ref{str}),}
respectively$,$

\newpage
\noindent
$$
J^{*}[\psi^{(k)}]_{T,t}^{s_l,\ldots,s_1}\ \stackrel{\rm def}{=}\ 
\prod_{p=1}^l {\bf 1}_{\{i_{s_p}=
i_{s_{p}+1}\ne 0\}}\ \times
$$
$$
\times
{\int\limits_t^{*}}^T
\psi_k(t_k)\ldots 
{\int\limits_t^{*}}^{t_{s_l+3}}
\psi_{s_l+2}(t_{s_l+2})
\int\limits_t^{t_{s_l+2}}\psi_{s_l}(t_{s_l+1})
\psi_{s_l+1}(t_{s_l+1}) \times
$$
$$
\times
{\int\limits_t^{*}}^{t_{s_l+1}}
\psi_{s_l-1}(t_{s_l-1})
\ldots
{\int\limits_t^{*}}^{t_{s_1+3}}
\psi_{s_1+2}(t_{s_1+2})
\int\limits_t^{t_{s_1+2}}\psi_{s_1}(t_{s_1+1})
\psi_{s_1+1}(t_{s_1+1}) \times
$$
$$
\times
{\int\limits_t^{*}}^{t_{s_1+1}}
\psi_{s_1-1}(t_{s_1-1})
\ldots 
{\int\limits_t^{*}}^{t_2}
\psi_1(t_1)
d{\bf w}_{t_1}^{(i_1)}\ldots d{\bf w}_{t_{s_1-1}}^{(i_{s_1-1})}
dt_{s_1+1}d{\bf w}_{t_{s_1+2}}^{(i_{s_1+2})}\ldots
$$

$$
\ldots\
d{\bf w}_{t_{s_l-1}}^{(i_{s_l-1})}
dt_{s_l+1}d{\bf w}_{t_{s_l+2}}^{(i_{s_l+2})}\ldots d{\bf w}_{t_k}^{(i_k)},
$$

\noindent
where 
$$
{\rm A}_{k,l}
=\bigl\{(s_l,\ldots,s_1):\
s_l>s_{l-1}+1,\ldots,s_2>s_1+1,\ s_l,\ldots,s_1=1,\ldots,k-1\bigr\},
$$
$$
(s_l,\ldots,s_1)\in{\rm A}_{k,l},\ \ \ 
l=1,\ldots,\left[k/2\right],\ \ \
i_s=0, 1,\ldots,m,\ \ \
s=1,\ldots,k,
$$

\noindent
$[x]$ is an integer part of a real number $x,$
${\bf 1}_A$ is the indicator of the set $A.$}

\vspace{2mm}

For example, from Proposition~4.3 for $k=1, 2, 3, 4$ we obtain
the following equalities w.~p.~1
$$
\int\limits_t^T\psi_1(t_1)d{\bf w}_{t_1}^{(i_1)}=
{\int\limits_t^{*}}^T\psi_1(t_1)d{\bf w}_{t_1}^{(i_1)},
$$
$$
\int\limits_t^T\psi_2(t_2)
\int\limits_t^{t_2}\psi_1(t_1)d{\bf w}_{t_1}^{(i_1)}
d{\bf w}_{t_2}^{(i_2)}=
{\int\limits_t^{*}}^T\psi_2(t_2)
{\int\limits_t^{*}}^{t_2}\psi_1(t_1)d{\bf w}_{t_1}^{(i_1)}
d{\bf w}_{t_2}^{(i_2)}-
$$
$$
-\frac{1}{2}{\bf 1}_{\{i_1=i_2\ne 0\}}
\int\limits_t^T\psi_2(t_2)\psi_1(t_2)dt_2,
$$
$$
\int\limits_t^T\psi_3(t_3)\ldots\hspace{-0.5mm}
\int\limits_t^{t_2}\psi_1(t_1)d{\bf w}_{t_1}^{(i_1)}\ldots
d{\bf w}_{t_3}^{(i_3)}=
\hspace{-1mm}
{\int\limits_t^{*}}^T\psi_3(t_3)\ldots
{\int\limits_t^{*}}^{t_2}\psi_1(t_1)d{\bf w}_{t_1}^{(i_1)}\ldots
d{\bf w}_{t_3}^{(i_3)}-
$$
$$
-\frac{1}{2}{\bf 1}_{\{i_1=i_2\ne 0\}}
{\int\limits_t^{*}}^T
\psi_3(t_3)
\int\limits_t^{t_3}\psi_2(t_2)\psi_1(t_2)dt_2
d{\bf w}_{t_3}^{(i_3)}-
$$
$$
-\frac{1}{2}{\bf 1}_{\{i_2=i_3\ne 0\}}
\int\limits_t^T\psi_3(t_3)\psi_2(t_3)
{\int\limits_t^{*}}^{t_3}
\psi_1(t_1)d{\bf w}_{t_1}^{(i_1)}
dt_3,
$$
$$
\int\limits_t^T\psi_4(t_4)\ldots\hspace{-0.5mm}
\int\limits_t^{t_2}\psi_1(t_1)d{\bf w}_{t_1}^{(i_1)}\ldots
d{\bf w}_{t_4}^{(i_4)}=
\hspace{-1mm}
{\int\limits_t^{*}}^T
\psi_4(t_4)\ldots
{\int\limits_t^{*}}^{t_2}\psi_1(t_1)d{\bf w}_{t_1}^{(i_1)}\ldots
d{\bf w}_{t_4}^{(i_4)}
-
$$
$$
-\frac{1}{2}{\bf 1}_{\{i_1=i_2\ne 0\}}
{\int\limits_t^{*}}^T
\psi_4(t_4)
{\int\limits_t^{*}}^{t_4}
\psi_3(t_3)\int\limits_t^{t_3}\psi_1(t_2)\psi_2(t_2)dt_2
d{\bf w}_{t_3}^{(i_3)}d{\bf w}_{t_4}^{(i_4)} -
$$
$$
-\frac{1}{2}{\bf 1}_{\{i_2=i_3\ne 0\}}
{\int\limits_t^{*}}^T\psi_4(t_4)\int\limits_t^{t_4}
\psi_3(t_3)\psi_2(t_3)
{\int\limits_t^{*}}^{t_3}
\psi_1(t_1)
d{\bf w}_{t_1}^{(i_1)}dt_3d{\bf w}_{t_4}^{(i_4)} -
$$
$$
- \frac{1}{2}{\bf 1}_{\{i_3=i_4\ne 0\}}
\int\limits_t^T\psi_4(t_4)\psi_3(t_4)
{\int\limits_t^{*}}^{t_4}
\psi_2(t_2)
{\int\limits_t^{*}}^{t_2}
\psi_1(t_1)
d{\bf w}_{t_1}^{(i_1)}d{\bf w}_{t_2}^{(i_2)}dt_4 +
$$
$$
+ \frac{1}{4}{\bf 1}_{\{i_1=i_2\ne 0\}}{\bf 1}_{\{i_3=i_4\ne 0\}}
\int\limits_t^T\psi_4(t_4)\psi_3(t_4)\int\limits_t^{t_4}
\psi_2(t_2)\psi_1(t_2)dt_2 dt_4.
$$

Further, using Proposition~4.3, we obtain for $r\ge 2$
$$
R({\bf x}_t,t)+
\sum_{q,k=1}^r \sum_{(\lambda_{k},\ldots,\lambda_1)
\in{\rm E}_{qk}}
\sum_{i_1=\lambda_1}^{m\lambda_1}
\ldots 
\sum_{i_k=\lambda_k}^{m\lambda_k}
Q_{\lambda_k}^{(i_k)}\ldots Q_{\lambda_1}^{(i_1)}
R({\bf x}_t,t){J}_{(\lambda_{k}\ldots \lambda_1)s,t}^{(i_k\ldots
i_1)}=
$$
\begin{equation}
\label{2024str201}
=R({\bf x}_t,t)+
\sum_{q,k=1}^r \sum_{(\lambda_{k},\ldots,\lambda_1)
\in{\rm E}_{qk}}
\sum_{i_1=\lambda_1}^{m\lambda_1}
\ldots 
\sum_{i_k=\lambda_k}^{m\lambda_k}
D_{\lambda_k}^{(i_k)}\ldots D_{\lambda_1}^{(i_1)}
R({\bf x}_t,t)
{J}_{(\lambda_{k}\ldots \lambda_1)s,t}^{*(i_k\ldots
i_1)}
\end{equation}

\noindent
w.~p.~1, where notations are the same as in 
(\ref{5.6.1rrr}), (\ref{5.6.1rrrx}).
Thus, (A) and (B) take place.

\section{The First Form of the Unified Taylor--Stratonovich Expansion}

In this section, we transform the right-hand side of (\ref{5.7.11xxx}) by
Theorem 3.1 and Lemma 4.2 to a
representation including the iterated 
Stratonovich stochastic integrals (\ref{a111}).
Moreover, we will use the remainder term $\left(D_{r+1}\right)_{s,t}$
of the form (\ref{5.7.12xxx}).

Denote
\begin{equation}
\label{opr1xxx}
~~~~I^{*(i_1\ldots i_k)} _ {{l_1 \ldots l_k}_{s, t}} 
={\int\limits_t^{*}}^{s}(t-t_k)^{l_{k}}\ldots\
{\int\limits_t^{*}}^{t _ {2}}(t-t _ {1}) ^ {l_ {1}} 
d{\bf w} ^ {(i_ {1})} _ {t_ {1}} \ldots
d{\bf w}_{t_ {k}}^{(i_ {k})}\ \ \ \hbox{for}\ \ \ k\ge 1
\end{equation}
and
$$
I^{*(i_1\ldots i_k)} _ {{l_1 \ldots l_k}_{s, t}}=1\ \ \ 
\hbox{for}\ \ \ k=0,
$$

\vspace{1mm}
\noindent
where $i_1,\ldots,i_k=1,\ldots,m.$ 

Futhermore, let
$$
{}^{(k)}I^{*}_{{l_1\ldots l_k}_{s,t}}=\Biggl\Vert 
I^{*(i_1\ldots i_k)} _ {{l_1 \ldots l_k}_{s, t}}\Biggr\Vert
_{i_1,\ldots,i_k=1}^{m},
$$

\begin{equation}
\label{a9x}
~~~~~~~\bar G_p^{(i)}\stackrel{\sf def}{=}\frac{1}{p}\left(
\bar G_{p-1}^{(i)}\bar L-\bar L\bar G_{p-1}^{(i)}\right),\ \ \
p=1, 2,\ldots,\ \ \ i=1,\ldots,m,
\end{equation}

\vspace{3mm}
\noindent
where $\bar G_0^{(i)}\stackrel{\sf def}{=}G_0^{(i)},$
$i=1,\ldots,m.$
The operators $\bar L$ and $G_0^{(i)},$ $i=1,\ldots,m$
are determined by the equalities
(\ref{2.3}), (\ref{2.4}), and (\ref{2.4a}). 

Denote

\vspace{-5mm}
$$
{A}_q\stackrel{\sf def}{=}
\left\{
(k,j,l_1,\ldots,l_k):\ k+j+\sum_{p=1}^k l_p=q;\ k,j,l_1,\ldots,l_k=0, 1,\ldots
\right\},
$$

\vspace{1mm}
$$
\Biggl\Vert 
\bar G_{l_1}^{(i_1)}\ldots \bar G_{l_k}^{(i_k)}\bar L^j
R({\bf x},t) \Biggr\Vert
_{i_1,\ldots,i_k=1}^
{m}\stackrel{\sf def}{=}{}^{(k)}
\bar G_{l_1}\ldots \bar G_{l_k}\bar L^j
R({\bf x},t),
$$

\vspace{1mm}
$$
\bar L^j R({\bf x},t)\stackrel{\sf def}{=}
\begin{cases}\underbrace{\bar L\ldots\bar L}_j
R({\bf x},t)\ &\hbox{for}\ j\ge 1\cr\cr
R({\bf x},t)\ &\hbox{for}\ j=0
\end{cases}.
$$

\vspace{3mm}

{\bf Theorem 4.3} \cite{Kuz} (also see
\cite{1}-\cite{12aa}, \cite{arxiv-24}, \cite{Kuz-3}).
{\it Suppose that sufficient conditions are satisfied under which the right-hand 
sides of {\rm(\ref{5.7.11xxx}), (\ref{5.7.12xxx})} exist. 
Then for any $s, t \in [0, T]$ such that $s>t$ 
and for any positive
integer $r,$ the following expansion takes place w.~p.~{\rm 1}
$$
R({\bf x}_s,s)=
R({\bf x}_t,t)+
$$

\vspace{-8mm}
\begin{equation}
\label{razl4x}
+
\sum_{q=1}^r
\sum_{(k,j,l_1,\ldots,l_k) \in {\rm A}_q}
\frac{(s-t)^j}{j!}
\sum_{i_1,\ldots,i_k=1}^m \bar G_{l_1}^{(i_1)}\ldots
\bar G_{l_k}^{(i_k)}\bar L^j R({\bf x}_t,t)
I_{{l_1\ldots l_k}_{s,t}}^{*(i_1\ldots i_k)}+
\left(D_{r+1}\right)_{s,t},
\end{equation}

\noindent
where $\left(D_{r+1}\right)_{s,t}$ is defined by {\rm (\ref{5.7.12xxx})}.
}

{\bf Proof.} We claim that
$$
\sum_{(\lambda_{q},\ldots,
\lambda_1)\in {M}_q}
{}^{(p_{q})}D_{\lambda_{q}}\ldots D_{\lambda_1}
R({\bf x}_t,t)\stackrel{p_q}{\cdot}
{}^{(p_{q})}J^{*}_{(\lambda_{q}\ldots \lambda_1)s,t}=
$$
\begin{equation}
\label{a22x}
~~~~~~=
\sum_{(k,j,l_1,\ldots,l_k) \in {\rm A}_q}
\frac{(s-t)^j}{j!}
\sum_{i_1,\ldots,i_k=1}^m \bar G_{l_1}^{(i_1)}\ldots
\bar G_{l_k}^{(i_k)}\bar L^j R({\bf x}_t,t)
I_{{l_1\ldots l_k}_{s,t}}^{*(i_1\ldots i_k)}
\end{equation}

\noindent
w.~p.~1. The equality (\ref{a22x}) is valid for $q = 1$. Assume 
that (\ref{a22x}) is valid for some $q > 1$. In 
this case using the induction hypothesis we obtain
$$
\sum_{(\lambda_{q+1},\ldots,
\lambda_1)\in {M}_{q+1}}\
{}^{(p_{q+1})}D_{\lambda_{1}}\ldots D_{\lambda_{q+1}}
R({\bf x}_t,t)\stackrel{p_{q+1}}{\cdot}
{}^{(p_{q+1})}J^{*}_{(\lambda_{1}\ldots \lambda_{q+1})s,t}=
$$
$$
=\sum_{\lambda_{q+1}\in\{1,\ 0\}}
{\int\limits_t^{*}}^s
\hspace{-3mm}\sum_{(\lambda_{q},\ldots,
\lambda_1)\in {M}_{q}}
\Biggl({}^{(p_{q+1})}D_{\lambda_{1}}
\ldots
D_{\lambda_{q+1}}
R({\bf x}_t,t)\stackrel{p_{q}}{\cdot}
{}^{(p_{q})}J^{*}_{(\lambda_{1}\ldots \lambda_{q})\theta,t}\Biggr)
\stackrel{\lambda_{q+1}}{\cdot}d{\bf w}_{\theta}=
$$
$$
=\sum_{\lambda_{q+1}\in\{1,\ 0\}}
{\int\limits_t^{*}}^s
\sum_{(k,j,l_1,\ldots,l_k)\in{A}_q}
\frac{(\theta-t)^j}{j!}\times
$$
$$
\times
\Biggl({}^{(k+\lambda_{q+1})}\bar G_{l_1}\ldots \bar G_{l_k}\bar L^j
D_{\lambda_{q+1}}R({\bf x}_t,t)
\stackrel{k}{\cdot}
{}^{(k)}I^{*}_{{l_1\ldots l_k}_{s,t}}\Biggr)
\stackrel{\lambda_{q+1}}{\cdot}d{\bf w}_{\theta}=
$$
$$
=\sum_{(k,j,l_1,\ldots,l_k)\in{A}_q}\left(
{}^{(k)}\bar G_{l_1}\ldots \bar G_{l_k}\bar L^{j+1}
R({\bf x}_t,t)
\stackrel{k}{\cdot}
\int\limits_t^s\frac{(\theta-t)^j}{j!}
{}^{(k)}I^{*}_{{l_1\ldots l_k}_{\theta,t}}d\theta+\right.
$$
\begin{equation}
\label{a30x}
~~~~~~\left.+\left(
{}^{(k+1)}\bar G_{l_1}\ldots \bar G_{l_k}\bar L^{j}\bar G_0
R({\bf x}_t,t)
\stackrel{k}{\cdot}
{\int\limits_t^{*}}^s
\frac{(\theta-t)^j}{j!}
{}^{(k)}I^{*}_{{l_1\ldots l_k}_{\theta,t}}\right)
\stackrel{1}{\cdot}d{\bf w}_{\theta}\right)
\end{equation}
w.~p.~1.

Using Lemma 4.2, we obtain
$$
\int\limits_t^s\frac{(\theta-t)^j}{j!}
{}^{(k)}I^{*}_{{l_1\ldots l_k}_{\theta,t}}d\theta=
$$
\begin{equation}
\label{a31x}
=\frac{1}{(j+1)!}
\begin{cases}(s-t)^{j+1}\ &\hbox{for}\ k=0 \cr \cr 
(s-t)^{j+1} \cdot {}^{(k)}I^{*}_{{l_1\ldots l_k}_{s,t}}-
(-1)^{j+1} \cdot
{}^{(k)}I^{*}_{{l_1\ldots l_{k-1}\ l_k+j+1}_{s,t}}\ &\hbox{for}\ k>0
\end{cases}
\end{equation}

\noindent
w.~p.~1. In addition (see (\ref{opr1xxx})) we get
\begin{equation}
\label{a32x}
{\int\limits_t^{*}}^s\frac{(\theta-t)^j}{j!}
I^{*(i_1\ldots i_k)}_{{l_1\ldots l_k}_{\theta,t}}
d{\bf w}_{\theta}^{(i_{k+1})}=
\frac{(-1)^j}{j!}I_{{l_1\ldots l_k j}_{s,t}}^{*(i_1\ldots i_k i_{k+1})}
\end{equation}
in the notations just introduced.
Substitute (\ref{a31x}) and (\ref{a32x}) 
into the formula (\ref{a30x}). Grouping summands in the 
obtained expression with
equal lower indices at iterated Stratonovich stochastic
integrals and using (\ref{a9x}) as well as the equality
\begin{equation}
\label{a33x}
~~~~~~~\bar G_p^{(i)}R({\bf x},t)=\frac{1}{p!}
\sum_{q=0}^p(-1)^q C_p^q\bar L^q\bar G_0^{(i)}\bar L^{p-q}
R({\bf x},t),\ \ \ C_p^q=\frac{p!}{q!(p-q)!}
\end{equation}

\noindent
(this equality follows from (\ref{a9x})), we note that the obtained 
expression equals to
$$
\sum_{(k,j,l_1,\ldots,l_k)\in{A}_{q+1}}
\frac{(s-t)^j}{j!}
{}^{(k)}\bar G_{l_1}\ldots \bar G_{l_k}\bar L^j\{\eta_t\} 
\stackrel{k}{\cdot}
{}^{(k)}I^{*}_{{l_1\ldots l_k}_{s,t}}
$$

\noindent
w.~p.~1. Summing the equalities (\ref{a22x}) for $q = 1, 2,\ldots,r$ 
and applying the formula (\ref{5.7.11xxx}), we obtain the expression
(\ref{razl4x}). The proof is completed.

Let us order terms of the expansion (\ref{razl4x}) according to 
their smallness orders as $s \downarrow t$ in the mean-square sense
$$
R({\bf x}_s,s)=
R({\bf x}_t,t)+
$$

\vspace{-5mm}
\begin{equation}
\label{t100x}
+\sum_{q=1}^r
\sum_{(k,j,l_1,\ldots,l_k) \in {\rm D}_q}
\frac{(s-t)^j}{j!}
\sum_{i_1,\ldots,i_k=1}^m \bar G_{l_1}^{(i_1)}\ldots
\bar G_{l_k}^{(i_k)}\bar L^j R({\bf x}_t,t)
I_{{l_1\ldots l_k}_{s,t}}^{*(i_1\ldots i_k)}+
\left(H_{r+1}\right)_{s,t}
\end{equation}

\noindent
w.~p.~1, where 
$$
\left(H_{r+1}\right)_{s,t}=
\sum_{(k,j,l_1,\ldots,l_k) \in {\rm U}_r}\ 
\frac{(s-t)^j}{j!}
\sum_{i_1,\ldots,i_k=1}^m \bar G_{l_1}^{(i_1)}\ldots
\bar G_{l_k}^{(i_k)}\bar L^j R({\bf x}_t,t)
I_{{l_1\ldots l_k}_{s,t}}^{*(i_1\ldots i_k)}+
$$
$$
+
\left(D_{r+1}\right)_{s,t},
$$

\vspace{-3mm}
\begin{equation}
\label{w1x}
{\rm D}_{q}=\left\{
(k, j, l_ {1},\ldots, l_ {k}): k + 2\left(j +
\sum\limits_{p=1}^k l_p\right)= q;\ k, j, l_{1},\ldots, 
l_{k} =0,1,\ldots\right\},
\end{equation}
$$
{\rm U}_{r}=\left\{
(k, j, l_ {1},\ldots, l_ {k}):
k + j +
\sum\limits_{p=1}^k l_p\le r,\right.
$$
\begin{equation}
\label{w2x}
\left.  k + 2\left(j +
\sum\limits_{p=1}^k l_p\right)\ge r+1;\ k, j, l_
{1},\ldots, l_ {k} =0,1,\ldots\right\},
\end{equation}

\noindent
and $\left(D_{r+1}\right)_{s,t}$ is defined by (\ref{5.7.12xxx}).
Note that the remainder term $\left(H_{r+1}\right)_{s,t}$ in 
(\ref{t100x}) has a higher order of smallness in the mean-square sense as
$s \downarrow t$ than the terms of the main part of 
the expansion (\ref{t100x}).

Further, using Proposition~4.3 analogously to (\ref{2024str201}), we obtain for $r\ge 2$
$$
R({\bf x}_t,t) + \sum_{q=1}^r
\sum_{(k,j,l_1,\ldots,l_k) \in {\rm D}_q}
\frac{(s-t)^j}{j!}
\sum_{i_1,\ldots,i_k=1}^m G_{l_1}^{(i_1)}\ldots
G_{l_k}^{(i_k)}L^j R({\bf x}_t,t) 
I_{{l_1\ldots l_k}_{s,t}}^{(i_1\ldots i_k)}=
$$
\begin{equation}
\label{useful40}
=R({\bf x}_t,t) + \sum_{q=1}^r
\sum_{(k,j,l_1,\ldots,l_k) \in {\rm D}_q}
\frac{(s-t)^j}{j!}
\sum_{i_1,\ldots,i_k=1}^m \bar G_{l_1}^{(i_1)}\ldots
\bar G_{l_k}^{(i_k)}\bar L^j R({\bf x}_t,t)
I_{{l_1\ldots l_k}_{s,t}}^{*(i_1\ldots i_k)}
\end{equation}
w.~p.~1.

\section{The Second Form of the Unified Taylor--Stratonovich Expansion}

Consider iterated Stratonovich stochastic integrals of the form
$$
J^{*(i_1\ldots i_k)} _ {{l_1 \ldots l_k}_{s, t}} 
={\int\limits_t^{*}}^{s}(s-t_k)^{l_{k}}\ldots\
{\int\limits_t^{*}}^{t _ {2}}(s-t _ {1}) ^ {l_ {1}} 
d{\bf w} ^ {(i_ {1})} _ {t_ {1}} \ldots
d{\bf w}_{t_ {k}}^{(i_ {k})}\ \ \ \hbox{for}\ \ \ k\ge 1
$$
and
$$
J^{*(i_1\ldots i_k)} _ {{l_1 \ldots l_k}_{s, t}}=1\ \ \ \hbox{for}\ \ \
k=0,
$$

\noindent
where $i_1,\ldots,i_k=1,\ldots,m.$ 

The additive property of stochastic integrals and the 
Newton binomial formula imply the following
equality
\begin{equation}
\label{70x}
~~~I^{*(i_1\ldots i_k)} _ {{l_1 \ldots l_k}_{s, t}}=
\sum_{j_1=0}^{l_1}\ldots \sum_{j_k=0}^{l_k}
\prod_{g=1}^k C_{l_g}^{j_g}(t-s)^{l_1+\ldots+l_k-j_1-\ldots-j_k}\
J^{*(i_1\ldots i_k)} _ {{j_1 \ldots j_k}_{s, t}}\ \ \ \hbox{w.~p.~1},
\end{equation}
where 
$$
C_l^k=\frac{l!}{k!(l-k)!}
$$

\noindent
is the binomial coefficient. Thus, the Taylor--Stratonovich expansion 
of the process $\eta_s = R({\bf x}_s, s)$, $s\in [0, T]$ 
can be constructed either using the iterated 
stochastic integrals 
$I^{*(i_1\ldots i_k)} _ {{l_1 \ldots l_k}_{s, t}}$
similarly to the
previous section or using the iterated stochastic integrals 
$J^{*(i_1\ldots i_k)} _ {{l_1 \ldots l_k}_{s, t}}.$
This is the main subject of this section.

Denote
$$
\Biggl\Vert 
J^{*(i_1\ldots i_k)} _ {{l_1 \ldots l_k}_{s, t}}\Biggr\Vert
_{i_1,\ldots,i_k=1}^{m}\stackrel{\sf def}
{=}{}^{(k)}J^{*}_{{l_1\ldots l_k}_{s,t}},
$$

$$
\Biggl\Vert 
\bar L^j \bar G_{l_1}^{(i_1)}\ldots \bar G_{l_k}^{(i_k)}
R({\bf x},t) \Biggr\Vert
_{i_1,\ldots,i_k=1}^
{m}\stackrel{\sf def}{=}{}^{(k)}
\bar L^j \bar G_{l_1}\ldots \bar G_{l_k}
R({\bf x},t).
$$

\vspace{1mm}

{\bf Theorem 4.4} \cite{Kuz} (also see
\cite{1}-\cite{12aa}, \cite{arxiv-24}, \cite{Kuz-3}).
{\it Suppose that sufficient conditions are satisfied under which the right-hand 
sides of {\rm(\ref{5.7.11xxx}), (\ref{5.7.12xxx})} exist. 
Then for any $s, t\in [0, T]$ such that $s>t$ and for any positive
integer $r,$ the following expansion is valid w.~p.~{\rm 1}
$$
R({\bf x}_s,s)=
R({\bf x}_t,t)+
$$

\vspace{-5mm}
\begin{equation}
\label{razl44x}
+\sum_{q=1}^r
\sum_{(k,j,l_1,\ldots,l_k) \in {\rm A}_q}
\frac{(s-t)^j}{j!}
\sum_{i_1,\ldots,i_k=1}^m \bar L^j\bar G_{l_1}^{(i_1)}\ldots
\bar G_{l_k}^{(i_k)}R({\bf x}_t,t)
J_{{l_1\ldots l_k}_{s,t}}^{*(i_1\ldots i_k)}+
\left(D_{r+1}\right)_{s,t},
\end{equation}

\noindent
where $\left(D_{r+1}\right)_{s,t}$ is defined by {\rm (\ref{5.7.12xxx})}.
}

{\bf Proof.} To prove the theorem, we check the equalities
$$
\sum_{(k,j,l_1,\ldots,l_k) \in {\rm A}_q}
\frac{(s-t)^j}{j!}
\sum_{i_1,\ldots,i_k=1}^m \bar L^j\bar G_{l_1}^{(i_1)}\ldots
\bar G_{l_k}^{(i_k)}R({\bf x}_t,t) 
J_{{l_1\ldots l_k}_{s,t}}^{*(i_1\ldots i_k)}=
$$
\begin{equation}
\label{f11x}
~~~=\sum_{(k,j,l_1,\ldots,l_k) \in {\rm A}_q}\ 
\frac{(s-t)^j}{j!}
\sum_{i_1,\ldots,i_k=1}^m \bar G_{l_1}^{(i_1)}\ldots
\bar G_{l_k}^{(i_k)}\bar L^j R({\bf x}_t,t)
I_{{l_1\ldots l_k}_{s,t}}^{*(i_1\ldots i_k)}\ \ \ \hbox{w.~p.~1}
\end{equation}

\vspace{1mm}
\noindent
for $q=1,2,\ldots,r.$ 
To check (\ref{f11x}), substitute the expression 
(\ref{70x}) 
into the right-hand side of (\ref{f11x}) and then use the formulas
(\ref{a9x}), (\ref{a33x}).

Let us order terms of 
the expansion (\ref{razl44x}) according to 
their smallness orders as $s \downarrow t$ in the mean-square sense
$$
R({\bf x}_s,s)=
R({\bf x}_t,t)+
$$

\vspace{-5mm}
$$
+
\sum_{q=1}^r
\sum_{(k,j,l_1,\ldots,l_k) \in {\rm D}_q}\ 
\frac{(s-t)^j}{j!}
\sum_{i_1,\ldots,i_k=1}^m \bar L^j \bar G_{l_1}^{(i_1)}\ldots
\bar G_{l_k}^{(i_k)}R({\bf x}_t,t)
J_{{l_1\ldots l_k}_{s,t}}^{*(i_1\ldots i_k)}+
\left(H_{r+1}\right)_{s,t}
$$

\vspace{1mm}
\noindent
w.~p.~1, where
$$
\left(H_{r+1}\right)_{s,t}=
\sum_{(k,j,l_1,\ldots,l_k) \in {\rm U}_r}
\frac{(s-t)^j}{j!}
\sum_{i_1,\ldots,i_k=1}^m \bar L^j\bar G_{l_1}^{(i_1)}\ldots
\bar G_{l_k}^{(i_k)}R({\bf x}_t,t) 
J_{{l_1\ldots l_k}_{s,t}}^{*(i_1\ldots i_k)}+
$$
$$
+
\left(D_{r+1}\right)_{s,t}.
$$

\vspace{2mm}

The remainder term $\left(D_{r+1}\right)_{s,t}$ 
is defined by (\ref{5.7.12xxx}); the sets ${\rm D}_q$ and ${\rm U}_r$ 
are defined by (\ref{w1x}) and (\ref{w2x}), respectively.
Finally, we note that the convergence w.~p.~1 of the 
truncated Taylor--Stratonovich expansion (\ref{5.7.11xxx}) (without
the remainder term $\left(D_{r+1}\right)_{s,t}$) to the process 
$R({\bf x}_s, s)$ as $r\to\infty$  
for all $s, t \in [0, T]$ such that $s>t$ and $T < \infty$
has been proved in \cite{Zapad-3} (Proposition 5.10.2).
Since the expansions (\ref{razl4x}) and 
(\ref{razl44x}) are obtained from the Taylor--Stratonovich 
expansion (\ref{5.7.11xxx}) without any additional
conditions, the truncated expansions (\ref{razl4x}) and 
(\ref{razl44x}) (without the 
reminder term $\left(D_{r+1}\right)_{s,t}$) under the
conditions of Proposition 5.10.2 \cite{Zapad-3}
converge to the 
process $R({\bf x}_s, s)$ w.~p.~1 as $r\to\infty$ for all
$s, t \in [0, T]$ such that $s>t$ and $T < \infty.$

\section{A Remark on Theorems~4.3 and 4.4}

Note that when proving Theorems~4.3 and 4.4
we established the following equalities w.~p.~1
$$
R({\bf x}_t,t)+
\sum_{k=1}^r \sum_{(\lambda_{k},\ldots,\lambda_1)
\in{M}_k}
\sum_{i_1=\lambda_1}^{m\lambda_1}
\ldots 
\sum_{i_k=\lambda_k}^{m\lambda_k}
D_{\lambda_k}^{(i_k)}\ldots D_{\lambda_1}^{(i_1)}
R({\bf x}_t,t){J}_{(\lambda_{k}\ldots \lambda_1)s,t}^{*(i_k\ldots i_1)}
=
$$
\begin{equation}
\label{2024str202}
=R({\bf x}_t,t)+
\sum_{q=1}^r
\sum_{(k,j,l_1,\ldots,l_k) \in {\rm A}_q}
\frac{(s-t)^j}{j!}
\sum_{i_1,\ldots,i_k=1}^m \bar G_{l_1}^{(i_1)}\ldots
\bar G_{l_k}^{(i_k)}\bar L^j R({\bf x}_t,t)
I_{{l_1\ldots l_k}_{s,t}}^{*(i_1\ldots i_k)},
\end{equation}
$$
R({\bf x}_t,t)+
\sum_{k=1}^r \sum_{(\lambda_{k},\ldots,\lambda_1)
\in{M}_k}
\sum_{i_1=\lambda_1}^{m\lambda_1}
\ldots 
\sum_{i_k=\lambda_k}^{m\lambda_k}
D_{\lambda_k}^{(i_k)}\ldots D_{\lambda_1}^{(i_1)}
R({\bf x}_t,t){J}_{(\lambda_{k}\ldots \lambda_1)s,t}^{*(i_k\ldots i_1)}
=
$$
\begin{equation}
\label{2024str203}
=R({\bf x}_t,t)+
\sum_{q=1}^r
\sum_{(k,j,l_1,\ldots,l_k) \in {\rm A}_q}
\frac{(s-t)^j}{j!}
\sum_{i_1,\ldots,i_k=1}^m \bar L^j\bar G_{l_1}^{(i_1)}\ldots
\bar G_{l_k}^{(i_k)}R({\bf x}_t,t)
J_{{l_1\ldots l_k}_{s,t}}^{*(i_1\ldots i_k)}.
\end{equation}

It is easy to see that by analogy with 
(\ref{2024str202}) and (\ref{2024str203})
the following equalities can be obtained w.~p.~1
$$
R({\bf x}_t,t)+
\sum_{q,k=1}^r \sum_{(\lambda_{k},\ldots,\lambda_1)
\in{\rm E}_{qk}}
\sum_{i_1=\lambda_1}^{m\lambda_1}
\ldots 
\sum_{i_k=\lambda_k}^{m\lambda_k}
D_{\lambda_k}^{(i_k)}\ldots D_{\lambda_1}^{(i_1)}
R({\bf x}_t,t)
{J}_{(\lambda_{k}\ldots \lambda_1)s,t}^{*(i_k\ldots
i_1)}=
$$
\begin{equation}
\label{2024str204}
=R({\bf x}_t,t)+
\sum_{q=1}^r
\sum_{(k,j,l_1,\ldots,l_k) \in {\rm D}_q}
\frac{(s-t)^j}{j!}
\sum_{i_1,\ldots,i_k=1}^m \bar G_{l_1}^{(i_1)}\ldots
\bar G_{l_k}^{(i_k)}\bar L^j R({\bf x}_t,t)
I_{{l_1\ldots l_k}_{s,t}}^{*(i_1\ldots i_k)},
\end{equation}
$$
R({\bf x}_t,t)+
\sum_{q,k=1}^r \sum_{(\lambda_{k},\ldots,\lambda_1)
\in{\rm E}_{qk}}
\sum_{i_1=\lambda_1}^{m\lambda_1}
\ldots 
\sum_{i_k=\lambda_k}^{m\lambda_k}
D_{\lambda_k}^{(i_k)}\ldots D_{\lambda_1}^{(i_1)}
R({\bf x}_t,t)
{J}_{(\lambda_{k}\ldots \lambda_1)s,t}^{*(i_k\ldots
i_1)}=
$$
\begin{equation}
\label{2024str205}
=R({\bf x}_t,t)+
\sum_{q=1}^r
\sum_{(k,j,l_1,\ldots,l_k) \in {\rm D}_q}\ 
\frac{(s-t)^j}{j!}
\sum_{i_1,\ldots,i_k=1}^m \bar L^j \bar G_{l_1}^{(i_1)}\ldots
\bar G_{l_k}^{(i_k)}R({\bf x}_t,t)
J_{{l_1\ldots l_k}_{s,t}}^{*(i_1\ldots i_k)}.
\end{equation}

Recall the equality (\ref{2024str201})
$$
R({\bf x}_t,t)+
\sum_{q,k=1}^r \sum_{(\lambda_{k},\ldots,\lambda_1)
\in{\rm E}_{qk}}
\sum_{i_1=\lambda_1}^{m\lambda_1}
\ldots 
\sum_{i_k=\lambda_k}^{m\lambda_k}
Q_{\lambda_k}^{(i_k)}\ldots Q_{\lambda_1}^{(i_1)}
R({\bf x}_t,t){J}_{(\lambda_{k}\ldots \lambda_1)s,t}^{(i_k\ldots
i_1)}=
$$
\begin{equation}
\label{2024str206}
=R({\bf x}_t,t)+
\sum_{q,k=1}^r \sum_{(\lambda_{k},\ldots,\lambda_1)
\in{\rm E}_{qk}}
\sum_{i_1=\lambda_1}^{m\lambda_1}
\ldots 
\sum_{i_k=\lambda_k}^{m\lambda_k}
D_{\lambda_k}^{(i_k)}\ldots D_{\lambda_1}^{(i_1)}
R({\bf x}_t,t)
{J}_{(\lambda_{k}\ldots \lambda_1)s,t}^{*(i_k\ldots
i_1)}
\end{equation}

\noindent
w.~p.~1, where $r\ge 2$.

Combining (\ref{2024str204})--(\ref{2024str206}), we obtain
$$
R({\bf x}_t,t)+
\sum_{q,k=1}^r \sum_{(\lambda_{k},\ldots,\lambda_1)
\in{\rm E}_{qk}}
\sum_{i_1=\lambda_1}^{m\lambda_1}
\ldots 
\sum_{i_k=\lambda_k}^{m\lambda_k}
Q_{\lambda_k}^{(i_k)}\ldots Q_{\lambda_1}^{(i_1)}
R({\bf x}_t,t){J}_{(\lambda_{k}\ldots \lambda_1)s,t}^{(i_k\ldots
i_1)}=
$$
$$
=R({\bf x}_t,t)+
\sum_{q=1}^r
\sum_{(k,j,l_1,\ldots,l_k) \in {\rm D}_q}
\frac{(s-t)^j}{j!}
\sum_{i_1,\ldots,i_k=1}^m \bar G_{l_1}^{(i_1)}\ldots
\bar G_{l_k}^{(i_k)}\bar L^j R({\bf x}_t,t)
I_{{l_1\ldots l_k}_{s,t}}^{*(i_1\ldots i_k)}=
$$
\begin{equation}
\label{2024str207}
=R({\bf x}_t,t)+
\sum_{q=1}^r
\sum_{(k,j,l_1,\ldots,l_k) \in {\rm D}_q}\ 
\frac{(s-t)^j}{j!}
\sum_{i_1,\ldots,i_k=1}^m \bar L^j \bar G_{l_1}^{(i_1)}\ldots
\bar G_{l_k}^{(i_k)}R({\bf x}_t,t)
J_{{l_1\ldots l_k}_{s,t}}^{*(i_1\ldots i_k)}
\end{equation}
w.~p.~1, where $r\ge 2$.

The equality (\ref{2024str207}) means that we have the following theorem
(see (\ref{5.6.1rrr}), (\ref{5.6.1rrrh})).

{\bf Theorem 4.5.}\ {\it Let conditions {\rm (i), (ii) (}see
Sect.~{\rm 4.3)}
be satisfied. Then for any $s, t \in [0, T]$ such that $s>t$ 
the following unified 
Taylor--Stratonovich expansions take place w.~p.~{\rm 1}
$$
R({\bf x}_s,s)=R({\bf x}_t,t)+
$$

\vspace{-7mm}
$$
+
\sum_{q=1}^r
\sum_{(k,j,l_1,\ldots,l_k) \in {\rm D}_q}
\frac{(s-t)^j}{j!}
\sum_{i_1,\ldots,i_k=1}^m \bar G_{l_1}^{(i_1)}\ldots
\bar G_{l_k}^{(i_k)}\bar L^j R({\bf x}_t,t)
I_{{l_1\ldots l_k}_{s,t}}^{*(i_1\ldots i_k)}+\left(H_{r+1}\right)_{s,t},
$$

\vspace{-1mm}
$$
R({\bf x}_s,s)=R({\bf x}_t,t)+
$$

\vspace{-7mm}
$$
+
\sum_{q=1}^r
\sum_{(k,j,l_1,\ldots,l_k) \in {\rm D}_q}\ 
\frac{(s-t)^j}{j!}
\sum_{i_1,\ldots,i_k=1}^m \bar L^j \bar G_{l_1}^{(i_1)}\ldots
\bar G_{l_k}^{(i_k)}R({\bf x}_t,t)
J_{{l_1\ldots l_k}_{s,t}}^{*(i_1\ldots i_k)}
+\left(H_{r+1}\right)_{s,t},
$$

\vspace{2mm}
\noindent
where $r\ge 2,$ the reainder term $\left(H_{r+1}\right)_{s,t}$
is definded by the relations {\rm(\ref{5.6.1rrrh})} and {\rm (\ref{5.7.12});}
another notations are the same as in Sect.~{\rm 4.3, 4.6--4.8.}}

\section{Comparison of the Unified Taylor--It\^{o} and
Taylor--Stratonovich Expansions
with the Classical Taylor--It\^{o} and Taylor--Stratonovich Expansions}

Note that the 
truncated 
unified 
Taylor--It\^{o} and Taylor--Stratonovich expansions contain 
the less number of various 
iterated It\^{o} 
and Stratonovich stochastic integrals (moreover, their major part 
will have 
less multiplicity) in comparison with 
the classical Taylor--It\^{o} and
Taylor--Stratonovich expansions \cite{KlPl1}.

It is easy to notice that the stochastic integrals from the sets
(\ref{re11}), (\ref{str11}) are connected by linear relations. 
However, the stochastic integrals from the sets 
(\ref{ll1}), (\ref{ll11}) cannot be connected by linear relations. 
This also holds for the 
stochastic integrals from the sets
(\ref{a111}), (\ref{a112}).
Therefore, we will call the sets 
(\ref{ll1})--(\ref{a112}) as the {\it stochastic 
bases.}

Let us 
call
the numbers ${\rm rank}_{\rm A}(r)$ and
${\rm rank}_{\rm D}(r)$ of various iterated 
It\^{o} and Stratonovich stochastic integrals, which are included in the 
sets (\ref{ll1})--(\ref{a112}) as the 
{\it ranks of stochastic bases} when 
summation in the stochastic expansions is performed using the 
sets 
${\rm A} _ {q}$ ($q=1,\ldots,r)$
and ${\rm D} _ {q}$ ($q=1,\ldots,r$) correspondingly.
Here $r$ is a fixed natural number.

At the beginning, let us analyze several examples
related to the Taylor--It\^{o} expansions (obviously,
the same conclusions will hold for the Taylor--Stratonovich expansions).

Assume that the summation in the unified  
Taylor--It\^{o} expansions is performed using 
the sets ${\rm D} _ {q}$ ($q=1,\ldots,r$).
It is easy to see that the 
truncated 
unified Taylor--It\^{o} expansion (\ref{t100}), 
where the summation is performed using the sets ${\rm D}_q$ when $r=3$  
includes 4 (${\rm rank}_{\rm D}(3)=4$) various iterated It\^{o}
stochastic 
integrals
$$
I_{0_{s,t}}^{(i_1)},\ \  I_{{00}_{s,t}}^{(i_1 i_2)},\ \ 
I_{1_{s,t}}^{(i_1)},\ \ 
I_{{000}_{s,t}}^{(i_1 i_2 i_3)}.
$$

The same truncated 
classical
Taylor--It\^{o} expansion (\ref{5.6.1rrr}) \cite{Zapad-3} contains 
5 various iterated It\^{o} stochastic integrals
$$
{J}_{(1)s,t}^{(i_1)},\ \  {J}_{(11)s,t}^{(i_1 i_2)},\ \ 
{J}_{(10)s,t}^{(i_1 0)},\ \  {J}_{(01)s,t}^{(0 i_1)},\ \ 
{J}_{(111)s,t}^{(i_1 i_2 i_3)}.
$$

For $r=4$ we have
7 (${\rm rank}_{\rm D}(4)=7$) stochastic integrals
$$
I_{0_{s,t}}^{(i_1)},\ \  I_{{00}_{s,t}}^{(i_1 i_2)},\ \ 
I_{1_{s,t}}^{(i_1)},\ \ 
I_{{000}_{s,t}}^{(i_1 i_2 i_3)},\ \ 
I_{{01}_{s,t}}^{(i_1 i_2)},\ \  I_{{10}_{s,t}}^{(i_1 i_2)},\ \ 
I_{{0000}_{s,t}}^{(i_1 i_2 i_3 i_4)}
$$
against 9 stochastic integrals
$$
{J}_{(1)s,t}^{(i_1)},\ \  {J}_{(11)s,t}^{(i_1 i_2)},\ \ 
{J}_{(10)s,t}^{(i_1 0)},\ \  {J}_{(01)s,t}^{(0 i_1)},\ \ 
{J}_{(111)s,t}^{(i_1 i_2 i_3)},\ \ 
{J}_{(101)s,t}^{(i_1 0 i_3)},\ \ 
{J}_{(110)s,t}^{(i_1 i_2 0)},\ \  {J}_{(011)s,t}^{(0 i_1 i_2)},\ \ 
{J}_{(1111)s,t}^{(i_1 i_2 i_3 i_4)}.
$$

For $r=5$ (${\rm rank}_{\rm D}(5)=12$) we get 
12 integrals against 17 integrals
and for 
$r=6$ and $r=7$ we have
20 against 29 and 33 against 50 
correspondingly.

We will obtain the same results when compare the unified 
Taylor--Stra\-to\-no\-vich 
expansions \cite{Kuz} (also see \cite{1}-\cite{12aa}, \cite{arxiv-24}, 
\cite{Kuz-3})
with their classical analogues \cite{Zapad-3}, \cite{KlPl1}
(see previous sections).

Note that the summation with respect to the sets ${\rm D}_q$ is usually used 
while constructing strong numerical methods (built according to 
the mean-square criterion of convergence) for It\^{o} SDEs
\cite{Zapad-1}, \cite{Zapad-3} (also see \cite{12}).
The summation 
with respect to the sets ${\rm A}_q$ is usually used when building 
weak numerical methods (built in accordance with the weak 
criterion of convergence) for It\^{o} SDEs
\cite{Zapad-1}, \cite{Zapad-3}.
For example, ${\rm rank}_{\rm A}(4)=15$ while the total number of various 
iterated It\^{o} stochastic integrals (included in the 
classical Taylor--It\^{o} 
expansion \cite{Zapad-3} when $r=4$) equals to 26.

Let us show that \cite{3}-\cite{12aa}, \cite{arxiv-24}
$$
{\rm rank}_{\rm A}(r)=2^r-1.
$$ 

Let $(l_1,\ldots,l_k)$ be an ordered set such that 
$l_1,\ldots,$ $l_k$ $=$ $0, 1,\ldots$ and
$k=1, 2,\ldots $ Consider $S(k)\stackrel{\sf def}{=}l_1+\ldots+l_k=p$
($p$ is a fixed natural number or zero).
Let $N(k,p)$ be a number of all ordered
combinations
$(l_1,\ldots,l_k)$ such that $l_1,\ldots,l_k=0, 1,\ldots,$
$k=1, 2,\ldots,$ and 
$S(k)=p$. First, let us show that
$$
N(k,p)=C_{p+k-1}^{k-1},
$$
where 
$$
C_n^m=\frac{n!}{m!(n-m)!}
$$ 

\noindent
is a binomial coefficient.

It is not difficult to see that
$$
N(1,p)=1=C_{p+1-1}^{1-1},
$$
$$
N(2,p)=p+1=C_{p+2-1}^{2-1},
$$
$$
N(3,p)=\frac{(p+1)(p+2)}{2}=C_{p+3-1}^{3-1}.
$$ 

Moreover,
$$
N(k+1,p)=\sum\limits_{l=0}^p N(k,l)=\sum_{l=0}^p C_{l+k-1}^{k-1}=C_{p+k}^{k},
$$
where we used the induction assumption and the well known 
property of binomial coefficients.

Then
$$
{\rm rank}_{A}(r)=
$$
$$
=N(1,0)+(N(1,1)+N(2,0))+
(N(1,2)+N(2,1)+N(3,0))+\ldots
$$
$$
\ldots
+(N(1,r-1)+N(2,r-2)+\ldots+N(r,0))=
$$
$$
=C_0^0+(C_1^0+C_1^1)+(C_2^0+C_2^1+C_2^2)+ \ldots
$$
$$
\ldots +(C_{r-1}^0+C_{r-1}^1+C_{r-1}^2+\ldots+C_{r-1}^{r-1})=
$$
$$
=2^0+2^1+2^2+\ldots+2^{r-1}=2^r-1.
$$

Let $n_{{\rm M}}(r)$ be the total number of various iterated stochastic 
integrals included in the 
classical
Taylor--It\^{o} expansion (\ref{5.7.11}) \cite{Zapad-3},
where summation is performed with respect to the set
$$
\bigcup\limits_{k=1}^r{\rm M}_k.
$$

If we exclude from the consideration 
the integrals, which are equal to
$$
\frac{(s-t)^j}{j!},
$$ 
then 
$$
n_{{M}}(r)=
$$
$$
=(2^1-1)+(2^2-1)+(2^3-1)+\ldots+(2^r-1)=
$$
$$
=2(1+2+2^2+\ldots+2^{r-1})-r=2(2^r-1)-r.
$$

It means that
$$
\lim\limits_{r\to\infty}
\frac{n_{{M}}(r)}{{\rm rank}_{A}(r)}=2.
$$

The numbers
$$
{\rm rank}_{\rm A}(r),\ \ \ n_{{\rm M}}(r),\ \ \ 
f(r)=n_{{\rm M}}(r)/{\rm rank}_{\rm A}(r)
$$
for various values $r$ are shown in Table 4.1.

Let us show that \cite{3}-\cite{12aa}, \cite{arxiv-24}
\begin{equation}
\label{gg1aa}
~~~~{\rm rank}_{\rm D}(r)=
\begin{cases}
\sum\limits_{s=0}^{r-1}~~
\sum\limits_{l=s}^{(r-1)/2+[s/2]}~ C_l^s\ \ \  &\hbox{for}\ \ \ 
r=1,\ 3,\ 5,\ldots\cr\cr
\sum\limits_{s=0}^{r-1}~~ \sum\limits_{l=s}^{
r/2-1+[(s+1)/2]}~ C_l^s\ \ \  &\hbox{for}\ \ \ 
r=2,\ 4,\ 6,\ldots
\end{cases},
\end{equation}

\noindent
where $[x]$ is an integer part of a real number $x$
and $C_n^m$ is a binomial coefficient.

For the proof of (\ref{gg1aa}) we rewrite
the condition
$$
k + 2(j + S(k))\le r,
$$ 
where 
$S(k)\stackrel{\sf def}{=}l_1+\ldots+l_k$
($k, j, l_1,\ldots, l_k =0,1,\ldots$)
in the form
$$
j+S(k)\le (r-k)/2
$$ 
and perform the consideration of all possible combinations
with respect to
$k=1,\ldots,r$. Moreover, we take into account the above 
reasoning.

Let us calculate the number
$n_{{E}}(r)$ of all different iterated It\^{o}
stochastic integrals from the classical Taylor--It\^{o}
expansion (\ref{5.6.1rrr}) \cite{Zapad-3} 
if the summation in this expansion is performed
with respect to the set
$$
\bigcup\limits_{q,k=1}^r{E}_{qk}.
$$

The summation condition can be rewritten in this case
in the form 
$$
0\le p+2q\le r,
$$ 
where $q$ 
is a total number of integrations
with respect to time while
$p$ is a total number of integrations
with respect to the Wiener processes in the selected iterated
stochastic integral from the Taylor--It\^{o} expansion 
(\ref{5.6.1rrr}) \cite{Zapad-3}.
At that the multiplicity of the mentioned
stochastic integral equals to $p+q$ and it is not more than
$r.$ Let us rewrite the above condition ($0\le p+2q\le r$)
in the form:
$0\le q\le (r-p)/2$ $\Leftrightarrow$ $0\le q\le [(r-p)/2]$,
where $[x]$ means an integer part of a real number $x$.
Then, performing the consideration of all possible combinations
with respect to $p=1,\ldots,r$ and using the
combinatorial reasoning, we come to 
the formula
\begin{equation}
\label{123u}
n_{{E}}(r)=\sum\limits_{s=1}^r~~ \sum\limits_{l=0}^{[(r-s)/2]}~
C_{[(r-s)/2]+s-l}^s,
\end{equation}

\noindent
where $[x]$ means an integer part of a real number $x$.

\begin{table}
\centering
\caption{Numbers ${\rm rank}_{\rm A}(r),$ $n_{{\rm M}}(r),$ 
$f(r)=n_{{\rm M}}(r)/{\rm rank}_{\rm A}(r)$}
\label{tab:4.1}      
\begin{tabular}{p{1.3cm}p{0.2cm}p{0.9cm}p{0.9cm}p{0.9cm}p{0.9cm}p{1cm}p{1cm}p{1cm}p{1cm}p{1cm}}
\hline\noalign{\smallskip}
$ r$ & 1 & 2& 3& 4 & 5 & 6 & 7 & 8 & 9 & 10\\
\noalign{\smallskip}\hline\noalign{\smallskip}
${\rm rank}_{\rm A}(r)$ &  1 & 3 & 7& 15 & 31 & 63& 127 & 255& 511& 1023\\
${n}_{\rm M}(r)$ &  1 & 4 & 11& 26 & 57 & 120& 247 & 502 & 1013& 2036\\
$f(r)$ & 1& 1.3333& 1.5714& 1.7333 & 1.8387 & 1.9048 & 1.9449 & 1.9686 & 
1.9824 & 1.9902\\
\noalign{\smallskip}\hline\noalign{\smallskip}
\end{tabular}
\end{table}

\begin{table}
\centering
\caption{Numbers ${\rm rank}_{\rm D}(r),$ $n_{{\rm E}}(r),$ 
$g(r)=n_{{\rm E}}(r)/{\rm rank}_{\rm D}(r)$}
\label{tab:4.2}      
\begin{tabular}{p{1.3cm}p{0.2cm}p{0.3cm}p{0.9cm}p{1cm}p{1cm}p{1cm}p{1cm}p{1cm}p{1cm}p{1cm}}
\hline\noalign{\smallskip}
$ r$ & 1 & 2& 3& 4 & 5 & 6 & 7 & 8 & 9 & 10\\
\noalign{\smallskip}\hline\noalign{\smallskip}
${\rm rank}_{\rm D}(r)$ &  1 & 2 & 4& 7 & 12 & 20& 33 & 54& 88 & 143\\
${n}_{\rm E}(r)$ &  1 & 2 & 5& 9 & 17 & 29& 50 & 83 & 138 & 261\\
$g(r)$ & 1& 1& 1.2500& 1.2857 & 1.4167 & 1.4500 & 1.5152 & 1.5370 & 
1.5682 & 1.8252\\
\noalign{\smallskip}\hline\noalign{\smallskip}
\end{tabular}
\vspace{5mm}
\end{table}

The numbers
$$
{\rm rank}_{\rm D}(r),\ \ \ n_{{\rm E}}(r),\ \ \
g(r)=n_{{\rm E}}(r)/{\rm rank}_{\rm D}(r)
$$
for various values $r$ are shown in Table 4.2.

\section{Application of First Form of the Unified 
Taylor--It\^{o} Expansion to the High-Order
Strong Numerical Methods for It\^{o} SDEs}

Let us rewrite (\ref{t100}) for all $s, t\in[0,T]$ such that $s>t$
in the following form
$$
R({\bf x}_s,s)=
R({\bf x}_t,t)+
$$
$$
+\sum_{q=1}^r\sum_{(k,j,l_1,\ldots,l_k)\in
{\rm D}_q}
\frac{(s-t)^j}{j!} \sum\limits_{i_1,\ldots,i_k=1}^m
G_{l_1}^{(i_1)}\ldots G_{l_k}^{(i_k)}
L^j R({\bf x}_t,t)
I^{(i_1\ldots i_k)}_{{l_1\ldots l_k}_{s,t}}+
$$
\begin{equation}
\label{15.001}
~~~~~+{\bf 1}_{\{r=2d-1,
d\in {\bf N}\}}\frac{(s-t)^{(r+1)/2}}{\left(
(r+1)/2\right)!}L^{(r+1)/2}R({\bf x}_t,t)
+\left(\bar H_{r+1}\right)_{s,t}\ \ \ \hbox{w.~p.~1},
\end{equation}

\noindent
where
$$
\left(\bar H_{r+1}\right)_{s,t}= \left(H_{r+1}\right)_{s,t}
-{\bf 1}_{\{r=2d-1,
d\in {\bf N}\}}\frac{(s-t)^{(r+1)/2}}{\left(
(r+1)/2\right)!}L^{(r+1)/2}R({\bf x}_t,t).
$$

\vspace{1mm}

Consider the partition 
$\{\tau_p\}_{p=0}^N$ of the interval
$[0,T]$ such that
$$
0=\tau_0<\tau_1<\ldots<\tau_N=T,\ \ \
\Delta_N=
\max\limits_{0\le j\le N-1}\left|\tau_{j+1}-\tau_j\right|.
$$

\vspace{-2mm}

From (\ref{15.001}) for $s=\tau_{p+1},$
$t=\tau_p$ we obtain the following representation of
explicit one-step strong numerical scheme for the It\^{o} SDE (\ref{1.5.2}),
which is based on first form of the unified Taylor--It\^{o}
expansion
$$
{\bf y}_{p+1}={\bf y}_{p}+
\sum_{q=1}^r\sum_{(k,j,l_1,\ldots,l_k)\in{\rm D}_q}
\frac{(\tau_{p+1}-\tau_p)^j}{j!} \sum\limits_{i_1,\ldots,i_k=1}^m
G_{l_1}^{(i_1)}\ldots G_{l_k}^{(i_k)}
L^j\hspace{0.4mm}{\bf y}_{p}\
{\hat I}^{(i_1\ldots i_k)}_{{l_1\ldots l_k}_{\tau_{p+1},\tau_p}}+
$$
\begin{equation}
\label{15.002}
+{\bf 1}_{\{r=2d-1,
d\in {\bf N}\}}\frac{(\tau_{p+1}-\tau_p)^{(r+1)/2}}{\left(
(r+1)/2\right)!}L^{(r+1)/2}\hspace{0.4mm}{\bf y}_{p},
\end{equation}
where
$$
\hat I^{(i_1\ldots i_k)}_{{l_1\ldots l_k}_{\tau_{p+1},\tau_p}}
$$

\noindent 
is an approximation of the following iterated It\^{o}
stochastic integral 
$$
I_{{l_1\ldots l_k}_{s,t}}^{(i_1\ldots i_k)}=
\int\limits_t^s
(t-t_{k})^{l_{k}}\ldots 
\int\limits_t^{t_2}
(t-t _ {1}) ^ {l_ {1}} d
{\bf w} ^ {(i_ {1})} _ {t_ {1}} \ldots 
d {\bf w} _ {t_ {k}} ^ {(i_ {k})}.
$$

Note that we understand the equality (\ref{15.002}) componentwise
with respect to the components ${\bf y}_p^{(i)}$ of the column
${\bf y}_p.$
Also for simplicity we put 
$\tau_p=p\Delta$,
$\Delta=T/N,$ $T=\tau_N,$ $p=0,1,\ldots,N.$

It is known \cite{Zapad-3} that under the appropriate conditions
the numerical scheme (\ref{15.002}) has strong order of convergence $r/2$
($r\in{\bf N}$).

Let $B_j({\bf x},t)$ is the $j$th column of the matrix
function $B({\bf x},t).$

Below we consider particular cases of the numerical scheme
(\ref{15.002}) for $r=2,3,4,5,$ and $6,$ i.e. 
ex\-plicit one-step strong numerical schemes for the
It\^{o} SDE (\ref{1.5.2})
with the convergence orders $1.0$, $1.5$, $2.0$, $2.5,$ and $3.0$.
At that for simplicity 
we will write ${\bf a},$ $L {\bf a},$ 
$B_i,$ $G_0^{(i)}B_{j},$ $\ldots$
instead of ${\bf a}({\bf y}_p,\tau_p),$ 
$L {\bf a}({\bf y}_p,\tau_p),$ 
$B_i({\bf y}_p,\tau_p),$ 
$G_0^{(i)}B_{j}({\bf y}_p,\tau_p),$ $\ldots$ correspondingly.
Moreover, the operators $L$ and $G_0^{(i)},$ $i=1,\ldots,m$
are determined by the equalities
(\ref{2.3}), (\ref{2.4}).

\vspace{11mm}

\centerline{\bf Scheme with strong order 1.0 (Milstein Scheme)}

\vspace{-1mm}
\begin{equation}
\label{al1}
~~~~~~~{\bf y}_{p+1}={\bf y}_{p}+\sum_{i_{1}=1}^{m}B_{i_{1}}
\hat I_{0_{\tau_{p+1},\tau_p}}^{(i_{1})}+\Delta{\bf a}
+\sum_{i_{1},i_{2}=1}^{m}G_0^{(i_{2})}
B_{i_{1}}\hat I_{00_{\tau_{p+1},\tau_p}}^{(i_{2}i_{1})}.
\end{equation}

\vspace{6mm}

\centerline{\bf Scheme with strong order 1.5}

\vspace{2mm}
$$
{\bf y}_{p+1}={\bf y}_{p}+\sum_{i_{1}=1}^{m}B_{i_{1}}
\hat I_{0_{\tau_{p+1},\tau_p}}^{(i_{1})}+\Delta{\bf a}
+\sum_{i_{1},i_{2}=1}^{m}G_0^{(i_{2})}
B_{i_{1}}\hat I_{00_{\tau_{p+1},\tau_p}}^{(i_{2}i_{1})}+
$$
$$
+
\sum_{i_{1}=1}^{m}\left[G_0^{(i_{1})}{\bf a}\left(
\Delta \hat I_{0_{\tau_{p+1},\tau_p}}^{(i_{1})}+
\hat I_{1_{\tau_{p+1},\tau_p}}^{(i_{1})}\right)
- LB_{i_{1}}\hat I_{1_{\tau_{p+1},\tau_p}}^{(i_{1})}\right]+
$$
\begin{equation}
\label{al2}
+\sum_{i_{1},i_{2},i_{3}=1}^{m} G_0^{(i_{3})}G_0^{(i_{2})}
B_{i_{1}}\hat I_{000_{\tau_{p+1},\tau_p}}^{(i_{3}i_{2}i_{1})}
+
\frac{\Delta^2}{2}L{\bf a}.
\end{equation}

\vspace{10mm}

\centerline{\bf Scheme with strong order 2.0}

\vspace{2mm}

$$
{\bf y}_{p+1}={\bf y}_{p}+\sum_{i_{1}=1}^{m}B_{i_{1}}
\hat I_{0_{\tau_{p+1},\tau_p}}^{(i_{1})}+\Delta{\bf a}
+\sum_{i_{1},i_{2}=1}^{m}G_0^{(i_{2})}
B_{i_{1}}\hat I_{00_{\tau_{p+1},\tau_p}}^{(i_{2}i_{1})}+
$$
$$
+
\sum_{i_{1}=1}^{m}\left[G_0^{(i_{1})}{\bf a}\left(
\Delta \hat I_{0_{\tau_{p+1},\tau_p}}^{(i_{1})}+
\hat I_{1_{\tau_{p+1},\tau_p}}^{(i_{1})}\right)
- LB_{i_{1}}\hat I_{1_{\tau_{p+1},\tau_p}}^{(i_{1})}\right]+
$$
$$
+\sum_{i_{1},i_{2},i_{3}=1}^{m} G_0^{(i_{3})}G_0^{(i_{2})}
B_{i_{1}}\hat I_{000_{\tau_{p+1},\tau_p}}^{(i_{3}i_{2}i_{1})}+
\frac{\Delta^2}{2} L{\bf a}+
$$
$$
+\sum_{i_{1},i_{2}=1}^{m}
\left[G_0^{(i_{2})} LB_{i_{1}}\left(
\hat I_{10_{\tau_{p+1},\tau_p}}^{(i_{2}i_{1})}-
\hat I_{01_{\tau_{p+1},\tau_p}}^{(i_{2}i_{1})}
\right)
-LG_0^{(i_{2})}
B_{i_{1}}\hat I_{10_{\tau_{p+1},\tau_p}}^{(i_{2}i_{1})}
+\right.
$$
$$
\left.+G_0^{(i_{2})}G_0^{(i_{1})}{\bf a}\left(
\hat I_{01_{\tau_{p+1},\tau_p}}
^{(i_{2}i_{1})}+\Delta \hat I_{00_{\tau_{p+1},\tau_p}}^{(i_{2}i_{1})}
\right)\right]+
$$
\begin{equation}
\label{al3}
+
\sum_{i_{1},i_{2},i_{3},i_{4}=1}^{m}G_0^{(i_{4})}G_0^{(i_{3})}G_0^{(i_{2})}
B_{i_{1}}\hat I_{0000_{\tau_{p+1},\tau_p}}^{(i_{4}i_{3}i_{2}i_{1})}.
\end{equation}

\vspace{8mm}

\centerline{\bf Scheme with strong order 2.5}

\vspace{2mm}

$$
{\bf y}_{p+1}={\bf y}_{p}+\sum_{i_{1}=1}^{m}B_{i_{1}}
\hat I_{0_{\tau_{p+1},\tau_p}}^{(i_{1})}+\Delta{\bf a}
+\sum_{i_{1},i_{2}=1}^{m}G_0^{(i_{2})}
B_{i_{1}}\hat I_{00_{\tau_{p+1},\tau_p}}^{(i_{2}i_{1})}+
$$
$$
+
\sum_{i_{1}=1}^{m}\left[G_0^{(i_{1})}{\bf a}\left(
\Delta \hat I_{0_{\tau_{p+1},\tau_p}}^{(i_{1})}+
\hat I_{1_{\tau_{p+1},\tau_p}}^{(i_{1})}\right)
- LB_{i_{1}}\hat I_{1_{\tau_{p+1},\tau_p}}^{(i_{1})}\right]+
$$
$$
+\sum_{i_{1},i_{2},i_{3}=1}^{m} G_0^{(i_{3})}G_0^{(i_{2})}
B_{i_{1}}\hat I_{000_{\tau_{p+1},\tau_p}}^{(i_{3}i_{2}i_{1})}+
\frac{\Delta^2}{2} L{\bf a}+
$$
$$
+\sum_{i_{1},i_{2}=1}^{m}
\left[G_0^{(i_{2})} LB_{i_{1}}\left(
\hat I_{10_{\tau_{p+1},\tau_p}}^{(i_{2}i_{1})}-
\hat I_{01_{\tau_{p+1},\tau_p}}^{(i_{2}i_{1})}
\right)
- LG_0^{(i_{2})}
B_{i_{1}}\hat I_{10_{\tau_{p+1},\tau_p}}^{(i_{2}i_{1})}
+\right.
$$
$$
\left.+G_0^{(i_{2})}G_0^{(i_{1})}{\bf a}\left(
\hat I_{01_{\tau_{p+1},\tau_p}}
^{(i_{2}i_{1})}+\Delta \hat I_{00_{\tau_{p+1},\tau_p}}^{(i_{2}i_{1})}
\right)\right]+
$$
$$
+
\sum_{i_{1},i_{2},i_{3},i_{4}=1}^{m}G_0^{(i_{4})}G_0^{(i_{3})}G_0^{(i_{2})}
B_{i_{1}}\hat I_{0000_{\tau_{p+1},\tau_p}}^{(i_{4}i_{3}i_{2}i_{1})}+
$$
$$
+\sum_{i_{1}=1}^{m}\Biggl[G_0^{(i_{1})} L{\bf a}\left(\frac{1}{2}
\hat I_{2_{\tau_{p+1},\tau_p}}
^{(i_{1})}+\Delta \hat I_{1_{\tau_{p+1},\tau_p}}^{(i_{1})}+
\frac{\Delta^2}{2}\hat I_{0_{\tau_{p+1},\tau_p}}^{(i_{1})}\right)\Biggr.+
$$
$$
+\frac{1}{2} L LB_{i_{1}}\hat I_{2_{\tau_{p+1},\tau_p}}^{(i_{1})}-
 LG_0^{(i_{1})}{\bf a}\Biggl.
\left(\hat I_{2_{\tau_{p+1},\tau_p}}^{(i_{1})}+
\Delta \hat I_{1{\tau_{p+1},\tau_p}}^{(i_{1})}\right)\Biggr]+
$$
$$
+
\sum_{i_{1},i_{2},i_{3}=1}^m\left[
G_0^{(i_{3})} LG_0^{(i_{2})}B_{i_{1}}
\left(\hat I_{100_{\tau_{p+1},\tau_p}}
^{(i_{3}i_{2}i_{1})}-\hat I_{010_{\tau_{p+1},\tau_p}}
^{(i_{3}i_{2}i_{1})}\right)
\right.+
$$
$$
+G_0^{(i_{3})}G_0^{(i_{2})} LB_{i_{1}}\left(
\hat I_{010_{\tau_{p+1},\tau_p}}^{(i_{3}i_{2}i_{1})}-
\hat I_{001_{\tau_{p+1},\tau_p}}^{(i_{3}i_{2}i_{1})}\right)+
$$

\vspace{-1mm}
$$
+G_0^{(i_{3})}G_0^{(i_{2})}G_0^{(i_{1})} {\bf a}
\left(\Delta \hat I_{000_{\tau_{p+1},\tau_p}}^{(i_{3}i_{2}i_{1})}+
\hat I_{001_{\tau_{p+1},\tau_p}}^{(i_{3}i_{2}i_{1})}\right)-
$$

\vspace{-1mm}
$$
\left.- LG_0^{(i_{3})}G_0^{(i_{2})}B_{i_{1}}
\hat I_{100_{\tau_{p+1},\tau_p}}^{(i_{3}i_{2}i_{1})}\right]+
$$
\begin{equation}
\label{al4}
+\sum_{i_{1},i_{2},i_{3},i_{4},i_{5}=1}^m
G_0^{(i_{5})}G_0^{(i_{4})}G_0^{(i_{3})}G_0^{(i_{2})}B_{i_{1}}
\hat I_{00000_{\tau_{p+1},\tau_p}}^{(i_{5}i_{4}i_{3}i_{2}i_{1})}+
\frac{\Delta^3}{6}LL{\bf a}.
\end{equation}

\vspace{10mm}

\centerline{\bf Scheme with strong order 3.0}

\vspace{2mm}

$$
{\bf y}_{p+1}={\bf y}_{p}+\sum_{i_{1}=1}^{m}B_{i_{1}}
\hat I_{0_{\tau_{p+1},\tau_p}}^{(i_{1})}+\Delta{\bf a}
+\sum_{i_{1},i_{2}=1}^{m}G_0^{(i_{2})}
B_{i_{1}}\hat I_{00_{\tau_{p+1},\tau_p}}^{(i_{2}i_{1})}+
$$
$$
+
\sum_{i_{1}=1}^{m}\left[G_0^{(i_{1})}{\bf a}\left(
\Delta \hat I_{0_{\tau_{p+1},\tau_p}}^{(i_{1})}+
\hat I_{1_{\tau_{p+1},\tau_p}}^{(i_{1})}\right)
- LB_{i_{1}}\hat I_{1_{\tau_{p+1},\tau_p}}^{(i_{1})}\right]+
$$
$$
+\sum_{i_{1},i_{2},i_{3}=1}^{m} G_0^{(i_{3})}G_0^{(i_{2})}
B_{i_{1}}\hat I_{000_{\tau_{p+1},\tau_p}}^{(i_{3}i_{2}i_{1})}+
\frac{\Delta^2}{2} L{\bf a}+
$$
$$
+\sum_{i_{1},i_{2}=1}^{m}
\left[G_0^{(i_{2})} LB_{i_{1}}\left(
\hat I_{10_{\tau_{p+1},\tau_p}}^{(i_{2}i_{1})}-
\hat I_{01_{\tau_{p+1},\tau_p}}^{(i_{2}i_{1})}
\right)
- LG_0^{(i_{2})}
B_{i_{1}}\hat I_{10_{\tau_{p+1},\tau_p}}^{(i_{2}i_{1})}
+\right.
$$
$$
\left.+G_0^{(i_{2})}G_0^{(i_{1})}{\bf a}\left(
\hat I_{01_{\tau_{p+1},\tau_p}}
^{(i_{2}i_{1})}+\Delta \hat I_{00_{\tau_{p+1},\tau_p}}^{(i_{2}i_{1})}
\right)\right]+
$$
\begin{equation}
\label{al5}
+
\sum_{i_{1},i_{2},i_{3},i_{4}=1}^{m}G_0^{(i_{4})}G_0^{(i_{3})}G_0^{(i_{2})}
B_{i_{1}}\hat I_{0000_{\tau_{p+1},\tau_p}}^{(i_{4}i_{3}i_{2}i_{1})}+
{\bf q}_{p+1,p}+{\bf r}_{p+1,p},
\end{equation}

\noindent
where
$$
{\bf q}_{p+1,p}=
\sum_{i_{1}=1}^{m}\Biggl[G_0^{(i_{1})} L{\bf a}\left(\frac{1}{2}
\hat I_{2_{\tau_{p+1},\tau_p}}
^{(i_{1})}+\Delta \hat I_{1_{\tau_{p+1},\tau_p}}^{(i_{1})}+
\frac{\Delta^2}{2}\hat I_{0_{\tau_{p+1},\tau_p}}^{(i_{1})}\right)\Biggr.+
$$
$$
+\frac{1}{2} L LB_{i_{1}}\hat I_{2_{\tau_{p+1},\tau_p}}^{(i_{1})}-
LG_0^{(i_{1})}{\bf a}\Biggl.
\left(\hat I_{2_{\tau_{p+1},\tau_p}}^{(i_{1})}+
\Delta \hat I_{1{\tau_{p+1},\tau_p}}^{(i_{1})}\right)\Biggr]+
$$
$$
+
\sum_{i_{1},i_{2},i_{3}=1}^m\left[
G_0^{(i_{3})} LG_0^{(i_{2})}B_{i_{1}}
\left(\hat I_{100_{\tau_{p+1},\tau_p}}
^{(i_{3}i_{2}i_{1})}-\hat I_{010_{\tau_{p+1},\tau_p}}
^{(i_{3}i_{2}i_{1})}\right)
\right.+
$$
$$
+G_0^{(i_{3})}G_0^{(i_{2})} LB_{i_{1}}\left(
\hat I_{010_{\tau_{p+1},\tau_p}}^{(i_{3}i_{2}i_{1})}-
\hat I_{001_{\tau_{p+1},\tau_p}}^{(i_{3}i_{2}i_{1})}\right)+
$$

\vspace{-1mm}
$$
+G_0^{(i_{3})}G_0^{(i_{2})}G_0^{(i_{1})} {\bf a}
\left(\Delta \hat I_{000_{\tau_{p+1},\tau_p}}^{(i_{3}i_{2}i_{1})}+
\hat I_{001_{\tau_{p+1},\tau_p}}^{(i_{3}i_{2}i_{1})}\right)-
$$

\vspace{-1mm}
$$
\left.- LG_0^{(i_{3})}G_0^{(i_{2})}B_{i_{1}}
\hat I_{100_{\tau_{p+1},\tau_p}}^{(i_{3}i_{2}i_{1})}\right]+
$$
$$
+\sum_{i_{1},i_{2},i_{3},i_{4},i_{5}=1}^m
G_0^{(i_{5})}G_0^{(i_{4})}G_0^{(i_{3})}G_0^{(i_{2})}B_{i_{1}}
\hat I_{00000_{\tau_{p+1},\tau_p}}^{(i_{5}i_{4}i_{3}i_{2}i_{1})}+
$$
$$
+
\frac{\Delta^3}{6}LL {\bf a},
$$

\noindent
and
$$
{\bf r}_{p+1,p}=\sum_{i_{1},i_{2}=1}^{m}
\Biggl[G_0^{(i_{2})}G_0^{(i_{1})} L {\bf a}\Biggl(
\frac{1}{2}\hat I_{02_{\tau_{p+1},\tau_p}}^{(i_{2}i_{1})}
+
\Delta \hat I_{01_{\tau_{p+1},\tau_p}}^{(i_{2}i_{1})}
+
\frac{\Delta^2}{2}
\hat I_{00_{\tau_{p+1},\tau_p}}^{(i_{2}i_{1})}\Biggr)+\Biggr.
$$
$$
+
\frac{1}{2} L LG_0^{(i_{2})}B_{i_{1}}
\hat I_{20_{\tau_{p+1},\tau_p}}^{(i_{2}i_{1})}+
$$

\vspace{-5mm}
$$
+G_0^{(i_{2})} LG_0^{(i_{1})} {\bf a}\left(
\hat I_{11_{\tau_{p+1},\tau_p}}
^{(i_{2}i_{1})}-\hat I_{02_{\tau_{p+1},\tau_p}}^{(i_{2}i_{1})}+
\Delta\left(\hat I_{10_{\tau_{p+1},\tau_p}}
^{(i_{2}i_{1})}-\hat I_{01_{\tau_{p+1},\tau_p}}^{(i_{2}i_{1})}
\right)\right)+
$$

\vspace{-2mm}
$$
+ LG_0^{(i_{2})} LB_{i_1}\left(
\hat I_{11_{\tau_{p+1},\tau_p}}
^{(i_{2}i_{1})}-\hat I_{20_{\tau_{p+1},\tau_p}}^{(i_{2}i_{1})}\right)+
$$
$$
+G_0^{(i_{2})} L LB_{i_1}\Biggl(
\frac{1}{2}\hat I_{02_{\tau_{p+1},\tau_p}}^{(i_{2}i_{1})}+
\frac{1}{2}\hat I_{20_{\tau_{p+1},\tau_p}}^{(i_{2}i_{1})}-
\hat I_{11_{\tau_{p+1},\tau_p}}^{(i_{2}i_{1})}\Biggr)-
$$
$$
\Biggl.- LG_0^{(i_{2})}G_0^{(i_{1})}{\bf a}\left(
\Delta \hat I_{10_{\tau_{p+1},\tau_p}}
^{(i_{2}i_{1})}+\hat I_{11_{\tau_{p+1},\tau_p}}^{(i_{2}i_{1})}\right)
\Biggr]+
$$
$$
+
\sum_{i_{1},i_2,i_3,i_{4}=1}^m\Biggl[
G_0^{(i_{4})}G_0^{(i_{3})}G_0^{(i_{2})}G_0^{(i_{1})}{\bf a}
\left(\Delta \hat I_{0000_{\tau_{p+1},\tau_p}}
^{(i_4i_{3}i_{2}i_{1})}+\hat I_{0001_{\tau_{p+1},\tau_p}}
^{(i_4i_{3}i_{2}i_{1})}\right)
+\Biggr.
$$
$$
+G_0^{(i_{4})}G_0^{(i_{3})} LG_0^{(i_{2})}B_{i_1}
\left(\hat I_{0100_{\tau_{p+1},\tau_p}}
^{(i_4i_{3}i_{2}i_{1})}-\hat I_{0010_{\tau_{p+1},\tau_p}}
^{(i_4i_{3}i_{2}i_{1})}\right)-
$$

\vspace{-2mm}
$$
- LG_0^{(i_{4})}G_0^{(i_{3})}G_0^{(i_{2})}B_{i_1}
\hat I_{1000_{\tau_{p+1},\tau_p}}
^{(i_4i_{3}i_{2}i_{1})}+
$$

\vspace{-2mm}
$$
+G_0^{(i_{4})} LG_0^{(i_{3})}G_0^{(i_{2})}B_{i_1}
\left(\hat I_{1000_{\tau_{p+1},\tau_p}}
^{(i_4i_{3}i_{2}i_{1})}-\hat I_{0100_{\tau_{p+1},\tau_p}}
^{(i_4i_{3}i_{2}i_{1})}\right)+
$$
$$
\Biggl.+G_0^{(i_{4})}G_0^{(i_{3})}G_0^{(i_{2})}LB_{i_1}
\left(\hat I_{0010_{\tau_{p+1},\tau_p}}
^{(i_4i_{3}i_{2}i_{1})}-\hat I_{0001_{\tau_{p+1},\tau_p}}
^{(i_4i_{3}i_{2}i_{1})}\right)\Biggr]+
$$
$$
+\sum_{i_{1},i_2,i_3,i_4,i_5,i_{6}=1}^m
G_0^{(i_{6})}G_0^{(i_{5})}
G_0^{(i_{4})}G_0^{(i_{3})}G_0^{(i_{2})}B_{i_{1}}
\hat I_{000000_{\tau_{p+1},\tau_p}}^{(i_6i_{5}i_{4}i_{3}i_{2}i_{1})}.
$$

\vspace{3mm}

It is well known \cite{Zapad-3} that under the standard conditions
the numerical schemes (\ref{al1})--(\ref{al5}) 
have strong orders of convergence 1.0, 1.5, 2.0, 2.5, and 3.0
correspondingly.
Among these conditions we consider only the condition
for approximations of iterated It\^{o} stochastic 
integrals from the numerical
schemes (\ref{al1})--(\ref{al5}) \cite{Zapad-3} (also see \cite{12})
\begin{equation}
\label{agentww3}
{\sf M}\left\{\Biggl(I_{{l_{1}\ldots l_{k}}_{\tau_{p+1},\tau_p}}
^{(i_{1}\ldots i_{k})} 
-\hat I_{{l_{1}\ldots l_{k}}_{\tau_{p+1},\tau_p}}^{(i_{1}\ldots i_{k})}
\Biggr)^2\right\}\le C\Delta^{r+1},
\end{equation}
where constant $C$ is independent of $\Delta$ and
$r/2$ are strong orders of convergence for the numerical schemes
(\ref{al1})--(\ref{al5}), i.e. $r/2=1.0, 1.5,$ $2.0, 2.5,$ and $3.0.$

As we mentioned above, the numerical schemes (\ref{al1})--(\ref{al5})
are unrealizable in practice without 
procedures for the numerical simulation 
of iterated It\^{o} stochastic integrals
from (\ref{15.001}).

In Chapter 5 
we give an extensive material on the mean-square 
approximation of specific iterated
It\^{o} stochastic integrals
from the numerical schemes (\ref{al1})--(\ref{al5}).
The mentioned material based on the results of Chapter 1.

\section{Application of First Form of the Unified 
Taylor--Stratonovich Expansion to the High-Order
Strong Numerical Methods for It\^{o} SDEs}

Let us rewrite (\ref{t100x}) for all $s, t\in[0,T]$ such that $s>t$
in the following from
$$
R({\bf x}_s,s)=
R({\bf x}_t,t)+
$$
$$
+\sum_{q=1}^r\sum_{(k,j,l_1,\ldots,l_k)\in
{\rm D}_q}
\frac{(s-t)^j}{j!} \sum\limits_{i_1,\ldots,i_k=1}^m
\bar G_{l_1}^{(i_1)}\ldots \bar G_{l_k}^{(i_k)}
\bar L^j R({\bf x}_t,t)
I^{*(i_1\ldots i_k)}_{{l_1\ldots l_k}_{s,t}}+
$$
\begin{equation}
\label{15.001x}
~~~~~+{\bf 1}_{\{r=2d-1,
d\in {\bf N}\}}\frac{(s-t)^{(r+1)/2}}{\left(
(r+1)/2\right)!}L^{(r+1)/2}R({\bf x}_t,t)
+\left(\bar H_{r+1}\right)_{s,t}\ \ \ \hbox{w.~p.~1},
\end{equation}

\noindent
where
$$
\left(\bar H_{r+1}\right)_{s,t}= \left(H_{r+1}\right)_{s,t}
-{\bf 1}_{\{r=2d-1,
d\in {\bf N}\}}\frac{(s-t)^{(r+1)/2}}{\left(
(r+1)/2\right)!}L^{(r+1)/2}R({\bf x}_t,t).
$$

\vspace{1mm}

Consider the partition 
$\{\tau_p\}_{p=0}^N$ of the interval
$[0,T]$ such that
$$
0=\tau_0<\tau_1<\ldots<\tau_N=T,\ \ \
\Delta_N=
\max\limits_{0\le j\le N-1}\left|\tau_{j+1}-\tau_j\right|.
$$

From (\ref{15.001x}) for $s=\tau_{p+1},$
$t=\tau_p$ we obtain the following representation of
explicit one-step strong numerical scheme for the It\^{o} SDE (\ref{1.5.2}),
which is based on first form of the unified Taylor--Stratonovich
expansion
$$
{\bf y}_{p+1}={\bf y}_{p}+
\sum_{q=1}^r\sum_{(k,j,l_1,\ldots,l_k)\in{\rm D}_q}
\frac{(\tau_{p+1}-\tau_p)^j}{j!} \sum\limits_{i_1,\ldots,i_k=1}^m
\bar G_{l_1}^{(i_1)}\ldots \bar G_{l_k}^{(i_k)}
\bar L^j\hspace{0.4mm}{\bf y}_{p}\
{\hat I}^{*(i_1\ldots i_k)}_{{l_1\ldots l_k}_{\tau_{p+1},\tau_p}}+
$$
\begin{equation}
\label{15.002x}
+{\bf 1}_{\{r=2d-1,
d\in {\bf N}\}}\frac{(\tau_{p+1}-\tau_p)^{(r+1)/2}}{\left(
(r+1)/2\right)!}L^{(r+1)/2}\hspace{0.4mm}{\bf y}_{p},
\end{equation}

\noindent
where
$$
\hat I^{*(i_1\ldots i_k)}_{{l_1\ldots l_k}_{\tau_{p+1},\tau_p}}
$$

\noindent 
is an approximation of the following iterated Stratonovich
stochastic integral 
$$
I_{{l_1\ldots l_k}_{s,t}}^{*(i_1\ldots i_k)}=
  {\int\limits_t^{*}}^s
(t-t_{k})^{l_{k}}\ldots 
{\int\limits_t^{*}}^{t_2}
(t-t _ {1}) ^ {l_ {1}} d
{\bf w} ^ {(i_ {1})} _ {t_ {1}} \ldots 
d {\bf w} _ {t_ {k}} ^ {(i_ {k})}.
$$

\vspace{1mm}

Note that we understand the equality (\ref{15.002x}) componentwise
with respect to the components ${\bf y}_p^{(i)}$ of the column
${\bf y}_p.$
Also for simplicity we put 
$\tau_p=p\Delta$,
$\Delta=T/N,$ $T=\tau_N,$ $p=0,1,\ldots,N.$

It is known \cite{Zapad-3} that under the appropriate conditions
the numerical scheme (\ref{15.002x}) has strong order of convergence $r/2$
($r\in{\bf N}$).

Denote
$$
\bar{\bf a}({\bf x},t)={\bf a}({\bf x},t)-
\frac{1}{2}\sum\limits_{j=1}^m G_0^{(j)}B_j({\bf x},t),
$$

\noindent
where $B_j({\bf x},t)$ is the $j$th column of the matrix
function $B({\bf x},t).$

It is not difficult to show that (see (\ref{2.4a}))
\begin{equation}
\label{2.4xxx}
~~~{\bar L}R({\bf x},t)=
\frac{\partial R}{\partial t}({\bf x},t)+
\sum\limits_{j=1}^n \bar {\bf a}^{(j)}({\bf x},t)
\frac{\partial R}{\partial {\bf x}^{(j)}}({\bf x},t),
\end{equation}
where $\bar {\bf a}^{(j)}({\bf x},t)$ is the $j$th component of the vector
function $\bar {\bf a}({\bf x},t).$

Below we consider particular cases of the numerical scheme
(\ref{15.002x}) for $r=2,3,4,5,$ and $6,$ i.e. 
explicit one-step strong numerical schemes for the It\^{o} SDE (\ref{1.5.2})
with the convergence orders $1.0,$ $1.5,$ $2.0,$ $2.5,$ and $3.0$.
At that, for simplicity 
we will write $\bar{\bf a},$ $\bar L\bar {\bf a},$ $L{\bf a},$
$B_i,$ $G_0^{(i)}B_{j},$ $\ldots$
instead of $\bar{\bf a}({\bf y}_p,\tau_p),$ 
$\bar L \bar {\bf a}({\bf y}_p,\tau_p),$ 
$L{\bf a}({\bf y}_p,\tau_p),$  $B_i({\bf y}_p,\tau_p),$ 
$G_0^{(i)}B_{j}({\bf y}_p,\tau_p),$ $\ldots$ correspondingly.
Moreover, the operators $\bar L$ and $G_0^{(i)},$ $i=1,\ldots,m$
are determined by the equalities
(\ref{2.3}), (\ref{2.4}), and (\ref{2.4xxx}).

\vspace{6mm}

\centerline{\bf Scheme with strong order 1.0}

\vspace{-1mm}
\begin{equation}
\label{al1x}
~~~~~~~{\bf y}_{p+1}={\bf y}_{p}+\sum_{i_{1}=1}^{m}B_{i_{1}}
\hat I_{0_{\tau_{p+1},\tau_p}}^{*(i_{1})}+\Delta\bar{\bf a}
+\sum_{i_{1},i_{2}=1}^{m}G_0^{(i_{2})}
B_{i_{1}}\hat I_{00_{\tau_{p+1},\tau_p}}^{*(i_{2}i_{1})}.
\end{equation}

\vspace{6mm}

\centerline{\bf Scheme with strong order 1.5}

\vspace{2mm}

$$
{\bf y}_{p+1}={\bf y}_{p}+\sum_{i_{1}=1}^{m}B_{i_{1}}
\hat I_{0_{\tau_{p+1},\tau_p}}^{*(i_{1})}+\Delta\bar{\bf a}
+\sum_{i_{1},i_{2}=1}^{m}G_0^{(i_{2})}
B_{i_{1}}\hat I_{00_{\tau_{p+1},\tau_p}}^{*(i_{2}i_{1})}+
$$
$$
+
\sum_{i_{1}=1}^{m}\left[G_0^{(i_{1})}\bar{\bf a}\left(
\Delta \hat I_{0_{\tau_{p+1},\tau_p}}^{*(i_{1})}+
\hat I_{1_{\tau_{p+1},\tau_p}}^{*(i_{1})}\right)
-\bar LB_{i_{1}}\hat I_{1_{\tau_{p+1},\tau_p}}^{*(i_{1})}\right]+
$$
\begin{equation}
\label{al2x}
+\sum_{i_{1},i_{2},i_{3}=1}^{m} G_0^{(i_{3})}G_0^{(i_{2})}
B_{i_{1}}\hat I_{000_{\tau_{p+1},\tau_p}}^{*(i_{3}i_{2}i_{1})}
+\frac{\Delta^2}{2}L{\bf a}.
\end{equation}

\vspace{10mm}

\centerline{\bf Scheme with strong order 2.0}

\vspace{3mm}

$$
{\bf y}_{p+1}={\bf y}_{p}+\sum_{i_{1}=1}^{m}B_{i_{1}}
\hat I_{0_{\tau_{p+1},\tau_p}}^{*(i_{1})}+\Delta\bar{\bf a}
+\sum_{i_{1},i_{2}=1}^{m}G_0^{(i_{2})}
B_{i_{1}}\hat I_{00_{\tau_{p+1},\tau_p}}^{*(i_{2}i_{1})}+
$$
$$
+
\sum_{i_{1}=1}^{m}\left[G_0^{(i_{1})}\bar{\bf a}\left(
\Delta \hat I_{0_{\tau_{p+1},\tau_p}}^{*(i_{1})}+
\hat I_{1_{\tau_{p+1},\tau_p}}^{*(i_{1})}\right)
-\bar LB_{i_{1}}\hat I_{1_{\tau_{p+1},\tau_p}}^{*(i_{1})}\right]+
$$
$$
+\sum_{i_{1},i_{2},i_{3}=1}^{m} G_0^{(i_{3})}G_0^{(i_{2})}
B_{i_{1}}\hat I_{000_{\tau_{p+1},\tau_p}}^{*(i_{3}i_{2}i_{1})}+
\frac{\Delta^2}{2}\bar L\bar{\bf a}+
$$
$$
+\sum_{i_{1},i_{2}=1}^{m}
\left[G_0^{(i_{2})}\bar LB_{i_{1}}\left(
\hat I_{10_{\tau_{p+1},\tau_p}}^{*(i_{2}i_{1})}-
\hat I_{01_{\tau_{p+1},\tau_p}}^{*(i_{2}i_{1})}
\right)
-\bar LG_0^{(i_{2})}
B_{i_{1}}\hat I_{10_{\tau_{p+1},\tau_p}}^{*(i_{2}i_{1})}
+\right.
$$
$$
\left.+G_0^{(i_{2})}G_0^{(i_{1})}\bar{\bf a}\left(
\hat I_{01_{\tau_{p+1},\tau_p}}
^{*(i_{2}i_{1})}+\Delta \hat I_{00_{\tau_{p+1},\tau_p}}^{*(i_{2}i_{1})}
\right)\right]+
$$
\begin{equation}
\label{al3x}
+
\sum_{i_{1},i_{2},i_{3},i_{4}=1}^{m}G_0^{(i_{4})}G_0^{(i_{3})}G_0^{(i_{2})}
B_{i_{1}}\hat I_{0000_{\tau_{p+1},\tau_p}}^{*(i_{4}i_{3}i_{2}i_{1})}.
\end{equation}

\vspace{10mm}

\centerline{\bf Scheme with strong order 2.5}

\vspace{3mm}

$$
{\bf y}_{p+1}={\bf y}_{p}+\sum_{i_{1}=1}^{m}B_{i_{1}}
\hat I_{0_{\tau_{p+1},\tau_p}}^{*(i_{1})}+\Delta\bar{\bf a}
+\sum_{i_{1},i_{2}=1}^{m}G_0^{(i_{2})}
B_{i_{1}}\hat I_{00_{\tau_{p+1},\tau_p}}^{*(i_{2}i_{1})}+
$$
$$
+
\sum_{i_{1}=1}^{m}\left[G_0^{(i_{1})}\bar{\bf a}\left(
\Delta \hat I_{0_{\tau_{p+1},\tau_p}}^{*(i_{1})}+
\hat I_{1_{\tau_{p+1},\tau_p}}^{*(i_{1})}\right)
-\bar LB_{i_{1}}\hat I_{1_{\tau_{p+1},\tau_p}}^{*(i_{1})}\right]+
$$
$$
+\sum_{i_{1},i_{2},i_{3}=1}^{m} G_0^{(i_{3})}G_0^{(i_{2})}
B_{i_{1}}\hat I_{000_{\tau_{p+1},\tau_p}}^{*(i_{3}i_{2}i_{1})}+
\frac{\Delta^2}{2}\bar L\bar{\bf a}+
$$
$$
+\sum_{i_{1},i_{2}=1}^{m}
\left[G_0^{(i_{2})}\bar LB_{i_{1}}\left(
\hat I_{10_{\tau_{p+1},\tau_p}}^{*(i_{2}i_{1})}-
\hat I_{01_{\tau_{p+1},\tau_p}}^{*(i_{2}i_{1})}
\right)
-\bar LG_0^{(i_{2})}
B_{i_{1}}\hat I_{10_{\tau_{p+1},\tau_p}}^{*(i_{2}i_{1})}
+\right.
$$
$$
\left.+G_0^{(i_{2})}G_0^{(i_{1})}\bar{\bf a}\left(
\hat I_{01_{\tau_{p+1},\tau_p}}
^{*(i_{2}i_{1})}+\Delta \hat I_{00_{\tau_{p+1},\tau_p}}^{*(i_{2}i_{1})}
\right)\right]+
$$
$$
+
\sum_{i_{1},i_{2},i_{3},i_{4}=1}^{m}G_0^{(i_{4})}G_0^{(i_{3})}G_0^{(i_{2})}
B_{i_{1}}\hat I_{0000_{\tau_{p+1},\tau_p}}^{*(i_{4}i_{3}i_{2}i_{1})}+
$$
$$
+\sum_{i_{1}=1}^{m}\Biggl[G_0^{(i_{1})}\bar L\bar{\bf a}\left(\frac{1}{2}
\hat I_{2_{\tau_{p+1},\tau_p}}
^{*(i_{1})}+\Delta \hat I_{1_{\tau_{p+1},\tau_p}}^{*(i_{1})}+
\frac{\Delta^2}{2}\hat I_{0_{\tau_{p+1},\tau_p}}^{*(i_{1})}\right)\Biggr.+
$$
$$
+\frac{1}{2}\bar L\bar LB_{i_{1}}\hat I_{2_{\tau_{p+1},\tau_p}}^{*(i_{1})}-
\bar LG_0^{(i_{1})}\bar{\bf a}\Biggl.
\left(\hat I_{2_{\tau_{p+1},\tau_p}}^{*(i_{1})}+
\Delta \hat I_{1{\tau_{p+1},\tau_p}}^{*(i_{1})}\right)\Biggr]+
$$
$$
+
\sum_{i_{1},i_{2},i_{3}=1}^m\left[
G_0^{(i_{3})}\bar LG_0^{(i_{2})}B_{i_{1}}
\left(\hat I_{100_{\tau_{p+1},\tau_p}}
^{*(i_{3}i_{2}i_{1})}-\hat I_{010_{\tau_{p+1},\tau_p}}
^{*(i_{3}i_{2}i_{1})}\right)
\right.+
$$
$$
+G_0^{(i_{3})}G_0^{(i_{2})}\bar LB_{i_{1}}\left(
\hat I_{010_{\tau_{p+1},\tau_p}}^{*(i_{3}i_{2}i_{1})}-
\hat I_{001_{\tau_{p+1},\tau_p}}^{*(i_{3}i_{2}i_{1})}\right)+
$$

\vspace{-1mm}
$$
+G_0^{(i_{3})}G_0^{(i_{2})}G_0^{(i_{1})}\bar {\bf a}
\left(\Delta \hat I_{000_{\tau_{p+1},\tau_p}}^{*(i_{3}i_{2}i_{1})}+
\hat I_{001_{\tau_{p+1},\tau_p}}^{*(i_{3}i_{2}i_{1})}\right)-
$$

\vspace{-1mm}
$$
\left.-\bar LG_0^{(i_{3})}G_0^{(i_{2})}B_{i_{1}}
\hat I_{100_{\tau_{p+1},\tau_p}}^{*(i_{3}i_{2}i_{1})}\right]+
$$
\begin{equation}
\label{al4x}
+\sum_{i_{1},i_{2},i_{3},i_{4},i_{5}=1}^m
G_0^{(i_{5})}G_0^{(i_{4})}G_0^{(i_{3})}G_0^{(i_{2})}B_{i_{1}}
\hat I_{00000_{\tau_{p+1},\tau_p}}^{*(i_{5}i_{4}i_{3}i_{2}i_{1})}
+
\frac{\Delta^3}{6}LL{\bf a}.
\end{equation}

\vspace{10mm}

\centerline{\bf Scheme with strong order 3.0}

\vspace{3mm}

$$
{\bf y}_{p+1}={\bf y}_{p}+\sum_{i_{1}=1}^{m}B_{i_{1}}
\hat I_{0_{\tau_{p+1},\tau_p}}^{*(i_{1})}+\Delta\bar{\bf a}
+\sum_{i_{1},i_{2}=1}^{m}G_0^{(i_{2})}
B_{i_{1}}\hat I_{00_{\tau_{p+1},\tau_p}}^{*(i_{2}i_{1})}+
$$
$$
+
\sum_{i_{1}=1}^{m}\left[G_0^{(i_{1})}\bar{\bf a}\left(
\Delta \hat I_{0_{\tau_{p+1},\tau_p}}^{*(i_{1})}+
\hat I_{1_{\tau_{p+1},\tau_p}}^{*(i_{1})}\right)
-\bar LB_{i_{1}}\hat I_{1_{\tau_{p+1},\tau_p}}^{*(i_{1})}\right]+
$$
$$
+\sum_{i_{1},i_{2},i_{3}=1}^{m} G_0^{(i_{3})}G_0^{(i_{2})}
B_{i_{1}}\hat I_{000_{\tau_{p+1},\tau_p}}^{*(i_{3}i_{2}i_{1})}+
\frac{\Delta^2}{2}\bar L\bar{\bf a}+
$$
$$
+\sum_{i_{1},i_{2}=1}^{m}
\left[G_0^{(i_{2})}\bar LB_{i_{1}}\left(
\hat I_{10_{\tau_{p+1},\tau_p}}^{*(i_{2}i_{1})}-
\hat I_{01_{\tau_{p+1},\tau_p}}^{*(i_{2}i_{1})}
\right)
-\bar LG_0^{(i_{2})}
B_{i_{1}}\hat I_{10_{\tau_{p+1},\tau_p}}^{*(i_{2}i_{1})}
+\right.
$$
$$
\left.+G_0^{(i_{2})}G_0^{(i_{1})}\bar{\bf a}\left(
\hat I_{01_{\tau_{p+1},\tau_p}}
^{*(i_{2}i_{1})}+\Delta \hat I_{00_{\tau_{p+1},\tau_p}}^{*(i_{2}i_{1})}
\right)\right]+
$$
\begin{equation}
\label{al5x}
+
\sum_{i_{1},i_{2},i_{3},i_{4}=1}^{m}G_0^{(i_{4})}G_0^{(i_{3})}G_0^{(i_{2})}
B_{i_{1}}\hat I_{0000_{\tau_{p+1},\tau_p}}^{*(i_{4}i_{3}i_{2}i_{1})}+
{\bf q}_{p+1,p}+{\bf r}_{p+1,p},
\end{equation}

\vspace{2mm}
\noindent
where
$$
{\bf q}_{p+1,p}=
\sum_{i_{1}=1}^{m}\Biggl[G_0^{(i_{1})}\bar L\bar{\bf a}\left(\frac{1}{2}
\hat I_{2_{\tau_{p+1},\tau_p}}
^{*(i_{1})}+\Delta \hat I_{1_{\tau_{p+1},\tau_p}}^{*(i_{1})}+
\frac{\Delta^2}{2}\hat I_{0_{\tau_{p+1},\tau_p}}^{*(i_{1})}\right)\Biggr.+
$$
$$
+\frac{1}{2}\bar L\bar LB_{i_{1}}\hat I_{2_{\tau_{p+1},\tau_p}}^{*(i_{1})}-
\bar LG_0^{(i_{1})}\bar{\bf a}\Biggl.
\left(\hat I_{2_{\tau_{p+1},\tau_p}}^{*(i_{1})}+
\Delta \hat I_{1{\tau_{p+1},\tau_p}}^{*(i_{1})}\right)\Biggr]+
$$
$$
+
\sum_{i_{1},i_{2},i_{3}=1}^m\left[
G_0^{(i_{3})}\bar LG_0^{(i_{2})}B_{i_{1}}
\left(\hat I_{100_{\tau_{p+1},\tau_p}}
^{*(i_{3}i_{2}i_{1})}-\hat I_{010_{\tau_{p+1},\tau_p}}
^{*(i_{3}i_{2}i_{1})}\right)
\right.+
$$
$$
+G_0^{(i_{3})}G_0^{(i_{2})}\bar LB_{i_{1}}\left(
\hat I_{010_{\tau_{p+1},\tau_p}}^{*(i_{3}i_{2}i_{1})}-
\hat I_{001_{\tau_{p+1},\tau_p}}^{*(i_{3}i_{2}i_{1})}\right)+
$$

\vspace{-1mm}
$$
+G_0^{(i_{3})}G_0^{(i_{2})}G_0^{(i_{1})}\bar {\bf a}
\left(\Delta \hat I_{000_{\tau_{p+1},\tau_p}}^{*(i_{3}i_{2}i_{1})}+
\hat I_{001_{\tau_{p+1},\tau_p}}^{*(i_{3}i_{2}i_{1})}\right)-
$$

\vspace{-1mm}
$$
\left.-\bar LG_0^{(i_{3})}G_0^{(i_{2})}B_{i_{1}}
\hat I_{100_{\tau_{p+1},\tau_p}}^{*(i_{3}i_{2}i_{1})}\right]+
$$
$$
+\sum_{i_{1},i_{2},i_{3},i_{4},i_{5}=1}^m
G_0^{(i_{5})}G_0^{(i_{4})}G_0^{(i_{3})}G_0^{(i_{2})}B_{i_{1}}
\hat I_{00000_{\tau_{p+1},\tau_p}}^{*(i_{5}i_{4}i_{3}i_{2}i_{1})}+
$$
$$
+
\frac{\Delta^3}{6}\bar L\bar L\bar {\bf a},
$$

\vspace{2mm}
\noindent
and
$$
{\bf r}_{p+1,p}=\sum_{i_{1},i_{2}=1}^{m}
\Biggl[G_0^{(i_{2})}G_0^{(i_{1})}\bar L\bar {\bf a}\Biggl(
\frac{1}{2}\hat I_{02_{\tau_{p+1},\tau_p}}^{*(i_{2}i_{1})}
+
\Delta \hat I_{01_{\tau_{p+1},\tau_p}}^{*(i_{2}i_{1})}
+
\frac{\Delta^2}{2}
\hat I_{00_{\tau_{p+1},\tau_p}}^{*(i_{2}i_{1})}\Biggr)+\Biggr.
$$
$$
+
\frac{1}{2}\bar L\bar LG_0^{(i_{2})}B_{i_{1}}
\hat I_{20_{\tau_{p+1},\tau_p}}^{*(i_{2}i_{1})}+
$$

\vspace{-6mm}
$$
+G_0^{(i_{2})}\bar LG_0^{(i_{1})}\bar {\bf a}\left(
\hat I_{11_{\tau_{p+1},\tau_p}}
^{*(i_{2}i_{1})}-\hat I_{02_{\tau_{p+1},\tau_p}}^{*(i_{2}i_{1})}+
\Delta\left(\hat I_{10_{\tau_{p+1},\tau_p}}
^{*(i_{2}i_{1})}-\hat I_{01_{\tau_{p+1},\tau_p}}^{*(i_{2}i_{1})}
\right)\right)+
$$

\vspace{-3mm}
$$
+\bar LG_0^{(i_{2})}\bar LB_{i_1}\left(
\hat I_{11_{\tau_{p+1},\tau_p}}
^{*(i_{2}i_{1})}-\hat I_{20_{\tau_{p+1},\tau_p}}^{*(i_{2}i_{1})}\right)+
$$
$$
+G_0^{(i_{2})}\bar L\bar LB_{i_1}\Biggl(
\frac{1}{2}\hat I_{02_{\tau_{p+1},\tau_p}}^{*(i_{2}i_{1})}+
\frac{1}{2}\hat I_{20_{\tau_{p+1},\tau_p}}^{*(i_{2}i_{1})}-
\hat I_{11_{\tau_{p+1},\tau_p}}^{*(i_{2}i_{1})}\Biggr)-
$$
$$
\Biggl.-\bar LG_0^{(i_{2})}G_0^{(i_{1})}\bar{\bf a}\left(
\Delta \hat I_{10_{\tau_{p+1},\tau_p}}
^{*(i_{2}i_{1})}+\hat I_{11_{\tau_{p+1},\tau_p}}^{*(i_{2}i_{1})}\right)
\Biggr]+
$$
$$
+
\sum_{i_{1},i_2,i_3,i_{4}=1}^m\Biggl[
G_0^{(i_{4})}G_0^{(i_{3})}G_0^{(i_{2})}G_0^{(i_{1})}\bar{\bf a}
\left(\Delta \hat I_{0000_{\tau_{p+1},\tau_p}}
^{*(i_4i_{3}i_{2}i_{1})}+\hat I_{0001_{\tau_{p+1},\tau_p}}
^{*(i_4i_{3}i_{2}i_{1})}\right)
+\Biggr.
$$
$$
+G_0^{(i_{4})}G_0^{(i_{3})}\bar LG_0^{(i_{2})}B_{i_1}
\left(\hat I_{0100_{\tau_{p+1},\tau_p}}
^{*(i_4i_{3}i_{2}i_{1})}-\hat I_{0010_{\tau_{p+1},\tau_p}}
^{*(i_4i_{3}i_{2}i_{1})}\right)-
$$

\vspace{-2mm}
$$
-\bar LG_0^{(i_{4})}G_0^{(i_{3})}G_0^{(i_{2})}B_{i_1}
\hat I_{1000_{\tau_{p+1},\tau_p}}
^{*(i_4i_{3}i_{2}i_{1})}+
$$

\vspace{-2mm}
$$
+G_0^{(i_{4})}\bar LG_0^{(i_{3})}G_0^{(i_{2})}B_{i_1}
\left(\hat I_{1000_{\tau_{p+1},\tau_p}}
^{*(i_4i_{3}i_{2}i_{1})}-\hat I_{0100_{\tau_{p+1},\tau_p}}
^{*(i_4i_{3}i_{2}i_{1})}\right)+
$$
$$
\Biggl.+G_0^{(i_{4})}G_0^{(i_{3})}G_0^{(i_{2})}\bar LB_{i_1}
\left(\hat I_{0010_{\tau_{p+1},\tau_p}}
^{*(i_4i_{3}i_{2}i_{1})}-\hat I_{0001_{\tau_{p+1},\tau_p}}
^{*(i_4i_{3}i_{2}i_{1})}\right)\Biggr]+
$$
$$
+\sum_{i_{1},i_2,i_3,i_4,i_5,i_{6}=1}^m
G_0^{(i_{6})}G_0^{(i_{5})}
G_0^{(i_{4})}G_0^{(i_{3})}G_0^{(i_{2})}B_{i_{1}}
\hat I_{000000_{\tau_{p+1},\tau_p}}^{*(i_6i_{5}i_{4}i_{3}i_{2}i_{1})}.
$$

\vspace{3mm}

It is well known \cite{Zapad-3} that under the standard conditions
the numerical schemes (\ref{al1x})--(\ref{al5x}) 
have strong orders of convergence 1.0, 1.5, 2.0, 2.5, and 3.0
correspondingly.
Among these conditions we consider only the condition
for approximations of iterated Stratonovich stochastic 
integrals from the numerical
schemes (\ref{al1x})--(\ref{al5x}) \cite{Zapad-3} (also see \cite{12})
$$
{\sf M}\left\{\Biggl(I_{{l_{1}\ldots l_{k}}_{\tau_{p+1},\tau_p}}
^{*(i_{1}\ldots i_{k})} 
-\hat I_{{l_{1}\ldots l_{k}}_{\tau_{p+1},\tau_p}}^{*(i_{1}\ldots i_{k})}
\Biggr)^2\right\}\le C\Delta^{r+1},
$$
where constant $C$ is independent of $\Delta$ and
$r/2$ are strong orders of convergence for the numerical schemes
(\ref{al1x})--(\ref{al5x}), i.e. $r/2=1.0, 1.5,$ $2.0, 2.5,$ and $3.0.$

As we mentioned above, the numerical schemes (\ref{al1x})--(\ref{al5x})
are unrealizable in practice without 
procedures for the numerical simulation 
of iterated Stratonovich stochastic integrals
from (\ref{15.001x}).

In Chapter 5 
we give an extensive material on the mean-square 
approximation of specific iterated
It\^{o} and Stratonovich stochastic integrals
from the numerical schemes (\ref{al1})--(\ref{al5}),
(\ref{al1x})--(\ref{al5x}).
The mentioned material based on the results of Chapters 1 and 2.

\vspace{-4mm}

\chapter{Mean-Square Approximation of Specific Iterated It\^{o} 
and Stratonovich
Sto\-chas\-tic Integrals of Multiplicities 1 to 6 
from the Taylor--It\^{o} and Taylor--Stra\-to\-no\-vich Expansions
Based on Theorems From Chapters 1 and 2}

\vspace{-4mm}

\section{Mean-Square Approximation of Specific Iterated It\^{o} and 
Stratonovich Stochastic Integrals of multiplicities 1 to 6 
Based on Legendre Polynomials}

This section is devoted to the extensive practical material on
expansions and mean-square approximations of specific iterated
It\^{o} and Stratonovich stochastic integrals
of multiplicities 1 to 6 
on the base of Theorems 1.1, 2.1--2.9, 2.33--2.36, 2.50, 2.51, 2.62, 2.63 and multiple 
Fourier--Legendre series.
The considered iterated It\^{o} and 
Stratonovich stochastic integrals are part of the Taylor--It\^{o} and 
Taylor--Stratonovich expansions.
Therefore, the results of this section can be useful for the numerical 
solution of It\^{o} SDEs with non-commutative noise.

Consider the following iterated
It\^{o} and Stratonovich stochastic integrals
\begin{equation}
\label{ito-ito}
J[\psi^{(k)}]_{T,t}=\int\limits_t^T\psi_k(t_k) \ldots \int\limits_t^{t_{2}}
\psi_1(t_1) d{\bf w}_{t_1}^{(i_1)}\ldots
d{\bf w}_{t_k}^{(i_k)},
\end{equation}
\begin{equation}
\label{str-str}
J^{*}[\psi^{(k)}]_{T,t}=
{\int\limits_t^{*}}^T
\psi_k(t_k) \ldots 
{\int\limits_t^{*}}^{t_{2}}
\psi_1(t_1) d{\bf w}_{t_1}^{(i_1)}\ldots
d{\bf w}_{t_k}^{(i_k)},
\end{equation}
where every $\psi_l(\tau)$ $(l=1,\ldots,k)$ is a continuous 
nonrandom function 
on $[t,T],$ ${\bf w}_{\tau}^{(i)}$
$(i=1,\ldots,m)$ are independent 
standard Wiener processes,
${\bf w}_{\tau}^{(0)}=\tau,$ 
$i_1,\ldots,i_k = 0, 1,\ldots,m.$

As we saw in Chapter 4, $\psi_l(\tau)\equiv 1$ $(l=1,\ldots,k)$ and
$i_1,\ldots,i_k = 0, 1,\ldots,m$ 
in (\ref{ito-ito}), (\ref{str-str}) if we consider the iterated
stochastic integrals from
the classical Taylor--It\^{o} and Taylor--Stratonovich expansions
\cite{Zapad-3}. At the same time
$\psi_l(\tau)\equiv (t-\tau)^{q_l}$ ($l=1,\ldots,k$,\
$q_1,\ldots,q_k=0, 1, 2,\ldots $) and $i_1,\ldots,i_k = 1,\ldots,m$ 
for the iterated stochastic integrals from the unified
Taylor--It\^{o} and Taylor--Stratonovich expansions
\cite{1}-\cite{12aa}, \cite{arxiv-24}, \cite{Kul-Kuz}, \cite{Kuz}.

Thus, in this section, we will 
consider the following  
collections 
of iterated It\^{o} and Stratonovich stochastic integrals
\begin{equation}
\label{k1000xxxx}
I_{(l_1\ldots l_k)T,t}^{(i_1\ldots i_k)}
=\int\limits_t^T(t-t_k)^{l_k} \ldots \int\limits_t^{t_{2}}
(t-t_1)^{l_1} d{\bf w}_{t_1}^{(i_1)}\ldots
d{\bf w}_{t_k}^{(i_k)},
\end{equation}
\begin{equation}
\label{k1001xxxx}
I_{(l_1\ldots l_k)T,t}^{*(i_1\ldots i_k)}
=
{\int\limits_t^{*}}^T
(t-t_k)^{l_k} \ldots 
{\int\limits_t^{*}}^{t_{2}}
(t-t_1)^{l_1} d{\bf w}_{t_1}^{(i_1)}\ldots
d{\bf w}_{t_k}^{(i_k)},
\end{equation}
where $i_1,\ldots, i_k=1,\dots,m,$  $l_1,\ldots,l_k=0, 1,\ldots$

The complete orthonormal system of Legendre polynomials in the 
space $L_2([t,T])$ looks as follows
\begin{equation}
\label{4009d}
~~~~~~\phi_j(x)=\sqrt{\frac{2j+1}{T-t}}P_j\left(\left(
x-\frac{T+t}{2}\right)\frac{2}{T-t}\right),\ \ \ j=0, 1, 2,\ldots,
\end{equation}
where 
\begin{equation}
\label{agentww1}
P_j(x)=\frac{1}{2^j j!} \frac{d^j}{dx^j}\left(x^2-1\right)^j
\end{equation}
is the Legendre polynomial.

Let us recall some properties of Legendre polynomials \cite{suet}
(see Sect.~2.1.2)
$$
P_j(1)=1,\ \ \ P_{j+1}(-1)=-P_j(-1),\ \ \ j=0, 1, 2,\ldots,
$$
$$
\frac{dP_{j+1}}{dx}(x)-\frac{dP_{j-1}}{dx}(x)=(2j+1)P_j(x),
$$
$$
xP_{j}(x)=\frac{(j+1)P_{j+1}(x)+jP_{j-1}(x)}{2j+1},\ \ \ j=1, 2,\ldots,
$$
$$
\int\limits_{-1}^1 x^k P_j(x)dx=0,\ \ \ k=0, 1,\ldots,j-1,
$$
$$
\int\limits_{-1}^1 P_k(x)P_j(x)dx=\left\{
\begin{matrix}
0\ \  &\hbox{if}\ \  k\ne j\cr\cr
2/(2j+1)\ \  &\hbox{if}\ \  k=j
\end{matrix}\ \ ,\right. 
$$

$$
P_n(x)P_m(x)=\sum_{k=0}^m K_{m,n,k} P_{n+m-2k}(x),
$$
where
$$
K_{m,n,k}=\frac{a_{m-k}a_k a_{n-k}}{a_{m+n-k}}\cdot
\frac{2n+2m-4k+1}{2n+2m-2k+1},\ \ \ a_k=\frac{(2k-1)!!}{k!},\ \ \ m\le n.
$$

\vspace{1mm}

Applying the above properties of 
the Legendre polynomial system
(\ref{4009d}) and Theorems 1.1, 2.1--2.9, 2.33--2.36, 2.50, 2.51, 2.62, 2.63,
we obtain the following expansions of iterated It\^{o} 
and Stratonovich stochastic integrals from the sets (\ref{k1000xxxx}),
(\ref{k1001xxxx})

\vspace{-4mm}
\begin{equation}
\label{4001}
I_{(0)T,t}^{(i_1)}=\sqrt{T-t}\zeta_0^{(i_1)},
\end{equation}

\vspace{-1mm}
\begin{equation}
\label{4002}
I_{(1)T,t}^{(i_1)}=-\frac{(T-t)^{3/2}}{2}\left(\zeta_0^{(i_1)}+
\frac{1}{\sqrt{3}}\zeta_1^{(i_1)}\right),
\end{equation}

\vspace{-1mm}
\begin{equation}
\label{4003}
I_{(2)T,t}^{(i_1)}=\frac{(T-t)^{5/2}}{3}\left(\zeta_0^{(i_1)}+
\frac{\sqrt{3}}{2}\zeta_1^{(i_1)}+
\frac{1}{2\sqrt{5}}\zeta_2^{(i_1)}\right),
\end{equation}

\vspace{-1mm}
\begin{equation}
\label{4004xxxx}
~~~~~I_{(00)T,t}^{*(i_1 i_2)}=
\frac{T-t}{2}\left(\zeta_0^{(i_1)}\zeta_0^{(i_2)}+\sum_{i=1}^{\infty}
\frac{1}{\sqrt{4i^2-1}}\left(
\zeta_{i-1}^{(i_1)}\zeta_{i}^{(i_2)}-
\zeta_i^{(i_1)}\zeta_{i-1}^{(i_2)}\right)\right),
\end{equation}

\vspace{-1mm}
\begin{equation}
\label{4004yes}
I_{(00)T,t}^{(i_1 i_2)}=
\frac{T-t}{2}\left(\zeta_0^{(i_1)}\zeta_0^{(i_2)}+\sum_{i=1}^{\infty}
\frac{1}{\sqrt{4i^2-1}}\left(
\zeta_{i-1}^{(i_1)}\zeta_{i}^{(i_2)}-
\zeta_i^{(i_1)}\zeta_{i-1}^{(i_2)}\right) - {\bf 1}_{\{i_1=i_2\}}\right),
\end{equation}

\vspace{-5mm}
$$
I_{(01)T,t}^{*(i_1 i_2)}=-\frac{T-t}{2}I_{(00)T,t}^{*(i_1 i_2)}
-\frac{(T-t)^2}{4}\Biggl(
\frac{1}{\sqrt{3}}\zeta_0^{(i_1)}\zeta_1^{(i_2)}+\Biggr.
$$

\vspace{-2mm}
\begin{equation}
\label{ud111}
~~~~~+\Biggl.\sum_{i=0}^{\infty}\Biggl(
\frac{(i+2)\zeta_i^{(i_1)}\zeta_{i+2}^{(i_2)}
-(i+1)\zeta_{i+2}^{(i_1)}\zeta_{i}^{(i_2)}}
{\sqrt{(2i+1)(2i+5)}(2i+3)}-
\frac{\zeta_i^{(i_1)}\zeta_{i}^{(i_2)}}{(2i-1)(2i+3)}\Biggr)\Biggr),
\end{equation}

\newpage
\noindent
$$
I_{(10)T,t}^{*(i_1 i_2)}=-\frac{T-t}{2}I_{(00)T,t}^{*(i_1 i_2)}
-\frac{(T-t)^2}{4}\Biggl(
\frac{1}{\sqrt{3}}\zeta_0^{(i_2)}\zeta_1^{(i_1)}+\Biggr.
$$
\begin{equation}
\label{4006}
~~~~~+\Biggl.\sum_{i=0}^{\infty}\Biggl(
\frac{(i+1)\zeta_{i+2}^{(i_2)}\zeta_{i}^{(i_1)}
-(i+2)\zeta_{i}^{(i_2)}\zeta_{i+2}^{(i_1)}}
{\sqrt{(2i+1)(2i+5)}(2i+3)}+
\frac{\zeta_i^{(i_1)}\zeta_{i}^{(i_2)}}{(2i-1)(2i+3)}\Biggr)\Biggr)
\end{equation}

\vspace{4mm}
\noindent
or
$$
I_{(01)T,t}^{*(i_1i_2)}
=\hbox{\vtop{\offinterlineskip\halign{
\hfil#\hfil\cr
{\rm l.i.m.}\cr
$\stackrel{}{{}_{p\to \infty}}$\cr
}} }
\sum_{j_1,j_2=0}^{p}
C_{j_2j_1}^{01}
\zeta_{j_1}^{(i_1)}\zeta_{j_2}^{(i_2)},
$$

$$
I_{(10)T,t}^{*(i_1i_2)}
=\hbox{\vtop{\offinterlineskip\halign{
\hfil#\hfil\cr
{\rm l.i.m.}\cr
$\stackrel{}{{}_{p\to \infty}}$\cr
}} }
\sum_{j_1,j_2=0}^{p}
C_{j_2j_1}^{10}
\zeta_{j_1}^{(i_1)}\zeta_{j_2}^{(i_2)},
$$

\vspace{2mm}
\noindent
where

\vspace{-4mm}
$$
C_{j_2j_1}^{01}
=\frac{\sqrt{(2j_1+1)(2j_2+1)}}{8}(T-t)^{2}\bar
C_{j_2j_1}^{01},
$$

\begin{equation}
\label{caaas1}
C_{j_2j_1}^{10}
=\frac{\sqrt{(2j_1+1)(2j_2+1)}}{8}(T-t)^{2}\bar
C_{j_2j_1}^{10},
\end{equation}
$$
\bar C_{j_2j_1}^{01}=-\int\limits_{-1}^{1}(1+y)P_{j_2}(y)
\int\limits_{-1}^{y}
P_{j_1}(x)dx dy,
$$
$$
\bar C_{j_2j_1}^{10}=-\int\limits_{-1}^{1}P_{j_2}(y)
\int\limits_{-1}^{y}
(1+x)P_{j_1}(x)dx dy;
$$

\vspace{3mm}
$$
I_{(10)T,t}^{(i_1 i_2)}=
I_{(10)T,t}^{*(i_1 i_2)}+
\frac{1}{4}{\bf 1}_{\{i_1=i_2\}}(T-t)^2\ \ \ {\rm w.~p.~1},
$$

$$
I_{(01)T,t}^{(i_1 i_2)}=
I_{(01)T,t}^{*(i_1 i_2)}+
\frac{1}{4}{\bf 1}_{\{i_1=i_2\}}(T-t)^2\ \ \ {\rm w.~p.~1},
$$

\vspace{3mm}

$$
I_{(01)T,t}^{(i_1 i_2)}=
-\frac{T-t}{2}
I_{(00)T,t}^{(i_1 i_2)}
-\frac{(T-t)^2}{4}\Biggl(
\frac{1}{\sqrt{3}}\zeta_0^{(i_1)}\zeta_1^{(i_2)}+\Biggr.
$$
\begin{equation}
\label{kol}
~~~~~+\Biggl.\sum_{i=0}^{\infty}\Biggl(
\frac{(i+2)\zeta_i^{(i_1)}\zeta_{i+2}^{(i_2)}
-(i+1)\zeta_{i+2}^{(i_1)}\zeta_{i}^{(i_2)}}
{\sqrt{(2i+1)(2i+5)}(2i+3)}-
\frac{\zeta_i^{(i_1)}\zeta_{i}^{(i_2)}}{(2i-1)(2i+3)}\Biggr)\Biggr),
\end{equation}

\newpage
\noindent
$$
I_{(10)T,t}^{(i_1 i_2)}=
-\frac{T-t}{2}I_{(00)T,t}^{(i_1 i_2)}
-\frac{(T-t)^2}{4}\Biggl(
\frac{1}{\sqrt{3}}\zeta_0^{(i_2)}\zeta_1^{(i_1)}+\Biggr.
$$
\begin{equation}
\label{kola}
~~~~~+\Biggl.\sum_{i=0}^{\infty}\Biggl(
\frac{(i+1)\zeta_{i+2}^{(i_2)}\zeta_{i}^{(i_1)}
-(i+2)\zeta_{i}^{(i_2)}\zeta_{i+2}^{(i_1)}}
{\sqrt{(2i+1)(2i+5)}(2i+3)}+
\frac{\zeta_i^{(i_1)}\zeta_{i}^{(i_2)}}{(2i-1)(2i+3)}\Biggr)\Biggr)
\end{equation}

\vspace{6mm}
\noindent
or
$$
I_{(01)T,t}^{(i_1 i_2)}=
\hbox{\vtop{\offinterlineskip\halign{
\hfil#\hfil\cr
{\rm l.i.m.}\cr
$\stackrel{}{{}_{p\to \infty}}$\cr
}} }
\sum_{j_1,j_2=0}^{p}
C_{j_2j_1}^{01}\Biggl(\zeta_{j_1}^{(i_1)}\zeta_{j_2}^{(i_2)}
-{\bf 1}_{\{i_1=i_2\}}
{\bf 1}_{\{j_1=j_2\}}\Biggr),
$$

$$
I_{(10)T,t}^{(i_1 i_2)}=
\hbox{\vtop{\offinterlineskip\halign{
\hfil#\hfil\cr
{\rm l.i.m.}\cr
$\stackrel{}{{}_{p\to \infty}}$\cr
}} }
\sum_{j_1,j_2=0}^{p}
C_{j_2j_1}^{10}\Biggl(\zeta_{j_1}^{(i_1)}\zeta_{j_2}^{(i_2)}
-{\bf 1}_{\{i_1=i_2\}}
{\bf 1}_{\{j_1=j_2\}}\Biggr),
$$

\vspace{5mm}

\begin{equation}
\label{good1}
I_{(000)T,t}^{*(i_1 i_2 i_3)}=
\hbox{\vtop{\offinterlineskip\halign{
\hfil#\hfil\cr
{\rm l.i.m.}\cr
$\stackrel{}{{}_{p\to \infty}}$\cr
}} }
\sum\limits_{j_1, j_2, j_3=0}^{p}
C_{j_3 j_2 j_1}\zeta_{j_1}^{(i_1)}\zeta_{j_2}^{(i_2)}\zeta_{j_3}^{(i_3)},
\end{equation}

\vspace{3mm}

$$
I_{(000)T,t}^{(i_1i_2i_3)}
=\hbox{\vtop{\offinterlineskip\halign{
\hfil#\hfil\cr
{\rm l.i.m.}\cr
$\stackrel{}{{}_{p\to \infty}}$\cr
}} }
\sum_{j_1,j_2,j_3=0}^{p}
C_{j_3j_2j_1}
\Biggl(
\zeta_{j_1}^{(i_1)}\zeta_{j_2}^{(i_2)}\zeta_{j_3}^{(i_3)}
-{\bf 1}_{\{i_1=i_2\}}
{\bf 1}_{\{j_1=j_2\}}
\zeta_{j_3}^{(i_3)}-
\Biggr.
$$
\begin{equation}
\label{zzz1}
\Biggl.
-{\bf 1}_{\{i_2=i_3\}}
{\bf 1}_{\{j_2=j_3\}}
\zeta_{j_1}^{(i_1)}-
{\bf 1}_{\{i_1=i_3\}}
{\bf 1}_{\{j_1=j_3\}}
\zeta_{j_2}^{(i_2)}\Biggr),
\end{equation}

\vspace{6mm}

$$
I_{(000)T,t}^{(i_1 i_1 i_1)}
=\frac{1}{6}(T-t)^{3/2}\biggl(
\left(\zeta_0^{(i_1)}\right)^3-3
\zeta_0^{(i_1)}\biggr)\ \ \ \hbox{w.~p.~1},
$$

\begin{equation}
\label{dest4}
I_{(000)T,t}^{*(i_1 i_1 i_1)}
=\frac{1}{6}(T-t)^{3/2}
\left(\zeta_0^{(i_1)}\right)^3\ \ \ \hbox{w.~p.~1},
\end{equation}

\vspace{3mm}
\noindent
where

\vspace{-4mm}
\begin{equation}
\label{zzz2}
C_{j_3j_2j_1}
=\frac{\sqrt{(2j_1+1)(2j_2+1)(2j_3+1)}}{8}(T-t)^{3/2}\bar
C_{j_3j_2j_1},
\end{equation}
\begin{equation}
\label{zzz3}
\bar C_{j_3j_2j_1}=\int\limits_{-1}^{1}P_{j_3}(z)
\int\limits_{-1}^{z}P_{j_2}(y)
\int\limits_{-1}^{y}
P_{j_1}(x)dx dy dz;
\end{equation}

\noindent
here and further in this section
$$
\zeta_{j}^{(i)}=
\int\limits_t^T \phi_{j}(s) d{\bf w}_s^{(i)}\ \ \ (i=1,\ldots,m,\
j=0,1,\ldots)
$$

\noindent
are independent standard Gaussian random variables
for various
$i$ or $j$;

$$
I_{(000)T,t}^{(i_1 i_2 i_3)}=I_{(000)T,t}^{*(i_1 i_2 i_3)}
+\frac{1}{2}{\bf 1}_{\{i_1=i_2\ne 0\}}I_{(1)T,t}^{(i_3)}-
$$
$$
-
\frac{1}{2}{\bf 1}_{\{i_2=i_3\ne 0\}}\left((T-t)
I_{(0)T,t}^{(i_1)}
+I_{(1)T,t}^{(i_1)}\right)\ \ \ \hbox{w.~p.~1},
$$

\vspace{2mm}

$$
I_{(02)T,t}^{*(i_1 i_2)}=-\frac{(T-t)^2}{4}I_{(00)T,t}^{*(i_1 i_2)}
-(T-t)I_{(01)T,t}^{*(i_1 i_2)}+
\frac{(T-t)^3}{8}\Biggl[
\frac{2}{3\sqrt{5}}\zeta_2^{(i_2)}\zeta_0^{(i_1)}+\Biggr.
$$
$$
+\frac{1}{3}\zeta_0^{(i_1)}\zeta_0^{(i_2)}+
\sum_{i=0}^{\infty}\Biggl(
\frac{(i+2)(i+3)\zeta_{i+3}^{(i_2)}\zeta_{i}^{(i_1)}
-(i+1)(i+2)\zeta_{i}^{(i_2)}\zeta_{i+3}^{(i_1)}}
{\sqrt{(2i+1)(2i+7)}(2i+3)(2i+5)}+
\Biggr.
$$

\vspace{1mm}
\begin{equation}
\label{leto1}
\Biggl.\Biggl.+
\frac{(i^2+i-3)\zeta_{i+1}^{(i_2)}\zeta_{i}^{(i_1)}
-(i^2+3i-1)\zeta_{i}^{(i_2)}\zeta_{i+1}^{(i_1)}}
{\sqrt{(2i+1)(2i+3)}(2i-1)(2i+5)}\Biggr)\Biggr],
\end{equation}

\vspace{6mm}

$$
I_{(20)T,t}^{*(i_1 i_2)}=-\frac{(T-t)^2}{4}I_{(00)T,t}^{*(i_1 i_2)}
-(T-t)I_{(10)T,t}^{*(i_1 i_2)}+
\frac{(T-t)^3}{8}\Biggl[
\frac{2}{3\sqrt{5}}\zeta_0^{(i_2)}\zeta_2^{(i_1)}+\Biggr.
$$
$$
+\frac{1}{3}\zeta_0^{(i_1)}\zeta_0^{(i_2)}+
\sum_{i=0}^{\infty}\Biggl(
\frac{(i+1)(i+2)\zeta_{i+3}^{(i_2)}\zeta_{i}^{(i_1)}
-(i+2)(i+3)\zeta_{i}^{(i_2)}\zeta_{i+3}^{(i_1)}}
{\sqrt{(2i+1)(2i+7)}(2i+3)(2i+5)}+
\Biggr.
$$

\vspace{1mm}
\begin{equation}
\label{leto2}
\Biggl.\Biggl.+
\frac{(i^2+3i-1)\zeta_{i+1}^{(i_2)}\zeta_{i}^{(i_1)}
-(i^2+i-3)\zeta_{i}^{(i_2)}\zeta_{i+1}^{(i_1)}}
{\sqrt{(2i+1)(2i+3)}(2i-1)(2i+5)}\Biggr)\Biggr],
\end{equation}

\vspace{11mm}

$$
I_{(11)T,t}^{*(i_1 i_2)}=-\frac{(T-t)^2}{4}I_{(00)T,t}^{*(i_1 i_2)}-
\frac{(T-t)}{2}\left(I_{(10)T,t}^{*(i_1 i_2)}+
I_{(01)T,t}^{*(i_1 i_2)}\right)+
$$

\vspace{-2mm}
$$
+\frac{(T-t)^3}{8}\Biggl[
\frac{1}{3}\zeta_1^{(i_1)}\zeta_1^{(i_2)}+
\sum_{i=0}^{\infty}\Biggl(
\frac{(i+1)(i+3)\left(\zeta_{i+3}^{(i_2)}\zeta_{i}^{(i_1)}
-\zeta_{i}^{(i_2)}\zeta_{i+3}^{(i_1)}\right)}
{\sqrt{(2i+1)(2i+7)}(2i+3)(2i+5)}+
\Biggr.
$$

\newpage
\noindent
\begin{equation}
\label{leto3}
\Biggl.\Biggl.+
\frac{(i+1)^2\left(\zeta_{i+1}^{(i_2)}\zeta_{i}^{(i_1)}
-\zeta_{i}^{(i_2)}\zeta_{i+1}^{(i_1)}\right)}
{\sqrt{(2i+1)(2i+3)}(2i-1)(2i+5)}\Biggr)\Biggr]
\end{equation}

\vspace{5mm}
\noindent
or
$$
I_{(02)T,t}^{*(i_1 i_2)}=
\hbox{\vtop{\offinterlineskip\halign{
\hfil#\hfil\cr
{\rm l.i.m.}\cr
$\stackrel{}{{}_{p\to \infty}}$\cr
}} }
\sum_{j_1,j_2=0}^{p}
C_{j_2j_1}^{02}\zeta_{j_1}^{(i_1)}\zeta_{j_2}^{(i_2)},
$$
$$
I_{(20)T,t}^{*(i_1 i_2)}=
\hbox{\vtop{\offinterlineskip\halign{
\hfil#\hfil\cr
{\rm l.i.m.}\cr
$\stackrel{}{{}_{p\to \infty}}$\cr
}} }
\sum_{j_1,j_2=0}^{p}
C_{j_2j_1}^{20}\zeta_{j_1}^{(i_1)}\zeta_{j_2}^{(i_2)},
$$
$$
I_{(11)T,t}^{*(i_1 i_2)}=
\hbox{\vtop{\offinterlineskip\halign{
\hfil#\hfil\cr
{\rm l.i.m.}\cr
$\stackrel{}{{}_{p\to \infty}}$\cr
}} }
\sum_{j_1,j_2=0}^{p}
C_{j_2j_1}^{11}\zeta_{j_1}^{(i_1)}\zeta_{j_2}^{(i_2)},
$$

\vspace{1mm}
\noindent
where

\vspace{-4mm}
$$
C_{j_2j_1}^{02}=
\frac{\sqrt{(2j_1+1)(2j_2+1)}}{16}(T-t)^{3}\bar
C_{j_2j_1}^{02},
$$

$$
C_{j_2j_1}^{20}=
\frac{\sqrt{(2j_1+1)(2j_2+1)}}{16}(T-t)^{3}\bar
C_{j_2j_1}^{20},
$$

$$
C_{j_2j_1}^{11}=
\frac{\sqrt{(2j_1+1)(2j_2+1)}}{16}(T-t)^{3}\bar
C_{j_2j_1}^{11}, 
$$

$$
\bar C_{j_2j_1}^{02}=
\int\limits_{-1}^{1}P_{j_2}(y)(y+1)^2
\int\limits_{-1}^{y}
P_{j_1}(x)dx dy,
$$
$$
\bar C_{j_2j_1}^{20}=
\int\limits_{-1}^{1}P_{j_2}(y)
\int\limits_{-1}^{y}
P_{j_1}(x)(x+1)^2 dx dy,
$$
$$
\bar C_{j_2j_1}^{11}=
\int\limits_{-1}^{1}P_{j_2}(y)(y+1)
\int\limits_{-1}^{y}
P_{j_1}(x)(x+1)dx dy,
$$

\vspace{4mm}

$$
I_{(11)T,t}^{*(i_1 i_1)}=\frac{1}{2}\left(I_{(1)T,t}^{(i_1)}
\right)^2\ \ \ \hbox{w.~p.~1,}
$$
\begin{equation}
\label{seg1}
I_{(02)T,t}^{(i_1 i_2)}=
I_{(02)T,t}^{*(i_1 i_2)}-
\frac{1}{6}{\bf 1}_{\{i_1=i_2\}}(T-t)^3\ \ \ \hbox{w.~p.~1},
\end{equation}
\begin{equation}
\label{seg2}
I_{(20)T,t}^{(i_1 i_2)}=
I_{(20)T,t}^{*(i_1 i_2)}-
\frac{1}{6}{\bf 1}_{\{i_1=i_2\}}(T-t)^3\ \ \ \hbox{w.~p.~1},
\end{equation}
$$
I_{(11)T,t}^{(i_1 i_2)}=
I_{(11)T,t}^{*(i_1 i_2)}-
\frac{1}{6}{\bf 1}_{\{i_1=i_2\}}(T-t)^3\ \ \ \hbox{w.~p.~1},
$$

\vspace{3mm}

$$
I_{(02)T,t}^{(i_1 i_2)}
=-\frac{(T-t)^2}{4}I_{(00)T,t}^{(i_1 i_2)}
-(T-t) I_{01_{T,t}}^{(i_1 i_2)}+
\frac{(T-t)^3}{8}\Biggl[
\frac{2}{3\sqrt{5}}\zeta_2^{(i_2)}\zeta_0^{(i_1)}+\Biggr.
$$
$$
+\frac{1}{3}\zeta_0^{(i_1)}\zeta_0^{(i_2)}+
\sum_{i=0}^{\infty}\Biggl(
\frac{(i+2)(i+3)\zeta_{i+3}^{(i_2)}\zeta_{i}^{(i_1)}
-(i+1)(i+2)\zeta_{i}^{(i_2)}\zeta_{i+3}^{(i_1)}}
{\sqrt{(2i+1)(2i+7)}(2i+3)(2i+5)}+
\Biggr.
$$

$$
\Biggl.\Biggl.+
\frac{(i^2+i-3)\zeta_{i+1}^{(i_2)}\zeta_{i}^{(i_1)}
-(i^2+3i-1)\zeta_{i}^{(i_2)}\zeta_{i+1}^{(i_1)}}
{\sqrt{(2i+1)(2i+3)}(2i-1)(2i+5)}\Biggr)\Biggr] - 
$$

\begin{equation}
\label{quak1}
-\frac{1}{24}{\bf 1}_{\{i_1=i_2\}}{(T-t)^3},
\end{equation}

\vspace{4mm}

$$
I_{(20)T,t}^{(i_1 i_2)}=-\frac{(T-t)^2}{4}
I_{(00)T,t}^{(i_1 i_2)}
-(T-t) I_{(10)T,t}^{(i_1 i_2)}+
\frac{(T-t)^3}{8}\Biggl[
\frac{2}{3\sqrt{5}}\zeta_0^{(i_2)}\zeta_2^{(i_1)}+\Biggr.
$$
$$
+\frac{1}{3}\zeta_0^{(i_1)}\zeta_0^{(i_2)}+
\sum_{i=0}^{\infty}\Biggl(
\frac{(i+1)(i+2)\zeta_{i+3}^{(i_2)}\zeta_{i}^{(i_1)}
-(i+2)(i+3)\zeta_{i}^{(i_2)}\zeta_{i+3}^{(i_1)}}
{\sqrt{(2i+1)(2i+7)}(2i+3)(2i+5)}+
\Biggr.
$$

$$
\Biggl.\Biggl.+
\frac{(i^2+3i-1)\zeta_{i+1}^{(i_2)}\zeta_{i}^{(i_1)}
-(i^2+i-3)\zeta_{i}^{(i_2)}\zeta_{i+1}^{(i_1)}}
{\sqrt{(2i+1)(2i+3)}(2i-1)(2i+5)}\Biggr)\Biggr] - 
$$
\begin{equation}
\label{quak2}
-
\frac{1}{24}{\bf 1}_{\{i_1=i_2\}}{(T-t)^3},
\end{equation}

\vspace{11mm}

$$
I_{(11)T,t}^{(i_1 i_2)}
=-\frac{(T-t)^2}{4}I_{(00)T,t}^{(i_1 i_2)}
-\frac{T-t}{2}\left(
I_{(10)T,t}^{(i_1 i_2)}+
I_{(01)T,t}^{(i_1 i_2)}\right)+
$$

\vspace{-4mm}
$$
+
\frac{(T-t)^3}{8}\Biggl[
\frac{1}{3}\zeta_1^{(i_1)}\zeta_1^{(i_2)}+\Biggr.
\sum_{i=0}^{\infty}\Biggl(
\frac{(i+1)(i+3)\left(\zeta_{i+3}^{(i_2)}\zeta_{i}^{(i_1)}
-\zeta_{i}^{(i_2)}\zeta_{i+3}^{(i_1)}\right)}
{\sqrt{(2i+1)(2i+7)}(2i+3)(2i+5)}+
\Biggr.
$$

\newpage
\noindent
$$
\Biggl.\Biggl.
+\frac{(i+1)^2\left(\zeta_{i+1}^{(i_2)}\zeta_{i}^{(i_1)}
-\zeta_{i}^{(i_2)}\zeta_{i+1}^{(i_1)}\right)}
{\sqrt{(2i+1)(2i+3)}(2i-1)(2i+5)}\Biggr)\Biggr] - 
$$

\begin{equation}
\label{quak3}
-
\frac{1}{24}{\bf 1}_{\{i_1=i_2\}}{(T-t)^3}
\end{equation}

\vspace{3mm}
\noindent
or
$$
I_{(02)T,t}^{(i_1 i_2)}=
\hbox{\vtop{\offinterlineskip\halign{
\hfil#\hfil\cr
{\rm l.i.m.}\cr
$\stackrel{}{{}_{p\to \infty}}$\cr
}} }
\sum_{j_1,j_2=0}^p
C_{j_2j_1}^{02}\Biggl(\zeta_{j_1}^{(i_1)}\zeta_{j_2}^{(i_2)}
-{\bf 1}_{\{i_1=i_2\}}
{\bf 1}_{\{j_1=j_2\}}\Biggr),
$$

\vspace{1mm}
$$
I_{(20)T,t}^{(i_1 i_2)}=
\hbox{\vtop{\offinterlineskip\halign{
\hfil#\hfil\cr
{\rm l.i.m.}\cr
$\stackrel{}{{}_{p\to \infty}}$\cr
}} }
\sum_{j_1,j_2=0}^{p}
C_{j_2j_1}^{20}\Biggl(\zeta_{j_1}^{(i_1)}\zeta_{j_2}^{(i_2)}
-{\bf 1}_{\{i_1=i_2\}}
{\bf 1}_{\{j_1=j_2\}}\Biggr),
$$

\vspace{1mm}
$$
I_{(11)T,t}^{(i_1 i_2)}=
\hbox{\vtop{\offinterlineskip\halign{
\hfil#\hfil\cr
{\rm l.i.m.}\cr
$\stackrel{}{{}_{p\to \infty}}$\cr
}} }
\sum_{j_1,j_2=0}^{p}
C_{j_2j_1}^{11}\Biggl(\zeta_{j_1}^{(i_1)}\zeta_{j_2}^{(i_2)}
-{\bf 1}_{\{i_1=i_2\}}
{\bf 1}_{\{j_1=j_2\}}\Biggr),
$$

\vspace{4mm}

\begin{equation}
\label{uruk}
~~~~~~~~~I_{(3)T,t}^{(i_1)}=-\frac{(T-t)^{7/2}}{4}\left(\zeta_0^{(i_1)}+
\frac{3\sqrt{3}}{5}\zeta_1^{(i_1)}+
\frac{1}{\sqrt{5}}\zeta_2^{(i_1)}+
\frac{1}{5\sqrt{7}}\zeta_3^{(i_1)}\right),
\end{equation}

\vspace{2mm}
$$
I_{(0000)T,t}^{*(i_1 i_2 i_3 i_4)}=
\hbox{\vtop{\offinterlineskip\halign{
\hfil#\hfil\cr
{\rm l.i.m.}\cr
$\stackrel{}{{}_{p\to \infty}}$\cr
}} }
\sum\limits_{j_1, j_2, j_3, j_4=0}^{p}
C_{j_4 j_3 j_2 j_1}\zeta_{j_1}^{(i_1)}\zeta_{j_2}^{(i_2)}\zeta_{j_3}^{(i_3)}
\zeta_{j_4}^{(i_4)},
$$

\vspace{7mm}

$$
I_{(0000)T,t}^{(i_1 i_2 i_3 i_4)}
=
\hbox{\vtop{\offinterlineskip\halign{
\hfil#\hfil\cr
{\rm l.i.m.}\cr
$\stackrel{}{{}_{p\to \infty}}$\cr
}} }
\sum_{j_1,j_2,j_3,j_4=0}^{p}
C_{j_4 j_3 j_2 j_1}\Biggl(
\zeta_{j_1}^{(i_1)}\zeta_{j_2}^{(i_2)}\zeta_{j_3}^{(i_3)}\zeta_{j_4}^{(i_4)}
-\Biggr.
$$
$$
-
{\bf 1}_{\{i_1=i_2\}}
{\bf 1}_{\{j_1=j_2\}}
\zeta_{j_3}^{(i_3)}
\zeta_{j_4}^{(i_4)}
-
{\bf 1}_{\{i_1=i_3\}}
{\bf 1}_{\{j_1=j_3\}}
\zeta_{j_2}^{(i_2)}
\zeta_{j_4}^{(i_4)}-
$$
$$
-
{\bf 1}_{\{i_1=i_4\}}
{\bf 1}_{\{j_1=j_4\}}
\zeta_{j_2}^{(i_2)}
\zeta_{j_3}^{(i_3)}
-
{\bf 1}_{\{i_2=i_3\}}
{\bf 1}_{\{j_2=j_3\}}
\zeta_{j_1}^{(i_1)}
\zeta_{j_4}^{(i_4)}-
$$
$$
-
{\bf 1}_{\{i_2=i_4\}}
{\bf 1}_{\{j_2=j_4\}}
\zeta_{j_1}^{(i_1)}
\zeta_{j_3}^{(i_3)}
-
{\bf 1}_{\{i_3=i_4\}}
{\bf 1}_{\{j_3=j_4\}}
\zeta_{j_1}^{(i_1)}
\zeta_{j_2}^{(i_2)}+
$$
$$
+
{\bf 1}_{\{i_1=i_2\}}
{\bf 1}_{\{j_1=j_2\}}
{\bf 1}_{\{i_3=i_4\}}
{\bf 1}_{\{j_3=j_4\}}+
{\bf 1}_{\{i_1=i_3\}}
{\bf 1}_{\{j_1=j_3\}}
{\bf 1}_{\{i_2=i_4\}}
{\bf 1}_{\{j_2=j_4\}}+
$$
\begin{equation}
\label{zzz10}
+\Biggl.
{\bf 1}_{\{i_1=i_4\}}
{\bf 1}_{\{j_1=j_4\}}
{\bf 1}_{\{i_2=i_3\}}
{\bf 1}_{\{j_2=j_3\}}\Biggr),
\end{equation}

\newpage
\noindent
$$
I_{(0000)T,t}^{(i_1i_1i_1i_1)}=
\frac{1}{24}(T-t)^2
\left(\left(\zeta_0^{(i_1)}\right)^4-
6\left(\zeta_0^{(i_1)}\right)^2+3\right)\ \ \ \hbox{w.~p.~1},
$$
\begin{equation}
\label{creat100}
I_{(0000)T,t}^{*(i_1i_1i_1i_1)}=
\frac{1}{24}(T-t)^2
\left(\zeta_0^{(i_1)}\right)^4\ \ \ \hbox{w.~p.~1},
\end{equation}

\vspace{3mm}
\noindent
where

\vspace{-3mm}
\begin{equation}
\label{zzz11}
C_{j_4j_3j_2j_1}=
\frac{\sqrt{(2j_1+1)(2j_2+1)(2j_3+1)(2j_4+1)}}{16}(T-t)^{2}\bar
C_{j_4j_3j_2j_1},
\end{equation}
\begin{equation}
\label{zzz12}
\bar C_{j_4j_3j_2j_1}=\int\limits_{-1}^{1}P_{j_4}(u)
\int\limits_{-1}^{u}P_{j_3}(z)
\int\limits_{-1}^{z}P_{j_2}(y)
\int\limits_{-1}^{y}
P_{j_1}(x)dx dy dz du;
\end{equation}

\vspace{4mm}

$$
I_{(001)T,t}^{*(i_1i_2i_3)}
=\hbox{\vtop{\offinterlineskip\halign{
\hfil#\hfil\cr
{\rm l.i.m.}\cr
$\stackrel{}{{}_{p\to \infty}}$\cr
}} }
\sum_{j_1,j_2,j_3=0}^{p}
C_{j_3 j_2 j_1}^{001}
\zeta_{j_1}^{(i_1)}\zeta_{j_2}^{(i_2)}\zeta_{j_3}^{(i_3)},
$$

$$
I_{(010)T,t}^{*(i_1i_2i_3)}
=\hbox{\vtop{\offinterlineskip\halign{
\hfil#\hfil\cr
{\rm l.i.m.}\cr
$\stackrel{}{{}_{p\to \infty}}$\cr
}} }
\sum_{j_1,j_2,j_3=0}^{p}
C_{j_3 j_2 j_1}^{010}
\zeta_{j_1}^{(i_1)}\zeta_{j_2}^{(i_2)}\zeta_{j_3}^{(i_3)},
$$

$$
I_{(100)T,t}^{*(i_1i_2i_3)}
=\hbox{\vtop{\offinterlineskip\halign{
\hfil#\hfil\cr
{\rm l.i.m.}\cr
$\stackrel{}{{}_{p\to \infty}}$\cr
}} }
\sum_{j_1,j_2,j_3=0}^{p}
C_{j_3 j_2 j_1}^{100}
\zeta_{j_1}^{(i_1)}\zeta_{j_2}^{(i_2)}\zeta_{j_3}^{(i_3)},
$$

\vspace{4mm}

$$
I_{(001)T,t}^{(i_1i_2i_3)}
=\hbox{\vtop{\offinterlineskip\halign{
\hfil#\hfil\cr
{\rm l.i.m.}\cr
$\stackrel{}{{}_{p\to \infty}}$\cr
}} }
\sum_{j_1,j_2,j_3=0}^{p}
C_{j_3j_2j_1}^{001}\Biggl(
\zeta_{j_1}^{(i_1)}\zeta_{j_2}^{(i_2)}\zeta_{j_3}^{(i_3)}
-{\bf 1}_{\{i_1=i_2\}}
{\bf 1}_{\{j_1=j_2\}}
\zeta_{j_3}^{(i_3)}-
\Biggr.
$$
\begin{equation}
\label{sss1}
\Biggl.
-{\bf 1}_{\{i_2=i_3\}}
{\bf 1}_{\{j_2=j_3\}}
\zeta_{j_1}^{(i_1)}-
{\bf 1}_{\{i_1=i_3\}}
{\bf 1}_{\{j_1=j_3\}}
\zeta_{j_2}^{(i_2)}\Biggr),
\end{equation}

\vspace{4mm}

$$
I_{(010)T,t}^{(i_1i_2i_3)}
=\hbox{\vtop{\offinterlineskip\halign{
\hfil#\hfil\cr
{\rm l.i.m.}\cr
$\stackrel{}{{}_{p\to \infty}}$\cr
}} }
\sum_{j_1,j_2,j_3=0}^{p}
C_{j_3j_2j_1}^{010}\Biggl(
\zeta_{j_1}^{(i_1)}\zeta_{j_2}^{(i_2)}\zeta_{j_3}^{(i_3)}
-{\bf 1}_{\{i_1=i_2\}}
{\bf 1}_{\{j_1=j_2\}}
\zeta_{j_3}^{(i_3)}-
\Biggr.
$$
\begin{equation}
\label{sss2}
\Biggl.
-{\bf 1}_{\{i_2=i_3\}}
{\bf 1}_{\{j_2=j_3\}}
\zeta_{j_1}^{(i_1)}-
{\bf 1}_{\{i_1=i_3\}}
{\bf 1}_{\{j_1=j_3\}}
\zeta_{j_2}^{(i_2)}\Biggr),
\end{equation}

\newpage
\noindent
$$
I_{(100)T,t}^{(i_1i_2i_3)}
=\hbox{\vtop{\offinterlineskip\halign{
\hfil#\hfil\cr
{\rm l.i.m.}\cr
$\stackrel{}{{}_{p\to \infty}}$\cr
}} }
\sum_{j_1,j_2,j_3=0}^{p}
C_{j_3j_2j_1}^{100}\Biggl(
\zeta_{j_1}^{(i_1)}\zeta_{j_2}^{(i_2)}\zeta_{j_3}^{(i_3)}
-{\bf 1}_{\{i_1=i_2\}}
{\bf 1}_{\{j_1=j_2\}}
\zeta_{j_3}^{(i_3)}-
\Biggr.
$$
\begin{equation}
\label{sss3}
\Biggl.
-{\bf 1}_{\{i_2=i_3\}}
{\bf 1}_{\{j_2=j_3\}}
\zeta_{j_1}^{(i_1)}-
{\bf 1}_{\{i_1=i_3\}}
{\bf 1}_{\{j_1=j_3\}}
\zeta_{j_2}^{(i_2)}\Biggr),
\end{equation}

\vspace{1mm}
\noindent
where

\vspace{-2mm}
$$
C_{j_3j_2j_1}^{001}
=\frac{\sqrt{(2j_1+1)(2j_2+1)(2j_3+1)}}{16}(T-t)^{5/2}\bar
C_{j_3j_2j_1}^{001},
$$

\vspace{1mm}
$$
C_{j_3j_2j_1}^{010}
=\frac{\sqrt{(2j_1+1)(2j_2+1)(2j_3+1)}}{16}(T-t)^{5/2}\bar
C_{j_3j_2j_1}^{010},
$$

\vspace{1mm}
$$
C_{j_3j_2j_1}^{100}
=\frac{\sqrt{(2j_1+1)(2j_2+1)(2j_3+1)}}{16}(T-t)^{5/2}\bar
C_{j_3j_2j_1}^{100},
$$

\vspace{1mm}
$$
\bar C_{j_3j_2j_1}^{100}=-
\int\limits_{-1}^{1}P_{j_3}(z)
\int\limits_{-1}^{z}P_{j_2}(y)
\int\limits_{-1}^{y}
P_{j_1}(x)(x+1)dx dy dz,
$$
$$
\bar C_{j_3j_2j_1}^{010}=-
\int\limits_{-1}^{1}P_{j_3}(z)
\int\limits_{-1}^{z}P_{j_2}(y)(y+1)
\int\limits_{-1}^{y}
P_{j_1}(x)dx dy dz,
$$
$$
\bar C_{j_3j_2j_1}^{001}=-
\int\limits_{-1}^{1}P_{j_3}(z)(z+1)
\int\limits_{-1}^{z}P_{j_2}(y)
\int\limits_{-1}^{y}
P_{j_1}(x)dx dy dz;
$$

\vspace{4mm}

$$
I_{(lll)T,t}^{(i_1i_1i_1)}=
\frac{1}{6}\left(\left(I_{(l)T,t}^{(i_1)}\right)^3-
3I_{(l)T,t}^{(i_1)}\Delta_{l(T,t)}\right)\ \ \ \hbox{w.~p.~1},
$$

$$
I_{(lll)T,t}^{*(i_1i_1i_1)}=
\frac{1}{6}\left(I_{(l)T,t}^{(i_1)}\right)^3\ \ \ \hbox{w.~p.~1},
$$

$$
I_{(llll)T,t}^{(i_1i_1i_1i_1)}=
\frac{1}{24}\left(\left(I_{(l)T,t}^{(i_1)}\right)^4-
6\left(I_{(l)T,t}^{(i_1)}\right)^2\Delta_{(l)T,t}+3
\left(\Delta_{(l)T,t}\right)^2\right)\ \ \ \hbox{w.~p.~1},
$$

$$
I_{(llll)T,t}^{*(i_1i_1i_1i_1)}=
\frac{1}{24}\left(I_{(l)T,t}^{(i_1)}\right)^4\ \ \ \hbox{w.~p.~1},
$$

\newpage
\noindent
where
\begin{equation}
\label{uruk1}
I_{(l)T,t}^{(i_1)}=\sum_{j=0}^l C_j^l \zeta_j^{(i_1)}\ \ \ \hbox{w.~p.~1},
\end{equation}
$$
\Delta_{l(T,t)}=\int\limits_t^T(t-s)^{2l}ds,\ \ \
C_j^l=\int\limits_t^T(t-s)^l\phi_j(s)ds;
$$ 

\vspace{2mm}

$$
I_{(00000)T,t}^{*(i_1 i_2 i_3 i_4 i_5)}=
\hbox{\vtop{\offinterlineskip\halign{
\hfil#\hfil\cr
{\rm l.i.m.}\cr
$\stackrel{}{{}_{p\to \infty}}$\cr
}} }
\sum\limits_{j_1, j_2, j_3, j_4, j_5=0}^{p}
C_{j_5j_4 j_3 j_2 j_1}
\zeta_{j_1}^{(i_1)}\zeta_{j_2}^{(i_2)}\zeta_{j_3}^{(i_3)}
\zeta_{j_4}^{(i_4)}\zeta_{j_5}^{(i_5)},
$$

\vspace{5mm}

$$
I_{(00000)T,t}^{(i_1 i_2 i_3 i_4 i_5)}
=
\hbox{\vtop{\offinterlineskip\halign{
\hfil#\hfil\cr
{\rm l.i.m.}\cr
$\stackrel{}{{}_{p\to \infty}}$\cr
}} }
\sum_{j_1,j_2,j_3,j_4,j_5=0}^p
C_{j_5 j_4 j_3 j_2 j_1}\Biggl(
\prod_{l=1}^5\zeta_{j_l}^{(i_l)}
-\Biggr.
$$
$$
-
{\bf 1}_{\{i_1=i_2\}}
{\bf 1}_{\{j_1=j_2\}}
\zeta_{j_3}^{(i_3)}
\zeta_{j_4}^{(i_4)}
\zeta_{j_5}^{(i_5)}-
{\bf 1}_{\{i_1=i_3\}}
{\bf 1}_{\{j_1=j_3\}}
\zeta_{j_2}^{(i_2)}
\zeta_{j_4}^{(i_4)}
\zeta_{j_5}^{(i_5)}-
$$
$$
-
{\bf 1}_{\{i_1=i_4\}}
{\bf 1}_{\{j_1=j_4\}}
\zeta_{j_2}^{(i_2)}
\zeta_{j_3}^{(i_3)}
\zeta_{j_5}^{(i_5)}-
{\bf 1}_{\{i_1=i_5\}}
{\bf 1}_{\{j_1=j_5\}}
\zeta_{j_2}^{(i_2)}
\zeta_{j_3}^{(i_3)}
\zeta_{j_4}^{(i_4)}-
$$
$$
-
{\bf 1}_{\{i_2=i_3\}}
{\bf 1}_{\{j_2=j_3\}}
\zeta_{j_1}^{(i_1)}
\zeta_{j_4}^{(i_4)}
\zeta_{j_5}^{(i_5)}-
{\bf 1}_{\{i_2=i_4\}}
{\bf 1}_{\{j_2=j_4\}}
\zeta_{j_1}^{(i_1)}
\zeta_{j_3}^{(i_3)}
\zeta_{j_5}^{(i_5)}-
$$
$$
-
{\bf 1}_{\{i_2=i_5\}}
{\bf 1}_{\{j_2=j_5\}}
\zeta_{j_1}^{(i_1)}
\zeta_{j_3}^{(i_3)}
\zeta_{j_4}^{(i_4)}
-{\bf 1}_{\{i_3=i_4\}}
{\bf 1}_{\{j_3=j_4\}}
\zeta_{j_1}^{(i_1)}
\zeta_{j_2}^{(i_2)}
\zeta_{j_5}^{(i_5)}-
$$
$$
-
{\bf 1}_{\{i_3=i_5\}}
{\bf 1}_{\{j_3=j_5\}}
\zeta_{j_1}^{(i_1)}
\zeta_{j_2}^{(i_2)}
\zeta_{j_4}^{(i_4)}
-{\bf 1}_{\{i_4=i_5\}}
{\bf 1}_{\{j_4=j_5\}}
\zeta_{j_1}^{(i_1)}
\zeta_{j_2}^{(i_2)}
\zeta_{j_3}^{(i_3)}+
$$
$$
+
{\bf 1}_{\{i_1=i_2\}}
{\bf 1}_{\{j_1=j_2\}}
{\bf 1}_{\{i_3=i_4\}}
{\bf 1}_{\{j_3=j_4\}}\zeta_{j_5}^{(i_5)}+
{\bf 1}_{\{i_1=i_2\}}
{\bf 1}_{\{j_1=j_2\}}
{\bf 1}_{\{i_3=i_5\}}
{\bf 1}_{\{j_3=j_5\}}\zeta_{j_4}^{(i_4)}+
$$
$$
+
{\bf 1}_{\{i_1=i_2\}}
{\bf 1}_{\{j_1=j_2\}}
{\bf 1}_{\{i_4=i_5\}}
{\bf 1}_{\{j_4=j_5\}}\zeta_{j_3}^{(i_3)}+
{\bf 1}_{\{i_1=i_3\}}
{\bf 1}_{\{j_1=j_3\}}
{\bf 1}_{\{i_2=i_4\}}
{\bf 1}_{\{j_2=j_4\}}\zeta_{j_5}^{(i_5)}+
$$
$$
+
{\bf 1}_{\{i_1=i_3\}}
{\bf 1}_{\{j_1=j_3\}}
{\bf 1}_{\{i_2=i_5\}}
{\bf 1}_{\{j_2=j_5\}}\zeta_{j_4}^{(i_4)}+
{\bf 1}_{\{i_1=i_3\}}
{\bf 1}_{\{j_1=j_3\}}
{\bf 1}_{\{i_4=i_5\}}
{\bf 1}_{\{j_4=j_5\}}\zeta_{j_2}^{(i_2)}+
$$
$$
+
{\bf 1}_{\{i_1=i_4\}}
{\bf 1}_{\{j_1=j_4\}}
{\bf 1}_{\{i_2=i_3\}}
{\bf 1}_{\{j_2=j_3\}}\zeta_{j_5}^{(i_5)}+
{\bf 1}_{\{i_1=i_4\}}
{\bf 1}_{\{j_1=j_4\}}
{\bf 1}_{\{i_2=i_5\}}
{\bf 1}_{\{j_2=j_5\}}\zeta_{j_3}^{(i_3)}+
$$
$$
+
{\bf 1}_{\{i_1=i_4\}}
{\bf 1}_{\{j_1=j_4\}}
{\bf 1}_{\{i_3=i_5\}}
{\bf 1}_{\{j_3=j_5\}}\zeta_{j_2}^{(i_2)}+
{\bf 1}_{\{i_1=i_5\}}
{\bf 1}_{\{j_1=j_5\}}
{\bf 1}_{\{i_2=i_3\}}
{\bf 1}_{\{j_2=j_3\}}\zeta_{j_4}^{(i_4)}+
$$
$$
+
{\bf 1}_{\{i_1=i_5\}}
{\bf 1}_{\{j_1=j_5\}}
{\bf 1}_{\{i_2=i_4\}}
{\bf 1}_{\{j_2=j_4\}}\zeta_{j_3}^{(i_3)}+
{\bf 1}_{\{i_1=i_5\}}
{\bf 1}_{\{j_1=j_5\}}
{\bf 1}_{\{i_3=i_4\}}
{\bf 1}_{\{j_3=j_4\}}\zeta_{j_2}^{(i_2)}+
$$
$$
+
{\bf 1}_{\{i_2=i_3\}}
{\bf 1}_{\{j_2=j_3\}}
{\bf 1}_{\{i_4=i_5\}}
{\bf 1}_{\{j_4=j_5\}}\zeta_{j_1}^{(i_1)}+
{\bf 1}_{\{i_2=i_4\}}
{\bf 1}_{\{j_2=j_4\}}
{\bf 1}_{\{i_3=i_5\}}
{\bf 1}_{\{j_3=j_5\}}\zeta_{j_1}^{(i_1)}+
$$
\begin{equation}
\label{sss4}
+\Biggl.
{\bf 1}_{\{i_2=i_5\}}
{\bf 1}_{\{j_2=j_5\}}
{\bf 1}_{\{i_3=i_4\}}
{\bf 1}_{\{j_3=j_4\}}\zeta_{j_1}^{(i_1)}\Biggr),
\end{equation}

\vspace{4mm}
$$         
I_{(00000)T,t}^{(i_1i_1i_1i_1i_1)}=
\frac{1}{120}(T-t)^{5/2}
\left(\left(\zeta_0^{(i_1)}\right)^5-
10\left(\zeta_0^{(i_1)}\right)^3+15\zeta_0^{(i_1)}\right)\ \ \ 
\hbox{w.~p.~1},
$$

\newpage
\noindent
$$
I_{(00000)T,t}^{*(i_1i_1i_1i_1i_1)}=
\frac{1}{120}(T-t)^{5/2}\left(\zeta_0^{(i_1)}\right)^5\ \ \ \hbox{w.~p.~1},
$$

\vspace{2mm}
\noindent
where

\vspace{-2mm}
$$
C_{j_5j_4 j_3 j_2 j_1}=
\frac{\sqrt{(2j_1+1)(2j_2+1)(2j_3+1)(2j_4+1)(2j_5+1)}}{32}(T-t)^{5/2}\bar
C_{j_5j_4 j_3 j_2 j_1},
$$

\vspace{-2mm}
$$
\bar C_{j_5j_4 j_3 j_2 j_1}=
\int\limits_{-1}^{1}P_{j_5}(v)
\int\limits_{-1}^{v}P_{j_4}(u)
\int\limits_{-1}^{u}P_{j_3}(z)
\int\limits_{-1}^{z}P_{j_2}(y)
\int\limits_{-1}^{y}
P_{j_1}(x)dx dy dz du dv;
$$

\vspace{5mm}

$$
I_{(0001)T,t}^{*(i_1i_2i_3)}
=\hbox{\vtop{\offinterlineskip\halign{
\hfil#\hfil\cr
{\rm l.i.m.}\cr
$\stackrel{}{{}_{p\to \infty}}$\cr
}} }
\sum_{j_1,j_2,j_3,j_4=0}^{p}
C_{j_4j_3 j_2 j_1}^{0001}
\zeta_{j_1}^{(i_1)}\zeta_{j_2}^{(i_2)}\zeta_{j_3}^{(i_3)}\zeta_{j_4}^{(i_4)},
$$

$$
I_{(0010)T,t}^{*(i_1i_2i_3)}
=\hbox{\vtop{\offinterlineskip\halign{
\hfil#\hfil\cr
{\rm l.i.m.}\cr
$\stackrel{}{{}_{p\to \infty}}$\cr
}} }
\sum_{j_1,j_2,j_3,j_4=0}^{p}
C_{j_4j_3 j_2 j_1}^{0010}
\zeta_{j_1}^{(i_1)}\zeta_{j_2}^{(i_2)}\zeta_{j_3}^{(i_3)}\zeta_{j_4}^{(i_4)},
$$

$$
I_{(0100)T,t}^{*(i_1i_2i_3)}
=\hbox{\vtop{\offinterlineskip\halign{
\hfil#\hfil\cr
{\rm l.i.m.}\cr
$\stackrel{}{{}_{p\to \infty}}$\cr
}} }
\sum_{j_1,j_2,j_3,j_4=0}^{p}
C_{j_4j_3 j_2 j_1}^{0100}
\zeta_{j_1}^{(i_1)}\zeta_{j_2}^{(i_2)}\zeta_{j_3}^{(i_3)}\zeta_{j_4}^{(i_4)},
$$

$$
I_{(1000)T,t}^{*(i_1i_2i_3)}
=\hbox{\vtop{\offinterlineskip\halign{
\hfil#\hfil\cr
{\rm l.i.m.}\cr
$\stackrel{}{{}_{p\to \infty}}$\cr
}} }
\sum_{j_1,j_2,j_3,j_4=0}^{p}
C_{j_4j_3 j_2 j_1}^{1000}
\zeta_{j_1}^{(i_1)}\zeta_{j_2}^{(i_2)}\zeta_{j_3}^{(i_3)}\zeta_{j_4}^{(i_4)},
$$

\vspace{6mm}

$$
I_{(0001)T,t}^{(i_1 i_2 i_3 i_4)}
=\hbox{\vtop{\offinterlineskip\halign{
\hfil#\hfil\cr
{\rm l.i.m.}\cr
$\stackrel{}{{}_{p\to \infty}}$\cr
}} }
\sum_{j_1,j_2,j_3,j_4=0}^{p}
C_{j_4 j_3 j_2 j_1}^{0001}\Biggl(
\zeta_{j_1}^{(i_1)}\zeta_{j_2}^{(i_2)}\zeta_{j_3}^{(i_3)}\zeta_{j_4}^{(i_4)}
-\Biggr.
$$
$$
-
{\bf 1}_{\{i_1=i_2\}}
{\bf 1}_{\{j_1=j_2\}}
\zeta_{j_3}^{(i_3)}
\zeta_{j_4}^{(i_4)}
-
{\bf 1}_{\{i_1=i_3\}}
{\bf 1}_{\{j_1=j_3\}}
\zeta_{j_2}^{(i_2)}
\zeta_{j_4}^{(i_4)}-
$$
$$
-
{\bf 1}_{\{i_1=i_4\}}
{\bf 1}_{\{j_1=j_4\}}
\zeta_{j_2}^{(i_2)}
\zeta_{j_3}^{(i_3)}
-
{\bf 1}_{\{i_2=i_3\}}
{\bf 1}_{\{j_2=j_3\}}
\zeta_{j_1}^{(i_1)}
\zeta_{j_4}^{(i_4)}-
$$
$$
-
{\bf 1}_{\{i_2=i_4\}}
{\bf 1}_{\{j_2=j_4\}}
\zeta_{j_1}^{(i_1)}
\zeta_{j_3}^{(i_3)}
-
{\bf 1}_{\{i_3=i_4\}}
{\bf 1}_{\{j_3=j_4\}}
\zeta_{j_1}^{(i_1)}
\zeta_{j_2}^{(i_2)}+
$$
$$
+
{\bf 1}_{\{i_1=i_2\}}
{\bf 1}_{\{j_1=j_2\}}
{\bf 1}_{\{i_3=i_4\}}
{\bf 1}_{\{j_3=j_4\}}+
{\bf 1}_{\{i_1=i_3\}}
{\bf 1}_{\{j_1=j_3\}}
{\bf 1}_{\{i_2=i_4\}}
{\bf 1}_{\{j_2=j_4\}}+
$$
$$
+\Biggl.
{\bf 1}_{\{i_1=i_4\}}
{\bf 1}_{\{j_1=j_4\}}
{\bf 1}_{\{i_2=i_3\}}
{\bf 1}_{\{j_2=j_3\}}\Biggr),
$$

\newpage
\noindent
$$
I_{(0010)T,t}^{(i_1 i_2 i_3 i_4)}
=\hbox{\vtop{\offinterlineskip\halign{
\hfil#\hfil\cr
{\rm l.i.m.}\cr
$\stackrel{}{{}_{p\to \infty}}$\cr
}} }
\sum_{j_1,j_2,j_3,j_4=0}^{p}
C_{j_4 j_3 j_2 j_1}^{0010}\Biggl(
\zeta_{j_1}^{(i_1)}\zeta_{j_2}^{(i_2)}\zeta_{j_3}^{(i_3)}\zeta_{j_4}^{(i_4)}
-\Biggr.
$$
$$
-
{\bf 1}_{\{i_1=i_2\}}
{\bf 1}_{\{j_1=j_2\}}
\zeta_{j_3}^{(i_3)}
\zeta_{j_4}^{(i_4)}
-
{\bf 1}_{\{i_1=i_3\}}
{\bf 1}_{\{j_1=j_3\}}
\zeta_{j_2}^{(i_2)}
\zeta_{j_4}^{(i_4)}-
$$
$$
-
{\bf 1}_{\{i_1=i_4\}}
{\bf 1}_{\{j_1=j_4\}}
\zeta_{j_2}^{(i_2)}
\zeta_{j_3}^{(i_3)}
-
{\bf 1}_{\{i_2=i_3\}}
{\bf 1}_{\{j_2=j_3\}}
\zeta_{j_1}^{(i_1)}
\zeta_{j_4}^{(i_4)}-
$$
$$
-
{\bf 1}_{\{i_2=i_4\}}
{\bf 1}_{\{j_2=j_4\}}
\zeta_{j_1}^{(i_1)}
\zeta_{j_3}^{(i_3)}
-
{\bf 1}_{\{i_3=i_4\}}
{\bf 1}_{\{j_3=j_4\}}
\zeta_{j_1}^{(i_1)}
\zeta_{j_2}^{(i_2)}+
$$
$$
+
{\bf 1}_{\{i_1=i_2\}}
{\bf 1}_{\{j_1=j_2\}}
{\bf 1}_{\{i_3=i_4\}}
{\bf 1}_{\{j_3=j_4\}}+
{\bf 1}_{\{i_1=i_3\}}
{\bf 1}_{\{j_1=j_3\}}
{\bf 1}_{\{i_2=i_4\}}
{\bf 1}_{\{j_2=j_4\}}+
$$
$$
+\Biggl.
{\bf 1}_{\{i_1=i_4\}}
{\bf 1}_{\{j_1=j_4\}}
{\bf 1}_{\{i_2=i_3\}}
{\bf 1}_{\{j_2=j_3\}}\Biggr),
$$

\vspace{4mm}

$$
I_{(0100)T,t}^{(i_1 i_2 i_3 i_4)}
=\hbox{\vtop{\offinterlineskip\halign{
\hfil#\hfil\cr
{\rm l.i.m.}\cr
$\stackrel{}{{}_{p\to \infty}}$\cr
}} }
\sum_{j_1,j_2,j_3,j_4=0}^{p}
C_{j_4 j_3 j_2 j_1}^{0100}\Biggl(
\zeta_{j_1}^{(i_1)}\zeta_{j_2}^{(i_2)}\zeta_{j_3}^{(i_3)}\zeta_{j_4}^{(i_4)}
-\Biggr.
$$
$$
-
{\bf 1}_{\{i_1=i_2\}}
{\bf 1}_{\{j_1=j_2\}}
\zeta_{j_3}^{(i_3)}
\zeta_{j_4}^{(i_4)}
-
{\bf 1}_{\{i_1=i_3\}}
{\bf 1}_{\{j_1=j_3\}}
\zeta_{j_2}^{(i_2)}
\zeta_{j_4}^{(i_4)}-
$$
$$
-
{\bf 1}_{\{i_1=i_4\}}
{\bf 1}_{\{j_1=j_4\}}
\zeta_{j_2}^{(i_2)}
\zeta_{j_3}^{(i_3)}
-
{\bf 1}_{\{i_2=i_3\}}
{\bf 1}_{\{j_2=j_3\}}
\zeta_{j_1}^{(i_1)}
\zeta_{j_4}^{(i_4)}-
$$
$$
-
{\bf 1}_{\{i_2=i_4\}}
{\bf 1}_{\{j_2=j_4\}}
\zeta_{j_1}^{(i_1)}
\zeta_{j_3}^{(i_3)}
-
{\bf 1}_{\{i_3=i_4\}}
{\bf 1}_{\{j_3=j_4\}}
\zeta_{j_1}^{(i_1)}
\zeta_{j_2}^{(i_2)}+
$$
$$
+
{\bf 1}_{\{i_1=i_2\}}
{\bf 1}_{\{j_1=j_2\}}
{\bf 1}_{\{i_3=i_4\}}
{\bf 1}_{\{j_3=j_4\}}+
{\bf 1}_{\{i_1=i_3\}}
{\bf 1}_{\{j_1=j_3\}}
{\bf 1}_{\{i_2=i_4\}}
{\bf 1}_{\{j_2=j_4\}}+
$$
$$
+\Biggl.
{\bf 1}_{\{i_1=i_4\}}
{\bf 1}_{\{j_1=j_4\}}
{\bf 1}_{\{i_2=i_3\}}
{\bf 1}_{\{j_2=j_3\}}\Biggr),
$$

\vspace{4mm}

$$
I_{(1000)T,t}^{(i_1 i_2 i_3 i_4)}
=\hbox{\vtop{\offinterlineskip\halign{
\hfil#\hfil\cr
{\rm l.i.m.}\cr
$\stackrel{}{{}_{p\to \infty}}$\cr
}} }
\sum_{j_1,j_2,j_3,j_4=0}^{p}
C_{j_4 j_3 j_2 j_1}^{1000}\Biggl(
\zeta_{j_1}^{(i_1)}\zeta_{j_2}^{(i_2)}\zeta_{j_3}^{(i_3)}\zeta_{j_4}^{(i_4)}
-\Biggr.
$$
$$
-
{\bf 1}_{\{i_1=i_2\}}
{\bf 1}_{\{j_1=j_2\}}
\zeta_{j_3}^{(i_3)}
\zeta_{j_4}^{(i_4)}
-
{\bf 1}_{\{i_1=i_3\}}
{\bf 1}_{\{j_1=j_3\}}
\zeta_{j_2}^{(i_2)}
\zeta_{j_4}^{(i_4)}-
$$
$$
-
{\bf 1}_{\{i_1=i_4\}}
{\bf 1}_{\{j_1=j_4\}}
\zeta_{j_2}^{(i_2)}
\zeta_{j_3}^{(i_3)}
-
{\bf 1}_{\{i_2=i_3\}}
{\bf 1}_{\{j_2=j_3\}}
\zeta_{j_1}^{(i_1)}
\zeta_{j_4}^{(i_4)}-
$$
$$
-
{\bf 1}_{\{i_2=i_4\}}
{\bf 1}_{\{j_2=j_4\}}
\zeta_{j_1}^{(i_1)}
\zeta_{j_3}^{(i_3)}
-
{\bf 1}_{\{i_3=i_4\}}
{\bf 1}_{\{j_3=j_4\}}
\zeta_{j_1}^{(i_1)}
\zeta_{j_2}^{(i_2)}+
$$
$$
+
{\bf 1}_{\{i_1=i_2\}}
{\bf 1}_{\{j_1=j_2\}}
{\bf 1}_{\{i_3=i_4\}}
{\bf 1}_{\{j_3=j_4\}}+
{\bf 1}_{\{i_1=i_3\}}
{\bf 1}_{\{j_1=j_3\}}
{\bf 1}_{\{i_2=i_4\}}
{\bf 1}_{\{j_2=j_4\}}+
$$
$$
+\Biggl.
{\bf 1}_{\{i_1=i_4\}}
{\bf 1}_{\{j_1=j_4\}}
{\bf 1}_{\{i_2=i_3\}}
{\bf 1}_{\{j_2=j_3\}}\Biggr),
$$

\noindent
where
\newpage
\noindent
$$
C_{j_4j_3j_2j_1}^{0001}
=\frac{\sqrt{(2j_1+1)(2j_2+1)(2j_3+1)(2j_4+1)}}{32}(T-t)^{3}\bar
C_{j_4j_3j_2j_1}^{0001},
$$

\vspace{-2mm}
$$
C_{j_3j_2j_1}^{0010}
=\frac{\sqrt{(2j_1+1)(2j_2+1)(2j_3+1)(2j_4+1)}}{32}(T-t)^{3}\bar
C_{j_4j_3j_2j_1}^{0010},
$$

\vspace{-2mm}
$$
C_{j_4j_3j_2j_1}^{0100}=
\frac{\sqrt{(2j_1+1)(2j_2+1)(2j_3+1)(2j_4+1)}}{32}(T-t)^{3}\bar
C_{j_3j_2j_1}^{0100},
$$

\vspace{-2mm}
$$
C_{j_4j_3j_2j_1}^{1000}
=\frac{\sqrt{(2j_1+1)(2j_2+1)(2j_3+1)(2j_4+1)}}{32}(T-t)^{3}\bar
C_{j_4j_3j_2j_1}^{1000},
$$

\vspace{-2mm}
$$
\bar C_{j_4j_3j_2j_1}^{1000}=-
\int\limits_{-1}^{1}P_{j_4}(u)
\int\limits_{-1}^{u}P_{j_3}(z)
\int\limits_{-1}^{z}P_{j_2}(y)
\int\limits_{-1}^{y}
P_{j_1}(x)(x+1)dx dy dz du,
$$
$$
\bar C_{j_4j_3j_2j_1}^{0100}=-
\int\limits_{-1}^{1}P_{j_4}(u)
\int\limits_{-1}^{u}P_{j_3}(z)
\int\limits_{-1}^{z}P_{j_2}(y)(y+1)
\int\limits_{-1}^{y}
P_{j_1}(x)dx dy dz du,
$$
$$
\bar C_{j_4j_3j_2j_1}^{0010}=-
\int\limits_{-1}^{1}P_{j_4}(u)
\int\limits_{-1}^{u}P_{j_3}(z)(z+1)
\int\limits_{-1}^{z}P_{j_2}(y)
\int\limits_{-1}^{y}
P_{j_1}(x)dx dy dz du,
$$
$$
\bar C_{j_4j_3j_2j_1}^{0001}=-
\int\limits_{-1}^{1}P_{j_4}(u)(u+1)
\int\limits_{-1}^{u}P_{j_3}(z)
\int\limits_{-1}^{z}P_{j_2}(y)
\int\limits_{-1}^{y}
P_{j_1}(x)dx dy dz du;
$$

\vspace{3mm}

$$
I_{(000000)T,t}^{*(i_1 i_2 i_3 i_4 i_5 i_6)}=
\hbox{\vtop{\offinterlineskip\halign{
\hfil#\hfil\cr
{\rm l.i.m.}\cr
$\stackrel{}{{}_{p\to \infty}}$\cr
}} }
\sum\limits_{j_1, j_2, j_3, j_4, j_5, j_6=0}^{p}
C_{j_6j_5j_4 j_3 j_2 j_1}
\zeta_{j_1}^{(i_1)}\zeta_{j_2}^{(i_2)}\zeta_{j_3}^{(i_3)}
\zeta_{j_4}^{(i_4)}\zeta_{j_5}^{(i_5)}\zeta_{j_6}^{(i_6)},
$$

\vspace{4mm}

$$
I_{(000000)T,t}^{(i_1 i_2 i_3 i_4 i_5 i_6)}
=\hbox{\vtop{\offinterlineskip\halign{
\hfil#\hfil\cr
{\rm l.i.m.}\cr
$\stackrel{}{{}_{p\to \infty}}$\cr
}} }\sum_{j_1,j_2,j_3,j_4,j_5,j_6=0}^{p}
C_{j_6 j_5 j_4 j_3 j_2 j_1}\Biggl(
\prod_{l=1}^6
\zeta_{j_l}^{(i_l)}
-\Biggr.
$$
$$
-
{\bf 1}_{\{j_1=j_6\}}
{\bf 1}_{\{i_1=i_6\}}
\zeta_{j_2}^{(i_2)}
\zeta_{j_3}^{(i_3)}
\zeta_{j_4}^{(i_4)}
\zeta_{j_5}^{(i_5)}-
{\bf 1}_{\{j_2=j_6\}}
{\bf 1}_{\{i_2=i_6\}}
\zeta_{j_1}^{(i_1)}
\zeta_{j_3}^{(i_3)}
\zeta_{j_4}^{(i_4)}
\zeta_{j_5}^{(i_5)}-
$$
$$
-
{\bf 1}_{\{j_3=j_6\}}
{\bf 1}_{\{i_3=i_6\}}
\zeta_{j_1}^{(i_1)}
\zeta_{j_2}^{(i_2)}
\zeta_{j_4}^{(i_4)}
\zeta_{j_5}^{(i_5)}-
{\bf 1}_{\{j_4=j_6\}}
{\bf 1}_{\{i_4=i_6\}}
\zeta_{j_1}^{(i_1)}
\zeta_{j_2}^{(i_2)}
\zeta_{j_3}^{(i_3)}
\zeta_{j_5}^{(i_5)}-
$$
$$
-
{\bf 1}_{\{j_5=j_6\}}
{\bf 1}_{\{i_5=i_6\}}
\zeta_{j_1}^{(i_1)}
\zeta_{j_2}^{(i_2)}
\zeta_{j_3}^{(i_3)}
\zeta_{j_4}^{(i_4)}-
{\bf 1}_{\{j_1=j_2\}}
{\bf 1}_{\{i_1=i_2\}}
\zeta_{j_3}^{(i_3)}
\zeta_{j_4}^{(i_4)}
\zeta_{j_5}^{(i_5)}
\zeta_{j_6}^{(i_6)}-
$$
$$
-
{\bf 1}_{\{j_1=j_3\}}
{\bf 1}_{\{i_1=i_3\}}
\zeta_{j_2}^{(i_2)}
\zeta_{j_4}^{(i_4)}
\zeta_{j_5}^{(i_5)}
\zeta_{j_6}^{(i_6)}-
{\bf 1}_{\{j_1=j_4\}}
{\bf 1}_{\{i_1=i_4\}}
\zeta_{j_2}^{(i_2)}
\zeta_{j_3}^{(i_3)}
\zeta_{j_5}^{(i_5)}
\zeta_{j_6}^{(i_6)}-
$$
$$
-
{\bf 1}_{\{j_1=j_5\}}
{\bf 1}_{\{i_1=i_5\}}
\zeta_{j_2}^{(i_2)}
\zeta_{j_3}^{(i_3)}
\zeta_{j_4}^{(i_4)}
\zeta_{j_6}^{(i_6)}-
{\bf 1}_{\{j_2=j_3\}}
{\bf 1}_{\{i_2=i_3\}}
\zeta_{j_1}^{(i_1)}
\zeta_{j_4}^{(i_4)}
\zeta_{j_5}^{(i_5)}
\zeta_{j_6}^{(i_6)}-
$$
$$
-
{\bf 1}_{\{j_2=j_4\}}
{\bf 1}_{\{i_2=i_4\}}
\zeta_{j_1}^{(i_1)}
\zeta_{j_3}^{(i_3)}
\zeta_{j_5}^{(i_5)}
\zeta_{j_6}^{(i_6)}-
{\bf 1}_{\{j_2=j_5\}}
{\bf 1}_{\{i_2=i_5\}}
\zeta_{j_1}^{(i_1)}
\zeta_{j_3}^{(i_3)}
\zeta_{j_4}^{(i_4)}
\zeta_{j_6}^{(i_6)}-
$$
$$
-
{\bf 1}_{\{j_3=j_4\}}
{\bf 1}_{\{i_3=i_4\}}
\zeta_{j_1}^{(i_1)}
\zeta_{j_2}^{(i_2)}
\zeta_{j_5}^{(i_5)}
\zeta_{j_6}^{(i_6)}-
{\bf 1}_{\{j_3=j_5\}}
{\bf 1}_{\{i_3=i_5\}}
\zeta_{j_1}^{(i_1)}
\zeta_{j_2}^{(i_2)}
\zeta_{j_4}^{(i_4)}
\zeta_{j_6}^{(i_6)}-
$$
$$
-
{\bf 1}_{\{j_4=j_5\}}
{\bf 1}_{\{i_4=i_5\}}
\zeta_{j_1}^{(i_1)}
\zeta_{j_2}^{(i_2)}
\zeta_{j_3}^{(i_3)}
\zeta_{j_6}^{(i_6)}+
$$
$$
+
{\bf 1}_{\{j_1=j_2\}}
{\bf 1}_{\{i_1=i_2\}}
{\bf 1}_{\{j_3=j_4\}}
{\bf 1}_{\{i_3=i_4\}}
\zeta_{j_5}^{(i_5)}
\zeta_{j_6}^{(i_6)}
+
{\bf 1}_{\{j_1=j_2\}}
{\bf 1}_{\{i_1=i_2\}}
{\bf 1}_{\{j_3=j_5\}}
{\bf 1}_{\{i_3=i_5\}}
\zeta_{j_4}^{(i_4)}
\zeta_{j_6}^{(i_6)}+
$$
$$
+
{\bf 1}_{\{j_1=j_2\}}
{\bf 1}_{\{i_1=i_2\}}
{\bf 1}_{\{j_4=j_5\}}
{\bf 1}_{\{i_4=i_5\}}
\zeta_{j_3}^{(i_3)}
\zeta_{j_6}^{(i_6)}
+
{\bf 1}_{\{j_1=j_3\}}
{\bf 1}_{\{i_1=i_3\}}
{\bf 1}_{\{j_2=j_4\}}
{\bf 1}_{\{i_2=i_4\}}
\zeta_{j_5}^{(i_5)}
\zeta_{j_6}^{(i_6)}+
$$
$$
+
{\bf 1}_{\{j_1=j_3\}}
{\bf 1}_{\{i_1=i_3\}}
{\bf 1}_{\{j_2=j_5\}}
{\bf 1}_{\{i_2=i_5\}}
\zeta_{j_4}^{(i_4)}
\zeta_{j_6}^{(i_6)}
+
{\bf 1}_{\{j_1=j_3\}}
{\bf 1}_{\{i_1=i_3\}}
{\bf 1}_{\{j_4=j_5\}}
{\bf 1}_{\{i_4=i_5\}}
\zeta_{j_2}^{(i_2)}
\zeta_{j_6}^{(i_6)}+
$$
$$
+
{\bf 1}_{\{j_1=j_4\}}
{\bf 1}_{\{i_1=i_4\}}
{\bf 1}_{\{j_2=j_3\}}
{\bf 1}_{\{i_2=i_3\}}
\zeta_{j_5}^{(i_5)}
\zeta_{j_6}^{(i_6)}
+
{\bf 1}_{\{j_1=j_4\}}
{\bf 1}_{\{i_1=i_4\}}
{\bf 1}_{\{j_2=j_5\}}
{\bf 1}_{\{i_2=i_5\}}
\zeta_{j_3}^{(i_3)}
\zeta_{j_6}^{(i_6)}+
$$
$$
+
{\bf 1}_{\{j_1=j_4\}}
{\bf 1}_{\{i_1=i_4\}}
{\bf 1}_{\{j_3=j_5\}}
{\bf 1}_{\{i_3=i_5\}}
\zeta_{j_2}^{(i_2)}
\zeta_{j_6}^{(i_6)}
+
{\bf 1}_{\{j_1=j_5\}}
{\bf 1}_{\{i_1=i_5\}}
{\bf 1}_{\{j_2=j_3\}}
{\bf 1}_{\{i_2=i_3\}}
\zeta_{j_4}^{(i_4)}
\zeta_{j_6}^{(i_6)}+
$$
$$
+
{\bf 1}_{\{j_1=j_5\}}
{\bf 1}_{\{i_1=i_5\}}
{\bf 1}_{\{j_2=j_4\}}
{\bf 1}_{\{i_2=i_4\}}
\zeta_{j_3}^{(i_3)}
\zeta_{j_6}^{(i_6)}
+
{\bf 1}_{\{j_1=j_5\}}
{\bf 1}_{\{i_1=i_5\}}
{\bf 1}_{\{j_3=j_4\}}
{\bf 1}_{\{i_3=i_4\}}
\zeta_{j_2}^{(i_2)}
\zeta_{j_6}^{(i_6)}+
$$
$$
+
{\bf 1}_{\{j_2=j_3\}}
{\bf 1}_{\{i_2=i_3\}}
{\bf 1}_{\{j_4=j_5\}}
{\bf 1}_{\{i_4=i_5\}}
\zeta_{j_1}^{(i_1)}
\zeta_{j_6}^{(i_6)}
+
{\bf 1}_{\{j_2=j_4\}}
{\bf 1}_{\{i_2=i_4\}}
{\bf 1}_{\{j_3=j_5\}}
{\bf 1}_{\{i_3=i_5\}}
\zeta_{j_1}^{(i_1)}
\zeta_{j_6}^{(i_6)}+
$$
$$
+
{\bf 1}_{\{j_2=j_5\}}
{\bf 1}_{\{i_2=i_5\}}
{\bf 1}_{\{j_3=j_4\}}
{\bf 1}_{\{i_3=i_4\}}
\zeta_{j_1}^{(i_1)}
\zeta_{j_6}^{(i_6)}
+
{\bf 1}_{\{j_6=j_1\}}
{\bf 1}_{\{i_6=i_1\}}
{\bf 1}_{\{j_3=j_4\}}
{\bf 1}_{\{i_3=i_4\}}
\zeta_{j_2}^{(i_2)}
\zeta_{j_5}^{(i_5)}+
$$
$$
+
{\bf 1}_{\{j_6=j_1\}}
{\bf 1}_{\{i_6=i_1\}}
{\bf 1}_{\{j_3=j_5\}}
{\bf 1}_{\{i_3=i_5\}}
\zeta_{j_2}^{(i_2)}
\zeta_{j_4}^{(i_4)}
+
{\bf 1}_{\{j_6=j_1\}}
{\bf 1}_{\{i_6=i_1\}}
{\bf 1}_{\{j_2=j_5\}}
{\bf 1}_{\{i_2=i_5\}}
\zeta_{j_3}^{(i_3)}
\zeta_{j_4}^{(i_4)}+
$$
$$
+
{\bf 1}_{\{j_6=j_1\}}
{\bf 1}_{\{i_6=i_1\}}
{\bf 1}_{\{j_2=j_4\}}
{\bf 1}_{\{i_2=i_4\}}
\zeta_{j_3}^{(i_3)}
\zeta_{j_5}^{(i_5)}
+
{\bf 1}_{\{j_6=j_1\}}
{\bf 1}_{\{i_6=i_1\}}
{\bf 1}_{\{j_4=j_5\}}
{\bf 1}_{\{i_4=i_5\}}
\zeta_{j_2}^{(i_2)}
\zeta_{j_3}^{(i_3)}+
$$
$$
+
{\bf 1}_{\{j_6=j_1\}}
{\bf 1}_{\{i_6=i_1\}}
{\bf 1}_{\{j_2=j_3\}}
{\bf 1}_{\{i_2=i_3\}}
\zeta_{j_4}^{(i_4)}
\zeta_{j_5}^{(i_5)}
+
{\bf 1}_{\{j_6=j_2\}}
{\bf 1}_{\{i_6=i_2\}}
{\bf 1}_{\{j_3=j_5\}}
{\bf 1}_{\{i_3=i_5\}}
\zeta_{j_1}^{(i_1)}
\zeta_{j_4}^{(i_4)}+
$$
$$
+
{\bf 1}_{\{j_6=j_2\}}
{\bf 1}_{\{i_6=i_2\}}
{\bf 1}_{\{j_4=j_5\}}
{\bf 1}_{\{i_4=i_5\}}
\zeta_{j_1}^{(i_1)}
\zeta_{j_3}^{(i_3)}
+
{\bf 1}_{\{j_6=j_2\}}
{\bf 1}_{\{i_6=i_2\}}
{\bf 1}_{\{j_3=j_4\}}
{\bf 1}_{\{i_3=i_4\}}
\zeta_{j_1}^{(i_1)}
\zeta_{j_5}^{(i_5)}+
$$
$$
+
{\bf 1}_{\{j_6=j_2\}}
{\bf 1}_{\{i_6=i_2\}}
{\bf 1}_{\{j_1=j_5\}}
{\bf 1}_{\{i_1=i_5\}}
\zeta_{j_3}^{(i_3)}
\zeta_{j_4}^{(i_4)}
+
{\bf 1}_{\{j_6=j_2\}}
{\bf 1}_{\{i_6=i_2\}}
{\bf 1}_{\{j_1=j_4\}}
{\bf 1}_{\{i_1=i_4\}}
\zeta_{j_3}^{(i_3)}
\zeta_{j_5}^{(i_5)}+
$$
$$
+
{\bf 1}_{\{j_6=j_2\}}
{\bf 1}_{\{i_6=i_2\}}
{\bf 1}_{\{j_1=j_3\}}
{\bf 1}_{\{i_1=i_3\}}
\zeta_{j_4}^{(i_4)}
\zeta_{j_5}^{(i_5)}
+
{\bf 1}_{\{j_6=j_3\}}
{\bf 1}_{\{i_6=i_3\}}
{\bf 1}_{\{j_2=j_5\}}
{\bf 1}_{\{i_2=i_5\}}
\zeta_{j_1}^{(i_1)}
\zeta_{j_4}^{(i_4)}+
$$
$$
+
{\bf 1}_{\{j_6=j_3\}}
{\bf 1}_{\{i_6=i_3\}}
{\bf 1}_{\{j_4=j_5\}}
{\bf 1}_{\{i_4=i_5\}}
\zeta_{j_1}^{(i_1)}
\zeta_{j_2}^{(i_2)}
+
{\bf 1}_{\{j_6=j_3\}}
{\bf 1}_{\{i_6=i_3\}}
{\bf 1}_{\{j_2=j_4\}}
{\bf 1}_{\{i_2=i_4\}}
\zeta_{j_1}^{(i_1)}
\zeta_{j_5}^{(i_5)}+
$$
$$
+
{\bf 1}_{\{j_6=j_3\}}
{\bf 1}_{\{i_6=i_3\}}
{\bf 1}_{\{j_1=j_5\}}
{\bf 1}_{\{i_1=i_5\}}
\zeta_{j_2}^{(i_2)}
\zeta_{j_4}^{(i_4)}
+
{\bf 1}_{\{j_6=j_3\}}
{\bf 1}_{\{i_6=i_3\}}
{\bf 1}_{\{j_1=j_4\}}
{\bf 1}_{\{i_1=i_4\}}
\zeta_{j_2}^{(i_2)}
\zeta_{j_5}^{(i_5)}+
$$
$$
+
{\bf 1}_{\{j_6=j_3\}}
{\bf 1}_{\{i_6=i_3\}}
{\bf 1}_{\{j_1=j_2\}}
{\bf 1}_{\{i_1=i_2\}}
\zeta_{j_4}^{(i_4)}
\zeta_{j_5}^{(i_5)}
+
{\bf 1}_{\{j_6=j_4\}}
{\bf 1}_{\{i_6=i_4\}}
{\bf 1}_{\{j_3=j_5\}}
{\bf 1}_{\{i_3=i_5\}}
\zeta_{j_1}^{(i_1)}
\zeta_{j_2}^{(i_2)}+
$$
$$
+
{\bf 1}_{\{j_6=j_4\}}
{\bf 1}_{\{i_6=i_4\}}
{\bf 1}_{\{j_2=j_5\}}
{\bf 1}_{\{i_2=i_5\}}
\zeta_{j_1}^{(i_1)}
\zeta_{j_3}^{(i_3)}
+
{\bf 1}_{\{j_6=j_4\}}
{\bf 1}_{\{i_6=i_4\}}
{\bf 1}_{\{j_2=j_3\}}
{\bf 1}_{\{i_2=i_3\}}
\zeta_{j_1}^{(i_1)}
\zeta_{j_5}^{(i_5)}+
$$
$$
+
{\bf 1}_{\{j_6=j_4\}}
{\bf 1}_{\{i_6=i_4\}}
{\bf 1}_{\{j_1=j_5\}}
{\bf 1}_{\{i_1=i_5\}}
\zeta_{j_2}^{(i_2)}
\zeta_{j_3}^{(i_3)}
+
{\bf 1}_{\{j_6=j_4\}}
{\bf 1}_{\{i_6=i_4\}}
{\bf 1}_{\{j_1=j_3\}}
{\bf 1}_{\{i_1=i_3\}}
\zeta_{j_2}^{(i_2)}
\zeta_{j_5}^{(i_5)}+
$$
$$
+
{\bf 1}_{\{j_6=j_4\}}
{\bf 1}_{\{i_6=i_4\}}
{\bf 1}_{\{j_1=j_2\}}
{\bf 1}_{\{i_1=i_2\}}
\zeta_{j_3}^{(i_3)}
\zeta_{j_5}^{(i_5)}
+
{\bf 1}_{\{j_6=j_5\}}
{\bf 1}_{\{i_6=i_5\}}
{\bf 1}_{\{j_3=j_4\}}
{\bf 1}_{\{i_3=i_4\}}
\zeta_{j_1}^{(i_1)}
\zeta_{j_2}^{(i_2)}+
$$
$$
+
{\bf 1}_{\{j_6=j_5\}}
{\bf 1}_{\{i_6=i_5\}}
{\bf 1}_{\{j_2=j_4\}}
{\bf 1}_{\{i_2=i_4\}}
\zeta_{j_1}^{(i_1)}
\zeta_{j_3}^{(i_3)}
+
{\bf 1}_{\{j_6=j_5\}}
{\bf 1}_{\{i_6=i_5\}}
{\bf 1}_{\{j_2=j_3\}}
{\bf 1}_{\{i_2=i_3\}}
\zeta_{j_1}^{(i_1)}
\zeta_{j_4}^{(i_4)}+
$$
$$
+
{\bf 1}_{\{j_6=j_5\}}
{\bf 1}_{\{i_6=i_5\}}
{\bf 1}_{\{j_1=j_4\}}
{\bf 1}_{\{i_1=i_4\}}
\zeta_{j_2}^{(i_2)}
\zeta_{j_3}^{(i_3)}
+
{\bf 1}_{\{j_6=j_5\}}
{\bf 1}_{\{i_6=i_5\}}
{\bf 1}_{\{j_1=j_3\}}
{\bf 1}_{\{i_1=i_3\}}
\zeta_{j_2}^{(i_2)}
\zeta_{j_4}^{(i_4)}+
$$
$$
+
{\bf 1}_{\{j_6=j_5\}}
{\bf 1}_{\{i_6=i_5\}}
{\bf 1}_{\{j_1=j_2\}}
{\bf 1}_{\{i_1=i_2\}}
\zeta_{j_3}^{(i_3)}
\zeta_{j_4}^{(i_4)}-
$$
$$
-
{\bf 1}_{\{j_6=j_1\}}
{\bf 1}_{\{i_6=i_1\}}
{\bf 1}_{\{j_2=j_5\}}
{\bf 1}_{\{i_2=i_5\}}
{\bf 1}_{\{j_3=j_4\}}
{\bf 1}_{\{i_3=i_4\}}-
$$
$$
-
{\bf 1}_{\{j_6=j_1\}}
{\bf 1}_{\{i_6=i_1\}}
{\bf 1}_{\{j_2=j_4\}}
{\bf 1}_{\{i_2=i_4\}}
{\bf 1}_{\{j_3=j_5\}}
{\bf 1}_{\{i_3=i_5\}}-
$$
$$
-
{\bf 1}_{\{j_6=j_1\}}
{\bf 1}_{\{i_6=i_1\}}
{\bf 1}_{\{j_2=j_3\}}
{\bf 1}_{\{i_2=i_3\}}
{\bf 1}_{\{j_4=j_5\}}
{\bf 1}_{\{i_4=i_5\}}-
$$
$$
-               
{\bf 1}_{\{j_6=j_2\}}
{\bf 1}_{\{i_6=i_2\}}
{\bf 1}_{\{j_1=j_5\}}
{\bf 1}_{\{i_1=i_5\}}
{\bf 1}_{\{j_3=j_4\}}
{\bf 1}_{\{i_3=i_4\}}-
$$
$$
-
{\bf 1}_{\{j_6=j_2\}}
{\bf 1}_{\{i_6=i_2\}}
{\bf 1}_{\{j_1=j_4\}}
{\bf 1}_{\{i_1=i_4\}}
{\bf 1}_{\{j_3=j_5\}}
{\bf 1}_{\{i_3=i_5\}}-
$$
$$
-
{\bf 1}_{\{j_6=j_2\}}
{\bf 1}_{\{i_6=i_2\}}
{\bf 1}_{\{j_1=j_3\}}
{\bf 1}_{\{i_1=i_3\}}
{\bf 1}_{\{j_4=j_5\}}
{\bf 1}_{\{i_4=i_5\}}-
$$
$$
-
{\bf 1}_{\{j_6=j_3\}}
{\bf 1}_{\{i_6=i_3\}}
{\bf 1}_{\{j_1=j_5\}}
{\bf 1}_{\{i_1=i_5\}}
{\bf 1}_{\{j_2=j_4\}}
{\bf 1}_{\{i_2=i_4\}}-
$$
$$
-
{\bf 1}_{\{j_6=j_3\}}
{\bf 1}_{\{i_6=i_3\}}
{\bf 1}_{\{j_1=j_4\}}
{\bf 1}_{\{i_1=i_4\}}
{\bf 1}_{\{j_2=j_5\}}
{\bf 1}_{\{i_2=i_5\}}-
$$
$$
-
{\bf 1}_{\{j_3=j_6\}}
{\bf 1}_{\{i_3=i_6\}}
{\bf 1}_{\{j_1=j_2\}}
{\bf 1}_{\{i_1=i_2\}}
{\bf 1}_{\{j_4=j_5\}}
{\bf 1}_{\{i_4=i_5\}}-
$$
$$
-
{\bf 1}_{\{j_6=j_4\}}
{\bf 1}_{\{i_6=i_4\}}
{\bf 1}_{\{j_1=j_5\}}
{\bf 1}_{\{i_1=i_5\}}
{\bf 1}_{\{j_2=j_3\}}
{\bf 1}_{\{i_2=i_3\}}-
$$
$$
-
{\bf 1}_{\{j_6=j_4\}}
{\bf 1}_{\{i_6=i_4\}}
{\bf 1}_{\{j_1=j_3\}}
{\bf 1}_{\{i_1=i_3\}}
{\bf 1}_{\{j_2=j_5\}}
{\bf 1}_{\{i_2=i_5\}}-
$$
$$
-
{\bf 1}_{\{j_6=j_4\}}
{\bf 1}_{\{i_6=i_4\}}
{\bf 1}_{\{j_1=j_2\}}
{\bf 1}_{\{i_1=i_2\}}
{\bf 1}_{\{j_3=j_5\}}
{\bf 1}_{\{i_3=i_5\}}-
$$
$$
-
{\bf 1}_{\{j_6=j_5\}}
{\bf 1}_{\{i_6=i_5\}}
{\bf 1}_{\{j_1=j_4\}}
{\bf 1}_{\{i_1=i_4\}}
{\bf 1}_{\{j_2=j_3\}}
{\bf 1}_{\{i_2=i_3\}}-
$$
$$
-
{\bf 1}_{\{j_6=j_5\}}
{\bf 1}_{\{i_6=i_5\}}
{\bf 1}_{\{j_1=j_2\}}
{\bf 1}_{\{i_1=i_2\}}
{\bf 1}_{\{j_3=j_4\}}
{\bf 1}_{\{i_3=i_4\}}-
$$
$$
\Biggl.-
{\bf 1}_{\{j_6=j_5\}}
{\bf 1}_{\{i_6=i_5\}}
{\bf 1}_{\{j_1=j_3\}}
{\bf 1}_{\{i_1=i_3\}}
{\bf 1}_{\{j_2=j_4\}}
{\bf 1}_{\{i_2=i_4\}}\Biggr),
$$

\vspace{5mm}

$$         
I_{(000000)T,t}^{(i_1i_1i_1i_1i_1i_1)}=
\frac{1}{720}(T-t)^{3}
\left(\left(\zeta_0^{(i_1)}\right)^6-
15\left(\zeta_0^{(i_1)}\right)^4+45\left(\zeta_0^{(i_1)}\right)^2-
15\right)\ \ \ 
\hbox{w.~p.~1},
$$

$$
I_{(000000)T,t}^{*(i_1i_1i_1i_1i_1i_1)}=
\frac{1}{720}(T-t)^{3}\left(\zeta_0^{(i_1)}\right)^6\ \ \ \hbox{w.~p.~1},
$$

\noindent
where

\newpage
\noindent
$$
C_{j_6j_5j_4 j_3 j_2 j_1}
=
$$

\vspace{-8mm}
$$
=\frac{\sqrt{(2j_1+1)(2j_2+1)(2j_3+1)
(2j_4+1)(2j_5+1)(2j_6+1)}}{64}(T-t)^{3}\bar
C_{j_6j_5j_4 j_3 j_2 j_1},
$$
$$
\bar C_{j_6j_5j_4 j_3 j_2 j_1}=
$$
$$
=
\int\limits_{-1}^{1}P_{j_6}(w)
\int\limits_{-1}^{w}P_{j_5}(v)
\int\limits_{-1}^{v}P_{j_4}(u)
\int\limits_{-1}^{u}P_{j_3}(z)
\int\limits_{-1}^{z}P_{j_2}(y)
\int\limits_{-1}^{y}
P_{j_1}(x)dx dy dz du dv dw.
$$

\vspace{2mm}

It should be noted that instead of the expansion (\ref{good1})
we can consider the following expansion, which is derived by direct 
calculation
$$
I_{(000)T,t}^{*(i_1 i_2 i_3)}=-\frac{1}{T-t}\left(
I_{(0)T,t}^{(i_3)}I_{(10)T,t}^{*(i_2 i_1)}+
I_{(0)T,t}^{(i_1)}I_{(10)T,t}^{*(i_2 i_3)}\right)+
\frac{1}{2}I_{(0)T,t}^{(i_3)}\left(
I_{(00)T,t}^{*(i_1 i_2)}-I_{(00)T,t}^{*(i_2 i_1)}\right)-
$$
\begin{equation}
\label{4004.ii}
~~-(T-t)^{3/2}\left(\frac{1}{6}\zeta_0^{(i_1)}\zeta_0^{(i_3)}
\left(\zeta_0^{(i_2)}+\sqrt{3}\zeta_1^{(i_2)}-\frac{1}{\sqrt{5}}
\zeta_2^{(i_2)}\right)
+\frac{1}{4}D^{(i_1 i_2 i_3)}_{T,t}\right),
\end{equation}

\vspace{2mm}
\noindent
where
$$
D^{(i_1 i_2 i_3)}_{T,t}=
\sum\limits_{\stackrel{i=1,\ j=0,\ k=i}{{}_{2i\ge k+i-j\ge -2;\
k+i-j\ - {\rm even}}}}^{\infty}
N_{ijk}K_{i+1,k+1,\frac{k+i-j}{2}+1}\zeta_i^{(i_1)}\zeta_j^{(i_2)}
\zeta_k^{(i_3)}+
$$
$$
+\sum\limits_{\stackrel{i=1,\ j=0,\ 1\le k\le i-1}{{}_{2k\ge k+i-j\ge -2;\
k+i-j\ - {\rm even}}}}^{\infty}
N_{ijk}K_{k+1,i+1,\frac{k+i-j}{2}+1}\zeta_i^{(i_1)}\zeta_j^{(i_2)}
\zeta_k^{(i_3)}-
$$
$$
-\sum\limits_{\stackrel{i=1,\ j=0,\ k=i+2}{{}_{2i+2\ge k+i-j\ge 0;\ k+i-j\ 
- {\rm even}}}}^{\infty}
N_{ijk}K_{i+1,k-1,\frac{k+i-j}{2}}\zeta_i^{(i_1)}\zeta_j^{(i_2)}
\zeta_k^{(i_3)}-
$$
$$
-\sum\limits_{\stackrel{i=1,\ j=0,\ 1\le k\le i+1}
{{}_{2k-2\ge k+i-j\ge 0;\ k+i-j\ - {\rm even}}}}^{\infty}
N_{ijk}K_{k-1,i+1,\frac{k+i-j}{2}}\zeta_i^{(i_1)}\zeta_j^{(i_2)}
\zeta_k^{(i_3)}-
$$
$$
-\sum\limits_{\stackrel{i=1,\ j=0,\ k=i-2, k\ge 1}
{{}_{2i-2\ge k+i-j\ge 0;\ k+i-j\ - {\rm even}}}}^{\infty}
N_{ijk}K_{i-1,k+1,\frac{k+i-j}{2}}\zeta_i^{(i_1)}\zeta_j^{(i_2)}
\zeta_k^{(i_3)}-
$$
$$
-\sum\limits_{\stackrel{i=1,\ j=0,\ 1\le k \le i-3} 
{{}_{2k+2\ge k+i-j\ge 0;\ k+i-j\ - {\rm even}}}}^{\infty}
N_{ijk}K_{k+1,i-1,\frac{k+i-j}{2}}\zeta_i^{(i_1)}\zeta_j^{(i_2)}
\zeta_k^{(i_3)}+
$$
$$
+\sum\limits_{\stackrel{i=1,\ j=0,\ k=i}
{{}_{2i\ge k+i-j\ge 2;\ k+i-j\ - {\rm even}}}}^{\infty}
N_{ijk}K_{i-1,k-1,\frac{k+i-j}{2}-1}\zeta_i^{(i_1)}\zeta_j^{(i_2)}
\zeta_k^{(i_3)}+
$$
$$
+\sum\limits_{\stackrel{i=1,\ j=0\ 1\le k\le i-1}
{{}_{2k\ge k+i-j\ge 2;\ k+i-j\ - {\rm even}}}}^{\infty}
N_{ijk}K_{k-1,i-1,\frac{k+i-j}{2}-1}\zeta_i^{(i_1)}\zeta_j^{(i_2)}
\zeta_k^{(i_3)},
$$

\noindent
where
$$
N_{ijk}=\sqrt{\frac{1}{(2k+1)(2j+1)(2i+1)}}~,
$$
$$
K_{m,n,k}=\frac{a_{m-k}a_k a_{n-k}}{a_{m+n-k}}\cdot
\frac{2n+2m-4k+1}{2n+2m-2k+1},\ \ \ a_k=\frac{(2k-1)!!}{k!},\ \ \ m\le n.
$$

\vspace{2mm}

However, as we will see further the expansion (\ref{zzz1})
is more convenient for the practical implementation then
(\ref{4004.ii}).

Also note the following relation between iterated
It\^{o} and Stratonovich stochastic integrals
$$
I_{(0000)T,t}^{(i_1 i_2 i_3 i_4)}=
I_{(0000)T,t}^{*(i_1 i_2 i_3 i_4)}+
\frac{1}{2}{\bf 1}_{\{i_1=i_2\}}I_{(10)T,t}^{*(i_3 i_4)}-
\frac{1}{2}{\bf 1}_{\{i_2=i_3\}}\biggl(
I_{(10)T,t}^{*(i_1 i_4)}-
I_{(01)T,t}^{*(i_1 i_4)}\biggr)-
$$
$$
-
\frac{1}{2}{\bf 1}_{\{i_3=i_4\}}\biggl(
(T-t) I_{(00)T,t}^{*(i_1 i_2)}+
I_{(01)T,t}^{*(i_1 i_2)}\biggr)
+\frac{1}{8}(T-t)^2{\bf 1}_{\{i_1=i_2\}}{\bf 1}_{\{i_3=i_4\}}\ \ \
\hbox{w.~p.~1}.
$$

\vspace{4mm}

Let us denote as
$$
I_{(l_1\ldots l_k)T,t}^{(i_1\ldots i_k)q},\ \ \
I_{(l_1\ldots l_k)T,t}^{*(i_1\ldots i_k)q}
$$
the approximations of iterated It\^{o} and Stratonovich
stochastic integrals 
$$
I_{(l_1\ldots l_k)T,t}^{(i_1\ldots i_k)},\ \ \
I_{(l_1\ldots l_k)T,t}^{*(i_1\ldots i_k)}
$$ 
defined by
(\ref{k1000xxxx}), (\ref{k1001xxxx}), i.e. we replace
$\infty$ on $q$ in the expansions of these sto\-chas\-tic integrals.
For example, $I_{(00)T,t}^{*(i_1 i_2)q}$ is the
approximation
of the iterated 
Stratonovich stochastic integral $I_{(00)T,t}^{*(i_1 i_2)}$ obtained from 
(\ref{4004xxxx}) by replacing $\infty$ on $q$, etc.

It is easy to prove that
\begin{equation}
\label{fff09}
~~~~~{\sf M}\biggl\{\left(I_{(00)T,t}^{*(i_1 i_2)}-
I_{(00)T,t}^{*(i_1 i_2)q}
\right)^2\biggr\}
=\frac{(T-t)^2}{2}\left(\frac{1}{2}-\sum_{i=1}^q
\frac{1}{4i^2-1}\right)\ \ \ (i_1\ne i_2).
\end{equation}

Moreover, using Theorem 1.3, we obtain for $i_1\ne i_2$
$$
{\sf M}\biggl\{\left(I_{(10)T,t}^{*(i_1 i_2)}-I_{(10)T,t}^{*(i_1 i_2)q}
\right)^2\biggr\}=
{\sf M}\biggl\{\left(I_{(01)T,t}^{*(i_1 i_2)}-
I_{(01)T,t}^{*(i_1 i_2)q}\right)^2\biggr\}=
$$
$$
=\frac{(T-t)^4}{16}\left(\frac{5}{9}-
2\sum_{i=2}^q\frac{1}{4i^2-1}-
\sum_{i=1}^q
\frac{1}{(2i-1)^2(2i+3)^2}
-\right.
$$
\begin{equation}
\label{fff09xxx}
\left.-\sum_{i=0}^q\frac{(i+2)^2+(i+1)^2}{(2i+1)(2i+5)(2i+3)^2}
\right).
\end{equation}

\vspace{2mm}

For the case $i_1=i_2$
using Theorem 1.3, we have
$$
{\sf M}\biggl\{\left(I_{(10)T,t}^{(i_1 i_1)}-
I_{(10)T,t}^{(i_1 i_1)q}
\right)^2\biggr\}=
{\sf M}\biggl\{\left(I_{(01)T,t}^{(i_1 i_1)}-
I_{(01)T,t}^{(i_1 i_1)q}\right)^2\biggr\}=
$$
\begin{equation}
\label{fff09xxxx}
=\frac{(T-t)^4}{16}\left(\frac{1}{9}-
\sum_{i=0}^{q}
\frac{1}{(2i+1)(2i+5)(2i+3)^2}
-2\sum_{i=1}^{q}
\frac{1}{(2i-1)^2(2i+3)^2}\right).
\end{equation}

\vspace{2mm}

In Tables 5.1--5.3 we have calculations according to the formulas 
(\ref{fff09})--(\ref{fff09xxxx}) for various values of $q.$ 
In the given tables $\varepsilon$
means the right-hand sides of these formulas. Obviously,
these results are consistent with the estimate (\ref{zsel1}).

Let us consider (\ref{ud111}), (\ref{4006}) for $i_1=i_2$
$$
I_{(01)T,t}^{*(i_1 i_1)}
=-\frac{(T-t)^2}{4}\Biggl(
\left(\zeta_0^{(i_1)}\right)^2+
\frac{1}{\sqrt{3}}\zeta_0^{(i_1)}\zeta_1^{(i_1)}+\Biggr.
$$
\begin{equation}
\label{leto1000}
\Biggl.
+\sum_{i=0}^{\infty}\Biggl(\frac{1}{\sqrt{(2i+1)(2i+5)}(2i+3)}
\zeta_i^{(i_1)}\zeta_{i+2}^{(i_1)}-
\frac{1}{(2i-1)(2i+3)}\left(\zeta_i^{(i_1)}\right)^2
\Biggr)\Biggr),
\end{equation}

$$
I_{(10)T,t}^{*(i_1 i_1)}
=-\frac{(T-t)^2}{4}\Biggl(
\left(\zeta_0^{(i_1)}\right)^2+
\frac{1}{\sqrt{3}}\zeta_0^{(i_1)}\zeta_1^{(i_1)}+\Biggr.
$$
\begin{equation}
\label{leto1001}
\Biggl.
+\sum_{i=0}^{\infty}\Biggl(-\frac{1}{\sqrt{(2i+1)(2i+5)}(2i+3)}
\zeta_i^{(i_1)}\zeta_{i+2}^{(i_1)}+
\frac{1}{(2i-1)(2i+3)}\left(\zeta_i^{(i_1)}\right)^2
\Biggr)\Biggr).
\end{equation}

\vspace{3mm}

From (\ref{leto1000}), (\ref{leto1001}), 
considering (\ref{4001}) and (\ref{4002}), we obtain
\begin{equation}
\label{leto1002}
I_{(10)T,t}^{*(i_1 i_1)}+I_{(01)T,t}^{*(i_1 i_1)}=
-\frac{(T-t)^2}{2}\left(\left(\zeta_0^{(i_1)}\right)^2+
\frac{1}{\sqrt{3}}\zeta_0^{(i_1)}\zeta_1^{(i_1)}\right)=
I_{(0)T,t}^{(i_1)}I_{(1)T,t}^{(i_1)}\ \ \ \hbox{w.~p.~1.}
\end{equation}

\begin{table}
\centering
\caption{Confirmation of the formula (\ref{fff09})}      
\label{tab:5.1} 

\begin{tabular}{p{2.3cm}p{1.7cm}p{1.7cm}p{1.7cm}p{2.3cm}p{2.3cm}p{2.3cm}}

\hline\noalign{\smallskip}

$2\varepsilon/(T-t)^2$&0.1667&0.0238&0.0025&$2.4988\cdot 10^{-4}$&$2.4999\cdot 10^{-5}$\\

\noalign{\smallskip}\hline\noalign{\smallskip}

$q$&1&10&100&1000&10000\\

\noalign{\smallskip}\hline\noalign{\smallskip}
\end{tabular}

\end{table}

\begin{table}
\centering
\caption{Confirmation of the formula (\ref{fff09xxx})}
\label{tab:5.2}      

\begin{tabular}{p{2.3cm}p{1.7cm}p{1.7cm}p{1.7cm}p{2.3cm}p{2.3cm}p{2.3cm}}

\hline\noalign{\smallskip}

$16\varepsilon/(T-t)^4$&0.3797&0.0581&0.0062&$6.2450\cdot 10^{-4}$&$6.2495\cdot 10^{-5}$\\

\noalign{\smallskip}\hline\noalign{\smallskip}

$q$&1&10&100&1000&10000\\

\noalign{\smallskip}\hline\noalign{\smallskip}
\end{tabular}
\end{table}

\begin{table}
\centering
\caption{Confirmation of the formula (\ref{fff09xxxx})}
\label{tab:5.3}      

\begin{tabular}{p{2.1cm}p{1.2cm}p{2.1cm}p{2.1cm}p{2.3cm}p{2.3cm}p{2.3cm}}

\hline\noalign{\smallskip}

$16\varepsilon/(T-t)^4$&0.0070&$4.3551\cdot 10^{-5}$&$6.0076\cdot 10^{-8}$&$6.2251\cdot 10^{-11}$&$6.3178\cdot 10^{-14}$\\

\noalign{\smallskip}\hline\noalign{\smallskip}

$q$&1&10&100&1000&10000\\

\noalign{\smallskip}\hline\noalign{\smallskip}
\end{tabular}
\end{table}

Obtaining (\ref{leto1002}) we supposed that the formulas
(\ref{ud111}), (\ref{4006})
are valid w.~p.~1. The complete proof of this fact is
given in Sect.~1.7.2 (Theorem 1.10).

Note that it is easy to obtain the equality (\ref{leto1002}) using 
the It\^{o} formula 
and standard relations between 
iterated It\^{o}  and 
Stratonovich stochastic integrals.

Using the It\^{o} formula, we obtain
\begin{equation}
\label{leto1010}
I_{(11)T,t}^{*(i_1 i_1)}=\frac{\left(I_{(1)T,t}^{(i_1)}
\right)^2}{2}\ \ \ \hbox{w.~p.~1.}
\end{equation}

In addition, using the It\^{o} formula, we have
\begin{equation}
\label{tqtq}
I_{(20)T,t}^{(i_1 i_1)}+I_{(02)T,t}^{(i_1 i_1)}=
I_{(0)T,t}^{(i_1)}I_{(2)T,t}^{(i_1)}-
\frac{(T-t)^3}{3}\ \ \ \hbox{w.~p.~1.}
\end{equation}

From (\ref{tqtq}), considering the formulas (\ref{seg1}), (\ref{seg2}), 
we obtain
\begin{equation}
\label{leto1012}
I_{(20)T,t}^{*(i_1 i_1)}+I_{(02)T,t}^{*(i_1 i_1)}=
I_{(0)T,t}^{(i_1)}I_{(2)T,t}^{(i_1)}\ \ \ \hbox{w.~p.~1.}
\end{equation}

Let us check whether the formulas (\ref{leto1010}), (\ref{leto1012}) 
follow
from (\ref{leto1})--(\ref{leto3}), if we suppose $i_1=i_2$ in the last ones.
From (\ref{leto1})--(\ref{leto3}) for
$i_1=i_2$ we get
$$
I_{(20)T,t}^{*(i_1 i_1)}+I_{(02)T,t}^{*(i_1 i_1)}=
-\frac{(T-t)^2}{2}I_{(00)T,t}^{*(i_1 i_1)}
-(T-t)\left(I_{(10)T,t}^{*(i_1 i_1)}+I_{(01)T,t}^{*(i_1 i_1)}\right)+
$$
\begin{equation}
\label{leto1020}
+\frac{(T-t)^3}{4}
\left(\frac{1}{3}\left(\zeta_0^{(i_1)}\right)^2+
\frac{2}{3\sqrt{5}}\zeta_2^{(i_1)}\zeta_0^{(i_1)}\right),
\end{equation}

\begin{equation}
\label{leto1021}
I_{(11)T,t}^{*(i_1 i_1)}=
-\frac{(T-t)^2}{4}I_{(00)T,t}^{*(i_1 i_1)}
-\frac{T-t}{2}\left(I_{(10)T,t}^{*(i_1 i_1)}+I_{(01)T,t}^{*(i_1 i_1)}\right)
+\frac{(T-t)^3}{24}
\left(\zeta_1^{(i_1)}\right)^2.
\end{equation}

It is easy to see that 
considering (\ref{leto1002}) and (\ref{4001})--(\ref{4004xxxx}),
we actually obtain the equalities  
(\ref{leto1010}) and (\ref{leto1012})
from (\ref{leto1020}) and (\ref{leto1021}). This fact indirectly con\-firms 
the correctness of 
the formulas (\ref{leto1})--(\ref{leto3}).

Obtaining (\ref{leto1010}), (\ref{leto1012}) we supposed that the formulas
(\ref{leto1})--(\ref{leto3})
are valid w.~p.~1. The complete proof of this fact is
given in Sect.~1.7.2 (Theorem 1.10).

On the basis of 
the presented 
expansions of 
iterated stochastic integrals we 
can see that increasing of multiplicities of these integrals 
or degree indices of their weight functions 
leads
to noticeable complication of formulas 
for the mentioned expansions. 

However, increasing of the mentioned parameters leads to increasing 
of orders of smallness with respect to $T-t$ in the mean-square sense 
for iterated stochastic integrals. This leads to a sharp decrease  
of member 
quantities
in expansions of iterated stochastic 
integrals, which are required for achieving the acceptable accuracy
of approximation. In this context, let us consider the approach 
to the approximation of iterated stochastic integrals, which 
provides a possibility to obtain the mean-square approximations of 
the required accuracy without using the 
complex expansions like (\ref{4004.ii}).

Let us analyze the following 
approximation of iterated It\^{o} stochastic integral of multiplicity 3
using (\ref{zzz1})
$$
I_{(000)T,t}^{(i_1i_2i_3)q_1}
=\sum_{j_1,j_2,j_3=0}^{q_1}
C_{j_3j_2j_1}\Biggl(
\zeta_{j_1}^{(i_1)}\zeta_{j_2}^{(i_2)}\zeta_{j_3}^{(i_3)}
-{\bf 1}_{\{i_1=i_2\}}
{\bf 1}_{\{j_1=j_2\}}
\zeta_{j_3}^{(i_3)}-
\Biggr.
$$
\begin{equation}
\label{sad001xxxx}
\Biggl.
-{\bf 1}_{\{i_2=i_3\}}
{\bf 1}_{\{j_2=j_3\}}
\zeta_{j_1}^{(i_1)}-
{\bf 1}_{\{i_1=i_3\}}
{\bf 1}_{\{j_1=j_3\}}
\zeta_{j_2}^{(i_2)}\Biggr),
\end{equation}

\noindent
where $C_{j_3j_2j_1}$ is defined by (\ref{zzz2}), (\ref{zzz3}).

In particular, from (\ref{sad001xxxx}) for 
$i_1\ne i_2$, 
$i_2\ne i_3$, $i_1\ne i_3$
we obtain
\begin{equation}
\label{38}
I_{(000)T,t}^{(i_1i_2i_3)q_1}=
\sum_{j_1,j_2,j_3=0}^{q_1}
C_{j_3j_2j_1}
\zeta_{j_1}^{(i_1)}\zeta_{j_2}^{(i_2)}\zeta_{j_3}^{(i_3)}.
\end{equation}

Furthermore, using Theorem 1.3 for $k=3$,
we get
$$
{\sf M}\left\{\left(
I_{(000)T,t}^{(i_1i_2 i_3)}-
I_{(000)T,t}^{(i_1i_2 i_3)q_1}\right)^2\right\}=
$$
\begin{equation}
\label{39}
=
\frac{(T-t)^{3}}{6}
-\sum_{j_1,j_2,j_3=0}^{q_1}
C_{j_3j_2j_1}^2\ \ \ (i_1\ne i_2, i_1\ne i_3, i_2\ne i_3),
\end{equation}

\vspace{2mm}

$$
{\sf M}\left\{\left(
I_{(000)T,t}^{(i_1i_2 i_3)}-
I_{(000)T,t}^{(i_1i_2 i_3)q_1}\right)^2\right\}=
$$
\begin{equation}
\label{39a}
=
\frac{(T-t)^{3}}{6}-\sum_{j_1,j_2,j_3=0}^{q_1}
C_{j_3j_2j_1}^2
-\sum_{j_1,j_2,j_3=0}^{q_1}
C_{j_2j_3j_1}C_{j_3j_2j_1}\ \ \ (i_1\ne i_2=i_3),
\end{equation}

\vspace{2mm}

$$
{\sf M}\left\{\left(
I_{(000)T,t}^{(i_1i_2 i_3)}-
I_{(000)T,t}^{(i_1i_2 i_3)q_1}\right)^2\right\}=
$$
\begin{equation}
\label{39b}
=
\frac{(T-t)^{3}}{6}-\sum_{j_1,j_2,j_3=0}^{q_1}
C_{j_3j_2j_1}^2
-\sum_{j_1,j_2,j_3=0}^{q_1}
C_{j_3j_2j_1}C_{j_1j_2j_3}\ \ \ (i_1=i_3\ne i_2),
\end{equation}

\vspace{2mm}

$$
{\sf M}\left\{\left(
I_{(000)T,t}^{(i_1i_2 i_3)}-
I_{(000)T,t}^{(i_1i_2 i_3)q_1}\right)^2\right\}=
$$
\begin{equation}
\label{39c}
=
\frac{(T-t)^{3}}{6}-\sum_{j_1,j_2,j_3=0}^{q_1}
C_{j_3j_2j_1}^2
-\sum_{j_1,j_2,j_3=0}^{q_1}
C_{j_3j_1j_2}C_{j_3j_2j_1}\ \ \ (i_1=i_2\ne i_3).
\end{equation}

\vspace{3mm}

From the other hand, from Theorem 1.4 for $k=3$ we obtain
\begin{equation}
\label{leto1041}
~~~{\sf M}\left\{\left(
I_{(000)T,t}^{(i_1i_2 i_3)}-
I_{(000)T,t}^{(i_1i_2 i_3)q_1}\right)^2\right\}\le
6\left(\frac{(T-t)^{3}}{6}-\sum_{j_1,j_2,j_3=0}^{q_1}
C_{j_3j_2j_1}^2\right),
\end{equation}
where $i_1, i_2, i_3=1,\ldots,m$.

We can act similarly with more complicated 
iterated stochastic integrals. For example, for the 
approximation of stochastic integral
$I_{(0000)T,t}^{(i_1 i_2 i_3 i_4)}$ 
we can write (see (\ref{zzz10}))
$$
I_{(0000)T,t}^{(i_1 i_2 i_3 i_4)q_2}=
\sum_{j_1,j_2,j_3,j_4=0}^{q_2}
C_{j_4 j_3 j_2 j_1}\Biggl(
\zeta_{j_1}^{(i_1)}\zeta_{j_2}^{(i_2)}\zeta_{j_3}^{(i_3)}\zeta_{j_4}^{(i_4)}
-\Biggr.
$$
$$
-
{\bf 1}_{\{i_1=i_2\}}
{\bf 1}_{\{j_1=j_2\}}
\zeta_{j_3}^{(i_3)}
\zeta_{j_4}^{(i_4)}
-
{\bf 1}_{\{i_1=i_3\}}
{\bf 1}_{\{j_1=j_3\}}
\zeta_{j_2}^{(i_2)}
\zeta_{j_4}^{(i_4)}-
$$
$$
-
{\bf 1}_{\{i_1=i_4\}}
{\bf 1}_{\{j_1=j_4\}}
\zeta_{j_2}^{(i_2)}
\zeta_{j_3}^{(i_3)}
-
{\bf 1}_{\{i_2=i_3\}}
{\bf 1}_{\{j_2=j_3\}}
\zeta_{j_1}^{(i_1)}
\zeta_{j_4}^{(i_4)}-
$$
$$
-
{\bf 1}_{\{i_2=i_4\}}
{\bf 1}_{\{j_2=j_4\}}
\zeta_{j_1}^{(i_1)}
\zeta_{j_3}^{(i_3)}
-
{\bf 1}_{\{i_3=i_4\}}
{\bf 1}_{\{j_3=j_4\}}
\zeta_{j_1}^{(i_1)}
\zeta_{j_2}^{(i_2)}+
$$
$$
+
{\bf 1}_{\{i_1=i_2\}}
{\bf 1}_{\{j_1=j_2\}}
{\bf 1}_{\{i_3=i_4\}}
{\bf 1}_{\{j_3=j_4\}}
+
{\bf 1}_{\{i_1=i_3\}}
{\bf 1}_{\{j_1=j_3\}}
{\bf 1}_{\{i_2=i_4\}}
{\bf 1}_{\{j_2=j_4\}}+
$$
\begin{equation}
\label{kol1}
+\Biggl.
{\bf 1}_{\{i_1=i_4\}}
{\bf 1}_{\{j_1=j_4\}}
{\bf 1}_{\{i_2=i_3\}}
{\bf 1}_{\{j_2=j_3\}}\Biggr),
\end{equation}

\noindent
where $C_{j_4 j_3 j_2 j_1}$ is defined by (\ref{zzz11}), (\ref{zzz12}).

Moreover, according to Theorem 1.4 for $k=4$, we get
$$
{\sf M}\left\{\left(
I_{(0000)T,t}^{(i_1i_2 i_3 i_4)}-
I_{(0000)T,t}^{(i_1i_2 i_3 i_4)q_2}\right)^2\right\}\le
24\left(\frac{(T-t)^{4}}{24}-\sum_{j_1,j_2,j_3,j_4=0}^{q_2}
C_{j_4j_3j_2j_1}^2\right),
$$
where
$i_1, i_2, i_3, i_4=1,\ldots,m$.

For pairwise different $i_1, i_2, i_3, i_4=1,\ldots,m$ from 
Theorem 1.3 we obtain
\begin{equation}
\label{r7}
~~~~~{\sf M}\left\{\left(
I_{(0000)T,t}^{(i_1i_2 i_3i_4)}-
I_{(0000)T,t}^{(i_1i_2 i_3i_4)q_2}\right)^2\right\}=
\frac{(T-t)^{4}}{24}
-\sum_{j_1,j_2,j_3,j_4=0}^{q_2}
C_{j_4j_3j_2j_1}^2.
\end{equation}

Using Theorem 1.3, we can calculate exactly the left-hand
side of (\ref{r7})
for any possible combinations
of $i_1, i_2, i_3, i_4$. These relations were obtained in 
Sect.~1.2. For example
$$
{\sf M}\left\{\left(
I_{(0000)T,t}^{(i_1i_2 i_3 i_4)}-
I_{(0000)T,t}^{(i_1i_2 i_3 i_4)q_2}\right)^2\right\}=
$$
$$
= \frac{(T-t)^{4}}{24} - \sum_{j_1,j_2,j_3,j_4=0}^{q_2}
C_{j_4j_3j_2j_1}\Biggl(\sum\limits_{(j_1,j_2)}\Biggl(
\sum\limits_{(j_3,j_4)}
C_{j_4j_3j_2j_1}\Biggr)\Biggr),
$$

\noindent
where $i_1=i_2\ne i_3=i_4$ and
$$
\sum\limits_{(j_1,j_2)}
$$
means the sum with respect to permutations $(j_1,j_2)$.

\begin{table}
\centering
\caption{Coefficients $\bar C_{0j_2j_1}$}
\label{tab:5.4}      

\begin{tabular}{p{1.3cm}p{1.3cm}p{1.3cm}p{1.3cm}p{1.3cm}p{1.3cm}p{1.3cm}p{1.3cm}}

\hline\noalign{\smallskip}

&$j_1=0$&$j_1=1$&$j_1=2$&$j_1=3$&$j_1=4$&$j_1=5$&$j_1=6$\\

\noalign{\smallskip}\hline\noalign{\smallskip}

$j_2=0$&$\frac{4}{3}$&$\frac{-2}{3}$&$\frac{2}{15}$&$0$&0&0&0\\

\noalign{\smallskip}

$j_2=1$&$0$&$\frac{2}{15}$&$\frac{-2}{15}$&$\frac{4}{105}$&0&0&0\\

\noalign{\smallskip}

$j_2=2$&$\frac{-4}{15}$&$\frac{2}{15}$&$\frac{2}{105}$&$\frac{-2}{35}$&
$\frac{2}{105}$&0&0\\

\noalign{\smallskip}

$j_2=3$&$0$&$\frac{-2}{35}$&$\frac{2}{35}$&$\frac{2}{315}$&
$\frac{-2}{63}$&$\frac{8}{693}$&0\\

\noalign{\smallskip}

$j_2=4$&$0$&$0$&$\frac{-8}{315}$&$\frac{2}{63}$&
$\frac{2}{693}$&$\frac{-2}{99}$&$\frac{10}{1287}$\\

\noalign{\smallskip}

$j_2=5$&$0$&$0$&$0$&$\frac{-10}{693}$&
$\frac{2}{99}$&$\frac{2}{1287}$&$\frac{-2}{143}$\\

\noalign{\smallskip}

$j_2=6$&$0$&$0$&$0$&$0$&$\frac{-4}{429}$&$\frac{2}{143}$&$\frac{2}{2145}$\\

\noalign{\smallskip}\hline\noalign{\smallskip}
\end{tabular}
\end{table}

\begin{table}
\centering
\caption{Coefficients $\bar C_{1j_2j_1}$}
\label{tab:5.5}      

\begin{tabular}{p{1.3cm}p{1.3cm}p{1.3cm}p{1.3cm}p{1.3cm}p{1.3cm}p{1.3cm}p{1.3cm}}

\hline\noalign{\smallskip}

&$j_1=0$&$j_1=1$&$j_1=2$&$j_1=3$&$j_1=4$&$j_1=5$&$j_1=6$\\

\noalign{\smallskip}\hline\noalign{\smallskip}

$j_2=0$&$\frac{2}{3}$&$\frac{-4}{15}$&$0$&$\frac{2}{105}$&0&0&0\\

\noalign{\smallskip}

$j_2=1$&$\frac{2}{15}$&$0$&$\frac{-4}{105}$&$0$&$\frac{2}{315}$&0&0\\

\noalign{\smallskip}

$j_2=2$&$\frac{-2}{15}$&$\frac{8}{105}$&$0$&$\frac{-2}{105}$&
0&$\frac{4}{1155}$&0\\

\noalign{\smallskip}

$j_2=3$&$\frac{-2}{35}$&0&$\frac{8}{315}$&0&
$\frac{-38}{3465}$&0&$\frac{20}{9009}$\\

\noalign{\smallskip}

$j_2=4$&$0$&$\frac{-4}{315}$&0&$\frac{46}{3465}$&
0&$\frac{-64}{9009}$&0\\

\noalign{\smallskip}

$j_2=5$&$0$&$0$&$\frac{-4}{693}$&0&
$\frac{74}{9009}$&0&$\frac{-32}{6435}$\\

\noalign{\smallskip}

$j_2=6$&$0$&$0$&$0$&$\frac{-10}{3003}$&$0$&$\frac{4}{715}$&$0$\\

\noalign{\smallskip}\hline\noalign{\smallskip}
\end{tabular}

\end{table}

\begin{table}
\centering
\caption{Coefficients $\bar C_{2j_2j_1}$}
\label{tab:5.6}      

\begin{tabular}{p{1.3cm}p{1.3cm}p{1.3cm}p{1.3cm}p{1.3cm}p{1.3cm}p{1.3cm}p{1.3cm}}

\hline\noalign{\smallskip}

&$j_1=0$&$j_1=1$&$j_1=2$&$j_1=3$&$j_1=4$&$j_1=5$&$j_1=6$\\

\noalign{\smallskip}\hline\noalign{\smallskip}

$j_2=0$&$\frac{2}{15}$&0&$\frac{-4}{105}$&0&$\frac{2}{315}$&0&0\\

\noalign{\smallskip}

$j_2=1$&$\frac{2}{15}$&$\frac{-4}{105}$&0&$\frac{-2}{315}$&0&$\frac{8}{3465}$&0\\

\noalign{\smallskip}

$j_2=2$&$\frac{2}{105}$&0&$0$&0&
$\frac{-2}{495}$&0&$\frac{4}{3003}$\\

\noalign{\smallskip}

$j_2=3$&$\frac{-2}{35}$&$\frac{8}{315}$&0&$\frac{-2}{3465}$&
0&$\frac{-116}{45045}$&0\\

\noalign{\smallskip}

$j_2=4$&$\frac{-8}{315}$&0&$\frac{4}{495}$&0&
$\frac{-2}{6435}$&0&$\frac{-16}{9009}$\\

\noalign{\smallskip}

$j_2=5$&$0$&$\frac{-4}{693}$&0&$\frac{38}{9009}$&
0&$\frac{-8}{45045}$&0\\

\noalign{\smallskip}

$j_2=6$&$0$&$0$&$\frac{-8}{3003}$&$0$&$\frac{118}{45045}$&$0$&$\frac{-4}{36465}$\\

\noalign{\smallskip}\hline\noalign{\smallskip}
\end{tabular}
\end{table}

\begin{table}
\centering
\caption{Coefficients $\bar C_{3j_2j_1}$}
\label{tab:5.7}      

\begin{tabular}{p{1.3cm}p{1.3cm}p{1.3cm}p{1.3cm}p{1.3cm}p{1.3cm}p{1.3cm}p{1.3cm}}

\hline\noalign{\smallskip}

&$j_1=0$&$j_1=1$&$j_1=2$&$j_1=3$&$j_1=4$&$j_1=5$&$j_1=6$\\

\noalign{\smallskip}\hline\noalign{\smallskip}

$j_2=0$&$0$&$\frac{2}{105}$&$0$&$\frac{-4}{315}$&$0$&$\frac{2}{693}$&0\\
\noalign{\smallskip}
$j_2=1$&$\frac{4}{105}$&0&$\frac{-2}{315}$&0&$\frac{-8}{3465}$&0&$\frac{10}{9009}$\\
\noalign{\smallskip}
$j_2=2$&$\frac{2}{35}$&$\frac{-2}{105}$&$0$&$\frac{4}{3465}$&
$0$&$\frac{-74}{45045}$&0\\
\noalign{\smallskip}
$j_2=3$&$\frac{2}{315}$&$0$&$\frac{-2}{3465}$&0&
$\frac{16}{45045}$&0&$\frac{-10}{9009}$\\
\noalign{\smallskip}
$j_2=4$&$\frac{-2}{63}$&$\frac{46}{3465}$&0&$\frac{-32}{45045}$&
0&$\frac{2}{9009}$&0\\
\noalign{\smallskip}
$j_2=5$&$\frac{-10}{693}$&0&$\frac{38}{9009}$&0&
$\frac{-4}{9009}$&0&$\frac{122}{765765}$\\
\noalign{\smallskip}
$j_2=6$&$0$&$\frac{-10}{3003}$&$0$&$\frac{20}{9009}$&$0$&$\frac{-226}{765765}$&$0$\\

\noalign{\smallskip}\hline\noalign{\smallskip}
\end{tabular}

\end{table}

\begin{table}
\centering
\caption{Coefficients $\bar C_{4j_2j_1}$}
\label{tab:5.8}      

\begin{tabular}{p{1.3cm}p{1.3cm}p{1.3cm}p{1.3cm}p{1.3cm}p{1.3cm}p{1.3cm}p{1.3cm}}

\hline\noalign{\smallskip}

&$j_1=0$&$j_1=1$&$j_1=2$&$j_1=3$&$j_1=4$&$j_1=5$&$j_1=6$\\

\noalign{\smallskip}\hline\noalign{\smallskip}

$j_2=0$&$0$&0&$\frac{2}{315}$&0&$\frac{-4}{693}$&0&$\frac{2}{1287}$\\
\noalign{\smallskip}
$j_2=1$&0&$\frac{2}{315}$&0&$\frac{-8}{3465}$&0&$\frac{-10}{9009}$&0\\
\noalign{\smallskip}
$j_2=2$&$\frac{2}{105}$&0&$\frac{-2}{495}$&0&
$\frac{4}{6435}$&0&$\frac{-38}{45045}$\\
\noalign{\smallskip}
$j_2=3$&$\frac{2}{63}$&$\frac{-38}{3465}$&0&$\frac{16}{45045}$&
0&$\frac{2}{9009}$&0\\
\noalign{\smallskip}
$j_2=4$&$\frac{2}{693}$&0&$\frac{-2}{6435}$&0&
0&0&$\frac{2}{13923}$\\
\noalign{\smallskip}
$j_2=5$&$\frac{-2}{99}$&$\frac{74}{9009}$&0&$\frac{-4}{9009}$&
0&$\frac{-2}{153153}$&0\\
\noalign{\smallskip}
$j_2=6$&$\frac{-4}{429}$&$0$&$\frac{118}{45045}$&$0$&$\frac{-4}{13923}$&$0$&
$\frac{-2}{188955}$\\
\noalign{\smallskip}\hline\noalign{\smallskip}
\end{tabular}

\end{table}

\begin{table}
\centering
\caption{Coefficients $\bar C_{5j_2j_1}$}
\label{tab:5.9}      

\begin{tabular}{p{1.3cm}p{1.3cm}p{1.3cm}p{1.3cm}p{1.3cm}p{1.3cm}p{1.3cm}p{1.3cm}}

\hline\noalign{\smallskip}

&$j_1=0$&$j_1=1$&$j_1=2$&$j_1=3$&$j_1=4$&$j_1=5$&$j_1=6$\\

\noalign{\smallskip}\hline\noalign{\smallskip}

$j_2=0$&$0$&0&0&$\frac{2}{693}$&0&$\frac{-4}{1287}$&0\\
\noalign{\smallskip}

$j_2=1$&0&0&$\frac{8}{3465}$&0&$\frac{-10}{9009}$&0&$\frac{-4}{6435}$\\
\noalign{\smallskip}
$j_2=2$&0&$\frac{4}{1155}$&0&$\frac{-74}{45045}$&
0&$\frac{16}{45045}$&0\\
\noalign{\smallskip}
$j_2=3$&$\frac{8}{693}$&0&$\frac{-116}{45045}$&0&
$\frac{2}{9009}$&0&$\frac{8}{58905}$\\
\noalign{\smallskip}
$j_2=4$&$\frac{2}{99}$&$\frac{-64}{9009}$&0&$\frac{2}{9009}$&
0&$\frac{4}{153153}$&0\\
\noalign{\smallskip}
$j_2=5$&$\frac{2}{1287}$&$0$&$\frac{-8}{45045}$&0&
$\frac{-2}{153153}$&0&$\frac{4}{415701}$\\
\noalign{\smallskip}
$j_2=6$&$\frac{-2}{143}$&$\frac{4}{715}$&$0$&$\frac{-226}{765765}$&
$0$&$\frac{-8}{415701}$&$0$\\
\noalign{\smallskip}\hline\noalign{\smallskip}
\end{tabular}
\end{table}

\begin{table}
\centering
\caption{Coefficients $\bar C_{6j_2j_1}$}
\label{tab:5.10}      

\begin{tabular}{p{1.3cm}p{1.3cm}p{1.3cm}p{1.3cm}p{1.3cm}p{1.3cm}p{1.3cm}p{1.3cm}}

\hline\noalign{\smallskip}

&$j_1=0$&$j_1=1$&$j_1=2$&$j_1=3$&$j_1=4$&$j_1=5$&$j_1=6$\\

\noalign{\smallskip}\hline\noalign{\smallskip}

$j_2=0$&$0$&0&0&$0$&$\frac{2}{1287}$&$0$&$\frac{-4}{2145}$\\
\noalign{\smallskip}
$j_2=1$&0&0&$0$&$\frac{10}{9009}$&$0$&$\frac{-4}{6435}$&$0$\\
\noalign{\smallskip}
$j_2=2$&0&$0$&$\frac{4}{3003}$&$0$&
$\frac{-38}{45045}$&$0$&$\frac{8}{36465}$\\
\noalign{\smallskip}
$j_2=3$&$0$&$\frac{20}{9009}$&0&$\frac{-10}{9009}$&0&
$\frac{8}{58905}$&0\\
\noalign{\smallskip}
$j_2=4$&$\frac{10}{1287}$&$0$&$\frac{-16}{9009}$&0&$\frac{2}{13923}$&
0&$\frac{4}{188955}$\\
\noalign{\smallskip}
$j_2=5$&$\frac{2}{143}$&$\frac{-32}{6435}$&0&
$\frac{122}{765765}$&0&$\frac{4}{415701}$&0\\
\noalign{\smallskip}
$j_2=6$&$\frac{2}{2145}$&$0$&$\frac{-4}{36465}$&$0$&$\frac{-2}{188955}$&
$0$&$0$\\
\noalign{\smallskip}\hline\noalign{\smallskip}
\end{tabular}
\end{table}

\begin{table}
\centering
\caption{Coefficients $\bar C_{00j_2j_1}$}
\label{tab:5.11}      

\begin{tabular}{p{1.3cm}p{1.3cm}p{1.3cm}p{1.3cm}}

\hline\noalign{\smallskip}

&$j_1=0$&$j_1=1$&$j_1=2$\\

\noalign{\smallskip}\hline\noalign{\smallskip}

$j_2=0$&$\frac{2}{3}$&$\frac{-2}{5}$&$\frac{2}{15}$\\
\noalign{\smallskip}
$j_2=1$&$\frac{-2}{15}$&$\frac{2}{15}$&$\frac{-2}{21}$\\
\noalign{\smallskip}
$j_2=2$&$\frac{-2}{15}$&$\frac{2}{35}$&$\frac{2}{105}$\\
\noalign{\smallskip}\hline\noalign{\smallskip}
\end{tabular}

\end{table}

\begin{table}
\centering
\caption{Coefficients $\bar C_{10j_2j_1}$}
\label{tab:5.12}      
\begin{tabular}{p{1.3cm}p{1.3cm}p{1.3cm}p{1.3cm}}
\hline\noalign{\smallskip}
&$j_1=0$&$j_1=1$&$j_1=2$\\
\noalign{\smallskip}\hline\noalign{\smallskip}
$j_2=0$&$\frac{2}{5}$&$\frac{-2}{9}$&$\frac{2}{35}$\\
\noalign{\smallskip}
$j_2=1$&$\frac{-2}{45}$&$\frac{2}{35}$&$\frac{-2}{45}$\\
\noalign{\smallskip}
$j_2=2$&$\frac{-2}{21}$&$\frac{2}{45}$&$\frac{2}{315}$\\
\noalign{\smallskip}\hline\noalign{\smallskip}
\end{tabular}
\vspace{3mm}
\end{table}

\begin{table}
\centering
\caption{Coefficients $\bar C_{02j_2j_1}$}
\label{tab:5.13}      
\begin{tabular}{p{1.3cm}p{1.3cm}p{1.3cm}p{1.3cm}}
\hline\noalign{\smallskip}
&$j_1=0$&$j_1=1$&$j_1=2$\\
\noalign{\smallskip}\hline\noalign{\smallskip}
$j_2=0$&$\frac{-2}{15}$&$\frac{2}{21}$&$\frac{-4}{105}$\\
\noalign{\smallskip}
$j_2=1$&$\frac{2}{35}$&$\frac{-4}{105}$&$\frac{2}{105}$\\
\noalign{\smallskip}
$j_2=2$&$\frac{4}{105}$&$\frac{-2}{105}$&$0$\\
\noalign{\smallskip}\hline\noalign{\smallskip}
\end{tabular}
\vspace{3mm}
\end{table}

\begin{table}
\centering
\caption{Coefficients $\bar C_{01j_2j_1}$}
\label{tab:5.14}      
\begin{tabular}{p{1.3cm}p{1.3cm}p{1.3cm}p{1.3cm}}
\hline\noalign{\smallskip}
&$j_1=0$&$j_1=1$&$j_1=2$\\
\noalign{\smallskip}\hline\noalign{\smallskip}
$j_2=0$&$\frac{2}{15}$&$\frac{-2}{45}$&$\frac{-2}{105}$\\
\noalign{\smallskip}
$j_2=1$&$\frac{2}{45}$&$\frac{-2}{105}$&$\frac{2}{315}$\\
\noalign{\smallskip}
$j_2=2$&$\frac{-2}{35}$&$\frac{2}{63}$&$\frac{-2}{315}$\\
\noalign{\smallskip}\hline\noalign{\smallskip}
\end{tabular}
\vspace{3mm}
\end{table}

\begin{table}
\centering
\caption{Coefficients $\bar C_{11j_2j_1}$}
\label{tab:5.15}      
\begin{tabular}{p{1.3cm}p{1.3cm}p{1.3cm}p{1.3cm}}
\hline\noalign{\smallskip}
&$j_1=0$&$j_1=1$&$j_1=2$\\
\noalign{\smallskip}\hline\noalign{\smallskip}
$j_2=0$&$\frac{2}{15}$&$\frac{-2}{35}$&$0$\\
\noalign{\smallskip}
$j_2=1$&$\frac{2}{105}$&$0$&$\frac{-2}{315}$\\
\noalign{\smallskip}
$j_2=2$&$\frac{-4}{105}$&$\frac{2}{105}$&$0$\\
\noalign{\smallskip}\hline\noalign{\smallskip}
\end{tabular}
\vspace{3mm}
\end{table}

\newpage

\begin{table}
\centering
\caption{Coefficients $\bar C_{20j_2j_1}$}
\label{tab:5.16}      
\begin{tabular}{p{1.3cm}p{1.3cm}p{1.3cm}p{1.3cm}}
\hline\noalign{\smallskip}
&$j_1=0$&$j_1=1$&$j_1=2$\\
\noalign{\smallskip}\hline\noalign{\smallskip}
$j_2=0$&$\frac{2}{15}$&$\frac{-2}{35}$&$0$\\
\noalign{\smallskip}
$j_2=1$&$\frac{2}{105}$&$0$&$\frac{-2}{315}$\\
\noalign{\smallskip}
$j_2=2$&$\frac{-4}{105}$&$\frac{2}{105}$&$0$\\
\noalign{\smallskip}\hline\noalign{\smallskip}
\end{tabular}
\vspace{3mm}
\end{table}

\begin{table}
\centering
\caption{Coefficients $\bar C_{21j_2j_1}$}
\label{tab:5.17}      
\begin{tabular}{p{1.3cm}p{1.3cm}p{1.3cm}p{1.3cm}}
\hline\noalign{\smallskip}
&$j_1=0$&$j_1=1$&$j_1=2$\\
\noalign{\smallskip}\hline\noalign{\smallskip}
$j_2=0$&$\frac{2}{21}$&$\frac{-2}{45}$&$\frac{2}{315}$\\
\noalign{\smallskip}
$j_2=1$&$\frac{2}{315}$&$\frac{2}{315}$&$\frac{-2}{225}$\\
\noalign{\smallskip}
$j_2=2$&$\frac{-2}{105}$&$\frac{2}{225}$&$\frac{2}{1155}$\\
\noalign{\smallskip}\hline\noalign{\smallskip}
\end{tabular}
\end{table}

\begin{table}
\centering
\caption{Coefficients $\bar C_{12j_2j_1}$}
\label{tab:5.18}      
\begin{tabular}{p{1.3cm}p{1.3cm}p{1.3cm}p{1.3cm}}
\hline\noalign{\smallskip}
&$j_1=0$&$j_1=1$&$j_1=2$\\
\noalign{\smallskip}\hline\noalign{\smallskip}
$j_2=0$&$\frac{-2}{35}$&$\frac{2}{45}$&$\frac{-2}{105}$\\
\noalign{\smallskip}
$j_2=1$&$\frac{2}{63}$&$\frac{-2}{105}$&$\frac{2}{225}$\\
\noalign{\smallskip}
$j_2=2$&$\frac{2}{105}$&$\frac{-2}{225}$&$\frac{-2}{3465}$\\
\noalign{\smallskip}\hline\noalign{\smallskip}
\end{tabular}
\end{table}

Assume that $q_1=6$. In Tables 5.4--5.10 we have the exact 
values of coefficients 
$\bar C_{j_3j_2j_1}$ ($j_1,j_2,j_3=0, 1,\ldots,6$). 
Here and further in this section the Fourier--Legendre coefficients 
have been calculated exactly using computer algebra
system Derive. Note that in \cite{Kuz-Kuz}, \cite{Mikh-1}
the database with 270,000 exactly
calculated Fourier--Legendre coefficients was described.
This database was used in the software package,
which is written in the Python programming language
for the implementation of the numerical schemes (\ref{al1})-(\ref{al5}),
(\ref{al1x})-(\ref{al5x}).

Calculating the value on the right-hand side of
(\ref{39}) for $q_1=6$ 
($i_1\ne i_2,$ $i_1\ne i_3,$ $i_3\ne i_2$),
we obtain the following  
approximate equality
$$
{\sf M}\left\{\left(
I_{(000)T,t}^{(i_1i_2 i_3)}-
I_{(000)T,t}^{(i_1i_2 i_3)q_1}\right)^2\right\}\approx
0.01956(T-t)^3.
$$

Let us choose, for example, $q_2=2.$ In Tables 5.11--5.19  
we have the exact values of coefficients 
$\bar C_{j_4j_3j_2j_1}$
($j_1,j_2,j_3,j_4=0, 1, 2$).  
In the case of pairwise different
$i_1, i_2, i_3, i_4$ we obtain from (\ref{r7}) the following  
approximate equality
\begin{equation}
\label{46000}
~~~{\sf M}\left\{\left(
I_{(0000)T,t}^{(i_1i_2i_3 i_4)}-
I_{(0000)T,t}^{(i_1i_2i_3 i_4)q_2}\right)^2\right\}\approx
0.0236084(T-t)^4.
\end{equation}

Let us analyze the following four approximations of the iterated It\^{o} 
stochastic integrals (see (\ref{sss1})--(\ref{sss4}))
$$
I_{(001)T,t}^{(i_1i_2i_3)q_3}
=
\sum_{j_1,j_2,j_3=0}^{q_3}
C_{j_3j_2j_1}^{001}\Biggl(
\zeta_{j_1}^{(i_1)}\zeta_{j_2}^{(i_2)}\zeta_{j_3}^{(i_3)}
-{\bf 1}_{\{i_1=i_2\}}
{\bf 1}_{\{j_1=j_2\}}
\zeta_{j_3}^{(i_3)}-
\Biggr.
$$
\begin{equation}
\label{r9}
\Biggl.
-{\bf 1}_{\{i_2=i_3\}}
{\bf 1}_{\{j_2=j_3\}}
\zeta_{j_1}^{(i_1)}-
{\bf 1}_{\{i_1=i_3\}}
{\bf 1}_{\{j_1=j_3\}}
\zeta_{j_2}^{(i_2)}\Biggr),
\end{equation}

$$
I_{(010)T,t}^{(i_1i_2i_3)q_3}
=
\sum_{j_1,j_2,j_3=0}^{q_3}
C_{j_3j_2j_1}^{010}\Biggl(
\zeta_{j_1}^{(i_1)}\zeta_{j_2}^{(i_2)}\zeta_{j_3}^{(i_3)}
-{\bf 1}_{\{i_1=i_2\}}
{\bf 1}_{\{j_1=j_2\}}
\zeta_{j_3}^{(i_3)}-
\Biggr.
$$
\begin{equation}
\label{r10}
\Biggl.
-{\bf 1}_{\{i_2=i_3\}}
{\bf 1}_{\{j_2=j_3\}}
\zeta_{j_1}^{(i_1)}-
{\bf 1}_{\{i_1=i_3\}}
{\bf 1}_{\{j_1=j_3\}}
\zeta_{j_2}^{(i_2)}\Biggr),
\end{equation}

$$
I_{(100)T,t}^{(i_1i_2i_3)q_3}
=
\sum_{j_1,j_2,j_3=0}^{q_3}
C_{j_3j_2j_1}^{100}\Biggl(
\zeta_{j_1}^{(i_1)}\zeta_{j_2}^{(i_2)}\zeta_{j_3}^{(i_3)}
-{\bf 1}_{\{i_1=i_2\}}
{\bf 1}_{\{j_1=j_2\}}
\zeta_{j_3}^{(i_3)}-
\Biggr.
$$
\begin{equation}
\label{r10a}
\Biggl.
-{\bf 1}_{\{i_2=i_3\}}
{\bf 1}_{\{j_2=j_3\}}
\zeta_{j_1}^{(i_1)}-
{\bf 1}_{\{i_1=i_3\}}
{\bf 1}_{\{j_1=j_3\}}
\zeta_{j_2}^{(i_2)}\Biggr),
\end{equation}

$$
I_{(00000)T,t}^{(i_1i_2i_3i_4 i_5)q_4}=
\sum_{j_1,j_2,j_3,j_4,j_5=0}^{q_4}
C_{j_5 j_4 j_3 j_2 j_1}\Biggl(
\prod_{l=1}^5\zeta_{j_l}^{(i_l)}
-\Biggr.
$$
$$
-
{\bf 1}_{\{i_1=i_2\}}
{\bf 1}_{\{j_1=j_2\}}
\zeta_{j_3}^{(i_3)}
\zeta_{j_4}^{(i_4)}
\zeta_{j_5}^{(i_5)}-
{\bf 1}_{\{i_1=i_3\}}
{\bf 1}_{\{j_1=j_3\}}
\zeta_{j_2}^{(i_2)}
\zeta_{j_4}^{(i_4)}
\zeta_{j_5}^{(i_5)}-
$$
$$
-
{\bf 1}_{\{i_1=i_4\}}
{\bf 1}_{\{j_1=j_4\}}
\zeta_{j_2}^{(i_2)}
\zeta_{j_3}^{(i_3)}
\zeta_{j_5}^{(i_5)}-
{\bf 1}_{\{i_1=i_5\}}
{\bf 1}_{\{j_1=j_5\}}
\zeta_{j_2}^{(i_2)}
\zeta_{j_3}^{(i_3)}
\zeta_{j_4}^{(i_4)}-
$$
$$
-
{\bf 1}_{\{i_2=i_3\}}
{\bf 1}_{\{j_2=j_3\}}
\zeta_{j_1}^{(i_1)}
\zeta_{j_4}^{(i_4)}
\zeta_{j_5}^{(i_5)}-
{\bf 1}_{\{i_2=i_4\}}
{\bf 1}_{\{j_2=j_4\}}
\zeta_{j_1}^{(i_1)}
\zeta_{j_3}^{(i_3)}
\zeta_{j_5}^{(i_5)}-
$$
$$
-
{\bf 1}_{\{i_2=i_5\}}
{\bf 1}_{\{j_2=j_5\}}
\zeta_{j_1}^{(i_1)}
\zeta_{j_3}^{(i_3)}
\zeta_{j_4}^{(i_4)}
-{\bf 1}_{\{i_3=i_4\}}
{\bf 1}_{\{j_3=j_4\}}
\zeta_{j_1}^{(i_1)}
\zeta_{j_2}^{(i_2)}
\zeta_{j_5}^{(i_5)}-
$$
$$
-
{\bf 1}_{\{i_3=i_5\}}
{\bf 1}_{\{j_3=j_5\}}
\zeta_{j_1}^{(i_1)}
\zeta_{j_2}^{(i_2)}
\zeta_{j_4}^{(i_4)}
-{\bf 1}_{\{i_4=i_5\}}
{\bf 1}_{\{j_4=j_5\}}
\zeta_{j_1}^{(i_1)}
\zeta_{j_2}^{(i_2)}
\zeta_{j_3}^{(i_3)}+
$$
$$
+
{\bf 1}_{\{i_1=i_2\}}
{\bf 1}_{\{j_1=j_2\}}
{\bf 1}_{\{i_3=i_4\}}
{\bf 1}_{\{j_3=j_4\}}\zeta_{j_5}^{(i_5)}+
{\bf 1}_{\{i_1=i_2\}}
{\bf 1}_{\{j_1=j_2\}}
{\bf 1}_{\{i_3=i_5\}}
{\bf 1}_{\{j_3=j_5\}}\zeta_{j_4}^{(i_4)}+
$$
$$
+
{\bf 1}_{\{i_1=i_2\}}
{\bf 1}_{\{j_1=j_2\}}
{\bf 1}_{\{i_4=i_5\}}
{\bf 1}_{\{j_4=j_5\}}\zeta_{j_3}^{(i_3)}+
{\bf 1}_{\{i_1=i_3\}}
{\bf 1}_{\{j_1=j_3\}}
{\bf 1}_{\{i_2=i_4\}}
{\bf 1}_{\{j_2=j_4\}}\zeta_{j_5}^{(i_5)}+
$$
$$
+
{\bf 1}_{\{i_1=i_3\}}
{\bf 1}_{\{j_1=j_3\}}
{\bf 1}_{\{i_2=i_5\}}
{\bf 1}_{\{j_2=j_5\}}\zeta_{j_4}^{(i_4)}+
{\bf 1}_{\{i_1=i_3\}}
{\bf 1}_{\{j_1=j_3\}}
{\bf 1}_{\{i_4=i_5\}}
{\bf 1}_{\{j_4=j_5\}}\zeta_{j_2}^{(i_2)}+
$$
$$
+
{\bf 1}_{\{i_1=i_4\}}
{\bf 1}_{\{j_1=j_4\}}
{\bf 1}_{\{i_2=i_3\}}
{\bf 1}_{\{j_2=j_3\}}\zeta_{j_5}^{(i_5)}+
{\bf 1}_{\{i_1=i_4\}}
{\bf 1}_{\{j_1=j_4\}}
{\bf 1}_{\{i_2=i_5\}}
{\bf 1}_{\{j_2=j_5\}}\zeta_{j_3}^{(i_3)}+
$$
$$
+
{\bf 1}_{\{i_1=i_4\}}
{\bf 1}_{\{j_1=j_4\}}
{\bf 1}_{\{i_3=i_5\}}
{\bf 1}_{\{j_3=j_5\}}\zeta_{j_2}^{(i_2)}+
{\bf 1}_{\{i_1=i_5\}}
{\bf 1}_{\{j_1=j_5\}}
{\bf 1}_{\{i_2=i_3\}}
{\bf 1}_{\{j_2=j_3\}}\zeta_{j_4}^{(i_4)}+
$$
$$
+
{\bf 1}_{\{i_1=i_5\}}
{\bf 1}_{\{j_1=j_5\}}
{\bf 1}_{\{i_2=i_4\}}
{\bf 1}_{\{j_2=j_4\}}\zeta_{j_3}^{(i_3)}+
{\bf 1}_{\{i_1=i_5\}}
{\bf 1}_{\{j_1=j_5\}}
{\bf 1}_{\{i_3=i_4\}}
{\bf 1}_{\{j_3=j_4\}}\zeta_{j_2}^{(i_2)}+
$$
$$
+
{\bf 1}_{\{i_2=i_3\}}
{\bf 1}_{\{j_2=j_3\}}
{\bf 1}_{\{i_4=i_5\}}
{\bf 1}_{\{j_4=j_5\}}\zeta_{j_1}^{(i_1)}+
{\bf 1}_{\{i_2=i_4\}}
{\bf 1}_{\{j_2=j_4\}}
{\bf 1}_{\{i_3=i_5\}}
{\bf 1}_{\{j_3=j_5\}}\zeta_{j_1}^{(i_1)}+
$$
\begin{equation}
\label{r11}
+\Biggl.
{\bf 1}_{\{i_2=i_5\}}
{\bf 1}_{\{j_2=j_5\}}
{\bf 1}_{\{i_3=i_4\}}
{\bf 1}_{\{j_3=j_4\}}\zeta_{j_1}^{(i_1)}\Biggr).
\end{equation}

\begin{table}
\centering
\caption{Coefficients $\bar C_{22j_2j_1}$}
\label{tab:5.19}      
\begin{tabular}{p{1.3cm}p{1.3cm}p{1.3cm}p{1.3cm}}
\hline\noalign{\smallskip}
&$j_1=0$&$j_1=1$&$j_1=2$\\
\noalign{\smallskip}\hline\noalign{\smallskip}
$j_2=0$&$\frac{2}{105}$&$\frac{-2}{315}$&$0$\\
\noalign{\smallskip}
$j_2=1$&$\frac{2}{315}$&$0$&$\frac{-2}{1155}$\\
\noalign{\smallskip}
$j_2=2$&$0$&$\frac{2}{3465}$&$0$\\
\noalign{\smallskip}\hline\noalign{\smallskip}
\end{tabular}
\end{table}

\begin{table}
\centering
\caption{Coefficients $\bar C_{0j_2j_1}^{001}$}
\label{tab:5.20}      
\begin{tabular}{p{1.3cm}p{1.3cm}p{1.3cm}p{1.3cm}}
\hline\noalign{\smallskip}
&$j_1=0$&$j_1=1$&$j_1=2$\\
\noalign{\smallskip}\hline\noalign{\smallskip}
$j_2=0$&$-2$&$\frac{14}{15}$&$\frac{-2}{15}$\\
\noalign{\smallskip}
$j_2=1$&$\frac{-2}{15}$&$\frac{-2}{15}$&$\frac{6}{35}$\\
\noalign{\smallskip}
$j_2=2$&$\frac{2}{5}$&$\frac{-22}{105}$&$\frac{-2}{105}$\\
\noalign{\smallskip}\hline\noalign{\smallskip}
\end{tabular}
\end{table}

\begin{table}
\centering
\caption{Coefficients $\bar C_{1j_2j_1}^{001}$}
\label{tab:5.21}      
\begin{tabular}{p{1.3cm}p{1.3cm}p{1.3cm}p{1.3cm}}
\hline\noalign{\smallskip}
&$j_1=0$&$j_1=1$&$j_1=2$\\
\noalign{\smallskip}\hline\noalign{\smallskip}
$j_2=0$&$\frac{-6}{5}$&$\frac{22}{45}$&$\frac{-2}{105}$\\
\noalign{\smallskip}
$j_2=1$&$\frac{-2}{9}$&$\frac{-2}{105}$&$\frac{26}{315}$\\
\noalign{\smallskip}
$j_2=2$&$\frac{22}{105}$&$\frac{-38}{315}$&$\frac{-2}{315}$\\
\noalign{\smallskip}\hline\noalign{\smallskip}
\end{tabular}
\end{table}

\begin{table}
\centering
\caption{Coefficients $\bar C_{2j_2j_1}^{001}$}
\label{tab:5.22}      
\begin{tabular}{p{1.3cm}p{1.3cm}p{1.3cm}p{1.3cm}}
\hline\noalign{\smallskip}
&$j_1=0$&$j_1=1$&$j_1=2$\\
\noalign{\smallskip}\hline\noalign{\smallskip}
$j_2=0$&$\frac{-2}{5}$&$\frac{2}{21}$&$\frac{4}{105}$\\
\noalign{\smallskip}
$j_2=1$&$\frac{-22}{105}$&$\frac{4}{105}$&$\frac{2}{105}$\\
\noalign{\smallskip}
$j_2=2$&$0$&$\frac{-2}{105}$&$0$\\
\noalign{\smallskip}\hline\noalign{\smallskip}
\end{tabular}
\vspace{2mm}
\end{table}

\begin{table}
\centering
\caption{Coefficients $\bar C_{0j_2j_1}^{100}$}
\label{tab:5.23}      
\begin{tabular}{p{1.3cm}p{1.3cm}p{1.3cm}p{1.3cm}}
\hline\noalign{\smallskip}
&$j_1=0$&$j_1=1$&$j_1=2$\\
\noalign{\smallskip}\hline\noalign{\smallskip}
$j_2=0$&$\frac{-2}{3}$&$\frac{2}{15}$&$\frac{2}{15}$\\
\noalign{\smallskip}
$j_2=1$&$\frac{-2}{15}$&$\frac{-2}{45}$&$\frac{2}{35}$\\
\noalign{\smallskip}
$j_2=2$&$\frac{2}{15}$&$\frac{-2}{35}$&$\frac{-4}{105}$\\
\noalign{\smallskip}\hline\noalign{\smallskip}
\end{tabular}
\vspace{2mm}
\end{table}

\begin{table}
\centering
\caption{Coefficients $\bar C_{1j_2j_1}^{100}$}
\label{tab:5.24}      
\begin{tabular}{p{1.3cm}p{1.3cm}p{1.3cm}p{1.3cm}}
\hline\noalign{\smallskip}
&$j_1=0$&$j_1=1$&$j_1=2$\\
\noalign{\smallskip}\hline\noalign{\smallskip}
$j_2=0$&$\frac{-2}{5}$&$\frac{2}{45}$&$\frac{2}{21}$\\
\noalign{\smallskip}
$j_2=1$&$\frac{-2}{15}$&$\frac{-2}{105}$&$\frac{4}{105}$\\
\noalign{\smallskip}
$j_2=2$&$\frac{2}{35}$&$\frac{-2}{63}$&$\frac{-2}{105}$\\
\noalign{\smallskip}\hline\noalign{\smallskip}
\end{tabular}
\vspace{2mm}
\end{table}

\begin{table}
\centering
\caption{Coefficients $\bar C_{2j_2j_1}^{100}$}
\label{tab:5.25}      
\begin{tabular}{p{1.3cm}p{1.3cm}p{1.3cm}p{1.3cm}}
\hline\noalign{\smallskip}
&$j_1=0$&$j_1=1$&$j_1=2$\\
\noalign{\smallskip}\hline\noalign{\smallskip}
$j_2=0$&$\frac{-2}{15}$&$\frac{-2}{105}$&$\frac{4}{105}$\\
\noalign{\smallskip}
$j_2=1$&$\frac{-2}{21}$&$\frac{-2}{315}$&$\frac{2}{105}$\\
\noalign{\smallskip}
$j_2=2$&$\frac{-2}{105}$&$\frac{-2}{315}$&$0$\\
\noalign{\smallskip}\hline\noalign{\smallskip}
\end{tabular}
\vspace{2mm}
\end{table}

\begin{table}
\centering
\caption{Coefficients $\bar C_{0j_2j_1}^{010}$}
\label{tab:5.26}      
\begin{tabular}{p{1.3cm}p{1.3cm}p{1.3cm}p{1.3cm}}
\hline\noalign{\smallskip}
&$j_1=0$&$j_1=1$&$j_1=2$\\
\noalign{\smallskip}\hline\noalign{\smallskip}
$j_2=0$&$\frac{-4}{3}$&$\frac{8}{15}$&$0$\\
\noalign{\smallskip}
$j_2=1$&$\frac{-4}{15}$&$0$&$\frac{8}{105}$\\
\noalign{\smallskip}
$j_2=2$&$\frac{4}{15}$&$\frac{-16}{105}$&$0$\\
\noalign{\smallskip}\hline\noalign{\smallskip}
\end{tabular}
\vspace{2mm}
\end{table}

\begin{table}
\centering
\caption{Coefficients $\bar C_{1j_2j_1}^{010}$}
\label{tab:5.27}      
\begin{tabular}{p{1.3cm}p{1.3cm}p{1.3cm}p{1.3cm}}
\hline\noalign{\smallskip}
&$j_1=0$&$j_1=1$&$j_1=2$\\
\noalign{\smallskip}\hline\noalign{\smallskip}
$j_2=0$&$\frac{-4}{5}$&$\frac{4}{15}$&$\frac{4}{105}$\\
\noalign{\smallskip}
$j_2=1$&$\frac{-4}{15}$&$\frac{4}{105}$&$\frac{4}{105}$\\
\noalign{\smallskip}
$j_2=2$&$\frac{4}{35}$&$\frac{-8}{105}$&$0$\\
\noalign{\smallskip}\hline\noalign{\smallskip}
\end{tabular}
\vspace{2mm}
\end{table}

\begin{table}
\centering
\caption{Coefficients $\bar C_{2j_2j_1}^{010}$}
\label{tab:5.28}      
\begin{tabular}{p{1.3cm}p{1.3cm}p{1.3cm}p{1.3cm}}
\hline\noalign{\smallskip}
&$j_1=0$&$j_1=1$&$j_1=2$\\
\noalign{\smallskip}\hline\noalign{\smallskip}
$j_2=0$&$\frac{-4}{15}$&$\frac{4}{105}$&$\frac{4}{105}$\\
\noalign{\smallskip}
$j_2=1$&$\frac{-4}{21}$&$\frac{4}{105}$&$\frac{4}{315}$\\
\noalign{\smallskip}
$j_2=2$&$\frac{-4}{105}$&$0$&$0$\\
\noalign{\smallskip}\hline\noalign{\smallskip}
\end{tabular}
\vspace{2mm}
\end{table}

\begin{table}
\centering
\caption{Coefficients $\bar C_{000j_2j_1}$}
\label{tab:5.29}      
\begin{tabular}{p{1.3cm}p{1.3cm}p{1.3cm}}
\hline\noalign{\smallskip}
&$j_1=0$&$j_1=1$\\
\noalign{\smallskip}\hline\noalign{\smallskip}
$j_2=0$&$\frac{4}{15}$&$\frac{-8}{45}$\\
\noalign{\smallskip}
$j_2=1$&$\frac{-4}{45}$&$\frac{8}{105}$\\
\noalign{\smallskip}\hline\noalign{\smallskip}
\end{tabular}
\vspace{2mm}
\end{table}

\begin{table}
\centering
\caption{Coefficients $\bar C_{010j_2j_1}$}
\label{tab:5.30}      
\begin{tabular}{p{1.3cm}p{1.3cm}p{1.3cm}}
\hline\noalign{\smallskip}
&$j_1=0$&$j_1=1$\\
\noalign{\smallskip}\hline\noalign{\smallskip}
$j_2=0$&$\frac{4}{45}$&$\frac{-16}{315}$\\
\noalign{\smallskip}
$j_2=1$&$\frac{-4}{315}$&$\frac{4}{315}$\\
\noalign{\smallskip}\hline\noalign{\smallskip}
\end{tabular}
\vspace{2mm}
\end{table}

\begin{table}
\centering
\caption{Coefficients $\bar C_{110j_2j_1}$}
\label{tab:5.31}      
\begin{tabular}{p{1.3cm}p{1.3cm}p{1.3cm}}
\hline\noalign{\smallskip}
&$j_1=0$&$j_1=1$\\
\noalign{\smallskip}\hline\noalign{\smallskip}
$j_2=0$&$\frac{8}{105}$&$\frac{-2}{45}$\\
\noalign{\smallskip}
$j_2=1$&$\frac{-4}{315}$&$\frac{4}{315}$\\
\noalign{\smallskip}\hline\noalign{\smallskip}
\end{tabular}
\vspace{2mm}
\end{table}

\begin{table}
\centering
\caption{Coefficients $\bar C_{011j_2j_1}$}
\label{tab:5.32}      
\begin{tabular}{p{1.3cm}p{1.3cm}p{1.3cm}}
\hline\noalign{\smallskip}
&$j_1=0$&$j_1=1$\\
\noalign{\smallskip}\hline\noalign{\smallskip}
$j_2=0$&$\frac{8}{315}$&$\frac{-4}{315}$\\
\noalign{\smallskip}
$j_2=1$&$0$&$\frac{2}{945}$\\
\noalign{\smallskip}\hline\noalign{\smallskip}
\end{tabular}
\vspace{2mm}
\end{table}

\begin{table}
\centering
\caption{Coefficients $\bar C_{001j_2j_1}$}
\label{tab:5.33}      
\begin{tabular}{p{1.3cm}p{1.3cm}p{1.3cm}}
\hline\noalign{\smallskip}
&$j_1=0$&$j_1=1$\\
\noalign{\smallskip}\hline\noalign{\smallskip}
$j_2=0$&$0$&$\frac{4}{315}$\\
\noalign{\smallskip}
$j_2=1$&$\frac{8}{315}$&$\frac{-2}{105}$\\
\noalign{\smallskip}\hline\noalign{\smallskip}
\end{tabular}
\vspace{2mm}
\end{table}

\begin{table}
\centering
\caption{Coefficients $\bar C_{100j_2j_1}$}
\label{tab:5.34}      
\begin{tabular}{p{1.3cm}p{1.3cm}p{1.3cm}}
\hline\noalign{\smallskip}
&$j_1=0$&$j_1=1$\\
\noalign{\smallskip}\hline\noalign{\smallskip}
$j_2=0$&$\frac{8}{45}$&$\frac{-4}{35}$\\
\noalign{\smallskip}
$j_2=1$&$\frac{-16}{315}$&$\frac{2}{45}$\\
\noalign{\smallskip}\hline\noalign{\smallskip}
\end{tabular}
\end{table}

\begin{table}
\centering
\caption{Coefficients $\bar C_{101j_2j_1}$}
\label{tab:5.35}      
\begin{tabular}{p{1.3cm}p{1.3cm}p{1.3cm}}
\hline\noalign{\smallskip}
&$j_1=0$&$j_1=1$\\
\noalign{\smallskip}\hline\noalign{\smallskip}
$j_2=0$&$\frac{4}{315}$&$0$\\
\noalign{\smallskip}
$j_2=1$&$\frac{4}{315}$&$\frac{-8}{945}$\\
\noalign{\smallskip}\hline\noalign{\smallskip}
\end{tabular}
\end{table}

\begin{table}
\centering
\caption{Coefficients $\bar C_{111j_2j_1}$}
\label{tab:5.36}      
\begin{tabular}{p{1.3cm}p{1.3cm}p{1.3cm}}
\hline\noalign{\smallskip}
&$j_1=0$&$j_1=1$\\
\noalign{\smallskip}\hline\noalign{\smallskip}
$j_2=0$&$\frac{2}{105}$&$\frac{-8}{945}$\\
\noalign{\smallskip}
$j_2=1$&$\frac{2}{945}$&$0$\\
\noalign{\smallskip}\hline\noalign{\smallskip}
\end{tabular}
\end{table}

Assume that 
$q_3=2,$ $q_4=1.$  In 
Tables 5.20--5.36 we have the exact values of 
Fo\-u\-ri\-er--Le\-gen\-dre
coefficients 
$\bar C_{j_3j_2j_1}^{001},$
$\bar C_{j_3j_2j_1}^{010},$
$\bar C_{j_3j_2j_1}^{100}$
($j_1,j_2,j_3=0, 1, 2),$
$\bar C_{j_5j_4j_3j_2j_1}$
($j_1,\ldots,j_5=0, 1$).

In the case of pairwise different 
$i_1, \ldots, i_5$ from 
Tables  
5.20--5.36 we obtain

\vspace{-2mm}
$$
{\sf M}\left\{\left(
I_{(100)T,t}^{(i_1i_2 i_3)}-
I_{(100)T,t}^{(i_1i_2 i_3)q_3}\right)^2\right\}=
$$
$$
=
\frac{(T-t)^{5}}{60}-\sum_{j_1,j_2,j_3=0}^{2}
\left(C_{j_3j_2j_1}^{100}\right)^2\approx 0.00815429(T-t)^5,
$$

\vspace{3mm}

$$
{\sf M}\left\{\left(
I_{(010)T,t}^{(i_1i_2 i_3)}-
I_{(010)T,t}^{(i_1i_2 i_3)q_3}\right)^2\right\}=
$$
$$
=
\frac{(T-t)^{5}}{20}-\sum_{j_1,j_2,j_3=0}^{2}
\left(C_{j_3j_2j_1}^{010}\right)^2\approx 0.01739030(T-t)^5,
$$

\vspace{3mm}

$$
{\sf M}\left\{\left(
I_{(001)T,t}^{(i_1i_2 i_3)}-
I_{(001)T,t}^{(i_1i_2 i_3)q_3}\right)^2\right\}=
$$
$$
=
\frac{(T-t)^5}{10}-\sum_{j_1,j_2,j_3=0}^{2}
\left(C_{j_3j_2j_1}^{001}\right)^2
\approx 0.02528010(T-t)^5,
$$

\vspace{3mm}

$$
{\sf M}\left\{\left(
I_{(00000)T,t}^{(i_1i_2i_3i_4 i_5)}-
I_{(00000)T,t}^{(i_1i_2i_3i_4 i_5)q_4}\right)^2\right\}=
$$
$$
=
\frac{(T-t)^5}{120}-\sum_{j_1,j_2,j_3,j_4,j_5=0}^{1}
C_{j_5j_4j_3j_2j_1}^2\approx 0.00759105(T-t)^5.
$$

\vspace{6mm}

Note that from Theorem 1.4 for $k=5$ we have

\vspace{-5mm}
$$
{\sf M}\left\{\left(
I_{(00000)T,t}^{(i_1i_2 i_3 i_4 i_5)}-
I_{(00000)T,t}^{(i_1i_2 i_3 i_4 i_5)q_4}\right)^2\right\}\le
120\left(\frac{(T-t)^{5}}{120}-\sum_{j_1,j_2,j_3,j_4,j_5=0}^{q_4}
C_{j_5j_4j_3j_2j_1}^2\right),
$$

\noindent
where $i_1, \ldots, i_5=1,\ldots,m$.

Moreover, from Theorem 1.4 we obtain the following useful estimates

\vspace{-1mm}
$$
\vspace{2mm}
{\sf M}\left\{\left(
I_{(01)T,t}^{(i_1i_2)}-
I_{(01)T,t}^{(i_1i_2)q}\right)^2\right\}\le
2\Biggl(\frac{(T-t)^{4}}{4}-\sum_{j_1,j_2=0}^{q}
\left(C_{j_2j_1}^{01}\right)^2\Biggr),
$$

\vspace{-3mm}
$$
{\sf M}\left\{\left(
I_{(10)T,t}^{(i_1i_2)}-
I_{(10)T,t}^{(i_1i_2)q}\right)^2\right\}\le
2\Biggl(\frac{(T-t)^{4}}{12}-\sum_{j_1,j_2=0}^{q}
\left(C_{j_2j_1}^{10}\right)^2\Biggr),
$$

\vspace{-3mm}
$$
{\sf M}\left\{\left(
I_{(100)T,t}^{(i_1i_2 i_3)}-
I_{(100)T,t}^{(i_1i_2 i_3)q}\right)^2\right\}\le
6\Biggl(\frac{(T-t)^{5}}{60}-\sum_{j_1,j_2,j_3=0}^{q}
\left(C_{j_3j_2j_1}^{100}\right)^2\Biggr),
$$

\vspace{-3mm}
$$
{\sf M}\left\{\left(
I_{(010)T,t}^{(i_1i_2 i_3)}-
I_{(010)T,t}^{(i_1i_2 i_3)q}\right)^2\right\}\le
6\Biggl(\frac{(T-t)^{5}}{20}-\sum_{j_1,j_2,j_3=0}^{q}
\left(C_{j_3j_2j_1}^{010}\right)^2\Biggr),
$$

\vspace{-3mm}
$$
{\sf M}\left\{\left(
I_{(001)T,t}^{(i_1i_2 i_3)}-
I_{(001)T,t}^{(i_1i_2 i_3)q}\right)^2\right\}\le
6\Biggl(\frac{(T-t)^5}{10}-\sum_{j_1,j_2,j_3=0}^{q}
\left(C_{j_3j_2j_1}^{001}\right)^2\Biggr),
$$

\vspace{-1mm}
$$
{\sf M}\left\{\left(
I_{(20))T,t}^{(i_1i_2)}-
I_{(20)T,t}^{(i_1i_2)q}\right)^2\right\}\le
2\Biggl(\frac{(T-t)^6}{30}-\sum_{j_2,j_1=0}^{q}
\left(C_{j_2j_1}^{20}\right)^2\Biggr),
$$

\vspace{-1mm}
$$
{\sf M}\left\{\left(
I_{(11)T,t}^{(i_1i_2)}-
I_{(11)T,t}^{(i_1i_2)q}\right)^2\right\}\le
2\Biggl(\frac{(T-t)^6}{18}-\sum_{j_2,j_1=0}^{q}
\left(C_{j_2j_1}^{11}\right)^2\Biggr),
$$

\vspace{-1mm}
$$
{\sf M}\left\{\left(
I_{(02)T,t}^{(i_1i_2)}-
I_{(02)T,t}^{(i_1i_2)q}\right)^2\right\}\le
2\Biggl(\frac{(T-t)^6}{6}-\sum_{j_2,j_1=0}^{q}
\left(C_{j_2j_1}^{02}\right)^2\Biggr),
$$

\vspace{-3mm}
$$
{\sf M}\left\{\left(
I_{(1000)T,t}^{(i_1i_2 i_3i_4)}-
I_{(1000)T,t}^{(i_1i_2 i_3i_4)q}\right)^2\right\}\le
24\Biggl(\frac{(T-t)^{6}}{360}-\sum_{j_1,j_2,j_3, j_4=0}^{q}
\left(C_{j_4j_3j_2j_1}^{1000}\right)^2\Biggr),
$$

\vspace{-3mm}
$$
{\sf M}\left\{\left(
I_{(0100)T,t}^{(i_1i_2 i_3i_4)}-
I_{(0100)T,t}^{(i_1i_2 i_3i_4)q}\right)^2\right\}\le
24\Biggl(\frac{(T-t)^{6}}{120}-\sum_{j_1,j_2,j_3, j_4=0}^{q}
\left(C_{j_4j_3j_2j_1}^{0100}\right)^2\Biggr),
$$

\vspace{-3mm}
$$
{\sf M}\left\{\left(
I_{(0010)T,t}^{(i_1i_2 i_3i_4)}-
I_{(0010)T,t}^{(i_1i_2 i_3 i_4)q}\right)^2\right\}\le
24\Biggl(\frac{(T-t)^6}{60}-\sum_{j_1,j_2,j_3, j_4=0}^{q}
\left(C_{j_4j_3j_2j_1}^{0010}\right)^2\Biggr),
$$

\vspace{-3mm}
$$
{\sf M}\left\{\left(
I_{(0001)T,t}^{(i_1i_2 i_3 i_4)}-
I_{(0001)T,t}^{(i_1i_2 i_3 i_4)q}\right)^2\right\}\le
24\Biggl(\frac{(T-t)^6}{36}-\sum_{j_1,j_2,j_3, j_4=0}^{q}
\left(C_{j_4j_3j_2j_1}^{0001}\right)^2\Biggr),
$$

\newpage
\noindent
$$
{\sf M}\left\{\left(
I_{(000000)T,t}^{(i_1 i_2 i_3 i_4 i_5 i_6)}-
I_{(000000)T,t}^{(i_1 i_2 i_3 i_4 i_5 i_6)q}\right)^2\right\}\le
$$
$$
\le
720\left(\frac{(T-t)^{6}}{720}-\sum_{j_1,j_2,j_3,j_4,j_5,j_6=0}^{q}
C_{j_6 j_5 j_4 j_3 j_2 j_1}^2\right).
$$

\vspace{4mm}

In addition, from Theorem 1.3 for $k=2$ we get

\vspace{-1mm}
$$
{\sf M}\biggl\{\left(I_{(10)T,t}^{(i_1 i_2)}-
I_{(10)T,t}^{(i_1 i_2)q}
\right)^2\biggr\}=
$$
$$
=\frac{(T-t)^4}{12}
-\sum_{j_1,j_2=0}^q
\left(C_{j_2j_1}^{10}\right)^2-
\sum_{j_1,j_2=0}^q C_{j_2j_1}^{10}C_{j_1j_2}^{10}\ \ \ (i_1=i_2),
$$

\vspace{1mm}

$$
{\sf M}\biggl\{\left(I_{(10)T,t}^{(i_1 i_2)}-
I_{(10)T,t}^{(i_1 i_2)q}
\right)^2\biggr\}=\frac{(T-t)^4}{12}
-\sum_{j_1,j_2=0}^q
\left(C_{j_2j_1}^{10}\right)^2\ \ \ (i_1\ne i_2),
$$

\vspace{1mm}

$$
{\sf M}\biggl\{\left(I_{(01)T,t}^{(i_1 i_2)}-
I_{(01)T,t}^{(i_1 i_2)q}
\right)^2\biggr\}=
$$
$$
=\frac{(T-t)^4}{4}
-\sum_{j_1,j_2=0}^q
\left(C_{j_2j_1}^{01}\right)^2-
\sum_{j_1,j_2=0}^q C_{j_2j_1}^{01}C_{j_1j_2}^{01}\ \ \ (i_1=i_2),
$$

\vspace{1mm}

$$
{\sf M}\biggl\{\left(I_{(01)T,t}^{(i_1 i_2)}-
I_{(01)T,t}^{(i_1 i_2)q}
\right)^2\biggr\}=\frac{(T-t)^4}{4}
-\sum_{j_1,j_2=0}^q
\left(C_{j_2j_1}^{01}\right)^2\ \ \ (i_1\ne i_2),
$$

\vspace{1mm}

$$
{\sf M}\biggl\{\left(I_{(20)T,t}^{(i_1 i_2)}-
I_{(20)T,t}^{(i_1 i_2)q}
\right)^2\biggr\}=
$$
$$
=\frac{(T-t)^6}{30}
-\sum_{j_1,j_2=0}^q
\left(C_{j_2j_1}^{20}\right)^2-
\sum_{j_1,j_2=0}^q C_{j_2j_1}^{20}C_{j_1j_2}^{20}\ \ \ (i_1=i_2),
$$

\vspace{1mm}

$$
{\sf M}\biggl\{\left(I_{(20)T,t}^{(i_1 i_2)}-
I_{(20)T,t}^{(i_1 i_2)q}
\right)^2\biggr\}=\frac{(T-t)^6}{30}
-\sum_{j_1,j_2=0}^q
\left(C_{j_2j_1}^{20}\right)^2\ \ \ (i_1\ne i_2),
$$

\vspace{1mm}

$$
{\sf M}\biggl\{\left(I_{(11)T,t}^{(i_1 i_2)}-
I_{(11)T,t}^{(i_1 i_2)q}
\right)^2\biggr\}=
$$
$$
=\frac{(T-t)^6}{18}
-\sum_{j_1,j_2=0}^q
\left(C_{j_2j_1}^{11}\right)^2-
\sum_{j_1,j_2=0}^q C_{j_2j_1}^{11}C_{j_1j_2}^{11}\ \ \ (i_1=i_2),
$$

\vspace{1mm}

$$
{\sf M}\biggl\{\left(I_{(11)T,t}^{(i_1 i_2)}-
I_{(11)T,t}^{(i_1 i_2)q}
\right)^2\biggr\}=\frac{(T-t)^6}{18}
-\sum_{j_1,j_2=0}^q
\left(C_{j_2j_1}^{11}\right)^2\ \ \ (i_1\ne i_2),
$$

\vspace{1mm}

$$
{\sf M}\biggl\{\left(I_{(02)}^{(i_1 i_2)}-
I_{(02)T,t}^{(i_1 i_2)q}
\right)^2\biggr\}=
$$
$$
=\frac{(T-t)^6}{6}
-\sum_{j_1,j_2=0}^q
\left(C_{j_2j_1}^{02}\right)^2-
\sum_{j_1,j_2=0}^q C_{j_2j_1}^{02}C_{j_1j_2}^{02}\ \ \ (i_1=i_2),
$$

\vspace{1mm}

$$
{\sf M}\biggl\{\left(I_{(02)T,t}^{(i_1 i_2)}-
I_{(02)T,t}^{(i_1 i_2)q}
\right)^2\biggr\}=\frac{(T-t)^6}{6}
-\sum_{j_1,j_2=0}^q
\left(C_{j_2j_1}^{02}\right)^2\ \ \ (i_1\ne i_2).
$$

\vspace{1mm}

Clearly, expansions for iterated Stratonovich stochastic integrals
(see Theorems 1.1, 2.1--2.9, 2.33--2.36, 2.50, 2.51, 2.62--2.65) are simpler 
than expansions for
iterated It\^{o} stochastic integrals (see Theorems 1.1, 1.2, 1.16 and
(\ref{a1})--(\ref{a7})). However, the calculation of the mean-square
approximation error for iterated Stratonovich
stochastic integrals turns out to be much more difficult than for 
iterated It\^{o} stochastic integrals.
Below we consider how we can estimate or calculate exactly
(for some particular cases)
the mean-square
approximation error for iterated Stratonovich
stochastic integrals.

Consider the iterated Stratonovich stochastic integral
of multiplicity 2
$$
J^{*}[\psi^{(2)}]_{T,t}=
{\int\limits_t^{*}}^T
\psi_2(t_2)
{\int\limits_t^{*}}^{t_{2}}
\psi_1(t_1) d{\bf w}_{t_1}^{(i_1)}
d{\bf w}_{t_2}^{(i_1)}\ \ \ (i_1=1,\ldots,m),
$$
where 
$\psi_1(\tau),$ $\psi_2(\tau)$ are 
continuously differentiable functions on $[t, T]$.

By Theorem 2.2 we have
$$
J^{*}[\psi^{(2)}]_{T,t}
=
\hbox{\vtop{\offinterlineskip\halign{
\hfil#\hfil\cr
{\rm l.i.m.}\cr
$\stackrel{}{{}_{p\to \infty}}$\cr
}} }
\sum_{j_1,j_2=0}^{p}
C_{j_2j_1}\zeta_{j_1}^{(i_1)}\zeta_{j_2}^{(i_1)}.
$$

Consider the following approximation of the stochastic integral 
$J^{*}[\psi^{(2)}]_{T,t}$
$$
J^{*}[\psi^{(2)}]_{T,t}^q
=\sum_{j_1,j_2=0}^{q}
C_{j_2j_1}\zeta_{j_1}^{(i_1)}\zeta_{j_2}^{(i_1)}.
$$

According to the standard relation between 
Stratonovich and It\^{o} stochastic integrals (see (\ref{uyes1}))
and (\ref{uyes2}), we obtain
$$
{\sf M}\left\{\left(J^{*}[\psi^{(2)}]_{T,t}-
J^{*}[\psi^{(2)}]_{T,t}^q\right)^2\right\}=
$$
$$
=
{\sf M}\left\{\left(J[\psi^{(2)}]_{T,t}
+\frac{1}{2}\int\limits_t^T\psi_1(s)\psi_2(s)ds-
\sum_{j_1,j_2=0}^{q}
C_{j_2j_1}\zeta_{j_1}^{(i_1)}\zeta_{j_2}^{(i_1)}\right)^2\right\}=
$$
$$
={\sf M}\left\{\left(J[\psi^{(2)}]_{T,t}
-J[\psi^{(2)}]_{T,t}^q
+\frac{1}{2}\int\limits_t^T\psi_1(s)\psi_2(s)ds-
\sum_{j_1=0}^{q} 
C_{j_1j_1}\right)^2\right\}=
$$
$$
={\sf M}\left\{\left(J[\psi^{(2)}]_{T,t}
-J[\psi^{(2)}]_{T,t}^q\right)^2\right\}+
\left(\frac{1}{2}\int\limits_t^T\psi_1(s)\psi_2(s)ds-
\sum_{j_1=0}^{q} 
C_{j_1j_1}\right)^2=
$$
$$
=\int\limits_{[t,T]^2}K^2(t_1,t_2)dt_1dt_2
-\sum_{j_1,j_2=0}^q
C_{j_2j_1}^2-
\sum_{j_1,j_2=0}^q
C_{j_2j_1}C_{j_1j_2}+
$$
$$
+
\left(\frac{1}{2}\int\limits_t^T\psi_1(s)\psi_2(s)ds-
\sum_{j_1=0}^{q} 
C_{j_1j_1}\right)^2,
$$

\noindent
where 
$$
J[\psi^{(2)}]_{T,t}^q
=
\sum_{j_1,j_2=0}^q
C_{j_2j_1}\zeta_{j_1}^{(i_1)}\zeta_{j_2}^{(i_1)}
-\sum_{j_1=0}^{q}
C_{j_1j_1}
$$

\vspace{2mm}
\noindent 
is the approximation (see (\ref{a2}))
of the iterated It\^{o} stochastic integral
$$
J[\psi^{(2)}]_{T,t}=
\int\limits_t^{T}\psi_2(t_2)\int\limits_t^{t_{2}}
\psi_1(t_1) d{\bf w}_{t_1}^{(i_1)}
d{\bf w}_{t_2}^{(i_1)}\ \ \ (i_1=1,\ldots,m).
$$

It is not difficult to see that the value
$$
{\sf M}\left\{\left(J^{*}[\psi^{(2)}]_{T,t}-
J^{*}[\psi^{(2)}]_{T,t}^q\right)^2\right\}
$$

\noindent
is greater than the value
$$
{\sf M}\left\{\left(J[\psi^{(2)}]_{T,t}
-J[\psi^{(2)}]_{T,t}^q\right)^2\right\}
$$

\noindent
by
$$
E^{(i_1)}_q=\left(\frac{1}{2}\int\limits_t^T\psi_1(s)\psi_2(s)ds-
\sum_{j_1=0}^{q} 
C_{j_1j_1}\right)^2.
$$

For some particular cases $E^{(i_1)}_q=0.$ For example,
for the case $\psi_1(\tau),$ $\psi_2(\tau)$ $\equiv 1$ 
($\{\phi_j(x)\}_{j=0}^{\infty}$ is a complete
orthonormal system of Legendre polynomials
or trigonometric functions in the space $L_2([t, T])$)
we have
$$
\sum_{j_1=0}^{q}C_{j_1j_1}=
\frac{1}{2}\sum_{j_1=0}^{q}\left(C_{j_1}\right)^2=
\frac{1}{2}\left(C_{0}\right)^2=\frac{1}{2}(T-t)=\frac{1}{2}
\int\limits_t^T ds.
$$

However, $E^{(i_1)}_q\ne 0$ in a general case.

Consider the following iterated Stratonovich stochastic integral
of multiplicity 3
$$
I_{(000)T,t}^{*(i_1i_2i_3)}=
{\int\limits_t^{*}}^T
{\int\limits_t^{*}}^{t_{3}}
{\int\limits_t^{*}}^{t_{2}}
d{\bf w}_{t_1}^{(i_1)}
d{\bf w}_{t_2}^{(i_2)}
d{\bf w}_{t_3}^{(i_3)}\ \ \ (i_1, i_2, i_3=1,\ldots,m).
$$

Taking into account the standard relations between 
It\^{o} and Stratonovich stochastic integrals (see (\ref{uyes3})) and
Theorem 1.1 (the case $k=3$) together with Theorem 2.8, we obtain
$$
{\sf M}\left\{\left(I_{(000)T,t}^{*(i_1i_2i_3)}-
I_{(000)T,t}^{*(i_1i_2i_3)q}\right)^2\right\}=
$$
$$
={\sf M}\hspace{-0.8mm}\left\{\hspace{-1.2mm}\left(
\hspace{-1mm}I_{(000)T,t}^{(i_1i_2i_3)}+
\frac{1}{2}{\bf 1}_{\{i_1=i_2\}}
\int\limits_t^T \hspace{-1mm}
\int\limits_t^{\tau}dsd{\bf w}_{\tau}^{(i_3)}
+\frac{1}{2}{\bf 1}_{\{i_2=i_3\}}
\int\limits_t^T \hspace{-1mm}\int\limits_t^{\tau}d{\bf w}_{s}^{(i_1)}d\tau-
I_{(000)T,t}^{*(i_1i_2i_3)q}\hspace{-0.5mm}\right)^{\hspace{-2mm}2}\right\}
$$
$$
={\sf M}\left\{\left(
I_{(000)T,t}^{(i_1i_2i_3)}-I_{(000)T,t}^{(i_1i_2i_3)q}+
I_{(000)T,t}^{(i_1i_2i_3)q}+
{\bf 1}_{\{i_1=i_2\}}
\frac{1}{2}\int\limits_t^{\stackrel{~}{T}}
\int\limits_{t}^{\stackrel{~}{\tau}}
dsd{\bf w}_{\tau}^{(i_3)}
\right.\right.+
$$
\begin{equation}
\label{tango3}
\left.\left.
+{\bf 1}_{\{i_2=i_3\}}
\frac{1}{2}\int\limits_t^T\int\limits_t^{\tau}d{\bf w}_{s}^{(i_1)}d\tau-
I_{(000)T,t}^{*(i_1i_2i_3)q}\right)^2\right\},
\end{equation}

\vspace{2mm}
\noindent
where the approximations $I_{(000)T,t}^{*(i_1i_2i_3)q},$
$I_{(000)T,t}^{(i_1i_2i_3)q}$ are defined by 
the relations (see (\ref{good1}), (\ref{zzz1})) 

\newpage
\noindent
$$
I_{(000)T,t}^{(i_1i_2i_3)q}=
\sum_{j_1,j_2,j_3=0}^{q}
C_{j_3j_2j_1}\Biggl(
\zeta_{j_1}^{(i_1)}\zeta_{j_2}^{(i_2)}\zeta_{j_3}^{(i_3)}
-{\bf 1}_{\{i_1=i_2\}}
{\bf 1}_{\{j_1=j_2\}}
\zeta_{j_3}^{(i_3)}-
\Biggr.
$$
\begin{equation}
\label{tango1}
\Biggl.
-{\bf 1}_{\{i_2=i_3\}}
{\bf 1}_{\{j_2=j_3\}}
\zeta_{j_1}^{(i_1)}-
{\bf 1}_{\{i_1=i_3\}}
{\bf 1}_{\{j_1=j_3\}}
\zeta_{j_2}^{(i_2)}\Biggr),
\end{equation}

\vspace{-2mm}
\begin{equation}
\label{tango2}
I_{(000)T,t}^{*(i_1i_2i_3)q}=
\sum_{j_1,j_2,j_3=0}^{q}
C_{j_3j_2j_1}
\zeta_{j_1}^{(i_1)}\zeta_{j_2}^{(i_2)}\zeta_{j_3}^{(i_3)}.
\end{equation}

\vspace{3mm}

Substituting (\ref{tango1}) and (\ref{tango2}) into (\ref{tango3}) yields
$$
{\sf M}\left\{\left(I_{(000)T,t}^{*(i_1i_2i_3)}-
I_{(000)T,t}^{*(i_1i_2i_3)q}\right)^2\right\}=
$$
$$
={\sf M}\left\{\left(I_{(000)T,t}^{(i_1i_2i_3)}-
I_{(000)T,t}^{(i_1i_2i_3)q}+
{\bf 1}_{\{i_1=i_2\}}
\left(\frac{1}{2}\int\limits_t^T
\int\limits_t^{\tau}dsd{\bf w}_{\tau}^{(i_3)}-
\sum_{j_1,j_3=0}^{q}
C_{j_3j_1j_1}
\zeta_{j_3}^{(i_3)}\right)+\right.\right.
$$
\begin{equation}
\label{dest10}
\left.\left.+{\bf 1}_{\{i_2=i_3\}}\left(
\frac{1}{2}\int\limits_t^T\int\limits_t^{\tau}d{\bf w}_{s}^{(i_1)}d\tau-
\hspace{-1.6mm}\sum_{j_1,j_3=0}^{q}
C_{j_3j_3j_1}
\zeta_{j_1}^{(i_1)}\right)
\hspace{-0.75mm}-{\bf 1}_{\{i_1=i_3\}}
\sum_{j_1,j_2=0}^{q}
C_{j_1j_2j_1}
\zeta_{j_2}^{(i_2)}\right)^{\hspace{-1.6mm}2}\right\}\hspace{-1mm}\le
\end{equation}
\begin{equation}
\label{tango4}
\le 4 \biggl({\sf M}\left\{\left(I_{(000)T,t}^{(i_1i_2i_3)}-
I_{(000)T,t}^{(i_1i_2i_3)q}\right)^2\right\}+
{\bf 1}_{\{i_1=i_2\}}F^{(i_3)}_q+
{\bf 1}_{\{i_2=i_3\}}G^{(i_1)}_q+
{\bf 1}_{\{i_1=i_3\}}H^{(i_2)}_q\biggr),
\end{equation}

\vspace{1mm}
\noindent
where 
\begin{equation}
\label{uyes10}
F^{(i_3)}_q={\sf M}\left\{\left(
\frac{1}{2}\int\limits_t^T\int\limits_t^{\tau}dsd{\bf w}_{\tau}^{(i_3)}-
\sum_{j_1,j_3=0}^{q}
C_{j_3j_1j_1}
\zeta_{j_3}^{(i_3)}\right)^2\right\},
\end{equation}
\begin{equation}
\label{uyes11}
G^{(i_1)}_q={\sf M}\left\{\left(
\frac{1}{2}\int\limits_t^T\int\limits_t^{\tau}d{\bf w}_{s}^{(i_1)}d\tau-
\sum_{j_1,j_3=0}^{q}
C_{j_3j_3j_1}
\zeta_{j_1}^{(i_1)}\right)^2\right\},
\end{equation}
\begin{equation}
\label{uyes12}
H^{(i_2)}_q={\sf M}\left\{\left(
\sum_{j_1,j_2=0}^{q}
C_{j_1j_2j_1}
\zeta_{j_2}^{(i_2)}\right)^2\right\}.
\end{equation}

In the cases of Legendre polynomials or trigonometric functions, 
we have (see Theorem 2.8) the equalities
$$
\lim\limits_{q\to\infty}F^{(i_3)}_q=0,\ \ \
\lim\limits_{q\to\infty} G^{(i_1)}_q=0,\ \ \
\lim\limits_{q\to\infty}H^{(i_2)}_q=0.
$$

However, in accordance with (\ref{tango4}) the value
$$
{\sf M}\left\{\left(I_{(000)T,t}^{*(i_1i_2i_3)}-
I_{(000)T,t}^{*(i_1i_2i_3)q}\right)^2\right\}
$$
with a finite $q$ can be estimated by
terms 
of a rather complex structure (see (\ref{uyes10})-(\ref{uyes12})). 
As is easily observed, this peculiarity will also apply to the iterated
Stratonovich stochastic
integrals of multiplicities $k\ge 4$ with the only difference that 
the number of additional terms like
(\ref{uyes10})-(\ref{uyes12})
will be considerably higher and their structure will be more complicated
(the exact calculation of the mean-square error of approximation 
for iterated Stratonovich stochastic integrals of multiplicities 1 to 4
is presented in Sect.~5.5, 5.6).

Therefore, the payment for a relatively simple approximation of 
the iterated Stratonovich
stochastic integrals (Theorems 2.1--2.9, 2.30, 2.33--2.36, 2.50, 2.51, 2.62--2.65) 
in comparison with the 
iterated It\^{o}
stochastic integrals (Theorems 1.1, 1.2, 1.16)
is a much more difficult calculation or estimation procedure 
of their mean-square
approximation errors.

As we mentioned above, on the basis of 
the presented 
approximations of 
iterated Stratonovich stochastic integrals we 
can see that increasing of multiplicities of these integrals 
leads to increasing 
of orders of smallness with respect to $T-t$
in the mean-square sense 
for iterated Stratonovich stochastic integrals
($T-t\ll 1$ 
because the length $T-t$ of integration interval $[t, T]$ 
of the iterated Stratonovich 
stochastic integrals 
plays the role of integration step for the numerical 
methods for It\^{o} SDEs, i.e. $T-t$ is already fairly small).
This leads to a sharp decrease  
of member 
quantities
in the appro\-xi\-ma\-ti\-ons of iterated Stratonovich stochastic 
integrals
which are required for achieving the acceptable accuracy
of approximation.

From (\ref{fff09}) $(i_1\ne i_2)$ we obtain
$$
{\sf M}\left\{\left(I_{(00)T,t}^{*(i_1 i_2)}-
I_{(00)T,t}^{*(i_1 i_2)q}
\right)^2\right\}=\frac{(T-t)^2}{2}
\sum\limits_{i=q+1}^{\infty}\frac{1}{4i^2-1}\le 
$$
\begin{equation}
\label{teacxx}
~~~~\le \frac{(T-t)^2}{2}\int\limits_{q}^{\infty}
\frac{1}{4x^2-1}dx
=-\frac{(T-t)^2}{8}{\rm ln}\left|
1-\frac{2}{2q+1}\right|\le C_1\frac{(T-t)^2}{q},
\end{equation}

\noindent
where constant $C_1$ does not depend on $q.$

It is easy to notice that for a sufficiently
small $T-t$ (recall that $T-t\ll 1$ since it is a step of integration
for the numerical schemes for It\^{o} SDEs) there 
exists a constant $C_2$ such that
\begin{equation}
\label{teac3xx}
~~~~~~~~{\sf M}\left\{\left(
I_{(l_1\ldots l_k)T,t}^{*(i_1\ldots i_k)}-
I_{(l_1\ldots l_k)T,t}^{*(i_1\ldots i_k)q}\right)^2\right\}
\le C_2 {\sf M}\left\{\left(I_{(00)T,t}^{*(i_1 i_2)}-
I_{(00)T,t}^{*(i_1 i_2)q}\right)^2\right\},
\end{equation}
where 
$I_{(l_1\ldots l_k)T,t}^{*(i_1\ldots i_k)q}$
is an approximation of the iterated Stratonovich stochastic integral 
$I_{(l_1\ldots l_k)T,t}^{*(i_1\ldots i_k)}.$

From (\ref{teacxx}) and (\ref{teac3xx}) we finally obtain
\begin{equation}
\label{teac4}
{\sf M}\left\{\left(
I_{(l_1\ldots l_k)T,t}^{*(i_1\ldots i_k)}-
I_{(l_1\ldots l_k)T,t}^{*(i_1\ldots i_k)q}\right)^2\right\}
\le C \frac{(T-t)^2}{q},
\end{equation}
where constant $C$ is independent of $T-t$.

The same idea can be found in \cite{Zapad-3} in the framework of 
the method of approximation of iterated Stratonovich stochastic 
integrals based
on the trigonometric expansion of the
Brownian bridge process.
Note that, in contrast to the estimate (\ref{teac4}), 
the constant $C$ in Theorems 2.38--2.40 does not depend on $p.$

We can get more information about the numbers $q$ (these
numbers are different for different iterated Stratonovich
stochastic integrals)
using the another approach.
Since for pairwise different $i_1,\ldots,i_k=1,\ldots,m$

\vspace{-3mm}
$$
J^{*}[\psi^{(k)}]_{T,t}=J[\psi^{(k)}]_{T,t}\ \ \ \hbox{w.~p.~1,}
$$

\vspace{2mm}
\noindent
where $J[\psi^{(k)}]_{T,t},$ $J^{*}[\psi^{(k)}]_{T,t}$
are defined by (\ref{ito-ito}) and (\ref{str-str}) correspondingly,
then 
for pairwise different 
$i_1,\ldots,i_6=1,\ldots,m$ from Theorem 1.3 we obtain
$$
{\sf M}\left\{\left(
I_{(01)T,t}^{*(i_1i_2)}-
I_{(01)T,t}^{*(i_1i_2)q}\right)^2\right\}=
\frac{(T-t)^{4}}{4}-\sum_{j_1,j_2=0}^{q}
\left(C_{j_2j_1}^{01}\right)^2,
$$
$$
{\sf M}\left\{\left(
I_{(10)T,t}^{*(i_1i_2)}-
I_{(10)T,t}^{*(i_1i_2)q}\right)^2\right\}=
\frac{(T-t)^{4}}{12}-\sum_{j_1,j_2=0}^{q}
\left(C_{j_2j_1}^{10}\right)^2,
$$
\begin{equation}
\label{dest3}
~~~~~~{\sf M}\left\{\left(
I_{(000)T,t}^{*(i_1i_2 i_3)}-
I_{(000)T,t}^{*(i_1i_2 i_3)q}\right)^2\right\}=
\frac{(T-t)^{3}}{6}-\sum_{j_3,j_2,j_1=0}^{q}
C_{j_3j_2j_1}^2,
\end{equation}
\begin{equation}
\label{destxyz1}
~~~~~~{\sf M}\left\{\left(
I_{(0000)T,t}^{*(i_1i_2 i_3 i_4)}-
I_{(0000)T,t}^{*(i_1i_2 i_3 i_4)q}\right)^2\right\}=
\frac{(T-t)^{4}}{24}-\sum_{j_1,j_2,j_3,j_4=0}^{q}
C_{j_4j_3j_2j_1}^2,
\end{equation}
$$
{\sf M}\left\{\left(
I_{(100)T,t}^{*(i_1i_2 i_3)}-
I_{(100)T,t}^{*(i_1i_2 i_3)q}\right)^2\right\}=
\frac{(T-t)^{5}}{60}-\sum_{j_1,j_2,j_3=0}^{q}
\left(C_{j_3j_2j_1}^{100}\right)^2,
$$
$$
{\sf M}\left\{\left(
I_{(010))T,t}^{*(i_1i_2 i_3)}-
I_{(010)T,t}^{*(i_1i_2 i_3)q}\right)^2\right\}=
\frac{(T-t)^{5}}{20}-\sum_{j_1,j_2,j_3=0}^{q}
\left(C_{j_3j_2j_1}^{010}\right)^2,
$$
$$
{\sf M}\left\{\left(
I_{(001)T,t}^{*(i_1i_2 i_3)}-
I_{(001)T,t}^{*(i_1i_2 i_3)q}\right)^2\right\}=
\frac{(T-t)^5}{10}-\sum_{j_1,j_2,j_3=0}^{q}
\left(C_{j_3j_2j_1}^{001}\right)^2,
$$
$$
{\sf M}\left\{\left(
I_{(00000)T,t}^{*(i_1 i_2 i_3 i_4 i_5)}-
I_{(00000)T,t}^{*(i_1 i_2 i_3 i_4 i_5)q}\right)^2\right\}=
\frac{(T-t)^{5}}{120}-\sum_{j_1,j_2,j_3,j_4,j_5=0}^{q}
C_{j_5 i_4 i_3 i_2 j_1}^2,
$$
$$
{\sf M}\left\{\left(
I_{(20)T,t}^{*(i_1i_2)}-
I_{(20)T,t}^{*(i_1i_2)q}\right)^2\right\}=
\frac{(T-t)^6}{30}-\sum_{j_2,j_1=0}^{q}
\left(C_{j_2j_1}^{20}\right)^2,
$$
$$
{\sf M}\left\{\left(
I_{(11)T,t}^{*(i_1i_2)}-
I_{(11)T,t}^{*(i_1i_2)q}\right)^2\right\}=
\frac{(T-t)^6}{18}-\sum_{j_2,j_1=0}^{q}
\left(C_{j_2j_1}^{11}\right)^2,
$$
$$
{\sf M}\left\{\left(
I_{(02)T,t}^{*(i_1i_2)}-
I_{(02)T,t}^{*(i_1i_2)q}\right)^2\right\}=
\frac{(T-t)^6}{6}-\sum_{j_2,j_1=0}^{q}
\left(C_{j_2j_1}^{02}\right)^2,
$$
$$
{\sf M}\left\{\left(
I_{(1000)T,t}^{*(i_1i_2 i_3i_4)}-
I_{(1000)T,t}^{*(i_1i_2 i_3i_4)q}\right)^2\right\}=
\frac{(T-t)^{6}}{360}-\sum_{j_1,j_2,j_3, j_4=0}^{q}
\left(C_{j_4j_3j_2j_1}^{1000}\right)^2,
$$
$$
{\sf M}\left\{\left(
I_{(0100)T,t}^{*(i_1i_2 i_3i_4)}-
I_{(0100)T,t}^{*(i_1i_2 i_3i_4)q}\right)^2\right\}=
\frac{(T-t)^{6}}{120}-\sum_{j_1,j_2,j_3, j_4=0}^{q}
\left(C_{j_4j_3j_2j_1}^{0100}\right)^2,
$$
$$
{\sf M}\left\{\left(
I_{(0010)T,t}^{*(i_1i_2 i_3i_4)}-
I_{(0010)T,t}^{*(i_1i_2 i_3 i_4)q}\right)^2\right\}=
\frac{(T-t)^6}{60}-\sum_{j_1,j_2,j_3, j_4=0}^{q}
\left(C_{j_4j_3j_2j_1}^{0010}\right)^2,
$$
$$
{\sf M}\left\{\left(
I_{(0001)T,t}^{*(i_1i_2 i_3 i_4)}-
I_{(0001)T,t}^{*(i_1i_2 i_3 i_4)q}\right)^2\right\}=
\frac{(T-t)^6}{36}-\sum_{j_1,j_2,j_3, j_4=0}^{q}
\left(C_{j_4j_3j_2j_1}^{0001}\right)^2,
$$
$$
{\sf M}\left\{\left(
I_{(000000)T,t}^{*(i_1 i_2 i_3 i_4 i_5 i_6)}-
I_{(000000)T,t}^{*(i_1 i_2 i_3 i_4 i_5 i_6)q}\right)^2\right\}=
\frac{(T-t)^{6}}{720}-\sum_{j_1,j_2,j_3,j_4,j_5,j_6=0}^{q}
C_{j_6 j_5 j_4 j_3 j_2 j_1}^2.
$$

\vspace{3mm}

Recall that the systems of iterated  stochastic integrals 
(\ref{ito-ito})--(\ref{k1001xxxx})
are part of the 
Taylor--It\^{o} and Taylor--Stratonovich expansions 
(see Chapter 4).

The function $K(t_1,\ldots,t_k)$ from 
Theorem 1.1 for
the set (\ref{k1000xxxx}) is defined by

\vspace{-6mm}
\begin{equation}
\label{leto7000}
~~~K(t_1,\ldots,t_k)=
(t-t_k)^{l_k}\ldots (t-t_1)^{l_1}\ {\bf 1}_{\{t_1<\ldots<t_k\}},\ \ \
t_1,\ldots,t_k\in[t, T],
\end{equation}
where ${\bf 1}_A$ is the indicator of the set $A$.

In particular, for the stochastic integrals 
$I_{(1)T,t}^{(i_1)},$ $I_{(2)T,t}^{(i_1)},$ $I_{(00)T,t}^{(i_1i_2)},$
$I_{(000)T,t}^{(i_1i_2i_3)},$ $I_{(01)T,t}^{(i_1i_2)},$
$I_{(10)T,t}^{(i_1i_2)},$
$I_{(0000)T,t}^{(i_1\ldots i_4)},$
$I_{(20)T,t}^{(i_1i_2)},$ $I_{(11)T,t}^{(i_1i_2)},$
$I_{(02)T,t}^{(i_1i_2)}$ $(i_1,\ldots, i_4=1,\ldots,m)$~
the ~func-\\

\vspace{-5mm}
\noindent
tions $K(t_1,\ldots,t_k)$ defined by (\ref{leto7000})
look as follows

\vspace{-6.5mm}
\begin{equation}
\label{aaa1}
K_1(t_1)=t-t_1,\ \ \ K_2(t_1)=(t-t_1)^2,\ \ \
K_{00}(t_1,t_2)={\bf 1}_{\{t_1<t_2\}},
\end{equation}

\vspace{-7.5mm}
\begin{equation}
\label{aaa2}
K_{000}(t_1,t_2,t_3)={\bf 1}_{\{t_1<t_2<t_3\}},\ \ \
K_{01}(t_1,t_2)=(t-t_2){\bf 1}_{\{t_1<t_2\}},
\end{equation}

\vspace{-7.5mm}
\begin{equation}
\label{aaa3}
K_{10}(t_1,t_2)=(t-t_1){\bf 1}_{\{t_1<t_2\}},\ \ \
K_{0000}(t_1,t_2)={\bf 1}_{\{t_1<t_2<t_3<t_4\}},
\end{equation}

\vspace{-7.5mm}
\begin{equation}
\label{aaa4}
K_{20}(t_1,t_2)=(t-t_1)^2{\bf 1}_{\{t_1<t_2\}},\ \ \
K_{11}(t_1,t_2)=(t-t_1)(t-t_2){\bf 1}_{\{t_1<t_2\}},
\end{equation}

\vspace{-4.5mm}
\begin{equation}
\label{aaa5}
K_{02}(t_1,t_2)=(t-t_2)^2{\bf 1}_{\{t_1<t_2\}},
\end{equation}

\noindent
where $t_1,\ldots, t_4\in [t, T].$

It is obviously that the most simple 
expansion 
for the polynomial of a finite degree
into the Fourier series using the 
complete orthonormal system of functions in the space $L_2([t, T])$
will be its Fourier--Legendre 
expansion (finite sum). The polynomial functions 
are included in the functions (\ref{aaa1})--(\ref{aaa5}) 
as their components if~ $l_1^2+\ldots+l_k^2>0$.
So, it is logical to expect that the most simple expansions 
for the functions
(\ref{aaa1})--(\ref{aaa5}) 
into generalized multiple Fourier series 
will 
be Fourier--Legendre expansions of these functions
when $l_1^2+\ldots+l_k^2>0$.
Note that the given assumption is confirmed completely 
(compare the formulas (\ref{4002}), (\ref{4003}) with the formulas
(\ref{420}), (\ref{460}) (see below) correspondently).
So, usage of Legendre polynomials for
the 
approximation of iterated It\^{o} and Stratonovich stochastic integrals is 
a step forward.

\vspace{3mm}

\section{Mean-Square Approximation of Specific Iterated  
Stra\-to\-no\-vich Stochastic Integrals of multiplicities 1 to 3 
Based on Trigonometric System of Functions}

\vspace{3mm}

In \cite{1}-\cite{12aa}, \cite{arxiv-22} 
on the base of Theorems 1.1, 2.2, 2.6, and 2.8
the author obtained
(also see
early publications \cite{old-art-1} (1997), \cite{old-art-2} (1998),
\cite{very-old-1} (1994), \cite{very-old-2} (1996))
the
following expansions of the iterated Stratonovich stochastic
integrals (\ref{k1001xxxx}) 
(independently from the papers
\cite{Zapad-1}-\cite{Zapad-4}, \cite{Zapad-8} 
excepting the method, in which the additional
random variables $\xi_q^{(i)}$ and $\mu_q^{(i)}$ 
are introduced)

\newpage
\noindent 
$$
I_{(0)T,t}^{*(i_1)}=\sqrt{T-t}\zeta_0^{(i_1)},
$$

\begin{equation}
\label{420}
~~~~~~ I_{(1)T,t}^{*(i_1)q}=-\frac{{(T-t)}^{3/2}}{2}
\Biggl(\zeta_0^{(i_1)}-\frac{\sqrt{2}}{\pi}\Biggl(\sum_{r=1}^{q}
\frac{1}{r}
\zeta_{2r-1}^{(i_1)}+\sqrt{\alpha_q}\xi_q^{(i_1)}\Biggr)
\Biggr),
\end{equation}

\vspace{1mm}
$$
I_{(00)T,t}^{*(i_1 i_2)q}=\frac{1}{2}(T-t)\Biggl(
\zeta_{0}^{(i_1)}\zeta_{0}^{(i_2)}
+\frac{1}{\pi}
\sum_{r=1}^{q}\frac{1}{r}\left(
\zeta_{2r}^{(i_1)}\zeta_{2r-1}^{(i_2)}-
\zeta_{2r-1}^{(i_1)}\zeta_{2r}^{(i_2)}+
\right.\Biggr.
$$
\begin{equation}
\label{430}
~~~~~ +\Biggl.\left.\sqrt{2}\left(\zeta_{2r-1}^{(i_1)}\zeta_{0}^{(i_2)}-
\zeta_{0}^{(i_1)}\zeta_{2r-1}^{(i_2)}\right)\right)
+\frac{\sqrt{2}}{\pi}\sqrt{\alpha_q}\left(
\xi_q^{(i_1)}\zeta_0^{(i_2)}-\zeta_0^{(i_1)}\xi_q^{(i_2)}\right)\Biggr),
\end{equation}

\vspace{5mm}
$$
I_{(000)T,t}^{*(i_1 i_2 i_3)q}=(T-t)^{3/2}\Biggl(\frac{1}{6}
\zeta_{0}^{(i_1)}\zeta_{0}^{(i_2)}\zeta_{0}^{(i_3)}+\Biggr.
\frac{\sqrt{\alpha_q}}{2\sqrt{2}\pi}\left(
\xi_q^{(i_1)}\zeta_0^{(i_2)}\zeta_0^{(i_3)}-\xi_q^{(i_3)}\zeta_0^{(i_2)}
\zeta_0^{(i_1)}\right)+
$$
$$
+\frac{1}{2\sqrt{2}\pi^2}\sqrt{\beta_q}\left(
\mu_q^{(i_1)}\zeta_0^{(i_2)}\zeta_0^{(i_3)}-2\mu_q^{(i_2)}\zeta_0^{(i_1)}
\zeta_0^{(i_3)}+\mu_q^{(i_3)}\zeta_0^{(i_1)}\zeta_0^{(i_2)}\right)+
$$
$$
+
\frac{1}{2\sqrt{2}}\sum_{r=1}^{q}
\Biggl(\frac{1}{\pi r}\left(
\zeta_{2r-1}^{(i_1)}
\zeta_{0}^{(i_2)}\zeta_{0}^{(i_3)}-
\zeta_{2r-1}^{(i_3)}
\zeta_{0}^{(i_2)}\zeta_{0}^{(i_1)}\right)+\Biggr.
$$
$$
\Biggl.+
\frac{1}{\pi^2 r^2}\left(
\zeta_{2r}^{(i_1)}
\zeta_{0}^{(i_2)}\zeta_{0}^{(i_3)}-
2\zeta_{2r}^{(i_2)}
\zeta_{0}^{(i_3)}\zeta_{0}^{(i_1)}+
\zeta_{2r}^{(i_3)}
\zeta_{0}^{(i_2)}\zeta_{0}^{(i_1)}\right)\Biggr)+
$$
$$
+
\sum_{r=1}^{q}
\Biggl(\frac{1}{4\pi r}\left(
\zeta_{2r}^{(i_1)}
\zeta_{2r-1}^{(i_2)}\zeta_{0}^{(i_3)}-
\zeta_{2r-1}^{(i_1)}
\zeta_{2r}^{(i_2)}\zeta_{0}^{(i_3)}-
\zeta_{2r-1}^{(i_2)}
\zeta_{2r}^{(i_3)}\zeta_{0}^{(i_1)}+
\zeta_{2r-1}^{(i_3)}
\zeta_{2r}^{(i_2)}\zeta_{0}^{(i_1)}\right)+\Biggr.
$$
$$
+
\frac{1}{8\pi^2 r^2}\left(
3\zeta_{2r-1}^{(i_1)}
\zeta_{2r-1}^{(i_2)}\zeta_{0}^{(i_3)}+
\zeta_{2r}^{(i_1)}
\zeta_{2r}^{(i_2)}\zeta_{0}^{(i_3)}-
6\zeta_{2r-1}^{(i_1)}
\zeta_{2r-1}^{(i_3)}\zeta_{0}^{(i_2)}+\right.
$$
\begin{equation}
\label{44xxxx}
\Biggl.\left.
~~~~~~~ +
3\zeta_{2r-1}^{(i_2)}
\zeta_{2r-1}^{(i_3)}\zeta_{0}^{(i_1)}-
2\zeta_{2r}^{(i_1)}
\zeta_{2r}^{(i_3)}\zeta_{0}^{(i_2)}+
\zeta_{2r}^{(i_3)}
\zeta_{2r}^{(i_2)}\zeta_{0}^{(i_1)}\right)\Biggr)
\Biggl.+D_{T,t}^{(i_1i_2i_3)q}\Biggr),
\end{equation}

\vspace{4mm}
\noindent
where

\vspace{-3mm}
$$
D_{T,t}^{(i_1i_2i_3)q}=
\frac{1}{2\pi^2}\sum_{\stackrel{r,l=1}{{}_{r\ne l}}}^{q}
\Biggl(\frac{1}{r^2-l^2}\biggl(
\zeta_{2r}^{(i_1)}
\zeta_{2l}^{(i_2)}\zeta_{0}^{(i_3)}-
\zeta_{2r}^{(i_2)}
\zeta_{0}^{(i_1)}\zeta_{2l}^{(i_3)}+\biggr.\Biggr.
$$
$$
\Biggl.+\biggl.
\frac{r}{l}
\zeta_{2r-1}^{(i_1)}
\zeta_{2l-1}^{(i_2)}\zeta_{0}^{(i_3)}-\frac{l}{r}
\zeta_{0}^{(i_1)}
\zeta_{2r-1}^{(i_2)}\zeta_{2l-1}^{(i_3)}\biggr)-
\frac{1}{rl}\zeta_{2r-1}^{(i_1)}
\zeta_{0}^{(i_2)}\zeta_{2l-1}^{(i_3)}\Biggr)+
$$
$$
+
\frac{1}{4\sqrt{2}\pi^2}\Biggl(
\sum_{r,m=1}^{q}\Biggl(\frac{2}{rm}
\left(-\zeta_{2r-1}^{(i_1)}
\zeta_{2m-1}^{(i_2)}\zeta_{2m}^{(i_3)}+
\zeta_{2r-1}^{(i_1)}
\zeta_{2r}^{(i_2)}\zeta_{2m-1}^{(i_3)}+
\right.\Biggr.\Biggr.
$$
$$
\left.+
\zeta_{2r-1}^{(i_1)}
\zeta_{2m}^{(i_2)}\zeta_{2m-1}^{(i_3)}-
\zeta_{2r}^{(i_1)}
\zeta_{2r-1}^{(i_2)}\zeta_{2m-1}^{(i_3)}\right)+
$$
$$
+\frac{1}{m(r+m)}
\left(-\zeta_{2(m+r)}^{(i_1)}
\zeta_{2r}^{(i_2)}\zeta_{2m}^{(i_3)}-
\zeta_{2(m+r)-1}^{(i_1)}
\zeta_{2r-1}^{(i_2)}\zeta_{2m}^{(i_3)}-
\right.
$$
$$
\Biggl.\left.
-\zeta_{2(m+r)-1}^{(i_1)}
\zeta_{2r}^{(i_2)}\zeta_{2m-1}^{(i_3)}+
\zeta_{2(m+r)}^{(i_1)}
\zeta_{2r-1}^{(i_2)}\zeta_{2m-1}^{(i_3)}\right)\Biggr)+
$$
$$
+
\sum_{m=1}^{q}\sum_{l=m+1}^{q}\Biggl(\frac{1}{m(l-m)}
\left(\zeta_{2(l-m)}^{(i_1)}
\zeta_{2l}^{(i_2)}\zeta_{2m}^{(i_3)}+
\zeta_{2(l-m)-1}^{(i_1)}
\zeta_{2l-1}^{(i_2)}\zeta_{2m}^{(i_3)}-
\right.\Biggr.
$$
$$
\left.
-\zeta_{2(l-m)-1}^{(i_1)}
\zeta_{2l}^{(i_2)}\zeta_{2m-1}^{(i_3)}+
\zeta_{2(l-m)}^{(i_1)}
\zeta_{2l-1}^{(i_2)}\zeta_{2m-1}^{(i_3)}\right)+
$$
$$
+
\frac{1}{l(l-m)}
\left(-\zeta_{2(l-m)}^{(i_1)}
\zeta_{2m}^{(i_2)}\zeta_{2l}^{(i_3)}+
\zeta_{2(l-m)-1}^{(i_1)}
\zeta_{2m-1}^{(i_2)}\zeta_{2l}^{(i_3)}-
\right.
$$
$$
\Biggl.
\Biggl.
\Biggl.
\left.
-\zeta_{2(l-m)-1}^{(i_1)}
\zeta_{2m}^{(i_2)}\zeta_{2l-1}^{(i_3)}-
\zeta_{2(l-m)}^{(i_1)}
\zeta_{2m-1}^{(i_2)}\zeta_{2l-1}^{(i_3)}\right)\Biggr)\Biggr),
$$

\vspace{7mm}
$$
I_{(10)T,t}^{*(i_1 i_2)q}=-(T-t)^{2}\Biggl(\frac{1}{6}
\zeta_{0}^{(i_1)}\zeta_{0}^{(i_2)}-\frac{1}{2\sqrt{2}\pi}
\sqrt{\alpha_q}\xi_q^{(i_2)}\zeta_0^{(i_1)}+\Biggr.
$$
$$
+\frac{1}{2\sqrt{2}\pi^2}\sqrt{\beta_q}\Biggl(
\mu_q^{(i_2)}\zeta_0^{(i_1)}-2\mu_q^{(i_1)}\zeta_0^{(i_2)}\Biggr)+
$$
$$
+\frac{1}{2\sqrt{2}}\sum_{r=1}^{q}
\Biggl(-\frac{1}{\pi r}
\zeta_{2r-1}^{(i_2)}
\zeta_{0}^{(i_1)}+
\frac{1}{\pi^2 r^2}\left(
\zeta_{2r}^{(i_2)}
\zeta_{0}^{(i_1)}-
2\zeta_{2r}^{(i_1)}
\zeta_{0}^{(i_2)}\right)\Biggr)-
$$
$$
-
\frac{1}{2\pi^2}\sum_{\stackrel{r,l=1}{{}_{r\ne l}}}^{q}
\frac{1}{r^2-l^2}\Biggl(
\zeta_{2r}^{(i_1)}
\zeta_{2l}^{(i_2)}+
\frac{l}{r}
\zeta_{2r-1}^{(i_1)}
\zeta_{2l-1}^{(i_2)}
\Biggr)+
$$
\begin{equation}
\label{9440}
+
\sum_{r=1}^{q}
\Biggl(\frac{1}{4\pi r}\left(
\zeta_{2r}^{(i_1)}
\zeta_{2r-1}^{(i_2)}-
\zeta_{2r-1}^{(i_1)}
\zeta_{2r}^{(i_2)}\right)+
\Biggl.\Biggl.
\frac{1}{8\pi^2 r^2}\left(
3\zeta_{2r-1}^{(i_1)}
\zeta_{2r-1}^{(i_2)}+
\zeta_{2r}^{(i_2)}
\zeta_{2r}^{(i_1)}\right)\Biggr)\Biggr),
\end{equation}

\newpage
\noindent
$$
I_{(01)T,t}^{*(i_1 i_2)q}=
(T-t)^{2}\Biggl(-\frac{1}{3}
\zeta_{0}^{(i_1)}\zeta_{0}^{(i_2)}-\frac{1}{2\sqrt{2}\pi}
\sqrt{\alpha_q}\left(\xi_q^{(i_1)}\zeta_0^{(i_2)}-
2\xi_q^{(i_2)}\zeta_0^{(i_1)}\right)
+\Biggr.
$$
$$
+\frac{1}{2\sqrt{2}\pi^2}\sqrt{\beta_q}\Biggl(
\mu_q^{(i_1)}\zeta_0^{(i_2)}-2\mu_q^{(i_2)}\zeta_0^{(i_1)}\Biggr)-
$$
$$
-\frac{1}{2\sqrt{2}}\sum_{r=1}^{q}
\Biggl(\frac{1}{\pi r}\left(
\zeta_{2r-1}^{(i_1)}
\zeta_{0}^{(i_2)}-
2\zeta_{2r-1}^{(i_2)}
\zeta_{0}^{(i_1)}\right)
-\frac{1}{\pi^2 r^2}\left(
\zeta_{2r}^{(i_1)}
\zeta_{0}^{(i_2)}-
2\zeta_{2r}^{(i_2)}
\zeta_{0}^{(i_1)}\right)\Biggr)+
$$
$$
+
\frac{1}{2\pi^2}\sum_{\stackrel{r,l=1}{{}_{r\ne l}}}^{q}
\frac{1}{r^2-l^2}\Biggl(
\frac{r}{l}\zeta_{2r-1}^{(i_1)}
\zeta_{2l-1}^{(i_2)}+
\zeta_{2r}^{(i_1)}
\zeta_{2l}^{(i_2)}
\Biggr)-
$$
\begin{equation}
\label{450}
-
\sum_{r=1}^{q}
\Biggl(\frac{1}{4\pi r}\left(
\zeta_{2r}^{(i_1)}
\zeta_{2r-1}^{(i_2)}-
\zeta_{2r-1}^{(i_1)}
\zeta_{2r}^{(i_2)}\right)-
\Biggl.\Biggl.
\frac{1}{8\pi^2 r^2}\left(
3\zeta_{2r-1}^{(i_1)}
\zeta_{2r-1}^{(i_2)}+
\zeta_{2r}^{(i_1)}
\zeta_{2r}^{(i_2)}\right)\Biggr)\Biggr),
\end{equation}

\vspace{3mm}

$$
I_{(2)T,t}^{*(i_1)q}=
(T-t)^{5/2}\Biggl(
\frac{1}{3}\zeta_0^{(i_1)}+\frac{1}{\sqrt{2}\pi^2}\Biggl(
\sum_{r=1}^{q}\frac{1}{r^2}\zeta_{2r}^{(i_1)}+
\sqrt{\beta_q}\mu_q^{(i_1)}\Biggr)-\Biggr.
$$
\begin{equation}
\label{460}
\Biggl.-
\frac{1}{\sqrt{2}\pi}\Biggl(\sum_{r=1}^q
\frac{1}{r}\zeta_{2r-1}^{(i_1)}+\sqrt{\alpha_q}\xi_q^{(i_1)}\Biggr)\Biggr),
\end{equation}

\vspace{3mm}
\noindent
where
$$
\xi_q^{(i)}=\frac{1}{\sqrt{\alpha_q}}\sum_{r=q+1}^{\infty}
\frac{1}{r}~\zeta_{2r-1}^{(i)},\ \ \
\alpha_q=\frac{\pi^2}{6}-\sum_{r=1}^q\frac{1}{r^2},\ \ \
\mu_q^{(i)}=\frac{1}{\sqrt{\beta_q}}\sum_{r=q+1}^{\infty}
\frac{1}{r^2}~\zeta_{2r}^{(i)},
$$
$$
\beta_q=\frac{\pi^4}{90}-\sum_{r=1}^q\frac{1}{r^4},\ \ \
\zeta_j^{(i)}=\int\limits_t^T\phi_j(s)d{\bf w}_s^{(i)},
$$
where
$\phi_j(s)$ is defined by (\ref{trig11}) and
$\zeta_0^{(i)},$ $\zeta_{2r}^{(i)},$
$\zeta_{2r-1}^{(i)},$ $\xi_q^{(i)},$ $\mu_q^{(i)}$ ($r=1,\ldots,q,$\
$i=1,\ldots,m$) are independent
standard Gaussian random variables
($i_1, i_2, i_3=1,\ldots,m$).

Note that (\ref{9440}), (\ref{450}) imply the following
\begin{equation}
\label{leto3000mil}
\sum\limits_{j=0}^{\infty}C_{jj}^{10}=
\sum\limits_{j=0}^{\infty}C_{jj}^{01}=-\frac{(T-t)^2}{4},
\end{equation}
where
$$
C_{jj}^{10}=\int\limits_{t}^{T}\phi_j(x)
\int\limits_{t}^{x}\phi_j(y)(t-y)
dy dx,
$$
$$
C_{jj}^{01}=\int\limits_{t}^{T}\phi_j(x)(t-x)
\int\limits_{t}^{x}\phi_j(y)
dy dx.
$$

Note that the formulas (\ref{leto3000mil}) are particular cases 
of the more general relation (\ref{5t}), which has been applied for the
proof 
of Theorems 2.1--2.3.

Let us consider the mean-square errors of approximations
(\ref{430})--(\ref{450}). From the relations (\ref{430})--(\ref{450})
when $i_1\ne i_2,$ $i_2\ne i_3,$ $i_1\ne i_3$ by direct 
calculation we obtain 
\begin{equation}
\label{801}
{\sf M}\left\{\left(I_{(00)T,t}^{*(i_1 i_2)}-
I_{(00)T,t}^{*(i_1 i_2)q}
\right)^2\right\}
=\frac{(T-t)^{2}}{2\pi^2}\Biggl(\frac{\pi^2}{6}-
\sum_{r=1}^q \frac{1}{r^2}\Biggr),
\end{equation}

$$
{\sf M}\left\{\left(I_{(000)T,t}^{*(i_1 i_2 i_3)}-
I_{(000)T,t}^{*(i_1 i_2 i_3)q}\right)^2\right\}
=(T-t)^{3}\Biggl(\frac{1}{4\pi^2}
\Biggl(\frac{\pi^2}{6}-
\sum_{r=1}^q \frac{1}{r^2}\Biggr)+
\Biggr.
$$
\begin{equation}
\label{802}
~~~~~~\Biggl.
+\frac{55}{32\pi^4}\Biggl(\frac{\pi^4}{90}-
\sum_{r=1}^q \frac{1}{r^4}\Biggr)
+\frac{1}{4\pi^4}
\Biggl(\sum_{\stackrel{r,l=1}{{}_{r\ne l}}}^{\infty}
-\sum_{\stackrel{r,l=1}{{}_{r\ne l}}}^{q}
\Biggr)
\frac{5l^4+4r^4-3l^2r^2}{r^2 l^2(r^2-l^2)^2}\Biggr),
\end{equation}

\vspace{2mm}

$$
{\sf M}\left\{\left(I_{(01)T,t}
^{*(i_1i_2)}-I_{(01)T,t}^{*(i_1i_2)q}\right)^2\right\}
=(T-t)^{4}\Biggl(\frac{1}{8\pi^2}
\Biggl(\frac{\pi^2}{6}-
\sum_{r=1}^q \frac{1}{r^2}\Biggr)+\Biggr.
$$
\begin{equation}
\label{804}
~~~~~~\Biggl.+\frac{5}{32\pi^4}\Biggl(\frac{\pi^4}{90}-
\sum_{r=1}^q \frac{1}{r^4}\Biggr)+
\frac{1}{4\pi^4}\Biggl(\sum_{\stackrel{k,l=1}{{}_{k\ne l}}}^{\infty}
-\sum_{\stackrel{k,l=1}{{}_{k\ne l}}}^{q}
\Biggr)\frac{l^2+k^2}{k^2(l^2-k^2)^2}\Biggr),
\end{equation}

\vspace{2mm}

$$
{\sf M}\left\{\left(I_{(10)T,t}^{*(i_1i_2)}
-I_{(10)T,t}^{*(i_1i_2)q}\right)^2\right\}
=(T-t)^{4}\Biggl(\frac{1}{8\pi^2}
\Biggl(\frac{\pi^2}{6}-
\sum_{r=1}^q \frac{1}{r^2}\Biggr)+\Biggr.
$$
\begin{equation}
\label{805}
~~~~~~+\Biggl.
\frac{5}{32\pi^4}\Biggl(\frac{\pi^4}{90}-
\sum_{r=1}^q \frac{1}{r^4}\Biggr)+
\frac{1}{4\pi^4}\Biggl(\sum_{\stackrel{k,l=1}{{}_{k\ne l}}}^{\infty}
-\sum_{\stackrel{k,l=1}{{}_{k\ne l}}}^{q}
\Biggr)\frac{l^2+k^2}{l^2(l^2-k^2)^2}\Biggr).
\end{equation}

\vspace{3mm}

It is easy to demonstrate that the relations 
(\ref{802}), (\ref{804}), and (\ref{805})
can be represented using Theorem 1.3 
in the following form
$$
{\sf M}\left\{\left(I_{(000)T,t}^{*(i_1 i_2 i_3)}-
I_{(000)T,t}^{*(i_1 i_2 i_3)q}\right)^2\right\}=
(T-t)^3\Biggl(\frac{4}{45}-\frac{1}{4\pi^2}\sum_{r=1}^q\frac{1}{r^2}-
\Biggl.
$$
\begin{equation}
\label{101.100}
\Biggl.-\frac{55}{32\pi^4}\sum_{r=1}^q\frac{1}{r^4}-
\frac{1}{4\pi^4}\sum_{\stackrel{r,l=1}{{}_{r\ne l}}}^q
\frac{5l^4+4r^4-3r^2l^2}{r^2 l^2 \left(r^2-l^2\right)^2}\Biggr),
\end{equation}

\vspace{2mm}

$$
{\sf M}\left\{\left(I_{(10)T,t}^{*(i_1 i_2)}-
I_{(10)T,t}^{*(i_1 i_2)q}\right)^2\right\}=
\frac{(T-t)^4}{4}\Biggl(\frac{1}{9}-
\frac{1}{2\pi^2}\sum_{r=1}^q \frac{1}{r^2}-\Biggr.
$$
\begin{equation}
\label{101.101}
\Biggl.-\frac{5}{8\pi^4}\sum_{r=1}^q \frac{1}{r^4}-
\frac{1}{\pi^4}\sum_{\stackrel{k,l=1}{{}_{k\ne l}}}^q
\frac{k^2+l^2}{l^2\left(l^2-k^2\right)^2}\Biggr),
\end{equation}

\vspace{2mm}

$$
{\sf M}\left\{\left(I_{(01)T,t}^{*(i_1 i_2)}-
I_{(01)T,t}^{*(i_1 i_2)q}\right)^2\right\}=
\frac{(T-t)^4}{4}\Biggl(\frac{1}{9}-
\frac{1}{2\pi^2}\sum_{r=1}^q \frac{1}{r^2}-\Biggr.
$$
\begin{equation}
\label{101.102}
\Biggl.-\frac{5}{8\pi^4}\sum_{r=1}^q \frac{1}{r^4}-
\frac{1}{\pi^4}\sum_{\stackrel{k,l=1}{{}_{k\ne l}}}^q
\frac{l^2+k^2}{k^2\left(l^2-k^2\right)^2}\Biggr).
\end{equation}

\vspace{3mm}

Comparing (\ref{101.100})--(\ref{101.102}) and
(\ref{802})--(\ref{805}), we note that
\begin{equation}
\label{101.103}
\sum_{\stackrel{k,l=1}{{}_{k\ne l}}}^{\infty}\frac{l^2+k^2}
{k^2\left(l^2-k^2\right)^2}=
\sum_{\stackrel{k,l=1}{{}_{k\ne l}}}^{\infty}\frac{l^2+k^2}
{l^2\left(l^2-k^2\right)^2}=\frac{\pi^4}{48},
\end{equation}
\begin{equation}
\label{daug1}
\sum_{\stackrel{r,l=1}{{}_{r\ne l}}}^{\infty}
\frac{5l^4+4r^4-3r^2 l^2}{r^2 l^2\left(r^2-l^2\right)^2}=
\frac{9\pi^4}{80}.
\end{equation}

Let us consider approximations of stochastic 
integrals $I_{(10)T,t}^{*(i_1i_1)},$
$I_{(01)T,t}^{*(i_1i_1)}$ and conditions for selecting 
number $q$ using the trigonometric system of functions
$$
I_{(10)T,t}^{*(i_1 i_1)q}=-(T-t)^{2}\Biggl(\frac{1}{6}
\left(\zeta_{0}^{(i_1)}\right)^2-\frac{1}{2\sqrt{2}\pi}
\sqrt{\alpha_q}\xi_q^{(i_1)}\zeta_0^{(i_1)}-\Biggr.
$$
$$
-\frac{1}{2\sqrt{2}\pi^2}\sqrt{\beta_q}
\mu_q^{(i_1)}\zeta_0^{(i_1)}
-\frac{1}{2\sqrt{2}}\sum_{r=1}^{q}
\Biggl(\frac{1}{\pi r}
\zeta_{2r-1}^{(i_1)}
\zeta_{0}^{(i_1)}+
\frac{1}{\pi^2 r^2}
\zeta_{2r}^{(i_1)}
\zeta_{0}^{(i_1)}\Biggr)-
$$
$$
-
\frac{1}{2\pi^2}\sum_{\stackrel{r,l=1}{{}_{r\ne l}}}^{q}
\frac{1}{r^2-l^2}\Biggl(
\zeta_{2r}^{(i_1)}
\zeta_{2l}^{(i_1)}+
\frac{l}{r}
\zeta_{2r-1}^{(i_1)}
\zeta_{2l-1}^{(i_1)}
\Biggr)+
$$
$$
\Biggl.+
\frac{1}{8\pi^2}\sum_{r=1}^{q}
\frac{1}{r^2}\left(
3\left(\zeta_{2r-1}^{(i_1)}\right)^2
+
\left(\zeta_{2r}^{(i_1)}\right)^2\right)\Biggr),
$$

\vspace{5mm}

$$
I_{(01)T,t}^{*(i_1 i_1)q}=(T-t)^{2}\Biggl(-\frac{1}{3}
\left(\zeta_{0}^{(i_1)}\right)^2+\frac{1}{2\sqrt{2}\pi}
\sqrt{\alpha_q}\xi_q^{(i_1)}\zeta_0^{(i_1)}-\Biggr.
$$
$$
-\frac{1}{2\sqrt{2}\pi^2}\sqrt{\beta_q}
\mu_q^{(i_1)}\zeta_0^{(i_1)}
+\frac{1}{2\sqrt{2}}\sum_{r=1}^{q}
\Biggl(\frac{1}{\pi r}
\zeta_{2r-1}^{(i_1)}
\zeta_{0}^{(i_1)}-
\frac{1}{\pi^2 r^2}
\zeta_{2r}^{(i_1)}
\zeta_{0}^{(i_1)}\Biggr)+
$$
$$
+
\frac{1}{2\pi^2}\sum_{\stackrel{r,l=1}{{}_{r\ne l}}}^{q}
\frac{1}{r^2-l^2}\Biggl(
\zeta_{2r}^{(i_1)}
\zeta_{2l}^{(i_1)}+
\frac{r}{l}
\zeta_{2r-1}^{(i_1)}
\zeta_{2l-1}^{(i_1)}
\Biggr)+
$$
$$
\Biggl.+
\frac{1}{8\pi^2}\sum_{r=1}^{q}
\frac{1}{r^2}\left(
3\left(\zeta_{2r-1}^{(i_1)}\right)^2
+
\left(\zeta_{2r}^{(i_1)}\right)^2\right)\Biggr).
$$

\vspace{3mm}

Furthermore, we have
$$
{\sf M}\left\{\left(I_{(01)T,t}
^{*(i_1i_1)}-I_{(01)T,t}^{*(i_1i_1)q}\right)^2\right\}=
{\sf M}\left\{\left(I_{(10)T,t}^{*(i_1i_1)}
-I_{(10)T,t}^{*(i_1i_1)q}\right)^2\right\}=
$$
$$
=\frac{(T-t)^{4}}{4}\Biggl(\frac{2}{\pi^4}\Biggl(\frac{\pi^4}{90}-
\sum_{r=1}^q \frac{1}{r^4}\Biggr)
+\frac{1}{\pi^4}
\Biggl(\frac{\pi^2}{6}-
\sum_{r=1}^q \frac{1}{r^2}\Biggr)^2+\Biggr.
$$
\begin{equation}
\label{101.104}
\Biggl.+
\frac{1}{\pi^4}\Biggl(\sum_{\stackrel{k,l=1}{{}_{k\ne l}}}^{\infty}
-\sum_{\stackrel{k,l=1}{{}_{k\ne l}}}^{q}
\Biggr)\frac{l^2+k^2}{k^2(l^2-k^2)^2}\Biggr).
\end{equation}

\vspace{2mm}

Considering (\ref{101.103}), we can rewrite the relation (\ref{101.104}) in 
the following form
$$
{\sf M}\left\{\left(I_{(01)T,t}
^{*(i_1i_1)}-I_{(01)T,t}^{*(i_1i_1)q}\right)^2\right\}=
{\sf M}\left\{\left(I_{(10)T,t}^{*(i_1i_1)}
-I_{(10)T,t}^{*(i_1i_1)q}\right)^2\right\}=
$$
$$
=\frac{(T-t)^{4}}{4}\Biggl(\frac{17}{240}-
\frac{1}{3\pi^2}
\sum_{r=1}^q \frac{1}{r^2}-\frac{2}{\pi^4}
\sum_{r=1}^q \frac{1}{r^4} +\Biggr.
$$
\begin{equation}
\label{daug}
\Biggl.+ \frac{1}{\pi^4}
\Biggl(
\sum_{r=1}^q \frac{1}{r^2}\Biggr)^2-
\frac{1}{\pi^4}
\sum_{\stackrel{k,l=1}{{}_{k\ne l}}}^{q}
\frac{l^2+k^2}{k^2(l^2-k^2)^2}\Biggr).
\end{equation}

\vspace{2mm}

In Tables 5.37--5.39 we confirm numerically the formulas 
(\ref{101.100})--(\ref{101.102}), 
(\ref{daug})
for various values of $q$. In Tables 5.37--5.39
the number $\varepsilon$  
means right-hand sides of the mentioned formulas.
Obviously, these results are consistent with 
the estimate (\ref{zsel1}).

The formulas (\ref{101.103}), (\ref{daug1}) appear to be interesting. 
Let us 
confirm numerically their correctness in Tables 5.40 and 5.41 
(the number $\varepsilon_q$ is an absolute 
deviation of multiple partial sums with 
the upper limit of summation $q$ for the 
series (\ref{101.103}), (\ref{daug1})
from the right-hand sides of the formulas (\ref{101.103}), (\ref{daug1});
convergence of multiple series is regarded here 
when $p_1=p_2=q\to\infty$, which is acceptable according to Theorems 
1.1, 2.2, 2.6, and 2.8).

\begin{table}
\centering
\caption{Confirmation of the formula (\ref{101.100})}
\label{tab:5.37}      
\begin{tabular}{p{2.1cm}p{1.7cm}p{1.7cm}p{2.1cm}p{2.3cm}p{2.3cm}p{2.3cm}}
\hline\noalign{\smallskip}
$\varepsilon/(T-t)^3$&0.0459&0.0072&$7.5722\cdot 10^{-4}$
&$7.5973\cdot 10^{-5}$&
$7.5990\cdot 10^{-6}$\\
\noalign{\smallskip}\hline\noalign{\smallskip}
$q$&1&10&100&1000&10000\\
\noalign{\smallskip}\hline\noalign{\smallskip}
\end{tabular}
\end{table}

\begin{table}
\centering
\caption{Confirmation of the formulas (\ref{101.101}), (\ref{101.102})}
\label{tab:5.38}      
\begin{tabular}{p{2.1cm}p{1.7cm}p{1.7cm}p{2.1cm}p{2.3cm}p{2.3cm}p{2.3cm}}
\hline\noalign{\smallskip}
$4\varepsilon/(T-t)^4$&0.0540&0.0082&$8.4261\cdot 10^{-4}$
&$8.4429\cdot 10^{-5}$&
$8.4435\cdot 10^{-6}$\\
\noalign{\smallskip}\hline\noalign{\smallskip}
$q$&1&10&100&1000&10000\\
\noalign{\smallskip}\hline\noalign{\smallskip}
\end{tabular}
\end{table}

\begin{table}
\centering
\caption{Confirmation of the formula (\ref{daug})}
\label{tab:5.39}      
\begin{tabular}{p{2.1cm}p{1.7cm}p{1.7cm}p{2.1cm}p{2.3cm}p{2.3cm}p{2.3cm}}
\hline\noalign{\smallskip}
$4\varepsilon/(T-t)^4$&0.0268&0.0034&$3.3955\cdot 10^{-4}$
&$3.3804\cdot 10^{-5}$&
$3.3778\cdot 10^{-6}$\\
\noalign{\smallskip}\hline\noalign{\smallskip}
$q$&1&10&100&1000&10000\\
\noalign{\smallskip}\hline\noalign{\smallskip}
\end{tabular}
\end{table}

\begin{table}
\centering
\caption{Confirmation of the formula (\ref{101.103})}
\label{tab:5.40}      
\begin{tabular}{p{1.3cm}p{1.8cm}p{1.8cm}p{1.8cm}p{1.8cm}p{2.3cm}p{1.8cm}}
\hline\noalign{\smallskip}
$\varepsilon_q$&2.0294&0.3241&0.0330
&0.0033&
$3.2902\cdot 10^{-4}$\\
\noalign{\smallskip}\hline\noalign{\smallskip}
$q$&1&10&100&1000&10000\\
\noalign{\smallskip}\hline\noalign{\smallskip}
\end{tabular}
\end{table}

\begin{table}
\centering
\caption{Confirmation of the formula (\ref{daug1})}
\label{tab:5.41}      
\begin{tabular}{p{1.1cm}p{1.5cm}p{1.5cm}p{1.5cm}p{1.5cm}p{1.5cm}p{1.5cm}}
\hline\noalign{\smallskip}
$\varepsilon_q$&10.9585&1.8836&0.1968
&0.0197&
0.0020\\
\noalign{\smallskip}\hline\noalign{\smallskip}
$q$&1&10&100&1000&10000\\
\noalign{\smallskip}\hline\noalign{\smallskip}
\end{tabular}
\end{table}

Using the trigonometric system of functions, let us consider
approximations of iterated stochastic integrals of the following form
$$
{J}_{(\lambda_{1}\ldots \lambda_k)T,t}^{*(i_1\ldots
i_k)}=
{\int\limits_t^{*}}^T
\ldots
{\int\limits_t^{*}}^{t_{2}}
d{\bf w}_{t_{1}}^{(i_1)}\ldots
d{\bf w}_{t_k}^{(i_k)},
$$
where $\lambda_l=1$ if $i_l=1,\ldots,m$ and
$\lambda_l=0$ if $i_l=0,$ $l=1,\ldots,k,$
${\bf w}_{\tau}^{(i)}$
$(i=1,\ldots,m)$ are independent standard Wiener processes, 
${\bf w}_{\tau}^{(0)}=\tau$.

It is easy to see that the approximations
$$
J_{(\lambda_1\lambda_2)T,t}^{*(i_1 i_2)q},\ \ \
J_{(\lambda_1\lambda_2\lambda_3)T,t}^{*(i_1 i_2 i_3)q}
$$ 
of the stochastic integrals  
$$
J_{(\lambda_1\lambda_2)T,t}^{*(i_1 i_2)},\ \ \
J_{(\lambda_1\lambda_2\lambda_3)T,t}^{*(i_1 i_2 i_3)}
$$
are defined by the right-hand sides of 
the formulas (\ref{430}), (\ref{44xxxx}),  where it is necessary to take
\begin{equation}
\label{123}
\zeta_j^{(i)}=\int\limits_t^T\phi_j(s)d{\bf w}_s^{(i)}
\end{equation}
and $i_1, i_2, i_3=0, 1,\ldots,m$.

Since
$$
\int\limits_t^T\phi_j(s)d{\bf w}_s^{(0)}=
\left\{
\begin{matrix}
\sqrt{T-t} &\hbox{if}\ j=0\cr\cr
0 &\hbox{if}\ j\ne 0
\end{matrix},\right.
$$ 
then
it is easy to get
from (\ref{430}) and (\ref{44xxxx}), considering
that in 
these equalities $\zeta_j^{(i)}$ is defined by (\ref{123})
and
$i_1, i_2, i_3=0, 1,\ldots,m$,
the following family of formulas
\begin{equation}
\label{df5}
~~~~~~~~~ J_{(10)T,t}^{*(i_1 0)q}=
\frac{1}{2}(T-t)^{3/2}\Biggl(
\zeta_0^{(i_1)}+\frac{\sqrt{2}}{\pi}
\Biggl(\sum_{r=1}^{q}\frac{1}{r}\zeta_{2r-1}^{(i_1)}+
\sqrt{\alpha_q}\xi_q^{(i_1)}\Biggr)\Biggr),
\end{equation}

\begin{equation}
\label{df6}
~~~~~~~~~ J_{(01)T,t}^{*(0 i_2)q}=
\frac{1}{2}(T-t)^{3/2}\Biggl(
\zeta_0^{(i_2)}-\frac{\sqrt{2}}{\pi}
\Biggl(\sum_{r=1}^{q}\frac{1}{r}\zeta_{2r-1}^{(i_2)}+
\sqrt{\alpha_q}\xi_q^{(i_2)}\Biggr)\Biggr),
\end{equation}

\newpage
\noindent
$$
J_{(001)T,t}^{*(00 i_3)q}=
(T-t)^{5/2}\Biggl(
\frac{1}{6}\zeta_0^{(i_3)}+\frac{1}{2\sqrt{2}\pi^2}\Biggl(
\sum_{r=1}^{q}\frac{1}{r^2}\zeta_{2r}^{(i_3)}+
\sqrt{\beta_q}\mu_q^{(i_3)}\Biggr)-\Biggr.
$$
$$
\Biggl.-
\frac{1}{2\sqrt{2}\pi}\Biggl(\sum_{r=1}^q
\frac{1}{r}\zeta_{2r-1}^{(i_3)}+\sqrt{\alpha_q}\xi_q^{(i_3)}\Biggr)\Biggr),
$$

\vspace{3mm}

$$
J_{(010)T,t}^{*(0 i_2 0)q}=
(T-t)^{5/2}\Biggl(
\frac{1}{6}\zeta_0^{(i_2)}-\frac{1}{\sqrt{2}\pi^2}
\Biggl(\sum_{r=1}^{q}\frac{1}{r^2}\zeta_{2r}^{(i_2)}+
\sqrt{\beta_q}\mu_q^{(i_2)}\Biggr)\Biggr),
$$

\vspace{3mm}
$$
J_{(100)T,t}^{*(i_1 0 0)q}=
(T-t)^{5/2}\Biggl(
\frac{1}{6}\zeta_0^{(i_1)}+\frac{1}{2\sqrt{2}\pi^2}\Biggl(
\sum_{r=1}^{q}\frac{1}{r^2}\zeta_{2r}^{(i_1)}+
\sqrt{\beta_q}\mu_q^{(i_1)}\Biggr)+\Biggr.
$$
$$
\Biggl.+
\frac{1}{2\sqrt{2}\pi}\Biggl(\sum_{r=1}^q
\frac{1}{r}\zeta_{2r-1}^{(i_1)}+\sqrt{\alpha_q}\xi_q^{(i_1)}\Biggr)\Biggr),
$$

\vspace{7mm}

$$
J_{(011)T,t}
^{*(0 i_2 i_3)q}=(T-t)^{2}\Biggl(\frac{1}{6}
\zeta_{0}^{(i_2)}\zeta_{0}^{(i_3)}-\frac{1}{2\sqrt{2}\pi}
\sqrt{\alpha_q}\xi_q^{(i_3)}\zeta_0^{(i_2)}+\Biggr.
$$
$$
+\frac{1}{2\sqrt{2}\pi^2}\sqrt{\beta_q}\Biggl(
\mu_q^{(i_3)}\zeta_0^{(i_2)}-2\mu_q^{(i_2)}\zeta_0^{(i_3)}\Biggr)+
$$
$$
+\frac{1}{2\sqrt{2}}\sum_{r=1}^{q}
\Biggl(-\frac{1}{\pi r}
\zeta_{2r-1}^{(i_3)}
\zeta_{0}^{(i_2)}+
\frac{1}{\pi^2 r^2}\left(
\zeta_{2r}^{(i_3)}
\zeta_{0}^{(i_2)}-
2\zeta_{2r}^{(i_2)}
\zeta_{0}^{(i_3)}\right)\Biggr)-
$$
$$
-
\frac{1}{2\pi^2}\sum_{\stackrel{r,l=1}{{}_{r\ne l}}}^{q}
\frac{1}{r^2-l^2}\Biggl(
\zeta_{2r}^{(i_2)}
\zeta_{2l}^{(i_3)}+
\frac{l}{r}
\zeta_{2r-1}^{(i_2)}
\zeta_{2l-1}^{(i_3)}
\Biggr)+
$$
$$
+
\sum_{r=1}^{q}
\Biggl(\frac{1}{4\pi r}\left(
\zeta_{2r}^{(i_2)}
\zeta_{2r-1}^{(i_3)}-
\zeta_{2r-1}^{(i_2)}
\zeta_{2r}^{(i_3)}\right)+\Biggr.
$$
\begin{equation}
\label{9000back}
\Biggl.\Biggl.+
\frac{1}{8\pi^2 r^2}\left(
3\zeta_{2r-1}^{(i_2)}
\zeta_{2r-1}^{(i_3)}+
\zeta_{2r}^{(i_3)}
\zeta_{2r}^{(i_2)}\right)\Biggr)\Biggr),
\end{equation}

\vspace{5mm}

$$
J_{(110)T,t}
^{*(i_1 i_2 0)q}=(T-t)^{2}\Biggl(\frac{1}{6}
\zeta_{0}^{(i_1)}\zeta_{0}^{(i_2)}+
\frac{1}{2\sqrt{2}\pi}
\sqrt{\alpha_q}\xi_q^{(i_1)}\zeta_0^{(i_2)}+\Biggr.
$$
$$
+\frac{1}{2\sqrt{2}\pi^2}\sqrt{\beta_q}\Biggl(
\mu_q^{(i_1)}\zeta_0^{(i_2)}-2\mu_q^{(i_2)}\zeta_0^{(i_1)}\Biggr)+
$$
$$
+\frac{1}{2\sqrt{2}}\sum_{r=1}^{q}
\Biggl(\frac{1}{\pi r}
\zeta_{2r-1}^{(i_1)}
\zeta_{0}^{(i_2)}+
\frac{1}{\pi^2 r^2}\left(\zeta_{2r}^{(i_1)}
\zeta_{0}^{(i_2)}
-2\zeta_{2r}^{(i_2)}
\zeta_{0}^{(i_1)}\right)\Biggr)+
$$
$$
+
\frac{1}{2\pi^2}\sum_{\stackrel{r,l=1}{{}_{r\ne l}}}^{q}
\frac{1}{r^2-l^2}\Biggl(\frac{r}{l}
\zeta_{2r-1}^{(i_1)}
\zeta_{2l-1}^{(i_2)}
+\zeta_{2r}^{(i_1)}
\zeta_{2l}^{(i_2)}
\Biggr)+
$$
$$
+
\sum_{r=1}^{q}
\Biggl(\frac{1}{4\pi r}\left(\zeta_{2r-1}^{(i_2)}
\zeta_{2r}^{(i_1)}
-\zeta_{2r-1}^{(i_1)}
\zeta_{2r}^{(i_2)}\right)+\Biggr.
$$
$$
\Biggl.\Biggl.+
\frac{1}{8\pi^2 r^2}\left(
3\zeta_{2r-1}^{(i_1)}
\zeta_{2r-1}^{(i_2)}+
\zeta_{2r}^{(i_1)}
\zeta_{2r}^{(i_2)}\right)\Biggr)\Biggr),
$$

\vspace{2mm}

$$
J_{(101)T,t}
^{*(i_1 0 i_3)q}=(T-t)^{2}\Biggl(\frac{1}{6}
\zeta_{0}^{(i_1)}\zeta_{0}^{(i_3)}
+\frac{1}{2\sqrt{2}\pi}\sqrt{\alpha_q}\left(
\xi_q^{(i_1)}\zeta_0^{(i_3)}-\xi_q^{(i_3)}\zeta_0^{(i_1)}\right)+
\Biggr.
$$
$$
+\frac{1}{2\sqrt{2}\pi^2}\sqrt{\beta_q}\left(
\mu_q^{(i_1)}\zeta_0^{(i_3)}+\mu_q^{(i_3)}\zeta_0^{(i_1)}\right)+
$$
$$
+
\frac{1}{2\sqrt{2}}\sum_{r=1}^{q}
\Biggl(\frac{1}{\pi r}
\left(\zeta_{2r-1}^{(i_1)}
\zeta_{0}^{(i_3)}-\zeta_{2r-1}^{(i_3)}
\zeta_{0}^{(i_1)}\right)+\Biggr.
$$
$$
\Biggl.+
\frac{1}{\pi^2 r^2}\left(
\zeta_{2r}^{(i_1)}
\zeta_{0}^{(i_3)}+
\zeta_{2r}^{(i_3)}
\zeta_{0}^{(i_1)}\right)\Biggr)-
\frac{1}{2\pi^2}\sum_{\stackrel{r,l=1}{{}_{r\ne l}}}^{q}
\frac{1}{rl}
\zeta_{2r-1}^{(i_1)}
\zeta_{2l-1}^{(i_3)}-
$$
$$
\Biggl.-
\sum_{r=1}^q\frac{1}{4\pi^2 r^2}\left(
3\zeta_{2r-1}^{(i_1)}
\zeta_{2r-1}^{(i_3)}+
\zeta_{2r}^{(i_1)}
\zeta_{2r}^{(i_3)}\right)\Biggr).
$$

\vspace{4mm}

\section{A Comparative Analysis of Efficiency of Using the 
Legendre Polynomials and Trigonometric Functions for 
the Numerical Solution of It\^{o} SDEs}

The section is devoted to comparative analysis of efficiency of 
application the Legendre polynomials and trigonometric functions for the 
numerical integration of It\^{o} SDEs in the 
framework of the method of approximation of iterated It\^{o} and Stratonovich 
stochastic integrals 
based on generalized multiple Fourier series (Theorems 1.1, 2.1--2.9, 2.30, 2.33--2.36,  
2.50, 2.51, 2.62--2.65). 
This section is written on the base of the papers \cite{art-4},
\cite{arxiv-12}, and \cite{arxiv-4} (Sect.~4).

Using the iterated It\^{o} stochastic integrals of 
multiplicities 1 to 3 appearing in the Taylor--It\^{o} 
expansion as an example, it is shown that their expansions 
obtained using multiple
Fourier--Legendre series are significantly simpler and 
less computationally costly than their analogues
obtained on the basis of multiple trigonometric Fourier series.

Let us consider the following set of iterated It\^{o} 
and Stratonovich stochastic integrals from the
classical Taylor--It\^{o} and Taylor--Stratonovich expansions
\cite{Zapad-3}
\begin{equation}
\label{ito1}
J_{(\lambda_1\ldots \lambda_k)T,t}^{(i_1\ldots i_k)}=
\int\limits_t^T \ldots \int\limits_t^{t_{2}}d{\bf w}_{t_1}^{(i_1)}\ldots
d{\bf w}_{t_k}^{(i_k)},
\end{equation}
\begin{equation}
\label{str1}
J_{(\lambda_1\ldots \lambda_k)T,t}^{*(i_1\ldots i_k)}=
{\int\limits_t^{*}}^T
\ldots 
{\int\limits_t^{*}}^{t_{2}}
d{\bf w}_{t_1}^{(i_1)}\ldots
d{\bf w}_{t_k}^{(i_k)},
\end{equation}

\vspace{1mm}
\noindent
where ${\bf w}_{\tau}^{(i)}$
$(i=1,\ldots,m)$ are independent standard Wiener processes,
${\bf w}_{\tau}^{(0)}=\tau,$
$i_1,\ldots,i_k = 0, 1,\ldots,m,$ $\lambda_l=0$ for $i_l=0$ and
$\lambda_l=1$ for $i_l=1,\ldots,m$ $(l=1,\ldots,k).$

In \cite{Zapad-1} Milstein G.N. obtained the following expansion of 
$J_{(11)T,t}^{(i_1i_2)}$ on the base of the
Karhunen--Lo\`{e}ve expansion of the Brownian bridge
process (we will discuss
the method \cite{Zapad-1} in detail in Sect.~6.2) 
$$
J_{(11)T,t}^{(i_1 i_2)}=\frac{1}{2}(T-t)\Biggl(
\zeta_{0}^{(i_1)}\zeta_{0}^{(i_2)}
+\frac{1}{\pi}
\sum_{r=1}^{\infty}\frac{1}{r}\left(
\zeta_{2r}^{(i_1)}\zeta_{2r-1}^{(i_2)}-
\zeta_{2r-1}^{(i_1)}\zeta_{2r}^{(i_2)}
\right.\Biggr.
+
$$
\begin{equation}
\label{43}
+\Biggl.\left.\sqrt{2}\left(\zeta_{2r-1}^{(i_1)}\zeta_{0}^{(i_2)}-
\zeta_{0}^{(i_1)}\zeta_{2r-1}^{(i_2)}\right)\right)
\Biggr),
\end{equation}
where the series converges in the mean-square sense,
$i_1\ne i_2,$ $i_1, i_2=1,\ldots,m,$
$$
\zeta_{j}^{(i)}=
\int\limits_t^T \phi_{j}(s) d{\bf w}_s^{(i)}
$$ 
are independent standard Gaussian random variables for various
$i$ or $j$  ($i=1,\ldots,m,\ j=0, 1,\ldots $),

\newpage
\noindent
\begin{equation}
\label{666.6}
\phi_j(s)=\frac{1}{\sqrt{T-t}}
\left\{
\begin{matrix}
1\ &\hbox{for}\ j=0\cr\cr
\sqrt{2}{\rm sin}(2\pi r(s-t)/(T-t))\ &\hbox{for}\ j=2r-1\cr\cr
\sqrt{2}{\rm cos}(2\pi r(s-t)/(T-t))\ &\hbox{for}\ j=2r
\end{matrix},\right.
\end{equation}

\noindent
where $r=1, 2,\ldots$

Moreover,
$$
J_{(1)T,t}^{(i_1)}=\sqrt{T-t}\zeta_0^{(i_1)},
$$
where $i_1=1,\ldots,m.$

In principle, for implementing the strong numerical method with the 
convergence order 1.0
(Milstein method \cite{Zapad-1}, see Sect.~4.10) 
for It\^{o} SDEs
we can
take the following approximations
\begin{equation}
\label{411}
J_{(1)T,t}^{(i_1)}=\sqrt{T-t}\zeta_0^{(i_1)},
\end{equation}
$$
J_{(11)T,t}^{(i_1 i_2)q}=\frac{1}{2}(T-t)\Biggl(
\zeta_{0}^{(i_1)}\zeta_{0}^{(i_2)}
+\frac{1}{\pi}
\sum_{r=1}^{q}\frac{1}{r}\left(
\zeta_{2r}^{(i_1)}\zeta_{2r-1}^{(i_2)}-
\zeta_{2r-1}^{(i_1)}\zeta_{2r}^{(i_2)}
\right.\Biggr.
+
$$
\begin{equation}
\label{431}
+\Biggl.\left.\sqrt{2}\left(\zeta_{2r-1}^{(i_1)}\zeta_{0}^{(i_2)}-
\zeta_{0}^{(i_1)}\zeta_{2r-1}^{(i_2)}\right)\right)
\Biggr),
\end{equation}
where $i_1\ne i_2,$ $i_1, i_2=1,\ldots,m.$

It is not difficult to show that
\begin{equation}
\label{801xxxx}
~~~~~~~~~~~ {\sf M}\left\{\left(J_{(11)T,t}^{(i_1 i_2)}-
J_{(11)T,t}^{(i_1 i_2)q}
\right)^2\right\}
=\frac{3(T-t)^{2}}{2\pi^2}\Biggl(\frac{\pi^2}{6}-
\sum_{r=1}^q \frac{1}{r^2}\Biggr).
\end{equation}

However, this approach has an obvious drawback. Indeed, we have too complex 
formulas for the stochastic integrals with Gaussian distribution
\begin{equation}
\label{100www}
J_{(01)T,t}^{(0 i_1)}=
\frac{{(T-t)}^{3/2}}{2}
\biggl(\zeta_0^{(i_1)}-\frac{\sqrt{2}}{\pi}\sum_{r=1}^{\infty}
\frac{1}{r}\zeta_{2r-1}^{(i_1)}\biggr),
\end{equation}
$$
J_{(001)T,t}^{(00 i_1)}=
(T-t)^{5/2}\Biggl(
\frac{1}{6}\zeta_0^{(i_1)}+\frac{1}{2\sqrt{2}\pi^2}
\sum_{r=1}^{\infty}\frac{1}{r^2}\zeta_{2r}^{(i_1)}-
\frac{1}{2\sqrt{2}\pi}\sum_{r=1}^{\infty}
\frac{1}{r}\zeta_{2r-1}^{(i_1)}\Biggr),
$$

\vspace{2mm}
$$
J_{(01)T,t}^{(0 i_1)q}=
\frac{{(T-t)}^{3/2}}{2}
\biggl(\zeta_0^{(i_1)}-\frac{\sqrt{2}}{\pi}\sum_{r=1}^{q}
\frac{1}{r}\zeta_{2r-1}^{(i_1)}\biggr),
$$
$$
J_{(001)T,t}^{(00 i_1)q}=
(T-t)^{5/2}\Biggl(
\frac{1}{6}\zeta_0^{(i_1)}+\frac{1}{2\sqrt{2}\pi^2}
\sum_{r=1}^{q}\frac{1}{r^2}\zeta_{2r}^{(i_1)}
-
\frac{1}{2\sqrt{2}\pi}\sum_{r=1}^q
\frac{1}{r}\zeta_{2r-1}^{(i_1)}\Biggr),
$$

\vspace{1mm}
\noindent
where the meaning of notations used in (\ref{431}) is retained.

In \cite{Zapad-1} Milstein G.N. proposed the following mean-square 
approximations on the base of
(\ref{43}), (\ref{100www})
\begin{equation}
\label{444www}
~~~~~~~~ J_{(01)T,t}^{(0 i_1)q}=\frac{{(T-t)}^{3/2}}{2}
\biggl(\zeta_0^{(i_1)}-\frac{\sqrt{2}}{\pi}\biggl(\sum_{r=1}^{q}
\frac{1}{r}
\zeta_{2r-1}^{(i_1)}+\sqrt{\alpha_q}\xi_q^{(i_1)}\biggr)
\biggr),
\end{equation}

$$
J_{(11)T,t}^{(i_1 i_2)q}=\frac{1}{2}(T-t)\Biggl(
\zeta_{0}^{(i_1)}\zeta_{0}^{(i_2)}
+\frac{1}{\pi}
\sum_{r=1}^{q}\frac{1}{r}\left(
\zeta_{2r}^{(i_1)}\zeta_{2r-1}^{(i_2)}-
\zeta_{2r-1}^{(i_1)}\zeta_{2r}^{(i_2)}+
\right.\Biggr.
$$
\begin{equation}
\label{555ccc}
~~~~~ +\Biggl.\left.\sqrt{2}\left(\zeta_{2r-1}^{(i_1)}\zeta_{0}^{(i_2)}-
\zeta_{0}^{(i_1)}\zeta_{2r-1}^{(i_2)}\right)\right)
+\frac{\sqrt{2}}{\pi}\sqrt{\alpha_q}\left(
\xi_q^{(i_1)}\zeta_0^{(i_2)}-\zeta_0^{(i_1)}\xi_q^{(i_2)}\right)\Biggr),
\end{equation}

\vspace{2mm}
\noindent
where $i_1\ne i_2$ in (\ref{555ccc}), and
\begin{equation}
\label{333}
\xi_q^{(i)}=\frac{1}{\sqrt{\alpha_q}}\sum_{r=q+1}^{\infty}
\frac{1}{r}~\zeta_{2r-1}^{(i)},\ \ \
\alpha_q=\frac{\pi^2}{6}-\sum_{r=1}^q\frac{1}{r^2},
\end{equation}
where
$\zeta_0^{(i)},$ $\zeta_{2r}^{(i)},$
$\zeta_{2r-1}^{(i)},$ $\xi_q^{(i)},$ $r=1,\ldots,q,$
$i=1,\ldots,m$
are independent standard Gaussian random variables.

Obviously, for the approximations (\ref{444www}) and (\ref{555ccc}) we obtain 
\cite{Zapad-1}
$$
{\sf M}\left\{\left(J_{(01)T,t}^{(0 i_1)}-
J_{(01)T,t}^{(0 i_1)q}
\right)^2\right\}=0,
$$
\begin{equation}
\label{8010}
~~~~~~~~~~ {\sf M}\left\{\left(J_{(11)T,t}^{(i_1 i_2)}-
J_{(11)T,t}^{(i_1 i_2)q}
\right)^2\right\}
=\frac{(T-t)^{2}}{2\pi^2}\Biggl(\frac{\pi^2}{6}-
\sum_{r=1}^q \frac{1}{r^2}\Biggr).
\end{equation}

This idea has been developed in \cite{Zapad-2}-\cite{Zapad-4}.
For example, the approximation $J_{(001)T,t}^{(0 0 i_1)q},$ 
which corresponds to (\ref{444www}), (\ref{555ccc})
is defined by \cite{Zapad-2}-\cite{Zapad-4}
$$
J_{(001)T,t}^{(00 i_1)q}=
(T-t)^{5/2}\Biggl(
\frac{1}{6}\zeta_0^{(i_1)}+\frac{1}{2\sqrt{2}\pi^2}\Biggl(
\sum_{r=1}^{q}\frac{1}{r^2}\zeta_{2r}^{(i_1)}+
\sqrt{\beta_q}\mu_q^{(i_1)}\Biggr)-\Biggl.
$$
\begin{equation}
\label{1970}
\Biggl.
-\frac{1}{2\sqrt{2}\pi}\Biggl(\sum_{r=1}^q
\frac{1}{r}\zeta_{2r-1}^{(i_1)}+\sqrt{\alpha_q}\xi_q^{(i_1)}\Biggr)\Biggr),
\end{equation}

\newpage
\noindent
where
$\xi_q^{(i)},$ $\alpha_q$ have the form (\ref{333}),
$$
\mu_q^{(i)}=\frac{1}{\sqrt{\beta_q}}\sum_{r=q+1}^{\infty}
\frac{1}{r^2}~\zeta_{2r}^{(i)},\ \ \
\beta_q=\frac{\pi^4}{90}-\sum_{r=1}^q\frac{1}{r^4},
$$

\noindent
$\phi_j(s)$ is defined by (\ref{666.6}), and
$\zeta_0^{(i)},$ $\zeta_{2r}^{(i)},$
$\zeta_{2r-1}^{(i)},$ $\xi_q^{(i)},$ $\mu_q^{(i)}$ ($r=1,\ldots,q,$
$i=1,\ldots,m$) are independent
standard Gaussian random variables.

Moreover,
$$
{\sf M}\left\{\left(J_{(001)T,t}^{(0 0 i_1)}-
J_{(0 0 1)T,t}^{(0 0 i_1)q}
\right)^2\right\}=0.
$$

Nevetheless, 
the expansions (\ref{444www}), (\ref{1970}) are too complex for
the approximation of two Gaussian random variables
$J_{(01)T,t}^{(0 i_1)},$ $J_{(001)T,t}^{(0 0 i_1)}.$

Further, we will see that the introducing of random 
variables $\xi_q^{(i)}$ and 
$\mu_q^{(i)}$ will sharply 
complicate the approximation of stochastic 
integral $J_{(111)T,t}^{(i_1 i_2 i_3)}$
($i_1,i_2,i_3=1,\ldots,m$). 
This is due to the fact that
the number $q$ is fixed for stochastic integrals 
included into the considered collection. However, it is clear that due 
to the smallness of $T-t$, the number $q$ for $J_{(111)T,t}^{(i_1 i_2 i_3)}$
could be taken significantly 
less than in the formula (\ref{555ccc}). 
This feature is also valid for the formulas (\ref{444www}), (\ref{1970}).

On the other hand, the following very simple formulas are well known
(see (\ref{4001})--(\ref{4003}))
\begin{equation}
\label{4}
J_{(1)T,t}^{(i_1)}=\sqrt{T-t}\zeta_0^{(i_1)},
\end{equation}
\begin{equation}
\label{5}
J_{(01)T,t}^{(0 i_1)}=\frac{(T-t)^{3/2}}{2}\biggl(\zeta_0^{(i_1)}+
\frac{1}{\sqrt{3}}\zeta_1^{(i_1)}\biggr),
\end{equation}
\begin{equation}
\label{6}
J_{(001)T,t}^{(0 0 i_1)}=\frac{(T-t)^{5/2}}{6}\biggl(\zeta_0^{(i_1)}+
\frac{\sqrt{3}}{2}\zeta_1^{(i_1)}+
\frac{1}{2\sqrt{5}}\zeta_2^{(i_1)}\biggr),
\end{equation}
where
$\zeta_0^{(i)},$ $\zeta_1^{(i)},$ $\zeta_2^{(i)}$ ($i=1,\ldots,m$) 
are indepentent standard Gaussian random variables.
Obviously, that the formulas (\ref{4})-(\ref{6}) 
are part of the method based on Theorem 1.1
(also see Sect.~5.1).

To obtain the Milstein expansion for the stochastic integral (\ref{str-str}) 
the truncated 
expansions of components of the 
Wiener process ${\bf w}_s$ must be
iteratively substituted in the single integrals in (\ref{str-str}), 
and the integrals
must be calculated starting from the innermost integral.
This is a complicated procedure that obviously does not lead to a general
expansion of (\ref{str-str}) valid for an arbitrary multiplicity $k.$
For this reason, only expansions of simplest single, double, and triple
integrals (\ref{str-str}) were obtained \cite{Zapad-1}-\cite{Zapad-4}, 
\cite{Zapad-8}, \cite{Zapad-9} by the Milstein approach \cite{Zapad-1}
based on the Karhunen--Lo\`{e}ve expansion of the Brownian bridge process.

At that, in \cite{Zapad-1}, \cite{Zapad-8} the case 
$\psi_1(s),$ $\psi_2(s)\equiv 1$ and
$i_1, i_2=0, 1,\ldots,m$ $(i_1\ne i_2)$ is considered. In 
\cite{Zapad-2}-\cite{Zapad-4}, \cite{Zapad-9} 
the attempt to consider the case 
$\psi_1(s),$ $\psi_2(s),$ $\psi_3(s)$ $\equiv 1$ and 
$i_1, i_2, i_3=0, 1,\ldots,m$ is realized.
Note that, generally speaking, the mean-square
convergence of $J_{(111)T,t}^{*(i_1 i_2 i_3)q}$ to
$J_{(111)T,t}^{*(i_1 i_2 i_3)}$ if $q\to\infty$
was not proved rigorously in
\cite{Zapad-2}-\cite{Zapad-4}, \cite{Zapad-9}
within the 
frames of the Milstein approach \cite{Zapad-1}
together with the Wong--Zakai approximation \cite{W-Z-1}-\cite{Watanabe}
(see discussions in Sect.~2.42, 2.43, 6.2).

\subsection{A Comparative Analysis of Efficiency of Using the 
Legendre Polynomials and Trigonometric Functions for 
the Integral $J_{(11)T,t}^{(i_1 i_2)}$}

Using Theorem 1.1 and complete orthonormal system of 
Legendre polynomials
in the space $L_2([t, T])$, we have (see (\ref{4004yes}))
\begin{equation}
\label{4004xxx}
J_{(11)T,t}^{(i_1 i_2)}=
\frac{T-t}{2}\Biggl(\zeta_0^{(i_1)}\zeta_0^{(i_2)}+\sum_{i=1}^{\infty}
\frac{1}{\sqrt{4i^2-1}}\left(
\zeta_{i-1}^{(i_1)}\zeta_{i}^{(i_2)}-
\zeta_i^{(i_1)}\zeta_{i-1}^{(i_2)}\right)-{\bf 1}_{\{i_1=i_2\}}\Biggr),
\end{equation}

\noindent
where series converges in the mean-square sense, $i_1, i_2=1,\ldots,m,$
$$
\zeta_{j}^{(i)}=
\int\limits_t^T \phi_{j}(s) d{\bf w}_s^{(i)}
$$ 
are independent standard Gaussian random variables
for various
$i$ or $j$,
\begin{equation}
\label{4009ququ}
~~~~~~~~~ \phi_j(x)=\sqrt{\frac{2j+1}{T-t}}P_j\left(\left(
x-\frac{T+t}{2}\right)\frac{2}{T-t}\right),\ \ \ j=0, 1, 2,\ldots,
\end{equation}
where $P_j(x)$ is the Legendre polynomial.

The formula (\ref{4004xxx}) has been derived 
for the first time in \cite{old-art-1} (1997)
with using Theorem 2.10.

Remind the formula (\ref{fff09}) \cite{old-art-1} (1997)
\begin{equation}
\label{400}
~~~~~~~~ {\sf M}\left\{\left(J_{(11)T,t}^{(i_1 i_2)}-
J_{(11)T,t}^{(i_1 i_2)q}
\right)^2\right\}
=\frac{(T-t)^2}{2}\Biggl(\frac{1}{2}-\sum_{i=1}^q
\frac{1}{4i^2-1}\Biggr),
\end{equation}
where
\begin{equation}
\label{401}
J_{(11)T,t}^{(i_1 i_2)q}=
\frac{T-t}{2}\Biggl(\zeta_0^{(i_1)}\zeta_0^{(i_2)}+\sum_{i=1}^{q}
\frac{1}{\sqrt{4i^2-1}}\left(
\zeta_{i-1}^{(i_1)}\zeta_{i}^{(i_2)}-
\zeta_i^{(i_1)}\zeta_{i-1}^{(i_2)}\right)-{\bf 1}_{\{i_1=i_2\}}\Biggr).
\end{equation}

Let us compare (\ref{401}) with (\ref{555ccc}) and (\ref{400}) with (\ref{8010}).
Consider minimal natural numbers $q_{\rm trig}$ and 
$q_{\rm pol},$ which satisfy to (see Table 5.42)
\begin{equation}
\label{uslov1}
\frac{(T-t)^2}{2}\Biggl(\frac{1}{2}-\sum_{i=1}^{q_{\rm pol}}
\frac{1}{4i^2-1}\Biggr)\le (T-t)^3,
\end{equation}
$$
\frac{(T-t)^{2}}{2\pi^2}\Biggl(\frac{\pi^2}{6}-
\sum_{r=1}^{q_{\rm trig}}\frac{1}{r^2}\Biggr)\le (T-t)^3.
$$

Thus, we have
$$
\frac{q_{\rm pol}}{q_{\rm trig}}\ \ \approx\ \ 
1.67,\ \ 2.22,\ \ 2.43,\ \ 2.36,\ \ 2.41,\ 
\ 2.43,\ \ 2.45,\ \ 2.45.
$$

From the other hand, the formula (\ref{555ccc}) 
includes $(4q+4)m$ independent
standard Gaussian random variables. At the same time the folmula
(\ref{401}) includes only $(2q+2)m$ independent
standard Gaussian random variables. Moreover, the formula
(\ref{401}) is simpler than the formula (\ref{555ccc}).
Thus, in this case we can talk about approximately equal computational costs
for the formulas (\ref{555ccc}) and (\ref{401}).

\begin{table}
\centering
\caption{Numbers $q_{\rm trig},$ $q_{\rm trig}^{*}$, $q_{\rm pol}$}
\label{tab:5.42}      
\begin{tabular}{p{1.1cm}p{1.1cm}p{1.1cm}p{1.1cm}p{1.1cm}p{1.1cm}p{1.1cm}p{1.1cm}p{1.1cm}}
\hline\noalign{\smallskip}
$T-t$&$2^{-5}$&$2^{-6}$&$2^{-7}$&$2^{-8}$&$2^{-9}$&$2^{-10}$&$2^{-11}$&$2^{-12}$\\
\noalign{\smallskip}\hline\noalign{\smallskip}
$q_{\rm trig}$&3&4&7&14&27&53&105&209\\
\noalign{\smallskip}\hline\noalign{\smallskip}
$q_{\rm trig}^{*}$&6&11&20&40&79&157&312&624\\
\noalign{\smallskip}\hline\noalign{\smallskip}
$q_{\rm pol}$&5&9&17&33&65&129&257&513\\
\noalign{\smallskip}\hline\noalign{\smallskip}
\end{tabular}
\end{table}

There is one important feature. 
As we mentioned above, further we will see that the introducing of random 
variables $\xi_q^{(i)}$ and 
$\mu_q^{(i)}$ will sharply 
complicate the approximation of 
stochastic integral $J_{(111)T,t}^{(i_1 i_2 i_3)}$
($i_1,i_2,i_3=1,\ldots,m$). 
This is due to the fact that
the number $q$ is fixed for all stochastic integrals, which 
included into the considered collection. However, it is clear that due 
to the smallness of $T-t$, the number $q$ for $J_{(111)T,t}^{(i_1 i_2 i_3)}$
could be chosen significantly 
less than in the formula (\ref{555ccc}). 
This feature is also valid for the formulas (\ref{444www}), (\ref{1970}).
However, for the case of Legendre polynomials we can choose different 
numbers $q$ for different stochastic integrals.

From the other hand, if we will not introduce the random 
variables $\xi_q^{(i)}$ and 
$\mu_q^{(i)},$ then the mean-square error of approximation of the
stochastic integral $J_{(11)T,t}^{(i_1 i_2)}$ will be three times larger
(see (\ref{801xxxx})). 
Moreover, in this case the stochastic integrals 
$J_{(01)T,t}^{(0 i_1)}$, $J_{(001)T,t}^{(00i_1)}$ (with Gaussian distribution)
will be approximated worse.

Consider minimal natural numbers $q_{\rm trig}^{*}$, 
which satisfy to (see Table 5.42)
$$
\frac{3(T-t)^{2}}{2\pi^2}\Biggl(\frac{\pi^2}{6}-
\sum_{r=1}^{q_{\rm trig}^{*}}\frac{1}{r^2}\Biggr)\le (T-t)^3.
$$

In this situation we can talk about
the advantage of Legendre polynomials ($q_{\rm trig}^{*} > q_{pol}$ and
(\ref{555ccc}) is more complex than (\ref{401})).

\subsection{A Comparative Analysis of Efficiency of Using the 
Legendre Polynomials and Trigonometric Functions for 
the Integrals $J_{(1)T,t}^{(i_1)},$
$J_{(11)T,t}^{(i_1 i_2)},$ $J_{(01)T,t}^{(0 i_1)},$ $J_{(10)T,t}^{(i_1 0)},$
$J_{(111)T,t}^{(i_1 i_2 i_3)}$}

It is well known \cite{Zapad-1}-\cite{Zapad-4}, \cite{Zapad-8} 
(also see \cite{12a}-\cite{12aa})
that for the numerical realization
of strong Taylor--It\^{o} numerical methods with the 
convergence order 1.5 for It\^{o} SDEs
we need to
approximate the following collection of iterated It\^{o} stochastic integrals
(see Sect.~4.10)
$$
J_{(1)T,t}^{(i_1)},\ \ \
J_{(11)T,t}^{(i_1 i_2)},\ \ \ J_{(01)T,t}^{(0 i_1)},\ \ \ J_{(10)T,t}^{(i_1 0)},\ \ \
J_{(111)T,t}^{(i_1 i_2 i_3)}.
$$

Using Theorem 1.1 for the system of trigonometric 
functions, we have (see Sect.~5.2)
\begin{equation}
\label{x1}
J_{(1)T,t}^{(i_1)}=\sqrt{T-t}\zeta_0^{(i_1)},
\end{equation}

$$
J_{(11)T,t}^{(i_1 i_2)q}=\frac{1}{2}(T-t)\Biggl(
\zeta_{0}^{(i_1)}\zeta_{0}^{(i_2)}
+\frac{1}{\pi}
\sum_{r=1}^{q}\frac{1}{r}\left(
\zeta_{2r}^{(i_1)}\zeta_{2r-1}^{(i_2)}-
\zeta_{2r-1}^{(i_1)}\zeta_{2r}^{(i_2)}+
\right.\Biggr.
$$
$$
+\left.\sqrt{2}\left(\zeta_{2r-1}^{(i_1)}\zeta_{0}^{(i_2)}-
\zeta_{0}^{(i_1)}\zeta_{2r-1}^{(i_2)}\right)\right)
+
$$
\begin{equation}
\label{x2}
\Biggl.+\frac{\sqrt{2}}{\pi}\sqrt{\alpha_q}\left(
\xi_q^{(i_1)}\zeta_0^{(i_2)}-\zeta_0^{(i_1)}\xi_q^{(i_2)}\right)-
{\bf 1}_{\{i_1=i_2\}}\Biggr),
\end{equation}

\newpage
\noindent
\begin{equation}
\label{x3}
~~~~~~~J_{(01)T,t}^{(0 i_1)q}=\frac{{(T-t)}^{3/2}}{2}
\biggl(\zeta_0^{(i_1)}-\frac{\sqrt{2}}{\pi}\biggl(\sum_{r=1}^{q}
\frac{1}{r}
\zeta_{2r-1}^{(i_1)}+\sqrt{\alpha_q}\xi_q^{(i_1)}\biggr)
\biggr),
\end{equation}

\begin{equation}
\label{x4}
~~~~~~~J_{(10)T,t}^{(i_1 0)q}=\frac{{(T-t)}^{3/2}}{2}
\biggl(\zeta_0^{(i_1)}+\frac{\sqrt{2}}{\pi}\biggl(\sum_{r=1}^{q}
\frac{1}{r}
\zeta_{2r-1}^{(i_1)}+\sqrt{\alpha_q}\xi_q^{(i_1)}\biggr)
\biggr),
\end{equation}

\vspace{1mm}

$$
J_{(111)T,t}^{(i_1 i_2 i_3)q}=(T-t)^{3/2}\Biggl(\frac{1}{6}
\zeta_{0}^{(i_1)}\zeta_{0}^{(i_2)}\zeta_{0}^{(i_3)}+\Biggr.
\frac{\sqrt{\alpha_q}}{2\sqrt{2}\pi}\left(
\xi_q^{(i_1)}\zeta_0^{(i_2)}\zeta_0^{(i_3)}-\xi_q^{(i_3)}\zeta_0^{(i_2)}
\zeta_0^{(i_1)}\right)+
$$
$$
+\frac{1}{2\sqrt{2}\pi^2}\sqrt{\beta_q}\left(
\mu_q^{(i_1)}\zeta_0^{(i_2)}\zeta_0^{(i_3)}-2\mu_q^{(i_2)}\zeta_0^{(i_1)}
\zeta_0^{(i_3)}+\mu_q^{(i_3)}\zeta_0^{(i_1)}\zeta_0^{(i_2)}\right)+
$$
$$
+
\frac{1}{2\sqrt{2}}\sum_{r=1}^{q}
\Biggl(\frac{1}{\pi r}\left(
\zeta_{2r-1}^{(i_1)}
\zeta_{0}^{(i_2)}\zeta_{0}^{(i_3)}-
\zeta_{2r-1}^{(i_3)}
\zeta_{0}^{(i_2)}\zeta_{0}^{(i_1)}\right)+\Biggr.
$$
$$
\Biggl.+
\frac{1}{\pi^2 r^2}\left(
\zeta_{2r}^{(i_1)}
\zeta_{0}^{(i_2)}\zeta_{0}^{(i_3)}-
2\zeta_{2r}^{(i_2)}
\zeta_{0}^{(i_3)}\zeta_{0}^{(i_1)}+
\zeta_{2r}^{(i_3)}
\zeta_{0}^{(i_2)}\zeta_{0}^{(i_1)}\right)\Biggr)+
$$
$$
+
\sum_{r=1}^{q}
\Biggl(\frac{1}{4\pi r}\left(
\zeta_{2r}^{(i_1)}
\zeta_{2r-1}^{(i_2)}\zeta_{0}^{(i_3)}-
\zeta_{2r-1}^{(i_1)}
\zeta_{2r}^{(i_2)}\zeta_{0}^{(i_3)}-
\zeta_{2r-1}^{(i_2)}
\zeta_{2r}^{(i_3)}\zeta_{0}^{(i_1)}+
\zeta_{2r-1}^{(i_3)}
\zeta_{2r}^{(i_2)}\zeta_{0}^{(i_1)}\right)+\Biggr.
$$
$$
+
\frac{1}{8\pi^2 r^2}\left(
3\zeta_{2r-1}^{(i_1)}
\zeta_{2r-1}^{(i_2)}\zeta_{0}^{(i_3)}+
\zeta_{2r}^{(i_1)}
\zeta_{2r}^{(i_2)}\zeta_{0}^{(i_3)}-
6\zeta_{2r-1}^{(i_1)}
\zeta_{2r-1}^{(i_3)}\zeta_{0}^{(i_2)}+\right.
$$
\begin{equation}
\label{x5}
\Biggl.\left.
~~~~~~~ +
3\zeta_{2r-1}^{(i_2)}
\zeta_{2r-1}^{(i_3)}\zeta_{0}^{(i_1)}-
2\zeta_{2r}^{(i_1)}
\zeta_{2r}^{(i_3)}\zeta_{0}^{(i_2)}+
\zeta_{2r}^{(i_3)}
\zeta_{2r}^{(i_2)}\zeta_{0}^{(i_1)}\right)\Biggr)
\Biggl.+D_{T,t}^{(i_1i_2i_3)q}\Biggr),
\end{equation}

\vspace{4mm}
\noindent
where in (\ref{x5})
we suppose that $i_1\ne i_2,$ $i_1\ne i_3,$ $i_2\ne i_3$,

\vspace{1mm}
$$
D_{T,t}^{(i_1i_2i_3)q}=
\frac{1}{2\pi^2}\sum_{\stackrel{r,l=1}{{}_{r\ne l}}}^{q}
\Biggl(\frac{1}{r^2-l^2}\biggl(
\zeta_{2r}^{(i_1)}
\zeta_{2l}^{(i_2)}\zeta_{0}^{(i_3)}-
\zeta_{2r}^{(i_2)}
\zeta_{0}^{(i_1)}\zeta_{2l}^{(i_3)}+\biggr.\Biggr.
$$
$$
\Biggl.+\biggl.
\frac{r}{l}
\zeta_{2r-1}^{(i_1)}
\zeta_{2l-1}^{(i_2)}\zeta_{0}^{(i_3)}-\frac{l}{r}
\zeta_{0}^{(i_1)}
\zeta_{2r-1}^{(i_2)}\zeta_{2l-1}^{(i_3)}\biggr)-
\frac{1}{rl}\zeta_{2r-1}^{(i_1)}
\zeta_{0}^{(i_2)}\zeta_{2l-1}^{(i_3)}\Biggr)+
$$
$$
+
\frac{1}{4\sqrt{2}\pi^2}\Biggl(
\sum_{r,m=1}^{q}\Biggl(\frac{2}{rm}
\left(-\zeta_{2r-1}^{(i_1)}
\zeta_{2m-1}^{(i_2)}\zeta_{2m}^{(i_3)}+
\zeta_{2r-1}^{(i_1)}
\zeta_{2r}^{(i_2)}\zeta_{2m-1}^{(i_3)}+
\right.\Biggr.\Biggr.
$$
$$
\left.+
\zeta_{2r-1}^{(i_1)}
\zeta_{2m}^{(i_2)}\zeta_{2m-1}^{(i_3)}-
\zeta_{2r}^{(i_1)}
\zeta_{2r-1}^{(i_2)}\zeta_{2m-1}^{(i_3)}\right)+
$$
$$
+\frac{1}{m(r+m)}
\left(-\zeta_{2(m+r)}^{(i_1)}
\zeta_{2r}^{(i_2)}\zeta_{2m}^{(i_3)}-
\zeta_{2(m+r)-1}^{(i_1)}
\zeta_{2r-1}^{(i_2)}\zeta_{2m}^{(i_3)}-
\right.
$$
$$
\Biggl.\left.
-\zeta_{2(m+r)-1}^{(i_1)}
\zeta_{2r}^{(i_2)}\zeta_{2m-1}^{(i_3)}+
\zeta_{2(m+r)}^{(i_1)}
\zeta_{2r-1}^{(i_2)}\zeta_{2m-1}^{(i_3)}\right)\Biggr)+
$$
$$
+
\sum_{m=1}^{q}\sum_{l=m+1}^{q}\Biggl(\frac{1}{m(l-m)}
\left(\zeta_{2(l-m)}^{(i_1)}
\zeta_{2l}^{(i_2)}\zeta_{2m}^{(i_3)}+
\zeta_{2(l-m)-1}^{(i_1)}
\zeta_{2l-1}^{(i_2)}\zeta_{2m}^{(i_3)}-
\right.\Biggr.
$$
$$
\left.
-\zeta_{2(l-m)-1}^{(i_1)}
\zeta_{2l}^{(i_2)}\zeta_{2m-1}^{(i_3)}+
\zeta_{2(l-m)}^{(i_1)}
\zeta_{2l-1}^{(i_2)}\zeta_{2m-1}^{(i_3)}\right)+
$$
$$
+
\frac{1}{l(l-m)}
\left(-\zeta_{2(l-m)}^{(i_1)}
\zeta_{2m}^{(i_2)}\zeta_{2l}^{(i_3)}+
\zeta_{2(l-m)-1}^{(i_1)}
\zeta_{2m-1}^{(i_2)}\zeta_{2l}^{(i_3)}-
\right.
$$
$$
\Biggl.
\Biggl.
\Biggl.
\left.
-\zeta_{2(l-m)-1}^{(i_1)}
\zeta_{2m}^{(i_2)}\zeta_{2l-1}^{(i_3)}-
\zeta_{2(l-m)}^{(i_1)}
\zeta_{2m-1}^{(i_2)}\zeta_{2l-1}^{(i_3)}\right)\Biggr)\Biggr),
$$

\vspace{3mm}
\noindent
where
$$
\xi_q^{(i)}=\frac{1}{\sqrt{\alpha_q}}\sum_{r=q+1}^{\infty}
\frac{1}{r}~\zeta_{2r-1}^{(i)},\ \ \
\alpha_q=\frac{\pi^2}{6}-\sum_{r=1}^q\frac{1}{r^2},
$$
$$
\mu_q^{(i)}=\frac{1}{\sqrt{\beta_q}}\sum_{r=q+1}^{\infty}
\frac{1}{r^2}~\zeta_{2r}^{(i)},\ \ \
\beta_q=\frac{\pi^4}{90}-\sum_{r=1}^q\frac{1}{r^4},\
$$

\noindent
and $\zeta_0^{(i)},$ $\zeta_{2r}^{(i)},$
$\zeta_{2r-1}^{(i)},$ $\xi_q^{(i)},$ $\mu_q^{(i)}$ 
($r=1,\ldots,q,$
$i=1,\ldots,m$) are independent
standard Gaussian random variables.

The mean-square errors of approximations (\ref{x2})--(\ref{x5})
are represented by the formulas
$$
{\sf M}\left\{\left(J_{(01)T,t}^{(0 i_1)}-
J_{(01)T,t}^{(0 i_1)q}
\right)^2\right\}=0,
$$

\vspace{-1mm}
$$
{\sf M}\left\{\left(J_{(10)T,t}^{(i_1 0)}-
J_{(10)T,t}^{(i_1 0)q}
\right)^2\right\}=0,
$$

\vspace{-4mm}
\begin{equation}
\label{y1ssss}
~~~~~~~ {\sf M}\left\{\left(J_{(11)T,t}^{(i_1 i_2)}-
J_{(11)T,t}^{(i_1 i_2)q}
\right)^2\right\}
=\frac{(T-t)^{2}}{2\pi^2}\Biggl(\frac{\pi^2}{6}-
\sum_{r=1}^q \frac{1}{r^2}\Biggr),
\end{equation}

\vspace{2mm}
$$
{\sf M}\left\{\left(J_{(111)T,t}^{(i_1 i_2 i_3)}-
J_{(111)T,t}^{(i_1 i_2 i_3)q}\right)^2\right\}=
(T-t)^3\Biggl(\frac{4}{45}-\frac{1}{4\pi^2}\sum_{r=1}^q\frac{1}{r^2}-
\Biggl.
$$
\begin{equation}
\label{x7}
\Biggl.-\frac{55}{32\pi^4}\sum_{r=1}^q\frac{1}{r^4}-
\frac{1}{4\pi^4}\sum_{\stackrel{r,l=1}{{}_{r\ne l}}}^q
\frac{5l^4+4r^4-3r^2l^2}{r^2 l^2 \left(r^2-l^2\right)^2}\Biggr),
\end{equation}

\vspace{1mm}
\noindent
where $i_1\ne i_2,$ $i_1\ne i_3,$ $i_2\ne i_3$.

\begin{table}
\centering
\caption{Confirmation of the formula (\ref{x7})}
\label{tab:5.43}      
\begin{tabular}{p{2.1cm}p{1.7cm}p{1.7cm}p{2.1cm}p{2.3cm}p{2.3cm}p{2.3cm}}
\hline\noalign{\smallskip}
$\varepsilon/(T-t)^3$&0.0459&0.0072&$7.5722\cdot 10^{-4}$
&$7.5973\cdot 10^{-5}$&
$7.5990\cdot 10^{-6}$\\
\noalign{\smallskip}\hline\noalign{\smallskip}
$q$&1&10&100&1000&10000\\
\noalign{\smallskip}\hline\noalign{\smallskip}
\end{tabular}
\end{table}

In Table 5.43 we can see the numerical confirmation of 
the formula (\ref{x7}) ($\varepsilon$ means the 
right-hand side of (\ref{x7})).

Note that the formulas (\ref{x1}), (\ref{x2}) have been obtained 
for the first time in \cite{Zapad-1}. Using 
(\ref{x1}), (\ref{x2}), we can realize numerically 
an explicit one-step strong numerical 
method with the convergence order 1.0 for It\^{o} SDEs
(Milstein method \cite{Zapad-1}; also see Sect.~4.10). 

An analogue of the formula (\ref{x5}) 
has been obtained for the first time in 
\cite{Zapad-2}, \cite{Zapad-3}. 

As we mentioned above, the Milstein expansion (i.e. expansion
based on the Karhunen--Lo\`{e}ve expansion of the Brownian bridge
process)
for iterated stochastic integrals leads to iterated application
of the operation of 
limit transition. An analogue of (\ref{x5}) for iterated
Stratonovich stochastic integrals has been derived in 
\cite{Zapad-2}, \cite{Zapad-3} on the base of the Milstein expansion
together with the Wong--Zakai approximation \cite{W-Z-1}-\cite{Watanabe}
(without rigorous proof).
It means that the authors in \cite{Zapad-2}, \cite{Zapad-3} formally
could not use the double sum with the upper limit $q$ in the analogue of
(\ref{x5}). From the other hand, the correctness of (\ref{x5}) follows 
directly from Theorem 1.1. Note that (\ref{x5}) has been
obtained  reasonably
for the first time in \cite{1}. The version of (\ref{x5}) but without
the introducing of random variables $\xi_q^{(i)}$ and $\mu_q^{(i)}$ 
can be found in \cite{old-art-1} (1997).

Note that the formula (\ref{y1ssss}) 
appears for the first time in \cite{Zapad-1}.
The mean-square error (\ref{x7}) has been obtained for the first time
in \cite{very-old-2} (1996) on the base of the simplified variant of 
Theorem 1.1 (the case of pairwise different $i_1,\ldots,i_k$).

The number $q$ as we noted above must be the same
in (\ref{x2})--(\ref{x5}). This is the main drawback of this approach,
because really the number $q$ in (\ref{x5})
can be chosen essentially smaller than in (\ref{x2}).

Note that in (\ref{x5}) we can replace 
$J_{(111)T,t}^{(i_1 i_2 i_3)q}$ with $J_{(111)T,t}^{*(i_1 i_2 i_3)q}$
and (\ref{x5}) then will be valid for any $i_1, i_2, i_3 = 0, 1,\ldots,m$
(see Theorems 2.6--2.8).

Consider now approximations of iterated stochastic integrals 
$$
J_{(1)T,t}^{(i_1)},\ \ \
J_{(11)T,t}^{(i_1 i_2)},\ \ \ J_{(01)T,t}^{(0 i_1)},\ \ \
J_{(10)T,t}^{(i_1 0)},\ \ \
J_{(111)T,t}^{(i_1 i_2 i_3)}\ \ \ \ (i_1,i_2,i_3=1,\ldots,m)
$$
on the base of Theorem 1.1 (the case of Legendre 
polynomials) \cite{1}-\cite{12aa}, \cite{arxiv-4}

\vspace{-1mm}
\begin{equation}
\label{y2}
J_{(1)T,t}^{(i_1)}=\sqrt{T-t}\zeta_0^{(i_1)},
\end{equation}

\vspace{-4mm}
\begin{equation}
\label{y3}
J_{(11)T,t}^{(i_1 i_2)q}=
\frac{T-t}{2}\Biggl(\zeta_0^{(i_1)}\zeta_0^{(i_2)}+\sum_{i=1}^{q}
\frac{1}{\sqrt{4i^2-1}}\left(
\zeta_{i-1}^{(i_1)}\zeta_{i}^{(i_2)}-
\zeta_i^{(i_1)}\zeta_{i-1}^{(i_2)}\right)-{\bf 1}_{\{i_1=i_2\}}\Biggr),
\end{equation}

\vspace{-4mm}
\begin{equation}
\label{y3a}
J_{(01)T,t}^{(0 i_1)}=\frac{(T-t)^{3/2}}{2}\biggl(\zeta_0^{(i_1)}+
\frac{1}{\sqrt{3}}\zeta_1^{(i_1)}\biggr),
\end{equation}

\begin{equation}
\label{y4a}
J_{(10)T,t}^{(i_1 0)}=\frac{(T-t)^{3/2}}{2}\biggl(\zeta_0^{(i_1)}-
\frac{1}{\sqrt{3}}\zeta_1^{(i_1)}\biggr),
\end{equation}

\vspace{3mm}
$$
J_{(111)T,t}^{(i_1 i_2 i_3)q_1}=
\sum_{j_1,j_2,j_3=0}^{q_1}
C_{j_3j_2j_1}
\Biggl(
\zeta_{j_1}^{(i_1)}\zeta_{j_2}^{(i_2)}\zeta_{j_3}^{(i_3)}
-{\bf 1}_{\{i_1=i_2\}}
{\bf 1}_{\{j_1=j_2\}}
\zeta_{j_3}^{(i_3)}-\Biggr.
$$
\begin{equation}
\label{y4}
~~~~~~~~\Biggl.
-{\bf 1}_{\{i_2=i_3\}}
{\bf 1}_{\{j_2=j_3\}}
\zeta_{j_1}^{(i_1)}-
{\bf 1}_{\{i_1=i_3\}}
{\bf 1}_{\{j_1=j_3\}}
\zeta_{j_2}^{(i_2)}\Biggr),\ \ \ q_1\ll q,
\end{equation}

\vspace{4mm}

$$
J_{(111)T,t}^{(i_1 i_1 i_1)}
=\frac{1}{6}(T-t)^{3/2}\left(
\left(\zeta_0^{(i_1)}\right)^3-3
\zeta_0^{(i_1)}\right),
$$

\vspace{2mm}
\noindent
where
$$
C_{j_3j_2j_1}=\int\limits_{t}^{T}\phi_{j_3}(z)
\int\limits_{t}^{z} \phi_{j_2}(y)
\int\limits_{t}^{y}
\phi_{j_1}(x) dx dy dz=
$$
$$
=\frac{\sqrt{(2j_1+1)(2j_2+1)(2j_3+1)}}{8}(T-t)^{3/2}\bar
C_{j_3j_2j_1},
$$
\begin{equation}
\label{ww}
\bar C_{j_3j_2j_1}=\int\limits_{-1}^{1}P_{j_3}(z)
\int\limits_{-1}^{z}P_{j_2}(y)
\int\limits_{-1}^{y}
P_{j_1}(x)dx dy dz,
\end{equation}

\vspace{1mm}
\noindent
$\phi_{j}(x)$ is defined by (\ref{4009ququ}) 
and $P_i(x)$ is the Legendre polynomial
$(i=0, 1, 2,\ldots)$.

The mean-square errors of approximations (\ref{y3}), (\ref{y4})
are represented by the 
formulas (see Theo\-rems 1.3 and 1.4; also see Sect.~5.1) 

\vspace{-4mm}
\begin{equation}
\label{sad004}
~~~~ {\sf M}\left\{\left(J_{(11)T,t}^{(i_1 i_2)}-
J_{(11)T,t}^{(i_1 i_2)q}
\right)^2\right\}
=\frac{(T-t)^2}{2}\Biggl(\frac{1}{2}-\sum_{i=1}^q
\frac{1}{4i^2-1}\Biggr)\ \ \ (i_1\ne i_2),
\end{equation}

\vspace{1mm}
$$
{\sf M}\left\{\left(
J_{(111)T,t}^{(i_1 i_2 i_3)}-
J_{(111)T,t}^{(i_1 i_2 i_3)q_1}\right)^2\right\}=
$$
\begin{equation}
\label{sad005}
~~~~~~~~ =\frac{(T-t)^{3}}{6}-
\sum_{j_3,j_2,j_1=0}^{q_1}
C_{j_3j_2j_1}^2\ \ \ (i_1\ne i_2, i_1\ne i_3, i_2\ne i_3),
\end{equation}

\vspace{4mm}

$$
{\sf M}\left\{\left(
J_{(111)T,t}^{(i_1 i_2 i_3)}-
J_{(111)T,t}^{(i_1 i_2 i_3)q_1}\right)^2\right\}=
\frac{(T-t)^{3}}{6}-\sum_{j_3,j_2,j_1=0}^{q_1}
C_{j_3j_2j_1}^2-
$$
\begin{equation}
\label{sad006}
-\sum_{j_3,j_2,j_1=0}^{q_1}
C_{j_2j_3j_1}C_{j_3j_2j_1}\ (i_1\ne i_2=i_3),
\end{equation}

\vspace{4mm}

$$
{\sf M}\left\{\left(
J_{(111)T,t}^{(i_1 i_2 i_3)}-
J_{(111)T,t}^{(i_1 i_2 i_3)q_1}\right)^2\right\}=
\frac{(T-t)^{3}}{6}
-\sum_{j_3,j_2,j_1=0}^{q_1}
C_{j_3j_2j_1}^2-
$$
\begin{equation}
\label{sad007}
-\sum_{j_3,j_2,j_1=0}^{q_1}
C_{j_3j_2j_1}C_{j_1j_2j_3}\ \ \ (i_1=i_3\ne i_2),
\end{equation}

\vspace{2mm}
$$
{\sf M}\left\{\left(
J_{(111)T,t}^{(i_1 i_2 i_3)}-
J_{(111)T,t}^{(i_1 i_2 i_3)q_1}\right)^2\right\}=
\frac{(T-t)^{3}}{6}-\sum_{j_3,j_2,j_1=0}^{q_1}
C_{j_3j_2j_1}^2-
$$
\begin{equation}
\label{sad008}
-\sum_{j_3,j_2,j_1=0}^{q_1}
C_{j_3j_1j_2}C_{j_3j_2j_1}\ \ \ (i_1=i_2\ne i_3),
\end{equation}

\begin{equation}
\label{sad009}
~~~~~~~~~~{\sf M}\left\{\left(
J_{(111)T,t}^{(i_1 i_2 i_3)}-
J_{(111)T,t}^{(i_1 i_2 i_3)q_1}\right)^2\right\}
\le
6\left(\frac{(T-t)^{3}}{6}-\sum_{j_3,j_2,j_1=0}^{q_1}
C_{j_3j_2j_1}^2\right),
\end{equation}

\vspace{1mm}
\noindent
where $i_1,i_2,i_3=1,\ldots,m$ in (\ref{sad009}).

Let us compare 
the efficiency of application of Legendre polynomials 
and tri\-go\-no\-met\-ric functions for the approximation
of iterated stochastic integrals
$J_{(11)T,t}^{(i_1 i_2)},$ $J_{(111)T,t}^{(i_1 i_2 i_3)}$.

Consider the following conditions $(i_1\ne i_2,\ 
i_1\ne i_3,\ i_2\ne i_3)$

\vspace{-3mm}
\begin{equation}
\label{z1eee}
\frac{(T-t)^2}{2}\Biggl(\frac{1}{2}-\sum_{i=1}^{q}
\frac{1}{4i^2-1}\Biggr)\le (T-t)^4,
\end{equation}

\begin{equation}
\label{z2ququ}
(T-t)^{3}\Biggl(\frac{1}{6}-\sum_{j_1,j_2,j_3=0}^{q_1}
\frac{\left(C_{j_3j_2j_1}\right)^2}{(T-t)^3}\Biggr)\le (T-t)^4,
\end{equation}

\vspace{3mm}
\begin{equation}
\label{z3}
\frac{(T-t)^{2}}{2\pi^2}\Biggl(\frac{\pi^2}{6}-
\sum_{r=1}^{p}\frac{1}{r^2}\Biggr)\le (T-t)^4,
\end{equation}

\vspace{-3mm}
\begin{equation}
\label{z4}
(T-t)^3\Biggl(\frac{4}{45}-\frac{1}{4\pi^2}\sum_{r=1}^{p_1}\frac{1}{r^2}
\Biggl.
\Biggl.-\frac{55}{32\pi^4}\sum_{r=1}^{p_1}\frac{1}{r^4}-
\frac{1}{4\pi^4}\sum_{\stackrel{r,l=1}{{}_{r\ne l}}}^{p_1}
\frac{5l^4+4r^4-3r^2l^2}{r^2 l^2 \left(r^2-l^2\right)^2}\Biggr)\le (T-t)^4,
\end{equation}

\vspace{2mm}
\noindent
where
$$
C_{j_3j_2j_1}=\frac{\sqrt{(2j_1+1)(2j_2+1)(2j_3+1)}}{8}(T-t)^{3/2}\bar
C_{j_3j_2j_1},
$$
$$
\bar C_{j_3j_2j_1}=\int\limits_{-1}^{1}P_{j_3}(z)
\int\limits_{-1}^{z}P_{j_2}(y)
\int\limits_{-1}^{y}
P_{j_1}(x)dx dy dz,
$$
where $P_i(x)$ is the Legendre polynomial.

In Tables 5.44 and 5.45 we can see the minimal numbers 
$q,$ $q_1,$ $p,$ $p_1,$ which satisfy the conditions
(\ref{z1eee})--(\ref{z4}). As we mentioned above, the
numbers $q,$ $q_1$ are different. At that $q_1\ll q$
(the case of Legendre polynomials). As we saw in the previous
sections, we cannot take different numbers 
$p,$ $p_1$ for the case of trigonometric functions. Thus, we should
choose $q=p$ in (\ref{x2})--(\ref{x5}). This leads
to huge computational costs (see the fairly complicated formula (\ref{x5})).

From the other hand, we can take different numbers $q$
in (\ref{x2})--(\ref{x5}). At that we should exclude
random variables $\xi_q^{(i)},$ $\mu_q^{(i)}$ from 
(\ref{x2})--(\ref{x5}). 
At this situation for the case $i_1\ne i_2,$ $i_2\ne i_3,$ 
$i_1\ne i_3$ we have

\vspace{-3mm}
\begin{equation}
\label{zzz3xx}
\frac{3(T-t)^{2}}{2\pi^2}\Biggl(\frac{\pi^2}{6}-
\sum_{r=1}^{p^{*}}\frac{1}{r^2}\Biggr)\le (T-t)^4,
\end{equation}

\vspace{2mm}
$$
(T-t)^3\Biggl(\frac{5}{36}-\frac{1}{2\pi^2}\sum_{r=1}^{p_1^{*}}\frac{1}{r^2}
-\frac{79}{32\pi^4}\sum_{r=1}^{p_1^{*}}\frac{1}{r^4}-\Biggr.
$$
\begin{equation}
\label{zzz4}
\Biggl.-
\frac{1}{4\pi^4}\sum_{\stackrel{r,l=1}{{}_{r\ne l}}}^{p_1^{*}}
\frac{5l^4+4r^4-3r^2l^2}{r^2 l^2 \left(r^2-l^2\right)^2}\Biggr)\le 
(T-t)^4,
\end{equation}

\vspace{1mm}
\noindent
where the left-hand sides of (\ref{zzz3xx}), (\ref{zzz4}) correspond
to (\ref{x2}), (\ref{x5}) but without $\xi_q^{(i)},$ $\mu_q^{(i)}$.
In Table 5.45 we can see minimal numbers 
$p^{*},$ $p_1^{*}$, which satisfy the conditions
(\ref{zzz3xx}), (\ref{zzz4}). 

\begin{table}
\centering
\caption{Numbers $q,$ $q_1$}
\label{tab:5.44}      
\begin{tabular}{p{1.1cm}p{2.1cm}p{2.1cm}p{2.1cm}p{2.1cm}p{2.1cm}p{2.1cm}}
\hline\noalign{\smallskip}
$T-t$&$0.08222$&$0.05020$&$0.02310$&$0.01956$\\
\noalign{\smallskip}\hline\noalign{\smallskip}
$q$&19&51&235&328\\
\noalign{\smallskip}\hline\noalign{\smallskip}
$q_1$&1&2&5&6\\
\noalign{\smallskip}\hline\noalign{\smallskip}
\end{tabular}
\end{table}

\begin{table}
\centering
\caption{Numbers $p,$ $p_1,$ $p^{*},$ $p_1^{*}$}
\label{tab:5.45}      
\begin{tabular}{p{1.1cm}p{2.1cm}p{2.1cm}p{2.1cm}p{2.1cm}p{2.1cm}p{2.1cm}}
\hline\noalign{\smallskip}
$T-t$&$0.08222$&$0.05020$&$0.02310$&$0.01956$\\
\noalign{\smallskip}\hline\noalign{\smallskip}
$p$&8&21&96&133\\
\noalign{\smallskip}\hline\noalign{\smallskip}
$p_1$&1&1&3&4\\
\noalign{\smallskip}\hline\noalign{\smallskip}
$p^{*}$&23&61&286&398\\
\noalign{\smallskip}\hline\noalign{\smallskip}
$p_1^{*}$&1&2&4&5\\
\noalign{\smallskip}\hline\noalign{\smallskip}
\end{tabular}
\end{table}

\begin{table}
\centering
\caption{Confirmation of the formula (\ref{zzz4})}
\label{tab:5.46}      
\begin{tabular}{p{2.1cm}p{1.7cm}p{1.7cm}p{2.1cm}p{2.3cm}p{2.3cm}p{2.3cm}}
\hline\noalign{\smallskip}
$\varepsilon/(T-t)^3$&$0.0629$&$0.0097$&$0.0010$&$1.0129\cdot 10^{-4}$&
$1.0132\cdot 10^{-5}$\\
\noalign{\smallskip}\hline\noalign{\smallskip}
$q$&1&10&100&1000&10000\\
\noalign{\smallskip}\hline\noalign{\smallskip}
\end{tabular}
\end{table}

Moreover,
$$
{\sf M}\left\{\left(J_{(01)T,t}^{(0 i_1)}-
J_{(01)T,t}^{(0 i_1)q}
\right)^2\right\}=
{\sf M}\left\{\left(J_{(10)T,t}^{(i_1 0)}-
J_{(10)T,t}^{(i_1 0)q}
\right)^2\right\}=
$$
\begin{equation}
\label{1001x}
=\frac{(T-t)^{3}}{2\pi^2}\Biggl(\frac{\pi^2}{6}-
\sum_{r=1}^q \frac{1}{r^2}\Biggr)\ne 0,
\end{equation}

\vspace{2mm}
\noindent
where
$J_{(01)T,t}^{(0 i_1)q}$, $J_{(10)T,t}^{(i_1 0)q}$ 
are defined by (\ref{x3}), (\ref{x4}) 
but without $\xi_q^{(i)}.$

It is not difficult to see that the numbers $q_{\rm trig}$ in Table 5.42
correspond to minimal numbers $q_{\rm trig}$, which satisfy 
the condition (compare with (\ref{1001x}))
$$
\frac{(T-t)^{3}}{2\pi^2}\Biggl(\frac{\pi^2}{6}-
\sum_{r=1}^{q_{\rm trig}} \frac{1}{r^2}\Biggr)\le (T-t)^4.
$$

From the other hand, the right-hand sides of (\ref{y3a}), (\ref{y4a}) 
include only 2 random variables.
In this situation we again can talk about
the advantage of Ledendre polynomials.

In Table 5.46 we can see the numerical confirmation of 
the formula (\ref{zzz4}) ($\varepsilon$ means the left-hand side 
of (\ref{zzz4})).

\subsection{A Comparative Analysis of Efficiency of Using the 
Legendre Polynomials and Trigonometric Functions for the Integral
$J_{(011)T,t}^{*(0 i_1 i_2)}$}

In this section, we compare computational costs for 
approximation of the iterated Stratonovich
stochastic integral $J_{(011)T,t}^{*(0 i_1 i_2)}$ 
$(i_1, i_2=1,\ldots,m)$ within the framework
of the method of generalized multiple Fourier series for 
the Legendre polynomial system
and the system of trigomomenric functions.

Using Theorem 2.1 for the case of trigonometric system of
functions, we obtain \cite{6}-\cite{12aa}, \cite{arxiv-12}
$$
J_{(011)T,t}^{*(0 i_1 i_2)q}=(T-t)^{2}\Biggl(\frac{1}{6}
\zeta_{0}^{(i_1)}\zeta_{0}^{(i_2)}-\frac{1}{2\sqrt{2}\pi}
\sqrt{\alpha_q}\xi_q^{(i_2)}\zeta_0^{(i_1)}+\Biggr.
$$
$$
+\frac{1}{2\sqrt{2}\pi^2}\sqrt{\beta_q}\Biggl(
\mu_q^{(i_2)}\zeta_0^{(i_1)}-2\mu_q^{(i_1)}\zeta_0^{(i_2)}\Biggr)+
$$
$$
+\frac{1}{2\sqrt{2}}\sum_{r=1}^{q}
\Biggl(-\frac{1}{\pi r}
\zeta_{2r-1}^{(i_2)}
\zeta_{0}^{(i_1)}+
\frac{1}{\pi^2 r^2}\left(
\zeta_{2r}^{(i_2)}
\zeta_{0}^{(i_1)}-
2\zeta_{2r}^{(i_1)}
\zeta_{0}^{(i_2)}\right)\Biggr)-
$$
$$
-
\frac{1}{2\pi^2}\sum_{\stackrel{r,l=1}{{}_{r\ne l}}}^{q}
\frac{1}{r^2-l^2}\Biggl(
\zeta_{2r}^{(i_1)}
\zeta_{2l}^{(i_2)}+
\frac{l}{r}
\zeta_{2r-1}^{(i_1)}
\zeta_{2l-1}^{(i_2)}
\Biggr)+
$$
$$
+
\sum_{r=1}^{q}
\Biggl(\frac{1}{4\pi r}\left(
\zeta_{2r}^{(i_1)}
\zeta_{2r-1}^{(i_2)}-
\zeta_{2r-1}^{(i_1)}
\zeta_{2r}^{(i_2)}\right)+\Biggr.
$$
\begin{equation}
\label{t1}
\Biggl.\Biggl.+
\frac{1}{8\pi^2 r^2}\left(
3\zeta_{2r-1}^{(i_1)}
\zeta_{2r-1}^{(i_2)}+
\zeta_{2r}^{(i_2)}
\zeta_{2r}^{(i_1)}\right)\Biggr)\Biggr).
\end{equation}

\vspace{3mm}

\begin{table}
\centering
\caption{Confirmation of the formula (\ref{t2})}
\label{tab:5.47}      
\begin{tabular}{p{2.1cm}p{1.7cm}p{1.7cm}p{2.1cm}p{2.3cm}p{2.3cm}p{2.3cm}}
\hline\noalign{\smallskip}
$4\varepsilon/(T-t)^4$&0.0540&0.0082&$8.4261\cdot 10^{-4}$
&$8.4429\cdot 10^{-5}$&
$8.4435\cdot 10^{-6}$\\
\noalign{\smallskip}\hline\noalign{\smallskip}
$q$&1&10&100&1000&10000\\
\noalign{\smallskip}\hline\noalign{\smallskip}
\end{tabular}
\end{table}

\begin{table}
\centering
\caption{Confirmation of the formula (\ref{t4})}
\label{tab:5.48}      
\begin{tabular}{p{2.1cm}p{1.7cm}p{1.7cm}p{2.1cm}p{2.3cm}p{2.3cm}p{2.3cm}}
\hline\noalign{\smallskip}
$16\varepsilon/(T-t)^4$&0.3797&0.0581&0.0062&$6.2450\cdot 10^{-4}$&$6.2495\cdot 10^{-5}$\\
\noalign{\smallskip}\hline\noalign{\smallskip}
$q$&1&10&100&1000&10000\\
\noalign{\smallskip}\hline\noalign{\smallskip}
\end{tabular}
\end{table}

For the case $i_1\ne i_2$ from Theorem 1.3 we get 
\cite{6}-\cite{art-1}, \cite{arxiv-3}, \cite{arxiv-12}

\vspace{-3mm}
$$
{\sf M}\left\{\left(J_{(011)T,t}^{*(0 i_1 i_2)}-
J_{(011)T,t}^{*(0 i_1 i_2)q}\right)^2\right\}=
\frac{(T-t)^4}{4}\Biggl(\frac{1}{9}-
\frac{1}{2\pi^2}\sum_{r=1}^q \frac{1}{r^2}-\Biggr.
$$
\begin{equation}
\label{t2}
\Biggl.-\frac{5}{8\pi^4}\sum_{r=1}^q \frac{1}{r^4}-
\frac{1}{\pi^4}\sum_{\stackrel{k,l=1}{{}_{k\ne l}}}^q
\frac{k^2+l^2}{l^2\left(l^2-k^2\right)^2}\Biggr).
\end{equation}

\vspace{3mm}

Analogues of the formulas (\ref{t1}), (\ref{t2}) for the case of 
Legendre polynomials will look as follows 
\cite{6}-\cite{art-1}, \cite{arxiv-3}, \cite{arxiv-12}

\vspace{-1mm}
$$
J_{(011)T,t}^{*(0 i_1 i_2)q}=\frac{T-t}{2}J_{(11)T,t}^{*(i_1 i_2)q}
+\frac{(T-t)^2}{4}\Biggl(
\frac{1}{\sqrt{3}}\zeta_0^{(i_2)}\zeta_1^{(i_1)}+\Biggr.
$$

\vspace{-2mm}
\begin{equation}
\label{t3}
~~~~~~~~ +\Biggl.\sum_{i=0}^{q}\Biggl(
\frac{(i+1)\zeta_{i+2}^{(i_2)}\zeta_{i}^{(i_1)}
-(i+2)\zeta_{i}^{(i_2)}\zeta_{i+2}^{(i_1)}}
{\sqrt{(2i+1)(2i+5)}(2i+3)}+
\frac{\zeta_i^{(i_1)}\zeta_{i}^{(i_2)}}{(2i-1)(2i+3)}\Biggr)\Biggr),
\end{equation}

\vspace{2mm}
\noindent
where
$$
J_{(11)T,t}^{*(i_1 i_2)q}=
\frac{T-t}{2}\biggl(\zeta_0^{(i_1)}\zeta_0^{(i_2)}+\sum_{i=1}^{q}
\frac{1}{\sqrt{4i^2-1}}\left(
\zeta_{i-1}^{(i_1)}\zeta_{i}^{(i_2)}-
\zeta_i^{(i_1)}\zeta_{i-1}^{(i_2)}\right)\biggr),
$$
$$
{\sf M}\left\{\left(J_{(011)T,t}^{*(0 i_1 i_2)}-J_{(011)T,t}^{*(0 i_1 i_2)q}
\right)^2\right\}=
\frac{(T-t)^4}{16}\left(\frac{5}{9}-
2\sum_{i=2}^q\frac{1}{4i^2-1}-\right.
$$
\begin{equation}
\label{t4}
~~~~~~~~ -
\sum_{i=1}^q
\frac{1}{(2i-1)^2(2i+3)^2}
-\left.\sum_{i=0}^q\frac{(i+2)^2+(i+1)^2}{(2i+1)(2i+5)(2i+3)^2}
\right),
\end{equation}

\vspace{3mm}
\noindent
where $i_1\ne i_2.$

In Tables 5.47 and 5.48 we can see the numerical confirmation of 
the formulas (\ref{t2}) and (\ref{t4}) ($\varepsilon$ means
the right-hand side of (\ref{t2}) or
(\ref{t4})).

Let us compare the complexity of the formulas (\ref{t1}) and (\ref{t3}).
The formula (\ref{t1}) includes the double sum

\vspace{-2mm}
$$
\frac{1}{2\pi^2}\sum_{\stackrel{r,l=1}{{}_{r\ne l}}}^{q}
\frac{1}{r^2-l^2}\Biggl(
\zeta_{2r}^{(i_1)}
\zeta_{2l}^{(i_2)}+
\frac{l}{r}
\zeta_{2r-1}^{(i_1)}
\zeta_{2l-1}^{(i_2)}
\Biggr).
$$

\vspace{2mm}

Thus, the formula (\ref{t1}) is more complex, than the formula (\ref{t3})
even if we take identical numbers $q$ in these formulas.
As we noted above, the number $q$ in (\ref{t1}) must be equal to the
number $q$ from the formula (\ref{x2}), so it is much 
larger than the number $q$ 
from the formula (\ref{t3}). As a result, we have 
obvious advantage of the formula (\ref{t3})
in computational costs. 

As we mentioned above,
if we will not introduce the random 
variables $\xi_q^{(i)}$ and 
$\mu_q^{(i)},$ then the number $q$ in (\ref{t1}) can be chosen smaller, but
the mean-square error of approximation of the
stochastic integral $J_{(11)T,t}^{(i_1 i_2)}$ will be three times larger
(see (\ref{801xxxx})). Moreover, in this case the stochastic integrals 
$J_{(01)T,t}^{(0 i_1)}$, $J_{(10)T,t}^{(i_1 0)},$ 
$J_{(001)T,t}^{(00i_1)}$ (with Gaussian distribution)
will be approximated worse. In this situation, we can again talk about
the advantage of Ledendre polynomials.

\subsection{Conclusions}

Summing up the results of previous sections we 
can come to the following conclusions.

1.\;We can talk about approximately equal computational costs
for the formulas (\ref{x2}) and (\ref{y3}). This means that
computational costs for realizing the Milstein scheme (explicit
one-step strong numerical method with the convergence order $\gamma=1.0$
for It\^{o} SDEs; see Sect.~4.10) 
for the case of Legendre polynomials and for the case 
of trigonometric functions are approximately the same.

2.\;If we will not introduce the random 
variables $\xi_q^{(i)}$ (see (\ref{x2})), then 
the mean-square error of approximation of the
stochastic integral $J_{(11)T,t}^{(i_1 i_2)}$ will be three times larger
(see (\ref{801xxxx})). 
In this situation, we can talk about
the advantage of Ledendre polynomials in
the Milstein method.
Moreover, in this case the stochastic integrals 
$J_{(01)T,t}^{(0 i_1)}$, $J_{(10)T,t}^{(i_1 0)},$ 
$J_{(001)T,t}^{(00i_1)}$ (with Gaussian distribution)
will be approximated worse.

3.\;If we talk about the explicit one-step strong numerical scheme 
with the convergence order $\gamma=1.5$
for It\^{o} SDEs (see Sect.~4.10), then 
the numbers $q,$ $q_1$ (see (\ref{y3}), (\ref{y4})) 
are different. At that $q_1\ll q$
(the case of Legendre polynomials). 
The number $q$ must be the same in (\ref{x2})--(\ref{x5}) 
(the case of trigonometric functions).
This leads to huge computational costs (see the fairly
complicated formula (\ref{x5})).
From the other hand, we can take different numbers $q$
in (\ref{x2})--(\ref{x5}). At that we should exclude
the random variables $\xi_q^{(i)},$ $\mu_q^{(i)}$ from 
(\ref{x2})--(\ref{x5}). This leads to another 
problems which we discussed above (see Conclusion 2).

4.\;In addition, 
the author supposes that effect described in Conclusion 
3 will be more impressive when 
analyzing more complex sets of iterated It\^{o} and Stratonovich 
stochastic
integrals (when $\gamma=$ $2.0,$ $2.5,$ $3.0,$ $\ldots $). 
This supposition is based on the fact that the polynomial 
system of functions has the significant advantage (in comparison with 
the trigonometric system) when approximating the iterated stochastic 
integrals for which not all weight functions are equal to 1
(see Sect 5.4.3 and conclusion at the end of Sect.~5.1).

\section{Optimization of the Mean-Square Approximation Procedures for 
Iterated It\^{o} Stochastic Integrals Based on Theorem 1.1 and Multiple Fourier--Le\-gend\-re Series}

This section is devoted to optimization of the mean-square approximation 
procedures for iterated It\^{o} stochastic integrals (\ref{k1000xxxx})
of multiplicities 1 to 4 based on 
Theorem 1.1 and multiple Fourier--Legendre series. 
The mentioned stochastic integrals are part of 
strong numerical methods with convergence orders 1.0, 1.5, and 2.0 for 
It\^{o} SDEs with multidimensional 
non-commutative noise (see (\ref{al1})--(\ref{al3})). 
We show that the lengths of sequences of independent standard 
Gaussian random variables required for the mean-square approximation of 
iterated It\^{o} stochastic integrals (\ref{k1000xxxx}) 
can be significantly reduced without the loss of 
the mean-square accuracy of approximation for these stochastic integrals.
This section is written on the base of paper \cite{Mikh-2aaaaa}.
An extension of the mentioned results to iterated It\^{o}
stochastic integrals of multiplicity 5 can be found in \cite{Mikh-2}.

Using Theorem 1.1 and
the system of Legendre polynomials, we obtain
the following approximations of iterated It\^{o}
stochastic integrals (\ref{k1000xxxx})
$$
I_{(0)T,t}^{(i_1)}=\sqrt{T-t}\zeta_0^{(i_1)},
$$
$$
I_{(1)T,t}^{(i_1)}=-\frac{(T-t)^{3/2}}{2}\left(\zeta_0^{(i_1)}+
\frac{1}{\sqrt{3}}\zeta_1^{(i_1)}\right),
$$
\begin{equation}
\label{nach1}
I_{(00)T,t}^{(i_1 i_2)q}=
\frac{T-t}{2}\left(\zeta_0^{(i_1)}\zeta_0^{(i_2)}+\sum_{i=1}^{q}
\frac{1}{\sqrt{4i^2-1}}\left(
\zeta_{i-1}^{(i_1)}\zeta_{i}^{(i_2)}-
\zeta_i^{(i_1)}\zeta_{i-1}^{(i_2)}\right) - {\bf 1}_{\{i_1=i_2\}}\right),
\end{equation}

\vspace{-3mm}
$$
I_{(000)T,t}^{(i_1i_2i_3)q_1}
=
\sum_{j_1,j_2,j_3=0}^{q_1}
C_{j_3j_2j_1}
\Biggl(
\zeta_{j_1}^{(i_1)}\zeta_{j_2}^{(i_2)}\zeta_{j_3}^{(i_3)}
-{\bf 1}_{\{i_1=i_2\}}
{\bf 1}_{\{j_1=j_2\}}
\zeta_{j_3}^{(i_3)}-
\Biggr.
$$
\begin{equation}
\label{901}
\Biggl.
-{\bf 1}_{\{i_2=i_3\}}
{\bf 1}_{\{j_2=j_3\}}
\zeta_{j_1}^{(i_1)}-
{\bf 1}_{\{i_1=i_3\}}
{\bf 1}_{\{j_1=j_3\}}
\zeta_{j_2}^{(i_2)}\Biggr),
\end{equation}

\vspace{3mm}
\begin{equation}
\label{902}
I_{(10)T,t}^{(i_1 i_2)q_2}=
\sum_{j_1,j_2=0}^{q_2}
C_{j_2j_1}^{10}\Biggl(\zeta_{j_1}^{(i_1)}\zeta_{j_2}^{(i_2)}
-{\bf 1}_{\{i_1=i_2\}}
{\bf 1}_{\{j_1=j_2\}}\Biggr),
\end{equation}
\begin{equation}
\label{903}
I_{(01)T,t}^{(i_1 i_2)\bar q_2}=
\sum_{j_1,j_2=0}^{\bar q_2}
C_{j_2j_1}^{01}\Biggl(\zeta_{j_1}^{(i_1)}\zeta_{j_2}^{(i_2)}
-{\bf 1}_{\{i_1=i_2\}}
{\bf 1}_{\{j_1=j_2\}}\Biggr),
\end{equation}

\vspace{3mm}
$$
I_{(0000)T,t}^{(i_1 i_2 i_3 i_4)q_3}
=
\sum_{j_1,j_2,j_3,j_4=0}^{q_3}
C_{j_4 j_3 j_2 j_1}\Biggl(
\zeta_{j_1}^{(i_1)}\zeta_{j_2}^{(i_2)}\zeta_{j_3}^{(i_3)}\zeta_{j_4}^{(i_4)}
-\Biggr.
$$
$$
-
{\bf 1}_{\{i_1=i_2\}}
{\bf 1}_{\{j_1=j_2\}}
\zeta_{j_3}^{(i_3)}
\zeta_{j_4}^{(i_4)}
-
{\bf 1}_{\{i_1=i_3\}}
{\bf 1}_{\{j_1=j_3\}}
\zeta_{j_2}^{(i_2)}
\zeta_{j_4}^{(i_4)}-
$$
$$
-
{\bf 1}_{\{i_1=i_4\}}
{\bf 1}_{\{j_1=j_4\}}
\zeta_{j_2}^{(i_2)}
\zeta_{j_3}^{(i_3)}
-
{\bf 1}_{\{i_2=i_3\}}
{\bf 1}_{\{j_2=j_3\}}
\zeta_{j_1}^{(i_1)}
\zeta_{j_4}^{(i_4)}-
$$
$$
-
{\bf 1}_{\{i_2=i_4\}}
{\bf 1}_{\{j_2=j_4\}}
\zeta_{j_1}^{(i_1)}
\zeta_{j_3}^{(i_3)}
-
{\bf 1}_{\{i_3=i_4\}}
{\bf 1}_{\{j_3=j_4\}}
\zeta_{j_1}^{(i_1)}
\zeta_{j_2}^{(i_2)}+
$$
$$
+
{\bf 1}_{\{i_1=i_2\}}
{\bf 1}_{\{j_1=j_2\}}
{\bf 1}_{\{i_3=i_4\}}
{\bf 1}_{\{j_3=j_4\}}+
{\bf 1}_{\{i_1=i_3\}}
{\bf 1}_{\{j_1=j_3\}}
{\bf 1}_{\{i_2=i_4\}}
{\bf 1}_{\{j_2=j_4\}}+
$$
\begin{equation}
\label{904}
+\Biggl.
{\bf 1}_{\{i_1=i_4\}}
{\bf 1}_{\{j_1=j_4\}}
{\bf 1}_{\{i_2=i_3\}}
{\bf 1}_{\{j_2=j_3\}}\Biggr),
\end{equation}

\vspace{2mm}
\noindent
where ${\bf 1}_A$ is the indicator of the set $A$,
$$
\zeta_{j}^{(i)}=
\int\limits_t^T \phi_{j}(s) d{\bf w}_s^{(i)}\ \ \ (i=1,\ldots,m,\
j=0,1,\ldots)
$$
are independent standard Gaussian random variables
for various
$i$ or $j$, $\{\phi_j(x)\}_{j=0}^{\infty}$ is a complete
orthonormal system of Legendre polynomials  
in the space $L_2([t, T])$  (see (\ref{4009d})),
\begin{equation}
\label{2001}
~~~~~~C_{j_3 j_2 j_1}=\frac{1}{8}L_{j_1j_2j_3}(T-t)^{3/2}\bar C_{j_3j_2j_1},\ \ \
C_{j_2j_1}^{01}=\frac{1}{8}L_{j_1j_2}(T-t)^{2}\bar C_{j_2j_1}^{01},
\end{equation}
\begin{equation}
\label{2002}
~~~~~C_{j_2j_1}^{10}=\frac{1}{8}L_{j_1j_2}(T-t)^{2}\bar C_{j_2j_1}^{10},\ \ \
C_{j_4j_3j_2j_1}=\frac{1}{16}L_{j_1 j_2 j_3 j_4}(T-t)^{2}\bar C_{j_4j_3j_2j_1},
\end{equation}

\vspace{-2mm}
$$
L_{j_1 j_2}=\sqrt{(2j_1+1)(2j_2+1)},\ \ \ 
L_{j_1 j_2 j_3}=\sqrt{(2j_1+1)(2j_2+1)(2j_3+1)}, 
$$

\vspace{-2mm}
$$
L_{j_1 j_2 j_3 j_4}=\sqrt{(2j_1+1)(2j_2+1)(2j_3+1)(2j_4+1)},
$$

\vspace{-2mm}
$$
\bar C_{j_3j_2j_1}=
\int\limits_{-1}^{1}P_{j_3}(z)
\int\limits_{-1}^{z}P_{j_2}(y)
\int\limits_{-1}^{y}
P_{j_1}(x)dx dy dz,
$$
$$
\bar C_{j_4j_3j_2j_1}=\int\limits_{-1}^{1}P_{j_4}(u)
\int\limits_{-1}^{u}P_{j_3}(z)
\int\limits_{-1}^{z}P_{j_2}(y)
\int\limits_{-1}^{y}
P_{j_1}(x)dx dy dz du,
$$
$$
\bar C_{j_2j_1}^{01}=-
\int\limits_{-1}^{1}(1+y)P_{j_2}(y)
\int\limits_{-1}^{y}
P_{j_1}(x)dx dy,
$$
$$
\bar C_{j_2j_1}^{10}=-
\int\limits_{-1}^{1}P_{j_2}(y)
\int\limits_{-1}^{y}
(1+x)P_{j_1}(x)dx dy,
$$
          
\noindent
$P_j(x)$ is the Legendre polynomial (see (\ref{agentww1})).

Combining the estimates (\ref{agentww3}) 
and (\ref{z1}) for $p_1=\ldots=p_k=p$, we obtain
\begin{equation}
\label{agentyz1}
~~~~ k!\left(\int\limits_{[t,T]^k}
K^2(t_1,\ldots,t_k)
dt_1\ldots dt_k -\sum_{j_1,\ldots,j_k=0}^{p}C^2_{j_k\ldots j_1}\right)\le 
C(T-t)^{r+1},
\end{equation}
 
\noindent
where $K(t_1,\ldots,t_k)$ is defined by (\ref{leto7000}) (see (\ref{aaa1})--(\ref{aaa3})),
$r/2$ is the strong convergence orders for the numerical schemes
(\ref{al1})--(\ref{al3}), i.e. $r/2=1.0, 1.5,$ and $2.0;$
constant $C$ is independent of $T-t$.

It is not difficult to see that the multiplier factor $k!$ on the left-hand 
side of the inequality (\ref{agentyz1}) 
leads to a significant increase of computational costs for approximation of iterated It\^{o} 
stochastic integrals. The mentioned problem can be overcome if we calculate 
the mean-square approximation error for iterated It\^{o} 
stochastic integrals exactly (see Theorem 1.3 and Sect.~1.2.3).
In this section, we discuss how to essentially minimize the numbers
$q,$ $q_1,$ $q_2,$ $\bar q_2,$ $q_3$
from (\ref{nach1})--(\ref{904}). At that we will use the results from Sect.~1.2.3.

Denote
\begin{equation}
\label{agentys1000}
E^{(l_1\ldots l_k)}_p\stackrel{\sf def}{=}
{\sf M}\left\{\left(
I_{(l_1\ldots l_k)T,t}^{(i_1\ldots i_k)}-
I_{(l_1\ldots l_k)T,t}^{(i_1\ldots i_k)p}\right)^2\right\},
\end{equation}

\noindent
where $I_{(l_1\ldots l_k)T,t}^{(i_1\ldots i_k)}$ is the 
iterated It\^{o} stochastic integral (\ref{k1000xxxx})
and $I_{(l_1\ldots l_k)T,t}^{(i_1\ldots i_k)p}$ is the
mean-square approximation of this stochastic integral. More precisely, 
the approximations 
$I_{(00)T,t}^{(i_1i_2)q},$ 
$I_{(000)T,t}^{(i_1i_2i_3)q_1},$
$I_{(10)T,t}^{(i_1i_2)q_2},$ $I_{(01)T,t}^{(i_1i_2)\bar q_2},$
$I_{(0000)T,t}^{(i_1i_2i_3i_4)q_3}$ are defined by
(\ref{nach1})-(\ref{904}).

The results of Sect.~1.2.3 give the following formulas for the case of Legendre polynomials
\begin{equation}
\label{prod1}
E^{(00)}_q
=
\frac{(T-t)^2}{2}\left(\frac{1}{2}-\sum_{i=1}^q
\frac{1}{4i^2-1}\right),\ \ i_1\ne i_2,
\end{equation}
\begin{equation}
\label{prod2}
~~~~~~~ E^{(000)}_{q_{1,1}} = (T-t)^3\left(\frac{1}{6}-\frac{1}{64}
\sum_{j_1,j_2,j_3=0}^{q_{1,1}}L_{j_1j_2j_3}^2
\left(\bar C_{j_3j_2j_1}\right)^2\right),
\end{equation}
where $i_1\ne i_2,\ i_1\ne i_3,\ i_2\ne i_3,$
\begin{equation}
\label{prod3}
E^{(000)}_{q_{1,2}} = (T-t)^3\left(\frac{1}{6}-\frac{1}{64}
\sum_{j_1,j_2,j_3=0}^{q_{1,2}}L_{j_1j_2j_3}^2
\left(\left(\bar C_{j_3j_2j_1}\right)^2+
\bar C_{j_3j_1j_2}\bar C_{j_3j_2j_1}\right)\right),
\end{equation}
where $i_1=i_2\ne i_3,$
\begin{equation}
\label{prod4}
E^{(000)}_{q_{1,3}} = (T-t)^3\left(\frac{1}{6}-\frac{1}{64}
\sum_{j_1,j_2,j_3=0}^{q_{1,3}}L_{j_1j_2j_3}^2
\left(\left(\bar C_{j_3j_2j_1}\right)^2+
\bar C_{j_2j_3j_1}\bar C_{j_3j_2j_1}\right)\right),
\end{equation}
where $i_1\ne i_2=i_3,$
\begin{equation}
\label{prod5}
E^{(000)}_{q_{1,4}} = (T-t)^3\left(\frac{1}{6}-\frac{1}{64}
\sum_{j_1,j_2,j_3=0}^{q_{1,4}}L_{j_1j_2j_3}^2
\left(\left(\bar C_{j_3j_2j_1}\right)^2+
\bar C_{j_3j_2j_1}\bar C_{j_1j_2j_3}\right)\right),
\end{equation}
where $i_1=i_3\ne i_2,$
\begin{equation}
\label{prod6}
E^{(10)}_{q_{2,1}}
=(T-t)^4\left(\frac{1}{12}-\frac{1}{64}
\sum_{j_1,j_2=0}^{q_{2,1}} L_{j_1j_2}^2
\left(\bar C_{j_2j_1}^{10}\right)^2
\right),\ \ i_1\ne i_2,
\end{equation}
\begin{equation}
\label{prod7}
E^{(10)}_{q_{2,2}}
=(T-t)^4\left(\frac{1}{12}-\frac{1}{64}
\sum_{j_1,j_2=0}^{q_{2,2}} L_{j_1j_2}^2
\bar C_{j_2j_1}^{10}\left(\sum\limits_{(j_1,j_2)}\bar C_{j_2j_1}^{10}
\right)\right),\ \ i_1=i_2,
\end{equation}

\vspace{-4mm}
\begin{equation}
\label{prod8}
E^{(01)}_{\bar q_{2,1}}
=(T-t)^4\left(\frac{1}{4}-\frac{1}{64}
\sum_{j_1,j_2=0}^{\bar q_{2,1}} L_{j_1j_2}^2
\left(\bar C_{j_2j_1}^{01}\right)^2
\right),\ \ i_1\ne i_2,
\end{equation}
\begin{equation}
\label{prod9}
E^{(01)}_{\bar q_{2,2}}
=(T-t)^4\left(\frac{1}{4}-\frac{1}{64}
\sum_{j_1,j_2=0}^{\bar q_{2,2}} L_{j_1j_2}^2
\bar C_{j_2j_1}^{01}\left(\sum\limits_{(j_1,j_2)}\bar C_{j_2j_1}^{01}
\right)\right),\ \ i_1=i_2,
\end{equation}

\vspace{-4mm}
\begin{equation}
\label{prod10}
E^{(0000)}_{q_{3,1}} = (T-t)^4\left(\frac{1}{24}-\frac{1}{256}
\sum_{j_1,\ldots,j_4=0}^{q_{3,1}}L_{j_1\ldots j_4}^2
\left(\bar C_{j_4\ldots j_1}\right)^2\right),
\end{equation}
where $i_1,\ldots,i_4$ are pairwise different,
\begin{equation}
\label{prod11}
E^{(0000)}_{q_{3,2}} = (T-t)^4\left(\frac{1}{24}-\frac{1}{256}
\sum_{j_1,\ldots,j_4=0}^{q_{3,2}}L_{j_1\ldots j_4}^2
\bar C_{j_4\ldots j_1}\left(
\sum\limits_{(j_1,j_2)}
\bar C_{j_4\ldots j_1}\right)\right),
\end{equation}
where $i_1=i_2\ne i_3, i_4;\ i_3\ne i_4,$
\begin{equation}
\label{prod12}
E^{(0000)}_{q_{3,3}} = (T-t)^4\left(\frac{1}{24}-\frac{1}{256}
\sum_{j_1,\ldots,j_4=0}^{q_{3,3}}L_{j_1\ldots j_4}^2
\bar C_{j_4\ldots j_1}\left(
\sum\limits_{(j_1,j_3)}
\bar C_{j_4\ldots j_1}\right)\right),
\end{equation}
where $i_1=i_3\ne i_2, i_4;\ i_2\ne i_4,$
\begin{equation}
\label{prod13}
E^{(0000)}_{q_{3,4}} = (T-t)^4\left(\frac{1}{24}-\frac{1}{256}
\sum_{j_1,\ldots,j_4=0}^{q_{3,4}}L_{j_1\ldots j_4}^2
\bar C_{j_4\ldots j_1}\left(
\sum\limits_{(j_1,j_4)}
\bar C_{j_4\ldots j_1}\right)\right),
\end{equation}
where $i_1=i_4\ne i_2, i_3;\ i_2\ne i_3,$
\begin{equation}
\label{prod14}
E^{(0000)}_{q_{3,5}}= (T-t)^4\left(\frac{1}{24}-\frac{1}{256}
\sum_{j_1,\ldots,j_4=0}^{q_{3,5}}L_{j_1\ldots j_4}^2
\bar C_{j_4\ldots j_1}\left(
\sum\limits_{(j_2,j_3)}
\bar C_{j_4\ldots j_1}\right)\right),
\end{equation}
where $i_2=i_3\ne i_1, i_4;\ i_1\ne i_4,$
\begin{equation}
\label{prod15}
E^{(0000)}_{q_{3,6}}= (T-t)^4\left(\frac{1}{24}-\frac{1}{256}
\sum_{j_1,\ldots,j_4=0}^{q_{3,6}}L_{j_1\ldots j_4}^2
\bar C_{j_4\ldots j_1}\left(
\sum\limits_{(j_2,j_4)}
\bar C_{j_4\ldots j_1}\right)\right),
\end{equation}
where $i_2=i_4\ne i_1, i_3;\ i_1\ne i_3,$
\begin{equation}
\label{prod16}
E^{(0000)}_{q_{3,7}} = (T-t)^4\left(\frac{1}{24}-\frac{1}{256}
\sum_{j_1,\ldots,j_4=0}^{q_{3,7}}L_{j_1\ldots j_4}^2
\bar C_{j_4\ldots j_1}\left(
\sum\limits_{(j_3,j_4)}
\bar C_{j_4\ldots j_1}\right)\right),
\end{equation}
where $i_3=i_4\ne i_1, i_2;\ i_1\ne i_2,$
\begin{equation}
\label{prod17}
E^{(0000)}_{q_{3,8}}= (T-t)^4\left(\frac{1}{24}-\frac{1}{256}
\sum_{j_1,\ldots,j_4=0}^{q_{3,8}}L_{j_1\ldots j_4}^2
\bar C_{j_4\ldots j_1}\left(
\sum\limits_{(j_1,j_2,j_3)}
\bar C_{j_4\ldots j_1}\right)\right), 
\end{equation}
where $i_1=i_2=i_3\ne i_4,$
\begin{equation}
\label{prod18}
E^{(0000)}_{q_{3,9}}= (T-t)^4\left(\frac{1}{24}-\frac{1}{256}
\sum_{j_1,\ldots,j_4=0}^{q_{3,9}}L_{j_1\ldots j_4}^2
\bar C_{j_4\ldots j_1}\left(
\sum\limits_{(j_2,j_3,j_4)}
\bar C_{j_4\ldots j_1}\right)\right),
\end{equation}
where $i_2=i_3=i_4\ne i_1$,
\begin{equation}
\label{prod19}
E^{(0000)}_{q_{3,10}}= (T-t)^4\left(\frac{1}{24}-\frac{1}{256}
\sum_{j_1,\ldots,j_4=0}^{q_{3,10}}L_{j_1\ldots j_4}^2
\bar C_{j_4\ldots j_1}\left(
\sum\limits_{(j_1,j_2,j_4)}
\bar C_{j_4\ldots j_1}\right)\right),
\end{equation}
where $i_1=i_2=i_4\ne i_3$,
\begin{equation}
\label{prod20}
E^{(0000)}_{q_{3,11}}= (T-t)^4\left(\frac{1}{24}-\frac{1}{256}
\sum_{j_1,\ldots,j_4=0}^{q_{3,11}}L_{j_1\ldots j_4}^2
\bar C_{j_4\ldots j_1}\left(
\sum\limits_{(j_1,j_3,j_4)}
\bar C_{j_4\ldots j_1}\right)\right),
\end{equation}
where $i_1=i_3=i_4\ne i_2$,
\begin{equation}
\label{prod21}
E^{(0000)}_{q_{3,12}}= (T-t)^4\left(\frac{1}{24}-\frac{1}{256}
\sum_{j_1,\ldots,j_4=0}^{q_{3,12}}L_{j_1\ldots j_4}^2
\bar C_{j_4\ldots j_1}\left(\sum\limits_{(j_1,j_2)}\left(
\sum\limits_{(j_3,j_4)}
\bar C_{j_4\ldots j_1}\right)\right)\right),
\end{equation}
where $i_1=i_2\ne i_3=i_4$,
\begin{equation}
\label{prod22}
E^{(0000)}_{q_{3,13}}= (T-t)^4\left(\frac{1}{24}-\frac{1}{256}
\sum_{j_1,\ldots,j_4=0}^{q_{3,13}}L_{j_1\ldots j_4}^2
\bar C_{j_4\ldots j_1}\left(\sum\limits_{(j_1,j_3)}\left(
\sum\limits_{(j_2,j_4)}
\bar C_{j_4\ldots j_1}\right)\right)\right),
\end{equation}
where $i_1=i_3\ne i_2=i_4$,
\begin{equation}
\label{prod23}
E^{(0000)}_{q_{3,14}}= (T-t)^4\left(\frac{1}{24}-\frac{1}{256}
\sum_{j_1,\ldots,j_4=0}^{q_{3,14}}L_{j_1\ldots j_4}^2
\bar C_{j_4\ldots j_1}\left(\sum\limits_{(j_1,j_4)}\left(
\sum\limits_{(j_2,j_3)}
\bar C_{j_4\ldots j_1}\right)\right)\right),
\end{equation}
where $i_1=i_4\ne i_2=i_3$.

Obviously, the conditions (\ref{prod1})--(\ref{prod23}) do not contain the multiplier factors 
$2!,$ $3!,$ and $4!$    
in contrast to the estimate (\ref{z1}) (see Theorem 1.4). 
However, the number of the mentioned conditions is quite 
large, which is inconvenient for practice. In this section, we propose the hypothesis 
\cite{Kuz-Kuz}-\cite{Mikh-2aaaaa}
that all the formulas (\ref{prod1})--(\ref{prod23}) can be replaced by 
the formulas (\ref{prod1}), (\ref{prod2}), (\ref{prod6}), (\ref{prod8}), (\ref{prod10})
in which we can suppose that $i_1, i_2, i_3, i_4=1,\ldots, m.$
At that we will not have a noticeable loss of the 
mean-square approximation accuracy of iterated It\^{o} stochastic integrals. 

It should be noted that unlike the method based on Theorem 1.1, existing approaches 
to the mean-square approximation of iterated stochastic integrals based on Fourier
series (see, for example, \cite{Zapad-1}-\cite{Zapad-4}, \cite{Zapad-8}, \cite{Zapad-11}) 
do not allow to choose different numbers $p$ (see (\ref{agentys1000})) for 
approximations of different iterated stochastic integrals with multiplicities 
$k=2,3,4,\ldots $
Moreover, the noted approaches exclude the possibility for obtaining of approximate 
and exact expressions similar to (\ref{tttr11}), (\ref{z1})
(see Theorems 1.3, 1.4). The detailed comparison of Theorem 1.1 with methods 
from \cite{Zapad-1}-\cite{Zapad-4}, \cite{Zapad-8}-\cite{Zapad-10}, \cite{Zapad-11},
\cite{Zapad-12a}, \cite{Zapad-12xxx} is given in Chapter 6 of this monograph.

Consider the following conditions
\begin{equation}
\label{agent36}
~~~~~~~ E^{(00)}_q\le (T-t)^4,\ \ \ E^{(000)}_{q_{1,i}} \le (T-t)^4,\ \ \ i=1,\ldots,4,
\end{equation}
and
\begin{equation}
\label{agent37}
~~~~~~~ E^{(00)}_q\le (T-t)^5,\ \ \ E^{(000)}_{q_{1,i}} \le (T-t)^5,\ \ \ 
E^{(10)}_{q_{2,j}}\le (T-t)^5,
\end{equation}
\begin{equation}
\label{agent37a}
E^{(01)}_{\bar q_{2,j}}\le (T-t)^5,\ \ \ 
E^{(0000)}_{q_{3,k}}\le (T-t)^5,
\end{equation}

\vspace{3mm}
\noindent
where $i=1,\ldots,4;$\ \ $j=1,2;$\ \ $k=1,\ldots,14.$

Let us show by numerical experiments that in most situations the following inequalities 
are fulfilled (under conditions (\ref{agent36}) and (\ref{agent37}),
(\ref{agent37a}))
\begin{equation}
\label{agentor10}
q_{1,1}\ge q_{1,i},\ \ \ i=2,3,4,
\end{equation}
\begin{equation}
\label{agentor11}
q_{2,1}\ge q_{2,2},\ \ \ 
\bar q_{2,1}\ge \bar q_{2,2},
\end{equation}

\vspace{-5mm}
\begin{equation}
\label{agentor12}
q_{3,1}\ge q_{3,k},\ \ \ k=2,\ldots,14,
\end{equation}

\noindent
where $q_{1,i}$, $q_{2,j},$ $\bar q_{2,j},$ $q_{3,k}$
($i=1,\ldots,4;$\ \ $j=1,2;$\ \ $k=1,\ldots,14$)
are minimal natural numbers satisfying 
the conditions (\ref{agent36}) and (\ref{agent37}),
(\ref{agent37a}).

In Tables 5.49--5.56 we can see the results of numerical experiments. 
These results confirm the hypothesis proposed earlier in this section.
Note that in Tables 5.54--5.56 we calculate the mean-square 
approximation errors of iterated It\^{o} stochastic integrals
in the case when 
$$
q_{1,i}=q_{1,1},\ \ \ i=2,3,4,
$$
$$
q_{2,2}=q_{2,1},\ \ \ 
\bar q_{2,2}=\bar q_{2,1},
$$

\vspace{-6mm}
$$
q_{3,k}=q_{3,1},\ \ \ k=2,\ldots,14,
$$

\noindent
where $q_{1,1}$, $q_{2,1},$ $\bar q_{2,1},$ $q_{3,1}$
are minimal natural numbers satisfying 
the conditions (\ref{agent36}) and (\ref{agent37}),
(\ref{agent37a}). In this case, there is no noticeable loss of the 
mean-square approximation accuracy of iterated It\^{o} stochastic integrals
(see Tables 5.54--5.56). 
This means that all the formulas (\ref{prod1})--(\ref{prod23}) can be re\-placed by 
the formulas (\ref{prod1}), (\ref{prod2}), (\ref{prod6}), (\ref{prod8}), (\ref{prod10})
in which we can suppose that $i_1, i_2, i_3, i_4=1,\ldots, m.$

Let $q_{1,1}$ and $q_{3,1}$ be minimal natural numbers satisfying the conditions
\begin{equation}
\label{agentor13}
E^{(000)}_{q_{1,1}} \le (T-t)^4,
\end{equation}
\begin{equation}
\label{agentor14}
E^{(0000)}_{q_{3,1}}\le (T-t)^5,
\end{equation}

\noindent
where the left-hand sides of 
these inequalities are defined by the formulas (\ref{prod2}) and (\ref{prod10}), respectively.

Let $p_{1,1}$ and $p_{3,1}$ be minimal natural numbers satisfying the conditions
\begin{equation}
\label{agentor15}
3! \cdot E^{(000)}_{p_{1,1}} \le (T-t)^4,
\end{equation}
\begin{equation}
\label{agentor16}
4! \cdot E^{(0000)}_{p_{3,1}}\le (T-t)^5,
\end{equation}

\noindent
where the values $E^{(000)}_{p_{1,1}}$ and $E^{(0000)}_{p_{3,1}}$ on the left-hand sides of 
these inequalities are defined by the formulas (\ref{prod2}) and (\ref{prod10}), respectively.

In Tables 5.57, 5.58 we can see the 
numerical comparison of the numbers $q_{1,1}$ and $q_{3,1}$ with the numbers
$p_{1,1}$ and $p_{3,1}$, respectively. Obviously, 
excluding of the multiplier factors\hspace{0.3mm} $3!$ and $4!$\hspace{0.3mm} 
essentially (in many times)\hspace{0.3mm} reduces 
the calculation

\begin{table}
\centering
\caption{Conditions $E^{(000)}_{q_{1,i}} \le (T-t)^4$,\ $i=1,\ldots,4.$}
\label{tab:5.49}      

\begin{tabular}{p{1.3cm}p{1.3cm}p{1.3cm}p{1.3cm}p{1.3cm}p{1.3cm}p{1.3cm}}

\hline\noalign{\smallskip}

$T-t$&0.011&0.008&0.0045&0.0035&0.0027&0.0025\\

\noalign{\smallskip}\hline\noalign{\smallskip}
$q_{1,1}$&12&16&28&36&47&50\\
\noalign{\smallskip}\hline\noalign{\smallskip}
$q_{1,2}$&6&8&14&18&23&25\\
\noalign{\smallskip}\hline\noalign{\smallskip}
$q_{1,3}$&6&8&14&18&23&25\\
\noalign{\smallskip}\hline\noalign{\smallskip}
$q_{1,4}$&12&16&28&36&47&51\\
\noalign{\smallskip}\hline\noalign{\smallskip}
\end{tabular}
\end{table}

\begin{table}
\centering
\caption{Conditions $E^{(0000)}_{q_{3,k}} \le (T-t)^5$,\ $i=1,\ldots,14.$}
\label{tab:5.50}      
\begin{tabular}{p{1.5cm}p{1.5cm}p{1.5cm}p{1.5cm}p{1.5cm}p{1.5cm}}
\hline\noalign{\smallskip}
$T-t$&$0.011$&$0.008$&$0.0045$&$0.0042$&0.0040\\
\noalign{\smallskip}\hline\noalign{\smallskip}
$q_{3,1}$&$6$&$8$&$14$&$15$&$16$\\
\noalign{\smallskip}\hline\noalign{\smallskip}
$q_{3,2}$&4&5&10&11&11\\
\noalign{\smallskip}\hline\noalign{\smallskip}
$q_{3,3}$&6&8&14&15&16\\
\noalign{\smallskip}\hline\noalign{\smallskip}
$q_{3,4}$&6&8&14&15&16\\
\noalign{\smallskip}\hline\noalign{\smallskip}
$q_{3,5}$&3&5&9&9&10\\
\noalign{\smallskip}\hline\noalign{\smallskip}
$q_{3,6}$&6&8&14&15&16\\
\noalign{\smallskip}\hline\noalign{\smallskip}
$q_{3,7}$&4&5&10&11&11\\
\noalign{\smallskip}\hline\noalign{\smallskip}
$q_{3,8}$&2&3&4&5&5\\
\noalign{\smallskip}\hline\noalign{\smallskip}
$q_{3,9}$&2&3&4&5&5\\
\noalign{\smallskip}\hline\noalign{\smallskip}
$q_{3,10}$&4&6&10&11&11\\
\noalign{\smallskip}\hline\noalign{\smallskip}
$q_{3,11}$&4&6&10&11&11\\
\noalign{\smallskip}\hline\noalign{\smallskip}
$q_{3,12}$&2&3&5&6&6\\
\noalign{\smallskip}\hline\noalign{\smallskip}
$q_{3,13}$&6&8&14&15&16\\
\noalign{\smallskip}\hline\noalign{\smallskip}
$q_{3,14}$&3&5&9&9&10\\
\noalign{\smallskip}\hline\noalign{\smallskip}
\end{tabular}
\end{table}

\begin{table}
\centering
\caption{The conditions (\ref{agent37}), (\ref{agent37a}).}
\label{tab:5.51}      
\begin{tabular}{p{1.8cm}p{1.8cm}p{1.8cm}p{1.8cm}}
\hline\noalign{\smallskip}
$T-t$&0.010&0.005&0.0025\\
\noalign{\smallskip}\hline\noalign{\smallskip}
$q_{2,1}$&4&8&16\\
\noalign{\smallskip}\hline\noalign{\smallskip}
$q_{2,2}$&1&1&1\\
\noalign{\smallskip}\hline\noalign{\smallskip}
$\bar q_{2,1}$&4&8&16\\
\noalign{\smallskip}\hline\noalign{\smallskip}
$\bar q_{2,2}$&1&1&1\\
\noalign{\smallskip}\hline\noalign{\smallskip}
\end{tabular}
\end{table}

\begin{table}
\centering
\caption{The condition (\ref{agent36}).}
\label{tab:5.52}      
\begin{tabular}{p{1.5cm}p{1.5cm}p{1.5cm}p{1.5cm}p{1.5cm}}
\hline\noalign{\smallskip}
$T-t$&$2^{-1}$&$2^{-3}$&$2^{-5}$&$2^{-8}$\\
\noalign{\smallskip}\hline\noalign{\smallskip}
$q$&1&8&128&8192\\
\noalign{\smallskip}\hline\noalign{\smallskip}
$q_{1,1}$&0&1&4&32\\
\noalign{\smallskip}\hline\noalign{\smallskip}
$q_{1,2}$&0&0&2&16\\
\noalign{\smallskip}\hline\noalign{\smallskip}
$q_{1,3}$&0&0&2&16\\
\noalign{\smallskip}\hline\noalign{\smallskip}
$q_{1,4}$&0&0&4&33\\
\noalign{\smallskip}\hline\noalign{\smallskip}
\end{tabular}
\end{table}

\begin{table}
\centering
\caption{The conditions (\ref{agent37}), (\ref{agent37a}).}
\label{tab:5.53}      
\begin{tabular}{p{2.5cm}p{1.3cm}p{1.3cm}p{1.3cm}p{1.3cm}}
\hline\noalign{\smallskip}
$T-t$&$2^{-1}$&$2^{-3}$&$2^{-5}$&$2^{-8}$\\
\noalign{\smallskip}\hline\noalign{\smallskip}
$q$&1&8&64&512\\
\noalign{\smallskip}\hline\noalign{\smallskip}
$q_{1,1}$&0&2&4&32\\
\noalign{\smallskip}\hline\noalign{\smallskip}
$q_{1,2}$&0&1&4&16\\
\noalign{\smallskip}\hline\noalign{\smallskip}
$q_{1,3}$&0&1&4&16\\
\noalign{\smallskip}\hline\noalign{\smallskip}
$q_{1,4}$&0&2&8&33\\
\noalign{\smallskip}\hline\noalign{\smallskip}
$\bar q_{2,1}$&0&0&1&1\\
\noalign{\smallskip}\hline\noalign{\smallskip}
$\bar q_{2,2}, q_{2,1}, q_{2,2}$&0&0&0&0\\
\noalign{\smallskip}\hline\noalign{\smallskip}
$q_{3,1}, \ldots, q_{3,14}$&0&0&0&0\\
\noalign{\smallskip}\hline\noalign{\smallskip}
\end{tabular}
\end{table}

\begin{table}
\centering
\caption{Values $E_{q_{1,i}}^{(000)} \cdot (T-t)^{-3}\stackrel{\sf def}{=}E_{q_{1,i}},$\ 
$i=1,\ldots,4.$}
\label{tab:5.54}      
\begin{tabular}{p{1.3cm}p{1.9cm}p{1.9cm}p{1.9cm}p{1.9cm}p{1.9cm}p{1.9cm}}
\hline\noalign{\smallskip}
$T-t$&$0.011$&$0.008$&$0.0045$&$0.0035$&$0.0027$&$0.0025$\\
\noalign{\smallskip}\hline\noalign{\smallskip}
$q_{1,1}$&$12$&$16$&$28$&$36$&$47$&$50$\\
\noalign{\smallskip}\hline\noalign{\smallskip}
$E_{q_{1,1}}$&0.010154&0.007681&0.004433&0.003456&0.002652&0.002494\\
\noalign{\smallskip}\hline\noalign{\smallskip}
$q_{1,2}$&$12$&$16$&$28$&$36$&$47$&$50$\\
\noalign{\smallskip}\hline\noalign{\smallskip}
$E_{q_{1,2}}$&0.005077&0.003841&0.002216&0.001728&0.001326&0.001247\\
\noalign{\smallskip}\hline\noalign{\smallskip}
$q_{1,3}$&$12$&$16$&$28$&$36$&$47$&$50$\\
\noalign{\smallskip}\hline\noalign{\smallskip}
$E_{q_{1,3}}$&0.005077&0.003841&0.002216&0.001728&0.001326&0.001247\\
\noalign{\smallskip}\hline\noalign{\smallskip}
$q_{1,4}$&$12$&$16$&$28$&$36$&$47$&$50$\\
\noalign{\smallskip}\hline\noalign{\smallskip}
$E_{q_{1,4}}$&0.010308&0.007787&0.004480&0.003488&0.002673&0.002513\\
\noalign{\smallskip}\hline\noalign{\smallskip}
\end{tabular}
\end{table}

\begin{table}
\centering
\caption{Values $E_{\bar q_{2,j}}^{(01)} \cdot (T-t)^{-4}\stackrel{\sf def}{=}E_{\bar q_{2,j}},$\
$E_{q_{2,j}}^{(10)} \cdot (T-t)^{-4}\stackrel{\sf def}{=}E_{q_{2,j}}$,\ $j=1,2.$}
\label{tab:5.55}      
\begin{tabular}{p{1.8cm}p{2.4cm}p{2.4cm}p{2.4cm}}
\hline\noalign{\smallskip}
$T-t$&$0.010$&$0.005$&$0.0025$\\
\noalign{\smallskip}\hline\noalign{\smallskip}
$\bar q_{2,1}$&4&8&16\\
\noalign{\smallskip}\hline\noalign{\smallskip}
$E_{\bar q_{2,1}}$&0.008950&0.004660&0.002383\\
\noalign{\smallskip}\hline\noalign{\smallskip}
$\bar q_{2,2}$&4&8&16\\
\noalign{\smallskip}\hline\noalign{\smallskip}
$E_{\bar q_{2,2}}$&0.000042&0.000006&0.000001\\
\noalign{\smallskip}\hline\noalign{\smallskip}
$q_{2,1}$&4&8&16\\
\noalign{\smallskip}\hline\noalign{\smallskip}
$E_{q_{2,1}}$&0.008950&0.004660&0.002383\\
\noalign{\smallskip}\hline\noalign{\smallskip}
$q_{2,2}$&4&8&16\\
\noalign{\smallskip}\hline\noalign{\smallskip}
$E_{q_{2,2}}$&0.000042&0.000006&0.000001\\
\noalign{\smallskip}\hline\noalign{\smallskip}
\end{tabular}
\end{table}

\begin{table}
\centering
\vspace{-7mm}
\caption{Values $E_{q_{3,k}}^{(0000)}\cdot (T-t)^{-4}\stackrel{\sf def}{=}E_{q_{3,k}},$\ 
$k=1,\ldots,14.$}
\label{tab:5.56}      
\begin{tabular}{p{1.6cm}p{2.2cm}p{2.2cm}p{2.2cm}p{2.2cm}}
\hline\noalign{\smallskip}
$T-t$&$0.011$&$0.008$&$0.0045$&$0.0042$\\
\noalign{\smallskip}\hline\noalign{\smallskip}
$q_{3,1}$&6&8&14&15\\
\noalign{\smallskip}\hline\noalign{\smallskip}
$E_{q_{3,1}}$&0.009636&0.007425&0.004378&0.004096\\
\noalign{\smallskip}\hline\noalign{\smallskip}
$q_{3,2}$&6&8&14&15\\
\noalign{\smallskip}\hline\noalign{\smallskip}
$E_{q_{3,2}}$&0.006771&0.005191&0.003041&0.002843\\
\noalign{\smallskip}\hline\noalign{\smallskip}
$q_{3,3}$&6&8&14&15\\
\noalign{\smallskip}\hline\noalign{\smallskip}
$E_{q_{3,3}}$&0.009722&0.007502&0.004424&0.004139\\
\noalign{\smallskip}\hline\noalign{\smallskip}
$q_{3,4}$&6&8&14&15\\
\noalign{\smallskip}\hline\noalign{\smallskip}
$E_{q_{3,4}}$&0.009641&0.007427&0.004379&0.004097\\
\noalign{\smallskip}\hline\noalign{\smallskip}
$q_{3,5}$&6&8&14&15\\
\noalign{\smallskip}\hline\noalign{\smallskip}
$E_{q_{3,5}}$&0.005997&0.004614&0.002720&0.002545\\
\noalign{\smallskip}\hline\noalign{\smallskip}
$q_{3,6}$&6&8&14&15\\
\noalign{\smallskip}\hline\noalign{\smallskip}
$E_{q_{3,6}}$&0.009722&0.007502&0.004424&0.004139\\
\noalign{\smallskip}\hline\noalign{\smallskip}
$q_{3,7}$&6&8&14&15\\
\noalign{\smallskip}\hline\noalign{\smallskip}
$E_{q_{3,7}}$&0.006771&0.005191&0.003041&0.002843\\
\noalign{\smallskip}\hline\noalign{\smallskip}
$q_{3,8}$&6&8&14&15\\
\noalign{\smallskip}\hline\noalign{\smallskip}
$E_{q_{3,8}}$&0.003095&0.002364&0.001379&0.001290\\
\noalign{\smallskip}\hline\noalign{\smallskip}
$q_{3,9}$&6&8&14&15\\
\noalign{\smallskip}\hline\noalign{\smallskip}
$E_{q_{3,9}}$&0.003095&0.002364&0.001379&0.001290\\
\noalign{\smallskip}\hline\noalign{\smallskip}
$q_{3,10}$&6&8&14&15\\
\noalign{\smallskip}\hline\noalign{\smallskip}
$E_{q_{3,10}}$&0.006885&0.005282&0.003090&0.002889\\
\noalign{\smallskip}\hline\noalign{\smallskip}
$q_{3,11}$&6&8&14&15\\
\noalign{\smallskip}\hline\noalign{\smallskip}
$E_{q_{3,11}}$&0.006885&0.005282&0.003090&0.002889\\
\noalign{\smallskip}\hline\noalign{\smallskip}
$q_{3,12}$&6&8&14&15\\
\noalign{\smallskip}\hline\noalign{\smallskip}
$E_{q_{3,12}}$&0.003690&0.002834&0.001663&0.001555\\
\noalign{\smallskip}\hline\noalign{\smallskip}
$q_{3,13}$&6&8&14&15\\
\noalign{\smallskip}\hline\noalign{\smallskip}
$E_{q_{3,13}}$&0.009756&0.007545&0.004457&0.004170\\
\noalign{\smallskip}\hline\noalign{\smallskip}
$q_{3,14}$&6&8&14&15\\
\noalign{\smallskip}\hline\noalign{\smallskip}
$E_{q_{3,14}}$&0.006010&0.004621&0.002722&0.002547\\
\noalign{\smallskip}\hline\noalign{\smallskip}
\end{tabular}
\end{table}

\begin{table}
\centering
\caption{Comparison of numbers $q_{1,1}$ and $p_{1,1}$.}
\label{tab:5.57}      
\begin{tabular}{p{1.9cm}p{1.3cm}p{1.3cm}p{1.3cm}p{1.3cm}p{1.3cm}p{1.3cm}}
\hline\noalign{\smallskip}
$T-t$&$2^{-1}$&$2^{-2}$&$2^{-3}$&$2^{-4}$&$2^{-5}$&$2^{-6}$\\
\noalign{\smallskip}\hline\noalign{\smallskip}
$q_{1,1}$&0&0&1&2&4&8\\
\noalign{\smallskip}\hline\noalign{\smallskip}
$(q_{1,1}+1)^3$&1&1&8&27&125&729\\
\noalign{\smallskip}\hline\noalign{\smallskip}
$p_{1,1}$&1&3&6&12&24&48\\
\noalign{\smallskip}\hline\noalign{\smallskip}
$(p_{1,1}+1)^3$&8&64&343&2197&15625&117649\\
\noalign{\smallskip}\hline\noalign{\smallskip}
\end{tabular}
\end{table}

\begin{table}
\centering
\caption{Comparison of numbers $q_{3,1}$ and $p_{3,1}$.}
\label{tab:5.58}      
\begin{tabular}{p{1.9cm}p{1.3cm}p{1.3cm}p{1.3cm}p{1.3cm}p{1.3cm}p{1.3cm}}
\hline\noalign{\smallskip}
$T-t$&$2^{-1}$&$2^{-2}$&$2^{-3}$&$2^{-4}$&$2^{-5}$&$2^{-6}$\\
\noalign{\smallskip}\hline\noalign{\smallskip}
$q_{3,1}$&0&0&0&0&0&0\\
\noalign{\smallskip}\hline\noalign{\smallskip}
$(q_{3,1}+1)^4$&1&1&1&1&1&1\\
\noalign{\smallskip}\hline\noalign{\smallskip}
$p_{3,1}$&3&4&6&9&12&17\\
\noalign{\smallskip}\hline\noalign{\smallskip}
$(p_{3,1}+1)^4$&256&625&2401&10000&28561&104976\\
\noalign{\smallskip}\hline\noalign{\smallskip}
\end{tabular}
\end{table}

\noindent
costs for the mean-square approximations of iterated It\^{o} stochastic 
integrals. Note that in this section we use the exactly calculated Fourier--Legendre 
coefficients using the Python programming language \cite{Kuz-Kuz}, \cite{Mikh-1}.

As we mentioned above, 
existing approaches to the mean-square approximation of iterated 
stochastic integrals based on Fourier
series (see, for example, \cite{Zapad-1}-\cite{Zapad-4}, \cite{Zapad-8}, \cite{Zapad-11}) 
do not allow to choose different numbers $p$  
(see Theorem 1.3) for approximations of different iterated stochastic integrals 
with multiplicities $k=2,3,4,\ldots $
and exclude the possibility for obtaining of approximate and exact expressions similar 
to the formulas (\ref{tttr11}), (\ref{z1}) (see Theorems 1.3, 1.4). 
This leads to unnecessary terms usage in the expansions of iterated 
It\^{o} stochastic integrals and, as a consequence, to essential increase of computational 
costs for the implementation of numerical methods for It\^{o} SDEs.

In this section
(also see \cite{Mikh-2}, \cite{Mikh-2aaaaa})
we have optimized the method based on Theorems 1.1 and 1.3, which makes it possible to correctly 
choose the lengths of sequences of standard Gaussian random variables required for the 
approximation of iterated It\^{o} stochastic integrals. Thus, the computational costs for 
the implementation of numerical methods for It\^{o} SDEs are significantly reduced.

On the base of the obtained results we recommend to use in practice 
the following conditions
(for any $i_1,\ldots,i_4=1,\ldots,m$) for 
correct choosing the minimal natural numbers $q,$ $q_1,$ $q_2,$ $\bar q_2,$ $q_3$

\vspace{-2mm}
$$
E^{(00)}_q\le C(T-t)^3
$$

\noindent
(for the Milstein scheme (\ref{al1})),  

$$
E^{(00)}_q\le (T-t)^4,\ \ \ E^{(000)}_{q_{1,1}} \le C(T-t)^4
$$

\vspace{3mm}
\noindent
(for the strong scheme (\ref{al2}) with order 1.5), and 

$$
E^{(00)}_q\le C(T-t)^5,\ \ \ E^{(000)}_{q_{1,1}} \le C(T-t)^5,\ \ \ 
E^{(10)}_{q_{2,1}}\le C(T-t)^5,
$$
$$
E^{(01)}_{\bar q_{2,1}}\le C(T-t)^5,\ \ \ E^{(0000)}_{q_{3,1}}\le C(T-t)^5
$$

\vspace{3mm}
\noindent
(for the strong scheme (\ref{al3}) with order 2.0). Here the 
left-hand sides of the above inequalities are defined by the relations (\ref{prod1}), 
(\ref{prod2}), (\ref{prod6}), 
(\ref{prod8}), (\ref{prod10})
and $C$ is a constant from the condition (\ref{agentww3}).

Taking into account the results of this section (also see \cite{Mikh-2}),
we recommend to use in practice the following condition (for any $i_1,\ldots,i_k=1,\ldots,m$)
on the mean-square approximation accuracy for iterated 
It\^{o} stochastic integrals 
$$
{\sf M}\left\{\left(
I_{(l_1\ldots l_k)T,t}^{(i_1\ldots i_k)}-
I_{(l_1\ldots l_k)T,t}^{(i_1\ldots i_k)p}\right)^2\right\}=
$$
$$
=
\int\limits_{[t,T]^k}
K^2(t_1,\ldots,t_k)
dt_1\ldots dt_k -\sum_{j_1,\ldots,j_k=0}^{p}C^2_{j_k\ldots j_1}\le 
C(T-t)^{r+1},
$$
 
\noindent
where $I_{(l_1\ldots l_k)T,t}^{(i_1\ldots i_k)}$ is the 
iterated It\^{o} stochastic integral (\ref{k1000xxxx}),
$I_{(l_1\ldots l_k)T,t}^{(i_1\ldots i_k)p}$ is the
mean-square approximation of this stochastic integral 
based on Theorem 1.1 and multiple Fourier--Legendre series,
$p$ and $k\in {\bf N},$ 
$$
K(t_1,\ldots,t_k)=
(t-t_k)^{l_k}\ldots (t-t_1)^{l_1}\ {\bf 1}_{\{t_1<\ldots<t_k\}},\ \ \
t_1,\ldots,t_k\in[t, T],
$$

\noindent
${\bf 1}_A$ is the indicator of the set $A$,\ \ $l_1,\ldots,l_k=0, 1,\ldots,$\ \ 
$C$ and $r$ have the same meaning as in the formula (\ref{agentww3}).

\section{Exact Calculation 
of the Mean-Square Approximation Errors
for Iterated Stra\-to\-no\-vich Stochastic Integrals
$I_{(0)T,t}^{*(i_1)},$
$I_{(1)T,t}^{*(i_1)},$
$I_{(00)T,t}^{*(i_1i_2)},$
$I_{(000)T,t}^{*(i_1i_2i_3)}$}

Consider the question on the exact calculation 
of the mean-square approximation errors
for the following iterated Stratonovich stochastic integrals
\begin{equation}
\label{dest1}
~~~~~~~~~I_{(0)T,t}^{*(i_1)},\ \ \ 
I_{(1)T,t}^{*(i_1)},\ \ \ 
I_{(00)T,t}^{*(i_1i_2)},\ \ \ 
I_{(000)T,t}^{*(i_1i_2i_3)},\ \ \ 
i_1, i_2, i_3=1,\ldots,m.
\end{equation}

We assume that the stochastic integrals (\ref{dest1})
are approximated using Theorems 1.1, 2.1, 2.8 and the Legendre
polynomial system. Since
$I_{(0)T,t}^{(i_1)}=I_{(0)T,t}^{*(i_1)},$
$I_{(1)T,t}^{(i_1)}=I_{(1)T,t}^{*(i_1)}$\ w.~p.~1,
we can use (\ref{4001}), (\ref{4002}) 
to approximate the stochastic integrals 
$I_{(0)T,t}^{*(i_1)},$
$I_{(1)T,t}^{*(i_1)}.$ In this case, we will have zero
mean-square approximation errors.

To approximate the iterated Stratonovich stochastic integral 
$I_{(00)T,t}^{*(i_1i_2)}$ 
we can use the formula (see (\ref{4004xxxx}))
\begin{equation}
\label{dest2}
~~~~I_{(00)T,t}^{*(i_1 i_2)q}=
\frac{T-t}{2}\left(\zeta_0^{(i_1)}\zeta_0^{(i_2)}+\sum_{i=1}^{q}
\frac{1}{\sqrt{4i^2-1}}\left(
\zeta_{i-1}^{(i_1)}\zeta_{i}^{(i_2)}-
\zeta_i^{(i_1)}\zeta_{i-1}^{(i_2)}\right)\right).
\end{equation}

The mean-square approximation error for (\ref{dest2})
will be determined by the formula (\ref{fff09}) $(i_1\ne i_2)$.
For the case $i_1=i_2$ we can use the formula (see (\ref{4.2.880}))
$$
I_{(00)T,t}^{*(i_1 i_1)}
=\frac{T-t}{2}\left(\zeta_0^{(i_1)}\right)^2\ \ \ \hbox{w.~p.~1.}
$$

Consider now the iterated Stratonovich stochastic integral
$I_{(000)T,t}^{*(i_1i_2i_3)}$ of multiplicity 3
$(i_1, i_2, i_3=$ $1,\ldots,m)$.
For the case of pairwise different 
$i_1, i_2, i_3$ we can use the formula (\ref{dest3}).
In the case $i_1=i_2=i_3,$ to approximate the stochastic integral
$I_{(000)T,t}^{*(i_1i_1i_1)},$ we use the formula (\ref{dest4}).

Thus, it remains to consider the following three cases
\begin{equation}
\label{dest5}
i_1=i_2\ne i_3,
\end{equation}
\begin{equation}
\label{dest6}
i_1\ne i_2=i_3,
\end{equation}

\vspace{-6mm}
\begin{equation}
\label{dest7}
i_1=i_3\ne i_2.
\end{equation}

Consider the case (\ref{dest5}). From (\ref{dest10}) we obtain
$$
{\sf M}\left\{\left(I_{(000)T,t}^{*(i_1i_2i_3)}-
I_{(000)T,t}^{*(i_1i_2i_3)q}\right)^2\right\}=
$$
\begin{equation}
\label{dest11}
={\sf M}\left\{\left(I_{(000)T,t}^{(i_1i_2i_3)}-
I_{(000)T,t}^{(i_1i_2i_3)q}+
\frac{1}{2}\int\limits_t^T
\int\limits_t^{\tau}dsd{\bf w}_{\tau}^{(i_3)}-
\sum_{j_1,j_3=0}^{q}
C_{j_3j_1j_1}
\zeta_{j_3}^{(i_3)}\right)^2\right\}.
\end{equation}

\vspace{2mm}

According to the formulas (\ref{yeee2}), (\ref{ttt2}), the quantity
$$
I_{(000)T,t}^{(i_1i_2i_3)}-
I_{(000)T,t}^{(i_1i_2i_3)q}
$$
includes only iterated It\^{o} stochastic integrals
of multiplicity 3. At the same time, the quantity
$$
\frac{1}{2}\int\limits_t^T
\int\limits_t^{\tau}dsd{\bf w}_{\tau}^{(i_3)}-
\sum_{j_1,j_3=0}^{q}
C_{j_3j_1j_1}
\zeta_{j_3}^{(i_3)}
$$
contains only iterated It\^{o} stochastic integrals
of multiplicity 1. This means that from (\ref{dest11}) we get
$$
{\sf M}\left\{\left(I_{(000)T,t}^{*(i_1i_2i_3)}-
I_{(000)T,t}^{*(i_1i_2i_3)q}\right)^2\right\}
={\sf M}\left\{\left(I_{(000)T,t}^{(i_1i_2i_3)}-
I_{(000)T,t}^{(i_1i_2i_3)q}\right)^2\right\}+
$$
\begin{equation}
\label{dest12}
+{\sf M}\left\{\left(\frac{1}{2}\int\limits_t^T
(\tau-t)d{\bf w}_{\tau}^{(i_3)}-
\sum_{j_1,j_3=0}^{q}
C_{j_3j_1j_1}
\zeta_{j_3}^{(i_3)}\right)^2\right\}.
\end{equation}

The relation (\ref{qq1}) implies that
$$
{\sf M}\left\{\left(I_{(000)T,t}^{(i_1i_2i_3)}-
I_{(000)T,t}^{(i_1i_2i_3)q}\right)^2\right\}=\frac{(T-t)^3}{6}
-
$$
\begin{equation}
\label{dest14}
-\sum_{j_1,j_2,j_3=0}^q C_{j_3j_2j_1}^2-
\sum_{j_1,j_2,j_3=0}^q C_{j_3j_1j_2}C_{j_3j_2j_1},
\end{equation}

\noindent
where $i_1=i_2\ne i_3.$

We have
$$
{\sf M}\left\{\left(\frac{1}{2}\int\limits_t^T
(\tau-t)d{\bf w}_{\tau}^{(i_3)}-
\sum_{j_1,j_3=0}^{q}
C_{j_3j_1j_1}
\zeta_{j_3}^{(i_3)}\right)^2\right\}=
\frac{1}{4}\int\limits_t^T
(\tau-t)^2 d\tau-
$$
\begin{equation}
\label{dest15}
~~~~~~~ -\sum_{j_1,j_3=0}^{q}
C_{j_3j_1j_1}\int\limits_t^T
(\tau-t)\phi_{j_3}(\tau)d\tau+
\sum_{j_3=0}^{q}\left(\sum_{j_1=0}^{q}
C_{j_3j_1j_1}\right)^2,
\end{equation}

\vspace{2mm}
\noindent
where $\phi_{j_3}(\tau)$ is the Legendre polynomial defined by (\ref{4009d}).

According to (\ref{ogo11}), we obtain
\begin{equation}
\label{dest16}
\int\limits_t^T
(\tau-t)\phi_{j_3}(\tau)d\tau=\frac{(T-t)^{3/2}}{2}
\left\{
\begin{matrix}
1,\ & j_3=0\cr\cr
1/\sqrt{3},\ & j_3=1\cr\cr
0,\ & j_3\ge 2
\end{matrix}
.\right.
\end{equation}

Combining (\ref{dest12})--(\ref{dest16}), we get
$$
{\sf M}\left\{\left(I_{(000)T,t}^{*(i_1i_2i_3)}-
I_{(000)T,t}^{*(i_1i_2i_3)q}\right)^2\right\}=
\frac{(T-t)^3}{4}
-
$$
$$
-\sum_{j_1,j_2,j_3=0}^q C_{j_3j_2j_1}^2-
\sum_{j_1,j_2,j_3=0}^q C_{j_3j_1j_2}C_{j_3j_2j_1}-
$$
$$
-\frac{(T-t)^{3/2}}{2}
\sum_{j_1=0}^{q}
\left(C_{0j_1j_1}+\frac{1}{\sqrt{3}}C_{1j_1j_1}\right)+
$$
\begin{equation}
\label{dest80}
+
\sum_{j_3=0}^{q}\left(\sum_{j_1=0}^{q}
C_{j_3j_1j_1}\right)^2,
\end{equation}

\noindent
where $i_1=i_2\ne i_3.$

Consider the case (\ref{dest6}). From (\ref{dest10}) we obtain
$$
{\sf M}\left\{\left(I_{(000)T,t}^{*(i_1i_2i_3)}-
I_{(000)T,t}^{*(i_1i_2i_3)q}\right)^2\right\}=
$$
$$
={\sf M}\left\{\left(I_{(000)T,t}^{(i_1i_2i_3)}-
I_{(000)T,t}^{(i_1i_2i_3)q}+
\frac{1}{2}\int\limits_t^T
\int\limits_t^{\tau}d{\bf w}_{s}^{(i_1)}d\tau-
\sum_{j_1,j_3=0}^{q}
C_{j_3j_3j_1}
\zeta_{j_1}^{(i_1)}\right)^2\right\}=
$$
$$
={\sf M}\left\{\left(I_{(000)T,t}^{(i_1i_2i_3)}-
I_{(000)T,t}^{(i_1i_2i_3)q}+
\frac{1}{2}\int\limits_t^T
(T-s)d{\bf w}_{s}^{(i_1)}-
\sum_{j_1,j_3=0}^{q}
C_{j_3j_3j_1}
\zeta_{j_1}^{(i_1)}\right)^2\right\}=
$$

$$
={\sf M}\left\{\left(I_{(000)T,t}^{(i_1i_2i_3)}-
I_{(000)T,t}^{(i_1i_2i_3)q}\right)^2\right\}+
$$
$$
+
{\sf M}\left\{\left(\frac{1}{2}\int\limits_t^T
(T-s)d{\bf w}_{s}^{(i_1)}-
\sum_{j_1,j_3=0}^{q}
C_{j_3j_3j_1}
\zeta_{j_1}^{(i_1)}\right)^2\right\}=
$$

$$
={\sf M}\left\{\left(I_{(000)T,t}^{(i_1i_2i_3)}-
I_{(000)T,t}^{(i_1i_2i_3)q}\right)^2\right\}+
$$
$$
+\frac{1}{4}\int\limits_t^T
(T-s)^2 ds-
\sum_{j_1,j_3=0}^{q}
C_{j_3j_3j_1}\int\limits_t^T
(T-s)\phi_{j_1}(s)ds+
$$
\begin{equation}
\label{dest27}
+\sum_{j_1=0}^{q}\left(\sum_{j_3=0}^{q}
C_{j_3j_3j_1}\right)^2,
\end{equation}

\noindent
where $\phi_{j_1}(\tau)$ is the Legendre polynomial defined by (\ref{4009d}).

The relation (\ref{dest30}) implies that
$$
{\sf M}\left\{\left(I_{(000)T,t}^{(i_1i_2i_3)}-
I_{(000)T,t}^{(i_1i_2i_3)q}\right)^2\right\}=\frac{(T-t)^3}{6}
-
$$
\begin{equation}
\label{dest31}
-\sum_{j_1,j_2,j_3=0}^q C_{j_3j_2j_1}^2-
\sum_{j_1,j_2,j_3=0}^q C_{j_2j_3j_1}C_{j_3j_2j_1},
\end{equation}

\noindent
where $i_1\ne i_2=i_3.$

Moreover,
\begin{equation}
\label{dest32}
~~~~~~~~\int\limits_t^T
(T-s)\phi_{j_1}(s)ds=\frac{(T-t)^{3/2}}{2}
\left\{
\begin{matrix}
1,\ & j_1=0\cr\cr
-1/\sqrt{3},\ & j_1=1\cr\cr
0,\ & j_1\ge 2
\end{matrix}
.\right.
\end{equation}

Combining (\ref{dest27})--(\ref{dest32}), we get
$$
{\sf M}\left\{\left(I_{(000)T,t}^{*(i_1i_2i_3)}-
I_{(000)T,t}^{*(i_1i_2i_3)q}\right)^2\right\}=
\frac{(T-t)^3}{4}
-
$$
$$
-\sum_{j_1,j_2,j_3=0}^q C_{j_3j_2j_1}^2-
\sum_{j_1,j_2,j_3=0}^q C_{j_2j_3j_1}C_{j_3j_2j_1}-
$$
$$
-\frac{(T-t)^{3/2}}{2}
\sum_{j_3=0}^{q}
\left(C_{j_3j_3 0}-\frac{1}{\sqrt{3}}C_{j_3j_3 1}\right)+
$$
\begin{equation}
\label{dest70}
+\sum_{j_1=0}^{q}\left(\sum_{j_3=0}^{q}
C_{j_3j_3j_1}\right)^2,
\end{equation}

\noindent
where $i_1\ne i_2=i_3.$

Consider the case (\ref{dest7}). From (\ref{dest10}) we obtain
$$
{\sf M}\left\{\left(I_{(000)T,t}^{*(i_1i_2i_3)}-
I_{(000)T,t}^{*(i_1i_2i_3)q}\right)^2\right\}=
$$
$$
={\sf M}\left\{\left(I_{(000)T,t}^{(i_1i_2i_3)}-
I_{(000)T,t}^{(i_1i_2i_3)q}-
\sum_{j_1,j_2=0}^{q}
C_{j_1j_2j_1}
\zeta_{j_2}^{(i_2)}\right)^{2}\right\}=
$$
$$
={\sf M}\left\{\left(I_{(000)T,t}^{(i_1i_2i_3)}-
I_{(000)T,t}^{(i_1i_2i_3)q}\right)^{2}\right\}+
{\sf M}\left\{\left(\sum_{j_1,j_2=0}^{q}
C_{j_1j_2j_1}
\zeta_{j_2}^{(i_2)}\right)^{2}\right\}=
$$
\begin{equation}
\label{dest49}
~~~~~~~ ={\sf M}\left\{\left(I_{(000)T,t}^{(i_1i_2i_3)}-
I_{(000)T,t}^{(i_1i_2i_3)q}\right)^{2}\right\}+
\sum_{j_2=0}^{q}
\left(\sum_{j_1=0}^{q}C_{j_1j_2j_1}\right)^2.
\end{equation}

\vspace{1mm}

The relation (\ref{dest40}) implies that
$$
{\sf M}\left\{\left(I_{(000)T,t}^{(i_1i_2i_3)}-
I_{(000)T,t}^{(i_1i_2i_3)q}\right)^2\right\}=\frac{(T-t)^3}{6}
-
$$
\begin{equation}
\label{dest50}
-\sum_{j_1,j_2,j_3=0}^q C_{j_3j_2j_1}^2-
\sum_{j_1,j_2,j_3=0}^q C_{j_3j_2j_1}C_{j_1j_2j_3},
\end{equation}

\vspace{0.4mm}
\noindent
where $i_1=i_3\ne i_2.$

Combining (\ref{dest49}) and (\ref{dest50}), we obtain
$$
{\sf M}\left\{\left(I_{(000)T,t}^{*(i_1i_2i_3)}-
I_{(000)T,t}^{*(i_1i_2i_3)q}\right)^2\right\}=
\frac{(T-t)^3}{6}-
$$
$$
-\sum_{j_1,j_2,j_3=0}^q C_{j_3j_2j_1}^2-
\sum_{j_1,j_2,j_3=0}^q C_{j_3j_2j_1}C_{j_1j_2j_3}+
$$
\begin{equation}
\label{dest60}
+\sum_{j_2=0}^{q}
\left(\sum_{j_1=0}^{q}C_{j_1j_2j_1}\right)^2,
\end{equation}

\vspace{0.4mm}
\noindent
where $i_1=i_3\ne i_2.$

Thus, the exact calculaton of the mean-square approximation error
for the iterated Stratonovich stochastic integral 
$I_{(000)T,t}^{*(i_1i_2i_3)}$ $(i_1,i_2,i_3=1,\ldots,m)$
is given by the formulas (\ref{dest3}),
(\ref{dest80}), (\ref{dest70}), and (\ref{dest60}).

\section{Exact Calculation 
of the Mean-Square Approximation Error
for Iterated Stra\-to\-no\-vich Stochastic Integral
$I_{(0000)T,t}^{*(i_1i_2i_3i_4)}$}

Consider now the iterated Stratonovich stochastic integral
$I_{(0000)T,t}^{*(i_1i_2i_3i_4)}$ of multiplicity 4
$(i_1, i_2, i_3, i_4=$ $1,\ldots,m)$.
For the case of pairwise different 
$i_1, i_2, i_3, i_4$ we can use the formula (\ref{destxyz1}).
In the case $i_1=i_2=i_3=i_4,$ to approximate the iterated stochastic integral
$I_{(0000)T,t}^{*(i_1i_1i_1i_1)},$ we use the formula (\ref{creat100}).

Thus, it remains to consider the following 13 cases
\begin{equation}
\label{casee1}
i_1=i_2\ne i_3, i_4;\ i_3\ne i_4,
\end{equation}

\vspace{-9.5mm}
\begin{equation}
\label{casee2}
i_1=i_3\ne i_2, i_4;\ i_2\ne i_4,
\end{equation}
\begin{equation}
\label{casee3}
i_1=i_4\ne i_2, i_3;\ i_2\ne i_3,
\end{equation}
\begin{equation}
\label{casee4}
i_2=i_3\ne i_1, i_4;\ i_1\ne i_4,
\end{equation}
\begin{equation}
\label{casee5}
i_2=i_4\ne i_1, i_3;\ i_1\ne i_3,
\end{equation}
\begin{equation}
\label{casee6}
i_3=i_4\ne i_1, i_2;\ i_1\ne i_2,
\end{equation}
\begin{equation}
\label{casee7}
i_1=i_2=i_3\ne i_4,
\end{equation}
\begin{equation}
\label{casee8}
i_2=i_3=i_4\ne i_1,
\end{equation}
\begin{equation}
\label{casee9}
i_1=i_2=i_4\ne i_3,
\end{equation}
\begin{equation}
\label{casee10}
i_1=i_3=i_4\ne i_2,
\end{equation}
\begin{equation}
\label{casee11}
i_1=i_2\ne i_3=i_4,
\end{equation}
\begin{equation}
\label{casee12}
i_1=i_3\ne i_2=i_4,
\end{equation}
\begin{equation}
\label{casee13}
i_1=i_4\ne i_2=i_3.
\end{equation}

\vspace{2mm}

By analogy with (\ref{dest10}) and using (\ref{uyes3may}), (\ref{a4}), we obtain
$$
{\sf M}\left\{\left(I_{(0000)T,t}^{*(i_1i_2i_3i_4)}-
I_{(0000)T,t}^{*(i_1i_2i_3i_4)q}\right)^2\right\}=
$$
$$
={\sf M}\left\{\Biggl(I_{(0000)T,t}^{(i_1i_2i_3i_4)}+
\frac{1}{2}{\bf 1}_{\{i_1=i_2\ne 0\}}
\int\limits_t^T\int\limits_t^{t_4}\int\limits_t^{t_3}dt_1
d{\bf w}_{t_3}^{(i_3)}
d{\bf w}_{t_4}^{(i_4)}+\Biggr.\right.
$$
$$
+\frac{1}{2}{\bf 1}_{\{i_2=i_3\ne 0\}}
\int\limits_t^T\int\limits_t^{t_4}\int\limits_t^{t_2}
d{\bf w}_{t_1}^{(i_1)}dt_2
d{\bf w}_{t_4}^{(i_4)}
+\frac{1}{2}{\bf 1}_{\{i_3=i_4\ne 0\}}
\int\limits_t^T\int\limits_t^{t_3}\int\limits_t^{t_2}
d{\bf w}_{t_1}^{(i_1)}
d{\bf w}_{t_2}^{(i_2)}dt_3+
$$
$$
+\frac{1}{4}{\bf 1}_{\{i_1=i_2\ne 0\}}
{\bf 1}_{\{i_3=i_4\ne 0\}}
\int\limits_t^T\int\limits_t^{t_2}dt_1dt_2-
I_{(0000)T,t}^{(i_1i_2i_3i_4)q}-
$$
$$
-{\bf 1}_{\{i_1=i_2\ne 0\}}\sum\limits_{j_4,j_3=0}^q
\sum\limits_{j_1=0}^q C_{j_4 j_3 j_1 j_1}\zeta_{j_3}^{(i_3)}\zeta_{j_4}^{(i_4)}-
{\bf 1}_{\{i_1=i_3\ne 0\}}\sum\limits_{j_4,j_2=0}^q
\sum\limits_{j_1=0}^q C_{j_4 j_1 j_2 j_1}\zeta_{j_2}^{(i_2)}\zeta_{j_4}^{(i_4)}-
$$
$$
-{\bf 1}_{\{i_1=i_4\ne 0\}}\sum\limits_{j_3,j_2=0}^q
\sum\limits_{j_1=0}^q C_{j_1 j_3 j_2 j_1}\zeta_{j_2}^{(i_2)}\zeta_{j_3}^{(i_3)}-
{\bf 1}_{\{i_2=i_3\ne 0\}}\sum\limits_{j_4,j_1=0}^q
\sum\limits_{j_2=0}^q C_{j_4 j_2 j_2 j_1}\zeta_{j_1}^{(i_1)}\zeta_{j_4}^{(i_4)}-
$$
$$
-{\bf 1}_{\{i_2=i_4\ne 0\}}\sum\limits_{j_3,j_1=0}^q
\sum\limits_{j_2=0}^q C_{j_2 j_3 j_2 j_1}\zeta_{j_1}^{(i_1)}\zeta_{j_3}^{(i_3)}-
{\bf 1}_{\{i_3=i_4\ne 0\}}\sum\limits_{j_2,j_1=0}^q
\sum\limits_{j_3=0}^q C_{j_3 j_3 j_2 j_1}\zeta_{j_1}^{(i_1)}\zeta_{j_2}^{(i_2)}+
$$
$$
+{\bf 1}_{\{i_1=i_2\ne 0\}}{\bf 1}_{\{i_3=i_4\ne 0\}}
\sum\limits_{j_3,j_1=0}^q C_{j_3 j_3 j_1 j_1}+
{\bf 1}_{\{i_1=i_3\ne 0\}}{\bf 1}_{\{i_2=i_4\ne 0\}}
\sum\limits_{j_2,j_1=0}^q C_{j_2 j_1 j_2 j_1}+
$$
\begin{equation}
\label{casee400}
\left.\Biggl.+{\bf 1}_{\{i_1=i_4\ne 0\}}{\bf 1}_{\{i_2=i_3\ne 0\}}
\sum\limits_{j_2,j_1=0}^q C_{j_1 j_2 j_2 j_1}\Biggr)^2\right\},
\end{equation}

\vspace{2mm}
\noindent 
where $I_{(0000)T,t}^{(i_1i_2i_3i_4)q}$ is defined by (\ref{904}).

Consider the case (\ref{casee1}). From (\ref{casee400}) we get
$$
{\sf M}\left\{\left(I_{(0000)T,t}^{*(i_1i_2i_3i_4)}-
I_{(0000)T,t}^{*(i_1i_2i_3i_4)q}\right)^2\right\}=
$$
$$
={\sf M}\left\{\Biggl(I_{(0000)T,t}^{(i_1i_2i_3i_4)}
-I_{(0000)T,t}^{(i_1i_2i_3i_4)q}+
\frac{1}{2}\int\limits_t^T\int\limits_t^{t_4}\int\limits_t^{t_3}dt_1
d{\bf w}_{t_3}^{(i_3)}
d{\bf w}_{t_4}^{(i_4)}-
\Biggr.\right.
$$
\begin{equation}
\label{casee401}
\left.\Biggl.-
\sum\limits_{j_4,j_3=0}^q
\sum\limits_{j_1=0}^q C_{j_4 j_3 j_1 j_1}\zeta_{j_3}^{(i_3)}\zeta_{j_4}^{(i_4)}
\Biggr)^2\right\}.
\end{equation}

\vspace{2mm}

Note that
\begin{equation}
\label{casee402}
\zeta_{j_3}^{(i_3)}\zeta_{j_4}^{(i_4)}=
\int\limits_t^T \phi_{j_4}(t_4)\int\limits_t^{t_4}\phi_{j_3}(t_3)
d{\bf w}_{t_3}^{(i_3)}
d{\bf w}_{t_4}^{(i_4)}+
\int\limits_t^T \phi_{j_3}(t_3)\int\limits_t^{t_3}\phi_{j_4}(t_4)
d{\bf w}_{t_4}^{(i_4)}
d{\bf w}_{t_3}^{(i_3)}
\end{equation}

\noindent
w.~p.~1, where $i_3\ne i_4.$

According to the formulas (\ref{yeee2}), (\ref{ttt2}), the quantity
$$
I_{(0000)T,t}^{(i_1i_2i_3i_4)}-
I_{(0000)T,t}^{(i_1i_2i_3i_4)q}
$$
includes only iterated It\^{o} stochastic integrals
of multiplicity 4. At the same time (see (\ref{casee402})), the quantity
$$
\frac{1}{2}\int\limits_t^T\int\limits_t^{t_4}\int\limits_t^{t_3}dt_1
d{\bf w}_{t_3}^{(i_3)}
d{\bf w}_{t_4}^{(i_4)}-
\sum\limits_{j_4,j_3=0}^q
\sum\limits_{j_1=0}^p C_{j_4 j_3 j_1 j_1}\zeta_{j_3}^{(i_3)}\zeta_{j_4}^{(i_4)}
$$
contains only iterated It\^{o} stochastic integrals
of multiplicity 2. This means that from (\ref{casee401}) we have
$$
{\sf M}\left\{\left(I_{(0000)T,t}^{*(i_1i_2i_3i_4)}-
I_{(0000)T,t}^{*(i_1i_2i_3i_4)q}\right)^2\right\}=
{\sf M}\left\{\left(I_{(0000)T,t}^{(i_1i_2i_3i_4)}
-I_{(0000)T,t}^{(i_1i_2i_3i_4)q}\right)^2\right\}
+
$$
$$
+{\sf M}\left\{\left(\frac{1}{2}\int\limits_t^T\int\limits_t^{t_4}(t_3-t)
d{\bf w}_{t_3}^{(i_3)}
d{\bf w}_{t_4}^{(i_4)}-
\sum\limits_{j_4,j_3=0}^q
\sum\limits_{j_1=0}^q C_{j_4 j_3 j_1 j_1}\zeta_{j_3}^{(i_3)}\zeta_{j_4}^{(i_4)}
\right)^2\right\}=
$$
$$
=
{\sf M}\left\{\left(I_{(0000)T,t}^{(i_1i_2i_3i_4)}
-I_{(0000)T,t}^{(i_1i_2i_3i_4)q}\right)^2\right\}
+\frac{1}{4}\int\limits_t^T\int\limits_t^{t_4}(t_3-t)^2
dt_3dt_4+
$$
$$
+\sum\limits_{j_4,j_3=0}^q
\left(\sum\limits_{j_1=0}^q C_{j_4 j_3 j_1 j_1}\right)^2- 
$$
$$
-
\sum\limits_{j_4,j_3=0}^q
\sum\limits_{j_1=0}^q C_{j_4 j_3 j_1 j_1}
\int\limits_t^T \phi_{j_4}(t_4)\int\limits_t^{t_4}\phi_{j_3}(t_3)(t_3-t)
dt_3 dt_4=
$$
$$
=
{\sf M}\left\{\left(I_{(0000)T,t}^{(i_1i_2i_3i_4)}
-I_{(0000)T,t}^{(i_1i_2i_3i_4)q}\right)^2\right\}
+\frac{(T-t)^4}{48}+
\sum\limits_{j_4,j_3=0}^q
\left(\sum\limits_{j_1=0}^q C_{j_4 j_3 j_1 j_1}\right)^2+
$$
\begin{equation}
\label{casee500}
+
\sum\limits_{j_4,j_3=0}^q
\sum\limits_{j_1=0}^q C_{j_4 j_3 j_1 j_1}C_{j_4 j_3}^{10},
\end{equation}

\noindent
where (see (\ref{caaas1}))
\begin{equation}
\label{casee700}
C_{j_4 j_3}^{10}=
\int\limits_t^T \phi_{j_4}(t_4)\int\limits_t^{t_4}\phi_{j_3}(t_3)(t-t_3)
dt_3 dt_4.
\end{equation}

Using (\ref{usl1}) and (\ref{casee500}), we finally get
$$
{\sf M}\left\{\left(I_{(0000)T,t}^{*(i_1i_2i_3i_4)}-
I_{(0000)T,t}^{*(i_1i_2i_3i_4)q}\right)^2\right\}=
\frac{(T-t)^4}{16}-
$$
$$
-
\sum_{j_1,j_2,j_3,j_4=0}^{q}
C_{j_4j_3j_2j_1}\Biggl(\sum\limits_{(j_1,j_2)}
C_{j_4j_3j_2j_1}\Biggr)
+
\sum\limits_{j_4,j_3=0}^q
\left(\sum\limits_{j_1=0}^q C_{j_4 j_3 j_1 j_1}\right)^2+
$$
\begin{equation}
\label{casee501}
+
\sum\limits_{j_4,j_3=0}^q
\sum\limits_{j_1=0}^q C_{j_4 j_3 j_1 j_1}C_{j_4 j_3}^{10},
\end{equation}

\vspace{1mm}
\noindent
where $i_1=i_2\ne i_3, i_4;\ i_3\ne i_4.$

Consider the cases (\ref{casee2}), (\ref{casee3}) by analogy
with the case (\ref{casee1}) using (\ref{usl2}), (\ref{usl3}).
We have
$$
{\sf M}\left\{\left(I_{(0000)T,t}^{*(i_1i_2i_3i_4)}-
I_{(0000)T,t}^{*(i_1i_2i_3i_4)q}\right)^2\right\}=
\frac{(T-t)^4}{24}-
$$
$$
-
\sum_{j_1,j_2,j_3,j_4=0}^{q}
C_{j_4j_3j_2j_1}\Biggl(\sum\limits_{(j_1,j_3)}
C_{j_4j_3j_2j_1}\Biggr)
+
\sum\limits_{j_4,j_2=0}^q
\left(\sum\limits_{j_1=0}^q C_{j_4 j_1 j_2 j_1}\right)^2,
$$

\vspace{1mm}
\noindent
where $i_1=i_3\ne i_2, i_4$ and $i_2\ne i_4;$
$$
{\sf M}\left\{\left(I_{(0000)T,t}^{*(i_1i_2i_3i_4)}-
I_{(0000)T,t}^{*(i_1i_2i_3i_4)q}\right)^2\right\}=
\frac{(T-t)^4}{24}-
$$
$$
-
\sum_{j_1,j_2,j_3,j_4=0}^{q}
C_{j_4j_3j_2j_1}\Biggl(\sum\limits_{(j_1,j_4)}
C_{j_4j_3j_2j_1}\Biggr)
+
\sum\limits_{j_3,j_2=0}^q
\left(\sum\limits_{j_1=0}^q C_{j_1 j_3 j_2 j_1}\right)^2,
$$

\vspace{1mm}
\noindent
where $i_1=i_4\ne i_2, i_3$ and $i_2\ne i_3.$

Consider the case (\ref{casee4}) by analogy
with the case (\ref{casee1}). 
We have
$$
{\sf M}\left\{\left(I_{(0000)T,t}^{*(i_1i_2i_3i_4)}-
I_{(0000)T,t}^{*(i_1i_2i_3i_4)q}\right)^2\right\}=
{\sf M}\left\{\left(I_{(0000)T,t}^{(i_1i_2i_3i_4)}
-I_{(0000)T,t}^{(i_1i_2i_3i_4)q}\right)^2\right\}
+
$$
$$
+{\sf M}\left\{\left(\frac{1}{2}\int\limits_t^T\int\limits_t^{t_4}\int\limits_t^{t_2}
d{\bf w}_{t_1}^{(i_1)}dt_2
d{\bf w}_{t_4}^{(i_4)}-
\sum\limits_{j_4,j_1=0}^q
\sum\limits_{j_2=0}^q C_{j_4 j_2 j_2 j_1}\zeta_{j_1}^{(i_1)}\zeta_{j_4}^{(i_4)}
\right)^2\right\}=
$$
$$
=
{\sf M}\left\{\left(I_{(0000)T,t}^{(i_1i_2i_3i_4)}
-I_{(0000)T,t}^{(i_1i_2i_3i_4)q}\right)^2\right\}
+
$$
$$
+{\sf M}\left\{\left(\frac{1}{2}\int\limits_t^T\int\limits_t^{t_4}(t_4-t_1)
d{\bf w}_{t_1}^{(i_1)}
d{\bf w}_{t_4}^{(i_4)}-
\sum\limits_{j_4,j_1=0}^q
\sum\limits_{j_2=0}^q C_{j_4 j_2 j_2 j_1}\zeta_{j_1}^{(i_1)}\zeta_{j_4}^{(i_4)}
\right)^2\right\}=
$$
$$
=
{\sf M}\left\{\left(I_{(0000)T,t}^{(i_1i_2i_3i_4)}
-I_{(0000)T,t}^{(i_1i_2i_3i_4)q}\right)^2\right\}
+\frac{(T-t)^4}{48}+
\sum\limits_{j_4,j_1=0}^q
\left(\sum\limits_{j_2=0}^q C_{j_4 j_2 j_2 j_1}\right)^2-
$$
$$
-
\sum\limits_{j_4,j_1=0}^q
\sum\limits_{j_2=0}^q C_{j_4 j_2 j_2 j_1}
\int\limits_t^T \phi_{j_4}(t_4)\int\limits_t^{t_4}\phi_{j_1}(t_1)(t_4-t_1)
dt_3 dt_4.
$$

\vspace{3mm}

Then using (\ref{usl4}), we obtain
$$
{\sf M}\left\{\left(I_{(0000)T,t}^{*(i_1i_2i_3i_4)}-
I_{(0000)T,t}^{*(i_1i_2i_3i_4)q}\right)^2\right\}=\frac{(T-t)^4}{16}-
$$
$$
-
\sum_{j_1,j_2,j_3,j_4=0}^{q}
C_{j_4j_3j_2j_1}\Biggl(\sum\limits_{(j_2,j_3)}
C_{j_4j_3j_2j_1}\Biggr)+
$$
$$
+
\sum\limits_{j_4,j_1=0}^q
\left(\sum\limits_{j_2=0}^q C_{j_4 j_2 j_2 j_1}\right)^2-
\sum\limits_{j_4,j_1=0}^q
\sum\limits_{j_2=0}^q C_{j_4 j_2 j_2 j_1}
\left(C_{j_4 j_1}^{10}-C_{j_4 j_1}^{01}\right),
$$

\vspace{1mm}
\noindent
where $i_2=i_3\ne i_1, i_4$ and $i_1\ne i_4;$
$C_{j_4 j_1}^{10}$ is defined by (\ref{casee700}) and
\begin{equation}
\label{casee700a}
C_{j_4 j_1}^{01}=
\int\limits_t^T \phi_{j_4}(t_4)(t-t_4)\int\limits_t^{t_4}\phi_{j_1}(t_1)
dt_1 dt_4.
\end{equation}

\vspace{2mm}

For the case (\ref{casee5}) by analogy
with the case (\ref{casee1}) and using (\ref{usl5}), we get
$$
{\sf M}\left\{\left(I_{(0000)T,t}^{*(i_1i_2i_3i_4)}-
I_{(0000)T,t}^{*(i_1i_2i_3i_4)q}\right)^2\right\}=
\frac{(T-t)^4}{24}-
$$
$$
-
\sum_{j_1,j_2,j_3,j_4=0}^{q}
C_{j_4j_3j_2j_1}\Biggl(\sum\limits_{(j_2,j_4)}
C_{j_4j_3j_2j_1}\Biggr)
+
\sum\limits_{j_3,j_1=0}^q
\left(\sum\limits_{j_2=0}^q C_{j_2 j_3 j_2 j_1}\right)^2,
$$

\vspace{1mm}
\noindent
where $i_2=i_4\ne i_1, i_3$ and $i_1\ne i_3.$

Consider the case (\ref{casee6}) by analogy
with the case (\ref{casee1}). 
We have (see Example~3.1 in Sect.~3.6)
$$
{\sf M}\left\{\left(I_{(0000)T,t}^{*(i_1i_2i_3i_4)}-
I_{(0000)T,t}^{*(i_1i_2i_3i_4)q}\right)^2\right\}=
{\sf M}\left\{\left(I_{(0000)T,t}^{(i_1i_2i_3i_4)}
-I_{(0000)T,t}^{(i_1i_2i_3i_4)q}\right)^2\right\}
+
$$
$$
+{\sf M}\left\{\left(\frac{1}{2}\int\limits_t^T(T-t_2)\int\limits_t^{t_2}
d{\bf w}_{t_1}^{(i_1)}
d{\bf w}_{t_2}^{(i_2)}-
\sum\limits_{j_2,j_1=0}^q
\sum\limits_{j_3=0}^q C_{j_3 j_3 j_2 j_1}\zeta_{j_1}^{(i_1)}\zeta_{j_2}^{(i_2)}
\right)^2\right\}=
$$
$$
=
{\sf M}\left\{\left(I_{(0000)T,t}^{(i_1i_2i_3i_4)}
-I_{(0000)T,t}^{(i_1i_2i_3i_4)q}\right)^2\right\}
+\frac{(T-t)^4}{48}+
\sum\limits_{j_2,j_1=0}^q
\left(\sum\limits_{j_3=0}^q C_{j_3 j_3 j_2 j_1}\right)^2-
$$
$$
-
\sum\limits_{j_2,j_1=0}^q
\sum\limits_{j_3=0}^q C_{j_3 j_3 j_2 j_1}
\int\limits_t^T (T-t_2)\phi_{j_2}(t_2)\int\limits_t^{t_2}\phi_{j_1}(t_1)
dt_1 dt_2.
$$

\vspace{3mm}

Then using (\ref{usl6}), we obtain
$$
{\sf M}\left\{\left(I_{(0000)T,t}^{*(i_1i_2i_3i_4)}-
I_{(0000)T,t}^{*(i_1i_2i_3i_4)q}\right)^2\right\}=\frac{(T-t)^4}{16}-
$$
$$
-
\sum_{j_1,j_2,j_3,j_4=0}^{q}
C_{j_4j_3j_2j_1}\Biggl(\sum\limits_{(j_3,j_4)}
C_{j_4j_3j_2j_1}\Biggr)+
$$
$$
+
\sum\limits_{j_2,j_1=0}^q
\left(\sum\limits_{j_3=0}^q C_{j_3 j_3 j_2 j_1}\right)^2-
\sum\limits_{j_2,j_1=0}^q
\sum\limits_{j_3=0}^q C_{j_3 j_3 j_2 j_1}
\left((T-t)C_{j_2 j_1}+C_{j_2 j_1}^{01}\right),
$$

\vspace{1mm}
\noindent
where $i_3=i_4\ne i_1, i_2$ and $i_1\ne i_2;$
$C_{j_2 j_1}^{01}$ is defined by (\ref{casee700a}) and
$$
C_{j_2 j_1}=
\int\limits_t^T \phi_{j_2}(t_2)\int\limits_t^{t_2}\phi_{j_1}(t_1)
dt_1 dt_2.
$$

\vspace{2mm}

Consider the case (\ref{casee7}). 
From (\ref{casee400}) we have 
$$
{\sf M}\left\{\left(I_{(0000)T,t}^{*(i_1i_1i_1i_4)}-
I_{(0000)T,t}^{*(i_1i_1i_1i_4)q}\right)^2\right\}=
{\sf M}\left\{\Biggl(I_{(0000)T,t}^{(i_1i_1i_1i_4)}+
\frac{1}{2}
\int\limits_t^T\int\limits_t^{t_4}\int\limits_t^{t_3}dt_1
d{\bf w}_{t_3}^{(i_1)}
d{\bf w}_{t_4}^{(i_4)}+\Biggr.\right.
$$
$$
+\frac{1}{2}
\int\limits_t^T\int\limits_t^{t_4}\int\limits_t^{t_2}
d{\bf w}_{t_1}^{(i_1)}dt_2
d{\bf w}_{t_4}^{(i_4)}-
I_{(0000)T,t}^{(i_1i_1i_1i_4)q}-
\sum\limits_{j_4,j_3=0}^q
\sum\limits_{j_1=0}^q C_{j_4 j_3 j_1 j_1}\zeta_{j_3}^{(i_1)}\zeta_{j_4}^{(i_4)}-
$$
\begin{equation}
\label{casee800}
~~~~~~~ \left.\Biggl.-
\sum\limits_{j_4,j_2=0}^q
\sum\limits_{j_1=0}^q C_{j_4 j_1 j_2 j_1}\zeta_{j_2}^{(i_1)}\zeta_{j_4}^{(i_4)}-
\sum\limits_{j_4,j_1=0}^q
\sum\limits_{j_2=0}^q C_{j_4 j_2 j_2 j_1}\zeta_{j_1}^{(i_1)}\zeta_{j_4}^{(i_4)}
\Biggr)^2\right\}.
\end{equation}

\vspace{2mm}

Furthermore,
$$
\int\limits_t^T\int\limits_t^{t_4}\int\limits_t^{t_3}dt_1
d{\bf w}_{t_3}^{(i_1)}
d{\bf w}_{t_4}^{(i_4)}+
\int\limits_t^T\int\limits_t^{t_4}\int\limits_t^{t_2}
d{\bf w}_{t_1}^{(i_1)}dt_2
d{\bf w}_{t_4}^{(i_4)}=
$$
$$
=\int\limits_t^T\int\limits_t^{t_4}(t_1-t)
d{\bf w}_{t_1}^{(i_1)}
d{\bf w}_{t_4}^{(i_4)}+
\int\limits_t^T\int\limits_t^{t_4}
(t_4-t_1)d{\bf w}_{t_1}^{(i_1)}
d{\bf w}_{t_4}^{(i_4)}=
$$
\begin{equation}
\label{casee801}
=\int\limits_t^T(t_4-t)\int\limits_t^{t_4}
d{\bf w}_{t_1}^{(i_1)}
d{\bf w}_{t_4}^{(i_4)}\ \ \ \hbox{w.~p.~1.}
\end{equation}

\vspace{1mm}

From (\ref{casee800}) and (\ref{casee801}) we get
$$
{\sf M}\left\{\left(I_{(0000)T,t}^{*(i_1i_1i_1i_4)}-
I_{(0000)T,t}^{*(i_1i_1i_1i_4)q}\right)^2\right\}=
{\sf M}\left\{\left(I_{(0000)T,t}^{(i_1i_1i_1i_4)}-
I_{(0000)T,t}^{(i_1i_1i_1i_4)q}\right)^2\right\}+
$$
$$
+{\sf M}\left\{\Biggl(
\frac{1}{2}\int\limits_t^T(t_4-t)\int\limits_t^{t_4}
d{\bf w}_{t_1}^{(i_1)}
d{\bf w}_{t_4}^{(i_4)}-\Biggr.\right.
$$
$$
\left.\Biggl.-
\sum\limits_{j_4,j_1=0}^q
\sum\limits_{j_2=0}^q \left(C_{j_4 j_1 j_2 j_2}+C_{j_4 j_2 j_1 j_2}+
C_{j_4 j_2 j_2 j_1}\right)\zeta_{j_1}^{(i_1)}\zeta_{j_4}^{(i_4)}
\Biggr)^2\right\}=
$$
$$
=
{\sf M}\left\{\left(I_{(0000)T,t}^{(i_1i_1i_1i_4)}-
I_{(0000)T,t}^{(i_1i_1i_1i_4)q}\right)^2\right\}+\frac{(T-t)^4}{16}+
$$
$$
+
\sum\limits_{j_4,j_1=0}^q
\left(\sum\limits_{j_2=0}^q \left(C_{j_4 j_1 j_2 j_2}+C_{j_4 j_2 j_1 j_2}+
C_{j_4 j_2 j_2 j_1}\right)\right)^2-
$$
\begin{equation}
\label{casee805}
-\sum\limits_{j_4,j_1=0}^q
\sum\limits_{j_2=0}^q \left(C_{j_4 j_1 j_2 j_2}+C_{j_4 j_2 j_1 j_2}+
C_{j_4 j_2 j_2 j_1}\right)
\int\limits_t^T (t_4-t)\phi_{j_4}(t_4)\int\limits_t^{t_4}\phi_{j_1}(t_1)
dt_1 dt_4.
\end{equation}

\vspace{1mm}

Using (\ref{usl7}) and (\ref{casee805}), we finally obtain
$$
{\sf M}\left\{\left(I_{(0000)T,t}^{*(i_1i_2i_3i_4)}-
I_{(0000)T,t}^{*(i_1i_2i_3i_4)q}\right)^2\right\}=
\frac{5(T-t)^4}{48}-
$$
$$
-
\sum_{j_1,j_2,j_3,j_4=0}^{q}
C_{j_4j_3j_2j_1}\Biggl(\sum\limits_{(j_1,j_2,j_3)}
C_{j_4j_3j_2j_1}\Biggr)+
$$
$$
+
\sum\limits_{j_4,j_1=0}^q
\left(\sum\limits_{j_2=0}^q \left(C_{j_4 j_1 j_2 j_2}+C_{j_4 j_2 j_1 j_2}+
C_{j_4 j_2 j_2 j_1}\right)\right)^2+
$$
$$
+\sum\limits_{j_4,j_1=0}^q
\sum\limits_{j_2=0}^q \left(C_{j_4 j_1 j_2 j_2}+C_{j_4 j_2 j_1 j_2}+
C_{j_4 j_2 j_2 j_1}\right)C_{j_4j_2}^{01},
$$

\vspace{1mm}
\noindent
where $i_1=i_2=i_3\ne i_4.$

Consider the case (\ref{casee8}). 
From (\ref{casee400}) we have 
$$
{\sf M}\left\{\left(I_{(0000)T,t}^{*(i_1i_2i_2i_2)}-
I_{(0000)T,t}^{*(i_1i_2i_2i_2)q}\right)^2\right\}=
{\sf M}\left\{\Biggl(I_{(0000)T,t}^{(i_1i_2i_2i_2)}+
\frac{1}{2}
\int\limits_t^T\int\limits_t^{t_4}\int\limits_t^{t_2}
d{\bf w}_{t_1}^{(i_1)}dt_2
d{\bf w}_{t_4}^{(i_2)}+\Biggr.\right.
$$
$$
+\frac{1}{2}
\int\limits_t^T\int\limits_t^{t_3}\int\limits_t^{t_2}
d{\bf w}_{t_1}^{(i_1)}
d{\bf w}_{t_2}^{(i_2)}dt_3-
I_{(0000)T,t}^{(i_1i_2i_2i_2)q}-
\sum\limits_{j_4,j_1=0}^q
\sum\limits_{j_2=0}^q C_{j_4 j_2 j_2 j_1}\zeta_{j_1}^{(i_1)}\zeta_{j_4}^{(i_2)}-
$$
\begin{equation}
\label{casee807}
~~~~~~~ \left.\Biggl.-
\sum\limits_{j_3,j_1=0}^q
\sum\limits_{j_2=0}^q C_{j_2 j_3 j_2 j_1}\zeta_{j_1}^{(i_1)}\zeta_{j_3}^{(i_2)}-
\sum\limits_{j_2,j_1=0}^q
\sum\limits_{j_3=0}^q C_{j_3 j_3 j_2 j_1}\zeta_{j_1}^{(i_1)}\zeta_{j_2}^{(i_2)}
\Biggr)^2\right\}.
\end{equation}

\vspace{2mm}

Moreover,
$$
\int\limits_t^T\int\limits_t^{t_4}\int\limits_t^{t_2}
d{\bf w}_{t_1}^{(i_1)}dt_2
d{\bf w}_{t_4}^{(i_2)}+
\int\limits_t^T\int\limits_t^{t_3}\int\limits_t^{t_2}
d{\bf w}_{t_1}^{(i_1)}
d{\bf w}_{t_2}^{(i_2)}dt_3=
$$
$$
=\int\limits_t^T\int\limits_t^{t_4}(t_4-t_1)
d{\bf w}_{t_1}^{(i_1)}
d{\bf w}_{t_4}^{(i_2)}+
\int\limits_t^T\int\limits_t^{t_4}
(T-t_4)d{\bf w}_{t_1}^{(i_1)}
d{\bf w}_{t_4}^{(i_2)}=
$$
\begin{equation}
\label{casee808}
=\int\limits_t^T\int\limits_t^{t_4}
(T-t_1)d{\bf w}_{t_1}^{(i_1)}
d{\bf w}_{t_4}^{(i_2)}\ \ \ \hbox{w.~p.~1.}
\end{equation}

\vspace{1mm}

From (\ref{casee807}) and (\ref{casee808}) we get
$$
{\sf M}\left\{\left(I_{(0000)T,t}^{*(i_1i_2i_2i_2)}-
I_{(0000)T,t}^{*(i_1i_2i_2i_2)q}\right)^2\right\}=
{\sf M}\left\{\left(I_{(0000)T,t}^{(i_1i_2i_2i_2)}-
I_{(0000)T,t}^{(i_1i_2i_2i_2)q}\right)^2\right\}+
$$
$$
+{\sf M}\left\{\Biggl(
\frac{1}{2}\int\limits_t^T\int\limits_t^{t_4}(T-t_1)
d{\bf w}_{t_1}^{(i_1)}
d{\bf w}_{t_4}^{(i_2)}-\Biggr.\right.
$$
$$
\left.\Biggl.-
\sum\limits_{j_4,j_1=0}^q
\sum\limits_{j_2=0}^q \left(C_{j_4 j_2 j_2 j_1}+C_{j_2 j_4 j_2 j_1}+
C_{j_2 j_2 j_4 j_1}\right)\zeta_{j_1}^{(i_1)}\zeta_{j_4}^{(i_2)}
\Biggr)^2\right\}=
$$
$$
=
{\sf M}\left\{\left(I_{(0000)T,t}^{(i_1i_2i_2i_2)}-
I_{(0000)T,t}^{(i_1i_2i_2i_2)q}\right)^2\right\}+\frac{(T-t)^4}{16}+
$$
$$
+
\sum\limits_{j_4,j_1=0}^q
\left(\sum\limits_{j_2=0}^q \left(C_{j_4 j_2 j_2 j_1}+C_{j_2 j_4 j_2 j_1}+
C_{j_2 j_2 j_4 j_1}\right)\right)^2-
$$
\begin{equation}
\label{casee809}
-\sum\limits_{j_4,j_1=0}^q
\sum\limits_{j_2=0}^q \left(C_{j_4 j_2 j_2 j_1}+C_{j_2 j_4 j_2 j_1}+
C_{j_2 j_2 j_4 j_1}\right)
\int\limits_t^T \phi_{j_4}(t_4)\int\limits_t^{t_4}(T-t_1)\phi_{j_1}(t_1)
dt_1 dt_4.
\end{equation}

\vspace{1mm}

Using (\ref{usl8}) and (\ref{casee809}), we finally obtain
$$
{\sf M}\left\{\left(I_{(0000)T,t}^{*(i_1i_2i_3i_4)}-
I_{(0000)T,t}^{*(i_1i_2i_3i_4)q}\right)^2\right\}=
\frac{5(T-t)^4}{48}-
$$
$$
-
\sum_{j_1,j_2,j_3,j_4=0}^{q}
C_{j_4j_3j_2j_1}\Biggl(\sum\limits_{(j_2,j_3,j_4)}
C_{j_4j_3j_2j_1}\Biggr)+
$$
$$
+
\sum\limits_{j_4,j_1=0}^q
\left(\sum\limits_{j_2=0}^q \left(C_{j_4 j_2 j_2 j_1}+C_{j_2 j_4 j_2 j_1}+
C_{j_2 j_2 j_4 j_1}\right)\right)^2-
$$
$$
-\sum\limits_{j_4,j_1=0}^q
\sum\limits_{j_2=0}^q \left(C_{j_4 j_2 j_2 j_1}+C_{j_2 j_4 j_2 j_1}+
C_{j_2 j_2 j_4 j_1}\right)\left((T-t)C_{j_4 j_1}+C_{j_4 j_1}^{10}\right),
$$

\vspace{1mm}
\noindent
where $i_2=i_3=i_4\ne i_1.$

For the cases (\ref{casee9}), (\ref{casee10}) by analogy with the case (\ref{casee8}) and using
(\ref{usl9}), (\ref{usl10}), we get
$$
{\sf M}\left\{\left(I_{(0000)T,t}^{*(i_1i_2i_3i_4)}-
I_{(0000)T,t}^{*(i_1i_2i_3i_4)q}\right)^2\right\}=
\frac{(T-t)^4}{16}-
$$
$$
-
\sum_{j_1,j_2,j_3,j_4=0}^{q}
C_{j_4j_3j_2j_1}\Biggl(\sum\limits_{(j_1,j_2,j_4)}
C_{j_4j_3j_2j_1}\Biggr)+
$$
$$
+
\sum\limits_{j_4,j_3=0}^q
\left(\sum\limits_{j_1=0}^q \left(C_{j_4 j_3 j_1 j_1}+C_{j_1 j_3 j_4 j_1}+
C_{j_1 j_3 j_1 j_4}\right)\right)^2+
$$
$$
+\sum\limits_{j_4,j_3=0}^q
\sum\limits_{j_1=0}^q \left(C_{j_4 j_3 j_1 j_1}+C_{j_1 j_3 j_4 j_1}+
C_{j_1 j_3 j_1 j_4}\right)C_{j_4 j_3}^{10},
$$

\vspace{1mm}
\noindent
where $i_1=i_2=i_4\ne i_3;$
$$
{\sf M}\left\{\left(I_{(0000)T,t}^{*(i_1i_2i_3i_4)}-
I_{(0000)T,t}^{*(i_1i_2i_3i_4)q}\right)^2\right\}=
\frac{(T-t)^4}{16}-
$$
$$
-
\sum_{j_1,j_2,j_3,j_4=0}^{q}
C_{j_4j_3j_2j_1}\Biggl(\sum\limits_{(j_1,j_3,j_4)}
C_{j_4j_3j_2j_1}\Biggr)+
$$
$$
+
\sum\limits_{j_4,j_2=0}^q
\left(\sum\limits_{j_1=0}^q \left(C_{j_4 j_1 j_2 j_1}+C_{j_1 j_4 j_2 j_1}+
C_{j_1 j_1 j_2 j_4}\right)\right)^2-
$$
$$
-\sum\limits_{j_4,j_2=0}^q
\sum\limits_{j_1=0}^q \left(C_{j_4 j_1 j_2 j_1}+C_{j_1 j_4 j_2 j_1}+
C_{j_1 j_1 j_2 j_4}\right)\left((T-t)C_{j_2 j_3}+C_{j_2 j_3}^{01}\right),
$$

\vspace{1mm}
\noindent
where $i_1=i_3=i_4\ne i_2.$

Let us consider the case (\ref{casee11}). Using (\ref{casee400}), we obtain
$$
{\sf M}\left\{\left(I_{(0000)T,t}^{*(i_1i_1i_3i_3)}-
I_{(0000)T,t}^{*(i_1i_1i_3i_3)q}\right)^2\right\}=
$$
$$
=
{\sf M}\left\{\Biggl(I_{(0000)T,t}^{(i_1i_1i_3i_3)}+
\frac{1}{2}
\int\limits_t^T\int\limits_t^{t_4}(t_3-t)
d{\bf w}_{t_3}^{(i_3)}
d{\bf w}_{t_4}^{(i_3)}+\Biggr.\right.
$$
$$
+\frac{1}{2}
\int\limits_t^T\int\limits_t^{t_3}\int\limits_t^{t_2}
d{\bf w}_{t_1}^{(i_1)}
d{\bf w}_{t_2}^{(i_1)}dt_3+ \frac{(T-t)^2}{8}-
I_{(0000)T,t}^{(i_1i_1i_3i_3)q}-
$$
$$
-\sum\limits_{j_4,j_3=0}^q
\sum\limits_{j_1=0}^q C_{j_4 j_3 j_1 j_1}\zeta_{j_3}^{(i_3)}\zeta_{j_4}^{(i_3)}-
\sum\limits_{j_2,j_1=0}^q
\sum\limits_{j_3=0}^q C_{j_3 j_3 j_2 j_1}\zeta_{j_1}^{(i_1)}\zeta_{j_2}^{(i_1)}+
\left.\Biggl.
\sum\limits_{j_3,j_1=0}^q C_{j_3 j_3 j_1 j_1}\Biggr)^2\right\}=
$$
$$
=
{\sf M}\left\{\Biggl(I_{(0000)T,t}^{(i_1i_1i_3i_3)}-
I_{(0000)T,t}^{(i_1i_1i_3i_3)q}+
\frac{1}{2}
\int\limits_t^T\int\limits_t^{t_4}(t_3-t)
d{\bf w}_{t_3}^{(i_3)}
d{\bf w}_{t_4}^{(i_3)}-\Biggr.\right.
$$
$$
-
\sum\limits_{j_4,j_3=0}^q
\sum\limits_{j_1=0}^q C_{j_4 j_3 j_1 j_1}\left(\zeta_{j_3}^{(i_3)}\zeta_{j_4}^{(i_3)}-
{\bf 1}_{\{j_3=j_4\}}\right)+
$$
$$
+\frac{1}{2}
\int\limits_t^T\int\limits_t^{t_3}\int\limits_t^{t_2}
d{\bf w}_{t_1}^{(i_1)}
d{\bf w}_{t_2}^{(i_1)}dt_3-
\sum\limits_{j_2,j_1=0}^q
\sum\limits_{j_3=0}^q C_{j_3 j_3 j_2 j_1}\left(\zeta_{j_1}^{(i_1)}\zeta_{j_2}^{(i_1)}-
{\bf 1}_{\{j_1=j_2\}}\right)+
$$
\begin{equation}
\label{caseu1}
\left.\Biggl.+\frac{(T-t)^2}{8}-
\sum\limits_{j_3,j_1=0}^q C_{j_3 j_3 j_1 j_1}\Biggr)^2\right\}.
\end{equation}

\vspace{2mm}

Note that
$$
\zeta_{j_3}^{(i_3)}\zeta_{j_4}^{(i_3)}-
{\bf 1}_{\{j_3=j_4\}}=
$$
\begin{equation}
\label{caseu2}
=\int\limits_t^T \phi_{j_4}(t_4)\int\limits_t^{t_4}\phi_{j_3}(t_3)
d{\bf w}_{t_3}^{(i_3)}
d{\bf w}_{t_4}^{(i_3)}+
\int\limits_t^T \phi_{j_3}(t_3)\int\limits_t^{t_3}\phi_{j_4}(t_4)
d{\bf w}_{t_4}^{(i_3)}
d{\bf w}_{t_3}^{(i_3)},
\end{equation}

$$
\zeta_{j_1}^{(i_1)}\zeta_{j_2}^{(i_1)}-
{\bf 1}_{\{j_1=j_2\}}=
$$
\begin{equation}
\label{caseu3}
=\int\limits_t^T \phi_{j_2}(t_2)\int\limits_t^{t_2}\phi_{j_1}(t_1)
d{\bf w}_{t_1}^{(i_1)}
d{\bf w}_{t_2}^{(i_1)}+
\int\limits_t^T \phi_{j_1}(t_1)\int\limits_t^{t_1}\phi_{j_2}(t_2)
d{\bf w}_{t_2}^{(i_1)}
d{\bf w}_{t_1}^{(i_1)}
\end{equation}

\noindent
w.~p.~1.

The relations (\ref{caseu1})--(\ref{caseu3}) and Example~3.1 in Sect.~3.6 imply the following 
$$
{\sf M}\left\{\left(I_{(0000)T,t}^{*(i_1i_1i_3i_3)}-
I_{(0000)T,t}^{*(i_1i_1i_3i_3)q}\right)^2\right\}=
{\sf M}\left\{\left(I_{(0000)T,t}^{(i_1i_1i_3i_3)}-
I_{(0000)T,t}^{(i_1i_1i_3i_3)q}\right)^2\right\}+
$$
$$
+
{\sf M}\left\{\Biggl(\frac{1}{2}
\int\limits_t^T\int\limits_t^{t_4}(t_3-t)
d{\bf w}_{t_3}^{(i_3)}
d{\bf w}_{t_4}^{(i_3)}-\Biggr.\right.
$$
$$
\left.\Biggl.-
\sum\limits_{j_4,j_3=0}^q
\sum\limits_{j_1=0}^q C_{j_4 j_3 j_1 j_1}\left(\zeta_{j_3}^{(i_3)}\zeta_{j_4}^{(i_3)}-
{\bf 1}_{\{j_3=j_4\}}\right)\Biggr)^2\right\}+
$$
$$
+{\sf M}\left\{\Biggl(\frac{1}{2}
\int\limits_t^T\int\limits_t^{t_3}\int\limits_t^{t_2}
d{\bf w}_{t_1}^{(i_1)}
d{\bf w}_{t_2}^{(i_1)}dt_3-\Biggr.\right.
$$
$$
\left.\Biggl.-
\sum\limits_{j_2,j_1=0}^q
\sum\limits_{j_3=0}^q C_{j_3 j_3 j_2 j_1}\left(\zeta_{j_1}^{(i_1)}\zeta_{j_2}^{(i_1)}-
{\bf 1}_{\{j_1=j_2\}}\right)\Biggr)^2\right\}+
$$
$$
+\left(\frac{(T-t)^2}{8}-
\sum\limits_{j_3,j_1=0}^q C_{j_3 j_3 j_1 j_1}\right)^2=
$$
$$
=
{\sf M}\left\{\left(I_{(0000)T,t}^{(i_1i_1i_3i_3)}-
I_{(0000)T,t}^{(i_1i_1i_3i_3)q}\right)^2\right\}+
{\sf M}\left\{\Biggl(\frac{1}{2}
\int\limits_t^T\int\limits_t^{t_4}(t_3-t)
d{\bf w}_{t_3}^{(i_3)}
d{\bf w}_{t_4}^{(i_3)}-\Biggr.\right.
$$
$$
\left.\Biggl.-
\sum\limits_{j_4,j_3=0}^q
\sum\limits_{j_1=0}^q C_{j_4 j_3 j_1 j_1}\left(\zeta_{j_3}^{(i_3)}\zeta_{j_4}^{(i_3)}-
{\bf 1}_{\{j_3=j_4\}}\right)\Biggr)^2\right\}+
$$
$$
+{\sf M}\left\{\Biggl(\frac{1}{2}
\int\limits_t^T(T-t_2)\int\limits_t^{t_2}
d{\bf w}_{t_1}^{(i_1)}
d{\bf w}_{t_2}^{(i_1)}-\Biggr.\right.
$$
$$
\left.\Biggl.-
\sum\limits_{j_2,j_1=0}^q
\sum\limits_{j_3=0}^q C_{j_3 j_3 j_2 j_1}\left(\zeta_{j_1}^{(i_1)}\zeta_{j_2}^{(i_1)}-
{\bf 1}_{\{j_1=j_2\}}\right)\Biggr)^2\right\}+
$$
$$
+\left(\frac{(T-t)^2}{8}-
\sum\limits_{j_3,j_1=0}^q C_{j_3 j_3 j_1 j_1}\right)^2=
$$
$$
=
{\sf M}\left\{\left(I_{(0000)T,t}^{(i_1i_1i_3i_3)}-
I_{(0000)T,t}^{(i_1i_1i_3i_3)q}\right)^2\right\}+
\frac{(T-t)^4}{48}
+\sum\limits_{j_4,j_3=0}^q
\sum\limits_{j_1=0}^q C_{j_4 j_3 j_1 j_1}\left(C_{j_3j_4}^{10}+C_{j_4j_3}^{10}\right)+
$$
$$
+
{\sf M}\left\{\left(
\sum\limits_{j_4,j_3=0}^q
\sum\limits_{j_1=0}^q C_{j_4 j_3 j_1 j_1}\left(\zeta_{j_3}^{(i_3)}\zeta_{j_4}^{(i_3)}-
{\bf 1}_{\{j_3=j_4\}}\right)\right)^2\right\}+
$$
$$
+
\frac{(T-t)^4}{48}
-\sum\limits_{j_2,j_1=0}^q
\sum\limits_{j_3=0}^q C_{j_3 j_3 j_2 j_1}\left((T-t)\left(C_{j_1j_2}+C_{j_2j_1}\right)+
C_{j_1j_2}^{01}+C_{j_2j_1}^{01}\right)+
$$
$$
+
{\sf M}\left\{\left(
\sum\limits_{j_2,j_1=0}^q
\sum\limits_{j_3=0}^q C_{j_3 j_3 j_2 j_1}\left(\zeta_{j_1}^{(i_1)}\zeta_{j_2}^{(i_1)}-
{\bf 1}_{\{j_1=j_2\}}\right)\right)^2\right\}+
$$
\begin{equation}
\label{caseu100}
+\left(\frac{(T-t)^2}{8}-
\sum\limits_{j_3,j_1=0}^q C_{j_3 j_3 j_1 j_1}\right)^2.
\end{equation}

Furthermore,
$$
{\sf M}\left\{\left(
\sum\limits_{j_4,j_3=0}^q
\sum\limits_{j_1=0}^q C_{j_4 j_3 j_1 j_1}\left(\zeta_{j_3}^{(i_3)}\zeta_{j_4}^{(i_3)}-
{\bf 1}_{\{j_3=j_4\}}\right)\right)^2\right\}=
$$
$$
={\sf M}\left\{\left(
\sum\limits_{j_4,j_3=0}^q
\sum\limits_{j_1=0}^q C_{j_4 j_3 j_1 j_1}\zeta_{j_3}^{(i_3)}\zeta_{j_4}^{(i_3)}
\right)^2\right\}-
2\left(\sum\limits_{j_3,j_1=0}^q C_{j_3 j_3 j_1 j_1}\right)^2+
$$
$$
+
\left(\sum\limits_{j_3,j_1=0}^q C_{j_3 j_3 j_1 j_1}\right)^2=
$$
\begin{equation}
\label{caseu101}
~~~~~~~ ={\sf M}\left\{\left(
\sum\limits_{j_4,j_3=0}^q
\sum\limits_{j_1=0}^q C_{j_4 j_3 j_1 j_1}\zeta_{j_3}^{(i_3)}\zeta_{j_4}^{(i_3)}
\right)^2\right\}-
\left(\sum\limits_{j_3,j_1=0}^q C_{j_3 j_3 j_1 j_1}\right)^2,
\end{equation}

\vspace{1mm}
$$
{\sf M}\left\{\left(
\sum\limits_{j_2,j_1=0}^q
\sum\limits_{j_3=0}^q C_{j_3 j_3 j_2 j_1}\left(\zeta_{j_1}^{(i_1)}\zeta_{j_2}^{(i_1)}-
{\bf 1}_{\{j_1=j_2\}}\right)\right)^2\right\}=
$$
\begin{equation}
\label{caseu102}
~~~~~~~ ={\sf M}\left\{\left(
\sum\limits_{j_2,j_1=0}^q
\sum\limits_{j_3=0}^q C_{j_3 j_3 j_2 j_1}\zeta_{j_1}^{(i_1)}\zeta_{j_2}^{(i_1)}
\right)^2\right\}-
\left(
\sum\limits_{j_1,j_3=0}^q C_{j_3 j_3 j_1 j_1}\right)^2.
\end{equation}

\vspace{2mm}

Using (\ref{otit321}), we get
$$
{\sf M}\left\{\left(
\sum\limits_{j_4,j_3=0}^q
\sum\limits_{j_1=0}^q C_{j_4 j_3 j_1 j_1}\zeta_{j_3}^{(i_3)}\zeta_{j_4}^{(i_3)}
\right)^2\right\}=
$$
$$
=\left(\sum\limits_{j_3,j_1=0}^q C_{j_3 j_3 j_1 j_1}\right)^2+
\sum\limits_{j_4=0}^q \sum\limits_{j_3=0}^{j_4-1}
\left(\sum\limits_{j_1=0}^q C_{j_3 j_4 j_1 j_1}+ \sum\limits_{j_1=0}^q C_{j_4 j_3 j_1 j_1}
\right)^2+
$$
\begin{equation}
\label{caseu103}
+2\sum\limits_{j_4=0}^q\left(\sum\limits_{j_1=0}^q C_{j_4 j_4 j_1 j_1}\right)^2.
\end{equation}

\vspace{2mm}

From (\ref{caseu101}) and (\ref{caseu103}) we have
$$
{\sf M}\left\{\left(
\sum\limits_{j_4,j_3=0}^q
\sum\limits_{j_1=0}^q C_{j_4 j_3 j_1 j_1}\left(\zeta_{j_3}^{(i_3)}\zeta_{j_4}^{(i_3)}-
{\bf 1}_{\{j_3=j_4\}}\right)\right)^2\right\}=
$$
\begin{equation}
\label{caseu104}
~~~~~ =\sum\limits_{j_4=0}^q \sum\limits_{j_3=0}^{j_4-1}
\left(\sum\limits_{j_1=0}^q C_{j_3 j_4 j_1 j_1}+ 
\sum\limits_{j_1=0}^q C_{j_4 j_3 j_1 j_1}\right)^2
+2\sum\limits_{j_4=0}^q\left(\sum\limits_{j_1=0}^q C_{j_4 j_4 j_1 j_1}\right)^2.
\end{equation}

\vspace{2mm}

By analogy with (\ref{caseu104}) we obtain
$$
{\sf M}\left\{\left(
\sum\limits_{j_2,j_1=0}^q
\sum\limits_{j_3=0}^q C_{j_3 j_3 j_2 j_1}\left(\zeta_{j_1}^{(i_1)}\zeta_{j_2}^{(i_1)}-
{\bf 1}_{\{j_1=j_2\}}\right)\right)^2\right\}=
$$
\begin{equation}
\label{caseu105}
~~~~~ =\sum\limits_{j_2=0}^q \sum\limits_{j_1=0}^{j_2-1}
\left(\sum\limits_{j_3=0}^q C_{j_3 j_3 j_1 j_2}+ 
\sum\limits_{j_3=0}^q C_{j_3 j_3 j_2 j_1}\right)^2
+2\sum\limits_{j_2=0}^q\left(\sum\limits_{j_3=0}^q C_{j_3 j_3 j_2 j_2}\right)^2.
\end{equation}

\vspace{2mm}

Combining (\ref{usl11}), (\ref{caseu100}), (\ref{caseu104}),
and (\ref{caseu105}), we finally have
$$
{\sf M}\left\{\left(I_{(0000)T,t}^{*(i_1i_2i_3i_4)}-
I_{(0000)T,t}^{*(i_1i_2i_3i_4)q}\right)^2\right\}=
\frac{(T-t)^4}{12}-
$$
$$
- \sum_{j_1,j_2,j_3,j_4=0}^{q}
C_{j_4j_3j_2j_1}\Biggl(\sum\limits_{(j_1,j_2)}\Biggl(
\sum\limits_{(j_3,j_4)}
C_{j_4j_3j_2j_1}\Biggr)\Biggr)
+\sum\limits_{j_4,j_3=0}^q
\sum\limits_{j_1=0}^q C_{j_4 j_3 j_1 j_1}\left(C_{j_3j_4}^{10}+C_{j_4j_3}^{10}\right)+
$$
$$
+\sum\limits_{j_4=0}^q \sum\limits_{j_3=0}^{j_4-1}
\left(\sum\limits_{j_1=0}^q C_{j_3 j_4 j_1 j_1}+ 
\sum\limits_{j_1=0}^q C_{j_4 j_3 j_1 j_1}\right)^2
+2\sum\limits_{j_4=0}^q\left(\sum\limits_{j_1=0}^q C_{j_4 j_4 j_1 j_1}\right)^2
-
$$
$$
-\sum\limits_{j_2,j_1=0}^q
\sum\limits_{j_3=0}^q C_{j_3 j_3 j_2 j_1}\left((T-t)C_{j_1}C_{j_2}+
C_{j_1j_2}^{01}+C_{j_2j_1}^{01}\right)+
$$
$$
+\sum\limits_{j_2=0}^q \sum\limits_{j_1=0}^{j_2-1}
\left(\sum\limits_{j_3=0}^q C_{j_3 j_3 j_1 j_2}+ 
\sum\limits_{j_3=0}^q C_{j_3 j_3 j_2 j_1}\right)^2
+2\sum\limits_{j_2=0}^q\left(\sum\limits_{j_3=0}^q C_{j_3 j_3 j_2 j_2}\right)^2+
$$
$$
+\left(\frac{(T-t)^2}{8}-
\sum\limits_{j_3,j_1=0}^q C_{j_3 j_3 j_1 j_1}\right)^2,
$$

\vspace{2mm}
\noindent
where $i_1=i_2\ne i_3=i_4$ and
$$
C_{j}=\int\limits_t^T
\phi_{j}(\tau)d\tau
=\left\{
\begin{matrix}
\sqrt{T-t},\ & j=0\cr\cr
0,\ & j\ne 0
\end{matrix}
.\right.
$$

\vspace{2mm}

Consider the case (\ref{casee12}) by analogy with the case (\ref{casee11}).
Using (\ref{casee400}), we obtain
$$
{\sf M}\left\{\left(I_{(0000)T,t}^{*(i_1i_2i_1i_2)}-
I_{(0000)T,t}^{*(i_1i_2i_1i_2)q}\right)^2\right\}=
{\sf M}\left\{\Biggl(I_{(0000)T,t}^{(i_1i_2i_1i_2)}-
I_{(0000)T,t}^{(i_1i_2i_1i_2)q}\Biggr.\right.-
$$
$$
-\sum\limits_{j_4,j_2=0}^q
\sum\limits_{j_1=0}^q C_{j_4 j_1 j_2 j_1}\zeta_{j_2}^{(i_2)}\zeta_{j_4}^{(i_2)}-
\sum\limits_{j_3,j_1=0}^q
\sum\limits_{j_2=0}^q C_{j_2 j_3 j_2 j_1}\zeta_{j_1}^{(i_1)}\zeta_{j_3}^{(i_1)}
\Biggl.\Biggl.+\sum\limits_{j_2,j_1=0}^q C_{j_2 j_1 j_2 j_1}\Biggr)^2\Biggr\}=
$$
$$
={\sf M}\left\{\Biggl(I_{(0000)T,t}^{(i_1i_2i_1i_2)}-
I_{(0000)T,t}^{(i_1i_2i_1i_2)q}\Biggr.\right.
-\sum\limits_{j_4,j_2=0}^q
\sum\limits_{j_1=0}^q C_{j_4 j_1 j_2 j_1}\left(\zeta_{j_2}^{(i_2)}\zeta_{j_4}^{(i_2)}-
{\bf 1}_{\{j_2=j_4\}}\right)-
$$
$$
-\sum\limits_{j_3,j_1=0}^q
\sum\limits_{j_2=0}^q C_{j_2 j_3 j_2 j_1}\left(\zeta_{j_1}^{(i_1)}\zeta_{j_3}^{(i_1)}-
{\bf 1}_{\{j_1=j_3\}}\right)
\Biggl.\Biggl.-\sum\limits_{j_2,j_1=0}^q C_{j_2 j_1 j_2 j_1}\Biggr)^2\Biggr\}=
$$
$$
={\sf M}\left\{\left(I_{(0000)T,t}^{(i_1i_2i_1i_2)}-
I_{(0000)T,t}^{(i_1i_2i_1i_2)q}\right)^2\right\}+
$$
$$
+{\sf M}\left\{\left(\sum\limits_{j_4,j_2=0}^q
\sum\limits_{j_1=0}^q C_{j_4 j_1 j_2 j_1}\left(\zeta_{j_2}^{(i_2)}\zeta_{j_4}^{(i_2)}-
{\bf 1}_{\{j_2=j_4\}}\right)\right)^2\right\}+
$$
$$
+{\sf M}\left\{\left(\sum\limits_{j_3,j_1=0}^q
\sum\limits_{j_2=0}^q C_{j_2 j_3 j_2 j_1}\left(\zeta_{j_1}^{(i_1)}\zeta_{j_3}^{(i_1)}-
{\bf 1}_{\{j_1=j_3\}}\right)\right)^2\right\}+
$$
\begin{equation}
\label{caseu700}
+\left(\sum\limits_{j_2,j_1=0}^q C_{j_2 j_1 j_2 j_1}\right)^2.
\end{equation}

\vspace{2mm}

Using (\ref{usl12}) and (\ref{caseu700}), we finally get
$$
{\sf M}\left\{\left(I_{(0000)T,t}^{*(i_1i_2i_3i_4)}-
I_{(0000)T,t}^{*(i_1i_2i_3i_4)q}\right)^2\right\}=
\frac{(T-t)^4}{24}-
$$
$$
- \sum_{j_1,j_2,j_3,j_4=0}^{q}
C_{j_4j_3j_2j_1}\Biggl(\sum\limits_{(j_1,j_3)}\Biggl(
\sum\limits_{(j_2,j_4)}
C_{j_4j_3j_2j_1}\Biggr)\Biggr)+
$$
$$
+\sum\limits_{j_4=0}^q \sum\limits_{j_2=0}^{j_4-1}
\left(\sum\limits_{j_1=0}^q C_{j_2 j_1 j_4 j_1}+ \sum\limits_{j_1=0}^q C_{j_4 j_1 j_2 j_1}
\right)^2
+2\sum\limits_{j_4=0}^q\left(\sum\limits_{j_1=0}^q C_{j_4 j_1 j_4 j_1}\right)^2+
$$
$$
+\sum\limits_{j_3=0}^q \sum\limits_{j_1=0}^{j_3-1}
\left(\sum\limits_{j_2=0}^q C_{j_2 j_1 j_2 j_3}+\sum\limits_{j_2=0}^q C_{j_2 j_3 j_2 j_1}
\right)^2
+2\sum\limits_{j_3=0}^q\left(\sum\limits_{j_2=0}^q C_{j_2 j_3 j_2 j_3}\right)^2+
$$
$$
+
\left(\sum\limits_{j_2,j_1=0}^q C_{j_2 j_1 j_2 j_1}\right)^2,
$$

\vspace{2mm}
\noindent
where $i_1=i_3\ne i_2=i_4$.

Consider the case (\ref{casee13}) by analogy with the cases (\ref{casee11}) and (\ref{casee12}).
Using (\ref{casee400}), we obtain
$$
{\sf M}\left\{\left(I_{(0000)T,t}^{*(i_1i_2i_2i_1)}-
I_{(0000)T,t}^{*(i_1i_2i_2i_1)q}\right)^2\right\}=
$$
$$
={\sf M}\left\{\Biggl(I_{(0000)T,t}^{(i_1i_2i_2i_1)}
+\frac{1}{2}
\int\limits_t^T\int\limits_t^{t_4}\int\limits_t^{t_2}
d{\bf w}_{t_1}^{(i_1)}dt_2
d{\bf w}_{t_4}^{(i_1)}
-
I_{(0000)T,t}^{(i_1i_2i_2i_1)q}-\Biggr.\right.
$$
$$
-\sum\limits_{j_3,j_2=0}^q
\sum\limits_{j_1=0}^q C_{j_1 j_3 j_2 j_1}\zeta_{j_2}^{(i_2)}\zeta_{j_3}^{(i_2)}-
\sum\limits_{j_4,j_1=0}^q
\sum\limits_{j_2=0}^q C_{j_4 j_2 j_2 j_1}\zeta_{j_1}^{(i_1)}\zeta_{j_4}^{(i_1)}+
\left.\Biggl.
\sum\limits_{j_2,j_1=0}^q C_{j_1 j_2 j_2 j_1}\Biggr)^2\right\}=
$$

\vspace{1mm}
$$
={\sf M}\left\{\Biggl(I_{(0000)T,t}^{(i_1i_2i_2i_1)}-
I_{(0000)T,t}^{(i_1i_2i_2i_1)q}+\Biggr.\right.
$$
$$
+\frac{1}{2}
\int\limits_t^T\int\limits_t^{t_4}\int\limits_t^{t_2}
d{\bf w}_{t_1}^{(i_1)}dt_2
d{\bf w}_{t_4}^{(i_1)}
-\sum\limits_{j_4,j_1=0}^q
\sum\limits_{j_2=0}^q C_{j_4 j_2 j_2 j_1}\left(\zeta_{j_1}^{(i_1)}\zeta_{j_4}^{(i_1)}-
{\bf 1}_{\{j_1=j_4\}}\right)-
$$
$$
\Biggl.\Biggl.-\sum\limits_{j_3,j_2=0}^q
\sum\limits_{j_1=0}^q C_{j_1 j_3 j_2 j_1}\left(\zeta_{j_2}^{(i_2)}\zeta_{j_3}^{(i_2)}-
{\bf 1}_{\{j_2=j_3\}}\right)
-
\sum\limits_{j_2,j_1=0}^q C_{j_1 j_2 j_2 j_1}\Biggr)^2\Biggr\}=
$$
$$
={\sf M}\left\{\left(I_{(0000)T,t}^{(i_1i_2i_2i_1)}-
I_{(0000)T,t}^{(i_1i_2i_2i_1)q}\right)^2\right\}+
{\sf M}\left\{\Biggl(\frac{1}{2}
\int\limits_t^T\int\limits_t^{t_4}(t_4-t_1)
d{\bf w}_{t_1}^{(i_1)}
d{\bf w}_{t_4}^{(i_1)}
-\Biggr.\right.
$$
$$
\left.\left.-\sum\limits_{j_4,j_1=0}^q
\sum\limits_{j_2=0}^q C_{j_4 j_2 j_2 j_1}\left(\zeta_{j_1}^{(i_1)}\zeta_{j_4}^{(i_1)}-
{\bf 1}_{\{j_1=j_4\}}\right)\right)^2\right\}+
$$
$$
+{\sf M}\left\{\left(\sum\limits_{j_3,j_2=0}^q
\sum\limits_{j_1=0}^q C_{j_1 j_3 j_2 j_1}\left(\zeta_{j_2}^{(i_2)}\zeta_{j_3}^{(i_2)}-
{\bf 1}_{\{j_2=j_3\}}\right)\right)^2\right\}+
$$
$$
+
\left(\sum\limits_{j_2,j_1=0}^q C_{j_1 j_2 j_2 j_1}\right)^2=
$$

$$
={\sf M}\left\{\left(I_{(0000)T,t}^{(i_1i_2i_2i_1)}-
I_{(0000)T,t}^{(i_1i_2i_2i_1)q}\right)^2\right\}+\frac{(T-t)^4}{48}-
$$
$$
-\sum\limits_{j_4,j_1=0}^q
\sum\limits_{j_2=0}^q C_{j_4 j_2 j_2 j_1}\left(
\int\limits_t^T \phi_{j_4}(t_4)\int\limits_t^{t_4}(t_4-t_1)\phi_{j_1}(t_1)
dt_1 dt_4+\right.
$$
$$
\left.
+\int\limits_t^T \phi_{j_1}(t_4)\int\limits_t^{t_4}(t_4-t_1)\phi_{j_4}(t_1)
dt_1 dt_4\right)+
$$
$$
+{\sf M}\left\{\left(\sum\limits_{j_4,j_1=0}^q
\sum\limits_{j_2=0}^q C_{j_4 j_2 j_2 j_1}\left(\zeta_{j_1}^{(i_1)}\zeta_{j_4}^{(i_1)}-
{\bf 1}_{\{j_1=j_4\}}\right)\right)^2\right\}+
$$
$$
+{\sf M}\left\{\left(\sum\limits_{j_3,j_2=0}^q
\sum\limits_{j_1=0}^q C_{j_1 j_3 j_2 j_1}\left(\zeta_{j_2}^{(i_2)}\zeta_{j_3}^{(i_2)}-
{\bf 1}_{\{j_2=j_3\}}\right)\right)^2\right\}+
$$
\begin{equation}
\label{caseu800}
+
\left(\sum\limits_{j_2,j_1=0}^q C_{j_1 j_2 j_2 j_1}\right)^2.
\end{equation}

\vspace{2mm}

Using (\ref{usl13}) and (\ref{caseu800}), we finally get
$$
{\sf M}\left\{\left(I_{(0000)T,t}^{*(i_1i_2i_3i_4)}-
I_{(0000)T,t}^{*(i_1i_2i_3i_4)q}\right)^2\right\}=\frac{(T-t)^4}{16}-
$$
$$
- \sum_{j_1,j_2,j_3,j_4=0}^{q}
C_{j_4j_3j_2j_1}\Biggl(\sum\limits_{(j_1,j_4)}\Biggl(
\sum\limits_{(j_2,j_3)}
C_{j_4j_3j_2j_1}\Biggr)\Biggr)-
$$
$$
-\sum\limits_{j_4,j_1=0}^q
\sum\limits_{j_2=0}^q C_{j_4 j_2 j_2 j_1}\left(C_{j_4j_1}^{10}+
C_{j_1j_4}^{10}-C_{j_4j_1}^{01}-C_{j_1j_4}^{01}\right)+
$$
$$
+\sum\limits_{j_4=0}^q \sum\limits_{j_1=0}^{j_4-1}
\left(\sum\limits_{j_2=0}^q C_{j_1 j_2 j_2 j_4} + \sum\limits_{j_2=0}^q C_{j_4 j_2 j_2 j_1}
\right)^2+
2\sum\limits_{j_4=0}^q\left(\sum\limits_{j_2=0}^q C_{j_4 j_2 j_2 j_4}\right)^2+
$$
$$
+\sum\limits_{j_3=0}^q \sum\limits_{j_2=0}^{j_3-1}
\left(\sum\limits_{j_1=0}^q C_{j_1 j_2 j_3 j_1} + \sum\limits_{j_1=0}^q C_{j_1 j_3 j_2 j_1}
\right)^2
+2\sum\limits_{j_3=0}^q\left(\sum\limits_{j_1=0}^q C_{j_1 j_3 j_3 j_1}\right)^2+
$$
$$
+
\left(\sum\limits_{j_2,j_1=0}^q C_{j_1 j_2 j_2 j_1}\right)^2,
$$

\vspace{2mm}
\noindent
where $i_1=i_4\ne i_2=i_3$.

\section{Optimization of the Mean-Square Approximation Procedures for 
Iterated Stra\-to\-no\-vich Stochastic Integrals Based on 
Theorems 2.2, 2.8 and Multiple Fourier--Le\-gend\-re Series}

This section is devoted to optimization of the mean-square approximation 
procedures for iterated Stratonovich stochastic integrals (\ref{k1001xxxx})
of multiplicities 1 to 3 based on 
Theorems 2.2, 2.8 and multiple Fourier--Legendre series \cite{matec1}${}^{1}$.

\footnotetext[1]{The results of this section were obtained 
jointly with Kuznetsov M.D., who is also a co-author of the publications 
\cite{123789000}, 
\cite{Kuz-Kuz}-\cite{Mikh-2aaaaa}, \cite{new-new-1}, \cite{new-new-3}, \cite{matec1}.}

The mentioned stochastic integrals are part of 
strong numerical methods with convergence orders 1.0 and 1.5 for 
It\^{o} SDEs with multidimensional 
non-commutative noise (see (\ref{al1x}), (\ref{al2x})).

We show that the lengths of sequences of independent standard 
Gaussian random variables required for the mean-square approximation of 
iterated Stratonovich stochastic integrals (\ref{k1001xxxx}) 
can be significantly reduced without the loss of 
the mean-square accuracy of approximation for these stochastic integrals.

Using Theorems 2.2, 2.8 and
the system of Legendre polynomials, we obtain
the following approximations of iterated Stra\-to\-no\-vich
stochastic integrals (\ref{k1001xxxx})
$$
I_{(0)T,t}^{*(i_1)}=\sqrt{T-t}\zeta_0^{(i_1)},
$$
$$
I_{(1)T,t}^{*(i_1)}=-\frac{(T-t)^{3/2}}{2}\left(\zeta_0^{(i_1)}+
\frac{1}{\sqrt{3}}\zeta_1^{(i_1)}\right),
$$
\begin{equation}
\label{dsds400}
~~~~~I_{(00)T,t}^{*(i_1 i_2)q}=
\frac{T-t}{2}\left(\zeta_0^{(i_1)}\zeta_0^{(i_2)}+\sum_{i=1}^{q}
\frac{1}{\sqrt{4i^2-1}}\left(
\zeta_{i-1}^{(i_1)}\zeta_{i}^{(i_2)}-
\zeta_i^{(i_1)}\zeta_{i-1}^{(i_2)}\right)\right),
\end{equation}
\begin{equation}
\label{dsds401}
I_{(000)T,t}^{*(i_1i_2i_3)q_1}
=
\sum_{j_1,j_2,j_3=0}^{q_1}
C_{j_3j_2j_1}
\zeta_{j_1}^{(i_1)}\zeta_{j_2}^{(i_2)}\zeta_{j_3}^{(i_3)},
\end{equation}

\vspace{2mm}
\noindent
where 
$$
\zeta_{j}^{(i)}=
\int\limits_t^T \phi_{j}(s) d{\bf w}_s^{(i)}\ \ \ (i=1,\ldots,m,\
j=0,1,\ldots)
$$
are independent standard Gaussian random variables
for various
$i$ or $j$, $\{\phi_j(x)\}_{j=0}^{\infty}$ is a complete
orthonormal system of Legendre polynomials  
in the space $L_2([t, T])$  (see (\ref{4009d})),

\vspace{-1mm}
$$
C_{j_3 j_2 j_1}=\frac{\sqrt{(2j_1+1)(2j_2+1)(2j_3+1)}}{8}(T-t)^{3/2}\bar C_{j_3j_2j_1},
$$
$$
\bar C_{j_3j_2j_1}=
\int\limits_{-1}^{1}P_{j_3}(z)
\int\limits_{-1}^{z}P_{j_2}(y)
\int\limits_{-1}^{y}
P_{j_1}(x)dx dy dz,
$$
          
\vspace{1mm}
\noindent
$P_j(x)$ is the Legendre polynomial (see (\ref{agentww1})).

Denote
$$
E^{*(l_1\ldots l_k)}_p\stackrel{\sf def}{=}
{\sf M}\left\{\left(
I_{(l_1\ldots l_k)T,t}^{*(i_1\ldots i_k)}-
I_{(l_1\ldots l_k)T,t}^{*(i_1\ldots i_k)p}\right)^2\right\},
$$

\noindent
where $I_{(l_1\ldots l_k)T,t}^{*(i_1\ldots i_k)}$ is the 
iterated Stra\-to\-no\-vich stochastic integral (\ref{k1001xxxx})
and $I_{(l_1\ldots l_k)T,t}^{*(i_1\ldots i_k)p}$ is the
mean-square approximation of this stochastic integral. More precisely, 
the approximations 
$I_{(00)T,t}^{*(i_1i_2)q},$ 
$I_{(000)T,t}^{*(i_1i_2i_3)q_1}$
are defined by
(\ref{dsds400}), (\ref{dsds401}).

Using (\ref{fff09}), (\ref{dest3}), (\ref{dest80}), (\ref{dest70}), (\ref{dest60}), we get

\vspace{-2mm}
\begin{equation}
\label{dsds402a}
E^{*(00)}_q
=
\frac{(T-t)^2}{2}\left(\frac{1}{2}-\sum_{i=1}^q
\frac{1}{4i^2-1}\right)\ \ (i_1\ne i_2),
\end{equation}

\vspace{-2mm}
\begin{equation}
\label{dsds402b}
~~~~~~~~E^{*(000)}_{q_{1,1}} = \frac{(T-t)^3}{6}-
\sum_{j_1,j_2,j_3=0}^{q_{1,1}}C_{j_3j_2j_1}^2
\ \ \ (i_1\ne i_2,\ i_1\ne i_3,\ i_2\ne i_3),
\end{equation}

$$
E^{*(000)}_{q_{1,2}}=
\frac{(T-t)^3}{4}
-\sum_{j_1,j_2,j_3=0}^{q_{1,2}} C_{j_3j_2j_1}^2-
\sum_{j_1,j_2,j_3=0}^{q_{1,2}} C_{j_3j_1j_2}C_{j_3j_2j_1}-
$$
\begin{equation}
\label{dsds402}
-\frac{(T-t)^{3/2}}{2}
\sum_{j_1=0}^{q_{1,2}}
\left(C_{0j_1j_1}+\frac{1}{\sqrt{3}}C_{1j_1j_1}\right)+
\sum_{j_3=0}^{q_{1,2}}\left(\sum_{j_1=0}^{q_{1,2}}
C_{j_3j_1j_1}\right)^2\ \ \ (i_1=i_2\ne i_3),
\end{equation}

\vspace{-2mm}
$$
E^{*(000)}_{q_{1,3}}=
\frac{(T-t)^3}{4}
-\sum_{j_1,j_2,j_3=0}^{q_{1,3}} C_{j_3j_2j_1}^2-
\sum_{j_1,j_2,j_3=0}^{q_{1,3}} C_{j_2j_3j_1}C_{j_3j_2j_1}-
$$
\begin{equation}
\label{dsds403}
-\frac{(T-t)^{3/2}}{2}
\sum_{j_3=0}^{q_{1,3}}
\left(C_{j_3j_3 0}-\frac{1}{\sqrt{3}}C_{j_3j_3 1}\right)+
\sum_{j_1=0}^{q_{1,3}}\left(\sum_{j_3=0}^{q_{1,3}}
C_{j_3j_3j_1}\right)^2\ \ \ (i_1\ne i_2=i_3),
\end{equation}

\vspace{-2mm}
$$
E^{*(000)}_{q_{1,4}}=
\frac{(T-t)^3}{6}
-\sum_{j_1,j_2,j_3=0}^{q_{1,4}} C_{j_3j_2j_1}^2-
\sum_{j_1,j_2,j_3=0}^{q_{1,4}} C_{j_3j_2j_1}C_{j_1j_2j_3}+
$$
\begin{equation}
\label{dsds404}
+\sum_{j_2=0}^{q_{1,4}}
\left(\sum_{j_1=0}^{q_{1,4}}C_{j_1j_2j_1}\right)^2\ \ \ (i_1=i_3\ne i_2).
\end{equation}

\vspace{2mm}

Note that 
the number of conditions (\ref{dsds402b})--(\ref{dsds404}) is quite 
large, which is inconvenient for practice. In this section, we propose the hypothesis 
that all the formulas (\ref{dsds402b})--(\ref{dsds404}) can be replaced by 
the formula (\ref{dsds402b})
in which we can suppose that $i_1, i_2, i_3=1,\ldots, m.$
At that we will not have a noticeable loss of the 
mean-square approximation accuracy of iterated Stra\-to\-no\-vich stochastic integrals.

Consider the following condition
\begin{equation}
\label{dsds405}
~~~~~~~ E^{*(00)}_q\le (T-t)^4,\ \ \ E^{*(000)}_{q_{1,i}} \le (T-t)^4,\ \ \ i=1,\ldots,4.
\end{equation}

\begin{table}
\centering
\caption{Conditions $E^{*(000)}_{q_{1,i}} \le (T-t)^4$,\ $i=1,\ldots,4.$}
\label{tab:5.59}      

\begin{tabular}{p{1.3cm}p{1.3cm}p{1.3cm}p{1.3cm}p{1.3cm}p{1.3cm}p{1.3cm}}

\hline\noalign{\smallskip}

$T-t$&0.011&0.008&0.0045&0.0035&0.0027&0.0025\\

\noalign{\smallskip}\hline\noalign{\smallskip}
$q_{1,1}$&12&16&28&36&47&50\\
\noalign{\smallskip}\hline\noalign{\smallskip}
$q_{1,2}$&6&8&14&18&23&25\\
\noalign{\smallskip}\hline\noalign{\smallskip}
$q_{1,3}$&6&8&14&18&23&25\\
\noalign{\smallskip}\hline\noalign{\smallskip}
$q_{1,4}$&12&16&28&36&47&51\\
\noalign{\smallskip}\hline\noalign{\smallskip}
\vspace{4mm}
\end{tabular}
\end{table}

\begin{table}
\centering
\caption{The condition (\ref{dsds405}).}
\label{tab:5.60}      
\begin{tabular}{p{1.5cm}p{1.5cm}p{1.5cm}p{1.5cm}p{1.5cm}}
\hline\noalign{\smallskip}
$T-t$&$2^{-1}$&$2^{-3}$&$2^{-5}$&$2^{-8}$\\
\noalign{\smallskip}\hline\noalign{\smallskip}
$q$&1&8&128&8192\\
\noalign{\smallskip}\hline\noalign{\smallskip}
$q_{1,1}$&0&1&4&32\\
\noalign{\smallskip}\hline\noalign{\smallskip}
$q_{1,2}$&0&0&2&16\\
\noalign{\smallskip}\hline\noalign{\smallskip}
$q_{1,3}$&0&0&2&16\\
\noalign{\smallskip}\hline\noalign{\smallskip}
$q_{1,4}$&0&0&4&33\\
\noalign{\smallskip}\hline\noalign{\smallskip}
\vspace{4mm}
\end{tabular}
\end{table}

\begin{table}
\centering
\caption{Values $E_{q_{1,i}}^{*(000)} \cdot (T-t)^{-3}\stackrel{\sf def}{=}E_{q_{1,i}}^{*},$\ 
$i=1,\ldots,4.$}
\label{tab:5.61}      
\begin{tabular}{p{1.3cm}p{1.9cm}p{1.9cm}p{1.9cm}p{1.9cm}p{1.9cm}p{1.9cm}}
\hline\noalign{\smallskip}
$T-t$&$0.011$&$0.008$&$0.0045$&$0.0035$&$0.0027$&$0.0025$\\
\noalign{\smallskip}\hline\noalign{\smallskip}
$q_{1,1}$&$12$&$16$&$28$&$36$&$47$&$50$\\
\noalign{\smallskip}\hline\noalign{\smallskip}
$E_{q_{1,1}}^{*}$&0.010154&0.007681&0.004433&0.003456&0.002652&0.002494\\
\noalign{\smallskip}\hline\noalign{\smallskip}
$q_{1,2}$&$12$&$16$&$28$&$36$&$47$&$50$\\
\noalign{\smallskip}\hline\noalign{\smallskip}
$E_{q_{1,2}}^{*}$&0.005102&0.003855&0.002221&0.001731&0.001328&0.001248\\
\noalign{\smallskip}\hline\noalign{\smallskip}
$q_{1,3}$&$12$&$16$&$28$&$36$&$47$&$50$\\
\noalign{\smallskip}\hline\noalign{\smallskip}
$E_{q_{1,3}}^{*}$&0.005102&0.003855&0.002221&0.001731&0.001328&0.001248\\
\noalign{\smallskip}\hline\noalign{\smallskip}
$q_{1,4}$&$12$&$16$&$28$&$36$&$47$&$50$\\
\noalign{\smallskip}\hline\noalign{\smallskip}
$E_{q_{1,4}}^{*}$&0.010407&0.007845&0.004500&0.003501&0.002680&0.002519\\
\noalign{\smallskip}\hline\noalign{\smallskip}
\end{tabular}
\end{table}

\begin{table}
\centering
\caption{Comparison of numbers $q_{1,1}$ and $p_{1,1}$.}
\label{tab:5.62}      
\begin{tabular}{p{1.9cm}p{1.3cm}p{1.3cm}p{1.3cm}p{1.3cm}p{1.3cm}p{1.3cm}}
\hline\noalign{\smallskip}
$T-t$&$2^{-1}$&$2^{-2}$&$2^{-3}$&$2^{-4}$&$2^{-5}$&$2^{-6}$\\
\noalign{\smallskip}\hline\noalign{\smallskip}
$q_{1,1}$&0&0&1&2&4&8\\
\noalign{\smallskip}\hline\noalign{\smallskip}
$(q_{1,1}+1)^3$&1&1&8&27&125&729\\
\noalign{\smallskip}\hline\noalign{\smallskip}
$p_{1,1}$&1&3&6&12&24&48\\
\noalign{\smallskip}\hline\noalign{\smallskip}
$(p_{1,1}+1)^3$&8&64&343&2197&15625&117649\\
\noalign{\smallskip}\hline\noalign{\smallskip}
\end{tabular}
\vspace{2mm}
\end{table}

Let us show by numerical experiments that in most situations the following inequality 
is fulfilled 
\begin{equation}
\label{dsds500}
q_{1,1}\ge q_{1,i},\ \ \ i=2,3,4,
\end{equation}

\noindent
where $q_{1,i}$
($i=1,\ldots,4$)
are minimal natural numbers satisfying 
the condition (\ref{dsds405}).

In Tables 5.59--5.61 we can see the results of numerical experiments. 
These results confirm the hypothesis proposed earlier in this section.
Note that in Table 5.61 we calculate the mean-square 
approximation errors of iterated Stra\-to\-no\-vich stochastic integrals
in the case when 
$$
q_{1,i}=q_{1,1},\ \ \ i=2,3,4,
$$

\noindent
where $q_{1,1}$ 
is the minimal natural number satisfying 
the condition
(\ref{dsds405}). In this case, there is no noticeable loss of the 
mean-square approximation accuracy of iterated Stra\-to\-no\-vich stochastic integrals
(see Table 5.61). 
This means that all the formulas (\ref{dsds402b})--(\ref{dsds404}) can be re\-placed by 
the formula (\ref{dsds402b})
in which we can suppose that $i_1, i_2, i_3=1,\ldots, m.$

Let $q_{1,1}$ be the minimal natural number satisfying the condition
\begin{equation}
\label{dsds501}
E^{(000)}_{q_{1,1}} \le (T-t)^4,
\end{equation}

\noindent
where the left-hand side of (\ref{dsds501})
is defined by the formula (\ref{dsds402b}).

Let $p_{1,1}$ be the minimal natural number satisfying the condition
\begin{equation}
\label{dsds503}
3! \cdot E^{(000)}_{p_{1,1}} \le (T-t)^4,
\end{equation}
where the value $E^{(000)}_{p_{1,1}}$ on the left-hand side of 
(\ref{dsds503}) is defined by the formula (\ref{dsds402b})
(recall that $3!$ is included in the inequality (\ref{agentyz1}) for the case $k=3$).

In Table 5.62 we can see the 
numerical comparison of numbers $q_{1,1}$  and $p_{1,1}$. Obviously, 
excluding of the multiplier factor\hspace{0.3mm} $3!$ 
essentially (in many times)\hspace{0.3mm} reduces 
the calculation 
costs for the mean-square approximations of iterated Stra\-to\-no\-vich stochastic 
integrals. Note that in this section we use the exactly calculated Fourier--Legendre 
coefficients using the Python programming language \cite{Kuz-Kuz}, \cite{Mikh-1}.

\chapter{Other Methods of Approximation of Specific Iterated
It\^{o} and Stra\-to\-no\-vich Sto\-chas\-tic Integrals of Multiplicities
1 to 4}

\section{New Simple Method for Obtainment an Expansion of Iterated
It\^{o} Stochastic integrals of Multiplicity 2 Based on the 
Wiener Process Expansion 
Using Legendre Polynomials and Trigonometric 
Functions}

This section is devoted to the expansion of 
iterated It\^{o} stochastic integrals of multiplicity 2
based on the Wiener process expansion 
using complete orthonormal systems of functions in $L_2([t, T])$.
The expansions of these stochastic integrals 
using Legendre polynomials and trigonometric functions are considered.
In contrast to the method of expansion of iterated It\^{o} stochastic integrals 
based on the Karhunen--Lo\`{e}ve expansion of the Brownian
bridge process
\cite{Zapad-1}-\cite{Zapad-3}, 
this method allows the use of different systems of basis functions, 
not only the trigonometric system of functions. 
The proposed method 
makes it possible to obtain expansions of 
iterated It\^{o} stochastic integrals of multiplicity
2 much 
easier than the method based on generalized multiple Fourier series
(see Chapters 1 and 2).                       
The latter involve the calculation of coefficients 
of multiple Fourier series, which is a time-consuming task.
However, the proposed method can be applied only to iterated
It\^{o} stochastic integrals of multiplicity 2.

It is well known that the idea of representing of the Wiener 
process as a functional series with random coefficients (that are 
independent standard Gaussian random variables) with using the
complete orthonormal system of trigonometric functions in $L_2([0, T])$ 
goes back to the works of Wiener \cite{7b} (1924) and L{\'e}vy
\cite{7c} (1951).
The specified series was used in \cite{7b} and \cite{7c} 
for construction of the Brownian motion process (Wiener process). 
A little later, It\^{o} and McKean in \cite{7d} (1965) used for 
this purpose the complete orthonormal system of
Haar functions in $L_2([0, T])$.

Let ${\bf w}_{\tau},$ $\tau\in[0, T]$ be an $m$-dimestional
standard Wiener process with independent components
${\bf w}_{\tau}^{(i)},$ $i=1,\ldots,m.$ 

We have w.~p.~1
$$
{\bf w}_s^{(i)}-{\bf w}_t^{(i)}=
\int\limits_t^s d{\bf w}_{\tau}^{(i)}=
\int\limits_t^T {\bf 1}_{\{\tau<s\}} d{\bf w}_{\tau}^{(i)},
$$
where 
$$
\int\limits_t^T {\bf 1}_{\{\tau<s\}} d{\bf w}_{\tau}^{(i)}
$$
is the It\^{o} stochastic integral,
$t\ge 0,$ and
$$
{\bf 1}_{\{\tau<s\}}=
\left\{
\begin{matrix}
1,\ &\ \tau<s\cr\cr
0,\ &\ \hbox{\rm otherwise}
\end{matrix},\ \ \ \tau, s\in [t, T].\right.
$$

Consider the Fourier expansion of ${\bf 1}_{\{\tau<s\}}\in L_2([t,T])$ at the 
interval $[t, T]$ (see, for example, \cite{7e})
\begin{equation}
\label{eee}
\sum_{j=0}^{\infty}\int\limits_t^T
{\bf 1}_{\{\tau<s\}}\phi_j(\tau)d\tau \phi_j(\tau)=
\sum_{j=0}^{\infty}\int\limits_t^s
\phi_j(\tau)d\tau \phi_j(\tau),
\end{equation}
where $\{\phi_j(\tau)\}_{j=0}^{\infty}$ is a complete 
orthonormal system of functions in the space $L_2([t, T])$
and the series (\ref{eee}) converges
in the mean-square sense, i.e.
$$
\int\limits_t^T\left({\bf 1}_{\{\tau<s\}}-
\sum_{j=0}^{q}\int\limits_t^s
\phi_j(\tau)d\tau \phi_j(\tau)\right)^2 d\tau\to 0\ \ \
\hbox{if}\ \ \ q\to\infty.
$$

Let ${\bf w}_{s,t}^{(i)q}$ be 
the mean-square approximation of the process
${\bf w}_s^{(i)}-{\bf w}_t^{(i)}$, 
which has the following form
\begin{equation}
\label{jq1}
~~~~{\bf w}_{s,t}^{(i)q}=
\int\limits_t^T
\left(\sum_{j=0}^{q}\int\limits_t^s
\phi_j(\tau)d\tau \phi_j(\tau)\right)
d{\bf w}_{\tau}^{(i)}=\sum_{j=0}^{q}\int\limits_t^s
\phi_j(\tau)d\tau \int\limits_t^T\phi_j(\tau)
d{\bf w}_{\tau}^{(i)}.
\end{equation}

Moreover,
$$
{\sf M}\left\{\biggl(
{\bf w}_s^{(i)}-{\bf w}_t^{(i)}-
{\bf w}_{s,t}^{(i)q}\biggr)^2\right\}=
$$
$$
={\sf M}\left\{\left(
\int\limits_t^T\left(
{\bf 1}_{\{\tau<s\}}-
\sum_{j=0}^{q}\int\limits_t^s
\phi_j(\tau)d\tau \phi_j(\tau)\right)
d{\bf w}_{\tau}^{(i)}\right)^2\right\}=
$$
\begin{equation}
\label{t1xx}
=\int\limits_t^T
\left({\bf 1}_{\{\tau<s\}}-
\sum_{j=0}^{q}\int\limits_t^s
\phi_j(\tau)d\tau \phi_j(\tau)\right)^2 d\tau\to 0\ \ \
\hbox{if}\ q\to\infty.
\end{equation}

In \cite{Zapad-5} it was proposed to use an expansion similar 
to (\ref{jq1}) for the expansion of iterated It\^{o} stochastic 
integrals 
\begin{equation}
\label{ito-kuku}
I_{(00)T,t}^{(i_1i_2)}=
\int\limits_t^T\int\limits_t^{t_{2}}
d{\bf w}_{t_1}^{(i_1)}d{\bf w}_{t_2}^{(i_2)}\ \ \ (i_1,i_2=1,\ldots,m).
\end{equation}

At that, to obtain the mentioned expansion of (\ref{ito-kuku}), the truncated 
expansions (\ref{jq1}) of components of the Wiener 
process ${\bf w}_s$ have been
iteratively substituted in the single integrals \cite{Zapad-5}. 
This procedure leads to the calculation
of coefficients of the double Fourier series, 
which is a time-consuming task for not too complex problem
of expansion of the iterated It\^{o} stochastic integral (\ref{ito-kuku}).
In \cite{Zapad-5} the expansions on 
the base of Haar functions and trigonometric functions have been
considered.

In contrast to \cite{Zapad-5} we substitute the expansion 
(\ref{jq1}) only one time and only into the innermost integral
in (\ref{ito-kuku}). 
This procedure leads to the simple calculation
of the coefficients 
$$
\int\limits_t^s
\phi_j(\tau)d\tau\ \ \ (j=0, 1, 2,\ldots)
$$
of the usual (not double) Fourier series.

Moreover, we use the Legendre polynomials \cite{arxiv-23}, \cite{new-art-1}
for the construction 
of the expansion of (\ref{ito-kuku}). For the first time the 
Legendre polynomials have been applied in the 
framework of the mentioned  problem in the author's papers 
\cite{old-art-1} (1997), \cite{old-art-2} (1998), \cite{old-art-3} (2000),
\cite{old-art-4} (2001) (also see \cite{1}-\cite{new-art-1})
while in the papers of other author's these polynomials 
have not been considered as the basis functions for the construction
of expansions of iterated It\^{o} or Stratonovich stochastic integrals.

{\bf Theorem 6.1} \cite{12a}-\cite{12aa}, \cite{arxiv-23}, \cite{new-art-1}. {\it Let
$\phi_j(\tau)$ $(j=0, 1, \ldots )$ be an arbitrary complete
orthonormal system of functions in the space $L_2([t, T]).$
Let 
\begin{equation}
\label{l1}
\int\limits_t^T
{\bf w}_{s,t}^{(i_1)q}
d{\bf w}_s^{(i_2)}=
\sum_{j=0}^{q}\int\limits_t^T\phi_j(\tau)
d{\bf w}_{\tau}^{(i_1)}
\int\limits_t^T\int\limits_t^s
\phi_j(\tau)d\tau d{\bf w}_{s}^{(i_2)}
\end{equation}
be an approximation of the iterated It\^{o} stochastic integral
{\rm (\ref{ito-kuku})}
for $i_1\ne i_2$. Then 
$$
I_{(00)T,t}^{(i_1i_2)}
=\hbox{\vtop{\offinterlineskip\halign{
\hfil#\hfil\cr
{\rm l.i.m.}\cr
$\stackrel{}{{}_{q\to \infty}}$\cr
}} }
\int\limits_t^T
{\bf w}_{s,t}^{(i_1)q}
d{\bf w}_s^{(i_2)}\ \ \ (i_1\ne i_2),
$$
where $i_1, i_2=1,\ldots,m$.}

{\bf Proof.} Using the standard properties of the It\^{o} stochastic
integral as well as (\ref{t1xx}) and the property of orthonormality
of functions 
$\phi_j(\tau)$ $(j=0, 1, \ldots )$ at the interval $[t, T],$
we obtain 
$$
{\sf M}\left\{\left(
\int\limits_t^T
\int\limits_t^s
d{\bf w}_{\tau}^{(i_1)}d{\bf w}_{s}^{(i_2)}-
\int\limits_t^T
{\bf w}_{s,t}^{(i_1)q}
d{\bf w}_s^{(i_2)}\right)^2\right\}=
$$
$$
=
\int\limits_t^T
{\sf M}\left\{\biggl(
{\bf w}_s^{(i_1)}-{\bf w}_t^{(i_1)}-
{\bf w}_{s,t}^{(i_1)q}\biggr)^2\right\}ds=
$$
$$
=
\int\limits_t^T\int\limits_t^T
\left({\bf 1}_{\{\tau<s\}}-
\sum_{j=0}^{q}\int\limits_t^s
\phi_j(\tau)d\tau \phi_j(\tau)\right)^2d\tau ds=
$$
\begin{equation}
\label{l5}
=\int\limits_t^T\left((s-t)-
\sum_{j=0}^{q}\left(\int\limits_t^s
\phi_j(\tau)d\tau\right)^2\right) ds.
\end{equation}

Using the continuity of the functions $u_q(s)$ (see below),
the nondecreasing property of
the functional sequence
$$
u_q(s)=\sum_{j=0}^{q}\left(\int\limits_t^s
\phi_j(\tau)d\tau\right)^2,
$$
and the continuity of the limit function
$u(s)= s-t$
according to Dini's 
Theorem,
we have the uniform convergence 
$u_q(s)$ to $u(s)$ at the interval $[t, T]$.

Then from this fact as well as from (\ref{l5}) we obtain
\begin{equation}
\label{rez1}
I_{(00)T,t}^{(i_1i_2)}=\hbox{\vtop{\offinterlineskip\halign{
\hfil#\hfil\cr
{\rm l.i.m.}\cr
$\stackrel{}{{}_{q\to \infty}}$\cr
}} }
\int\limits_t^T
{\bf w}_{s,t}^{(i_1)q}
d{\bf w}_s^{(i_2)}.
\end{equation}

Note that we could also use Lebesgue's Dominated Convergence Theorem 
in (\ref{l5}) to obtain (\ref{rez1}).
Theorem~6.1 is proved.

Let $\{\phi_j(\tau)\}_{j=0}^{\infty}$ be a complete 
orthonormal system of Legendre polynomials in the space $L_2([t, T]),$
which has the form (\ref{4009d}).
Then
\begin{equation}
\label{l2}
~~~~~~~~\int\limits_t^s
\phi_j(\tau)d\tau=
\frac{T-t}{2}\left(\frac{\phi_{j+1}(s)}{\sqrt{(2j+1)(2j+3)}}-
\frac{\phi_{j-1}(s)}{\sqrt{4j^2-1}}\right)\ \ \hbox{for}\ \ j\ge 1.
\end{equation}

Let us denote 
$$
\zeta_j^{(i)}=\int\limits_t^T\phi_j(\tau)
d{\bf w}_{\tau}^{(i)}\ \ \ (i=1,\ldots,m).
$$

From (\ref{l1}) and (\ref{l2}) we obtain
$$
\int\limits_t^T
{\bf w}_{s,t}^{(i_1)q}
d{\bf w}_s^{(i_2)}=
\frac{1}{\sqrt{T-t}}\zeta_0^{(i_1)}
\int\limits_t^T(s-t){\bf w}_s^{(i_2)}+
$$
$$
+
\frac{T-t}{2}\sum_{j=1}^q
\zeta_j^{(i_1)}
\left(\frac{1}{\sqrt{(2j+1)(2j+3)}}\zeta_{j+1}^{(i_2)}-
\frac{1}{\sqrt{4j^2-1}}\zeta_{j-1}^{(i_2)}\right)=
$$

$$
=
\frac{T-t}{2}\zeta_0^{(i_1)}\left(
\zeta_0^{(i_2)}+\frac{1}{\sqrt{3}}\zeta_1^{(i_2)}\right)
+
$$
$$
+
\frac{T-t}{2}\sum_{j=1}^q
\zeta_j^{(i_1)}
\left(\frac{1}{\sqrt{(2j+1)(2j+3)}}\zeta_{j+1}^{(i_2)}-
\frac{1}{\sqrt{4j^2-1}}\zeta_{j-1}^{(i_2)}\right)=
$$
$$
=\frac{T-t}{2}\left(\zeta_0^{(i_1)}\zeta_0^{(i_2)}+\sum_{j=1}^{q}
\frac{1}{\sqrt{4j^2-1}}\left(
\zeta_{j-1}^{(i_1)}\zeta_{j}^{(i_2)}-
\zeta_j^{(i_1)}\zeta_{j-1}^{(i_2)}\right)\right)+
$$
\begin{equation}
\label{l6}
+
\frac{T-t}{2}\zeta_q^{(i_1)}\zeta_{q+1}^{(i_2)}
\frac{1}{\sqrt{(2q+1)(2q+3)}}.
\end{equation}

\vspace{2mm}

Then from (\ref{rez1}) and (\ref{l6}) we get
$$
I_{(00)T,t}^{(i_1i_2)}
=\hbox{\vtop{\offinterlineskip\halign{
\hfil#\hfil\cr
{\rm l.i.m.}\cr
$\stackrel{}{{}_{q\to \infty}}$\cr
}} }
\int\limits_t^T
{\bf w}_{s,t}^{(i_1)q}
d{\bf w}_s^{(i_2)}=
$$
\begin{equation}
\label{l10}
~~~~~~ =
\frac{T-t}{2}\left(\zeta_0^{(i_1)}\zeta_0^{(i_2)}+\sum_{j=1}^{\infty}
\frac{1}{\sqrt{4j^2-1}}\left(
\zeta_{j-1}^{(i_1)}\zeta_{j}^{(i_2)}-
\zeta_j^{(i_1)}\zeta_{j-1}^{(i_2)}\right)\right).
\end{equation}

It is not difficult to see that the relation (\ref{l10}) 
has been obtained in Sect.~5.1 (see (\ref{4004yes})).

Let $\{\phi_j(\tau)\}_{j=0}^{\infty}$ be a complete 
orthonormal system of trigonomertic functions in the space $L_2([t, T]),$
which has the form (\ref{666.6}).

We have
\begin{equation}
\label{rre11}
~~~~~~\int\limits_t^s
\phi_j(\tau)d\tau=
\frac{T-t}{2\pi r}\left\{
\begin{matrix}
\phi_{2r-1}(s),\ &\ j=2r\cr\cr
\sqrt{2}\phi_0(s)-\phi_{2r}(s),\ &\ j=2r-1
\end{matrix}
\right.\ ,
\end{equation}
where $j\ge 1$ and $r=1, 2,\ldots$

From (\ref{l1}) and (\ref{rre11}) we obtain
$$
\int\limits_t^T
{\bf w}_{s,t}^{(i_1)q}
d{\bf w}_s^{(i_2)}=
\frac{1}{\sqrt{T-t}}\zeta_0^{(i_1)}
\int\limits_t^T(s-t){\bf w}_s^{(i_2)}+
$$
$$
+
\frac{T-t}{2}\sum_{r=1}^q
\frac{1}{\pi r}
\left(\left(\zeta_{2r}^{(i_1)}\zeta_{2r-1}^{(i_2)}-
\zeta_{2r-1}^{(i_1)}\zeta_{2r}^{(i_2)}\right)+
\sqrt{2}\zeta_{0}^{(i_2)}\zeta_{2r-1}^{(i_1)}\right)=
$$
$$
=
\frac{1}{\sqrt{T-t}}\zeta_0^{(i_1)}
\frac{(T-t)^{3/2}}{2}\left(\zeta_0^{(i_2)}-
\frac{\sqrt{2}}{\pi}\sum_{r=1}^{\infty}\frac{1}{r}
\zeta_{2r-1}^{(i_2)}\right)+
$$
$$
+
\frac{T-t}{2}\sum_{r=1}^q
\frac{1}{\pi r}
\left(\left(\zeta_{2r}^{(i_1)}\zeta_{2r-1}^{(i_2)}-
\zeta_{2r-1}^{(i_1)}\zeta_{2r}^{(i_2)}\right)+
\sqrt{2}\zeta_{0}^{(i_2)}\zeta_{2r-1}^{(i_1)}\right)=
$$
$$
=\frac{1}{2}(T-t)\Biggl(
\zeta_{0}^{(i_1)}\zeta_{0}^{(i_2)}
+\frac{1}{\pi}
\sum_{r=1}^{q}\frac{1}{r}\biggl(
\zeta_{2r}^{(i_1)}\zeta_{2r-1}^{(i_2)}-
\zeta_{2r-1}^{(i_1)}\zeta_{2r}^{(i_2)}+\biggr.\Biggr.
$$
\begin{equation}
\label{zq1}
~~~~~~+\biggl.\Biggl.
\sqrt{2}\left(\zeta_{2r-1}^{(i_1)}\zeta_{0}^{(i_2)}-
\zeta_{0}^{(i_1)}\zeta_{2r-1}^{(i_2)}\right)\biggr)\Biggr)-
\frac{T-t}{\pi \sqrt{2}}
\zeta_{0}^{(i_1)}\sum_{r=q+1}^{\infty}\frac{1}{r}
\zeta_{2r-1}^{(i_2)}.
\end{equation}

\vspace{2mm}

From (\ref{zq1}) and (\ref{rez1}) we get
$$
I_{(00)T,t}^{(i_1i_2)}=\hbox{\vtop{\offinterlineskip\halign{
\hfil#\hfil\cr
{\rm l.i.m.}\cr
$\stackrel{}{{}_{q\to \infty}}$\cr
}} }
\int\limits_t^T
{\bf w}_{s,t}^{(i_1)q}
d{\bf w}_s^{(i_2)}=
\frac{1}{2}(T-t)\Biggl(
\zeta_{0}^{(i_1)}\zeta_{0}^{(i_2)}
+\Biggr.
$$
\begin{equation}
\label{e987}
~~~~~+\frac{1}{\pi}
\sum_{r=1}^{\infty}\frac{1}{r}\biggl(
\zeta_{2r}^{(i_1)}\zeta_{2r-1}^{(i_2)}-
\zeta_{2r-1}^{(i_1)}\zeta_{2r}^{(i_2)}+\biggr.\Biggr.
\biggl.
\sqrt{2}\left(\zeta_{2r-1}^{(i_1)}\zeta_{0}^{(i_2)}-
\zeta_{0}^{(i_1)}\zeta_{2r-1}^{(i_2)}\right)\biggr)\Biggr),
\end{equation}

\noindent
where $i_1\ne i_2.$

It is obvious that 
(\ref{e987}) is consistent with (\ref{430})
for $i_1\ne i_2$ (we consider here 
(\ref{430}) without the random variables
$\xi_q^{(i)}$).

\section{Milstein method
of Expansion of Iterated It\^{o} and Stratonovich Stochastic Integrals}

The method that is considered in this section was proposed by 
Milstein G.N. \cite{Zapad-1} (1988) and probably until the mid-2000s 
remained one of the most famous methods 
for strong approximation of iterated 
stochastic integrals (also see \cite{Zapad-2}-\cite{Zapad-4},
\cite{Zapad-7b}-\cite{Zapad-9}, \cite{Zapad-11}, \cite{Zapad-12a}). 
However, in light of the results of Chapters 1 and 2
as well as Sect.~5.1 and 5.3, it can be argued that 
the method based on Theorem 1.1 is more general and effective.

The mentioned Milstein method \cite{Zapad-1}
is based on the expansion of the Brownian bridge 
process into 
the trigonometric Fourier series with random coefficients
(version of the so-called Karhunen--Lo\`{e}ve expansion). 

Let us consider the Brownian bridge process
\begin{equation}
\label{6.5.1}
{\bf w}_t-\frac{t}{\Delta}{\bf w}_{\Delta},\ \ \
t\in[0,\Delta],\ \ \ \Delta>0,
\end{equation}
where ${\bf w}_t$ is a standard
Wiener process with independent components
${\bf w}^{(i)}_t,$ $i=1,\ldots,m.$

The componentwise Karhunen--Lo\`{e}ve expansion of the process (\ref{6.5.1}) 
has the following form
\begin{equation}
\label{6.5.2}
~~~~~~~~{\bf w}_t^{(i)}-\frac{t}{\Delta}{\bf w}_{\Delta}^{(i)}=
\frac{1}{2}a_{i,0}+\sum_{r=1}^{\infty}\left(
a_{i,r}{\rm cos}\frac{2\pi rt}{\Delta} + b_{i,r}{\rm sin}
\frac{2\pi rt}{\Delta}\right),
\end{equation}
where the series converges in the mean-square sense and
$$
a_{i,r}=\frac{2}{\Delta} \int\limits_0^{\Delta}
\left({\bf w}_s^{(i)}-\frac{s}{\Delta}{\bf w}_{\Delta}^{(i)}\right)
{\rm cos}\frac{2\pi rs}{\Delta}ds,
$$
$$
b_{i,r}=\frac{2}{\Delta} \int\limits_0^{\Delta}
\left({\bf w}_s^{(i)}-\frac{s}{\Delta}{\bf w}_{\Delta}^{(i)}\right)
{\rm sin}\frac{2\pi rs}{\Delta}ds,
$$
$r=0, 1,\ldots,$ $i=1,\ldots,m.$ 

It is easy to demonstrate \cite{Zapad-1} that the random variables
$a_{i,r}, b_{i,r}$ 
are Gaussian ones and they satisfy the following relations
$$
{\sf M}\left\{a_{i,r}b_{i,r}\right\}=
{\sf M}\left\{a_{i,r}b_{i,k}\right\}=0,\ \ \
{\sf M}\left\{a_{i,r}a_{i,k}\right\}=
{\sf M}\left\{b_{i,r}b_{i,k}\right\}=0,
$$
$$
{\sf M}\left\{a_{i_1,r}a_{i_2,r}\right\}=
{\sf M}\left\{b_{i_1,r}b_{i_2,r}\right\}=0,\ \ \
{\sf M}\left\{a_{i,r}^2\right\}=
{\sf M}\left\{b_{i,r}^2\right\}=\frac{\Delta}{2\pi^2 r^2},
$$

\vspace{2mm}
\noindent
where $i, i_1, i_2=1,\ldots,m,$ $ r\ne k,$ $i_1\ne i_2.$

According to (\ref{6.5.2}), we have
\begin{equation}
\label{6.5.7}
~~~~~~~~{\bf w}_t^{(i)}={\bf w}_{\Delta}^{(i)}\frac{t}{\Delta}+
\frac{1}{2}a_{i,0}+
\sum_{r=1}^{\infty}\left(
a_{i,r}{\rm cos}\frac{2\pi rt}{\Delta}+b_{i,r}{\rm sin}
\frac{2\pi rt}{\Delta}\right),
\end{equation}

\noindent
where the series
converges in the mean-square sense.

Note that the trigonometric functions are the eigenfunctions of the 
covariance 
operator
of the Brownian bridge process. That is why the basis functions are the
trigonometric functions in the considered approach.

Using the relation (\ref{6.5.7}), it is easy to get the 
following expansions \cite{Zapad-1}-\cite{Zapad-3} 

\vspace{-4mm}
\begin{equation}
\label{6.5.8}
~~~~~~~ \int\limits_0^{t}
d{\bf w}_{\tau}^{(i)}=\frac{t}{\Delta}{\bf w}_{\Delta}^{(i)}+
\frac{1}{2}a_{i,0}+\sum_{r=1}^{\infty}\left(
a_{i,r}{\rm cos}\frac{2\pi rt}{\Delta} +b_{i,r}{\rm sin}
\frac{2\pi rt}{\Delta}\right),
\end{equation}

\vspace{2mm}
$$
\int\limits_0^t\int\limits_0^{\tau}
d{\bf w}_{\tau_1}^{(i)}d\tau=\frac{t^2}{2\Delta}{\bf w}_{\Delta}^{(i)}+
\frac{t}{2}a_{i,0}+
$$
\begin{equation}
\label{6.5.9}
+\frac{\Delta}{2\pi}\sum_{r=1}^{\infty}
\frac{1}{r}\left(
a_{i,r}{\rm sin}\frac{2\pi rt}{\Delta} 
- b_{i,r}\left({\rm cos}
\frac{2\pi rt}{\Delta}-1\right)\right),
\end{equation}

\newpage
\noindent
$$
\int\limits_0^{t}\int\limits_0^{\tau}
d\tau_1 d{\bf w}_{\tau}^{(i)}=
t\int\limits_0^{t}
d{\bf w}_t^{(i)}-\int\limits_0^{t}\int\limits_0^{\tau}
d{\bf w}_{\tau_1}^{(i)}d\tau=
\frac{t^2}{2\Delta}{\bf w}_{\Delta}^{(i)}+
$$
$$
+ t \sum_{r=1}^{\infty}
\left(
a_{i,r}{\rm cos}\frac{2\pi rt}{\Delta}
+b_{i,r}{\rm sin}\frac{2\pi rt}{\Delta}\right)-
$$
\begin{equation}
\label{6.5.10}
-\frac{\Delta}{2\pi}\sum_{r=1}^{\infty}
\frac{1}{r}\left(
a_{i,r}{\rm sin}\frac{2\pi rt}{\Delta}-
b_{i,r}\left({\rm cos}
\frac{2\pi rt}{\Delta}-1\right)\right),
\end{equation}

\vspace{8mm}

$$
\int\limits_0^{t}\int\limits_0^{\tau}
d{\bf w}_{\tau_1}^{(i_1)}d{\bf w}_{\tau}^{(i_2)}=
\frac{1}{\Delta}{\bf w}_{\Delta}^{(i_1)}
\int\limits_0^{t}\int\limits_0^{\tau}
d\tau_1 d{\bf w}_{\tau}^{(i_2)} +
\frac{1}{2}a_{i_1,0}
\int\limits_0^{t}
d{\bf w}_{\tau}^{(i_2)}
+
$$
$$
+\frac{t\pi}{\Delta}\sum_{r=1}^{\infty} r
\left(a_{i_1,r}b_{i_2,r}-
b_{i_1,r}a_{i_2,r}\right)+
$$
$$
+
\frac{1}{4}\sum_{r=1}^{\infty} \Biggl(
\left(a_{i_1,r}a_{i_2,r}-b_{i_1,r}b_{i_2,r}\right)
\left(1-{\rm cos}\frac{4\pi rt}{\Delta}\right)+\Biggr.
$$
$$
+\left(a_{i_1,r}b_{i_2,r}+
b_{i_1,r}a_{i_2,r}\right)
{\rm sin}\frac{4\pi rt}{\Delta}+
$$
$$
\Biggl.+\frac{2}{\pi r}{\bf w}_{\Delta}^{(i_2)}\left(
a_{i_1,r}{\rm sin}\frac{2\pi rt}{\Delta} +b_{i_1,r}\left(
{\rm cos}\frac{2\pi rt}{\Delta}-1\right)\right)\Biggr)+
$$
$$
+\sum_{k=1}^{\infty}\sum_{r=1 (r\ne k)}^{\infty}
k\left(a_{i_1,r}a_{i_2,k}
\left(\frac{{\rm cos}\left(\frac{2\pi(k+r)t}{\Delta}\right)}{2(k+r)}+
\frac{{\rm cos}\left(\frac{2\pi(k-r)t}{\Delta}\right)}{2(k-r)}-
\frac{k}{k^2-r^2}\right)+\right.
$$
$$
+a_{i_1,r}b_{i_2,k}
\left(\frac{{\rm sin}\left(\frac{2\pi(k+r)t}{\Delta}\right)}{2(k+r)}+
\frac{{\rm sin}\left(\frac{2\pi(k-r)t}{\Delta}\right)}{2(k-r)}\right)+
$$
$$
+b_{i_1,r}b_{i_2,k}
\left(\frac{{\rm cos}\left(\frac{2\pi(k-r)t}{\Delta}\right)}{2(k-r)}-
\frac{{\rm cos}\left(\frac{2\pi(k+r)t}{\Delta}\right)}{2(k+r)}-
\frac{r}{k^2-r^2}\right)+
$$
\begin{equation}
\label{6.5.11}
\left.+\frac{\Delta}{2\pi}b_{i_1,r}a_{i_2,k}
\left(\frac{{\rm sin}\left(\frac{2\pi(k+r)t}{\Delta}\right)}{2(k+r)}-
\frac{{\rm sin}\left(\frac{2\pi(k-r)t}{\Delta}\right)}{2(k-r)}\right)\right)
\end{equation}

\newpage
\noindent
converging in the mean-square sense, where we suppose that $i_1\ne i_2$
in (\ref{6.5.11}).

It is necessary to pay a special attention to the fact that the 
double series in (\ref{6.5.11}) should be understood 
as the iterated one, and not as a multiple series 
(as in Theorem 1.1), 
i.e. as the iterated passage to the limit for the sequence 
of double partial sums.
So, the Milstein method of approximation of iterated 
stochastic integrals \cite{Zapad-1} leads to iterated application of the
limit transition (in contrast with the method 
of generalized multiple Fourier series (Theorem 1.1), for which the 
limit transition is implemented only once) 
starting at least from the second or third multiplicity of 
iterated stochastic integrals (we mean at least double or 
triple integration with respect to components
of the Wiener process).
Multiple series are more 
preferential for approximation than the iterated ones, since the
partial sums of multiple series converge for any possible case of joint 
converging to infinity of their upper limits of summation (let us denote 
them as $p_1,\ldots, p_k$). 
For example,
when $p_1=\ldots=p_k=p\to\infty$. 
For iterated series, the condition $p_1=\ldots=p_k=p\to\infty$ obviously 
does not guarantee the convergence of this series.
However, as we will see further in this section 
in \cite{Zapad-2} (pp.~438-439), \cite{Zapad-3}
(Sect.~5.8, pp.~202--204), \cite{Zapad-4} (pp.~82-84), 
\cite{Zapad-9} (pp.~263-264)
the authors use (without rigorous proof) the condition 
$p_1=p_2=p_3=p\to\infty$
within the frames of the Milstein method \cite{Zapad-1}
together with the Wong--Zakai approximation \cite{W-Z-1}-\cite{Watanabe}
(also see discussions in Sect.~2.42, 2.43).
Furthermore,
in order to obtain the Milstein expansion for
iterated stochastic integral, the truncated 
expansions (\ref{6.5.7}) of components of the Wiener 
process ${\bf w}_t$ must be
iteratively substituted in the single integrals, and the integrals
must be calculated, starting from the innermost integral.
This is a complicated procedure that obviously does not lead to the
expansion of iterated stochastic integral 
of multiplicity $k$ $(k\in {\bf N})$.

Assume that $t=\Delta$ in  (\ref{6.5.8})--(\ref{6.5.11}) 
(at that 
double 
partial sums 
of iterated series in (\ref{6.5.11}) will become zero).
As a result, we get the expansions
\begin{equation}
\label{6.5.12}
\int\limits_0^{\Delta}
d{\bf w}_{\tau}^{(i)}={\bf w}_{\Delta}^{(i)},
\end{equation}
\begin{equation}
\label{6.5.13}
\int\limits_0^{\Delta}\int\limits_0^{\tau}
d{\bf w}_{\tau_1}^{(i)}d\tau=\frac{1}{2}\Delta\left(
{\bf w}_{\Delta}^{(i)}+a_{i,0}\right),
\end{equation}
\begin{equation}
\label{6.5.14}
\int\limits_0^{\Delta}\int\limits_0^{\tau}
d\tau_1 d{\bf w}_{\tau}^{(i)}=
\frac{1}{2}\Delta\left(
{\bf w}_{\Delta}^{(i)}-a_{i,0}\right),
\end{equation}

\newpage
\noindent
$$
\int\limits_0^{\Delta}\int\limits_0^{\tau}
d{\bf w}_{\tau_1}^{(i_1)}d{\bf w}_{\tau}^{(i_2)}=
\frac{1}{2}{\bf w}_{\Delta}^{(i_1)}{\bf w}_{\Delta}^{(i_2)}-
\frac{1}{2}\left(a_{i_2,0}{\bf w}_{\Delta}^{(i_1)}-
a_{i_1,0}{\bf w}_{\Delta}^{(i_2)}\right)+
$$
\begin{equation}
\label{6.5.15}
+
\pi\sum_{r=1}^{\infty} r
\left(a_{i_1,r}b_{i_2,r}-b_{i_1,r}a_{i_2,r}\right)
\end{equation}

\noindent
that converge
in the mean-square sense,
where we suppose that $i_1\ne i_2$ in (\ref{6.5.15}).

Deriving (\ref{6.5.12})--(\ref{6.5.15}),
we used the relation
\begin{equation}
\label{6.5.16}
a_{i,0}=-2\sum_{r=1}^{\infty}a_{i,r},
\end{equation}
which results from (\ref{6.5.2}) when $t=\Delta.$

Let us compare expansions of some iterated stochastic
integrals of 
first and second multiplicity obtained by Milstein method \cite{Zapad-1}
and method based on generalized multiple Fourier series (Theorem 1.1).

Let us denote
\begin{equation}
\label{df0}
~~~~~~~ \zeta_{2r-1}^{(i)}=\sqrt{\frac{2}{\Delta}}\int\limits_0^{\Delta}
{\rm sin}\frac{2\pi r s}{\Delta}d{\bf w}_{s}^{(i)},\ \ \
\zeta_{2r}^{(i)}=\sqrt{\frac{2}{\Delta}}\int\limits_0^{\Delta}
{\rm cos}\frac{2\pi r s}{\Delta}d{\bf w}_{s}^{(i)},
\end{equation}
\begin{equation}
\label{df01}
\zeta_{0}^{(i)}=\frac{1}{\sqrt{\Delta}}\int\limits_0^{\Delta}
d{\bf w}_{s}^{(i)},
\end{equation}
where $r=1, 2,\ldots,$ $i=1,\ldots,m.$

Using the It\^{o} formula, it is not difficult to show that
\begin{equation}
\label{df1wwww}
a_{i,r}=-\frac{1}{\pi r}\sqrt{\frac{\Delta}{2}}\zeta_{2r-1}^{(i)},\ \ \
b_{i,r}=\frac{1}{\pi r}\sqrt{\frac{\Delta}{2}}\zeta_{2r}^{(i)}\ \ \
\hbox{w.~p.~1.}
\end{equation}

From (\ref{6.5.16}) we get
\begin{equation}
\label{df2wwww}
a_{i,0}=\frac{\sqrt{2\Delta}}{\pi }\sum\limits_{r=1}^{\infty}
\frac{1}{r}\zeta_{2r-1}^{(i)}.
\end{equation}

After substituting (\ref{df1wwww}), (\ref{df2wwww}) into 
(\ref{6.5.12})--(\ref{6.5.15}) and taking into account
(\ref{df0}), (\ref{df01}), we have
\begin{equation}
\label{41a}
\int\limits_0^{\Delta}
d{\bf w}_{\tau}^{(i_1)}=\sqrt{\Delta}\zeta_0^{(i_1)},
\end{equation}
\begin{equation}
\label{42a}
\int\limits_0^{\Delta}\int\limits_0^{\tau}
d\tau_1 d{\bf w}_{\tau}^{(i_1)}
=\frac{{\Delta}^{3/2}}{2}
\Biggl(\zeta_0^{(i_1)}-\frac{\sqrt{2}}{\pi}\sum_{r=1}^{\infty}
\frac{1}{r}\zeta_{2r-1}^{(i_1)}\Biggr),
\end{equation}
\begin{equation}
\label{42aa}
\int\limits_0^{\Delta}\int\limits_0^{\tau}
d{\bf w}_{\tau_1}^{(i_1)}d\tau
=\frac{{\Delta}^{3/2}}{2}
\Biggl(\zeta_0^{(i_1)}+\frac{\sqrt{2}}{\pi}\sum_{r=1}^{\infty}
\frac{1}{r}\zeta_{2r-1}^{(i_1)}\Biggr),
\end{equation}

$$
\int\limits_0^{\Delta}\int\limits_0^{\tau}
d{\bf w}_{\tau_1}^{(i_1)}d{\bf w}_{\tau}^{(i_2)}=\frac{\Delta}{2}\Biggl(
\zeta_{0}^{(i_1)}\zeta_{0}^{(i_2)}
+\frac{1}{\pi}
\sum_{r=1}^{\infty}\frac{1}{r}\left(
\zeta_{2r}^{(i_1)}\zeta_{2r-1}^{(i_2)}-
\zeta_{2r-1}^{(i_1)}\zeta_{2r}^{(i_2)}+
\right.\Biggr.
$$
\begin{equation}
\label{43xx0}
+\left.\sqrt{2}\left(\zeta_{2r-1}^{(i_1)}\zeta_{0}^{(i_2)}-
\zeta_{0}^{(i_1)}\zeta_{2r-1}^{(i_2)}\right)\right)
\Biggl.
\Biggr).
\end{equation}

Obviously, the formulas (\ref{41a})--(\ref{43xx0}) 
are consistent with the formulas (\ref{4001}), (\ref{430}),
(\ref{df5}), (\ref{df6}).
It testifies that at least for 
the considered iterated stochastic integrals and 
trigonometric system of functions, the Milstein method and 
the method based on generalized multiple Fourier series 
(Theorem 1.1) give the same result 
(it is an interesting fact, although it is rather expectable).

Further, we will 
discuss the usage of Milstein method
for the iterated sto\-chas\-tic integrals of third multiplicity.

First, we note that the 
authors of the monograph \cite{Zapad-3}
based on the results of Wong E. and Zakai M. \cite{W-Z-1}, \cite{W-Z-2}
(also see \cite{Watanabe}) 
concluded (without rigorous proof) 
that 
the expansions of iterated stochastic
integrals on the basis of (\ref{6.5.7}) (the case $i_1,i_2,i_3=1,\ldots,m$)
converge to the iterated Stratonovich stochastic integrals
(see discussions in Sect.~2.42, 2.43).
It is obvious that this conclusion is consistent with the 
results given above in this section for the case $i_1\ne i_2.$

As we mentioned before, the technical peculiarities 
of the Milstein method \cite{Zapad-1} may result 
to the iterated series 
of products
of standard Gaussian random variables
(in contradiction to multiple series 
as in Theorem 1.1).
In the case of simplest 
stochastic integral of second multiplicity this problem was avoided 
as we saw above.  However, the situation is not 
the same for the simplest stochastic integrals of third multiplicity.

Let us denote
$$
{J}_{(\lambda_{1}\ldots \lambda_k)T,t}^{*(i_1\ldots
i_k)}=
{\int\limits_t^{*}}^T
\ldots
{\int\limits_t^{*}}^{t_{2}}
d{\bf w}_{t_{1}}^{(i_1)}\ldots
d{\bf w}_{t_k}^{(i_k)},
$$
where $\lambda_l=1$ if 
$i_l=1,\ldots,m$ and 
$\lambda_l=0$ if $i_l=0,$ $l=1,\ldots,k,$
${\bf w}_{\tau}^{(i)}$
$(i=1,\ldots,m)$ are independent standard Wiener processes,
${\bf w}_{\tau}^{(0)}=\tau$.

Let us consider the expansion of iterated Stratonovich stochastic
integral of third multiplicity obtained in \cite{Zapad-2}-\cite{Zapad-4}, 
\cite{Zapad-9}
by the Milstein method \cite{Zapad-1}

\vspace{-1mm}
$$
J_{(111)\Delta,0}^{*(i_1 i_2 i_3)}
=\frac{1}{\Delta}
J_{(1)\Delta,0}^{*(i_1)}J_{(011)\Delta,0}^{*(0 i_2 i_3)}
+
$$

$$
+\frac{1}{2}a_{i_1,0}J_{(11)\Delta,0}^{*(i_2 i_3)}+
\frac{1}{2\pi}b_{i_1}J_{(1)\Delta,0}^{(i_2)}
J_{(1)\Delta,0}^{*(i_3)}-\Delta J_{(1)\Delta,0}^{*(i_2)}B_{i_1 i_3}+
$$

\begin{equation}
\label{starr}
+\Delta J_{(1)\Delta,0}^{*(i_3)}
\left(\frac{1}{2}A_{i_1i_2}-C_{i_2i_1}\right)
+\Delta^{3/2}D_{i_1i_2i_3},
\end{equation}

\vspace{3mm}
\noindent
where
$$
J_{(011)\Delta,0}^{*(0 i_2 i_3)}=\frac{1}{6}
J_{(1)\Delta,0}^{*(i_2)}
J_{(1)\Delta,0}^{*(i_3)}-\frac{1}{\pi}\Delta
J_{(1)\Delta,0}^{*(i_3)}b_{i_2}+
$$

$$
+\Delta^2 B_{i_2i_3}-\frac{1}{4}\Delta a_{i_3,0}
J_{(1)\Delta,0}^{*(i_2)}+
\frac{1}{2\pi}\Delta b_{i_3}
J_{(1)\Delta,0}^{*(i_2)}
+\Delta^2 C_{i_2i_3}+\frac{1}{2}\Delta^2 A_{i_2i_3},
$$

\vspace{2mm}

$$
A_{i_2i_3}=\frac{\pi}{\Delta}\sum_{r=1}^{\infty}
r\left(a_{i_2,r}b_{i_3,r}-
b_{i_2,r}
a_{i_3,r}\right),
$$

$$
C_{i_2i_3}=-\frac{1}{\Delta}\sum_{l=1}^{\infty}
\sum_{r=1 (r\ne l)}^{\infty}
\frac{r}{r^2-l^2}\left(ra_{i_2,r}a_{i_3,l}+
lb_{i_2,r}
b_{i_3,l}\right),
$$

$$
B_{i_2i_3}=\frac{1}{2\Delta}\sum_{r=1}^{\infty}
\left(a_{i_2,r}a_{i_3,r}+
b_{i_2,r}
b_{i_3,r}\right),\
b_i=\sum_{r=1}^{\infty}\frac{1}{r}b_{i,r},
$$

\vspace{3mm}

$$
D_{i_1i_2i_3}=-\frac{\pi}{2\Delta^{3/2}}
\sum_{l=1}^{\infty}\sum_{r=1}^{\infty}
l
\Biggl(a_{i_2,l}
\left(a_{i_3,l+r}b_{i_1,r}-a_{i_1,r}
b_{i_3,l+r}\right)
+\Biggr.
$$
$$
\Biggl.
+
b_{i_2,l}
\left(a_{i_1,r}a_{i_3,r+l}+b_{i_1,r}
b_{i_3,l+r}\right)\Biggr)+
$$
$$
+\frac{\pi}{2\Delta^{3/2}}
\sum_{l=1}^{\infty}\sum_{r=1}^{l-1}
l
\Biggl(a_{i_2,l}
\left(a_{i_1,r}b_{i_3,l-r}+a_{i_3,l-r}
b_{i_1,r}\right)
-\Biggr.
$$
$$
\Biggl.
-b_{i_2,l}
\left(a_{i_1,r}a_{i_3,l-r}-b_{i_1,r}
b_{i_3,l-r}\right)\Biggr)+
$$
$$
+\frac{\pi}{2\Delta^{3/2}}
\sum_{l=1}^{\infty}\sum_{r=l+1}^{\infty}
l
\Biggl(a_{i_2,l}
\left(a_{i_3,r-l}b_{i_1,r}-a_{i_1,r}
b_{i_3,r-l}\right)
+\Biggr.
$$
$$
\Biggl.
+
b_{i_2,l}
\left(a_{i_1,r}a_{i_3,r-l}+b_{i_1,r}
b_{i_3,r-l}\right)\Biggr).
$$

\vspace{1mm}

From 
the expansion
(\ref{starr}) and expansion 
of the stochastic integral
$J_{(011)\Delta,0}^{*(0i_2 i_3)}$
we can conclude
that they include 
iterated (double)
series. 
Moreover, for approximation of the                      
stochastic integral $J_{(111)\Delta,0}^{*(i_1 i_2 i_3)}$ in the works
\cite{Zapad-2} (pp.~438-439), \cite{Zapad-3}
(Sect.~5.8, pp.~202--204), \cite{Zapad-4} (pp.~82-84), 
\cite{Zapad-9} (pp.~263-264)
it is proposed to put upper limits of summation 
by equal $q$
(on the base of the Wong--Zakai approximation \cite{W-Z-1}-\cite{Watanabe}
but without rigorous proof;
also see discussions in Sect.~2.42, 2.43).

For example, the value $D_{i_1i_2i_3}$ is approximated 
in \cite{Zapad-2}-\cite{Zapad-4}, \cite{Zapad-9}
by the double sums of the form
$$
D_{i_1i_2i_3}^{q}=-\frac{\pi}{2\Delta^{3/2}}
\sum_{l=1}^{q}\sum_{r=1}^{q}
l
\Biggl(a_{i_2,l}
\left(a_{i_3,l+r}b_{i_1,r}-a_{i_1,r}
b_{i_3,l+r}\right)
+\Biggr.
$$
$$
\Biggl.
+
b_{i_2,l}
\left(a_{i_1,r}a_{i_3,r+l}+b_{i_1,r}
b_{i_3,l+r}\right)\Biggr)+
$$
$$
+\frac{\pi}{2\Delta^{3/2}}
\sum_{l=1}^{q}\sum_{r=1}^{l-1}
l
\Biggl(a_{i_2,l}
\left(a_{i_1,r}b_{i_3,l-r}+a_{i_3,l-r}
b_{i_1,r}\right)
-\Biggr.
$$
$$
\Biggl.
-b_{i_2,l}
\left(a_{i_1,r}a_{i_3,l-r}-b_{i_1,r}
b_{i_3,l-r}\right)\Biggr)+
$$
$$
+\frac{\pi}{2\Delta^{3/2}}
\sum_{l=1}^{q}\sum_{r=l+1}^{2q}
l
\Biggl(a_{i_2,l}
\left(a_{i_3,r-l}b_{i_1,r}-a_{i_1,r}
b_{i_3,r-l}\right)
+\Biggr.
$$
$$
\Biggl.
+
b_{i_2,l}
\left(a_{i_1,r}a_{i_3,r-l}+b_{i_1,r}
b_{i_3,r-l}\right)\Biggr).
$$

\vspace{2mm}

We can avoid the mentioned problem (iterated application
of the operation of limit transition)  
using the method based on Theorems 1.1, 2.1--2.9, 2.33--2.36, 2.50, 2.51, 2.62, 2.63.

From the other hand, if we prove that the members of 
the expansion (\ref{starr}) coincide 
with the members of its analogue
obtained using Theorem 1.1,
then we can replace the iterated series in (\ref{starr}) by 
the multiple 
series (see Theorems 1.1, 2.1--2.9, 2.33--2.36, 2.50, 2.51, 2.62, 2.63) 
as was made formally 
in \cite{Zapad-2}-\cite{Zapad-4}, \cite{Zapad-9}. 
However, it 
requires the separate argumentation.

\section{Usage of Integral Sums for Approximation of 
Iterated It\^{o} Stochastic Integrals}

It should be noted that there is an approach to 
the mean-square approximation of iterated stochastic
integrals based on multiple integral sums (see, for example, 
\cite{Zapad-1}, \cite{Zapad-8},
\cite{Zapad-10}, \cite{rupp}). This method implies 
the partitioning
of the integration interval $[t, T]$ of the iterated stochastic 
integral under consideration; this interval is the
integration step of the numerical methods used to solve It\^{o} SDEs
(see Chapter 4); therefore, it is already fairly small
and does not need to be partitioned. 
Computational experiments \cite{1} (also see below in this section) 
show that the application of the
method \cite{Zapad-1}, \cite{Zapad-8},
\cite{Zapad-10}, \cite{rupp}
to stochastic integrals with multiplicities $k\ge 2$ 
leads to unacceptably high computational
cost and accumulation of computation errors.

As we noted in the introduction to this book, considering the modern state 
of question on the approximation of iterated stochastic integrals, 
the method analyzed in this section is hardly important
for practice. 
However, we will consider this method in order to get the overall view.
In this section, we will analyze one of the simplest 
modifications of the mentioned method.

Let the functions $\psi_l(\tau),$ $l=1,\ldots,k$ satisfy 
the Lipschitz condition at the interval $[t, T]$ with constants
$C_l$
\begin{equation}
\label{321.001}
|\psi_l(\tau_1)-\psi_l(\tau_2)|\le C_l|\tau_1-\tau_2|\ \ \
\hbox{for all}\ \ \ \tau_1,\tau_2\in [t,T].
\end{equation}

Then, according to Lemma 1.1 (see Sect.~1.1.3), 
the following 
equality is correct
$$
J[\psi^{(k)}]_{T,t}=
\hbox{\vtop{\offinterlineskip\halign{
\hfil#\hfil\cr
{\rm l.i.m.}\cr
$\stackrel{}{{}_{N\to \infty}}$\cr
}} }
\sum_{j_k=0}^{N-1} 
\ldots \sum_{j_1=0}^{j_{2}-1}
\prod_{l=1}^k
\psi_l(\tau_{j_l})\Delta{\bf w}_{\tau_
{j_l}}^{(i_l)}\ \ \ \hbox{w.~p.~1,}
$$

\noindent
where notations are the same as in 
(\ref{30.30}).

Let us consider the following approximation 
\begin{equation}
\label{321.002}
J[\psi^{(k)}]_{T,t}^N=
\sum_{j_k=0}^{N-1} 
\ldots \sum_{j_1=0}^{j_{2}-1}
\prod_{l=1}^k
\psi_l(\tau_{j_l})\Delta{\bf w}_{\tau_{j_l}}^{(i_l)}
\end{equation}
of the iterated
It\^{o} stochastic integral $J[\psi^{(k)}]_{T,t}$. 
The relation (\ref{321.002}) can be re\-writ\-ten as
\begin{equation}
\label{321.003}
J[\psi^{(k)}]_{T,t}^N=
\sum_{j_k=0}^{N-1} 
\ldots \sum_{j_1=0}^{j_{2}-1}
\prod_{l=1}^k
\sqrt{\Delta \tau_{j_l}}
\psi_l(\tau_{j_l}) {\bf u}_{j_l}^{(i_l)},
\end{equation}

\noindent
where ${\bf u}_j^{(i)}\stackrel{\rm def}{=}\left({\bf w}_{\tau_{j+1}}^{(i)}-
{\bf w}^{(i)}_{\tau_j}\right)/\sqrt{\Delta\tau_j}$,\ \  
$i=1,\ldots,m$ are independent
standard Gaussian random variables
for various $i$ or $j$,
${\bf u}_j^{(0)}=\sqrt{\Delta\tau_j}.$

Assume that
\begin{equation}
\label{321.005}
\tau_j=t+j\Delta,\ \ \ j=0,1,\ldots,N,\ \ \ \tau_N=T,\ \ \ \Delta>0.
\end{equation}

Then 
\begin{equation}
\label{321.006}
J[\psi^{(k)}]_{T,t}^N=\Delta^{k/2}
\sum_{j_k=0}^{N-1}
\ldots \sum_{j_1=0}^{j_{2}-1}
\prod_{l=1}^k
\psi_l(t+j_l\Delta) {\bf u}_{j_l}^{(i_l)},
\end{equation}

\noindent
where
${\bf u}_j^{(i)}\stackrel{\rm def}{=}\left({\bf w}_{t+(j+1)\Delta}^{(i)}-
{\bf w}^{(i)}_{t+j\Delta}\right)/\sqrt{\Delta}$,\ \ 
$i=1,\ldots,m,$ ${\bf u}_j^{(0)}=\sqrt{\Delta}.$

{\bf Lemma 6.1.} {\it Suppose that the functions 
$\psi_l(\tau),$ $l=1,\ldots k$
satisfy the Lipschitz condition
{\rm (\ref{321.001})} and 
$\left\{\tau_j\right\}_{j=0}^{N}$ is a partition 
of the interval $[t, T],$ which satisfies the condition
{\rm (\ref{321.005})}.
Then
for 
a sufficiently small value $T-t$
there 
exists 
a constant $H_k<\infty$ such that

\vspace{-2mm}
$$
{\sf M}\left\{\left(J[\psi^{(k)}]_{T,t}-J[\psi^{(k)}]_{T,t}^N\right)^2
\right\}\le
\frac{H_k(T-t)^2}{N}.
$$
}

\par
{\bf Proof.} It is easy to see that in the case of a sufficiently small 
value $T-t$ there 
exists 
a constant $L_k$ such that
$$
{\sf M}\left\{\left(J[\psi^{(k)}]_{T,t}-J[\psi^{(k)}]_{T,t}^N\right)^2
\right\}\le L_k
{\sf M}\left\{\left(J[\psi^{(2)}]_{T,t}-J[\psi^{(2)}]_{T,t}^N\right)^2
\right\},
$$
where
$$
J[\psi^{(2)}]_{T,t}-J[\psi^{(2)}]_{T,t}^N=
\sum_{j=1}^3 S_j^N,
$$
$$
S_1^N=\sum_{j_1=0}^{N-1}\int\limits_{\tau_{j_1}}^{\tau_{j_1+1}}
\psi_2(t_2)\int\limits_{\tau_{j_1}}^{t_2}\psi_1(t_1)d{\bf w}_{t_1}^{(i_1)}
d{\bf w}_{t_2}^{(i_2)},
$$
$$
S_2^N=\sum_{j_1=0}^{N-1}\int\limits_{\tau_{j_1}}^{\tau_{j_1+1}}
\left(\psi_2(t_2)-\psi_2(\tau_{j_1})\right)d{\bf w}_{t_2}^{(i_2)}
\sum_{j_2=0}^{j_1-1}
\int\limits_{\tau_{j_2}}^{\tau_{j_2+1}}
\psi_1(t_1)d{\bf w}_{t_1}^{(i_1)},
$$
$$
S_3^N=
\sum_{j_1=0}^{N-1}
\psi_2(\tau_{j_1}) \Delta{\bf w}_{\tau_{j_1}}^{(i_2)}
\sum_{j_2=0}^{j_1-1}\int\limits_{\tau_{j_2}}^{\tau_{j_2+1}}
\left(\psi_1(t_1)-\psi_1(\tau_{j_2})\right)d{\bf w}_{t_1}^{(i_1)}.
$$

Therefore, according to the Minkowski inequality, we have
$$
\left({\sf M}\left\{\left(J[\psi^{(2)}]_{T,t}-J[\psi^{(2)}]_{T,t}^N\right)^2
\right\}\right)^{1/2}\le
\sum_{j=1}^3 \biggl({\sf M}\left\{\left(S_j^N\right)^2\right\}
\biggr)^{1/2}.
$$

Using standard  moment properties of stochastic integrals 
(see (\ref{99.010}), (\ref{99.010a})), let us estimate
the values 
${\sf M}\left\{
\left(S_j^N\right)^2\right\},$ $j=1, 2, 3.$

Let us consider four cases.

Case 1. $i_1, i_2\ne 0:$
$$
{\sf M}\left\{\left(
S_1^N\right)^2\right\}\le
\frac{\Delta}{2}(T-t)
\max\limits_{s\in[t,T]}
\psi^2_2(s)\psi^2_1(s),
$$
$$
{\sf M}\left\{\left(
S_2^N\right)^2\right\}\le
\frac{\Delta^2}{6}(T-t)^2 \left(C_2\right)^2
\max\limits_{s\in[t,T]}\psi^2_1(s),
$$
$$
{\sf M}\left\{\left(
S_3^N\right)^2\right\}\le
\frac{\Delta^2}{6}(T-t)^2 \left(C_1\right)^2
\max\limits_{s\in[t,T]}\psi^2_2(s).
$$

Case 2. $i_1\ne 0,$ $i_2=0:$
$$
{\sf M}\left\{\left(
S_1^N\right)^2\right\}\le
\frac{\Delta}{2}(T-t)^2
\max\limits_{s\in[t,T]}\psi^2_2(s)\psi^2_1(s),
$$
$$
{\sf M}\left\{\left(
S_2^N\right)^2\right\}\le
\frac{\Delta^2}{3}(T-t)^3 \left(C_2\right)^2
\max\limits_{s\in[t,T]}\psi^2_1(s),
$$
$$
{\sf M}\left\{\left(
S_3^N\right)^2\right\}\le
\frac{\Delta^2}{3}(T-t)^3 \left(C_1\right)^2
\max\limits_{s\in[t,T]}\psi^2_2(s).
$$

Case 3. $i_2\ne 0,$ $i_1=0:$
$$
{\sf M}\left\{\left(
S_1^N\right)^2\right\}\le
\frac{\Delta^2}{3}(T-t)
\max\limits_{s\in[t,T]}\psi^2_2(s)\psi^2_1(s),
$$
$$
{\sf M}\left\{\left(
S_2^N\right)^2\right\}\le
\frac{\Delta^2}{3}(T-t)^3 \left(C_2\right)^2
\max\limits_{s\in[t,T]}\psi^2_1(s),
$$
$$
{\sf M}\left\{\left(
S_3^N\right)^2\right\}\le
\frac{\Delta^2}{8}(T-t)^3 \left(C_1\right)^2
\max\limits_{s\in[t,T]}\psi^2_2(s).
$$

Case 4. $i_1=i_2=0:$
$$
{\sf M}\left\{\left(
S_1^N\right)^2\right\}\le
\frac{\Delta^2}{4}(T-t)^2
\max\limits_{s\in[t,T]}\psi^2_2(s)\psi^2_1(s),
$$
$$
{\sf M}\left\{\left(
S_2^N\right)^2\right\}\le
\frac{\Delta^2}{4}(T-t)^4 \left(C_2\right)^2
\max\limits_{s\in[t,T]}\psi^2_1(s),
$$
$$
{\sf M}\left\{\left(
S_3^N\right)^2\right\}\le
\frac{\Delta^2}{16}(T-t)^4 \left(C_1\right)^2
\max\limits_{s\in[t,T]}\psi^2_2(s).
$$

\vspace{2mm}

According to the obtained estimates, we have
$$
{\sf M}\left\{\left(J[\psi^{(k)}]_{T,t}-J[\psi^{(k)}]_{T,t}^N\right)^2\right\}
\le H_k(T-t)\Delta=\frac{H_k(T-t)^2}{N},
$$
where $H_k<\infty.$ Lemma 6.1 is proved. 

It is easy to check that the following relation is correct
\begin{equation}
\label{321.017}
{\sf M}\left\{\left(I_{(00)T,t}^{(i_1i_2)}-
I_{(00)T,t}^{(i_1i_2)N}\right)^2\right\}=\frac{(T-t)^2}{2N},
\end{equation}
where $i_1,i_2=1,\ldots,m$
and $I_{(00)T,t}^{(i_1i_2)N}$
is the approximation of the iterated stochastic 
integral $I_{(00)T,t}^{(i_1i_2)}$
(see (\ref{ito-kuku})) obtained according to the formula 
(\ref{321.006}).

Finally, we will demonstrate that 
the 
method based on generalized multiple Fourier series (Theorem 1.1) is 
signficiantly better,
than the 
method based on multiple integral sums
in the sense of computational costs 
on modeling of iterated stochastic integrals.
   
Let us consider the approximations of iterated It\^{o} stochastic integrals
obtained using the method based on multiple integral sums
\begin{equation}
\label{bob30}
I_{(0)T,t}^{(1)q}=\sqrt{\Delta}\sum_{j=0}^{q-1}\xi_{j}^{(1)},
\end{equation}
\begin{equation}
\label{bob31}
I_{(00)T,t}^{(12)q}=\Delta\sum_{j=0}^{q-1}
\left(\sum_{i=0}^{j-1}\xi_i^{(1)}\right) \xi_{j}^{(2)},
\end{equation}
where 
$$
\xi_j^{(i)}=\left({\bf w}^{(i)}_{t+(j+1)\Delta}-{\bf w}^{(i)}_{t+j\Delta}
\right)/
\sqrt{\Delta},\ \ \ i=1, 2
$$
are independent standard  
Gaussian random variables, $\Delta=(T-t)/q,$
$I_{(00)T,t}^{(12)q},$ $I_{(0)T,t}^{(1)q}$ are approximations
of the iterated It\^{o} stochastic 
integrals $I_{(00)T,t}^{(12)}$ (see (\ref{ito-kuku})),
$I_{(0)T,t}^{(1)}={\bf w}_T^{(1)}-{\bf w}_t^{(1)}$.

Let us choose the number $q$ 
(see (\ref{bob30}), (\ref{bob31}))
from the condition 
$$
{\sf M}\left\{\left(I_{(00)T,t}^{(12)}-I_{(00)T,t}^{(12)q}\right)^2
\right\}=\frac{(T-t)^2}{2q}\le
(T-t)^3.
$$

Let us implement 200 independent 
numerical modelings
of the collection of iterated It\^{o} stochastic 
integrals $I_{(00)T,t}^{(12)},$ $I_{(0)T,t}^{(1)}$ using 
the formulas (\ref{bob30}), (\ref{bob31})
for $T-t=2^{-j},$ $j=5, 6, 7.$ 
We denote by $T_{\rm sum}$ the computer time 
which is necessary for performing this task.

Let us repeat the above experiment for the case
when the approximations of iterated It\^{o}
stochastic integrals $I_{(00)T,t}^{(12)},$ $I_{(0)T,t}^{(1)}$
are defined by (\ref{y2}), (\ref{y3}) and the number $q$
is chosen from the condition 
(\ref{uslov1}) (method based on Theorem 1.1, the case of Legendre
polynomials).
Let $T_{\rm pol}$ be the computer time 
which is necessary for performing this task.

\begin{table}
\centering
\caption{Values $T_{\rm sum}/T_{\rm pol}$.}
\label{tab:6.1}      
\begin{tabular}{p{1.7cm}p{1.7cm}p{1.7cm}p{1.7cm}}
\hline\noalign{\smallskip}
$T-t$&$2^{-5}$&$2^{-6}$&$2^{-7}$\\
\noalign{\smallskip}\hline\noalign{\smallskip}
$T_{\rm sum}/T_{\rm pol}$&8.67&23.25&55.86\\
\noalign{\smallskip}\hline\noalign{\smallskip}
\end{tabular}
\end{table}

Considering the results from Table 6.1, we come to conclusion that the method 
based on multiple integral sums even when $T-t=2^{-7}$ is more than 50 times 
worse in terms of computer time for 
modeling
the collection 
of iterated It\^{o} stochastic integrals
$I_{00_{T,t}}^{(12)},$ $I_{0_{T,t}}^{(1)}$, than the method based on 
generalized multiple Fourier series.

It is not difficult to see that this effect will be
more essential if we consider iterated stochastic integrals
of multiplicities $3, 4, \ldots$ or choose value $T-t$
smaller than $2^{-7}.$

\section{Iterated It\^{o} Stochastic Integrals as Solutions of Systems 
of Linear It\^{o} SDEs}

Milstein G.N. \cite{Zapad-1} (also see \cite{Zapad-12xxx}))
proposed an approach to numerical 
modeling
of iterated It\^{o}
stochastic integrals based on their representation in the 
form of systems of linear It\^{o} SDEs.  
Let us consider this approach using the following set of 
iterated It\^{o} stochastic integrals
\begin{equation}
\label{13.001}
I_{(0)s,t}^{(i_1)}=\int\limits_{t}^s d{\bf w}_{t_1}^{(i_1)},\ \ \
I_{(00)s,t}^{(i_1i_2)}=\int\limits_t^s
\int\limits_t^{t_2}d{\bf w}_{t_1}^{(i_1)}d{\bf w}_{t_2}^{(i_2)},
\end{equation}
where $i_1, i_2=1,\ldots,m,$ $0\le t<s\le T,$
${\bf w}_{s}^{(i)}$ ($i=1,\ldots,m$) are independent
standard Wiener processes.

Obviously, we have the following representation
\begin{equation}
\label{13.002}
~~~~~~~~  d\left(
\begin{matrix}
I_{(0)s,t}^{(i_1)}\cr\cr
I_{(00)s,t}^{(i_1i_2)}
\end{matrix}\right)=
\left(
\begin{matrix}
0\ \ \ 0\cr 1\ \ \ 0
\end{matrix}\right)
\left( 
\begin{matrix}
I_{(0)s,t}^{(i_1)}\cr\cr
\vspace{1mm}
I_{(00)s,t}^{(i_1i_2)}
\end{matrix}\right)d{\bf w}_s^{(i_2)}+
\left(
\begin{matrix}
1\ \ \ 0\cr 0\ \ \ 0
\end{matrix}\right)
d\left(
\begin{matrix}
{\bf w}_s^{(i_1)}\cr \cr {\bf w}_s^{(i_2)}
\end{matrix}\right).
\end{equation}

It is well known \cite{Zapad-1}, \cite{Zapad-3} that the solution of 
system (\ref{13.002}) has the following integral form
\begin{equation}
\label{13.003}
\left(
\begin{matrix}I_{(0)s,t}^{(i_1)}\cr\cr
I_{(00)s,t}^{(i_1i_2)}
\end{matrix}\right)=\int\limits_t^s
e^{{\small\left(
\begin{matrix}
0\ \ \ 0\cr 1\ \ \ 0
\end{matrix}\right)}\left({\bf w}_s^{(i_2)}-
{\bf w}_{\theta}^{(i_2)}\right)}
\left(
\begin{matrix}1\ \ \ 0\cr 0\ \ \ 0
\end{matrix}\right)
d\left(
\begin{matrix}
{\bf w}_{\theta}^{(i_1)}\cr \cr {\bf w}_{\theta}^{(i_2)}
\end{matrix}\right),
\end{equation}
where $e^A$ is a matrix exponent
$$
e^A\stackrel{\rm def}{=}\sum\limits_{k=0}^{\infty}\frac{A^k}{k!},
$$
$A$ is a square matrix, and $A^0\stackrel{\rm def}{=}I$ 
is a unity
matrix.

Numerical modeling of the right-hand side of (\ref{13.003}) is 
unlikely
simpler task than the jointly 
numerical modeling of the 
collection of stochastic integrals (\ref{13.001}). We have to 
perform numerical modeling of (\ref{13.001}) 
within the frames
of the considered approach by numerical integration of the system of linear 
It\^{o} SDEs (\ref{13.002}). This procedure can be 
realized using the Euler (Euler--Ma\-ru\-yama) method \cite{Zapad-1}. Note that the 
expressions of more 
accurate 
numerical methods for the 
system (\ref{13.002}) (see Chapter 4) 
contain the iterated It\^{o} stochastic integrals (\ref{13.001}) and 
therefore they useless in our situation.

Let $\{\tau_j\}_{j=0}^N$ be the partition of 
$[t,s]$ such that
$$
\tau_j=t+j\Delta,\ \ \ j=0, 1,\ldots, N,\ \ \
\tau_N=s.
$$

Let us consider the Euler method
for the system of 
linear It\^{o} SDEs (\ref{13.002})
\begin{equation}
\label{13.004}
~~~~~\left(
\begin{matrix}{\bf y}_{p+1}^{(i_1)}\cr\cr 
{\bf y}_{p+1}^{(i_1i_2)}
\end{matrix}\right)=
\left(
\begin{matrix}
{\bf y}_{p}^{(i_1)}\cr\cr {\bf y}_{p}^{(i_1i_2)}
\end{matrix}\right)+
\left(
\begin{matrix}
\Delta{\bf w}_{\tau_p}^{(i_1)}\cr\cr
{\bf y}_p^{(i_1)}\Delta{\bf w}_{\tau_p}^{(i_2)}
\end{matrix}\right),\ \ \
{\bf y}_0^{(i_1)}=0,\ \ \ {\bf y}_0^{(i_1i_2)}=0,
\end{equation}
where 
$$
{\bf y}_{\tau_p}^{(i_1)}\stackrel{\rm def}{=}{\bf y}_p^{(i_1)},\ \ \
{\bf y}_{\tau_p}^{(i_1i_2)}\stackrel{\rm def}{=}{\bf y}_p^{(i_1i_2)}
$$ 
are approximations of the iterated It\^{o}
stochastic integrals
$I_{(0)\tau_p,t}^{(i_1)},$ $I_{(00)\tau_p,t}^{(i_1i_2)}$
obtained using the numerical scheme
(\ref{13.004}), $\Delta{\bf w}_{\tau_p}^{(i)}={\bf w}_{\tau_{p+1}}^{(i)}-
{\bf w}_{\tau_p}^{(i)},$ $i=1,\ldots,m.$

Iterating the expression 
(\ref{13.004}), we have
\begin{equation}
\label{13.005}
{\bf y}_N^{(i_1)}=\sum_{l=0}^{N-1}\Delta{\bf w}_{\tau_l}^{(i_1)},\ \ \
{\bf y}_N^{(i_1i_2)}=\sum_{q=0}^{N-1}
\sum_{l=0}^{q-1}\Delta{\bf w}_{\tau_l}^{(i_1)}\Delta{\bf w}_{\tau_q}^{(i_2)},
\end{equation}
where $\sum\limits_{\emptyset}\stackrel{\rm def}{=}0.$

Obviously, the formulas (\ref{13.005}) are formulas 
for approximations of the iterated It\^{o} stochastic 
integrals (\ref{13.001})
obtained using the method based on multiple integral sums
(see (\ref{bob30}), (\ref{bob31})).
   
Consequently, the efficiency of methods for the approximation 
of iterated It\^{o} stochastic integrals based on multiple integral sums 
and numerical integration of systems of linear It\^{o} SDEs 
on the base of the Euler method turns out to be equivalent.

\section{Combined Method of the Mean-Square
Approximation of Iterated It\^{o} Stochastic Integrals}

This section is written of the base of the work \cite{comb} 
(also see \cite{12aa})
and devoted to the combined method of approximation of
iterated It\^{o} stochastic integrals based on Theorem 1.1 and
the method of multiple integral sums (see Sect.~6.3).
 
The combined method of approximation of iterated It\^{o} stochastic 
integrals provides a possibility to minimize significantly the total 
number of the Fourier--Legendre
coefficients which are necessary 
for the approximation of iterated It\^{o} stochastic integrals. 
However, in this connection the 
computational costs for approximation of the mentioned stochastic integrals 
are become bigger.

Using the additive property of the It\^{o} stochastic integral, we have
\begin{equation}
\label{4.my}
I_{(0)T,t}^{(i_1)}=\sqrt{\Delta}\sum_{k=0}^{N-1}
\zeta_{0,k}^{(i_1)}\ \ \ \hbox{w.~p.~1,}
\end{equation}
\begin{equation}
\label{5.my}
I_{(1)T,t}^{(i_1)}=\sum_{k=0}^{N-1}
\left(I_{(1)\tau_{k+1},\tau_k}^{(i_1)}-
\Delta^{3/2}k\zeta_{0,k}^{(i_1)}\right)\ \ \ \hbox{w.~p.~1,}
\end{equation}
\begin{equation}
\label{6.my}
I_{(00)T,t}^{(i_1 i_2)}=\Delta\sum_{k=0}^{N-1}\sum_{l=0}^{k-1}
\zeta_{0,l}^{(i_1)}\zeta_{0,k}^{(i_2)}+
\sum_{k=0}^{N-1}I_{(00)\tau_{k+1},\tau_k}^{(i_1 i_2)}\ \ \ \hbox{w.~p.~1,}
\end{equation}

\vspace{1mm}
$$
I_{(000)T,t}^{(i_1 i_2 i_3)}=\Delta^{3/2}\sum_{k=0}^{N-1}\sum_{l=0}^{k-1}
\sum_{q=0}^{l-1}
\zeta_{0,q}^{(i_1)}\zeta_{0,l}^{(i_2)}\zeta_{0,k}^{(i_3)}+
$$
\begin{equation}
\label{7.my}
+\sqrt{\Delta}\sum_{k=0}^{N-1}\sum_{l=0}^{k-1}\left(
I_{(00)\tau_{l+1},\tau_l}^{(i_1 i_2)}\zeta_{0,k}^{(i_3)}
+\zeta_{0,l}^{(i_1)}I_{(00)\tau_{k+1},\tau_k}^{(i_2 i_3)}\right)
+
\sum_{k=0}^{N-1}I_{(000)\tau_{k+1},\tau_k}^{(i_1 i_2 i_3)}\ \ \
\hbox{w.~p.~1,}
\end{equation}

\vspace{2mm}
\noindent
where stochastic integrals 
$$
I_{(0)T,t}^{(i_1)},\ \ \
I_{(1)T,t}^{(i_1)},\ \ \ I_{(00)T,t}^{(i_1 i_2)},\ \ \
I_{(000)T,t}^{(i_1 i_2 i_3)}
$$
have the form (\ref{k1000xxxx}),
$i_1,i_2,i_3=1,\ldots,m,$
$T-t=N\Delta,$ $\tau_k=t+k\Delta,$ 
$$
\zeta_{0,k}^{(i)}\stackrel{\rm def}{=}
\frac{1}{\sqrt{\Delta}}\int\limits_{\tau_k}^{\tau_{k+1}}d{\bf w}_s^{(i)},
$$
$k=0, 1,\ldots,N-1,$ 
the sum with respect to the empty set is equal to zero.

Substituting the relation
$$
I_{(1)\tau_{k+1},\tau_k}^{(i_1)}=-\frac{\Delta^{3/2}}{2}\left(
\zeta_{0,k}^{(i_1)}+\frac{1}{\sqrt{3}}\zeta_{1,k}^{(i_1)}\right)\ \ \
\hbox{w.~p.~1}
$$
into (\ref{5.my}), where
$\zeta_{0,k}^{(i_1)},$ $\zeta_{1,k}^{(i_1)}$ are
independent standard Gaussian random variables,
we get
\begin{equation}
\label{8.my}
~~~~~~~ I_{(1)T,t}^{(i_1)}=-\Delta^{3/2}\sum_{k=0}^{N-1}
\left(\left(\frac{1}{2}+k\right)\zeta_{0,k}^{(i_1)}+
\frac{1}{2\sqrt{3}}\zeta_{1,k}^{(i_1)}\right)\ \ \ \hbox{w.~p.~1.}
\end{equation}

Consider approximations of the following 
iterated It\^{o} stochastic integrals using the 
method based on multiple Fourier--Legendre series (Theorem 1.1)
$$
I_{(00)\tau_{k+1},\tau_k}^{(i_1 i_2)},\ \ \
I_{(00)\tau_{k+1},\tau_k}^{(i_2 i_3)},\ \ \
I_{(000)\tau_{k+1},\tau_k}^{(i_1 i_2 i_3)}.
$$

As a result, we get 
\begin{equation}
\label{9.my}
~~I_{(00)T,t}^{(i_1 i_2)N,q}=\Delta\sum_{k=0}^{N-1}\sum_{l=0}^{k-1}
\zeta_{0,l}^{(i_1)}\zeta_{0,k}^{(i_2)}+
\sum_{k=0}^{N-1}I_{(00)\tau_{k+1},\tau_k}^{(i_1 i_2)q},
\end{equation}

$$
I_{(000)T,t}^{(i_1 i_2 i_3)N,q_1,q_2}=
\Delta^{3/2}\sum_{k=0}^{N-1}\sum_{l=0}^{k-1}
\sum_{q=0}^{l-1}
\zeta_{0,q}^{(i_1)}\zeta_{0,l}^{(i_2)}\zeta_{0,k}^{(i_3)}+
$$
\begin{equation}
\label{10.my}
~~~~~ +\sqrt{\Delta}\sum_{k=0}^{N-1}\sum_{l=0}^{k-1}\left(
I_{(00)\tau_{l+1},\tau_l}^{(i_1 i_2)q_1}\zeta_{0,k}^{(i_3)}
+\zeta_{0,l}^{(i_1)}I_{(00)\tau_{k+1},\tau_k}^{(i_2 i_3)q_1}\right)
+\sum_{k=0}^{N-1}I_{(000)\tau_{k+1},\tau_k}^{(i_1 i_2 i_3)q_2},
\end{equation}

\vspace{2mm}
\noindent
where we suppose that the approximations
$$
I_{(00)\tau_{k+1},\tau_k}^{(i_1 i_2)q},\ \ \
I_{(00)\tau_{k+1},\tau_k}^{(i_1 i_2)q_1},\ \ \
I_{(000)\tau_{k+1},\tau_k}^{(i_1 i_2 i_3)q_2}
$$ 
are obtained using Theorem 1.1 (the case of Legendre polynomials).

In particular, when $N=2$, the formulas 
(\ref{4.my}), (\ref{8.my})-(\ref{10.my})
will look as follows

\vspace{-6mm}
\begin{equation}
\label{roar1}
I_{(0)T,t}^{(i_1)}=\sqrt{\Delta}\left(
\zeta_{0,0}^{(i_1)}+\zeta_{0,1}^{(i_1)}\right)\ \ \ \hbox{w.~p.~1,}
\end{equation}

\vspace{-6mm}
\begin{equation}
\label{roar2}
~~~~~~~I_{(1)T,t}^{(i_1)}=-\Delta^{3/2}\left(\frac{1}{2}\zeta_{0,0}^{(i_1)}+
\frac{3}{2}\zeta_{0,1}^{(i_1)}+\frac{1}{2\sqrt{3}}\left(
\zeta_{1,0}^{(i_1)}+\zeta_{1,1}^{(i_1)}\right)\right)\ \ \ \hbox{w.~p.~1,}
\end{equation}

\begin{equation}
\label{roar3}
I_{(00)T,t}^{(i_1 i_2)2,q}=\Delta\left(
\zeta_{0,0}^{(i_1)}\zeta_{0,1}^{(i_2)}+
I_{(00)\tau_{1},\tau_0}^{(i_1 i_2)q}+
I_{(00)\tau_{2},\tau_1}^{(i_1 i_2)q}\right),
\end{equation}

$$
I_{(000)T,t}^{(i_1 i_2 i_3)2,q_1,q_2}=
\sqrt{\Delta}\left(
I_{(00)\tau_{1},\tau_0}^{(i_1 i_2)q_1}\zeta_{0,1}^{(i_3)}
+\zeta_{0,0}^{(i_1)}I_{(00)\tau_{2},\tau_1}^{(i_2 i_3)q_1}\right)
+
$$

\vspace{-2mm}
\begin{equation}
\label{roar4}
+I_{(000)\tau_{1},\tau_0}^{(i_1 i_2 i_3)q_2}+
I_{(000)\tau_{2},\tau_1}^{(i_1 i_2 i_3)q_2},
\end{equation}

\vspace{2mm}
\noindent
where $\Delta=(T-t)/2,$ $\tau_k=t+k\Delta,$ $k=0, 1, 2.$

Note that if  $N=1$, then (\ref{4.my}), (\ref{8.my})-(\ref{10.my}) 
are the formulas for numerical modeling of the mentioned stochastic 
integrals using the method based on Theorem 1.1.

Further, we will demonstrate that 
modeling 
of the iterated  It\^{o} stochastic integrals
$$
I_{(0)T,t}^{(i_1)},\ \ \
I_{(1)T,t}^{(i_1)},\ \ \
I_{(00)T,t}^{(i_1 i_2)},\ \ \
I_{(000)T,t}^{(i_1 i_2 i_3)}
$$

\noindent
using the formulas (\ref{roar1})--(\ref{roar4}) results in 
abrupt decrease of the total number of Fourier--Legendre 
coefficients, which are necessary for approximation of these 
stochastic integrals using the method based on Theorem 1.1.

From the other hand, the formulas 
(\ref{roar3}), (\ref{roar4})
include two approximations of iterated It\^{o} stochastic integrals of
second 
and third multiplicity, and each one of them should be obtained 
using the method based on Theorem 1.1. Obviously, this leads to an
increase in computational costs for the approximation.

Let us calculate the mean-square approximation errors
for the formulas (\ref{9.my}),
(\ref{10.my}). We have

\vspace{-3mm}
$$
E^{q}_N\stackrel{\rm def}{=}
{\sf M}\left\{\biggl(I_{(00)T,t}^{(i_1 i_2)}-
I_{(00)T,t}^{(i_1 i_2)N,q}\biggr)^2\right\}=
\sum_{k=0}^{N-1}{\sf M}\left\{\biggl(
I_{(00)\tau_{k+1},\tau_k}^{(i_1 i_2)}-
I_{(00)\tau_{k+1},\tau_k}^{(i_1 i_2)q}\biggr)^2\right\}=
$$

\begin{equation}
\label{15.my}
~~~~=N\frac{\Delta^2}{2}\Biggl(\frac{1}{2}-\sum_{l=1}^q
\frac{1}{4l^2-1}\Biggr)
=\frac{(T-t)^2}{2N}
\Biggl(\frac{1}{2}-\sum_{l=1}^q
\frac{1}{4l^2-1}\Biggr),
\end{equation}

\vspace{6mm}

$$
E^{q_1,q_2}_N\stackrel{\rm def}{=}
{\sf M}\left\{\biggl(I_{(000)T,t}^{(i_1 i_2 i_3)}-
I_{(000)T,t}^{(i_1 i_2 i_3)N,q_1,q_2}\biggr)^2\right\}=
$$

\vspace{1mm}
$$
={\sf M}\Biggl\{\Biggl(\sum_{k=0}^{N-1}\Biggl(
\sqrt{\Delta}\sum_{l=0}^{k-1}\biggl(\zeta_{0,k}^{(i_3)}\biggl(
I_{(00)\tau_{l+1},\tau_l}^{(i_1 i_2)}-
I_{(00)\tau_{l+1},\tau_l}^{(i_1 i_2)q_1}\biggr)\biggr.\Biggr.\Biggr.\Biggr.+
$$

$$
\Biggl.\Biggl.\Biggl.\biggl.+
\zeta_{0,l}^{(i_1)}\biggl(
I_{(00)\tau_{k+1},\tau_k}^{(i_2 i_3)}-
I_{(00)\tau_{k+1},\tau_k}^{(i_2 i_3)q_1}\biggr)\biggr)+
I_{(000)\tau_{k+1},\tau_k}^{(i_1 i_2 i_3)}-
I_{(000)\tau_{k+1},\tau_k}^{(i_1 i_2 i_3)q_2}\Biggr)\Biggr)^2\Biggr\}=
$$

\vspace{2mm}

$$
=\sum_{k=0}^{N-1}{\sf M}\Biggl\{\Biggl(
\sqrt{\Delta}\sum_{l=0}^{k-1}\biggl(\zeta_{0,k}^{(i_3)}\biggl(
I_{(00)\tau_{l+1},\tau_l}^{(i_1 i_2)}-
I_{(00)\tau_{l+1},\tau_l}^{(i_1 i_2)q_1}\biggr)\biggr.\Biggr.\Biggr.+
$$

$$
\Biggl.\Biggl.\biggl.+
\zeta_{0,l}^{(i_1)}\biggl(
I_{(00)\tau_{k+1},\tau_k}^{(i_2 i_3)}-
I_{(00)\tau_{k+1},\tau_k}^{(i_2 i_3)q_1}\biggr)\biggr)+
I_{(000)\tau_{k+1},\tau_k}^{(i_1 i_2 i_3)}-
I_{(000)\tau_{k+1},\tau_k}^{(i_1 i_2 i_3)q_2}\Biggr)^2\Biggr\}=
$$

\newpage
\noindent
$$
=\sum_{k=0}^{N-1}\left(\Delta{\sf M}\left\{\left(
\zeta_{0,k}^{(i_3)}\sum_{l=0}^{k-1}\biggl(
I_{(00)\tau_{l+1},\tau_l}^{(i_1 i_2)}-
I_{(00)\tau_{l+1},\tau_l}^{(i_1 i_2)q_1}\biggr)\right)^2\right\}\right.+
$$
$$
\left.+\Delta{\sf M}\left\{\left(\biggl(
I_{(00)\tau_{k+1},\tau_k}^{(i_2 i_3)}-
I_{(00)\tau_{k+1},\tau_k}^{(i_2 i_3)q_1}\biggr)\sum_{l=0}^{k-1}
\zeta_{0,l}^{(i_1)}\right)^2\right\}+H_{k,q_2}^{(i_1i_2i_3)}\right)=
$$

\vspace{2mm}

$$
=\sum_{k=0}^{N-1}\left(\Delta\sum_{l=0}^{k-1}{\sf M}\Biggl\{\biggl(
I_{(00)\tau_{l+1},\tau_l}^{(i_1 i_2)}-
I_{(00)\tau_{l+1},\tau_l}^{(i_1 i_2)q_1}\biggr)^2\Biggr\}+\right.
$$
$$
+\left. k\Delta{\sf M}\Biggl\{\biggl(
I_{(00)\tau_{k+1},\tau_k}^{(i_2 i_3)}-
I_{(00)\tau_{k+1},\tau_k}^{(i_2 i_3)q_1}\biggr)^2\Biggr\}+
H_{k,q_2}^{(i_1i_2i_3)}\right)=
$$

\vspace{2mm}

$$
=\sum_{k=0}^{N-1}\Biggl(2k\Delta{\sf M}\Biggl\{\biggl(
I_{(00)\tau_{k+1},\tau_k}^{(i_1 i_2)}-
I_{(00)\tau_{k+1},\tau_k}^{(i_1 i_2)q_1}\biggr)^2\Biggr\}+
H_{k,q_2}^{(i_1i_2i_3)}\Biggr)=
$$

$$
=\sum_{k=0}^{N-1}\Biggl(2k\Delta
\frac{\Delta^2}{2}\Biggl(\frac{1}{2}-\sum_{l=1}^{q_1}
\frac{1}{4l^2-1}\Biggr)+H_{k,q_2}^{(i_1i_2i_3)}\Biggr)=
$$

$$
=\Delta^3\frac{N(N-1)}{2}
\Biggl(\frac{1}{2}-\sum_{l=1}^{q_1}
\frac{1}{4l^2-1}\Biggr)+\sum_{k=0}^{N-1}H_{k,q_2}^{(i_1i_2i_3)}=
$$

$$
=\frac{1}{2}(T-t)^3\Biggl(\frac{1}{N}-\frac{1}{N^2}\Biggr)
\Biggl(\frac{1}{2}-\sum_{l=1}^{q_1}
\frac{1}{4l^2-1}\Biggr)+
$$

\begin{equation}
\label{16.my}
+\sum_{k=0}^{N-1}H_{k,q_2}^{(i_1i_2i_3)},
\end{equation}

\noindent
where
$$
H_{k,q_2}^{(i_1i_2i_3)}=
{\sf M}\Biggl\{\biggl(
I_{(000)\tau_{k+1},\tau_k}^{(i_1 i_2 i_3)}-
I_{(000)\tau_{k+1},\tau_k}^{(i_1 i_2 i_3)q_2}\biggr)^2\Biggr\}.
$$

\noindent
\par
Moreover, we suppose that $i_1\ne i_2$ in (\ref{15.my}) 
and not all indices
$i_1, i_2, i_3$
in (\ref{16.my}) are 
equal. Otherwise there are  
simple relationships for modeling 
the integrals
$I_{(00)T,t}^{(i_1i_2)},$ $I_{(000)T,t}^{(i_1i_2i_3)}$
$$
I_{(00)T,t}^{(i_1 i_1)}
=\frac{1}{2}(T-t)\biggl(
\left(\zeta_0^{(i_1)}\right)^2-1\biggr)\ \ \ \hbox{w.~p.~1},
$$
$$
I_{(000)T,t}^{(i_1 i_1 i_1)}
=\frac{1}{6}(T-t)^{3/2}\biggl(
\left(\zeta_0^{(i_1)}\right)^3-3
\zeta_0^{(i_1)}\biggr)\ \ \ \hbox{w.~p.~1},
$$
where
$$
\zeta_0^{(i_1)}=
\frac{1}{\sqrt{T-t}}\int\limits_t^T d{\bf w}_s^{(i_1)}
$$
is a standard Gaussian random variable.

For definiteness, assume that 
$i_1, i_2, i_3$ are  pairwise
different in (\ref{16.my})  (other cases are represented by 
(\ref{39a})--(\ref{leto1041})). 
Then from Theorem 1.3 we have
\begin{equation}
\label{17.my}
H_{k,q_2}^{(i_1i_2i_3)}=\Delta^3\Biggl(\frac{1}{6}-
\sum_{j_1,j_2,j_3=0}^{q_2}\frac{C_{j_3j_2j_1}^2}{\Delta^3}\Biggr),
\end{equation}

\noindent
where
$$
C_{j_3j_2j_1}=\frac{\sqrt{(2j_1+1)(2j_2+1)(2j_3+1)}}{8}
\Delta^{3/2}\bar C_{j_3j_2j_1},
$$
$$
\bar C_{j_3j_2j_1}=\int\limits_{-1}^1 P_{j_3}(z)
\int\limits_{-1}^z P_{j_2}(y)
\int\limits_{-1}^y P_{j_1}(x)dx dy dz,
$$

\noindent
and $P_i(x)$ ($i=0, 1, 2,\ldots $) is the Legendre polynomial.

Substituting (\ref{17.my}) into (\ref{16.my}), we obtain
$$
E_{N}^{q_1,q_2}
=\frac{1}{2}(T-t)^3\left(\frac{1}{N}-\frac{1}{N^2}\right)
\Biggl(\frac{1}{2}-\sum_{l=1}^{q_1}
\frac{1}{4l^2-1}\Biggr)+
$$
\begin{equation}
\label{18.my}
~~~~~~~ +\frac{(T-t)^3}{N^2}\Biggl(\frac{1}{6}-
\sum_{j_1,j_2,j_3=0}^{q_2}\frac{(2j_1+1)(2j_2+1)(2j_3+1)}{64}
\bar C_{j_3j_2j_1}^2\Biggr).
\end{equation}

Note that for $N=1$ the formulas (\ref{15.my}), (\ref{18.my}) pass into the 
corresponding formulas for the mean-square 
approximation errors of the iterated It\^{o} stochastic integrals 
$I_{(00)T,t}^{(i_1i_2)},$ $I_{(000)T,t}^{(i_1i_2i_3)}$ 
(see Theorem 1.3).

Let us consider modeling
the integrals $I_{(0)T,t}^{(i_1)},$ $I_{(00)T,t}^{(i_1i_2)}$.
To do it we can use the relations 
(\ref{4.my}), (\ref{9.my}).  At that, the 
mean-square approximation error for the 
integral $I_{(00)T,t}^{(i_1i_2)}$ is defined by the formula 
(\ref{15.my}) for the case of 
Legendre polynomials. Let us calculate the value 
$E_{N}^q$ for various
$N$ and $q$ 
\begin{equation}
\label{19.my}
E_{3}^2\approx 0.0167(T-t)^2,\ \ \
E_{2}^3\approx 0.0179(T-t)^2,
\end{equation}
\begin{equation}
\label{20.my}
E_{1}^{6}\approx 0.0192(T-t)^2.
\end{equation}

Note that the combined method (see (\ref{19.my})) requires 
calculation of a significantly smaller number of the Fourier--Legendre
coefficients than the method based on Theorem 1.1
(see (\ref{20.my})).

\begin{table}
\centering
\caption{$T-t=0.1.$}
\label{tab:6.2}      
\begin{tabular}{p{1.7cm}p{1.7cm}p{1.7cm}p{1.7cm}p{1.7cm}}
\hline\noalign{\smallskip}
$N$&$q$&$q_1$&$q_2$&$M$\\
\noalign{\smallskip}\hline\noalign{\smallskip}
$1$&13&--&1&21\\
$2$&6&0&0&7\\
$3$&4&0&0&5\\
\noalign{\smallskip}\hline\noalign{\smallskip}
\end{tabular}
\end{table}

\begin{table}
\centering
\caption{$T-t=0.05.$}
\label{tab:6.3}      
\begin{tabular}{p{1.7cm}p{1.7cm}p{1.7cm}p{1.7cm}p{1.7cm}}
\hline\noalign{\smallskip}
$N$&$q$&$q_1$&$q_2$&$M$\\
\noalign{\smallskip}\hline\noalign{\smallskip}
$1$&50&--&2&77\\
$2$&25&2&0&26\\
$3$&17&1&0&18\\
\noalign{\smallskip}\hline\noalign{\smallskip}
\end{tabular}
\end{table}

\begin{table}
\centering
\caption{$T-t=0.02.$}
\label{tab:6.4}      
\begin{tabular}{p{1.7cm}p{1.7cm}p{1.7cm}p{1.7cm}p{1.7cm}}
\hline\noalign{\smallskip}
$N$&$q$&$q_1$&$q_2$&$M$\\
\noalign{\smallskip}\hline\noalign{\smallskip}
$1$&312&--&6&655\\
$2$&156&4&2&183\\
$3$&104&6&0&105\\
\noalign{\smallskip}\hline\noalign{\smallskip}
\end{tabular}
\end{table}

Assume that the mean-square approximation error
of the iterated It\^{o} stochastic integrals
$I_{(00)T,t}^{(i_1i_2)},$ $I_{(000)T,t}^{(i_1i_2i_3)}$ 
equals to  $(T-t)^4.$

In Tables  6.2--6.4 we can see the values 
$N, q, q_1, q_2,$ which 
satisfy the system of inequalities
\begin{equation}
\label{21.my}
\left\{
\begin{matrix}
E_{N}^{q}\le(T-t)^4\cr\cr
E_{N}^{q_1,q_2}\le(T-t)^4
\end{matrix}\right.
\end{equation}

\noindent
as well as the total number $M$ of the Fourier--Legendre coefficients,
which are nesessary for approximation of the
iterated It\^{o} stochastic integrals  
$I_{(00)T,t}^{(i_1i_2)},$ $I_{(000)T,t}^{(i_1i_2i_3)}$
when $T-t=0.1, 0.05, 0.02$ (the numbers $q, q_1, q_2$
were taken in such a manner that the number $M$ was
the smallest one).

From Tables  6.2--6.4 it is clear that the combined method 
with the small $N$ $(N=2)$ provides a possibility to decrease significantly 
the total number of the Fourier--Legendre coefficients, 
which are necessary for the approximation of the 
iterated It\^{o} stochastic integrals
$I_{(00)T,t}^{(i_1i_2)},$ $I_{(000)T,t}^{(i_1i_2i_3)}$  
in comparison with the method based on Theorem 1.1
$(N=1).$ However, as we noted before, as a result the computational
costs for the approximation are increased. The approximation accuracy 
of iterated It\^{o}
stochastic integrals for the combined method and the method 
based on Theorem 1.1 was taken $(T-t)^4$.

\section{Representation of Iterated It\^{o} Stochastic Integrals 
of Multiplicity $k$
with Respect to the Scalar Standard Wiener Process Based on Hermite Polynomials}

In Chapters 1, 2, and 5 we analyzed the general 
theory of the approximation of iterated It\^{o} and Stratonovich 
stochastic integrals
with respect to components of the multidimensional Wiener process.
However, in some narrow special cases we can get exact 
expressions for
iterated It\^{o} and Stratonovich stochastic integrals in the 
form of polynomials of finite degrees from one standard Gaussian 
random variable. This and next sections will be devoted to this 
question. The results described in them can be found, 
for example, in \cite{Ch} (also see \cite{Zapad-1}, \cite{Zapad-3}).

Let us consider the set of 
polynomials 
$H_n(x,y),$ $n=0, 1,\ldots$ defined by
$$
\Biggl.H_n(x,y)=\left(\frac{d^n}{d\alpha^n} 
e^{\alpha x-\alpha^2 y/2}\right)
\Biggr|_{\alpha=0}.
$$

It is well known that polynomials $H_n(x,y)$ are connected with 
the Hermite polynomials
$$
h_n(x)=(-1)^n e^{x^2} \frac{d^n}{dx^n}\left(e^{-x^2}\right)
$$
by the formula 
$$
H_n(x,y)=\left(\frac{y}{2}\right)^{n/2}
h_n\left(\frac{x}{\sqrt{2y}}\right)=y^{n/2}H_n\left(\frac{x}{\sqrt{y}}\right),
$$
where $H_n(x)$ is the Hermite polynomial (\ref{ziko500}).

Using the recurrent formulas
$$
\frac{dh_n}{dz}(z)=2nh_{n-1}(z),\ \ \ n=1,\ 2,\ldots,
$$
$$
h_n(z)=2zh_{n-1}(z)-2(n-1)h_{n-2}(z),\ \ \ n=2,\ 3,\ldots,
$$
it is easy to get the following recurrent relations 
for polynomials $H_n(x,y)$
\begin{equation}
\label{4.2.11}
\frac{\partial H_n}{\partial x}(x,y)=nH_{n-1}(x,y),\ \ \ n=1, 2,\ldots,
\end{equation}
\begin{equation}
\label{4.2.12}
~~~~~~~~~~~~\frac{\partial H_n}{\partial y}(x,y)=\frac{n}{2y}
H_{n}(x,y)-\frac{nx}{2y}H_{n-1}(x,y),\ \ \ n=1, 2,\ldots,
\end{equation}

\begin{equation}
\label{4.2.13}
~~~~~\frac{\partial H_n}{\partial y}(x,y)
=-\frac{n(n-1)}{2}H_{n-2}(x,y),\ \ \ n=2, 3,\ldots
\end{equation}

From (\ref{4.2.11}) -- (\ref{4.2.13}) it follows that
\begin{equation}
\label{4.2.14}
\frac{\partial H_n}{\partial y}(x,y)
+\frac{1}{2}\frac{\partial^2H_n}{\partial x^2}(x,y)=0,\ \ \ n=2, 3,\ldots
\end{equation}

Using the It\^{o} formula, we have
\begin{equation}
\label{4.2.15}
H_n(w_t,t)-H_n(0,0)=\int\limits_0^t
\frac{\partial H_n}{\partial x }(w_s,s)dw_s
+\hspace{-0.7mm}
\int\limits_0^t \hspace{-1.3mm}\left(\hspace{-0.4mm}\frac{\partial H_n}{\partial y}(w_s,s)+
\frac{1}{2}\frac{\partial^2H_n}{\partial x^2}(w_s,s)\hspace{-0.7mm}\right)\hspace{-0.5mm}ds
\end{equation}
w.~p.~1, where $t\in [0, T]$ and $w_t$ is a scalar standard Wiener
process.

Note that $H_n(0,0)=0,$ $n=2, 3,\ldots$ Then
from (\ref{4.2.14})
and (\ref{4.2.15}) we get
\begin{equation}
\label{ziko9001}
~~~~~~ H_n(w_t,t)=\int\limits_0^t nH_{n-1}(w_s,s)dw_s,\ \ \ \hbox{w.~p.~1}\ \ \ 
(n=2, 3,\ldots).
\end{equation}

Furthermore, by
induction it is easy to get the following relation (see (\ref{ziko9001}))
\begin{equation}
\label{4.2.17}
~~~~~~\int\limits_0^t \ldots\int\limits_0^{t_{2}}dw_{t_1}\ldots
dw_{t_n}=\frac{H_n(w_t,t)}{n!}\ \ \ \hbox{w.~p.~1}\ \ \ 
(n=1, 2,\ldots).
\end{equation}

Let us consider one generalization of the formula (\ref{4.2.17}) \cite{Ch}
\begin{equation}
\label{4.2.18}
~~~~~ J_t^{(n)}\stackrel{\sf def}{=}
\int\limits_0^t \psi(t_n)\ldots
\int\limits_0^{t_{2}}\psi(t_1)dw_{t_1}\ldots
dw_{t_n}
=\frac{H_n(x_t,y(t))}{n!}\ \ \ \hbox{w.~p.~1},
\end{equation}
where $t\in [0, T]$, $n=1, 2,\ldots$\ \hspace{-1mm}, and
$$
x_t\stackrel{\rm def}{=}\int\limits_0^t\psi(s)dw_{s},\ \ \
y(t)\stackrel{\rm def}{=}\int\limits_0^t \psi^2(s)ds,
$$
where $\psi(s)\in L_2([0,T])$.

To prove the equality (\ref{4.2.18}), we apply the It\^{o} formula.
Using the It\^{o} formula and (\ref{4.2.11}), (\ref{4.2.14}), we obtain w.~p.~1
($H_n(0,0)=0,$ $n=2, 3,\ldots$)
$$
H_n(x_t,y(t))-H_n(0,0)=\int\limits_0^t
\psi(s)\frac{\partial H_n}{\partial x}(x_s,y(s))dw_s+
$$
$$
+\int\limits_0^t\left(
\frac{\partial H_n}{\partial s}(x_s,y(s))+\frac{1}{2}\psi^2(s)
\frac{\partial^2 H_n}{\partial x^2}(x_s,y(s))\right)ds=
$$
$$
=\int\limits_0^t
\psi(s)\frac{\partial H_n}{\partial x}(x_s,y(s))dw_s+
$$
$$
+\int\limits_0^t\left(
\frac{\partial H_n}{\partial y(s)}(x_s,y(s))y'(s)+\frac{1}{2}\psi^2(s)
\frac{\partial^2 H_n}{\partial x^2}(x_s,y(s))\right)ds=
$$
$$
=\int\limits_0^t
\psi(s)\frac{\partial H_n}{\partial x}(x_s,y(s))dw_s+
$$
$$
+\int\limits_0^t\psi^2(s)\left(
\frac{\partial H_n}{\partial y(s)}(x_s,y(s))+\frac{1}{2}
\frac{\partial^2 H_n}{\partial x^2}(x_s,y(s))\right)ds=
$$
$$
=\int\limits_0^t
\psi(s)\frac{\partial H_n}{\partial x}(x_s,y(s))dw_s=
\int\limits_0^t
\psi(s)nH_{n-1}(x_s,y(s))dw_s=
$$
$$
=\int\limits_0^t
\psi(s)\int\limits_0^s \psi(\tau)n(n-1)H_{n-2}(x_{\tau},y(\tau))dw_{\tau}dw_s=\ldots
$$
\begin{equation}
\label{ziko2221}
\ldots = n!\int\limits_0^t \psi(t_n)\ldots
\int\limits_0^{t_{2}}\psi(t_1)dw_{t_1}\ldots
dw_{t_n}.
\end{equation}

From (\ref{ziko2221}) we get (\ref{4.2.18}).

It is easy to check that first eight formulas from 
the set (\ref{4.2.18})
have the following form  
$$
J_t^{(1)}=\frac{1}{1!}x_t,
$$
$$
J_t^{(2)}=\frac{1}{2!}\left((x_t)^2-y(t)\right),
$$
$$
J_t^{(3)}=\frac{1}{3!}\left((x_t)^3-3x_t y(t)\right),
$$
$$
J_t^{(4)}=\frac{1}{4!}\left((x_t)^4-6(x_t)^2 y(t)
+3y^2(t)\right),
$$
$$
J_t^{(5)}=\frac{1}{5!}\left((x_t)^5-10(x_t)^3 y(t)
+15x_t y^2(t)\right),
$$
$$
J_t^{(6)}=\frac{1}{6!}\left((x_t)^6-15(x_t)^4 y(t)
+45(x_t)^2 y^2(t) -15 y^3(t)\right),
$$
$$
J_t^{(7)}=\frac{1}{7!}\left((x_t)^7-21(x_t)^5 y(t)
+105(x_t)^3 y^2(t) -105x_t y^3(t)\right),
$$
$$
J_t^{(8)}=\frac{1}{8!}\left((x_t)^8-28(x_t)^6 y(t)
+210(x_t)^4 y^2(t) -420(x_t)^2 y^3(t)+
105 y^4(t)\right)
$$

\noindent 
w.~p.~1.
As follows from the results of Sect.~1.1.6, 
for the case $\psi_1(\tau),\ldots,\psi_k(\tau)$ $\equiv \psi(\tau)$
and $i_1=\ldots=i_k=1,\ldots,m$
the formula (\ref{leto6000}) transforms into (\ref{4.2.18}).

\section{Representation of Iterated Stratonovich Stochastic Integrals of
Multiplicity $k$ with Respect to the Scalar Standard Wiener Process}

Let us prove the following relation for 
iterated Stratonovich stochastic integrals (see, for example, \cite{Zapad-3})
\begin{equation}
\label{4.2.880}
{\int\limits_0^{*}}^t\ldots
{\int\limits_0^{*}}^{t_{2}}
dw_{t_1}\ldots dw_{t_n}=\frac{\bigl(w_{t}\bigr)^n}{n!}\ \ \ 
\hbox{w.~p.~1},
\end{equation}

\noindent
where $t\in [0, T].$

At first, we will consider the case $n=2$. Using Theorem 2.12, we obtain
\begin{equation}
\label{4.2.881}
{\int\limits_0^{*}}^t
{\int\limits_0^{*}}^{t_{2}}
dw_{t_1}dw_{t_2}=\int\limits_0^t
\int\limits_0^{t_{2}}
dw_{t_1}dw_{t_2}+\frac{1}{2}
\int\limits_0^t dt_1\ \ \ \hbox{w.~p.~1}.
\end{equation}

\vspace{1mm}

From the relation (\ref{4.2.17}) for $n=2$ it follows that 
\begin{equation}
\label{4.2.882}
\int\limits_0^t
\int\limits_0^{t_{2}}
dw_{t_1}dw_{t_2}=\frac{\bigl(w_{t}\bigr)^2}{2!}
-\frac{1}{2}\int\limits_0^t dt_1\ \ \ \hbox{w.~p.~1}.
\end{equation}

\vspace{1mm}

Substituting (\ref{4.2.882}) into (\ref{4.2.881}),
we have
$$
{\int\limits_0^{*}}^t
{\int\limits_0^{*}}^{t_{2}}
dw_{t_1}dw_{t_2}=
\frac{\bigl(w_{t}\bigr)^2}{2!}\ \ \ \hbox{w.~p.~1}.
$$ 

\vspace{1mm}

So, the formula (\ref{4.2.880}) is correct 
for $n=2.$
Using the induction assumption and (\ref{d11}), we obtain
\begin{equation}
\label{ura900}
{\int\limits_0^{*}}^t\ldots
{\int\limits_0^{*}}^{t_{2}}
dw_{t_1}\ldots dw_{t_{n+1}}= 
{\int\limits_0^{*}}^t \frac{\left(w_{\tau}\right)^n}{n!}dw_{\tau}=
\int\limits_0^t \frac{\left(w_{\tau}\right)^n}{n!}dw_{\tau}+
\frac{1}{2}
\int\limits_0^t \frac{\left(w_{\tau}\right)^{n-1}}{(n-1)!}d\tau
\end{equation}
w.~p.~1. From the other hand, using the It\^{o} formula, we get
\begin{equation}
\label{ura901}
\frac{\bigl(w_{t}\bigr)^{n+1}}{(n+1)!}=
\int\limits_0^t\frac{\left(w_{\tau}\right)^{n-1}}{2(n-1)!}d\tau
+\int\limits_0^t\frac{\left(w_{\tau}\right)^{n}}{n!}dw_{\tau}\ \ \
\hbox{w.~p.~1.}
\end{equation}

From (\ref{ura900}) and (\ref{ura901}) we obtain (\ref{4.2.880}).
It is easy to see that the formula (\ref{4.2.880}) admits the 
following generalization
\begin{equation}
\label{ziko7771}
~~~~~~{\int\limits_0^{*}}^t\psi(t_{n})\ldots
{\int\limits_0^{*}}^{t_{2}}\psi(t_1)
dw_{t_1}\ldots dw_{t_n}=
\frac{1}{n!}\left(\int\limits_0^t\psi(\tau)dw_{\tau}\right)^n\ \ \ 
\hbox{w.~p.~1},
\end{equation}
where $t\in [0, T]$ and $\psi(\tau)$ is a continuous 
nonrandom function at the 
interval $[0, T]$.

To prove the equality (\ref{ziko7771}),
first consider the case $n=2$. Using Theorem 2.12, we get
\begin{equation}
\label{4.2.881ziko}
{\int\limits_0^{*}}^t\psi(t_2)
{\int\limits_0^{*}}^{t_{2}}\psi(t_1)
dw_{t_1}dw_{t_2}=\int\limits_0^t\psi(t_2)
\int\limits_0^{t_{2}}\psi(t_1)
dw_{t_1}dw_{t_2}+\frac{1}{2}
\int\limits_0^t \psi^2(s)ds\ \ \ \hbox{w.~p.~1}.
\end{equation}

From the relation (\ref{4.2.18}) for $n=2$ it follows that 
\begin{equation}
\label{4.2.882ziko}
~~~~\int\limits_0^t\psi(t_2)
\int\limits_0^{t_{2}}\psi(t_1)
dw_{t_1}dw_{t_2}=\frac{1}{2!}\left(\int\limits_0^t\psi(s)dw_s\right)^2
-\frac{1}{2}\int\limits_0^t \psi^2(s)ds\ \ \ \hbox{w.~p.~1}.
\end{equation}

Substituting (\ref{4.2.882ziko}) into (\ref{4.2.881ziko}),
we obtain
$$
{\int\limits_0^{*}}^t\psi(t_2)
{\int\limits_0^{*}}^{t_{2}}\psi(t_1)
dw_{t_1}dw_{t_2}=
\frac{1}{2!}\left(\int\limits_0^t\psi(s)dw_s\right)^2\ \ \ \hbox{w.~p.~1}.
$$

Thus the formula (\ref{ziko7771}) is proved
for $n=2.$
Applying the induction assumption and (\ref{d11}), we have 
$$
{\int\limits_0^{*}}^t\psi(t_{n+1}) \ldots
{\int\limits_0^{*}}^{t_{2}}\psi(t_1)
dw_{t_1}\ldots dw_{t_{n+1}}= 
{\int\limits_0^{*}}^t \psi(\tau)\frac{1}{n!}
\left(\int\limits_0^{\tau}\psi(s)dw_s\right)^n
dw_{\tau}=
$$
\begin{equation}
\label{ura900ziko}
=
\int\limits_0^t \psi(\tau) 
\frac{1}{n!}
\left(\int\limits_0^{\tau}\psi(s)dw_s\right)^n
dw_{\tau}
+\frac{1}{2}
\int\limits_0^t \psi^2(\tau)
\frac{1}{(n-1)!}
\left(\int\limits_0^{\tau}\psi(s)dw_s\right)^{n-1}
d\tau
\end{equation}
w.~p.~1. Applying the It\^{o} formula, we obtain
$$
\frac{1}{(n+1)!}
\left(\int\limits_0^{t}\psi(s)dw_s\right)^{n+1}
=
\int\limits_0^t\psi^2(\tau)\frac{1}{2(n-1)!}\left(\int\limits_0^{\tau}\psi(s)dw_s\right)^{n-1}d\tau
+
$$
\begin{equation}
\label{ura901ziko}
+\int\limits_0^t\psi(\tau)
\frac{1}{n!}\left(\int\limits_0^{\tau}\psi(s)dw_s\right)^{n}dw_{\tau}\ \ \
\hbox{w.~p.~1.}
\end{equation}

From (\ref{ura900ziko}) and (6.84) we get (\ref{ziko7771}).

\section{Weak Approximation of Iterated It\^{o} Stochastic Integrals
of Multiplicity 1 to 4}

In the previous chapters of the book and previous sections 
of this chapter we analyzed in detail the methods of so-called strong
or mean-square approximation of iterated stochastic
integrals. For numerical integration of It\^{o} SDEs 
the so-called weak approximations of iterated It\^{o} 
stochastic integrals from the Taylor--It\^{o} expansions
(see Chapter 4) are also interesting.

Let $(\Omega,$ ${\rm F},$ ${\sf P})$ be a complete probability space, let 
$\{{\rm F}_t, t\in[0,T]\}$ be a nondecreasing right-con\-ti\-nu\-o\-us family 
of $\sigma$-algebras of ${\rm F},$
and let ${\bf w}_t$ be a standard $m$-dimensional Wiener 
process, which is
${\rm F}_t$-measurable for all $t\in[0, T].$ We suppose that the components
${\bf w}_{t}^{(i)}$ $(i=1,\ldots,m)$ of this process are independent. 

Let us consider
an It\^{o} SDE in the integral form
\begin{equation}
\label{1.5.2xxa}
~~~~~~~ {\bf x}_t={\bf x}_0+\int\limits_0^t {\bf a}({\bf x}_{\tau},\tau)d\tau+
\int\limits_0^t B({\bf x}_{\tau},\tau)d{\bf w}_{\tau},\ \ \
{\bf x}_0={\bf x}(0,\omega),
\end{equation}
where ${\bf x}_t$ is some $n$-dimensional stochastic process 
satisfying to the It\^{o} SDE (\ref{1.5.2xxa}),
the nonrandom functions ${\bf a}: {\bf R}^n\times[0, T]\to{\bf R}^n$,
$B: {\bf R}^n\times[0, T]\to{\bf R}^{n\times m}$
guarantee the existence and uniqueness up to stochastic 
equivalence of a solution
of (\ref{1.5.2xxa}) \cite{Gih1},
${\bf x}_0$ is an $n$-dimensional random variable, which is 
${\rm F}_0$-measurable and 
${\sf M}\bigl\{\left|{\bf x}_0\right|^2\bigr\}<\infty$,
${\bf x}_0$ and ${\bf w}_t-{\bf w}_0$ are independent for $t>0.$

Let us consider the iterated It\^{o} stochastic integrals
from the classical Taylor--It\^{o} expansion (see Chapter 4)
$$
{J}_{(\lambda_{1}\ldots \lambda_k)s,t}^{(i_1\ldots i_k)}=
\int\limits_t^s\ldots
\int\limits_t^{\tau_{2}}
d{\bf w}_{t_{1}}^{(i_1)}\ldots
d{\bf w}_{t_k}^{(i_k)}\ \ \ (k\ge 1),
$$
where ${\bf w}_{{\tau}}^{(i)}$ $(i=1,\ldots,m)$ are independent standard
Wiener processes,
${\bf w}_{{\tau}}^{(0)}=\tau$,
$i_l=0$ if $\lambda_l=0$ 
and $i_l=1,\ldots,m$ if $\lambda_l=1$ $(l=1,\ldots,k)$.
Moreover, let
$$
{\rm M}_k=\biggl\{(\lambda_1,\ldots,\lambda_k):\
\lambda_l=0\ \hbox{or}\ 1,\ l=1,\ldots,k\biggr\}.
$$

Weak approximations of iterated It\^{o} stochastic stochastic integrals 
are formed or selected from the specific moment conditions 
\cite{Zapad-1}, \cite{Zapad-3}, \cite{Zapad-4}, 
\cite{Zapad-8}, \cite{Zapad-9} (see below) and 
they are significantly simpler than their mean-square analogues.
However, weak approximations are focused on the numerical 
solution of other problems \cite{Zapad-1}, \cite{Zapad-3}, \cite{Zapad-4}, 
\cite{Zapad-8}, \cite{Zapad-9} connected with 
It\^{o} SDEs in comparison with mean-square 
approximations.

We will say that the set of weak approximations
$$
{\hat J}_{(\lambda_{1}\ldots \lambda_k)s,t}^{(i_1\ldots i_k)}
$$ 
of the iterated It\^{o} stochastic integrals
$$
{J}_{(\lambda_{1}\ldots \lambda_k)s,t}^{(i_1\ldots, i_k)}
$$ 
from the Taylor--It\^{o} expansion (\ref{5.7.11}) has the order $r$, if 
\cite{Zapad-1}, \cite{Zapad-3}
for $t\in [t_0, T]$ and
$r\in {\bf N}$ there exists a constant $K\in(0,\infty)$ such that
the condition
\begin{equation}
\label{1uuuu}
\left|{\sf M}\left\{
\prod_{g=1}^l J_{\left(\lambda_{1}^{(g)}\ldots~\lambda_{k_g}^{(g)}\right)
t,t_0}^{\left(i_{1}^{(g)}\ldots~i_{k_g}^{(g)}\right)}-
\prod_{g=1}^l \hat J_{\left(\lambda_{1}^{(g)}
\ldots~\lambda_{k_g}^{(g)}\right)t,t_0}^{\left(i_{1}^{(g)}
\ldots~i_{k_g}^{(g)}\right)}\biggr|
{\rm F}_{t_0}\right\}\right|
\le K(t-t_0)^{r+1}\ \ \ \hbox{w.~p.~1}
\end{equation}

\noindent
is satisfied for all 
$\bigl(\lambda_{1}^{(g)}\ldots\lambda_{k_g}^{(g)}\bigr)\in {\rm M}_{k_g},$
$i_1^{(g)},\ldots,i_{k_g}^{(g)}=0, 1,\ldots,m,$
$k_g\le r,$ $g=1,\ldots,l,$ $l=1, 2,\ldots,2r+1.$

If we talk about the unified Taylor--It\^{o} expansion (\ref{razl4}), then
we will say that the set of weak approximations
$$
{\hat I}_{(l_1\ldots l_k)s,t}^{(i_1\ldots i_k)}
$$
of the iterated It\^{o} stochastic integrals
$$
I_{(l_1\ldots l_k)s,t}^{(i_1\ldots i_k)}
=\int\limits_t^s (t-t_k)^{l_k} \ldots \int\limits_t^{t_{2}}
(t-t_1)^{l_1} d{\bf w}_{t_1}^{(i_1)}\ldots
d{\bf w}_{t_k}^{(i_k)}\ \ \ (i_1,\ldots,i_k=1,\ldots,m)
$$
has the order $r$, if 
for $t\in [t_0,T]$ and
$r\in {\bf N}$ there exists 
a constant $K\in(0,\infty)$ such that the condition
$$
\left|{\sf M}\left\{
\prod_{g=1}^l
\frac{(t-t_0)^{j_g}}{j_g !}\
I_{\left(l_{1}^{(g)}\ldots~ l_{k_g}^{(g)}\right)t,t_0}^{\left(i_{1}^{(g)}\ldots~
i_{k_g}^{(g)}\right)}-
\prod_{g=1}^l 
\frac{(t-t_0)^{j_g}}{j_g !}\
\hat I_{\left(l_{1}^{(g)}\ldots~ l_{k_g}^{(g)}\right)t,t_0}
^{\left(i_{1}^{(g)}\ldots~
i_{k_g}^{(g)}\right)}\biggr|{\rm F}_{t_0}\right\}\right|\le
$$
\begin{equation}
\label{2uuuu}
\le K (t-t_0)^{r+1}
\end{equation}

\vspace{2mm}
\noindent
is satisfied w.~p.~1 for all $\bigl(k_g,j_g,l_1^{(g)},\ldots,~l_{k_g}^{(g)}\bigr)
\in {\rm A}_{q_g},$
$i_1^{(g)},\ldots,i_{k_g}^{(g)}=1,\ldots,m,$
$q_g\le r,$ $g=1,\ldots,l,$ $l=1, 2,\ldots,2r+1,$ 
where
$$
{\rm A} _ {q} = \left\{
(k, j, l_ {1},\ldots, l_ {k}): k + j +
\sum\limits_{p=1}^k l_ {p}= q;\ k, j, l_{1},\ldots, l_ {k} 
= 0,1,\ldots \right\}.
$$

The theory of weak approximations of iterated It\^{o} stochastic  
integrals is not so rich as the theory of mean-square approximations. 
On the one hand, it is connected with the sufficiency for practical needs 
of already found approximations \cite{Zapad-1}, \cite{Zapad-3}, 
\cite{Zapad-8}, and on the other hand, it is 
connected with the complexity of their formation owing to the necessity 
to satisfy a lot of moment conditions.

Let us consider the basic results in this area.

In \cite{Zapad-3} (also see \cite{Zapad-1}) the authors found
the weak approximations 
with the orders $r=1, 2$ when $m, n\ge 1$ 
as well as with the order $r=3$ when $m=1,$ $n\ge 1$
for iterated It\^{o} stochastic integrals
$$
J_{(\lambda_1\ldots \lambda_k)t,t_0}^{(i_1\ldots i_k)}.
$$

Recall that
$n$ is a dimension of the It\^{o} process ${\bf x}_t$, which is 
a solution of the It\^{o} SDE (\ref{1.5.2xxa}) 
and $m$ is a dimension of the Wiener process in (\ref{1.5.2xxa}).

Further, we will consider the
mentioned weak approximations as well as weak approximations 
with the order $r=4$ when $m=1,$ $n\ge 1$ \cite{2000} (2000) for
iterated It\^{o} stochastic integrals
$$
I_{{l_{1}\ldots l_{k}}_{t,t_0}}^{(i_{1}\ldots i_{k})}.
$$

In order to shorten the record let us write
\begin{equation}
\label{010.008}
{\sf M}\left\{\prod_{g=1}^l 
J_{\left(\lambda_{1}^{(g)}\ldots~
\lambda_{k_g}^{(g)}\right)t_0+\Delta,t_0}^{\left(i_{1}^{(g)}\ldots~
i_{k_g}^{(g)}\right)}\biggr|{\rm F}_{t_0}\right\}
\stackrel{\rm def}{=}
{\sf M}'\left\{\prod_{g=1}^l 
J_{\left(\lambda_{1}^{(g)}\ldots~\lambda_{k_g}^{(g)}
\right)}^{\left(i_{1}^{(g)}
\ldots~ i_{k_g}^{(g)}\right)}\right\},
\end{equation}
where $\Delta\in [0,T-t_0],$ 
$\bigl(\lambda_{1}^{(g)}\ldots\lambda_{k_g}^{(g)}\bigr)\in{\rm M}_{k_g},$
$k_g\le r,$ $g=1,\ldots,l.$

Note that the conditional expectation 
(\ref{010.008}) is equal w.~p.~1 to the corresponding expectation without the condition ${\rm F}_{t_0}$.
Further, equalities and inequalities for 
conditional expectations are understood w.~p.~1.
As before, ${\bf 1}_A$ means the indicator of the set $A$.

Let us consider the exact values of conditional expectations 
(\ref{010.008}) calculated in \cite{Zapad-1}, \cite{Zapad-3}
and necessary to 
form weak approximations 
$$
\hat J_{(\lambda_{1}\ldots\lambda_k)t_0+\Delta,t_0}^{(i_{1}\ldots i_k)}
$$ 
of the orders $r=1, 2$ when $m, n\ge 1$

\vspace{-1mm}
\begin{equation}
\label{010.009}
{\sf M}'\left\{J_{(1)}^{(i_1)}J_{(1)}^{(i_2)}\right\}=
\Delta {\bf 1}_{\{i_1=i_2\}},
\end{equation}

\vspace{-1mm}
\begin{equation}
\label{010.010}
{\sf M}'\left\{J_{(1)}^{(i_1)}J_{(01)}^{(0i_2)}\right\}=
{\sf M}'\left\{J_{(1)}^{(i_1)}J_{(10)}^{(i_2 0)}\right\}=
\frac{1}{2}\Delta^2 {\bf 1}_{\{i_1=i_2\}},
\end{equation}

\begin{equation}
\label{010.011}
{\sf M}'\left\{J_{(11)}^{(i_1 i_2)}J_{(11)}^{(i_3 i_4)}\right\}=
\frac{1}{2}\Delta^2 {\bf 1}_{\{i_1=i_3\}}{\bf 1}_{\{i_2=i_4\}},
\end{equation}

\newpage
\noindent
\begin{equation}
\label{010.012}
{\sf M}'\left\{J_{(1)}^{(i_1)}J_{(1)}^{(i_2)}J_{(11)}^{(i_3 i_4)}\right\}=
\left\{
\begin{matrix}
\Delta^2\ &\hbox{when}\ i_1=\ldots=i_4 \cr\cr
\Delta^2/2\ &
\begin{matrix}
&\hbox{when}\ i_3\ne i_4,\ i_1=i_3,\ i_2=i_4 \cr
&\hbox{or}\ i_3\ne i_4,\ i_1=i_4,\ i_2=i_3
\end{matrix} \cr\cr
0\ &\hbox{otherwise}
\end{matrix},\right.
\end{equation}

\vspace{4mm}

\begin{equation}
\label{010.013}
{\sf M}'\left\{J_{(1)}^{(i_1)}J_{(1)}^{(i_2)}
J_{(1)}^{(i_3)}J_{(1)}^{(i_4)}\right\}=
\left\{
\begin{matrix}3\Delta^2\ &\hbox{when}\ i_1=\ldots=i_4 \cr\cr
\Delta^2\ &
\begin{matrix}
&\hbox{if among}\ i_1,\ldots,i_4\ \hbox{there are}\cr
&\hbox{two pairs of identical numbers}
\end{matrix} \cr\cr
0\ &\hbox{otherwise}
\end{matrix},\right.
\end{equation}

\vspace{4mm}

\begin{equation}
\label{010.014}
{\sf M}'\left\{J_{(10)}^{(i_1 0)}J_{(01)}^{(0 i_2)}\right\}=
\frac{1}{6}\Delta^3 {\bf 1}_{\{i_1=i_2\}},
\end{equation}

\begin{equation}
\label{010.015}
{\sf M}'\left\{J_{(10)}^{(i_1 0)}J_{(10)}^{(i_2 0)}\right\}=
{\sf M}'\left\{J_{(01)}^{(0 i_1)}J_{(01)}^{(0 i_2)}\right\}=
\frac{1}{3}\Delta^3 {\bf 1}_{\{i_1=i_2\}},
\end{equation}

\vspace{4mm}

$$
{\sf M}'\left\{J_{(01)}^{(0 i_1)}J_{(1)}^{(i_2)}
J_{(1)}^{(i_3)}J_{(1)}^{(i_4)}\right\}=
{\sf M}'\left\{J_{(10)}^{(i_1 0)}J_{(1)}^{(i_2)}
J_{(1)}^{(i_3)}J_{(1)}^{(i_4)}\right\}=
$$

\begin{equation}
\label{010.016}
=\left\{
\begin{matrix}
3\Delta^3/2\ &\hbox{when}\ i_1=\ldots=i_4 \cr\cr
\Delta^3/2\ &
\begin{matrix}&\hbox{if among}\ i_1,\ldots,i_4\ 
\hbox{there are}\cr
&\hbox{two pairs of identical numbers}
\end{matrix} \cr\cr
0\ &\hbox{otherwise}
\end{matrix},\right.
\end{equation}

\vspace{4mm}

\begin{equation}
\label{010.017}
{\sf M}'\left\{J_{(01)}^{(0i_1)}J_{(1)}^{(i_2)}J_{(11)}^{(i_3i_4)}\right\}=
\frac{1}{6}\Delta^3 {\bf 1}_{\{i_1=i_3\}}{\bf 1}_{\{i_2=i_4\}},
\end{equation}

\begin{equation}
\label{010.018}
{\sf M}'\left\{J_{(10)}^{(i_1 0)}J_{(1)}^{(i_2)}J_{(11)}^{(i_3i_4)}\right\}=
\frac{1}{3}\Delta^3 {\bf 1}_{\{i_1=i_3\}}{\bf 1}_{\{i_2=i_4\}},
\end{equation}

\newpage
\noindent
\begin{equation}
\label{010.019}
{\sf M}'\left\{J_{(1)}^{(i_1)}\ldots J_{(1)}^{(i_6)}\right\} =
\left\{
\begin{matrix}
15\Delta^3\ &\hbox{when}\ i_1=\ldots=i_6 \cr\cr
3\Delta^3\ &
\begin{matrix}
&\hbox{if among}\ i_1,\ldots,i_6\ 
\hbox{there is a pair}\cr
&\hbox{and a quad of identical numbers}
\end{matrix} \cr\cr
\Delta^3\ &
\begin{matrix}
&\hbox{ if among}\ i_1,\ldots,i_6\ 
\hbox{there are three}\cr
&\hbox{ pairs of identical numbers}
\end{matrix} \cr\cr
0\ &\hbox{otherwise}
\end{matrix},\right.
\end{equation}

\vspace{2mm}
$$
{\sf M}'\left\{J_{(11)}^{(i_1 i_2)}
J_{(11)}^{(i_3i_4)}J_{(11)}^{(i_5i_6)}\right\}=
$$
$$
=
\frac{1}{6}\Delta^3\Biggl(
{\bf 1}_{\{i_2=i_4\}}\left({\bf 1}_{\{i_1=i_5\}}
{\bf 1}_{\{i_3=i_6\}}+{\bf 1}_{\{i_1=i_6\}}{\bf 1}_{\{i_3=i_5\}}\right)+
\Biggr.
$$
$$
+{\bf 1}_{\{i_2=i_6\}}\left({\bf 1}_{\{i_1=i_3\}}
{\bf 1}_{\{i_4=i_5\}}+{\bf 1}_{\{i_1=i_4\}}{\bf 1}_{\{i_3=i_5\}}\right)+
$$
\begin{equation}
\label{010.020}
\Biggl.+{\bf 1}_{\{i_4=i_6\}}\left({\bf 1}_{\{i_1=i_3\}}
{\bf 1}_{\{i_2=i_5\}}+{\bf 1}_{\{i_2=i_3\}}{\bf 1}_{\{i_1=i_5\}}\right)
\Biggr),
\end{equation}

\vspace{3mm}
$$
{\sf M}'\left\{J_{(11)}^{(i_1 i_2)}
J_{(11)}^{(i_3i_4)}J_{(1)}^{(i_5)}J_{(1)}^{(i_6)}\right\}=
$$
$$
=
\frac{1}{2}\Delta^3
{\bf 1}_{\{i_1=i_3\}}{\bf 1}_{\{i_2=i_4\}}{\bf 1}_{\{i_5=i_6\}}+
$$
$$
+\frac{1}{6}\Delta^3\Biggl(
2\cdot{\bf 1}_{\{i_1=i_3\}}\left({\bf 1}_{\{i_2=i_5\}}
{\bf 1}_{\{i_4=i_6\}}+{\bf 1}_{\{i_2=i_6\}}{\bf 1}_{\{i_4=i_5\}}\right)+
\Biggr.
$$
$$
+{\bf 1}_{\{i_2=i_3\}}\left({\bf 1}_{\{i_1=i_5\}}
{\bf 1}_{\{i_4=i_6\}}+{\bf 1}_{\{i_1=i_6\}}{\bf 1}_{\{i_4=i_5\}}\right)+
$$
$$
+{\bf 1}_{\{i_1=i_4\}}\left({\bf 1}_{\{i_3=i_5\}}
{\bf 1}_{\{i_2=i_6\}}+{\bf 1}_{\{i_3=i_6\}}{\bf 1}_{\{i_2=i_5\}}\right)+
$$
\begin{equation}
\label{010.021}
\Biggl.+2\cdot{\bf 1}_{\{i_2=i_4\}}\left({\bf 1}_{\{i_1=i_5\}}
{\bf 1}_{\{i_3=i_6\}}+{\bf 1}_{\{i_3=i_5\}}{\bf 1}_{\{i_1=i_6\}}\right)
\Biggr),
\end{equation}

$$
{\sf M}'\left\{J_{(11)}^{(i_1 i_2)}
J_{(1)}^{(i_3)}\ldots J_{(1)}^{(i_6)}\right\}=
$$
\begin{equation}
\label{010.022}
~~~~~~ =\frac{1}{2}\left({\sf M}'\left\{
J_{(1)}^{(i_1)}\ldots J_{(1)}^{(i_6)}\right\}-
\Delta{\bf 1}_{\{i_1=i_2\}}{\sf M}'\left\{
J_{(1)}^{(i_3)}\ldots J_{(1)}^{(i_6)}\right\}\right).
\end{equation}

\vspace{3mm}

Let us explain the formula 
(\ref{010.021}).
From the following equality 
$$
J_{(1)}^{(i_5)}J_{(1)}^{(i_6)}=
J_{(11)}^{(i_5i_6)}+J_{(11)}^{(i_6i_5)}+\Delta{\bf 1}_{\{i_5=i_6\}}\ \ \
\hbox{w.~p.~1}
$$

\noindent
we obtain
$$
{\sf M}'\left\{J_{(11)}^{(i_1 i_2)}
J_{(11)}^{(i_3i_4)}J_{(1)}^{(i_5)}J_{(1)}^{(i_6)}\right\}=
{\sf M}'\left\{J_{(11)}^{(i_1 i_2)}
J_{(11)}^{(i_3i_4)}J_{(11)}^{(i_5i_6)}\right\}+
$$
\begin{equation}
\label{010.02200}
~~~~~~~~ +
{\sf M}'\left\{J_{(11)}^{(i_1 i_2)}
J_{(11)}^{(i_3i_4)}J_{(11)}^{(i_6i_5)}\right\}+
\Delta{\bf 1}_{\{i_5=i_6\}}
{\sf M}'\left\{J_{(11)}^{(i_1 i_2)}
J_{(11)}^{(i_3i_4)}\right\}.
\end{equation}

\vspace{2mm}

Applying (\ref{010.011}), (\ref{010.020})
to the right-hand side of (\ref{010.02200}) 
gives (\ref{010.021}).
It is necessary to note \cite{Zapad-1}, \cite{Zapad-3} that
$$
{\sf M}'\left\{\prod\limits_{g=1}^l 
J_{\left(\lambda_{1}^{(g)}\ldots~ \lambda_{k_g}^{(g)}\right)}
^{\left(i_{1}^{(g)}\ldots~ i_{k_g}^{(g)}\right)}\right\}=0
$$
if the number of units included
in all multi-indices
$\bigl(\lambda_{1}^{(g)}\ldots\lambda_{k_g}^{(g)}\bigr)$
is odd ($k_g\le r,$ $g=1,\ldots,l$).
In addition \cite{Zapad-1}, \cite{Zapad-3}
$$
\left|{\sf M}'\left\{\prod\limits_{g=1}^l J_{\left(\lambda_{1}^{(g)}
\ldots~ \lambda_{k_g}^{(g)}\right)}
^{\left(i_{1}^{(g)}\ldots~ i_{k_g}^{(g)}\right)}
\right\}\right|\le K\Delta^{\gamma_l},
$$
where $\gamma_l=\delta_l/2+\rho_l,$
$\delta_l$ is a number of units and $\rho_l$ is a number of zeros
included in all multi-indices
$\bigl(\lambda_{1}^{(g)}\ldots\lambda_{k_g}^{(g)}\bigr),$
$k_g\le r,$ $g=1,\ldots,l,$
$K\in (0,\infty)$ is a constant.

In the case $n, m\ge 1$ and $r=1$ we can put
\cite{Zapad-1}, \cite{Zapad-3}
$$
\hat J_{(1)}^{(i)}=\Delta\tilde{\bf w}^{(i)}\ \ \
(i=1,\ldots,m),
$$
where $\Delta\tilde{\bf w}^{(i)},$ $i=1,\ldots,m$ are independent 
discrete random variables
for which
$$
{\sf P}\left\{\Delta\tilde{\bf w}^{(i)}=\pm\sqrt{\Delta}\right\}=
\frac{1}{2}.
$$

It is not difficult to see that the approximation
$$
\hat J_{(1)}^{(i)}=\sqrt{\Delta}\zeta_0^{(i)}\ \ \ (i=1,\ldots,m)
$$
also satisfies the condition
(\ref{1uuuu}) when $r=1.$
Here $\zeta_0^{(i)}$ are independent standard 
Gaussian random variables.

In the case $n, m\ge 1$ and $r=2$ as the approximations
$\hat J_{(1)}^{(i_1)},$ 
$\hat J_{(11)}^{(i_1 i_2)},$ 
$\hat J_{(10)}^{(i_1 0)},$ 
$\hat J_{(01)}^{(0 i_1)}$ 
are taken the following ones \cite{Zapad-3}
\begin{equation}
\label{010.023}
\hat J_{(1)}^{(i_1)}=\Delta\tilde{\bf w}^{(i_1)},\ \ \
\hat J_{(10)}^{(i_1 0)}=\hat J_{(01)}^{(0 i_1)}=
\frac{1}{2}\Delta\cdot\Delta\tilde{\bf w}^{(i_1)},
\end{equation}
\begin{equation}
\label{010.024}
\hat J_{(11)}^{(i_1 i_2)}=\frac{1}{2}\left(
\Delta\tilde{\bf w}^{(i_1)}\Delta\tilde{\bf w}^{(i_2)}+
V^{(i_1i_2)}\right),
\end{equation}

\noindent
where $\Delta\tilde{\bf w}^{(i)}$
are independent Gaussian random variables with zero  
expectation and variance $\Delta$
or independent discrete random variables for which the following 
conditions are fulfilled
$$
{\sf P}\left\{\Delta\tilde{\bf w}^{(i)}=\pm\sqrt{3\Delta}\right\}=
\frac{1}{6},
$$
$$
{\sf P}\left\{\Delta\tilde{\bf w}^{(i)}=0\right\}=\frac{2}{3},
$$
$V^{(i_1i_2)}$ are independent discrete random variables
satisfying the conditions
$$
{\sf P}\left\{V^{(i_1i_2)}=\pm\Delta\right\}=\frac{1}{2}\ \ \
\hbox{when}\ \ \ i_2<i_1,
$$
$$
V^{(i_1i_1)}=-\Delta,\ \ \ V^{(i_1i_2)}=-V^{(i_2i_1)}\ \ \ 
\hbox{when}\ \ \ i_1<i_2,
$$
where $i_1, i_2=1,\ldots,m.$

Let us consider the case $r=3$ and $m=1, n\ge 1.$
In this situation in addition to 
the formulas (\ref{010.009})--(\ref{010.02200}) 
we need a number of formulas for the conditional   
expectations (\ref{010.008}) when $m=1.$

We have \cite{Zapad-1}, \cite{Zapad-3} 
$$
{\sf M}'\{J_{(1)}J_{(111)}\}={\sf M}'\{J_{(01)}J_{(111)}\}=
{\sf M}'\{J_{(10)}J_{(111)}\}=0,
$$
$$
{\sf M}'\{J_{(011)}J_{(11)}\}={\sf M}'\{J_{(101)}J_{(11)}\}=
{\sf M}'\{J_{(110)}J_{(11)}\}=\frac{1}{6}\Delta^3,
$$
$$
{\sf M}'\{J_{(001)}J_{(1)}\}={\sf M}'\{J_{(010)}J_{(1)}\}=
{\sf M}'\{J_{(100)}J_{(1)}\}=\frac{1}{6}\Delta^3,
$$
$$
{\sf M}'\{J_{(100)}J_{(10)}\}={\sf M}'\{J_{(001)}J_{(01)}\}=
\frac{1}{8}\Delta^4,\ \ \ {\sf M}'\{J_{(111)}J_{(11)}\}=0,\
$$
$$
{\sf M}'\{J_{(010)}J_{(10)}\}={\sf M}'\{J_{(010)}J_{(01)}\}=
\frac{1}{6}\Delta^4,\ \ \
{\sf M}'\left\{\left(J_{(111)}\right)^2\right\}=\frac{1}{6}\Delta^3,
$$
$$
{\sf M}'\{J_{(100)}J_{(01)}\}={\sf M}'\{J_{(001)}J_{(10)}\}=
\frac{1}{24}\Delta^4,
$$
$$
{\sf M}'\{J_{(110)}J_{(10)}\}={\sf M}'\{J_{(110)}J_{(01)}\}=
{\sf M}'\{J_{(101)}J_{(10)}\}=0,
$$
$$
{\sf M}'\{J_{(101)}J_{(01)}\}={\sf M}'\{J_{(011)}J_{(10)}\}=
{\sf M}'\{J_{(011)}J_{(01)}\}=0,
$$
$$
{\sf M}'\left\{J_{(011)}\left(J_{(1)}\right)^2\right\}=
{\sf M}'\left\{J_{(101)}\left(J_{(1)}\right)^2\right\}=
{\sf M}'\left\{J_{(110)}\left(J_{(1)}\right)^2\right\}=
\frac{1}{6}\Delta^3,
$$
$$
{\sf M}'\left\{J_{(111)}\left(J_{(1)}\right)^3\right\}=\Delta^3,\ \ \
{\sf M}'\left\{J_{(111)}J_{(11)}J_{(1)}\right\}=\frac{1}{2}\Delta^3,
$$
where
$$
J_{(\lambda_1\ldots\lambda_k)}\stackrel{\rm def}{=}
\int\limits_{t_0}^{t_0+\Delta}\ldots
\int\limits_{t_0}^{t_{2}}d{w}_{t_1}^{(\lambda_1)}
\ldots d{w}_{t_k}^{(\lambda_k)},
$$
${w}^{(0)}_{t}\stackrel{\rm def}{=}t,$
${w}_{t}^{(1)}\stackrel{\rm def}{=}{w}_{t}$
is standard scalar Wiener process,
$\lambda_l=0$ or $\lambda_l=1$, $l=1,\ldots,k.$

In \cite{Zapad-1}, \cite{Zapad-3} using the given moment relations  
the authors proposed the following weak approximations of iterated It\^{o} 
stochastic integrals for $r=3$ when $m=1,$ $n\ge 1$
\begin{equation}
\label{aa.1}
\hat J_{(1)}=\Delta\tilde{w},
\end{equation}
\begin{equation}
\label{bb.1}
\hat J_{(10)}=\Delta\hat w,\ \ \ \hat J_{(01)}=\Delta\cdot
\Delta\tilde w-\Delta\hat w,
\end{equation}
$$
\hat J_{(11)}
=\frac{1}{2}\left(\left(\Delta\tilde w\right)^2-\Delta\right),\ \ \
\hat J_{(001)}=\hat J_{(010)}=\hat J_{(100)}=
\frac{1}{6}\Delta^2\cdot\Delta\tilde w,
$$
$$
\hat J_{(110)}=\hat J_{(101)}=\hat J_{(011)}=
\frac{1}{6}\Delta\left(\left(\Delta\tilde w\right)^2-
\Delta\right),
$$
$$
\hat J_{(111)}=\frac{1}{6}\Delta\tilde w
\left(\left(\Delta\tilde w\right)^2-
3\Delta\right),
$$

\noindent
where 
$$
\Delta\tilde w \sim {\rm N}(0,\Delta),\ \ \ \Delta\hat w \sim
{\rm N}\left(0,\frac{1}{3}\Delta^3\right),\ \ \ 
{\sf M}\left\{\Delta\tilde w\Delta\hat w\right\}=
\frac{1}{2}\Delta^2.
$$
Here ${\rm N}(0,\sigma^2)$ is a Gaussian distribution with zero
expectation and variance $\sigma^2.$

Finally, we will form the weak approximations of 
iterated It\^{o} stochastic integrals for $r=4$ when $m=1,$ $n\ge 1$ 
\cite{1}-\cite{12aa}.

The truncated Taylor--It\^{o} expansion
(\ref{5.7.11})
when $r=4$ and $m=1$ includes
26 various iterated It\^{o} stochastic integrals. 
The formation of weak approximations for these stochastic integrals
satisfying the condition (\ref{2uuuu}) when $r=4$ is extremely 
difficult due to the necessity to consider a lot of moment conditions.
However, this problem can be simplified if we consider the truncated 
unified Taylor--It\^{o} expansion (\ref{razl4}) 
when $r=4$ and $m=1,$ since this expansion includes only
15 various iterated It\^{o} stochastic integrals
$$
I_0,\   I_{1},\   I_{00},\ 
I_{000},\   I_{2},\   I_{10},\   I_{01},\    I_{3},\ 
I_{11},\   I_{20},\  
I_{02},\   I_{100},\   I_{010},\   I_{001},\   I_{0000},
$$
where
$$
I_{l_1\ldots l_k}\stackrel{\rm def}{=}
\int\limits_{t_0}^{t_0+\Delta}(t_0-t_k)^{l_k}
\ldots\int\limits_{t_0}^{t_2}(t_0-t_1)^{l_1}dw_{t_1}\ldots
dw_{t_k}\ \ \ (k\ge 1)
$$
and $w_{t}$ is standard scalar Wiener process.

It is not 
difficult to notice that the condition (\ref{2uuuu}) will be satisfied
for $r=4$ and $i_ 1=\ldots=i_ 4$ if the following more 
strong condition is fulfilled
\begin{equation}
\label{010.026}
~~~~~~~~~ \left|{\sf M}\left\{
\prod_{g=1}^l I_{l_{1}^{(g)}\ldots~ l_{k_g}^{(g)}}-
\prod_{g=1}^l \hat I_{l_{1}^{(g)}\ldots~ l_{k_g}^{(g)}}
\biggr|{\rm F}_{t_0}\right\}\right|
\le K(t-t_0)^{5}\ \ \ \hbox{w.~p.~1}
\end{equation}

\noindent
for all $l_1^{(g)}\ldots~ l_{k_g}^{(g)}\in A,$
$k_g\le 4,$ $g=1,\ldots,l,$ $l=1, 2,\ldots,9,$
where $K\in(0,\infty)$ and
$$
A=\biggl\{0, 1, 00, 000, 2, 10, 01, 3, 11, 20,
02, 100, 010, 001, 0000\biggr\}
$$
is the set of multi-indices. 

Let (see Sect.~5.1 and 6.6) \cite{12a}-\cite{12aa}, \cite{2000}
\begin{equation}
\label{16.00000}
\hat I_0=\sqrt{\Delta}\zeta_0,\ \ \ \hat I_{00}=\frac{1}{2}\Delta
\left(\bigl(\zeta_0\bigr)^2-1\right),
\end{equation}
\begin{equation}
\label{17.00000}
~~~~~~~~~ \hat I_1=-\frac{\Delta^{3/2}}{2}\left(\zeta_0+
\frac{1}{\sqrt{3}}\zeta_1\right),\ \ \
\hat I_{000}=\frac{\Delta^{3/2}}{6}\left(\bigl(\zeta_0\bigr)^3-
3\zeta_0\right),
\end{equation}
\begin{equation}
\label{18.00000}
\hat I_{0000}=\frac{\Delta^{2}}{24}\left(\bigl(\zeta_0\bigr)^4-
6\bigl(\zeta_0\bigr)^2+3\right).
\end{equation}

\noindent
Here and further
$$
\zeta_0\stackrel{\rm def}{=}\frac{1}{\sqrt{\Delta}}
\int\limits_{t_0}^{t_0+\Delta}
dw_s,\ \ \ 
\zeta_1\stackrel{\rm def}{=}\frac{2\sqrt{3}}{\Delta^{3/2}}
\int\limits_{t_0}^{t_0+\Delta}\left(s-t_0-\frac{\Delta}{2}\right)
dw_s,
$$
where $w_s$ is scalar standard Wiener process.

It is not difficult to see that $\zeta_0,$ $\zeta_1$ 
are independent standard Gaussian random variables.
In addition, the approximations (\ref{16.00000})--(\ref{18.00000}) 
equal w.~p.~1 to the iterated It\^{o} stochastic integrals
corresponding to these approximations. 
This implies that all products 
$$
\prod\limits_{g=1}^l\hat I_{l_1^{(g)}\ldots~ l_{k_g}^{(g)}},
$$
which contain only the approximations 
(\ref{16.00000})--(\ref{18.00000})
will convert the left-hand side of (\ref{010.026}) to zero w.~p.~1, i.e. 
the condition (\ref{010.026}) will be fulfilled automatically.

For forming the approximations 
$$
\hat I_{100},\
\hat I_{010},\ \hat I_{001},\ \hat I_{10},\ \hat I_{01},\
\hat I_{11},\ \hat I_{20},\ \hat I_{02},\ \hat I_{2},\ \hat I_{3}
$$
it is necessary to calculate several conditional 
expectations
\begin{equation}
\label{fas1}
{\sf M}\left\{\prod\limits_{g=1}^l 
I_{l_1^{(g)}\ldots~ l_{k_g}^{(g)}}\biggr|
{\rm F}_{t_0}\right\},
\end{equation}  
where $l_1^{(g)}\ldots~ l_{k_g}^{(g)}\in A.$
Note that the conditional expectation 
(\ref{fas1}) is equal w.~p.~1 to the corresponding expectation without the condition ${\rm F}_{t_0}$.
We will denote (\ref{fas1}) (as before) as follows
$$
{\sf M}'\left\{\prod\limits_{g=1}^l 
I_{l_1^{(g)}\ldots~ l_{k_g}^{(g)}}\right\}.
$$

We have
$$
{\sf M}'\{I_3\}={\sf M}'\{I_3(I_0)^2\}={\sf M'}\{I_3I_{00}\}=0,\ \ \
{\sf M}'\{I_3I_0\}=-\frac{\Delta^4}{4},
$$
$$
{\sf M}'\{I_2(I_0)^2\}={\sf M}'\{I_2I_{00}\}={\sf M}'\{I_2I_{000}\}=
{\sf M}'\{I_2I_{0000}\}=0,
$$

\vspace{-3mm}
$$
{\sf M}'\{I_2(I_{00})^2\}={\sf M}'\{I_2(I_0)^4\}={\sf M}'\{I_2I_{000}I_0\}=0,
$$

\vspace{-7mm}
$$
{\sf M}'\{I_2I_{00}(I_0)^2\}={\sf M}'\{I_2I_{10}\}=
{\sf M}'\{I_2I_{01}\}={\sf M}'\{I_2I_1I_0\}={\sf M}'\{I_2\}=0,
$$
$$
{\sf M}'\{I_2I_0\}=\frac{\Delta^3}{3},\ \ \ {\sf M}'\{I_2(I_0)^3\}=\Delta^4,\
\ \
{\sf M}'\{I_2I_{00}I_0\}=\frac{\Delta^4}{3},\ 
$$
$$
{\sf M}'\{I_2I_1\}=
-\frac{\Delta^4}{4},\ \ \
{\sf M}'\{I_{\mu}\}={\sf M}'\{I_{\mu}I_0\}={\sf M}'\{I_{\mu}I_{000}\}=
{\sf M}'\{I_{\mu}(I_0)^3\}=0,
$$
$$
{\sf M}'\{I_{\mu}I_{00}I_0\}={\sf M}'\{I_{\mu}I_1\}=0,\ \ \
{\sf M}'\{I_{20}(I_0)^2\}=\frac{\Delta^4}{6},\ \ \
{\sf M}'\{I_{20}I_{00}\}=\frac{\Delta^4}{12},
$$
$$
{\sf M}'\{I_{11}(I_0)^2\}=\frac{\Delta^4}{4},\ \ \
{\sf M}'\{I_{11}I_{00}\}=\frac{\Delta^4}{8},\ \ \
{\sf M}'\{I_{02}(I_0)^2\}=\frac{\Delta^4}{2},
$$
$$
{\sf M}'\{I_{02}I_{00}\}=\frac{\Delta^4}{4},\ \ \
{\sf M}'\{I_{\lambda}\}={\sf M}'\{I_{\lambda}I_0\}=
{\sf M}'\{I_{\lambda}(I_0)^2\}={\sf M}'\{I_{\lambda}I_{00}\}=0,
$$
$$
{\sf M}'\{I_{\lambda}I_1\}={\sf M}'\{I_{\lambda}I_{0000}\}=
{\sf M}'\{I_{\lambda}(I_{00})^2\}={\sf M}'\{I_{\lambda}(I_0)^4\}=0,
$$

\vspace{-5mm}
$$
{\sf M}'\{I_{\lambda}I_{000}I_0\}={\sf M}'\{I_{\lambda}I_{00}(I_0)^2\}=
{\sf M}'\{I_{\lambda}I_{10}\}=0,\
$$

\vspace{-5mm}
$$
{\sf M}'\{I_{\lambda}I_{01}\}={\sf M}'\{I_{\lambda}I_1I_0\}=0,
$$
$$
{\sf M}'\{I_{100}I_{000}\}=-\frac{\Delta^4}{24},\ \ \
{\sf M}'\{I_{100}(I_{0})^3\}=-\frac{\Delta^4}{4},\ \ \
{\sf M}'\{I_{100}I_{00}I_0\}=-\frac{\Delta^4}{8},
$$
$$
{\sf M}'\{I_{010}I_{000}\}=-\frac{\Delta^4}{12},\ \ \
{\sf M}'\{I_{010}(I_{0})^3\}=-\frac{\Delta^4}{2},\ \ \
{\sf M}'\{I_{010}I_{00}I_0\}=-\frac{\Delta^4}{4},
$$
$$
{\sf M}'\{I_{001}I_{000}\}=-\frac{\Delta^4}{8},\ \ \
{\sf M}'\{I_{001}(I_{0})^3\}=-\frac{3\Delta^4}{4},\ \ \
{\sf M}'\{I_{001}I_{00}I_0\}=-\frac{3\Delta^4}{8},
$$
$$
{\sf M}'\{I_{\rho}I_0\}={\sf M}'\{I_{\rho}I_{000}\}=
{\sf M}'\{I_{\rho}(I_0)^3\}={\sf M}'\{I_{\rho}I_{00}I_0\}=0,
$$

\vspace{-5mm}
$$
{\sf M}'\{I_{\rho}I_1\}={\sf M}'\{I_{\rho}I_{0000}\}=
{\sf M}'\{I_{\rho}(I_0)^5\}={\sf M}'\{I_{\rho}(I_{00})^2I_0\}=0,
$$

\vspace{-5mm}
$$
{\sf M}'\{I_{\rho}I_{00}(I_0)^3\}={\sf M}'\{I_{\rho}I_{000}(I_0)^2\}=
{\sf M}'\{I_{\rho}I_{0000}I_0\}=0,
$$

\vspace{-5mm}
$$
{\sf M}'\{I_{\rho}I_{000}I_{00}\}={\sf M}'\{I_{\rho}I_{100}\}=
{\sf M}'\{I_{\rho}I_{010}\}=0,
$$

\vspace{-5mm}
$$
{\sf M}'\{I_{\rho}I_{001}\}={\sf M}'\{I_{\rho}I_2\}=
{\sf M}'\{(I_{\rho})^2I_0\}={\sf M}'\{I_{\rho}I_{00}I_1\}=0,
$$

\vspace{-5mm}
$$
{\sf M}'\{I_{10}I_{01}I_0\}={\sf M}'\{I_{\rho}\}=
{\sf M}'\{I_{\rho}I_1(I_0)^2\}=0,
$$
$$
{\sf M}'\{I_{10}(I_0)^2\}=-\frac{\Delta^3}{3},\ \ \
{\sf M}'\{I_{10}I_{00}\}=-\frac{\Delta^3}{6},\ \ \
{\sf M}'\{I_{10}(I_{00})^2\}=-\frac{\Delta^4}{3},
$$
$$
{\sf M}'\{I_{10}(I_{0})^4\}=-2\Delta^4,\ \ \
{\sf M}'\{I_{10}I_{000}I_0\}=-\frac{\Delta^4}{6},
$$
$$
{\sf M}'\{I_{10}I_{00}(I_0)^2\}=-\frac{5\Delta^4}{6},
$$
$$
{\sf M}'\{(I_{10})^2\}=\frac{\Delta^4}{12},\ \ \
{\sf M}'\{I_{10}I_{01}\}=\frac{\Delta^4}{8},\ \ \
{\sf M}'\{I_{10}I_{1}I_0\}=\frac{5\Delta^4}{24},
$$
$$
{\sf M}'\{I_{01}(I_0)^2\}=-\frac{2\Delta^3}{3},\ \ \
{\sf M}'\{I_{01}I_{00}\}=-\frac{\Delta^3}{3},\ \ \
{\sf M}'\{I_{01}(I_{00})^2\}=-\frac{2\Delta^4}{3},
$$
$$
{\sf M}'\{I_{01}(I_{0})^4\}=-4\Delta^4,\ \ \
{\sf M}'\{I_{01}I_{000}I_0\}=-\frac{\Delta^4}{3},
$$
$$
{\sf M}'\{I_{01}I_{00}(I_0)^2\}=-\frac{5\Delta^4}{3},\ \ \
{\sf M}'\{(I_{01})^2\}=\frac{\Delta^4}{4},\ \ \
{\sf M}'\{I_{01}I_{1}I_0\}=\frac{3\Delta^4}{8},
$$

\noindent
where 
$$
\mu=02,\  11,\  20,\ \ \
\lambda=100,\  010,\  001,\ \ \
\rho=10,\  01
$$
(these recordings should be understood as sequences of digits).

The above relations are obtained using the standard properties
of the It\^{o} stochastic integral and the following  
equalities resulting from the It\^{o} formula
$$
(I_0)^4=24 I_{0000}+12\Delta I_{00}+3\Delta^2,\ \ \ 
(I_{00})^2=6I_{0000}+2\Delta I_{00}+\frac{\Delta^2}{2},
$$
$$
I_{00}(I_0)^2=12 I_{0000}+5\Delta I_{00}+\Delta^2,\ \ \
I_1I_0=I_{10}+I_{01}-\frac{\Delta^2}{2},\
$$
$$
I_{00}(I_0)^3=60I_{00000}+27\Delta I_{000}+6\Delta^2 I_0,
$$

\vspace{-5mm}
$$
(I_{0})^5=120I_{00000}+60\Delta I_{000}+15\Delta^2 I_0,\
$$
$$
(I_{00})^2 I_0=30I_{00000}+12\Delta I_{000}+\frac{10\Delta^2}{4}I_0,
$$
$$
I_{000}(I_0)^2=20I_{00000}+7\Delta I_{000}+\Delta^2 I_0,\ \ \
I_{0000}I_0=5I_{00000}+\Delta I_{000},
$$
$$
I_{000}I_{00}=10I_{00000}+3\Delta I_{000}+\frac{\Delta^2}{2}I_0,\ \ \
I_{00}I_1=I_{001}+I_{010}+I_{100}-\frac{\Delta^2}{2}I_0,
$$
$$
(I_0)^3=6I_{000}+3\Delta I_0,\ \ \
I_{00}I_0=3I_{000}+\Delta I_0,\
$$
$$
I_{10}I_0=I_{010}+I_{100}+\Delta I_1+I_2,\ \ \
I_{000}I_0=4 I_{0000}+\Delta I_{00},\ \ \ 
(I_0)^2=2I_{00}+\Delta,\
$$
$$
I_{01}I_0=2I_{001}+I_{010}-\frac{1}{2}\left(I_2+\Delta^2 I_0\right)
$$
w.~p.~1.

Using the given before moment relations, we can form the weak
approximations
$\hat I_{100},$
$\hat I_{010},$ $\hat I_{001},$ $\hat I_{10},$ $\hat I_{01},$
$\hat I_{11},$ $\hat I_{20},$ $\hat I_{02},$ $\hat I_{2},$ $\hat I_{3}$
\cite{12a}-\cite{12aa}, \cite{2000}
\begin{equation}
\label{19.00000}
~~~~~~~~\hat I_{100}=
-\frac{\Delta^{5/2}}{24}\left(\bigl(\zeta_0\bigr)^3-3\zeta_0\right),\ \ \
\hat I_{010}=
-\frac{\Delta^{5/2}}{12}\left(\bigl(\zeta_0\bigr)^3-3\zeta_0\right),
\end{equation}
\begin{equation}
\label{20.00000}
~~~~~~~~\hat I_{001}=
-\frac{\Delta^{5/2}}{8}\left(\bigl(\zeta_0\bigr)^3-3\zeta_0\right),\ \ \
\hat I_{11}=\frac{\Delta^3}{8}\left(\bigl(\zeta_0\bigr)^2-1\right),
\end{equation}
\begin{equation}
\label{21.00000}
\hat I_{20}=\frac{\Delta^3}{12}\left(\bigl(\zeta_0\bigr)^2-1\right),\ \ \
\hat I_{02}=\frac{\Delta^3}{4}\left(\bigl(\zeta_0\bigr)^2-1\right),
\end{equation}
\begin{equation}
\label{22.00000}
\hat I_3=-\frac{\Delta^{7/2}}{4}\zeta_0,\ \ \
\hat I_2=\frac{\Delta^{5/2}}{3}\left(\zeta_0+\frac{\sqrt{3}}{2}\zeta_1\right),
\end{equation}
\begin{equation}
\label{23.00000}
~~~~~~\hat I_{10}=\Delta^2\Biggl(-\frac{1}{6}\left(\bigl(\zeta_0\bigr)^2-1\right)
-\frac{1}{4\sqrt{3}}\zeta_0\zeta_1\pm \frac{1}{12\sqrt{2}}
\left(\bigl(\zeta_1\bigr)^2-1\right)\Biggr),
\end{equation}
\begin{equation}
\label{24.00000}
~~~~~~\hat I_{01}=\Delta^2\Biggl(-\frac{1}{3}\left(\bigl(\zeta_0\bigr)^2-1\right)
-\frac{1}{4\sqrt{3}}\zeta_0\zeta_1\mp \frac{1}{12\sqrt{2}}
\left(\bigl(\zeta_1\bigr)^2-1\right)\Biggr),
\end{equation}

\noindent
where $\zeta_0,$ $\zeta_1$ are the same random variables as in 
(\ref{16.00000})--(\ref{18.00000}). 

It is easy to check that the approximations 
(\ref{16.00000})--(\ref{18.00000}), 
(\ref{19.00000})--(\ref{24.00000})
satisfy the condition (\ref{010.026})
for $r=4$ and $m=1,$ $n\ge 1$, 
i.e. they are weak approximations
of the order $r=4$ for the case $m=1,$ $n\ge 1$.

\chapter{Approximation of Iterated Stochastic Integrals
with Respect to the $Q$-Wiener Process. Application
to the High-Order Strong Numerical Methods
for Non-Commutative Semilinear SPDEs with Nonliear Multiplicative
Trace Class Noise}

\section{Introduction}

There exists a lot of publications on the subject of 
numerical integration of stochastic partial differential
equations (SPDEs) (see, for example 
\cite{1zzzzz}-\cite{27zzzzz}).

One of the perspective approaches to the construction of high-order
strong numerical methods (with respect to the temporal 
discretization) for semilinear SPDEs is based
on the Taylor formula in Banach spaces and exponential formula for 
the mild solution of 
SPDEs \cite{8zzzzz}, \cite{10zzzzz}-\cite{13zzzzz}.
A significant step in this direction 
was made in \cite{12zzzzz} (2015), \cite{13zzzzz} (2016), 
where the exponential 
Milstein and 
Wagner--Platen methods for 
semilinear SPDEs with
nonlinear multiplicative trace class noise
were constructed.
Under the appropriate conditions \cite{12zzzzz}, \cite{13zzzzz} 
these methods
have strong orders of convergence $1.0-\varepsilon$ and $1.5-\varepsilon$
correspondingly
with respect to the temporal variable (where $\varepsilon$
is an arbitrary small posilive real number).
It should be noted that in \cite{18zzzzz} (2007) the convergence 
with strong
order $1.0$
of the 
exponential Milstein scheme 
for semilinear SPDEs 
was proved
under additional smoothness assumptions.

An important feature of the mentioned numerical methods is the presence 
in them the so-called iterated stochastic integrals with respect 
to the infinite-dimensional $Q$-Wiener process \cite{20zzzzz}.
Approximation of these stochastic integrals is a complex problem.
The problem of numerical modeling of these stochastic integrals
with multiplicities 1 to 3
was solved in \cite{12zzzzz}, \cite{13zzzzz} for the case when special 
commutativity conditions for semilinear SPDE with
nonlinear multiplicative trace class noise 
are fulfilled.

If the mentioned commutativity conditions are not fulfilled,
which often corresponds to SPDEs in numerous applications, the 
numerical modeling 
of iterated stochastic integrals with respect 
to the infinite-dimensional $Q$-Wiener process becomes 
much more difficult. 
Note that the exponential Milstein scheme \cite{12zzzzz} contains the
iterated stochastic integrals of multiplicities 1 and 2 with respect 
to the infinite-dimensional $Q$-Wiener process
and the exponential Wagner--Platen scheme \cite{13zzzzz}
contains the mentioned stochastic integrals of multiplicities
1 to 3 (see Sect.~7.2).

In \cite{26zzzzz} (2017), \cite{27zzzzz} (2018) two methods of 
the mean-square approximation of simplest
iterated (double) stochastic integrals from the 
exponential Milstein scheme for semilinear SPDEs with
nonlinear multiplicative trace class noise and
without the commutativity conditions
are considered and theorems on the convergence of 
these methods are given. 
At that, the basic idea (first of the mentioned methods \cite{26zzzzz}, 
\cite{27zzzzz}) about the Karhunen--Lo\`{e}ve
expansion of the Brownian bridge process was taken from the monograph
\cite{Zapad-1} (Milstein approach, see Sect.~6.2).
The second of the mentioned methods \cite{26zzzzz}, 
\cite{27zzzzz} is based on the results 
of Wiktorsson M. \cite{Zapad-6}, \cite{Zapad-7} (2001).

Note that the mean-square error of approximation
of iterated stochastic integrals with respect 
to the infinite-dimensional $Q$-Wiener process 
consists of two components \cite{26zzzzz}, \cite{27zzzzz}.
The first component is related with the finite-dimentional
approximation of the infinite-dimentional $Q$-Wiener process
while the second one is connected with the approximation
of iterated It\^{o} stochastic integrals with respect to the
scalar standard Brownian motions.

It is important to note that the approximation of iterated 
stochastic integrals with respect to the infinite-dimensional 
$Q$-Wiener process can be reduced to the approximation of iterated 
It\^{o} stochastic integrals with respect to 
the finite-dimensional 
Wiener process. In a lot of author's publications  
\cite{1}-\cite{new-new-6} (see Chapters 1, 2, and 5)
an effective method of the mean-square
approximation of 
iterated It\^{o} (and Stratonovich) stochastic 
integrals with respect to
the finite-dimensional Wiener 
process was proposed and developed. 
This method is based on the generalized multiple Fourier series,
in particular, on the multiple Fourier--Legendre series (see Sect.~5.1). 

The purpose of this chapter is an adaptation of the method 
\cite{1}-\cite{new-new-6}
for the mean-square approximation of iterated stochastic integrals 
with respect to the infinite-dimensional $Q$-Wiener process.
In the author's publications \cite{art-7}, \cite{arxiv-20} (see Sect.~7.3)
the problem of the mean-square approximation 
of iterated stochastic integrals 
with respect to the infinite-dimensional $Q$-Wiener process
in the sense of the second component of approximation error (see above)
has been solved 
for arbitraty multiplicity $k$ ($k\in{\bf N}$)
of stochastic integrals and without the assumptions of commutativity
for SPDE.
More precisely, in \cite{art-7}, \cite{arxiv-20} the method of 
generalized multiple Fourier series (Theorems 1.1, 1.2, 1.16) for 
the 
approximation of iterated It\^{o} stochastic integrals with respect to the 
scalar standard Brownian motions was adapted for  
iterated stochastic integrals with respect to the 
infinite-dimensional $Q$-Wiener process 
(in the sense of the second component of approximation error).

In Sect.~7.4 (also see \cite{OK}, \cite{arxiv-21}), 
we extend the method \cite{26zzzzz}, \cite{27zzzzz}
and estimate the first 
component of approximation error
for iterated stochastic integrals 
of multiplicities 1 to  3 with respect to the 
infinite-dimensional $Q$-Wiener process. 
In addition, we combine the obtained 
results with the results from \cite{art-7}, \cite{arxiv-20}
(see Sect.~7.3). Thus, the results 
of this chapter can be applied to the implementation of exponential Milstein
and Wagner--Platen schemes for semilinear SPDEs with
nonlinear multiplicative trace class noise and
without the commutativity conditions.

Let $U,$ $H$ be separable ${\bf R}$-Hilbert spaces and
$L_{HS}(U,H)$ be a space of Hilbert--Schmidt operators mapping 
from $U$ to $H$.
Let $(\Omega, {\bf F}, \sf{P})$ be a probability space 
with a normal filtration $\{{\bf F}_t, t\in [0, \bar{T}]\}$
\cite{20zzzzz}, let ${\bf W}_t$ be an $U$-valued $Q$-Wiener process 
with respect to $\{{\bf F}_t, t\in [0, \bar{T}]\},$  
which has a covariance trace class operator $Q\in L(U)$. 
Here and further $L(U)$ denotes all bounded linear operators
on $U$. Let $U_0$ be an ${\bf R}$-Hilbert space 
defined as $U_0=Q^{1/2}(U).$ At that, a scalar product in $U_0$ is
given by the relation \cite{13zzzzz}
$$
\left\langle u, w\right\rangle_{U_0}=\left\langle 
Q^{-1/2}u, Q^{-1/2}w\right\rangle_{U}
$$ 
for all $u, w\in U_0.$

Consider the semilinear parabolic 
SPDE with nonlinear multiplicative trace class noise
\begin{equation}
\label{xx1}
~~~~~~ dX_t = \left(A X_t + F(X_t)\right)dt + B(X_t)d{\bf W}_t,\ \ \
X_0=\xi,\ \ \ t\in [0, \bar{T}],
\end{equation}

\noindent
where nonlinear operators $F,$ $B$ ($F:$ $H\rightarrow H$, 
$B:$ $H\rightarrow L_{HS}(U_0,H)$), the linear operator
$A:$ $D(A)\subset H\rightarrow H$
as well as the initial value $\xi$ 
are assumed to satisfy the conditions
of existence and uniqueness of
the SPDE mild solution
(see \cite{13zzzzz}, Assumptions A1--A4).

It is well known \cite{23zzzzz} that Assumptions A1--A4 \cite{13zzzzz}
guarantee the existence and uniqueness (up to modifications) of 
the mild solution 
$X_t: [0, \bar T]\times \Omega \rightarrow H $
of the SPDE (\ref{xx1}) 
\begin{equation}
\label{mild}
X_t={\rm exp}(At)\xi +\int\limits_0^t {\rm exp}(A(t-\tau))F(X_\tau)d\tau+
\int\limits_0^t {\rm exp}(A(t-\tau))B(X_\tau)d{\bf W}_{\tau}
\end{equation}

\noindent
w.~p.~1 for all $t\in[0, \bar T],$ where
${\rm exp}(At),$ $t\ge 0$ is the semigroup generated by the operator $A$.

As we mentioned earlier, numerical methods of high orders of accuracy 
(with respect to the temporal discretization)
for approximating the mild solution of the SPDE (\ref{xx1}),
which are based on the Taylor formula for operators and an 
exponential formula for the mild solution of 
SPDEs, contain iterated stochastic integrals with respect to the $Q$-Wiener 
process \cite{8zzzzz}, \cite{10zzzzz}-\cite{13zzzzz}, \cite{18zzzzz}.

Note that the exponential Milstein type numerical scheme 
\cite{12zzzzz}
and the exponential Wagner-Platen type numerical scheme \cite{13zzzzz}
contain, for example, the following 
iterated stochastic integrals (see Sect.~7.2)
\begin{equation}
\label{xx2}
\int\limits_{t}^{T}B(Z) d{\bf W}_{t_1},\ \ \
\int\limits_{t}^{T}B'(Z) \left(
\int\limits_{t}^{t_2}B(Z) d{\bf W}_{t_1} \right) d{\bf W}_{t_2},\
\end{equation}
\begin{equation}
\label{xx3}
~~~ \int\limits_{t}^{T}B'(Z) \left(
\int\limits_{t}^{t_2}F(Z) dt_1 \right) d{\bf W}_{t_2},\ \ \
\int\limits_{t}^{T}
F'(Z)\left(\int\limits_{t}^{t_2} B(Z) d{\bf W}_{t_1} \right) dt_2,\ 
\end{equation}
\begin{equation}
\label{xx3a}
\int\limits_{t}^{T}B'(Z) \left(
\int\limits_{t}^{t_3}B'(Z) \left(
\int\limits_{t}^{t_2}B(Z)
d{\bf W}_{t_1} \right) d{\bf W}_{t_2} \right) d{\bf W}_{t_3},\
\end{equation}
\begin{equation}
\label{xx3aaa}
\int\limits_{t}^{T}B''(Z) \left(
\int\limits_{t}^{t_2}B(Z) d{\bf W}_{t_1},
\int\limits_{t}^{t_2}B(Z) d{\bf W}_{t_1} 
\right) d{\bf W}_{t_2},\
\end{equation}

\noindent
where $0\le t<T \le \bar{T},$
$Z: \Omega \rightarrow H$ is an ${\bf F}_t/{\cal B}(H)$-measurable mapping
and $F',$ $B',$ $B''$ denote 
Fr\^{e}chet derivatives.
At that, the 
exponential Milstein type scheme \cite{12zzzzz} contains integrals 
(\ref{xx2}) while the exponential Wagner--Platen 
type scheme \cite{13zzzzz} 
contains integrals (\ref{xx2})--(\ref{xx3aaa}) (see Sect.~7.2).

It is easy to notice that the numerical schemes for SPDEs with higher 
orders of convergence (with respect to the temporal discretization) 
in contrast with the numerical schemes from 
\cite{12zzzzz}, \cite{13zzzzz} will include iterated stochastic
integrals (with respect to the $Q$-Wiener process) 
with multiplicities $k>3$ \cite{11zzzzz} (2011). 
So, this chapter is partially devoted to the approximation of iterated 
stochastic integrals of the form
\begin{equation}
\label{xx4}
I[\Phi^{(k)}(Z)]_{T,t}=\int\limits_{t}^{T}\Phi_k(Z) \left( \ldots \left(
\int\limits_{t}^{t_3}\Phi_2(Z) \left(
\int\limits_{t}^{t_2}\Phi_1(Z)
d{\bf W}_{t_1} \right)
d{\bf W}_{t_2} \right) \ldots  \right) d{\bf W}_{t_k},
\end{equation}

\vspace{2mm}
\noindent
where $0\le t<T \le \bar{T},$
$Z: \Omega \rightarrow H$ is an ${\bf F}_t/{\cal B}(H)$-measurable 
mapping
and 
$\Phi_k(v)(\ \ldots (\Phi_2(v)(\Phi_1(v)) \ldots\ ))$
is a $k$-linear Hilbert--Schmidt operator
mapping from
$\underbrace{U_0\times \ldots \times U_0}_{\small{\hbox{$k$ times}}}$ 
to $H$
for all $v\in H$.

\vspace{1mm}

In Sect.~7.3.1 we consider the approximation of more general
iterated stochastic integrals than (\ref{xx4}). 
In Sect.~7.3.2 and 7.3.3 some other
types of ite\-ra\-ted stochastic integrals of multiplicities
2--4 with respect to the $Q$-Wiener process will be considered.

Note that the stochastic integral (\ref{xx3aaa}) is not a special case 
of the stochastic integral (\ref{xx4}) for $k=3$.
Nevertheless, the extended representation for approximation of the 
stochastic 
integral (\ref{xx3aaa}) is similar
to (\ref{xx405}) (see below) for $k=3$. Moreover,
the mentioned
representation for approximation of the stochastic integral (\ref{xx3aaa})
contains the same iterated It\^{o} stochastic integrals 
of third multiplicity as in (\ref{xx405}) for $k=3$
(see Sect.~7.3.2).
These conclusions mean that one of the main results
of this chapter (Theorem 7.1, Sect.~7.3.1) 
for $k=3$ can be reformulated naturally for the stochastic integral 
(\ref{xx3aaa}) (see Sect.~7.3.2).

It should be noted that by developing the 
approach from the work \cite{13zzzzz}, 
which uses the Taylor formula for operators and a formula for 
the mild solution of 
the SPDE (\ref{xx1}),
we obviously obtain a number of other iterated stochastic 
integrals. For example, the following stochastic integrals
$$
\int\limits_{t}^{T}B'''(Z) \left(
\int\limits_{t}^{t_2}B(Z) d{\bf W}_{t_1},
\int\limits_{t}^{t_2}B(Z) d{\bf W}_{t_1},
\int\limits_{t}^{t_2}B(Z) d{\bf W}_{t_1}
\right) d{\bf W}_{t_2},\
$$
$$
\int\limits_{t}^{T}B'(Z) \left(
\int\limits_{t}^{t_3}B''(Z) \left(
\int\limits_{t}^{t_2}B(Z) d{\bf W}_{t_1},
\int\limits_{t}^{t_2}B(Z) d{\bf W}_{t_1} 
\right) d{\bf W}_{t_2}\right)d{\bf W}_{t_3},
$$
$$
\int\limits_{t}^{T}B''(Z) \left(
\int\limits_{t}^{t_3}B(Z)d{\bf W}_{t_1},
\int\limits_{t}^{t_3}B'(Z)\left(
\int\limits_{t}^{t_2}B(Z) d{\bf W}_{t_1}
\right) d{\bf W}_{t_2}\right)d{\bf W}_{t_3},\
$$
$$
\int\limits_{t}^{T}F'(Z)\left(
\int\limits_{t}^{t_3}B'(Z) \left(
\int\limits_{t}^{t_2}B(Z) d{\bf W}_{t_1}
\right) d{\bf W}_{t_2} 
\right) dt_3,
$$
$$
\int\limits_{t}^{T}F''(Z) \left(
\int\limits_{t}^{t_2}B(Z) d{\bf W}_{t_1},
\int\limits_{t}^{t_2}B(Z) d{\bf W}_{t_1} 
\right) dt_2,
$$
$$
\int\limits_{t}^{T}B''(Z) \left(
\int\limits_{t}^{t_2}F(Z) dt_1,
\int\limits_{t}^{t_2}B(Z) d{\bf W}_{t_1} 
\right) d{\bf W}_{t_2}
$$
will be considered in Sect.~7.3.3. 
Here 
$Z: \Omega \rightarrow H$ is an ${\bf F}_t/{\mathcal{B}}(H)$-measurable 
mapping and 
$B'$, $B'',$ $B''',$ $F'$, $F''$ are
Fr\^{e}chet derivatives.

Consider eigenvalues $\lambda_i$ and 
eigenfunctions $e_i(x)$ of the covariance operator $Q,$ where
$i=(i_1,\ldots,i_d)\in J,$\ 
$x=(x_1,\ldots,x_d),$\ and\  
$J=\bigl\{i:\ i\in {\bf N}^d\ \hbox{and}\ \lambda_i>0\bigr\}$.

The series representation of 
the $Q$-Wiener process has the form \cite{20zzzzz}
$$
{\bf W}(t,x)=\sum\limits_{i\in J} e_i(x)\sqrt{\lambda_i}
{\bf w}_t^{(i)},\ \ \ t\in[0, \bar{T}]
$$
or in the shorter notations
$$
{\bf W}_t=\sum\limits_{i\in J} e_i\sqrt{\lambda_i}
{\bf w}_t^{(i)},\ \ \ t\in[0, \bar{T}],
$$
where ${\bf w}_t^{(i)},$ $i\in J$ are independent standard
Wiener processes.

Note that eigenfunctions $e_i,$ $i\in J$ form
an orthonormal basis of $U$ \cite{20zzzzz}.

Consider the finite-dimensional approximation of ${\bf W}_t$ \cite{20zzzzz}
\begin{equation}
\label{yy.1}
{\bf W}^M_t=\sum\limits_{i\in J_M} e_i\sqrt{\lambda_i}
{\bf w}_t^{(i)},\ \ \ t\in[0, \bar{T}],
\end{equation}
where 
\begin{equation}
\label{sat1}
J_M=\bigl\{i:\ 1\le i_1,\ldots,i_d\le M\ \hbox{and}\ 
\lambda_i>0\bigr\}.
\end{equation}

Using (\ref{yy.1}) and the relation \cite{20zzzzz}
\begin{equation}
\label{xx.1}
{\bf w}_t^{(i)}=\frac{1}{\sqrt{\lambda_i}}\langle e_i,{\bf W}_t\rangle_U,\ \ \
i\in J,
\end{equation}
we obtain

\vspace{-5mm}
\begin{equation}
\label{yy.12}
{\bf W}^M_t=\sum\limits_{i\in J_M} 
e_i\ \langle e_i,{\bf W}_t\rangle_U,\ \ \ t\in[0, \bar{T}],
\end{equation}

\vspace{2mm}
\noindent
where $\langle \cdot,\cdot\rangle_U$ is a scalar product in $U.$

Taking into account (\ref{xx.1}) and (\ref{yy.12}), we note that
the approximation $I[\Phi^{(k)}(Z)]_{T,t}^M$ of the ite\-ra\-ted stochastic 
integral
$I[\Phi^{(k)}(Z)]_{T,t}$ (see (\ref{xx4}))
can be written w.~p.~1 in the following form

\vspace{-2mm}
$$
I[\Phi^{(k)}(Z)]_{T,t}^M=
$$

\vspace{-4mm}
$$
=\int\limits_{t}^{T}\Phi_k(Z) \left( \ldots \left(
\int\limits_{t}^{t_3}\Phi_2(Z) \left(
\int\limits_{t}^{t_2}\Phi_1(Z)
d{\bf W}_{t_1}^M\ \right) 
d{\bf W}_{t_2}^M\ \right) \ldots  \right) d{\bf W}_{t_k}^M=
$$

\vspace{3mm}
$$
=
\sum_{r_1,\ldots,r_k\in J_M}
\Phi_k(Z)\left(\ldots \left(\Phi_2(Z) \left(\Phi_1(Z)
e_{r_1} \right)
e_{r_2} \right) \ldots \right) e_{r_k} \times
$$
$$
\times
\int\limits_t^T \ldots \int\limits_t^{t_{3}}\int\limits_t^{t_{2}}
d\langle e_{r_1},{\bf W}_{t_1}
\rangle_U\
d\langle e_{r_2},{\bf W}_{t_2}
\rangle_U\ \ldots \
d\langle e_{r_k},{\bf W}_{t_k}
\rangle_U =
$$

\vspace{2mm}
$$
=\sum_{r_1,\ldots,r_k\in J_M}
\Phi_k(Z)\left(\ldots \left(\Phi_2(Z) \left(\Phi_1(Z)
e_{r_1} \right) 
e_{r_2} \right) \ldots \right) e_{r_k}
\sqrt{\lambda_{r_1}\lambda_{r_2}\ldots \lambda_{r_k}}\times
$$
\begin{equation}
\label{xx405}
\times
\int\limits_t^T \ldots \int\limits_t^{t_3}\int\limits_t^{t_{2}}
d{\bf w}_{t_1}^{(r_1)}d{\bf w}_{t_2}^{(r_2)}\ldots
d{\bf w}_{t_k}^{(r_k)},
\end{equation}

\noindent
where $0\le t<T \le \bar{T}.$

{\bf Remark 7.1.}\ {\it Obviously, without loss of generality,
we can write $J_M=\{1, 2,\ldots, M\}.$}

As we mentioned before, when special conditions of 
commutativity for the SPDE (\ref{xx1}) 
be fulfilled, it is 
proposed to simulate 
numerically
the stochastic integrals (\ref{xx2})--(\ref{xx3aaa}) using 
the simple formulas \cite{12zzzzz}, \cite{13zzzzz}. 
In this case, the numerical simulation of the mentioned stochastic 
integrals requires the use of increments of the $Q$-Wiener process only.
However, if these commutativity conditions are not fulfilled
(which often corresponds to SPDEs in numerous applications), the 
numerical simulation of the stochastic integrals (\ref{xx2})--(\ref{xx3aaa})
becomes 
much more difficult. 
Recall that in \cite{26zzzzz}, \cite{27zzzzz}
two methods for the mean-square approximation of simplest
iterated (double) 
stochastic integrals defined by (\ref{xx2})
are proposed. In this chapter, we consider 
a substantially more general and effective method (based on 
the results of Chapters 1 and 5)
for the mean-square approximation of 
iterated stochastic integrals of multiplicity $k$ $(k\in{\bf N})$
with respect to 
the $Q$-Wiener process. 
The convergence analysis in the transition from $J_M$ to $J$, i.e., 
from the finite-dimensional Wiener process to the infinite-dimensional 
one will be carried out in Sect.~7.4 for integrals of multiplicities 1 to 3
similar to the 
proof of Theorem 1 \cite{27zzzzz}.

\section{Exponential Milstein and Wagner--Platen Numerical Schemes
for Non-Com\-mu\-ta\-tive Semilinear SPDEs}

Let assumptions of Sect.~7.1 are fulfilled.
Let $\Delta>0,$ $\tau_p=p\Delta$ $(p=0, 1,\ldots, N),$ 
and $N\Delta=\bar T.$ Consider the 
exponential Milstein numerical scheme \cite{12zzzzz} 
$$
Y_{p+1}={\rm exp}\left(A\Delta\right)\left(
Y_p+\Delta F(Y_p)+\int\limits_{\tau_p}^{\tau_{p+1}}B(Y_p)d{\bf W}_s
+\right.
$$
\begin{equation}
\label{fff100}
\left.+\int\limits_{\tau_p}^{\tau_{p+1}}B'(Y_p)
\left(\int\limits_{\tau_p}^{s}B(Y_p)d{\bf W}_{\tau}\right)
d{\bf W}_s\right)
\end{equation}

\vspace{4mm}
\noindent
and the exponential Wagner--Platen numerical scheme \cite{13zzzzz} 

\vspace{-2mm}
$$
Y_{p+1}={\rm exp}\left(\frac{A\Delta}{2}\right)\left(
{\rm exp}\left(\frac{A\Delta}{2}\right)Y_p+\Delta F(Y_p)+
\int\limits_{\tau_p}^{\tau_{p+1}}B(Y_p)d{\bf W}_s+\right.
$$
$$
+
\int\limits_{\tau_p}^{\tau_{p+1}}B'(Y_p)
\left(\int\limits_{\tau_p}^{s}B(Y_p)d{\bf W}_{\tau}\right)
d{\bf W}_s+
$$
$$
+
\frac{\Delta^2}{2}
F'(Y_p)\biggl(AY_p+F(Y_p)\biggr)+
\int\limits_{\tau_p}^{\tau_{p+1}}
F'(Y_p)\left(\int\limits_{\tau_p}^s
B(Y_p)d{\bf W}_{\tau}\right)ds+
$$
$$
+\frac{\Delta^2}{4}
\sum\limits_{i\in J}\lambda_i F''(Y_p)\biggl(B(Y_p)e_i,B(Y_p)e_i\biggr)+
$$
$$
+
A\left(\int\limits_{\tau_p}^{\tau_{p+1}}
\int\limits_{\tau_p}^{s}
B(Y_p)d{\bf W}_{\tau}ds-\frac{\Delta}{2}
\int\limits_{\tau_p}^{\tau_{p+1}}B(Y_p)d{\bf W}_s\right)+
$$
$$
+\Delta
\int\limits_{\tau_p}^{\tau_{p+1}}B'(Y_p)
\biggl(AY_p+F(Y_p)\biggr)d{\bf W}_s-
\int\limits_{\tau_p}^{\tau_{p+1}}\int\limits_{\tau_p}^{s}B'(Y_p)
\biggl(AY_p+F(Y_p)\biggr)d{\bf W}_{\tau}ds+
$$
$$
+
\frac{1}{2}\int\limits_{\tau_p}^{\tau_{p+1}}B''(Y_p)
\left(\int\limits_{\tau_p}^{s}B(Y_p)d{\bf W}_{\tau},
\int\limits_{\tau_p}^{s}B(Y_p)d{\bf W}_{\tau}
\right)
d{\bf W}_s+
$$
\begin{equation}
\label{fffguk}
~~~~~~~\left.+
\int\limits_{\tau_p}^{\tau_{p+1}}B'(Y_p)
\left(\int\limits_{\tau_p}^{s}B'(Y_p)
\left(\int\limits_{\tau_p}^{\tau}B(Y_p)
d{\bf W}_{\theta}\right)
d{\bf W}_{\tau}
\right)
d{\bf W}_s\right)
\end{equation}

\vspace{5mm}

\noindent 
for the SPDE (\ref{xx1}),
where $Y_p$ is an approximation of $X_{\tau_p}$ (mild solution 
(\ref{mild})
at the time moment  $\tau_p$), $p=0, 1,\ldots,N,$ and
$B'$, $B'',$ $F'$, $F''$ are
Fr\^{e}chet derivatives. 

Note that in addition to the temporal discretization, 
the implementation of numerical schemes 
(\ref{fff100}) and (\ref{fffguk}) also requires a 
discretization of the infinite-dimensional
Hilbert space $H$ (approximation with respect to the space domain)
and a finite-dimensional approximation of the $Q$-Wiener process.
Let us focus on the 
approximation connected with the $Q$-Wiener process.

Consider the following iterated It\^{o} stochastic integrals
\begin{equation}
\label{x111}
J_{(1)T,t}^{(r_1)}=\int\limits_t^Td{\bf w}_{t_1}^{(r_1)},\ \ \
J_{(10)T,t}^{(r_1 0)}=
\int\limits_t^T\int\limits_t^{t_2} d{\bf w}_{t_1}^{(r_1)}
dt_2,\ \ \
J_{(01)T,t}^{(0 r_2)}=\int\limits_t^T\int\limits_t^{t_2} dt_1
d{\bf w}_{t_2}^{(r_2)},
\end{equation}
\begin{equation}
\label{x112}
~~~J_{(11)T,t}^{(r_1 r_2)}=
\int\limits_t^T\int\limits_t^{t_2} d{\bf w}_{t_1}^{(r_1)}
d{\bf w}_{t_2}^{(r_2)},\ \ \
J_{(111)T,t}^{(r_1 r_2 r_3)}=
\int\limits_t^T\int\limits_t^{t_3}\int\limits_t^{t_2} 
d{\bf w}_{t_1}^{(r_1)}
d{\bf w}_{t_2}^{(r_2)}
d{\bf w}_{t_3}^{(r_3)},
\end{equation}

\noindent
where $r_1, r_2, r_3\in J_M,$ $0\le t<T\le \bar T,$
and $J_M$ is defined by (\ref{sat1}).

Let us replace the infinite-dimensional
$Q$-Wiener process in the iterated stochastic 
integrals from (\ref{fff100}), 
(\ref{fffguk}) by its finite-dimensional approximation (\ref{yy.1}).
Then we have w.~p.~1 
\begin{equation}
\label{ggg0}
\int\limits_{\tau_p}^{\tau_{p+1}}B(Y_p)d{\bf W}_s^M=
\sum\limits_{r_1\in J_M}
B(Y_p)e_{r_1}\sqrt{\lambda_{r_1}}
J_{(1)\tau_{p+1},\tau_p}^{(r_1)},
\end{equation}

\vspace{3mm}
$$
A\left(\int\limits_{\tau_p}^{\tau_{p+1}}
\int\limits_{\tau_p}^{s}
B(Y_p)d{\bf W}_{\tau}^M ds-\frac{\Delta}{2}
\int\limits_{\tau_p}^{\tau_{p+1}}B(Y_p)d{\bf W}_s^M\right)=
$$
$$
=
A\int\limits_{\tau_p}^{\tau_{p+1}}
B(Y_p)\left(\frac{\tau_{p+1}}{2}-s+\frac{\tau_p}{2}\right)
d{\bf W}_{s}^M=
$$
\begin{equation}
\label{ggg1}
=\sum\limits_{r_1\in J_M}
AB(Y_p)e_{r_1}\sqrt{\lambda_{r_1}}
\left(\frac{\Delta}{2}J_{(1)\tau_{p+1},\tau_p}^{(r_1)}-
J_{(01)\tau_{p+1},\tau_p}^{(0 r_1)}\right),
\end{equation}

\vspace{3mm}
$$
\Delta
\int\limits_{\tau_p}^{\tau_{p+1}}B'(Y_p)
\biggl(AY_p+F(Y_p)\biggr)d{\bf W}_s^M-
\int\limits_{\tau_p}^{\tau_{p+1}}\int\limits_{\tau_p}^{s}B'(Y_p)
\biggl(AY_p+F(Y_p)\biggr)d{\bf W}_{\tau}^Mds=
$$
$$
=\int\limits_{\tau_p}^{\tau_{p+1}}B'(Y_p)
\int\limits_{\tau_p}^{s}
\biggl(AY_p+F(Y_p)\biggr)d\tau d{\bf W}_{s}^M=
$$

\begin{equation}
\label{ggg2}
=
\sum\limits_{r_1\in J_M}
B'(Y_p)\biggl(AY_p+F(Y_p)\biggr)
e_{r_1}\sqrt{\lambda_{r_1}}J_{(01)\tau_{p+1},\tau_p}^{(0 r_1)},
\end{equation}

\vspace{6mm}

$$
\int\limits_{\tau_p}^{\tau_{p+1}}
F'(Y_p)\left(\int\limits_{\tau_p}^s
B(Y_p)d{\bf W}_{\tau}^M\right)ds=
$$

\begin{equation}
\label{ggg3}
=\sum\limits_{r_1\in J_M}
F'(Y_p)B(Y_p)e_{r_1}\sqrt{\lambda_{r_1}}
\left(\Delta J_{(1)\tau_{p+1},\tau_p}^{(r_1)}-
J_{(01)\tau_{p+1},\tau_p}^{(0 r_1)}\right),
\end{equation}

\newpage
\noindent
$$
\int\limits_{\tau_p}^{\tau_{p+1}}
B'(Y_p)\left(\int\limits_{\tau_p}^s
B(Y_p)d{\bf W}_{\tau}^M\right)d{\bf W}_s^M=
$$

\begin{equation}
\label{uspel1}
=\sum\limits_{r_1,r_2\in J_M}
B'(Y_p)\left(B(Y_p)e_{r_1}\right)
e_{r_2}\sqrt{\lambda_{r_1}\lambda_{r_2}}J_{(11)\tau_{p+1},\tau_p}^{(r_1 r_2)},
\end{equation}

\vspace{5mm}
$$
\int\limits_{\tau_p}^{\tau_{p+1}}
B'(Y_p)\left(\int\limits_{\tau_p}^s
B'(Y_p)
\left(\int\limits_{\tau_p}^{\tau}
B(Y_p)
d{\bf W}_{\theta}^M\right)d{\bf W}_{\tau}^M\right)d{\bf W}_{s}^M=
$$

\begin{equation}
\label{uspel2}
~~~~~ =\sum\limits_{r_1,r_2,r_3\in J_M}
B'(Y_p)\left(B'(Y_p)\left(B(Y_p)e_{r_1}\right)
e_{r_2}\right)e_{r_3}
\sqrt{\lambda_{r_1}\lambda_{r_2}\lambda_{r_3}}
J_{(111)\tau_{p+1},\tau_p}^{(r_1 r_2 r_3)},
\end{equation}

\vspace{5mm}
$$
\int\limits_{\tau_p}^{\tau_{p+1}}B''(Y_p)
\left(\int\limits_{\tau_p}^{s}B(Y_p)d{\bf W}_{\tau}^M,
\int\limits_{\tau_p}^{s}B(Y_p)d{\bf W}_{\tau}^M
\right)
d{\bf W}_s^M=
$$

$$
=
\sum_{r_1,r_2,r_3\in J_M}B''(Y_p)
\left(B(Y_p)e_{r_1}, B(Y_p)e_{r_2}\right)e_{r_3}
\sqrt{\lambda_{r_1}\lambda_{r_2}\lambda_{r_3}}\times
$$

\begin{equation}
\label{ggg4}
\times
\int\limits_{\tau_p}^{\tau_{p+1}}
\left(\int\limits_{\tau_p}^s d{\bf w}_{\tau}^{(r_1)}   
\int\limits_{\tau_p}^s d{\bf w}_{\tau}^{(r_2)}\right)
d{\bf w}_s^{(r_3)}.
\end{equation}

\vspace{5mm}

Note that in (\ref{ggg1})--(\ref{ggg3}) we used the It\^{o} formula.
Moreover, using the It\^{o} formula we obtain
\begin{equation}
\label{da12ccc}
~~~~~\int\limits_{\tau_p}^s d{\bf w}_{\tau}^{(r_1)}   
\int\limits_{\tau_p}^s d{\bf w}_{\tau}^{(r_2)}
=J_{(11)s,\tau_p}^{(r_1 r_2)}+J_{(11)s,\tau_p}^{(r_2 r_1)}
+
{\bf 1}_{\{r_1=r_2\}}(s-\tau_p)\ \ \ \hbox{w.~p.~1,}
\end{equation}

\noindent
where ${\bf 1}_A$ is the indicator of the set $A.$
From (\ref{da12ccc}) we have w.~p.~1 
\begin{equation}
\label{jjj2}
\int\limits_{\tau_p}^{\tau_{p+1}}
\left(\int\limits_{\tau_p}^s d{\bf w}_{\tau}^{(r_1)}   
\int\limits_{\tau_p}^s d{\bf w}_{\tau}^{(r_2)}\right)
d{\bf w}_s^{(r_3)}=J_{(111)\tau_{p+1},\tau_p}^{(r_1 r_2 r_3)}
+J_{(111)\tau_{p+1},\tau_p}^{(r_2 r_1 r_3)}+
{\bf 1}_{\{r_1=r_2\}}J_{(01)\tau_{p+1},\tau_p}^{(0 r_3)}.
\end{equation}

After substituting (\ref{jjj2}) into (\ref{ggg4}), we obtain w.~p.~1 
$$
\int\limits_{\tau_p}^{\tau_{p+1}}B''(Y_p)
\left(\int\limits_{\tau_p}^{s}B(Y_p)d{\bf W}_{\tau}^M,
\int\limits_{\tau_p}^{s}B(Y_p)d{\bf W}_{\tau}^M
\right)
d{\bf W}_s^M=
$$
$$
=
\sum_{r_1,r_2,r_3\in J_M}B''(Y_p)\left(B(Y_p)
e_{r_1}, B(Y_p)e_{r_2}\right)e_{r_3}
\sqrt{\lambda_{r_1}\lambda_{r_2}\lambda_{r_3}}\times
$$

\begin{equation}
\label{da1quka}
\times \left(J_{(111)\tau_{p+1},\tau_p}^{(r_1 r_2 r_3)}
+J_{(111)\tau_{p+1},\tau_p}^{(r_2 r_1 r_3)}+
{\bf 1}_{\{r_1=r_2\}}
J_{(01)\tau_{p+1},\tau_p}^{(0 r_3)}\right).
\end{equation}

\vspace{4mm}

Thus, for the implementation of numerical schemes 
(\ref{fff100}) and (\ref{fffguk})
we need to approximate the following collection of iterated 
It\^{o} stochastic integrals
\begin{equation}
\label{ha100}
J_{(1)T,t}^{(r_1)},\ \ \
J_{(01)T,t}^{(0 r_1)},\ \ \ 
J_{(11)T,t}^{(r_1 r_2)},\ \ \
J_{(111)T,t}^{(r_1 r_2 r_3)},
\end{equation}
where $r_1, r_2, r_3\in J_M,$\ \  $0\le t<T\le \bar T.$

The problem of the mean-square approximation of 
iterated It\^{o} stochastic integrals (\ref{ha100}) is considered completely 
in Chapters 1 and 5.

\section{Approximation of Iterated Stochastic Integrals
of Multiplicity $k$ ($k\in{\bf N}$)
with Respect to the Finite-Dimensional Approximation ${\bf W}_t^M$ of
the $Q$-Wiener Process}

In this section, we consider a method for the approximation 
of iterated stochastic 
integrals of multiplicity $k$ $(k\in {\bf N})$
with respect to the 
finite-dimensional approximation ${\bf W}_t^M$ of the $Q$-Wiener process 
${\bf W}_t$
using the mean-square
approximation method of iterated It\^{o} stochastic integrals based on
Theorems 1.1, 1.2, 1.16.

\subsection{Theorem on the Mean-Square Approximation
of Iterated Sto\-chas\-tic Integrals
of Multiplicity $k$ ($k\in{\bf N}$)
with Respect to the Finite-Dimensional Approximation ${\bf W}_t^M$ of
the $Q$-Wiener Process}

Consider the iterated stochastic integral with respect to 
the $Q$-Wiener process in the following form
\newpage
\noindent
$$
I[\Phi^{(k)}(Z), \psi^{(k)}]_{T,t}
=\int\limits_{t}^{T}\Phi_k(Z) \left( \ldots \left(
\int\limits_{t}^{t_3}\Phi_2(Z)\times\right.\right.
$$

\vspace{-2mm}
\begin{equation}
\label{ssss11}
~~~~~\left.\left.\times \left(
\int\limits_{t}^{t_2}\Phi_1(Z)
\psi_1(t_1)d{\bf W}_{t_1} \right)
\psi_2(t_2)d{\bf W}_{t_2} \right) 
\ldots\right) 
\psi_k(t_k)d{\bf W}_{t_k},
\end{equation}

\vspace{3mm}
\noindent
where $Z: \Omega \rightarrow H$ is 
an ${\bf F}_t/{\mathcal{B}}(H)$-measurable mapping, 
${\bf W}_{\tau}$ is the $Q$-Wiener process,
$\Phi_k(v)(\ \ldots (\Phi_2(v)(\Phi_1(v))) \ldots\ )$ 
is a $k$-linear Hilbert--Schmidt operator
mapping from
$\underbrace{U_0\times \ldots \times U_0}_{\small{\hbox{$k$ times}}}$ to $H$
for all $v\in H,$ and $\psi_1(\tau),\ldots,\psi_k(\tau)\in L_2([t, T])$.

\vspace{4mm}

Let $I[\Phi^{(k)}(Z), \psi^{(k)}]_{T,t}^M$ be the approximation
of the iterated stochastic integral (\ref{ssss11})

\vspace{-3mm}
$$
I[\Phi^{(k)}(Z), \psi^{(k)}]_{T,t}^M
=\int\limits_{t}^{T}\Phi_k(Z) \left( \ldots \left(
\int\limits_{t}^{t_3}\Phi_2(Z)\times\right.\right.
$$

$$
\left.\left.\times \left(
\int\limits_{t}^{t_2}\Phi_1(Z)
\psi_1(t_1)d{\bf W}_{t_1}^M \right)
\psi_2(t_2)d{\bf W}_{t_2}^M \right) 
\ldots\right) 
\psi_k(t_k)d{\bf W}_{t_k}^M=
$$

\vspace{-3mm}
$$
=\sum_{r_1,r_2,\ldots,r_k\in J_M}
\Phi_k(Z)\left(\ldots \left(\Phi_2(Z) \left(\Phi_1(Z)
e_{r_1} \right) 
e_{r_2} \right) \ldots \right) e_{r_k}
\left(\prod_{l=1}^k\lambda_{r_l}\right)^{1/2}\times
$$

\begin{equation}
\label{xx605}
\times
J[\psi^{(k)}]_{T,t}^{(r_1 r_2\ldots r_k)},
\end{equation}

\vspace{6mm}
\noindent
where $0\le t<T\le \bar{T}$  and

\vspace{-3mm}
$$
J[\psi^{(k)}]_{T,t}^{(r_1 r_2\ldots r_k)}=
\int\limits_t^T \psi_k(t_k)\ldots \int\limits_t^{t_3}
\psi_2(t_2)\int\limits_t^{t_{2}}\psi_1(t_1)
d{\bf w}_{t_1}^{(r_1)}d{\bf w}_{t_2}^{(r_2)}\ldots
d{\bf w}_{t_k}^{(r_k)}
$$

\vspace{2mm}
\noindent
is the iterated It\^{o} stochastic integral (\ref{ito}).

Let $I[\Phi^{(k)}(Z), \psi^{(k)}]_{T,t}^{M,p_1,\ldots, p_k}$ be the
approximation of the iterated stochastic integral (\ref{xx605})

\newpage
\noindent
$$
I[\Phi^{(k)}(Z), \psi^{(k)}]_{T,t}^{M,p_1,\ldots, p_k}
=
$$
$$
=\sum_{r_1,r_2,\ldots,r_k\in J_M}
\Phi_k(Z)\left(\ldots \left(\Phi_2(Z) \left(\Phi_1(Z)
e_{r_1} \right) 
e_{r_2} \right) \ldots \right) e_{r_k}
\left(\prod_{l=1}^k\lambda_{r_l}\right)^{\hspace{-2.2mm}1/2}\times
$$

\vspace{-1mm}
\begin{equation}
\label{xx705}
\times J[\psi^{(k)}]_{T,t}^{(r_1 r_2\ldots r_k)p_1,\ldots ,p_k},
\end{equation}

\vspace{4mm}
\noindent
where
$J[\psi^{(k)}]_{T,t}^{(r_1 r_2\ldots r_k)p_1, \ldots, p_k}$
is defined as the expression before passing to the limit
on the right-hand side of (\ref{razzar1})

\vspace{-2mm}
$$
J[\psi^{(k)}]_{T,t}^{(r_1 r_2\ldots r_k)p_1,\ldots ,p_k}=
\sum\limits_{j_1=0}^{p_1}\ldots
\sum\limits_{j_k=0}^{p_k}
C_{j_k\ldots j_1}\Biggl(
\prod_{l=1}^k\zeta_{j_l}^{(r_l)}+\sum\limits_{m=1}^{[k/2]}
(-1)^m \times
\Biggr.
$$

\vspace{-2mm}
\begin{equation}
\label{f1}
\times
\sum_{\stackrel{(\{\{g_1, g_2\}, \ldots, 
\{g_{2m-1}, g_{2m}\}\}, \{q_1, \ldots, q_{k-2m}\})}
{{}_{\{g_1, g_2, \ldots, 
g_{2m-1}, g_{2m}, q_1, \ldots, q_{k-2m}\}=\{1, 2, \ldots, k\}}}}
\prod\limits_{s=1}^m
{\bf 1}_{\{r_{g_{{}_{2s-1}}}=~r_{g_{{}_{2s}}}\ne 0\}}
\Biggl.{\bf 1}_{\{j_{g_{{}_{2s-1}}}=~j_{g_{{}_{2s}}}\}}
\prod_{l=1}^{k-2m}\zeta_{j_{q_l}}^{(r_{q_l})}\Biggr).
\end{equation}

\vspace{4mm}

Let $U,$ $H$ be separable
${\bf R}$-Hilbert spaces,
$U_{0}=Q^{1/2}(U)$, and
$L(U, H)$ be the space of linear and bounded operators mapping
from $U$ to $H$. Let

\vspace{-4mm}
$$
L(U, H)_{0}=\bigl\{T \vert_{U_{0}}:\
T\in L(U, H)\bigr\},
$$

\noindent
where $T \vert_{U_{0}}$ is the restriction of operator $T$ to
the space $U_0$. It is known \cite{20zzzzz} that
$L(U, H)_{0}$
is a dense subset of the space of
Hilbert--Schmidt operators $L_{HS}(U_{0}, H)$.

\vspace{2mm}

{\bf Theorem 7.1} \cite{12a}-\cite{12aa}, \cite{art-7}, \cite{arxiv-20}. {\it 
Suppose that
$\psi_1(\tau),\ldots,\psi_k(\tau)\in L_2([t, T])$ and
$\{\phi_j(x)\}_{j=0}^{\infty}$ is an arbitrary complete orthonormal system  
of functions in the space $L_2([t,T]).$
Furthermore$,$ let the following conditions be satisfied$:$

{\rm 1}. $Q\in L(U)$ is a nonnegative and symmetric 
trace class operator {\rm (}$\lambda_i$ and $e_i$ $(i\in J)$ are
its eigenvalues and eigenfunctions {\rm (}which form
an orthonormal basis of $U${\rm )}
correspondingly{\rm )} and ${\bf W}_{\tau},$ $\tau\in [0, \bar T]$
is an $U$-valued $Q$-Wiener process.

{\rm 2}. $Z: \Omega \rightarrow H$ is 
an ${\bf F}_t/{\mathcal{B}}(H)$-measurable mapping.

{\rm 3}. $\Phi_1\in L(U, H)_{0},$ 
$\Phi_2\in L(H,L(U,H)_0),$ and 
$\Phi_k(v)(\ \ldots (\Phi_2(v)(\Phi_1(v))) \ldots\ )$ 
is a $k$-linear Hilbert--Schmidt operator mapping from
$\underbrace{U_0\times \ldots \times U_0}_{\small{\hbox{$k$ times}}}$ to $H$
for all $v\in H$
such that

\newpage
\noindent
$$
\Biggl\Vert \Phi_k(Z)\left(\ldots \left(\Phi_2(Z) \left(\Phi_1(Z)
e_{r_1} \right) 
e_{r_2} \right) \ldots \right) e_{r_k}\Biggr\Vert_H^2
\le L_k<\infty
$$

\vspace{3mm}
\noindent
w.~p.~{\rm 1}\ for all $r_1,r_2,\ldots,r_k\in J_M$, $M\in\mathbb{\bf N}$.

Then

$$
{\sf M}\left\{
\Biggl\Vert
I[\Phi^{(k)}(Z), \psi^{(k)}]_{T,t}^M-
I[\Phi^{(k)}(Z), \psi^{(k)}]_{T,t}^{M,p_1,\ldots, p_k}
\Biggr\Vert_H^2\right\}\le 
$$

\begin{equation}
\label{zzz1ee}
\le L_k (k!)^2
\left({\rm tr}\ Q\right)^k
\left(I_k-\sum_{j_1=0}^{p_1}\ldots
\sum_{j_k=0}^{p_k}C^2_{j_k\ldots j_1}\right),
\end{equation}

\vspace{3mm}
\noindent
where 
$$
{\rm tr}\ Q=\sum\limits_{i\in J}\lambda_i<\infty,
$$

$$
I_k=\left\Vert K\right\Vert_{L_2([t,T]^k)}^2=
\int\limits_{[t,T]^k} K^2(t_1,\ldots,t_k)
dt_1\ldots dt_k,
$$

$$
C_{j_k\ldots j_1}=\int\limits_{[t,T]^k}
K(t_1,\ldots,t_k)\prod_{l=1}^{k}\phi_{j_l}(t_l)dt_1\ldots dt_k
$$

\vspace{2mm}
\noindent
is the Fourier coefficient,

\vspace{-4mm}
$$
K(t_1,\ldots,t_k)=
\left\{\begin{matrix}
\psi_1(t_1)\ldots \psi_k(t_k),\ &t_1<\ldots<t_k\cr\cr
0,\ &\hbox{\rm otherwise}
\end{matrix}
\right.\ \ \
=\ \
\prod\limits_{l=1}^k
\psi_l(t_l)\ \prod\limits_{l=1}^{k-1}{\bf 1}_{\{t_l<t_{l+1}\}},
$$

\vspace{4mm}
\noindent
where $t_1,\ldots,t_k\in [t, T]$ $(k\ge 2)$ and 
$K(t_1)\equiv\psi_1(t_1)$ for $t_1\in[t, T]$ 
{\rm (}${\bf 1}_A$ denotes the indicator of the set $A${\rm )}. 
}

\vspace{2mm}

{\bf Remark 7.2.}\ {\it It should be noted that the right-hand side 
of the inequality {\rm (\ref{zzz1ee})} is independent of $M$
and tends to zero if $p_1,\ldots,p_k\to\infty$ due to the
Parseval equality.}

{\bf Remark 7.3.}\ {\it Recall the estimate {\rm (\ref{z1new}),} which
we will use in the proof of Theorem {\rm 7.1}
$$
{\sf M}\left\{\biggl(
J[\psi^{(k)}]_{T,t}^{(r_1 r_2\ldots r_k)}-
J[\psi^{(k)}]_{T,t}^{(r_1 r_2\ldots r_k)p_1,\ldots ,p_k}\biggr)^2\right\}
\le 
$$

\vspace{-1mm}
$$
\le k!
\left(I_k-\sum_{j_1=0}^{p_1}\ldots
\sum_{j_k=0}^{p_k}C^2_{j_k\ldots j_1}\right),
$$

\vspace{3mm}
\noindent
where $J[\psi^{(k)}]_{T,t}^{(r_1 r_2\ldots r_k)}$
is defined by {\rm (\ref{ito})} and 
$J[\psi^{(k)}]_{T,t}^{(r_1 r_2\ldots r_k)p_1,\ldots ,p_k}$
is defined by {\rm (\ref{f1})}.  
}

{\bf Proof.} Using (\ref{z1}), we obtain

\vspace{-2mm}
$$
{\sf M}\left\{
\Biggl\Vert
I[\Phi^{(k)}(Z), \psi^{(k)}]_{T,t}^M-
I[\Phi^{(k)}(Z), \psi^{(k)}]_{T,t}^{M,p_1,\ldots ,p_k}
\Biggr\Vert_H^2\right\}=
$$

\vspace{1mm}
$$
={\sf M}\left\{
\Biggl\Vert
\sum_{r_1,r_2,\ldots,r_k\in J_M}
\Phi_k(Z)\left(\ldots \left(\Phi_2(Z) \left(\Phi_1(Z)
e_{r_1} \right) 
e_{r_2} \right) \ldots \right) e_{r_k}
\left(\prod_{l=1}^k\lambda_{r_l}\right)^{\hspace{-2.2mm}1/2}
\times\right.\Biggr.
$$
\begin{equation}
\label{zz1ququ}
\times
\left.\Biggl.\Biggl(J[\psi^{(k)}]_{T,t}^{(r_1 r_2\ldots r_k)}-
J[\psi^{(k)}]_{T,t}^{(r_1 r_2\ldots r_k)p_1,\ldots, p_k}\Biggr)
\Biggr\Vert_H^2\right\}=
\end{equation}

\vspace{4mm}
$$
=\Biggl|
{\sf M}\Biggl\{\sum_{r_1,r_2,\ldots,r_k\in J_M}\
\sum_{(r_1^{,},r_2^{,},\ldots,r_k^{,}):\ 
\{r_1^{,},r_2^{,},\ldots,r_k^{,}\}=\{r_1,r_2,\ldots,r_k\}}
\left(\prod_{l=1}^k\lambda_{r_l}\right)^{\hspace{-2.2mm}1/2}
\left(\prod_{l=1}^k\lambda_{r_l^{,}}\right)^{\hspace{-2.2mm}1/2}
\times\Biggr.
$$
$$
\hspace{-30mm}\times
\Biggl\langle \Phi_k(Z)\left(\ldots \left(\Phi_2(Z) \left(\Phi_1(Z)
e_{r_1} \right) 
e_{r_2} \right) \ldots \right) e_{r_k}\ ,\Biggr.
$$
$$
~~~~~~~~~~~~~~~\Biggl.
\Phi_k(Z)\left(\ldots \left(\Phi_2(Z) \left(\Phi_1(Z)
e_{r_1^{,}} \right) 
e_{r_2^{,}} \right) \ldots \right) e_{r_k^{,}} \Biggr\rangle_H
\times
$$
$$
\hspace{-30mm}\times\
{\sf M}\Biggl\{\Biggl(J[\psi^{(k)}]_{T,t}^{(r_1 r_2\ldots r_k)}-
J[\psi^{(k)}]_{T,t}^{(r_1 r_2\ldots r_k)p_1,\ldots ,p_k}\Biggr)\times\Biggr.
$$
\begin{equation}
\label{zz100}
\Biggl.
\Biggl.
\Biggl.
~~~~~~~~~~~~~~~~~~~\times
\Biggl(J[\psi^{(k)}]_{T,t}^{(r_1^{,} r_2^{,}\ldots r_k^{,})}-
J[\psi^{(k)}]_{T,t}^{(r_1^{,} r_2^{,}\ldots r_k^{,})p_1,\ldots, p_k}\Biggr)
\biggl.\biggr|{\bf F}_t\Biggr\}
\Biggr\}\Biggr|\le
\end{equation}

\newpage
\noindent
$$
\le
\sum_{r_1,r_2,\ldots,r_k\in J_M}\ 
\sum_{(r_1^{,},r_2^{,},\ldots,r_k^{,}):\ 
\{r_1^{,},r_2^{,},\ldots,r_k^{,}\}=\{r_1,r_2,\ldots,r_k\}}
\left(\prod_{l=1}^k\lambda_{r_l}\right)^{\hspace{-2.2mm}1/2}
\left(\prod_{l=1}^k\lambda_{r_l^{,}}\right)^{\hspace{-2.2mm}1/2}\times
$$

$$
\hspace{-20mm}
\times
{\sf M}\Biggl\{
\Biggl\Vert \Phi_k(Z)\left(\ldots \left(\Phi_2(Z) \left(\Phi_1(Z)
e_{r_1} \right) 
e_{r_2} \right) \ldots \right) e_{r_k} \Biggr\Vert_H\times \Biggr.
$$
$$
~~~~~~~~~~~~~~~~\times
\Biggl\Vert \Phi_k(Z)\left(\ldots \left(\Phi_2(Z) \left(\Phi_1(Z)
e_{r_1^{,}} \right) 
e_{r_2^{,}} \right) \ldots \right) e_{r_k^{,}}\ \Biggr\Vert_H
\times
$$
$$
\hspace{-30mm}\times
\left|{\sf M}\Biggl\{\Biggl(J[\psi^{(k)}]_{T,t}^{(r_1 r_2\ldots r_k)}-
J[\psi^{(k)}]_{T,t}^{(r_1 r_2\ldots r_k)p_1,\ldots, p_k}\Biggr)\times
\Biggr.\right.
$$
$$
\Biggl.\left.\Biggl.
~~~~~~~~~~~~~~~~\times
\Biggl(J[\psi^{(k)}]_{T,t}^{(r_1^{,} r_2^{,}\ldots r_k^{,})}-
J[\psi^{(k)}]_{T,t}^{(r_1^{,} r_2^{,}\ldots r_k^{,})p_1,\ldots, p_k}\Biggr)
\biggl.\biggr|{\bf F}_t\Biggr\}
\right|
\Biggr\}\le
$$

\vspace{5mm}
$$
\le L_k
\sum_{r_1,r_2,\ldots,r_k\in J_M}\ \ 
\sum_{(r_1^{,},r_2^{,},\ldots,r_k^{,}):\ 
\{r_1^{,},r_2^{,},\ldots,r_k^{,}\}=\{r_1,r_2,\ldots,r_k\}}\ \
\left(\prod_{l=1}^k\lambda_{r_l}\right)^{\hspace{-2.2mm}1/2}
\left(\prod_{l=1}^k\lambda_{r_l^{,}}\right)^{\hspace{-2.2mm}1/2}\times
$$

$$
\hspace{-20mm}
\times{\sf M}\Biggl\{\left|\Biggl(J[\psi^{(k)}]_{T,t}^{(r_1 r_2\ldots r_k)}-
J[\psi^{(k)}]_{T,t}^{(r_1 r_2\ldots r_k)p_1, \ldots, p_k}\Biggr)\times
\right.\Biggr.
$$
$$
~~~~~~~~~\Biggl.\left.\times
\Biggl(J[\psi^{(k)}]_{T,t}^{(r_1^{,} r_2^{,}\ldots r_k^{,})}-
J[\psi^{(k)}]_{T,t}^{(r_1^{,} r_2^{,}\ldots r_k^{,})p_1,\ldots ,p_k}\Biggr)
\right|\Biggr\}\le
$$

\vspace{5mm}
$$
\le L_k
\sum_{r_1,r_2,\ldots,r_k\in J_M}\ \ 
\sum_{(r_1^{,},r_2^{,},\ldots,r_k^{,}):\ 
\{r_1^{,},r_2^{,},\ldots,r_k^{,}\}=\{r_1,r_2,\ldots,r_k\}}\ \
\left(\prod_{l=1}^k\lambda_{r_l}\right)^{\hspace{-2.2mm}1/2}
\left(\prod_{l=1}^k\lambda_{r_l^{,}}\right)^{\hspace{-2.2mm}1/2}\times
$$
$$
\hspace{-3mm}\times
\left({\sf M}\left\{\Biggl(J[\psi^{(k)}]_{T,t}^{(r_1 r_2\ldots r_k)}-
J[\psi^{(k)}]_{T,t}^{(r_1 r_2\ldots r_k)p_1, \ldots ,p_k}\Biggr)^2\right\}
\right)^{1/2}\times
$$
$$
~~~~~~~~~~~~~~~~~~~\times
\left({\sf M}\left\{
\Biggl(J[\psi^{(k)}]_{T,t}^{(r_1^{,} r_2^{,}\ldots r_k^{,})}-
J[\psi^{(k)}]_{T,t}^{(r_1^{,} r_2^{,}\ldots r_k^{,})p_1,\ldots ,p_k}\Biggr)^2
\right\}
\right)^{1/2}\le
$$

\newpage
\noindent
$$
\le L_k
\sum_{r_1,r_2,\ldots,r_k\in J_M}\ \ 
\sum_{(r_1^{,},r_2^{,},\ldots,r_k^{,}):\ 
\{r_1^{,},r_2^{,},\ldots,r_k^{,}\}=\{r_1,r_2,\ldots,r_k\}}\ \
\left(\prod_{l=1}^k\lambda_{r_l}\right)^{\hspace{-2.2mm}1/2}
\left(\prod_{l=1}^k\lambda_{r_l^{,}}\right)^{\hspace{-2.2mm}1/2}\times
$$

$$
\times
\left(k!\left(I_k-\sum_{j_1=0}^{p_1}\ldots
\sum_{j_k=0}^{p_k}C^2_{j_k\ldots j_1}\right)\right)^{1/2}
\left(k!\left(I_k-\sum_{j_1=0}^{p_1}\ldots
\sum_{j_k=0}^{p_k}C^2_{j_k\ldots j_1}\right)\right)^{1/2}\le
$$

\vspace{1mm}
$$
\le L_k
\sum_{r_1,r_2,\ldots,r_k\in J_M} k!\ 
\lambda_{r_1}\lambda_{r_2}\ldots \lambda_{r_k}
\left(k!\left(I_k-\sum_{j_1=0}^{p_1}\ldots
\sum_{j_k=0}^{p_k}C^2_{j_k\ldots j_1}\right)\right)=
$$

\vspace{1mm}
$$
=L_k \left(k!\right)^2
\sum_{r_1,r_2,\ldots,r_k\in J_M}
\lambda_{r_1}\lambda_{r_2}\ldots \lambda_{r_k}
\left(I_k-\sum_{j_1=0}^{p_1}\ldots
\sum_{j_k=0}^{p_k}C^2_{j_k\ldots j_1}\right)\le
$$

\vspace{1mm}
$$
\le L_k \left(k!\right)^2
\left({\rm tr}\ Q\right)^k
\left(I_k-\sum_{j_1=0}^{p_1}\ldots
\sum_{j_k=0}^{p_k}C^2_{j_k\ldots j_1}\right),
$$

\vspace{5mm}
\noindent
where $\langle \cdot,\cdot\rangle_H$ is a scalar product in $H,$ and
$$
\sum_{(r_1^{,},r_2^{,},\ldots,r_k^{,}):\ 
\{r_1^{,},r_2^{,},\ldots,r_k^{,}\}=\{r_1,r_2,\ldots,r_k\}}
$$
means the sum with respect to all
possible permutations
$(r_1^{,},r_2^{,},\ldots,r_k^{,})$ such that
$\{r_1^{,},r_2^{,},\ldots,r_k^{,}\}=\{r_1,r_2,\ldots,r_k\}$.

The transition from (\ref{zz1ququ}) to (\ref{zz100}) 
is based on the following theorem.

{\bf Theorem 7.2} \cite{12a}-\cite{12aa}, \cite{art-7}, \cite{arxiv-20}. {\it 
Suppose that
$\psi_1(\tau),\ldots,\psi_k(\tau)\in L_2([t, T])$ and
$\{\phi_j(x)\}_{j=0}^{\infty}$ is an arbitrary complete orthonormal system  
of functions in the space $L_2([t,T]).$
Then$,$ the following 
equality is true
$$
\hspace{-40mm}
{\sf M}
\Biggl\{\Biggl(
J[\psi^{(k)}]_{T,t}^{(r_1\ldots r_k)}-
J[\psi^{(k)}]_{T,t}^{(r_1\ldots r_k)p_1,\ldots, p_k}
\Biggr)\times\Biggr.
$$
\begin{equation}
\label{uuu2}
~~~~~~~~~~~~~~~~~\Biggl.\times
\Biggl(
J[\psi^{(k)}]_{T,t}^{(m_1\ldots m_k)}-
J[\psi^{(k)}]_{T,t}^{(m_1\ldots m_k)p_1,\ldots, p_k}
\Biggr)\biggl.\biggr|{\bf F}_t\Biggr\}=0
\end{equation}

\noindent
w.~p.~{\rm 1} for all $r_1,\ldots,r_k,m_1,\ldots,m_k\in 
J_M$ $(M\in{\bf N})$
such that $\{r_1,\ldots,r_k\}\ne 
\{m_1,\ldots,m_k\}$.}

{\bf Proof.} Using the standard moment properties 
of the It\^{o} stochastic integral,
we obtain
\begin{equation}
\label{uuu3}
{\sf M}
\biggl\{J[\psi^{(k)}]_{T,t}^{(r_1\ldots r_k)}
J[\psi^{(k)}]_{T,t}^{(m_1\ldots m_k)}\biggl.\biggr|{\bf F}_t
\biggr\}=0
\end{equation}

\noindent
w.~p.~1 for all $r_1,\ldots,r_k,m_1,\ldots,m_k\in J_M$ $(M\in{\bf N})$
such that $(r_1,\ldots,r_k)\ne 
(m_1,\ldots,m_k)$.

Using (\ref{chain102}), (\ref{chain7771}), (\ref{chain7878}), and (\ref{chain104}), we obtain

\vspace{-4mm}
\begin{equation}
\label{ttt1}
~~~~~~J[\psi^{(k)}]_{T,t}^{(m_1\ldots m_k)p_1,\ldots,p_k}=
\sum\limits_{j_1=0}^{p_1}\ldots
\sum\limits_{j_k=0}^{p_k}
C_{j_k\ldots j_1}J'[\phi_{j_1}\ldots \phi_{j_k}]_{T,t}^{(m_1\ldots m_k)},
\end{equation}

\vspace{1mm}
\noindent
where
$$
J'[\phi_{j_1}\ldots \phi_{j_k}]_{T,t}^{(m_1\ldots m_k)}
=
$$

\vspace{-4mm}
\begin{equation}
\label{ttt2www}
~~~~=\sum\limits_{(j_1,\ldots,j_k)}
\int\limits_t^T \phi_{j_k}(t_k)
\ldots
\int\limits_t^{t_{2}}\phi_{j_{1}}(t_{1})
d{\bf w}_{t_1}^{(m_1)}\ldots
d{\bf w}_{t_k}^{(m_k)}\ \ \ {\rm w.~p.~1,}
\end{equation}

\vspace{2mm}
\noindent
and 
$$
\sum\limits_{(j_1,\ldots,j_k)}
$$ 

\noindent
means the sum with respect to all
possible permutations
$(j_1,\ldots,j_k)$. At the same time if 
$j_r$ swapped  with $j_q$ in the permutation $(j_1,\ldots,j_k)$,
then $m_r$ swapped  with $m_q$ in the permutation
$(m_1,\ldots,m_k).$
Another notations are the same as in Theorems 1.1, 1.2, 1.16 
($J'[\phi_{j_1}\ldots \phi_{j_k}]_{T,t}^{(m_1\ldots m_k)}$ is defined by (\ref{WiI})).

Then w.~p.~1
$$
{\sf M}
\biggl\{J[\psi^{(k)}]_{T,t}^{(r_1\ldots r_k)}
J[\psi^{(k)}]_{T,t}^{(m_1\ldots m_k)p_1,\ldots ,p_k}
\biggl.\biggr|{\bf F}_t\biggr\}
=\sum_{j_1=0}^{p_1}\ldots\sum_{j_k=0}^{p_k}
C_{j_k\ldots j_1}\times
$$
$$
\times{\sf M}\left\{J[\psi^{(k)}]_{T,t}^{(r_1\ldots r_k)}
\sum\limits_{(j_1,\ldots,j_k)}
\int\limits_t^T \phi_{j_k}(t_k)
\ldots
\int\limits_t^{t_{2}}\phi_{j_{1}}(t_{1})
d{\bf w}_{t_1}^{(m_1)}\ldots
d{\bf w}_{t_k}^{(m_k)}
\biggl.\biggr|{\bf F}_t\right\}.
$$

\vspace{3mm}

From the standard moment properties of the It\^{o} stochastic integral 
it follows that
$$
{\sf M}\left\{J[\psi^{(k)}]_{T,t}^{(r_1\ldots r_k)}
\sum\limits_{(j_1,\ldots,j_k)}
\int\limits_t^T \phi_{j_k}(t_k)
\ldots
\int\limits_t^{t_{2}}\phi_{j_{1}}(t_{1})
d{\bf w}_{t_1}^{(m_1)}\ldots
d{\bf w}_{t_k}^{(m_k)}
\biggl.\biggr|{\bf F}_t\right\}=0
$$

\noindent
w.~p.~1 for all $r_1,\ldots,r_k,m_1,\ldots,m_k\in J_M$ $(M\in{\bf N})$
such that $\{r_1,\ldots,r_k\}\ne 
\{m_1,\ldots,m_k\}$.

Then 
\begin{equation}
\label{uuu4}
{\sf M}
\biggl\{J[\psi^{(k)}]_{T,t}^{(r_1\ldots r_k)}
J[\psi^{(k)}]_{T,t}^{(m_1\ldots m_k)p_1,\ldots, p_k}
\biggl.\biggr|{\bf F}_t\biggr\}=0
\end{equation}

\noindent
w.~p.~1 for all $r_1,\ldots,r_k,m_1,\ldots,m_k\in J_M$ $(M\in{\bf N})$
such that $\{r_1,\ldots,r_k\}\ne 
\{m_1,\ldots,m_k\}$.

Using (\ref{ttt1}), (\ref{ttt2www}), we have
$$
{\sf M}\biggl\{J[\psi^{(k)}]_{T,t}^{(r_1\ldots r_k)p_1,\ldots, p_k}
J[\psi^{(k)}]_{T,t}^{(m_1\ldots m_k)p_1,\ldots, p_k} 
\biggl.\biggr|{\bf F}_t\biggr\}=
$$
$$
=\sum_{j_1=0}^{p_1}\ldots\sum_{j_k=0}^{p_k}
C_{j_k\ldots j_1}
\sum_{{q}_1=0}^{p_1}\ldots\sum_{{q}_k=0}^{p_k}
C_{{q}_k\ldots {q}_1}\times
$$
$$
\hspace{-40mm}\times
{\sf M}\left\{\left(
\sum\limits_{(j_1,\ldots,j_k)}
\int\limits_t^T \phi_{j_k}(t_k)
\ldots
\int\limits_t^{t_{2}}\phi_{j_{1}}(t_{1})
d{\bf w}_{t_1}^{(r_1)}\ldots
d{\bf w}_{t_k}^{(r_k)}\right) \times\right.
$$
\begin{equation}
\label{uuu5}
~~~~~~~~~~~~~~~~~~~~\times
\left.\left(
\sum\limits_{(q_1,\ldots,q_k)}
\int\limits_t^T \phi_{q_k}(t_k)
\ldots
\int\limits_t^{t_{2}}\phi_{q_{1}}(t_{1})
d{\bf w}_{t_1}^{(m_1)}\ldots
d{\bf w}_{t_k}^{(m_k)}\right) \biggl.\biggr|{\bf F}_t\right\}=0
\end{equation}

\noindent
w.~p.~1 for all $r_1,\ldots,r_k,m_1,\ldots,m_k\in J_M$ $(M\in{\bf N})$
such that $\{r_1,\ldots,r_k\}\ne 
\{m_1,\ldots,m_k\}$.

From (\ref{uuu3}), (\ref{uuu4}), and (\ref{uuu5}) 
we obtain (\ref{uuu2}).
Theorem 7.2 is proved.

\vspace{1mm}

{\bf Corollary 7.1} \cite{12a}-\cite{12aa}, \cite{art-7}, \cite{arxiv-20}. 
{\it Suppose that $\{\phi_j(x)\}_{j=0}^{\infty}$ is an arbitrary complete orthonormal system  
of functions in the space $L_2([t,T])$
and
$\psi_1(\tau),\ldots,\psi_k(\tau)\in L_2([t, T]).$ 
Then$,$ the following equality is true
$$
\hspace{-40mm}
{\sf M}
\Biggl\{\Biggl(
J[\psi^{(k)}]_{T,t}^{(r_1\ldots r_k)}-
J[\psi^{(k)}]_{T,t}^{(r_1\ldots r_k)p_1,\ldots, p_k}
\Biggr)\times\Biggr.
$$
$$
~~~~~~~~~~~~~~~~~\Biggl.\times
\Biggl(
J[\psi^{(l)}]_{T,t}^{(m_1\ldots m_l)}-
J[\psi^{(l)}]_{T,t}^{(m_1\ldots m_l)q_1,\ldots, q_l}
\Biggr)\biggl.\biggr|{\bf F}_t\Biggr\}=0
$$

\noindent
w.~p.~{\rm 1}\ for all $l=1, 2, \ldots, k-1$ and 
$r_1,\ldots,r_k,$ $m_1,\ldots,m_l$ $\in $ $J_M,$ $p_1,\ldots,p_k,$
$q_1,\ldots,q_l$
$=$ $0,1,2,\ldots $}

\subsection{Approximation of Some Iterated Stochastic Integrals 
of Miltiplicities 2 and 3 with
Respect to the Finite-Dimensional Approximation ${\bf W}_t^M$
of the $Q$-Wiener Process}

This section is devoted to the approximation of iterated
stochastic integrals of the following form (see Sect.~7.1)
\begin{equation}
\label{abc1}
\hspace{-10mm}I_0[B(Z),F(Z)]_{T,t}^M=
\int\limits_{t}^{T}B'(Z) \left(
\int\limits_{t}^{t_2}F(Z) dt_1 \right) d{\bf W}_{t_2}^M,
\end{equation}
\begin{equation}
\label{kkk1}
\hspace{-10mm}I_1[B(Z),F(Z)]_{T,t}^M=\int\limits_{t}^{T}F'(Z)
\left(\int\limits_{t}^{t_2} B(Z) d{\bf W}_{t_1}^M \right) dt_2,
\end{equation}
\begin{equation}
\label{kkk2}
~~~~~I_2[B(Z)]_{T,t}^M=\int\limits_{t}^{T}B''(Z) \left(
\int\limits_{t}^{t_2}B(Z) d{\bf W}_{t_1}^M,
\int\limits_{t}^{t_2}B(Z) d{\bf W}_{t_1}^M 
\right) d{\bf W}_{t_2}^M.\
\end{equation}

\vspace{1mm}

Let Conditions 1, 2 of Theorem 7.1 be fulfilled.
Let $B''(v)(B(v),B(v))$ 
be a 3-linear Hilbert--Schmidt operator mapping from
$U_0\times U_0\times U_0$ to $H$ for all $v\in H$.
Then we have w.~p.~1 (see (\ref{xx605}))

\vspace{-1mm}
\begin{equation}
\label{abc2}
I_0[B(Z),F(Z)]_{T,t}^M=
\sum_{r_1\in J_M}B'(Z)F(Z)e_{r_1}\sqrt{\lambda_{r_1}}J_{(01)T,t}^{(0 r_1)},
\end{equation}

\vspace{-1mm}
\begin{equation}
\label{da0}
I_1[B(Z),F(Z)]_{T,t}^M=
\sum_{r_1\in J_M}F'(Z)(B(Z)e_{r_1})\sqrt{\lambda_{r_1}}J_{(10)T,t}^{(r_1 0)},
\end{equation}

$$
I_2[B(Z)]_{T,t}^M=
\sum_{r_1,r_2,r_3\in J_M}B''(Z)\left(B(Z)e_{r_1}, B(Z)e_{r_2}\right)e_{r_3}
\sqrt{\lambda_{r_1}\lambda_{r_2}\lambda_{r_3}}\times
$$
\begin{equation}
\label{jjj1}
\times \int\limits_t^T
\left(\int\limits_t^s d{\bf w}_{\tau}^{(r_1)}   
\int\limits_t^s d{\bf w}_{\tau}^{(r_2)}\right)
d{\bf w}_s^{(r_3)}.
\end{equation}

\vspace{2mm}

Using the It\^{o} formula, we obtain
\begin{equation}
\label{da12eee}
~~~~~~ \int\limits_t^s d{\bf w}_{\tau}^{(r_1)}   
\int\limits_t^s d{\bf w}_{\tau}^{(r_2)}
=
J_{(11)s,t}^{(r_1 r_2)}+
J_{(11)s,t}^{(r_2 r_1)}+
{\bf 1}_{\{r_1=r_2\}}(s-t)\ \ \ \hbox{w.~p.~1.}
\end{equation}

From (\ref{da12eee}) we have
\begin{equation}
\label{jjj2sss}
\int\limits_t^T
\left(\int\limits_t^s d{\bf w}_{\tau}^{(r_1)}   
\int\limits_t^s d{\bf w}_{\tau}^{(r_2)}\right)
d{\bf w}_s^{(r_3)}=
J_{(111)T,t}^{(r_1 r_2 r_3)}+
J_{(111)T,t}^{(r_2 r_1 r_3)}+
{\bf 1}_{\{r_1=r_2\}}J_{(01)T,t}^{(0 r_3)}\ \ \ \hbox{w.~p.~1.}
\end{equation}

Note that in 
(\ref{abc2}), (\ref{da0}), (\ref{da12eee}), and (\ref{jjj2sss}) 
we use the notations from Sect.~7.2 (see (\ref{x111}), (\ref{x112})).
After substituting (\ref{jjj2sss}) into (\ref{jjj1}), we have

\vspace{-5mm}
$$
I_2[B(Z)]_{T,t}^M=
\sum_{r_1,r_2,r_3\in J_M}B''(Z)\left(B(Z)e_{r_1}, B(Z)e_{r_2}\right)e_{r_3}
\sqrt{\lambda_{r_1}\lambda_{r_2}\lambda_{r_3}} \times
$$
\begin{equation}
\label{da1}
\times
\left(J_{(111)T,t}^{(r_1 r_2 r_3)}+
J_{(111)T,t}^{(r_2 r_1 r_3)}+
{\bf 1}_{\{r_1=r_2\}}J_{(01)T,t}^{(0 r_3)}\right)\ \ \ \hbox{w.~p.~1.}
\end{equation}

\vspace{3mm}

Taking into account (\ref{y3a}), (\ref{y4a}), we put for $q=1$

\vspace{-3mm}
\begin{equation}
\label{opp3}
~~~~~~~~ J_{(01)T,t}^{(0 r_3)q}=
J_{(01)T,t}^{(0 r_3)}=\frac{(T-t)^{3/2}}{2}\biggl(\zeta_0^{(r_3)}+
\frac{1}{\sqrt{3}}\zeta_1^{(r_3)}\biggr)\ \ \ \hbox{w.~p.~1},
\end{equation}

\vspace{-4mm}
\begin{equation}
\label{opp4}
~~~~~~~~ J_{(10)T,t}^{(r_1 0)q}=
J_{(10)T,t}^{(r_1 0)}=\frac{(T-t)^{3/2}}{2}\biggl(\zeta_0^{(r_1)}-
\frac{1}{\sqrt{3}}\zeta_1^{(r_1)}\biggr)\ \ \ \hbox{w.~p.~1},
\end{equation}

\vspace{4mm}
\noindent
where 
$J_{(01)T,t}^{(0 r_3)q},$
$J_{(10)T,t}^{(r_1 0)q}$ denote the approximations of
corresponding
iterated It\^{o} stochastic integrals.

Denote by
$I_0[B(Z),F(Z)]_{T,t}^{M,q},$ $I_1[B(Z),F(Z)]_{T,t}^{M,q},$ 
$I_2[B(Z)]_{T,t}^{M,q}$ the approximations of iterated stochastic integrals
(\ref{abc2}), (\ref{da0}), (\ref{da1})

\vspace{-1mm}
\begin{equation}
\label{abc3}
I_0[B(Z),F(Z)]_{T,t}^{M,q}=
\sum_{r_1\in J_M}B'(Z)F(Z)e_{r_1}\sqrt{\lambda_{r_1}}J_{(01)T,t}^{(0 r_1)q},
\end{equation}

\vspace{-1mm}
\begin{equation}
\label{da3}
I_1[B(Z),F(Z)]_{T,t}^{M,q}=
\sum_{r_1\in J_M}F'(Z)(B(Z)e_{r_1})\sqrt{\lambda_{r_1}}J_{(10)T,t}^{(r_1 0)q},
\end{equation}

\vspace{-2mm}
$$
I_2[B(Z)]_{T,t}^{M,q}=
\sum_{r_1,r_2,r_3\in J_M}B''(Z)\left(B(Z)e_{r_1}, B(Z)e_{r_2}\right)e_{r_3}
\sqrt{\lambda_{r_1}\lambda_{r_2}\lambda_{r_3}} \times
$$
\begin{equation}
\label{da4www}
\times
\left(J_{(111)T,t}^{(r_1 r_2 r_3)q}+
J_{(111)T,t}^{(r_2 r_1 r_3)q}+
{\bf 1}_{\{r_1=r_2\}}J_{(01)T,t}^{(0 r_3)q}\right),
\end{equation}

\vspace{3mm}
\noindent
where $q=1$ in (\ref{abc3}), (\ref{da3}) 
and the approximations $J_{(111)T,t}^{(r_1 r_2 r_3)q},$
$J_{(111)T,t}^{(r_2 r_1 r_3)q}$ are defined by 
(\ref{sad001}) for some $q\ge 1$.

From (\ref{abc2}), (\ref{da0}), (\ref{da1}),
(\ref{abc3})--(\ref{da4www}) we have

\vspace{-2mm}
$$
I_0[B(Z),F(Z)]_{T,t}^{M}-I_0[B(Z),F(Z)]_{T,t}^{M,q}=0\ \ \ \hbox{w.~p.~1,}
$$

\vspace{-2mm}
$$
I_1[B(Z),F(Z)]_{T,t}^{M}-I_1[B(Z),F(Z)]_{T,t}^{M,q}=0\ \ \ \hbox{w.~p.~1,}
$$

\vspace{-2mm}
$$
I_2[B(Z)]_{T,t}^{M}-I_2[B(Z)]_{T,t}^{M,q}=
$$

\vspace{-3mm}
$$
=
\sum_{r_1,r_2,r_3\in J_M}B''(Z)\left(B(Z)e_{r_1}, B(Z)e_{r_2}\right)e_{r_3}
\sqrt{\lambda_{r_1}\lambda_{r_2}\lambda_{r_3}} \times
$$
$$
\times
\left(\left(J_{(111)T,t}^{(r_1 r_2 r_3)}-J_{(111)T,t}^{(r_1 r_2 r_3)q}\right)+
\left(J_{(111)T,t}^{(r_2 r_1 r_3)}-J_{(111)T,t}^{(r_2 r_1 r_3)q}
\right)
\right)\ \ \ \hbox{w.~p.~1.}
$$

\vspace{4mm}

Repeating with an insignificant 
modification the proof 
of Theorem 7.1 for the case $k=3$, we obtain

\vspace{-2mm}
$$
\hspace{-65mm}
{\sf M}\left\{\biggl\Vert
I_2[B(Z)]_{T,t}^{M}-I_2[B(Z)]_{T,t}^{M,q}\biggr\Vert_H^2\right\}\le
$$
$$
~~~~~~~~~~~~~~~~~~~~~~~~~~~\le 4C(3!)^2
\left({\rm tr}\ Q\right)^3
\left(\frac{(T-t)^{3}}{6}-\sum_{j_1,j_2,j_3=0}^{q}
C_{j_3j_2j_1}^2\right),
$$

\vspace{3mm}
\noindent
where here and further constant $C$ has the same meaning 
as constant $L_k$ in Theorem 7.1
($k$ is the multiplicity of the iterated stochastic integral),
and

\vspace{-1mm}
$$
C_{j_3j_2j_1}=\frac{\sqrt{(2j_1+1)(2j_2+1)(2j_3+1)}}{8}(T-t)^{3/2}\bar
C_{j_3j_2j_1},
$$
$$
\bar C_{j_3j_2j_1}=\int\limits_{-1}^{1}P_{j_3}(z)
\int\limits_{-1}^{z}P_{j_2}(y)
\int\limits_{-1}^{y}
P_{j_1}(x)dx dy dz,
$$

\vspace{2mm}
\noindent
where $P_j(x)$ is the Legendre polynomial.

\subsection{Approximation of Some Iterated Stochastic Integrals
of Miltiplicities 3 and 4 with
Respect to the Finite-Dimensional Approximation ${\bf W}_t^M$ of 
the $Q$-Wiener Process}

In this section, we consider the approximation of iterated
stochastic integrals of the following form (see Sect.~7.1)
$$
I_3[B(Z)]_{T,t}^M=\int\limits_{t}^{T}B'''(Z) \left(
\int\limits_{t}^{t_2}B(Z) d{\bf W}_{t_1}^M,
\int\limits_{t}^{t_2}B(Z) d{\bf W}_{t_1}^M,
\int\limits_{t}^{t_2}B(Z) d{\bf W}_{t_1}^M 
\right) d{\bf W}_{t_2}^M,\
$$
$$
I_4[B(Z)]_{T,t}^M=
$$
$$
=\int\limits_{t}^{T}B'(Z) \left(
\int\limits_{t}^{t_3}B''(Z) \left(
\int\limits_{t}^{t_2}B(Z) d{\bf W}_{t_1}^M,
\int\limits_{t}^{t_2}B(Z) d{\bf W}_{t_1}^M 
\right) d{\bf W}_{t_2}^M\right)d{\bf W}_{t_3}^M,\
$$
$$
I_5[B(Z)]_{T,t}^M=
$$
$$
=\int\limits_{t}^{T}B''(Z) \left(
\int\limits_{t}^{t_3}B(Z)d{\bf W}_{t_1}^M,
\int\limits_{t}^{t_3}B'(Z)\left(
\int\limits_{t}^{t_2}B(Z) d{\bf W}_{t_1}^M
\right) d{\bf W}_{t_2}^M\right)d{\bf W}_{t_3}^M,\
$$
$$
I_6[B(Z),F(Z)]_{T,t}^M=\int\limits_{t}^{T}F'(Z)\left(
\int\limits_{t}^{t_3}B'(Z) \left(
\int\limits_{t}^{t_2}B(Z) d{\bf W}_{t_1}^M
\right) d{\bf W}_{t_2}^M 
\right) dt_3,
$$
$$
I_7[B(Z),F(Z)]_{T,t}^M=\int\limits_{t}^{T}F''(Z) \left(
\int\limits_{t}^{t_2}B(Z) d{\bf W}_{t_1}^M,
\int\limits_{t}^{t_2}B(Z) d{\bf W}_{t_1}^M 
\right) dt_2,
$$
$$
I_8[B(Z),F(Z)]_{T,t}^M=\int\limits_{t}^{T}B''(Z) \left(
\int\limits_{t}^{t_2}F(Z) dt_1,
\int\limits_{t}^{t_2}B(Z) d{\bf W}_{t_1}^M 
\right) d{\bf W}_{t_2}^M.
$$

Consider the stochastic integral $I_3[B(Z)]_{T,t}^M.$
Let Conditions 1, 2 of Theorem 7.1 be fulfilled.
Let 
$B'''(v)(B(v),B(v),B(v))$
be a 4-linear Hilbert--Schmidt operator mapping from
$U_0\times U_0\times U_0\times U_0$ to $H$ for all $v\in H$.

We have (see (\ref{xx605}))
$$
I_3[B(Z)]_{T,t}^M=
\sum_{r_1,r_2,r_3,r_4\in J_M}
B'''(Z)\left(B(Z)e_{r_1}, B(Z)e_{r_2},B(Z)e_{r_3}\right)e_{r_4}
\sqrt{\lambda_{r_1}\lambda_{r_2}\lambda_{r_3}\lambda_{r_4}}\times
$$
\begin{equation}
\label{jjjk1}
\times\int\limits_t^T
\left(\int\limits_t^s d{\bf w}_{\tau}^{(r_1)}   
\int\limits_t^s d{\bf w}_{\tau}^{(r_2)}
\int\limits_t^s d{\bf w}_{\tau}^{(r_3)}\right)
d{\bf w}_s^{(r_4)}\ \ \ \hbox{w.~p.~1.}
\end{equation}

\vspace{3mm}

By analogy with (\ref{s1s}) or using the It\^{o} formula,
we obtain

\vspace{-3mm}
$$
J_{(1)s,t}^{(r_1)}
J_{(1)s,t}^{(r_2)}
J_{(1)s,t}^{(r_3)}
=
J_{(111)s,t}^{(r_1r_2r_3)}+
J_{(111)s,t}^{(r_1r_3r_2)}+
J_{(111)s,t}^{(r_2r_1r_3)}+
J_{(111)s,t}^{(r_2r_3r_1)}+
J_{(111)s,t}^{(r_3r_1r_2)}+
J_{(111)s,t}^{(r_3r_2r_1)}+
$$

\vspace{-3mm}
$$
+{\bf 1}_{\{r_1=r_2\}}\left(
J_{(10)s,t}^{(r_3 0)}+J_{(01)s,t}^{(0 r_3)}\right)+
{\bf 1}_{\{r_1=r_3\}}\left(
J_{(10)s,t}^{(r_2 0)}+J_{(01)s,t}^{(0 r_2)}\right)+
$$

\vspace{-1mm}
$$
+{\bf 1}_{\{r_2=r_3\}}\left(
J_{(10)s,t}^{(r_1 0)}+J_{(01)s,t}^{(0 r_1)}\right)=
$$

\vspace{-5mm}
\begin{equation}
\label{q12}
=\sum\limits_{(r_1,r_2,r_3)}
J_{(111)s,t}^{(r_1r_2r_3)}+
(s-t)\left({\bf 1}_{\{r_2=r_3\}}J_{(1)s,t}^{(r_1)}+
{\bf 1}_{\{r_1=r_3\}}J_{(1)s,t}^{(r_2)}+
{\bf 1}_{\{r_1=r_2\}}J_{(1)s,t}^{(r_3)}\right)
\end{equation}

\vspace{1mm}
\noindent
w.~p.~1, where 
$$
\sum\limits_{(r_1,r_2,r_3)}
$$ 
means the sum with respect to all
possible permutations
$(r_1,r_2,r_3).$ We also use the notations  
from Sect.~7.2 (see (\ref{x111}), (\ref{x112})).

After substituting (\ref{q12}) into (\ref{jjjk1}), we obtain

\vspace{-4mm}
$$
I_3[B(Z)]_{T,t}^M=
\sum_{r_1,r_2,r_3,r_4\in J_M}
B'''(Z)\left(B(Z)e_{r_1}, B(Z)e_{r_2},B(Z)e_{r_3}\right)e_{r_4}
\sqrt{\lambda_{r_1}\lambda_{r_2}\lambda_{r_3}\lambda_{r_4}}\times
$$
\begin{equation}
\label{jjjk2}
\times 
\left(\sum\limits_{(r_1,r_2,r_3)}
J_{(1111)T,t}^{(r_1r_2r_3r_4)}-
{\bf 1}_{\{r_1=r_2\}}I_{(01)T,t}^{(r_3r_4)}-
{\bf 1}_{\{r_1=r_3\}}I_{(01)T,t}^{(r_2r_4)}-
{\bf 1}_{\{r_2=r_3\}}I_{(01)T,t}^{(r_1r_4)}\right)
\end{equation}
w.~p.~1, 
where $J_{(1111)T,t}^{(r_1r_2r_3r_4)}$ is defined by (\ref{ito1}) and
\begin{equation}
\label{ros1}
I_{(01)T,t}^{(r_1r_2)}=
\int\limits_t^T
(t-s)\int\limits_t^s d{\bf w}_{\tau}^{(r_1)}   
d{\bf w}_s^{(r_2)}.
\end{equation}

Denote by
$I_3[B(Z)]_{T,t}^{M,q}$ the approximation of the iterated stochastic integral
(\ref{jjjk2}),
which has the following form

\newpage
\noindent
$$
I_3[B(Z)]_{T,t}^{M,q}=
\hspace{-2mm}\sum_{r_1,r_2,r_3,r_4\in J_M}
B'''(Z)\left(B(Z)e_{r_1}, B(Z)e_{r_2},B(Z)e_{r_3}\right)e_{r_4}
\sqrt{\lambda_{r_1}\lambda_{r_2}\lambda_{r_3}\lambda_{r_4}}\times
$$
\begin{equation}
\label{12346}
\times 
\left(\sum\limits_{(r_1,r_2,r_3)}
J_{(1111)T,t}^{(r_1r_2r_3r_4)q}-
{\bf 1}_{\{r_1=r_2\}}I_{(01)T,t}^{(r_3r_4)q}-
{\bf 1}_{\{r_1=r_3\}}I_{(01)T,t}^{(r_2r_4)q}-
{\bf 1}_{\{r_2=r_3\}}I_{(01)T,t}^{(r_1r_4)q}\right),
\end{equation}

\noindent
where 
the approximations $J_{(1111)T,t}^{(r_1r_2r_3r_4)q},$
$I_{(01)T,t}^{(r_1r_2)q}$
are based on Theorem 1.1 and Le\-gen\-dre polynomials (see (\ref{kol})
and (\ref{kol1})).

For example, from (\ref{kol}) we have (here we use the notation 
$I_{(01)T,t}^{(r_1 r_2)}$ 
from the formula (\ref{kol}))
$$
I_{(01)T,t}^{(r_1 r_2)q}=
-\frac{T-t}{2}
J_{(11)T,t}^{(r_1 r_2)q}
-\frac{(T-t)^2}{4}\Biggl(
\frac{1}{\sqrt{3}}\zeta_0^{(r_1)}\zeta_1^{(r_2)}+\Biggr.
$$
\begin{equation}
\label{vini0}
~~~~~~~ +\Biggl.\sum_{i=0}^{q}\Biggl(
\frac{(i+2)\zeta_i^{(r_1)}\zeta_{i+2}^{(r_2)}
-(i+1)\zeta_{i+2}^{(r_1)}\zeta_{i}^{(r_2)}}
{\sqrt{(2i+1)(2i+5)}(2i+3)}-
\frac{\zeta_i^{(r_1)}\zeta_{i}^{(r_2)}}{(2i-1)(2i+3)}\Biggr)\Biggr),
\end{equation}

\begin{equation}
\label{vini}
J_{(11)T,t}^{(r_1 r_2)q}=
\frac{T-t}{2}\left(\zeta_0^{(r_1)}\zeta_0^{(r_2)}+\sum_{i=1}^{q}
\frac{1}{\sqrt{4i^2-1}}\biggl(
\zeta_{i-1}^{(r_1)}\zeta_{i}^{(r_2)}-
\zeta_i^{(r_1)}\zeta_{i-1}^{(r_2)}\biggr)-{\bf 1}_{\{r_1=r_2\}}
\right),
\end{equation}

\noindent
where notations are the same as in Theorem 1.1.
For $r_1\ne r_2$ we get (see (\ref{fff09xxx}))
$$
{\sf M}\left\{\left(I_{(01)T,t}^{(r_1 r_2)}-
I_{(01)T,t}^{(r_1 r_2)q}\right)^2\right\}=
\frac{(T-t)^4}{16}\left(\frac{5}{9}-
2\sum_{i=2}^q\frac{1}{4i^2-1}-
\right.
$$
\begin{equation}
\label{987}
~~~~~~~ \left.-\sum_{i=1}^q
\frac{1}{(2i-1)^2(2i+3)^2}-
\sum_{i=0}^q\frac{(i+2)^2+(i+1)^2}{(2i+1)(2i+5)(2i+3)^2}
\right).
\end{equation}

\vspace{2mm}

From (\ref{z1}) and (\ref{987}) we obtain
$$
{\sf M}\left\{\left(I_{(01)T,t}^{(r_1 r_2)}-
I_{(01)T,t}^{(r_1 r_2)q}\right)^2\right\}\le
\frac{(T-t)^4}{8}\left(\frac{5}{9}-
2\sum_{i=2}^q\frac{1}{4i^2-1}-
\right.
$$
$$
\left.-\sum_{i=1}^q
\frac{1}{(2i-1)^2(2i+3)^2}-
\sum_{i=0}^q\frac{(i+2)^2+(i+1)^2}{(2i+1)(2i+5)(2i+3)^2}
\right),
$$

\noindent
where $r_1,r_2=1,\ldots,M.$

From (\ref{jjjk2}) and (\ref{12346})
it follows that
$$
I_3[B(Z)]_{T,t}^{M}-
I_3[B(Z)]_{T,t}^{M,q}=
$$
$$
=
\sum_{r_1,r_2,r_3,r_4\in J_M}
B'''(Z)\left(B(Z)e_{r_1}, B(Z)e_{r_2},B(Z)e_{r_3}\right)e_{r_4}
\sqrt{\lambda_{r_1}\lambda_{r_2}\lambda_{r_3}\lambda_{r_4}}\times
$$

\vspace{-8mm}
$$
\times 
\Biggl(\sum\limits_{(r_1,r_2,r_3)}
\left(J_{(1111)T,t}^{(r_1r_2r_3r_4)}-J_{(1111)T,t}^{(r_1r_2r_3r_4)q}\right)-
{\bf 1}_{\{r_1=r_2\}}\left(I_{(01)T,t}^{(r_3r_4)}-I_{(01)T,t}^{(r_3r_4)q}
\right)-\Biggr.
$$

\vspace{-11mm}
\begin{equation}
\label{12341}
~~~~~\Biggl.-{\bf 1}_{\{r_1=r_3\}}
\left(I_{(01)T,t}^{(r_2r_4)}-I_{(01)T,t}^{(r_2r_4)q}
\right)-
{\bf 1}_{\{r_2=r_3\}}\left(I_{(01)T,t}^{(r_1r_4)}-
I_{(01)T,t}^{(r_1r_4)q}\right)\Biggr)\ \ \ \hbox{w.~p.~1.}
\end{equation}

Repeating with an insignificant 
modification the proof 
of Theorem 7.1 for the cases $k=2$ and $k=4$, we obtain
$$
\hspace{-85mm}{\sf M}\Biggl\{\biggl\Vert
I_3[B(Z)]_{T,t}^{M}-I_3[B(Z)]_{T,t}^{M,q}\biggr\Vert_H^2\Biggr\}\le
$$
$$
~~~~~~~~~~~~~~\le C\left({\rm tr}\ Q\right)^4
\left(6^2(4!)^2
\left(\frac{(T-t)^{4}}{24}-\sum_{j_1,j_2,j_3,j_4=0}^{q}
C_{j_4j_3j_2j_1}^2\right)+3^2(2!)^2E_q\right),
$$

\vspace{1mm}
\noindent
where 
$E_q$ is the right-hand side  of (\ref{987}) and

\vspace{-2mm}
\begin{equation}
\label{tq1}
~~~~~~~C_{j_4j_3j_2j_1}
=\frac{\sqrt{(2j_1+1)(2j_2+1)(2j_3+1)(2j_4+1)}}{16}(T-t)^{2}\bar
C_{j_4j_3j_2j_1},
\end{equation}
$$
\bar C_{j_4j_3j_2j_1}=\int\limits_{-1}^{1}P_{j_4}(u)
\int\limits_{-1}^{u}P_{j_3}(z)
\int\limits_{-1}^{z}P_{j_2}(y)
\int\limits_{-1}^{y}
P_{j_1}(x)dx dy dz du,
$$

\vspace{2mm}
\noindent
where $P_j(x)$ is the Legendre polynomial.

Consider the stochastic integral $I_4[B(Z)]_{T,t}^M.$
Let Conditions 1, 2 of Theorem 7.1 be fulfilled.
Let 
$B'(v)(B''(v)(B(v),B(v)))$
be a 4-linear Hilbert--Schmidt operator mapping from
$U_0\times U_0\times U_0\times U_0$ to $H$ for all $v\in H$.

We have (see (\ref{xx605}))

\vspace{-5mm}
$$
I_4[B(Z)]_{T,t}^M=
\sum_{r_1,r_2,r_3,r_4\in J_M}
B'(Z)\left(B''(Z)
\left(B(Z)e_{r_1}, B(Z)e_{r_2}\right)e_{r_3}\right)e_{r_4}
\times
$$
\begin{equation}
\label{jjjk5}
~~~~~\times
\sqrt{\lambda_{r_1}\lambda_{r_2}\lambda_{r_3}\lambda_{r_4}}
\int\limits_t^T
\int\limits_t^s \left(
\int\limits_t^{\tau} d{\bf w}_{u}^{(r_1)}
\int\limits_t^{\tau} d{\bf w}_{u}^{(r_2)}\right)
d{\bf w}_{\tau}^{(r_3)}
d{\bf w}_s^{(r_4)}\ \ \ \hbox{w.~p.~1.}
\end{equation}

\vspace{4mm}

From (\ref{jjj2sss}) and (\ref{jjjk5}) we obtain

\vspace{-5mm}
$$
I_4[B(Z)]_{T,t}^M=
\sum_{r_1,r_2,r_3,r_4\in J_M}
B'(Z)\left(B''(Z)
\left(B(Z)e_{r_1}, B(Z)e_{r_2}\right)e_{r_3}\right)e_{r_4}
\times
$$
\begin{equation}
\label{jjjk6}
~~~~~~\times 
\sqrt{\lambda_{r_1}\lambda_{r_2}\lambda_{r_3}\lambda_{r_4}}
\left(J_{(1111)T,t}^{(r_1r_2r_3r_4)}+
J_{(1111)T,t}^{(r_2r_1r_3r_4)}-
{\bf 1}_{\{r_1=r_2\}}I_{(10)T,t}^{(r_3r_4)}\right)\ \ \ \hbox{w.~p.~1,}
\end{equation}
\noindent

\vspace{4mm}
\noindent
where 
\begin{equation}
\label{ros2}
I_{(10)T,t}^{(r_3r_4)}=
\int\limits_t^T
\int\limits_t^s (t-\tau)d{\bf w}_{\tau}^{(r_3)}   
d{\bf w}_s^{(r_4)}.
\end{equation}

Denote by
$I_4[B(Z)]_{T,t}^{M,q}$ the approximation of the iterated stochastic integral
(\ref{jjjk6}), 
which has the following form

\vspace{-5mm}
$$
I_4[B(Z)]_{T,t}^{M,q}=
\sum_{r_1,r_2,r_3,r_4\in J_M}
B'(Z)\left(B''(Z)
\left(B(Z)e_{r_1}, B(Z)e_{r_2}\right)e_{r_3}\right)e_{r_4}
\times
$$
\begin{equation}
\label{jjjk7}
~~~~~~\times
\sqrt{\lambda_{r_1}\lambda_{r_2}\lambda_{r_3}\lambda_{r_4}}
\left(J_{(1111)T,t}^{(r_1r_2r_3r_4)q}+
J_{(1111)T,t}^{(r_2r_1r_3r_4)q}
-{\bf 1}_{\{r_1=r_2\}}I_{(10)T,t}^{(r_3r_4)q}\right)\ \ \ \hbox{w.~p.~1,}
\end{equation}

\vspace{2mm}
\noindent
where 
the approximations $J_{(1111)T,t}^{(r_1r_2r_3r_4)q},$
$I_{(10)T,t}^{(r_1r_2)q}$
are based on Theorem 1.1 and Legendre polynomials.

For example, from (\ref{kola}) we have (here we use the notation 
$I_{(10)T,t}^{(r_1 r_2)}$ 
from the formula (\ref{kola}))
$$
I_{(10)T,t}^{(r_1 r_2)q}=
-\frac{T-t}{2}J_{(11)T,t}^{(r_1 r_2)q}
-\frac{(T-t)^2}{4}\Biggl(
\frac{1}{\sqrt{3}}\zeta_0^{(r_2)}\zeta_1^{(r_1)}+\Biggr.
$$
\begin{equation}
\label{4006x}
~~~~~~ +\Biggl.\sum_{i=0}^{q}\Biggl(
\frac{(i+1)\zeta_{i+2}^{(r_2)}\zeta_{i}^{(r_1)}
-(i+2)\zeta_{i}^{(r_2)}\zeta_{i+2}^{(r_1)}}
{\sqrt{(2i+1)(2i+5)}(2i+3)}+
\frac{\zeta_i^{(r_1)}\zeta_{i}^{(r_2)}}{(2i-1)(2i+3)}\Biggr)\Biggr),
\end{equation}

\vspace{4mm}
\noindent
where the approximation $J_{(11)T,t}^{(r_1 r_2)q}$
is defined by (\ref{vini}).

Moreover, 
\begin{equation}
\label{roza1}
{\sf M}\left\{\left(I_{(10)T,t}^{(r_1 r_2)}-
I_{(10)T,t}^{(r_1 r_2)q}\right)^2\right\}=E_q\ \ \ (r_1\ne r_2),
\end{equation}

\noindent
where 
$E_q$ is the right-hand side of (\ref{987}) 
(see (\ref{fff09xxx})).

From (\ref{jjjk6}), (\ref{jjjk7})
we have
$$
I_4[B(Z)]_{T,t}^{M}-
I_4[B(Z)]_{T,t}^{M,q}
=
$$
$$
=\sum_{r_1,r_2,r_3,r_4\in J_M}
B'(Z)\left(B''(Z)
\left(B(Z)e_{r_1}, B(Z)e_{r_2}\right)e_{r_3}\right)e_{r_4}
\sqrt{\lambda_{r_1}\lambda_{r_2}\lambda_{r_3}\lambda_{r_4}}\times
$$
$$
\times 
\Biggl(\left(J_{(1111)T,t}^{(r_1r_2r_3r_4)}
-J_{(1111)T,t}^{(r_1r_2r_3r_4)q}\right)
+
\left(J_{(1111)T,t}^{(r_2r_1r_3r_4)}
-J_{(1111)T,t}^{(r_2r_1r_3r_4)q}\right)-\Biggr.
$$
$$
\Biggl.-
{\bf 1}_{\{r_1=r_2\}}\left(I_{(10)T,t}^{(r_3r_4)}
-I_{(10)T,t}^{(r_3r_4)q}\right)\Biggr)\ \ \ \hbox{w.~p.~1.}
$$

Repeating with an insignificant 
modification the proof 
of Theorem 7.1 for the cases $k=2$ and $k=4$, we obtain
$$
\hspace{-85mm}{\sf M}\Biggl\{\biggl\Vert
I_4[B(Z)]_{T,t}^{M}-I_4[B(Z)]_{T,t}^{M,q}\biggr\Vert_H^2\Biggr\}\le
$$
$$
~~~~~~~~~~~~~~\le
C
\left({\rm tr}\ Q\right)^4
\Biggl(2^2(4!)^2
\Biggl(\frac{(T-t)^{4}}{24}-\sum_{j_1,j_2,j_3,j_4=0}^{q}
C_{j_4j_3j_2j_1}^2\Biggr) +(2!)^2E_q\Biggr),
$$

\vspace{2mm}
\noindent
where 
$E_q$ is the  right-hand side of (\ref{987}) and
$C_{j_4j_3j_2j_1}$ is defined by (\ref{tq1}).

Consider the stochastic integral $I_5[B(Z)]_{T,t}^M.$
Let Conditions 1, 2 of Theorem 7.1 be fulfilled.
Let 
$B''(v)(B(v),B'(v)(B(v)))$
be a 4-linear Hilbert--Schmidt operator mapping from
$U_0\times U_0\times U_0\times U_0$ to $H$ for all $v\in H$.

We have (see (\ref{xx605}))

\vspace{-5mm}
$$
I_5[B(Z)]_{T,t}^M=
\hspace{-1mm}\sum_{r_1,r_2,r_3,r_4\in J_M}
B''(Z)(B(Z)e_{r_3},
B'(Z)(B(Z)e_{r_2})e_{r_1})e_{r_4}
\sqrt{\lambda_{r_1}\lambda_{r_2}\lambda_{r_3}\lambda_{r_4}}\times
$$
\begin{equation}
\label{jjj500}
\times\int\limits_t^T\left(
\int\limits_t^s d{\bf w}_{\tau}^{(r_3)}
\int\limits_t^{s} 
\int\limits_t^{\tau} d{\bf w}_{u}^{(r_2)}
d{\bf w}_{\tau}^{(r_1)}
\right)
d{\bf w}_s^{(r_4)}\ \ \ \hbox{w.~p.~1.}
\end{equation}

\vspace{1mm}

Using the theorem on replacement of the integration order in iterated 
It\^{o} stochastic integrals (see Theorem 3.1 and Example 3.1) or 
the It\^{o} formula, we obtain
$$
\int\limits_t^T\left(
\int\limits_t^s d{\bf w}_{\tau}^{(r_3)}
\int\limits_t^{s} 
\int\limits_t^{\tau} d{\bf w}_{u}^{(r_2)}
d{\bf w}_{\tau}^{(r_1)}
\right)
d{\bf w}_s^{(r_4)}
=
$$

\vspace{-1mm}
$$
=J_{(1111)T,t}^{(r_2r_1r_3r_4)}+
J_{(1111)T,t}^{(r_2r_3r_1r_4)}+
J_{(1111)T,t}^{(r_3r_2r_1r_4)}+
$$

\vspace{-3mm}
\begin{equation}
\label{jjj501}
+
{\bf 1}_{\{r_1=r_3\}}\left(I_{(10)T,t}^{(r_2r_4)}-
I_{(01)T,t}^{(r_2r_4)}\right)-
{\bf 1}_{\{r_2=r_3\}}I_{(10)T,t}^{(r_1r_4)}\ \ \ \hbox{w.~p.~1,}
\end{equation}

\vspace{3mm}
\noindent
where 
we use the notations from Sect.~7.2 (see (\ref{x112}))
and $I_{(01)T,t}^{(r_1r_2)}$, $I_{(10)T,t}^{(r_1r_2)}$
are defined by (\ref{ros1}), (\ref{ros2}).

After substituting (\ref{jjj501}) into (\ref{jjj500}), we obtain

\vspace{-4mm}
$$
I_5[B(Z)]_{T,t}^M=
\sum_{r_1,r_2,r_3,r_4\in J_M}
B''(Z)(B(Z)e_{r_3},
B'(Z)(B(Z)e_{r_2})e_{r_1})e_{r_4}\times
$$
$$
\times\sqrt{\lambda_{r_1}\lambda_{r_2}\lambda_{r_3}\lambda_{r_4}}
\Biggl(J_{(1111)T,t}^{(r_2r_1r_3r_4)}+
J_{(1111)T,t}^{(r_2r_3r_1r_4)}+
J_{(1111)T,t}^{(r_3r_2r_1r_4)}+\Biggr.
$$
\begin{equation}
\label{ros3}
~~~~~~~~\Biggl.+
{\bf 1}_{\{r_1=r_3\}}\left(I_{(10)T,t}^{(r_2r_4)}-
I_{(01)T,t}^{(r_2r_4)}\right)-
{\bf 1}_{\{r_2=r_3\}}I_{(10)T,t}^{(r_1r_4)}\Biggr)\ \ \ \hbox{w.~p.~1.}
\end{equation}

\vspace{2mm}

Denote by
$I_5[B(Z)]_{T,t}^{M,q}$ the approximation of the iterated stochastic integral
(\ref{ros3}), 
which has the following form
$$
I_5[B(Z)]_{T,t}^{M,q}=
\sum_{r_1,r_2,r_3,r_4\in J_M}
B''(Z)(B(Z)e_{r_3},
B'(Z)(B(Z)e_{r_2})e_{r_1})e_{r_4}\times
$$
$$
\times
\sqrt{\lambda_{r_1}\lambda_{r_2}\lambda_{r_3}\lambda_{r_4}}
\Biggl(J_{(1111)T,t}^{(r_2r_1r_3r_4)q}+
J_{(1111)T,t}^{(r_2r_3r_1r_4)q}+
J_{(1111)T,t}^{(r_3r_2r_1r_4)q}+\Biggr.
$$
\begin{equation}
\label{ros4}
~~~~~~~~\Biggl.+
{\bf 1}_{\{r_1=r_3\}}\left(I_{(10)T,t}^{(r_2r_4)q}-
I_{(01)T,t}^{(r_2r_4)q}\right)-
{\bf 1}_{\{r_2=r_3\}}I_{(10)T,t}^{(r_1r_4)q}\Biggr)\ \ \ \hbox{w.~p.~1,}
\end{equation}

\vspace{2mm}
\noindent
where 
the approximations $J_{(1111)T,t}^{(r_1r_2r_3r_4)q},$
$I_{(01)T,t}^{(r_1r_2)q},$ and $I_{(10)T,t}^{(r_1r_2)q}$
are based on Theorem 1.1 and Legendre polynomials.

From (\ref{ros3}), (\ref{ros4}) it follows that
$$
I_5[B(Z)]_{T,t}^{M}-I_5[B(Z)]_{T,t}^{M,q}=
$$
$$
=
\sum_{r_1,r_2,r_3,r_4\in J_M}
B''(Z)(B(Z)e_{r_3},
B'(Z)(B(Z)e_{r_2})e_{r_1})e_{r_4}
\sqrt{\lambda_{r_1}\lambda_{r_2}\lambda_{r_3}\lambda_{r_4}}\times
$$
$$
\times\hspace{-0.5mm}\Biggl(\hspace{-0.7mm}\left(\hspace{-0.4mm}
J_{(1111)T,t}^{(r_2r_1r_3r_4)}\hspace{-0.2mm}-\hspace{-0.2mm}
J_{(1111)T,t}^{(r_2r_1r_3r_4)q}\right)+
\left(\hspace{-0.4mm}
J_{(1111)T,t}^{(r_2r_3r_1r_4)}\hspace{-0.2mm}-
\hspace{-0.2mm}J_{(1111)T,t}^{(r_2r_3r_1r_4)q}\right)+
\left(\hspace{-0.4mm}
J_{(1111)T,t}^{(r_3r_2r_1r_4)}\hspace{-0.2mm}-\hspace{-0.2mm}
J_{(1111)T,t}^{(r_3r_2r_1r_4)q}\right)\hspace{-1.2mm}+\Biggr.
$$
$$
+
{\bf 1}_{\{r_1=r_3\}}\left(\left(I_{(10)T,t}^{(r_2r_4)}-
I_{(10)T,t}^{(r_2r_4)q}\right)-
\left(I_{(01)T,t}^{(r_2r_4)}-I_{(01)T,t}^{(r_2r_4)q}\right)\right)-
$$
$$
\Biggl.-
{\bf 1}_{\{r_2=r_3\}}
\left(I_{(10)T,t}^{(r_1r_4)}-I_{(10)T,t}^{(r_1r_4)q}\right)
\hspace{-0.6mm}\Biggr)\ \ \ \hbox{w.~p.~1.}
$$

Repeating with an insignificant 
modification the proof 
of Theorem 7.1 for the cases $k=2$ and $k=4$ and taking into account
(\ref{roza1}),
we obtain
$$
\hspace{-85mm}{\sf M}\Biggl\{\biggl\Vert
I_5[B(Z)]_{T,t}^{M}-I_5[B(Z)]_{T,t}^{M,q}\biggr\Vert_H^2\Biggr\}\le
$$
$$
~~~~~~~~~~~\le
C
\left({\rm tr}\ Q\right)^4
\Biggl(3^2(4!)^2
\Biggl(\frac{(T-t)^{4}}{24}-\sum_{j_1,j_2,j_3,j_4=0}^{q}
C_{j_4j_3j_2j_1}^2\Biggr) +3^2(2!)^2E_q\Biggr),
$$

\vspace{2mm}
\noindent
where 
$E_q$ is the right-hand side of (\ref{987}) and
$C_{j_4j_3j_2j_1}$ is defined by (\ref{tq1}).

Consider the stochastic integral $I_6[B(Z),F(Z)]_{T,t}^M.$
Let Conditions 1, 2 of Theorem 7.1 be fulfilled.
We have (see (\ref{xx605}))

\vspace{-6mm}
$$
I_6[B(Z),F(Z)]_{T,t}^M=
\sum_{r_1,r_2\in J_M}
F'(Z)(B'(Z)(B(Z)e_{r_1})e_{r_2})
\sqrt{\lambda_{r_1}\lambda_{r_2}}\times
$$
\begin{equation}
\label{jjjk51}
\times\int\limits_t^T
\int\limits_t^s 
\int\limits_{t}^{\tau} d{\bf w}_{u}^{(r_1)}
d{\bf w}_{\tau}^{(r_2)}
ds\ \ \ \hbox{w.~p.~1.}
\end{equation}

Using the theorem on replacement of the integration order in iterated 
It\^{o} stochastic integrals (see Theorem 3.1 and Example 3.1)
or the It\^{o} formula, we obtain
\begin{equation}
\label{ff11}
~~~~~~~~~\int\limits_t^T
\int\limits_t^s 
\int\limits_{t}^{\tau} d{\bf w}_{u}^{(r_1)}
d{\bf w}_{\tau}^{(r_2)}
ds=(T-t)J_{(11)T,t}^{(r_1r_2)}+I_{(01)T,t}^{(r_1r_2)}\ \ \ \hbox{w.~p.~1.}
\end{equation}

After substituting (\ref{ff11}) into (\ref{jjjk51}),
we have

\vspace{-4mm}
$$
I_6[B(Z),F(Z)]_{T,t}^M=
\sum_{r_1,r_2\in J_M}
F'(Z)(B'(Z)(B(Z)e_{r_1})e_{r_2})
\sqrt{\lambda_{r_1}\lambda_{r_2}}\times
$$
\begin{equation}
\label{jjjk101}
\times\left(
(T-t)J_{(11)T,t}^{(r_1r_2)}+I_{(01)T,t}^{(r_1r_2)}
\right)\ \ \ \hbox{w.~p.~1.}
\end{equation}

\vspace{3mm}

Denote by
$I_6[B(Z),F(Z)]_{T,t}^{M,q}$ the approximation of the
iterated stochastic integral
(\ref{jjjk101}), 
which has the following form

\vspace{-4mm}
$$
I_6[B(Z),F(Z)]_{T,t}^{M,q}=                     
\sum_{r_1,r_2\in J_M}
F'(Z)(B'(Z)(B(Z)e_{r_1})e_{r_2})
\sqrt{\lambda_{r_1}\lambda_{r_2}}\times
$$

\vspace{-4mm}
\begin{equation}
\label{jjjk102}
\times\left(
(T-t)J_{(11)T,t}^{(r_1r_2)q}+I_{(01)T,t}^{(r_1r_2)q}
\right),
\end{equation}

\vspace{3mm}
\noindent
where 
the approximations $I_{(01)T,t}^{(r_1r_2)q}$, $J_{(11)T,t}^{(r_1r_2)q}$
are defined by (\ref{vini0}), (\ref{vini}).

From (\ref{jjjk101}), (\ref{jjjk102}) we get

\vspace{-4mm}
$$
I_6[B(Z),F(Z)]_{T,t}^M-I_6[B(Z),F(Z)]_{T,t}^{M,q}=
$$

\vspace{-2mm}
$$
=                     
\sum_{r_1,r_2\in J_M}
F'(Z)(B'(Z)(B(Z)e_{r_1})e_{r_2})
\sqrt{\lambda_{r_1}\lambda_{r_2}}\times
$$
$$
~~\times\left(
(T-t)\left(J_{(11)T,t}^{(r_1r_2)}-J_{(11)T,t}^{(r_1r_2)q}\right)
+\left(I_{(01)T,t}^{(r_1r_2)}-I_{(01)T,t}^{(r_1r_2)q}\right)
\right)\ \ \ \hbox{w.~p.~1.}
$$

\vspace{4mm}

Repeating with an insignificant 
modification the proof 
of Theorem 7.1 for the case $k=2$, we obtain
$$
\hspace{-35mm}
{\sf M}\left\{\biggl\Vert
I_6[B(Z),F(Z)]_{T,t}^M-I_6[B(Z),F(Z)]_{T,t}^{M,q}\biggr\Vert_H^2\right\}\le
$$
$$
~~~~~~~~~~~~~~~~~~~~~~~~~~~~~~~~~~~~~~~\le
2C(2!)^2\left({\rm tr}\ Q\right)^2
\biggl((T-t)^2 G_q+E_q\biggr),
$$

\vspace{2mm}
\noindent
where $G_q$ and $E_q$ are the right-hand sides of (\ref{fff09}) and
(\ref{987}) correspondingly.

Consider the stochastic integral $I_7[B(Z),F(Z)]_{T,t}^M.$
Let Conditions 1, 2 of Theorem 7.1 be fulfilled.

Then we have (see (\ref{xx605}))
$$
I_7[B(Z),F(Z)]_{T,t}^M=
\sum_{r_1,r_2\in J_M}F''(Z)\left(B(Z)e_{r_1}, B(Z)e_{r_2}\right)
\sqrt{\lambda_{r_1}\lambda_{r_2}}\times
$$
\begin{equation}
\label{da11}
\times \int\limits_t^T
\left(\int\limits_t^s d{\bf w}_{\tau}^{(r_1)}   
\int\limits_t^s d{\bf w}_{\tau}^{(r_2)}\right)
ds\ \ \ \hbox{w.~p.~1.}
\end{equation}

\vspace{4mm}

From (\ref{da12eee}) and (\ref{ff11}) we get w.~p.~1
$$
\int\limits_t^T
\left(\int\limits_t^s d{\bf w}_{\tau}^{(r_1)}   
\int\limits_t^s d{\bf w}_{\tau}^{(r_2)}\right)ds=
$$
$$
=
\int\limits_t^T
J_{(11)s,t}^{(r_1 r_2)}ds+
\int\limits_t^T
J_{(11)s,t}^{(r_2 r_1)}ds+{\bf 1}_{\{r_1=r_2\}}\frac{(T-t)^2}{2}=
$$
$$
=(T-t)\left(J_{(11)T,t}^{(r_1 r_2)}+
J_{(11)T,t}^{(r_2 r_1)}\right)+
I_{(01)T,t}^{(r_1 r_2)}+I_{(01)T,t}^{(r_2 r_1)}+
{\bf 1}_{\{r_1=r_2\}}\frac{(T-t)^2}{2}=
$$

\vspace{-1mm}
$$
=
(T-t)\left(J_{(1)T,t}^{(r_1)}J_{(1)T,t}^{(r_2)}-
{\bf 1}_{\{r_1=r_2\}}(T-t)\right)+
$$

\vspace{-1mm}
$$
+
I_{(01)T,t}^{(r_1 r_2)}+I_{(01)T,t}^{(r_2 r_1)}+
{\bf 1}_{\{r_1=r_2\}}\frac{(T-t)^2}{2}=
$$

\vspace{-4mm}
\begin{equation}
\label{da13}
~~~~~~~~=
(T-t)J_{(1)T,t}^{(r_1)}J_{(1)T,t}^{(r_2)}+
I_{(01)T,t}^{(r_1 r_2)}+I_{(01)T,t}^{(r_2 r_1)}
-
{\bf 1}_{\{r_1=r_2\}}\frac{(T-t)^2}{2}.
\end{equation}

\vspace{5mm}

After substituting (\ref{da13}) into (\ref{da11}), we obtain

\vspace{-4mm}
$$
I_7[B(Z),F(Z)]_{T,t}^M=
\sum_{r_1,r_2\in J_M}F''(Z)\left(B(Z)e_{r_1}, B(Z)e_{r_2}\right)
\sqrt{\lambda_{r_1}\lambda_{r_2}}\times
$$
\begin{equation}
\label{da14}
~~\times
\left(
(T-t)J_{(1)T,t}^{(r_1)}J_{(1)T,t}^{(r_2)}+
I_{(01)T,t}^{(r_1 r_2)}+I_{(01)T,t}^{(r_2 r_1)}-
{\bf 1}_{\{r_1=r_2\}}\frac{(T-t)^2}{2}\right)\ \ \ \hbox{w.~p.~1.}
\end{equation}

\vspace{3mm}

Denote by
$I_7[B(Z),F(Z)]_{T,t}^{M,q}$ the approximation of 
the
iterated stochastic integral
(\ref{da14}),
which has the following form

\newpage
\noindent
$$
I_7[B(Z),F(Z)]_{T,t}^{M,q}=                     
\sum_{r_1,r_2\in J_M}F''(Z)\left(B(Z)e_{r_1}, B(Z)e_{r_2}\right)
\sqrt{\lambda_{r_1}\lambda_{r_2}}\times
$$
\begin{equation}
\label{da20}
~~~~~~\times
\left(
(T-t)J_{(1)T,t}^{(r_1)}J_{(1)T,t}^{(r_2)}+
I_{(01)T,t}^{(r_1 r_2)q}+I_{(01)T,t}^{(r_2 r_1)q}-
{\bf 1}_{\{r_1=r_2\}}\frac{(T-t)^2}{2}\right),
\end{equation}

\vspace{3mm}
\noindent
where 
the approximation $I_{(01)T,t}^{(r_1r_2)q}$
is defined by (\ref{vini0}).

From (\ref{da14}), (\ref{da20}) it follows that 

\vspace{-4mm}
$$
I_7[B(Z),F(Z)]_{T,t}^M-I_7[B(Z),F(Z)]_{T,t}^{M,q}=                    
\sum_{r_1,r_2\in J_M}F''(Z)\left(B(Z)e_{r_1}, B(Z)e_{r_2}\right)\times
$$
$$
\times\sqrt{\lambda_{r_1}\lambda_{r_2}}
\left(
\left(I_{(01)T,t}^{(r_1 r_2)}-I_{(01)T,t}^{(r_1 r_2)q}\right)
+\left(I_{(01)T,t}^{(r_2 r_1)}-I_{(01)T,t}^{(r_2 r_1)q}\right)
\right)\ \ \ \hbox{w.~p.~1.}
$$

\vspace{4mm}

Repeating with an insignificant 
modification the proof 
of Theorem 7.1 for the case $k=2$, we obtain
$$
{\sf M}\left\{\biggl\Vert
I_7[B(Z),F(Z)]_{T,t}^M-I_7[B(Z),F(Z)]_{T,t}^{M,q}\biggr\Vert_H^2\right\}\le
4C (2!)^2\left({\rm tr}\ Q\right)^2
E_q,
$$

\vspace{2mm}
\noindent
where $E_q$ is the right-hand side of (\ref{987}).

Consider the stochastic integral $I_8[B(Z),F(Z)]_{T,t}^M.$
Let Conditions 1, 2 of Theorem 7.1 be fulfilled.

Then we have w.~p.~1 (see (\ref{xx605}))

\vspace{-4mm}
\begin{equation}
\label{da110}
I_8[B(Z),F(Z)]_{T,t}^M=-
\sum_{r_1,r_2\in J_M}B''(Z)\left(F(Z), B(Z)e_{r_1}\right)e_{r_2}
\sqrt{\lambda_{r_1}\lambda_{r_2}}I_{(01)T,t}^{(r_1r_2)}.
\end{equation}

\vspace{2mm}

Denote by
$I_8[B(Z),F(Z)]_{T,t}^{M,q}$ the approximation of the
iterated stochastic integral
(\ref{da110}),
which has the following form

\vspace{-4mm}
\begin{equation}
\label{da200}
I_8[B(Z),F(Z)]_{T,t}^{M,q}=                     
-\sum_{r_1,r_2\in J_M}B''(Z)\left(F(Z), B(Z)e_{r_1}\right)e_{r_2}
\sqrt{\lambda_{r_1}\lambda_{r_2}}I_{(01)T,t}^{(r_1r_2)q},
\end{equation}

\vspace{1mm}
\noindent
where 
the approximation $I_{(01)T,t}^{(r_1r_2)q}$
is defined by (\ref{vini0}).

From (\ref{da110}), (\ref{da200}) we get

\newpage
\noindent
$$
I_8[B(Z),F(Z)]_{T,t}^M-I_8[B(Z),F(Z)]_{T,t}^{M,q}=
$$

\vspace{-5mm}
$$
=-\sum_{r_1,r_2\in J_M}B''(Z)\left(F(Z), B(Z)e_{r_1}\right)e_{r_2}
\sqrt{\lambda_{r_1}\lambda_{r_2}}
\left(I_{(01)T,t}^{(r_1r_2)}-
I_{(01)T,t}^{(r_1r_2)q}\right)\ \ \ \hbox{w.~p.~1.}
$$

Repeating with an insignificant 
modification the proof 
of Theorem 7.1 for the case $k=2$, we obtain
$$
{\sf M}\left\{\biggl\Vert
I_8[B(Z),F(Z)]_{T,t}^M-I_8[B(Z),F(Z)]_{T,t}^{M,q}\biggr\Vert_H^2\right\}\le
C (2!)^2\left({\rm tr}\ Q\right)^2
E_q,
$$

\noindent
where $E_q$ is the right-hand side of (\ref{987}).

\section{Approximation of Iterated Stochastic Integrals 
of Mil\-ti\-pli\-ci\-ties 1 to 3 
with Respect to the Infinite-Di\-men\-si\-onal $Q$-Wiener Process}

This section is devoted to 
the application of Theorem 1.1 and multiple Fourier--Legendre series for  
the approximation of iterated stochastic integrals
of multiplicities 1 to 3 with respect to the infinite-dimensional
$Q$-Wiener process. These iterated stochastic integrals are part of
the exponential Milstein and 
Wagner--Platen numerical methods for semilinear SPDEs
with nonlinear multiplicative trace class 
noise (see Sect.~7.2). 
Theorem 7.3 (see below) on the mean-square convergence
of approximations of iterated stochastic integrals
of multiplicities 2 and 3 with respect to the infinite-dimensional
$Q$-Wiener process is formulated and proved.
The results of this section can be applied to the
implementation of high-order strong numerical methods 
for non-commutative semilinear SPDEs
with nonlinear multiplicative trace class 
noise.

\subsection{Formulas for the Numerical Modeling of 
Iterated Stochastic Integrals 
of Miltiplicities 1 to 3 
with Respect to the Infinite-Dimensional $Q$-Wiener Process
Based on Theorem 1.1 and Legendre Polynomials}

This section is devoted to the formulas 
for numerical modeling of iterated
stochastic integrals 
from the Milstein type scheme (\ref{fff100})
and the Wagner--Platen type scheme (\ref{fffguk}) for non-commutative 
semilinear SPDEs. These integrals have the following form 
(below we introduce new notations for the stochastic integrals 
(\ref{ret4})-(\ref{ccc3}) and their approximations)
\begin{equation}
\label{ret1}
J_1[B(Z)]_{T,t}=\int\limits_{t}^{T}B(Z)d{\bf W}_{t_1},
\end{equation}
\begin{equation}
\label{ret2}
~~~~~~J_2[B(Z)]_{T,t}=A\left(\int\limits_{t}^{T}
\int\limits_{t}^{t_2}
B(Z)d{\bf W}_{t_1}dt_2-\frac{(T-t)}{2}
\int\limits_{t}^{T}B(Z)d{\bf W}_{t_1}\right),
\end{equation}
$$
J_3[B(Z),F(Z)]_{T,t}
=(T-t)\int\limits_{t}^{T}B'(Z)\biggl(AZ+F(Z)\biggr)d{\bf W}_{t_1}-
$$
\begin{equation}
\label{ret3}
-
\int\limits_{t}^{T}\int\limits_{t}^{t_2}
B'(Z)\biggl(AZ+F(Z)\biggr)d{\bf W}_{t_1}dt_2,
\end{equation}
\begin{equation}
\label{ret4}
J_4[B(Z),F(Z)]_{T,t}=\int\limits_{t}^{T}F'(Z)
\left(\int\limits_{t}^{t_2} B(Z) d{\bf W}_{t_1} \right) dt_2,
\end{equation}
\begin{equation}
\label{ccc1}
I_1[B(Z)]_{T,t}=\int\limits_{t}^{T}
B'(Z)\left(\int\limits_{t}^{t_2}
B(Z)d{\bf W}_{t_1}\right)d{\bf W}_{t_2},
\end{equation}

\vspace{-4mm}
\begin{equation}
\label{ccc2}
~~~~~I_2[B(Z)]_{T,t}=\int\limits_{t}^{T}B'(Z) \left(
\int\limits_{t}^{t_3}B'(Z)\left(
\int\limits_{t}^{t_2}B(Z) d{\bf W}_{t_1} 
\right) d{\bf W}_{t_2} \right) d{\bf W}_{t_3},
\end{equation}

\vspace{-4mm}
\begin{equation}
\label{ccc3}
~~I_3[B(Z)]_{T,t}=\int\limits_{t}^{T}B''(Z) \left(
\int\limits_{t}^{t_2}B(Z) d{\bf W}_{t_1},
\int\limits_{t}^{t_2}B(Z) d{\bf W}_{t_1} 
\right) d{\bf W}_{t_2},
\end{equation}

\vspace{1mm}
\noindent
where $Z: \Omega \rightarrow H$ is 
an ${\bf F}_t/{\mathcal{B}}(H)$-measurable mapping,
$0\le t<T\le \bar T.$

Note that according to
(\ref{ggg0})--(\ref{ggg3}), (\ref{4001}), (\ref{y3a}), and (\ref{y4a}),
we can write the following relatively simple formulas
for numerical modeling \cite{OK}, \cite{arxiv-21}

\vspace{-2mm}
$$
J_1[B(Z)]_{T,t}^M=\int\limits_{t}^{T}B(Z)d{\bf W}_s^M=
$$
$$
=
(T-t)^{1/2}\sum\limits_{r_1\in J_M}
B(Z)e_{r_1}\sqrt{\lambda_{r_1}}
\zeta_0^{(r_1)},
$$

$$
J_2[B(Z)]_{T,t}^M=
A\left(\int\limits_{t}^{T}
\int\limits_{t}^{t_2}
B(Z)d{\bf W}_{t_1}^M dt_2-\frac{(T-t)}{2}
\int\limits_{t}^{T}B(Z)d{\bf W}_{t_1}^M\right)=
$$
\begin{equation}
\label{sd11}
=
-\frac{(T-t)^{3/2}}{2\sqrt{3}}\sum\limits_{r_1\in J_M}
AB(Z)e_{r_1}\sqrt{\lambda_{r_1}}
\zeta_1^{(r_1)},
\end{equation}

\vspace{4mm}
$$
J_3[B(Z),F(Z)]_{T,t}^M=
(T-t)\int\limits_{t}^{T}B'(Z)\biggl(AZ+F(Z)\biggr)d{\bf W}_{t_1}^M-
$$
$$
-
\int\limits_{t}^{T}\int\limits_{t}^{t_2}
B'(Z)\biggl(AZ+F(Z)\biggr)d{\bf W}_{t_1}^M dt_2=
$$
\begin{equation}
\label{sd12}
~~~~~~ =
\frac{(T-t)^{3/2}}{2}\sum\limits_{r_1\in J_M}
B'(Z)\biggl(AZ+F(Z)\biggr)
e_{r_1}\sqrt{\lambda_{r_1}}\left(\zeta_0^{(r_1)}+
\frac{1}{\sqrt{3}}\zeta_1^{(r_1)}\right),
\end{equation}

\vspace{5mm}
$$
J_4[B(Z),F(Z)]_{T,t}^M=\int\limits_{t}^{T}F'(Z)
\left(\int\limits_{t}^{t_2} B(Z) d{\bf W}_{t_1}^M \right) dt_2=
$$
\begin{equation}
\label{sd13}
~~~~~~~ =\frac{(T-t)^{3/2}}{2}\sum\limits_{r_1\in J_M}
F'(Z)B(Z)e_{r_1}\sqrt{\lambda_{r_1}}
\left(\zeta_0^{(r_1)}-
\frac{1}{\sqrt{3}}\zeta_1^{(r_1)}\right),
\end{equation}

\vspace{5mm}
\noindent
where $\zeta_0^{(r_1)},$ $\zeta_1^{(r_1)}$ $(r_1\in J_M)$ are 
independent standard Gaussian random variables.

Further, consider the stochastic
integrals (\ref{ccc1})--(\ref{ccc3}), which are more
complicate, in detail.

Let $I_1[B(Z)]_{T,t}^M,$ $I_2[B(Z)]_{T,t}^M,$ 
$I_3[B(Z)]_{T,t}^M$ be approximations
of the stochastic integrals (\ref{ccc1})--(\ref{ccc3}),
which have the following form
(see (\ref{uspel1}), (\ref{uspel2}), and (\ref{da1quka}))

\vspace{-3mm}
$$
I_1[B(Z)]_{T,t}^M=\int\limits_{t}^{T}
B'(Z)\left(\int\limits_{t}^{t_2}
B(Z)d{\bf W}_{t_1}^M\right)d{\bf W}_{t_2}^M=
$$

\newpage
\noindent
\begin{equation}
\label{vvv1}
=\sum\limits_{r_1,r_2\in J_M}
B'(Z)\left(B(Z)e_{r_1}\right)
e_{r_2}\sqrt{\lambda_{r_1}\lambda_{r_2}}J_{(11)T,t}^{(r_1 r_2)},
\end{equation}

$$
I_2[B(Z)]_{T,t}^M=
\int\limits_{t}^{T}B'(Z) \left(
\int\limits_{t}^{t_3}B'(Z)\left(
\int\limits_{t}^{t_2}B(Z) d{\bf W}_{t_1}^M 
\right) d{\bf W}_{t_2}^M \right) d{\bf W}_{t_3}^M=
$$

\vspace{-2mm}
\begin{equation}
\label{da4xx}
~~~~=\sum_{r_1,r_2,r_3\in J_M}B'(Z)\left(B'(Z)\left(B(Z)e_{r_1}\right)
e_{r_2}\right)e_{r_3}
\sqrt{\lambda_{r_1}\lambda_{r_2}\lambda_{r_3}}
J_{(111)T,t}^{(r_1 r_2 r_3)},
\end{equation}

\vspace{1mm}

$$
I_3[B(Z)]_{T,t}^M=\int\limits_{t}^{T}B''(Z) \left(
\int\limits_{t}^{t_2}B(Z) d{\bf W}_{t_1}^M,
\int\limits_{t}^{t_2}B(Z) d{\bf W}_{t_1}^M 
\right) d{\bf W}_{t_2}^M=
$$

\vspace{-2mm}
$$
=\sum_{r_1,r_2,r_3\in J_M}B''(Z)\left(B(Z)e_{r_1}, B(Z)e_{r_2}\right)e_{r_3}
\sqrt{\lambda_{r_1}\lambda_{r_2}\lambda_{r_3}} \times
$$

\vspace{-2mm}
\begin{equation}
\label{da4x}
\times
\left(J_{(111)T,t}^{(r_1 r_2 r_3)}+
J_{(111)T,t}^{(r_2 r_1 r_3)}+
{\bf 1}_{\{r_1=r_2\}}J_{(01)T,t}^{(0 r_3)}\right).
\end{equation}

\vspace{6mm}

Let $I_1[B(Z)]_{T,t}^{M,q},$
$I_2[B(Z)]_{T,t}^{M.q},$ $I_3[B(Z)]_{T,t}^{M,q}$ be approximations
of the stochastic integrals (\ref{vvv1})--(\ref{da4x}),
which look as follows

$$
I_1[B(Z)]_{T,t}^{M,q}
=\sum\limits_{r_1,r_2\in J_M}
B'(Z)\left(B(Z)e_{r_1}\right)
e_{r_2}\sqrt{\lambda_{r_1}\lambda_{r_2}}J_{(11)T,t}^{(r_1 r_2)q},
$$

\begin{equation}
\label{da4xxx}
I_2[B(Z)]_{T,t}^{M,q}=
\sum_{r_1,r_2,r_3\in J_M}B'(Z)\left(B'(Z)\left(B(Z)e_{r_1}\right)
e_{r_2}\right)e_{r_3}
\sqrt{\lambda_{r_1}\lambda_{r_2}\lambda_{r_3}}
J_{(111)T,t}^{(r_1 r_2 r_3)q},
\end{equation}

\vspace{-4mm}
$$
I_3[B(Z)]_{T,t}^{M,q}=
\sum_{r_1,r_2,r_3\in J_M}B''(Z)\left(B(Z)e_{r_1}, B(Z)e_{r_2}\right)e_{r_3}
\sqrt{\lambda_{r_1}\lambda_{r_2}\lambda_{r_3}} \times
$$
$$
\times
\left(J_{(111)T,t}^{(r_1 r_2 r_3)q}+
J_{(111)T,t}^{(r_2 r_1 r_3)q}+
{\bf 1}_{\{r_1=r_2\}}J_{(01)T,t}^{(0 r_3)}\right),
$$

\vspace{4.5mm}
\noindent
where the approximations 
$J_{(11)T,t}^{(r_1 r_2)q},$ $J_{(111)T,t}^{(r_1 r_2 r_3)q},$
$J_{(111)T,t}^{(r_2 r_1 r_3)q}$ 
of the stochastic integrals (\ref{x112})
are defined by 
(\ref{y3}), (\ref{y4}) and $J_{(01)T,t}^{(0 r_3)}$
has the form (\ref{y3a}), $q\ge 1.$

\subsection{Theorem on the Mean-Square Approximation
of Iterated Sto\-chas\-tic Integrals
of Multiplicities 2 and 3
with Res\-pect to the Ininite-Dimensional 
$Q$-Wiener Process}

Recall that $L_{HS}(U_0,H)$
is a space of Hilbert--Schmidt operators mapping from $U_0$ to $H.$
Moreover, let $L_{HS}^{(2)}(U_0, H)$ and
$L_{HS}^{(3)}(U_0, H)$
be spaces of bilinear and 3-linear
Hilbert--Schmidt operators 
mapping 
from $U_0 \times U_0$ to $H$ and from $U_0 \times U_0 \times U_0$ to $H$
correspondingly.
Furthermore, let
$$
\left\Vert \cdot \right\Vert_{L_{HS}(U_0,H)},\ \ \
\left\Vert \cdot \right\Vert_{L_{HS}^{(2)}(U_0, H)},\ \ \
\left\Vert \cdot \right\Vert_{L_{HS}^{(3)}(U_0, H)}
$$ 
be operator norms in these spaces.

\vspace{2mm}

{\bf Theorem 7.3} \cite{12a}-\cite{12aa}, \cite{OK}, \cite{arxiv-21},
\cite{Kuzh-1}, \cite{new-new-1}. {\it 
Let Conditions {\rm 1, 2} of Theorem {\rm 7.1} 
be fulfilled. Furthermore, let

\vspace{-3mm}
$$
B(v)\in L_{HS}(U_0,H),\ \
B'(v)(B(v))\in L_{HS}^{(2)}(U_0, H),
$$

\vspace{-4mm}
$$
B'(v)(B'(v)(B(v))),\ \ B''(v)(B(v),B(v)) \in
L_{HS}^{(3)}(U_0, H)
$$

\vspace{2.7mm}
\noindent
for all $v\in H$ {\rm (}we suppose that Fr\^{e}chet derivatives 
$B',$ $B''$ exist{\rm ;} see Sect.~{\rm 7.1)}. 
Moreover, 
let there exists a constant $C$ such that w.~p.~{\rm 1}

\vspace{-4mm}
$$
\biggl\Vert B(Z)Q^{-\alpha}\biggr\Vert_{L_{HS}(U_0,H)}<C,\ \ \
\biggl\Vert B'(Z)(B(Z))Q^{-\alpha}\biggr\Vert_{L_{HS}^{(2)}(U_0, H)}<C,
$$
$$
\biggl\Vert B'(Z)(B'(Z)(B(Z)))Q^{-\alpha}
\biggr\Vert_{L_{HS}^{(3)}(U_0, H)}<C,
$$
$$
\biggl\Vert B''(Z)(B(Z), B(Z))Q^{-\alpha}
\biggr\Vert_{L_{HS}^{(3)}(U_0, H)}<C
$$

\vspace{3mm}
\noindent
for some $\alpha>0.$
Then

$$
{\sf M}\left\{\Biggl\Vert
I_1[B(Z)]_{T,t}-
I_1[B(Z)]_{T,t}^{M,p}\Biggr\Vert_H^2\right\}\le
$$
\begin{equation}
\label{rock6}
~~~~\le (T-t)^2\left(C_0
\left({\rm tr}\ Q\right)^2
\left(\frac{1}{2}-\sum_{j=1}^{p}
\frac{1}{4j^2-1}\right)+
K_Q\left(\sup\limits_{i\in J\backslash J_M}\lambda_i\right)^{2\alpha}
\right),
\end{equation}

\newpage
\noindent
$$
{\sf M}\left\{\Biggl\Vert
I_2[B(Z)]_{T,t}-
I_2[B(Z)]_{T,t}^{M,q}\Biggr\Vert_H^2\right\}\le
$$
\begin{equation}
\label{ppp6}
~~\le (T-t)^3\left(C_1
\left({\rm tr}\ Q\right)^3
\Biggl(\frac{1}{6}-\sum_{j_1,j_2,j_3=0}^{q}
{\hat C}_{j_3j_2j_1}^2\Biggr)+
L_Q\left(\sup\limits_{i\in J\backslash J_M}\lambda_i\right)^{2\alpha}
\right),
\end{equation}

\vspace{2mm}
$$
{\sf M}\left\{\Biggl\Vert
I_3[B(Z)]_{T,t}-
I_3[B(Z)]_{T,t}^{M,q}\Biggr\Vert_H^2\right\}\le
$$
\begin{equation}
\label{uuu1}
~~\le (T-t)^3\left(C_2
\left({\rm tr}\ Q\right)^3
\Biggl(\frac{1}{6}-\sum_{j_1,j_2,j_3=0}^{q}
{\hat C}_{j_3j_2j_1}^2\Biggr)+
M_Q\left(\sup\limits_{i\in J\backslash J_M}\lambda_i\right)^{2\alpha}
\right),
\end{equation}

\vspace{6mm}
\noindent
where $p, q\in{\bf N},$ $C_0, C_1, C_2, K_Q, L_Q, M_Q<\infty,$ and 

\vspace{2mm}
$$
{\hat C}_{j_3j_2j_1}=\frac{\sqrt{(2j_1+1)(2j_2+1)(2j_3+1)}}{8}\bar
C_{j_3j_2j_1},
$$

\vspace{-3mm}
$$
\bar C_{j_3j_2j_1}=\int\limits_{-1}^{1}P_{j_3}(z)
\int\limits_{-1}^{z}P_{j_2}(y)
\int\limits_{-1}^{y}
P_{j_1}(x)dx dy dz,
$$

\vspace{2mm}
\noindent
where $P_j(x)$ $(j=0, 1, 2,\ldots)$ is the Legendre polynomial.
}

{\bf Remark 7.4.}\ {\it Note that the estimate similar to
{\rm (\ref{rock6})} has been derived in {\rm \cite{26zzzzz}}, 
{\rm \cite{27zzzzz}} 
{\rm (}also see {\rm \cite{12zzzzz}}{\rm )}
with the 
difference connected with the first term
on the right-hand side of {\rm (\ref{rock6})}.
In {\rm \cite{27zzzzz}} the authors used the Karhunen--Lo\`{e}ve
expansion of the Brownian bridge process for the
approximation of iterated It\^{o} stochastic integrals
with respect to the finite-dimensional Wiener process
{\rm (}Milstein approach, see Sect.~{\rm 6.2)}.
In this section, we apply Theorem {\rm 1.1} and the
system of Legendre polynomials to obtain
the first term on the right-hand side of {\rm (\ref{rock6}).}}

{\bf Remark 7.5.}\ {\it If we assume that
$\lambda_i\le C' i^{-\gamma}$ $(\gamma>1, C'<\infty)$ for
$i\in J$, then the parameter 
$\alpha>0$ obviously increases with decreasing of ${\gamma}$
{\rm \cite{26zzzzz}}.}

{\bf Proof.}\ The estimate (\ref{rock6}) follows directly 
from (\ref{zzz1ee}) for $k=2$ 
(the first term on the right-hand side of (\ref{rock6}))
and Theorem 1 from \cite{27zzzzz} (the second term on the 
right-hand side of (\ref{rock6})).
Further $C_3,$ $C_4,$ $\ldots $ denote various constants.

Let us prove the estimates (\ref{ppp6}), (\ref{uuu1}).
Using Theorem 7.1, we obtain

\vspace{-4mm}
$$
{\sf M}\left\{\Biggl\Vert
I_2[B(Z)]_{T,t}-
I_2[B(Z)]_{T,t}^{M,q}\Biggr\Vert_H^2\right\}\le 
2{\sf M}\left\{\Biggl\Vert
I_2[B(Z)]_{T,t}-
I_2[B(Z)]_{T,t}^{M}\Biggr\Vert_H^2\right\}+
$$

\vspace{-3mm}
$$
+2
{\sf M}\left\{\Biggl\Vert
I_2[B(Z)]_{T,t}^M-
I_2[B(Z)]_{T,t}^{M,q}\Biggr\Vert_H^2\right\}\le
$$

$$
\hspace{-60mm}\le 
2{\sf M}\left\{\Biggl\Vert
I_2[B(Z)]_{T,t}-
I_2[B(Z)]_{T,t}^{M}\Biggr\Vert_H^2\right\}+
$$

\vspace{-3mm}
\begin{equation}
\label{ppp4}
~~~~~~~~~~~~~~~~~~~~~~~~~~~~~~~~~+C_3(T-t)^3
\left({\rm tr}\ Q\right)^3
\Biggl(\frac{1}{6}-\sum_{j_1,j_2,j_3=0}^{q}
{\hat C}_{j_3j_2j_1}^2\Biggr),
\end{equation}

\vspace{2mm}

$$
{\sf M}\left\{\Biggl\Vert
I_3[B(Z)]_{T,t}-
I_3[B(Z)]_{T,t}^{M,q}\Biggr\Vert_H^2\right\}
\le 2
{\sf M}\left\{\Biggl\Vert
I_3[B(Z)]_{T,t}-
I_3[B(Z)]_{T,t}^{M}\Biggr\Vert_H^2\right\}+
$$

\vspace{-3mm}
\begin{equation}
\label{hhh}
+2{\sf M}\left\{\Biggl\Vert
I_3[B(Z)]_{T,t}^M-
I_3[B(Z)]_{T,t}^{M,q}\Biggr\Vert_H^2\right\}.
\end{equation}

\vspace{4mm}

Repeating with an insignificant 
modification the proof 
of Theorem 7.1 for the case $k=3$, 
we have (also see Sect.~7.3.2)

$$
\hspace{-70mm}{\sf M}\left\{\biggl\Vert
I_3[B(Z)]_{T,t}^{M}-I_3[B(Z)]_{T,t}^{M,q}\biggr\Vert_H^2\right\}\le
$$
\begin{equation}
\label{hhh1}
~~~~~~~~~~~~~~~~~~~~~~~~~\le 4{\tilde C}(3!)^2
\left({\rm tr}\ Q\right)^3 (T-t)^3
\Biggl(\frac{1}{6}-\sum_{j_1,j_2,j_3=0}^{q}
{\hat C}_{j_3j_2j_1}^2\Biggr),
\end{equation}

\vspace{4mm}
\noindent
where constant ${\tilde C}$ has the same meaning 
as constant $L_k$ in Theorem 7.1 
($k$ is the multiplicity of the iterated stochastic integral).

Combining (\ref{hhh}) and (\ref{hhh1}), we obtain

\newpage
\noindent
$$
{\sf M}\left\{\Biggl\Vert
I_3[B(Z)]_{T,t}-
I_3[B(Z)]_{T,t}^{M,q}\Biggr\Vert_H^2\right\}\le 
2{\sf M}\left\{\Biggl\Vert
I_3[B(Z)]_{T,t}-
I_3[B(Z)]_{T,t}^{M}\Biggr\Vert_H^2\right\}+
$$
\begin{equation}
\label{sss1aaa}
+C_4(T-t)^3
\left({\rm tr}\ Q\right)^3
\Biggl(\frac{1}{6}-\sum_{j_1,j_2,j_3=0}^{q}
{\hat C}_{j_3j_2j_1}^2\Biggr).
\end{equation}

\vspace{4mm}

Let us estimate the right-hand sides of (\ref{ppp4}) and (\ref{sss1aaa}).
Using the ele\-men\-tary inequality $(a+b+c)^2\le 3(a^2+b^2+c^2)$, we obtain
          
\vspace{-3.5mm}
\begin{equation}
\label{ppp0}
~~~~~~{\sf M}\left\{\Biggl\Vert
I_2[B(Z)]_{T,t}-
I_2[B(Z)]_{T,t}^{M}\Biggr\Vert_H^2\right\}\le
3\left(E^{1,M}_{T,t}+E^{2,M}_{T,t}+E^{3,M}_{T,t}\right),
\end{equation}

\begin{equation}
\label{ppp0x}
~~~~~~{\sf M}\left\{\Biggl\Vert
I_3[B(Z)]_{T,t}-
I_3[B(Z)]_{T,t}^{M}\Biggr\Vert_H^2\right\}\le
3\left(G^{1,M}_{T,t}+G^{2,M}_{T,t}+G^{3,M}_{T,t}\right),
\end{equation}

\vspace{4.5mm}
\noindent 
where

\vspace{-4mm}
$$
E^{1,M}_{T,t}=
$$
$$
=
{\sf M}\left\{\left\Vert
\int\limits_{t}^{T}B'(Z) \left(
\int\limits_{t}^{t_3}B'(Z)\left(
\int\limits_{t}^{t_2}B(Z) d\left({\bf W}_{t_1}-{\bf W}_{t_1}^M\right) 
\right) d{\bf W}_{t_2} \right) d{\bf W}_{t_3}\right\Vert_H^2\right\},
$$

\vspace{1.5mm}
$$
E^{2,M}_{T,t}=
$$
$$
=
{\sf M}\left\{\left\Vert
\int\limits_{t}^{T}B'(Z) \left(
\int\limits_{t}^{t_3}B'(Z)\left(
\int\limits_{t}^{t_2}B(Z) d{\bf W}_{t_1}^M
\right) d\left({\bf W}_{t_2}-{\bf W}_{t_2}^M\right) \right) 
d{\bf W}_{t_3}\right\Vert_H^2\right\},
$$

\vspace{1.5mm}
$$
E^{3,M}_{T,t}=
$$
$$
=
{\sf M}\left\{\left\Vert
\int\limits_{t}^{T}B'(Z) \left(
\int\limits_{t}^{t_3}B'(Z)\left(
\int\limits_{t}^{t_2}B(Z) d{\bf W}_{t_1}^M
\right) d{\bf W}_{t_2}^M\right) 
d\left({\bf W}_{t_3}-{\bf W}_{t_3}^M\right)\right\Vert_H^2\right\},
$$

\newpage
\noindent
$$
G^{1,M}_{T,t}=
{\sf M}\left\{\left\Vert
\int\limits_{t}^{T}B''(Z) \left(
\int\limits_{t}^{t_2}B(Z)d{\bf W}_{t_1},
\int\limits_{t}^{t_2}B(Z)d\left({\bf W}_{t_1}-{\bf W}_{t_1}^M\right)
\right) 
d{\bf W}_{t_2}\right\Vert_H^2\right\},
$$

\vspace{1mm}
$$
G^{2,M}_{T,t}=
{\sf M}\left\{\left\Vert
\int\limits_{t}^{T}B''(Z) \left(
\int\limits_{t}^{t_2}B(Z)d\left({\bf W}_{t_1}-{\bf W}_{t_1}^M\right),
\int\limits_{t}^{t_2}B(Z)d{\bf W}_{t_1}^M
\right) 
d{\bf W}_{t_2}\right\Vert_H^2\right\},
$$

$$
G^{3,M}_{T,t}=
{\sf M}\left\{\left\Vert
\int\limits_{t}^{T}B''(Z) \left(
\int\limits_{t}^{t_2}B(Z)d{\bf W}_{t_1}^M,
\int\limits_{t}^{t_2}B(Z)d{\bf W}_{t_1}^M
\right) 
d\left({\bf W}_{t_2}-{\bf W}_{t_2}^M\right)\right\Vert_H^2\right\}.
$$

\vspace{6mm}

We have

\vspace{-5mm}
$$
E^{1,M}_{T,t}=
$$

\vspace{-4mm}
$$
=
\int\limits_{t}^{T}\hspace{-0.5mm}{\sf M}\left\{\left\Vert
B'(Z) \hspace{-1mm}\left(
\int\limits_{t}^{t_3} \hspace{-1mm} B'(Z) \hspace{-1mm}\left(
\int\limits_{t}^{t_2}B(Z) d\left({\bf W}_{t_1}-{\bf W}_{t_1}^M\right) 
\hspace{-1mm} \right) \hspace{-0.5mm}d{\bf W}_{t_2} \hspace{-1mm}
\right)\right\Vert_{L_{HS}(U_0,H)}^2
\hspace{-0.2mm}\right\}\hspace{-0.2mm}dt_3\hspace{-0.2mm}\le
$$

\vspace{-1mm}
$$
\le C_5 
\int\limits_{t}^{T}{\sf M}\left\{\left\Vert
\int\limits_{t}^{t_3}B'(Z)\left(
\int\limits_{t}^{t_2}B(Z) d\left({\bf W}_{t_1}-{\bf W}_{t_1}^M\right) 
\right) d{\bf W}_{t_2}\right\Vert_H^2
\right\}dt_3=
$$

\vspace{-1mm}
$$
=
C_5 
\int\limits_{t}^{T}\int\limits_{t}^{t_3}{\sf M}\left\{\left\Vert
B'(Z)\left(
\int\limits_{t}^{t_2}B(Z) d\left({\bf W}_{t_1}-{\bf W}_{t_1}^M\right) 
\right)\right\Vert_{L_{HS}(U_0,H)}^2 \right\}dt_2
dt_3\le
$$

\vspace{-1mm}
\begin{equation}
\label{22u}
~~~~~~~\le 
C_6 
\int\limits_{t}^{T}\int\limits_{t}^{t_3}{\sf M}\left\{\left\Vert
\int\limits_{t}^{t_2}B(Z) d\left({\bf W}_{t_1}-{\bf W}_{t_1}^M\right) 
\right\Vert_H^2 \right\}dt_2
dt_3\le
\end{equation}

\vspace{-1mm}
\begin{equation}
\label{23u}
~~\le 
C_6 \left(\sup\limits_{i\in J\backslash J_M}\lambda_i\right)^{2\alpha}
\int\limits_{t}^{T}\int\limits_{t}^{t_3}
\int\limits_{t}^{t_2}
{\sf M}\left\{\biggl\Vert
B(Z)Q^{-\alpha}\biggr\Vert_{L_{HS}(U_0,H)}^2 \right\}dt_1 dt_2
dt_3 \le 
\end{equation}

\newpage
\noindent
\begin{equation}
\label{ppp1aaa}
\le C_7 
\left(\sup\limits_{i\in J\backslash J_M}\lambda_i\right)^{2\alpha}
(T-t)^3.
\end{equation}

\vspace{3mm}

Note that the transition from (\ref{22u}) to (\ref{23u}) was made
by analogy with the proof of Theorem 1 in \cite{27zzzzz} (also see 
\cite{12zzzzz}). 
More precisely, taking into account the relation $Q^{\alpha}e_i=
\lambda_i^{\alpha}e_i$, we have 
(see \cite{27zzzzz}, Sect.~3.1)

\vspace{-1mm}
$$
{\sf M}\left\{\left\Vert
\int\limits_{t}^{t_2}B(Z) d\left({\bf W}_{t_1}-{\bf W}_{t_1}^M\right) 
\right\Vert_H^2\right\}=
$$

$$
=
{\sf M}\left\{\left\Vert \sum\limits_{i\in J\backslash J_M}
\sqrt{\lambda_i}
\int\limits_{t}^{t_2}B(Z)e_i d{\bf w}_{t_1}^{(i)} 
\right\Vert_H^2\right\}=
$$

$$
=\sum\limits_{i\in J\backslash J_M}\lambda_i
\int\limits_{t}^{t_2}{\sf M}\left\{\biggl\Vert
B(Z)Q^{-\alpha}Q^{\alpha}e_i 
\biggr\Vert_H^2\right\}dt_1=
$$

$$
=\sum\limits_{i\in J\backslash J_M}\lambda_i^{1+2\alpha}
\int\limits_{t}^{t_2}{\sf M}\left\{\biggl\Vert
B(Z)Q^{-\alpha}e_i 
\biggr\Vert_H^2\right\}dt_1=
$$

$$
=\left(\sup\limits_{i\in J\backslash J_M}\lambda_i\right)^{2\alpha}
\int\limits_{t}^{t_2}{\sf M}\left\{
\sum\limits_{i\in J\backslash J_M}\lambda_i
\biggl\Vert
B(Z)Q^{-\alpha}e_i 
\biggr\Vert_H^2\right\}dt_1\le
$$

$$
\le 
\left(\sup\limits_{i\in J\backslash J_M}\lambda_i\right)^{2\alpha}
\int\limits_{t}^{t_2}{\sf M}\left\{
\sum\limits_{i\in J}\lambda_i
\biggl\Vert
B(Z)Q^{-\alpha}e_i 
\biggr\Vert_H^2\right\}dt_1=
$$

\begin{equation}
\label{24u}
~~~~~~~~ =\left(\sup\limits_{i\in J\backslash J_M}\lambda_i\right)^{2\alpha}
\int\limits_{t}^{t_2}{\sf M}\left\{
\biggl\Vert
B(Z)Q^{-\alpha} 
\biggr\Vert_{L_{HS}(U_0,H)}^2\right\}dt_1.
\end{equation}

\vspace{5mm}

Further,
we also will use the estimate similar to (\ref{24u}).

We have

\vspace{-4mm}
$$
E^{2,M}_{T,t}=
$$

\newpage
\noindent
$$
=
\int\limits_{t}^{T}\hspace{-1mm}{\sf M}\left\{\left\Vert
B'(Z) \hspace{-1mm}\left(
\int\limits_{t}^{t_3} \hspace{-1mm} B'(Z) \hspace{-1mm}\left(
\int\limits_{t}^{t_2} \hspace{-1mm} B(Z) d{\bf W}_{t_1}^M
\right) \hspace{-0.3mm}
d\left({\bf W}_{t_2}-{\bf W}_{t_2}^M \right) \hspace{-1mm}
\right)\right\Vert_{L_{HS}(U_0,H)}^2
\right\}\hspace{-1mm}dt_3\hspace{-0.3mm}\le
$$

$$
\le C_8
\int\limits_{t}^{T}{\sf M}\left\{\left\Vert
\int\limits_{t}^{t_3}B'(Z)\left(
\int\limits_{t}^{t_2}B(Z) d{\bf W}_{t_1}^M
\right) d\left({\bf W}_{t_2}-{\bf W}_{t_2}^M\right) 
\right\Vert_H^2
\right\}dt_3\le
$$

$$
\le C_8 \hspace{-0.3mm}\left(\sup\limits_{i\in J\backslash J_M}
\lambda_i\right)^{\hspace{-1mm}2\alpha}
\int\limits_{t}^{T}\int\limits_{t}^{t_3}{\sf M}\left\{\left\Vert
B'(Z)\hspace{-1mm}\left(
\int\limits_{t}^{t_2}\hspace{-1mm}B(Z) d{\bf W}_{t_1}^M
\right) \hspace{-1mm} Q^{-\alpha}
\right\Vert_{L_{HS}(U_0,H)}^2
\right\}dt_2 dt_3\le  
$$

$$
\le C_{8} \left(\sup\limits_{i\in J\backslash J_M}\lambda_i\right)^{2\alpha}
\int\limits_{t}^{T}\int\limits_{t}^{t_3}
{\sf M}\left\{\biggl\Vert B'(Z)\left(B(Z)\right)Q^{-\alpha}
\biggr\Vert^2_{L_{HS}^{(2)}(U_0,H)}\right\}(t_2-t)dt_2
dt_3 \le
$$

\vspace{-1mm}
\begin{equation}
\label{ppp2}
\le C_9 
\left(\sup\limits_{i\in J\backslash J_M}\lambda_i\right)^{2\alpha}
(T-t)^3.
\end{equation}

\vspace{6mm}

Moreover,

\vspace{-3mm}
$$
E^{3,M}_{T,t}\le
\left(\sup\limits_{i\in J\backslash J_M}\lambda_i\right)^{2\alpha}\times
$$

\vspace{-1mm}
$$
\times
\int\limits_{t}^{T}{\sf M}\left\{\left\Vert
B'(Z) \left(
\int\limits_{t}^{t_3}B'(Z)\left(
\int\limits_{t}^{t_2}B(Z) d{\bf W}_{t_1}^M
\right) d{\bf W}_{t_2}^M\right) 
Q^{-\alpha}\right\Vert_{L_{HS}(U_0,H)}^2\right\}dt_3\le
$$

\vspace{-1mm}
$$
\le C_{10}
\left(\sup\limits_{i\in J\backslash J_M}\lambda_i\right)^{2\alpha}\times
$$

\vspace{-1mm}
$$
\times\int\limits_{t}^{T}
{\sf M}\left\{\biggl\Vert B'(Z)\left(B'(Z)\left(B(Z)\right)
\right)Q^{-\alpha}
\biggr\Vert^2_{L_{HS}^{(3)}(U_0,H)}\right\}
\frac{(t_3-t)^2}{2}dt_3
\le
$$

\newpage
\noindent
\begin{equation}
\label{ppp3}
\le C_{11}\left(\sup\limits_{i\in J\backslash J_M}\lambda_i\right)^{2\alpha}
(T-t)^3.
\end{equation}

\vspace{4mm}

Combining (\ref{ppp4}), (\ref{ppp0}), (\ref{ppp1aaa}), 
(\ref{ppp2}), and (\ref{ppp3}),
we obtain (\ref{ppp6}).

We have

\vspace{-4mm}
$$
G^{1,M}_{T,t}=
$$
$$
=\int\limits_{t}^{T}
{\sf M}\left\{\left\Vert
B''(Z) \left(
\int\limits_{t}^{t_2}B(Z)d{\bf W}_{t_1},
\int\limits_{t}^{t_2}B(Z)d\left({\bf W}_{t_1}-{\bf W}_{t_1}^M\right)
\right)\right\Vert_{L_{HS}(U_0,H)}^2\right\}dt_3\le
$$

\vspace{-1mm}
$$
\le C_{12}\int\limits_{t}^{T}
{\sf M}\left\{\left\Vert
\int\limits_{t}^{t_2}B(Z)d{\bf W}_{t_1}
\right\Vert_H^2
\left\Vert
\int\limits_{t}^{t_2}B(Z)d\left({\bf W}_{t_1}-{\bf W}_{t_1}^M\right)
\right\Vert_H^2\right\}dt_3\le
$$

$$
\hspace{-40mm}\le C_{12}\int\limits_{t}^{T}
\left({\sf M}\left\{\left\Vert
\int\limits_{t}^{t_2}B(Z)d{\bf W}_{t_1}
\right\Vert_H^4\right\}\right)^{1/2}\times
$$

\vspace{-6mm}
$$
~~~~~~~~~~~~~~~~~~~~~~~\times
\left({\sf M}\left\{\left\Vert
\int\limits_{t}^{t_2}B(Z)d\left({\bf W}_{t_1}-{\bf W}_{t_1}^M\right)
\right\Vert_H^4\right\}\right)^{1/2}dt_3\le
$$
$$
\hspace{-40mm}\le C_{13}\int\limits_{t}^{T}
\int\limits_{t}^{t_2}\left({\sf M}\left\{\biggl\Vert B(Z)
\biggr\Vert_{L_{HS}(U_0,H)}^4\right\}\right)^{1/2}
dt_1\times
$$
$$
~~~~~~~~~~~~~~~~~~~~~~~\times\left({\sf M}\left\{\left\Vert
\int\limits_{t}^{t_2}B(Z)d\left({\bf W}_{t_1}-{\bf W}_{t_1}^M\right)
\right\Vert_H^4\right\}\right)^{1/2}dt_3\le
$$
\begin{equation}
\label{lll1}
~~~~~~~~~~~\le C_{14}\int\limits_{t}^{T}
(t_2-t)
\left({\sf M}\left\{\left\Vert
\int\limits_{t}^{t_2}B(Z)d\left({\bf W}_{t_1}-{\bf W}_{t_1}^M\right)
\right\Vert_H^4\right\}\right)^{1/2}dt_3.
\end{equation}

\vspace{4mm}

Let us estimate the right-hand side of (\ref{lll1}).
If $s>t,$ then for fixed $M\in{\bf N}$ and 
for some $N>M$ ($N\in{\bf N}$) we have

\newpage
\noindent
$$
{\sf M}\left\{\left\Vert
\int\limits_{t}^{s}B(Z)d\left({\bf W}_{t_1}^N-{\bf W}_{t_1}^M\right)
\right\Vert_H^4\right\}=
$$
$$
\hspace{-40mm}={\sf M}\left\{
\Biggl\langle 
\sum\limits_{j\in J_N\backslash J_M}\sqrt{\lambda_j}
B(Z)e_j\left({\bf w}_{s}^{(j)}-{\bf w}_{t}^{(j)}\right),\Biggr.\right.
$$

\vspace{-9mm}$$                     
\left.\Biggl.~~~~~~~~~~~~~~~~~~~~~~~~~\sum\limits_{j'\in J_N\backslash J_M}
\sqrt{\lambda_{j'}}
B(Z)e_{j'}\left({\bf w}_{s}^{(j')}-{\bf w}_{t}^{(j')}\right)
\Biggr\rangle_H^2\right\}=
$$
$$
=
\sum\limits_{j, j', l, l'\in J_N\backslash J_M}
\hspace{-2mm}\sqrt{\lambda_j\lambda_{j'}\lambda_l\lambda_{l'}}\ \
{\sf M}\Biggl\{
\biggl\langle 
B(Z)e_j,
B(Z)e_{j'}
\biggr\rangle_H\biggl\langle 
B(Z)e_l,
B(Z)e_{l'}
\biggr\rangle_H\times\Biggr.
$$

\vspace{-2mm}
$$
\Biggl.\times
{\sf M}\left\{\left({\bf w}_{s}^{(j)}-{\bf w}_{t}^{(j)}\right)
\left({\bf w}_{s}^{(j')}-{\bf w}_{t}^{(j')}\right)
\left({\bf w}_{s}^{(l)}-{\bf w}_{t}^{(l)}\right)
\left({\bf w}_{s}^{(l')}-{\bf w}_{t}^{(l')}\right)
\biggl.\biggr|{\bf F}_t\right\}\Biggr\}
=
$$

$$
=3(s-t)^2\sum\limits_{j\in J_N\backslash J_M}
\lambda_j^{2}
{\sf M}\left\{\biggl\Vert B(Z)e_j\biggr\Vert_H^4\right\}+
$$

\vspace{3mm}
$$
\hspace{-30mm}+(s-t)^2
\sum\limits_{j, j'\in J_N\backslash J_M (j\ne j')}
\lambda_j\lambda_{j'}\Biggl(
{\sf M}\left\{\biggl\Vert B(Z)e_j\biggr\Vert_H^2
\biggl\Vert B(Z)e_{j'}\biggr\Vert_H^2\right\}+\Biggr.
$$

\vspace{-6mm}
$$
~~~~~~~~~~~~~~~~~~~~~~~~~~~~~~~~~~~~~~~~~~~~~~~~~~~~~~~~~~~~\Biggl.+
2\biggl\langle 
B(Z)e_j,
B(Z)e_{j'}
\biggr\rangle_H^2\Biggr)\le
$$

\vspace{-6mm}
$$
\hspace{-70mm}\le 3(s-t)^2\left(\sum\limits_{j\in J_N\backslash J_M}
\lambda_j^{2}
{\sf M}\left\{\biggl\Vert B(Z)e_j\biggr\Vert_H^4\right\}+\right.
$$
$$
~~~~~~~~~~~~~~~~~~~~~~~~~~\left.+
\sum\limits_{j, j'\in J_N\backslash J_M (j\ne j')}
\lambda_j\lambda_{j'}
{\sf M}\left\{\biggl\Vert B(Z)e_j\biggr\Vert_H^2
\biggl\Vert B(Z)e_{j'}\biggr\Vert_H^2\right\}\right)=
$$

$$
=
3(s-t)^2
{\sf M}\left\{\left(
\sum\limits_{j\in J_N\backslash J_M}\lambda_j
\biggl\Vert B(Z)e_j\biggr\Vert_H^2\right)^2\right\}\le
$$
$$
\le
3(s-t)^2\left(\sup\limits_{i\in J_N\backslash J_M}\lambda_i\right)^{4\alpha}
{\sf M}\left\{\left(
\sum\limits_{j\in J_N\backslash J_M}\lambda_j
\biggl\Vert B(Z)Q^{-\alpha}e_j\biggr\Vert_H^2\right)^2\right\}\le
$$

\newpage
\noindent
\begin{equation}
\label{hhhh1}
~~~~~~~~\le
C_{15}(s-t)^2
\left(\sup\limits_{i\in J_N\backslash J_M}\lambda_i\right)^{4\alpha}
{\sf M}\left\{
\biggl\Vert B(Z)Q^{-\alpha}\biggr\Vert_{L_{HS}(U_0,H)}^4\right\}.
\end{equation}

\vspace{5mm}

Using (\ref{hhhh1}), we obtain
$$
{\sf M}\left\{\left\Vert
\int\limits_{t}^{s}B(Z)d\left({\bf W}_{t_1}^N-{\bf W}_{t_1}^M\right)
\right\Vert_H^4\right\}\to 0
$$

\noindent
if $N, M\to \infty\  (N>M)$. This means that
$$
\int\limits_{t}^{s}B(Z)d{\bf W}_{t_1}^N
$$
is a Cauchy sequence in the $L_4$--space of $H$--valued 
stochastic processes.

It is well known \cite{Pis} that $L_p$--space ($1\le p <\infty $) of Banach space valued 
stochastic processes is a Banach space, i.e. a complete space. Then, carrying out the passage to the 
limit $\lim\limits_{N\to\infty}$ in (\ref{hhhh1}), we get
$$
{\sf M}\left\{\left\Vert
\int\limits_{t}^{s}B(Z)d\left({\bf W}_{t_1}-{\bf W}_{t_1}^M\right)
\right\Vert_H^4\right\}
=
$$
$$
=
\lim\limits_{N\to\infty}{\sf M}\left\{\left\Vert
\int\limits_{t}^{s}B(Z)d\left({\bf W}_{t_1}^N-{\bf W}_{t_1}^M\right)
\right\Vert_H^4\right\}\le
$$
\begin{equation}
\label{hhhh1xxx}
~~~~~~~~\le
C_{15}(s-t)^2
\left(\sup\limits_{i\in J\backslash J_M}\lambda_i\right)^{4\alpha}
{\sf M}\left\{
\biggl\Vert B(Z)Q^{-\alpha}\biggr\Vert_{L_{HS}(U_0,H)}^4\right\}.
\end{equation}

\vspace{4mm}

Combining 
(\ref{lll1}) and (\ref{hhhh1xxx}), we obtain
\begin{equation}
\label{aaa1ddd}
G^{1,M}_{T,t}\le C_{16} 
\left(\sup\limits_{i\in J\backslash J_M}\lambda_i\right)^{2\alpha}
(T-t)^3.
\end{equation}

\vspace{2mm}

Absolutely analogously we get
\begin{equation}
\label{plan1000}
G^{2,M}_{T,t}\le C_{17} 
\left(\sup\limits_{i\in J\backslash J_M}\lambda_i\right)^{2\alpha}
(T-t)^3.
\end{equation}

\vspace{4mm}

Let us estimate $G^{3,M}_{T,t}.$ We have

$$
G^{3,M}_{T,t}\le
\left(\sup\limits_{i\in J\backslash J_M}\lambda_i\right)^{2\alpha}\times
$$
$$
\hspace{-10mm}\times
\int\limits_{t}^{T}
{\sf M}\left\{\left\Vert
B''(Z) \left(
\int\limits_{t}^{t_2}B(Z)d{\bf W}_{t_1}^M,
\int\limits_{t}^{t_2}B(Z)d{\bf W}_{t_1}^M
\right)Q^{-\alpha} 
\right\Vert_{L_{HS}(U_0,H)}^2\right\}dt_2\le
$$
$$
\hspace{-50mm}\le 
\left(\sup\limits_{i\in J\backslash J_M}\lambda_i\right)^{2\alpha}
\sum\limits_{i\in J}
\sum\limits_{j,l\in J_M}\lambda_i \lambda_j
\lambda_l\times
$$
$$
~~~~~~~~~~~~~~~~~~~~~~~~\times
\int\limits_t^T(t_2-t)^2
\left({\sf M}\Biggl\{\biggl\Vert B''(Z)(B(Z)e_j,B(Z)e_l)Q^{-\alpha}e_i
\biggr\Vert_H^2\right\}+\Biggr.
$$
$$
\hspace{-3mm}
+{\sf M}\Biggl\{\biggl\Vert B''(Z)(B(Z)e_j,B(Z)e_j)Q^{-\alpha}e_i
\biggr\Vert_H
\biggl\Vert B''(Z)(B(Z)e_l,B(Z)e_l)Q^{-\alpha}e_i
\biggr\Vert_H
\Biggr\}
+
$$
$$
\hspace{-25mm}+
{\sf M}\Biggl\{\biggl\Vert B''(Z)(B(Z)e_j,B(Z)e_l)Q^{-\alpha}e_i
\biggr\Vert_H\times\Biggr.
$$
$$
~~~~~~~~~~~~~~~~~~~~~\Biggl.\Biggl.\times
\biggl\Vert B''(Z)(B(Z)e_l,B(Z)e_j)Q^{-\alpha}e_i
\biggr\Vert_H
\Biggr\}\Biggr)dt_2\le
$$
\begin{equation}
\label{kkk1www}
\le C_{18}
\left(\sup\limits_{i\in J\backslash J_M}\lambda_i\right)^{2\alpha}
(T-t)^3.
\end{equation}

\vspace{4mm}

Combining (\ref{sss1aaa}), (\ref{ppp0x}), 
and (\ref{aaa1ddd})--(\ref{kkk1www}),
we obtain (\ref{uuu1}). Theorem 7.3 is proved.

Let us consider the convergence analysis
for the stochastic integrals
(\ref{ret2})--(\ref{ret4}) (convergence for the 
stochastic integral (\ref{ret1}) follows from 
(\ref{24u}) (see Theorem 1 in \cite{27zzzzz} or \cite{12zzzzz})). 

Using the It\^{o} formula, we obtain w.~p.~1 \cite{13zzzzz}

\vspace{-2mm}
$$
J_2[B(Z)]_{T,t}=\int\limits_{t}^{T}\left(\frac{T}{2}-s+\frac{t}{2}\right)
AB(Z)d{\bf W}_{s},
$$
$$
J_3[B(Z),F(Z)]_{T,t}
=\int\limits_{t}^{T}(s-t)B'(Z)\biggl(AZ+F(Z)\biggr)d{\bf W}_{s}.
$$

\vspace{2mm}

Suppose that

$$
{\sf M}\left\{\biggl\Vert B'(Z)\biggl(AZ+F(Z)\biggr)Q^{-\alpha}
\biggr\Vert_{L_{HS}(U_0,H)}^2\right\}<\infty,
$$

\vspace{1mm}
$$
{\sf M}\left\{\biggl\Vert AB(Z)Q^{-\alpha}\biggr\Vert_{L_{HS}(U_0,H)}^2\right\}
<\infty
$$

\vspace{2mm}
\noindent
for some $\alpha>0.$

Then by analogy with (\ref{24u}) we get

\vspace{-3mm}
$$
{\sf M}\left\{\Biggl\Vert
J_2[B(Z)]_{T,t}-
J_2[B(Z)]_{T,t}^{M}\Biggr\Vert_H^2\right\}\le
$$
$$
\le C_{19} (T-t)^3
\left(\sup\limits_{i\in J\backslash J_M}\lambda_i\right)^{2\alpha},
$$

$$
{\sf M}\left\{\Biggl\Vert
J_3[B(Z),F(Z)]_{T,t}-
J_3[B(Z),F(Z)]_{T,t}^{M}\Biggr\Vert_H^2\right\}\le
$$
$$
\le C_{20} (T-t)^3
\left(\sup\limits_{i\in J\backslash J_M}\lambda_i\right)^{2\alpha},
$$

\vspace{5mm}
\noindent
where
$J_2[B(Z)]_{T,t}^{M},$
$J_3[B(Z),F(Z)]_{T,t}^{M}$
are defined by 
(\ref{sd11}), (\ref{sd12}).

Moreover, under the conditions of Theorem 7.3 we obtain for some $\alpha>0$

\vspace{-1mm}
$$
{\sf M}\left\{\Biggl\Vert
J_4[B(Z),F(Z)]_{T,t}-
J_4[B(Z),F(Z)]_{T,t}^{M}\Biggr\Vert_H^2\right\}=
$$
$$
={\sf M}\left\{\left\Vert
\int\limits_t^T F'(Z)\left(\int\limits_t^{t_2}
B(Z)
d\left({\bf W}_{t_1}-{\bf W}_{t_1}^M\right)\right)
dt_2\right\Vert_H^2\right\}\le
$$
$$
\le 
(T-t)
\int\limits_t^T
{\sf M}\left\{\left\Vert
F'(Z)\left(\int\limits_t^{t_2}
B(Z)
d\left({\bf W}_{t_1}-{\bf W}_{t_1}^M\right)\right)
\right\Vert_H^2\right\}dt_2\le
$$
$$
\le 
C_{21}(T-t)
\int\limits_t^T
{\sf M}\left\{\left\Vert
\int\limits_t^{t_2}
B(Z)
d\left({\bf W}_{t_1}-{\bf W}_{t_1}^M\right)
\right\Vert_H^2\right\}dt_2\le
$$

\vspace{-3mm}
$$
\le C_{21}(T-t)
\left(\sup\limits_{i\in J\backslash J_M}\lambda_i\right)^{2\alpha}
\int\limits_t^T
\int\limits_t^{t_2}
{\sf M}\left\{\biggl\Vert
B(Z)Q^{-\alpha}
\biggr\Vert_{L_{HS}(U_0,H)}^2\right\} dt_1dt_2\le
$$
$$
\le 
C_{22}(T-t)^3
\left(\sup\limits_{i\in J\backslash J_M}\lambda_i\right)^{2\alpha},
$$

\vspace{5mm}
\noindent
where 
$J_4[B(Z),F(Z)]_{T,t}^{M}$ is defined by 
(\ref{sd13}).

\newpage

\addcontentsline{toc}{chapter}{Epilogue}

\newpage

\vspace{30mm}
$$
{}
$$
$$
{}
$$

\noindent
{\Huge {\bf Epilogue}}

\vspace{10mm}

The results presented in this book were developed \cite{Kuz-Kuz}, \cite{Mikh-1}
(also see \cite{123789000}, Chapter~16)
in the form
of a software package in the Python programming language
that implements the numerical methods (\ref{al1})--(\ref{al5}), 
(\ref{al1x})--(\ref{al5x}) (see Chapter 4) 
with the orders $1.0,$ $1.5,$ $2.0,$ $2.5$, and $3.0$
of strong convergence
based on the unified Taylor--It\^{o} and Taylor--Stratonovich 
expansions. At that for the numerical simulation of iterated 
It\^{o} and Stratonovich stochastic integrals of multiplicities 
1 to 6 we used \cite{Kuz-Kuz}, \cite{Mikh-1} the formulas from Sect.~5.1, i.e. method
based on Theorem 1.1 and multiple Fourier--Legendre series.
Note that in \cite{Kuz-Kuz}, \cite{Mikh-1} we used
the database with 270,000 exactly calculated
Fourier--Legendre coefficients.

Using computational experiments
it was shown in \cite{Mikh-2}, \cite{Mikh-2aaaaa} (also see Sect.~5.4) that
in most cases all the exact formulas from Sect.~1.2.3
for the mean-square approximation errors of iterated It\^{o}
stochastic integarls can be replaced by the formula
(\ref{pochti}) for $k=1,\ldots,5$. This allows us to neglect the 
multiplier factor $k!$ (see the formula (\ref{z1})).
As a result, the computational costs for the approximation 
of iterated It\^{o} stochastic integrals are significantly reduced.
For the same reason, we can replace the multiplier factor
$(k!)^2$ by $k!$ in the formula (\ref{zzz1ee}) in practical calculations.

Iterated stochastic integrals are a fundamental tool for describing and 
studying the dynamics of various types of stochastic equations.
In recent years and decades, numerical methods 
of high orders of accuracy have been 
constructed using iterated stochastic integrals
not only for It\^{o} SDEs, 
but also for SDEs with jumps \cite{Zapad-9},
SPDEs with multiplicative trace class noise
\cite{12zzzzz}, \cite{13zzzzz}, \cite{18zzzzz},
McKean SDEs \cite{koms1}, 
SDEs with switchings
\cite{koms2}, mean-field SDEs \cite{koms3}, 
It\^{o}--Volterra SDEs \cite{18zzzzz}, etc.

\newpage

\addcontentsline{toc}{chapter}{Bibliography}

{

}

\end{document}